\errorstopmode
\documentclass[10pt,CMpaper]{CMbook}

\usepackage{graphicx,color}
\usepackage{eurosym}
\usepackage{%amsmath,
amsfonts,amssymb}
\usepackage{mathrsfs}
\usepackage{oldstyle}
\usepackage[T1]{fontenc}
\usepackage[utf8]{inputenc}
\usepackage[french]{babel}
\usepackage[all]{xy}
\usepackage{pb-diagram}
\usepackage[french,nohints]{minitoc}
\usepackage{multicol}
\usepackage{CMmacros2}
\usepackage{makeidx}\makeindex\makeglossary
\usepackage{makehints}\makehints

\usepackage{lmodern}
\usepackage{microtype}
\usepackage{array}
\usepackage{stmaryrd}
\usepackage{float}
\usepackage{wasysym}

\usepackage{bm}

\DeclareFontFamily{U}{wasy}{}
\DeclareFontShape{U}{wasy}{m}{n}{<5> <6> <7> <8> <9> gen * wasy
      <10> <10.95> <12> <14.4> <17.28> <20.74> <24.88>wasy10  }{}
\DeclareFontShape{U}{wasy}{b}{n}{ <-10> ssub * wasy/m/n
 <10> <10.95> <12> <14.4> <17.28> <20.74> <24.88>wasyb10 }{}
\DeclareFontShape{U}{wasy}{bx}{n}{<5> <6> <7> <8> <9> gen * wasy
 <10> <10.95> <12> <14.4> <17.28> <20.74> <24.88>wasyb10}{}
\DeclareSymbolFont{wasy}{U}{wasy}{m}{n}
\SetSymbolFont{wasy}{bold}{U}{wasy}{bx}{n}
\DeclareFontFamily{U}{lasy}{}
\DeclareFontShape{U}{lasy}{m}{n}{ <5> <6> <7> <8> <9> gen * lasy
      <10> <10.95> <12> <14.4> <17.28> <20.74> <24.88>lasy10  }{}
\DeclareFontShape{U}{lasy}{b}{n}{ <5> <6> <7> <8> <9> gen * lasy
      <10> <10.95> <12> <14.4> <17.28> <20.74> <24.88>lasyb10  }{}
\DeclareFontFamily{U}{stmry}{}
\DeclareFontShape{U}{stmry}{m}{n}
   {  <5> <6> <7> <8> <9> <10> gen * stmary
      <10.95><12><14.4><17.28><20.74><24.88>stmary10%
   }{}
\DeclareFontShape{U}{stmry}{b}{n}
   {  <5> <6> <7> <8> <9> <10> gen * stmary
      <10.95><12><14.4><17.28><20.74><24.88>stmary10%
   }{}

\ifpdf
\makeatletter
\def\toclevel@chapternb{0}
\makeatother
\usepackage[bookmarksopen=false,pdftex=true,breaklinks=true,backref=page,%
    pagebackref=true,plainpages=false,%
    hyperindex=true,pdfstartview=FitH,
    pdfpagelabels=true,colorlinks=false,citecolor=red,colorlinks=false]%
    {hyperref}
\else
\newcommand{\texorpdfstring}[2]{#1}
\fi

\newcommand\sinotenglish[1]{#1}

\title{Alg\`ebre commutative.  M\'ethodes constructives.
2\ieme \'edition.
}
\author{Henri Lombardi \&\
Claude Quitté}
\date{\today }
\begin{document}
\overfullrule=0cm
\hfuzz14pt

\let\oldref\ref
\renewcommand{\ref}[1]{\hbox{\oldref{#1}}}

\makeatletter
\def\U@@@ref#1\relax#2\relax{\let\showsection\relax#2}
\newcommand{\U@@ref}[5]{\U@@@ref#1}
\newcommand{\U@ref}[1]{\NR@setref{#1}\U@@ref{#1}}
\newcommand{\@iref}{\@ifstar\@refstar\U@ref}
\newcommand{\iref}[1]{\hbox{\@iref{#1}}}
\def\l@chapter#1#2{\ifnum\c@tocdepth>\m@ne\addpenalty{-\@highpenalty}\vskip1.0em\@plus\p@
\setlength\@tempdima{2.8em}\begingroup\parindent\z@\rightskip\@pnumwidth\parfillskip-\@pnumwidth\penalty-2000\leavevmode\bfseries
\advance\leftskip\@tempdima\hskip-\leftskip\hskip0.8em{\boldmath#1}\nobreak\hfil\nobreak\hb@xt@\@pnumwidth{\hss}\par\penalty\@highpenalty\endgroup\fi}
\def\l@chapterbis#1#2{\ifnum\c@tocdepth>\m@ne\addpenalty{-\@highpenalty}\vskip1.0em\@plus\p@
\setlength\@tempdima{2em}\begingroup\parindent\z@\rightskip\@pnumwidth\parfillskip-\@pnumwidth\penalty-2000\leavevmode\bfseries
\advance\leftskip\@tempdima\hskip-\leftskip{\boldmath#1}\nobreak\hfil\nobreak\hb@xt@\@pnumwidth{\hss#2}\par\penalty\@highpenalty\endgroup\fi}
\def\toclevel@chapterbis{0}
\def\l@section{\@dottedtocline{1}{1em}{0em}}
\def\l@subsection{\@dottedtocline{2}{1.75em}{0.25em}}
\def\l@subsubsection{\@dottedtocline{3}{2.5em}{0.25em}}
\def\numberline#1{\hb@xt@\@tempdima{\hss#1\ \hfil}}
\def\contentsline#1#2#3#4{\ifx\\#4\\\csname
l@#1\endcsname{#2}{#3}\else\ifHy@linktocpage\csname
l@#1\endcsname{{#2}}{\hyper@linkstart{link}{#4}{#3}\hyper@linkend }
\else\csname l@#1\endcsname{\hyper@linkstart{link}{#4}{#2}\hyper@linkend
}{#3}\fi\fi}
\makeatother

\setcounter{minitocdepth}{3}
\dominitoc
%\pageblanche0
%\maketitle
%%%%%%%%%%%%%%%%%%%%%%%%%%%%%%
%!TEX root =  ACMC-1-copie.tex
%!TEX encoding =  UTF-8 Unicode

\ifx\mesMacrosDejaChargees\undefined\let\chargeMesMacros\relax
\else\let\chargeMesMacros \fi
\chargeMesMacros

\gdef\mesMacrosDejaChargees{}

%%%%%%%%%%%%%%%%%%%%%%%%%%%%%%%%%%%%%%%%%%%%%%%%%%%%%%%%%%%%%

\def\smalltimes{\raise0.75pt\hbox{\boldmath$\mathsurround0pt\scriptstyle\times$}}
\def\ssmalltimes{\raise0.5pt\hbox{$\mathsurround0pt\scriptscriptstyle\times$}}

\FrenchFootnotes

\newcommand{\DBxk}{{\Der \gk\gB\xi}}%
\newcommand{\DAxk}{{\Der \gk\gA\xi}}%
\newcommand{\DkXxk}{{\Der \gk\kuX\xi}}%
\newcommand \STR {\MA{\rm SR}}
\newcommand \Str {\MA{\mathsf{Sr}}}
\newcommand \Bdim {\MA{\mathsf{Bdim}}}
\newcommand \Sdim {\MA{\mathsf{Sdim}}}
\newcommand \Cdim {\MA{\mathsf{Cdim}}}
\newcommand \Gdim {\MA{\mathsf{Gdim}}}
\newcommand \DAbul {\rD_{\!\Abul}}

\newcommand \dotdiv {\; {{\,.\,} \over {} } \;}

\renewcommand{\labelenumii}{\theenumii.}

\newcommand\perso[1]{\sibrouillon{\marginpar{\hspace*{1em}
\begin{minipage}{10em}{\begin{flushleft}\footnotesize #1

\end{flushleft}}
\end{minipage}}}}

\newcommand \entrenous[1]{\sibrouillon{{\sf Entrenous. #1 \hfill \hbox{FinEntrenous}\par}}}
\newcommand \hum[1]{\sibrouillon{{\par\small\tt Hum: #1 \eoe\par}}}

\newcommand \incertain[1]{}
\newcommand \Oui[1]{}

\newcounter{bidon}
\newcommand{\rdb}{\refstepcounter{bidon}}

%:  index
\newcommand \ix[1] {\index{#1}\emph{#1}}
\newcommand \ixc[2] {\index{#1!#2}\emph{#1}}
\newcommand \ixx[2] {\index{#1!#2}\emph{#1 #2}}
\newcommand \ixy[2] {\index{#2!#1}\emph{#1 #2}}
\newcommand \ixd[2] {\index{#2!#1}\emph{#1}}
\newcommand \ixe[2] {\index{#2@#1}\emph{#1}}
\newcommand \ixf[3] {\index{#3@#2!#1}\emph{#1}}
\newcommand \ixg[3] {\index{#3@#1!#2}\emph{#1}}

%:  Grandcadre
\newcommand \Grandcadre[1]{%
\begin{center}
\begin{tabular}{|c|}
\hline
~\\[-3mm]
#1\\[-3mm]
~\\ 
\hline
\end{tabular}
\end{center}
}

%:  lecteur, lectrice, masculin feminin
% compteur impair: masculin, compteur pair: feminin

\newcounter{MF}
\newcommand\stMF{\stepcounter{MF}}

\newcommand{\lec}{\stMF\ifodd\value{MF}lecteur \else lectrice \fi}
\newcommand{\lecz}{\stMF\ifodd\value{MF}lecteur\else lectrice\fi}

\newcommand{\lecs}{\stMF\ifodd\value{MF}lecteurs \else lectrices \fi}
\newcommand{\lecsz}{\stMF\ifodd\value{MF}lecteurs\else lectrices\fi}

\newcommand{\alec}{\stMF\ifodd\value{MF}au lecteur \else%
\`a la lectrice \fi}
\newcommand{\alecz}{\stMF\ifodd\value{MF}au lecteur\else%
\`a la lectrice\fi}

\newcommand{\dlec}{\stMF\ifodd\value{MF}du lecteur \else%
de la lectrice \fi}
\newcommand{\dlecz}{\stMF\ifodd\value{MF}du lecteur\else%
de la lectrice\fi}

\newcommand{\llec}{\stMF\ifodd\value{MF}le lecteur \else la lectrice \fi}
\newcommand{\llecz}{\stMF\ifodd\value{MF}le lecteur\else la lectrice\fi}

\newcommand{\Llec}{\stMF\ifodd\value{MF}Le lecteur \else La lectrice \fi}
\newcommand{\Llecz}{\stMF\ifodd\value{MF}Le lecteur\else La lectrice\fi}

% les suivants ne changent pas le genre
\newcommand{\lui}{\ifodd\value{MF}lui \else
elle \fi}
\newcommand{\luiz}{\ifodd\value{MF}lui\else
elle\fi}

\newcommand{\celui}{\ifodd\value{MF}celui \else
celle \fi}

\newcommand{\ceux}{\ifodd\value{MF}ceux \else
celles \fi}

\newcommand{\er}{\ifodd\value{MF}er \else
ère \fi}

\newcommand{\eux}{\ifodd\value{MF}eux \else
elles \fi}

\newcommand{\eUx}{\ifodd\value{MF}eux \else
euse \fi}

\newcommand{\leux}{\ifodd\value{MF}leux \else
leuse \fi}

\newcommand{\il}{\ifodd\value{MF}il \else
elle \fi}

\newcommand{\Sil}{\ifodd\value{MF}S'il \else
Si elle \fi}

\newcommand{\e}{\ifodd\value{MF} \else e \fi}
\newcommand{\ez}{\ifodd\value{MF}\else e\fi}

\newcommand{\n}{\ifodd\value{MF}n \else nne \fi}
\newcommand{\nz}{\ifodd\value{MF}n\else nne\fi}

%  \la* fait le ou la sans changer le genre,
%  \la fait le ou la en changeant de genre
\makeatletter
\newcommand{\la}{\@ifstar{\ifodd\value{MF}le\else
la\fi}{\stMF\ifodd\value{MF}le \else la \fi}}
\makeatother

%:  newtheorems  

%\theoremstyle{plain}
\CMnewtheorem{theorem}{Théorème}{\itshape}   
\CMnewtheorem{thdef}{Théorème et définition}{\itshape}
\CMnewtheorem{plcc}{Principe local-global concret}{\itshape}
\CMnewtheorem{prcf}{Principe de recouvrement fermé}{\itshape}
\CMnewtheorem{prmf}{Principe de recollement fermé}{\itshape}
\CMnewtheorem{prvq}{Principe de recouvrement par quotients}{\itshape}
\CMnewtheorem{prmq}{Principe de recollement de quotients}{\itshape}
\CMnewtheorem{plca}{Principe local-global abstrait$\mathbf{^*}$}{\itshape}
\CMnewtheorem{plcd}{Principe local-global dynamique}{\itshape}
\CMnewtheorem{proposition}{Proposition}{\itshape}
\CMnewtheorem{propdef}{Proposition et définition}{\itshape}
\CMnewtheorem{lemma}{Lemme}{\itshape}
\CMnewtheorem{corollary}{Corolaire}{\itshape}
\CMnewtheorem{fact}{Fait}{\itshape}
\CMnewtheorem{theoremc}{Th\'{e}or\`{e}me\etoz}{\itshape}
\CMnewtheorem{lemmac}{Lemme\etoz}{\itshape}
\CMnewtheorem{corollaryc}{Corolaire\etoz}{\itshape}
\CMnewtheorem{proprietec}{Propri\'{e}t\'{e}\etoz}{\itshape}
\CMnewtheorem{propositionc}{Proposition\etoz}{\itshape}
\CMnewtheorem{factc}{Fait\etoz}{\itshape}

\CMnewtheorem{remark}{Remarque}{}
\CMnewtheorem{remarks}{Remarques}{}
\CMnewtheorem{comment}{Commentaire}{}
\CMnewtheorem{comments}{Commentaires}{}
\CMnewtheorem{example}{Exemple}{}
\CMnewtheorem{examples}{Exemples}{}
\CMnewtheorem{definition}{Définition}{}
\CMnewtheorem{definitions}{Définitions}{}
\CMnewtheorem{definota}{Définition et notation}{}
\CMnewtheorem{definotas}{Définitions et notations}{}
\CMnewtheorem{convention}{Convention}{}
\CMnewtheorem{notation}{Notation}{} 
\CMnewtheorem{notations}{Notations}{} 
\CMnewtheorem{question}{Question}{}
%\CMnewtheorem{exercise}{Exercice}{\itshape}
%\CMnewtheorem{problem}{Problème}{}
\CMnewtheorem{algorithm}{Algorithme}{}

\newcounter{exercise}[chapter]
\newenvironment{exercise}{\ifhmode\par\fi
\vskip-\lastskip\vskip1.5ex\mou\penalty-300 \relax
\everypar{}\noindent
\refstepcounter{exercise}{\bfseries Exercice \theexercise.}\relax
\itshape\ignorespaces}{\par\vskip-\lastskip\vskip1em}

\newcounter{problem}[chapter]
\newenvironment{problem}{\ifhmode\par\fi
\vskip-\lastskip\vskip1em\mou\penalty-300 \relax
\everypar{}\noindent
\refstepcounter{problem}{\bfseries Problème \theproblem.}\relax
\itshape\ignorespaces}{\par\vskip-\lastskip\vskip1em}

\newcommand
\CHAP[1]{ 
\goodbreak\vskip4mm \mou \noindent  {\bf #1}\par\nobreak
\vskip1mm \mou \nobreak
}

\newcommand {\junk}[1]{}
\newcommand {\eop}{\hbox{$\square$}}

\newcommand\DebP{\raisebox{1pt}{\large\bf \Rightcircle}\,}
\newcommand\FinP{\raisebox{1pt}{\large\bf \Leftcircle$\!\!$}}

\def\thefootnote{\arabic{footnote}}

\newenvironment{proof}{\ifhmode\par\fi\vskip-\lastskip\vskip0.5ex\global\insidedemotrue
\everypar{}\noindent{\setbox0=\hbox{$\!\!$\DebP }\global\wdTitreEnvir\wd0\box0}
%: \DebP  a la place de \demoname  espace apres plus petit 
\ignorespaces}{\enddemobox\par\vskip.5em}
\def\enddemobox{\ifinsidedemo
\ifmmode\hbox{$\square$}\else
\ifhmode\unskip\else\noindent\fi\nobreak\null\nobreak\hfill
\nobreak$\square$\fi\fi\global\insidedemofalse}

%:  espace apres plus petit 
\newenvironment{Proof}[1]{\ifhmode\par\fi\vskip-\lastskip\vskip0.5ex\global\insidedemotrue
\everypar{}\noindent{\setbox0=\hbox{\it #1}\global\wdTitreEnvir\wd0\box0}\ignorespaces}{\enddemobox\par\vskip.5em}
\def\enddemobox{\ifinsidedemo
\ifmmode\hbox{$\square$}\else
\ifhmode\unskip\else\noindent\fi\nobreak\null\nobreak\hfill
\nobreak$\square$\fi\fi\global\insidedemofalse}

\newcommand\facile{\begin{proof}
La \dem est laissée \alecz.
\end{proof}
}

\newcommand\dsp{\displaystyle}
\newcommand\ndsp{\textstyle}

\newcommand\mapright[1]{\smash{\mathop{\longrightarrow}\limits^{#1}}} 
\newcommand\maprightto[1]{\smash{\mathop{\longmapsto}\limits^{#1}}} 
\def\mapdown#1{\downarrow\rlap{$\vcenter{\hbox{$\scriptstyle 
#1$}}$}}

\newcommand\oge{\leavevmode\raise.3ex\hbox{$\scriptscriptstyle\langle\!\langle\,$}}
\newcommand\feg{\leavevmode\raise.3ex\hbox{$\scriptscriptstyle\,\rangle\!\rangle$}}
\DeclareRobustCommand{\guig}{\mbox{{\usefont{U}{lasy}%
{\if b\expandafter\@car\f@series\@nil b\else m\fi}{n}%
\char40\kern-0.20em\char40}~}}
\DeclareRobustCommand{\guid}{\mbox{~\usefont{U}{lasy}%
{\if b\expandafter\@car\f@series\@nil b\else m\fi}{n}%
\char41\kern-0.20em\char41}}
\newcommand\gui[1]{\oge{#1}\feg}

%:  subsection* avec table des matières
%%%%%%%%%%%%%%%%%%%%%%%%%%%%%%%%%%%%%%%%%
\newcommand\subsec[1]{\goodbreak\rdb\addcontentsline{toc}{subsection}{#1}\subsection*{#1}}
%%%%%%%%%%%%%%%%%%%%%%%%%%%%%%%%%%%%%%%%%

%%%%%%%%%%%%%%%%%%%%%%%%%%%%%%%%%%%%%%%%%
\newcommand\subsect[2]{\goodbreak\rdb\addcontentsline{toc}{subsection}{#2}
\subsection*{#1}}
%%%%%%%%%%%%%%%%%%%%%%%%%%%%%%%%%%%%%%%%%

%%%%%%%%%%%%%%%%%%%%%%%%%%%%%%%%%%%%%%%%%
\newcommand\subsubs[1]{

\goodbreak\rdb\medskip 

{\bf #1}

\smallskip }
%%%%%%%%%%%%%%%%%%%%%%%%%%%%%%%%%%%%%%%%%

%%%%%%%%%%%%%%%%%%%%%%%%%%%%%%%%%%%%%%%%%
\newcommand\subsubsec[1]{\goodbreak\rdb\addcontentsline{toc}{subsubsection}{#1}\subsubsection*{#1}}
%%%%%%%%%%%%%%%%%%%%%%%%%%%%%%%%%%%%%%%%%

%%%%%%%%%%%%%%%%%%%%%%%%%%%%%%%%%%%%%%%%%
\newcommand\subsubsect[2]{\goodbreak\rdb\addcontentsline{toc}{subsubsection}{#2}\subsubsection*{#1}}
%%%%%%%%%%%%%%%%%%%%%%%%%%%%%%%%%%%%%%%%%

% ci apres a utiliser si on veut des subsubsection* directement apres section
% sans passer par subsection
%%%%%%%%%%%%%%%%%%%%%%%%%%%%%%%%%%%%%%%%%
\newcommand\Subsubsec[1]{\goodbreak\rdb\addcontentsline{toc}{subsection}{#1}\subsubsection*{#1}}
%%%%%%%%%%%%%%%%%%%%%%%%%%%%%%%%%%%%%%%%%

%%%%%%%%%%%%%%%%%%%%%%%%%%%%%%%%%%%%%%%%%
\newcommand\Subsubsect[2]{\rdb\addcontentsline{toc}{subsection}{#2}\subsubsection*{#1}}
%%%%%%%%%%%%%%%%%%%%%%%%%%%%%%%%%%%%%%%%%

%%%%%%%%%%%%%%%%%%%%%%%%%%%%%%%%%%%%%%%%%
\newcommand\Subsec[1]{\goodbreak\rdb\addcontentsline{toc}{section}{#1}\subsection*{#1}}
%%%%%%%%%%%%%%%%%%%%%%%%%%%%%%%%%%%%%%%%%

%%%%%%%%%%%%%%%%%%%%%%%%%%%%%%%%%%%%%%%%%
\newcommand\Exercices{\rdb\addcontentsline{toc}{section}{Exercices et \pbsz}%
\markright{Exercices et \pbsz}\section*{Exercices et \pbsz}\pagestyle{CMExercicesheadings}\small }
%%%%%%%%%%%%%%%%%%%%%%%%%%%%%%%%%%%%%%%%%

%%%%%%%%%%%%%%%%%%%%%%%%%%%%%%%%%%%%%%%%%
\newcommand\sol{%
\rdb\normalsize\addcontentsline{toc}{subsection}{Solutions d'exercices}%
\markright{Solutions d'exercices}\subsection*{Quelques solutions, ou esquisses de solutions}%
\pagestyle{CMExercicesheadings}\small}
%%%%%%%%%%%%%%%%%%%%%%%%%%%%%%%%%%%%%%%%%

%: exo corr
\newcommand\exer[1]{{\smallskip \goodbreak\noindent \textbf{Exercice \ref{#1}.} }}

\newcommand\prob[1]{{\smallskip \goodbreak\noindent \textbf{Problème \ref{#1}.} }}

%%%%%%%%%%%%%%%%%%%%%%%%%%%%%%%%%%%%%%%%%
\newcommand\Biblio{\goodbreak\rdb\normalsize\addcontentsline{toc}{section}{Commentaires bibliographiques}%
\markright{Commentaires bibliographiques}%
\pagestyle{CMExercicesheadings}\section*{Commentaires bibliographiques} 

}
%%%%%%%%%%%%%%%%%%%%%%%%%%%%%%%%%%%%%%%%%

%%%%%%%%%%%%%%%%%%%%%%%%%%%%%%%%%%%%%%%%%
\newcommand\Intro{\addcontentsline{toc}{section}{Introduction}\markright{Introduction}%
\thispagestyle{CMExercicesheadings}\subsection*{Introduction} }
%%%%%%%%%%%%%%%%%%%%%%%%%%%%%%%%%%%%%%%%%

%:  definition alternative avec meme numero
\newcommand\defa[2]{\mni\rdb\textbf{Définition alternative \ref{#1}
\it #2 
}}

%:  def  commandes a parametres

%: \renewcommand\matrix
\renewcommand\matrix[1]{{\begin{array}{ccccccccccccccccccccccccc} #1 \end{array}}}

\newcommand\MA[1]{\mathop{#1}\nolimits}

\newcommand{\vect}[1]{\mathchoice{\overrightarrow{\strut#1}}%
{\overrightarrow{\textstyle\strut#1}}{\overrightarrow{\scriptstyle#1}}{\overrightarrow{\scriptscriptstyle#1}}}

\newcommand\vab[2]{[\,#1\;#2\,]}

\newcommand\abs[1]{\left|{#1}\right|}
\newcommand\abS[1]{\big|{#1}\big|}
\newcommand\aqo[2]{#1\sur{\gen{#2}}\!}
\newcommand\Aqo[2]{#1\sur{\big\langle{#2}\big\rangle}\!}
\newcommand\bloc[4]{\left[\matrix{#1 & #2 \cr #3 & #4}\right]}

\newcommand\carray[2]{{\left[\begin{array}{#1} #2 \end{array}\right]}}
\newcommand\cmatrix[1]{\left[\matrix{#1}\right]}
\newcommand\clmatrix[1]{{\left[\begin{array}{lllllll} #1 \end{array}\right]}}
\newcommand\crmatrix[1]{{\left[\begin{array}{rrrrrrr} #1 \end{array}\right]}}
\newcommand\dmatrix[1]{\abs{\matrix{#1}}}
\newcommand\Cmatrix[2]{\setlength{\arraycolsep}{#1}\left[\matrix{#2}\right]}
\newcommand\Dmatrix[2]{\setlength{\arraycolsep}{#1}\left|\matrix{#2}\right|}

\newcommand{\Dpp}[2]{{{\partial #1}\over{\partial #2}}}

\newcommand\Dlu[2]{{\rm Dl}_{#1}(#2)}
\newcommand\dessus[2]{{\textstyle {#1} \atop \textstyle {#2}}}
\newcommand\eqdf[1]{\buildrel{#1}\over =}
\newcommand\formule[1]{{\left\{ {\arraycolsep2pt\begin{array}{lll} #1 \end{array}}\right.}}
\newcommand\formul[2]{{\left\{ \begin{array}{#1} #2 \end{array}\right.}}
\newcommand\gen[1]{\left\langle{#1}\right\rangle}
\newcommand\geN[1]{\big\langle{#1}\big\rangle}
\newcommand\impdef[1]{\buildrel{#1}\over \Longrightarrow}

\newcommand\eqdefi{\eqdf{\rm def}}
\newcommand\eqdef{\buildrel{{\rm def}}\over \Longleftrightarrow }

\newcommand{\Kr}[2]{#1\lrb{#2}}

\newcommand\lra[1]{\langle{#1}\rangle}
\newcommand\lrb[1] {\llbracket #1 \rrbracket}
\newcommand\lrbn {\lrb{1..n}}
\newcommand\lrbzn {\lrb{0..n}}
\newcommand\lrbl {\lrb{1..\ell}}
\newcommand\lrbm {\lrb{1..m}}
\newcommand\lrbk {\lrb{1..k}}
\newcommand\lrbh {\lrb{1..h}}
\newcommand\lrbp {\lrb{1..p}}
\newcommand\lrbq {\lrb{1..q}}
\newcommand\lrbr {\lrb{1..r}}
\newcommand\lrbs {\lrb{1..s}}

\newcommand\fraC[2] {{{#1}\over {#2}}}
\newcommand\nsup[1] {\Vert#1\Vert_\infty}
\newcommand\meck[2] {\{#1, #2\}}
\newcommand\sat[1] {#1^{\rm sat}}
\newcommand\satu[2] {#1^{\rm sat_{#2}}}
\newcommand\scp[2] {\gen{#1\,|\,#2}\!}
\newcommand\sur[1]{\!\left/#1\right.}
\newcommand\so[1]{\left\{{#1}\right\}}
\newcommand\sO[1]{\big\{{#1}\big\}}
\newcommand\sotq[2]{\so{\,#1\,\vert\,#2\,}}
\newcommand\sotQ[2]{\sO{\,#1\;\big\vert\;#2\,}}
\newcommand\frt[1]{\!\left|_{#1}\right.\!}
\newcommand\Frt[2]{\left.#1\right|_{#2}\!}
\newcommand\sims[1]{\buildrel{#1}\over \sim}
\newcommand\tra[1]{{\,^{\rm t}\!#1}}
\newcommand\tralst[1] {\tra\,{\lst{#1}}}
\newcommand\Al[1]{\Vi^{\!#1}}
\newcommand\Ae[1]{\gA^{\!#1}}

%:HHH je change snic snac snuc
\newcommand\Snic[1]{$$\nds #1$$}
\newcommand\snic[1]{

{\centering$#1$\par}

}
\newcommand\snac[1]{

{\small\centering$#1$\par}

}
\newcommand\snuc[1]{

{\footnotesize\centering$#1$\par}
}

\newcommand\snucc[1]{

{\footnotesize$$\preskip0pt\postskip0pt\textstyle#1$$}}

\newcommand\snicc[1]{
{$$\preskip0pt\postskip0pt\textstyle#1$$}}

\newcommand\snif[3]{\vspace{#1}\noindent\centerline{$#3$}
\vspace{#2}}
%exemple   \snif{.1cm}{.2cm}{ma formule}

\newcommand\env[2]{{{#2}_{#1}^{\mathrm{e}}}} % algebre enveloppante
\newcommand\Om[2]{\Omega_{{#2}/{#1}}}
\newcommand\Der[3]{{\rm Der}_{{#1}}({#2},{#3})}

\newcommand \isA[1] {_{#1/\!\gA}}
\newcommand\OmA[1]{\Omega\isA{#1}}

\newcommand \bra[1] {\left[{#1}\right]}
\newcommand \bu[1] {{{#1}\bul}}
\newcommand \ci[1] {{{#1}^\circ}}
\newcommand \wi[1] {\widetilde{#1}}
\newcommand \wh[1]{{\widehat{#1}}}
\newcommand \ov[1] {\overline{#1}}
\newcommand \und[1] {\underline{#1}}
\newcommand \uci[1]{{\buildrel{\circ}\over{#1}}}

%:H il semble que le package variroref pose pb 
% a cause de notre redefinition de \ref
\newcommand{\pref}[1]{\textup{\hbox{\normalfont(\ref{#1})}}}
\newcommand \vref[1] {\ref{#1}}
\newcommand \vpageref[1] {\paref{#1}}
\newcommand \thref[1] {\thoz~\ref{#1}}
\newcommand \paref[1] {page~\pageref{#1}}
\newcommand \thrf[1] {\thoz~\ref{#1}}
\newcommand \thrfs[2] {\thosz~\ref{#1} et~\ref{#2}}
\renewcommand \rref[1] {\ref{#1} \paref{#1}}
\newcommand \egrf[1] {\egtz~(\ref{#1})}
\newcommand \egref[1] {\egtz~(\ref{#1})  \paref{#1}}
\newcommand \eqrf[1] {équation~(\ref{#1})}
\newcommand \eqvrf[1] {équation (\ref{#1}) \paref{#1}}
\newcommand \prirf[1] {prin\-cipe~\ref{#1}}
\newcommand \plgrf[1] {\plgz~\ref{#1}}
\newcommand \plgref[1] {\plgz~\ref{#1}}
\newcommand \plgrfs[2] {\plgs \ref{#1} et~\ref{#2}}

\newcommand\VRT[1]{\rotatebox{90}{\hbox{$#1$}}}
\newcommand\VRTsubseteq{\VRT{\subseteq}}
\newcommand\VRTsupseteq{\VRT{\supseteq}}
\newcommand\VRTlongrightarrow{\VRT{\longrightarrow}}
\newcommand\VRTlongleftarrow{\VRT{\longleftarrow}}

\newcommand \rC[1]{\MA{{\rm C}_{#1}}}
\newcommand \rF[1]{\MA{{\rm F}_{\!#1}}}
\newcommand \rR[1]{\MA{{\rm R}_{#1}}}
\newcommand \rRs[1]{\MA{\Rs_{#1}}}
\newcommand \ep[1]{^{(#1)}}

\newcommand \bal[1] {^\rK_{#1}}
\newcommand \ul[1] {_\rK^{#1}}
\newcommand \SNw[1] {P_{#1}}
\newcommand \gBtst {\gB[[t]]^{\!\times}}

% indice d'un sous groupe
\newcommand \idg[1] {\hbox{$|\,#1\,|$}}
\newcommand \idG[1] {\hbox{$\big|\,#1\,\big|$}}

% dimension d'une extension
\newcommand \dex[1] {\hbox{$[\,#1\,]$}}
\newcommand \deX[1] {\hbox{$\big[\,#1\,\big]$}}

% liste
\newcommand \lst[1] {\hbox{$[\,#1\,]$}}
\newcommand \lsT[1] {\hbox{$\big[\,#1\,\big]$}}
%:h2013 rajout
\newcommand \brk[1] {[#1]}

\newcommand \THo[2]{\rdb
\mni{\bf Théorème {#1}~} {\it #2

}}

\newcommand \PLCC[2]{\rdb
\mni{\bf Principe \lgb concret \ref{#1} bis~} {\it #2

}}

\newcommand \THO[2]{\rdb
\mni{\bf Théorème \ref{#1} bis~} {\it #2

}}

%:   fleches longues
\newcommand{\mt}{\mapsto}
\newcommand \lmt{\longmapsto}

\newcommand{\llongrightarrow}{\relbar\joinrel\mkern-1mu\longrightarrow}
\newcommand{\lllongrightarrow}{\relbar\joinrel\mkern-1mu\llongrightarrow}
\newcommand{\llllongrightarrow}{\relbar\joinrel\mkern-1mu\lllongrightarrow}
\newcommand{\lllllongrightarrow}{\relbar\joinrel\mkern-1mu\llllongrightarrow}
\newcommand \lora {\longrightarrow}
\newcommand \llra {\llongrightarrow}
\newcommand \lllra {\lllongrightarrow}
\newcommand \llllra {\llllongrightarrow}
\newcommand \lllllra {\lllllongrightarrow}
\newcommand\simarrow{\vers{_\sim}}
\newcommand\Simarrow{\buildrel{_\sim}\over \longleftrightarrow}
\newcommand\isosim {\simarrow}
\newcommand\vers[1]{\buildrel{#1}\over \lora }
\newcommand\vvers[1]{\buildrel{#1}\over \llra }
\newcommand\vvvers[1]{\buildrel{#1}\over \lllra }
\newcommand\vvvvers[1]{\buildrel{#1}\over \llllra }
\newcommand\vvvvvers[1]{\buildrel{#1}\over \lllllra }

%:  leqslant
\renewcommand \le{\leqslant}
\renewcommand \leq{\leqslant}
\renewcommand \preceq{\preccurlyeq}
\renewcommand \ge{\geqslant}
\renewcommand \geq{\geqslant}
\renewcommand \succeq{\succurlyeq}
\newcommand \cuvu {\curlyvee}
\newcommand \cuvi {\curlywedge}

%:  tableau deux colonnes
\newcommand \DeuxCol[2]{%
\sni\mbox{\parbox[t]{.475\textwidth}{#1}%
\hspace{.05\textwidth}%
\parbox[t]{.475\textwidth}{#2}}}

\newcommand \Deuxcol[4]{%
\sni\mbox{\parbox[t]{#1\textwidth}{#3}%
\hspace{.05\textwidth}%
\parbox[t]{#2\textwidth}{#4}}}

%:  Algorithmes

\floatstyle{boxed}
\floatname{agc}{Algorithme}
\newfloat{agc}{ht}{lag}[section]

\floatstyle{boxed}
\floatname{agC}{Algorithme}
\newfloat{agC}{H}{lag}[section]

\newenvironment{algor}[1][]%
{\par\smallskip\begin{agc}
\vskip 1mm
\begin{algorithm}{\bfseries#1}
\upshape\sffamily
}
{\end{algorithm}
\end{agc}
}

\newenvironment{algoR}[1][]%
{\par\smallskip\begin{agC}
\vskip 1mm
\begin{algorithm}{\bfseries#1}
\upshape\sffamily
}
{\end{algorithm}
\end{agC}
}

\newcommand\Vrai{\mathsf{Vrai}}
\newcommand\Faux{\mathsf{Faux}}
\newcommand\ET{\mathsf{ et }}
\newcommand\OU{\mathsf{ ou }}
\newcommand\pour[3]{\textbf{pour } $#1$ \textbf{ de } $#2$
           \textbf{ \`a } $#3$ \textbf{ faire }}
\newcommand\pur[2]{\textbf{pour } $#1$ \textbf{ dans } $#2$
            \textbf{ faire }}
\newcommand\por[3]{\textbf{pour } $#1$ \textbf{ de } $#2$
           \textbf{ \`a } $#3$  }
\def\sialors#1{\textbf{si } $#1$ \textbf{ alors }}
\def\tantque#1{\textbf{tant que } $#1$ \textbf{ faire }}
\newcommand\finpour{\textbf{fin pour}}
\newcommand\sinon{\textbf{sinon }}
\newcommand\finsi{\textbf{fin si }}
\newcommand\fintantque{\textbf{fin tant que }}
\newcommand\aff{\leftarrow }
\newcommand\Debut{\\[1mm] \textbf{Début }}
\newcommand\Fin{\textbf{\\ Fin.}}
\newcommand\Entree{\\[1mm] \textbf{Entrée : }}
\newcommand\Sortie{\\ \textbf{Sortie : }}
\newcommand\Varloc{\\ \textbf{Variables locales : }}
\newcommand\Repeter{\textbf{Répéter }}
\newcommand\jusqua{\textbf{jusqu\`a ce que }}
\newcommand\hsz{\\ }
\newcommand\hsu{\\ \hspace*{4mm}}
\newcommand\hsd{\\ \hspace*{8mm}}
\newcommand\hst{\\ \hspace*{1,2cm}}
\newcommand\hsq{\\ \hspace*{1,6cm}}
\newcommand\hsc{\\ \hspace*{2cm}}
\newcommand\hsix{\\ \hspace*{2,4cm}}
\newcommand\hsept{\\ \hspace*{2,8cm}}

\newcommand \legendre {\overwithdelims()}
\newcommand \legendr[2] {\Big(\frac {#1}{#2} \Big)}
\newcommand \som {\sum\nolimits}
\newcommand \ds {\displaystyle}

\newcommand\et{\;\;\hbox{ et }\;\;}

%:  treillis distributifs, divisibilité

\newcommand \divi {\mid}
\def \nedivi {\not\kern 2.5pt\mid}

\newcommand \vu {\vee} % sup dans les treillis
\newcommand \vi {\wedge} % inf dans les treillis
\newcommand \Vu {\bigvee\nolimits}
\newcommand \Vi {\bigwedge\nolimits}

\newcommand \vda {\,\vdash\,}
\newcommand \vdi[1] {\,\vdash_{#1}\,}
\newcommand \im {\rightarrow} % fleche des treillis implicatifs
\newcommand \dar[1] {\MA{\downarrow \!#1}}
\newcommand \uar[1] {\MA{\uparrow \!#1}}

\newcommand \Un {\mathbf{1}}
\newcommand \Deux {\mathbf{2}}
\newcommand \Trois {\mathbf{3}}
\newcommand \Quatre {\mathbf{4}}
\newcommand \Cinq {\mathbf{5}}

%:  logique

\newcommand \Pf {{\rm P}_{{\rm f}}}
\newcommand \Pfe {{\rm P}_{{\rm fe}}}

\newcommand \Ex {{\exists}}
\newcommand \Tt {{\forall}}

\newcommand \Lst{\mathsf{Lst}}
\newcommand \Irr{\mathsf{Irr}}
\newcommand \Prim{\mathsf{Prim}}
\newcommand \Rec{\mathsf{Rec}}

%: tsbf

\newcommand \tsbf[1]{\textbf{\textsf{#1}}}
\newcommand \FAN{\tsbf{FAN} }
\newcommand \FANz{\tsbf{FAN}}
\newcommand \KL{\FAN}
\newcommand \KLz{\FANz}
\newcommand \KLp{\tsbf{KL}$_2$ }
\newcommand \KLpz{\tsbf{KL}$_2$}
\newcommand \kl{\tsbf{KL}$_1$ }
\newcommand \klz{\tsbf{KL}$_1$}
\newcommand \HC{\tsbf{HC} }
\newcommand \HCz{\tsbf{HC}}
\newcommand \LLPO{\tsbf{LLPO} }
\newcommand \LLPOz{\tsbf{LLPO}}
\newcommand \LPO{\tsbf{LPO} }
\newcommand \LPOz{\tsbf{LPO}}
\newcommand \MP{\tsbf{MP} }
\newcommand \MPz{\tsbf{MP}}
\newcommand \TEM{\tsbf{PTE} }
\newcommand \TEMz{\tsbf{PTE}}
\newcommand \UC{\tsbf{UC} }
\newcommand \UCz{\tsbf{UC}}
\newcommand \UCp{\tsbf{UC}$^+$ }
\newcommand \UCpz{\tsbf{UC}$^+$}
\newcommand \Mini{\tsbf{Min} }
\newcommand \Minip{\tsbf{Min}$^+$ }
\newcommand \Minipz{\tsbf{Min}$^+$}
\newcommand \Minim{\tsbf{Min}$^-$ }
\newcommand \Miniz{\tsbf{Min}}

%:  en exposant ou en indice
\newcommand \bul{^\bullet}
\newcommand \eci{^\circ}
\newcommand \esh{^\sharp}
\newcommand \efl{^\flat}
\newcommand \epr{^\perp}
\newcommand \eti{^\times}
\newcommand \etl{^* }
\newcommand \eto{$^*\!$ }
\newcommand \etoz{$^*\!$}
\newcommand \sta{^\star}
\newcommand \ist{_\star}
\newcommand \eo {^{\mathrm{op}}}

\newcommand\Abul {\gA\!\bul}
\newcommand\Ati {\gA^{\!\times}}
\newcommand\Asta {\gA^{\!\star}}
\newcommand\Atl {\gA^{\!*}}
\newcommand\Bti {\gB^{\times}}
\newcommand\Bst {\Bti}
\newcommand\te  {\otimes}

\newcommand \iBA {_{\gB/\!\gA}}
\newcommand \iBk {_{\gB/\gk}}
\newcommand \iBK {_{\gB/\gK}}
\newcommand \iAk {_{\gA/\gk}}
\newcommand \iAK {_{\gA/\gK}}
\newcommand \iCk {_{\gC/\gk}}

\newcommand \tgaBG {\gB\{G\}}  %% twisted group algebra of G
%% first cohomology group
\newcommand \zcoho {Z^1(G, \Bti)}
\newcommand \bcoho {B^1(G, \Bti)}
\newcommand \hcoho {H^1(G, \Bti)}

\newcommand \vep{{\varepsilon}}

%:HHH espace rétréci autour de equidef
\newcommand\equidef{\buildrel{{\rm def}}\over{\;\Longleftrightarrow\;}}

\newcommand \E{\mathaccent19}    %accent aigu math
\newcommand \aigu{\mathaccent19}    %accent aigu math

%:  noindent
\newcommand \noi {\noindent}
\renewcommand \ss {\smallskip}
\newcommand \sni {\ss\noi}
\newcommand \snii {}
\newcommand \ms {\medskip}
\newcommand \mni {\ms\noi}
\newcommand \bs {\bigskip}
\newcommand \bni {\bs\noi}
\newcommand \ce{\centering}
\newcommand \alb{\allowbreak}

%:  souligné tout pret
\newcommand \ua  {{\underline{a}}}
\newcommand \ub  {{\underline{b}}}
\newcommand \ual {{\underline{\alpha}}}
\newcommand \ube {{\underline{\beta}}}
\newcommand \uc  {{\underline{c}}}
\newcommand \ud  {{\underline{d}}}
\newcommand \udel{{\underline{\delta}}}
\newcommand \ue  {{\underline{e}}}
\newcommand \uf  {{\underline{f}}}
\newcommand \uF  {{\underline{F}}}
\newcommand \ug  {{\underline{g}}}
\newcommand \uh  {{\underline{h}}}
\newcommand \uH  {{\underline{H}}}
\newcommand \uga {{\underline{\gamma}}}
\newcommand \ur  {{\underline{r}}}
\newcommand \us  {{\underline{s}}}
\newcommand \ut  {{\underline{t}}}
\newcommand \uu  {{\underline{u}}}
\newcommand \ux  {{\underline{x}}}
\newcommand \uxi {{\underline{\xi}}}
\newcommand \uy  {{\underline{y}}}
\newcommand \uP  {{\underline{P}}}
\newcommand \uS  {{\underline{S}}}
\newcommand \uT  {{\underline{T}}}
\newcommand \uU  {{\underline{U}}}
\newcommand \uX  {{\underline{X}}}
\newcommand \uY  {{\underline{Y}}}
\newcommand \uZ  {{\underline{Z}}}
\newcommand \uz  {{\underline{z}}}
\newcommand \uzeta  {{\underline{\zeta}}}
\newcommand \uze {{\underline{0}}}

%:  séquences x_1,\ldots,x_m toutes pretes
\newcommand \am {a_1,\ldots,a_m}
\newcommand \an {a_1,\ldots,a_n}
\newcommand \bn {b_1,\ldots,b_n}
\newcommand \bbm {b_1,\ldots,b_m}
\newcommand \br {b_1,\ldots,b_r}
\newcommand \azn {a_0,\ldots,a_n}
\newcommand \bzn {b_0,\ldots,b_n}
\newcommand \czn {c_0,\ldots,c_n}
\newcommand \cq {c_1,\ldots,c_q}
\newcommand \gq {g_1,\ldots,g_q}
\newcommand \un {u_1,\ldots,u_n}
\newcommand \xk {x_1,\ldots,x_k}
\newcommand \Xk {X_1,\ldots,X_k}
\newcommand \xl {x_1,\ldots,x_\ell}
\newcommand \xm {x_1,\ldots,x_m}
\newcommand \xn {x_1,\ldots,x_n}
\newcommand \xp {x_1,\ldots,x_p}
\newcommand \yp {y_1,\ldots,y_p}
\newcommand \xq {x_1,\ldots,x_q}
\newcommand \xzk {x_0,\ldots,x_k}
\newcommand \xzn {x_0,\ldots,x_n}
\newcommand \xhn {x_0:\ldots:x_n}
\newcommand \Xn {X_1,\ldots,X_n}
\newcommand \Xzn {X_0,\ldots,X_n}
\newcommand \Xm {X_1,\ldots,X_m}
\newcommand \Xr {X_1,\ldots,X_r}
\newcommand \xr {x_1,\ldots,x_r}
\newcommand \Yr {Y_1,\ldots,Y_r}
\newcommand \Yn {Y_1,\ldots,Y_n}
\newcommand \ym {y_1,\ldots,y_m}
\newcommand \Ym {Y_1,\ldots,Y_m}
\newcommand \yk {y_1,\ldots,y_k}
\newcommand \yr {y_1,\ldots,y_r}
\newcommand \yn {y_1,\ldots,y_n}
\newcommand \zn {z_1,\ldots,z_n}
\newcommand \Zn {Z_1,\ldots,Z_n}

\newcommand \Sun {$S_1$, $\dots$, $S_n$ }
\newcommand \Sunz {$S_1$, $\dots$, $S_n$}

\newcommand \xpn {x'_1,\ldots,x'_n}
\newcommand \uxp  {{\underline{x'}}}
\newcommand \ypm {y'_1,\ldots,y'_m}
\newcommand \uyp  {{\underline{y'}}}

\newcommand \aln {\alpha_1,\ldots,\alpha_n}
\newcommand \gan {\gamma_1,\ldots,\gamma_n}
\newcommand \xin {\xi_1,\ldots,\xi_n}
\newcommand \xihn {\xi_0:\ldots:\xi_n}

\newcommand \lfs {f_1,\ldots,f_s}

\newcommand \Cin{C^{\infty}}
\newcommand \Ared {\gA\red}
\newcommand \Aqim {\gA\qim}
\newcommand \Amin {\Aqim}
\newcommand \Aqi {\gA_\mathrm{qi}}

%:  anneaux de polynomes
\newcommand \AT {{\gA[T]}}
\newcommand \Ax {{\gA[x]}}
\newcommand \AX {{\gA[X]}}
\newcommand \AY {{\gA[Y]}}
\newcommand \ArX {{\gA\lra X}}
\newcommand \ArY {{\gA\lra Y}}
\newcommand \BX {{\gB[X]}}
\newcommand \BY {{\gB[Y]}}
\newcommand \kx {{\gk[x]}}
\newcommand \kX {{\gk[X]}}
\newcommand \kT {{\gk[T]}}
\newcommand{\Kfi}{{\gK[\varphi]}}
\newcommand \KT {{\gK[T]}}
\newcommand \Kx {{\gK[x]}}
\newcommand \KX {{\gK[X]}}
\newcommand \KY {{\gK[Y]}}
\newcommand \QQX {{\QQ[X]}}
\newcommand \VX {{\gV[X]}}
\newcommand \ZZX {{\ZZ[X]}}
\newcommand \ZZx {{\ZZ[x]}}

\newcommand \AXY {\AuX \lra Y}

\newcommand \AXm {{\gA[\Xm]}}

\newcommand \kGa {{\gk[\Gamma]}}

\newcommand \kXn {{\gk[\Xn]}}
\newcommand \lXn {{\gl[\Xn]}}
\newcommand \AXk {{\gA[\Xk]}}
\newcommand \AXn {{\gA[\Xn]}}
\newcommand \BXn {{\gB[\Xn]}}
\newcommand \CXn {{\gC[\Xn]}}
\newcommand \KXn {{\gK[\Xn]}}
\newcommand \LXn {{\gL[\Xn]}}
\newcommand \RXn {{\gR[\Xn]}}
\newcommand \QQXk {{\QQ[\Xk]}}
\newcommand \QQXn {{\QQ[\Xn]}}
\newcommand \RXzn {{\gR[\Xzn]}}
\newcommand \Rxzn {{\gR[\xzn]}}
\newcommand \Ryn {{\gR[\yn]}}
\newcommand \VXn {{\gV[\Xn]}}
\newcommand \ZZXn {{\ZZ[\Xn]}}
\newcommand \ZZXk {{\ZZ[\Xk]}}

\newcommand \kXr {{\gk[\Xr]}}
\newcommand \lXr {{\gl[\Xr]}}
\newcommand \AXr {{\gA[\Xr]}}
\newcommand \KXr {{\gK[\Xr]}}
\newcommand \LXr {{\gL[\Xr]}}
\newcommand \RXr {{\gR[\Xr]}}
\newcommand \VXr {{\gV[\Xr]}}

\newcommand \kYr {{\gk[\Yr]}}
\newcommand \kyr {{\gk[\yr]}}
\newcommand \kyn {{\gk[\yn]}}
\newcommand \AYr {{\gA[\Yr]}}
\newcommand \KYr {{\gK[\Yr]}}
\newcommand \kXm {{\gk[\Xm]}}
\newcommand \KXm {{\gK[\Xm]}}
\newcommand \KYm {{\gK[\Ym]}}
\newcommand \kYm {{\gk[\Ym]}}
\newcommand \BYm {{\gB[\Ym]}}
\newcommand \lYr {{\gl[\Yr]}}
\newcommand \LYr {{\gL[\Yr]}}
\newcommand \RYr {{\gR[\Yr]}}

\newcommand \Axr {{\gA[\xr]}}
\newcommand \Kxr {{\gK[\xr]}}
\newcommand \kxr {{\gk[\xr]}}
\newcommand \Ayr {{\gA[\yr]}}
\newcommand \Kyr {{\gK[\yr]}}
\newcommand \Ayn {{\gA[\yn]}}
\newcommand \AYn {{\gA[\Yn]}}

\newcommand \kxm {{\gk[\xm]}}
\newcommand \kxn {{\gk[\xn]}}
\newcommand \lxn {{\gl[\xn]}}
\newcommand \Axn {{\gA[\xn]}}
\newcommand \Bxn {{\gB[\xn]}}
\newcommand \Cxn {{\gC[\xn]}}
\newcommand \Kxn {{\gK[\xn]}}
\newcommand \Kyn {{\gK[\yn]}}
\newcommand \Lxn {{\gL[\xn]}}
\newcommand \Rxn {{\gR[\xn]}}

\newcommand \Aux {{\gA[\ux]}}
\newcommand \Auy {{\gA[\uy]}}
\newcommand \Bux {{\gB[\ux]}}
\newcommand \Kuy {{\gK[\uy]}}
\newcommand \Kuu {{\gK[\uu]}}
\newcommand \Kux {{\gK[\ux]}}
\newcommand \kux {{\gk[\ux]}}
\newcommand \kuy {{\gk[\uy]}}
\newcommand \Rux {{\gR[\ux]}}
\newcommand \Ruy {{\gR[\uy]}}

\newcommand \ZZuS {{\ZZ[\uS]}}

\newcommand \AuX {{\gA[\uX]}}
\newcommand \BuX {{\gB[\uX]}}
\newcommand \KuX {{\gK[\uX]}}
\newcommand \kuX {{\gk[\uX]}}
\newcommand \LuX {{\gL[\uX]}}
\newcommand \RuX {{\gR[\uX]}}
\newcommand \ZZuX {{\ZZ[\uX]}}
\newcommand \QQuX {{\QQ[\uX]}}
\newcommand \AuY {{\gA[\uY]}}
\newcommand \BuY {{\gB[\uY]}}
\newcommand \KuY {{\gK[\uY]}}
\newcommand \kuY {{\gk[\uY]}}

%:  Foncteurs grasmanniennes, GLn,SLn ...
\newcommand \Gn  {\gG_n}
\newcommand \Gnk {\gG_{n,k}}
\newcommand \Gnr {\gG_{n,r}}
\newcommand \cGn {{\cG_n}}
\newcommand \cGnk{{\cG_{n,k}}}

\newcommand \GGn {\GG_{n}}
\newcommand \GGnk{{\GG_{n,k}}} % Grasmannienne ``projective'' usuelle
\newcommand \GGnr{{\GG_{n,r}}}
\newcommand \GA  {\mathbb{GA}}
\newcommand \GAn {\GA_{n}}  % Grasmannienne ``affine'': matrices de proj
\newcommand \GAq {\GA_{q}}
\newcommand \GAnk{\GA_{n,k}}
\newcommand \GAnr{\GA_{n,r}}

\newcommand \Mm {{\MM_{m}}}
\newcommand \Mn {{\MM_{n}}}
\newcommand \Mk {{\MM_{k}}}
\newcommand \Mq {{\MM_{q}}}
\newcommand \Mr {{\MM_{r}}}
\newcommand \MMn {{\MM_{n}}}

\newcommand \Bo{\BB\mathrm{o}}

\newcommand \GL {\mathbb{GL}}
\newcommand \GLn {{\GL_n}}
\newcommand \Gl {\mathbf{GL}}
\newcommand \Gln {{\Gl_n}}
\newcommand \PGL {\mathbb{PGL}}
\newcommand \SL {\mathbb{SL}}
\newcommand \SLn {{\SL_n}}
\newcommand \EE {\mathbb{E}}
\newcommand \En {\EE_n}
\newcommand \Pn {\PP^n}
\newcommand \An {\AA^{\!n}}
\newcommand \Am {\AA^{\!m}}
\newcommand \Sl {\mathbf{SL}}
\newcommand \Sln {{\Sl_n}}

%:   mathrm
\newcommand \I  {\mathrm{I}}
\newcommand \G  {\mathrm{G}}

\newcommand \rA {\mathrm{A}}
\newcommand \rD {\mathrm{D}}
\newcommand \rc {\mathrm{c}}
\newcommand \rE {\mathrm{E}}
\newcommand \rG {\mathrm{G}}
\newcommand \rH {\mathrm{H}}
\newcommand \rI {\mathrm{I}}
\newcommand \rJ {\mathrm{J}}
\newcommand \rK {\mathrm{K}}
\newcommand \rL {\mathrm{L}}
\newcommand \rM {\mathrm{M}}
\newcommand \rN {\mathrm{N}}
\newcommand \rP {\mathrm{P}}
\newcommand \rQ {\mathrm{Q}}
\newcommand \rS {\mathrm{S}}
\newcommand \rT {\mathrm{T}}
\newcommand \rZ {\mathrm{Z}}
\newcommand \rd {\mathrm{d}}

\newcommand \ini {{\mathrm{ini}_\preceq}}

%:   mathbb
\renewcommand \AA{\mathbb{A}}
\newcommand \BB{\mathbb{B}}
\newcommand \CC{\mathbb{C}}
\newcommand \FF{\mathbb{F}}
\newcommand \FFp{\FF_p}
\newcommand \FFq{\FF_q}
\newcommand \GG{\mathbb{G}}
\newcommand \MM{\mathbb{M}}
\newcommand \NN{\mathbb{N}}
\newcommand \PP{\mathbb{P}}
\newcommand \QQ{\mathbb{Q}}
\newcommand \UU{\mathbb{U}}
\newcommand \ZZ{\mathbb{Z}}
\newcommand \ZB{\mathbb{ZB}}
\newcommand \RR{\mathbb{R}}
\newcommand \ASL {\mathbb{ASL}}
\newcommand \AGL {\mathbb{AGL}}

\newcommand \Z{\mathbb{Z}} % a modifier eventuellement voir sous anneau premier chap3

%:   mathcal
\newcommand \cA {\mathcal{A}}
\newcommand \cB {\mathcal{B}}
\newcommand \cC {\mathcal{C}}
\newcommand \cD {\mathcal{D}}
\newcommand \cE {\mathcal{E}}
\newcommand \Diff {\mathcal{D}}
\newcommand \cG {\mathcal{G}}
\newcommand \cI {\mathcal{I}}
\newcommand \cJ {\mathcal{J}}
\newcommand \cF {\mathcal{F}}
\newcommand \cK {\mathcal{K}}
\newcommand \cL {\mathcal{L}}
\newcommand \cO {\mathcal{O}}
\newcommand \cP {\mathcal{P}}
\newcommand \cR {\mathcal{R}}
\newcommand \cM {\mathcal{M}}
\newcommand \cN {\mathcal{N}}
\newcommand \cS {\mathcal{S}}
\newcommand \cV {\mathcal{V}}
\newcommand \cZ {\mathcal{Z}}

\newcommand \SK {\cS^\rK}
\newcommand \IK {\cI^\rK}
\newcommand \JK {\cJ^\rK}
\newcommand \IH {\cI^\rH}
\newcommand \JH {\cJ^\rH}

%:   mathbf
\newcommand \ga{\mathbf{a}}
\newcommand \gb{\mathbf{b}}
\newcommand \gc{\mathbf{c}}
\newcommand \bfe{\mathbf{e}}  %% pas \ge !!
\newcommand \gh{\mathbf{h}}
\newcommand \gk{\mathbf{k}}
\newcommand \gl{\mathbf{l}}
\newcommand \gs{\mathbf{s}}
\newcommand \gv{\mathbf{v}}
\newcommand \gw{\mathbf{w}}
\newcommand \gA{\mathbf{A}}
\newcommand \gB{\mathbf{B}}
\newcommand \gC{\mathbf{C}}
\newcommand \gD{\mathbf{D}}
\newcommand \gE{\mathbf{E}}
\newcommand \gF{\mathbf{F}}
\newcommand \gG{\mathbf{G}}
\newcommand \gK{\mathbf{K}}
\newcommand \gL{\mathbf{L}}
\newcommand \gM{\mathbf{M}}
\newcommand \gQ{\mathbf{Q}}
\newcommand \gR{\mathbf{R}}
\newcommand \gS{\mathbf{S}}
\newcommand \gT{\mathbf{T}}
\newcommand \gU{\mathbf{U}}
\newcommand \gV{\mathbf{V}}
\newcommand \gx{\mathbf{x}}
\newcommand \gy{\mathbf{y}}
\newcommand \gX{\mathbf{X}}
\newcommand \gW{\mathbf{W}}
\newcommand \gZ{\mathbf{Z}}

%:   mathfrak
\newcommand\fa{\mathfrak{a}}
\newcommand\fb{\mathfrak{b}}
\newcommand\fc{\mathfrak{c}}
\newcommand\fd{\mathfrak{d}}
\newcommand\fA{\mathfrak{A}}
\newcommand\fB{\mathfrak{B}}
\newcommand\fC{\mathfrak{C}}
\newcommand\fD{\mathfrak{D}}
\newcommand\fI{\mathfrak{i}}
\newcommand\fII{\mathfrak{I}}
\newcommand\fj{\mathfrak{j}}
\newcommand\fJ{\mathfrak{J}}
\newcommand\ff{\mathfrak{f}}
\newcommand\ffg{\mathfrak{g}}
\newcommand\fF{\mathfrak{F}}
\newcommand\fh{\mathfrak{h}}
\newcommand\fl{\mathfrak{l}}
\newcommand\fm{\mathfrak{m}}
\newcommand\fp{\mathfrak{p}}
\newcommand\fq{\mathfrak{q}}
\newcommand\fM{\mathfrak{M}}
\newcommand\fN{\mathfrak{N}}
\newcommand\fP{\mathfrak{P}}
\newcommand\fQ{\mathfrak{Q}}
\newcommand\fR{\mathfrak{R}}
\newcommand\fx{\mathfrak{x}}
\newcommand\fV{\mathfrak{V}}
\newcommand\fZ{\mathfrak{Z}}

%:   mots mathematiques roman

\newcommand \Ig {\mathrm{Ig}}
\newcommand \Id {\mathrm{Id}}
\newcommand \In {{\rI_n}}
\newcommand \J {\MA{\mathrm{Jac}}}
\newcommand \JJ {\MA{\mathrm{JAC}}}
\newcommand \Sn {{\mathrm{S}_n}}

\newcommand \DA {\rD_{\!\gA}}
\newcommand \DB {\rD_\gB}
\newcommand \DV {\rD_\gV}
\newcommand \JA {\rJ_\gA}

\newcommand \red {_{\mathrm{red}}}
\newcommand \qim {_{\mathrm{min}}}
\newcommand \rja {\mathrm{Ja}}

\newcommand \ide {\mathrm{e}}

\newcommand \Adj {\MA{\mathrm{Adj}}}
\newcommand \adj {\MA{\mathrm{adj}}}
\newcommand \Adu {\MA{\mathrm{Adu}}}
\newcommand \Ann {\mathrm{Ann}}
\newcommand \Aut {\MA{\mathrm{Aut}}}
\newcommand \BZ {\mathrm{BZ}}
\newcommand \car {\MA{\mathrm{car}}}
\newcommand \Cl {\MA{\mathrm{Cl}}}
\newcommand \Coker {\MA{\mathrm{Coker}}}
\renewcommand \det {\MA{\mathrm{det}}}
\renewcommand \deg {\MA{\mathrm{deg}}}
\newcommand \Diag {\MA{\mathrm{Diag}}}
\newcommand \di {\MA{\mathrm{di}}}
\newcommand \disc {\MA{\mathrm{disc}}}
\newcommand \Disc {\MA{\mathrm{Disc}}}
\newcommand \Div {\MA{\mathrm{Div}}}
\newcommand \dv {\MA{\mathrm{div}}}
\newcommand \ev {{\mathrm{ev}}}
\newcommand \End {\MA{\mathrm{End}}}
\newcommand \Fix {\MA{\mathrm{Fix}}}
\newcommand \Frac {\MA{\mathrm{Frac}}}
\newcommand \gr {\MA{\mathrm{gr}}}
\newcommand \Gal {\MA{\mathrm{Gal}}}
\newcommand \Gfr {\MA{\mathrm{Gfr}}}
\newcommand \Gram {\MA{\mathrm{Gram}}}
\newcommand \gram {\MA{\mathrm{gram}}}
\newcommand \hauteur {\mathrm{hauteur}}
\renewcommand \ker {\Ker}
\newcommand \Ker {\MA{\mathrm{Ker}}}
\newcommand \Hom {\MA{\mathrm{Hom}}}
\newcommand \Ifr {\MA{\mathrm{Ifr}}}
\newcommand \Icl {\MA{\mathrm{Icl}}}
\renewcommand \Im {\MA{\mathrm{Im}}}
\newcommand \Lin {\mathrm{L}}
\newcommand \LIN {\mathrm{Lin}}
\newcommand \Mat {\MA{\mathrm{Mat}}}
\newcommand \Mip {\mathrm{Min}}
\newcommand \md {\mathrm{md}}
\newcommand \mod {\;\mathrm{mod}\;}
\newcommand \Mor {\MA{\mathrm{Mor}}}
\newcommand \poles {\hbox {\rm p\^oles}}
\newcommand \pgcd {\MA{\mathrm{pgcd}}}
\newcommand \ppcm {\MA{\mathrm{ppcm}}}
\newcommand \Rad {\MA{\mathrm{Rad}}}
\newcommand \Reg {\MA{\mathrm{Reg}}}
\newcommand \rg{\MA{\mathrm{rg}}}
\newcommand \Res {\mathrm{Res}}
\newcommand \Rs {\MA{\mathrm{Rs}}}
\newcommand \rPr{\MA{\mathrm{Pr}}}
\newcommand \Rv {\mathrm{Rv}}
\newcommand \Syl {\mathrm{Syl}}
\newcommand \Stp {\MA{\mathrm{Stp}}}
\newcommand \St {\mathrm{St}}
\newcommand \Tri {\MA{\mathrm{Tri}}}
\newcommand \Tor {\MA{\mathrm{Tor}}}
\newcommand \tr {\MA{\mathrm{tr}}}
\newcommand \Tr {\MA{\mathrm{Tr}}}
\newcommand \Tsc {\MA{\mathrm{Tsch}}}
\newcommand \Um {\MA{\mathrm{Um}}}
\newcommand \val {\MA{\mathrm{val}}}
\newcommand \Suslin{{\rm Suslin}}

\newcommand \rgl {\rg^\lambda}
\newcommand \rgg {\rg^\gamma}

%:  gras math avec \bm 
\newcommand{\bma}{\bm{a}}
\newcommand{\bmb}{\bm{b}}
\newcommand{\bmc}{\bm{c}}
\newcommand{\bmd}{\bm{d}}
\newcommand{\bme}{\bm{e}}
\newcommand{\bmf}{\bm{f}}
\newcommand{\bmu}{\bm{u}}
\newcommand{\bmv}{\bm{v}}
\newcommand{\bmw}{\bm{w}}
\newcommand{\bmy}{\bm{y}}
\newcommand{\bmx}{\bm{x}}
\newcommand{\bmz}{\bm{z}}

%:   mots mathematiques  sf

\newcommand \sfP {\mathsf{P}}
\newcommand \sfC {\mathsf{C}}

\newcommand \Gr {\MA{\mathsf{Gr}}}
\newcommand \GK {\MA{\mathsf{GK}}}
\newcommand \GKO {\MA{\mathsf{GK}_0}}
\newcommand \HO {\MA{\mathsf{H}_0}}
\newcommand \HOp {\MA{\mathsf{H}_0^+}}
\newcommand \Hdim {\MA{\mathsf{Hdim}}}
\newcommand \HeA {{\Heit\gA}}
\newcommand \Heit {\MA{\mathsf{Heit}}}
\newcommand \Hspec {\MA{\mathsf{Hspec}}}
\newcommand \Jdim {\MA{\mathsf{Jdim}}}
\newcommand \jdim {\MA{\mathsf{jdim}}}
\newcommand \Jspec {\MA{\mathsf{Jspec}}}
\newcommand \jspec {\MA{\mathsf{jspec}}}
\newcommand \KO {\MA{\mathsf{K}_0}}
\newcommand \KOp {\MA{\mathsf{K}_0^+}}
\newcommand \KTO {\MA{\wi{\mathsf{K}}_0}}
\newcommand \Kdim {\MA{\mathsf{Kdim}}}
\newcommand \Max {\MA{\mathsf{Max}}}
\newcommand \Min {\MA{\mathsf{Min}}}
\newcommand \OQC {\MA{\mathsf{Oqc}}}
\newcommand \Pic {\MA{\mathsf{Pic}}}
\newcommand \Spec {\MA{\mathsf{Spec}}}
\newcommand \SpecA {\Spec\gA}
\newcommand \SpecT {\Spec\gT}
\newcommand \Vdim {\MA{\mathsf{Vdim}}}
\newcommand \Zar {\MA{\mathsf{Zar}}}
\newcommand \ZF {\MA{\mathsf{ZF}}}
\newcommand \ZarA {{\Zar\gA}}

%:   mots mathematiques  frak

\newcommand \fRes {\MA{\mathfrak{Res}}}

\newcommand \SPEC {\MA{\mathfrak{spec}}}
\newcommand \SPECK {\SPEC_\gK}
\newcommand \SPECk {\SPEC_\gk}

\newcommand \Ap {{\gA_\fp}}
\newcommand \zg {\ZZ[G]}

\newcommand \num {{n$^{\mathrm{ o}}$}}
%\newcommand \num {{n\o}}

%:  Commentaires, remarques, exemples, problemes
\newcommand\comm{\rdb
\noi{\it Commentaire. }}

\newcommand\COM[1]{\rdb
\noi{\it Commentaire #1. }}

\newcommand\comms{\rdb
\noi{\it Commentaires. }}

\newcommand\rem{\rdb
\noi{\it Remarque. }}

\newcommand\REM[1]{\rdb
\noi{\it Remarque#1. }}

\newcommand\rems{\rdb
\noi{\it Remarques. }}

\newcommand\exl{\rdb
\noi{\bf Exemple. }}

\newcommand\EXL[1]{\rdb
\noi{\bf Exemple: #1. }}

\newcommand\exls{\rdb
\noi{\bf Exemples. }}

\newcommand\Pb{\rdb
\noi{\bf Problème. }}

\newcommand\PB[1]{\rdb
\noi{\bf Problème #1. }}

\newcommand \eoe {\hbox{}\nobreak\hfill
\vrule width 1.4mm height 1.4mm depth 0mm \par \smallskip}

\newcommand\eoq{\hbox{}\nobreak
\vrule width 1.4mm height 1.4mm depth 0mm}

%:   cad hdr spdg propeq ...
\newcommand \recu {récur\-rence }
\newcommand \recuz {récur\-rence}
\newcommand \hdr {hypo\-thèse de \recu }
\newcommand \hdrz {hypo\-thèse de \recuz}
\newcommand \Cad {C'est-\`a-dire }
\newcommand \cad {c'est-\`a-dire }
\newcommand \cadz {c'est-\`a-dire}
\newcommand \cade {c'est-\`a-dire en\-co\-re }
\newcommand \ssi {si, et seu\-le\-ment si, }
\newcommand \ssiz {si, et seu\-le\-ment si,~}
\newcommand \cnes {con\-di\-tion néces\-sai\-re et suf\-fi\-san\-te }
\newcommand \spdg {sans per\-te de géné\-ra\-lité }
\newcommand \spdgz {sans per\-te de géné\-ra\-lité}
\newcommand \Spdg {Sans per\-te de géné\-ra\-lité }
\newcommand \Spdgz {Sans per\-te de géné\-ra\-lité}

\newcommand \Propeq {Les pro\-pri\-é\-tés sui\-van\-tes sont
équi\-va\-len\-tes.}
\newcommand \propeq {les pro\-pri\-é\-tés sui\-van\-tes sont
équi\-va\-len\-tes.}

%:  Kev Alg etc...
\newcommand \Kev {$\gK$-\evc }
\newcommand \Kevs {$\gK$-\evcs }
\newcommand \Kevz {$\gK$-\evcz}
\newcommand \Kevsz {$\gK$-\evcsz}

\newcommand \Lev {$\gL$-\evc }
\newcommand \Levs {$\gL$-\evcs }
\newcommand \Levz {$\gL$-\evcz}
\newcommand \Levsz {$\gL$-\evcsz}

\newcommand \Qev {$\QQ$-\evc }
\newcommand \Qevs {$\QQ$-\evcs }
\newcommand \Qevz {$\QQ$-\evcz}
\newcommand \Qevsz {$\QQ$-\evcsz}

\newcommand \kev {$\gk$-\evc }
\newcommand \kevs {$\gk$-\evcs }
\newcommand \kevz {$\gk$-\evcz}
\newcommand \kevsz {$\gk$-\evcsz}

\newcommand \lev {$\gl$-\evc }
\newcommand \levs {$\gl$-\evcs }
\newcommand \levz {$\gl$-\evcz}
\newcommand \levsz {$\gl$-\evcsz}

\newcommand \Alg {$\gA$-\alg}
\newcommand \Algs {$\gA$-\algs}
\newcommand \Algz {$\gA$-\algz}
\newcommand \Algsz {$\gA$-\algsz}

\newcommand \Blg {$\gB$-\alg}
\newcommand \Blgs {$\gB$-\algs}
\newcommand \Blgz {$\gB$-\algz}
\newcommand \Blgsz {$\gB$-\algsz}

\newcommand \Clg {$\gC$-\alg}
\newcommand \Clgs {$\gC$-\algs}
\newcommand \Clgz {$\gC$-\algz}
\newcommand \Clgsz {$\gC$-\algsz}

\newcommand \klg {$\gk$-\alg}
\newcommand \klgs {$\gk$-\algs}
\newcommand \klgz {$\gk$-\algz}
\newcommand \klgsz {$\gk$-\algsz}

\newcommand \llg {$\gl$-\alg}
\newcommand \llgs {$\gl$-\algs}
\newcommand \llgz {$\gl$-\algz}
\newcommand \llgsz {$\gl$-\algsz}

\newcommand \Klg {$\gK$-\alg}
\newcommand \Klgs {$\gK$-\algs}
\newcommand \Klgz {$\gK$-\algz}
\newcommand \Klgsz {$\gK$-\algsz}

\newcommand \Llg {$\gL$-\alg}
\newcommand \Llgs {$\gL$-\algs}
\newcommand \Llgz {$\gL$-\algz}
\newcommand \Llgsz {$\gL$-\algsz}

\newcommand \QQlg {$\QQ$-\alg}
\newcommand \QQlgs {$\QQ$-\algs}
\newcommand \QQlgz {$\QQ$-\algz}
\newcommand \QQlgsz {$\QQ$-\algsz}

\newcommand \Rlg {$\gR$-\alg}
\newcommand \Rlgs {$\gR$-\algs}
\newcommand \Rlgz {$\gR$-\algz}
\newcommand \Rlgsz {$\gR$-\algsz}

\newcommand \RRlg {$\RR$-\alg}
\newcommand \RRlgs {$\RR$-\algs}
\newcommand \RRlgz {$\RR$-\algz}
\newcommand \RRlgsz {$\RR$-\algsz}

\newcommand \ZZlg {$\ZZ$-\alg}
\newcommand \ZZlgs {$\ZZ$-\algs}
\newcommand \ZZlgz {$\ZZ$-\algz}
\newcommand \ZZlgsz {$\ZZ$-\algsz}

\newcommand \Amo {$\gA$-mo\-du\-le }
\newcommand \Amos {$\gA$-mo\-du\-les }
\newcommand \Amoz {$\gA$-mo\-du\-le}
\newcommand \Amosz {$\gA$-mo\-du\-les}

\newcommand \Amrc {\Amo \prc}
\newcommand \Amrcs {\Amos \prcs}
\newcommand \Amrcz {\Amo \prcz}
\newcommand \Amrcsz {\Amos \prcsz}

\newcommand \ZZmo {$\ZZ$-mo\-du\-le }
\newcommand \ZZmos {$\ZZ$-mo\-du\-les }
\newcommand \ZZmoz {$\ZZ$-mo\-du\-le}
\newcommand \ZZmosz {$\ZZ$-mo\-du\-les}

\newcommand \Bmo {$\gB$-mo\-du\-le }
\newcommand \Bmos {$\gB$-mo\-du\-les }
\newcommand \Bmoz {$\gB$-mo\-du\-le}
\newcommand \Bmosz {$\gB$-mo\-du\-les}

\newcommand \Cmo {$\gC$-mo\-du\-le }
\newcommand \Cmos {$\gC$-mo\-du\-les }
\newcommand \Cmoz {$\gC$-mo\-du\-le}
\newcommand \Cmosz {$\gC$-mo\-du\-les}

\newcommand \Dmo {$\gD$-mo\-du\-le }
\newcommand \Dmos {$\gD$-mo\-du\-les }
\newcommand \Dmoz {$\gD$-mo\-du\-le}
\newcommand \Dmosz {$\gD$-mo\-du\-les}

\newcommand \kmo {$\gk$-mo\-du\-le }
\newcommand \kmos {$\gk$-mo\-du\-les }
\newcommand \kmoz {$\gk$-mo\-du\-le}
\newcommand \kmosz {$\gk$-mo\-du\-les}

\newcommand \Kmo {$\gK$-mo\-du\-le }
\newcommand \Kmos {$\gK$-mo\-du\-les }
\newcommand \Kmoz {$\gK$-mo\-du\-le}
\newcommand \Kmosz {$\gK$-mo\-du\-les}

\newcommand \kmrc {\kmo \prc}
\newcommand \kmrcs {\kmos \prcs}
\newcommand \kmrcz {\kmo \prcz}
\newcommand \kmrcsz {\kmos \prcsz}

\newcommand \Lmo {$\gL$-mo\-du\-le }
\newcommand \Lmos {$\gL$-mo\-du\-les }
\newcommand \Lmoz {$\gL$-mo\-du\-le}
\newcommand \Lmosz {$\gL$-mo\-du\-les}

\newcommand \Ali {\apl $\gA$-\lin }
\newcommand \Alis {\apls $\gA$-\lins }
\newcommand \Aliz {\apl $\gA$-\linz}
\newcommand \Alisz {\apls $\gA$-\linsz}

\newcommand \kli {\apl $\gk$-\lin }
\newcommand \klis {\apls $\gk$-\lins }
\newcommand \kliz {\apl $\gk$-\linz}
\newcommand \klisz {\apls $\gk$-\linsz}

\newcommand \Kli {\apl $\gK$-\lin }
\newcommand \Klis {\apls $\gK$-\lins }
\newcommand \Kliz {\apl $\gK$-\linz}
\newcommand \Klisz {\apls $\gK$-\linsz}

\newcommand \Bli {\apl $\gB$-\lin }
\newcommand \Blis {\apls $\gB$-\lins }
\newcommand \Bliz {\apl $\gB$-\linz}
\newcommand \Blisz {\apls $\gB$-\linsz}

\newcommand \Cli {\apl $\gC$-\lin }
\newcommand \Clis {\apls $\gC$-\lins }
\newcommand \Cliz {\apl $\gC$-\linz}
\newcommand \Clisz {\apls $\gC$-\linsz}

%:  abbréviations a

\newcommand \abi {\apl \bil}
\newcommand \abis {\apls \bils}
\newcommand \abiz {\apl \bilz}
\newcommand \abisz {\apls \bilsz}

\newcommand \ac{\agqt clos }
\newcommand \acz{\agqt clos}

\newcommand \acl {an\-neau \icl}
\newcommand \acls {an\-neaux \icl}
\newcommand \aclsz {an\-neaux \iclz}
\newcommand \aclz {an\-neau \iclz}

\newcommand \adk {an\-neau de Dedekind }
\newcommand \adks {an\-neaux de Dedekind }
\newcommand \adkz {an\-neau de Dedekind}
\newcommand \adksz {an\-neaux de Dedekind}

\newcommand \adp {an\-neau de Pr\"u\-fer }
\newcommand \adps {an\-neaux de Pr\"u\-fer }
\newcommand \adpsz {an\-neaux de Pr\"u\-fer}
\newcommand \adpz {an\-neau de Pr\"u\-fer}

\newcommand \adpc {\adp \coh }
\newcommand \adpcs {\adps \cohs }
\newcommand \adpcz {\adp \cohz}
\newcommand \adpcsz {\adps \cohsz}

\newcommand \adu {\alg de \dcn \uvle }
\newcommand \adus {\algs de \dcn \uvle }
\newcommand \aduz {\alg de \dcn \uvlez}
\newcommand \adusz {\algs de \dcn \uvlez}

\newcommand \adv {an\-neau de va\-lu\-ation }
\newcommand \advs {an\-neaux de va\-lu\-ation }
\newcommand \advsz {an\-neaux de va\-lu\-ation}
\newcommand \advz {an\-neau de va\-lu\-ation}

\newcommand \advd {\adv discrète }
\newcommand \advds {\advs discrète }
\newcommand \advdz {\adv discrète}
\newcommand \advdsz {\advs discrète}

\newcommand \aG {\alg galoisienne }
\newcommand \aGs {\algs galoisiennes }
\newcommand \aGz {\alg galoisienne}
\newcommand \aGsz {\algs galoisiennes}

\newcommand \AgB {Algèbre de Boole }
\newcommand \agB {\alg de Boole }
\newcommand \agBs {\algs de Boole }
\newcommand \agBz {\alg de Boole}
\newcommand \agBsz {\algs de Boole}

\newcommand \agH {\alg de Heyting }
\newcommand \agHs {\algs de Heyting }
\newcommand \agHz {\alg de Heyting}
\newcommand \agHsz {\algs de Heyting}

\newcommand \agsp {\agq \spl }
\newcommand \agsps {\agqs \spls }
\newcommand \agspz {\agq \splz}
\newcommand \agspsz {\agqs \splsz}

\newcommand \agq{algé\-brique }
\newcommand \agqs{algé\-briques }
\newcommand \agqz{algé\-brique}
\newcommand \agqsz{algé\-briques}

\newcommand \agqt{\agqz\-ment }

\newcommand \alg {algèbre }
\newcommand \algs {algèbres }
\newcommand \algz {algèbre}
\newcommand \algsz {algèbres}

\newcommand \Algo{Algo\-rithme }
\newcommand \algo{algo\-rithme }
\newcommand \algos{algo\-rithmes }
\newcommand \algoz{algo\-rithme}
\newcommand \algosz{algo\-rithmes}

\newcommand \algq{algo\-rith\-mique }
\newcommand \algqs{algo\-rith\-miques }
\newcommand \algqz{algo\-rith\-mique}
\newcommand \algqsz{algo\-rith\-miques}

\newcommand \algqt{algo\-rith\-mi\-quement }

\newcommand \ali {appli\-ca\-tion \lin }
\newcommand \alis {appli\-ca\-tions \lins }
\newcommand \aliz {appli\-ca\-tion \linz}
\newcommand \alisz {appli\-ca\-tions \linsz}

\newcommand \alo {an\-neau lo\-cal }
\newcommand \aloz {an\-neau lo\-cal}
\newcommand \alos {an\-neaux lo\-caux }
\newcommand \alosz {an\-neaux lo\-caux}

\newcommand \algb {an\-neau \lgb}
\newcommand \algbs {an\-neaux \lgbs}
\newcommand \algbz {an\-neau \lgbz}
\newcommand \algbsz {an\-neaux \lgbsz}

\newcommand \alrd {\alo \dcd }
\newcommand \alrds {\alos \dcds }
\newcommand \alrdz {\alo \dcdz}
\newcommand \alrdsz {\alos \dcdsz}

\newcommand \anar {anneau \ari}
\newcommand \anars {anneaux \aris}
\newcommand \anarsz {anneaux \arisz}
\newcommand \anarz {anneau \ariz}

\newcommand \anor {an\-neau nor\-mal }
\newcommand \anors {an\-neaux nor\-maux }
\newcommand \anorsz {an\-neaux nor\-maux}
\newcommand \anorz {an\-neau nor\-mal}

\newcommand \apf {\alg \pf}
\newcommand \apfs {\algs \pf}
\newcommand \apfz {\alg \pfz}
\newcommand \apfsz {\algs \pfz}

\newcommand \apG {\alg pré\-galoisienne }
\newcommand \apGs {\algs pré\-galoisiennes }
\newcommand \apGz {\alg pré\-galoisienne}
\newcommand \apGsz {\algs pré\-galoisiennes}

\newcommand \apl {appli\-cation }
\newcommand \apls {appli\-cations }
\newcommand \aplz {appli\-cation}
\newcommand \aplsz {appli\-cations}

\newcommand \areg{\apl régu\-lière }
\newcommand \aregs{\apls régu\-lières }
\newcommand \aregz{\apl régu\-lière}
\newcommand \aregsz{\apls régu\-lières}

\newcommand \ari{arith\-mé\-tique }
\newcommand \ariz{arith\-mé\-tique}
\newcommand \aris{arith\-mé\-tiques }
\newcommand \arisz{arith\-mé\-tiques}

\newcommand \ase {\alg \ste}
\newcommand \ases {\algs \stes}
\newcommand \asez {\alg \stez}
\newcommand \asesz {\algs \stesz}

\newcommand \asf {\alg \stfe}
\newcommand \asfs {\algs \stfes}
\newcommand \asfz {\alg \stfez}
\newcommand \asfsz {\algs \stfesz}

\newcommand \atf {\alg \tf}
\newcommand \atfs {\algs \tf}
\newcommand \atfz {\alg \tfz}
\newcommand \atfsz {\algs \tfz}

\newcommand \auto {auto\-mor\-phis\-me }
\newcommand \autos {auto\-mor\-phis\-mes }
\newcommand \autoz {auto\-mor\-phis\-me}
\newcommand \autosz {auto\-mor\-phis\-mes}

%:  abbréviations b c

\newcommand \bb {bien bon }
\newcommand \bbs {bien bons }
\newcommand \bbz {bien bon}

\newcommand \bdp {base de \dcn partielle }
\newcommand \bdpz {base de \dcn partielle}

\newcommand \bdf {base de \fap }
\newcommand \bdfz {base de \fapz}

\newcommand \bdg {base de Gr\"obner }
\newcommand \bdgs {bases de Gr\"obner }
\newcommand \bdgz {base de Gr\"obner}
\newcommand \bdgsz {bases de Gr\"obner}

\newcommand \bil {bili\-né\-aire }
\newcommand \bils {bili\-né\-aires }
\newcommand \bilz {bili\-né\-aire}
\newcommand \bilsz {bili\-né\-aires}

\newcommand \cac{corps \ac }
\newcommand \cacz{corps \acz}

\newcommand \calf{calcul formel }
\newcommand \calfz{calcul formel}

\newcommand \cara{carac\-té\-ris\-tique }
\newcommand \caras{carac\-té\-ris\-tiques }
\newcommand \caraz{carac\-té\-ris\-tique}
\newcommand \carasz{carac\-té\-ris\-tiques}

\newcommand \care{carac\-té\-risé }
\newcommand \caree{carac\-té\-risée }
\newcommand \cares{carac\-té\-risés }
\newcommand \carees{carac\-té\-risées }

\newcommand \Carn{Carac\-té\-ri\-sation }
\newcommand \carn{carac\-té\-ri\-sation }
\newcommand \carnz{carac\-té\-ri\-sation}
\newcommand \Carns{Carac\-té\-ri\-sations }
\newcommand \carns{carac\-té\-ri\-sations }

\newcommand \carar{carac\-té\-riser }

\newcommand \carf{de carac\-tère fini }
\newcommand \carfz{de carac\-tère fini}

\newcommand \cdac{\cdi \ac }
\newcommand \cdacs{\cdis \ac }
\newcommand \cdacz{\cdi \acz}
\newcommand \cdacsz{\cdis \acz}

\newcommand \cdi{corps discret }
\newcommand \cdis{corps discrets }
\newcommand \cdiz{corps discret}
\newcommand \cdisz{corps discrets}

\newcommand \cdr{corps de racines }
\newcommand \cdrz{corps de racines}

\newcommand \cdv{changement de variables }
\newcommand \cdvs{changements de variables }
\newcommand \cdvz{changement de variables}
\newcommand \cdvsz{changements de variables}

\newcommand \cli{cl\^o\-ture inté\-grale }
\newcommand \cliz{cl\^o\-ture inté\-grale}

\newcommand \cpb {compa\-tible } 
\newcommand \cpbs {compa\-tibles } 
\newcommand \cpbz {compa\-tible} 
\newcommand \cpbsz {compa\-tibles} 

\newcommand \cpbt {compa\-tibi\-lité } 
\newcommand \cpbtz {compa\-tibi\-lité} 

\newcommand \coe {coef\-fi\-cient }
\newcommand \coes {coef\-fi\-cients }
\newcommand \coez {coef\-fi\-cient}
\newcommand \coesz {coef\-fi\-cients}

\newcommand \coh {cohé\-rent }
\newcommand \cohs {cohé\-rents }
\newcommand \cohz {cohé\-rent}
\newcommand \cohsz {cohé\-rents}

\newcommand \cohc {co\-hé\-rence }
\newcommand \cohcz {co\-hé\-rence}

\newcommand \coli {combi\-naison \lin }
\newcommand \colis {combi\-naisons \lins }
\newcommand \coliz {combi\-naison \linz}
\newcommand \colisz {combi\-naisons \linsz}

\newcommand \com {co\-maxi\-maux }
\newcommand \comz {co\-maxi\-maux}
\newcommand \come {co\-maxi\-ma\-les }
\newcommand \comez {co\-maxi\-ma\-les}

\newcommand \como {module \coh}
\newcommand \comoz {module \cohz}
\newcommand \comos {modules \cohs}
\newcommand \comosz {modules \cohsz}

\newcommand \coo {coor\-don\-née }
\newcommand \coos {coor\-don\-nées }
\newcommand \cooz {coor\-don\-née}
\newcommand \coosz {coor\-don\-nées}

\newcommand \cop {complé\-men\-taire }
\newcommand \cops {complé\-men\-taires }
\newcommand \copz {complé\-men\-taire}
\newcommand \copsz {complé\-men\-taires}

\newcommand \cor {coré\-guliers }
\newcommand \core {coré\-gulières }
\newcommand \corez {coré\-gulières}
\newcommand \corz {coré\-guliers}

\newcommand \cori {anneau \coh}
\newcommand \coriz {anneau \cohz}
\newcommand \coris {anneaux \cohs}
\newcommand \corisz {anneaux \cohsz}

\newcommand \crc{carac\-tère }
\newcommand \crcs{carac\-tères }
\newcommand \crcz{carac\-tère}
\newcommand \crcsz{carac\-tères}

\newcommand \cro {croi\-sée }
\newcommand \cros {croi\-sées }
\newcommand \croz {croi\-sée}
\newcommand \crosz {croi\-sées}

%:  abbréviations d

\newcommand \dcd {\rdt dis\-cret }
\newcommand \dcdz {\rdt dis\-cret}
\newcommand \dcds {\rdt dis\-crets }
\newcommand \dcdsz {\rdt dis\-crets}

\newcommand \dcn {dé\-com\-po\-sition }
\newcommand \dcns {dé\-com\-po\-sitions }
\newcommand \dcnz {dé\-com\-po\-sition}
\newcommand \dcnsz {dé\-com\-po\-sitions}

\newcommand \dcnb {\dcn bornée }
\newcommand \dcnbz {\dcn bornée}

\newcommand \dcnc {\dcn complète }
\newcommand \dcncz {\dcn complète}

\newcommand \dcnp {\dcn partielle }
\newcommand \dcnps {\dcn partielles }
\newcommand \dcnpz {\dcn partielle}

\newcommand \dcp {décom\-posa\-ble }
\newcommand \dcps {décom\-posa\-bles }
\newcommand \dcpz {décom\-posa\-ble}
\newcommand \dcpsz {décom\-posa\-bles}

\newcommand \ddk {dimen\-sion de Krull }
\newcommand \ddkz {dimen\-sion de Krull}

\newcommand \dDk {do\-mai\-ne de Dedekind }
\newcommand \dDks {do\-mai\-nes de Dedekind }
\newcommand \dDkz {do\-mai\-ne de Dedekind}
\newcommand \dDksz {do\-mai\-nes de Dedekind}

\newcommand \ddi {de di\-men\-sion infé\-rieure ou égale \`a~}

\newcommand \ddp {domaine de Pr\"u\-fer }
\newcommand \ddps {domaines de Pr\"u\-fer }
\newcommand \ddpz {domaine de Pr\"u\-fer}
\newcommand \ddpsz {domaines de Pr\"u\-fer}

\newcommand \ddv {domaine de valu\-ation }

\newcommand \dgn {dégé\-néré }
\newcommand \dgne {dégé\-nérée }
\newcommand \dgns {dégé\-nérés }
\newcommand \dgnes {dégé\-nérées }
\newcommand \dgnz {dégé\-néré}
\newcommand \dgnez {dégé\-nérée}
\newcommand \dgnsz {dégé\-nérés}
\newcommand \dgnesz {dégé\-nérées}

\newcommand \Demo{Démons\-tra\-tion }

\newcommand \dem{démons\-tra\-tion }
\newcommand \demz{démons\-tra\-tion}
\newcommand \dems{démons\-tra\-tions }
\newcommand \demsz{démons\-tra\-tions}

\newcommand \denb{dénom\-brable }
\newcommand \denbs{dénom\-brables }
\newcommand \denbz{dénom\-brable}
\newcommand \denbsz{dénom\-brables}

\newcommand \deno{déno\-mi\-nateur }
\newcommand \denos{déno\-mi\-nateurs }
\newcommand \denoz{déno\-mi\-nateur}
\newcommand \denosz{déno\-mi\-nateurs}

\newcommand \deter {déter\-minant }
\newcommand \deters {déter\-minants }
\newcommand \deterz {déter\-minant}
\newcommand \detersz {déter\-minants}

\newcommand \dig{diago\-na\-li\-sable }
\newcommand \digs{diago\-na\-li\-sables }
\newcommand \digz{diago\-na\-li\-sable}
\newcommand \digsz{diago\-na\-li\-sables}

\newcommand \dil{diffé\-rentiel }
\newcommand \dils{diffé\-rentiels }
\newcommand \dile{diffé\-ren\-tielle }
\newcommand \diles{diffé\-ren\-tielles }
\newcommand \dilz{diffé\-rentiel}
\newcommand \dilsz{diffé\-rentiels}
\newcommand \dilez{diffé\-ren\-tielle}
\newcommand \dilesz{diffé\-ren\-tielles}

\newcommand \din{diago\-na\-li\-sation }
\newcommand \dinz{diago\-na\-li\-sation}
\newcommand \dins{diago\-na\-li\-sations }
\newcommand \dinsz{diago\-na\-li\-sations}

\newcommand \dit{distri\-bu\-ti\-vité }
\newcommand \ditz{distri\-bu\-ti\-vité}

\newcommand \iDKM {\index{Dedekind-Mertens}}
\newcommand \DKM {\iDKM Dedekind-Mertens }
\newcommand \DKMz {\iDKM Dedekind-Mertens}

\newcommand \dfn{défi\-nition }
\newcommand \dfns{défi\-nitions }
\newcommand \dfnz{défi\-nition}
\newcommand \dfnsz{défi\-nitions}

\newcommand \dlg{d'élar\-gis\-sement }
\newcommand \dlgz{d'élar\-gis\-sement}

\newcommand \discri{discri\-minant }
\newcommand \discris{discri\-minants }
\newcommand \discriz{discri\-minant}
\newcommand \discrisz{discri\-minants}

\newcommand \dok {domaine de~Dedekind }
\newcommand \doks {domaines de~Dedekind }
\newcommand \dokz {domaine de~Dedekind}
\newcommand \doksz {domaines de~Dedekind}

\newcommand \dvn {déri\-vation }
\newcommand \dvns {déri\-vations }
\newcommand \dvnz {déri\-vation}
\newcommand \dvnsz {déri\-vations}

\newcommand \dvz {di\-viseur de zé\-ro }
\newcommand \dvzs {di\-viseurs de zé\-ro }
\newcommand \dvzz {di\-viseur de zé\-ro}
\newcommand \dvzsz {di\-viseurs de zé\-ro}

\newcommand \dve {divi\-si\-bi\-lité }
\newcommand \dvez {divi\-si\-bi\-lité}

\newcommand \dvee {\`a \dve explicite }
\newcommand \dveez {\`a \dve explicite}

%:  abbréviations e

\newcommand \eco {\elts \com}
\newcommand \ecoz {\elts \comz}

\newcommand \Eds {Exten\-sion des sca\-laires }
\newcommand \Edsz {Exten\-sion des sca\-laires}
\newcommand \eds {exten\-sion des sca\-laires }
\newcommand \edsz {exten\-sion des sca\-laires}

\newcommand \egmt {éga\-le\-ment }
\newcommand \egmtz {éga\-le\-ment}

\newcommand \egt {éga\-lité }
\newcommand \egts {éga\-lités }
\newcommand \egtz {éga\-lité}
\newcommand \egtsz {éga\-lités}

\newcommand \eli{élimi\-nation }
\newcommand \eliz{élimi\-nation}

\newcommand \elr{élé\-men\-taire }
\newcommand \elrs{élé\-men\-taires }
\newcommand \elrz{élé\-men\-taire}
\newcommand \elrsz{élé\-men\-taires}

\newcommand \elrt{élé\-men\-tai\-rement }
\newcommand \elrtz{élé\-men\-tai\-rement}

\newcommand \elt{élé\-ment }
\newcommand \elts{élé\-ments }
\newcommand \eltz{élé\-ment}
\newcommand \eltsz{élé\-ments}

\def \endo {en\-do\-mor\-phisme }
\def \endos {en\-do\-mor\-phismes }
\def \endoz {en\-do\-mor\-phisme}
\def \endosz {en\-do\-mor\-phismes}

\newcommand \entrel {rela\-tion impli\-ca\-tive }
\newcommand \entrelz {rela\-tion impli\-ca\-tive}
\newcommand \entrels {rela\-tions impli\-ca\-tives }
\newcommand \entrelsz {rela\-tions impli\-ca\-tives}

\newcommand \enum {énu\-mé\-rable }
\newcommand \enums {énu\-mé\-rables }
\newcommand \enumz {énu\-mé\-rable}
\newcommand \enumsz {énu\-mé\-rables}

\newcommand \eqn{équa\-tion }
\newcommand \eqns{équa\-tions }
\newcommand \eqnz{équa\-tion}
\newcommand \eqnsz{équa\-tions}

\newcommand \eqv {équi\-valent }
\newcommand \eqve {équi\-valente }
\newcommand \eqvs {équi\-valents }
\newcommand \eqves {équi\-valentes }
\newcommand \eqvz {équi\-valent}
\newcommand \eqvez {équi\-valente}
\newcommand \eqvsz {équi\-valents}
\newcommand \eqvesz {équi\-valentes}

\newcommand \eqvc {équi\-va\-lence }
\newcommand \eqvcs {équi\-va\-lences }
\newcommand \eqvcz {équi\-va\-lence}
\newcommand \eqvcsz {équi\-va\-lences}

\newcommand \ecr {\elts \cor}
\newcommand \ecrz {\elts \corz}
\newcommand \Erg {$E$-\ndz }
\newcommand \Erge {$E$-\ndze }
\newcommand \Ergs {$E$-\ndzs }
\newcommand \Erges {$E$-\ndzes }
\newcommand \Ergz {$E$-\ndzz }
\newcommand \Ergez {$E$-\ndzez }
\newcommand \Ergsz {$E$-\ndzsz }
\newcommand \Ergesz {$E$-\ndzesz }

\newcommand \Frg {$F$-\ndz }
\newcommand \Frge {$F$-\ndze }
\newcommand \Frgs {$F$-\ndzs }
\newcommand \Frges {$F$-\ndzes }
\newcommand \Frgz {$F$-\ndzz }
\newcommand \Frgez {$F$-\ndzez }
\newcommand \Frgsz {$F$-\ndzsz }
\newcommand \Frgesz {$F$-\ndzesz }

\newcommand \Grg {$G$-\ndz }
\newcommand \Grge {$G$-\ndze }
\newcommand \Grgs {$G$-\ndzs }
\newcommand \Grges {$G$-\ndzes }
\newcommand \Grgz {$G$-\ndzz }
\newcommand \Grgez {$G$-\ndzez }
\newcommand \Grgsz {$G$-\ndzsz }
\newcommand \Grgesz {$G$-\ndzesz }

\newcommand \Hrg {$H$-\ndz }
\newcommand \Hrge {$H$-\ndze }
\newcommand \Hrgs {$H$-\ndzs }
\newcommand \Hrges {$H$-\ndzes }
\newcommand \Hrgz {$H$-\ndzz }
\newcommand \Hrgez {$H$-\ndzez }
\newcommand \Hrgsz {$H$-\ndzsz }
\newcommand \Hrgesz {$H$-\ndzesz }

\newcommand\evc{espa\-ce vecto\-riel }
\newcommand\evcs{espa\-ces vecto\-riels }
\newcommand\evcz{espa\-ce vecto\-riel}
\newcommand\evcsz{espa\-ces vecto\-riels}

\newcommand\evd{\evn dyna\-mique }
\newcommand\evdz{\evn dyna\-mique}

\newcommand\evn{éva\-lua\-tion }
\newcommand\evnz{éva\-lua\-tion}
\newcommand\evns{éva\-lua\-tions }
\newcommand\evnsz{éva\-lua\-tions}

%:  abbréviations f

\newcommand \fab {\fcn bornée }
\newcommand \fabz {\fcn bornée}

\newcommand \fac {\fcn totale }
\newcommand \facz {\fcn totale}

\newcommand \fap {\fcn par\-tiel\-le }
\newcommand \faps {\fcns par\-tiel\-les }
\newcommand \fapz {\fcn par\-tiel\-le}
\newcommand \fapsz {\fcns par\-tiel\-les}

\newcommand \fcn {facto\-ri\-sation }
\newcommand \fcns {facto\-ri\-sations }
\newcommand \fcnz {facto\-ri\-sation}
\newcommand \fcnsz {facto\-ri\-sations}

\newcommand \fit {fidè\-lement }

\newcommand \fip {filtre pre\-mier }
\newcommand \fips {filtres pre\-miers }
\newcommand \fipz {filtre pre\-mier}
\newcommand \fipsz {filtres pre\-miers}

\newcommand \fima {filtre maxi\-mal }
\newcommand \fimas {filtres maxi\-maux }
\newcommand \fimaz {filtre maxi\-mal}
\newcommand \fimazs {filtres maxi\-maux}

\newcommand \fdi {for\-te\-ment dis\-cret }
\newcommand \fdis {for\-te\-ment dis\-crets }
\newcommand \fdisz {for\-te\-ment dis\-crets}
\newcommand \fdiz {for\-te\-ment dis\-cret}

\newcommand \fmt {formel\-lement }

\newcommand \fnt {\fmt nette }
\newcommand \fnts {\fmt nettes }
\newcommand \fntz {\fmt nette}
\newcommand \fntsz {\fmt nettes}

\newcommand \fpt {\fit plat }
\newcommand \fpte {\fit plate }
\newcommand \fpts {\fit plats }
\newcommand \fptes {\fit plates }
\newcommand \fptz {\fit plat}
\newcommand \fptez {\fit plate}
\newcommand \fptsz {\fit plats}
\newcommand \fptesz {\fit plates}

\newcommand \frg{fonction régu\-lière }
\newcommand \frgs{fonctions régu\-lières }
\newcommand \frgz{fonction régu\-lière}
\newcommand \frgsz{fonctions régu\-lières}

\newcommand \ftr {forme trace }
\newcommand \ftrz {forme trace}

%:  abbréviations g

\newcommand\gaq{\gmt \agq}
\newcommand\gaqz{\gmt \agqz}

\newcommand\gmt{géo\-métrie }
\newcommand\gmts{géo\-métries }
\newcommand\gmtz{géo\-métrie}
\newcommand\gmtsz{géo\-métries}

\newcommand\gmq{géo\-mé\-trique }
\newcommand\gmqs{géo\-mé\-triques }
\newcommand\gmqz{géo\-mé\-trique}
\newcommand\gmqsz{géo\-mé\-triques}

\newcommand\gmqt{géo\-mé\-tri\-quement }
\newcommand\gmqtz{géo\-mé\-tri\-quement}

\newcommand \Gmqtz{Géo\-mé\-tri\-quement}

\newcommand\gne{gé\-né\-ra\-li\-sé }
\newcommand\gnee{gé\-né\-ra\-li\-sée }
\newcommand\gnes{gé\-né\-ra\-li\-sés }
\newcommand\gnees{gé\-né\-ra\-li\-sées }
\newcommand\gnez{gé\-né\-ra\-li\-sé}
\newcommand\gneez{gé\-né\-ra\-li\-sée}
\newcommand\gnesz{gé\-né\-ra\-li\-sés}
\newcommand\gneesz{gé\-né\-ra\-li\-sées}

\newcommand\gnl{gé\-né\-ral }
\newcommand\gnle{gé\-né\-ra\-le }
\newcommand\gnls{gé\-né\-raux }
\newcommand\gnles{gé\-né\-ra\-les }
\newcommand\gnlz{gé\-né\-ral}
\newcommand\gnlez{gé\-né\-ra\-le}
\newcommand\gnlsz{gé\-né\-raux}
\newcommand\gnlesz{gé\-né\-ra\-les}

\newcommand\gnlt{géné\-ra\-le\-ment }
\newcommand\gnltz{gé\-né\-ra\-le\-ment}

\newcommand\gnn{géné\-ra\-li\-sa\-tion }
\newcommand\gnns{géné\-ra\-li\-sa\-tions }
\newcommand\gnnz{géné\-ra\-li\-sa\-tion}
\newcommand\gnnsz{géné\-ra\-li\-sa\-tions}

\newcommand\gnq{gé\-né\-rique }
\newcommand\gnqs{gé\-né\-riques }
\newcommand\gnqz{gé\-né\-rique}
\newcommand\gnqsz{gé\-né\-riques}

\newcommand\gnr{géné\-ra\-liser }

\newcommand\gns{géné\-ra\-lise }

\newcommand\gnt{géné\-ra\-lité }
\newcommand\gnts{géné\-ra\-lités }
\newcommand\gntz{géné\-ra\-lité}
\newcommand\gntsz{géné\-ra\-lités}

\newcommand\grl{groupe réti\-culé }
\newcommand\grls{groupes réti\-culés }
\newcommand\grlz{groupe réti\-culé}
\newcommand\grlsz{groupes réti\-culés}

\newcommand\gtr{géné\-rateur }
\newcommand\gtrs{géné\-rateurs }
\newcommand\gtrz{géné\-rateur}
\newcommand\gtrsz{géné\-rateurs}

%:  abbréviations h i

\newcommand \homo {homo\-mor\-phisme }
\newcommand \homos {homo\-mor\-phismes }
\newcommand \homosz {homo\-mor\-phismes}
\newcommand \homoz {homo\-mor\-phisme}

\newcommand \hmg {homo\-gène }
\newcommand \hmgs {homo\-gènes }
\newcommand \hmgz {homo\-gène}
\newcommand \hmgsz {homo\-gènes}

\newcommand \icl {inté\-gra\-le\-ment clos }
\newcommand \iclz {inté\-gra\-le\-ment clos}

\newcommand \id {idéal }
\newcommand \ids {idéaux }
\newcommand \idz {idéal}
\newcommand \idsz {idéaux}

\newcommand \ida {\idt \agq }
\newcommand \idas {\idts \agqs }
\newcommand \idasz {\idts \agqsz}
\newcommand \idaz {\idt \agqz}

\newcommand \idc  {\idt de Cramer }
\newcommand \idcs {\idts de Cramer }
\newcommand \idcsz {\idts de Cramer}
\newcommand \idcz {\idt de Cramer}

\newcommand \idd {\id déter\-mi\-nantiel }
\newcommand \idds {\ids déter\-mi\-nantiels }
\newcommand \iddz {\id déter\-mi\-nantiel}
\newcommand \iddsz {\ids déter\-mi\-nantiels}

\newcommand \idema {\id maxi\-mal }
\newcommand \idemas {\ids maxi\-maux }
\newcommand \idemaz {\id maxi\-mal}
\newcommand \idemasz {\ids maxi\-maux}

\newcommand \idep {\id pre\-mier }
\newcommand \idepz {\id pre\-mier}
\newcommand \ideps {\ids pre\-miers }
\newcommand \idepsz {\ids pre\-miers}

\newcommand \idemi {\idep minimal }
\newcommand \idemis {\ideps minimaux }
\newcommand \idemiz {\idep minimal}
\newcommand \idemisz {\ideps minimaux}

\newcommand \idf {idéal de Fitting }
\newcommand \idfs {idéaux de Fitting }
\newcommand \idfz {idéal de Fitting}
\newcommand \idfsz {idéaux de Fitting}

\newcommand \idm {idem\-po\-tent }
\newcommand \idms {idem\-po\-tents }
\newcommand \idmz {idem\-po\-tent}
\newcommand \idmsz {idem\-po\-tents}

\newcommand \idme {idem\-po\-tente }
\newcommand \idmes {idem\-po\-tentes }
\newcommand \idmez {idem\-po\-tente}
\newcommand \idmesz {idem\-po\-tentes}

\newcommand \idp {idé\-al prin\-cipal }
\newcommand \idps {idé\-aux prin\-cipaux }
\newcommand \idpsz {idé\-aux prin\-cipaux}
\newcommand \idpz {idé\-al prin\-cipal}

\newcommand\idst{\idm de \spt}
\newcommand\idstz{\idm de \sptz}

\newcommand \idt {iden\-ti\-té }
\newcommand \idts {iden\-ti\-tés }
\newcommand \idtz {iden\-ti\-té}
\newcommand \idtsz {iden\-ti\-tés}

\newcommand \idtr {indé\-ter\-minée }
\newcommand \idtrs {indé\-ter\-minées }
\newcommand \idtrz {indé\-ter\-minée}
\newcommand \idtrsz {indé\-ter\-minées}

\newcommand \ifr {idéal frac\-tion\-naire }
\newcommand \ifrs {idéaux frac\-tion\-naires }
\newcommand \ifrz {idéal frac\-tion\-naire}
\newcommand \ifrsz {idéaux frac\-tion\-naires}

\newcommand \imd {immé\-diat }
\newcommand \imde {immé\-diate }
\newcommand \imds {immé\-diats }
\newcommand \imdes {immé\-diates }
\newcommand \imdz {immé\-diat}
\newcommand \imdez {immé\-diate}
\newcommand \imdsz {immé\-diats}
\newcommand \imdesz {immé\-diates}

\newcommand \imdt {immé\-dia\-te\-ment }
\newcommand \imdtz {immé\-dia\-te\-ment}

\newcommand \ing {in\-ver\-se \gne }
\newcommand \ings {in\-ver\-ses \gnes }
\newcommand \ingz {in\-ver\-se \gnez}
\newcommand \ingsz {in\-ver\-ses \gnesz}

\newcommand \iMP {in\-ver\-se de Moo\-re-Pen\-ro\-se }
\newcommand \iMPz {in\-ver\-se de Moo\-re-Pen\-ro\-se}
\newcommand \iMPs {in\-ver\-ses de Moo\-re-Pen\-ro\-se }
\newcommand \iMPsz {in\-ver\-ses de Moo\-re-Pen\-ro\-se}

\newcommand \ipp {\idep poten\-tiel }
\newcommand \ipps {\ideps poten\-tiels }
\newcommand \ippz {\idep poten\-tiel}
\newcommand \ippsz {\ideps poten\-tiels}

\newcommand \iso {iso\-mor\-phisme }
\newcommand \isos {iso\-mor\-phismes }
\newcommand \isosz {iso\-mor\-phismes}
\newcommand \isoz {iso\-mor\-phisme}

\newcommand \ird {irré\-duc\-tible }
\newcommand \irds {irré\-duc\-tibles }
\newcommand \irdz {irré\-duc\-tible}
\newcommand \irdsz {irré\-duc\-tibles}

\newcommand \itf {idé\-al \tf}
\newcommand \itfs {idé\-aux \tf}
\newcommand \itfz {idé\-al \tfz}
\newcommand \itfsz {idé\-aux \tfz}

\newcommand \iv {inver\-sible }
\newcommand \ivs {inver\-sibles }
\newcommand \ivz {inver\-sible}
\newcommand \ivsz {inver\-sibles}

\newcommand \ivt {inver\-si\-bi\-lité }
\newcommand \ivtz {inver\-si\-bi\-lité}

%: indexs
\newcommand \iJG{\index{Lemme de Gauss-Joyal}}

\newcommand \ihi {\index{Hilbert} }
\newcommand \ihiz {\index{Hilbert}}
\newcommand \imlg {\index{machinerie locale-globale \elrz} }
\newcommand \imlgz {\index{machinerie locale-globale \elrz}}
\newcommand \imlb {\index{machinerie locale-globale de base (\`a \idepsz)} }
\newcommand \imlbz {\index{machinerie locale-globale de base (\`a \idepsz)}}
\newcommand \imla {\index{machinerie locale-globale des \anarsz} }
\newcommand \imlaz {\index{machinerie locale-globale des \anarsz}}
\newcommand \imlma {\index{machinerie locale-globale \`a \idemasz} }
\newcommand \imlmaz {\index{machinerie locale-globale \`a \idemasz}}
\newcommand \iplg {\index{principe local-global de base} }
\newcommand \iplgz {\index{principe local-global de base}}

\newcommand \KRA {\index{Kronecker!astuce de ---}Kronecker }
\newcommand \KRO {\index{Kronecker!\tho de ---  (1)}Kronecker }
\newcommand \KRN {\index{Kronecker!\tho de ---  (2)}Kronecker }
\newcommand \KRAz {\index{Kronecker!astuce de ---}Kronecker}
\newcommand \KROz {\index{Kronecker!\tho de --- (1)}Kronecker}
\newcommand \KRNz {\index{Kronecker!\tho de --- (2)}Kronecker}

%:  abbréviations l

\newcommand \lgb {local-global }
\newcommand \lgbe {locale-globale }
\newcommand \lgbes {locales-globales }
\newcommand \lgbs {local-globals }
\newcommand \lgbz {local-global}
\newcommand \lgbez {locale-globale}
\newcommand \lgbsz {local-globals}

\newcommand \lin {liné\-aire }
\newcommand \lins {liné\-aires }
\newcommand \linz {liné\-aire}
\newcommand \linsz {liné\-aires}

\newcommand \lint {liné\-ai\-rement }

\newcommand \lmo {\lot mono\-gène }
\newcommand \lmos {\lot mono\-gènes }
\newcommand \lmoz {\lot mono\-gène}
\newcommand \lmosz {\lot mono\-gènes}

\newcommand \lnl {\lot \nl}
\newcommand \lnls {\lot \nls}
\newcommand \lnlz {\lot \nlz}
\newcommand \lnlsz {\lot \nlsz}

\newcommand \lot {loca\-lement }
\newcommand \lotz {loca\-lement}

\newcommand \Lon {Loca\-li\-sation }
\newcommand \lon {loca\-li\-sation }
\newcommand \lons {loca\-li\-sations }
\newcommand \lonz {loca\-li\-sation}
\newcommand \lonsz {loca\-li\-sations}

\newcommand \lop {\lot prin\-ci\-pal }
\newcommand \lops {\lot prin\-ci\-paux }
\newcommand \lopsz {\lot prin\-ci\-paux}
\newcommand \lopz {\lot prin\-ci\-pal}

\newcommand \lsdz {\lot \sdz}
\newcommand \lsdzz {\lot \sdzz}

%:  abbréviations m
\newcommand \mdi {mo\-dule des \diles }
\newcommand \mdiz {mo\-dule des \dilesz}

\newcommand \mlm {mo\-dule \lmo}
\newcommand \mlms {mo\-dules \lmos}
\newcommand \mlmz {mo\-dule \lmoz}
\newcommand \mlmsz {mo\-dules \lmosz}

\newcommand \mlmo {ma\-tri\-ce de \lon mono\-gène }
\newcommand \mlmos {ma\-tri\-ces de \lon mono\-gène }
\newcommand \mlmoz {ma\-tri\-ce de \lon mono\-gène}
\newcommand \mlmosz {ma\-tri\-ces de \lon mono\-gène}

\newcommand \mlp {ma\-tri\-ce de \lon prin\-ci\-pa\-le }
\newcommand \mlpz {ma\-tri\-ce de \lon prin\-ci\-pa\-le}
\newcommand \mlps {ma\-tri\-ces de \lon prin\-ci\-pa\-le }
\newcommand \mlpsz {ma\-tri\-ces de \lon prin\-ci\-pa\-le}

\newcommand \mlr {mani\-pu\-lation \elr}
\newcommand \mlrs {mani\-pu\-lations \elrs}
\newcommand \mlrz {mani\-pu\-lation \elrz}
\newcommand \mlrsz {mani\-pu\-lations \elrsz}

\newcommand \mo {mo\-no\"{\i}de }
\newcommand \mos {mo\-no\"{\i}des }
\newcommand \mosz {mo\-no\"{\i}des}
\newcommand \moz {mo\-no\"{\i}de}

\newcommand \mon {\unt }
\newcommand \mons {\unts }
\newcommand \monz {\untz }
\newcommand \monsz {\untsz }

\newcommand \moco {\mos\com}
\newcommand \mocoz {\mos\comz}

\newcommand \molo {morphisme de \lon }
\newcommand \molos {morphismes de \lon }
\newcommand \moloz {morphisme de \lonz }
\newcommand \molosz {morphismes de \lonz }

\newcommand \mom {mo\-n\^o\-me }
\newcommand \moms {mo\-n\^o\-mes }
\newcommand \momz {mo\-n\^o\-me}
\newcommand \momsz {mo\-n\^o\-mes}

\newcommand \mor{mor\-phisme }
\newcommand \mors{mor\-phismes }
\newcommand \morz{mor\-phisme}
\newcommand \morsz{mor\-phismes}

\newcommand \mpf {mo\-dule \pf}
\newcommand \mpfs {mo\-dules \pf}
\newcommand \mpfz {mo\-dule \pfz}
\newcommand \mpfsz {mo\-dules \pfz}

\newcommand \mpl {mo\-dule plat }
\newcommand \mpls {mo\-dules plats }
\newcommand \mplz {mo\-dule plat}
\newcommand \mplsz {mo\-dules plats}

\newcommand \mpn {ma\-trice de \pn }
\newcommand \mpns {ma\-trices de \pn }
\newcommand \mpnz {ma\-trice de \pnz}
\newcommand \mpnsz {ma\-trices de \pnz}

\newcommand \mprn {ma\-trice de \prn }
\newcommand \mprns {ma\-trices de \prn }
\newcommand \mprnz {ma\-trice de \prnz}
\newcommand \mprnsz {ma\-trices de \prnz}

\newcommand \mptf {mo\-dule \ptf}
\newcommand \mptfs {mo\-dules \ptfs}
\newcommand \mptfz {mo\-dule \ptfz}
\newcommand \mptfsz {mo\-dules \ptfsz}

\newcommand \mrc {mo\-dule \prc }
\newcommand \mrcz {mo\-dule \prcz}
\newcommand \mrcs {mo\-dules \prcs }
\newcommand \mrcsz {mo\-dules \prcsz}

\newcommand \mtf {mo\-du\-le \tf}
\newcommand \mtfs {mo\-du\-les \tf}
\newcommand \mtfz {mo\-du\-le \tfz}
\newcommand \mtfsz {mo\-du\-les \tfz}

%:  abbréviations n

\newcommand \imN {\index{Newton!méthode de ---}\index{methode de N@méthode de Newton}}
\newcommand \isN {\index{Newton!sommes de ---}}

\newcommand \ncr{né\-ces\-sai\-re }
\newcommand \ncrs{né\-ces\-sai\-res }
\newcommand \ncrz{né\-ces\-sai\-re}
\newcommand \ncrsz{né\-ces\-sai\-res}

\newcommand \ncrt{né\-ces\-sai\-re\-ment }
\newcommand \ncrtz{né\-ces\-sai\-re\-ment}

\newcommand \ndz {régu\-lier }
\newcommand \ndzs {régu\-liers }
\newcommand \ndzz {régu\-lier}
\newcommand \ndzsz {régu\-liers}

\newcommand \ndze {régu\-lière }
\newcommand \ndzes {régu\-lières }
\newcommand \ndzez {régu\-lière}
\newcommand \ndzesz {régu\-lières}

\newcommand \nl {sim\-ple }
\newcommand \nlz {sim\-ple}
\newcommand \nls {sim\-ples }
\newcommand \nlsz {sim\-ples}

\newcommand \Noe {Noether }
\newcommand \Noez {Noether}

\newcommand \inoe {\index{Noether@Noether!position de ---}}
\newcommand \iNoe {\inoe\Noe}
\newcommand \iNoez {\inoe\Noez}

\newcommand \noe {noethé\-rien }
\newcommand \noes {noethé\-riens }
\newcommand \noee {noethé\-rienne }
\newcommand \noees {noethé\-riennes }
\newcommand \noez {noethé\-rien}
\newcommand \noesz {noethé\-riens}
\newcommand \noeez {noethé\-rienne}
\newcommand \noeesz {noethé\-riennes}

\newcommand \noet {noethé\-ria\-nité }
\newcommand \noetz {noethé\-ria\-nité}

\newcommand \inst {\index{Null\-stellen\-satz}}
\newcommand \nst {\inst Null\-stellen\-satz }
\newcommand \nstz {\inst Null\-stellen\-satz}
\newcommand \nsts {\inst Null\-stellen\-s\"atze }
\newcommand \nstsz {\inst Null\-stellen\-s\"atze}

%:  abbréviations o

\newcommand \op{opé\-ra\-tion }
\newcommand \ops{opé\-ra\-tions }
\newcommand \opz{opé\-ra\-tion}
\newcommand \opsz{opé\-ra\-tions}
\newcommand \opari{\op\ari}
\newcommand \oparis{\ops\aris}
\newcommand \opariz{\op\ariz}
\newcommand \oparisz{\ops\arisz}
\newcommand \oparisv{\ops\arisv}

\newcommand \oqc {ouvert \qc}
\newcommand \oqcs {ouverts \qcs}
\newcommand \oqcz {ouvert \qcz}
\newcommand \oqcsz {ouverts \qcsz}

\newcommand \ort{ortho\-gonal }
\newcommand \orte{ortho\-go\-nale }
\newcommand \orts{ortho\-gonaux }
\newcommand \ortes{ortho\-go\-nales }
\newcommand \ortz{ortho\-gonal}
\newcommand \ortez{ortho\-go\-nale}
\newcommand \ortsz{ortho\-gonaux}
\newcommand \ortesz{ortho\-go\-nales}

%:  abbréviations p

\newcommand \pa {couple saturé }
\newcommand \pas {couples saturés }
\newcommand \paz {couple saturé}
\newcommand \pasz {couples saturés}

\newcommand \paral{pa\-ral\-lè\-le }
\newcommand \parals{pa\-ral\-lè\-les }
\newcommand \paralz{pa\-ral\-lè\-le}
\newcommand \paralsz{pa\-ral\-lè\-les}

\newcommand \paralm{pa\-ral\-lè\-le\-ment }

\newcommand \pb{pro\-blè\-me }
\newcommand \pbs{pro\-blè\-mes }
\newcommand \pbz{pro\-blè\-me}
\newcommand \pbsz{pro\-blè\-mes}

\newcommand \pf {de \pn finie }
\newcommand \pfz {de \pn finie}

\newcommand \plc {\rdt \zed }
\newcommand \plcs {\rdt \zeds }
\newcommand \plcz {\rdt \zedz}
\newcommand \plcsz {\rdt \zedsz}

\newcommand \plg {prin\-cipe \lgb }
\newcommand \plgs {prin\-cipes \lgbs }
\newcommand \plgz {prin\-cipe \lgbz}
\newcommand \plgsz {prin\-cipes \lgbsz}

\newcommand \plga {\plg abs\-trait }
\newcommand \plgas {\plgs abs\-traits }
\newcommand \plgaz {\plg abs\-trait}
\newcommand \plgasz {\plgs abs\-traits}

\newcommand \plgc {\plg con\-cret }
\newcommand \plgcz {\plg con\-cret}
\newcommand \plgcs {\plgs con\-crets }
\newcommand \plgcsz {\plgs con\-crets}

\newcommand \pn {présen\-ta\-tion }
\newcommand \pns {présen\-ta\-tions }
\newcommand \pnz {présen\-ta\-tion}
\newcommand \pnsz {présen\-ta\-tions}

\newcommand \pog {\pol \hmg }
\newcommand \pogs {\pols \hmgs }
\newcommand \pogz {\pol \hmgz}
\newcommand \pogsz {\pols \hmgsz}

\newcommand \Pol {Poly\-n\^ome }
\newcommand \Pols {Poly\-n\^omes }

\newcommand \pol {poly\-n\^ome }
\newcommand \pols {poly\-n\^omes }
\newcommand \polz {poly\-n\^ome}
\newcommand \polsz {poly\-n\^omes}

\newcommand \polu {\pol \unt }
\newcommand \polus {\pols \unts }
\newcommand \poluz {\pol \untz }
\newcommand \polusz {\pols \untsz }

\newcommand \polcar {\pol \cara }
\newcommand \polcarz {\pol \caraz}
\newcommand \polcars {\pols \caras }
\newcommand \polcarsz {\pols \carasz}

\newcommand \polfon {\pol fon\-da\-men\-tal }
\newcommand \polfonz {\pol fon\-da\-men\-tal}

\newcommand \poll{poly\-nomial }
\newcommand \polls{poly\-nomiaux }
\newcommand \pollsz{poly\-nomiaux}
\newcommand \pollz{poly\-nomial}
\newcommand \polle{poly\-no\-miale }
\newcommand \polles{poly\-no\-miales }
\newcommand \pollesz{poly\-no\-miales}
\newcommand \pollez{poly\-no\-miale}

\newcommand \pollt{poly\-no\-mia\-le\-ment }
\newcommand \polltz{poly\-no\-mia\-le\-ment}

\newcommand \polmin {\pol mini\-mal }
\newcommand \polmins {\pols mini\-maux }
\newcommand \polminsz {\pols mini\-maux}
\newcommand \polminz {\pol mini\-mal}

\newcommand \polmu {\pol rang }
\newcommand \polmus {\pols rang }
\newcommand \polmuz{\pol rang}

\newcommand \ppv {primi\-tif par valeurs }
\newcommand \ppvs {primi\-tifs par valeurs }
\newcommand \ppvz {primi\-tif par valeurs}
\newcommand \ppvsz {primi\-tifs par valeurs}

\newcommand \prc {\pro de rang constant }
\newcommand \prcs {\pros de rang constant }
\newcommand \prcz {\pro de rang constant}
\newcommand \prcsz {\pros de rang constant}

\newcommand \prca {prin\-ci\-pe de \rca }
\newcommand \prcc {prin\-ci\-pe de \rcc }
\newcommand \prce {prin\-ci\-pe de \rce }

\newcommand \prf {prin\-ci\-pe de recou\-vre\-ment fer\-mé }
\newcommand \prfz {prin\-ci\-pe de recou\-vre\-ment fer\-mé}

\newcommand \prmt {préci\-sé\-ment }
\newcommand \prmtz {préci\-sé\-ment}
\newcommand \Prmt {Préci\-sé\-ment }
\newcommand \Prmtz {Préci\-sé\-ment}

\newcommand \prn {pro\-jec\-tion }
\newcommand \prns {pro\-jec\-tions }
\newcommand \prnz {pro\-jec\-tion}
\newcommand \prnsz {pro\-jec\-tions}

\newcommand \pro {pro\-jec\-tif }
\newcommand \pros {pro\-jec\-tifs }
\newcommand \proz {pro\-jec\-tif}
\newcommand \prosz {pro\-jec\-tifs}

\newcommand \Prof {Profondeur }
\newcommand \prof {profondeur }
\newcommand \profz {profondeur}
\newcommand \profs {profondeurs }
\newcommand \profsz {profondeurs}

\newcommand \prr {pro\-jec\-teur }
\newcommand \prrs {pro\-jec\-teurs }
\newcommand \prrz {pro\-jec\-teur}
\newcommand \prrsz {pro\-jec\-teurs}

\newcommand \Prt {Propriété }
\newcommand \Prts {Propriétés }
\newcommand \prt {pro\-pri\-été }
\newcommand \prts {pro\-pri\-étés }
\newcommand \prtz {pro\-pri\-été}
\newcommand \prtsz {pro\-pri\-étés}

\newcommand \ptf {\pro \tf }
\newcommand \ptfz {\pro \tfz}
\newcommand \ptfs {\pros \tf }
\newcommand \ptfsz {\pros \tfz}

%:  abbréviations q

\newcommand \qc {quasi-compact }
\newcommand \qcs {quasi-compacts }
\newcommand \qcz {quasi-compact}
\newcommand \qcsz {quasi-compacts}

\newcommand \qi {qua\-si in\-tè\-gre }
\newcommand \qis {qua\-si in\-tè\-gres }
\newcommand \qisz {qua\-si in\-tè\-gres}
\newcommand \qiz {qua\-si in\-tè\-gre}

\newcommand \qiv {qua\-si inver\-se }
\newcommand \qivs {qua\-si inver\-ses }
\newcommand \qivz {qua\-si inver\-se}
\newcommand \qivsz {qua\-si inver\-sez}

\newcommand \qnl {quasi \nl}
\newcommand \qnls {quasi \nls}
\newcommand \qnlz {quasi \nlz}
\newcommand \qnlsz {quasi \nlsz}

\newcommand \qreg {quasi-régu\-lière }
\newcommand \qregz {quasi-régu\-lière}
\newcommand \qregs {quasi-régu\-lières }
\newcommand \qregsz {quasi-régu\-lières}

\newcommand \qtf {quanti\-fi\-cateur }
\newcommand \qtfs {quanti\-fi\-cateurs }
\newcommand \qtfz {quanti\-fi\-cateur}
\newcommand \qtfsz {quanti\-fi\-cateurs}

\newcommand \qtn {quanti\-fi\-cation }
\newcommand \qtnz {quanti\-fi\-cation}

%:  abbréviations r

\newcommand \RC {R-compa\-tible }
\newcommand \RCs {R-compa\-tibles }
\newcommand \RCz {R-compa\-tible}
\newcommand \RCsz {R-compa\-tibles}

\newcommand \rcc {\rcm con\-cret }
\newcommand \rca {\rcm abs\-trait }
\newcommand \rce {\rcc des éga\-li\-tés }

\newcommand \rcm {recol\-le\-ment }
\newcommand \rcmz {recol\-le\-ment}
\newcommand \rcms {recol\-le\-ments }
\newcommand \rcmsz {recol\-le\-ments}

\newcommand \rcv {recou\-vrement } 
\newcommand \rcvs {recou\-vrements } 
\newcommand \rcvz {recou\-vrement} 
\newcommand \rcvsz {recou\-vrements} 

\newcommand \rde {rela\-tion de dépen\-dance }
\newcommand \rdes {rela\-tions de dépen\-dance }
\newcommand \rdesz {rela\-tions de dépen\-dance}
\newcommand \rdez {rela\-tion de dépen\-dance}

\newcommand \rdi {\rde inté\-grale }
\newcommand \rdis {\rdes inté\-grale }
\newcommand \rdiz {\rde inté\-grale}
\newcommand \rdisz {\rdes inté\-grale}

\newcommand \rdl {\rde \lin }
\newcommand \rdls {\rdes \lin }
\newcommand \rdlsz {\rdes \linz}
\newcommand \rdlz {\rde \linz}

\newcommand \rdt {rési\-duel\-lement }
\newcommand \rdtz {rési\-duel\-lement}

\newcommand \reg {régu\-lier }   % morphisme régulier
\newcommand \regs {régu\-liers }
\newcommand \regz {régu\-lier}
\newcommand \regsz {régu\-liers}

\newcommand \rpf {réduite-de-\pnz-finie }
\newcommand \rpfs {réduites-de-\pnz-finie }
\newcommand \rpfz {réduite-de-\pnz-finie}
\newcommand \rpfsz {réduites-de-\pnz-finie}

%:  abbréviations s

\newcommand \scf {schéma spectral }
\newcommand \scfs {schémas spectraux }
\newcommand \scfz {schéma spectral}
\newcommand \scfsz {schémas spectraux}

\newcommand \scl {schéma \elr}
\newcommand \scls {schémas \elrs}
\newcommand \sclz {schéma \elrz}
\newcommand \sclsz {schémas \elrsz}

\newcommand \sdo {\sdr \orte }
\newcommand \sdos {\sdrs \ortes }
\newcommand \sdoz {\sdr \ortez }
\newcommand \sdosz {\sdrs \ortesz }

\newcommand \sdr {somme directe }
\newcommand \sdrs {sommes directes }
\newcommand \sdrz {somme directe}
\newcommand \sdrsz {sommes directes}

\newcommand \sdz {sans \dvz}
\newcommand \sdzz {sans \dvzz}

\newcommand \seco {\sex courte }
\newcommand \secos {\sexs courtes }
\newcommand \secoz {\sex courte}
\newcommand \secosz {\sexs courtes}

\newcommand \seqreg {suite régu\-lière }
\newcommand \seqregz {suite régu\-lière}
\newcommand \seqregs {suites régu\-lières }
\newcommand \seqregsz {suites régu\-lières}

\newcommand \seqqreg {suite \qreg}
\newcommand \seqqregz {suite \qregz}
\newcommand \seqqregs {suites \qregs}
\newcommand \seqqregsz {suites \qregsz}

\newcommand \sex {suite exacte }
\newcommand \sexs {suites exactes }
\newcommand \sexz {suite exacte}
\newcommand \sexsz {suites exactes}

\newcommand \sfio {\sys fonda\-mental d'\idms \orts }
\newcommand \sfios {\syss fonda\-mentaux d'\idms \orts }
\newcommand \sfioz {\sys fonda\-mental d'\idms \ortsz}
\newcommand \sfiosz {\syss fonda\-mentaux d'\idms \ortsz}

\newcommand \sgr {\sys \gtr }
\newcommand \sgrs {\syss \gtrs }
\newcommand \sgrz {\sys \gtrz}
\newcommand \sgrsz {\syss \gtrsz}

\newcommand \sing {singu\-lière }
\newcommand \sings {singu\-lières }
\newcommand \singz {singu\-lière}
\newcommand \singsz {singu\-lières}

\newcommand \slgb {semi-local }
\newcommand \slgbs {semi-locaux }
\newcommand \slgbz {semi-local}
\newcommand \slgbsz {semi-locaux}

\newcommand \sli {\sys \lin }
\newcommand \slis {\syss \lins }
\newcommand \slisz {\syss \linsz}
\newcommand \sliz {\sys \linz}

\newcommand \sml {semi-local strict }
\newcommand \smls {semi-locaux stricts }
\newcommand \smlz {semi-local strict}
\newcommand \smlss {semi-locaux stricts}

\newcommand \smq {symé\-trique }
\newcommand \smqs {symé\-triques }
\newcommand \smqz {symé\-trique}
\newcommand \smqsz {symé\-triques}

\newcommand \spb {sépa\-rable }  % algebres
\newcommand \spbs {sépa\-rables }
\newcommand \spbz {sépa\-rable}
\newcommand \spbsz {sépa\-rables}

\newcommand \spl {sépa\-rable }  % polynomes
\newcommand \spls {sépa\-rables }
\newcommand \splz {sépa\-rable}
\newcommand \splsz {sépa\-rables}

\newcommand\spt{sépa\-ra\-bi\-lité }
\newcommand\sptz{sépa\-ra\-bi\-lité}

\newcommand \sErg {suite $E$-\ndze }
\newcommand \sErgz {suite $E$-\ndzez }
\newcommand \sErgs {suites $E$-\ndzes }
\newcommand \sErgsz {suites $E$-\ndzesz }

\newcommand \srg {suite régu\-lière }
\newcommand \srgs {suites régu\-lières }
\newcommand \srgz {suite régu\-lière}
\newcommand \srgsz {suites régu\-lières}

\newcommand{\SSOz} {Splitting Off de Serre\index{Serre!Splitting Off de ---}}
\newcommand{\SSO} {Splitting Off de Serre\index{Serre!Splitting Off de ---} }

\newcommand \ste {stric\-te\-ment éta\-le }
\newcommand \stes {stric\-te\-ment éta\-les }
\newcommand \stez {stric\-te\-ment éta\-le}
\newcommand \stesz {stric\-te\-ment éta\-les}

\newcommand \stf {stric\-te\-ment fini }
\newcommand \stfs {stric\-te\-ment finis }
\newcommand \stfz {stric\-te\-ment fini}
\newcommand \stfsz {stric\-te\-ment finis}
\newcommand \stfe {stric\-te\-ment finie }
\newcommand \stfes {stric\-te\-ment finies }
\newcommand \stfez {stric\-te\-ment finie}
\newcommand \stfesz {stric\-te\-ment finies}

\newcommand \stl {stable\-ment libre }
\newcommand \stls {stable\-ment libres }
\newcommand \stlz {stable\-ment libre}
\newcommand \stlsz {stable\-ment libres}

\newcommand \styc {\sys tracique de \coos }
\newcommand \stycz {\sys tracique de \coosz}

\newcommand \sul {supplé\-men\-taire } % un supplement d'information
\newcommand \suls {supplé\-men\-taires }
\newcommand \sulz {supplé\-men\-taire}
\newcommand \sulsz {supplé\-men\-taires}

\newcommand \supl {supplé\-men\-taire } % sous modules
\newcommand \supls {supplé\-men\-taires } % sous modules
\newcommand \suplz {supplé\-men\-taire} % sous modules
\newcommand \suplsz {supplé\-men\-taires} % sous modules

\newcommand \Sus {Suslin\index{Suslin} }
\newcommand \Susz {Suslin\index{Suslin}}

\newcommand \susi {suite \sing }
\newcommand \susis {suites \sings }
\newcommand \susiz {suite \singz }
\newcommand \susisz {suites \singsz }

\newcommand \Sut {Support }
\newcommand \Suts {Supports }
\newcommand \sut {support }
\newcommand \suts {supports }
\newcommand \sutz {support}
\newcommand \sutsz {supports}

\newcommand \syc {\sys de coordon\-nées }
\newcommand \sycs {\syss de coordon\-nées }
\newcommand \sycz {\sys de coordon\-nées}
\newcommand \sycsz {\syss de coordon\-nées}

\newcommand \syp {\sys \poll }
\newcommand \syps {\syss \polls }
\newcommand \sypz {\sys \pollz}
\newcommand \sypsz {\syss \pollsz}

\newcommand \Sys {Système }
\newcommand \sys {sys\-tème }
\newcommand \syss {sys\-tèmes }
\newcommand \sysz {sys\-tème}
\newcommand \syssz {sys\-tèmes}

\newcommand \syzy {syzygie }
\newcommand \syzys {syzygies }
\newcommand \syzyz {syzygie}
\newcommand \syzysz {syzygies}

%:  abbréviations t

\newcommand \tf {de ty\-pe fi\-ni }
\newcommand \tfz {de ty\-pe fi\-ni}

\newcommand \Tho {Théo\-rème }
\newcommand \Thos {Théo\-rèmes }
\newcommand \tho {théo\-rème }
\newcommand \thos {théo\-rèmes }
\newcommand \thoz {théo\-rème}
\newcommand \thosz {théo\-rèmes}

\newcommand \thoc {\thoz\etoz~}

\newcommand \trdi {treillis distri\-butif }
\newcommand \trdis {treillis distri\-butifs }
\newcommand \trdiz {treillis distri\-butif}
\newcommand \trdisz {treillis distri\-butifs}

\newcommand \trel {trans\-for\-mation \elr}
\newcommand \trels {trans\-for\-mations \elrs}
\newcommand \trelz {trans\-for\-mation \elrz}
\newcommand \trelsz {trans\-for\-mations \elrsz}

%:  abbréviations u

\newcommand \umd {unimo\-du\-laire }
\newcommand \umds {unimo\-du\-laires }
\newcommand \umdz {unimo\-du\-laire}
\newcommand \umdsz {unimo\-du\-laires}

\newcommand \unt {uni\-taire }
\newcommand \unts {uni\-taires }
\newcommand \untz {uni\-taire}
\newcommand \untsz {uni\-taires}

\newcommand \uvl {uni\-versel }
\newcommand \uvle {uni\-ver\-selle }
\newcommand \uvls {uni\-versels }
\newcommand \uvles {uni\-ver\-selles }
\newcommand \uvlz {uni\-versel}
\newcommand \uvlez {uni\-ver\-selle}
\newcommand \uvlsz {uni\-versels}
\newcommand \uvlesz {uni\-ver\-selles}

%:  abbréviations v

\newcommand \vgq {\vrt \agq}
\newcommand \vgqs {\vrts \agqs}
\newcommand \vgqz {\vrt \agqz}
\newcommand \vgqsz {\vrts \agqsz}

\newcommand \vrt {varié\-té }
\newcommand \vrts {varié\-tés }
\newcommand \vrtz {varié\-té}
\newcommand \vrtsz {varié\-tés}

\newcommand \vfn {véri\-fi\-cation }
\newcommand \vfns {véri\-fi\-cations }
\newcommand \vfnz {véri\-fi\-cation}
\newcommand \vfnsz {véri\-fi\-cations}

\newcommand \vmd {vec\-teur \umd}
\newcommand \vmds {vec\-teurs \umds}
\newcommand \vmdz {vec\-teur \umdz}
\newcommand \vmdsz {vec\-teurs \umdsz}

\newcommand \zed {zéro-dimen\-sionnel }
\newcommand \zedz {zéro-dimen\-sionnel}
\newcommand \zede {zéro-dimen\-sion\-nelle }
\newcommand \zedez {zéro-dimen\-sion\-nelle}
\newcommand \zeds {zéro-dimen\-sionnels }
\newcommand \zedsz {zéro-dimen\-sionnels}
\newcommand \zedes {zéro-dimen\-sion\-nelles }
\newcommand \zedesz {zéro-dimen\-sion\-nelles}

\newcommand \zedr {\zed réduit }
\newcommand \zedrs {\zeds réduits }
\newcommand \zedrz {\zed réduit}
\newcommand \zedrsz {\zeds réduits}

\newcommand \zmt {\tho de Zariski-Grothen\-dieck }
\newcommand \zmtz {\tho de Zariski-Grothen\-dieck}

%:  ------- maths constructives

\newcommand \cof {cons\-truc\-tif }
\newcommand \cofs {cons\-truc\-tifs }
\newcommand \cofz {cons\-truc\-tif}
\newcommand \cofsz {cons\-truc\-tifs}

\newcommand \cov {cons\-truc\-ti\-ve }
\newcommand \covz {cons\-truc\-ti\-ve}
\newcommand \covsz {cons\-truc\-ti\-ves}
\newcommand \covs {cons\-truc\-ti\-ves }

\newcommand \coma {\maths\covs}
\newcommand \comaz {\maths\covsz}
\newcommand \clama {\maths clas\-si\-ques }
\newcommand \clamaz {\maths clas\-si\-ques}

\renewcommand \cot {cons\-truc\-ti\-vement }
\newcommand \cotz {cons\-truc\-ti\-vement}

\newcommand \matn {mathé\-ma\-ti\-cien }
\newcommand \matne {mathé\-ma\-ti\-cienne }
\newcommand \matns {mathé\-ma\-ti\-ciens }
\newcommand \matnes {mathé\-ma\-ti\-ciennes }

\newcommand \maths {mathé\-ma\-tiques }
\newcommand \mathsz {mathé\-ma\-tiques}
\newcommand \mathe {mathé\-ma\-tique }
\newcommand \mathz {mathé\-ma\-tique}

\newcommand \prco {preuve \cov}
\newcommand \prcos {preuves \covs}
\newcommand \prcoz {preuve \covz}
\newcommand \prcosz {preuve \covsz}

%: Macros compliquees
\newcommand\cm{cm}
\makeatletter

%:  \\ retabli
\DeclareRobustCommand\\{%
  \let \reserved@e \relax
  \let \reserved@f \relax
  \@ifstar{\let \reserved@e \vadjust \let \reserved@f \nobreak
             \@xnewline}%
          \@xnewline}
\makeatother

%:             blocs  (2 X 2 blocs)
\newcommand{\blocs}[8]{%
{\setlength{\unitlength}{.0833\textwidth}
\tabcolsep0pt\renewcommand{\arraystretch}{0}%
\begin{tabular}{|c|c|}
\hline
\parbox[t][#3\cm][c]{#1\cm}{\begin{minipage}[c]{#1\cm}
\centering#5
\end{minipage}}&
\parbox[t][#3\cm][c]{#2\cm}{\begin{minipage}[c]{#2\cm}
\centering#6
\end{minipage}}\\
\hline
\parbox[t][#4\cm][c]{#1\cm}{\begin{minipage}[c]{#1\cm}
\centering#7
\end{minipage}}&
\parbox[t][#4\cm][c]{#2\cm}{\begin{minipage}[c]{#2\cm}
\centering#8
\end{minipage}}\\
\hline
\end{tabular}
}}

% exemple
% \blocs{.8}{.6}{.8}{.6}{$\I_k$}{$0$}{$0$}{$0$}

%:             tri (triangle)
\newcommand\tri[7]{
$$\quad\quad\quad\quad
\vcenter{\xymatrix@C=1.5cm
{
#1 \ar[d]_{#2} \ar[dr]^{#3} \\
{#4} \ar[r]_{{#5}}   & {#6} \\
}}
\quad\quad \vcenter{\hbox{\small {#7}}\hbox{~\\[1mm] ~ }}
$$
}

%:             carre (carré)

\newcommand\carre[8]{
$$
\xymatrix @C=1.2cm{
#1\,\ar[d]^{#4}\ar[r]^{#2}   & \,#3\ar[d]^{#5}   \\
#6\,\ar[r]    ^{#7}    & \,#8  \\
}
$$
}

%:             pun  fragile
\newcommand\pun[7]{
$$\quad\quad\quad\quad
\vcenter{\xymatrix@C=1.5cm
{
#1 \ar[d]_{#2} \ar[dr]^{#3} \\
{#4} \ar@{-->}[r]_{{#5}\,!}   & {#6} \\
}}
\quad\quad \vcenter{\hbox{\small {#7}}\hbox{~\\[1mm] ~ }}
$$
}

\newcommand\puN[8]{
$$\hspace{#8}
\vcenter{\xymatrix@C=1.5cm
{
#1 \ar[d]_{#2} \ar[dr]^{#3} \\
{#4} \ar@{-->}[r]_{{#5}\,!}   & {#6} \\
}}
\quad\quad \vcenter{\hbox{\small {#7}}\hbox{~\\[1mm] ~ }}
$$
}

%:              Pun et autres
\newcommand\Pun[8]{
$$\quad\quad\quad\quad
\vcenter{\xymatrix@C=1.5cm
{
#1 \ar[d]_{#2} \ar[dr]^{#3} \\
{#4} \ar@{-->}[r]_{{#5}\,!}   & {#6} \\
}}
\quad\quad
\vcenter{\hbox{\small {#7}}
\hbox{~\\[3.5mm] ~ }
\hbox{\small {#8}}
\hbox{~\\[-3.5mm] ~ }}
$$
}

%\def\PUn{\PUN}

%%%%%%%%%%%%%%%%%%%%%%%%%%%%%%%%%%%%%%%%%
\newcommand\PUN[9]{
$$\quad\quad
\vcenter{\xymatrix@C=1.5cm
{
#1 \ar[d]_{#2} \ar[dr]^{#3} \\
{#4} \ar@{-->}[r]_{{#5}\,!}   & {#6} \\
}}
\quad\quad
\vcenter{
\hbox{\small {#7}}
\hbox{~\\[-3mm] ~}
\hbox{\small {#8}}
\hbox{~\\[-3mm] ~}
\hbox{\small {#9}}
\hbox{~\\[-3.5mm] ~ }}
$$
}

%%%%%%%%%%%%%%%%%%%%%%%%%%%%%%%%%%%%%%%%%
\newcommand\Pnv[9]{
$$\quad\quad\quad\quad
\vcenter{\xymatrix@C=1.5cm
{
#1 \ar[d]_{#2} \ar[dr]^{#3} \\
{#4} \ar@{-->}[r]_{{#5}\,!}   & {#6} \\
}}
\quad\quad
\vcenter{
\hbox{\small {#7}}
\hbox{~\\[1mm] ~}
\hbox{\small {#8}}
\hbox{~\\[-1mm] ~}
\hbox{\small {#9}}
\hbox{~\\[0mm] ~ }}
$$
}

%%%%%%%%%%%%%%%%%%%%%%%%%%%%%%%%%%%%%%%%%
\newcommand\PNV[9]{
$$\quad\quad\quad\quad
\vcenter{\xymatrix@C=1.5cm
{
#1 \ar[d]_{#2} \ar[dr]^{#3} \\
{#4} \ar@{-->}[r]_{{#5}\,!}   & {#6} \\
}}
\quad\quad
\vcenter{
\vspace{4mm}
\hbox{\small {#7}}
\hbox{~\\[-1.7mm] ~}
\hbox{\small {#8}}
\hbox{~\\[-1.7mm] ~}
\hbox{\small {#9}}
\hbox{~\\[2mm] ~ }
}
$$
}

%%%%%%%%%%%%%%%%%%%%%%%%%%%%%%%%%%%%%%%%%

%:     SCO  (suites complementaires)
\newdimen\xyrowsp
\xyrowsp=3pt
\newcommand{\SCO}[6]{
\xymatrix @R = \xyrowsp {
                                  &1 \ar@{-}[dl] \ar@{-}[dr] \\
#3 \ar@{-}[ddr]                   &   & #6 \ar@{-}[ddl] \\
                                  &\bullet\ar@{-}[d] \\
                                  &\bullet   \\
#2 \ar@{-}[ddr] \ar@{-}[uur]      &   & #5 \ar@{-}[ddl] \ar@{-}[uul] \\
                                  &\bullet \ar@{-}[d] \\
                                  &\bullet  \\
#1 \ar@{-}[uur]                   &   & #4 \ar@{-}[uul] \\
                                  & 0 \ar@{-}[ul] \ar@{-}[ur] \\
}
}

%%%%%%%%%%%%%%%%%%%%%%%%%%%%%%%%%%%%%%%%%%%%%%%%%%%%%%%%%%%%%%%%%%%

\makeatletter
\newif\if@borderstar
\def\bordercmatrix{\@ifnextchar*{%
  \@borderstartrue\@bordercmatrix@i}{\@borderstarfalse\@bordercmatrix@i*}%
}
\def\@bordercmatrix@i*{\@ifnextchar[{%
  \@bordercmatrix@ii}{\@bordercmatrix@ii[()]}
}
\def\@bordercmatrix@ii[#1]#2{%
  \begingroup
    \m@th\@tempdima.875em\setbox\z@\vbox{%
      \def\cr{\crcr\noalign{\kern 2\p@\global\let\cr\endline}}%
      \ialign {$##$\hfil\kern.2em\kern\@tempdima&\thinspace%
      \hfil$##$\hfil&&\quad\hfil$##$\hfil\crcr\omit\strut%
      \hfil\crcr\noalign{\kern-\baselineskip}#2\crcr\omit%
      \strut\cr}}%
    \setbox\tw@\vbox{\unvcopy\z@\global\setbox\@ne\lastbox}%
    \setbox\tw@\hbox{\unhbox\@ne\unskip\global\setbox\@ne\lastbox}%
    \setbox\tw@\hbox{%
      $\kern\wd\@ne\kern-\@tempdima\left\@firstoftwo#1%
        \if@borderstar\kern.2em\else\kern -\wd\@ne\fi%
      \global\setbox\@ne\vbox{\box\@ne\if@borderstar\else\kern.2em\fi}%
      \vcenter{\if@borderstar\else\kern-\ht\@ne\fi%
        \unvbox\z@\kern-\if@borderstar2\fi\baselineskip}%
\if@borderstar\kern-2\@tempdima\kern.4em\else\,\fi\right\@secondoftwo#1 $%
    }\null\;\vbox{\kern\ht\@ne\box\tw@}%
  \endgroup
}
\makeatother

%:  hyphenation

\hyphenation{
abs-trait abs-traite
abs-traits abs-traites
ad-di-tion ad-di-tions
af-fec-ta-tion af-fec-ta-tions an-neau an-neaux
ap-par-tien-nent ap-par-tient
ap-pli-ca-tion ap-pli-ca-tions
ap-proxi-ma-tif ap-proxi-ma-tifs
ap-proxi-ma-ti-ve ap-proxi-ma-ti-ves
ap-proxi-ma-ti-ve-ment
ar-bi-trai-re ar-bi-trai-re-ment ar-bi-trai-res
a-symp-to-ti-que a-symp-to-ti-ques a-symp-to-ti-que-ment
au-cu-ne au-cun au-cu-nes au-cuns
auto-mor-phi-sme auto-mor-phi-smes
auxi-liai-re auxi-liai-res
axio-me axio-mes
calcu-ler cal-cul cal-culs
calcu-la-ble calcu-la-bles clas-se clas-ses
cano-ni-ques cano-ni-que cano-ni-que-ment
cepen-dant
chan-ge-ment chan-ge-ments
cher-che cher-cher
clai-re-ment
clas-si-que clas-si-ques clas-si-que-ment
co-lon-nes co-lon-ne
com-me com-ment com-men-tai-re com-men-tai-res
com-mu-ta-tif com-mu-ta-tifs com-mu-ta-ti-ve com-mu-ta-ti-ves
com-por-te com-por-tent com-por-ter com-por-te-ment
com-po-san-te com-po-san-tes
con-jec-tu-re con-jec-tu-res con-jec-tu-rer
con-gru-en-tiel con-gru-en-tiels
cons-tant cons-tants cons-tan-te cons-tan-tes
cons-trui-re cons-trui-te cons-trui-tes
cons-truit cons-truits cons-trui-sons
cons-truc-tion cons-truc-tions
conve-na-ble conve-na-bles
cor-rect cor-rec-te cor-rects cor-rec-tes
cor-res-pond cor-res-pon-dre cor-res-pon-dent
cor-res-pon-dant cor-res-pon-dante
cor-res-pon-dants cor-res-pon-dantes
cou-ram-ment cou-rant cou-rante cou-rants cou-rantes
cyclo-to-mi-que cyclo-to-mi-ques
des-cen-dant des-cen-dants des-cen-dan-te des-cen-dan-tes
des-crip-tion des-crip-tions des-crip-ti-ble des-crip-ti-bles
de-ve-nir de-vien-nent
dia-go-nal dia-go-na-le dia-go-na-les dia-go-naux
di-mi-nu-tion
di-vi-sion di-vi-sions
do-mi-nant do-mi-nan-te do-mi-nants do-mi-nan-tes
en-sem-ble en-sem-bles
es-sen-tiel es-sen-tiels es-sen-tiel-le es-sen-tiel-les
es-sen-tiel-le-ment
exac-te-ment exac-te exac-tes
exa-men exa-mi-ne exa-mi-ner
exem-ples exem-ple
ex-cep-tion-nel ex-cep-tion-nels
ex-cep-tion-nel-le ex-cep-tion-nel-les
ex-cep-tion-nel-le-ment
ex-pli-ci-te ex-pli-ci-tes ex-pli-ci-ter ex-pli-ci-tons
ex-po-sant ex-po-sants ex-po-sons ex-plo-sion
exis-te exis-tent exis-tant exis-tan-te exis-tants
exis-tan-tes exis-ten-ce
feuil-le feuil-les
fonc-tion fonc-tions
fonc-tion-ne fonc-tion-ner fonc-tion-nent
fon-da-men-tal fon-da-men-tale fon-da-men-tale-ment
for-me for-mes for-mel for-mel-le for-mels for-mel-les
for-mel-le-ment for-ma-tion
grou-pe grou-pes
ha-bi-tuel-le-ment ha-bi-tuel-le ha-bi-tuel-les
 ha-bi-tuel  ha-bi-tuels
i-den-ti-que i-den-ti-ques i-den-ti-que-ment
i-den-ti-fi-ca-tion i-den-ti-fi-ca-tions
i-den-ti-fie i-den-ti-fier
im-por-tant im-por-tan-te im-por-tants im-por-tan-tes
in-for-ma-tion in-for-ma-tions
ins-truc-tion ins-truc-tions
in-ter-po-la-tion in-ter-po-la-tions
in-ter-vien-nent in-ter-vient
in-ver-se in-ver-ses in-ver-si-ble in-ver-si-bles
in-ver-se-ment in-ver-sion in-ver-sions
le-quel les-quels la-quel-le les-quel-les
le-xi-co-gra-phi-que
lo-cal loca-li-sation loca-li-sations
loca-li-sant loca-li-ser loca-li-se
lors-que main-te-nant
ma-jo-rant ma-jo-rants ma-jo-ra-tion ma-jo-ra-tions
ma-ni-pu-la-tion ma-ni-pu-la-tions
ma-tri-ce ma-tri-ces ma-tri-ciel ma-tri-ciels
 ma-tri-ciel-le  ma-tri-ciel-les ma-xi-mum mi-ni-mum
 meil-leur meil-leurs meil-leu-re meil-leu-res
mi-no-rant mi-no-rants mi-no-ra-tion mi-no-ra-tions
mo-du-le mo-du-les
na-tu-rel na-tu-rels na-tu-rel-le na-tu-rel-les
na-tu-rel-le-ment nor-me nor-mes
nor-mal nor-male nor-maux nor-males nor-ma-le-ment
no-tam-ment
op-ti-mal op-ti-male
or-tho-go-nal or-tho-go-na-le or-tho-go-naux or-tho-go-na-les
pa-ra-gra-phe pa-ra-gra-phes
par-ti-cu-lier par-ti-cu-liers
per-mu-ta-tion per-mu-ta-tions
pi-vot pi-vots
po-si-tif po-si-tifs po-si-ti-ve po-si-ti-ves  po-si-ti-ve-ment
preu-ve preu-ves prio-ri prin-ci-pe
pro-ces-seur pro-ces-seurs pro-ces-sus
pro-gram-me pro-gram-mes
pro-po-si-tion pro-po-si-tions
qua-dra-ti-que qua-dra-ti-ques
quel-le quel-les
rai-son-na-ble rai-son-na-bles rai-son-na-ble-ment
rai-son-ne-ment rai-son-ne-ments
ra-re ra-res ra-re-ment
ra-tion-nel ra-tion-nels ra-tion-nel-le ra-tion-nel-les
recol-lement
rec-tan-gu-lai-re rec-tan-gu-lai-res
re-gis-tre re-gis-tres
rela-tif rela-ti-ve rela-tifs rela-ti-ves
rela-tion rela-tions
ren-sei-gne-ment ren-sei-gne-ments
res-pec-tif res-pec-tifs res-pec-ti-ve res-pec-ti-ves
 res-pec-ti-ve-ment
 res-trein-dre res-tric-tion res-tant
seu-le-ment
si-gni-fi-ca-tif si-gni-fi-ca-ti-ve si-gni-fi-ca-tifs
si-gni-fi-ca-ti-ves si-gni-fi-ca-ti-ve-ment
si-mu-ler si-mu-la-tion si-mu-la-tions
som-me som-mes som-ma-tions som-ma-tion
sous-trac-tion sous-trac-tions
subs-ti-tuer subs-ti-tu-tion subs-ti-tu-tions
suc-ces-sifs suc-ces-si-ves suc-ces-si-ve-ment
suf-fi-sant suf-fi-sants suf-fi-san-te suf-fi-san-tes
suf-fi-sam-ment
sui-vant sui-van-te sui-vants sui-van-tes
sup-po-ser sup-po-sons Sup-po-sons
sym-bo-li-que sym-bo-li-ques sym-bo-li-que-ment
ten-seur ten-seurs ten-so-riel-le
tou-te tou-tes tou-te-fois
tra-ce tra-ces
trans-for-ma-tions trans-for-ma-tion
trans-por-teur trans-por-teurs
tri-an-gu-lai-re tri-an-gu-lai-res
tri-an-gu-la-tions tri-an-gu-la-tion
tri-an-gu-la-ri-sa-tion
u-ni-for-me u-ni-for-mes
u-ni-mo-du-lai-re u-ni-mo-du-lai-res
u-ni-que u-ni-ques u-ni-que-ment u-ni-tai-re
 u-ti-le u-ti-les
u-ti-li-sa-tion u-ti-li-ser u-ti-li-se u-ti-li-sons
vala-ble vala-bles
va-ria-ble va-ria-bles va-rian-te va-rian-tes
va-ria-tion va-ria-tions
}

\newcommand\sibrouillon[1]{}

\newcounter{nbtotalexos}
\newcounter{nbtotalprob}
\newcommand*{\incrementeexosetprob}{%
\addtocounter{nbtotalexos}{\value{exercise}}
\addtocounter{nbtotalprob}{\value{problem}}}
\newcommand*{\finincrementeexosetprob}{%
\addtocounter{nbtotalexos}{-1}\refstepcounter{nbtotalexos}\label{nombreexos}%
\addtocounter{nbtotalprob}{-1}\refstepcounter{nbtotalprob}\label{nombreprob}%
}
\romanpagenumbers

%!TEX encoding =  UTF-8 Unicode
%!TEX root =  ACMC-A.tex
\begingroup

\pagestyle{CMcadreseul}

\def\entree#1.{\setbox0=\hbox{#1. --- }\hangindent \wd0\hangafter1\noindent\box0\ignorespaces}

\setpagenumber1
\begin{center}\Large
\sc mathématiques en devenir
\end{center}

\vspace{2cm}
%:     center
\begin{center} \label{center}

Texte pour la 2\ieme\,édition française  du livre 

\medskip 
\bf {\Large Algèbre commutative 

Méthodes constructives}

\medskip 
{\large Modules projectifs de type fini}
\end{center}

\rm
\medskip 
\Grandcadre{version actualisée le \today.}

\medskip 
Nous avons corrigé des erreurs, ajouté des solutions d'exercices ainsi que quelques compléments, et le nombre de pages a augmenté d'une centaine.

\medskip 
Aucune numérotation n'a changé, sauf le \plg XII-7.13  devenu XII-7.14.

Toutes précisions utiles sur le site:

\url{http://hlombardi.free.fr/publis/LivresBrochures.html}

\bigskip 
\entrenous{

Ceci est en fait la version commentée, avec des notes en marge\perso{voici une note en marge} et des commentaires, signalés comme celui-ci.

}
%%%%%%%%%%%%%%%%%%%%%%%%%%%%%%%%%%%%%%%%%%%%%%%%%%%%%%%%%%%%%%%%%%%%%%%%%%%

\vfill\eject

\noindent{\large\it Mathématiques en devenir}
\vskip1em

\entree 101. Jacques Faraut. {\sl Analyse sur les groupes de Lie. Une introduction}

\entree 102. Patrice Tauvel. {\sl Corps commutatifs et théorie de Galois}

\entree 103. Jean Saint Raymond. {\sl Topologie, calcul différentiel et 
variable complexe}

\entree 104. Clément de Seguins Pazzis. {\sl Invitation aux formes quadratiques}

\entree 105. Bruno Ingrao. {\sl Coniques projectives, affines et métriques}

\entree 106. Wolfgang Bertram. {\sl Calcul différentiel topologique élémentaire}

\entree 107. Henri Lombardi \&\ Claude Quitté. {\sl Algèbre commutative. Méthodes constructives. Modules projectifs de type fini}

\entree 108. Frédéric Testard. {\sl Analyse mathématique. La maîtrise de l'implicite}

\entree 109. Grégory Berhuy. {\sl Modules: théorie, pratique\dots\ et un peu d'arithmétique}

\entree 110. Bernard Candelpergher. {\sl Théorie des probabilités. Une introduction élémentaire}

\entree 111. Philippe Caldero et Jérôme Germoni. {\sl Histoires hédonistes de groupes et de géométries. Deux tomes.}

\entree 112. Gema-Maria Díaz-Toca, Henri Lombardi \&\ Claude Quitté. {\sl Modules sur les anneaux commutatifs.}

\cleardoublepage

\begin{center}\Large\bf
\CMauthor 
\end{center}
\vskip1.7cm
\begin{center}
%\Huge\bf \CMtitle
\bf \Huge 
 Algèbre commutative\\ méthodes constructives
\\[3mm]
\LARGE Modules \ptfs  
\\[5mm]
\large
Cours et exercices
\\[7mm]
\rm
2\ieme\,édition
\\[3mm]
\normalsize \rm
Dernière mise à jour, \today

\vskip0pt plus1filll
\large\bf
Calvage \& Mounet
\\[3mm]

\end{center}
\break

\thispagestyle{CMcadreseul}
\begingroup
\parskip0pt

\begin{flushleft}

\noindent {\sc Henri Lombardi. }Maître de Conférences à l'Université de Franche-Comté et membre de
l'\'Equipe de Mathématique de Besan\c con (UMR 6623).
Ses recherches concernent les mathématiques constructives, l'algèbre réelle
et la complexité algorithmique. \\
Il est l'un des initiateurs du groupe international M.A.P. (Mathematics, Algorithms, Proofs), créé en 2003: voir le site \url{http://map.disi.unige.it/}

%\smallskip 
Il a publié les ouvrages suivants. 
\begin{itemize}
\item \emph{Modules sur les anneaux commutatifs}, Calvage\&Mounet, 2014,
en collaboration avec Gema Díaz-Toca et Claude Quitté.
\item \emph{\'Epistémologie mathématique}, Ellipse, 2011.
\item \emph{Méthodes matricielles. Introduction à la complexité algébrique}, Springer, 2003, en collaboration avec Jounaïdi Abdeljaoued.
\item \emph{Géométries élémentaires (tome 1)}, Presses Universitaires de Franche-Comté. 1999.
\end{itemize}

\hfill{\tt henri.lombardi@univ-fcomte.fr}
\null\hfill \url {http://hlombardi.free.fr}

\smallskip   {\sc Claude Quitté. }Maître de conférences à l'Université de Poitiers 
et membre du Laboratoire de Mathématiques et 
Applications de l'Université de Poitiers (UMR 6086).
Ses recherches concernent l'algèbre commutative 
effective et le calcul formel. 
Il a enseigné à tous les niveaux (en particulier dans la préparation à l'agrégation), et il est intervenu dans des enseignements combinant
mathématiques et informatique. 
Il a programmé en {\tt Magma} de très nombreux algorithmes en relation directe
avec le présent ouvrage (cours et/ou exercices).\\
En collaboration avec Patrice Naudin, il a publié l'ouvrage
\emph{Algorithmique algébrique}, Masson, 1991.\\
Avec Henri Lombardi, il a participé à la rédaction de l'ouvrage 
collectif \emph{Mathématiques L3 Algèbre.} Pearson Education, 2009.
\\
Il a publié \emph{Modules sur les anneaux commutatifs}, Calvage\&Mounet, 2014,
en collaboration avec Gema Díaz-Toca et Henri Lombardi.

\hfill{\tt claude.quitte@math.univ-poitiers.fr}

\vfill
\begingroup
\noindent 
Mathematics Subject Classification (2010) 

-- Primary: 13 ~ Commutative Algebra.

-- Secondary:

\leftskip4ex
03F ~ Proof theory and constructive mathematics. 

\leftskip4ex
06D ~ Distributive lattices.

\leftskip4ex
14Q ~ Computational aspects of algebraic geometry.

\vfill
\vfill

\end{flushleft}

\endgroup

%%%%%%% ISBN
%\vspace{1em}
%\noindent
%\kern0.5pt\raise-2pt\hbox{\includegraphics{C8}} \kern0.5pt\ 
%Imprimé sur papier permanent\hfill
%\vtop to0em{\kern-35pt\hbox{\ISBN 0000000000000  \barheight=0,9cm \EAN 0000000000000 %\www\ (Lombardi-Quitt-Tome 1)
%}\vss}
%\vskip0.2truecm
%\noindent\copyright\ \ Calvage \&\ Mounet, Paris, 201?
% \endgroup

\newpage
\thispagestyle{CMcadreseul}
\null\vfill\vfill
\hfill{\it à  James Brewer}
\vfill
\null
\newpage

%%%%%%%%%%%%%%%%%%%%%%%%%%%%%%
\let\oldshowchapter\showchapter
\let\oldshowsection\showsection
\def\showchapter#1{\edef\temp{\thechapter}\def\ttemp{#1}\ifx\ttemp\temp\relax\def\showsection##1{##1}\else\let\showsection\oldshowsection\oldshowchapter{#1}\fi}
%%%%%%%%%%%%%%%%%%%%%%%%%%%%%%
\chapter*{Préface de la première édition}
\markboth{Préface}{Préface}

Ce livre est un cours d'introduction à l'\alg commutative de base, avec un accent
particulier mis sur les \mptfsz, qui constituent la version \agq des
fibrés vectoriels en \gmt \dilez.

Nous adoptons le point de vue \cofz, avec lequel tous les \thos
d'existence ont un contenu \algq explicite.
 En particulier, lorsqu'un \tho affirme l'existence d'un objet, solution d'un \pbz,
un \algo de construction de l'objet peut toujours
être extrait de la \dem qui est donnée.

Nous revisitons avec un regard nouveau et souvent simplificateur
plusieurs théories classiques abstraites.  
En particulier, nous revenons sur des théories qui n'avaient 
pas de contenu \algq  dans leur cadre naturel \gnlz, comme la théorie de Galois, celle des \adksz, 
celle des \mptfs ou celle de la \ddkz.  

L'\alg \cov est en fait une vieille discipline,
développée entre autres par Gauss et Kronecker.
Nous nous situons dans la lignée de la \gui{bible} moderne sur le 
sujet, qu'est le livre {\it  A Course in
Constructive Algebra} de Ray Mines, Fred Richman et Wim Ruitenburg, paru en 1988.
Nous le citerons sous forme abrégée \cite{MRR}.

L'ouvrage correspond à un niveau de Master 2, du moins 
jusqu'au chapitre~\ref{chapNbGtrs}, mais ne réclame comme prérequis que les notions de base concernant la théorie des groupes, l'\alg \lin sur les corps,
les \detersz,
les modules sur les anneaux commutatifs, ainsi que la \dfn 
des anneaux quotients et localisés.
Une familiarité avec les anneaux de \polsz, les \prts \aris de $\ZZ$ et des anneaux euclidiens
est \egmt souhaitable.

Signalons enfin que nous considérons les  exercices et  \pbs (un peu plus de 320 en tout) 
comme  une partie essentielle de l'ouvrage. 

Nous essaierons de publier le maximum de corrigés manquants, ainsi que des exercices supplémentaires, sur la page web de l'un des auteurs:\\ 
\url{http://hlombardi.free.fr/publis/LivresBrochures.html}.

\subsection*{Remerciements.} Nous remercions tou(te)s les collègues
qui nous ont encouragés dans notre projet, nous ont apporté quelques sérieux coups de main ou fourni de précieuses informations. Et tout particulièrement  MariEmi Alonso, Thierry Coquand, Gema Díaz-Toca, Lionel Ducos,  M'hammed El Kahoui, Marco Fontana, Sarah Glaz, Laureano Gonz\'alez-Vega, Emmanuel Hallouin,  Hervé Perdry, Jean-Claude Raoult, Fred Richman, Marie-Fran\c{c}oise Roy,
  Peter Schuster et Ihsen Yengui. 
Last but not least, une mention toute spéciale pour notre expert Latex, Fran\c{c}ois Pétiard.

\goodbreak
Enfin, nous ne saurions oublier le Centre International de Recherches Mathématiques à Luminy et 
le Mathematisches Forschungsinstitut Oberwolfach, 
qui nous ont accueillis pour des séjours de recherche
pendant la préparation de ce livre, 
nous offrant des conditions de travail inappréciables. 

\vspace{5mm}
\begin{flushright}
Henri Lombardi, Claude Quitté\\
Août 2011
\end{flushright}

%\newpage
%
%~
%\vspace{3em}
%\markboth{Préface}{Préface}

\chapter*{Préface de la deuxième édition}

Dans cette deuxième édition, nous avons corrigé les erreurs
que nous avons débusquées ou qui nous ont été signalées.

Nous avons ajouté des solutions d'exercices (notamment celle de l'exercice VII-3, légèrement modifié) ainsi que quelques compléments. La plupart des compléments sont des corrections d'exercices ou de nouveaux exercices ou \pbsz.

Les ajouts dans le cours sont les suivants. Un paragraphe sur les tenseurs nuls ajouté à la fin de la section \ref{secStabPf}. 
Le paragraphe sur les quotients de modules plats à la fin de la section \ref{secPlatDebut} a été étoffé. La section \ref{secAnor} a été rajoutée pour discuter un \pb intéressant de décryptage des \dems classiques,
insensibles à la distinction entre anneaux \sdz et anneaux intègres,
pertinente du point de vue \cofz.
Enfin, on a rajouté  deux sections \iref{secPlgcor} et \iref{secPlgprof2} dans le chapitre \ref{chapPlg}
consacré aux \plgsz.\\
Notons aussi que nous avons en général remplacé l'expression \gui{\rdlz} par le terme plus court et plus usuel aujourd'hui \gui{syzygie}.  

Aucune numérotation n'a changé, sauf le \plg XII-7.13  devenu \ref{plcc.ddk}. Le nombre de pages a  augmenté d'une centaine.

Il y a maintenant  \iref{nombreexos} exercices et \iref{nombreprob} \pbsz.

L'édition anglaise chez Springer en 2015 correspond à très peu près à cette version corrigée et augmentée française. Il manque cependant
dans l'édition anglaise les sections \ref{secconstrdi} et \ref{secAnor}
ainsi que quelques nouveaux exercices ou solutions d'anciens exercices.

 Toutes précisions utiles sur le site:

\url{http://hlombardi.free.fr/publis/LivresBrochures.html}

\vspace{5mm}
\begin{flushright}
Henri Lombardi, Claude Quitté\\
\today
\end{flushright}

%%%%%%%%%%%%%%%%%%%%%%%%%%%%%%
\cleardoublepage
\setcounter{tocdepth}{1}
\tableofcontents
%%%%%%%%%%%%%%%%%%%%%%%%%%%%%%
%%%%%%%%%%% AVANT PROPOS  %%%%%%%%%%%%%%%%%

%\chapter*{Avant-Propos}
~

\vskip1cm
\rdb
%\chapter*{\protect\LARGE Avant-propos}
{\LARGE \bf Avant-propos}
\pagestyle{CMExercicesheadings}
\thispagestyle{CMchapitre}%\mtcaddchapter
\addcontentsline{toc}{chapterbis}{Avant-propos}
\markboth{Avant-propos}{Avant-propos}
%\perso{compilé le \today}

\vskip1cm

\begin{flushright}
{\small
 Quant  à moi, je proposerais de s'en tenir
 aux règles
suivantes:

%-----------------begin enum------------------

\item 1. Ne jamais envisager que des objets susceptibles
d'être définis
\\
en un nombre fini de mots;
\\
\item 2. Ne jamais perdre de vue que toute proposition
 sur l'infini doit
\\
être la traduction, l'énoncé abrégé
de propositions sur le
fini;
\\
\item 3. Éviter les classifications et les définitions
non-prédicatives.

%-----------------end enum------------------

\medskip \rm Henri Poincaré,

dans  {\it La logique de l'infini }
(Revue de Métaphysique et de Morale, 1909).
\\
Réédité dans  {\it
Dernières pensées}, Flammarion.
}
\end{flushright}

\vspace{.3cm}

Ce livre est un cours d'introduction à l'\alg commutative de base, avec un accent
particulier mis sur les \mptfsz, qui constituent la version \agq des
fibrés vectoriels en \gmt \dilez.

Comme indiqué dans la préface, nous adoptons la méthode \covz, avec laquelle tous les \thos
d'existence ont un contenu \algq explicite.
Les \coma peuvent \^{e}tre regardées comme la partie la plus théorique du \calf
(computer algebra en anglais), qui s'occupe des \maths qui \gui{tournent
sur ordinateur}. 
Notre cours se distingue cependant des cours de \calf usuels sous deux aspects
essentiels.

Tout d'abord, nos \algos sont le plus souvent seulement implicites,
sous-jacents à la \demz, et ne sont en aucune manière
optimisés pour s'exécuter le plus rapidement possible,
comme il est naturel lorsque l'on vise une implémen\-tation efficace.

Ensuite, notre approche théorique est entièrement \covz, alors que les cours
de \calf usuels se préoccupent peu de cette question.
La philosophie n'est donc pas ici, comme il est d'usage
\gui{blanc ou noir, le bon chat est  celui qui attrape la
souris{\footnote{Proverbe chinois.}}},
mais plutôt
la suivante \gui{le moyen fait partie
de la recherche de la vérité, aussi bien que le résultat.
Il faut que la recherche de la vérité soit elle-même vraie;
la recherche vraie, c'est la vérité
déployée,  dont les membres épars se réunissent
dans le résultat{\footnote{Karl Marx, Remarques à propos de la récente instruction prussienne sur la censure, 1843 (cité par Georges Perec
dans \emph{Les Choses}).}}}.

Nous sommes amenés à parler
souvent des deux points de vue, classique et  \cofz, sur un m\^{e}me
sujet.
En particulier, nous avons mis une étoile pour signaler les
énoncés
(théor\`{e}mes, lemmes \dots) qui sont vrais en \clamaz, mais
dont nous ne donnons pas de \dem \covz, et qui souvent ne peuvent
pas en avoir.
Ces énoncés \gui{étoilés} ne seront donc 
probablement jamais
implémentés sur machine, mais ils sont bien souvent utiles
comme guides pour l'intuition, et au moins pour faire
le lien avec les exposés usuels écrits dans le
style des \clamaz.

 \stMF
Pour ce qui concerne les \dfnsz, nous donnons \gnlt en
premier une variante \covz, \llecz{\footnote{La personne qui lit ce
livre subit la règle inexorable de l'alternance des sexes.
Espérons que les lecteurs n'en seront pas plus affectés que
les lectrices. En tout cas, cela nous économisera bien des \gui{ou}
et bien des \gui{(e)}.}} voudra bien nous le pardonner,
quitte à montrer en \clama l'équivalence avec la
définition
usuelle. \Llec constatera que dans les démonstrations
\gui{étoilées} nous utilisons librement le lemme de Zorn et le
principe du tiers exclu, tandis que les autres \dems ont toujours
une traduction directe sous forme d'\algoz.

L'\alg \cov est en fait une vieille discipline,
développée en particulier par Gauss et Kronecker.
Comme précisé également dans la préface, nous nous situons dans la lignée de la \gui{bible} moderne sur le 
sujet, qu'est le livre {\it  A Course in
Constructive Algebra} de Ray Mines, Fred Richman et Wim Ruitenburg, paru en 1988.
Nous le citerons sous forme abrégée \cite{MRR}.
Notre ouvrage est cependant autocontenu et nous ne réclamons pas \cite{MRR}
comme prérequis.
Les livres de Harold M.~Edwards de \coma \cite{Ed,Ed2} et celui de Ihsen Yengui \cite{Yengui} sont aussi à recommander.

%%%%%%%%%%%%%
\subsection*{Le contenu de l'ouvrage}

Nous commençons par un bref commentaire sur les choix qui ont été faits
concernant les thèmes traités. 

La théorie des \mptfs est un des thèmes unificateurs de l'ouvrage. 
Nous voyons cette théorie sous forme abstraite comme une théorie \agq des fibrés vectoriels, et sous forme concrète comme celle des matrices \idmesz.
La comparaison des deux points de vue est esquissée dans le chapitre introductif.

La théorie des \mptfs proprement dite est traitée dans les 
chapitres~\ref{chap ptf0} (premières \prtsz),~\ref{chap AlgStricFi} (\algs qui sont des \mptfsz),~\ref{chap ptf1} (théorie du rang et exemples),~\ref{chapNbGtrs} (splitting off de Serre) et~\ref{ChapMPEtendus} (\mptfs étendus).

Un autre thème unificateur est fourni par les \plgsz, comme dans \cite{Kun} par exemple. Il s'agit d'un cadre conceptuel très efficace, même s'il est un peu vague. 
D'un point de vue \cofz, on remplace la \lon en un \idep arbitraire 
par un nombre fini de \lons en des  \mocoz. 
Les notions qui respectent le \plg sont \gui{de bonnes notions}, en ce sens qu'elles
sont mûres pour le passage des anneaux commutatifs aux schémas de Grothendieck,
que nous ne pourrons malheureusement pas aborder dans l'espace restreint de cet
ouvrage.

Enfin, un dernier thème récurrent est donné par la méthode, tout à fait familière en calcul formel, dite de \emph{l'\evn paresseuse}, ou  dans sa forme la plus aboutie, la méthode de \emph{l'\evn dynamique}. Cette méthode est indispensable lorsque l'on veut mettre en place un traitement \algq des questions qui requièrent a priori la solution d'un \pb de \fcnz. 
Cette méthode  a \egmt permis la mise au point des
machineries \covs \lgbes que l'on trouve dans les chapitres~\ref{chap mpf} 
et~\ref{chapPlg}, ainsi que celle de 
la théorie \cov de la \ddk (chapitre~\ref{chapKrulldim}), 
 avec d'importantes applications
dans les derniers chapitres.

Nous passons maintenant
à une description plus détaillée du contenu de l'ouvrage.

%: 1
\smallskip Dans le chapitre~\ref{chapMotivation}, 
nous expliquons les liens étroits que l'on peut établir entre 
les notions de fibrés vectoriels
en \gmt \dile et de \mptf en algèbre commutative. 
Ceci fait partie du processus \gnl d'algébrisation en \mathsz,
processus qui permet souvent de simplifier, d'abstraire et de généraliser 
de manière surprenante des concepts provenant de théories particulières.

%: 2
\smallskip Le chapitre~\ref{chapSli} est consacré aux \slis sur un 
anneau commutatif, traités sous forme \elrz. 
Il ne requiert presqu'aucun appareillage théorique,
mis à part la question de la \lon en un \moz, dont nous donnons 
un rappel dans la section~\ref{secPrelimCh2}.
Nous entrons ensuite dans notre sujet en mettant
  en place le \plgc pour la résolution des \slis
(section~\ref{secPLGCBasic}),
un outil simple et efficace qui sera repris et diversifié
sans cesse. D'un point de vue \cofz, la résolution des \slis fait \imdt apparaître comme central
le concept d'anneau \coh que nous traitons dans la 
section~\ref{secAnneauxCoherents}.
Les anneaux \cohs sont ceux pour lesquels on a une prise
minimale sur la solution des \slis homogènes.
De manière très étonnante, ce concept n'apparaît pas dans les
traités classiques d'\alg commutative. 
C'est qu'en général celle notion est complètement occultée par celle d'anneau \noez. Cette occultation n'a pas lieu en \coma o\`u la \noet n'implique
pas \ncrt la \cohcz. Nous développons  dans la section~\ref{sec sfio} 
la question des produits
finis d'anneaux, avec la notion de \sfio  et le \tho
des restes chinois.
La longue section~\ref{secCramer} est consacrée à de nombreuses
variations sur le thème des \detersz. Enfin, la section~\ref{secPLGCBasicModules}
revient sur le \plg de base, dans une version un peu plus \gnle consacrée
aux suites exactes de modules. 

%: 3
\smallskip Le chapitre~\ref{chapGenerique} développe la méthode des \coes
indéterminés, développée par Gauss. 
De très nombreux \thos d'existence en \alg commutative reposent sur
des \gui{\idas sous conditions} et donc sur des appartenances %du type 
$g\in\gen{\lfs}$ dans un anneau~$\ZZ[c_1,\ldots,c_r,\Xn]$, o\`u les $X_i$ sont les variables et les $c_j$ les paramètres du \tho considéré. En ce sens, on peut considérer que l'\alg commutative est une vaste théorie des \idasz, qui trouve son cadre naturel dans la méthode des \coes indéterminés, \cad la méthode
dans laquelle les paramètres du \pb à traiter sont pris comme des \idtrsz.
Forts de cette certitude, nous sommes, autant que faire se pouvait, 
systématiquement \gui{partis à la chasse des \idasz}, ceci non seulement dans les chapitres~\ref{chapSli} et~\ref{chapGenerique} 
\gui{purement calculatoires}, mais dans tout l'ouvrage. 
En bref, plutôt que d'affirmer en filigrane d'un \tho d'existence 
\gui{il existe une \ida qui certifie cette existence}, nous avons t\^aché 
de donner chaque fois l'\ida elle-même. 

%: 3 suite
Ce chapitre~\ref{chapGenerique} peut être considéré comme un cours 
d'algèbre de base avec les méthodes du 19\ieme\ siècle.
Les sections~\ref{secAnnPols},~\ref{secLemArtin} et~\ref{secThKro} donnent 
quelques  \gnts sur les  \polsz,
avec notamment l'\algo de \fcn partielle, la \gui{théorie des \idasz}
(qui explique la méthode des \coes indéterminés), les \pols \smqs \elrsz, le lemme de \DKM et le \tho de \KROz.  Ces deux 
derniers résultats sont des outils de base qui donnent des informations  précises sur les \coes du produit de deux \polsz;
ils sont souvent utilisés dans le reste de l'ouvrage.
La section~\ref{sec0adu} introduit l'\adu d'un \pol \unt sur un anneau
commutatif arbitraire, qui est un substitut efficace au corps des racines
d'un \pol sur un corps. La section~\ref{secDisc} est consacrée au discriminant
et explique en quel sens précis une matrice \gnq est \digz.
Avec ces outils en mains, on peut traiter la théorie de Galois de base
dans la section~\ref{secGaloisElr}.
La théorie \elr de l'élimination via le résultant est donnée dans 
la section~\ref{secRes}. On peut alors donner les bases de la théorie \agq
des nombres, avec le \tho de \dcn unique en facteurs premiers pour un
\itf d'un corps de nombres 
(section~\ref{secApTDN}). La section~\ref{secChap3Nst} donne le \nst de Hilbert comme application du résultant.
Enfin, la section~\ref{secNewton} sur la méthode de Newton en \alg termine 
ce chapitre.

%: 4
\smallskip Le chapitre~\ref{chap mpf} est consacré à l'étude des \prts
\elrs des \mpfsz. 
Ces modules jouent un peu le m\^{e}me rôle pour les anneaux que les \evcs de 
dimension finie pour les corps: la théorie des \mpfs est une manière un peu
 plus abstraite, et souvent profitable, d'aborder la question des \slisz.
Les sections~\ref{sec pf chg} à~\ref{secStabPf} donnent les \prts de stabilité de base ainsi que l'exemple important de l'\id d'un zéro pour un \syp (sur
un anneau commutatif arbitraire). On s'intéresse ensuite au \pb de classification
des \mpfs sur un anneau donné. Sur le chemin des anneaux principaux,
pour lesquels le \pb de classification est complètement résolu (section~\ref{secBézout}), nous rencontrons les anneaux \qis (section~\ref{subsecAnneauxqi}),
qui sont les anneaux o\`u l'annulateur d'un \elt est toujours engendré 
par un \idmz. 
C'est l'occasion de mettre en place une \emph{machinerie \lgbe \elrz}  
qui permet de passer d'un résultat
établi \cot pour les anneaux intègres au même résultat, convenablement reformulé, pour les anneaux \qisz\imlgz. 
Cette machinerie de transformation
de preuves est \elrz, car fondée sur la \dcn d'un anneau en produit fini d'anneaux.
La chose intéressante est que cette \dcn est obtenue par relecture
de la \dem \cov écrite dans le cas intègre: on voit ici qu'en \coma la \dem
est souvent encore plus importante que le résultat.  
De la même manière, on a une machinerie \lgbe \elr qui permet de passer d'un résultat
établi \cot pour les \cdis au même résultat, convenablement reformulé, pour les anneaux \zedrs (section~\ref{secKrull0dim}). Les anneaux \zedsz, ici définis
de manière \elrz, constituent une clé  importante de l'\alg commutative,
comme étape intermédiaire pour \gnr certains résultats des \cdis aux anneaux commutatifs arbitraires.  Dans la littérature classique, ils apparaissent
souvent sous leur forme \noeez, \cad celle des anneaux artiniens.
La section~\ref{sec Fitt} introduit les
invariants très importants que sont les \idfs d'un \mpfz. Enfin,  la section
\ref{subsecIdealResultant} applique cette notion pour introduire l'\id résultant
d'un \itf dans un anneau de \pols  quand l'\id en
question contient un \pol \untz, et démontrer un  \tho d'\eli \agq sur un anneau arbitraire.

%: 5
\smallskip Le chapitre~\ref{chap ptf0} est une première approche de la théorie des \mptfsz. Les sections~\ref{subsecPropCarPTFS} à~\ref{secPtfCoNat} donnent les \prts  de base ainsi que l'exemple important des anneaux \zedsz.
La section~\ref{secMPTFlocLib} donne le \tho de structure locale: un module est \ptf \ssi il devient
libre après \lon en des \eco convenables. Sa \dem \cov est une relecture
d'un résultat établi dans le chapitre~\ref{chapSli} pour les \slis \gui{bien conditionnés} (\thref{theoremIFD}).
La section~\ref{Idpp} développe l'exemple des modules  \pros \lmosz. La section
\ref{subsec det ptf} introduit le \deter d'un \endo d'un \mptfz. 
Ceci donne accès à la \dcn d'un tel module en 
somme directe de ses composants de rang constant.
Enfin, la section~\ref{secSalutFini}, que l'on ne savait pas bien o\`u mettre dans l'ouvrage,
héberge quelques considérations \suls sur les \emph{\prts \carfz}, une notion
introduite au chapitre~\ref{chapSli} pour discuter les rapports entre
\plgs concrets et \plgs abstraits.

%: 6
\smallskip Le chapitre~\ref{chap AlgStricFi} est consacré pour l'essentiel
aux \algs qui sont des \mptfs sur leur anneau de base.  
Nous les appelons
des \asfsz. Elles constituent
une \gnn naturelle pour les anneaux commutatifs
de la notion d'\alg finie sur un corps. Comme cas important, cerise sur le g\^{a}teau, les \aGsz,
qui généralisent les extensions galoisiennes de \cdis aux anneaux commutatifs.

La section~\ref{secEtaleSurCD} traite le cas o\`u l'anneau de base est un \cdiz.
Elle donne des versions \covs pour les \thos de structure obtenus en \clamaz.
Le cas des \algs étales (quand le \discri est \ivz) est particulièrement éclairant. On découvre que les \thos classiques supposent toujours implicitement que l'on sache factoriser les \pols \spls sur le corps de base.
La \dem \cov du \tho de l'\elt primitif \ref{thEtalePrimitif} est significative par son
écart avec la \dem classique.
La section~\ref{sec2GaloisElr} applique les résultats précédents
pour terminer la théorie de Galois de base commencée dans la section~\ref{secGaloisElr} en caractérisant les extensions galoisiennes de \cdis comme les extensions
étales et normales.
La section~\ref{sec1Apf} est une brève introduction aux \apfsz, en insistant
sur le cas des \algs entières, avec un \nst faible et le lemme lying over.
La section~\ref{subsecAlgStfes} introduit les \algs \stfes sur un anneau arbitraire. 
Dans les sections~\ref{secAlgSte} et~\ref{secAlgSpb}, sont  introduites
les notions voisines d'\alg \ste et d'\alg \spb
qui généralisent la notion d'\alg étale sur un \cdiz. 
Dans la  section~\ref{secAGTG}, on donne un exposé \cof des bases
de la théorie des \aGs
pour les anneaux commutatifs.
Il s'agit en fait d'une théorie d'Artin-Galois, puisqu'elle reprend
l'approche qu'Artin avait développée pour le cas des corps, en partant directement d'un groupe fini d'\autos d'un corps, le corps de base n'apparaissant
que comme un sous-produit des constructions qui s'ensuivent.

%: 7
\smallskip Dans le chapitre~\ref{ChapGalois},
la méthode dynamique, une pierre angulaire des méthodes modernes
en \alg \covz,  est 
mise en {\oe}uvre  pour traiter d'un point de vue \cof le corps des racines d'un \pol et la théorie de Galois dans le cas \spbz, lorsque la proie s'échappe pour laisser place à son ombre, \cad lorsque l'on ne sait pas factoriser les \pols sur le corps de base que l'on considère.
\`A titre d'entraînement, la section~\ref{secNstSCA} commence par établir des résultats sous forme \cov pour le \nst lorsque
l'on  ne sait pas factoriser les \pols sur le corps de base.
Des considérations d'ordre \gnl sur la méthode dynamique sont développées
dans la section~\ref{subsecDyna}.
Plus de détails sur le déroulement des festivités  sont donnés dans l'introduction du chapitre.

%: 8
\smallskip Le chapitre~\ref{chap mod plats} est une brève introduction
aux modules plats et aux \algs plates et \fptesz.
En langage intuitif, une \Alg $\gB$ est plate  lorsque les
\slis sur $\gA$ sans second membre n'ont \gui{pas plus} de solutions dans 
$\gB$ que dans $\gA$, et elle est fidèlement plate si cette affirmation
est vraie \egmt des \slis avec second membre. 
Ces notions cruciales de l'\alg commutative ont été
introduites par Serre dans~\cite[GAGA,1956]{Serregaga}. Nous ne donnons que les
résultats vraiment fondamentaux. C'est \egmt l'occasion d'introduire les
notions d'anneau \lsdzz, de module sans torsion (pour un anneau arbitraire),
d'\anar et d'anneau de Pr\"ufer. Nous insistons comme toujours sur
le \plg quand il s'applique.

%: 9
\smallskip Le chapitre~\ref{chap Anneaux locaux} parle des anneaux locaux et de
quelques \gnnsz. La section~\ref{secAloc1} introduit la terminologie \cov
pour quelques notions classiques usuelles, dont la notion importante de radical
de Jacobson. Une notion connexe est celle d'anneau \plc (un anneau $\gA$
tel que $\gA/\Rad\gA$ est \zedz). C'est une notion robuste, 
qui n'utilise jamais les \idemasz,
et la plupart des \thos de la littérature concernant
les anneaux semi-locaux (en \clama ce sont les anneaux qui n'ont
qu'un nombre fini d'\idemasz) s'appliquent aux anneaux \plcsz.
La section~\ref{secAloc2} répertorie quelques résultats qui montrent
que sur un \alo on ramène la solution de certains \pbs au cas des corps.
Les sections~\ref{secLoc1+fa} et~\ref{secExlocGeoAlg} établissent
sur des exemples \gmqs (\cad concernant l'étude de \sypsz) 
un lien entre la notion d'étude locale au sens
intuitif topologique et l'étude de certaines \lons d'anneaux (dans le cas d'un \cdi
à la base, ces \lons sont des \alosz). On introduit notamment
les notions d'espaces tangent et cotangent en un zéro d'un \sypz.
La section~\ref{secRelIdm} fait une brève étude des anneaux \dcpsz,
dont un cas particulier en \clama sont les anneaux décomposés (produits 
finis d'anneaux locaux), 
qui jouent un rôle important dans la théorie des \alos henséliens.
Enfin la section~\ref{secAlocglob} traite la notion d'anneau \lgbz,
qui \gns à la fois celle d'\alo et celle d'anneau \zedz. Ces anneaux vérifient des \prts \lgbes très fortes,
par exemple les \mrcs  sont toujours libres, et ils sont stables
par extensions entières.

%: 10
\smallskip Le chapitre~\ref{chap ptf1}
poursuit l'étude des \mptfs commencée dans  le chapitre~\ref{chap ptf0}.
Dans la section
\ref{sec ptf loc lib}, nous reprenons la question de la
\carn  des \mptfs comme modules \lot libres, \cad du \tho de structure locale.
Nous en donnons une version matricielle (\thref{th matproj}), 
qui résume et précise les différents énoncés du \thoz.
La section~\ref{subsecCalRang} est consacrée à l'anneau des rangs sur~$\gA$.
En \clamaz, le rang d'un \mptf est défini
comme une fonction localement constante sur le spectre de Zariski.
Nous donnons ici une théorie \elr du rang qui ne fait pas appel aux \idepsz.
Dans la section~\ref{secAppliLocPtf},
nous donnons quelques applications simples du \tho
de structure locale.
La section~\ref{secGrassman} est une introduction aux grassmanniennes.
Dans la section~\ref{subsecClassifMptfs}, nous introduisons le \pb \gnl de
la classification complète des \mptfs sur un anneau $\gA$ fixé.
Cette classification est un \pb fondamental et difficile, qui n'admet pas de
solution \algq \gnlez.
La section~\ref{secAppliIdenti}
présente un exemple non trivial pour
lesquels cette classification peut être obtenue.

%: 11
\smallskip Le chapitre~\ref{chapTrdi} est consacré aux \trdis et groupes réticulés. Les deux premières sections décrivent ces structures \agqs ainsi que leurs
\prts de base. Ces structures sont importantes en \alg commutative pour plusieurs raisons.

D'une part, la théorie de la \dve a comme \gui{modèle idéal} la théorie de la \dve des entiers naturels. La structure du \mo 
multiplicatif~\hbox{$(\NN\etl,\times ,1)$}
en fait la partie positive d'un \grlz. Ceci se généralise en \alg commutative dans deux directions. La première \gnn est la théorie des anneaux intègres
 dont les \itfs forment un \trdiz, appelés des \ddpsz,
que nous étudierons dans le chapitre~\ref{ChapAdpc}:
leurs \itfs non nuls forment  la partie positive d'un \grlz.
La deuxième est la théorie des anneaux à pgcd 
que nous étudions dans la section~\ref{secAnnPgcd}.
Signalons la première apparition de la dimension de Krull $\leq 1$
dans le \thref{propGCDDim1}: un anneau à pgcd intègre de dimension $\leq 1$
est un anneau de Bézout.

D'autre part, les \trdis interviennent comme la contrepartie \cov
des espaces spectraux divers et variés qui se sont imposés comme
des outils puissants de l'\alg abstraite.
Les rapports entre \trdis et espaces spectraux seront abordés 
dans la section~\ref{secEspSpectraux}.
Dans la section~\ref{secZarAcom}, nous mettons en place le treillis de Zariski d'un anneau commutatif $\gA$, qui est la contrepartie \cov du fameux spectre de Zariski. Notre but ici est d'établir le parallèle entre la 
construction de la clôture \zede réduite d'un anneau (notée $\Abul$) 
et celle de l'\agB
engendrée par un \trdi (qui fait l'objet du \thref{thZedGenEtBoolGen}).
L'objet  $\Abul$ ainsi construit contient
essentiellement la même information que le produit des anneaux
$\Frac(\gA\sur\fp)$ pour tous les \ideps $\fp$ de $\gA$(\footnote{Ce produit n'est pas accessible en \comaz, $\Abul$ en est un substitut \cof tout à fait efficace.}).
Ce résultat est en relation étroite avec le fait que le treillis de Zariski de $\Abul$ est l'\agB engendrée par
le treillis de Zariski de $\gA$.

Une troisième raison de s'intéresser aux \trdis est la logique \cov
(ou intuitionniste). Dans cette logique, l'ensemble des valeurs de vérité de la logique classique, à savoir $\so{\Vrai,\Faux}$,
qui est une \agB à deux \eltsz, est remplacé par un \trdi
assez mystérieux.
La logique \cov sera abordée de manière informelle dans l'annexe.
Dans la section~\ref{secEntRelAgH}, 
nous mettons en place les outils qui servent
de cadre à une étude \agq formelle de la logique \covz:
les \entrels et les \agHsz.
Par ailleurs, \entrels et \agHs ont leur utilité propre dans l'étude \gnle des
\trdisz. Par exemple, le treillis
de Zariski d'un anneau \noe \coh  est une \agH (proposition~\ref{propNoetAgH}).

%: 12
\smallskip Le chapitre~\ref{ChapAdpc} traite les \anarsz, les \adps et
les \adksz. Les \anars sont les anneaux dont le treillis des \itfs est distributif.
Un \adp est un \anar réduit et il est caractérisé par le fait que tous ses \ids sont plats. Un \adpc est la même chose qu'un \anar \qiz. Il est caractérisé par les fait que ses \itfs sont projectifs. Un \adk est un \adpc \noe et \fdi
(en \clama avec le principe du tiers exclu tout anneau est \fdi et tout anneau \noe est \cohz). Ces anneaux sont apparus tout d'abord avec les anneaux d'entiers de corps de nombres. Le paradigme dans le cas intègre est la \dcn unique en facteurs premiers de tout \itf non nul. Les \prts \aris du \mo
multiplicatif des \itfs sont pour l'essentiel vérifiées par les \anarsz.
Pour les \prts les plus subtiles concernant la \fcn des \itfsz, et notamment la \dcn en facteurs premiers, une hypothèse \noeez, ou au moins de dimension~$\leq 1$, est indispensable.
Dans ce chapitre, nous avons voulu montrer la progression des \prts satisfaites
par les anneaux au fur et à mesure que l'on renforce les hypothèses,
depuis les \anars jusqu'aux \adks à \fcn totale. 
Nous insistons sur le caractère
\algq simple des \dfns dans le cadre \cofz. 
Certaines \prts ne dépendent que de la dimension $\leq 1$, et nous avons voulu rendre justice aux anneaux \qis \ddi1. Nous avons \egmt fait une étude
du \pb de la \dcn en facteurs premiers plus progressive et plus fine que
dans les exposés qui s'autorisent le principe du tiers exclu.
Par exemple, les \thosz~\ref{thAESTE} et~\ref{lemthAESTE} donnent des versions \covs précises du \tho concernant
les extensions finies normales d'\adksz, avec ou sans la \prt de \fcn totale.

Le chapitre commence par quelques remarques d'ordre épistémologique
sur l'intérêt intrinsèque d'aborder les \pbs de \fcn avec le \tho de \fcn partielle plut\^{o}t qu'avec celui de \fcn totale. Pour avoir une bonne idée du déroulement des festivités, il suffit de se reporter à la 
table des matières en tête
du chapitre \paref{ChapAdpc} et à la table des \thos \paref{tdtChapAdpc}.

%: 13
\smallskip Le chapitre~\ref{chapKrulldim} est consacré à la \ddk des anneaux commutatifs, à celle de leurs morphismes, à celle des \trdis et à la dimension valuative des anneaux commutatifs. 

Plusieurs notions importantes de dimension en algèbre commutative classique sont des dimensions d'espaces spectraux. Ces espaces
topologiques très particuliers jouissent de la \prt d'être entièrement
décrits (au moins en \clamaz) par leurs \oqcsz, qui forment un \trdiz. Il s'avère
que le \trdi correspondant a en \gnl une interprétation simple, sans recours aucun aux espaces spectraux. En 1974, Joyal a montré comment définir \cot la \ddk d'un \trdiz. Depuis ce jour faste, la théorie de la dimension qui semblait baigner dans des espaces éthérés, invisibles lorsque l'on ne fait pas confiance à l'axiome du choix, est devenue (au moins en principe) une théorie de nature \elrz, sans plus aucun mystère.

La section~\ref{secEspSpectraux} décrit l'approche de la \ddk en \clamaz.
Elle explique aussi 
comment on peut interpréter la \ddk d'un tel espace en terme
du \trdi de ses \oqcsz. 
La section~\ref{secDefConsDimKrull} donne la \dfn \cov de la \ddk
d'un anneau commutatif, notée $\Kdim\gA$, et en tire quelques conséquences.
La section~\ref{secKrullElem} donne quelques \prts plus avancées, et notamment
le \plg et le principe de recouvrement fermé pour la \ddkz. 
La section~\ref{secDDKExtEn} traite la \ddk des extensions entières et la section~\ref{secDimGeom} 
celle des anneaux \gmqs (correspondant aux \sypsz) sur les \cdisz.
La section~\ref{secDDKTRDIS} donne la \dfn \cov de la \ddk d'un \trdi et
montre que la \ddk d'un anneau commutatif et celle de son treillis de Zariski
co\"{\i}ncident.
La section~\ref{secKdimMor} est consacrée à la dimension 
des morphismes entre anneaux
commutatifs. La \dfn utilise la clôture \zed réduite de l'anneau source du morphisme. Pour démontrer la formule qui majore $\Kdim \gB$ à partir de
$\Kdim \gA$ et $\Kdim\rho$ (lorsque l'on a un morphisme $\rho:\gA\to\gB$),
nous devons introduire la clôture \qi minimale d'un anneau commutatif.
Cet objet est une contrepartie \cov du produit de tous les $\gA/\fp$,
lorsque $\fp$ parcourt les \idemis de~$\gA$.
La section~\ref{secValdim} introduit la dimension valuative d'un anneau commutatif
et utilise cette notion notamment pour démontrer le résultat important suivant:
pour un \anar non nul $\gA$, on a $\Kdim\AXn=n+\Kdim\gA$.
La section~\ref{secGoingUp} donne des versions \covs des \thos Going up et Going down.

%: 14
\smallskip Dans le chapitre~\ref{chapNbGtrs}, intitulé \emph{Nombre de \gtrs d'un module}, on établit la version \elrz,  non
\noee et \cov de \gui{grands} \thos d'\alg
commutative, dus dans leur version originale
à \KRNz, Bass, Serre, Forster et Swan.  Ces résultats concernent
le nombre de \gtrs radicaux d'un \itfz,
le nombre de \gtrs d'un module, la possibilité de produire
un sous-module libre en facteur direct dans un module,
et la possibilité de simplifier des \isosz,
dans le style suivant: si $M\oplus N\simeq M'\oplus N$ alors $M\simeq M'$.
Ils font intervenir la \ddk ou d'autres dimensions plus sophistiquées,
introduites par R. Heitmann ainsi que par les auteurs de cet ouvrage et T. Coquand.

La section~\ref{secKroBass} est consacrée au \tho de Kronecker et à ses extensions
(la plus aboutie, non \noeez, est due à R. Heitmann \cite{Hei84}).
Le \tho de \KRN est usuellement énoncé
sous la forme suivante: une
\vgq dans $\CC^n$ peut toujours \^{e}tre définie par $n+1$
\eqnsz. La forme due à Heitmann est que dans un anneau de dimension de Krull inférieure ou égale à $n$, pour tout \itf $\fa$ il existe un \id $\fb$ engendré par au plus $n+1$ \elts 
de $\fa$ tel que $\sqrt \fb=\sqrt\fa$.
La \dem donne aussi le \tho de Bass, dit \gui{stable range}. 
Ce dernier \tho a été amélioré en faisant intervenir des dimensions \gui{meilleures} que la \ddkz. Ceci fait l'objet de la section~\ref{subsecDimHeit},
o\`u est définie la \emph{dimension de Heitmann}, 
découverte en lisant avec attention
les \dems de Heitmann (Heitmann utilise une autre dimension,
a priori un peu moins bonne, que nous expliquons \egmt en termes \cofsz).
Dans la section \ref{secSOSFSa}, nous expliquons quelles sont les \prts 
matricielles d'un anneau qui permettent de faire fonctionner
les \thos de Serre (splitting off), de Forster-Swan 
(contrôle du nombre de \gtrs d'un \mtf en fonction du nombre
de \gtrs local) et le \tho de simplification de Bass.
La section~\ref{secSUPPORTS} introduit les notions de support (une application d'un anneau dans un \trdi vérifiant certains axiomes) et de $n$-stabilité.
Cette dernière notion a été définie   par Thierry Coquand,
après avoir analysé une \dem de Bass qui établit que les
\mptfs sur un anneau~$\gV[X]$, o\`u~$\gV$ est un anneau de valuation
de \ddk finie, sont libres.
Dans la dernière section, on démontre que la \prt matricielle cruciale
introduite dans la section~\ref{secSOSFSa} est satisfaite, d'une part, par les anneaux $n$-stables, d'autre part par les anneaux de dimension de Heitmann $<n$.

%: 15
\smallskip Le chapitre~\ref{chapPlg} est consacré au \plg et à ses variantes.
La section~\ref{subsecMoco} introduit la notion de recouvrement d'un \mo
par une famille finie de \mosz, ce qui \gns la notion de \mocoz. 
Le lemme de recouvrement~\ref{lemRecouvre} sera décisif 
dans la section~\ref{secMachLoGlo}.
La section~\ref{subsec loc glob conc} donne des 
\plgcsz. Il s'agit de dire que certains \prts
sont vraies globalement
dès qu'elles le sont localement.
Ici, \og localement\fg\/ est pris au sens \cofz: après \lon en
un nombre fini de \mocoz. La plupart des résultats
 ont été
établis dans les chapitres précédents. Leur regroupement 
fait voir la portée très générale de ces principes.
La section~\ref{subsec loc glob abs} reprend certains de
ces principes sous forme de \plgasz. Ici, \og localement\fg\/ est pris au sens
abstrait, \cad après \lon en n'importe quel \idepz. 
C'est surtout la comparaison avec les \plgcs correspondants qui nous
intéresse.
La section~\ref{secColleCiseaux} explique la construction d'objets
\gui{globaux} à partir d'objets de même nature définis uniquement
de manière locale, comme il est usuel en \gmt \dilez.
C'est l'impossibilité de cette construction lorsque l'on cherche
à recoller certains anneaux qui est à l'origine des schémas de Grothendieck.
En ce sens, les sections~\ref{subsec loc glob conc} et~\ref{secColleCiseaux}
constituent la base à partir de laquelle on peut développer la théorie
des schémas dans un cadre complètement \cofz.

Les sections suivantes sont d'une autre nature. D'ordre méthodologique, elles
sont consacrées au décryptage de différentes variantes 
du \plg en \clamaz.
Par exemple, la \lon en tous les \idepsz, le passage au quotient par tous les \idemas ou la \lon en tous les \idemisz, qui s'appliquent
chacune dans des situations particulières.
Un tel décryptage  présente un caractère certainement déroutant
dans la mesure o\`u il prend pour point de départ une \dem classique
qui utilise des \thos en bonne et due forme, mais o\`u le décryptage
\cof de cette \dem n'est pas seulement donné par l'utilisation 
de \thos \cofs en bonne et due forme. 
Il faut aussi regarder ce que fait la \dem classique
avec ses objets purement idéaux (des \idemas par exemple) pour 
comprendre comment elle
nous donne le moyen de construire un nombre fini d'\elts qui
vont être impliqués dans un \tho \cof (un \plgc par exemple)
pour aboutir au résultat souhaité.
En décryptant une telle \demz, nous utilisons la méthode dynamique \gnle 
exposée au chapitre~\ref{ChapGalois}. Nous décrivons ainsi 
des \emph{machineries locales-globales} nettement moins \elrs
que celles du chapitre~\ref{chap mpf}: la machinerie \lgbe \cov de base  \gui{à \idepsz} (section~\ref{secMachLoGlo}), la  machinerie \lgbe \cov à \idemas  (section~\ref{subsecLGIdeMax}) et la  machinerie \lgbe \cov à 
\idemis  (section~\ref{subsecLGIdepMin})\imlbz.  
En réalisant \gui{le programme de Poincaré} cité en exergue de cet
avant-propos, nos machineries \lgbes prennent en compte une remarque essentielle
de Lakatos, à savoir que la chose la plus intéressante et robuste dans un \thoz, 
c'est toujours sa \demz,
même si elle est critiquable à certains égards (voir \cite{La}).     

Dans les sections \ref{secPlgcor} et \ref{secPlgprof2}, nous examinons
dans quelle mesure certains \plgs restent valides lorsque l'on remplace 
dans les énoncés les listes d'\eco 
par des listes de profondeur $\geq 1$ ou de profondeur $\geq 2$.

%: 16
\smallskip Dans le chapitre~\ref{ChapMPEtendus}, nous traitons la question
des \mptfs sur les anneaux de \polsz. La question décisive est d'établir
pour quelles classes d'anneaux les modules \ptfs sur un anneau de \pols
proviennent par \eds d'un module \ptf sur l'anneau lui-même (éventuellement
en posant certaines restrictions sur les \mptfs considérés ou sur le nombre de variables dans l'anneau de \polsz).
Quelques \gnts sur les modules étendus sont données dans la section~\ref{sec etendus}.
Le cas des \mrcsz~$1$, complètement éclairci par le \tho de
Traverso-Swan-Coquand, est traité dans la section~\ref{sec.Traverso}.
La \dem \cov de Coquand utilise de manière cruciale la  machinerie \lgbe \cov à 
\idemisz.
La section~\ref{subsecQPatch} traite les \thos 
de recollement de Quillen (Quillen patching) et Vaserstein,
qui disent que certains objets sont obtenus par \eds 
(depuis l'anneau de base à un anneau de \polsz) \ssi cette \prt 
est vérifiée \lotz. 
Nous donnons aussi une sorte de réciproque du Quillen patching,
due à Roitman, sous forme \covz. 
La section~\ref{sec.Horrocks} est consacrée aux \thos de Horrocks.
La \dem \cov du \tho de Horrocks global fait subir à la \dem du \tho de Horrocks local la machinerie \lgbe de base
et se conclut avec le Quillen patching \cofz. 
La section~\ref{sec.QS} donne plusieurs \prcos du \tho de Quillen-Suslin
(les \mptfs sur un anneau de \pols sur un \cdi sont libres), fondées 
sur différentes \dems classiques.
La section~\ref{sec.Etendus.Valuation} établit le \tho de Lequain-Simis
 (les \mptfs sur un anneau de \pols sur un \anar sont étendus). La \dem utilise la méthode dynamique exposée au
chapitre~\ref{ChapGalois}, cela permet d'établir le \tho d'induction de Yengui,
 une variante \cov de l'induction de Lequain-Simis.

%: 17
\smallskip Dans le chapitre~\ref{ChapSuslinStab},
nous démontrons le \gui{Suslin Stability Theorem} dans le cas particulier
des \cdisz.
Ici aussi, pour obtenir une \dem \covz, nous utilisons la machinerie \lgbe
 de base, exposée au chapitre~\ref{chapPlg}.

%: Annexe
\smallskip L'annexe décrit la théorie des ensembles \cov à la Bishop.
Elle peut être vue comme une introduction à la logique \covz.
On y explique la sémantique de Brouwer-Heyting-Kolmogorov pour les
connecteurs et quantificateurs. On discute certaines formes faibles
du principe du tiers exclu ainsi que plusieurs principes problématiques
en \comaz.

%%%%%%%%%%%%%%%%%%%%%%%%%%%%%%%%%%%%%%%%%%%%%%%%%%%%%%%%%%%%%%%%%%%%%%%%%%%
\subsection*{Quelques remarques d'ordre épistémologique}

{\spaceskip 0.6ex plus 0.3ex minus 0.1ex
Nous espérons dans cet ouvrage montrer que des livres classiques d'\alg
commutative comme  \cite{Atiyah}, \cite{Eis}, \cite{Gil}, \cite{Gla},  \cite{Kapl}, \cite{Knap1}, \cite{Knap2}, \cite{Kun}, \cite{Laf},  \cite{Lam06}
(dont la lecture est vivement recommandée), \cite{Mat},
\cite{Nor}, ou même \cite{Bou} 
%:h2014 rajout Stacks
et le remarquable ouvrage disponible sur le réseau \cite{Stacks},
pourront entièrement \^{e}tre récrits
avec un point de vue \cofz, dissipant le voile de mystère qui
entoure les \thos d'existence non explicites des \clamaz.
Naturellement, nous espérons que les \lecs  profiteront de notre ouvrage pour
jeter un regard nouveau sur les livres de calcul formel classiques,
comme par exemple \cite{CLS}, \cite{KrRo}, \cite{Singular}, \cite{EnHe}, \cite{ElkMo}, \cite{Mora}, \cite{CCS} ou\linebreak \cite{vzGaGe}.

}

\smallskip 
Dans la mesure o\`{u} nous voulons un traitement \algq
de l'algèbre commutative, nous ne pouvons pas utiliser
toutes les facilités que donnent l'usage systématique du
lemme de Zorn et du principe du tiers exclu en \clamaz.
Sans doute,  \llec comprend bien qu'il est difficile
d'implé\-menter le lemme de Zorn en calcul formel.
Le refus du principe du tiers exclu doit par contre lui sembler plus
dur à avaler.
Ce n'est de notre part qu'une constatation pratique.
Si dans une \dem classique, vous trouvez un raisonnement qui conduit
à un calcul du type:  \gui{si $x$  est inversible, faire ceci,
sinon faire cela}, il est bien clair que cela ne se traduit
directement sous forme d'un \algo que dans le cas o\`{u} l'on dispose
d'un test d'inversibilité dans l'anneau en question.
C'est pour insister sur cette difficulté, que nous devons
contourner en permanence, que nous sommes amenés à parler
souvent des deux points de vue, classique et  \cofz, sur un m\^{e}me
sujet.

On peut discuter indéfiniment
pour savoir si les \coma sont une partie des \clamaz, la partie
qui s'occupe exclusivement des aspects explicites des choses, ou si au
contraire ce sont les \clama qui sont une partie des \comaz, la
partie dont les \thos sont \gui{étoilés}, \cad qui rajoutent
systématiquement dans leurs hypothèses le principe du tiers exclu
et l'axiome du choix.
Un de nos objectifs est de faire pencher la balance dans la deuxième
direction, non par le débat philosophique, mais par la pratique.

Signalons enfin deux traits marquants de cet ouvrage par 
rapport aux ouvrages classiques d'\alg commutative.

Le premier est la mise au second plan de la
\noetz.
L'expérience prouve en effet que la \noet est bien souvent une
hypothèse trop forte, qui cache la vraie nature \algq des choses.
Par exemple, tel \tho habituellement énoncé pour les anneaux \noes et les \mtfsz,
lorsque l'on met sa \dem à plat pour en extraire un \algoz,
s'avère \^{e}tre un \tho sur les anneaux \cohs et les \mpfsz.
Le \tho habituel n'est qu'un corolaire du bon \thoz, mais avec deux
arguments non \cofs qui permettent de déduire en \clama la \cohc et
la présentation finie de
la \noet et du type fini.
Une \dem dans le cadre plus satisfaisant de la \cohc et
des \mpfs se trouve bien souvent  déjà publiée
dans des articles de recherche, quoique rarement sous forme
entièrement \covz, mais
\gui{le bon énoncé} est en \gnl absent dans les ouvrages de
référence{\footnote{Cette déformation professionnelle \noee a produit un travers
linguistique dans la littérature anglaise qui consiste à prendre \gui{local ring}
dans le sens de \gui{Noetherian local ring}.}}.

Le deuxième trait marquant de l'ouvrage est l'absence presque totale de la négation
dans les énoncés \cofsz. 
Par exemple, au lieu d'énoncer que pour un anneau $\gA$  non trivial, 
deux modules libres de rang $m$ et $n$ avec~$m>n$ ne peuvent pas être isomorphes, 
nous préférons dire,
sans aucune hypothèse sur l'anneau, que si ces modules sont isomorphes, alors l'anneau est trivial (proposition~\ref{propDimMod1}). Cette nuance peut sembler bien mince au premier abord, mais elle
a une importance \algqz. Elle va permettre de remplacer une \dem en \clama
utilisant un anneau
$\gA=\gB\sur\fa$,  qui conclurait que~$1\in\fa$ au moyen d'un raisonnement par l'absurde, par une \dem pleinement
\algq qui construit~$1$ en tant qu'\elt de l'\id $\fa$ à partir d'un \iso entre $\Ae m$ et~$\Ae n$.

Pour une présentation \gnle des idées qui ont conduit aux nouvelles méthodes utilisées en \alg \cov dans cet ouvrage, on pourra lire
l'article de synthèse \cite[Coquand\&Lombardi, 2006]{CL05}.

\vspace{5mm}
\begin{flushright}
Henri Lombardi, Claude Quitté\\
Mai 2014
\end{flushright}

%\rdb
%\newpage \thispagestyle{CMcadreseul}
%\cleardoublepage %\thispagestyle{CMcadreseul}
\newpage
%\ifodd\arabic{page}\relax\else\null\thispagestyle{CMcadreseul}\newpage\fi
\ifodd\arabic{page}\null\thispagestyle{CMcadreseul}\newpage\else\relax\fi

%:HHH organigrame
\null\vfil
\centerline{\includegraphics[width=8cm]{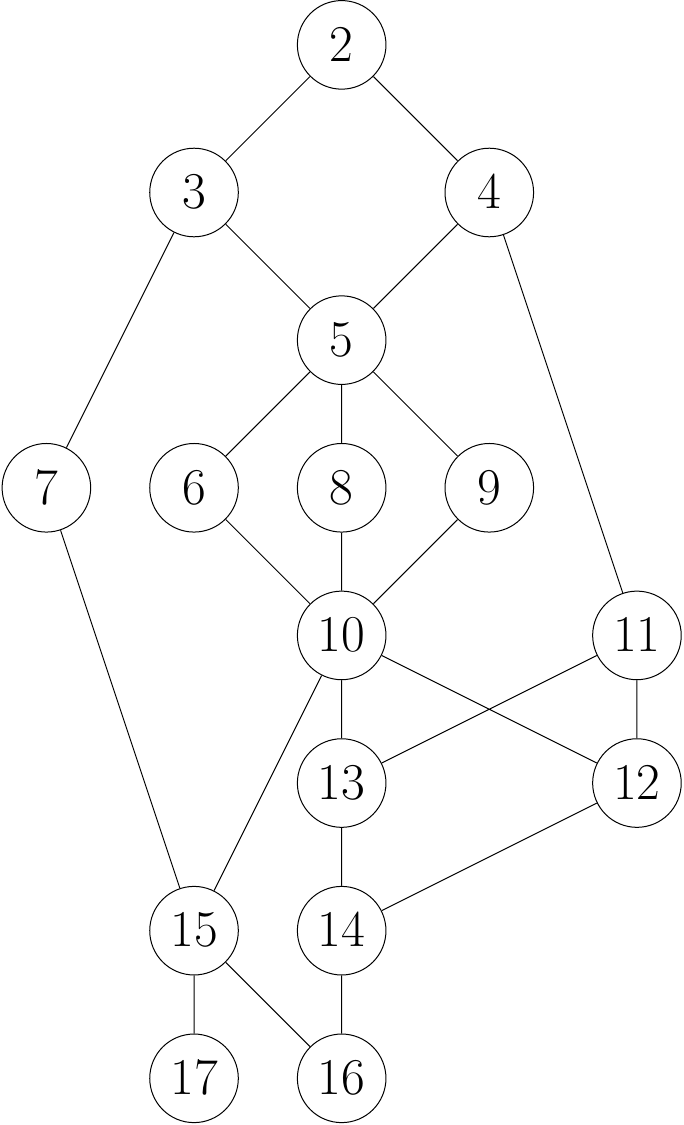}}
\kern3.4truemm

\newpage

L'organigramme de la page précédente donne les liens de dépendance entre les différents chapitres

\bigskip 
\begin{enumerate}\setcounter{enumi}{1}
\item \nameref{chapSli}\\
Anneaux et modules \cohsz. Un peu d'\alg extérieure.
\item \nameref{chapGenerique}\\
Lemme de Dedekind-Mertens et \tho de Kronecker. 
%Algèbre de \dcn \uvlez. Discriminant et résultant. 
Théorie de Galois de base. \nst classique.  
\item \nameref{chap mpf}\\
Catégorie des \mpfsz. Anneaux \zedsz. 
Machineries \lgbes \elrsz. Idéaux de Fitting.
\item \nameref{chap ptf0}\\
\Tho de structure locale. Déterminant. Rang.
\item \nameref{chap AlgStricFi}
\item \nameref{ChapGalois}\\
\nst \gnl (sans clôture \agqz). Théorie de Galois \gnle (sans \algo 
de \fcnz).
\item \nameref{chap mod plats}\\
Algèbres plates et \fptesz.
\item \nameref{chap Anneaux locaux}\\
Anneau décomposable. Anneau \lgbz.
\item \nameref{chap ptf1}
\item \nameref{chapTrdi}\\
Anneaux à pgcd. Treillis de Zariski d'un anneau commutatif. Relations implicatives.
\item \nameref{ChapAdpc}\\
Extensions entières. Dimension $\leq 1$. Factorisation d'\itfsz.
\item \nameref{chapKrulldim}\\
Dimension de Krull. Dimension des morphismes. Dimension valuative.
Dimension des extensions entières et \pollesz.
\item \nameref{chapNbGtrs}\\
Théorèmes de Kronecker, Bass et Forster-Swan. Splitting off de Serre. Dimension de Heitmann.
\item \nameref{chapPlg}
\item \nameref{ChapMPEtendus}\\
Théorèmes de Traverso-Swan-Coquand, Quillen-Suslin, Bass-Lequain-Simis.
\item \nameref{ChapSuslinStab}
%
%\item Résolutions libres finies
%%
%\item Théorie des diviseurs
%
\end{enumerate}

\newpage\thispagestyle{CMcadreseul}
\incrementeexosetprob 

\cleardoublepage
%%%%%%%%%%%%%%%%%%%%%%%%%%%%%%
\pagestyle{CMheadings}
\arabicpagenumbers
\setpagenumber1

%%%%%%%%%%%%%%%%%%%%%%%%%%%%%%%%%%%%%%%%%%%%%%%%%%%%%%%%%%%%%%%%%%%%%%%%%%
%!TEX root =  ACMC-A.tex
%!TEX encoding =  UTF-8 Unicode

%---- CHAP 1 {Exemples}    -------------  
\chapter{Exemples} 
\label{chapMotivation}
%------------------------------- 
%--------------------
\minitoc
\subsection*{Introduction} 
\addcontentsline{toc}{section}{Introduction}
%-----------------------------------------

Dans tout l'ouvrage, sauf mention expresse du
contraire, les anneaux sont commutatifs et \untsz, et un \homo d'anneaux $\varphi:\gA\to\gB$ doit vérifier~$\varphi(1_\gA)=1_\gB$. 

Soit $\gA$ un anneau.
On dit qu'un \Amo $M$ est \ixd{libre de rang fini}{module} 
lorsqu'il est isomorphe à un module $\Ae n$. 
On dit qu'il est \ixd{projectif de type fini}{module} lorsqu'il 
 existe un \Amo $N$ tel que  $M\oplus N$ est libre de rang fini. 
 Il revient au même 
de dire que 
$M$ est isomorphe à l'image d'une \ixx{matrice}{de projection} 
 (une matrice $P$ telle que $P^2=P$). Il s'agit de la matrice de la \prn sur $M$
 \paralm à $N$, définie \prmt comme suit: 
 $$ 
 M\oplus N \lora M\oplus N, \quad  x+y\longmapsto x \qquad \hbox{pour $x\in M$ et $y\in N$}. 
 $$ 
Une \mprn est encore appelée un \ix{projecteur}.

Lorsque l'on a un \iso
 $M\oplus \gA^\ell\simeq \gA^k$, le \mptf $M$ est dit \ixd{stablement 
libre}{module}. 
\index{projectif de type fini!module ---} 

Alors que sur un corps ou sur un anneau principal les  \mptfs 
 sont libres (sur un corps ce sont des \evcs de dimension finie), 
 sur un anneau commutatif
\gnlz, la classification des \mptfs est un \pb à la fois important et difficile.

En théorie des nombres Kronecker et Dedekind ont démontré qu'un \itf non nul dans l'anneau d'entiers d'un corps de nombres est toujours \iv (donc \ptfz), mais qu'il est rarement libre (\cad principal). Il s'agit d'un phénomène fondamental, qui est à l'origine du développement moderne de la théorie des nombres.   

\smallskip 
Dans ce chapitre nous essayons d'expliquer pourquoi la 
notion de \mptf est importante, en donnant des exemples significatifs en \gmt \dilez.

La donnée d'un fibré vectoriel sur une variété compacte lisse $V$ est en effet \eqve à la donnée d'un \mptf sur l'anneau~\hbox{$\gA= \Cin(V)$} des fonctions
lisses sur $V$: à un fibré vectoriel, on associe le \Amo des sections
du fibré, ce \Amo est toujours \ptfz, mais il n'est libre que lorsque le fibré est trivial. 

Le fibré tangent correspond à un module que l'on construit par un procédé purement formel à partir de l'anneau $\gA$. 
Dans le cas où la variété $V$ est une sphère, le module des sections du fibré tangent
est \stlz. Un résultat important concernant la sphère est qu'il n'existe pas de champ de vecteurs lisse partout non nul. Cela équivaut au fait
que le module des sections du fibré tangent n'est pas libre.    

\medskip 
Nous essayons d'être le plus explicite possible, mais dans ce 
chapitre de motivation, nous utilisons librement les raisonnements de 
\clama sans nous soucier d'être totalement rigoureux d'un point de 
vue \cofz.

%\hum{Faire un vrai laius introductif}

\penalty-5000
%--- Section{subsecChampVect}---------  
\section{Fibrés vectoriels sur une \vrt compacte lisse} 
\label{subsecChampVect}
%--------------------

Ici, on donne quelques motivations pour les \mptfs et la localisation en 
expliquant l'exemple  des fibrés vectoriels sur une \vrt lisse 
compacte. Deux cas particuliers importants sont les fibrés tangents et 
cotangents correspondants aux champs de vecteurs et aux formes 
\dilesz~$\Cin$. 

Nous utiliserons le terme \gui{lisse} comme synonyme de 
\gui{de classe~$\Cin$}.

Nous allons voir que le fait que la sphère ne peut pas être 
peignée admet une interprétation purement \agqz.

\ss Dans cette section, on considère une \vrt différentiable 
réelle lisse  $V$  et l'on note $\gA= \Cin(V)$ l'\alg
réelle des fonctions lisses globales sur la \vrtz.

%:      Subsubsection*{Quelques localisées de l'\alg des fonctions}------- 
\Subsubsect{Quelques localisées de l'\alg des fonctions continues}{Quelques localisations} 
%-----------------------------------------

Considérons tout d'abord un \elt $f\in\gA$ ainsi que l'ouvert 

\snic{U=\sotq{x\in V}{f(x)\neq 0}} 

%\sni
et regardons comment on peut 
interpréter l'\alg $\gA[1/f]$: deux \eltsz~$g/f^k$ et $h/f^k$  sont 
égaux dans $\gA[1/f]$ \ssi pour un exposant $\ell$ on a 
$gf^\ell=hf^\ell$ ce qui signifie exactement~$g\frt{U}=h\frt{U}.$

Il s'ensuit que l'on peut interpréter  $\gA[1/f]$ 
comme une sous-\alg de l'\alg des fonctions lisses  sur $U$: 
cette sous-\alg a pour \elts les fonctions
qui peuvent s'écrire sous la forme $(g\frt{U})/(f\frt{U})^k$ (pour un 
certain exposant $k$)
avec $g\in\gA$, ce qui introduit a priori certaines restrictions sur le 
comportement de la fonction au bord de~$U$.

Pour  ne pas avoir à gérer ce \pb délicat, on utilise le lemme 
suivant.
%--- Lemma{lemlocvar1}---------------- 
\begin{lemma} 
\label{lemlocvar1} 
Soit  $U'$ un ouvert contenant le support de $f$. Alors, l'application 
naturelle (par restriction),

\snic{\hbox{de  }\Cin(V)[1/f]=\gA[1/f]\hbox{  vers   }\Cin(U')[1/f\frt{U'}],}

%\sni
est un \isoz. 
\end{lemma}
%--- end-lemma-----------------------------------------
%-----------------begin proof------------------
\begin{proof}
Rappelons que le support de $f$ est l'adhérence de l'ouvert 
$U.$ On a un \homo de restriction $h\mapsto h\frt{U'}$ de $\Cin(V)$ vers $\Cin(U')$
qui induit un \homo $\varphi:\Cin(V)[1/f]\to\Cin(U')[1/f\frt{U'}]$.
Nous voulons montrer que $\varphi$ est un \isoz.
Si $g\in \Cin(U')$, la fonction $gf$, qui est nulle sur $U'\setminus \ov U$, 
peut se prolonger en une fonction lisse à $V$  tout entier, en la 
prenant nulle en dehors de $U'$. Nous la notons encore $gf$. 
Alors, l'\iso réciproque de $\varphi$ est donné par $g\mapsto gf/f$
et $g/f^m\mapsto gf/f^{m+1}$.
\end{proof}
%-----------------end proof------------------

Un \emph{germe
de fonction lisse en un point  $p$}  de la \vrt  
$V$ est donné par
un couple $(U,f)$ où  $U$ est un ouvert contenant $p$ et $f$ est une 
fonction \hbox{lisse $U\rightarrow \RR$}.
Deux couples $(U_1,f_1)$ et $(U_2,f_2)$ définissent le même germe 
s'il existe un ouvert $U\subseteq U_1\cap U_2$ contenant $p$ tel que 
$f_1\frt{U}=f_2\frt{U}$. Les germes de fonctions lisses au point $p$ 
forment une \RRlg que l'on
note~$\gA_{p}.$

On a alors le petit \gui{miracle \agqz} suivant.

%--- Lemma{lemlocvar2}--------------
\begin{lemma} 
\label{lemlocvar2} 
L'\alg $\gA_{p}$ est
naturellement isomorphe au localisé $\gA_{S_p}$, où $S_p$ est la
partie multiplicative des fonctions non nulles au point $p$.
\end{lemma}
%--- end-lemma-----------------------------------------
%-----------------begin proof------------------
\begin{proof}
Tout d'abord, on a une application naturelle $\gA\to\gA_p$ qui à une 
fonction définie sur $V$  associe son germe en $p$. 
Il est immédiat 
que l'image de $S_p$ est contenue dans les inversibles 
de $\gA_p$.
Donc, on a une factorisation de l'application naturelle ci-dessus qui 
fournit un \homoz~$\gA_{S_p}\to\gA_p$.

Ensuite, on définit un \homo  $\gA_p\to\gA_{S_p}$.
Si $(U,f)$ définit le germe~$g$ considérons une fonction
$h\in \gA$ qui est égale à $1$ sur un  ouvert~${U'}$ contenant $p$ 
avec $\overline{U'}\subseteq U$ et qui est nulle en dehors de $U$ 
(dans une carte on pourra prendre pour $U'$ une boule ouverte de centre~$p$).
Alors, chacun des trois couples
$(U,f)$, $(U',f\frt{U'})$ et $(V,fh)$ définit le même germe $g.$ 
Maintenant,~$fh$ définit un \elt de $\gA_{S_p}$. Il reste à 
vérifier que la correspondance que l'on vient d'établir
produit bien un \homo de l'\alg $\gA_{p}$ sur l'\alg  
$\gA_{S_p}$: quelle que soit la manière de représenter le germe sous la forme $(U,f)$, l'\elt $fh/1$ de  $\gA_{S_p}$ ne dépend que du 
germe~$g$.  

Enfin, on vérifie que les deux \homos 
de \RRlgs  que l'on a définis sont bien des \isos inverses l'un de 
l'autre.
\end{proof}
%-----------------end proof------------------

Bref, nous venons d'algébriser la notion de germe de fonction  lisse.
\`A ceci près que le \mo $S_p$ est défini à partir de la 
\vrt $V$, pas seulement à partir de l'\algz~$\gA$. 

Mais si $V$ est compacte, les \mos $S_p$ sont exactement les
complé\-men\-taires des \idemas de $\gA$. En effet, d'une part, que 
$V$ soit ou non compacte, l'ensemble des $f\in \gA$ nulles en $p$ 
constitue toujours un \idema $\fm_p$ de corps résiduel égal à 
$\RR$. D'autre part, si $\fm$ est  un \idema de $\gA$ l'intersection des 
$Z(f)=\sotq{x\in V}{f(x)=0}$ pour les $f\in \fm$ est un 
compact non vide (notez que $Z(f)\cap Z(g)=Z(f^2+g^2)$). Comme l'idéal 
est maximal, ce compact est \ncrt réduit à un point $p$ 
et l'on obtient ensuite~$\fm=\fm_p$. 

%:      Subsubsection*{Fibrés vectoriels}------- 
\Subsubsect{Fibrés vectoriels et \mptfsz}{Fibrés vectoriels} 
%-----------------------------------------

 Rappelons maintenant la notion de \emph{fibré vectoriel} au dessus de 
$V$. 
\\
Un fibré vectoriel est donné par une \vrt lisse $W$, 
une application surjective lisse $\pi:W\rightarrow V$,  et une structure 
d'espace vectoriel de dimension finie sur chaque fibre 
$\pi^{-1}(p)$. En outre, localement, tout ceci doit être 
difféomorphe à la situation simple suivante, dite triviale:  
$$
\pi_1: (U\times{\RR}^m) \rightarrow U,\;(p,v)\mapsto p,
$$ 
avec $m$ qui peut dépendre de $U$ si $V$ n'est pas connexe. Cela 
signifie que la structure d'espace vectoriel  (de dimension finie) sur 
la fibre au dessus de $p$ doit dépendre \gui{convenablement} de~$p$.

Un tel ouvert $U$, qui trivialise le fibré, est appelé un 
\emph{ouvert distingué.}

Une \emph{section} du fibré vectoriel  $\pi:W\rightarrow V$ est par 
\dfn une application~\hbox{$\sigma:V\rightarrow W$} telle que 
$\pi\circ\sigma=\Id_V$. On notera $\Gamma(W)$ l'ensemble des sections  
lisses  de ce fibré. Il est muni d'une structure naturelle de \Amoz. 

Supposons maintenant la \vrt $V$ compacte.  
Comme le fibré est localement trivial il existe un recouvrement fini 
de $V$ par des ouverts distingués $U_i$ et une partition de l'unité 
$(f_i)_{i\in \lrbs}$ \emph{subordonnée à ce recouvrement}: le 
support de $f_i$ est un compact $K_i$ contenu dans~$U_i$.

On remarque d'après le lemme \ref{lemlocvar1} que les \algs 
$\gA[1/f_i]=\Cin(V)[1/f_i]$ et $\Cin(U_i)[1/f_i]$ sont naturellement 
isomorphes.

  Si on localise l'anneau $\gA$ et le module  $M=\Gamma(W)$ en rendant 
$f_i$ inversible, on obtient l'anneau $\gA_i=\gA[1/f_i]$ et le module 
$M_i$. Notons $W_i=\pi^{-1}(U_i)$. Alors, $W_i\to U_i$ est \gui{isomorphe} à 
$\RR^{m_i}\times U_i\to U_i$. Il revient donc au même de se donner 
une section du fibré $W_i$, ou de se donner les $m_i$ 
fonctions~$U_i\to\RR$ qui fabriquent une section du fibré  
$\RR^{m_i}\times U_i\to U_i$.   Autrement dit, le module des sections de $W_i$ est libre 
de rang~$m$. 

Vu qu'un module qui 
devient libre après \lon en un nombre fini d'\eco est \ptf 
(\plgref{plcc.cor.pf.ptf}), on obtient  alors la partie directe 
(point \emph{1}) du \tho suivant.

%--- Theorem{thFVMPTF}---------
\begin{theorem} 
\label{thFVMPTF} 
Soit $V$ une \vrt compacte lisse, on \hbox{note  $\gA= \Cin(V)$}.
\begin{enumerate}
\item Si $W\vvers\pi V$ est un fibré vectoriel sur $V$, le \Amo des sections lisses 
de $W$ est \ptfz.
\item Réciproquement, tout \Amo \ptf est isomorphe au module des sections lisses
d'un fibré vectoriel sur $V$.
\end{enumerate}
\end{theorem}
%--- end-theorem-----------------------------------------

\'Evoquons la partie réciproque du \thoz: si l'on se donne un \Amo
\ptf $M$, on peut construire un fibré vectoriel $W$ au dessus de $V$
dont le module des sections est isomorphe à $M$. On procède comme suit.
On considère une \mprn $F=(f_{ij})\in\Mn(\gA)$ telle que~\hbox{$\Im F\simeq M$}
et l'on pose 
$$
W=\sotq{(x,h)\in V \times \RR^n}{h\in\Im F\frt x},
$$
où $F\frt x$ désigne la matrice $(f_{ij}(x))$.
\Llec pourra montrer alors que~$\Im F$ s'identifie au module des sections $\Gamma(W)$:
à l'\elt $s\in\Im F$ on fait correspondre la section $\wi{s}$ définie par $x\mapsto \wi{s}(x)=(x,s\frt x)$.
Par ailleurs, dans le cas où $F$ est la \mprn standard 
$$
\I_{k,n}=\blocs{.8}{.6}{.8}{.6}{$\I_k$}{$0$}{$0$}{$0_r$}\quad (k+r=n),
$$ 
 il est clair que $W$
est trivial: il est égal à $V\times \left(\RR^k\times \so{0}^{r}\right)$.
Enfin, un \mptf devient libre après \lon en des \eco convenables (\thref{prop Fitt ptf 2}, point~\emph{3}, 
ou \thref{propPTFDec}, forme matricielle plus précise).
En conséquence, le fibré $W$ défini ci-dessus est localement trivial: c'est bien un fibré vectoriel.

%:      Subsubsection*{Vecteurs tangents et \dvns}------- 
\Subsubsec{Vecteurs tangents et \dvns} 
%-----------------------------------------

Un exemple décisif de fibré vectoriel est le fibré tangent, dont 
les \elts sont les couples $(p,v)$, où $p\in V$ et $v$ est un 
vecteur tangent au point $p$. 

Lorsque la \vrt $V$ est une \vrt plongée dans un espace 
$\RR^n$, un vecteur tangent $v$ au point $p$ peut être identifié 
à la \dvn au point $p$ dans la direction de $v$. 

Lorsque la \vrt $V$ n'est pas une \vrt plongée dans un 
espace $\RR^n,$ un vecteur tangent $v$ peut être \emph{défini} 
comme une \emph{\dvn au point $p$,} 
\cad comme une forme $\RR$-\lin $v:\gA\rightarrow \RR$ qui vérifie la 
règle de Leibniz%
%:\omega index
\index{derivation@dérivation!en un point d'une \vrtz}
%--------------------begin equation---------------
\begin{equation}\label{Leibniz1}
v(fg)=f(p)v(g)+g(p)v(f).
\end{equation}
%---------------------end equation--------------

On vérifie par quelques calculs que les vecteurs tangents à $V$ 
forment bien un fibré vectoriel $\mathrm{T}_V$ au dessus de~$V$.

\ss \`A un fibré vectoriel $\pi:W\rightarrow V,$ est associé le 
$\gA$-module $\Gamma (W)$ formé par les sections lisses du fibré.
Dans le cas du fibré tangent, $\Gamma({\mathrm{T}_V})$ n'est rien 
d'autre que le  $\gA$-module des champs de vecteurs (lisses) usuels.

De même qu'un vecteur tangent au point $p$ est identifié à une 
déri\-vation au point $p$, qui peut être définie en termes \agqs 
(\eqrf{Leibniz1}), de même, un champ (lisse) de vecteurs tangents peut 
être identifé à un \elt du~\emph{\hbox{$\gA$-module} des 
\dvns de la \RRlg $\gA$}, défini comme suit. 

Une \dvn d'une
\RRlg $\gB$ dans un \Bmo $M$ est une appli\-cation~$\RR$-\lin $v:\gB\rightarrow M$ qui vérifie la règle de~Leibniz%
%:HHH index
\index{derivation@dérivation!module des ---}
%--------------------begin equation---------------
\begin{equation}\label{Leibniz2}
v(fg)=f\,v(g)+g\,v(f).
\end{equation}
%---------------------end equation--------------
Le \Bmo des \dvns de $\gB$ dans $M$ est noté 
$\Der{\RR}{\gB}{M}$.%
%:HHH index, puis saut de ligne
\index{derivation@dérivation!d'une \alg dans un module}
\\
Une \dvn d'une
\RRlg $\gB$  \gui{tout court} est une \dvn à valeurs dans 
$\gB$. Lorsque le contexte est clair nous noterons $\mathrm{Der}(\gB)$ 
comme une abréviation pour~$\Der{\RR}{\gB}{\gB}$.%
%:HHH index 
\index{derivation@dérivation!d'une \algz}  

%:HHH rajout
Les \dvns au point $p$ sont donc les \elts de $\Der{\RR}{\gA}{\RR_p}$
où $\RR_p=\RR$ muni de la structure de \Amo donnée par 
l'\homo $f\mapsto f(p)$ de $\gA$ dans $\RR$.
Ainsi $\Der{\RR}{\gA}{\RR_p}$ est une version \agq abstraite de
l'espace tangent au point $p$ à la \vrt $V$. 
\index{tangent!espace ---}\index{espace tangent}

\ss Une \vrt lisse est dite \emph{parallélisable} si elle 
possède
un champ (lisse) de bases ($n$ sections lisses du fibré tangent qui en 
tout point donnent une base). Cela revient à dire que le fibré 
tangent est trivial, ou encore que le~\hbox{\Amoz} des sections de ce fibré, 
le module $\mathrm{Der}(\gA)$ des \dvns de~$\gA,$ est libre.

%:      Subsubsection*{Différentielles et fibré cotangent}------- 
\Subsubsec{Différentielles et fibré cotangent} 
%-----------------------------------------

Le fibré dual du fibré tangent, appelé fibré cotangent, admet 
pour sections les formes \diles sur la \vrtz~$V$.

Le $\gA$-module correspondant, appelé \mdiz, 
peut être défini \emph{par \gtrs et relations} de la manière 
suivante.  

\rdb\label{dilesK}
De manière \gnlez, si $(f_i)_{i\in I}$ est une famille 
d'\elts qui engendre \hbox{une \RRlgz} $\gB$, le \emph{\Bmo des  \diles 
(de K\"ahler)}%
\index{module!des diff@des \diles (de K\"ahler)}%
\index{differentielle (de@\dile (de K\"ahler)} 
 de~$\gB$, no\-té~$\Om{\RR}{\gB}$, est engendré par 
les $\rd f_i$ 
(purement formels) soumis aux relations \gui{dérivées} des relations 
qui lient les $f_i$: si  $P\in \RR[z_1,\ldots ,z_n]$ et \hbox{si 
$P(f_{i_1},\ldots ,f_{i_n})=0$}, la relation dérivée est 
$$
\som_{k=1}^n \Dpp P {z_k}(f_{i_1},\ldots , f_{i_n})\rd f_{i_k} =0.
$$
On dispose en outre de l'application canonique $\rd:\gB\to\Om{\RR}{\gB}$, 
définie \hbox{par $\rd f=$ la classe de $f$}
(si $f=\som \alpha_if_i$, avec $\alpha_i\in\RR$, $\rd f=\som \alpha_i\rd f_i$), qui est une \dvnz%
\footnote{Pour plus de précisions sur ce sujet voir les \thosz~\ref{thDerivUniv} et \ref{thDerivUnivPF}.}%
.

On montre alors que, 
pour toute  \RRlg $\gB$, le \Bmo des \dvns de $\gB$ est le dual 
du \Bmo des \diles de K\"ahler.

Dans le cas où le \Bmo des \diles de $\gB$ est \ptf 
(par exemple si $\gB=\gA$), alors 
il est lui-même le dual du 
\Bmo des \dvns de~$\gB$.

%:      Subsubsection*{Cas \agq lisse}------- 
\Subsubsect{Cas des \vrts compactes \agqs lisses}{Cas \agq lisse} 
%-----------------------------------------

Dans le cas d'une \vrt \emph{\agqz} réelle compacte lisse $V$, 
l'\alg $\gA$ des fonctions lisses sur $V$ admet comme sous-\alg celle 
des fonctions \pollesz, notée~$\RR[V]$.

Les modules des champs de vecteurs et des formes \diles 
peuvent être définis comme ci-dessus au niveau de l'\alg $\RR[V]$.

Tout \mptf $M$ sur $\RR[V]$ correspond à un fibré vec\-toriel~\hbox{$W\to V$} 
que l'on qualifie de \emph{fortement \agqz}. 
Les sections lisses de ce fibré vectoriel forment un
\Amo qui est (isomorphe au) le module obtenu à partir de $M$ en 
étendant les scalaires à~$\gA$.

Alors, le fait que la \vrt est parallélisable peut être 
testé au niveau le plus \elrz, celui du module~$M$.  

En effet l'affirmation concernant le cas lisse \gui{le $\gA$-module des 
sections lisses de $W$ est libre} équivaut à l'affirmation 
correspondante de même nature pour le cas \agq \gui{le $\RR[V]$-module 
$M$ est libre}. Esquisse d'une preuve: le \tho d'approximation de 
Weierstrass permet d'approcher une section lisse par une section \polle 
et un champ de bases lisse ($n$ sections lisses du fibré qui en tout 
point donnent une base), par un champ de bases \pollz.

\ms Examinons maintenant le cas des surfaces compactes lisses. 
Une telle surface est parallélisable \ssi elle est orientable et 
possède un champ de vecteurs partout non nul. 
De manière imagée cette deuxième condition se lit: la surface peut 
être peignée.
Les courbes intégrales du champ de vecteurs forment alors une 
\emph{belle famille de courbes}, \cad une famille de courbes localement 
rectifiable.

Donc pour une surface \agq compacte $V$ lisse orientable \propeq
%-----------------begin enum------------------
\begin{enumerate}
\item Il existe un champ de vecteurs partout non nul.
\item Il existe une belle famille de courbes.
\item La \vrt est parallélisable.
\item Le \mdi de K\"ahler de $\RR[V]$ est libre.
\end{enumerate}
%-----------------end enum------------------

Comme expliqué précédemment, la dernière condition relève de 
l'\alg pure (voir aussi la section~\ref{subsecFDVL}).
 
D'où la possibilité d'une preuve \gui{\agqz} du fait que la sphère 
ne peut pas être peignée.

Il semble qu'une telle preuve ne soit pas encore disponible sur le 
marché. 

%:      Subsubsection*{Dérivations d'une \alg \pfz}------- 
\Subsubsect{Module des \diles et module des \dvns 
d'une \alg \pfz}{Dérivations d'une \alg \pfz} 
%--------------------

Soit $\gR$  un anneau commutatif. Pour une \Rlg \pf  
$$\gA=\aqo{\gR[\Xn]}{\lfs}=\gR[\xn],$$ 
les \dfns du module des \dvns et du module des 
\diles s'actualisent comme suit. 

On note $\pi:\gR[\Xn]\to\gA,\;g(\uX)\mapsto g(\ux)$ la projection 
canonique.

On considère la matrice jacobienne du \sys $\lfs$, 
$$
J(\uX)=
\cmatrix{ 
 \Dpp{f_1}{X_1}(\uX)  & \cdots & 
\Dpp{f_1}{X_n}(\uX)      \cr 
 \vdots   &  & 
\vdots       \cr 
 \Dpp{f_s}{X_1}(\uX)  & \cdots & 
\Dpp{f_s}{X_n}(\uX)       
}.
$$
La matrice $J(\ux)$ définit une \Ali $\gA^n\to\gA^s$.
Alors, on a deux \isos naturels $\Om{\gR}{\gA}\simeq 
\Coker\tra{J(\ux)}$
et  $\mathrm{Der}(\gA)\simeq \Ker J(\ux)$. Le premier \iso résulte de 
la \dfn du \mdiz. Le deuxième peut s'expliciter
 comme suit: si $u=(u_1,\ldots ,u_n)\in\Ker J(\ux)$, on lui 
associe \gui{la \dvn partielle dans la direction du vecteur 
tangent~$u$} (en fait c'est plutôt un champ de vecteurs) définie par 
$$
\delta_u: \gA\to\gA,\;\pi(g)\mapsto \som_{i=1}^n u_i 
\Dpp{g}{X_i}(\ux).
$$
Alors, $u\mapsto \delta_u$ est l'\iso en question.

%--- Exercise{exoDerDifPF}----------
\begin{exercise} 
\label{exoDerDifPF} 
Démontrer l'affirmation qui vient d'être faite concernant le 
module des \dvnsz. Confirmer ensuite à partir de cela le fait 
que  $\mathrm{Der}(\gA)$ est le dual de $\Om{\gR}{\gA}$: si $\varphi 
:E\to F$ est une \ali entre modules libres de rang fini on a toujours
$\Ker \varphi \simeq (E\sta/\Im \tra{\varphi})\sta$.
\end{exercise}
%--- end-exercise-----------------------------------------

Nous nous intéressons dans la suite au cas lisse, dans lequel les 
notions purement \agqs coïncident avec les notions analogues en 
géométrie différentielle.

%--- Section{subsecFDVL}--------------  
\section[Formes \diles sur une \vrt affine 
lisse]{Formes \diles à \coes polynomiaux sur une 
\vrt affine lisse} 
\label{subsecFDVL}
%--------------------

\vspace{4pt}
%:  Subsection*{Sphere}----------- 
\subsect{Le module des formes \diles à \coes \\
\polls sur  la sphère}{Le cas de la sphère} \label{subsecFDSphere}
%--------------------

Soit $S=\sotq{(\alpha,\beta,\gamma)\in\RR^3}{
\alpha^2+\beta^2+\gamma^2=1}$. 
L'anneau des fonctions \pols sur $S$ est la \RRlg 
$$
\gA=\aqo{\RR[X,Y,Z]}{X^2+Y^2+Z^2-1}=\RR[x,y,z].
$$ 
Le \Amo des formes \diles à \coes polynomiaux sur $S$ est
%-----------------begin $$----------------
$$ \Om{\RR}{\gA}= (\gA\ {\rd}x\oplus \gA\ {\rd}y\oplus \gA\ 
{\rd}z)/\gen{x{\rd}x+y{\rd}y+z{\rd}z }\simeq \gA^3/\gA v,
$$
%-----------------end $$------------------
où $v$ est le vecteur colonne $\tra{ [\,x\;y\;z\,] }$.\\ 
Ce vecteur est \emph{unimodulaire} (cela signifie que ses \coos
sont des \eco de $\gA$) puisque $[\,x\;y\;z\,]\cdot v=1$. Donc, la matrice 
$$
P=v\cdot[\,x\;y\;z\,]=
\cmatrix{ 
 x^2   &  xy   &   xz   \cr 
  xy  &  y^2   &  yz    \cr
 xz   &  yz   & z^2
}
%---------------------end pmatrix--------------
$$
vérifie $P^2=P$, $P\cdot v=v$, $\Im(P)=\gA v$ de sorte qu'en posant 
$Q=\I_3-P$ on obtient  
$$
\Im(Q)\simeq \gA^3/\Im(P)\simeq \Om{\RR}{\gA}, 
\;\;\mathrm{et}\;\; \Om{\RR}{\gA}\oplus \Im(P)\simeq \Om{\RR}{\gA}
\oplus \gA\simeq \gA^3.
$$
Ceci met en évidence que $\Om{\RR}{\gA}$ est un \Amo \pro de rang 2, 
\stlz.

Les considérations précédentes auraient fonctionné en 
remplaçant $\RR$ par un corps de \cara $\neq 2$ ou même par un 
anneau commutatif $\gR$ où 2 est inversible.

Un \pb intéressant qui se pose est de savoir pour quels anneaux 
$\gR$ exactement le \Amo $\Om{\gR}{\gA}$ est libre.

%:  Subsection{Le cas d'une \vrt lisse}-------
\subsect{Le module des formes \diles à \coes \\
\polls sur une \vrt \agq lisse}{Le cas d'une \vrt \agq
lisse} 
 
\vspace{2pt}   
%:      Subsubsection*{Hypersurface}------- 
\subsubsec{Cas d'une hypersurface lisse} 
%--------------------
Soit $\gR$ un anneau commutatif, et $f(\Xn)\in \RXn=\RuX$. On considère 
la $\gR$-\alg 
$$\gA=\aqo{\gR[\Xn]}{f}=\Rxn=\Rux .$$ 
On dira que \emph{l'hypersurface $S$ définie par $f=0$ est lisse} si,
pour tout corps~$\gK$ \gui{extension de $\gR$}~({\footnote{Dans ce 
chapitre introductif, quand nous utilisons l'expression incantatoire
imagée \emph{corps $\gK$ \gui{extension de $\gR$}}, nous entendons 
simplement que $\gK$ est un corps muni d'une
structure de \Rlgz. Cela revient à dire qu'un
sous-anneau de $\gK$ est isomorphe à un quotient
(intègre) de $\gR$, et que l'\iso est donné.
En conséquence les \coes de $f$ peuvent être \gui{vus}
dans $\gK$ et le discours qui suit l'expression incantatoire
a bien un sens \agq précis. Dans le chapitre \ref{chapGenerique}
nous définirons une extension d'anneaux comme un \homo
\emph{injectif}. Cette \dfn est en conflit direct avec l'expression 
imagée utilisée ici si $\gR$ n'est pas un corps. Ceci explique les guillemets
utilisés dans le chapitre présent.}}) et pour tout point $\uxi 
=(\xin)\in \gK^n$ vérifiant~\hbox{$f(\uxi )=0$} on a une des \coos 
$(\partial f/\partial X_i)(\uxi )$ qui est non nulle.
Par le \nst formel, cela équivaut à l'existence de %\pols
$F$, $B_1$, \ldots, $B_n$ dans~$\RuX$ vérifiant
%-----------------begin $$----------------
$$ 
Ff+B_1 \Dpp{f}{X_1}+\cdots +B_n{\partial f\over 
\partial X_n}=1.
$$
%-----------------end $$------------------
Notons $b_i=B_i(\ux )$ l'image de $B_i$ dans $\gA$ et
$\partial_if=(\partial f/\partial X_i)(\ux )$. On a donc dans~$\gA$
%-----------------begin $$----------------
$$ 
b_1 \,\partial_1f+\cdots +
b_n\, \partial_nf=1.
$$
%-----------------end $$------------------
Le \Amo des formes \diles à \coes polynomiaux sur $S$ est
%-----------------begin $$----------------
$$ \Om{\gR}{\gA}= \aqo{(\gA\ \rd x_1\oplus \cdots \oplus \gA\ \rd 
x_n)}{\rd f}\simeq \gA^n/\gA v,
$$
%-----------------end $$------------------
où $v$ est le vecteur colonne 
$\tra{[\,\partial_1f\;\cdots\;\partial_nf\,]}$. 
Ce vecteur est \umd puisque $[\,b_1\;\cdots\;b_n\,]\cdot v=1$. Alors, 
la matrice 
$$P=v\cdot[\,b_1\;\cdots\;b_n\,]=
\cmatrix{ 
 b_1\partial_1f   &  \ldots   &   b_n\partial_1f   \cr 
  \vdots  &           &  \vdots    \cr
 b_1\partial_nf   &  \ldots   &   b_n\partial_nf
}
%---------------------end pmatrix--------------
$$
vérifie $P^2=P$, $P\cdot v=v$, $\Im(P)=\gA v$ de sorte qu'en posant 
$Q=\I_n-P$ on obtient  
$$\Im(Q)\simeq \gA^n/\Im(P)\simeq \Om{\gR}{\gA} \;\;\mathrm{et}\;\; 
\Om{\gR}{\gA}\oplus \Im(P)\simeq \Om{\gR}{\gA}\oplus \gA\simeq \gA^n.$$
Ceci met en évidence que $\Om{\gR}{\gA}$ est un \Amo \pro de rang 
$n-1$, \stlz.

%:      Subsubsect   intersection complète lisse 
\subsubsect{Cas d'une intersection complète lisse}{Cas d'une intersection complète} 
%--------------------

Nous traitons le cas de deux \eqns qui définissent une 
intersection complète lisse. La \gnn avec un nombre 
quelconque d'\eqns est immédiate.

Soit $\gR$ un anneau commutatif, et $f(\uX)$, $g(\uX)\in \gR[\Xn]$. On 
considère 
la $\gR$-\alg 
$$\gA=\aqo{\gR[\Xn]}{f,g}=\gR[\xn]=\gR[\ux ].$$ 
La matrice jacobienne du système $(f,g)$ est
$$J(\uX)=
\cmatrix{ 
 {\partial f\over \partial X_1}(\uX)  & \cdots & 
{\partial f\over \partial X_n}(\uX)      \cr 
 {\partial g\over \partial X_1}(\uX)  & \cdots &
{\partial g\over \partial X_n}(\uX)      
}.
$$ 
On dira que \emph{la \vrt \agq $S$ définie par $f=g=0$ est 
lisse de codimension $2$} si, pour tout corps $\gK$ \gui{extension de 
$\gR$} et pour tout point~\hbox{$(\uxi) =(\xin)\in \gK^n$} 
vérifiant~\hbox{$f(\uxi )=g(\uxi )=0,$} on a l'un des mineurs~\hbox{$2\times 2$} de la matrice 
jacobienne, $J_{k,\ell}(\uxi ),$ où
$$J_{k,\ell}(\uX)=
\left\vert\matrix{ 
 {\partial f\over \partial X_k}(\uX)  &  
{\partial f\over \partial X_\ell}(\uX)      \cr 
 {\partial g\over \partial X_k}(\uX)  & 
{\partial g\over \partial X_\ell}(\uX)      
}\right\vert,
$$ 
qui est non nul.

\noi Par le \nst formel, cela équivaut à l'existence de \pols
$F,\,G$ et $(B_{k,\ell})_{1\leq k<\ell\leq n}$ dans $\RuX$ qui 
vérifient
%-----------------begin $$----------------
$$ Ff+Gg+\som_{ 1\leq k<\ell\leq n} {B_{k,\ell}(\uX) J_{k,\ell}(\uX) }=1.
$$
%-----------------end $$------------------
Notons $b_{k,\ell}=B_{k,\ell}(\ux )$ l'image de $B_{k,\ell}$ dans $\gA$ 
et
$j_{k,\ell}=J_{k,\ell}(\ux )$. On a donc dans~$\gA$
%-----------------begin $$----------------
$$ \som_{1\leq k<\ell\leq n}{b_{k,\ell}\,j_{k,\ell}}=1.
$$
%-----------------end $$------------------
Le \Amo des formes \diles à \coes polynomiaux sur $S$ est
%-----------------begin $$----------------
$$ 
\Om{\gR}{\gA}= \aqo{(\gA\ \rd x_1\oplus \cdots \oplus \gA\ \rd x_n)}
{\rd f,\rd g}\simeq \gA^n/\Im\tra J,
$$
%-----------------end $$------------------
où $\tra J$ est la transposée de la matrice jacobienne (vue 
dans~$\gA$):
$$\tra J=\tra J(\ux )=
%--------------------begin pmatrix---------------
\cmatrix{ 
\partial_1f    &  \partial_1g   \cr 
    \vdots     &  \vdots        \cr 
\partial_nf    &  \partial_ng    
}
%---------------------end pmatrix--------------
.$$ 
La matrice jacobienne 
$J(\ux )$ définit une \ali surjective 
parce que~\hbox{$\sum_{1\leq k<\ell\leq n}{b_{k,\ell}\,j_{k,\ell}}=1$}, et sa 
transposée définit une \ali injective: 
plus \prmtz, si l'on pose 
$$T_{k,l}(\ux )=
\cmatrix{ 
0 &\cdots&0& \partial_\ell g&0&\cdots &0& -\partial_k g&0&\cdots   &0 \cr 
0 &\cdots&0& -\partial_\ell f&0&\cdots &0& \partial_k f&0&\cdots   &0  
}
%---------------------end pmatrix--------------
$$
et $T=\sum_{1\leq k<\ell\leq n}{b_{k,\ell}T_{k,l}}$, alors $T\cdot \tra 
J=\I_2=J\cdot \tra T$ et la matrice $P=\tra J\cdot T$ vérifie 
$$
P^2=P,\;P\cdot \tra J=\tra J,\;\Im P=\Im\tra J\simeq \gA^2,
$$ 
de sorte qu'en posant $Q=\I_n-P$ on obtient  
$$
\Im Q \simeq \gA^n/\Im P\simeq \Om{\gR}{\gA} \;\;\mathrm{et}\;\; 
\Om{\gR}{\gA}\oplus \Im P\simeq \Om{\gR}{\gA}\oplus \gA^2\simeq 
\gA^n.
$$
Ceci met en évidence que $\Om{\gR}{\gA}$ est un \Amo \pro de rang
 $n-2$, \stlz.

%:      Subsubsection*{Le cas général}------- 
\subsubsect{Le cas \gnl}{Cas \gnlz} 
%-------------------- 

Nous traitons le cas de $m$ \eqns qui définissent une 
\vrt lisse de codimension~$r$.

Soit $\gR$ un anneau commutatif, et $f_i(\uX)\in \RXn$, 
$i=1,\ldots,m$. On considère la $\gR$-\alg 
$$\gA=\aqo{\RXn}{f_1,\ldots ,f_m}=\Rxn=\gR[\ux].$$ 
La matrice jacobienne du système $(f_1,\ldots ,f_m)$ est
$$J(\uX)=
\cmatrix{ 
 {\partial f_1\over \partial X_1}(\uX)  & \cdots & 
{\partial f_1\over \partial X_n}(\uX)      \cr 
 \vdots   &  & 
\vdots       \cr 
 {\partial f_m\over \partial X_1}(\uX)  & \cdots & 
{\partial f_m\over \partial X_n}(\uX)       
}.
$$ 
On dira que \emph{la \vrt \agq $S$ définie par 
$f_1=\cdots=f_m=0$ est lisse de codimension $r$} si la matrice 
jacobienne
vue dans $\gA$ est \gui{de rang $r$}, \cad: 
\Grandcadre{tous les mineurs d'ordre 
$r+1$ sont nuls, \\et les mineurs d'ordre $r$ sont \comz.}

Ceci implique que pour tout corps $\gK$ \gui{extension de $\gR$} et en 
tout point \hbox{$(\uxi)\in \gK^n$} de la \vrt des zéros des $f_i$ dans 
$\gK^n$, l'espace tangent est de codimension~$r$. Si l'anneau $\gA$ est 
réduit, cette condition \gui{géométrique} est d'ailleurs 
suffisante (en \clamaz).  

Notons $J_{k_1,\ldots,k_r}^{i_1,\ldots,i_r}(\uX)$  le mineur $r\times r$ 
extrait sur les lignes $i_1,\ldots ,i_r$ et sur les colonnes 
$k_1,\ldots,k_r$ de $J(\uX)$, %Ce mineur vu dans $\gA$ est noté
et vu dans $\gA$: $j_{k_1,\ldots,k_r}^{i_1,\ldots,i_r}= 
J_{k_1,\ldots,k_r}^{i_1,\ldots,i_r} (\ux )$. 

La condition sur les mineurs $r\times r$ signifie l'existence 
d'\elts $b_{k_1,\ldots,k_r}^{i_1,\ldots,i_r}$  de~$\gA$
tels que
%-----------------begin $$----------------
$$ \sum_{1\leq k_1<\cdots<k_r\leq n, 1\leq i_1<\cdots<i_r\leq m} 
{b_{k_1,\ldots,k_r}^{i_1,\ldots,i_r}\,j_{k_1,\ldots,k_r}^{i_1,\ldots,i_r
}}=1.
$$
%-----------------end $$------------------

Le \Amo des formes \diles à \coes polynomiaux sur $S$ est
%-----------------begin $$----------------
$$ 
\Om{\gR}{\gA}= \aqo{(\gA\ \rd x_1\oplus \cdots \oplus \gA\ \rd x_n)}{\rd 
f_1,\ldots ,\rd f_m}\simeq \gA^n/\Im\tra J,
$$
%-----------------end $$------------------
où $\tra J=\tra J(\ux )$ est la transposée de la matrice jacobienne 
(vue dans~$\gA$).

Nous allons voir que $\Im\tra J$ est l'image d'une matrice de 
projection de \hbox{rang $n-r$}.
Ceci mettra en évidence que $\Om{\gR}{\gA}$ est un \Amo \pro de \hbox{rang 
$n-r$}, (mais a priori il n'est pas \stlz).

Pour cela il suffit de calculer une matrice $H$ de $\gA^{m\times n}$ 
telle que  $\tra J\,H\,\tra J=\tra J$, car alors la matrice $P=\tra J\,H$ 
est la matrice de projection recherchée.

On est donc ramené à résoudre un \sli dont les inconnues sont les 
\coes de la matrice $H$. Or la solution d'un \sli est essentiellement 
une affaire locale, et si on localise en rendant un mineur d'ordre~$r$ 
inversible, la solution n'est pas trop difficile à trouver, en sachant 
que tous les mineurs d'ordre $r+1$ sont nuls.  

Voici par exemple comment cela peut fonctionner.

%--- Exercise{exorecolle}------- 
\begin{exercise} 
\label{exorecolle} 
{\rm  Dans cet exercice, on fait un recollement de la manière la plus 
naïve qui soit. Soit $A \in \Ae{n \times m}$ une matrice de rang $r$,
on cherche à construire une matrice $B\in \Ae{m \times n}$ telle que $ABA=A$.
On notera que si l'on a une solution pour une matrice $A$ on a ipso facto une solution pour toute matrice \eqvez.
\begin{enumerate}
\item [\emph{1.}] Traiter le cas où $A=\I_{r,n,m}=\blocs{1}{1.6}{1}{.8}{$\I_{r}$}{$0$}{$0$}{$0$}$
\item [\emph{2.}] Traiter le cas où $PAQ=\I_{r,n,m}$ avec $P$ et $Q$ \ivsz.
\item [\emph{3.}] Traiter le cas où $A$ possède un mineur d'ordre $r$ \ivz.
\item [\emph{4.}] Traiter le cas \gnlz.
\end{enumerate}
}
\end{exercise}
%--- end -exercise-----------------------------------------
%
\begin{Proof}{\bf Solution. } \emph{1.} On prend $B=\tra A$.

%\sni
 \emph{2.} On prend $B=Q\tra{(PAQ)}P$.

%\sni
\emph{3.} On suppose \spdg que le mineur \iv est dans le coin nord-ouest.
On pose $s=n-r$, $t=m-r$.
On écrit $\delta_1=\det R,$
$$  A=\blocs{1.0}{1.2}{1.0}{.8}{$R$}{$-V$}{$-U$}{$W$}\,,\;
L=\blocs{1.0}{.8}{1.0}{.8}{$\I_r$}{$0$}{$U\wi R$}{$\delta_1\I_s$}\,,
\; 
C=\blocs{1.0}{1.2}{1.0}{1.2}{$\I_r$}{$\wi R V$}{$0$}{$\delta_1\I_t$}.
$$
On obtient
$LA=\blocs{1.0}{1.2}{1.0}{.8}{$R$}{$-V$}{$0$}{$W'$}$ avec $W'=-\delta_1 U\wi R V + W$,
puis  

\snic{LAC=\blocs{1.0}{1.2}{1.0}{.8}{$R$}{$0$}{$0$}{$\delta_1 W'$}\,.}

%\sni
Or les mineurs d'ordre $r+1$  de $A$, donc de  $LAC$, sont nuls, donc $\delta_1^{2} W'=0$. 
 Avec~$M=\blocs{1.0}{.8}{1.0}{1.2}{$\wi R$}{$0$}{$0$}{$0$}\,$,
on a $(LAC)M(LAC)=\blocs{1.0}{1.2}{1.0}{.8}{$\delta_1R$}{$0$}{$0$}{$0$}=\delta_1 LAC$. 
\\[1mm]
Avec $B_1=CML$ cela donne $$LAB_1AC=(LAC)M(LAC)=\delta_1 LAC,$$ donc
en multipliant à gauche par  $\wi L$ et à droite par $\wi C$ 
$$ \delta_1^{s+t}AB_1A=\delta_1^{s+t+1}A.
$$
D'où la solution $B=B_1/\delta_1$ puisqu'on a supposé $\delta_1$ \ivz.

%\sni
\emph{4.} La précalcul qui a été fait avec le mineur $\delta_1$ 
n'a pas demandé qu'il soit \ivz. Il peut être fait avec chacun des mineurs
$\delta_\ell$ d'ordre $r$ de $A$. 
Cela donne autant d'\egts $\delta_\ell^{s+t}AB_\ell A=\delta_\ell^{s+t+1}A$. 
\\
 Une \coli $\som_\ell a_\ell\delta_\ell=1$, élevée à une puissance suffisante, 
donne une \egt $\som_\ell b_\ell\delta_\ell^{s+t+1}=1$, d'où $ABA=A$
pour $B=\som_\ell b_\ell \delta_\ell^{s+t} B_\ell$.
\end{Proof}

\rems ~\\
1) Nous reviendrons sur l'\egt $ABA = A$ en utilisant une formule
magique à la Cramer, cf. le \tho \ref{propIGCram}.

%\sni
2) Dans le dernier exemple, nous nous sommes directement inspirés du \gui{\tho du rang}
qui affirme que si une application lisse $\varphi:U\to \RR^k$ a un rang constant
$r$
en tous les points de $V=\sotq{x\in U}{\varphi(x)=0}$, alors $V$
est une sous-variété lisse de codimension $r$ de l'ouvert $U\subseteq \RR^n$. 
Il s'avère qu'en fait l'analogue que nous avons développé ici ne fonctionne pas toujours correctement. Par exemple avec $\gR=\FF_2$, $
f_1 = X^2 + Y$ et  $f_2 = Y^2$, la \vrtz~$V$ est réduite à un point, l'origine (même si l'on passe à la clôture \agq de $\FF_2$), en lequel la matrice jacobienne est 
de rang  $1$: $\cmatrix{0&1\cr0&0}$. Mais $V$ n'est pas une courbe, c'est un point multiple.
Cela signifie que le \tho du rang pose quelques \pbs en \cara non nulle.
Notre \dfn est donc abusive lorsque l'anneau $\gR$ n'est pas une \QQlgz.      
\eoe

%: entrenous
\entrenous{ on pourrait mettre un exemple de ce genre, à condition de mieux expliquer le rapport entre le fait que les \itfs sont \pros et le fait que tout se passe bien dans les calculs

%:--- Section{subsecThN}--------------- 
%\section{Les pgcds idéaux en théorie des nombres} 
\label{subsecThN}
%--------------------
\rm

\bni{\large\bf 1.3 Les pgcds idéaux en théorie des nombres}

La notion d'idéal a été inventée en théorie des nombres pour 
pallier au défaut d'existence de vrai pgcd pour des nombres entiers 
algébriques.
Par exemple si l'on considère $\alpha=1+2\sqrt{-5}$ 
et $\beta=1-2\sqrt{-5}$ on a
\hum{très insatisfaisant pour le moment,
si on le prend, on le développe}
%--------------------begin equation---------------
\begin{equation}\label{eq1gcdideal}
\alpha \times \beta = 7 \times 3
\end{equation}
%---------------------end equation--------------
mais cette \egt est mal expliquée
dans l'anneau $\ZZ[\alpha]$ engendré par $\alpha$. Par la relation de Bézout, $7$ est étranger à $3$ et l'on devrait 
donc avoir $\alpha=\pgcd(\alpha,7)\times \pgcd(\alpha,3)$ à un 
inversible près.
D'où l'idée d'introduire des pgcds idéaux $\gen{\alpha,7}$,
 $\gen{\alpha,3}$, $\gen{\beta,7}$, $\gen{\beta,3}$ et d'écrire
%--------------------begin equation---------------
\begin{equation}\label{eq2gcdideal}
\gen{\alpha}=\gen{\alpha,7}\times \gen{\alpha,3},\;
\gen{\beta}=\gen{\beta,7}\times \gen{\beta,3}
\end{equation}
%---------------------end equation--------------
et 
%--------------------begin equation---------------
\begin{equation}\label{eq3gcdideal}
\gen{7}=\gen{\alpha,7}\times \gen{\beta,7},\;
\gen{3}=\gen{\alpha,3}\times \gen{\beta,3}
\end{equation}
%---------------------end equation--------------
ce qui explique (3)%(\ref{eq1gcdideal})
. 

\'{E}videmment, pour que ceci soit plus qu'un jeu formel, il faut montrer 
que l'on peut développer un calcul cohérent 
sur la base de ces idées.

Une manière efficace de développer ce calcul, et de donner la preuve qu'il fonctionne bien, est la suivante.

On considère un anneau $\gZ$ contenant les entiers \agqs qui
nous intéres\-sent.
On interprète le pgcd idéal d'\elts de $\gZ$ comme 
l'idéal de $\gZ$ (au sens moderne) engendré par ces \eltsz.

Dans l'exemple ci-dessus,
comme $\gen{7,3}=\gen{1}$, le \tho des restes chinois (cf. \paref{restes chinois}) nous dit que les \egts (4) %(\ref{eq2gcdideal}) 
est automatiquement
vérifiée. De même $\gen{\alpha,\beta}=\gen{1}$ (l'\id contient $2$ et $21$), ce qui justifie les   \egtsz~(5).

Dans la situation \gnle où l'on considère un anneau engendré par des entiers \agqsz, les choses ne se passent pas toujours aussi bien. La solution mise au point par Dedekind:
en présence d'entiers \agqs $(\alpha_i)_{1\leq i\leq k}$ le 
bon anneau à considérer est l'anneau $\gZ$, clôture intégrale de 
$\ZZ[(\alpha_i)_{1\leq i\leq k}]$  dans son corps de fractions.

Le fait que tous les calculs (pgcd, ppcm, produits) fonctionnent bien 
(\cad se comportent comme les calculs analogues dans $\ZZ$) est 
directement relié au fait suivant: tout \itf de $\gZ$ est localement 
libre, (i.e. localement principal).  En effet, comme dans un anneau 
commutatif l'\alg \lin est essentiellement une affaire locale, on est 
ramené à la situation
des entiers usuels: $\gZ$ se comporte localement comme $\ZZ$.
Enfin, pour un \itf d'un anneau intègre être localement libre 
signifie être un \mptfz.

%:section: Exercices   
%\Exercices{ 

\hum{Pas d'exos ni de biblio-comment en fin de chap 1}

%}% fin des exos

%:   ---- Section*{references}-----------    
%\Biblio

\sf
}%fin d'entrenous

\newpage \thispagestyle{CMcadreseul}
\incrementeexosetprob

%:        %%%%%%%%%%%%%%%%%%%%%%%%%%%%%%%%%%%%
%:        %%%%%%%%%%%%%%%%%%%%%%%%%%%%%%%%%%%%
%---- Chapitre 2 {Resolution des slis}------------
\chapter{Principe local-global de base et \slisz}
\label{chapSli}
\minitoc

Dans ce chapitre, comme dans tout l'ouvrage sauf mention expresse du
contraire, les anneaux sont commutatifs et \untsz, et les \homos entre anneaux respectent les $1$. En particulier, un sous-anneau a le même $1$ que l'anneau.

%-------------------- 
\subsection*{Introduction}
\addcontentsline{toc}{section}{Introduction}
%-----------------------------------------
\vskip 0pt\vskip-2pt La théorie de la résolution des \slis  est un thème omniprésent en
\alg commutative (sa forme la plus évoluée est l'algèbre
homologique).
Nous donnons dans ce chapitre un rappel de quelques
résul\-tats classiques sur ce sujet.
Nous y reviendrons souvent.

Nous insistons particulièrement sur le \plg de base, sur la notion de module 
\coh et sur les variations autour de la formule de Cramer.\iplg

%--- Sec{Preliminaires}  secPrelimCh2
\section[Quelques faits concernant les \lonsz]{Quelques faits concernant les quotients et les \lonsz}
\label{secPrelimCh2}
%-----------------------------------------

Commençons par un rappel sur les quotients.
Soit $\fa$  un \id de $\gA$. En cas de besoin, on notera l'application
canonique par $\pi_{\gA,\fa}:\gA\to\gA\sur{\fa}$.

L'anneau quotient $(\gA\sur{\fa},\pi_{\gA,\fa})$
est caractérisé,
\emph{à \iso unique près},
par la \prt universelle suivante.%
\index{anneau!quotient par l'\id $\fa$}\label{NOTAgAfa}

%--- Fact{factUnivQuot}---------------
\begin{fact}
\label{factUnivQuot} \emph{(Propriété  caractéristique du quotient
par l'\id $\fa$)}\\
Un \homo d'anneaux  $\psi : \gA\to\gB$ se factorise par $\pi_{\gA,\fa}$
\ssi
$\fa\subseteq \Ker \psi$, \cade si $\psi (\fa)\subseteq \so{0_\gB}$. Dans ce cas, la factorisation est
unique.
\pun{\gA}{\pi_{\gA,\fa}}{\psi}{\gA\sur{\fa}}{\theta}{\gB}{\homos
nuls sur $\fa$.}

\vspace{-3mm}
\noindent \emph{Explication concernant la figure.} Dans une figure du type ci-dessus,
tout est donné, sauf le morphisme $\theta$ correspondant à la flèche en
traits tiretés. Le point d'exclamation signifie que $\theta$ fait commuter le diagramme et qu'il est \emph{l'unique}
morphisme possédant cette \prtz.
\end{fact}
 
On note $M\sur{\fa M}$ le $\gA\sur{\fa}$-module quotient du \Amo $M$ par
le sous-module engendré par les $ax$ pour $a\in\fa$ et $x\in M$.
Ce module peut aussi être défini par extension des scalaires
à $\gA\sur{\fa}$ du \Amo $M$ (voir \paref{pageChgtBase}, et l'exercice \ref{exoEdsQuot}).

\medskip Passons aux \lonsz, qui sont très analogues aux
quotients (nous reviendrons plus en détail sur cette analogie, en   \paref{secIDEFIL}). Dans cet ouvrage, lorsque l'on parle d'un
\emph{monoïde} contenu dans un anneau, on entend toujours une partie contenant $1$
et stable pour la multiplication.\index{monoide@monoïde!dans un anneau}

\rdb
Pour un anneau $\gA$ nous noterons $\Ati$ le groupe multiplicatif
des \elts inversibles, encore appelé \ix{groupe des unités}.
\label{NOTAAst}

Si $S$ est un \moz, on note~$\gA_S$ ou $S^{-1}\gA$ le localisé de $\gA$
en  $S$.  Tout \elt de $\gA_S$ s'écrit $x/s$ avec $x\in \gA$
et $s\in S$.
\\
Par \dfn on \hbox{a $x_1/s_1=x_2/s_2$} s'il existe $s\in S$ tel que
$ss_2x_1=ss_1x_2$.
En cas de besoin, on notera~\hbox{$j_{\gA,S}:\gA\to\gA_S$} l'application
canonique $x\mapsto x/1$.

Le localisé $(\gA_S,j_{\gA,S})$ est caractérisé,
\emph{à \iso unique près},
par la \prt universelle suivante.
\index{anneau!localisé en $S$}

%:--- Fact{factUnivLoc}---------------
\begin{fact}
\label{factUnivLoc} \emph{(Propriété caractéristique de la
\lon en $S$)}%
\index{localisation!propriété caractéristique}\index{localisation!en un monoïde}
\\
Un \homo d'anneaux  $\psi : \gA\to\gB$ se factorise par $j_{\gA,S}$
\ssi
$\psi (S)\subseteq \gB\eti$, et dans ce cas la factorisation est
unique.
\puN{\gA}{j_{\gA,S}}{\psi}{S^{-1}\gA}{\theta}{\gB}{\homos
qui envoient $S$ dans $\gB\eti$.}{8mm}
\end{fact}
%--- end-fact-----------------------------------------

De même, on note $M_S=S^{-1}M$ le $\gA_S$-module localisé du \Amo $M$
en $S$. Tout \elt de $M_S$ s'écrit $x/s$ avec $x\in M$ et $s\in S$.
Par \dfnz, on a $x_1/s_1=x_2/s_2$ s'il existe $s\in S$ tel que
$ss_2x_1=ss_1x_2$.
Ce module $M_S$ peut aussi être défini par extension des scalaires
à $\gA_S$ du \Amo $M$ (voir \paref{pageChgtBase}, et l'exercice \ref{exoEdsQuot}).
\index{module!localisé en $S$}

\smallskip Un \mo $S$ dans un anneau $\gA$ est dit
\emph{saturé} lorsque l'implication
\index{monoide@monoïde!saturé}
\index{sature@saturé!monoide@monoïde ---}
%-----------------begin $$------------------
$$\forall s, t \in \gA \;\;( st\in S \;\Rightarrow\; s\in S)
$$
%-----------------end $$------------------
est satisfaite. Un monoïde saturé est \egmt appelé un \ixc{filtre}{d'un anneau commutatif}.
Nous appellerons \emph{filtre principal} un filtre engendré par un \eltz:
il est constitué de l'ensemble des diviseurs d'une puissance de cet \eltz.%
\index{principal!filtre --- d'un anneau commutatif}%
\index{filtre!principal d'un anneau commutatif}

\rdb\label{NotaSatmon} On  note $\sat{S}$ le saturé du \mo $S$; il est obtenu en rajoutant tous les \elts qui divisent un \elt de $S$.
Si l'on sature un \mo $S$, on ne change pas la \lonz{\footnote{En fait,
selon la construction précise que l'on choisit pour définir une
\lonz, on aura ou bien \egtz, ou bien \iso canonique, entre les
deux localisés.}}.
Deux \mos $S_1$ et $S_2$ sont dits \emph{\eqvsz} s'ils
ont  le même saturé. Dans ce cas, on 
écrit $\gA_{S_1}=\gA_{S_2}$.%
\index{monoide@monoïde!equiva@équivalents}%
\index{equiva@équivalents!monoïdes ---}

\vspace{-2mm}
\Grandcadre{Nous  gardons la possibilité de localiser en un \mo qui contient
$0$. \\
Le résultat est alors l'anneau \emph{trivial} (rappelons qu'un anneau est trivial\\ s'il est réduit à un seul \eltz, \cad encore si $1=0$).}%
\index{trivial!anneau ---}%
\index{anneau!trivial}

\rdb
Si $S$ est engendré par $s\in \gA$, \cad si
$S=s^\NN\eqdefi\sotq{s^k}{k\in\NN}$, on note~$\gA_s$ ou $\gA[1/s]$ le
localisé
$S^{-1}\gA$, qui est isomorphe à
$\aqo{\gA[T]}{sT-1}$.\label{NOTAA[1/s]}\relax

 \rdb
\medskip
Dans un anneau le \ixc{transporteur}{d'un idéal dans un autre}
d'un \id $\fa$ dans un \id $\fb$
est l'\id 

\snic{(\fb:\fa)_\gA=\sotq{a\in\gA}{a \fa\subseteq \fb}.}

%\sni
Plus \gnltz, si $N$ et $P$  sont deux sous-modules d'un 
\Amoz~$M$,
on définit le \ixc{transporteur}{d'un sous-module dans un autre}
de $N$ dans $P$ comme l'\id\index{ideal@idéal!transporteur}\label{NOTATransp} 

\snic{(P:N)_\gA=\sotq{a\in\gA}{aN\subseteq P}.}

\rdb%\sni
Rappelons aussi que \emph{l'annulateur} d'un \elt $x$ d'un \Amo $M$ est
l'\id $\Ann_\gA(x)=(\gen{0_\gA}:\gen{x})=\sotq{a\in\gA}{ax=0}$.

L'\emph{annulateur du module $M$} est l'\id $\Ann_\gA(M)=(\gen{0_M}:M)_\gA$.
Un module ou un idéal est \ixc{fidèle}{module ---} 
si son annulateur est réduit à $0$.%
\index{annulateur!d'un élément}%
\index{annulateur!d'un module}%
\index{module!fidèle}%
\label{NOTAAnn}%
\index{fidèle!idéal ---}\index{ideal@idéal!fidèle}

\rdb
Les notations suivantes sont \egmt utiles pour un sous-module $N$ de~$M$.

\snic{
(N:\fa)_M=\sotq{x\in M}{\,x\,
\fa\subseteq N}.\label{NOTAAnn2}
}

\snic{
(N:\fa^\infty)_M=\sotq{x\in M}{\exists n,\,x\,
\fa^n\subseteq N}.
}

%\sni
Ce dernier module s'appelle le \emph{saturé} de $N$ par $\fa$.
\index{saturé!module --- de $N$ par $\fa$}

\smallskip  Nous disons qu'un \elt $x$ d'un \Amo $M$ est {\em  \ndzz}
(si $M=\gA$ on dit aussi que $x$ est \ix{non diviseur de zéro}, en un
seul mot)
si la suite

\snic{
0\vers{}\gA\vers{.x}M}

%\sni
est exacte, autrement dit si $\Ann(x)=0$. Si $0_\gA$ est non diviseur de
zéro dans~$\gA$, l'anneau est trivial.%
\index{regulier@régulier!element@\elt ---}

\smallskip
En \gnl pour alléger les notations précédentes
concernant les transporteurs on omet l'indice $\gA$ ou $M$
lorsqu'il est clair d'après le contexte.

\rdb
L'\ixx{anneau}{total des fractions} de $\gA$,
que nous notons $\Frac\gA$, est l'anneau localisé $\gA_S$, où $S$
est le \mo des \elts \ndzs de $\gA$, que nous notons $\Reg \gA $.\label{NOTATotFrac}

%--- Fact{factKerAAsMMs}-----------
\goodbreak
\begin{fact}
\label{factKerAAsMMs}~
%-----------------begin enum------------------
\begin{enumerate}
\item Le noyau de l'\homo naturel $j_{\gA,s}:\gA \to \gA_s=\gA[1/s]$ est
l'\id $(0:s^\infty)_\gA$. Il est réduit à $0$ \ssi $s$ est \ndzz.

\item De même le noyau de l'\homo naturel de $M$ dans $M_s=M[1/s]$ est
le sous-\Amo $(0:s^\infty)_M.$

\item L'\homo naturel $\gA\to\Frac\gA$ est injectif.
\end{enumerate}
%-----------------end enum------------------
\end{fact}
%--- end-fact-----------------------------------------

%--- Fact{fact.bilocal} --------
\begin{fact}
\label{fact.bilocal}\relax
Si $S\subseteq S'$ sont deux \mos de $\gA$ et $M$ un \Amoz, on a des
identifications canoniques $(\gA_S)_{S'}\simeq \gA_{S'}$  et
$(M_S)_{S'}\simeq M_{S'}$.
\end{fact}
%--- end-fact--------------------

%--- Sec{Prin local-global de base}
\section{Principe local-global de base}
\label{secPLGCBasic}
%-----------------------------------------
Nous étudierons le fonctionnement
\gnl du \plg
en \alg commutative
dans le chapitre~\ref{chapPlg}.
Nous le rencontrerons cependant à tous les détours de notre chemin
sous des formes particulières, adaptées à chaque situation.
Une instance essentielle de ce principe
est donnée dans cette section  parce qu'elle est tellement
simple qu'il serait bête de se priver plus longtemps de ce petit
plaisir et de cette machinerie si efficace.

Le \plg affirme que certaines \prts sont vraies
\ssi elles sont vraies après des \lons \gui{en quantité
suffisante}. En \clama on invoque souvent la \lon en tous les \ids
maximaux. C'est beaucoup, et un peu mystérieux, surtout d'un point de
vue \algqz. Nous utiliserons des versions plus simples, et moins
effrayantes, dans lesquelles seulement un nombre fini de \lons
sont mises en {\oe}uvre.

%:--- SUBsection{subsecLocCom}--------
\subsec{Localisations comaximales et principe local-global}

La \dfn qui suit correspond à l'idée intuitive que certains systèmes
de localisés d'un anneau $\gA$ sont \gui{en quantité suffisante} pour
récupérer à travers eux toute l'information contenue dans $\gA$.

%:--- definition --def.moco0----------
%\goodbreak
\begin{definition}\label{def.moco0}\relax ~
\begin{enumerate}
\item   Des \elts $s_1$, $\dots$, $s_n$ sont dits
\ixc{comaximaux}{elements@\elts ---} si
$\gen{1} = \gen{s_1,\dots,s_n}$. Deux \elts \com sont aussi appelés
\emph{étrangers}. \index{etrangers@étrangers!elements@\elts ---}
\item   Des \mos $S_1$, $\dots$, $S_n$ sont dits \ixc{comaximaux}{monoides@\mos ---}
si chaque fois que $s_1\in S_1$, \dots, $s_n\in S_n$, les $s_i$  sont \comz.
\end{enumerate}
\end{definition}

\rdb
\noindent \textbf{Deux exemples fondamentaux.} \label{explfonda}

\noindent 
1) Si $s_1$, $\dots$, $s_n$ sont \comz, les \mos qu'ils engendrent sont \comz.
En effet  considérons des $s_i^{m_i}$ ($m_i\geq 1$) dans les \mos $s_i^{\NN}$, en élevant l'\egt $\sum_{i=1}^na_is_i=1$ à la puissance $1-n+\sum_{i=1}^nm_i$, on obtient, en regroupant convenablement les termes de la somme
obtenue, l'\egt souhaitée $\sum_{i=1}^nb_is_i^{m_i}=1$.

\noindent 
2)
Si $a=a_1\cdots a_n\in\gA$, alors
les \mos $a^\NN$, $1+a_1\gA,$ \ldots, $1+a_n\gA$ sont \comz.
Prenons en effet un \elt $b_i=1-a_ix_i$ dans chaque \mo $1+a_i\gA$ et un \elt $a^{m}$ dans le \mo $a^{\NN}$. 
On doit montrer que l'\id $\fm=\gen{a^{m},b_1,\dots,b_n}$ contient $1$.
Or, modulo  $\fm$  on a $1=a_ix_i$, donc $1=a\prod_ix_i=ax$, et enfin~$1=1^{m}=a^{m}x^{m}=0$.
\eoe

\smallskip  Voici une \carn en \clamaz.
%--- Fact{factMoco}-------------
\begin{factc}
 \label{factMoco}
 {\rm  Soient  des  \mos $S_1$, $\dots$, $S_n$ dans un anneau non trivial
 $\gA$ (i.e., $1\neq_\gA0$).
Les  \mos $S_i$  sont \com \ssi pour tout \idep (resp. pour tout \idemaz) $\fp$ l'un des $S_i$ est contenu
dans~$\gA\setminus\fp$.
 } \end{factc}
%--- end-fact-----------------------------------------
%
\begin{proof}
Soit $\fp$ un \idepz.  Si aucun des $S_i$ n'est
contenu dans $\gA\setminus\fp$, il existe, pour chaque $i$, un $s_i \in S_i
\cap \fp$; alors $s_1, \ldots, s_n$ ne sont pas \comz.
\\
 Inversement, supposons que pour tout  \idema $\fm$ l'un des $S_i$ est contenu
dans~$\gA\setminus\fm$ et soient $s_1\in S_1$, $\ldots$, $s_n\in S_n$ alors l'\id
$\gen{s_1,\ldots,s_n}$ n'est contenu dans aucun \idemaz, donc il contient $1$. 
\end{proof}

\rdb
\medskip
Nous notons $\Ae{m\times p}$ ou $\MM_{m,p}(\gA)$
le \Amo des matrices à $m$ lignes et
$p$ colonnes à \coes dans $\gA$, et $\Mn(\gA)$ désigne $\MM_{n,n}(\gA)$.
Le groupe formé par les matrices \ivs est noté $\GLn(\gA)$,
le sous-groupe des matrices de \deter $1$ est noté $\SLn(\gA)$.
Le sous-ensemble de  $\Mn(\gA)$ formé par les \mprns
(\cad les matrices $F$ telles que $F^2=F$) est noté $\GAn(\gA)$.
L'explication des acronymes est la suivante: $\GL$ pour groupe \linz,
$\SL$ pour groupe \lin spécial et $\GA$ pour grassmannienne affine.
\label{NOTAmatrices}

%:    Plgc {plcc.basic}
\begin{plcc}\label{plcc.basic}\relax   
\emph{(Principe \lgb de base\iplgz,
recollement concret de solutions d'un \sliz)}
\\
Soient $S_1$, $\dots$, $S_n$ des \moco de $\gA$,  $B$
une matrice \hbox{de $\Ae{m\times p}$} et $C$ un vecteur colonne de $
\Ae{m}$.
Alors \propeq
\begin{enumerate}
\item  {Le \sli $BX=C$ admet une solution dans $\gA^{p}$}.
\item  {Pour $ i\in\lrbn$
le \sli $BX=C$ admet une solution dans~$\gA_{S_i}^{p}$}.
\end{enumerate}
Ce principe vaut \egmt pour les \slis à \coes dans un \Amo $M$.
\end{plcc}
%--- end-plgc-----------------------------------------
\begin{proof}
\emph{1} $\Rightarrow$ \emph{2.} Clair.
\\
\emph{2} $\Rightarrow$ \emph{1.}
Pour chaque $i$, on a  $Y_i \in \gA^{p}$ et $s_i \in S_i$ tels que
$B  (Y_i/s_i) =  C$ dans~$\gA_{S_i}^m$. 
Ceci signifie que l'on a un  $t_i \in S_i$
tel que  $t_i\,B  Y_i = s_i t_i\,C$ dans~$\Ae{m}$.
En utilisant %une \coli 
$\sum_i  a_is_i t_i =1$, on a une solution %du \sli 
dans
$\gA$: $X=\sum_i a_i t_i Y_i$.
\end{proof}

\rem Quant au fond, ce \plgc se ramène à la remarque suivante
dans le cas d'un anneau intègre (un anneau est dit \emph{intègre} si tout \elt  est nul ou \ndzz\footnote{La notion est discutée plus en détail  \paref{subsecAnneauxqi}.}). Si les $s_i$ sont \ndzs et si
\index{anneau!intègre}\index{integre@intègre!anneau ---}
$$
\frac{x_1}{s_1}=\frac{x_2}{s_2}
=\cdots=\frac{x_n}{s_n},
$$
la valeur commune de cette fraction, lorsque $\sum_is_iu_i=1$,
est aussi égale à
$$
\frac{x_1u_1+\cdots+x_nu_n}{s_1u_1+\cdots+s_nu_n}
=x_1u_1+\cdots+x_nu_n.
$$
Ce principe pourrait donc s'appeler aussi \gui{l'art de chasser astucieusement les dénominateurs}.
La chose la plus remarquable est sans doute que cela fonctionne
en toute \gntz, même si l'anneau n'est pas intègre.
Merci donc à Claude
Chevalley d'avoir introduit les \lons arbitraires.
Dans certains ouvrages savants, on trouve
la même chose formulée ainsi
(au prix d'une perte d'information sur le caractère très
concret du résultat):
\hbox{le \Amoz} $\bigoplus_\fm\!\gA_{1+\fm}$ (où $\fm$ parcourt tous les \idemas de
$\gA$) est  fidèlement plat.
\eoe

%:    Corollary{corplcc.basic}--------
\begin{corollary}
\label{corplcc.basic}\relax
Soient $S_1$, $\dots$, $S_n$ des \moco de $\gA$,  $x\in\gA$
 et $\fa, \fb$ deux \itfs de $\gA$.
Alors, on a les \eqvcs suivantes.
%
%
%\vspace{-.15em} 
%\begin{enumerate}\itemsep=.1em \mou
\begin{enumerate}
\item  $x=0$ dans $\gA$ \ssi pour $i\in\lrbn,$
$x=0$ dans~$\gA_{S_i}$.
\item  $x$ est \ndz dans $\gA$ \ssi pour $i\in \lrbn,$
$x$ est \ndz dans~$\gA_{S_i}$.
\item
{$\fa=\gen{1}$ dans $\gA$} \ssi pour $i\in \lrbn,$
$\fa=\gen{1}$ dans~$\gA_{S_i}$.
\item
{$\fa\subseteq\fb$ dans $\gA$} \ssi pour $i\in \lrbn,$
$\fa\subseteq\fb$ dans~$\gA_{S_i}$.
\end{enumerate}
\end{corollary}
%--- end-corollary------------------------------------

%
\facile

\rem En fait, comme nous le verrons
dans le \plgrf{plcc.basic.modules}, les \ids n'ont pas besoin d'être \tfz.
\eoe

%--- SUBSUBsection*{Exemples}---
\subsubsection*{Exemples}

Donnons des exemples simples d'application du \plgc de base.
Un cas d'application typique du premier exemple (fait \ref{factExl1Plg})
est celui où le
module $M$ dans l'énoncé  est un \id non nul d'un anneau de
Dedekind.\iplg
Un module $M$ est dit \emph{\lmoz} si, après chaque
\lon en des \moco $S_1,\ldots ,S_n$, il
est engendré par un seul \eltz.%
\index{module!localement monogène}%
\index{localement!module --- monogène}

%:     Fact{factExl1Plg}
\begin{fact}\label{factExl1Plg}
Soit $M=\gen{a,b}=\gen{c,d}$ un module avec deux \sgrsz. 
On suppose que ce module est fidèle
 et \lmoz.
Alors, il existe une matrice $A\in\SL_2(\gA)$  telle que
$\vab a b \,A=\vab c d $.
\end{fact}
\begin{proof}
Si $A=\cmatrix{x&y\cr z&t}$, la matrice cotransposée doit être égale à

\snic{B=\Adj A=\cmatrix{t&-y\cr -z&x}.}

%\snii
En particulier, on cherche à
résoudre le \sli suivant:
$$
\vab a b \,A=\vab c d , \quad \vab c d \,B=\vab a b  \eqno (*)
$$
dont les inconnues sont $x$, $y$, $z$, $t$. Notons que $A\,B=\det(A)\;\I_2$.\\
Inversement si ce \sli est résolu,
on aura $\vab a b =\vab a b \,A\,B$, donc $\big(1-\det(A)\big)\vab a b =0$, et puisque le module
est fidèle, $\det(A)=1$.\\
On a des \moco $S_i$ tels que $M_{S_i}$ est engendré par $g_i/1$ pour un $g_i\in M$.
Pour résoudre le \sli il suffit de le résoudre après
\lon en chacun des $S_i$. 
\\
Dans l'anneau $\gA_{S_i}$,
on a les \egts $a=\alpha_ig_i$,  $b=\beta_ig_i$, $g_i=\mu_ia+\nu_ib$, \hbox{donc
$\big(1-(\alpha_i\mu_i+\beta_i\nu_i)\big)\,g_i=0$}. 
\\
Le module $M_{S_i}=\gen{g_i}$
 reste fidèle, donc $1=\alpha_i\mu_i+\beta_i\nu_i$ dans $\gA_{S_i}$.
Ainsi:

\snic{\vab a b\, E_i=\vab{g_i}0 \;$  avec  $\; E_i=\cmatrix{\mu_i&-\beta_i\cr \nu_i&\alpha_i}$  et  $\det(E_i)=1.}

%\snii
De même on obtiendra $\vab c d \,C_i=\vab{g_i}0$ avec une matrice $C_i$ de
\deterz~1 dans $\gA_{S_i}$. En prenant $A_i=E_i\,\Adj(C_i)$ 
on obtient $\vab a b\, A_i=\vab c d$ et $\det(A_i)=1$ dans $\gA_{S_i}$. 
Ainsi le \sli $(*)$ admet une solution dans~$\gA_{S_i}$.
\end{proof}

\rdb
Notre deuxième exemple est donné par le lemme de Gauss-Joyal:
le point~\emph{\ref{i1lemPrimitf}}
 dans le lemme suivant est prouvé en application
du \plg de base.\iplg
Nous devons d'abord rappeler quelques \dfnsz.
\rdb

Un \elt $a$ d'un anneau est dit \ix{nilpotent}
si $a^n=0$ pour un entier $n\in\NN$.
Les \elts nilpotents dans un anneau $\gA$ forment un \id appelé
\emph{nilradical}, ou encore \ixx{radical}{nilpotent} de l'anneau.
Un anneau est \emph{réduit} si son nilradical est égal à $0$.
Plus \gnlt le nilradical d'un \id $\fa$ de $\gA$ est l'\id formé par les
$x\in\gA$ dont une puissance est dans $\fa$. Nous le noterons $\sqrt{\fa}$ ou
$\DA(\fa)$. Nous notons aussi $\DA(x)$ pour $\DA(\gen{x})$.
Un \id $\fa$ est appelé \emph{un idéal radical} lorsqu'il
est égal à son nilradical.
L'anneau $\gA/\DA(0)=\gA\red$ est \emph{l'anneau réduit associé à $\gA$}.%
\index{anneau!réduit}\label{NOTADA}%
\index{reduit@réduit!anneau ---}%
\index{nilradical!d'un anneau}%
\index{nilradical!d'un idéal}%
\index{radical!idéal ---}%
\index{ideal@idéal!radical}

Pour un  \pol $f$ de $\AXn=\AuX$, on appelle
\emph{contenu} de $f$ et l'on note $\rc_{\gA,\uX}(f)$ ou  $\rc(f)$
l'\id engendré par les \coes de $f$.
Le \pol  $f$ est dit \ixc{primitif}{poly@\pol ---} (en $\uX$) lorsque  $\rc_{\gA,\uX}(f)=\gen{1}$.%
\index{polynome@\pol!primitif}%
\index{contenu!d'un \polz}

Lorsqu'un \pol $f$ de $\AX$ est donné sous la forme
$f(X)=\sum_{k=0}^{n}a_kX^k$, on dit que $n$ est le \ix{degré formel}
de $f$, et $a_n$ est son \emph{\coe formellement dominant}.
Enfin, si $f$ est donné comme nul, son degré formel est $-1$.%
\index{formellement dominant!\coe ---}

%:     Lemma{lemGaussJoyal}
\begin{lemma}\label{lemGaussJoyal}\iJG~
\begin{enumerate}
\item \label{i1lemPrimitf} \emph{(Gauss-Joyal du pauvre)}
Le produit de deux \pols primitifs est un \pol primitif.
\item \label{i2lemPrimitf}  \emph{(Gauss-Joyal)}
Pour $f$, $g\in\AuX$, il existe un entier $p\in\NN$ tel que $\big(\rc(f)\rc(g)\big)^p\subseteq\rc(fg)$.%
\item \label{i3lemPrimitf} \emph{(\'Eléments nilpotents dans $\AuX$)}
Un \elt $f$ de $\AuX$ est nilpotent \ssi tous ses \coes sont nilpotents.
Autrement dit, on a l'\egt  $(\AuX)\red=\Ared[\uX]$.
\item \label{i4lemPrimitf} \emph{(\'Eléments inversibles dans $\AuX$)}
Un \elt $f$ de $\AuX$ est \iv \ssi $f(\uze)$ est \iv et $f-f(\uze)$
est nilpotent. Autrement dit, $\AuX\eti = \Ati+\DA(0)[\uX]$
et en particulier $(\Ared[\uX])\eti=(\Ared)\eti$.
\end{enumerate}
\end{lemma}

\begin{proof}
Notez que l'on a a priori l'inclusion
$\rc(fg)\subseteq \rc(f)\rc(g)$.

\noindent 
\emph{\ref{i1lemPrimitf}. Pour des \pols $f,\,g\in\AX$ en une
 variable.} On a $\rc(f)=\rc(g)=\gen{1}$.
On considère l'anneau quotient
$\gB=\gA\sur{\DA\big(\rc(fg)\big)}$. On doit démontrer que cet anneau est trivial.
Il suffit de le faire après localisation en des \ecoz, par exemple
les \coes de $f$. Autrement dit, on peut supposer qu'un \coe de $f$
est \ivz.
Faisons la preuve sur un exemple suffisamment \gnlz, en supposant que

\snic{f(X)=a+bX+X^2+cX^3+\dots$ et $g(X)=g_0+g_1X+g_2X^2+\dots}

%\sni
Dans l'anneau %réduit 
$\gB$ on a $ag_0=0$, $ag_1+bg_0=0$, $ag_2+bg_1+g_0=0$, \hbox{donc  $bg_0^2=0$}, puis $g_0^3=0$, donc $g_0=0$.
On a alors $g=Xh$ et $c(fg)=c(fh)$, et puisque le degré formel de $h$
est plus petit que celui de $g$, on peut conclure par \recu sur le
degré formel que $g=0$. Comme $\rc(g)=\gen{1}$, l'anneau est trivial. 

\noindent 
\emph{\ref{i2lemPrimitf}. Pour des \pols en une variable.} On considère un \coe $a$ de $f$ et un \coe
$b$ de $g$. Montrons que $ab$ est nilpotent dans %l'anneau
$\gB=\gA\sur{\rc(fg)}$. Ceci revient à démontrer que $\gC=\gB[1/(ab)]$ est trivial.
Or dans  $\gC$, $f$ et $g$ sont primitifs,
donc le point \emph{\ref{i1lemPrimitf}}
implique que $\gC$ est trivial.

\emph{\ref{i2lemPrimitf}} et \emph{\ref{i1lemPrimitf}. Cas général.} 
Le point  \emph{\ref{i2lemPrimitf}} se démontre par \recu sur le nombre
de variables à partir du cas univarié.
En effet, pour $f\in \AX[Y]$ on a l'\egt  

\snic{\rc_{\gA,X,Y}(f)=\gen{\rc_{\gA,X}(h)\mid h\in\rc_{\AX,Y}(f)}.}

%\sni
Ensuite on en déduit
le point \emph{\ref{i1lemPrimitf}.}

\emph{\ref{i3lemPrimitf}.} On note que $f^2=0$ implique
$\rc(f)^p=0$ pour un certain $p$ d'après le point~\emph{\ref{i2lemPrimitf}.}

\emph{\ref{i4lemPrimitf}.}
La condition est suffisante: dans un anneau si $x$ est nilpotent,
$1-x$ est inversible parce que $(1-x)(1+x+\cdots+x^n)=1-x^{n+1}$,
donc si $u$ est \iv et $x$ nilpotent, $u+x$ est \ivz.
Pour voir que la condition est \ncr il suffit de traiter le cas en une variable (on conclut par \recu sur le nombre de variables).
\'Ecrivons $fg=1$ avec $f=f(0)+XF(X)$ \hbox{et  $g=g(0)+XG(X)$}. On a $f(0)g(0)=1$.
Soit $n$ le degré formel de $F$ et $m$ celui de $G$.
On doit montrer que $F$ et $G$ sont nilpotents.
\\
Si $n=-1$ ou $m=-1$, le résultat est clair. On raisonne par \recu sur $n+m$
en supposant $n$, $m\geq0$, $F_n$ et $G_m$ étant les \coes \fmt
dominants. Par \hdr le résultat est obtenu pour les anneaux $(\aqo{\gA}{F_n})[X]$
et $(\aqo{\gA}{G_m})[X]$. Puisque $F_nG_m=0$, on peut conclure par le lemme qui suit.

\noindent NB: on donne des précisions dans l'exercice~\ref{exoNilIndexInversiblePol}.
\end{proof}
%

%:     Lemma{lemNilpotProd}
\begin{lemma}\label{lemNilpotProd}
Soient $a$, $b$, $c\in\gA$. Si $c$ est nilpotent modulo $a$ et modulo~$b$, et si $ab=0$,
alors $c$ est nilpotent.
\end{lemma}
\begin{proof}
On a $c^n=xa$ et $c^m=yb$ donc $c^{n+m}=xyab=0$.
\end{proof}

\rem
On peut formuler ce lemme de manière plus structurelle en considérant pour
deux \ids $\fa,\, \fb$ le morphisme canonique $\gA \to \gA\sur{\fa}
\times \gA\sur{\fb}$ de noyau $\fa\cap\fb$. Si un \elt de $\gA$ est nilpotent
modulo $\fa$ et modulo $\fb$, il l'est modulo $\fa\cap\fb$, donc aussi
modulo $\fa\fb$, car $(\fa\cap\fb)^2 \subseteq \fa\fb$.
On touche ici au  \gui {principe de recouvrement fermé}, voir \paref{prcf1}.
\eoe

%:--- SUBsection{Propriétés de caractère fini}---

\subsec{Propriétés de caractère fini}

Le \plgc de base peut être reformulé comme un \gui{principe de
transfert}.

%:    Principe de Transfert  PrTransfertBasic
\CMnewtheorem{ptfb}{Principe de transfert de base}{\itshape}
\begin{ptfb}\label{lem}
\label{PrTransfertBasic}\index{principe de transfert}\iplg\\
Pour un \sli dans un anneau $\gA$ les \elts $s$ tels que
le \sli ait une solution dans $\gA[1/s]$ forment un \id de $\gA$. 
\end{ptfb}
%--------- fin principe transfert de base ------------------------------- 

%\medskip\noindent
%{\bf Principe de transfert de base.}\label{PrTransfertBasic}\index{principe de transfert}\iplg\\
%{\it Pour un \sli dans un anneau $\gA$ les \elts $s$ tels que
%le \sli ait une solution dans $\gA[1/s]$ forment un \id de $\gA$.
% }
%%--- end-Thm---------------------------
%
%\medskip
Nous proposons tout d'abord
 \alec de faire la \dem que ce principe de transfert est \eqv au \plgc de base.

Nous faisons maintenant une analyse détaillée de ce qui se passe.
L'\eqvc
repose en fait sur la notion suivante.

%:     Definition{defiPropCarFini}
\begin{definition}\label{defiPropCarFini}
Une \prt $\sfP $ concernant les anneaux commutatifs
et les modules est dite \emph{de caractère fini}
si elle est conservée par \lon
(par passage de $\gA$ à $S^{-1}\gA$) et si, lorsqu'elle est
vérifiée avec  $S^{-1}\gA$, alors elle est vérifiée avec
$\gA[1/s]$ pour un certain $s\in S$.\index{propriété de caractère fini}
\end{definition}
%
%:       Fact{fact1PropCarFin}
\begin{fact}
\label{fact1PropCarFin}
Soit  $\sfP $ une \prt \carfz.
Alors, le \plgc pour $\sfP $ est \eqv  au principe de transfert pour~$\sfP$.
Autrement dit  les principes suivants sont \eqvsz.
\begin{enumerate}
\item  Si la \prt $\sfP $ est vraie   après \lon en
 une famille de \mocoz, alors elle est vraie.
\item L'ensemble des \elts $s$ de l'anneau pour lesquels la \prt
$\sfP$ est vraie après localisation en $s$ forme un \idz.
\end{enumerate}
\end{fact}
%--- end-fact-----------------------------------------
%
\begin{proof}
Soit $\gA$ un anneau qui fournit le contexte pour la \prt $\sfP$.
Considé\-rons alors l'ensemble
$I=\sotq{s\in\gA}{\sfP  \mathrm{\;est\; vraie\; pour\;} \gA_s}$.

\noindent
\emph{1. $\Rightarrow$ 2}. Supposons \emph{1}.
Soient $s,t\in I$,
$a,b\in\gA$ et $u=as+bt$.
Les \eltsz~$s$ et $t$ sont \com dans $\gA_u$.
Puisque  $\sfP $ est  stable par \lonz,
$\sfP $ est vraie pour $(\gA_u)_s=(\gA_s)_u$ et $(\gA_u)_t=(\gA_t)_u$.
En appliquant \emph{1}, $\sfP $ est vraie pour $\gA_u$, i.e., $u=as+bt\in I$.

\noindent
\emph{2. $\Rightarrow$ 1}. 
Supposons \emph{2} et soit $(S_i)$ la famille de \moco
considérée. Puisque la \prt
est \carfz, on trouve dans chaque~$S_i$ un \elt $s_i$ tel que $\sfP $ soit vraie
après \lon en $s_i$. Puisque les $S_i$ sont \com les $s_i$ sont des \ecoz.
En appliquant \emph{2}, on obtient $I=\gen{1}$. Et la \lon en $1$ donne la réponse.
\end{proof}

La plupart des \plgcs que nous considérerons dans cet ouvrage s'appliquent
pour des \prts \carfz.
Si \llec le préfère, \il a tout le loisir de remplacer alors le \plgc
par le prin\-cipe de transfert correspondant.

\smallskip  
En \clama
on a pour les \prts \carf  l'\eqvc  de deux notions,
l'une concrète et l'autre abstraite, comme expliqué dans le fait suivant.
On utilisera la version concrète dans les chapitres \ref{chapPlg} et~\ref{ChapSuslinStab}.
%que l'on peut qualifier \gui{de nature quasi globale}.
%En \clama
%on a pour les \prts \carf  l'\eqvc  de deux notions,
%l'une concrète et l'autre abstraite,
%que l'on peut qualifier \gui{de nature quasi globale}.

%--- Fact{fact2PropCarFin}--------------
\begin{factc}
\label{fact2PropCarFin}
Soit  $\sfP $ une \prt \carfz.
Alors, en \clama  \propeq
\begin{enumerate}
\item Il existe des \moco tels que la \prt $\sfP $ soit vraie après \lon
en chacun des \mosz.
\item La \prt $\sfP $ est vraie après \lon en
tout \idemaz.
\end{enumerate}
\end{factc}
%--- end-fact-----------------------------------------
%-----------------begin proof------------------
\begin{proof}
\emph{1. $\Rightarrow$ 2}.
Soit $(S_i)$ la famille de \moco considérée. Puisque la \prt
est \carfz, on trouve dans chaque $S_i$ un \elt $s_i$ tel que $\sfP $ soit vraie
après \lon en $s_i$. Puisque les $S_i$ sont \com les~$s_i$ sont des \ecoz.
Soit $\fm$ un \idemaz. L'un des~$s_i$ n'est pas dans $\fm$. La \lon en $1+\fm$
 est une \lon de la \lon en $s_i$.
Donc $\sfP $ est vraie après \lon en $1+\fm$.
\\
\emph{2. $\Rightarrow$ 1}.  Pour chaque \idema $\fm$ sélectionnons un
$s_\fm\notin\fm$ tel que la \prt $\sfP $ soit vraie après \lon en
$s_\fm$. L'ensemble des $s_\fm$ engendre un \id qui n'est contenu dans aucun
\idemaz, donc c'est l'\idz~$\gen{1}$. Une famille finie de certains
de ces $s_\fm$ est donc un \sys d'\ecoz. La famille des \mos
engendrés par ces \elts convient.
\end{proof}
%-----------------end proof------------------

On a le corolaire \imd suivant.

%--- Fact{factPropCarFin}--------------
\begin{factc}
\label{factPropCarFin}
Soit  $\sfP $ une \prt \carfz.
Alors, le \plgc pour $\sfP $ est \eqv (en \clamaz) au \plga pour~$\sfP $.
Autrement dit  les principes suivants sont \eqvsz.
\begin{enumerate}
\item Si la \prt $\sfP $ est vraie   après \lon
en une famille de \mocoz, alors elle est vraie.
\item Si la \prt $\sfP $ est vraie   après \lon en
tout \idemaz, alors elle est vraie.
\end{enumerate}
\end{factc}
%--- end-fact-----------------------------------------
%%-----------------begin proof------------------
%\begin{proof}
%Résulte \imdt des faits \ref{fact1PropCarFin} et \ref{fact2PropCarFin}.

%\end{proof}
%%-----------------end proof------------------

\rem Donnons une \dem
directe de l'\eqvc en \clama du principe de transfert
et du \plga  pour la \prt $\sfP $ (supposée \carfz).

\noindent
\emph{Transfert $\Rightarrow$ Abstrait}. Supposons la \prt vraie après \lon en
tout \idemaz.  L'idéal donné par le principe de transfert ne peut pas être
strict\footnote{Rappelons qu'un \id est dit \emph{strict} lorsqu'il ne contient pas $1$. Nous ferons usage de cette notion essentiellement dans nos commentaires au sujet des \clamaz.}\index{strict!ideal@idéal}\index{ideal@idéal!strict} car sinon il serait contenu dans un \id maximal
$\fm$, ce qui est en contradiction avec le fait que la \prt est
vraie  après \lon en un~$s\notin\fm$.

\noindent
\emph{Abstrait $\Rightarrow$ Transfert}. Pour chaque \idema $\fm$ sélectionnons un
$s_\fm\notin\fm$ tel que la \prt $\sfP $ soit vraie après \lon en
$s_\fm$. L'ensemble des $s_\fm$ engendre un \id qui n'est contenu dans aucun
\idemaz, donc c'est l'\id $\gen{1}$. On peut conclure par le principe de transfert:
la \prt est vraie après \lon en $1$!
\eoe

\medskip \comm \label{comabstraitconcret}
L'avantage de la \lon en un \idep est que le localisé est un
anneau local, lequel a de très bonnes \prts (voir le chapitre~\ref{chap Anneaux locaux}). 
Le désavantage est que les preuves
qui utilisent un \plga en lieu et place du \plgc correspondant sont
non \covs dans la mesure où le seul accès que l'on a (dans une situation
\gnlez) aux \ideps est donné par le lemme de Zorn. En outre
même le fait \ref{factMoco} est obtenu au moyen d'un raisonnement par
l'absurde qui enlève tout caractère \algq à la \gui{construction}
correspondante.\\
Certains \plgcs n'ont pas de correspondant abstrait, parce que la
\prt concernée n'est pas de caractère fini.
Ce sera le cas des \plgcs \vref{plcc.tf} pour les \mtfs et \vref{plcc.coh} pour les anneaux \cohsz.
\\
Nous ferons un usage systématique efficace et \cof du \plg concret de
base et de ses conséquences.\iplg Souvent, nous nous inspirerons d'une
\dem  d'un \plga en \clamaz.
Dans le chapitre~\ref{chapPlg} nous mettrons au point une
machinerie locale-globale \gnle pour exploiter à fond de manière \cov les
preuves classiques de type local-global.
\eoe

%
%--- SUBSUBsection*{Version abstraite du \plg de base}---
\subsubsection*{Version abstraite du \plg de base}

Vu que la \prt considérée est
\carfz, on obtient en \clama  la version abstraite suivante
pour le \plg de base.\iplg

%:    Plga {plca.basic}
\begin{plca}\label{plca.basic}\relax %~\\
\emph{(Principe \lgb abstrait de base: \rca de
solutions d'un \sliz)}~\\
Soient  $B$ une matrice $\in \Ae{m\times p}$ et
$C$ un vecteur colonne de $\Ae{m}$.
Alors \propeq
\begin{enumerate}
\item  Le \sli $BX=C$ admet une solution dans~$\gA^{p}$.
\item  Pour tout \idema $\fm$
le \sli $BX=C$ admet une solution dans~$(\gA_{1+\fm})^p$.
\end{enumerate}
\end{plca}
%--- end-plga-----------------------------------------
%-----------------begin proof------------------

%:--- SUBsection*{Rendre des \elts \com par force}---
\subsec{Rendre des \elts \com par force}

La \lon en un \elt $s\in\gA$ est une opération fondamentale en
\alg commutative pour rendre $s$ inversible par force. 

Il arrive que l'on ait besoin de rendre \com des \elts $a_1$, \ldots
, $a_n$ d'un anneau $\gA$. \`A cet effet on introduit l'anneau
$$\gB={\AXn}\big/\geN{1-\som_ia_iX_i}=\Axn.$$

%--- Lemma{lemKerCom}--------------
\begin{lemma}
\label{lemKerCom}
Le noyau de l'\homo naturel $\psi :\gA\to\gB$ est \hbox{l'\id
$(0:\fa^\infty)$}, où
$\fa=\gen{a_1,\ldots ,a_n}$. En particulier, l'\homo est injectif
\ssi $\Ann\,\fa=0$.
\end{lemma}
%--- end-lemma-----------------------------------------
%-----------------begin proof------------------
\begin{proof}
Soit $c$  un \elt du noyau, vu l'\iso $\aqo{\gB}{(x_j)_{j\neq i}}\simeq
\gA[1/a_i]$, on a $c=_{\gA[1/a_i]}0$, donc $c\in (0:a_i^\infty)$. On en
déduit $c\in(0:\fa^\infty)$. Inversement si $c\in(0:\fa^\infty)$, il
existe un $r$ tel que $ca_i^r=0$ pour chaque $i$, et \hbox{donc
$\psi(c)=\psi(c)(\sum a_ix_i)^{nr}= 0$}.
\end{proof}
%-----------------end proof------------------

\penalty-2500
%--- Sec{Anneaux et modules cohérents}-----
\section{Anneaux et modules cohérents }
\label{secAnneauxCoherents}
%-----------------------------------------

\subsection*{Une notion fondamentale}
\addcontentsline{toc}{subsection}{Une notion fondamentale}

Un anneau $\gA$ est dit \ixc{cohérent}{anneau ---} \index{anneau!cohérent}
si toute équation \lin 

\snic{LX=0 \;\hbox{ avec }\; L\in \Ae{1{\times}n}\;\hbox{ et }\;X\in \Ae{n{\times}1}}

%\sni
admet pour solutions les
\elts d'un sous-\Amo \tf de $\Ae{n{\times}1}$. Autrement dit:
%---  equation eqAnCoh --------
\begin{equation}\label{eqAnCoh}
\formule{\forall n\in\NN,\, \forall L\in \Ae{1{\times}n},\, \exists
m\in\NN,
\, \exists G\in \Ae{n{\times}m},\,\forall X\in \Ae{n{\times}1} \,,\\[2mm]
\quad \quad LX=0\quad \Longleftrightarrow\quad
 \exists Y\in \Ae{m{\times}1},\; X=GY\;.}
\end{equation}
%---------------------end equation--------------
Cela signifie que l'on maîtrise un peu l'ensemble des solutions
du \sli homogène \hbox{$LX=0$.}

Il est clair qu'un produit fini d'anneaux est \coh \ssi chaque facteur est
\cohz.

Plus \gnlt si $V=(v_1,\ldots ,v_n)\in M^n$, où $M$ est un \Amoz, on
appelle \emph{module des relations} \emph{entre les $v_i$} le sous-\Amo
de $\Ae{n}$ noyau de l'\ali
$$
\breve{V}:\Ae{n}\to M,\quad   (\xn)\mapsto\som_ix_iv_i.
$$
On dira aussi de manière plus précise qu'il s'agit du \emph{module des
relations pour (le vecteur) $V$} , ou encore \emph{du module des \syzys pour (le vecteur) $V$}. Un \elt $(\xn)$ de ce noyau est
appelé une \emph{\rdlz}, ou encore une \emph{\syzyz}  entre les  $v_i$.%
\index{relation de dépendance!linéaire, syzygie}%
\index{syzygie, relation de dépendance linéaire}%
\index{module!des syzygies pour un vecteur}\index{module!des relations pour un vecteur}\index{syzygies!module des --- (pour un vecteur)}

Par abus de langage on parle indifféremment de \emph{la relation}
$\som_ix_iv_i=0$
ou de \emph{la relation} $(\xn)\in \Ae{n}$.
Le \Amo $M$ est  dit \ixc{cohérent}{module ---}
\index{module!cohérent} si pour tout $V\in M^n$,
le module des \syzys est \tfz, autrement dit si l'on~a:
%---  equation eqMoCoh --------
\begin{equation}\label{eqMoCoh}
\formule{
\forall n\in\NN,\,\forall V\in M^{n{\times}1},\, \exists m\in\NN
,\, \exists G\in \Ae{m{\times}n},\,\forall X\in \Ae{1{\times}n}\,, \\[2mm]
\quad \quad
XV=0\quad \Longleftrightarrow\quad
\exists Y\in \Ae{1{\times}m},\; X=YG\;.}
\end{equation}
%---------------------end equation--------------
Un anneau $\gA$ est donc \coh \ssi il est \coh en tant que \Amoz.

Notez que nous avons utilisé dans la formule (\ref{eqMoCoh}) une notation
transposée
par rapport à la formule (\ref{eqAnCoh}). 
C'est pour ne pas avoir la somme $\som_ix_iv_i$
écrite sous forme $\som_iv_ix_i$ avec $v_i\in M$ et $x_i\in\gA$.
Dans la suite, nous ne ferons \gnlt plus cette transposition,
car il nous semble préférable de
garder la forme usuelle $AX=V$ pour un \sliz, même si les matrices $A$
et $V$ sont à \coes dans $M$.

%--- Proposition{propCoh1}------
\begin{proposition}
\label{propCoh1}  Soit $M$  un \Amo \cohz.\\
Tout \sli {\em sans second membre $BX=0$ ($B\in
M^{k{\times}n},\;X\in \Ae{n{\times}1}$)} admet pour solutions les \elts
d'un
sous-\Amo \tf de~$\Ae{n{\times}1}$.
\end{proposition}
%--- end-proposition----------------------------------------
%-----------------begin proof------------------
\begin{proof}
Faisons la \dem par exemple pour $k=2$ (la \dem \gnle fonctionne par
\recu de la même manière). Le principe est le suivant:
%\gui
{on résout la première \eqn et l'on porte la solution \gnle
dans la seconde}. Voyons ceci plus \prmtz. La matrice $B$ est
constituée des lignes $L$ et $L'$. On a une matrice $G$ telle~que

\snic{LX=0\quad \Longleftrightarrow\quad  \exists Y\in \Ae{m{\times}1},\;
X=GY.}

%\sni
Il reste à résoudre $L'GY=0$ qui équivaut à
l'existence d'un vecteur colonne~$Z$ tel que $Y=G'Z$ pour une matrice $G'$
convenable. Donc $BX=0$ \ssi $X$ peut s'écrire sous forme $GG'Z$.
\end{proof}
%-----------------end proof------------------

La proposition précédente est particulièrement
importante pour les \slis sur $\gA$ (\cad lorsque~$M=\gA$).

%H
% -- Commentaire
\medskip \comm La notion 
d'anneau \coh est donc fondamentale du point de vue \algq en algèbre
commutative.
Dans les traités usuels, cette notion est rarement mise en avant
parce que l'on préfère la notion d'anneau \emph{\noez}\footnote{Nous donnons
après ce commentaire une \dfn \cov de cette notion.}.
En \clama tout anneau \noe $\gA$ est \coh parce que tous les sous-modules
de $\Ae{n}$ sont \tfz, et tout \mtf est \coh pour la même raison.
En outre, on a le \tho de Hilbert\ihi qui dit que \emph{si $\gA$ est \noez, toute \Alg \tf est \egmt un anneau \noez,} tandis que la même affirmation est en
défaut si l'on remplace \gui{\noez} par \gui{\cohz}.

D'un point de vue \algq cependant, il semble
impossible de trouver une formulation \cov
satisfaisante de la \noet qui implique la \cohc
(voir l'exercice \ref{exo.quo.coh}).
Et la \cohc est souvent la \prt la plus
importante du point de vue \algqz.
Comme conséquence, la \cohc ne peut pas
être sous-entendue (comme c'est le cas en \clamaz)
lorsque l'on parle d'un anneau ou d'un module \noez.

Le \tho classique disant que sur un anneau \noe
tout \Amo \tf est \noe est souvent avantageusement
remplacé par le \tho \cof 
suivant\footnote{Pour la version non-\noee voir le \thrf{propCoh2},
et pour la version \noeez, voir \cite[corolaire 3.2.8 p.~83]{MRR}.}.

 \emph{Sur un anneau \coh (resp. \noe \cohz)
tout \Amo \pf est  \coh (resp. \noe \cohz)}.

En fait, comme le montre cet exemple, la
\noet est souvent une hypothèse inutilement forte.
\eoe

%H-------------------------
\medskip
La \dfn suivante d'un module \noe est \eqve à la \dfn usuelle en \clamaz,
mais elle est beaucoup mieux adaptée à l'\alg \cov
(seul l'anneau trivial satisfait \cot la \dfn usuelle).
%:     Definition{definoetherien}
\begin{definition}\label{definoetherien}\label{noetherien} \emph{(Noethérianité à la Richman-Seidenberg, \cite{ric74,sei74b})}\\
Un \Amo est dit \emph{\noez} s'il vérifie la \emph{condition de chaîne
ascendante} suivante: toute suite croissante de sous-\mtfs possède deux termes consécutifs égaux. Un anneau $\gA$ est dit \emph{\noez} s'il est \noe en tant que \Amoz.
\index{module!noeth@\noez}
\index{anneau!noeth@\noez}
\index{noetherien@\noez!module ---}
\index{noetherien@\noez!anneau ---}
\end{definition}

\smallskip Voici un corolaire de la proposition \ref{propCoh1}.

%--- Corollary{corpropCoh1}-------
\begin{corollary}\emph{(Transporteurs et \cohcz)}
\label{corpropCoh1} ~\\
Soit $\gA$ un anneau \cohz. Alors, le transporteur
d'un \itf  dans un autre est un \itfz.
Plus \gnltz, si $N$ et $P$  sont deux sous-\mtfs d'un \Amo  \cohz, alors  $(P:N)$ est un \itfz.
\end{corollary}
%--- end-corollary------------------------------------

%--- Theorem{propCoh4}----------
\begin{theorem}\label{propCoh4}
Un \Amo $M$ est \coh \ssi sont véri\-fiées les deux conditions
suivantes.
%-----------------begin enum------------------
\begin{enumerate}
\item  L'intersection de deux sous-\mtfs arbitraires est un \mtfz.
\item  L'annulateur d'un \elt arbitraire est un \itfz.
\end{enumerate}
%-----------------end enum------------------
\end{theorem}
%--- end-theorem-----------------------------------------
\begin{proof}
\emph{La première condition est \ncrz.}
Soient $g_1$, \ldots, $g_n$ des \gtrs du premier sous-module
et $g_{n+1}$, \ldots, $g_{m}$ des \gtrs du second.
Se donner un \elt de l'intersection revient à se donner une relation
$\som_{i=1}^{m}\alpha_ig_i=0$ entre les $g_i$:  à une telle relation
 $\alpha =(\alpha_1,\ldots ,\alpha_m)\in\Ae{m}$, correspond l'\elt
$\varphi(\alpha ) = \alpha_1g_1+\cdots +\alpha_ng_n =-(\alpha
_{n+1}g_{n+1}+\cdots+\alpha_mg_{m})$ dans l'intersection. Donc si $S$ est
un \sys \gtr pour les  relations entre les $g_i$,  $\varphi(S)$
engendre l'intersection des deux
sous-modules.

\noindent \emph{La deuxième condition est \ncrz} par \dfnz.

\noindent \emph{Les deux conditions mises ensemble sont suffisantes.} 
Nous donnons
l'idée essentielle de la \dem et laissons les détails \alecz.
Nous considérons le module des relations pour un $L\in M^n$.
On raisonne par \recu sur  $n$. Pour $n=1$ la deuxième condition
s'applique et donne un \sys \gtr pour les relations liant l'unique \elt de
$L$.
\\
 Supposons que le module des relations pour tout $L\in M^{n}$
soit \tf et considérons un  $L'\in M^{n+1}$. Soit  un entier $k\in \lrbn$,
on écrit
$L'=L_1\bullet L_2$ où $L_1=(a_1,\ldots ,a_{k})$ et
$L_2=(a_{k+1},\ldots ,a_{n+1})$. Posons
$M_1=\gen{a_1,\ldots ,a_{k}}$ \hbox{et $M_2=\gen{a_{k+1},\ldots ,a_{n+1}}$}.
Se donner une relation $\sum_{i=1}^{n+1}\alpha_ia_i=0$ revient à se
donner un \elt de l'intersection  $M_1\cap M_2$ (comme ci-dessus). On
obtiendra donc un \sgr pour les relations entre les $a_i$ en prenant la
réunion des trois
\syss de relations suivants:  celui des relations entre les \elts de $L_1$,
celui des relations entre les \elts de $L_2$, et celui qui provient du \sys
\gtr de l'intersection $M_1\cap M_2$.
\end{proof}
%-----------------end proof------------------

En particulier, \emph{un anneau est \coh \ssi d'une part l'intersection de
deux \itfs est toujours un \itfz, et d'autre part l'annulateur d'un \elt
est toujours un \itfz}.

\medskip \exls
Si $\gK$ est un corps discret, toute \alg \pf sur $\gK$ 
est un anneau \coh (\thref{thpolcohfd}).
Il est clair aussi que tout anneau de Bézout intègre
(cf. \paref{secBézout}) est un anneau \cohz.
\eoe

\subsection*{Caractère local de la \cohcz}
\addcontentsline{toc}{subsection}{Caractère local de la \cohcz}

%:     Prc loc glob conc {plcc.coh}  -
La \cohc est une notion locale, au sens suivant.

\begin{plcc}
\label{plcc.coh}\relax
{\em  (Modules \cohsz)}\\
On considère un anneau  $\gA$, $S_1$, $\ldots$, $ S_n$  des \moco et $M$ 
un~\hbox{\Amoz}.
\begin{enumerate}
\item Le module $M$ est \coh \ssi chacun des $M_{S_i}$ est \cohz.
\item L'anneau $\gA$  est \coh \ssi chacun des $\gA_{S_i}$ est
\cohz.
\end{enumerate} 
\end{plcc}
%--- end-plcc-----------------------------------------
%-----------------begin proof------------------
\begin{proof}
Soit  $a=(a_1,\ldots ,a_m)\in M^m$, et $N\subseteq\Ae{m}$ le module des
relations pour~$a$.
Nous constatons que pour n'importe quel \mo $S$, $N_S$ est le module des
relations pour $a$ dans $M_S$. Ceci nous ramène à démontrer
le \plgc qui suit.
\end{proof}
%-----------------end proof------------------

%:     Prc loc glob conc {plcc.tf} -
\begin{plcc}
\label{plcc.tf}\relax
{\em  (Modules \tfz)}\\
Soient $S_1$, $\ldots$, $ S_n$  des \moco de $\gA$ et $M$ un \Amoz.
Alors,~$M$ est \tf \ssi chacun des $M_{S_i}$ est \tfz.
\end{plcc}
%--- end-plcc-----------------------------------------
%-----------------begin proof------------------
\begin{proof}
Supposons que  $M_{S_i}$ soit un $\gA_{S_i}$-\mtf pour chaque $i$.
Montrons que $M$ est \tfz.
Soient $g_{i,1}$, \ldots, $g_{i,q_i}$ des \elts de $M$ qui engendrent
$M_{S_i}$. Soit $x\in M$ arbitraire. Pour chaque $i$ on a un $s_i\in S_i$
et des $a_{i,j}\in \gA$ %convenables
tels que:
%----begin $$----------
$$ s_ix=a_{i,1}g_{i,1}+\cdots+a_{i,q_i}g_{i,q_i} \quad {\rm  dans} \quad M.
$$
%----end $$----------
En écrivant $\sum_{i=1}^{n} b_i s_i =1$, on voit que $x$ est \coli des
$g_{i,j}$.
\end{proof}
%-----------------end proof------------------

\rdb\label{remplcc.tf}\relax
\rem
Considérons le sous-$\ZZ$-module $M$ de $\QQ$ engendré par les
\eltsz~$1/p$ où $p$ parcourt l'ensemble des nombres premiers.
On vérifie facilement que $M$ n'est pas \tf mais qu'il
devient \tf après \lon en n'importe quel \idepz.
Cela signifie que le \plgc \ref{plcc.tf} n'admet pas de version
\gui{abstraite} correspondante, dans laquelle la \lon en des
\moco serait remplacée par la \lon en tous les \idepsz.
En fait la \prt $\,\sfP\,$ pour un module d'être \tf n'est pas
une \prt \carfz, comme on peut le voir avec le module
$M$ ci-dessus et les \mos $\ZZ\setminus\so0$ ou $1+p\ZZ$.
La \prt vérifie par ailleurs le principe de transfert,
mais en l'occurrence, cela n'est d'aucune utilité.
\eoe

%:  ---SUBsection Au sujet du test d'\egt
\subsec{Au sujet du test d'\egt et du test d'appartenance}
\label{subsecTestDEgalite}
Nous introduisons maintenant quelques notions \covs
relatives au test d'\egt et au test d'appartenance.

\medskip Un ensemble $E$ est bien défini lorsque l'on a indiqué comment
construire ses \elts et lorsque l'on a construit
une relation d'\eqvc  qui définit l'\egt de deux \elts dans l'ensemble.
On note $x=y$ l'\egt dans $E$, ou $x=_Ey$ si \ncrz.
L'ensemble $E$ est appelé
\ixc{discret}{ensemble ---} lorsque l'axiome suivant est vérifié
\index{ensemble!discret}
$$
\forall x,y\in E \qquad x=y \;\; \hbox{ou} \;\;  \lnot (x=y).
$$

Classiquement, tous les ensembles sont discrets,
car le \gui{ou} présent dans la \dfn est compris
de manière \gui{abstraite}.
Constructivement, le \gui{ou} présent dans la \dfn est compris selon la
signification du langage usuel: une des deux alternatives au moins
doit avoir lieu de manière certaine.
Il s'agit donc d'un \gui{ou} de nature algorithmique.
En bref un ensemble est discret si l'on a un test
pour l'\egt de deux \elts arbitraires de cet ensemble.

Si l'on veut être plus précis et expliquer en détail ce qu'est
un test d'\egt dans l'ensemble $E$, on dira qu'il s'agit d'une construction
qui,
à partir de deux \elts de $E$ donnés en tant que tels, fournit une
réponse \gui{oui} ou \gui{non} à la question posée
(ces \elts sont-ils égaux?). Mais on ne pourra guère aller
plus loin.
En \coma les notions  de nombre entier et de construction sont des
concepts de base. Elles peuvent être expliquées et commentées,
mais pas à proprement parler \gui{définies}. La signification \cov du
\gui{ou} et
celle du \gui{il existe} sont ainsi directement dépendantes de la notion
de construction\footnote{En \clama on peut vouloir définir
la notion de construction à partir de la notion de \gui{programme
correct}.
Mais ce que l'on définit ainsi est plutôt la notion de \gui{construction
mécanisable}. Et surtout dans la notion de \gui{programme correct}, il y
a le fait que le programme doit s'arrêter après un nombre fini
d'étapes.
Ceci cache un \gui{il existe}, qui en \coma renvoie de manière
irréductible à la notion de construction. Voir à ce sujet la section \ref{AnnexeCalculsMec} de l'Annexe.}, que l'on ne tente pas de
définir.

\smallskip Un \ixx{corps}{discret} (en un seul mot)
est un anneau où est vérifié l'axio\-me suivant:
\index{discret!corps ---}\index{corps}
%---  equation eqDefCodi --------
\begin{equation}\label{eqDefCodi}
\forall x\in \gA \qquad x=0 \;\; {\rm  ou} \;\;  x\in\Ati
\end{equation}
%---------------------end equation--------------
L'anneau trivial est un corps discret.

\medskip
\rem
La méthode chinoise du pivot (souvent appelée méthode du pivot de
Gauss)
fonctionne de façon \algq avec les corps discrets.
Ceci signifie que l'\alg \lin de base est explicite sur les corps
discrets.\eoe

\medskip
Notons qu'un corps discret $\gA$
est un ensemble discret \ssi
le test \gui{$1=_\gA0$?} est explicite\footnote{La notion \gnle de corps
en \coma sera définie \paref{corpsdeHeyting}. Nous verrons que si un corps est un ensemble discret c'est un corps discret.}.
Il arrive cependant que l'on sache qu'un anneau construit au cours d'un \algo
est un corps discret sans savoir s'il est trivial ou non.

Si $\gA$ est un corps discret non trivial, l'affirmation \gui{$M$ est un
\evc libre de dimension finie} est plus précise que l'affirmation \gui{$M$ est
un \evc \tfz}, car dans le dernier cas, savoir extraire une base du
\sys \gtr revient à disposer d'un test d'indépendance \lin
dans~$M$.

\rdb
\medskip
Une partie $P$ d'un ensemble $E$ est dite \ixe{détachable}{detachable}
lorsque la \prt suivante est vérifiée:\label{detachable}

\snic{\forall x\in E \qquad x\in P \;\; {\rm  ou} \;\;
   \lnot(x\in P).}

\smallskip
Il revient au même de se donner une partie détachable de $E$ ou sa fonction
\cara $\chi_P:E\to\so{0,1}$.

En \coma on considère que si deux ensembles
$E$ et $F$ sont correctement définis, il en va de même pour l'\ixx{ensemble}
{des fonctions de $E$ vers $F$}, que l'on note $F^E$. En conséquence
l'\ixx{ensemble}{des parties détachables} d'un ensemble $E$ est lui-même correctement défini car il s'identifie à l'ensemble $\so{0,1}^E$
des fonctions \caras de source $E$.

%:  ---SUBsection Anneaux et modules coh ftt dis
\subsec{Anneaux et modules cohérents  fortement discrets}

Un anneau (resp. un module) est dit \ixc{fortement discret}{anneau,
module}
\index{module!fortement discret}
\index{anneau!fortement discret}
lorsque les \itfs (resp. les sous-\mtfsz) sont détachables, \cade si les
quotients par les \itfs (resp. par les sous-\mtfsz) sont discrets.

Cela revient à dire que l'on a un test pour décider si une équation
\linz~\hbox{$LX=c$} a ou non une solution, et en calculer une en cas de réponse
positive.

Un résultat essentiel pour l'\alg \cov et le \calf affirme que
$\ZZ[\Xn]$ est un anneau \coh \fdiz. 
\\
Plus \gnltz, on a la version \cov
suivante du \tho de Hilbert\ihi (voir \cite{MRR,Lou}).

\emph{Si $\gA$ est un anneau \noe \coh  \fdiz,
 il en va de même pour toute \Alg \pfz.}

La proposition suivante se démontre comme la proposition
\ref{propCoh1}.

%--- Proposition{propCohfd1}------
\begin{proposition}
\label{propCohfd1}
Sur un module \coh fortement discret $M$, tout \sli  $BX=C$ ($B\in
M^{k{\times}n},\;C\in M^{k{\times}1},\;X\in \Ae{n{\times}1}$) peut
être testé.
En cas de réponse positive, une solution particulière $X_0$ peut
être calculée.
En outre les solutions $X$ sont tous les \elts de $X_0+N$ où $N$ est un
sous-\Amo \tf de $\Ae{n{\times}1}$.
\end{proposition}
%--- end-proposition----------------------------------------

%--- Sec{sec sfio}-------------
\section{Systèmes fondamentaux d'\idms \ortsz} \label{sec sfio}\relax
%-------------------
Un \elt $e$ d'un anneau est dit \ix{idempotent} si $e^2=e$.
 Dans ce cas, $1-e$ est aussi un \idmz, appelé
l'\emph{idempotent complémentaire de $e$}, ou encore le \emph{complément de $e$}.
\index{idempotent!complémentaire}\index{complement@complément!d'un \idmz}
Pour deux \idms $e_1$ et $e_2$, on a
$$\gen{e_1}\cap\gen{e_2}=\gen{e_1e_2}, \quad \gen{e_1}+\gen{e_2}=\gen{e_1,e_2}=\gen{e_1+e_2-e_1e_2},
$$
avec $e_1e_2$ et $e_1+e_2-e_1e_2$ \idmsz.
Deux \idms $e_1$ et $e_2$ sont dits \ixc{orthogonaux}{idempotents ---}
lorsque $e_1e_2=0$. On a alors $\gen{e_1}+\gen{e_2}=\gen{e_1+e_2}$.

Un anneau est dit \ixc{connexe}{anneau ---} si tout \idm est égal à $0$ ou
$1$. \index{anneau!connexe}

Dans la suite nous utilisons implicitement le fait évident suivant:
%\label{lem ide}\relax
pour un idempotent $e$ et un \elt $x$, $e$
divise  $x$ \ssi $x=ex$.

La présence d'un \idm $\neq 0$, $1$ signifie que l'anneau
$\gA$ est isomorphe à un produit de deux anneaux $\gA_1$ et
$\gA_2$, et que tout calcul dans $\gA$ peut être scindé
en deux calculs \gui{plus simples} dans $\gA_1$ et $\gA_2$.
On décrit  cette situation comme suit.  

%:     Fact{lemCompAnnComm}
\begin{fact}\label{lemCompAnnComm}
Pour tout \iso $\lambda:\gA\to\gA_1\times \gA_2$, il existe un unique
\elt $e\in\gA$
satisfaisant les \prts suivantes.
\begin{enumerate}
\item L'\elt $e$ est \idm  (on note son complément $f=1-e$).
\item L'\homo % composé
$\gA\to\gA_1$
identifie $\gA_1$ avec $\gA\sur{\gen{e}}$  et avec $\gA[1/f]$.
\item L'\homo % composé
$\gA\to\gA_2$
identifie $\gA_2$ avec $\gA\sur{\gen{f}}$ et avec $\gA[1/e]$.
\end{enumerate}
Réciproquement, si $e$ est un \idm et $f$ son complément, l'\homo
canonique $\gA\to\aqo{\gA}{e} \times \aqo{\gA}{f}$ est un \isoz.
\end{fact}
\begin{proof}
L'\elt $e$ est défini par $\lambda(e)=(0,1)$.
\end{proof}

On peut apporter quelques précisions souvent utiles.
%: --- Fact{fact.loc.idm}---------
\begin{fact}\label{fact.loc.idm}\label{fact.loc.idm1}Soit $e$ un \idm de $\gA$, $f=1-e$ et $M$ un \Amoz.
%-----------------begin enum------------------
\begin{enumerate}
\item Les \mos $e^\NN=\so{1,e}$
et $1+f\gA$ ont le même saturé.

\item En tant que \Amoz, $\gA$ est somme directe de $\gen{e}=e\gA$ et
$\gen{f}=f\gA$.
L'idéal $e\gA$ est un anneau si l'on prend $e$ comme \elt neutre pour la
multiplication. On a alors trois anneaux isomorphes
$$
\gA[1/e]=(1+f\gA)^{-1}\gA\;\simeq\;\aqo{\gA}{f}\;\simeq\; e\gA.
$$
Ces \isos proviennent des trois applications canoniques
\[\arraycolsep2pt
\begin{array}{lllllllllllllllll} 
\gA&\rightarrow&\gA[1/e]&\ {\string:}\ \ & x&\mapsto& x/1,  \\ 
\gA&\rightarrow&\aqo{\gA}{f}&\ {\string:}\ & x&\mapsto& x\,\mod\gen{f},  \\ 
\gA&\rightarrow& e\gA&\ {\string:}\ & x&\mapsto&   e\,x, 
 \end{array}
\]  
qui sont surjectives et ont même noyau.

%\item On a l'\iso d'anneaux fondamental
%$\gA\simeq \aqo{\gA}{e}\times \aqo{\gA}{f}$.

\item On a trois \Amos isomorphes $
M[1/e]\simeq M/fM\simeq eM.$
Ces \isos proviennent des trois applications canoniques  
\[ \arraycolsep2pt
\begin{array}{lllllllllllllllll} 
M &\rightarrow&M [1/e]  &\ {\string:}\ \ & x&\mapsto& x/1,  \\ 
M &\rightarrow&{M }/{fM}  &\ {\string:}\ & x&\mapsto& x\,\mod\gen{f},  \\ 
M &\rightarrow& eM   &\ {\string:}\ & x&\mapsto&   e\,x, 
 \end{array}
\]  
 qui sont surjectives et ont même noyau.
\end{enumerate}
%-----------------end enum------------------
\end{fact}
%--- end-fact-----------------------------------------

Par ailleurs, il faut prendre garde que l'\id $e\gA$, qui est un anneau
avec
$e$ pour \elt neutre, n'est pas un sous-anneau de $\gA$ (sauf si $e=1$).

\smallskip
Dans un anneau $\gA$ un \emph{\sfioz}
%(que nous abrègerons en \ix{sfio})
est une liste
$(e_1,\ldots ,e_n)$ d'\elts de $\gA$ qui satisfait les \egts suivantes:%
\index{systeme fon@\sfioz}
$$\preskip.4em \postskip.4em\ndsp 
e_ie_j=0 \;\hbox{  pour   }\; i\not= j, \quad
 \hbox{ et }\quad \som_{i=1}^n\, e_i=1  . 
$$
Ceci implique que les $e_i$ sont \idmsz. Nous ne réclamons pas qu'ils
soient tous non nuls{\footnote{C'est beaucoup plus confortable pour
obtenir des énoncés uniformes. En outre c'est pratiquement
indispensable lorsque l'on ne sait pas tester l'\egt à zéro
des \idms dans l'anneau avec lequel on travaille.}}.

%Le \tho de structure ci-après et le lemme qui suit sont des \gnns
%des deux faits précédents.

%:--- Theorem{fact.sfio} --------
\begin{theorem}\label{fact.sfio}\relax \emph{(Système fondamental d'\idms \ortsz)}\\
Soit $(e_1,\ldots ,e_n)$ un \sfio d'un
anneau $\gA$, et  $M$  un \Amoz.
Notons $\gA_i=\aqo{\gA}{1-e_i}\simeq\gA[1/e_i]$. Alors:
\[\preskip.3em \postskip.5em\arraycolsep2pt
\begin{array}{rcl}
\gA  & \simeq  &  \gA_1\times\cdots
\times \gA_n, \\
M  & =  &  {e_1}M\oplus\cdots\oplus {e_n}M .
\end{array}
\]
\end{theorem}
%---end-theorem-----------------

Notez que ${e_1}M$ est un \Amo et un $\gA_1$-module, mais que
ce n'est pas un $\gA_2$-module (sauf s'il est nul). 

%:2012 
Le lemme suivant donne une réciproque du \thref{fact.sfio}

%--- Lemma{lemfacile}-------------
\begin{lemma}
\label{lemfacile}
Soient $(\fa_i)_{i\in\lrbn}$  des idéaux de $\gA$. On a
$\gA=\bigoplus_{i\in\lrbn}\fa_i$ \ssi il existe un \sfio
$(e_i)_{i\in\lrbn}$ tel que $\fa_i=\gen{e_i}$ pour $i\in\lrbn$.
Dans ce cas le \sfio est déterminé de manière unique.
\end{lemma}
%--- end-lemma-----------------------------------------
%-----------------begin proof------------------
\begin{proof}
Supposons  $\gA=\bigoplus_{i\in\lrbn}\fa_i$. On a des $e_i\in\fa_i$
tels que $\som_ie_i=1$, et comme $e_ie_j\in\fa_i\cap\fa_j=\{0\}$ pour
$i\neq j$, on obtient bien un \sfioz.\\ 
En outre si $x\in\fa_j$, on a
$x=x\som_ie_i=xe_j$ et donc $\fa_j=\gen{e_j}$. \\
L'implication réciproque est \imdez.
L'unicité résulte de celle d'une écriture d'un \elt dans une
somme directe.
\end{proof}
%-----------------end proof------------------

Voici maintenant deux lemmes très utiles.

%:--- Lemme{lem2ide.idem}---------
\begin{lemma}
\label{lem2ide.idem}\relax \emph{(Lemme de l'\id engendré par un \idmz)}\\
Un \id  $\fa$ est engendré par un \idm \ssi

\snic{\fa+\Ann\,\fa=\gen{1}.}
\end{lemma}
%--------------
%
\begin{proof}
Tout d'abord, si $e$ est \idmz, on a $\Ann\,\gen{e}=\gen{1-e}$.
Pour l'implication réciproque, soit $e\in\fa$ tel que $1-e\in\Ann\,\fa$. Alors $e(1-e)=0$, donc~$e$ est \idmz. Et pour tout $y\in\fa$, $y=ye$,
donc $\fa\subseteq \gen{e}.$
\end{proof}
%

%:--- Lemme{lem.ide.idem}---------
\begin{lemma}
\label{lem.ide.idem}\relax \emph{(Lemme de l'\itf \idmz)}\\
Si $\fa$ est un \itf \idm (i.e., $\fa =\fa^2 $)
dans  $\gA$, alors $\fa=\gen{e}$ avec $e^2=e$ entièrement
déterminé par
$\fa$.%
\index{Lemme de l'ideal de type@Lemme de l'\itf \idmz}
\end{lemma}
%--------------
%----begin{proof----------
\begin{proof} On utilise le truc du \deterz.
%(cf. \cite{Nor} chap. 4 exercice 11, p. 129)
On considère un \sgrz~\hbox{$(a_1,\ldots a_q)$} de $\fa$ et le
vecteur colonne ${\ua}=\tra{[\,a_1\;\cdots\; a_q\,]}$.
\\
Puisque $a_j\in\fa^2$ pour $j\in\lrbq$, il existe $C\in\MM_q(\fa)$ telle 
\hbox{que ${\ua}=C\,{\ua}$}, \hbox{donc
$(\I_q-C)\,{\ua}=\uze$} et $\det(\I_q-C)\,\ua=\uze$. Or $\det(\I_q-C)=1-e$
\hbox{avec $e \in \fa$}. Donc $(1-e)\fa=0$, et l'on applique le lemme~\ref{lem2ide.idem}.\\
Enfin, l'unicité de $e$ est déjà dans le lemme~\ref{lemfacile}.
\end{proof}
%----end{proof----------

%:    Thm{restes chinois}-------

Rappelons enfin le {\tho chinois}, outil très efficace, qui cache un \sfioz.
Des \ids $\fb_1$, \ldots, $\fb_\ell$ d'un anneau~$\gA$ sont dit \ixc{comaximaux}{ideaux@\ids ---}
lorsque $\fb_1+\cdots+\fb_\ell=\gen{1}$.

\CMnewtheorem{Thresteschi}{Théorème  des restes chinois}{\itshape}
%:     Theorem{Thresteschi}
\begin{Thresteschi} ~\label{restes chinois}\\
Soient dans $\gA$ des \ids  $(\fa_i)_{i\in\lrbn}$
deux à deux \com  et
$\fa=\bigcap_i \fa_i$.
\begin{enumerate}
\item On a l'\egt $\fa=\prod_i \fa_i$,
\item l'application canonique $\gA/\fa\to\prod_i
\gA/\fa_i$ est un \isoz,
\item il existe  $e_1$, $\dots$, $e_n$
dans  $\gA$ tels que
$\fa_i=\fa+\gen{1-e_i}$ et les $\pi_{\gA,\fa}(e_i)$ forment un \sfio de~$\gA/\fa$.
\end{enumerate}
\end{Thresteschi}
%--------- fin theorem ---------------------------------------------- 

 Comme corolaire on obtient le lemme des noyaux.

%:    Lemma{lemDesNoyaux}
\begin{lemma}\label{lemDesNoyaux} \emph{(Lemme des noyaux)}\\
Soit $P=P_1\cdots P_\ell\in\AX$ et une \Ali $\varphi:M\to M$
 vérifiant $P(\varphi)=0$. On suppose que
les $P_i$ sont deux à deux \comz. Notons $K_i=\Ker \big(P_i(\varphi)\big)$, $Q_i=\prod_{j\neq i}P_j$.
%H  index modifié
Alors on a:\index{Lemme des noyaux}
\[ 
\begin{array}{c} \dsp
K_i=\Im \big(Q_i(\varphi)\big),\;  M =\bigoplus_{j=1}^\ell K_j \hbox{ et }    %\\[1mm] 
\Im \big(P_i(\varphi)\big)=\Ker \big(Q_i(\varphi)\big)=\bigoplus_{j\neq i}K_i.    
\end{array}
\]
\end{lemma}
\begin{proof}
On considère l'anneau $\gB=\aqo{\AX}{P}$. Le module $M$ peut être vu comme
un \Bmo pour la loi $(Q,y)\mapsto Q\cdot y=Q(\varphi)(y)$. On applique alors
le \tho des restes chinois et le \tho de structure~\vref{fact.sfio}.
\\
%Ceci est une manière abstraite de présenter le calcul plus classique
%suivant. 
%\\
Cette \dem résume le calcul plus classique suivant.
\`A partir des \egts $U_{ij}P_i+U_{ji}P_j=1$, on obtient des
\egts $U_iP_i+V_iQ_i=1$ ainsi qu'une \egt $\sum_iW_iQ_i=1$. Notons $p_i=P_i(\varphi)$, $q_i=Q_i(\varphi)$ etc.\\
 Alors, tous les \endos obtenus commutent et l'on obtient des \hbox{\egts $p_iq_i=0$}, $u_ip_i+v_iq_i=\Id_M$,  $\sum_iw_iq_i=\Id_M$.
Le lemme en découle facilement.
\end{proof}
%

%--- Sec{Un peu d'algèbre extérieure}------
\section{Un peu d'algèbre extérieure}
\label{secCramer}
%-------------------------------

\begin{flushright}
{\em Qu'un \sli \hmg de $n$ équations à $n$ inconnues \\
admette (sur un corps discret) une solution non triviale\\
 \ssi le \deter du \sys est nul, \\
voilà un fait d'une importance capitale dont on \\
n'aura
jamais fini de mesurer la portée.}\\
Anonyme \\
\end{flushright}

\begin{flushright}
{\em \'Eliminons, éliminons, éliminons \\
les éliminateurs de l'élimination!}\\
Poème \mathe (extrait)\\
S. Abhyankar \\
\end{flushright}

Quelques exemples simples illustrent ces idées dans la section
présente.

%:--- SUBSec Sous-modlibres en facteur direct
\subsec{Sous-modules libres en facteur direct (splitting off)}

Soit $k\in\NN$. Un \emph{module} \ixc{libre de rang $k$}{module ---} est par
\dfn un \Amo isomorphe à $\Ae{k}$. \index{module!libre de rang $k$}
Si $k$ n'est pas précisé, on dira \emph{module} \ixc{libre de rang
fini}{module ---}.%
\index{module!libre de rang fini}\index{rang!d'un module libre}

%:h2014 rajout ci-dessous
Lorsque $\gA$ est un \cdi on parle indifféremment d'\emph{\evc de dimension finie} ou \emph{de rang fini}.\index{espace vectoriel!libre de dimension finie}%
\index{libre de dimension finie!espace vectoriel}%
\index{dimension!d'un espace vectoriel}

Les modules dont la structure est la plus simple sont les modules libres
de rang fini. On est donc
intéressé par la possibilité d'écrire un module arbitraire $M$
sous la forme $L\oplus N$ où $L$ est un module libre de rang fini. Une
réponse (partielle) à cette question est
donnée par l'algèbre extérieure.

% Proposition{propSplittingOffAlgExt}
\begin{proposition}\label{propSplittingOffAlgExt}
\emph{(Splitting off)} \\
 Soient $a_1$, \ldots, $a_k$ des \elts d'un \Amo $M$, alors \propeq
%-----------------begin enum------------------
\begin{enumerate}
\item Le sous-module
$L=\gen{a_1,\ldots ,a_k}$ de $M$
est libre de base $(a_1,\ldots ,a_k)$ et il est facteur direct
de $M$.
\item Il existe une forme $k$-\lin alternée $\varphi:M^k\to\gA$
qui satisfait l'\egt $\varphi(a_1,\ldots ,a_k)$ $=1$.
\end{enumerate}
%-----------------end enum------------------
\end{proposition}
\begin{proof}
\emph{1} $\Rightarrow$ \emph{2}. Si $L\oplus N=M$, si $\pi:M\to L$ est la \prn \paralm à
$N$, et \hbox{si $\theta_j:L\to\gA$} est la $j$-ième forme \coo pour la
base $(a_1,\ldots ,a_k)$, on définit
$$\preskip -.4em \postskip.4em 
\varphi(\xk)=\det\big(\big(\theta_j(\pi(x_i))\big)_{i,j\in \lrbk}\big). 
$$
\emph{2} $\Rightarrow$ \emph{1}.  On définit l'\ali $\pi:M\to M$ par
$$\preskip.2em \postskip.4em 
\pi(x)=\som_{j=1}^k\,\varphi(\underbrace{a_1,\ldots ,x,\ldots ,a_k}_{(x
\mathrm{\,\,est\,\, en\,\, position\,\,} j)})\,a_j
.
$$
On a \imdt $\pi(a_i)=a_i$ et $\Im\pi\subseteq
L:=\gen{a_1,\ldots ,a_k}$, donc $\pi^2=\pi$ et $\Im\pi=L$.
Enfin, si $x=\sum_j\lambda_ja_j=0$, alors $\varphi(a_1,\ldots ,x,\ldots
,a_k)=\lambda_j=0$ (avec $x$ en position $j$).
\end{proof}
Cas particulier: pour $k=1$ on dit que l'\elt $a$ de $M$ est
\ixc{unimodulaire}{elem@\elt --- d'un module} lorsqu'il existe une
forme \lin $\varphi:M\to\gA$ tel que $\varphi(a)=1$.
%:h2014 rajout ``unimodulaire''
Dire que le vecteur $b=(\bn)\in\gA^{n}$ est \umd revient à dire que les~$b_i$ sont \comz.
On dit aussi dans ce cas que la suite $(\bn)$ est 
\emph{\umdz}.\index{unimodulaire!suite (ou vecteur) ---}\index{suite!unimodulaire}\index{vecteur unimodulaire}

%:--- SUBSec Rang d'un module libre
\subsec{Le rang d'un module libre}

Comme nous allons le voir,
le rang d'un module libre est un entier bien déterminé si l'anneau
n'est pas trivial.
Autrement dit, deux \Amosz~\hbox{$M\simeq \Ae{m}$} et $P\simeq \gA^{p}$ avec $m\neq p$
ne peuvent être isomorphes que si $1=_\gA0$.\label{Nota1rang}

Nous utiliserons la notation $\rg_\gA(M)=k$ (ou $\rg(M)=k$ si
$\gA$ est clair d'après le contexte) pour indiquer qu'un module
(supposé libre) est de rang~$k$.

Une \dem savante consiste à dire que, si $m>p$, la puissance
extérieure
$m$-ième de $P$ est $\{0\}$ tandis que celle de $M$ est isomorphe à
$\gA$
(c'est pour l'essentiel la preuve faite dans le
corolaire~\ref{corprop inj surj det}).

La même \dem peut être présentée de façon plus \elr
comme suit.
Rappelons tout d'abord la formule de Cramer de base.
Si~$B$ est une matrice carrée d'ordre $n$, nous notons $\wi{B}$ ou
$\Adj B$ la matrice \emph{cotransposée} (on dit parfois, \emph{adjointe}). La forme
\elr des
\idcs s'écrit alors:
\label{NOTACotrans}
\index{cotransposée!matrice ---}\index{adjointe!matrice ---}
\index{matrice!adjointe (cotransposée)}
%------begin equation-----------
\begin{equation}\label{eqIDC1}
A\; \Adj( A)=\Adj( A)\; A=\det( A)\;\I_n.
\end{equation}
%---------------------end equation--------------
Cette formule, jointe à la formule du produit \gui{$\det(AB)=\det (A)\det (B)$},
implique qu'une matrice carrée $A$ est \iv \ssi son \deter est \ivz,
ou encore si elle est \iv d'un seul coté, et que son inverse est alors égal à $(\det A)^{-1}\Adj A$.
\\
On considère maintenant deux \Amos $M\simeq \Ae{m}$ et $P\simeq \gA^{p}$ avec $m\geq
p$  et une  \ali surjective $\varphi: P\rightarrow M$.
Il existe donc une
\ali $\psi :M\rightarrow P$ telle que $\varphi \circ \psi =\Id_M$.
Ceci correspond à deux matrices
$A\in\Ae{m\times p}$ et $B\in\Ae{p\times m}$  avec $AB=\I_m$. Si $m=p$,
la matrice  $A$ est \iv d'inverse $B$  et $\varphi$ et $\psi$ sont des \isos
réciproques.
Si~$m>p$, on a  $AB=A_1B_1$  avec $A_1$  et $B_1$ carrées obtenues à
partir de~$A$ et~$B$ en complétant par des zéros ($m-p$ colonnes pour
$A_1$, $m-p$ lignes pour~$B_1$).

\vspace{-1.2em}
{\small$$
A_1 =\blocs{.5}{1.3}{0}{1.8}{}{}{$0$ \\[1mm] $~\vdots~$ \\[1mm]$0$}{$A$}\,,
\qquad
B_1 =\blocs{0}{1.8}{.5}{1.3}{}{$\begin{array}{ccc}0 & \cdots & 0\end{array}$}{}{$B$}\,,
\qquad
A_1 B_1 = {\rm I}_m.
$$}

\vspace{-.7em}
Ainsi  $1=\det\I_m=\det (AB)=\det(A_1B_1)=\det(A_1)\det(B_1)=0.$

Dans cette \dem on voit clairement apparaître la commutativité de
l'anneau (qui est vraiment \ncrz). Résumons.

%--- Proposition{propDimMod1}--
\begin{proposition}
\label{propDimMod1}
Soient  deux \Amos $M\simeq \Ae{m}$ et $P\simeq \gA^{p}$ et une  \ali
surjective $\varphi: P\rightarrow M$.
%-----------------begin enum------------------
\begin{enumerate}
\item  Si $m=p$, alors $\varphi$ est un \isoz. Autrement dit, dans un
module~$\Ae{m}$ tout \sys \gtr de $m$ \elts est une base.
\item  Si $m>p$, alors $1=_\gA0$. Et si l'anneau n'est pas
trivial, $m>p$ est impossible.
\end{enumerate}
%-----------------end enum------------------
\end{proposition}
%--- end-proposition----------------------------------------

Dans la suite ce \tho de classification important apparaîtra souvent
comme corolaire de \thos plus subtils, comme par exemple le
\thref{prop unicyc} ou le \thoz~\ref{prop quot non iso}.

%:--- SUBSection{subsecPuissExt}---
\subsec{Puissances extérieures d'un module}
\label{subsecPuissExt}
%-----------------------------------------

\noi{\bf Terminologie.} Rappelons que l'on appelle \ix{mineur} d'une matrice
$A$ tout \deter d'une matrice carrée extraite de $A$ sur certaines lignes
et certaines colonnes. On parle de \emph{mineur d'ordre $k$}
lorsque la matrice carrée extraite est dans $\Mk(\gA)$.
Lorsque $A$ est une matrice carrée, un \emph{mineur principal} est un mineur correspondant à une matrice extraite pour le même
ensemble d'indices sur les lignes et sur les colonnes.
Par exemple \hbox{si $A\in\Mn(\gA)$}, le \coe de $X^k$ dans le \pol $\det(\In+XA)$
est la somme des mineurs principaux d'ordre $k$ de $A$.
Enfin, on appelle \emph{mineur principal dominant} un mineur principal
en position nord-ouest, \cad obtenu en extrayant la matrice sur les premières lignes et les premières colonnes.\eoe%
\index{mineur!d'ordre $k$}%
\index{mineur!principal}%
\index{mineur!principal dominant}%
\index{principal!mineur ---}

\medskip
Soit $M$ un \Amoz.
Une application $k$-\lin alternée $\varphi :M^k\to P$ est
appelée une \ixc{puissance extérieure}{d'un module} $k$-ième 
du \Amo $M$  si
toute \ali alternée $\psi :M^k\to R$ s'écrit de
manière unique sous la forme~\hbox{$\psi=\theta\circ\varphi$}, où $\theta$
est une \Ali de $P$ vers~$R$.
\Pnv{M^k}{\varphi}{\psi}{P}{\theta}{R}{}{applications  $k$-linéaires alternées}{\alisz.}
\label{PuissExtMod}

Il est clair que  $\varphi :M^k\to P$ est unique au sens catégorique,
\cad que pour toute autre puissance extérieure  $\varphi' :M^k\to P'$ il
y a 
une
\ali unique $\theta:P\to P'$ qui rend le diagramme convenable commutatif,
et que $\theta$ est un \isoz.

On note alors $\Al{k}\!M$ ou $\Vi_\Ae{k} M$ pour $P$ et $\lambda_{k}(x_{1},\ldots,x_{k})$ ou $x_1\vi\cdots\vi x_k$ pour~$\varphi(\xk)$.

L'existence d'une puissance extérieure $k$-ième pour tout module $M$
résulte de considérations \gnles analogues à celles que nous
détaillerons pour le produit tensoriel  \paref{ProdTens} de la
section~\ref{secStabPf}.

La théorie la plus simple des puissances extérieures, analogue à la
théorie \elr du \deterz, démontre que si $M$ est un module libre
ayant une base de $n$ \elts $(a_1,\ldots ,a_n),$ alors  $\Al{k}M$ est nul
si $k>n$, et sinon c'est un module libre ayant pour base les $n \choose k$
$k$-vecteurs $a_{i_1}\vi \cdots\vi a_{i_k}$, où $(i_1,\ldots ,i_k)$
parcourt l'ensemble des $k$-uplets strictement croissants d'\elts de
$\lrbn$. \\
En particulier,  $\Al{n}M$ est libre de rang~1 avec pour base
$a_{1}\vi \cdots\vi a_{n}$.

\`A toute \Ali $\alpha :M\to N$ correspond une unique \Ali $\Al{k}\alpha :
\Al{k}M\to\Al{k}N$ vérifiant l'\egt 

\snic{\big(\Al{k}\alpha \big)(x_1\vi\cdots\vi
x_k)=\alpha(x_1)\vi \cdots \vi \alpha(x_k)}

%\sni
pour tout $k$-vecteur
$x_1\vi\cdots\vi x_k$ de $\Al{k}M$. L'\ali $\Al{k}\alpha$ s'appelle
la \ixc{puissance extérieure}{d'une \aliz} $k$-ième de l'\aliz~$\alpha.$

En outre on a  $\big(\Al{k}\alpha\big) \circ
\big(\Al{k}\beta\big) =\Al{k}(\alpha \circ \beta )$ quand $\alpha
\circ\beta $ est défini.
En bref, chaque $\Al{k}(\bullet)$ est un foncteur.

Si $M$ et $N$ sont libres de bases respectives  $(a_1,\ldots ,a_n)$ et
 $(b_1,\ldots ,b_m)$,  et si~$\alpha$ admet la matrice $H$ sur ces bases,
alors  $\Al{k}\alpha$ admet la matrice notée~$\Al{k}H$ sur les bases
correspondantes de $\Al{k}M$ et  $\Al{k}N$. Les \coes de cette
matrice sont tous les mineurs d'ordre $k$ de la matrice~$H$.

\penalty-5000
%:--- Subsec{Ideaux determinantiels}------
\subsec{Idéaux déterminantiels}
\label{secIdd}
%-------------------------------

%--- Defi{defIdDet}-------------
\begin{definition}\label{defIdDet}
Soient $G\in\Ae{n\times m}$ et $ k\in\lrb{1..\min(m,n)}$,
 {\em  l'\idd d'ordre $k$ de la matrice $G$}  est l'\idz,
noté $\cD_{\gA,k}(G)$ ou $\cD_k(G)$,  engendré  par  les   mineurs d'ordre  $k$
de   $G$.
\index{ideaux determ@idéaux déterminantiels!d'une matrice}
Pour $k\leq 0$ on pose par convention~\hbox{$\cD_k(G)=\gen{1}$}, et pour
$k> \min(m,n)$,  $\cD_k(G)=\gen{0}$.
\end{definition}

Ces conventions sont naturelles car elles permettent d'obtenir en toute
\gnt les \egts suivantes.
\begin{itemize}
\item Si $H=\blocs{.6}{.9}{.6}{.5}{$\I_r$}{$0$}{$0$}{$G$}$\ , pour tout $k\in\ZZ$
on a $\cD_k(G)=\cD_{k+r}(H)$.
\item Si $H=\blocs{.6}{.9}{.7}{.5}{$0$}{$0$}{$0$}{$G$}$\ , pour tout $k\in\ZZ$ on a
$\cD_k(H)=\cD_{k}(G)$.
\end{itemize}
%\begin{itemize}
%%
%\item Si $H=\cmatrix{\I_r&0\cr0&G}$, pour tout $k$
%on a $\cD_k(G)=\cD_{k+r}(H)$.
%%
%\item Si $H=\cmatrix{0&G}$, ou  $H=\cmatrix{0\cr G}$, pour tout $k$ on a
%$\cD_k(G)=\cD_{k}(H)$.
%%
%\end{itemize}
%

%--- Fact{fact.idd inc}-----------
\begin{fact}
\label{fact.idd inc}\relax
Pour toute matrice $G$ de type $n\times m$ on a les inclusions
%---  equation eqfact.idd inc ---
\begin{equation}\label{eqfact.idd inc}\relax
\{0\}=\cD_{1+\min(m,n)}(G) \subseteq\cdots\subseteq\cD_1(G)
\subseteq\cD_0(G)=\gen{1}=\gA
%---------------------end equation--------------
\end{equation}
Plus \prmt pour tout $k,r\in\NN$ on a une inclusion
%----- equation--eqfact.idd inc2-----------
\begin{equation}\label{eqfact.idd inc2}\relax
\cD_{k+r}(G)\subseteq\cD_{k}(G)\,\cD_{r}(G)
\end{equation}
%---------------------end equation--------------
\end{fact}
%--------------

En effet, tout mineur d'ordre $h+1$ s'exprime comme \coli de mineurs
d'ordre $h$. Et l'inclusion (\ref{eqfact.idd inc2}) s'obtient avec le
développement de Laplace du \deterz. 
\hum{Un exo à mettre pour
le développement de Laplace du \deterz?}

%:     Fact{fact.idd.sousmod}
\begin{fact}\label{fact.idd.sousmod}\label{factlnlImage}
Soient $G_1\in\Ae{n\times m_1}$,  $G_2\in\Ae{n\times m_2}$ et $H\in \gA^{p\times n}$.
\begin{enumerate}
\item Si $\Im G_1\subseteq\Im G_2$, alors pour tout entier $k$ on a
$\cD_k(G_1) \subseteq \cD_k(G_2)$.
\item Pour tout entier $k$, on a $\cD_k(HG_1) \subseteq \cD_k(G_1)$.
\item Les \idds d'une matrice $G\in\Ae{n\times m}$ ne dépendent que de la classe d'\eqvc
du sous-module image de $G$ (i.e., ils ne dépendent que de $\Im G$,
à \auto près du module $\Ae{n}$).
\item  En particulier, si  $\varphi$ est une \ali entre modules libres de rangs finis,
les \idds  d'une matrice de $\varphi$ ne dépendent pas des bases choisies.
On les note $\cD_k(\varphi)$ et on les appelle les
{\em  \idds de l'\ali $\varphi$}.
\end{enumerate}
\index{ideaux determ@idéaux déterminantiels!d'une \ali (modules libres)}
\index{ideal@idéal!déterminantiel}
\end{fact}
\begin{proof}
\emph{1.} Chaque colonne de $G_1$ est une \coli de colonnes de $G_2$. On conclut par
la multilinéarité du \deterz. \\
\emph{2.} Même raisonnement en remplaçant les colonnes par les lignes.\\
Enfin, \emph{3} implique \emph{4} et résulte des deux points précédents.
\end{proof}

\rem Un \idd est donc attaché essentiellement à un sous-\mtf $M$ d'un
module libre $L$. Mais c'est la structure de l'inclusion $M\subseteq L$
et non pas seulement la structure de $M$ qui intervient pour déterminer
les \iddsz. Par exemple $M=3\ZZ\times 5\ZZ$ est un \hbox{sous-\ZZmo} libre de
$L=\ZZ^2$ et ses \idds sont~\hbox{$\cD_1(M)=\gen{1}$}, $\cD_2(M)=\gen{15}$.
Si l'on remplace $3$ et $5$ par $6$ et $10$ par exemple, on obtient un autre sous-module
libre, mais la structure de l'inclusion est différente puisque les \idds sont
maintenant~$\gen{2}$ et~$\gen{60}$.
\eoe

%--- Fact{fact.idd prod}----------
\begin{fact}
\label{fact.idd prod}\relax
Si $G$ et $H$ sont des matrices telles que $GH$ est définie, alors, pour
tout $n\geq 0$ on a
%--------------------begin equation---------------
\begin{equation}\label{eqfact.idd prod}\relax
\cD_n(GH)\subseteq \cD_n(G)\,\cD_n(H)
\end{equation}
%---------------------end equation--------------
\end{fact}
%--------------

%
\begin{proof}
Le résultat est clair pour $n=1$. 
Pour $n>1$, on se ramène
au \hbox{cas $n=1$} en notant que les mineurs d'ordre~$n$ de~$G$,~$H$ et $GH$
représentent les \coes des matrices \gui{puissance extérieure
$n$-ième de $G$, $H$ et $GH$}  \big(en tenant compte de l'\egt $\Al{n}(\varphi\psi )=\Al{n}\varphi
\circ  \Al{n}\psi$\big).
\end{proof}

L'\egt suivante est immédiate.
%------begin equation-------
\begin{equation}\label{eqIDDSDIR}
\cD_n(\varphi\oplus\psi)=\som_{k=0}^n\cD_{k}(\varphi)\,\cD_{n-k}(\psi)
\end{equation}
%---------------------end equation--------------

%:--- Subsec{Rang d'une matrice}------
\subsec{Rang d'une matrice}
\label{secRangmat}
%-------------------------------

%--- Definition{defRangk}------
\begin{definition}
\label{defRangk} ~\\
Une \ali $\varphi$ entre modules libres de
rangs finis est dite 
\begin{itemize}
\item \emph{de rang $\leq k$} si $\cD_{k+1}(\varphi)=0$,
\item \emph{de rang $\geq k$} \hbox{si $\cD_{k}(\varphi)=\gen{1},$} 
\item  \emph{de rang $k$} si elle est à la fois de \hbox{rang
$\geq k$} et de rang~$\leq k$.
\end{itemize}
\index{rang!d'une matrice}%
\index{rang!d'une application linéaire}%
\index{matrice!de rang $\geq k$}%
\index{matrice!de rang $\leq k$}%
\index{matrice!de rang $k$}
\end{definition}
%--- end-definition------------------------------------

Nous utiliserons les notations $\rg(\varphi)\geq k$ et $\rg(\varphi)\leq k$,
conformément à la \dfn précédente, sans présupposer que $\rg(\varphi)$ soit défini.
Seule l'écriture~\hbox{$\rg(\varphi)= k$} signifiera que le rang est défini.

Nous généra\-li\-serons plus loin cette \dfn au cas d'\alis entre
\mptfsz: voir la notation \ref{notaRgfi} ainsi que les exercices~\hbox{\ref{exoIDDPTF1},
\ref{exoIDDPTF2}} et \ref{exoLocSimpPtf}.

\medskip \comm
\Llec doit prendre garde qu'il n'existe pas de \dfn universellement
acceptée
pour \gui{matrice de rang $k$} dans la
litté\-ra\-ture.
En lisant un autre ouvrage, \il doit d'abord s'assurer de la \dfn
adoptée par l'auteur. Par exemple dans le cas d'un anneau intègre
$\gA$, on trouve souvent
le rang défini comme celui de la matrice vue dans le corps de fractions
de $\gA$. Néanmoins une matrice de rang $k$ au sens de la \dfnz~\ref{defRangk} est \gnlt de rang $k$  au sens des autres auteurs.
\eoe

\medskip
Le \plgc suivant est une conséquence immédiate du \plg de base.

%--- Prin loc glob conc plccRangMat}
\begin{plcc}
\label{plccRangMat}\emph{(Rang d'une matrice)}\\
Soient $S_1$, $\dots$, $S_n$ des \moco de $\gA$ et  $B$
une matrice $\in \Ae{m\times p}$.
Alors \propeq
\begin{enumerate}
\item  La matrice est de rang $\leq k$ (resp. de  rang $\geq k$) sur
$\gA$.
\item Pour $i\in\lrbn,$
la matrice est de rang $\leq k$ (resp. de  rang $\geq k$)
sur~$\gA_{S_i}$.
\end{enumerate}
\end{plcc}
%--- end-plcc-----------------------------------------

%%%%%%%%%%%%%%%%%%%%%%%%%%%%%%%%%%%%%%%%%%%%%%%%%%%%%%%%%%%%%%%%%%%%%%%%%%%
%:--- Subsec{Pivot de Gauss}------
\subsec{Méthode du pivot \gneez}
\label{secPivotdeGauss}
%-------------------------------

\rdb
\vspace{3pt}

\noi{\bf Terminologie.} 

1) Deux matrices  sont dites
\emph{\eqvesz}
lorsque l'on  passe de l'une à l'autre en multipliant à droite et à
gauche par des matrices inversibles.

\noindent 2)
Deux matrices carrées dans $\Mn(\gA)$ sont dites \ixd{semblables}{matrices}
lorsqu'elles repré\-sentent le même \endo
de $\Ae n$ sur deux bases (distinctes ou non), autrement dit lorsqu'elles
sont conjuguées pour l'action $(G,M)\mapsto GMG^{-1}$ de $\GLn(\gA)$ sur 
$\Mn(\gA)$.

\noindent 3)
Une \emph{manipulation \elr de lignes} sur une matrice de $n$ lignes
consiste en le remplacement d'une ligne $L_i$ par la ligne $L_i+\lambda L_j$
avec~\hbox{$ i\neq j$}. On la note  aussi $L_i\aff L_i+\lambda L_j$. Cela correspond à la multiplication à gauche
par une matrice, dite \emph{\elrz},  notée $\rE^{(n)}_{i,j}(\lambda)$
(ou, si le contexte le permet, $\rE_{i,j}(\lambda)$). Cette matrice est
obtenue à partir de $\In$ par la même manipulation \elr
de lignes.
\\
La multiplication à droite par la
même matrice $\rE_{i,j}(\lambda)$ correspond, elle,
à la \emph{manipulation \elr de colonnes} (pour une
matrice qui possède $n$ colonnes) qui
transforme la matrice $\In$ en $\rE_{i,j}(\lambda)$: $C_j\aff C_j+\lambda C_i$.

\noindent 4)
Le sous-groupe de $\SLn(\gA)$
engendré par les matrices \elrs est appelé le \emph{groupe \elrz} et il est noté
$\En(\gA)$. Deux matrices sont dites \emph{\elrt \eqvesz}
lorsque l'on peut passer de l'une à l'autre
par des manipulations \elrs de lignes et de colonnes.% 
\index{matrices!equiva@équivalentes}%
\index{equivalentes@équivalentes!matrices ---}%
\index{elementairement@élémentairement équi\-valentes!matrices ---}%
\index{matrices!elem@élémentairement équi\-valentes}%
\index{matrice!elem@élémentaire}%
\index{groupe elem@groupe élémentaire}%
\index{manipulation!elem@élémentaire}%
\index{elementaire@élémentaire!matrice ---}%
\index{elementaire@élémentaire!groupe ---}%
\index{elementaire@élémentaire!manipulation --- de lignes}%
\label{NOTAEn}%
\eoe

%:  Lem du mineur inv {lem.min.inv}
\CMnewtheorem{lemmininv}{Lemme du mineur inversible}{\itshape}
\begin{lemmininv}\label{lem.min.inv}\index{Lemme du mineur inversible}
\emph{(Pivot généralisé)}\\
Si une matrice $G\in\gA^{q\times m}$ possède un mineur d'ordre $k\leq \min(m,q)$
inversible% (où $k\leq \min(m,q)$)
, elle est \eqve à une matrice

\snic{
\cmatrix{
    \I_{k}   &0_{k,m-k}      \cr
    0_{q-k,k}&  G_1      },}

%\sni
avec $\cD_r(G_1)=\cD_{k+r}(G)$ pour tout $r\in\ZZ$. 
\end{lemmininv}
%--------- fin lemma ---------------------------------------------- 

%----begin{proof---------------
\begin{proof}
En permutant éventuellement les lignes et les colonnes on ramène le
mineur inversible en haut à gauche. Puis en multipliant  à droite  (ou
à gauche) par une  matrice inversible   on se ramène  à la forme

\snic{G'\;=\;
\cmatrix{
\I_k & A     \cr
  B &  C
}, }

puis par des manipulations \elrs de lignes et de colonnes, on
obtient
$$
%\snic{
G''\;=\;
\cmatrix{
   \I_{k}   &0_{k,m-k}      \cr
    0_{q-k,k}&     G_1}.
%    }
$$
%\sni
Enfin $\cD_r(G_1)=\cD_{k+r}(G'')=\cD_{k+r}(G)$ pour tout $r\in\ZZ$.
\end{proof}
%----end{proof------------------

Comme conséquence immédiate on obtient le lemme de la liberté.

%--- Lem de la liberte{lem pf libre}
\CMnewtheorem{lemli}{Lemme de la liberté}{\itshape}
%:     Lemma{lem pf libre}
\begin{lemli}\label{NOTAIkqm} \label{lem pf libre}\index{Lemme de la liberté}
Considérons une matrice $G\in\gA^{q\times m}$ de \hbox{rang $\leq k$} avec $1\leq k\leq \min(m,q)$. Si la
matrice~$G$ possède un mineur d'ordre~$k$ inversible, alors elle est
\eqve à la matrice

\snic{\I_{k,q,m}\;=\;
\cmatrix{
    \I_{k}   &0_{k,m-k}      \cr
    0_{q-k,k}&     0_{q-k,m-k}      }.}

%\sni
Dans ce cas, l'image, le noyau et le conoyau de $G$ sont libres,
respectivement de rangs $k$, $m-k$ et $q-k$. En outre l'image et le noyau
possèdent des supplémen\-taires libres.
\\ 
Si $i_1$, $\ldots$, $i_k$ (resp. $j_1$, $\ldots$, $j_k$) sont
les numéros de lignes (resp. de colonnes) du mineur inversible, alors
les colonnes $j_1$, $\ldots$, $j_k$ forment une base du module $\Im G$, et
$\Ker G$ est le sous-module défini par l'annulation des formes \lins
correspondant aux lignes $i_1$, $\ldots$, $i_k$.
\end{lemli}

%----begin{proof---------------
\begin{proof}
Avec les notations du lemme précédent on a  $\cD_1(G_1)=\cD_{k+1}(G)=\gen{0}$,
\hbox{donc $G_1=0$}.   Le reste est laissé \alecz.
\end{proof}
%----end{proof------------------

La matrice $\I_{k,q,m}$ est appelée une \emph{matrice simple standard}.
On note $\I_{k,n}$ pour $\I_{k,n,n}$ et on l'appelle une
\emph{\mprn standard}.%
\index{matrice!simple standard}%
\index{matrice!de projection standard}

%--- Definition{defnl} ---------
\begin{definition}
\label{defnl}
Une \ali entre modules libres de rangs finis est dite
\emph{simple} si elle peut être représentée par une matrice
 $\I_{k,q,m}$ sur des bases convenables. De même une matrice
est dite \emph{\nlz} lorsqu'elle est \eqve à une
matrice $\I_{k,q,m}$.  
\index{simple!matrice ---}%
\index{simple!application \lin ---}%
\index{matrice!simple}%
\end{definition}
%--- end-definition------------------------------------

%:--- SUBSection{Formule de Cramer generalisee
\subsec{Formule de Cramer généralisée}

Nous étudions dans ce paragraphe quelques \gnns des formules
de Cramer usuelles. Nous les exploiterons dans les paragraphes suivants.

Pour une matrice $A\in\Ae{m{\times}n}$ nous notons
$A_{\alpha,\beta}$ la matrice extraite sur les lignes
$\alpha=\{\alpha_1,\ldots ,\alpha_r\}\subseteq\lrbm$ et les
colonnes
$\beta=\{\beta_1,\ldots ,\beta_s\}\subseteq\lrbn$.
\label{NOTAextraite}

Supposons la matrice $A$ de rang $\leq k$.
Soit $V\in\Ae{m{\times}1}$ un
vecteur colonne tel que la matrice bordée $[\,A\,|\,V\,]$
soit aussi de rang $\leq k$.
Appelons $A_j$ la $j$-ième colonne de $A$.
Soit $\mu_{\alpha,\beta}=\det(A_{\alpha,\beta})$ le mineur
d'ordre $k$ de la matrice $A$ extrait sur les lignes
$\alpha=\{\alpha_1,\ldots ,\alpha_k\}$ et les colonnes
$\beta=\{\beta_1,\ldots ,\beta _k\}$.
Pour $j\in\lrbk$ soit $\nu_{\alpha,\beta,j}$ le \deter
de la même matrice extraite, à ceci près que la
colonne $j$ a été remplacée par la colonne extraite de $V$
sur les lignes $\alpha$.
Alors, on obtient pour chaque couple $(\alpha,\beta)$ de
multi-indices une \idt de Cramer:
%---- equation {eqMPC1} ----
\begin{equation}\label{eqMPC1}
\qquad \mu_{\alpha,\beta}\;V=\som_{j=1}^k
\nu_{\alpha,\beta,j}\,A_{\beta_j}\qquad
\end{equation}
%---------------------end equation--------------
due au fait que le rang de la matrice bordée
$[\,A_{1..m,\beta}\,|\,V\,]$ est
$\leq k$. Ceci peut se relire comme
suit:
\begin{eqnarray}\preskip.0em \postskip.4em
\qquad \mu_{\alpha,\beta}\;V&=& \left[
\begin{array}{ccc}
    A_{\beta_1} & \ldots  & A_{\beta_k}
\end{array}
\right] \cdot  \left[
\begin{array}{c}
      \nu_{\alpha,\beta,1} \\
    \vdots  \\
     \nu_{\alpha,\beta,k}
\end{array}
\right]= \nonumber \\[.1em]
 &=&\left[
\begin{array}{ccc}
    A_{\beta_1} & \ldots  & A_{\beta_k}
\end{array}
\right] \cdot
\Adj(A_{\alpha,\beta})
 \cdot
 \left[
 \begin{array}{c}
     v_{\alpha_1}  \\
     \vdots   \\
     v_{\alpha_k}
 \end{array} \right]=\nonumber\\[.2em]
\label{eqMPC2}
 &=&A \cdot (\I_n)_{1.. n,\beta} \cdot
\Adj(A_{\alpha,\beta}) \cdot (\I_m)_{\alpha,1.. m}\cdot V
\end{eqnarray}

Ceci nous conduit à introduire la notation suivante.
%--- Notation{notaAdjalbe} -----
\begin{notation}
\label{notaAdjalbe}
{\rm Nous notons $\cP_{\ell}$ l'ensemble des parties  de
$\lrb{1..\ell}$ et $\cP_{k,\ell}$ l'ensemble des parties à $k$ \elts de
$\lrb{1..\ell}$. 
\\
Pour  $A\in\Ae{m{\times}n}$ et $\alpha\in
\cP_{k,m},\,\beta\in \cP_{k,n}$ nous notons
$$\preskip.2em \postskip.1em
\Adj_{\alpha,\beta}(A):=
(\I_n)_{1.. n,\beta} \cdot
\Adj(A_{\alpha,\beta}) \cdot (\I_m)_{\alpha,1.. m}.$$
}
\end{notation}
%--- end-notation-----------------------------------------

Par exemple avec la matrice
$$A=\crmatrix{ 5&-5&7&4\cr
9&-1&2&7\cr13&3&-3&10},$$
et les parties $\alpha=\{1,2\}$ et
$\beta=\{2,3\}$, on obtient
{\small
$$A_{\alpha,\beta}=\crmatrix{
-5&7\cr-1&2},\;\Adj(A_{\alpha,\beta})=\crmatrix{
2&-7\cr1&-5}\; \hbox{et}\;
\Adj_{\alpha,\beta}(A)=\crmatrix{ 0&0&0\cr2&-7&0\cr
   1&-5&0\cr0&0&0}.
$$
}

L'\egrf{eqMPC2} s'écrit comme suit, sous l'hypothèse
que $\cD_{k+1}([\,A\,|\,V\,])=0$.
%---- equation {eqGema} (formule miracle)
\begin{equation}\label{eqGema}
\mu_{\alpha,\beta}\;V\;=\;A \cdot \Adj_{\alpha,\beta}(A) \cdot V
\end{equation}
%---------------------end equation--------------

On obtient donc l'\egt ci-dessous, 
sous l'hypothèse que $A$ est de rang~$\leq k$.
%---- equation {eqIGCram} ---
\begin{equation}\label{eqIGCram}
\mu_{\alpha,\beta}\;A\;=\;A \cdot \Adj_{\alpha,\beta}(A) \cdot A
\end{equation}
%---------------------end equation--------------

Les \idcs (\ref{eqGema}) et (\ref{eqIGCram}) fournissent des congruences
qui ne sont soumises à aucune hypothèse: il suffit par exemple de lire
(\ref{eqGema}) dans l'anneau quotient $\gA/\cD_{k+1}([\,A\,|\,V\,])$ pour obtenir la
congruence~(\ref{eqGema2}).

%--- Lemma{lemCram}---
\begin{lemma}
\label{lemCram}  \emph{(Formule de Cramer \gneez)}\\
Sans aucune hypothèse sur la matrice $A$ ou le vecteur $V$, on a \hbox{pour
$\alpha\in \cP_{k,m}$} \hbox{et $\beta\in \cP_{k,n}$} les congruences suivantes.
%---- equation {eqGema2,eqCGCram,eqCGCram2} ---
\begin{eqnarray}\label{eqGema2}
\mu_{\alpha,\beta}\;V&\equiv&A \cdot \Adj_{\alpha,\beta}(A) \cdot V
\qquad \mathrm{mod}\quad \cD_{k+1}([\,A\,|\,V\,])
\\[.1em] \label{eqCGCram}
\mu_{\alpha,\beta}\;A&\equiv&A \cdot \Adj_{\alpha,\beta}(A) \cdot A
\qquad \mathrm{mod}\quad \cD_{k+1}(A).
\end{eqnarray}
%---------------------end equation--------------
\end{lemma}
%--- end-lemma-----------------------------------------

Un cas particulier simple est le suivant avec $k=m\leq n$.
%---- equation {eqGema3} ---
\begin{equation}\label{eqGema3}\preskip.2em \postskip.3em
\mu_{1..m,\beta}\;\I_m\;=\;A \cdot \Adj_{1..m,\beta}(A)\quad (\beta\in
\cP_{m,n})
\end{equation}
%---------------------end equation--------------
Cette \egt est d'ailleurs une conséquence directe de
l'\idc de base (\ref{eqIDC1}). De la même manière on obtient
%---- equation {eqGema4} ---
\begin{equation}\label{eqGema4}\preskip.2em \postskip.2em
\mu_{\alpha,1..n}\;\I_n\;=\; \Adj_{\alpha,1..n}(A)  \cdot A
\quad (\alpha \in \cP_{n,m},\,n\leq m)
\end{equation}
%---------------------end equation--------------

%:--- SUBSection{Une formule magique}----
\subsec{Une formule magique}

Une conséquence immédiate de l'\idc (\ref{eqIGCram})
est l'\idt (\ref{eqIGCram2}) moins usuelle donnée dans le \tho suivant.
De même les \egts
(\ref{eqIGCram3}) et (\ref{eqIGCram4}) résultent facilement de
(\ref{eqGema3}) et (\ref{eqGema4}).

%--- Theorem{propIGCram}----
\begin{theorem}
\label{propIGCram}
Soit  $A\in\Ae{m\times n}$ une matrice de rang $ k$.
On a donc une \egt
%
%\snic{
$\som_{\alpha\in \cP_{k,m},\beta\in \cP_{k,n}}
c_{\alpha,\beta}\,\mu_{\alpha,\beta}=1.$
%}
%
%\sni 
Posons 

\snic{B\;=\;\som_{\alpha\in \cP_{k,m},\beta\in
\cP_{k,n}}\,c_{\alpha,\beta}\,\Adj_{\alpha,\beta}(A).}
%-----------------begin enum------------------
\begin{enumerate}
\item  On a
%---- equation {eqIGCram2}
\begin{equation}\label{eqIGCram2}\preskip.0em \postskip.3em
A\cdot B\cdot A=A.
\end{equation}
En conséquence $A\, B$ est une \prn de rang $k$ et le sous-module 
$\Im A=\Im AB$
est facteur direct dans~$\Ae{m}$.

\item  Si  $k=m$, alors
%---- equation {eqIGCram3}
\begin{equation}\label{eqIGCram3}\preskip.0em \postskip.3em
A\cdot B=\I_m.
\end{equation}

\item  Si $k=n$, alors
%---- equation {eqIGCram4}
\begin{equation}\label{eqIGCram4}\preskip.0em \postskip.3em
B\cdot A=\I_n.
\end{equation}
%---------------------end equation--------------
\end{enumerate}
%-----------------end enum------------------
 \end{theorem}
%--- end-theorem----------------------------------------

L'\idt suivante, que nous n'utiliserons pas dans la suite, est encore plus
miraculeuse.
%--- Proposition{propIGCram2}---
\begin{proposition}
\label{propIGCram2} \emph{(Prasad et Robinson)}\\
Avec  les  hypothèses et les notations du \thref{propIGCram}, 
si l'on a 

\snic{\forall \alpha,\alpha'\in \cP_{k,m},$
$\forall\beta,\beta '\in \cP_{k,n}$
$\;c_{\alpha,\beta}\,c_{\alpha',\beta'}=
c_{\alpha,\beta'}\,c_{\alpha',\beta},}

%\sni
alors
%---- equation {eqIGCramPraRo}
\begin{equation}\label{eqIGCramPraRo}
B\cdot A\cdot B=B.
\end{equation}
%---------------------end equation--------------
\end{proposition}
%--- end-proposition----------------------------------------

%:--- SUBSection{subsecInvGen}--------
\subsec{Inverses généralisés et applications \lnlsz}
 \label{subsecInvGen}
%-----------------------------------------

Soient $E$ et $F$ deux \Amosz,  et une \ali $\varphi:E\rightarrow F$.
On peut voir cette donnée comme une sorte de \sli \gne (un \sli usuel
correspond au cas de modules libres de rang fini). 
Informellement un tel \sli est considéré comme \gui{bien
conditionné} s'il y a une façon systématique de trouver une
solution à l'équation en $x$, $\varphi(x)=y$, à partir de la
donnée $y$, lorsqu'une telle solution existe.
Plus \prmtz, on se demande s'il existe une \aliz~\hbox{$\psi:F\rightarrow E$} vérifiant  $\varphi\big(\psi(y)\big)=y$ chaque fois qu'il
existe une solution~$x$. Cela revient à demander
$\varphi\big(\psi\big(\varphi(x)\big)\big)=\varphi(x)$ pour tout $x\in E$.

Ceci éclaire l'importance de l'\eqrf{eqIGCram2} et conduit à la notion
d'\ingz.

La terminologie concernant les \ings ne semble pas entièrement fixée.
Nous adoptons celle de \cite{Lan}. \\
Dans le livre \cite{Bha}, l'auteur utilise le
terme
\gui{reflexive g-inverse}.

%:    Definition{defIng} et \lnl -----
\begin{definition}
\label{defIng}
Soient $E$ et $F$ deux \Amosz,  et une \ali $\varphi:E\rightarrow F$.
Une \ali $\psi :F\rightarrow E$ est appelée un \ix{inverse
généralisé} de $\varphi$ si l'on a
%---  equation eqdefIng --------
\begin{equation}\label{eqdefIng}
\varphi \circ\psi \circ\varphi =\varphi \quad \mathrm{et}  \quad \psi
\circ\varphi \circ\psi =\psi.
\end{equation}
%---------------------end equation--------------
Une \ali est dite \emph{\lnlz} lorsqu'elle possède un
\ingz. \index{{localement!application linéaire --- simple}}
\end{definition}
%--- end-definition------------------------------------

Le fait suivant est immédiat.
%--- Fact{factIng0}------------------
\begin{fact}
\label{factIng0}
Lorsque $\psi$ est un \ing de $\varphi$, on a:
%-----------------begin item------------------
\begin{enumerate}
\item [--] $\varphi\, \psi$ et $\psi\, \varphi$ sont des \prnsz,
\item [--] $\Im\varphi =\Im\varphi\, \psi$,  $\Im\psi  =\Im\psi\,
\varphi$,  $\Ker\varphi =\Ker\psi \,\varphi$, $\Ker\psi
=\Ker\varphi\,\psi$,
\item [--] $E=\Ker\varphi \oplus \Im\psi$ et
 $F=\Ker\psi  \oplus \Im\varphi$,
\item [--] $\Ker\varphi \simeq \Coker\psi$ et
 $\Ker\psi \simeq \Coker\varphi$.
\end{enumerate}
En outre $\varphi$ et $\psi$ donnent par restriction des \isos réciproques
$\varphi_1$ et $\psi_1$ entre $\Im\psi$ et $\Im\varphi$. Matriciellement on obtient:
\vspace{-1mm}
$$
\bordercmatrix [\lbrack\rbrack]{
           & \Im\psi   %&\oplus   
                       & \Ker\varphi \cr
\Im\varphi & \varphi_1 %& 
                       &   0    \cr
%~~\oplus &   & &        \cr
\Ker\psi   &     0     %& 
                       &0
}~=~\varphi,
\qquad
\bordercmatrix [\lbrack\rbrack]{
           & \Im\varphi      & \Ker\psi \cr
\Im\psi & \psi_1  &   0    \cr
\Ker\varphi   &  0 & 0
}~=~\psi.
$$
%-----------------end item------------------
\end{fact}
%--- end-fact-----------------------------------------

%--- Remarque
\rems\\
1) Si l'on a une \ali $\psi_0$ vérifiant comme dans le \thref{propIGCram} l'\egtz~\hbox{$\varphi\, \psi_0\,\varphi =\varphi$}, on obtient un
\ing de $\varphi$ en posant~\hbox{$\psi=\psi_0\,\varphi \,\psi_0$}.
Autrement dit, une \ali $\varphi$ est \lnl \ssi il existe $\psi$
vérifiant
{$\varphi \,\psi \,\varphi =\varphi$}.

\noindent  2) Une \ali \nl entre modules libres de rangs finis est \lnl
(\vfn immédiate).

\noindent  3) Le \thrf{propIGCram} nous dit qu'une \ali qui possède un
rang $k$ au sens de la \dfnz~\ref{defRangk} est \lnlz.
\eoe

%--- Fact{factInvGenCrois}-----------
\begin{fact}
\label{factInvGenCrois} Soit une \ali $\varphi:\Ae{n}\rightarrow \Ae{m}.$
\Propeq
\vspace{-2pt}
%-----------------begin enum------------------
\begin{enumerate}\itemsep0pt
\item L'\ali $\varphi$ est \lnlz.
\item Il existe  $\varphi\bul :\Ae{m}\rightarrow \Ae{n}$   telle que

\snic{\Ae{n}=\Ker\varphi\oplus\Im\varphi\bul $ et  $\Ae{m}=\Ker\varphi\bul
\oplus\Im\varphi.}
\item Le sous-module $\Im\varphi$ est facteur direct dans $\Ae{m}$.
\end{enumerate}
%-----------------end enum------------------
\end{fact}
%--- end-fact-----------------------------------------
%-----------------begin proof------------------
\begin{proof}
\emph{1} $\Rightarrow$ \emph{2.} Si $\psi$ est un \ing de $\varphi$, on peut prendre
$\varphi\bul=\psi$.
\\
\emph{2} $\Rightarrow$ \emph{3.}  \'Evident.
\\
\emph{3} $\Rightarrow$ \emph{1.}
Si $\Ae{m}=P\oplus\Im\varphi$, notons $\pi:\Ae{m}\to\Ae{m}$ la projection
sur $\Im\varphi$ parallèlement à $P$. Pour chaque vecteur $e_i$ de
la base canonique de $\Ae{m}$ il existe un \elt $a_i$ de $\Ae{n}$ tel que
$\varphi(a_i)=\pi(e_i)$. On définit $\psi :\Ae{m}\rightarrow \Ae{n}$  par
$\psi(e_i)=a_i$. Alors, $\varphi \circ  \psi =\pi$
et  $\varphi\circ \psi\circ \varphi =\pi\circ \varphi=\varphi.$
Et $\psi \circ \varphi\circ  \psi $ est un
\ing de $\varphi$.
\end{proof}
%-----------------end proof------------------

%:    Prin loc glob conc fact.lnl.loc
La notion d'\ali \lnl est une notion locale au sens suivant.

\rdb
\begin{plcc}
\label{fact.lnl.loc}\relax
\emph{(Applications \lins \lot sim\-ples)}
Soient $S_1$, $\ldots$, $S_n$ des \moco d'un anneau $\gA$.
Soit une \ali  $\varphi : \Ae{m}\rightarrow \gA^{q}$. Si les $\varphi _{S_i}:
\gA_{S_i}^m\rightarrow \gA_{S_i}^q$ sont \nlsz, alors $\varphi $ est \lnlz.
Plus \gnlt  $\varphi $ est \lnl \ssi les  $\varphi _{S_i}$ sont \lnlsz.
\end{plcc}
%--- end-plcc-----------------------------------------

%-----------------begin proof------------------
\begin{proof}
Voyons la deuxième affirmation. Montrer que $\varphi $ est \lnl revient
à trouver $\psi$ vérifiant $\varphi \,\psi \,\varphi =\varphi $. Ceci
est un \sli en les \coes de la matrice de $\psi$ et l'on peut donc appliquer
le
\plgc de base (\prirf{plcc.basic}).\iplg
\end{proof}
%-----------------end proof------------------

La terminologie d'\ali \lnl est justifiée par le  \plg
précédent et par la réciproque donnée au
point \emph{\ref{IFDg}}  du
\thrf{theoremIFD} (voir aussi le lemme de l'application \lnl
dans le cas des anneaux locaux, \paref{lelnllo}).

%:--- Subsection{Grassmanniennes}
\penalty-2500
\subsec{Grassmanniennes}
\label{subsecGrassmanniennes}

Le \tho suivant sert d'introduction aux \vrts grassmanniennes. Il résulte  du fait \ref{factInvGenCrois} et du \thref{propIGCram}.

%:     theorem{propFactDirRangk}
\begin{theorem}\label{propFactDirRangk}
\emph{(Sous-\mtfs en facteur direct dans un module libre)}
Soit $M=\gen{C_1,\ldots,C_m}$ un sous-\mtf de $\Ae{n}$ et~\hbox{$C=[\,C_1\;\cdots\;C_m\,]\in\Ae {n\times m}$}
la matrice correspondante.
\begin{enumerate}
\item \Propeq
\begin{enumerate}
\item La matrice $C$ est \lnlz.
\item Le module $M$ est en facteur direct dans $\Ae{n}$.
\item Le module $M$ est l'image d'une matrice $F\in\GAn(\gA)$.
\end{enumerate}
\item \Propeq
\begin{enumerate}
\item La matrice $C$ est de rang $k$.
\item \label{i1propFactDirRangk}Le module $M$ est l'image d'une matrice $F\in\GAn(\gA)$ de rang $k$.
%
%\item \label{i2propFactDirRangk}Le module $M$ est l'image d'une matrice $F\in\GAn(\gA)$ telle que
%$\det(\In+TF)=(1+T)^k$.
%
\end{enumerate}
\end{enumerate}
\end{theorem}

La \gui{variété} des droites vectorielles dans un \Kev de dimension~$n+1$
est, intuitivement, de dimension~$n$, car une droite dépend pour l'essentiel
de $n$ paramètres (un vecteur non nul, à une constante multiplicative près, cela fait~$(n+1)-1$
paramètres indépendants).
On appelle cette variété l'espace \pro de dimension $n$ sur $\gK$.
\\
Par ailleurs, en passant d'un corps $\gK$ à un anneau arbitraire
$\gA$, la bonne \gnn
d'une \gui{droite vectorielle dans $\gK^{n+1}$} est
\gui{l'image d'une \mprn de rang $1$
dans $\Ae {n+1}$}.
Ceci conduit aux \dfns suivantes.
%

%:     Definition{defiGrassmanniennes}
\begin{definition}\label{defiGrassmanniennes}~
\begin{enumerate}
\item On définit l'espace $\GAnk(\gA)\subseteq\GAn(\gA)$ comme l'ensemble
des  \mprns de rang $k$ et $\GGnk(\gA)$ comme l'ensemble
des sous-modules
de $\Ae n$ qui sont images de matrices de $\GAnk(\gA)$.
\item L'espace $\GG_{n+1,1}(\gA)$ est encore noté $\Pn(\gA)$ et on l'appelle
l'\emph{espace projectif de dimension $n$ sur $\gA$}.
\item On note $\GGn(\gA)$ l'espace de tous les sous-modules
en facteur direct dans~$\Ae n$ (i.e., images d'une \mprnz).
\end{enumerate}
\index{espace projectif de dimension $n$ sur un anneau}
\end{definition}

Naturellement la \dfn ci-dessus est peu satisfaisante, dans la mesure où on n'explique pas comment est structuré l'ensemble $\GGnk(\gA)$. Seule cette structure lui fait mériter son nom  d'\gui{espace}.

Une réponse partielle
est donnée par la constatation que $\GGnk$ est un foncteur.
Plus \prmtz,
à tout \homo $\varphi:\gA\to\gB$ on associe une application naturelle
$\GGnk(\varphi):\GGnk(\gA)\to\GGnk(\gB)$, avec notamment

\snic{\GGnk(\Id_\gA)=\Id_{\GGnk(\gA)}$,\, et \,$\GGnk(\psi\circ \varphi)=\GGnk(\psi)\circ \GGnk(\varphi),}

%\sni
lorsque $\psi\circ \varphi$ est défini.

%Nous reviendrons plus en détail par la suite sur la \gui{structure}
%des variétés grassmanniennes $\Pn(\gA)$, $\GGn(\gA)$ et
%  $\GGnk(\gA)$.

%:--- SUBSection{subsecCritInjSurj}---
\subsec{Critères d'injectivité et de surjectivité}
\label{subsecCritInjSurj}
%-----------------------------------------

Deux propositions célèbres sont contenues dans le \tho suivant.
%:   Theorem{prop inj surj det}------
\begin{theorem}
\label{prop inj surj det}\relax
Soit $\varphi:\Ae{n}\rightarrow \Ae{m}$ une \ali de matrice~$A$.
%---------begin item----------
\begin{enumerate}
\item L'application $\varphi$ est surjective \ssi $\varphi$ est de rang $m$,
\cad ici $\cD_m(\varphi)=\gen{1}$ (on dit alors que $A$ est \emph{\umdz}).
\item %:h2014 rajout du titre  et de l'index ci-dessous
 \emph{(\Tho de McCoy)}\index{McCoy!theor@\tho de ---}
L'application $\varphi$ est injective \ssi l'\id $\cD_n(\varphi)$ est fidèle, c.-à-d. si 
$\Ann_\gA\big(\cD_n(\varphi)\big)=0$.\index{unimodulaire!matrice ---}
\end{enumerate}
%---------end item----------
\end{theorem}
%--- end-theorem-----------------------------------------

%---------begin proof----------
\begin{proof}
\emph{1.} Si $\varphi$ est surjective, elle admet une inverse à droite $\psi$, et le fait \ref{fact.idd prod} donne
$\gen{1}=\cD_m(\I_m)\subseteq \cD_m(\varphi)\cD_m(\psi)$, donc $\cD_m(\varphi)=\gen{1}$.
Réciproquement, si $A$ est de rang $m$, l'\eqrf{eqIGCram3} montre que $A$
admet une inverse à droite, et
$\varphi$ est surjective.

\noindent  \emph{2.} Supposons que $\cD_n(A)$ est fidèle. D'après l'\egt
(\ref{eqGema4}), si $AV=0$, \hbox{alors $\mu_{\alpha,1..n}V=0$} pour tous les
\gtrs $\mu_{\alpha,1..n}$ de  $\cD_n(A)$, et donc $V=0$.
\\
Pour la réciproque\footnote{Voir aussi la \dem alternative donnée en \ref{propInjIdd}.}, nous montrons par récurrence sur $k$ la
\prt suivante: \emph{si $k$ vecteurs colonnes $\xk$ sont
\lint indé\-pendants, alors l'annulateur du vecteur $x_1\land\cdots\land
x_k$ est réduit à $0$.} Pour $k=1$ c'est trivial. Pour passer de $k$
à $k+1$ nous raisonnons comme suit. Soit $z$ un scalaire annulant
$x_1\land\cdots\land x_{k+1}$. Pour $\alpha\in \cP_{k,m} $, nous notons
$d_\alpha(\yk)$ le mineur extrait sur les lignes indices de
$\alpha$ pour les vecteurs colonnes~$\yk$ de~$\Ae{m}$. 
Puisque~$z (x_1\land\cdots\land x_{k+1})=0$, et vues les formules de Cramer, 
on a l'\egt

\snic{z \,\big(d_\alpha(\xk)x_{k+1}-
d_\alpha(x_1,\ldots, x_{k-1},x_{k+1} )x_k+ \cdots\big) = 0,
}

%\sni
donc $z \,d_\alpha (\xk)=0$. \\
Comme ceci est vrai pour tout
$\alpha$, cela donne~\hbox{$z (x_1\land\cdots\land x_{k})=0$}. Et par
l'\hdrz, $z=0$.
\end{proof}
%---------end proof----------

\rem Le \thrf{prop inj surj det} peut  se relire sous la forme suivante.
%-----------------begin item------------------
\begin{enumerate}
\item L'\ali $\varphi :\Ae{n}\to\Ae{m}$ est surjective \ssi
l'application  $\Al{m}\varphi
:\gA^{n\choose m}\to \gA$ est surjective.
\item L'\ali $\varphi :\Ae{n}\to\Ae{m}$ est injective \ssi  
l'application  $\Al{n}\varphi
:\gA  \to\gA^{m\choose n}$ est injective. \eoe
\end{enumerate}
%-----------------end item------------------

%--- Corollary{corprop inj surj det}-
\begin{corollary}
\label{corprop inj surj det}\relax
Soit une \Ali $\varphi:\Ae{n}\rightarrow \Ae{m}$.
%-----------------begin enum------------------
\begin{enumerate}
\item  Si $\varphi$ est surjective et  $n< m$, l'anneau est trivial.
\item  Si $\varphi$ est injective et  $n> m$, l'anneau est trivial.
\end{enumerate}
%-----------------end enum------------------
\end{corollary}
%--- end-corollary------------------------------------

\rem Une formulation plus positive, \eqvez, mais sans doute encore
plus dérou\-tante,
pour les résultats du corolaire précédent est la suivante.
\begin{enumerate}
\item Si $\varphi$ est surjective, alors $X^m$ divise $X^n$ dans $\gA[X]$.
\item Si $\varphi$ est injective, alors $X^n$ divise $X^m$ dans $\gA[X]$.
\end{enumerate}
D'une certaine manière, cela se rapproche plus de la formulation en \clamaz:
si l'anneau est non trivial, alors  $m\leq n$ dans le premier cas
(resp.   $n\leq m$ dans le deuxième cas).\\
L'avantage de nos formulations est qu'elles
fonctionnent dans tous les cas, sans avoir besoin de
présupposer que l'on sache décider si l'anneau est trivial ou pas. 
\eoe

%--- Corollary{corInjPuisExt}--------
\begin{corollary}
\label{corInjPuisExt}
Si $\varphi :\Ae{n}\rightarrow \Ae{m}$ est injective, il en va de même
pour toute puissance extérieure de $\varphi$.
\end{corollary}
%--- end-corollary------------------------------------
%-----------------begin proof------------------
\begin{proof}
L'annulateur de $\cD_n(\varphi)$ est réduit à $0$ par le \tho
précédent. Il existe un anneau $\gB\supseteq\gA$ tel que
les \gtrs de  $\cD_n(\varphi)$ deviennent \com dans $\gB$ (lemme
\ref{lemKerCom}). L'\Bli $\varphi_1 :\gB^n\rightarrow \gB^m$ obtenue en étendant $\varphi$ à $\gB$, est donc de
rang $n$ et admet un inverse à gauche $\psi$ (point~\emph{3} du \thref{propIGCram}), \cad
$\psi \circ\varphi_1 =\Id_{\gB^n}$. Par suite 

\snic{\Al{k} \psi \,\circ\,\Al{k}
\varphi_1 =\Id_{\Al{k} \!\gB^n}.}

%\sni
Ainsi la matrice de $\Al{k}\! \varphi_1$,  est 
injective.  Et puisque c'est la même matrice que celle de $\Al{k}\!\varphi$,
l'\ali  $\Al{k}\! \varphi$  est injective.
\end{proof}
%-----------------end proof------------------

%:--- Subsec{Caracterisation des applications \lnlsz}
\subsec{Caractérisation des applications \lnlsz} \label{CarLnl}
%-------------------

Le lemme suivant met en correspondance bijective
les \sfios et les suites d'\idms croissantes pour la \dvez.

%--- Lemme{lem ide-div}---
\begin{lemma}\label{lem ide-div}\relax
Soit une liste d'\idms $(e_{q+1} =0$, $
e_q$, $\ldots$, $e_{1}$, $e_{0}=1)$  telle que  $e_{i}$ divise $e_{i+1}$ pour
$i=0,\ldots ,q$. Alors, les \elts 
$r_i:=e_i-e_{i+1}$ pour $i\in\lrb{0..q}$,
forment un \sfioz. 
Réciproquement, tout \sfioz~\hbox{$(r_{0},\ldots,r_{q})$}
définit une telle liste d'\idms en posant 

\centerline{$e_j=\sum_{k\geq j}r_{k}$ pour $j\in\lrb{0..q+1}$.}
\end{lemma}
%--------------
%----begin proof ---------
\begin{proof} Il est clair que  $\som_ir_i=1$.
%D'après le lemme \ref{lem ide}, p
Pour
$0 \leq i <  q$, on a
$e_{i+1}=e_{i}e_{i+1}$.\\
D'où~$(e_i-e_{i+1})e_{i+1}=0$, \cad
$(r_q+\cdots+ r_{i+1})\cdot r_{i}=0$. On en déduit facilement que $r_ir_j=0$ pour $j>i$.
\end{proof}
%----end{proof----------

On note $\Diag(a_1,\ldots ,a_n)$ la matrice diagonale d'ordre $n$
dont le \coe en position $(i,i)$ est l'\elt $a_i$. \label{NOTADiag}

Dans le \tho qui suit certains des \idms $r_i$ dans le \sfio peuvent très
bien être nuls. Par exemple si l'anneau est connexe et non trivial
ils sont tous nuls
sauf un.

%:     Theorem{theoremIFD}--------
\begin{theorem}
\label{theoremIFD} \emph{(Matrice \lnlz)}%
\index{{localement!matrice --- simple}}%
\index{{matrice!localement simple}}  
Soit $G\in \Ae{m\times n}$ la matrice  d'une \ali 
$\varphi:\Ae{n}\rightarrow \Ae{m}$
\hbox{et  $q=\inf(m,n)$}.
\Propeq
%-----------------begin enum------------------
\begin{enumerate}\itemsep=1pt
\item \label{IFDa} %\label{e}\label{f}
 L'\ali $\varphi$ est \lnlz.

\item \label{IFDb} %\label{a}
Le sous-module  $\Im\varphi$  est   facteur direct dans  $\Ae{m}$.

\item \label{IFDc} %\label{c}
  $\Im\varphi$  est   facteur direct dans  $\Ae{m}$ et  $\Ker\varphi$
  est facteur direct dans  $\Ae{n}$.

\item \label{IFDd} %\label{d}
  Il existe une \ali  $\varphi\bul :\Ae{m}\rightarrow \Ae{n}$
  avec $\Ae{n}=\Ker\varphi\oplus\Im\varphi\bul $
  et  $\Ae{m}=\Ker\varphi\bul \oplus\Im\varphi$.

\item \label{IFDe} %\label{g}
  Chaque \idd $\cD_k(\varphi)$ est \idmz.

\item \label{IFDf1} Il existe un (unique) \sfio \hbox{$(r_0, r_1, \dots , r_q)$}
 tel que sur chaque localisé $\gA[1/r_k]$ 
 l'application $\varphi$ est de rang~$k$.

\item \label{IFDf}
  Chaque \idd $\cD_k(\varphi)$ est engendré par un \idm $e_k$.
  Soit  alors $r_k=e_k-e_{k+1}$. Les $r_k$ forment un \sfioz.
  Pour tout mineur $\mu$ d'ordre $k$ de $G$,
  sur le localisé $\gA[1/(r_k\,\mu)]$  l'\ali $\varphi$
  devient \nl de rang~$k$.

\item \label{IFDg}
  L'\ali $\varphi$ devient \nl après \lon en des
  \eco convenables.

\item \label{IFDi}
Chaque \idd $\cD_k(\varphi)$ est engendré par un \idm $e_k$ et  la
matrice de $\varphi$ devient \eqve à la matrice
$\Diag(e_1,e_2, \ldots ,e_q)$, éventuellement complétée par des
lignes ou colonnes nulles,  après \lon en des  \eco convenables.

\item $\!\!\!\sta$ \label{IFDh}
L'\ali $\varphi$ devient \nl après \lon en n'importe quel \id maximal.
\end{enumerate}
%-----------------end enum------------------
\end{theorem}
%------------------------------

%-----------------begin proof------------------
\begin{proof}
L'\eqvc  des points \emph{\ref{IFDa}}, \emph{\ref{IFDb}}, \emph{\ref{IFDc}}, \emph{\ref{IFDd}}
est déjà claire
(voir les faits \ref{factIng0} et~\ref{factInvGenCrois}).
Par ailleurs, on a trivialement \emph{\ref{IFDf}} $\Rightarrow$ \emph{\ref{IFDf1}}
$\Rightarrow$ \emph{\ref{IFDe}} et
\emph{\ref{IFDi}} $\Rightarrow$ \emph{\ref{IFDe}.}

\noindent  Puisque $q=\inf(m,n)$, on a  $\cD_{q+1}(\varphi)=0$.

\noindent  \emph{\ref{IFDa}} $\Rightarrow$ \emph{\ref{IFDe}.} 
On a $GHG=G$ pour une certaine
matrice $H$ et l'on applique le fait~\ref{fact.idd prod}.

\noindent  \emph{\ref{IFDe}} $\Rightarrow$ \emph{\ref{IFDf}.}  
Le fait que chaque
$\cD_k(\varphi)$ est engendré par un \idm $e_k$ résulte du fait
\ref{lem.ide.idem}.
Le fait que $(r_0,\ldots,r_q)$ est un \sfio résulte du lemme
\ref{lem ide-div} (et du fait \ref{fact.idd inc}). \\
Comme
$r_ke_{k+1}=0$,  sur l'anneau $\gA[1/r_k]$, et donc sur l'anneau $\gA[1/(\mu r_k)]$, où $\mu$ est un mineur d'ordre $k$, tous les mineurs d'ordre $k+1$
de la matrice~$G$ sont nuls.
Donc, par le lemme de la liberté, $G$  est
\nl de rang $k$.

\noindent  \emph{\ref{IFDf}} $\Rightarrow$ \emph{\ref{IFDi}.} Sur $\gA[1/r_k]$ et donc
sur $\gA[1/(\mu r_k)]$ ($\mu$ un mineur d'ordre $k$), 
on~a~$\Diag(e_1,\ldots ,e_q)=\Diag(1,\ldots ,1,0,\ldots ,0)$ avec $k$ fois 1.

\noindent  \emph{\ref{IFDf}} $\Rightarrow$ \emph{\ref{IFDg}.}
Notons $t_{k,j}$ les mineurs
d'ordre $k$ de $G$. Les \lons sont celles en les $t_{k,j}r_k$.
Nous devons vérifier qu'elles sont \comez.
Chaque  $e_k$ s'écrit sous
forme $\sum t_{k,j}v_{k,j}$,
donc $\sum_{k,j}v_{k,j}(t_{k,j}r_k)= \sum_ke_kr_k =\sum r_k=1$.

\noindent \emph{\ref{IFDg}} $\Rightarrow$ \emph{\ref{IFDa}.} Par application du
\plgrf{fact.lnl.loc} puisque toute application \nl est \lnlz.

\noindent \emph{\ref{IFDg}}   $\Rightarrow$ \emph{\ref{IFDh}.}
(En \clamaz.)
 Parce que le complémentaire
d'un \idema contient toujours au moins un \elt dans un \sys d'\eco (on
peut supposer l'anneau non trivial).

\noindent \emph{\ref{IFDh}} $\Rightarrow$ \emph{\ref{IFDg}.}
 (En \clamaz.) Pour chaque \idema
$\fm$ on obtient un $s_\fm\notin\fm$ et une matrice $H_\fm$ tels que l'on
ait $GH_\fm G=G$ dans~$\gA[1/s_\fm]$. L'\id engendré par les $s_\fm$
n'est contenu dans aucun \idema donc c'est l'\id $\gen{1}$. Il y a donc un
nombre fini de ces $s_\fm$ qui sont \comz.

\noindent Terminons en donnant une preuve directe pour
l'implication \emph{\ref{IFDf1}} $\Rightarrow$ \emph{\ref{IFDa}.}\\
 Sur l'anneau
$\gA[1/r_k]$ la matrice $G$ est de rang $k$ donc il existe une matrice~$B_k$
vérifiant $GB_kG=G$
(\thrf{propIGCram}). Ceci signifie sur l'anneau~$\gA$ que l'on a une 
matrice~$H_k$ dans $\Ae {n\times m}$  vérifiant $r_kH_k=H_k$ et $r_kG=GH_kG$. On
prend alors $H=\sum_k H_k$ et l'on obtient $G=GHG$.
\end{proof}
%-----------------end proof------------------

L'\eqvc  des points \emph{\ref{IFDa}} à \emph{\ref{IFDi}} a été établie
de manière \covz, tandis que le point \emph{\ref{IFDh}} implique les
précédents uniquement en \clamaz.

%:--- Subsection{Trace, norme, discriminant, transitivité}
\subsec{Trace, norme, discriminant, transitivité}
\label{subsecTransDet}

Nous notons $\Tr(\varphi)$ et $\rC{\varphi}(X)$ la trace et
le \emph{\polcarz} d'un \endo
$\varphi$ d'un module libre de rang fini
(nous prenons pour \polcar d'une matrice $F\in \Mn(\gA)$
le \pol $\det(X\I_n-F)$, qui a l'avantage d'être  \untz).
\index{polynome@\pol!caractéristique d'un endomorphisme}
\label{NOTA1Polcar}

%:    Notation{notaCTrN}
\begin{notation}\label{notaCTrN}~
\begin{enumerate}
\item [--] Si $\gA\subseteq\gB$ et si $\gB$ est un \Amo libre de rang fini,
on note $\dex{\gB:\gA}$ pour $\rg_\gA(\gB)$.
\item [--]
Pour $a\in\gB$ on note alors $\Tr_{\gB/\!\gA}(a)$,   $\rN_{\gB/\!\gA}(a)$
et  $\rC{\gB/\!\gA}(a)(X)$
la trace, le \deter et le \polcar de la multiplication par~$a$,
vue comme \endo du \Amo $\gB$.
\end{enumerate}
\index{polynome@\pol!caractéristique d'un élément}
\end{notation}

%%%%%%%%%%%%%%%%%%%%%%%%%%%%%%%%%%%%%%%%%
\begin {lemma}\label{lem1TransDet}
 Supposons que
  $\gA\subseteq\gB$ et que $\gB$ est un \Amo libre de rang fini $m$.
\begin{enumerate}
\item  Soit $E$ un
\Bmo libre de rang fini $n$.  Si $\ue = (e_i)_{i\in \lrbm}$ est une base
de $\gB$ sur $\gA$ et $\uf = (f_j)_{j\in \lrbn}$ une base de $E$ sur $\gB$, alors
$(e_if_j)_{i,j}$ est une base de $E$ sur $\gA$.
En conséquence, $E$ est libre
sur $\gA $ et

\snic{\rg_\gA( E) = \rg_\gB  (E) \times \rg_\gA (\gB) .}
\item Si $\gB\subseteq\gC$ et si $\gC$ est un \Bmo libre de rang fini,
on a

\snic{\dex{\gC:\gA}=\dex{\gC:\gB}\,\dex{\gB:\gA}.}
\end{enumerate}
\end {lemma}

%%%%%%%%%%%%%%%%%%%%%%%%%%%%%%%%%%%%%%%%%
\rem
 Soit $\gC = \aqo{\gA[Y]}{Y^3}=\gA[y]$, c'est une \Alg libre de rang~3. Puisque $y^4 = 0$, $\gB  = \gA
\oplus \gA y^2$ est une sous-\alg de $\gC$ libre sur $\gA$ dont le rang
(égal à 2)
ne divise pas le rang de $\gC$ (égal à 3). 
L'\egtz~\hbox{$\dex{\gC :  \gA} = \dex{\gC :  \gB} \dex{\gB  :\gA}$} ne s'applique pas car $\gC$ n'est pas libre sur $\gB$. 
\eoe

%%%%%%%%%%%%%%%%%%%%%%%%%%%%%%%%%%%%%%%%%

%:    theorem Th.transitivity
\begin{theorem}\label{Th.transitivity}\relax \emph{(Formules de transitivité pour la trace, le \deter et le \polcarz)}
Sous les mêmes hypothèses soit $u_\gB  :  E \to E$ une \Bliz. On
note $u_\gA$ cette application considérée comme une \Aliz. On
a alors les \egtsz:
\begin{enumerate}
\item [ ] $\det(u_\gA) = \rN_{\gB/\!\gA}\big(\det(u_\gB )\big) $,
 $\;\Tr(u_\gA) = \Tr_{\gB /\!\gA  } \big( \Tr(u_\gB )\big)$,
\item [ ] $\rC{u_\gA}(X) = \rN_{\gB [X]/\!{\gA[X]}} \big(
\rC{u_\gB}(X)\big).$
\end{enumerate}
\end {theorem}

%%%%%%%%%%%%%%%%%%%%%%%%%%%%%%%%%%%%%%%%%
\begin {proof}  On prend les notations du lemme \ref{lem1TransDet}.
Soient $u_{kj}$ les \elts de $\gB$ définis par $u(f_j) =
\sum_{k=1}^n u_{kj} f_k$.  Alors, la matrice $M$ de $u_\gA$ sur la base $(e_if_j)_{i,j}$ s'écrit comme une matrice par blocs

\snic{ M = \cmatrix {
M_{11} &  \cdots & M_{1n} \cr
\vdots &        &    \vdots \cr
M_{n1} &  \cdots & M_{nn} \cr
},}

%\sni
où $M_{kj}$ représente l'\Ali $b \mapsto b u_{kj}$ 
de $\gB $ dans $\gB$ sur la base~$\ue$.
Cela fournit la relation
sur la trace puisque:

\snic{ \arraycolsep2pt
\begin{array}{ccccccc}
 \Tr(u_\gA) & =  & \ds \som_{i=1}^n \Tr(M_{ii}) &=&\ds \som_{i=1}^n \Tr_{\gB /\gA}(u_{ii})\\[.6em]
  &  = & \ds \Tr_{\gB /\gA}  \big(\som_{i=1}^n u_{ii}\big) &=& \Tr_{\gB /\gA}\big(\Tr(u_\gB)\big).\
\end{array}
}

%\sni

Quant à l'\egt pour le \deterz, remarquons que les matrices $M_{ij}$
commutent deux à deux ($M_{ij}$ est la matrice de la multiplication par
$u_{ij}$). On peut donc appliquer le lemme \ref{lemdeterblocs} qui suit, ce qui nous donne:

\snic{ 
\det(M) = \det(\Delta) \qquad \hbox {avec} \qquad
\Delta = \som_{\sigma \in \Sn} \varepsilon(\sigma)
M_{1\sigma_1} M_{2\sigma_2} \ldots M_{n\sigma_n}.
}

%\sni
Or  $\Delta$ n'est autre que la matrice de la multiplication par l'\elt

\snic{\sum_{\sigma \in \Sn} \varepsilon(\sigma) u_{1\sigma_1} u_{2\sigma_2}
\ldots u_{n\sigma_n},
}

%\sni
i.e.,  par
$\det(u_\gB )$, donc:

\snic{ 
\det(u_\gA) = \det(M) = \rN_{\gB /\gA}\big(\det(u_\gB )\big).
}

%\sni
Enfin, l'\egt sur le \polcar se déduit
de celle sur les \deters en utilisant le fait que $\rC{u_\gA}(X)$ est le
\deter de l'\endo $X\Id_{E[X]} - u_\gA $  du $\gA[X]$-module $E[X]$ tandis
que $\rC{u_\gB}(X)$ est celui de la même application vue
comme \endo  du $\gB[X]$-module~$E[X]$.
%  $X\Id_E - u_\gB  :  E[X]
%\to E[X]$.
\end {proof}

Dans un anneau non commutatif, deux \elts $a$ et $b$ sont dits \emph{permutables} si $ab=ba$.

%%%%%%%%%%%%%%%%%%%%%%%%%%%%%%%%%%%%%%%%%
\begin {lemma}\label{lemdeterblocs}
Soit $(N_{ij})_{i,j}$ une famille de $n^2$
matrices carrées $\in\MM_m(\gA)$, deux à deux permutables,
et $N$ la matrice carrée d'ordre $mn$:
$$
N = \cmatrix {
N_{11} &  \cdots & N_{1n} \cr
\vdots &        &    \vdots \cr
N_{n1} &  \cdots & N_{nn} \cr
}.
$$
Alors:
$
\det(N) = \det\big( \som_{\sigma \in \Sn} \varepsilon(\sigma)
\,N_{1\sigma_1} N_{2\sigma_2} \ldots N_{n\sigma_n} \big).
$
\end {lemma}

\begin {proof}
  Notons $\Delta$ la matrice $n \times n$ définie par $\Delta = \sum_{\sigma
\in \Sn} \varepsilon(\sigma) N_{1\sigma_1} N_{2\sigma_2} \ldots
N_{n\sigma_n}$. Il faut donc démontrer que $\det(N) = \det(\Delta)$.

\noindent Traitons les cas particuliers
$n = 2$ puis $n = 3$. On remplace $\gA$ par $\gA[Y]$ et $N_{ii}$ par
$N_{ii} + Y
\rI_m$, ce qui a l'avantage de rendre certains \deters \ndzs
dans $\gA[Y]$.
Il suffit d'établir les \egts avec ces nouvelles matrices,
car on termine en faisant $Y=0$.

\noindent Le point-clef de la \dem pour $n = 2$ réside
dans l'\egt suivante:
$$
\cmatrix {N_{11} & N_{12} \cr N_{21} & N_{22} \cr}
\cmatrix {N_{22} & 0      \cr -N_{21} & \I_m \cr} =
\cmatrix {N_{11}N_{22} - N_{12}N_{21} & N_{12} \cr 0 & N_{22} \cr}.
$$
On considère ensuite le \deter des deux membres:
$$
\det(N) \det(N_{22}) = \det(N_{11}N_{22} - N_{12}N_{21}) \det(N_{22}),
$$
  puis on simplifie par $\det(N_{22})$ (qui est \ndzz) pour
obtenir le résultat.

\noindent Le cas $n = 3$ passe par l'\egtz:

\snac{
\cmatrix {
N_{11} & N_{12} & N_{13} \cr
N_{21} & N_{22} & N_{23} \cr
N_{31} & N_{32} & N_{33} \cr}
\cmatrix {
N_{22}N_{33} - N_{23}N_{32} & 0   & 0   \cr
N_{31}N_{23} - N_{21}N_{33} & \I_m & 0   \cr
N_{21}N_{32} - N_{22}N_{31} & 0   & \I_m \cr} =
\cmatrix {
\Delta & N_{12} & N_{13} \cr
0      & N_{22} & N_{23} \cr
0      & N_{32} & N_{33} \cr},
}  

%\sni
qui conduit à

\snic{\det(N) \det(N_{22}N_{33} - N_{23}N_{32}) = \det(\Delta) \det
\cmatrix {N_{22} & N_{23}\cr N_{32} & N_{33}\cr}.}

%\sni
Le cas $n = 2$ 
fournit $\det(N_{22}N_{33} - N_{23}N_{32}) = \det \cmatrix {N_{22} & N_{23}\cr
N_{32} & N_{33}\cr}$. On simplifie  par ce  \deter et on obtient $\det(N) = \det(\Delta)$.

\noindent   Le cas \gnl
  est laissé \alec (voir l'exercice~\ref{exolemdeterblocs}).
\end {proof}

\begin {corollary}
Soient $\gA  \subseteq \gB  \subseteq \gC$ trois anneaux
avec $\gC$ libre de rang fini sur $\gB $ et $\gB $ libre de rang fini
sur $\gA $.  On
a les \egts suivantes:
\begin{enumerate}
\item [ ] $\rN_{\gC/\!\gA} = \rN_{\gB /\!\gA} \circ  \rN_{\gC/\gB }  $,
 $\quad\Tr_{\gC/\!\gA} = \Tr_{\gB /\!\gA}   \circ \Tr_{\gC/\gB } $,
\item [ ] $\rC{\gC/\!\gA}(c)(X)= \rN_{\gB [X]/\!\gA[X]}\!  \big(\rC{\gC/\gB}(c)( X)\big)$ pour~$c \in \gC$.
\end{enumerate}
\end {corollary}

%%%%%%%%%%%%%
\penalty-2500
\subsubsection*{Déterminants de Gram et discriminants}

%:     Definition{defiGram}
\begin{definition}\label{defiGram}
Soit $M$ un \Amoz, $\varphi:M\times M\to\gA$ une forme bi\lin \smq et
$\ux=\xk$ une liste d'\elts de $M$. On appelle \emph{matrice de Gram
de  $(\xk)$ pour $\varphi$} la matrice
$$\Gram_\gA(\varphi,\ux)\eqdefi\big(\varphi(x_i,x_j)\big)_{i,j\in\lrbk}.$$
Son \deter est appelé le \emph{\deter de Gram de $(\xk)$ pour $\varphi$},
il est noté $\gram_\gA(\varphi,\ux)$.%
\index{matrice!de Gram}%
\index{determinant@\deterz!de Gram}%
\index{Gram!matrice de ---}%
\index{Gram!déterminant de ---}%
\end{definition}

Si $\gA y_1 + \cdots + \gA y_k \subseteq \gA x_1 + \cdots + \gA x_k$ on a
une \egt

\snic{
\gram(\varphi,\yk)=\det(A)^2\gram(\varphi,\xk),}

%\sni
 où $A$ est une matrice $k\times k$
qui exprime les $y_j$ en fonction des $x_i$.

\smallskip Nous introduisons maintenant un cas important de \deter de Gram,  le \discriz.
Rappelons que deux \elts $a$, 
$b$ d'un anneau $\gA$ sont dits \emph{associés}
s'il existe $u\in\Ati$ tels que $a=ub$.%
\index{associés!elements@\elts --- dans un anneau}

%:     Propdef{defiDiscTra}
\begin{propdef}\label{defiDiscTra}
Soit $\gC\supseteq\gA$ une \Alg qui est \hbox{un \Amo libre} de rang fini
et $x_1$, \dots, $x_k$, $y_1$, \dots, $y_k\in\gC$.
\begin{enumerate}
\item On appelle \ix{discriminant} de $(\xk)$
le \deter de la matrice\index{discriminant!d'une famille finie dans une algèbre libre de rang fini}

\snic{
\big(\Tr_{\gC/\!\gA}(x_ix_j)\big)_{i,j\in\lrbk}.
}

%\sni
 On le note $\disc_{\gC/\!\gA}(\xk)$ ou $\disc(\xk)$.
\item Si $\gA y_1 + \cdots + \gA y_k \subseteq \gA x_1 + \cdots + \gA x_k$ on a

\snic{
\disc(\yk)=\det(A)^2\disc(\xk),}

%\sni
 où $A$ est une matrice $k\times k$
qui exprime les $y_j$ en fonction des $x_i$.
\item En particulier, si $(\xn)$ et $(\yn)$ sont deux bases de $\gC$ comme \Amoz,
 les \elts
$\disc(\xn)$ et $\disc(\yn)$  sont congrus multiplicativement
modulo les carrés
de $\Ati$.
On appelle \emph{\discri de l'extension} $\gC\sur\gA$
la  classe d'\eqvc correspondante. On le note $\Disc_{\gC/\!\gA}$.%
\index{discriminant!d'une algèbre libre de rang fini}
\item Si $\Disc_{\gC/\!\gA}$ est \ndz et $n=\dex{\gC:\gA}$, un \sys $\un$ dans $\gC$
est une $\gA$-base de $\gC$ \ssi $\disc(\un)$ et $\Disc_{\gC/\!\gA}$
sont associés.

\end{enumerate}
\end{propdef}

Par exemple dans le cas où $\gA=\ZZ$ le \discri de  l'extension
est un entier bien défini,
tandis que si  $\gA=\QQ$, le \discri est \care
d'une part par son signe, d'autre part par la liste des
nombres premiers qui y figurent
avec une puissance impaire.

%:     Proposition{propTransDisc}
\begin{proposition}\label{propTransDisc}
Soient $\gB$ et $\gC$ deux \Algs libres de rangs $m$ et~$n$.
Soit l'\alg produit $\gB \times \gC $.
\'Etant données une liste~\hbox{$(\ux)=(\xm)$} d'\elts de
$\gB$ et une liste $(\uy)=(\yn)$ d'\elts de $\gC$, on a:
$$
\disc_{(\gB\times\gC)/\!\gA}(\ux,\uy) =
\disc_{\gB/\!\gA}(\ux) \times \disc_{\gC /\!\gA}(\uy).
$$
En particulier,
$\Disc_{(\gB\times \gC)/\!\gA}=\Disc_{\gB/\!\gA}\times\Disc_{\gC /\!\gA}.$
\end {proposition}
\facile

%:     proposition{prop1TransDisc}
\begin {proposition}\label{prop1TransDisc}
Soit $\gB\supseteq\gA$ une \Alg libre de rang fini $p$.  \\
On considère

\vspace{-.25em} 
\begin{itemize}\itemsep=0pt%\partopsep=0pt
\item un \Bmo $E$,
\item une forme $\gB$-bi\lin \smq
$\varphi_\gB: E \times E \to \gB$,
\item  une base $(\ub) = (b_i)_{i \in\lrbp}$  de $\gB$ sur $\gA$, et
\item   une famille $(\ue) = (e_j)_{j \in\lrbn}$ de $n$
\elts de $E$.
\end{itemize}
Notons $(\ub  \star  \ue)$ la famille $(b_i e_j)$ de
$np$ \elts de $E$ et $\varphi_\gA: E \times E \to \gA$
la forme $\gA$-bi\lin \smq définie par:
$$\preskip.3em \postskip.3em
\varphi_\gA(x, y) = \Tr_{\gB/\!\gA}\big(\varphi_\gB(x, y)\big).
$$
On a alors la formule de
transitivité suivante:
$$\preskip.3em \postskip.4em
\gram(\varphi_\gA, \ub \star \ue) =
\disc_{\gB/\!\gA}(\ub)^n \times \rN_{\gB/\!\gA} \big(\gram(\varphi_\gB, \ue)\big).
$$
\end {proposition}
\begin {proof}
Dans la suite les indices $i$, $i'$, $k$, $j$, $j'$ satisfont
à $i$, $i'$, $k\in\lrbp$ \hbox{et $j$, $j'\in\lrbn$}.  Convenons de
ranger $\ub \star \ue$ dans l'ordre

\snic{\ub \star \ue = 
b_1 e_1, \ldots, b_p e_1,
b_1 e_2, \ldots, b_p e_2, \ldots,
b_1 e_n, \ldots, b_p e_n .
}

%\sni
Pour $x \in \gB$, notons $\mu_x:  \gB \to \gB$ la multiplication par $x$ et $m(x)$
la matrice de $\mu_x$ dans la base $(b_i)_{i\in\lrbp}$ de $\gB$ sur $\gA$. On
définit ainsi un \isoz~$m$ de l'anneau $\gB$ vers un sous-anneau
commutatif de $\MM_p(\gA)$. Si l'on note $m_{ki}(x)$ les \coes de la matrice $m(x)$,
on a donc:
$$\preskip.4em \postskip.4em\ndsp 
\mu_x(b_i) = b_i x = \som_{k=1}^p m_{ki}(x) b_k, 
$$
avec $\rN_{\gB/\!\gA}(x) = \det\big(m(x)\big)$.  En posant
$\varphi_{jj'} = \varphi_\gB(e_j, e_{j'}) \in \gB$, on a

\snic{\varphi_\gA( b_i e_j
b_{i'} e_{j'} ) = \Tr_{\gB/\!\gA}\big(\varphi_\gB( b_i e_j b_{i'} e_{j'} )\big) = \Tr_{\gB/\!\gA} (b_i b_{i'} \varphi_{jj'} ).}

%\sni
  En utilisant l'\egt $b_{i'} \varphi_{jj'} =
\som_{k=1}^p m_{ki'}(\varphi_{jj'})\, b_k$, il vient avec 
$\Tr=\Tr_{\gB/\!\gA}$:
{\small
$$
\Tr ( b_i b_{i'} \varphi_{jj'} ) =
\Tr \big( \som_{k = 1}^p b_i\, m_{ki'}(\varphi_{jj'})\, b_k \big) =
\som_{k = 1}^p \Tr(b_i b_k)\, m_{ki'}(\varphi_{jj'}).\eqno(*)
$$
}

\vspace{-12pt}\noindent On définit $\beta \in \MM_p(\gA)$ par $\beta_{ik} = \Tr_{\gB/\!\gA}(b_i
b_k)$. La somme de droite dans~$(*)$ n'est autre que le \coe d'un produit de
matrices: $\big(\beta \cdot m(\varphi_{jj'})\big)_{ii'}$.  Le \deter de Gram de
 $\ub \star \ue$
pour $\varphi_\gA$  est donc une matrice $np \times np$
constituée de $n^2$ blocs de matrices $p \times p$. Voici cette matrice
en notant~\hbox{$\phi_{jj'}=m(\varphi_{jj'})$} pour alléger l'écriture:
{
\small
$$\!\!
\cmatrix {
\beta \phi_{11} &\beta  \phi_{12} & \ldots &\beta  \phi_{1n} \cr
\beta \phi_{21} &\beta  \phi_{22} & \ldots &\beta  \phi_{2n} \cr
      \vdots       &                    &        &       \vdots       \cr
\beta \phi_{n1} &\beta  \phi_{n2} & \ldots &\beta  \phi_{nn} \cr
} =
\cmatrix {
\beta  & 0     & \ldots & 0      \cr
0      & \beta & \ldots & \vdots \cr
\vdots &       & \ddots & 0      \cr
0      &       & \ldots & \beta  \cr
}
\cmatrix {
\phi_{11} & \phi_{12} & \ldots & \phi_{1n} \cr
\phi_{21} & \phi_{22} & \ldots & \phi_{2n} \cr
\vdots       &              &        & \vdots       \cr
\phi_{n1} & \phi_{n2} & \ldots & \phi_{nn} \cr
}\!.
$$
}

\vspace{-10pt}
\noindent 
En prenant les \deters on obtient
$$
\gram(\varphi_\gA, \ub \star \ue) = \det(\beta)^n \cdot
\det \cmatrix {
\phi_{11} & \phi_{12} & \ldots & \phi_{1n} \cr
\phi_{21} & \phi_{22} & \ldots & \phi_{2n} \cr
\vdots       &              &        & \vdots       \cr
\phi_{n1} & \phi_{n2} & \ldots & \phi_{nn} \cr
}.
$$

En utilisant le fait que les matrices $\phi_{jl}$ commutent deux à
deux, on obtient que le \deter de droite est égal à
$$\mathrigid1mu
\det\big( \sum_{\sigma \in \Sn} \varepsilon(\sigma)
\phi_{1\sigma_1} \phi_{2\sigma_2} \ldots
\phi_{n\sigma_n} \big) \,=\;
\det m\big(\det(\varphi_{jl})\big) \,=\; \rN_{\gB/\!\gA} \big(\gram(\varphi_\gB, \ue)\big),
$$
ce qui démontre le résultat.
\end {proof}
%

%:     theorem{thTransDisc}
\begin{theorem}\label{thTransDisc} \emph{(Formule de transitivité pour les \discrisz)}\\
 Soient $\gA \subseteq \gB \subseteq \gC$, avec $\gB$ libre  sur $\gA$,  $\gC$
 libre  sur $\gB$, \hbox{$\dex{\gC: \gB} = n$} \hbox{et $\dex{\gB: \gA} = m$}.
Soit $(\ub) = (b_i)_{i\in\lrbm}$ une base de $\gB$ sur~$\gA$,~\hbox{$(\uc) = (c_j)_{j\in\lrbn}$} une base
de~$\gC$ sur $\gB$ et notons $(\ub \star \uc)$ la base $(b_i c_j)$ de $\gC$ sur $\gA$.
Alors:
\[\arraycolsep2pt
\begin{array}{rcl}
\disc_{\gC/\!\gA}(\ub \star \uc)  &  = & \disc_{\gB/\!\gA}(\ub)^{\dex{\gC: \gB}}\ \rN_{\gB/\!\gA}\big(\disc_{\gC/\gB}(\uc)\big) , \\ [1.5mm]
\hbox{et donc}\quad \Disc_{\gC/\!\gA}  &  = & \Disc_{\gB/\!\gA}^{~~~\dex{\gC: \gB}}\
 \rN_{\gB/\!\gA}(\Disc_{\gC/\gB}).  \end{array}
\]
\end {theorem}

\begin {proof}
Application directe de la proposition~\ref{prop1TransDisc}.
\end {proof}

\penalty-5000
%--- Sec{Prin local-global de base modules}
\section{Principe local-global de base pour les modules}
\label{secPLGCBasicModules}
%-----------------------------------------

Les résultats de cette section ne seront pas utilisés avant le chapitre
\ref{chap ptf0}. 

Nous allons donner une version un peu plus \gnle du
\plg de base \vref{plcc.basic}, version qui concerne des \Amos et des
\alis arbitraires,
tandis que le principe de base peut être considéré comme le cas
particulier où les modules sont libres de rang fini.
La preuve est essentiellement la même que
celle du principe de base.\iplg
\\
Auparavant nous commençons par un bref rappel sur les suites exactes
et nous établissons quelques \prts \elrs de la \lon pour les modules.

%:--- SUBsection{Suites exactes}---------
\subsec{Complexes et suites exactes}
\index{suite exacte!d'\alisz}

Lorsque l'on a des \alis successives 

\snic{M\vers{\alpha}N\vers{\beta}P\vers{\gamma}Q\;,}

on dit qu'elles forment un \ix{complexe} si la composée de deux applications qui se suivent est nulle. 

On dit que la suite est \emph{exacte en $N$} si $\Im\alpha=\Ker\beta$. La suite toute entière est dite exacte si elle est exacte en $N$ et $P$. Ceci s'étend à des suites de longueur arbitraire.

%:HHH surjection scindée :  defi ramenée ici d'un peu plus loin
%       et developpée 
\rdb
Une \ali $\varphi:E\to F$ est appelée une \emph{surjection scindée}, s'il
 existe une \ali
$\psi:F\to E$ avec $\varphi\circ \psi=\Id_F$.
Dans ce cas on dit que $\psi$ est une \emph{section} de $\varphi$, et l'on a $E=\Ker\varphi\oplus \psi(F)\simeq \Ker\varphi\oplus F.$%
\index{scindee@scindée!surjection ---}\index{surjection scindée}%
\index{section!d'une surjection scindée}

Une \emph{suite exacte courte} est une suite exacte du type%
\index{suite exacte!courte}\index{suite exacte}
$$\preskip.2em \postskip.4em
{0\to M\vers{\alpha}N\vers{\beta}P\to 0}
$$
Dans ce cas, le module $M$ s'identifie à un sous-module $M'$ de $N$, et $P$ s'identifie~\hbox{à $N/M'$}.\\
Une suite exacte courte est dite \emph{scindée} si sa surjection est scindée. 
\index{scindee@scindée!suite exacte courte ---}%
\index{suite exacte!courte scindée}

Ce langage \gui{abstrait} a une contrepartie \imde en termes de \slis lorsque l'on a affaire à des modules libres de rang fini. Par exemple si $N=\Ae n$, $P=\Ae m$
et si l'on a une suite exacte   

\snic{0\to M\vers{\alpha}N\vers{\beta}P\vers{\gamma}Q\to 0\;,}

l'\ali $\beta$ est représentée par une matrice associée à un \sli
de $m$ \eqns à $n$ inconnues, le module $M$, isomorphe à $\Ker \beta$,
représente le défaut d'injectivité  de $\beta$ et le module $Q$, 
isomorphe à $\Coker \beta$,
représente son défaut de surjectivité.

\rdb
Il est important de noter qu'une suite exacte du type \label{sexaseco}
\[
\begin{array}{ccccccccccccccccccccccccccccccc} 
 0 & \to  & M_m&\vvers{u_m}& M_{m-1}&\lora& \cdots\cdots\cdots& \vvers{u_1}& M_0 &\to & 0\end{array}
\]
(avec $m\geq 3$) \gui{se décompose}
en $m-1$ \emph{suites exactes courtes} selon le schéma suivant.
\[ 
\begin{array}{ccccccccccccccccccccccccccccccc} 
 0 & \to  & E_2 & \vvers{\iota_2} & M_1 & \vvers{u_1} & M_0 &\to & 0  \\[1mm] 
 0 & \to  & E_3 & \vvers{\iota_3} & M_2 & \vvers{v_2} & E_2 &\to & 0\\[1mm] 
   & \vdots   &   &   &   &  &   &\vdots \\[1mm] 
 0 & \to  & E_{m-1} & \vvers{\iota_{m-1}} & M_{m-2} & \vvers{v_{m-2}} & E_{m-2} &\to & 0\\[1mm] 
 0 & \to   & M_m & \vvers{u_m} & M_{m-1} & \vvers{v_{m-1}} & E_{m-1} &\to & 0\\%[1mm] 
 \end{array}
\]
avec $E_i=\Im u_{i}\subseteq M_{i-1}$ pour $i\in\lrb{2..m-1}$, 
les $\iota_k$ des injections canoniques, et les $v_k$ obtenus à partir des $u_k$ en restreignant le module image à $\Im u_k$.

\smallskip  Un thème important de l'\alg commutative est fourni par les transformations qui conservent, ou ne conservent pas, les suites exactes.

Nous allons donner deux exemples de base, qui utilisent les modules d'\alisz.

Nous notons  $\Lin_\gA(M,P)$ le \Amo des \Alis de $M$ dans~$P$ et $\End_\gA(M)$
désigne  $\Lin_\gA(M,M)$ (avec sa structure d'anneau \gnlt non commutatif).
Le \emph{module dual} de $M$, $\Lin_\gA(M,\gA)$ sera en \gnl noté~$M\sta$.%
\label{NOTAAlis}\index{module!dual}

%:     Fact{fact1HomEx}
\begin{fact}\label{fact1HomEx}
Si $0\to M\vers{\alpha}N\vers{\beta}P$ est une suite exacte de \Amosz, et si~$F$ est un \Amoz, alors la suite 
$$\preskip.4em \postskip.0em 
0\to \Lin_\gA(F,M)\lora \Lin_\gA(F,N)\lora \Lin_\gA(F,P) 
$$
est exacte.
\end{fact}
%--------- fin fact ---------------------------------------------- 
%
\begin{proof} \emph{Exactitude en $\Lin_\gA(F,M)$.} 
Soit $\varphi\in\Lin_\gA(F,M)$ telle que $\alpha\circ \varphi=0$. 
Alors, puisque la première suite est exacte en $M$, pour tout $x\in F$, $\varphi(x)=0$, donc~$\varphi=0$.

\noindent \emph{Exactitude en $\Lin_\gA(F,N)$.}  
Soit $\varphi\in\Lin_\gA(F,N)$ telle
que $\beta\circ \varphi=0$.
Alors, puisque la première suite est exacte en $N$, pour tout $x\in F$, $\varphi(x)\in\Im\alpha$. Soient $\alpha_1:\Im\alpha\to M$ la bijection réciproque de
$\alpha$ (lorsque l'on  regarde~$\alpha$ comme à valeurs dans $\Im\alpha$) et
$\psi=\alpha_1\, \varphi$.
\\
 On obtient alors les \egts
%Alors,
$\Lin_\gA(F,\alpha)(\psi)=\alpha\, \alpha_1\, \varphi=\varphi$.
\end{proof}
%
%:     Fact{fact2HomEx}
\begin{fact}\label{fact2HomEx}
Si $N\vers{\beta}P\vers{\gamma}Q\to 0$ est une suite exacte de \Amos et si~$F$ est un \Amoz, alors la suite 

\snic{0\to \Lin_\gA(Q,F)\lora \Lin_\gA(P,F)\lora \Lin_\gA(N,F)}

%\sni
est exacte. 
\end{fact}
%--------- fin fact ---------------------------------------------- 
%
\begin{proof} \emph{Exactitude en $\Lin_\gA(Q,F)$.} Si $\varphi\in\Lin_\gA(Q,F)$ vérifie
$\varphi\circ \gamma=0$, alors, puisque~$\gamma$ est surjective, $\varphi=0$.

\noindent \emph{Exactitude en $\Lin_\gA(P,F)$.} Si $\varphi:P\to F$ vérifie
$\varphi\circ \beta=0$, alors $\Im\beta\subseteq\Ker\varphi$ et $\varphi$ se factorise par $P\sur{\,\Im\beta}\simeq Q$,
i.e. $\varphi=\psi\circ \gamma$ pour une \ali $\psi:Q\to F$, \cad $\varphi\in\Im\Lin_\gA(\gamma,F)$.
\end{proof}
%

%:     Fact{factDualReflexif}
\begin{fact} \label{factDualReflexif}
Soit $\beta:N\to P$ une \ali et $\gamma:P\to\Coker\beta$ la \prn canonique.  
\begin{enumerate}
\item L'application canonique $\tra \gamma:(\Coker\beta)\sta \to P\sta$ induit un \iso de $(\Coker\beta)\sta$ sur $\Ker\!\tra{\beta}$.
\item Si les \alis canoniques $N\to N^{\star\star}$ et $P\to P^{\star\star}$
sont des \isosz, alors la surjection canonique de $N\sta$ dans $\Coker\!\tra{\beta}$
 fournit par dualité un \iso de
 $(\Coker\! \tra{\beta})\sta $  sur $\Ker\beta$.
\end{enumerate}

\end{fact}

\begin{proof}
\emph{1.} On applique le fait \ref{fact2HomEx} avec $F=\gA$.

 \emph{2.} 
On applique le point \emph{1} à l'\ali $\tra\beta$
en identifiant  $N$ et~$N^{\star\star}$, ainsi que $P$ et~$P^{\star\star}$, 
et donc aussi $\beta$ et $^{\rm t}{(\!\tra{\beta})}$.
\end{proof}

\rem Il est possible d'affaiblir légèrement l'hypothèse en demandant
pour l'\ali $P\to P^{\star\star}$ qu'elle soit injective.
\eoe

%:--- SUBsection{subsecFaitsLoc}---------
\subsec{Localisation et suites exactes}

%--- Fact{fact.sexloc} ---------
\begin{fact}\label{fact.sexloc}\relax Soit $S$ un \mo d'un anneau $\gA$.
%-----------------begin item------------------
\begin{enumerate}
\item
 Si $M$ est un sous-module de $N$, on a l'identification
canonique de $M_S$ avec un sous-module de $N_S$ et de $(N/M)_S$ avec
$N_S/M_S$.\\
En particulier,  pour tout \id $\fa$ de $\gA$, le \Amo $\fa_S$ s'identifie
canoniquement avec l'\id $\fa\gA_S$ de $\gA_S$.
\item   Si $\varphi:M\rightarrow N$ est une \Aliz, alors:
%-----------------begin item------------------
\begin{enumerate}
\item    ${\rm  Im}(\varphi_S)$  s'identifie canoniquement à
$\big({\rm  Im}(\varphi)\big)_S$\,,
\item    ${\rm  Ker}(\varphi_S)$  s'identifie canoniquement à
$\big({\rm  Ker}(\varphi)\big)_S$\,,
\item    ${\rm  Coker}(\varphi_S)$  s'identifie canoniquement à
$\big({\rm  Coker}(\varphi)\big)_S$\,.
\end{enumerate}
%-----------------end item------------------
%
\item   Si l'on a une   suite exacte  de \Amos
$$\preskip.2em \postskip.1em
M\vers{\varphi}N\vers{\psi}P\;,
$$
alors la suite  de $\gA_S$-modules
$$\preskip.0em \postskip.0em
M_S\vers{\varphi_S}N_S\vers{\psi_S}P_S
$$
est \egmt exacte.
\end{enumerate}
%-----------------end item------------------
\end{fact}
%---------

%:     Fact{fact.LocIntersect}
\begin{fact}\label{fact.LocIntersect}\relax
Si $M_1$, $\ldots$, $M_r$ sont des sous-modules de $N$ et $M=\bigcap_{i=1}^rM_i$,
alors en identifiant les modules $(M_i)_S$ et $M_S$ à des sous-modules de $N_S$ on
obtient $M_S=\bigcap_{i=1}^r(M_i)_S$.
\end{fact}

%--- Fact{fact.transporteur}----
\begin{fact}
\label{fact.transporteur}\relax
Soient $M$ et $N$ deux sous-modules d'un \Amo $P$, \emph{avec $N$ \tfz}.
Alors, l'idéal transporteur $(M_S:N_S)$ s'identifie à   $(M:N)_S$, via
les applications naturelles de $(M:N)$ dans  $(M_S:N_S)$ et   $(M:N)_S$.
\end{fact}
%--- end-fact-----------------------------------------

Ceci s'applique en particulier pour l'annulateur d'un \itfz.

%:--- SUBsection{Principe \lgb pour les suites exactes}---------
\subsec{Principe \lgb pour les suites exactes de modules}

%--- Pr loc glog concret{plcc.basic.modules}---------
\begin{plcc}
\label{plcc.basic.modules}\emph{(Pour les suites exactes)}\\
Soient $S_1$, $\ldots$, $ S_n$  des \moco de $\gA$, $M$, $N$, $P$ des \Amos et deux
\alis $\varphi:M\to N$, $\psi:N\to P$. On note $\gA_i$ \hbox{pour $\gA_{S_i}$},
$M_i$ pour $M_{S_i}$ etc. \Propeq
\begin{enumerate}
\item  La suite
$M\vers{\varphi}N\vers{\psi}P$ est exacte.
\item  Pour chaque $ i\in\lrbn,$
la suite
$M_i\vers{\varphi_i}N_i\vers{\psi_i}P_i$ est exacte.
\end{enumerate}
Comme conséquence, $\varphi$ est injective (resp. surjective) \ssi pour
chaque $ i\in\lrbn,$  $\varphi_i$ est injective (resp. surjective)
\end{plcc}
%--- end-plcc-----------------------------------------
%-----------------begin proof------------------
\begin{proof}
Nous avons vu que \emph{1} $\Rightarrow$ \emph{2} dans le fait \ref{fact.sexloc}.\\
Supposons \emph{2}. Notons $\mu_i:M\to M_i,$  $\nu_i:N\to N_i,$  $\pi_i:P\to
P_i$ les \homos canoniques. Soit $x\in M$ et $z=\psi\big(\varphi(x)\big)$, on a
$$\preskip.4em \postskip.4em 
0=\psi_i\big(\varphi_i(\mu_i(x))\big)=\pi_i\big(\psi(\varphi(x))\big)=\pi_i(z), 
$$
donc pour un $s_i\in S_i$, $s_iz=0$ dans $P$. On conclut que $z=0$ en
utilisant la comaximalité des $S_i$:
$\som_i u_is_i=1$. Soit maintenant $y\in N$ tel \hbox{que $\psi(y)=0$}. Pour
chaque $i$ il existe un $x_i\in M_i$ tel que $\varphi_i(x_i)=\nu_i(y)$. \\
On
écrit $x_i=_{M_i}a_i/s_i$ avec $a_i\in M$ et $s_i\in S_i$. L'\egt
 $\varphi_i(x_i)=\nu_i(y)$ signifie que pour un certain $t_i\in S_i$
on a $t_i\varphi(a_i)=t_is_iy$ dans~$N$. Si $\som_i v_it_is_i=1$,
on en déduit que $\varphi\big(\som_iv_it_ia_i\big)=y$. Ainsi  $\Ker\psi$ est bien inclus dans  $\Im\varphi$.
\end{proof}
%-----------------end proof------------------

%--- Pr loc glog abstrait{plca.basic.modules}---------
\begin{plca}
\label{plca.basic.modules}\emph{(Pour les suites exactes)}\\
Soient $M$, $N$, $P$ des \Amosz, et deux
\alis $\varphi:M\to N$ \hbox{et $\psi:N\to P$}.  \Propeq
\begin{enumerate}
\item  La suite
$M\vers{\varphi}N\vers{\psi}P$ est exacte.
\item  Pour tout \idema $\fm$ la suite
$M_\fm\vers{\varphi_\fm}N_\fm\vers{\psi_\fm}P_\fm$ est exacte.
\end{enumerate}
Comme conséquence, $\varphi$ est injective (resp. surjective) \ssi pour
 tout \idema $\fm$,  $\varphi_\fm$ est injective (resp. surjective)
\end{plca}
%--- end-plcc-----------------------------------------
%
\begin{proof}
La \prt $x=0$ pour un \elt $x$ d'un module est une \prt \carfz.
De même pour la \prt $y\in\Im\varphi$. Ainsi, même si la \prt \gui{la suite est exacte} n'est pas \carfz, c'est une conjonction de \prts \carfz, et l'on peut
 appliquer le fait\etoz~\ref{fact2PropCarFin} pour déduire le \plg abstrait du \plg concret.
\end{proof}

Signalons enfin un \plgc pour les \mosz.

%--- Pr loc glog concret{plcc.basic.monoides}---------
\begin{plcc}
\label{plcc.basic.monoides}\emph{(Pour les \mosz)}\\
Soient $S_1$, $\ldots$, $ S_n$  des \moco de $\gA$, $V$ un \moz. \Propeq
\begin{enumerate}
\item  Le \mo $V$ contient $0$.
\item  Pour $ i\in\lrbn,$ le \mo $V$ vu dans $\gA_{S_i}$ contient $0$.
\end{enumerate}
\end{plcc}
%--- end-plcc-----------------------------------------
%
\begin{proof}
Pour chaque $i$ on a un $v_i\in V$ et un $s_i\in S_i$ tels que $s_iv_i=0$.
On pose~$v=\prod_i v_i\in V$. Alors, $v$ est nul dans les  $\gA_{S_i}$, donc dans
$\gA$.
\end{proof}
%

%%%%%%%%%%%%%%%%%%%%%%%%%%%%%%%%%%%%%%%%%%%%%%%%%%%%%%%%%%%%%%%%%%%%%%%%%%%
%%%%%%%%%%%%%%%%%%%%%%%%%%%%%%%%%%%%%%%%%%%%%%%%%%%%%%%%%%%%%%%%%%%%%%%%%%%
%%%%%%%%%%%%%%%%%%%%                                     %%%%%%%%%%%%%%%%%%
%%%%%%%%%%%%%%%%%%%%            EXERCICES                %%%%%%%%%%%%%%%%%%
%%%%%%%%%%%%%%%%%%%%                                     %%%%%%%%%%%%%%%%%%
%%%%%%%%%%%%%%%%%%%%%%%%%%%%%%%%%%%%%%%%%%%%%%%%%%%%%%%%%%%%%%%%%%%%%%%%%%%
%%%%%%%%%%%%%%%%%%%%%%%%%%%%%%%%%%%%%%%%%%%%%%%%%%%%%%%%%%%%%%%%%%%%%%%%%%%
%:section: Exercices
\Exercices{

%--- Exercise{exo2Lecteur}-------------
\begin{exercise}
\label{exo2Lecteur}
{\rm  Il est recommandé de faire les \dems non données, esquissées,
laissées \alecz,
etc\ldots
\, On pourra notamment traiter les cas suivants.
\begin{itemize}
\item \label{exofactKerAAsMMs} Vérifier les affirmations des faits \ref{factUnivLoc}
à \ref{fact.bilocal}. 
\item  \label{exocorplcc.basic}\relax
 Démontrer le corolaire \ref{corplcc.basic}.
\item  \label{exolemGaussJoyal}
Dans le lemme \ref{lemGaussJoyal} calculer des exposants
 convenables dans les points \emph{\ref{i2lemPrimitf}},
  \emph{\ref{i3lemPrimitf}},  \emph{\ref{i4lemPrimitf}}, 
  en explicitant complètement la \demz.
\item  \label{exocorpropCoh1}
  Démontrer le corolaire \ref{corpropCoh1}.
\label{exopropCoh4}
  Donner une preuve plus détaillée
du \thrf{propCoh4}.
\label{exoplcc.coh}\relax
  Vérifiez les détails dans la preuve du \plgrf{plcc.coh}.
\label{exopropCohfd1}
  Démontrer  la proposition \ref{propCohfd1}.
\item \label{exofact.transporteur}\relax 
Vérifier les affirmations des faits \ref{fact.sexloc}
à \ref{fact.transporteur}. Pour le
fait \ref{fact.LocIntersect} on utilisera la suite exacte
$0\to M\to N\to \bigoplus_{i=1}^rN/M_i$ qui est préservée par \lonz.

\end{itemize}
}
\end{exercise}
%--- end -exercise-----------------------------------------

%--- Exercise{exoNilpotentChap2}-------------
\begin{exercise}\label{exoNilpotentChap2} 
{\rm (Voir aussi l'exercice \ref{exoNilIndexInversiblePol})
\begin{enumerate}\itemsep0pt
\item  \emph{(Inversibles dans $\gB[T]$, cf. lemme \ref{lemGaussJoyal})}\\
Soient deux \pols $f= \sum_{i=0}^n a_i T^i$, $g= \sum_{j=0}^m b_j T^j$
 avec  $fg = 1$.
Montrer que les \coes $a_i$, $i \ge 1$, $b_j$, $j \ge 1$ sont nilpotents
et que $a_n^{m+1} = 0$.

\item  \emph{(\Pol \cara d'une matrice nilpotente)}\\
Soit $A \in \Mn(\gB)$ une matrice nilpotente et 
$\rC{A}(T) = T^{n} + \sum_{k=0}^{n-1}a_kT^k
%a_{n-1}T^{n-1} + \cdots + a_1T + a_0
$ son \polcarz.
\begin{enumerate}
\item Montrer que les
\coes $a_i$ sont nilpotents.
\item Précisément, si $A^e = 0$, alors  $\Tr(A)^{(e-1)n + 1} = 0$ et

\snic{%\dsp
 a_i^{e_i} = 0  \quad \hbox {avec} \quad
e_i = (e-1) {n \choose i} + 1\quad(i=0,\ldots,n-1).}
\end{enumerate}
\end{enumerate}
}
\end{exercise}
%--- end -exercise-----------------------------------------

%--- Exercise{exolemUMD}-------------
\begin{exercise}
 \label{exolemUMD}
 {\rm  On considère un vecteur
$x=(\xn)\in\Ae{n}$ et  $s\in\gA$.
\begin{enumerate}\itemsep0pt
\item Si $x$ est \umd  dans
$\aqo{\gA}{s}$ et dans $\gA[1/s]$, il est \umd dans~$\gA$.
\item Soient $\fb$ et $\fc$ deux \ids de $\gA$, si $x$ est \umd modulo
$\fb$ et modulo $\fc$, il l'est modulo $\fb\fc$.
\end{enumerate}
} \end{exercise}
%--- end-exercise-----------------------------------------

%--- Exercise{exoPlgb1}-----------
\begin{exercise} % \label{exo2.2.0}\relax
\label{exoPlgb1} 
(Une application typique du \plg de base)\\
{\rm   Soit
$x=(\xn)\in\Ae{n}$, \emph{\umdz}.
Pour $d\geq 1$, on note $\AXn_d$ le sous-\Amo des \pogs de degré $d$ et 

\snic{I_{d,x}=\sotq{f\in\AuX_d}{f(x)=0}, \hbox{ sous-\Amo de  }\AuX.\label{NOTAAXd}}
%-----------------begin item------------------
\begin{enumerate}\itemsep0pt
\item Si  $x_1\in\Ati$,  tout $f\in I_{d,x}$
est \coli \hbox{des $x_1X_j-x_jX_1$} avec pour \coes des \pogs
de degré~$d-1$. %On pourra commencer par le cas $x_1=1$.
\item En  \gnlz,  tout $f\in I_{d,x}$
est une \coli des $(x_kX_j-x_jX_k)$ avec pour \coes des \pogs
de degré~$d-1$.
\item Soit $I_x=\bigoplus_{d\geq 1}I_{d,x}$. Montrer
que $I_x=\sotq{F}{F(tx)=0}$ (où $t$ est une nouvelle
\idtrz). Montrer que $I_x$ est \emph{saturé}, i.e., si
$X_j^mF\in I_x$ pour un $m$ et pour chaque $j$, alors $F\in
I_x$.
\end{enumerate}
%-----------------end item------------------
}
\end{exercise}
%--- end-exercise-----------------------------------------

%--- Exercise{exoGaussJoyal}-------------
\begin{exercise}
\label{exoGaussJoyal} (Variations sur le \iJG lemme de Gauss-Joyal \ref{lemGaussJoyal})
%:HHH index \iJG
\\
{\rm  Montrer que les affirmations suivantes sont \eqves
(chacune des affirmations est \uvlez, i.e., valable pour tous \pols et
tout anneau commutatif $\gA$):
%-----------------begin item------------------
\begin{enumerate}\itemsep0pt
\item  $\rc(f)=\rc(g)=\gen{1}\; \Rightarrow\;  \rc(fg)=\gen{1}$,
\item  $(\exists i_0,j_0\;\;f_{i_0}=g_{j_0}=1)\;\Rightarrow\;\rc(fg)=\gen{1}$,
\item  $\exists p\in\NN, \; \; \big(\rc(f)\rc(g)\big)^p\subseteq \rc(fg)$,
\item  \emph{(Gauss-Joyal)} $\; \DA\big({\rc(f)\rc(g)}\big) = \DA\big({\rc(fg)}\big)$.
\end{enumerate}
%-----------------end item------------------
}
\end{exercise}
%--- end -exercise-----------------------------------------

%--- Exercise{exoNormPrimitivePol}-------------
\begin{exercise}\label{exoNormPrimitivePol}
{(Norme d'un \pol primitif via l'utilisation d'un anneau nul)}\\
{\rm  
Soient $\gB$ une \Alg libre de dimension finie, $\uX = (\Xn)$ des \idtrsz, $Q \in \gB[\uX]$ et $P = \rN_{\gB[\uX]/\gA[\uX]}(Q) \in \gA[\uX]$.  
Montrer que si $Q$ est primitif, alors $P$ l'est aussi.  
\emph{Indication:} vérifier que $\gA\cap \rc_\gB(P) = \rc_\gA(P)$, 
considérer le sous-anneau $\gA' = \gA/\rc_\gA(P)$ de $\gB' = \gB/\rc_\gB(P)$ et l'application $\gA'$-\lin \gui{multiplication par $Q$},
$m_Q : \gB'[\uX] \to \gB'[\uX],\;R\mapsto QR$.

}

\end {exercise}
%--- end -exercise-----------------------------------------

%--- Exercise{exoCohfd1}-------
\begin{exercise}
\label{exoCohfd1}
{\rm  Montrer qu'un anneau $\gA$ \coh est \fdi \ssi le test
\gui{$1\in\gen{a_1,\ldots ,a_n}$?}
est explicite pour toute suite finie $(a_1,\ldots ,a_n)$ dans $\gA$.
}
\end{exercise}
%--- end-exercise-----------------------------------------

%--- Exercise{exo.quo.coh}-------
\begin{exercise}
\label{exo.quo.coh} (Un exemple  d'anneau \noe \coh avec un quotient
\emph{non} \cohz)\relax \\
{\rm
On considère l'anneau $\ZZ$ et un \id $\fa$ engendré par une suite
infinie d'\eltsz, tous nuls sauf éventuellement un, qui est alors égal
à $3$
(par exemple on met un~$3$ la première fois, si cela arrive, qu'un
zéro de la fonction zéta
de Riemann n'a pas sa partie réelle égale à $1/2$).
Si l'on est capable de donner un \sys fini de \gtrs pour l'annulateur de $3$
dans $\ZZ/\fa$, on est capable de dire si la suite infinie
est identiquement nulle ou pas.
Ceci signifierait qu'il existe une méthode sûre
pour résoudre les conjectures du type de celle de Riemann.

\comm
Comme toute \dfn \cov raisonnable de la \noet semble réclamer qu'un
quotient
d'un anneau \noe reste \noez, et vu le \gui{contre-exemple} précédent,
on ne peut espérer avoir une preuve \cov du \tho
de \clama qui affirme que tout anneau \noe est \cohz.
\eoe
}
\end{exercise}
%--- end-exercise-----------------------------------------

%--- Exercise{exoIdempAX}-------------
\begin{exercise}\label{exoIdempAX} (Idempotents de $\AX$)\\
{\rm  
Montrer que tout \idm de $\AX$ est un \idm de $\gA$.
}
\end {exercise}
%--- end -exercise-----------------------------------------

%--- Exercise{exoIdmsSupInf}-------
\begin{exercise}
\label{exoIdmsSupInf}
{\rm  Soient $u$ et $v$ deux \idms et $x$
un \elt de $\gA$. \\
L'\elt $1-(1-u)(1-v)=u+v-uv$ est noté $u\vu v$.
\begin{enumerate}\itemsep0pt
% 1
\item Montrer que $x\in u\gA\,\Leftrightarrow\,ux=x$.
En particulier, $u\gA=v\gA \,\Leftrightarrow\, u=v$.
% 2
\item  L'\elt $uv$
est \und{le} plus petit commun multiple de $u$ et $v$ parmi les \idms de $\gA$
(i.e., si $w$ est un \idmz, $w\in u\gA\cap v\gA \,\Leftrightarrow\,
w\in uv\gA$). En fait, on a même $u\gA\cap v\gA=uv\gA$. On note $u\vi v=uv$.
% 3
\item  Démontrer l'\egt $u\gA+ v\gA=(u\vu v)\gA$.
En déduire que $u\vu v$ est \und{le} plus grand commun diviseur de $u$ et $v$ parmi les \idms de $\gA$
(en fait un \elt arbitraire de $\gA$ divise $u$ et $v$ \ssi il divise $u\vu v$).
\item Par une suite de manipulations \elrsz, transformer la matrice
$\Diag(u,v)$ en la matrice $\Diag(u\vu v,u\vi v)$.\\
 En déduire que les deux
\Amos $u\gA\oplus v\gA$ et   ${(u\vu v)}\gA\oplus {(u\vi v)}\gA$
sont isomorphes.
\item Montrer que les deux anneaux $\aqo{\gA}{u}\times \aqo{\gA}{v}$
et $\aqo{\gA}{u\vu v}\times \aqo{\gA}{u\vi v}$
sont isomorphes.
\end{enumerate}
}
\end{exercise}
%--- end-exercise-----------------------------------------

%--- Exercise{exoFracSfio}-------------
\begin{exercise}
\label{exoFracSfio}
{\rm Soit $\gA$ un anneau  et $(e_1,\ldots,e_n)$ un \sfio de $\Frac\gA=\gK$.
On écrit $e_i=a_i/d$ avec $a_i\in\gA$ et $d\in\Reg\gA$.
On a alors $a_ia_j=0$ pour $i\neq j$ et $\sum_ia_i$ \ndzz.
\\
\emph{1.} 
\'Etablir une réciproque.  
\\
\emph{2.} Montrer que $\gK[1/e_i]\simeq \Frac\big(\gA\sur{\Ann_\gA(a_i)}\!\big)$ et $\gK\simeq\prod_i\Frac\big(\gA\sur{\Ann_\gA(a_i)}\!\big)$. 
}
\end{exercise}
%--- end -exercise-----------------------------------------

%--- Exercise{exoFracZed}-------------
\begin{exercise}
\label{exoFracZed} (Séparer les composantes \irdsz)
\\
{\rm 
\emph{1.} Soit $\gA=\QQ[x,y,z]=\aqo{\QQ[X,Y,Z]}{XY,XZ,YZ}$ et $\gK=\Frac\gA$.
Quels sont les zéros de $\gA$ dans $\QQ^3$ (i.e.  $(x,y,z)\in\QQ^3$
tels que $xy=yz=zx=0$)?
Donner une forme réduite pour les \elts de $\gA$. Montrer que $x+y+z\in\Reg\gA$. Montrer que
les \elts $\frac{x}{x+y+z}$, $\frac{y}{x+y+z}$ et $\frac{z}{x+y+z}$ forment un \sfio  dans $\gK$. 
Montrer que $\gK\simeq\QQ (X) \times\QQ (Y) \times\QQ (Z) $.

\emph{2.} Soit $\gB=\QQ[u,v,w]=\aqo{\QQ[U,V,W]}{UVW}$ et $\gL=\Frac\gB$.\\
Quels sont les zéros de $\gB$ dans $\QQ^3$?
Donner une forme réduite pour les \elts de $\gB$. 
%Montrer que $uv+vw+wu$ est \ndzz. 
Montrer que $\gL\simeq\QQ(U,V)\times\QQ(V,W)\times\QQ(W,U)$.  
%\\
%NB: voir aussi les exercices \ref{}, \ldots
 
}
\end{exercise}
%--- end -exercise-----------------------------------------

%--- Exercise{exoIdempotentE2}-------
\begin {exercise} \label {exoIdempotentE2}
       (Idempotent et groupe \elrz)\\
{\rm
Soit $a \in \gA$ un \idmz. Pour $b \in \gA$, expliciter
une matrice $A \in \EE_2(\gA)$  et un \elt $d\in\gA$ tels que $A \cmatrix {a\cr b\cr} = \cmatrix {d\cr 0\cr}$. En particulier, $\gen {a,b} = \gen {d}$.\\
En outre, si $b$ est \ndz (resp.\ \ivz) modulo $a$, alors $d$ est \ndz (resp.\ \ivz). Enfin si $b$ est \idmz,
$d=a\vu b=a+b-ab$.
}

\end {exercise}

%--- Exercise{exoSfio}-----------
\begin{exercise}\label{exoSfio} %~\\
{\rm  Soit $(r_1,\ldots,r_m)$  une famille finie d'\idms
dans un anneau $\gA$. Posons $s_i=1-r_i$ et,
pour une partie $I$ de $\lrbm$, notons
$r_I=\prod_{i\in I}r_i\prod_{i\notin I}s_i$.

\noindent \emph{1.} 
Montrer que la matrice diagonale $D=\Diag(r_1,\ldots,r_m)$ est
semblable à une matrice  $D'=\Diag(e_1,\ldots,e_m)$ où les $e_i$ sont
des \idms qui vérifient: $e_i$ divise $e_j$ si $j> i$. On pourra
commencer par le cas $n=2$ et utiliser l'exercice
\ref{exoIdmsSupInf}. Montrer que  $\gen{e_k}=\cD_k(D)$ pour tout $k$.

\noindent \emph{2.} 
Montrer que l'on peut écrire $D'=PDP^{-1}$ avec $P$  une
\ixx{matrice}{de permutation
généralisée}, \cad une matrice 
qui s'écrit $\sum_jf_jP_j$ où les
$f_j$ forment un \sfio et chaque $P_j$ est une matrice de permutation.
Suggestions:
%-----------------begin item------------------
\begin{itemize}
\item  Les $r_I$ forment un \sfioz.
La matrice diagonale $r_ID$ a pour
 \coe en position $(i,i)$ l'\elt $r_I$ si $i\in I$ et $0$ sinon.
La matrice $P_I$ correspond alors à une
permutation ramenant les \coes $r_I$ en tête de la liste.
Enfin $P=\sum_Ir_IP_I$. Notez que le test \gui{$r_I=0$?} n'est pas \ncrz!
\item On peut aussi traiter le cas $m=2$: on trouve $P=e
\cmatrix{1&0\cr0&1}+ f\cmatrix{0&1\cr1&0}$ avec $f=r_2s_1$, $e=1-f$,
et $D'=\Diag(r_1\vu r_2,r_1\vi r_2)$.
\\
Ensuite on  traite
le cas  $m>2$ de proche en proche.
\end{itemize}
%-----------------end item------------------

}
\end{exercise}
%--- end-exercise-----------------------------------------

%--- Exercise{exoChinois}-------
\begin{exercise}
\label{exoChinois}
{\rm  Rappeler une preuve du \tho des restes chinois (\paref{restes
chinois}) et expliciter les \idmsz.
}
\end{exercise}
%--- end-exercise-----------------------------------------

%%%%%%%%%%%  exoFacileGrpElem1  %%%%%%%%%%%%%%
\begin {exercise} \label {exoFacileGrpElem1}
       (Groupe \elrz: premiers pas) {\rm Cas de $\MM_2(\gA)$.\\
 \emph{1.} 
Soit $a \in \gA$. 
%expliciter une suite de transformations
%\elrs pour $\cmatrix {a\cr 0\cr} \mapsto \cmatrix {0\cr a\cr}$
%et donner 
Déterminer une matrice $P \in \EE_2(\gA)$
%, produit de deux matrices \elrsz,
telle que $P\cmatrix {a\cr 0\cr} = \cmatrix {0\cr
a\cr}$. Même chose pour $\cmatrix {\varepsilon a\cr 0\cr}
\mapsto \cmatrix {a\cr 0\cr}$ où $\varepsilon \in \Ati$.

\noindent \emph{2.} 
\'Ecrire comme \elts de $\EE_2(\gA)$ les matrices $
\crmatrix {0 & -1\cr 1 & 0\cr}$  et $\crmatrix {-1 & 0\cr 0 & -1\cr} .
$

\noindent \emph{3.} 
Toute matrice triangulaire  de $\SL_2(\gA)$ est dans $\EE_2(\gA)$.

\noindent \emph{4.} 
Soient $u = \cmatrix {x\cr y}$, $v = \cmatrix {y\cr x}$, $w = \crmatrix {-y\cr
x}$ avec $x,y\in\gA$.  
Montrer \linebreak que
$v\in \GL_2(\gA)\cdot u$ et $w\in \EE_2(\gA)\cdot u$, mais pas \ncrt
$v \in \SL_2(\gA)\cdot u$.  Par exemple, si $x$, $y$ sont deux \idtrs sur un anneau
$\gk$, $\gA=\gk[x,y]$ et $v = Au$, avec $A \in\GL_2(\gA)$, alors
$\big(\det(A)\big)(0,0)=-1$. 
En conséquence on~a~$\det (A) \in -1 +\rD_\gk(0) \gen{x,y}$ 
(lemme \ref{lemGaussJoyal}), donc $\det(A) = -1$ si $\gk$
est réduit. De~plus, si $\det(A) =1$, alors $2=0$ dans~$\gk$.
Par suite, $v \in \SL_2(\gA)\cdot u$ \ssiz$2=0$ dans~$\gk$.
}
\end {exercise}

%%%%%%%%%%%%%%%%%%%%%%%%%%%%%%%%%%%%%%%%%%%%%%%%%%%%%%%%%%%%%%%%%%%%%%%%%%%%%%%

%%%%%%%%%%%  exoFacileGrpElem2  %%%%%%%%%%%%%%
%\penalty-2500
\begin {exercise} \label {exoFacileGrpElem2}
       (Groupe élémentaire : deuxièmes pas)\\
{\rm
\emph{1.} Soit $A \in \MM_{n,m}(\gA)$  avec
un \coe \iv et $(n,m)\neq(1,1)$. Déterminer des matrices $P \in \EE_n(\gA)$ \hbox{et $Q \in\EE_m(\gA)$} telles que
$
P A\, Q = \cmatrix {1 & 0_{1,m-1}\cr 0_{n-1,1}& A'}
$.
Exemple: avec $a\in\Ati$ donner $P$ pour
$P\,\cmatrix {a \cr 0\cr} = \cmatrix {1 \cr 0\cr}$
 (exercice~\ref{exoFacileGrpElem1} point~\emph{1}).

\noindent \emph{2.} 
Soit $A \in \MM_2(\gA)$ avec un \coe \ivz. Calculer
des matrices $P $ \hbox{et $Q \in \EE_2(\gA)$} telles que:
$
P A\, Q = \cmatrix {1 & 0\cr 0 & \delta\cr}
$  avec $\delta = \det(A)
$.\\
Toute matrice $A \in \SL_2(\gA)$ ayant un \coe
\iv appartient à $\EE_2(\gA)$. Expliciter  les cas suivants:
$$
\cmatrix {a & 0\cr 0 & a^{-1}\cr},\qquad
\cmatrix {0 & a\cr -a^{-1} & 0\cr},\qquad
\hbox {avec $a \in \Ati$.}
$$
\'Ecrire les matrices suivantes (avec $a \in
\Ati$) dans  $\EE_2(\gA)$:
$$
\cmatrix {a & b\cr 0 & a^{-1}\cr},\qquad
\cmatrix {a & 0\cr b & a^{-1}\cr},\qquad
\cmatrix {0 & a\cr -a^{-1} & b\cr},\qquad
\cmatrix {b & a\cr -a^{-1} & 0\cr}.
$$

\noindent \emph{3.} 
Si $A = \Diag(a_1, a_2, \ldots, a_n) \in \SLn(\gA)$,
alors $A\in\EE_n(\gA)$.

\noindent \emph{4.} 
Toute matrice triangulaire $A \in \SL_n(\gA)$
appartient à $\EE_n(\gA)$.
}
\end {exercise}

%%%%%%%%%%%%%%%%%%%%%%%%%%%%%%%%%%%%%%%%%%%%%%%%%%%%%%%%%%%%%%%%%%%%%%%%%%%%%%%
%%%%%%%%%%%%%%%%%%%%%%%%%%%%%%%%%%%%%%%%%%%%%%%%%%%%%%%%%%%%%%%%%%%%%%%%%%%%%%%

%%%%%%%%%%%  exoFacileGrpElem4  %%%%%%%%%%%%%%
\begin {exercise} \label {exoFacileGrpElem4}
                 (Les matrices de division $D_q$ de \deter $1$)\\
{\rm
\noindent Une \gui {division générale} $a = bq - r$ peut s'écrire
matriciellement:
$$
\crmatrix {0 & 1\cr -1& q\cr} \cmatrix {a\cr b\cr} = \cmatrix {b\cr r\cr}.
$$
Ceci conduit à introduire la matrice $D_q = \crmatrix {0 & 1\cr -1& q\cr} \in
\SL_2(\gA)$.

\noindent Montrer que $\EE_2(\gA)$ est le \mo engendré par les matrices
$D_q$.

}
\end {exercise}

%%%%%%%%%%%%%%%%%%%%%%%%%%%%%%%%%%%%%%%%%%%%%%%%%%%%%%%%%%%%%%%%%%%%%%%%%%%%%%%

%%%%%%%%%%%  exoFacileGrpElem5  %%%%%%%%%%%%%%
\begin {exercise} \label {exoFacileGrpElem5} %~
%
%\noindent 
{\rm
 Soient $\gA$ un anneau  et $A$, $B
\in \MMn(\gA)$. On suppose que l'on  a un~$i\in\gA$ avec
$i^2 = -1$ et que $2\in\Ati$. Montrer que les matrices de $\MM_{2n}(\gA)$

\snic{M = \cmatrix {A & -B\cr B & A\cr}\; $ et $\;
M' = \cmatrix {A+iB & 0\cr 0 & A-iB\cr}}

%\sni
sont \emph{\elrt semblables}, (i.e., $\Ex P \in
\EE_{2n}(\gA),\;P M P^{-1} = M'$).\\
\emph{Indication}: traiter d'abord
le cas $n=1$% et chercher $P$ comme produit de deux matrices \elrsz
.
}

\end {exercise}

%%%%%%%%%%%%%%%%%%%%%%%%%%%%%%%%%%%%%%%%%%%%%%%%%%%%%%%%%%%%%%%%%%%%%%%%%%%%%%%

%%%%%%%%%%%  exoFacileGrpElem6  %%%%%%%%%%%%%%
\begin {exercise} \label {exoFacileGrpElem6}
{\rm
\noindent
Pour $d \in \Ati$ et $\lambda \in \gA$  calculer la matrice

\snic{
\Diag(1, \dots, d, \dots, 1)\cdot \rE_{ij}(\lambda)\cdot
\Diag(1, \dots, d^{-1}, \dots, 1).}

%\sni

Montrer que le sous-groupe des matrices diagonales de $\GL_n(\gA)$
normalise~$\En(\gA)$.
}
\end {exercise}
%%%%%%%%%%%%%%%%%%%%%%%%%%%%%%%%%%%%%%%%%%%%%%%%%%%%%%%%%%%%%%%%%%%%%%%%%%%%%%%

%--- Exercise{exoTroisiemeLemmeLiberte}-------------
\begin{exercise}\label{exoTroisiemeLemmeLiberte} (Un lemme de liberté, ou un splitting off, au choix \dlecz) 
 \\
{\rm
Soit $F \in \GAn(\gA)$ un \prr possédant un mineur principal
d'ordre $k$ \ivz. Montrer que $F$ est semblable à une matrice
$\cmatrix {\I_k & 0\cr 0 & F'\cr}$ où $F' \in \GA_{n-k}(\gA)$.

\noindent
Le \mptf $P\eqdefi\Im F\subseteq \Ae n$ admet un facteur direct libre ayant pour base $k$ colonnes de $F$.
}
\end {exercise}
%--- end -exercise-----------------------------------------

%%%%%%%%%%%%%%%%%%%%%%%%%%%%%%%%%%%%%%%%%%%%%%%%%%%%%%%%%%%%%%%%%%%%%%%%%%%

%%%%%%%%%%%%%     exoABArang1       %%%%%%%%%%%%%%%%%%%
\begin {exercise}\label{exoABArang1}
{\rm
Soit $A \in \Ae {n \times m}$ de rang $1$.
 Construire $B \in \Ae {m \times n}$ telle que $ABA = A$
et vérifier que $AB$ est un \prr de rang $1$. Comparez votre
solution à celle qui résulterait de la preuve du \thrf{propIGCram}.
}

\end {exercise}

%%%%%%%%%%%%%%%%%%%%%%%%%%%%%%%%%%%%%%%%%%%%%%%%%%%%%%%
%--- Exercise{exoABAabstrait}-------------
\begin{exercise}
\label{exoABAabstrait}
{\rm
 \emph{Cet exercice constitue une abstraction des calculs qui ont mené
au  \thrf{propIGCram}.} 
On considère un \Amo $E$ \gui {ayant assez de formes \linsz},
i.e.  si $x \in E$ vérifie $\mu(x) = 0$ pour tout
$\mu\in E\sta$, alors $x = 0$. Ceci signifie que l'application
canonique de $E$ dans son bidual, $E \to E{\sta}{\sta}$, est injective. Cette
condition est vérifiée si $E$ est un module \emph{réflexif},
 i.e. $E \simeq
E{\sta}{\sta}$, e.g. un \mptfz%
, ou un module libre de rang fini%
.

\noindent Pour $x_1$, $\ldots$, $x_n \in E$, on note $\Vi_r(x_1, \ldots, x_n)$
l'\id de $\gA$ engendré par les évaluations de toutes les formes
$r$-\lins alternées de~$E$ en tous les $r$-uplets
d'\elts de $\{x_1, \ldots, x_n\}$.

\noindent \emph{On suppose que $1 \in
\Vi_r(x_1, \ldots, x_n)$ et $\Vi_{r+1}(x_1, \ldots, x_n) = 0$.}

\noindent On veut montrer
que le sous-module $\som \gA x_i$ est facteur direct dans~$E$ en explicitant un projecteur $\pi : E \to E$ dont l'image est 
ce sous-module.
\begin{enumerate}\itemsep0pt
\item \emph{(Formules de Cramer)}
Soit $f$ une forme $r$-\lin alternée sur~$E$.  Montrer, pour $y_0$, $
\ldots$, $y_r \in \som \gA x_i$, que
$$
\som_{i=0}^r
(-1)^i f(y_0, \ldots, y_{i-1}, \widehat {y_i}, y_{i+1}, \ldots, y_r)\, y_i = 0
.$$
Ou encore, pour $y$, $y_1$, $\ldots$, $y_r \in \som \gA x_i$
$$
f(y_1, \ldots, y_r)\,y =
\som_{i=1}^r f(y_1, \ldots, y_{i-1}, y, y_{i+1}, \ldots, y_r)\, y_i
.$$

\item
Donner $n$ formes \lins $\alpha_i \in E\sta$ telles que
l'\ali 

\snic{\pi : E \to E$, $\; x\mapsto \som_i \alpha_i(x) x_i}

%\sni
  soit un projecteur d'image
$\som \gA x_i$. \\
 On notera $\psi : \Ae{n} \to E$ définie par $e_i \mapsto x_i$,
et $\varphi : E \to \Ae{n}$ définie par $\varphi(x) = \big(\alpha_1(x), \ldots,
\alpha_n(x)\big)$.
On s'arrangera pour que $\pi = \psi \circ \varphi$ et $\pi \circ \psi = \psi$,
ce qui donne $\psi \circ \varphi \circ \psi = \psi$.

\item \emph{(Nouvelle \dem du \thrf{propIGCram})}
Soit $A \in \Ae {m \times n}$ une matrice de rang~$r$.
Montrer qu'il existe $B \in \Ae {n \times m}$ telle que~$A\,B\,A = A$.

\end{enumerate}

}
\end{exercise}
%--- end -exercise-----------------------------------------
%%%%%%%%%%%%%%%%%%%%%%%%%%%%%%%%%%%%%%%%%%%%%%%%%%%%%%%%%%%%%%%%%%%%%%%%%%%%%%%

%--- Exercise{exoPrepBinetCauchy}-------------
\begin{exercise}
\label{exoPrepBinetCauchy}
{\rm  
Soient $A\in \Ae {n\times m}$ et $B\in\Ae {m\times n}$.

\noindent  \emph{1.}    On a la formule de commutativité suivante:
$\det(\I_m+XBA)=\det(\I_n+XAB).$

\noindent   \emph{Première \demz}. Traiter d'abord le cas où $m=n$, par exemple par la méthode des \coes indéterminés.
Si $m\neq n$, on peut compléter $A$ et $B$ par des lignes et des colonnes de $0$
pour en faire des matrices carrées $A_1$ et $B_1$ de taille $q=\max(m,n)$ comme dans la \dem donnée \paref{eqIDC1}. On vérifie alors que $\det(\I_m+XBA)=\det(\I_q+XB_1A_1)$
et $\det(\I_n+XAB)=\det(\I_q+XA_1B_1)$.

\noindent  \emph{Deuxième \demz}. On considère une indéterminée $X$ et les matrices

\snic{B'=\cmatrix{XB&\I_m\cr\I_n&0_{n,m}}\quad \mathrm{et}\quad
A'=\cmatrix{A&\I_n\cr\I_m&-XB}.}

%\sni
Calculer $A'B'$ et $B'A'$ et conclure.

\noindent  \emph{2.} Qu'en déduit-on pour les \polcars de $A\,B$ et $B\,A$?
}
\end{exercise}
%--- end -exercise-----------------------------------------

%--- Exercise{exoBinetCauchy}-------
\begin{exercise}\label{exoBinetCauchy} (Formule de Binet-Cauchy)
\index{Binet-Cauchy!formule de --- }\\
 {\rm  On utilise les notations \paref{notaAdjalbe}.
    Si $A\in \Ae {n\times m}$ et $B\in\Ae {m\times n}$ sont deux matrices
    de formats transposés, on a la formule de Binet-Cauchy:
$$
\det(BA)=\som_{\alpha \in \cP_{m,n}
}\det(B_{1..m,\alpha})\det(A_{\alpha,1..m}).
$$

\noindent  \emph{Première \demz}.  On utilise la formule
$\det(\I_m+XBA)=\det(\I_n+XAB)$  (exercice \ref{exoPrepBinetCauchy}).
On considère alors le \coe de $X^m$ dans chacun des \pols
$\det(\I_m+XBA)$ et $\det(\I_n+XAB)$.

\noindent  \emph{Deuxième \demz}. Les matrices $A$ et $B$ représentent
des \alis $u :\Ae{m}\to\Ae{n}$ et $v :\Ae{n}\to\Ae{m}$. \\
On considère
alors les matrices de $\Al mu $, $\Al mv $ et $\Al m(v\circ u )$
sur les bases naturellement associées
aux bases canoniques \hbox{de $\Ae{n}$} et~$\Ae{m}$.
\\
On conclut en écrivant que $\Al m(v\circ u )=\Al m v \circ\Al  mu $.

\noindent  \emph{Troisième \demz}. Dans le produit $BA$
on intercale entre $B$ et $A$ une matrice
diagonale $D$ ayant pour \coes des \idtrs $\lambda_i$,
et l'on regarde quel est le \coe de $\lambda_{i_1}\cdots\lambda_{i_m}$
dans le \pol $\det(BDA)$ (pour cela on prend $\lambda_{i_1}=\cdots=\lambda_{i_m}=1$
et les autres nuls). On conclut en prenant tous les $\lambda_i$ égaux à $1$.
}
\end{exercise}
%--- end-exercise-----------------------------------------

%%%%%%%%%%%%%%%%%%%%%%%%%%%%%%%%%%%%%%%%%%%%%%%%%%%%%%%%%%%%%%%%%%%%%%%%%%%
%--- Exercise{exoDetExtPower}-------------
\begin{exercise}
\label {exoDetExtPower}
\noindent {\rm
Soit $u \in \End_\gA(\Ae{n})$. Pour $ k \in\lrb{0..n}$, on note $u_k =
\Al k(u)$. \\
Montrer que $\det(u_k) =
\det(u)^{n-1 \choose k-1}$ et que

\snic{\det(u_k)
\det(u_{n-k}) = \det(u)^{n \choose k}\ .}
}
\end{exercise}
%--- end -exercise-----------------------------------------
%%%%%%%%%%%%%%%%%%%%%%%%%%%%%%%%%%%%%%%%%%%%%%%%%%%%%%%%%%%%%%%%%%%%%%%%%%%%%%%

%--- Exercise{exoMatInjLocSimple}-------------
\begin{exercise}
 \label{exoMatInjLocSimple}
 {\rm  Pour $A\in \Ae {n\times r}$ \propeq
 \begin{enumerate}
\item La matrice $A$ est injective et \lnlz.
\item Il existe une matrice $B\in \Ae {r\times n}$ telle que $B\,A=\I_r$.
\item L'\idd $\cD_r(A)=\gen{1}$.
\end{enumerate}
Indication: voir les \thos \vref{propIGCram}, \vref{prop inj surj det} et \vref{theoremIFD}.
 } \end{exercise}
%--- end-exercise-----------------------------------------
%%%%%%%%%%%%%%%%%%%%%%%%%%%%%%%%%%%%%%%%%%%%%%%%%%%%%%%%%%%%%%%%%%%%%%%%%%%%%%%

%--- Exercise{exolemdeterblocs}--------
\begin{exercise}
\label{exolemdeterblocs}
{\rm  Traiter le cas \gnl dans la \dem du lemme \ref{lemdeterblocs}.
}
\end{exercise}
%--- end-exercise------------------
%%%%%%%%%%%%%%%%%%%%%%%%%%%%%%%%%%%%%%%%%%%%%%%%%%%%%%%%%%%%%%%%%%%%%%%%%%%%%%%

%--- Exercise{exoGram}-------------
\begin{exercise}
\label{exoGram}  ~\\
{\rm Si $\gram_\gA(\varphi,\xn)$ est \ivz, le sous-module
$\gA x_1+\cdots+\gA x_n$ est libre avec $(\xn)$ pour base.
}
\end{exercise}
%--- end -exercise-----------------------------------------

%%%%%%%%%%%%%%%%%%%%%%%%%%%%%%%%%%%%%%%%%%%%%%%%%%%%%%%%%%%%%%%%%%%%%%%%%%%
%:sinotenglish
\sinotenglish{
%:--- Exercise{exoABegal0}-------------
\begin{exercise} ~
\label{exoABegal0}
{\rm Soient $A\in\Ae{m\times n}$, $B\in\Ae{n\times p}$, et $r$, $s$ avec $r+s>n$.
\begin{enumerate}
\item Si  $AB=0$ 
alors $\cD_r(A)\cD_s(B)=0$.
\item En \gnlz, $\cD_r(A)\cD_s(B)\subseteq \cD_1(AB)$.
\item Plus \gnltz, si $r+s\geq n+q$,
alors pour tout mineur $\mu$ d'ordre $r$ de $A $ on a l'inclusion $\mu^q\cD_s(B)\subseteq \cD_{q}(AB)$.
\end{enumerate}
}
\end{exercise}
%--- end -exercise-----------------------------------------

%--- Exercise{exosexcNoether1}-------------
\begin{exercise}
\label{exosexcNoether1}  {\rm
 On considère un \Amo $M$ et deux sous-\Amos $N_1$ et $N_2$.
 On a une \sex courte:

\snic{0\lora  N_1\cap N_2 \vvers{j} N_1 \times  N_2 \vvers \pi  N_1+N_2 \lora 0}

avec $j(x)=(x,-x)$ et $\pi(y,z)=y+z$.
\begin{enumerate}
\item Qu'est-ce que cela donne en termes de dimensions d'\evcs 
lorsque~$\gA$ est un \cdi et $M$ un \evc de dimension finie?
\item \'Etudier la signification du caractère scindé de cette suite exacte.
\end{enumerate}
}
\end{exercise}
%--- end -exercise-----------------------------------------

%--- Exercise{exosexcNoether2}-------------
\begin{exercise}
\label{exosexcNoether2}  {\rm
 On considère un \Amo $M$ et deux sous-\Amos $N_1$ et $N_2$.
 On définit un complexe comme suit:

\snic{0\lora M/(N_1\cap N_2)\vvers{j}M/N_1 \times M/N_2 \vvers \pi M/(N_1+N_2)\lora 0}

avec $j(\wh x)=(\wi x,-\uci x)$ et $\pi(\wi y,\uci z)=\ov{y+z}$.
\begin{enumerate}
\item Montrer qu'il s'agit d'une suite exacte.
\item Qu'est-ce que cela donne en termes de dimensions d'\evcs 
lorsque~$\gA$ est un \cdi et $M$ un \evc de dimension finie?
\item Donner des exemples où cette suite exacte est scindée et d'autres
où elle ne l'est pas.
\end{enumerate}
}
\end{exercise}
%--- end -exercise-----------------------------------------

%--- Exercise{exoNoeSeco}-------------
\begin{exercise}
\label{exoNoeSeco} {\rm
 On considère  deux sous-modules $E$ et~$F'$ d'un \Amo $F$.
\\ On note $E'=E\cap F'$, $G=F/E $, $ G'=F'/E'$, $S=E+F'$,
    $E''=E/E'$,    $F''=F/F'$ et  $G''=F/S$.  
\begin{enumerate}
\item 
Montrer que l'on a un diagramme commutatif comme ci-dessous dans lequel 
\begin{itemize}
\item $\iota$, $\iota'$, $\iota_E$ et $\iota_F$ sont les injections canoniques,
\item  $\pi$, $\pi'$, $\pi_E$ et $\pi_F$ sont les surjections canoniques,
\item  et toutes les suites horizontales et verticales sont exactes. 
\end{itemize}

\smallskip 
{\small\centerline{$
\xymatrix@C=3.2em@R=1.5em {
&0\ar[d]  &  0\ar[d]  & 0\ar[d]
\\
0\ar[r]&E'\ar[d]_{\iota_E} \ar[r]^{\iota '} 
&F'\ar[d]_{\iota_F}\ar[r]^{\pi ' } &
G'\ar[d]^{\iota_G}\ar[r]&0
\\
0\ar[r]&E\ar[d]_{\pi_E} \ar[r]^{\iota} 
&F\ar[d]_{\pi_F}\ar[r]^{\pi} &
G\ar[d]^{\pi_G}\ar[r]&0
\\
0\ar[r]&E''\ar[d] \ar[r]_{\iota ''} &F''\ar[d]\ar[r]_{\pi'' } &G''\ar[d]\ar[r]&0
\\
&0  &  0  & 0
\\
}
$
}}

Faites le lien avec les \thos de Noether concernant les quotients de sous-modules.
\item Le diagramme construit est-il 
le seul diagramme {commutatif} satisfaisant les conditions 
requises au point \emph{1}?
\end{enumerate}
}
\end{exercise}

%--- Exercise{exoNoeSeco2}-------------
\begin{exercise}
\label{exoNoeSeco2}{\rm 
\begin{enumerate}
\item On  considère un diagramme commutatif comme ci-dessous dans lequel toutes les suites horizontales et verticales
sont supposées exactes

\smallskip 
{\small\centerline{$
\xymatrix@C=3.2em@R=1.5em {
&0\ar[d]  &  0\ar[d]  & 
\\
0\ar[r]&E_0\ar[d]_{\jmath_E} \ar[r]^{\iota _0} 
&F'\ar[d]_{\iota_F}\ar[r]^{\pi _0 } &
G_0\ar[r]&0
\\
0\ar[r]&E \ar[r]^{\iota} 
&F\ar[r]^{\pi} &
G\ar[r]&0
}
$
}}

\smallskip \`A \isos et renommages près, on peut supposer que $E$, $F'$ et~$E_0$
sont des sous-modules de $F$ et que toutes les injections et surjections sont canoniques (donc $G_0=F'/E_0$ et $G=F/E$). Nous le supposons désormais
et nous notons $E'=E\cap F'$. 
\begin{enumerate}
\item Montrer qu'il existe une unique \ali ${\jmath_G}:G_0\to G$ qui rend le diagramme commutatif (i.e., telle que $\jmath_G\circ \pi_0=\pi\circ \iota_F$).
\item Montrer que l'image de $\jmath_G$ est le sous-module $(E+F')/E=S/E$
de~$G=F/E$.   
\item Montrer que $\jmath_G$ est injective \ssi $E_0=E'$.
Dans ce cas~$\jmath_G$ réalise un \iso de~$F'/E'$ sur $S/E$.
Et l'on est ramené à la situation de l'exercice \ref{exoNoeSeco}.
\end{enumerate}

\item On étudie maintenant la situation \gui{duale} de celle du point \emph{1.} 
\\
\Prmtz, on suppose que l'on a un diagramme commutatif comme ci-dessous dans lequel
les suites horizontales et verticales sont exactes

\smallskip 
{\small\centerline{$
\xymatrix@C=3.2em@R=1.5em {
0\ar[r]&E%\ar[d]_{\pi_E} 
        \ar[r]^{\iota} &F\ar[d]_{\pi_F}\ar[r]^{\pi} & G\ar[d]^{\theta_G}\ar[r] & 0
\\
0\ar[r]&E_3%\ar[d] 
         \ar[r]_{\iota _3} &F''\ar[d]\ar[r]_{\pi_3 }    & G_3\ar[d]\ar[r] & 0
\\
&%0  
&  0  & 0
\\
}
$
}}

\smallskip \`A \iso et renommage près, on peut supposer que $E$
est un sous-module de $F$ et que l'injection $\iota$ est canonique.
Nous le supposons désormais et nous notons $F'=\Ker\pi_F$.
\\
Notons $S_3$ le noyau de l'\ali $\theta_G\circ \pi=\pi_3\circ \pi_F$.
On a donc une inclusion $S_3\supseteq \Ker\pi+\Ker\pi_F=E+F'$.
\begin{enumerate}
\item Montrer qu'il existe une unique \ali 
${\beta:E\to E_3}$
qui rend le diagramme commutatif (i.e., telle que $\iota_3\circ \beta=\pi_F\circ \iota$).
\item Montrer que $\Ker\beta=E\cap F'$.
\item  Montrer que $\beta$ est surjective \ssi $S_3=E+F'$. \\
Préciser dans ce cas en quoi on retrouve  la situation correspondant à l'exercice \ref{exoNoeSeco}.   
\end{enumerate}
 
%
%\item  
%Faites le lien entre les résultats précédents et les \thos de Noether concernant les quotients de sous-modules.
%
\end{enumerate}
}
\end{exercise}
%--- end -exercise-----------------------------------------

%--- Exercise{exoNoeSeco3}-------------
\begin{exercise}
\label{exoNoeSeco3}  {\rm
On suppose que dans le diagramme {commutatif} ci-dessous, les suites verticales sont exactes, et que la suite $0\to E\to F \to G \to 0$ est exacte.

\smallskip 
{\small\centerline{$
\xymatrix@C=3.2em@R=1.5em {
&0\ar[d]  &  0\ar[d]  & 0\ar[d]
\\
0\ar[r]&E_1\ar[d]_{\iota_E} \ar[r]^{\iota _1} 
&F_1\ar[d]_{\iota_F}\ar[r]^{\pi _1 } &
G_1\ar[d]^{\iota_G}\ar[r]&0
\\
0\ar[r]&E\ar[d]_{\pi_E} \ar[r]^{\iota} 
&F\ar[d]_{\pi_F}\ar[r]^{\pi} &
G\ar[d]^{\pi_G}\ar[r]&0
\\
0\ar[r]&E_2\ar[d] \ar[r]_{\iota _2} &F_2\ar[d]\ar[r]_{\pi_2 } &G_2\ar[d]\ar[r]&0
\\
&0  &  0  & 0
\\
}
$
}}
\begin{enumerate}
\item Montrer que la suite $0\to E_1\to F_1 \to G_1 \to 0$ est exacte \ssi la 
suite~\hbox{$0\to E_2\to F_2 \to G_2 \to 0$}  exacte. 
\item Dans ce cas, à renommages et \isos près, on retrouve le diagramme de 
l'exercice~\ref{exoNoeSeco}. En particulier, si les injections $\iota$, $\iota_1$, $\iota_E$  et~$\iota_F$
sont canoniques (ce qui n'est pas restrictif), on a $E_1= E\cap F_1$ \hbox{et
$\Ker(\pi_G\circ \pi)=E+F_1$} (donc $G_2\simeq F/(E+F_1)$).
\end{enumerate}
}
\end{exercise}
%--- end -exercise-----------------------------------------

%--- Exercise{exoDualite}-------------
\begin{exercise}
\label{exoDualite} (Dualité exacte et \endo cotransposé)\\
{\rm
 On considère un \abi $\Psi:E\times F\to G$ où $G$ est un \kmo libre de rang $1$. On dit que \emph{$\Psi$ est une dualité exacte entre $E$ et $F$} lorsque les \alis correspondantes 
 
\snic{E\lora \Lin(F,G),\,x\lmt\big(y\mt \varphi(x,y)\big)\;\hbox{ et }\;F\lora \Lin(E,G),\,y\lmt\big(x \mt \varphi(x,y)\big)}
 
sont des \isosz. On en déduit que $E\sta\simeq F$ et $F\sta\simeq E$.

\emph{1.} Lorsque l'on a une dualité exacte $\Psi$ entre $E$ et $F$, pour tout  $\varphi\in\End(F)$ on a une $\Psi$-transposée $\varphi^{\star_\Psi}:E\to E$ qui est l'unique \kli satisfaisant 

\snic{\Psi\big(\varphi^{\star_\Psi}(x),y\big)=\Psi\big(x,\varphi(y)\big)$ pour tous $x\in E$ et $y\in F.}

 On a comme pour la transposition usuelle
$\varphi_1^{\star_\Psi}\circ \varphi_2^{\star_\Psi}=(\varphi_2\circ \varphi_1)^{\star_\Psi}$. Et aussi, avec la \dfn \smq et une notation légèrement ambivalente $(\varphi^{\star_\Psi})^{\star_\Psi}=\varphi$.

\emph{2.} Soit $E$ un \kmo libre de rang $n$ et $k\in\lrb{1..n-1}$. 
\\
Montrer qu'une dualité exacte entre $\Al k E$ et $\Al {n-k} E$ est donnée par

\snic{\Psi_k:\Al k E\times \Al {n-k} E\lora \Al n E,\;\; (x,y)\lmt x\vi y.}

Montrer aussi que le cotransposé d'un \endo $\varphi\in\Lin_\gk(E)$
au sens usuel est égal à $\big(\Al {n-1}\varphi\big)^{\star_{\Psi_1}}$.
Ainsi est expliqué le fait que la matrice de $\wi\varphi$ sur une base donnée est
la transposée de la matrice des cofacteurs. Ceci donne aussi une \gui{bonne} raison pour laquelle l'\endo cotransposé est intrinsèque.}
\end{exercise}}%
% 
%: fin sinotenglish 

%%%%%%%%%%%%%%%%%%%%%%%%%%%%%%%%%%%%%%%%%%%%%%%%%%%%%%%%%%%%%%%%%%%%%%%%%%%
%%%%%%%%%%%%%%%%%%%%%%%%%%%%%%%%%%%%%%%%%%%%%%%%%%%%%%%%%%%%%%%%%%%%%%%%%%%
%:  problemes

%--- problem{exoRationaliteLineaire}-------------
\begin{problem}\label{exoRationaliteLineaire} {(Pivot de Gauss, $A\,B\,A=A$, et rationalité \linz)}\\
{\rm  
Soit $\gK$ un \cdiz. Si $x \in \gK^n$ est un vecteur non nul, son \emph {indice
pivot} $i$ est le plus petit indice $i$ tel que $x_i \ne 0$. On dit que le
\coe $x_i$ est le \emph {pivot} de~$x$.  La \emph {hauteur} $h(x)$ de $x$ est
l'entier $n-i+1$ et l'on convient que $h(0) = 0$. Par exemple, pour $n = 4$ et
$x = \cmatrix {0\cr 1\cr *\cr *\cr}$, l'indice pivot de $x$ est $i = 2$, et
$h(x) = 3$. Les notions d'échelonnement qui suivent sont relatives à cette
hauteur $h$.  \\
On dit qu'une matrice $A \in
\MM_{n,m}(\gK)$ est \emph {échelonnée en colonnes} si les colonnes non
nulles de $A$ ont des hauteurs distinctes; on dit qu'elle est \emph
{strictement échelonnée} si de plus les lignes passant par les indices
pivot sont des vecteurs de la base canonique de $\gK^m$ (ces vecteurs sont
\ncrt distincts). Voici une matrice strictement échelonnée
($0$ a été remplacé par un point):

\snic {
\cmatrix {
\cdot      &\cdot      &\cdot      &1      &\cdot      &\cdot\cr
\cdot      &\cdot      &\cdot      &a_{24} &\cdot      &\cdot\cr
\cdot      &\cdot      &1      &\cdot      &\cdot      &\cdot\cr
\cdot      &\cdot      &a_{43} &a_{44} &\cdot      &\cdot\cr
1      &\cdot      &\cdot      &\cdot      &\cdot      &\cdot\cr
\cdot      &1      &\cdot      &\cdot      &\cdot      &\cdot\cr
a_{71} &a_{72} &a_{73} &a_{74} &\cdot      &\cdot\cr
\cdot      &\cdot      &\cdot      &\cdot      &1      &\cdot\cr
a_{91} &a_{92} &a_{93} &a_{94} &a_{95} &\cdot\cr
}.}

\emph {1.}
Soit $A \in \MM_{n,m}(\gK)$ strictement échelonnée; on définit $\ov {A}
\in \MM_{n,m}(\gK)$ en annulant les \coes non pivots (les $a_{ij}$ dans
l'exemple ci-dessus) et $B = \tra {\,\ov A} \in \MM_{m,n}(\gK)$.  Vérifier
que $ABA = A$. \\
Décrire les \prrs $AB$, $BA$ et la \dcn $\gK^n = \Im AB
\oplus \Ker AB$.

\emph {2.}
Soit $A \in \MM_{n,m}(\gK)$ une matrice quelconque. Comment obtenir $Q \in
\GL_m(\gK)$ telle que $A' = AQ$ soit strictement échelonnée?
Comment calculer $B \in \MM_{m,n}(\gK)$ vérifiant $ABA = A$?

\emph {3.}
Soient $A \in \MM_{n,m}(\gK)$ et $y \in \gK^n$. On suppose que le
\sli $Ax = y$ admet une solution $x$ sur un sur-anneau
de $\gK$. Montrer qu'il admet une solution sur $\gK$.

\emph {4.}
Soient $\gK_0 \subseteq \gK$ un sous-corps et $E$, $F$ deux sous-\evcs 
 \suls de $\gK^n$. On suppose que $E$ et $F$
sont engendrés par des vecteurs à composantes dans
$\gK_0$. Montrer que $\gK_0^n = (E\cap \gK_0^n) \oplus (F\cap \gK_0^n)$.

Soit $E \subseteq \gK^n$ un sous-$\gK$-\evcz.  On dit que \emph{$E$ est
$\gK_0$-rationnel} s'il est engendré par des vecteurs à composantes dans $\gK_0$.

\emph {5.}
Soit $F$ un \sul de $E$ dans $\gK^n$ engendré par des
vecteurs de la base canonique de $\gK^n$: $\gK^n = E \oplus F$
et $\pi : \gK^n \twoheadrightarrow E$ la \prn associée.
\begin {itemize}
\item [\emph {a.}]
Montrer que $E$ est $\gK_0$-rationnel \ssi $\pi(e_j) \in \gK_0^n$ pour
tout vecteur $e_j$ de la base canonique.

\item [\emph {b.}]
En déduire l'existence d'un plus petit corps de rationalité pour $E$.

\item [\emph {c.}]
Quel est le corps de rationalité de l'image dans $\gK^n$ d'une
matrice strictement échelonnée en colonnes?
\end {itemize}

}

\end {problem}
%--- end -problem-----------------------------------------

%--- Problem{exoPlgb2}--------
\begin{problem}
\label{exoPlgb2} ~\\
{\rm  
\emph{1.} \emph{Algorithme de \fcn partielle.}
 \'{E}tant donnés deux entiers $a$ et
$b$ montrer que l'on peut calculer \gui{rapidement} une famille finie
d'entiers positifs $p_i$ premiers entre eux deux à deux tels que
$a=\pm\prod_{i=1}^np_i^{\alpha_i}$ et  $b=\pm\prod_{i=1}^np_i^{\beta_i}$.

\emph{2.} On considère un \sli $AX=B$ dans $\ZZ$ qui admet une infinité de
solutions dans $\QQ^m$.
Pour savoir s'il admet une solution dans $\ZZ^m$ on peut essayer une
méthode \lgbez. On commence par déterminer une solution dans
$\QQ$, qui est un vecteur $X\in\QQ^{m}$. On trouve un entier $d$ tel que
$dX\in\ZZ^{m}$,
de sorte que $X$ est à \coes dans $\ZZ[1/d]$.
Il suffit ensuite de construire une solution dans chaque localisé
$\ZZ_{1+p\ZZ}$ pour les $p$ premiers qui divisent $d$ et d'appliquer le
\plgc \ref{plcc.basic}. Pour savoir s'il y a une solution dans $\ZZ_{1+p\ZZ}$
et en construire une, on peut utiliser la méthode du pivot, à
condition
de prendre pour pivot un \elt de la matrice (ou plutôt de la partie
restant à traiter de la matrice) qui divise tous les autres \coesz, \cad
un \coe
dans lequel $p$ figure avec un exposant minimum.\\
L'inconvénient de cette méthode est qu'elle nécessite
de factoriser $d$, ce qui peut la rendre impraticable. \\
Cependant, on peut légèrement modifier la méthode de façon à
ne pas avoir à factoriser complètement $d$. On utilisera l'\algo de
factorisation partielle. On commence par faire comme si $d$ était un
nombre premier. Plus \prmt on travaille  avec
l'anneau $\ZZ_{1+d\ZZ}$.
On cherche si un \coe de la matrice est étranger à
$d$. Si l'on en trouve un, on le choisit comme pivot.
Dans le cas contraire aucun \coe de la matrice n'est étranger à $d$ et
(en utilisant si \ncr l'\algo de factorisation partielle) on est dans l'un des trois cas suivants:
%-----------------begin item------------------
\begin{itemize}
\item $d$ divise tous les \coes de la matrice, auquel cas, ou bien il
divise
aussi les \coes de $B$ et l'on est ramené à un \pb plus simple, ou bien
il ne divise pas un \coe de $B$ et le \sli n'admet pas de solution,
\item  $d$ s'écrit sous forme d'un produit de facteurs deux à deux
étrangers $d=d_1\cdots d_k$ avec $k\geq 2$, auquel cas on peut travailler
ensuite avec les \lons
en les \mos  $(1+d_1\ZZ)$, \ldots,   $(1+d_k\ZZ)$,
\item   $d$ s'écrit comme une puissance pure d'un $d'$ divisant $d$,
ce qui nous ramène, avec~$d'$ à la place de $d$, à un
\pb du même type mais plus simple.
\end{itemize}
%-----------------end item------------------
Vérifier que l'on peut exploiter récursivement l'idée
exprimée ci-dessus. \'Ecrire un \algo et l'expérimenter.
Examiner si l'\algo obtenu s'exécute en temps
raisonnable.
}
\end{problem}
%--- end-problem-----------------------------------------
%%%%%%%%%%%%%%%%%%%%%%%%%%%%%%%%%%%%%%%%%%%%%%%%%%%%%%%%%%%%%%%%%%%%%%%%%%%%%%%

}% fin des exos
%:  solutions
\penalty-2500
\sol{

%%%%%%%%%%%%%%%% solution d'un Fait laisse au lecteur, non retenue %%%%%%%%%%%%%
%\exer{exo2Lecteur}{\Demo du fait \ref{factMoco}. 
%Supposons que les $S_i$ ne sont pas \comz.
%Prenons dans chaque $S_i$ un $x_i$ de manière que
%$1\notin\gen{\xn}$. Soit $\fp$ un \idep
%(resp. un \idemaz) contenant $\gen{\xn}$: aucun des  $S_i$
%n'est contenu dans $\gA\setminus\fp$.}

%%%%%%%%%%%%%%%%%%%%%%%%%%%%%%%%%%%%%%%%%
\exer{exoNilpotentChap2}{
\emph{1.} On suppose \spdg $a_0 = b_0 = 1$.
Lorsque l'on écrit que $fg = 1$, il vient 

\snic{0 = a_n b_m$, $0 = a_n b_{m-1}
+ a_{n-1} b_m$, $0 = a_n b_{m-2} + a_{n-1} b_{m-1} + a_{n-2} b_{m},}

%\sni
et ainsi
de suite jusqu'au degré~$1$.\\
On montre alors par \recu sur~$j$ que $\deg(a_n^j g) \leq m-j$. 
\\
En particulier, pour $j = m+1$, on obtient $\deg(a_n^{m+1} g) \leq -1$, i.e.
$a_n^{m+1} g = 0$. D'où $a_n^{m+1} = 0$.
Enfin en raisonnant modulo $\DB(0)$, on obtient $a_j$ nilpotent
successivement pour $j=n-1$, $\ldots$, $1$.

\noindent
\emph{2a.} On considère les \pols sur l'anneau commutatif $\gB[A]$:

\snic{
f(T) = \det(\In - TA) \; \hbox{ et } \;
g(T) = \det(\In + TA + T^2A^2 + \cdots + T^{e-1}A^{e-1}).}

%\sni
On a $f(T) g(T) = \det(\In - T^eA^e)=1$.
Le \coe de degré~$n-i$ de $f$ est $\pm a_i$. On applique  \emph{1}.
}

\noindent
\emph{2b.} Il suffit de montrer que $\Tr(A)^{(e-1)n + 1} = 0$, car $a_i = \pm
\Tr\big(\bigwedge^{n-i}(A)\big)$.\\
%:HHH rajout correction fin 2.3b
 On considère le \deter défini par rapport à une base fixée
$\cB$ de $\Ae n$. Si l'on prend la base canonique formée par les $e_i$,
on a une \egt évidente 
$$
\Tr(f)= \det_\cB(f(e_1),e_2,\dots,e_n)+\cdots+\det_\cB\big(e_1,e_2,\dots,f(e_n)\big).
$$
Elle peut être vue sous la forme suivante:
$$
\Tr(f)\det_\cB(e_1,\dots,e_n)= \det_\cB(f(e_1),e_2,\dots,e_n)+\cdots+
\det_\cB\big(e_1,e_2,\dots,f(e_n)\big).
$$
Sous cette forme on peut remplacer les $e_i$ par n'importe quel système de $n$
vecteurs de $\Ae{n}$: les deux membres sont des formes $n$-\lins alternées
(en les $e_i$) sur $\Ae n$, donc sont égales parce qu'elles coïncident
sur une base.
\\
 Ainsi, multiplier un \deter par $\Tr(f)$ revient à
le remplacer par une somme de \deters dans lesquels on a fait opérer $f$ 
sur chacun des vecteurs.
\\
 On en déduit que  l'expression $\Tr(f)^{n(e-1)+1}\det_\cB(e_1,\dots,e_n)$ est égale à une somme 
dont chaque terme est un \deter de la forme 
$$
 \det_\cB\big(f^{m_1}(e_1),f^{m_2}(e_2),\dots,,f^{m_n}(e_n)\big),
$$
avec $\som_i m_i=n(e-1)+1 $,
donc au moins l'un des exposants $m_i$ est $\geq e$.

\noindent \rem Cette solution pour la borne $n(e-1)+1$ est due à Gert Almkvist. Voir à ce sujet:
{\sc Zeilberger D.} {\it Gert Almkvist's generalization of a mistake of Bourbaki.}
Contemporary Mathematics {\bf 143} (1993), p.~609--612. \eoe

%%%%%%%%%%%%%%%   exolemUMD    %%%%%%%%%%%%%%%%%%%%%%%%
\exer{exolemUMD}{\emph{1.} Posons $\fa=\gen{x_1,\ldots,x_n}$. On obtient $s^r\in\fa$ (pour un certain
$r$), \hbox{et $1-as\in\fa$} (pour un certain $a$). On écrit $1=a^rs^r+(1-as)(1+as+\cdots) \in \fa$.\\
\emph{2.}
$\fa+\fb=\gen{1}$, $\fa+\fc=\gen{1}$ et $(\fa+\fb)(\fa+\fc)\subseteq \fa+\fb\fc $,
donc
$\fa+\fb\fc=\gen{1}$.

}

%%%%%%%%%%%%%%   exoPlgb1    %%%%%%%%%%%%%%%%%%%%%%%%%%
\exer{exoPlgb1}{
\emph{1.}
Puisque $f$ est \hmgz, on a $f(tx) = 0$ pour une nouvelle
\idtr $t$. D'où des $U_i \in \gA[X_1,
\ldots, X_n, t]$ tels que $f = \sum_{i=1}^n (X_i - tx_i)U_i$.
\\ En
faisant $t := x_1^{-1}X_1$, on obtient des $v_i \in \gA[X_1, \ldots,
X_n]$ tels que 

\snic{f = \sum_{i=2}^n (x_1X_i - x_iX_1)v_i.}

%\sni
Enfin, puisque $f$ est
\hmg de degré $d$, on peut remplacer $v_i$ par sa composante
\hmg de degré $d-1$.
\\
\emph{2.}
Considérons l'\egt $f = \sum_{k,j}(x_k X_j - x_jX_k)u_{kj}$, où les
$u_{kj}$ sont des \pogs de degré $d-1$.  Il s'agit d'un
\sli en les \coes des~$u_{kj}$.  Puisque ce \sys
admet une solution sur chaque localisé~$\gA_{x_i}$ et que les~$x_i$ sont
\comz, il admet une solution sur~$\gA$.
\\
\emph{3.}
Si $F = \sum_d F_d$ est la \dcn de $F \in \gA[X_1, \ldots, X_n]$ en
composantes \hmgsz, on a $F(tx) = 0$ \ssi $F_d(x) = 0$ pour
tout $d$, d'où le premier point de la question.  Pour la saturation, montrons que si $X_iF \in I_x$ pour tout $i$, alors $F \in I_x$.
Or on a $x_iF(tx) = 0$ donc, par comaximalité des $x_i$, 
on obtient $F(tx) = 0$,
i.e. $F \in I_x$.
}
%%%%%%%%%%%%%%%%%%%%%%%%%%%%%%%%%%%%%%%%%%%%%%%%%%%%%%%%%%%%%%%%%%%%%%%%%%%

\exer{exoNormPrimitivePol} 
Le \pol $Q$, vu comme un \pol à \coes dans $\gB'$, reste primitif
donc \ndz (Gauss-Joyal, point \iref{i2lemPrimitf} du lemme
\ref{lemGaussJoyal}). 
Puisque $m_Q$ est
injective, son \deter $\det(m_Q) = P \in \gA'[\uX]$ est \ndz (\thref {prop inj surj
  det}, point~{\emph 2}).  Mais $P$ est \egmt nul dans $\gA'[\uX]$. Donc  $\gA'$ est l'anneau nul, autrement dit $1 \in \rc_\gA(P)$.

%%%%%%%%%%%%%%%%%%%%%%%%%%%%%%%%%%%%%%%%%%%%%%%%%%%%%%%%%%%%%%%%%%%%%%%%%%%
\exer{exoIdempAX}
Soit $f(X)$ un \idm de $\AX$. Il est clair que $e=f(0)$ est \idmz.
On veut montrer que $f=e$. Pour cela on peut raisonner séparément modulo $e$ et modulo $1-e$.
\\
  Si $e=0$, alors $f=Xg$. On a $(Xg)(1-Xg)=0$, or $1-Xg$ est \ndzz, donc~$g=0$.
\\
 Si $e=1$, on considère l'\idm $1-f$ et l'on est ramené au cas précédent.

%%%%%%%%%%%    exoIdmsSupInf    %%%%%%%%%%%%%%%%%%%%%%%%%%%%%
\exer{exoIdmsSupInf}
Pour la question \emph{5} on montre d'abord le résultat lorsque $uv=0$.
Dans la situation \gnlez, on note $u'=1-u$ et $v'=1-v$. On a alors un \sfio
$(uv,uv',u'v,u'v')$ et en appliquant le cas particulier précédent on voit que les deux anneaux sont isomorphes au \linebreak 
produit $\aqo{\gA}{uv'}\times (\aqo{\gA}{uv})^2\times \aqo{\gA}{u'v}$.

%%%%%%%%%%%    exoFracSfio    %%%%%%%%%%%%%%%%%%%%%%%%%%%%%

\exer{exoFracSfio}\emph{2.} On a $\gK[1/e_i]\simeq \gK\sur{\Ann_\gK(e_i)}$
et $\Ann_\gK(e_i)=\Ann_\gA(a_i)\gK$. Pour un \elt $x$  
de~$\gA$, on écrit  $dx=\sum_{i\in\lrbn}x_i$ dans $\gK$,
avec $x_i=e_idx=a_ix$. La \dcn est donc entièrement dans~$\gA$. 
Et $dx\equiv x_i\mod \Ann_\gA(a_i)$, donc la composante $\gK\sur{\Ann_\gK(e_i)}$ du produit, quand on la voit comme l'\id $e_i\gK$, est formée des 
\elts de la forme~$a_ix/y$ avec $x\in\gA$ et~$y$ \ndz dans~$\gA$. 
Mais~$y$ est \ndz dans~$\gA$ \ssi chaque
$y_i=a_iy$ est \ndz modulo $\Ann_\gA(a_i)$, 
de sorte $\gK\sur{\Ann_\gK(e_i)}$ s'identifie à  $\Frac(\gA\sur{\Ann_\gA(a_i)}\!)$.

%%%%%%%%%%%    exoFracZed    %%%%%%%%%%%%%%%%%%%%%%%%%%%%%

\exer{exoFracZed}
\emph{1.} Les zéros de $\gA$ sont les trois \gui{axes de coordonnées}. 
\\
Tout \elt de $\gA$ s'écrit de manière unique
sous forme 

\snic{u=a+xf(x)+yg(y)+zh(z),}

%\sni
avec $f$, $g$, $h\in\QQ[T]$. Ceci implique  que $x+y+z$ est \ndz
car 

\snic{(x+y+z)u=x\big(a+xf(x)\big)+y\big(a+yg(y)\big)+z\big(a+zh(z)\big).}

%\sni
Donc les \elts $\frac{x}{x+y+z}$, $\frac{y}{x+y+z}$ et $\frac{z}{x+y+z}$ forment un \sfio de $\gK$. 
On conclut avec l'exercice \ref{exoFracSfio} en notant que $\Ann_\gA(x)=\gen{y,z}$, 
et donc que

\snic{\gA\sur{\Ann_\gA(x)}\simeq \QQX.}

\sni
\emph{2.} Les zéros de $\gB$ sont les trois \gui{plans de coordonnées}. Le \sfio dans $\gL$ est  donné par
 $\frac{uv}{uv+vw+wu}$, $\frac{vw}{uv+vw+wu}$ et $\frac{wu}{uv+vw+wu}$.

%%%%%%%%%%%   exoIdempotentE2    %%%%%%%%%%%%%%%%%%%%%%%%%%%%%%
\vspace{2pt}
\exer{exoIdempotentE2}{
Il suffit de résoudre la question modulo $a$ et modulo $1-a$.\\
Modulo $a$: $\cmatrix{a\cr b}=\cmatrix{0\cr b}\mapsto \cmatrix{b\cr b}\mapsto \cmatrix{b\cr 0}$.\\          
Modulo $1-a$, $\cmatrix{a\cr b}=\cmatrix{1\cr b}\mapsto  \cmatrix{1\cr 0}$.
En recollant: $d=(1-a)b+a$ avec par exemple la matrice
$A=A_2A_1,$ où
$$\preskip.4em \postskip.4em  A_1=(1-a)\bloc{1}{1}{0}{1}+ a \bloc{1}{0}{0}{1}=
\bloc{1}{1-a}{0}{1} , $$
$$\preskip.0em \postskip.4em
A_2=(1-a)\bloc{\phantom-1}{0}{-1}{1}+ a \bloc{\phantom-1}{0}{-b}{1}=\bloc{1}{0}{a-ab-1}{1} \et $$
$$\preskip.0em \postskip.4em
A=\bloc{1}{1-a}{a-ab-1}{a}.
$$
}

%%%%%%%%%%%%%%%%%%%%%%%%%%%%%%%%%%%%%%%%%

%: \exer{exoSfio}
\exer{exoSfio}
\emph{1} et \emph{2.} Cas $m=2$. On a par manipulations \elrs avec 

\snic{e_1=r_1\vu r_2=r_1+s_1r_2=r_2+s_2r_1,}

\snii
en notant que $e_1r_2=r_2$ et $-r_2(r_2s_1)=-r_2s_1=r_1r_2-r_2$.
\[\!\!\! \!\!
\begin{array}{rcl} 
\cmatrix{r_1&0\cr 0 & r_2}  &  \longmapsto & \cmatrix{r_1&0\cr r_2 & r_2} \;  \longmapsto \;  \cmatrix{r_1+r_2s_1&r_2s_1\cr r_2 & r_2}=\cmatrix{e_1&r_2s_1\cr r_2 & r_2}  \\[4mm] 
&  \longmapsto &\cmatrix{e_1&0\cr 0 & r_1r_2}.  
 \end{array}
\]
En outre en posant $f=r_2s_1$, $e=1-f$ et $P=\cmatrix{e&f\cr f & e}$
on a $P^2=\I_2$, $er_1=r_1$, $er_2=r_1r_2$ et
\[ 
\begin{array}{rcl} 
P\,  \cmatrix{r_1&0\cr 0 & r_2}\,P  & =    & \cmatrix{er_1&fr_2\cr fr_1 & er_2}\,P    \;  =  \;
\cmatrix{r_1&f\cr 0 & r_1r_2}\,\cmatrix{e&f\cr f & e}  \;  = \\[5mm]  
\cmatrix{r_1e+ f&0\cr 0 &r_1r_2e}&  = & \cmatrix{e_1&0\cr 0 &r_1r_2}.
 \end{array}
\]

%%%%%%%%%%%%%%%%%%%%%%%%%%%%%%%%%%%%%%%%%

%: \exer{exoChinois}
\exer{exoChinois} On pose $\fb_i=\prod_{j:j\neq i}\fa_j$.
On note  $\varphi:\gA\to\prod_{k=1}^n \gA\sur{\fa_k}$ l'application
canonique.
\'Ecrivons $a_{ij}+a_{ji}=1$ pour $i\neq j$ avec $a_{ij}\in \fa_i$, $a_{ji}\in \fa_j$.
On écrit 
$$
1=\prod\nolimits_{k:k\neq i}(a_{ik}+a_{ki})=\big(\prod\nolimits_{k:k\neq i}a_{ki}\big)+b_i=e_i+b_i \eqno(\#)
$$
avec $b_i\in \fa_i$ et $e_i\in \fb_i$, donc
$$ 
e_i\equiv 0 \mod \fb_i \et e_i\equiv 1 \mod \fa_i \eqno(+)
$$
En conséquence, pour $\xn\in\gA$ 
$$
\varphi \big(\som_{i=1}^ne_ix_i\big) =(x_1 \mod \fa_1,\ldots,x_n \mod \fa_n)
%\eqno(*)
$$
ce qui montre que $\varphi$ est surjective. Le \tho de \fcn donne alors l'\iso 
$\theta:\gA/\fa\to\prod_i \gA/\fa_i$ car on a évidemment $\Ker \varphi=\bigcap\nolimits_{k=1}^n \fa_k=\fa$. Les congurences $(+)$ montrent
que les $\pi(e_i)\in\gA/\fa$ donnent
par $\theta$ le \sfio associé à la structure de produit $\prod_i\gA/\fa_i$.
Vus dans ce produit, les \elts de $\fa_1$ sont ceux dont la première coordonnée est nulle: ils forment donc bien l'\id engendré par $\varphi(1-e_1)$.
Autrement dit en remontant dans $\gA/\fa$, $\pi(\fa_1)=\pi(\gen{1-e_1})$,
et en remontant dans $\gA$, $\fa_1=\fa+\gen{1-e_1}$.

\sni  
L'\egt  $\bigcap\nolimits_{k=1}^n \fa_k=\prod\nolimits_{k=1}^n \fa_k$ se démontre par \recu sur $n$ pour $n\geq 2$ en notant que $(\#)$ implique
que $\fa_i$ et $\fb_i$ sont \comz. Voyons l'initialisation, \cad le cas $n=2$: si $x\in \fa_1\cap \fa_2$ 
et si $a+b=1$ avec $a\in \fa_1$ et $b\in \fa_2$, alors $x=ax+bx$, avec $ax\in \fa_1 \fa_2$
parce que $x\in \fa_2$ et  $bx\in \fa_1 \fa_2$
parce que $x\in \fa_1$, donc $x\in \fa_1\fa_2$.

%%%%%%%%%%%%%%%%%%%%%%%%%%%%%%%%%%%%%%%%%

\exer{exoFacileGrpElem4}
La matrice $D_0 = \crmatrix {0 & -1\cr 1& 0\cr}$ transforme $\cmatrix {x\cr
y\cr}$ en $\crmatrix {-y\cr x\cr}$, donc $D_0^2 = -\I_2$ et $D_0^3 = -D_0
= D_0^{-1}$. \\
On a aussi $D_0=\rE_{12}(1)\rE_{21}(-1)\rE_{12}(1)$,
$D_0 D_q= -\rE_{12}(q)$ et $ D_q D_0=- \rE_{21}(q)$.

%%%%%%%%%%%%%%%%%%%%%%%%%%%%%%%%%%%%%%%%%

\exer{exoTroisiemeLemmeLiberte}
Notons $(e_1, \ldots, e_n)$ la base canonique de~$\Ae n$ et $(f_1, \ldots, f_n)$
les $n$ colonnes de~$F$.  On peut supposer que le mineur  principal \iv
 est en position nord-ouest
%s'appuie sur les $k$ premières lignes et les $k$ premières
%colonnes
de sorte que $(f_1, \ldots, f_k, e_{k+1}, \ldots, e_n)$ est une base
de~$\Ae n$. \\
Puisque $F(f_j) = f_j$, la matrice de~$F$ dans cette base est $G
\eqdf {\rm def} \cmatrix {\I_k & *\cr 0 & *\cr}$.

La matrice $G$ est \idme ainsi
que sa transposée $G'$.
On applique au \prr $G'$ l'opération que l'on vient de faire
subir  à~$F$.\\
%avec la base $g'_1, \ldots, g'_k, e_{k+1}, \ldots, e_n$~;
Puisque $G'(e_j) \in \bigoplus%\limits
_{i \ge k+1} \gA e_i$ pour $j \ge k+1$,
la matrice de $G'$ dans la  nouvelle base est de la forme $H = \cmatrix {\I_k &
0\cr 0 & *\cr}$, d'où le résultat car $F$ est semblable \hbox{à $\tra{H}$}.

%%%%%%%%%%%%%%%%%%%%%%%%%%%%%%%%%%%%%%
\exer{exoABArang1}{
On a des $b_{ji} \in \gA$
tels que $1 = \sum_{i,j} b_{ji} a_{ij}$.  Soit $B \in \Ae {m \times n}$
définie par $B = (b_{ji})$. Vérifions que $ABA = A$: $(ABA)_{ij} =
\sum_{l,k} a_{il} b_{lk} a_{kj}$.  \\
Mais $\dmatrix {a_{il} & a_{ij}\cr
a_{kl} & a_{kj}\cr} = 0$, donc $(ABA)_{ij} = \sum_{l,k} a_{ij} a_{kl}
b_{lk} = a_{ij} \sum_{l,k} a_{kl} b_{lk} = a_{ij}$.   En
conséquence, $AB$ est un \prrz.  \\
Montrons que
$AB$ est de rang~$1$. On a
$\Tr(AB) = \sum_i (AB)_{ii} = \sum_{i,j} a_{ij} b_{ji} = 1$,
\hbox{donc $\cD_1(AB)=1$}.
Par ailleurs,
  $\cD_2(AB)\subseteq\cD_2(A)=0$.
}

%%%%%%%%%%%%%%%%%%%%%%%%%%%%%%%%%%%%%%%%%
\exer{exoABAabstrait} {\emph{1.}
Fixons une forme \lin $\mu$. L'application $E^{r+1} \to \gA$ définie~par

\snic{(y_0, \ldots, y_r) \mapsto \som_{i=0}^r
(-1)^i f(y_0, \ldots, y_{i-1}, \widehat {y_i}, y_{i+1}, \ldots, y_r)\, \mu(y_i),}

%\sni
où $\widehat {y_i}$ symbole de l'omission de l'\eltz, est une forme
$(r+1)$-\lin alternée.

  D'après l'hypothèse $\Vi_{r+1}(x_1,
\ldots, x_n) = 0$ et l'injectivité de $E \mapsto E{\sta}{\sta}$,
on obtient

\snic{
\som_{i=0}^r
(-1)^i f(y_0, \ldots, y_{i-1}, \widehat {y_i}, y_{i+1}, \ldots, y_r)\, y_i = 0.}

%\sni
Notons $y$ au lieu de $y_0$ et réalisons l'opération suivante: dans
l'expression

\snic{
(-1)^i f(y, \ldots, y_{i-1}, \widehat {y_i}, y_{i+1}, \ldots, y_r),}

%\sni
amenons $y$ entre $y_{i-1}$ et $y_i$. La permutation ainsi réalisée
nécessite une multiplication par $(-1)^{i-1}$. On obtient alors la
deuxième \egt dans laquelle tous les signes \gui {ont
disparu}. Par exemple avec $r = 4$, l'expression
$$
\begin {array} {c}
f(\widehat {y}, y_1, y_2, y_3, y_4)y -
f(y, \widehat {y_1}, y_2, y_3, y_4)y_1 +
f(y, y_1, \widehat {y_2}, y_3, y_4)y_2 -
\\
f(y, y_1, y_2, \widehat {y_3}, y_4)y_3 +
f(y, y_1, y_2, y_3, \widehat {y_4})y_4 =
\\
f(y_1, y_2, y_3, y_4)y -
f(y, y_2, y_3, y_4)y_1 +
f(y, y_1, y_3, y_4)y_2 -
\\
f(y, y_1, y_2, y_4)y_3 +
f(y, y_1, y_2, y_3)y_4
\\
\end {array}
$$
n'est autre que
$$
\begin {array} {c}
f(y_1, y_2, y_3, y_4)y -
f(y, y_2, y_3, y_4)y_1 -
f(y_1, y, y_3, y_4)y_2 -
\\
f(y_1, y_2, y, y_4)y_3 -
f(y_1, y_2, y_3, y)y_4.
\end {array}
$$
Une preuve plus expéditive: on applique une forme \lin
$\mu$ à la dernière expression ci-dessus, on vérifie
que l'application obtenue $(y, y_1, y_2, y_3, y_4) \mapsto \mu(\ldots)$ est
$5$-\lin alternée donc nulle d'après les hypothèses.

\noindent  \emph{2.}
Traitons le cas $r = 3$. On a une hypothèse
$$
1 = \som_{ijk} \alpha_{ijk} f_{ijk}(x_i,x_j,x_k), \qquad
\hbox {$f_{ijk}$ $3$-\lin alternée sur $E$.}
$$
On définit $\pi : E \to E$ par:
$$
\pi(x) = \som_{ijk} \alpha_{ijk}
[ f_{ijk}(x,x_j,x_k)x_i +  f_{ijk}(x_i,x,x_k)x_j +  f_{ijk}(x_i,x_j,x) x_k].
$$
Il est clair que l'image de $p$ est contenue dans le sous-module $\som \gA x_i$.
De plus, pour $ {x \in \som \gA x_i}$, on a
$$
f_{ijk}(x,x_j,x_k)x_i +  f_{ijk}(x_i,x,x_k)x_j +  f_{ijk}(x_i,x_j,x) x_k =
f_{ijk}(x_i,x_j,x_k)x.
$$
D'où $\pi(x) = x$: l'endomorphisme $\pi : E \to E$ est un \prr d'image
$\som \gA x_i$. On voit que $p$ est de la forme $\pi(x) = \som_i \alpha_i(x)
x_i$ i.e. $\pi = \psi \circ \varphi$ et que $\pi \circ \psi = \psi$.

\noindent  \emph{3.}
Le module $E$ en question est $\Ae{m}$ et les vecteurs $x_1$, $
\ldots$, $x_n$ sont les colonnes de~$A$. On a $\psi = A : \Ae{n} \to \Ae{m}$, et si
l'on note $B \in \Ae {n \times m}$ la matrice de $\varphi : \Ae{m} \to \Ae{n}$, on
a bien~$ABA = A$.  Alors, l'\ali $AB : \Ae{m} \to \Ae{m}$ est un projecteur de même
image que~$A$.

}

%%%%%%%%%%%%%%%%%%%%%%%%%%%%%%%%%%%%%%%%%
\exer{exoDetExtPower}{
Voyons d'abord le cas où $u=\Diag(\lambda_1, \ldots, \lambda_n)$. On dispose d'une base
$(e_I)$ de $\Vi^k(\Ae{n})$ indexée par les parties $I \subseteq \{1,
\ldots, n\}$ de cardinal~$k$:

\snic{e_I = e_{i_1} \wedge \cdots \wedge e_{i_k} \qquad
I = \so{i_1 < \cdots < i_k }.
}

%\sni
Alors, $u_k$ est diagonale dans la base $(e_I)$: $u_k(e_I) =
\lambda_I e_I$ avec $ {\lambda_I = \prod_{i \in I} \lambda_i}$.
Il s'ensuit que $ {\det(u_k) = \prod_{\#I = k} \prod_{i \in I}
\lambda_i}$.  Reste à déterminer, pour un $j$ donné dans~$\lrbn$, le nombre d'occurrences de $\lambda_j$ dans le produit ci-dessus. Autrement
dit, combien de parties $I$, de cardinal~$k$, contenant~$j$? Autant que de
parties de cardinal $k-1$ contenues dans $\{1, \cdots, n\} \setminus \{j\}$,
i.e.  ${n-1 \choose k-1}$. Le résultat  est démontré pour
une matrice générique. Donc il est vrai pour une matrice quelconque.  Le
deuxième point résulte des \egts

\snic{\dsp{n-1 \choose k-1} + {n-1 \choose n-k-1} =
{n-1 \choose k-1} + {n-1 \choose k} = {n \choose k}.
}
}

%%%%%%%%%%%%%%%%%%%%%%%%%%%%%%%%%%%%%%%%%
\exer{exolemdeterblocs}{Le cas \gnl se traite par \recu sur $n$.
On considère l'anneau de
\pols $\ZZ[(x_{ij})]$ à $n^2$ \idtrs et la matrice \uvle
$A=(x_{ij})$ à \coes dans cet anneau. Notons $\Delta_{1k} \in \ZZ[(x_{ij})]$ le  cofacteur  de $x_{1k}$ dans $A$. Ces cofacteurs vérifient les
identités:
$$
\som_{j = 1}^n x_{1j} \Delta_{1j} =
\det A
%= \sum_{\sigma \in \Sn} \varepsilon(\sigma)
%x_{1\sigma_1} x_{2\sigma_2} \cdots x_{n\sigma_n}
,
\;\;\;\;
%:HHH \Delta_{1j} ci-dessous
\som_{j = 1}^n x_{ij} \Delta_{1j} = 0 \quad \hbox {pour } i > 1  .
$$
  Puisque les $N_{kl}$ commutent deux à deux, la spécialisation
  $x_{kl} \mapsto N_{kl}$ est légitime. Notons $N'_{1j} =
\Delta_{1j}(x_{kl} \mapsto N_{kl})$, alors on a

\snic{
N'_{11} =
\sum_{\sigma \in \rS_{n-1}} \varepsilon(\sigma)
N_{2\sigma_2} N_{3\sigma_3} \ldots N_{n\sigma_n}.}

%\sni
Définissons $N'$ par:
$$
N' = \cmatrix {
N'_{11}      &  0    & \cdots & 0      \cr
N'_{12}      &  \rI_m  &        & \vdots \cr
\vdots       &  \vdots     & \ddots & 0      \cr
N'_{1n}      &  0    & \cdots       & \rI_m    \cr}\,,
\; \hbox {d'où} \;\;
N N' = \cmatrix {
\Delta  & N_{12} & \cdots & N_{1n} \cr
0       & N_{22} & \cdots & N_{2n} \cr
\vdots  &        &        & \vdots \cr
0       & N_{n2} & \cdots & N_{nn} \cr
}.
$$
En prenant les \detersz, on obtient
$$
\det(N) \det(N'_{11}) = \det(\Delta)
\det \cmatrix {
N_{22} & \cdots & N_{2n} \cr
\vdots &        & \vdots \cr
N_{n2} & \cdots & N_{nn}  }.
$$
L'\hdr fournit les \egts
$$
\det \cmatrix {
N_{22} & \cdots & N_{2n} \cr
\vdots &        & \vdots \cr
N_{n2} & \cdots & N_{nn} \cr} =
\det\Big( \sum_{\sigma \in \rS_{n-1}} \varepsilon(\sigma)
N_{2\sigma_2} N_{3\sigma_3} \cdots N_{n\sigma_n} \Big)
= \det(N'_{11}).
$$
La simplification par l'\elt \ndz $\det(N'_{11})$ donne
l'\egt $\det(N) = \det(\Delta)$.
}
%%%%%%%%%%%%%%%%%%%%%%%%%%%%%%%%%%%%%%%%%%%%%%%%%%%%%%%%%%%%%%%%%%%%%%%%%%%

%: sinotenglish
\sinotenglish{
%: \exer{exoABegal0}
\exer{exoABegal0} ~
\\ \emph{1.} On peut supposer $r\leq n$.
On considère un mineur $\mu$ d'ordre $r$ de $A$, \spdg on suppose 
la matrice carrée extraite correspondante $A_1$ située dans le coin
nord-ouest. On écrit
$$
A =\blocs{1.3}{1}{1.3}{1.8}{$A_1$}{$X$}{$Y$}{$Z$}\,,
\;
A' =\blocs{1.3}{1.8}{1.3}{1}{$\wi{A_1}$}{$0$}{$0$}{$0$}\,,
\;
A' A =\blocs{1.3}{1}{1.3}{1}{$\mu\,\I_r$}{$0$}{$0$}{$0$} 
$$
avec $\mu=\det(A_1)$. On partage $B$ en $B_1\in\Ae{r\times p}$ et $B_2\in\Ae{(n-r)\times p}$ 
$$
B =\blocs{4}{0}{1.3}{1}{$B_1$}{}{$B_2$}{}\,,
\quad
A' A B =\blocs{4}{0}{1.3}{1}{$\mu\, B_1$}{}{$0$}{}=0.
$$
Un mineur $\nu$ d'ordre $s$ de $B$ est le \deter d'une matrice carrée extraite
$C$
dont au moins une ligne est dans $B_1$. On exprime ce mineur $\nu$ au moyen d'un développement de Laplace en partageant $C$ en deux parties, l'une correspondant aux lignes empruntées à $B_1$, l'autre, éventuellement vide,
correspondant aux lignes empruntées à $B_2$. On voit que $\mu\,\nu$
est dans l'\id engendré par les \coes de~$\mu\,B_1$, donc $\mu\,\nu=0$. Ce qu'il fallait démontrer.  
 
 \emph{2.} Il suffit d'appliquer le point \emph{1} avec l'anneau
$\gA\sur{\cD_1\!(AB)}$.

 \emph{3.} Supposons par exemple $r+s \geq n+2$.
Le même calcul que dans le point \emph{1} donne cette fois-ci $\mu^2\,\cD_2(B_1)=\cD_2(\mu\,B_1)\subseteq \cD_2(AB)$. On utilise le développement de Laplace pour exprimer $\nu=\det(C)$, la matrice $C$ a maintenant 
au moins deux lignes empruntées à $B_1$, on obtient donc  $\mu^2\,\nu\in \cD_2(AB)$.

%%%%%%%%%%%%%%%%%%%%%%%%%%%%%%%%%%%%%%%%%%%%%%%%%%%%%%%%%%%%%%%%%%%%%%%%%%%

\exer{exosexcNoether1}
Notez que l'on est dans le  contexte usuel des \thos de \fcn de \Noez, ici les deux sous-modules sont notés $N_1$ et $N_2$, ce qui montre mieux la symétrie de la situation.

\emph{1.} Si l'anneau est un corps $\gK$ avec 

\snic{  \dim_\gK(N_i)=n_i$,    $\dim_\gK(N_1+N_2)=n$ et  $\dim_\gK(N_1\cap N_2)=n',}

la suite exacte est automatiquement scindée (\tho de la base incomplète)
et l'on obtient l'\egt classique
$n+n'=n_1+n_2$.

\emph{2.} Notons $N=N_1\cap N_2$. La suite exacte est scindée si l'on a une section 

\snic{\sigma:N_1+N_2\lora N_1\times N_2.}

On écrit $\sigma(x_1+x_2)=\sigma_1(x_1)+\sigma_2(x_2)$, avec $\sigma_1(x_1)=\big(x_1-\alpha_1(x_1),\alpha_1(x_1)\big)$ 
(en effet~\hbox{$\pi(\sigma_1(x_1))=x_1$}) et $\sigma_2(x_2)=\big(\alpha_2(x_2),x_2-\alpha_2(x_2)\big)$. \\
On a donc $\alpha_1:N_1\to N$ \hbox{et $\alpha_2:N_2\to N$}. L'application  $\sigma$ est bien définie \ssi pour \hbox{tout $y\in N$}, on a $\sigma_1(y)=\sigma_2(y)$, i.e. $\alpha_1(y)+\alpha_2(y)=y$.
\\
En résumé, la suite est scindée \ssi on peut trouver $\alpha_1:N_1\to N$ \hbox{et $\alpha_2:N_2\to N$} vérifiant $\alpha_1(y)+\alpha_2(y)=y$ pour $y\in N$.
\\
Cette condition est un peu mystérieuse.
Elle est satisfaite par exemple si $N$ est facteur direct dans $N_1$ en prenant $\alpha_2=0$ et $\alpha_1$ une \prn de $N_1$ sur $N$. 
Mais en général, le critère n'est pas très parlant.
\\
Prenons par exemple avec un anneau à pgcd $\gA$, les sous-modules $N_1=a_1\gA$ et  $N_2=a_2\gA$ du module $\gA$. Si $g$ est le pgcd et $m$ le ppcm, on a $N=m \gA$ \hbox{avec $a_1=gc_1$, $a_2=gc_2$, $m=a_1c_2=a_2c_1$}. Pour obtenir une section il nous faut des $\alpha_i:N_i\to N$.
On a alors $\alpha_1(a_1)=mx$, $\alpha_2(a_2)=my$,
ce qui donne 

\snic{\alpha_1(m)=c_2mx$, \,\, $\alpha_2(m)=c_1my,}

et l'\egt $\alpha_1(m)+\alpha_2(m)=m$ signifie $c_2x+c_1y=1$.
En bref les deux \eltsz~$c_1$ et~$c_2$ premiers entre eux doivent engendrer l'\id $\gen{1}$.
Ainsi, la suite sera toujours scindée si $\gA$ est un domaine de Bézout,
mais pas toujours scindée dans le cas contraire.

%%%%%%%%%%%%%%%%%%%%%%%%%%%%%%%%%%%%%%%%%%%%%%%%%%%%%%%%%%%%%%%%%%%%%%%%%%%

\exer{exosexcNoether2}
 On étudie le complexe:

\snic{0\lora M/(N_1\cap N_2)\vvers{j}M/N_1 \times M/N_2 \vvers \pi M/(N_1+N_2)\lora 0}

où $j(\wh x)=(\wi x,-\uci x)$ et $\pi(\wi y,\uci z)=\ov y+\ov z$.

 \emph{1.} \emph{Le complexe est exact.} Tout d'abord $j(\wh x)=0$ \ssi 
$\wi x$ et $\uci x$ sont nuls, i.e. $x\in N_1\cap N_2$.
Ceci donne l'exactitude en $M/(N_1\cap N_2)$.
Ensuite $\pi(\wi y,0)=\ov y$ donc $\pi$ est surjective.
Ceci donne l'exactitude en $M/(N_1+ N_2)$. 
\\
Soit maintenant un \elt
arbitraire $(\wi y,\uci z)\in\Ker\pi$, 
i.e. $y+z\in N_1+N_2$. \\
On écrit $y+z=y_1+z_2$ avec $y_1\in N_1$ et
$z_2\in N_2$, d'où $(y-y_1)=-(z-z_2)$. 
\\
Alors $\wi y=\wi{y-y_1}$, $\uci z =\uci{z-z_2}$ donc $(\wi y,\uci z)=j(\wh x)$ pour $x=y-y_1$.
\\
Ceci donne l'exactitude au milieu.

\emph{2.} Si l'anneau est un corps $\gK$ avec 

\snic{\dim_\gK(M)=m$,  $\dim_\gK(N_i)=n_i$,  $\dim_\gK(N_1+N_2)=n$ et  $\dim_\gK(N_1\cap N_2)=n',}

la suite exacte est automatiquement scindée (\tho de la base incomplète)
et l'on obtient  $(m-n)+(m-n')=(m-n_1)+(m-n_2)$, \cad l'\egt classique
$n+n'=n_1+n_2$.   

\emph{3.} On prend {les mêmes modules} $M=\gA$, $N_1=a_1\gA\subseteq M$, $N_2=a_2\gA\subseteq M$ que pour la fin de la solution de l'exercice \ref{exosexcNoether1}.
On suppose que $\gA$ est un anneau à pgcd, on reprend les mêmes notations.
\\
Une section $\sigma:\aqo\gA{a_1,a_2}\to \aqo\gA{a_1}\times \aqo\gA{a_2}$
est a priori donnée par deux \alis $\sigma_i:\aqo\gA{a_1,a_2}\to \aqo\gA{a_i}$. On les définit par 

\centerline{$\sigma_1(\ov 1)=\wi{uc_2}$ et $\sigma_2(\ov 1)=\uci{vc_1}.$}

Pour que $\pi(\sigma(\ov 1))=\ov 1$, il faut que $uc_2+vc_1\equiv 1\mod\gen{a_1,a_1}$, ce qui signifie~\hbox{$\gen{c_1,c_2}=\gen{1}$} dans $\gA$.
Ainsi on retrouve que la suite est scindée si $\gA$ est de Bézout, et qu'elle peut ne pas être scindée dans le cas contraire. 
%%%%%%%%%%%%%%%%%%%%%%%%%%%%%%%%%%%%%%%%%%%%%%%%%%%%%%%%%%%%%%%%%%%%%%%%%%%

% \exer{exoNoeSeco}
\exer{exoNoeSeco}
\emph{1.} Rappelons le diagramme que l'on souhaite: les deux premières lignes et les deux premières colonnes sont des \secos canoniques et~$G''=F/(E+F')$ 

\smallskip 
{\small\centerline{$
\xymatrix@C=3.2em@R=1.5em {
&0\ar[d]  &  0\ar[d]  & 0\ar[d]
\\
0\ar[r]&E'\ar[d]_{\iota_E} \ar[r]^{\iota '} 
&F'\ar[d]_{\iota_F}\ar[r]^{\pi ' } &
G'\ar[d]^{\iota_G}\ar[r]&0
\\
0\ar[r]&E\ar[d]_{\pi_E} \ar[r]^{\iota} 
&F\ar[d]_{\pi_F}\ar[r]^{\pi} &
G\ar[d]^{\pi_G}\ar[r]&0
\\
0\ar[r]&E''\ar[d] \ar[r]_{\iota ''} &F''\ar[d]\ar[r]_{\pi'' } &G''\ar[d]\ar[r]&0
\\
&0  &  0  & 0
\\
}
$
}}

Le \tho de \fcn de \Noe donne un \iso naturel 

\snic{G'=F'/(E\cap F')\vers j (E+F')/E, \; \hbox{ avec }(E+F')/E\subseteq F/E.}

Cet \iso est
défini par $j\big(\pi'(x)\big)=\pi(x)=\pi\big(\iota_F(x)\big)$. Cela nous dit que l'on a une \ali injective $\iota_G:G'\to F/E$
qui vérifie 

\snic{\iota_G\big(\pi'(x)\big)=\pi\big(\iota_F(x)\big),}

\cad qui rend le diagramme commutatif.
\\
Comme $\pi'$ est surjective, $\iota_G$ est même l'unique \Ali  qui rend le diagramme commutatif.
\\
De même on a une unique \ali  

\snic{E''=E/(E\cap F')\vers{\iota''} F''=F/F'}

qui rend le diagramme commutatif, et $\iota''$ est injective, d'image $(E+F')/F'$.
\\
La surjection canonique $\theta:F\to F/(E+F')$ se factorise de manière unique via~$\pi$
parce que $\Ker\pi=E\subseteq E+F'=\Ker\theta$ et l'on obtient ainsi $\pi_G:G\to G''$ satisfaisant $\pi_G\circ \pi=\theta$.\\
De même on obtient une \ali surjective $\pi'':F''\to G''$ satisfaisant
l'\egt $\pi''\circ \pi_F=\theta$.\\
On a ainsi obtenu un diagramme commutatif complet. Il reste à voir que 
la troisième suite verticale et la  troisième suite horizontale sont exactes.\\
Or $\Ker\pi_G=\pi(\Ker\theta)=(E+F')/E=S/E$ et $\Im\iota_G=\Im j=S/E$. Ceci montre que la suite verticale est exacte, et l'on vient de redécouvrir le
\tho de Noether qui établit
un \iso naturel 

\snic{G/\Ker\pi_G=(F/E)/(S/E)\vers \alpha F/S=G''}

\snii
satisfaisant $\alpha\big(\ov{\pi(x)}\big)=\pi(x)$ pour tout $x\in F$.
\\
Symétriquement la troisième suite horizontale est exacte.

\emph{2.}
On a déjà vu que la commutativité du diagramme sur les deux premières lignes (resp. colonnes) impose l'\ali $\iota_G$  (resp. $\iota''$).
Il nous reste à voir si l'affirmation analogue concernant $\pi_G$ et $\pi''$ est correcte.
On  suppose que l'on a des \alis $\lambda_G$ et $\lambda''$ qui satisfont l'\egt $\lambda_G\circ \pi=\lambda''\circ \pi_F$ et que toutes les lignes et colonnes sont exactes. On obtient donc $\Ker\lambda_G =\Im\iota_G=E+F'$, mais ceci ne force
pas l'\egt $\lambda_G\circ \pi=\theta$. Par exemple, si $\beta$ est un \auto arbitraire de $G''$, on peut prendre $\lambda_G=\beta\circ \pi_G$
et $\lambda''=\beta\circ \pi''$.

\rems\\ 
i. Le point \emph{2} montre une certaine absence de symétrie (regrettable)
dans la situation. Cela sera élucidé d'une certaine manière dans l'exercice \ref{exoNoeSeco2}.
\\
On peut néanmoins conclure cet exercice comme suit. 
\\
\emph{Supposons que:
\begin{itemize}
\item les deux premières suites horizontales et les deux premières suites
verticales sont des \secos canoniques,
\item  et $\theta=\pi_G\circ \pi=\pi''\circ \pi_F$.
\end{itemize}
Alors il y a un unique diagramme commutatif de la forme annoncée, 
et il rend les
troisièmes suites horizontale et verticale exactes.}
\\
Notons que ceci constitue une forme particulièrement précise des \thos de Noether dans la mesure où les \isos de Noether sont ici complètement explicites et déterminés de manière unique.

ii.
Remarquons aussi que l'hypothèse selon laquelle certaines injections et surjections sont canoniques  est un peu artificielle
dans la mesure où  $\iota_G$ et $\iota''$ ne sont pas des injections canoniques et $\pi_G$ et $\pi''$ ne sont pas des surjections canoniques.
Voir à ce sujet la remarque à la fin du corrigé de l'exercice \ref{exoNoeSeco2}.
\eoe

%%%%%%%%%%%%%%%%%%%%%%%%%%%%%%%%%%%%%%%%%%%%%%%%%%%%%%%%%%%%%%%%%%%%%%%%%%%

% \exer{exoNoeSeco2}
\exer{exoNoeSeco2}
\emph{1.} Rappelons le diagramme donné en hypothèse dans le cas (non vraiment restrictif) où
les injections et les surjections sont toutes canoniques.

\smallskip 
{\small\centerline{$
\xymatrix@C=3.2em@R=1.5em {
&0\ar[d]  &  0\ar[d]  & 
\\
0\ar[r]&E_0\ar[d]_{\jmath_E} \ar[r]^{\iota _0} 
&F'\ar[d]_{\iota_F}\ar[r]^(.3){\pi _0 } &
G_0=F'/E_0\ar[r]&0
\\
0\ar[r]&E \ar[r]^{\iota} 
&F\ar[r]^(.35){\pi} &
G=F/E\ar[r]&0
}
$
}}

\smallskip
\emph{1a.} L'\ali $\jmath_G$ doit être obtenue en factorisant $\pi\circ \iota_F$.
Si elle existe, elle est unique, et elle existe \ssi $\Ker\pi_0\subseteq \Ker(\pi\circ \iota_F)$. Or $\Ker\pi_0=\Im\iota_0$.  La condition  équivaut donc à  $\pi\circ \iota_F\circ \iota_0=0$.
Or  $\pi\circ \iota_F\circ \iota_0= \pi\circ \iota\circ \jmath_E$
et $\pi\circ \iota=0$.
 
\emph{1b.} Puisque $\pi_0$ est surjective, on a $\Im\jmath_G=\Im(\jmath_G\circ \pi_0)=\Im(\pi\circ \iota_F)=\pi(F')$. Enfin $\pi^{-1}\big(\pi(F')\big)=E+F'=S$. Ainsi $\Im\jmath_G=S/E\subseteq F/E$. 

\emph{1c.} L'application $\jmath_G$ est injective \ssi le noyau de $\jmath_G\circ\pi_0 $ est égal au noyau de  $\pi_0$, qui est  $E_0$, \cade si $\Ker(\pi\circ \iota_F)=E_0$. 
\\
Or $\Ker(\pi\circ \iota_F)=\iota_F^{-1}(E)=E\cap F'$. Ainsi, la condition est bien $E_0=E'$.
\\
Ceci nous ramène à la situation
de l'exercice \ref{exoNoeSeco}.

\smallskip \emph{2.} Rappelons le diagramme donné en hypothèse dans lequel on peut supposer que $\iota$ est une injection canonique et $\pi_F$ une surjection canonique~\hbox{$F\to F''=F/F'$} \hbox{avec $F'=\Ker\pi_F$}.

\smallskip 
{\small\centerline{$
\xymatrix@C=3.2em@R=1.5em {
%&  & 0\ar[d]            & 0\ar[d]
%\\
%&\Ker\pi_E\ar[d]_{\iota_E}  
%&\Ker\pi_F\ar[d]_{\iota_F} & \Ker\pi_3\ar[d]^{\iota_3}
%\\
0\ar[r]&E%\ar[d]_{\pi_E} 
        \ar[r]^{\iota} &F\ar[d]_{\pi_F}\ar[r]^{\pi} & G\ar[d]^{\theta_G}\ar[r] & 0
\\
0\ar[r]&E_3%\ar[d] 
         \ar[r]_{\iota _3} &F''\ar[d]\ar[r]_{\pi_3 }    & G_3\ar[d]\ar[r] & 0
\\
&%0  
&  0  & 0
\\
}
$
}}

\smallskip \emph{2a.} Comme $\iota$ est une injection canonique, l'\ali  $\beta:E\to E_3$ que l'on veut définir doit satisfaire pour $x\in E$ l'\egt  $\iota_3\big(\beta(x)\big)=\pi_F(x)$, ce qui est possible si $\pi_F(E)\subseteq \iota_3(E_3)$, \cad si $\pi_F(E)\subseteq \Ker\pi_3$,
\cad encore~\hbox{$\pi_3\circ \pi_F\circ \iota=0$}. 
Or $\pi_3\circ \pi_F\circ \iota=\theta_G\circ \pi \circ \iota$ et $\pi \circ \iota=0$. 
\\
Ainsi $\beta$ est bien définie, et elle est unique parce que $\iota_3$ est injective.

\emph{2b.} Puisque $\iota_3$ est injective, on a 
 
\snic{\Ker\beta=\Ker(\iota_3\circ \beta)=\Ker(\pi_F\circ \iota)=\iota^{-1}(\Ker\pi_F)=\iota^{-1}(F')=E\cap F'.}

\emph{2c.} L'\ali $\beta$ est surjective \ssi $\Im(\iota_3\circ \beta)\supseteq \Ker\pi_3$, \cade si %$\Im(\pi_F\circ \iota)=E_3$, i.e. 
$\pi_F(E)\supseteq \Ker\pi_3$, ou aussi $\pi_F^{-1}\big(\pi_F(E)\big)\supseteq \pi_F^{-1}\big(\Ker\pi_3\big)$, i.e. enfin $E+\Ker\pi_F\supseteq S_3$.  
Ainsi $\beta$ est surjective \ssi $E+F'= S_3$. 

Dans ce cas, en prenant $E'=E\cap F'$ on retrouve à \isos près la situation de l'exercice \ref{exoNoeSeco}. Voyons ceci \prmtz.\\
Tout d'abord puisque $\beta$ est surjective de noyau $E'$, on a 
une unique \ali $\alpha_E=E/E'\to E_3$ qui satisfait $\alpha_E\circ \pi_E=\beta$ ($\pi_E:E\to E/E'$ surjection canonique). Et $\alpha_E$ est un \isoz.
\\
Ensuite, puisque $F''=F/F'$ et $\Ker(\pi_3\circ \pi_F)=E'+F$, on a un unique \ali $\alpha_G:F/S=G''\to G_3$ qui satisfait $\alpha_G\circ \pi''=\pi_3$
(avec $\pi''$ comme dans l'exercice \ref{exoNoeSeco}), et $\alpha_G$ est un \isoz.
\\
Enfin, vu la commutativité du diagramme, on a $\alpha_G\circ \pi_G=\theta_G$
(avec $\pi_G$ comme dans l'exercice \ref{exoNoeSeco}).
\\
Ainsi, on retrouve bien, modulo les \isos $\alpha_E$ et $\alpha_G$, les deux dernières lignes du diagramme de l'exercice \ref{exoNoeSeco}. 

\rem
En se libérant de l'hypothèse rajoutée un peu artificiellement pour faciliter la \demz, selon laquelle les injections et surjections données au départ sont canoniques, le point \emph{1} donnerait l'énoncé suivant,  
\\
\emph{Un diagramme commutatif de suites exactes du type $(I)$  
peut être complété  en un diagramme commutatif complet $(C)$ de suites exactes \ssi  on a l'\egt $\Ker(\iota\circ \iota_E)=\Ker \iota_E\cap \Ker \iota'$.
Et dans ce cas $(C)$ est essentiellement unique.}

\smallskip 
{\small\centerline{$
\xymatrix@C=3.2em@R=1.5em {
&&0\ar[d]  &  0\ar[d]&~&~~   
\\
(I)&0\ar[r]&E'\ar[d]_{\iota_E} \ar[r]^{\iota '} 
&F'\ar[d]_{\iota_F}\\
&0\ar[r]&E \ar[r]^{\iota} 
&F
\\
}
$
}}

\medskip 
{\small\centerline{$
\xymatrix@C=3.2em@R=1.5em {
&&0\ar[d]  &  0\ar[d]  & 0\ar[d]
\\
&0\ar[r]&E'\ar[d]_{\iota_E} \ar[r]^{\iota '} 
&F'\ar[d]_{\iota_F}\ar[r]^{\pi ' } &
G'\ar[d]^{\iota_G}\ar[r]&0
\\
(C)&0\ar[r]&E\ar[d]_{\pi_E} \ar[r]^{\iota} 
&F\ar[d]_{\pi_F}\ar[r]^{\pi} &
G\ar[d]^{\pi_G}\ar[r]&0
\\
&0\ar[r]&E''\ar[d] \ar[r]_{\iota ''} &F''\ar[d]\ar[r]_{\pi'' } &G''\ar[d]\ar[r]&0
\\
&&0  &  0  & 0
\\
}
$
}}

\medskip 
En outre, le point \emph{2} fournit un \tho dual.\\
\emph{Un diagramme commutatif de suites exactes du type $(I')$  
peut être complété  en un diagramme commutatif complet $(C)$ de suites exactes \ssi  on a l'\egt $\Ker(\pi_G\circ \pi)=\Ker\pi+\Ker\pi_F$.
Et dans ce cas $(C)$ est essentiellement unique.}

\smallskip {\small\centerline{$
\xymatrix@C=3.2em@R=1.5em {
&~~&~~ 
&F\ar[d]_{\pi_F}\ar[r]^{\pi} &
G\ar[d]^{\pi_G}\ar[r]&0
\\
(I')& &  &F''\ar[d]\ar[r]_{\pi'' } &G''\ar[d]\ar[r]&0
\\
&&  &  0  & 0
\\
}
$
}}

\smallskip 
\rem La dualité qui apparaît ici entre les points \emph{1} et \emph{2} frise maintenant la perfection. Elle a donné lieu à une abstraction qui permet de mieux la comprendre: la théorie des catégories abéliennes. La catégorie opposée d'une catégorie abélienne étant elle-ême abélienne, un énoncé du style de \emph{1} prouvé dans une catégorie abélienne fournit ipso facto un énoncé correct tel que \emph{2.} 
\eoe

%%%%%%%%%%%%%%%%%%%%%%%%%%%%%%%%%%%%%%%%%%%%%%%%%%%%%%%%%%%%%%%%%%%%%%%%%%%

% \exer{exoNoeSeco3}
\exer{exoNoeSeco3}  
On rappelle le diagramme

\smallskip 
{\small\centerline{$
\xymatrix@C=3.2em@R=1.5em {
&0\ar[d]  &  0\ar[d]  & 0\ar[d]
\\
0\ar[r]&E_1\ar[d]_{\iota_E} \ar[r]^{\iota _1} 
&F_1\ar[d]_{\iota_F}\ar[r]^{\pi _1 } &
G_1\ar[d]^{\iota_G}\ar[r]&0
\\
0\ar[r]&E\ar[d]_{\pi_E} \ar[r]^{\iota} 
&F\ar[d]_{\pi_F}\ar[r]^{\pi} &
G\ar[d]^{\pi_G}\ar[r]&0
\\
0\ar[r]&E_2\ar[d] \ar[r]_{\iota _2} &F_2\ar[d]\ar[r]_{\pi_2 } &G_2\ar[d]\ar[r]&0
\\
&0  &  0  & 0
\\
}
$
}}

On suppose \spdg que $\iota$, $\iota_1$, $\iota_E$ et $\iota_F$ sont des injections canoniques, et $\pi$, $\pi_1$, $\pi_E$ et $\pi_F$ sont des surjections canoniques.\\ 
Notons $S_2$ le noyau de l'\ali $\pi_2\circ \pi_F=\pi_G\circ \pi$.

\emph{1.} Supposons tout d'abord la suite $0\to E_1\to F_1 \to G_1 \to 0$
exacte. Alors d'après le point \emph{1} de l'exercice~\ref{exoNoeSeco2},
on a $E_1=E\cap F_1$, donc $E_2= E/(E\cap F_1)$, \hbox{et $\Im\iota_G=(E+F_1)/E\subseteq F/E$}, 
donc $G_2\simeq F/(E+F_1)$.\\
Ceci implique que la troisième ligne est exacte.

Supposons la suite $0\to E_2\to F_2 \to G_2 \to 0$
exacte. Alors d'après le point \emph{2} de l'exercice~\ref{exoNoeSeco2},
on a $\Ker\pi_E=E\cap F_1$, donc $E_1=E\cap F_1$, et le noyau de l'\ali $F\to G_2$ doit être égal à $E+F_1$, ce qui implique $\Ker\pi_G=(E+F_1)/E$.\\
Comme $F_1/(E\cap F_1)\simeq (E+F_1)/E$, ceci implique que la première ligne est exacte.

\emph{2.} Déjà démontré

%%%%%%%%%%%%%%%%%%%%%%%%%%%%%%%%%%%%%%%%%%%%%%%%%%%%%%%%%%%%%%%%%%%%%%%%%%%

\exer{exoDualite} Laissé \alecz.

}
%: fin sinotenglish

%%%%%%%%%%%%%%%%%%%%%%%%%%%%%%%%%%%%%%%%%%%%%%%%%%%%%%%%%%%%%%%%%%%%%%%%%%%
%:prob{exoRationaliteLineaire}
\prob{exoRationaliteLineaire} 
\emph {1.}
Si $A_j$ est une colonne non nulle de $A$, on a $BA_j = e_j$ donc $ABA_j =
A_j$; ainsi $AB$ est l'identité sur $\Im A$ donc $ABA = A$. La matrice $AB$
est triangulaire inférieure, et ses \coes diagonaux sont $0,1$.  La matrice
$BA$ est diagonale et ses \coes diagonaux sont $0,1$.

\snic {
B = \cmatrix {
\cdot&   \cdot&   \cdot&   \cdot&   1&  \cdot&  \cdot&  \cdot&  \cdot\cr
\cdot&   \cdot&   \cdot&   \cdot&   \cdot&  1&  \cdot&  \cdot&  \cdot\cr
\cdot&   \cdot&   1&   \cdot&   \cdot&  \cdot&  \cdot&  \cdot&  \cdot\cr
1&   \cdot&   \cdot&   \cdot&   \cdot&  \cdot&  \cdot&  \cdot&  \cdot\cr
\cdot&   \cdot&   \cdot&   \cdot&   \cdot&  \cdot&  \cdot&  1&  \cdot\cr
\cdot&   \cdot&   \cdot&   \cdot&   \cdot&  \cdot&  \cdot&  \cdot&  \cdot\cr
},
\qquad
BA = \cmatrix {
1&  \cdot&  \cdot&  \cdot&  \cdot&  \cdot\cr
\cdot&  1&  \cdot&  \cdot&  \cdot&  \cdot\cr
\cdot&  \cdot&  1&  \cdot&  \cdot&  \cdot\cr
\cdot&  \cdot&  \cdot&  1&  \cdot&  \cdot\cr
\cdot&  \cdot&  \cdot&  \cdot&  1&  \cdot\cr
\cdot&  \cdot&  \cdot&  \cdot&  \cdot&  \cdot\cr
},}

\snic {
AB = \cmatrix {
1&      \cdot&       \cdot&   \cdot&       \cdot&       \cdot&   \cdot&       \cdot&  \cdot\cr
a_{24}& \cdot&       \cdot&   \cdot&       \cdot&       \cdot&   \cdot&       \cdot&  \cdot\cr
\cdot&      \cdot&       1&   \cdot&       \cdot&       \cdot&   \cdot&       \cdot&  \cdot\cr
a_{44}& \cdot&  a_{43}&   \cdot&       \cdot&       \cdot&   \cdot&       \cdot&  \cdot\cr
\cdot&      \cdot&       \cdot&   \cdot&       1&       \cdot&   \cdot&       \cdot&  \cdot\cr
\cdot&      \cdot&       \cdot&   \cdot&       \cdot&       1&   \cdot&       \cdot&  \cdot\cr
a_{74}& \cdot&  a_{73}&   \cdot&  a_{71}&  a_{72}&   \cdot&       \cdot&  \cdot\cr
\cdot&      \cdot&       \cdot&   \cdot&       \cdot&       \cdot&   \cdot&       1&  \cdot\cr
a_{94}& \cdot&  a_{93}&   \cdot&  a_{91}&  a_{92}&   \cdot&  a_{95}&  \cdot\cr
}.}

%\sni
Le \sul $\Ker AB$ de $\Im A = \Im AB$ dans $\gK^n$ admet comme base les~$e_i$ pour les indices $i$ de lignes ne
contenant pas un indice pivot.\\
Dans l'exemple, $(e_2, e_4, e_7, e_9)$ est une
base de $\Ker AB$.

\noindent\emph {2.}
On obtient $(Q, A')$ par la méthode (classique) d'échelonnement de
Gauss. Si la matrice $B' \in M_{n,m}(\gK)$ vérifie $A'B'A' = A'$, alors $AQB'AQ = AQ$,
donc la matrice $B = QB'$ vérifie $ABA = A$.

\noindent\emph {3.}
Considérons une matrice $B \in \MM_{m,n}(\gK)$ telle que $ABA = A$. Alors, si $y = Ax$
pour un $m$-vecteur à \coes dans un sur-anneau de $\gK$, on a $A(By) = y$,
d'où l'existence d'une solution sur $\gK$, à savoir $By$.

\noindent\emph {4.}
Soient $(u_1, \ldots, u_r)$ un \sgr du $\gK$-\evc $E$, constitué de
vecteurs de $\gK_0^n$; idem pour $(v_1, \ldots, v_s)$ et $F$.  Soit $z \in
\gK_0^n$, que l'on cherche à écrire sous la forme $z = x_1 u_1 + \cdots + x_r u_r + y_1
v_1 + \cdots + y_s v_s$ avec les~$x_i, y_j \in \gK_0$. On obtient ainsi un \sys
$\gK_0$-\lin en les inconnues $x_i, y_j$ qui admet une solution sur $\gK$, donc
\egmt sur $\gK_0$.

\noindent\emph {5.a.}
Si tous les $\pi(e_j)$ sont dans $\gK_0^n$, alors le sous-espace $E$,
engendré par les $\pi(e_j)$,
est \hbox{$\gK_0$-rationnel}. 
Réciproquement, si $E$ est \hbox{$\gK_0$-rationnel},
comme $F$ l'est aussi, on~a, d'après la question précédente, 
$\pi(e_j) \in \gK_0^n$ pour tout $j$.

\noindent\emph {5b.}
Facile maintenant: $\gK_0$ est le sous-corps engendré  par les composantes des vecteurs~$\pi(e_j)$.

\noindent\emph {5c.}
Le corps de rationalité d'une matrice strictement échelonnée est le
sous-corps engendré par les \coes de la
matrice. Par exemple avec $E = \Im A \subset \gK^5$:

\snic{
A = \bordercmatrix [\lbrack\rbrack]{
    & w_1 & w_2 & w_3 \cr
e_1 & 1   & 0   & 0 \cr
e_2 & a   & 0   & 0 \cr
e_3 & 0   & 1   & 0 \cr
e_4 & 0   & 0   & 1 \cr
e_5 & b   & c   & d \cr
},}

%\sni
on obtient $E = \gK w_1 \oplus \gK w_2 \oplus \gK w_3$ et l'on a
$\gK^5 = E \oplus F$ avec $F = \gK e_2 \oplus \gK e_5$. 
Puisque 

\snic {
e_1 - w_1 \in F, \quad e_3 - w_2 \in F, \quad e_4 - w_3 \in F, 
}

%\sni
on a $\pi(e_1) = w_1$, $\pi(e_3) = w_2$, $\pi(e_4) = w_3$
et $\pi(e_2) = \pi(e_5) = 0$. Le corps de rationalité de $E$
est $\gK_0 = \gk(a,b,c,d)$, où $\gk$ est le sous-corps premier de $\gK$.

}% fin des solutions d'exos

%:  ---- Section*{references}-----------
\Biblio

Le \iJG lemme de Gauss-Joyal est dans  \cite{Es}, qui lui donne son nom de baptême.
Sur le sujet \gnl de la comparaison entre
les \ids $\rc(f)\rc(g)$ et $\rc(fg)$ on peut consulter \cite{CDLQ02,Glaz2,Nor2} et, dans cet ouvrage, les sections \ref{secLemArtin} et \ref{secThKro}
et la proposition \ref{propLG}.

Concernant le traitement \cof de la \noet on peut consulter
\cite{MRR,JL,Per1,Per2,ric74,sei74b,sei84,Ten}.

L'ensemble de la section \ref{secCramer} se trouve plus ou moins
dans \cite{Nor}.
Par exemple la formule (\ref{eqIGCram}) \paref{eqIGCram}
se trouve sous une forme voisine dans le \tho 5 page~10.
De même notre formule magique à la Cramer (\ref{eqIGCram2}) \paref{eqIGCram2} est très proche du \thoz~6 page 11: Northcott attache une importance centrale à l'équation matricielle $A\,B\,A=A$.
Sur ce sujet, voir aussi \cite{RM} et \cite[Díaz-Toca\&al.]{DiGLQ}.

La proposition \ref{propIGCram2} se trouve dans \cite{Bha} \thoz~5.5.

Concernant le \thrf{theoremIFD}: dans \cite{Nor} le \tho 18
page 122 établit l'\eqvc  des points  \emph{\ref{IFDa}} et \emph{\ref{IFDe}}
par une méthode qui n'est pas entièrement \covz,
mais le \tho 5 page 10
permettrait de donner une formule explicite pour l'implication \emph{\ref{IFDe}}
$\Rightarrow$ \emph{\ref{IFDa}}.

\newpage \thispagestyle{CMcadreseul}
\incrementeexosetprob

%:        %%%%%%%%%%%%%%%%%%%%%%%%%%%%%%%%%%%%
%:        %%%%%%%%%%%%%%%%%%%%%%%%%%%%%%%%%%%%
%---- Chapitre  {coefficients indéterminés}
\chapter{La méthode des coefficients indéterminés}
\label{chapGenerique}
%--------------------
\minitoc

%: Intro
\Intro
\pagestyle{CMExercicesheadings}

\begin{flushright}
{\em Weil Gauss ein echter Prophet der Wissenschaft ist,\\
deshalb reichen die Begriffe,\\
die er aus der Tiefe der Wissencshaft sch\"opft,\\
weit hinaus über den Zweck, \\
zu welchem sie aufgestellt wurden.
}\\
Kronecker\\
Vorlesungen Sommersemester 1891. Leçon 11 \cite{BoSc} \\[1mm]
Trad. approx.\\
{\em Parce que Gauss est un vrai Prophète de la Science,\\
les concepts qu'il puise aux profondeurs de la Science\\
vont au delà du but pour lequel ils ont été établis.
}
\end{flushright}

En 1816 Gauss publie un article fondamental \cite{Gauss} dans lequel il rectifie 
 (sans la citer) la \dem du \tho fondamental de l'\alg donnée par Laplace
 quelques années auparavant.
 La \dem de Laplace est elle-même remarquable en ce qu'elle est 
 \gui{purement \agqz}: elle ne réclame pour les nombres réels que deux 
 \prts très \elrsz: l'existence de la racine carrée d'un nombre 
 $\geq0$ et celle d'un zéro pour un \pol de degré impair.
 
 L'objectif de Gauss est de traiter ce \tho sans faire appel à un corps 
 de nombres imaginaires, hypothétique, sur lequel se décomposerait en facteurs \lins un \pol réel arbitraire. 
 La \dem de Laplace suppose implicitement l'existence d'un tel corps~$\gK$ contenant~$\CC=\RR[i]$, et montre que la \dcn en produit de facteurs \lins
 a lieu en fait dans~$\CC[X]$.
 
 La \dem de Gauss s'affranchit de l'hypothèse du corps~$\gK$
 et constitue
 un tour de force qui montre que l'on peut traiter les choses de manière purement
 formelle. Il prouve l'existence du pgcd de deux \pols par l'\algo d'Euclide
 ainsi que la relation de Bézout correspondante.
 Il démontre que tout \pol \smq s'écrit de manière unique comme
 \pol en les fonctions \smqs \elrs (en introduisant un ordre lexicographique sur les \momsz). Il définit le \discri d'un \polu de manière purement formelle. Il démontre (sans recours aux racines)
 que tout \pol se décompose en produit de \pols de \discri non nul.
 Il démontre (sans recours aux racines) qu'un \pol admet un facteur carré
 \ssi son \discri est nul (il est en \cara nulle). Il fait enfin fonctionner
 la \dem de Laplace de façon purement formelle, sans recours à un \cdrz,
 en utilisant uniquement résultants et \discrisz. 
 
 En bref il établit une \gui{méthode \gnle des \coes indéterminés}
 sur une base ferme, qui sera systématiquement reprise,  notamment
 par Leopold Kronecker, Richard Dedekind, Jules Drach, Ernest Vessiot\ldots

\medskip Dans ce chapitre nous introduisons la méthode des \coes indéterminés et nous en donnons quelques applications.

Nous commençons par quelques \gnts sur les anneaux de \polsz. Le lemme de \DKM et le \tho de \KRO sont deux outils de base qui donnent des informations  précises sur les \coes du produit de deux \polsz. Ces deux résultats seront souvent utilisés dans le reste de l'ouvrage.
 
Nous étudions  les \prts \elrs du discriminant et du résultant et nous introduisons l'outil fondamental qu'est l'\adu d'un \pol unitaire.
Celle-ci permet de simplifier des preuves purement formelles
à la Gauss en donnant un substitut formel au \gui{\cdrz} du \polz. 

Tout ceci est très uniforme et fonctionne avec des anneaux commutatifs
arbitraires. \Llec ne verra apparaître les corps qu'à partir de la section~\ref{secGaloisElr}.

Les applications que nous traitons concernent la théorie de Galois de base,
les premiers pas en théorie algébrique des nombres, et le \nst de Hilbert.
Nous avons \egmt consacré une section à la méthode de Newton en \algz.

%%%%%%%%%%%%%%%%%%%%%%%%%%%%%%%%%%%%%%%%%%%%%%%%%%%%%%%%%%%%%%%%%%%%%%%%%%%
%:    ensembles finis
\subsec{Deux mots sur les ensembles finis}\label{Deux mots}\relax

  Un ensemble~$E$ est dit \ixc{fini}{ensemble ---}
lorsque l'on a explicitement une bijection
entre~$E$ et un segment initial~$\sotq{x\in\NN}{x<n}$ de~$\NN$.
Il est dit \emph{finiment énumérable}
lorsque l'on  a explicitement une surjection d'un ensemble fini~$F$
sur~$E$.%
\index{finiment énumérable!ensemble ---}%
\index{ensemble!fini}%
\index{ensemble!finiment énumérable}
%-% PERSO
\hum{La \dfn \emph{finiment énumérable}
exclut des parties qui hésiteraient entre le vide
et un \eltz, mais elle n'exclut pas des parties
qui hésitent entre 1 et 2 \eltsz.
C'est légèrement ennuyeux, mais éviter cette dissymétrie
conduit à une \dfn qui est un peu dure à avaler.}%
%-% Fin PERSO

\rdb \label{NOTAPfPfe}
En \gnl le contexte est suffisant pour faire
la distinction entre les deux notions. Parfois, on a intérêt
à être très précis. Nous ferons
la distinction si nécessaire en utilisant la notation~$\Pf$ ou~$\Pfe$:
nous noterons~$\Pf(S)$ \emph{l'ensemble des parties finies} de l'ensemble~$S$
et~$\Pfe(S)$ \emph{l'ensemble des parties finiment énumérables}.
En \coma lorsque~$S$ est discret (resp. fini), on a l'\egt
$\Pf(S)=\Pfe(S)$ et c'est un ensemble discret  (resp. fini)%
\footnote{En \coma
on s'abstient en  \gnl de considérer l'\gui{ensemble de toutes les parties
d'un ensemble}, même fini, car ce n'est pas un ensemble \gui{raisonnable}:
il ne semble pas possible de donner une \dfn claire de ses \elts
(voir la discussion \paref{P(X)}).
Quand nous avons utilisé la notation~$\cP_{\ell}$ pour \gui{l'ensemble des parties  de~$\so{1,\ldots ,\ell}$}, \paref{notaAdjalbe}, il s'agissait
en fait de l'ensemble des parties finies de~$\so{1,\ldots ,\ell}$.
}.
Lorsque~$S$ \emph{n'est pas} discret,~$\Pf(S)$ \emph{n'est pas} égal
à~$\Pfe(S)$.\index{ensemble!des parties finies}%
\index{ensemble!des parties finiment énumérées}

Notons aussi que lorsque~$S$ est un ensemble fini toute partie détachable
(cf. \paref{subsecTestDEgalite}) est finie: l'ensemble des parties finies est alors égal à l'ensemble des parties
détachables.

Les parties finiment énumérées sont omniprésentes dans le discours
\mathe usuel. Par exemple lorsque l'on parle d'un \itf on veut dire
un \id engendré par une partie finiment énumérée et non
par une partie finie. 
De même
quand nous parlons d'une \ixx{famille}{finie}~$(a_i)_{i\in I}$
dans l'ensemble~$E$,
nous entendons que~$I$ est un ensemble fini, donc la partie~\hbox{$\sotq{a_i}{i\in I}\subseteq E$} est  finiment énumérée.

Enfin un ensemble non vide~$E$ est dit \emph{énumérable}
si l'on a une application surjective~$\NN\to E$.%
\index{enumerable@énumérable!ensemble ---}

\pagestyle{CMheadings}
%%%%%%%%%%%%%
\section{Anneaux de \pols}
\label{secAnnPols}

%:   subsec    Algorithme de factorisation partielle
\subsec{Algorithme de factorisation partielle}

Nous supposons \llec familier avec l'\algo d'Euclide étendu qui permet
de calculer le pgcd \mon de deux \polus dans~$\KX$ lorsque~$\gK$ est un \cdi
(voir par exemple le \pbz~\ref{exoAnneauEuclidien}).
%:     Lemma{lemPartialDec}
\begin{lemma}\label{lemPartialDec}
Si~$\gK$ est un \cdiz, on dispose d'un \algo de \emph{\fapz}
pour les familles finies de \polus dans~$\KX$: une \fap pour une
famille finie~$(g_1,\ldots,g_r)$ est donnée par une famille finie~$(\lfs)$
de \polus deux à deux étrangers et
par  l'écriture de chaque~$g_i$ sous la forme
$$\preskip.4em \postskip.4em \ndsp
g_i=\prod\nolimits_{k=1}^sf_k^{m_{k,i}}\; (m_{k,i}\in\NN).
$$
%\sni
La famille~$(\lfs)$ s'appelle une \emph{\bdfz}
pour la famille~$(g_1,\ldots,g_r)$.%
\index{algorithme de factorisation partielle}% 
\index{factorisation!partielle}\index{factorisation partielle!base de ---} 
\end{lemma}
%--------- fin lemma ---------------------------------------------- 
%
\begin{proof}
Si les~$g_i$ sont deux à deux étrangers, il n'y a rien à faire.
Sinon, supposons par exemple que~$\pgcd(g_1,g_2)=h_0$, $g_1=h_0h_1$
et~$g_2=h_0h_2$ avec~$\deg(h_0)\geq1$. On remplace la famille~$(g_1,\ldots,g_r)$
par la famille~$(h_0,h_1,h_2,g_3,\ldots,g_r)$. On note que la somme des degrés a diminué. On note aussi que l'on peut supprimer dans la liste les \pols égaux à~$1$, ou les occurrences multiples d'un même \polz.
On termine par \recu sur la somme des degrés.
 Les détails sont laissés \alecz.
\end{proof}
\vspace{-1pt}

%:   subsec    Propriété \uvle
\subsec{Propriété \uvle des anneaux de \polsz}

\vspace{1pt}
Un anneau de \pols $\AXn$ vérifie la \prt \uvle
qui le définit comme \emph{l'anneau commutatif librement engendré par~$\gA$
et~$n$ nouveaux \eltsz.} C'est la \prt décrite au moyen de l'\homo d'évaluation dans les termes suivants.

%:     Proposition{propApolLibre}
\begin{proposition}\label{propApolLibre}
\'Etant donnés deux anneaux commutatifs~$\gA$ et~$\gB$, un \homo~$\rho:\gA\to\gB$ et~$n$ \elts $b_1$, $\ldots$, $b_n\in\gB$ il existe un unique \homo
$\varphi:\AXn=\AuX\to\gB$ qui prolonge~$\rho$ et qui envoie les~$X_i$ sur les~$b_i$.

\vspace{-1.1em}
\Pun{\gA}{j}{\rho}{\AuX}{\varphi}{\gB}{~}{$\varphi(X_i)=b_i$, $i\in\lrbn$. \qquad~}
\end{proposition}

\vspace{-1.1em}
Cet \homo $\varphi$ s'appelle \emph{l'\homo d'\evnz} (des~$X_i$ en les~$b_i$).
Si~$P\in\AuX$ a pour image~$P^{\rho}$ dans~$\BXn$, on obtient l'\egt $\varphi(P)=P^{\rho}(b_1,\ldots,b_n)$. L'\homo d'\evn s'appelle encore une \emph{spécialisation}, et
l'on dit  que~$\varphi(P)$ est obtenu en \emph{spécialisant} les~$X_i$
en les~$b_i$.
Lorsque~$\gA\subseteq\gB$, les \elts  $b_1$, $\ldots$, $b_n\in\gB$ sont dits
\emph{\agqt indépendants sur~$\gA$} si l'\homo d'évaluation correspondant
est injectif.%
\index{algebriquement@\agqt indépendants!elements@\elts --- sur un sous-anneau}%
\index{homomorphisme d'évaluation}%
\index{specialisation@spécialisation}

D'après la proposition \ref{propApolLibre} tout calcul fait dans~$\AuX$ se transfère
dans~$\gB$ au moyen de l'\homo d'évaluation.

Il est clair que~$\Sn$ agit comme groupe d'\autos
de~$\AuX$ par permutation des \idtrsz: $(\sigma,Q)\mapsto Q(X_{\sigma1},\ldots, X_{\sigma n})$.

Le corolaire suivant résulte \imdt de la proposition \ref{propApolLibre}.

%Comme conséquence \imde de la proposition \ref{propApolLibre} on obtient
%le corolaire suivant.

%:     Corollary{propZXnLibre}
\begin{corollary}\label{propZXnLibre}
\'Etant donnés $n$ \elts $b_1$, $\ldots$, $b_n$ dans un anneau commutatif  $\gB$, il existe un unique \homo
$\varphi:\ZZXn\to\gB$ qui  envoie les~$X_i$ sur les~$b_i$.
\end{corollary}

%:  subsec    Identités \agqsz
\subsec{Identités \agqsz}

Une \ida est une \egt entre deux \elts de $\ZZXn$ définis de manière
différente. Elle se transfère automatiquement dans tout anneau commutatif au moyen du corolaire précédent.

Comme l'anneau  $\ZZXn$   a des \prts particulières, il arrive
que des \idas soient plus faciles à démontrer sur  $\ZZXn$
que dans \gui{un anneau~$\gB$ arbitraire}.
En conséquence, si la structure d'un \tho se ramène à une famille d'\idasz, ce qui est très fréquent en \alg commutative, on a souvent
intérêt à utiliser un anneau de \pols à \coes dans~$\ZZ$ en prenant comme \idtrs les \elts pertinents dans l'énoncé du \thoz.

Les \prts des anneaux~$\ZZuX$ qui peuvent s'avérer utiles sont nombreuses.
La première est qu'il s'agit d'un anneau intègre.
 Donc il se plonge dans son corps de fractions~$\QQ(\Xn)$ qui offre toutes les facilités des \cdisz.

La deuxième est qu'il s'agit d'un anneau infini et intègre.
En conséquence, \gui{on peut faire disparaître les cas ennuyeux mais rares}. Un cas est rare quand il correspond à l'annulation d'un \polz~$Q$ non identiquement nul.
Il suffit de vérifier l'\egt correspondant à l'\ida lorsque celle-ci est
évaluée pour les points de~$\ZZ^n$ qui n'annulent pas~$Q$.
En effet, si l'\ida à démontrer est~$P=0$, 
on obtient que le \polz~$PQ$ définit la
fonction identiquement nulle sur~$\ZZ^n$, ceci implique~$PQ=0$ et donc~$P=0$
puisque~$Q\neq0$ et~$\ZZuX$ est intègre.
Ceci est parfois appelé le \gui{principe
de prolongement des \idasz}.%
\index{principe de prolongement des \idasz}

D'autres \prts remarquables de~$\ZZ[\uX]$ pourront parfois être utilisées,
comme le fait que c'est un anneau factoriel, \noe \coh \fdi de \ddk finie.

%:   subsec   Exemple d'application
\penalty-2500
\subsubsection*{Un exemple d'application}

%:     Lemma{lemPrincipeIdentitesAlgebriques}
\begin{lemma}\label{lemPrincipeIdentitesAlgebriques}
Pour $A$, $B \in \Mn(\gA)$, on a les résultats suivants.
\begin{enumerate}
\item $\wi {AB} = \wi{B}\wi{A}$.
\item $\rC{AB} = \rC{BA}$. %(\egt des \polcarsz),
\item  $\wi {PAP^{-1}} = P\wi {A}P^{-1}$ pour $P \in \GLn(\gA)$.
\item  $\wi{\wi{A}} = \det(A)^{n-2} A$ si $n \ge 2$.
\item  \emph{(\Tho de Cayley-Hamilton)} $\rC A(A) = 0$.
\item
  Si $\Gamma_A(X):=(-1)^{n+1}\big(\rC A(X)-\rC A(0)\big)\sur{X}$,
on a $\wi{A} = \Gamma_A(A)$ ($n\geq2$).\\
On a aussi
$\Tr \big(\wi{A}\big)=(-1)^{n+1}\Gamma_A(0)$.
\item  \emph{(Identités de Sylvester)} Soient  $r\geq 1$, $s\geq 2$ avec $n=r+s$. Soient les matrices $C\in\MM_r(\gA)$, $F\in\MM_s(\gA)$,
$D\in\MM_{r,s}(\gA)$ et~$E\in\MM_{s,r}(\gA)$ extraites de la matrice $A$
comme ci-dessous
$$\preskip.4em \postskip.4em
\;\;\;A=\blocs{1}{.6}{1}{.6}{$C$}{$D$}{$E$}{$F$}\;.
$$
Notons $\alpha_i=\so{1,\ldots,r,r+i}$ et $\mu_{i,j}=\det(A_{\alpha_i,\alpha_j})$
pour~$i,j\in\lrb{1..s}$. Alors:
$$\preskip-.4em \postskip.4em
\det(C)^{s-1}\det(A)=\det\big((\mu_{i,j})_{i,j\in\lrbs} \big).
$$
\item  Si $\det A = 0$, alors $\Al{2}\wi A=0$.
\end{enumerate}%
\index{Cayley-Hamilton}%
\index{Sylvester!identités de ---}
\end{lemma}
\begin{proof}
 On peut prendre toutes les matrices à \coes indéterminés
sur~$\ZZ$ et localiser en~$\det P$.
Dans ce cas~$A$, $B$, $C$ et~$P$ sont \ivs dans le corps des
fractions de l'anneau~$\gB=\ZZ[(a_{ij}),(b_{ij}),(p_{ij})]$.
Par ailleurs, la matrice~$\wi A$
vérifie l'\egt \smash{$\wi A A=\det(A)\,\In$}, ce qui la caractérise
puisque~$\det A$ est \ivz. Ceci fournit
le point \emph{1} via l'\egt $\det(AB)=\det (A)\det (B)$, les points
\emph{3} et~\emph{4},
 et le point \emph{6}
via le point \emph{5} et l'\egt $\rC A(0)=(-1)^n\det A$.
Pour le point~\emph{2} on note que~$AB=A(BA)A^{-1}$. \\
Pour le \tho de Cayley-Hamilton, on traite d'abord
le cas de la \emph{matrice compagne
d'un \poluz} $f=T^n-\sum_{k=1}^na_kT^{n-k}$:  \label{matrice.compagne}\relax
\index{matrice!compagne d'un \polz}%
\index{compagne!matrice --- d'un \polz}
$$
P =\cmatrix{
0&\cdots&\cdots&\cdots&0&a_n\cr
1&0&&&\vdots&a_{n-1}\cr
0&\ddots&\ddots&&\vdots&\vdots\cr
\vdots&\ddots&\ddots&\ddots&\vdots&\vdots\cr
\vdots& & \ddots&1&0&a_2\cr
0&\cdots&\cdots&0&1&a_1
}.
$$
Il s'agit de la matrice de %l'\Ali 
la \gui{multiplication par $t$}, $\mu_t:y\mapsto ty$  (où $t$ est la classe de~$T$) dans l'anneau quotient
$\aqo{\AT}{f(T)}=\gA[t]$,
exprimée sur la base des \moms ordonnés par degrés croissants.
En effet, d'une part un calcul direct montre que~$\rC P(T)=f(T)$. D'autre
part~$f(\mu_t)=\mu_{f(t)}=0$, donc~$f(P)=0$.\\
Par ailleurs,
dans le cas de la matrice \gnqz, le \deter de la famille~$(e_1, Ae_1, \ldots,
A^{n-1}e_1)$ est \ncrt non nul, donc la matrice \gnq est semblable
à la matrice compagne de son \polcar sur le corps de fractions
de~$\ZZ[(a_{ij})]$.\\
\emph{7.} Puisque $C$ est \ivz, on peut utiliser le pivot de Gauss \gnez,
par multiplication à gauche par une matrice
\smashbot{$\blocs{.8}{.4}{.8}{.4}{$C^{-1}$}{$0$}{$E'$}{$\I_s$}\,,$}
ceci nous ramène au cas où~$C=\I_r$ et~$E=0$. \\
Enfin le point \emph{8} résulte de l'identité de Sylvester (point \emph{7})
avec~$s=2$.
\end{proof}

\rem Le point \emph{3} permet de définir l'\emph{\endo cotransposé}
d'un \endo d'un module libre de rang fini, à partir de
la matrice cotransposée.%
\index{cotransposé!endomorphisme ---}

%  subsec     %%%%%%%%%%%%
\subsubsection*{Poids, \pols \hmgs}

On dit que l'on a défini un poids sur une \alg de \pols $\AXk$
lorsque l'on attribue à chaque \idtr $X_i$ un poids~$w(X_i)\in\NN$.
On définit ensuite le poids du \mom $\uX^{\underline{m}}=X_1^{m_1}\cdots X_k^{m_k}$ par

\snic{w(\uX^{\underline{m}})=\sum_im_iw(X_i),}

%\sni
de sorte que~$w(\uX^{\underline{m}+\underline{m'}})=
w(\uX^{\underline{m}})+w(\uX^{\underline{m'}})$. Le degré d'un \polz~$P$
pour ce poids, noté en \gnl $w(P)$, est le plus grand des poids des \moms
apparaissant avec un \coe non nul.
Ceci n'est bien défini que si l'on dispose d'un test
d'\egt à~$0$ dans~$\gA$. Dans le cas contraire on se contente de définir la phrase
\gui{$w(P)\leq r$}.

Un \pol est dit \ixc{homogène}{polynome@\pol ---} (pour un poids~$w$) si tous ses \moms ont même poids.

Lorsque l'on dispose d'une \ida et d'un poids, chaque composante \hmg de l'\ida
fournit une \ida particulière.

On peut aussi définir des poids à valeurs dans des \mos ordonnés plus compliqués
que~$(\NN,0,+,\geq)$. On demande alors que ce \mo soit
la partie positive d'un produit de groupes abéliens totalement ordonnés,
ou plus \gnlt un \mo à pgcd (cette notion sera introduite au chapitre~\ref{chapTrdi}).

%%%%%%%%%%%%%%%%%%%%%%%%%%%%%%%%%%%%%%%%%%%%%%%%%%%%%%%%%%%%%%%%%%%%%%%%%%%
%:subsec   Polynomes symetriques
\subsec{Polynômes \smqs}

On fixe $n$ et $\gA$ et l'on note  $S_1$, $\ldots$, $S_n$ les \emph{\pols \smqs \elrs
en les~$X_i$}  dans~$\AXn$. Ils sont définis par l'\egt
$$
T^n+S_1T^{n-1}+S_2T^{n-2}+\cdots+S_n=\prod\nolimits_{i=1}^n(T+X_i).
$$
On a $S_1=\sum_iX_i$, $S_n=\prod_iX_i$,
$S_k=\sum_{J\in\cP_{k,n}}\prod_{i\in J}X_i$.
Rappelons le \tho bien connu suivant (une \dem est suggérée en
exercice~\ref{exothSymEl}).%
\index{polynome@\pol!symétrique \elrz}

%:     Theorem{thSymEl}
\begin{theorem}\label{thSymEl} \emph{(\Pols \smqs \elrsz)}
\begin{enumerate}
\item Un \polz~$Q\in\AXn=\AuX$, invariant par les permutations de variables, s'écrit de
manière unique comme un \pol en les fonctions \smqs \elrs $S_1$, \ldots, $S_n$. En d'autres termes
\begin{itemize}
  \item le sous-anneau des points fixes de $\AuX$ par l'action du groupe
  \smq $\Sn$
est l'anneau $\gA[S_1,\ldots,S_n]$ engendré par~$\gA$
et les~$S_i$, et
  \item les $S_i$ sont \agqt indépendants sur $\gA$.
\end{itemize}
\item Notons $d(P)$ le degré total de $P\in\AuX$ lorsque chaque $X_i$ est affecté du poids~$1$, et~$d_1(P)$ son degré en~$X_1$.
Notons~$\delta(Q)$ le degré total de~$Q\in\gA[S_1,\ldots,S_n]$ lorsque chaque variable~$S_i$ est affectée du poids~$i$ et~$\delta_1(Q)$ son degré total  lorsque chaque variable~$S_i$ est affectée du poids~$1$.
Supposons que~$Q(S_1,\ldots,S_n)$ s'évalue en~$P(\uX)$.
\begin{itemize}
\item [a.] $d(P)=\delta(Q)$, et si $Q$ est $\delta$-\hmgz, alors $P$ est $d$-\hmgz.
%
%\item [(b)] Donc dans tous les cas $d(P)=\delta(Q)$.
%
\item [b.] % Dans tous les cas
$d_1(P)=\delta_1(Q)$.
\end{itemize}

\item $\AXn$ est un module libre de rang $n!$ sur $\gA[S_1,\ldots,S_n]$ et une base est formée
par les \moms $X_1^{k_1}\cdots X_{n-1}^{k_{n-1}}$ tels que~$k_i\in\lrb{0.. n-i}$ pour chaque~$i$.
\end{enumerate}
\end{theorem}

%:     Corollary{corthSymE1}
\begin{corollary}\label{corthSymE1}
Sur un anneau $\gA$ on considère le \pol \gnq %de degré $n$

\snic{f=T^n+f_1T^{n-1}+f_2T^{n-2}+\cdots+f_n,}

%\sni
où les~$f_i$ sont des \idtrsz.
\\
On a un \homo injectif $j:\gA[f_1,\ldots,f_n]\to\AXn$ tel que les~$(-1)^{k}j(f_k)$
sont les \pols \smqs \elrs en les~$X_i$.
\end{corollary}

En bref on peut toujours se ramener au cas où~$f(T)=\prod_i(T-X_i)$, où les~$X_i$ sont d'autres \idtrsz.

%:     Corollary{corthSymE2}
\begin{corollary}\label{corthSymE2}
Sur un anneau $\gA$ on considère le \pol \gnq %de degré $n$
$$\preskip.4em \postskip.2em 
f=f_0T^n+f_1T^{n-1}+f_2T^{n-2}+\cdots+f_n, 
$$
où les~$f_i$ sont des \idtrsz. 
\\
On a un \homo injectif $j:\gA[f_0,\ldots,f_n]\to\gB{=}\gA[F_0,\Xn]$, avec
dans~$\gB[T]$ l'\egt suivante.
$$\preskip.2em \postskip.4em 
\dsp j(f_0)\,T^n+j(f_1)\,T^{n-1}+\cdots+j(f_n)=F_0\,\prod\nolimits_i\,(T-X_i)\,. 
$$
\end{corollary}

En bref, on peut toujours se ramener au cas où~$f(T)=f_0\prod_i(T-X_i)$, 
avec des \idtrs $f_0$, $X_1$, \dots, $X_n$.

\begin{proof}
Il suffit de voir que si $f_0$, $g_1$, $\ldots$, 
$g_n\in\gB%\supseteq\gA
$ sont \agqt
indépen\-dants sur~$\gA$, alors il en va de même pour $f_0$, $f_0g_1$, $\ldots$, $f_0g_n$.
Il suffit de vérifier que $f_0g_1$, $\ldots$, $f_0g_n$ sont \agqt
indépendants sur~$\gA[f_0]$. Cela résulte de ce que
$f_0$ est \ndz et de ce que $g_1$, $\ldots$, $g_n$ sont \agqt
indépendants sur~$\gA[f_0]$.
\end{proof}
%

%%%%%%%%%%%%%
\section{Lemme de Dedekind-Mertens}
\label{secLemArtin}

Rappelons que pour un  \polz~$f$ de $\AXn=\AuX$, on appelle
\gui{contenu de~$f$} et l'on note~$\rc_{\gA,\uX}(f)$ ou~$\rc(f)$
 l'\id engendré par les \coes de~$f$.

Notons que l'on a toujours $\rc(f)\rc(g){\supseteq}\rc(fg)$
et donc $\rc(f)^{k+1}\rc(g){\supseteq}\rc(f)^k\rc(fg)$ pour tout~$k\geq0$.
Pour~$k$ assez grand cette inclusion devient une \egtz.

%:     Lemma{lemdArtin}
\CMnewtheorem{lemDKME}{Lemme de \DKMz}{\itshape}

\begin{lemDKME}\label{lemdArtin}% 
\index{Lemme de Dedekind-Mertens} ~\\
Pour $f,g\in\AT$ avec $m\geq\deg g$ on a

\snic{\rc(f)^{m+1}\rc(g)=\rc(f)^m\rc(fg).}
\end{lemDKME}
\begin{proof}
Tout d'abord on remarque que les produits $f_ig_j$ sont les \coes du
\polz~$f(Y)g(X)$. Pareillement, pour des \idtrs
$Y_0,\ldots ,Y_m$, le contenu du \polz~$f(Y_0)\cdots f(Y_m)g(X)$ est égal
à~$\rc(f)^{m+1}\rc(g)$.
\\
 Notons~$h=fg$.
Imaginons que dans l'anneau~$\gB=\gA[X,Y_0,\ldots ,Y_m]$ on puisse montrer
l'appartenance du \pol  $f(Y_0)\cdots f(Y_m)g(X)$ à l'idéal

\snic{\som_{j=0}^m \big(h(Y_j)\,\prod\nolimits_{k,k\neq j} \gen{f(Y_k)}\big).}

%\sni
On en déduirait immédiatement que $\rc(f)^{m+1}\rc(g)\subseteq \rc(f)^m\rc(h)$.
\\
\`A quelque chose près, c'est ce qui va arriver.
On chasse les \denos dans la formule d'interpolation de Lagrange
(on a besoin d'au moins~$1+\deg g$ points d'interpolation):

\snic{g(X)=\sum_{j=0}^m\; \; \fraC{\prod\nolimits_{k,k\neq j}(X-Y_k)}
{\prod\nolimits_{k,k\neq j}(Y_j-Y_k)}\;\,g(Y_j)\,.}

%\sni
On obtient dans l'anneau $\gB$, en posant
$\Delta=\prod_{j\neq k}(Y_j-Y_k)$:

\snic{\Delta \cdot g(X)\;\in \;\som_{j=0}^m \gen{g(Y_j)}.}

%\sni
Donc en multipliant par  $f(Y_0)\cdots f(Y_m)$:

\snic{\Delta \cdot f(Y_0)\cdots f(Y_m)\cdot g(X)\;\in\; \som_{j=0}^m h(Y_j)\,
\prod\nolimits_{k,k\neq j} \gen{f(Y_k)}.
}

%\sni
Si l'on montre que pour n'importe quel {\mathrigid 4mu $Q\in \gB$ on a
$\rc(Q)=\rc(\Delta \cdot Q)$}, l'appartenance précédente donne 
$\rc(f)^{m+1}\rc(g)\subseteq \rc(f)^m\rc(h).$
\\
 On note que $\rc(Y_i\,Q)=\rc(Q)$ et surtout que
%-----------------begin $$----------------
$$ \rc\big(Q(Y_0\pm Y_1,Y_1,\ldots ,Y_m)\big)\subseteq \rc\big(Q(Y_0,Y_1,\ldots ,Y_m)\big).
$$
%-----------------end $$------------------
Donc, en faisant $Y_0=(Y_0\pm Y_1)\mp Y_1$, 
$\rc\big(Q(Y_0\pm Y_1,Y_1,\ldots,Y_m)\big)= \rc(Q).$  
 Les \pols suivants ont donc tous même contenu:

\snac{Q,\;Q(Y_0+ Y_1,Y_1,\ldots ,Y_m),\;Y_0\,Q(Y_0+ Y_1,Y_1,\ldots
,Y_m),\;(Y_0-Y_1)\,Q(Y_0,Y_1,\ldots ,Y_m).}

%\sni
D'où ensuite  $\rc(Q)=\rc(\Delta \cdot Q)$.
\end{proof}

On en déduit les corolaires suivants.

%:     corollary} \label {corLDM}
\begin {corollary} \label {corLDM}
Si $f_1$, $\ldots$, $f_d$ sont $d$ \pols (à une \idtrz) 
de degré~$\le \delta$,
on a, avec~$e_i = 1 + (d-i)\delta$:
$$\preskip.3em \postskip.4em 
\rc(f_1)^{e_1} \rc(f_2)^{e_2} \cdots \rc(f_d)^{e_d}  \subseteq \rc(f_1f_2\cdots f_d). 
$$
\end {corollary}
%---------  end corollary ---------------------------------------------

%
\begin {proof}
Soient $f = f_1$ et $g = f_2 \cdots f_d$. On a $\deg g \le (d-1)\delta$ et~$e_1 = 1 + (d-1)\delta$. Le lemme de \DKM donne donc:

\snac {\rc(f)^{e_1} \rc(g) = \rc(f)^{(d-1)\delta} \rc(fg) \subseteq \rc(fg), $ i.e. $
\rc(f_1)^{e_1} \rc(f_2 \cdots f_d) \subseteq \rc(f_1f_2\cdots f_d).
}

%\sni
On termine par \recu sur $d$.
\end {proof}
%

%:     Corollary{corlemdArtin}
\begin{corollary}\label{corlemdArtin}  Soient $f$ et $g\in \AT.$
\begin{enumerate}
\item \label{i1corlemdArtin} Si $\Ann_\gA\big(\rc(f)\big)=0$, alors $\Ann_\AT(f)=0$ \emph{(Lemme de McCoy)}.\index{Lemme de McCoy}
\index{McCoy!lemme de ---}
\item \label{i2corlemdArtin} Si $\gA$ est réduit, alors $\Ann_\AT(f)= \Ann_\gA\big(\rc(f)\big)[T]$.
\item \label{i3corlemdArtin} Le \polz~$f$ est nilpotent \ssi chacun de ses \coes est nilpotent.
\item Si $\rc(f)=1$, alors $\rc(fg)=\rc(g)$.
\end{enumerate}
\end{corollary}
\begin{proof}
Soit $g\in \Ann_\AT(f)$  et $m \ge \deg(g)$. Le lemme de \DKM  implique:
$$\preskip-.4em \postskip.3em 
\qquad\qquad\rc(f)^{1+m}g=0.\eqno (*) 
$$

%\sni
\emph{1.}
Donc,  $\Ann_\gA \rc(f)=0$  implique  $g=0$.
\\
\emph{2.}
Puisque l'anneau est réduit, $(*)$ implique $\rc(f)g=0$. Ainsi tout \polz~$g$ annulé par
$f$ est annulé par~$\rc(f)$.\\
Par ailleurs, $\Ann_\gA\big(\rc(f)\big)=\gA\cap \Ann_\AT(f)$
et donc  l'inclusion
$$\preskip.3em \postskip.3em 
\Ann_\AT(f)\supseteq \Ann_\gA\big(\rc(f)\big)[T] 
$$
est toujours vraie (que $\gA$ soit réduit ou non).\\
\emph{3.}  Si $f^2=0$, le lemme de \DKM implique $\rc(f)^{2+\deg f}=0$.\\
\emph{4.} Immédiat d'après $\rc(f)^{m+1}\rc(g)=\rc(f)^m\rc(fg)$.
\end{proof}
%

%%%%%%%%%%%%%
\section{Un \tho de Kronecker}
\label{secThKro}

\vspace{3pt}
%:  subsec    Algèbres et \elts entiers
\subsec{Algèbres et \elts entiers}

Nous introduisons tout d'abord
la terminologie des \Algsz. Les \algs que nous considérons dans cet ouvrage sont
associatives, commutatives et unitaires, sauf précision contraire.

%:    Definition{def0Alg}---------
\begin{definition}
\label{def0Alg} ~
%-----------------begin enum------------------
\begin{enumerate}
\item Une \emph{\Algz} est un anneau commutatif $\gB$ avec un \homo d'anneaux
commutatifs $\rho:\gA\to\gB$. Cela fait de~$\gB$ un \Amoz. Lorsque
$\gA\subseteq\gB$, ou plus \gnlt si~$\rho$ est injectif, on dira que~$\gB$ est une \ix{extension} de~$\gA$.%
\index{algèbre!sur un anneau}

\item Un \emph{morphisme} de l'\Alg $\gA\vers{\rho}\gB$ vers \hbox{l'\Alg $\gA\vers{\rho'}\gB'$} est un \homo d'anneaux~$\gB\vers{\varphi}\gB'$ vérifiant
$\varphi\circ\rho=\rho'$. L'ensemble des morphismes d'\Algs de~$\gB$ vers~$\gB'$ sera noté~$\Hom_\gA(\gB,\gB')$.

\vspace{-1.2em}
\tri{\gA}{\rho}{\rho'}{\gB}{\varphi}{\gB'}{}
\end{enumerate}
%-----------------end enum------------------
\end{definition}
%--- end-definition------------------------------------

\vspace{-.5em}
\rems ~\\
1) Nous n'avons pas voulu réserver la terminologie \gui{extension} au cas des
corps. Ceci nous obligera par la suite à utiliser dans le cas des corps des
phrases comme: $\gL$ est une extension de corps de $\gK$, ou: $\gL$ est un corps,
extension de~$\gK$.

2) Tout anneau est une \ZZlg de manière unique et tout \homo d'anneaux est un morphisme des \ZZlgs correspondantes. La catégorie des anneaux
commutatifs peut donc être vue comme un cas particulier parmi les catégories
d'\algs définies ci-dessus.
\eoe \rdb

\mni
{\bf Notation:} Si $b\in\gB$ et $M$ est un \Bmoz, on note $\mu_{M,b}$ ou $\mu_b$ la multiplication par~$b$ dans~$M$: $y\mapsto by,\, M\to M$. Ceci peut être vu comme une \Bliz, ou, si~$\gB$ est une \Algz, comme une \Ali pour la structure de
\Amo de~$M$. \label{NOTAmux}

%  subsec     %%%%%%%%%%%%
%\subsubsection*{\'Eléments entiers}
%:    Definition{defEntierAnn0}--------
\begin{definition}
\label{defEntierAnn0} Soit  $\gA\subseteq\gB$ des anneaux.
\begin{enumerate}
\item Un \elt $x\in\gB$  est dit \ixc{entier}{element@élément --- sur un anneau} sur~$\gA$ s'il existe un entier~$k\geq 1$ tel \hbox{que $x^k=a_1x^{k-1}+a_2 x^{k-2}+\cdots +a_k$}
avec  les~$a_h\in\gA$. Si~$\gA$ est un \cdiz, on dit aussi que~$x$
est \emph{algébrique sur~$\gA$}.\index{algebrique@\agq!element@\elt --- sur un corps discret}
\item Dans ce cas, le \polu $P=X^k-(a_1X^{k-1}+a_2 X^{k-2}+\cdots +a_k)$ est appelé une \emph{\rdiz} de  $x$ sur~$\gA$.
En fait, par abus de langage on dit aussi que l'\egt $P(x)=0$ est une  \emph{\rdiz}. Si~$\gA$ est un \cdiz,~on parle aussi de \emph{\rde \agqz}.%
\index{relation de dépendance!intégrale}%
\index{relation de dépendance!algébrique}
\item L'anneau $\gB$ est dit \emph{entier} sur $\gA$ si tout \elt de~$\gB$ est entier
sur~$\gA$. On dira aussi que l'\Alg $\gB$ est \emph{entière}.
Si~$\gA$ et~$\gB$ sont des \cdisz,~on dit que~$\gB$ est  \emph{\agqz} sur~$\gA$.%
\index{algèbre!entière}%
\index{anneau!entier sur un sous-anneau}%
\index{entier!anneau --- sur un sous-anneau}%
\index{algebrique@\agq!corps --- sur un sous-corps}
\item
Si  $\rho:\gC\to\gB$ est une \Clg avec~$\rho(\gC)=\gA$, on dira
que l'\algz~$\gB$ est entière sur~$\gC$ si elle est entière sur~$\gA$.
\end{enumerate}
\end{definition}
%--- end-definition-------------------------

%:  subsec    Le \tho
\penalty-2500
\subsec{Le \tho}

%:     Theorem{thKro}
\begin{theorem}\label{thKro} \emph{(\Tho de \KROz) \cite{Kro1}}\\
Soit dans $\gB[T]$ les \pols
%{\small
$$\preskip.3em \postskip.3em
f=\sum_{i=0}^n(-1)^if_iT^{n-i},\;  g=\sum_{j=0}^m(-1)^j g_j T^{m-j} \; \mathit{et} \; h=fg=\sum_{r=0}^{p}(-1)^r h_rT^{p-r},
$$
%}
où $p=m+n$. Soit $\gA=\gZ[h_0,\ldots,h_p]$ le sous-anneau
engendré par les \coes de $h$ ($\gZ$ est le sous-anneau de $\gB$ engendré par $1_\gB$). 
\begin{enumerate}
\item Chaque $f_ig_j$ est entier sur $\gA$.
\item Dans le cas où on prend pour $f_i$ et $g_j$ des \idtrs
sur l'anneau $\ZZ$,
on trouve une \rdi sur $\gA$ pour $z_{i,j}=f_ig_j$ qui est \hmg
pour différents systèmes de poids attribués aux \momsz:
\begin{enumerate}
\item les poids respectifs de $z_{k,\ell} $ et $h_r$ sont $k+\ell$ et $ r$.
\item les poids respectifs de $z_{k,\ell} $ et $h_r$ sont $p-k-\ell$ et $p-r$.
\item les poids de $z_{k,\ell} $ et $h_r$ sont  $w(z_{k,\ell})=w(h_r)=1$.
\end{enumerate}
Naturellement ces \rdis s'appliquent ensuite dans tout anneau.
\end{enumerate}
\end{theorem}
\begin{proof}
Il suffit de traiter le point \emph{2.}
\\
 Voyons d'abord un cas \gnq intermédiaire.
Nous prenons $f_0=g_0=1$ et pour les autres $f_i$ et $g_j$
des \idtrs sur $\ZZ$. Les \pols $f$ et $g$ sont donc des \polus
dans $\gB[T]$ avec $\gB=\ZZ[f_1,\ldots,f_n,g_1,\ldots,g_m]$,
et $\gA=\ZZ[h_1,\ldots,h_p]$.\\
On suppose \spdg \hbox{que $\gB\subseteq\gC=\ZZ[\xn,\alb\ym]$}, où les $x_i$ et $y_j=x_{n+j}$ sont des \idtrsz,
les $f_i$ sont les \pols \smqs \elrs en les $x_i$, et les $g_j$ sont les \pols \smqs \elrs en les $y_j$ (appliquer 
deux fois le corolaire~\ref{corthSymE1}).
Si nous attribuons à $x_i$ et $y_j$ le poids 1, les $z_{k,\ell}$ et $ h_r$ sont \hmgs et obtiennent
les poids décrits en \emph{2a}.
Pour calculer une \rdi pour $f_ig_j$ (avec éventuellement $i$ ou $j=0$)
sur~$\gA$, on considère le sous-groupe $H_{i,j}$ \hbox{de $\rS_p$} formé par les
 $\sigma$  qui vérifient $\sigma(f_ig_j)=f_ig_j$
(ce sous-groupe contient au moins toutes les
permutations qui stabilisent $\lrbn$).
On considère alors le \pol
$$
P_{i,j}(T)=\prod\nolimits_{\tau\in\rS_p/H_{i,j}}\big(T-\tau(f_ig_j)\big),\eqno(*)
$$
où $\tau\in\rS_p/H_{i,j}$ signifie que l'on prend exactement un $\tau$ dans chaque classe
à gauche modulo $H_{i,j}$. Alors, $P_{i,j}
$ est \hmg pour les poids $w_a$ décrits \hbox{en \emph{2a}} ($i,j$ étant fixés, on note $w_a$ les poids
\emph{2a}, avec $w_a(T)=w_a(z_{i,j})$).
En outre,~$P_{i,j}$ est \smq en les $x_k$ ($k\in\lrbp$). Il s'écrit
donc de manière unique
comme un \polz~$Q_{i,j}(\uh,T)$ en les $h_r$ et $T$, et $Q_{i,j}$ 
 est $w_a$-\hmg (\tho \ref{thSymEl} points \emph{1} et \emph{2a}).
Le degré en $T$ de $Q_{i,j}$ \hbox{est $d_{i,j}=(\rS_p:H_{i,j})$}.
Pour $R\in\gC[T]$, notons $\delta(R)$ l'entier $\deg_{x_1}(R)+\deg_T(R)$. On voit que~$\delta$ est un poids, et que $\delta(f_ig_j)=w(f_ig_j)\leq1$,
$\delta(h_r)=w(h_r)\leq1$ (\hbox{avec $w(h_r)=1$} si $i$, $j$, $r\geq1$). 
En outre, chaque facteur
de $P_{i,j}$ dans~$(*)$ est de poids~$1$ (mais pas forcément \hmg
car on peut avoir $\delta(\sigma\big(f_ig_j)\big)=0$).
Ceci donne $\delta(Q_{i,j})=d_{i,j}$ lorsque le  \pol
est évalué dans $\gC[T]$. En outre, d'après le \tho \ref{thSymEl} point
\emph{2b}, lorsque l'on écrit un \pol \smq  en $(\xp)$, disons $S(\ux)$, comme \polz~$S_1(\uh)$ en les $h_i$,
\hbox{on a $\delta(S)=w(S_1)$}.
Ainsi  $w(Q_{i,j})=d_{i,j}$.
\\
Pour traiter le point \emph{2} proprement dit il suffit d'\gui{homogénéiser}.
En effet, si l'on pose $\wi f_i=f_i/f_0$
et $\wi g_j=g_j/g_0$, ce qui est légitime d'après 
le corolaire~\ref{corthSymE2}, on retombe pour les $\wi f_i$ et $\wi g_j$ sur
la situation précédente pour ce qui concerne les poids \emph{2a}.
On obtient une \rdi \hmg pour $\wi z_{i,j}=\wi f_i\wi g_j$ sur le sous-anneau
engendré par les~$\wi h_{r}$:
$$\preskip.3em \postskip.3em
Q_{i,j}(\wi h_{1},\ldots,\wi h_{p},\wi z_{i,j})=0,
$$
avec $\wi z_{i,j}=f_ig_j/h_0$ et $\wi h_{r}=h_{r}/h_0$.
\\
On multiplie l'\ida
obtenue par $h_0^{d_{i,j}}$ de manière à obtenir un
\polu en $z_{i,j}$.
\\
Tous les \denos ont disparu parce que
$w(Q_{i,j})=d_{i,j}$. On obtient:
$$\preskip.4em \postskip.4em
R_{i,j}(h_0,\ldots,h_p,f_ig_j)=0,
$$
où $R_{i,j}(h_0,\ldots,h_p,T)$ est unitaire en $T$ et homogène pour les poids $w_a$ et~$w$.
\\
Reste la question de l'homogénéité pour les poids $w_b$ en \emph{2b}:
il suffit de noter que l'on a pour tout $R\in\AT$ l'\egt
%
%\snic{
$w_a(R)+w_b(R)=pw(R).$
%}
\end{proof}
\exl Dans le cas
où $m=n=2$, le calcul indiqué donne les résultats suivants.
\\
Lorsque $f_0=g_0=1$ le \coe $g_1$ annule le \polz:
$$\arraycolsep2pt\begin{array}{rcl}
p_{01}(t)&=& t^6 - 3 h_1 t^5 + (3 h_1^2 + 2 h_2) \,t^4 + (-h_1^3 - 4 h_1 h_2) \,t^3 \\[1mm]
&&+(2 h_1^2 h_2 + h_1 h_3 + h_2^2 - 4 h_4) \,t^2 + (-h_1^2 h_3 - h_1 h_2^2 + 4 h_1 h_4) \,t
\\[1mm]
&&-h_1^2 h_4 + h_1 h_2 h_3 - h_3^2\,,
\end{array}$$
donc dans le cas \gnl $f_0g_1$
annule le \polz:
$$\arraycolsep2pt\begin{array}{rcl}
q_{01}(t)&=& t^6 - 3 h_1 t^5 + (3 h_1^2 + 2h_0 h_2) \,t^4 + (-h_1^3 - 4 h_0 h_1 h_2) \,t^3 \\[1mm]
&&+(2 h_0 h_1^2 h_2 + h_0^2 h_1 h_3 +h_0^2  h_2^2 - 4h_0^3 h_4) \,t^2 
\\[1mm]
&&+(-h_0^2 h_1^2 h_3 - h_0^2 h_1 h_2^2 + 4 h_0^3h_1 h_4) \,t
-h_0^3h_1^2 h_4 + h_0^3h_1 h_2 h_3 - h_0^4h_3^2\,.
\end{array}$$
Lorsque $f_0=g_0=1$ le \coe $g_2$  annule le \polz:
$$\arraycolsep2pt\begin{array}{rcl}
p_{02}(t)      & =& t^6 - h_2 t^5 + (h_1 h_3 - h_4) \,t^4 + (-h_1^2 h_4 + 2 h_2 h_4 - h_3^2) \,t^3   \\[1mm]
      &  &+ (h_1 h_3 h_4 - h_4^2) \,t^2 - h_2 h_4^2 \,t + h_4^3
\,,\end{array}$$
donc  $f_0g_2$
annule le \polz:
$$\arraycolsep2pt\begin{array}{rcl}
q_{02}(t)      & =& t^6 - h_2 t^5 + (h_1 h_3 - h_0 h_4) \,t^4 + (-h_1^2 h_4 + 2 h_0 h_2 h_4 - h_0 h_3^2) \,t^3   \\[1mm]
      &  & +(h_0 h_1 h_3 h_4 - h_0^2 h_4^2) \,t^2 - h_0^2 h_2 h_4^2 \,t + h_0^3 h_4^3\,.
\end{array}$$
Lorsque $f_0=g_0=1$ le \coe
$f_1  g_1$ annule le \polz:
$$\arraycolsep2pt\begin{array}{rcl}
p_{11}(t)&=&t^3 - 2 h_2 t^2 + (h_1 h_3 + h_2^2 - 4 h_4) \,t +
      h_1^2 h_4 - h_1 h_2 h_3 +
h_3^2\,.
\end{array}
$$
Lorsque $f_0=g_0=1$ le \coe $f_1  g_2$ annule le \polz:
$$\arraycolsep2pt\begin{array}{rcl}
p_{12}(t)&=&t^6 - 3 h_3 t^5 + (2 h_2 h_4 + 3 h_3^2) \,t^4 + (-4 h_2 h_3 h_4 - h_3^3) \,t^3
 \\[1mm]
      &  &+(h_1 h_3 h_4^2 + h_2^2 h_4^2 + 2 h_2 h_3^2 h_4 - 4 h_4^3) \,t^2  \\[1mm]
      &  &+(-h_1 h_3^2 h_4^2 -
h_2^2 h_3 h_4^2 + 4 h_3 h_4^3) \,t - h_1^2 h_4^4 + h_1 h_2 h_3 h_4^3 - h_3^2 h_4^3\,.
\end{array}$$

%:     Corollary{corthKro}
\begin{corollary}\label{corthKro}
 \emph{(\Tho de \KRO en plusieurs variables) }\\
Soit dans $\gB[\Xk]$ les \pols

\snic{f=\som_{\alpha}f_\alpha X^{\alpha},\;\;  g=\som_{\beta}b_\beta X^{\beta} \;\;\hbox{ et }\;\; h=fg=\som_{\gamma}h_\gamma X^{\gamma},}

%\sni
(ici, $\alpha$, $\beta$, $\gamma$ sont des multi-indices, et si $\alpha=(\alpha_1,\ldots,\alpha_k)$,
$X^{\alpha}$ est une notation pour $X_1^{\alpha_1}\cdots X_k^{\alpha_k}$).
Soit $\gA=\gZ[(h_\gamma)]$ le sous-anneau
engendré par les \coes de $h$ ($\gZ$ est le sous-anneau de $\gB$ engendré par $1_\gB$). Alors, chaque $f_\alpha g_\beta$ est entier sur $\gA$.
\end{corollary}
\begin{proof}
On applique ce qu'il est convenu d'appeler l'\emph{astuce de \KRAz}:
on pose~$X_j=T^{n^j}$ avec $n$ assez grand. Ceci transforme $f$, $g$ et $h$
en des  \pols $F(T)$, $G(T)$, $H(T)$  dont les \coes sont
respectivement ceux de~$f$,~$g$ et~$h$.
\end{proof}
%

%\newpage

%%%%%%%%%%%%%%%%%%%%%%%%%%%%%%%%%%%%%%%%%%%%%%%%%%%%%%%%%%%%%%%%%%%%%%%%%%%
\section[L'\adu (1)]{L'\adu \\ pour un \polu sur un anneau
commutatif (1)} \label{sec0adu}

\emph{Avertissement.}
Dans un contexte où l'on manipule des \algs il est parfois préférable
de garder l'intuition qu'à la base, on a envie d'avoir un corps, même
si c'est seulement un anneau commutatif. Dans un tel cas nous choisissons
de donner un nom comme $\gk$
à l'anneau de base.
C'est ce que nous ferons dans cette section dédiée à l'\aduz.
\\
Lorsque nous avons vraiment affaire à un corps discret,
nous utiliserons plutôt une écriture telle que $\gK$.

Nous procédons maintenant à l'opération inverse de
celle qui passait de l'anneau des \pols au sous-anneau des \pols \smqsz.
\\
En présence d'un \polu $f=T^n+\sum_{k=1}^n(-1)^{k} s_kT^{n-k} \in \kT$ sur un anneau $\gk$,
nous voulons disposer d'une extension de $\gk$ où
le \pol se décompose en facteurs \linsz.
Une telle extension peut être construite de manière
purement formelle. Elle s'appelle l'\aduz.

%:     Definota{definotaAdu}
\begin{definota}\label{definotaAdu}
Soit $f=T^n+\sum_{k=1}^n(-1)^{k} s_kT^{n-k} \in \kT$ un \polu de
degré~$n$.
On note
$\Adu_{\gk,f}$ l'\emph{\adu de $f$ sur $\gk$}
définie comme suit:%
\index{algebre de@\aduz!de $f$ sur $\gk$}%
\index{ideal@idéal!des relateurs \smqsz}%
\index{relateur!symétrique}
$$\preskip.2em \postskip.4em
\Adu_{\gk,f}=\kXn/\cJ(f)=\kxn
,$$
où $\cJ(f)$ est l'\emph{\id des relateurs \smqsz} \ncr pour identifier
 dans le quotient $\prod_{i=1}^n(T-x_i)$ avec
$f(T)$.
Précisément, si  $S_1$, $ S_2$, \ldots,  $S_n$
sont les fonctions \smqs \elrs des~$X_i$,
 l'\id   $\cJ(f)$ est donné par:
$$\preskip.2em \postskip.4em
\cJ(f)=\gen{S_1-s_1,\;S_2-s_2,\ldots ,\;S_n-s_n}.$$
\end{definota}

L'\adu $\gA=\Adu_{\gk,f}$ peut être \caree par la \prt suivante.

%:    Fact{factEvident}---------
\begin{fact}
\label{factEvident} \emph{(Algèbre de \dcn \uvlez, \prt \caraz)}
\begin{enumerate}
\item Soit $\gC$ une \klg pour laquelle $f(T)$ se décompose en produit de
facteurs $T-z_i$. Alors, il existe un unique \homo de \klgs \hbox{de $\gA$} \hbox{vers $\gC$} qui
envoie les $x_i$ sur les $z_i$.
\item Ceci caractérise l'\adu $\gA=\Adu_{\gk,f}$, à \iso unique près.
\item Si en outre $\gC$ est engendrée (comme \klgz) par les $z_i$,
elle est isomorphe à un quotient de $\gA$%
%, la classe de $x_i$ correspondant à $z_i$
.
\end{enumerate}
\end{fact}
%--- end-fact-----------------------------------------
%
\begin{proof}
Pour le point \emph{1} on utilise la proposition
\ref{propApolLibre} qui décrit les \algs de \pols
comme des \algs librement engendrées par les \idtrs et le fait
\ref{factUnivQuot} qui décrit les anneaux quotients comme ceux qui permettent
de factoriser de manière unique certains \homosz. Le point \emph{2} résulte
de la constatation qu'un objet qui résout un \pb universel est toujours
unique à \iso unique près.
\end{proof}

Et en prenant $\gC=\gA$  on obtient que toute permutation de  $\so{1,\ldots,n}$ produit
un (unique) $\gk$-\auto de $\gA$.
\\
Dit autrement: le groupe $\Sn$ des permutations de $\so{\Xn}$
agit \hbox{sur $\kXn$} et fixe l'\id $\cJ(f)$, donc l'action passe
au quotient et ceci définit  $\Sn$ comme groupe d'\autos de l'\aduz.

  Pour étudier l'\adu
on introduit les \ix{modules de Cauchy} qui sont les \pols suivants:
{\small
%--------------------begin array---------------
$$
{\arraycolsep2pt\begin{array}{rcl}
f_1(X_1)& =  & f(X_1)   \\
f_2(X_1,X_2)& =  &  \big(f_1(X_1)-f_1(X_2)\big)\big/(X_1-X_2)    \\
 & \vdots  &  \\
f_{k+1}(X_1,\ldots ,X_{k+1})& =  &  {\dsp \frac{f_k(X_1,\ldots ,X_{k-1},X_k)-
f_k(X_1,\ldots ,X_{k-1},X_{k+1})}{X_k-X_{k+1}}}    \\
 & \vdots  &       \\
f_n(X_1,\ldots ,X_n)& =  &{\dsp \frac{f_{n-1}(X_1,\ldots ,X_{n-2},X_{n-1})-f_{n-
1}(X_1,\ldots ,X_{n-2},X_{n})}{X_{n-1}-X_{n}}} .    \\
\end{array}}$$
%---------------------end array--------------
}

\medskip Le fait suivant résulte de la propriété \cara des \adusz.

%:     Fact{factAduAdu}
\begin{fact}\label{factAduAdu}
Avec les notations précédentes pour les modules de Cauchy, \hbox{soit $\gk_1=\gk[x_1]$} et $g_2(T)=f_2(x_1,T)$.
Alors, l'application $\gk_1$-\lin canonique~\hbox{$\Adu_{\gk,f}\to\Adu_{\gk_1,g_2}$} (qui envoie chaque
$x_i$ ($i\geq2$)
de $\Adu_{\gk,f}$ sur le~$x_i$  \hbox{de $\Adu_{\gk_1,g_2}$}) est un \isoz.
\end{fact}

\mni \exls (Modules de Cauchy)\\
 Avec $n=4$:
{\small
%--------------------begin array---------------
$$\arraycolsep2pt\begin{array}{rcl}
f_1(x)& =  & x^4-s_1x^3+s_2x^2-s_3x+s_4   \\
f_2(x,y)& =  &  (y^{3}+ y^{2}x+ yx^{2} + x^{3})- s_1(y^{2}+yx+x^2)  +
s_2(y+x)  - s_3 \\
&=&y^{3} +y^{2}(x- s_1) + y(x^{2}- s_1x+ s_2) + (x^{3}- s_1x^{2}+ s_2x-
s_3)  \\
f_{3}(x,y,z)& =  &     (z^2+ y^2+ x^2+ zy + zx+ yx) - s_1(z+y+x)   + s_2
\\
&=& z^{2} +z (y+x- s_1)   + \big((y^{2} + yx + x^{2})- s_1(y+ x)  + s_2\big)
\\
f_4(x,y,z,t)& =  &t+z+y+x-s_1 .
\end{array}$$
%---------------------end array--------------
}
Pour $f(T)=T^6$:

\snac{
\arraycolsep2pt\begin{array}{rcl}
f_2(x,y)& =  &  y^{5}+ y^{4}x+ y^3x^{2} + y^2x^{3}+yx^4+x^5 \\
f_{3}(x,y,z)& =  &   (z^{4}+ y^{4}+ x^{4})+ ( z^{2}y^{2} + z^{2} x^{2}+
y^{2}x^{2}) +  \\
&&( zy^{3} + zx^{3}+ y z^{3}+  y x^{3}+ x z^{3} +x y^{3})  +\\
&& ( zy x^{2}  + zx  y^{2} + yx  z^{2}    )             \\
f_4(x,y,z,t)& =  & (t^{3}+ z^{3}+ y^{3} +x^{3}) + (tzy+tyx+tzx+zyx)+ \\
&&t^{2}(z+y+x)+ z^{2}(t+y+x)+ \\
&&y^{2}(t+z+x)+ x^{2} (t+z+y)       \\
f_5(x,y,z,t,u)&=&  (u^{2} + t^{2}+ z^{2}+ y^{2} +x^{2}) +\\
&& (x u + x t + x z +
x y + t u + z u + z t + y u + y t + y z)\\
f_6(x,y,z,t,u,v)&=& v+u+t+z+y+x .
\end{array}}

\smallskip 
Plus \gnltz, pour $f(T)=T^n$, $f_k(t_1,\ldots ,t_k)$ est la somme de
tous les \moms de degré $n+1-k$ en $t_1,\ldots ,t_k$.\\
Ceci permet par linéarité d'obtenir une description précise explicite
des modules de Cauchy pour un \pol arbitraire.
\eoe

\medskip
D'après la remarque qui suit le dernier exemple,
le \polz~$f_i$ est \smq en les variables $X_1$, $\ldots$, $ X_i$, \mon
en $X_i$, de degré total $n-i+1$.

Le fait \ref{factEvident} implique que l'idéal $\cJ(f)$  est égal à
l'idéal engendré par les modules de Cauchy. En effet, le quotient
par ce dernier \id réalise clairement la même \prt \uvle que le quotient
par $\cJ(f)$.

Donc l'\adu est un \kmo
libre de rang~$n!$. Plus \prmtz, on obtient le résultat suivant.

%--- Fact{factBase}-------------
\begin{fact}
\label{factBase}
Le \kmo $\gA=\Adu_{\gk,f}$ est libre et une base est formée par les \gui{\momsz}
$x_1^{d_1}\cdots x_{n-1}^{d_{n-1}}$ tels que pour $k=1,\ldots ,n-1$ on ait
$d_k\leq n-k$.
\end{fact}
%--- end-fact-----------------------------------------

%:     Corollary{corfactBase}
\begin{corollary}\label{corfactBase}
En considérant l'\adu du \polu \gnq $f(T)=T^n+\sum_{k=1}^n(-1)^kS_kT^{n-k}$, où les $S_i$ sont des \idtrsz, on obtient
une \alg de \pols $\kxn$ avec les $S_i$ qui s'identifient aux \pols \smqs \elrs 
en les~$x_i$.
\end{corollary}

\stMF
\comm (Pour \ceux qui connaissent les bases de Gr\"obner)\\
Dans le cas où $\gk$ est un \cdiz, les modules de Cauchy peuvent
être vus comme une base de Gr\"obner de l'\id $\cJ(f)$, pour l'ordre monomial
lexicographique avec $X_1<X_2<\cdots< X_n$.
\\
En fait, même si $\gk$ n'est pas un \cdiz, les modules de Cauchy
fonctionnent comme une base de Gr\"obner: tout \pol en les $x_i$ se
réécrit sur la base de \moms précédente par divisions successives
par les modules de Cauchy. On divise tout d'abord
par $f_n$ par rapport à la variable~$X_n$, ce qui la fait
disparaître. Ensuite on divise par $f_{n-1}$ par rapport à la
variable  $X_{n-1}$, ce qui la ramène en degré $\leq 1$, et ainsi de
suite.
\eoe

%%%%%%%%%%%%%%%%%%%%%%%%%%%%%%%%%%%%%%%%%%%%%%%%%%%%%%%%%%%%%%%%%%%%%%%%%%%%%%%%%%%%%%%
\section{Discriminant, \dinz}
\label{secDisc}
\vspace{4pt}
%  subsec  Définition du \discri d'un \polu
\subsec{Définition du \discri d'un \polu}

On définit le \ix{discriminant} d'un \polu $f$ en une variable sur un anneau commutatif $\gA$ en commençant par le cas où $f$ est
le \polu \gnq de degré $n$:\index{discriminant!d'un polynôme unitaire}
$$
\mathrigid1mu f(T)\,=\,T^n-S_1T^{n-1}+S_2T^{n-2}+\cdots+(-1)^nS_n\;\in\,\ZZ[S_1,\ldots,S_n][T]\,=\,\ZZuS[T].
$$
On peut écrire $f(T)=\prod_i(T-X_i)$ dans $\ZZXn$  
(corolaire \ref{corthSymE1}), et l'on pose
\begin{equation}\preskip.0em \postskip.4em
\label{eqDiscri}
\disc_T(f)=  (-1)^{n(n-1)/2}\, \prod\nolimits_{i=1}^n f'(X_i)
=\prod\nolimits_{1\leq i< j\leq n}(X_i-X_j)^2.
\end{equation}
Comme manifestement ce \pol en les $X_i$ est invariant par permutation
des variables il existe un unique \pol en les $S_i$, $D_n(S_1,\ldots,S_n)\in\ZZuS$
qui est égal à $\disc_T(f)$. En bref, les variables auxiliaires $X_i$
peuvent bien disparaître.

Ensuite pour un \pol \gui{concret}
$$
g(T)=T^n-s_1T^{n-1}+s_2T^{n-2}+\cdots+(-1)^ns_n\in\AT,
$$
 on définit $\disc_T(g)=D_n(s_1,\ldots,s_n)$.

Naturellement, s'il arrive que $g(T)=\prod_{i=1}^n(T-b_i)$ dans un anneau
$\gB\supseteq\gA$, on obtiendra $\disc_T(g)=\prod\nolimits_{1\leq i< j\leq n}(b_i-b_j)^2$
en évaluant la formule (\ref{eqDiscri}).
En particulier, en utilisant l'\adu  on pourrait définir directement le \discri
par cette formule.

Un \polu est dit \emph{séparable} lorsque son \discri est
\ivz.%
\index{separable@séparable!\polu ---}%
\index{polynome@\pol!unitaire séparable}

%:  subsec  Diagonalisation de matrices
\subsec{Diagonalisation de matrices sur un anneau}

Rappelons d'abord que si $f\in\AT$, un \emph{zéro de $f$ dans une \Alg
$\gB$} (donnée par un \homo $\varphi:\gA\to\gB$) est un $y\in\gB$ qui annule
le \polz~$f^{\varphi}$, image de $f$ dans $\gB[T]$.

Le  zéro $y$ est dit \emph{simple} si en outre $f'(y)\in\gB\eti$ (on dit aussi que c'est une \emph{racine simple} de $f$).%
\index{simple!zéro ---}%
\index{simple!racine ---}%
\index{racine!simple}%
\index{zero@zéro!simple d'un \polz}%
\index{zero@zéro!d'un \polz}

Nous nous intéressons ici aux \dins de matrices sur un anneau commutatif
arbitraire, lorsque le \polcar est \emph{\splz}.

Nous avons tout d'abord le classique \gui{lemme des noyaux}~\ref{lemDesNoyaux}.
\\
Voici ensuite une \gnn du \tho qui affirme (dans le cas d'un \cdiz) qu'un zéro simple du \polcar
définit un sous-espace propre de dimension $1$.

%:     Lemma{lemValProp}
\begin{lemma}\label{lemValProp}
Soit $n\geq2$, $a\in\gA$ et $A\in\Mn(\gA)$ une matrice dont le \polcar
$f(X)=\rC A(X)$ admet $a$ comme zéro simple. \hbox{Soit  $g=f/(X-a)$},
$h=X-a$,
 $K=\Ker h(A)$ et $I=\Im h(A)$.
\begin{enumerate}
\item On a $K=\Im g(A)$,  $I=\Ker g(A)$ et
 $\Ae n=I\oplus K$.
\item La matrice $g(A)$ est de rang $1$, et $h(A)$ de rang $n-1$.
\item Si un \polz~$R(X)$ annule $A$,
alors $R(a)=0$, \cad $R$ est multiple de $X-a$.
\item Les mineurs principaux d'ordre $n-1$ de $A-a\In$ sont \comz.
Quand on localise
en inversant un tel mineur, la matrice
$g(A)$ devient simple de rang $1$, les modules $I$ et $K$ deviennent libres de
rangs $n-1$ et $1$.
\end{enumerate}
\end{lemma}
\begin{proof}
On suppose \spdg que $a=0$.
\\
Alors, $f(X)=Xg(X)$, $h(A)=A$,
$g(A)=\pm\wi A$, $\Tr\big(g(A)\big)=g(0)$ (lemme~\ref{lemPrincipeIdentitesAlgebriques}
 point \emph{6}), et $g(0)=f'(0)\in\Ati$.\\
 \emph{1.}
 On écrit $g(X)=Xk(X)+g(0)$: cela montre que les \pols $g(X)$ et $X$ sont \comz.
Vu le \tho de Cayley-Hamilton, le lemme des noyaux s'applique
et donne le point~\emph{1.}

\emph{2.} Notons $\mu_1$, $\ldots$, $\mu_n$ les mineurs principaux d'ordre $n-1$ de $A$. \\
Puisque $g(A)=\pm\wi A$, on obtient
$g(0)=\Tr \big(g(A)\big)=\pm \Tr\wi A =\pm\sum_i\mu_i$. Ceci montre que
$\rg\big(h(A)\big)=n-1$ et $\rg\big(g(A)\big)\geq1$.
Enfin, on sait que $\rg(\wi A)\leq1$ d'après
le lemme  \ref{lemPrincipeIdentitesAlgebriques} point~\emph{8}.

\emph{3.} Supposons $R(A)=0$.
En multipliant par $\wi{A}$, on obtient $R(0)\wi{A} = 0$
   (puisque $\wi{A} A = 0$). En prenant la trace,
   $R(0)\Tr(\wi{A}) = 0$ donc $R(0) = 0$.
%Pour $u\in K$, on a $Au=0$ et donc
% $0u=R(A)(u)=R(0)u$. Donc $R(0)g(A)=0$,  $R(0)\Tr(g(A))=0$ et $R(0)=0$.
\\
Remarquons que le point \emph{3} résulte aussi du point \emph{4.}

\emph{4.} On a déjà vu que les $\mu_i$ sont \comz.
Après \lon en un~$\mu_i$, la matrice $g(A)$
devient simple de rang $1$ en vertu du lemme de
la liberté~\ref{lem pf libre}. Donc $I$ et $K$ deviennent libres de
rangs $n-1$ et $1$.
\end{proof}
%

%:     Proposition{propDiag}
\begin{proposition}\label{propDiag} \emph{(Diagonalisation d'une matrice dont le \polcar
est \splz)}  %\\
Soit $A\in\Mn(\gA)$ une matrice dont le \polcar $\rC A(X)$ est \splz,
et un anneau $\gA_1 \supseteq\gA$ sur lequel on peut écrire
$\rC A(X)=\prod_{i=1}^n(X-x_i)$
(par exemple, $\gA_1=\Adu_{\gA,f}$). \\
Soit $K_i=\Ker(A-x_i\In)\subseteq \gA_1^n$.
\begin{enumerate}
\item $\gA_1^n=\bigoplus_i K_i$.
\item Chaque $K_i$ est l'image d'une matrice de rang $1$.
\item Tout \polz~$R$ qui annule $A$ est multiple de $\rC A$.
\item Après \lon en des \eco de $\gA_1$ la matrice est \digz, semblable à $\Diag(\xn)$.
\end{enumerate}
\end{proposition}
NB:  si $\alpha\in\End_{\gA_1}(\gA_1^n)$ a pour matrice  $A$,
on a $\alpha\frt{K_i}=x_i\,\Id_{K_i}$ pour chaque~$i$.
\begin{proof}
Conséquence \imde du lemme des noyaux et du lemme~\ref{lemValProp}.
Pour rendre la matrice \dig il suffira d'inverser un produit $\nu_1\cdots\nu_n$
où chaque $\nu_i$ est un mineur principal d'ordre $n-1$ de la matrice $A-x_i\In$ (ce qui fait a priori $n^n$ \lons \comez).
\end{proof}

\rem Un résultat analogue concernant une matrice qui annule un \polz~$\prod_{i}(X-x_i)$ \spl est donné en exercice~\ref{exo2Diag}. La preuve est \elrz.
\eoe

%:  subsec  La matrice \gnq est \dig
\penalty-2500
\subsec{La matrice \gnq est \dig}

Considérons  $n^2$ \idtrs
$(a_{i,j})_{i,j\in \lrbn}$ et notons $A$ la matrice
correspondante (elle est à \coes dans
$\gA=\ZZ[(a_{i,j})]$).

%Il y a deux manières d'interpréter le titre du paragraphe présent.
%La première, très contraignante, est la suivante.

%:     Proposition{propMatGenDiag}
\begin{proposition}\label{propMatGenDiag}
La matrice \gnq $A$ est \dig sur un anneau~$\gB$ contenant $\ZZ[(a_{i,j})]=\gA$.
\end{proposition}
\begin{proof}
Soit
$f(T)=T^n-s_1T^{n-1}+\cdots+(-1)^ns_n$ le \polcar de $A$.
Alors les \coes $s_i$ sont \agqt indépendants sur $\ZZ$.\linebreak 
Il suffit pour s'en rendre compte de spécialiser $A$ en la matrice compagne d'un \polu \gnqz.
\\
En particulier, le \discri $\Delta=\disc(f)$ est non nul dans l'anneau
intègre~$\gA$.
On  considère alors
l'anneau $\gA_1=\gA[1/\Delta]\supseteq\gA$ puis l'\adu $\gC=\Adu_{\gA_1,f}$. Notons $x_i$
les \elts de $\gC$ tels que $f(T)=\prod_i(T-x_i)$.
\\
On applique enfin la proposition  \ref{propDiag}. Si l'on veut
aboutir à une matrice \digz, on inverse par exemple $a=\prod_i\det\big((A-x_i\In)_{1..n-1,1..n-1}\big)$. Il s'agit d'un \elt de $\gA$ et il suffit de se convaincre qu'il n'est pas nul en exhibant une matrice particulière,
par exemple la matrice compagne du \pol $X^n-1$.
\\
En définitive on considère $\gA_2=\gA[1/(a\Delta)]\supseteq\gA$
et l'on prend

\snic{\gB=\Adu_{\gA_2,f}\supseteq\gA_2.}
\end{proof}

La force du résultat précédent, \gui{qui simplifie considérablement la vie}
est illustrée dans les deux paragraphes qui suivent.

%:  subsec Identité concernant les \polcars
\subsec{Identité concernant les \polcarsz}

%:     Proposition{prop1tschir}
\begin{proposition}\label{prop1tschir}
Soient $A$ et $B\in \Mn(\gA)$ deux matrices qui ont le même \polcarz,
et soit $g\in\AT$. Alors les matrices $g(A)$ et $g(B)$ ont même \polcarz.
\end{proposition}

%:     Corollary{corprop1tschir}
\begin{corollary}\label{corprop1tschir}\label{lemPolCar}~
\begin{enumerate}
\item Si $A$ est une matrice de \polcar $f$, et si l'on peut écrire
$f(T)=\prod_{i=1}^n(T-x_i)$ sur un anneau $\gA_1\supseteq \gA$, 
alors le \polcar de $g(A)$
est égal au produit $\prod_{i=1}^n\big(T-g(x_i)\big)$.
\item Soit $\gB$ une \Alg libre de rang fini $n$ et
$x\in\gB$.
On suppose que dans   $\gB_1\supseteq \gB$,
on a $\rC{\gB/\!\gA}(x)(T) = (T - x_1) \cdots (T - x_n)$. Alors, pour
tout $g \in \gA[T]$, on a les \egts suivantes:
$$\preskip.3em \postskip.3em 
\rC{\gB/\!\gA}\big(g(x)\big)(T) = \big(T - g(x_1)\big) \cdots \big(T - g(x_n)\big),  
$$
$\Tr\iBA \big(g(x)\big) = \som_{i=1}^n g(x_i)$   et 
$\rN\iBA \big(g(x)\big) = \prod_{i=1}^n g(x_i)$.
\end{enumerate}
\end{corollary}
\begin{Proof}{\Demo de la proposition et du corolaire.} \\
Point \emph{1} du corolaire. On considère la matrice
$\Diag(\xn)$ qui a même \polcar que $A$ et on applique la proposition
avec l'anneau~$\gA_1$.\\
Inversement, si le point \emph{1} du corolaire est démontré pour
$\gA_1=\Adu_{\gA,f}$,
il implique la proposition \ref{prop1tschir} car le \pol
$\prod_{i=1}^n\big(T-g(x_i)\big)$ calculé \hbox{dans $\Adu_{\gA,f}$}
ne dépend que de $f$ et $g$.\\
On note maintenant que la structure de l'énoncé du corolaire, point \emph{1},
lorsque l'on prend $\gA_1=\Adu_{\gA,f}$, est une famille d'\idas
avec pour \idtrs les \coes de la matrice~$A$.
Il suffit donc de le démontrer pour la matrice \gnqz.
Or elle est \dig sur un suranneau (proposition~\ref{propMatGenDiag}),
et pour une matrice \dig le résultat est clair.
\\
 Enfin, le point \emph{2} du corolaire est une conséquence \imde du point~\emph{1.}
\end{Proof}
%

%%%%%%%%%%%%%%%%%%%%%%%%%%%%%%%%%%%%%%%%%%%%%%%%%%%%%%%%%%%%%%%%%%%%%%%%%%%
%:  subsec  Identité concernant les puissances extérieures
\subsec{Identité concernant les puissances extérieures}

Les résultats suivants,
analogues à la proposition~\ref{prop1tschir}
et au corolaire~\ref{corprop1tschir} peuvent être démontrés
en suivant exactement les mêmes lignes.

%:     Proposition{propPolCarPuissExt}
\begin{proposition}\label{propPolCarPuissExt}
Si $\varphi$ est un \endo d'un \Amo libre de rang fini, le \polcar
de ${\Al{k}\!\varphi}$
ne dépend  que de l'entier $k$ et du \polcar de $\varphi$.
\end{proposition}

%:     Corollary{corpropPolCarPuissExt}
\begin{corollary}\label{corpropPolCarPuissExt}
Si $A\in\Mn(\gA)$ est une matrice de \polcar $f$, et si
$f(T)=\prod_{i=1}^n(T-x_i)$ dans un suranneau de $\gA$,
 alors le \polcar de $\Al{k}\!A$
est égal au produit 

\snic{\prod_{J\in \cP_{k,n}}(T-x_J)$, où $x_J=\prod_{i\in J}x_i.}
\end{corollary}

%:  subsec    Tschirnhaus
\subsec{Transformation de Tschirnhaus}

%:     Definition{defiTschir}
\begin{definition}\label{defiTschir}
Soient $f$ et $g\in\AT$ avec $f$ \mon de degré $p$.
On considère l'\Alg $\gB=\aqo{\AT}{f}$,
qui est un \Amo libre de rang~$p$.
On~définit le \emph{transformé de Tschirnhaus de
$f$ par $g$},  noté
$\Tsc_{\gA,g}(f)$ ou~$\Tsc_g(f)$, par l'\egt

\snic{\Tsc_{\gA,g}(f)=\rC{\gB/\!\gA}(\ov g),\quad (\ov g $ est la classe de $g$ dans $\gB).}%
\index{polynome@\pol!transformé de Tschirnhaus}
\end{definition}

%Rappelons que $\rC{\gB/\!\gA}(g)$ est le  \polcar de $\mu_g$,
%où $\mu_g$ la multiplication par
%(la classe de) $g$ dans $\gB$.

La proposition \ref{prop1tschir} et le
corolaire \ref{corprop1tschir}, donnent le résultat suivant.

%:     Proposition{prop2tschir}
\begin{proposition}\label{prop2tschir}
Soient $f$ et $g\in\AT$ avec $f$ \mon de degré $p$.
\begin{enumerate}
\item Si $A$ est une matrice telle que $f(T)=\rC A(T)$, on a
$$\preskip.2em \postskip.4em 
\Tsc_g(f)(T)=\rC{g(A)}(T). 
$$
\item Si $f(T)=\prod_i(T-x_i)$ sur un anneau  qui contient $\gA$, on a
$$\preskip.3em \postskip.4em\ndsp 
\Tsc_g(f)(T)=\prod_i\big(T-g(x_i)\big), 
$$
en particulier, on obtient avec $\gB=\aqo{\AT}{f}$
$$\preskip.3em \postskip.4em\ndsp 
\rN\iBA (g)=\prod_ig(x_i) \hbox{ et }
\Tr\iBA (g)=\sum_ig(x_i). 
$$
\end{enumerate}
\end{proposition}

\smallskip 
\rem On peut aussi écrire $\Tsc_{\gA,g}(f)(T)=\rN_{\gB[T]/\!\AT}(T-\ov g)$.
En fait pour une notation entièrement non ambigu\"{e} on devrait noter
$\Tsc(\gA,f,g,T)$ au lieu de $\Tsc_{\gA,g}(f)$.
Une ambiguïté analogue se trouve d'ailleurs dans la nota\-tion $\rC{\gB/\!\gA}(g)$.
\eoe

%:  subsubsection*{Calcul du transformé de Tschirnhaus}
\rdb
\subsubsection*{Calcul du transformé de Tschirnhaus}
Rappelons que la matrice $C$ de l'\endo $\mu_t$ de multiplication par $t$
(la classe de $T$ dans $\gB$)
est appelée la matrice compagne de $f$ (voir \paref{matrice.compagne}).
Alors  la  matrice (sur la même base) de $\mu_{\ov g}=g(\mu_t)$ est la matrice
$g(C)$. Donc~$\Tsc_g(f)$ est le \polcarz\footnote{Le calcul rapide des
\deters et \polcars suscite un grand intérêt en calcul formel.
On pourra par exemple consulter \cite{Jou}.
Une autre formule que l'on peut utiliser pour le calcul du transformé
de Tschirnhaus est $\Tsc_g(f)=\Res_X(f(X),T-g(X))$ (voir le lemme~\ref{lemResultant})} de $g(C)$.

%:   subsec  Nouvelle version du discriminant
\subsec{Nouvelle version du discriminant}

Rappelons (\dfn \ref{defiDiscTra})
que lorsque $\gC\supseteq \gA$ est une \Alg libre de rang fini et $x_1$, \dots, $x_k\in\gC$,
on appelle discriminant de $(\xk)$
le \deter de la matrice
$\big(\Tr_{\gC/\!\gA}(x_ix_j)\big)_{i,j\in \lrbk}$. On le note $\disc_{\gC/\!\gA}(\xk)$.\\ 
En outre, si $(\xk)$ est une $\gA$-base de $\gC$,
on note $\Disc_{\gC/\!\gA}$ la classe multiplicative
de $\disc_{\gC/\!\gA}(\xk)$ modulo les carrés
de $\Ati$. On l'appelle  le {\discri de l'extension} $\gC\sur\gA$.

Nous faisons dans ce paragraphe le lien entre le \discri
des \algs libres de rang fini et le \discri des \polusz.
\\
Insistons sur le \crc remarquable de l'implication \emph{1a} $\Rightarrow$
\emph{1b} dans la proposition suivante.

%%%%%%%%%%%%%%%%%%%%%%%%%%%%%%%%%%%%%%%%%%%%%%%%%%%%%%%%%%%%%%%%

%:     Proposition{propdiscTra}
\begin{proposition}\label{propdiscTra}
 \emph{(Discriminant d'un \polu et forme trace)}%
 \index{discriminant!et forme trace}\\
Soit $\gB$ une \Alg libre de rang fini $n$,
$x \in \gB$ et $f=\rC{\gB/\!\gA}(x)(T)$. On~a:
$$\preskip.1em \postskip.3em
\disc(1, x, \ldots, x^{n-1}) = \disc(f) =
(-1)^{n(n-1) \over 2} \rN\iBA \big(f'(x)\big).
$$
On dit que $f'(x)$ est la \emph{différente de $x$}.%
\index{differente@différente!d'un \elt dans une \alg libre finie}
Les résultats suivants en découlent.
\begin{enumerate}
\item \Propeq
\begin{enumerate} \itemsep0pt
\item $\disc(f)\in\Ati$.
\item  $\Disc\iBA \in\Ati$ et $(1,x,\ldots,x^{n-1})$ est une $\gA$-base de $\gB$.
\item $\Disc\iBA \in\Ati$ et $\gB=\gA[x]$.
\end{enumerate}
\item Si $\Disc\iBA $ est \ndzz, \propeq
\begin{enumerate} \itemsep0pt
\item $\Disc\iBA $ et $\disc(f)$ sont associés.
\item $(1,x,\ldots,x^{n-1})$ est une $\gA$-base de $\gB$.
\item $\gB=\gA[x]$.
\end{enumerate}
\item Le \discri d'un \polu $g\in\AT$ représente
(modulo les carrés de
$\Ati$) le \discri de l'extension $\aqo{\AT}{g}$ de $\gA$.
\\
On a $\disc_T(g)\in\Ati$ \ssi $\gen{g(T),g'(T)}=\gA.$
\end{enumerate}
\end {proposition}
\begin{proof}
Dans un sur-anneau $\gB'$ de $\gB$, on peut écrire 
$f(T) = (T - x_1)\cdots (T - x_n)$.
\\
 Pour un
$g \in \gA[T]$, en appliquant le corolaire~\ref{lemPolCar}, on obtient les \egts

\snic{\Tr\iBA \big(g(x)\big) = g(x_1) + \cdots +
g(x_n)$  et $\rN\iBA \big(g(x)\big) = g(x_1) \cdots g(x_n).}

%\sni
 On note $M \in \Mn(\gA)$ la matrice intervenant dans le calcul du
\discri de $(1, x, \ldots, x^{n-1})$:

\snic{M=\big((a_{ij})_{i,j\in\lrb{0..n-1}}\big),\qquad a_{ij} = \Tr\iBA (x^{i+j}) = x_1^{i+j} + \cdots + x_n^{i+j}\,.}

%\sni
Soit $V \in \Mn(\gB')$  la matrice de Vandermonde ayant pour lignes $[\,x_1^i\; \ldots\;
x_n^i\,]$ (où $i\in\lrb{0.. n-1}$). Alors $M = V\tra {V}$.  On en déduit:

\snic{\det(M) = \det(V)^2 = \prod\nolimits_{i < j} (x_i - x_j)^2 = \disc(f)\,.}

%\sni
Ceci démontre la première \egtz. Puisque
$\rN\iBA \big(f'(x)\big) = f'(x_1) \cdots f'(x_n)$ et
$
f'(x_i) =\prod_{j \mid j \ne i} (x_i - x_j)
$,
il vient:
$$\preskip.1em \postskip.4em\ndsp
\rN\iBA \big(f'(x)\big) = \prod\nolimits_{(i,j) \mid j \ne i} (x_i - x_j) =
(-1)^{n(n-1) \over 2} \prod\nolimits_{i < j} (x_i - x_j)^2.
$$
La \dem des conséquences est laissée \alec (utiliser la proposition~\ref{defiDiscTra}).
\end {proof}

%:   subsec  Discriminant d'une \adu
\subsect{Discriminant d'une \alg de \dcn \\ \uvle}{Discriminant d'une \aduz}

L'\egt du \discri \gui{tracique} et du \discri \gui{\pollz}, jointe à la formule de transitivité (\thrf{thTransDisc}) nous permet le calcul suivant.

%:     Fact{factDiscriAdu}
\begin{fact}\label{factDiscriAdu} \emph{(Discriminant d'une \aduz)}\\
Soit $f$ un \polu de degré $n\geq2$ de $\kT$ et $\gA=\Adu_{\gk,f}$.
\\ Alors
$\Disc\iAk =\big(\disc_T(f)\big)^{n!/2}$.
\end{fact}
\begin{proof}
On reprend les notations de la section \ref{sec0adu}.
On raisonne par \recu sur~$n$, le cas $n=2$ étant clair.
On a $\gA=\gk_1[x_2,\ldots,x_n]$
avec 

\snic{\gk_1=\gk[x_1]\simeq\aqo{\gk[X_1]}{f(X_1)}.}

%\sni
En outre, $\gA\simeq\Adu_{\gk_1,g_2}$
où 

\snic{g_2(T)=f_2(x_1,T)
=\big(f(T)-f(x_1)\big)\big/(T-x_1)\in\gk_1[T]\subseteq\gA[T].}

%\sni
La formule de transitivité des \discris
donne alors les \egts suivantes.
$$
{\mathrigid1mu \Disc\iAk =
\Disc_{\gk_1/\gk}^{~\dex{\gA: \gk_1}}\
 \rN_{\gk_1/\gk}(\Disc_{\gA/\gk_1})=
 (\disc f)^{(n-1)!}\ \rN_{\gk_1/\gk}(\Disc_{\gA/\gk_1}).}
$$
%\sni
En utilisant l'\hdr on obtient l'\egt

\snic{\Disc_{\gA/\gk_1}=(\disc g_2)^{(n-1)!/2}=
\big(\prod\nolimits_{2\leq i<j\leq n}(x_i-x_j)^2\big)^{(n-1)!/2}.}

%\sni
Pour $i \in \lrb {2..n}$, notons $\tau_i$ la transposition $(1,i)$; pour $z
\in \gk_1$, d'après le corolaire~\ref{lemPolCar}, 
$\rN_{\gk_1/\gk}(z) = z
\prod_{i=2}^n \tau_i(z)$.  
Appliqué à $z = \prod\nolimits_{2\leq i<j\leq n}(x_i-x_j)^2$, cela donne
$$\preskip.2em \postskip.2em
\rN_{\gk_1/\gk}(z) = (\disc f)^{n-2},
\; \hbox { d'où } \,
\rN_{\gk_1/\gk}(\Disc_{\gA/\gk_1}) = (\disc f)^{(n-2)\cdot(n-1)!/2}\,,
$$
puis
$$\preskip-.2em \postskip.4em
\Disc\iAk = (\disc f)^{(n-1)! + (n-2)\cdot(n-1)!/2} = (\disc f)^{n!/2}.
$$
%\sni

NB: un examen détaillé du calcul précédent montre que l'on a en fait calculé le \discri de la base \gui{canonique} de l'\adu décrite dans le fait~\ref{factBase}.
\end{proof}
%

%Le même type de calcul donne:

%:    Lemma{lemPolCarAdu}
\begin{lemma}
\label{lemPolCarAdu} (Mêmes hypothèses que pour le fait \ref{factDiscriAdu})\\
Pour tout  $z\in\gA$ on a:
$$\preskip.2em \postskip.4em 
\rC{\gA/\gk}(z)(T)=\prod\nolimits_{\sigma\in\Sn}\big(T-\sigma(z)\big). 
$$
En particulier, $\Tr\iAk (z)=\som_{\sigma\in\Sn}\!\sigma(z)$ et $\rN\iAk (z)=\prod_{\sigma\in\Sn}\!\sigma(z)$.
\end{lemma}
%--- end-lemma-----------------------------------------
%
\begin{proof}
Il suffit de montrer la formule pour la norme, car on obtient ensuite celle pour le \polcar en remplaçant $\gk$ par $\kT$ 
(ce qui remplace~$\gA$ par $\AT$).
La formule pour la norme se prouve par \recu sur le nombre de variables
en utilisant le fait \ref{factAduAdu}, la formule de transitivité pour les normes
et le corolaire~\ref{lemPolCar}.
\end{proof}
%

%%%%%%%%%%%%%%%%%%%%%%%%%%%%%%%%%%%%%%%%%%%%%%%%%%%%%%%%%%%%%%%%%%%%%%%%%%%
\newpage	
\section{Théorie de Galois de base (1)}
\label{secGaloisElr}

\Grandcadre{Dans la section  \ref{secGaloisElr},
$\gK$ désigne un \cdi non trivial.}

\vspace{-8pt}
%%%%%%%%%%%%%%%%%%%%%%%%%%%%%%%%%%%%%%%%%%%%%%%%%%%%%%%%%%%%%%%%%%%%%%%
%:subsec{Factorisation et zéros}
\subsec{Factorisation et zéros}

Rappelons qu'un anneau  est intègre  si tout \elt  est nul ou \ndzz\footnote{La notion est discutée plus en détail  \paref{subsecAnneauxqi}.}.
Un sous-anneau d'un anneau intègre est intègre.
Un \cdi est un anneau intègre. Un anneau $\gA$ est intègre \ssi son anneau
total de fractions $\Frac\gA$ est un \cdiz. On dit alors que $\Frac\gA$
est le \ixc{corps de fractions}{d'un anneau intègre} de~$\gA$.

%:     Proposition{propZerFactPol}
\begin{proposition}\label{propZerFactPol}
Soient $\gA\subseteq \gB$ des anneaux et $f\in\AT$ un \polu
de degré~$n$.
\begin{enumerate}
\item Si $z$ est un zéro de $f$ dans  $\gB$, $f(T)$
est divisible par $T-z$ dans~$\gB[T]$. 
\item On suppose désormais que $\gB$ est intègre et non trivial\footnote{On pourrait se passer de l'hypothèse négative \gui{non trivial} en lisant l'hypothèse que les $z_i$ sont \gui{distincts} comme signifiant que les $z_i-z_j$ sont \ndzsz.}. Si $z_1$, $\dots$, $z_k$ sont des zéros de $f$ deux à deux distincts
dans   $\gB$, le \pol $f(T)$ est divisible par
$\prod_{i=1}^k(T-z_i)$ dans~$\gB[T]$.
\item Si en outre $k=n$, alors $f(T)=\prod_{i=1}^n(T-z_i)$, et les $z_i$ sont les seuls zéros de $f$ dans $\gB$ et dans toute extension intègre de~$\gB$.
\end{enumerate}
\end{proposition}
%--------- fin proposition ---------------------------------------------- 
%
\begin{proof}
La \dem est \imdez, certains résultats plus précis sont dans l'exercice \ref{exoLagrange} consacré à l'interpolation de Lagrange.
\end{proof}
%

%%%%%%%%%%%%%%%%%%%%%%%%%%%%%%%%%%%%%%%%%%%%%%%%%%%%%%%%%%%%%%%%%%%%%%%
%:subsec{Algèbres \stfes sur un \cdiz}
\subsec{Algèbres \stfes sur un \cdiz}

%:     Definition{defiSTF}
\begin{definition}\label{defiSTF}~\\
Une \Klg $\gA$ est dite \emph{\stfez} si c'est
un \Kev libre de dimension finie.%
\index{strictement finie!algèbre --- sur un corps discret}%
\index{algèbre!strictement finie sur un corps discret}
\end{definition}

Autrement dit, on connaît une base finie de $\gA$ comme \Kevz.
Dans ce cas, pour un $x\in\gA$,  la trace, la norme, le \polcar de (la
multiplication par) $x$, ainsi que le \polmin de $x$ sur $\gK$,
noté~$\Mip_{\gK,x}(T)$ ou $\Mip_{x}(T)$, peuvent se
calculer par les méthodes standards de l'\alg \lin sur un \cdiz.
De même toute sous-\Klg finie de~$\gA$ est \stfe et
l'intersection de deux sous-\algs \stfes est \stfez.

%:     Lemma{lemEntReduitConnexe}
\begin{lemma}\label{lemEntReduitConnexe}
Soit $\gB\supseteq\gK$  un anneau entier sur $\gK$. \Propeq
\begin{enumerate}
\item $\gB$ est un \cdiz.
\item $\gB$ est \sdzz~: $xy=0\Rightarrow (x=0$ \emph{ou} $y=0)$.
\item $\gB$ est connexe et réduit.
\end{enumerate}
En conséquence si $\gB$ est un \cdiz, toute sous-\Klg finie de $\gB$
est un \cdiz.
 \end{lemma}
\begin{proof}
 Les implications \emph{1} $\Rightarrow$ \emph{2}
$\Rightarrow$ \emph{3} sont claires.
\\
\emph{3} $\Rightarrow$ \emph{1}.
Soit $x\in\gB$, il annule un \pol non nul de $\KX$ que l'on peut supposer
de la forme $X^k\big(1-XR(X)\big)$. Alors $x\big(1-xR(x)\big)$ est nilpotent donc nul.
L'\elt $e=xR(x)$ est \idm et $x=ex$. Si $e=0$, alors~$x=0$. Si $e=1$,
alors $x$ est \ivz.
\end{proof}
%

%:     Lemma{lemdefiSTF}
\begin{lemma}\label{lemdefiSTF}
Soient $\gK\subseteq\gL\subseteq\gA$ avec $\gA$ et $\gL$ \stfs sur $\gK$.
Si~$\gL$ est un \cdiz, alors $\gA$ est \stfe sur $\gL.$
\end{lemma}
\begin{proof}
\Demo laissée \alec (ou voir le fait \ref{fact1Etale} point \emph{\iref{i3fact1Etale}}).
\end{proof}

Si $g$ est un \pol \ird de $\gK[T]$, l'\alg quotient
$\aqo{\gK[T]}{g}$ est un \cdi \stf sur $\gK$.
En fait, comme corolaire des deux lemmes précédents on obtient que toute extension \stfe de \cdis s'obtient en itérant cette construction.

%:     Fact{factStrucExtFiC}
\begin{fact}\label{factStrucExtFiC} \emph{(Structure d'une extension \stfe de \cdisz)}\\
Soit $\gL=\gK[\xm]$ un \cdi \stf sur $\gK$.\\ 
Pour $ k\in\lrb{1.. m+1}$, notons $\gK_k=\gK[(x_i)_{i<k}]$ et 
$f_k=\Mip_{\gK_{k},{x_k}}(T)$, de sorte que $\gK_1=\gK$, et pour $k\in \lrbm$, $\gK_{k+1}\simeq\aqo{\gK_k[X_k]}{f_k(X_k)}$. \\
Alors, pour $k<\ell$ dans $\lrb{1.. m+1}$,  l'inclusion $\gK_k\to\gK_\ell$ est une extension \stfe de \cdisz, avec 
$$\preskip.4em \postskip.4em\ndsp 
\dex{\gK_\ell:\gK_k} = \prod\nolimits_{k\leq i< \ell}\dex{\gK_{i+1}:\gK_i}
  = \prod\nolimits_{k\leq i< \ell}\deg_T(f_i). 
$$
En outre, si $F_k\in\gK[X_1,\dots,X_{k}]$ est un \polu en $X_k$ pour lequel on a $F_k\big((x_i)_{i<k},X_k\big)=f_k(X_k)$, on obtient par \fcn de l'\homo d'\evnz, un \iso
$$\preskip.3em \postskip.4em 
       \aqo{\gK[\Xm]}{F_1,\dots,F_m} \; \simarrow \; \gL. 
$$
\end{fact}
%--------- fin fact ---------------------------------------------- 

%:     Definition{defCorpsdesRacines}--
\begin{definition}
\label{defCorpsdesRacines}
Soit $g\in\gK[T]$ un \poluz, on appelle \emph{\cdr de $g$ au dessus de $\gK$}  un \cdi
$\gL$ extension  de $\gK$  dans lequel~$g$ se décompose complètement
 et qui est engendré comme \Klg par les zéros de~$g$.%
\index{corps de racines!d'un \polz}
\end{definition}
%--- end-definition------------------------------------

Notez que $\gL$ est fini sur $\gK$ mais que l'on ne demande pas que $\gL$ soit \stf sur $\gK$
(d'ailleurs, il n'y a pas de \dem \cov qu'un tel \cdr doive être \stf sur $\gK$). Ceci nécessite quelques subtilités dans le \tho qui suit.

%:     theorem{propUnicCDR}
\begin{theorem} \emph{(Unicité du \cdr dans le cas \stfz)}
\label{propUnicCDR}\\
Soit $f\in\gK[T]$ un \poluz. On suppose qu'il existe un
\cdr  $\gL$ de $f$ au dessus de $\gK$.
\begin{enumerate}
\item Soit $\gM\supseteq\gK$ un \cdi \stf sur $\gK$,  engendré
par $\gK$ et des zéros de $f$ dans $\gM$.
 Le corps $\gM$ est isomorphe à un
sous-corps de~$\gL$.   
\item  Supposons qu'il existe un \cdr pour $f$  \stf sur~$\gK$. Alors,
 tout \cdr de $f$ au dessus de $\gK$ est isomorphe à $\gL$ (qui est donc \stf sur $\gK$).
\item Soient $\gK_1$, $\gK_2$ deux corps discrets non triviaux, $\tau : \gK_1 \to
\gK_2$ un \isoz, $f_1 \in \gK_1[T]$ un \poluz, $f_2 = f_1^\tau \in
\gK_2[T]$.  Si~$\gL_i$ est un corps de racines \stf pour $f_i$ sur $\gK_i$ ($i = 1,
2$), alors~$\tau$ se prolonge en un \iso de $\gL_1$ sur $\gL_2$.
\end{enumerate}
\end{theorem}
\begin{proof}
On montre seulement le point \emph{1} dans un cas particulier (suffisamment \gnlz). Le reste est laissé \alecz.\\
On écrit $f(T)=\prod_{i=1}^n(T-x_i)$ dans~$\gL[T]$. 
Supposons aussi
\hbox{que $\gM=\gK[y,z]$} avec $y\neq z$ \hbox{et $f(y)=f(z)=0$}. \\
On a donc dans $\gM[T]$ l'\egt $f(T)=(T-y)f_1(T)=(T-y)(T-z)f_2(T)$
(proposition \ref{propZerFactPol}). \\ 
Puisque $f(y)=0$, le \polmin $g(Y)$ de $y$
sur $\gK$ divise $f(Y)$ dans~$\gK[Y]$. Donc $\prod_{i=1}^ng(x_i)=0$ dans $\gL$
qui est un \cdiz, et l'un des $x_i$, disons $x_1$, annule $g$.  On obtient ici

\snic{
\gK[y] \simeq \aqo{\gK[Y]}{g(Y)}\simeq\gK[x_1] \subseteq\gL.
}

%\sni
Le \cdi $\gK[y]$ est \stf sur $\gK$ et  $\gM$ est \stf sur $\gK[y]$
(lemme \ref{lemdefiSTF}). Soit alors $h\in\gK[Y,Z]$ un \pol \unt en $Z$ tel que
$h(y,Z)$ soit le \polmin de $z$ sur $\gK[y]$. \\
Puisque $f_1(z)=0$, 
le \polz~$h(y,Z)$ divise $f_1(Z)=f(Z)/(Z-y)$ dans~$\gK[y][Z]$, 
donc son image $h(x_1,Z)$ dans~$\gK[x_1][Z]$ est un \pol \ird qui divise
$f(Z)/(Z-x_1)$.  Donc $h(x_1,Z)$ admet pour zéro un des $x_i$ 
pour $i\in\lrb{2..n}$, disons $x_2$, et $h(x_1,Z)$ est le \polmin de~$x_2$
sur $\gK[x_1]$. On
obtient donc les \isos

\snic{
\gK[y,z]\simeq\aqo{\gK[y][Z]}{h(y,Z)}\simeq\aqo{\gK[x_1][Z]}{h(x_1,Z)} \simeq
\gK[x_1,x_2]\subseteq\gL.
}

%\sni
Notons que l'on a aussi $\gK[y,z]\simeq\aqo{\gK[Y,Z]}{g(Y),h(Y,Z)}$.
\end{proof}

\rem Une inspection détaillée de la \dem précédente conduit à la conclusion que si $\gL$ est un \cdr \stf sur $\gK$, le groupe des 
$\gK$-\autos de $\gL$
est un groupe fini ayant au plus $[\gL:\gK]$ \eltsz. 
Si l'on ne suppose pas $\gL$ \stf sur $\gK$, on obtiendra seulement qu'il est absurde de supposer que ce groupe contient plus que~$[\gL:\gK]$ \eltsz.
\eoe

%:subsec{Le cas \elr}
\subsec{Le cas \elr de la théorie de Galois}
\label{subsecGaloisElr}

%:     DefiNota{NOTAStStp}
\begin{definota}\label{NOTAStStp}
Nous utiliserons les notations suivantes lorsqu'un groupe $G$
 opère sur un ensemble~$E$.
\begin{enumerate}\itemsep0pt
\item [---]
Pour $x\in E$, $\St_G(x)=\St(x)\eqdefi\sotq{\sigma\in G}{\sigma(x)=x}$
désigne le \ix{stabilisateur} de $x$.
\item [---] $G.x$ désigne l'orbite de $x$ sous $G$, et l'écriture
$G.x=\so{x_1,\ldots,x_k}$ est une abréviation pour: \emph{$(x_1,\ldots,x_k)$
est une énumération sans répétition de $G.x$, avec $x_1=x$.}
\item [---]
 Pour $F\subseteq E$, $\Stp_G(F)$ ou $\Stp(F)$ désigne le
stabilisateur point par point de~$F$.
\item [---]
 Si $H$ est un sous-groupe de $G$,
 \begin{enumerate}\itemsep0pt
\item [--] on note $\idg{G:H}$ l'indice de $H$ dans $G$,%
\index{indice!d'un sous-groupe dans un groupe}
\item [--] on note $\Fix_E(H)=\Fix(H)=E^H=\sotq{x\in E}{\forall \sigma\in H,\;\sigma(x)=x}$,
\item [--]  l'écriture $\sigma\in G/H$ signifie que l'on prend un \elt $\sigma\in G$
dans chaque classe à gauche modulo $H$.
\end{enumerate}
 \end{enumerate}
Lorsque $G$ est un groupe fini %d'\autos d'un
opérant sur un anneau
$\gB$
on note pour $b\in \gB$:
\[\preskip.3em \postskip.1em
\Tr_G(b)   =   \sum_{\sigma\in G}\sigma(b), \;
 \rN_G(b)  =   \prod_{\sigma\in G}\sigma(b), \hbox{ et }
\rC{G}(b)(T)   =    \prod_{\sigma\in G}\big(T-\sigma(b)\big).
\]
Et si $G.b=\so {b_1,\ldots ,b_k}$, (les $b_i$ deux à deux distincts), on note:
$$\preskip.2em \postskip.4em\ndsp 
\Rv_{G,b}(T)=\prod_{i=1}^k(T-b_i). 
$$
Ce \pol est appelé la \emph{résolvante} de $b$ (relativement à $G$).
Il est clair que $\Rv_{G,b}^r=\rC{G}(b)$ avec $r=\idG{G:\St_G(b)}$.%
\index{resolvante@résolvante}
\end{definota}

\rdb
\'Etant donnée une \Alg $\gB$ on note $\Aut_\gA(\gB)$ le groupe des $\gA$-\autos de $\gB$.
\label{NOTAAutAB}

%:     Definition{defiGalGal}
\begin{definition}\label{defiGalGal}
\label{i4defiCorGal}
Si $\gL$ est une extension \stfe de $\gK$, et un
\cdr  pour un \polu \spl sur $\gK$,
on dit que~$\gL$ est une \ix{extension galoisienne} de $\gK$, on note 
alors $\Gal(\gL/\gK)$ au lieu de~$\Aut_\gK(\gL)$ et l'on dit que c'est le \ix{groupe de Galois} de l'extension $\gL/\gK$. 
\end{definition}
%--------- fin definition ---------------------------------------------- 

Notez bien que dans la \dfn d'une extension galoisienne $\gL/\gK$, est compris
le fait que $\gL$ est \stf (et non pas seulement fini) sur $\gK$.

%\penalty-2500
%:     propdef{defiCorGal}
\begin{propdef}\label{defiCorGal} \emph{(Correspondance galoisienne)}\\
Soit $\gL\supseteq\gK$ un corps \stf sur~$\gK$.
\begin{enumerate}
\item  \label{i1defiCorGal}
  Le groupe
 $\Aut_\gK(\gL)$
est un sous-groupe détachable de $\GL_\gK(\gL)$.
\perso{Je pense qu'il n'y a pas de \dem \cov que $\Aut_\gK(\gL)$ soit fini. Sinon par exemple on aurait un \algo pour décider si une extension de degré 3 est galoisienne et en caractéristique 0 cela revient à décider si un discriminant est un carré?}
  Si $H$ est un sous-groupe de $\Aut_\gK(\gL)$, le sous-corps
$\gL^H$ s'appelle le \emph{corps fixe de~$H$}.
\item \label{i2defiCorGal}  On appelle \ix{correspondance galoisienne} les deux applications $\Fix$
et $\Stp$ entre les deux ensembles suivants. D'une part $\cG=\cG_{\gL/\gK}$ est l'ensemble des sous-groupes finis de $\Aut_\gK(\gL)$. D'autre part $\cK=\cK_{\gL/\gK}$
est l'ensemble des sous-extensions \stfes de~$\gL$.
\item \label{i3defiCorGal} Dans la correspondance galoisienne chacune des deux applications est décroissante.
En outre, $H\subseteq \Stp(\gL^H)$ pour tout~$H\in\cG$,
$\gM\subseteq \gL^{\Stp(\gM)}$ pour tout~$\gM\in\cK$,
$\Stp\circ\Fix\circ\Stp=\Stp$ et  $\Fix\circ\Stp\circ\Fix=\Fix$.
%
%\item 
%
\end{enumerate}
\end{propdef}
\begin{proof}
Dans le point \emph{\ref{i1defiCorGal}} il faut montrer que le sous-groupe est détachable
et dans le point \emph{\ref{i2defiCorGal}} que $\Fix$ et $\Stp$ agissent bien sur les deux ensembles tels qu'ils sont décrits. Ceci est basé sur l'\alg \lin en dimension finie sur les \cdisz. Nous laissons les détails \alecz.
\end{proof}

\rem
Bien que l'on puisse décider si un \elt
donné de $\GL_\gK(\gL)$ est
dans $\Aut_\gK(\gL)$, et bien qu'il soit facile de borner le nombre d'\elts de $\Aut_\gK(\gL)$,
il n'y a pas de méthode \gnle s\^ure pour calculer ce nombre.
\eoe

\medskip 
On a comme conséquence du \thrf{propUnicCDR} le corolaire suivant.

%:    Theorem{thPrIs}  prolongement des \isos
\begin{theorem}
\label{thPrIs} \emph{(\Tho de prolongement des \isosz)}\\
Soit  $\gL/\gK$ une extension galoisienne et $\gM$ une sous-$\gK$-extension finie de~$\gL$. Tout $\gK$-\homo
$\tau:\gM\to\gL$ se prolonge en un \eltz~$\wi\tau$ de~$\Gal(\gL/\gK)$.
\end{theorem}
%--- end-theorem-----------------------------------------
%
\begin{proof}
$\gL$ est le \cdr d'un \pol \spl $g\in\KT$.
On remarque que puisque $\gL$ est \stf sur $\gK$, $\gM$ est \stf sur $\gK$
et $\gL$ \stf sur $\gM$.
  Notons $\gM'$ l'image de $\tau$. C'est un corps \stf sur $\gK$, donc $\gL$
  est \stf sur $\gM'$.
Ainsi $\gL$ est un corps de racines
pour $g$  \stf sur $\gM$ et sur $\gM'$. D'après le \tho \ref{propUnicCDR} (point \emph{3}), on
peut prolonger $\tau$ en un $\gK$-\iso $\wi\tau : \gL \to \gL$.
%  On peut \gui{prolonger la photocopie} $\gM\vers{\tau_1}\gM'$ (où $\tau_1(x)=\tau(x)$ pour tout $x\in \gM$) en une
%\gui{photocopie} $\tau':\gL\to\gL'$, où $\gL'$ est une $\gM'$-\alg
%dont la structure recopie la structure de $\gL$ comme $\gM$-\algz.
%Alors $\gL$ et $\gL'$ sont deux \cdr pour $g$ sur $\gM'$,  \stfsz. Le \tho
%\ref{propUnicCDR} s'applique et fournit un $\gM'$-\iso $\sigma:\gL'\to\gL$.
%En composant on obtient un $\gK$-\auto convenable de $\gL$:
%$\wi\tau=\sigma\circ \tau'$.
\end{proof}

\rdb \label{NOTAGalKf}
 Lorsqu'un \pol \spl sur $\gK$ possède un \cdr $\gL$ \stf sur $\gK$,
 le groupe $\Gal(\gL/\gK)$ peut aussi être noté
$\Gal_\gK(f)$ dans la mesure où le \thrf{propUnicCDR} donne l'unicité
de $\gL$ (à $\gK$-\auto près).

\medskip
\rem
En \coma on a les résultats suivants (triviaux en \clamaz).
Pour un sous-groupe $H$ d'un groupe fini  \propeq
%-----------------begin item------------------
\begin{itemize}
\item  $H$  est fini.
\item  $H$  est \tfz.
\item  $H$  est détachable.
\end{itemize}
%-----------------end item------------------
De même pour un sous-\Kev $M$ d'un \Kev de dimension finie \propeq
%-----------------begin item------------------
\begin{itemize}
\item  $M$  est de dimension finie.
\item  $M$  est \tf (i.e., l'image d'une matrice).
\item  $M$  est le noyau d'une matrice.\eoe
\end{itemize}
%-----------------end item------------------

%:     Proposition{propGaloiselr}
\begin{propdef}\label{propGaloiselr}\emph{(Situation galoisienne \elrz)}\\
Soient deux anneaux $\gA\subseteq\gB$.
Une \emph{situation galoisienne \elrz} est définie comme suit.
\begin{enumerate}
\item [{i.}] On a un \polz~$Q\in\AT$  \mon \spl de degré $d$ et
des \elts $y_1$, $y_2$, $\ldots$, $y_d$  de $\gB$ tels que:

\snic{Q(T)=\prod_{i=1}^d (T-y _i).}

\item [{ii.}]  On note $y=y_1$. On suppose pour chaque $i$ que
$\gB=\gA[y _i]$ et que  $\gen{Q}$ est le noyau de l'\homo
 de \Algs $\AT\to\gB$ qui envoie $T$ en $y_i$
(d'où $\gB=\gA[y]=\gA[y _i]\simeq\aqo{\AT}{Q}$).
Pour chaque $i$ il existe donc un unique $\gA$-\auto $\sigma_i$ de $\gB$
vérifiant $\sigma_i(y)=y_i$.

\item [{iii.}]  On suppose que ces \autos forment un groupe,
que l'on note~$G$. En particulier, $\abs{G}=d=\dex{\gB:\gA}$.
\end{enumerate}
Dans une situation galoisienne \elr on a les résultats suivants.
\begin{enumerate}
\item
\begin{enumerate}
\item \label{i1propGaloiselr}
$\Fix_\gB(G)=\gA$.
\item \label{i2propGaloiselr}
Pour tout $z\in\gB$, $\rC{\gB/\gA}(z)(T)=\rC G(z)(T)$.
\end{enumerate}
\item \label{i3propGaloiselr}
Soit $H$ un sous-groupe détachable de $G$, $\gA'=\gB^H$ et

\snic{Q_H(T)=\prod_{\sigma\in H} \big(T-\sigma(y)\big).}

%\sni
Alors, on retrouve la situation galoisienne
\elr avec $\gA'$, $\gB$,  $Q_H$ et $\big(\sigma(y)\big)_{\sigma\in H}$.
En particulier, $\gB=\gA'[y]$ est un $\gA'$-module libre de rang $\abs{H}=\dex{\gB:\gA'}$.
En outre, $H$ est égal à $\Stp_G(\gA')$.
\end{enumerate}
\end{propdef}
\begin{proof}
\emph{1a.}
Considérons un $x=\sum_{k=0}^{d-1} \xi_ky ^k$ dans $\gB$
(avec les $\xi_k\in\gA$) invariant par l'action de
$G=\so{\sigma_1,\ldots,\sigma_{d}}$.
On a donc pour tout $\sigma\in G$, $x=\sum_{k=0}^{d-1} \xi_k\sigma(y )^k$.
Si~$V \in \Mn(\gB)$ est la matrice de Vandermonde
 $$V=\cmatrix {1& y_1  &y_1^2& \cdots&y_1^{d-1}\cr
 \vdots&&&&\vdots\cr
 \vdots&&&&\vdots\cr
 1& y_d&y_d^2& \cdots&y_d^{d-1}
 },$$
  on obtient
$$
V\ \left[\begin{array}{c} \xi _0 \\ \xi _1 \\ \vdots \\ \xi _{d-1} \end{array}\right]
=\left[\begin{array}{c} x \\ x \\ \vdots \\ x \end{array}\right]=
V\ \left[\begin{array}{c} x \\ 0 \\ \vdots \\ 0 \end{array}\right].
$$
Or, $\det(\!\tra{V}V)=\disc_T(Q)\in\Ati$. Donc
$[\,\xi_0 \ \xi_1 \ \cdots\ \xi_{d-1}\,]=[\,x \ 0\ \cdots\ 0\,]$, \hbox{et $x =\xi _0\in\gA$}.
\\
\emph{1b.}
Puisque $\gB\simeq\aqo{\AT}{Q}$, le corolaire \ref{corprop1tschir} donne,
pour $g\in\AY$ \hbox{et~$z=g(y_1)$}, les \egts
$$\preskip.0em \postskip.4em
\rC{\gB/\gA}(z)(T)=\prod\nolimits_i\big(T-g(y_i)\big) =
\prod\nolimits_{\sigma\in G}\big(T-\sigma(g(y_1))\big) = \rC G(z)(T).
$$
\emph{2.} Il est clair que $\gB=\gA'[\sigma(y)]$ pour chaque $\sigma\in H$
et que $Q_H$ est un \pol \spl de $\gA'[T]$. Il reste à voir
que tout \polz~$P\in\gA'[T]$ qui annule un $y_i=\sigma_i(y)$ ($\sigma_i\in H$)
est multiple de $Q_H$. Pour tout $\sigma\in H$, puisque~$\sigma$ est un
$\gA'$-\auto de $\gB$, on a~$P\big(\sigma(y_i)\big)=\sigma\big(P(y_i)\big)=0$. Ainsi~$P$ est divisible
par chacun des~$T-\sigma(y)$, pour $\sigma\in H$. Comme ces \pols sont
deux à deux \comz, $P$ est multiple de leur produit~$Q_H$. 
Enfin, si $\sigma_j\in G$ est un $\gA'$-\auto de $\gB$, $\sigma_j(y)=y_j$ doit être
un zéro de $Q_H$. Mais puisque $Q$ est \splz, les seuls $y_i$ qui
annulent $Q_H$ sont les $\sigma(y)$ pour $\sigma\in H$. Donc $\sigma_j\in H$.
\end{proof}

\rems 1) Dans la situation galoisienne \elr rien ne dit que les $y_i$
sont les seuls zéros de $Q$ dans $\gB$, ni que les $\sigma_i$ soient les seuls
\hbox{$\gA$-\autosz} de~$\gB$. Prenons par exemple $\gB=\gK^3$, et $a$, $b$, $c$  distincts
dans le \cdi $\gK$. Le \polz~$Q=(T-a)(T-b)(T-c)$ admet~27 zéros dans $\gB$,
dont six qui ont $Q$ pour \polminz, ce qui fait \hbox{six $\gK$-\autosz} de $\gB$.
\\
En outre, si l'on prend $z_1=(a,b,c)$, $z_2=(b,a,b)$ et $z_3=(c,c,a)$, 
on voit que~$Q=(T-z_1)(T-z_2)(T-z_3)$, ce qui montre que la première condition n'implique
pas la seconde. Par contre, avec $y_1=(a,b,c)$, $y_2=(b,c,a)$ \hbox{et $y_3=(c,a,b)$},
on est dans la situation galoisienne \elrz.

2) Concernant la condition \emph{iii} pour définir la situation
galoisienne \elrz, on voit facilement qu'elle
équivaut au fait que chaque $\sigma_i$ permute les $y_j$.
Cette condition n'est pas conséquence des deux premières comme le
prouve l'exemple qui suit. Considérons le carré latin $5 \times 5$ suivant
(dans chaque ligne et chaque colonne, les entiers sont différents), qui
n'est pas la table d'un groupe:
$$
\cmatrix {
1 & 2 & 3 & 4 & 5\cr
2 & 4 & 1 & 5 & 3\cr
3 & 5 & 4 & 2 & 1\cr
4 & 1 & 5 & 3 & 2\cr
5 & 3 & 2 & 1 & 4\cr
}.
$$
Chaque ligne définit une permutation $\sigma_i \in S_5$; ainsi $\sigma_1 =
\Id$, $\sigma_2 = (12453)$, \dots, $\sigma_5 = (154)(23)$.  Les $\sigma_i$ ne
forment pas un groupe (qui serait d'ordre $5$) car $\sigma_5$ est d'ordre $6$.
Posons $\gB = \gK^5$ où $\gK$ est un corps ayant au moins $5$ \elts
$a_1$, $\ldots$, $a_5$, $y = (a_1, \ldots, a_5) \in \gB$, $y_i = \sigma_i(y)$
et

\snic{Q(T) = \prod\nolimits_i (T - y_i) = \prod\nolimits_i (T - a_i) \in \gK[T].}

%\sni
Alors, dans \ref {propGaloiselr}, les deux premières conditions \emph{i},
\emph{ii} sont vérifiées mais pas la condition~\emph{iii}.\\
Fort heureusement les choses sont plus simples
dans le cas des corps. \eoe

%:     Lemma{lemGaloiselr}
\begin{lemma}\label{lemGaloiselr}
Soit $\gL=\gK[y]$ un \cdi \stf sur $\gK$.
Soit $Q$ le \polmin de $y$ sur $\gK$.
Si $Q$ est \spl et se factorise complètement dans $\gL[T]$, on se trouve dans la situation galoisienne \elr et
le groupe $G$ correspondant est le groupe $\Gal(\gL/\gK)$ de tous \hbox{les
$\gK$-\autosz} de~$\gL$.
\end{lemma}
\begin{proof}
Notons $y=y_1$, $\ldots$, $y_d$ les zéros de $Q$ (de degré $d$)
dans $\gL$.
Chaque $y_i$ annule $Q$ et $Q$ est \ird dans $\KT$, donc
$Q$ est le \polmin de $y_i$ sur $\gK$ et $\gK[y_i]$ est un sous-\Kev de $\gL$,
libre et de même dimension $d$, donc égal à $\gL$.
Enfin,
puisque $\gL$ est intègre, les $y_i$ sont les seuls zéros de $Q$
dans $\gL$, donc tout $\gK$-\auto de $\gL$ est un $\sigma_i$, et les $\sigma_i$
forment donc bien un groupe: le groupe de Galois $G=\Gal(\gL/\gK)$.
\end{proof}
%
%:     Theorem{thGaloiselr}
\begin{theorem} \emph{(Correspondance galoisienne, le cas \elrz)}
\label{thGaloiselr} \\
Soit $\gL=\gK[y]$ un \cdi \stf sur $\gK$.
Soit $Q$ le \polmin de $y$ sur $\gK$.
On suppose que $Q$ est \spl et se factorise
complètement dans $\gL[T]$. En particulier, $\gL$ est une extension
galoisienne de~$\gK$.
On a les résultats suivants.
\begin{enumerate}
\item  Les deux applications de la correspondance galoisienne sont
deux bijections réci\-proques l'une de l'autre.
\item  Pour tout $\gM\in \cK_{\gL/\gK}$, $\gL/\gM$ est une extension
galoisienne de groupe de Galois $\Fix(\gM)$ et $\dex{\gL:\gM}=\abs{\Fix(\gM)}$.
\item Si $H_1, H_2\in\cG_{\gL/\gK}$ et $\gM_i=\Fix(H_i)\in \cK_{\gL/\gK}$, alors:
\begin{itemize}
\item  $H_1\cap H_2$
correspond à la sous-\Klg engendrée par $\gM_1$ et $\gM_2$,
\item  $\gM_1\cap \gM_2$
correspond au sous-groupe engendré par $H_1$ et $H_2$.
\end{itemize}
\item Si $H_1\subseteq H_2$ dans $\cG_{\gL/\gK}$ et $\gM_i=\Fix(H_i)$, alors
$\gM_1\supseteq\gM_2$ et on a l'\egt 
$\idg{H_2:H_1}=\dex{\gM_1:\gM_2}$.
\item Pour tout $z\in\gL$, $\rC{\gL/\gK}(z)(T)=\rC{\Gal(\gL/\gK)}(z)(T)$.
\end{enumerate}
\end{theorem}
\begin{proof}
Il suffit de prouver le premier point.
D'après la proposition \ref{propGaloiselr} on a l'\egt $\Stp\circ \Fix=\Id_{\cG_{\gL/\gK}}$.\\
Soit maintenant $\gM\in \cK_{\gL/\gK}$. Puisque $\gL=\gK[y]$, on a $\gL=\gM[y]$.
Comme~$\gL$ est \stf sur~$\gM$, on peut calculer
le \polminz~$P$ de~$y$ sur~$\gM$. Il divise~$Q$
donc il est \splz. Il se factorise complètement dans~$\gL[T]$.
Ainsi, avec $\gM$, $\gL=\gM[y]$ et~$P$, on est dans les hypothèses du lemme
\ref{lemGaloiselr},
donc dans la situation galoisienne \elrz.  
Les $\gM$-\autos de~$\gL$ sont des $\gK$-\autos donc ce sont exactement les \elts du stabilisateur $H=\Stp_G(\gM)$
(où $G=\Gal(\gL/\gK)$).
Dans cette situation le point \emph{1b} de la proposition
\ref{propGaloiselr} nous dit que $\Fix(H)=\gM$.
\end{proof}

Nous venons d'établir que la correspondance galoisienne
est bijective,
ce qui est le \tho fondamental de la théorie de
Galois, dans le cas \elrz. Mais  ce cas est
en fait le cas \gui{\gnlz}: chaque fois que l'on a une extension galoisienne
on peut se ramener à la situation \elr (\thref{thResolUniv} et \tho de l'\elt primitif
\vref{thEtalePrimitif}).

%%%%%%%%%%%%%%%%%%%%%%%%%%%%%%%%%%%%%%%%%%%%%%%%%%%%%%%%%%%%%%%%%%%%%%%%%%%
%:subsec  Corps de racines
\subsect{Construction d'un \cdr au
moyen d'une résol\-vante de Galois, théorie de Galois de base}{Construction d'un \cdrz}
\label{secResolUniv}

\Grandcadre{Dans ce paragraphe
$f\in\gK[T]$ est un \polu \spl \\ de degré $n$
et
$\gA=\Adu_{\gK,f}$ avec $f(T)=\prod_i(T-x_i)$ dans $\gA$.}

Le but du paragraphe présent est de montrer le résultat suivant: si $\gK$ est infini, et si l'on sait factoriser les \polus \spls dans $\KT$, alors on sait construire un \cdr pour n'importe quel \polu \splz, et l'extension obtenue rentre dans le cadre
\elr du \thoz~\ref{thGaloiselr}.

Nous construisons ce \cdr par une méthode \gui{uniforme}.
Comme le \cdr construit est \stfz, le \thref{propUnicCDR} nous dit que 
tout autre \cdr lui est isomorphe.

%:     Theorem{thResolUniv}
\begin{theorem}\label{thResolUniv}
On introduit des \idtrs $u_1$,  \dots, $u_n$.
Pour $\sigma\in\Sn$ on note $u_\sigma=\sum_iu_ix_{\sigma i}$.
On pose
$$\preskip.3em \postskip.3em
R(\uu,T):=\prod\nolimits_{\sigma\in\Sn}(T-u_\sigma)\in\gK[\uu,T],
$$
et $D(\uu):=\disc_T(R)\in\Kuu$.
\begin{enumerate}
\item \label{i1thResolUniv} Un des \coes de~$D$ est égal à $\pm\disc(f)^{(n-2)!(n!-1)}$.
\end{enumerate}
Dans la suite on suppose que l'on spécialise les $u_i$
en des \elts $a_i\in\gK$
et que $D(\ua)\neq0$ (c'est toujours possible si $\gK$ est infini).
\begin{enumerate}\setcounter{enumi}{1}
\item \label{i2thResolUniv} Pour n'importe quel $\sigma\in\Sn$, l'\elt
 $a_\sigma=\sum_ia_ix_{\sigma i}$ admet le \pol $R(\ua,T)\in\KT$
pour \polminz, de sorte que
$$\preskip.4em \postskip.4em
\gA=\gK[a_\sigma]\simeq\aqo{\KT}{R(\ua,T)}.
$$
On note $a=a_\Id=\sum_i a_ix_i$.
\item \label{i3thResolUniv} Les seuls \elts de $\gA$ fixés par $\Sn$
sont les \elts de $\gK$.
\item  Supposons que
l'on sache décomposer $R(\ua,T)$ en produit de facteurs \irds
dans $\KT$:
$R(\ua,T)=\prod_{j=1}^\ell Q_j$.
\begin{enumerate}
\item \label{i4thResolUniv} Si $\ell=1$, $\gA$ est un corps,
l'extension $\gA/\gK$ est un \cdr pour le \pol $f$, ainsi
que pour $R(\ua,T)$, et la situation relève du \thrf{thGaloiselr}. En particulier,
$\Gal(\gA/\gK)\simeq\Sn.$
\item \label{i5thResolUniv} Si $\ell>1$, alors
 $\gA\simeq\prod_j\gK_j$ où 
$$\preskip.4em \postskip.4em 
\gK_j=\gK[\pi_j(a)]=\aqo{\gA}{Q_j(a)}\simeq\aqo{\KT}{Q_j}. 
$$
($\pi_j:\gA\to\gK_j$ est
la \prn canonique.) \\
Soit $H_j$ le sous-groupe  de $\Sn$
qui stabilise l'\id $\gen{Q_j(a)}_{\!\gA}$. Alors:
\begin{enumerate}
\item [--] $\Sn$ opère transitivement sur les \ids $\gen{Q_j(a)}_{\!\gA}$,
de sorte que les $Q_j$ ont tous même degré,
$\abs{H_j}=\deg(Q_j)=\dex{\gK_j:\gK}$,
et  les $\gK_j$ sont des \cdis deux à deux isomorphes,
\item [--] l'extension $\gK_1/\gK$ est un \cdr pour $f$, ainsi
que pour chacun des $Q_j$, et  la situation relève du \tho \ref{thGaloiselr}, en particulier,
$H_1=\Gal(\gK_1/\gK)$.
\end{enumerate}
\end{enumerate}
%  le truc suivant peut etre en exo ?
%\item \label{i6thResolUniv} Réciproquement supposons que l'on connaisse un corps
%$\gL\supseteq\gK$ tel que:
%\begin{itemize}
%%
%\item $f(T)=\prod_{i=1}^n(T-\xi_i)$ dans $\gL[T]$,
%%
%\item $\gL$ est un \Kev de dimension finie,
%%
%\item $\gL=\gK[\xin]$.
%\end{itemize}
%Alors on sait factoriser $R(\ua,T)$ dans $\gK[T]$ et $\gL\simeq\gK_1$
%(avec la notation du point précédent).
\end{enumerate}
\end{theorem}
\begin{proof}
\emph{\ref{i1thResolUniv}.} Le discriminant $D$ est égal (au signe près)
au produit des $u_\sigma-u_\tau$ pour $\sigma\neq\tau\in\Sn$.
Chaque  $u_\sigma-u_\tau$ est une somme d'\elts $u_i(x_{\sigma i}-x_{\tau i})$:
chaque $u_i$ a pour \coe $0$ ou un $x_j-x_k$ ($j\neq k$).
Le premier \mom pour l'ordre lexicographique qui apparaît
dans le produit~$D$ est le \mom
$$\preskip.3em \postskip.4em
u_1^{n!(n!-(n-1)!)}u_2^{n!\big((n-1)!-(n-2)!\big)}\cdots u_{n-1}^{n!(2!-1!)},
$$
avec pour \coe un produit d'\elts du type $x_i-x_j$  ($i\neq j$). Plus \prmt si
$\delta=\disc(f)$, le \coe en question sera, au signe près,

\snic{ \delta^{(n-2)!(n!-1)}.}

%\sni
\emph{\ref{i2thResolUniv}.} On  utilise la proposition \ref{propdiscTra} puisque $R(\ua,T)$ est le \polcar de  $a$ (lemme \ref{lemPolCarAdu}).
\\
\emph{\ref{i3thResolUniv}.}
Voir le point \emph{1b} de la proposition
\ref{propGaloiselr}.\\
\emph{4a.} C'est clair.\\
\emph{4b.} Le fait que $\gA\simeq\prod_j\gK_j$
résulte du \tho des restes chinois.\\
L'\egt $\prod_jQ_j(T)=\prod_{\sigma}(T-a_{\sigma})$ dans $\AT$
 reste valable dans
$\gK_1[T]$. \\
Donc, il existe pour tout $j$ un $\sigma_j$ tel que
$Q_j\big(\pi_1(a_{\sigma_j})\big)=0$, autrement 
\linebreak 
dit, $Q_j(a_{\sigma_j})\in \gen{Q_1(a)}_{\!\gA}$. D'autre part, dans
$\gA$ on a $Q_j(a_{\sigma_j})=\sigma_j\big(Q_j(a)\big)$ parce que
$Q_j\in\KT$. Donc
$\sigma_j\big(\gen{Q_j(a)}_{\!\gA}\big)\subseteq \gen{Q_1(a)}_{\!\gA}$.
\\
Ceci nous donne une surjection
$\sigma_j : \aqo\gA {Q_j(a)}   \to   \aqo\gA {Q_1(a)}$,
i.e. une surjection  $\aqo\KT {Q_j} \to   \aqo\KT {Q_1}$.
Il en résulte $\deg Q_1\leq\deg Q_j$, et par symétrie $\deg Q_j=\deg Q_1$,
d'où $\sigma_j\big(\gen{Q_j(a)}_{\!\gA}\big)= \gen{Q_1(a)}_{\!\gA}$.
\\
Ainsi $\Sn$ opère transitivement sur les \ids $\gen{Q_j(a)}_{\!\gA}$
et les $\gK_j$ sont deux à deux isomorphes.
\end{proof}

\rem La construction du \cdr suggérée ici
est en fait à peu près
impraticable dès que le degré $n$ de $f$ est supérieur ou égal à 7, car elle nécessite de factoriser un \pol de degré $n!$.
Nous proposons dans le chapitre \ref{ChapGalois} une méthode dynamique moins
brutale qui a l'avantage supplémentaire de ne pas réclamer de savoir factoriser les \pols \spls de $\KT$.
La contrepartie de cette absence de factorisation sera que, bien que l'on sache calculer dans \gui{un}
\cdrz, on n'est a priori jamais certain de le connaître
de manière complète (au sens où on connaîtrait sa dimension comme
\Kevz). En outre, le même manque de précision se retrouve pour ce qui concerne le groupe de Galois.
\eoe

\medskip
\exl On considère le \pol {\tt p(T)}$\in\QQ[T]$ ci-dessous.
On demande à {\tt Magma} de prendre au hasard une \coli {\tt z} des {\tt xi} (les zéros \hbox{de {\tt p(T)}}  dans l'\adu $\gA=\Adu_{\QQ,p}$), de calculer $\Mip_{\QQ,z}(T)$, puis de le factoriser.
Le logiciel donne rapidement le \polmin {\tt pm} de degré $720$
et le décompose en un produit de $30$ facteurs de degré $24$ (le tout
en une ou deux minutes).
Un de ces facteurs est le \pol {\tt q}. Comme {\tt q} est très encombrant,
on demande à   {\tt Magma} de calculer une \bdg de l'idéal
engendré par les modules de Cauchy d'une part, et par {\tt q(z)} d'autre
part, ce qui fournit une description plus claire du \cdr $\aqo{\gA}{q(z)}$:
{\tt x6} est annulé par {\tt p}, {\tt x5} est annulé par un \pol
de degré $4$ à \coes dans $\QQ[x_6]$, {\tt x1}, \ldots, {\tt x4}
s'expriment en fonction de {\tt x5} et {\tt x6}. Le calcul de la \bdg prend
plusieurs heures. {\tt Magma} peut ensuite calculer le groupe de Galois,
qui est donné par deux \gtrsz.
Voici les résultats:

%\vspace{-.2cm}
%:     verbatim
{\small \def\baselinestretch{1}\label{exemple1Galois}
\begin{verbatim}
p:=T^6 - 3*T^5 + 6*T^4 - 7*T^3 + 2*T^2 + T - 1;
z:=x1 + 2*x2 + 13*x3 - 24*x4 + 35*x5 - 436*x6;
pm:=T^720 + 147240*T^719 + 10877951340*T^718 + 537614218119000*T^717 +
    19994843992714365210*T^716 + 596880113924932859498208*T^715 +
    14896247531385087685472255280*T^714 + ...
q:= T^24 + 4908*T^23 + 13278966*T^22 + 25122595960*T^21 +
    36160999067785*T^20 + 41348091425849608*T^19 +
    38304456918334801182*T^18 + 28901611463650323108996*T^17 +...
//on annule  q(z): description du corps des racines;
Affine Algebra of rank 6 over Rational Field
Variables: x1, x2, x3, x4, x5, x6
Quotient relations:
x1 + 18/37*x5^3*x6^5 - 45/37*x5^3*x6^4 + 104/37*x5^3*x6^3 - 3*x5^3*x6^2
    + 36/37*x5^3*x6 - 1/37*x5^3 - 27/37*x5^2*x6^5 + 135/74*x5^2*x6^4 -
    156/37*x5^2*x6^3 + 9/2*x5^2*x6^2 - 54/37*x5^2*x6 + 3/74*x5^2 +
    91/37*x5*x6^5 - 455/74*x5*x6^4 + 460/37*x5*x6^3 - 25/2*x5*x6^2 +
    108/37*x5*x6 + 31/74*x5 - 41/37*x6^5 + 205/74*x6^4 - 204/37*x6^3 +
    11/2*x6^2 - 45/37*x6 - 53/74,
x2 + x6 - 1,
x3 + x5 - 1,
x4 - 18/37*x5^3*x6^5 + 45/37*x5^3*x6^4 - 104/37*x5^3*x6^3 + 3*x5^3*x6^2
    - 36/37*x5^3*x6 + 1/37*x5^3 + 27/37*x5^2*x6^5 - 135/74*x5^2*x6^4 +
    156/37*x5^2*x6^3 - 9/2*x5^2*x6^2 + 54/37*x5^2*x6 - 3/74*x5^2 -
    91/37*x5*x6^5 + 455/74*x5*x6^4 - 460/37*x5*x6^3 + 25/2*x5*x6^2 -
    108/37*x5*x6 - 31/74*x5 + 41/37*x6^5 - 205/74*x6^4 + 204/37*x6^3 -
    11/2*x6^2 + 45/37*x6 - 21/74,
x5^4 - 2*x5^3 + x5^2*x6^2 - x5^2*x6 + 4*x5^2 - x5*x6^2 + x5*x6 - 3*x5 +
     x6^4 - 2*x6^3 + 4*x6^2 - 3*x6 - 1,
x6^6 - 3*x6^5 + 6*x6^4 - 7*x6^3 + 2*x6^2 + x6 - 1
// le groupe de Galois;
Permutation group acting on a set of cardinality 6
Order = 24 = 2^3 * 3
    (1, 4)(2, 5)(3, 6)
    (1, 2, 4, 6)
\end{verbatim}
}

%\vspace{-.2cm}
On notera que $\disc_T(p)=2^4\times 37^3$, ce qui n'est pas sans rapport
avec les \denos apparaissant dans la \bdgz. L'exemple sera repris
\paref{exemple2Galois} avec la méthode dynamique.
\eoe

\medskip
\rem
Nous interrompons ici le traitement de la théorie de Galois de base. Nous
reprendrons le fil de ces idées dans les sections \ref{secEtaleSurCD}
et \ref{sec2GaloisElr}
qui peuvent être lues directement ici
(les résultats des chapitres intermédiaires ne seront pas utilisés).
Dans le chapitre \ref{ChapGalois} nous aborderons une théorie plus sophistiquée
qui s'avère \ncr lorsque l'on ne dispose d'aucun \algo de factorisation des \pols \spls sur le corps de base.
\eoe

%%%%%%%%%%%%%%%%%%%%%%%%%%%%%%%%%%%%%%%%%%%%%%%%%%%%%%%%%%%%%%%%%%%%%%%%%%%
\section{Le résultant}
\label{secRes}

Le résultant est l'outil de base de la théorie de l'\eliz.
Ceci est basé sur le lemme d'\eli de base \ref{LemElimAffBasic} 
qui s'applique avec  n'importe quel anneau et sur son corolaire~\ref{cor2LemElimAffBasic}
pour le cas  \gmqz.

%:  subsec   La théorie de l'\eli
\subsec{La théorie de l'\eliz}

La théorie de l'\eli s'intéresse aux \syss d'équations polynomiales (ou \emph{\sypsz}).\index{systeme polynomial@\sypz}\index{elimination@\eliz!théorie de l'---}

Un tel \sys $(\lfs)$ dans $\kXn=\kuX$, où $\gk$ est un \cdiz, peut admettre des
zéros dans $\gk^{n}$, ou dans $\gL^{n}$, avec $\gL$ un surcorps de $\gk$,
ou même $\gL$  une \klg arbitraire.
Les zéros dépendent seulement de l'\id $\fa=\gen{\lfs}$ de $\kuX$ engendré par 
les $f_i$. Aussi on les appelle \emph{les zéros de l'\id $\fa$}.

Soit $\pi:\gL^{n}\to\gL^{r}$ la projection qui oublie les $n-r$ dernières \coosz.
Si $V\subseteq \gL^{n}$ est l'ensemble des zéros de $\fa$ sur $\gL$, on est intéressé par une description aussi précise que possible de la projection  $W=\pi(V)$. Si possible comme zéros d'un \syp en les variables $(X_1,\dots,X_{r})$.

Ici intervient de manière naturelle l'\emph{\id d'\eliz} (\eli des variables $X_{r+1}$, \dots, $X_n$ pour le \syp considéré), qui est défini par $\fb=\fa\cap\kXr$.
En effet tout \elt de $W$ est clairement un zéro de $\fb$.

La réciproque n'est pas toujours vraie (et de toute manière pas du tout évidente), mais elle est vraie dans certains bons cas: si $\gL$ est un \cac et si l'\id est en position de \Noe  
(\thref{thNstfaibleClass}).

Un fait rassurant, et facile à établir par des considérations d'\alg \lin sur les \cdisz, est que l'\id d'\eli $\fb$ \gui{ne dépend pas}
du corps de base $\gk$ considéré. Plus \prmtz, si $\gk_1$ est un surcorps
de $\gk$, on a les résultats suivants.
\begin{itemize}
\item L'\id $\gen{\lfs}_{\gk_1[\Xn]}$ ne dépend que de l'\id $\fa$: \\
c'est l'\id $\fa_1$ de~\hbox{$\gk_1[\Xn]$} engendré par $\fa$. 
\item L'\id d'\eli $\fb_1=\fa_1\cap\gk_1[\Xr]$ ne dépend que de $\fb$: 
\\
c'est l'\id  de $\gk_1[\Xr]$ engendré par $\fb$. 
\end{itemize}

\medskip 
La théorie \elr de l'\eli se heurte à deux obstacles.

Le premier est la difficulté de calculer $\fb$ à partir de $\fa$,
\cad de calculer un \sgr fini de $\fb$ à partir du \syp $(\lfs)$.
Ce calcul est rendu possible par la théorie des \bdgsz, que nous n'aborderons pas dans l'ouvrage. En outre ce calcul n'est pas uniforme, contrairement aux calculs liés à la théorie du résultant.

Le deuxième obstacle, c'est que les choses ne se passent de manière vraiment satisfaisante qu'avec les \syps \hmgsz.
L'exemple de base qui montre ceci est le \deterz. 
On considère un \sli générique $(f_1,\dots,f_n)$ de $\gk[\ua][\uX]$,
où les variables $a_{ij}$ dans $\ua$ représentent les~$n^{2}$ \coes des $n$ formes \lins $f_i$, et les $X_j$ sont les inconnues.
Alors l'\id $\gen{\det(\ua)}$ de $\gk[\ua]$ est bien l'\id d'\eli des variables $X_j$ pour le \sys $(f_1,\dots,f_n)$, mais à condition de ne prendre en compte que les zéros du \sys distincts de 
$\uze=(0,\dots,0)$, \cad de se situer dans un cadre entièrement \hmgz.
\\
La simplicité du résultat est à mettre en regard avec la complication
de la discussion, dans le cadre non \hmgz, pour les \syss où les $f_i$ sont des formes affines.
\\
D'autre part,  bien que les zéros de l'\id $\gen{\det(\ua)}$ 
correspondent effectivement aux \syss qui admettent un zéro  $\neq \uze$,
cet \id n'est pas exactement égal à $\gen{f_1,\dots,f_n}\cap \gk[\ua]$, il faut d'abord \emph{saturer} $\fa=\gen{f_1,\dots,f_n}$ par rapport aux variables \hmgs $X_j$, \cad lui rajouter tous les $g$ tels que, pour chaque 
$j\in\lrbn$, $gX_j^{N}\in\fa$ pour un $N$ assez grand. 
Dans le cas présent, ce saturé est l'\id $\fa+\det(\ua)\gk[\ua][\uX]$,  chaque
$\det(\ua)X_j$ est dans $\fa$, et l'intersection du saturé avec $\gk[\ua]$
est bien $\gen{\det(\ua)}$.    

Nous retiendrons de cette petite introduction à la théorie de l'\eli  une \dfnz:
soient $\gk$ un anneau commutatif, $\fa$ un \id de $\kXn$ \hbox{et $r\in\lrb{0..n-1}$}, on définit  l'\emph{\id d'\eli des variables $X_{r+1}$, \dots, $X_n$ pour l'\id $\fa$} comme étant l'\id $\fb=\fa\cap\kXr$.%
\index{elimination@\eliz!ideal@\id d'---}\index{ideal@idéal!d'elimi@d'\eliz}

On prendra garde au fait que si $\gk$ est un anneau arbitraire, l'\id $\fa$
peut très bien être \tf sans que $\fb$ le soit.

%:  subsec   la matrice de Sylvester
\subsec{La matrice de Sylvester}

Dans ce qui suit, on ne suppose pas l'anneau $\gA$ discret,
si bien que le degré
d'un \pol de $\AX $ n'est pas \ncrt connu de manière exacte.
Du point de vue du calcul, on doit en \gnl
prendre les \pols dans~$\AX $
sous forme de \emph{\pols formels} \cad des couples $(f,p)$ où $f$  est un \pol
et $p$ majore son degré. Cette notion est \egmt utile en cas de
changement d'anneau de base, car un \pol peut voir par exemple son degré baisser sans que l'on sache le tester
(lors d'un passage dans un quotient par exemple).%
\index{polynome@\pol!formel}

Rappelons la \dfn de la matrice de Sylvester et du résultant de deux \pols (des \pols formels
de degrés $p$ et $q\geq0$):
$$\arraycolsep2pt
\begin{array}{lcl} 
f & = & a_pX^p+ a_{p-1}X^{p-1}+\cdots+a_0, \\[1mm] 
g & = &b_{q}X^{q}+b_{q-1}X^{q-1}+\cdots+b_0 .
\end{array}
$$
La \emph{matrice de Sylvester} de $f$ et $g$ (en degrés $p$ et $q$) est la
matrice
suivante\index{matrice!de Sylvester}\index{Sylvester!matrice de ---}\rdb
$$\label{SylvMat}
\Syl_X(f,p,g,q)=
\underbrace{
\left[
\begin{array}{ccccccccccc}
 a_p&\cdots&\cdots&\cdots&\cdots&a_0&& \\[.5mm]
&\ddots& & & & &\ddots& \\[.5mm]
 &&a_p&\cdots&\cdots&\cdots&\cdots&a_0\\[.8mm]
b_{q}&\cdots& \cdots &b_0&&&& \\[.5mm]
 &\ddots& & &\ddots&&& \\[3.5mm]
 &&&\ddots& & &\ddots& \\[.5mm]
 &&&&b_{q}&\cdots&\cdots& b_0  \\[.5mm]
\end{array}
 \right]
}_{p+q}
\matrix{ \left.\matrix{\cr \cr \cr \cr}\right\}&{q} \cr
 \left.\matrix{\cr \cr \cr \cr \cr\cr}\right\}&{p}
\cr}
$$
Cette matrice peut être vue comme la matrice dont les lignes
sont les \coos des \pols $X^{q-1}f$, $\ldots $, $Xf$, $f$, $X^{p-1}g$, $\ldots $, $Xg$, $g$
sur la \hbox{base
$(X^{p+q-1},X^{p+q-2},\ldots ,X,1)$}.

\rdb

Le \emph{résultant de $f$ et $g$ (en degrés $p$ et $q$)}, noté
$\Res_{X}(f,p,g,q)$, est le \deter de cette
matrice de Sylvester:\index{resulant@résultant!de deux \polsz}
%------begin equation--eqResultant-----------
\begin{equation}\label{eqResultant}
\Res_{X}(f,p,g,q)\eqdefi\det\big(\Syl_X(f,p,g,q)\big).
\end{equation}
%---------------------end equation--------------
Si le contexte est clair, on note aussi $\Res_X(f,g)$
ou $\Res(f,g)$.
On a
\begin{equation}
\label{eqres1.0}\relax
{\Res_{X}(f,p,g,q)=(-1)^{pq}\Res_{X}(g,q,f,p)},
\end{equation}
et aussi, pour $a,b\in\gA$,
\begin{equation}
\label{eqres1.1}\relax
\Res_{X}(af,p,bg,q) = a^qb^p\Res_{X}(f,p,g,q).
\end{equation}

Si $p=q=0$, on obtient le \deter d'une matrice vide, \cad $1$. 

%:     Fact{fact1Res}
\begin{fact}\label{fact1Res} 
Si $p\geq 1$ ou $q\geq1$, alors
$\Res_X(f,p,g,q)\in\gen{f,g}_{\AX }\cap\gA$. Plus \prmtz,
pour chaque $n\in\lrb{0..p+q-1}$, il existe
 $u_n$ et $v_n\in\AX$ tels que $\deg u_n<q$, $\deg v_n<p$ et
%------begin equation--eqResultant2-----------
\begin{equation}\label{eqResultant2}
X^n\,\Res_X(f,g)=u_n(X)f(X)+v_n(X)g(X).
\end{equation}
%---------------------end equation--------------
\end{fact}
\begin{proof}
Soit  $S$ la transposée de $\Syl_X(f,p,g,q)$. Les colonnes de $S$
expriment des \pols $X^kf$ ou $X^\ell g$ sur la base des \moms de degrés $<p+q$. En utilisant la formule de Cramer

\snic{S \,\wi S=\det S\cdot\rI_{p+q}\;,}

%\sni
on voit que chaque  $X^n\Res(f,g)$ (qui correspond à l'une des colonnes de la matrice du second membre)  est une \coli
des colonnes de~$S$.
\end{proof}

\rem On peut aussi voir l'\egrf{eqResultant2} dans le cas $n=0$ comme exprimant le
\deter de la matrice ci-dessous développé selon la dernière colonne
(il s'agit de la matrice de Sylvester dans laquelle on a remplacé dans la dernière
colonne chaque \coe par le \gui{nom} de sa ligne):
$$
\left[
\begin{array}{ccccccccccc}
 a_p&\cdots&\cdots&\cdots&\cdots&a_0&& X^{q-1}f\\[.5mm]
&\ddots& & & & &\ddots& \\[.5mm]
 &&a_p&\cdots&\cdots&\cdots&\cdots&f\\[2mm]
b_q&\cdots& \cdots &b_0&&&&X^{p-1}g \\[.5mm]
 &\ddots& &&\ddots&&& \\[5.5mm]
 &&&b_q&\cdots&\cdots&b_0& Xg\\[1.5mm]
 &&&&b_{q}&\cdots&\cdots& g  %\\[.5mm]
\end{array}
 \right]. $$
\eoe

%:     Corollary{lem0Resultant}
\begin{corollary}\label{lem0Resultant}
Soient $f$, $g\in\AX$ et $a\in\gB\supseteq\gA$, avec $f(a)=g(a)=0$,
et~$p\geq1$ ou~$q\geq1$, alors $\Res_X(f,p,g,q)=0$.
\end{corollary}

Notez que si les deux degrés sont surévalués le résultant s'annule,
et l'interprétation intuitive est que les deux \pols ont un zéro en commun \gui{à l'infini}.
Tandis que si $a_p=1$, le résultant (pour $f$ en degré $p$)
est le même quel que soit le degré formel  choisi pour $g$.
Ceci permet alors de passer sans ambiguïté à la notation $\Res(f,g)$,
comme dans le lemme suivant.

%:     Lemma{lemResultant}----------
\begin{lemma}
\label{lemResultant}
Soient $f$ et $g\in\AX$ avec $f$ \mon de degré $p$.
\begin{enumerate}
\item On note $\gB=\aqo{\AX}{f}$ et $\mu_g$ la multiplication par
(la classe de) $g$ dans~$\gB$, qui est un \Amo libre de rang $p$.
Alors:
\begin{equation}
\label{eq0lemResultant}
\rN\iBA (g)=\det\mu_g=\Res(f,g).
\end{equation}
\item Par suite
%------begin equation--eqResultant4-----------
\begin{eqnarray}\label{eqResultant4}
\Res(f,gh)&=&\Res(f,g)\,\Res(f,h), \\ %\quad \mathrm{et}\quad
\Res(f,g+fh)&=&\Res(f,g).\label{eqResultant5}
\end{eqnarray}
%---------------------end equation--------------
%
\item Pour toute matrice carrée $A\in\MM_p(\gA)$
dont le \polcar est égal à $f$, on a
\begin{equation}
\label{eqlemResultant}
\Res(f,g)=\det\big(g(A)\big).
\end{equation}
\item Si l'on écrit  $f=\prod_{i=1}^p(X-x_i)$
dans une extension de $\gA$, on obtient
%------begin equation--eqResultant3-----------
\begin{equation}\label{eqResultant3}
\Res(f,g)= \prod\nolimits_{i=1}^pg(x_i).
\end{equation}
%---------------------end equation--------------
%
\end{enumerate}
\end{lemma}
%--- end-lemma-----------------------------------------
%-----------------begin proof------------------
\begin{proof}
\emph{1.} Par manipulations \elrs de lignes,
la matrice de Sylvester
$$
\Syl_X(f,p,g,q)=
\left[
\begin{array}{ccccccccccc}
 1&a_{p-1}&\cdots&\cdots&\cdots&a_0&& \\[.5mm]
&\ddots&\ddots& & & &\ddots& \\[.5mm]
 &&1&a_{p-1}&\cdots&\cdots&\cdots&a_0\\[2mm]
b_{q}&\cdots& \cdots &b_0&&&& \\[.5mm]
 &\ddots& & &\ddots&&& \\[3.5mm]
 &&&\ddots& & &\ddots& \\[.5mm]
 &&&&b_{q}&\cdots&\cdots& b_0
\end{array}
 \right]
$$
est transformée en la matrice qui est visualisée ci-après,
dans laquelle les 
lignes~\hbox{$q+1$, $\ldots$, $q+p$} contiennent maintenant les restes
de la division par $f$ des \pols $X^{p-1}g$, $\ldots$, $Xg$, $g$.
Ainsi la matrice $p\times p$ dans le coin sud-est, est exactement
la matrice transposée de la matrice de l'\endo $\mu_g$ de $\gB$
sur la base des \momsz.
Et son \deter est égal à celui de la matrice de Sylvester.
$$
\left[
\begin{array}{ccccccccccc}
 1&a_{p-1}&\cdots&\cdots&\cdots&a_0&& \\[.5mm]
&\ddots&\ddots& & & &\ddots& \\[.5mm]
 &&1&a_{p-1}&\cdots&\cdots&\cdots&a_0\\[2mm]
0&\cdots& 0 &\times &\cdots&\cdots&\cdots& \times \\[.5mm]
\vdots &&\vdots&\vdots& &&&\vdots \\[3.5mm]
\vdots &&\vdots&\vdots& &&&\vdots \\[.5mm]
0&\cdots& 0 &\times &\cdots&\cdots&\cdots& \times
\end{array}
 \right]
$$

\emph{2.} Résulte du point \emph{1.}
\\
\emph{3}
et \emph{4.} Résultent de la proposition \ref{prop2tschir} via le point~\emph{1.}

On peut aussi donner les preuves directes suivantes.

\emph{4.} Tout d'abord, de l'\eqrf{eqResultant4} on déduit la formule \smq
\[
 \Res(f_1f_2,g)=\Res(f_1,g)\;\Res(f_2,g)
\] 
pour $f_1$ et $f_2$ \mons (utiliser
les équations (\ref{eqres1.0}) et (\ref{eqres1.1}) et le fait que
dans le cas où les \coes de $g$ sont des \idtrs on peut supposer $g=b_qg_1$
avec $g_1$ \monz). Ensuite, un calcul direct donne $\Res(X-a,g)=g(a)$.

\emph{3.} On doit démontrer $\Res(\rC A,g)=\det\big(g(A)\big)$ pour un \polz~$g$ et une matrice $A$ arbitraires. Il s'agit d'une \ida concernant les \coes de~$A$
et de~$g$. On peut donc
se limiter au cas où la matrice $A$ est la matrice générique. Alors, elle se diagonalise dans un suranneau et l'on conclut 
en appliquant le point~\emph{4.}
\end{proof}
%-----------------end proof------------------

\rem \label{remreciproquenonnegligeable} Le point \emph{4} offre une réciproque non négligeable au corolaire~\ref{lem0Resultant}:
si $\gA$ est intègre et si $f$ et $g$ sont deux \polus de~$\AT$
qui se factorisent complètement dans un anneau intègre contenant~$\gA$, ils ont un zéro commun
\ssi leur résultant est nul. \eoe

Dans le cas d'un \cdi non trivial $\gK$ on a un peu mieux. 

%:     Fact{factResNulPGCD}
\begin{fact}\label{factResNulPGCD}
Soient $f$ et $g\in\KX$ de degrés $p$ et $q\geq1$, avec $\Res(f,g)=0$. Alors, $f$ et $g$ ont un pgcd de degré $\geq1$.
\end{fact}
\begin{proof}
L'\Kli $(u,v)\mapsto uf+vg$ où $\deg u<q$ et $\deg v<p$ admet pour matrice
sur les bases de \moms la transposée de la matrice de Sylvester.
Soit donc $(u,v)\neq(0,0)$ dans le noyau. Le \polz~$uf=-vg$ est de degré
$<p+q$. Donc $\deg \!\big(\ppcm(f,g)\big)<p+q$, et cela implique $\deg\! \big(\pgcd(f,g)\big)>0$.
\end{proof}

\comm La \dem ci-dessus suppose que l'on connaisse la théorie \elr de la \dve (via l'\algo d'Euclide) dans les anneaux du type $\KX$.
Cette théorie montre l'existence d'un pgcd et d'un ppcm avec la relation

\snic{
\qquad\qquad\ppcm(f,g)\, \pgcd(f,g)=\alpha fg, \qquad (\alpha\in\gK\eti).}

%\sni
Une autre preuve consisterait à dire que dans un \cdi $\gL$
qui est une extension de $\gK$, les \pols $f$ et $g$ sont scindés
(i.e., se décomposent en facteurs de degré $1$) ce qui implique,
vu la remarque qui était faite juste avant, que $f$ et $g$
ont un zéro commun et donc un facteur commun de \hbox{degré $>0$}.
Il faut ensuite terminer en remarquant que le pgcd se calcule par l'\algo
d'Euclide et ne dépend donc pas du corps de base choisi (qui doit seulement
contenir les \coes de $f$ et $g$). Néanmoins cette seconde \demz,
qui en quelque sorte donne \gui{la vraie raison du \thoz}
suppose l'existence de
$\gL$ (qui n'est pas garantie d'un point de vue \cofz) et elle n'évite
nullement la théorie de la \dve dans~$\KX$ via l'\algo d'Euclide.
\eoe

%:  --- Lemme{LemElimAffBasic}--------
\CMnewtheorem{lemelibas}{Lemme d'\eli de base}{\itshape}
%      Lemma{LemElimAffBasic}
\begin{lemelibas}\label{LemElimAffBasic}%
\index{Lemme d'elimination de base@Lemme d'\eli de base}%
\index{elimination@\eliz!lemme d'--- de base}%
\index{elimination@\eliz!ideal@\id d'---}\\
Soient $f$ et $g\in\AX $ avec $f$ \mon de degré $p$. Alors,
$R=\Res_X(f,g)$ est bien défini et l'\id d'\eli  $\fa=\gen{f,g}_{\AX}\cap\gA$
vérifie

\snic{
\fa^p\; \subseteq\; \Res_X(f,g)\gA\;\subseteq\; \fa.
}

%\sni
En particulier:
\begin{enumerate}
\item $R$ est \iv \ssi $1\in\gen{f,g}$,
\item $R$ est \ndz
\ssi
$\fa$ est fidèle, et
\item $R$ est nilpotent
\ssi  $\fa$ est nilpotent.
\end{enumerate}
 
\end{lemelibas}
%--------- fin lemma ---------------------------------------------- 

%-----------------begin proof------------------
\begin{proof}
On sait déjà que $\Res_X(f,g)\in\gen{f,g}_{\AX}$.\\
On reprend les notations du lemme \ref{lemResultant}, point \emph{1.}
On note $x$ la classe de~$X$ dans~$\gB=\aqo{\AX}{f}$.
Une base de $\gB$ sur $\gA$ est $(1,x,\dots,x^{p-1})$.
Soient $(\gamma_i)_{i\in\lrbp} $ des \elts de $\fa$. Les \elts $\gamma_1$, $\gamma_2x$, $\ldots$, $\gamma_px^{p-1}$ sont dans~$\Im\mu_g$,
 donc la matrice $D=\Diag(\gamma_1,\ldots,\gamma_p)$ peut s'écrire
 sous la forme $GB$, où~$G$ est la matrice de $\mu_g$
 sur la base des \momsz.
 Par suite
 $$
 \prod\nolimits_{k=1}^p\gamma_k\,=\,\det D\,=\,\det G\,\det B\,=\,\Res(f,g)\,\det B.
 $$
 Ainsi l'\elt $\prod\nolimits_{k=1}^p\gamma_k$
 de  $\fa^p$ appartient à $\gen{\Res(f,g)}_\gA$.
\end{proof}
%-----------------end proof------------------

Le lemme d'\eli de base sera \gne plus loin (lemme \ref{lemElimPlusieurs}
et lemme d'\eli \gnl \ref{LemElimAff}).

L'appellation \gui{\id d'\eliz} correspond aux faits suivants
qui résultent du lemme précédent et du
lemme~\ref{lemResultant}:
%:     Corollary{corLemElimAffBasic}
\begin{corollary}\label{corLemElimAffBasic}
Soit $\gA$ un anneau intègre et $f$, $g\in\AX$. Si $f$ est \mon et se factorise 
complètement, \propeq
\begin{enumerate}
\item L'\id d'\eli $\gen{f,g}_{\AX}\cap\gA$ est nul.
\item Le résultant $\Res_X(f,g)=0$.
\item Les \pols $f$ et $g$ ont une racine commune.
\end{enumerate}

\end{corollary}
%--------- fin corollary ---------------------------------------------- 

 Un \cdi $\gK$ est dit \emph{\agqt clos} si tout \polu de~$\KX$ se décompose
 en produit de facteurs $X-x_i$ ($x_i\in\gK$).%
\index{corps!algebriquement@\agqt clos}%
\index{algebriquement clos@\agqt clos!corps discret ---}

%:     Corollary{cor2LemElimAffBasic}
\begin{corollary}\label{cor2LemElimAffBasic}
Soit $\gK$ un \cdacz.
\\  
On pose $\gA=\KYm$. Soient $f$ et $g\in\AX$ avec $f$ \mon en~$X$. 
Pour un \elt arbitraire $\uzeta=(\zeta_1,\ldots,\zeta_m)$ de $\gK^m$, \propeq
\begin{enumerate}
\item $\uzeta$ annule tous les \pols de l'\id d'\eli $\gen{f,g}\cap{\gA}$.%
\index{elimination@\eliz!ideal@\id d'---}\index{ideal@idéal!d'elimi@d'\eliz}
\item $\Res_X\big(f(\uzeta,X),g(\uzeta,X)\big)=0$.
\item $f(\uzeta,X)$ et $g(\uzeta,X)$ ont une racine commune.
\end{enumerate}
En conséquence si $V$ est l'ensemble des zéros communs à $f$ et $g$ dans 
$\gK^{m+1}$, et si $\pi:\gK^{m+1}\to\gK^{m}$ est la projection qui oublie la dernière coordonnée, alors $\pi(V)$ est l'ensemble des zéros de $\Res_X(f ,g)\in\KYm.$
\end{corollary}
%--------- fin corollary ---------------------------------------------- 

%:subsec   Retour sur le \discri  
\subsec{Retour sur le \discri}

Lorsque $g=\prod_{i=1}^n(X-y_i)$,
le lemme \ref{lemResultant} nous donne $\Res_X(g,g')=\prod_{i=1}^ng'(y_i)$
et donc
\begin{equation}
\label{eqDiscriRes}
\disc(g)=(-1)^{n(n-1)/2}\Res_X(g,g').
\end{equation}
Comme l'\egt $g(X)=\prod_{i=1}^n(X-y_i)$ peut toujours être
réalisée dans l'\adu si $g$ est \monz, on obtient que
l'\egrf{eqDiscriRes} est valable pour tout \poluz, sur
tout anneau commutatif.

Le fait suivant résulte donc du lemme d'\eli
 de base.

%:    Fact{factDiscUnit}-----------
\begin{fact}
\label{factDiscUnit}\index{discriminant!quand le --- est inversible}
On considère un \polu $g\in\AX$.
\vspace{-1pt}
\begin{enumerate}
\item [--] On a $\gen{g(X),g'(X)}=\gen{1}$  \ssi  $\disc g$ est \ivz.
\item [--]  L'\id $\gen{g(X),g'(X)}\cap\gA$ est fidèle \ssi $\disc g$ est un
\elt \ndz de $\gA$.
\end{enumerate}
\end{fact}
%--- end-fact-----------------------------------------

%:    Fact{factDiscProd}
\begin{fact} \label{factDiscProd}\index{discriminant!d'un produit}
Si $f=gh\in\AX$ avec $g$, $h$ \monsz, on a l'\egt suivante:
\vspace{-2pt}
\begin{equation}
\label{eqfactDiscProd}
\disc(f)=\disc(g)\disc(h)\Res(g,h)^2
\end{equation}
\end{fact}
\begin{proof}
Cela résulte \imdt des \eqns  (\ref{eqResultant4}), (\ref{eqResultant5})
\paref{eqResultant5} et (\ref{eqDiscriRes}).
\end{proof}

%:     Corollary{corfactDiscProd}
\begin{corollary}\label{corfactDiscProd} 
Soit $f\in\AX$  \mon et $\gB=\gA[x]=\aqo{\AX}{f}$.
\begin{enumerate}
\item Si  $f$ possède un facteur carré, $\disc f =0$. Inversement,
si  $\disc f =0$ et si $f(X)=\prod(X-x_i)$ dans un
anneau intègre  contenant $\gA$, deux des zéros $x_i$ sont égaux.
\item \label{icorfactDiscProd} Supposons $f$ \spl et $f=gh$ ($g$
et $h$ \monsz). 
\begin{enumerate}
\item Les \pols $g$ et $h$ sont \spls et \comz.
\item
Il existe un \idm $e$ de $\gB$ tel que $\gen{e}=\gen{\pi(g)}$.\\
On a $\gB\simeq\aqo{\gB}{g}\times \aqo{\gB}{h}$.
\end{enumerate}
\item Supposons $\disc f$ \ndz et $f=gh$ ($g$
et $h$ \monsz). \\
Alors, les \elts $\disc g $, $\disc h $ et
$\Res(g,h)$ sont \ndzsz.
\end{enumerate}
\end{corollary}
\begin{proof}
Tout ceci résulte du fait \ref{factDiscProd}, sauf peut-être
l'\idm $e$ dans le point~\emph{\ref{icorfactDiscProd}}. Si~$gu+hv=1$,
il faut prendre $e=\ov{gu}$.
\end{proof}
%

%:     Corollary{corcorfactDiscProd}
\begin{corollary}\label{corcorfactDiscProd}
Soient $\gK$ un \cdiz, $f\in\KX$ un \polu \spl et $\gB=\aqo{\KX}{f}$.
Dans le point \ref{icorfactDiscProd}. du corolaire précédent,
on associe à tout diviseur $g$ de $f$ l'\idm $e$ tel que
$\gen{\ov g}=\gen{e}$. Ceci établit une bijection entre les diviseurs
\mons  de $f$ et les \idms de~$\gB$. Cette bijection respecte la \dvez.
\end{corollary}
\begin{proof}
La bijection réciproque est donnée par $e=\ov v\mapsto\pgcd(v,f)$.
\end{proof}

Nous introduisons maintenant les notions de sous-corps premier et de
\cara d'un \cdiz. \rdb

Plus \gnltz, si $\gA$ est un anneau arbitraire, nous notons $\Z_\gA$ le \emph{sous-anneau premier de} $\gA$ défini comme suit:
$$
\Z_\gA=\sotq{n \cdot (m \cdot 1_A)^{-1}}{n,m\in \ZZ,\;m\cdot1_\gA\in\Ati}
.\label{NOTAZA}
$$
Si $\rho:\ZZ\to\gA$ est l'unique \homo d'anneaux de $\ZZ$ dans $\gA$,
le sous-anneau premier est donc isomorphe à $S^{-1}\ZZ\sur{\Ker\rho}$, où $S=\rho^{-1}(\Ati)$. Un anneau peut être appelé premier s'il est égal
à son sous-anneau premier.
En fait la terminologie n'est usuelle que dans le cas des corps.\rdb

Lorsque $\gK$ est un \cdiz, le sous-anneau premier
est un sous-corps, appelé \emph{sous-corps premier de~$\gK$}.
Pour un $m>0$ on dira que \emph{$\gK$ est de \cara $>m$},
et nous écrivons \gui{$\car(\gK)>m$}
si pour tout %entier 
$n\in\lrbm$, l'\elt $n\cdot1_\gK$ est \ivz.\label{NOTACarK}%
\index{premier!sous-anneau --- d'un anneau}%
\index{caractéristique!d'un corps}%
\index{premier!sous-corps --- d'un corps}%
\index{corps!premier}

Lorsque $\gK$ est non trivial, s'il existe un $m>0$ tel que $m\cdot1_\gK=0$, alors il en existe un minimum,
qui est un nombre premier $p$, et l'on dit que le corps \emph{est de \cara $p$}.
Lorsque le sous-corps premier de $\gK$ est isomorphe à $\QQ$, la tradition
est de parler de \emph{\cara nulle}, mais nous utiliserons aussi la terminologie  de \emph{\cara infinie} dans les contextes où cela est utile
pour rester cohérent avec la notation précédente, par exemple dans le fait~\ref{factPolSepFC}.

On peut concevoir\footnote{Il peut aussi s'en présenter à nous comme résultat d'une construction compliquée dans une \dem subtile.} des \cdis non triviaux dont la \cara n'est pas bien
définie du point de vue \cofz. Par contre pour un \cdi l'affirmation \gui{$\car(\gK)>m$} est toujours décidable.

%:     Fact{factPolSepFC}
\begin{fact}\label{factPolSepFC}
Soit $\gK$ un \cdi et $f\in\KX$ un \poluz.
Si~$\disc f=0$ et $\car(\gK)>\deg f$, $f$ possède un facteur carré
de degré $\geq1$.
\end{fact}
\begin{proof}
Soit $n=\deg f$. Le \polz~$f'$ est de degré $n-1$.
Soit $g=\pgcd(f,f')$, on a $\deg g \in\lrb{1..n-1}$ (fait \ref{factResNulPGCD}).
On écrit $f=gh$ donc

\snic{\disc(f)=\Res(g,h)^2\disc(g)\disc(h).}

%\sni
Ainsi, $\Res(g,h)=0$, ou $\disc(g)=0$, ou $\disc(h)=0$. Dans le premier
cas, les \pols $g$ et $h$ ont un pgcd $k$ de degré $\geq1$ 
et $k^2$ divise $f$.
Dans les deux autres cas, puisque $\deg g<\deg f$ et
$\deg h<\deg f$, on peut terminer par \recu sur le degré, en notant
que si $\deg f=1$, alors $\disc f\neq0$, ce qui assure l'initialisation.
\end{proof}
%

%:entrenous
\entrenous{%  subsec     %%%%%%%%%%%%
%\subsubsection*{Polynomes sous-résultants}
écrire qqch sur les sous-résultants?
}% finentrenous

%%%%%%%%%%%%%
\penalty-2500
\section[Théorie \agq des nombres, premiers pas]{Théorie \agq des nombres, \\ \hspace*{20pt}premiers pas}
\label{secApTDN}

Nous donnons ici quelques applications \gnlesz, en théorie des nombres \elrz,
des résultats précédemment obtenus dans ce chapitre. Pour entrevoir les multiples
facettes passionnantes de la théorie des nombres, \llec
pourra consulter le merveilleux ouvrage \cite {IR}.

%:   subsec{Algèbres finies, entières}
\subsec{Algèbres finies, entières}

Nous donnons quelques
précisions par rapport à la \dfnz~\ref{defEntierAnn0}.

%--- Definition{def05Alg}---------
\begin{definition}
\label{def05Alg}~
\begin{enumerate}
\item
Une \Alg $\gB$ est dite \ixd{finie}{algèbre}
si $\gB$ est un \Amo \tfz.  On dit aussi: \emph{$\gB$ est finie sur $\gA$}. 
Dans le cas
d'une extension, on parlera d'\emph{extension finie} de $\gA$.%
\index{algèbre!finie}
\item  Supposons  $\gA\subseteq\gB$. L'anneau
$\gA$ est dit \ix{intégralement clos} dans $\gB$ si tout \elt de
$\gB$ entier sur $\gA$  est dans~$\gA$.\index{anneau!intégralement clos dans
\ldots }
\end{enumerate}
%-----------------end enum------------------
\end{definition}
%--- end-definition------------------------------------

%--- Fact{factEntiersAnn}------------
\begin{fact}
\label{factEntiersAnn}
Soit  $\gA\subseteq\gB$  et $x\in\gB.$ \Propeq
%-----------------begin item------------------
%-----------------begin enum------------------
\begin{enumerate}
\item L'\elt $x$ est entier sur $\gA.$
\item La sous-\alg
$\gA[x]$ de $\gB$ est finie.
\item Il existe un \Amo fidèle et \tf $M\subseteq\gB$ tel que
$xM\subseteq M$.
\end{enumerate}
%-----------------end enum------------------
\end{fact}
%--- end-fact-------------------------------
%-----------------begin proof------------------
\begin{proof}
\emph{3} $\Rightarrow$~\emph{1} (a fortiori \emph{2} $\Rightarrow$~\emph{1.})
On considère une
matrice $A$ à \coes dans~$\gA$ qui repré\-sente $\mu_{x,M}$ (la multiplication par $x$ dans $M$) sur
un \sgr fini de $M$. Si $f$ est le \polcar de $A$, on a par le \tho de
Cayley-Hamilton  $0=f(\mu_{x,M})=\mu_{f(x),M}$
et puisque le module est fidèle,
$f(x)=0$. \\
Nous laissons le reste \alecz.
\end{proof}
%-----------------end proof------------------

On a aussi facilement le fait suivant.
%--- Fact{factEntiersAnn2}-----------
\begin{fact}
\label{factEntiersAnn2}
Soit  $\gB$  une \Alg et   $\gC$ une \Blgz.
%-----------------begin enum------------------
\begin{enumerate}
\item Si $\gC$ est finie sur $\gB$ et $\gB$ finie sur
$\gA$, alors $\gC$ est finie sur $\gA$.
\item Une \Alg engendrée par un nombre fini d'\elts entiers sur $\gA$ est
finie.
\item Les \elts de $\gB$  entiers sur $\gA$ forment un anneau
\icl dans~$\gB$.
On l'appelle la \emph{clôture (ou fermeture) intégrale de $\gA$ dans $\gB$}.%
\index{cloture@clôture intégrale!de $\gA$ dans $\gB\supseteq\gA$}
\end{enumerate}
%-----------------end enum------------------
\end{fact}
%--- end-fact-------------------------------

%:     Lemma{lemPolEnt}
\begin{lemma}\label{lemPolEnt}
Soient $\gA \subseteq \gB$ et $f \in \BuX$. Le \pol
$f$ est entier sur~$\AuX$ \ssi chaque \coe de $f$ est entier sur $\gA$.
\end{lemma}
\begin{proof} La condition est suffisante, d'après le point \emph{3}
du lemme précédent. Dans l'autre sens
on considère une \rdi  $P(f)=0$ pour $f$
(avec $P\in\AuX[T]$, \monz). On a dans $\gB[\uX,T]$ une \egt
$$
P(\uX,T)=\,\bigr(T - f(\uX) \bigl)\,
\bigr(T^n + u_{n-1}(\uX) T^{n-1} + \cdots + u_{0}(\uX) \bigl)
.$$
Puisque le \coe de $T^n$ dans le deuxième facteur est $1$, le \tho de
\KRO en plusieurs variables \ref{corthKro} implique que chaque \coe de $f$ est entier sur~$\gA$.
\end{proof}
%

%:     Lemma{lemPolCarInt}
\begin{lemma}\label{lemPolCarInt}
Soit $\gA \subseteq \gB$, $L$ un \Bmo libre de rang fini et
$u \in \End_\gB(L)$ entier sur~$\gA$.  Alors,  les \coes du
\polcar de~$u$ sont entiers sur~$\gA$. En
particulier, $\det( u)$ et $\Tr( u)$ sont entiers sur~$\gA$.
\end{lemma}
\begin{proof}
Démontrons d'abord que $\det(u)$ est entier sur~$\gA$.
Soit $\cE=(e_1,\ldots,e_n)$  une base fixée de $L$.
Le \Amo $\gA[u]$ est un \Amo
\tfz, et donc le module
$$
E = \som_{i\in\lrbn, k\geq 0} \gA u^k(e_i) \subseteq L
$$
est un \Amo  \tfz, avec $u(E) \subseteq E$. Introduisons le \Amo
$$
D = \som_{\ux \in E^n}
\gA\det_{\cE} (\ux) \subseteq \gB.
$$
Puisque $E$ est un \Amo \tfz, $D$ est un \Amo \tfz,
et il est fidèle: $1 \in D$ car $\det_{\cE} (\cE)=1$.
Enfin, l'\egtz
$$
\det(u) \det_{\cE} (x_1, \ldots, x_n) =
\det_{\cE} \big(u(x_1), \ldots, u(x_n)\big)
$$
et le fait que $u(E) \subseteq E$ montrent que $\det(u) D \subseteq D$.\\
Considérons ensuite $\AX\subseteq\BX$
et le $\BX$-module $L[X]$.\\ On a $X\Id_{L[X]} - u \in \End_{\gB[X]}(L[X])$.
Si $u$ est
entier sur~$\gA$, $X\Id_{L[X]} - u$ est entier sur~$\gA[X]$ donc $\rC{u}(X) =
\det(X\Id_{L[X]} - u)$ est entier sur $\gA[X]$. On conclut avec le lemme~\ref{lemPolEnt}.
\end{proof}

%:     corollary{corlemPolcarEntier}
\begin{corollary}\label{corlemPolcarEntier}
Soit $\gA\subseteq\gB\subseteq\gC$ avec $\gC$ une \Blg qui est un \Bmo
libre de rang fini.
Soit $x\in\gC$  entier sur $\gA$. Alors,   $\Tr_{\gC/\gB}(x)$,
 $\rN_{\gC/\gB}(x)$ et tous les \coes de  $\rC{\gC/\gB}(x)$
sont entiers sur $\gA$.
Si en plus $\gB$ est un \cdiz, les \coes du \polmin  $\Mip_{\gB,x}$
sont entiers sur~$\gA$.
\end{corollary}
\begin{proof}
On applique le lemme précédent avec $L=\gC$ et $u=\mu_x$.
Pour la dernière affirmation, on utilise le \tho de \KRO
et le fait que le \polmin divise le \polcarz.
\end{proof}
%

%
%  subsec   Anneaux intégralement clos
\penalty-2500
\subsubsection*{Anneaux intégralement clos}

%:     Definition{defiIntClos}
\begin{definition}\label{defiIntClos}
Un anneau intègre $\gA$ est dit \ix{intégralement clos} s'il
est \icl dans son corps de fractions.%
\index{anneau!intégralement clos}
\end{definition}

%:     Fact{fact.loc.entier}---------
\begin{fact}
\label{fact.loc.entier}\relax
Soit $\gA\subseteq\gB$,
$S$ un \mo de $\gA$, $x\in\gB$ et $s\in S$.
%-----------------begin item------------------
\begin{enumerate}
\item   L'\elt $x/s\in\gB_S$ est entier sur $\gA_S$ \ssi il existe
$u\in S$ tel que
$xu$ est entier sur  $\gA$.
\item  Si $\gC$ est la \cli de $\gA$ dans $\gB$,
alors $\gC_{S}$ est la \cli de $\gA_S$
dans~$\gB_{S}$.
\item  Si $\gA$ est \iclz, alors $\gA_S$ \egmtz.
\end{enumerate}
%-----------------end item------------------
\end{fact}
%--- end-fact-----------------------------------------
%
\begin{proof}
Il suffit de montrer le point \emph{1.}
 Supposons d'abord $x/s$ entier sur $\gA_S$.
On a par exemple une \egt dans $\gB$

\snic{u(x^3+a_2sx^2+a_1s^2x+a_0s^3)=0,}

%\sni
avec $u\in S$ et les $a_i\in\gA$. En multipliant par $u^2$ on obtient

\snic{(ux)^3+a_2us(ux)^2+a_1u^2s^2(ux)+a_0u^3s^3=0}

%\sni
 dans $\gB$. Inversement supposons $xu$ entier sur $\gA$ avec $u\in S$.
On a par exemple une \egt

\snic{(ux)^3+a_2(ux)^2+a_1(ux)+a_0=0}

%\sni
dans $\gB$, donc dans $\gB_S$:

\snic{x^3+(a_2/u)x^2+(a_1/u^2)x+(a_0/u^3)=0.}

%\sni
\vspace{-5pt}
\end{proof}
%

%:    Principe local global concret{plcc.entier}---
\begin{plcc}
\label{plcc.entier}\relax \emph{(\'Eléments entiers)}\\
Soient $S_1$, $\ldots$, $S_n$  des \moco d'un anneau $\gA\subseteq\gB$ et $x\in\gB$.
On a les \eqvcs suivantes.
%-----------------begin enum------------------
\begin{enumerate}
\item  L'\elt $x$ est entier sur $\gA$ \ssi il est entier sur chacun des~$\gA_{S_i}$.
\item  Supposons $\gA$ intègre: $\gA$ est \icl \ssi chacun des~$\gA_{S_i}$ est \iclz.
\end{enumerate}
%-----------------end enum------------------
\end{plcc}
%--- end-plcc-----------------------------------------
%-----------------begin proof------------------
\begin{proof}
Il faut montrer dans le point \emph{1} que si la condition est réalisée
localement, elle l'est globalement. On considère donc un $x\in\gB$
qui vérifie pour chaque $i$ une relation  $(s_ix)^k=a_{i,1}(s_ix)^{k-1}+ a_{i,2} (s_ix)^{k-2}+ \cdots +a_{i,k}$ avec
les~$a_{i,h}\in\gA$ et les $s_i\in S_i$ (on peut supposer \spdg que les degrés sont les mêmes).
On utilise alors une relation $\sum s_i^ku_i=1$
pour obtenir une \rdi de $x$ sur $\gA$.
\end{proof}
%-----------------end proof------------------

Le \tho de \KRO implique facilement le lemme qui suit.

%:     Lemma{lem0IntClos}
\begin{lemma}\label{lem0IntClos} 
\emph{(\Tho de Kronecker, cas d'un anneau intègre)}\\
Soit $\gA$ intégralement clos, de corps de fractions $\gK$.
Si l'on~a~$f=gh$ dans~$\KT$ avec $g$, $h$ \mons et $f\in\AT$, alors $g$ et $h$
sont aussi dans~$\AT$.
\end{lemma}

%:     Lemma{lemZintClos}
\begin{lemma}\label{lemZintClos}
L'anneau $\ZZ$ ainsi que l'anneau $\KX$ lorsque $\gK$ est un \cdiz,
sont intégralement clos.
\end{lemma}
\begin{proof}
En fait cela fonctionne avec tout anneau à pgcd intègre $\gA$ 
(voir la section~\ref{secGpReticules}).
Soient $f(T)=T^n-\sum_{k=0}^{n-1}f_kT^k$ et $a/b$ une fraction réduite dans le corps de fractions de $\gA$ avec $f(a/b)=0$.
En multipliant par $b^n$ on obtient  

\snic{a^n=b\;\sum_{k=0}^{n-1}f_ka^kb^{n-1-k}.}

%\sni
Puisque $\pgcd(a,b)=1$,  $\pgcd(a^n,b)=1$. Mais $b$ divise $a^n$, donc $b$ est \ivz, et
 $a/b\in\gA$.
\end{proof}
%

%:     Theorem{thIntClosStab}
\begin{theorem}\label{thIntClosStab}
Si $\gA$ est  intégralement clos, il en est de même pour~$\AX$.
\end{theorem}
\begin{proof}
Posons $\gK=\Frac\gA$. Si  un \elt $f$ de $\gK(X)$ est entier sur $\AX$, il est entier sur $\KX$,
donc dans $\KX$ car $\KX$ est intégralement clos.
On conclut avec le lemme~\ref{lemPolEnt}: tous les \coes du \polz~$f$ sont entiers sur $\gA$,
donc dans $\gA$.
\end{proof}

Un corolaire intéressant du \tho de \KRO est  la
\prt suivante (avec les mêmes notations que dans le 
\thref{thKro}).

%:    Proposition{propArm}--------------
\begin{proposition}
\label{propArm} Soient $f,g\in\AuX$.
Supposons que $\gA$  est intégralement clos, et que $a\in \gA$
divise tous les \coes de $h=fg$, alors $a$ divise tous les $f_{\!\alpha\,} g_\beta$.
Autrement dit
%-----------------begin $$----------------
$$ \rc(fg)\equiv 0\;\mod\;a\;\;\iff \;\;
\rc(f)\rc(g)\equiv 0\;\mod\;a.
$$
%-----------------end $$------------------
\end{proposition}
%--- end-proposition----------------------------------------
%-----------------begin proof------------------
\begin{proof}
En effet, en considérant les \pols $f/a$ et $g$ à \coes dans le corps
des fractions de $\gA$,  le \tho de \KRO implique que $f_{\!\alpha\,} g_\beta/a$
est entier sur $\gA$ car les $h_\gamma/a$ sont dans $\gA$.
\end{proof}
%-----------------end proof------------------

%  subsec   Décomposition de \pols
\penalty-2500
\subsubsection*{Décomposition de \pols en produits de
facteurs\\ \irdsz}

%:     Lemma{lemKXfactor}
\begin{lemma}\label{lemKXfactor}
Soit $\gK$ un \cdiz. Les \pols de $\KX$
se décom\-posent en produits de facteurs
\irds \ssi on a un \algo pour le calcul des zéros
dans $\gK$ d'un \pol arbitraire de~$\KX$.
\end{lemma}
\begin{proof}
La deuxième condition est a priori plus faible puisqu'elle revient à
déterminer les facteurs de degré 1 pour un \pol  de $\KX$. Supposons cette condition
vérifiée. Pour savoir s'il existe une \dcn $f=gh$
avec~$g$ et~$h$ \mons de degrés $>0$ fixés, on applique le \tho de \KROz. On obtient pour chaque \coe de $g$ et $h$ un nombre fini
de possibilités (ce sont des zéros de \polus que l'on peut
expliciter en fonction des \coes de $f$).
\end{proof}
%

%:     Proposition{propZXfactor}
\begin{proposition}\label{propZXfactor}
Dans $\ZZ[X]$ et $\QQ[X]$ les \pols se décomposent en produits de facteurs
\irdsz. Un \pol non constant de $\ZZ[X]$ est irréductible dans $\ZZ[X]$ \ssi il est
primitif et irréductible dans~$\QQ[X]$.
\end{proposition}
\begin{proof}
Pour $\QQ[X]$ on applique le lemme \ref{lemKXfactor}. Il faut donc montrer que l'on sait déterminer les zéros rationnels d'un \polu $f$ à \coes rationnels. On peut même supposer les \coes de $f$ entiers.
La théorie \elr
de la \dve dans $\ZZ$ montre alors que si $a/b$ est un zéro de~$f$,
$a$ doit diviser le \coe dominant et $b$ le \coe constant de $f$: il n'y a donc qu'un nombre fini de tests à faire.\\
Pour $\ZZ[X]$, un \pol primitif $f$ étant donné, on cherche à savoir s'il existe une \dcn   $f=gh$ avec $g$ et $h$  de degrés $>0$ fixés. On peut supposer $f(0)\neq0$.
On applique le \tho de \KROz. Un produit $g_0h_j$ par exemple
doit être un zéro dans $\ZZ$ d'un \polu $q_{0,j}$ de $\ZZ[T]$ que l'on peut calculer.
En particulier, $g_0h_j$ doit diviser  $q_{0,j}(0)$, ce qui ne laisse qu'un nombre fini de possibilités pour $h_j$.\\
Enfin pour le dernier point si un \polz~$f$ primitif dans $\ZZ[X]$
se décompose sous la forme $f=gh$ dans $\QQ[X]$ on peut supposer que $g$ est primitif dans~$\ZZ[X]$; soit alors $a$ un \coe de $h$, tous les $ag_j$ sont dans $\ZZ$ (\tho de \KROz),
et une relation de Bézout $\sum_jg_ju_j=1$ montre \hbox{que $a\in\ZZ$}.
\end{proof}
%

%:  subsec  Corps de nombres
%\penalty-2500
\vspace{5pt}
\subsec{Corps de nombres}

On appelle \emph{corps de nombres} un \cdi  $\gK$
\stf sur $\QQ$.

%  subsubsec cloture galoisienne
\vspace{-1pt}
\subsubsection*{Clôture galoisienne}

%:     Theorem{thClotAlgQ}
\begin{theorem}\label{thClotAlgQ} \emph{(Corps de racines, \tho de l'\elt primitif)}
\begin{enumerate}
\item
Si $f$ est un \polu \spl de $\QQ[X]$
il existe un corps de nombres $\gL$ sur lequel on peut écrire
$f(X)=\prod_i(X-x_i)$.
En outre, avec un $\alpha\in\gL$ on a:
$$\gL=\QQ[\xn]=\QQ[\alpha]\simeq\aqo{\QQ[T]}{Q},$$
où $Q(\alpha)=0$ et  le \polu $Q$
est \ird dans $\QQ[T]$ et se décom\-pose complètement dans $\gL[T]$.
\\
En particulier, l'extension $\gL/\QQ$ est galoisienne et le \thrf{thGaloiselr} s'applique.
\item
Tout corps de nombres $\gK$ est contenu dans une extension galoisienne
du type précédent. En outre, il existe un $x\in\gK$ tel que $\gK=\QQ[x]$.
\end{enumerate}

\end{theorem}
\begin{proof}
\emph{1.} Cela résulte du \thrf{thResolUniv} et de la proposition \ref{propZXfactor}.\\
\emph{2.} Un corps de nombres est engendré par un nombre fini d'\elts
qui sont  \agqs sur $\QQ$.
Chacun de ces \elts admet un \polmin qui est \ird sur
$\QQ$ donc \spl (fait \ref{factPolSepFC}).
En prenant le ppcm~$f$ de ces \pols on obtient
un \pol \splz. En appliquant le point~\emph{1} à $f$
et en utilisant le \thrf{propUnicCDR}, on voit que $\gK$ est isomorphe
à un sous-corps de $\gL$. Enfin comme la correspondance galoisienne est
bijective et comme le groupe de Galois $\Gal(\gL/\QQ)$ est fini,
le corps $\gK$ ne contient qu'un nombre fini, explicite, de sous-corps
$\gK_i$ \stfs sur $\QQ$. Si l'on choisit $x\in\gK$ en dehors de la réunion
de ces sous-corps (qui sont des sous-\Qevs stricts), on a \ncrt
$\QQ[x]=\gK$: c'est un sous-corps de $\gK$ \stf sur $\QQ$ et distinct de
tous les $\gK_i$.
\end{proof}

%%%%%%%%%%%%%%%%%%%%%%%%%%%%%%%%%%%%%%%%%%%%%%%%%%%%%%%%%%%%%%%%%%%%%%%%%%%
%  subsubsec \'Elément cotransposé
\subsubsection*{\'Elément cotransposé} \rdb

Si $\gB$ est une \Alg libre de rang fini,
on peut identifier $\gB$ à une sous-\alg
commutative de $\End_\gA(B)$, où $B$ désigne le \Amo $\gB$ privé de sa structure multiplicative, au moyen de
l'\homo $x\mapsto\mu_{\gB,x}$, \hbox{où $\mu_{\gB,x}=\mu_x$} est la multiplication par $x$ dans $\gB$.
Alors, puisque $\wi \mu_x=G(\mu_x)$ pour un \polz~$G$ de $\AT$
(lemme \ref{lemPrincipeIdentitesAlgebriques} point~\emph{6}), on peut définir
$\wi x$ par l'\egt $\wi x=G(x)$, ou ce qui revient au même
$\wi {\mu_x}=\mu_{\wi x}$.
Si plus de précision est \ncrz,
on utilisera la notation $\Adj\iBA (x)$.
Cet \elt $\wi x$ s'appelle \emph{l'\elt cotransposé de x}.
On a alors l'\egt importante:%
\index{cotransposé!element@\elt --- (dans une \alg libre)}
\begin{equation}
\label{eqelt0cotransp}
x\ \wi x=x\ \Adj\iBA(x)=\rN\iBA(x).
\end{equation}

\rdb
\rem
\label{factNormeRationnelle}
Notons aussi que les applications \gui{norme de} et \gui{\elt cotransposé de} jouissent de propriétés de \gui{$\gA$-rationalité}, qui résultent directement de leurs \dfnsz:
si $P\in\gB[\Xk]$,
alors en prenant les $x_i$ dans $\gA$,
 $\rN\iBA \big(P(\xk)\big)$ et $\Adj\iBA \big(P(\xk)\big)$ sont donnés par des \pols de $\gA[\Xk]$.
\\
En fait $\BuX$ est libre sur $\AuX$ avec la même base que celle de $\gB$
sur~$\gA$ et  $\rN\iBA \big(P(\ux)\big)$ est donné par l'évaluation
en $\ux$ de $\rN_{\BuX/\!\AuX}\big(P(\uX)\big)$ (même chose pour l'\elt cotransposé). On utilisera par abus la
notation $\rN\iBA \big(P(\uX)\big)$.
\\
   En outre, si $\dex{\gB:\gA}=n$ et si $P$ est
\hmg de degré $d$, alors
 $\rN\iBA \big(P(\uX)\big)$ est \hmg de degré
$nd$  et
  $\Adj\iBA \!\big(P(\uX)\big)$ est \hmg de degré
$(n-1)\,d$.
\eoe

%:  subsec  Anneau d'entiers
\vspace{-7pt}
\subsec{Anneau d'entiers d'un corps de nombres}

Si $\gK$ est  un corps de nombres son \ixc{anneau d'entiers}{d'un corps de nombres}
est la \cli de  $\ZZ$ dans $\gK$.

%:     Propdef{propAECDN}
\begin{propdef}\label{propAECDN} \emph{(Discriminant d'un corps de nombres)}
 \\
Soit $\gK$  un corps de nombres et $\gZ$
son anneau d'entiers.
\begin{enumerate}
\item  Un \elt $y$ de $\gK$ est dans $\gZ$ \ssi $\Mip_{\QQ,y}(X)\in\ZZX$. 
\item  On a $\gK=(\NN\etl)^{-1}\gZ$.  
\item  
Supposons que $\gK=\QQ[x]$ avec $x\in\gZ$. Soit $f(X)=\Mip_{\QQ,x}(X)$ dans
$\ZZ[X]$ et $\Delta^2$ le plus grand facteur carré de $\disc_X f$. Alors,
$\ZZ[x]\subseteq\gZ\subseteq  {1 \over \Delta} \ZZ[x]$.

\item L'anneau $\gZ$ est un \ZZmo libre de rang $\dex{\gK:\gQ}$.
\item  L'entier $\Disc_{\gZ/\ZZ}$
est bien défini, on l'appelle le \emph{discriminant du corps de nombres $\gK$}.%
\index{discriminant!d'un corps de nombres} 
\end{enumerate}
 
\end{propdef}
%--------- fin proposition ---------------------------------------------- 
%
\begin{proof}
\emph{1.} Résulte du lemme \ref{lem0IntClos} (\tho de \KROz). 
 
\emph{2.} Soit $y\in\gK$ et $g(X)\in\ZZX$ un \pol non nul qui annule $y$. 
Si $a$ est le \coe dominant de $g$, $ay$ est entier sur $\ZZ$.
 
\emph{3.}  
Posons $\gA=\ZZ[x]$ et $n=\dex{\gK:\gQ}$.  
Soit $z\in\gZ$, que l'on écrit $h(x)/\delta$ avec~$\delta\in\NN\etl$,
$\gen{\delta}+\rc(h)=\gen{1}$ et $\deg h < n$.  On a $\gA+\ZZ z\subseteq
{1 \over \delta} \gA$ et il suffit donc de montrer que $\delta^2$ divise~$\disc_X(f)$.
%
%\\
%
 L'anneau $\gA$ est un \ZZmo libre de rang $n$, avec la base
$\cB_0=(1,x,\ldots,x^{n-1})$.  La proposition \ref{propdiscTra} donne:

\snic{\Disc_{\gA/\ZZ}=\disc_{\gA/\ZZ}(\cB_0)=\disc_{\gK/\QQ}(\cB_0)=\disc_X f .}

%\sni

Le \ZZmo $M =\gA+\ZZ z$ est \egmt libre de rang $n$, avec une base $\cB_1$,
et l'on obtient les \egts
 
\snic{\disc_X f=\disc_{\gK/\QQ}(\cB_0)= \disc_{\gK/\QQ}(\cB_1)\times d^2,}

%\sni
 où $d$ est le \deter de la matrice de $\cB_0$ sur $\cB_1$ (proposition~\ref{defiDiscTra}~\emph{2}).
\\
Enfin $d=\pm\delta$ d'après le lem\-me \ref{lemSousLibre}
qui suit.  Et l'on peut conclure.
 
\emph{4.} On se place \spdg dans la situation du point \emph{3.}
Il n'y a qu'un nombre fini de \ZZmos \tf entre $\ZZ[x]$ et ${1\over\Delta}\ZZ[x]$. 
Et pour chacun d'entre eux on peut tester s'il est contenu dans $\gZ$.
Le plus grand possible est \ncrt égal à $\gZ$.
\end{proof}
\rems ~\\
1) Comme corolaire, on voit que dans la situation du point~\emph{3}, si
$\disc_X(f)$ est sans facteur carré, alors $\gZ=\ZZ[x]$.

2) La \dem du point \emph {4} ne donne pas de moyen {pratique} pour calculer une
$\ZZ$-base de $\gZ$. Pour quelques informations plus précises voir le
\pbz~\ref{exoLemmeFourchette} (lemme de la fourchette). 
En fait on ne connaît pas d'\algo \gnl
\gui{en temps \pollz} pour calculer une $\ZZ$-base de $\gZ$.
\eoe

%:     Lemma{lemSousLibre}
\begin{lemma}\label{lemSousLibre}
Soient $N \subseteq M$ deux \Amos libres de même \hbox{rang $n$}
\hbox{avec $M = N + \gA z$}. On suppose que pour un \elt \ndz $\delta \in \gA$, on 
\hbox{a~$\delta z \in N$} et $\delta z = a_1 e_1 + \cdots + a_n e_n$,
où $(e_1, \ldots, e_n)$ est une base de $N$. Alors, le \deter $d$
d'une matrice d'une base de $N$ sur une base de $M$ vérifie:
\begin{equation}
\label{eqlemSousLibre}
d \,\gen {\delta, \an} = \gen {\delta}
\end{equation}
En particulier, $\gen {\delta, \an}$ est un \idpz, et si $\delta$, $a_1$, \dots, $a_n$ sont \comz, alors $\gen {d} = \gen {\delta}$. En outre, $M/N \simeq 
\aqo {\gA}{d}$.  
\end{lemma}
%--------- fin lemma ---------------------------------------------- 
%
\begin{proof} L'\egrf{eqlemSousLibre}  est laissée \alec (voir 
%:2012 reference inutile supprimée
%le fait \ref{fact.idd.sousmod} et 
l'exercice \ref{exolemSousLibre}).
\\
Il nous reste à montrer $M/N \simeq \aqo {\gA}{d}$.
En notant $\ov z$ la classe de $z$ dans~$M/N$, puisque $M/N\simeq\gA \ov z$, on doit montrer que $\Ann_\gA(\ov z) = \gen {d}$,
\cad \hbox{que $bz \in N \Leftrightarrow b \in \gen {d}$}. Il est clair que $dz \in N$.\\
Si $bz \in N$, alors $b\delta z \in \delta N$, donc en écrivant
$\delta z = a_1e_1 + \cdots + a_ne_n$, il \hbox{vient $ba_i \in \gen{\delta}$}, puis
$b\gen {\delta,\an} \subseteq \gen {\delta}$. En 
multipliant par $d$ et en simplifiant par~$\delta$, on obtient $b \in \gen {d}$.
\end{proof}
%

%%%%%%%%%%%%%%%%%%%%%%%%%%%%%%%%%%%%%%%%%%%%%%%%%%%%%%%%%%%%%%%%%%%%%%%%%%%
%  subsec  Théorie multiplicative des idéaux d'un corps de nombres
\subsubsection*{Théorie multiplicative des idéaux d'un corps de nombres}

%:    Definition{defiiv}--------------
\begin{definition}
\label{defiiv}
Un \id $\fa$ d'un anneau $\gA$ est dit \emph{\ivz}
s'il existe un \id $\fb$ et un
\elt \ndz $a$ tels que $\fa\,\fb=\gen{a}$.%
\index{inversible!idéal}%
\index{ideal@idéal!inversible}
\end{definition}
%--- end-definition------------------------------------

%:    Fact{factdefiiv}------------
\begin{fact}
\label{factdefiiv}
Soit $\fa$ un \id \iv d'un anneau $\gA$.
\begin{enumerate}
\item  L'\id $\fa$ est \tfz.
\item Si $\fa$ est engendré par $k$ \elts et si $\fa\,\fb=\gen{a}$
avec $a$ \ndzz, alors  $\fb$ est engendré par $k$ \eltsz. En outre 
$\fb=(\gen{a}:\fa)$.
\item  On a la règle  $\fa\,\fc \subseteq \fa\,\fd \;\Rightarrow\;\fc \subseteq \fd $
pour tous \idsz~$\fc $ et~$\fd $.
\item Si $\fc \subseteq \fa$ il existe un unique $\fd $ tel que $\fd \,\fa=\fc $, à savoir $\fd=(\fc:\fa)$. 
\\ Et si $\fc $ est \tfz, il en va de même pour $\fd $. 
\end{enumerate}
\end{fact}
%--- end-fact-----------------------------------------
%-----------------begin proof------------------
\begin{proof}
 \emph{3.} Si $\fa\,\fc \subseteq \fa\,\fd $ en multipliant par $\fb$ on obtient $a\,\fc \subseteq a\,\fd $.
Et puisque $a$ est \ndzz, cela implique  $\fc \subseteq \fd $.
\\
\emph{1.}
Si  $\fa\,\fb=\gen{a}$,  on trouve deux \itfs $\fa_1\subseteq\fa$ et $\fb_1\subseteq\fb$
tels \hbox{que $a\in\fa_1\,\fb_1$} et donc $\fa\,\fb=\gen{a}\subseteq\fa_1\,\fb_1\subseteq\fa\,\fb_1\subseteq \fa\,\fb$.
On en déduit les \egts $\fa_1\,\fb_1=\fa\,\fb_1=\fa\,\fb$.
D'où $\fb=\fb_1$ d'après le point \emph{3}. De même, $\fa=\fa_1$.
\\
\emph{2.} Si $\fa=\gen{a_1,\dots,a_k}$, on trouve $b_1$, \dots, $b_k\in\fb$
tels que $\som_i a_ib_i=a$. \\
En raisonnant comme au point \emph{1}
avec $\fa_1=\fa$ et $\fb_1=\gen{b_1,\dots,b_k}$ on obtient 
l'\egt $\fb=\gen{b_1,\dots,b_k}$. %\\
Puisque $\fa\,\fb=\gen{a}$, on a $\fb\subseteq (\gen{a}:\fa)$. Réciproquement, si $x\fa\subseteq \gen{a}$, alors 
%
%\snic
${x\gen{a}=x\,\fa\,\fb\subseteq a\,\fb},$
%
%\snii
donc $ax=ab$ pour un $b\in\fb$ et $x\in\fb$ car $a$ est \ndzz. 
\\
\emph{4.} De $\fa \,\fb=\gen{a}$ on déduit $\fc \,\fb\subseteq \gen{a}$.
Tous les \elts de $\fc \,\fb$ étant multiples de $a$, en les divisant
par $a$ on obtient un \id $\fd $, que l'on  note $\fraC 1 a \,\fc \,\fb$, et
avec lequel on obtient l'\egt $\fa\,\fd=\fraC 1 a \,\fc \,\fb\,\fa =\fraC 1 a \,\fc \,\gen{a}=\fc $ car $a$ est \ndzz.
\\
Si $\fc $ est \tfz, $\fd $ est engendré par les \elts obtenus en divisant chaque \gtr de $\fc \,\fb$ par $a$.
\\
L'unicité de $\fd$ résulte du point \emph{3.}
\\
Il reste à montrer que $\fd=(\fc:\fa)$. L'inclusion $\fd\subseteq (\fc:\fa)$
est \imdez.   
Réciproquement, si $x\fa\subseteq \fc$, alors $x\gen{a}\subseteq \fc\,\fb$, 
donc $x\in\fraC 1 a \,\fc \,\fb=\fd$.
\end{proof}
%-----------------end proof------------------

Le \tho suivant est le \tho clé dans la théorie multiplicative des \ids de corps de nombres. Nous en donnons deux \demsz.
Auparavant nous convions \llec à visiter le \pb \ref{exoPetitKummer} 
qui donne le petit \tho de Kummer, lequel résout à moindres frais la question pour \gui{presque tous} les \itfs
des corps de nombres.% 
\index{Kummer!petit \tho de ---}
Le \pbz~\ref{exoCyclotomicRing} est \egmt instructif car il donne une preuve 
directe de l'inversi\-bilité de tous les \itfs non nuls 
ainsi que de leur \dcn unique
en produit de \gui{facteurs premiers} pour l'anneau $\ZZ[\root n \of 1\,]$. 

%:     Theorem{th1IdZalpha}
\begin{theorem} \emph{(Inversibilité des \ids d'un corps de nombres)}\label{th1IdZalpha}\\
Tout \itf  non nul de l'anneau d'entiers $\gZ$ d'un corps de nombres~$\gK$ est \ivz.
\end{theorem}
\begin{proof}
\emph{Première \demz.} (à la Kronecker\footnote{En fait Kronecker n'utilise
pas l'\emph{\elt cotransposé} de $\alpha+\beta X+\gamma X^2$ (selon la \dfn que nous avons donnée), mais le produit de tous les conjugués  de $\alpha X+\beta Y+\gamma Z$ dans une extension galoisienne.
Ceci introduit une légère variation dans la \demz.})
\\
Prenons par exemple $\fa=\gen{\alpha,\beta,\gamma}$.
Notons $\gA=\QQ[X]$ et $\gB=\KX$.
L'\alg $\gB$ est libre sur $\gA$ avec la même base que celle de $\gK$
sur $\QQ$.
On considère le \polz~$g=\alpha+\beta X+\gamma X^2$
qui vérifie $\rc_\gZ(g)=\fa$.
Puisque $\alpha$, $\beta$, $\gamma$ sont entiers sur $\ZZ$, $g$ est entier
sur $\ZZ[X]$.
Soit $h(X)=\Adj\iBA (g)$ l'\elt cotransposé de $g$.
On sait que $h$ s'exprime comme un \pol en $g$ et en les \coes du
\polcar de $g$. En appliquant le corolaire~\ref{corlemPolcarEntier}
on en déduit que $h$ est à \coes dans $\gZ$.
 Notons $\fb$ l'\itf de $\gZ$
engendré par les \coes de $h$.
On a $gh=\rN\iBA (g)\in\gZ[X]\cap\QQ[X]=\ZZ[X]$.
Soit $d$ le pgcd des \coes de $gh$. La proposition \ref{propArm}
nous dit qu'un \elt arbitraire  de $\gZ$ divise $d$ \ssi il
divise tous les \elts de $\fa\,\fb$. En particulier, $d\gZ\supseteq\fa\,\fb$. Vu la relation de Bézout qui exprime $d$ en fonction des \coes de $gh$
on a aussi $d\in\fa\,\fb$. Donc $d\gZ=\fa\,\fb$.
\\
\emph{Deuxième \demz.} (à la Dedekind)\index{Dedekind!inversion d'un idéal à la ---}
\\
 Tout d'abord on remarque qu'il suffit de savoir inverser
les \ids à deux \gtrs en vertu de la remarque suivante. Pour trois \ids
arbitraires~$\fa$,~$\fb$,~$\fc$ dans un anneau on a toujours l'\egt

\snic{(\fa+\fb)(\fb+\fc)(\fc+\fa)=(\fa+\fb+\fc)(\fa\fb+\fb\fc+\fa\fc),}

%\sni
donc, si l'on sait inverser les \ids à $k$ \gtrs ($k\geq2$), on sait \egmt
inverser les \ids à $k+1$ \gtrsz.\\
On considère donc un \id $\gen{\alpha,\beta} $ avec $\alpha\neq0$.
Comme  $\alpha$ est
entier sur $\ZZ$, on peut trouver $\ov{\alpha} \in \gZ$ tel que
$\ov{\alpha}\alpha \in \ZZ \setminus \so{0}$.
Ainsi, quitte à remplacer $(\alpha,\beta )$ par $(\ov{\alpha}\alpha,\ov{\alpha}\beta )$,
on se restreint à l'étude d'un \id
$\gen{a,\beta}$ avec $(a,\beta)\in\ZZ\times\gZ$.
\\
 Soit $f \in \ZZ[X]$ un \polu s'annulant en $\beta$. On écrit

\snic{f(X)=(X-\beta)h(X)$,
où $h \in \gZ[X]\,.}

%\sni
On a donc $f(a X)=(a X-\beta )h(a X)$, que l'on réécrit
$f_1=g_1h_1$.
Soit alors $d$ le pgcd  des \coes de $f_1$ dans~$\ZZ$.
 Avec  $\fb=\rc_\gZ(h_1)$ et
$\fa=\rc_\gZ(g_1)=\gen{a,\beta}$,
on a clairement $d \in \fa \fb$.
Par ailleurs, la proposition \ref{propArm}
nous dit qu'un \elt arbitraire de~$\gZ$ divise tous les \elts de $\rc_\gZ(f_1)=\gen{d}$ \ssi il
divise tous les \elts de l'\id produit $\fa\,\fb$. En particulier,~$d$ 
divise tous les \elts de $\fa\,\fb$.
Ainsi $\fa \fb=\gen{d}$.
\end{proof}

Le \tho suivant montre que les \itfs d'un corps de nombres se comportent
vis à vis des opérations \elrs (somme, intersection, produit, division exacte)  de manière essentiellement équivalente aux \ids principaux de $\ZZ$, lesquels traduisent
de façon très précise la théorie de la \dve pour les
entiers naturels.
\\
Rappelons que dans 
la bijection $n \mapsto n\ZZ$ ($n\in\NN$, $n\ZZ$
\itf de~$\ZZ$), le produit correspond au produit, la \dve à l'inclusion, le pgcd à la somme, le ppcm à l'intersection et la division exacte au transporteur.

%:     Theorem{th2IdZalpha}
\begin{theorem} \emph{(Les \itfs d'un corps de nombres)}\label{th2IdZalpha}\\
Soit $\gK$  un corps de nombres et $\gZ$ son anneau d'entiers.
\begin{enumerate}
\item \label{i2th2IdZalpha}
Si $\fb$  et $\fc$ sont deux \ids arbitraires, et si  $\fa$ est un \itf
non nul de $\gZ$, on a l'implication:
$$
\fa\, \fb\subseteq \fa\, \fc\ \ \Rightarrow\ \ \fb\subseteq  \fc\,.
$$
\item \label{i3th2IdZalpha}
 Si $\fb\subseteq\fc$ sont deux \itfsz, il existe un \itf  $\fa$
tel que $\fa\, \fc=\fb$.
\item \label{i4th2IdZalpha}
 L'ensemble des \itfs de $\gZ$ est stable par intersections finies et l'on a
les \egts suivantes (où $\fa$, $\fb$, $\fc$ désignent des \itfs de~$\gZ$):
\vspace{-8pt}
\[\arraycolsep2pt\def\mathit#1{\llap{{#1}\hskip1.5cm}}
\begin{array}{lrclc}
\mathit{a.} &(\fa\cap\fb)(\fa+\fb)  &  = & \fa\fb \,,&\qquad\qquad \qquad \\[1mm]
\mathit{b.} &\fa\cap(\fb+\fc)  &  = & (\fa\cap\fb)+(\fa\cap\fc)  \,, \\[1mm]
\mathit{c.} &\fa+(\fb\cap\fc)  &  = & (\fa+\fb)\cap(\fa+\fc)   \,,\\[1mm]
\mathit{d.} &\fa  (\fb\cap\fc)  &  = & (\fa \fb)\cap(\fa \fc)  \,,\\[1mm]
\mathit{e.} &(\fa+\fb)^n  &  = & \fa^n+\fb^n\quad(n\in\NN) \,.
\end{array}
\]
\item \label{i5th2IdZalpha} Si  $\fa$ est un \itf
non nul de $\gZ$ l'anneau $\gZ\sur\fa$ est fini. 
\\
En particulier,
on a des tests pour décider:
\begin{itemize}
\item si un $x\in\gZ$ est dans $\fa$,
\item si un $x\in\gZ$ est \iv modulo $\fa$, 
\item si $\fa$ est contenu dans un autre \itf $\fb$, 
\item  si  $\gZ\sur\fa$ est un \cdi (on dit alors que $\fa$
est un \idema détachable).
\end{itemize}

\item \label{i6th2IdZalpha} Tout \itf distinct de $\gen{0}$ et $\gen{1}$ est égal à un produit d'\idemas
\ivs détachables, et cette \dcn est unique à l'ordre près des facteurs.
\end{enumerate}
\end{theorem}
\begin{proof}
\emph{\ref{i2th2IdZalpha}} et \emph{\ref{i3th2IdZalpha}.}
D'après le fait \ref{factdefiiv}.
 
\emph{\ref{i4th2IdZalpha}.}
Si l'un des \itfs est nul tout est clair. On les suppose dans la suite non nuls.
 
\emph{\ref{i4th2IdZalpha}a.}
Soit $\fc$ tel que $\fc(\fa+\fb)=\fa\fb$.
Puisque $(\fa\cap\fb)(\fa+\fb)\subseteq\fa\fb$, on obtient 
l'inclusion $\fa\cap\fb\subseteq\fc$ (simplification par $\fa+\fb$). Inversement, $\fc\fa \subseteq \fa\fb$, \hbox{donc $\fc \subseteq \fb$}
   (simplification par $\fa$). De même $\fc \subseteq \fa$.
 
\emph{\ref{i4th2IdZalpha}c.}
On multiplie les deux membres par $\fa+\fb+\fc=(\fa+\fb)+(\fa+\fc)$.
Le membre de droite donne $(\fa+\fb)(\fa+\fc)$.\\
Le membre de gauche donne $\fa(\fa+\fb+\fc)+ \fa(\fb\cap\fc)+(\fb+\fc)(\fb\cap\fc)$. \\
Dans les deux cas cela fait $\fa(\fa+\fb+\fc)+\fb\fc$.
 
\emph{\ref{i4th2IdZalpha}b.}
Les \itfs forment pour l'inclusion un treillis (le sup est la somme
et le inf est l'intersection). On vient de voir qu'une des lois est distributive par rapport à l'autre. Il est classique dans un treillis que cela implique l'autre \dit (voir \paref{DistriTrdi}).
 
\emph{\ref{i4th2IdZalpha}d.}
L'application $\fx\mapsto\fa\,\fx$ (de l'ensemble des \itfs
vers l'ensemble des \itfs multiples de $\fa$) est un \iso pour la structure d'ordre d'après le point \emph{\ref{i2th2IdZalpha}} Ceci implique qu'elle transforme
le inf en le inf.  Il suffit donc d'établir que $\fa\fb\cap\fa\fc$
est multiple de $\fa$. Cela résulte du point \emph{\ref{i3th2IdZalpha}.}
 
\emph{\ref{i4th2IdZalpha}e.}
Par exemple avec $n=3$, $(\fa+\fb)^3=\fa^3+\fa^2\fb+\fa\fb^2+\fb^3$.
\\
En multipliant  $(\fa+\fb)^3$ et $\fa^3+\fb^3$ par $(\fa+\fb)^2$
on trouve dans les deux \hbox{cas $\fa^5+\fa^4\fb+\cdots+\fa\fb^4+\fb^5$}.
 
\emph{\ref{i5th2IdZalpha}.} 
On regarde $\gZ$ comme un \ZZmo libre de rang $n=\dex{\gK:\QQ}$. On se convainc facilement qu'un \itf $\fa$ contenant l'entier $m\neq0$ peut être explicité comme un sous
\ZZmo \tf de $\ZZ^n$ contenant $m\ZZ^n$. 
 
\emph{\ref{i6th2IdZalpha}.}
Soit $\fa$ un \itf $\neq\gen{0},\gen{1}$. 
Les \idemas \tf de $\gZ$ contenant $\fa$ sont obtenus 
en déterminant les \idemas \tf de $\gZ\sur\fa$ (ce qui est possible parce que
l'anneau $\gZ\sur\fa$ est fini). 
Si $\fp$ est un \idema \tf contenant $\fa$, on peut écrire $\fa=\fb\,\fp$.
En outre, 
on a l'\egt $\idg{\gZ:\fa}=\idg{\gZ:\fb} \,\idg{\fb:\fa}$. On obtient alors la \dcn en produit d'\idemas \tf par \recu sur $\idg{\gZ:\fa}$. L'unicité résulte du fait que si un \idema \tf $\fp$ contient un produit d'\idemas \tfz, il est forcément égal à l'un d'entre eux, car sinon il serait comaximal avec le produit. 
\end{proof}

Nous terminons cette section par quelques généralités concernant \emph{les \ids qui évitent le conducteur}. 
La situation en théorie des nombres est la suivante.
On a un corps de nombres $\gK=\QQ[\alpha]$ avec $\alpha$ entier sur $\ZZ$.
On note $\gZ$ l'anneau des entiers de $\gK$, \cad la \cli de $\ZZ$
dans $\gK$. Bien que ce soit en principe possible, il n'est pas toujours facile d'obtenir une base de $\gZ$ comme \ZZmoz, ni d'étudier la structure du \mo (multiplicatif) des \itfs de~$\gZ$.\index{monoide@monoïde!des ideaux@des \itfs}

On suppose que l'on dispose d'un anneau $\gZ'$ qui constitue une approximation de $\gZ$
en ce sens que $\ZZ[\alpha]\subseteq\gZ'\subseteq\gZ$. Par exemple en un premier
\hbox{temps $\gZ'=\ZZ[\alpha]$}. On est intéressé par la structure multiplicative
du groupe des \ifrs de $\gZ$(\footnote{Un \ifr de $\gZ$ est un sous-$\gZ$-module de $\gK$ de la forme $\fraC{1}m\,\fa$ pour \hbox{un $m\in\ZZ\sta$} et un \itf $\fa$ de $\gZ$, cf. \paref{NOTAIfr}.}), et l'on veut s'appuyer sur celle de $\gZ'$
pour l'étudier en détail.

Le \tho qui suit dit que \gui{cela marche très bien pour la plupart des \idsz, \cad pour tous ceux qui évitent le conducteur de $\gZ'$ dans~$\gZ$}.  

%:     Definition{defiConducteur}
\begin{definition}\label{defiConducteur}~
\\
Soient deux anneaux $\gA \subseteq \gB$, $\fa$ un \id de $\gA$
et $\fb$ un \id de $\gB$. 
\begin{enumerate}
\item Le \emph{conducteur de $\gA$ dans
$\gB$} est  $(\gA:\gB) = \sotq {x \in \gB} {x\gB \subseteq \gA}$.% 
\index{conducteur!d'un anneau dans un sur-anneau}
\item L'\emph{extension de $\fa$} est l'\id $\fa\gB$ de $\gB$.%
\index{extension!d'un idéal dans un sur-anneau}
\item La \emph{contraction de $\fb$}  est l'\id $\gA \cap \fb$ de $\gA$.%
\index{contraction!d'un idéal dans un sous-anneau}
\end{enumerate} 
\end{definition}
%--------- fin definition ---------------------------------------------- 
\rem La terminologie concernant le conducteur est flottante.
Des auteurs disent \gui{conducteur de $\gB$ dans $\gA$} là où nous disons
\gui{conducteur de~$\gA$ dans $\gB$}. Pour eux, conducteur est synonyme de
transporteur. En théorie des nombres, Dedekind a introduit la notion
de conducteur en tant qu'\id attaché au \gui{petit anneau} (un sous-anneau $\gA$ de l'anneau d'entiers~$\gZ$ d'un corps de nombres, avec même corps de fractions).
\eoe

%%%%%%%%%%%%%%%%%%%%%%%%%%%%%%%%%%%%%%%%%%%%%%%%%%%%%%%%%%%%%%%%%%%%%%%%%%%
%:     theorem{propEvitementConducteur}
\begin{theorem}\label{propEvitementConducteur} 
\emph{(\Tho de Dedekind, \ids qui évitent le conducteur)}\index{Dedekind!idéaux qui évitent le conducteur}
Soient $\gA \subseteq \gB$ deux anneaux et $\ff$ le conducteur de $\gA$ dans~$\gB$. 
\begin{enumerate}
\item L'\id $\ff$ est
l'annulateur du \Amo $\gB\sur\gA$. C'est à la
fois un \id de $\gA$ et un \id de $\gB$,  et c'est le plus grand \id pour cette
\prtz.
\end{enumerate}
 On note $\cA$ (resp. $\cB$) la classe des \ids de $\gA$ (resp. de $\gB$)
\comz~à~$\ff$.
\begin{enumerate}\setcounter{enumi}{1}
\item Pour $\fa \in \cA$, on a $\gA\sur\fa \simeq \gB\sur{\fa\gB}$ et
pour $\fb \in \cB$, on a $\gB\sur\fb \simeq \gA\sur{\gA\cap\fb}$.

\item $\cA$ est stable par multiplication,  somme, intersection et vérifie: 

\snic{\fa \in \cA, \,\fa'
\supseteq \fa \;\;\;\Longrightarrow\;\;\; \fa' \in \cA.}

%\sni
En particulier, $\fa_1\fa_2 \in \cA$
\ssi $\fa_1$ et $\fa_2 \in \cA$. 
Les mêmes \prts sont valables pour~$\cB$.
\item L'extension et la contraction, restreintes respectivement à $\cA$ et $\cB$,
sont deux correspondances réciproques l'une de
l'autre. Elles préservent la multiplication, l'inclusion, l'intersection et le \crc \tfz.
\item On suppose $\gB$ intègre.  Alors, un \id $\fa\in\cA$ est \iv dans
$\gA$ \ssi $\fa\gB$ l'est dans $\gB$. De même, un \id $\fb\in\cB$ est \iv dans $\gB$ \ssi $\gA \cap \fb$ l'est dans~$\gA$.

\end{enumerate}
\end{theorem}
%--------- fin theorem ---------------------------------------------- 
 
%
\begin{proof}
On montre seulement quelques \prtsz. 
Remarquons que l'on a toujours les inclusions $\fa \subseteq \gA\cap
\fa\gB$ et $(\gA\cap\fb)\gB \subseteq \fb$.

Soit $\fa \in \cA$, donc $1 = a + f$ avec $a \in \fa$ et $f \in \ff$;
a fortiori, $1 \in \fa\gB + \ff$. Montrons que $\gA\cap\fa\gB = \fa$.
On prend $x \in \gA\cap\fa\gB$ et l'on écrit 

\snic{x = xf + xa \in \fa\gB\ff + \fa \subseteq
\fa\gA + \fa = \fa\,.}

%\sni
D'où le résultat. On voit aussi que $\gB = \gA + \fa\gB$,
donc le morphisme composé $\gA \to \gB\sur{\fa\gB}$ est surjectif de noyau $\fa$, ce qui donne 
un \iso $\gA\sur\fa \simeq \gB\sur{\fa\gB}$.

Soit $\fb \in \cB$, donc $1 = b + f$ avec $b \in \fb$, $f \in \ff$.
Puisque $\ff \subseteq \gA$, on a $b \in \gA \cap \fb$ donc 
$1 \in \gA\cap \fb + \ff$. Montrons que $(\gA \cap\fb)\gB = \fb$.
\\
Si $x \in \fb$, alors:

\snac {
x = (b+f)x = bx + xf \in (\gA \cap\fb)\gB + \fb\ff \subseteq 
(\gA \cap\fb)\gB + \gA\cap\fb \subseteq (\gA \cap\fb)\gB.
}

%\sni
Ainsi $\fb \subseteq (\gA \cap\fb)\gB$ puis $\fb = (\gA \cap\fb)\gB$.
De plus, puisque $\gB = \fb + \ff \subseteq \fb + \gA$, le morphisme composé
$\gA \to \gB\sur\fb$  est
surjectif, de noyau $\gA\cap \fb$, ce qui donne un \iso $\gA\sur{\gA\cap
\fb} \simeq \gB\sur\fb$.

L'extension est multiplicative, donc la contraction (restreinte à $\cB$)
qui est son inverse, est \egmt multiplicative. La contraction est compatible
avec l'intersection, donc l'extension (restreinte à $\cA$) qui est son
inverse, est \egmt compatible avec l'intersection.

Soit $\fb=\gen{b_1, \ldots, b_n}_\gB \in \cB$.  Montrons que
$\gA\cap\fb$ est  \tfz.  \\
On écrit $1 = a + f^2$ avec $a \in
\fb$, $f \in \ff$. Puisque $f \in \gA$, on a $a \in \gA\cap\fb$. Montrons que
$(a, fb_1, \ldots, fb_n)$ est un \sgr de $\gA\cap\fb$.  \\
Soit $x \in \gA\cap\fb$
que l'on écrit $x = \sum_i y_ib_i$ avec $y_i \in \gB$, alors:

\snic {
x = \sum_i (y_i a + y_i f^2) b_i =  
xa + \sum_i (y_if)fb_i\in \gen {a,fb_1, \ldots, fb_n}_\gA.
}

%\sni
Pour un \id $\fb\in\cB$ (non \ncrt \tfz), on a en fait montré le
résultat suivant: si $1 = a + f^2$ avec $a \in \fb$ et $f \in \ff$, alors
$\gA\cap\fb = \gA a + f(f\fb)$ (et $f\fb$ est un \id de $\gA$).

Soit $\fb \in \cB$ un \id\ivz, montrons que $\fa = \gA\cap\fb$ est un
\id\ivz. On écrit $1 = a + f$ avec $a \in \fb$ et $f \in \ff$, de sorte que
$a \in \fa$. \\
Si $a = 0$, alors $1 = f \in \ff$, donc $\gA = \gB$ et
il n'y a rien à montrer. Sinon, $a$ est \ndz et il existe un \id
$\fb'$ de $\gB$ tel que $\fb\fb' = a\gB$.  \\
Puisque les \ids $a\gB$, $\fb$ et $\fb'$
sont comaximaux à $\ff$, on peut appliquer le \crc multiplicatif de
la contraction à l'\egt $\fb\fb' = a\gB$ pour obtenir 
l'\egt $\fa\fa' = a\gA$ avec
$\fa' = \gA\cap\fb'$.
\end{proof}

%: entrenous
\entrenous{ \par
1) Donner des exemples pour lesquels $\gA \cap \fa\gB \ne \fa$ et
$(\gA \cap \fb)\gB \ne \fb$.

2) Il faudrait aussi peut-être faire le lien avec certains exos.
En fait le remarquable \tho précédent est insuffisamment exploité.

3) 
%:     Proposition{prop0ITFSCDN}
%\begin{proposition}\label{prop0ITFSCDN}
\textbf{Proposition } (\tho un et demi)~~
Encombrant à faire ici, voir le corolaire \ref{corpropZerdimLib}, mais c'est fait \gui{à la Hilbert} dans le livre Modules
%\end{proposition}
} % finentrenous
%%%%%%%%%%%%%%%%%%%%%%%%%%%%%%%%%%%%%%%%%%%%%%%%%%%%%%%%%%%%%%%%%%%%%%%%%%%

%%%%%%%%%%%%%%%%%%%%%%%%%%%%%%%%%%%%%%%%%%%%%%%%%%%%%%%%%%%%%%%%%%%%%%%%%%%
%  sec{Le \nst} {secChap3Nst} 
\section{Le \nst de Hilbert}\label{secChap3Nst}\ihi   
%-----------------------------------------

Nous illustrons dans cette section l'importance du résultant en montrant comment on peut en déduire le \nst de Hilbert. Nous utiliserons une \gnn du  lemme d'\eli de base.

%:    subsec{Clôture algébrique}
\subsec{Clôture algébrique de $\QQ$ et des corps finis}

Soient $\gK\subseteq\gL$ des \cdisz, on dit que \emph{$\gL$ 
est une clôture \agq
 de~$\gK$} si~$\gL$ est  \agq sur $\gK$ et \agqt clos.%
\index{cloture@clôture!algébrique}

\Llec nous accordera que $\QQ$ et les corps $\FFp$ possèdent une clôture \agqz.
Ceci sera discuté plus en détail dans la section \ref{secEtaleSurCD},
avec notamment le \thref{thClsep}.

%:    subsec{Le \nst classique (cas \agqt clos)}
\subsec{Le \nst classique (cas \agqt clos)}

Le \nst est un \tho qui concerne les \syss d'\eqns \polles sur un corps discret.
De manière très informelle sa signification peut être décrite comme suit: une affirmation de nature \gmq 
possède \ncrt un certificat \agqz. Ou encore: une \dem en \alg 
commutative peut (presque)
toujours, si elle est suffisamment \gnlez, être résumée 
par de simples \idasz.

Si l'on a des \cdis $\gK\subseteq\gL$, et si $(\uf)=(\lfs)$ est un \sys de \pols dans $\KXn=\KuX$,
on dit que $(\xin)=(\uxi)$ est un \emph{zéro de $(\uf)$ dans $\gL^n$}, ou encore un \emph{zéro de
$(\uf)$ à \coos dans $\gL$}, si les \eqns $f_i(\uxi)=0$ sont satisfaites.
Notons $\ff=\gen{\lfs}_\KuX$. Alors, tous les \pols $g\in\ff$ s'annulent en un tel $(\uxi).$
On parle donc aussi bien de~$(\uxi)$ comme \emph{zéro de l'\id $\ff$ dans $\gL^n$},
ou \emph{à \coos dans $\gL$}.

 Nous commençons par un fait presque évident.   

%:     Fact{factGCDDeg}
\begin{fact}\label{factGCDDeg} Soit $\gk$ un anneau commutatif
 et $h\in\kX$ un \polu de degré $\geq1$. 
 \begin{itemize}
\item Si un multiple de $h$
est dans $\gk$, ce multiple est nul.
\item  Soient $f$ et $g\in\kX$ de degrés formels $p$ et  $q$. Si $h$   divise~$f$ et~$g$, alors~$\Res_X(f,p,g,q)=0$.
\end{itemize}
\end{fact}
%--------- fin fact ---------------------------------------------- 

Voici maintenant une \gnn  du lemme  d'\eli de base \ref{LemElimAffBasic}.

%:     Lemma{lemElimPlusieurs}
\begin{lemma}\label{lemElimPlusieurs} \emph{(\'Elimination d'une variable entre plusieurs \polsz)}
\label{lemElimParametre}%
\index{elimination@\eliz!d'une variable}
\\
Soient $f$,  $g_1$, $\ldots$, $g_r$  $\in\kX$ ($r\geq1$), avec $f$  \mon de degré $d$. \\
On pose $\ff=\gen{f,g_1,\ldots,g_r}$ et $\fa=\ff\cap\gk$ (c'est l'\id d'\eli de la variable $X$ dans $\ff$).
On pose aussi:
\[ 
\begin{array}{c} 
%\hbox{les \pols}  &  
g(T,X)=g_1+Tg_2+\cdots+T^{r-1}g_r\in\gk[T,X],    \\[2mm] 
%\hbox{et}  &  
R(T)=R(f,g_1,\ldots,g_r)(T)=\Res_X\big(f,g(T,X)\big)\in\gk[T],    \\[1.3mm] 
%\hbox{et l'\id}  &  
\fb=\fR(f,g_1,\ldots,g_r)\eqdefi\rc_{\gk,T}\big(R(f,g_1,\ldots,g_r)(T)\big)\subseteq\gk.  
 \end{array}
\]
\begin{enumerate}
% 1
\item  L'\id $\fb$ est engendré par $d(r-1) + 1$ \elts et
l'on a les inclusions: 
\begin{eqnarray}
\label{eq0lemElimPlusieurs} %\fa^{d}\subseteq
\fb\subseteq\fa\subseteq\sqrt\fb=\sqrt\fa\,.
\end{eqnarray}
Plus \prmtz, posons $e_i = 1 + (d-i)(r-1)$, $i \in \lrb{1..d}$, alors pour
des \elts arbitraires $a_1$, $\ldots$, $a_d \in \fa$, on a:
$$
a_1^{e_1} a_2^{e_2} \cdots a_d^{e_d}\, \in \,\fR(f, g_1,\ldots, g_r)\,.
$$
%\sni
En particulier, on a les \eqvcs suivantes. 
\begin{eqnarray}
\label{eq01lemElimPlusieurs} 1\in \fb\; \iff\; 1\in\fa \;\iff\; 1\in\ff\,.~\;\; 
\end{eqnarray}

% 2
\item Si $\gk$ est un \cdi contenu dans un \cac discret $\gL$, notons $h$ le pgcd \mon de $f$, $g_1$, $\ldots$, $g_r$ et $V$ l'ensemble des zéros de $\ff$
dans $\gL^n$. Alors, on a les \eqvcs suivantes: 
\begin{eqnarray}
\label{eqlemElimPlusieurs} 1\in \fb\; \iff\; 1\in\fa \;\iff\; 1\in\ff \;\iff\;h=1\;\iff\;
V=\emptyset~\;\;~\;\;
\end{eqnarray}

%
%\item 
%
\end{enumerate}
\end{lemma}
%--------- fin lemma ---------------------------------------------- 
%
\begin{proof}
\emph{1.} On sait que $R(T)$
s'écrit 

\snic{u(T,X)f(X)+v(T,X)g(T,X),}

%\sni
donc chaque \coe de $R(T)$ est une \coli de $f$ et des~$g_i$ dans~$\kX$. Ceci
donne l'inclusion $\fb\subseteq\fa$. L'in\egt $\deg_T(R)\leq d(r-1)$
donne la majoration pour le nombre de \gtrs de $\fb$.
\\
Si $f_1$, $\ldots$, $f_d$ sont $d$ \pols (à une \idtrz) de degré $< r$, on déduit 
du lemme de \DKM (voir corolaire \ref{corLDM}) l'inclusion suivante.

\snic {
(\star) \qquad\qquad\qquad
\rc(f_1)^{e_1} \rc(f_2)^{e_2} \cdots \rc(f_d)^{e_d}  \subseteq \rc(f_1f_2\cdots f_d).
\qquad\qquad\qquad}

%\sni
Supposons $f(X) = (X-x_1) \cdots (X-x_d)$. On pose alors, pour $i \in \lrb{1..d}$,

\snic {
f_i(T) = g_1(x_i) + g_2(x_i)T + \cdots + g_r(x_i)T^{r-1},
}

%\sni
de sorte que $f_1 f_2 \cdots f_d = \Res_X(f, g_1 + g_2T + \cdots +
g_rT^{r-1})$.  
\\
Ainsi, pour $a_j \in \fa= \gen {f, g_1, \ldots, g_r}_{\kX} \cap \gk$, en évaluant
en $x_i$, on obtient~\hbox{$a_j \in \gen {g_1(x_i), \ldots, g_r(x_i)} =
\rc(f_i)$}.  En appliquant l'inclusion $(\star)$ on obtient l'appartenance
$a_1^{e_1} a_2^{e_2} \cdots a_d^{e_d}\in\fb$.
\\
Passons au cas \gnlz. On considère l'\aduz~\hbox{$\gk'=\Adu_{\gk,f}$}.
Le calcul précédent vaut pour $\gk'$. Comme $\gk'=\gk\oplus E$
en tant que \kmoz, on a l'\egt
$(\fb\gk')\cap \gk=\fb$. Pour des $a_j\in\fa$, cela permet de conclure
que $a_1^{e_1} a_2^{e_2} \cdots a_d^{e_d}\in\fb$, car le produit est dans
$(\fb\gk')\cap \gk$.

\emph{2.} 
Par \dfn du pgcd, on a $\ff = \gen {h}$. Par ailleurs, $h=1\Leftrightarrow V=\emptyset$. Donc tout est clair d'après le point \emph{1.} 
\\
 \emph{Voici cependant pour ce cas particulier
 une \dem plus directe, qui donne l'origine de la
\dem magique du point {1.}}
\\
Supposons que $h$ soit égal à $1$; alors dans ce cas $1\in\ff$ et
$1\in\fa$.  Supposons ensuite que $h$ soit de degré $\geq 1$; alors
$\fa=\gen{0}$.  On a donc obtenu les \eqvcs $\;1\in\fa \iff 1\in\ff\iff
\deg(h)=0\;$ et $\;\fa=\gen{0} \iff \deg(h) \ge 1$.
\\
Montrons maintenant l'\eqvc $\;\deg(h) \ge 1 \iff \fb = \gen{0}$. \\
 Si
$\deg(h)\geq1$, alors $h(X)$ divise $g(T,X)$, donc $R(f,g_1,\ldots,g_r)(T)=0$ (fait~\ref{factGCDDeg}), i.e. $\fb=\gen{0}$.
\\ Inversement, supposons $\fb=\gen{0}$. Alors, pour toute valeur du
para\-mètre~$t\in\gL$, les \pols $f(X)$ et $g(t,X)$ ont un zéro en commun
dans $\gL$ ($f$ est \mon et le résultant des deux \pols est nul).
\\
Considérons les
zéros $\xi_1$, $\ldots$, $\xi_d \in \gL$ de $f$. En prenant $d(r-1)+1$ valeurs distinctes de
$t$, on trouvera un $\xi_\ell$ tel que $g(t, \xi_\ell) = 0$ pour au moins
$r$ valeurs de $t$.  Ceci implique que $g(T,\xi_\ell)$ est identiquement
nul, i.e. $\xi_\ell$ annule tous les $g_i$, et que $h$ est multiple de
$X-\xi_\ell$ donc $\deg(h) \ge 1$.
\end{proof}

Le point \emph{2} du lemme \ref{lemElimPlusieurs} donne le corolaire suivant.
%:     Corollary{corlemElimPlusieurs}
\begin{corollary}\label{corlemElimPlusieurs}
Soit $\gK$  un \cdi 
non trivial contenu dans
un \cac $\gL$. Reprenons les hypothèses du lemme \ref{lemElimPlusieurs},
avec l'anneau $\gk=\gK[X_1,\ldots,X_{n-1}]$. Alors, pour  $\alpha=(\alpha_1,\ldots,\alpha_{n-1})\in\gL^{n-1}$ 
\propeq
\begin{enumerate}
\item Il existe $\xi\in\gL$ tel que $(\alpha,\xi)$ annule $(f,g_1,\ldots,g_r)$.
\item $\alpha$ est un zéro de l'\id $\fb=\fR(f,g_1,\ldots,g_r)\subseteq\gk$.
\end{enumerate}
Précisions: si le degré total des \gtrs de $\ff$ est majoré par $d$,
on obtient comme \gtrs de $\fb$, $d(r-1)+1$ \pols de degré total majoré
 par~$2d^2$.  
\end{corollary}
%--------- fin corollary ---------------------------------------------- 

\rem Le corolaire précédent a la structure voulue pour enchaîner une 
\recu qui nous permet une description des zéros de $\ff$ dans $\gL^n$. 
\\
En effet, en partant de l'\itf $\ff\subseteq\KXn$ on produit un \itf $\fb\subseteq\gk$ avec la \prt suivante: \emph{les
zéros de $\ff$ dans $\gL^n$ se projettent exactement sur les zéros de $\fb$ 
dans $\gL^{n-1}$}. Plus \prmt au dessus de chaque
zéro de $\fb$ dans $\gL^{n-1}$ se trouve un nombre fini, non nul,
majoré par $\deg_{X_n}(f)$
de zéros de $\ff$ dans $\gL^n$.
\\
 Donc
ou bien tous les \gtrs de $\fb$ sont nuls et le processus décrivant les zéros de $\ff$ est terminé, 
ou bien un des \gtrs de $\fb$
est non nul et l'on est près à faire à $\fb\subseteq\gK[X_1,\ldots,\alb X_{n-1}]$ ce que l'on a fait à
$\ff\subseteq\KXn$ \emph{à condition toutefois de trouver} un \pol unitaire en $X_{n-1}$ dans l'\id $\fb$.
\\
 Cette dernière question est réglée par le lemme de \cdv suivant. 
\eoe

%:    Lemma{lemCDV}---------------
\begin{lemma}
\label{lemCDV} \emph{(Lemme de \cdvsz)}\\
Soit $\gK$  un \cdi infini  et  $g\neq0$ dans $\KuX=\KXn$ de degré~$d$.
Il existe $(a_1, \ldots, a_{n-1}) \in \gK^{n-1}$ tel que le \pol

\snic {
g(X_1 + a_1X_n, \ldots, X_{n-1} + a_{n-1} X_n, X_n) 
}

%\sni
soit de la forme $aX_n^d + h$ avec $a \in \gK^\times$ et $\deg_{X_n} h < d$.
\end{lemma}
%--- end-lemma-----------------------------------------
%-----------------begin proof------------------
\begin{proof}
Soit $g_d$ la composante \hmg de degré $d$ de $g$. Alors:

\snic {
g(X_1 + a_1X_n, \ldots, X_{n-1} + a_{n-1} X_n, X_n) =
g_d(a_1, \ldots, a_{n-1}, 1) X_n^d + h, 
}

%\sni
avec $\deg_{X_n} h < d$. Comme $g_d(X_1, \ldots, X_n)$ est \hmg non nul,
le \polz~$g_d(X_1, \ldots, X_{n-1}, 1)$ est non nul. Il existe donc $(a_1,
\ldots, a_{n-1}) \in \gK^{n-1}$ tel que $g_d(a_1, \ldots, a_{n-1}, 1) \ne 0$.
\end{proof}
%-----------------end proof------------------

On obtient maintenant un \gui{\nst faible} 
(\cad l'\eqvc entre~$\,V=\emptyset\,$  et~$\,\gen{\lfs}=\gen{1}\,$ dans le \thoz)
et une \gui{mise en position de \iNoez} qui donne une description de $V$ dans le cas non vide.

%:     Theorem{thNstfaibleClass}
\begin{theorem}\label{thNstfaibleClass}\emph{(\nst faible et mise en position de \iNoez)}\\
Soit $\gK$  un \cdi infini contenu dans
un \cac $\gL$ et $(\lfs)$ un \syp dans $\KXn$.%
\index{systeme polynomial@\sypz}
\\
 Notons $\ff=\gen{\lfs}_\KuX$ et $V$ la \vrt des zéros de $\ff$ dans $\gL^n$.
\begin{enumerate}
\item Ou bien  $\gen{\lfs}=\gen{1},$ et $V=\emptyset$.
\item Ou bien $V\neq\emptyset$. Alors, il existe un entier $r\in\lrb{0..n}$, un \cdv $\gK$-\lin (les nouvelles variables sont notées $\Yn$), et des \itfs $\ff_j\subseteq\gK[Y_1,\ldots,Y_j]$
$(j\in\lrb{r..n})$, qui satisfont les \prts suivantes.
\begin{itemize}
\item [$\bullet$] On a $\ff\,\cap\,\gK[\Yr]=0$. Autrement dit, l'anneau  $\gK[\Yr]$
s'iden\-tifie à un sous-anneau du quotient $\KuX\sur\ff$.
\item [$\bullet$] Chaque $Y_j$ $(j\in\lrb{r+1..n})$ est entier sur $\gK[\Yr]$ modulo $\ff$.
Autrement dit l'anneau $\KuX\sur\ff$ est entier sur le sous-anneau $\gK[\Yr]$. 
\item [$\bullet$] On a les inclusions $\gen{0}=\ff_r\subseteq\ff_{r+1} \subseteq \ldots \subseteq\ff_{n-1}\subseteq\ff$ et pour chaque~\hbox{$j\in\lrb{r..n}$}
l'\egt $\sqrt {\ff\,} \cap \gK[Y_1,\ldots,Y_j]=\sqrt {\ff_j}$.
\item [$\bullet$] Pour les nouvelles \coos correspondant aux $Y_i$,
soit  $\pi_j$ la \prn $\gL^n\to\gL^j$
qui oublie les dernières \coos ($j\in\lrbn$).
Pour chaque $j\in\lrb{r..n-1}$ 
la \prn  de la \vrt $V\subseteq\gL^n$ sur~$\gL^j$ 
est exactement la \vrt $V_j$ des
zéros de $\ff_j$. En outre, pour chaque \elt~$\alpha$ de~$V_j$, la fibre 
$\pi_j^{-1}(\alpha)$ est finie, non vide, avec un nombre d'\elts uniformément borné.
\end{itemize}
\end{enumerate}
En particulier:
\begin{itemize}
\item Ou bien $V$ est vide (et l'on peut convenir de~$r = -1$).
\item Ou bien $V$ est finie non vide, $r=0$
et les \coos des points de $V$ sont \agqs sur~$\gK$. 
\item Ou bien $r\geq1$ et la \prn $\pi_r$ envoie surjectivement la \vrt $V$ sur $\gL^r$
(donc $V$ est infinie).
Dans ce cas, si $\alpha\in\gK^r$, les  \coos des points de~$\pi_r^{-1}(\alpha)$
sont  \agqs sur~$\gK$.
\end{itemize}
\end{theorem}
%--------- fin theorem ---------------------------------------------- 
%
\begin{proof}
On raisonne comme on l'a indiqué dans la remarque précédant le lemme de \cdvz. Notons que la première étape du processus n'a lieu que si le \syp de départ est non nul, auquel cas la première opération consiste en un \cdv \lin qui rend l'un des $f_i$ \mon en~$Y_n$.
\end{proof}

\rems ~\\
 1) Le nombre $r$ ci-dessus correspond au nombre maximum d'\idtrs
pour un anneau de \pols $\gK[Z_1,\ldots,Z_r]$ qui soit isomorphe à une sous-\Klg
de $\aqo\KuX{\lfs}$. Ceci est relié à la théorie de la \ddk
qui sera exposée au chapitre~\ref{chapKrulldim} (voir notamment le \thref{thDKAG}).

2) Supposons que les degrés des $f_j$ soient majorés par $d$.
\\
 En s'appuyant sur le résultat
énoncé à la fin du corolaire \ref{corlemElimPlusieurs}, on peut donner des précisions dans le \tho précédent en calculant a priori,
uniquement en fonction des entiers $n$, $s$, $j$ et $d$,
\begin{itemize}
\item d'une part une majoration pour le nombre de \gtrs pour chaque \idz~$\ff_j$,
\item d'autre part une majoration pour les degrés de ces \gtrsz.
\end{itemize}

3) Le calcul des \ids $\ff_j$ ainsi que toutes les affirmations du \tho 
qui ne concernent pas la variété $V$ sont valables même lorsque l'on ne connaît pas de \cac $\gL$ contenant $\gK$.
On utilise pour cela seulement les lemmes \ref{lemElimPlusieurs} 
et~\ref{lemCDV}.
On reprendra ceci dans les \thosz~\ref{thNstfaibleClassSCA} et \ref{thNstNoe}.
\eoe

\medskip
La restriction introduite par l'hypothèse \gui{$\gK$ est infini} 
va disparaître dans le \nst classique en raison de la constatation suivante.

%:     Fact{factNstRationnel}
\begin{fact}\label{factNstRationnel}
Soit $\gK\subseteq\gL$ des \cdis et $h$, $f_1$, \dots, $f_s\in\KXn$, alors 
$h\in\gen{\lfs}_\KXn\iff h\in\gen{\lfs}_\LXn$.  
\end{fact}
%--------- fin fact ---------------------------------------------- 
%
\begin{proof}
En effet, une \egt $h=\sum_ia_if_i$, une fois les degrés des $a_i$ fixés,
peut être vue comme un \sli dont les inconnues sont les \coes des $a_i$.
Le fait qu'un \sli admet une solution ne dépend pas du corps dans lequel on cherche la solution, du moment qu'il contient les \coes du \sliz: la méthode du pivot est un processus entièrement rationnel.
\end{proof}

Comme corolaire  du \nst faible et du fait précédent on obtient le \nst classique. 

%:     Theorem{thNstClass}
\begin{theorem}\label{thNstClass} \emph{(\nst classique)}
\\
Soit $\gK$  un \cdi contenu dans
un \cac $\gL$ et des \pols $g$, $f_1$, \dots, $f_s$  dans $\KXn$. 
Notons  $V$ la \vrt des zéros de $(\lfs)$ dans $\gL^n$. 
\begin{enumerate}
\item Ou bien il existe un point $\xi$ de $V$ tel que $g(\xi)\neq0$.
\item Ou bien il existe un entier $N$ tel que $g^N\in\gen{\lfs}_\KuX$.
\end{enumerate}
\end{theorem}
%--------- fin theorem ---------------------------------------------- 
%
\begin{proof} Le cas $g=0$ est clair, on suppose $g\neq 0$.
On applique \emph{l'astuce de Rabinovitch}\index{Rabinovitch!astuce de ---}, \cad on introduit une \idtr supplémentaire
$T$ et l'on remarque que $g$ s'annule aux zéros de $(\lfs)$ \ssi le système
$(1-gT,\lfs)$ n'admet pas de solution. On applique alors le \nst faible à ce nouveau \sypz, avec $\gL$ (qui est infini) à la place de $\gK$.
On obtient dans $\KuX[T]$ (gr\^ace au fait \ref{factNstRationnel}) une \egt
$$
\big(1-g(\uX)T\big)a(\uX,T) + f_1(\uX)b_1(\uX,T)+\cdots+f_s(\uX)b_s(\uX,T)=1.
$$
Dans le localisé $\KuX[1/g]$, on réalise la substitution $T = 1/g$.
Plus \prmtz, en restant dans $\gK[\uX,T]$, si $N$ est le plus grand des degrés en $T$ des~$b_i$, on multiplie l'\egt précédente par $g^N$ et l'on remplace dans  $g^Nb_i(\uX,T)$ chaque $g^NT^k$ par $g^{N-k}$
modulo $(1-gT)$. On obtient alors une \egt
$$
\big(1-g(\uX)T\big)a_1(\uX,T) + f_1(\uX)c_1(\uX)+\cdots+f_s(\uX)c_s(\uX)=g^N,
$$
dans laquelle \ncrt $a_1=0$, puisque, si l'on regarde $a_1$ dans~$\KuX[T]$, son \coe \fmt dominant en $T$ est nul.
\end{proof}

\rem On notera que la séparation des différents cas dans les \thosz~\ref{thNstfaibleClass}
et~\ref{thNstClass} est explicite.
\eoe

%:     Corollary{corthNstClass}
\begin{corollary}\label{corthNstClass}~
Soit $\gK$  un \cdi  contenu dans
un \cac $\gL$ et deux \itfs $\fa=\gen{\lfs}$, $\fb$ de $\KXn$.
Soit $\gK_0$ le sous-corps de $\gK$ engendré par les \coes des $f_i$.
\\
\Propeq
\begin{enumerate}
\item $\fb\subseteq\rD_{\KuX}(\fa)$.
\item $\fb \subseteq\rD_{\LuX}(\fa)$.
\item  Tout zéro de $\fa$ dans $\gL^n$ est un zéro de $\fb$.
\item  Pour tout sous-corps $\gK_1$ de $\gL$ fini sur $\gK_0$, 
tout zéro de $\fa$  dans $\gK_1^n$ est un zéro de~$\fb$.
\end{enumerate}
En particulier,
 $\rD_{\KuX}(\fa)=\rD_{\KuX}(\fb)$ \ssi $\fa$ et $\fb$ ont les mêmes zéros dans $\gL^n$.
\end{corollary}
%--------- fin corollary ---------------------------------------------- 
%
\begin{proof}
Conséquence \imde du \nstz. 
\end{proof}
%

%:    subsec{Le \nst formel}
\subsec{Le \nst formel}
Nous passons maintenant à un \emph{\nst formel}, formel en ce sens qu'il
s'applique  (en \clamaz) à un \id arbitraire sur un anneau arbitraire. 
Néanmoins pour avoir un énoncé \cof nous nous contenterons d'un anneau de \pols $\ZZuX$ pour notre anneau arbitraire et d'un \itf  pour notre \id arbitraire. 

Cela peut sembler très restrictif, mais la pratique montre que ce n'est pas le cas, en raison du fait que l'on peut (presque) toujours appliquer la méthode des \coes indé\-ter\-minés à un \pb d'\alg commutative, méthode qui ramène le \pb à un \pb \poll sur $\ZZ$. Une illustration en sera donnée ensuite.

Notons que pour lire l'énoncé, lorsque l'on parle d'un zéro d'un $f_i\in\ZZuX$ sur un anneau $\gA$, il faut d'abord voir $f_i$ modulo $\Ker\varphi$, 
où $\varphi$
est l'unique \homo $\ZZ\to\gA$, d'image $\gA_1\simeq\ZZ\sur{\Ker\varphi}$. 
On se ramène ainsi à un \polz~$\ov {f_i}$ de $\gA_1[\uX]\subseteq\gA[\uX]$.

%:     Theorem{thNSTsurZ}
\begin{theorem}\label{thNSTsurZ}
\emph{(\nst sur $\ZZ$, \nst formel)} 
\\
On écrit $\ZZuX=\ZZXn$
On considère  $g$, $f_1$, \dots, $f_s$ dans $\ZZuX$
\begin{enumerate}
\item Pour le \sys $(\lfs)$ \propeq
\begin{enumerate}
\item $1\in\gen{\lfs}$.
%
%\item Si le \sys $(\lfs)$ admet un zéro dans un anneau, cet anneau est trivial.
%
\item Le \sys  n'admet de zéro sur aucun \cdi non trivial.
\item Le \sys n'admet de zéro sur aucun corps  fini et sur
aucune extension finie de $\QQ$.
\item Le \sys n'admet de zéro sur aucun corps  fini.
\end{enumerate}
\item \Propeq
\begin{enumerate}
\item $\exists N \in \NN,\;g^N\in\gen{\lfs}$.
%
%\item Même chose si l'on évalue les \idtrs dans un anneau  arbitraire.
%
\item Le \polz~$g$ s'annule aux zéros du \sys $(\lfs)$  sur n'importe quel \cdiz.
\item Le \polz~$g$ s'annule aux zéros du \sys $(\lfs)$  sur tout corps  fini et sur
toute extension finie de $\QQ$.
\item Le \polz~$g$ s'annule aux zéros du \sys $(\lfs)$  sur tout corps  fini.
\end{enumerate}
\end{enumerate}
\end{theorem}
%--------- fin theorem ---------------------------------------------- 

%
\begin{proof}
Il suffit de démontrer la version faible \emph{1}, car on passe ensuite à la version \gnle \emph{2} en appliquant l'astuce de Rabinovitch.
Pour ce qui concerne la version faible, 
la chose difficile est l'implication \emph{d} $\Rightarrow
$ \emph{a}.

Voyons d'abord \emph{c} $\Rightarrow
$ \emph{a}. On applique le \nst faible
en considérant~\hbox{$\ZZ\subseteq\QQ$}. 
Cela donne une appartenance:
$$
m \in \gen {f_1, \ldots, f_s}_{\ZZuX}
\quad\hbox{ avec } m\in\ZZ\setminus\so0\eqno(\star_\QQ)
.$$ 
En appliquant le \nst faible avec une clôture \agq $\gL_p$ de $\FFp$
on obtient aussi pour chaque nombre premier $p \divi m$ une appartenance:
$$
1 \in \gen {f_1, \ldots, f_s}_{\ZZuX} +p\ZZuX \eqno(\star_{\FFp})
.$$ 
Or dans n'importe quel anneau, pour trois \ids quelconques $\fa, \fb, \fc$,
on a l'inclusion $(\fa + \fb)(\fa + \fc) \subseteq \fa + \fb\fc$.
En écrivant le $m$ ci-dessus dans $(\star_\QQ)$ sous forme $\prod_jp_j^{k_j}$
avec les $p_j$ premiers, on obtient donc
$$
1 \in \gen {f_1, \ldots, f_s}_{\ZZuX} + m\ZZuX.
$$
Cette appartenance, jointe à  $(\star_\QQ)$, fournit 
$1 \in \gen {f_1, \ldots, f_s}_{\ZZuX}$.

\emph{d} $\Rightarrow$ \emph{c.} Nous montrons qu'un zéro $(\uxi)$ du  \sys $(\lfs)$ dans une extension finie de $\QQ$  donne lieu
à un zéro de $(\lfs)$ dans une extension finie de $\FFp$ pour tous les nombres premiers, à l'exception d'un nombre fini d'entre eux.
\\
 En effet, soit $\gQ=\QQ[\alpha]\simeq \aqo{\QQX}{h(X)}$ 
(avec $h$ \mon \ird
dans~$\ZZX$) une extension finie de $\QQ$ et $(\uxi)\in\gQ^n$ un zéro de 
$(\lfs)$.
Si~\hbox{$\xi_j=q_j(\alpha)$} avec $q_j\in\QQX$ pour $j\in\lrbn$, cela signifie que 
$$
f_i(q_1,\ldots,q_n)\equiv 0\mod h \quad\hbox{dans }\QQX, \;i\in\lrbs.
$$ 
Ceci reste vrai dans $\FFp[X]$ dès qu'aucun des \denos figurant dans les $q_j$
n'est multiple de $p$, à condition de prendre les fractions dans $\FFp$:
$$
\ov{f_i}(\ov{q_1},\ldots,\ov{q_n})\equiv 0\mod \ov{h} \quad\hbox{dans }\FFp[X], \;i\in\lrbs.
$$ 
Pour un tel $p$, on prend un diviseur \mon \ird
$h_p(X)$ de $\ov{h}(X)$ dans~$\FFp[X]$ et l'on considère le corps fini 
$\gF=\aqo{\FFp[X]}{h_p(X)}$ avec $\alpha_p$ la classe de~$X$. Alors, 
$\big(q_1(\alpha_p),\ldots,q_n(\alpha_p)\big)$ est un zéro de $(\lfs)$ dans~$\gF^n$.
\end{proof}

On a le corolaire \imd suivant, avec des \itfsz. 

%:     corollary{corthNSTsurZ}
\begin{corollary}\label{corthNSTsurZ}
\emph{(\nst sur $\ZZ$, \nst formel, 2)} 
\\
On écrit $\ZZuX=\ZZXn$.
Pour deux \itfs $\fa$, $\fb$ de $\ZZuX$ \propeq
\begin{enumerate}
\item $\rD_{\ZZuX}(\fa)\subseteq\rD_{\ZZuX}(\fb)$.
\item  
 $\rD_{\gK}\big(\varphi(\fa)\big)\subseteq\rD_{\gK}\big(\varphi(\fb)\big)$ pour tout
\cdi $\gK$ et tout
\homo $\varphi:\ZZuX\to\gK$.
%
%\item $\fa$ et $\fb$ ont les mêmes zéros sur n'importe quel anneau réduit.
%
\item Même chose en se limitant aux 
 extensions \agqs de $\QQ$ et aux corps finis.
\item Même chose en se limitant aux 
 corps finis.
\end{enumerate}
\end{corollary}
%--------- fin corolaire ---------------------------------------------- 

%:    subsec{Exemples d'applications}
\mni\textbf{Un exemple d'application}

\mni Nous considérons %tout d'abord 
le résultat suivant déjà
démontré dans le lemme \ref{lemGaussJoyal}:
\emph{Un \elt $f$ de $\AuX$ est \iv \ssi $f(\uze)$ est \iv et $f-f(\uze)$
est nilpotent. Autrement dit $\AuX\eti = \Ati+\DA(0)[\uX]$.}
\\
On peut supposer que $fg=1$ avec $f=1+Xf_1$ et $g=1+Xg_1$. On considère les
\coes de $f_1$ et $g_1$ comme des \idtrsz. On est ramené à montrer le résultat suivant.
\\
\emph{Une \egt $\;f_1+g_1+Xf_1g_1=0\;\;(*)\;$ implique que les \coes de~$f_1$ sont nilpotents.}
\\
Or lorsque les \idtrs sont évaluées dans un corps, les \coes de~$f_1$
s'annulent aux zéros du \syp en les \idtrs donné par l'\egtz~$(*)$.
On conclut par le \nst formel.

 Si l'on compare à la \dem donnée pour le point \emph{4}
du lemme~\ref{lemGaussJoyal}, on peut constater que celle donnée ici
est à la fois plus simple (pas besoin de trouver un calcul un peu subtil)
et plus savante (utilisation du \nst formel).

\emph{Note.} D'autre exemples sont donnés dans  
%: sinotenglish
\sinotenglish{
l'exercice~\ref{exoPgcdNst} et} 
le \pb \ref{exoChasserIdeauxPremiers1}.
\eoe

%%%%%%%%%%%%%%%%%%%%%%%%%%%%%%%%%%%%%%%%%%%%%%%%%%%%%%%%%%%%%%%%%%%%%%%%%%%
\section{La méthode de Newton en \algz }  \label{secNewton}

Soit $\gk$ un anneau et $f_1$, $\ldots$, $f_s\in\kuX= \gk[\Xn]$.
La  \emph{matrice jacobienne} du \sys est la matrice%
\index{matrice!jacobienne}\imN
%-----------------begin $$----------------
$$ \JJ_{\Xn}(f_1,\ldots,f_s) =
%\left[
\Big( \frac{\partial f_i}{\partial X_j}
\Big)_{i\in \lrbs,j\in \lrbn }
%\right]
\in \kuX^{s\times n}.
$$
%-----------------end $$------------------
Celle-ci est encore notée $\JJ_{\uX}(\uf)$ ou $\JJ(\uf)$.
On la visualise comme ceci:
$$
\bordercmatrix [\lbrack\rbrack]{
    & X_1                     & X_2                     &\cdots  & X_n \cr
f_1 & \Dpp {f_1}{X_1} &\Dpp {f_1}{X_2}  &\cdots  &\Dpp {f_1}{X_n} \cr
f_2 & \Dpp {f_2}{X_1} &\Dpp {f_2}{X_2}  &\cdots  &\Dpp {f_2}{X_n} \cr
f_i & \vdots                  &                         &        & \vdots              \cr
 & \vdots                  &                         &        & \vdots              \cr
f_s & \Dpp {f_s}{X_1} &\Dpp {f_s}{X_2}  &\cdots  &\Dpp {f_s}{X_n} \cr
}
.$$
Si $s=n$, on note $\J_{\uX}(\uf)$
ou $\J_{\Xn}(f_1,\ldots,f_n)$ ou $\J(\uf)$ 
le \ixc{jacobien}{d'un \sypz} du \sys $(\uf)$, 
\cad le \deter de la matrice jacobienne.

En analyse la méthode de Newton pour approcher un zéro d'une fonction
différentiable $f:\RR\to\RR$ est la suivante. On part d'un point $x_0$ qui est
\gui{proche  d'une racine}, en lequel la dérivée est \gui{loin de 0} et l'on
construit une suite $(x_m)_{m\in\NN}$ par \recu en posant
$$
x_{m+1}=x_m-\frac{f(x_m)}{f'(x_m)}.
$$
La méthode se généralise pour un système de $p$ équations à $p$
inconnues. Une solution d'un tel système est un zéro d'une
fonction~\hbox{$f:\RR^p\to\RR^p$}. On applique \gui{la même formule} que ci-dessus:
$$
x_{m+1}=x_m-f'(x_m)^{-1}\cdot f(x_m).
$$
où $f'(x)$ est la  \dile (la matrice jacobienne)
de $f$ au point $x\in\RR^p$, qui
doit être \iv dans un voisinage de $x_0$.

Cette méthode, et d'autres méthodes du calcul infinitésimal,
s'appliquent dans certains cas \egmt en algèbre, en remplaçant les
infiniment petits leibniziens par des \elts nilpotents.

Si par exemple $\gA$ est une \QQlg et $x\in\gA$ est nilpotent, la série formelle

\snic{1+x+x^2/2+x^3/6+\ldots }

%\sni
qui définit $\exp(x)$ n'a qu'un nombre fini de termes non nuls dans
$\gA$ et définit donc un \elt $1+y$ avec $y$ nilpotent. Comme l'\egt

\snic{\exp(x+x')=\exp(x)\exp(x'),}

%\sni
 parce qu'elle a lieu en analyse,
 valide la même formule au niveau des séries
formelles sur $\QQ$, on obtiendra lorsque $x$ et $x'$ sont nilpotents dans $\gA$
la même \egt dans $\gA$. De même la série formelle

\snic{y-y^2/2+y^3/3-\ldots ,}

%\sni
 qui définit $\log(1+y)$, n'a qu'un nombre fini de termes dans $\gA$
lorsque $y$ est nilpotent et permet de définir  $\log(1+y)$ comme un \elt
nilpotent de $\gA$. En outre, pour $x$ et $y$ nilpotents, on obtient les \egts
$\log\big(\exp(x)\big)=x$ et~$\exp\big(\log(1+y)\big)=1+y$ comme conséquences des \egts
correspondantes pour les séries formelles.

Dans le même style on obtient facilement, en utilisant la série formelle
inverse de $1-x$, le résultat suivant.
%--- Lemma{lemMatInvNil}-------
\begin{lemma}
\label{lemMatInvNil} \emph{(Lemme des \elts \rdt \ivsz)}
%-----------------begin enum------------------
\begin{enumerate}
\item Si $ef\equiv1$ modulo le nilradical, alors $e$ est \iv et

\snic{e^{-1}=f\,\sum_{k\geq 0}(1-ef)^k .}
\item Une matrice carrée $E\in\Mn(\gA)$ \iv modulo le nilradical
est \ivz. Supposons que $d\det(E)\equiv1$ modulo le nilradical.\\
Posons $F=d\wi{E}$ (où $\wi{E}$ est la matrice cotransposée de $E$).
 Alors,  $E^{-1}$  est dans le sous-anneau de $\Mn(\gA)$ engendré par les \coes du \polcar de $E$, $d$  et $E$.\\ 
 Plus \prmtz,
 la matrice $\In-EF = \big(1 - d\det(E)\big)\In$ est nilpotente~et

 \snic{E^{-1}=F\,\sum_{k\geq 0}\big(1 - d\det(E)\big)^k.}
\end{enumerate}
%-----------------end enum------------------
\end{lemma}
%--- end-lemma-----------------------------------------

Passons à la méthode de Newton. 

%:   Theorem{thNewtonLin}------------
\begin{theorem}
\label{thNewtonLin} \emph{(Méthode de Newton \linz\imN)}\\
Soient $\fN$ un \id d'un anneau $\gA$,
 $\uf=\tra{[\,f_1\;\cdots\;f_n\,]}$ un vecteur dont les \coos sont des \pols dans  $\AXn$, et $\ua=\tra{(a_1,\ldots ,a_n)}$
 dans~$\Ae n$
un \emph{zéro simple approché} du \sys au sens suivant.
%-----------------begin item------------------
\begin{enumerate}
\item [--] La matrice jacobienne $J({\ua})$ de $\uf$ au point $\ua$ est
\iv modu\-lo~$\fN$; soit $U\in\Mn(\gA)$ un tel inverse.
\item [--] Le vecteur $\uf(\ua)$ est nul modulo~$\fN$.
\end{enumerate}
%-----------------end item------------------
Considérons la suite $(\ua^{(m)})_{m\geq 1}$ dans $\gA^n$  définie par
l'itération de Newton linéaire:
$$
\ua^{(1)}=\ua, \,\,\,\,\ua^{(m+1)}=\ua^{(m)}- U \cdot \uf(\ua^{(m)})
.$$
%-----------------begin item------------------
\begin{enumerate}
\item [a.]
Cette suite satisfait les exigences $\fN$-adiques suivantes:
$$\ua^{(1)}\equiv\ua\,\,\,\mod\,\fN, \,\mathrm{\,et\,} \Tt m,\,\,\,
\ua^{(m+1)}\equiv\ua^{(m)} \,\, \mathrm{et}\,\,
\uf(\ua^{(m)})\equiv 0 \,\,\,\,\mod\,\fN^m.
$$
\item [b.] Cette suite est unique au sens suivant, si $\ub^{(m)}$ est une autre
suite vérifiant les  exigences du point  a, alors pour tout $m$,
$\ua^{(m)}\equiv\ub^{(m)}\,\,\,\mod\,\fN^m.$
\item [c.] Soit $\gA_1$ le sous-anneau engendré par les
\coes des $f_i$, par  ceux
de $U$ et par les \coos de $\ua$.
Dans cet anneau soit $\fN_1$ l'\id
engendré par les \coes de $\In-UJ(\ua)$ et les \coos de $\ua$.
Si les \gtrs de $\fN_1$ sont nilpotents, la suite converge en un nombre fini
d'étapes vers un vrai zéro du \sys $\uf$, et c'est l'unique zéro du \sys
congru à $\ua$ modulo~$\fN_1$.
\end{enumerate}
%-----------------end item------------------
\end{theorem}
%--- end-theorem-----------------------------------------

Sous les mêmes hypothèses, on a la méthode quadratique suivante.

%:   Theorem{thNewtonQuad}-----------
\begin{theorem}
\label{thNewtonQuad} \emph{(Méthode de Newton quadratique\imN)}\\
Définissons les suites $(\ua^{(m)})_{m\geq
0}$ dans $\gA^n$  et $(U^{(m)})_{m\geq 0}$ dans $\Mn(\gA)$,
par l'itération de
Newton quadratique suivante:
%--------------------begin array---------------
$$\begin{array}{lcl}
\ua^{(0)}=\ua,&  \quad \quad  & \ua^{(m+1)}=\ua^{(m)}- U^{(m)} \cdot
\uf(\ua^{(m)}),    \\[1mm]
U^{(0)}=U,&   & U^{(m+1)}=U^{(m)}\,\left(2\I_n-J(\ua^{(m+1)})U^{(m)}\right).
\end{array}$$
%---------------------end array--------------
Alors, on obtient pour tout $m$ les congruences suivantes:
%--------------------begin array---------------
$$
\begin{array}{lcll}
 \ua^{(m+1)}\equiv\ua^{(m)} & \,\,\mathrm{et}\,\, &
  U^{(m+1)}\equiv U^{(m)} &\,\mod \,\fN^{2^m}    \\
  \uf(\ua^{(m)})\equiv 0 &  \,\,\mathrm{et}\,\,  &
   U^{(m)}\,J(\ua^{(m)})\equiv \In  &\,\mod \,\fN^{2^m}.
\end{array}
$$
%---------------------end array--------------
\end{theorem}
%--- end-theorem-----------------------------------------

Nous laissons les démonstrations \alec (cf. \cite{ValuHallouin}) en remarquant
que l'itération concernant l'inverse de la matrice jacobienne peut être
justifiée par la méthode de Newton linéaire ou par le calcul suivant dans un
anneau non \ncrt commutatif:
$$
(1-ab)^2=1-ab'\quad \mathrm{avec}\quad b'=b(2-ab).
$$

%-% entrenous
\entrenous{Après le \tho 10.3 sur la méthode de Newton quadratique

 \rem Le résultat précédent très simple est en quelque sorte trop
simple, par son cadre trop restrictif.
Il mérite d'être \gne en \gaqz, par exemple à d'autres fonctions
que les fonctions \polsz, et avec un autre radical que le nilradical.
\eoe

Le lemme de Hensel
multidimen\-sion\-nel et le Main Theorem de Zariski couvent là-dessous?
Nous reviendrons là dessus? 

La remarque serait utile à condition de la développer.
%-% Fin entrenous
}

%--- Corollary{corIdmNewton}---------
\begin{corollary}
\label{corIdmNewton} \emph{(Lemme des \idms résiduels)}
%-----------------begin enum------------------
\begin{enumerate}
\item Pour tout anneau commutatif $\gA$:
\begin{enumerate}
\item deux \idms égaux modulo $\DA(0)=\sqrt{\gen{0}}$
sont égaux,
\item tout \idm $e$ modulo un idéal $\fN$  se
relève de manière unique en un \idm $e'$ modulo $\fN^2$. L'itération de
Newton quadratique est donnée par $e\mapsto 3e^2-2e^3$.
\end{enumerate}
\item De même toute matrice $E\in\Mn(\gA)$ \idme modulo~$\fN$  se
relève en une matrice $F$ \idme modulo $\fN^2$. Le \gui{relèvement} $F$ est
unique si l'on exige que $F\in\gA[E]$. L'itération de Newton quadratique est
donnée par $E\mapsto 3E^2-2E^3$.
\end{enumerate}
%-----------------end enum------------------
\end{corollary}
%--- end-corollary------------------------------------
%-----------------begin proof------------------
\begin{proof}
\emph{1a.} Laissé \alecz. Une version plus forte est démontrée dans le lemme~\ref{lemIdmIsoles}. 
\\
\emph{1b.} Considérer le \hbox{\polz~$T^2-T$}, et noter
que  $2e-1$ est \iv modulo~$\fN$ puisque $(2e-1)^2=1$ modulo $\fN$.
%que
%$2e-1$ est \iv si $e$ est un \idm résiduel puisque $(2e-1)^2=1$
%modulo $\fN$.
\\
\emph{2.} On applique le point \emph{1} avec l'anneau commutatif
$\gA[E]\subseteq\End(\gA^n)$.
\end{proof}
%-----------------end proof------------------

%%%%%%%%%%%%%%%%%%%%%%%%%%%%%%%%%%%%%%%%%
%:section: Exercices
\Exercices{

%--- Exercise{exoLagrange}-------------
\begin{exercise}
\label{exoLagrange} (Interpolation de Lagrange)~
{\rm  Soit $\gA$ un anneau commutatif.\index{interpolation de Lagrange}%
\index{Lagrange!interpolation de ---}
\begin{enumerate}\itemsep=0pt
\item Soient $f$, $g\in\AX$ et $a_1$, $\ldots$, $a_k$ des \elts de $\gA$
tels que $a_i-a_j\in\Reg\gA$ pour $i\neq j$.
\begin{enumerate}
\item Si les $a_i$ sont des zéros  de $f$, 
$f$ est multiple de $(X-a_1)\cdots(X-a_k)$.
\item Si $f(a_i)=g(a_i)$ pour $i\in\lrbk$ et si $\deg(f-g)<k$, alors $f=g$.
\end{enumerate}  

\item Si $\gA$ est intègre et infini, l'\elt $f$ de $\AX$ est \care par la fonction \polle qu'il définit sur $\gA$. 

\item \emph{(\Pol d'interpolation de Lagrange)}
Soient $(\xzn)$ dans $\gA$ tels que les $x_i-x_j\in\Ati$
(pour $i\neq j$).
Alors, pour  $(y_0, \ldots, y_n)$ dans $\gA$ il existe exactement un \polz~$f$ de degré $\leq n$ tel que pour chaque $j\in\lrb{0..n}$ on ait $f(x_j)=y_j$.
\\
Plus \prmtz, le \polz~$f_i$ de degré $\leq n$ tel que $f_i(x_i)=1$
et $f_i(x_j)=0$ pour $j\neq i$ est égal à
$$
\preskip.0em \postskip.4em 
f_i= \frac{\prod\nolimits_{j\in\lrb{0..n},j\neq i}(X-x_j)}{\prod\nolimits_{j\in\lrb{0..n},j\neq i}(x_i-x_j)}, 
$$
et le \pol d'interpolation $f$ ci-dessus est égal à $\som_{i\in\lrbzn}y_if_i$.  

\item Avec les mêmes hypothèses, en posant $h=(X-x_0)\cdots(X-x_n)$,
on obtient un \iso  d'\Algsz: $\aqo\AX h \to \gA^{n+1},\;\ov g\mapsto\big(g(x_0),\ldots,g(x_n)\big)$.
\item Interprétez les résultats précédents avec l'\alg \lin (matrice et \deter de Vandermonde) et avec le \tho des restes chinois (utilisez les \ids deux à deux \com $\gen{X-x_i}$).
\end{enumerate}
 
}
\end{exercise}
%--- end -exercise-----------------------------------------

%--- Exercise{exoGensIdealEnsFini}-------------
\begin{exercise}
\label {exoGensIdealEnsFini} 
       {(Générateurs de l'idéal d'un ensemble fini)} {\rm Voir aussi l'exercice~\ref{exoGensPolIdeal}.\\ 
Soit $\gK$ un corps discret et $V \subset \gK^n$ un ensemble fini. On va
montrer que l'\idz~$\fa(V) = \sotq{ f \in \Kux}{ \forall\ w \in V,\ f(w) = 0}$ 
est engendré par $n$ \elts (notez que cette borne ne dépend pas de $\#V$
et que le résultat est clair pour $n=1$).
On note~$\pi_n : \gK^n \to \gK$ la $n$-ième \prn et pour chaque~$\xi
\in \pi_n(V)$, 

\snic{V_\xi = \sotq {(\xi_1, \ldots, \xi_{n-1}) \in \gK^{n-1}} 
         {(\xi_1, \ldots, \xi_{n-1}, \xi) \in V}.}

\begin {enumerate}
\item
Soit $U \subset \gK$ une partie finie et pour chaque $\xi \in U$, un \pol

\snic{Q_\xi
\in \gK[x_1, \ldots, x_{n-1}].}

%\sni
Expliciter un
\polz~$Q \in \Kux$ vérifiant $Q(x_1, \ldots, x_{n-1},\xi) = Q_\xi$ 
pour tout~$\xi \in U$.
\item
Soit $V \subset \gK^n$ une partie telle que $\pi_n(V)$ soit finie. On suppose
que pour chaque $\xi \in \pi_n(V)$, l'idéal $\fa(V_\xi)$ est engendré par $m$
\polsz. Montrer que~$\fa(V)$ est engendré par $m+1$ \polsz.
Conclure.
\end {enumerate}
}
\end {exercise}
%--- end -exercise-----------------------------------------

%--- Exercise{exothSymEl}-------------
\begin{exercise}
\label{exothSymEl}
(\Demo détaillée du \thref{thSymEl})\\
{\rm   On considère l'anneau  $\AXn=\AuX$ et l'on note $S_1$, $\ldots$, $S_n$ 
les fonctions \smqs \elrs de $\uX$. Tous les \pols considérés sont des \pols formels, car on ne suppose pas 
que $\gA$ est discret.  
On introduit un autre   
jeu d'\idtrsz, $(\underline{s}) = (s_1, \ldots, s_n)$, et  sur l'anneau 
$\gA[\underline{s}]$ on définit le
poids $\delta$ par~$\delta(s_i) = i$ (un \pol \fmt non nul a un poids formel bien défini). 
\\
On note $\varphi:\gA[\underline{s}]\to\AuX$
l'\homo d'\evn défini par $\varphi(s_i)=S_i$.
\\
On considère sur les \moms de $\AuX = \gA[X_1,
\ldots, X_n]$ l'ordre  {\tt deglex} 
pour lequel deux \moms sont d'abord comparés selon leur degré total, puis
ensuite selon l'ordre lexicographique avec $X_1 > \cdots > X_n$. Ceci fournit pour un $f \in \AuX$ (\fmt non nul) une notion de \emph{\mom \fmt dominant} que l'on note $\md(f)$. Cet \gui{ordre monomial} est clairement isomorphe à $(\NN,\leq)$.

\smallskip 
\emph{0.} Vérifier que tout \pol  \smq (i.e. invariant par l'action de $\Sn$) de~$\AuX$
est égal à un \pol \fmt \smqz, i.e. invariant par l'action de~$\Sn$
en tant que \pol formel.

\emph{1. (Injectivité de $\varphi$)}\,\\
Soit $\alpha = (\alpha_1, \ldots, \alpha_n)$ un exposant décroissant
($\alpha_1 \ge \cdots \ge \alpha_n)$. \\
On pose $\beta_i=\alpha_i-\alpha_{i+1}$
($i\in\lrb{1..n-1}$). Montrer que

\snic{\md(S_1^{\beta_1} S_2^{\beta_2} \cdots
S_{n-1}^{\beta_{n-1}} S_n^{\alpha_n}) =
X_1^{\alpha_1} X_2^{\alpha_2} \cdots X_n^{\alpha_n}.}

%\sni
En déduire que $\varphi$ est injectif.

\emph{2. (Fin de la \dem des points 1 et 2 du \tho \ref{thSymEl})}\,
Soit $f \in \AuX$ un \pol \fmt \smqz, \fmt non nul, et
$\uX^\alpha = \md(f)$. 

\begin{itemize}
\item Montrer que $\alpha$ est décroissant.
  En déduire un \algo pour écrire tout \pol \smq
de $\AuX$ comme un \pol en $(S_1, \ldots, S_n)$ à \coes
dans $\gA$, i.e. dans l'image de $\varphi$. La terminaison de l'\algo peut 
être prouvée par \recu sur l'ordre monomial, isomorphe à $\NN$.
\item \`A titre d'exemple, écrire le symétrisé du \mom $X_1^4X_2^2X_3$
dans $\gA[X_1,\ldots,X_4]$ comme \pol en les $S_i$.

\end{itemize}

\goodbreak
\emph{3.} 
\emph{(\Demo du point 3 du \thoz)} 

\begin{itemize}
\item
Soit $g(T) \in \gB[T]$ un \polu de degré $n \ge 1$. 
Montrer que
$\gB[T]$ est un $\gB[g]$-module libre de base $(1, T, \ldots, T^{n-1})$.\\
En déduire que $\gA[S_1,\ldots,S_{n-1}][X_n]$ est un module libre \hbox{sur
$\gA[S_1,\ldots,S_{n-1}
][S_n]$},   de base $(1, X_n, \ldots, X_n^{n-1})$.
\item
On note $\uS' = (S'_1, \ldots, S'_{n-1})$ les fonctions \smqs \elrs
des variables~$(X_1, \ldots,\alb X_{n-1})$. 
Montrer que $\gA[\uS', X_n] = \gA[S_1, \ldots, S_{n-1}, X_n]$.
\item
Déduire des deux points précédents que $\gA[\uS',X_n]$ est
un $\gA[\uS]$-module libre de base $(1, X_n, \ldots, X_n^{n-1})$.
\item Conclure par \recu sur $n$ que la famille

\snic{\sotq{X^\alpha}{\alpha=(\aln)\in\NN^n,\,\forall k\in\lrbn, \,\alpha_k<k}}

%\sni
forme une base de $\AuX$ sur $\gA[\uS]$.
\end{itemize}

\emph{4.} 
\emph{(Une autre \dem  du point 3 du \thoz, et même plus, après avoir lu la section \ref{sec0adu})}\,\,  
Montrer que $\AuX$ est canoniquement isomorphe à l'\adu du \polz~$t^n+\sum_{k=1}^{n}
(-1)^ks_kt^{n-k}$ sur l'anneau~$\gA[s_1,\ldots,s_n]$.
}
\end{exercise}

%--- Exercise{exoPolSym1}-------------
\begin{exercise}
\label{exoPolSym1}
{\rm
On note $S_1$, \ldots, $S_n \in \AuX = \gA[X_1, \ldots, X_n]$
les $n$ fonctions \smqs \elrsz.
\begin{enumerate}\itemsep=0pt
\item
Pour $n = 3$, vérifier que $X_1^3 + X_2^3 + X_3^3 = S_1^3 - 3S_1S_2 + 3S_3$.
En déduire que 
pour tout $n$, $\sum_{i = 1}^n X_i^3 = S_1^3 - 3S_1S_2 + 3S_3$.
\item
En utilisant une méthode analogue à la question précédente, exprimer
les \pols $\sum_{i \ne j} X_i^2 X_j$, $\sum_{i \ne j} X_i^3 X_j$, $\sum_{i
< j} X_i^2 X_j^2$ à l'aide des fonctions \smqs \elrs.

\item
\'Enoncer un résultat \gnlz.
\end {enumerate}
}
\end{exercise}
%--- end -exercise-----------------------------------------

%--- Exercise{exoNewtonSum}-------------
\begin{exercise}\label{exoNewtonSum}
 {(Les sommes de Newton et les fonctions \smqs complètes)}
{\rm
On note $S_i \in \AuX = \gA[X_1, \ldots, X_n]$ les fonctions
\smqs \elrs en convenant de prendre $S_i = 0$ pour $i > n$ et  $S_0 = 1$.
\\
Pour $r\ge 1$, on définit les \ix{sommes de Newton} $\SNw{r} = X_1^r + \cdots +
X_n^r$. On travaille dans l'anneau des séries formelles~$\gA[\uX][[t]]$ 
et l'on introduit les séries\isN

\snic {
P(t) = \sum_{r \ge 1} \SNw{r}\, t^r  \quad \hbox{ et } \quad
E(t) = \sum_{r \ge 0} S_r \, t^r.
}
\begin {enumerate}\itemsep=0pt
\item
Vérifier l'\egt
%
%\snic {
$P(t) = \sum_{i = 1}^n {X_i \over 1 - X_i t}$.
%}
%
\item
Si $u \in \gB[[t]]$ est inversible, on introduit sa dérivée logarithmique

\snic{D_{\rm log}(u) = u' u^{-1}.}

%\sni
On obtient ainsi un morphisme de groupes $D_{\rm log} :
(\gBtst, \times) \to (\gB[[t]], +)$.
\item
En utilisant la dérivée logarithmique, montrer la
{\it relation de Newton}:

\snic {
P(-t)= {E'(t) \over E(t)}, \;\;\hbox { ou encore }  P(-t)E(t) = E'(t).
}

\item
Pour $d \ge 1$, en déduire la formule de Newton:

\snic {
\sum_{r = 1}^d (-1)^{r-1} \SNw{r}\, S_{d-r} = d\,S_d.
}
\end {enumerate}
Pour $r \ge 0$, on définit la \ix{fonction symé\-trique complète de degré~$r$} par

\snic{H_r = \som_{|\alpha| = r} \uX^\alpha.}

%\sni
Ainsi $H_1 = S_1$, $H_2 = \sum_{i \le
j} X_iX_j$, $H_3 = \sum_{i \le j \le k} X_iX_jX_k$. On définit la série:

\snic {
H(t) = \sum_{r \ge 1} H_r\, t^r.
}
\begin {enumerate} \setcounter {enumi}{4}
\item
Montrer l'\egt
%
%\snic {
$H(t) = \sum_{i = 1}^n {1 \over 1 - X_i \,t}$.
%}
%%
\item
En déduire l'\egt $H(t)\,E(-t) = 1$, puis  pour $d \in \lrb {1..n}$:

\snic {
\sum_{r = 0}^{d} (-1)^r S_r\, H_{d-r} = 0, \quad
H_d \in \gA[S_1, \ldots, S_d], \quad
S_d \in \gA[H_1, \ldots, H_d].
}

\item
On considère l'\homo  $\varphi : \gA[S_1,\ldots, S_n]\to\gA[S_1,\ldots, S_n]$
défini \hbox{par $\varphi(S_i) = H_i$}. Montrer que $\varphi(H_d) = S_d$
pour $d \in \lrb {1..n}$. 
Ainsi: 
\begin{itemize}
\item $\varphi \circ \varphi =
\I_{\gA[\uS]}$,
\item  $ H_1$, \ldots, $H_n$ sont \agqt
indépendants sur $\gA$,
\item  $\gA[\uS] = \gA[\uH]$, et 
l'expression de~$S_d$ en fonction de $H_1$, $\ldots$, $H_d$
est la même que celle de~$H_d$ en fonction de $S_1$,~$\ldots$,~$S_d$.
\end{itemize}
\end {enumerate}
}
\end {exercise}

%--- end -exercise-----------------------------------------

%--- Exercise{exolemArtin}-------------
\begin{exercise}
\label{exolemArtin} (Formes \eqves du lemme de \DKMz)\\
{\rm   Les affirmations suivantes sont \eqves
(chacune des affirmations est \uvlez, i.e., valable pour tous \pols et
tous anneaux commutatifs):
%-----------------begin item------------------

\emph{1.} $\rc(f)=\gen{1}\; \Longrightarrow\;  \rc(g)=\rc(fg)$.

\emph{2.} $\exists p\in\NN \; \; \rc(f)^{p}\rc(g)\subseteq \rc(fg)$.

\emph{3.} \emph{(Dedekind-Mertens, forme affaiblie)} 
$\;\;\;\;\;\exists p\in\NN \; \;
\rc(f)^{p+1}\rc(g)=\rc(f)^p\rc(fg)$.

\emph{4.} $\Ann\big(\rc(f)\big)=0\; \; \Longrightarrow\; \;
\Ann\big(\rc(fg)\big)=\Ann\big(\rc(g)\big)$.

\emph{5.} \emph{(McCoy)} $\;\;\;\;\;\big(\Ann(\rc(f)\big)=0, \; fg=0)\; \Longrightarrow\;
g=0$.

\emph{6.} $(\rc(f)=\gen{1}, \; fg=0)\; \Longrightarrow\;  g=0$.
}
\end{exercise}
%--- end -exercise-----------------------------------------

%-% ENTRE NOUS
\entrenous{ si possible un exo sur
%--- Exercise{exo2lemArtin}-------------
%\begin{exercise}
%\label{exo2lemArtin} 
\DKMz, version Jouanolou
%{\rm
%
%}
%\end{exercise}
%--- end -exercise-----------------------------------------
}
%-% Fin ENTRENOUS

%--- Exercise{exoMcCoy}-------------
\begin{exercise}
\label{exoMcCoy}
{\rm
Soit $\fc = \rc(f)$ le contenu de $f \in \gA[T]$. Le lemme de \DKM
 donne:
$\Ann_\gA(\fc)[T] \subseteq \Ann_{\gA[T]}(f) \subseteq \DA(\Ann_\gA\big(\fc)\big)[T]$.
Donner un exemple pour lequel il n'y a pas \egtz.
}
\end {exercise}
%--- end -exercise-----------------------------------------

%--- Exercise{exothKro}-------------
\begin{exercise}
\label{exothKro}
{\rm Déduire le \tho de \KRO \ref{thKro} du lemme de \DKMz.
 }
\end{exercise}
%--- end -exercise-----------------------------------------

%--- Exercise{exoModCauBase}----------
\begin{exercise}
\label{exoModCauBase} (Modules de Cauchy) 
{\rm  On peut donner une explication très précise pour le fait que
l'idéal $\cJ(f)$  (\dfn \ref{definotaAdu}) est égal à l'idéal engendré par les modules de
Cauchy. Cela fonctionne avec une belle formule.
Introduisons une nouvelle variable $T$.
Démontrer les résultats suivants.
%-----------------begin enum------------------
\begin{enumerate}\itemsep=0pt
\item Dans $\gA[X_1,\ldots ,X_n,T]=\gA[\uX,T]$,
on a
%---  equation eqModCau --------
\begin{equation}\label{eqModCau}
\arraycolsep2pt\begin{array}{rcl}
f(T)&   =&f_1(X_1)+(T-X_1)f_2(X_1,X_2)+
\\[1mm]
&&\quad (T-X_1)(T-X_2)f_3(X_1,X_2,X_3)+\cdots+
\\[1mm]
&   &\quad\quad
(T-X_1)\cdots (T-X_{n-1})f_n(X_1,\ldots ,X_n)+\\[1mm]
&&\quad\quad\quad (T-X_1)\cdots (T-X_{n})
\end{array}
\end{equation}
%---------------------end equation--------------
\item Dans le sous-$\gA[\uX]$-module de  $\gA[\uX,T]$ 
formé par les \pols de degré $\leq n$ en $T$,
le \polz~$f(T)-(T-X_1)\cdots (T-X_{n})$ possède deux expressions
différentes:
%-----------------begin item------------------
\begin{itemize}
\item D'une part, sur la base $(1,T,T^2,\ldots ,T^n)$,
il a pour \coos 

\centerline{$\big((-1)^n(s_n-S_n),\ldots ,(s_2-S_2),
-(s_1-S_1),0\big)$.}
\item D'autre part, sur la base 

\centerline{$\big(1,(T-X_1),(T-X_1)(T-X_2),\ldots ,(T-
X_1)\cdots (T-X_{n})\big)$,}

il a pour \coos $(f_1,f_2,\ldots ,f_{n},0)$.
\end{itemize}
%-----------------end item------------------
%\item
En conséquence sur l'anneau $\AXn$,
chacun des deux vecteurs

\snic{\big((-1)^n(s_n-S_n),\ldots ,(s_2-S_2),-(s_1-S_1)\big)
\quad \hbox{ et }\quad (f_1,\ldots ,f_{n-1},f_n)\qquad}

%\sni
s'exprime en fonction de l'autre
au moyen d'une matrice unipotente (triangulaire avec des 1 sur la
diagonale).
\end{enumerate}
%-----------------end enum------------------

}
\end{exercise}
%--- end-exercise-----------------------------------------

%--- Exercise{exoPrimePowerRoot}-------------
\begin{exercise}\label{exoPrimePowerRoot} {(Le \polz~$X^p - a$)}
{\rm
Soit $a \in \Ati$ et $p$ un nombre premier. On suppose que le \polz~$X^p - a$
possède dans $\gA[X]$ un diviseur \mon non trivial. Montrer que $a$ est une
puissance $p$-ième dans $\gA$.
}
\end {exercise}
%--- end -exercise-----------------------------------------

%--- Exercise{exoPrincipeIdentitesAlgebriques}-------------
\begin{exercise}
\label{exoPrincipeIdentitesAlgebriques} (Avec le principe
de prolongement des \idasz)\\
{\rm
Notons $S_n(\gA)$ le sous-module de $\Mn(\gA)$ constitué des matrices
\smqsz.  \\
Pour
$A \in S_n(\gA)$, notons $\varphi_{\!A}$ l'\endo de $S_n(\gA)$
défini par $S \mapsto \tra {A} S A$.  Calculer $\det(\varphi_{\!A})$ en
fonction de $\det(A)$. Montrer que $\rC{\varphi_{\!A}}$ ne dépend que
de~$\rC {\!A}$.
}
\end {exercise}
%--- end -exercise-----------------------------------------

%--- Exercise{exoFreeFracTransfert}-------------
\begin{exercise}
\label{exoFreeFracTransfert}
{\rm
Soit $\gB\supseteq\gA$ une \Alg  intègre, libre de rang $n$,
 $\gK = \Frac(\gA)$ et~$\gL = \Frac(\gB)$. Montrer que toute base de $\gB/\gA$
est une base de~$\gL/\gK$.
}
\end{exercise}
%--- end -exercise-----------------------------------------

%--- Exercise{exoResultant}-------------
\begin{exercise}
\label{exoResultant}
{\rm
 Soient $f\in\AX$, $g\in\AY$, $h\in\gA[X,Y]$. Démontrer que

\snic{\Res_Y\big(g,\Res_X(f,h)\big)=\Res_X\big(f,\Res_Y(g,h)\big).}

\hum{utile dans le cours?}
}
\end{exercise}
%--- end -exercise-----------------------------------------

%-% ENTRE NOUS
\entrenous{ si possible\\
%--- Exercise{exoTschirnhaus}-------------
%\begin{exercise}
%\label{exoTschirnhaus}
\textbf{Exercice}  Différentes expressions pour la transformation de Tschirnhaus.

%\end{exercise}
%--- end -exercise-----------------------------------------

%--- Exercise{exo2Tschirnhaus}-------------
%\begin{exercise}
%\label{exo2Tschirnhaus}
\textbf{Exercice}  Généralisations de la transformation de Tschirnhaus.
%\end{exercise}
%%--- end -exercise-----------------------------------------
}
%-% Fin ENTRENOUS

%--- Exercise{exoSommesNewton}-------------
\begin{exercise}\label{exoSommesNewton}
{(\isN Sommes de Newton et $\Tr(A^k)$)}
{\rm
~ Soit une matrice $A \in \Mn(\gB)$. On pose $\rC{A}(X)=X^n+\sum_{j=1}^n (-1)^j s_j X^{n-j}$, $s_0=1$ et  $p_k=\Tr(A^k)$.
 
\emph{1.}
Montrer que les $p_k$ et $s_j$ sont reliés par les
formules de Newton pour les sommes des puissances $k$-ièmes
(exercice \ref{exoNewtonSum}): $\sum_{r = 1}^d (-1)^{r-1} p_{r} s_{d-r} = ds_d$
($d\in \lrbn$).
 
\emph{2.}
Si $\Tr(A^k) = 0$ pour $k \in \lrb {1..n}$, et si $n!$ est \ndz dans
$\gB$, alors $\rC{A}(X) = X^n$.
\\
NB: cet exercice peut être considéré comme une variation sur
le thème de la proposition~\ref{prop2tschir}.
}
\end {exercise}
%--- end -exercise-----------------------------------------

%--- Exercise{exoCorpsFiniEltPrimitif}-------------
\begin{exercise}
\label{exoCorpsFiniEltPrimitif}
{\rm
Soient $\gK \subseteq \gL$ deux corps finis, $q = \#\gK$
et $n = \dex{\gL:\gK}$. Le sous-anneau de $\gK$ engendré par $1$ est un corps $\FF_p$
où $p$ est un nombre premier, \hbox{et $q=p^r$} pour un entier $r>0$.
L'\ix{automorphisme de Frobenius} de (la $\gK$-extension)~$\gL$ est donné par
$\sigma : \gL \to \gL,\,\sigma(x) = x^q$.
 
\emph{1.}
Soit $R$ la réunion des racines dans $\gL$ des $X^{q^d} - X$ avec
$1 \le d < n$. \\
Montrer que $\#R < q^n$ et que pour $x \in \gL
\setminus R$, $\gL = \gK[x]$.

\emph{2.}
Ici $\gK = \FF_2$ et $\gL = \FF_2[X]/\gen {\Phi_5(X)} = \FF_2[x]$ où
$\Phi_5(X)$ est le \pol cyclotomique $X^4 + X^3 + X^2 + X + 1$.
Vérifier que $\gL$ est bien un corps; $x$ est un \elt primitif de $\gL$ sur
$\gK$ mais n'est pas un \gtr du groupe multiplicatif~$\gL\eti$.

\emph{3.}
Pour $x \in \gL\eti$, notons  $o(x)$ son ordre  dans le groupe multiplicatif
$\gL\eti$. Montrer que $\gL = \gK[x]$ \ssi l'ordre de $q$ dans le groupe
$\left(\aqo\ZZ{o(x)}\right)\eti$ est~$n$.
}
\end {exercise}
%--- end -exercise-----------------------------------------

%--- Exercise{exoRacinesUniteCyclique}-------------
\begin{exercise}
\label{exoRacinesUniteCyclique}
{\rm  Le but de cet exercice est de montrer que dans un \cdi le groupe des racines $n$-ièmes de l'unité est cyclique. En conséquence le
groupe multiplicatif d'un corps fini est cyclique. On montre un résultat à peine plus \gnlz.
\\
 Montrer que dans un anneau commutatif
non trivial $\gA$, si des \elts $(x_i)_{i\in\lrbn}$
forment un groupe $G$ pour la multiplication, et si $x_i-x_j$ est \ndz pour tout couple $i,j$ ($i\neq j$), alors $G$ est cyclique.
\\
Suggestion: d'après le \tho de structure
des groupes abéliens finis, un groupe abélien fini, noté additivement, dans lequel toute \eqn $dx=0$ admet au plus~$d$ solutions est cyclique.
Utilisez aussi l'exercice~\ref{exoLagrange}.

}
\end{exercise}
%--- end -exercise-----------------------------------------

%--- Exercise{exoCorpsFinis}-------------
\begin{exercise}
\label{exoCorpsFinis} (Structure des corps finis, automorphisme de Frobenius){\rm
  
\emph{1.}
Démontrer que deux corps finis qui ont le même ordre sont isomorphes.

\emph{2.} Si $\gF\supseteq\FF_p$ est un corps fini d'ordre $p^r$, montrer que
$\tau:x\mapsto x^p$ définit un \auto de $\gF$. On l'appelle l'\emph{automorphisme de Frobenius}. Montrer le groupe des \autos de $\gF$ est un groupe cyclique
d'ordre $r$ engendré par $\tau$.  

\emph{3.} Dans le cas précédent, $\gF$ est une extension galoisienne de
$\FF_p$. Préciser la correspondance galoisienne.
 
NB. On note souvent $\FF_q$ un corps fini d'ordre $q$, tout en sachant
qu'il s'agit d'une notation légèrement ambigu\"{e} si $q$ n'est pas premier. 

}
\end{exercise}
%--- end -exercise-----------------------------------------

%--- Exercise{exo2CorpsFinis}-------------
\begin{exercise}
\label{exo2CorpsFinis} (Clôture \agq de $\FF_p$) 
{\rm

\emph{1.}  Pour chaque entier $r>0$ construire un corps $\FF_{p^{r!}}$ d'ordre $p^{r!}$. En procédant par \recu on peut
avoir une inclusion $\imath_r:\FF_{p^{r!}}\hookrightarrow \FF_{p^{(r+1)!}}$.
 
\emph{2.} 
Construire un corps $\FF_{p^{\infty}}$ en prenant la réunion 
des $\FF_{p^{r!}}$ via les inclusions $\imath_r$.
Montrer que $\FF_{p^{\infty}}$ est  un corps \agqt clos qui contient une
copie (unique) de chaque corps fini de \cara $p$.
}
\end {exercise}
%--- end -exercise-----------------------------------------

%--- Exercise{exoPpcmPolsSeparables}-------------
\begin{exercise}\label{exoPpcmPolsSeparables} {(Ppcm de  \pols séparables)}
{\rm
  
\emph{1.} Soient $x$, $x'$, $y$, $y'\in \gB$. Montrer
que 
$\gen {x,x'} \gen {y,y'} \gen {x,y}^2
\subseteq \gen {xy, x'y + y'x}.$
\\ En déduire que le produit de deux \polus \spls
et \com  dans $\gA[T]$ est un \pol \splz.

\emph{2.} 
Si $\gA$ est un \cdiz, le 
ppcm de plusieurs \pols \spls est \splz.

}

\end {exercise}
%--- end -exercise-----------------------------------------

%--- Exercise{exolemSousLibre}-------------
\begin{exercise}
\label{exolemSousLibre} (Indice d'un sous-module \tf dans un module libre)
{\rm  

\emph{1.}
Soit $A\in\Ae{m\times n}$ et $E\subseteq \Ae m$ le sous-module image de $A$.
Montrer que $\cD_m(A)$ ne dépend que de $E$. 
On appelle cet \id l'\ixc{indice}{d'un sous-module dans un libre} \emph{de $E$
 dans $L=\Ae m$}, et on le note~\hbox{$\idg{L:E}_\gA$} \hbox{(ou $\idg{L:E}$)}. Remarquez que cet indice est nul dès
que $E$ ne s'approche pas suffisamment de $L$, par exemple si $n<m$.
\\
Vérifier que, dans le cas où $\gA=\ZZ$, on retrouve l'indice usuel du sous-groupe d'un groupe pour deux groupes abéliens libres de même rang. 
 
 \emph{2.} Si $E\subseteq F$ sont des sous-modules \tf de $L\simeq \Ae m$,
on a  $\idg{L:E}\subseteq \idg{L:F}$.
 
 \emph{3.} Si en outre $F$ est libre de rang $m$, on a la formule de transitivité

\snic {
\idg{L:E} = \idg{L:F}\,\idg{F:E}. 
}

\emph{4.} Si $\delta$ est un \elt \ndz de $\gA$, on a $\idg{\delta L:\delta E}=\idg{L:E}$. 
En déduire l'\egref{eqlemSousLibre} annoncée dans le lemme \ref{lemSousLibre}.
}
\end{exercise}
%--- end -exercise-----------------------------------------

%--- Exercise{exoNbGtrsInverse}-------------
\begin{exercise}
\label{exoNbGtrsInverse} (Précision sur le fait \ref{factdefiiv})
{\rm  Soient dans un anneau $\gA$ deux \idsz~$\fa$ et~$\fb$ tels que $\fa\,\fb=\gen{a}$
avec $a$ \ndzz. Montrer que si $\fa$ est engendré par $k$ \eltsz, on peut trouver dans $\fb$ un \sgr de $k$ \eltsz.
 
}
\end{exercise}
%--- end -exercise-----------------------------------------

%--- Exercise{exoDecompIdeal}-------------
\begin{exercise}
\label{exoDecompIdeal} (Décomposition d'un \id en produit d'\idemas \ivsz)
{\rm  On considère un anneau intègre non trivial $\gA$ \emph{à \dve explicite}\footnote{On dit qu'un anneau arbitraire est à \dve explicite lorsque l'on a un \algo qui, pour $a$ et $b\in\gA$, teste si $\exists x,\, a=bx$, et en cas de réponse positive, fournit un $x$ convenable.}.%
\index{divsibilité explicite!anneau à ---}\index{anneau!a divi@à divisbilité explicite}  
\vspace{-4pt}
\begin{enumerate}
\item Si $\fa$ est un \id \iv et si $\fb$ est un \itfz, il y a un test pour~$\fb \subseteq \fa$.
\end{enumerate}
\vspace{-4pt}
Soient $\fq_1$, $\ldots$,
$\fq_n$ des \idemas (dans le sens que les quotients $\gA\sur{\fq_k}$ dont des \cdis non triviaux), un \itf $\fb$ et un \elt $a$ \ndz de $\gA$ vérifiant
$a\gA=\fq_1 \cdots \fq_n \subseteq \fb$.
\vspace{-4pt}
\begin{enumerate} \setcounter{enumi}{1}
\item  Les $\fq_i$ sont \ivs et $\fb$ est le produit de certains des
$\fq_i$ (et par suite il est \ivz).  En outre, cette \dcn de  $\fb$ en produit d'\idemas
\tf est unique à l'ordre près des facteurs.
\end{enumerate}

}
\end{exercise}
%--- end -exercise-----------------------------------------
%--- Exercise{exoSymboleLegendre}-------------
\begin{exercise}\label{exoSymboleLegendre}
{(Symbole de Legendre)}\index{symbole de Legendre}\\
{\rm  
Soit $\gk$ un corps fini de cardinal $q$ impair; on définit
le \emph {symbole de Legendre}

\snic{\dsp
{ \legendr \bullet \gk} :
\gk\eti \lora \{\pm 1\}, \;
x \;\longmapsto \;\formule{
\phantom-1 \hbox{ si $x$ est un carré dans } \gk\eti ,
\\  
-1 \hbox{ sinon.}
}
}

%\sni
Montrer que ${. \legendre \gk}$ est un morphisme de groupes
et que ${x \legendre \gk} = x^{q-1 \over 2}$.\\
En particulier,~$-1$ est un carré dans $\gk\eti$
\ssi $q \equiv 1 \bmod 4$.\\
NB: si $p$ est un nombre premier impair et $x$ un entier étranger à $p$
on retrouve sous la forme ${{\ov x} \legendre \FFp}$
 le symbole  ${x \legendre p} $ défini par Legendre.}
\end {exercise}
%--- end -exercise-----------------------------------------

%--- Exercise{exoRabinovitchTrick}-------------
\begin{exercise}\label{exoRabinovitchTrick}
 {(L'astuce de Rabinovitch)}\\
{\rm  
Soit $\fa\subseteq \gA$ un \id et $x \in \gA$.
On considère l'\id suivant de $\AT$:

\snic {
\fb = \gen {\fa, 1-xT} = \fa[T] + \gen {1-xT}_{\gA[T]}.
}

%\sni
Montrer l'\eqvc $x \in \sqrt\fa \iff 1 \in \fb$.

}
\end {exercise}
%--- end -exercise-----------------------------------------

%--- Exercise{exoDunford}-------------
\begin{exercise}
\label{exoDunford} (Décomposition de Jordan-Chevalley-Dunford)\index{Dunford!décomposition de Jordan-Chevalley- ---}
\\
 {\rm  
Soit  $M\in\Mn(\gA)$. On suppose que le \polcar de $M$ 
divise une puissance  d'un \pol \spl $f$. 

\emph{1.} Montrer qu'il existe $D$, $N\in\Mn(\gA)$ tels que:
%-----------------begin item------------------
\begin{itemize}
\item  $D$ et $N$ sont des \pols en $M$ (à \coes dans $\gA$).
\item  $M=D+N$.
\item  $f(D)=0$.
\item  $N$ est nilpotente.
\end{itemize}
%-----------------end item------------------
%
\emph{2.} 
Montrer l'unicité de la décomposition ci-dessus. Y compris en 
affaiblissant la première contrainte: en demandant seulement $DN=ND$.

}
\end{exercise}
%--- end -exercise-----------------------------------------

%--- Exercise{exoDunfordBis}-------------
\begin{exercise}\label{exoDunfordBis}
{(\'Eléments \spbz{ment} entiers)}
\\{\rm  
Soit $\gA \subseteq \gB$. On  dira que $z \in \gB$ est
\emph{\spbz{ment} entier} sur $\gA$ si $z$ est racine d'un \polu\spl de $\AT$.
On cherche ici un exemple pour lequel la somme de deux
\elts \spbz{ment} entiers est  un \elt nilpotent, non nul et non \spbz{ment} entier.
\\
Soient $\gB = \gA[x] = \aqo{\gA[X]}{X^2 + bX + c}$. On suppose que  $\Delta = b^2 - 4c$ est une unité de $\gA$. 
Pour $a\in\gA$, calculer le \polcar de $ax$ sur $\gA$ et son \discriz.  
En déduire
un exemple comme annoncé lorsque $\DA(0)\neq 0$.

}

\end {exercise}
%--- end -exercise-----------------------------------------

%: sinotenglish
\sinotenglish{
%--- Exercise{exoAdj2Minors}-------------
\begin{exercise} \label{exoAdj2Minors}
 {(Mineurs d'ordre 2 de la matrice cotransposée)} \\
 {\rm Cet exercice explicite le point \emph{8} du lemme \ref{lemPrincipeIdentitesAlgebriques}. 
Etant donnés $i$, $i' \in \lrbn$ distincts, on note $\vep_{i,i'}$
la signature de la permutation $(I,i,i')$ où $I$ est la
suite $\lrbn$ privée de $(i,i')$. Ou encore si $(e_1, \ldots,
e_n)$ est la base canonique de $\Ae n$, $\vep_{i,i'}$ est défini par

\snic {
e_I \wedge e_i \wedge e_{i'} = \vep_{i,i'}\, e_1 \wedge \cdots \wedge e_n.
}

\snii
On a $\vep_{i,i'} = (-1)^{i+i'+1}$ si $i < i'$ et $\vep_{i,i'} = (-1)^{i+i'}$
sinon.

\snii
Soient $\wi {A}$ la cotransposée de $A \in \Mn(\gA)$, 
$i \ne i'$, $j \ne j'$, $I = \lrbn \setminus \{i,i'\}$,
$J = \lrbn \setminus \{j,j'\}$. Montrer que:
$$
\left| \matrix {
\wi{A}_{i,j} &\wi{A}_{i,j'}\cr \wi{A}_{i',j} &\wi{A}_{i',j'}\cr} \right| =
\vep_{i,i'}\vep_{j,j'} \det(A) \det(A_{J,I}).
$$
En déduire que si $\det(A) = 0$, alors $\wi {A}$ est de rang
inférieur ou égal à 1.
}

\end {exercise}

%%%%%%%%%%%%%%%%%%%%%%%%%%%%%%%%%%%%%%%%%%%%%%%%%%%%%%%%%%%%%%%%%%%%%%%%%%%
%--- Exercise{exoResInversible}-------------
\begin{exercise}
\label{exoResInversible}
{\rm  Soient $f$, $g\in\AX$ de degrés formels $p$, $q\geq 1$, \hbox{et
$f_p$, $g_q$} leurs \coes \fmt dominants. \\
Soient $F(X,Y)=Y^pf(\fraC X Y)$ et $G(X,Y)=Y^qg(\fraC X Y)$ les \pols homognéisés en leurs degrés formels.
Montrer que \propeq
\begin{enumerate}
\item  \label{i1exoResInversible} $\Res_{X}(f,p,g,q)\in\Ati$.
\item \label{i2exoResInversible} $f$ et $g$ sont \com dans $\AX$, et $f_p$ et $g_q$ sont \com dans $\gA$.
\item \label{i3exoResInversible} Il existe un $k\in\NN$ tel que $X^k$ et $Y^k\in \gen{F,G}$ dans $\gA[X,Y]$.
\item \label{i4exoResInversible} On a $\gen{X,Y}^{p+q-1}\subseteq  \gen{F,G}$ dans $\gA[X,Y]$.
\end{enumerate} 
}
\end{exercise}
%--- end -exercise-----------------------------------------

%%%%%%%%%%%%%%%%%%%%%%%%%%%%%%%%%%%%%%%%%%%%%%%%%%%%%%%%%%%%%%%%%%%%%%%%%%%

%--- Exercise{exoPgcdNst}-------------
\begin{exercise}\label{exoPgcdNst}
{(Un exemple d'application du \nst formel)} \\
{\rm
Il s'agit ici de \gnr au cas d'un anneau commutatif arbitraire un résultat
utile en théorie des \cdisz: \emph{si l'on divise un \pol $f(x)$ par le pgcd de~$f$ et $f'$, on obtient un \pol \splz.} \\
Pour un anneau arbitraire, on devra supposer que le pgcd de~$f$ et $f'$
existe en un sens fort. 

\emph{1.} Soit $\gK$ un corps discret, 
$x$ une \idtrz, $f\in\Kx$ un \pol non nul de degré $n\geq 0$, \hbox{$h=\pgcd(f,f')$} et
${f_1}=f/h$. \\
On suppose que $f$ se décompose en un produit de facteurs \lins dans
un corps discret contenant $\gK$. 
Montrer que $\Res_x({f_1},{f_1}')\in\gK\eti$, ou, ce qui revient au même, que $1\in\gen{f_1,f_1'}\subseteq \Kx$.
\\
Si en outre $\deg(f)=n$ et  $n!\in\gK\eti$, alors $f$ divise ${f_1}^{n}$.

\emph{2.} Démontrer les mêmes résultats sans faire d'hypothèse de \fcn concernant~$f$.

\emph{3.} Soit $\gk$ un anneau commutatif et $f\in\kx$ primitif de degré
formel $n\geq 2$. On suppose que l'\id $\gen{f,f'}$ est engendré par un \pol $h$ (\ncrt primitif). %
\begin{enumerate}
\item [\emph{a.}] Montrer qu'il existe des \pols  $u$, $v$, $f_2$, ${f_1}\in\kx$,  satisfaisant
les \egts 

\snic{u{f_1}+vf_2=1 \;\;\hbox{ et }\;\;\cmatrix{u&v\cr-f_2&{f_1}}\cmatrix{f\cr f'}=\cmatrix{h\cr 0}.}

\item [\emph{b.}] En utilisant le \nst formel, montrer que $1\in\gen{f_1,f_1'}\subseteq \kx$.
\item [\emph{c.}] Si en outre  $n!\in\gk\eti$, alors
peut-on montrer que $f$ divise ${f_1}^{n}$?
\end{enumerate}

\emph{4.} Question subsidiaire.
Donner une \dem directe du point \emph{3}
qui n'utilise pas le \nst formel.

}

\end {exercise}
%--- end -exercise-----------------------------------------
}
%: fin sinotenglish

%: problemes
%--- problem{exoDiscriminantsUtiles}-------------
\begin{problem}\label{exoDiscriminantsUtiles} {(Quelques résultants et discriminants utiles)}\\
{\rm
\emph {1.}
Montrer que $\disc(X^n + c) = (-1)^{n(n-1) \over 2} n^n c^{n-1}$.
Plus
généralement, montrer pour~$n \ge 2$ l'\egt

\snic {
\disc(X^n + bX + c) = (-1)^{n(n-1) \over 2}
\bigl(n^n c^{n-1} + (1 - n)^{n-1} b^n \bigr).
}

%\sni
\emph {2.}
Pour $n, m \in \NN^*$, en posant $d = \pgcd(n,m)$, $n_1=\frac n d$ et $m_1=\frac m d$ montrer l'\egt:

\snic {
\Res(X^n - a, X^m - b) = (-1)^n (b^{n_1} - a^{m_1})^d.
}

%\sni
Plus \gnlt

\snic {
\Res(\alpha X^n - a,n, \beta X^m - b,m) =
(-1)^n (\alpha^{m_1}b^{n_1} - \beta^{n_1}a^{m_1})^d.
}

%\sni
\emph {3.}
Notations comme au point \emph{2}, avec $1 \le m \le n-1$.  Alors:

\snic {
\disc(X^n + bX^m + c) = (-1)^{n(n-1) \over 2} c^{m-1}
\bigl(n^{n_1} c^{n_1 - m_1} - 
(n - m)^{n_1 - m_1} m^{m_1} (-b)^{n_1}\bigr)^d.
}

%\sni
\emph {4.}
Pour $n \in \NN^*$, on note $\Phi_n$ le polynôme cyclotomique de niveau $n$
(voir le \pbz~\ref{exoPolCyclotomique}).
Alors, pour  $p$ premier $\ge 3$\index{polynome@\pol!cyclotomique}

\snic {
\disc(\Phi_p) = (-1)^{p-1 \over 2} p^{p-2}.
}

%\sni
\emph {5.}
Soient $p$ premier et $k\ge 1$. Alors, $\Phi_{p^k}(X) = \Phi_p(X^{p^{k-1}})$
et:

\snic {
\disc(\Phi_{p^k}) = (-1)^{\varphi(p^k) \over 2} p^{(k(p-1)-1) p^{k-1}}
\qquad (p,k) \ne (2,1),
}

%\sni
avec pour $p \ne 2$, $(-1)^{\varphi(p^k)\over 2} = (-1)^{p-1\over 2}$. 
Pour $p=2$, on obtient $\disc(\Phi_4) = -4$ 
et~$\disc(\Phi_{2^k}) = 2^{(k-1) 2^{k-1}}$ pour $k \ge 3$.
Par ailleurs,  $\disc(\Phi_2) = 1$.

\emph {6.}
Soit $n \ge 1$ et $\zeta_n$ une racine primitive $n$-ième de l'unité.
\\
Si $n$ n'est pas la puissance d'un nombre premier, alors $\Phi_n(1) = 1$,
et $1 - \zeta_n$ est \iv dans $\ZZ[\zeta_n]$. 
\\
Si $n = p^k$ avec
$p$ premier, $k \ge 1$, alors $\Phi_n(1) = p$. Enfin $\Phi_1(1) = 0$.

\emph {7.}
Soit $\Delta_n = \disc(\Phi_n)$.  Pour $n$, $m$ premiers entre eux, on a la
formule de multiplicativité $\Delta_{nm} = \Delta_n^{\varphi(m)}
\Delta_m^{\varphi(n)}$ et l'\egt
$$
\preskip.4em \postskip.4em\ndsp 
\Delta_n = (-1)^{\varphi(n) \over 2}
\frac{n^{\varphi(n)}} { \prod_{p \mid n} p^{\varphi(n) \over p-1}} \qquad
\hbox {pour }n \ge 3. 
$$
}
\end {problem}
%--- end -problem-----------------------------------------

\vspace{-2em}
%--- problem{exoAnneauEuclidien}-------------
\begin{problem}
\label{exoAnneauEuclidien} (Anneaux euclidiens, l'exemple $\ZZ[i]$)
\\
{\rm  Un \emph{stathme euclidien} est une application $\varphi:\gA\to\NN$ 
qui vérifie
les \prts suivantes\footnote{Dans la littérature on trouve parfois  un \gui{stathme euclidien} défini
comme une application
$\varphi:\gA\to\NN\cup\so{-\infty}$, ou $\varphi:\gA\to\NN\cup\so{-1}$
(la valeur minimum étant toujours égale à $\varphi(0)$).} (grosso modo, on recopie la division euclidienne dans $\NN$)
\begin{itemize}
\item $\varphi(a)=0\iff a=0$.
\item $\forall a,b\neq0,\;\exists q,r,\;\;\;a=bq+r \hbox{ et } \varphi(r)<\varphi(b)$. 
%
%\item  
\end{itemize}
 Un \emph{anneau euclidien} est un anneau intègre non trivial donné avec un stathme euclidien. Notez que l'anneau est discret. On peut alors faire avec la \gui{division} qui est donnée par le stathme la même chose que l'on fait dans 
$\ZZ$ avec la division euclidienne.\index{anneau!euclidien}%
\index{euclidien!anneau ---}
\\
Les exemples les plus connus sont les suivants.
\begin{itemize}
\item $\ZZ$, avec $\varphi(x)=\abs{x}$,
\item  $\KX$ ($\gK$ un corps discret), avec $\varphi(P)=1+\deg(P)$ pour $P\neq0$,
\item $\ZZ[i]\simeq\aqo\ZZX{X^2+1}$, avec $\varphi(m+in)=m^2+n^2$, 
\item $\ZZ[i\sqrt 2]\simeq\aqo\ZZX{X^2+2}$, avec $\varphi(m+i\sqrt 2 n)=m^2+2n^2$. 
\end{itemize}
Dans ces exemples on a en outre l'\eqvcz: $\;x\in\Ati \iff \varphi(x)=1$.

\begin{enumerate}
\item \label{i1exoAnneauEuclidien} \emph{(Algorithme d'Euclide étendu)}
Pour tous $a,b$, il existe $u,v,a_1,b_1,g$ tels que 
$$ \cmatrix{g\cr 0}=\bloc{u}{v}{-b_1}{a_1} \cmatrix{a\cr b} \;\hbox{ et } \; ua_1+vb_1=1
.$$
En particulier, $\gen{a,b}=\gen{g}$ et $g$ est un pgcd de $a$ et $b$. 
Si $(a,b)\neq(0,0)$, $\frac{ab}g$ est un ppcm de  $a$ et $b$.
\item \label{i2exoAnneauEuclidien}
\begin{enumerate}
\item Montrer que l'anneau $\gA$ est principal.
\item   Faisons les hypothèses suivantes.
 \begin{itemize}
\item $\Ati$ est une partie détachable de $\gA$.
\item On dispose d'un test de primalité pour les \elts de $\gA\setminus\Ati$ au sens suivant: étant donné $a\in \gA\setminus\Ati$ on sait décider si $a$
est \irdz, et, en cas de réponse négative, écrire $a$ sous la forme $bc$
\hbox{avec $b$, $c\in \gA\setminus\Ati$}. 
\end{itemize}
Montrer qu'alors
$\gA$ vérifie le \gui{\tho fondamental de l'\ariz}
(\dcn unique en facteurs premiers, à association près).
\end{enumerate} 
\end{enumerate}
\emph{L'exemple $\ZZ[i]$.}
 On rappelle que $z = m+in \mapsto \ov z = m-in$ est un \auto \hbox{de $\ZZ[i]$} et
  que la norme $\rN=\rN_{\ZZ[i]/\ZZ}$  ($\rN(z)=z\ov z$) est un stathme euclidien: 
  prendre pour $q$ ci-dessus un \elt de
$\ZZ[i]$ proche de $a/b\in\QQ[i]$ et vérifier \hbox{que $\rN(r)\leq\rN(b)/2$}.\\
Pour connaître les \elts \irds de $\ZZ[i]$, il suffit de savoir décomposer tout nombre premier $p$. Ceci revient
à déterminer les \ids contenant $p\ZZ[i]$, \cade les
\ids de $\gZ_p:=\aqo{\ZZ[i]}p$. Or $\gZ_p\simeq \aqo{\FFp[X]}{X^2+1}$. On est donc ramené à trouver les diviseurs de $X^2+1$, donc à factoriser $X^2+1$, \hbox{dans $\FFp[X]$}.  
\begin{enumerate}\setcounter{enumi}{2}
\item \label{i3exoAnneauEuclidien} 
Montrer qu'a priori trois cas qui peuvent se présenter.
\begin{itemize}
\item $X^2+1$ est \ird dans $\FFp[X]$, et $p$ est \ird dans $\ZZ[i]$.
\item $X^2+1=(X+u)(X-u)$ dans $\FFp[X]$ avec $u\neq -u$, et alors: 

\snic{ \gen{p}=\gen{i+u,p}\gen{i-u,p}=\gen{m+in}\gen{m-in}$
 et $p=m^2+n^2.}

\vspace{2pt}
\item \label{i4exoAnneauEuclidien}
$X^2+1=(X+u)^2$ dans $\FFp[X]$, et alors $\gen{p}=\gen{i+u}^2$. Ceci se produit 
uniquement pour $p=2$, avec  $2=(-i)(1+i)^2$ (où $-i\in\ZZ[i]\eti$). 
\end{itemize}

\item \label{i5exoAnneauEuclidien}
Si $p\equiv 3\mod 4$, alors $-1$ n'est pas un carré dans $\FFp$.
Si $p\equiv 1\mod 4$, alors~$-1$ est un carré dans~$\FFp$. Dans ce cas donner un \algo
rapide pour écrire $p$ sous forme $m^2+n^2$ dans $\NN$.
\item \label{i6exoAnneauEuclidien} Soit $z\in\ZZ[i]$. On peut écrire $z=m(n+qi)$ avec $m,n,q\in\NN$ $\pgcd(n,q)=1. $ Donner un \algo rapide pour décomposer $z$ en facteurs premiers dans $\ZZ[i]$ connaissant une
\dcn  en facteurs premiers de $\rN(z)=m^2(n^2+q^2)$ dans $\NN$. 
\\
 Connaissant une \dcn en facteurs premiers de $s\in\NN$, décrire sous quelle condition $s$ est une somme de deux
carrés, ainsi que le nombre d'écritures~\hbox{$s=a^2+b^2$} avec $0<a\leq b$ dans $\NN$.

\item \label{i7exoAnneauEuclidien} Dire dans quels cas (relativement rares)
on peut généraliser la démarche précédente pour décomposer en produit de facteurs premiers les \itfs d'un anneau $\ZZ[\alpha]$, lorsque $\alpha$ est un entier \agqz.
\end{enumerate}
 
}
\end{problem}
%--- end -problem-----------------------------------------

%--- problem{exoPetitKummer}-------------
\begin{problem}
\label{exoPetitKummer} (Petit \tho de Kummer)
\\
{\rm  Le \pb \ref{exoAnneauEuclidien} peut se généraliser pour des anneaux d'entiers principaux de la forme $\ZZ[\alpha]$, mais ce cas est relativement rare. Bien au contraire, 
le petit \tho de Kummer donne la \dcn d'un nombre premier (dans $\NN$) en produits d'\idemas $2$-engendrés pour presque tous les nombres premiers, dans tous les 
anneaux d'entiers. Ceci montre la supériorité intrinsèque des \gui{nombres idéaux} introduits par Kummer. En outre, l'argument est extrêmement simple
et ne nécessite que le \tho chinois. Cependant les nombres premiers qui ne tombent
pas sous la coupe du petit \tho de Kummer constituent en fait le c{\oe}ur de la
théorie \agq des nombres, c'est eux qui ont nécessité une mise au point fine de la théorie (selon deux méthodes distinctes par Kronecker et Dedekind), sans laquelle tout progrès décisif e\^ut été impossible.

 On considère un zéro  $\alpha$ d'un \polu \irdz~\hbox{$f(T)\in\ZZ[T]$},
de sorte \hbox{que $\ZZ[\alpha]\simeq\aqo{\ZZ[T]}{f(T)}$}. On note $\Delta=\disc(f)$.
\begin{enumerate}
\item Soit $p$ un nombre premier qui ne divise pas $\Delta$.
\begin{itemize}
\item  Montrer que $f(T)$ est \spl dans $\FFp[T]$. 
\item  On  décompose
$f(T)$ dans $\FFp[T]$ sous forme $\prod_{k=1}^{\ell}Q_k(T)$ avec les $Q_k$ \irds \mons distincts. On pose $q_k=Q_k(\alpha)$ (a vrai dire ce n'est défini que modulo $p$, mais on peut relever $Q_k$ dans $\ZZ[T]$).
Montrer que dans $\ZZ[\alpha]$ on a $\gen{p}=\prod_{k=1}^{\ell}\gen{p,q_k}$
et que les \ids $\gen{p,q_k}$ sont maximaux, distincts et inversibles.
En particulier, si $\ell=1$, $\gen{p}$ est maximal. 
\item  Montrer que cette \dcn reste valable dans tout anneau $\gA$ tel que 
$\ZZ[\alpha]\subseteq\gA\subseteq\gZ$, où $\gZ$ est l'anneau
des entiers de $\QQ[\alpha]$.
\end{itemize}
\item Soit $a\in\ZZ[\alpha]$ tel que $A=\rN_{\ZZ[\alpha]/\ZZ}(a)$ soit étranger
à $\Delta$. Soit $\fa=\gen{b_1,\ldots,b_r}$ un \itf de $\ZZ[\alpha]$ contenant $a$. Montrer que dans $\ZZ[\alpha]$ l'\id $\fa$ est \iv et 
se décompose 
en produits d'\idemas qui divisent les facteurs premiers de $A$. 
Enfin, cette \dcn est unique à l'ordre près des facteurs et tout ceci reste
valable dans tout anneau $\gA$ tel que ci-dessus.
\end{enumerate} 
}
\end{problem}
%--- end -problem-----------------------------------------

%--- problem{exoPolCyclotomique}-------------
\begin{problem}\label{exoPolCyclotomique}  {(Le \pol cyclotomique $\Phi_n$)}\index{polynome@\pol!cyclotomique}\\
{\rm
Dans $\gA[X]$, le \polz~$X^n - 1$ est \spl \ssi $n \in \Ati$% (dérivée)
.\\
Notons $\gQ_n$ un \cdr au dessus de $\QQ$ pour ce \polz.
Soit $\UU_n$ le groupe des racines $n$-ièmes de l'unité dans  $\gQ_n$.
C'est un groupe cyclique
d'ordre~$n$, qui possède donc $\varphi(n)$ \gtrs (racines
primitives $n$-ièmes de l'unité). On définit $\Phi_n(X) \in \gQ_n[X]$ par
$\Phi_n(X) = \prod_{o(\xi) = n} (X - \xi)$. C'est un \polu
de degré $\varphi(n)$.  On a l'\egt fondamentale

\snic {
X^n - 1 = \prod_{d \divi n} \Phi_d(X),
}

%\sni
qui permet de démontrer par \recu sur $n$ que $\Phi_n(X) \in \ZZ[X]$.

\begin {enumerate}
\item
On va montrer que $\Phi_n(X)$ est \ird dans $\ZZ[X]$
(donc dans $\QQ[X]$: proposition \ref{propZXfactor}).
Soient  $f$, $g$ deux \polus de $\ZZ[X]$
avec  $\Phi_n = fg$ et $\deg f \ge 1$; il faut prouver que $g = 1$.

\begin {itemize}
\item [a.]
Il suffit de prouver que $f(\xi^p) = 0$ pour tout premier
$p \nedivi n$ et pour tout zéro $\xi$ de $f$ dans $\gQ_n$.
\item [b.]
On suppose que $g(\xi^p) = 0$ pour un zéro $\xi$ de $f$ dans $\gQ_n$.
Examiner ce qui se passe dans $\FF_p[X]$ et conclure.
\end {itemize}

\item
Fixons une racine  $\xi_n$ de $\Phi_n$ dans $\gQ_n$.
\\
Montrer que $\gQ_n = \QQ(\xi_n)$
et qu'avec $(\QQ,\gQ_n,\Phi_n)$, on est dans la situation
galoisienne \elr du lemme \ref {lemGaloiselr}.
\\
  Décrire des
\isos explicites de groupes:

\snic {
\Aut(\UU_n)  \simeq (\ZZ/n\ZZ)^{\!\times} \simeq \Gal(\gQ_n/\QQ).
}

\item
Soit $\gK$ un corps de caractéristique $0$. Que peut-on dire d'un
\cdrz~$\gL$ de $X^n-1$ au dessus de $\gK$?

\end {enumerate}
}
\end {problem}
%--- end -problem-----------------------------------------

%--- problem{exoCyclotomicRing}-------------
\begin{problem}\label{exoCyclotomicRing}
{(L'anneau $\ZZ[\!\root n \of 1]$: \ddpz, \fcn des \idsz)}
\\
 {\rm  
Soit $\Phi_n(X) \in \ZZ[X]$ le \pol cyclotomique d'ordre $n$, \ird sur
$\QQ$. 
On note $\gQ_n  = \QQ(\zeta_n)\simeq \aqo {\QQ[X]} {\Phi_n}$.  Le groupe
multiplicatif $\UU_n$
 engendré par $\zeta_n$ (racine primitive $n$-ième de l'unité) est cyclique d'ordre~$n$.
\\
On va démontrer 
entre autres
que l'anneau $\gA = \ZZ[\UU_n] = \ZZ[\zeta_n]\simeq\aqo{\ZZ[X]}{\Phi_n}$ est 
un \ixx{domaine}{de Prüfer}\index{Prufer@Prüfer!domaine de ---}: 
un anneau intègre dont les \itfs non nuls
sont \ivs (cf. section \ref{secIplatTf} et chapitre~\ref{ChapAdpc}).

 \emph {1.}
Soit $p \in \NN$ un nombre premier. On montre ici que $\sqrt {p\gA}$ est un
\idp et on l'explicite comme un produit fini d'\idemas \ivs 2-engendrés.  On
considère les facteurs \irds distincts de $\Phi_n$ modulo $p$ que l'on
relève en des \polus $f_1, \ldots, f_k \in \ZZ[X]$. On note $g = f_1
\cdots f_k$ (de sorte que $\ov g$ est la partie sans facteur carré de
$\Phi_n$ modulo $p$) et $\fp_i = \gen {p, f_i(\zeta_n)}$ pour $i \in \lrbk$.
\begin {enumerate}\itemsep0pt
\item [\emph {a.}]
Montrer que $\fp_i$ est un \idema et que

\snic {
\sqrt {p\gA} = \gen {p, g(\zeta_n)} = \fp_1 \cdots\ \fp_k
}

\item [\emph {b.}]
Si $p$ ne divise pas $n$, montrer que~$\ov g = \ov {\Phi_n}$, donc $\sqrt
{p\gA} = \gen {p}$ est un \idpz.

\item [\emph {c.}]
On suppose que $p$ divise $n$ et l'on écrit $n = mp^k$ avec $k \ge 1$,
$\pgcd(m, p) = 1$.\\  
En étudiant la factorisation de $\Phi_n$ modulo $p$,
montrer que~$\ov g = \ov {\Phi_m}$. En déduire que~$\sqrt {p\gA} = \gen {p,
\Phi_m(\zeta_n)}$.  Montrer ensuite que $p \in \gen {\Phi_m(\zeta_n)}$, et donc
que $\sqrt {p\gA} = \gen {\Phi_m(\zeta_n)}$ est un \idpz.

\item [\emph {d.}]
En déduire que $p\gA$ est un produit de la forme
$\fp_1^{e_1} \ldots \fp_k^{e_k}$.
\end {enumerate}
 \emph {2.}
Soit $a \in \ZZ \setminus \{0\}$; montrer que $a\gA$ est un produit d'\idemas
\ivs à deux \gtrsz.  En déduire que dans $\gA$ 
 tout \itf non nul se décompose en un produit d'\idemas
\ivs $2$-engendrés et que la \dcn est unique à l'ordre près des facteurs.

}

\end {problem}
%--- end -problem-----------------------------------------

%--- problem{exoSommeGauss}-------------
\begin{problem}\label{exoSommeGauss}
{(Une \prt \elr des sommes de Gauss)}
\\ 
{\rm  
On désigne par $\gk$ un corps fini de cardinal $q$ et $\gA$ un
anneau intègre. On considère:
\begin{itemize}
\item un \gui {caractère
multiplicatif} $\chi : \gk\eti \to \Ati$,  i.e. un morphisme de groupes multiplicatifs,
\item un \gui{caractère additif}  $\psi : \gk
\to \Ati$,  i.e.  un morphisme de groupes 

\centerline{$\psi :
(\gk, +) \to (\Ati, \times)$.}
\end{itemize}
   On suppose que ni $\chi$ ni $\psi$ ne sont
triviaux et l'on prolonge $\chi$ à $\gk$ tout entier \hbox{via $\chi(0) =
0$}. Enfin, on définit la somme de Gauss de $\chi$, relativement
à $\psi$, par:

\snic {
G_\psi(\chi) = \sum_{x \in \gk} \chi(x) \psi(x) = 
\sum_{x \in \gk\eti} \chi(x) \psi(x).
}

%\sni
On va montrer que

\snic {
G_\psi(\chi)G_\psi(\chi^{-1}) = q\chi(-1) 
,}

%\sni

et donner des applications \aris de ce résultat (question \emph{4}).

%\sni
\emph {1.}
Soit  $G$  un groupe fini et $\varphi : G \to \Ati$ un \homo non trivial.  Montrer que $\sum_{x \in G} \varphi(x) = 0$. 

 \emph {2.}
Montrer que:

\snic {
\sum_{x + y = z} \chi(x) \chi^{-1}(y) = \cases {
-\chi(-1) &si $z \ne 0$,\cr
(q-1) \chi(-1) &sinon.\cr}
}

%\sni
\emph {3.}
En déduire que $G_\psi(\chi)G_\psi(\chi^{-1}) = q\chi(-1)$.

 \emph {4.}
On considère $\gk = \FF_p$ où $p$ est un nombre premier impair,
 $\gA =\QQ(\!\root p \of 1)$, et $\zeta$  une racine primitive $p$-ième de l'unité
dans $\gA$. Les caractères $\psi$ et $\chi$ sont définis par:

\snic {
\psi(i \bmod p) = \zeta^i, \qquad  \chi(i \bmod p) = {i \legendre p}
\quad \hbox {(symbole de Legendre)}
}

\begin {enumerate}%\itemsep0pt
\item [\emph {a.}]
Alors, $\chi = \chi^{-1}$, les sommes de Gauss $G_\psi(\chi)$, 
$G_\psi(\chi^{-1})$ sont égales à

\snic {
\tau \eqdefi \sum_{i \in \FF_p^*} {i \legendre p} \zeta^i,
}

%\sni
et en posant $p^* = (-1)^{p-1 \over 2} p$ (de sorte que $p\etl\equiv1\mod 4$), on obtient:

\snic {
\tau^2 = p^*,  \quad \hbox {en particulier,} \quad
\QQ(\sqrt {p^*}) \subseteq \QQ(\!\root p \of 1)
}

\item [\emph {b.}]
On définit $\tau_0 = \sum_{i \in \FF_p^{\times 2}} \zeta^i, \, 
\tau_1 = \sum_{i \in \FF_p\eti \setminus \FF_p^{\times 2}} \zeta^i$
  de sorte que  
$\tau  = \tau_0 - \tau_1.$
Montrer que $\tau_0$ et $\tau_1$ sont les racines de $X^2
+ X + {1 - p^* \over 4}$ et que l'anneau~\hbox{$\ZZ[\tau_0] = \ZZ[\tau_1]$} est
l'anneau des entiers de $\QQ(\sqrt {p^*})$.

\end {enumerate}
}
\end {problem}
%--- end -problem-----------------------------------------

%--- problem{exoDedekindPolynomial}-------------
\begin{problem} 
\label{exoDedekindPolynomial}\index{Dedekind!poly@\pol de ---} 
{(Le \pol de Dedekind $f(X) = X^3 + X^2 -2X + 8$)}
\\ 
{\rm
Le but de ce \pb est de fournir un exemple d'anneau $\gA$ d'entiers
de corps de nombres qui n'est pas une $\ZZ$-algèbre monogène\footnote{Une \Alg $\gB$ est dite \emph{monogène} lorsqu'elle est engendrée en tant qu'\Alg par un unique \elt $x$. Ainsi, $\gB=\gA_1[x]$, où $\gA_1$ est l'image de $\gA$ dans $\gB$.}.

\begin {enumerate}\itemsep0pt\mou
\item
Montrer que $f$ est \ird dans $\ZZ[X]$ et que $\disc(f) = -2\,012 = -2^2
\times 503$.

\item
Soit $\alpha$ une racine de $f(X)$.  Montrer que $\beta =
4\alpha^{-1}$ est entier sur~$\ZZ$, que  

\centerline{$\gA = \ZZ \oplus \ZZ\alpha \oplus \ZZ\beta$} 

est l'anneau des entiers de $\QQ(\alpha)$ et
que $\Disc_{\gA/\ZZ} = -503$.

\item
Montrer que le nombre premier $p = 2$ est totalement décomposé dans
$\gA$, autrement dit que $\gA/2\gA \simeq \FF_2 \times \FF_2 \times \FF_2$.
En déduire que $\gA$ n'est pas une \ZZlg monogène.

\item
\emph{(\'Evitement du conducteur, Dedekind)} Soient $\gB \subseteq \gB'$
deux anneaux, $\ff$ un idéal de $\gB$ vérifiant $\ff\gB' \subseteq
\gB$; a fortiori $\ff\gB' \subseteq\gB'$ et $\ff$ est aussi un \id
de $\gB'$. Alors, pour tout idéal $\fb$ de $\gB$ tel que $1 \in \fb+\ff$, en
posant $\fb' = \fb\gB'$, le morphisme canonique $\gB/\fb \to \gB'/\fb'$ est un
\isoz.

\item
En déduire que $2$ est un \emph {diviseur essentiel} de $\gA$: on entend par là que $2$ divise l'indice $\idg{\gA:\ZZ[x]}$ quelque soit l'\elt
primitif~$x$ de~$\QQ(\alpha)/\QQ$ entier sur~$\ZZ$.
\end {enumerate}
}
\end {problem}
%--- end -problem-----------------------------------------

%--- problem{exoGaloisNormIdeal}-------------
\begin{problem}\label{exoGaloisNormIdeal} {(Norme d'un idéal en terrain quasi-galoisien)}
 \\
{\rm
Soit $(\gB, \gA, G)$ où $G \subseteq \Aut(\gB)$ est un groupe fini,
 et $\gA = \gB^G=\Fix_\gB(G)$.   Si~$\fb$ est un idéal de $\gB$, on note $\rN'_G(\fb) =
\prod_{\sigma \in G} \sigma(\fb)$ (idéal de $\gB$) et $\rN_G(\fb) =
\gA \cap \rN'_G(\fb)$ (idéal de $\gA$).
\begin {enumerate}\itemsep0pt\mou
\item %1
Montrer que $\gB$ est entier sur $\gA$.
\item %2
Soient $\gB = \ZZ[\sqrt d]$ où $d \in \ZZ$ n'est pas un carré, $\tau$
l'\auto (noté aussi~\hbox{$z
\mapsto \ov z$})  défini par $\sqrt d \mapsto -\sqrt d$,  et $G = \gen {\tau}$. Donc $\gA = \ZZ$. On suppose que~\hbox{$d \equiv 1
\bmod 4$} et l'on pose $\fm = \gen {1+\sqrt d, 1-\sqrt d}$.
\begin {enumerate}
\item [\emph{a.}]
On a $\fm = \ov\fm$, $\rN'_G(\fm) = \fm^2 = 2\fm$ et $\rN_G(\fm) = 2\ZZ$.
En déduire que $\fm$ n'est pas \iv et que l'on
n'a pas $\rN'_G(\fm) = \rN_G(\fm)\gB$.
\item [\emph{b.}]
Montrer que $\ZZ[\sqrt d]/\fm \simeq \FF_2$; donc $\fm$ est d'indice $2$ dans
$\ZZ[\sqrt d]$ mais 2 n'est pas le pgcd des $\rN_G(z)$, $z \in \fm$.
Vérifier \egmt que 
$\fb \mapsto \idg{\gB:\fb}$
n'est pas multiplicative
sur les \ids non nuls de $\gB$.
\end {enumerate}

\item %3
On suppose $\gB$ \icl et $\gA$ de Bézout. Soit $\fb \subseteq \gB$ un \id de type~fini.
\begin {enumerate}
\item [\emph{a.}]
Donner un $d \in \gA$ tel que $\rN'_G(\fb) = d\gB$.  En particulier, si $\fb$
est non nul, il est \ivz.
Ainsi, $\gB$ est un \ddpz.
\item [\emph{b.}]
Montrer que $\rN_G(\fb) = d\gA$, donc $\rN'_G(\fb) = \rN_G(\fb)\gB$.
\item [\emph{c.}]
On suppose que 
%$\gA$ est un \emph{anneau de Smith} (voir la \dfn \paref{exoSmith}) et que $\fb$ est non nul.
%On note $a_1$, $\ldots$, $a_k \in \gA$ les  \emph{facteurs invariants}  du $\gA$-module $\gB/\fb$ (voir le \thref{prop unicyc}). En particulier,
le \Amo $\gB/\fb$ est isomorphe à $\aqo{\gA}{a_1} \times \cdots \times \aqo{\gA}{a_k}$.
\\
Montrer que $\rN_G(\fb) = \gen {a_1 \cdots a_k}_\gA$.

\item [\emph{d.}]
On suppose $\#G = 2$.  Expliciter, en fonction d'un \sgr fini de $\fb$, des
\elts $z_1$, \dots, $z_m \in \fb$ tels que $\rN_G(\fb) = \gen {\rN(z_i), i \in
\lrb {1..m}}_\gA$.
\end {enumerate}
\end {enumerate}

}
\end {problem}
%--- end -problem-----------------------------------------

%%%%%%%%%%%%%%%%%%%%%%%%%%%%%%%%%%%%%%%%%
%--- problem{exoLemmeFourchette}-------------
\begin{problem}\label{exoLemmeFourchette} 
 {(Lemme de la fourchette)} \index{Lemme de la fourchette}
{\rm
\begin{enumerate}
\item Soit $\gA$ un \acl de corps de fractions $\gk$, $\gK$ une extension finie
séparable de $\gk$ de degré $n$, $\gB$ la \cli de $\gA$
dans~$\gK$.
Montrer qu'il existe une base $(\ue) = (e_1, \ldots, e_n)$ de $\gL/\gK$
contenue dans~$\gB$. On note~\hbox{$\Delta = \disc(\ue)$} et
$(\ue') = (e'_1, \ldots, e'_n)$ la base traciquement duale de $(\ue)$.
Montrer les inclusions:

\snic {
\bigoplus_{i=1}^n \gA e_i \;\subseteq\; \gB \;\subseteq\;
\bigoplus_{i=1}^n \gA e'_i \;\subseteq\;
\Delta^{-1} \bigoplus_{i=1}^n \gA e_i.
}

\end{enumerate}

Dans la suite $\gA=\ZZ$ et $\gk=\QQ$; $\gK$ est donc un corps de nombres
et $\gB=\gZ$ est son anneau d'entiers.   
On considère un $x\in\gZ$ tel que $\gK=\QQ[x]$. 
\\
Soit $f(X)=\Mip_{\QQ,x}(X)\in\ZZ[X]$ et $\delta^2$ le plus grand facteur carré de $\disc_X( f)$. 
\\
D'après la proposition~\ref{propAECDN},
 $\gZ$ est un $\ZZ$-module libre de rang $n = \dex{\gL : \QQ}$, et l'on \hbox{a $\ZZx\subseteq\gZ\subseteq \fraC 1 \delta \ZZx$}.
  Ceci est légèrement plus précis que le résultat du point \emph{1.}
\\
On considère une \ZZlg  \tf $\gB$ intermédiaire entre $\ZZx$ et $\gZ$.
Comme c'est un \ZZmo \tfz, $\gB$ est \egmt un \ZZmo libre de rang $n$.
Le cas le plus important est celui où $\gB=\gZ$.
\\ 
Le but du \pb est de préciser  une $\ZZ$-base de $\gB$ de la forme
$$\ndsp
\cB=\big(\fraC{g_0}{d_0},\fraC{g_1(x)}{d_1},\fraC{g_2(x)}{d_2},\dots,\fraC{g_{n-1}(x)}{d_{n-1}}\big)  
$$   
avec
$g_k\in\ZZ[X]$  de  degré  $k$  pour tout  $k$, et les $d_k>0$ les plus petits possibles. On va établir ce résultat avec des \polus $g_k$
 et $1=d_0\mid d_1\mid d_2\mid\cdots\mid d_{n-1}$.
\\
Le corps $\gK$ est un \Qev de base $(1,x,\dots,x^{n-1})$ et \hbox{pour $ k \in \lrb{0..n-1}$}, on note $\pi_k : \gK \to \QQ$ la forme \lin 
composante \hbox{sur $x^k$} et:
$$
Q_k = \bigoplus\nolimits_{i=0}^k \QQ\, x^i,\; Z_k = \frac1\delta \bigoplus\nolimits_{i=0}^k \ZZ\, x^i, \quad \hbox{et }{  F_k = Q_k \cap \gB=  Z_k \cap \gB}.
$$
Il est clair que $Q_0 = \QQ$, $Q_{n-1} = \gK$, $F_0 = \ZZ$  et $F_{n-1} = \gB$.
\begin{enumerate}\setcounter{enumi}{1}
\item Montrer que le \ZZmo $F_k$ est libre de rang $k+1$.  
\\ Le \ZZmo $\pi_k(F_k)$ est un sous-\ZZmo \tf de $\fraC1\delta\ZZ$.
Montrer qu'il est de la forme $\fraC1{d_k}\ZZ$ pour un $d_k$
qui divise $\delta$. NB:  a $d_0=1$. 
\item  On note $y_k$ un \elt de $F_k$ tel que $\pi_k(y_k)=\fraC1{d_k}$.
\\
On écrit $y_k$ sous la forme $f_k(x)/d_k$, avec $f_k\in\QQ[X]$  \mon et de degré~$k$. Il est clair que $y_0=1$. Mais les autres $y_i$ ne sont pas déterminés de manière unique. 
Montrer que $(1,y_1,\dots,y_k)$ est une $\ZZ$-base
de $F_k$.
\item Montrer que si $i+j\leq n-1$, on a $d_id_j\mid d_{i+j}$. 
En particulier $d_i$ divise $d_k$ \hbox{si $1\leq i<k\leq n-1$}.
En déduire aussi \hbox{que $d_1^{n(n-1)/2}$} divise~$\delta$. 
\item Montrer que $d_ky_k\in\ZZ[x]$ 
pour chaque $ k \in \lrb{0..n-1}$.  
En déduire\hbox{ que $f_k\in\ZZ[X]$} et que $\big( 1,  f_1(x),\dots,  f_{n-1}(x)\big)$
est une $\ZZ$-base de $\ZZx$.  
\item Montrer que   
$\cB={\big( 1, \fraC 1 {d_1} f_1(x),\dots, \fraC 1 {d_{n-1}} f_{n-1}(x)\big)}
$ 
est une $\ZZ$-base de $\gB$ adaptée à l'inclusion
$\ZZx\subseteq \gB$. 
Les $d_i$ sont donc les facteurs invariants de cette inclusion, et $\prod_{i=1}^{n-1}d_i$ est égal à 
l'indice $\idG{\gB:\ZZx}$ qui divise $\delta$. 
\end{enumerate}

}
\end {problem}
%--- end -problem-----------------------------------------

%--- problem{exoPolynomialAutomorphism}-------------
\begin{problem}
\label {exoPolynomialAutomorphism} 
  {(Changements de variables, automorphismes polynomiaux
    et \imN méthode de Newton)}
    \\
{\rm
Soient $F = (F_1, \ldots, F_n)$ avec $F_i \in \gA[X] = \gA[X_1, \ldots,
X_n]$ et $\theta_F : \gA[X] \to \gA[X]$ le morphisme de $\gA$-\algs
réalisant $X_i \mapsto F_i$; on a donc $\theta_F(g) = g(F)$. On suppose que
$\gA[X] = \gA[F]$: il existe donc $G_i \in \gA[X]$ vérifiant $X_i =
G_i(F)$, ce que l'on écrit classiquement (avec quelques abus) $X = G(F)$ et
parfois $X = G \circ F$ (au sens des applications de $\gA[X]^n$
dans $\gA[X]^n$). \\
On notera le renversement $\theta_F \circ \theta_G =
\I_{\gA[X]}$.
\\
  On va montrer ici que $\theta_G \circ \theta_F =
\I_{\gA[X]}$ ou encore $X = F(G)$.
\\
 En conséquence (cf. question 1) $G$ est déterminé de
manière unique, $\theta_F$ est un automorphisme de $\gA[X]$ et $F_1,
\ldots, F_n$ sont \agqt indépendants sur $\gA$.
\\
  L'idée consiste à utiliser l'anneau des
séries formelles $\gA[[X]]$ ou
du moins les quotients $\gA[X]/\fm^d$ où $\fm = \gen {X_1, \ldots,
X_n}$. Soit $F = (F_1, \ldots, F_n) \in \gA[[X]]^n$; on étudie à quelle
condition il existe $G = (G_1, \ldots, G_n)$, $G_i \in \gA[[X]]$ sans terme
constant, vérifiant $F(G) = X$. On a alors $F(0) = 0$, et en posant $J_0 =
\JJ(F)(0)$, on obtient $J_0 \in \GLn(\gA)$ (puisque $\JJ(F)(0) \circ \JJ(G)(0) =
\I_{\Ae n}$). \\
On va montrer la réciproque: dans le cas où $F(0) = 0$ et $J_0 \in
\GLn(\gA)$, il \hbox{existe $G = (G_1, \ldots, G_n)$}, avec les $G_i \in \gA[[X]]$, $G_i(0)=0$, et $F(G) = X$.
\begin {enumerate}
\item
En admettant cette réciproque, montrer que $G$ est unique et que
$G(F) = X$.
\item
Soit $\gS \subset \gA[[X]]$ l'ensemble des séries formelles sans terme
constant; $\gS^n$ est, pour la loi de composition, un \mo dont $X$ est le
neutre.  On rappelle le principe de la méthode de Newton pour résoudre en
$z$ une équation $P(z) = 0$: introduire l'itérateur $\Phi : z \mapsto z -
P'(z)^{-1}P(z)$ et la suite $z_{d+1} = \Phi(z_d)$ avec un $z_0$ adéquat; ou
bien une variante $\Phi : z \mapsto z - P'(z_0)^{-1}P(z)$. Pour résoudre en
$G$, $F(G) - X = 0$, vérifier que cela conduit à l'itérateur sur
$\gS^n$:

\snic{
\Phi : G \mapsto G - J_0^{-1} \cdot (F(G) - X)
}

\item % 3
On introduit $\val : \gA[[X]] \to \NN \cup \{\infty\}$: $\val(g) = d$ signifie
que $d$ est le degré (total) minimum des \moms de $g$, en convenant que
$\val(0) = +\infty$. On a donc $\val(g) \ge d$ \ssi $g \in \fm^d$. On note
pour $g$, $h \in \gA[[X]]$ \hbox{et $G$, $H \in \gA[[X]]^n$}:

\snic {
d(f,g) = {1 \over 2^{\val(f-g)}}, \qquad
d(F,G) = \max_i d(F_i, G_i)
}

%\sni
Montrer que $\Phi$ est contractant: $d\big(\Phi(G), \Phi(H)\big) \le d(G, H)/2$. En
déduire que~$\Phi$  admet un unique point fixe $G \in \gS^n$, unique
solution de $F(G) = X$.
\item
Résoudre le \pb initial relatif aux \polsz.
\item
Vérifier que les \syss suivants sont des \cdvs et
expliciter leurs inverses (dans $\ZZ[X,Y,Z]$ puis dans $\ZZ[X_1, X_2, X_3,
X_4, X_5]$):
\[ 
\begin{array}{ccc} 
 (X - 2fY - f^2 Z,\ Y + fZ,\ Z) \;\hbox { avec }\, f = XZ + Y^2 ,  \\[1mm] 
 (X_1 + 3X_2X_4^2 - 2X_3X_4X_5,\  X_2 + X_4^2  X_5,\  X_3 + X_4^3,\
X_4 + X_5^3,\   X_5). 
\end{array}
\]
\end {enumerate}
}
\end {problem}
%--- end -problem-----------------------------------------

%%%%%%%%%%%%%%%%%%%%%%%%%%%%%%%%%%%%%%%%%%%%%%%%%%%%%%%%%%%%%%%%%%%%%%%%%%%

%: sinotenglish
\sinotenglish{
%:--- Problem  exoFinitudeClassesIdeaux
\begin{problem} \label{exoFinitudeClassesIdeaux} {(Finitude de l'ensemble des classes d'\ids d'un 
anneau de nombres)}
\\
{\rm  
Un \emph{anneau de nombres} $\gA$ est un anneau intègre qui en tant que groupe additif est un
$\ZZ$-module libre de rang fini $n>0$. 
Son corps de fractions $\gK = \Frac(\gA)$
est un corps de nombres de degré $n$. On montre ici que l'ensemble
des classes d'\itfs de $\gA$ est fini (deux \itfs $\fa$ et $\fb$ sont dans la même classe
s'il existe  $x$ et $y$ non nuls dans $\gA$ tels que $x\fa=y\fb$).

\snii 
Pour $x =(\xn) \in \QQ^n$, on note $\nsup{x} = \max_i |x_i|$.

\snii \emph {1.}
Soient $n \in \NN^*$ et un $K > 0$. Montrer qu'il
existe $d \in \NN^*$ tel que pour \hbox{tout $x \in \QQ^n$},
il existe $y \in \ZZ^n$ et $m \in \lrb{1..d}$ vérifiant
$\nsup{mx - y} < K$.
\\
Idée: pour
 $N \in \NN^*$, montrer que pour tout $x \in \QQ^n$, il existe $y \in
\ZZ^n$ et $m \in \lrb{1..N^n}$ vérifiant $\nsup{mx - y} < 1/N$.

\snii \emph {2.}
On note $\rN = \rN_{\gK\sur\QQ}$.
On fixe une $\QQ$-base $(e_1, \ldots, e_n)$ de $\gK$
constituée d'\elts de $\gA$. On définit alors
$\nsup {\ } : \gK \to \QQ_+$  par 
$\nsup{x} = \max_i |x_i|$ \hbox{où $x=\sum_i x_ie_i$}.
Montrer qu'il existe $C > 0$ tel que
$|\rN(x)| \le C \nsup {x}^n$,  $\forall x \in \gK$.

\snii \emph {3.}
Montrer qu'il existe $d \in \NN^*$,
attaché uniquement à $\gA$,
tel que pour tout
$x \in \gK$, il y ait $m \in \lrb{1..d}$ et $q \in \gA$
vérifiant $|\rN(mx - q)| < 1$. En déduire que pour
$a \in \gA$ \hbox{et $b \in \gA\setminus \{0\}$}, il y a un
$m \in \lrb{1..d}$ tel que $|\rN(ma-bq)| < |\rN(b)|$.

\snii \emph {4.}
On pose $D = d!$.  Soit $\fb$ un \itf non nul  de $\gA$. On veut montrer
qu'il existe $b \in \fb\setminus \{0\}$ tel que $D\fb \subset \gen {b}$.

\begin {enumerate}
\item [\emph {a.}]
Soit $b \in \fb\setminus\{0\}$ tel que $|N(b)|$ soit minimum (ce qui du point
de vue des \clama ne pose pas de problème puisque pour $b \in
\fb\setminus\{0\}$, on a $|N(b)| \in \NN^*$).  Montrer que $b$ convient.

\item [\emph {b.}]
Fournir une preuve constructive de l'existence d'un tel $b$.

\item [\emph {c.}]
En déduire qu'il existe un \id $\fa$ de $\gA$ associé à $\fb$ tel que $D
\in \fa$.
\end {enumerate}

\noindent \emph {5.}
Conclure.

}
\end{problem}

%%%%%%%%%%%%%%%%%%%%%%%%%%%%%%%%%%%%%%%%%%%%%%%%%%%%%%%%%%%%%%%%%%%%%%%%%%%

%:--- Problem{pbMatrixAndIdealClasses}-------------
\begin{problem}\label{pbMatrixAndIdealClasses}  {(Classes de similitude de matrices et classes d'\idsz)}
\\
{\rm  
Soit $\gA$ un anneau, $\gK = \Frac(\gA)$ son anneau total des
fractions,~\hbox{$f \in \gA[X]$} un \pol\unt de degré $n \ge 1$ et $\gB = \gA[x] =
\aqo{\AX}{f}$. On note $C_f \in \Mn(\gA)$ la matrice compagne de $f$.  
Deux \ids $\fb$ et $\fb'$ d'un anneau $\gC$ sont dits  \emph{associés}, ou \emph{dans la même classe} s'il
existe deux \elts \ndzs $b$, $b'\in \gC$ tels \hbox{que $b'\fb = b\,\fb'$}.%
\index{classe!d'ideaux@d'\idsz}  

\noindent On étudie
dans ce problème une correspondance bijective entre 

\begin{itemize}
\item les classes de
similitude (sur $\gA$) de matrices de $\Mn(\gA)$ semblables sur $\gK$ à $C_f$
d'une part,
et
\item  les classes de certains \ids de $\gB$, à savoir ceux
qui sont des \Amo libres de rang $n$, d'autre part. 
\end{itemize}

\noindent \emph {1.}
Montrer que $C_f$ est semblable sur $\gA$ à sa transposée.  Montrer
que l'on peut prendre pour $Q$ satisfaisant $\tra {C_f} = Q^{-1} C_f Q$ une matrice de Hankel supérieure\footnote{Une matrice $R=((r_{i,j})_{i,j})$ est \emph{de Hankel} lorsque l'on a: $i+j=k+\ell\Rightarrow
r_{ij}=r_{k\ell}$. Elle est \emph{de Hankel supérieure} lorsque les \coes en dessous de l'antidiagonale sont nuls.} de \deter
$(-1)^{\lfloor n/2\rfloor}$.

\snii \emph {2.}
Montrer l'\idt $QS = \Adj(x\In-C_f)$, où   $S=(x^{i+j})_{0\le i,j \le n-1}$.

\snii \emph {3.}
Soit $M \in \Mn(\gA)$ semblable à $C_f$ sur $\gK$ et
$P \in \Mn(\gA)$ telle que $PM = C_fP$ et $\det(P)$ \ndzz. On définit $\vep_1, \ldots, \vep_n \in \gB$ par

\snic {
[\,\vep_1\; \ldots\; \vep_n\,] = [\,1\; x\; \ldots\; x^{n-1}\,]\cdot P
\;.}

\vspace{-.8em}
\begin {enumerate}
\item [\emph {a.}]
Montrer que $x\,[\,\vep_1\; \ldots\; \vep_n\,] = [\,\vep_1\; \ldots\; \vep_n\,]\cdot M$, puis que $\fb \eqdefi \gA\vep_1 + \cdots + \gA\vep_n$ est un \id de
$\gB$ et un \Amo libre de rang $n$ de base $(\vep_1, \ldots,\vep_n)$. 
Vérifier \egmt que $\fb$ contient un \elt de $\gA$
\ndz dans $\gA$ (donc dans $\gB$) et que la matrice de la multiplication par $x$
dans cette base est $M$.
\item [\emph {b.}]
Montrer que la classe d'\eqvc de $\fb$ ne dépend pas du choix de $P$
puis qu'elle ne dépend que de la classe de similitude de $M$ sur $\gA$.
\item [\emph {c.}]
Montrer que la suite ci-dessous est exacte:

\snic {
\gB^n\vvvvvers{x\In - M} \gB^n \vvvvvers{[\,\vep_1\; \ldots\; \vep_n\,]} 
\fb = \gA\vep_1 \oplus \cdots \oplus \gA\vep_n \to 0.
}

\snii 
On a donc un \iso de $\gB$-modules $\fb \simeq \gB^n\sur{\Im(x\In-M)}$.
\\
Et la matrice $x\In-M$ est  de rang $n-1$ \ssi $\fb$ est un \id projectif de
rang 1 (donc \iv car il contient un \elt
\ndzz).
\end {enumerate}

\snii \emph {4.}
Réciproquement, soit $\fb$ un \id de $\gB$ qui possède une $\gA$-base $(\vep_1, \ldots, \vep_n)$.  Montrer que $\fb$ contient un
\elt \ndz de $\gA$ (donc \ndz dans $\gB$).  Soit~$M$ la matrice de la
multiplication par $x$ dans cette base.  Montrer que $M$ est semblable à
$C_f$ sur $\gK$ et que la classe d'\ids associée à $M$ est celle de $\fb$.

%%%%%
\snii \emph {5.}
On reprend le contexte de la question \emph {3.}
\begin {enumerate}
\item [\emph {a.}]
Donner une matrice  $P' \in \Mn(\gA)$  telle que $P'\tra {M} = C_fP'$ et $\det(P')$ \ndzz.
\item [\emph {b.}]
On définit  $\vep'_1$, \ldots, $\vep'_n \in \gB$ par
l'\egt $[\,\vep'_1\; \ldots\; \vep'_n\,] = [\,1\; x\; \ldots\; x^{n-1}\,]\cdot P'$,
\hbox{et $\fb' = \gA\vep'_1 \oplus \cdots \oplus \gA\vep'_n$}.
Montrer que:

\snic {
\cmatrix {\vep'_1\cr \vdots\cr \vep'_n\cr} [\,\vep_1\; \ldots\; \vep_n\,] = 
\det(P) \Adj(x\In-M).
}

\snii 
En déduire l'\egt d'\ids $\fb'\fb = \det(P)\cD_{n-1}(x\In-M)$.

\item [\emph {c.}]
On suppose $x\In-M$ de rang $n-1$. \\
Alors $\fb$ est un
\id\ivz, et plus \prmt $\fb\fb'=\gen{\det(P)}$. \\
La suite $\gB^n\vvvvers{x\In-M} \gB^n
\vvvvvers{\Adj(x\In-M)} \gB^n$ est exacte et  $\fb \simeq \Im\Adj(x\In-M)$,
cette dernière matrice étant de rang~1.
\end {enumerate}

\noindent \emph {6.} (Exemple)
Soient $a \in \gA$ un \elt\ndzz, $f = X^n + \cdots + ab_0\in \AX$ un \pol\unt
de degré $n$ dont le \coe constant est multiple de $a$ et $\fb = \gen {a,x}
\subseteq \gB = \gA[x]$. Montrer que $\fb$ est un \Amo libre de rang $n$, en
expliciter une $\gA$-base, les matrices $M$, $P$, $P'$ correspondantes
(notations de la question précédente) ainsi que l'\id $\fb'$.  Si $a$ et
$b_0$ sont \comz, montrer \hbox{que $1 \in \cD_{n-1}(x\In-M)$}, donc $\fb\fb' = a\gB$
et l'\id $\fb$ est \ivz.

\snii \emph {7.} (Eisenstein)
Soient $a \in \gA$ un \elt\ndzz, $f = X^n + \sum_{i=0}^{n-1} a_iX^i
\in \AX$ un \pol \gui{Eisenstein en $a$}, i.e.  tel que

\snic{a_i \equiv 0
\bmod a$ pour $i\in\lrb{0..n-1}\;$ et $\;a_0 = b_0a$ avec $b_0$ \iv modulo $a.}

\snii 
Soit  $\fb := \gen {a,x}
\subseteq \gB = \gA[x]$. Montrer, pour $k \in \lrb{1..n}$, que $\fb^k$ est \hbox{un
\Amoz} libre de rang $n$. Plus \prmtz:

\snic {
\fb^k = \bigoplus_{i=0}^{k-1} \gA ax^i \oplus \bigoplus_{j=k}^{n-1} \gA x^j
= \gen {a, x^k}.
}

\snii  En particulier $\fb^n = a\gB$ et $\fb$ est \ivz.  \\
Montrer
\egmt que les \Amos $\gB\sur\fb$ et $\aqo{\gA}{a}$ sont isomorphes.

}

\end {problem}
%--- end -problem-----------------------------------------

}
%: fin sinotenglish

}% fin des exos

%:  Solutions
\penalty-2500
\sol{

%%%%%%%%%%%%%%%%%%%%%%%%%%%%%%%%%%%%%%%%%%%%%%%%%%%%%%%%%%%%%%%%%%%%%%%%%%%

\exer{exoGensIdealEnsFini}
\emph{1.}
Interpolation de Lagrange:
${
Q = \sum_{\xi \in U} \left(\prod_{\zeta \in U \setminus \{\xi\}} 
{x_n-\zeta \over \xi-\zeta}\right) Q_\xi
}.$
\\  \emph{2.}
Supposons que chaque $\fa(V_\xi) \subset \gK[x_1, \ldots, x_{n-1}]$ (pour $\xi
\in \pi_n(V)$) soit engendré par~$m$ \polsz:

\snic {
\fa(V_\xi) =  \gen{f^\xi_j, j \in \lrbm}, \quad 
f^\xi_j \in \gK[x_1, \ldots, x_{n-1}].
}

%\sni
D'après le point \emph{1}, il existe $f_j \in \Kux$ 
vérifiant $f_j(x_1, \ldots, x_{n-1},\xi) = f_j^\xi$ pour tout~$\xi \in \pi_n(V)$. On montre alors en s'appuyant sur le point \emph{1}
que:

\snic {
\fa(V) = \langle P, f_1, \ldots, f_m \rangle \qquad
\hbox {avec} \  P = \prod_{\xi \in \pi_n(V)} (x_n - \xi).
}

%\sni
On conclut par \recu sur~$n$.

%%%%%%%%%%%%%%%%%%%%%%%%%%%%%%%%%%%%%%%%%%%%%%%%%%%%%%%%%%%%%%%%%%%%%%%%%%%

\exer{exothSymEl}
 \emph{4.} 
Considérons l'anneau de \pols $\gB=\gA[s_1,\ldots,s_n]$
où les $s_i$ sont des \idtrsz, puis le \pol 
%
%\snic{
$
f(t)=t^n+\sum_{k=1}^{n}(-1)^ks_kt^{n-k}\in\gB[t].$
%}
\\ 
Considérons aussi l'\adu

\snic{\gC=\Adu_{\gB,f}=\Bxn=\Axn,}

%\sni
avec dans $\gC[t]$, l'\egt $f(t)=\prod_{i=1}^n(t-x_i)$.
\\
 Soient $\rho:\AXn\to\Axn$ et  $\varphi:\gA[s_1,\ldots,s_n]\to\gA[S_1,\ldots,S_n]$  les \homos d'\evn $X_i\mapsto x_i$ et $s_i\mapsto S_i$. 

\vspace{-7pt}
\Deuxcol{.7}{.15}
% premiere colonne
{
On a clairement $\rho(S_i)=s_i$. Donc, en notant  $\rho_1$ la restriction de $\rho$ à $\gA[\uS]$ et $\gA[\und{s}]%=\gB
$, on a $\varphi\circ \rho_1=\Id_{\gA[\uS]}$ \hbox{et $\rho_1\circ \varphi=\Id_{\gA[\und{s}]}$}.
Ceci montre que les $S_i$ sont \agqt indépendants sur $\gA$
et l'on peut identifier $\gA[\uS]$ et~$\gA[\und{s}]=\gB$.
}
% deuxieme colonne
{
\xymatrix @R = 10pt{
\gA[\uX]\ar@<0.7ex>[r]^\rho &
\gA[\ux]\ar@<0.7ex>[l]^\psi 
\\ 
\gA[\uS]\ar@{^{(}->}[u]\ar@<0.7ex>[r]^{\rho_1} &
   \ar@<0.7ex>[l]^\varphi \gA[\und s]\ar@{^{(}->}[u]
\\ 
}
}

Par la \prt \uvle de l'\aduz, il existe un 
(unique) $\gB$-\homo $\psi:\gC\to\AuX$
qui envoie $x_i$ sur $X_i$. Il s'ensuit que $\rho$ et $\psi$ sont deux \isos
réciproques l'un de l'autre. Ainsi les $x_i$ sont \agqt indépendants sur $\gA$
et $\AuX$ est libre de rang $n!$ sur $\gA[\uS]=\gB$, avec la base prescrite.
\\
 NB: cette \dem ne semble pas pouvoir donner  de manière simple le fait que les \pols \smqs de $\AuX$ sont dans $\gA[\uS]$.

%%%%%%%%%%%%%%%%%%%%%%%%%%%%%%%%%%%%%%%%%%%%%%%%%%%%%%%%%%%%%%%%%%%%%%%%%%%
\exer{exoPolSym1}
\emph{1.}
Soit $f = (X_1^3 + X_2^3 + \cdots + X_n^3) - (S_1^3 - 3S_2S_2 + 3S_3)$.
C'est
un \pol \smq \hmgz,  donc $f =
g(S_1, \ldots, S_n)$ où $g = g(Y_1, \ldots, Y_n)$ est \hmg en poids, de poids 3
pour le poids $\alpha_1 + 2\alpha_2 + \cdots + n\alpha_n$.\\
L'\egt $\alpha_1 + 2\alpha_2 + \cdots + n\alpha_n = 3$ implique $\alpha_i =
0$ pour~$i > 3$, donc $g$ ne dépend que de $Y_1$, $Y_2$, $Y_3$, disons $g = g(Y_1,
Y_2, Y_3)$. Dans l'\egt 

\snic{(X_1^3 + X_2^3 + \cdots + X_n^3) - (S_1^3 -
3S_2S_2 + 3S_3) = g(S_1, S_2, S_3),}

%\sni
on réalise $X_i := 0$ pour~$i > 3$; on
obtient $g(S'_1, S'_2, S'_3) = 0$ où $S'_1, S'_2, S'_3$ sont les  
fonctions \smqs \elrs de $X_1, X_2, X_3$.  On en déduit
que $g = 0$ puis~$f = 0$.
 
\emph{2.}
Pour le premier, on peut supposer $n = 3$; on trouve
$S_1S_2 - 3S_3$. Pour les deux autres qui sont \smqs
\hmgs de degré $4$, on travaille avec 4 \idtrs
et l'on obtient $S_1^2 S_2 - 2S_2^2 - S_1S_3 + 4S_4$ et
$S_2^2 - 2S_1S_3 + 2S_4$.

 \emph{3.}
Soient $n > d$ et $f(X_1, \ldots, X_{n})$ un \pol \smq
\hmg de degré $d$. Soit $h\in\gA[X_1,\ldots,X_d]=f(X_1,\ldots,X_d,0,\ldots,0)$.
Si $h=0$, alors $f = 0$.
\\
On peut traduire ce résultat en disant que l'on a des isomorphismes de
$\gA$-modules au niveau des composantes \smqs \hmgs de degré
$d$:

\snuc{\cdots \to
\gA[X_1, \ldots, X_{d+2}]^{\rm sym.}_d
\vvvers{X_{d+2} := 0}
\gA[X_1, \ldots, X_{d+1}]^{\rm sym.}_d
\vvvers{X_{d+1} := 0}
\gA[X_1, \ldots, X_{d}]^{\rm sym.}_d
\;.}

%%%%%%%%%%%%%%%%%%%%%%%%%%%%%%%%%%%%%%%%%%%%%%%%%%%%%%%%%%%%%%%%%%%%%%%%%%%
\exer{exoMcCoy}
On pose $\gA = \ZZ[U,V]/\gen {U^2, V^2} = \ZZ[u,v] = \ZZ \oplus \ZZ u \oplus
\ZZ v \oplus \ZZ uv$. \\
\emph{a.} On prend $f = uT + v$ donc $\fc = \gen {u,v}$.  On a alors

\snuc{\Ann(u) = \gA u$, $\Ann(v) = \gA v$, $\Ann(\fc) = \Ann(u) \cap \Ann(v) = \gA
uv$ et $\rD\big(\Ann(\fc)\big) = \fc.}

\sni
\emph{b.} On pose $g = uT - v$.  On a $fg = 0$ mais $g
\notin \Ann(\fc)[T]$; on a $u \in \rD\big(\Ann(\fc)\big)$ mais $u \notin
\Ann_{\gA[T]}(f)$ (idem pour $v$).

%%%%%%%%%%%%%%%%%%%%%%%%%%%%%%%%%%%%%%%%%%%%%%%%%%%%%%%%%%%%%%%%%%%%%%%%%%%
\exer{exoModCauBase}
{Il suffit de prouver le point \emph{1.}
On a $$f(T)=f(X_1)+(T-X_1)f_2(X_1,T)$$ par \dfn de $f_1=f$ et $f_2$.
De même $$f_2(X_1,T)=f_2(X_1,X_2)+(T-X_2)f_3(X_1,X_2,T)$$ par \dfn de
$f_3$. Donc $$f(T)=
%f(X_1)+(T-X_1)(f_2(X_1,X_2)+(T-X_2)f_3(X_1,X_2,T))=
f(X_1)+(T-X_1)f_2(X_1,X_2)+(T-X_1)(T-X_2)f_3(X_1,X_2,T).$$
On continue jusqu'à
%--------------------begin array---------------
$$\arraycolsep2pt
\begin{array}{rcl}
f_{n-1}(X_1,\ldots ,X_{n-2},T)&=&f_{n-1}(X_1,\ldots ,X_{n-2},X_{n-1})
\;+\\[1mm]
&   &
(T-X_{n-1})f_n(X_1,\ldots ,X_{n-1},T),
\end{array}
$$
%---------------------end array--------------
ce qui donne
%--------------------begin array---------------
$$
\arraycolsep2pt\begin{array}{rcl}
f(T)&   =&f_1(X_1)+(T-X_1)f_2(X_1,X_2)+(T-X_1)(T-X_2)f_3(X_1,X_2,X_3)
\\[1mm]
&   &
+\,\cdots\,+\,(T-X_1)\cdots (T-X_{n-1})f_n(X_1,\ldots ,X_{n-1},T).
\end{array}
$$
%---------------------end array--------------
Enfin $f_n(X_1,\ldots ,X_{n-1},T)$ est \mon de degré 1 en $T$
donc

\snic{f_n(X_1,\ldots ,X_{n-1},T)=f_n(X_1,\ldots ,X_{n-1},X_n)+(T-X_n).}

%\sni
Notez que ceci prouve en particulier que $f_n=S_1-s_1$.
}

%%%%%%%%%%%%%%%%%%%%%%%%%%%%%%%%%%%%%%%%%%%%%%%%%%%%%%%%%%%%%%%%%%%%%%%%%%%
\exer{exoPrimePowerRoot}
Soit $f \in \gA[X]$ \mon de degré $d$, avec $f
\divi X^p - a$ et $1 \le d \le p-1$. Dans un anneau $\gB \supseteq \gA$, on 
écrit $f(X) = \prod_{i=1}^d (X - \alpha_i)$, donc $\alpha_i^p = a$
et~$\prod_i \alpha_i = b$ avec~$b = (-1)^d f(0) \in \gA$. En élevant
à la puissance $p$, $a^d = b^p$. Mais $\pgcd(d,p) = 1$, donc~$1 = ud + vp$, puis $a = a^{ud} a^{vp} = (b^u a^v)^p$.

%%%%%%%%%%%%%%%%%%%%%%%%%%%%%%%%%%%%%%%%%%%%%%%%%%%%%%%%%%%%%%%%%%%%%%%%%%%
\exer{exoPrincipeIdentitesAlgebriques}
Notons $e_{ij}$ la matrice de $\Mn(\gA)$ ayant un seul \coe
non nul, le \coe en position $(i,j)$, égal à $1$.
Le module $S_n(\gA)$ est libre et une base est formée par les $e_{ii}$ pour
$i \in \lrbn$ et les $e_{ij} + e_{ji}$ pour $1 \leq i < j < n$.  Il suffit de
traiter le cas où %$A$ est diagonale,
$A = \Diag(\lambda_1, \ldots,
\lambda_n)$.
Alors, $\varphi_A= \Diag(\lambda_1^{2}, \ldots,
\lambda_n^{2})$, \hbox{et
$\varphi_A(e_{ij}+e_{ji})=\lambda_i\lambda_j (e_{ij}+e_{ji})$}.
D'où $\det(\varphi_A)=(\det A)^{n+1}$.

%%%%%%%%%%%%%%%%%%%%%%%%%%%%%%%%%%%%%%%%%%%%%%%%%%%%%%%%%%%%%%%%%%%%%%%%%%%
\exer{exoFreeFracTransfert}
Soit $\ue = (e_1, \ldots, e_n)$ une base de $\gB/\gA$. Il est clair
que $\ue$ est une famille $\gK$-libre. Soit $x = b/b' \in \gL$
avec $b \in \gB$, $b' \in \gB \setminus \{0\}$; on écrit:

\snic{
x = (b\wi{b'}) / (b'\wi{b'}) =
b\wi{b'} / \rN\iBA(b') \in
\gK e_1 + \cdots + \gK e_n.}

%%%%%%%%%%%%%%%%%%%%%%%%%%%%%%%%%%%%%%%%%%%%%%%%%%%%%%%%%%%%%%%%%%%%%%%%%%%

\exer{exoSommesNewton}
\emph{1.} Il suffit de le prouver pour la matrice \gnq $(a_{ij})_{i,j\in\lrbn}$ à \coes
dans $\gA=\ZZ[(a_{ij})_{i,j\in\lrbn}]$:
cette matrice  est  \dig dans un suranneau de~$\gA$.

\emph{2.} Résulte \imdt du \emph{1}.

%%%%%%%%%%%%%%%%%%%%%%%%%%%%%%%%%%%%%%%%%%%%%%%%%%%%%%%%%%%%%%%%%%%%%%%%%%%
\exer{exoCorpsFiniEltPrimitif}
\emph{1.} On a $\#R \le \sum_{d=1}^{n-1} q^d < 1 + q + \cdots + q^{n-1} 
= {q^n - 1\over q-1}$.
\\
 A fortiori, $\#R < q^n-1 < q^n$. Soit $x \in \gL\setminus R$
et $d = \dex{\gK[x] : \gK}$. On a $x^{q^d} = x$, et comme $x \notin R$,
c'est que $d = n$.

 \emph{2.}
Le \pol cyclotomique est \ird dans $\FF_2[X]$. En effet, le seul
\pol irréductible de degré $2$ de $\FF_2[X]$ est $X^2 + X + 1$, $\Phi_5$
est sans racine dans~$\FF_2$, et $\Phi_5 \ne (X^2 + X + 1)^2$.  On a $\#\gL =
2^4 = 16$, $\#\gL\eti = 15$, mais $x^5 = 1$.

 \emph{3.}
Soit $\sigma : \gL \to \gL$ l'\auto de Frobenius de $\gL/\gK$, i.e.
$\sigma(x) = x^q$. On vérifie facilement que $\gL = \gK[x]$ \ssi
les $\sigma^i(x)$, $i \in \lrb{0..n-1}$, sont deux à deux distincts.
Cette condition équivaut à $\,\sigma^k(x) = x \;\Rightarrow\; k \equiv 0
\bmod n\,$, i.e. $\,x^{q^k} = x \;\Rightarrow\; k \equiv 0 \bmod n\,$. Mais

\snic {
x^{q^k} = x \iff x^{q^k-1} = 1 \iff o(x) \divi q^k - 1 \iff
q^k \equiv 1 \bmod o(x).
}

%\sni
On en déduit, pour $x \in \gL\eti$, que $\gL=\gK[x]$ \ssi
l'ordre de $q$ dans le groupe des \ivs modulo $o(x)$
est exactement $n$.

%%%%%%%%%%%%%%%%%%%%%%%%%%%%%%%%%%%%%%%%%%%%%%%%%%%%%%%%%%%%%%%%%%%%%%%%%%%

\exer{exoPpcmPolsSeparables}
 \emph{1.}
On a $\gen {g,g'} \gen {g,h} \subseteq \gen {g, g'h} = \gen {g, g'h + gh'}$.
\\
De même, $\gen {h,h'} \gen {g,h} \subseteq \gen {h, g'h + gh'}$.
En faisant le produit, il vient:

\snic {
\gen {g,g'} \gen {h,h'} \gen {g,h}^2 \subseteq 
\gen {g, g'h + h'g} \gen {h, g'h + h'g} \subseteq \gen {gh, g'h + h'g}.
} 

%\sni
Pour le deuxième point de la question on applique le résultat que l'on vient d'établir et le fait 
\ref{factDiscUnit}. NB: cela résulte aussi de l'\eqnz~(\ref{eqfactDiscProd}),  fait~\ref{factDiscProd}.

 \emph{2.} Il suffit de traiter le cas de deux \pols \spls
$f$, $g \in \gA[T]$. 
\\
Soit
$h = \pgcd(f,g)$. On a $f = hf_1$, $g = hg_1$, avec
$\pgcd(f_1, g_1) = 1$.% et $\ppcm(f, g) = hf_1g_1$.
\\
Puisque $g$ est \splz, $\pgcd(h, g_1) = 1$, donc $\pgcd(hf_1, g_1) = 1=\pgcd(f, g_1)$.  \\
Les  \polsz~$f$,~$g_1$ sont \splsz,
\comz, donc leur produit $\ppcm(f, g)$ est \splz.
%%%%%%%%%%%%%%%%%%%%%%%%%%%%%%%%%%%%%%%%%%%%%%%%%%%%%%%%%%%%%%%%%%%%%%%%%%%

\exer{exolemSousLibre} 
 \emph{1} et \emph{2.} Ce sont des cas particuliers de ce qui est affirmé dans le fait \ref{fact.idd.sousmod}.

 \emph{3.} On suppose $L=\Ae m$. Si $A\in\Mm(\gA)$ est une matrice dont les colonnes forment une base de $F$, elle est injective et son \deter est \ndzz. Si $B$ est une matrice correspondant à l'inclusion $F\subseteq E$, on a 

\snic{\idg{L:F}=\gen{\det A}$, $\;\idg{F:E}=\cD_m(B)\;$ et $\;\idg{L:E}=\cD_m(AB),}

%\sni
d'où l'\egt souhaitée.

 \emph{4.} 
On a $\idg{N:\delta N}  = \gen {\delta^n}$.
On a aussi $\idg{N:\delta M} = \delta^{n-1} \gen {\delta, \an}$: prendre pour \sgr de $\delta M$ la famille $\delta e_1, \dots, \delta
e_n, \delta z$ où $e_1, \ldots, e_n$ est une base de $N$ (on utilise
$M = N + \gA z$), et calculer l'\idd d'ordre $n$ d'une matrice de type
suivant (pour $n = 3$): 
$\cmatrix {
\delta & 0     & 0      & a_1 \cr
0      &\delta & 0      & a_2 \cr
0      &0      & \delta & a_3 \cr
}.$
\\
Alors:

\snic {
\idg{N:\delta N} \,=\, \idg{N:\delta M}\, \idg{\delta M: \delta N} \, =\,  \idg{N:\delta M}\,\idg{M:N}, 
}

%\sni
i.e. $\gen {\delta^n} = \idg{M:N}\; \delta^{n-1} \gen {\delta, \an}$.
\\
En simplifiant par $\delta^{n-1}$ on obtient l'\egt $\gen{\delta}=d\,\gen {\delta, \an}$.

%%%%%%%%%%%%%%%%%%%%%%%%%%%%%%%%%%%%%%%%%%%%%%%%%%%%%%%%%%%%%%%%%%%%%%%%%%%

\exer{exoDecompIdeal} \emph{1.} Si $\fa\fa'=a\gA$ avec $a$ \ndzz, alors
$\fb\subseteq\fa$ équivaut à $\fb\fa'\subseteq a\gA$. On note que le test fournit
un \itf $\fc=\fb\fa'/a$ tel que $\fa\fc=\fb$ en cas de réponse positive, et un \elt $b\notin\fa$ parmi les \gtrs de $\fb$ en cas de réponse négative.

 \emph{2.}
Il est clair que les $\fq_i$ sont \ivs (et donc \tfz).
\\
On fait les tests $\fb\subseteq\fq_i$. Si une réponse est positive,
par exemple $\fb\subseteq\fq_1$, on écrit~$\fc\fq_1=\fb$, d'où 
$\fq_2\cdots\fq_n\subseteq\fc$, et l'on termine par \recuz.
\\
 Si tous les tests sont négatifs, on a des $x_i\in\fb$ et $y_i\in\gA$
tels que $1-x_iy_i\in\fq_i$ (on suppose ici que les quotients~$\gA\sur{\fq_i}$
sont des \cdisz), d'où en faisant le produit $1-b\in \fq_1\cdots\fq_n\subseteq\fb$
avec $b\in\fb$, donc $1\in\fb$.
\\
Voyons enfin la question de l'unicité. Supposons que $\fb=\fq_1\cdots\fq_k$.
\\
Il suffit de montrer que si un \idema $\fq$ \tf contient $\fb$, il est égal à l'un des $\fq_i$ ($i\in\lrbk$). 
\\
 Puisque l'on peut tester $\fq\subseteq\fq_i$, si chacun des tests était négatif
on aurait explicitement $1\in\fq+\fq_i$ pour chaque $i$ et donc $1\in\fq+\fb$.
\\
NB: si l'on ne suppose pas $\fb$ \tf et $\gA$ à \dve explicite,
la \dem du petit \tho de Kummer nécessiterait que l'on sache au moins tester~$\fq\subseteq\fb$ pour tout
\gui{sous-produit} $\fq$ de  $\fq_1\cdots\fq_n$.

%%%%%%%%%%%%%%%%%%%%%%%%%%%%%%%%%%%%%%%%%%%%%%%%%%%%%%%%%%%%%%%%%%%%%%%%%%%%%%

%%%%%%%%%%%%%%%%%%%%%%%%%%%%%%%%%%%%%%%%%%%%%%%%%%%%%%%%%%%%%%%%%%%%%%%%%%%
\exer{exoRabinovitchTrick} 
Supposons $x \in \sqrt\fa$; comme $\fa\subseteq\fb$, dans $\gA[T]\sur\fb$, $\ov
x$ est nilpotent et inversible (puisque $\ov x \ov T = 1$), donc
$\gA[T]\sur\fb$ est l'anneau nul, i.e. $1 \in \fb$.
\\
Inversement, supposons $1 \in \fb$ et raisonnons dans l'anneau $\gA[T]/\fa[T] =
(\gA\sur\fa)[T]$. Puisque $1 \in \fb$, $1-xT$ est \iv dans cet anneau,
donc $x$ est nilpotent dans~$\gA\sur\fa$, i.e. $x \in \sqrt\fa$.

%%%%%%%%%%%%%%%%%%%%%%%%%%%%%%%%%%%%%%%%%%%%%%%%%%%%%%%%%%%%%%%%%%%%%%%%%%%

\exer{exoDunford} \emph{(Décomposition de Jordan-Chevalley-Dunford)} \\ 
 \emph{Existence.}
On cherche un zéro $D$ de $f$, \gui{voisin de $M$}, (i.e., avec $M-D$ nilpotent), dans l'anneau commutatif $\gK[M]$.
On a par hypothèse $f(M)^{k}=0$ pour un $k\leq n$, et si $uf^{k}+vf'=1$, on obtient
${v(M)f'(M)=\In.}$

En conséquence, la méthode de Newton\imN, démarrant \hbox{avec $x_0=M$}, 
donne la solution dans $\gK[M]$ en $\lceil\log_2(k)\rceil$ itérations.

\emph{Unicité.}
La solution est unique, sous la condition $f(D)=0$, dans tout anneau commutatif contenant~$\gK[M]$,
par exemple dans $\gK[M,N]$ si le couple $(D,N)$ résout le \pb posé.\\
Lorsque l'on suppose seulement que le \polmin de $D$ est \splz, l'unicité est plus délicate.\\
Un solution serait de démontrer directement que le \polcar de~$D$ est \ncrt égal à celui de $M$, mais ce n'est pas si simple\footnote{En caractéristique nulle, une astuce consiste à récupérer le \polcar d'une matrice $A$ à partir des $\Tr(A^{k})$ en suivant la méthode de Le Verrier.}.\\
Appelons $(D_1,N_1)$ la solution dans $\gK[M]$ donnée par la méthode de Newton.
Puisque $D$ et $N$ commutent, elles commutent avec $M=D+N$ et donc avec~$D_1$ et~$N_1$ car ils appartiennent à $\gK[M]$. On en déduit que  $D-D_1$ est nilpotente car elle est égale à $N_1-N$ avec~$N$ et~$N_1$ nilpotentes qui commutent. Or l'\algz~$\gK[D,D_1]$ est étale
d'après le \thref{corlemEtaleEtage}, donc elle est réduite, et $D=D_1$.

%%%%%%%%%%%%%%%%%%%%%%%%%%%%%%%%%%%%%%%%%%%%%%%%%%%%%%%%%%%%%%%%%%%%%%%%%%%

\exer{exoDunfordBis} 
On a $\gB = \gA[x] = \gA \oplus \gA x$ avec $x$ \splz{ment} entier sur $\gA$.
Notons $\,z \mapsto \wi z\;$ l'\auto de l'\Alg $\gB$ qui échange $x$ et $-b-x$.
\\
Pour~$z \in \gB$, on a ${\rm C}\iBA(z)(T) = (T - z)(T - \wi z)$.\\ 
Ainsi
${\rm C}\iBA(ax)(T) = T^2 + abT + a^2c$, et son \discri est égal à $a^2\Delta$.
\\
Soit $\vep\in \gA$ nilpotent non nul et posons $y = (\vep-1)x$. Alors, $y$
est \spbz{ment} entier sur $\gA$ car $(\vep-1)^2\Delta$ est \ivz.
En outre, l'\elt $z=x + y = \vep x$ est nilpotent non nul. Supposons que $\vep^2=0$ et soit $g\in\AX$ un \polu qui annule $z$, on va montrer que $g$ n'est pas \splz.\\
\'Ecrivons $g(X)=u+vX+X^2h(X)$, alors $z^2=0$, donc $u+vz=0$.\\
Puisque $\gB=\gA\oplus\gA x$,
on obtient $u=v\vep=0$, puis \hbox{$g(X)=X\ell(X)$} \hbox{avec $\ell(0)=v$} non \iv (sinon, $\vep=0$). Enfin, $\disc(g)=\disc(\ell)\ \Res(X,\ell)^2=\disc(\ell)\ v^2$ est non \ivz.

%: sinotenglish
\sinotenglish{
\exer{exoAdj2Minors}
\\
On va se ramener à $i=j=n-1$, $i'=j'=n$ à l'aide de matrices de
permutations~$P_\sigma$. On rappelle tout d'abord que:

\snic {
(P_\tau A P_\sigma)_{i,j} = A_{\tau^{-1}(i),\sigma(j)}, \qquad
\widetilde {P_\tau AP_\sigma} = 
\vep(\sigma)\vep(\tau) P_{\sigma^{-1}} \widetilde A P_{\tau^{-1}}
}

\snii
Notons $G(A,i,j,i',j')$ le membre gauche de l'\egt à démontrer. Alors:

\snic {
G(P_\tau AP_\sigma,i,j,i',j') = G(A,\sigma(i),\tau^{-1}(j),\sigma(i'),\tau^{-1}(j'))
}

\snii
On a la même \prt d'invariance pour le membre droit. On peut donc supposer,
quitte à remplacer $A$ par $P_\tau AP_\sigma$, que $i=j=n-1$, $i'=j'=n$. On
applique alors le point \emph {7} du
lemme III-\ref{lemPrincipeIdentitesAlgebriques}).

%%%%%%%%%%%%%%%%%%%%%%%%%%%%%%%%%%%%%%%%%%%%%%%%%%%%%%%%%%%%%%%%%%%%%%%%%%% 
\exer{exoResInversible}

\noindent   
\emph{\ref{i2exoResInversible}} $\Rightarrow$ \emph{\ref{i1exoResInversible}.}
Si $f_p=1$ ou $g_q=1$ l'implication est claire d'après le lemme d'\eli de base.
L'implication est donc valide après \lon en $f_p$ ou $g_q$. 
On conclut par le \plg de base.

\snii  
\emph{\ref{i1exoResInversible}} $\Rightarrow$ \emph{\ref{i2exoResInversible}.} En développant le \deter de la matrice de Sylvester selon la première colonne
on obtient que $\Res_{X}(f,p,g,q)\in \gen{f_p,g_q}$.

\noindent 
Par ailleurs d'après le fait \ref{fact1Res}, $\Res_{X}(f,p,g,q)$ s'écrit toujours sous forme 

\snic{u(X)f(X)+v(X)g(X).}

\snii  \emph{\ref{i1exoResInversible}} $\Rightarrow$ \emph{\ref{i4exoResInversible}.} 
%Par homogénésiation, une \egt $u(X)f(X)+v(X)g(X)=1$ donne une \egt
%\hmg $U(X,Y)f(X,Y)+v(X,Y)g(X,Y)=Y^k$ (en fait $k=p+q-1$).
%En considérant les \pols réciproques $F(1,Y)$
%et $G(1,Y)$ qui donnent le même résultant au signe près,
%on obtient de la même manière que $X^k\in\gen{F,G}$.
%
%\snii  \emph{\ref{i3exoResInversible}.} $\Rightarrow$ \emph{\ref{i2exoResInversible}.} 
%Par déshomogénésiation, une \egt 
%
%\snic{U(X,Y)f(X,Y)+v(X,Y)g(X,Y)=Y^k}
%
%\snii
%donne une \egt $u(X)f(X)+v(X)g(X)=1$.
%A fortiori $1\in\gen{f(0),g(0)}$.\\
%En appliquant cette dernière remarque aux \pols réciproques 
%on obtient l'appartenance $1\in\gen{f_p,g_q}$.
La matrice de Sylvester $\Syl_X(f,p,g,q)$ 
\paref{SylvMat}
%page~119 
peut être vue comme la matrice dont les lignes
sont les \coos des \pols 

\snic{X^{q-1}F$, $\ldots $, $XY^{q-2}F$, $Y^{q-1}F$, $X^{p-1}G$, $\ldots $, $XY^{p-2}G$, $Y^{p-1}G}

\snii
sur la base

\snic{(X^{p+q-1},X^{p+q-2}Y,\ldots ,XY^{p+q-2},X^{p+q-1})}

du module des \pogs de degré $p+q-1$. Si son \deter est \ivz, c'est que la 
matrice de l'application  $(U,V)\mapsto UF+VG$ (pour les modules convenables à la source et au but) est surjective.
%%%%%%%%%%%%%%%%%%%%%%%%%%%%%%%%%%%%%%%%%%%%%%%%%%%%%%%%%%%%%%%%%%%%%%%%%%%

%: solution exer{exoPgcdNst}
\exer{exoPgcdNst} ~\\
\emph{1.} Si $n=\deg(f)=1$ ou $0$, on a $f_1=f$ et les conclusions sont satisfaites. Si $n\geq 2$ on raisonne comme suit. 
On suppose que $f=\prod_k(x-a_k)^{m_k}$,
avec les $a_k-a_\ell$ \ivs si $k\neq \ell$.
On écrit $f=(x-a_k)^{m_k}f_k$ et l'on a $\gen{f_k,x-a_k}=\gen{1}$.
\hbox{On a:}
$$\preskip.3em \postskip.0em
{f'=m_k(x-a_k)^{m_k-1} g_k+ (x-a_k)^{m_k} g'_k =(x-a_k)^{m_k-1}\,u_k.}
$$
avec
$$\preskip-.2em \postskip.4em
u_k=m_kg_k+ (x-a_k) g'_k.
$$
Si $m_k=0$ dans $\gK$, on obtient que le pgcd $h$ de $f$ et $f'$ est divisible par $(x-a_k)^{m_k}$. \\
Si $m_k\in\gK\eti$, alors $\gen{u_k,x-a_k}=\gen{1}$, et l'on obtient  que $f'$ est divisible par $(x-a_k)^{m_k-1}$ mais pas par $(x-a_k)^{m_k}$.
En fin de compte, on obtient l'\egt 
$$\preskip.4em \postskip.4em
{f_1=\prod\nolimits_{k:m_k\in\gK\eti}(x-a_k).}
$$
Dans tous les cas, $h$ est \splz. \\
Et si tous les $m_k\in\gK\eti$,
(par exemple \hbox{si $n!\in\gK\eti$}), alors $f$ divise $f_1^{n}.$

\emph{2.}  Cette question ne se comprend que d'un point de vue \cofz, car en \clama
tout corps possède une clôture \agqz, et il suffit alors de se reporter au point \emph{1.}

Sans doute il faut attendre d'avoir lu le chapitre VII pour se convaincre que l'on peut toujours \gui{faire comme si} l'on disposait d'un \cdr
pour le \pol $f$. On considère d'abord les zéros $x_k$ dans l'\adu $\gA$ de $f$.
Si $x_1-x_2\in\Ati$, $f$ est \splz, $h=1$ \hbox{et $f_1=f$}. \\
Sinon, on remplace~$\gA$ par un quotient de Galois $\gB$ de $\gA$. Dans ce quotient $\gB$, on a par exemple $x_1=x_2$.
On considère ensuite $x_1-x_3$ dans $\gB$. S'il est nul ou \ivz, tout est OK (et on continue en comparant les autres paires de racines).
Sinon, il faut considèrer un quotient de Galois plus poussé.
En fin de compte, après avoir renuméroté les $x_i$ on est certain d'obtenir dans un quotient de Galois~$\gC$ de l'\adu une \egt 
$$\preskip.4em \postskip.4em
f(x)=\prod\nolimits_{k=1}^{\ell}(x-x_k)^{m_k}, \hbox{ avec les } x_k-x_j\in\gC\eti
\hbox{ pour } k\neq j.
$$
La \dem donnée au point \emph{1} fonctionne alors dans ce nouveau cadre.
On obtient en effet 
$$\preskip-.4em \postskip.4em
f_1=\prod\nolimits_{k:k\leq \ell,\,m_k\in\gK\eti}(x-x_k).
$$  
Puis  $\Res_x(f_1,f_1')$ peut être calculé dans~$\gC$, où il est égal
à un sous-produit de {$\pm\prod_{j,k:j<k\leq \ell}(x_j-x_k)^{2}$.} 
C'est donc un \elt
de $\gK\cap \gC\eti=\gK\eti$.\\
Et si tous les $m_k\in\gK\eti$,
(par exemple \hbox{si $n!\in\gK\eti$}), alors $f$ divise $f_1^{n}.$ 

\emph{3a.} Si $\gen{f,f'}=\gen{h}$, on a des \pols $u$, $v$, $f_2$ et ${f_1}$ tels que 

\snic{uf+vf'=h$, $hf_1=f$ et $hf_2=f'.}

\snii
Cela donne $h(uf_1+vf_2)=h$.
Puisque $h$ divise $f$, il est primitif donc \ndzz, \hbox{d'où $uf_1+vf_2=1$}.
Et l'\egt matricielle du point \emph{3a} est bien satisfaite. 

\snic{\cmatrix{u&v\cr-f_2&{f_1}}\cmatrix{f\cr f'}=\cmatrix{h\cr 0}.}

\smallskip
\emph{3b.}
On considère l'anneau $\ZZ[(c_i)_{i\in\lrbl}]$, où les $c_i$ sont 
d'une part des \idtrs que l'on prend pour les \coes  des \pols  \hbox{$f$, $h$, $u$, $v$, $f_2$} et~${f_1}$,
et d'autre part des \idtrs pour obtenir une \coli des \coes de $f$ égale à~$1$. 
On a choisi pour les \pols en $x$ les degrés formels correspondant aux \eqns que l'on a par hypothèse dans~$\kx$.
\\
On considère le \syp sur les \idtrs $(c_i)$ qui correspond aux \eqns suivantes  dans $\ZZ[(c_i)][x]$:
$$\preskip.4em \postskip.4em
{f \hbox{ est primitif},\; uf_1+vf_2=1,\; hf_1=f,\; hf_2=f'.}
$$
Soit alors $\fa$  l'\id de $\ZZ[(c_i)]$ engendré par ce \syp  de~$\ZZ[(c_i)]$. On obtient ainsi l'anneau \gui{\gnqz} de la situation considérée: 
$$\preskip.4em \postskip.4em
\fbox{$\gA=\ZZ[(c_i)]/\fa$.}
$$
Il est clair que tout se passe dans le sous-anneau $\gk'$ (de~$\gk$)   quotient de l'anneau générique~$\gA$, obtenu en spécialisant les \idtrs $c_i$ en leurs valeurs dans $\gk$.
\\
  Si l'on évalue cette situation dans un corps fini $\gF$, 
  i.e. si l'on considère un \homo $\gA\to \gF$, on a $1\in \gen{f_1(x),f_1'(x)}\subseteq \gF[x]$ en vertu du 
point~\emph{1.}
En effet, comme on a forcé $f(x)$  à être primitif, son image dans $\gF[x]$
est un \pol non nul, et l'on peut appliquer le point \emph{1},
en notant que  tout corps fini possède une clôture \agqz.  
Notons que si au contraire $f(x)$ était le \pol nul de $\gF[x]$, les \eqns que l'on impose n'interdiraient pas d'avoir $f_1=0$. 
\\
Notons $\fb=\fa+\gen{f_1(x),f_1'(x)}\subseteq \ZZ[x,(c_i)_{i\in\lrbl}]$.
On a donc obtenu que le \syp qui correspond à l'\id $\fb$ n'a de solution
dans aucun corps fini. 
\\
Par le \nst formel on en déduit que $1\in \fb$, 
ce qui veut aussi dire que $1\in \gen{f_1(x),f_1'(x)}\subseteq \gA[x]$.
Et ceci implique que $1\in \gen{f_1(x),f_1'(x)}\subseteq \gk[x]$, car
cette appartenance est déjà certifiée avec le sous-anneau $\gk'$ de $\gk$ qui est un quotient de~$\gA$.

Notez que l'on n'a pas besoin de démontrer le point \emph{2} pour obtenir le résultat \gnl du point \emph{3b} (qui contient le point \emph{2} comme cas particulier).

\emph{3c.} Voyons la dernière question: \emph{si en outre  $n!\in\gk\eti$, alors $f$ divise ${f_1}^{n}$?} Ici $n$ est a priori le degré formel de $f$,
qui peut être son vrai degré  si on le connaît.\\
On doit introduire une \idtr \sul $z$ pour l'inverse de $n!$.

\emph{Première solution partielle.}\\
Notons $R$ le reste de la division de ${f_1}^{n}$ par $f$ (que l'on suppose ici unitaire).
Alors on sait que pour tout zéro de l'\id $\fc=\fa+\gen{zn!-1}$ dans un corps fini,
les \coes de $R$ sont nuls.  
Le \nst formel nous dit alors que les \coes de $R$ sont dans le nilradical $\sqrt\fc$
de $\fc$. \\
En conclusion, dans un anneau $\gk$ tel que $n!\in\gk\eti$, si les hypothèses
sont satisfaites, et si $f$ est \untz, on peut affirmer que les \coes de $R$ sont nilpotents. Comme conséquence, une certaine puissance de $R={f_1}^{n}-fq$ est nulle,
et donc $f$ divise une puissance de $f_1$.

\emph{Deuxième solution partielle.} \\
On va obtenir la même conclusion finale sans supposer $f$ \untz.\\
On introduit une \idtr  $z$ et 
l'on considère l'\id

\snic{\fd=\fa+\gen{f(x),zn!-1}\subseteq \ZZ[x,z,(c_i)_{i\in\lrbl}].}

\snii
D'après le point \emph{1}, le \pol $f_1(x)$ s'annule en tout zéro de $\fd$ dans tout corps fini. Le \nst formel implique qu'une puissance de $f_1$ est dans $\fd$,
d'où il suit que dans $\gk[x]$, $f$ divise une puissance de $f_1$.  

\emph{4.} Merci \alec qui résoudra cette question, et qui éclaircira
complètement le dernier point de la question \emph{3.}

}
%: fin sinotenglish

%%%%%%%%%%%%%%%%%%%%%%%%%%%%%%%%%%%%%%%%%%%%%%%%%%%%%%%%%%%%%%%%%%%%%%%%%%%

%:  sol des pbs
\prob{exoDiscriminantsUtiles}~\\
\emph {1.}
Soit $f(X) = X^n + c = (X-x_1) \cdots (X-x_n)$. Alors, $f' = nX^{n-1}$ et

\snic {
\Res(f,f') =  f'(x_1) \cdots f'(x_n) =
n^n (x_1 \cdots x_n)^{n-1} = n^n \bigl((-1)^n c\bigr)^{n-1} =
n^n c^{n-1}.
}

%\sni

Variante:

\snic {
\Res(f', f) = n^n \Res(X^{n-1}, f) = n^n \prod_{i=1}^{n-1} f(0) 
= n^n c^{n-1}.
}

%%%%%%%%%%%%%%%%%%%%%%%%%%%%%%%%%%%%%%%%%%%%%

%\sni
\emph {2.}
Soit $f(X) = X^n + bX + c = (X - x_1) \cdots (X - x_n)$;

\snic {
\disc(f) = (-1)^{n(n-1)\over 2} \prod_{i=1}^n y_i 
\quad \hbox {avec}\quad  y_i = f'(x_i) = n x_i^{n-1} + b\,. 
}

%\sni
Pour calculer le produit des $y_i$, on calcule le produit $P$ des $x_i y_i$
(celui des $x_i$ vaut~$(-1)^n c$). On a $x_iy_i = nx_i^n + bx_i = ux_i + v$,
avec $u = (1-n)b$, $v = -nc$.  On utilise les fonctions \smqs
\elrs $S_j(x_1, \ldots, x_n)$ (presque toutes nulles):

\snic {
\prod_{i=1}^n (ux_i + v) = \sum_{j=0}^n u^j S_j(x_1, \ldots, x_n) v^{n-j}.
}

%\sni
Il vient:

\snic {
P = v^n + u^n S_n + u^{n-1} S_{n-1} v =
v^n + u^n (-1)^n c + u^{n-1} (-1)^{n-1} b v\,,
}

%\sni
\cadz, en remplaçant $u$ et $v$ par leurs valeurs:

\snic {\arraycolsep2pt
\begin {array}{rcl}
P &=&
 (-1)^n n^n c^n + (n-1)^n b^n c - n (n-1)^{n-1} b^n c \\
&=&
 (-1)^n n^n c^n + b^n c \bigl( (n-1)^n - n (n-1)^{n-1} \bigr) \\
&=& (-1)^n n^n c^n - b^n c  (n-1)^{n-1}.
\\
\end {array}
}

%\sni
En divisant par $(-1)^n c$, on obtient le produit des~$y_i$ puis la formule
annoncée.
%%%%%%%%%%%%%%%%%%%%%%%%%%%%%%%%%%%%%%%%%%%%%%
 
\emph {3.}
Laissé à la sagacité \dlec qui pourra consulter \cite{Swan62}.
 
\emph {4.}
En notant $\Delta_p = \disc(\Phi_p)$, on a l'\egt

\snic {
\disc(X^p - 1) = \Res(X-1,\Phi_p)^2 \disc(X-1)\Delta_p =
\Phi_p(1)^2 \Delta_p  = p^2 \Delta_p. 
}

%\sni

En utilisant  $\disc(X^n - 1) = (-1)^{n(n-1) \over 2} n^n (-1)^{n-1}$, 
on obtient:
$$
%\snic {
\Delta_2 = 1, \qquad  \Delta_p =(-1)^{p-1\over 2} p^{p-2}
\quad\hbox {pour $p \ge 3$}.
%}
$$

%\sni
\emph {5.}
Soit $q = p^{k-1}$; montrons d'abord que $r := \Res(X^q - 1, \Phi_{p^k}) =
p^q$.\\
  Avec $X^q - 1 = \prod_{i=1}^{q} (X-\zeta_i)$, on a $r = \prod_{i=1}^{q}
\Phi_{p^k}(\zeta_i)$. Par ailleurs:

\snic {
\Phi_{p^k}(X) = {Y^p - 1 \over Y-1} = Y^{p-1} + \cdots + Y + 1
\quad \hbox {avec} \quad Y = X^q.
}

%\sni
En faisant $X := \zeta_i$, on doit faire $Y := 1$, on obtient $\Phi_{p^k}(\zeta_i) = p$,
puis $r = p^q$.
\\
On note $D_k = \disc(X^{p^k} - 1)$. Puisque $X^{p^k} - 1 = (X^q - 1)
\Phi_{p^k}(X)$, on a:

\snic {
D_k = \Res(X^q-1, \Phi_{p^k})^2 D_{k-1} \disc(\Phi_{p^k}) =
p^{2q} D_{k-1} \disc(\Phi_{p^k}).
}

%\sni
On utilise $\disc(X^n - 1) = (-1)^{n(n-1) \over 2} n^n (-1)^{n-1}$
pour $n = p^k$ et $q$: 

\snic {
D_k/D_{k-1} = \varepsilon\, p^N, \quad
\varepsilon = \pm 1, \quad N = kp^k - (k-1)q = \big(k(p-1)+1\big)\,q\,.
}

%\sni

Pour obtenir $\disc(\Phi_{p^k})$, il faut diviser $D_k/D_{k-1}$ par
$p^{2q}$, ce qui remplace l'exposant~$N$ par $N - 2q =
(k(p-1)-1) q$. Quant au signe $\varepsilon$, pour $p$ impair:

\snic {
\varepsilon = (-1)^{{p^k - 1} \over 2} (-1)^{{q - 1} \over 2} =
(-1)^{{p^k - q} \over 2} = (-1)^{{p - 1} \over 2} 
.}

%\sni
Pour $p = 2$, $\varepsilon = 1$ pour $k \ge 3$ ou $k = 1$ et $\varepsilon 
= -1$ pour $k = 2$.
 
\emph {6.}
Si $n$ n'est pas la puissance d'un nombre premier, on peut écrire
$n = mp^k$ avec $p$ premier, $\pgcd(m,p) = 1$, $k \ge 1$ et $m \ge 2$.
Alors, $\Phi_n(X) = \Phi_m(X^{p^k}) /\Phi_m(X^{p^{k-1}})$, \egt dans
laquelle on réalise $X = 1$ pour obtenir $\Phi_n(1) = 1$.
Les autres points sont faciles.
 
\emph {7.}
Soient $f$, $g$ deux \polusz, avec $d = \deg f$, $e = \deg g$ et $d, e \ge
1$. Notons $\gA[x] = \aqo{\gA[X]}{f(X)}$, $\gA[y] = \aqo{\gA[Y]}{g(Y)}$.
Notons $f
\otimes g$ le \polcar de $x \otimes y$ dans $\gA[x]\te_\gA\gA[y] =
\aqo {\gA[X,Y]} {f(X),g(Y)}$. C'est un \polu de degré $d\,e$. Lorsque
 $f(X) = \prod_i (X - x_i)$, $g(Y) = \prod_j(Y - y_j)$,
on obtient $(f\otimes g)(T) = \prod_{i,j} (T - x_iy_j)$. On voit
facilement que

\snic {
\disc(f \otimes g) = \prod_{(i,j) < (i',j')} (x_iy_j - x_{i'}y_{j'})^2 =
\disc(f)^e \disc(g)^d f(0)^e g(0)^d \pi\,,
}

%\sni
où $\pi \in \gA$ est le produit  
$\prod_{i\ne i' ,\, j\ne j'} (x_iy_j - x_{i'}y_{j'})$.
\\
Soient $n$, $m \ge 2$ avec $\pgcd(n,m) = 1$, $\zeta_n$, $\zeta_m$, $\zeta_{nm}$
des racines de l'unité d'ordres respectifs $n$, $m$, $nm$. Par le \tho chinois, on obtient $\Phi_{nm} = \Phi_n \otimes \Phi_m$. \\
Comme
$\Phi_n(0) = \Phi_m(0) = 1$ (car $n, m \ge 2$), on a l'\egt

\snic {
\Delta_{nm} = \Delta_n^{\varphi(m)}\Delta_m^{\varphi(n)} \,\pi\,,
}

%\sni
où $\pi \in \ZZ$ est le produit suivant.

\snic{\prod_{i\ne i' ,\, j\ne j'}
(\zeta_n^i \zeta_m^j - \zeta_n^{i'}\zeta_m^{j'})$, pour $i,i'\in(\ZZ\sur{n\ZZ})\eti$ et $j, j'\in(\ZZ\sur{m\ZZ})\eti.}

%\sni
Soit $C \subset
(\ZZ\sur{nm\ZZ})\eti \times (\ZZ\sur{nm\ZZ})\eti$ l'ensemble des couples
$(a,b)$ avec $a, b$ \ivs modulo $nm$, $a \not\equiv b \bmod n$, $a
\not\equiv b \bmod m$. Le \tho chinois nous donne

\snic {
\pi = \prod_{(a, b) \in C}(\zeta_{nm}^a - \zeta_{nm}^b). 
}

%\sni

Soit $z \mapsto \ov z$  la conjugaison complexe. Alors, $\pi$ est de la
forme $z\ov z$, donc  $\pi\in\NN\etl$.
\\
En effet, $(a,b) \in C \Rightarrow (-a,-b) \in C$ avec $(a,b) \ne
(-a,-b)$.\\ 
Par ailleurs, pour $c \in \ZZ$
non multiple de $n$, ni de $m$, considérons l'\elt $\zeta_{nm}^c$ qui est
d'ordre $nm/\pgcd(c,nm) = n'm'$ avec $n' = n/\pgcd(c,n) > 1$, $m' > 1$ et $\pgcd(n',m') = 1$. Donc $n'm'$ n'est pas la puissance d'un
nombre premier, et, d'après la question précédente, $1 - \zeta_{nm}^c$ est
\iv dans~$\ZZ[\zeta_{nm}^c]$, a fortiori dans~$\ZZ[\zeta_{nm}]$. On en déduit
que $\pi$ est inversible dans~$\ZZ[\zeta_{nm}]$, donc dans $\ZZ$. 
\\
Bilan: $\pi =
1$, et $\Delta_{nm} = \Delta_n^{\varphi(m)}\Delta_m^{\varphi(n)}$.
\\
Enfin, si la formule qui donne le \discri cyclotomique est vérifiée
pour deux entiers $n$, $m$ étrangers entre eux, elle est vérifiée pour
le produit $nm$ (utiliser le premier point). Or elle est vraie pour des entiers
puissances d'un premier d'après la question \emph {5}, donc elle est vraie
pour tout entier $\ge 3$.

%%%%%%%%%%%%%%%%%%%%%%%%%%%%%%%%%%%%%%%%%%%%%%%%%%%%%%%%%%%%%%%%%%%%%%%%%%%%%%

\prob{exoAnneauEuclidien}
\emph{\ref{i5exoAnneauEuclidien}.} Considérons $p\equiv1\mod 4$. Le \polz~$Y^{p-1 \over 2} - 1 \in \FF_p[Y]$ est de degré $< \#\FF_p\eti$.
Il existe donc $y \in \FF_p\eti$ non racine de ce \pol; on pose $x = y^{p-1
\over 4}$ de sorte que $x^2 = y^{p-1 \over 2} \ne 1$~; mais $x^4 = 1$ donc $x^2 =
-1$. En fait,  pour la moitié des $y
\in \FF_p\eti$, on a $y^{p-1 \over 2} = 1$ (les carrés), et pour l'autre
moitié (les non-carrés), on~a~$y^{p-1 \over 2} = -1$.
\\
Voyons la question de l'\algo rapide. On entend par là que le temps d'exécution a pour ordre de grandeur une petite puissance du nombre de chiffres de~$p$.
\\
On détermine d'abord un $x \in \FF_p$ tel que $x^2 = -1$. Pour cela on  tire au hasard des entiers $y$ sur $\lrb{2..(p-1)/2}$ et l'on calcule $y^{p-1
\over 4}$ dans $\FFp$ (on utilise pour cela un \algo rapide d'exponentiation modulo $p$). La probabilité d'échec (lorsque le résultat est $\pm1$) est de
$1/2$ à chaque tirage.  
\\
Une fois trouvé un tel $x$, il reste à calculer
$\pgcd(x+i,p)$ avec l'\algo d'Euclide. Comme la norme est divisée par au moins $2$
à chaque étape, l'\algo est rapide.
\\
NB: la méthode brutale qui consisterait à dire, \gui{puisque $p\equiv 1\mod 4$,
il possède un facteur de la forme $m+in$, et il ne reste qu'à essayer tous les
$m<p$}, s'avère rapidement impraticable dès que $p$ devient grand. 

\emph{\ref{i6exoAnneauEuclidien}.} La \dcn des diviseurs premiers de $m$ est traitée dans le point précédent. Il reste à décomposer $n+qi$. 
\\
 Pour ce qui concerne la \dcn de $n^2+q^2$, on sait déjà que les seuls nombres premiers  y figurant sont  $2$ (avec l'exposant $1$) ou des
 $p\equiv 1\mod 4$. \\
 Si $u+vi$ est facteur d'un $p$
qui divise $n^2+q^2$, alors
 $u+vi$ ou $u-vi$ \hbox{divise $n+qi$}. Si $p$ figure avec l'exposant $k$ 
 dans  $n^2+q^2$, et si $u+vi$ divise $n+qi$, alors $u+vi$ 
 figure avec l'exposant $k$  dans~$n+qi$.
\\
 Si $s=2^k\prod_ip_i^{m_i}\prod_jq_j^{n_j}$ avec les $p_i\equiv 3\mod 4$
 et les $q_j\equiv 1 \mod 4$, alors la condition pour que $s$ soit somme de deux carrés est que les $m_i$ soient tous pairs.
\\
 On note qu'à une écriture $s=a^2+b^2$ avec $0<a\leq b$ correspondent 
  deux \elts conjugués $a\pm i b$ définis à association près
 (par exemple multiplier par $i$ revient à permuter $a$ et $b$).
 Il s'ensuit que dans le cas où $s$ est somme de deux carrés,  
  le nombre d'écritures de $s$ comme somme de deux carrés est égal 
  à $(1/2)\prod_j(1+n_j)$ sauf si les $n_j$ sont tous pairs, auquel cas on rajoute  ou retranche $1/2$  selon que l'on considère qu'une écriture
 $a^2+0^2$ est ou n'est pas légitime comme somme de deux carrés. 
\\
 Par exemple avec $5=\rN(a)$, $a=2+i$ et $13=\rN(b)$, $b=3+2i$ on obtient:
{\footnotesize 
$$\!\!{\begin{array}{rcl} 
5=\rN(a)  &\hbox{donne}& 5=2^2+1^2,  \\[1mm] 
10=\rN\big(a(1+i)\big)=\rN(1+3i)  &\hbox{donne}&  10=1^2+3^2, \\[1mm] 
5^3=\rN(a^3)=\rN(5a)  &\hbox{donne}& 125=2^2+11^2=10^2+5^2,  \\[1mm] 
5^4=\rN(a^4)=\rN(5a^2)=\rN(25)  &\hbox{donne}& 625=7^2+24^2=15^2+20^2=25^2+0,  \\[1mm] 
5^2\times 13=\rN(a^2 b)=\rN(a^2 \ov b)=\rN(5 b)  &\hbox{donne}& 325=18^2+1=17^2+6^2=15^2+10^2.   
 \end{array}}$$}

\vspace{-8pt}
De même 
$5^3\times 13=\rN(a^3b)=\rN(a^3\ov b)=\rN(5ab) = \rN(5a\ov b)$
donne

\snic{1625=16^2+37^2=28^2+29^2= 20^2+35^2= 40^2+5^2.}

%\sni
Et un calcul analogue donne 

\snic{1105=5\times 13\times 17=9^2+32^2=33^2+4^2=23^2+24^2= 31^2+12^2.}

%%%%%%%%%%%%%%%%%%%%%%%%%%%%%%%%%%%%%%%%%%%%%%%%%%%%%%%%%%%%%%%%%%%%%%%%%%%

\prob{exoPetitKummer} \emph{1.} Le discriminant se spécialise et $\Delta$ 
est \iv modulo $p$.
\\
 Ensuite on note que
$\aqo{\ZZ[\alpha]}p\simeq\FFp[t]:=\aqo{\FFp[T]}{f(T)}$. 
Ceci implique déjà que les \ids $\gen{q_k,p}$ sont maximaux dans $\ZZ[\alpha]$.
Pour $j\neq k$, $\gen{Q_j(t)}+\gen{Q_k(t)}=\gen{1}$
\hbox{dans $\FFp[t]$}, \hbox{donc $\gen{q_j}+\gen{q_k}+\gen{p}=\gen{1}$} 
\hbox{dans $\ZZ[\alpha]$}. D'où $\gen{q_j,p}+\gen{q_k,p}=\gen{1}$. 
\\
Par le \tho chinois, le produit
des $\gen{q_k,p}$ est donc égal à leur intersection, qui est égale à
$\gen{p}$ parce que l'intersection des $\gen{Q_j(t)}$ dans $\FFp[t]$
est égale à leur produit, qui est nul. 
\\
 Notons que l'\egt $\gen{p}=\prod_{k=1}^{\ell}\gen{p,Q_k(\alpha)}$
se maintient dans tout anneau contenant $\ZZ[\alpha]$. Même chose
pour le caractère comaximal des \idsz.
\\
 Si l'on passe de $\ZZ[\alpha]$ à $\gA$, la seule chose qui reste donc 
à vérifier est que \hbox{les  $\gen{p,q_k}$} restent bien des \idemas (stricts). C'est bien le cas et les corps quotients sont isomorphes. En effet, tout \elt de $\gA$ s'écrit $a/m$ avec $a\in\ZZ[\alpha]$ et~$m^2$
qui divise $\Delta$ (proposition \ref{propAECDN}), donc qui est étranger à $p$. Et 
l'\homo naturel $\aqo{\ZZ[\alpha]}{p,q_k}\to\aqo\gA{p,q_k}$ est un \isoz. 

\emph{2.} On applique l'exercice \ref{exoDecompIdeal}.

%%%%%%%%%%%%%%%%%%%%%%%%%%%%%%%%%%%%%%%%%%%%%%%%%%%%%%%%%%%%%%%%%%%%%%%%%%%

\prob{exoPolCyclotomique}
\emph{1a.}
On en déduit pour des premiers $p_1$, $p_2,$ \ldots\, ne divisant pas $n$,
que $f(\xi^{p_1p_2\ldots}) = 0$, i.e. $f(\xi^m) = 0$ pour tout $m$
tel que $\pgcd(n,m) = 1$, ou encore que~$f(\xi') = 0$ pour tout $\xi'$, racine
primitive $n$-ième de l'unité. Donc $f = \Phi_n$.
\\
 \emph{1b.}
Soit $h(X) = \pgcd_{\QQ[X]}\big(f(X),g(X^p)\big)$. 
Par le \tho de \KRO $h\in\ZZ[X]$.  On a $h(\xi) = 0$ donc $\deg
h \ge 1$. Raisonnons modulo $p$.  On a $g(X^p) = g(X)^p$, \hbox{donc $\ov h
\divi \ov f$} \hbox{et $\ov h \divi \ov g^p$}. Si $\pi$ est un facteur irréductible de $\ov h$,
 $\pi^2$ est un facteur carré de $X^n -\ov 1$, mais $X^n - \ov 1$ 
 est \spl dans~$\FF_p[X]$.
\\
Note: le \discri du \polz~$X^n + c$ est $(-1)^{\frac{n(n-1)}2}n^n c^{n-1}$,
en particulier celui de $X^n-1$ est $(-1)^{\frac{(n+2)(n+3)}2}n^n$.
 
 \emph{2.}
Si $G$ un groupe cyclique d'ordre~$n$, 
%Pour $m \in\ZZ\sur{n\ZZ}$, on définit $e_m : G \to G$ par $e_m(x) = x^m$. 
%On obtient
%ainsi tous les \endos du groupe $G$. En outre, $e_m \in \Aut(G)$ \ssi $m \in
%(\ZZ\sur{n\ZZ})^{\!\times}$. 
on a les \isos  classiques 

\snuc{\End(G)\simeq  \ZZ/n\ZZ\hbox{ (comme anneaux) et }\Aut(G)\simeq\big((\ZZ/n\ZZ)^{\!\times},\times \big) \hbox{ (comme groupes)}.}

%\sni
D'où des \isos canoniques 
 $\Aut(\UU_n)\simeq(\ZZ/n\ZZ)^{\!\times}\simeq\Gal(\gQ_n/\QQ)$.\\ 
Si $m \in (\ZZ/n\ZZ)^{\!\times}$, on obtient l'\auto $\sigma_m$ de $\gQ_n$  
 défini par $\sigma_m(\zeta)=\zeta^{m}$ pour $\zeta\in \UU_n$.
 
\emph{3.}
Supposons connaître un corps de racines $\gL$ en tant qu'extension \stfe de $\gK$. L'application $\sigma
\mapsto \sigma\frt{\UU_n}$ un morphisme injectif de $\Aut_\gK(\gL)$ vers $\Aut(\UU_n)$. En particulier, $\Aut_\gK(\gL)$ est  isomorphe à un sous-groupe de $(\ZZ/n\ZZ)^{\!\times}$.  Par
ailleurs, pour toute racine primitive $n$-ième de l'unité $\xi$ dans
$\gL$, on a $\gL = \gK(\xi)$. Donc, tous les facteurs \irds de
$\Phi_n(X)$ dans $\gK[X]$ ont même degré $\dex{\gL : \gK}$.  Mais il n'est pas
évident a priori de préciser quel type d'opération  sur
$\gK$  est \ncr pour factoriser $\Phi_n(X)$ dans $\gK[X]$.
Voici un exemple où  l'on peut déterminer
de manière certaine $\dex{\gL : \gK}$:  $p$ est
premier $\ge 3$, $p^* = (-1)^{p-1 \over 2}p$ et $\gK = \QQ(\sqrt {p^*})$. Alors
$\gK \subseteq \gQ_p$ (Gauss), la seule racine $p$-ième de l'unité contenue
dans $\gK$ est $1$ et $\Phi_p(X)$ se factorise dans $\gK[X]$ en produit de deux
\pols\irds de même degré~$p-1 \over 2$.

%%%%%%%%%%%%%%%%%%%%%%%%%%%%%%%%%%%%%%%%%%%%%%%%%%%%%%%%%%%%%%%%%%%%%%%%%%%

\prob{exoCyclotomicRing} 
\emph {1a.}
D'une part $\gA\sur\fp_i \simeq \aqo {\FF_p[X]}{\ov {f_i}}$ donc $\fp_i$
est maximal.  D'autre part, on  a un \iso $\ov\gA
= \gA\sur{p\gA} \Simarrow \aqo {\FF_p[X]} {\ov {\Phi_n}}$: 

\snic{h(\zeta_n) \bmod p \;\Simarrow\; \ov h \bmod \ov
{\Phi_n} \quad (\hbox{pour tout }h \in \ZZ[X]).}

En notant $\pi : \gA \twoheadrightarrow \ov\gA$ la surjection
canonique, on a $\sqrt {p\gA} = \pi^{-1} \big(\rD_{\ov\gA}(0)\big)$ et dans
l'isomorphisme ci-dessus:

\snic {
\rD_{\ov\gA}(0) \Simarrow \aqo {\gen {\ov g}} {\ov {\Phi_n}} 
= \bigcap_i \aqo {\gen {\ov {f_i}}}{\ov {\Phi_n}} 
= \prod_i \aqo {\gen {\ov {f_i}}}{\ov {\Phi_n}} 
,}

d'où le résultat.

\emph {1b.}
Résulte du fait que $\Phi_n$ est \spb modulo $p$.

\emph {1c.}
On vérifie facilement les \egts suivantes dans $\ZZ[X]$:
$$
\Phi_n(X) = \Phi_{mp}(X^{p^{k-1}}) =
{\Phi_{m}(X^{p^{k}}) \over \Phi_{m}(X^{p^{k-1}}) }
,$$
et donc dans $\FF_p[X]$, en notant $\varphi$ l'indicateur d'Euler:
$$
\Phi_n(X) = {\Phi_{m}(X)^{p^{k}} \over \Phi_{m}(X)^{p^{k-1}}} =
\Phi_{m}(X)^{\varphi(p^{k})} \qquad \bmod p
.$$
Le \polz~$\Phi_m$ est \spb modulo $p$, donc la partie sans facteur carré
de $\Phi_n$ modulo $p$ est $\ov g = \ov {\Phi_m}$; d'où
$\sqrt {p\gA} = \gen {p, \Phi_m(\zeta_n)}$.
\\
Montrons que $p \in \gen {\Phi_m(\zeta_n)}$. Si $\zeta_p \in \UU_n$ est une
racine primitive $p$-ième de l'unité, on a l'\egt

\snic {
\Phi_p(X) = \sum_{i = 0}^{p-1} X^i = \prod_{j=1}^{p-1} (X - \zeta_p^j),
}

%\sni
d'où, en faisant $X := 1$:

\snic {
p = \prod_{j=1}^{p-1} (1 - \zeta_p^j) \in \gen {1 - \zeta_p}.
}

%\sni

En appliquant cela à $\zeta_p = \zeta_n^{mp^{k-1}}$, on obtient $p \in 
\big\langle{1 - \zeta_n^{mp^{k-1}}}\big\rangle$. Mais $X^{mp^{k-1}} - 1$ est un multiple de $\Phi_m$
dans $\ZZ[X]$, donc $\zeta_n^{mp^{k-1}} - 1$ est un multiple de
$\Phi_m(\zeta_n)$ dans $\gA$, d'où $p \in \gen {\Phi_m(\zeta_n)}$.

\emph {1d.}
Comme $\sqrt {p\gA} = \fp_1 \cdots \fp_k = \gen {\Phi_m(\zeta_n)}$ est de type
fini, il y a un exposant $e$ tel que $(\fp_1 \cdots \fp_k)^e \subseteq p\gA$
et l'on applique l'exercice \ref{exoDecompIdeal}.  \\
Note: on peut prendre $e =
\varphi(p^k) = p^k - p^{k-1}$.

\emph {2.}
Le premier point est immédiat. Ensuite, si $\fa$ est un \itf non nul
de $\gA$, il contient un \elt $z$ non nul. Alors, $a = \rN_{\gQ_n\sur\QQ}
(z) = z\wi z$ est un entier non nul appartenant à $\fa$. On
écrit $a\gA \subseteq \fa$ comme produit d'\idemas\ivs et l'on 
applique de nouveau à l'\id $\fa$ l'exercice~\ref{exoDecompIdeal}.

%%%%%%%%%%%%%%%%%%%%%%%%%%%%%%%%%%%%%%%%%%%%%%%%%%%%%%%%%%%%%%%%%%%%%%%%%%%
\prob{exoSommeGauss}
\emph {1.}
Soit $x_0 \in G$ tel que $\varphi(x_0) \ne 1$.  \\
On écrit $\sum_{x \in G}
\varphi(x) = \sum_{x \in G} \varphi(xx_0)$, donc $S
\varphi(x_0) = S$ avec $S = \sum_{x \in G} \varphi(x)$, \cad $\big(1 -
\varphi(x_0)\big) S = 0$, d'où $S = 0$.
\\
\emph{2.}
Remarquons d'abord que $\chi^{-1}(-1) = \chi(-1)$ puisque
 $\chi(-1)^2 = \chi\big((-1)^2\big) = 1$. On écrit:

\snic {
\sum_{x+y = z} \chi(x) \chi^{-1}(y) = 
\sum_{x \ne 0,z} \chi\left( {x \over z-x} \right)
.}

%\sni
Si $z \ne 0$, l'application 
$x \mapsto {x \over z-x}$ est une bijection de $\gk \cup
\{\infty\}$ sur $\gk \cup \{\infty\}$ qui transforme $z$ en $\infty$, $\infty$
en $-1$, $0$ en $0$, ce qui donne une bijection de $\gk\eti \setminus \{z\}$
sur~$\gk\eti \setminus \{-1\}$. On peut donc écrire:

\snic {
\sum_{x+y = z} \chi(x) \chi^{-1}(y) = 
\sum_{v \in \gk\eti  \setminus \{-1\}} \chi(v) =
\sum_{v \in \gk\eti} \chi(v) - \chi(-1) = 0 -\chi(-1).
}

%\sni
Si $z = 0$ on a l'\egt

\snic {
\sum_{x+y=z} \chi(x) \chi^{-1}(y) = \sum_{x \ne 0} \chi(-1) = (q-1)\chi(-1).
}

%\sni
\emph{3.}
On écrit:

\snic {
G_\psi(\chi)G_\psi(\chi^{-1}) = \sum_{x,y} \chi(x) \chi^{-1}(y) \psi(x+y) =
\sum_{z \in \gk} S(z) \psi(z), 
}

%\sni

avec $S(z) = \sum_{x+y = z} \chi(x) \chi^{-1}(y)$.  D'où:

\snic {\arraycolsep2pt
\begin {array}  {rcl}
G_\psi(\chi)G_\psi(\chi^{-1}) &=& 
(q-1) \chi(-1) - \chi(-1) \sum_{z \ne 0} \psi(z) 
\\[1mm]
&=&
q\chi(-1) - \chi(-1) \sum_{z \in \gk} \psi(z) = q\chi(-1).
\\
\end {array}
}

%\sni
\emph{4.}
Le premier point est immédiat.  
%En utilisant $\tau_0 + \tau_1 = -1$ et
%$4\tau_0\tau_1 = (\tau_0 + \tau_1)^2 - (\tau_0 - \tau_1)^2$, %
On a facilement
$\tau_0 \tau_1 = {1 - p^* \over 4}$.   Le reste suit.

%%%%%%%%%%%%%%%%%%%%%%%%%%%%%%%%%%%%%%%%%%%%%%%%%%%%%%%%%%%%%%%%%%%%%%%%%%%
\prob{exoDedekindPolynomial}
\emph{1.}  %%%%%%%%%%%%%%%%%%%%%%%%%%%%%%%%%%%%%%%%%
Si $g(x) = 0$, avec $x \in \ZZ$ et $g(X) \in \ZZ[X]$ \monz, alors $x \divi
g(0)$.  Ici $\pm1, \pm2, \pm4, \pm8$ ne sont pas racines de $f(X)$, donc ce
\pol est \irdz. Le \discri du \polz~$X^3 + aX^2 + bX + c$ est:

\snic {
18abc - 4a^3c + a^2b^2 - 4b^3 - 27c^2 ,\quad
\hbox {d'où le résultat pour $a = 1$, $b = -2$, $c = 8$.}
}

\emph{2.} %%%%%%%%%%%%%%%%%%%%%%%%%%%%%
L'\elt $\beta = 4\alpha^{-1} \in \QQ(\alpha)$ est entier sur~$\ZZ$
puisque:

\snuc {
\alpha^3 + \alpha^2 - 2\alpha + 8 = 0 \buildrel /\alpha^3 \over \Longrightarrow
    1 + \alpha^{-1} - 2\alpha^{-2} + 8\alpha^{-3} = 0
\buildrel \times 8 \over \Longrightarrow
8 + 2\beta - \beta^2 + \beta^3 = 0.
}

%\sni
Pour vérifier que $\gA = \ZZ \oplus \ZZ\alpha \oplus \ZZ\beta$ est un
anneau, il suffit de voir que $\alpha^2, \alpha\beta, \beta^2 \in
\gA$.  C'est clair pour $\alpha\beta = 4$. On a $\alpha^2 + \alpha - 2 +
2\beta = 0$, donc $\alpha^2 = 2 - \alpha - 2\beta$.
Et~$\beta^3 - \beta^2 +
2\beta + 8 = 0$ donc $\beta^2 = \beta - 2 - 8\beta^{-1} = \beta - 2 -
2\alpha$.  L'expression de~$(1, \alpha, \alpha^2)$ sur la base $(1, \alpha, \beta)$ est
fournie par:

\snic {
\matrix {
 \phantom{\J\limits_{i}}                 &
\matrix {1 & \alpha & \alpha^2}
                                   \cr
~~\matrix{1\cr\alpha\cr\beta}         &
\crmatrix{   1 & 0      &  2 \cr
            0 & 1      & -1 \cr
            0 & 0      & -2 \cr }
}.}

L'anneau $\ZZ[\alpha]$ est donc d'indice~2 dans~$\gA$; or

\snic {
\Disc_{\ZZ[\alpha]/\ZZ} =  \idg{\gA : \ZZ[\alpha]} ^2 \cdot \Disc_{\gA/\ZZ}
\quad \hbox {donc} \quad \Disc_{\gA/\ZZ} = -503.
}

%\sni
Le \discri de~$\gA$ étant sans facteur carré, $\gA$ est l'anneau des
entiers de $\QQ(\alpha)$.

\emph{3.} %%%%%%%%%%%%%%%%%%%%%%%%%%%%%
Montrons que $\alpha$, $\beta$ et $\gamma := 1 + \alpha + \beta$ forment, modulo~2,
un \sfio:
$$
%\snic {
\alpha + \alpha^2 = 2 - 2\beta,\qquad \beta^2 - \beta = 2 - 2\alpha,\qquad
\alpha\beta = 4,
%}
$$
%\penalty-2500
d'où modulo~2:

\snic {
\alpha \equiv \alpha^2,\quad\ \beta \equiv \beta^2,\quad
\gamma^2 \equiv \gamma,\quad \alpha + \beta + \gamma \equiv 1,\quad
\alpha \beta \equiv 0, \quad \alpha \gamma \equiv 0,\quad
\beta \gamma \equiv 0.
}

%\sni
On a donc $\gA/2\gA = \FF_2 \ov \alpha \oplus \FF_2 \ov \beta
\oplus \FF_2 \ov \gamma$. Si l'on veut calculer la \fcn de~2
dans~$\gA$, on remarque que $(\alpha, \beta, \gamma)$ est une $\ZZ$-base
de~$\gA$ et qu'en désignant par~$\pi$ le morphisme de réduction modulo~$2$, $\pi :\gA \to \gA/2\gA$,  les \ideps de~$\gA$ au dessus de $2$, sont les images réciproques des \ideps de  $\gA/2\gA$. Par exemple:~$\fa = \pi^{-1}(\{0\} \oplus \FF_2 \ov \beta
\oplus \FF_2 \ov \gamma) = \gen {2\alpha, \beta, \gamma}$.  
Ainsi en posant~$\fb = \gen {\alpha, 2\beta, \gamma}$ et~$\fc = \gen {\alpha, \beta, 2\gamma}$,
on~a~$\gA/\fa \simeq \gA/\fb \simeq \gA/\fc \simeq \FF_2$ et~$2\gA = \fa\fb\fc =
\fa \cap \fb \cap \fc$.
\\
De manière \gnlez, soit $\gK$ un corps de nombres vérifiant: $
\dex{\gK : \QQ} \ge 3$ et~$2$ est totalement décomposé
dans l'anneau d'entiers $\gZ_\gK$. Alors,~$\gZ_\gK$  n'est pas
monogène, i.e.  il n'existe aucun $x \in \gZ_\gK$ tel que $\gZ_\gK =
\ZZ[x]$.  En effet, $\gZ_\gK/2\gZ_\gK \simeq \FF_2^n$ et $\FF_2^n$ n'admet pas
d'\elt primitif sur $\FF_2$  si $n \ge 3$.

\emph{4.} %%%%%%%%%
En multipliant $1 \in \ff + \fb$ par $\gB'$, on obtient $\gB' \subseteq
\ff\gB' + \fb' \subseteq \gB + \fb'$, ce qui montre que $\gB \to \gB'/\fb'$ est
surjective. Montrons que $\gB \to \gB'/\fb'$ est injective, i.e.~$\fb' \cap
\gB = \fb$. En multipliant $1 \in \ff+\fb$ par $\fb'\cap\gB$ on obtient les inclusions

\snic{\fb' \cap \gB \subseteq (\fb' \cap \gB)\ff + (\fb' \cap \gB)\fb
\subseteq \fb\gB'\ff + \fb \subseteq \fb\gB + \fb \subseteq \fb.}

%\sni
\emph{5.} %%%%%%%%%
Dans le contexte précédent, soit $x \in \gZ_\gK$ de degré $n = \dex{\gK :\QQ}$. \\
Notons 
$d = \idg{\gZ_\gK : \ZZ[x]}$. On a $d\gZ_\gK \subseteq \ZZ[x]$ et
$d$ peut servir de conducteur de $\ZZ[x]$ dans $\gZ_\gK$.  Si $2 \nedivi d$,
par l'évitement de Dedekind, $\gZ_\gK/2\gZ_\gK \simeq \ZZ[x]/2\ZZ[x] =
\FF_2[\ov x]$. Or $\gZ_\gK/2\gZ_\gK \simeq \FF_2^n$ n'admet pas d'\elt
primitif sur $\FF_2$ pour $n \ge 3$.

%%%%%%%%%%%%%%%%%%%%%%%%%%%%%%%%%%%%%%%%%%%%%%%%%%%%%%%%%%%%%%%%%%%%%%%%%%%
\prob{exoGaloisNormIdeal}
\emph{1.} $z \in \gB$ est racine de $\prod_{\sigma \in G}(T - z)$, \polu
à \coes dans $\gA$.

\emph{2.} 
$\ov\fm = \fm$ est clair. Calculons $\fm^2$ en écrivant $d = 4q+1$,
donc $1+d = 2(2q+1)$:

\snic {\arraycolsep2pt
\begin {array} {rcl}
\fm^2 &=& \gen {1 + 2\sqrt d + d, 1-d, 1 - 2\sqrt d + d} \\[1mm]
&=& 2 \gen {2q+1 + \sqrt d, 2q, 2q+1 - \sqrt d} =
2 \gen {1 + \sqrt d, 1 - \sqrt d} = 2\fm. 
\end {array}
}

%\sni
Par ailleurs, comme $\ZZ$-module, $\fm = \ZZ(1+\sqrt d) \oplus \ZZ(1-\sqrt d)
= 2\ZZ \oplus \ZZ(1\pm\sqrt d)$. On ne peut pas simplifier $\fm^2 = 2\fm$ par
$\fm$ (car $\fm \neq 2\gB$ vu que $1 \pm \sqrt d \notin 2\gB$), donc $\fm$
n'est pas \ivz.  On a $\rN_G(\fm) = 2\ZZ$ donc $\rN_G(\fm)\gB = 2\gB \ne
\rN'_G(\fm)$.
\\
L'application canonique $\ZZ \to \gB/\fm$ est surjective (puisque $x+y\sqrt d
\equiv x+y \bmod\fm$), de noyau $2\ZZ$, donc $\FF_2 \simeq \gB/\fm$. Ou encore:
$x + y\sqrt d \to (x+y) \bmod 2$ définit un morphisme surjectif \emph
{d'anneaux} $\gB \twoheadrightarrow \FF_2$, de noyau $\fm$.
\\
Notons $\rN(\fb) = \#(\gB/\fb)$ pour $\fb$ non nul.  Si $z = x(1+\sqrt d) +
y(1-\sqrt d) \in\fm$ avec~$x$,~$y \in \ZZ$, alors $\rN_G(z) = (x+y)^2 - d(x-y)^2
\equiv 4xy \bmod 4$. \\
Donc $\rN_G(z) \in 4\ZZ$ pour $z \in \fm$, mais
$\rN(\fm) = 2$. On a $\rN(\fm^2) = \rN(2\fm) = 4\rN(\fm) = 8$, mais~$\rN(\fm)^2
= 4$.

\emph{3.}
Soit $\fb = \gen {b_1, \ldots, b_n}$ et $n$ \idtrs $\uX = (X_1,
\ldots, X_n)$.  Introduisons le \pol normique $h(\uX)$:

\snic {
h(\uX) = \prod\nolimits_{\sigma \in G} h_\sigma(\uX)
\quad \hbox {avec} \quad
h_\sigma(\uX) = \sigma(b_1) X_1 + \cdots +  \sigma(b_n) X_n.
}

%\sni
On a $h(\uX) \in \gA[\uX]$. Notons $d$ un \gtr de $\rc(h)_\gA$.
Comme $\gB$ est \icl et $\rc(h)_\gB = d\gB$ principal, on peut appliquer la
proposition \ref{propArm}: on a alors~$\prod_\sigma \rc(h_\sigma)_\gB =
\rc(h)_\gB = d\gB$, \cad  $\rN'_G(\fb) = d\gB$.
\\
On va utiliser $\gA$ \icl (car $\gA$ est de Bézout).  Soit $a \in \gA \cap d\gB$.
Alors l'\elt $a/d \in \Frac(\gA)$ est entier sur $\gA$ (car $a/d \in \gB$) donc $a/d
\in \gA$, i.e. $a \in d\gA$. \\
Bilan: $\gA \cap d\gB = d\gA$ i.e. $\rN_G(\fb) =
d\gA$.
\\
Par \dfnz, les \evns du \pol normique $h$ sur $\gB^n$ sont les
normes d'\elts de l'idéal $\fb$; elles appartiennent à l'idéal de $\gA$
engendré par les \coes du \pol normique, cet idéal
de $\gA$ étant $\rN_G(\fb)$.
\\
Si $\#G = 2$, le \coe de $X_1X_2$ dans $h$ est:

\snic {
h(1,1, \ldots, 0) - h(1, 0, \ldots,0) - h(0, 1, \ldots,0) =
\rN_{G}(b_1+b_2) - \rN_{G}(b_1) - \rN_{G}(b_2).
}

%\sni
Ceci revient d'ailleurs à écrire $b_1\ov{b_2} + b_2\ov{b_1} =
\rN_{G}(b_1+b_2) - \rN_{G}(b_1) - \rN_{G}(b_2)$.  De même, le \coe de
$X_iX_j$ dans $h$ est, pour $i \ne j$, $\rN_{G}(b_i+b_j) -
\rN_{G}(b_i) - \rN_{G}(b_j)$. En conséquence, l'idéal de $\gA$ engendré
par les normes $\rN_G(b_i)$ et $\rN_G(b_i + b_j)$ contient tous les
\coes de $h(\uX)$; c'est donc l'idéal $\rN_G(\fb)$.

%%%%%%%%%%%%%%%%%%%%%%%%%%%%%%%%%%%%%%%%%%%%%%%%%%%%%%%%%%%%%%%%%%%%%%%%%%%

\prob{exoLemmeFourchette} 
 \emph{(Lemme de la fourchette)}\index{Lemme de la fourchette}
\\
\emph{1.}
Pour $x \in \gL$, on a $x = \sum_j \Tr_{\gL/\gK}(xe_j)e'_j$.
\\
Si $x \in \gB$, alors
$\Tr_{\gL/\gK}(xe_j)$ est un \elt de $\gK$ entier sur $\gA$ donc
dans $\gA$. Ceci démontre l'inclusion du milieu.
\\
En écrivant $e_i = \sum_j \Tr_{\gL/\gK}(e_ie_j)e'_j$, on
obtient 

\snic{\tra {\ue} = A \tra {\ue'}$ où $A = \big(\Tr_{\gL/\gK}(e_ie_j)\big)
\in \Mn(\gA)$, avec $\det(A)=\Delta,}

d'où l'inclusion de droite.

\emph{2.} Le \ZZmo $F_k$ est l'intersection de $\gB$ et $Z_k$, qui sont deux sous-\mtfs de $Z_{n-1}$, libre de rang $n$. C'est donc un \ZZmo libre de rang fini.
Et les deux inclusions $\delta Z_k\subseteq F_k\subseteq Z_k$
montrent que $F_k$ est de rang $k+1$.   
\\
 Le \ZZmo $\pi_k(F_k)$ est un sous-\ZZmo \tf de $\fraC1\delta\ZZ$.
Donc il est engendré par  $a_k/\delta$ (où $a_k$ est le pgcd des
numérateurs des \gtrsz). 
\\
Enfin, comme $1=\pi_k(x^{k})$, $a_k$ doit diviser $\delta$ et l'on écrit $\fraC {a_k}\delta=\fraC1{d_k}$.

\emph{3.} Soit $k\geq 1$ et $z\in F_k$, si $\pi_k(z)=a/d_k$ (avec $a\in\ZZ$) on a $\pi_k(z-ay_k)=0$.
\\ 
Donc
$z-ay_k\in F_{k-1}$. Ainsi $F_k=\ZZ y_k\oplus F_{k-1}$ et l'on conclut par \recu sur $k$
que $z\in\bigoplus_{i=0}^{k}\ZZ y_k$.

\emph{4.} On a $y_iy_j\in F_{i+j}$ donc $\fraC 1 {d_id_j}=\pi_{i+j}(y_iy_j)\in \fraC 1 {d_{i+j}}
\ZZ$, autrement dit $d_{i+j}$ est multiple de $d_id_j$.

\emph{5} et \emph{6.} Montrons tout d'abord  que $d_kF_k\subseteq \ZZ[x]$ par \recu sur $k$. 
L'initialisation~\hbox{$k=0$} est claire.
On utilise ensuite le fait que $xy_{k-1} \in F_k$ \hbox{et $\pi_k(xy_{k-1})=\fraC1{d_{k-1}}$}, donc  

\snic{xy_{k-1} = \fraC  {d_k} {d_{k-1}} y_k + w_{k-1} \quad \hbox {avec } w_{k-1} \in F_{k-1}.}

Il vient $d_ky_k = xd_{k-1}y_{k-1} - d_{k-1}w_{k-1}$ et le second membre
est  dans $\ZZ[x]$ par \hdrz. Donc $d_ky_k\in\ZZx$ et 

\snic{d_kF_k=d_k(\ZZ y_k \oplus F_{k-1})=\ZZ d_ky_k \oplus d_k F_{k-1}\subseteq \ZZx + d_{k-1} F_{k-1}\subseteq \ZZx.}
 
\smallskip 
On a défini $f_k(X)$ unitaire de degré $k$ dans $\QQ[X]$
par l'\egt $f_k(x)=d_ky_k$. 
\\
Comme  $(1,\dots,x^{n-1})$
est aussi bien une $\ZZ$-base de $\ZZx$, qu'une $\QQ$-base de $\QQ[x]$,
et comme $d_ky_k\in\ZZ[X]$ on obtient {$f_k\in\ZZ[X]$}. 

Tout le reste suit facilement.

%%%%%%%%%%%%%%%%%%%%%%%%%%%%%%%%%%%%%%%%%%%%%%%%%%%%%%%%%%%%%%%%%%%%%%%%%%%

\prob{exoPolynomialAutomorphism}
\emph{1.} Si $F(G) = X$, on a
$\JJ(F)(0) \circ \JJ(G)(0) = \I_{\Ae n}$. \\
Comme $\JJ(G)(0)$ est inversible,
on applique le résultat à $G$.
On a $H \in \gS^n$
\hbox{avec $G(H) = X$}. Alors $F=F\circ G\circ H=  H$. Donc $F$, $G$ sont inverses l'un de l'autre
(comme transformations de $\gS^n$).

\emph{2.}
Immédiat. Et l'on peut vérifier a posteriori  $\Phi(\gS^n)
\subseteq \gS^n$ ainsi que l'\eqvcz:\\
\centerline{ $\Phi(G) = G \iff F(G) = X$.}

\emph{3.}
On écrit $F(X) = J_0 \cdot X + F_2(X)$, où le vecteur $F_2(X)$ est de degré $\ge 2$ en~$X$. Alors, $J_0^{-1} \cdot \big(F(G) - F(H)\big) =
G - H + J_0^{-1} \cdot \big(F_2(G) - F_2(H)\big)$.
\\
Puis $\Phi(G) - \Phi(H) = -J_0^{-1}
\cdot \big(F_2(G) - F_2(H)\big)$. Supposons $G_i - H_i \in \fm^d$ ($d \ge
1$), et montrons que chaque composante de $\Phi(G) - \Phi(H)$ appartient à
$\fm^{d+1}$; il en résultera l'inégalité voulue. Une telle composante
est une combinaison $\gA$-\lin de $G^\alpha - H^\alpha$ avec $\alpha \in
\NN^n$ et $|\alpha| \ge 2$. Pour simplifier les notations, faisons $n = 3$
et écrivons:

\snuc {
G^{\alpha} - H^{\alpha} =
(G_1^{\alpha_1} - H_1^{\alpha_1}) G_2^{\alpha_2} G_3^{\alpha_3} +
(G_2^{\alpha_2} - H_2^{\alpha_2}) H_1^{\alpha_1} G_3^{\alpha_3} +
(G_3^{\alpha_3} - H_3^{\alpha_3}) H_1^{\alpha_1} H_2^{\alpha_2}.
}

%\sni
Comme les $H_i$, $G_i$ sont sans terme constant, on a  $G^\alpha - H^\alpha \in \fm^{d+1}$, sauf peut-être pour
$(\alpha_2, \alpha_3) = (0,0)$ ou $(\alpha_1, \alpha_3) = (0,0)$ ou
$(\alpha_1, \alpha_3) = (0,0)$.  Il reste
à voir les cas particuliers, par exemple $\alpha_2 = \alpha_3 = 0$. Dans ce
cas, puisque $\alpha_1 - 1 \ge 1$:

\snic {
G^\alpha - H^\alpha =
G_1^{\alpha_1} - H_1^{\alpha_1} =
(G_1 - H_1) \sum_{i+j = \alpha_1-1} G_1^i H_1^j \in \fm^{d+1}.
}

%\sni
On a donc établi $d\big(\Phi(G), \Phi(H)\big) \le d(G, H)/2$. Ceci assure en
particulier qu'il existe au plus un point fixe de $\Phi$.  Soient $G^{(0)}\in
\gS^n$, par exemple $G^{(0)} = 0$, et la suite $G^{(d)}$ définie par
\recu au moyen de $G^{(d+1)} = \Phi(G^{(d)})$. \\
Pour $d \ge 1$, chaque composante
de $G^{(d)} - G^{(d-1)}$ est dans $\fm^d$, ce qui permet de
définir $G \in \gS^n$ par
$G = \sum_{d \ge 1} \big(G^{(d)} - G^{(d-1)}\big).$
\\
Alors, $G$ est la limite des $G^{(d)}$ pour $d \mapsto \infty$,
c'est un point fixe de $\Phi$, i.e. $F(G) = X$.

\emph{4.}
Supposons $G(F) = X$, donc $G\big(F(0)\big) = 0$. \\
On pose $\wi F = F - F(0)$, $\wi G =
G\big(X + F(0)\big)$. Alors, $\wi F(0) = \wi G(0) = 0$ \hbox{et $\wi G(\wi F) = X$}. D'où $\wi F(\wi G) = X$, puis $F(G) = X$.

\emph{5.}
On vérifie dans les deux cas que $\J(F) = 1$. Pour le premier, on obtient $G$
(de même degré maximum que $F$) en itérant $\Phi$ quatre fois:

\snic {
G = (-X^2Z^3 - 2XY^2Z^2 + 2XYZ + X - Y^4Z + 2Y^3,\
    -XZ^2 - Y^2Z + Y,\  Z).
}

%\sni
Pour le second, on obtient $G = (G_1, \ldots, G_5)$ en itérant $\Phi$ quatre fois:

\snic {\arraycolsep2pt
\begin {array} {rclrcl}
G_1 &=& X_1 - 3X_2X_4^2 + 6X_2X_4X_5^3 - 3X_2X_5^6 + 2X_3X_4X_5 - 2X_3X_5^4 +
\\[1mm]
&&X_4^4X_5 - 4X_4^3X_5^4 + 6X_4^2X_5^7 - 4X_4X_5^{10} + X_5^{13},
\\[1mm]
G_2&=&X_2 - X_4^2X_5 + 2X_4X_5^4 - X_5^7,
\\[1mm]
G_3&=&X_3 - X_4^3 + 3X_4^2X_5^3 - 3X_4X_5^6 + X_5^9,
\\[1mm]
G_4&=&X_4 - X_5^3,
\quad \quad 
G_5\;=\;X_4.
\end {array}
}

%\sni
On notera que le degré maximum de $G$ est 13 alors que celui de $F$
est $3$.
%%%%%%%%%%%%%%%%%%%%%%%%%%%%%%%%%%%%%%%%%%%%%%%%%%%%%%%%%%%%%%%%%%%%%%%%%%%

%: sinotenglish
\sinotenglish{

\prob{exoFinitudeClassesIdeaux} \emph{(Finitude de l'ensemble des classes d'\ids d'un 
anneau de nombres)}
\\ \emph {1a.}
Partageons l'intervalle semi-ouvert $[0,1[$ en $N$ sous-intervalles $[i/N,
(i+1)/N[$ de longueur $1/N$, pour $0 \le i \le N-1$. L'hyper-cube $[0,1[^n$
est une réunion de~$N^n$ petits hyper-cubes.  Pour $x \in \QQ^n$, notons
$\lfloor x\rfloor$ le vecteur de $\ZZ^n$ dont la $i$-ième composante est la
partie entière $\lfloor x_i\rfloor$ de la composante $x_i$.  Pour $0 \le k
\le N^n$, on considère les $N^n + 1$ vecteurs $kx - \lfloor kx\rfloor \in
[0,1[^n$. \\
D'après le principe des tiroirs, il en existe deux qui sont dans
le même petit hyper-cube, i.e. il existe $h$, $k$ distincts, $0 \le h<k \le
N^n$ avec:

\snic {
\nsup {(kx - \lfloor kx\rfloor) - (hx - \lfloor hx\rfloor)} < 1/N
}

\snii 
On pose  $m = k-h\in \lrb{1..N^n}$, $y = \lfloor kx\rfloor - \lfloor hx\rfloor
\in \ZZ^n$. On a bien $\nsup{mx - y} < 1/N$.

\snii \emph {1b.}
Prendre $N \in \NN^*$ tel que $N \ge 1/K$ et poser $d = N^n$.

\snii \emph {2.}
Soient $c_{ij}^k \in \QQ$ les constantes de structure définies par $e_ie_j =
\sum_k c_{ij}^k e_k$. \\
Pour $x \in \gK$, $x = x_1 e_1 + \cdots + x_n e_n$, on
veut calculer le \coe $a_{ij}$ de la matrice de la multiplication par $x$ dans
la base $(e_1, \ldots, e_n)$. On a:
$$
xe_j = \som_{k,i} x_k c_{kj}^i e_i = 
\som_i \big( \som_{k} x_k c_{kj}^i \big) e_i 
\quad \hbox {donc} \quad
a_{ij} = \som_{k} x_k c_{kj}^i.
$$ 
Soit $M = \max_{i,j} \sum_k |c_{kj}^i|$. Alors $|a_{ij}| \le M \abs{x}$,
donc $|\rN(x)| \le M^n \abs{x}^n$, et l'on peut prendre $C = M^n$.

\snii \emph {3.}
Soit $C > 0$ la constante de la question précédente. On prend 
$K > 0$ vérifiant~\hbox{$CK^n < 1$} et on lui applique la question
\emph {1}, ce qui nous fournit un $d \in \NN^*$. \hbox{Pour
$x \in \gK$}, il y a $m \in \lrb {1..d}$ et $q \in \gA$
vérifiant $\nsup{mx-q} < K$ donc:

\snic {
|\rN(mx-q)| \le C\nsup{mx-q}^n \le CK^n < 1.
}

\snii 
Pour le deuxième point de la question, on considère $x = a/b$.

\snii \emph {4a.}
Pour $a \in \fb$, il y a $m \in \lrb{1..d}$ et $q \in \gA$ tels que
$|\rN(ma-bq)| < |\rN(b)|$. La minimalité de $|\rN(b)|$ fait que $ma-bq = 0$;
donc $Da \in \gen {b}$ puis $D\fb \subseteq \gen {b}$.

\snii \emph {4b.}
Si $\fb=\gen{b_1, \ldots, b_k}$, pour $b \in \fb \setminus \{0\}$, il
existe $m_i \in \lrb{1..d}$ et $q_i \in \gA$ tels que~\hbox{$|\rN(m_ib_i - bq_i)| <
|\rN(b)|$}. \\
Ou bien $m_ib_i - bq_i = 0$ pour~\hbox{$i\in\lrbk$}, auquel
cas $D\fb \subset \gen {b}$.  \\
Ou bien il y a un $b' = m_ib_i - bq_i
\ne 0$; ce $b'$ appartient à~\hbox{$\fb \setminus \{0\}$ et $|\rN(b')| <
|\rN(b)|$}.  On recommence alors avec $b'$ à la place de~$b$. Au
départ, on peut utiliser pour $b$ l'un des $b_i$ non nul; ce processus
s'arrête car les valeurs absolues $|\rN(b)|$ des normes sont des entiers $>
0$.

\snii \emph {4c.}
On prend $\fa = (D/b)\,\fb$. C'est un \id (entier) de $\gA$, qui est associé à
$\fb$ et qui contient~$D$ (car $\fb$ contient $b$).

\snii \emph {5.}
Tout \itf non nul de $\gA$ est donc associé à un
\id $\fa$ contenant $D$. Or il n'y a qu'un nombre
fini de tels \ids $\fa$ car leur ensemble s'identifie
à celui des \ids de l'anneau fini $\gA\sur{D\gA}$
(de cardinal $D^n$).

%%%%%%%%%%%%%%%%%%%%%%%%%%%%%%%%%%%%%%%%%%%%%%%%%%%%%%%%%%%%%%%%%%%%%%%%%%%

%: solution prob  pbMatrixAndIdealClasses 
\prob{pbMatrixAndIdealClasses} \emph{(Classes de similitude de matrices et classes d'\idsz)} \\
On note $\gB = \Ax$.
\\
\emph {1.}
Notons $C$ au lieu de $C_f$.  Soit $(e_1,\dots,e_n)$  la base
canonique de $\Ae n$. On vérifie que les $\tra{C}^j e_n$, $j
\in \lrb {0..n-1}$, forment une base de $\Ae n$ dans laquelle 
la matrice de~$\tra{C}$ est $C$.  
Si l'on note $Q' \in \Mn(\gA)$ la matrice ayant pour
colonnes les $\tra{C}^j e_n$, on a donc $Q'C = \tra{C}Q'$.  Vu que $Q'$est une
matrice de Hankel inférieure, cf. ci-dessous,  son inverse $Q$ est  de Hankel supérieure et vérifie~\hbox{$CQ
= Q\tra{C}$}. \\
Si l'on utilise pour les lignes une numérotation de $0$ à $n-1$
(au lieu de $1$ à $n$), la matrice $Q$ a pour \coe d'indice $(i,j)$ le
coefficient de $f$ en $X^{i+j}$.  \\
Par exemple, si $f = X^5 + a_4X^4 + a_3X^3 +
a_2X^2 + a_1X + a_0$, voici $Q'$:

\snuc {
Q' = \cmatrix {
0 &0    &0          &0                 &1 \cr
0 &0    &0          &1                 &-a_4 \cr
0 &0    &1          &-a_4              &-a_3+a_4^2 \cr
0 &1    &-a_4       &-a_3+a_4^2        &-a_2+2a_3a_4-a_4^3 \cr
1 &-a_4 &-a_3+a_4^2 &-a_2+2a_3a_4-a_4^3 &-a_1+2a_2a_4 + a_3^2-3a_3a_4^2 +a_4^4\cr}
}

\snii  et son inverse $Q$:

\snic {
Q = \cmatrix {
a_1 & a_2 & a_3 & a_4 & 1 \cr
a_2 & a_3 & a_4 &   1 & 0 \cr
a_3 & a_4 &   1 &   0 & 0 \cr
a_4 &   1 &   0 &   0 & 0 \cr
1   &   0 &   0 &   0 & 0 \cr}
}

\snii  \rem Introduisons, en prenant $a_n = 1$, pour $i \in
\lrb{0..n}$, les \pols de Horner $f_i(X) = \sum_{j=i}^n a_jX^{j-i}$ (en
particulier
$f_0(X) = f(X)$ et $f_n(X) = 1$). \\
Alors $\big(f_1(x), \ldots, f_n(x)\big)$ est une
$\gA$-base de $\gB$. De plus $xf_{i+1}(x) = f_i(x) - a_i$, donc la
matrice de la multiplication par $x$ dans cette base est $\tra{C_f}$. Et la
matrice $Q$ ci-dessus est celle qui exprime la base $\big(f_1(x), \ldots, f_n(x)\big)$
dans la base $(1, \ldots, x^{n-1})$.\eoe

\snii \emph {2.}
Dans $\gA[X]$ la matrice $D = QS - \Adj(X\In-C_f)$ est triangulaire supérieure, de
diagonale nulle et $D_{i,j} = X^{j-(i+1)}f$ pour $j > i$. \\
Donc dans $\gA[x]$,
$QS - \Adj(x\In-C_f) = 0$.

\snii \emph{3a.}
D'une part $x\fb \subseteq \fb$,
donc $\fb$ est bien un \id de $\gB$.  
\\
D'autre part, en multipliant à
droite par $\wi P$ l'\egt définissant les $\vep_i$, on obtient
\hbox{que $\fb\ni\det(P)$}, \elt\ndz de $\gA$.

\snii \emph{3b.}
Clair.

\snii \emph{3c.}
Voir l'exercice \ref{exoAXmodule}

\snii \emph{4.}
Si $P \in \Mn(\gA)$ exprime une $\gA$-base de $\fb$ dans la $\gA$-base $(1, x,
\ldots, x^{n-1})$ de~$\gB$, alors $\det(P)\in\Reg(\gA)\cap\fb$ \hbox{(cf.
\emph {3a})}. 
\\
Remarque: si $b \in \gB$ est un \elt\ndz de $\gB$, alors
$\rN\iBA(b) = b\wi b$ est un \elt\ndz de $\gA$ (\thref{prop inj surj det}),
multiple de $b$ dans $\gB$.  \\
Dans $\gK[x]$, $\gK\vep_1 \oplus \cdots \oplus
\gK\vep_n$ est un \id de $\gK[x]$ contenant un \elt\ndz de $\gA$, donc cet \id est $\gK[x]$.  On dispose ainsi de deux $\gK$-bases
de~$\gK[x]$:

\snic {
\gK[x] = \gK.\vep_1 \oplus \cdots \oplus \gK.\vep_n = 
\gK.1 \oplus \gK.x \oplus \cdots \oplus \gK.x^{n-1}.
}

\snii 
Donc la matrice $M \in \Mn(\gA)$ est conjuguée sur $\gK$ de la matrice
compagne de $f$.

\snii \emph{5a.}
On a sur $\gK$:

\snic {
\tra {M} = \tra{P}\, \tra{C_f}\, \tra{P}^{-1} = 
\tra{P} Q^{-1} \ C_f\  Q\, \tra{P}^{-1},
}

\snii 
donc $\tra {M} = P'^{-1} C_f P'$ avec $P' = Q\tra {\wi {P}}$.

\snii \emph{5b.}
Par \dfn:

\snic {
N \eqdefi
\cmatrix {\vep'_1\cr \vdots\cr \vep'_n\cr} [\vep_1, \ldots, \vep_n] = 
\tra{P'} \cmatrix{1\cr \vdots\cr x^{n-1}\cr} [1, \ldots, x^{n-1}] P = 
\wi P \tra{Q} S P,
}

\snii  puis 

\snic { 
\begin{array}{ccccc} 
N  & =  & \det(P)P^{-1}\Adj(x\In-C_f)P  \\[1mm] 
  & =  &  \det(P)\Adj(x\In-P^{-1}C_fP)&=& \det(P)\Adj(x\In-M). \end{array}
}

\snii 
Cette \egt matricielle implique l'\egt d'\ids 

\snic{\fb'\fb =\det(P)\cD_{n-1}(x\In-M)}

\snii  
car les \coes de $\Adj(x\In-M)$ sont les mineurs
d'ordre $n-1$ de $x\In-M$.

\snii \emph{5c.}
Il s'agit de montrer de manière \gnlez, que pour $A \in \Mn(\gA)$
de rang $n-1$, on a $\Ker\wi A = \Im A$; d'après le \plg
de base, on peut supposer que $A$ possède un mineur d'ordre
$n-1$ \iv et d'après le lemme du mineur inversible \ref{lem.min.inv}
que $A = \Diag(1, \ldots, 1,0)$; dans ce cas,
$\wi A = \Diag(0, \ldots, 0, 1)$ et le résultat est clair.
\\
Une fois obtenue l'\egt $\Ker\wi A = \Im A$, on a
$\Im\wi A \simeq \Ae n/\Ker\wi A =
\Ae n/\Im A$.

%%%%%%
\snii \emph{6.}
L'\id $\fb$ a pour $\gA$-base $(a, x, x^2, \ldots, x^{n-1})$. La matrice $M
\in \Mn(\gA)$ de la multiplication par $x$ dans cette base est
à quelque chose près la matrice $C_f$ sauf \hbox{que $M_{1n} = -b_0$}
(au lieu de $-a_0$) et $M_{2,1} = a$ (au lieu de $1$). On
peut prendre pour~$P$ la matrice diagonale $\Diag(a, 1, \ldots, 1)$,
donc $\wi P = \Diag(1, a, \ldots, a)$. \\
Déterminons
l'\id $\fb'$ qui vérifie $\fb\fb' = a\cD_{n-1}(x\In-M)$;
on introduit, en convenant de $a_n = 1$, $f_i(X) = \sum_{j=i}^n
a_jX^{j-i}$ (donc $f_n(X) = 1$).\\
Alors $\big(f_1(x), af_2(x), \ldots,af_n(x)\big)$ est une $\gA$-base de $\fb'$, donc $\fb' = \gen {a,f_1(x)}$.\\
Pour $n=4$ par exemple, voici la matrice $x\In-M$:

\snic {
x\In-M = \cmatrix {
x  & 0  &  0 & b_0 \cr
-a &  x &  0 & a_1 \cr
0  & -1 &  x & a_2 \cr
0  &  0 & -1 & x+a_3 \cr
}.}

\snii 
Ainsi $\cD_{n-1}(x\In-M)$ contient $a$ et $b_0$; donc,
si $1 \in \gen {a,b_0}$, \hbox{alors $1 \in \cD_{n-1}(x\In-M)$}
\hbox{et $\fb\fb' = a\gB$}.

\snii \emph{7.}
Montrons d'abord que $a \in \fb^k$ pour $k \in \lrb {1.. n}$.  On va
utiliser à plusieurs reprises que $a_i \equiv 0 \bmod a$. On~a

\snic{
b_0a = -(a_1x + \cdots + a_{n-1}x^{n-1} + x^n) \in x\fb \subseteq \fb^2,
}

\snii 
i.e. $b_0a \equiv 0 \bmod\fb^2$; mais $b_0$ est \iv modulo $a$ donc modulo
$\fb$ (car $\fb$ contient~$a$) donc modulo toutes les puissances de $\fb$
d'o\'u $a \equiv 0 \bmod\fb^2$.  Ensuite, en utilisant que $a \in \fb^2$, on
voit que 

\snic{b_0a = -(a_1x + \cdots + a_{n-1}x^{n-1} + x^n) \in x\fb^2 \subseteq
\fb^3,}

\snii 
et par un raisonnement analogue, $a \in \fb^3$. De proche en proche,
on voit ainsi que $a \in \fb^k$ pour tout $k \le n$.

\snii  
Puisque $x^n \in a\gB$, le \Amo $\gA a + \gA x + \cdots + \gA x^{n-1}$ est
stable par $x$ donc c'est un \id de $\gB$, égal à $\fb$ et $(a, x,
\cdots, x^{n-1})$ est une $\gA$-base de $\fb$.  \\
De même, le \Amo $\gA a +\gA ax + \gA x^2+ \cdots + \gA x^{n-1}$ est stable par $x$ donc c'est un \idz.
Il contient $a^2$, $ax$, $x^2$ et d'autre part, il est contenu dans
$\fb^2$, donc c'est $\fb^2$.  \\
Même chose pour les autres puissances de $\fb$.

}
%: fin sinotenglish

}% fin des solutions d'exos

%:   ---- Section*{references}-----------
%\newpage	
\Biblio

\vspace{8pt}
La preuve du lemme de \DKM  \rref{lemdArtin} est prise dans Northcott \cite{Nor2} (il l'attribue à Artin).\perso{Chercher
plus de références pour \DKMz, Thierry propose de consulter
les oeuvres complètes de Dedekind.}

Le \tho de \KRO \rref{thKro} se trouve dans \cite[Kronecker]{Kro1}. Il est \egmt démontré par
Dedekind \cite{Ddk1} et Mertens \cite{Mer}.

Concernant les résultants et sous-résultants en une variable,
un livre de référence est \cite{AJ}.
On regrettera cependant l'absence de bibliographie:
même si les résultats sont soit très anciens soit complètement nouveaux,
on ne voit pas l'utilité de cacher les sources exactes.
\\
Un autre livre important pour les questions \algqs sur le sujet est l'ouvrage
\cite{BPR}.

La construction d'un \cdr abstrait pour un \pol \spl donnée dans le \thrf{thResolUniv} est (à très peu près) celle décrite par Jules Drach
dans \cite{Drach}, qui semble être celui qui introduit l'\adu comme outil fondamental
pour étudier les extensions \agqs de corps.

La preuve télégraphique du \thrf{thIntClosStab} nous a été
suggérée  par Thierry Coquand.

L'approche de Kronecker concernant la théorie des \ids de corps de nombres
fait l'objet d'un survey historique dans \cite[Fontana\&Loper]{FL}.

La \dem du \nst donnée dans la section \ref{secChap3Nst}
est  inspirée de celle dans \cite{BPR}, elle même inspirée 
d'une \dem de van der Waerden.

\newpage \thispagestyle{CMcadreseul}
\incrementeexosetprob

%:        %%%%%%%%%%%%%%%%%%%%%%%%%%%%%%%%%%%%
%:        %%%%%%%%%%%%%%%%%%%%%%%%%%%%%%%%%%%%

%---- Chapitre  {Modules \pf}------------
\chapter{Modules \pfz}
\label{chap mpf}
%--------------------
\minitoc

\subsection*{Introduction}
\addcontentsline{toc}{section}{Introduction}
%-----------------------------------------

Sur un anneau les \mpfs jouent un peu le même rôle que les \evcs de dimension
finie sur un corps: la théorie des \mpfs est une manière un peu plus
abstraite,
et souvent profitable, d'aborder la question des \slisz.

\smallskip Dans les premières sections du chapitre, 
on donne les bases de la théorie des \mpfsz.

Dans la section \ref{secBézout}, on traite l'exemple des \mpfs sur les anneaux principaux, et dans la section \ref{secKrull0dim} celui des \mpfs sur les anneaux \zedsz.

Enfin la section \ref{sec Fitt} est consacrée aux invariants importants
que sont les \idfsz, et la section \ref{subsecIdealResultant} introduit
l'\id résultant, comme application directe des \idfsz.

%--- Sec{Chgt syst gen}---  \label{sec pf chg}---
\newpage\section{Définition, changement de \sgrz} \label{sec pf chg}
%--------------------

 Un module \ixd{de présentation finie}{module} est un \Amo $M$ donné
par un nombre fini de \gtrs et de relations. C'est donc un \mtf
avec un \sgr possédant un module des relations \tfz.
De manière \eqvez, c'est un module $M$ isomorphe au conoyau d'une \ali
%-----------------begin $$------------------
$$
\gamma:  \Ae m\longrightarrow  \gA^q.
$$
%-----------------end $$------------------
La matrice $G\in \gA^{q\times m}$ de $\gamma$ a pour colonnes un \sgr
du module des syzygies entre les \gtrs $g_i$ qui sont les images de la base
canonique de $\gA^q$ par la
surjection $\pi: \gA^q\rightarrow M$.
Une telle matrice s'appelle une
\index{matrice!de présentation}
\emph{\mpn du module $M$ pour le \sgr $(\gq)$}.
Cela se traduit par:
%-----------------begin item------------------
\begin{itemize}
\item  $[\,g_1\;\cdots\;g_q\,]\,G=0$, et
\item  toute syzygie entre les $g_i$ est une \coli des colonnes de~$G$,
i.e.: si $[\,g_1\;\cdots\;g_q\,]\,C=0$ avec $C\in \gA^{q\times 1}$, il existe
 $C'\in \Ae {m\times 1}$  tel que~$C=G\,C'$.
\end{itemize}
%-----------------end item------------------

\smallskip \exls \rdb 
1) Un module libre de rang $k$ est un \mpf présenté par
une matrice colonne formée de $k$ zéros{\footnote{Si l'on considère qu'une
matrice est donnée par deux entiers  $q,m\geq 0$  et une famille d'\elts de
l'anneau indexée par les couples $(i,j)$ avec
$ i\in\lrbq,\;  j \in \lrbm$,
on peut accepter une matrice vide de type
$k\times 0$, qui serait la matrice canonique pour présenter un module libre de
rang $k$.}}.
Plus \gnlt toute matrice simple est la \mpn d'un module libre de rang fini.

 2) \label{exl1pf}
 Rappelons qu'un \mptf est un module $P$ isomorphe à l'image
d'une \mprn $F\in\Mn(\gA)$ pour un certain entier~$n$. Puisque $\Ae n=\Im (F)\oplus \Im(\In-F)$, on obtient $P\simeq \Coker(\In-F)$. Ceci montre que tout \mptf
est \pfz.

 3)
Soit $\varphi :V\to V$ un \endo d'un \evc de dimension finie sur un corps
discret $\gK$.
Considérons $V$ comme un $\KX$-module avec la loi externe suivante:
$$
\formule{\gK[X]\times V&\to& V\\[.2em] (P,u)&\mapsto& P\cdot u:=P(\varphi)(u).}
$$
Soit $(\un)$ une base de $V$ comme \Kev et $A$ la matrice de $\varphi$
sur cette base. Alors, on peut montrer qu'une \mpn de $V$ comme
$\KX$-module pour le \sgr  $(\un)$ est
la matrice $X\,\In-A$ (voir l'exercice~\ref{exoAXmodule}).  \eoe

%\smallskip
\setcounter{subsection}{-1}
%:     lemma factchangesgrmpf
\begin{lemma} \label{factchangesgrmpf}
Lorsque l'on  change de \sgr fini pour un \mpfz, les \syzys entre
les nouveaux \gtrs forment de nouveau un module \tfz.
\end{lemma}
\begin{proof}
Supposons en effet, avec $M\simeq \Coker G$, qu'un autre \sgr
de $M$ soit $(h_1,\ldots,h_r)$.
On a donc des matrices~$H_1\in\gA^{q\times r}$ et~$H_2\in\Ae {r\times q}$ 
telles que

\snic{[\,g_1\;\cdots\;g_q\,]\,H_1=[\,h_1\;\cdots\;h_r\,]$ et
$[\,h_1\;\cdots\;h_r\,]\,H_2=[\,g_1\;\cdots\;g_q\,].}

%\sni
Alors, le module des syzygies entre les $h_j$ est engendré
par les colonnes de~$H_2G$  et celles de $\I_r-H_2H_1$.
En effet, d'une part on a clairement

\snic{[\,h_1\;\cdots\;h_r\,]\,H_2\,G=0\;$ et
 $\;[\,h_1\;\cdots\;h_r\,]\,(\I_r-H_2H_1)=0.}

%\sni
D'autre part, si l'on a une \syzy $[\,h_1\;\cdots\;h_r\,]\,C=0$, on
en déduit

\snic{[\,g_1\;\cdots\;g_q\,]\,H_1C=0,} 

donc $H_1C=GC'$ pour un certain
vecteur colonne $C'$ et
$$\preskip.3em \postskip.2em
C=\big((\I_r-H_2H_1)+H_2H_1\big)C=(\I_r-H_2H_1)C+H_2GC'=HC'',
$$
où $H=\lst{\I_r-H_2H_1\mid H_2G}$ et $C''=\Cmatrix{1pt}{C\cr C'}$.
\end{proof}

\rdb
 La possibilité de remplacer un \sgr par un autre tout
en gardant un nombre fini de relations est un phénomène
extrêmement \gnlz. Il s'applique à toutes formes de structures
\agqs qui peuvent être définies par \gtrs et relations.
Par exemple, pour les structures dont tous les axiomes sont des \egts universelles.
Voici comment cela fonctionne (il suffira de vérifier dans chaque cas que le
raisonnement s'applique bien).\label{nouveausgr}

Supposons  que l'on  a des \gtrs  $g_1, $\ldots$, g_n$ et des relations

\snic{R_1(g_1,\alb\ldots,g_n)$, \ldots,  $R_s (g_1,\ldots,g_n),}

%\sni
qui \gui{présentent} une structure  $M$.
\\
Si l'on a d'autres \gtrs
 $h_1$, $\ldots$, $h_m$,
on les exprime en fonction des $g_j$ sous forme
  $h_i = H_i(g_1,\ldots,g_n)$. Notons  $ S_i(h_i,g_1,\ldots,g_n)$  cette relation.
\\
On exprime pareillement les $g_j$   en fonction des  $ h_i$:
  $g_j = G_j(h_1,\ldots, h_m)$.
Notons~$ T_j(g_j, h_1,\ldots, h_m)$  cette relation.

La structure ne change pas si l'on remplace la \pn
$$\preskip.3em \postskip-.2em 
(g_1,\ldots,g_n\; ; \alb\; R_1, \ldots, R_s) 
$$
    par
$$\preskip.0em \postskip.4em 
(g_1,\ldots,g_n, h_1,\ldots, h_m\; ; \; R_1,\alb \ldots,\alb R_s,
S_1,\ldots,S_m). 
$$
Comme les relations $T_j$  sont satisfaites, elles sont conséquences
des relations $R_1$, $\ldots$, $R_s,$ $S_1$, $\ldots$, $S_m$,
donc la structure est toujours la même avec la \pn suivante:

\snic{(g_1,\ldots,g_n, h_1,\ldots, h_m\; ; \; R_1, \ldots, R_s, S_1,\ldots,S_m,
T_1,\ldots,T_n).}

%\sni
Maintenant, dans chacune des relations $R_k$  et $S_\ell,$  on peut remplacer
chaque~$g_j$  par son expression en fonction des $h_i$  
(qui est donnée dans~$T_j$)
et cela ne change toujours pas la structure présentée.
On obtient
$$\preskip.3em \postskip.4em 
(g_1,\ldots,g_n, h_1,\ldots, h_m\; ; \; R'_1, \ldots, R'_s, S'_1,\ldots,S'_m,
T_1,\ldots,T_n). 
$$
Enfin, si l'on enlève un à un les couples $(g_j;T_j)$,  il est
clair que la structure ne change pas non plus,
donc on obtient la \pn finie

\snic{ (h_1,\ldots, h_m\; ; \; R'_1, \ldots, R'_s, S'_1,\ldots,S'_m).}

%\sni
On peut reprendre ce raisonnement sous une forme matricielle dans le cas des
modules \pfz. Voici ce que cela donne.

Tout d'abord on constate que l'on  ne change pas la structure de $M$ lorsque l'on  fait
subir à
la \mpn $G$ une des transformations suivantes.\rdb

%-----------------begin enumerate -------
\begin{enumerate}\label{ManipMpns}
\item   Ajout d'une colonne nulle (ceci ne change pas le module des syzygies
entre des \gtrs fixés).
\item  Suppression d'une colonne nulle, sauf à obtenir une matrice vide.
\item   Remplacement de $G$, de type $q\times m$, par $G'$ de type
$(q+1)\times (m+1)$ obtenue à partir de $G$ en rajoutant une ligne nulle en
dessous puis une colonne à droite avec $1$ en position $(q+1,m+1)$, (ceci
revient à rajouter un vecteur parmi les \gtrsz, en indiquant sa dépendance par
rapport aux \gtrs précédents):

\snic{G\;\mapsto \;G'\;=\;
\cmatrix{
    G      &C           \cr
    0_{1,m}&1
}.}
\item   Opération inverse de la précédente, sauf à aboutir à une matrice
vide.
\item   Ajout à une colonne d'une \coli des autres colonnes (ceci ne change pas
le module des syzygies entre des \gtrs fixés).
\item   Ajout à une ligne d'une \coli des autres lignes,
(par exemple si nous notons $L_i$ la $i$-ième ligne, le remplacement de 
$L_1$ par  $L_1+\gamma L_2$ revient à remplacer le \gtr $g_2$ par $g_2-\gamma g_1$).
\item   Permutation de colonnes ou de lignes.
%\item   multiplication d'une colonne ou d'une
%ligne par un \elt \iv  (facultatif).
\end{enumerate}
%-----------------end enumerate --------

On voit ensuite que si $G$ et $H$ sont deux matrices de \pn d'un même module
$M$, on peut passer de l'une à l'autre au moyen des transformations décrites
ci-dessus. Un peu mieux: on voit que pour tout \sgr fini de $M$, on peut
construire à partir de $G$, en utilisant ces transformations, une \mpn de $M$
pour le nouveau \sgrz. Notez qu'en conséquence, un changement de base
de $\gA^q$ ou $\Ae m,$ qui correspond à la multiplication de $G$ (à gauche ou
à droite) par une matrice \ivz,  peut être réalisé par les opérations
décrites précédemment.
\\
Précisément, on obtient le résultat suivant.
%--- Lemma{lem pres equiv}------
\begin{lemma}
\label{lem pres equiv}
Soient deux matrices $G\in \gA^{q\times m}$ et $H\in \Ae {r\times n}$.
Alors \propeq
%-----------------begin item------------------
\begin{enumerate}\itemsep0pt
\item  Les matrices $G$ et $H$ présentent \gui{le même} module, \cad leurs conoyaux sont
isomorphes.
\item  Les deux matrices  de la figure ci-dessous sont
\elrt \eqvesz.\perso{Ne devrait-on pas mettre en
corolaire une bonne variante du lemme de Schanuel?}
\item  Les deux matrices  de la figure ci-dessous sont
\eqvesz.
\end{enumerate}
\vspace{-10pt}
%-----------------end item------------------
\begin{figure}[htbp]
{\small
\begin{center}
\[
\begin{array}{c|p{35pt}|p{15pt}|p{25pt}|p{20pt}|}
\multicolumn{1}{c}{} & \multicolumn{1}{c}{m} & \multicolumn{1}{c}{r} &
\multicolumn{1}{c}{q} & \multicolumn{1}{c}{n} \\
\cline{2-5}
\vrule height20pt depth13pt width0pt q\; & \hfil G \hfil &\hfil
0\hfil &\hfil 0\hfil &\hfil 0\hfil \\
\cline{2-5}
\vrule height15pt depth8pt width0pt r\; & \hfil 0  &\hfil
${\rm I}_{r}$\hfil & \hfil 0\hfil &\hfil0\hfil \\
\cline{2-5}
\end{array}
\]
\vspace{6pt}
\[
\begin{array}{c|p{35pt}|p{15pt}|p{25pt}|p{20pt}|}
\cline{2-5}
\vrule height20pt depth13pt width0pt q\; & \hfil 0 \hfil&\hfil
0\hfil &\hfil $\I_{q}$\hfil &\hfil 0\hfil \\
\cline{2-5}
\vrule height15pt depth8pt width0pt r\; &\hfil 0 \hfil&\hfil
0\hfil & \hfil 0\hfil & \hfil H\hfil \\
\cline{2-5}
\end{array}\]

\end{center}}
\caption{\emph{Les deux matrices}}
\label{fig}
\end{figure}

\end{lemma}
%--- end-lemma--------------------

\vspace{-4pt}
Comme conséquence du lemme \ref{factchangesgrmpf}, on obtient une reformulation plus abstraite de la
\cohc comme suit.

%--- Fact{factCohFA}----------------
\begin{fact}
\label{factCohFA}
Un anneau est \coh \ssi tout \itf est \pf (en tant que \Amoz).
Un \Amo est \coh \ssi tout sous-\mtf est \pfz.
\end{fact}
%--- end-fact-----------------------------------------

%:   subsec{Digression sur le calcul \agqz}
\subsec{Digression sur le calcul \agqz} \label{DigCalcAlg}
Outre leur rapport direct avec la résolution des \slis une autre raison de
l'importance des \mpfs est la suivante.

Chaque fois qu'un calcul \agq aboutit à un
\gui{résultat intéressant} dans un \Amo $M$, ce calcul n'a fait intervenir
qu'un nombre fini d'\elts $x_1$, \dots, $x_n$ de $M$ et un nombre fini de relations
entre les $x_j$, de sorte qu'il existe un  \mpf $P=\Ae{n}\sur{R}$ et un
\homo surjectif $\theta:P\to x_1\gA+\cdots+x_n\gA\subseteq M$ qui envoie les $e_j$ sur les $x_j$ (où $e_j$ désigne la classe modulo $R$ du $j$-ième
vecteur de la base canonique de~$\Ae{n}$),
et tel que le \gui{résultat intéressant} avait déjà lieu dans $P$
pour les~$e_j$.

En langage plus savant on exprime cette idée comme ceci.
\emph{Tout \Amo est limite inductive filtrante  d'\Amos \pfz.\label{factLimIndFiltPf}
}
\\
Mais cet énoncé
nécessite un traitement un peu subtil en \comaz, et nous ne faisons donc que
signaler son existence.

%%%%%%%%%%%%%%%%%%%%%%%%%%%%%%%%%%%%%%%%%%%%%%%%%%%%%%%%%%%%%%%%%%%%%%%%%%%
%%%%%%%%%%%%%%%%%%%%%%%%%%%%%%%%%%%%%%%%%%%%%%%%%%%%%%%%%%%%%%%%%%%%%%%%%%%
%%%%%%%%%%%%%%%%%%%%%%%%%%%%%%%%%%%%%%%%%%%%%%%%%%%%%%%%%%%%%%%%%%%%%%%%%%%
%%%%%%%%%%%%%%%%%%%%%%%%%%%%%%%%%%%%%%%%%%%%%%%%%%%%%%%%%%%%%%%%%%%%%%%%%%%
% section{Idéaux \pfz, et un exemple en \gmt}
\section{Idéaux \pfz} \label{secIdPf}

\vspace{3pt}
On considère un anneau $\gA$ et un \sgr $(\an)=(\ua)$ pour un \itf $\fa$ de $\gA$.
On s'intéresse à la structure de \Amo de $\fa$. 

%:--- SUBsection Relations triviales et suites reg {secRelTrivSeqReg}
\subsec{Syzygies triviales %et suites régulières
} \label{secRelTrivSeqReg}

Parmi les syzygies entre les $a_i$ figurent ce que l'on  appelle les \emph{syzygies triviales} (ou \emph{relateurs triviaux} si on les voit
comme des \rdes \agqs sur $\gk$ lorsque $\gA$ est une \klgz): 

\snic{a_ia_j-a_ja_i=0\;$ pour $\;i\neq j.}

%\sni
Si $\fa$ est \pfz, on pourra toujours prendre une \mpn de $\fa$ pour le \sgr $(\ua)$
sous la forme

\snic{W=[\,R_\ua\mid U\,],}

%\sni
où $R_\ua$ est \gui{la} \emph{matrice des syzygies triviales} (l'ordre des colonnes est sans importance),
 de format
${n\times  n(n-1)/2}$. Par exemple, pour $n = 5$%
\index{syzygie!triviale}%
\index{matrice!des syzygies triviales}
$$
\preskip.4em \postskip.4em 
R_{\ua} = \Cmatrix{.3em} {
a_2 & a_3  &  0  & a_4  &  0  &   0 & a_5 & 0 & 0 & 0\cr
-a_1&   0 &  a_3 &   0 &  a_4 &   0& 0 & a_5 & 0 & 0\cr
 0 & -a_1 & -a_2 &   0 &   0 &  a_4&0  &   0 & a_5&   0\cr
 0 &   0 &   0 & -a_1 & -a_2 & -a_3&  0  &   0 & 0& a_5\cr
 0 &   0 &   0 & 0 &   0 &   0 & -a_1 & -a_2 & -a_3& -a_4}. 
$$
 
%:     Lemma{lemDnRz}
\begin{lemma}\label{lemDnRz} \emph{(Idéaux déterminantiels de la matrice des syzygies triviales)}
Avec les notations ci-dessus, on a les résultats suivants.
\begin {enumerate}
\item
$\cD_n(R_\ua)  = \{0\}$.
\item Si $1\leq r <n$, alors
$\cD_r(R_\ua) = \fa^r$ et 
$$\preskip.4em \postskip.4em 
\fa^r+\cD_r(U)\;\subseteq\; \cD_r(W)  \;\subseteq\;\fa+ \cD_r(U). 
$$
En particulier, on a l'\eqvc
$$\preskip.4em \postskip.4em 
1\in\cD_{\gA,r}(W)\iff 1\in\cD_{\gA/\!\fa,r}(\ov U), \;\;\hbox{  où } \;\ov U=U\mod\fa. 
$$
\item
$\cD_n(W)  = \cD_n(U)$.
\end {enumerate}
 \end{lemma}
%--------- fin lemma ---------------------------------------------- 
%
\begin{proof}
\emph {1.} Il s'agit d'\idas et l'on peut prendre pour $a_1$, \dots,~$a_n$
des \idtrs sur $\ZZ$. 
Comme $[\,a_1\;\cdots\;a_n\,] \cdot R_\ua = 0$, on obtient l'\egtz~\hbox{$\cD_n(R_\ua) \,[\,a_1\;\cdots\;a_n\,]  = 0$}. On conclut puisque $a_1$ est \ndzz.

\emph {2.} L'inclusion $\cD_r(R_\ua)\subseteq\fa^r$ est 
évidente pour tout $r\geq0$. 
Pour l'inclusion réciproque prenons par exemple $r=4$ et $n\geq5$ et 
montrons que
$$\preskip.4em \postskip.5em
\so{a_1^4,\, a_1^3a_2,\, a_1^2a_2^2,\, a_1^2a_2a_3,\, a_1a_2a_3a_4}\subseteq\cD_4(R_\ua).
$$
Il suffit de considérer les matrices dessinées ci-après 
(nous avons supprimé les $0$ et remplacé~$\pm a_i$ par~$i$ 
pour mieux voir la structure)
extraites de $R_\ua$, et les mineurs extraits sur les $4$
dernières lignes.
\[\preskip.4em \postskip.4em 
\begin{array}{ccccc} 
\cmatrix{ 2 &  3  &  4  &  5  \cr
         1 &  &  &  \cr  
        &  1 &   &  \cr 
        &  &  1 &   \cr 
        &  &  &  1},\,  
&   
\cmatrix{ 2 &  3  & 4   &   \cr
         1 &  &  &  5 \cr  
        &  1 &   &  \cr 
        &  &  1 &   \cr 
        &  &  &  2},\,  
& 
\cmatrix{ 2 &  3  &   &   \cr
         1 &  &  4 &  5 \cr  
             &  1 &   &  \cr 
             &      &  2 &   \cr 
             &      &      &  2}, \, 
\\[10mm]   
\cmatrix{ 2 &  3  &   &   \cr
         1 &  &  4 &  \cr  
             &  1 &   &  5 \cr 
             &      &  2 &   \cr 
             &      &      &  3},  
&   
\cmatrix{ 2 &   &   &   \cr
         1 &  3 &  &  \cr  
        &  2 &  4  &  \cr 
        &  &  3 &  5  \cr 
        &  &  &  4}.  
\end{array}
\]
L'inclusion $\fa^r+\cD_r(U)\subseteq \cD_r(W)$ résulte de
$\cD_r(R_\ua)+\cD_r(U)\subseteq \cD_r(W)$ et de l'\egt $\cD_r(R_\ua) = \fa^r$.
 L'inclusion $\cD_r(W)  \subseteq\fa+ \cD_r(U)$ est \imdez. 
Enfin l'\eqvc
finale résulte des inclusions précédentes et de l'\egt 
$$
\cD_{\gA/\!\fa,r}(\ov U)=\pi_{\gA,\fa}^{-1}\big(\fa+\cD_r(U)\big).
$$

\emph {3.} On doit montrer que si une matrice $A\in\Mn(\gA)$  extraite de $W$ 
contient une colonne dans $R_\ua$, alors $\det A=0$.
Prenons par exemple la première colonne de $A$ égale à la première colonne de $R_\ua$: $\tra{[\,a_2\;-a_1\;0\;\cdots\;0\,]}$. 
Le lemme \ref{lem acheval} ci-après implique,  lorsque $z_i=a_i$,
$\det A=0$, car les $s_j$ sont~nuls.
\end{proof}
%

%:2015 ci-dessous
\rdb\label{NOTAProdScal}
Rappelons que $A_{\alpha,\, \beta}$ est la sous-matrice de $A$ extraite
sur les lignes $\alpha$ et les colonnes $\beta$.
Introduisons aussi la notation \gui{produit scalaire}

\snic{\scp{x}{y}\,\eqdefi\,\sum_{i=1}^n x_iy_i}

%\sni
pour deux vecteurs colonnes $x$ et $y$.

%:     Lemma{lem acheval}
\begin {lemma} \label {lem acheval} 
Soient $A %= (a_{ij}) 
\in \Mn(\gA)$,  $A_j=A_{1..n,j}$,
et $z = \tra{[\,z_1\; \cdots\; z_n\,]}\in \Ae{n\times 1}$
avec $A_1 = \tra{[\,z_2\; -z_1\; 0\;\cdots\;0\,]}$.
En posant $s_j = \scp {z} {A_j}$ pour $j\in \lrb {2..n}$, on~a
%$$

\snic{
\dsp\det A  = \som_{j=2}^n (-1)^j \,s_j \,\det(A_{3..n,\, 2..n \setminus \{j\}}).
}
%\sni 
%$$

En particulier, $\det A  \in \gen {s_2, \ldots, s_n}$.

\end {lemma}
\begin{proof}
Notons $B=A_{3..n,2..n}$, $B_j=A_{3..n,j}$ et $B_{\hat \jmath}=A_{3..n,\, 2..n
\setminus \{j\}}$.  Le développement de Laplace du \deter de $A$ selon les deux
premières lignes donne l'\egtz:
$$\preskip.2em \postskip.2em \ndsp
\mathrigid1mu 
\det A  = \sum\limits_{j=2}^n (-1)^j \,\dmatrix{\phantom-z_2&a_{1j}\cr -z_1&a_{2j}} \,\det(B_{\hat \jmath}) 
= \sum\limits_{j=2}^n (-1)^j \,(z_1a_{1j}+z_2a_{2j}) \,\det(B_{\hat \jmath}). 
$$
 L'écart entre cette \egt et l'\egt voulue est

\snic{
\qquad \qquad 
\som_{j=2}^n (-1)^j \,(z_3a_{3j}+\cdots+z_na_{nj}) \,\det(B_{\hat \jmath}).\qquad \qquad(*)
}
%$$

%\sni
La syzygie de Cramer entre les colonnes d'une matrice  avec $m=n_2$,  donne pour $B$ les \egts  
$$\preskip.4em \postskip.2em\ndsp 
\som_{j=2}^n (-1)^j \det(B_{\hat \jmath})\, B_j=0, \hbox{ a fortiori } \som_{j=2}^n (-1)^j \scp{y}{B_j}\,\det(B_{\hat \jmath})=0,
 $$ 
  pour n'importe quel vecteur $y\in\Ae{(n-2)\times 1}$.  
En prenant $y=\tra{[\,z_3\; \cdots\; z_n\,]}$, on
voit que l'écart $(*)$ est nul.
\end{proof}
%

%:subsec{Suites régulières}
\subsec{Suites régulières}

%:  Definition{defSeqReg}------
\begin{definition}
\label{defSeqReg} ~
Une suite $(a_1, \ldots, a_k)$ dans  un anneau $\gA$
est
\ixg{régulière}{suite ---}{reguliere}
si chaque $a_i$ est \ndz dans l'anneau
$\aqo{\gA}{a_j\,;\,j<i}$.%
\index{suite!reguliere@régulière}
\end{definition}
%--- end-definition------------------------------------

\rem Nous avons retenu ici la \dfn de Bourbaki. La plupart des auteurs
réclament en outre que l'\id $\gen{a_1,\ldots ,a_k}$ ne contienne pas $1$. 
Mais l'expérience montre que cette négation introduit des complications inutiles dans les énoncés et les \demsz. \eoe

\smallskip Comme premier exemple, pour tout anneau $\gk$, la suite $(X_1, \ldots , X_k)$  est régulière dans $\gk[X_1,\ldots ,X_k]$.

 Notre but est de monter qu'un \id engendré par une suite 
régulière est un module \pfz.

 Nous établissons d'abord un petit lemme et une proposition.
\\
Rappelons qu'une matrice $M = (m_{ij}) \in \Mn(\gA)$ est dite
\emph {alternée} si c'est la matrice d'une forme bi\lin alternée,
i.e. $m_{ii} = 0$ et $m_{ij} + m_{ji} = 0$ pour~$i$,~$j\in\lrbn$.%
\index{matrice!alternée}\index{alternée!matrice ---}

Le \Amo des matrices alternées
 est libre de rang  ${n(n-1)}\over 2$ et admet une base naturelle.
Par exemple, pour $n = 3$,
$$\preskip.4em \postskip.4em
\Cmatrix{.3em} {0 & a & b \cr -a & 0 & c \cr -b & -c & 0 \cr} =
a\cmatrix {0 & 1 &  0 \cr -1 & 0 & 0 \cr 0 & 0 & 0 \cr} +
b\cmatrix {0 & 0 & 1 \cr 0 & 0 & 0 \cr -1 & 0 & 0 \cr} +
c\cmatrix {0 & 0 & 0 \cr 0 & 0 & 1 \cr 0 & -1 & 0 \cr}.
$$

%:    lemma{PetitLemmeAlterne}
\begin {lemma}\label{PetitLemmeAlterne}
Soit $a = \tra{[\,\ua\,]} = \tra{[\,a_1\;\cdots\;a_n\,]} \in \Ae {n\times 1}$.
\begin {enumerate}
\item
Soit $M  \in \Mn(\gA)$ une matrice alternée;
on a $\scp {Ma}{a} = 0$.

\item
Un $u \in \Ae {n\times 1}$ est dans $\Im R_\ua$ \ssi il
existe une matrice alternée $M\in \Mn(\gA)$
telle que $u = Ma$.
\end {enumerate}
\end {lemma}
%%%%%%%%%%%%%%%%%%%%%%%%%%%%%%%%%%%%%%%%%
\begin {proof}\emph{1.} En effet, $\scp {Ma}{a} = \varphi(a,a)$, où $\varphi$
est une forme bi\lin alternée.

\emph{2.} Par exemple, pour la première colonne de $R_\ua$ avec $n=4$, on~a:
$$\preskip.4em \postskip.2em 
\cmatrix {0 & 1 & 0 & 0\cr -1 & 0 & 0 & 0\cr 0 & 0 & 0 & 0\cr 0 & 0 & 0 & 0\cr}
\crmatrix {a_1\cr a_2\cr a_3\cr a_4\cr} =
\Cmatrix{.2em} {a_2\cr -a_1\cr 0\cr 0\cr}, 
$$
et les ${n(n-1)}\over 2$ colonnes de $R_\ua$ correspondent ainsi
aux ${n(n-1)}\over 2$ matrices alternées formant la base naturelle
du \Amo des matrices alternées de~$\Mn(\gA)$.
\end {proof}
%%%%%%%%%%%%%%%%%%%%%%%%%%%%%%%%%%%%%%%%%

%:     Proposition{prop1sregpf}
\begin{proposition}\label{prop1sregpf}
Soit $(z_1, \ldots, z_n)=(\uz)$ une suite régulière
d'\elts de~$\gA$ et $z = \tra{ [\,z_1 \;\cdots\; z_n\,]} \in \Ae{n\times 1}$.  Si
$\scp {u}{z} = 0$, il existe une matrice alternée~$M  \in
\Mn(\gA)$ telle que $u = Mz$, et donc $u\in\Im R_\uz$. 
\end{proposition}
%--------- fin proposition ---------------------------------------------- 
%\rem La \dem n'utilise pas la deuxième condition dans la \dfn 
% des suites régulières.
%\eoe
\begin{proof}
On raisonne par \recu sur $n$. Pour $n = 2$, on part de $u_1z_1 + u_2z_2 =
0$. Donc $u_2z_2 = 0$ dans $\aqo{\gA}{z_1}$, et puisque $z_2$ est
\ndz modulo $z_1$, on~a~$u_2 = 0$ dans $\aqo{\gA}{z_1}$, disons
$u_2 = -az_1$ dans $\gA$. Il vient $u_1z_1 -az_2z_1 = 0$, et comme~$z_1$ est
\ndzz, $u_1 = az_2$, ce qui s'écrit $\Cmatrix{.3em} {u_1\cr u_2\cr} =
\crmatrix {0 & a \cr -a & 0\cr} \cmatrix {z_1\cr z_2\cr}$.

 Pour $n+1$ ($n\geq2$), 
 on part de $u_1z_1 + \cdots + u_{n+1}z_{n+1} = 0$.  En utilisant
le fait que $z_{n+1}$ est \ndz modulo $\gen{\zn}$,
on obtient $u_{n+1} \in \gen{\zn}$, ce que l'on écrit
$a_1z_1 + \cdots + a_{n}z_{n} + u_{n+1} = 0$. D'où:

\snic{
(u_1 - a_1z_{n+1}) z_1 + \cdots + (u_{n} - a_{n}z_{n+1}) z_{n} = 0.}

%\sni
Par \hdrz, 
on sait construire $M \in \MM_{n}(\gA)$ alternée
avec:
$$\preskip.2em \postskip.4em
\mathrigid1mu 
\cmatrix {u_1 - a_1z_{n+1}\cr \vdots\cr u_{n} - a_{n} z_{n+1}\cr} =
M \cmatrix {z_1\cr \vdots\cr z_{n}\cr},
\hbox { i.e. }
\cmatrix {u_1\cr \vdots\cr u_{n}\cr} =
M \cmatrix {z_1\cr \vdots\cr z_{n}\cr} +
z_{n+1}\cmatrix {a_1\cr \vdots\cr a_{n}\cr}.
$$
Et l'on obtient le résultat voulu:
$$\preskip.4em \postskip.2em
\Cmatrix{.3em} {u_1\cr \vdots\cr u_{n}\cr u_{n+1}\cr} =
\cmatrix {
     &       &          & a_1  \cr
     & M     &          & \vdots \cr
     &       &          & a_{n} \cr
-a_1 & \dots & -a_{n} & 0 \cr}
\cmatrix {z_1\cr \vdots\cr z_{n}\cr z_{n+1}\cr}.
$$

\vspace{-10pt}
\end {proof}

%:     theorem{propsregpf}
\begin {theorem}\label{propsregpf}
Si $(z_1, \ldots, z_n)$ est une suite régulière
d'\elts de $\gA$, l'\id $\gen{\zn}$ est un \Amo \pfz.
Plus \prmtz, on a la suite exacte

\snic{\Ae {n(n-1)/2} \vvvers {R_\uz}  \Ae n
\vvvvvers {(z_1, \ldots, z_n)}  \gen{\zn} \lora 0.}
\end {theorem}
%%%%%%%%%%%%%%%%%%%%%%%%%%%%%%%%%%%%%

\rem
Les objets définis ci-dessus constituent une introduction au premier étage
du \emph{complexe de Koszul} descendant de $(\zn)$.
\eoe

\begin{proof}
Cela résulte de la proposition \ref{prop1sregpf} et du lemme \ref{PetitLemmeAlterne}.
\end{proof}
%

%:--- subsection{Un exemple en \gmt}
%\penalty-5000
\subsec{Un exemple en \gmt} \label{secExempleGeo}

Voici pour commencer une évidence fort utile.
%:     PropDef{prdfCaracAlg}
\begin{propdef}\label{prdfCaracAlg}
\emph{(Caractères d'une \algz)} 
\\
 Soit $\imath:\gk\to\gA$  une \algz.
\begin{itemize}
\item Un \homo de \klgs $\varphi:\gA\to\gk$ est appelé un 
\ixc{caractère}{d'une algèbre}.
\item  Si $\gA$ 
possède un caractère $\varphi$, alors   $\varphi\circ \imath=\Id_\gk$, $\imath\circ \varphi$ est un \prr et $\gA=\gk.1_\gA\oplus\Ker\varphi$. 
En particulier, on peut identifier $\gk$ et $\gk.1_\gA$.
\end{itemize}
\end{propdef}
\facile

Soit maintenant $(\uf)=(f_1, \ldots, f_s)$ un \syp sur un anneau~$\gk$,
avec les $f_i \in\kuX
= \kXn$.  On note 
$$\gA = \kxn = \aqo{\kuX}{\uf}.
$$

Dans ce paragraphe, de façon informelle, nous dirons que $\gA$
est  l'\emph{anneau de la \vrt affine $\uf = \uze$}. 

Pour l'\alg $\gA$, les caractères $\varphi:\gA\to\gk$ sont donnés par les zéros
dans~$\gk^n$ du \syp $(\lfs)$: 
$$
\uxi=(\xin)=\big(\varphi(x_1),\ldots,\varphi(x_n)\big),\quad\uf(\uxi)=\uze.
$$
Dans ce cas, on dit que $\uxi\in\gk^ n$ est un point de la \vrt $\uf = \uze$.

L'\id 
$$
\fm_\uxi\eqdefi\gen{x_1-\xi_1,\ldots,x_n-\xi_n}_\gA
$$ 
est appelé 
l'\emph{\id du point $\uxi$ dans la \vrtz}. 
On a alors comme cas particulier de la proposition \ref{prdfCaracAlg}: $\gA=\gk\oplus\fm_\uxi$, avec $\fm_\uxi=\Ker\varphi$.
\index{ideal@idéal!d'un point}

{Dans ce paragraphe on montre que l'\id $\fm_\uxi$  est un \Amo \pf en
explicitant une \mpn
pour le \sgr $(x_1-\xi_1, \ldots, x_n-\xi_n)$.}

 Par translation, il suffit de traiter le cas où $\uxi=\uze$, ce que nous supposons désormais.
\\
 Le cas le plus simple, celui pour lequel il n'y a aucune
\eqnz, a déjà été traité dans le \thrf{propsregpf}.
\\
Observons que tout $f \in \kuX$ tel que $f(\uze) =
0$ s'écrit, de plusieurs manières, sous la forme
$$\preskip.0em \postskip.3em 
f = X_1u_1 + \cdots + X_n u_n, \qquad  u_i \in \kuX. 
$$
Si $X_1v_1 + \cdots + X_nv_n$ est une autre écriture de $f$, on obtient
par soustraction une syzygie entre les $X_i$ dans $\kuX$, et donc:

\snic{
\tra{ [\, v_1\;\cdots\; v_n\,]} - \tra{[\, u_1\;\cdots\; u_n\,]} \in \Im {R_{\uX}}.}

%\sni

Pour le \syp $(f_1, \ldots, f_s)$, on définit ainsi (de manière non unique) une famille de \pols $(u_{ij})_{i \in\lrbn, j \in\lrbs}$, avec $f_j = \sum_{i=1}^n
X_iu_{ij}$. Ceci donne  une
matrice $U(\uX) = (u_{ij})$ et  son image $U(\ux) = \big(u_{ij}(\ux)\big) \in \Ae{n\times s}$.

%:     theorem{thidptva}
\begin {theorem} \label{thidptva}
Pour un \syp sur un anneau $\gk$ et un zéro~$\uxi\in\gk^n$, l'\id $\fm_\uxi$
du point $\uxi$ est un \Amo \pfz.%
\index{systeme polynomial@\sypz}
\\
Plus \prmtz, avec les notations précédentes, pour le cas $\uxi=\uze$   la
matrice $W = [\,R_\ux\,|\,U(\ux)\,]% \in \Ae{n\times m}
$  est une \mpn
de l'\id $\fm_\uze$ pour le \sgr $(\xn)$. Autrement dit on a une suite exacte

\snic{
\Ae{m} \vvvvvers {[\,R_\ux\,|\,U\,]}
 \Ae n
\vvvvvers {(x_1, \ldots, x_n)}  \fm_\uze\lora 0 \quad\quad(m = {n(n-1) \over 2} + s).
}
\end {theorem}
%----- end theorem -----------------------------------------
\begin {proof}
   Prenons par exemple $n = 3$, $s = 4$, $X = \tra {[\,X_1 \;X_2\;X_3 \,]}
 $
et pour économiser les indices écrivons  $f_1 =
X_1a_1 + X_2a_2 + X_3a_3$,  et  $f_2$, $f_3$, $f_4$ en utilisant les lettres $b$, $c$, $d$. On prétend avoir la \mpn suivante
pour le \sgr $(x_1,x_2,x_3)$ de $\fm_\uze$:
$$\preskip.4em \postskip.4em
\cmatrix {
x_2  & x_3  & 0    & a_1(\ux) & b_1(\ux) & c_1(\ux) & d_1(\ux)\cr
-x_1 & 0    & x_3  & a_2(\ux) & b_2(\ux) & c_2(\ux) & d_2(\ux)\cr
0    & -x_1 & -x_2 & a_3(\ux) & b_3(\ux) & c_3(\ux) & d_3(\ux)\cr}.
$$
 On définit $A = \tra {[\,a_1 \;a_2\;a_3 \,]}$  dans $\kuX^3$
 (ainsi que $B,C,D$) de sorte que:

 \snic{f_1 =\scp {{A}}{X},\; f_2 = \scp {{B}}{X} \;\;\ldots.}

Considérons une \syzy $v_1(\ux)x_1 + v_2(\ux)x_2 + v_3(\ux)x_3 =0$
dans~$\gA$. On la remonte dans $\kuX$:

 \snic{v_1X_1 + v_2X_2 + v_3X_3 \equiv 0 \mod \gen{\uf}.}

Ce que l'on écrit

 \snic{
v_1X_1 + v_2X_2 + v_3X_3 = \alpha f_1 + \beta f_2 + \gamma f_3 + \delta f_4,
\qquad \alpha,\, \beta,\, \gamma,\, \delta \in \kuX.
}

Donc, avec $V = \tra {[\,v_1 \;v_2\;v_3 \,]}$,
$V - (\alpha{A} + \beta{B} + \gamma{C} + \delta{D})$
est une \syzy pour~$(X_1, X_2, X_3)$, ce qui implique par la proposition~\ref{prop1sregpf}

 \snic{
V - (\alpha{A} + \beta{B} + \gamma{C} + \delta{D}) \in
\Im {R_{\uX}}.
}

 %\sni
Ainsi, $V \in \Im \,[\,R_{\uX}\,|\,U(\uX)\,]$, et 
$\tra {[\, v_1(\ux)\;v_2(\ux)\;v_3(\ux)\,] } \in \Im \,[\,R_\ux\,|\,U(\ux)\,]$.
\end {proof}

%%%%%%%%%%%%%%%%%%%%%%%%%%%%%%%%%%%%%%%%%%%%%%%%%%%%%%%%%%%%%%%%%%%%%%%%%%%
%%%%%%%%%%%%%%%%%%%%%%%%%%%%%%%%%%%%%%%%%%%%%%%%%%%%%%%%%%%%%%%%%%%%%%%%%%%
%%%%%%%%%%%%%%%%%%%%%%%%%%%%%%%%%%%%%%%%%%%%%%%%%%%%%%%%%%%%%%%%%%%%%%%%%%%
%%%%%%%%%%%%%%%%%%%%%%%%%%%%%%%%%%%%%%%%%%%%%%%%%%%%%%%%%%%%%%%%%%%%%%%%%%%
%--- Sec{Categorie}--secCatMpf------
\section{Catégorie des \mpfs} \label{secCatMpf}
%--------------------

La catégorie des \mpfs sur $\gA$ peut être construite à partir de la
catégorie des modules libres de rang fini sur $\gA$ par un procédé purement
catégorique.
%-----------------begin enum------------------
\begin{enumerate}
\item  Un \mpf $M$ est décrit par un triplet 

\snic{(\rK_M,\alb\rG_M,\rA_M),}

%\sni
où
$\rA_M$ est une \ali entre les modules libres
de rangs finis~$\rK_M$ et~$\rG_M$.
On a~$M\simeq \Coker\rA_M$ et~$\pi_M:\rG_M\rightarrow M$ est l'\ali surjective
de noyau $\Im\rA_M$.
La matrice de l'\aliz~$\rA_M$ est une \mpn de~$M$.

\item  Une \ali $\varphi$ du module $M$ (décrit par $(\rK_M,\rG_M,\rA_M)$) vers
le module  $N$ (décrit par $(\rK_N,\rG_N,\rA_N)$)
est décrite par deux \alis
$\rK_\varphi:\rK_M\rightarrow \rK_N$ et $\rG_\varphi:\rG_M\rightarrow \rG_N$
soumises à la relation de commutation
$\rG_\varphi\circ\rA_M=\rA_N\circ\rK_\varphi$.

\centerline{\xymatrix {
\rK_M \ar[r]^{\rA_M} \ar[d]_{\rK_\varphi} & \rG_M \ar[d]^{\rG_\varphi}
\ar@{->>}[r]^{\pi_M}
                   & M \ar[d]^{\varphi} \\
\rK_N \ar[r]_{\rA_N}                    & \rG_N \ar@{->>}[r]_{\pi_N}
                   & N \\}
}
\item  La somme de deux \alis $\varphi$ et $\psi$ de $M$ vers $N$
repré\-sentées par $(\rK_\varphi,\rG_\varphi)$ et  $(\rK_\psi,\rG_\psi)$ est
repré\-sentée par $(\rK_\varphi+\rK_\psi,\alb \rG_\varphi+\rG_\psi)$.
L'\ali  $a\varphi$ est repré\-sentée par $(a\rK_\varphi,\alb a\rG_\varphi)$.

\item  Pour représenter la composée de deux \alisz, on compose leurs
représentations.

\item  Enfin l'\ali  $\varphi$ de $M$ vers $N$ représentée par
$(\rK_\varphi,\rG_\varphi)$ est nulle
\ssi il existe $Z_\varphi:\rG_M\rightarrow \rK_N$ 
vérifiant~$\rA_N\circ Z_\varphi=\rG_\varphi$.
\end{enumerate}
%-----------------end enum------------------

Ceci montre que les problèmes concernant les \mpfs peuvent toujours être
interprétés comme des \pbs à propos de matrices, et se ramènent
\label{exMpf} souvent à des problèmes de résolution de \slis sur $\gA$.
\\
Par exemple, si l'on donne $M$, $N$ et $\varphi$ et si l'on cherche
une \ali $\sigma:N\rightarrow M$ vérifiant $\varphi\circ\sigma=\Id_N$,
il faut trouver des \alis $\rK_\sigma:\rK_N\rightarrow \rK_M$,
$\rG_\sigma:\rG_N\rightarrow \rG_M$  et $Z:\rG_N\rightarrow \rK_N$ qui 
vérifient:
%-----------------begin $$----------------
$$ 
\rG_\sigma\circ\rA_N=\rA_M\circ\rK_\sigma \quad {\rm et}\quad
\rA_N\circ Z=\rG_\varphi\circ\rG_\sigma-\Id_{\rG_N}.
$$
%-----------------end $$------------------
\perso{mettre un dessin}
Ceci n'est autre qu'un \sli ayant pour inconnues les \coes des matrices
des \alis $\rG_\sigma$, $\rK_\sigma$ et~$Z$.

De manière analogue,
si l'on se donne $\sigma:N\to M$ et si l'on se pose la question de savoir
s'il existe $\varphi :M\to N$ vérifiant  $\varphi\circ\sigma=\Id_N$,
on devra résoudre un \sli dont les inconnues sont les  \coes des matrices
des \alis $\rG_\varphi$, $\rK_\varphi$ et~$Z$.

De même, si l'on se donne $\varphi:M\to N$  et si l'on se pose la question de savoir
si $\varphi$ est \lnlz, on doit savoir s'il existe $\sigma:N\to M$ 
vérifiant~$\varphi \circ\sigma\circ\varphi =\varphi $, et l'on obtient un \sli ayant pour
inconnues les \coes des matrices
de $\rG_\sigma$, $\rK_\sigma$ et~$Z$.

On en déduit des \plgs correspondants.

%--- Principe local global concret{plcc.scinde}----
\begin{plcc}
\label{plcc.scinde} 
\emph{(Pour certaines \prts des \alis entre \mpfsz)}\\
Soient $S_1$, $\dots$, $S_n$ des \moco de $\gA$,  $\varphi :M\to N$
une \ali entre \mpfsz.
Alors \propeq
\begin{enumerate}
\item  L'\Ali $\varphi$ admet un inverse à gauche (resp. admet un inverse à
droite, resp. est \lnlz).
\item Pour $ i\in\lrbn,$
l'application $\gA_{S_i}$-\lin $\varphi_{S_i}:M_{S_i}\to N_{S_i}$ admet un inverse à  
gauche (resp. un
admet un inverse à droite, resp. est \lnlz).
\end{enumerate}
\end{plcc}
%--- end-plcc-----------------------------------------

%--- Sec{Stabilité}--secStabPf-
\penalty-2500
\section{Propriétés de stabilité} \label{secStabPf}
%--------------------

%--- Proposition{propPfInter}--------
\begin{proposition}
\label{propPfInter}
Soient $N_1$ et $N_2$  deux sous-\Amos \tf \linebreak 
d'un \Amo  $M$.
Si $N_1+N_2$ est \pfz, alors  $N_1\cap N_2$ est \tfz.
\end{proposition}
%--- end-proposition----------------------------------------
%
\begin{proof}
On peut reprendre presque mot pour mot la \dem du point
\emph{1} du \thref{propCoh4}
(condition \ncrz).
\end{proof}
%

%:2012 le point 3 est plus précis
%--- Proposition{propPfSex}----------
\begin{proposition}
\label{propPfSex}
Soit $N$ un sous-\Amo de $M$ et $P=M/N$.
%-----------------begin enum------------------
\begin{enumerate}
\item Si $M$ est \pf et   $N$  \tfz,  $P$ est \pfz.
\item Si $M$ est \tf et  $P$  \pfz,   $N$ est \tfz.
\item Si $P$ et $N$ sont \pfz, $M$ est \pfz. Plus \prmtz, si $A$ et
$B$ sont des \mpns pour $N$ et $P$, on a une \mpn 
%$D=\cmatrix{A&C\cr 0&B}$
$D=\blocs{1}{.8}{.7}{.6}{$A$}{$C$}{$0$}{$B$}$
pour $M$.
\end{enumerate}
%-----------------end enum------------------
\end{proposition}
%--- end-proposition----------------------------------------
%
\begin{proof}
\emph{1.} On peut supposer que $M=\gA^p/F$ avec $F$ de type fini. 
Si $N$ est \tfz, il s'écrit  $N=(F'+F)/F$ où $F'$ est \tfz, donc $P\simeq \gA^p/(F+F')$.
\\
 \emph{2.} On écrit $M=\gA^p/F$,
et $N=(F'+F)/F$. On a $P\simeq \gA^p/(F'+F)$, \linebreak 
donc $F'+F$ est \tf
(section \ref{sec pf chg}), 
et $N$ \egmtz.
\\
 \emph{3.} Soient $x_1$, \ldots, $x_m$ des
\gtrs de $N$ et $x_{m+1}$, \ldots, $x_n$ des \elts de $M$ dont les classes modulo $N$
engendrent $P$. Toute syzygie sur $(\ov{x_{m+1}}, \ldots, \ov{x_n})$ dans $P$
donne une syzygie sur $(\xn)$ dans $M$. De même, toute syzygie sur $(\xn)$ dans $M$ donne une syzygie sur $(\ov{x_{m+1}}, \ldots, \ov{x_n})$ dans $P$.
\linebreak 
Si $A$ est une \mpn
de $N$ pour   $(\xm)$ et si  $B$ est une \mpn
de $P$ pour  $(\ov{x_{m+1}},\ldots ,\ov{x_n})$, on obtient
donc une \mpn $D$ de $M$ pour $(\xn)$ du format voulu.
\end{proof}
On notera que dans la \dem du point \emph{2} 
les sous-modules $F$ et $F'$ ne sont pas \ncrt \tfz.

%--- Sec{Cohérence et \pn finie}
\subsection*{Cohérence et \pn finie}
\addcontentsline{toc}{subsection}{Cohérence et \pn finie}

Les propositions \ref{propCoh1} et \ref{propCohfd1}
(lorsque l'on  prend $\gA$  comme \Amoz~$M$)
se relisent sous la forme du \tho suivant.

%:   Theorem{propCoh2}-------
\begin{theorem}
\label{propCoh2}
Sur un \cori tout  \mpf est \cohz.
Sur un \cori \fdi tout  \mpf est \coh \fdiz.
\end{theorem}
%--- end-theorem-----------------------------------------

%:     Proposition{propCohpfKer}
\begin{proposition}\label{propCohpfKer}
Soit $\gA$ un \cori et $\varphi:M\to N$ une \ali entre \Amos \pfz,
alors $\Ker\varphi$, $\Im\varphi$ et~$\Coker\varphi$ sont des \mpfsz.
\end{proposition}

%--- Proposition{propCohSex}---------
\begin{proposition}
\label{propCohSex}
Soit $N$ un sous-\Amo \tf de $M$.
%-----------------begin enum------------------
\begin{enumerate}
\item Si $M$ est \cohz, $M/N$ est \cohz.
\item Si $M/N$ et $N$ sont \cohsz, $M$ est \cohz.
\end{enumerate}
%-----------------end enum------------------
\end{proposition}
%--- end-proposition----------------------------------------
%-----------------begin proof------------------
\begin{proof}
\emph{1.} On considère un sous-\mtf $P=\gen{\ov{x_1},\ldots
,\ov{x_\ell}}$
de $M/N$. Alors $P\simeq (\gen{\xl}+N)/N$. On conclut par la
proposition \ref{propPfSex} qu'il est \pfz.

\emph{2.} Soit $Q$ un sous-\mtf de $M$.
Le module $(Q+N)/N$ est \tf dans $M/N$ donc \pfz.
Puisque $(Q+N)/N$ et $N$ sont \pfz, $Q+N$  \egmt (proposition \ref{propPfSex}).
Donc $Q\cap N$ est \tf (proposition \ref{propPfInter}). Puisque $N$ est
\cohz, $Q\cap N$ est \pfz. Puisque $Q/(Q\cap N)\simeq  (Q+N)/N$ et  $Q\cap N$
sont \pfz, $Q$ est \pf (proposition~\ref{propPfSex}).
\end{proof}
%-----------------end proof------------------

%:--- Sec{Produit tensoriel}
\subsect{Produit tensoriel, puissances extérieures,
puissances\\ symé\-triques}{Produit tensoriel, puissances extérieures,
puissances symé\-triques}
 \label{ProdTens}

Soient $M$ et $N$ deux \Amosz.
Une application bi\lin $\varphi :M\times N\to P$ est
appelée
un \ix{produit tensoriel} des \Amos $M$ et $N$
si toute application bi\lin $\psi
:M\times N\to R$ s'écrit de manière unique sous la forme
$\psi=\theta\circ\varphi$, où $\theta$ est une \Ali de $P$ vers $R$.

\vspace{-1.1em}
\Pnv{M\times N}{\varphi}{\psi}{P}{\theta}{R}{}{applications bi\lins}{\alisz.}

\vspace{-1em}
Il est alors clair que  $\varphi :M\times N\to P$ est unique au sens
catégorique,
\cad que pour tout autre produit tensoriel  $\varphi' :M\times N\to P'$ il y a une
\ali unique $\theta:P\to P'$ qui rend le diagramme convenable commutatif,
et que $\theta$ est un \isoz.

Si $(\ug)$ est un \sgr de $M$ et $(\uh)$ un \sgr de~$N$, une application bi\lin
$\lambda : M\times N\rightarrow P$ est connue à partir de ses valeurs sur les
\elts de $\ug\times \uh$.
En outre, les valeurs $\lambda(x,y)$ sont liées par certaines contraintes, qui
proviennent des syzygies entre \elts de~$\ug$ dans $M$ et des  syzygies entre
\elts de $\uh$ dans $N$. \\
Par exemple, si l'on a une syzygie $a_1x_1+a_2x_2+a_3x_3=_M
0$ entre des \elts $x_i$ de $\ug$, avec les~$a_i$ dans~$\gA$, cela fournit
pour chaque $y\in \uh$ la syzygie suivante dans $P$:
$a_1\lambda(x_1,y)+a_2\lambda(x_2,y)+a_3\lambda(x_3,y)=0$.

En fait: \gui{ce sont les seules contraintes indispensables, et cela montre qu'un
produit tensoriel peut être construit}.

Plus \prmtz, notons $x\te y$ à la place de $(x,y)$
un \elt arbitraire de~$\ug\times \uh$. 
Considérons alors le \Amo $P$  engendré par
les $x\te y$, liés par les syzygies décrites ci-dessus
(%, c'est la  syzygie  
 $a_1(x_1\te y)+a_2(x_2\te
y)+a_3(x_3\te y)=_P 0$ pour l'exemple donné).

%--- Proposition{propPftens}---
\begin{proposition}
\label{propPftens} 
(Avec les notations ci-dessus)
%-----------------begin enum------------------
\begin{enumerate}
\item Il existe une unique application bi\lin $\varphi :M\times N\to P$
telle que pour tout $(x,y)\in\ug\times \uh$, on ait $\varphi(x,y)=x\te y$.
\item Cette  application bi\lin fait de $P$ un produit tensoriel des modules~$M$ et~$N$. 
En particulier, si $M$ et $N$ sont libres de bases $(\ug)$ et $(\uh)$, 
le module~$P$ est
libre de base $(\ug\te\uh):=({x\te y})_{x\in \ug,\,y\in\uh}$.

\end{enumerate}
%-----------------end enum------------------
\end{proposition}
%--- end-proposition-------------------
%-----------------begin proof------------------
\facile
%-----------------end proof------------------

Ainsi, le produit tensoriel de deux \Amos existe et peut toujours
être défini à partir de \pns de ces modules. Il est noté
$M\te_\gA N$.

Le fait qui suit est plus ou moins une paraphrase de la proposition précédente,
mais il ne peut être énoncé qu'une fois que l'on  sait que les produits tensoriels existent.

%:     Fact{factpropPftens}
\begin{fact}\label{factpropPftens}
\begin{enumerate}
\item  Si  deux modules
sont \tf (resp. \pfz) leur produit tensoriel l'est \egmtz.

\item Si $M$ est libre de base $(g_i)_{i\in I}$ et $N$ est libre
de base $(h_j)_{j\in J}$,
alors $M\otimes N$ est libre de base $(g_i\otimes h_j)_{(i,j)\in I\times J}$.

\item Si $M\simeq \Coker \alpha  $ et
$N\simeq \Coker \beta $, avec $\alpha  :L_1\rightarrow L_2$ et
$\beta :L_3\rightarrow L_4$, les modules $L_i$ étant libres, alors
l'\Ali

\snic{(\alpha \te\Id_{L_4})\oplus (\Id_{L_2}\te\beta ):
(L_1\otimes L_4)\oplus ( L_2\otimes L_3)\to L_2\otimes L_4}

%\sni
a pour conoyau un produit tensoriel de $M$ et $N$.
\end{enumerate}

\end{fact}

\comms~ \\
1) Il y a des raisons profondes, données dans la théorie qui a pour nom
\emph{\alg universelle}, qui font que la construction du produit tensoriel
\emph{ne peut pas ne pas marcher}. Mais cette théorie générale est un peu
trop lourde pour être exposée dans cet ouvrage, et il vaut mieux s'imbiber de
ce genre de choses par imprégnation sur des exemples.

\rdb%  
2) \Llec habitué\e aux \clama n'aura pas lu sans
appré\-hen\-sion notre \gui{présentation} du produit tensoriel de $M$ et $N$,
qui est un module
construit à partir de \pns de $M$ et $N$.\label{CommUnivPasUniv}
\Sil a lu Bourbaki, \il aura remarqué que notre construction est la même
que celle de l'illustre \matn multicéphale, à ceci près que Bourbaki se
limite à une \pn \gui{naturelle et universelle}: tout module est engendré par
\underline{tous} ses \elts liés par \underline{toutes} leurs syzygies. Si la
\gui{\pnz} de Bourbaki a le mérite de l'universalité, elle a l'inconvénient
de la lourdeur
de l'hippopotame.

En fait, en \comaz, on n'a pas la même \gui{théorie des ensembles} sous-jacente
qu'en \clamaz. Une fois que l'on  a \emph{donné} un module $M$ au moyen
d'une \pn $\alpha :L_1\to L_2$, on ne s'empresse pas d'oublier $\alpha$ comme on
fait semblant de le faire en \clamaz\footnote{Une inspection détaillée de l'objet $M$
construit selon la théorie des ensembles des \clama montrerait d'ailleurs que
ces dernières ne l'oublient pas non~plus.}.
Bien au contraire, du point de vue \cofz, le module $M$ n'est rien d'autre 
qu'\gui{un codage de
l'\ali $\alpha$}
(par exemple sous forme d'une matrice si la \pn est finie), avec l'information
complémentaire qu'il s'agit de la \pn d'un module.
D'autre part, un \gui{ensemble quotient} n'est pas vu comme un ensemble de classes
d'\eqvcz, mais comme \gui{le même préensemble muni d'une relation
d'\egt moins fine}: l'ensemble quotient de $(E,=_E)$ par la relation
d'\eqvc $\sim$ est simplement l'ensemble $(E,\sim)$.
En conséquence, notre
construction du produit tensoriel, conforme à son implémentation sur machine,
est entièrement \gui{naturelle et universelle} dans le cadre de la théorie
constructive des ensembles (\llec pourra consulter le simple et
génial chapitre 3
de \cite{B67}, ou l'un des autres ouvrages de référence classiques
pour les \coma \cite{Be,BB85,BR,MRR}).

3) Pour construire le produit tensoriel de deux modules \emph{non} discrets
$M$ et~$N$,
nous avons besoin a priori de la notion de module librement engendré par
un ensemble \emph{non} discret. Pour la \dfn \cov de ce type de modules libres, voir l'exercice \ref{propfreeplat}. On peut cependant contourner la
difficulté en ne faisant pas appel, dans la construction, à des \sgrs de $M$ et $N$. Les \elts du produit tensoriel $M\otimes_\gA N$ sont seulement des sommes formelles $\sum_{i=1}^{n}x_i\otimes y_i$ pour des familles finiment énumérées dans~$M$ et $N$. Le tout est de bien définir la relation d'\eqvc qui donne par passage au quotient le module $M\otimes_\gA N$.
Les détails sont laissés \alecz.
\eoe

%:     Fact{factSuitExTens}
D'après sa \dfn même, le produit tensoriel est \gui{fonctoriel}, i.e. si l'on a deux \Alis $f:M\to M'$ et $g:N\to N'$,
alors il existe une unique \ali $h:M\te_\gA N
\to M'\te_\gA N'$ vérifiant les \egts $h(x\te y)=f(x)\te g(y)$ pour $x\in M$ et $y\in N$.
Cette \ali est naturellement notée $h=f\te g$.

On a aussi des \isos canoniques 

\snic{M\te_\gA N\simarrow N\te_\gA M\;\hbox{  et  }\;
      M\te_\gA (N\te_\gA P)\simarrow (M\te_\gA N)\te_\gA P,}

%\sni
ce que l'on  exprime en
disant que le produit tensoriel est commutatif et associatif.

Le fait suivant résulte \imdt de la description du produit tensoriel
par \gtrs et relations.

%%%%%%%%%%%%%%%%%%%%%%%%%%%%%%%%%%%%%%%%%
\begin{fact}\label{factSuitExTens}
Pour toute suite exacte de \Amos $M\vers{f} N\vers{g} P\to 0$ et pour tout \Amo $Q$
la suite

\snic{M\te_\gA Q\vvvvers{f\te \Id_Q} N\te_\gA Q\vvvvers{g\te \Id_Q} P\te_\gA Q\to 0}

%\sni
est exacte.
\end{fact}
%%%%%%%%%%%%%%%%%%%%%%%%%%%%%%%%%%%%%%%%%

On exprime ce fait en disant que \gui{le foncteur $\bullet\otimes Q$
est exact à droite}.

 Nous ne rappellerons pas en détail l'énoncé des \pbs universels que
résolvent les puissances extérieures (déjà donné \paref{PuissExtMod}), les \emph{puissances \smqsz}\index{puissance symétrique!d'un module} et
l'\emph{\alg extérieure}\index{algèbre!extérieure d'un module}\index{extérieure!algèbre ---  d'un module} d'un \Amoz.
Voici néanmoins les \gui{petits diagrammes} correspondants pour les
deux derniers.

\vspace{-1.4em}
\Pnv{M^{k}}{\gs_\gA^k}{\psi}{\gS_\gA^k M}{\theta}{N}{}{applications
multi\lins \smqs~~~~~~~}{\alisz.}

\vspace{-2.6em}
\PNV{M}{ \lambda_\gA}{\psi}{\Vi_\gA\! M}{\theta}{\gB}%
{\Amos}{$\psi(x)\,\smalltimes\, \psi(x)=0$ pour tout $x\in M$~~~~~~~~~~~~~~}%
{\Algs associatives.}

\vspace{-1.4em}
Comme corolaire de la proposition \ref{propPftens} on obtient
la proposition qui suit.

%--- Proposition{propPfPex}---
\begin{proposition}
\label{propPfPex}
Si $M$ est un \Amo \pfz,  alors
il en va de même pour $\Al k_\gA\! M$
et pour les puissances \smqs  $\gS_\gA^k M$ ($k\in\NN$).\\
Plus \prmtz, si $M$ est engendré par le \sys $(x_1, \dots,x_n)$ soumis à des
syzygies~\hbox{$r_j\in\Ae n$}, on obtient les résultats suivants.
%-----------------begin enum------------------
\begin{enumerate}
\item Le module $\Al k_\gA \!M$ est engendré par les $k$-vecteurs 

\snic{x_{i_1}\vi\cdots\vi
x_{i_k}$ pour
$1\leq i_1<\cdots <i_k\leq n,}

%\sni
soumis aux syzygies obtenues en faisant le produit
extérieur des syzygies~$r_j$ par les $(k-1)$-vecteurs
$x_{i_1}\vi\cdots\vi x_{i_{k-1}}$.
\item Le module $\gS_\gA^k M$ est engendré par les tenseurs $k$-\smqs

\snic{\gs(x_{i_1},\dots, x_{i_k})$ pour
$1\leq i_1\leq \cdots \leq i_k\leq n,}

%\sni
soumis aux syzygies obtenues en faisant le
produit des syzygies $r_j$ par les tenseurs $(k-1)$-\smqs
$\gs(x_{i_1},\dots, x_{i_{k-1}})$.
\end{enumerate}
%-----------------end enum------------------
\end{proposition}
%--- end-proposition-------------------

Par exemple, avec $n=4$ et $k=2$ une syzygie $a_1x_1+\cdots +a_4x_4=0$
dans~$M$  donne lieu à 4 syzygies dans $\Al2_\gA M$:
$$\arraycolsep2pt
\begin{array}{rcl}
a_2 \,(x_1\vi x_2) \,+\,a_3 \,(x_1\vi x_3) \,+\,a_4 \,(x_1\vi x_4)&  = &  0    \\
a_1 \,(x_1\vi x_2)  \,-\,a_3 \,(x_2\vi x_3) \,-\,a_4 \,( x_2\vi x_4)&  = &  0 \\
a_1 \,(x_1\vi x_3) \,+\,a_2 \,(x_2\vi x_3) \,-\,a_4 \,(x_3\vi x_4)&  = &  0    \\
a_1 \,(x_1\vi x_4) \,+\,a_2 \,(x_2\vi x_4) \,+\,a_3 \,(x_3\vi x_4)&  = &  0
\end{array}
$$
et à 4 syzygies dans $\gS_\gA^2 M$:
$$\arraycolsep2pt
\begin{array}{rcl}
a_1\,\gs(x_1,x_1) \,+\, a_2\,\gs(x_1,x_2) \,+\,a_3\,\gs(x_1,x_3)
\,+\,a_4\,\gs(x_1,x_4)&  = &  0    \\
a_1\,\gs(x_1,x_2) \,+\, a_2\,\gs(x_2,x_2) \,+\,a_3\,\gs(x_2,x_3)
\,+\,a_4\,\gs(x_2,x_4)&  = &  0    \\
a_1\,\gs(x_1,x_3) \,+\, a_2\,\gs(x_2,x_3) \,+\,a_3\,\gs(x_3,x_3)
\,+\,a_4\,\gs(x_3,x_4)&  = &  0    \\
a_1\,\gs(x_1,x_4) \,+\, a_2\,\gs(x_2,x_4) \,+\,a_3\,\gs(x_3,x_4)
\,+\,a_4\,\gs(x_4,x_4)&  = &  0
\end{array}
$$

%:HHH rajout d'une remarque
\sni
\rem
De manière plus \gnlez, pour toute suite exacte:
$$\preskip.4em \postskip.1em 
K \vers {u} G \vers {p} M \to 0 
$$
on a une suite exacte:
$$\preskip.0em \postskip.4em\ndsp 
K \otimes \Al{k-1} G \vers {u'} \Al{k} G \vvvers {\Al{k} p} \Al{k} M \to 0 
$$
avec $u'(z \otimes y) = u(z) \vi y$ pour $z \in K$, $y \in \Al{k-1} G$. \\ 
\`A
droite, la surjectivité est immédiate et il est clair que $(\Al{k}
p) \circ u' = 0$, ce qui permet de définir $p' : \Coker u' \to \Al{k} M$ par
passage au quotient. Il reste à prouver que $p'$ est un \isoz. Pour cela, il
suffit de construire une \aliz~\hbox{$q' : \Al{k}
M \to \Coker u'$} qui soit l'inverse de~$p'$. On n'a pas le choix: pour
$x_1$, \ldots, $x_k \in M$ avec des antécédents $y_1$, \ldots, $y_k \in G$ par
$p$

\snic {
q'(x_1 \vi \cdots \vi x_k) = y_1 \vi \cdots \vi y_k \bmod \Im u'.
}

%\sni
On laisse le soin \alec de vérifier que $q'$ est bien définie
et convient.
\\
Le résultat analogue vaut pour les puissances symétriques.
\eoe

\sni\exl \label{belexemple}
Soit $\gB$ l'anneau des \pols $\gA[x,y]$ en les \idtrs $x$ et $y$
sur un anneau $\gA$ non trivial. On considère l'\id $\fb=\gen{x,y}$
de $\gB$, et on le regarde comme un \Bmo que l'on  note $M$.
Alors, $M$ admet le \sgr $(x,y)$ pour lequel une  \mpn est égale à
$\Cmatrix{3pt}{\phantom-y\cr-x}$. On en déduit que $M\te_\gB M$ admet 
$(x\te x,x\te y,y\te x,y\te y)$ pour \sgrz, avec une  \mpn  égale à:
$$\preskip.4em \postskip.4em
\crmatrix{y&0&0&y\cr-x&0&y&0\cr0&y&0&-x\cr0&-x&-x&0}
\leqno{ 
\begin{array}{c} 
x\te x\\ x\te y\\ y\te x\\ y\te y \\[-8.5pt] ~
\end{array}}
%\matrix{x\te x\cr x\te y\cr y\te x\cr y\te y}
$$
On en déduit les annulateurs suivants:
$$
\begin{array}{c}
\Ann_\gB(x\te y-y\te x)=\fb, \quad \Ann_\gB(x\te y+y\te x)=\Ann_\gA(2)\,\fb,   
\\[1mm]
\Ann_\gB(x\te x)=\Ann_\gB(x\te y)=\Ann_\gB(y\te x)=\Ann_\gB(y\te y)=0.
\end{array}
$$
Le dual $M\sta=\Lin_\gB(M,\gB)$ de $M$ est libre de rang 1, engendré par la forme
\label{NOTAdual}
$$
\alpha :M\longrightarrow \gB,\quad z \longmapsto z,
$$
ce qui donne seulement une information partielle sur la structure de $M$.
Par exemple, pour toute forme \lin $\beta  :M\to\gB$ on a $\beta (M)\subseteq\fb$
et donc  $M$ ne possède pas de facteur direct
libre de rang 1 (cf. proposition~\ref{propSplittingOffAlgExt}).\\
De même, le dual $(M\te_\gB M)\sta$ de $M\te_\gB M$ est libre de rang 1,
engendré par la forme
$$
\varphi  :M\te_\gB M\longrightarrow \gB,\quad z\te z' \longmapsto zz',
$$
et  $M\te_\gB M$  ne possède pas de facteur direct
libre de rang 1.\\
Concernant $\gS_\gB^2 M$, on trouve qu'il  admet
un \sgr égal à~$\big(\gs(x,x),\alb\gs(x,y),\gs(y,y)\big)$,
avec la \mpn
$$\preskip.4em \postskip.4em
\crmatrix{y&0\cr-x&y\cr0&-x}.
\leqno{
\begin{array}{c} 
\gs(x,x)\\ \gs(x,y)\\ \gs(y,y)\\[-8.5pt] ~
\end{array}} 
$$
Concernant $\Al2_\gB M$, on trouve qu'il est engendré par $x\vi y$
avec la \mpn $\lst{x\;y}$ ce qui donne
$$
\Al2_\gB M\,\simeq\,\gB/\fb\,\simeq\,\gA.
$$
Mais attention au fait que $\gA$ comme \Bmo est un quotient et non un
sous-module de $\gB$.
\eoe

\rdb
%--- Sec{Changement d'anneau}
\subsection*{Changement d'anneau de base}
\addcontentsline{toc}{subsection}{Changement d'anneau de base}
\label{pageChgtBase}

Soit  $\rho:\gA\rightarrow \gB$ une \algz.
Tout \Bmo $P$ peut être muni d'une structure de \Amo via $\rho$
en posant $a.x\eqdefi \rho(a)x$.

%--- Definition{defAliExtScal}-------
\begin{definition}
\label{defAliAliExtScal} Soit $\gA\vers{\rho}\gB$ une \Algz.
\begin{enumerate}
\item Soit $M$ un \Amoz. Une application \Ali $\varphi :M\to P$, où $P$
est un \Bmoz, est appelée
un \emph{morphisme d'extension des scalaires} (de $\gA$ à $\gB$ pour $M$),
ou encore un \ix{changement d'anneau de base} (de $\gA$ à $\gB$ pour $M$), 
si la \prt \uvle suivante est satisfaite.
\\
\emph{Pour tout \Bmo $R$,  
toute  \Ali $\psi :M\to R$ s'écrit de manière unique sous
la forme $\psi=\theta\circ\varphi$, où $\theta\in\Lin_\gB(P,R)$.}%

\vspace{-20pt}

\PNV{M}{\varphi}{\psi}{P}{\theta}{R}{\Amos}{\Alis}{\Bmosz, \Blis}

\vspace{-20pt}

\item Un \Bmo $P$ tel qu'il existe un \Amo $M$ et un morphisme
d'\eds  $\varphi :M\to P$ est dit
\emph{étendu} \emph{depuis~$\gA$}.
On dira aussi que $P$ \emph{provient du \Amo $M$ par \edsz}.
\end{enumerate}%
\index{morphisme!d'extension des scalaires}%
\index{extension des scalaires}%
\index{module!eten@étendu}%
\index{eten@étendu!module ---}
\end{definition}
%--- end-definition------------------------------------

Il est  clair qu'un morphisme d'\eds
$\varphi :M\to P$ est unique au sens catégorique,
\cad que pour tout autre morphisme d'\eds
$\varphi' :M\to P'$, il y a un
 unique    $\theta\in\Lin_\gB(P,P')$ qui rend le diagramme convenable
commutatif, et que $\theta$ est un \isoz.

Si $(\ug)$ est un \sgr de $M$ et $P$ un \Bmo arbitraire, une \Ali $\lambda :
M\rightarrow P$ est connue à partir de ses valeurs sur les \elts $x$ de $\ug$.
En outre, les valeurs $\lambda(x)$ sont liées par certaines contraintes, qui
proviennent des syzygies entre \elts de $\ug$ dans $M$. Par exemple, si l'on a une
syzygie $a_1x_1+a_2x_2+a_3x_3=_M 0$ entre des \eltsz~$x_i$ de~$\ug$, avec
les~$a_i$ dans $\gA$, cela fournit la syzygie suivante entre les $\lambda(x_i)$ dans~$P$:
$\rho(a_1)\lambda(x_1)+\rho(a_2)\lambda(x_2)+\rho(a_3)\lambda(x_3)=0$.

En fait \gui{ce sont les seules contraintes indispensables, et cela montre qu'une
\eds peut être construite}.

Plus \prmtz, notons $\rho\ist(x)$ à la place de $x$
(un \elt arbitraire de~$\ug$). Considérons alors le \Bmo $M_{1}$  engendré par les
$\rho\ist(x)$, liés par les syzygies décrites ci-dessus
(%, c'est la syzygie 
$\rho(a_1)\rho\ist(x_1)+\rho(a_2)\rho\ist(x_2)+\rho(a_3)\rho\ist(x_3)=_P 0$
pour l'exemple donné).

%--- Proposition{propPfExt}-----
\begin{proposition}
\label{propPfExt}
 (Avec les notations ci-dessus)
%-----------------begin enum------------------
\begin{enumerate}
\item \label{i1propPfExt}
%-----------------begin enum------------------
\begin{enumerate}
\item Il existe une unique \Ali $\varphi :M\to M_{1}$
telle que pour tout $x\in\ug$, on ait $\varphi(x)=\rho\ist(x)$.

\item Cette  \Ali fait de $M_{1}$ une \eds de $\gA$ à $\gB$
pour $M$. On  notera $M_{1}=\rho\ist(M)$.\perso{l'inconvénient est que $M_{1}$ n'est pas en \gnl égal à l'image de $M$ par $\varphi=\rho\ist$.}
\label{NOTArhosta}

\item Dans le cas d'un \mpfz, si $M$ est (isomorphe au) conoyau d'une matrice
$F=(f_{i,j})\in\gA^{q\times m}$, alors $M_{1}$
est (isomorphe au) conoyau de
{\em la même matrice vue dans $\gB$}, \cad la matrice  $F^\rho=\big(\rho(f_{i,j})\big)$.
En particulier, si $M$ est libre de base $(\ug)$, $M_{1}$ est libre de base
$\rho\ist(\ug)$.
\end{enumerate}
%-----------------end enum------------------
%
\item \label{i2propPfExt}
En conséquence l'\eds de $\gA$ à $\gB$ pour un \Amo
arbitraire existe et peut toujours être définie à partir d'une \pn de ce
module. Si le  module
est \tf (resp. \pfz) l'\eds l'est \egmtz.
\item \label{i3propPfExt}
 Sachant que les extensions de scalaires existent, on
peut décrire la construction précédente (de manière non circulaire) comme
suit: \\
 si~$M\simeq \Coker \alpha  $  avec $\alpha  :L_1\rightarrow L_2$, les
modules $L_i$ étant libres, alors le module
$M_{1}= \Coker \big(\rho\ist(\alpha)\big)$
est une \eds de $\gA$ à $\gB$ pour le module $M$.
\item \label{i4propPfExt}
 L'\eds est transitive: si $\gA\vers\rho\gB\vers{\rho'}\gC$
sont deux \algs \gui{successives} et si $\rho''=\rho'\circ \rho$ définit l'\alg
\gui{composée}, l'\Cli canonique $\rho''\ist(M)\to \rho'\ist\big(\rho\ist(M)\big)$
est un \isoz.
\item  \label{i5propPfExt}
L'\eds et le produit tensoriel commutent:
si~$M$,~$N$ sont des \Amos et $\rho:\gA\to\gB$ un \homo d'anneaux,
l'\Bli naturelle $\rho\ist(M\te_\gA N)\to\rho\ist(M)\te_\gB\rho\ist(N)$ est un
\isoz.
\item  \label{i6propPfExt}
De même l'\eds commute avec la construction
des puissances extérieures, des puissances \smqs et de l'\alg
extérieure.
\item  \label{i7propPfExt}
Vu comme \Amoz,  $\rho\ist(M)$ est isomorphe (de manière unique)
au produit tensoriel $\gB\te_\gA M$
 ($\gB$ est ici muni de sa structure de \Amo via $\rho$). En outre,
la \gui{loi externe} $\gB\times \rho\ist(M)\to\rho\ist(M),$ qui définit
la structure de \Bmo de $\rho\ist(M)$, s'interprète via l'\iso précédent
comme l'\Ali 
$$\preskip.4em \postskip.4em 
\pi\te_\gA \Id_M:\gB\te_\gA\gB\te_\gA M\lora  \gB\te_\gA M, 
$$
obtenue à partir de l'\Ali   $\pi:\gB\te_\gA\gB\to\gB$ \gui{produit dans $\gB$} ($\pi(b\te c)=bc$).
\item \label{i8propPfExt}
Pour toute suite exacte de \Amos $M\vers{f} N\vers{g} P\to 0$
la suite
$$\preskip.3em \postskip.0em 
\rho\ist(M)\vvvvers{\rho\ist(f)} \rho\ist(N) \vvvvers{\rho\ist(g)} \rho\ist(P)\to 0 
$$
est exacte.
\end{enumerate}
%-----------------end enum------------------
\end{proposition}
%--- end-proposition-------------------
%-----------------begin proof------------------
\facile
%-----------------end proof------------------

Ainsi, un \Bmo $P$ est étendu depuis
$\gA$ \ssi il est isomorphe à un module $\rho\ist(M)$.
On prendra garde cependant au fait qu'un
 \Bmo étendu peut provenir de plusieurs \Amos non isomorphes:
 par exemple lorsque
 l'on étend un $\ZZ$-module à $\QQ$, \gui{on tue la torsion},
 et $\ZZ$ et $\ZZ\oplus \aqo{\ZZ}{3}$ donnent tous deux par extension
 des scalaires un \Qev de dimension~$1$.

\smallskip
\rem Avec la notation tensorielle du point \emph{\ref{i7propPfExt}} l'\iso
canonique donné au point \emph{\ref{i5propPfExt}} s'écrit:

\snic{
\gC\otimes_\gA M \vvers{\varphi} \gC\otimes_\gB(\gB \otimes_\gA M)\simeq
(\gC\otimes_\gB\gB) \otimes_\gA M,}

avec $\varphi(c\te x)=c\te(1_\gB\te x)$. Nous reviendrons sur ce type
d'\gui{associativité} dans la remarque qui suit le corolaire \ref{corPlatEds} \paref{corPlatEds}.
\eoe

%%%%%%%%%%%%%%%%%%%%%%%%%%%%%%%%%%%%%%%%%%%%%%%%%%%%%%%%%%%%%%%%%%%%%%%%%%%
%:--- Subsec{Modules d'alin}
\subsec{Modules d'\alis}

%--- Proposition{propAliCoh}--------
\begin{proposition}
\label{propAliCoh}
  Si $M$ et $N$ sont des \mpfs sur un \cori $\gA$, alors $\Lin_\gA(M,N)$ est \pfz.
\end{proposition}
%--- end-proposition-----------------------------------------
%-----------------begin proof------------------
\begin{proof}
On reprend les notations de la section \ref{secCatMpf}.\\ 
Donner
un \elt $\varphi$ de $
\Lin_\gA(M,N)$ revient à donner les matrices
de $\rG_\varphi$ et~$\rK_\varphi$
qui satisfont la condition $\rG_\varphi\,\rA_M=\rA_N\,\rK_\varphi$.
\\
Puisque l'anneau est  \cohz,
les solutions de ce \sli forment \linebreak 
un \Amo \tfz, engendré par exemple par les solutions
correspondant à des \alis $\varphi_1$, \ldots, $\varphi_\ell$ données par
des couples de matrices $(\rG_{\varphi_1},\rK_{\varphi_1})$, \ldots, $(\rG_{\varphi_\ell},\rK_{\varphi_\ell})$.
Donc  $\Lin_\gA(M,N)=\gen{\varphi_1,\ldots ,\varphi_\ell}$.\\
Par ailleurs, une syzygie $\som_i a_i\varphi_i=0$ est vérifiée \ssi 
on %peut trouver 
a une
\ali
 $Z_\varphi:\rG_M\rightarrow \rK_N$ vérifiant
$\rA_N\, Z_\varphi=\som_ia_i\rG_{\varphi_i}$. En prenant le \sli correspondant,
dont les inconnues sont les $a_i$ d'une part et les \coes de la matrice de
$Z_\varphi$ d'autre part, on constate que le module des syzygies pour le \sgr
$(\varphi_1,\ldots ,\varphi_\ell)$ est bien \tfz.
\end{proof}
%-----------------end proof------------------

%--- Sec{Caractere local}
\subsection*{Le caractère local des \mpfsz}
\addcontentsline{toc}{subsection}{Le caractère local des \mpfsz}

Le fait qu'un \Amo est \pf est une notion locale, au sens suivant.

%--- Principe loc glo conc{plcc.pf}--
\begin{plcc}
\label{plcc.pf} \emph{(Modules \pfz)}\\
Soient $S_1$, $\ldots$, $S_n$  des \moco d'un anneau $\gA$, et $M$ 
\linebreak un \Amoz.
Alors,  $M$ est \pf \ssi chacun des~$M_{S_i}$ est un $\gA_{S_i}$-\mpfz.
\end{plcc}
%--- end-plcc-----------------------------------------
%-----------------begin proof------------------
\begin{proof}
Supposons que  $M_{S_i}$ soit un $\gA_{S_i}$-\mpf pour chaque~$i$.
Montrons que  $M$ est \pfz.
\\
D'après le \plgrf{plcc.tf}, $M$ est \tfz. Soit $(g_1,\ldots ,g_q)$ un \sgr de $M$.
\\
Soient $(a_{i,h,1}, $\ldots$, a_{i,h,q})\in \gA_{S_i}^q$ des syzygies entre les
$g_j/1\in M_{S_i}$ (autrement dit, $\sum_j \,a_{i,h,j}g_j=0$ dans $M_{S_i}$) pour
$h=1,$ \ldots$, k_i$, qui engendrent \linebreak 
le  $\gA_{S_i}$-module  des syzygies entre les $g_j/1$. 
\\
On suppose \spdg que les
$a_{i,h,j}$ sont de la forme $a'_{i,h,j}/1$, avec $a'_{i,h,j}\in \gA$. Il existe
alors un  $s_i\in S_i$ convenable tel que les vecteurs
$$\preskip.4em \postskip.4em 
s_i\,(a'_{i,h,1},\ldots,a'_{i,h,q}) = (b_{i,h,1},\ldots,b_{i,h,q}) 
$$
soient des $\gA$-syzygies entre les~$g_j\in M$.
\\
Montrons que les  syzygies ainsi construites entre les~$g_j$
engendrent toutes les syzygies. Considérons pour cela une syzygie arbitraire
$(\cq)$  entre les~$g_j$.
Regardons la comme une syzygie entre les $g_j/1\in M_{S_i}$, et écrivons
la  comme combinaison $\gA_{S_i}$-\lin des vecteurs
$(b_{i,h,1},\ldots,b_{i,h,q})$ dans $\gA_{S_i}^q$. Après multiplication par un
$s'_i\in S_i$ convenable on obtient une \egt dans~$\gA^q$:
%---------begin $$----------
$$\preskip.2em \postskip.4em
s'_i(\cq)=e_{i,1}(b_{i,1,1},\ldots,b_{i,1,q})+
\cdots+e_{i,k_i}(b_{i,k_i,1},\ldots,b_{i,k_i,q}).
$$
%---------end $$----------
On écrit $\sum_{i=1}^{n} u_i s'_i =1$. On voit que $(\cq)$ est
combinaison $\gA$-linéaire des  $(b_{i,h,1},\ldots,b_{i,h,q})$.
\end{proof}
%-----------------end proof------------------

%: subsec{Tenseurs nuls}
\subsec{Tenseurs nuls}

Soient $M$ et $N$ deux \Amos arbitraires,  et $t=\sum_{i\in \lrbn}x_i\te y_i\in M\te N$.\\
L'\egt $\sum_{i}x_i\te y_i=0$
ne dépend pas seulement a priori de la connaissance des sous-modules $\sum_i\gA x_i\subseteq M$ et $\sum_i\gA y_i\subseteq N$. 
\\
En conséquence la notation $\sum_{i}x_i\te y_i$
est en \gnl lourde d'ambigüité, et dangereuse. On devrait la préciser comme suit: $\sum_{i}x_i\otimes_{\gA,M,N} y_i$, ou au moins écrire  les \egts sous la forme 

\snic{\sum_{i}x_i\te y_i=_{M\te_\gA N} \dots}

Cette précaution ne devient inutile que dans le cas où les deux 
modules~$M$ et $N$ sont plats (voir le chapitre VIII), par exemple lorsque l'anneau~$\gA$ est un \cdiz.

%--- Lem du tenseur nul {lem-tenul}
\CMnewtheorem{lemtenul}{Lemme du tenseur nul}{\itshape}
%:     Lemma{lem-tenul}
\begin{lemtenul}\label{lem-tenul}\index{Lemme du tenseur nul}
Soit $M=\gA x_1+\cdots+\gA x_n$ un \mtfz, $N$ un autre module et $t=\sum_{i\in \lrbn}x_i\te y_i\in M\te_\gA N$. \\
Avec $X=[\,x_1\,\cdots\,x_n\,]\in M^{1\times n}$ et $Y=\tra[\,y_1\,\cdots\,y_n\,]\in N^{n\times 1}$,
on utilise la notation  $t=X\odot Y$.
\Propeq%
\begin{enumerate}
\item $t=_{M\te_\gA N}0$.
\item On a un $Z\in N^{m\times 1}$ et une matrice $G\in \Ae {n\times m}$ qui vérifient:
%---  equation eqdef.lem-tenul -------
\begin{equation}\label{eqdef.lem-tenul} \preskip.3em \postskip.4em
XG=_{M^ m} 0\quad {\rm et } \quad GZ=_{N^{n}}Y\,.
\end{equation}
%---------------------end equation--------------
%
\end{enumerate}
\end{lemtenul}
%--------- fin lem du tenseur nul ------------ 
%
\begin{proof}
\emph{2} $\Rightarrow$ \emph{1.} De manière \gnle l'\egt $X\odot G Z=XG\odot Z$ est assurée pour toute matrice $G$ à \coes dans $\gA$ parce que $x\otimes \alpha z=\alpha x\otimes z$ lorsque~$x\in M$, $z\in N$ et~$\alpha\in\gA$.

\emph{1} $\Rightarrow$ \emph{2.}\\
L'\egt $t=_{M\te N}0$ provient d'un nombre fini de syzygies
à l'intérieur des modules $M$ et $N$. Il existe donc un sous-module $N'$ tel que
$$\preskip.4em \postskip.4em 
\gA y_1+\cdots+\gA y_n\subseteq N'=\gA z_1+\cdots+\gA z_m\subseteq N, 
$$
et $X\odot Y=_{M\otimes N'} 0$. On note $Z=\tra[\,z_1\,\cdots\,z_m\,]$.
On a alors une suite exacte
$$\preskip.4em \postskip.4em 
K \vers{a}L\vers\pi N'\to 0 
$$
où $L$ est libre de base $(\ell_1,\dots,\ell_m)$ et $\pi(\ell_j)=z_j$, qui donne une suite exacte
$$\preskip.4em \postskip.4em 
M\te K \vvvers{\rI\te a}M\te L\vers{\rI\te \pi} M\te N'\to 0 
$$
Si $U\in M^{1\times m}$ vérifie $U\odot Z=_{M\te N'}0$, cela signifie que $U$ vu comme \elt de $M\te L\simeq M^{n}$, \cad vu comme $\sum_j u_j\te_{M\otimes L} \ell_j$,  est dans  le sous-module $\Ker(\rI\te \pi)=\Im(\rI\te a)$, autrement dit
$$\preskip.4em \postskip.4em 
\som_j u_j\te_{M\otimes L} \ell_j=\som_ix_i\otimes \som_{ij}a_{ij}\ell_j=\som_j \big(\som_i a_{ij}x_i\big)\te \ell_j
 $$
pour des $a_{ij}\in\gA$ qui vérifient $\sum_{j}a_{ij}z_j=0$.
Autrement dit $U=XA$ pour une matrice $A$ vérifiant $AZ=0$.
\\ 
Si l'on écrit $Y=HZ$ avec $H\in\Ae{n\times m}$, on a $XH\odot Z=0$, ce qui donne une \egt $XH=XA$ avec une matrice $A$ vérifiant $AZ=0$. \\
On pose alors
$G=H-A$ et l'on a $XG=0$ et $GZ=HZ=Y$. 
\end{proof}
% 

%--- Sec{classification}-----
\section[Problèmes de classification]{Problèmes de classification des \mpfsz}
%--------------------

Le premier \tho de classification concerne les \Amos libres de rangs finis: deux
\Amos $M\simeq \Ae m$ et $P\simeq \gA^p$ avec $m\neq p$ ne peuvent être
isomorphes que si $1=_\gA0$ (proposition~\ref{propDimMod1}).

\smallskip \rem
Notez que nous utilisons l'expression \gui{$M$ est un module libre de rang $k$} pour
signifier que $M$ est isomorphe à $\Ae k$, même dans le cas où nous
ignorons si l'anneau $\gA$ est trivial ou non. Cela n'implique donc pas
toujours a priori que l'entier $k$ est bien déterminé. 
\eoe

\smallskip 
Rares sont les anneaux pour lesquels on dispose d'une
classification complète \gui{satisfaisante} des \mpfsz. Le cas des corps
discrets est bien connu: tout \mpf est libre (cela résulte du pivot chinois ou
du lemme de la liberté).
Dans cet ouvrage nous traiterons quelques \gnns de ce cas \elrz:
 les anneaux de valuation, les anneaux principaux et
les anneaux \zeds réduits
(sections \ref{secBézout} et~\ref{secKrull0dim}), et certains anneaux de Prüfer
(proposition \ref{propAriCohZed} et \thref{thMpfPruCohDim}).

\smallskip
Concernant la classification des  \mtfsz, nous signalons les deux
résultats d'uni\-cité importants suivants.

%:--- Subsec{Type fini}---------
\subsec{Deux résultats concernant les \mtfsz}

%--- Theorem{prop unicyc}-----------
\begin{theorem}
\label{prop unicyc}
Soient $\fa_1\subseteq \cdots\subseteq\fa_n$ et $\fb_1\subseteq
\cdots\subseteq\fb_m$ des idéaux de~$\gA$ avec $n \le m$. 
Si un \Amo $M$ est isomorphe  
\hbox{à $\gA/\fa_1\oplus\cdots\oplus \gA/\fa_n$} et \hbox{à $\gA/\fb_1\oplus\cdots\oplus
\gA/\fb_m$}, alors:
%-----------------begin enum------------------
\begin{enumerate}
\item  on a $\fb_k=\gA$ pour $n < k \le m$,
\item  et $\fb_k=\fa_k$ pour $1 \le k \le n$.
\end{enumerate}
%-----------------end enum------------------
On dit que   $(\fa_1, \dots,\fa_n)$ est la liste des \emph{facteurs invariants\footnote{On notera que la liste donnée ici peut être raccourcie ou rallongée sur la fin par de termes $\fa_j=\gen{1}$ lorsque l'on n'a pas de test pour l'\egt en question. Un peu comme la liste des \coes d'un \pol qui peut être raccourcie ou rallongée par des $0$ lorsque l'anneau n'est pas discret.}} du module $M$.%
\index{facteurs invariants!d'un module}\index{invariant!facteur}
\end{theorem}
%--- end-theorem-----------------------------------------

%-----------------begin proof------------------
\begin{proof} 
\emph {1.} Il suffit de montrer que si $n < m$, alors $\fb_m=\gA$,
autrement dit que l'anneau $\gB := \gA/\fb_m$ est nul. En notant
$M=\gA/\fa_1\oplus\cdots\oplus \gA/\fa_n$, on a 
$$\ndsp
\preskip.4em \postskip.4em 
\gB^m=\bigoplus_{j=1}^m\gA/(\fb_j+\fb_m)\simeq M/\fb_mM\simeq
\bigoplus_{i=1}^n\gA/(\fa_i+\fb_m). 
$$
Or chaque $\gA/(\fa_i+\fb_m)$ est un
quotient de $\gB$, donc il existe une \ali surjective de $\gB^n$ sur $\gB^m$
et par suite $\gB$ est nul (proposition~\ref{propDimMod1}).
On suppose désormais \spdg que $m=n$.

\emph {2.} Il suffit de montrer que $\fb_k\subseteq\fa_k$ pour $k\in\lrbn$.  Remarquons que pour un \id $\fa$ et un \elt $x$ de $\gA$, le noyau de
l'\aliz~$y\mapsto yx \mod\fa$, de $\gA$ vers $x(\gA/\fa)$ est l'\id
$(\fa:x)$, et donc que
$$
\preskip.4em \postskip.4em 
\gA/\fa)\simeq \gA/(\fa:x). 
$$
Soit maintenant  $x \in \fb_k$. Pour $j\in\lrb{k..n}$, on a
$(\fb_j : x) = \gA$, et donc:
$$
\preskip.4em \postskip.4em \ndsp 
xM \simeq \bigoplus_{j=1}^n\gA/(\fb_j:x) =
\bigoplus_{j=1}^{k-1}\gA/(\fb_j:x),
\; \hbox { et } \;
xM \simeq \bigoplus_{i=1}^n\gA/(\fa_i:x). 
$$
En appliquant le point \emph {1} au module $xM$ 
avec les entiers $k-1$ et $n$, nous
obtenons $(\fa_k:x)=\gA$, i.e. $x\in\fa_k$.
\end{proof}
%-----------------end proof------------------

Notez que dans le \tho précédent, on n'a fait aucune hypothèse concernant
les idéaux (il n'est pas \ncr qu'ils soient \tf ou détachables  pour que le
résultat soit valide \cotz).
%--- Th {prop quot non iso}
\begin{theorem}
\label{prop quot non iso}
Soit $M$ un \Amo \tf et $\varphi\,:\,M\rightarrow M$ une \ali surjective.
Alors, $\varphi$ est un \iso et son inverse est un \pol en $\varphi$. Si un
quotient $M/N$ de $M$ est isomorphe à $M$, \linebreak alors~$N=0$.
\end{theorem}
%--- end-theorem-------------------
%--- Corollary{corInvDInv}-----
\begin{corollary}
\label{corInvDInv} Si  $M$ est un \mtfz, tout \elt $\varphi$ \iv à droite
dans $\End_\gA(M)$ est \ivz, et son inverse est un \pol en $\varphi$.
\end{corollary}
%--- end-corollary------------------------------------
%-----------------begin proof------------------
\begin{Proof}{Démonstration du \thrf{prop quot non iso}. }\\
Soit $(\xn)$ un \sgr de~$M$,
$\gB=\gA[\varphi]\subseteq\End_\gA(M)$, et $\fa=\gen{\varphi}$
l'idéal de $\gB$ engendré par $\varphi$. L'anneau~$\gB$ est
commutatif et l'on regarde $M$ comme un \Bmoz. Puisque
l'\ali $\varphi$ est surjective, il existe $P\in \Mn(\fa)$
avec
$P\,\tra{\lst{x_1\;\cdots\;x_n}}=\tra{\lst{x_1\;\cdots\;x_n}}$, c.-à-d.
$$\preskip.4em \postskip.4em 
(\In-P)\,\tra{\lst{x_1\;\cdots\;x_n}}= \tra{\lst{0\;\cdots\;0}}. 
$$
(où $\In=(\In)_\gB$ est la matrice identité de $\Mn(\gB)$),
et donc
$$\preskip.4em \postskip.4em 
\det(\In-P)\,\tra{\lst{x_1\;\cdots\;x_n}}=   
\wi{(\In-P)}\,(\In-P)\,\tra{\lst{x_1\;\cdots\;x_n}}      
= \tra{\lst{0\;\cdots\;0}}. 
$$
Donc $\det(\In-P)=0_\gB$,  or $\det(\In-P)=1_\gB-\varphi\, \psi$ avec
$\psi\in \gB$ (puisque~$P$ est à \coes dans $\fa=\varphi\, \gB$). Ainsi,
$\varphi\, \psi= \psi\,\varphi=\alb 1_\gB=\alb \Id_M$: $\varphi$
est \iv dans $\gB$.
\end{Proof}
%-----------------end proof------------------

%--- Sec{Anneaux Bézout}---------
\section{Anneaux \qisz}
\label{subsecAnneauxqi}
%-----------------------------------------

Dans la \dfn suivante, nous modifions de manière infinitésimale la notion
d'anneau intègre usuellement donnée en \comaz, non par plaisir, mais parce que  notre \dfn correspond mieux aux algorithmes mettant en \oe{}uvre les anneaux
intègres.

\begin{definition}
\label{defqi}
Un anneau est dit \emph{intègre} si tout \elt  est nul ou \ndzz\footnote{Un anneau intègre est aussi appelé un \emph{domaine d'intégrité} dans la  littérature classique. Nous préférons garder \gui{anneau intègre}
et n'utiliser le mot \gui{domaine} que dans des expressions composées comme \gui{domaine de Bézout}. Notons aussi qu'en \coma on établit une claire distinction entre les anneaux intègres et les anneaux \sdzz.
 Voir la \dfn \paref{eqSDZ}.}.
Un anneau $\gA$ est dit \ixe{\qiz}{quasi integre} lorsque tout \elt  admet pour
annulateur un (idéal engendré par un)
\idmz.\index{anneau!intègre}\index{anneau!quasi intègre}%
\index{integre@intègre!anneau ---}%
\index{integre@intègre!anneau quasi ---}%
\end{definition}
%--- end-definition---------------

Comme d'habitude, le \gui{ou} dans la \dfn précédente doit être lu comme
un ou explicite.
En conséquence, un anneau intègre est  un ensemble discret \ssi
en outre il est trivial ou non trivial.
Nos anneaux intègres non triviaux sont donc exactement  les \gui{discrete domains} de~\cite{MRR}.

Dans cet ouvrage, il arrive que l'on parle d'un \gui{\elt non nul} 
dans un anneau intègre, mais on devrait en fait dire 
\gui{\elt \ndzz} pour ne pas exclure le cas de l'anneau trivial.

%:     Fact{factQIRed}
\begin{fact}\label{factQIRed}
Un anneau \qi est réduit.
\end{fact}
%--------- fin fact ---------------------------------------------- 
%
\begin{proof}
Si $e$ est l'\idm annulateur de $x$ et si $x^2=0$, alors $x\in\gen{e}$,
\linebreak  donc $x=ex=0$.
\end{proof}

Un \cdi est un anneau intègre. Un anneau $\gA$ est intègre \ssi son anneau
total de fractions $\Frac\gA$ est un \cdiz.
Un produit fini d'anneaux \qis est \qiz.

Un anneau est intègre \ssi il est \qi et connexe.

Dans la littérature, un anneau \qi est parfois
appelé un \ixx{anneau}{de Baer}
ou encore, en anglais, un
\ix{pp-ring} (principal ideals are projective, cf. section~\ref{Idpp}).

%:--- Subsubsec*{Definition equationnelle}
\subsec{Définition équationnelle des anneaux \qisz}
%-------------------
Dans un anneau \qiz, pour $a\in\gA$,  notons $e_a$ l'unique \idm tel que
$\Ann(a)=\gen{1-e_a}$. On a %un \iso naturel 
$\gA\simeq \gA[1/e_a]\times \aqo{\gA}{e_a}$.\\
Dans l'anneau $\gA[1/e_a]$, l'\elt $a$ est \ndzz, et dans $\aqo{\gA}{e_a}$, $a$ est nul.
\\
On a alors $e_{ab}=e_a e_b$, $e_aa=a$ et $e_0=0$.
%\junk{Inversement supposons qu'un anneau commutatif soit muni d'une loi 
%unaire $a\mapsto
%\ci{a}$ qui vérifie les quatre axiomes:
%%------begin equation--eqaqis-----------
%\begin{equation}\label{eqaqis}
%(\ci{a})^2=\ci{a},\quad
%\ci{a}\,a=a,\quad
%\ci{(ab)}=\ci{a}\,\ci{b},\quad
%\ci{0}=0
%\end{equation}
%%---------------------end equation--------------
%Alors l'anneau est \qi et pour tout $a\in\gA$ on a $\Ann(a)=\gen{1-\ci{a}}$.
%En effet, $(1-\ci{a})a = 0$ et si $ax=0$ alors
%$\ci{a}\,x=\ci{a}\,\ci{x}\,x=\ci{(ax)}\,x=0$, donc $x=(1-\ci{a})x$.}
\\
Inversement, supposons qu'un anneau commutatif soit muni d'une loi \linebreak unaire $a\mapsto
\ci{a}$ qui vérifie les trois axiomes suivants:
%------begin equation--eqaqis-----------
\begin{equation}\label{eqaqis}\preskip.3em \postskip.3em
%(\ci{a})^2=\ci{a},\quad
\ci{a}\,a=a,\quad
\ci{(ab)}=\ci{a}\,\ci{b},\quad
\ci{0}=0
\end{equation}
%---------------------end equation--------------
Alors,   pour tout $a\in\gA$, on a $\Ann(a)=\gen{1-\ci{a}}$, et $\ci{a}$
est \idmz, de sorte que l'anneau  est \qiz.
En effet, tout d'abord $(1-\ci{a})a = 0$, et si $ax=0$, alors
$$\preskip.0em \postskip.4em 
\ci{a}\,x=\ci{a}\,\ci{x}\,x=\ci{(ax)}\,x=\ci{0}\,x=0, 
$$
 donc $x=(1-\ci{a})x$:
ainsi $\Ann(a)=\gen{1-\ci{a}}$.
Voyons ensuite que $\ci{a}$ est \idmz. Appliquons 
le résultat précédent à $x = 1-\ci{a}$ qui vérifie $ax = 0$
(d'après le premier axiome): l'\egt $x=(1-\ci{a})x$ %\linebreak
  donne $x=x^2 $,
i.e. l'\elt $1 - \ci{a}$ est \idmz.

%%--- Subsubsec*{Lemme de scindage \qi}
%\subsubsection*{Lemme de scindage \qiz}

Le lemme de scindage suivant est à peu près immédiat.

%:     Lemma{thScindageQi}
\CMnewtheorem{lemScindageQi}{Lemme  de scindage \qiz}{\itshape}

\begin{lemScindageQi}\label{thScindageQi}
Soit $n$ \elts $x_1$, \dots, $x_n$ dans un anneau \qiz~$\gA$.
Il existe un \sfio $(e_j)$ de cardinal~$2^n$ tel que dans chacune des
composantes~\hbox{$\gA[1/e_j]$}, chaque~$x_i$ est nul ou \ndzz.
\end{lemScindageQi}
\begin{proof}
Soit $r_i$ l'\idm tel que $\gen{r_i}=\Ann(x_i)$, et $s_i=1-r_i$.
En développant le produit $1=\prod_{i=1}^n(r_i+s_i)$ on obtient le
\sfio indexé par $\cP_n$: $e_J = \prod_{j\in J}r_j \prod_{k\notin J}s_k$.
On peut  supprimer certains \elts
de ce \sys quand on sait qu'ils sont nuls.
\end{proof}

%
%:--- Subsubsec*{Machinerie locale-glo elr 1}
\subsect{Des anneaux intègres aux anneaux \qisz}{Machinerie \lgbe \elr \num1: des anneaux intègres aux anneaux \qisz}

Le fait de pouvoir scinder systématiquement en deux composantes un anneau \qi
conduit à la méthode \gnle suivante.
La différence essentielle avec le lemme de scindage précédent
est que l'on  ne connaît pas a priori la famille finie d'\elts qui
va provoquer le scindage.

%:   Machinerie locale-globale \elr \num1
\rdb
\mni {\bf Machinerie locale-globale \elr \num1.}\label{MethodeQI}\imlg
{\it La plupart des \algos qui fonctionnent avec les anneaux intègres non
triviaux peuvent être modifiés de manière à fonctionner avec les anneaux
\qisz, en scindant l'anneau
en deux composantes chaque fois que l'\algo écrit pour les anneaux intègres
utilise le test
\gui{cet \elt est-il nul ou \ndzz?}. Dans la première composante l'\elt en question
est nul, dans la seconde il est \ndzz.}

\medskip Un premier exemple d'application de cette machinerie \lgbe
sera donné
\paref{exlBezQi}.
Mais déjà le corolaire
\ref{corlemQI} ci-dessous pourrait être obtenu à partir du cas intègre,
où il est évident, en appliquant cette machinerie \lgbez.
\\
Expliquons pourquoi nous parlons ici de machinerie locale-globale
\elrz. De manière générale un \plg dit qu'une \prtz~$\sfP$ est vraie
\ssi elle est vraie \gui{après \lon en des \mocoz}. Dans le cas
présent, les \moco sont engendrés par des \elts $1-r_i$ où les $r_i$ forment
un \sfioz. En conséquence l'anneau est simplement isomorphe au produit des
localisés, et la situation est donc tout à fait simple, \elrz.

\smallskip
\rem \label{remMLGE1}
\Llec aura remarqué la formulation très informelle
que nous avons donnée pour cette machinerie locale-globale: \gui{La plupart
des \algos \ldots}. C'est qu'il nous a paru bien difficile de donner
par avance des conditions très précises requises pour que la méthode indiquée fonctionne. On pourrait imaginer un \algo qui fonctionne
pour tout anneau intègre, mais de façon pas du tout uniforme, ce
qui ferait que l'arbre correspondant que l'on  construit dans le cas \qi
ne serait pas fini. Par exemple, dans le cas intègre,
 une situation de départ donnée exigerait trois tests (pour terminer le calcul) si les réponses sont $0,0,0$, mais quatre tests si les réponses
sont $0,0,1,0$, puis cinq tests si ce sont les réponses $0,0,1,1,0$,
puis six tests  si ce sont les réponses $0,0,1,1,1,1$, puis sept tests  si ce sont les réponses $0,0,1,1,1,0,1$,  etc.  Naturellement,
on peut mettre en doute qu'un tel \algo puisse exister sans qu'existe
en même temps un anneau intègre qui le mette en défaut.
Autrement dit, un \algo qui n'est pas suffisamment uniforme n'est sans doute
pas un \algoz. Mais nous ne préjugeons de rien.\\
Même si nous n'avons pour le moment rencontré aucun exemple du
type ci-dessus où la machinerie \lgbe \elr ne s'appliquerait pas, nous ne pouvons exclure a priori une telle possibilité.  \eoe

%:--- Subsubsec*{Annulateur des itfs dans \qi}
\subsect{Annulateurs des \itfs dans les anneaux\\ \qisz}
{Annulateurs des \itfsz}
%-------------------

Le lemme suivant peut être considéré comme une variante
économique du lemme de scindage \ref{thScindageQi}.

%--- lem {lemQI} ------------
\begin{lemma}\label{lemQI}
Soient $x_1$, \dots, $x_n$ des \elts d'un \Amoz.
\\
 Si l'on a $\Ann(x_i) = \gen{r_i}$ où $r_i$ est un \idm  ($ i
\in\lrbn$), posons 

\snic{s_i=1-r_i$,
$t_1=s_1$, $t_2=r_1s_2$, $t_3=r_1r_2s_3$, $\dots$,
$t_{n+1}=r_1r_2\cdots r_n.}

%\sni
Alors, $(t_1, \dots, t_{n+1})$ est un \sfio
et l'\elt $x=x_1+t_2x_2+\cdots +t_nx_n$ vérifie

\snic{\Ann(x_1,\dots,x_n) = \Ann(x) = \gen{t_{n+1}}.}

%\sni
NB: dans la composante $t_k = 1$ ($k\in\lrbn$), on a 
   $x_k$ \ndz et $x_j = 0$ pour $j < k$, et dans la composante $t_{n+1}=1$,
   on a $x_1=\cdots=x_n=0$.
\end{lemma}

%--- Corollary{corlemQI}------
\begin{corollary}
\label{corlemQI}
Sur un anneau \qi $\gA$ tout sous-\mtf $M$ d'un module libre a pour
annulateur un idéal $\gen{r}$ avec $r$ \idmz,
et $M$ contient un \elt
$x$ ayant le même annulateur.
Ceci s'applique en particulier à un \itf de $\gA$.
\end{corollary}
%--- end-corollary------------------------------------

\begin{Proof}{\Demo du lemme \ref{lemQI}. }
On a  $t_1 x_1 = x_1$  et

\snic{\arraycolsep2pt
\begin{array}{rl}
1~ = & s_1 + r_1 = s_1 + r_1(s_2 + r_2) ~=~  s_1 + r_1s_2 +
r_1r_2(s_3+r_3)
~=~ \cdots
\\
= & s_1 + r_1s_2 + r_1r_2s_3 +\cdots +r_1r_2\cdots r_{n-1}s_n +
   r_1r_2\cdots r_n
\end{array}}

%\sni
donc $t_1,\dots,t_{n+1}$ est un \sfio et
$x = t_1 x_1 + t_2 x_2 + \cdots + t_nx_n$.
Il est clair que

\snic{\gen{t_{n+1}} \subseteq \Ann(x_1,\dots,x_n) \subseteq \Ann(x).}

%\sni
Inversement,
soit $z\in \Ann(x)$. Alors $zx=0$, donc $z t_i x_i = z t_i x=0$\linebreak 
 pour $i \in\lrbn$. Ainsi,
$z t_i \in \Ann(x_i) = \gen{r_i}$ et $z t_i = z t_i r_i = 0$.
Enfin, puisque $z = \sum_{i=1}^{n+1} z t_i$, on a $z = z t_{n+1} \in
\gen{t_{n+1}}$.
\end{Proof}

%%%%%%%%%%%%%%%%%%%%%%%%%%%%%%%%%%%%%%%%%%%%%%%%%%%%%%%%%%%%%%%%%%%%%%%%%%%

%:--- Subsub{Principe \lgb concret pour les anneaux \qi}
\subsect{Principe \lgb concret pour les anneaux\\ \qisz}{Principe \lgbz}

La \prt pour un anneau d'être \qi est locale au sens suivant.

%--- Principe local global concret{plcc.aqi}----
\begin{plcc}
\label{plcc.aqi} 
\emph{(Anneaux \qisz)}\\
Soient $S_1$, $\dots$, $S_n$ des \moco de $\gA$.
\Propeq
\begin{enumerate}
\item  L'anneau $\gA$ est \qiz.
\item Pour $ i=1$, $\dots$, $n,$
chaque anneau $\gA_{S_i}$ est \qiz.
\end{enumerate}
\end{plcc}
%--- end-plcc-----------------------------------------
%
\begin{proof}
Soit $a\in\gA$. Pour tout \mo $S$ de $\gA$ 
on a $\Ann_{\gA_S}(a)=\big(\Ann_{\gA}(a)\big)_S$.
Donc l'annulateur $\fa$ de $a$ est \tf \ssi il l'est après \lon en les $S_i$
(\plgref{plcc.tf}).
Ensuite l'inclusion $\fa\subseteq\fa^2$ relève du \plgc de base (\rref{plcc.basic}).\iplg
\end{proof}
%

%--- Sec{Anneaux Bézout}---------
\section{Anneaux de Bézout}
\label{secBézout}
%-----------------------------------------

Un anneau $\gA$ est appelé un \ixx{anneau}{de Bézout} lorsque tout \itf est
principal.\index{Bézout!anneau de ---}
Il revient au même de dire que tout \id avec deux \gtrs est principal:
%------begin equation--eqBézout-----------
\begin{equation}\label{eqBézout}
\forall a,\,b\;\;\exists u,\,v,\,g,\,a_1,\,b_1\;\;
\left(au+bv=g,\,a=ga_1,\,b=gb_1\right).
\end{equation}
%---------------------end equation--------------

Un anneau de Bézout est \fdi \ssi la relation de divisibilité y est explicite.

 Un \ixx{anneau}{local} est un anneau $\gA$ où est vérifié
l'axiome suivant:\index{local!anneau ---}
%------begin equation--eqAloc-----------
\begin{equation}\label{eqAloc}
\forall x,\,y\in \gA \qquad x+y \in\Ati\;\Longrightarrow \; ( x\in\Ati\;\mathrm{ou}\;
y\in\Ati )\,.
\end{equation}
%---------------------end equation--------------
Il revient au même de demander:
$$
\forall x\in \gA \qquad x \in\Ati\;\; {\rm  ou} \;\;  1-x \in\Ati\,.
$$

Notez que selon cette \dfn l'anneau trivial est local.
Par ailleurs, les \gui{ou} doivent être compris dans leur sens \cofz:
l'alternative doit être explicite.
La plupart des anneaux locaux avec lesquels on travaille usuellement en \clama
vérifient en fait la \dfn précédente si on les regarde d'un point de vue
\cofz.

Tout quotient d'un anneau local est local.
Un \cdi est un \aloz.

%:     Lemma{lemBezloc}----------
\begin{lemma}
\label{lemBezloc} \emph{(Bézout toujours trivial pour un \aloz)}\\
Un anneau est un anneau de Bézout local \ssi il vérifie la \prt
suivante: $\forall a,b\in\gA$, $a$ divise $b$ ou $b$ divise~$a$.
\end{lemma}
%--- end-lemma-----------------------------------------
%-----------------begin proof------------------
%:H2018   ajout preuve que la condition est suffisante
\begin{proof}
\emph{La condition est suffisante}. Tout d'abord il est clair que l'anneau est de Bézout. Supposons $x+y$  inversible. Si $x$ divise $y$, il divise aussi~\hbox{$x+y$} qui divise~$1$, donc~$x$ est \ivz. Symétriquement, si $y$ divise $x$, $y$ est \ivz.
Ainsi l'anneau est local. 
\\
\emph{La condition est \ncrz}.
On~a~\hbox{$g(1-ua_1-vb_1)=0$}. Puisque $1=ua_1+vb_1+(1-ua_1-vb_1)$, l'un des trois termes
dans la somme est \ivz. Si $1-ua_1-vb_1$ est \ivz,  
alors $g=a=b=0$. Si $ua_1$ est \ivz, alors
$a_1$ \egmtz, et $a$ divise~$g$  qui divise $b$. Symétriquement si~$vb_1$ est \ivz, alors $b$ divise~$a$.
\end{proof}
%-----------------end proof------------------

\rdb
Les anneaux de Bézout locaux sont donc les \gui{anneaux de valuation} au sens de
Kaplansky. Nous préférerons  la définition aujourd'hui usuelle:
un \ixx{anneau}{de valuation} est  un anneau de Bézout local réduit.%
\index{valuation!anneau de ---}

%:--- Subsec{Anneaux de valuation}
\subsect{Modules \pf sur les anneaux \\de valuation}{Modules \pf sur les anneaux de valuation} \label{secpfval}
%-------------------

Une matrice $B=(b_{i,j})\in\Ae {m\times n}$ est dite \ixd{en forme de
Smith}{matrice} si  tout \coe hors de la diagonale principale est nul, et si pour
$1\leq i< \inf(m,n)$, le \coe diagonal $b_{i,i}$ divise le suivant $b_{i+1,i+1}$.%
\index{Smith!matrice en forme de ---}

%--- Proposition{propPfVal} -
\begin{proposition}
\label{propPfVal}
Soit  $\gA$ un anneau de Bézout local.
%-----------------begin enum------------------
\begin{enumerate}
\item  Toute matrice de $\Ae{m\times n}$ est \elrt \eqve à une matrice
en forme de Smith.
\item  Tout \Amo \pf $M$ est isomorphe à une somme directe de modules
$\aqo{\gA}{a_i}$:  $M\simeq\bigoplus_{i=1}^p \aqo{\gA}{a_i}$, avec de plus, pour
chaque $i< p$, $a_{i+1}$ divise $a_{i}$.
\end{enumerate}
%-----------------end enum------------------
\end{proposition}
%--- end-proposition----------------------------------------
%-----------------begin proof------------------
\begin{proof}
\emph{1.} On utilise la méthode du pivot de Gauss en choisissant comme
premier pivot un \coe de la matrice qui divise tous les autres. On termine par
\recuz. \\
\emph{2.} Conséquence directe du point~\emph{1.}
\end{proof}
%-----------------end proof------------------

\rem
Ce résultat se complète par le \tho d'uni\-cité (\thrf{prop unicyc}) comme suit.
\begin{itemize}\itemsep.0pt
\item [\emph{1.}] Dans la réduite en forme de Smith les \ids $\gen{b_{i,i}}$ sont
déterminés de manière unique.
\item [\emph{2.}] Dans la décomposition $\bigoplus_{i=1}^p \aqo{\gA}{a_i}$, les \ids
$\gen{a_{i}}$ sont déterminés de manière unique, à ceci près  que des
\ids  en surnombre peuvent être égaux à $\gen{1}$: on peut supprimer les
termes correspondants, mais ceci ne se fait à coup sûr que lorsque l'on  a un
test d'inversibilité dans l'anneau.
\eoe
\end{itemize}

\rdb
 Un anneau $\gA$ est appelé un \ixx{anneau}{de Bézout strict} lorsque tout
\linebreak  vecteur $\vab u v\in\gA^2$
peut être transformé en un vecteur $\vab h 0$ par multiplication par une
matrice $2\times 2$ \ivz.%
\index{Bézout!anneau de --- strict}

Voici maintenant un exemple d'utilisation de la machinerie locale-globale \elr
\num1 (décrite \paref{MethodeQI}).\imlg

\ms\exl \label{exlBezQi}
On va montrer que \emph{tout anneau de Bézout \qi est un anneau de Bézout strict}.
\\
Commençons par le cas intègre. Soient $u$, $v\in\gA$:

\snic{\exists\,h,a,b,u_1,v_1\quad (h=au+bv,\,u=hu_1,\,v=hv_1).}

%\sni
Si $\Ann(v)=1$, alors $v=0$ et $\vab u 0 =\vab u v\;\I_2$. \\
Si $\Ann(v)=0$,
alors  $\Ann(h)=0$, $h(au_1+bv_1)=h$, puis $au_1+bv_1=1$.
Enfin, $\vab h 0=\vab u v\crmatrix{a & -v_1 \cr b  &u_1}$
et la matrice est de \deterz~1.
\\
Appliquons maintenant la machinerie locale-globale \elr \num1 expliquée
\paref{MethodeQI}. On considère l'\idm $e$ tel que 

\snic{\Ann(v)=\gen{e}$
et $f=1-e.}

%\sni
Dans $\gA[1/e]$, on a  $\vab u 0=\vab u v\;\I_2$. \\
Dans $\gA[1/f]$, on a  $\vab h 0=\vab u v\Cmatrix{.3em}{a & -v_1 \cr b  &u_1}$.

Donc dans $\gA$, on a  \smashbot{$\vab{ue+hf} 0 =\vab u v\Cmatrix{.3em}{fa+e & -fv_1 \cr fb  &fu_1+e}$},
et  la matrice est de \deterz~1.
\eoe

%:--- SUBsec{Anneaux principaux}
\subsect{Modules \pf sur les anneaux \\ principaux}{Modules \pf sur les anneaux principaux}

Supposons que  $\gA$ est  un anneau de Bézout strict. Si $a$ et $b$ sont deux
\elts sur une même ligne (resp. colonne) dans une matrice $M$ à \coes 
\linebreak 
dans $\gA$, on peut  postmultiplier (resp. prémultiplier) $M$
par une matrice \ivz, ce qui modifiera les colonnes (resp. les lignes)
où se trouvent les \coesz~$a$ et~$b$, lesquels seront remplacés par~$c$ et~$0$.
Pour parler de cette transformation de matrices, nous parlerons de
\emph{manipulations de Bézout}.
Les manipulations \elrs peuvent être vues comme des cas particuliers de
manipulations de Bézout.%
\index{manipulation!de Bézout}

Un anneau intègre est dit \ixc{principal}{anneau ---}\index{anneau!principal}
lorsqu'il est de Bézout et lorsque toute suite croissante d'idéaux principaux
admet deux termes consécutifs égaux (cf. \cite{MRR}). 
Autrement dit un anneau principal est un domaine de Bézout \noe 
(cf. \dfnz~\ref{definoetherien})
C'est par exemple le cas de $\ZZ$
ou de l'anneau de \pols $\gK[X]$ lorsque $\gK$ est un \cdiz.

%:     Proposition{propPfPID}-- forme de Smith
\begin{proposition}
\label{propPfPID} \emph{(\Tho de la base adaptée)}
Soit  $\gA$ un anneau principal.
%-----------------begin enum------------------
%
\begin{enumerate}
\item  Toute matrice $A\in\Ae {m\times n}$ est \eqve à une matrice
en forme de Smith. En notant $b_i$ les \coes diagonaux de la réduite, les \idps $\gen{b_1}\supseteq\cdots\supseteq\gen{b_q}$  ($q=\inf(m,n)$) sont des invariants de la {matrice} $A$ à \eqvc près. 
Une base~\hbox{$(e_1,\dots,e_m)$} de~$\Ae{m}$ telle \hbox{que $\Im(A)=\sum_{i=1}^{m}\gen{b_i}\, e_i$}
 est appelée une \emph{base adaptée au sous-module~\hbox{$\Im(A)$}.}%
\index{base adaptée!à une inclusion} 

\item  Pour tout \Amo \pf $M$, il existe  $r$, $p\in\NN$ et des \elts \ndzs $a_{1},$ \ldots, $a_{p}$, avec $a_{i}$ divise $a_{i+1}$  pour
$i< p$, tels que~$M$ est isomorphe à la somme directe 
$ \bigl(\bigoplus_{i=1}^p \aqo{\gA}{a_i}\bigr)\oplus\Ae r$.
\end{enumerate}
Si en outre $\gA$ est \fdi non trivial, on peut
demander dans le point 2. qu'aucun $\gen{a_i}$ ne soit égal à $\gen{1}$.
Dans ce cas, on appelle \emph{facteurs invariants du module $M$} les \elts
de la liste 
$\big(a_1,\ldots,a_p,\MA{\underline{0,\ldots,0}}\limits_{r\;{\it fois}}\big)$. 
Et la liste des facteurs invariants de $M$ est bien définie\footnote{On retrouve la \dfn donnée dans le \thref{prop unicyc}. On notera cependant que l'ordre est renversé et que l'on a remplacé ici les \idps par leurs \gtrsz, tout ceci pour se conformer à la terminologie la plus fréquente.} \gui{à association près}.%
\index{facteurs invariants!d'un module}\index{invariant!facteur}
\end{proposition}
%--- end-proposition----------------------------------------
%-----------------begin proof------------------
\begin{Proof}{Idée de la démonstration. }
Par des manipulations de Bézout sur les colonnes,
on remplace la première ligne
par un vecteur $(g_1,0,\ldots,0)$.
Par des manipulations de Bézout sur les lignes,
on remplace la première colonne
par un vecteur $(g_2,0,\ldots,0)$.
On continue le processus jusqu'à ce que pour un indice  $k$, on ait
$g_{k}\gA=g_{k+1}\gA$.
Par exemple, avec $k$ impair cela signifie que les dernières opérations de
lignes au moyen de manipulations de Bézout ont été faites à tort,
puisque $g_{k}$ divisait la première colonne.
On revient une étape en arrière, et l'on utilise $g_{k}$ comme pivot de Gauss.
On obtient ainsi une matrice de la forme
$$
\blocs{.5}{1.9}{.5}{1.9}{$g$}{$\begin{array}{ccc}0 & \cdots & 0\end{array}$}{$0$ \\[1mm] $~\vdots~$ \\[1mm]$0$}{$B$}
$$
Par \recu on obtient une réduite \gui{diagonale}. On vérifie enfin que l'on 
peut passer, par manipulations de Bézout et manipulations \elrsz, d'une matrice
$\left[\begin{array}{cc}a & 0 \\0 & b\end{array}\right]$
à une matrice $\left[\begin{array}{cc}c & 0 \\0 & d\end{array}\right]$
avec $c$ divise $d$.
\\
Le point \emph{2.} est une conséquence directe du point~\emph{1}.
\end{Proof}
%-----------------end proof------------------

\rems~\\
 1) Un \algo  plus simple peut être écrit si $\gA$ est \fdiz.
 
2)
On ne sait toujours pas (en 2018) si la conclusion de la proposition précédente est vraie sous la seule hypothèse que~$\gA$ est un anneau de Bézout intègre. On ne dispose ni de \demz, ni de contre exemple.\\
%:2018 petit ajout et la date 2008 est remplacee par 2011
On sait par contre le résultat vrai pour les anneaux de Bézout intègres
de dimension $\leq 1$: voir la remarque qui suit le \thref{thMpfPruCohDim}.
\eoe

%--- Sec{Anneaux zero-dim}---------
\section{Anneaux zéro-dimensionnels}
\label{secKrull0dim}
%-----------------------------------------

\vspace{4pt}
On dira qu'un anneau est \emph{\zedz} lorsqu'il vérifie l'axiome suivant:%
\index{zero-dimensionnel@\zedz!anneau ---}\index{anneau!zéro-dimensionnel}
%---  equation eqZed --------
\begin{equation}\label{eqZed}
\forall x\in \gA~\exists a\in\gA~\exists k\in
\NN\quad \quad x^{k}=ax^{k+1}
\end{equation}
%---------------------end equation--------------

Un anneau est dit \ixc{artinien}{anneau} s'il est \zedz, \coh et \noez.%
\index{anneau!artinien}

%:--- SUBsec{Propriétés de base}-
\subsec{Propriétés de base}
\label{subsecBasicsZerdim}

%--- Fact{factExZed}----------------
\begin{fact}
\label{factExZed} ~
\begin{enumerate}
\item [--] Tout anneau fini, tout \cdi est \zedz.
\item [--] Tout quotient et tout localisé d'un anneau \zed est \zedz.
\item [--] Tout produit fini d'anneaux \zeds est un anneau \zedz.
\item [--] Une \alg de Boole (cf. section \ref{secBoole}) est un anneau \zedz.
\end{enumerate}
\end{fact}
%--- end-fact-----------------------------------------

%:     lemme:idempotentDimension0
\begin{lemma}
\label{lemme:idempotentDimension0}
\Propeq
%-----------------begin enum------------------
\begin{enumerate}
\item  \label{LID001} $\gA$ est \zedz.
\item  \label{LID002} $\forall x\in \gA~\exists s\in\gA~\exists  d \in \NN^*$ tels
que
$\gen{x^d} = \gen{s}$ et $s$ \idmz.
\item  \label{LID003}
Pour tout \itf $\fa$ de $\gA$, il existe $d \in \NN^*$ tel que
$\fa^d = \gen{s}$ où $s$ est un \idmz, et en particulier,
$\Ann(\fa^d)=\gen{1-s}$ et $\fa^e=\fa^d$ pour $e \geq d$.
\end{enumerate}
%-----------------end enum------------------
\end{lemma}
\begin{proof}
\emph{\ref{LID001}} $\Rightarrow$ \emph{\ref{LID002}.} 
%D'après l'\eqrf{eqZed}, p
Pour tout  $x \in \gA$, il existe
$a\in \gA$  
et $k \in \NN$ tels que
$x^k = ax^{k+1}$.\\
 Si $k=0$ on a $\gen{x}=\gen{1}$, on prend $s=1$  
et $d=1$.
\\
Si $k\geq1$, on prend $d=k$: en multipliant $k$ fois par $ax$, on
obtient les \egts $x^k = a x^{k+1} = a^2 x^{k+2} = \cdots = a^k x^{2k}$.
Donc l'\elt $s=a^k x^k$ est un \idmz,~$x^k = sx^k$, et~$\gen{x^k} = \gen{s}$.

 \emph{\ref{LID002}} $\Rightarrow$ \emph{\ref{LID001}.}
On a $s=bx^d$ et $x^d s=x^d$. Donc, en posant $a=bx^{d-1}$,
on obtient les \egts $x^d = bx^{2d} = ax^{d+1}$.

\emph{\ref{LID002}} $\Rightarrow$ \emph{\ref{LID003}.}
Si $\fa = x_1\gA + \cdots + x_n\gA$,  il existe des \idms
$s_1,\dots,s_n\in \gA$ et des entiers $d_1,\dots,d_n \geq 1$
tels que $x_i^{d_i}\gA = s_i\gA$.
Soit 

\snic{s = 1-(1-s_1)\cdots(1-s_n),}

%\sni
de sorte que $s\gA = s_1\gA + \cdots + s_n\gA$.
Il est clair que l'\idm $s$ appartient à~$\fa$,
donc à toutes les puissances de~$\fa$.
D'autre part, si~$d \geq d_1 + \cdots + d_n - (n-1)$ on a

\snic{
\fa^d \subseteq x_1^{d_1}\gA + \cdots + x_n^{d_n}\gA = s_1\gA + \cdots + s_n\gA = s
\gA.}

%\sni
En conclusion $\fa^d = s\gA$.

Enfin, \emph{\ref{LID003}} implique clairement \emph{\ref{LID002}}.
\end{proof}

%--- Corollary{corZedReg}------
\begin{corollary}
\label{corZedReg}
%:HHH fidele et non pas regulier
Si $\fa$ est un \itf fidèle d'un anneau \zedz, alors $\fa=\gen{1}$.
En particulier, dans un anneau \zedz, tout \elt \ndz est \ivz.
\end{corollary}
%--- end-corollary------------------------------------
%-----------------begin proof------------------
\begin{proof}
Pour $d$ assez grand l'\id $\fa^d$ est engendré par un \idm $s$. Cet \id est
\ndzz, donc l'\idm $s$ est égal à 1.
\end{proof}
%-----------------end proof------------------

\begin{lemma}\label{lemZerloc}
\emph{(Anneaux \zeds locaux)}\\
\Propeq
%-----------------begin enum------------------
\begin{enumerate}\itemsep.5pt
\item   $\gA$ est \zed local.
\item   Tout \elt de $\gA$ est inversible ou nilpotent.
\item   $\gA$ est \zed et connexe.
\end{enumerate}
%-----------------end enum------------------
\end{lemma}

En conséquence un \cdi peut aussi être défini comme un \alo \zed
réduit.

\penalty-2500
%:--- SUBsec{Anneaux zed red}-
\subsec{Anneaux \zeds réduits}
%\addcontentsline{toc}{subsubsection}{Anneaux \zeds réduits}
%-------------------
\label{subsecAzedred}

\vspace{2pt}
%:--- Subsubsec*{Proprietes carac}
\subsubsec{Propriétés \carasz}

Les \eqvcs du lemme suivant sont faciles (voir la \dem du lemme
analogue \ref{lemme:idempotentDimension0}).
%-------------------
\begin{lemma}\label{lemZerRed}
\emph{(Anneaux \zeds réduits)}\\
\Propeq
%-----------------begin enum------------------
\begin{enumerate}%\itemsep1pt
\item  \label{i1lemZerRed}
L'anneau $\gA$ est \zed réduit.
%:2015  simplifie les points 2 et 4 qui suivent
\item  \label{i2lemZerRed}
Tout \idp est \idm
(i.e., $\forall a\in \gA,\;a\in\gen{a^{2}}$).
\item  \label{i3lemZerRed}
   Tout \idp est engendré par un \idmz.
\item  \label{i4lemZerRed}
   Tout \itf  est engendré par un \idmz.
\item  \label{i5lemZerRed}
 Pour toute liste finie  $(a_1,\ldots ,a_k)$ d'\elts de $\gA$, il existe des \idms \orts $(e_1,\ldots ,e_{k})$
tels que pour $j\in\lrbk$
$$\preskip.4em \postskip.4em 
\gen{a_1,\ldots ,a_j}=\gen{a_1e_1+\cdots +a_je_j}=\gen{e_1+\cdots+e_j}. 
$$
\item  \label{i6lemZerRed}
   Tout idéal  est \idmz.
\item  \label{i6bislemZerRed}
   Le produit de deux \ids est toujours égal à leur intersection.
\item  \label{i7lemZerRed}
   L'anneau $\gA$ est \qi et tout \elt \ndz est \ivz.
\end{enumerate}
%-----------------end enum------------------
\end{lemma}

%On en déduit le fait suivant.

%:     Fact{factQoQiZed}
\begin{fact}\label{factQoQiZed} 
\begin{enumerate}%\itemsep0pt
\item Soit $\gA$ un anneau arbitraire. Si $\Ann(a) = \gen{\vep}$  avec $\vep$ \idm,
alors l'\elt $b = a+\vep$ est \ndz et $ab = a^2$.
\item Si $\gA$ est \qiz, $\Frac\gA$ est \zedr et
tout \idm de $\Frac\gA$ est dans $\gA$.
\end{enumerate}
\end{fact}
\begin{proof}
\emph {1.}
Regarder modulo $\vep$ et modulo $1-\vep$.

\emph {2.}
Pour un $a\in\gA$, on doit trouver $x\in\Frac\gA$ tel que $a^2x=a$.\\
On pose $b=a+(1-e_a)\in\Reg\gA$ donc $ab=a^{2}$, et l'on prend $x=b^{-1}$.
\\
Soit maintenant $a/b$ un \idm de $\Frac\gA$. On a $a^2=ab$.\\
\hspace*{.1em} --- Modulo ${1-e_a}$, on a $b=a$  et $a/b=1=e_a$ (car $a$ est \ndzz).\\ 
\hspace*{.1em} --- Modulo ${e_a}$, on a $a/b=0=e_a$ (car $a=0$). En résumé $a/b=e_a$. 
\end{proof}

%--- Fact{factZerRedCoh}-------
\begin{fact}
\label{factZerRedCoh}
Un anneau \zedr est \cohz.
Il est \fdi \ssi il y a un test d'\egt à zéro pour les \idmsz.
\end{fact}
%--- end-fact-----------------------------------------

\exl Soit $\PP$ l'ensemble des nombres premiers. L'anneau {\mathrigid 2mu $\gA=\prod_{p\in\PP} \ZZ/ \!\!\gen{p}$} est \zedr mais il \emph{n'est pas} dicret.
\eoe

\smallskip On obtient aussi facilement les \eqvcs suivantes.

%--- Fact{factZerRedConnexe}-------
\begin{fact}
\label{factZerRedConnexe}
Pour un anneau \zed  $\gA$ \propeq
\begin{enumerate}
  \item $\gA$ est connexe (resp. $\gA$ est connexe et réduit).
  \item $\gA$ est local (resp. $\gA$ est local et réduit).
  \item $\Ared$ est intègre (resp. $\gA$ est intègre).
  \item $\Ared$ est un \cdi (resp. $\gA$ est un \cdiz).
\end{enumerate} 
\end{fact}
%--- end-fact----------------------

%:--- Subsubsec*{Definition equationnelle}
\subsubsect{Définition équationnelle des anneaux \zeds réduits}
{Définition équationnelle}
%-------------------

%\vspace{2pt}
Un anneau non \ncrt commutatif vérifiant
$$\preskip.4em \postskip.4em
\forall x\,\exists a\quad xax=x
$$
est
souvent qualifié de \ix{Von Neumann régulier}. Dans le cas commutatif ce sont
les anneaux \zeds réduits. On les appelle encore \emph{anneaux absolument
plats}, parce qu'ils sont \egmt \cares par
la \prt suivante: tout \Amo est plat
(voir la
proposition~\ref{propEvcPlat}).\index{anneau!absolument plat}

Dans un anneau commutatif, deux \elts $a$ et $b$ sont dits \emph{quasi inverses}
si l'on~a:
%------begin equation--eqQuasiInv-----------
\begin{equation}\label{eqQuasiInv}\preskip.0em \postskip.4em
a^2b=a,\quad \quad  b^2a=b
\end{equation}
%---------------------end equation--------------

On dit aussi que $b$ est \und{le} \ix{quasi inverse} de $a$. On vérifie en effet
qu'il est unique: si  $a^2b=a=a^2c$,
 $b^2a=b$  et  $c^2a=c$, alors
$$\preskip.4em \postskip.4em
%\snic{
c-b=a(c^2-b^2)=a(c-b)(c+b)=a^2(c-b)(c^2+b^2)=0,
%}
$$
%\sni
puisque $ab=a^2b^2$, $ac=a^2c^2$ et   $a^2(c-b)=a-a=0$.

Par ailleurs, si $x^2y=x$, on vérifie que $xy^2$ est \qiv de $x$.

Ainsi: \emph{un anneau est \zedr \ssi tout \elt admet un \qivz.}

Comme le \qiv est unique, un anneau \zedr peut être vu comme un
anneau muni d'une loi unaire supplémentaire, $a\mapsto a\bul$
soumise à l'axiome (\ref{eqQuasiInv}) avec $a\bul$ à la place de $b$.\\
Notons que $(a\bul)\bul=a$ et $(a_1a_2)\bul=a_1\bul a_2\bul$.

%--- Fact{factQuasiInv}---------------
\begin{fact}
\label{factQuasiInv}
Un anneau \zedr $\gA$ est \qiz, avec l'\idm $e_a=aa\bul$:
$\Ann(a)=\gen{1-e_a}$.
On a $\gA\simeq \gA[1/e_a]\times \aqo{\gA}{e_a}$.
Dans $\gA[1/e_a]$, $a$ est inversible, et dans $\aqo{\gA}{e_a}$, $a$ est nul.
\end{fact}
%--- end-fact-----------------------------------------

% --- Subsubsec*{Lemme de scindage \zed}
\subsubsection*{Lemme de scindage \zedz}

Le lemme de scindage suivant est à peu près immédiat,
il se démontre comme le lemme de scindage \qi \ref{thScindageQi}.

%     Lemma{thScindageZed}
\begin{lemma}\label{thScindageZed}
Soit $(x_i)_{i\in I}$ une famille finie d'\elts dans un anneau \zed $\gA$.
Il existe  un \sfio $(e_1, \ldots, e_n)$ tel que dans chaque
composante $\gA[1/e_j]$, chaque $x_i$ est nilpotent ou \ivz.
\end{lemma}
%

%:--- Subsubsec*{Machinerie locale-glo  elr 2}
\subsubsect{Des \cdis aux anneaux \zeds réduits}{Machinerie locale-globale \elr \num2:
des \cdis aux anneaux \zeds réduits}

Les anneaux \zeds réduits ressemblent beaucoup à des produits finis de
\cdisz, et cela se manifeste \prmt comme suit.

%:   Machinerie locale-globale \elr \num2
\medskip 
{\bf Machinerie locale-globale \elr \num2.
%des \cdis aux anneaux \zeds réduits.
}\label{MethodeZedRed}\imlg
{\it La plupart des \algos qui fonctionnent avec les \cdis non triviaux
peuvent être modifiés de manière à fonctionner avec les anneaux \zeds
réduits, en scindant l'anneau
en deux composantes chaque fois que l'\algo écrit pour les \cdis utilise
le test
\gui{cet \elt est-il nul ou \ivz?}.
Dans la première composante l'\elt en question
est nul, dans la seconde il est \ivz.}

\medskip
\rems \label{remzedred1}
1) Nous avons mis \gui{la plupart} plutôt que \gui{tous} dans la mesure où
l'énoncé du résultat de l'\algo pour les \cdis doit être écrit
sous une forme où n'apparaît pas qu'un \cdi est connexe.

2) Par ailleurs, la même remarque que celle que nous avons faite
\paref{remMLGE1} concernant la
machinerie \lgbe \elr \num1 s'applique encore. L'\algo donné dans le cas des \cdis doit être suffisamment uniforme pour ne pas conduire à un arbre infini lorsque l'on  veut le transformer
en un \algo pour les anneaux \zeds réduits.
\eoe

\smallskip Nous donnons tout de suite un exemple d'application de cette machinerie.

%--- Proposition{propZedBez}-----
\begin{proposition}\label{exlZedBez}
\label{propZedBez} Pour un anneau $\gA$ \propeq
%-----------------begin enum------------------
\begin{enumerate}
\item $\gA$ est un anneau \zedrz.
\item $\gA[X]$ est un anneau de Bézout strict et \qiz.
\item $\gA[X]$ est un anneau de Bézout.
\end{enumerate}
%-----------------end enum------------------
\end{proposition}
%--- end-proposition----------------------------------------
%-----------------begin proof------------------
\begin{proof}
\emph{1 $\Rightarrow$ 2.}
Le fait est classique pour un \cdiz: on utilise l'\algo
d'Euclide étendu pour calculer sous la forme $g(X)=a(X)u(X)+b(X)v(X)$
un pgcd de $a(X)$ et $b(X)$. En outre, on obtient une matrice $\crmatrix{u& -b_1\cr
v&a_1}$ de \deter $1$
qui transforme $[\,a \;\;b\,]$ en $[\,g \;\;0\,]$.
Cette matrice est le produit de matrices $\cmatrix{0&-1\cr 1&-q_i}$
où les $q_i$ sont les quotients successifs.\\
Passons maintenant au cas d'un anneau \zedrz, donc \qiz.
Tout d'abord $\gA[X]$ est \qi car l'annulateur d'un \pol est l'intersection
des annulateurs de ses \coes 
(voir le corolaire~\ref{corlemdArtin}~\emph{\iref{i2corlemdArtin}}), donc engendré par le produit des \idms
correspondants.
Concernant le caractère \gui{Bézout strict} l'\algo qui vient d'être
expliqué pour les \cdis bute a priori sur l'obstacle de la non
inversibilité des \coes dominants dans les divisions successives. Mais cet
obstacle est à chaque fois
contourné par la considération d'un \idmz~$e_i$ convenable, l'annulateur du
\coe à inverser. Dans $\gA_i[1/e_i]$, (\hbox{où $\gA_i=\gA[1/u_i]$} est l'anneau
\gui{en cours} avec un certain \idm $u_i$) le \pol diviseur a un degré plus
petit que prévu et l'on recommence avec le \coe suivant. Dans $\gA_i[1/f_i]$,
($f_i=1-e_i$ dans $\gA_i$), le \coe dominant du diviseur est inversible et la
division peut être exécutée.
On obtient de cette façon un arbre de calcul aux feuilles duquel
on a le résultat souhaité. \`A chaque feuille le résultat est obtenu dans un
localisé $\gA[1/h]$ pour un certain \idm $h$. Et les $h$  aux feuilles de
l'arbre forment un \sfioz. Ceci permet
de recoller toutes les \egts en une seule\footnote{Pour une preuve
plus directe, voir l'exercice~\ref{exoZerRedBez}}.

 \emph{3 $\Rightarrow$ 1.}
Cela résulte du lemme suivant.

\textbf{Lemme.} \textit{Pour un anneau arbitraire $\gA$, si l'\id $\gen{a,X}$ est un \idp de $\gA[X]$, alors $\gen{a}=\gen{e}$ pour un certain \idm $e$.}

 On suppose que  $\gen{a,X}=\gen{p(X)}$ avec $p(X)q(X)=X$.
On a donc 

\snic{\gen{a}=\gen{p(0)}$,  $\;p(0)q(0)=0\;$
et $\;1=p(0)q'(0)+p'(0)q(0),}

%\sni
d'où $p(0)=p(0)^2q'(0)$. Ainsi, $e=p(0)q'(0)$ est \idm
et $\gen{a}=\gen{e}$.
\end{proof}
%-----------------end proof------------------

\rem La notion d'anneau \zedr peut être vue comme l'analogue
non-\noe de la notion de \cdiz, puisque si l'\agB des \idms est infinie, la \noet est perdue. Illustrons ceci sur l'exemple du \nstz, pour lequel ce n'est pas clair a priori si la \noet est un ingrédient essentiel ou un simple accident.
Un énoncé \cof  précis du \nst de Hilbert\ihi (forme faible) 
se formule comme suit.

%:2018  ajout d'un titre
\smallskip \noindent \textbf{\nst faible classique.} \emph{Soit $\gk$ un \cdi non trivial, $(\lfs)$ une liste d'\elts de $\kuX$,
et $\gA=\aqo\kuX\uf$ l'\alg quotient.
Alors, ou \hbox{bien $1\in\gen{\lfs}$}, ou bien il existe un quotient  
de $\gA$ qui est un \kev non nul de dimension finie}. \eoe

\smallskip  Comme la preuve est donnée par un \algo uniforme (pour plus de précisions voir le \thref{thNstNoe} et l'exercice \ref{exothNst1-zed})
 on obtient par application de la machinerie locale-globale \elr \num2 
 le résultat ci-après, sans disjonction, qui implique le \nst
  précédent pour un \cdi non trivial (cet exemple illustre aussi la première remarque \paref{remzedred1}).
  Un \Amo $M$ est dit \emph{quasi libre}\index{quasi libre!module ---}\index{module!quasi libre} s'il est isomorphe à une somme directe finie d'idéaux $\gen{e_i}$ avec les $e_i$ \idmsz.
On peut alors en outre imposer que  $e_ie_j=e_j$ si $j>i$, car  pour deux \idms
$e$ et $f$, on a

\snic{\gen{e}\oplus\gen{f}\simeq\gen{e\vu f}\oplus\gen{e\vi f}$, où
$e\vi f=ef$ et $e\vu f=e+f-ef.}

%:2018  ajout d'un titre, ajout de fidèle, i.e.
\smallskip \noindent \textbf{\nst faible pour un anneau \zedrz.}\\
\emph{Soit $\gk$
 un anneau \zedrz, $(\lfs)$ une liste d'\elts  de $\kXn$ et $\gA$ l'\alg quotient. Alors, l'\id $\gen{\lfs}\,\cap\, \gk$ est engendré par un idempotent ${e}$, et en posant $\gk_1=\aqo{\gk}{e}$,  
il existe un quotient~$\gB$ de $\gA$ qui est un $\gk_1$-module quasi libre  fidèle, i.e. avec l'\homo naturel $\gk_1\to\gB$ injectif}.
\eoe

%: Subsubsec*{Mpfs zed red}
\subsubsect{Modules \pf sur les anneaux \zeds réduits}{Modules \pfz}
%-------------------

%:     Theorem{thZerDimRedLib}
\begin{theorem}\label{thZerDimRedLib}
\emph{(Le paradis des anneaux \zedrsz)}\\
Soit $\gA$ un anneau \zedrz.
\begin{enumerate}
\item Toute matrice est \eqve à une matrice en forme de Smith
avec des \idms sur la diagonale principale.
\item Tout \mpf est quasi libre.
\item
Tout sous-module de type fini d'un \mpf est facteur
direct.
\end{enumerate}
\end{theorem}
\begin{proof}
Les résultats sont classiques
pour le cas des \cdis (une \prco peut être
basée sur la méthode du pivot).
La machinerie \lgbe \elr \num2 fournit alors (pour chacun des trois points)
le résultat séparément dans chacun des $\gA[1/e_j]$,
après avoir scindé l'anneau en un produit de localisés
$\gA[1/e_j]$ pour un 
\sfio $(e_1, $\ldots$, e_k)$.
Mais le résultat est justement formulé de façon
à être vrai globalement dès qu'il est vrai dans chacune
des composantes.\imlg
\end{proof}

%%%%%%%%%%%%%%%%%%%%%%%%%%%%%%%%%%%%%%%%%%%%%%%%%%%%%%%%%%%%%%%%%%%%%%%%%%
%%%%%%%%%%%%%%%%%%%%%%%%%%%%%%%%%%%%%%%%%%%%%%%%%%%%%%%%%%%%%%%%%%%%%%%%%%
%%%%%%%%%                                                   %%%%%%%%%%%%%%
%%%%%%%%%             subsec{Systèmes \polls \zedsz}      %%%%%%%%%%%%%%
%%%%%%%%%                                                   %%%%%%%%%%%%%%
%%%%%%%%%%%%%%%%%%%%%%%%%%%%%%%%%%%%%%%%%%%%%%%%%%%%%%%%%%%%%%%%%%%%%%%%%%
%%%%%%%%%%%%%%%%%%%%%%%%%%%%%%%%%%%%%%%%%%%%%%%%%%%%%%%%%%%%%%%%%%%%%%%%%%

%:HHH   section rajoutee (supprimee en chap 9))	
%:--- SUBsec{Systèmes \polls \zedsz}-
\subsec{Systèmes \polls \zedsz}
\label{subsecSypZerdim}

%%: subsubsec{Rappel sur la mise en position de \Noez}
%\subsubsec{Rappel sur la mise en position de \Noez}%\label{EntSurZdim}

Nous étudions dans ce paragraphe un exemple particulièrement important d'anneau \zedz, fourni
par les \algs quotients associées aux \syps \zeds sur les \cdisz.

Rappelons le contexte étudié dans la section \ref{secChap3Nst} consacrée au \nst de Hilbert.
Si $\gK\subseteq\gL$ sont des \cdisz, et si $(\lfs)$ est un \syp dans $\KXn=\KuX$,
on dit que $(\xin)=(\uxi)$ est un {zéro de $\uf$ dans $\gL^n$} si les \eqns $f_i(\uxi)=0$ sont satisfaites.
\\
L'étude de la \vrt des zéros du \sys est étroitement reliée à celle de l'\emph{\alg quotient associée au \sypz}, à savoir 

\snic{\gA=\aqo\KuX\uf=\Kux \quad (x_i \hbox{ est la classe de } X_i \hbox{ dans }\gA).}

%\sni
En effet, il revient au même de se donner un zéro $(\uxi)$ du \syp dans $\gL^n$ ou de se donner un homomorphisme  de \Klgs  $\psi:\gA\to\gL$ ($\psi$~est défini par $\psi(x_i)=\xi_i$ pour $i\in\lrbn$).
Pour $h\in\gA$, on note $h(\uxi)=\psi(h)$ l'\evn de $h$ en $\uxi$.

%Dans la suite nous notons $\ff=\gen{\uf}$.
% et $x_i$ la classe de $X_i$ dans $\gA$.
Lorsque $\gK$ est infini, le \thref{thNstfaibleClass} nous donne une mise en position de \iNoe
par un \cdv \linz, et un entier $r\in\lrb{-1..n}$ satisfaisant les \prts suivantes\footnote{Dans l'hypothèse du \thref{thNstfaibleClass}, on a mis
que $\gK$ est contenu dans un \cac $\gL$, mais la \dem des \prts que nous signalons ici, basée sur le lemme \ref{lemElimPlusieurs} et sur le lemme de changement de variables \ref{lemCDV}, n'utilise pas l'existence du corps $\gL$.}
(nous ne changeons pas le nom de variables, ce qui est un léger abus). 
\begin{enumerate}
\item Si $r=-1$, alors $\gA=0$, \cad $1\in\gen{\uf}$.
\item Si $r=0$, chaque $x_i$ est entier sur $\gK$, et $\gA\neq 0$.
\item Si $0<r<n$,  alors $\KXr\,\cap\,\gen{\uf}=0$ 
et les $x_{i}$ pour ${i\in\lrb{r+1..n}}$ sont entiers sur  $\Kxr$ (qui est isomorphe à $\KXr$).
\item Si $r=n$, $\gen{\uf}=0$ et $\gA=\KuX$
\end{enumerate}

%--- SUBsubsec*{Algèbres \stfes sur un \cdiz} -
%\subsubsection*{Algèbres \stfes sur un \cdiz}
%-----------------------------------------

%:     Lemma{lemSypZD1}
\begin{lemma}\label{lemSypZD1}
\emph{(Précisions sur le \thref{thNstfaibleClass})}%
\begin{enumerate}
\item Dans le cas où $r=0$, l'\alg quotient $\gA$ est finie sur $\gK$.
\item Si l'\alg quotient $\gA$ est finie sur $\gK$, elle est \stfe sur~$\gK$,
et c'est un anneau \zedz.
On dit alors que \emph{le \syp est \zedz}.

\item Si  l'anneau  $\gA$ est  \zedz, alors $r\leq 0$.
\item Toute \alg \stfe sur le \cdi $\gK$ peut être regardée
comme (est isomorphe à) l'\alg quotient d'un \syp \zed sur $\gK$.
\end{enumerate}\index{zero-dimensionnel@\zedz!systeme@\syp ---}\index{systeme polynomial@\sypz!zéro-dimensionnel}
\end{lemma}
%--------- fin lemma ----------------------------------------------
%
\begin{proof}
\emph{1.} 
En effet, si $x_i$ annule le \polu $p_i\in\KT$, l'\alg $\gA$ est un quotient
de $\gB=\KuX\big/\geN{\big(p_i(X_i)\big)\phantom{\!\!)}_{i\in\lrbn}}$, qui est un \Kev de dimension finie.

\emph{2.} On commence comme au point \emph{1} Pour obtenir l'\alg $\gA$, il nous suffit de prendre le quotient de $\gB$ par l'\id $\gen{f_1(\uz),\dots,f_s(\uz)}$ (où les $z_i$ sont les classes des $X_i$ dans $\gB$). On voit facilement que cet \id est un sous-\evc de type fini de $\gB$, donc le quotient est de nouveau un \Kev
de dimension finie. Ainsi, $\gA$ est \stfe sur $\gK$.\\
Montrons que $\gA$ est \zedz. Tout $x\in\gA$ annule son \polminz, disons $f(T)$, de sorte que  l'on a une \egt $x^k\big(1+xg(x)\big)=0$
(multiplier $f$ par l'inverse du \coe de plus bas degré non nul).

\emph{4.} L'\alg $\gA$ est engendrée par un nombre fini d'\elts $x_i$  
(par exemple une base comme \Kevz), qui annulent chacun leur \polminz,
disons $p_i(T)$. Ainsi $\gA$ est un quotient d'une \alg 

\snic{\gB=\KuX\big/\geN{\big(p_i(X_i)\big)\phantom{\!\!)}_{i\in\lrbn}}=\gA[\zn].}

%\sni
Le morphisme surjectif correspondant, de $\gB$ sur $\gA$, est une \ali
dont on peut calculer le noyau (puisque $\gA$ et $\gB$ sont de dimension finie), par exemple en précisant un \sgr $\big(g_1(\uz),\dots,g_\ell(\uz)\big)$.
En conclusion,
 l'\alg $\gA$ est isomorphe à l'\alg quotient associée au \syp
$\big(p_1(X_1),\dots,p_n(X_n),g_1(\uX),\dots,g_\ell(\uX)\big)$.

\emph{3.}  Ce point résulte des deux lemmes qui suivent.
\end{proof}

\rem Traditionnellement, on réserve l'appellation de \syp \zed au cas $r=0$, mais l'\alg quotient est \zede aussi lorsque $r=-1$.
\eoe

%:     Lemma{lemnonZR}
\begin{lemma}\label{lemnonZR}
Si l'anneau $\gC[\Xr]$ est \zed avec $r>0$, alors l'anneau $\gC$ est trivial.  
\end{lemma}
%--------- fin lemma ---------------------------------------------- 
%
\begin{proof}
On écrit $X_1^m\big(1-X_1P(\Xr)\big)=0$. Le \coe de $X_1^m$ est à la fois égal à $0$ et à $1$.
\end{proof}
%

%--- Lemma{lemZrZr2}----------------
\begin{lemma}
\label{lemZrZr2}
Soit~$\gk\subseteq \gA$,~$\gA$ entière sur~$\gk$. Si
$\gA$ est un anneau \zedz,
$\gk$ est un anneau \zedz.
\end{lemma}
%--- end-lemma-----------------------------------------
%-----------------begin proof------------------
\begin{proof}
Soit $x\in\gk$, on a un $y\in\gA$ tel que $x^{k}=yx^{k+1}$.
Supposons par exemple que $y^3+b_2y^2+b_1y+b_0=0$ avec les $b_i\in\gk$.
\\
Alors,  $x^{k}=yx^{k+1}=y^2x^{k+2}=y^3x^{k+3}$, et donc

\snic{\arraycolsep2pt\begin{array}{rclcl}
0& =  &  (y^3+b_2y^2+b_1y+b_0)x^{k+3} &   &   \\[1mm]
& =   & x^{k}+b_2x^{k+1}+b_1x^{k+2}+b_0x^{k+3}  & =  & x^{k}\big(1+x(b_2+b_1x+b_0x^2)\big).
\end{array}}

\vspace{-1em}
\end{proof}
%-----------------end proof------------------

%:   Theorem{thSPolZed}----------
\begin{theorem}
\label{thSPolZed} \emph{(\Sys \zed sur un \cdiz)}\\
Soit~$\gK$  un \cdi et $(\lfs)$ dans $\KXn=\KuX$.\\
 Notons~$\gA=\aqo\KuX\uf$ l'\alg quotient associée à ce \sypz.\\
\Propeq
%-----------------begin enum------------------
\begin{enumerate}
\item $\gA$ est finie sur~$\gK$.
\item $\gA$ est \stfe sur~$\gK$.
\item $\gA$ est un anneau \zedz.
\end{enumerate}
%-----------------end enum------------------
Si~$\gK$ est contenu dans un \cdacz~$\gL$, ces \prts sont aussi
\eqves aux suivantes.
\begin{enumerate} \setcounter{enumi}{3}
\item Le \syp a un nombre fini de zéros dans~$\gL^n$.
\item Le \syp a un nombre borné de zéros dans~$\gL^n$.
\end{enumerate}%
\index{zero-dimensionnel@\zedz!systeme@\syp ---}\index{systeme polynomial@\sypz!zéro-dimensionnel}
\end{theorem}
%--- end-theorem-----------------------------------------
%-----------------begin proof------------------
\begin{proof} Lorsque $\gK$ est infini, on obtient les \eqvcs
en appliquant le lemme~\ref{lemSypZD1} et le \thref{thNstfaibleClass}.
\\
Dans le cas \gnlz, on peut \egmt obtenir une
 mise en position de \iNoe en utilisant 
un \cdv \gnl non \ncrt \lin comme celui donné en~\ref{lemNoether}
(voir le \thref{thNstNoe}).  
\end{proof}
%-----------------end proof------------------

%:HHH rajout
Une variation sur  le \tho précédent est donnée au \thref{thNst0} 

\medskip 
\rem Plutôt que d'utiliser un \cdv non \lin comme proposé dans la \dem précédente, on peut recourir à la technique de \gui{changement de corps de base}. Cela fonctionne comme suit. On considère un corps infini $\gK_1\supseteq \gK$, par exemple $\gK_1=\gK(t)$, ou $\gK_1$ un \cac contenant $\gK$ si l'on sait en construire un. 
Alors l'\eqvc des points \emph{1}, \emph{2} et \emph{3} est assurée
pour l'\alg $\gA_1$ pour le même \syp vu sur $\gK_1$. L'\alg $\gA_1$
est obtenue à partir de $\gA$ par \eds de $\gK$ à $\gK_1$.
Il reste à voir que chacun des trois points est vérifié pour $\gA$
\ssi il est vérifié pour $\gA_1$. Ce que nous laisserons \alecz\footnote{On pourra voir à ce sujet les \thos \ref{thSurZedFidPlat}, \ref{propFidPlatTf} et~\ref{propFidPlatPrAlg}}. 
\eoe

%--- SUBsubsec*{Le \tho de Stickelberger} -
%\subsubsection*{Le \tho de Stickelberger}
%-----------------------------------------

%:    Theorem{thStickelberger}
\begin{theorem}\label{thStickelberger}  \emph{(\Tho de Stickelberger)}\\
Même contexte que le \thref{thSPolZed}, avec maintenant $\gK$  \cacz.%
\index{Stickelberger!\tho de ---}
\begin{enumerate}
\item Le \syp admet un nombre fini de zéros sur $\gK$.  
\\
On les note
$\uxi_1$, \ldots, $\uxi_\ell$. 
\item Pour chaque $\uxi_k$ il existe un \idm $e_k\in\gA$
satisfaisant $e_k(\uxi_j)=\delta_{j,k}$ (symbole de Kronecker)
pour tous  $j\in\lrb{1..\ell}$. 
\item Les \idms $(e_1,\dots,e_\ell)$ forment un \sfio de $\gA$.
\item Chaque \alg $\gA[1/e_k]$ est un \alo \zed (tout \elt est \iv ou nilpotent).
\item
Notons $m_k$   la dimension du \Kev $\gA[1/e_k]$. \\
On a
$[\gA:\gK]=\som_{k=1}^\ell m_k$ et
pour tout $h\in\gA$ on a

\snic{\rC{\gA\sur\gK}(h)(T)=\prod\nolimits_{k=1}^\ell\big(T-h(\uxi_k)\big)^{m_k}.
}

%\sni
En particulier,   $\Tr\iAK(h)=\som_{k=1}^\ell m_kh(\uxi_k)$ et
$\rN\iAK(h)=\prod_{k=1}^\ell h(\uxi_k)^{m_k}$.
\item Notons $\pi_k:\gA\to\gK,\;h\mapsto h(\uxi_k)$ l'\evn en $\uxi_k$,
et $\fm_k=\Ker\pi_k$. Alors $\gen{e_k-1}=\fm_k^{m_k}$ et $\fm_k=\sqrt{\gen{e_k-1}}$.

\end{enumerate}
\end{theorem}
%--------- fin theorem ---------------------------------------------- 
%
\begin{proof} On note $V=\big\{\uxi_1,\dots,\uxi_\ell\big\}$ 
la \vrt des zéros du \sys dans $\gK^n$. 

\emph{2} et \emph{3.}
On a des interpolants de Lagrange en plusieurs variables. Ce sont des \pols $L_k\in\KuX$ 
qui vérifient \smash{$L_k(\uxi_j)=\delta_{j,k}$}. On considère les~$L_k$ comme des \elts de
$\gA$.\\ 
Puisque $\gA$ est \zedz, il existe un entier $d$ et un \idm $e_k$ avec $\gen{e_k}=\gen{L_k}^d$, donc $e_kL_k^d=L_k^d$ et $L_k^db_k=e_k$ pour un certain $b_k$.
Ceci implique que $e_k(\uxi_j)=\delta_{j,k}$. \\
Pour $j\neq k$, $e_je_k$ est nul sur $V$, donc par le \nstz, $e_je_k$
est nilpotent dans $\gA$. Comme c'est un \idmz, $e_je_k=0$.\\
La somme des $e_j$ est donc un \idm $e$. Cet \elt ne s'annule nulle part, \cad qu'il a les mêmes zéros que $1$. Par le \nstz, on obtient $1\in\sqrt{\gen{e}}$. Ainsi $e=1$ car
c'est un \idm \iv de $\gA$.

\emph{4.} La \Klg $\gA_k=\gA[1/e_k]=\aqo\gA{1-e_k}$ est l'\alg quotient associée au \syp $(\lfs,1-e_k)$ qui admet $\uxi_k$ pour seul zéro.
On considère un \elt arbitraire $h\in\gA_k$. En raisonnant comme au point précédent, on obtient par le \nst que
si~$h(\uxi_k)=0$, alors~$h$ est nilpotent, et si~$h(\uxi_k)\neq 0$, 
alors~$h$ est \ivz.

\emph{5.} Puisque $\gA\simeq\prod_{k=1}^\ell\gA_k$, il suffit de démontrer que pour $h\in\gA_k$, on a l'\egt $\rC{\gA_k\sur\gk}(h)(T)=\big(T-h(\uxi_k)\big)^{m_k}$. On identifie $\gK$ à son image dans~$\gA_k$. L'\elt $h_k=h-h(\uxi_k)$ s'annule en $\uxi_k$, donc il est nilpotent. Si~$\mu$ désigne la multiplication par $h_k$ dans $\gA_k$, $\mu$ est un \endo nilpotent.
Sur une base convenable, sa matrice est triangulaire stricte inférieure, et celle de la multiplication par $h$ est triangulaire avec des $h(\uxi_k)$ sur la diagonale, donc son \polcar est $\big(T-h(\uxi_k)\big)^{m_k}$.

\emph{6.} On a clairement $e_k-1\in\fm_k$. Si $h\in\fm_k$, l'\elt $e_kh$
est partout nul sur $V$, donc nilpotent. Donc $h^Ne_k=0$ pour un certain $N$ et $h\in\sqrt{\gen{e_k-1}}$. Pour voir que $\fm_k^{m_k}=\gen{e_k-1}$,
on peut se situer dans $\gA_k$, où $\gen{e_k-1}=0$. Dans cet anneau, l'\id $\fm_k$ est un \Kev de dimension~$m_k-1$. Les puissances successives de $\fm_k$
forment alors une suite décroissante de sous-\Kevs de dimensions finies, qui stationne dès que deux termes consécutifs sont égaux.
Ainsi $\fm_k^{m_k}$ est un \itf \idm strict, donc nul. 
\end{proof}

\rems~\\ 1) Le fait que le \syp est \zed résulte d'un calcul rationnel
dans le corps des \coes (mise en position de \iNoe ou calcul de \bdgz). 

 2) Le point \emph{5} 
du \tho de Stickelberger permet de  calculer toutes les informations utiles sur les zéros du \sys
en se basant sur la seule forme trace. En outre, la forme trace peut être 
calculée dans le corps des \coes des \pols du \sysz.
Ceci a des applications importantes en calcul formel (voir par exemple \cite{BPR}).
\eoe

\medskip  Pour des exemples, on pourra consulter l'exercice \ref{exoFreeAlgebraPresentation}
et le \pb \ref{exoDimZeroXcYbZa}. Pour une étude purement locale des zéros isolés, voir la section \ref{secExlocGeoAlg}.

%--- Sec Idéaux de Fitting {sec Fitt}--------
\section{Idéaux de Fitting}
\label{sec Fitt}

La théorie des \idfs des \mpfs est une machinerie calculatoire extrêmement
efficace d'un point de vue théorique \cofz. Elle a un coté \gui{théorie de
l'élimination} dans la mesure où elle est entièrement basée sur des
calculs de \detersz, et elle a pendant un temps plus ou moins disparu 
de la littérature sous l'influence de l'idée qu'il fallait \gui{éliminer
l'élimination} pour se sortir du bourbier de calculs dont la signification ne
semblait pas claire.

Les \idfs redeviennent à la mode et c'est tant mieux.
Pour plus de détails, on pourra consulter \cite{Nor}.

%: subsec {sec Fitt pres fin}--
\subsec{Idéaux de Fitting d'un \mpfz}
\label{sec Fitt pres fin}

%--- Defi{ideaux Fitting}----
\begin{definition}\label{def ide fit}\\
Si $G\in\gA^{q\times m}$ est une \mpn d'un \Amo $M$ donné par~$q$ \gtrsz, les \emph{\idfs de} $M$ sont les \ids

\snic{\cF_{\gA,n}(M)=\cF_{n}(M):= \cD_{\gA,q-n}(G)}

%\sni
où $n$ est un entier arbitraire.\index{ideal@idéal!de Fitting}%
\index{Fitting!idéal de ---}
\end{definition}
%--------------

Cette \dfn est légitimée par le lemme facile mais fondamental suivant.

%--- Lemma{lemFitInv}-----------
\begin{lemma}
\label{lemFitInv}
Les \idfs du \mpf $M$ sont bien définis, autrement dit ces idéaux ne
dépendent pas de la \pn choisie $G$ pour $M$.
\end{lemma}
%--- end-lemma-----------------------------------------
\begin{proof}
Pour prouver ce lemme  il faut essentiellement voir que les idéaux  $\cD_{q-
n}(G)$ ne changent pas,
%-----------------begin item------------------
\begin{enumerate}
\item d'une part, lorsque l'on  ajoute une nouvelle syzygie, \coli des syzygies déjà
présentes,
\item d'autre part, lorsque l'on  ajoute un nouvel \elt à un \sgrz, avec une syzygie qui
exprime ce nouvel \elt en fonction des anciens \gtrsz.
\end{enumerate}
%-----------------end item------------------
Les détails sont laissés \alecz.
\end{proof}

\perso{Consulter le livre de Sarah Glaz sur les anneaux commutatifs cohérents
p97 et suivantes, il y a des choses intéressantes sur les Fitting.}
On a immédiatement les faits suivants.

%--- Fact{fact.idf inc}-----------
\begin{fact}\label{fact.idf inc}
Pour tout \mpf  $M$ avec $q$  \gtrsz, on a les inclusions

\snic{\gen{0} = \cF_{-1}(M) \subseteq \cF_{0}(M) \subseteq \cdots \subseteq
\cF_q(M)= \gen{1}.}

%\sni
Si $N$ est un \mpf quotient de $M$, on a les \linebreak 
inclusions $\cF_{k}(M) \subseteq
\cF_{k}(N)$ pour tout $k\geq 0$.
\end{fact}
%--------------

\rem En particulier, si $\cF_r(M)\neq \gen{1}$ le module $M$ ne peut pas être
engendré par $r$ \eltsz. On verra (lemme du nombre de \gtrs local
\paref{lemnbgtrlo}) que la signification
de l'\egt $\cF_r(M)= \gen{1}$ est que le module est
\emph{localement} engendré par
$r$ \eltsz.
\eoe

%--- Fact{fact.idf libre}---------
\begin{fact}\label{fact.idf libre}
Soit $M$ un \Amo libre de rang $k$. Alors,

\snic{\cF_{0}(M) = \cdots = \cF_{k-1}(M)  =  \gen{0} \subseteq
\cF_k(M)=\gen{1}.}

%\sni
Plus \gnltz, si $M$ est quasi libre isomorphe à $\bigoplus_{i=1}^{k}\gen{f_i}$, où les~$f_i$ sont des \idms tels que  $f_if_j=f_j$ si $j>i$,
alors   $\cF_k(M)=\gen{1}$ et~$\cF_i(M)=\gen{1-f_{i+1}}$ pour $ i\in\lrb{0..k-1}$.
\end{fact}
%--------------

Remarquez que ceci fournit une preuve savante du fait que si un  module
est libre avec deux rangs distincts, l'anneau est trivial.

\ms\exls ~
\\ 1. Pour un groupe abélien fini $H$ considéré comme $\ZZ$-module, l'\id $\cF_0(H)$ est engendré par l'ordre
du groupe tandis que l'annulateur est
engendré par son exposant. 
En outre, la structure du groupe est entièrement \caree par
ses \idfsz. Une \gnn est donnée dans l'exercice~\ref{exoFitt0}.

\noi 2. Reprenons le \Bmo $M$ de l'exemple \paref{belexemple}.
Le calcul donne les résultats suivants.
%-----------------begin item------------------
\begin{itemize}
\item Pour $M$: $\cF_0(M)=0$, $\cF_1(M)=\fb$ et $\cF_2(M)=\gen{1}$,
\item pour  $M'=M\te M$: $\cF_0(M')=0$,  $\cF_1=\fb^3$,  $\cF_2=\fb^2$,  $\cF_3=\fb$ et $\cF_4=\gen{1}$,
\item pour  $M''=\gS^2(M)$: $\cF_0 (M'')=0$, $\cF_1=\fb^2$,
$\cF_2=\fb$ et
$\cF_3=\gen{1}$,
\item pour  $\Al2M$:
 $\cF_0(\Al2M)=\fb$ et  $\cF_1(\Al2M)=\gen{1}$.
 \eoe
\end{itemize}
%-----------------end item------------------

%--- Fact{fact.idf.change}---
\begin{fact}\label{fact.idf.change}\label{fact.idf loc}
\emph{(Changement d'anneau de base)}\\
Soit $M$  un \Amo \pfz,  $\rho:\gA\rightarrow \gB$  un \homo d'anneaux,
et  $\rho\ist(M)$ le $\gB$-module obtenu par \eds à $\gB$.
On a pour tout entier $n\geq 0$ l'\egt:
$\gen{\rho\big(\cF_{n}(M)\big)} =  \cF_{n}\big(\rho\ist(M)\big).$
\\
En particulier, si $S$ est un \moz, on a
$\cF_{n}(M_S) =  \big(\cF_{n}(M)\big)_S .$
\end{fact}
%--------------

Les deux faits suivants sont moins évidents.
%--- Lemma{fact.idf.ann}---
\begin{lemma} \emph{(Annulateur et premier \idfz)}\label{fact.idf.ann}
\\
Soit $M$  un \Amo \pfz, engendré par $q$ \eltsz, on~a:

\snic{\Ann(M)^q\subseteq \cF_{0}(M) \subseteq \Ann(M).}
\end{lemma}
%--------------
%-----------------begin proof------------------
\begin{proof}
Soit  $(\xq)$ un \sgr de $M$, $X={[\,x_1\;\cdots\;x_q\,]}$ et $G$ une
\mpn associée à $X$.
Soient $a_1$, \ldots, $a_q\in\Ann(M)$. Alors, la matrice diagonale
$\Diag(a_1,\ldots ,a_q)$ a pour colonnes des \colis des colonnes de $G$, donc son
\deter $a_1\cdots a_q$ appartient à~$\cF_{0}(M)$. Ceci prouve la première inclusion. \\
Soit $\delta$ un mineur d'ordre $q$ extrait de $G$. On va montrer que
$\delta \in\Ann(M)$, d'où la deuxième inclusion. Si $\delta$
correspond à une sous-matrice $H$ de $G$ on~a~$X\,H=0$, donc $\delta X=0$, et cela
signifie bien $\delta \in\Ann(M)$.
\end{proof}
%-----------------end proof------------------

%--- Fact{fact.idf.sex}---
\begin{fact}\label{fact.idf.sex} 
\emph{(Idéaux de Fitting et suites exactes)}\\
Soit $0\rightarrow N\rightarrow M\rightarrow P\rightarrow 0$  une suite exacte de
\mpfsz. Pour tout $p\geq 0$ on~a

\snic{\cF_p(M)\;\supseteq\;\sum_{r\geq0,s\geq0,r+s= p}\cF_r(N)\cF_{s}(P) ,}

%\sni
et si $M\simeq N \oplus P$, l'inclusion est une \egtz.
\end{fact}
%--------------
%-----------------begin proof------------------
\begin{proof}
On peut considérer que $N\subseteq M$ et $P=M/N$. On reprend les notations du point \emph{3} de la proposition \ref{propPfSex}. On a une \mpn $D$ de~$M$  qui
s'écrit \gui{sous forme triangulaire}
$$
D=\cmatrix{A&C\cr 0&B}.
$$
Alors, tout produit d'un mineur d'ordre $k$ de $A$ et d'un
mineur d'ordre $\ell$ de $B$ est égal à un mineur d'ordre $k+\ell$ de $D$.
Ceci implique le résultat annoncé pour les \idfsz. 
\\
Le deuxième cas est clair, avec $C=0$.
\end{proof}
%-----------------end proof------------------

%:h2014  ajout d'un exemple
\exl Sur l'anneau de \pols $\gA=\ZZ[a,b,c,d]$, considérons le module
$M=\gA g_1+ \gA g_2=\Coker F$ où $F=\cmatrix{a&b\cr c&d}$. Ici $g_1$ et $g_2$ sont les images de la base naturelle $(e_1,e_2)$ de $\Ae2$. Notons $\delta=\det(F)$. \\
On voit facilement que  $\delta\, e_1$ est une base du
sous-module  $\Im F\cap e_1\gA$
de $\Ae2$. 
\\
Soient $N=\gA g_1$ et $P=M/N$.
Alors le nodule $N$ admet la \mpn $[\,\delta\,]$ pour le \sgr $(g_1)$ et $P$ admet la \mpn $[\,c\;d\,]$ pour le \sgr $(\ov{g_2})$. 
\\
En conséquence, on obtient $\cF_0(M)=\cF_0(N)=\gen{\delta}$ et $\cF_0(P)=\gen{c,d}$, et l'inclusion
$\cF_0(N)\cF_0(P)\subseteq \cF_0(M)$ est stricte.   
\eoe

%:--- SUBsec {sec Fitt tf}--------
\subsec{Idéaux de Fitting d'un \mtfz}
%\addcontentsline{toc}{subsection}{Les \idfs d'un \mtfz}
\label{sec Fitt tf}

On peut généraliser
la \dfn des \idfs à un \mtf arbitraire $M$ comme suit.
Si $(\xq)$
est un \sgr de $M$ et si $X=\tra{[\,x_1\;\cdots\;x_q\,]}$,
on définit $\cF_{q-k}(M)$ comme
l'\id engendré par tous les mineurs d'ordre $k$ de toutes les matrices
$G\in\Ae {k\times q}$ vérifiant~$GX=0$.
Une \dfn alternative est que chaque $\cF_j(M)$ est la somme de tous les $\cF_j(N)$
où $N$ parcourt les
\mpfs qui s'envoient surjectivement sur $M$.

Ceci montre que les \ids ainsi définis ne dépendent pas du \sgr considéré.

\ss La remarque suivante est souvent utile.

%--- Fact{facttfpf}--------------
\begin{fact}
\label{facttfpf}
Soit $M$ un \Amo \tfz.
%-----------------begin enum------------------
\begin{enumerate}
\item Si $\cF_k(M)$ est un \itfz, $M$ est le quotient d'un \mpf $M'$
pour lequel  $\cF_k(M')=\cF_k(M)$.
\item Si tous les \idfs sont \tfz,  $M$ est le quotient d'un \mpf $M'$
ayant les mêmes \idfs que~$M$.
\end{enumerate}
%-----------------end enum------------------
\end{fact}
%--- end-fact-----------------------------------------

% --- section  Idéal résultant {subsecIdealResultant}-----
\section{Idéal résultant}
\label{subsecIdealResultant}
%-----------------------------------------

Dans ce qui suit, on considère un anneau $\gk$ que l'on ne suppose pas discret.
%on utilise donc des \emph{\pols formels} \cad des couples $(f,p)$ où $f$ 
% est un \pol qui s'écrit sous forme $\sum_{k=0}^pa_kX^k$.
%\medskip 
Le résultant de deux \pols est à la base de la théorie de l'\eliz.
Si $f$, $g\in\kX$ avec $f$ unitaire,
le lemme d'\eli de base \paref{LemElimAffBasic} 
peut être lu dans l'\alg $\gB=%\gk[x]=
\aqo{\kX}{f}$ en écrivant:
$$ \rD_\gB(\ov g)\,\cap\, \gk=\rD_\gk\big(\Res_X(f,g)\big).$$

Il se \gns alors avec le résultat suivant, 
 que l'on peut voir comme une formulation très précise du 
 lemme \gui{lying over} (voir le lemme~\ref{lemLingOver}).

%--- Lemme{LemElimAff}--------
\CMnewtheorem{lemeligen}{Lemme d'\eli \gnlz}{\itshape}
%:     Lemma{lem}
\begin{lemeligen}\label{LemElimAff}%
\index{Lemme d'elimination general@Lemme d'\eli général}%
\index{elimination@\eliz!lemme d'--- général}%
\index{elimination@\eliz!ideal@\id d'---}~
\begin{enumerate}
\item Soient $\gk\vers{\rho}\gC$ une \alg qui est un \kmo engendré par
$m$ \eltsz, $\fa=\cF_{\gk,0}(\gC)$ son premier \idf et
$\fc=\Ker\rho$.
Alors:
\begin{enumerate}
\item $\fc=\Ann_\gk(\gC)$,
\item \fbox{$\fc^m\subseteq\fa\subseteq\fc$} et donc \fbox{$\rD_\gk(\fc)=\rD_\gk(\fa)$},
\item si par une \eds $\varphi:\gk\to\gk'$ on obtient l'\alg $\rho':\gk'\to\gC'$,
alors l'\id $\fa':=\cF_0(\gC')$ est égal à $\varphi(\fa)\gk'$ et en tant que $\gk'$-module, il est isomorphe à $\gk'\otimes_\gk\fa\simeq\varphi\ist(\fa)$.
\end{enumerate}
\item Soit $\gB\supseteq\gk$ une \klg qui est un \kmo libre de rang
$m$, et $\fb$  un \itf de~$\gB$.
\begin{enumerate}
\item L'\id d'\eli $\fb\,\cap\,\gk$ est le noyau de l'\homo canonique
$\rho:\gk\to\gB\sur\fb$, i.e. l'annulateur du \kmoz~$\gB\sur\fb $.

\item Le \kmo $\gB\sur\fb $ est \pf et l'on~a:

{\fbox{$(\fb\,\cap\,\gk)^m\subseteq\cF_0(\gB\sur\fb )\subseteq\fb\,\cap\,\gk$}  
\hbox{ et }  
\fbox{$\rD_\gB(\fb)\,\cap\, \gk=\rD_\gk\big(\cF_0(\gB\sur\fb )\big)$}.}

On  note $\fRes(\fb):=\cF_{\gk,0}(\gB\sur\fb )$, on l'appelle l'\emph{\id résultant de $\fb$}.
\end{enumerate}
\end{enumerate}
\end{lemeligen}
%--------- fin lemma ---------------------------------------------- 

%----------------- proof ------------------
\begin{proof}
\emph{1a} et \emph{1b.} En effet, $a\in\gk$ annule  $\gC$ \ssi
il annule $1_{\gC}$, \ssi $\rho(a)=0$. 
La double inclusion recherchée est donc donnée 
par le lemme  \ref{fact.idf.ann} (valable aussi pour les \mtfsz).

 \emph{1c.} Les \idfs se comportent bien par \edsz.

 \emph{2.} On applique le point \emph{1} avec $\gC=\gB\sur\fb $.
\end{proof}
%-----------------end proof------------------

\rems  
 1) L'\id résultant dans le point \emph{2} peut être décrit \prmt
 comme suit. Si $\fb=\gen{b_1,\ldots ,b_s}$ on considère l'\emph{application de Sylvester \gneez}
$$\preskip-.2em \postskip.3em \ndsp
\psi:\gB^s\to\gB,\quad (y_1,\ldots ,y_s)\mapsto\psi(\uy)=\som_iy_ib_i.
$$
%:2018 \sni remplacé par $$ ci-dessus  
C'est une \kli  entre \kmos libres de rangs  $ms$ et $m$.
Alors, on a $\fRes(\fb)=\cD_m(\psi)$. \rdb

2) Cela fait beaucoup de \gtrs pour l'\id $\fRes(\fb)$.
En fait il existe diverses techniques pour diminuer le nombre de \gtrs
en remplaçant~$\fRes(\fb)$ par un \itf 
ayant nettement moins de \gtrs mais ayant le
même nilradical.
Voir à ce sujet: le traitement donné en section~\ref{secChap3Nst}
avec notamment le lemme~\ref{lemElimParametre},  
les résultats du chapitre~\ref{chapKrulldim} sur le nombre de \gtrs radicaux d'un \id radicalement \tf (\thref{thKroH}),
et l'article~\cite{DiGLQ}.
\eoe%
\index{application de Sylvester!généralisée}%
\index{Sylvester!application de --- généralisée}%

\smallskip 
Voici maintenant un cas particulier du lemme d'\eli \gnlz.
Ce \tho complète le lemme \ref{lemElimPlusieurs}. 

%:2018  theorem  transformé en thdef
%:   Theorem{thElimAff}--------------
\begin{thdef}
\label{thElimAff} \label{corLemmeElim} \emph{(\Tho d'\eli \agqz: l'\id résultant)}%
\index{ideal@idéal!résultant}%
\index{resultant@résultant!idéal ---}%
\index{elimination@\eliz!theoreme@\tho d'--- \agqz} 
Soit $(f,g_1,\ldots,g_r)$ des \pols de $\kX$ avec $f$ \mon de degré $m$.
On pose 
$$\preskip.2em \postskip.4em
\ff=\gen{f,g_1,\ldots,g_r}\subseteq\kX \;\hbox{et}\;\gB=\aqo{\kX}{f}.
$$
Notons $\psi:\gB^r\to\gB$ l'\emph{application de Sylvester \gneez} définie par:
$$\preskip.4em \postskip.4em\ndsp
(y_1,\ldots ,y_r)\mapsto \psi(\uy)=\som_iy_i\ov{g_i}.
$$
Il s'agit d'une \kli entre \kmos libres de rangs respectifs $mr$ et $m$.
Notons $\fa$ l'\idd $\cD_m(\psi)$. 
\begin{enumerate}
\item  $\fa=\cF_{\gk,0}(\kX\sur\ff),$ 
 et l'on a
$$\preskip.3em \postskip.3em
(\ff\,\cap\,\gk)^m\subseteq\fa\subseteq\ff\,\cap\,\gk,\quad \mathit{et\;donc}\quad
\rD_\kX(\ff)\,\cap\, \gk=\rD_\gk(\fa).
$$
\item Supposons que $\gk=\gA[Y_1,\ldots,Y_q]$ et que $f$ et les $g_i$ soient
de degré total~$\leq d$ dans $\gA[\uY,X]$.  Alors, les \gtrs de $\cD_m(\psi)$
sont de degré total~$\leq d^2$ dans $\gA[\uY]$. 
\item L'\id $\fa$ ne dépend que de 
$\ff$ (sous la seule hypothèse que $\ff$ contienne un \poluz).  
Nous l'appelons \emph{l'\id résultant de $\ff$
par rapport à l'\idtr $X$} et nous le notons $\fRes_X(f,g_1,\ldots,g_r)$
ou $\fRes_{X}(\ff)$, ou~$\fRes(\ff)$. 
\item Si par une \eds $\theta:\gk\to\gk'$ on obtient l'\id $\ff'$ de~$\gk'[X]$,
alors l'\id $\fRes_{X}(\ff')\subseteq\gk'$ est égal à $\theta\big(\fRes_{X}(\ff)\big)\gk'$,
et en tant que module il est isomorphe à $\gk'\otimes_\gk\fRes_{X}(\ff)\simeq\theta\ist\big(\fRes_{X}(\ff)\big)$.

\end{enumerate}
\end{thdef}
%--------- fin corollary ---------------------------------------------- 
NB. Considérons la base $\cE=(1,\ldots ,X^{m-1})$  de $\gB$ sur $\gk$.
Soit  $F\in\gk^ {m\times mr}$ la matrice de  $\psi$ 
pour les bases déduites de  $\cE$. 
Ses colonnes sont les $X^jg_k \mod f$ pour $j\in\lrb{0..m-1}$, 
$k\in\lrbr$ écrits sur la base $\cE$.
On dit que $F$ est une \emph{matrice de Sylvester \gneez}. Par \dfn on a $\fRes_{X}(\ff)=\cD_m(F)$.
\index{matrice!de Sylvester \gneez}
\eoe

\begin{proof} Posons $\fb=\ff \mod f=\gen{\ov{g_1},\ldots,\ov{g_r}}\subseteq\gB$.
On applique les points \emph{2} et~\emph{1c} du 
lemme d'\eli \gnl en remarquant que
$\kX\sur\ff\simeq\gB\sur\fb$, avec~$\ff\cap\gk=\fb\cap\gk$.
\end{proof}

\rem 
Ainsi, le \tho \ref{thElimAff} établit un lien très étroit entre \id d'\eli
et \id résultant. 
Les avantages que présente l'\id résultant sur l'\id
d'\eli sont les suivants\index{ideal@idéal!d'elimi@d'\eliz}
\begin{itemize}
\item l'\id résultant est \tfz,
\item  son  calcul est \emph{uniforme},
\item il se comporte bien par \edsz. 
\end{itemize}
Notons que dans le cas où $\gk=\gK[Y_1,\ldots,Y_q]$ avec 
$\gK$ un \cdiz, l'\id d'\eli est aussi \tf mais son calcul, 
par exemple via les \bdgsz, n'est pas uniforme.
\\
Cependant l'\id résultant n'est défini que lorsque $\ff$ contient un \polu
et ceci limite la portée du \thoz.
\eoe

%:section: --- Exercices
\Exercices{

%--- Exercise{exo3Lecteur}-------------
\begin{exercise}
\label{exo3Lecteur}
{\rm  Il est recommandé de faire les \dems non données, esquissées,
laissées \alecz,
etc\ldots
\, On pourra notamment traiter les cas suivants.
\begin{itemize}
\item \label{exo lem pres equiv}
Donner une preuve détaillée du lemme \ref{lem pres equiv}.
\item \label{exopropCoh1.2}
Expliquez pourquoi les propositions \ref{propCoh1} et \ref{propCohfd1}
(lors\-que l'on  prend $\gA$  comme \Amo $M$)
se relisent sous la forme du \thrf{propCoh2}.

\item \label{exopropCohpfKer}
  Démontrer les propositions \ref{propPfInter} et \ref{propCohpfKer}.
Donner une preuve  détaillée des propositions
 \ref{propPftens} et~\ref{propPfExt}. Montrer que
$\gA\sur{\fa}\te_\gA \gA\sur{\fb}\simeq \gA\sur{(\fa+\fb)}.$

\item \label{exobelexemple}
Justifier les affirmations contenues dans l'exemple \paref{belexemple}.

\item \label{exoZerRed}
Démontrer les lemmes ou faits  \ref{lemZerloc}, \ref{lemZerRed},
\ref{factZerRedCoh} et
\ref{factZerRedConnexe}.

\item \label{exothZerDimRedLib}
Donnez des \algos pour les trois points du \thrf{thZerDimRedLib}.

\item \label{exoFittsex}  \label{exoFitt5}
Prouver le fait \ref{facttfpf}.
\end{itemize}
}
\end{exercise}
%--- end -exercise-----------------------------------------

%--- Exercise{exoptfpf}--------
\begin{exercise}
\label{exoptfpf}
{\rm  Soient $M\subseteq N$ des \Amos avec $M$ en facteur direct dans $N$.
Si~$N$ est \tf (resp. \pfz), alors $M$ \egmtz.
 }
\end{exercise}
%--- end-exercise-----------------------------------------

%:HHH exo rajouté

%:--- Exercise{exoAXmodule}-------------
\begin{exercise}\label{exoAXmodule} {(Structure de $\AX$-module sur $\Ae n$ associée à
un $A\in\Mn(\gA)$)}\\ 
{\rm  
Soit $\gA$ un anneau commuatif. On munit $\Ae n$ d'une structure de $\AX$-module en
posant 
$$\preskip.0em \postskip.4em 
Q \cdot x = Q(A) \cdot x\hbox{  pour  }Q \in \AX\hbox{ et }x \in \Ae n. 
$$
On va donner une \mpn pour ce $\AX$-module. Ceci généralise l'exemple 3)
\paref{exl1pf} donné au début de la section \ref{sec pf chg}, avec $\gA$ un \cdiz. 

%\sni
Soit $\theta_A : \AX^n \twoheadrightarrow \Ae n$ l'unique $\AX$-morphisme qui
transforme la base canonique de $\AX^n$ en celle de $\Ae n$. En notant du 
même nom $(e_1, \ldots, e_n)$ ces deux bases canoniques, $\theta_A$
transforme donc $Q_1 e_1 + \cdots + Q_n e_n$ en $Q_1(A)\cdot e_1 + 
\cdots + Q_n(A)\cdot e_n$.  On va montrer que la suite ci-dessous est exacte:

\snic {
\AX^n\vvvvvers{X\In-A} {\AX^n} \vvers{\theta_A} \Ae n \to 0
}

%\sni
Autrement dit $\Ae n$ est un $\AX$-\mpf et $X\In-A$ est une \mpn pour le
\sgr $(e_1, \ldots, e_n)$.

%\sni
\emph {1.}
Montrer que l'on a une somme directe de $\gA$-modules $\AX^n = \Im(X\In-A)
\oplus \Ae n$.

%\sni
\emph {2.}
Conclure.

}
\end{exercise}
%--- end-exercise-----------------------------------------

%--- Exercise{exoTenseursNuls}-------
\begin{exercise}\label{exoTenseursNuls}
 (Description des tenseurs nuls)\\
{\rm  
Soient $M$ et $N$ deux \Amos arbitraires,  et $z=\sum_{i\in \lrbn}x_i\te y_i\in M\te N$.

\emph{1.}
 Montrer que $z=0$ \ssi il existe
un sous-\mtfz~$M_1$  de $M$ 
tel que l'on ait $\sum_{i\in\lrbn}x_i\te y_i=_{M_1\te N}0$.

\emph{2.} On écrit $M_1=\gA x_1+\cdots+\gA x_p$ avec $p\geq n$.
On pose $y_k=_N0$ pour $n<k\leq p$. Utiliser le lemme du tenseur nul
avec l'\egt $\sum_{i\in\lrbp}x_i\te y_i=_{M_1\te N}0$ pour donner une \carn
des tenseurs nuls dans la situation \gnlez.
}
\end{exercise}
%--- end-exercise-----------------------------------------

%--- Exercise{exoEdsQuot}-------------
\begin{exercise}
\label{exoEdsQuot}
{\rm  Soit $M$ un \Amoz, $\fa$ un \id et $S$ un \mo de $\gA$.

\emph{1.} Montrer que l'\ali canonique $M\to M\sur{\fa M}$ résout le \pb \uvl de l'\eds
pour l'\homo $\gA\to\gA\sur\fa$ (i.e., selon la \dfn \ref{defAliAliExtScal}, cette \ali est un morphisme d'\eds de $\gA$
à $\gA\sur\fa$ pour $M$). 
En déduire que l'\ali natu-\linebreak relle $\gA\sur\fa\otimes _\gA M\to M\sur{\fa M}$ est un \isoz.

\emph{2.} Montrer que l'\ali canonique $M\to M_S$ résout le \pb \uvl de l'\eds
pour l'\homo $\gA\to\gA_S$. En déduire que l'\ali naturelle $\gA_S\otimes _\gA M\to M_S$ est un \isoz.
 
}
\end{exercise}
%--- end -exercise-----------------------------------------

%--- Exercise{exoBézoutstrict}-------
\begin{exercise}
\label{exoBézoutstrict}
{\rm
Toute matrice sur un anneau de Bézout intègre est
\eqve à une matrice de la forme $\cmatrix{T&0\cr0&0}$
avec $T$ triangulaire et les \elts sur la diagonale de $T$ non nuls
(naturellement, les lignes ou colonnes indiquées nulles peuvent être
absentes).
Cette \eqvc peut être obtenue par des manipulations de Bézout.\\
Généraliser aux anneaux \qis en utilisant la méthode générale
expliquée \paref{MethodeQI}.  
}
\end{exercise}
%--- end-exercise-----------------------------------------

%--- Exercise{exoAnneauBézoutStrict}-------------
\begin{exercise}\label{exoAnneauBézoutStrict}
{(Anneaux de Bézout strict)}\\
{\rm 
\emph {1.}
Pour un anneau $\gA$, montrer que \propeq
\vspace{-2pt}
\begin{itemize}%\itemsep=0pt
\item [\emph {a.}]
Si $A \in \Ae{n \times m}$, il existe $Q \in \GL_m(\gA)$
telle que $AQ$ soit triangulaire inférieure.

\item [\emph {b.}]
Même chose qu'en \emph {a} avec $(n,m) = (1,2)$, i.e. $\gA$
est un anneau de Bézout strict.

\item [\emph {c.}]
Pour $a$, $b \in \gA$, il existe $x$, $y \in \gA$ \com tels que $ax + by = 0$.

\item [\emph {d.}]
Pour $(\ua) = (a_1, \ldots, a_n)$ dans $\gA$, il existe $d\in \gA$ et
$(\underline {a'}) = (a'_1, \ldots, a'_n)$ \com vérifiant
$(\ua) = d(\underline {a'})$; on a alors $\gen {\ua} = \gen {d}$.

\item [\emph {e.}]
Même chose qu'en \emph {d} avec $n=2$.
\end{itemize}

\emph{2.}
Montrer que la classe des anneaux de Bézout strict est stable
par produit fini, par quotient et par localisation.

Dans la suite, on suppose que $\gA$ est un anneau de Bézout strict.

\emph {3.}
Soient $a$, $b$, $d_2 \in \gA$ tels que $\gen {a,b} = \gen {d_2}$. Montrer qu'il
existe $a_2$, $b_2 \in \gA$ \com tels que $(a,b) = d_2(a_2,b_2)$. On pourra
considérer $d_1$, $a_1$, $b_1$, $u_1$, $v_1$ avec $(a,b) = d_1(a_1,b_1)$, $1 =
u_1a_1 + v_1b_1$ et introduire:
$$\preskip.4em \postskip.4em 
(\star) \quad
\Cmatrix{.2em} {a_2\cr b_2} = \crmatrix {v_1 & a_1\cr -u_1 & b_1}
\cmatrix {\vep \cr k_{12}} \; \hbox {où} \;
d_1 = k_{12}d_2,\ d_2 = k_{21}d_1,\ \vep = k_{12}k_{21} - 1. 
$$

\emph {4.}
Même chose que dans le point précédent mais avec un nombre quelconque
d'\elts \cad pour $(\ua) = (a_1, \ldots, a_n)$ dans $\gA$ et $d$ donnés vérifiant
$\gen {\ua} = \gen {d}$, il existe $(\underline {a'}) = (a'_1, \ldots, a'_n)$, \comz,
tels que $\ua = d\underline {a'}$.

\emph {5.}
Montrer que toute matrice  {diagonale} $\Diag(\an)$ est $\SL_n$-\eqve
à une autre matrice diagonale $\Diag(\bn)$ avec $b_1 \divi b_2 \divi \cdots \divi b_n$.\\
 De plus, si l'on pose $\fa_i = \gen {a_i}$, $\fb_i = \gen {b_i}$, on a
$\fb_i = S_i(\fa_1, \ldots, \fa_n)$ où $S_i$ est la \gui {$i$-ème fonction
\smq \elr de $\fa_1, \ldots, \fa_n$} obtenue en remplaçant le produit par
l'intersection. Par exemple:

\snic {
S_2(\fa_1, \fa_2, \fa_3) = (\fa_1\cap\fa_2) + (\fa_1\cap\fa_3) + (\fa_2\cap\fa_3).
}

%\sni
En particulier, $\fb_1 = \sum_i \fa_i$, $\fb_n = \bigcap_i \fa_i$.
De plus $\prod_i \gA\sur{\fa_i} \simeq \prod_i \gA\sur{\fb_i}$.
\\
Ce dernier résultat sera \gne aux \anars (corolaire
\ref{corthAnar}). 
\\
D'autres \gui {vraies} fonctions \smqs \elrs d'\ids 
interviennent dans l'exercice~\ref{exoFitt0}.

}
\end{exercise}
%--- end -exercise-----------------------------------------

%:HHH petits changements dans l'exo suivant
%--- Exercise{exoSmith}-------------
\begin{exercise}\label{exoSmith}\index{Smith!anneau de ---}
{(Anneaux de Smith, ou elementary divisor rings)}\\
{\rm Définissons un \ixx{anneau}{de Smith}
comme un anneau sur lequel toute matrice admet une forme réduite de Smith,
(cf. la section \ref{secBézout} \paref{secpfval}).
Un tel anneau est de Bézout strict
(cf. 
l'exercice~\ref{exoAnneauBézoutStrict}).
Puisque sur un anneau de Bézout strict, toute matrice
carrée diagonale est \eqve à une matrice de Smith
(exercice~\ref{exoAnneauBézoutStrict}, question~\emph {5}),
un anneau est de Smith \ssi toute matrice est \eqve
à une matrice \gui{diagonale}, sans condition de
\dve sur les \coesz.
Ces anneaux ont été en particulier étudiés
par Kaplansky dans \cite{Kap}, y compris pour le cas non commutatif,
puis par Gillman \& Henriksen dans \cite {GillmanHenriksen}.
Nous nous limitons ici au cas commutatif.
\\
Montrer que \propeq
\begin{itemize}
  \item [\emph{1.}] $\gA$ est un anneau de Smith.
  \item [\emph{2.}] $\gA$ est un anneau  de Bézout strict et toute matrice
  triangulaire dans $\MM_2(\gA)$
est \eqve à une matrice diagonale.
  \item [\emph{3.}] L'anneau $\gA$ est de Bézout strict, et si $1 \in \gen{a,b,c}$, 
  alors il existe $(p,q)$, $(p',q')$  tels que $1=pp'a+qp'b+qq'c$.
  \item [\emph{4.}] L'anneau $\gA$ est de Bézout strict, et si $\gen{a,b,c} = \gen{g}$, 
  alors il existe $(p,q)$, $(p',q')$  tels que $g=pp'a+qp'b+qq'c$.
\end{itemize}
Ceci donne un joli \tho de structure pour les
\mpfsz, en tenant compte pour l'unicité du \thrf{prop unicyc}.
Notez aussi que ce \tho implique l'unicité de la réduite de Smith
d'une matrice $A$ (en considérant le module conoyau) au sens suivant: 
en notant $b_i$ les \coes diagonaux de la réduite, les \idps $\gen{b_1}\supseteq\cdots\supseteq\gen{b_q}$  ($q=\inf(m,n)$) sont des invariants de la {matrice} $A$ à \eqvc près. 
\\
En termes de modules, ces \idps caractérisent, 
à un \auto près de $\Ae{m}$, le morphisme d'inclusion
$P=\Im(A)\to \Ae{m}$. 
\\
Une base~\hbox{$(e_1,\dots,e_m)$} de $\Ae{m}$ telle que   $P=b_1\gA\, e_1+\cdots + b_m\gA\, e_m$
 est appelée une \emph{base de $\Ae{m}$ adaptée au sous-module $P$.}%
\index{base adaptée!à une inclusion} 
\\
Posons $b_r=0$ si $m\geq r>n$, on a $\gen{b_1}\supseteq\cdots\supseteq\gen{b_r}$.
Les \idpsz~\hbox{$\neq \gen{1}$} 
de cette liste sont  les
{facteurs invariants} du module $M=\Coker(A)$.
Le \thrf{prop unicyc} nous dit que cette liste caractérise la structure
du module~$M$.   
\\
Notons enfin que les anneaux de Smith sont stables par produit fini, \lon et
passage au quotient.
}

\end {exercise}
%--- end -exercise-----------------------------------------

%--- Exercise{exoUnitsOfSomeRings}-------------
\begin{exercise}\label{exoUnitsOfSomeRings}
{(Exemple \elr de détermination du groupe des \ivsz)}

{\rm 
\emph {1.}
Soit $\gk$ un anneau réduit et $\gA = \aqo{\gk[Y,Z]}{YZ}
=\gk[y,z]$ avec $yz = 0$. Montrer en utilisant une mise en position de
\iNoe de $\gA$ sur $\gk$, que $\Ati = \gk\eti$.

\emph {2.}
Soit $\gA = \aqo{\ZZ[a,b,X,Y]}{X-aY,Y-bX} = \ZZ[\alpha,\beta,x,y]$
avec $x = \alpha y$ et $y = \beta x$.
Montrer que $\Ati = \{\pm 1\}$;
on a donc $\gA x = \gA y$ mais $y \notin \Ati x$.

}

\end {exercise}
%--- end -exercise-----------------------------------------

%--- Exercise{exoUAtoUB}-------------
\begin{exercise}\label{exoUAtoUB}  
{(Conditions suffisantes pour la surjectivité de $\Ati \to (\gA/\fa)\eti$)}
\\
{\rm 
%:HHH 11 reference a un autre exo
Voir aussi l'exercice \ref{exoLgb2}.\\  
Pour un \id $\fa$ d'un anneau $\gA$, on considère la \prt $(\star)$:

\snic {
(\star)\qquad\qquad\qquad
\Ati \to (\gA/\fa)\eti \hbox { est surjectif}
,}

%\sni
i.e., pour $x \in \gA$ \iv modulo $\fa$, il existe
$y \in \Ati$ tel que $y \equiv x \bmod \fa$, ou encore: si 
$\gA x + \fa$ rencontre $\Ati$, alors  $x + \fa$ rencontre $\Ati$.

\emph{1.}
 Montrer que $(\star)$ est vérifiée quand $\gA$ est
\zedz.

\emph{2.}
 Si $(\star)$ est vérifiée pour tout \id principal
$\fa$, elle l'est pour tout \idz~$\fa$.

\emph{3.}
 On suppose $(\star)$ vérifiée. Soient $x$, $y$ deux \elts
d'un \Amo tels \linebreak que $\gA x = \gA y$; montrer que $y = ux$
pour un  $u \in \Ati$.\\
NB: l'exercice \ref{exoUnitsOfSomeRings}
fournit un exemple d'anneau $\gA$ avec $x$, $y \in \gA$
et $\gA x = \gA y$, \linebreak mais $y \notin \Ati x$.

\emph{4.}
 Soient $\gA' = \gA\sur{\Rad \gA}$, $\pi : \gA \twoheadrightarrow \gA'$
la surjection canonique et $\fa' = \pi(\fa)$. Montrer que si $(\star)$ est
vérifiée pour $(\gA',\fa')$, elle l'est pour $(\gA, \fa)$.
}

\end {exercise}
%--- end -exercise-----------------------------------------

%--- Exercise{exoCalculT(M)}-------------
\begin{exercise}
\label{exoCalculT(M)} (Calcul d'un sous-module de torsion)
\\
 {\rm  
Soit  $\gA$ un \cori intègre et $M$ un \Amo \pfz. 
 Alors le sous-module
de torsion de $M$ est un \mpfz.
\\
%:2012 il manquait ``une'' avant \mpn
 Plus \prmt si l'on a une \mpn $E$ pour $M$ avec une suite exacte 
$$\preskip-.2em \postskip.4em 
\Ae n \vvers{E} \Ae\ell \vvers\pi M\to 0 
$$
et si $F$ est une matrice telle que l'on ait une suite exacte
$$\preskip.2em \postskip.4em 
\Ae m \vvers{F} \Ae\ell \vvers {\tra E\,} \Ae n 
$$
(l'existence de la matrice $F$ résulte du fait que $\gA$
est \cohz) alors le sous-module de torsion $\rT(M)$ de $M$ est 
égal à $\pi(\Ker{\tra F})$ et isomorphe
à $\Ker{\tra F}/\Im E$.

 Montrer aussi que le résultat se généralise au cas où
$\gA$ est \coh et \qiz. 
}
\end{exercise}
%--- end -exercise-----------------------------------------

%--- Exercise{exoZerRedBez}-------------
\begin{exercise}
\label{exoZerRedBez} (L'\algo d'Euclide dans le cas \zedrz)\\
 {\rm On donne ici une version plus uniforme de la \dem de
la proposition \ref{propZedBez} et on la généralise. On considère un
anneau \zedr $\gA$.

\emph{1.} 
Soient $\gB$ un anneau quelconque et $b\in\gB$ tel que $\gen{b}$ soit engendré
par un \idmz.  Pour  $a\in\gB$ donner une matrice $M\in\EE_2(\gB)$ et
$d\in \gB$ vérifiant l'\egt $M\cmatrix{a\cr b}=\cmatrix{d\cr 0}$. En particulier
$\gen{a,b}=\gen{d}$.

\emph{2.}  Donner un \algo d'Euclide \gui{uniforme} pour deux \pols
de $\gA[X]$.

\emph{3.}  
L'anneau $\AX$ est un anneau de Smith:   donner un \algo qui
réduit toute matrice sur $\AX$ à une forme de Smith  au moyen de manipulations \elrs de lignes et de colonnes.  
}

\end{exercise}
%--- end-exercise-----------------------------------------

%--- Exercise{exoZDpiv}--------------
\begin{exercise}
\label{exoZDpiv} (Dépendance linéaire en dimension $0$)\\
%:HHH 11 petit commentaire en plus enonce restructure
{\rm On donne ici la \gnn du \tho selon lequel $n+1$ vecteurs de $\gK^n$ sont \lint dépendants, du cas des \cdis à celui des anneaux \zedrsz.
Notez que la syzygie, pour être digne de ce nom, doit avoir des
\coes \comz. 

Soit $\gK$ un anneau \zedrz, et  $y_1$, \ldots, $y_{n+1}\in\gK^{n}$.

\emph{1}. Construire un \sfio $(e_j)_{j\in \lrb{1..n+1}}$ de façon à ce que, dans
chaque composante $\gK[1/e_j]$, le vecteur $y_j$ soit \coli des $y_i$ qui le précèdent.

\emph{2}. En déduire qu'il existe un \sys d'\elts \com $(a_1, \ldots, a_{n+1})$ dans~$\gK$ 
tel que $\sum_i a_iy_i = 0$.

\sni
\rems 1) On rappelle la convention selon laquelle on accepte que certains \elts d'un
\sfio soient nuls: on voit sur cet exemple que l'énoncé de la \prt
désirée en est grandement facilité.

 2)
On pourra au choix, ou bien donner un traitement adéquat de la matrice des 
$y_i$ par des
manipulations \elrs en s'appuyant sur le lemme \ref{lemQI}, ou bien traiter le cas
des \cdis puis utiliser la machinerie locale globale \elr \num2
\paref{MethodeZedRed}.\imlgz\eoe
}
\end{exercise}
%--- end-exercise-----------------------------------------

%--- Exercise exoZedLG ----------
\begin{exercise}
\label{exoZedLG}
{\rm Soient $S_1$, $\dots$, $S_n$ des \moco de $\gA$. Montrer que $\gA$ est \zed \ssi
chacun des $\gA_{S_i}$ est \zedz.
}
\end{exercise}
%--- end-exercise-----------------------------------------

%--- Exercise{exoFreeAlgebraPresentation}-------------
\begin{exercise}\label{exoFreeAlgebraPresentation}
(Présentation d'une \alg qui est libre finie comme module)
\\
{\rm
Soit~$\gB$ une \Alg libre de rang $n$ de base
$\ue = (e_1, \ldots, e_n)$. On note

\snic{\varphi : \AuX = \gA[X_1, \ldots, X_n]\twoheadrightarrow \gB}

%\sni
l'\homo (surjectif) d'\Algs qui réalise
$X_i \mapsto e_i$. On note $c_{ij}^k$ les constantes de structure
définies par $e_ie_j = \sum_k c_{ij}^k e_k$.
On considère $a_1$, \ldots, $a_n \in \gA$ définis
par $1 = \sum_k a_k e_k$ et l'on pose:

\snic{
R_0 = 1 - \som_k a_k X_k, \qquad
R_{ij} = X_iX_j - \sum c_{ij}^k X_k.}

%\sni
On note $\fa = \gen {R_0, R_{ij}, i \leq j}$.
Montrer que  tout $f \in \AuX$ est congru modulo $\fa$ à un
\pol \hmg de degré 1. En déduire
que $\Ker\varphi = \fa$.
}
\end{exercise}
%--- end -exercise-----------------------------------------

%--- Exercise{exoFitt0}-------
\begin{exercise}\label{exoFitt0}
(Quelques calculs d'\idfsz)\\
{\rm
\emph{1.} 
Déterminer les \idfs d'un \Amo présenté par une matrice en forme de
Smith.

\emph{2.}   
Déterminer les \idfs de $\gA\sur\fa$.

\emph{3.}   
Soit $E$ un \Amo \tf et $\fa$ un \idz. Montrer que

\snic{\cF_k(E\oplus \gA\sur{\fa})= \cF_{k-1}(E)+\cF_{k}(E)\,\fa.}

\emph{4.}   
Déterminer les \idfs du \Amo  $M=\gA\sur{\fa_1}\oplus\cdots\oplus
\gA\sur{\fa_n}$ dans le cas \hbox{où $\fa_1\subseteq\fa_2\subseteq \cdots\subseteq\fa_n$}. 

\emph{5.} 
 Déterminer les \idfs du \Amo
$M=\gA\sur{\fa_1}\oplus\cdots\oplus$ $\gA\sur{\fa_n}$
sans faire d'hypothèse d'inclusion pour les \ids $\fa_k$.
\\
Comparer $\cF_0(M)$ et $\Ann(M)$.
}
\end{exercise}
%--- end-exercise-----------------------------------------

%--- Exercise{exoFitt4}-------
\begin{exercise}\label{exoFitt4}
(Les \idfs d'un \Amo \tfz)\\
{\rm
Montrer que les faits \ref{fact.idf inc}, \ref{fact.idf.change},
 \ref{fact.idf.sex} et le lemme \ref{fact.idf.ann} restent valables avec les
\mtfsz.
}
\end{exercise}
%--- end-exercise-----------------------------------------

%--- Exercise{exoFitt6}-------
\begin{exercise}
\label{exoFitt6}
{\rm  Une des \prts \caras des \emph{\adpsz} (qui seront étudiés au chapitre 
\ref{ChapAdpc}) est la
suivante: si $A\in\Ae {n\times m}$, $B\in\Ae {n\times 1}$, et si les \idds
de $A$ et $[\,A\,|\,B\,]$ sont les mêmes, 
alors le \sliz~$AX=B$ admet une solution.

\emph{1.} Soit  un \mtf $M$ sur un \adp et $N$ un quotient de~$M$. Montrer que
si~$M$ et $N$ ont les mêmes \idfsz, alors $M=N$.

\emph{2.} Montrer que si un \mtf $M$ sur un \adp a ses \idfs qui sont \tfz, c'est un
\mpfz.
}
\end{exercise}
%--- end-exercise-----------------------------------------

%--- Exercise{exoAutresIdF}----------
\begin{exercise}\label{exoAutresIdF}\index{ideal@idéal!de Kaplansky}\index{Kaplansky!idéal de}
(Idéaux de Kaplansky) \\
{\rm Pour un \Amo $M$ et un entier $r$  on note $\cK_r(M)$
l'\id somme de tous les
transporteurs $\big(\gen{m_1,\ldots m_r}:M\big)$ pour tous les
\syss $(m_1,\ldots m_r)$ dans  $M$.
On l'appelle \emph{l'idéal de Kaplansky d'ordre $r$  du module $M$}.
Ainsi,  $\cK_0(M)=\Ann (M)$, et si $M$ est engendré par $q$ \eltsz,
 on a $\cK_q(M)=\gen{1}$.
%-----------------begin item------------------
\begin{itemize}
\item Montrer que si  $\cK_q(M)=\gen{1}$, $M$ est \tfz.
\item Montrer que si $M$ est \tfz, alors pour tout entier $r$ on a les inclusions

\snic{\cF_r(M) \subseteq  \cK_r(M) \subseteq \sqrt{\cF_r(M)}= \sqrt{\cK_r(M)}.}
\end{itemize}
%-----------------end item------------------
NB: voir aussi 
%:HHH ref changée en dessous
l'exercice~\ref{exoVariationLocGenerated}.
}
\end{exercise}
%--- end-exercise-----------------------------------------

%--- Exercise{exoPetitsExemplesElim}-------------
\begin{exercise}\label{exoPetitsExemplesElim} 
{(Un exemple \elr d'\ids résultants)} \\
{\rm  
Soient $f$, $g_1$, \ldots, $g_r \in \AX$, $f$ \mon de degré $d \ge 1$ et
$\ff = \gen {f, g_1, \ldots, g_r} \subseteq \AX$.  On va comparer  l'\id 

\snic{\fa = \fR(f, g_1,\ldots,
g_r) = \rc_T\big(\Res(f, g_1 + g_2T + \cdots + g_rT^{r-1})\big)}

%\sni
(section \ref {secChap3Nst}),  et l'\id résultant 
$\fb = \fRes(\ff) = {\cF_{\gA,0}(\AX\sur\ff)}$ (voir le lemme d'\eli \gnl de la section \ref{subsecIdealResultant}).

\emph {1.}
On pose $\fa' = \rc_\uT\big(\Res(f, g_1T_1 + g_2T_2 + \cdots +
g_rT_r)\big)$. Montrer les inclusions:

\snic {
\fa \subseteq \fa' \subseteq \fb \subseteq \ff \cap \gA
.}

%\sni
\emph {2.}
Soient $\gA = \ZZ[a,b,c]$ où $a$, $b$, $c$ sont trois \idtrsz,
$f = X^d$, $g_1=a$, $g_2=b$ et $g_3=c$. Déterminer les 
\ids $\ff \cap \gA$, $\fa$, $\fa'$, $\fb$ et vérifier qu'ils
sont distincts. \linebreak 
Vérifier \egmt que $\fR(f, g_1, g_2, g_3)$ dépend de l'ordre des $g_i$. \\
Est-ce que l'on $(\ff\cap\gA)^d \subseteq \fa$?

}

\end {exercise}
%--- end -exercise-----------------------------------------

%--- Exercise{exoRelateursViaElimIdeal}-------------
\begin{exercise}\label{exoRelateursViaElimIdeal} 
{(Relateurs et \id d'élimination)}\index{elimination@\eliz!ideal@\id d'---}\index{ideal@idéal!d'elimi@d'\eliz}
\\
 {\rm  
Soient $f_1(\uX), \ldots, f_s(\uX) \in \kuX = \gk[\Xn]$ ($\gk$ est
un anneau commutatif). Notons $\fa \subseteq \kuY = \gk[Y_1, \ldots, Y_s]$
l'\id des relateurs sur $\gk$ de $(f_1, \ldots, f_s)$, \cad
$\fa = \ker\varphi$, où $\varphi : \kuY \to \kuX$ est le
morphisme d'\evn $Y_i \mapsto f_i$. On note $g_i = f_i(\uX) - Y_i \in
\gk[\uY, \uX]$ et $\ff = \gen {g_1, \ldots, g_s}$. \\
Montrer que $\fa = \ff \cap \kuY$.
Ainsi, $\fa$ est  l'\id d'\eli des variables $X_j$ dans le \syp des $g_i$.
}

\end {exercise}
%--- end -exercise-----------------------------------------

%: sinotenglish
\sinotenglish{

%:--- Exercise{exoZerdimloc}-------------
\begin{exercise}
\label{exoZerdimloc} 
{\rm \emph{(Localisation d'un anneau \zedz)}\\
Tout quotient ou localisé d'un anneau \zed est \zedz.
On s'intéresse ici aux \lons obtenues en inversant un seul \eltz. 

\noindent Soit $\gA$ un anneau \zed et $a\in\gA$.
\begin{enumerate}
\item L'anneau $\gA[1/a]$ s'identifie à un facteur de $\gA$. Plus \prmt
si $e$ est l'\idm  tel que $\gen{e}=\gen{a^d}$ pour $d$ assez grand on a un 
unique \isoz~\hbox{$\gA[1/e]\to \gA[1/a]$} qui factorise les \homos  $\gA\to \gA[1/e]$ \hbox{et $\gA\to \gA[1/a]$}.
\item Si $\gA$ est une \klg \stfe sur un \cdi $\gk$, il en est de même pour
 $\gA[1/a]$.
\item Dans le cas où $\gA=\gk[x]=\aqo{\kX}{f}$, 
\begin{enumerate}
\item l'\alg $\gA[1/g(x)]$ est isomorphe à une \klg $\aqo{\kX}{h}$ où $h$ est un diviseur de $f$ (le \pol $h$ est \gui{la partie de $f$ étrangère à $g$}, avec possiblement $h=1$),
\item l'\alg $\gA[1/f'(x)]$ est isomorphe à une \klg $\aqo{\kX}{f_1}$ où $f_1$ est un \pol \spl ($f_1$ est \gui{la partie sans carré} de $f$, \cad la partie étrangère à $f'$).
\end{enumerate}
\end{enumerate}
}
\end{exercise}
%--- end -exercise-----------------------------------------

%--- Exercise{exo}-------------
\begin{exercise}
\label{exoFittAnn}
{\rm  Si $M$ est un \Amo \pfz, pour tout $k\geq 1$ on a l'inclusion
{\fbox{$\Ann_\gA(M)\,\cF_k(M)\subseteq \cF_{k-1}(M)$}.}

}
\end{exercise}
%--- end -exercise-----------------------------------------

%--- Exercise{exoQuentelZerdimFrac}-------------
\begin{exercise}
\label{exoQuentelZerdimFrac} (Anneau réduit dont l'anneau total des fractions est \zedz)
{\rm (D'après  
\emph{Quentel Y., \emph{Sur la compacité du spectre minimal d'un anneau}. Bull. Soc. Math. France, {\bf 99}, (1971), 265--272}. Voir aussi le \pbz~\ref{exoAnneauNoetherienReduit})
\\ Soit $\gA$ un anneau réduit et $\gK$ son anneau total de fractions.
Alors \propeq 
\begin{enumerate}
\item $\gK$ est \zedz. 
\item \label{i2exoQuentelZerdimFrac}
Pour tout $x\in\gA$ il existe un  $y\in\gA$ tel que $xy=0$ et $\gen{x,y}$
contient un \elt \ndzz. 
\item \label{i3exoQuentelZerdimFrac}
Les deux conditions suivantes sont réalisées:
\begin{enumerate}
\item \label{i3aexoQuentelZerdimFrac} tout \itf fidèle de $\gA$ contient un \elt \ndzz,
\item \label{i3bexoQuentelZerdimFrac} pour tout $x\in\gA$ il existe un \itf $\fa$ tel que $x\fa=0$ et $\gen{x}+\fa$ est fidèle. 
\end{enumerate} 
\end{enumerate}
}
\end{exercise}
%--- end -exercise-----------------------------------------

%:--- Exercise{exoReduiteFrob}-------------
\begin{exercise}
\label{exoReduiteFrob}
 {(Forme réduite de Frobenius pour un \endo d'un \Kev de dimension finie, $\gK$ corps discret non trivial)} 
{\rm  Soit $\varphi:V\to V$ un \endo d'un \Kev $V$ de dimension $n$, et $F$
la matrice de $\varphi$ sur une base donnée $(e_1, \ldots, e_n)$.
On notera $\nu_\varphi$ et $\chi_\varphi$ le \polmin et le \polcar de $\varphi$. 
On munit $V$ d'une structure de $\KX$-module
en définissant la loi externe suivante 

\snic{\KX\times V\to V,\quad  (P,v)\mapsto P\cdot_\varphi v= P(\varphi)(v).}

En particulier  $X\cdot_\varphi v =\varphi(v)$. On note $V_\varphi$ le $\KX$-module ainsi obtenu. Les sous-modules de $V_\varphi$ sont exactement les sous-\Kevs stables par $\varphi$.

\emph{1.}  Pour tout $y$ non nul de $ V_\varphi$, le sous-module $\KX\cdot_\varphi y$ est le plus petit sous-\evc $\varphi$-stable de $V$ contenant $y$. 
On le notera aussi  
 $\Kfi \cdot y$.\\
Il admet une $\gK$-base de la forme  $\cB_{y,\varphi}=\big(y,\varphi(y),\dots,\varphi^{k-1}(y)\big)$, où 
$$\preskip.3em \postskip.4em
\ndsp
{\varphi^k(y)=a_0y+\som_{j=1}^{k-1}a_j\varphi^j(y)}
$$
est la première relation de dépendance $\gK$-\lin 
qui se présente entre~$y$ et ses transformés successifs par $\varphi$. On appelle \emph{\polmin de $y$ pour $\varphi$} et l'on note $\nu_{y,\varphi}(X)$ le \pol $X^k-\sum_{j=0}^{k-1}a_jX^j$. On a
alors  les résultats suivants.
\begin{itemize}
\item En tant que $\KX$-module, le sous-espace $\Kfi\cdot y$ est isomorphe  à $\aqo\KX {\nu_{y,\varphi}}$. L'\iso est donné par

\snic{\aqo\KX {\nu_{y,\varphi}}\,\lora\, \Kfi\cdot y,\quad \ov g\,\lmt\, g\cdot_\varphi y .}

 En particulier, en notant $x=\ov X$, l'image de  $(1,x,\dots,x^{k-1})$ 
($\gK$-base de $\aqo\KX {\nu_{y,\varphi}}$), est la base $\cB_{y,\varphi}$.

\item La matrice de la restriction de $\varphi$ à ${\Kfi\cdot y}$ sur la base $\cB_{y,\varphi}$ est la matrice compagne  du \pol $\nu_{y,\varphi}(X)$:% 
$$
 C_{\nu_{y,\varphi}} =\cmatrix{
0&\cdots&\cdots&\cdots&0&a_0\\[-3pt]
1&\ddots&&&\vdots&a_{1}\\[-1pt]
0&\ddots&\ddots&&\vdots&\vdots\\
\vdots&\ddots&\ddots&\ddots&\vdots&\vdots\\[-3pt]
\vdots& & \ddots&\ddots&0&a_{k-2}\\[3pt]
0&\cdots&\cdots&0&1&a_{k-1}
}.
$$

\end{itemize}

\emph{2.} On note encore $\varphi$ l'extension de $\varphi$ en un \endo de $V[X]=\KX\otimes_\gK V$. On note $\theta_\varphi:V[X]\to V_\varphi$ l'unique application $\KX$-\lin qui donne l'identité sur $V$.
 L'exercice \ref{exoAXmodule} montre que la suite ci-dessous est exacte:

\snic {
V[X]\vvvvvers{X\Id_V-\varphi} {V[X]} \vvers{\theta_\varphi} V \lora 0.
}

\smallskip 
Autrement dit $V_\varphi$ est un $\KX$-\mpf et $X\In-F$ est une \mpn de $V_\varphi$ pour le
\sgr $(e_1, \ldots, e_n)$. 

Comme $\KX$ est un anneau principal (à \dve explicite donc \fdiz)
on peut appliquer le \tho de structure donné par la proposition \ref{propPfPID}.
\begin{itemize}
\item Expliciter ce que ce \tho de structure donne en terme de \dcn de~$V$ en somme directe de sous-\evcs stables (par $\varphi$).
\item Donner une forme réduite de la matrice de $\varphi$ qui résulte du \tho de structure. Comparer les facteurs invariants \hbox{du $\KX$-module} $V_\varphi$ avec le \polcar $\chi_\varphi$ et le \polmin $\nu_\varphi$ de $\varphi$. 
\item Montrer que les facteurs invariants \hbox{du $\KX$-module} $V_\varphi$ caractérisent la classe de similitude de $\varphi$ comme endomorphisme de $V$.  On les appelle les \ixx{invariants de similitude}{de l'\endo $\varphi$}.
\item Les invariants de similitude peuvent être calculés dans le sous-corps $\gK_1$ engendré par les \coe de $F$. Ils ne changent pas si l'on étend les scalaires à n'importe quel surcorps de $\gK$.
\item 
Expliquer ce que donne le lemme des noyaux \ref{lemDesNoyaux} lorsque l'on a une \dcn du \polmin $\nu_\varphi$ de $\varphi$ en un produits de \pols deux à deux
étrangers. Par exemple on peut considérer une \bdf obtenue à partir 
la famille des invariants de similitude. 

\end{itemize}
 
}
\end{exercise}
%--- end -exercise-----------------------------------------

%:--- Exercise{exoEndosemisimples}-------------
\begin{exercise}
\label{exoEndosemisimples}
 {(Endomorphismes semi-simples)} Suite de l'exercice \ref{exoReduiteFrob}.
\\
{\rm  Un \endo $\varphi$ de $V\simeq \gK^{n}$ est dit \ixc{semi-simple}{endomorphisme} si 
tout sous-\evc stable est \sul d'un sous-espace stable.
\begin{enumerate}
\item Si $f=g^{2}h\in\KX$, de degré $r$, l'\endo  

\snic{\varphi: V\lora V,\;\ov g\mapsto x\ov{g},\quad  \hbox{ où } V=\aqo\KX f \hbox{ et }x=\ov X,}

qui est représenté sur la base $(1,x,\dots,x^{r-1})$ par la matrice compagne de $f$,
n'est pas semi simple.
\item Si le \polcar  $\chi_\varphi$  se décompose sous forme $\prod_i h_i^{m_i}$  avec les $h_i$ \irds \unts deux à deux distincts, l'\endo $\varphi$ est semi-simple \ssi son \polmin est sans facteur carré, \cad  égal à $\prod_i h_i$.
\item Supposons $\gK$  \agqt clos. Un \endo
est semi-simple \ssi il est diagonalisable.
\item Soit $S$ un sous-espace stable de $V$.
On peut certifier de manière \algq  une des deux alternatives
de la disjonction suivante:\begin{itemize}
\item le sous-espace $S$ admet un \sul stable $T$,
 OU
\item le \polmin $\nu_\varphi$ admet un facteur carré.
\end{itemize}
Dans le premier cas, on a décomposé $\nu_\varphi$ en un produit
$\prod_{j=1}^{m}f_j$ de \pols 
deux à deux étrangers, et l'on obtient, en posant 

\snic{K_j=\Ker\big(f_j(\varphi)\big)$
et $\gK_j=\aqo\KX{f_j},}

\begin{itemize}
\item  $S=\bigoplus_{j=1}^{m} S_j$, où $S_j=S\cap K_j$,
\item  $T=\bigoplus_{j=1}^{m} T_j$, où $T_j=T\cap K_j$,
\item  $S_j$ et $T_j$ sont des $\gK_j$-modules libres de rang fini,
\item  pour chaque $j$, $K_j=S_j\oplus T_j$.
\end{itemize}
\item Si sur le corps $\gK$ on sait décomposer le \polmin de $\varphi$
en produit de \pols \splsz, on sait tester si $\varphi$ est semi-simple,
et expliciter ce caractère semi-simple lorsque la réponse est positive. 
\end{enumerate}
}
\end{exercise}
%--- end -exercise-----------------------------------------

}
%: fin sinotenglish

%: problemes

%--- problem{exoDimZeroXcYbZa}-------------
\begin{problem}\label{exoDimZeroXcYbZa}
 {(Un exemple  de \sys \zedz)}\\
{\rm  
Soient $\gk$ un anneau et $a$, $b$, $c \in \NN^*$ avec $a \le b \le c$ et au moins une in\egt stricte. On définit trois \pols $f_i \in \gk[X,Y,Z]$:

\snic {
f_1 = X^c+Y^b+Z^a,\quad  f_2 = X^a+Y^c+Z^b,\quad  f_3 = X^b+Y^a+Z^c
.}

%\sni

Il s'agit d'étudier le \sys défini par ces trois \polsz. On note 
$\gA = \gk[x,y,z]$ \hbox{la \klg} $\aqo{\gk[X,Y,Z]}{f_1,f_2,f_3}$.

\emph {1.}
Pour un anneau quelconque $\gk$, $\gA$ est-elle libre finie
sur $\gk$? Si oui, calculer une base et donner la dimension.

\emph {2.}
\'{E}tudier de manière détaillée le \sys pour $\gk = \QQ$ et $(a,b,c) =
(2,2,3)$: déterminer tous les zéros du \sys dans une certaine extension
finie de $\QQ$ (à préciser), leur nombre et leurs multiplicités.

\emph {3.}
L'\alg localisée $\gA_{1+\gen{x,y,z}}$ est-elle libre sur $\gk$?  Si oui, donner une
base.

}

\end {problem}
%--- end -problem-----------------------------------------

%--- problem{exoIdealResultantGenerique}-------------
\begin{problem}
\label{exoIdealResultantGenerique} (L'\id résultant générique)\\
{\rm 
Soient $d$, $r$ deux entiers fixés avec $d\ge 1$. On étudie dans cet exercice
l'idéal résultant générique $\fb = \fRes(f, g_1, \ldots, g_r)$ où $f$
est \mon de degré $d$, \hbox{et $g_1$, \ldots, $g_r$} sont de degré $d-1$, les
\coes de ces \pols étant des \idtrs sur $\ZZ$. L'anneau de base est
donc $\gk =\ZZ[(a_i)_{i\in \lrb{1..d}}, (b_{ji})_{j\in \lrb{1..r}, i\in
\lrb{1..d}}]$ avec

\snic{
f = X^d + \sum_{i=1}^{d} a_iX^{d-i} \quad\hbox{et}\quad  g_j = \sum_{i=1}^{d} b_{ji}X^{d-i}.
}

%\sni
\emph{1.}  
Mettre des poids sur les $a_i$ et $b_{ij}$ de façon à ce que $\fb$ soit
un \id \hmgz.

\emph{2.}
Si $S$ est la matrice de Sylvester généralisée de $(f, g_1, \ldots, g_r)$,
préciser le poids des \coes de $S$ et ceux de ses mineurs
d'ordre $d$.

\emph{3.}
A l'aide d'un \sys de Calcul Formel, étudier le nombre minimal de
\gtrs de $\fb$. On pourra remplacer $\ZZ$ par $\QQ$, introduire
l'\id $\fm$ de $\gk$ engendré par toutes les \idtrs et considérer
$E = \fb\sur{\fm\fb}$ qui est un $\gk\sur{\fm} = \QQ$-\evc de dimension
finie.

}
\end{problem}
%--- end -problem-----------------------------------------

%:--- problem{exoNakayamaHomogeneRegularSequence}-------------
\begin{problem}\label{exoNakayamaHomogeneRegularSequence}
{(Nakayama \hmg et suites \ndzesz)}
\\
{\rm  
\noindent\emph {1.} \emph{(Suite \ndze et indépendance \agqz)}
Soit $(\an)$ une \seqreg d'un anneau $\gA$ et $\gk \subseteq \gA$
un sous-anneau tel que $\gk \cap \gen {\an} = \{0\}$. Montrer
que $a_1$, \dots, $a_n$ sont algébriquement indépendants sur~$\gk$.

\emph {2.} {\it(Nakayama \hmgz)}
Soient $\gA = \gA_0 \oplus \gA_1 \oplus \gA_2 \oplus \dots$
un anneau gradué \hbox{et $E = E_0 \oplus E_1 \oplus E_2 \oplus \dots$}
un \Amo gradué. \\
On note $\gA_+$ l'\id $\gA_1 \oplus \gA_2 \oplus \dots$,
si bien que $\gA/\gA_+ \simeq \gA_0$.
\vspace{-2pt}
\begin {itemize}\itemsep=0pt
\item [\emph {a.}]
Si $\gA_+ E = E$, alors $E = 0$.

\item [\emph {b.}]
Soit $(e_i)_{i \in I}$ une famille d'\elts \hmgs de $E$. Si les $e_i$ engendrent le $\gA_0$-module $E/\gA_+E$, alors ils engendrent le \Amoz~$E$.
\end {itemize}
%\vspace{-2pt}
Notez que l'on n'a pas supposé $E$ \tfz.

\emph {3.}
Soit $\gB = \gB_0 \oplus \gB_1 \oplus \gB_2 \oplus \dots$ un anneau gradué
et $h_1$, \ldots, $h_d$ des \elts \hmgs de l'\id $\gB_+$. On note $\fb = \gen {h_1, \ldots, h_d}$ et $\gA = \gB_0[h_1,
\ldots, h_d]$. On a donc $\gB_0 \cap \fb = \{0\}$, et $\gA$ est un sous-anneau
gradué de $\gB$.  Enfin, soit $(e_i)_{i \in I}$ une famille d'\elts \hmgs de
$\gB$ qui engendrent le $\gB_0$-module $\gB/\fb$.
%
%\vspace{-2pt}
\begin {itemize}%\itemsep=0pt
\item [\emph {a.}]
Vérifier que $\gA_0 = \gB_0$ et que $\fb = \gA_+ \gB$ puis
montrer que les $e_i$ forment un \sgr du \Amo $\gB$.
\item [\emph {b.}]
On suppose que $(h_1, \ldots, h_d)$ est une \seqreg et que les $e_i$
forment une base du $\gB_0$-module $\gB/\fb$.  Montrer que $ h_1$, \ldots, $h_d$
sont \agqt indépendants sur~$\gB_0$ et que les $e_i$ forment une
base du \Amo $\gB$.
\end {itemize}
\emph {Bilan}:
Soit $\gB = \gB_0 \oplus \gB_1 \oplus \gB_2 \oplus \dots$ un anneau gradué
et $(h_1, \ldots, h_d)$ une \seqreg \hmg de l'\id $\gB_+$.  Si
$\aqo{\gB}{h_1, \ldots,h_d}$ est un $\gB_0$-module libre, alors $\gB$ est un $\gB_0[h_1, \ldots,
h_d]$-module libre et $\gB_0[h_1, \ldots, h_d]$ est un anneau de \pols en
$(h_1, \ldots, h_d)$.

\emph {4.}
En guise de réciproque.
Soit $\gB = \gB_0 \oplus \gB_1 \oplus \gB_2 \oplus \dots$ un anneau gradué
\hbox{et $h_1$, \ldots, $h_d$} des \elts \hmgs de l'\id $\gB_+$, \agqt indépendants
sur $\gB_0$. Si $\gB$ est un $\gB_0[h_1, \ldots, h_d]$-module libre,
alors la suite $(h_1, \ldots, h_d)$ est \ndzez.

}
\end{problem}

}% fin des exos

%\vspace{-10pt}
%:   solutions d'exos
\sol

%%%%%%%%%%%%%%     exoptfpf      %%%%%%%%%%%%%%
\exer{exoptfpf} {Il suffit d'appliquer la proposition \ref{propPfSex}. Directement:
on considère un \prr $\pi:N\to N$ ayant pour image $M$. 
Si $X$ est un \sgr de~$N$,~$\pi(X)$ est un \sgr de $M$. Si $N$ est \pfz, le module de syzygies
pour $\pi(X)$ est obtenu en prenant les syzygies pour $X$ dans~$N$ et les syzygies
$\pi(x)=x$ pour chaque \elt $x$ de~$X$.
}

%%%%%%%%%%%%%     exoAXmodule   %%%%%%%%%%%%%%%%%%

\exer{exoAXmodule} On note pour commencer que $\theta_A \circ (X\In - A) = 0$
et que $\theta_A$ est l'identité sur $\Ae n$.

%\sni
\emph {1.} 
Montrons que
$\Im(X\In-A) \cap \Ae n = 0$.  Soit $x \in \Im(X\In-A) \cap \Ae n$,
les calculs préliminaires donnent $\theta_A(x)=x$ et $\theta_A(x) = 0$.
Montrons que $\AX^n = \Im(X\In-A) + \Ae n$.
Il suffit de voir que $X^k e_i \in \Im(X\In-A) + \Ae n$ pour $k\geq 0$ et $i\in\lrbn$. Si $k=0$ c'est clair, pour $k>0$ on écrit:

\snic {
X^k\In-A^k = (X\In-A) \, \sum_{j+\ell = k-1} X^j A^\ell.
}

%\sni 
En appliquant cette \egt à $e_i$, on obtient $X^k e_i - A^k e_i \in \Im(X\In-A)$,
\hbox{donc $X^k e_i \in \Im(X\In-A) + A^k e_i \subseteq \Im(X\In-A) + \Ae n$}.

%\sni
\emph {2.}
Soit $y \in \Ker\theta_A$.  On écrit $y = z + w$ avec $z
\in \Im(X\In-A)$ et $w \in \Ae n$. \\
Donc $0 =\theta_A(y) = \theta_A(z) + \theta_A(w) = 0 + w$ et $y = z \in \Im(X\In-A)$.

%%%%%%%%%%%%%     exoAnneauBézoutStrict   %%%%%%%%%%%%%%%%%%

\exer{exoAnneauBézoutStrict}
\emph{1} et \emph{2.} laissés \alecz.

\emph{3.}
Par construction, $\vep$ annule $d_2$ (i.e. annule $a,b$). On a donc les \egts

\snuc{
d_2\cmatrix {a_2\cr b_2} = 
\crmatrix {v_1 & a_1\cr -u_1 & b_1} \cmatrix {d_2\vep \cr d_2k_{12}} =
\crmatrix {v_1 & a_1\cr -u_1 & b_1} \cmatrix {0 \cr d_1} =
d_1\cmatrix {a_1\cr b_1} = \cmatrix {a\cr b} 
}

%\sni
Reste à voir que $1 \in \gen {a_2, b_2}$. En inversant
la matrice $2\times 2$ dans $(\star)$ (de \deterz~1), on
voit que l'\id $\gen {a_2, b_2}$ contient $\vep$ et $k_{12}$, donc
il contient $1 = k_{12}k_{21} - \vep$.

\emph{4.}
Par \recu sur $n$, $n=2$ étant la question précédente.
Supposons $n \ge 3$. Par \recuz, il existe $b_1$, \ldots, $b_{n-1}$
\com et $d$ tels que 

\snic{(a_1, \ldots, a_{n-1}) = d(b_1, \ldots, b_{n-1}),\hbox{ donc }\gen{\ua} = \gen {d, a_n}.}

%\sni
Le point \emph {3}
donne~$u$,~$v$ \com et~$\delta$ tels que $(d,a_n) = \delta(u,v)$.
\\
Alors
%\snic {
$
(\an) = (db_1, \ldots, db_{n-1}, \delta v) =
\delta (ub_1, \ldots, ub_{n-1}, v).
$
%}
%
%\sni
\\ Et $\gen {ub_1, \ldots, ub_{n-1}, v} = \gen {u,v} = \gen {1}$.

\emph{5.}
D'abord pour $n=2$ avec $(a,b)$. Il y a $d$ avec $(a,b) = d(a',b')$
et $1 = ua' + vb'$. On pose $m = da'b' = ab' = ba' \in \gen {a}\cap\gen{b}$;
on a $\gen {a} \cap \gen{b} = \gen{m}$ car si $x \in \gen {a} \cap \gen{b}$,
on écrit $x=x(ua' + vb') \in \gen {ba'} + \gen {ab'} = \gen{m}$.
La $\SL_2(\gA)$-\eqvc est fournie par l'\egt ci-dessous:

\snic {
\cmatrix {1 &-1 \cr vb' & ua'} \cmatrix {a &0 \cr 0 & b} =
\cmatrix {d & 0 \cr 0 & m} \cmatrix {a' &-b' \cr v & u}
.}

%\sni
Pour $n \ge 3$. En utilisant le cas $n=2$ pour les positions
$(1,2)$, $(1,3)$, \dots, $(1,n)$, on obtient $\Diag(a_1, a_2, \ldots, a_n)
\sim \Diag(a'_1, a'_2, \ldots, a'_n)$ avec $a'_1 \divi a'_i$ pour $i \ge 2$.
\\
Par \recuz, $\Diag(a'_2, \ldots, a'_n) \sim \Diag(b_2, \ldots, b_n)$
avec $b_2 \divi b_3 \cdots \divi b_n$. On vérifie alors
que $a'_1 \divi b_2$ et l'on pose $b_1 = a'_1$.
\Llec scrupu\leux vérifiera la \prt concernant les
fonctions \smqs \elrsz.

%%%%%%%%%%%%%     exoSmith   %%%%%%%%%%%%%%%%%%

\exer{exoSmith} \emph{(Anneaux de Smith, ou elementary divisor rings)}
\\
Calcul préliminaire avec $A = \cmatrix {a& b\cr 0 &c}$
et $B$ de la forme:

\snic {
B = \cmatrix {p' & q'\cr *  &* } A \cmatrix {p & * \cr q &* }
.}

%\sni
Le \coe $b_{11}$ de $B$ est égal à $b_{11} = p'(pa + qb) + q'qc$.

$\emph {2} \Rightarrow \emph {3.}$ 
La matrice $A$ est \eqve à une matrice
diagonale $\Diag(g,h)$, ce qui donne $(p,q)$ et $(p',q')$ \com avec $g
= p'(pa + qb) + q'qc$ (calcul préliminaire).  Et l'on a $\gen {a,b,c} = \gen
{g,h}$.  Comme $\gA$ est de Bézout strict, on peut supposer $g \divi h$ et
puisque $1 \in \gen {a,b,c}$, $g$ est \iv et donc $1$ s'écrit comme voulu.

$\emph {3} \Rightarrow \emph {4.}$ 
Attention, ici $g$ est imposé. Mais d'après
la question \emph {4} de l'exercice~\ref{exoAnneauBézoutStrict},
%:2012 ci-dessous  g\,(a',b',c')
on peut écrire $(a,b,c) = g\,(a',b',c')$ avec $(a',b',c')$
\comz. On applique le point~\emph {3} à $(a',b',c')$
et on multiplie le résultat obtenu par $g$.

$\emph {4} \Rightarrow \emph {2.}$ 
Soit $A \in \MM_2(\gA)$ triangulaire, $A = \cmatrix {a& b\cr 0 &c}$.  Avec les
paramètres du point~\emph {4}, on construit (calcul préliminaire) une
matrice $B$ \eqve à $A$ de \coe $b_{11} = g$.  Comme $g$ divise tous les
\coes de $B$, on a $B \,\sims{\EE_2(\gA)}\, \Diag(g,h)$.

\emph {1} $\Leftrightarrow$ \emph {2.} Laissé  \alec (qui pourra
consulter l'article de Kaplansky).

%%%%%%%%%%%%%     exoUnitsOfSomeRings   %%%%%%%%%%%%%%%%%%
\exer{exoUnitsOfSomeRings} 
\emph {1.}
Soit $s = y+z$; alors $\gk[s]$ est un anneau de \pols en $s$, et~$y$,~$z$ sont
entiers sur $\gk[s]$, car zéros de $(T-y)(T-z) = T(T-s)
\in \gk[s][T]$. On vérifie facilement que $\gA$ est libre sur $\gk[s]$ avec
 $(1,y)$ pour base.  Pour $u$, $v \in \gk[s]$, la norme sur $\gk[s]$ de 
$u + vy$ est:

\snic {
\rN_{\gA/\gk[s]}(u+vy) = (u+vy)(u+vz) = u^2 + suv = u(u + sv)
.}

%\sni
L'\elt $u + vy$ est \iv dans $\gA$ \ssi $u(u + sv)$ est \iv
dans $\gk[s]$. Comme $\gk$ est réduit, $(\gk[s])\eti = \gk\eti$;
donc $u \in \gk\eti$ et $v = 0$.

\emph {2.}
On a $\gA = \ZZ[\alpha,\beta,y] = \aqo{\ZZ[a,b,Y]}{(ab-1)Y}$ avec
$y(\alpha\beta-1) = 0$. Soit $t$ une \idtr sur $\ZZ$ et
$\gk = \ZZ[t,t^{-1}]$. Considérons la \klg $\gk[y,z]$ avec la seule relation $yz= 0$.
On a un morphisme $\gA \to \gk[y,z]$ qui réalise

\snic{\alpha \mapsto t(z+1)$, $\beta \mapsto t^{-1}$, $y \mapsto y,}

%\sni
et l'on vérifie que c'est une injection.\\
 Alors un \elt $w \in \Ati$
est aussi dans~$\gk[y,z]^{\times}$, et comme~$\gk$ est réduit, $w \in \gk\eti$.
Enfin, les \ivs de $\gk = \ZZ[t,t^{-1}]$
sont les $\pm t^k$ avec $k \in \ZZ$, donc $w = \pm 1$.

%%%%%%%%%%%%%     exoUAtoUB   %%%%%%%%%%%%%%%%%%

\exer{exoUAtoUB} \emph {1.}
On sait que $\gen {x^n} = \gen {e}$.  On cherche $y \in \Ati$ tel que $y
\equiv x \bmod \fa$ sur les composantes $\gA_e$ et $\gA_{1-e}$; d'abord, on a
$x^n (1-ax) = 0$ et $x$ \iv modulo $\fa$, donc $ax \equiv 1 \bmod \fa$ puis $e
\equiv 1 \bmod\fa$, i.e. $1-e \in \fa$.  \\
Dans la composante $\gA_e$, $x$ est
\ivz, donc on peut prendre $y=x$. Dans la composante $\gA_{1-e}$, $1\in \fa$,
donc on peut prendre $y=1$.  Globalement, on propose donc $y = ex + 1-e$ qui est
bien \iv (d'inverse $ea^{n}x^{n-1} + 1-e$) et qui vérifie $y \equiv x \bmod
\fa$.  Remarque: $y = ex + (1-e)u$ avec $u \in \Ati$ convient aussi.

\emph {2.}
Soit $x$ inversible modulo $\fa$ donc $1-ax \in \fa$ pour un certain 
$a \in \gA$. \\
Alors, $x$ est \iv modulo l'\id principal $\gen {1-ax}$,
donc il existe $y \in \Ati$ tel que $y \equiv x \bmod \gen {1-ax}$,
a fortiori $y \equiv x \bmod \fa$.

\emph {3.}
On écrit $y=bx$, $x=ay$ donc $(1-ab)x = 0$; $b$ est \iv modulo 
$\gen {1-ab}$ donc il existe $u \in \Ati$ tel que $u \equiv b \bmod
\gen{1-ab}$ d'où $ux = bx = y$.

\emph {4.}
Soit $x$ \iv modulo $\fa$. Alors $\pi(x)$ est \iv modulo $\fa'$, d'où
 $y \in \gA$ tel que~$\pi(y)$ soit \iv dans $\gA'$ et $\pi(y) \equiv
\pi(x) \bmod \fa'$. Alors, $y$ est \iv dans $\gA$ et $y-x \in \fa + \Rad\gA$, i.e. 
$y = x + a + z$ avec $a \in \fa$ et $z \in \Rad\gA$.
Ainsi, l'\elt $y-z$ est \iv dans $\gA$, et $y-z \equiv x \bmod \fa$.

%%%%%%%%%%%%%     exoCalculT(M)   %%%%%%%%%%%%%%%%%%
%:HHH corrige nouvel exo 
\exer{exoCalculT(M)}  
Appelons $\gA_1$ le corps de fractions de $\gA$ et mettons un indice $1$
pour indiquer que nous faisons une \eds de $\gA$ à $\gA_1$. Ainsi $M_1$ est le $\gA_1$-\evc correspondant à la suite exacte

\snic{\gA^n_1 \vers{E_1} \gA^\ell_1 \vers{\pi_1} M_1\to 0}

et le sous-module $\rT(M)$ de $M$ est le noyau de l'\Ali naturelle de $M$ vers $M_1$,
i.e. le module $\pi(\Ae\ell \cap \Ker \pi_1)$, ou encore 
le module $\pi(\Ae\ell \cap \Im E_1)$
 (en regardant $\Ae\ell$
comme un sous-module de $\gA^\ell_1$).
\\
 La suite exacte
% 
%\vspace{-2pt}  
%
%\snic{
$\Ae m \vers{F} \Ae\ell \vers {\tra E} \Ae n$
%}
%\sni
donne par \lon la suite exacte
 
\vspace{-1pt}  

\snic{
\gA^ m_1 \vers{F_1} \gA^\ell_1 \vers {\tra {E_1}} \gA^n_1}

%\sni
%\penalty-2500
et puisque $\gA_1$ est un \cdi cela donne par dualité la suite exacte
\vspace{-1pt}

\snic{\gA^n_1 \vers{E_1} \gA^\ell_1 \vers {\tra {F_1}} \gA^ m_1.}

%\sni
Ainsi $\Im E_1=\Ker \!{\tra {F_1}}$, donc $\Ae\ell \cap \Im E_1=\Ae\ell \cap \Ker\! {\tra {F_1}}$. Enfin on a 
l'\egtz~\hbox{$\Ae\ell \cap \Ker\! {\tra {F_1}}=\Ker\! {\tra F}$}
car le morphisme naturel $\gA\to \gA_1$ est injectif.
\\
Conclusion: $\rT(M)$ est égal à $\pi(\Ker \!{\tra F})$, isomorphe à
$\Ker\!{\tra F}/\Im E$, donc \pf (parce que $\gA$ est \cohz).

Si $\gA$ est \coh \qi l'anneau total des fractions $\gA_1=\Frac\gA$
est \zedrz, et tous les arguments donnés dans le cas intègre fonctionnent
pareillement.

%%%%%%%%%%%%%     exoZerRedBez   %%%%%%%%%%%%%%%%%%

\exer{exoZerRedBez} 
Tous les résultats peuvent être obtenus à partir du cas des \cdisz, cas
pour lequel les \algos sont classiques, en utilisant la machinerie \lgbe \elr
des anneaux \zedrsz. On va préciser ici un peu cette affirmation de
caractère très \gnlz.\imlg
\\
Faisons deux remarques préliminaires pour un anneau quelconque $\gA$.
\\
Premièrement, soit $e$ \idm et $E$ une matrice \elr modulo~${1-e}$.\linebreak  Si l'on  relève $E$ en une matrice $F\in \Mn(\gA)$, alors
la matrice $(1-e)\In+eF\in\En(\gA)$ est \elrz, elle agit comme $E$ dans la composante $\aqo\gA {1-e}$, et elle ne fait
rien dans la composante $\aqo\gA {e}$. Ceci permet de comprendre comment on
peut récupérer les résultats souhaités sur $\gA$ en utilisant des
résultats analogues modulo les \idms $1-e_i$ lorsque l'on dispose d'un \sfio
$(e_1,\ldots,e_k)$ (fourni par l'\algo que l'on construit).\\ 
Deuxièmement, si
$g \in \AX$ est \mon de degré $m \ge 0$, pour tout $f \in \AX$, on peut
diviser $f$ par $g$: $f = gq + r$ avec $r$ de degré formel $m-1$.

 \emph {1.}
Soit $e$ l'\idm tel que $\gen {e} = \gen {b}$. Il suffit de résoudre la
question modulo $e$ et $1-e$.  Dans la branche $e = 1$, $b$ est \ivz, $\gen
{a,b} = \gen {1}$ et le \pb est résolu (pivot de Gauss). 
Dans la branche $e = 0$, $b$ est nul et le \pb est résolu.
Si $e = bx$, on trouve $d = e + (1-e)a$ et

\snic {
M = \rE_{21}(-be) \rE_{12}\big(ex(1-a)\big) =
\cmatrix {1 & ex(1-a)\cr -eb & ae + (1-e)\cr}
.}

%\sni
\emph {2.} 
On part de deux \pols $f$ et $g$.  On va construire un \pol $h$ et une matrice
$M\in\EE_2(\AX)$ telles que $M\cmatrix{f\cr g}=\cmatrix{h\cr 0}$. A fortiori
$\gen{f,g}=\gen{h}$.
 \\
On procède par \recu sur $m$, degré formel de $g$, de \coe formellement
dominant $b$.  Si l'on amorce la \recu à $m=-1$, $g = 0$ et
$\I_2\cmatrix{f\cr g}=\cmatrix{f\cr 0}$.  On peut traiter $m=0$,
avec $g \in \gA$ et utiliser le point \emph {1}   ($\gB = \AX$, $a=f$,
$b=g$). Mais il est inutile de traiter ce cas à part (et donc on
n'utilise plus le point~\emph {1}).  En effet, si $e$ l'\idm tel
que $\gen {e} = \gen {b}$, il suffit de résoudre la question modulo $e$ et
$1-e$ et ce qui suit est valide pour tout $m \ge 0$.
\\
Dans la branche $e = 1$, $b$ est \ivz, et puisque $m \ge 0$, on peut 
réaliser une division euclidienne classique de $f$ par $g$: $f = qg -r$
avec le degré formel de $r$ égal à $m-1$.  Ce qui donne une matrice
$N\in\EE_2(\AX)$ telle que $N\cmatrix{f\cr g}=\cmatrix{g\cr r}$, à~savoir
$N = \crmatrix {0 &1\cr -1 &q\cr}$.  On peut alors appliquer l'\hdrz.
\\
Dans la branche $e = 0$, $g$ est de degré formel $m-1$ et l'\hdr s'applique.

 Dans la suite, on utilise le point \emph{2} en disant que l'on passe de~$\tra{\vab f  g}$ à~$\tra{\vab  h 0}$ au moyen de
\gui{manipulations de Bézout}.

 \emph{3.} En s'appuyant sur le résultat du point \emph{2}  on
s'inspire de la \dem de la proposition \ref{propPfPID} (un anneau principal
est un anneau de Smith).  Si l'on était sur un \cdi non trivial, l'\algo
terminerait en un nombre fini d'étapes qui peut être borné directement
en fonction de $(D,m,n)$, où $D$ est le degré maximum des \coes de
la matrice $M\in \AX ^{m\times n}$ que l'on désire réduire à la forme de
Smith. Il s'ensuit que lorsque $\gA$ est \zedr le nombre de scindages produits
par les calculs de pgcd (comme au point \emph{2}) est lui aussi borné en
fonction de $(D,m,n)$, où $D$ est maintenant le degré formel maximum 
des \coes de la matrice. Ceci montre que l'\algo complet, compte tenu de la
remarque préliminaire, termine lui aussi en un nombre d'étapes borné en
fonction de $(D,m,n)$.

%\sni
\rem  Les \algos ne demandent pas que $\gA$ soit discret.
\eoe

%%%%%%%%%%%%%%%%%%%%%%%%%%%%%%%%%%%%%%%%%%%%%%%%%%%%%%%%%%%%%%%%%%%%%%%%%%%
\exer{exoFreeAlgebraPresentation}
Il est clair que $\fa\subseteq\Ker\varphi$.
Soit $\cE \subseteq \AuX$ l'ensemble des \polsz~$f$ congrus modulo $\fa$ à un
\pol \hmg de degré 1.\\ On a $1 \in \cE$ et $f \in \cE \Rightarrow
X_if \in \cE$ car si $f \equiv \sum_j \alpha_j X_j \mod {\fa}$, alors:

\snic{
X_if \equiv \som_j \alpha_j X_iX_j \equiv \som_{j,k}
\alpha_j c_{ij}^k X_k \mod {\fa}.}

%\sni
Donc $\cE = \AuX$. Soit $f \in \Ker\varphi$; on écrit $f \equiv \sum_k \alpha_k X_k \mod {\fa}$.
\\
Alors $\varphi(f)=0 = \sum_k \alpha_k e_k$,
donc $\alpha_k = 0$, puis $f \in \fa$.
%%%%%%%%%%%%%%%%%%%%%%%%%%%%%%%%%%%%%%%%%%%%%%%%%%%%%%%%%%%%%%%%%%%%%%%%%%%

\exer{exoFitt0}~\\
\emph{2.} Si $\fa$ est \tf une \mpn du module $M=\gA\sur\fa$ est une matrice ligne $L$
ayant pour \coes des \gtrs de l'\idz. On en déduit que $\cD_1(L)=\fa$. Donc $\cF_{-1}(M)=0\subseteq\cF_0(M)=\fa\subseteq\cF_1(M)=\gen{1}$. 
Le résultat se \gns à un \id $\fa$ arbitraire.

\emph{3.} Résulte de \emph{2} et du fait \ref{fact.idf.sex}.

\emph{4} et \emph{5.} Dans le cas \gnl en appliquant \emph{2} et \emph{3} on trouve: 

\snic{\cF_0(M)=\prod_{i=1}^n \fa_i$, $\cF_{n-1}(M)=\som_{i=1}^n \fa_i,}

%\sni
et pour les \ids intermédiaires les \gui{fonctions \smqsz} 

\snic{\cF_{n-k}(M)=\som_{1\leq i_1<\ldots< i_k\leq n} \prod_{\ell=1}^k \fa_{i_\ell}.}

%\sni
Par ailleurs, $\Ann (M)=\fa_1\cap\cdots\cap\fa_n$.
\\
Lorsque $\fa_1\subseteq\fa_2\subseteq \cdots\subseteq\fa_n$ le résultat est un peu plus simple 

\snic{\cF_{n-1}(M)=\fa_n$, $\cF_{n-2}(M)=\fa_n\,\fa_{n-1}$, \ldots\, $\cF_{n-k}(M)=\fa_n\cdots \fa_{n-k+1}.}

%\sni
On retrouve alors pour le point \emph{1} le résultat du calcul direct donné par
les \idds d'une matrice en forme de Smith.

%%%%%%%%%%%%%%     exoFitt6      %%%%%%%%%%%%%%%%%
\exer{exoFitt6}{
Montrons le point \emph{1} (après, on peut appliquer le fait \ref{facttfpf}).\\
Prenons $M=\gen{g_1,\ldots ,g_q}$. On considère une syzygie $\som_i\alpha_ig_i=_N0$.
Le but est de montrer que le vecteur colonne $V=(\alpha_1,\ldots ,\alpha_q)$ est
une syzygie dans $M$.

\emph{Premier cas, $M$ est \pfz.} \\
La colonne $V$, ajoutée à une \mpn $F$ de $M$ pour $(g_1,\ldots
,g_q)$ ne change pas les \idds de cette matrice, donc $V$ est une \coli des
colonnes de $F$. 

\emph{Deuxième cas, $M$ est \tfz.} \\
Puisque $\cD_1(V)\subseteq \cF_{q-1}(M)$, il existe
une matrice $F_1$ de syzygies  pour $(g_1,\ldots ,g_q)$ dans $M$ avec
$\cD_1(V)\subseteq \cD_{1}(F_1)$. Puisque  $\cD_2(V|F_1)\subseteq \cF_{q-2}(M)$, il
existe une matrice $F_2$ de syzygies  pour $(g_1,\ldots ,g_q)$ dans $M$ avec
$\cD_2(V|F_1)\subseteq \cD_{2}(F_1|F_2)$,
mais aussi bien sûr $\cD_1(V|F_1)\subseteq \cD_{1}(F_1|F_2)$. Et ainsi de
suite jusqu'à: il existe une matrice $F=[\,F_1\mid\cdots \mid F_q\,]$ de syzygies  pour
$(g_1,\ldots ,g_q)$ dans $M$ telle que les \idds de $[V|F]$ soient contenus dans
ceux de $F$.
Donc $V$ est une \coli des colonnes de $F$.
}

%%%%%%%%%%%%%%%%%%%%%%%%%%%%%%%%%%%%%%%%%%%%%%%%%%%%%%%%%%%%%%%%%%%%%%%%%%%
\exer{exoAutresIdF}
Si un \id de Kaplansky est égal à $1$, cela implique que le module est \tfz, car le module est \tf pour des  $\gA[1/a_i]$  avec des  $a_i$  \comz.
\\
Morale: les \ids de Kaplansky sont un peu plus \gnlsz, mais apparemment sans
utilité dans le cas où le module n'est pas \tfz.
Notons quand même que les \ids de Kaplansky présentent l'avantage sur
les \idfs de
permettre de \carar les \mtfsz.

Pour le deuxième point, voici ce qui se passe.
\\
Si $a$ est un \gtr typique de  $\cK_r(M)$
et si $M$ est engendré par  $(g_1,...,g_q)$,
on sait qu'il y a  $(h_1,...,h_r)$ dans  $M $
tels que $a M $ est contenu dans   $\gen{h_1,...,h_r}$.
\\
Une matrice de syzygies  pour le \sgr
   $(g_1,...,g_q,h_1,...,h_r)$
est alors de la forme suivante
$\cmatrix{   a \I_q \cr    B}$
avec    $B$  de format  $r\times q$.
On a simplement écrit que l'on  peut exprimer
    $a g_j$    en fonction des    $h_i$.
Donc dans l'\idf d'ordre  $r$  du module
il y a un \gtr typique qui est le \deter de  $a \I_q$  \cadz~$a^q$.
Ainsi, tout \gtr typique du Kaplansky est dans le nilradical
du Fitting correspondant. Notez que l'exposant qui intervient ici est
simplement le nombre de \gtrs du module.
\\
 Si maintenant  $a$  est un \gtr typique de  $\cF_r(M)$
on obtient  $a$  comme mineur  d'ordre  $q-r$  pour une
matrice de syzygies entre  $q$ \gtrs   $(g_1, ..., g_q)$.
Quitte à renuméroter les \gtrsz,
cette matrice peut s'écrire
\smashbot{$\cmatrix{N\cr D}$}
avec  $D$  carrée d'ordre  $q-r$,   $N$  de type  $r\times (q-r)$,
et   $\det  D  = a$.
\\
Par \colis des colonnes
(\prmt en faisant le produit à droite
 par la matrice cotransposée de $D$)
on obtient d'autres syzygies  pour les mêmes \gtrs sous la forme
$\cmatrix{N'
\cr a \I_{q-r}}$
et cela implique exactement que les  $q-r$ derniers \gtrs
multipliés par  $a$  tombent dans le module engendré par les
$r$ premiers.
Bref tout \gtr typique du Fitting est aussi un \gtr
typique du Kaplansky correspondant.

%%%%%%%%%%%%%%%%%%%%%%%%%%%%%%%%%%%%%%%%%%%%%%%%%%%%%%%%%%%%%%%%%%%%%%%%%%%
\exer{exoPetitsExemplesElim}

\emph {2.}
On a $\ff \cap\gA = \gen {a,b,c}$ (si $x \in \gA$ vérifie
$x \in \gen {X^d, a, b, c}_\AX$, faire $X := 0$),
et aussi~$\fb = \gen {a,b,c}^d$.  L'\id $\fa$ est le contenu en $T$ du
\pol $(a + bT + cT^2)^d$ tandis que $\fa'$ est le contenu en $\uT$ du
\pol $(aT_1 + bT_2 + cT_3)^d$. Par exemple pour $d = 2$:

\snic {
\fa = \gen {a^2, ab^2, 2ab, 2ac+b^2, b^3, b^2c, 2bc, c^2}, \quad
\fa' = \gen{  a^2, 2ab, 2ac, b^2, 2bc, c^2}.
}

%\sni
On a $\fa \subsetneq \fa' \subsetneq \fb \subsetneq \ff \cap \gA$ et $\fb =
(\ff\cap\gA)^d$. On voit aussi que $\fa$ n'est pas symétrique en $a$, $b$, $
c$. Toujours pour $d = 2$, on a $(\ff\cap \gA)^4 \subsetneq \fa$
 et $(\ff\cap \gA)^3 \not\subset \fa'$.  Pour $d$ quelconque, il
semblerait que $(\ff\cap \gA)^{3d-2} \subseteq \fa$.
%%%%%%%%%%%%%%%%%%%%%%%%%%%%%%%%%%%%%%%%%%%%%%%%%%%%%%%%%%%%%%%%%%%%%%%%%%%

\exer{exoRelateursViaElimIdeal} 
Soit $\wi\varphi : \gk[\uX,\uY] \to \kuX$ le morphisme d'\evn
 $Y_i \mt f_i$, l'anneau de base étant $\kuX$.  On~a:

\snic {
\ker\wi\varphi = \gen{Y_1 - f_1, \ldots, Y_s - f_s} =
\gen{g_1, \ldots, g_s},
}

%\sni
et puisque $\wi\varphi$ prolonge $\varphi$, $\ker\varphi =
\kuY \cap \ker\wi\varphi$, ce qu'il fallait démontrer.

%: sinotenglish
\sinotenglish{

%: correction exer   exoZerdimloc 
\exer{exoZerdimloc} \emph{1.} Puisque $\gen{e}=\gen{a^d}$ les filtres engendrés
par $e$ et $a$ sont les mêmes.

 \emph{2.}  $\gA[1/a]$ est isomorphe à $\aqo\gA{1-e}$ et $\gen{1-e}$
est un sous-\kev \tf de $\gA$.

 \emph{3a.} On a dans $\kX$ l'\egt $f=hf_0$ avec $f_0$ divise une puissance de $g$ et $\gen{h,g}=1$. Donc dans $\gA[1/g(x)]$, $h(x)=0$ et dans $\aqo\kX h$, $g(x)$ et $f_0(x)$ sont \ivsz.
Ceci montre que $\gA[1/g(x)]$ est isomorphe à $\aqo\kX h$.

 \emph{3b.} Mêmes notations avec $g=f'$ et $f_1=h$. Une relation 
 $f_1u+f'v=1$ donne l'\egt $f_1(u+f'_0v)+f'_1(f_0v)=1$, donc $f_1$ est bien \splz. \\
 Notez que si $f$ se factorise sous forme $\prod_k(x-a_k)^{m_k}$,
 alors $(x-a_k)$ divise $f'$ \ssi $m_k\geq 2$ (ceci indépendamment de la \cara du \cdiz~$\gk$), donc $f_1$ est le produit des $(x-a_k)$ qui ne figurent pas au carré dans $f$.

%: correction exer exoFittAnn
\exer{exoFittAnn}
Soit $F\in\MM_{q,n}$ une \mpn de $M$ pour un \sgr $(\xq)$.
Dire que $a\in\Ann_\gA(M)$, c'est dire que les $ax_i$ sont nuls, \cad que la matrice $a\,\I_q$ a pour colonnes des \colis des colonnes de $F$, donc que la matrice $F$ peut être élargie avec la matrice $a\,\I_q$, sans changer ses \iddsz. 
\\
Avec cette nouvelle \mpn \smashtop{$F'\,=\;$\blocs{1.3 }{.8 }{.8 }{0}{$F$}{$a\,\I_q$}{}{}\,}, il est clair qu'en multipliant un mineur d'ordre $r<q$ par $a$ on obtient un mineur d'ordre~\hbox{$r+1$}.
\\
Donc avec $k=q-r$, on obtient l'inclusion $a\,\cF_k(M)\subseteq \cF_{k-1}(M)$.

\smallskip \rem Puisque $\cF_q(M)=\gen{1}$, ceci donne une nouvelle \dem
de l'inclusion: $\Ann(M)^{q}\subseteq \cF_0(M)$.
\eoe

%: correction exer   exoQuentelZerdimFrac 
\exer{exoQuentelZerdimFrac}
On voit tout de suite que le point \emph{1} signifie:
$$\preskip.3em \postskip.3em 
\forall x\in\gA,\,\exists b\in\gA,\,z\in\Reg(\gA)\hbox{ 
tels que } x(z-bx)=0 
$$
Donc \emph{1} $\Leftrightarrow$ \emph{2.} Pour l'implication directe, prendre $y=z-bx$. 

\emph{2} $\Rightarrow$ \emph{3b}. Prendre $\fa=\gen{y}$.
 
\emph{1} $\Rightarrow$ \emph{3a.}  Car dans $\gK$ un \itf fidèle est égal à $\gen{1}$. 

\emph{3} $\Rightarrow$
\emph{1.} Prendre pour $z$ un \elt \ndz de
$\gen{x}+\fa$. 

\exer{exoReduiteFrob}
 \emph{(Forme réduite de Frobenius pour un \endo d'un \Kev de dimension finie, $\gK$ corps discret)}
 
\emph{1.}  Laissé \alecz. Voir aussi \paref{lemPrincipeIdentitesAlgebriques}.

\emph{2.}  On obtient les résultats précis suivants.
\begin{itemize}
\item La réduction de Smith de la matrice  
$X\In-F$ est du type
$$\preskip.3em \postskip.3em
\!\!\!\!\! {\, L\, (X\In-F) \, C =\Diag(1,\dots,1,f_1,\dots,f_k),\, k\in\NN\etl,\,L,\,C\in\GL_n(\KX)}
$$
avec pour $f_i$ des \pols \unts de degrés $>0$  vérifiant $f_1\mid\cdots\mid f_k$.
\item Le $\KX$-module $V_\varphi$ est somme directe de sous-espaces stables $\Kfi\cdot v_i$ avec $\nu_{v_i,\varphi}=f_i$. Il est isomorphe à 
$$\preskip.3em \postskip.3em 
\aqo\KX{f_1}\oplus\cdots\oplus \aqo\KX{f_k}.
$$
\item La matrice  $F$ est semblable à une {matrice diagonale par blocs}
dont les blocs diagonaux sont les matrices compagnes des \pols $f_i$.
Cette forme réduite de la matrice de $\varphi$ est appelée \emph{forme de Frobenius}.\index{forme de Frobenius!d'une matrice carrée sur un \cdiz}
\item Le \pol $f_k$ est égal au \polmin $\nu_\varphi$ de $\varphi$. Le \polcar $\chi_\varphi$ de $\varphi$ est égal au produit des $f_i$.
\item  Si $\chi_\varphi=\nu_\varphi$, alors
$V_\varphi=\Kfi\cdot y$ pour un $y\in V$,
 et la forme réduite de Frobenius de $F$ est la matrice compagne
de~$\chi_\varphi$.
\item  Considérons une \bdf $(g_1,\dots,g_r)$ pour les
invariants de similitude $(f_1,\dots,f_k)$ de $\varphi$, avec $\nu_{\varphi}=\prod_{i=1}^{r} g_i^{\ell_i}$.
Le lemme des noyaux donne la somme directe $V=\bigoplus_{i=1}^{r}\Ker g_i^{\ell_i}$, chaque sous-espace  $V_i=\Ker g_i^{\ell_i}$ est stable et si l'on note $\varphi_{i}$ la restriction de $\varphi_i$ à $V_i$, on a $\nu_{\varphi_i}=g_i^{\ell_i}$, et les invariants de similitude de $\varphi_i$ sont tous des puissances de $g_i$. 
\begin{itemize}
\item La description précédente reste valable
si l'on a une \dcn de $\nu_\varphi$ plus fine que celle donnée
par la \bdf des invariants de similitude.  
\item  Si $\ell_i=1$ pour un $i$,
tous les invariants de similitude de $\varphi_i$ sont égaux à $g_i$ 
et~$V_i$ est un
$\aqo\KX{g_i}$-module libre. 
\item Si $\gK$ est \acz, on retrouve la forme réduite de Jordan classique
et décomposant $\nu_\varphi$ en produit de facteurs du premier degré. 
\end{itemize}
\end{itemize}

%: correction exer   exoEndosemisimples 
\exer{exoEndosemisimples}
 \emph{(Endomorphismes semi-simples)} Suite de l'exercice \ref{exoReduiteFrob}.\\
On traite seulement le point \emph{4.}

On commence par calculer les invariants de similitude de~$\varphi$ et ceux de $\psi=\varphi\frt S$.
On calcule ensuite une \bdfz~\hbox{$(g_1,\dots,g_s)$} pour cette famille de \polsz. 
Comme $\nu_\psi$ divise $\nu_\varphi$, tous les $g_i$ divisent~$\nu_\varphi$.
\\
Si l'un des $g_i$ apparaît avec un exposant $>1$ dans $\nu_\varphi$, on a trouvé un facteur carré dans $\nu_\varphi$, on a terminé.
\\
Sinon, chacun des invariants de similitude de $\varphi$ et de $\psi$ est un
produit $\prod_{i\in J} g_i$ (sans exposant) pour une partie $J$ de $\lrbs$.
\\
On applique le lemme des noyaux 
correspondant à la \dcn   $\nu_\varphi=\prod_{i=1}^{s}g_i$.
Notons $K_i=\Ker g_i(\varphi)$, $S_i=S\cap K_i$, $\varphi_i:K_i\to K_i$   \hbox{et  $\psi_i:S_i\to S_i$} les  restrictions de $\varphi$, et  $\gK_i=\aqo\KX{g_i}=\gK[x_i]$ ($x_i$ est donc la classe de $X$ modulo~$g_i$).
\\
Il suffit de d'établir l'un des deux termes de l'alternative demandée pour chacun des 
triplets~\hbox{$(K_i,S_i,\varphi_i)$}. 
\\
La liste des invariants de similitude de $\varphi_i$ est formée de $v_i$ \pols tous égaux à $g_i$. Cela signifie que $K_i$ est un $\gK_i$-module libre de rang $v_i$.\\
La liste des invariants de similitude de $\psi_i$ est la même liste, mais plus courte, disons de longueur $u_i$. Le sous-espace  $S_i$ est un $\gK_i$-module libre de rang  $u_i$.\\
On regarde désormais $K_i$ sous la forme~$\gK_i^{v_i}$.\\
Si $u_i=0$ on pose $T_i=K_i$, et si $u_i=v_i$ on pose $T_i=0$.
En dehors de ces cas simples, on fait l'étude suivante. 
\\
Tout \elt $x\in K_i$ donne  un vecteur colonne de $\gK_i^{v_i}$.
Chacune des \coos ainsi obtenue est un \elt $h(x_i)$ de $\gK_i$.
En calculant une \bdf pour $h$ et $g_i$, on peut décider que l'une des alternatives
suivantes a certainement lieu:
\begin{itemize}
\item $h(x_i)=0$, OU
\item $h(x_i)$ est \ivz, OU 
\item  $g_i$ se décompose en un produit de plusieurs facteurs.
\end{itemize}
En bref, ou bien on découvre un facteur carré de $g_i$ (et l'\algo se termine),
ou bien $g_i$ se décompose en un produit de facteurs deux à deux étrangers, ce qui ramène le \pb à un \pb \gui{plus simple}
(les degrés des facteurs sont plus petits, et ne descendront jamais en dessous de~$1$), ou bien $\gK_i$ se comporte au cours du calcul comme s'il était un corps.\\ 
Ce que nous venons de dire à propos d'un \elt arbitraire de $K_i$,
nous l'appliquons pour une base de
$S_i$ comme $\gK_i$-module.
Nous traitons la matrice obtenue par la méthode du pivot, du moins si $\gK_i$
veut bien se comporter comme un corps au cours du calcul,
et nous calculons ainsi un \sul stable~$T_i$, isomorphe \hbox{à $\gK_i^{v_i-u_i}$}, de~$S_i$ dans $K_i$. 

\rem Pour plus de détails sur les deux exercices présédents
on peut consulter le chapitre 7 de \cite{DiLQ14}
.
\eoe
}
%: fin sinotenglish

%%%%%%%%%%%%%%%%%%%%%%%%%%%%%%%%%%%%%%%%%%%%%%%%%%%%%%%%%%%%%%%%%%%%%%%%%%%
%: problemes

\prob{exoDimZeroXcYbZa} 
D'abord, le cycle $\sigma = (1,2,3)$ réalise $\sigma(f_1) =
f_2$, $\sigma(f_2) = f_3$ et~\hbox{$\sigma(f_3) = f_1$}. Donc $C_3 = \gen
{\sigma}$ opère sur $\gA = \gk[x,y,z]$.
Si de plus $a=b$ ou $b=c$, alors $\{f_1, f_2, f_3\}$ est
invariant par  $\rS_3$.
Enfin, remarquons que l'origine est un zéro du \sysz, mais
aussi l'existence de solutions avec $x = y = z\neq 0$ (dans
une extension de $\gk$).

\emph {1.}
Il y a deux cas de figure: le cas $a \le b < c$, le plus facile
à étudier (cas~I), et le cas $a < b=c$ (cas~II).

$\bullet$ cas~I ($b < c$).
\\
On considère la relation d'ordre {\tt deglex} sur les \moms
de $\gk[X,Y,Z]$ (voir l'exercice \ref{exothSymEl}). 
Montrons que $\gA = \sum_{p,q,r:\max(p,q,r) < c} \gk\, x^py^qz^r$. \\
Soit $m = x^iy^j z^k$ vérifiant $\max(i,j,k) \ge c$.
Si $i \ge c$, on remplace dans $m$, $x^c$ \hbox{par $x^{i-c}x^c =
-x^{i-c}(y^b+z^a)$}. Même chose si $j \ge c$ ou si $k \ge c$. Il vient alors

\snic {
m = -(m_1 + m_2) \hbox { avec }  m_1, m_2 = \cases {
x^{i-c}y^{b+j}z^k,\ x^{i-c}y^jz^{a+k} & si $i \ge c$,\cr
x^{a+i}y^{j-c}z^k,\ x^iy^{j-c}z^{b+k} & si $j \ge c$,\cr
x^{b+i}y^jz^{k-c},\ x^iy^{a+j}z^{k-c} & si $k \ge c$.\cr
}}

%\sni
On voit alors que $m_1 < m$ et $m_2 < m$; on termine par \recuz.
\Llec vérifiera que les $x^py^qz^r$ avec $p$, $q$, $r < c$ forment
une $\gk$-base de $\gA$.  Pour ceux qui connaissent~: lorsque
$\gk$ est un \cdiz, $(f_1, f_2, f_3)$ est une base de Gr\"obner
pour la relation d'ordre {\tt deglex}.
Bilan: $\dim_\gk\gA = c^3$.

$\bullet$ cas~II ($a < b=c$). Ce cas est plus difficile. 
On suppose d'abord que $2$ est inversible dans~$\gk$. On introduit:
$$\preskip.4em \postskip.0em 
\arraycolsep2pt
\begin {array} {rcl}
g_1 &=& -f_1 + f_2 + f_3 = 2Z^c + X^a + Y^a - Z^a, \cr
g_2 &=& f_1 - f_2 + f_3 = 2X^c - X^a + Y^a + Z^a, \cr
g_3 &=& f_1 + f_2 - f_3  = 2Y^c + X^a - Y^a + Z^a .
\end {array} 
$$
On a alors:
$$\preskip.0em \postskip.4em 
2f_1 = g_2 + g_3, \quad  2f_2 = g_1 + g_3,\quad 2f_3 = g_1 + g_2, 
$$
de sorte que $\gen {f_1, f_2, f_3} = \gen {g_1, g_2, g_3}$.  On peut alors
opérer avec les $g_j$ comme on a fait avec les $f_i$ dans le cas~I. Si $\gk$ est un \cdiz,
$(g_1, g_2, g_3)$ est une base de Gr\"obner pour la relation d'ordre gradué
lexicographique {\tt deglex}.  \\
Bilan: $\dim_\gk\gA = c^3$ et les $x^py^qz^r$ avec
$p$, $q$, $r < c$ forment une $\gk$-base de $\gA$.  

 $\bullet$
Le cas II avec un \cdi $\gk$ de \cara $2$ est laissé à la sagacité
\dlecz.  L'anneau $\gA$ n'est pas toujours \zedz! 
Ceci arrive par exemple pour $\gk =
\FF_2$ et $(a,b) = (1,3)$, $(1,7)$, $(2,6)$, $(3,9)$.
Quand il est \zedz, il semble que $\dim_\gk\gA < c^3$.

\emph {2.}
Pour $(a,b,c) = (2,2,3)$, on sait que $\dim_\gk \gk[x,y,z] = 3^3 = 27$.  On
utilise le \tho de Stickelberger  \ref {thStickelberger}, sauf que
l'on ne connaît pas les zéros du \sysz.  On vérifie, à l'aide d'un \sys
de Calcul Formel, que le \polcar de $x$ sur~$\gk$ se factorise en
\pols \irds $(\gk = \QQ$):

\snic {
\rC{x} = t^8 (t+2) (t^3-t^2+1)^2  
       (t^4 - 2t^3 + 4t^2 - 6t + 4)  (t^4 + t^3 + t^2 - t + 2)^2
,}

%\sni

mais la \fcn de $\rC{x+2y}$ est du type $1^8 \cdot 1^1 \cdot 4^1
\cdot 4^1 \cdot 4^1 \cdot 6^1$. En conséquence, la \prn $(x,y,z) \mapsto
x$ ne sépare pas les zéros du \sysz, tandis que la \prn $(x,y,z) \mapsto
x+2y$ le fait. De plus, on voit que l'origine est le seul zéro
avec multiplicité (égale à $8$). Aidé de la \fcn
de $\rC{x}$ et en réalisant quelques petits calculs \sulsz,
on obtient:

\begin {itemize}
\item
Un (autre) zéro défini sur $\gk$, $(x,y,z) = (-2,-2,-2)$ et il est simple. 
\item
Si $\alpha$, $\beta$, $\gamma$ sont les trois 
racines distinctes de $t^3 - t^2 +1$, on obtient 6 zéros simples en faisant agir le groupe 
$\rS_3$ sur le zéro $(\alpha, \beta, \gamma)$. Si $s_1$, $s_2$, $s_3$ sont les fonctions \smqs \elrs de
$(X,Y,Z)$, alors, sur $\QQ$, on a l'\egt d'\ids $\gen {f_1, f_2, f_3, s_1-1} =
\gen {s_1-1, s_2,s_3+1}$, i.e. l'\alg de ces 6 zéros est l'\adu du \pol $t^3
- t^2 + 1$.
\item
Soit $\delta_i$ une racine de $t^4 + t^3 + t^2 - t + 2$ ($i\in\lrb{1..4}$).\\
En posant $y = x=\delta_i$ et
$z = 2/(x+1) = -(x^3 + x - 2)/2$, on obtient un zéro du
\sysz. Le \polmin de $z$ sur $\QQ$ est celui que l'on voit dans la
\fcn de $\rC{x}$: $t^4 - 2t^3 + 4t^2 - 6t + 4$.
On obtient ainsi quatre zéros simples du \sysz.
\item
On peut faire agir $\rA_3$ sur les quatre zéros précédents.
\end {itemize}
On a donc obtenu $1 + 6 + 3\times 4 = 19$ zéros simples et un zéro de
multiplicité~$8$. Le compte est bon.
\\
Remarque: alors que $\dim_\gk \gk[x,y,z] = 27$, on~a:

\snic {
\begin {array} {c}
\dim_\gk \gk[x] = \dim_\gk \gk[y] = \dim_\gk \gk[z] = 14,
\\[1mm]
\dim_\gk \gk[x,y] = \dim_\gk \gk[x,z] = \dim_\gk \gk[y,z] = 23.
\end {array}
}

%\sni
Ainsi, ni $\gk[x,y]$, ni $\gk[x,y,z]$ ne sont libres sur $\gk[x]$,
et $\gk[x,y,z]$ n'est pas libre sur~$\gk[x,y]$.

%\sni
\emph {3.}
Si $\gk$ est un \cdiz, dans le cas~I en \cara $\ne 2$, on trouve, de manière
expérimentale, que l'\alg locale de l'origine est $\aqo{\gk[X,Y,Z]}{X^a,
Y^a, Z^a}$ et donc la multiplicité de l'origine serait $a^3$.  Quant au cas
II, cela semble bien mystérieux.

%%%%%%%%%%%%%%%%%%%%%%%%%%%%%%%%%%%%%%%%%%%%%%%%%%%%%%%%%%%%%%%%%%%%%%%%%%%

\prob{exoIdealResultantGenerique}
\emph{1.}
On met les poids suivants sur $\gk[X]$: $X$ est de poids 1, et le
poids de $a_i$ et $b_{ji}$ est $i$. Ainsi $f$ et $g_j$ sont
\hmgs de poids $d$. On vérifie facilement pour tout $k \ge 0$
que $(X^k g_j) \bmod f$ est \hmg de poids $d+k$.

\emph{2.}
On indexe les $d$ lignes de $S$ par $1, \ldots, d$, la ligne $i$ correspondant
au poids $i$ via $i \leftrightarrow X^{d-i} \leftrightarrow a_i$. La matrice
$S$ est la concaténation horizontale de $r$ matrices carrées d'ordre $d$,
la $j$-ième matrice carrée étant celle de la multiplication par $g_j$
modulo $f$ dans la base $(X^{d-1}, \ldots, X, 1)$.  
Si l'on numérote par $(0, 1,
\ldots, d-1)$, les colonnes de la première sous-matrice carrée d'ordre $d$
de $S$ (correspondant à $g_1$), alors le \coe d'indice $(i,j)$ est
\hmg de poids $i+j$. Idem pour les autres \coes avec des conventions
analogues.\\
 Par exemple, pour $d = 3$, si $f = X^3 + a_1X^2 + a_2X + a_3$, $g =
b_1X^2 + b_2X + b_3$, la matrice de la multiplication par $g \bmod f$ est:
$$\preskip.3em \postskip.5em 
\arraycolsep.3em\mathrigid 2mu
\bordercmatrix [\lbrack\rbrack]{
                          &g    &Xg \bmod f      &X^2g \bmod f  \cr
X^{d-1}\leftrightarrow 1  &b_1  & -a_1b_1 + b_2  &a_1^2b_1 - a_1b_2 - a_2b_1 + b_3 \cr
X^{d-2}\leftrightarrow 2  &b_2  & -a_2b_1 + b_3  &a_1a_2b_1 - a_2b_2 - a_3b_1 \cr
X^{d-3}\leftrightarrow 3  &b_3  & -a_3b_1        &a_1a_3b_1 - a_3b_2 \cr
}
\; \hbox { de poids } \;
\cmatrix {1 & 2 & 3\cr 2 & 3 & 4\cr 3 &4 &5 \cr}. 
$$
Soit $M$ une sous-matrice d'ordre $d$ de $S$, $(k_1, \ldots, k_d)$ les
exposants de $X$ correspondant à ses colonnes ($k_i \in\lrb{0.. d-1}$,
et les colonnes sont des $X^{k_i} g_j \bmod f$).\\
Alors, $\det(M)$ est \hmgz, et son poids est la somme des poids des \coes
diagonaux, \cad

\snic {
(1 + k_1) + (2 + k_2) + \cdots + (d + k_d) = d(d+1)/2 + \sum_{i=1}^d k_i.
}

%\sni
Par exemple, le poids du premier mineur d'ordre $d$ de $S$ (correspondant
à la multiplication par $g_1$) est $d(d+1)/2 + \sum_{k=0}^{d-1}k = d^2$.  
\\
Le
poids de chacun des $rd\choose d$ mineurs est minoré par $d(d+1)/2$ (borne
obtenue pour $k_i=0$) et majoré par~$ d(3d-1)/2$, (borne obtenue pour
$k_i = d-1$). Ces bornes sont atteintes si $r \ge d$.

\emph{3.}
Le nombre $\dim_\QQ E$ minore le cardinal de n'importe quel
\sgr de~$\fb$. On trouve de manière expérimentale, pour des petites
valeurs de $r$ et $d$, que~$\dim_\QQ E = r^d$.  Mais on a mieux. En effet, la
considération d'objets gradués permet d'affirmer le résultat suivant
(Nakayama \hmgz, \pb \ref{exoNakayamaHomogeneRegularSequence}): toute famille graduée de $\fb$ dont l'image dans $E$ est
un \sgr\hmg du $\QQ$-\evc gradué $E$ est un
\sgr (\hmgz) de $\fb$. En~particulier, il existe un \sgr\hmg de $\fb$
de cardinal $\dim_\QQ E$, de manière conjecturale, $r^d$. On peut
aller plus loin en examinant les poids des \sgrs\hmgs minimaux de $\fb$:
ceux-ci sont uniques et fournies par la série (finie) du 
$\QQ$-\evc gradué $E$. Par exemple, pour $d=5$, $r=2$, cette
série est

\snic {
6t^{25} + 4t^{24} + 6t^{23} + 6t^{22} + 6t^{21} + 2t^{20} + 2t^{19},
}

%\sni
ce qui signifie que dans n'importe 
quel \sgr\hmg minimal de~$\fb$,
il y a exactement 6 \pols de poids~$25$, 
4 \pols de poids~$24$,
\dots, 2 \pols de poids 19 (avec $6 + 4 + \cdots + 2 = 32 = 2^5 = r^d$).
Dans cet exemple, le nombre ${rd \choose d}$ de mineurs d'ordre
$d$ de $S$ est $252$.

De manière conjecturale, il semblerait que $\fb$ soit engendré
par des \pols\hmgs de poids $\le d^2$, avec $d+r-1 \choose r-1$
\pols de poids $d^2$ exactement.

%: prob   exoNakayamaHomogeneRegularSequence
\prob{exoNakayamaHomogeneRegularSequence} ~\\
\emph {1.}
On fait une \dem par \recu sur $n$. 
\\ 
\emph{Cas $n=0$}:  résultat trivial.
%\emph{Cas $n=1$}:  suite réduite à un \elt
%\ndz $a$. Soit $\lambda_d a^d + \cdots + \lambda_1 a + \lambda_0 = 0$ avec
%$\lambda_i \in \gk$. On a $\lambda_0 \in \gk \cap \gen {a} = \{0\}$, donc
%$(\lambda_d a^{d-1} + \cdots + \lambda_1) a = 0$. On simplifie par $a$
%\ndz et on termine par \recu sur $d$.  
\\
\emph{Pour $n \ge 1$}, on considère
$\gA' = \aqo{\gA}{a_1}$. On a $\gk \hookrightarrow \gA'$ car $\gk \cap \gen
{a_1} = \{0\}$. La suite $(\overline {a_2}, \cdots, \overline {a_n})$ dans $\gA'$
vérifie les bonnes hypothèses pour la \recu sur $n$.  Supposons $f(a_1,
\ldots, a_n) = 0$ avec $f\in\kXn$ et $\deg_{X_1}(f)\le d$. \\
On
écrit $f = X_1q(X_1, \ldots, X_n) + r(X_2, \ldots, X_n)$ avec
$q$, $r$ à \coes dans $\gk$ et~$q$ de degré $\le d-1$ en $X_1$. 
Dans $\gA'$, on a
$r(\overline {a_2}, \ldots, \overline {a_n}) = 0$. 
Par \recu sur~$n$, on  a $r=0$.
Puisque $a_1$ est \ndzz, $q(a_1, \ldots, a_n)= 0$. Par
\recu sur $d$, on obtient $q=0$, donc $f=0$.

\emph {2a.}
Par \dfnz, $\gA_+E \subseteq E_1 \oplus E_2 \oplus \dots$; et comme $\gA_+E = E$,
c'est que $E_0 = 0$. Alors $\gA_+E \subseteq E_2 \oplus E_3 \oplus \dots$, et en
utilisant de nouveau $\gA_+E = E$, il vient $E_1 = 0$. Et ainsi de suite,
$E_n = 0$ pour tout $n$, donc $E = 0$.

\emph {2b.}
Soit $F$ le sous-\Amo de $E$ engendré par les $e_i$. C'est un sous-module
gradué car les $e_i$ sont \hmgsz. L'hypothèse équivaut à
$F + \gA_+ E = E$ ou encore $\gA_+(E/F) = E/F$. D'après la question \emph {2a},
on a $E/F = 0$ i.e. $E = F$: les~$e_i$ engendrent le \Amo~$E$.

\emph {3a.}
Il est clair que $\gA_0 = \gB_0$ et  $\fb = \gA_+ \gB$. En appliquant la
question précédente \hbox{au \Amoz} gradué $\gB$ et aux $e_i$,
on obtient que les $e_i$ forment un \sgr du \Amo $\gB$.

\emph {3b.}
Notons $S = \sum_i \gB_0 e_i$ (en fait, c'est une somme directe).\\ 
Montrons que $\gen {h_1, \ldots, h_d} \cap S = \{0\}$. Si $s = \sum_i \lambda_i
e_i \in \gen {h_1, \ldots, h_d}$ avec $\lambda_i \in \gB_0$, alors en
réduisant modulo $\gen {h_1, \ldots, h_d}$, il vient $\sum_i \lambda_i
\overline {e_i} = 0$ donc $\lambda_i = 0$ pour tout $i$ et $s = 0$.
\\
Pour $\alpha = (\alpha_1, \ldots, \alpha_d) \in \NN^d$, notons $h^\alpha =
h_1^{\alpha_1} \cdots h_d^{\alpha_d}$. Montrons que:
$$\preskip.4em \postskip.4em
\som_\alpha s_\alpha h^\alpha = 0 \hbox { avec } s_\alpha \in S
\;\Longrightarrow\; s_\alpha = 0 \;\; \hbox { pour tout }\alpha.
\leqno (\star)
$$
Pour cela, on va prouver par \recu descendante sur $i$, que:

\snic {
\bigl( f \in S[X_i, \ldots, X_d] \hbox { et }
f(h_i, \ldots, h_d) \equiv 0 \bmod \gen {h_1, \ldots, h_{i-1}} \bigr)
\;\Longrightarrow\; f = 0.
}

%\sni
\emph{D'abord pour $i = d$}. L'hypothèse est $s_m h_d^m + \cdots + s_1 h_d +
s_0 \equiv 0 \bmod \gen {h_1, \ldots, h_{d-1}}$ et on veut $s_k = 0$ pour tout
$k$.  On a $s_0 \in S \cap \gen {h_1, \ldots, h_d} = \{0\}$. On peut
simplifier la congruence par $h_d$ (qui est \ndz modulo $\gen {h_1, \ldots,
h_{d-1}}$) pour obtenir $s_m h_d^{m-1} + \cdots + s_1 \equiv 0 \bmod \gen {h_1,
\ldots, h_{d-1}}$.  En itérant le procédé, on obtient que tous les $s_k$
sont nuls.
\\
\emph{Passons de $i+1$ à $i$}. \\
Soit $f \in S[X_i, \ldots, X_d]$ de degré
$\le m$ avec $f(h_i, \ldots, h_d) \equiv 0 \bmod \gen {h_1, \ldots,
h_{i-1}}$. On écrit $f = X_i q(X_i, \ldots X_d) + r(X_{i+1}, \ldots, X_d)$
avec $q$, $r$ à \coes dans~$S$ et~$q$ de degré $\le m-1$.  On a donc
$r(h_{i+1}, \ldots, h_d) \equiv 0 \bmod \gen {h_1, \ldots, h_i}$, d'où par
\recu sur $i$, $r = 0$. On peut  simplifier la congruence par $h_i$ (qui
est \ndz modulo $\gen {h_1, \ldots, h_{i-1}}$) pour obtenir $q(h_i, \ldots,
h_d) \equiv 0 \bmod \gen {h_1, \ldots, h_{i-1}}$. \hbox{Donc $q = 0$} par \recu sur
$m$, puis $f = 0$.
\\
Bilan: on a donc le résultat pour $i = 1$ et ce résultat n'est
autre que $(\star)$.
\\
Une fois $(\star)$ prouvé, on peut montrer que les $e_i$ sont \lint
indépendants sur~$\gA$. Soit $\sum_i a_i e_i = 0$ avec $a_i \in \gA$; on
écrit $a_i = \sum_{\alpha} \lambda_{\alpha,i} h^\alpha$ et
$$\preskip.4em \postskip.4em
\som_i a_ie_i = \som_{i,\alpha} \lambda_{i,\alpha} h^\alpha e_i =
\som_{\alpha} s_\alpha h^\alpha 
\quad \hbox {avec} \quad s_\alpha = \som_i \lambda_{i,\alpha} e_i \in S.
$$
%\sni
Donc $s_\alpha = 0$ pour tout $\alpha$, puis $\lambda_{i,\alpha} = 0$ pour 
tout $i$, et  $a_i = 0$.

\emph {4.}
De manière \gnlez, si $(a_1, \ldots, a_d)$ est une \seqreg
d'un anneau $\gA$, elle est $L$-\ndze pour tout $\gA$-module 
libre $L$ (laissé \alecz). On applique ceci à
l'anneau $\gA = \gB_0[h_1, \ldots, h_d]$, à la suite
$(h_1, \ldots, h_d)$ (qui est bien une \seqreg de~$\gA$)
et à $L = \gB$ (qui est un \Amo libre par hypothèse).

%%%%%%%%%%%%%%%%%%%%%%%%%%%%%%%%%%%%%%%%%%%%%%%%%%%%%%%%%%%%%%%%%%%%%%%%%%%

% fin des solutions d'exos

%%%%%%%%%%%%%%%%%%%%%%%%%%%%%%%%%%%%%%%%%
%\exer{}{
%}

%:  ---- Section*{references}-----------
\Biblio
Bourbaki (Algèbre, chapitre X, ou Algèbre commutative chapitre I) appelle module \emph{pseudo cohérent} ce que nous appelons \como (conformément à l'usage le plus répandu, notamment dans la littérature anglaise), et \emph{\comoz} ce que nous appelons \como \pfz. Ceci est naturellement à relier aux \gui{Faisceaux Algébriques Cohérents} \hbox{de J.-P. Serre} (précurseurs des faisceaux de modules sur un schéma de Grothendieck) qui sont localement donnés par des modules \pf \cohsz. 
Signalons aussi que \cite{Stacks} adopte la \dfn de Bourbaki pour les modules \cohsz.

\smallskip 
Le \thrf{prop unicyc} est recopié de \cite{MRR} chap.\ V, th.\ 2.4.
Le \thrf{prop quot non iso} est recopié de \cite{MRR}  chap.\ III, exercice 9 p.\ 80.

La référence standard pour les \idfs est \cite{Nor}.

Pour ce qui concerne les structures \agqs purement équationnelles et l'\alg
\uvle on peut consulter~\cite{BuSa}. 

Une première introduction aux catégories
se trouve dans~\cite{Cohn}. 
\\ Des livres consacrés au sujet que l'on peut
recommander sont \cite{MACL} et~\cite{LaRo}.

Les \ids de Kaplansky d'un module $M$ étudiés dans l'exercice \ref{exoAutresIdF} sont utilisés dans~\cite[Chap.~IV]{Kun}
et~\cite[Chap. 9]{IRa}.

Les anneaux de Bézout strict (exercice
\ref {exoAnneauBézoutStrict}) et les anneaux de Smith 
(exercice~\ref{exoSmith}) ont été étudiés par Kaplansky 
dans \cite{Kap} dans un cadre plus \gnl d'anneaux non \ncrt commutatifs. Il les appelle respectivement des \gui{Hermite rings} et
des \gui{elementary divisor rings}. Mais cette terminologie n'est pas fixée.
Dans \cite{Lam06}, où l'exercice
\ref {exoAnneauBézoutStrict} trouve sa source, Lam utilise~\hbox{\emph{K-Hermite ring}} pour anneau de Bézout strict. Cela est
à distinguer de \emph{Hermite ring}: aujourd'hui un anneau $\gA$ est appelé 
\emph{anneau de Hermite} si tout \Amo stablement libre est libre, \cad encore
si tout \vmd %de $\Ae n$ 
est complétable
(voir chapitre \ref{chap ptf0}, section~\ref{ModStabLibre}).
Quant aux \gui{diviseurs \elrsz} ils sont désormais souvent utilisés dans
un sens plus particulier. Par exemple, il est fréquent de trouver écrit
que le \ZZmo $$\aqo\ZZ{900}\oplus\aqo\ZZ{10}\simeq \aqo\ZZ{25}\oplus\aqo\ZZ{5}\oplus\aqo\ZZ{4}\oplus\aqo\ZZ{2}\oplus\aqo\ZZ{9}$$ admet pour facteurs invariants   
la liste $(10,900)$ et pour diviseurs \elrs la liste non ordonnée $(25,5,4,2,9)$.%
\index{anneau!de Hermite}

L'exercice \ref{exoCalculT(M)} nous a été communiqué par Thierry Coquand.

\newpage \thispagestyle{CMcadreseul}
\incrementeexosetprob
%:        %%%%%%%%%%%%%%%%%%%%%%%%%%%%%%%%%%%%
%:        %%%%%%%%%%%%%%%%%%%%%%%%%%%%%%%%%%%%
%---- Chapitre  {Modules \ptf, 1}------------
\chapter{Modules \ptfsz, 1}
\label{chap ptf0}\relax
%--------------------
\minitoc

\section{Introduction}
%-----------------------------------------

Rappelons qu'un \mptf est un module isomorphe
à un facteur direct dans un \Amo
libre de rang fini. Cette notion s'avère être la \gnn naturelle,
pour les modules sur un anneau commutatif,
de la notion d'\evc de dimension finie sur un \cdiz.
Ce chapitre développe la théorie de base de ces modules.

\smallskip  
Une des motivations initiales de ce livre était de comprendre \emph{en termes
concrets} les \thos suivants concernant les \mptfsz.

%: --- Theorem{th.ptf.loc}--------
\begin{theorem}\label{th.ptf.loc}\relax {\em  (\Tho de structure locale des \mptfsz)}
Un \Amo $P$    est \ptf \ssi il est {\em  localement libre}
au sens suivant. Il existe des \ecoz~$s_1$,~$\ldots$,~$s_{\ell}$ dans~$\gA$ tels que
les modules~$P_{s_i}$ obtenus à partir de $P$ en étendant 
les scalaires aux anneaux~$\gA_{s_i}=\gA[1/s_i]$ sont libres.
\end{theorem}
%----------------------

%: --- Theorem{th.ptf.Fitting} ---
\begin{theorem}\label{th.ptf.Fitting}\relax {\em  (\Carn des \mptfs par
leurs idéaux de Fitting)} Un \Amo \pf est
\pro \ssi ses \idfs sont (des \idps engendrés par des)
idempotents.
\end{theorem}
%----------------------

%H  ''A-module'' remplace ''module sur A'' dans 3 des theoremes

%: --- Theorem{th.ptf.idpt}-------
\begin{theorem}\label{th.ptf.idpt}\relax {\em  (Décomposition d'un \mptf en somme
directe de modules de rang constant)} Si $P$ est un \Amo \ptf engendré par $n$
\eltsz, il existe un \sfio $(r_0, r_1, \ldots, r_n)$ (certains
éventuellement nuls) tel que chaque  $r_kP$  soit un module projectif de rang
$k$ sur l'anneau~$\aqo{\gA}{1-r_k}$. Alors, $P=\bigoplus_{k>0}r_kP$ et 
$\Ann(P)=\gen{r_0}$.
\end{theorem}
%----------------------

Dans cette somme directe on peut naturellement se limiter aux indices $k>0$ tels
que $r_k\neq 0$.

%: --- Theorem{th.ptf.plat} ---
\begin{theorem}\label{th.ptf.plat}\relax {\em  (\Carn des \mptfs par la
platitude)} Un \Amo \pf  est projectif \ssi il est
plat.
\end{theorem}
%----------------------

%H le commentaire suivant a changé
Dans ce chapitre nous démontrerons les trois premiers de ces \thosz.
Ils seront repris avec de nouvelles \dems dans le chapitre~\ref{chap ptf1}.
Le quatrième sera démontré dans le  
chapitre~\ref{chap mod plats} consacré
aux modules plats.

D'autres \thos importants concernant les \mptfs seront démontrés
dans les chapitres \ref{chap ptf1}, \ref{chapNbGtrs} et \ref{ChapMPEtendus}.
La théorie des \algs qui sont des \mptfs (nous les appelons des \asfsz) est développée dans le chapitre~\ref{chap AlgStricFi}.

%--- subsection{subsecPropCarPTFS}--------
\section{Généralités}
\label{subsecPropCarPTFS}
%-----------------------------------------

Rappelons qu'un \mptf est \pf (exemple 2, \paref{exl1pf}).

%:    subsection{subsecPropCarPTFS}--------
\subsec{Propriétés caractéristiques}

Lorsque $M$ et $N$ sont deux \Amosz, on a une \Ali naturelle
$\theta_{M,N}:M\sta\te
N\to \Lin_\gA(M,N)$ donnée par
%H j'ai numéroté ceci
\begin{equation}
\label{NOTAthetaMN}
\theta_{M,N}(\alpha \te y)= \big(x\mapsto \alpha (x)y\big)
\end{equation}

On note aussi $\theta_M$ pour $\theta_{M,M}$.\label{NOTAthetaM}

\medskip \rem On note parfois $\alpha \te y$  pour $\theta_{M,N}(\alpha \te y)$
mais ce n'est certai\-nement pas recommandé
lorsque $\theta_{M,N}$ n'est pas injective.
\eoe

\medskip Le \tho suivant donne quelques \prts immédiatement
équi\-valentes.
%: --- Theorem{propdef ptf}-------
\begin{theorem} {\em  (Modules \ptfsz)}\label{propdef ptf}\\
Pour un  \Amo  $P$,  \propeq
\begin{enumerate}
\item  [$(a)$] $P$  est un \ixx{module}{projectif de type fini}, i.e. il existe un
entier $n$, un \Amo $N$ et un \iso de $P\oplus N$ sur~$\Ae n$.
\index{projectif de type fini!module ---}

\item  [$(b1)$] Il existe un entier  $n$,
des \elts  $(g_i)_{i\in\lrbn}$  de
$P$  et des formes \lins  $(\alpha_i)_{i\in\lrbn}$  sur  $P$  telles que
pour tout $x\in P,   \;    x = \sum_i \alpha_i (x)\,  g_i$.

\item  [$(b2)$] $P$  est de type fini, et pour tout \sys fini de \gtrs
$(h_i)_{i\in\lrbm}$  de  $P$  il existe des formes \lins
$(\beta_i)_{i\in\lrbm}$  sur  $P$  telles que pour 
tout~$x\in P,   \;    x = \sum_i  \beta_i (x) \, h_i$.

\item  [$(b3)$] L'image de $P\sta\te_\gA P$ dans $\Lin_\gA(P,P)$ par l'\homo
cano\-nique~$\theta_P$  contient $\Id_P$.

\item  [$(c1)$] Il existe un entier  $n$  et deux \alis   
$\varphi : P\rightarrow\Ae n$  
\linebreak et~$\psi : \Ae n\rightarrow P$,  telles que  $\psi\circ \varphi= \Id_P$. On a alors~$\Ae n=\Im(\varphi )\oplus\Ker(\psi)$  et~$P\simeq\Im(\varphi\circ\psi )$.

\item  [$(c2)$] Le module $P$  est \tfz, et pour toute \ali 
surjec-\linebreak tive   
$\psi :\Ae m\rightarrow P$, il existe une \aliz~$\varphi : P\rightarrow \Ae m$  telle que~$\psi\circ \varphi= \Id_P$.  
On a alors~$\Ae m=\Im(\varphi )\oplus\Ker(\psi)$  et~$P\simeq\Im(\varphi\circ\psi )$.

\item  [$(c3)$] Comme $(c2)$ mais en remplaçant $\Ae m$ par un \Amo $M$
arbitraire: le module $P$  est \tfz, et pour toute \ali surjec-\linebreak tive $\psi :M\rightarrow P$,
$\Ker(\psi)$ est facteur direct.

\item  [$(c4)$]
Le module $P$  est \tf et le foncteur $\Lin_\gA(P,\bullet)$ transforme les \alis
surjectives en applications surjectives. \\
Autrement dit: 
pour tous \Amos $M$, $N$, pour toute \ali surjective   
$\psi : M\rightarrow N$ et
toute \ali $\Phi:P\rightarrow N$, il existe une \ali $\varphi : P\rightarrow M$
telle que
$\psi\circ \varphi= \Phi$.

\snic{
\xymatrix {
                                         & M\ar@{>>}[d]^{\psi} \\
P\ar@{-->}[ur]^{\varphi} \ar[r]_{\Phi} & N \\
}}

\end{enumerate}
\end{theorem}
%---------
%----begin{proof------------------
\begin{proof}
Le point $(b1)$ (resp. $(b2)$) n'est qu'une reformulation de $(c1)$ (resp. $(c2)$).
\\
Le point $(b3)$ n'est qu'une reformulation de~$(b1)$.
\\
On a trivialement $(c3) \Rightarrow (c2)\Rightarrow  (c1)$.

 $(a) \Rightarrow(c1)$ Considérer les applications canoniques 
 
\snic{P\rightarrow P\oplus N$ et  $P\oplus N\rightarrow P.}

%\sni
 $(c1) \Rightarrow(a)$ 
Considérer $\pi=\varphi\circ\psi$. On a $\pi^2=\pi$. 
Ceci définit une \prn de $\Ae n$ sur $\Im\pi=\Im \varphi\simeq P$ \paralm
à~$N=\Ker\pi=\Ker\psi$.

  $(b1)\Rightarrow (c4)$ Si $\Phi(g_i)=\psi(y_i)$ ($ i\in\lrbn$), on pose
$\varphi(x)=\sum \alpha_i(x)\,y_i$. On a alors pour tout $x\in P$:
$$\ndsp\Phi(x) =  \Phi\big(\sum \alpha_i(x)\,g_i\big) = \sum \alpha_i(x)\,
\psi(y_i) =\psi \big(\sum \alpha_i(x)\,y_i\big) = \psi \big(\varphi(x)\big).$$

 $(c4)\Rightarrow (c3)$ On prend $N=P$ et $\Phi=\Id_P$.
\end{proof}
%----end{proof------------------

On a aussi directement $(b1)\Rightarrow (b2)$ comme suit: en exprimant les $g_i$
comme \colis des $h_j$ on obtient les $\beta_j$ à partir des~$\alpha_i$.

 En pratique, conformément à la \dfn initiale, nous
considérerons un \mptf comme (copie par \iso  de l') image d'une
\mprnz~$F$. Une telle matrice, ou l'\ali qu'elle représente,
est encore appelée un \emph{projecteur}. Plus \gnltz, tout \endo
idempotent d'un module $M$ est appelé un \ix{projecteur}.

Lorsque l'on voit un \mptf selon la \dfn $(c1)$, la matrice de \prn est celle de
l'\ali $\varphi\circ\psi$. De même,
si l'on utilise la \dfn $(b1)$, la
\mprn est celle ayant pour \coes les $\alpha_i(g_j)$ en position $(i,j)$.

Un \sys $\big((g_1,\ldots,g_n),(\aln)\big)$ qui vérifie $(b1)$ est
appelé un \ixe{\sycz}{systeme de coordonnees} pour le module \pro $P$.
Certains auteurs parlent d'une \emph{base} du \mptfz, 
mais nous ne les suivrons~pas.

%--- Fact{factDualPTF}---------
\begin{fact}
\label{factDualPTF}\label{lemDual} \emph{(Dual d'un \mptfz, 1)}\\
Soit  $\big((g_1,\ldots,g_n),(\aln)\big)$  un \syc pour un \mptf
$P$.  Alors:
\begin{enumerate}
\item  [--]   les  $g_i$ engendrent $P$,
\item  [--]   les $\alpha_j$ engendrent $\Lin(P,\gA)=P\sta$,
\item  [--]   le module $P\sta$ est \ptfz,
\item  [--]   le module $(P\sta)\sta$ est canoniquement isomorphe à $P$,
\item  [--]   via cette identification canonique,
$\big((\aln),(g_1,\ldots,g_n)\big)$ est un \syc pour $P\sta$.
\end{enumerate}
En particulier, si $P$ est (isomorphe à) l'image d'une matrice de \prnz~$F$,
le module dual $P\sta$  est (isomorphe à) l'image de la matrice de \prnz~$\tra{F}$.
\end{fact}
%--- end-fact-----------------------------------------
%-----------------begin proof------------------
\begin{proof}
Le premier point est clair. Tout le reste est clair à partir du moment où on
montre que $\lambda  =\sum \lambda(g_i)\, \alpha_i $ pour tout $\lambda \in
P\sta$. Or cette \egt se démontre en évaluant
les deux membres en un \elt $x$
arbitraire de~$P$:

\snic{\lambda(x)= \lambda\bigl(\sum\, \alpha_i (x)\,  g_i\bigr)=
\sum\, \alpha_i (x)\,  \lambda(g_i)=
\bigl(\sum\, \lambda(g_i) \, \alpha_i \bigr)(x).}
\vspace{-.5em}\end{proof}
%-----------------end proof------------------

%: --- Theorem{prop pf ptf}---
\begin{theorem}
\label{prop pf ptf}\relax
Soit $\Ae m\vers{\psi} \gA^q\vers{\pi}P\to 0$ une \pn d'un module $P$.
Alors, $P$ est \ptf \ssi $\psi$ est \lnlz.
\end{theorem}
%--- end-Theorem-------------------
 Rappelons que \gui{$\psi$ est \lnlz} signifie qu'il
existe $\varphi :\gA^q\rightarrow \Ae m$
vérifiant $\psi\,\varphi\,\psi=\psi$. Par ailleurs, d'après le
\thrf{propIGCram} toute \ali qui a un rang au sens de la \dfnz~\ref{defRangk}
est \lnlz.

%-----------------begin proof------------------
\begin{proof} Si $\psi$ est \lnlz, le fait \ref{factInvGenCrois}
nous dit que $\Im\psi$ est facteur direct, et $\Coker\psi$
est isomorphe à un \sul de $\Im\psi$.
Réciproquement, si le module $P:=\Coker\psi$ est projectif, on applique à la \prnz~$\pi:\gA^q\rightarrow P$ la \prtz~$(c2)$ du \thref{propdef ptf}.
On obtient~$\tau:P\rightarrow \gA^q$ avec~$\pi\circ\tau=\Id_P$, de sorte que
$\gA^q=\Im \tau \oplus \Im\psi$. Donc $\Im\psi$ est \ptf et l'on peut
appliquer à $\psi :\Ae m\rightarrow \Im\psi$ la \prt $(c2)$, ce qui
nous donne $\varphi$ sur la composante $\Im\psi$ (et l'on prend par exemple $0$
sur~$\Im \tau$).
\end{proof}
%-----------------end proof------------------

%:    subsec{Principe local-global}
\subsec{Principe local-global}

Le fait qu'un \Amo est \ptf est une notion locale au sens suivant.

%: --- Principe locglobplcc.cor.pf.ptf
\begin{plcc}
\label{plcc.cor.pf.ptf}\relax
\emph{(Modules \ptfsz)}\\
Soient $S_1$, $\ldots$, $S_n$ des \moco de $\gA$ et
 $P$ un \Amoz. 
 \\
 Si les~$P_{S_i}$ sont libres,~$P$ est \ptfz. 
 \\
 Plus \gnltz, le module $P$
est \ptf \ssi les  $P_{S_i}$ sont des $\gA_{S_i}$-\mptfsz.
\end{plcc}
%--- end-plcc-----------------------------------------
%-----------------begin proof------------------
\begin{proof}
Cela résulte du \thref{prop pf ptf}, du \plgref{plcc.pf} pour les \mpfs et du
\plgrf{fact.lnl.loc} pour les \alis \lnlsz.
\end{proof}
%-----------------end proof------------------

Le \plgref{plcc.cor.pf.ptf} établit l'implication \gui{si} dans le
\thref{th.ptf.loc}. La réciproque \gui{seulement si} a de
fait été démontrée au \thref{theoremIFD} ce qui
nous donnera le \thref{prop Fitt ptf 2}.
Nous  donnerons pour cette réciproque un énoncé plus précis et une \dem
plus conceptuelle avec le \thref{th ptf loc libre}.

%:    subsec{Modules projectifs et lemme de Schanuel}
\subsec{Modules projectifs et lemme de Schanuel}

La notion de module \pro peut être définie pour des modules
qui ne sont pas  \tfz.
Dans la suite nous utiliserons  rarement de tels modules, 
mais il est cependant utile de
donner quelques précisions sur ce~sujet.

%:     Definition{defiMPRO}
\begin{definition}\label{defiMPRO}
Un \Amo $P$ (non \ncrt \tfz) est dit \ixc{projectif}{module ---}
s'il vérifie la \prt suivante.
\\
 Pour tous \Amos $M,\,N$, pour toute \ali surjec-\linebreak tive
$\psi : M\rightarrow N$ et  toute \ali $\Phi:P\rightarrow N$, il existe une \ali $\varphi : P\rightarrow M$   telle que
$\psi\circ \varphi= \Phi$.

\snic{
\xymatrix {
                                       & M\ar@{>>}[d]^{\psi} \\
P\ar@{-->}[ur]^{\varphi} \ar[r]_{\Phi} & N \\
}}
\end{definition}
%

%:HHH rajout une phrase
Ainsi, vue la \carn (c4) dans le \thref{propdef ptf}, un \Amo est \ptf \ssi il est \pro et \tfz.
Dans le fait suivant, la dernière \prt est
comme l'implication  \hbox{$(c4)$ $\Rightarrow$  $(c3)$} dans ce \thoz.

%:     Fact{factdefiMPRO}
\begin{fact}\label{factdefiMPRO} ~
\begin{enumerate}
\item Un module libre ayant pour base un ensemble en bijection avec $\NN$ est \proz.
Par exemple l'anneau des \pols $\AX$ est un \Amo \proz.
\item Tout module en facteur direct dans un module \pro est \proz.
\item 
%:HHH enonce raccourci
Si $P$ est  \proz, toute suite exacte courte
$0\to N\to M\to P\to 0$ est scindée.
%Si $P$ est un module \proz, pour tout \Amo $M$ et toute \ali 
%surjective $\psi :M\rightarrow P$,
%$\Ker(\psi)$ est facteur direct. Autrement dit, toute suite exacte courte
%$0\to N\to M\to P\to 0$ est scindée.
\end{enumerate}
\end{fact}

\comm En \coma les modules libres \emph{ne sont pas toujours} \prosz.
En outre, il semble impossible de réaliser tout module
comme quotient d'un module libre et \proz.
De même il semble impossible de mettre
tout module \pro en facteur direct dans un module libre et \proz.
Pour plus de détails sur ce sujet on peut consulter 
l'exercice~\ref{propfreeplat} et le livre~\cite{MRR}.
\eoe
\perso{peut-être donner un exo? cela risque d'être démoralisant.}

%:     Lemma{lemScha}
\begin{lemma}\label{lemScha}
On considère deux \Alis surjectives de même image
$$\preskip-.2em \postskip.1em 
P_1  \vers{\varphi_1}   M  \rightarrow  0  \;\;\;\hbox{et} \;\;\;
   P_2  \vers{\varphi_2}   M  \rightarrow  0
$$
avec $P_1$ et $P_2$ \prosz.
\begin{enumerate}
\item  Il existe des \isos réciproques
$$\preskip.2em \postskip-.4em
\alpha ,\,  \beta:P_1\oplus P_2\to P_1\oplus P_2
$$ 
tels que
$$\preskip.0em \postskip-.2em
(\varphi_1\oplus 0_{P_2})\circ \alpha=0_{P_1} \oplus \varphi_2\;\hbox{
et}\; \varphi_1\oplus 0_{P_2}=(0_{P_1} \oplus \varphi_2)\circ \beta.
$$
\item Notons $K_1=\Ker\varphi_1$ et $K_2=\Ker\varphi_2$. 
Par restriction de~$\alpha$ et~$\beta$ on obtient des \isos réciproques entre
$K_1\oplus P_2$ et $P_1\oplus K_2$.
\end{enumerate}
%On considère deux \Alis surjectives de même image
%$ P_1  \vers{\varphi_1}   M  \rightarrow  0  $, $
%   P_2  \vers{\varphi_2}   M  \rightarrow  0$
%avec les modules $P_1$ et $P_2$ \prosz.
%\begin{enumerate}
%%
%\item  Il existe des \isos réciproques
%$\alpha$, $\beta:P_1\oplus P_2\to P_1\oplus P_2$ et
%tels que
%$(\varphi_1\oplus 0_{P_2})\circ \alpha=0_{P_1} \oplus \varphi_2$
%et  $\varphi_1\oplus 0_{P_2}=(0_{P_1} \oplus \varphi_2)\circ \beta$.
%%
%\item Si l'on note $K_1=\Ker\varphi_1$ et $K_2=\Ker\varphi_2$, on obtient
%par restriction de~$\alpha$ et~$\beta$ des \isos réciproques entre
%$K_1\oplus P_2$ et $P_1\oplus K_2$.
%\end{enumerate}
\end{lemma}
\begin{proof}
Il existe $u:P_1\to P_2$ tel que $\varphi_2\circ u=\varphi_1$ et
 $v:P_2\to P_1$ tel que~$\varphi_1\circ v=\varphi_2$.
 $$
\xymatrix @R = 0.3cm{
P_1 \ar@{-->}[dd]_{u} \ar[rd]^{\varphi_1} \\
                                        & M \\
P_2 \ar@{>>}[ru]_{\varphi_2} \\
}
\qquad\qquad
\xymatrix @R = 0.3cm{
P_1 \ar@{>>}[rd]^{\varphi_1} \\
                                        & M \\
P_2 \ar@{-->}[uu]_{v}\ar[ru]_{\varphi_2} \\
}
\qquad\qquad
\xymatrix @R = 0.4cm{
P_1 \oplus P_2 \ar@/^/@{-->}[dd]^{\beta}
                 \ar[rd]^{~\varphi_1 \oplus 0_{P_2}} \\
                                        & M \\
P_1 \oplus P_2 \ar@/^/@{-->}[uu]^{\alpha}
    \ar[ru]_{~0_{P_1} \oplus \varphi_2} \\
}
$$
On vérifie que
 $\alpha$ et $\beta$ définis par les matrices ci-dessous conviennent.
 $$
 \alpha=\cmatrix{\Id_{P_1}-vu&  v\cr -u& \Id_{P_2} }
\quad\quad
 \beta=\cmatrix{\Id_{P_1}& -v\cr u& \Id_{P_2}-uv }
  .$$
NB: la matrice $\beta$ est une variante sophistiquée de ce que serait
la matrice  cotransposée de  $\alpha$ si $\Id_{P_1}$, $\Id_{P_2}$, $u$ et $v$
étaient des scalaires.
\end{proof}

%

%:     Corollary{corlemScha}
\begin{corollary}\label{corlemScha} \emph{(Lemme de Schanuel)}
On considère deux  suites exactes:
\vspace{-1mm}
\[\preskip.2em \postskip.2em
\begin{array}{ccccccccc}
0 &\rightarrow& K_1& \vers{j_1} & P_1& \vers{\varphi_1} & M& \rightarrow& 0   \\
0 &\rightarrow& K_2& \vers{j_2} & P_2& \vers{\varphi_2} & M& \rightarrow& 0
 \end{array}
\]
avec les modules $P_1$ et $P_2$ \prosz. Alors, $K_1\oplus P_2\simeq K_2\oplus P_1.$
\end{corollary}

%:  subsec{Catégorie des mptfs}
\subsec{Catégorie des \mptfsz}
\label{secCatMptf}
%-----------------------------------------
\vspace{3pt}
%--- SUBSUBsection{Cart cat}---------------
\subsubsection*{Une construction purement catégorique}

La catégorie des \mptfs sur $\gA$ peut être construite à partir de la
catégorie des modules libres de rang fini sur $\gA$ par un procédé purement
catégorique.
%-----------------begin enum------------------
\begin{enumerate}
\item  Un \mptf $P$ est décrit par un couple $(\rL_P,\rPr_P)$ 
où~$\rL_P$ est un module libre
de rang fini et $\rPr_P\in\End(\rL_P)$ est un \prrz.
On a $P\simeq \Im\rPr_P\simeq \Coker(\Id_{\rL_P}-\rPr_P)$.
%La matrice de $\rPr_P$ est une matrice de \prn pour $P$.
\item  Une \ali $\varphi$ du module $P$ (décrit par $(\rL_P,\rPr_P)$) vers
le module  $Q$ (décrit par $(\rL_Q,\rPr_Q)$)
est décrite par une \ali
$\rL_\varphi:\rL_P\rightarrow \rL_Q$  soumise aux relations de commutation

\snic{\rPr_Q\circ\,\rL_\varphi=\rL_\varphi=\rL_\varphi\circ\rPr_P.}

%\sni
En d'autres termes $\rL_\varphi$ est nulle sur $\Ker(\rPr_P)$ et son image est
contenue dans $\Im(\rPr_Q)$.
\item  L'identité de $P$ est représentée par $\rL_{\Id_P}=\rPr_P$.
\item  La somme de deux \alis $\varphi$ et $\psi$ de $P$ vers $Q$
repré\-sentées par $\rL_\varphi$ et  $\rL_\psi$ est repré\-sentée par
$\rL_\varphi+\rL_\psi$.
L'\ali  $a\varphi$ est repré\-sentée par $a\rL_\varphi$.
\item  Pour représenter la composée de deux \alisz, on compose leurs
repré\-sen\-tations.
\item  Enfin, une \ali  $\varphi$ de $P$ vers $Q$ représentée par
$\rL_\varphi$ est nulle
\ssi  $\rL_\varphi=0$.
\end{enumerate}
%-----------------end enum------------------

Ceci montre que les problèmes concernant les \mptfs peuvent toujours être
interprétés comme des \pbs à propos de matrices de \prnz, et se ramènent
\label{exMptf} souvent à des problèmes de résolution de \slis sur $\gA$.

\medskip Une catégorie \eqvez, mieux adaptée aux calculs, est la
catégorie dont les objets sont les \mprns à \coes dans
$\gA$, un morphisme
de $F$ vers $G$ étant une matrice $H$ de format convenable 
vérifiant les \egts 

\snic{GH=H=HF.}

%--- SUBSUBsection{Sycs}---------------
\subsubsection*{Avec des \sycsz}

%:HHH rajout une phrase
Le fait suivant reprend les affirmations du paragraphe précédent
lorsque l'on prend le point de vue des \sycsz.

%:     Fact{factMatriceAlin}
\begin{fact}\label{factMatriceAlin}
Soient $P$ et $Q$ deux \mptfs avec des \sycs 

\snic{\big((\xn),(\aln)\big)$ et $\big((\ym),(\beta_1,\ldots,\beta_m)\big),}

%\sni
et soit $\varphi:P\to Q$ une \Aliz.\\
Alors, on peut coder $P$ et $Q$ par les matrices 

\snic{F\eqdefi \big(\alpha_i(x_j)\big)_{i,j\in\lrbn}\quad\et\quad G\eqdefi \big(\beta_i(y_j)\big)_{i,j\in\lrbm}.}

%\sni
Précisément, on a les \isos

\snic{
\begin{array}{rcl} 
\pi_1: P\to\Im F\,,  &   & x\mapsto \tra[\,\alpha_1(x) \;\cdots\; \alpha_n(x)\,],  
 \\[1mm] 
\pi_2: Q\to\Im G\,,  &   & y\mapsto \tra[\,\beta_1(y) \;\cdots\; \beta_m(y)\,].   \end{array}
}

%\sni
Quant à l'\ali $\varphi$, elle est codée par la matrice 

\snic{H\eqdefi \big(\beta_i(\varphi(x_j))\big)_{i\in\lrbn,j\in\lrbm}}

%\sni
qui vérifie $GH=H=HF$. La matrice $H$ est celle de l'\ali 

\snic{\Ae n\to \Ae m, \quad\pi_1(x)+z\;\mapsto\; \pi_2 \big(\varphi(x)\big)\quad$ si $ x\in P $ et $z\in\Ker F.}

%\sni
 On dira que \emph{la matrice $H$ représente l'\ali $\varphi$ dans les \sycs $\big((\ux),(\ual)\big)$
et $\big((\uy),(\und{\beta})\big)$}.%
\index{matrice!d'une appli@d'une \ali dans des \sycsz}
%:HHH   la phrase précédente a été rappatriée ici depuis plus loin

\end{fact}
%:HHH  suppression de :  \facile

%--- SUBSUBsection{Iso ptfs}---------------
\subsubsection*{Application: les \isos entre modules \pros \\de
 type~fini}

%:2018
%Le \tho suivant dit que, pour $F\in\GA_m(\gA)$ et  $G\in\GA_n(\gA)$, si 
Le lemme d'élargissement \ref{propIsoIm} dit que, pour $F\in\GA_m(\gA)$ et  $G\in\GA_n(\gA)$, si 
$\Im F$ et~$\Im G$ sont isomorphes, quitte à \gui{agrandir}
les matrices $F$ et $G$, elles peuvent être supposées semblables.

%:2018
%Dans le lemme qui suit, 
Dans la suite, 
nous employons la notation $\Diag(M_1,\dots,M_k)$ d'une manière plus large que celle utilisée jusqu'ici. Au lieu d'une liste d'\elts
de l'anneau, nous considérons 
pour $(M_1,\dots,M_k)$ une liste de matrices carrées. La matrice représentée ainsi
est usuellement appelée une \ixx{matrice}{diagonale par blocs}.\label{Notadiagblocs}

%--- lemme{propIsoIm}----------
\begin{lemma}
\label{propIsoIm} \emph{(Lemme \dlgz)}\\
On considère le codage matriciel de la catégorie des \mptfsz.
Si un \iso $\varphi$ de $\Im F$ sur $\Im G$ est codé par $U$ et son inverse codé par
$U',$ on obtient une matrice $A\in\EE_{n+m}(\gA)$
$$A=\cmatrix{\I_m-F&-U'\cr U&\In-G}=\cmatrix{\I_m&0 \cr U &\In }\,
\cmatrix{ \I_m&-U' \cr 0& \In}\,\cmatrix{\I_m &0 \cr U&\In },$$
avec
%------begin equation--eqpropIsoIm-----------
\begin{equation}\label{eqpropIsoIm}
\bloc{0_m}{0}{0}{G}=A\,\bloc{F}{0}{0}{0_n}\,A^{-1}.
\end{equation}
%---------------------end equation--------------
Réciproquement, une
conjugaison entre $\Diag(0_m,G)$ et $\Diag(F,0_n)$ fournit un \iso entre
$\Im F$ et $\Im G$.
\end{lemma}
%--- end-theorem----------------------------------------
%-----------------begin proof------------------
\begin{proof} La matrice suivante
$$\bordercmatrix [\lbrack\rbrack]{
       & \Im F & \Ker F   & \Im G & \Ker G \cr
\Im F &    0    &  0 &       -\varphi^{-1}   &  0\cr
\Ker F &    0    &    \Id   &    0  &  0\cr
\Im G &    \varphi     &    0   &    0  &  0\cr
\Ker G &    0    &    0   &    0  &  \Id
}  ,$$
une fois remplacé $\Im F \oplus \Ker F$ par $\Ae m$ et $\Im G \oplus \Ker G$
par $\Ae n$, donne la matrice $A$. La présence du signe $-$ est due
à la \dcn classique en produit de matrices \elrs
$$
\cmatrix{0&-a^{-1}\cr a&0} =\cmatrix{1&0\cr a&1} \,\cmatrix{1&-a^{-1}\cr 0&1}
\,\cmatrix{1&0\cr a&1}.
$$
\end{proof}
%-----------------end proof------------------

%  SUBsubsection{Quand l'image est libre}-
\subsubsection*{Quand l'image d'une matrice de \prn est libre}

Si un \prr $P\in\GAn(\gA)$  a pour image un module libre
de rang $r$, son noyau n'est
pas automatiquement libre, et la matrice n'est donc pas à tout coup semblable
à la matrice standard~$\I_{r,n}$.

Il est intéressant de savoir \carar simplement
le fait que l'image est~libre.

%--- Proposition{propImProjLib}------
\begin{proposition}
\label{propImProjLib} \emph{(Matrices de \prn dont l'image est libre)}
\\
Soit $P\in\Mn(\gA)$. La matrice
$P$ est \idme et d'image libre de rang~$r$  \ssi
il existe deux matrices $X\in\Ae {n\times r}$ et  $Y\in\Ae {r\times n}$ telles que~$YX=\I_r$ et~$P=XY$. On donne en outre les précisions suivantes.
%-----------------begin enum------------------
\begin{enumerate}
\item  $\Ker P=\Ker Y$, $\Im P=\Im X\simeq \Im Y$, et les colonnes de $X$ 
forment une base de~$\Im P$.
\item  Pour toutes matrices $X'$, $Y'$ de mêmes formats que $X$ et $Y$ et
telles que~$P=X'Y'$, il existe
une unique matrice~$U\in\GL_r(\gA)$ telle que 

\snic{X'=X\,U\,\hbox{ et }\,Y=U\,Y'.}

%\sni
En fait,
$U=YX'$, $U^{-1}=Y'X$ et $Y'X'=\I_r$.
\end{enumerate}
%-----------------end enum------------------
\end{proposition}
%--- end-proposition----------------------------------------
%-----------------begin proof------------------
\begin{proof}
Supposons que $P$ est \idme d'image libre de rang $r$.
Pour colonnes de $X$ on prend une base de $\Im P$. 
Alors, il existe une unique matrice $Y$ telle que $P=XY$. 
Puisque $PX=X$ (car $\Im X\subseteq\Im P$ et $P^2=P$), on obtient~$XYX=X$.\\
Puisque les colonnes de $X$ sont indépendantes et que $X(\I_r-YX)=0$, on obtient~$\I_r=YX$.
\\
Réciproquement, supposons $YX=\I_r$ et $P=XY$. 
Alors:

\snic{P^2=XYXY=X\I_rY=XY=P\;$ et $\;PX=XYX=X.}

%\sni
Donc $\Im P=\Im X$. En outre, les colonnes de $X$ sont
indépendantes car~$XZ=0$ implique $Z=YXZ=0$.
 
 \emph{1.} La suite $\Ae n\vvvers{\In-P}\Ae n\vers{Y}\Ae r$ est exacte. En effet,
$Y(\In-P)=0$, et si~$YZ=0$, alors $PZ=0$, donc $Z=(\In-P)Z$. Ainsi:

\snic{\Ker Y= \Im(\In-P)=\Ker P$, et $\Im Y\simeq
\Ae n\sur{\Ker Y}=\Ae n\sur{\Ker P\simeq \Im P}.}

\sni
\emph{2.} Si maintenant  $X'$, $Y'$ sont de mêmes formats que $X$, $Y$, et si $P=X'Y'$,
on pose $U=YX'$ et $V=Y'X$.
Alors: 
\begin{itemize}
\item $UV=YX'Y'X=YPX=YX=\I_r$,
\item  $X'V=X'Y'X=PX=X$, donc $X'=XU$,
\item  $UY'=YX'Y'=YP=Y$, donc $Y'=VY$.
\end{itemize}
 Enfin $Y'X'=VYXU=VU=\I_r$.
\end{proof}
%-----------------end proof------------------

%%%%%%%%%%%%%%%%%%%%%%%%%%%%%%%%%%%%%%%%%%%%%%%%%%%%%%%%%%%%%%%%%%%%%%%%%%%
%--- subsec Mptfs sur les anneaux \zeds
\section[Sur les anneaux \zedsz]{Modules \ptfs sur les anneaux \zedsz}
\label{secMPTFzed}
%-----------------------------------------

Le \tho suivant  généralise le
\thrf{thZerDimRedLib}.
\perso{les \algbs vérifient les points 2,3,4, mais les
\alos ne vérifient pas en \gnl les points 5 et 6}

%:     Theorem{propZerdimLib}
\begin{theorem}
\label{propZerdimLib}
Soit $\gA$ un anneau \zedz.
%-----------------begin enum------------------
\begin{enumerate}
\item \label{ite1propZerdimLib} Si $\gA$ est réduit  tout \mpf $M$ est quasi libre, et tout
sous-\mtf de $M$ est en facteur direct.

\item  \label{ite2propZerdimLib} \emph{(Lemme de la liberté \zedez)}\\
Tout \Amo \ptf  est quasi libre.

\item    \label{ite3propZerdimLib}
Toute matrice  $G\in\gA^{q\times m}$ de
rang $\geq k$ est \eqve à une matrice

\snic{
\cmatrix{
    \I_{k}   &0_{k,m-k}      \cr
    0_{q-k,k}&  G_1      }}

%\sni
avec $\cD_r(G_1)=\cD_{k+r}(G)$ pour
tout $r\geq 0$. En particulier, toute matrice de rang $k$ est \nlz.

\item    \label{ite4propZerdimLib}
Tout \mpf $M$ tel que $\cF_r(M)=\gen{1}$
(\cad \lot engendré par $r$ \eltsz, cf. la \dfn \ref{deflocgenk})
est engendré par $r$ \eltsz.

\item \label{ite5propZerdimLib}
\emph{(\Tho de la base incomplète)}\\
Si un sous-module $P$ d'un \mptf $Q$ est \ptfz,
il possède un \sulz.
Si $Q$ est libre de rang $q$ et $P$ libre de rang $p$, tout \sul est
libre de rang $q-p$.

\item  \label{ite6propZerdimLib}
Soit  $Q$  un \Amo \ptf et $\varphi :Q\to Q$ un \endoz. \Propeq
%-----------------begin enum------------------
\begin{enumerate}\itemsep0pt
\item $\varphi$ est injectif.
\item $\varphi$ est surjectif.
\item $\varphi$ est un \isoz.
\end{enumerate}
%-----------------end enum------------------
\end{enumerate}
%-----------------end enum------------------
\end{theorem}
%--- end-theorem----------------------------

\begin{proof}
Le point \emph{\ref{ite1propZerdimLib}.} est un rappel du \thrf{thZerDimRedLib}.

\emph{\ref{ite2propZerdimLib}.}
On considère une matrice de \pn $A$ du module et
on commence par remarquer que puisque le module est \proz,  $\cD_1(A)=\gen{e}$
avec $e$ \idmz. On peut considérer que la
première étape du calcul se fait au niveau de l'anneau $\gA_{e}=\gA[1/e]$.
On est ramené au cas où $\cD_1(A)=e=1$, ce que nous supposons désormais.
On applique le point \emph{\ref{ite3propZerdimLib}} avec $k=1$ et l'on
termine par \recuz.

Le point \emph{\ref{ite3propZerdimLib}} est une sorte de lemme du mineur \iv (\ref{lem.min.inv}) sans mineur \iv dans l'hypothèse.  
%:2012  le numero du lemme plutot que sa page 
%(voir \paref{lem.min.inv}).
On applique avec  l'anneau  $\Ared$ le point~\emph{1} du
\thrf{thZerDimRedLib}.
On obtient alors la matrice voulue, mais seulement modulo $\DA(0)$.
\\
On remarque que la matrice $\rI_k+R$ avec $R\in \Mk \big(\DA(0)\big)$ a un \deter \ivz, ce qui permet d'appliquer le lemme du mineur \ivz.

\emph{\ref{ite4propZerdimLib}.} Résulte du point \emph{\ref{ite3propZerdimLib}} appliqué à une matrice
de \pn du module.

\emph{\ref{ite5propZerdimLib}.}
Voyons d'abord le deuxième cas.  Considérons la matrice $G$
dont les vecteurs colonnes forment une base du sous-module $P$. 
Puisque $G$ est la matrice d'une \ali injective, 
son \idd d'ordre~$p$ est \ndzz, donc égal à
$\gen{1}$ (corolaire \ref{corZedReg}).
Il reste à appliquer le point~\emph{\ref{ite3propZerdimLib}}.
Dans le cas \gnlz, si $P$ est engendré par $p$ \eltsz, considérons un
module~$P'$ tel que $P\oplus P'\simeq\gA^p$.
Le module $Q\oplus P'$ est \ptfz, donc en facteur direct dans un
module  $L\simeq \Ae n$. Alors,
d'après le deuxième cas,~$P\oplus P'$ est en facteur direct dans $L$.
On en déduit que $P$ est l'image d'une \prn $\pi:L\to L$. Enfin, la restriction
de $\pi$ à $Q$ est une \prn qui a pour image $P$.

\emph{\ref{ite6propZerdimLib}}. 
On sait déjà que \emph{b} et \emph{c} sont \eqvs parce que $Q$
est \tf (\thrf{prop quot non iso}). 
Pour démontrer que \emph{a} implique \emph{b}, on peut
supposer que $Q$ est libre (quitte à considérer $Q'$ tel que $Q\oplus Q'$ est
libre). Alors,~$\varphi$ est représenté par une matrice dont le \deter est \ndz
donc \ivz.
\end{proof}
%-----------------end proof------------------

Le \tho précédent admet un corolaire important en théorie des nombres.

%:     Corollary{corpropZerdimLib}
\begin{corollary} \label{corpropZerdimLib} \emph{(\Tho un et demi)}
\begin{enumerate}
\item Soit $\fa$ un \id  de $\gA$.
On suppose que c'est un \Amo \pf avec $\cF_1(\fa)=\gen{1}$ et
qu'il existe $a\in\fa$
tel que l'anneau $\gB=\aqo{\gA}{a}$ soit \zedz. Alors, il existe~\hbox{$c\in\fa$} tel que~\hbox{$\fa=\gen{a,c}=\gen{a^m,c}$} pour tout $m\geq1$.
\item  Soit $\gZ$ l'anneau d'entiers d'un corps de nom\-bres~$\gK$ et~$\fa$ un \itf non nul de $\gZ$. Pour tout $a\neq0$ dans $\fa$ il existe~$c\in\fa$ tel \hbox{que
$\fa=\gen{a,c}=\gen{a^m,c}$} pour tout $m\geq1$.%
%:HHH index
\index{un et demi!\Tho ---}
\end{enumerate}
\end{corollary}
\begin{proof}
\emph{1.} 
Le \Bmo $\fa/a\fa$ est obtenu à partir du \Amo $\fa$ par \eds de $\gA$ à $\gB$,
donc son premier \idf reste égal à $\gen{1}$. On applique le point \emph{\ref{ite4propZerdimLib}} du \tho \ref{propZerdimLib}: il existe un $c\in\fa$
tel que $\fa/a\fa=\gen{c}$ en tant que \Bmoz. Cela signifie $\fa=c\gA+a\fa$ et donne le résultat souhaité.

 \emph{2.} Si $\fa=\gen{\xn}$ est un \itf de $\gZ$, il existe un \itf $\fb$
tel que $\fa\fb=\gen{a}$ (\thrf{th1IdZalpha}).
Notons $\ux=[\,x_1 \;\cdots\;x_n \,]$. Il existe donc $y_1$, \dots, $y_n$ dans $\fb$ tels que
$\ux\tra{\uy}=\sum_i x_iy_i=a$. Si $y_ix_j=\alpha_{ij}a$,\linebreak  on~a~$\alpha_{ii}x_k=\alpha_{ki}x_i$. Donc, l'\id $\fa$ devient principal dans $\gZ[1/\alpha_{ii}]$,
égal à~$\gen{x_i}$, qui est libre de rang $1$ (on peut supposer les $x_i$ non nuls). \\
Puisque~$\sum_i \alpha_{ii}=1$, les $\alpha_{ii}$
sont \comz, donc $\fa$ est \ptf et~$\cF_1(\fa)=\gen{1}$ (ceci est vrai localement donc globalement).\\
Pour appliquer le point \emph{1} il reste à vérifier que $\aqo{\gZ}{a}$ est \zedz.
L'\elt $a$ annule un \polu $P\in\ZZ[X]$ de \coe constant
non nul, ce que l'on écrit $aQ(a)=r\neq 0$.
Donc, $\aqo{\gZ}{a}$ est un quotient de $\gC=\aqo{\gZ}{r}$. Il suffit de montrer que $\gC$ est \zedz. 
Notons~\hbox{$\gA=\aqo{\ZZ}{r}$}.
Soit $\ov u \in\gC$, $u$ annule un \polu $R\in\ZZ[T]$ de degré $n$. Donc l'anneau $\gA[\ov u]$ est un quotient de l'anneau  $\aqo{\gA[T]} {\ov R(T)}$, 
qui est un \Amo libre de rang $n$, et donc est fini. Donc on peut trouver explicitement  $k\geq 0$ et $\ell\geq 1$ tels que ${\ov u}^{k}(1-{\ov u}^{\ell})=0$.  
\end{proof}

\rem La matrice $A=(\alpha_{ij})$ vérifie les \egts suivantes:
$$
\tra{\uy}\,\ux=aA,\;A^2=A,\;\cD_2(A)=0,\;\Tr(A)=1,\;\ux A=\ux.
$$
On en déduit que $A$ est une \mprn de rang $1$. 
En outre, on~a~$\ux(\In-A)=0$. Et si $\ux \tra\uz=0$,
 alors $\tra{\uy}\,\ux \tra\uz=0=aA\tra\uz$, donc $A\tra\uz=0$ et~$\tra\uz=(\In-A)\tra\uz$. Ceci montre que $\In-A$ est une \mpn de~$\fa$ (sur le \sgrz~$(\xn)$).
 Donc,~$\fa$ est isomorphe comme~$\gZ$-module à~$\Im A\simeq\Coker(\In-A)$.
\eoe

%%%%%%%%%%%%%%%%%%%%%%%%%%%%%%%%%%%%%%%%%%%%%%%%%%%%%%%%%%%%%%%%%%%%%%%%%%%
%  Section{Modules stablement libres}-
\section{Modules stablement libres} 
\label{ModStabLibre}
\index{module!stablement libre}
\index{stablement libre!module ---}

Rappelons qu'un module $M$ est dit \emph{\stlz} s'il \emph{est \sul d'un libre dans un libre}, autrement dit s'il existe
un \iso entre~$\Ae n$ et~$M\oplus \Ae r$ pour deux entiers $r$ et $n$.
\\ On dira alors\footnote{Cette notion de rang sera \gneez, \dfns \ref{def ptf rang constant} puis \ref{defiRang}, et \llec pourra constater qu'il s'agit bien de \gnnsz.}
que $M$ est \emph{de rang} $s=n-r$.
Le rang d'un module stablement libre sur un anneau non trivial est bien 
défini. En effet, si $M\oplus \Ae r\simeq\Ae {n}$ et~$M\oplus \Ae {r'}\simeq\Ae {n'}$, alors
 on a $\Ae r\oplus \Ae {n'} \simeq \Ae {r'}\oplus \Ae n$ par le lemme de Shanuel~\ref{corlemScha}.
\`A partir d'un \iso $M\oplus \Ae r\to\Ae n$, 
on obtient la \prn $\pi:\Ae n\to\Ae n$
sur $\Ae r$ \paralm à $M$. Cela donne aussi une \Ali surjective
$\varphi:\Ae n\to\Ae r$ avec $\Ker\pi=\Ker\varphi\simeq M$: il suffit de poser
$\varphi(x)=\pi(x)$ pour tout $x\in \Ae n$.

Inversement, si l'on a une \ali $\varphi : \Ae {n}\to\Ae r$  
surjective,
il existe   $\psi : \Ae r\to\Ae n$ telle que $\varphi \circ \psi =\Id_{\Ae r}$.
Alors, $\pi=\psi\circ \varphi:\Ae {n}\to\Ae n$ est une \prnz, avec
$\Ker\pi=\Ker\varphi$, $\Im\pi=\Im\psi$ et
$\Ker\pi\oplus \Im\pi =\Ae n$.\linebreak 
 Et puisque $\Im\pi\simeq\Im\varphi=\Ae r$, le module

\snic{M=\Ker\varphi=\Ker\pi\simeq\Coker\pi=\Coker\psi
}

%\sni
est \stlz, et isomorphe à  $\Im(\Id_{\Ae n}-\pi)$.
%:HHH  la seconde référence suffit
Rappelons que d'après 
%les \thrfs{propIGCram}{prop inj surj det},
le \thrf{prop inj surj det},
dire que
$\varphi:\Ae n \to \Ae r$ est surjective revient à dire que~$\varphi$ est de rang $r$,
\cad ici que $\cD_r(\varphi)=\gen{1}$.

Enfin, si l'on part d'une \ali injective $\psi : \Ae r\to\Ae n$,
dire qu'il existe  $\varphi : \Ae {n}\to\Ae r$ telle que
$\varphi \circ \psi =\Id_{\Ae r}$ revient à dire que  $\cD_r(\psi)=\gen{1}$
(\thrf{propIGCram}).
Résumons la discussion précédente.

%:     Fact{factStablib}
\begin{fact}\label{factStablib} Pour un module $M$ \propeq
\begin{enumerate}
\item $M$ est \stlz.
\item $M$ est isomorphe au noyau d'une matrice surjective.
\item $M$ est isomorphe au conoyau d'une matrice injective
de  rang maximum.
\end{enumerate}
\end{fact}

Ce résultat peut permettre de définir un nouveau codage, spécifique
pour les modules \stlsz. Un tel module sera codé par
les matrices des \alis $\varphi$ et $\psi$.
Concernant le dual de $M$ il sera codé par les matrices transposées,
comme indiqué dans le fait suivant.

%:     Fact{factStablibDual}
\begin{fact}\label{factStablibDual}
Avec les notations précédentes $M\sta$ est \stlz,
canoniquement isomorphe à $\Coker\!\tra{\varphi} $ et à $\Ker\!\tra{\psi}.$
\end{fact}

Ceci est un cas particulier du résultat plus \gnl suivant
(voir aussi le fait~\ref{factDualReflexif}).

%:     Proposition{propDualSurjScind}
\begin{proposition}\label{propDualSurjScind}
Soit $\varphi:E\to F$ une surjection scindée et $\psi:F\to E$ une section de $\varphi$.
\\
Notons $\pi:E\to E $ la \prn $\psi\circ \varphi$, et  $j:\Ker\varphi\to E$ l'injection canonique.
\begin{enumerate}
\item $E=\Im\psi\oplus\Ker\varphi$, $\Ker\varphi=\Ker\pi\simeq\Coker\pi=\Coker\psi$.
\item $\Ker\!\tra{j}=\Im\!\tra{\varphi}$ et $\tra{j}$ est surjective,
ce qui donne par \fcn
un \iso canonique  $\Coker \tra{\varphi} \simarrow (\Ker\varphi)\sta$.
\end{enumerate}
\end{proposition}
\begin{proof}
L'\ali $\psi$ est un \ing de $\varphi$ (\dfn \ref{defIng}).
On a donc $E=\Im\psi\oplus\Ker\varphi $, et $\psi$ et $\varphi$ définissent des \isos réciproques entre $F$ et $\Im \psi$. La proposition suit facilement  (voir le fait~\ref{factIng0}).
\end{proof}
%

%:HHH un Grandcadre remplace par une \subsec

\penalty-2500
\subsec{Quand un module \stl est-il libre?}

On obtient alors les résultats suivants, formulés en termes de noyau
d'une matrice surjective.

%\subsubsection*{Quand un module \stlz, donné comme le noyau\\ 
%d'une matrice surjective $R\in\Ae {r\times n}$, est-il libre?}

%On obtient alors les résultats suivants.

%:     Proposition{propStabliblib}
\begin{proposition}\label{propStabliblib}\emph{(Quand un module \stl est libre, 1)}\\
Soit  $n=r+s$ et $R\in\Ae {r\times n}$. \Propeq

\begin{enumerate}
\item $R$ est surjective et le noyau de $R$ est libre.
\item Il existe une matrice $S\in\Ae {s\times n}$ telle que la matrice
$\cmatrix{S\cr R}$ est \ivz.
\end{enumerate}
En particulier, tout module \stl de rang $1$ est libre.

\end{proposition}
\begin{proof}
\emph{1 $\Rightarrow$ 2.} Si $R$ est surjective, il existe $R'\in\Ae {n\times r}$ avec $RR'=\I_r$.
\\
Les matrices $R$ et $R'$ correspondent aux \alis $\varphi$
 et $\psi$ dans la discussion préliminaire. En particulier,
 on a $\Ae n=\Ker R\oplus \Im R'$.
 On considère une matrice $S'$
dont les vecteurs colonnes constituent
une base du noyau de $R$. Puisque $\Ae n=\Ker R\oplus \Im R'$, la
matrice $A'=\cmatrix{S'\mid  R'}$ a pour colonnes une base de $\Ae n$.
Elle est \iv et son inverse est de la forme $\cmatrix{S\cr R}$ 
car $R$ est la seule matrice qui vérifie 
$R\,A'=\lst{0_{r,n-r}\mid \I_r}$.
 
\emph{2 $\Rightarrow$ 1.}  Notons $A=\cmatrix{S\cr R}$ et posons $A'=A^{-1}$, que nous écrivons sous la forme $\cmatrix{S'\,|\, R'}$. On a $RS'=0_{r,n-r}$, donc 

\snic{\Im S'\subseteq \Ker R\qquad(\alpha),}

%\sni
et  $RR'=\I_r$.
Donc 

\snic{\qquad\Ker R\oplus \Im R'=\Ae n=\Im S'\oplus \Im R'\qquad (\beta).}

%\sni
Enfin, $(\alpha)$
et $(\beta)$ impliquent $\Im S'= \Ker R$.

Si $M$ est un module \stl de rang $1$, c'est le noyau d'une matrice
surjective $R\in\Ae{(n-1)\times n}$. Puisque la matrice est surjective,
%\linebreak 
on 	\hbox{obtient $1\in\cD_{n-1}(R)$}, et cela donne la ligne $S$ pour compléter $R$ en une matrice \iv
(développer le \deter selon la première ligne).
\end{proof}

%:     Corollary{prop2Stabliblib}
\begin{corollary}\label{prop2Stabliblib}\emph{(Quand un module \stl est libre, 2)}\\
On considère $R\in\Ae {r\times n}$ et $R'\in\Ae {n\times r}$ avec $RR'=\I_r$,
$s:=n-r$. Alors, les modules $\Ker R$  et $\Coker R'$ sont isomorphes
et \propeq
\begin{enumerate}
\item Le noyau de $R$  est libre.
\item Il existe une matrice $S'\in\Ae {s\times n}$ telle que
$\cmatrix{S'\mid R'}$ est \ivz.
\item Il existe une matrice $S'\in\Ae {s\times n}$ et  une matrice $S\in\Ae {s\times n}$ telles que
$$
\blocs{1.4}{0}{.6}{.8}{$S$}{}{$R$}{}\;\;
\blocs{.6}{.8}{1.4}{0}{$S'$}{$R'$}{}{}
=\In.
$$
\end{enumerate}
\end{corollary}

Rappelons qu'un vecteur $x\in\gA^q$ est dit \emph{\umdz}
lorsque ses \coos sont des \ecoz. Il est dit \emph{complétable}
s'il est le premier vecteur (ligne ou colonne) d'une matrice inversible.
\index{unimodulaire!vecteur ---}
\index{completa@complétable!vecteur ---}

%:     Corollary{corpropStabliblib}
\begin{proposition}\label{corpropStabliblib}
\Propeq \perso{peut on l'améliorer avec $=m$ et $q=m+1$?}
\begin{enumerate}
\item Tout \Amo \stl de rang $\geq m$ est libre.
\item Tout vecteur \umd dans $\gA^{{q}\times 1}$
avec $q>m$ est complétable.
\item Tout vecteur \umd dans $\gA^q$
avec $q>m$ engendre un sous-module \sul d'un module libre dans~$\gA^q$.
\end{enumerate}
\end{proposition}
\begin{proof}
Les points \emph{2} et \emph{3} sont clairement \eqvsz.

\emph{1} $\Rightarrow$ \emph{3.} Soit un vecteur \umd $x \in \gA^q$ avec $q > m$.
Alors, on peut écrire $\gA^q = M \oplus \gA x$, et  $M$ est \stl de rang
$q-1 \ge m$, donc libre.

\emph{3} $\Rightarrow$ \emph{1.} Soit $M$ un \Amo \stl de rang $n \ge m$.
On peut écrire $L=M\oplus \gA
x_1\oplus\cdots\oplus\gA x_r$, où $L\simeq \Ae{n+r}$.
Si $r=0$, il n'y a rien à faire.
Sinon, $x_r$ est un vecteur \umd dans $L$,
 donc par hypothèse~$\gA x_r$ admet un
\sul libre dans $L$. Ainsi, $L\sur{\gA x_r}\simeq \Ae{n+r-1}$, 
et il en est de même de $M\oplus \gA x_1\oplus\cdots\oplus\gA
x_{r-1}$, qui est isomorphe à $L\sur{\gA x_r}$.
On peut donc conclure par \recu sur $r$ que $M$ est libre.
\end{proof}
%

%: subsec{Le stable range de Bass}
\subsec{Le stable range de Bass}

La notion de stable range est liée aux manipulations \elrs (de lignes ou de colonnes) et permet dans une certaine mesure de contrôler les modules \stlsz.

%:     Definition{defiStableRange}
\begin{definition}\label{defiStableRange}
Soit $n\geq 0$. On dit qu'un anneau $\gA$ est de \ix{stable range} \emph{(de Bass) inférieur ou égal à} $n$ lorsque l'on peut 
\gui{raccourcir} les \vmds de longueur $n+1$ au sens suivant:

\snic{
\hbox { $1 \in \gen {a,\an}\;\Longrightarrow\;\exists \,\xn,\;1 \in \gen {a_1 + x_1a, \ldots, a_n + x_na}$.}
}

%\sni
Dans ce cas on écrit \gui{$\Bdim \gA < n$}.  
\end{definition}
%--------- fin definition ---------------------------------------------- 

Dans l'acronyme $\Bdim$, $\mathsf{B}$ fait allusion à \gui{Bass}.

La notation $\Bdim \gA < n$ est légitimée d'une part par le point \emph{1} dans le fait suivant, et d'autre part par des résultats à venir qui comparent  la $\Bdim$ à des dimensions naturelles en \alg commutative\footnote{Voir par exemple les résultats dans le chapitre \ref{chapNbGtrs} qui établissent une comparaison avec les dimensions de Krull et de Heitmann.}.
 
Le point \emph{3} utilise l'\id $\Rad \gA $ qui sera défini au chapitre
\ref{chap Anneaux locaux}. La chose à savoir est qu'un \elt de $\gA$ est \iv \ssi il est \iv modulo $\Rad\gA$.

%:     Fact{factStableRange}
\begin{fact}\label{factStableRange} Soit $\gA$ un anneau et $\fa$
un \idz.
\begin{enumerate}
\item Si $\Bdim\gA< n$ et $n<m$ alors $\Bdim\gA< m$.  
\item Pour tout $n\geq 0$, on a $\Bdim\gA< n\Rightarrow\Bdim\gA/\fa < n$.
En abrégé, on écrit cette implication sous forme: $\Bdim\gA/\fa \le \Bdim\gA$.
\item On a $\Bdim(\gA/\!\Rad \gA) =\Bdim\gA$ (en utilisant la même abréviation).
\end{enumerate}
\end{fact}
%--------- fin fact ---------------------------------------------- 

\begin{proof} \emph{1.} On prend $m=n+1$.
Soient $(a, a_0,\ldots, a_n)$ avec $1 \in \gen {a, a_0, \ldots, a_n}$.  \\
On a $1 = ua + va_0 + \dots $, donc $1 \in \gen {a', a_1,\ldots,
a_n}$ avec $a' = ua + va_0$. \\
Donc on a $x_1$, \dots, $x_n$ dans $\gA$  avec $1 \in \gen {a_1 + x_1a',
\ldots, a_n + x_na'}$,  et par suite ${1 \in \gen {a_0 + y_0a, \dots, a_n + y_na}}$
avec $y_0 = 0$ et $y_i = x_iu$ \hbox{pour $i \ge 1$}.

\emph{2} et \emph{3.} Laissé \alec
\end{proof}
%

%--- Fact{corBass}-----------
\begin{fact}
\label{corBass} \emph{(Vecteurs \umds et \trelsz)}\\
Soit $n\geq 0$. Si $\Bdim \gA <  n$ et $V\in\Ae{n+1}$ est \umdz, il peut
être transformé en le vecteur $(1,0\ldots ,0)$
par des \mlrsz.
\end{fact}
%--- end-corollary------------------------------------
%-----------------begin proof------------------
\begin{proof}
Posons $V=(v_0,v_1,\ldots ,v_{n})$, avec $1\in\gen{v_0,v_1,\ldots ,v_{n}}$. On
applique la \dfn avec $a=v_0$, on obtient $x_1$, \dots, $x_n$ tels que
$$1\in  \gen{v_1+x_1v_0,\alb\ldots ,\alb v_n+x_nv_0}.$$
Le vecteur $V$ peut être transformé par \mlrs en le vecteur
$V'=(v_0,v_1+x_1v_0,\ldots ,v_n+x_nv_0)=(v_0,\alb v'_1,\alb \ldots ,\alb v'_{n})$,
et l'on a des $y_i$ tels que $\sum_{i=1}^ny_iv'_i=1$.
Par \mlrsz, on peut transformer~$V'$ en $(1,v'_1,\ldots ,v'_{n})$,
puis en $(1,0,\ldots ,0)$.
\end{proof}
%-----------------end proof------------------

La proposition \ref{corpropStabliblib} 
et le fait \ref{corBass} donnent le \gui{\tho de Bass} suivant.
En fait, le vrai \tho de Bass est plutôt la conjonction du \tho qui suit avec
un \tho qui fournit une condition suffisante pour avoir $\Bdim\gA<n$.
Nous réaliserons différentes variantes dans les \thrfs{Bass0}{Bass} et le fait~\ref{factGdimBdim}.

%:  --- Corollary{corBass2}---\Tho de Bass-Heitmann
\begin{theorem}
\label{corBass2} \emph{(\Tho de Bass, modules stablement libres)}\\
Si $\Bdim\gA<n$, tout \Amo stablement libre de rang $\geq n$ est libre.
\end{theorem}
%--- end-corollary------------------------------------

\penalty-2500
%--- Section{Constructions natu}
\section{Constructions naturelles}\label{secPtfCoNat}

%--- Proposition{propPtfExt}----
\begin{proposition}
\label{propPtfExt} \emph{(Changement d'anneau de base)}\\
Si $P$  est un \Amo \ptf et si $\rho:\gA\rightarrow \gB$ est un \homo d'anneaux,
alors le $\gB$-module $\rho\ist(P)$ obtenu par \eds à $\gB$
est \ptfz. Si $P$ est isomorphe à l'image d'une matrice de \prn $F=(f_{i,j})$,
$\rho\ist(P)$
est isomorphe à l'image de
{\em la même matrice vue dans $\gB$}, \cad la matrice de \prn
$F^\rho= \big(\rho(f_{i,j})\big)$.
\end{proposition}
%--- end-proposition-------------------
%-----------------begin proof------------------
\begin{proof}
Le changement d'anneau de base conserve les sommes directes et les \prnsz.
\end{proof}
%-----------------end proof------------------

Dans la proposition qui suit, on peut a priori prendre pour ensembles d'indices
 $I=\lrbm$ et
$J=\lrbn$, mais de toute manière $I\times J$,
qui sert d'ensemble d'indices pour la matrice carrée qui définit
le \ix{produit de Kronecker} des deux matrices $F$ et $G$ n'est pas égal
à $\lrb{1.. mn}$. Ceci est un argument important en faveur
de la \dfn des matrices à la Bourbaki, \cad avec des ensembles finis d'indices
(pour les lignes et les colonnes) qui ne sont pas \ncrt du type $\lrbm$.

%--- Proposition{propTensptf}----
\begin{proposition}
\label{propTensptf} \emph{(Produit tensoriel)}\\
Si $P$ et $Q$ sont des modules \pros représentés par les \mprns  $F=(p_{i,j})_{i,j\in
I}\in \gA^{I\times I}$ et $G=(q_{k,\ell})_{k,\ell\in J}\in \gA^{J\times J}$,
alors le produit tensoriel
 $P\te Q$ est un \mptf représenté par le produit de Kronecker

\snic{F\te G=(r_{(i,k),(j,\ell)})_{(i,k),(j,\ell)\in I\times J},}

%\sni
où $r_{(i,k),(j,\ell)}=p_{i,j}q_{k,\ell}$.
\end{proposition}
%--- end-proposition----------------------------------------
%-----------------begin proof------------------
\begin{proof}
Supposons $P\oplus P'=\Ae m$ et $Q\oplus Q'=\Ae n$. La matrice  $F$ (resp.\,$G$) représente la \prn sur $P$ (resp.\,$Q$)
\paralm à $P'$ (resp $Q'$).
Alors, la matrice produit de Kronecker $F\otimes G$ représente la \prn de $\Ae m\otimes \Ae n$ sur  $P\otimes Q$,
\paralm au sous-espace  $(P'\otimes Q)\oplus (P\otimes Q')\oplus
(P'\otimes Q')$.
\end{proof}
%-----------------end proof------------------

%--- Proposition{propDual}----
\begin{proposition}
\label{propDual} \emph{(Dual d'un \mptfz, 2)}\\
Si $P$ est représenté par la \mprn  $F=(p_{i,j})_{i,j\in
I}\in \gA^{I\times I}$,
alors le dual de $P$ est  un \mptf représenté par la matrice 
transposée
de~$F$. Si $x$ est un vecteur colonne dans $\Im F$ et $\alpha$ un vecteur colonne dans l'image de $\tra{F}$, le scalaire
$\alpha(x)$ est l'unique \coe de la matrice $\tra{\alpha}\,x$.
\end{proposition}
%--- end-proposition----------------------------------------
%-----------------begin proof------------------
\begin{proof}
Ce point résulte du fait \ref{factDualPTF}.
\end{proof}
%-----------------end proof------------------

%:     Proposition{propAliPtfs}--
\begin{proposition}
\label{propAliPtfs} \emph{(Modules d'\alisz)}
%-----------------begin enum------------------
\begin{enumerate}
\item Si $P$ ou $Q$ est \ptfz, l'\homo naturel (\paref{NOTAthetaMN})  

\snic{\theta_{P,Q}:P\sta\te Q\to\Lin_\gA(P,Q)}

%\sni
est un \isoz.
\item Si $P$ et $Q$ sont \ptfsz, le module $\Lin_\gA(P,Q)$ est  un \mptf
canoniquement
 isomorphe à $P\sta\te Q$, représenté par la matrice $\tra{F}\te G$.
\item Un \Amo $P$  est \ptf \ssi l'\homo naturel $\theta_{P}$ est un \isoz.
\end{enumerate}
%-----------------end enum------------------
\end{proposition}
%--- end-proposition----------------------------------------
%-----------------begin proof------------------
\begin{proof}
\emph{1.}
Supposons $P\oplus P'=\Ae m$.
On a des \isos
%--------------------begin array---------------
$$\arraycolsep2pt\begin{array}{rcl}
\Lin_\gA(\Ae m,Q)& \simeq& \Lin_\gA(P,Q)\oplus \Lin_\gA(P',Q),     \\[1mm]
(\Ae m)\sta\te Q&    \simeq&  (P\oplus P')\sta \te Q \\
&    \simeq&  (P\sta\oplus (P')\sta) \te Q \\
&    \simeq&  \big(P\sta\te Q \big)\oplus \big((P')\sta\te Q \big).
\end{array}$$
%---------------------end array--------------
Ces \isos sont compatibles avec les \homos naturels
%--------------------begin array---------------
$$\arraycolsep2pt\begin{array}{rcl}
Q^m\simeq (\Ae m)\sta\te Q&  \longrightarrow &
\Lin_\gA(\Ae m,Q)\simeq Q^m,   \\
P\sta\te Q&  \longrightarrow &\Lin_\gA(P,Q),     \\
(P')\sta\te Q&  \longrightarrow &\Lin_\gA(P',Q).
\end{array}$$
%---------------------end array--------------
Comme le premier est un \isoz, les autres le sont \egmtz.\\
Le cas où $Q$ est \ptf se traite de façon analogue.
%:% PERSO junk
\perso{il semble inutile
de répéter, détails dans le source (junk)}
\junk{Supposons $Q\oplus Q'=\Ae n$. On a
%--------------------begin array---------------
$$\begin{array}{rcl}
\Lin_\gA(P,\Ae n)& \simeq& \Lin_\gA(P,Q)\oplus  \Lin_\gA(P,Q') \\[1mm]
P\sta\te \Ae n&    \simeq&  P\sta \te (Q\oplus Q') \\
&    \simeq&  \big(P\sta\te Q \big)\oplus  \big(P\sta\te Q' \big).
\end{array}$$
%---------------------end array--------------
Et tout ceci est compatible avec les \homos naturels
%--------------------begin array---------------
$$\begin{array}{rcl}
(P\sta)^n\simeq\big(P\sta\te \Ae n\big)&  \longrightarrow &  \Lin_\gA(P,\Ae n)
\simeq(P\sta)^n   \\
P\sta\te Q &  \longrightarrow &\Lin_\gA(P,Q)     \\
P\sta\te Q' &  \longrightarrow &\Lin_\gA(P,Q').
\end{array}$$
%---------------------end array--------------
Comme le premier est un \isoz, les autres le sont \egmtz.}
%:% fin PERSO et fin junk

\emph{2.} Cas particulier du point \emph{1}. 

\emph{3.} Résulte du point \emph{1.}
et du fait que $P$ est \ptf si l'image de $\theta_P$ contient $\Id_P$
(\thref{propdef ptf}~$(b3)$).
\end{proof}
%-----------------end proof------------------

En utilisant la commutation de l'\eds avec le produit tensoriel on obtient alors le corolaire suivant.
%:     Corollary{corpropAliPtfs}
\begin{corollary}\label{corpropAliPtfs}
 Si $P$ ou $Q$ est \ptf (sur $\gA$), et si $\gA\vers{\rho}\gB$ est une \algz, l'\homo naturel 
 
\snic{\rho\ist \big(\Lin_\gA(P,Q)\big)\to
\Lin_\gB \big(\rho\ist(P),\rho\ist(Q)\big)}

%\sni
est un \isoz.
\end{corollary}

%--- section{\Tho de structure locale}
\section{\Tho de structure locale}
\label{secMPTFlocLib}
%-----------------------------------------

Dans cet ouvrage, nous donnons plusieurs \dems du \tho de structure locale des \mptfsz. Il y a la voie ouverte par
les \idfsz, qui solde la question le plus rapidement. C'est l'objet de cette
section.

Il y a une méthode éclair basée sur une sorte de formule magique donnée
en exercice \ref{exo7.1}. Cette solution miracle est en fait directement
inspirée par une autre approche du \pbz, basée sur la \gui{relecture dynamique} du lemme de la liberté locale
\paref{lelilo}.
Cette relecture dynamique est expliquée \paref{quasiglobaldynamique}
dans la section~\ref{secMachLoGlo}.

Nous considérons cependant que ce qui éclaire le mieux la situation
est une voie d'accès entièrement basée
sur les matrices de projection et sur des explications
plus structurelles qui font intervenir
l'usage systématique du \deter des \endos des \mptfsz.
Ceci sera fait au chapitre~\ref{chap ptf1}.

%\medskip
%Les \thrfs{prop pf ptf}{theoremIFD} ont comme corolaire immédiat
%le résultat important suivant.
%:HHH je prefere donner la demonstration
%: --- Theorem{prop Fitt ptf 2}-----
\begin{theorem} 
\label{prop Fitt ptf 2}\relax \label{prop Fitt ptf 1}\relax
\emph{(Structure locale et \idfs d'un \mptfz, 1)}
%-----------------begin enum------------------
\begin{enumerate}
\item Un \Amo $P$ \pf est \ptf \ssi ses \idfs sont (engendrés par des)
idempotents.
\item Plus \prmt pour la réciproque, supposons qu'un \Amo $P$ \pf ait
ses \idfs \idmsz, et \hbox{que~$G\in\gA^{q\times n}$} soit une matrice de \pn de $P$,
correspondant à un \sys de~$q$ \gtrsz. 
\\
Notons $f_h$ l'\idm qui engendre $\cF_h(P)$, \hbox{et~$r_h:=f_h-f_{h-1}$}. 
%-----------------begin enum------------------
\begin{enumerate}
\item $(r_0,\ldots,r_q)$ est un \sfioz.
\item Soient $t_{h,j}$ un mineur d'ordre $q-h$ de $G$, et $s_{h,j}:=t_{h,j}r_h$. Alors, le~$\gA[1/{s_{h,j}}]$-module $P[1/{s_{h,j}}]$ est libre de rang~$h$.
\item Les \elts  $s_{h,j}$ sont \comz.
\item On a $r_k=1$ \ssi la matrice $G$ est de rang $q-k$.
\item Le module $P$ est \ptfz.
\end{enumerate}
%-----------------end enum------------------
\item  En particulier, un \mptf devient libre après \lon en un nombre
fini d'\ecoz.
\end{enumerate}
%-----------------end enum------------------
\end{theorem}
%--- end-theorem-----------------------------------------
%
\begin{proof}
Le \thref{prop pf ptf}
nous dit que le module $P$ présenté par la matrice $G$ est \pro \ssi
la matrice $G$ est \lnlz. On applique ensuite la \carn des matrices \lnls par leurs \idds donnée dans le \thref{theoremIFD}, ainsi que la description précise de la structure des matrices \lnls donnée dans ce \tho
(points \emph{5} et  \emph{7} du \thoz).

Note: Le point \emph{3} peut être obtenu plus directement en appliquant le \thref{theoremIFD} à
une matrice \idme (donc \lnlz) dont l'image est isomorphe au module $P$. 
\end{proof}
Ainsi, les \mptfs sont  localement libres,  au sens fort
donné dans le \thrf{th.ptf.loc}.

\smallskip Dans la section \ref{sec ptf loc lib} nous donnerons une preuve
alternative pour le \thref{th.ptf.loc}, plus intuitive et plus éclairante
que celle que nous venons de fournir. En outre, les \eco qui fournissent
des \lons libres seront moins nombreux.

\medskip \rem  %--- Remark{rem.test.projectif}--
\label{rem.test.projectif}\relax
On peut donc tester si un \mpf  est projectif
ou non lorsque l'on sait tester si ses \idfs sont idempotents ou non.
Ceci est possible si l'on sait tester l'appartenance
$x\in \gen{a_1,\ldots ,a_h}$ pour tout \sys $(x,a_1,\ldots,a_h)$
d'\elts de $\gA$, c.-à-d. si l'anneau est fortement discret.
%Pour plus de précisions voir le chapitre~\ref{ChapMptfCalcForm}.
On pourra comparer à \cite{MRR} chap. III exer\-cice 4 p. 96.
\eoe

%--- subsubsec{Annulateur d'un mptf
\penalty-2500
\rdb
\subsection*{Annulateur d'un \mptfz}\label{subsubsecAnnMptf}

%:     Lemma{lemAnnMptf}
\begin{lemma}\label{lemAnnMptf}
L'annulateur d'un \mptf $P$ est égal à son premier \idf
$\cF_0(P)$, il est engendré par un \idmz.
\end{lemma}
\begin{proof}
On sait que les \idfs sont  engendrés par des \idmsz.
On sait aussi que $\cF_{0}(P)\subseteq\Ann(P)$ (lemme \ref{fact.idf.ann}).\\
Voyons l'inclusion contraire. Le fait \ref{fact.transporteur} implique que l'annulateur d'un \mtf se comporte bien par \lonz, donc pour tout \mo $S$, on a $\Ann_{\gA_S}(P_S)= \big(\Ann_\gA(P)\big)_S$. On sait qu'il en est de même pour les \idfs d'un \mpfz. 
Par ailleurs, pour prouver une inclusion d'\idsz, 
on peut localiser en des \ecoz. 
On choisit donc des \eco qui rendent le module $P$ libre, 
auquel cas le résultat est évident.
\end{proof}

La \dem précédente illustre la force du \tho de structure locale
(point \emph{3} du \thrf{prop Fitt ptf 2}).
La section suivante en est une autre illustration.

%--- subsec Ideaux projectifs
\section[Modules \lmos \pros]{Modules \lmos \pros et idéaux \ptfsz}
\label{secIdProj}\label{Idpp}
%-----------------------------------------

\vspace{4pt}
%:--- Sec {sec mlm}--------
\subsec{Modules localement monogènes}
\label{sec mlm}\relax

Un \Amo $M$ est dit \ixc{monogène}{module ---} ou \ixc{cyclique}{module ---}
s'il est engendré par un seul \eltz: $M=\gA{a}$. Autrement dit, s'il est isomorphe
à un quotient~$\gA/\fa$.
%H j'ai un doute ...
%-- Definition{defmlm}-------

\ms En \clama un module est dit {\em \lmoz} s'il devient monogène après
\lon en n'importe quel \idepz.
Il semble difficile de fournir un énoncé \eqv qui fasse sens en \comaz.
Rappelons aussi que la remarque \paref{remplcc.tf} montre que la notion ne semble pas
pertinente lorsque le module n'est pas supposé \tfz.
Néanmoins lorsque l'on se limite aux \mtfs il n'y a pas de problème.
%:  defi  defmlm
La \dfn suivante a déjà été donnée avant le fait~\ref{factExl1Plg}.
\begin{definition}
\label{defmlm} Un \Amo \emph{\tfz} $M$  est dit {\em \lmoz}
s'il existe des \moco $S_1$, \ldots, $S_n$ de $\gA$ tels que chaque
$M_{S_j}$ est monogène comme $\gA_{S_j}$-module. Dans le cas d'un
\id on parle d'\emph{\id \lopz}.%
\index{localement!module --- monogène}\index{module!localement monogène}%
\index{localement!idéal --- principal}\index{ideal@idéal!localement principal}
\end{definition}
%--- end defi --------------------------------------

Notez que la \prt \gui{locale concrète} dans la \dfn
précé\-dente,
sans l'hypothèse que $M$ est \tfz, implique que $M$ est \tf
(\plgrf{plcc.tf}).

\ss Nous aurons besoin de la remarque suivante.
%-- Fact{factLocCas}----------------
\begin{fact} \emph{(Lemme des \lons successives, 1)}
\index{Lemme des \lons successives, 1}
\label{factLocCas} \\
Si $s_1$, \ldots, $s_n$ sont des \eco de $\gA$ et si pour chaque $i$,
on a des \elts
$
s_{i,1},\; \ldots ,\;
s_{i,k_i},
$
 \com dans $\gA[1/s_i]$,
alors les $s_{i}s_{i,j}$ sont \com dans~$\gA$.
\end{fact}
%--- end-fact-----------------------------------------

Voici maintenant, en point \emph{3} du \tho suivant, une machinerie
calculatoire efficace pour les modules \lmosz.

%:-- theorem{propmlm}--
\begin{theorem}
\label{propmlm} {\em (Modules \tf \lmosz)}\\
Soit $M=\gA x_1+\cdots+\gA x_n$ un \mtfz.
\Propeq
\begin{enumerate}
\item [1.\phantom{*}] Le module $M$ est \lmoz.
\item [2.\phantom{*}] Il existe $n$  \eco $s_i$ de $\gA$
tels que pour chaque $i$ on ait l'\egt  $M=_{\gA_{s_i}}\gen{ x_i}$.
\item [3.\phantom{*}] Il existe
une matrice $A = (a_{ij})\in\Mn(\gA)$
qui vérifie:
%---- equation {eqmlp} ----
\begin{equation}\preskip.2em \postskip.4em
\label{eqmlm}
\left\{\arraycolsep2pt
\begin{array}{rcl}
\sum a_{ii}      &=&  1\\[.2em]
a_{\ell j}x_{i} & =& a_{\ell i}x_{j} \qquad \forall i, j,
\ell \in \lrbn,
\end{array}
\right.
\end{equation}
%---------------------end equation--------------
autrement dit, pour chaque ligne $\ell$, la matrice suivante
est formellement de rang $\leq 1$ (ses mineurs d'ordre 2 sont nuls)
$$\preskip.4em \postskip-.2em
\cmatrix{
a_{\ell 1}&\cdots &a_{\ell n}\cr
x_1&\cdots &x_n
}
.$$
\item [4.\phantom{*}]$\Vi_\Ae 2(M)=0$.
\item [5.\phantom{*}] $\cF_1(M)=\gen{1}$.
\item [6*.] Après \lon en n'importe quel \idepz,  $M$ est monogène.
\item [7*.] Après \lon en n'importe quel \idemaz,  $M$ est monogène.
\end{enumerate}
\end{theorem}
%--- end-theorem----------------------------------------
%-----------------begin proof------------------
\begin{proof}
\emph{3} $\Rightarrow$ \emph{2} $\Rightarrow$ \emph{1.} Clair, avec $s_i=a_{ii}$ dans
le point~\emph{2}.

Montrons qu'un module monogène vérifie la condition \emph{3}.
\\
Si $M = \gen{g}$, on~a~$g = \sum_{i=1}^{n} u_i x_i$ et  
%:2015 en dessous
$x_i = y_ig$. 
Posons $b_{ij} = u_iy_j$.\\ 
Alors, pour tous
   $i$, $j$, $\ell \in \lrbn$, on a
   $b_{\ell j} x_{i} = u_\ell  y_i y_j g= b_{\ell i} x_{j}$. 
En outre:

\hfil
$%\displaystyle
g = \som_{i=1}^{n} u_i x_i = \som_{i=1}^{n} u_i y_i g
= \big( \som_{i=1}^{n} b_{ii} \big) g.
$

 Posons $s=1-\som_{i=1}^{n} b_{ii}$.
On a $sg=0$, et donc $sx_k=0$ pour tout $k$.
\\
Prenons $a_{ij} = b_{ij}$ pour $(i, j) \neq (n,n)$ et $a_{nn} = b_{nn} + s $.
Alors, la matrice~$(a_{ij})$ vérifie bien les équations~(\ref{eqmlm}).

\emph{1} $\Rightarrow$  \emph{3.}
La \prt \emph{3} peut être vue comme l'existence d'une
solution pour un \sli dont les \coes s'expriment en fonction
des \gtrsz~$x_i$.
Or un module monogène vérifie la \prtz~\emph{3.}
On peut donc appliquer le \plg de base.\iplg

Ainsi, \emph{1} $ \Leftrightarrow$ \emph{2} $\Leftrightarrow$ \emph{3.}

\emph{1} $\Rightarrow$  \emph{4} et \emph{1} $\Rightarrow$ 
\emph{5.} Parce que les foncteurs $\Vi_\Ae 2\bullet$
et $\cF_1(\bullet)$ se comportent bien par \lonz.

\emph{5} $\Rightarrow$ \emph{1}. $M$ est le quotient d'un \mpf $M'$ tel que
$\cF_1(M')=\gen{1}$, on peut donc supposer \spdg que~$M$ est \pf avec une
\mpn $B\in\Ae{n\times m}$.
Par hypothèse, les mineurs d'ordre $n-1$ de la matrice $B$ sont \comz. 
Lorsque l'on inverse l'un de ces mineurs, par le lemme du mineur inversible
\paref{lem.min.inv}, la matrice $B$ est \eqve à une matrice
%--------------------begin pmatrix---------------

\snic{
\cmatrix{
   \I_{n-1}   &0_{n-1,m-n+1}      \cr
    0_{1,n-1}&     B_1},}

%\sni
et la matrice $B_1\in\Ae {1\times (m-n+1)}$ est aussi une
 \mpn de~$M$.
\\
Supposons \emph{4}  et $n\geq 2$, et montrons que $M$ est,
après \lon en des \eco convenables,
engendré par $n-1$ \eltsz.
Cela suffira à montrer
(en utilisant une \recu sur $n$) que  \emph{4} $\Rightarrow$ \emph{1},
en utilisant le fait~\ref{factLocCas}.
Le module $\Al2_\gA(M)$ est engendré par
les $v_{j,k}=x_j\vi x_k$ ($1\leq j<k\leq n$)
et  les \syzys entre les $v_{j,k}$ sont toutes obtenues
à partir des \syzys entre les $x_i$. Donc si  $\Al2_\gA(M)=0$, $M$ est le
quotient d'un \mpf $M'$ tel que $\Al2_\gA(M')=0$. On suppose alors \spdg que $M$
est \pf avec une \mpn $A=(a_{ij})$. Une \mpn $B$ pour $\Al2_\gA(M)$ avec
les \gtrs $v_{j,k}$ est obtenue comme indiqué dans la proposition
\ref{propPfPex}. C'est une matrice de format ${{n(n-1)}\over{2}}\times m$
(pour un $m$ convenable), et chaque
\coe de $B$ est nul ou égal à un $a_{ij}$. Cette matrice est surjective
donc $\cD_{n(n-1)/2}(B)=\gen{1}$ et les  $a_{ij}$ sont \comz. Or lorsque l'on
passe de $\gA$ à $\gA[1/a_{ij}]$, $x_i$ devient \coli des $x_k$ ($k\neq i$)
et $M$ est engendré par $n-1$ \eltsz.

\emph{1} $\Rightarrow$ \emph{6*} $\Rightarrow$ \emph{7*}. \'Evident.

La preuve que \emph{7*} 
implique \emph{3} est non \covz: on remplace dans la preuve que \emph{1} implique \emph{3}
l'existence d'une solution pour un \sli en vertu du \plg de base, par l'existence
d'une solution en vertu du \plg abstrait correspondant.\iplg
\end{proof}
%-----------------end proof------------------

%\newpage

Dans la suite, nous appellerons \ixx{matrice}{de localisation monogène pour
le $n$-uplet $(\xn)$} une matrice $(a_{ij})$ qui vérifie les
équations~(\ref{eqmlm}). Si les $x_i$ sont des \elts de
$\gA$, ils engendrent un \id \lop et nous parlerons de
\ixx{matrice}{de localisation principale}.%
\index{localisation!matrice de --- monogène}\index{localisation!matrice de --- principale}

\medskip
\rem Dans le cas d'un module engendré par 2 \elts $M=\gA x+ \gA y$,
les équations~(\ref{eqmlm}) sont très simples et
une \mlmo pour $(x,y)$ est une matrice $\cmatrix{1-u&-b\cr-a&u}$ qui vérifie:
%%%%%%%%%%%%%%%  equation   eqmlm2gen  %%%%%%%%%%
\begin{equation}\label{eqmlm2gen}
\dmatrix{1-u&-b\cr x&y}= \dmatrix{-a&u\cr x&y}=0,\,\mathrm{i.e.}\quad(1-u)y=bx \;\et\;
ux=ay
\end{equation}
%%%%%%%%%%%%%%%%%%%%%%%%%%%%%%%%%%%%%%%%%
\vspace{-6pt}
\eoe

\vspace{3pt}
%:-- Proposition pmlm ------------
\begin{proposition}\label{pmlm}
Soit $M = \gA x_1+\cdots+\gA x_n $ un \Amo \tfz.
\begin{enumerate}
\item %1
 Si $M$ est \lmo
et si $A= (a_{ij})$ est une \mlmo pour $(x_1,\dots,x_n)$, nous avons les résultats suivants.
\begin{enumerate}
\item  %1a
$\lst{x_1 \;\cdots \; x_n}\;
A =
\lst{x_1 \;\cdots \; x_n}
$.
\item  %1b 
    Les \ids  $\cD_2(A)$ et $\cD_1(A^2 - A)$ annulent $M$.
\item   %1c
On a {\mathrigid 4mu $a_{ii}M\subseteq \gA x_i$, et sur l'anneau $\gA_i = \gA[\fraC1{a_{ii}}]$, on a $M =_{\gA_i} 
\gA_i {x_i}$}.
%Sur l'anneau $\gA_i = \gA[1/a_{ii}]$, on a $M =_{\gA_i} \gA_i {x_i}$.
%
\item   %1d
$ \gen{a_{1j},\dots,a_{nj}} M = \gA{x_j}$.
\item   %1e
Plus \gnltz, pour n'importe quel \elt $y = \sum\alpha_{i}x_i$
de $M$, si l'on pose
$\alpha = \tra {\lst{\alpha_1\;\cdots\;\alpha_n}} $ et $\beta  = A\, \alpha$, alors
$y=\sum_i\beta_ix_i$ et l'on obtient
une \egt de matrices carrées à \coes dans~$M$:
%--------------------begin equation---------------
\begin{equation}\label{eqpmlm}
\preskip.4em \postskip.4em
\beta x = \cmatrix{\beta_1\cr \vdots\cr \beta_n}\lst{x_1 \;\cdots \; x_n}=
A\,y,\;\;\hbox{i.e.}\quad \forall i,j\;\beta_ix_j=a_{ij}y
\end{equation}
%---------------------end equation--------------
En particulier,
${\gen{\beta _{1},\dots,\beta _{n}} M = \gA y.}$
\end{enumerate}
\item %2
\Propeq
\begin{enumerate}
\item [--] $M$  est isomorphe à l'image d'une matrice
de \prn de rang $1$.
\item [--] $M$ est fidèle (i.e. $\Ann(M)=0$) et \lmoz.
\end{enumerate} 
Dans ce cas, soit $A$ une \mlmo pour $(\xn)$. On obtient:
\begin{enumerate}
\item [--] $A$ est une  \mprn
de rang 1,
\item [--] la suite ci-après est exacte:
$
\Ae n\vvvvers{\In-A} {\Ae n}\vvvvvers{[\,x_1\;\cdots\;x_n\,]}M\to 0,
$
\item [--]  $M\simeq\Im A$.
\end{enumerate}
\end{enumerate}
\end{proposition}

\begin{proof}
\emph{1.} Le point \emph{1c} est clair, et \emph{1d} est un cas
particulier de \emph{1e}.

\emph{1a.} La $j$-ème \coo de $\lst{x_1 \;\cdots \; x_n}\, A$ s'écrit:
$${
\som_{i=1}^{n}a_{ij}x_i = \som_{i=1}^{n}a_{ii}x_j = x_j}.$$

\emph{1b.} Montrons que tout mineur d'ordre 2 de $A$ annule $x_i$: on
considère la matrice suivante
$$
\cmatrix{
a_{ji}  &  a_{j\ell}  &  a_{jh}\cr
a_{ki}  &  a_{k\ell}  &  a_{kh}\cr
x_i     &  x_\ell     &  x_h
}.
$$
Son \deter est nul (en développant par rapport à la première ligne)
et le développement par rapport à la première colonne fournit
$$\preskip.4em \postskip.4em
  (a_{j\ell} a_{kh} - a_{jh} a_{k\ell})x_i =0.
$$
Montrons que $A^2 = A$ modulo $\Ann(M)$. Ce qui suit est écrit modulo cet
annulateur. On vient de montrer que les mineurs d'ordre 2 de $A$ sont nuls.
Ainsi,  $A$ est une \mlmo pour chacune de ses \hbox{lignes $L_i$}.
D'après le point \emph{1a} appliqué à $L_i$, on a $L_i\, A=L_i$, et 
\hbox{donc $A^2=A$}.

\emph{1e.}  Notons $\ux=\lst{x_1 \;\cdots \; x_n}$.
D'une part 
$$
\preskip.2em \postskip.0em 
\som_i\beta_ix_i=\ux\beta=\ux A\alpha=\ux\alpha=\som_i\alpha_ix_i. 
$$
D'autre part,
$$
\preskip.4em \postskip.4em 
\beta_i x_{j} = \som_{k}\alpha_{k}\, a_{ik} \,x_j=
\som_{k}\alpha_{k}\, a_{ij} \,x_k =
a_{ij}\,\big(\som_{k}\alpha_{k} x_k\big) = a_{ij}\, y. 
$$
Ceci montre l'\egt (\ref{eqpmlm}) et l'on en déduit
$\gen{\beta _{1},\dots,\beta _{n}} M =  \gA y$.

\sni \emph{2.}
Supposons tout d'abord que   $M$  est isomorphe à l'image d'une
 \mprn  $A$ de rang $1$. Notons $x_i$ la $i$-ème colonne de $A$.
Comme~$\cD_2(A)=0$, on a les \egts
$a_{\ell j}x_{i}  = a_{\ell i}x_{j}$ pour $i, j,
\ell \in \lrbn$. Ceci implique que
sur l'anneau $\gA[1/a_{\ell j}]$,  $M$ est engendré par $x_j$, et
puisque  $\cD_1(A)=\gen{1}$, le module est \lmoz. Enfin, soit~$b\in\Ann(M)$,
 \hbox{alors~$bA=0$}, et~$\cD_1(A)=\gen{1}$ implique $b=0$: le module est fidèle.
\\
 Supposons maintenant que $M$ est \lmoz, et que $A$
est une \mlmo pour un \sys \gtr
$(\xn)$.
Si $M$ est fidèle, vu  \emph{1b}, on a $\cD_2(A)=0$ et $A^2=A$ donc $A$
est une \mprn de rang $\leq1$.
Puisque $\Tr(A)=1$, $A$ est de rang $1$. Vu~\emph{1a}, la matrice $\I_n-A$ est une
matrice de \syzys pour $(\xn)$. Soit maintenant
$\sum_{i=1}^n\alpha_ix_i=0$ une \syzy arbitraire en les $x_i$.
Posons comme dans \emph{1e}
$$
\preskip.3em \postskip.4em 
\beta=\tralst{\beta_1\;\cdots\;\beta_n} = A\;
\tralst{\alpha_1\;\cdots\;\alpha_n}, 
$$
on obtient $\gen{\beta_{1},\dots,\beta_{n}} M =  0$ et, 
puisque $M$  est fidèle, $\beta=0$.\\
Ainsi,
$\,A\,\tralst{\alpha_1\;\cdots\;\alpha_n}=0$
et $(\In-A)\,\tralst{\alpha_1\;\cdots\;\alpha_n}=\tralst{\alpha_1\;\cdots\;\alpha_n}$:
toute  \syzy pour $(\xn)$ est
une \coli des colonnes de  $\I_n-A$. Ceci montre que
 $\I_n-A$ est une  \mpn de $M$  pour le \sgr $(\xn)$.
Puisque $A^2=A$, on a $M\simeq\Coker(\I_n-A)\simeq\Im A$.
\end{proof}

%:--- subsec Modules monogenes projectifs
%\rdb
%\penalty-2500
\subsec{Modules monogènes  \prosz}%\label{Idpp}

La description suivante s'applique en particulier pour les \idps \prosz.
%--- Lemma{lemIdpPtf}-----------
\begin{lemma}
\label{lemIdpPtf}
Pour un module monogène $M$, \propeq
%-----------------begin item------------------
\begin{enumerate}
\item $M$ est un \Amo \ptfz.
\item $\Ann(M)=\gen{s}$ avec $s$ \idmz.
\item $M\simeq \gen{r}$ avec $r$ \idmz.
\end{enumerate}
%-----------------end item------------------
\end{lemma}
%--- end-lemma--------------------
%-----------------begin proof------------------
\begin{proof}
Les implications \emph{2 $\Rightarrow$ 3
$\Rightarrow$ 1} sont évidentes, et
l'implication \emph{1 $\Rightarrow$ 2} est donnée dans le lemme~\ref{lemAnnMptf}.

%-% ENTRE NOUS
\entrenous{Preuve indépendante de \emph{1 $\Rightarrow$ 2}: Voyons que  \emph{1 $\Rightarrow$ 2}. On a l'\ali surjective
$\psi : \gA\rightarrow x\gA,\;a\mapsto ax$. Donc
(\tho \ref{propdef ptf}~(c3)) $\Ker\psi=\Ann(x)$ est facteur direct dans $\gA$.
On conclut avec le lemme~\ref{lemfacile}.
}
%-% Fin ENTRENOUS
\end{proof}
%-----------------end proof------------------

On en déduit qu'un anneau $\gA$ est \ixe{\qiz}{quasi integre}
\ssi tout \idp est \proz,
ce qui justifie la terminologie anglaise de
\ix{pp-ring} (principal ideals are projective).\index{anneau!quasi intègre}

%:--- subsec{Modules localement monogènes proj
%\rdb
\subsec{Modules \lmos \prosz}\label{subsubseclmoproj}

Le lemme suivant généralise l'\eqvc donnée dans la proposition
\ref{pmlm} entre module \lmo fidèle et image d'une \mprn de rang $1$.
%:     Lemma{lemLmoProj}
\begin{lemma}\label{lemLmoProj} 
\Propeq
\begin{enumerate}
\item $M$ est \lmo et $\Ann (M)$ est engendré par un \idmz.
\item $M$ est \ptf et \lmoz.
\item $M$ est isomorphe à l'image d'une \mprn de rang~$\leq1$.
\end{enumerate}
\end{lemma}
\begin{proof}
\emph{1} $\Rightarrow$ \emph{2.} On localise en des \eco qui rendent le module monogène
et l'on applique le lemme~\ref{lemIdpPtf}.

Dans \emph{2} et \emph{3} on note $F$ une \mprnz, carrée d'ordre $n$, ayant~$M$ pour image.  Après \lon en des \eco elle devient semblable à  une \mprn standard $\I_{k,n}$, $k$ dépendant de la \lonz.

\emph{2} $\Rightarrow$ \emph{3.} Si $k>1$,
on obtient dans la \lon correspondante $\cF_1(M)=\gen{0}$.
Comme on a déjà $\cF_1(M)=\gen{1}$, la \lon est triviale. Le rang de $F$ est donc $\leq1$ dans toutes les \lonsz.
%-% PERSO
\perso{un peu bizarre que ce soit si compliqué,
on doit pouvoir raisonner en disant $\Al2F=0$ \ssi $\Al2M=0$, ce
qui doit être vrai chaque fois que $M$ est le module image de la matrice~$F$
}
%-% Fin PERSO

\emph{3} $\Rightarrow$ \emph{1.}  Après \lonz, comme
la matrice  est
de rang $\leq1$, on a $k\leq1$. Le module devient donc monogène.
Par ailleurs, d'après le lemme~\ref{lemAnnMptf}, $\Ann (M)$ est engendré par un \idmz.
 \end{proof}

%:--- subsubsec{Ideaux proj
\subsec{Idéaux \ptfs}\label{subsubsecIdproj}
%/HHH \index{ideal@idéal!régulier}  plutot: fidele
 Rappelons qu'un \id $\fa$ est dit \emph{fidèle} s'il est
fidèle comme \Amoz.%
\index{fidèle!idéal ---}\index{ideal@idéal!fidèle}

\smallskip 
\rem Dans la terminologie la plus répandue,
 un idéal est appelé régulier s'il contient un \elt \ndzz.
 A fortiori c'est un \id fidèle.
Nous n'utiliserons pas cette terminologie car elle nous semble prêter à confusion. \eoe

%:     Lemma{lemIdproj}
\begin{lemma}\label{lemIdproj}~
\begin{enumerate}
\item \label{i1lemIdproj} Si $\fa\subseteq\fb$ avec $\fa$ \tf et $\fb$ \lopz, il existe un \itf
 $\fc$ tel que
$\fb\fc=\fa$.

\item \label{i2lemIdproj} Un \id $\fa$ est \ptf \ssi il est \lop et son annulateur
est engendré par un \idmz.

\item \label{i20lemIdproj} Un \id $\fa$ est quasi libre \ssi il est
principal et son  annulateur
est engendré par un \idmz.

\item \label{i3lemIdproj} Soient $\fa_1$ et $\fa_2$ des \ids et $\fb$
un \id \ptf fidèle.
\hbox{Si $\fb\fa_1=\fb\fa_2$,} alors $\fa_1=\fa_2$.

\item \label{i6lemIdproj}
Un \id est \iv \ssi il est \lop et il contient un \elt \ndzz.
\end{enumerate}
\index{inversible!idéal}
\index{ideal@idéal!inversible}
\end{lemma}
%--- end{lemma} --------------------------------------------
%
\begin{proof}
\emph{\ref{i1lemIdproj}.} Il suffit de montrer que pour un $a\in\fb$ arbitraire on a un \itf $\fc$ tel que $\fb\,\fc=\gen{a}$. 
Ceci est donné par le point~\emph{1e} de la proposition~\ref{pmlm} lorsque
$M=\fb$.

\emph{\ref{i2lemIdproj}.} L'implication directe utilise
le corolaire~\ref{corprop inj surj det}:  si une \ali $\Ae k\to\gA$ est injective avec $k>1$, l'anneau est trivial. Donc dans chaque \lonz, l'\id $\fa$ est non seulement libre mais monogène.
L'implication réciproque est dans le lemme~\ref{lemLmoProj}.

\emph{\ref{i20lemIdproj}.} Pour l'implication directe, on écrit
$\fa\simeq \bigoplus_{i\in \lrbn} \gen{e_i}$, où les $e_i$ sont des \idms
avec $e_{i+1}$ multiple de $e_i$. On veut montrer que si $n> 1$,  $e_2=0$.
On localise l'injection $\fa \to \gA$ en $e_2$ et l'on obtient une injection

\snic{
\gA_{e_2}\oplus \gA_{e_2}\simeq e_1\gA_{e_2}\oplus e_2\gA_{e_2}\hookrightarrow \bigoplus e_i\gA_{e_2}\simeq \fa\gA_{e_2} \hookrightarrow \gA_{e_2},}

%\sni
donc $\gA_{e_2}$ est nul (corolaire~\ref{corprop inj surj det}).

\emph{\ref{i3lemIdproj}.} L'\id $\fb$ devient
libre (après \lonz), et monogène d'après le point~\emph{\ref{i2lemIdproj}.} 
Si en plus il est fidèle,
son annulateur est nul, et le \gtr est un \elt \ndzz.

\emph{\ref{i6lemIdproj}.} Le point \emph{\ref{i1lemIdproj}} implique
qu'un \id \lop qui contient un \elt \ndz est \ivz. Réciproquement,
soit $\fa=\gen{a_1,\ldots,a_n}$ un \id \ivz. Il existe $c$ \ndz dans $\fa$ et
un \id $\fb$ tels que $\fa\,\fb=\gen{c}$.\linebreak 
 Soient $b_1$, \ldots, $b_n\in \fb$
avec $\sum_ia_ib_i=c$. On a pour chaque $i$, $j\in\lrbn$ \hbox{un $c_{ij}\in\gA$}
tel que $b_ia_j=c\,c_{ij}$. En utilisant le fait que $c$ est \ndz on vérifie sans peine que la matrice $(c_{ij})_{1\leq i,j\leq n}$ est une \mlp
pour $(a_1,\ldots,a_n)$.
\end{proof}
%

%  Section{subsec det ptf} ---------
%\penalty-2500
\section[Déterminant, \polfon et \polmuz]{Déterminant, \polcarz, \\
\hspace*{18pt}\polfon et \polmuz}
\label{subsec det ptf}\relax
%------------------

\vspace{4pt}
Si $M$ est un \Amoz, nous notons $M[X]$ le $\AX$-module obtenu par \edsz.

Lorsque $\gA$ est un anneau intègre, si $P$ est un module \ptfz,
isomorphe à l'image d'un \prr $F\in\GAn(\gA)$,
on obtient par passage au corps des
fractions un espace vectoriel $P'$ de dimension finie, disons $k$.
On en déduit que le \polcar de la matrice $F$  est égal à
$(X-1)^kX^{n-k}$.
De manière plus simple encore, le \deter de la multiplication par $X$ dans
$P'[X]$ est égal à
$X^k$, i.e.:

\snic{\det \big((\In -F)+XF \big)=X^k.}

\smallskip
Lorsque $\gA$ est un anneau arbitraire, nous allons voir que l'on peut définir l'analogue du \pol $X^k$ ci-dessus.
Tout d'abord, nous introduisons le \deter d'un \endo d'un \mptfz.

%:  subsec  {Deter, polcar et endo cotranspose}
\subsec{Le \deterz, le \polcar et l'\endo cotransposé}

%:     Thdef{propdef det ptf}-----
\begin{thdef} 
\label{propdef det ptf} Soit $P$ un  \mptfz.
%:HHH enonce simplifié
%---------begin item----------
\begin{enumerate}
\item  Soit $\varphi\in\End(P)$. Supposons que
$P\oplus Q_1$ soit isomorphe à un module libre et 
notons
$\varphi_1=\varphi\oplus \Id_{Q_1}$. %cet \endo de $P\oplus Q$. Alors
%---------begin item----------
\begin{enumerate}
\item  Le \deter de $\varphi_1$ ne dépend que de $\varphi$. 
Le scalaire ainsi
défini est appelé le {\em  \deter de l'endomorphisme} $\varphi$%
\index{determinant@\deterz!d'un endomorphisme}.
On le note $\det(\varphi)$ ou $\det \varphi$.\label{NOTAdet}
\item  Le \deter de l'\endo
$X\Id_{P[X]}-\varphi$ de $P[X]$  est appelé le {\em
\polcar de l'\endoz}~$\varphi$. \\
On le note
$\rC\varphi(X)$; et l'on a $\rC{-\varphi}(0)=\det\varphi$.\label{NOTAPolcar}%
\index{polynome@\pol!caractéristique d'un endomorphisme}
\item  Considérons l'\endo cotransposé $\Adj(\varphi_1)=\wi{\varphi_1}$ de $\varphi_1$. Il
opère sur~$P$, et l'\endo de $P$ ainsi défini ne dépend que de~$\varphi$. On
l'appelle {\em  l'\endo cotransposé de} $\varphi$%
\index{cotransposé!endomorphisme ---} et on le 
note~$\wi{\varphi}$ ou $\Adj(\varphi)$.\label{NOTACotransp}
\item Soit $\rho:\gA\to\gB$ un morphisme. Par \eds de~$\gA$ à~$\gB$,
on obtient un \mptf $\rho\ist(P)$ avec un \endo $\rho\ist(\varphi)$.
 Alors les objets définis précédemment se comportent \gui{fonctoriellemnt}, \cad \prmt 
 
\vspace{-.3em}
 \[ 
 \begin{array}{ccccccc} 
 \det\big(\rho\ist(\varphi)\big)    =  \rho\big(\det(\varphi)\big),~&&\rC{\rho\ist(\varphi)}(X)    =  \rho \big(\rC\varphi(X)\big),  \\[.3em] 
  \Adj\big(\rho\ist(\varphi)\big)   =  \rho\ist\big(\Adj(\varphi)\big). 
  \end{array}
 \] 
\end{enumerate}
%---------end enumerate----------
%
\item 
On a $\det(\Id_p)=1$, et si $\psi : P\rightarrow P$ est un autre \endo de $P$, on~a:
$$\preskip-.4em \postskip-.4em 
\det(\varphi\circ\psi)=\det(\varphi)\det(\psi). 
$$
\vspace{-1mm}
\item  Si $P'$ est un autre \mptf et
si \smashtop{$\psi=\blocs{.7}{.8}{.7}{.8}{$\varphi$}{$\gamma$}{$0$}{$\varphi'$}$} est un \endo
de $P\oplus P'$ \gui{triangulaire par blocs},
on~a
$$\preskip.1em \postskip.3em 
\det(\psi)=d d' \;\hbox{ et }\;
\wi{\psi}=\blocs{.7}{.8}{.7}{.8}{$d'\wi{\varphi}$}{$\eta$}{$0$}{$d\wi{\varphi'}$}
\,,\hbox{ où } d=\det(\varphi),\,d'=\det(\varphi') . 
$$
%4
\item  Si $\varphi : P\rightarrow P$  et $\varphi' : P'\rightarrow P'$ sont des
\endos de \mptfsz, et si $\alpha \circ \varphi =\varphi'\circ \alpha  $ pour un \iso  $\alpha : P\rightarrow P'$,  alors $\det(\varphi )=\det(\varphi' )$.
\item  L'\ali $\varphi$ est un \iso (resp. est injective) \ssi $\det(\varphi)$ est
 inversible (resp. est \ndzz).
\item  L'\egt classique suivante est satisfaite:
$$\preskip.2em \postskip.4em 
\wi{\varphi}\circ \varphi \,=\, \varphi\circ\wi{\varphi}
\,=\,\det(\varphi)\,\Id_P. 
$$
%7
\item  Le \tho de Cayley-Hamilton s'applique:
$\rC\varphi(\varphi)=0$.

\item 
On note
$$\preskip-.0em \postskip.4em 
\begin{array}{rcccl}
\Gamma_{\varphi}(X)& :=  &  -\frac{\rC{-\varphi}(-X)-\rC{-\varphi}(0)}{X} & = &
\frac{-\rC{-\varphi}(-X)+\det(\varphi)}{X}\,,
\end{array} 
$$
%\sni
de sorte que
$ \rC{-\varphi}(-X) =
-X\Gamma_{\varphi}(X)+\det(\varphi)$.
Alors  $\widetilde{\varphi}=\Gamma_{\varphi}(\varphi)$.
\end{enumerate}
%---------end item----------
\end{thdef}
%--- end-thdef-----------------

%---------begin proof----------
\begin{proof}  On remarque que les \dfns données dans le point \emph{1} redonnent bien les
objets usuels de même nom dans le cas où le module est libre.
De même, la formule du point \emph{8} donne, lorsque $\varphi$ est un \endo d'un module libre,  le même $\Gamma_\varphi$
que la formule du lemme~\ref{lemPrincipeIdentitesAlgebriques}. Il n'y a donc pas de conflit de notations. 

\emph{1a.} Supposons que $\Ae m\simeq P\oplus Q_1$ et $\Ae n\simeq P\oplus Q_2$,  et considérons la somme directe
$$\preskip-.5em \postskip.4em 
\Ae {m+n}\simeq P\,\oplus\, Q_1 \oplus P\,\oplus\, Q_2. \eqno(*)
$$
On pose 
$\varphi_1=\varphi\,\oplus\, \Id_{Q_1}$ \hbox{et $\varphi_2=\varphi\,\oplus \,\Id_{Q_2}$}. On doit démontrer l'\egtz~\hbox{$\det\varphi_1=\det\varphi_2$}.
On considère l'\endo 
$$\preskip.3em \postskip.3em 
\phi=\varphi\,\oplus\,\Id_{Q_1}\oplus\,\Id_{P}\oplus\,\Id_{Q_2} 
$$ 
de $\Ae {m+n}$,
de sorte que $\phi$ est conjugué de $\varphi_1\oplus \Id_{\Ae n}$ et de $\varphi_2\oplus \Id_{\Ae m}$. D'où~\hbox{$\det\phi=\det\varphi_1$} \hbox{et $\det\phi=\det\varphi_2$}.

\emph{1c.} On procède pour ce point de la même
manière.
La cotransposition des \endos vérifie le
point~\emph{3} dans le cas des modules libres, donc $\wi\phi$
opère sur $P\oplus Q_1$ et se restreint en $\wi{\varphi_1}$.
En outre, comme $\wi\phi=\Gamma_\phi(\phi)$, $\wi\phi$ opère sur chaque composante dans la somme directe $(*)$. 
De même  $\wi\phi$
opère \hbox{sur $P\oplus Q_2$} et se restreint en $\wi{\varphi_2}$.
Ainsi,
$\wi{\varphi_1}$  et~$\wi{\varphi_2}$ opèrent tous deux sur $P$ de la même
manière que $\wi\phi$. Notez que $\wi\varphi=\Gamma_\phi(\varphi)$.

\emph{1d.} 
Ce point résulte directement des \dfnsz.

Tous les points restants du \tho résultent du cas libre (où le résultat est clair), du \tho de structure locale et du point \emph{1d}. En effet les
énoncés peuvent être certifiés en les vérifiant après \lon est \ecoz, et les \mptfs concernés deviennent simultanément libres après
\lon en un \sys d'\ecoz. Nous donnons néanmoins des \dems plus directes.

 Les affirmations \emph{2}, \emph{3}, \emph{4} et \emph{5} résultent facilement des
\dfnsz, sachant que les résultats sont vrais dans le cas libre.

\emph{6.} On a défini $\wi{\varphi}$  comme
la restriction de
 $\wi{\varphi_1}$ à $P$. Et puisque $\varphi_1$ est un \endo d'un module
libre, on~a
$$\preskip.2em \postskip.3em 
\wi{\varphi_1}\circ \varphi_1 =  \det(\varphi_1)\,\Id_{P\oplus Q_1}, 
$$
ce qui donne par restriction à $P$ l'\egt souhaitée $\wi{\varphi}\circ \varphi =  \det(\varphi)\,\Id_{P}$, car $\det\varphi=\det\varphi_1$.

\emph{7.}
Nous pouvons reproduire la preuve suivante, classique dans
le cas des modules libres. Considérons l'\endo 
$$\preskip.4em \postskip.4em 
\psi=X\Id_{P[X]}-\varphi \in\End_{\AX}(P[X]).
$$ 
D'après le point \emph{6} on~a
$$\preskip.0em \postskip.4em 
~\qquad\qquad\wi{\psi}\psi=\psi\wi{\psi}=\rC \varphi (X)\,\Id_{P[X]}.\eqno (+) 
$$
%\sni
En outre,  \smash{$\wi{\psi}$} est un \pol en $X$ à \coes dans
$\gA[\varphi]$. Autrement dit $\wi{\psi}=\som_{k\geq 0}\phi_kX^k$, avec les
$\phi_k:P\rightarrow P$ qui sont dans $\gA[\varphi]$. En posant $\rC \varphi
(X) = \som_{k\geq 0}a_kX^k$ et en identifiant les deux membres de 
l'\egt $(+)$ on obtient
(en convenant de $\phi_{-1} = 0$)
$$\preskip.4em \postskip.4em 
\phi_{k-1}-\phi_k\varphi =a_k\Id_P\hbox{ pour tout }k \ge 0. 
$$
Ainsi, $\rC \varphi (\varphi ) =\som_{k\ge 0}(\phi_{k-1}-
\phi_k\varphi)\varphi^k=0$.

\emph{8.}
Le \pol $\Gamma_\varphi$ a été défini de façon à vérifier
$$\preskip.4em \postskip.4em 
\rC{-\varphi}(-X) = -X\Gamma_{\varphi}(X)+\det(\varphi). 
$$
En évaluant $X := \varphi$, on
obtient $\varphi\Gamma_\varphi(\varphi) =\det(\varphi)\Id_P$ (\tho de Cayley-Hamilton), d'où $\varphi\Gamma_\varphi(\varphi) =\varphi\wi\varphi$.
En remplaçant
$\varphi$ par  $\theta := T\Id_{P[T]} + \varphi$, on
obtient $\theta\Gamma_\theta(\theta) = \theta\wi\theta$, puis $\Gamma_\theta(\theta)=\wi\theta$,
car $\theta$ est un \elt \ndz de $\gA[T,\varphi] = \gA[\varphi][T]$.  On
termine en faisant $T := 0$.
\end{proof}
%---------end proof----------
\rem
\label{rem det}\relax
Le \deter de l'application identique de tout \mptfz, y compris le module réduit
à $\so{0}$, est égal à 1 (en suivant la \dfn ci-dessus).
\eoe

%--- Corollary{corDetAnn}-------
\begin{corollary}
\label{corDetAnn}
Soit $\varphi:P\rightarrow P$ un \endo d'un
\mptfz, et $x\in P$ vérifiant $\varphi(x)=0$, alors $\det(\varphi)x=0$.
\end{corollary}
%--- end-corollary------------
%---------begin proof----------
\begin{proof} Résulte de
$\wi{\varphi}\circ \varphi=\det(\varphi)\Id_P$.
\perso{cela aurait été beau de
déduire Cayley-Hamilton directement de ce corolaire. Mais la preuve de Richman
ne semble pas passer, et ta tentative ne m'a pas convaincu}
\end{proof}
%---------end proof----------

%
%:\entrenous
\entrenous{Pour un \endo d'un \mrc $k+1$ on aimerait une identification entre
$\wi{\varphi}$ et $\Al k({^t\varphi})$, comme dans l'exercice \ref{exoDualite}
pour le cas libre.
Cela semble coincer parce que $\Al {k+1}(P)$ est de rang 1, mais pas libre en \gnlz. 
Il faudrait éclaircir la chose.
On doit ensuite pouvoir \gnr à $\wi{\varphi}$ et $\Al{(\rg\varphi)-1}({^t\varphi})$. 
Si on comprend ce qui se passe, à mettre en exo quelque part?}

%:  subsec {Le \polfon et le \polmuz}-
\subsec{Le \polfon et le \polmuz}

Nous sommes intéressés par le \polcar de l'identité sur un \mptfz. Il est
cependant plus simple d'introduire un autre \pol qui lui est directement relié
et qui est l'analogue du \pol $X^k$ dont nous parlions au début de la
section~\ref{subsec det ptf}.

%--- Notation{nota RM}-------
\begin{definotas}
\label{nota RM}\relax
Soit $P$ un \Amo \ptf et $\varphi$ un \endo de $P$. On considère le
$\AX $-module  $P[X]$ et l'on définit les \pols
$\rF{\gA,\varphi}(X)$ et $\rR {\gA,P}(X)$ (ou $\rF{\varphi}(X)$ et $\rR {P}(X)$ si le contexte est clair)
par les \egts suivantes:

\snic{\rF{\varphi}(X)= \det(\Id_{P[X]}+X\varphi) \quad {\rm  et} \quad
\rR{P}(X)=\det(X\Id_{P[X]}).
}

%\sni
Donc $\rR{P}(1+X)=\rF{\,\Id_{P}}(X)$.
\begin{itemize}
\item Le \pol $\rF{\varphi}(X)$ est appelé
le {\em \polfonz} de l'\endo $\varphi$.
\item   Le \coe de $X$ dans le \polfon est appelé
la \ixc{trace}{d'un \endo d'un \mptfz} de~$\varphi$ et 
est noté $\Tr_P(\varphi)$.
\item  Le \pol $\rR{P}(X)$ est appelé le {\em \polmuz} du module
$P${\rm(}{\footnote{Cette terminologie est justifiée par le fait que pour un
module libre de rang $k$ le \polmu est égal à~$X^k$, ainsi que par  le
\thrf{th ptf sfio}.}}{\rm)}.
\label{NOTAPolfon}
\end{itemize}
\index{rang!polynome@\pol --- d'un \mptfz}
\label{NOTAPolmu}
\index{polynome@\pol!fondamental}
\index{polynome@\pol!rang}
\end{definotas}
%--- end-definotas------------

On notera que

\snic{
\rF{\varphi}(0)=1=\rR{P}(1),\;\;\rC{\varphi}(0)=\det(-\varphi),\;\;\hbox{et}\;\;
\rF{a\varphi}(X)=\rF{\varphi}(aX),}

%\sni
mais $\rC{\varphi}(X)$ n'est pas toujours unitaire
(voir l'exemple \paref{ex ideal
idempotent}).

On notera également que pour tout $a\in \gA$ on obtient:
%%%%%%%%%%%%%%%%%%%%%%%%%%%%%%%%%%%%%%%%%
\begin{equation}
\label{eq2polmu}
\det(a\varphi)=\det(a\,\Id_P)\det(\varphi)=
\rR{P}(a)\det(\varphi).
\end{equation}
%%%%%%%%%%%%%%%%%%%%%%%%%%%%%%%%%%%%%%%%%

On en déduit les \egts suivantes
$$ 
\begin{array}{c}%\arraycolsep2pt 
 \rR{P}(0)  =    \det(0_{\End_\gA(P)}),\,\,\,\,\,\\[1mm]  
\mathrigid1mu \rC{-\varphi}(-X)\,=\,\det(\varphi -X\Id_{P[X]})  \, = \,   \det \big(-(X\Id_{P[X]}-\varphi)\big)\, =\, \rR{P}(-1)\rC{\varphi}(X),\\[1mm]
\det(\varphi)=\rR{P}(-1)\,\rC{\varphi}(0).\,\,\,\,
  \end{array}
$$
La dernière \egt remplace l'\egt $ \det(\varphi)=(-1)^k\rC{\varphi}(0)$ valable pour les
modules libres de rang $k$.

%:     Theorem{th ptf sfio}-------
%Si $P$ est un \mptfz, il lui correspond un  \sfio de manière canonique.
\smallskip Un \pol $R(X)$ est dit \ixc{multiplicatif}{polynome@\pol ---} si {\mathrigid 2mu $R(1)=1$ et $R(XY)=R(X)R(Y)$}.\index{polynome@\pol!multiplicatif}

\begin{theorem}
\label{th ptf sfio}\relax \emph{(Le \sfio associé à un \mptfz)}
\index{systeme fon@\sfioz!associé à un \mptfz}
%-----------------begin enum------------------
\begin{enumerate}
\item Si $P$ est un \mptf sur un anneau $\gA$, son \polmu
 $\rR{P}$ est multiplicatif.

\item \label{NOTAide} Autrement dit, les \coes de  $\rR{P}(X)$ forment un \sfioz.
Si
$\rR{P}(X)=r_0+r_1X+\cdots+r_nX^n$, on note $\ide_h(P):=r_h$: il est appelé {\em
l'\idm  associé à l'entier $h$ et au module $P$}  (si $h>n$ on pose
$\ide_h(P):=0$).

\item \label{th ptf sfio item reg}\relax Tout \polmu $\rR{P}(X)$ est un \elt \ndz de
$\AX $.

\item \label{remRang} Une \gnn de l'\egt $\,\rg(P\oplus Q)= \rg(P)+\rg(Q)\,$
concernant les rangs des modules libres
 est donnée pour les \mptfs par:
$$
%\snic{
\rR{P\oplus Q}(X)=\rR{P}(X)\,\rR{Q}(X)\,.
%}
$$
\item \label{remRang2} Si $P\oplus Q\simeq\Ae n$ et
$\rR{P}(X)=\sum _{k=0}^nr_kX^k$, alors $\rR{Q}(X)=\sum _{k=0}^nr_kX^{n-k}$.

\item \label{remRang3} L'\egt $\rR{P}(X)=1$ caractérise, parmi les \mptfsz, le module $P=\so{0}$. Elle
équivaut aussi à $\ide_0(P)=\rR{P}(0)=1$.
\end{enumerate}
%-----------------end enum------------------
\end{theorem}
%--- end-lemma--------

%---------begin proof----------
\begin{proof}
\emph{1}  et \emph{2.} Si $\mu_a$ désigne la multiplication par $a$ dans $P[X,Y]$,
on a clairement l'\egt
 $\mu_X\mu_Y=\mu_{XY}$, donc $\rR{P}(X)\,\rR{P}(Y)=\rR{P}(XY)$
(\thref{propdef det ptf}.\emph{2}). Puisque $\rR{P}(1)=\det(\Id_P)=1$, on en
déduit que les \coes de $\rR{P}(X)$ forment un \sfioz.

\emph{3.} Résulte du lemme de McCoy (corolaire~\ref{corlemdArtin}).
On pourrait aussi le démontrer en utilisant le \plg de base (en localisant en les~$r_i$).\iplg

\emph{4.} Résulte du point \emph{3} dans le \thoz~\ref{propdef det ptf}.
\perso{On serait heureux d'avoir aussi une preuve simple et directe pour le
\polmu d'un produit tensoriel.
Naturellement il existe une preuve très facile par \lonz,
mais normalement, on ne voudrait pas utiliser ICI cet argument. Cf. la discussion
qui précède le \thref{corth.ptf.sfio}.
Bref, on a repoussé la chose plus loin.}

\emph{5.} Résulte des points \emph{3} et \emph{4} puisque
$\big(\som_{k=0}^nr_kX^k\big)\big(\som_{k=0}^nr_{n-k}X^k\big)=X^n$.

\emph{6.} On a $r_0=\det(0_{\End(P)})$. Puisque les $r_i$ forment un \sfioz,
les \egts $\rR P=1$ et $r_0=1$ sont \eqvesz.\\
  %\\
   Si $P=\so{0}$, alors $0_{\End(P)}=\Id_P$,
donc $r_0=\det(\Id_P)=1$. 
\\
Si $r_0=1$, \hbox{alors $0_{\End(P)}$} est \ivz, donc $P=\so{0}$.
\end{proof}
%---------end proof----------

Si $P$ est un \Amo libre de rang $k$, on a $\rR{P}(X)=X^k$, la \dfn
suivante est donc une extension légitime des modules libres aux \mptfsz.

%--- Defini{def ptf rang constant}
\begin{definition}
\label{def ptf rang constant}\relax
Un \mptf $P$ est dit {\em  de rang égal à $k$} si $\rR{P}(X)=X^k$. Si l'on ne
précise pas la valeur du rang, on dit simplement que le module est {\em  de rang
constant}.
Nous utiliserons la notation $\rg(M)=k$ pour indiquer qu'un module
(supposé \prcz) est de rang $k$. \label{Nota2rang}
\index{rang!module de --- constant}
\end{definition}
%--- end-definition------------

Notons que d'après la proposition~\ref{propAnnul}, tout module \pro de rang
$k>0$ est fidèle.

%:     Fact{factPolCarRangConstant}
\begin{fact}\label{factPolCarRangConstant}
Le \polcar d'un \endo d'un \mrc $k$ est \mon de degré $k$.
\end{fact}
\begin{proof} On peut donner une élégante \dem directe
(voir l'exercice~\ref{exoPrecisionsDet3}).
On pourrait aussi éviter toute fatigue et utiliser un argument de \lonz,
en s'appuyant sur le \tho de structure locale et sur le fait~\ref{fact.det loc}, qui affirme que tout se passe bien pour le \polcar par \lonz.
\end{proof}

La convention dans la remarque
suivante permet une formulation plus uniforme des \thos et des
preuves dans la suite.

%--- rem{conven rgc}---------

\medskip \rem
\label{conven rgc}\relax
Lorsque l'anneau $\gA$ est réduit à $\so{0}$, tous les \Amos sont
 triviaux. Néanmoins, conformément à la \dfn
ci-dessus, le module nul sur l'anneau nul est un \mrc égal à $k$, pour n'importe quelle
valeur de l'entier $k\geq 0$.
Par ailleurs, il est immédiat que si un \mptf $P$ a deux rangs constants distincts, alors l'anneau est trivial: on a
$\rR{P}(X)=1_\gA X^h=1_\gA X^k$ avec $h\neq k$ donc le
\coe de $X^h$ est égal à la fois à $1_\gA$ et à $0_\gA$.
\eoe

%%%%%%%%%%%%%%%%%%%%%%%%%%%%%%%%%%%%%%%%%%%%%%%%%%%%%%%%%%%%%%%%%%%%%%%%%%%
%:  subsec {Calculs explicites}-

\subsec{Quelques calculs explicites}

\vspace{3pt}
%--- Lemma{lem calculs}--------
Le \polfon d'un \endo $\varphi$ est plus facile à utiliser que le \polcarz.
Cela tient à ce que le
\polfon est invariant lorsque l'on ajoute \gui{en somme directe} un \endo nul à
$\varphi$. Ceci permet de ramener systématiquement et facilement le calcul d'un
\polfon au cas où le module \pro est libre. De manière précise, on
pourra calculer les \pols précédemment définis en suivant le lemme
ci-après.

%:  lemme lem calculs   :  Calcul explicite du \deterz, 
\begin{lemma}
\label{lem calculs}\relax {\rm (Calcul explicite du \deterz, du \polfonz, du \polcarz,
du \polmu et de l'\endo cotransposé)} \\
Soit un \Amo $P\simeq \Im F$ avec  $F\in\GAn(\gA)$. Notons $Q=\Ker(F)$, de sorte
que $P\oplus Q\simeq \Ae n$, \hbox{et $\I_n-F$} est la matrice de la \prnz~$\pi_Q$ sur~$Q$
parallèlement à~$P$.
Un \endo   $\varphi$ de $P$ est \care par la matrice~$H$ de l'\endo
$\varphi_0= \varphi\oplus 0_Q$ de $\Ae n$. Une telle matrice $H$ est soumise à
l'unique restriction
$F\cdot H\cdot F=H$. On pose $G=\I_n-F+H$.

%---------begin item----------
\begin{enumerate}\itemsep0pt
\item  Calcul du \deterz:
$$
\det(\varphi)=\det(\varphi\oplus\Id_Q)=\det(G).
$$
\item  Donc aussi
%---------begin $$--------
%--------------------begin array---------------
$$\arraycolsep2pt\begin{array}{rcl}
\det(X\Id_{P[X,Y]}+Y\varphi)&=&  \det \big((X\Id_{P[X,Y]}+Y\varphi)\oplus\Id_Q \big)\,=\\[1mm]
\det(\I_n-F+XF+YH)&  = &  \det(\I_n+(X-1)F+YH).
\end{array}$$
%---------------------end array--------------
%---------end $$----------
\item  Calcul du  \polmu de $P$:
$$\rR{P}(1+X)=\det \big((1+X)\Id_{P[X]} \big)=\det(\I_n+XF),
$$ 
en particulier,

\snic{ \rR{P}(0)=\det(\I_n-F),}

\smallskip
et $\rR{P}(1+X)=1+u_1X+\cdots+u_nX^n$, où $u_h$ est  la somme des mineurs
principaux d'ordre $h$ de la matrice $F$.

\item  Calcul du  \polfon de $\varphi$:
$$
\rF{\varphi}(Y)=\det(\Id_{P[Y]}+Y\varphi) =\det(\I_n+YH)=
1+\som_{k=1}^n v_k Y^k,
$$
où $v_k$ est la somme des mineurs principaux d'ordre $k$ de la matrice~$H$.
En particulier, $\Tr_P(\varphi)=\Tr(H)$.

\item  Calcul du  \polcar de $\varphi$:
$$
\rC{\varphi}(X)=\det(X\Id_{P[X]}-\varphi) =  \det(\I_n-H+(X-1)F).
$$
\item  Calcul de l'\endo cotransposé $\wi{\varphi}$ de $\varphi$: il est
défini par la matrice
$$
\wi{G}\cdot F =  F \cdot \wi{G}=  \wi{G}-\det(\varphi)(\In-F).
$$
\end{enumerate}
%---------end item----------
\end{lemma}
%--- end-lemma-----------------

Pour le dernier point on applique le point \emph{3} du \thref{propdef det ptf} avec
$\varphi$ et $\Id_Q$ en remarquant
que $G$ est la matrice de $\psi=\varphi\oplus \Id_Q=\varphi_0+\pi_Q$.

Notez que le \polcar de $\Id_P$ est égal à $\rR{P}(X-1)$.

%--- Fact{fact.det loc}---------
 Le fait suivant est une conséquence immédiate de la proposition
\ref{propPtfExt} et du lemme précédent.

\begin{fact}
\label{fact.det loc}\relax
Le \deterz, l'\endo cotransposé, le \polcarz, le \polfon et le \polmu se
comportent bien par \eds via un \homo $\gA\rightarrow \gB$. \\
En particulier, si $\varphi :P\rightarrow P$ est un \endo d'un $\gA$-\mptf et $S$
un \mo de $\gA$, alors $\det(\varphi)_S=\det(\varphi_S)$ (ou, si l'on préfère,
$\det(\varphi)/1=_{\gA_S}\det(\varphi_S)$).
 La même chose vaut pour l'\endo cotransposé, le \polfonz, le \polcar et le
\polmuz.
\end{fact}
%--- end-fact-----------------

%--- exl{ex ideal idempotent}--
\exl
\label{ex ideal idempotent}\relax
Soit  $e$ un \idm de $\gA$ et $f=1-e$. Le module $\gA$ est somme directe des
sous-modules $e\gA$  et $f\gA$ qui sont donc \ptfsz.
La matrice $1\times 1$ ayant pour unique \coe $e$ est une matrice $F$ dont l'image
est $P=e\gA$. Pour $a\in \gA$ considérons $\mu_a=\mu_{P,a}\in\End_\gA(P)$.
La matrice $H$ %correspondante
a pour unique \coe $ea$.
On a alors en appliquant les formules précédentes:
%--------------------begin array---------------
$$\begin{array}{c}
\det(0_{e\gA})=f,\; \rR{e\gA}(X)=f+eX,\;\rC{\Id_{e\gA}}(X)=f-e+eX,\\[1mm]
\; \det(\mu_a)=f+ea,\\[1mm]
\rF{\mu_a}(X)=1+eaX,\;\rC{\mu_a}(X)=1-ea+e(X-1)=f-ea+eX.
\end{array}$$
%---------------------end array--------------
%:HHH  $e\neq 1,0$
Notez que le \polcar de $\mu_a$ n'est pas unitaire \hbox{si $e\neq 1,0$}.
Et l'on a bien le \tho de Cayley-Hamilton:

\snic{\rC{\mu_a}(\mu_a)=(f-ea)\Id_{e\gA}+e\mu_a= (f-ea+ea) \Id_{e\gA} = f\Id_{e\gA} =
0_{e\gA}.  }
\vskip-2em\eoe

%--- SUBSUBsection{Syc}---------------
\subsubsection*{Avec un \sycz}
Lorsque l'on utilise un \syc le lemme~\ref{lem calculs} conduit au résultat suivant.

%:     Fact{factMatriceEndo}
\begin{fact}\label{factMatriceEndo}
Soit $P$ un \mptf avec un \syc $\big((\xn),(\aln)\big)$ et $\varphi$ un \endo de $P$.\\
 On rappelle (fait \ref{factMatriceAlin})  que  l'on peut coder $P$ par la matrice 

\snic{F\eqdefi \big(\alpha_i(x_j)\big)_{i,j\in\lrbn}}

%\sni
($P$ est  isomorphe à  $\Im F\subseteq\Ae n$ au moyen de 
$x\mapsto \pi( x) =\tra[\,\alpha_1(x) \;\cdots\; \alpha_n(x)\,]$). En outre l'\endo $\varphi$ est représenté par la matrice 

\snic{H\eqdefi \big(\alpha_i(\varphi(x_j))\big)_{i,j\in\lrbn}}

%\sni
qui vérifie $H=HF=FH$. %Alors on obtient les calculs suivants.
\begin{enumerate}
%
%\item 
%:HHH  : transféré beaucoup plus tot:
% On dira que \emph{la matrice $\big(\alpha_i(\varphi(x_j))\big)_{i,j\in\lrbn}$ représente l'\endo $\varphi$ dans le \syc $\big((\ux),(\ual)\big)$}.
%\index{matrice!d'un \endo dans un \sycz}
%
\item On a $\rF\varphi(X)=\det (\In+XH)$ et $\Tr(\varphi)=\Tr(H)=\sum_i\alpha_i \big(\varphi(x_i)\big)$.
\item 
Pour $\nu \in P\sta$ et $x$, $y \in P$, rappelons que $\theta_P(\nu \te x)(y)=\nu(y)x$.
%(voir \paref{NOTAthetaMN}). 
La trace de cet \endo est donnée par $\Tr_P \big(\theta_P(\nu \te x)\big)=\nu(x)$.
\end{enumerate}
\end{fact}
\begin{proof}
La matrice $H$
est aussi celle de l'\Ali $\varphi_0$ introduite dans le lemme~\ref{lem calculs}: 

\snic{\pi( x) + y\mapsto \pi \big(\varphi(x)\big)$ avec $\pi( x)\in \Im F$ et $y\in \Ker F.}

%\sni
Le point \emph{2} résulte donc du lemme  \ref{lem calculs}.

\emph{3.} D'après le point \emph{2}, on a:

\snic {\mathrigid1mu
\Tr \big(\theta_P(\nu \te x)\big) = %\sum_i \alpha_i(\varphi(x_i)) = 
\sum_i \alpha_i(\nu(x_i)x) = \sum_i \nu(x_i) \alpha_i(x) = 
\nu\big(\sum_i \alpha_i(x)x_i\big)= \nu(x). 
}
\end{proof}
%

%:     Lemma{lemTraceProT}
\begin{lemma}\label{lemTraceProT}
Soient $M$, $N$ deux \kmos \ptfs et des \endos $\varphi\in\alb\End_\gk(M)$ et
 $\psi\in\End_\gk(N)$. \\
 Alors, $\Tr_{M\te N}(\varphi\otimes \psi)=\Tr_M(\varphi)\Tr_N(\psi)$.
\end{lemma}
\begin{proof}
On considère des \sycs pour $M$ et $N$ et l'on applique la formule pour la trace
des \endos
(fait~\ref{factMatriceEndo}).
\end{proof}
%

%%%%%%%%%%%%%%%%%%%%%%%%%%%%%%%%%%%%%%%%%%%%%%%%%%%%%%%%%%%%%%%%%%%%%%%%%%%
%:  subsec{L'annulateur d'un \mptf}-
\subsec{L'annulateur d'un \mptfz}
\label{subsecAnnul}
%-----------------

Nous avons déjà établi certains résultats concernant cet annulateur
 en nous appuyant sur le \tho de structure locale
 des \mptfsz, démontré en utilisant les \idfs (voir le lemme~\ref{lemAnnMptf}).

Ici nous donnons quelques précisions \suls en utilisant
une \dem qui ne s'appuie pas sur le \tho de structure locale.
%:     Proposition{propAnnul}-----
\begin{proposition}
\label{propAnnul}
Soit $P$ un \Amo \ptfz. On considère
l'\id $J_P= \gen{ \alpha (x) \,|\,\alpha \in P\sta,\; x\in P  }$.  
On note $r_0=\rR{P}(0)=\ide_0(P)$.
\begin{enumerate}
\item $\gen{r_0}=\Ann(P)=\Ann(J_P)$.
%:HHH reformulé
%\item $s_0=1-r_0$ est un \idm et $\gen{s_0}=J_P$.
\item $J_P=\gen{s_0}$, où $s_0$ est l'\idm $1-r_0$.
\end{enumerate}
\end{proposition}
%--- end-proposition----------------
%---------begin proof----------
\begin{proof} On a évidemment $\Ann(P)\subseteq \Ann(J_P)$.
Soit $\big((x_i)_{i\in\lrbn},(\alpha_i)_{i\in\lrbn} \big)$ un \syc
sur $P$. Alors:

\snic{J_P=
\gen{ \alpha_i(x_j)\,;\,{i,j\in\lrbn}},
}

%\sni
et la \mprn $F= \big(\alpha_i(x_j)\big)_{i,j\in\lrbn}$ a une
image isomorphe à $P$. Par \dfnz, $r_0$ est l'\idm
$r_0=\det(\I_n-F)$. Puisque $(\I_n - F)F = 0$, on a $r_0F = 0$, i.e. $r_0P=0$. \\
Donc $\gen{r_0}\subseteq
\Ann(P)\subseteq \Ann(J_P)$ et $J_P\subseteq
\Ann(r_0)$.\\
 Par ailleurs, on a $\I_n-F\equiv \I_n$ modulo $J_P$, donc en
prenant les \detersz, on~a~$r_0\equiv 1$ modulo $J_P$, \cad $s_0\in J_P$, puis
$\Ann(J_P)\subseteq\Ann(s_0)$.
\\
On peut donc conclure:

\snac{\gen{r_0}\subseteq \Ann(P)\subseteq\Ann(J_P)\subseteq\Ann(s_0)=\gen{r_0}
\,\hbox{ et }\,\gen{s_0}\subseteq J_P\subseteq \Ann(r_0) = \gen{s_0}.}
\end{proof}
%---------end proof----------

%%%%%%%%%%%%%%%%%%%%%%%%%%%%%%%%%%%%%%%%%%%%%%%%%%%%%%%%%%%%%%%%%%%%%%%%%%%
%:  subsec {subsec decomp ptf}-------
\subsec{Décomposition canonique d'un module \proz}%{Décomposition canonique}
\label{subsec decomp ptf}\relax
%------------------

%:     Definition{defi comp ptf}
\begin{definition}\label{defi comp ptf}\relax
Soit $P$ un \Amo \ptf et $h\in\NN$.\\
 Si $r_h=\ide_h(P)$,
on note $P\ep{h}$
le sous-\Amo
$r_hP$. %($r_h=0$ si $h>n$)
Il est appelé le {\em  composant du module $P$ en rang~$h$}.
\perso{la notation $P^{[1]}$  entrait en conflit avec le rang
\gne $[1]$ j'ai remplacé par  $P^{(1)}$, cela peut être facilement modifié
car c'est contrôlé par une macro.}
\label{NOTAep}
\end{definition}

Rappelons que pour un \idm $e$ et un \Amo $M$, le module obtenu par \eds à $\gA[1/e]\simeq\aqo{\gA}{1-e}$
s'identifie au sous-module~$eM$, lui-même isomorphe au quotient ${M}\sur{(1-e)M}$.

%:     theorem{propdef comp ptf}--
\begin{theorem}
\label{propdef comp ptf}\label{th decomp ptf}\relax
Soit $P$ un $\gA$-\mptfz.
%---------begin item----------
\begin{enumerate}
\item  Le module $r_hP=P\ep{h}$ est un
$\gA[1/r_h]$-module \pro de rang $h$.
\item Le module $P$ est  somme directe des $P\ep{h}$.
\item L'idéal $\gen{r_0}$ est l'annulateur du \Amo $P$.
\item  Pour $h>0$,   $P\ep{h}= \so{0}$ implique $r_h=  0.$
\end{enumerate}
%---------end item----------
\end{theorem}
%--- end-propdef-----------------
%---------begin proof----------
\begin{proof}
\emph{1.} Localiser en $r_h$: on obtient
$\; \rR{P\ep{h}}(X)=_{\gA[1/r_h]}\rR{P}(X) =_{\gA[1/r_h]}X^h$.

 \emph{2.} Car les $r_h$ forment un \sfioz.

 \emph{3.} Déjà vu (proposition~\ref{propAnnul}).

 \emph{4.} Résulte \imdt du point \emph{3.}
%---------end proof----------
\end{proof}

Notez que, sauf si $r_h=1$ ou $h=0$, le module $r_hP$ n'est pas de rang constant
en tant que \Amoz.

Le \tho précédent donne une \dem \gui{structurelle} du \thrf{th.ptf.idpt}.

%--- Remark{rem th decomp ptf}------
\medskip
\rem
\label{rem th decomp ptf}\relax
Si $P$ est (isomorphe à) l'image d'une \mprn $F$ les \idms
$r_k=\ide_k(P)$ attachés au module $P$ peuvent être reliés au \polcar de  la
matrice $F$ comme suit:

\snic{\det(X\In-F)=\som_{k=0}^n r_kX^{n-k}(X-1)^k.
}

%\sni
(Notez que les $X^{n-k}(X-1)^k$ forment une base du module des \pols de degré~$\leq n$,
triangulaire par rapport à la base usuelle.) \eoe

%%%%%%%%%%%%%%%%%%%%%%%%%%%%%%%%%%%%%%%%%%%%%%%%%%%%%%%%%%%%%%%%%%%%%%%%%%%
%:  subsec{Comparaison des deux rangs}-
\subsec{\Pol rang et \idfsz}%{Comparaison des deux rangs}
\label{subsecRangRang}

%:--- Theorem{corth.ptf.sfio}---

La démonstration du \tho \ref{corth.ptf.sfio} qui suit
s'appuie sur le \tho
\ref{prop Fitt ptf 2}, qui affirme qu'un \mptf devient libre après \lon en des
\ecoz.

Nous avons placé ce \tho ici car il répond aux questions que l'on se pose
naturellement après le \thrf{th ptf sfio}.
D'abord, vérifier qu'une \mprn est de rang $k$ \ssi son image est un \mrc $k$. 
Plus \gnltz, caractériser le \sfio
qui intervient dans le \polmu en termes des \idfs du
module.

En fait, on peut donner
une \dem alternative du \tho \ref{corth.ptf.sfio} sans passer par un argument de
\lonz, en s'appuyant sur les puissances extérieures (voir la proposition~\ref{prop puissance ext}).
%H La phrase suivante est améliorée
%\Llec pourra vérifier que cette \dem alternative
%ne fait pas appel au \tho \ref{prop Fitt ptf 2} ni à ses conséquences.

%H L'ordre des theoremes ayant changé j'ai du changer l'explication qui suit
Signalons que pour un \mpf $M$ l'\egt $\cF_h(M)=\gen{1}$ signifie que $M$ est
localement engendré par $h$ \elts (on a vu cela pour le cas~$h=1$ dans le
\thrf{propmlm}, dans le cas \gnlz,
voir le lemme du nombre de \gtrs local
\paref{lemnbgtrlo} et la \dfnz~\ref{deflocgenk})
%H: rappel inutile:  et qu'une matrice $A$ est de rang
%$\geq h$ lorsque $\cD_h(A)=\gen{1}$.

%:  theorem   corth.ptf.sfio
%:HHH 11 titre changé
\begin{theorem}
\label{corth.ptf.sfio}\relax %\emph{(\Pol rang et \idfsz)}
\emph{(Structure locale et \idfs d'un \mptfz, 2)}\\
Soient $F\in\GAq(\gA)$, $P\simeq\Im  F$ et
$\rR{P}(X)=\sum_{i=0}^qr_iX^i$.
%-----------------begin enum------------------
\begin{enumerate}
\item  Posons  $S(X)=\rR{P}(1+X)=1+u_1X+\cdots+u_qX^q$ ($u_h$ est
la somme des mineurs principaux d'ordre $h$ de la matrice $F$).\\
On a pour tout $h\in\lrb{0..q}$:

\snic{
\formule{\cD_h(F)=\gen{r_h+\cdots+r_q}=\gen{r_h,\ldots ,r_q}= \gen{u_h,\ldots
,u_q}\\[1mm]
\cF_h(P)= \gen{r_0+\cdots +r_h}= \gen{r_0,\ldots ,r_h}}
}

\item  En particulier:%
%H l'endroit ou se trouve la definition ne me semble
% pas devoir etre rappelé  (voir \dfnz~\ref{defRangk})
%:HHH  equivalences avec des  <=>  plutot que des phrases
%-----------------begin enum------------------
\begin{enumerate}
\item  $\rg(F)=h\iff \rg(P)=h$, 
%Le rang de $F$ est $h$ \ssi le module $P$ est de rang $h$ 
%(i.e., $r_h=1$),  \ssi $\cF_{h-1}(P)=\gen{0}$ et
%$\cF_{h}(P)=\gen{1}$.
\item  $\rg(F)\leq h\iff \deg \rR{P}\leq h,$
%Le rang de $F$ est $\leq h$ \ssi $\deg \rR{P}\leq h,$ \ssi
%$\cF_{h}(P)=\gen{1}$.
\item  $\rg(F)> h\iff r_0=\cdots=r_h=0 \iff \cF_{h}(P)=0.$
%Le rang de $F$ est $> h$ \ssi  $r_0=\cdots=r_h=0,$  \ssi
%$\cF_{h}(P)=\gen{0}$.
\end{enumerate}
%-----------------end enum------------------
\end{enumerate}
%-----------------end enum------------------
\end{theorem}
%--- end-theorem------------------------------------
%-----------------begin proof------------------
\begin{proof}
L'\egt $ \gen{u_h,\ldots,u_q}=\gen{r_h,\ldots ,r_q}$ résulte
 des \egts

\snic{S(X)=\rR{P}(1+X)\;\hbox{ et }\;\rR{P}(X)=S(X-1).}

%\sni
Pour vérifier les \egts 
$\cD_h(F)=\gen{r_h+\cdots+r_q}=\gen{r_h,\ldots ,r_q}$ et

\snic{\cD_{q-h}(\I_q-F)= \gen{r_0+\cdots +r_h}= \gen{r_0,\ldots ,r_h},}

%\sni
il suffit de le faire après \lon
en des \ecoz. Or le noyau et l'image de $F$  deviennent libres après \lon en des \eco (\thref{theoremIFD} ou \thref{prop Fitt ptf 2}), et la matrice devient
donc semblable à une \mprn standard.  
%\\
% Les \egts $\gen{r_h,\ldots ,r_q}=\gen{r_h+\cdots+r_q}$ et
%$\gen{r_0+\cdots +r_h}= \gen{r_0,\ldots ,r_h}$
%résul\-tent du fait que les $r_j$ sont des \idms \ortsz. Il est clair que
%$u_h,\ldots ,u_q\in\cD_h(F)$.
%Pour montrer l'inclusion réciproque $\cD_h(F)\subseteq\gen{r_h,\ldots ,r_q}$
%on peut localiser en des \eco de façon à ce que les modules $\Im F$
%et $\Ker F$ deviennent
% libres dans chacun
%des localisés (\thrf{prop Fitt ptf 2}). 
%Alors la matrice $F$ devient semblable
%à une \mprn standard et le résultat est clair.\\
%De même l'\egt $\cF_h(P)= \gen{r_0,\ldots ,r_h}$
%(qui signifie  $\gen{r_0,\ldots ,r_h}=\cD_{q-h}(\I_q-F)$) est claire lorsque
%la matrice est une \mprn standard. Or
%il suffit de la vérifier après \lon en des \ecoz.
\end{proof}
%-----------------end proof------------------

%%%%%%%%%%%%%
\section{Propriétés de caractère fini}
\label{secSalutFini}\index{propriété de caractère fini}

Cette section  veut illustrer
l'idée que les bons concepts en algèbre sont ceux qui
sont contrôlables par des procédures finies. 

Nous avons en vue de mettre en évidence des \gui{bonnes
\prtsz}. Il y a naturellement celles qui acceptent de se soumettre au principe local-global: pour que la \prt soit vraie il faut et suffit qu'elle le soit après \lon en des \mocoz. C'est un phénomène que nous avons déjà beaucoup rencontré, et qui continuera par la suite.

Rappelons qu'une \prt est dite \gui{\carfz}
si elle est conservée par \lon
(par passage de $\gA$ à $S^{-1}\gA$) et si, lorsqu'elle est
vérifiée après \lon en  $S$, alors elle est vérifiée après \lon en $s$ pour un certain $s\in S$.

Dans le  fait\eto \ref{factPropCarFin} nous avons démontré
en \clama que pour les \prts
 \carfz, le \plgc (\lon en des \mocoz)
est \eqv au \plga (\lon en tous les \idemasz).
Par contre, une \prco du \plgc contient une information plus précise
a priori qu'une preuve classique du \plgaz.

%:     Proposition{propFiniBon1}
\begin{proposition}\label{propFiniBon1} Soit  $S$ un \mo de $\gA$.
\perso{à compléter éventuellement en ajoutant d'autres items}
\begin{enumerate}
  \item Soit $AX=B$ un \sli sur $\gA$. Alors, s'il admet une
  solution dans $\gA_{S}$, il existe $s\in S$ tel qu'il admette une solution
  dans~$\gA_{s}$.
  \item Soient $M$ et $N$ deux sous-\Amos d'un même module,
\emph{avec $M$ \tfz}. Alors, si $M_{S}\subseteq N_{S}$, il existe $s\in S$ tel que
$M_{s}\subseteq N_{s}$.
% 3
  \item \label{propFiniBon1-3}\relax
  Soient $\gA$ un  \emph{\coriz}, $M$, $N$, $P$ des \Amos \emph{\pfz}, et deux \alis $\varphi:M\to N$, $\psi:N\to P$. \\
Si la suite
$M\vers{\varphi}N\vers{\psi}P$ devient exacte après \lon en $S$  il existe $s\in S$ tel que la suite devienne exacte après \lon en $s$.
  \item \label{propFiniBon1-4}\relax Soient $M$ et $N$ deux \Amos \emph{\pfz}. \hbox{Si $M_{S}\simeq N_{S}$,} il existe $s\in S$ tel que
$M_{s}\simeq N_{s}$.
  \item Soit $M$ un \Amo  \pfz. Si $M_S$ est libre, il existe
un $s \in S$ tel que $M_s$ soit libre. De même, si $M_S$ est \stlz,
il existe un $s \in S$ tel que $M_s$ soit \stlz.
  \item Si un \mpf devient \pro après \lon en
$S$, il  devient \pro après \lon en
un \elt $s$ de $S$.
%
%   \item
\end{enumerate}
\end{proposition}
\begin{proof}
Montrons le point \emph{\ref{propFiniBon1-3}}.
On trouve d'abord un $u\in S$ tel que $u\,\psi \big(\varphi(x_{j})\big)=0$
pour des \gtrs $x_{j}$ de $N$. On en déduit que $\psi\circ\varphi$ devient nul après \lon en $u$. Par ailleurs,  les hypothèses assurent  que
$\Ker\psi$ est \tfz. Soient $\yn$ des \gtrs de  $\Ker\psi$.
Pour chacun d'eux on trouve un $z_{j}$ dans $N$ et un $s_{j}\in S$ tels
que $s_{j}(\varphi(z_{j})-y_{j})=0$. On prend pour $s$ le produit de $u$
et des $s_{j}$.\\
Montrons le point \emph{\ref{propFiniBon1-4}}. 
Soient $G$ et $H$ des \mpns pour $M$ et $N$.
Notons $G_{1}$ et $H_{1}$ les deux matrices données dans le lemme \ref{lem pres equiv}.
Par hypothèse il existe deux matrices carrées $Q$ et $R$ à \coes dans $\gA$
telles que $v=\det(Q)\det(R)\in S$ et
$
Q\,G_{1}=_{\gA_{S}} H_{1}\,R
$.
Ceci signifie que l'on a sur $\gA$ une \egt
$$
w\,(Q\,G_{1}- H_{1}\,R)=0,\quad w \in S
.$$
Il suffit donc de prendre $s=vw$.
\end{proof}

On a vu que l'\eds se comporte bien par rapport aux produits
tensoriels, aux puissances extérieures et aux puissances symétriques.
Pour le foncteur $\Lin_\gA$ les choses ne se passent pas toujours aussi bien. Des
résultats importants pour la suite sont les suivants.
%:  proposition{fact.hom egaux}
\begin{proposition}
\label{fact.hom egaux}\relax
Soient $f:M\rightarrow N$ et $g:M\rightarrow N$ deux \alis entre \Amosz, {\em
avec $M$ de type fini}.  Alors, $f_S=g_S$ \ssi il existe $s\in S$ tel que
$sf=sg$.
En d'autres termes, l'application canonique
$\big(\Lin_\gA(M,N)\big)_S\rightarrow \Lin_{\gA_S}(M_S,N_S)$
est injective.
\end{proposition}
%--- end-fact--------------------------
%:   proposition{fact.homom loc pf} --------
\begin{proposition}\label{fact.homom loc pf}\relax
Soient $M$ et $N$ deux \Amos  et $\varphi:M_S\rightarrow N_S$ une \Aliz. On suppose que  $M$ est  \pfz, ou
que~$\gA$ est intègre, $M$ \tf et $N$ sans torsion (i.e. $a\in\gA$, $x\in N$, $ax=0$ impliquent $a=0$ ou $x=0$).
\\ 
Alors, il existe une \Ali $\phi:M\rightarrow N$ et un $s\in S$ tels que
%----begin $$------------------
$$\forall x\in M\quad  \varphi(x/1)  = \phi(x)/s ,$$
%----end $$------------------
et l'application canonique
$\big(\Lin_\gA(M,N)\big)_S\rightarrow \Lin_{\gA_S}(M_S,N_S)$ est bijective.
\end{proposition}
%---------
%----begin{proof------------------
\begin{proof}
Le deuxième cas, facile, est laissé \alecz. Pour suivre la \dem du premier cas il faut
regarder la figure ci-après.
\begin{figure}[htbp]
\centerline{
\xymatrix{
\Ae m \ar[d]_{\displaystyle g} \ar[rrr]^{\displaystyle j_m} & &
&\gA_S^m\ar[d]^{\displaystyle {g_S}} &
\\
\Ae q \ar[d]_{\displaystyle \pi}
   \ar@{-->}[rdd]^(.4){\displaystyle \Psi}\ar[rrr]^{\displaystyle j_q} & &
&\gA_S^q\ar[d]_{\displaystyle\pi_S}\ar[rdd]^{\displaystyle \psi} &
\\
M\ar[rrr]
% |!{[u];[rd]}\hole
^{\displaystyle j_M}\ar@{-->}[rd]_{\displaystyle \phi}
& & &M_S\ar[rd]_{\displaystyle \varphi} &
\\
 & N\ar[rrr]^{\displaystyle j_N} & & &N_S}
}
\caption{Localisation des \homos}
\label{fig1}
\end{figure}
Supposons que $M$ est le conoyau de l'\ali $g:\Ae m\rightarrow \gA^q$ avec une
matrice $G=(g_{i,j})$ par rapport aux bases canoniques, alors d'après le fait
\ref{fact.sexloc} le module~$M_S$ est le conoyau de l'\ali $g_S:\gA_S^m\rightarrow
\gA_S^q$, représentée par la matrice $G_S=(g_{i,j}/1)$ sur les bases canoniques.
On note 

\snac{\Ae m\vers{j_m} \gA_S^m,\, \gA^q\vers{j_q} \gA_S^q$,
$M\vers{j_M} M_S$, $N\vers{j_N} N_S$, $\gA^q\vers{\pi} M$,
$\gA_S^q\vers{\pi_S} M_S,}

%\sni
les applications canoniques. Soit $\psi:=\varphi
\circ \pi_S$, de sorte que $\psi\circ g_S=0$.\linebreak 
Donc $\psi\circ g_S\circ j_m=0=\psi \circ j_q \circ g$. 
Il existe un  $s\in S$,
dénominateur commun pour les images par $\psi$ des vecteurs de la base canonique. D'où une \ali $\Psi:\gA^q\rightarrow N$ avec 
$(s\psi)\circ j_q=j_N\circ \Psi$.\\ 
Ainsi,~$j_N\circ\Psi\circ  g= s(j_m\circ g_S\circ \psi)=0$.
D'après la proposition~\ref{fact.hom egaux} 
appliquée à~$\Psi\circ g$,  l'\egt
$j_N\circ(\Psi\circ  g)= 0$ dans $N_S$ implique qu'il existe $s'\in S$ tel que
$s'(\Psi\circ  g)=0$. Donc $s' \Psi$ se factorise sous forme $\phi\circ \pi$.
On obtient alors 

\snic{(ss'\varphi)\circ j_M\circ \pi=ss'(\varphi \circ\pi_S\circ j_q)
= ss'\psi\circ j_q=s'j_N\circ\Psi =
j_N\circ \phi\circ \pi,}

%\sni
 et puisque $\pi$ est surjective,
$ss'\varphi \circ j_M=j_N\circ\phi$. Ainsi, pour tout $x\in M$,
on~a~$\varphi(x/1)=\phi(x)/ss'$.
\end{proof}
%----end{proof------------------

%:     Corollary{corfinchapMPF}
\begin{corollary}\label{corfinchapMPF}
Supposons que $M$ et $N$ sont \pfz, ou qu'ils sont \tf sans torsion et
que $\gA$ est intègre. Si $\varphi : M_S \to N_S$ est un \isoz, il existe $s \in S$
et un \iso $\psi : M_s \to N_s$ tel que $\psi_S = \varphi$.
\end{corollary}
\begin{proof}
Soit $\varphi' : N_S \to M_S$ l'inverse de $\varphi$.  D'après la
proposition précédente,
il existe $\phi : M \to N$, $\phi' : N \to M$, $s \in
S$, $s' \in S$ tel que $\varphi = \phi_S/s$,
$\varphi' = \phi'_S/s'$.
Posons $t = ss'$ et définissons $\psi = \phi_t/s : M_t \to N_t$,
$\psi' = \phi'_t/s' : N_t \to
M_t$.  Alors, $(\psi' \circ \psi)_S$ est l'identité sur $M_S$,
et $(\psi \circ\psi')_S$ est l'identité sur $N_S$. On en déduit l'existence d'un $u \in S$
tel que $(\psi' \circ \psi)_{tu}$ est l'identité sur $M_{tu}$,
et~$(\psi \circ\psi')_{tu}$ est l'identité sur $N_{tu}$. En conséquence, $\psi_{tu} : M_{tu} \to N_{tu}$ est un \iso tel que $(\psi_{tu})_S = \varphi$.
\end{proof}

%:section: --- Exercices
\Exercices{

%--- Exercise{exoptf0Lecteur}-------------
\rdb\begin{exercise}
\label{exoptf0Lecteur}
{\rm  Il est recommandé de faire les \dems non données, esquissées,
laissées \alecz,
etc\ldots\,
 On pourra notamment traiter les cas suivants.
\begin{itemize}
\item Montrer les faits \ref{factdefiMPRO} et \ref{factMatriceAlin}.
\item Vérifier les détails du lemme \ref{lem calculs}.
\item \label{exofact.homegaux}\relax
Montrer le fait \ref{fact.hom egaux} ainsi que le
deuxième cas dans la proposition \ref{fact.homom loc pf}.
\end{itemize}
}
\end{exercise}
%--- end -exercise-----------------------------------------

%--- Exercise{exoProjmemeImage}-------------
\rdb\begin{exercise}
 \label{exoProjmemeImage} (Projecteurs ayant même image)\\
 {\rm  Soient $a$, $c$ dans un anneau $\gB$ non \ncrt commutatif.
 \Propeq
\begin{itemize}
\item  $ac=c$ et $ca=a$.
\item $a^2=a$, $c^2=c$ et $a\gB=c\gB$.
\end{itemize}
Dans un tel cas on pose $h=c-a$ et $x=1+h$.
Montrer les résultats suivants. 

\snuc{ha=hc=0$, $ah=ch=h$, $h^2=0$, $x\in\gB\eti$, $ax=c$,
$xa=x^{-1}a=a$ et $\fbox{$x^{-1}ax=c$}.}

%\sni
On notera en passant que l'\egt $ax = c$
%:2012 ci-dessous \gB et non \gA
redonne l'\egt $a\gB = c\gB$.
\\
Cas particulier.  $\gA$ un anneau commutatif, 
 $M$ un \Amoz, et $\gB=\End_\gA(M)$: deux  \prrs qui ont même image sont semblables.

} \end{exercise}
%--- end-exercise-----------------------------------------

%--- Exercise{exo2.4.1}-----------
\rdb\begin{exercise}
\label{exo2.4.1} (Deux \prrs \eqvs sont sembla\-bles) \relax\\
{\rm  Soient dans un anneau $\gB$ non \ncrt commutatif, deux \idms  \emph{\eqvsz} ($a^2= a$, $b^2=b$, $\exists p,\,q\in\gB\eti,\, b=paq$).
On va montrer qu'il sont \emph{conjugués} ($\exists d\in\gB\eti, \;dad^{-1}=b$).
\begin {itemize}
\item
Dans cette question, $a$, $b \in \gB$ sont \eqvs ($b=paq$), mais ne sont pas supposés
\idmsz. %si $p, q$ \ivs sont tels que $b=paq$,
Montrer que l'\elt $c=p^{-1}bp$ vérifie
$a\gB = c\gB$.

\item En particulier, si $b$ est \idmz, $c$ est un \idm
conjugué de $b$ qui vérifie $a\gB = c\gB$.
Conclure en utilisant l'exercice précédent.
\end {itemize}
Cas particulier.  $\gA$ un anneau commutatif, 
 $M$ un \Amoz, et $\gB=\End_\gA(M)$: deux
 \prrs de $M$ \eqvs  
%:HHH  \mprns \eqves  (M n'est pas forcement libre)
 sont semblables.
}
\end{exercise}
%--- end-exercise-----------------------------------------

%--- Exercise{exoSchanuelVariation}-------------
\rdb\begin{exercise}\label{exoSchanuelVariation} {(Une cons\'equence importante du lemme de Schanuel \ref{corlemScha})}
\\
{\rm
\emph{1.} On consid\`ere deux suites exactes:
$$\preskip.4em \postskip.2em
\arraycolsep2pt
\begin{array}{ccccccccccccccc}
0 &\rightarrow& K& \to & P_{n-1}& \to & \cdots & \to & P_1 & \vers{u} & P_0
& \to & M& \rightarrow& 0
\\
0 &\rightarrow& K'& \to & P'_{n-1}& \to & \cdots & \to & P'_1 & \vers{u'} & P'_0
& \to & M&  \rightarrow& 0
\end {array}
$$
avec les modules $P_i$, $P'_i$ \prosz. Alors, on obtient un \isoz:
$$\preskip.4em \postskip.4em\ndsp 
K \oplus \bigoplus\limits_{i \equiv n-1 \mod 2} P'_i \oplus
\bigoplus\limits_{j \equiv n \mod 2} P_j \quad \simeq\quad 
K' \oplus \bigoplus\limits_{k \equiv n-1 \mod 2} P_k \oplus
\bigoplus\limits_{\ell \equiv n \mod 2} P'_\ell 
$$

\sni
\emph{2.} En déduire que si l'on a une suite exacte
où les $P_i$, $i\in \lrbn$ sont \prosz

\snic {
0 \rightarrow P_n \to  P_{n-1} \to  \cdots  \to  P_1 \to P_0 \to M
\rightarrow 0,
}

%\sni
alors, pour toute suite exacte

\snic {
0 \to K' \to  P'_{n-1} \to  \cdots  \to  P'_1  \to  P'_0 \to M
\to 0,
}

%\sni
 où les $P'_i$ sont \prosz, le module $K'$ est \egmt \proz.
}
\end {exercise}
%--- end-exercise-----------------------------------------

%--- Exercise{exoSuitExPtfs}-------------
\rdb\begin{exercise}
\label{exoSuitExPtfs}
{\rm
On considère une suite exacte entre modules projectifs de type fini

\snic{
0  \lora
P_n \vvers {u_n} 
P_{n-1} \vvvers {u_{n-1}}  P_{n-2}
\lora \cdots \lora
P_2 \vvers {u_2} 
P_{1} \lora   0}

%\sni
Montrer que
${  {\bigoplus\limits_{i\rm\ impair} P_i \;\simeq\;
\bigoplus\limits_{j\rm\ pair} P_j}}.$

  En déduire que si  les $P_i$ pour $i \ge
2$ sont \stlsz, il en est de même de~$P_1$.
}
\end{exercise}
%--- end -exercise-----------------------------------------

%--- Exercise{exoZerDiR}----------
\rdb\begin{exercise}
\label{exoZerDiR}
{\rm  Montrer que %pour un anneau $\gA$ 
\propeq
%-----------------begin enum------------------
\begin{itemize}
\item L'anneau $\gA$  est \zedrz.
\item Les $\gA$-\mpfs sont toujours \ptfsz.
\item Tout module $\aqo{\gA}{a}$ est \ptfz.
\end{itemize}
%-----------------end enum------------------
(autrement dit, montrer la réciproque pour le point \emph{1} dans le
\thrf{propZerdimLib}).
}
\end{exercise}
%--- end-exercise-----------------------------------------

%--- Exercise{exoProjRang1}-------------
\rdb\begin{exercise}\label{exoProjRang1} {(Projecteurs de rang $1$, voir la proposition \ref{pmlm})}
\\
{\rm
Soit $A = (a_{ij}) \in \Mn(\gA)$. On étudie différents \syps en les  $a_{ij}$ dont l'annulation définit la sous-\vrt $\GA_{n,1}(\gA)$ de~$\Mn(\gA)$.
On note $\cD'_2(A)$ l'\id engendré par
les mineurs ayant au moins l'un des \gui {quatre coins} sur la diagonale (à ne pas confondre avec les mineurs principaux, sauf si $n=2$).

\emph {1.}
Si $A$ est un \prr de rang $\le 1$, alors $\Ann\, A$ est
engendré par  $1-\Tr A$ (\idmz). En particulier, un \prr de rang $1$ est de
trace~$1$.

\emph {2.}
Les \egts $\Tr A = 1$ et $\cD'_2(A) = 0$ impliquent $A^2 = A$ et $\cD_2(A) = 0$.
Dans ce cas, $A$ est un \prr de rang $1$ 
(mais on peut avoir $\Tr A = 1$ et $A^2 = A$ sans avoir
$\cD_2(A) = 0$, par exemple pour un \prr de rang $3$ sur un anneau \hbox{où $2=0$}). 
En conséquence, pour une matrice~$A$ quelconque on a
$$\preskip.3em \postskip.25em 
\gen {1-\Tr A} + \cD_1(A^2 - A) \subseteq
\gen {1-\Tr A} + \cD'_2(A) = \gen {1-\Tr A} + \cD_2(A) 
$$
sans avoir \ncrt l'\egt à gauche.

\emph {3.}
On considère le \pol  $\det\!\big
(\In+(X-1)A\big)$ (si $A\in\GA_n(\gA)$, c'est le \pol rang du module $P=\Im A$) et l'on note $r_1(A)$ son
\coe en $X$. On a alors l'\egt des trois \ids suivants, qui définissent
la sous-\vrt $\GA_{n,1}(\gA)$ de~$\Mn(\gA)$:
$$\preskip.2em \postskip.3em 
\gen {1-\Tr A} + \cD'_2(A) = \gen {1-\Tr A} + \cD_2(A) =
\geN {1-r_1(A)} + \cD_1(A^2 - A). 
$$
Préciser le cardinal de chaque \sgrz.
%On obtient
%ainsi trois \syss d'\eqns en $n^2$ variables qui définissent la
%\vrt $\GA_{n,1}$ des \prrs de rang 1.

}
\end {exercise}
%--- end -exercise-----------------------------------------

%--- Exercise{exoMatProjRang1CoeffReg}-------------
\rdb\begin{exercise}\label{exoMatProjRang1CoeffReg}
{(Projecteur de rang 1 ayant un \coe \ndzz)}\\
{\rm 
Soit $A = (a_{ij}) \in \GA_n(\gA)$ un \prr de rang $1$,
$L_i$ sa ligne $i$, $C_j$ sa colonne~$j$.

 \emph {1.}
Fournir une preuve directe de l'\egt matricielle $C_j \cdot L_i = a_{ij}A$.
En remarquant que $L_i\cdot C_j = a_{ij}$, en déduire l'\egt d'\ids $\gen
{L_i}\gen{C_j} = \gen {a_{ij}}$.

 \emph {2.}
On suppose $a_{ij}$ \ndz; donc $\gen {L_i}$ et $\gen{C_j}$ sont des \ids
\ivsz, inverses l'un de l'autre.  Fournir une preuve directe de l'exactitude
au milieu de la suite:

\snic {
\Ae n\vvvers{\In-A} {\Ae n}\vvers{L_i} \gen {L_i}\to 0
}

%\sni
et par conséquent $\gen {L_i} \simeq \Im A$.

 \emph {3.} Montrer que la matrice $A$ est entièrement déterminée
par $L_i$ et $C_j$. Plus \prmtz, si l'anneau $\gA$ est à \dve explicite:
\begin{itemize}
\item calculer la matrice $A$,
\item en déduire à quelle condition une ligne $L$
et une colonne $C$ peuvent être la ligne $i$ et la colonne $j$ d'une \mprn
de rang $1$ (on suppose que le \coe commun en position $(i,j)$ est \ndzz). 
%
%\item  
\end{itemize}
 
\emph {4.}
Soit $C \in \Im A$, $\tra L \in \Im \tra A$ et $a = L\cdot C$.
Montrer l'\egt matricielle $C \cdot L = aA$ et en déduire l'\egt d'\ids $\gen
{L}\lra{\tra C\,} = \gen {a}$. Si $a$ est \ndzz, les \ids $\gen {L}$ et~$\lra {\tra C\,}$
sont \ivsz, inverses l'un de l'autre, $\gen {L} \simeq \Im A$ et
$\lra {\tra C\,} \simeq \Im \tra A$.

}
\end {exercise}
%--- end -exercise-----------------------------------------

%--- Exercise{exoFittTfPtf}----------
\rdb\begin{exercise}
\label{exoFittTfPtf}
{\rm Si un \Amo \emph{\tfz} a ses \idfs engendrés par des \idmsz, il est \ptfz.
}
\end{exercise}
%--- end-exercise-----------------------------------------

%:--- Exercise{exoRelateursCourts}-------------
\rdb\begin{exercise}
\label{exoRelateursCourts}  (Syzygies courtes)\\
{\rm
\emph{Notations, terminologie.}
On note $(e_1,\ldots,e_n)$  la base canonique de $\Ae n$.\\
Soient $x_1$, \dots, $x_n$ des \elts d'un \Amoz. On note $x=\tra[\,x_1\; \cdots\; x_n\,]$ \hbox{et $x^\perp := \Ker (\tra {x})\subseteq \Ae n$} le
module des \syzys entre les~$x_i$.
\\
 On dira d'une \syzy $z \in x^\perp$
qu'elle est \gui {courte} si elle possède au plus deux \coos non
nulles, i.e. si $z \in \gA e_i \oplus \gA e_j$ $(1\leq i\neq j\leq n)$.
\begin{enumerate}\itemsep=0pt
\item 
Soit $z \in x^\perp$. Montrer que la condition  \gui{$z$ est somme de
\syzys courtes} est une condition \linz. En conséquence, si
$z$ est \gui{localement} somme de \syzys courtes, elle l'est globalement.

\item 
En déduire que si $M=\sum \gA x_i$ est un \mlmz, alors tout \elt de
 $x^{\perp}$ est somme de \syzys courtes.

\item 
Si toute \syzy entre trois \elts de $\gA$
est somme de \syzys courtes, alors $\gA$ est un \emph{\anarz},
i.e., tout \id $\gen{x,y}$ est \lopz.

\item 
Dans la question \emph{2} donner une solution
globale en utilisant une \mlmo $A = (a_{ij}) \in \Mn(\gA)$
 pour $x$.
\end{enumerate}
}
\end{exercise}
%--- end -exercise-----------------------------------------

%--- Exercise{exoPetitsRelateurs}-------------
\rdb\begin{exercise}
\label{exoPetitsRelateurs} (Syzygies triviales)\\
{\rm On utilise les notations de l'exercice \ref{exoRelateursCourts}. Maintenant
$x_1$, \ldots, $x_n \in \gA$. \\
Pour $z \in \Ae n$ on note 
$\scp{z}{x} \eqdefi \sum z_i x_i$. Le module
des \syzys $x^\perp$  contient
les \gui {\syzys triviales} $x_j e_i - x_i e_j$ (qui sont un cas particulier
de \syzys courtes).
\\
Dans les deux premières questions, on montre que si $x$
est \umdz, alors $x^\perp$ est engendré par ces \syzys triviales. On
fixe $y \in \Ae n$ tel que $\scp {x}{y} = 1$.
\begin{enumerate}\itemsep=0pt
\item 
Rappeler pourquoi $\Ae n = \gA y \oplus x^\perp$.

\item 
Pour $1 \leq i < j \leq n$, on définit $\pi_{ij} \,: \,\Ae n \to \Ae n$
par

\snic{\pi_{ij}(z) = (z_i y_j - z_j y_i) (x_j e_i - x_i e_j),}

si bien que $\Im \pi_{ij} \subseteq x^\perp \cap (\gA e_i \oplus \gA
e_j)$. Montrer que $\pi = \sum_{i < j} \pi_{ij}$  est la projection
sur $x^\perp$ \paralm à $\gA y$. En déduire le résultat
sur les \syzys triviales. Voir aussi l'exercice~\ref{exoPlgb1}.

\end{enumerate}
On ne suppose plus que $x$ est \umdz.  Soit $M  \in \Mn(\gA)$ une
matrice alternée.
\begin{enumerate} \setcounter{enumi}{2}\itemsep=0pt
\item 
Montrer qu'en
posant $z = Mx$, on a $\scp {x}{z} = 0$. 
\item 
En quel sens, une matrice alternée est-elle \gui {somme de petites matrices
alternées}? Faire le lien avec la \dfn de $\pi_{ij}$ dans
la question \emph{2}.

\end{enumerate}
}
\end{exercise}
%--- end -exercise-----------------------------------------

%:--- Exercise{exoProjImLibre}-------------
\rdb\begin{exercise}
 \label{exoProjImLibre} (Matrices de \prn qui ont une image libre)\\
 {\rm  Soit $P \in \GAn(\gA)$ un projecteur dont l'image est libre de rang $r$;
d'après la proposition \ref{propImProjLib} il existe
$X \in \Ae {n \times r}$, $Y \in \Ae {r \times n}$
vérifiant $YX = \I_r$ et $P = XY$.

\emph{1.}
On demande d'expliciter le lemme \dlg (lemme \ref{propIsoIm}), autrement dit
de calculer $A \in \SL_{n+r}(\gA)$ (et son inverse) telle que
$$A^{-1} \Diag(0_r, P) A = \I_{r,n+r}.\eqno(*)$$

\emph{2.}
On suppose que $X=\tra Y$ (donc $P$ est
symétrique). \\
Vérifier que l'on peut imposer à $A$ d'être \gui{orthonormale}
i.e. $\tra A = A^{-1}$. \\
 Réciproquement, si $A \in \SL_{n+r}(\gA)$ est  orthonormale et vérifie~$(*)$, alors on
peut écrire $P = X\tra{X}$ avec $X \in \Ae {n \times r}$ et~$\tra{X} X =
\I_r$ (la matrice $P$ est donc symétrique).
} \end{exercise}
%--- end-exercise-----------------------------------------

%--- Exercise{exoStabLibRang1}-------------
\rdb\begin{exercise}
 \label{exoStabLibRang1} (Modules \stls de rang $1$)\\
 {\rm   \Demo directe que tout module \stl de rang $1$
 est libre (proposition \ref{propStabliblib}), 
 en utilisant la formule de Binet-Cauchy (exercice~\ref{exoBinetCauchy}).\\
On considère deux matrices $R\in\Ae {(n-1)\times n}$ et $R'\in\Ae {n\times (n-1)}$ avec $RR'=\I_{n-1}$. Montrer que $\Ker R$ est un module libre.
Conclure.
} \end{exercise}
%--- end-exercise-----------------------------------------

%--- Exercise{exoStablib}-------------
\rdb\begin{exercise}
 \label{exoStablib} (Vecteurs \umdsz, modules $M$ vérifiant $M\oplus\gA\simeq\Ae n$)\\ 
 {\rm
 Soient $x$, $y \in \Ae n$ deux vecteurs et $A \in \Mn(\gA)$ une matrice
de première colonne $x$. On construit une matrice $B \in \Mn(\gA)$
de la manière suivante: sa première ligne  est $\tra {y}$ et ses $n-1$
dernières lignes sont les $n-1$ dernières lignes de $\wi {A}$, matrice
adjointe de $A$.
\begin{enumerate}\itemsep=0pt
\item Montrer que $\det(B) = \det(A)^{n-2} \scp {x}{y}$ et que les $n-1$
dernières lignes de $B$ appartiennent à $x^\perp := \ker \tra {x}$.
\end{enumerate}
On suppose désormais que $\scp {x}{y} = 1$. On sait alors que les deux
modules \stls $x^\perp$ et $y^\perp$ sont duaux l'un de l'autre
(faits \ref{factStablib} et \ref{factStablibDual});
on détaille cette \prt de manière matricielle dans le cas où
$y^\perp$ est libre.
\begin{enumerate}\itemsep=0pt
\setcounter{enumi}{1}
\item  Rappeler pourquoi $\Ae n = \gA x \oplus y^\perp$ et
$\Ae n = \gA y \oplus x^\perp$.

\item On suppose que $\gA x$ possède un \sul libre dans
$\Ae n$. Montrer de manière matricielle qu'il en est de même de
$\gA y$ en construisant une matrice inversible $n \times n$ \gui{adaptée}
à la décomposition $\Ae n = \gA y \oplus x^\perp$.
\end{enumerate}
} \end{exercise}
%--- end-exercise-----------------------------------------

%--- Exercise{exoMatLocSym}-------------
\rdb\begin{exercise}\label{exoMatLocSym}  {(Matrice de \lon principale \smqz)}
\\
{\rm
Soit $(x_1, \ldots, x_n) \in \Ae n$ possédant une \mlp $A \in \Mn(\gA)$ \smqz.  On pose $\fa = \gen {x_1, \ldots,x_n}$; 
en utilisant l'\egt (\ref{eqpmlm}) de la proposition \ref{pmlm},
montrer que $\fa^2$ est principal et \prmtz:
$\fa^2 = \gen {x_1^2, \cdots, x_n^2} = \gen {x_1^2 + \cdots + x_n^2}$.
}
\end {exercise}
%--- end -exercise-----------------------------------------

%%%%%%%%%%%  exoPgcdPpcm  %%%%%%%%%%%
\begin {exercise}\label{exoPgcdPpcm} (Au sujet de $\gA/\fa \oplus \gA/\fb
\simeq \gA/(\fa \cap \fb) \oplus \gA/(\fa + \fb)$)\\
{\rm Voir aussi les exercices \ref{exogcdlcm} et \ref{exoSECSci}, et le corolaire \ref{corthAnar}.

 \emph{1.}
 Soient $\fa$, $\fb$ deux idéaux de $\gA$ vérifiant $1 \in (\fa : \fb) +
(\fb : \fa)$.
Expliciter $\theta \in \GL_2(\gA)$ vérifiant
$\theta(\fa \oplus \fb) = (\fa \cap \fb) \oplus (\fa + \fb)$.
En déduire que $\gA/\fa \oplus \gA/\fb
$ est isomorphe à~$\gA/(\fa \cap \fb) \oplus \gA/(\fa + \fb)$.

 \emph{2.}
Soient $a$, $b \in \gA$, $\fa = \gen {a}$, $\fb = \gen {b}$.  On suppose qu'il
existe $A \in \GL_2(\gA)$ telle que $A\cmatrix {a\cr b\cr} = \cmatrix
{*\cr0\cr}$. Montrer que $1 \in (\fb : \fa) + (\fa : \fb)$. Expliciter $d$ et $m$
tels que~$\fa \cap \fb = \gen {m}$, $\fa + \fb = \gen {d}$, ainsi qu'une
\eqvc matricielle entre~$\Diag (a,b)$ et~$\Diag (m,d)$.

 \emph{3.}
Soient $a$, $b \in \gA$ avec $a\in \gen{a^2}$. Montrer que $a$, $b$
vérifient les conditions de la question~\emph{2.}

 \emph{4.}
Soient $\fa$, $\fb$ deux \itfs tels que $\fa + \fb$ soit
\lopz.  Montrer: $1 \in (\fa : \fb) + (\fb :
\fa)$, $\fa \cap \fb$ est \tf et  $\fa\fb = (\fa
\cap \fb)(\fa + \fb)$.

}

\end {exercise}
%%%%%%%%%%%%%%%%%%%%%%%%%%%%%%%%%%%%%%

Les exercices qui suivent apportent quelques précisions sur le \deterz,
le \polcar et le \polfonz.

%--- Exercise{exoPrecisionsDet1}-----
\rdb\begin{exercise}
\label{exoPrecisionsDet1}
{\rm
Soient $M$ un $\gA$-\mptfz,  $e$ un \idm de~$\gA$,
$f=1-e$ et $\varphi$ un \endo de $M$.
Alors $M=eM\oplus fM$, de sorte que $eM$ et $fM$ sont \ptfsz.
On a aussi $\varphi(eM)\subseteq eM$, et en notant $\varphi_e:eM\rightarrow eM$
l'\endo induit par $\varphi$, on~a:
%--------------------begin array---------------

\snic{\begin{array}{c}
\det(\varphi_e)=f+e\det(\varphi)   \quad {\rm et } \quad
\det(e\varphi)=r_0f+e\det(\varphi)     \\[1mm]
\rF{e\varphi}(X)=\rF{\varphi}(eX)=\rF{\varphi_e}(X)=
f+e\rF{\varphi}(X)     \\[1mm]
\rC{\varphi_e}(X)=f+e\rC{\varphi}(X)     \\[1mm]
\rR{eM}(X)=f+e\rR{M}(X)     \\[1mm]
\end{array}}
%---------------------end array--------------

%\sni
En outre,  $e\,\det(\varphi)$ est le \deter de $\varphi_e$ en tant 
qu'\endo \linebreak du $\gA[1/e]$-module $eM$.
}
\end{exercise}
%--- end-exercise-----------------------------------------

%--- Exercise{exoPeticalculPolrang}-------------
\rdb\begin{exercise}
 \label{exoPeticalculPolrang}
 {\rm On considère le module quasi libre
$M=\bigoplus _{k\in\lrbn}(r_k\gA)^{k}$, 
 où les $r_k$ sont des \idms \ortsz. 
On a $M\simeq e_1\gA \oplus  
\cdots \oplus e_n\gA$ avec $e_k=\sum_{j=k}^nr_j$, et $e_k \divi e_{k+1}$ pour $k\in\lrb{1..n-1}$ (cf. lemme \ref{lem ide-div},
et  exercices \ref{exoIdmsSupInf} et  \ref{exoSfio}).
\\
On pose $r_0=1-\sum_{i=1}^nr_i$ et $s_k=1-r_k$.
\begin{itemize}
  \item [--] Rappeler pourquoi  $\rR{r_k\gA}(X)=s_k+r_kX$.
  \item [--] Montrer que $\rR{M}(X)=r_0+r_1X+\cdots +r_nX^n = \prod\nolimits_{k=1}^n (s_k+r_kX)^k$.

  \item [--] Vérifier cette \egt par un calcul direct.
\end{itemize}
}
\end{exercise}
%--- end-exercise-----------------------------------------

%--- Exercise{exoPrecisionsDet2}----
\rdb\begin{exercise}
\label{exoPrecisionsDet2} (Le \deterz, composante par composante)\\
{\rm
Soit $\varphi$  un \endo d'un \mptf $M$ ayant~$n$ \gtrsz. Soient $r_h=\ide_h(M)$
(pour $h\in\lrb{0..n}$) et $d=\det(\varphi)$. Notons~$\varphi\ep{h}$ l'\endo du
\Amo $M\ep{h}$  induit par $\varphi$,  
$d_h=r_hd$, $\delta_h=\det(\varphi\ep{h})$ et $s_h=1-r_h$.
% (ici~$\varphi\ep{h}$ est vu comme \endo du~\Amo $M\ep{h}$).

\emph{1.}  On a les \egts suivantes:  

\snic{d_0=r_0 ,\; \delta_0=1,\; \delta_h=s_h+d_h \;  \hbox{ et }  \; d=d_0+d_1+\cdots+d_n=\delta_1 \smalltimes
\cdots \smalltimes\, \delta_n.}

%\sni
\emph{2.} En outre, $d_h$ est le \deter de $\varphi\ep{h}$ dans $\gA[1/{r_h}]$
lorsque l'on voit $\varphi\ep{h}$ comme un \endo du~$\gA[1/{r_h}]$-module
$M\ep{h}$.

\emph{3.}  De la même manière, on a:

\snic{\rF{\varphi\ep{h}}(X)=s_h+r_h\rF{\varphi}(X) \quad {\rm et }\quad
\rC{\varphi\ep{h}}(X)=s_h+
r_h\rC{\varphi}(X).}

}
\end{exercise}
%--- end-exercise-----------------------------------------

%--- Exercise{exoPrecisionsDet3}-----
\rdb\begin{exercise}
\label{exoPrecisionsDet3} (\Pol \cara et \polfon en cas de rang constant)
{\rm
Soit $\varphi$  un \endo d'un mo\-du\-le  $M$ de rang cons\-tant $h$.
On montrera les résultats suivants.
\\
 Le \polcar de $\varphi$ est unitaire de degré $h$ et le
\polfon de $\varphi$
est de degré $\leq h$. Les homogénéises en degré $h$ de
$\rC{\varphi}(X)$ et $\rF{\varphi}(X)$ sont égaux respectivement à
$\det(X\Id_M-Y\varphi)$ et $\det(Y\Id_M+X\varphi)$. Autrement dit on a les deux \egts

\snic{ \rC{\varphi}(X)=X^h\rF{\varphi}(-1/X)\quad {\rm et }\quad
\rF{\varphi}(X)=(-X)^h\rC{\varphi}(-1/X).}

%\sni
En outre,  $\det(\varphi)=(-1)^h\rC{\varphi}(0)$ est égal au \coe en $X^h$
de $\rF{\varphi}(X)$.
}
\end{exercise}
%--- end-exercise-----------------------------------------

%--- Exercise{exoPrecisionsDet4}-----
\rdb\begin{exercise}
\label{exoPrecisionsDet4} (\Pol \cara et \polfonz, cas \gnlz)\\
{\rm
Soit $\varphi$  un \endo d'un \mptf $M$. Notons

\snic{\rF{\varphi}(X)=1+v_1X+\cdots+v_nX^n$  et
$\rR{M}(X)=r_0+r_1X+\cdots+r_nX^n.}

%\sni
Alors,  on a les \egts suivantes.

\snic{\arraycolsep2pt\begin{array}{rcl}
  r_hv_k&=&0\;\hbox{ pour }0\leq h< k\leq n,\\[1mm]
   \rC{\varphi}(X)   &=& r_0\,+\,\sum\nolimits_{1\leq h\leq n}
r_h\,X^h\,\rF{\varphi}(-1/X),   \\[1mm]
   \rF{\varphi}(-X)   &=&  r_0\,+\,\sum\nolimits_{1\leq h\leq n}
r_h\,X^h\,\rC{\varphi}(1/X),  \\[1mm]
    \det(\varphi -X\Id_M)  &=&   \rR{M}(-1)\, \rC{\varphi}(X), \\[1mm]
   \det(\varphi)   &=&   r_0+r_1v_1+\cdots+r_nv_n=\rR{M}(-1)\,\rC{\varphi}(0).
\end{array}}
}
\end{exercise}
%--- end-exercise-----------------------------------------

%: problemes

%--- problem{exoSuslinCompletableLemma}-------------
\rdb\begin{problem}\label{exoSuslinCompletableLemma}
 {(Complétion de vecteurs \umds: un résultat d\^u à \Susz)}\\
{\rm  
Un vecteur de $\Ae n$ est dit \emph {complétable} s'il est égal à la première
colonne d'une matrice de~$\GLn(\gA)$. Il est
alors \umdz. On veut montrer le résultat suivant.\index{completable@complétable!\vmd ---} 
 
\emph{Soient
$b\in\gA$ et $(\an)\in\Ae n$ tels que $(\ov{a_1}, \ldots, \ov{a_n})$ soit
complétable sur $\gA\sur{b\gA}$, alors $(\an, b^n)$ est complétable (sur
$\gA$)}. 
 
Par hypothèse, on a  $A$, $D \in \Mn(\gA)$ vérifiant $
A\, D \equiv \In \mod b$, avec $\vab{a_1}{\cdots\;a_n}$ pour première ligne de $A$.
On veut trouver une matrice de~$\GL_{n+1}(\gA)$ dont la première ligne soit 
$\vab{a_1}{\cdots\;a_n\;b^n}$.  Notons~$a = \det(A)$.

\emph {1.}
Montrer qu'il existe $C\in \Mn(\gA)$ telle que $\cmatrix {A &b\,\In\cr C &D\cr}
\in \GL_{2n}(\gA)$.

 Il s'agit maintenant de transformer le
coin haut-droit $b\,\In$ de la matrice ci-dessus en $B' := \Diag(b^n, 1,\dots,1)$.

\emph {2.}
Montrer que l'on peut écrire $B' = bE + aF$ avec $E \in \En(\gA)$ et
$F \in \Mn(\gA)$.

\emph {3.}
Vérifier que
$\cmatrix {A &b\,\In\cr C &D\cr} \cmatrix {\In &\wi AF\cr 0 &E\cr} =
\cmatrix {A &B'\cr C &D'\cr}$ avec $D' \in \Mn(\gA)$.

\emph {4.}
Montrer que $\cmatrix {A &B'\cr C &D'\cr}$ est \eqve
à une matrice $\cmatrix {A &B'\cr C &D''\cr}$ où
$D''$ a \hbox{ses $n-1$} dernières
colonnes nulles. En déduire l'existence d'une matrice
\iv dont la première ligne est  $\vab{a_1}{\cdots\;a_n\;b^n}$.  

\emph {5.}
Exemple (Krusemeyer). Si $(x,y,z) \in \Ae 3$ est \umdz,
 $(x,y,z^2)$ est complétable. Plus \prmtz, si
$ux + vy + wz = 1$, la matrice ci-dessous convient:

\snic {
\cmatrix {x&y&z^2 \cr v^2 &w-uv &-x-2vz\cr -w-uv &u^2&-y+2uz}.
}

%\sni
Quel est son \deter (indépendamment du fait que
$ux+vy+wz=1$)?

\emph {6.}
Plus \gnltz, on a le résultat suivant (Suslin): si $(a_0, \an)$
est \umdz, alors $(a_0, a_1, a_2^2, \ldots, a_n^n)$ est complétable.

\emph {7.}
%:h2013 \Sus
Montrer le résultat suivant (\Susz's $n!$ theorem): si $(a_0, \an)$ est
\umdz, alors pour des exposants $e_0$, $e_1$, \ldots, $e_n$ tels que $n!$ divise
$e_0\cdot e_1\cdots e_n$, le vecteur $(a_0^{e_0}, a_1^{e_1}, \ldots,
a_n^{e_n})$ est complétable.

}

\end {problem}
%--- end -problem----------

%--- problem{exo2SphereCompletableYengui}-------------
\rdb\begin{problem}\label{exoSphereCompletableYengui} {(La $n$-sphère quand $-1$ est une somme de $n$ carrés,  avec I.~Yengui)}\\
{\rm  
\emph {1.}
Soit $\gA$ un anneau dans lequel $-1$ est une somme de 2 carrés et $x_0$,
$x_1$, $x_2 \in \gA$ vérifiant $x_0^2 + x_1^2 + x_2^2 = 1$.
\begin {enumerate}\itemsep=0pt
\item [\emph {a.}]
Montrer que le vecteur $(x_0,x_1,x_2)$ est complétable en
considérant une matrice $M = \cmatrix {x_0 & u & a\cr x_1 & v & b\cr x_2 & 0
& c\cr}$ où $u$, $v$ sont des formes linéaires en $x_0$, $x_1$, $x_2$ \linebreak 
et $a$, $b$, $c$  sont
des constantes.

\item [\emph {b.}]
Donner des exemples d'anneaux $\gA$ dans lesquels $-1$ est une somme de 2 carrés.
\end {enumerate}

\emph {2.}
On suppose que $-1$ est une somme de $n$ carrés dans l'anneau $\gA$.
\begin {enumerate}\itemsep=0pt
\item [\emph {a.}]  On utilise la  {\mathrigid 2mu notation \smash{$A\sims{\cG }B$} \paref{NOTAfGg}.
Soit $(x_i)_{i\in\lrbzn}$ avec $x_0^2 + \cdots + x_n^2 = 1$.} Montrer que
$$\preskip-.6em \postskip.4em 
\tra {[\,x_0\;x_1\; \cdots\; x_n\,]} \sims{\EE_{n+1}} \tra {[\,1\; 0\; \cdots\; 0\,]}. 
$$
En particulier, $\tra {[\,x_0\;x_1\; \cdots\; x_n\,]}$ est complétable.

\item [\emph {b.}]
Soit $m \ge n$, $x_0$, $x_1$, \ldots, $x_m$ et $y_{n+1}$, \ldots, $y_m$ vérifiant
\smashbot{$\sum\limits_{i=0}^n x_i^2 + \sum\limits_{j=n+1}^m y_jx_j = 1$}. 
Montrer que $\tra {[\,x_0\;x_1\; \cdots\; x_m\,]} \sims{\EE_{m+1}} \tra {[\,1\; 0\; \cdots\; 0\,]}$. 
\end {enumerate}

\emph {3.}
On suppose qu'il existe $a \in \gA$ tel que $1 + a^2$ soit nilpotent. C'est le
cas si $-1$ est un carré dans $\gA$, ou si $2$ est nilpotent.
\begin {enumerate}\itemsep=0pt
\item [\emph {a.}]
Soient $x_0$, $x_1 \in \gA$ avec $x_0^2 + x_1^2 = 1$. Montrer
que $\crmatrix {x_0 & -x_1\cr x_1 & x_0} \in \EE_2(\gA)$.

\item [\emph {b.}]
Soient $x_0$, $x_1$, \dots,$x_n$ et $y_2$, \ldots, $y_n$ dans $\gA$ tels que $x_0^2 + x_1^2 +
\sum_{i=2}^n x_i y_i = 1$.  Montrer que $\tra {[\,x_0\;x_1\; \cdots\; x_n\,]}
\sims{\EE_{n+1}} \tra {[\,1\; 0\; \cdots\; 0\,]}$.

\item [\emph {c.}]
Soit $\gk$ un anneau, $\gk[\uX,\uY] = \gk[X_0, \Xn, Y_2, \ldots, Y_n]$ et 

\snic{f=1 - \big(X_0^2 + X_1^2 + \sum_{i=2}^n X_i Y_i\big).}

%\sni
On pose $\gA_n =
\gk[x_0, \xn, y_2, \ldots, y_n] = \aqo {\gk[\uX,\uY]}{f}$. Donner des exemples
pour lesquels, pour tout $n$, $\tra {[\,x_0\;x_1\; \cdots\; x_n\,]}$ est complétable 
sans que $-1$  soit un carré dans $\gA_n$.
\end {enumerate}
}

\end {problem}
%--- end -problem-----------------------------------------

}% fin des exos
%:   solutions d'exos
\penalty-2500
\sol{

%%%%%%%%%%%%%%%%%%%%%%%%%%%%%%%%%%%%%%%%%%%%%%%%%%%%%%%%%%%%%%%%%%%%%%%%%%%
\exer{exoSchanuelVariation}
\emph{1.} Par \recu sur $n$, le cas $n=1$ étant exactement le lemme de Schanuel
(corolaire \ref{corlemScha}). \`A partir de chaque suite exacte, on en
construit une autre de longueur un de moins
\snic {\arraycolsep4pt
\begin{array}{ccccccccccccc}
0 &\rightarrow& K& \to & P_{n-1}& \to & \cdots & \to &
P_1 \oplus P'_0 &   \smash{\vvvers{u \oplus \I_{{P_0}'}}}
& \Im u \oplus P'_0 & \to & 0
\\
0 &\rightarrow& K'& \to & P'_{n-1}& \to & \cdots & \to &
P_1 \oplus P'_0 &  \vvvers{u' \oplus \I_{P_0}}
& \Im u' \oplus P_0 & \to & 0
\end {array}
}

%\sni
Mais on a $\Im u \oplus P'_0 \simeq \Im u' \oplus P_0$ d'après le lemme de
Schanuel appliqué aux deux suites exactes courtes:

\snic {\arraycolsep4pt
\begin{array}{ccccccccc}
0 &\to& \Im u & \to & P_0 & \to & M & \to & 0 \\
0 &\to& \Im u' & \to & P'_0 & \to & M & \to & 0 \\
\end {array}
}

%\sni
On peut donc appliquer la \recu (aux deux longues suites exactes de
longueur un de moins), ce qui fournit le résultat demandé. 

 \emph{2.} Conséquence \imde de \emph{1.}

%%%%%%%%%%%%%%%%%%%%%%%%%%%%%%%%%%%%%%%%%%%%%%%%%%%%%%%%%%%%%%%%%%%%%%%%%%%
\exer{exoSuitExPtfs}
{
Montrons par \recu sur $i$ que $\Im u_i$ est un \mptfz. C'est vrai pour $i = 1$. Supposons le vrai pour $i\geq 1$; on a donc une
\ali surjective $P_i \vers{u_i}  \Im u_i$ où $\Im
u_i$ est \ptf et par conséquent $P_i \simeq \Ker  u_i \oplus \Im
u_i$. Mais $\Ker u_i = \Im u_{i+1}$ donc $\Im u_{i+1}$ est \ptfz.
De plus $P_i \simeq \Im u_{i} \oplus \Im u_{i+1}$. Ensuite
$$\arraycolsep4pt
\begin {array} {rcl}
P_1 \oplus P_3 \oplus P_5
\oplus \cdots&\simeq
& 
(\Im u_1 \oplus \Im u_2) \oplus (\Im u_3 \oplus \Im u_4) \oplus
 \cdots
\\
&\simeq & 
\! \Im u_1 \oplus (\Im u_2 \oplus \Im u_3) \oplus
(\Im u_4  \oplus \Im u_5) \oplus  \cdots
\\
&\simeq & 
\! P_2 \oplus P_4 \oplus P_6 \oplus \cdots
\end {array}
$$
}

%%%%%%%%%%%%%%%%%%%%%%%%%%%%%%%%%%%%%%%%%%%%%%%%%%%%%%%%%%%%%%%%%%%%%%%%%%%

\exer{exoProjRang1}
On note $A_1$, \ldots, $A_n$ les colonnes de $A$ et $t = \Tr A = \sum_i a_{ii}$.

\emph {1.} On suppose que $A$ est un \prr de rang $\leq 1$.
Vérifions d'abord $tA_j = A_j$:

\snic {
\left| \matrix {a_{ii} & a_{ij} \cr a_{ki} & a_{kj} \cr} \right| = 0
\hbox { et } A^2 = A, \;\hbox{impliquent}\;
ta_{kj} = \sum_i a_{ii} a_{kj} =  \sum_i a_{ki} a_{ij} =  a_{kj}.
}

%\sni
Donc $(1-t)A = 0$, puis $(1-t)t= 0$, i.e. $t$ \idmz. De plus, si $a A=0$,
\hbox{alors $a t = 0$}, i.e. $a = a(1-t)$.

\emph {2.}
Sur le localisé en $a_{ii}$, deux colonnes quelconques $A_j$, $A_k$ sont
multiples de~$A_i$ donc $A_j \wedge A_k = 0$. D'où globalement $A_j \wedge
A_k = 0$, et donc $\cD_2(A) = 0$. Par ailleurs, en utilisant $\left| \matrix
{a_{ik} & a_{ij}\cr a_{kk} & a_{kj}\cr}\right| = 0$, on obtient $\sum_k a_{ik}a_{kj}
= \sum_k a_{ij}a_{kk} = a_{ij} \Tr A = a_{ij}$, i.e. $A^2 = A$. 
%Sur $\FF_p$,
%la matrice $A = \I_{p+1}$ vérifie $A^2 = A$ et $\Tr A = 1$, mais 
%\hbox{pas $\cD_2(A) = 0$}.

\emph {3.}
Le \sys de droite est de cardinal $1 + n^2$, celui du milieu de cardinal
$1 + {n \choose 2}^2$. Pour obtenir celui de gauche, il faut compter les
mineurs sans coin sur la diagonale. Supposons $n \ge 3$, il y en a ${n \choose
2} {n-2 \choose 2}$, il en reste ${n \choose 2}^2 - {n \choose 2} {n-2 \choose
2} = (2n-3){n \choose 2}$ d'où le cardinal $1 + (2n-3){n \choose 2}$.
Pour $n = 3$, chaque \sys est de cardinal $10$.
Pour $n>3$,   $1 + n^2$ est strictement plus petit que les deux autres.
%%%%%%%%%%%%%%%%%%%%%%%%%%%%%%%%%%%%%%%%%%%%%%%%%%%%%%%%%%%%%%%%%%%%%%%%%%%

\exer{exoMatProjRang1CoeffReg} 
\emph {1.}
On a $\dmatrix {a_{i\ell} & a_{ij}\cr a_{k\ell} & a_{kj}} =
0$, i.e. $a_{kj} a_{i\ell} = a_{ij}a_{k\ell}$.\\
 C'est l'\egt $C_j \cdot
L_i = a_{ij} A$. Quant à $L_i \cdot C_j$, c'est le \coe en
position $(i,j)$ de $A^2 = A$, i.e. $a_{ij}$.

\emph {2.}
On a $L_i\cdot A = L_i$ donc $L_i\cdot(\In-A) = 0$.  Réciproquement, pour $u
\in \gA^n$ tel que $\scp {L_i}{u} = 0$, il faut montrer que $u = (\In-A)(u)$,
i.e. $Au = 0$, i.e.  $\scp {L_k}{u} = 0$. Mais $a_{ij}L_k = a_{kj}L_i$ et comme
$a_{ij}$ est \ndzz, c'est immédiat.

\emph {3.} L'\egt $a_{kj}a_{i\ell}=a_{ij}a_{k\ell}$ montre que 
$C\cdot L=a_{ij}A$. En outre, si $\gA$
est à \dve explicite, on peut  calculer  $A$ à partir de $L$ et $C$. 
\\
Si l'on se donne une ligne $L$ dont les \coes sont appelés
$a_{i\ell}$ ($\ell\in\lrbn$) et une colonne $C$ dont les \coes sont appelés 
$a_{kj}$ ($k\in\lrbn$), avec l'\elt commun $a_{ij}$ \ndzz,
les conditions sont les suivantes:
\begin{itemize}
\item chaque \coe de $C\cdot L$ doit être divisible par $a_{ij}$, d'où $A=\frac1{a_{ij}}C\cdot L$,
\item on doit avoir $\Tr(A)=a_{ij}$,
i.e., $L\cdot C =a_{ij}$.
\end{itemize}
Naturellement, ces conditions sont directement reliées à 
l'inver\-si\-bi\-lité de l'\id
engendré par les \coes de $L$.

\emph {4.}
Dans l'\egt matricielle $C \cdot L = (L\cdot C)\,A$ à démontrer, chaque membre
est bi\lin en $(L, C)$. Or, l'\egt est vraie si $\tra L$ est une colonne de $\tra A$ et
$C$ une colonne de $A$, donc elle reste vraie pour $\tra L \in \Im\tra A$ et $C \in \Im
A$.  Le reste est facile.

%%%%%%%%%%%%%%%%%%%%%%%%%%%%%%%%%%%%%%%%%%%%%%%%%%%%%%%%%%%%%%%%%%%%%%%%%%%
\exer{exoFittTfPtf}
$M$ est le quotient d'un \mptf $P$ qui a les mêmes \idfs que lui.
Si $P\oplus N=\Ae n$, $M\oplus N$ est un quotient de~$\Ae n$ avec les mêmes
\idfsz. Donc il n'y a pas de relation non nulle entre les \gtrs de $\Ae n$
dans le quotient $M\oplus N$. 
Donc 

\snic{M\oplus N=\Ae n \;\hbox{ et }\;P/M\simeq (P\oplus N)/(M\oplus N)=0.}

%%%%%%%%%%%%%%%%%%%%%%%%%%%%%%%%%%%%%%%%%%%%%%%%%%%%%%%%%%%%%%%%%%%%%%%%%%%
\exer{exoRelateursCourts}{
\emph{1.} 
Une \syzy $z=\sum z_k e_k$ est somme de \syzys courtes \ssi il existe des
\syzys $z_{ij} \in \gA e_i \oplus \gA e_j$ telles que $z = \sum_{i < j}
z_{ij}$. Ceci se relit comme suit:

\snic{ \exists\alpha_{ij}, \beta_{ij} \in
\gA,\; z_{ij} = \alpha_{ij} e_i + \beta_{ij} e_j, \;\scp {z_{ij}}
{x} = 0 \et z = \sum_{i < j} z_{ij}.}

%\sni
Cela équivaut  à $z_k = \sum_{k < j} \alpha_{kj} + \sum_{i < k}
\beta_{ik}$ ($k\in\lrbn$) et $\alpha_{ij} x_i + \beta_{ij} x_j = 0$
\hbox{(pour $i<j$)}. Il s'agit bien d'un
\sli en les \gui{inconnues} $\alpha_{ij}$,~$\beta_{ij}$.

\emph{2.} 
En raisonnant localement, on peut supposer que les $x_i$ sont
multiples de $x_1$, ce que l'on écrit $b_i x_1 + x_i = 0$. D'où
les \syzys $r_i = b_i e_1 + e_i$ pour $i\in\lrb{2..n}$.  
\\
Soit $z \in
x^\perp$. Posons $y = z - (z_2r_2 + \cdots + z_nr_n)$, on a $y_i = 0$
pour $i \ge 2$, et donc~$y$ est une \syzy (très) courte.
Ainsi,
$z = y + \sum_{i=2} z_i r_i$ est une somme de \syzys courtes.

\emph{3.} 
Soient $x$, $y \in \gA$. On cherche $s$, $t$
avec $s+t = 1$,  $sx \in \gA y$ et $ty \in \gA x$. On écrit
la \syzy $(-1,-1,1)$ entre $(x, y, x+y)$  
comme somme de \syzys courtes

\snic{
(-1,-1,1) = (0,a,a') + (b,0,b') + (c,c',0).
}

%\sni
En particulier, $a'+b'=1$,
 et l'on conclut.

\emph{4.} 
Par \dfn $\sum_i a_{ii} = 1$ et 
$\left| 
\matrix {a_{ij} & a_{ik} \cr
x_{j} & x_{k}} \right| = 0$. Ceci fournit
de nombreuses \syzys courtes $a_{ij} e_k - a_{ik} e_j$.  On retient les
$r_{ik} = a_{ii} e_k - a_{ik} e_i$, i.e. celles correspondant à un \gui{mineur
diagonal} $\left| \matrix {a_{ii} & a_{ik} \cr x_{i} & x_{k}}
\right|$. Pour $z \in \Ae n$, on pose

\snic{y = Az \et  z' = \sum_{i,k} z_k r_{ik} = \sum_{i,k} z_k (a_{ii} e_k - a_{ik} e_i).
}

Alors $z=z'+y$: en effet, le \coe
de $e_j$ dans $z'$ est

\snic{
\big(\sum_i a_{ii}\big) z_j - \sum_{k} a_{jk} z_k = z_j - (Az)_j.
}

Puisque $A^2-A$ annule $M$, $z-y\in x^{\perp}$,
donc $z \in x^\perp \Rightarrow y  \in x^\perp$.
Chaque $y_ie_i$ est un \syzy (très) courte puisque $y_ix_i = 0$.
Donc  $z = z' + y=z'+\sum y_i e_i$
est une somme de \syzys courtes.
}

%%%%%%%%%%%%%%%%%%%%%%%%%%%%%%%%%%%%%%%%%%%%%%%%%%%%%%%%%%%%%%%%%%%%%%%%%%%
\exer{exoPetitsRelateurs}
\emph{1.}
On écrit $z \in \Ae n$ sous la forme

\snic{z = \scp {x}{z}.y +
(z - \scp {x}{z}.y),}

%\sni
ce qui fournit la \dcn
$\Ae n = \gA y \oplus x^\perp$.

\emph{2.}
Pour $i \leq j$, définissons $z_{ij} \in \gA$ par $z_{ij} = z_i y_j - z_j
y_i$ et posons
$$ \ndsp
z' = \som_{i < j} z_{ij} (x_j e_i - x_i e_j) =
\som_{i \leq j} z_{ij} (x_j e_i - x_i e_j)
.$$
Pour $k$ fixé, le \coe de $e_k$ dans la somme de droite est
\[
\begin{array}{r}
%\dsp
\sum_{j \ge k} z_{kj} x_j - \sum_{i < k} z_{ik} x_i  \;  =  \;
%\dsp
\sum_{j \ge k} (z_ky_j - z_jy_k) x_j - \sum_{i < k} (z_iy_k - z_ky_i) x_i
 \\[2mm]
    =  \; %\dsp
  z_k \sum_{j = 1}^n y_jx_j - y_k \sum_{j = 1}^n z_jx_j\;=\;
   z_k - \scp {z}{x} y_k\,.
  \end{array}
\]
Ce qui signifie $z' = z - \scp {z}{x} y$ et prouve le résultat demandé.

\emph{3.} Si $\psi$ est la forme bi\lin alternée associée à la matrice,
l'\egt $\scp {x}{z} = 0$ signifie simplement que $\psi(x,x)=0$.

\emph{4.}  
On peut
écrire une matrice alternée $n \times n$ comme somme de $n(n-1) \over 2$
petites matrices alternées. 
Pour $n = 3$, voici la matrice alternée permettant de faire le lien avec la
question~\emph{2} ($y$ est fixé et c'est $z$ qui varie) :
$$
M_z = \cmatrix {
0                & z_1y_2 - z_2y_1  & z_1y_3-z_3y_1 \cr
-z_1y_2 + z_2y_1 & 0                & z_2y_3-z_3y_2 \cr
-z_1y_3 + z_3y_1 & -z_2y_3 + z_3y_2 & 0 }.
$$
La \dcn de $M_z$ en petites matrices alternées fournit les
$\pi_{ij}$. Il faut noter que
$z\mapsto M_z,\;\Ae n \to \Mn(\gA)$ est une \ali
 et que $\pi(z)= M_z x$.

%%%%%%%%%%%%%%%%%%%%%%%%%%%%%%%%%%%%%%%%%%%%%%%%%%%%%%%%%%%%%%%%%%%%%%%%%%%
\exer{exoProjImLibre} \emph{1.} On suit la méthode du cours. Elle conduit à poser:
$$
A = \cmatrix {0_r & -Y\cr X & \I_n - P\cr},
\quad 
A' = \cmatrix {0_r & Y\cr -X & \I_n - P\cr}.
$$
Ces matrices vérifient
$$
A\cmatrix {\I_r & 0\cr 0 & 0_n\cr} = \cmatrix {0_r & 0\cr 0 & P\cr} A,
\
AA' = \I_{n+r}.
$$
\emph{2.} Immédiat puisque l'on a les formules sous les yeux.

%%%%%%%%%%%%%%%%%%%%%%%%%%%%%%%%%%%%%%%%%%%%%%%%%%%%%%%%%%%%%%%%%%%%%%%%%%%
\exer{exoStabLibRang1}{
La formule de Binet-Cauchy donne $1=\det(RR')=\sum \delta_i \delta'_i$
avec 

\snic{\delta_i=\det(R_{1..n-1,1..n\setminus i})$
et $\delta'_i=\det(R'_{1..n\setminus i,1..n-1}).}

%\sni
On pose $S=[\,\delta'_1\;-\delta'_2\;\cdots\;(-1)^{n-1}\delta'_n\,]$.
On vérifie que la matrice carrée $A=\cmatrix{S\cr R}$ est de \deter $1$.
Ceci montre que $\Ker R$ est libre (proposition \ref{propStabliblib}).
\\
En fait soit $S'=\tra{[\,\delta_1\;-\delta_2\;\cdots\;(-1)^{n-1}\delta_n\,]}$ et $A'=[S'\,R']$. Alors $AA'=\In$, et ceci montre que $S'\in\Ae n$
est une base de $\Ker R$.
}

%%%%%%%%%%%%%%%%%%%%%%%%%%%%%%%%%%%%%%%%%%%%%%%%%%%%%%%%%%%%%%%%%%%%%%%%%%%

\exer{exoStablib}{
\emph{1.} On considère la matrice $B\,A$. Par \dfn de $B$, $B\,A$
est triangulaire supérieure, de diagonale $(\scp {x}{y}, \delta, \ldots,
\delta)$ où $\delta = \det(A)$.  En prenant le déterminant, 
on obtient $\det(B)
\det(A) = \scp {x}{y} \,\delta^{n-1}$. Les \idas annoncées sont donc vraies lorsque $\delta$ est \ivz. Puisqu'il s'agit d'\idasz, elles sont toujours vraies.
Le deuxième point de la question
est immédiat.

 \emph{2.} 
On écrit $z \in \Ae n$ sous la forme $z = \scp {y}{z}\,x +
(z - \scp {y}{z}\,x)$, ce qui fournit la décom\-position
$\Ae n = \gA x \oplus y^\perp$

 \emph{3.} 
L'hypothèse revient à dire que $x$ est la première colonne d'une
matrice inversible~$A$. Donc $y$ est la première ligne de la matrice
inversible $B$ ci-dessus. La matrice~$\tra {B}$ est adaptée à la
décomposition $\Ae n = \gA y \oplus x^\perp$.
}

%%%%%%%%%%%%%%%%%%%%%%%%%%%%%%%%%%%%%%%%%
\exer{exoMatLocSym}
Notons $x= [\,x_1 \;\cdots\;x_n \,]$. Pour $\alpha = \tra {[\alpha_1, \ldots, \alpha_n]}$ et $\beta = A\alpha$,
l'\egt en question est:

\snic {
\beta \,x=\cmatrix{\beta_1\cr \vdots\cr \beta_n}\cmatrix{x_1 \;\cdots \; x_n}=
{x}\,{\alpha}\,A\;\;$ avec $\;\;{x}\,{\beta}=x \,\alpha.
}

%\sni
Prenons $\alpha_i = x_i$. Puisque $x A = x$ et $A$ est \smqz,
on obtient $A \tra {x} = \tra {x}$, \hbox{i.e.  $\beta = \tra {x}$}. D'où
$\tra {x}\, x = {x}\tra{x}\,A = (x_1^2 + \cdots + x_n^2)A$. \\
Finalement, 
 $x_ix_j \in \gen {x_1^2 + \cdots + x_n^2}$ ($i,j\in\lrbn$).

%%%%%%%%%%%%%%%%%%%%%%%%%%%%%%%%%%%%%%%%%
\exer{exoPgcdPpcm}{
\emph{1.}
Soit $\alpha \in (\fb : \fa)$ et $\beta \in (\fa : \fb)$ vérifiant $1 =
\alpha + \beta$. 
Alors, la matrice $\theta = \crmatrix {\alpha & \beta\cr -1 & 1\cr}$,
de déterminant $1$ et
d'inverse $\theta^{-1} = \crmatrix {1 & -\beta\cr 1 & \alpha\cr}$, convient.  En
effet:

\snic{\crmatrix {\alpha & \beta\cr -1 & 1\cr} \cmatrix {\fa\cr \fb\cr} \subseteq
\cmatrix {\fa\cap \fb\cr \fa+\fb\cr}
\quad \hbox {et} \quad
\crmatrix {1 & -\beta\cr 1 & \alpha\cr} \cmatrix {\fa\cap \fb\cr \fa+\fb\cr}
\subseteq \cmatrix {\fa\cr \fb\cr}.}

%\sni
\`A gauche, l'inclusion haute vient du fait que $\alpha \fa + \beta \fb \subseteq
\fa \cap \fb$, celle du bas est triviale. 
\`A droite, l'inclusion haute vient
du fait que \hbox{$\fa \cap \fb + \beta(\fa + \fb) \subseteq \fa$}, 
 et celle du bas vient
du fait que $\fa \cap \fb + \alpha (\fa + \fb) \subseteq \fb$. Bilan: on a
\hbox{$\theta(\fa \oplus \fb) = (\fa \cap \fb) \oplus (\fa +
\fb)$} avec $\theta = \crmatrix {1 & \beta\cr
0 & 1\cr} \crmatrix {1 & 0\cr -1 & 1\cr} \in \EE_2(\gA)$.

 \emph{2.}
On peut prendre $A$ de la forme $A = \cmatrix {u & v\cr -b' & a'}$ avec $ua' +
vb' = 1$ et $a'b = b'a$. Posons $m = a'b = b'a$ et $d = ua
+ vb$. L'\egt $\cmatrix {a\cr b} =A^{-1} \cmatrix {d\cr 0}$ donne $a =
da'$ et~$b = db'$.  Il est clair que $\fa \cap \fb = \gen {m}$ et $\fa + \fb =
\gen {d}$. On a $a' \in (\fa : \fb)$ et $b' \in (\fb : \fa)$.  Donc $1 =
\alpha + \beta$ avec $\alpha = vb' \in (\fb : \fa)$, $\beta = ua' \in (\fa :
\fb)$.  Pour expliciter une \eqvc matricielle, il suffit d'utiliser une
matrice $\theta$ de la question précédente:

\snuc{\theta \cmatrix {a\cr 0} = \cmatrix {vm\cr -a} =
v\cmatrix {m\cr 0} - a'\cmatrix {0\cr d},
\;
\theta \cmatrix {0\cr b} = \cmatrix {um\cr b} =
u\cmatrix {m\cr 0} + b'\cmatrix {0\cr d}
.}

%\sni
D'où l'\eqvc matricielle: $\theta \cmatrix {a & 0\cr 0 &b\cr} =
\cmatrix {m & 0\cr 0 &d\cr}\cmatrix {v & u\cr -a' &b'\cr}$.

 \emph{3.}
L'hypothèse est $a = a^2x$ pour un certain $x$. Alors, l'\elt $e = ax$ est \idm et
$\gen{a}=\gen{e}$. On
doit résoudre $a'b = b'a$, $1 = ua' + vb'$, qui est un \sli
en $a', b', u, v$. \\
Modulo $1 - e$, on a $ax = 1$, on prend $a' = a$, $b' = b$,
$u = x$, $v = 0$. \\
Modulo $e$, on a $a = 0$, on  prend $a' = a$, $b' = 1$, $u
= 0$, $v = 1$. 
\\ Donc globalement:

\snic{a' = a, \quad b' = axb + (1-ax)1 = 1 - ax + axb, \quad
u = ax^2, \quad v = 1-ax.}

\sni \emph{4.}
Soit $\fa = \langle x_1, \ldots, x_n \rangle$ et $\fb = \langle y_1, \ldots,
y_m \rangle$. 
\\
On écrit $\fa + \fb = \langle z_1, \ldots, z_{n+m}\rangle$
avec $z_1 = x_1$, \dots, $z_{n+m} = y_m$.  Soient  $s_1, \ldots, s_{n+m}$
\com tels que sur $\gA_{s_i}$, on ait $\fa
+\fb= \gen{z_i}$.  \\
Dans chaque localisé on a  $\fa\subseteq \fb$ ou $\fb\subseteq \fa$,
d'où $\so{\fa+\fb, \fa\cap\fb} =\so{\fa,\fb}$, et:

\snic{1 \in (\fa : \fb) + (\fb : \fa)$,  $\,\fa\fb = (\fa\cap\fb)(\fa \oplus\fb)\,$ et  $\,\fa \cap \fb$ est \tfz.$ }

}

%%%%%%%%%%%%%%%%%%%%%%%%%%%%%%%%%%%%%%%%%%%%%%%%%%%%%%%%%%%%%%%%%%%%%%%%%%%
\exer{exoPrecisionsDet1} {
Nous reprenons les notations du lemme
\ref{lem calculs}.
\\
Voyons les \deters de $e\varphi$ et $\varphi_e$.
On a 
$$\preskip.3em \postskip.1em\mathrigid1mu
\det(\varphi)=\det(\I_n-F+H),\, \det(e\varphi)=\det(\I_n-F+eH)
\;\hbox{et}\;\det(\varphi_e)=\det(\I_n-eF+eH).
$$
 On en déduit 
$$\preskip.0em \postskip.4em 
\begin{array}{ccc} 
e\det(\varphi_e)=\det(e\I_n-eF+eH)=e\det(\varphi) \hbox{ et }
 \\[.3em] 
f\det(\varphi_e)=\det(f\I_n-feF+feH)= \det(f\I_n)=f. 
\end{array}
$$ 
Donc  $\det(\varphi_e)=f\det(\varphi_e)+ e\det(\varphi_e)
=f+e\det(\varphi)$.
\\
 De même $e\det(e\varphi)=\det(e\I_n-eF+eH)=e\det(\varphi)$  et
$$\preskip.4em \postskip.4em 
f\det(e\varphi)=\det(f\I_n-fF+feH)=f\det(\I_n-F)=f\rR{M}(0)=f\ide_0(M). 
$$
 En appliquant $\det(\varphi_e)=f+e\det(\varphi)$ aux \endos
$\Id+X\varphi$, $X\Id-\varphi$ et $X\Id$ du $\AX $-module $M[X]$ on obtient
$\rF{\varphi_e}(X)= f+e\rF{\varphi}(X)$, $\rC{\varphi_e}(X)=f+e\rC{\varphi}(X)$ et
$\rR{eM}(X)=f+e\rR{M}(X).$

 Par ailleurs,
la matrice $eH$ représente à la fois l'\endo $e\varphi$ de $M$ et l'\endo
$\varphi_e$ de $eM$. On a donc  $\rF{\varphi_e}(X)=\rF{e\varphi}(X)=\det(\I_n+eXH)
=
\rF{\varphi}(eX)$.

 En ce qui concerne la dernière affirmation: on doit regarder
$\det(\varphi_e)$ dans $\aqo{\gA}{f}$, on obtient $e\det(\varphi)$ modulo $f\gA$,
et cela correspond à l'\elt $e\det(\varphi)$ de $e\gA$.
}

%%%%%%%%%%%%%%%%%%%%%%%%%%%%%%%%%%%%%%%%%%%%%%%%%%%%%%%%%%%%%%%%%%%%%%%%%%%
\exer{exoPrecisionsDet2}
{\emph{1.} On a $\varphi\ep{h}=\varphi_{r_h}$ en appliquant la
notation de l'exercice \ref{exoPrecisionsDet1}. \\
Donc
$\delta_h=s_h+d_h$. On a $\delta_0=1$ parce que $M\ep{0}=\so{0}$, et puisque
$\delta_0=s_0+d_0$, cela donne $d_0=r_0$.
\\
 L'\egt $d=d_0+d_1+\cdots+d_n$ est triviale.
\\ 
 L'\egt $d=\delta_1 \smalltimes \cdots \smalltimes\, \delta_n$ résulte du point \emph{3}
du \thrf{propdef det ptf}.
On peut aussi démontrer $d_0+d_1+\cdots+d_n=\delta_1 \smalltimes \cdots \smalltimes\,
\delta_n$ par un calcul direct.

 \emph{2} et \emph{3.} 
Déjà vus dans l'exercice \ref{exoPeticalculPolrang}.
}

%%%%%%%%%%%%%%%%%%%%%%%%%%%%%%%%%%%%%%%%%%%%%%%%%%%%%%%%%%%%%%%%%%%%%%%%%%%
\exer{exoPrecisionsDet3}{
Rappel: pour $a\in \gA$ on a
$\det(a\varphi)=\rR{M}(a)\det(\varphi)=a^h\det(\varphi)$.
\\
 On se place alors sur l'anneau $\gA[X,1/X]$ et l'on considère le module
$M[X,1/X]$, on obtient
%---------begin $$--------
$$\preskip.4em \postskip.4em
X^h\rF{\varphi}(-1/X)=
\det \big(X(\Id_M-(1/X)\varphi)\big)=\det(X\Id_M-\varphi)= \rC{\varphi}(X)
.$$
%---------end $$----------
En remplaçant $X$ par $-1/X$ dans
$\rC{\varphi}(X)=X^h\rF{\varphi}(-1/X)$ on obtient l'autre \egtz. Les deux \pols
sont donc de degrés $\leq h$. Comme le \coe constant de $\rF{\varphi}$ est
égal à 1, on obtient aussi que $\rC{\varphi}$ est unitaire.
\\
 Pour les homogénéisés, le même calcul fonctionne.
\\
 Pour le \deter on remarque que $\det(-\varphi)=\rC{\varphi}(0)$.
}

%%%%%%%%%%%%%%%%%%%%%%%%%%%%%%%%%%%%%%%%%%%%%%%%%%%%%%%%%%%%%%%%%%%%%%%%%%%
\exer{exoPrecisionsDet4}{
On se place sur l'anneau $\gA_{r_h}$ et l'on considère le module $r_hM$ et l'\endo
$\varphi\ep{h}$. On obtient un module de rang constant $h$. Donc
$r_h\rF{\varphi}(X)$ \hbox{et $r_h\rC{\varphi}(X)$} sont de degrés $\leq h$, et 
$r_h \big(X^h\rF{\varphi}(-1/X)\big)=r_h\rC{\varphi}(X)$.
Il reste à faire la somme des \egts ainsi obtenues pour ${ h\in\lrbn}$.
\\
 Même calcul pour la deuxième \egtz.
 Les deux dernières \egts étaient déjà connues, sauf pour
$\det(\varphi)= r_0+r_1v_1+\cdots+r_nv_n$ qui peut se démontrer comme la
première.
}

%%%%%%%%%%%%%%%%%%%%%%%%%%%%%%%%%%%%%%%%%
%:  sol des pb

\prob{exoSuslinCompletableLemma}
\emph {1.}
Soient $C, U \in \Mn(\gA)$ telles que $AD = \In + bU$, $DA = \In + bC$.
Alors: 
$$\preskip.0em \postskip.4em 
\cmatrix {A &b\In\cr C &D\cr} \cmatrix {D &-b\In\cr -U &A\cr} =
\cmatrix {\In &0\cr * &\In\cr} \in \GL_{2n}(\gA). 
$$

\emph {2.}
On travaille modulo $a$ en remarquant que $b$ est \iv modulo $a$. On
peut donc, sur $\gA\sur{a\gA}$ considérer $b^{-1}B'$: c'est une matrice
diagonale de \deterz~1, donc elle appartient
à $\En(\gA\sur{a\gA})$ (cf exercice \ref {exoFacileGrpElem2}), on la remonte
en une \hbox{matrice $E \in \En(\gA)$} et l'on obtient $B' \equiv bE \bmod a$.

\emph {3.}
Immédiat.

\emph {4.}
Il suffit d'utiliser la sous-matrice $\I_{n-1}$ qui figure dans $B'$
pour tuer les coefficients des $n-1$ dernières colonnes de $D'$.
La sous-matrice carrée d'ordre $n+1$ obtenue
à partir de $\cmatrix {A &B'\cr C &D''\cr}$ en supprimant les
lignes $2$ à $n$ et les $n-1$ dernières colonnes est \iv
de première ligne $[\,a_1\;\cdots\;a_n\; b^n\,]$.

\emph {5.}
Modulo $z$, le vecteur $\vab x y$ est complétable en $A := \crmatrix {x &y\cr
-v & u\cr}$. \\
On a $\det(A)= a := ux + vy\equiv1 \mod z$ et l'on peut prendre $D=\wi
A$. \\
On écrit $DA = a\,\I_2 = \I_2 - wz\I_2$, donc $C = -w\I_2$. 
La matrice $\cmatrix {A &z\I_2\cr C &D\cr}$ est de \deter
$(a + wz)^2$. Pour trouver $E$, on utilise l'\egt 

\snic {
\cmatrix {z &0\cr 0&z^{-1}} = \rE_{21}(-1) \rE_{12}(1-z^{-1}) 
\rE_{21}(z) \rE_{12} \big(z^{-1}(z^{-1}-1)\big)
}

%\sni
et le fait que modulo $a$,  $zw\equiv 1$.
L'auteur de l'exercice a obtenu une matrice $G$ plus compliquée que celle de 
Krusemeyer. Avec $p = (y+u)w - u$, $q = (x-v)w + v$:

\snic {
G = \cmatrix {x&y&z^2\cr p(w-1)v-w & -p(w-1)u & y+u(z+1)\cr
            -q(w-1)v & q(w-1)u-w & -x+v(z+1)\cr}.
}

%\sni
On a $\det(G)=1 + (xu+yv+zw-1) (wz+1) (yq-xp+1)$ tandis
que la matrice de Krusemeyer est de \deter $(ux + vy + wz)^2$!

\emph {6.}
Immédiat par \recuz.

%%%%%%%%%%%%%%%%%%%%%%%%%%%%%%%%%%%%%%%%%%%%%%%%%%%%%%%%%%%%%%%%%%%%%%%%%%%

\prob{exoSphereCompletableYengui} 
\emph {1a.}
On a $\det(M) = -(cx_1-bx_2)u + (cx_0-ax_2)v$.\\
Avec 
$u = -(cx_1+bx_2)$, $v = cx_0+ax_2$, on obtient
$\det(M) = cx_0^2 + cx_1^2 - (a^2 + b^2)x_2^2$. Il suffit
de prendre $c = 1$ et $a$, $b \in \gA$ tels que $-1 = a^2 + b^2$.

\emph {1b.}
Montrons que $-1$ est une somme de deux carrés si $\gA$ contient un corps
fini. On peut supposer que $\gA$ est un corps de cardinal impair $q$. 
\\
On considère
les ensembles $A =
\sotQ {a^2}{a \in \gA}$ et $B = \sotQ {-1-b^2}{b \in \gA}$. 
\\
Ils ont 
$(q+1)/2$ \eltsz, donc $A \cap B \ne \emptyset$, ce qui donne le résultat.

Voici maintenant un résultat plus \gnlz:
 si $n \not\equiv 0 \mod 4$, alors $-1$ est une somme
de deux carrés dans $\ZZ/n\ZZ$. L'hypothèse peut
s'écrire $\pgcd(n, 4) = 1,2$ donc $2 \in n\ZZ + 4\ZZ$,
$2 = nu + 4v$. On pose $m = -1 + nu = -4v + 1$; puisque
$\pgcd(4n, m) = 1$, la progression arithmétique $4n\NN + m$
contient un nombre premier~$p$ (Dirichlet), qui vérifie
$p \equiv m \equiv -1 \mod n$ et
$p \equiv m \equiv 1 \mod 4$. \\
D'après cette dernière congruence,
$p$ est une somme de deux carrés, $p = a^2 + b^2$,
donc $-1 = a^2 + b^2$ dans $\ZZ/n\ZZ$.

On en déduit que si $n.1_\gA = 0$ avec $n \not\equiv 0 \mod 4$
(c'est le cas si $n$ est un nombre premier), 
alors $-1$ est une somme de deux carrés dans $\gA$.

\emph {2a.}
Soient $a_1$, \dots, $a_n$ tels que $-1 = \sum_{i=1}^n a_i^2$. On va utiliser:

\snic {
\sum_{i=1}^n (x_i-a_ix_0) (x_i + a_ix_0) = 
\sum_{i=1}^n x_i^2 -x_0^2 \sum_{i=1}^n a_i^2 = 1
}

%\sni
\penalty-2500
On a:
$$\preskip.4em \postskip.4em 
\cmatrix {x_0\cr x_1\cr \vdots\cr x_n\cr}  \sims{\EE_{n+1}}
\cmatrix {x_0\cr x_1 + a_1x_0\cr \vdots\cr x_n + a_nx_0\cr}  \sims{\EE_{n+1}}
\cmatrix {x_0 + h\cr x_1 + a_1x_0\cr \vdots\cr x_n + a_nx_0\cr}
\quad
\begin {array}{c}
\hbox {avec}\\
%% h = \lambda_1 (x_1+a_1x_0) + \cdots + \lambda_n (x_n+a_nx_0)\\
h = \sum_{i=1}^n \lambda_i (x_i+a_ix_0)\\
\end {array}
$$
En prenant $\lambda_i = (1-x_0)(x_i - a_ix_0)$, on obtient $h = 1-x_0$ donc
$x_0 + h = 1$. 
\\
Il est ensuite clair que 

\snic{\tra {[1, x_1+a_1x_0, \ldots,
x_n+a_nx_0]} \sims{\EE_{n+1}} \tra {[1, 0, \ldots, 0]}.}

%\sni
 De manière
explicite, en numérotant les $n+1$ lignes de $0$ à $n$ (au lieu
de $1$ à $n+1$) et en posant:

\snuc {
\begin {array} {c}
N = \prod_{i=1}^n \rE_{i,0} \big(-(x_i + a_ix_0)\big)\
\prod_{i=1}^n \rE_{0,i} \big((1-x_0)(x_i - a_ix_0)\big)\  \prod_{i=1}^n \rE_{i,0}(a_i) 
\\[1mm]
M = N^{-1} = \prod_{i=1}^n \rE_{i,0}(-a_i) \
\prod_{i=1}^n \rE_{0,i} \big((x_0-1)(x_i - a_ix_0)\big)\ \prod_{i=1}^n \rE_{i,0}(x_i + a_ix_0),
\end {array}
}

%\sni
on obtient une matrice $M \in \EE_{n+1}(\gA)$ de première colonne
$\tra {[\,x_0\; \cdots\;x_n\,]}$.

\emph {2b.}
On utilise $\gB = \aqo {\gA}{x_{n+1}, \ldots, x_m}$. 
Les morphismes $\EE_r(\gA) \twoheadrightarrow \EE_r(\gB)$ sont
surjectifs.  On obtient d'abord:

\snic {
\tra {[x_0, \ldots, x_n]} \sims{\EE_{n+1}(\gB)}  \tra {[1, 0, \ldots, 0]},
}

%\sni
donc des $x'_0, \ldots, x'_n \in \gA$ avec en particulier 
$x'_0 \equiv 1 \mod \gen {x_{n+1}, \ldots, x_m}$ tels que

\snic {
\tra {[x_0, \ldots, x_n, x_{n+1}, \ldots, x_m]} \sims{\EE_{m+1}(\gA)} 
\tra {[x'_0, \ldots, x'_n, x_{n+1}, \ldots, x_m].}
}

%\sni
On en déduit facilement:

\snuc {
\tra {[x'_0, \ldots, x'_n, x_{n+1}, \ldots, x_m]} \sims{\EE_{m+1}(\gA)} 
\tra {[1, \ldots, x'_n, x_{n+1}, \ldots, x_m]} \sims{\EE_{m+1}(\gA)} 
\tra {[1, 0, \ldots, 0].}
}

%\sni
\emph {3a.}
On a $\crmatrix {x_0 & -x_1\cr x_1 & x_0\cr} \sims{\EE_{2}(\gA)} B = 
\cmatrix {x_0 + ax_1& -x_1 + ax_0\cr x_1 & x_0\cr}$.  En utilisant le fait que
$1+a^2$ est nilpotent, on voit que $x_0 + ax_1$ est \iv car

\snic {
(x_0 + ax_1)(x_0 - ax_1) = x_0^2 + x_1^2 - (1+a^2)x_1^2 = 1 -(1+a^2)x_1^2.
}

%\sni
La matrice $B \in \SL_2(\gA)$ possède un \coe\iv donc elle est dans
$\EE_2(\gA)$.

\emph {3b.}
Raisonner d'abord modulo $\gen {x_2, \ldots, x_n}$, puis comme dans la
question \emph {2b}.

\emph {3c.}
On peut prendre $\gk = \ZZ/2^e\ZZ$ avec $e \ge 2$: $-1$ n'est pas un
carré dans $\gk$. \hbox{Et $-1$} n'est pas non plus un carré dans $\gA_n$ 
puisqu'il y a des morphismes $\gA_n \to \gk$, par exemple le
morphisme d'\evn en $x_0 = 1$, $x_i = 0$ pour $i \ge 1$, 
$y_j = 0$ \hbox{pour $j \ge 2$}.

%%%%%%%%%%%%%%%%%%%%%%%%%%%%%%%%%%%%%%%%%%%%%%%%%%%%%%%%%%%%%%%%%%%%%%%%%%%

}% fin des solutions d'exos

%:  ---- Section*{references}-----------
%\newpage	
\Biblio

Concernant le \thrf{prop Fitt ptf 2} et la \carn des \mptfs par leurs \idfs voir \cite{Nor} \tho 18 p.~122 et exercice 7 p.~49.
\perso{l'exercice 7 page 49 établit l'\eqvc entre conoyau \pro et
application \lnlz.}
Notons cependant que la preuve de Northcott n'est pas entièrement \covz,
puisqu'il fait appel à un \prca des \mptfsz.

Nous avons défini le \deter d'un \endo d'un \mptf comme dans
\cite[Goldman]{Gold}. La différence réside dans le fait que
nos \dems sont \covsz.

Une étude sur la faisabilité du \tho de structure locale des \mptfs se trouve dans \cite[Díaz-Toca\&Lombardi]{DiL09}.

La proposition \ref{fact.homom loc pf} concernant $\big(\Lin_\gA(M,N)\big)_S$
est un résultat crucial que l'on trouve par exemple dans  \cite{Nor},
exercice~9
p.~50, et dans \cite{Kun} (chapitre~IV, proposition~1.10).
Ce résultat sera \gne dans la proposition~\ref{propPlateHom}.

Le \pb \ref{exoSuslinCompletableLemma} est d\^u à \Sus \cite{Sus77b}.

\newpage \thispagestyle{CMcadreseul}
\incrementeexosetprob

%:        %%%%%%%%%%%%%%%%%%%%%%%%%%%%%%%%%%%%
%:        %%%%%%%%%%%%%%%%%%%%%%%%%%%%%%%%%%%%
%---- Chapitre  {chap AlgStricFi}------------
\chapter{Algèbres \stfes et \aGsz}
\label{chap AlgStricFi}\relax
\label{chap AlgEtales} 
%--------------------
\minitoc

\Intro

Ce chapitre est consacré  à une \gnn naturelle pour les anneaux 
commutatifs de la notion d'\alg finie sur un corps. 
En \comaz, pour obtenir les conclusions dans le cas des corps, il est 
souvent \ncr de supposer non
seulement que l'\alg est un \evc \tfz, mais plus \prmt que le corps est discret et que l'on connaît une
base de l'\evcz. C'est ce qui nous a amené à introduire la notion
d'\alg \stfe sur un \cdiz.

La \gnn pertinente de cette notion aux anneaux commutatifs
est donnée par les \algs qui sont des \mptfs
sur l'anneau de base. Nous les appelons donc des \algs \stfesz.

Les sections~\ref{secEtaleSurCD} et~\ref{sec2GaloisElr} qui ne concernent que les
\algs sur les \cdis peuvent être lues directement
après la section~\ref{secGaloisElr}. Même chose pour la section
\ref{secAGTG} si l'on prend à la base un \cdi (certaines \dems sont alors simplifiées).

La section \ref{sec1Apf} est une brève introduction aux \apfsz, en insistant
sur le cas des \algs entières. 

Le reste du chapitre est consacré aux \algs \stfes proprement dites.

Dans les sections \ref{secAlgSte} et \ref{secAlgSpb} sont  introduites
les notions voisines d'\alg \ste et d'\alg \spbz, 
qui généralisent la notion
d'\alg étale sur un \cdiz.

Dans la  section \ref{secAGTG} on donne un exposé \cof des bases
de la théorie des \aGs
pour les anneaux commutatifs.
Il s'agit en fait d'une théorie d'Artin-Galois, puisqu'elle reprend
l'approche qu'Artin avait développée pour le cas des corps en partant directement d'un groupe fini d'\autos d'un corps, le corps de base n'apparaissant
que comme un sous-produit des constructions qui s'ensuivent.

%%%%%%%%%%%%%
\section{Algèbres étales sur un \cdiz}
\label{secEtaleSurCD}

\Grandcadre{Dans les sections \ref{secEtaleSurCD} et \ref{sec2GaloisElr},~$\gK$ désigne un \cdi non trivial}

Rappelons qu'une \Klgz~$\gB$ est dite finie (resp. \stfez) si elle est \tf en tant que \Kev (resp. si~$\gB$ est  un \Kev de dimension finie).
Si~$\gB$ est une \Klg finie,
cela n'implique pas que l'on sache déterminer une
base de~$\gB$ comme \Kevz, ni même que
$\gB$ soit discrète. Si elle est \stfez, au contraire,
on connaît une base finie de~$\gB$ comme \Kevz.
Dans ce cas, pour un~$x\in\gB$,  la trace, la norme, le \polcar de (la
multiplication par) $x$, ainsi que le \polmin de $x$ sur~$\gK$ peuvent se
calculer par les méthodes standards de l'\alg \lin sur un \cdiz.
De même toute sous-\Klg finie de~$\gB$ est \stfe et
l'intersection de deux sous-\algs \stfes est \stfez.

%:     Definition{defi1Etale}
\begin{definition}\label{defi1Etale}
\index{etale@étale!algèbre --- sur un corps discret}
\index{algèbre!etale@étale sur un corps discret}
%\index{algèbre!separable@séparable}
\index{algebrique@algébrique séparable!algèbre --- sur un corps discret}
\index{algebrique@algébrique séparable!element@\elt --- sur un corps discret}
\index{algèbre!algebrique@algébrique séparable sur un corps discret}
%\index{separable@séparable!algèbre --- }
Soit~$\gL$ un \cdi et~$\gA$ une \Llgz.
\begin{enumerate}
\item L'\algz~$\gA$ est dite  \emph{étale (sur~$\gL$)} si
elle est \stfe  et si le \discri $\Disc_{\gA/\gL}$ est \ivz.
\item Un \elt de~$\gA$ est dit  \emph{\agsp  (sur~$\gL$)} s'il annule un \pol \splz.
\item L'\algz~$\gA$ est dite  \emph{\agsp  (sur~$\gL$)} si tout \elt de
$\gA$ est \agsp sur~$\gL$.
\end{enumerate}
\end{definition}

%En particulier, 
Lorsque~$f$ est un \polu de~$\gL[X]$, l'\alg 
quotient~$\aqo{\gL[X]}{f}$ est étale \ssiz$f$ est \spl (proposition \ref{propdiscTra}).

%: subsec{Théreme de structure}
\penalty-2500
\subsec{Théorèmes de structure des \algs étales}

\vspace{3pt}
La proposition~\ref{propdiscTra} donne le lemme qui suit.

%:     Lemma{lem1EtaleCD}
\begin{lemma}\label{lem1EtaleCD}
Soit~$\gA$  une \Klg \stfe et  $a\in\gA$. Si le \polcar $\rC{\gA/\gK}(a)(T)$ est \splz,
alors l'\alg est étale et~$\gA=\gK[a]$.
\end{lemma}

Dans le fait \ref{fact1Etale}, les points \emph{\ref{i1fact1Etale}}
et \emph{\ref{i2fact1Etale}} précisent certains points du lemme~\ref{lemZerRed}
et du fait~\ref{factZerRedConnexe} (concernant les anneaux \zeds réduits \gnlsz),
dans le cas d'une \Klg \stfe réduite.

Des résultats \gnls sur les extensions entières d'anneaux \zeds
sont donnés dans la section~\ref{sec1Apf} \paref{EntSurZdim} et suivantes. 

%:     Fact{fact1Etale}
\begin{fact}\label{fact1Etale} Soit~$\gB\supseteq\gK$ une \alg \stfez.
\begin{enumerate}
\item \label{i1fact1Etale} L'\algz~$\gB$ est \zedez. Si elle est réduite, pour tout $a\in\gB$ il existe un unique \idm $e\in\gK[a]$ tel que~$\gen{a}=\gen{e}$. En outre, lorsque $e=1$, \cad lorsque $a$ est \ivz, $a^{-1}\in\gK[a]$.
\item \label{i2fact1Etale} \Propeq
\begin{enumerate}
\item $\gB$ est un \cdiz.
\item $\gB$ est sans diviseur de zéro: $xy=0\Rightarrow (x=0$ \emph{ou} $y=0)$.
\item $\gB$ est connexe et réduite.
\item Le \polmin sur~$\gK$ de n'importe quel \elt de~$\gB$ est
\irdz.
\end{enumerate}
\item \label{i3fact1Etale} Si~$\gK\subseteq\gL\subseteq\gB$ et~$\gL$ est 
un \cdi \stf sur~$\gK$, alors~$\gB$ est \stfe sur~$\gL$. En outre,~$\gB$ 
est étale sur~$\gK$ \ssi elle est étale
sur~$\gL$ et~$\gL$ est étale sur~$\gK$.
\item \label{i4fact1Etale} Si $(e_1,\ldots,e_r)$ est un \sfio de~$\gB$, 
$\gB$ est étale sur
$\gK$ \ssi chacune des compo\-santes~$\gB [1/e_i]$
est étale sur~$\gK$.
\item \label{i5fact1Etale} Si~$\gB$ est étale elle est réduite.
\item \label{i6fact1Etale} Si $\car(\gK) >\dex{\gB:\gK}$
et si~$\gB$ est réduite, elle est étale.
\end{enumerate}
\end{fact}
% -----  end{fact} -------------------------------------------------------
\begin{proof}
\emph{\ref{i1fact1Etale}.}
L'\elt $a$ de~$\gB$ est annulé par un \polu de $\KT$ que l'on 
écrit $uT^k\big(1 - T\,h(T)\big)$
avec $u\in\gK\eti$, $k\geq0$.  Donc~$\gB$ est \zedez. Si elle est réduite, $a\big(1 - ah(a)\big) = 0$.
Alors, $e=ah(a)$ vérifie $a(1-e)=0$ et a fortiori $e(1-e) = 0$. 
Ce qui permet de conclure.

\emph{\ref{i2fact1Etale}.} L'\eqvc de \emph{a}, \emph{b} et \emph{c} est un cas particulier du lemme~\ref{lemEntReduitConnexe}. %\\
L'implication \emph{d} $\Rightarrow$ \emph{c} est claire.
%\\
Voyons \emph{b} $\Rightarrow$ \emph{d.}
Soit $x$ dans~$\gB$ et $f(X)$ son \polmin sur~$\gK$. Si $f=gh$, avec
$g$, $h$ \monsz, alors
$g(x)h(x)=0$ donc $g(x)=0$ ou $h(x)=0$.
Par exemple $g(x)=0$, et puisque~$f$ est le \polminz,~$f$ divise $g$, et $h=1$.
 
\emph{\ref{i3fact1Etale}.} Soit $(f_1,\ldots,f_s)$
une~$\gK$-base de~$\gL$.
On peut calculer une~$\gL$-base de~$\gB$
 comme suit. La base commence avec $e_1=1$. Supposons avoir calculé
des \elts $e_1$, \ldots, $e_r$ de~$\gB$ \lint indépendants sur~$\gL$.
Les %sous-\Kevs 
$\gL e_i$ sont en somme directe dans~$\gB$ et l'on  a une~$\gK$-base~$(e_if_1,\ldots,e_if_s)$ pour chaque~$\gL e_i$. Si $rs=\dex{\gB:\gK}$, on a terminé.  
Dans le cas contraire, on
peut trouver~$e_{r+1}\in\gB$ qui n'est pas dans %le sous-\Kev 
$F_r=\gL e_1\oplus \cdots\oplus \gL e_r$. \\
Alors,~$\gL e_{r+1}\cap F_r=\so{0}$ (sinon, on exprimerait $e_{r+1}$
comme~$\gL$-\coli  de $(e_1,\ldots,e_r)$).
Et l'on itère le processus en remplaçant $(e_1, \ldots, e_r)$
par~$(e_1, \ldots, e_{r+1})$.
\\
Une fois que l'on dispose d'une base de~$\gB$ comme \Levz, il reste à utiliser la formule de transitivité des \discris (\thrf{thTransDisc}).
 
\emph{\ref{i4fact1Etale}.} On utilise
le \tho de structure \rref{fact.sfio} pour les \sfios
et la formule du discriminant d'une \alg produit direct d'\algs
(proposition~\ref{propTransDisc}).
 
\emph{\ref{i5fact1Etale}.} Soit $b$ un \elt nilpotent de~$\gB$. Pour tout $x\in\gB$ la multiplication par~$bx$ est un \endo nilpotent $\mu_{bx}$ de~$\gB$. On peut alors trouver une~\hbox{$\gK$-base}
de~$\gB$ dans laquelle la matrice de $\mu_{bx}$ est strictement triangulaire, donc
$\Tr(\mu_{bx})=\Tr_{\gB/\!\gK}(bx)=0$. \\
Ainsi $b$ est dans le noyau de l'\Kli 

\snic{tr:\gB\to\Lin_\gK(\gB,\gK),\;\;b\mapsto(x\mapsto \Tr_{\gB/\!\gK}(bx).}

%\sni
Enfin, $tr$ est un \iso  puisque $\Disc_{\gB/\!\gK}$ est \ivz, donc~$b=0$.
 
\emph{\ref{i6fact1Etale}.} Avec la notation précédente,
on suppose~$\gB$ réduite et l'on veut montrer que l'\Kli $tr$ est un \isoz.
\\
Il suffit de montrer que $\Ker tr=0$. Supposons $tr(b)=0$, alors  $\Tr_{\gB/\!\gK}(bx)=0$ pour tout $x$ et en particulier $\Tr_{\gB/\!\gK}(b^n)=0$
pour tout $n>0$. Donc l'\endo $\mu_b$ de multiplication par $b$ vérifie
$\Tr(\mu_b^n)=0$ pour tout $n>0$. Les formules qui relient les \isN sommes de Newton
aux fonctions \smqs \elrs
montrent alors que le \polcar de $\mu_b$ est égal
à $T^{\dex{\gB:\gK}}$  (cf. exercice~\ref{exoSommesNewton}).
Le \tho de Cayley-Hamilton et le fait que~$\gB$ est réduite permettent de conclure que $b=0$.
\end{proof}
%

%:     Theorem{th1Etale}
\begin{theorem} \emph{(\Tho de structure des \Klgs étales, 1)}
\label{th1Etale}
\\
Soit %$\gK$ un \cdi non trivial et
$\gB$ une \Klg étale.
\begin{enumerate}
\item \label{i1th1Etale} Tout \idz~$\gen{b_1,\ldots,b_r}_\gB$  est engendré par un \idm $e$ qui appartient à~$\gen{b_1,\ldots,b_r}_{\gK[b_1,\ldots,b_r]}$.
Et l'\alg quotient est étale sur~$\gK$.
\item \label{i2th1Etale} Soit~$\gA$ une sous-\Klg \tf de~$\gB$.
\begin{enumerate}
\item $\gA$ est une \Klg étale.
\item Il existe un entier $r\geq1$ et un \sfio $(e_1,\ldots,e_r)$ de~$\gA$
tel que, pour chaque $i\in \lrbr$,~$\gB[1/e_i]$ est un module libre de rang fini
sur~$\gA[1/e_i]$. En d'autres termes,~$\gB$ est un module quasi libre sur~$\gA$.
\end{enumerate}
\item \label{i3th1Etale} $\gB$ est  \agsp sur~$\gK$.
\item \label{i4th1Etale}
Pour tout $b\in\gB$, le \polcar $\rC{\gB/\gK}(b)$ est un produit de \pols \splsz.
\end{enumerate}
\end{theorem}
\begin{proof}
\emph{\ref{i1th1Etale}.} Si l'\id est principal cela résulte du fait \ref{fact1Etale} point~\emph{\ref{i1fact1Etale}.}
Par ailleurs, pour deux \idms $e_1$, $e_2$, on~a~$\gen{e_1,e_2}=\gen{e_1+e_2-e_1e_2}$.
Enfin l'\alg quotient est elle même étale sur~$\gK$ d'après la formule du discriminant d'une \alg produit direct.

\emph{\ref{i2th1Etale}.} Il suffit de démontrer le point \emph{b},
car alors on conclut en utilisant la formule
de transitivité des \discris pour chaque
$\gK\subseteq\gA[1/e_i]\subseteq\gB[1/e_i]$ et la formule du \discri d'une \alg produit direct.
\\
Pour démontrer le point \emph{b}, on essaie de calculer une base de~$\gB$ sur
$\gA$ en utilisant la méthode indiquée dans le cas où~$\gA$ est un \cdi
dont on connaît une~$\gK$-base, donnée dans le fait~\ref{fact1Etale}~\emph{\ref{i3fact1Etale}}.
Le point où l'\algo risque d'achopper est lorsque $e_{r+1}\gA\cap F_r$ n'est pas réduit à $\so{0}$. On a alors une \egt $\alpha_{r+1}e_{r+1}=\sum_{i=1}^r\alpha_ie_i$
avec tous les $\alpha_i$ dans~$\gA$, et $\alpha_{r+1}\neq0$ mais non \iv dans~$\gA$.
Ceci implique (point \emph{\ref{i1th1Etale}}) que l'on trouve un \idm $e\neq0,1$
dans~$\gK[\alpha_{r+1}]\subseteq\gA$.
On recommence alors avec les deux \lons en $e$ et $1-e$. Enfin, on remarque que le nombre de scindages ainsi opérés est a priori borné par~$\dex{\gB:\gK}$.

\emph{\ref{i3th1Etale}} et \emph{\ref{i4th1Etale}.} Résultent facilement de \emph{\ref{i2th1Etale}.}
\end{proof}

\rem Une \gnn du point \emph{\ref{i1th1Etale}}
du \tho précédent se trouve dans les lemmes~\ref{lemKlgEntiere} et~\ref{lemZrZr1}.
\eoe

\medskip
On peut construire des \Klgs étales de proche en proche en vertu du lemme suivant, qui prolonge le lemme~\ref{lem1EtaleCD}.

%:     Lemma{lemEtaleEtage}
\begin{lemma}\label{lemEtaleEtage}
Soit~$\gA$ une \Klg étale et $f\in\AT$ un \polu \splz.
Alors, $\aqo{\AT}{f}$ est une \Klg étale.
\end{lemma}
\begin{proof}
On regarde d'abord $\aqo{\AT}{f}$ comme une \Alg libre de rang $\deg f$.
On a $\Disc\iBA =\disc(f)$ (proposition~\ref{propdiscTra} point \emph{3}). On conclut par la formule de transitivité des \discrisz.
\end{proof}

Les deux \thos qui suivent sont des corolaires.

%:     theorem{cor1lemEtaleEtage}
\begin{theorem}\label{cor1lemEtaleEtage}
Soit~$\gB$ une \Klgz. Les \elts de~$\gB$ \agsps sur~$\gK$ forment une
sous-\algz~$\gA$. En outre, tout \elt de~$\gB$ qui annule
un \polu \spl de~$\gA[T]$ est dans~$\gA$.
\end{theorem}
\begin{proof}
Montrons d'abord que si $x$ est \agsp sur~$\gK$ et $y$ annule un \polu \spl $g$ de~$\gK[x][Y]$, alors tout \elt de~$\gK[x,y]$ est \agsp sur~$\gK$.
Si $f\in\KX$ \spl annule $x$, alors la sous-\algz~$\gK[x,y]$ est un quotient
$\aqo{\gK[X,Y]}{f(X),g(X,Y)}$. Cette \Klg est étale d'après le
lemme~\ref{lemEtaleEtage}.\\
En raisonnant par \recuz, on peut itérer la construction précédente.
On obtient le résultat souhaité en notant qu'une \Klg étale est
 \agsp sur~$\gK$, et que tout quotient d'une telle \alg est encore \agsp
 sur~$\gK$.
\end{proof}

Voici une variante \gui{\stfez}. Nous redonnons la \dem
car les variations, bien que simples, sont significatives
des précautions à prendre dans le cas \stfz.
%:     Theorem{corlemEtaleEtage}
\begin{theorem} \emph{(\Carn des \Klgs étales)}\label{corlemEtaleEtage}
\\
Soit~$\gB$ une \Klg \stfe donnée sous la forme $\Kxn$.
\Propeq
\begin{enumerate}
\item $\gB$ est étale sur~$\gK$.
\item  Le \polmin sur~$\gK$ de chacun des
$x_i$ est  \splz.
\item $\gB$ est \agsp sur~$\gK$.
\end{enumerate}
En particulier, un corps~$\gL$ qui est une extension galoisienne de~$\gK$
est étale sur~$\gK$.
\end{theorem}
\begin{proof}
\emph{1} $\Rightarrow$ \emph{3.}
 D'après le \thrf{th1Etale}.

\emph{2} $\Rightarrow$ \emph{1.}
Traitons d'abord le cas d'une \Klg \stfez~$\gA[x]$ où~$\gA$ est étale sur~$\gK$ et où  le \polminz~$f$ de $x$ sur~$\gK$
est \splz.
On a alors un \homo surjectif de la \Klg \stfez~$\aqo{\gA[T]}{f}$
sur~$\gA[x]$ et le noyau de cet \homo (qui se calcule comme noyau
d'une \ali entre \Kevs de dimensions finies) est \tfz, donc engendré par un
\idm $e$. La \Klgz~$\gC=\aqo{\gA[T]}{f}$ est étale d'après le lemme~\ref{lemEtaleEtage}. On en déduit  que~$\gA[x]\simeq\aqo{\gC}{e}$ est étale sur~$\gK$.
\\
On peut alors terminer par \recu sur $n$.
\end{proof}
%

%:     Corollary{corcorlemEtaleEtage}
\begin{corollary}\label{corcorlemEtaleEtage}
Soit $f\in\KT$ un \poluz. L'\adu
$\Adu_{\gK,f}$ est étale \ssiz$f$ est \splz.
\end{corollary}

\rem On avait déjà ce résultat par calcul direct du \discri de l'\adu (fait~\ref{factDiscriAdu}).\eoe

%:     Theorem{thEtalePrimitif}
\begin{theorem} \emph{(\Tho de l'\elt primitif)}
\label{thEtalePrimitif}
\\
Soit~$\gB$ une \Klg étale.
\begin{enumerate}
\item  Si~$\gK$ est infini ou si~$\gB$ est un \cdiz,~$\gB$ est une \alg monogène, \prmt de la forme~$\gK[b]\simeq\aqo{\KT}{f}$ pour un $b\in \gB$ et un $f\in\KT$ \splz. 
\\ 
Ceci s'applique en particulier pour un corps~$\gL$ qui est une extension galoisienne de~$\gK$, de sorte que l'extension~$\gL/\gK$ relève
du cas \elr étudié dans le \thref{thGaloiselr}.
\item $\gB$ est un produit fini de \Klgs étales monogènes.
\end{enumerate}
\end{theorem}
\begin{proof}
\emph{1.} Il suffit de traiter le cas d'une
\alg à deux \gtrsz~$\gB=\gK[x,z]$. On va chercher un \gtr de~$\gB$
de la forme $\alpha x+\beta z$ avec $\alpha$, $\beta\in\gK$.
On note~$f$ et $g$ les \polmins de $x$ et $z$ sur~$\gK$.
On sait qu'ils sont \splsz.
On note~$\gC=\aqo{\gK[X,Z]}{f(X),g(Z)}=\gK[\xi,\zeta]$.
Il suffit de trouver $\alpha$, $\beta\in\gK$ tels que~$\gC=\gK[\alpha\xi+\beta\zeta]$.
Pour avoir ce résultat, il suffit que le \polcar de $\alpha\xi+\beta\zeta$ soit \splz, car on peut alors appliquer le lemme~\ref{lem1EtaleCD}.
On introduit deux \idtrs $a$ et~$b$, et l'on note~$h_{a,b}(T)$ le \polcar de la
multiplication par~$a\xi+b\zeta$ dans~$\gC[a,b]$ vue comme~$\gK[a,b]$-\alg
libre de rang fini.
En fait: 
$$\preskip.4em \postskip.3em 
\gC[a,b]\simeq \aqo{\gK[a,b][X,Z]}{f(X),g(Z)}. 
$$
On note $d(a,b)=\disc_T(h_{a,b})$.
On fait
un calcul dans une \gui{double \aduz} sur~$\gC[a,b]$, dans laquelle
on factorise séparément~$f$ et $g$: 

\snic{f(X)=\prod_{i\in\lrbn}(X-x_i) \;\hbox{ et }\;g(Z)=\prod_{i\in\lrbk}(Z-z_j).}

%\sni
 On obtient

\snac{
\pm d(a,b)=\!\!\!\!\!\!\prod\limits_{(i,j)\neq(k,\ell)}\!\!\!\!(a(x_i-x_k)+b\big(z_j-z_\ell)\big)=
(a^{n^2-n}\! \disc f)^{p^2}  (b^{p^2-p} \!\disc g)^n + \ldots
}

%\sni

Dans le membre le plus à droite des \egts ci-dessus on a indiqué le terme de plus haut degré lorsque l'on ordonne les \moms en $a$, $b$ selon un ordre lexicographique. Ainsi le \pol $d(a,b)$ a au moins
un \coe \ivz.
Il suffit de choisir  $\alpha$, $\beta$ de façon que
$d(\alpha,\beta)\in\gK\eti$ pour obtenir un \elt $\alpha\xi+\beta\zeta$
de~$\gC$ dont le \polcar est \splz.
Ceci achève la \dem pour le cas où~$\gK$ est infini.
\\
Dans le cas où~$\gB$ est un \cdi on énumère les entiers de~$\gK$
jusqu'à obtenir $\alpha$, $\beta$ dans~$\gK$ avec 
$d(\alpha,\beta)\in\gK\eti$, ou à conclure que la \cara est égale 
à un nombre premier $p$. On énumère ensuite
les puissances des \coes de~$f$ et de $g$ jusqu'à obtenir 
suffisamment d'\elts
dans~$\gK$, ou à conclure que le corps~$\gK_0$ engendré 
par les \coes de~$f$ et $g$ est un corps fini. Dans ce cas,~$\gK_0[x,z]$
est lui même un corps fini et il est engendré par un \gtrz~$\gamma$
de son groupe multiplicatif, donc~$\gK[x,z]=\gK[\gamma]$.

\emph{2.} On reprend la preuve qui vient d'être donnée pour le cas où
$\gB$ est un \cdiz. Si l'on n'arrive pas à la conclusion, c'est que
la preuve a achoppé à un endroit précis, qui manifeste que~$\gB$ n'est pas un \cdiz. Puisque l'on  est avec une \Klg \stfez, cela nous fournit\footnote{Pour plus de précisions voir la solution de l'exercice~\ref{exothEtalePrimitif}}
un \idm $e\neq0,1$ dans~$\gB$.
%Cela se produit sous la forme que le corps~$\gK_0$ décrit
%ci-dessus est fini, mais que~$\gK_0[x,z]$ n'est pas un corps fini.
%On peut donc calculer un \idm $\neq0,1$ dans~$\gK_0[x,z]$, de sorte
%que
Ainsi~$\gB\simeq\gB[1/e]\times \gB[1/(1-e)]$. On peut alors conclure par
\recu sur $\dex{\gB:\gK}$.
\end{proof}
%

%%%%%%%%%%%%%%%%%%%%%%%%%%%%%%%%%%%%%%%%%%%%%%%%%%%%%%%%%%%%%%%%%%%%%%%%%%%
%: subsec{Algèbres étales sur un corps \splz ment factoriel}
%\penalty-5000
\subsec{Algèbres étales sur un corps \splz ment factoriel}

\vspace{3pt}
Lorsque tout  \pol \spl sur~$\gK$ se décompose en un
produit de facteurs irréductibles, le corps~$\gK$ est dit
\emph{sépara\-blement factoriel}.
\index{separablement@séparablement factoriel!corps discret ---}
\index{corps!separablement@séparablement factoriel}

%:     Lemma{lemsepfactcorps}
\begin{lemma}\label{lemsepfactcorps}
Un corps~$\gK$ est  \splz ment factoriel  \ssi on a un test pour l'existence d'un zéro
dans~$\gK$ pour un \pol \spl  arbitraire de $\KT$.
\end{lemma}
%
%\begin{proof}
%La preuve est à peu près la même que pour le lemme~\ref{lemKXfactor}.
%Ici, lorsque l'on tente une \fcn $f=gh$ avec $g$ et $h$
%\mons de degrés fixés, on se situe dans $\Adu_{\gK,f}$.
%Alors, chaque \coe de $g$ ou $h$ annule un produit de \pols \spls
%de $\KT$, d'après le \thrf{th1Etale} (point \emph{\ref{i4th1Etale}.})
%et le corolaire~\ref{corcorlemEtaleEtage}.
%\end{proof}
%
\begin{proof}
La deuxième condition est a priori plus faible puisqu'elle revient à
déterminer les facteurs de degré 1 pour un \pol \spl  de $\KX$. Supposons cette condition vérifiée. La preuve est à peu près la même que pour le lemme~\ref{lemKXfactor},
mais demande quelques détails \sulsz. On note $f(T)=T^{n}+\sum_{j=0}^{n-1}a_jY^{j}$, on fixe un entier $k\in\lrb{2..n-2}$ et l'on cherche les \pols {$g=T^{k}+\sum_{j=0}^{k-1}b_jT^{j}$}
qui divisent $f$. On va montrer qu'il n'y a qu'un nombre fini de possibilités, explicites, pour chacun des $b_j$.
La \dem du \tho de \KRO utilise des
\pols \uvls  $Q_{n,k,r}(a_0,\dots,a_{n-1},X)\in\ZZ[\ua,X]$, \unts en $X$, tels que $Q_{n,k,r}(\ua,b_r)=0$. 
\\
Ces \pols peuvent être calculés dans l'\adu $\gA=\Adu_{\gK,f}$ comme suit. On pose 
$$\preskip.4em \postskip.4em \ndsp
G(T)=\prod_{i=1}^{k}(T-x_i)=T^{k}+\sum_{j=0}^{k-1}g_jT^{j}. 
$$
On considère l'orbite $(g_{r,1},\dots,g_{r,\ell})$ de $g_r$ sous l'action de $\Sn$, et l'on obtient  
$$\preskip.4em \postskip.0em \ndsp
Q_{n,k,r}(\ua,X)=\prod_{i=1}^{\ell}(X-g_{r,i}). 
$$ 
On en déduit que 
$$\preskip-.2em \postskip.3em \ndsp
\prod\nolimits_{\sigma\in\Sn}\big(W-\sigma(g_r)\big) =Q_{n,k,r}^{n!/\ell}. 
$$ 
Donc, d'après le lemme \ref{lemPolCarAdu}, 
$\rC{\gA/\gk}(z)(X)=Q_{n,k,r}^{n!/\ell}(X)$.   
Enfin, comme~$\gA$ est étale sur $\gK$ (corolaire~\ref{corcorlemEtaleEtage}), le \polcar de $g_r$  annule un produit de \pols \spls de $\KT$ d'après le  \thrf{th1Etale}~\emph{\ref{i4th1Etale}}.
\\
Ainsi, $b_r$ doit être cherché parmi les zéros d'un nombre fini de \pols \splsz: il y a un nombre fini de possibilités, toutes explicites.
\end{proof}

%:     Theorem{th2Etale}
\begin{theorem} \emph{(\Tho de structure des \Klgs étales, 2)}
\label{th2Etale}\\
Supposons~$\gK$  \splz ment factoriel.
Une \Klgz~$\gB$ est étale \ssi elle est isomorphe à un produit fini de corps étales sur~$\gK$.
\end{theorem}
\begin{proof}
Conséquence du \tho de l'\elt primitif
(\thref{thEtalePrimitif}).
\end{proof}
%

%:     Corollary{corth1Etale}
\begin{corollary}\label{corth1Etale}
Si~$\gL$ est un corps étale sur~$\gK$ et si~$\gK$ est \splz ment factoriel, il en va de même pour~$\gL$.
\end{corollary}
\begin{proof}
Soit $f\in\gL[T]$ un \polu \splz. La \Llgz~$\gB=\aqo{\gL[T]}{f}$ est étale, donc c'est aussi une \Klg étale. On peut donc trouver un \sfio
tel que chaque composante correspondante de~$\gB$ est connexe. Cela revient à factoriser~$f$ en produit de facteurs irréductibles.
\end{proof}
%

%:    corollary{propIdemMini}----
\begin{corollary}
\label{propIdemMini}
\Propeq
\begin{enumerate}
\item  Toute \Klg étale est isomorphe à un produit de corps étales
sur~$\gK$.
\item  Le corps~$\gK$ est séparablement factoriel.
\item  Tout \pol \spl possède un corps de racines qui est
une extension \stfe (donc galoisienne) de~$\gK$.
\item  Tout \pol \spl possède un corps de racines qui est
étale sur~$\gK$.
\end{enumerate}
\end{corollary}
%--- end-corollary----------------------------------------
%
\begin{proof}
Pour \emph{2} $\Rightarrow$ \emph{4}, on utilise le fait que l'\adu
pour un \pol \spl est étale (corolaire~\ref{corcorlemEtaleEtage}) et l'on
applique le \thrf{th2Etale}.
\end{proof}
%

%:     Corollary{corth3Etale}
\begin{corollary}\label{corth3Etale}
Si~$\gK$ est \splz ment factoriel et si $(\gK_i)$ est une famille finie
de corps étales sur~$\gK$, il existe une extension galoisienne~$\gL$
de~$\gK$ qui contient une copie de chacun des~$\gK_i$.
\end{corollary}
%--------- fin corollary ---------------------------------------------- 
%
\begin{proof}
Chaque~$\gK_i$ est isomorphe à un $\aqo\KT {f_i}$
avec $f_i$ \ird \splz. On considère le
ppcm~$f$ des $f_i$ puis un \cdr de~$f$.
\end{proof}
%

%: subsec{Corps parfaits}
\subsec{Corps parfaits, clôture \spl et clôture \agqz}

\vspace{3pt}
Pour un corps~$\gK$ de \cara finie $p$ l'application $x\mapsto x^p$ est un \homo injectif.

En \clama un corps~$\gK$ est dit \ixc{parfait}{corps ---} s'il est de \cara infinie, ou 
si, étant de \cara finie $p$,  le morphisme~$x\mapsto x^p$
%$\gK\to\gK,\;x\mapsto x^p$ 
est un \isoz.

En \coma pour éviter la disjonction sur la \cara dans le \gui{ou} ci-dessus (qui peut ne pas être explicite), on formule la chose comme suit: 
\emph{si $p$ est un nombre premier tel que $p.1_\gK=0_\gK$, alors l'\homoz~$\gK\to\gK,\;x\mapsto x^p$ est surjectif.}
 
Le corps des rationnels $\QQ$ et les corps finis  (dont le corps trivial) sont parfaits.

Soit~$\gK$ un corps de \cara finie $p$.
Un surcorps~$\gL\supseteq\gK$ est appelé une \emph{clôture parfaite} de~$\gK$
si c'est un corps parfait et si tout \elt de~$\gL$, élevé à une certaine puissance $p^k$ est un \elt de~$\gK$.

\index{cloture@clôture!parfaite}
 
%:     Lemma{lemClotParf}
\begin{lemma}\label{lemClotParf}
Un \cdiz~$\gK$ de \cara finie $p$ possède une clôture parfaite~$\gL$,
unique à \iso unique près. \\
En outre,~$\gK$ est une partie détachable de~$\gL$
\ssi il existe un test pour \gui{$\exists x\in\gK,\; y=x^p$?} (avec extraction de la racine $p$-ième de $y$ quand elle existe). 
\end{lemma}
%--------- fin lemma ----------------------------------------------
%
\begin{Proof}{Idée de la \demz. } Un \elt de~$\gL$
 est codé par un
couple $(x,k)$, où $k\in\NN$ et $x\in\gK$. Ce code représente la racine 
$p^k$-ième de $x$. 
\\
L'\egt dans~$\gL$, $(x,k)=_\gL(y,\ell)$, est définie par $x^{p^{\ell}}=y^{p^{k}}$ (dans~$\gK$), de sorte \hbox{que $(x^p,k+1)=_\gL(x,k)$}. 
\end{Proof}
%

%:     Lemma{lemSqfDec}
\begin{lemma}\label{lemSqfDec} \emph{(\Algo de \fcn sans carrés)}\index{algorithme de factorisation sans carrés}% 
\index{factorisation!sans carrés} 
\\
Si~$\gK$ est un \cdi parfait, on dispose d'un \algo de \emph{\fcn sans carrés}
des listes de \pols de $\KX$ au sens suivant. Une \fcn sans carrés d'une
famille  $(g_1,\ldots,g_r)$ est donnée par:

\begin{itemize}
\item une famille  $(\lfs)$
de \pols \spls deux à deux étrangers,
\item  l'écriture de chaque $g_i$ sous
forme

\snic{g_i=\prod_{k=1}^sf_k^{m_{k,i}}\; (m_{k,i}\in\NN) .}
\end{itemize}%
\end{lemma}
%--------- fin lemma ---------------------------------------------- 
%

\begin{Proof} {Idée de la \demz. } On commence par calculer une
\bdf pour la famille $(g_i)_{i\in\lrbr}$ (voir le lemme~\ref{lemPartialDec}). 
Si certains des \pols dans la base sont de la forme $h(X^p)$,
on sait les écrire sous forme $g(X)^p$, on remplace alors $h$ par $g$.
On itère ce processus jusqu'à ce que tous les \pols de la famille
aient une dérivée non nulle. On introduit alors les dérivées des
\pols de la famille. Pour cette nouvelle famille on calcule une nouvelle
\bdfz. \\
On itère le processus d'ensemble jusqu'à ce que l'objectif de départ soit
atteint. Les détails sont laissés \alecz. 
\end{Proof}

 Un \cdiz~$\gK$ est dit \emph{\splz ment clos} 
 si tout \polu \spl de $\KX$ se décompose
 en produit de facteurs $X-x_i$ ($x_i\in\gK$).
 
 Soient~$\gK\subseteq\gL$  des \cdisz. On dit que \emph{$\gL$ est une clôture \spl
 de~$\gK$} si~$\gL$ est \splz ment clos et \agq \spl sur~$\gK$.%
\index{corps!separablement@séparablement clos}%
\index{cloture@clôture!séparable}%
\index{separablement clos@séparablement clos!corps discret ---}
 
%:     Lemma{lemSepParfait}
\begin{lemma}\label{lemSepParfait}~
\begin{enumerate}
\item Un \cdi est \agqt clos \ssi il est parfait et \splz ment clos.
\item Si un \cdiz~$\gK$ est parfait, tout corps étale sur~$\gK$ est parfait.
\item Si un \cdi parfait possède une clôture \spbz, c'est aussi
une clôture \agqz.
\end{enumerate}
 \end{lemma}
%--------- fin lemma ---------------------------------------------- 
%
\begin{proof}
\emph{1.} Résulte du lemme \ref{lemSqfDec} et \emph{3} résulte de 
\emph{1} et \emph{2.}

\emph{2.} On considère~$\gL$ étale sur~$\gK$. 
On note $\sigma:\gL\to\gL:z\mapsto z^p$. \\
On sait que~$\gL=\gK[x]\simeq\aqo\KX f$
avec~$f$ le \polmin de $x$ sur~$\gK$. L'\eltz~$y=x^p$ est zéro du \pol $f^{\sigma}$, qui est \spl et \ird sur~$\gK$ parce que $\sigma$ est un \auto de~$\gK$. 
On obtient donc un \iso $\aqo\KX{f^\sigma}\to \gK[y]\subseteq\gL$. Ainsi~$\gK[y]$ et~$\gL$ sont des \Kevs de même dimension, donc~$\gK[y]=\gL$ et $\sigma$ est surjectif. 
\end{proof}
%
 
%:     Theorem{thClsep}
\begin{theorem}\label{thClsep}
Soit~$\gK$ un \cdi \splz ment factoriel et dénombrable. 
\begin{enumerate}
\item $\gK$ possède une clôture \splz~$\gL$,
et toute clôture \spl de~$\gK$ est~\hbox{$\gK$-isomorphe} à~$\gL$.

\item Ceci s'applique pour~$\gK=\QQ$,  
$\QQ(\Xn)$, $\FF_p$ ou  $\FF_p(\Xn)$.

\item Si en outre~$\gK$ est parfait, alors~$\gL$ est une clôture \agq
de~$\gK$ et toute clôture \agq de~$\gK$ est~$\gK$-isomorphe à~$\gL$.
\end{enumerate}
\end{theorem}
%--------- fin theorem ----------------------------------------------
%
\begin{proof}
Nous donnons seulement une esquisse de \dem du point \emph{1.}
\\
Rappelons tout d'abord le point \emph{2} du \thref{propUnicCDR}: si un \cdr 
pour $f\in\KX$ existe et est \stf sur~$\gK$, alors tout autre \cdr pour~$f$
sur~$\gK$ est isomorphe au premier.
\\
Admettons un moment que l'on sache construire un \cdr \stf pour tout \pol \spl sur~$\gK$.
On énumère tous les \polus \spls de $\KX$ en une suite infinie $(p_n)_{n\in\NN}$.
On appelle $f_n$ le ppcm des \pols $p_0$, \ldots, $p_n$.
On construit des \cdr successifs~$\gK_0$, \ldots,~$\gK_i$, $\ldots$ pour ces $f_i$.
\\
En raison du résultat évoqué précédemment, on sait construire des \homos injectifs de \Klgsz, 
$$\gK_0\vers{\jmath_1}\gK_1\vers{\jmath_2}\cdots\cdots\vers{\jmath_n}\gK_n\vvers{\jmath_{n+1}}\cdots$$
La clôture \spl de~$\gK$ est alors la limite inductive du \sys ainsi construit.
\\
 Il reste à voir pourquoi on sait construire un \cdr \stf pour tout \pol \splz~$f$ sur~$\gK$. Si le corps est infini cela est donné par le \thref{thResolUniv}. Dans le cas d'un corps fini, l'étude des corps finis montre directement comment construire un \cdrz.
Dans le cas le plus \gnlz, on peut de toute manière construire un \cdr
par force brute,  
en ajoutant des racines l'une après l'autre: on considère un facteur irréductible $h$ de~$f$ et le corps~$\gK[\xi_1]=\aqo{\KX}{h}$. 
Sur le nouveau corps~$\gK[\xi_1]$, on considère un facteur irréductible 
$h_1(X)$ de 
$f_1(X)={f(X)\over X-\xi_1} $ ce qui permet de construire~$\gK[\xi_1,\xi_2]$ etc \ldots\ 
Ce processus est possible en vertu du 
corolaire~\ref{corth1Etale} car les corps successifs~$\gK[\xi_1]$,~$\gK[\xi_1,\xi_2]$ \ldots\ restent \splz ment factoriels.   
\end{proof}

\rem Il existe de nombreuses manières de construire une clôture \agq de $\QQ$.
Celle qui est proposée dans le \tho précédent dépend de l'énumération que l'on choisit pour les \polus \spls de $\QQX$ et elle manque de pertinence \gmqz. 
De ce point de vue, la limite inductive que l'on construit présente en fait nettement moins d'intérêt que les \cdr particuliers
 que l'on peut construire chaque fois que le besoin s'en fait sentir.
\\
Il existe d'autres constructions, de nature \emph{\gmqz}, 
de  clôtures \agqs de $\QQ$
qui, elles, sont intéressantes en tant qu'objets globaux. 
La plus connue est celle
via le corps des nombres réels \agqs auquel on ajoute
un \elt $i=\sqrt{-1}$. 
\\
Pour chaque nombre premier $p$, une autre clôture \agq de $\QQ$ \egmt très pertinente est obtenue en passant par le corps 
intermédiaire formé par les nombres \agqs $p$-adiques.
\eoe

%%%%%%%%%%%%%%%%%%%%%%%%%%%%%%%%%%%%%%%%%%%%%%%%%%%%%%%%%%%%%%%%%%%%%%%%%%%
%      sec{Théorie de Galois de base}
\section{Théorie de Galois de base (2)} \label{sec2GaloisElr}

Cette section complète la section \ref{secGaloisElr}
(voir aussi les \thosz~\ref{corlemEtaleEtage} et~\ref{thEtalePrimitif}).

\smallskip 
\emph{Quelques rappels.} Une extension galoisienne de~$\gK$ est 
définie comme un corps \stf
sur~$\gK$ qui est un \cdr pour un \pol \spl de $\KT$.
D'après le \thref{thEtalePrimitif}
 une extension galoisienne de~$\gK$ relève toujours
du cas \elr étudié dans le \thref{thGaloiselr}.
Enfin, d'après le \thref{propUnicCDR}, un tel  \cdr est unique à un \iso  près. \eoe

%--- Definition{defNormale}-----
\begin{definition}
\label{defNormale}
Un surcorps~$\gL$ de~$\gK$ est dit \emph{normal} (sur~$\gK$) si tout $x\in\gL$
annule un \polu de~$\gK[T]$
qui se décompose en produit de facteurs \lins dans~$\gL[T]$.
%Il est dit
%\emph{strictement normal}, si pour tout $x\in\gL$,~$\gK[x]$ est une
%extension \stfe de~$\gK$, et le \polmin de $x$ sur~$\gK$ (qui existe donc) a
%toutes ses racines dans~$\gL$.
\end{definition}
%--- end-definition------------------------------------
\index{normal!surcorps ---}
\perso{extension normale: la \dfn est celle de Richman \cite{MRR}.}

\rem Notez que si~$\gL$ est une extension \stfe de~$\gK$ ou plus \gnlt si~$\gL$
possède une base discrète comme \Kevz,
alors  le \polmin d'un \elt arbitraire de~$\gK$ existe. Si la condition de la \dfn ci-dessus est vérifiée, le \polmin lui-même se décompose en
facteurs \lins dans~$\gL[T]$.
\eoe

%:     Fact{factNormal1}
\begin{fact}\label{factNormal1}
Soit $f(T)\in\gK[T]$ un \polu et~$\gL\supseteq\gK $  un corps de racines pour~$f$.
%~$\gK$  engendré par les racines
%d' (i.e.,~$\gL=\Kxn$ où $f(T)=\prod_{i=1}^n(T-x_{i})$).
Alors,~$\gL$ est une extension normale de~$\gK$.
\end{fact}
\begin{proof}
On a~$\gL=\Kxn$ où $f(T)=\prod_{i=1}^n(T-x_{i})$. Soit $y=h(\xn)$ un \elt arbitraire de~$\gL$.
On pose

\snic{g(\Xn,T)=\prod_{\sigma\in\Sn}\big(T-h^{\sigma}(\uX)\big).}

%\sni
On a clairement $g(\ux,y)=0$. En outre, $g(\ux,T)\in\gK[T]$, car chacun des
\coes de $g(\uX)(T)$ dans $\KuX$ est un \pol symétrique en
les $X_{i}$, donc un \pol en les fonctions symétriques \elrsz,
qui se spécialisent en des \elts de~$\gK$ (les \coes de~$f$)
par le~$\gK$-\homo $\uX\mapsto\ux$.
\end{proof}

%:     theorem{thNormalEtaleGalois}
\begin{theorem} \emph{(\Carn des extensions galoisiennes)}\\
Soit~$\gL$ un corps \stf sur~$\gK$. 
\Propeq\label{thNormalEtaleGalois}
\begin{enumerate}
\item \label{i1thNormalEtaleGalois}
$\gL$ est une extension galoisienne de~$\gK$.
\item \label{i2thNormalEtaleGalois} $\gL$ est étale et normal sur~$\gK$.
\item \label{i3thNormalEtaleGalois} $\Aut_\gK(\gL)$ est fini et
la correspondance galoisienne est bijective.
\item \label{i4thNormalEtaleGalois} Il existe un groupe
fini $G\subseteq\Aut_\gK(\gL)$ dont le corps fixe est~$\gK$.
%-% PERSO
\perso{NB: \gui{$\Aut_\gK(\gL)$ est fini} \emph{ne résulte pas} de
 \gui{$\gL$  \stf sur~$\gK$} (exemple de Richman, on construit~$\gK\supseteq\QQ$
 qui contient \emph{peut-être} une racine cubique de 1
 et l'on prend~$\gL=\gK[\alpha]$ où $\alpha^3=2$)}
%-% Fin PERSO
%
\end{enumerate}
Dans ce cas,
 dans le point \ref{i4thNormalEtaleGalois}, on a \ncrt $G=\Gal(\gL/\gK)$.
\end{theorem}
\begin{proof}
\emph{\ref{i1thNormalEtaleGalois}} $\Rightarrow$  \emph{\ref{i2thNormalEtaleGalois}.}
 C'est le fait
\ref{factNormal1}.
 
\emph{\ref{i2thNormalEtaleGalois}} $\Rightarrow$  \emph{\ref{i1thNormalEtaleGalois}}
et \emph{\ref{i3thNormalEtaleGalois}.}
Par le \tho de l'\elt primitif,~$\gL=\gK[y]$ pour un $y$ dans~$\gL$.
Le \polminz~$f$ de~$y$ sur~$\gK$ est \splz, et~$f$
se factorise complètement dans~$\gL[T]$ parce que~$\gL$ est normal sur~$\gK$.
Donc~$\gL$ est un \cdr pour~$f$. En outre, le \thrf{thGaloiselr} s'applique.

\emph{\ref{i4thNormalEtaleGalois}} $\Rightarrow$  \emph{\ref{i2thNormalEtaleGalois}.} Il suffit de montrer que tout $x\in\gL$
annule un \pol \spl de $\KT$ qui se factorise complètement dans~$\gL[T]$,
car alors l'extension est normale (par \dfnz)
et étale (\thrf{corlemEtaleEtage}).
Posons

\snic{P(T)=\Rv_{G/H,x}(T)=\prod_{\sigma\in G/H}\big(T-\sigma(x)\big)\;$ où $\;H=\St(x).}

%\sni
En indice, l'expression $\sigma\in G/H$ signifie que l'on prend un $\sigma$ dans chaque classe à gauche modulo $H$. Le \pol $P$ est fixé par $G$, donc $P\in\KT$.
Par ailleurs,   $\disc (P)=\prod_{i,j\in\lrbk,i< j}(x_i-x_j)^2$ est \ivz.

Enfin, vu que la correspondance galoisienne est bijective, et puisque le corps fixe de $G$ est~$\gK$, dans le point \emph{\ref{i4thNormalEtaleGalois}},
on a \ncrt $G=\Gal(\gL/\gK)$.
\end{proof}
%

%
%:    Theorem{thSG}  correspondance galoisienne, complément
\begin{theorem}
\label{thSG} \emph{(Correspondance galoisienne, complément)}
\\
Soit~$\gL/\gK$ une extension galoisienne de groupe de Galois $G=\Gal(\gL/\gK)$.
Soient $H$ un sous-groupe détachable de~$G$, $\sigma$ un \elt de~$G$,~%
%.  Posons
$H_\sigma=\sigma H\sigma^{-1}$.
\begin{enumerate}
\item Le corps $\sigma(\gL^H)$ est égal à~$\gL^{H_\sigma}$.
\item $\gL^H$ est une extension galoisienne de~$\gK$ \ssi $H$ est normal
dans $G$. Dans ce cas le groupe de Galois $\Gal(\gL^H/\gK)$ 
est canoniquement isomorphe à $G/H$.
\end{enumerate}
\end{theorem}
%--- end-theorem-----------------------------------------
%-----------------begin proof------------------
\begin{proof}
\emph{1.} Calcul immédiat.

\emph{2.} Posons~$\gM=\gL^H$. 
Par le \tho de l'\elt primitif écrivons~$\gM=\gK[y]$, de sorte que $H=\St(y)$.
Le corps~$\gM$ est normal sur~$\gK$ \ssi pour  chaque $\tau\in G$, on a $\tau(y)\in\gM$,
\cad
$\tau(\gM)=\gM$. Et d'après le point  \emph{1} cela signifie
$\tau H \tau^{-1}=H$.
\end{proof}
%-----------------end proof------------------

Nous reprenons le \thrf{thGaloiselr} en apportant quelques précisions.

%:    Theorem{thCGSynthese}  correspondance galoisienne
\begin{theorem}
\label{thCGSynthese} \emph{(Correspondance galoisienne, synthèse)}
\\
Soit~$\gL/\gK$ une extension galoisienne.
La correspondance galoisienne fonctionne comme suit.
\begin{enumerate}
\item Pour tout~$\gM\in \cK_{\gL/\gK}$,~$\gL/\gM$ est une extension
galoisienne de groupe de Galois $\Fix(\gM)$ et $\dex{\gL:\gM}=\#{\Fix(\gM)}$.
\item Si $H_1, H_2\in\cG_{\gL/\gK}$ et~$\gM_i=\Fix(H_i)\in \cK_{\gL/\gK}$, alors:
\begin{itemize}
\item  $H_1\cap H_2$
correspond à la sous-\Klg engendrée par~$\gM_1$ et~$\gM_2$,
\item  $\gM_1\cap \gM_2$
correspond au sous-groupe engendré par $H_1$ et $H_2$.
\end{itemize}
\item Si $H_1\subseteq H_2$, %dans $\cG_{\gL/\gK}$ et~$\gM_i=\Fix(H_i)$
alors:
\begin{itemize}
\item $\gM_1\supseteq\gM_2$ et $(H_2:H_1)=\dex{\gM_1:\gM_2}$,
\item $\gM_1/\gM_2$ est une extension galoisienne \ssi $H_1$ est normal dans $H_2$.
Dans ce cas le groupe $\Gal(\gM_1/\gM_2)$ est naturellement isomorphe à $H_2/H_1$.
%
%\item
\end{itemize}
\end{enumerate}
\end{theorem}
%--- end-theorem-----------------------------------------
%

%%%%%%%%%%%%%%%%%%%%%%%%%%%%%%%%%%%%%%%%%%%%%%%%%%%%%%%%%%%%%%%%%%%%%%%%%%%
%     subsec{Algèbres \pfz}
\section{Algèbres \pfz} \label{sec1Apf}

\vspace{4pt}
%:subsec{Généralités}   %%%%%%%%%%%%%%%%%%%%%
\subsec{Généralités}

Les \algs \pf sont aux \syss d'équations polynomiales (ou
\emph{\sypsz}) ce que sont les \mpfs aux \slisz.\index{systeme polynomial@\sypz}

Nous présentons ici quelques faits \gnls de base concernant ces \algsz.

Les \algs que nous considérons dans cette section sont
associatives, commutatives et unitaires.

%--- Definition{defAlg}---------
\begin{definition}
\label{defAlg}\label{defi2STF} Soit~$\gA$ une \klgz.
%-----------------begin enum------------------
\begin{enumerate}

\item L'\algz~$\gA$ est dite \emph{\tfz} si elle est engendrée par une famille finie en tant que~\hbox{\klgz}. Ceci
revient à dire qu'elle est isomorphe à une \alg quotient $\kXn\sur{\fa} $. On note alors~$\gA=\kxn$, où $x_i$ est l'image de $X_i$ dans
$\gA$.
Cette notation ne sous-entend pas que~$\gA$ est une extension de~$\gk$.%
\index{algèbre!de type fini}

\item L'\algz~$\gA$ est dite \emph{\pfz} si elle est \pf en tant que \klgz. Ceci
qui revient à dire qu'elle est isomorphe à une \alg 
$\kXn\sur{\fa}$, avec un \itf $\fa=\gen{\lfs}$.%
\index{algèbre!de présentation finie}

\item L'\algz~$\gA$ est dite \emph{\rpfz} (en un seul mot) si elle est \pf en
tant que \klg réduite. Autrement dit si elle est isomorphe à une
\alg quotient $\kXn/\sqrt{\fa}$ avec un \itfz~$\fa$.%
\index{algèbre!réduite-de-présentation-finie}

\item L'\alg $\gA$ est dite \emph{strictement finie}
si~$\gA$ est un \kmo \ptfz. On dit aussi que \emph{$\gA$ est \stfe sur~$\gk$}.
Dans le cas d'une extension, on parlera d'\emph{extension \stfez} de~$\gk$.%
\index{algèbre!strictement finie}%
\index{strictement finie!algèbre --- }

\item Si~$\gA$ est \stfe   on note

\snic{\Tr\iAk (x),\;\; \rN\iAk (x), \;\;
\rF{\gA/\gk}(x)(T) \;\;\hbox{et}\;\; \rC{\gA/\gk}(x)(T),}

%\sni
la trace, le
\deterz, le \polfon et le \polcar de l'\kli $\mu_{\gA,x}\in\End_\gk(\gA)$.
En outre, en notant $g(T)=\rC{\gA/\gk}(x)(T)$, l'\elt $g'(x)$ est appelé \emph{la
différente de~$x$}.%
\index{differente@différente!d'un \elt dans une \alg strictement finie}
\end{enumerate}
%-----------------end enum------------------
\end{definition}
%--- end-definition------------------------------------

Notez que dans le cas où~$\gk$ est un \cdiz, on retrouve bien
la notion d'\alg \stfe donnée dans la \dfnz~\ref{defiSTF}.

%:     Fact{factPropUnivAPF}-----------
\begin{fact}
\label{factPropUnivAPF} \emph{(\Prt universelle d'une \apfz)}
\\
L'\alg \pf $\aqo{\kXn}{\lfs}=\kxn$ est \caree par la \prt
suivante: si une \klgz~$\gk\vers{\varphi}\gA$ contient des \elts $\yn$
tels que les $f_i^\varphi(\yn)$ sont nuls, il existe un unique
\homo de \klgs $\kxn\to\gA$ qui envoie les $x_i$ sur les~$y_i$.
\end{fact}
%--- end-fact-----------------------------------------

%%%%%%%%%%%%%%%%%%%%%%%%%%%%%%%%%%%%%%%%%%%%%%%%%%%%%%%%%%%%%%%%%%%%%%%%%%%
\subsubsection*{Changement de \sgrz}

%:     Fact{factChscalg}-------------
\begin{fact}
\label{factChscalg}
Lorsque l'on change de \sys \gtr pour une \alg \pfz~$\gA$ les relations entre les
nouveaux \gtrs forment de nouveau un \itfz.
\end{fact}
%--- end-fact-----------------------------------------

Vous pouvez vous reporter à la section \ref{sec pf chg} et vérifier
que ce qui a été expliqué un peu informellement \paref{nouveausgr}
fonctionne
bien dans le cas présent.

%%%%%%%%%%%%%%%%%%%%%%%%%%%%%%%%%%%%%%%%%%%%%%%%%%%%%%%%%%%%%%%%%%%%%%%%%%%
%\penalty-2500
\subsubsection*{Transitivité (\apfsz)}

%\vspace{-10pt}

%--- Fact{factTransAPF}-------------
\begin{fact}
\label{factTransAPF}
Si~$\gk\vers \lambda \gA$ et~$\gA\vers\rho \gC$  sont deux \algs \pfz, alors
$\gC$ est une \klg \pfz.
\end{fact}
%--- end-fact-----------------------------------------
%-----------------begin proof------------------
\begin{proof}
\'Ecrivons
$\gA=\kuy\simeq\aqo{\kuY}{g_1,\ldots ,g_t}$
et~$\gC=\Aux\simeq\aqo{\AuX}{\lfs}$.
Soient $F_{1},\ldots,F_{s}\in\gk[\uY,\uX]$ des \pols tels que $F_{i}(\uy,\uX)=f_{i}(\uX)$.
\\
Alors,~$\gC=\gk[\rho(\uy),\ux]\simeq\aqo{\gk[\uY,\uX]}{g_1,\ldots ,g_t,F_{1},\ldots,F_{s}}$.
\end{proof}
%-----------------end proof------------------

%%%%%%%%%%%%%%%%%%%%%%%%%%%%%%%%%%%%%%%%%%%%%%%%%%%%%%%%%%%%%%%%%%%%%%%%%%%
\subsubsection*{Sous-algèbres}

%\vspace{-10pt}

%--- Fact{factSousAlg}-------------
\begin{fact}
\label{factSousAlg}
Soient~$\gA\subseteq\gC$ deux \klgs \tfz.
Si~$\gC$ est une \klg \pf c'est aussi une \Alg \pf
(avec \gui{la même} \pnz, lue dans~$\gA$).
\end{fact}
%--- end-fact-----------------------------------------
%-----------------begin proof------------------
\begin{proof}
\'Ecrivons \spdg $\gC=\kxn\simeq\aqo{\kuX}{\uf}$ et~$\gA=\kxr$.
On a~$\gA\simeq \kXr/\ff$ avec $$\ff=\gen{\lfs}\cap\kXr.$$
Notons $\pi:\kXr\to\gA$ le passage au quotient et pour $h\in\kXn$, 
$h^\pi\in\gA[X_{r+1},\ldots ,X_n]$ son image:
 $$h^\pi=h(\xr,X_{r+1},\ldots ,X_n).$$
On considère l'\homo
$$
\gamma:
\begin{array}{l}
 \aqo{\gA[X_{r+1},\ldots ,X_n]}{f_1^\pi,\ldots ,f_s^\pi} \simeq \\[1mm]
 \aqo{\AXn}{X_1-x_1,\ldots ,X_r-x_r,f_1^\pi,\ldots ,f_s^\pi}
\end{array}
\to\gC.
$$
C'est l'\homo qui fixe~$\gA$ et envoie $X_k$ sur $x_k$ pour $k\in\lrb{r+1..n}$.
Il suffit de montrer que~$\gamma$ est injectif. Tout \elt
$g$ de~$\gA[X_{r+1},\ldots ,X_n]$ peut s'écrire $g=G^\pi$
avec $G\in\kXn$.
\\
Supposons que $g$ modulo~$\gen{f_1^\pi,\ldots ,f_s^\pi}$ soit dans $\Ker\gamma$.
On a alors
$$\preskip.4em \postskip.4em 
g(x_{r+1},\ldots ,x_n)=G(\xn)=0. 
$$
Donc $G\in\gen{\lfs}$, ce qui donne 
$g\in\gen{f_1^\pi,\ldots ,f_s^\pi}$ (après transformation par $\pi$). Ce que nous voulions.
\end{proof}
%-----------------end proof------------------

\rem La condition~$\gA\subseteq\gC$ est indispensable pour le bon fonctionnement
de la preuve. Par ailleurs, il faut noter que l'\id $\ff$
n'est pas \ncrt de type fini.
\eoe

%%%%%%%%%%%%%%%%%%%%%%%%%%%%%%%%%%%%%%%%%%%%%%%%%%%%%%%%%%%%%%%%%%%%%%%%%%%
%:subsec{Les zéros d'un \sypz}
%\penalty-2500
\subsec{Les zéros d'un \sypz}  \label{secZerosSyp}
Considérons un \syp $(\uf)=(\lfs)$ dans $\kXn,$ et une
 \klgz~$\rho:\gk\to\gB$. 
 
%:     Definition{defiZeroSyp}
\begin{definition}\label{defiZeroSyp}
Un \emph{zéro du \sys $(\uf)$ sur~$\gB$} est un
$n$-uplet 

\snic{(\uxi)=(\xin)\in\gB^n}

%\sni
vérifiant $f_i^\rho(\uxi)=0$  pour chaque $i$.
L'ensemble des zéros de $(\uf)$ sur~$\gB$ est souvent appelé, de
manière imagée, la
\ixc{variété des zéros}{d'un \sypz} sur~$\gB$ du \sypz,
et pour cela, on le note $\cZ_\gk(\uf,\gB)$ ou $\cZ(\uf,\gB)$.
\index{zero@zéro!d'un \sypz, sur une \algz}
\end{definition}
%--------- fin definition ---------------------------------------------- 

Certains zéros sont plus intéressants que d'autres: plus l'\algz~$\gB$ est
proche de~$\gk$ et plus le zéro est intéressant. On est particulièrement
attentif aux zéros sur~$\gk$, ou à défaut sur des \klgs finies.

Deux zéros sont a priori
particulièrement décevants. Celui fourni par l'\alg
triviale, et le zéro $(\xn)$ sur
l'\ixc{algèbre quotient}{pour un \sypz} associée au \sypz, \cad
\Grandcadre{$\gA=\kxn = \aqo{\kXn}{\lfs}.$}

Néanmoins cette dernière \alg
joue un rôle central pour notre \pb en raison de deux
constatations. La première est la suivante.

%--- Fact{factZeros}-----------------
\begin{fact}
\label{factZeros}
Pour toute \klgz~$\gB$ l'ensemble des zéros de $(\uf)$ sur~$\gB$ s'identifie
naturellement à l'ensemble des morphismes de \klgs de~$\gA$ vers~$\gB$.
En particulier, les zéros sur~$\gk$ s'identifient aux carac\-tères
de l'\algz~$\gA$.
\end{fact}
%--- end-fact-----------------------------------------

%
\begin{Proof}{\Demo sur un exemple. }
Posons $\QQ[x,y]=\aqo{\QQ[X,Y]}{X^2+Y^2-1}$. 
Il revient au même de se donner un
point réel $(\alpha,\beta)$ du cercle $X^2+Y^2=1$ ou un morphisme
$\rho:\QQ[x,y]\vers{}\RR$
(celui qui envoie $x$ et $y$ sur $\alpha$ et~$\beta$).
\end{Proof}

\rdb 
On a donc une identification cruciale, que nous écrivons comme une \egtz:
\Grandcadre{$\Hom_{\gk} (\gA,\gB)=\cZ_{\gk}(\uf,\gB)\subseteq \gB^n.$
\label{ZerosCrucial}}

En bref l'\alg quotient~$\gA$ résume \emph{de manière intrinsèque} les informations
pertinentes contenues dans le \syp $(\uf)$.
Ce pourquoi on dit aussi que $\cZ_\gk(\uf,\gB)$ est la 
\ixc{variété des zéros}{d'une \alg sur une autre} de~$\gA$ sur~$\gB$.

\medskip
La seconde constatation (étroitement reliée à la précédente
 d'ailleurs) est la suivante.

Du point de vue \gmq \emph{deux \syps $(\uf)$ et $(\ug)$ dans~$\kuX$
 qui ont les mêmes zéros},
sur n'importe quelle \klgz, doivent être considérés comme \emph{\eqvsz.}
Si tel est le cas, notons~$\gA_{1}=\kux$ et~$\gA_{2}=\kuy$ les deux \algs quotients
(on ne donne pas le même nom aux classes des~$X_i$
dans les deux quotients).
Considérons le zéro canonique~$(\xn)$ de~$(\uf)$ dans~$\gA_1$.
 Puisque $\cZ(\uf,\gA_1)=\cZ(\ug,\gA_1)$,
 on doit avoir $g_j(\ux)=0$ pour chaque~$j.$ Cela signifie
que $g_j(\uX)$ est nul modulo~$\gen{\uf}$%
.
De même, chaque $f_i$ doit être dans~$\gen{\ug}$.

Résumons cette deuxième constatation.

%--- Fact{fact1Zeros}-----------------
\begin{fact}
\label{fact1Zeros}
Deux \syps $(\uf)$ et $(\ug)$ dans $\kuX$
admettent les mêmes zéros,
sur n'importe quelle \klgz, \ssi ils définissent
la même \alg quotient.
\end{fact}
%--- end-fact-----------------------------------------

\medskip
%  exl
\exl
Les cercles
$x^2+y^2-3=0$ et $x^2+y^2-7=0$
ne peuvent pas être distingués par leurs points rationnels, ils n'en ont
pas (puisque sur $\ZZ$, la congruence $a^2 + b^2 \equiv 3c^2 \bmod 4$ 
 entraîne $a, b, c$ pairs), mais les \QQlgs quotients
sont non isomorphes, et l'on peut constater sur $\QQ[\sqrt{3},\sqrt{7}]$
qu'ils ont des zéros distincts \gui{quelque part}.
\eoe

\smallskip
Lorsque~$\gk$ est réduit et si l'on s'intéresse particulièrement aux zéros
sur les \klgs réduites, l'\algz~$\gA=\aqo{\kuX}{\uf}$  doit être
remplacée par sa variante réduite, qui est une \alg \rpfz:

\snic{\gA\sur{\DA(0)}=\kXn/\sqrt{\gen{\lfs}}\,.}

%\sni
     Nous poursuivrons cette discussion \paref{subsecNstMorphismes} 
      dans le paragraphe
      \gui{\nst et équivalence de deux catégories}.

%%%%%%%%%%%%%%%%%%%%%%%%%%%%%%%%%%%%%%%%%%%%%%%%%%%%%%%%%%%%%%%%%%%%%%%%%%%
%   subsec{Digression sur le calcul \agqz}
\subsubsection*{Digression sur le calcul \agqz}
Outre leur rapport direct avec la résolution des \syps une autre raison de
l'importance des \algs \pf est la suivante.
Chaque fois qu'un calcul \agq aboutit à un
\gui{résultat intéressant} dans une \klgz~$\gB$, ce calcul n'a fait intervenir
qu'un nombre fini d'\elts $\yn$ de~$\gB$ et un nombre fini de relations
entre les $y_i$, de sorte qu'il existe une \klg \pfz~$\gC=\kxn$ et un
morphisme surjectif $\theta:\gC\to\kyn\subseteq\gB$ qui envoie les $x_i$ sur les $y_i$
et tel que le \gui{résultat intéressant} avait déjà lieu dans~$\gC$
pour les $x_i$.
En langage plus savant\footnote{\Llec notera que le paragraphe présent
est directement recopié du paragraphe analogue pour les \mpfsz, \paref{DigCalcAlg}.}: toute \klg est une limite inductive filtrante de \klgs
\pfz.

%:subsec{Produit tensoriel de deux \klgs}   %%%%%%%%%%%%%%%%%%%%%
%\penalty-2500
\subsec{Produit tensoriel de deux \klgs}

\vspace{3pt}
La \emph{somme directe} de deux \klgsz~$\gA$ et~$\gB$ dans la catégorie des \klgs
est donnée par la solution du \pb universel suivant (\gui{morphisme} signifie ici \gui{\homo de \klgsz}).%
\index{somme directe!dans une catégorie}
\\
\emph{Trouver une \klgz~$\gC$ et deux morphismes $\alpha:\gA\to\gC$ et
$\lambda:\gB\to\gC$ tels que, pour toute \klgz~$\gD$ et pour tout couple de
 morphismes $\varphi:\gA\to\gD$ et $\psi:\gB\to\gD$,
il existe un unique morphisme~$\gamma:\gC\to\gD$ tel que
$\varphi=\gamma\circ\alpha$ et $\psi=\gamma\circ\lambda$.}

\vspace{-1.8em}
$$
\xymatrix @R=15pt{
 & \gA\ar[dr]_{\alpha} \ar[drrr]^{\varphi}
\\
\gk\ar@{->}[ru]^{\beta} \ar@{->}[rd]_{\rho} &&\gC\ar@{-->}[rr]^{\gamma !~~~~~~} && \gD
\\
 & \gB\ar[ur]^{\lambda} \ar[urrr]_{\psi}
\\
}
$$

Notons que  dans la catégorie des anneaux commutatifs,
 la \prt universelle ci-dessus signifie que~$\gC$,
avec les deux morphismes $\alpha$ et $\lambda$,
 est la \ixc{somme amalgamée}{de deux flèches de même source dans une catégorie}
des deux flèches $\beta:\gk\to\gA$ et $\rho:\gk\to\gB$.
En anglais on dit que~$\gC$ est le push-out de $\beta$ et $\rho$.
En français on dit encore que l'on a un \emph{carré cocartésien}, formé avec les 4 flèches $\beta$, $\rho$,~$\alpha$ et~$\lambda$.

%:     theorem  factSDIRKlg
\begin{theorem}\label{factSDIRKlg}%\label{factChsgalg}
On considère deux \klgsz~$\gk\vers{\rho}\gB$ et~$\gk\vers{\beta}\gA$.

 A. \emph{(Somme directe dans la catégorie des \klgsz)}\\
Les \algs $\gA$ et $\gB$
admettent une somme directe~$\gC$ dans la catégorie \hbox{des \klgsz}.
En voici différentes descriptions possibles:
\begin{enumerate}
 \item Si~$\gA=\aqo{\kXn}{\lfs}$,~$\gC=\aqo{\gB[\Xn]}{f_1^\rho,\ldots,f_s^\rho}$
 avec les deux morphismes naturels~$\gA\to\gC$ et~$\gB\to\gC$.
 \item Si en outre~$\gB=\kyr\simeq\aqo{\kYr}{g_1,\ldots ,g_t}$ est elle-même
une \klg \pfz, on obtient
$$\gC\simeq\aqo{\gk[\Xn,\Yr]}{\lfs,g_1,\ldots ,g_t}.$$
 \item En \gnlz, on peut considérer le \kmoz~$\gC=\gB\otimes_{\gk}\gA$.
 Il est muni d'une structure d'anneau commutatif en définissant
 le produit par

\snic{(x\otimes a)\cdot(y\otimes b)=xy\otimes ab.}

\snii
On obtient une structure de \klg
 et l'on a deux morphismes naturels~$\gB\to\gC,\,x\mapsto x\otimes 1$ et~$\gA\to\gC,\,a\mapsto 1\otimes a$. Ceci fait de~$\gC$ la somme directe
 de~$\gB$ et~$\gA$.
\item \label{SdirQuo} Si~$\gB=\gk\sur{\fa}$, on obtient~$\gC\simeq\gA/\fb$ où $\fb=\beta(\fa)\gA$.
\item Si~$\gB=S^{-1}\gk$, on obtient~$\gC\simeq U^{-1}\gA$ où $U=\beta(S)$.

\end{enumerate}

 B. \emph{(Extension des scalaires)}
\\
On peut voir~$\gC$ comme une \Blgz, on dit alors
 que~$\gC$ est la \Blg obtenue à partir de~$\gA$ par
\ix{changement d'anneau de base}, ou encore par
\ix{extension des scalaires}. Il est logique alors de la noter $\rho\ist(\gA)$.
\end{theorem}
\facile

On prendra garde au fait que~$\gk\subseteq \gA$ n'implique pas en \gnl
$\gB\subseteq \gC$, en particulier dans le cas~\emph{\ref{SdirQuo}}.

Notons aussi que la tradition est de parler de \emph{produit tensoriel de \klgsz} plutôt que de somme directe. \index{produit tensoriel!d'algèbres}

%:     Fact{factBimodule}
\begin{fact}\label{factBimodule}
Si~$\gA$ et~$\gB$ sont deux \klgsz, et $(M,+)$ est un groupe additif,
 se donner une structure de~$\gA\te_\gk\gB$-module sur 
$M$  revient à se donner
une loi externe de \Amoz~$\,\gA\times M\to M\,$
et une loi externe de~\Bmoz~$\,\gB\times M\to M\,$ qui commutent, et qui \gui{coïncident sur~$\gk$}. On dit aussi que $M$ est muni d'une \emph{structure de}
$(\gA,\gB)$-\ix{bimodule}.
\end{fact}
\begin{proof}
L'explication est la suivante avec~$\gk\vers{\rho}\gB$,~$\gk\vers{\alpha}\gA$.
\\
Si l'on a une structure de~$\gA\te_\gk\gB$-module sur $M$, on a les deux lois externes 

\snic{\gB\times M\to M,\, (c,m)\mapsto c\cdot m=(1\te c)m$, et $~}

\snic{\gA\times M\to M,\, (b,m)\mapsto b\star m=(b\te 1)m.}

%\sni
Puisque $b\te c=(b\te 1)(1\te c)=(1\te c)(b\te 1)$, on doit avoir $b\star(c\cdot m)=c\cdot(b\star m)$. Si $a\in\gk$, $a(1\te 1)=\alpha(a)\te1=1\te \rho(a)$ donc on doit avoir
 $\rho(a)\cdot m=\alpha(a)\star m$. Ainsi les deux lois commutent et
 coïncident sur~$\gk$. 
\\
 Inversement, à partir de deux lois externes qui commutent et coïncident sur~$\gk$, on peut définir $(b\te c)m$ par $b\star(c\cdot m)$.
\end{proof}

Voici un fait important, et facile, concernant l'\edsz.
%:     factEdsAlg
\begin{fact}\label{factEdsAlg}
On considère deux \klgsz~$\gk\vers{\rho}\gk'$ et~$\gk\vers{\alpha}\gA$
et l'on note~$\gA'=\rho\ist(\gA)$.
Si la \klgz~$\gA$ est \tf (resp. \pfz, finie, entière, \stfez)
il en va de même pour la~$\gk'$-\algz~$\gA'$.
\end{fact}
\facile

%: subsec{Algèbres entières}
\subsec{Algèbres entières}
\vspace{4pt}
%: subsubsec{Le lemme lying over}
\subsubsec{Le lemme lying over}

Dans ce paragraphe et le suivant 
nous complétons ce qui a déjà été dit sur les \algs entières dans la
section~\ref{secApTDN}. 

\smallskip  Le lemme qui suit exprime le contenu \cof du lemme de \clamaz, appelé
\gui{lying over},
qui affirme que si~$\gB$ est un anneau entier sur un sous-anneau~$\gA$, il y a
toujours un \idep de~$\gB$ au dessus d'un \idep donné de~$\gA$. \index{lying
over}

Rappelons que nous notons $\DA(\fa)$
le nilradical de l'\id $\fa$ de~$\gA$.
%--- Lemma{lemLingOver}--------------
\begin{lemma}\label{lemLingOver} \emph{(Lying over)}\\
 Soit $\gA\subseteq\gB$ avec~$\gB$ entier sur~$\gA$
et $\fa$ un idéal de~$\gA$, alors $\fa\gB\,\cap\,\gA\subseteq\DA(\fa)$,
ou ce qui revient au même

\snic{\rD_\gB(\fa\gB)\,\cap\,\gA=\DA(\fa).}

%\sni
En particulier, $1\in\fa\;\Leftrightarrow\; 1\in\fa\gB$.
\end{lemma}
%--- end-lemma-----------------------------------------
%-----------------begin proof------------------
\begin{proof}
Si $x\in\fa\gB$ on a
$x=\sum a_ib_i$, avec $a_i\in \fa,\;b_i\in \gB$.
Les $b_i$  engendrent une sous-\Algz~$\gB'$ qui est finie.
Soit $G$ un \sys \gtr fini (avec $\ell$ \eltsz) du \Amoz~$\gB'$.
Soit $B_i\in\MM_{\ell}(\gA)$ une matrice qui exprime la multiplication par
$b_i$ sur $G$.
La multiplication par
$x$ est exprimée par la matrice $\sum a_iB_i$, qui est à \coes dans $\fa$.  Le
\polcar de cette matrice, qui annule $x$ (parce que~$\gB'$ est un \Amo fidèle), a
donc tous ses \coes (sauf le \coe dominant) dans~$\fa$. 
Lorsque~$x\in\gA$, ceci implique $x^{\ell}\in\fa$.
\perso{En fait tout $x\in\fa\gB$
est entier sur l'\id $\fa$. Mais la terminologie n'est pas encore introduite.
Voir le lemme~\ref{lemLingOver2}}
\end{proof}
%-----------------end proof------------------

\rem \label{remlyingover}
Indiquons comment on en déduit le lying over classique en \clamaz.
On considère le cas où $\fa$ est un \idep et l'on note $S=\gA\setminus\fa$.
Alors, $\fa\gB\,\cap\, S = (\fa\gB\,\cap\,\gA)\,\cap\, S$ est vide d'après le lemme~\ref{lemLingOver}. D'après le \emph{lemme de Krull}, il existe donc un \idepz~$\fp$
de~$\gB$ tel que $\fa\gB\subseteq\fp$ et $\fp\,\cap\, S=\emptyset$, ce qui implique
$\fp\,\cap\,\gA=\fa$. Il serait \egmt facile de déduire en \clama le lemme
\ref{lemLingOver} du lying over classique.\eoe
\index{Lemme de Krull}

\medskip \exl
On montre ici que la condition \gui{$\gB$ entier sur~$\gA$} est
cruciale dans le lying over. On considère~$\gA=\ZZ$,~$\gB=\ZZ[1/3]$ et
$\fa=3\ZZ$. Alors, on obtient $\fa\gB=\gen{1}$, mais $\fa\neq \gen{1}$.\eoe

%%%%%%%%%%%%%%%%%%%%%%%%%%%%%%%%%%%%%%%%%%%%%%%%%%%%%%%%%%%%%%%%%%%%%%%%%%%
%: subsubsec{Algebres entières  sur un anneau \zedz}
\subsubsec{Algèbres entières  sur un anneau \zedz}\label{EntSurZdim}

Nous  examinons ici le cas particulier
des \algs sur un anneau \zedz. 
%Nous  reviendrons sur ce sujet dans la section \ref{secStfeSurZed}.

Les \algs entières sur les \cdis
sont un exemple important d'anneaux \zedsz.
Dans cette situation, on précise le point \emph{\iref{LID003}}
du lemme~\ref{lemme:idempotentDimension0} comme suit
(voir aussi le \thrf{th1Etale}).

%--- Lemma{lemKlgEntiere}------------
\begin{lemma}
\label{lemKlgEntiere}
Une  \alg entière~$\gA$ sur un \cdi~$\gK$ est  \zedez.
Plus \prmtz, soit
\hbox{$\fa=\gen{a_1,\ldots ,a_n}=\gen{\ua}$}  un \itfz. Il existe un entier $d$ et un \idm $s\in a_1\gK[\ua]+\cdots +a_n\gK[\ua]$ 
tels \hbox{que $\fa^d=\gen{s}$}.
\end{lemma}
%--- end-lemma----------------------------
%-----------------begin proof------------------
\penalty-2500
\begin{proof}
Un \elt $x$ de~$\gA$ est annulé par un \polu de $\KX$ que l'on écrit $uX^k\big(1 - X\,h(X)\big)$
avec $u\in\gK\eti$, $k\geq0$ 
et donc $x^k\big(1-xh(x)\big)=0$. 

L'\idm $e_x$ tel que
$\gen{e_x}=\gen{x}^d$ pour $d$ assez grand est alors égal à~$\big(xh(x)\big)^k$,
et~$d$ est \gui{assez grand} dès que $d\geq k$.

Dans le cas de l'\itf $\fa=\gen{a_1,\ldots ,a_n}$,
chaque \idmz~$e_{a_i}$ est un \elt de $a_i\gK[a_i]$. Donc leur pgcd,
qui est l'\idm $s$ dans l'énoncé du lemme, est dans $a_1\gK[\ua]+\cdots
+a_n\gK[\ua]$ (car le pgcd de deux \idms $e$ et $f$ est $e\vu f=e+f-ef$).
\end{proof}
%-----------------end proof------------------

%--- Lemma{lemZrZr1}----------
\begin{lemma}
\label{lemZrZr1}
Soit~$\gk$  un anneau \zed et~$\gA$ une \klg entière sur~$\gk$.
\begin{enumerate}\itemsep0pt
\item L'anneau~$\gA$ est  \zedz.
\item Plus \prmtz, si % pour un \itf
$\fa=\gen{a_1,\ldots ,a_n}%=\gen{\ua}
$, il existe un entier $d$ et un \idm $s\in
a_1\gk[\ua]+\cdots +a_n\gk[\ua]$ tel que~\hbox{$\fa^d=\gen{s}$}.
\item En particulier, on obtient pour chaque $a\in\gA$ une \egt
$$\preskip.3em \postskip.2em
a^d\big(1-af(a)\big)=0,
$$
%\sni
 avec un $f(X)\in\kX$ (donc, $\big(af(a) \big)^d$ est \idmz).
\end{enumerate}
\emph{NB:} on ne réclame pas que $\rho:\gk\to\gA$ soit injectif.
\end{lemma}
%--- end-lemma-----------------------------------------
%-----------------begin proof------------------
\begin{proof}
Il suffit de montrer le point \emph{2.}\\
En appliquant la machinerie \lgbe \elr \paref{MethodeZedRed}, on étend le
résultat du lemme~\ref{lemKlgEntiere} au cas où~$\gk$ est \zed réduit.\imlg Ensuite
on ramène le cas \zed au  cas \zed réduit en passant au quotient par le
nilradical et en utilisant \gui{la \imN méthode de Newton en \algz}
(section~\ref{secNewton}). Plus \prmtz, posons  $\fN=\DA(0)$. D'après le cas \zed réduit, il existe $x_1$, \dots, $x_n$ $\in\gk[\ua]$ tels que
$$
s=a_1x_1+\cdots +a_nx_n,\hbox{ avec }s^2\equiv s  \mod  \fN\hbox{ et }sa_i\equiv a_i
\mod \fN.
$$
%\sni
L'\elt $s$ est congru modulo $\fN$  à un unique \idm $s_1$, lequel
s'écrit~$sp(s)$ avec $p(T)\in\ZZ[T]$
(corolaire~\ref{corIdmNewton}). Puisque~$s\in\gk[\ua]$, cela donne
une %\linebreak
 \egt $s_1=a_1y_1+\cdots +a_ny_n$ avec $y_1$, \dots, $y_n\in\gk[\ua]$. 
En outre,~$s_1a_i\equiv sa_i\equiv a_i$ modulo $\fN$ pour chaque $i$. 
Puisque $(1-s_1)a_i\in\fN$%et $1-s_1$ est \idmz
, il existe $k_i$ 
tel que $(1-s_1)a_i^{k_i}=0$.
Finalement avec $k=k_1+\cdots +k_n$, on obtient $\fa^k=\gen{s_1}$.
\end{proof}
%-----------------end proof------------------

%:HHH  le lemme a migré chap 4

Rappelons que le lemme \ref{lemZrZr2} établit la réciproque suivante.

\emph{Soient~$\gk\subseteq \gA$, avec $\gA$ entier sur~$\gk$. Si
$\gA$ est un anneau \zedz, alors
$\gk$ est un anneau \zedz.
}

%:HHH (\Sys \zed sur un \cdiz a migre chap 4

%%%%%%%%%%%%%%%%%%%%%%%%%%%%%%%%%%%%%%%%%%%%%%%%%%%%%%%%%%%%%%%%%%%%%%%%%%%
%: subsubsec{Un \nst faible}
\subsubsec{Un \nst faible}

Le \tho suivant, pour l'implication \emph{2} $\Rightarrow$ \emph{3} limitée au cas où $\gA$ est un \cdiz, est souvent appelé \gui{\nst faible} dans la littérature, car il peut servir de préliminaire au \nst (en \clamaz). C'est à distinguer des autres \nsts faibles déjà envisagés dans cet ouvrage.

%:   Theorem{thNst0}------
\begin{theorem}
\label{thNst0} \emph{(Un \nst faible)}
\\
Soit~$\gK$  un anneau \zed réduit et~$\gA$ une \Klg \tfz. 
Pour les \prts suivantes, on a 1 $\Rightarrow$ 2 $\Leftrightarrow$ 3.
\begin{enumerate}
\item $\gA$ est  un anneau local.
\item $\gA$ est  \zedz.
\item $\gA$ est fini sur~$\gK$.
\end{enumerate}
\end{theorem}
%--- end-theorem-----------------------------------------
NB: on ne suppose pas que $\rho:\gK\to\gA$ est injective.
%-----------------begin proof------------------
\begin{proof} On sait déjà que \emph{3} implique \emph{2.} Voyons que \emph{1} ou \emph{2} implique \emph{3}.

On peut remplacer~$\gK$ par $\rho(\gK)$ qui est \egmt
\zed réduit, on a alors~$\gK\subseteq\gA=\Kxn=\gK[\ux]$. Nous faisons une
preuve par \recu sur $n$.
Le cas $n=0$ est trivial. Passons de $n-1$  à $n$. \\
Si~$\gA$ est \zedz, il
existe un \pol $R\in\KXn$ et un entier~$\ell$ tel que $x_n^\ell\big(1-x_nR(\ux)\big)=0$.
Le \pol $X_n^\ell\big(1-X_nR(\uX)\big)$ a l'un de ses \coes égal à 1 et est donc primitif.\\
Si $\gA$ est local, $x_n$ ou $1+x_n$ est \ivz. \Spdg on suppose que $x_n$
est \ivz. Il existe un \pol $R\in\KXn$ tel que $1+x_nR(\ux)=0$.
Le \pol $1+X_nR(\uX)$ a l'un de ses \coes égal à 1 et est donc primitif.
\\
Dans les deux cas, on peut faire un \cdv comme dans le lemme~\ref{lemCDV} (cas d'un \cdi infini) ou~\ref{lemNoether} (cas \gnlz). On a alors~$\gA=\Kyn$, et~$\gA$ est finie sur~$\gA_1=\gK[y_1,\ldots ,y_{n-1}]\subseteq\gA$. \\
Si $\gA$ est \zedz, le lemme~\ref{lemZrZr2}
implique que~$\gA_1$ est \zed et l'on peut donc appliquer l'\hdrz.\\
Si $\gA$ est local, le point \emph{3} du \thref{thJacplc}
implique que~$\gA_1$ est local et l'on peut donc appliquer l'\hdrz.
\end{proof}
%-----------------end proof------------------

%
\rem Ce qui est nouveau 
%:HHH rajout une ligne ci-dessous
pour l'implication \emph{2} $\Rightarrow$ \emph{3} 
dans le \thref{thNst0} par rapport au \thref{thSPolZed},
qui utilise la mise en position de \iNoez,
c'est donc le fait que l'on  ne suppose pas l'\alg \pf mais seulement \tfz.
Les deux  \dems sont en définitive basées sur le lemme~\ref{lemZrZr2}
et sur un lemme de \cdvz.
\eoe

%%%%%%%%%%%%%%%%%%%%%%%%%%%%%%%%%%%%%%%%%%%%%%%%%%%%%%%%%%%%%%%%%%%%%%%%%%%
%: subsubsec{Algebres entières  sur un anneau \qiz}
\penalty-2500
\subsubsec{Algèbres entières  sur un anneau \qiz}

\rdb \label{NOTAReg}
On note $\Reg\gA$ le filtre des \elts \ndzs de l'anneau~$\gA$, de sorte que 
l'anneau total des fractions $\Frac\gA$ est égal à $(\Reg\gA)^{-1}\!\gA$.

%:     fact{factReduitEntierQi}
\begin{fact}\label{factReduitEntierQi} ~\\
Soient~$\gA$ un anneau \qiz, 
$\gK=\Frac \gA$,~$\gL\supseteq\gK$ une \Klg entière réduite et~$\gB$ la \cli de~$\gA$ dans~$\gL$.\\
Alors,~$\gB$  est \qi et  $\Frac\gB=\gL=(\Reg\gA)^{-1}\gB$. 
\end{fact}
%%%%%%%%%%%%%%%%%%%%%%%%%%%%%%%%%%%%%%%%%
%
\begin{proof} $\gK$ est \zedr parce que~$\gA$ est \qi (fait~\ref{factQoQiZed}).
L'anneau~$\gL$ est \zed parce qu'il est entier sur~$\gK$. 
Comme il est réduit, il est \qiz. Comme~$\gB$ est \icl dans 
$\gL$, tout \idm de~$\gL$ est dans~$\gB$, donc~$\gB$ est \qiz.
\\
Considérons un $x\in\gL$ et un \polu  $f\in\KX$ qui annule $x$. En {chassant les \denosz} on obtient un \pol 

\snic{g(X)=a_mX^m+a_{m-1}X^{m-1}+\cdots+a_0\in\AX}

%\sni
qui annule $x$,
avec $a_m\in\Reg\gA$. Alors, $y=a_mx$, entier sur~$\gA$, est  dans~$\gB$ et
$x\in(\Reg\gA)^{-1}\gB$.
\end{proof}

%: subsubsec{Algèbres qui sont des \mpfsz}
\subsubsec{Algèbres qui sont des \mpfsz}

%:      Theorem{propAlgFinPresfin}------
\begin{theorem}
\label{propAlgFinPresfin} \emph{(Quand une \klg est un \kmo \pfz)}
%-----------------begin enum------------------
\begin{enumerate}
\item Pour une \klgz~$\gA$ \propeq
%-----------------begin item------------------
\begin{enumerate}
\item $\gA$ est  un \kmo \pfz.
\item $\gA$ est finie et c'est une \klg \pfz.
\item $\gA$ est entière et \pf sur~$\gk$.
\end{enumerate}
%-----------------end item------------------
\item Si ces conditions sont vérifiées et si en outre~$\gk$ est \coh (resp.
\coh \fdiz), alors~$\gA$ est \coh (resp. \coh \fdiz).
\end{enumerate}
%-----------------end enum------------------
\end{theorem}
%--- end-theorem-----------------------------------------
%-----------------begin proof------------------
\begin{proof}
\emph{1a} $\Rightarrow$ \emph{1b.}
Soit $\gA=\som_{i=1}^{m}b_i\gk$  un \kmo \pfz.
On doit donner une présentation finie de~$\gA$ comme~\klgz. 
On considère le
\sgr $(b_1,\ldots ,b_m)$.
D'une part, on prend les relations~$\gk$-\lins données par la \pn de~$\gA$ comme
\kmoz. D'autre part on écrit chaque $b_ib_j$ comme combinaison~$\gk$-\lin des $b_k$.
Modulo ces dernières relations, tout \pol en les $b_i$ à \coes dans~$\gk$ 
se réécrit comme
une combinaison~$\gk$-\lin des $b_i$. Il s'évalue donc à~0 dans~$\gA$ \ssi 
(en tant que \polz)
il est dans l'idéal engendré par toutes les relations que l'on  a données.
\\
\emph{1b} $\Longleftrightarrow$ \emph{1c.} Clair.
\\
\emph{1b} $\Rightarrow$ \emph{1a.}
On suppose que~$\gA$ est finie sur~$\gk$ avec

\snic{\gA=\kxn
=\aqo{\kuX}{\uf}.}

%\sni
Pour chaque $i$, soit $t_i(X_i)\in\gk[X_i]$ un \pol unitaire tel que $t_i(x_i)=0$, et
$\delta_i=\deg\,t_i$. On a 
$$\preskip.2em \postskip.4em 
\;\;\gA=\aqo{\kuX}{t_1,\ldots ,t_n,h_1,\ldots ,h_s}, 
$$
où les $h_j$ sont les $f_j$ réduits modulo~$\gen{t_1,\ldots,t_n }$.
\\
Les \gui{\momsz} $\ux^\ud=x_1^{d_1}\cdots x_n^{d_n}$ où $d_1<\delta_1,\ldots
,\, d_n<\delta_n$ (ce que nous notons $\ud<\udel$) forment une base de l'\alg
$\aqo{\gk[\uX]}{\ut}$ et
un \sys \gtr $G$ du \kmoz~$\gA$.
Une relation arbitraire entre ces \gtrs est obtenue lorsque l'on  écrit
$\som_{j=1}^sg_j(\ux)h_j(\ux)=0$,
à condition de l'exprimer comme une combinaison~$\gk$-\lin d'\elts de $G$.
On peut naturellement se limiter
aux $g_j$ qui sont de degré $<\delta_i$ en chaque variable $X_i$.
Si l'on fixe un indice $j\in\lrb{1..s}$ et un \mom $\ux^\ud$ avec
$\ud<\udel$,
on obtient une relation de dépendance~$\gk$-\lin entre les \elts de $G$ en
réécrivant
$\uX^\ud\, h_j(\uX)$ modulo~$\gen{t_1,\ldots,t_n }$ et en disant que la \coli des
\elts de~$G$ ainsi obtenue est nulle.
Ces relations engendrent le \kmo des \syzys entre les \elts de~$G$.

\emph{2}. Si~$\gk$ est \coh (resp. \coh \fdiz), alors
on sait que~$\gA$ est \coh (resp. \coh \fdiz) en tant que \kmo (puisqu'il est
\pfz).
Soit $(b_i)_{i=1}^m$ un \sys \gtr de~$\gA$ comme~\hbox{\kmo} et $v=(v_1,\ldots
,v_n)\in\Ae n$.
L'\idz~$\gen{v_1,\ldots ,v_n}$ est le~\hbox{\kmo} \tf engendré par les $v_ib_j$, donc
il est détachable si~$\gk$ est \fdiz.\\
En outre, une~$\gA$-\syzy pour $v$ peut se réécrire comme
une~$\gk$-\syzy entre les $v_ib_j$. Donc un \sys \gtr du \kmo des
$\gk$-\syzys entre les $v_ib_j$ donne  par relecture un \sys \gtr du \Amo 
des~\hbox{$\gA$-syzygies} entre les~$v_i$.

\end{proof}
%-----------------end proof------------------

%:  subsubsec{Algèbre entière  sur un anneau \iclz}
\subsubsec{Algèbre entière  sur un anneau \iclz}

On généralise ici la proposition \ref{propAECDN}.
%:   theorem \label{Thextent}
\begin{theorem}
\label{Thextent} 
Soit $\gA$ un anneau \iclz, $\gK$ son corps de fractions, $\gL$
un surcorps \stf sur $\gK$ et $\gB$ la clôture intégrale de~$\gA$ dans~$\gL$.
Pour $z\in\gL$, notons $\mu_{\gL,z}\in \End_\gK(\gL)$ la multiplication par $z$, puis
 $\nu_z(X)$ et $\chi_z(X)$
le \polmin et le \polcar de~$\mu_{\gL,z}$ (ce sont des \elts de $\gK[X]$). %
\begin{enumerate}
\item \label{i1Thextent} Pour  $z\in\gL$, on a: $z\in\gB\Longleftrightarrow\nu_z\in \AX\Longleftrightarrow\chi_z\in\AX$. En particulier, pour $z\in\gB$,
$\rN_{\gL/\gK}(z)$ et  $\Tr_{\gL/\gK}(z)\in\gA.$  
\end{enumerate} 
On suppose maintenant que $\gL$ est étale sur $\gK$, i.e. que 
$\Disc_{\gL/\gK}\in\gK\eti$.
\begin{enumerate} \setcounter{enumi}{1}
\item \label{i2Thextent}
Soit $x$ un \elt de $\gB$ tel que $\Kx=\gL$. On note  $\Delta_x=\disc(\chi_x)$.
\begin{enumerate}
\item $\gA[x]\simeq \aqo\AX{\chi_x}$, \Amo libre de rang $\dex{\gL:\gK}$.
\item   On a \hbox{$\gA[x][1/\Delta_x]=\gB[1/\Delta_x]$},
 anneau  \iclz. 
\item Si $\gA$ est un  anneau à pgcd,
si  $\Delta_x=d^{2}b$  et $b$
est sans facteur carré, alors $\gA[x][1/d]=\gB[1/d]$
et c'est un anneau est \iclz.
\end{enumerate} 
\item \label{i3Thextent}
Soit $\cB$  une base de $\gL$ sur $\gK$ contenue dans $\gB$ et $M\subseteq \gB$ le \Amo de base~$\cB$.
\begin{enumerate}
\item L'\elt $\Delta=\disc_{\gL/\gK}(\cB)$ est dans $\gA$.
\item Pour tout $x\in\gB$,  $\Delta x\in M$, autrement dit \hbox{$M\subseteq \gB\subseteq \fraC 1\Delta M$}. 
\item Si $\gA$ est un anneau à pgcd, pour tout $x\in\gB$,
il existe $\delta\in\gA$, tel que~$\delta^{2}$ divise $\Delta$ \hbox{et  $\delta x\in M$}.
\\
 Si en outre $\Delta=d^{2}b$ avec $b$ sans facteur carré, $M\subseteq \gB\subseteq \fraC 1 d M$. 
\end{enumerate}
\end{enumerate}
\end{theorem}
%--- end-theorem-----------------------------------------
%
\begin{proof}
\emph{\ref{i1Thextent}.} 
Si $z\in\gB$, il annule un \pol \unt $h(X)\in\AX$, et le \polz~$\nu_z$
divise~$h$ dans $\KX$.
 Comme $\nu_z$ est \unt et $\gA$ \iclz, on obtient
$\nu_z\in\AX$ par le lemme \ref{lem0IntClos}. 

Par ailleurs dans $\KX$, 
$\nu_z$ divise  $\chi_z$ et $\chi_z$
divise une puissance de $\nu_z$, donc,
toujours par le lemme \ref{lem0IntClos},  $\nu_z\in\AX$ équivaut à $\chi_z\in\AX$. 

\emph{\ref{i2Thextent}a.} Clair: $(1,x,\dots,x^{\dex{\gL:\gK}-1})$ est 
à la fois une base de $\gA[x]$ sur $\gA$ et de $\gL$ sur $\gK$.
Notez que par hypothèse $\chi_x=\nu_x$.

\emph{\ref{i2Thextent}b.}  Considérons le cas particulier de \emph{\ref{i3Thextent}b} où $M=\gA[x]$. 
\\
 On obtient $\gB[1/\Delta_x]=\gA[x][1/\Delta_x]$, et  comme $\gB$ est \iclz, il en va de même pour
$\gB[1/\Delta_x]$. 

\emph{\ref{i2Thextent}c.}  Cas particulier de \emph{\ref{i3Thextent}c} avec $M=\gA[x]$, en raisonnant 
comme en~\emph{3c.} 

\emph{\ref{i3Thextent}a.} Conséquence \imde de \emph{\ref{i1Thextent}.}

\emph{\ref{i3Thextent}b.} \'Ecrivons $\cB=(\bn)$ et $x=\sum_{i}x_ib_i$ avec les $x_i\in\gK$.
Considérons par exemple le \coe $x_1$, supposé non nul.
Le $n$-uplet  $\cB'=(x,b_2,\dots,b_n)$ est une $\gK$-base de $\gL$ contenue dans $\gB$.
La matrice de $\cB'$ sur $\cB$ a pour \deter $x_1$. 
Donc  $x_1^{2}\Delta=x_1^{2}\disc(\cB)=\disc(\cB')\in\gA$. A fortiori $(x_1\Delta)^{2}\in\gA$, et puisque $\gA$ est \iclz, $x_1\Delta\in\gA$.
Ainsi  toutes les \coos sur $\cB$ de $\Delta x$ sont dans $\gA$.

\emph{\ref{i3Thextent}c.} Lorsque $\gA$ est un anneau à pgcd, on écrit l'\elt $x_1$ comme une fraction réduite $x_1=a_1/\delta_1$. Alors, puisque $x_1^{2}\Delta\in\gA$, $\delta_1^{2}$ divise
$a_1^{2}\Delta$, et comme $\pgcd(a_1,\delta_1)=1$, l'\elt $\delta_1^{2}$ divise
$\Delta$. On fait de même pour chacun des $x_i=a_i/\delta_i$. Si $\delta$ est le ppcm des $\delta_i$, $\delta^{2}$ est le ppcm des $\delta_i^2$, donc
il divise $\Delta$, et $\delta x\in M$. 
\end{proof}
%
%\newpage
%%%%%%%%%%%%%%%%%%%%%%%%%%%%%%%%%%%%%%%%%%%%%%%%%%%%%%%%%%%%%%%%%%%%%%%%%%%
\penalty-2500
\section{Algèbres \stfesz}
\label{subsecAlgStfes}

%\vspace{4pt}
%%%%%%%%%%%%%%%%%%%%%%%%%%%%%%%%%%%%%%%%%%%%%%%%%%%%%%%%%%%%%%%%%%%%%%%%%%%
%:  subsec{Le module dual et la trace}
\subsec{Le module dual et la trace}

Si $P$ et $Q$ sont des \kmos \ptfs on a un \iso canonique $\theta_{P,Q}:P\sta\te_\gk Q\to \Lin_\gk(P,Q)$.

Lorsque le contexte est clair on peut identifier $\alpha\te x\in P\sta\te_\gk Q$ avec l'\kli correspondante $y\mapsto \alpha(y)x$. 

En particulier, un \syc de $P$, $\big((\xn),(\aln) \big)$, est caractérisé par l'\egtz:
\begin{equation}
\label{eq1syc}
\som_{i=1}^n\alpha_i\te x_i=\Id_P.
\end{equation}
De manière duale on a, modulo l'identification de $P$ avec $(P\sta)\sta$:  
\begin{equation}
\label{eq2syc}
\som_{i=1}^nx_i\te \alpha_i=\Id_{P\sta}.
\end{equation}
Cette \eqn signifie que pour tout~$\gamma\in P\sta$ on a~$\gamma=\sum_{i=1}^n\gamma(x_i)\alpha_i$.

%:     Definition{definotaAsta}
\begin{definota}\label{definotaAsta}
Soit~$\gA$ une \klgz. 
\\
Le dual $\Asta$  du
\kmoz~$\gA$ est muni d'une structure de \Amo via
la loi externe $(a,\alpha)\mapsto a\centerdot\alpha\eqdefi \alpha\circ \mu_a$, i.e.
$(a\centerdot\alpha)(x) = \alpha(ax)$.  
 
\end{definota}

%:HHH 11 petite modif ci dessous
Les faits \ref{factMatriceAlin} et/ou~\ref{factMatriceEndo} donnent le résultat suivant.

%:     Fact{factTraceSyc}
\begin{fact}\label{factTraceSyc}
Soit $\big((\xn),(\aln) \big)$  un \syc pour la \klg \stfez~$\gA$, alors 
 l'\kli $\mu_{\gA,a}$ est représentée dans ce \sys par la matrice 
${\big(\alpha_i(ax_j)\big)_{i,j\in\lrbn}}$
et l'on a
\begin{equation}
\label{eqfactTraceSyc}
\Tr\iAk=\som_{i=1}^nx_i\centerdot\alpha_i,\quad \big(\mathit{i.e.,}\, \forall a\in\gA,\;\Tr\iAk(a)=\som_{i=1}^n\alpha_i(ax_i)\big).
\end{equation}
\end{fact}

%%%%%%%%%%%%%%%%%%%%%%%%%%%%%%%%%%%%%%%%%%%%%%%%%%%%%%%%%%%%%%%%%%%%%%%%%%%
%:  subsec{Norme et \elt cotransposé}
\subsec{Norme et \elt cotransposé}

Nous introduisons  la notion d'\emph{\elt cotransposé}
dans une \alg \stfez. Il suffit de reprendre ce qui a été dit
dans le cas d'une \alg libre de rang fini \paref{eqelt0cotransp}.
Si~$\gA$ est \stfe sur~$\gk$ on peut identifier~$\gA$ à une sous-\klg
commutative de $\End_\gk(A)$, où $A$ désigne le \kmoz~$\gA$ privé de sa structure multiplicative, au moyen de
l'\homo de multiplication $x\mapsto\mu_{A,x}=\mu_x$.
Alors, puisque $\wi \mu_x=G(\mu_x)$ pour un \pol $G$ de $\kT$
(point \emph{8} du \thref{propdef det ptf}), on peut définir
$\wi x$ par l'\egt $\wi x=G(x)$, ou ce qui revient au même $\wi {\mu_x}=\mu_{\wi x}$. Si plus de précision est \ncr
on utilisera la notation $\Adj\iAk (x)$.
Cet \elt $\wi x$ s'appelle \emph{l'\elt cotransposé de x}.
L'\egt $\wi {\mu_x}\,\mu_x=\det(\mu_x)\Id_\gA$ donne alors:
\index{cotransposé!element@\elt --- (dans une \alg strictement finie)}
%%%%%%%%%%%%%%%%%%%%%%%%%%%%%%%%%%%%%%%%%
\begin{equation}
\label{eqeltcotransp}
x\ \Adj\iAk (x)=\rN\iAk (x)
\end{equation}

%:     Lemma{lemIRAdu}--------------
\begin{lemma}
\label{lemIRAdu} Soit~$\gk\vers{\rho}\gA$ une \alg \stfez, $x\in\gA$ et $y\in\gk$.
%-----------------begin enum------------------
\begin{enumerate}
\item On a $x\in\Ati$ \ssi $\rN\iAk (x)\in\Ati$.
\\
Dans ce cas $x^{-1}=\wi x/\rN\iAk (x)$.
\item $x$ est \ndz dans~$\gA$ \ssi $\rN\iAk (x)$ est \ndz dans
$\gk$. Dans ce cas $\wi x$ est \egmt \ndzz.
\item $\rho(\gk)$ est facteur direct dans~$\gA$.
\end{enumerate}
Notons $e=\ide_0(\gA)$ (de sorte que~$\gen{e}_\gk=\Ann_\gk(\gA)$).
\begin{enumerate}\setcounter{enumi}{3}
\item On a $\rho(y)\in\Ati$ \ssi $y\in(\aqo{\gk}{e})\eti$.
\item $\rho(y)$ est \ndz dans~$\gA$ \ssi $y$ est \ndz dans $\aqo{\gk}{e}$.
\end{enumerate}
%-----------------end enum------------------
\end{lemma}
%--- end-lemma-----------------------------------------
NB. Si~$\gA$ est un \kmo fidèle, i.e., si $\rho$ est injectif,
on identifie~$\gk$ à $\rho(\gk)$.
Alors,~$\gk$ est facteur direct dans~$\gA$, et un \elt $y$ de~$\gk$ est \iv (resp. \ndzz) dans~$\gk$ \ssi il est \iv (resp. \ndzz) dans~$\gA$.
%-----------------begin proof------------------
\begin{proof}
\emph{1.} Dans un \mptf un \endo (ici $\mu_x$) est
bijectif \ssi son \deter est \ivz.

\emph{2.} Dans un \mptf un \endo est injectif \ssi son \deter est \ndzz.

Les points \emph{3}, \emph{4} et \emph{5} peuvent être démontrés
après \lon en des \eco de~$\gk$. D'après le \tho de structure locale
on est ramené au cas où~$\gA$ est libre de rang fini, disons $k$.
Si $k=0$, alors~$e=1$, donc~$\gA$ et~$\aqo{\gk}{e}$ sont triviaux et tout est clair (même si c'est un peu troublant). Examinons le cas où $k\geq1$, d'où $e=0$
et identifions~$\gk$ à $\rho(\gk)$.
\\
Les points \emph{4} et \emph{5} découlent alors des points \emph{1} et \emph{2}
parce que $\rN\iAk (y)=y^{k}$.
\\
Pour le point \emph{3}, on considère une base $(b_1,\ldots,b_k)$ de~$\gA$
sur~$\gk$ et des \elts $a_1$, \ldots, $a_k\in\gk$ tels que $\sum_ia_ib_i=1$.
On a $\rN\iAk (\sum_ia_ib_i)=1$. 
Par ailleurs,
pour~$y_1$, \dots, $y_k\in\gk$, $\rN\iAk (\sum_i y_ib_i)$ s'écrit
comme un \pog de degré $k$ dans $\kuy$
(voir la remarque \paref{factNormeRationnelle}), et donc

\snic{\rN\iAk (\sum_ia_ib_i)=\sum_ia_i\beta_i=1}

%\sni
pour des $\beta_i\in\gk$
convenables. On considère l'\elt $\beta\in\End_\gk(\gA)$ défini par $\beta(\sum_ix_i{b_i})=\sum_ix_i\beta_i$. \\
Alors, $\beta(1)=1$, donc $\beta(z)=z$ pour $z\in\gk$,
$\Im \beta=\gk$ et $\beta\circ \beta=\beta$.
\end{proof}
%-----------------end proof------------------

%%%%%%%%%%%%%%%%%%%%%%%%%%%%%%%%%%%%%%%%%%%%%%%%%%%%%%%%%%%%%%%%%%%%%%%%%%%
%:  subsec{Transitivité et rang}
\subsec{Transitivité et rang}\label{subsecTransRang}
\vspace{3pt}
Lorsque~$\gA$ est de rang constant $n$, nous noterons $\dex{\gA:\gk}=n$. Ceci généralise
la notation déjà définie dans le cas des \algs libres, et cela sera \gne
au chapitre~\ref{chap ptf1} (notation~\ref{notaTraceDetCarAlg}).
Dans ce paragraphe, $m$ et $n$ désignent des entiers.

%:     Fact{factTransptf}
\begin{fact}\label{factTransptf}
%:HHH    $\gB$ une \Alg \stfe  mis au depart
Soit~$\gA$ une \klg \stfez, $M$  un \Amo \ptf et $\gB$ une \Alg \stfez.
\begin{enumerate}
\item  $M$ est aussi un  \kmo \ptfz.
\item  Supposons $\rg_\gA M=m$  et notons $f(T)=\rR\gk(\gA)\in\BB(\gk)[T]$ le \polmu de~$\gA$ comme \kmoz,
alors $\rR\gk(M)=f^m(T)=f(T^m)$.
\item $\gB$ est \stfe sur~$\gk$ et $\Tr\iBk=\Tr\iAk\circ \Tr\iBA$.
\end{enumerate}
\end{fact}
\begin{proof} \emph{1.} Supposons que~$\gA\oplus E\simeq \gk^r$
(\kmosz) et $M\oplus N\simeq \Ae s$ (\Amosz).
Alors~$M\oplus N\oplus E^s\simeq \gk^{rs}$ (\kmosz).
 On peut redire ceci avec des \sycs sous la forme suivante: 
 si $\big((\xn),(\aln) \big)$ est un \syc pour le \kmoz~$\gA$ et $\big((\ym),(\beta_1,\ldots,\beta_m) \big)$ un \syc pour le \Amo $M$, alors $\big((x_iy_j),(\alpha_i\circ \beta_j) \big)$ est un \syc pour le \kmo $M$. 
\\
\emph{2.} Laissé \alec (qui peut s'appuyer sur la description précédente du \sycz, ou consulter la \dem du lemme~\ref{lem1TransPtf}).
\\
\emph{3.} On travaille avec des \sycs comme dans le point \emph{1} et l'on applique le fait~\ref{factTraceSyc} concernant la trace.
\end{proof}
%

%:     Theorem{propTraptf}
\begin{theorem}\label{propTraptf}
Soient~$\gk\subseteq\gA\subseteq\gB$ des anneaux.
Supposons que~$\gB$ est  \stf sur~$\gA$.
\begin{enumerate}
\item L'anneau $\gB$ est \stf sur~$\gk$ \ssiz$\gA$ est \stf sur~$\gk$.\perso{Lorsque~$\gA$ et~$\gB$ sont \stes sur~$\gk$  
alors~$\gB$ est \stfe sur~$\gA$? C'est en tout cas ce que semblent dire
Ferrero et Paques.}
\item Si $\dex{\gA:\gk}=n$  et
$\dex{\gB:\gA}=m$, alors $\dex{\gB:\gk}=mn$.
\item Si $\dex{\gB:\gk}=mn$ et
$\dex{\gB:\gA}=m$, alors $\dex{\gA:\gk}=n$.
\end{enumerate}
\end{theorem}
\begin{proof}
\emph{1.} Si~$\gB$ est \stf sur~$\gk$, alors
$\gA$ est \stf sur~$\gk$: cela résulte de ce que~$\gA$ est facteur direct dans
$\gB$ (lemme~\ref{lemIRAdu} point \emph{3}), qui est un~$\gk$-\mptfz.
\\
L'implication réciproque est dans le lemme~\ref{factTransptf}.
 
\emph{2} et \emph{3.} Résultent du fait \ref{factTransptf} en notant que le
seul \pol multiplicatif~$f$ de $\kT$ qui vérifie $f^m(T)=T^{mn}$ est $f=T^n$ puisque
$f^m(T)=f(T^{m})$.
\end{proof}

\rem Des formules de transitivité plus \gnles (dans le cas de rang non constant)
sont données en section~\ref{secAppliLocPtf} dans le paragraphe 
\gui{\nameref{subsecTransPtf}} \paref{subsecTransPtf} (voir notamment le corolaire \ref{lem2TransPtf} et le \thoz~\ref{corthTransPtf}).
\eoe

%\newpage
%%%%%%%%%%%%%%%%%%%%%%%%%%%%%%%%%%%%%%%%%%%%%%%%%%%%%%%%%%%%%%%%%%%%%%%%%%%
%  subsec{Formes \lins dualisantes, \algs \stes}
%\penalty-2500
\section{Formes \lins dualisantes, \algs \stesz}\label{secAlgSte}

%:      defdualisante
\begin{definition}\label{defdualisante}
\emph{(Forme bi\lin \smq non \dgnez, forme \lin dualisante, \asez)}
\\
Soit $M$ un \kmo et~$\gA$ une \klgz.

\begin{enumerate}

\item
Si $\phi : M \times M \to \gk$ est une forme bi\lin \smqz, on lui
associe l'\kli $\varphi : M \to M\sta$ définie par $\varphi(x) = \phi(x, \bullet) =
\phi(\bullet, x)$. \\
On dit que $\phi$ est \emph{non \dgnez} si $\varphi$
est un \isoz.
\index{non dege@non \dgnez!forme bilinéaire}
\index{forme bilinéaire!non \dgnez}

\item
Si $\lambda \in \Lin_\gk(\gA,\gk)=\Asta$, on lui associe la forme~$\gk$-bi\lin
\smq sur~$\gA$, notée $\Phi_{\gA/\gk,\lambda}=\Phi_\lambda$ et définie par
$\Phi_\lambda(x,y) = \lambda(xy)$.
\\
On dit que la  forme \lin $\lambda$ est \ix{dualisante} si
$\Phi_\lambda$ est non \dgnez. 
\\
On appelle \emph{algèbre de Frobenius}
une \alg pour laquelle il existe une forme \lin dualisante.%
\index{forme linéaire!dualisante}%
\index{algèbre!de Frobenius}
\index{Frobenius!algèbre de ---}

\item
Si~$\gA$ est \stfe sur~$\gk$ on appelle \ix{forme trace} (ou \emph{forme bi\lin tracique}) la forme $\Phi_{\Tr\iAk}$.

\item L'\alg $\gA$ est dite \emph{\stez} sur~$\gk$ si elle est \stfe et si
la trace est dualisante, i.e. la \ftr est non \dgnez.
\index{strictement etale@strictement étale!algèbre --- }
\index{algèbre!strictement etale@strictement étale}

\end {enumerate}
\end {definition}

\rem Si~$\gA$ est libre de base $(\ue)=(e_1,\ldots,e_n)$ sur~$\gk$,
la matrice de $\phi$ et celle de $\varphi$ coïncident (pour les bases convenables). En outre, $\phi$ est non \dgne \ssi $\Disc\iAk =\disc\iAk (\ue)$
est \ivz.
On note que lorsque~$\gk$ est un \cdi on retrouve la \dfn
\ref{defi1Etale} pour une \alg étale\footnote{Nous n'avons pas donné
la \dfn \gnle d'une \alg étale. Il se trouve que les \algs étales sur les \cdis sont toujours \stes (au moins en \clamaz, c'est en rapport avec le \thref{thSepProjFi}),
mais que ce n'est plus le cas pour un anneau commutatif arbitraire,
d'où la nécessité d'introduire ici la terminologie \gui{\stez}.}.
\eoe

%%%%%%%%%%%%%%%%%%%%%%%%%%%%%%%%%%%%%%%%%%%%%%%%%%%%%%%%%%%%%%%%%%%%%%%%%%%
%:  subsec{Formes dualisantes}
\subsec{Formes dualisantes}

%:    theorem  \label{factCarDua}
\begin {theorem}\label{factCarDua} \emph{(\Carn des formes dualisantes dans le cas \stfz)} \\
Soit~$\gA$ une \klg et $\lambda\in\Asta$.
 Pour $x\in\gA$, notons  
$x\sta=x\centerdot\lambda \in\Asta.$
\begin {enumerate}
\item 
Si~$\gA$ est  \stfe et si $\lambda$ est dualisante, alors pour tout \sgr $(x_i)_{i\in\lrbn}$, il existe
un \sys $(y_i)_{i\in\lrbn}$ tel que l'on ait 
\begin{equation}\ndsp\preskip.4em \postskip.4em
\label{eqDua}
\som_{i=1}^ny_i\sta\te x_i=\Id_\gA, \quad i.e. \quad\forall x\in\gA,\;
x = \som_{i=1}^n \lambda(xy_i)x_i
\end{equation}
En outre, si~$\gA$ est fidèle, $\lambda$ est surjective.

\item
Réciproquement, s'il existe deux \syss $(x_i)_{i\in\lrbn}$, $(y_i)_{i\in\lrbn}$ tels que  $\sum_iy_i\sta\te x_i=\Id_\gA$, alors:
\begin{itemize}
\item $\gA$ est  \stfez, 
\item la forme $\lambda$ est dualisante,  
\item   et l'on a l'\egt $\sum_ix_i\sta\te y_i=\Id_\gA$.
\end{itemize}

\item Si~$\gA$ est \stfez, \propeq
\begin{enumerate}
\item  $\lambda$ est dualisante.
\item  $\lambda$ est une base du \Amo $\Asta$ 
(qui est donc libre de rang $1$). 
\item  $\lambda$ engendre le \Amo $\Asta$, i.e.~$\gA\centerdot\lambda=\Asta$. 
\end{enumerate}

\end {enumerate}
\end {theorem}
%%%%%%%%%%%%%%%%%%%%%%%%%%%%%%%%%%%%%%%%%
\begin {proof}
\emph {1.}
D'une part $y \mapsto y\sta$ est un \iso
de~$\gA$ sur $\Asta$, et d'autre part
 tout \sgr est la première composante d'un \sycz.
Voyons la surjectivité. Comme $\gA$ est fidèle on peut supposer que $\gk\subseteq \gA$. Soit $\fa$ l'\id de~$\gk$
engendré par les $\lambda(y_i)$. L'\egrf{eqDua} donne
l'appartenance $1=\sum_i\lambda(y_i)x_i\in\fa\gA$. Comme~$\gA$ est entière sur~$\gk$, le lying over (lemme~\ref{lemLingOver})
montre que $1\in\fa$.

 \emph {2.} L'\egrf{eqDua} donne %par dualité 
 $\alpha = \sum_i \alpha(x_i) y_i\sta$ pour $\alpha\in\Asta$. Ceci prouve que $y \mapsto y\sta$ est surjective. Par
ailleurs, si $x\sta=0$, alors  \hbox{on a $\lambda(xy_i) = 0$}, \hbox{puis $x = 0$}. %\linebreak  
Ainsi $\lambda$ est dualisante. 
\\
Enfin,  l'\egt $\alpha = \sum_i \alpha(x_i) y_i\sta$ donne avec $\alpha = x\sta$:  $x\sta = \sum_i \lambda(x_i x) y_i\sta$. \linebreak 
Et puisque $z\mapsto z\sta$ est un~$\gk$-\isoz, 
 $x= \sum_i \lambda(x_i x)y_i$.

%\penalty-2500 
 \emph {3a} $\Leftrightarrow$ \emph{3b.} \gui{$\lambda$ est dualisante} signifie que $x\mapsto x\sta$ est un \isoz, i.e. que
$\lambda$ est une~$\gA$-base de $\Asta$. 
L'implication \emph {c} $\Rightarrow$ \emph{a} résulte du point~\emph{2}
car un \syc s'écrit $\big((x_i),(y_i\sta)\big)$.
\end {proof}

\exls Voir les exercices~\ref{exoGroupAlgebra} à 
\ref{exoFrobeniusAlgExemples} et le \pbz~\ref{exoBuildingFrobAlgebra}.

1) Si $f\in\kX$ est \monz, l'\algz~$\gk[x]=\aqo{\gk[X]}{f(X)}$ est une \alg de Frobenius
(exercice~\ref{exo1Frobenius}).

2) L'\algz~$\gk[x,y]=\aqo{\gk[X,Y]}{X^2,Y^2,XY}$ n'est pas une \alg de Frobenius
(exercice~\ref{exoFrobeniusAlgExemples}).
\eoe

\penalty-2500
\subsubsection*{Extension des scalaires}

%:     Fact{factEdsDualisante}
\begin{fact}\label{factEdsDualisante}\emph{(Stabilité des formes dualisantes
par \edsz)}\\
On considère deux \klgsz~$\gk'$ et~$\gA$
et l'on note~$\gA'=\gk'\te_\gk\gA$.\\ 
Si la forme $\alpha\in\Lin_\gk(\gA,\gk)$
est dualisante, il en va de même pour la \linebreak forme 
$\alpha'\in\Lin_{\gk'}(\gA',\gk')$
obtenue par \edsz.
\\
Comme conséquence, une \alg de Frobenius reste une \alg de Frobenius
par \edsz. 
\end{fact}

%\penalty-2500
\vspace{4pt}
\subsubsection*{Transitivité pour les formes dualisantes}
\vspace{-2pt}

%:     Fact{factDuaTrans}
\begin{fact}\label{factDuaTrans}
Soient~$\gA$ une \klg \stfez, $\gB$ une \Alg \stfez, 
 $\beta\in\Lin_\gA(\gB,\gA)$  et $\alpha\in\Lin_\gk(\gA,\gk)$.
\begin{enumerate}
\item Si $\alpha$ et $\beta$ sont dualisantes, il en va de même pour $\alpha\circ \beta$.
\item Si $\alpha\circ \beta$ est  dualisante et $\beta$ surjective (par exemple
$\gB$ est fidèle et $\beta$ est dualisante), alors $\alpha$ est dualisante.
%
%\item 
%
\end{enumerate}
\end{fact}
\begin{proof}  Si $\big((a_i),(\alpha_i) \big)$ est un \syc de~$\gA\sur\gk$
et $\big((b_j),(\beta_j) \big)$ un \syc de~$\gB\sur\gA$,
alors $\big((a_ib_j),(\alpha_i\circ\beta_j) \big)$ est un \syc 
de~$\gB\sur\gk$.

 \emph{1.}
Pour $a \in \gA$, $b \in \gB$, $\eta\in\Lin_\gk(\gA,\gk)$
et $\epsilon\in\Lin_\gA(\gB,\gA)$
 on vérifie facilement que $ab\centerdot(\eta \circ
\epsilon) = (a\centerdot\eta) \circ (b\centerdot\epsilon)$. 
\\
Puisque $\alpha$ est dualisante, on a des $u_i\in\gA$ tels que $u_i\centerdot\alpha=\alpha_i$ pour  $i\in\lrbn$. 
 Puisque $\beta$ est dualisante, on a des $v_j\in\gB$ tels que $v_j\centerdot\beta=\beta_j$ pour  $j\in\lrbm$. 
Alors, $u_iv_j\centerdot (\alpha \circ\beta)=\alpha_i\circ\beta_j$,
et ceci montre que $\alpha \circ\beta$ est dualisante.

 \emph{2.} Soit $\alpha'\in\Lin_\gk(\gA,\gk)$ que l'on cherche à l'écrire sous la forme
$a\centerdot\alpha$.  Remarquons que pour tout $b_0\in\gB$, on a $\big(b_0\centerdot(\alpha'\circ\beta)\big)\frt\gA =
\beta(b_0)\centerdot\alpha'$; en particulier, si $\beta(b_0)=1$, alors  $\big(b_0\centerdot(\alpha'\circ\beta)\big)\frt\gA =\alpha'$.  Puisque $\alpha\circ
\beta$ est dualisante, il existe $b\in\gB$ tel que $\alpha'\circ \beta =
b\centerdot(\alpha\circ \beta)$.  En multipliant cette \egt par $b_0
\centerdot$, on obtient, par restriction à~$\gA$, $\alpha' = \big((b_0b)
\centerdot (\alpha\circ\beta) \big)\frt\gA = \beta(b_0b) \centerdot \alpha$.
\end{proof}
%

%%%%%%%%%%%%%%%%%%%%%%%%%%%%%%%%%%%%%%%%%%%%%%%%%%%%%%%%%%%%%%%%%%%%%%%%%%%
%:  subsec{Algèbres \stesz}
%\penalty-2500
\subsec{Algèbres \stesz}

Le \tho suivant est un corolaire \imd du  \thref{factCarDua}.
%:    fact  \label{factCarAste}
\begin {theorem}\label{factCarAste} \emph{(\Carn des \asesz)} 
Soit~$\gA$ une \klg \stfez. Pour $x\in\gA$, on note $x\sta=x\centerdot\Tr\iAk\in\Asta$. 
\begin {enumerate}
\item 
Si~$\gA$ est \stez, alors pour tout \sgr  $(x_i)_{i\in\lrbn}$, il existe
un \sys $(y_i)_{i\in\lrbn}$ tel que l'on ait 
\begin{equation}\preskip.4em \postskip.4em
\label{eqaste}
\som_{i=1}^ny_i\sta\te x_i=\Id_\gA, \quad i.e. \quad\forall x\in\gA,\;
x = \som_{i=1}^n \Tr\iAk(xy_i)x_i
\end{equation}
Un tel couple $\big((x_i),(y_i) \big)$  est appelé un \emph{\stycz}.
\\
Si en outre~$\gA$ est fidèle, $\Tr\iAk$ est surjective.
\index{systeme trac@système tracique de \coosz}

\item
Réciproquement, si l'on a un couple $\big((x_i)_{i\in\lrbn}, (y_i)_{i\in\lrbn}\big)$
qui vérifie~\pref{eqaste}, alors~$\gA$ est  \stez,  
et l'on a  $\sum_ix_i\sta\te y_i=\Id_\gA$.

\item \Propeq
\begin{enumerate}
\item  $\Tr\iAk$ est dualisante (i.e. $\gA$ est \stez).
\item  $\Tr\iAk$ est une base du \Amo $\Asta$ 
(qui est donc libre de rang~$1$). 
\item  $\Tr\iAk$ engendre le \Amo $\Asta$. 
\end {enumerate}
\end {enumerate}
\end {theorem}

%%%%%%%%%%%%%%%%%%%%%%%%%%%%%%%%%%%%%%%%%%%%%%%%%%%%%%%%%%%%%%%%%%%%%%%%%%%
\subsubsection*{Extension des scalaires}

Le fait qui suit prolonge les faits \ref{factEdsAlg} et \ref{factEdsDualisante}.
%:     Fact{factEdsEtale}
\begin{fact}\label{factEdsEtale}
On considère deux \klgsz~$\gk'$ et~$\gA$
et l'on note~$\gA'=\gk'\te_\gk\gA$. %
\begin{enumerate}
\item Si~$\gA$ est  \ste sur~$\gk$,
alors~$\gA'$ est  \ste sur~$\gk'$.
\item Si~$\gk'$ est \stfe et contient~$\gk$, et si~$\gA'$
est \ste sur~$\gk'$, alors~$\gA$ est \ste sur~$\gk$.
\end{enumerate}
\end{fact}
\begin{proof}
\emph{1.} Laissé \alecz.
\\
\emph{2.} Supposons d'abord~$\gA$ libre sur~$\gk$. Soit $\Delta=\Disc\iAk =\disc\iAk (\ue)\in\gk$ pour une base $\ue$ de~$\gA$ sur~$\gk$. Par \eds on obtient 
l'\egtz~$\Delta=\Disc_{\gA'/\gk'}\in\gk'$. Si $\Delta$ est \iv dans~$\gk'$
il est \iv dans~$\gk$ d'après le lemme~\ref{lemIRAdu}.
Dans le cas \gnl on se ramène au cas précédent par \lon en des \eco
de~$\gk$.\perso{la \dem fonctionne aussi bien avec~$\gk'$ \fpte sur~$\gk$}
\end{proof}
%

%%%%%%%%%%%%%%%%%%%%%%%%%%%%%%%%%%%%%%%%%%%%%%%%%%%%%%%%%%%%%%%%%%%%%%%%%%%
\subsubsection*{Transitivité pour les \asesz}

%:     Fact{propTrEta}
\begin{fact}\label{propTrEta} 
Soit~$\gA$ une \klg \stfe et~$\gB$ une \Alg \stez. 
\begin{enumerate}
\item Si~$\gA$ est \ste sur~$\gk$, alors~$\gB$ est \ste sur~$\gk$.
\item Si~$\gB$ est \ste sur~$\gk$ et fidèle sur~$\gA$, alors
$\gA$ est \ste sur~$\gk$. 
%
%\item 
%
\end{enumerate}
\end{fact}

\begin{proof}
Résulte de \ref{factDuaTrans} et \ref{factTransptf}.
\end{proof}

\subsubsection*{Séparabilité et nilpotence}

%:     theorem{prop2EtaleReduit}
\begin{theorem}\label{prop2EtaleReduit} \label{propEtaleReduit}
Soit~$\gA$ une \klg \stez.
\begin{enumerate}
\item Si~$\gk$ est réduit,~$\gA$ \egmtz.
\item L'\id $\DA(0)$ est engendré par l'image de ${\rD_\gk(0)}$ dans~$\gA$.
\item Si~$\gk'$ est une \klg réduite, ${\gA'} =\gk'\otimes_\gk \gA$
est réduite.
\perso{ce serait super d'avoir la réciproque, si~$\gA$ est \stfe et si pour tout~$\gk'$ réduit \ldots\ldots alors \ldots \stez}

\end{enumerate}
\end{theorem}

\begin{proof}
\emph{1.} On raisonne à peu près comme dans le cas où~$\gk$ est un \cdi
(fait~\ref{fact1Etale}). On suppose d'abord que~$\gA$ est libre
sur~$\gk$.
 Soit $a\in\DA(0)$. Pour tout $x\in\gA$ la multiplication par
$ax$ est un \endo nilpotent $\mu_{ax}$ de~$\gA$. Sa matrice est nilpotente
donc les \coes de son \polcar sont nilpotents (voir par exemple l'exercice~\ref{exoNilpotentChap2}), donc nuls puisque~$\gk$ est réduit.
En particulier, $\Tr\iAk(ax)=0$. Ainsi $a$ est dans le noyau de l'\kli $tr:a\mapsto\big(x\mapsto \Tr\iAk(ax)\big)$.
Or $tr$ est un \iso par hypothèse
donc $a=0$.
\\
Dans le cas \gnl on se ramène au cas où~$\gA$ est libre sur~$\gk$ par le \tho de structure locale
des \mptfs (en tenant compte du fait~\ref{factEdsEtale}~\emph{1}).

 Le  point \emph{3} résulte de \emph{1}
et du fait~\ref{factEdsEtale}~\emph{1.}
Le  point \emph{2} résulte de  \emph{3}, en considérant~$\gk'=\gk\red$.
\end{proof}

La même technique prouve  le lemme suivant.
%:     Lemma{lemNilpNilp}
\begin{lemma}\label{lemNilpNilp}
Si~$\gA$ est \stfe sur~$\gk$ et si $a\in\gA$ est nilpotent, les \coes de
$\rF{\gA/\gk}(a)(T)$ sont nilpotents (hormis le \coe constant).
\end{lemma}

%:  subsec{Produits tensoriels}
\subsec{Produits tensoriels}

Si $\phi$ et $\phi'$ sont deux formes bi\lins
\smqs sur $M$ et $M'$, on définit une forme
bi\lin \smq sur $M\otimes_\gk M'$, notée $\phi\otimes\phi'$ par:
$$
(\phi\otimes\phi')(x\otimes x', y\otimes y') = \phi(x,y)\phi'(x',y')
\label{NOTAptfb}
.$$

%:     Proposition{propProTNdg}----------
\begin{proposition}
\label{propProTNdg}
 \emph{(Produit tensoriel de deux formes non \dgnesz)}
\\
Soient $M$, $M'$ deux \kmos \ptfs et~$\gA$,~${\gA'} $ deux
\klgs \stfesz.
\begin{enumerate}
\item
Si $\phi$ sur $M$ et $\phi'$ sur $M'$ sont deux formes bi\lins
\smqs non \dgnesz, il en est de m\^ eme de $\phi \otimes
\phi'$.
\item
Si $\lambda \in \Asta$ et $\lambda' \in {\gA'}\sta$ sont deux formes
$\gk$-\lins dualisantes, il en est de même de
$\lambda\otimes\lambda' \in (\gA\otimes_\gk\gA')\sta$.
\end {enumerate}
\end {proposition}
\begin{proof}
\emph {1.}
L'\kli canonique $M\sta \otimes_\gk {M'}\sta \to
(M\otimes_\gk M')\sta$ est un \iso puisque $M, M'$ sont \ptfsz.
Soit $\varphi : M \to M\sta$ l'\iso associé à $\phi$, et
$\varphi' : M' \to {M'}\sta$ celui associé à $\phi'$. Le morphisme
associé à $\phi\otimes\phi'$ est composé de deux \isos donc est
un \iso:
$$
\xymatrix {
M\otimes_\gk M' \ar@{->}[rr]\ar@{->}[dr]_{\varphi\otimes\varphi'}
 && (M \otimes_\gk M')\sta
\\
 & M\sta \otimes_\gk {M'}\sta \ar@{->}[ur]_{\rm iso.\ can.}
\\
}
$$
\emph {2.}
Résulte de $\Phi_{\lambda\otimes\lambda'} =
\Phi_{\lambda}\otimes\Phi_{\lambda'}$.
\end{proof}

La proposition précédente et le lemme~\ref{lemTraceProT} donnent le résultat suivant.
%:     Corollary{corlemTraceProT}
\begin{corollary}\label{corlemTraceProT}
Soient~$\gA$ et~$\gC$ deux \klgs \stfesz.
Alors:

\snic {\Phi_{\Tr_{(\gA\otimes_\gk\gC)/\gk}}=\Phi_{\Tr\iAk}\otimes \Phi_{\Tr_{\gC/\!\gk}}.}

%\sni
En particulier,~$\gA\otimes_\gk\gC$ est \ste si~$\gA$ et~$\gC$
sont \stesz. 
(Pour  le calcul précis du \discriz, voir 
l'exercice~\ref{exoTensorielDiscriminant}.)
\end{corollary}

%%%%%%%%%%%%%%%%%%%%%%%%%%%%%%%%%%%%%%%%%%%%%%%%%%%%%%%%%%%%%%%%%%%%%%%%%%%
%:  subsec {Idempotents, \dinz}
\subsec {\'Eléments entiers, \idmsz, \dinz} \label{secEtaleIdmDin}

%-% ENTRE NOUS
\entrenous{Ce paragraphe est surtout intéressant si l'on démontre le
(projet de)
\thref{thEtaleDin}. Néanmoins le \thrf{thIdmEtale} est utilisé dans le
chapitre~\ref{ChapGalois}
section Galois sur \cdi pour expliquer où se trouvent les \idmsz.

Je suis un peu surpris de ne pas arriver à obtenir le résultat 
qui suit en supposant seulement~$\gA$ \stfe sur~$\gk$.
}
%-% Fin ENTRENOUS

Le \tho suivant est une conséquence subtile du remarquable lem-
\linebreak 
me~\ref{lemPolCarInt}. Il sera utile en théorie de Galois pour le \thref{thZsuffit}. 

%:     Theorem{thIdmEtale}
\begin{theorem}\label{thIdmEtale}
Soit $\rho:\gk\to\gk'$ un \homo injectif d'anneaux avec~$\gk$  \icl dans~$\gk'$, et~$\gA$  une \klg  \stez. Par \eds on obtient ${\gA'}=\rho\ist(\gA)\simeq\gk'\otimes_\gk \gA$  \ste sur~$\gk'$. 
\begin{enumerate}
\item L'\homoz~$\gA\to\gA'$ est injectif.
\item L'anneau~$\gA$ est \icl dans~$\gA'$.
\item Tout \idm de~$\gA'$ est dans~$\gA$.
\end{enumerate}
\end{theorem}
\begin{proof}
Le point \emph{3} est un cas particulier du point \emph{2.}

\emph{1.} On applique le \tho de structure locale des \mptfs et le \plgrf{plcc.basic.modules} pour les suites exactes.

\emph{2.}
On peut
identifier~$\gk$ à un sous-anneau
de~$\gk'$ et~$\gA$ à un sous-anneau de~$\gA'$.
Rappelons que~$\gA$ est finie, donc entière, sur~$\gk$. 
Il suffit de traiter le cas où~$\gA$ est libre sur~$\gk$  (\tho de structure locale des \mptfs  et \plgrf{plcc.entier} pour les \elts entiers).\\
Soit $\ue=(e_1,\ldots,e_n)$ une base de~$\gA$ sur~$\gk$ et $\uh$ la base duale relativement à la \ftrz. 
Si $n=0$ ou $n=1$ le résultat est évident. On suppose $n\geq2$.
\\
Notons que  $\ue$ est aussi une base de~$\gA'$ sur~$\gk'$.
 En outre, puisque, pour~$a\in\gA$, les \endos $\mu_{\gA,a}$ et $\mu_{\gA',a}$
ont la même matrice sur $\ue$, la \ftr sur~$\gA'$ est une extension de la \ftr sur~$\gA$ et $\uh$ reste la base duale relativement à la \ftr
dans~$\gA'$. Soit $x=\sum_ix_ie_i$ un \elt de~$\gA'$ entier sur~$\gA$ ($x_i\in\gk'$).
On doit montrer que les $x_i$ sont dans~$\gk$, ou ce qui revient au même,
entiers sur~$\gk$. Or  $xh_i$ est entier sur~$\gk$.  La matrice de~$\mu_{\gA',xh_i}$ est donc un \elt de $\Mn(\gk')$ entier sur~$\gk$.
Par suite son \polcar a ses \coes entiers sur~$\gk$ (lemme~\ref{lemPolCarInt}),
donc dans~$\gk$, et
en particulier $x_i=\Tr_{\gA'/\gk'}(xh_i)\in\gk$.
\end{proof}
%

%:     Lemma{lemDiag}
\begin{lemma}\label{lemDiag}
La \klg produit~$\gk^n$ est \stez, 
le \discri de la base canonique est égal à $1$.
Si~$\gk$ est un anneau connexe non trivial, 
cette \klg  possède exactement $n$ \crcs et $n!$
\autos (ceux que l'on voit au premier coup d'oeil). 
\end{lemma}
%--------- fin lemma ---------------------------------------------- 
%
\begin{proof}  L'affirmation concernant le \discri est claire (proposition~\ref{propTransDisc}).\\
On a évidemment comme \crcs les $n$ \prns naturelles
$\pi_i:\gk^n\to\gk$ sur chacun des facteurs, et comme~$\gk$-\autos les $n!$ \autos obtenus par permutation des \coosz. 
Soit $e_i$ l'\idm défini  par $\Ker\pi_i=\gen{1-e_i}$. 
 Si $\pi:\gk^n\to\gk$ est un \crcz, les $\pi(e_i)$ forment un \sfio de~$\gk$.
Puisque~$\gk$ est connexe non trivial, ils sont tous nuls sauf un, $\pi(e_j)=1$ par exemple. \linebreak 
Alors,  $\pi=\pi_j$, car ce sont des \klis qui coïncident sur les~$e_i$.
Enfin, comme conséquence tout~$\gk$-\auto de~$\gk^n$ permute les~$e_i$. 
\end{proof}
%
%

%:     Definition{defiDiagAlg}
\begin{definition} \emph{(Algèbres diagonales)}  \label{defiDiagAlg} 
\begin{enumerate}
\item Une \klg  est dite \emph{diagonale} si elle est isomorphe à
une \alg produit~$\gk^n$ pour  un $n\in\NN$. En particulier, elle est \stez.

\item Soit~$\gA$ une \klg \stfe et~$\gL$ une \klgz. \\
On dit que
\emph{$\gL$ diagonalise~$\gA$} si~$\gL\otimes _\gk\gA$ est une
\Llg diagonale.
\end{enumerate}
\index{diagonaliser}
\end{definition}

%:     Fact{factDiagAlg}
\begin{fact}\label{factDiagAlg} \emph{(Algèbres diagonales monogènes)}\\
Soit  $f\in\kX$ un \polu de degré $n$ et~$\gA=\aqo\kX{f}$.
\begin{enumerate}
\item La \klgz~$\gA$ est diagonale \ssiz$f$ est \spl et se décom\-pose en facteurs \lins dans $\kX$.
\item  Dans ce cas, si~$\gk$ est connexe non trivial,~$f$ admet exactement $n$ zéros dans~$\gk$, et la \dcn de~$f$ est unique à l'ordre des facteurs près. 
\item Une \klgz~$\gL$ diagonalise~$\gA$ \ssiz$\disc(f)$ est \iv dans 
$\gL$ et~$f$ se décom\-pose en facteurs \lins dans~$\gL[X]$. 
\end{enumerate}

\end{fact}
\begin{proof} \emph{1.}
Si~$f$ est \spl et se factorise complètement, on a un \isoz~$\gA\simeq\gk^n$ par le \tho d'interpolation de Lagrange (exercice~\ref{exoLagrange}).
\\
Montrons la réciproque. 
Tout \crc $\kX\to\gk$ est un \homo d'\evnz, donc  tout
\crcz~$\gA\to\gk$ est l'\evn en un zéro de~$f$ dans~$\gk$.
Ainsi l'\iso donné en hypothèse est  
de la forme

\snic{\ov g\mapsto \big(g(x_1),\ldots,g(x_n)\big)\quad  (x_i\in\gk\hbox{ et }f(x_i) = 0).}

%\sni

Soit alors $g_i$ vérifiant   $g_i(x_i) = 1$ et, pour $j \ne i$, $g_i(x_j) = 0$. Pour $j \ne i$, l'\elt $x_i - x_j$ divise $g_i(x_i) - g_i(x_j) = 1$, donc $x_i - x_j$ est \ivz. Ceci implique que $f=\prod_{i=1}^n(X-x_i)$ (toujours par Lagrange).

 \emph{2.} Avec les notations précédentes on doit montrer que les seuls zéros de~$f$ dans~$\gk$ sont les $x_i$. Un zéro de~$f$ correspond à un \crc 
$\pi:\gA\to\gk$. On doit donc démontrer que~$\gk^n$ n'admet pas d'autre \crc
que les \prns sur chaque facteur. Or cela a été démontré
dans le lemme~\ref{lemDiag}.  

 \emph{3.} Appliquer le point \emph{1} à la \Llgz~$\gL\te_\gk\gA\simeq\aqo{\gL[X]} f$.
\end{proof}

\rems~\\
 1) Le point \emph{2} nécessite~$\gk$ connexe.

 2) (Exercice laissé \alecz) Si~$\gk$ est un \cdi et si $A$ est une
matrice de $\Mn(\gk)$, dire que~$\gL$ diagonalise
$\gk[A]$ signifie que cette
matrice est \gui{\digz} dans $\Mn(\gL)$, au sens (faible) que~$\gL^n$
est somme directe des sous-espaces propres de $A$.

\rdb%
3) La \dcn d'un anneau~$\gA$ en produit fini d'anneaux connexes non nuls, quand elle est possible, est unique à l'ordre des facteurs près. 
Chaque facteur connexe, isomorphe à un localisé~$\gA[1/e]$,  correspond en effet à un \idm $e$ \emph{indécomposable}\footnote{L'\idm $e$ est dit indécomposable si une \egt $e=e_1+e_2$ avec $e_1$, $ e_2$ \idms et $e_1e_2=0$ implique $e_1=0$ ou $e_2=0$}. Ceci peut
 être compris 
comme conséquence de la structure des \agBs finies
(voir le \thref{factagb}). On peut aussi obtenir le résultat en raisonnant avec un \sfio comme dans la \dem du lemme~\ref{lemDiag}.%
\index{indecomp@indécomposable!idempotent ---}

4)  Dans le point \emph{2}, l'hypothèse \gui{non trivial} donne un énoncé plus usuel. Sans cette hypothèse on aurait dit dans la première partie de la phrase: tout zéro de~$f$ est donné par l'un des $x_i$
correspondant à la \dcn supposée de~$f$ en facteurs \linsz.

5) Pour l'essentiel le fait précédent est une reformulation
plus abstraite du \tho d'interpolation de Lagrange. 
\eoe

%:     Proposition{propEtaleCdi}
\begin{proposition}\label{propEtaleCdi}
Soit~$\gK$ un \cdi \splz ment factoriel et~$\gB$ une \Klg \stfez.
Alors,~$\gB$ est étale \ssi elle est diagonalisée par un surcorps de~$\gK$
étale sur~$\gK$.
\end{proposition}
\begin{proof}
Supposons~$\gB$ étale. Elle est isomorphe à un produit de corps~$\gK_i$
étales sur~$\gK$ (\thref{th2Etale}) 
et il existe un corps~$\gL$ étale sur~$\gK$,
extension galoisienne qui contient une copie de chaque~$\gK_i$ 
(corolaire~\ref{corth3Etale}).
On voit facilement que~$\gL$ diagonalise~$\gB$.
\\
Supposons qu'un corps~$\gL$ étale sur~$\gK$ diagonalise~$\gB$.
Alors, $\Disc_{\gB/\gK}$ est \iv dans~$\gL$ donc dans~$\gK$:~$\gB$ est étale.
\end{proof}

%-% ENTRE NOUS
\entrenous{
Dans le cas d'un \cdi non \spbz ment factoriel, on doit pouvoir diagonaliser une
\Klg étale au moyen d'une autre \Klg étale (qui lui ressemble).

En fait on serait intéressé par un \tho beaucoup plus \gnl
du type suivant.

NB: Le \tho suivant est un peu prématuré si l'on veut vraiment
$\gB$ galoisienne.

%:     Theorem{thEtaleDin}
\begin{theorem} 
\label{thEtaleDin}
Toute \klg \stez~$\gA$ de rang constant peut être diagonalisée par une \klg
fidèle et \stez~$\gB$, et même par une \aG $(\gk,\gB,G)$. 
\end{theorem}

S'il y avait une solution \gnq à ce \pb cela serait une jolie \gnn de l'\aduz. Lorsque~$\gA$ est une \klg monogène
$\aqo{\kT}{f}$
(avec~$f$ \splz) le candidat le plus naturel pour~$\gB$ est $\Adu_{\gk,f}$.

Et il semble bien qu'après \eds fidèle, on puisse rendre toute \alg
\ste produit d'\algs \stes monogènes (voir Ferrand et sa notion de ``\lotz'')

Tant qu'à faire on pourrait
avoir une \aG \gui{assez universelle} modulo la donnée d'\homos
séparants~$\gA\to\gB$ (correspondant aux \homos $t\mapsto x_i$
dans le cas monogène).

Si cela marche cela doit s\^urement exister dans la littérature,
mais où???
}
%-% Fin ENTRENOUS
%

%%%%%%%%%%%%%%%%%%%%%%%%%%%%%%%%%%%%%%%%%%%%%%%%%%%%%%%%%%%%%%%%%%%%%%%%%%%
%  sec{Idempotent de séparabilité}
\penalty-2500
\section[Algèbres \spbsz]{Algèbres \spbsz, \idstz} %secAlgSpb
\label{secAlgSpb}

Les résultats de cette section seront utilisés dans la section
\ref{secAGTG} consacrée aux \aGsz, mais uniquement pour le \thref{thCorGalGen}
qui établit la correspondance galoisienne dans le cas connexe.

Par ailleurs, ils sont également très utiles %plus tard 
%\ref{chap Differentielles} consacré aux 
pour l'étude du \mdiz. Nous nous limiterons ici à parler de \dvnsz.

%--- Notation{NotaAGenv}-------
\begin{definotas}
\label{NotaAGenv}
Soit une \klgz~$\gA$.
\begin{enumerate}
\item  L'\algz~$\gA\otimes_\gk\gA$,
appelée \emph{\alg enveloppante} de~$\gA$, est notée $\env{\gk}{\gA}$.
\item  
Cette \klg possède deux structures naturelles de \Algz, données respectivement par les \homos $g\iAk:a\mapsto a\te 1$
(structure à gauche) et $d\iAk:a\mapsto 1\te a$ (structure à droite).
\index{algèbre!enveloppante}
 Nous utiliserons la
notation abrégée suivante pour les deux structures de
\Amo correspondantes: pour $a\in\gA$ et~$\gamma\in\env\gk\gA$,

\snic{
a\cdot \gamma\,=\, g\iAk(a)\gamma\,=\,(a\te 1) \gamma \,\et\,
 \gamma\cdot a\,=\, d\iAk(a)\gamma\,=\, \gamma(1\te a).
}

\item Nous noterons  $\rJ\iAk$  (ou $\rJ$ si le contexte est clair) l'\id 
de $\env\gk\gA$ engendré par les \elts de la forme
$a\te1-1\te a=a\cdot 1_{\env\gk\gA}-1_{\env\gk\gA}\cdot a$. 
\item Nous introduisons aussi  les \klis suivantes:
\begin{eqnarray}\preskip.0em \postskip.2em
\label{eqDiAk}
\Delta\iAk:\gA\to\rJ\iAk,\quad a\mapsto a\te1-1\te a.\qquad\\
\label{eqmuAk}
\mu\iAk \,:\,\env\gk\gA\rightarrow \gA,\,\,a\otimes b \mapsto ab
\quad\mathrm{(multiplication)}
\end{eqnarray}
\item Dans le cas où $\gA$ est une $\gk$-\atfz, $\gA=\gk[\xn]$, il en est de même
de $\env{\gk}{\gA}$ et l'on a la description suivante possible des objets précédents.
\begin{itemize}
\item $\env{\gk}{\gA}=\gk[\yn,\zn]=\gk[\uy,\uz]$ avec
$y_i = x_i\te 1$, $z_i = 1\te x_i$.
\item Pour $a=a(\ux)\in\gA$, et  $h(\uy,\uz)\in\gk[\uy,\uz]$, on a:
\begin{itemize}
\item 
$g\iAk(a)=a(\uy)$,  $d\iAk(a)=a(\uz)$,
\item 
$a\cdot h=a(\uy)h(\uy,\uz)$, $h\cdot a= a(\uz)h(\uy,\uz)$,
\item $\Delta\iAk(a)=a(\uy)-a(\uz)$,
\item et $\mu\iAk(h)=h(\ux,\ux)$.
\end{itemize} 
\item  $\rJ\iAk$ est l'\id de $\gk[\uy,\uz]$ engendré par les $y_i-z_i$.
\end{itemize}
%:HHH ajout d'un point
\item Enfin, dans le cas où $\gA=\aqo{\kXn}{\lfs}=\gk[\ux]$, autrement dit lorsque $\gA$ est 
une \klg \pfz, il en est de même de $\env{\gk}{\gA}$ (voir le \thref{factSDIRKlg}).

\snic {
\env{\gk}{\gA}=\aqo{\gk[\Yn,\Zn]}{\uf(\uY),\uf(\uZ)}=\gk[\uy,\uz]
}

\end{enumerate}
\end{definotas}
%--- end-notation-----------------------------------------
Notons que 
$\mu\iAk(a\cdot\gamma)=a\mu\iAk(\gamma)=\mu\iAk(\gamma\cdot a)$ pour~$\gamma\in\env\gk\gA$ et $a\in\gA$.

%%%%%%%%%%%%%%%%%%%%%%%%%%%%%%%%%%%%%%%%%%%%%%%%%%%%%%%%%%%%%%%%%%%%%%%%%%%
%:  subsec idst

\subsec{Vers l'\idm de \sptz}

%:     Fact{fact2OmAbsrait}
\begin{fact}\label{fact2OmAbsrait}~

\begin {enumerate}
\item
L'application $\mu\iAk$ est un \crc de \Algs
(pour les deux structures).
\item
On a $\rJ\iAk =\Ker(\mu\iAk )$. Donc~$\gA\simeq\env\gk\gA\sur{\rJ\iAk}$
et
$$
\env\gk\gA = (\gA \otimes 1) \oplus \rJ\iAk =
(1 \otimes \gA) \oplus \rJ\iAk 
,$$
et $\rJ\iAk$ est le \Amo à gauche (ou à droite) engendré par
$\Im\Delta\iAk$.
%:HHH ajout du point 3
\item Dans le cas où $\gA=\aqo{\kXn}{\lfs}=\gk[\ux]$ on obtient
$$
\gk[\uy,\uz]=\gk[\uy]\oplus \gen{y_1-z_1,\dots,y_n-z_n}=\gk[\uz]\oplus \gen{y_1-z_1,\dots,y_n-z_n}.
$$
\end {enumerate}
\end{fact}
\begin{proof}
L'inclusion $\rJ\iAk\subseteq\Ker(\mu\iAk)$ est claire. En notant $\Delta$
pour $\Delta\iAk$, on~a:

\snac
{
\som_ia_i\te b_i = \big(\som_i a_ib_i\big)\te 1 - \som_i a_i \cdot \Delta(b_i) =
1\te \big(\som_i a_ib_i\big) - \som_i \Delta(a_i) \cdot b_i
.}

%\sni
On en déduit $\Ker(\mu\iAk)$ est le \Amo (à droite ou à gauche)
engendré par $\Im\Delta$ et donc qu'il est contenu dans $\rJ\iAk$.  On 
conclut avec~\ref{prdfCaracAlg}.
\end{proof}
%

%%%%%%%%%%%%%%%%%%%%%%%%%%%%%%%%%%%%%%%%%%%%%%%%%%%%%%%%%%%%%%%%%%%%%%%%%%%
\exl Pour~$\gA=\kX$, on a $\env\gk\gA\simeq \gk[Y,Z]$ avec les
\homos 

\snic{h(X)\mapsto h(Y)\hbox{  (à gauche) et } h(X)\mapsto h(Z)\hbox{  (à droite)}}

%\sni

donc $h\cdot g=h(Y)g$ et $g\cdot h=h(Z)g$. On a aussi

\snic{\Delta\iAk(h)=h(Y)-h(Z)$,
 $\mu\iAk\big(g(Y,Z)\big)=g(X,X)$
et $\rJ\iAk=\gen{Y-Z}.}

%\sni
On voit que $\rJ\iAk$ est libre de base $Y-Z$ sur $\env\gk\gA$,
et comme \Amo à gauche, il est libre de base $\big((Y-Z)Z^n \big)_{n\in\NN}$. 
\eoe

%:     Fact{factOmAbsrait}
\begin{fact}\label{factOmAbsrait} On note $\Delta$ pour $\Delta\iAk$.
\begin{enumerate}
% 1
\item Pour $a,b\in\gA$ on a $\Delta(ab)=\Delta(a)\cdot b+a\cdot\Delta(b)$. Plus
\gnltz,
\vspace{.5mm}

\snic{
\arraycolsep2pt
\begin{array}{rcl}
\Delta(a_1\cdots a_n)  & =  & \Delta(a_1)\cdot a_2\cdots a_{n}+a_1\cdot \Delta(a_2)\cdot a_3\cdots a_n +
\cdots \\[.8mm]
  & &
\quad{}+
a_1\cdots a_{n-2}\cdot\Delta(a_{n-1})\cdot a_n
+
a_1\cdots a_{n-1}\cdot\Delta(a_n).
\end{array}
}
%2
\item Si~$\gA$ est une \klg \tfz, engendrée par $(\xr)$, $\rJ\iAk$ est un \itf
de $\env\gk\gA$,  engendré par $(\Delta(x_1),\ldots,\Delta(x_r))$.
\item  Sur l'\id $\Ann(\rJ\iAk)$, les deux structures
de \Amos à gauche et à droite coïncident.
De plus, pour $\alpha \in \Ann(\rJ\iAk)$ et~$\gamma \in \env\gk\gA$,
on~a:
\begin{equation}
\label{eqa2struct}
\gamma\alpha =  \mu\iAk(\gamma) \cdot \alpha =
\alpha \cdot\mu\iAk(\gamma)
\end{equation}

\end{enumerate}

\end{fact}
%-------- end fact --------------------------------
%
\begin{proof} \emph{1.} Calcul \imdz. Le point \emph{2} en résulte
puisque $\rJ\iAk$ est l'\id engendré par l'image de $\Delta$, et que
pour tout \gui{\momz} en les \gtrsz, par exemple $x^3y^4z^2$, $\Delta(x^3y^4z^2)$
est égal à une \coli (à \coes dans $\env\gk\gA$) des images
des \gtrs $\Delta(x)$, $\Delta(y)$ et $\Delta(z)$.

 \emph{3.} L'\id $\fa=\Ann(\rJ\iAk)$ est un $\env\gk\gA$-module, 
donc il est stable pour les deux lois de \Amoz. 
Montrons que ces deux structures coïncident.
Si~$\alpha\in\fa$, pour tout $a\in\gA$ on a 
$0=\alpha(a\cdot 1-1\cdot a)=a\cdot \alpha-\alpha\cdot a$.
\\
 L'\egt (\ref{eqa2struct}) découle du fait que
$\gamma - \mu\iAk(\gamma) \cdot 1 $ et $\gamma - 1 \cdot \mu\iAk(\gamma)$
 sont \hbox{dans $\Ker\mu\iAk=\rJ\iAk$}.
\end{proof}
%

%:      lemmeIdempotentAnnJ
\begin{lemma}\label{lemmeIdempotentAnnJ} 
L'idéal $\rJ\iAk$ est engendré par un \idm si, et seulement~si, 
$$\preskip-.1em \postskip.4em 
1 \in
\mu\iAk\big(\Ann(\rJ\iAk)\big). 
$$
De plus, si $1 = \mu\iAk(\vep)$  avec  $\vep
\in \Ann(\rJ\iAk)$, alors $\vep$ est un \idmz, et l'on~a
$$\preskip-.1em \postskip.4em 
\Ann(\rJ\iAk) = \gen {\vep}\hbox{  et  }\rJ\iAk = \gen {1 - \vep}, 
$$
de sorte que $\vep$ est déterminé de manière unique.
\end {lemma}

\begin{proof} On omet les~$\gA\sur\gk$ en indice.
Si $\rJ = \gen {\vep}$ avec un \idm $\vep$, on obtient les \egts $\Ann(\rJ) = \gen {1 - \vep}$ et
$\mu(1 - \vep) = 1$. \\
Réciproquement, supposons $1 = \mu(\vep)$ avec
$\vep \in \Ann(\rJ)$. Alors $\mu(1 - \vep) = 0$, donc $1 - \vep \in \rJ$, puis
$(1-\vep)\vep = 0$, i.e. $\vep$ est \idmz. \\
Et l'\egt $1 =
(1-\vep) + \vep$ implique que $\Ann(\rJ) = \gen {\vep}$ et $\rJ = \gen {1 -
\vep}$.
\end{proof}

\subsubsection*{Matrice bezoutienne d'un \syp} 

Soient $f_1$, \ldots, $f_s \in \gk[\Xn]=\kuX$.\\ 
On définit la \ixx{matrice} {bezoutienne} du \sys $(\uf)=(f_1, \ldots, f_s)$
en les variables $(\Yn,\Zn)$
par%
\index{besoutienne!matrice}
\[ 
\begin{array}{ccc} 
  \BZ_{\uY,\uZ}(\uf)=(b_{ij})_{i\in\lrbs, j \in \lrbn},\;\hbox{ où }   \\[2mm] 
b_{ij} = 
{f_i(Z_{1..j-1}, Y_j, Y_{j+1..n}) - f_i(Z_{1..j-1}, Z_j, Y_{j+1..n}) \over
Y_j - Z_j} \,.\end{array}
\]

Ainsi pour $n = 2$, $s=3$:

\snic {
\BZ_{\uY, \uZ}(f_1, f_2, f_3) = 
\cmatrix {
{f_1(Y_1,Y_2) - f_1(Z_1,Y_2) \over Y_1-Z_1} & 
{f_1(Z_1,Y_2) - f_1(Z_1,Z_2) \over Y_2-Z_2} \cr
\noalign {\smallskip}
{f_2(Y_1,Y_2) - f_2(Z_1,Y_2) \over Y_1-Z_1} & 
{f_2(Z_1,Y_2) - f_2(Z_1,Z_2) \over Y_2-Z_2} \cr
\noalign {\smallskip}
{f_3(Y_1,Y_2) - f_3(Z_1,Y_2) \over Y_1-Z_1} & 
{f_3(Z_1,Y_2) - f_3(Z_1,Z_2) \over Y_2-Z_2} \cr
}.}

%\sni
Pour $n=3$, la ligne $i$ de la matrice bezoutienne est:

\snic {
\big[{\,
{f_i (Y_1,Y_2,Y_3) - f_i (Z_1,Y_2,Y_3) \over Y_1-Z_1} \;
{f_i (Z_1,Y_2,Y_3) - f_i (Z_1,Z_2,Y_3) \over Y_2-Z_2} \;
{f_i (Z_1,Z_2,Y_3) - f_i (Z_1,Z_2,Z_3) \over Y_3-Z_3}}\,\big].
}

%\sni
On a l'\egt:
$$
\BZ_{\uY,\uZ}(\uf) \cdot \cmatrix {Y_1 - Z_1\cr \vdots\cr Y_n - Z_n\cr}
= \cmatrix {f_1(\uY) - f_1(\uZ)\cr \vdots\cr f_s(\uY) - f_s(\uZ)\cr}
\eqno(\star)
$$
De plus $\BZ_{\uX,\uX}(\uf)=\JJ_\uX(\uf)$, la matrice jacobienne de 
$(f_1, \ldots, f_s)$.

Considérons maintenant une \klg \tf 

\snic{\gA =
\gk[\xn] = \gk[\ux],}

%\sni
avec des \pols $f_i$ vérifiant $f_i(\ux) = 0$
pour tout $i$. 
Son \alg
enveloppante est $\env\gk\gA = \gk[\yn,\zn]$ (notations de début de section).

Alors la matrice $\BZ_{\uy,\uz}(\uf)\in\MM_{s,n}(\env{\gk}{\gA})$  a pour
image par $\mu\iAk$
la matrice jacobienne $\JJ_\ux(f_1, \ldots, f_s)\in\MM_{s,n}(\gA)$.

Pour $D $ mineur d'ordre $n$ de  $\BZ_{\uy,\uz}(\uf)$,  
l'\egt  $(\star)$ montre que
$D \,(y_j - z_j) = 0$  pour $j\in\lrbn$. Autrement dit $D  \in \Ann(\rJ\iAk)$.  La matrice bezoutienne
nous permet donc de
construire des \elts de l'\id $\Ann(\rJ\iAk)$.  
\\
En outre, $\delta :=
\mu\iAk(D )$ est le mineur correspondant dans 
$\JJ_\ux(\uf)$.

Donnons  une application lorsque la transposée
$\tra{\JJ_{\ux}(\uf)} :\gA^s \to \gA^n$ est surjective, i.e.
 $1\in\cD_n(\JJ_{\ux}(\uf))$. On a  donc une
\egt $1 = \sum_{I\in\cP_n} u_I\delta_I$ dans~$\gA$, où $\delta_I$
est le mineur de la matrice extraite de $\JJ_{\ux}(\uf)$ sur les lignes d'indices $i\in I$.  En posant $\vep =
\sum_{I\in\cP_n} u_ID _I \in \env{\gk}{\gA}$, on obtient $\mu\iAk(\vep) =
1$ avec $\vep \in \Ann(\rJ\iAk)$.
\\
Bilan: $\vep$ est ce que l'on appelle l'\idst de $\gA$, et $\gA$ est
une \alg\spbz, notions définies plus loin (\dfnz~\ref{defiSpb}).\\
Donc, si~$\gA$ est une \klg \pf $\aqo{\gk[\uX]}{\uf}$ et 
si l'\ali $\tra{\JJ_{\ux}(\uf)} :\gA^s \to \gA^n$ est surjective,  
alors $\gA$ est \spbz.
\\
Plus \gnltz, pour une \apf $\gA = \aqo{\gk[\uX]}{\uf}$, on
va voir que  $\Coker\big(\tra{\JJ_{\ux}(\uf)}\big)$ et $\rJ\iAk \big/ \rJ\iAk^2$
sont des \Amos isomorphes (\thref{thDerivUnivPF}).

%:subsec{Dérivations} 
\subsec{Dérivations}

%:     Definition{defiDeriv}
\begin{definition}\label{defiDeriv}
Soit $\gA$ une \klg et $M$ un \Amoz. \\
On appelle \emph{$\gk$-dérivation de $\gA$
dans $M$}, une \kli $\delta$ qui vérifie l'\egt de Leibniz
$$
\delta(ab)=a\delta(b)+b\delta(a).
$$%
\index{derivation@dérivation!d'une \alg dans un module}%
On note $\Der\gk\gA M$ le \Amo des $\gk$-\dvns de $\gA$ dans $M$.%
\index{derivation@dérivation!module des ---}
Une \dvn de $\gA$  \gui{tout court} est une \dvn à valeurs dans 
$\gA$. Lorsque le contexte est clair,  $\mathrm{Der}(\gA)$ 
est une abréviation pour~$\Der{\gk}{\gA}{\gA}$.%
\index{derivation@dérivation!d'une \algz}   
\end{definition}
%--------- fin definition ----------------------------------------------

Notez que $\delta(1)=0$ car $1^{2}=1$, et donc $\delta\frt\gk=0.$ 

%:     Theorem{thDerivUniv}
\begin{thdef}\label{thDerivUniv} \emph{(Dérivation \uvle de K\"ahler)}%
\index{derivation@dérivation!universelle}
\\
Le contexte est celui de la \dfn \ref{NotaAGenv}.
\begin{enumerate}
\item Sur $\rJ/\rJ^{2}$ les deux structures de \Amo (à gauche et à droite)
coïncident.
\item L'application composée $\rd:\gA\to\rJ/\rJ^{2}$, définie par $\rd(a)=\ov{\Delta(a)}$, est \linebreak 
une $\gk$-\dvnz.
\item C'est une $\gk$-\dvn \uvle au sens suivant. %
\\
Pour tout \Amo $M$ et
toute $\gk$-\dvn $\delta:\gA\to M$, il existe une unique \Ali $\theta:\rJ/\rJ^{2}\to M$
telle que $\theta\circ \rd=\delta$.
\Pnv{\gA}{\rd}{\delta}{\rJ/\rJ^{2}}{\theta}{M}{}{$\gk$-\dvnsz}{\Alisz.}
\end{enumerate}
Le \Amo $\rJ/\rJ^{2}$, noté  $\Om{\gk}{\gA}$, est appelé le \emph{module des  \diles 
(de K\"ahler) de~$\gA$.} 
\end{thdef}%
\index{module!des diff@des \diles (de K\"ahler)}%
\index{differentielle (de@\dile (de K\"ahler)}
%--------- fin theorem ---------------------------------------------- 
%
\begin{proof}
Les points \emph{1} et \emph{2} sont laissés \alecz.\\
\emph{3.} L'unicité est claire, montrons l'existence. \\
On définit l'\kli
$\tau:\gA\te_\gk\gA\to M$
\tri{\gA}{\Delta}{\delta}{\gA\te_\gk\gA}{\tau}{M}{$\tau(a\te b)=-a\,\delta(b)$.}

\vspace{-5mm}
Le diagramme commute et $\tau$ est $\gA$-\lin à gauche.\\
Il reste à voir que $\tau(\rJ^{2})=0$, car $\theta$ est alors définie
par restriction et passage au quotient de $\tau$. %\\ 
 On vérifie que $\tau(\Delta(a)\Delta(b))=b\delta(a)+a\delta(b)-\delta(ab)=0$.
\end{proof}

On considère maintenant le cas d'une \apf 

\snic{\gA=\aqo{\kXn}{\lfs}=\gk[\ux].}

%\sni
On va utiliser les notations \ref{NotaAGenv}. 
On rappelle que la matrice jacobienne du \syp est définie comme
$$ \JJ_{\uX}(\uf)\; =\;
\bordercmatrix [\lbrack\rbrack]{
    & X_1                     & X_2                     &\cdots  & X_n \cr
f_1 & \Dpp {f_1}{X_1} &\Dpp {f_1}{X_2}  &\cdots  &\Dpp {f_1}{X_n} \cr
f_2 & \Dpp {f_2}{X_1} &\Dpp {f_2}{X_2}  &\cdots  &\Dpp {f_2}{X_n} \cr
f_i & \vdots                  &                         &        & \vdots              \cr
 & \vdots                  &                         &        & \vdots              \cr
f_s & \Dpp {f_s}{X_1} &\Dpp {f_s}{X_2}  &\cdots  &\Dpp {f_s}{X_n} \cr
}
.$$
Dans le \tho qui suit, on note $\rja=\tra{\JJ_{\uX}(\uf)}:\Ae s\to \Ae n$
l'\ali définie par la matrice transposée, et  $(e_1,\dots,e_n)$ la base canonique 
de~$\Ae n$.
On définit 
\[ 
\begin{array}{rcl} 
\delta:\gA\to\Coker(\rja)  &:   & g(\ux)\mapsto   \som_{i=1}^{n} \Dpp {g}{X_i}(\ux)\,\ov{e_i},\\[1.5mm] 
\lambda:\Ae n \to \rJ/\rJ^{2}  &:   &  e_i\mapsto\rd(x_i)=\ov{y_i-z_i}. 
 \end{array}
\]

%:     Theorem{thDerivUnivPF}
\begin{theorem}\label{thDerivUnivPF}  \emph{(Dérivation \uvle via la jacobienne)}
\begin{enumerate}
\item L'application $\delta$ est une $\gk$-\dvn avec $\delta(x_i)=\ov{e_i}$.
\item L'\Ali $\lambda$  induit par passage au quotient
un \iso $\ov \lambda$ de $\Coker(\rja)$ sur $\rJ/\rJ^{2}$.
\end{enumerate}
En conséquence, $\delta$ est \egmt une \dvn \uvlez.
$$
\xymatrix @C=.75cm@R=.4cm{ 
         &&\Coker(\rja) \ar[dd]^{\ov\lambda}&\delta(x_i)\ar@{<->}[dd]\\
\gA\ar[urr]^{\delta}\ar[drr]_\rd\\
         &&\rJ/\rJ^{2}& \rd(x_i)\\ 
}
$$ 
%--------- fin theorem ----------------------------------------------
\end{theorem}
\begin{proof}
\emph{1.} Laissé \alecz.

\emph{2.} 
On commence par montrer l'inclusion $\Im(\rja)\subseteq \Ker\lambda$, i.e. pour chaque~$k$:

\snic{\lambda\big(\som_{i=1}^{n}\Dpp{f_k}{X_i}(\ux)\,e_i\big)= 0.}

%\sni
Pour $g\in\kuX$ on utilise une formule de Taylor à l'ordre 1:

\snic{g(\uy)\equiv g(\uz)+\som_{i=1}^{n}\Dpp{g}{X_i}(\uz)\,(y_i-z_i) \mod \rJ^{2}.
}

%\sni
Pour $g=f_k$ on a $f_k(\uy)=f_k(\uz)=0$, donc $\som_{i=1}^{n}\Dpp{f_k}{X_i}(\uz)\,(y_i-z_i)\in\rJ^{2}$. Ceci démontre l'\egt ci-dessus
en tenant compte de la loi de \Amo \hbox{sur $\rJ/\rJ^{2}$}.
Cela montre que $\lambda$ passe au quotient, avec 

\snic{\ov\lambda:\delta(x_i)=\ov{e_i}\mapsto\rd(x_i)=\ov{y_i-z_i}.}

%\sni
Par ailleurs, puisque $\delta$ est une $\gk$-\dvnz, la \prt \uvle de la \dvn $\rd:\gA\to\rJ/\rJ^{2}$ nous donne une \fcn
$\gA$-\lin 
$$\preskip.4em \postskip.4em 
\rJ/\rJ^{2} \to \Coker(\rja)\;:\;\rd(x_i)\mapsto \delta(x_i) .
$$
Il est clair que les deux applications sont réciproques l'une de l'autre.
\end{proof}
% 

%:subsec{Idempotent de \spt d'une \asez} 
\subsect{Idempotent de \spt d'une \alg \\ \stez}{Idempotent de \spt d'une \asez} 

Soit~$\gA$  une \klg\stfez.  Pour $a\in\gA$, notons $a\sta =
a\centerdot\Tr\iAk$% (i.e. $ a\sta(x)= \Tr\iAk(ax)$)
.
On a une \kli canonique $\env{\gk}{\gA}
\to \End_\gk(\gA)$, composée de l'\ali $\env{\gk}{\gA} \to
\Asta\te_\gk\gA$, $a \otimes b \mapsto a\sta \otimes b$, 
et de l'\iso naturel $\Asta\te_\gk\gA \to \End_\gk(\gA)$.

Si~$\gA$ est \ste ces \alis sont toutes des \isosz.
Alors, si $\big((x_i), (y_i) \big)$ est un \stycz,
l'\elt $\sum_i x_i\otimes y_i$ est
indépendant du choix du \sys car son image  dans $\End_\gk(\gA)$
est $\Id_\gA$.  En particulier, $\sum_i x_i\otimes y_i = \sum_i y_i\otimes x_i$.  

Le \tho suivant dégage les
\prts \caras de cet élé-\linebreak ment $\sum_i x_i\otimes y_i$.
Ces \prts  conduisent à la notion d'\alg\spbz.

%:     theorem   thAlgSteIdmSpb
\begin {theorem}\label{thAlgSteIdmSpb} \emph{(Idempotent  de \spt d'une \asez)}
Soit~$\gA$ une \klg\ste et $\big((x_i), (y_i) \big)$ un \styc de~$\gA$.  Alors, l'\elt
$\vep = \sum_i x_i \te y_i \in \env\gk\gA$
vérifie les conditions du lemme~\ref{lemmeIdempotentAnnJ}. En particulier,
$\vep$ est \idm et l'on a
$$
\som_i x_iy_i = 1, \quad  
a\cdot\vep = \vep\cdot a \quad \forall a \in \gA.
$$
\end {theorem}
NB: on démontre la réciproque (pour les \asfsz) un peu plus loin (\thref{thAlgStfSpbSte}).
\begin {Proof} {Preuve dans le cas galoisien (à lire après le \thref{thA}). }\\
Soit $(\gk,\gA,G)$ une \aGz. Puisque le résultat à prouver est
indépendant du \stycz, on peut supposer que les familles $(x_i)$ et $(y_i)$ sont deux \syss
d'\elts de~$\gA$ vérifiant les conditions du point \emph {2} du
\tho d'Artin~\ref{thA}.

Dire que $\mu(\vep) = 1$ consiste à dire que $\sum_i x_i y_i = 1$,
ce que vérifie $\big((x_i), (y_i) \big)$.  Pour montrer que $\sum_i ax_i \otimes
y_i = \sum_i x_i \otimes ay_i$, il suffit d'appliquer $\psi_G$; on note
$(g_\sigma)_\sigma$ l'image du membre gauche, $(d_\sigma)_\sigma$ l'image du
membre droit. On obtient, en notant $\delta$ le symbole de Kronecker:

\snic {
g_\sigma = \sum_i ax_i\sigma(y_i) = a \delta_{\sigma, \Id}, 
\qquad
d_\sigma = \sum_i x_i\sigma(ay_i) = \sigma(a) \delta_{\sigma, \Id}
.}

%\sni
On a bien l'\egt puisque les composantes des deux familles $(d_\sigma)$ et
$(g_\sigma)$ sont nulles sauf en l'indice $\sigma = \Id$, indice pour lequel
leur valeur (commune) est $a$.

Remarquons que $\vep$ est égal à l'\elt $\vep_{\Id}$ introduit dans le lemme
\ref{lemArtin}. Son image par $\varphi_G$ est l'\idm $e_\Id$, ce qui confirme
que $\vep$ est \idmz.
\end {Proof}
\begin {Proof} {Preuve (\gnlez) dans le cas \stez. }\\
On note $\Tr$ pour $\Tr_{\gA\sur\gk}$ et $m_\vep : \env\gk\gA \to
\env\gk\gA$ la multiplication par $\vep$.  On a:
$$ 
\Tr(ab) = \som_i \Tr(ay_i) \Tr(bx_i), \qquad a, b \in \gA.\eqno(\star)
.$$
 En effet, cela découle facilement de l'\egt $a = \sum_i \Tr(ay_i)x_i$.

On réécrit $(\star)$  comme l'\egt des deux formes~$\gk$-\linsz,
$\env\gk\gA \to \gk$:
$$
\Tr_{\gA\sur\gk} \circ\, \mu\iAk = 
\Tr_{\env\gk\gA\sur\gk} \circ\, m_\vep.\eqno(*)
$$
Montrons que $\vep \in \Ann(\rJ)$. Soient $z \in \env\gk\gA$,
$z' \in \rJ$. En évaluant l'\egt $(*)$ en $zz'$, on obtient:
$$
\Tr_{\gA\sur\gk}\big(\mu\iAk(zz')\big) = 
\Tr_{\env\gk\gA\sur\gk}(\vep zz').
$$
Mais $\mu\iAk(zz') = \mu\iAk(z)\mu\iAk(z') = 0$ car $z' \in \rJ =
\Ker\mu\iAk$. 
On en déduit que $\Tr_{\env\gk\gA\sur\gk}(\vep zz') =
0$ pour tout $z \in \env\gk\gA$. Comme $\Tr_{\env\gk\gA\sur\gk}$ est non \dgne
%(car égale à $\Tr_{\gA\sur\gk} \otimes \Tr_{\gA\sur\gk}$), on en déduit que 
on obtient $\vep z' = 0$. 
Ainsi $\vep \in \Ann(\rJ)$.

Il reste à montrer que $\mu\iAk(\vep) = 1$, i.e. $s=\sum_i x_iy_i = 1$. 
\\
L'\egt 
$\Tr(x) = \sum_i \Tr(xx_iy_i)$ (fait~\ref{factMatriceEndo})
%:HHH   ci dessus reference simplifiee
 exprime que $\Tr \big((1-s)x \big) = 0$
pour tout $x \in \gA$ donc $s = 1$.
\end {Proof}

%%%%%%%%%%%%%%%%%%%%%%%%%%%%%%%%%%%%%%%%%%%%%%%%%%%%%%%%%%%%%%%%%%%%%%%%%%%
%:  subsec{Algèbres \spbsz}  subsecAlgSpb
\subsec{Algèbres \spbsz} \label{subsecAlgSpb}

%:     Theorem{thSepIversen}
\begin{theorem}\label{thSepIversen}
Pour une \klgz~$\gA$ \propeq
\begin{enumerate}
\item \label{i1thSepIversen}
$\gA$ est projectif comme $\env\gk\gA$-module.
\item \label{i2thSepIversen}
$\rJ\iAk$ est engendré par un \idm de $\env\gk\gA$.
\item \label{i4thSepIversen}
$\rJ\iAk$ est \tf et \idmz.
\item \label{i5thSepIversen}
$1 \in\mu\iAk\big(\Ann(\rJ\iAk)\big)$.
\item \label{i6thSepIversen}
Il existe $n\in\NN$, et $x_1$, \dots,  $x_n$,  $y_1$, \dots,  $y_n\in\gA$ tels que
$\sum_ix_iy_i=1$ et pour tout $a\in\gA$, $\sum_iax_i\otimes y_i=\sum_ix_i\otimes ay_i$. 
\end{enumerate}
Dans le cas où $\gA$ est  \tfz, on a aussi l'\eqvc avec:
\begin{enumerate}\setcounter{enumi}{5}
\item \label{i7thSepIversen}  
$\Om{\gk}{\gA}=0$.\footnote{Lorsque $\gA=\aqo{\kXn}{\lfs}=\gk[\ux]$, le \thref{thDerivUnivPF} donne la condition suivante pour $\Om{\gk}{\gA}=0$: 
la matrice $\rja(\ux)$, transposée de la matrice jacobienne du \sysz, est surjective, i.e. $1\in\cD_n(\JJ_{\uX}(\uf)(\ux))$.}
\end{enumerate}
On note alors $\vep\iAk$ l'unique \idm qui engendre l'\id $\Ann(\rJ\iAk).$
%\\
%Enfin, si~$\gA$ est une \klg \tfz, alors $\rJ\iAk$
%est un \itf de~$\env\gk\gA$, donc la condition~\ref{i4thSepIversen}. 
%se réduit à: $\rJ\iAk$ est \idmz, i.e. $\Om{\gk}{\gA}=0$.

\end{theorem}
\begin{proof}
Puisque~$\gA\simeq\env\gk\gA\sur{\rJ\iAk}$,
les points~\emph{\ref{i1thSepIversen}} et~\emph{\ref{i2thSepIversen}}
sont \eqvs en application du lemme~\ref{lemIdpPtf} concernant
les modules monogènes \prosz. Les points~\emph{\ref{i2thSepIversen}} et~\emph{\ref{i4thSepIversen}}
sont \eqvs en application du lemme~\ref{lem.ide.idem} sur les \ids \idms \tfz.
Le lemme~\ref{lemmeIdempotentAnnJ} donne l'\eqvc de \emph{\ref{i2thSepIversen}} et
\emph{\ref{i5thSepIversen}.}\\
\emph{\ref{i4thSepIversen}} $\Leftrightarrow$
\emph{\ref{i7thSepIversen}.} Si~$\gA$ est une \klg \tfz, alors $\rJ\iAk$
est un \itf de~$\env\gk\gA$, donc la condition~\emph{\ref{i4thSepIversen}} 
se réduit à: $\rJ\iAk$ est \idmz, \hbox{i.e. $\Om{\gk}{\gA}=0$}.\\ 
Enfin  \emph{\ref{i6thSepIversen}} est la forme concrète de \emph{\ref{i5thSepIversen}.}
\end{proof}
%

%:     Definition{defiSpb}
\begin{definition}\label{defiSpb}
On appelle \emph{\alg \spbz} une \alg qui vérifie les \prts \eqves
énoncées au \thrf{thSepIversen}. L'\idm $\vep\iAk\in\env\gk\gA$ s'appelle
\emph{l'\idstz} de~$\gA$.
\index{idempotent!de séparabilité}
\index{algèbre!separable@séparable}
\index{separable@séparable!algèbre --- }
\end{definition}

\comm
Il faut noter que Bourbaki utilise une notion d'extension
\spb pour les corps assez différente de la \dfn ci-dessus.
En \clama les \algs sur un corps~$\gK$  \gui{\spbs au sens de la \dfnz~\ref{defiSpb}}
sont les \algs \gui{finies et \spbs au sens de Bourbaki} (voir le \tho \vref{thSepProjFi}). Beaucoup d'auteurs suivent Bourbaki
au moins pour les extensions \agqs de corps, qu'elles soient finies ou pas.
Dans le cas d'une \Klg \agq sur un corps discret~$\gK$, la \dfn à la Bourbaki signifie que tout \elt de l'\alg annule un \polu \spl de~$\KT$.
\eoe

%: perso
\perso{Le commentaire qui suit pourrait peut-être mis quelque part.

\comm
 Il est bien dommage que la terminologie se soit fixée ainsi.
La notion la plus naturelle est celle d'\alg \fntz.
Dans \gui{\alg \spbz} il y a en plus une condition de finitude un peu \gui{hors sujet}.
En fait nous avons en vue surtout les \algs \tf et les deux notions coïncident dans ce cas.
\\
Nous aurions donc préféré utiliser une terminologie moins conflictuelle,
et par exemple dire \gui{finiment \spbz} au lieu de \gui{\spbz}, mais
nous n'avons pas osé. \eoe
}
%-% Fin perso

%:     Fact{factSpbEds}
\begin{fact}\label{factSpbEds} \emph{(Stabilité des \algs \spbs par \edsz)}\\
Soit $\imath:\gk\to\gA$ et $\rho:\gk\to\gk'$ deux \klgs et~$\gA'=\rho\ist(\gA)$. 
\\
On a un \iso canonique $\rho\ist(\env\gk\gA)\to 
{\gA'_{\,\gk'}}^{\!\!\mathrm{e}}$
et le diagramme ci-dessous commute

\vspace{-2em} 
\carre{\env\gk\gA}{\env\gk\rho}{{\gA'_{\,\gk'}}^{\!\!\mathrm{e}}}
{\mu\iAk}{\mu_{\gA'\!/\gk'}}{\gA}{\rho}{\gA'}

\vspace{-.8em} 
En particulier, une \alg \spb reste \spb par \edsz.
\end{fact}
%--------- fin fact ---------------------------------------------- 
\facile

\medskip 
Nous montrons maintenant la réciproque du \thref{thAlgSteIdmSpb},
ce qui nécessite un lemme préparatoire.

%:    lemme   lemTracAlgEnv
\begin {lemma}\label{lemTracAlgEnv}
Soit~$\gA$ une \klg \stfe et $\env\gk\gA $ son
\alg enveloppante.
\begin {enumerate}
\item
$\env\gk\gA$ est une~\Alg à gauche \stfe dont la trace
est donnée par~$\gamma_g\circ (\Id_\gA \te \Tr\iAk)$
(où~$\gamma_g:\gA\te_\gk\gk\to\gA$ est l'\iso canonique), i.e. pour $\alpha = \sum_i a_i\te b_i$:

\snic {
\Tr_{(\env\gk\gA\sur\gA)_{\rm g}}(\alpha) = 
\sum_i a_i\Tr\iAk(b_i).
}

%\sni
De même, 
$\env\gk\gA$ est une~$\gA$-\alg à droite \stfe dont la trace
est donnée par~$\gamma_d\circ (\Tr\iAk\te\Id_\gA)$, i.e.
$
\Tr_{(\env\gk\gA\sur\gA)_{\rm d}}(\alpha) = 
\sum_i \Tr\iAk(a_i) b_i.
$

\item
Sur $\Ann(\rJ\iAk)$, les  formes~$\gA$-\lins
$\Tr_{(\env\gk\gA\sur\gA)_{\rm g}}$, $\Tr_{(\env\gk\gA\sur\gA)_{\rm d}}$
et $\mu\iAk$ coïncident, i.e. si $\alpha = \sum_i a_i\te b_i \in
\Ann(\rJ\iAk)$:

\snic {
\sum_i a_ib_i = \sum_i a_i\Tr\iAk(b_i) = \sum_i \Tr\iAk(a_i) b_i.
}
\end {enumerate}

\end {lemma}
%------  end {lemma} --------------------------------------------------

%
\begin {proof}
\emph{1.}
Il s'agit d'un résultat structurel \gnlz: la trace se conserve par \eds (voir fait~\ref{fact.det loc}).
Autrement dit si~$\gk'$
est une \klgz,~$\gk' \te_\gk \gA$ est une~$\gk'$-\alg
\stfe dont la trace est~$\gamma\circ (\Id_{\gk'} \otimes \Tr\iAk)$
où~$\gamma:\gk'\te_\gk \gk\to\gk'$ est l'\iso canonique.

\emph{2.}
De manière \gnlez, sous les hypothèses que $E$ est un \Amo \ptfz, $x \in E$, $\nu \in E\sta$ et
$u=\theta_{E}(\nu\te x) \in \End_\gA(E)$, on obtient \linebreak l'\egt $\Tr_E(u)= \nu(x)$
(voir le fait~\ref{factMatriceEndo}). 
\\
On applique
ceci à $E = \env\gk\gA$, $x = \alpha\in E$ et $\nu = \mu\iAk\in E\sta$, en
notant qu'alors $u = \theta_{E}(\nu\te \alpha) = \mu_{\env\gk\gA,\alpha}$. En effet,
d'après le point \emph {3}  du fait~\ref {factOmAbsrait}, on~a pour
$\gamma \in \env\gk\gA$,
$\gamma\alpha=\mu\iAk(\gamma)\cdot\alpha=\theta_{E}(\nu\te \alpha)(\gamma)$.
\end {proof}

%:    Theorem  thAlgStfSpbSte
\begin {theorem}\label{thAlgStfSpbSte}\emph{(Algèbres \stes et \algs \spbsz)}\\
Toute \klg\spb et \stfez~$\gA$ est \stez. Plus \prmtz, si $\vep\iAk = \sum x_i\te y_i \in \env\gk\gA $ est l'\idst de~$\gA$, alors \smash{$\big((x_i), (y_i) \big)$} est un
\styc de~$\gA\sur\gk$.
\\
En résumé une \alg \stfe est \spb \ssi elle est \stez.
\end {theorem}
NB: \prmtz, le lien entre les deux notions est obtenu par la relation qui lie l'\idst et les \sycsz, comme cela ressort  du \tho direct~\ref{thAlgSteIdmSpb} et du \tho réciproque \ref{thAlgStfSpbSte}.
%------  end theorem --------------------------------------------------

\begin {proof}
Soit $x \in \gA$, alors $(x\te1)\vep\iAk = \sum_i xx_i \te y_i$ est dans
$\Ann(\rJ\iAk)$, donc d'après le lemme~\ref{lemTracAlgEnv}, on a $\sum_i xx_iy_i =
\sum_i \Tr\iAk(x_ix)y_i$.\\
 Et comme $\sum_i x_iy_i = 1$, cela donne $x = \sum_i
\Tr\iAk(x_ix)y_i$. On conclut par la \carn des \ases donnée dans le fait~\ref{factCarAste}.
\end {proof}

Le \tho suivant renforce le \tho précédent et montre que l'existence d'un \idst est une condition de finitude très forte.

%:     Theorem{thSepProjFi}
\begin{theorem}\label{thSepProjFi}
%:HHH  index pour systeme de coordonnees
Soit~$\gA$ une \klg \spbz. \\
On suppose que~$\gA$  possède un \ixe{\sycz}{systeme de coordonnees}, dans le sens suivant: on a un ensemble discret $I$, une famille $(a_i)_{i\in I}$
dans~$\gA$ et une famille~$(\alpha_i)_{i\in I}$ dans le \kmo dual $\Asta=\Lin_\gk(\gA,\gk)$, telles que
pour tout $x\in\gA$ on ait
$$\ndsp\preskip.3em \postskip.2em 
x=\sum_{i\in J_x}\alpha_i(x)a_i, 
$$
ici $J_x$ est une partie finie de $I$, et tous les $\alpha_i(x)$ pour $i\in I\setminus J_x$ sont nuls.
\\
Alors,~$\gA$ est \stfez, donc \stez.%
\\
C'est par exemple le cas si $\gk$ est un \cdi et si $\gA$ est une \klg \pfz.
\end{theorem}
%------  end theorem --------------------------------------------------
%
\begin{proof} 
Concernant le cas particulier, l'\alg quotient possède une base finie ou dénombrable de \momsz, d'après la théorie des \bdgsz.
\\
Soit $\vep=\sum_{k=1}^rb_k\te c_k$ l'\idm de séparabilité.
On a $\vep\cdot x= x\cdot \vep$ pour tout $x\in\gA$,  et $\sum_{k=1}^rb_kc_k=1$.
\\
Pour $\alpha \in \Asta$ et $x \in \gA$, en appliquant 
$1\te \alpha$ à~$x \cdot \vep  = \vep \cdot x$ 
on obtient: 

\snic{\som_k x b_k \alpha(c_k) = \som_k b_k \alpha(xc_k).
}

%\sni
En notant $J$ la partie finie $J = \bigcup J_{c_k}$, on obtient pour chaque~$k$

\snic{
c_k = \som_{i \in J} \alpha_i(c_k) a_i.
}

%\sni
On écrit alors:

\snac{x = \som_{k \in \lrbr} xb_kc_k = \som_{k \in \lrbr, i \in J} 
x b_k \alpha_i(c_k) a_i = \som_{i\in J, k\in\lrbr} \alpha_i(c_kx)b_ka_i.}

%\sni
Ce qui donne maintenant un \syc fini pour $\gA$, avec les \elts~$b_ka_i$ et 
les formes $x\mapsto \alpha_i(c_kx)$ pour $(i,k)\in J\times \lrbr$.
\end{proof}

\comm
Notons que lorsque l'on a un \sys de \coos pour un module, le module est
\pro au sens usuel.
La \dfn d'un \sys de \coos pour un module $M$ revient à dire
que $M$ est isomorphe à un facteur direct du module $\Ae{(I)}$. Ce dernier module, librement engendré par $I$, est \pro parce que $I$ est discret.
\\
\emph{En \clamaz}, tout module \pro possède un \sycz, parce
que tous les ensembles sont discrets, donc le \tho précédent s'applique: 
toute \klg \spb qui est un \kmo \pro est \stfez. Par la même occasion \emph{toute
\alg \spb sur un \cdi ou sur un anneau \zedr est \stfez}. 
\eoe

%-% ENTRE NOUS
\entrenous{Le \tho 6.14, on se demande bien s'il ne devrait
pas fonctionner, avec la conclusion que l'\alg est finie,
 sous l'hypothèse que~$\gA$ est plate sur~$\gk$, qui remplacerait l'existence d'un \sycz.
}
%-% Fin ENTRENOUS

%:HHH on insiste sur la cas des corps discrets

Dans le cas d'une \apf sur un \cdiz, les \thos \ref{thSepIversen} et \ref{thSepProjFi} donnent le résultat suivant.

%:     Corollary{corthSepProjFicdi}
\begin{corollary}\label{corthSepProjFicdi}
Pour $f_1$, \dots, $f_s\in\kXn$ lorsque $\gk$ est un \cdiz, \propeq
\begin{enumerate}
\item L'\alg quotient $\gA=\kux$ est \stez.
\item L'\alg quotient est \spbz.
\item La matrice $\rja(\ux)$, transposée de la matrice jacobienne du \sypz, est surjective.
\end{enumerate} 
\end{corollary}
%--------- fin corollary ---------------------------------------------- 

Nous allons maintenant 
montrer qu'une \alg \spb ressemble beaucoup à une \alg diagonale,
y compris dans le cas où l'anneau de base est quelconque.
\\
Considérons  la \klg
diagonale~$\gk^n$. Notons $(e_1, \ldots, e_n)$ sa base canonique et~$p_i : \gk^n \to \gk$ la forme \coo relative à~$e_i$.  Alors on a: 

\snic{e_i\in\BB(\gk^n)$, $p_i \in \Hom_\gk(\gk^n,
\gk)$,  $p_i(e_i) = 1$ et $xe_i = p_i(x)e_i \;\Tt x \in \gk^n.}

%\sni
  Il s'agit en quelque sorte de généraliser cela aux
\algs\spbsz.

%%%%%%%%%%%%%%%%%%%%%%%%%%%%%%%%%%%%%
%:     lemme  lemIdmHomSpb

\begin {lemma}\label{lemIdmHomSpb}\emph{(Caractères d'une \alg \spbz)}\\
Soit~$\gA$ une \klg\spb avec~$\gk\subseteq\gA$.
\begin {enumerate}
\item Notons $\imath:\gk\to\gA$ l'injection canonique.
Si $\varphi \in \Hom_\gk(\gA,\gk)$, 
%on a pour tout $x\in\gA$
%$x = (x - \varphi(x).1) + \varphi(x).1$,
  $\imath\circ\varphi$ est un \prr d'image~$\gk.1$,
donc

\snic{\gA= \gk.1 \oplus\Ker\varphi \;\hbox{ et }\; \Im(\Id_\gA - \imath\circ \varphi)=\Ker\varphi.}

%\sni
En fait l'\id $\Ker\varphi$ est engendré par un 
\idm  de~$\gA$, on notera~$\vep_\varphi$ l'\idm \copz.
\item
Pour $\varphi$, $\varphi' \in \Hom_\gk(\gA,\gk)$, on a
$\varphi'(\vep_\varphi) = \varphi(\vep_{\varphi'})$.
\\
Cet \eltz, noté~$e_{\{\varphi,\varphi'\}}$ est un \idm de~$\gk$ et l'on~a:
\vspace{-1mm}
\[ 
\begin{array}{ccc} 
\vep_\varphi\vep_{\varphi'} =
e_{\{\varphi,\varphi'\}} \vep_\varphi =
e_{\{\varphi,\varphi'\}} \vep_{\varphi'}=\varphi(\vep_\varphi\vep_{\varphi'})=
\varphi'(\vep_\varphi\vep_{\varphi'}),     \\[2mm] 
\gen {\Im(\varphi - \varphi')}_\gk = \gen {1-e_{\{\varphi,\varphi'\}}}_\gk\; \;\mathit{et}
\;\;\Ann_\gk(\varphi - \varphi') = \gen {e_{\{\varphi,\varphi'\}}}_\gk.
 \end{array}
\]

\item
En conséquence on a les \eqvcsz:

\snic {
e_{\{\varphi,\varphi'\}} = 1 \iff
\vep_\varphi = \vep_{\varphi'} \iff \varphi = \varphi',\,
$  et $\,
e_{\{\varphi,\varphi'\}} = 0 \iff \vep_\varphi \vep_{\varphi'} = 0.
}

\vspace{1mm}
\item
Si~$\gk$ est connexe, deux \idms $\vep_\varphi$, (pour $\varphi \in
\Hom_\gk(\gA, \gk)$), sont égaux ou \ortsz.
\end {enumerate}
\end {lemma}
%---------- end {lemma} ------------------------------------------
%
\begin {proof}
Soit $\vep\iAk = \sum x_i \otimes y_i$. 
On sait que $a\cdot\vep\iAk =
\vep\iAk\cdot a$ pour tout $a \in \gA$, que $\sum x_i \otimes y_i = \sum y_i \otimes x_i$ et que
$\sum_i x_iy_i = 1$.

 \emph{1.} La première affirmation est valable pour tout \crc de toute
\algz~$\gA$ (proposition \ref{prdfCaracAlg}).
Il reste à voir que $\Ker\varphi$ est engendré par un \idmz.
On considère l'\homo de \klgs {\mathrigid 2mu $\nu=\mu\iAk\circ(\varphi\te\Id_\gA):\env\gk\gA\to\gA$},
et   l'\elt $\vep=\nu(\vep\iAk)$. Ainsi $\vep=\som_i \varphi(x_i)y_i$ est un \idm et 
l'on obtient les \egts

\snic{\varphi(\vep)=\som_i \varphi(x_i)\varphi(y_i)=\varphi(\som_ix_iy_i)=\varphi(1)=1.}

%\sni
Donc $1-\vep\in\Ker\varphi$.
\\
En appliquant
$\nu$ à l'\egt $\som_i ax_i \otimes y_i = \som_i x_i
\otimes ay_i$, on obtient $\varphi(a)\vep = a\vep$.
Donc $a\in\Ker\varphi$ implique $a=(1-\vep )a$, et $\Ker\varphi=\gen{1-\vep}$.

 \emph{2.}
On a pour $a \in \gA$:
$$\preskip.3em \postskip.4em
\varphi'(a) \varphi'(\vep_\varphi) = \varphi'(a \vep_\varphi) = 
\varphi'(\varphi(a) \vep_\varphi) = \varphi(a)\varphi'(\vep_\varphi).
\eqno (\star)
$$
Pour $a = \vep_{\varphi'}$, on obtient $\varphi'(\vep_\varphi) =
\varphi(\vep_{\varphi'})\varphi'(\vep_\varphi)$. 
\\
Par symétrie,
$\varphi(\vep_{\varphi'}) = \varphi'(\vep_\varphi)$. Notons $e$
cet \idm de~$\gk$. Par \dfnz, on a $a\vep_\varphi = \varphi(a)
\vep_\varphi$. En faisant $a = \vep_{\varphi'}$, on obtient
$\vep_{\varphi'}\vep_\varphi = e\vep_\varphi$.
\\ 
Enfin, notons $\fa = \gen {\Im(\varphi - \varphi')}$. La relation $(\star)$
montre que $\fa e = 0$. D'autre part $1-e = (\varphi - \varphi')(\vep_\varphi)
\in \fa$. Donc $\fa = \gen {1-e}_\gk$ et $\Ann_\gk(\fa) = \gen {e}_\gk$.

 \emph{3} et \emph{4.}
Découlent du point précédent.
\end {proof}

%%%%%%%%%%%%%%%%%%%%%%%%%%%%%%%%%%%%%%%%%%%%%%%%%%%%%%%%%%%%
%:     lemme  lemSspbdiag
\begin {lemma} \label{lemSspbdiag}
 \emph{(Sous-algèbre \spb d'une extension diagonale)}
\\
Soit~$\gk$ un anneau \emph {connexe non trivial},~$\gB=\gk^n$, $p_i:\gB\to\gk$ la $i$-ième
\prn canonique, $e_i$ l'\idm défini par $\Ker p_i=\gen{1-e_i}$ $(i\in\lrbn)$.
Pour une partie finie $I$ de $\lrbn$ on note $e_I = \sum_{i \in I} e_i$.\\
Soit~$\gA$ une \klg\spb avec
$\gk\subseteq\gA\subseteq\gk^n$ et $\pi_i$ la restriction de $p_i$ à~$\gA$
pour $i \in \lrbn$.  
\begin {enumerate}
\item On considère la relation d'équivalence sur $\lrbn$ définie par $\pi_i = \pi_j$.
La partition correspondante $\cP$ est un ensemble fini de parties finies de~$\lrbn$. 
Pour $J\in\cP$ on note $\pi_J$ la valeur commune des $\pi_j$ pour $j \in J$.
\item
$\gA$ est un \kmoz-libre de base $\sotq{e_J}{J \in \cP}$.
\item $\Asta$ est un
\kmoz-libre de base $\sotq{\pi_J}{J \in \cP}=\Hom_\gk(\gA, \gk)$.
\end {enumerate}
\end {lemma}
%---------- end {lemma} ------------------------------------------

%
\begin {proof}
\emph {1.} Comme~$\gk$ est connexe non trivial, tout \idm de~$\gB$ est de la forme $e_I$ pour une unique partie finie
$I$ de $\lrbn$.
\\
 Soit $i \in \lrbn$. 
D'après le lemme~\ref{lemIdmHomSpb} 
il existe un et un seul \idmz~$\vep_i$ de $\gA$ tel que
$\pi_i(\vep_i) = 1$ et $a \vep_i = \pi_i(a)\vep_i$ pour tout $a\in\gA$. 
Cet \idm est aussi un \idm de~$\gB$ donc de la forme $e_{J_i}$ pour une partie finie~$J_i$ de $\lrbn$. 
Puisque $\pi_i(\vep_i)=p_i(e_{J_i}) = 1$, on a $i \in
J_i$, et la réunion des~$J_i$ est~$\lrbn$. 
Deux $J_i$ distincts sont disjoints
d'après le dernier point du lemme~\ref{lemIdmHomSpb}.
Les~$J_i$ forment donc une partition finie formée de parties finies de~$\lrbn$.
\\
Si  $\pi_i=\pi_j$, alors $\vep_i=\vep_j$ donc $J_i=J_j$.
Si $J_i=J_j$, alors $\vep_i=\vep_j$ et $\pi_i(\vep_j)=1$.
Le point \emph{2} du lemme~\ref{lemIdmHomSpb} donne $1\in\Ann_\gA(\pi_i-\pi_j)$,
donc $\pi_i=\pi_j$.

\emph{2.} Résulte du point \emph{1.}

\emph{3.} 
Soit $\varphi \in \Hom_\gk(\gA,\gk)$. Les $\varphi(e_J)$ sont des \idms de
$\gk$.  Comme~$\gk$ est connexe, on a $\varphi(e_J) = 0$ ou $1$. Mais les $(e_J)_{J
\in \cP}$ forment un \sfioz, donc il y a un seul $J \in \cP$ pour 
lequel~$\varphi(e_J) = 1$ et par suite $\varphi = \pi_J$.  Le reste est \imdz.\end {proof}

%%%%%%%%%%%%%%%%%%%%%%%%%%%%%%%%%%%%%%%%%%%%%%%%%%%%%%%%%%%%%%%%%%%%%%%%%%%
\section{Algèbres galoisiennes, théorie générale}\label{secAGTG}
%-----------------------------------------

Dans la théorie mise au point par Artin,
on considère un groupe fini $G$ d'\autos d'un \cdi
$\gL$, on appelle~$\gK$ le sous-corps des points fixes de $G$ et l'on démontre
que~$\gL$ est une extension galoisienne de~$\gK$, avec $G$ pour groupe de Galois.

Dans la section présente on donne
la \gnn de la théorie d'Artin pour des anneaux commutatifs au lieu de \cdisz.
Une bonne idée de \gui{comment cela peut fonctionner} est déjà donnée par
le petit exemple significatif suivant, qui montre que l'hypothèse \gui{\cdiz} n'est pas
essentielle.

\mni
\textbf{Un petit exemple pour commencer}
\\
Soit~$\gA$ un anneau commutatif,  $\sigma \in \Aut(\gA)$ un \auto d'ordre 3,
$G$ le groupe qu'il engendre. Supposons qu'il existe $x\in\gA$  tel que
$\sigma(x) - x \in \gA^{\times}$. Posons~$\gk =
\gA^{G}$ le sous-anneau des points fixes.  Alors, $(1, x,
x^2)$ est une base de~$\gA$ sur~$\gk$.
En effet, soit $V$ la matrice de Vandermonde
%%%%%%%%%%%%%%%%%%%%%%%%%%%%%%%%%%%%%%%%%
$$
V = \cmatrix {1 & x & x^2 \cr
1 & \sigma(x) & \sigma(x^2) \cr
1 & \sigma^2(x) & \sigma^2(x^2) \cr}
= \cmatrix {1 & x_0 & x_0^2 \cr
1 & x_1 & x_1^2 \cr
1 & x_2 & x_2^2 \cr}
\quad \hbox {avec} \quad x_i = \sigma^i(x)
.$$
%%%%%%%%%%%%%%%%%%%%%%%%%%%%%%%%%%%%%%%%%
On pose $\varepsilon = \sigma(x) - x$.
Alors, $\det(V) = (x_1-x_0)(x_2-x_1)(x_2-x_0)$ est \ivz:

\snic{\det(V) =
\big(\sigma(x)-x\big) \cdot \sigma\big(\sigma(x)-x\big) \cdot \sigma^2\big(x-\sigma(x)\big)
= - \varepsilon \sigma(\varepsilon) \sigma^2(\varepsilon).}

%\sni
Pour $y \in \gA$, on cherche à écrire $y=\lambda_0+\lambda_1x+\lambda_2x^2$
avec les $\lambda_i\in \gk$. On a alors \ncrtz:
%%%%%%%%%%%%%%%%%%%%%%%%%%%%%%%%%%%%%%%%%
$$
\cmatrix {y\cr \sigma(y)\cr \sigma^2(y)\cr} =
\cmatrix {1 & x & x^2 \cr
1 & \sigma(x) & \sigma(x^2) \cr
1 & \sigma^2(x) & \sigma^2(x^2) \cr}
\cmatrix {\lambda_0\cr \lambda_1\cr \lambda_2\cr}
.$$
%%%%%%%%%%%%%%%%%%%%%%%%%%%%%%%%%%%%%%%%%
Or  le \sli ci-dessus a  une et une seule solution  dans~$\gA$.
Puisque la solution est unique, $\sigma(\lambda_i) = \lambda_i$, i.e.
$\lambda_i \in \gk$ ($i=0,1,2$). \\
Finalement, $(1, x, x^2)$ est bien une~$\gk$-base
de~$\gA$.
\eoe

%:  subsec{Correspondance galoisienne, faits évidents}
\subsec{Correspondance galoisienne, faits évidents}

\vspace{3pt}
Ceci peut être considéré comme une reprise de la proposition~\ref{defiCorGal}.
\vspace{-1pt}

%:     Fact{lemGal1bis}-------
\begin{fact}
\label{lemGal1bis} \emph{(Correspondance galoisienne, faits évidents)}\\
On considère un groupe fini $G$ d'\autos d'un anneau~$\gA$.
On utilise les notations définies en~\ref{NOTAStStp}: en particulier,
$\Ae H=\Fix_\gA(H)$ pour un sous-groupe $H$ de $G$. On pose $\gk=\gA^{G}$.
\begin{enumerate}
\item  Si $H\subseteq H'$ sont deux sous-groupes de $G$, alors $\Ae H\supseteq\Ae {H'}$, et si $H$ est le sous-groupe engendré par $H_1\cup H_2$, alors
$\Ae {H}= \Ae {H_1}\cap \Ae {H_2}$.
\item
$H\subseteq \Stp(\Ae H)$ pour tout sous-groupe $H$ de $G$.
\item Si $\sigma\in G$ et $H$ est un sous-groupe de $G$
alors

\snic{\sigma(\Ae H)=\Ae {\sigma H\sigma^{-1}}.}
\item  Si~$\gC \subseteq \gC'$ sont deux sous-\klgs de~$\gA$, alors
$\Stp(\gC)\supseteq\Stp(\gC')$,
et si~$\gC$ est la sous-\klg engendrée par~$\gC_1\cup \gC_2$,
alors 

\snic{\Stp(\gC)= \Stp(\gC_1)\cap \Stp(\gC_2).}

\item $\gC\subseteq \Ae {\Stp(\gC)}$ pour toute sous-\klgz~$\gC$ de~$\gA$.
\item Après un aller-retour-aller on retombe sur l'arrivée du premier aller:

\snic{    \Ae H = \Ae {\Stp(\Ae H)} \; \et \; \Stp(\gC)= \Stp\big(\Ae {\Stp(\gC)}\big).}
\end{enumerate}
%-----------------begin enum------------------
%-----------------end enum------------------
\end{fact}
%--- end-fact----------------------------------------
%-----------------begin proof------------------
\begin{proof}
Le dernier point est une conséquence directe des précédents,
qui sont \imdsz. Comme dans toutes les \gui{dualités} de ce type.
\end{proof}
%-----------------end proof------------------

%:  subsec{Une \dfn naturelle}
\subsec{Une \dfn naturelle}

Notons $\cG=\cG_G$ l'ensemble des sous-groupes finis (i.e., détachables) de $G$,
et~$\cA=\cA_G$ l'ensemble des sous-anneaux de~$\gA$ qui sont de la forme $\Fix(H)$
pour un $H\in\cG$. Considérons
les restrictions de $\Fix$
et $\Stp$ aux ensembles $\cG$
et $\cA$. Nous sommes intéressés pour déterminer
dans quelles conditions on obtient ainsi deux bijections réciproques
l'une de l'autre entre $\cG$ et $\cA$, et à donner une \carn
agréable des sous-\algs appartenant à $\cA$.

Dans le cas où~$\gA$ est un corps discret, la théorie d'Artin
montre que l'on  se trouve dans une situation galoisienne classique:~$\gA$ est une extension galoisienne du sous-corps~$\gk=\Ae G$,
 $G$ est le groupe de Galois de cette extension et $\cA$ est l'ensemble de
 toutes les sous-extensions \stfes de~$\gA$.

Cette théorie \gui{d'Artin-Galois} a ensuite été \gnee à un anneau
 commutatif arbitraire~$\gA$, à condition d'imposer certaines contraintes au groupe $G$ et aux sous-\klgs de~$\gA$.

En fait, on veut que la notion correspondante d'\aG soit suffisamment stable.
En particulier, on souhaite que lorsque l'on  remplace~$\gk$ par un quotient
non trivial~$\gk\sur\fa$ et
$\gA$ par~$\gA\sur{\fa\gA}$, on maintienne la notion d'\aGz. Il ne faut donc pas
que deux \autos de~$\gA$ présents dans $G$ puissent devenir un seul \auto
en passant au quotient.

Ceci conduit à la \dfn suivante.

%:     Definition{defaG}----defi alg galoisienne-
\begin{definition}\label{defaG}
\emph{(Applications bien séparées, \autos séparants, \aGsz)}
\index{algèbre!galoisienne}
\index{bien séparées!applications ---}
\begin{enumerate}
\item
Deux applications  $\sigma$, $\sigma'$ d'un ensemble $E$ dans
un anneau~$\gA$ sont dites \emph{bien
séparées}  si

\snic{\gen{\,\sigma(x)-\sigma'(x)\;;\;x\in E\,}_{\!\gA}=\gen{1}.}

\item Un \auto $\tau$ de~$\gA$ est dit \ixg{séparant}{autom@\auto ---}{separant} s'il est bien séparé de  $\Id_\gA$.

\item Un groupe fini $G$
qui opère sur~$\gA$ est dit \ixg{séparant}{groupe --- d'\autosz}{separant},
si les \elts $\sigma\neq 1_G$ de
$G$ sont séparants (il revient au même de dire que toute paire d'\elts distincts de $G$ donne deux \autos  bien séparés).
\\
On dira aussi que $G$ opère \emph{de façon séparante} sur~$\gA$.
\item Une \emph{\aGz} est par définition un triplet
$(\gk,\gA,G)$,
où~$\gA$ est un anneau,
$G$  est un groupe fini opérant sur~$\gA$ de façon séparante, et~$\gk=\Fix(G)$.
\end{enumerate}
\end{definition}

\comms 
\\
1) Pour ce qui concerne la \dfn d'une \aGz,
nous n'avons pas voulu interdire un groupe fini opérant sur l'anneau
trivial%
\footnote{L'unique \auto  de l'anneau trivial est séparant,
et tout groupe fini opère sur l'anneau trivial de manière à en faire
une \aGz.},
 et en conséquence nous ne définissons pas $G$
comme un groupe d'\autos de~$\gA$, mais comme un groupe fini
opérant sur~$\gA$.
En fait, la \dfn implique que $G$ opère toujours de manière fidèle sur $\gA$ (et donc s'identifie à un sous-groupe de $\Aut(\gA)$) sauf dans le cas où l'anneau est trivial.
 Ceci présente plusieurs avantages.
\\
D'une part, une \aG reste galoisienne, \emph{avec le même groupe~$G$}, pour toute
\edsz: il arrive que l'on  ne sache pas si une \edsz~$\gk\to\gk'$, qui débarque au cours d'une \demz, est triviale ou non. \\
D'autre part, le fait de ne pas changer de groupe
est de toute manière plus confortable, pour n'importe quelle \edsz.
 
2) Nous avons imposé la condition $\gk\subseteq \gA$, qui n'est pas dans
le style catégorique usuel.  \Llec qui le désire pourra rétablir une \dfn plus
catégorique, en disant que le morphisme  $\gk\to\gA$ établit un \iso entre
$\gk$ et $\gA^G$. Cela serait parfois \ncrz, par  exemple dans le point~\emph{2} du fait~\ref{factEdsAG}.
\eoe

\penalty-2500 
\exls \\
1) Soit~$\gL/\gK$  une extension galoisienne de \cdisz.\\
Alors  $\big(\gK,\gL,\Gal(\gL/\gK)\big)$ est une \aGz.

2) Nous démontrerons plus loin (\thrf{thAduAGB}) que
pour un \polu $f\in\kT$ \splz,
le triplet $(\gk,\Adu_{\gk,f},\Sn)$ est une \aGz.

3) Un \auto $\sigma$ d'un \alo $\gA$ est séparant \ssi il existe un $x\in\gA$
tel que $x-\sigma(x)$ est \ivz. 
\eoe

\medskip
Les notions d'\auto séparant et d'\aG ont été mises au point de manière à vérifier les faits fondamentaux suivants.

%:     Fact{factEdsAG}
\begin{fact}\label{factEdsAG}~
\begin{enumerate}
\item Un \auto séparant $\sigma$ d'un anneau~$\gA$ fournit par \eds $\rho:\gA\to\gB$ un \auto séparant $\rho\ist(\sigma)$ de~$\gB$.
\item Si $(\gk,\gA,G)$ est une \aG  et si $\rho:\gk\to\gk'$ est un \homo d'anneaux,
alors $(\gk',\rho\ist(\gA),G)$ est une \aGz.
\end{enumerate}

\end{fact}
\begin{proof}
Le point \emph{1}, ainsi que le point \emph{2} dans le cas d'une \eds par \lonz, sont faciles et laissés \alecz. 
\\
La \dem du cas \gnl pour le point \emph{2} devra attendre
le \thref{corAGextsca}.
\end{proof}

%:     PrincipeLocGlob{plcc.aGs}
\begin{plcc}\label{plcc.aGs}\relax \emph{(Algèbres galoisiennes)}
\\
Soit $G$ un groupe fini opérant sur une \klgz~$\gA$ avec $\gk\subseteq \gA$. 
\\
 Soient   $S_1$, $\dots$, $S_n$  des \moco de~$\gk$.
\\
 Alors, $(\gk,\gA,G)$ est une \aG \ssi chaque \linebreak 
 triplet $(\gk_{S_i},\gA_{S_i},G)$
 est une \aGz.
\end{plcc}

\facile

%:  subsec{Lemme de Dedekind}
\subsec{Lemme de Dedekind}
\label{subsecLDAC}\index{Dedekind!Lemme de ---}

Soit~$\gA$ un anneau commutatif. Considérons   l'\Alg puissance $m$-ième~$\Ae m$. Ses \elts seront vus comme des
vecteurs colonnes et les lois sont les lois produit:

\snic{\cmatrix{a_1 \cr\vdots\cr  a_m}\star\cmatrix{b_1 \cr\vdots\cr  b_m} =
\cmatrix{a_1\star b_1 \cr\vdots\cr  a_m\star b_m},\quad a \cmatrix{a_1 \cr\vdots\cr  a_m}=\cmatrix{aa_1 \cr\vdots\cr  aa_m}.}

%\smallskip
%L'anneau  $\Ae m$ est une \Alg via le morphisme diagonal $x\mapsto \tra{[\,x \;\cdots\; x\,]}$.% de~$\gA$ dans $\Ae m$.

%--- Lemma{lemDA}------------------------------------------
\begin{lemma}
\label{lemDA}
Soit $C$ une partie finie de $\Ae m$ qui \gui{sépare les lignes}:
i.e.,~$\gen{x_i- x_j \,;\,x\in C}_\gA=\gen{1}$ (pour $i\neq j\in\lrbm$).
Alors, l'\Alg engendrée par $C$ est égale à $\Ae m$.
\end{lemma}
%--- end-lemma-----------------------------------------
%-----------------begin proof------------------
\begin{proof}
La remarque fondamentale est que
dans le \Amo engendré par $1_{\Ae m}$ et  $x=\tra{[\,x_1 \,\cdots\, x_m\,]}$ il
y a les vecteurs

\snic{x-x_2 \,1_{\Ae m}=\tra{[\,x_1-x_2\;0\;* \,\cdots\, *\,]}$
et 
$-x+x_1 \,1_{\Ae m}=\tra{[\,0\;x_1-x_2\;* \,\cdots\, *\,]}.}

%\sni
Donc lorsque l'on  suppose que
l'idéal  engendré par les $x_1-x_2$ contient $1$, cela implique
que dans le \Amo engendré par $C$ il y a un vecteur $g^{1,2}$ du
type $\tra{[\,1\;0\;g^{1,2}_3 \,\cdots\,  g^{1,2}_m\,]}$ et un vecteur $g^{2,1}$ du
type $\tra{[\,0\;1\;g^{2,1}_3 \,\cdots\,  g^{2,1}_m\,]}$. Même chose en
remplaçant $1$ et $2$ par deux entiers $i\neq j\in\lrbm$.
\\
On en déduit que $\tra{[\,1\;0\;0 \,\cdots\, 0\,]}=
g^{1,2}\cdot g^{1,3}\cdots g^{1,m}$, est dans l'\Alg  engendrée par $C$.
De même, chaque vecteur de la base canonique
de $\Ae m$ sera dans l'\Alg  engendrée  par $C$.
On  obtient en fait que $\Ae m$ est l'image d'une matrice dont les colonnes
sont les produits d'au plus $m$ colonnes dans~$C$.
\end{proof}
%-----------------end proof------------------
%--- Notations{notas2.1}----------
\begin{notations}
\label{notas2.1}  \emph{(Contexte du lemme de Dedekind)}
%-----------------begin item------------------
\begin{enumerate}\itemsep0pt
\item [--] $\gA$ est un anneau commutatif.
\item [--] $(M,\cdot,1)$ est un \moz.\index{monoide@monoïde!(lemme de Dedekind)}
\item [--] $\tau=(\tau_1,\tau_2,\ldots,\tau_m)$ est une liste de $m$ \homosz,
deux à deux bien séparés, de $(M,\cdot,1)$ dans $(\gA,\cdot,1)$.
\item [--] Pour $z\in M$ on note $\tau(z)$ l'\elt de $\Ae m$ défini par

\snic{\tau(z)=\tra{[\,\tau_1(z) \,\cdots\, \tau_m(z)\,]}.}
\end{enumerate}
%-----------------end item------------------
\end{notations}
%--- end-notation-----------------------------------------

%:     Theorem{thDA}  Lemme de Dedekind 1
\begin{theorem}
\label{thDA}
\emph{(Lemme de Dedekind)}\\
Avec les notations~\ref{notas2.1} il existe $y_1$, \ldots, $y_r\in M$ tels que la matrice

\snic{[\,\tau(y_1) \mid \cdots \mid \tau(y_r)\,]= \big(\tau_i(y_j)\big)_{i\in\lrbm,j\in\lrbr}}

%\sni
 est surjective. En particulier,
$\tau_1$, \ldots, $\tau_m$ sont~$\gA$-\lint indépendants.
\end{theorem}
%--- end-theorem-----------------------------------------
%
\begin{proof}
Se déduit du lemme \ref{lemDA} en remarquant que,
puisque $\tau(xy) = \tau(x)\tau(y)$, l'\Alg
engendrée par les $\tau(x)$ coïncide avec
le \Amo engendré par les~$\tau(x)$.
\end{proof}

\rems\\
 1) Posons $F=\big(\tau_i(y_j)\big)_{ij}\in\Ae{m\times r}$.
L'indépendance \lin des lignes signifie que $\cD_m(F)$
est fidèle, tandis que la surjectivité de $F$ signifie que $\cD_m(F)$
contient $1$. Parfois, le lemme de Dedekind est appelé \gui{lemme d'indépendance des \homosz}, lorsqu'on a en vue le cas où $\gA$ est un \cdiz. En fait, c'est seulement lorsque $\gA$ est un anneau \zed que 
l'on peut déduire \gui{$\cD_m(F)=\gen{1}$} de \gui{$\cD_m(F)$ fidèle}.

2) L'entier $r$  peut être contrôlé à partir des
données du problème.
\eoe

%%%%%%%%%%%%%%%%%%%%%%%%%%%%%%%%%%%%%%%%%%%%%%%%%%%%%%%%%%%%%%%%%%%%%%%%%%%
%: subsec{\Tho d'Artin
%\penalty-2500
\subsec{\Tho d'Artin et premières conséquences}

%une \aG $(\gk,\gA,G)$ %(en fait $(\gA,G)$ est suffisant)
%est donnée par un anneau~$\gA$ non trivial
%et un groupe fini $G$ d'\autos séparants
%de~$\gA$,~$\gk$ étant le sous-anneau des points fixes de $G$.

%--- Notation{notaAGAL}-------
\begin{definota}
%\label{notadebut}\label{notadebutADU}\label{defAGAL}
\label{NotaAGAL}
Soit une \klgz~$\gA$ avec $\gk\subseteq \gA$.
%-----------------begin item------------------
%
\begin{enumerate}
\item On peut munir le \kmo $\Lin_\gk(\gA,\gA)$ d'une structure de \Amo  par
la loi externe
$$\preskip.2em \postskip.2em \ndsp
(y,\varphi)\mapsto \big(x\mapsto y\varphi(x)\big), \quad\gA\times
\Lin_\gk(\gA,\gA)\to\Lin_\gk(\gA,\gA)\,.
$$
On note alors $\LIN_\gk(\gA,\gA)$ ce \Amoz.
\end{enumerate}
%:HHH  finiment enumere remplacé par fini, avec contorsion
Soit $G=\so{\sigma_1=\Id,\sigma_2, \ldots ,\sigma_n}$ un groupe fini opérant (par~$\gk$-\autosz) sur~$\gA$.
\begin{enumerate}\setcounter{enumi}{1}
\item
L'\Ali $\iota_G:\prod_{\sigma\in G}\gA \to \LIN_\gk(\gA,\gA)$
est définie par 
$$\preskip.2em \postskip.2em \ndsp
\iota_G \big((a_\sigma)_{\sigma\in G}\big)=\sum_{\sigma\in G} a_\sigma
\sigma\,. 
$$
\item
L'\kli
$\psi_G:\env\gk\gA \to \prod_{\sigma\in G}\gA$  est définie par
$$\preskip.2em \postskip.2em 
\psi_G(a\otimes b)=\big(a\sigma(b)\big)_{\sigma\in G}\,.
$$
C'est un \homo d'\Algs
(à gauche).
\perso{Il y a aussi une action de $G$ qui est conservée dans l'\isoz, à
préciser si vraiment \ncrz.}
\end{enumerate}
%-----------------end enum------------------
\end{definota}
%--- end-notation-----------------------------------------

%:     Fact{factdefAGAL}
\begin{fact}\label{factdefAGAL} Avec les notations ci-dessus, et la structure à gauche pour\linebreak  le~\Amo $\env\gk\gA$, on a les résultats suivants.
\begin{enumerate}
\item Dire que $\iota_G$ est un \iso signifie que $\LIN_\gk(\gA,\gA)$ est un \Amo libre dont $G$ est une base.
\item Si~$\gA$ est \ste de rang constant sur~$\gk$, dire
que $\env{\gk}{\gA}$ est un~\Amo libre
de rang fini signifie que~$\gA$ se diagonalise
elle-même.
\item Dire que $\psi_G$ est un \iso signifie \prmt la chose suivante.
Le \Amo $\env{\gk}{\gA}$ est  libre
de rang $\#{G}$, avec une base $\cB$ telle que,
après \eds de~$\gk$ à~$\gA$, l'\ali $\mu_{\gA,a}$,
qui est devenue $\mu_{\env{\gk}{\gA},1\te a}$,
est maintenant diagonale sur la base $\cB$,
avec pour matrice~$\Diag\big(\sigma_1(a),\sigma_2(a),\ldots,\sigma_n(a)\big)$,
ceci pour n'importe quel $a\in\gA$.\perso{Le 3 mérite d'être expliqué?
}
\end{enumerate}
 \end{fact}

%%%%%%%%%%%%%%%%%%%%%%%%%%%%%%%%%%%%%%%%%%%%%%%%%%%%%%%%%%%%%%%%%%%%%%%%%%%

%:    lemma} \label{lemArtin}
\begin{lemma}\label{lemArtin}~\\
Soit $G=\so{\sigma_1=\Id,\sigma_2, \ldots ,\sigma_n}$ un groupe fini
opérant sur un anneau~$\gA$ et~$\gk=\Ae G$. 
Pour $y\in\gA$ notons $y\sta$ l'\elt de $\Asta$ défini par
$x\mapsto \Tr_G(xy)$. \Propeq
\begin{enumerate}
\item $(\gk,\gA,G)$ est une \aGz.
\item  Il existe  $x_1$, \ldots, $x_r$, $y_1$, \ldots, $y_r$ dans~$\gA$
tels que pour tout $\sigma\in G$ on ait
\begin{equation}\label{thAeq4}
\som_{i=1}^r x_i \sigma(y_i)= \formule{1\;\;\hbox{si}\;\sigma=\Id\\
0\;\;\hbox{sinon}\,.}
\end{equation}
\end{enumerate}
Dans ce cas on a les résultats suivants. 
\begin{enumerate} \setcounter {enumi}{2}
\item 
Pour $z\in\gA$, on a $z=\sum_{i=1}^r \Tr_G(zy_i)\, x_i = 
\sum_{i=1}^r \Tr_G(zx_i)\, y_i$. 
\\
Autrement dit, $\gA$ est un \kmo \ptf
et 
$$\preskip.2em \postskip.2em 
\big((\xr),(y_1\sta,\ldots,y_r\sta)\big)\;\hbox{  et  }\;\big((\yr),(x_1\sta,\ldots,x_r\sta)\big) 
$$
sont des \sycsz.
\item 
La forme
$\Tr_G: \gA \rightarrow \gk$ est dualisante, surjective.
\item 
Pour $\sigma \in G$, on pose $\varepsilon_\sigma = \sum_i \sigma(x_i) \otimes
y_i \in \env\gk\gA $. Alors,
$(\varepsilon_\sigma)_{\sigma \in G}$ est une~$\gA$-base \gui {à gauche}
de $\env\gk\gA$. De plus, pour $a$, $b \in \gA$, on a
$$\preskip.2em \postskip.2em \ndsp
b \otimes a = \sum_\sigma  b\sigma(a) \varepsilon_\sigma
,
$$
et l'image de cette base $(\varepsilon_\sigma)_\sigma$ par $\psi_G :
\env\gk\gA \to \prod_{\tau \in G}\gA$ est la~$\gA$-base canonique
$(e_\sigma)_{\sigma\in G}$ de $\prod_{\tau \in G}\gA$.
En conséquence, $\psi_G$ est un \iso d'\Algsz.

\end{enumerate}
\end{lemma}

%%%%%%%%%%%%%%%%%%%%%%%%%%%%%%%%%%%%%%%

\begin {proof}
 \emph{1} $\Rightarrow$ \emph{2.} 
D'après le lemme de Dedekind, il existe
un entier $r$ et des \eltsz\hbox{ $x_1$, \ldots, $x_r$, $y_1$, \ldots, $y_r\in\gA$} 
tels que 
%-----------------begin $$----------------
$$ \sum_{i=1}^r \,x_i\,\cmatrix{\sigma_1(y_i) \cr\sigma_2(y_i) \cr\vdots 
\cr\sigma_n(y_i)}=
\cmatrix{1\cr0\cr\vdots\cr0},
$$
%-----------------end $$------------------
c'est-à-dire exactement, pour $\sigma \in G$, les \eqnsz~(\ref{thAeq4}).

 \emph{2} $\Rightarrow$ \emph{1.}
Pour $\sigma \ne \Id$, on a $\sum_{i=1}^r x_i \big(y_i - \sigma(y_i)\big) = 1$, ce qui
prouve que $\sigma$ est séparant.

 \emph{3.}
Pour $z \in \gA$, on a les \egts
\[ \arraycolsep2pt
\begin{array}{rccccc} 
\sum_{i=1}^r \Tr_G(zy_i) \, x_i&=&\sum_{i=1}^r \sum_{j=1}^n \sigma_j(zy_i) x_i 
&=& \\[2mm] 
 \sum_{j=1}^n \sigma_j(z) \big( \sum_{i=1}^r \sigma_j(y_i)x_i \big)&=&
\sigma_1(z)\cdot 1+\sum_{j=2}^n\sigma_j(z)\cdot 0&=&z. 
\end{array}
\]

 \emph{3} $\Rightarrow$ \emph{4.} D'après le point \emph{1} du \thref{factCarDua}.
%La surjectivité en résulte. 
%:HHH  la dem alternative a été repoussée en remarque
%Voici une \dem alternative.
%\\
%Pour $z = 1$, $1 = \sum_{i=1}^r t_i x_i$ avec
%$t_i = \Tr_G(y_i) \in \Tr_G(\gA) \subseteq \gk$. Introduisons le \pol 
%\gui {normique} $N(T_1, \ldots, T_r)$:
%
%\snic{N(T_1, \ldots, T_r) =\rN_G \big(\som_{i=1}^r T_i x_i\big) = 
%\prod\nolimits_{\sigma \in G} \big(T_1\sigma(x_1) + \cdots + T_r\sigma(x_r)\big).
%}
%
%\sni
%C'est un \pol\hmg de degré $n \ge 1$, invariant par $G$, donc à \coes dans
%$\gk$: $N(\und{T}) = \sum_{|\alpha|=n} \lambda_\alpha \und{T}^\alpha$ avec
%$\lambda_\alpha \in \gk$. En conséquence, pour $u_1$, \ldots, $u_r \in \gk$,
%on a $N(u_1, \ldots, u_r) \in \gk u_1 + \cdots + \gk u_r$.
%En particulier: 
%
%\snic{1 =\rN_G(1)= \rN_G\big(\som_{i=1}^r t_i x_i\big)= N(t_1, \ldots, t_r) 
%\in \gk t_1 + \cdots + \gk t_r \subseteq \Tr_G(\gA).}

\sni
\emph {5.}
On a $\psi_G(\varepsilon_\sigma) = \big(\sum_i
\sigma(x_i)\tau(y_i)\big)_\tau = e_\sigma$. Montrons maintenant
l'\egt relative à $b\otimes a$. Vue la structure choisie de \Amo à
gauche, on peut supposer $b = 1$. Alors:

\snic {\arraycolsep2pt
\begin {array} {rcl}
\sum_\sigma  \sigma(a) \varepsilon_\sigma &=&
\sum_\sigma  \sigma(a) \sum_i \sigma(x_i) \otimes y_i =
\sum_i \Tr_G(ax_i) \otimes y_i
\\[1mm] 
&=& \sum_i 1 \otimes \Tr_G(ax_i)y_i = 
     1 \otimes \sum_i \Tr_G(ax_i)y_i = 1 \otimes a.
\\
\end {array}
}

Ceci montre que $(\varepsilon_\sigma)_\sigma$ est un \sgr du \Amo $\env\gk\gA$.
Comme son image par $\psi_G$ est la~$\gA$-base canonique de $\prod_{\tau\in
G}\gA$, ce \sys est libre sur~$\gA$. Le reste en découle.
\end {proof}

%:HHH remarque avec demo alternative
\rem Voici une \dem alternative de la surjectivité de la trace (point \emph{4}). Pour $z = 1$, $1 = \sum_{i=1}^r t_i x_i$ avec
$t_i = \Tr_G(y_i) \in \Tr_G(\gA) \subseteq \gk$. Introduisons le \pol 
\gui {normique} $N(T_1, \ldots, T_r)$:

\snic{N(T_1, \ldots, T_r) =\rN_G \big(\som_{i=1}^r T_i x_i\big) = 
\prod\nolimits_{\sigma \in G} \big(T_1\sigma(x_1) + \cdots + T_r\sigma(x_r)\big).
}

%\sni
C'est un \pol\hmg de degré $n \ge 1$, invariant par $G$, donc à \coes dans
$\gk$: $N(\und{T}) = \sum_{|\alpha|=n} \lambda_\alpha \und{T}^\alpha$ avec
$\lambda_\alpha \in \gk$. En conséquence, pour $u_1$, \ldots, $u_r \in \gk$,
on a $N(u_1, \ldots, u_r) \in \gk u_1 + \cdots + \gk u_r$.
En particulier: 

\smallskip \centerline{~~~~{\small $1 =\rN_G(1)= \rN_G\big(\som_{i=1}^r t_i x_i\big)= N(t_1, \ldots, t_r) 
\in \gk t_1 + \cdots + \gk t_r \subseteq \Tr_G(\gA).$} \eoe}

%%%%%%%%%%%%%%%%%%%%%%%%%%%%%%%%%%%%%%%%%%%%%%%%%%%%%%%%%%%%%%%%%%%%%%%%%%%%

%:     Theorem{thA}-  Artin  --
\begin{theorem}
\label{thA}\emph{(\Tho d'Artin, version  \aGsz)}\index{Artin!Théorème d'---} \\
Soit $(\gk,\gA,G)$  une \aG (notations \ref{NotaAGAL}).
%-----------------begin enum------------------
\begin{enumerate}
\item Le \kmoz~$\gA$ est \prc $\#G$, et~$\gk$ est facteur
direct dans~$\gA$.
\item
Il existe $x_1$, \ldots, $x_r$ et $y_1$, \ldots, $y_r$ tels que pour tous
$\sigma$, $\tau\in G$ on ait
\begin{equation}\label{thAeq5}
\forall\sigma,\tau \in G\qquad \som_{i=1}^r \tau(x_i) \sigma(y_i)= 
\formule{1\;\;\mathrm{si}\;\sigma=\tau\\
0\;\;\mathrm{sinon}.}
\end{equation}
\item La forme $\Tr_G$ est dualisante.
\item
L'application $\psi_G : \env{\gk}{\gA} \to \prod_{\sigma\in G} \gA$ est un \iso
d'\Algsz.  En particulier,~$\gA$ se diagonalise elle-même.
\item 
\begin{enumerate}
\item $\rC{G}(x)(T)=\rC{\gA\sur\gk}(x)(T)$,
 $\Tr_G=\Tr\iAk $ et
 $\rN_G=\rN\iAk $, 
\item $\gA$ est \ste sur~$\gk$. 
\end{enumerate}
\item Si~$\gA$ est un \cdiz, c'est une extension galoisienne de~$\gk$,
et l'on~a~$G=\Gal(\gA\sur\gk)$.
\end{enumerate}
%-----------------end enum------------------
\end{theorem}
%--- end-theorem-----------------------------------------
%-----------------begin proof------------------
\begin{proof}
Dans cette preuve, pour $x\in\gA$, on note $\Tr(x)=\Tr_G(x)$, et
$x\sta$ est la forme~$\gk$-\lin $z\mapsto\Tr(zx)$.

%:HHH  phrase mieux ecrite
Le lemme \ref{lemArtin} prouve les points \emph{1} (mis
à part la question du rang),  \emph{3} et  \emph{4.} Il
prouve aussi le point \emph {2}, car~(\ref{thAeq5}) résulte clairement de (\ref{thAeq4}).
\\
Notons que~$\gk$ est en facteur direct dans~$\gA$
d'après le lemme~\ref{lemIRAdu}~\emph{3}\footnote{Ou plus directement, d'après la
surjectivité de la trace (qui résulte du \thref{factCarAste}~\emph{1}).
Soit en effet $x_0 \in \gA$
tel que $\Tr(x_0) = 1$, on a~$\gA = \gk\cdot 1 \oplus \Ker x_0\sta$, car tout $y \in \gA$\linebreak  s'écrit
$y = x_0\sta(y)\cdot 1 + (y - x_0\sta(y)\cdot 1)$ avec $y - x_0\sta(y)\cdot 1\in\Ker x_0\sta$.}

Voyons que~$\gA$ est bien de rang constant $n$.
Le point \emph{4} montre que, après \eds de~$\gk$ à~$\gA$, le \kmoz~$\gA$
devient libre de rang $\#{G}$.
Ainsi~$\gA$ est bien de rang constant $n$ sur~$\gk$: le \polmu du \kmoz~$\gA$  \gui{ne change pas} par l'\edsz%
\footnote{En fait, ses \coes sont transformés en \gui{eux-mêmes}, vus dans~$\gA$.}~$\gk\to\gA$ (injective), il est donc lui-même égal à $T^n$.

\emph{5a.} (et donc \emph{5b}) Puisque $\psi_G$ est un \iso d'\Algs (point \emph{4}),~$\gA$ se diagonalise elle-même.
On déduit  alors  du fait~\ref{factdefAGAL} point~\emph{3},
l'\egt 

\snic{\rC{G}(x)(T)=\rC{\gA/\gk}(x)(T).}

%\sni
Elle est vraie pour les \pols vus dans $\AT$,
donc aussi dans $\kT$.

\emph{6.} Tout d'abord, l'anneau~$\gk$ est \zed d'après le lemme~\ref{lemZrZr2},
c'est donc un \cdiz, car il est connexe et réduit. L'extension est %clairement 
étale. Elle est normale, car tout
$x\in\gA$ annule %le \pol 
$\rC G(x)(T)$, et  ce \pol se décompose en produit de
facteurs \lins dans $\AT$.
%
%\sni
%Montrons enfin que $\iota_G$ est un \isoz.
%
\end{proof}
%
%-----------------end proof------------------

\rem Le calcul qui suit peut éclairer les choses, bien qu'il n'ait pas été \ncrz.
\\
On note que d'après le point \emph{3} du lemme~\ref{lemArtin}, le \kmoz~$\gA$ est image de la \mprn
$$\preskip.2em \postskip.4em
P=(p_{ij})_{i,j\in \lrbr}=\big(y_i\sta(x_j)\big)_{i,j\in \lrbr}=
\big(\Tr(y_ix_j)\big)_{i,j\in \lrbr}\,.
$$
Rappelons aussi l'\eqrf{thAeq5}:
${\sum_{i=1}^r \tau(x_i) \sigma(y_i)=
\formule{1\;\;\mathrm{si}\;\sigma=\tau\\
0\;\;\mathrm{sinon}}.
}$

%\sni
Posons alors: 
\[\preskip-.4em \postskip.4em 
\begin{array}{cccccc} 
X&=&\cmatrix{\sigma_1(x_1)&\sigma_1(x_2)&\cdots& \sigma_1(x_r) \cr
\sigma_2(x_1)&\sigma_2(x_2)&\cdots &\sigma_2(x_r) \cr\vdots&\vdots&  &\vdots 
\cr
\sigma_n(x_1)&\sigma_n(x_2)&\cdots &\sigma_n(x_r) }
&\,\hbox{ et }  \\[2.5em] 
Y&=&~\cmatrix{\sigma_1(y_1)&\sigma_1(y_2)&\,\cdots\, &\sigma_1(y_r) \cr
\sigma_2(y_1)&\sigma_2(y_2)&\cdots &\sigma_2(y_r) \cr\vdots&\vdots&  &\vdots 
\cr
\sigma_n(y_1)&\sigma_n(y_2)&\cdots &\sigma_n(y_r) }.& 
\end{array}
\]
D'après l'\eqn (\ref{thAeq5}), on a $X\tra{Y}=\In$ et $P=\tra{Y}X$. 
\\
D'après la proposition \ref{propImProjLib}, ceci signifie que le \kmoz~$\gA$, 
%vu  comme  $\Im P \subseteq\gk^r$, 
devient libre
de rang~$n$, avec pour base les $n$ lignes de $Y$,
après \eds de~$\gk$ à~$\gA$. Autrement dit, le \Amo $\env\gk\gA$, vu comme image de
la matrice~$P$ \gui{à \coes dans~$\gA$} est un sous-\Amo libre de
rang~$n$ de $\Ae r$, en facteur direct.  
\eoe

%%%%%%%%%%%%%%%%%%%%%%%%%%%%%%%%%%%%%%%%%%%%%%%%%%%%%%%%%%%%%%%%%%%%%%%%%%%

%%%%%%%%%%%%%%%%%%%%%%%%%%%%%%%%%%%%%%%%%%%%%%%%%%%%%%%%%%%%%%%%%%%%%%%%%%%
%:     Corollary{corAGlibres}
\begin{corollary}\label{corAGlibres} \emph{(Algèbre galoisienne libre)}\\
Soit $(\gk,\gA,G)$ une \aG libre, \hbox{et $n=\#{G}$}.
Si $\ub= (\bn)$ dans $\gA$, on définit
$M_\ub \in \Mn(\gA)$ par 

\snic{M_\ub = \big(\sigma_i(b_j)\big)_{i,j\in\lrbn}.}

%\sni
Alors, pour deux \syss $\ub$, $\ub'$ de $n$ \elts de~$\gA$ on obtient:

\snic {
\tra M_\ub\, M_{\ub'} = \Tr_G(b_ib'_j)_{i,j\in\lrbn}.
}

%\sni
En conséquence, on obtient les résultats suivants.
\begin{itemize}
\item $\det(M_\ub)^2 = \disc(b_1, \ldots, b_n)$.
\item Le \sys $(\bn)$
est une~$\gk$-base de~$\gA$ \ssi la matrice $M_\ub$ est \ivz.
\item  Dans ce
cas, si $\ub'$ est la base duale de $\ub$ relativement
à la forme bi\lin tracique, alors les matrices $M_\ub$ et $M_{\ub'}$
sont inverses l'une de l'autre.
\end{itemize}    
\end{corollary}

\rem Dans la situation où $\gA$ est un \cdiz, le {lemme de Dedekind} dans sa forme originale affirme que la \gui{matrice de Dedekind} $M_\ub$ est \iv lorsque $(\ub)$ est une base de $\gA$ comme \kevz.\eoe

%%%%%%%%%%%%%%%%%%%%%%%%%%%%%%%%%%%%%%%%%%%%%%%%%%%%%%%%%%%%%%%%%%%%%%%%%%%
%-% ENTRE NOUS
\entrenous{Je pense que l'inversibilité de la matrice des $\sigma_i(x_j)$
était la forme originelle du lemme de Dedekind dans le cas des corps.
Si cela se confirme cela pourrait faire l'objet d'un commentaire: comment
Artin occulte ce fait capital, qui donne concrètement  $G$ comme base
de $\LIN_\gk(\gA,\gA)$, en faisant d'une part une preuve super élégante de
l'indépendance \lin des \homosz, et en concluant d'autre part
que $G$ est une base par un argument de dimension. Cela va bien dans le cas
des corps, au prix d'une perte d'information concrète, mais cela ne va plus du tout dans le cas des anneaux, où indépendance \lin n'implique plus base,
même avec un argument de dimension.
}
%-% Fin ENTRENOUS

%%%%%%%%%%%%%%%%%%%%%%%%%%%%%%%%%%%%%%%%%%%%%%%%%%%%%%%%%%%%%%%%%%%%%%%%%%%

%:     Theorem{corAGextsca}
\begin{theorem}\label{corAGextsca}\emph{(\Eds pour les \aGsz)}\\
Soit $(\gk,\gA,G)$ une \aGz,  $\rho:\gk\to\gk'$ une \alg  et~$\gA'=\rho\ist(\gA)$.
\begin{enumerate}
\item Le groupe $G$ opère de façon naturelle sur~$\gA'$
et $(\gk',\gA',G)$ est une \aGz.
\item La \gui{théorie de Galois} de $(\gk',\gA',G)$ se déduit
par \eds de celle de $(\gk,\gA,G)$, au sens suivant: pour chaque sous-groupe fini $H$ de~$G$,
l'\homo naturel $\rho\ist(\Ae H)\to \gA'^H$ est un \isoz.
\end{enumerate}
\end{theorem}
\begin{proof} \emph{1.}
On voit facilement que $G$ agit sur~$\gA'$ de façon séparante. 
Il reste à montrer que~$\gk'$ est le sous-anneau des \elts
$G$-invariants de~$\gA'$.  \\
Notons $\Tr = \Tr_G$. Nous voyons $\Tr$ comme un
$\gk$-\endo de~$\gA$, qui par \eds donne le $\gk'$-\endo  $\Id_{\gk'} \otimes \Tr$
de $\gA'$.\\
puisque $y$ est~$G$-invariant, on a l'\egt 

\snic{(\Id_{\gk'} \otimes \Tr)(yz) =
y\,(\Id_{\gk'} \otimes \Tr)(z).}

%\sni
En prenant $z_0 = 1_{\gk'} \otimes x_0$, où~$x_0 \in
\gA$ vérifie $\Tr(x_0) = 1$, on obtient l'appartenance souhaitée:

\snic {
y = (\Id_{\gk'} \otimes \Tr) (yz_0) \in
\gk'\otimes_\gk \gk = \gk'
.}

\sni \emph{2.} Résulte du point \emph{1}. En effet, considérons l'\aG $(\Ae H,\gA,H)$ et l'\eds $\varphi:\Ae H\to \gk'\te_\gk \Ae H=\rho\ist(\Ae H)$. On obtient l'\egt $\varphi\ist(\gA)=\gA'$,
d'où l'\aG $\big(\rho\ist(\Ae H),\gA',H\big)$.\\
 Ainsi~$\gA'^H=\rho\ist(\Ae H)$. 
\end{proof}
% 	

%%%%%%%%%%%%%%%%%%%%%%%%%%%%%%%%%%%%%%%%%%%%%%%%%%%%%%%%%%%%%%%%%%%%%%%%%%%
Dans le \tho qui suit, on aurait pu exprimer l'hypothèse en disant que
le groupe fini $G$ opère sur l'anneau $\gA$, et que $\gk $ est un sous-anneau de $\Ae H$. 

%:     Theorem{thAGACar}
\begin{theorem}\label{thAGACar}\emph{(\Carn des \aGsz)}\\
Soit $G$ un groupe fini opérant sur une \klgz~$\gA$ avec $\gk\subseteq \gA$. \Propeq
\begin{enumerate}
\item $(\gk,\gA,G)$ est une \aG (en particulier,~$\gk = \Ae G$).
\item $\gk=\Ae G$, et il existe  $x_1$, \dots, $x_r$, $y_1$, \dots, $y_r$ dans~$\gA$
tels que l' on ait pour tout $\sigma\in G$

\snic{\sum_{i=1}^r x_i \sigma(y_i)=
\formule{1\;\;\mathrm{si}\;\sigma=\Id\\
0\;\;\mathrm{sinon}.}
}
\item $\gk=\Ae G$,~$\gA$ est finie sur~$\gk$, et pour tout \sgr fini $(a_j)_{j\in J}$ de~$\gA$ comme
\kmoz, il existe une famille $(b_j)_{j \in J}$ dans~$\gA$ tel que l'on ait
pour tous $\sigma$, $\tau\in G$

\snic{\sum_{j \in J} \tau(a_j) \sigma(b_j)=
\formule{1\;\;\mathrm{si}\;\sigma=\tau\\
0\;\;\mathrm{sinon}.}
}
\item $\gk=\Ae G$, et
$\psi_G : \env{\gk}{\gA} \to \prod_{\sigma\in G} \gA$ est un \iso d'\Algsz.
\item $\gA$ est \stfe sur $\gk$, et $G$
est une base de~$\LIN_\gk(\gA,\gA)$.
%(i.e. $\iota_G : \prod_{\sigma \in G}\gA \to
%\LIN_\gk(\gA,\gA)$ est un \iso de \Amosz).
%
\end{enumerate}
\end{theorem}

\begin{proof}
On a déjà vu  \emph{1} $\Leftrightarrow$ \emph{2}
et \emph{1} $\Rightarrow$ \emph{4} (lemme~\ref{lemArtin}).
\\
L'implication \emph{3} $\Rightarrow$ \emph{2} est claire.

 \emph{2} $\Rightarrow$ \emph{3.}
On exprime $x_i$ en fonction des $a_j$: $x_i = \sum_j u_{ij} a_j$
avec $u_{ij} \in \gk$. Alors,

\snic{\som_j \sigma\big(\som_i u_{ij} y_i\big) a_j =
\som_{j,i} u_{ij} \sigma(y_i) a_j =
\som_i \sigma(y_i) x_i = \delta_{\Id, \sigma},}

%\sni
d'où le résultat en prenant $b_j = \sum_i u_{ij} y_i$.

 \emph{2} $\Rightarrow$ \emph{5.}
Notons d'abord que si $\varphi \in \LIN_\gk(\gA,\gA)$ s'écrit $\varphi
= \sum_\sigma a_\sigma \sigma$, alors en évaluant en $y_i$, en multipliant
par $\tau(x_i)$ et en sommant sur les $i$, il vient:

\snic{\som_i \varphi(y_i) \tau(x_i) =
\som_{i,\sigma} a_\sigma \sigma(y_i) \tau(x_i) = a_\tau.
}

%\sni
Ceci montre d'une part que $G$ est~$\gA$-libre. D'autre part, cela conduit à
penser que tout $\varphi \in \LIN_\gk(\gA,\gA)$ s'écrit $\varphi =
\sum_\sigma a_\sigma \sigma$ avec $a_\sigma =
\sum_i \varphi(y_i) \sigma(x_i)$. Vérifions-le
en évaluant $\varphi' := \sum_\sigma a_\sigma\sigma$ en $x \in \gA$:

\snac{\varphi'(x) = \som_{i,\sigma} \varphi(y_i) \sigma(x_i) \sigma(x) =
\som_i \Tr_G(x_ix) \varphi(y_i) =
\varphi(\sum_i \Tr_G(x_ix) y_i) = \varphi(x).}

\sni \emph{5} $\Rightarrow$ \emph{2.}
Puisque~$\gk \subseteq \gA$, on a une inclusion $\Asta \hookrightarrow
\LIN_\gk(\gA,\gA)$.  Montrons d'abord que $\Ae G \subseteq \gk$ (on aura alors
l'\egtz). Chaque $\sigma \in G$ est $\Ae G$-\lin donc, puisque $G$ engendre
$\LIN_\gk(\gA,\gA)$ comme \Amoz, chaque \eltz~$\varphi$ de~$\LIN_\gk(\gA,\gA)$ est
$\Ae G$-\linz. En particulier, chaque $\alpha \in \Asta$ \hbox{est $\Ae G$-\linz}.
Soit $\big((x_i), (\alpha_i)\big)$ un \syc du
\kmoz~$\gA$. Comme~$\gA$ est
un \kmo fidèle, d'après la proposition~\ref {propAnnul}, il existe une famille $(z_i)$ dans~$\gA$ telle que $1 = \sum_i
\alpha_i(z_i)$. Alors, si $x \in \Ae G$, on obtient les \egtsz~$x = \sum_i
\alpha_i(z_i) x = \sum_i \alpha_i(z_ix)$: $x$ appartient à~$\gk$.

Montrons ensuite que pour chaque $\alpha \in \Asta$, il existe un unique~$a \in \gA$ tel que $\alpha = \sum_{\sigma \in G} \sigma(a)\sigma$, i.e. tel que $\alpha$ soit la forme~$\gk$-\linz~$x \mapsto
\Tr_G(ax)$. Puisque $G$ est une~$\gA$-base de
$\LIN_\gk(\gA,\gA)$, on a $\alpha = \sum_\sigma a_\sigma \sigma$ avec
des $a_\sigma \in \gA$. Posons $a = a_\Id$. En écrivant, pour $\tau \in G$, $\tau \circ \alpha =
\alpha$, on obtient $\tau(a_\sigma) = a_{\sigma\tau}$, en particulier $a_\tau
= \tau(a)$, d'où l'\egt souhaitée $\alpha = \sum_{\sigma \in G} \sigma(a)\sigma$. En passant, on vient de prouver que l'\kli 

\snic{\gA \to
\Asta$, $a \mapsto \Tr_G(a\bullet)}

%\sni
est un \iso de \kmosz. 
On peut donc définir un
\sys $(y_i)$ par les \egts $\alpha_i = \Tr_G(y_i\bullet)$. Alors, pour $x \in \gA$ on obtient

\snic{
x = \som_i \alpha_i(x) x_i =
\som_{i,\sigma} \sigma(y_ix) x_i =
\som_{\sigma} \big(\sum_i x_i\sigma(y_i)\big)\sigma(x),
}

%\sni
\cad $\Id = \sum_{\sigma} \big(\sum_i x_i\sigma(y_i)\big)\sigma$.
Mais comme $G$ est~$\gA$-libre, l'écriture de $\Id \in G$ est
réduite à $\Id$, donc $\sum_i x_i\sigma(y_i) = 1$
si $\sigma = \Id$, $0$ sinon.

NB. Puisque $\som_i x_iy_i = 1$, on a les \egts

\snic {
\Tr(x) = \sum_i\alpha_i(x_ix) =
\sum_{i,\sigma} \sigma(x_iy_i) \sigma(x) =
\sum_\sigma \sum_i \sigma(x_iy_i) \sigma(x) = \Tr_G(x).
}

\mni \emph{4} $\Rightarrow$ \emph{2.}
Soit $z=\sum_i x_i \otimes y_i$ l'\elt de $\env\gk\gA $
défini par: $\psi_G(z)$ est l'\elt de $\prod_{\sigma \in G} \gA$ dont toutes
les composantes sont nulles, sauf celle d'indice~$\Id$ qui vaut~1. Cela
signifie exactement que $\sum_i x_i\sigma(y_i) = 1$ si $\sigma = \Id$, $0$
sinon.
\end{proof}
%
%%%%%%%%%%%%%%%%%%%%%%%%%%%%%%%%%%%%%%%%%%%%%%%%%%%%%%%%%%%%%%%%%%%%%%%%%%%

Le cas des \aGs libres est décrit dans le corolaire suivant, qui est une conséquence \imde des résultats plus \gnls précédents.

%:     corollary{corAGAlibreCar}
\begin{corollary}\label{corAGAlibreCar}\emph{(\Carn des \aGs libres)}\\
Soit $G$ un groupe fini opérant sur une \klgz~$\gA$ avec $\gk\subseteq \gA$. \\
On suppose que $\gA$ est libre sur $\gk$, de rang $n=\abs G$, avec $\ux=(\xn)$ pour base.
\Propeq
\begin{enumerate}
\item $(\gk,\gA,G)$ est une \aG (en particulier,~$\gk = \Ae G$).
\item La matrice $M_\ux=\big(\sigma_i(x_j)\big)_{i,j\in\lrbn}$ est \iv (on a indexé le groupe~$G$ par $\lrbn$). 
\item La forme $\Tr_G$ est dualisante.
\item $\gk=\Ae G$, et il existe   $y_1$, \dots, $y_n$ dans~$\gA$
tels que l' on ait pour tout $\sigma\in G$

\snic{\sum_{i=1}^n x_i \sigma(y_i)=
\formule{1\;\;\mathrm{si}\;\sigma=\Id\\
0\;\;\mathrm{sinon}.}
}
\item  Le groupe $G$ est une~$\gA$-base de~$\LIN_\gk(\gA,\gA)$.
\item $\gk=\Ae G$, et
$\psi_G : \env{\gk}{\gA} \to \prod_{\sigma\in G} \gA$ est un \iso d'\Algsz.
\end{enumerate}
Dans ce cas on a les résultats suivants. 
\begin{enumerate} \setcounter {enumi}{6}
\item 
Dans les points 4 et 3,
\begin{itemize}
\item on obtient les $y_i$ comme la solution de
$M_\ux . \tra {[\,y_1\;\cdots\; y_n\,]} = \tra {[\,1\;0\cdots\; 0\,]}$,
où $M_\ux$ est définie comme dans le point 2, avec $\sigma_1 = Id$,
\item  $(y_1\sta,\ldots,y_n\sta)$ est la base duale
de $(x_1,\dots,x_n)$.
\end{itemize}
\item Le point 6 peut être précisé comme suit.\\
Pour $\sigma \in G$, on pose $\varepsilon_\sigma = \sum_i \sigma(x_i) \otimes
y_i \in \env\gk\gA $. Alors,
$(\varepsilon_\sigma)_{\sigma \in G}$ est une~$\gA$-base \gui {à gauche}
de $\env\gk\gA$. De plus, pour $a$, $b \in \gA$, on a

\snic {
b \otimes a = \sum_\sigma  b\sigma(a) \varepsilon_\sigma
,}

%\sni
et l'image de cette base $(\varepsilon_\sigma)_\sigma$ par $\psi_G :
\env\gk\gA \to \prod_{\tau \in G}\gA$ est la~$\gA$-base canonique
$(e_\sigma)_{\sigma\in G}$ de $\prod_{\tau \in G}\gA$.
\end{enumerate}
Enfin, on souligne les points suivants, dans lesquels on ne suppose pas que~$\gA$ est libre sur $\Ae G$. 
\begin{itemize}
\item Lorsque $\gA$ est un \cdi (cadre historique du \tho d'Artin), si un groupe $G$ opère fidèlement
sur $\gA$, l'\alg $(\Ae G,\gA,G)$ est toujours galoisienne,
 $\Ae G$ est un \cdi et $\gA$ est libre de rang~$n$ sur $\Ae G$.
%
%\item Lorsque $\gA$ est  \zedz, l'\alg est  galoisienne \ssi
%pour chaque $\sigma\in G$ distinct de $1_G$, $0$ et $1$ sont les seuls \idms de~$\gk$ fixés par
%$\sigma$.
%
\item Lorsque $\gA$ est un \alrdz, l'\alg $(\Ae G,\gA,G)$ est galoisienne \ssi $G$
opère fidèlement sur le corps résiduel~$\gA/\Rad\gA$.  
Dans ce cas, $\Ae G$ est un \alrd et $\gA$ est libre de rang~$n$ sur $\Ae G$.    
\end{itemize} 

\end{corollary}

Naturellement, nous encourageons vivement \llec à donner une \dem plus directe et plus courte
du corolaire précédent. Il est \egmt envisageable de déduire 
les résultats \gnls des résultats particuliers énoncés dans le cas où $\gA$ est un \alrdz, qui pourraient eux-mêmes se déduire du cas des \cdisz.

%%%%%%%%%%%%%%%%%%%%%%%%%%%%%%%%%%%%%%%%%%%%%%%%%%%%%%%%%%%%%%%%%%%%%%%%%%%

%:     Theorem{corDAexplicite2}----     
\begin{theorem} 
\label{corDAexplicite2} \emph{(La correspondance galoisienne pour une \aGz)}
Soit $(\gk,\gA,G)$  une \aG non triviale, et $H$ un sous-groupe fini de $G$.
%:HHH ajout: non triviale, H fini est mis en facteur, notamment pour le point 5
%-----------------begin enum------------------
\begin{enumerate}
\item 
Le triplet $(\Ae H,\gA,H)$ est une \aGz, $\Ae H$ 
est \ste sur~$\gk$, de rang constant $[\Ae H:\gk]=\idg{G:H}$. 

\item Si $ H'\supseteq H$ est un sous-groupe fini de $G$, 
 $\Ae H$ est  \stfe sur~$\Ae {H'}$, de rang constant 
 $[\Ae H:\Ae {H'}]=\idg{H':H}$.
\item  
On a $H= \Stp(\Ae H)$.
\item  
L'application $\Fix_\gA$ 
restreinte aux sous-groupes finis de $G$ est injective.
\item 
Si $H$ est  normal dans $G$, $(\gk,\Ae H,G/H)$ est une \aGz. 
\end{enumerate}
%-----------------end enum------------------
\end{theorem}
%--- end-corollary------------------------------------
%-----------------begin proof------------------
\begin{proof}
\emph{1.}  
 Puisque $H$ est un groupe séparant
d'\autos de~$\gA$, $(\Ae H,\gA,H)$ est une \aGz. Donc~$\gA$ est une $\Ae H$-\alg \stfe
de rang constant  
$\#H$. 
Donc $\Ae H$ est \stfe sur~$\gk$, de rang constant égal à~$\idg{G:H}$
(\thref{propTraptf}). En outre, elle est \ste d'après le fait~\ref{propTrEta}.

\emph{2.} On applique le \thoz~\ref{propTraptf}. 
%obtient $\Ae {H'}$ facteur direct dans $\Ae H$  et l'on conclut.

\emph{3.} L'inclusion  $H\subseteq \Stp(\Ae H)$ est évidente. Soient 
$\sigma\in \Stp(\Ae H)$ et~$H'$ le sous-groupe engendré par $H$ et $\sigma$.
On~a~$\idg{H':H}=[\Ae H:\Ae {H'}]$, or $\Ae H=\Ae {H'}$,
donc~$H'=H$ et~$\sigma\in H$.

\emph{4.} 
Résulte de \emph{3.}

\emph{5.} Tout d'abord, pour $\sigma\in G$, on a $\sigma(\Ae H) = \Ae H$. Si
l'on note $\ov\sigma$ la restriction de $\sigma$ à $\Ae H$, on
obtient un morphisme de groupes $G \to \Aut_\gk(\Ae H)$, 
$\sigma \mapsto \ov\sigma$, dont le noyau est $H$ d'après le
point~\emph {3}. Le groupe quotient $G/H$ se réalise donc comme sous-groupe de $\Aut_\gk(\Ae H)$.
\\
Soient  un $x \in \gA$ vérifiant $\Tr_H(x) =
1$, un \sgr $(a_1, \ldots, a_r)$ de~$\gA$ comme \kmoz, et des \elts $b_1$, \ldots, $b_r$ tels que pour tous $\sigma$, $\tau\in G$ on ait
${\sum_{i=1}^r \tau(a_i) \sigma(b_i)= 
\formule{1\;\;\mathrm{si}\;\sigma=\tau\\
0\;\;\mathrm{sinon}.}
}$.
On définit alors, pour $i \in \lrb{1..r}$, les \elts de $\Ae H$,
${
a'_i = \Tr_H(xa_i), \hbox{ et } b'_i = \Tr_H(b_i).
}$
\\
On vérifie facilement que pour $\sigma\in G$ on a

\snic{\sum_{i=1}^r a'_i \sigma(b'_i)= 
\formule{1\;\;\mathrm{si}\;\sigma \in H\\
0\;\;\mathrm{sinon}.}
}

%\sni
Ainsi, en
application du point~\emph{2} du \thref{thAGACar},
 $(\gk, \Ae H, \G/H)$ est une \aGz. 
\end{proof}
%-----------------end proof------------------

Le \thref{corDAexplicite2} qui précède établit la correspondance galoisienne entre sous-groupes finis de $G$
d'une part et \gui{certaines} sous-\klgs \stes de~$\gA$ d'autre part.
Une correspondance bijective exacte va être établie dans le paragraphe suivant
lorsque~$\gA$ est connexe.

Mais auparavant nous donnons quelques précisions \sulsz.

%:     Proposition{propAdiagAH}
\begin{proposition}\label{propAdiagAH}
Soit $(\gk,\gA,G)$  une \aG  et $H$ un sous-groupe fini de $G$.
\begin{enumerate}
\item $\gA$ diagonalise $\Ae H$.
\item Pour  $b\in\Ae H$, le \polcar de $b$ (sur~$\gk$, dans $\Ae H$)
est donné par

\snic {
\rC{\Ae H\!/\gk}(b)(T) = \prod_{\sigma \in G/H}\big(T - \sigma(b)\big)
.}

%\sni
(Ici $G/H$ désigne un \sys de représentants des classes à gauche,
et l'on note que $\sigma(b)$ ne dépend pas du représentant $\sigma$ choisi% 
%dans une classe
.) 
\end{enumerate} 
\end{proposition}
%--------- fin proposition ---------------------------------------------- 
%
\begin{proof}
Rappelons que~$\gA$ se diagonalise elle-même, comme le montre l'\iso
$\psi_G:\env\gk\gA \to \prod_{\sigma\in G}\gA$. Nous regardons ce produit comme
l'\alg de fonctions  $\cF(G,\gA)$. Il est muni d'une action naturelle de $G$ à gauche, de la façon suivante:

\snic{\sigma\in G,\, w\in \cF(G,\gA): \sigma \cdot w\in \cF(G,\gA)
\hbox{ définie par } \tau \mapsto w(\tau\sigma).  
}

%\sni
De même $G$ agit à gauche sur l'\Alg  $\env\gk\gA = \gA\te_\gk\gA$ via
$\Id \te G$. On vérifie alors que $\psi_G$ est un $G$-morphisme, i.e. 
que pour $\tau \in G$, le diagramme suivant commute:
$$
\preskip-.4em \postskip.2em 
\hspace{.03cm}\xymatrix {
\gA\te_\gk\gA \ar[d]^{\Id \te\tau} \ar[r]^(0.35){\psi_G} &
       \cF(G,\gA)=\prod_{\sigma\in G}\gA \ar[d]^{w \mapsto \tau\cdot w}\\
\gA\te_\gk\gA \ar[r]^(0.35){\psi_G} &\cF(G,\gA)=\prod_{\sigma\in G}\gA\\
} 
$$

\emph{1.}
Considérons
le diagramme commutatif:
$$\preskip.4em \postskip.4em 
\hspace{.05cm}\xymatrix {
\gA\te_\gk\Ae H \ar[d] \ar[r]^(0.33){\varphi_H} &
       \cF(G/H,\gA)=\prod_{\sigma\in G/H}\gA \ar[d]\\
\gA\te_\gk\gA \ar[r]^(0.36){\psi_G}_(0.36){\sim} &\cF(G,\gA)=\prod_{\sigma\in G}\gA\\
}
$$

\`A droite,  la flèche verticale est  injective, 
et elle identifie $\cF(G/H,\gA)$ à la partie~$\cF(G,\gA)^H$ de~$\cF(G,\gA)$ (fonctions constantes sur les classes à gauche de $G$ modulo $H$).  
\\
\`A gauche, la flèche verticale (correspondant
à l'injection $\Ae H \hookrightarrow \gA$) est aussi une injection car
$\Ae H$ est facteur direct dans~$\gA$ en tant que $\Ae H$-module. 
\\
Enfin,  $\varphi_H$ est défini par $a\te b \mapsto
\big(a\sigma(b)\big)_{\sigma \in G/H}$. 
%La notation $\sigma(b)$ a bien un sens: puisque
%$b \in \Ae H$, $\sigma(b)$ ne dépend pas du représentant choisi~$\sigma$.
\\
Alors, $\varphi_H$ est un \iso d'\Algsz. En effet, $\varphi_H$ est injective,
et pour la surjectivité, il suffit de voir que
$(\gA\te_\gk\gA)^{\Id\te H} = \gA\te_\gk\Ae H$: 
%ceci résulte de la
%technique de la trace surjective qui intervient dans la preuve du 
%corolaire~\ref{corAGextsca}.
ceci est donné par le \thref{corAGextsca} pour l'\aG $(\Ae H,\gA,H)$
et l'\eds $\Ae H\hookrightarrow \gA$. 

 \emph{2.} Ceci résulte du point \emph{1} et du lemme suivant.  
\end{proof}
%

%:     Lemma{lemPolCarDiag}
\begin{lemma}\label{lemPolCarDiag}
%:HHH  légèrement reformulé
Soient~$\gA$ et~$\gB$ deux \klgsz,~$\gB$ \stfe de rang constant $n$. On suppose que $\gA$ diagonalise~$\gB$ au moyen d'un \iso  

\snic{\psi:
\gA\te_\gk \gB \lora \Ae n}

%\sni
donné par des \gui{\coosz} notées $\psi_i : \gB \to \gA$.
\\
Alors, pour $b \in \gB$, on a une \egt 

\snic{\rC{\gB/\gk}(b)(T)=\prod_{i=1}^n \big(T - \psi_i(b)\big),}

%\sni
si l'on transforme le membre de gauche (qui est un \elt de $\kT$)  en un
\elt de~$\gA[T]$ via~$\gk \to \gA$.
\end{lemma}
%--------- fin lemma ---------------------------------------------- 
%
\begin{proof}
Immédiat d'après le calcul du \polcar d'un \elt dans une \alg diagonale.
\end{proof}
%

%%%%%%%%%%%%%%%%%%%%%%%%%%%%%%%%%%%%%%%%%%%%%%%%%%%%%%%%%%%%%%%%%%%%%%%%%%%
%: subsect La correspondance galoisienne dans le cas connexe
\subsect{La correspondance galoisienne dans le cas où~$\gA$ est connexe}{La correspondance galoisienne dans le cas connexe}

\Llec est invité\e à revoir le lemme~\ref{lemSspbdiag}.

%:     Theorem{thCorGalGen}
\begin{theorem}
\label{thCorGalGen}
Si $(\gk,\gA,G)$ est une \aG non triviale et si~$\gA$ est \emph{connexe},
la correspondance galoisienne
établit une bijection décroissante entre
\begin{itemize}
\item  d'une part, l'ensemble des
sous-groupes détachables de $G$,
\item  et d'autre part, l'ensemble des sous-\klgs
de~$\gA$ qui sont \spbsz.
\end{itemize}
Ce dernier ensemble est \egmt celui des sous-\algs de~$\gA$ qui sont \stes
sur~$\gk$.
\end{theorem}
%--- end-theorem-----------------------------------------

\begin {proof}
Soit~$\gk\subseteq\gA'\subseteq\gA$ avec~$\gA'$ \spbz. En posant $H =
\Stp(\gA')$, nous devons montrer que~$\gA' = \Ae H$. On a bien s\^ur~$\gA'
\subseteq \Ae H$. \\
Considérons l'\Alg produit~$\gC = \prod_{\sigma \in
G}\gA \simeq \gA^n$ avec $n = \#G$. 
\\
Soit $p_\sigma : \gC \to \gA$ la \prn
définie par $p_\sigma\big((a_\tau)_\tau\big) = a_\sigma$.  
%Comme dans le théorème d'Artin, 
Rappelons  
l'\iso  d'\Algs $\psi_G : \gA\otimes_\gk\gA \to \gC,\;a \otimes b \mapsto \big(a\sigma(b)\big)_{\sigma \in G}$.
\\
Puisque~$\gA$ est un \kmo\ptfz, le morphisme canonique~$\gA\otimes_\gk\gA'
\to \gA\otimes_\gk\gA$ est injectif. En le composant avec $\psi_G$,
nous obtenons un morphisme injectif d'\Algs 
$\gA\otimes_\gk\gA' \to \gC$. Dans les notations,
nous identifierons~$\gA\otimes_\gk\gA'$ à son image~$\gB$ dans~$\gC \simeq \gA^n$.
\\
Puisque~$\gA'$ est une \klg\spbz,~$\gB$ est une \Alg\spbz. Nous pouvons donc
appliquer le lemme~\ref{lemSspbdiag}. Si l'on note $\pi_\sigma$ la restriction 
de~$p_\sigma$ à~$\gB$, il faut identifier la relation d'équivalence sur $G$
définie par $\pi_\sigma = \pi_{\sigma'}$. 
Pour $a' \in \gA'$, $1 \otimes a'$
correspond par $\psi_G$ à $\big(\tau(a')\big)_\tau$, 
donc $\pi_\sigma(1 \otimes a') = \sigma(a')$. 
En conséquence, $\pi_\sigma =
\pi_{\sigma'}$ \ssi $\sigma $ et $\sigma'$ coïncident sur~$\gA'$ ou encore,
par \dfn de~$H$, \ssi $\sigma^{-1}\sigma' \in H$, i.e.
$\sigma H = \sigma' H$. 
On en déduit que les classes d'équivalence 
 sont les classes à gauche modulo~$H$. Avec les notations
du lemme~\ref{lemSspbdiag}, on a donc~$\gB = \bigoplus_J \gA e_J$, où~$J$
décrit $G/H$. 
En utilisant la~$\gA$-base $(e_J)_J$ de~$\gB$, on voit alors que~$\gB = \gC^H$.
\\
Reste à \gui {redescendre} à~$\gA$. Par image réciproque par $\psi_G$,
on a 
$$
(\gA\otimes_\gk\gA)^{\Id\otimes H} = \gA\otimes_\gk\gA'\,.
$$
%\sni
En particulier,~$\gA\otimes_\gk\Ae H \subseteq \gA\otimes_\gk\gA'$.
En appliquant $\Tr_G \otimes \Id_\gA$ à cette inclusion et en
utilisant le fait que $\Tr_G : \gA \to \gk$ est surjective, on obtient
l'inclusion

\snic{\gk\otimes_\gk\Ae H \subseteq \gk\otimes_\gk\gA'$, i.e.~$\;\Ae H \subseteq
\gA'\,.}

%\sni
Ainsi~$\Ae H = \gA'$, ce qu'il fallait démontrer.

 Enfin, puisque les \klgs~$\Ae H$ sont \stes et que les \algs
\stes sont \spbsz, il est clair que les sous-\klgs \spbs de~$\gA$
coïncident avec les sous-\klgs \stesz.
\end {proof}

\rem La théorie des \aGs ne réclame pas vraiment le recours aux \algs \spbsz,
même pour le \tho précédent que l'on peut énoncer avec les seules
 sous-\algs \stes de~$\gA$. Pour une \dem de ce \tho sans recours aux
 \algs \spbsz, voir les exercices~\ref{lemPaquesFerrero} et~\ref{exoCorGalste}.
 Néanmoins la théorie des \algs \spbsz, remarquable pour elle-même,
  apporte un bel éclairage à la situation galoisienne.
 \eoe

%%%%%%%%%%%%%%%%%%%%%%%%%%%%%%%%%%%%%%%%%%%%%%%%%%%%%%%%%%%%%%%%%%%%%%%%%%%
%: subsec{Quotient d'\aGs
\subsec{Quotients d'\aGsz}

%:     Proposition}\label{prop1GalQuo}
\begin {proposition}\label{prop1GalQuo}
\emph{(Quotient d'une \aG par un \id invariant)} Soit $(\gk,\gC,G)$ une \aGz, $\fc$ un \id $G$-invariant de~$\gC$
 et~$\fa = \fc\cap\gk$. 
\begin{enumerate}
\item Le triplet $(\gk\sur\fa, \gC\sur\fc, G)$
est une \aGz. 
\item Cette \aG est naturellement isomorphe à celle obtenue à partir
de~$(\gk,\gC,G)$ au moyen de l'\edsz~$\gk\to\gk\sur\fa$. 

%
%\item 
%
\end{enumerate}
 \end {proposition}

\begin {proof} \emph{1.} Le groupe $G$ opère sur~$\gC\sur\fc$ parce que $\fc$ est (globalement) invariant.
Montrons que l'\homo injectif naturel~$\gk\sur\fa \to
(\gC\sur\fc)^G$ est surjectif.  Si $x \in \gC$ est $G$-invariant modulo
$\fc$, on doit trouver un \elt de~$\gk$ égal à $x$ modulo $\fc$.  On
considère $x_0 \in \gC$ vérifiant $\Tr_G(x_0) = 1$; alors $\Tr_G(xx_0)$ convient:

\snic{x = \som_{\sigma \in G} x\sigma(x_0) \equiv 
\som_{\sigma \in G}  \sigma(x)\sigma(x_0) = \Tr_G(xx_0) \mod \fc\,.}

%\sni
Ainsi $(\gC\sur\fc)^G =\gk\sur\fa$. Enfin il est clair que $G$ opère de façon séparante sur~$\gC\sur\fc$.

\emph{2.} L'\edsz~$\gk\to\gk\sur\fa$ donne 
$(\gk\sur\fa\!,\gC\sur{\fa\gC}\!,G)$ (\aGz), avec $\fa\gC\subseteq\fc$. On doit vérifier que $\fc=\fa\gC$. \\
La \prn  $\pi:\gC\sur{\fa\gC}\to\gC\sur\fc$ est une 
application~$\gk\sur\fa$-\lin
surjective entre modules \prosz, donc~$\gC\sur{\fa\gC}\simeq \gC\sur\fc\oplus\Ker\pi$. Comme les deux modules ont même rang constant $\# G$, le \polmu
de $\Ker\pi$ est égal à~$1$, donc $\Ker\pi=0$ (\thref{th ptf sfio}).
\end {proof}

Dans la \dfn qui suit, on n'a pas besoin de supposer que $(\gk,\gC,G)$
est une \aGz.

%:     Definition{defQuoDeGal}----
\begin{definition}
\label{defQuoDeGal}\label{defIdmGalAdu}\index{ideal@idéal!galoisien}
\index{quotient de Galois!d'une \algz}
Soit $G$ un groupe fini qui opère sur une \klgz~$\gC$.
\begin{enumerate}
\item Un \idm  de~$\gC$ est dit \ixc{galoisien}{idempotent --- dans une \alg munie 
d'un groupe fini d'\autos} si son orbite sous $G$ est un \sfio (cela requiert que cette orbite soit un ensemble fini, ou, de manière \eqvez, que 
le sous-groupe~$\St_G(e)$ soit détachable).
\item Un \id de~$\gC$ est dit
\ixc{galoisien}{idéal ---}
lorsqu'il est engendré par l'\idm
\cop d'un \idm galoisien $e$.
\item Dans ce cas, le groupe $\St_G(e)$
opère sur l'\algz~$\gC[1/e]\simeq\aqo{\gC}{1-e}$,
et~$\big(\gk,\gC[1/e],\St_G(e)\big)$ est appelé
un \emph{quotient de Galois} de $(\gk,\gC,G)$.
\end{enumerate} 
\end{definition}
%--- end-definition------------------------------------

%:     Fact{factQuoDeGal}
\begin{fact}\label{factQuoDeGal}
Avec les hypothèses de la \dfn \ref{defQuoDeGal}, si $\so{e_1,\ldots,e_r}$
est l'orbite de $e$, l'\kli naturelle~$\gC\to \prod_{i=1}^r \gC[1/ {e_i}]$
est un \iso de \klgsz. En outre, les $\St_G(e_i)$ sont deux à deux conjugués
par des \elts de $G$ qui permutent les \klgsz~$\gC[1/ {e_i}]$ 
(elles sont donc deux à deux isomorphes). En particulier~$\gC\simeq\gC[1/e]^r$.
\end{fact}
\facile

Le \tho qui suit est le point \emph{7} du \thref{thADG1Idm} consacré aux quotients de Galois des \apGsz.

%:     proposition{corAGQuo}
\begin{theorem}\label{corAGQuo} \emph{(Quotients de Galois)}\\
Tout quotient de Galois d'une \aG est une \aGz.
\end{theorem}
\incertain{

%%%%%%%%%%%%%%%%%%%%%%%%%%%%%%%%%%%%%%%%%%%%%%%%%%%%%%%%%%%%%%%%%%%%%%%%%%%
%: subsec{Algèbres galoisiennes particulières
\subsec{Algèbres galoisiennes particulières}

%:     Proposition{propAlgGalSitElem}
\begin{proposition}\label{propAlgGalSitElem}
Dans la situation galoisienne \elr décrite dans la proposition~\ref{propGaloiselr}, le triplet $(\gA,\gB,G)$ est une \aGz.
\end{proposition}
%--------- fin proposition ---------------------------------------------- 

%
%:     proposition{corAGcdi}
\begin{proposition}\label{corAGcdi}
Toute \aG sur un \cdi infini~$\gK$ est un quotient de Galois d'une \adu pour un \pol \splz.
\end{proposition}
\begin{proof}
Voir le fichier \texttt{AlgGalPart.tex} 
\end{proof}
%

%-% ENTRE NOUS
\entrenous{on ne peut sans doute pas se débarrasser de l'hypothèse~$\gK$ infini:
faire opérer le groupe cyclique d'ordre $r$ sur $\FFp^r$ pour $r$ grand. 
}
%-% Fin ENTRENOUS

%:     proposition{corAGicl}
\begin{proposition}\label{corAGicl}
Toute \aG $(\gk,\gA,G)$ sur un  \acl infini est un quotient de Galois d'une \adu pour un \pol \splz.
\end{proposition}
\begin{proof}
%-% ENTRE NOUS
\entrenous{La situation n'est pas claire}
%-% Fin ENTRENOUS
\end{proof}

}
%:-% Fin Incertain

%%%%%%%%%%%%%%%%%%%%%%%%%%%%%%%%%%%%%%%%%%%%%%%%%%%%%%%%%%%%%%%%%%%%%%%%%%%
%:    Exercices
\Exercices

%--- Exercise{exochapAlgStricFiLecteur}-------------
\begin{exercise}
\label{exochapAlgStricFiLecteur}
{\rm  Il est recommandé de faire les \dems non données, esquissées,
laissées \alecz,
etc\ldots\, On pourra notamment traiter les cas suivants. \perso{à compléter à la fin}
\begin{itemize}
\item  Montrer le  \thref{factSDIRKlg}.
\item  Montrer le fait \rref{factEdsAlg}.
\item  Montrer le \plgref{plcc.aGs} pour les \aGsz.
\item  Vérifier le fait \rref{factdefAGAL}. 
\end{itemize}
}
\end{exercise}
%--- end -exercise-----------------------------------------

%--- Exercise{exothEtalePrimitif}-------------
\begin{exercise}
\label{exothEtalePrimitif}
{\rm  Donner une \dem détaillée du point \emph{2} du \thref{thEtalePrimitif}.
 
}
\end{exercise}
%--- end -exercise-----------------------------------------

%--- Exercise{exoKnEltPrimitif}-------------
\begin{exercise}
\label{exoKnEltPrimitif}
{\rm
On considère l'\Alg produit~$\gB = \Ae n$.
 
\emph{1.} \`A quelle condition un $x \in \gB$ vérifie-t-il~$\gB = \gA[x]$? \\
Dans ce cas, montrer que
$(1, x, \ldots, x^{n-1})$ est une~$\gA$-base de~$\gB$.
 
\emph{2.} Si~$\gA$ est un \cdiz, à quelle condition~$\gB$ admet-elle un \elt
primitif?
}
\end {exercise}
%--- end -exercise-----------------------------------------

%--- Exercise{exoIdmChangBase}-------------
\begin{exercise}
\label{exoIdmChangBase}
{\rm Soient~$\gK$ un \cdi non trivial,
$\gB$ une \Klg \stfe réduite et $v$ une \idtrz.\\
On considère la \Llgz~$\gB(v)\eqdefi \gK(v)\otimes_\gK\gB$.
Montrer les résultats suivants.

\emph{1.}~$\gB(v)$ est \stfe sur~$\gK(v)$.
 
\emph{2.}   Si~$\gB$
est étale sur~$\gK$,~$\gB(v)$ est étale sur~$\gK(v)$.
 
\emph{3.} Tout \idm de~$\gB(v)$ est en fait dans~$\gB$.

%-% ENTRE NOUS
\entrenous{l'exo \ref{exoIdmChangBase} rappelle des choses dans le cours:
\thrf{thIdmEtale}. Mais la preuve directe
pour le cas des \cdis peut être un bon entraînement;
peut-être il faudrait demander
de \gnr lorsque~$\gK$ est un anneau arbitraire: la situation
est contrôlée par le fait que~$\gK$ est \icl dans~$\gL$,
et il semble que l'on  a besoin de~$\gB$ étale sur~$\gK$.
}
%-% Fin ENTRENOUS
}
\end{exercise}
%--- end -exercise-----------------------------------------

%--- Exercise{exo1SepFact}-------------
\begin{exercise}
\label{exo1SepFact}
{\rm  Si~$\gK$ est un \cdi \splz ment factoriel, il en
va de même pour~$\gK(v)$, où $v$ est une \idtrz.
\perso{L'exo \ref{exo1SepFact} semble redoutable}
\\
NB: on ne suppose pas que~$\gK$ est infini, ni non plus qu'il est fini.
}
\end{exercise}
%--- end -exercise-----------------------------------------

%:HHH exoFreeAlgebraPresentation ramene dans le chap 4
%--- Exercise{exoBiquadratique}-------------
\begin{exercise}\label{exoBiquadratique}  {(Les anneaux d'entiers de l'extension
$\QQ(\sqrt a) \subset \QQ(\sqrt a,\sqrt 2)$)}
\\
{\rm
Soient~$\gK \subseteq \gL$ deux corps de nombres et~$\gA \subseteq \gB$ leurs
anneaux d'entiers; on donne ici un exemple élémentaire où~$\gB$
n'est pas un~$\gA$-module libre.

 \emph {1.}
Soit $d \in \ZZ$ sans facteur carré; déterminer l'anneau
des entiers de $\QQ(\sqrt d)$.

Soit $a \in \ZZ$ sans facteur carré avec $a \equiv 3 \bmod 4$. %\\
On pose~$\gK
= \QQ(\sqrt a)$,~$\gL = \gK(\sqrt 2)$, et~$\beta = \sqrt 2 {1 + \sqrt a \over
2}$. 
On définit $\sigma \in \Aut(\gL/\gK)$  et $\tau
\in \Aut\big(\gL/\QQ(\sqrt 2)\big)$, 
 par:\\
\centerline{$\sigma  (\sqrt 2) = -\sqrt 2\;$  et $\;\tau ( \sqrt a) = -\sqrt a$.}

 \emph {2.}
Vérifier que $\beta \in \gB$ et calculer $(\sigma\tau)(\beta)$.

 \emph {3.}
On veut montrer que $(1, \sqrt 2, \sqrt a, \beta)$ (qui est une $\QQ$-base de
$\gL$) est une $\ZZ$-base de~$\gB$.  Soit $z = r + s\sqrt 2 + t\sqrt a +
u\beta \in \gB$ avec $r$, $s$, $t$, $u \in \QQ$. \\
En considérant $(\sigma\tau)(z)$,
montrer que $u \in \ZZ$ puis que $r$, $s$, $t \in \ZZ$.

 \emph {4.}
Expliciter~$\gB$ comme \Amo \ptfz. Vérifier qu'il est isomorphe
à son dual.
}
\end {exercise}
%--- end -exercise-----------------------------------------

%--- Exercise{exoTensorielDiscriminant}-------------
\begin{exercise}\label{exoTensorielDiscriminant}
 {(Discriminant du produit tensoriel)}\\
{\rm
Soient $\gA$, $\gA'$ deux
\klgs libres de rangs $n$, $n'$, $(\ux) = (x_i)$ une famille
de $n$ \elts de~$\gA$, $(\ux') = (x'_j)$ une famille
de $n'$ \elts de~$\gA'$. On note~$\gB = \gA\otimes_\gk\gA'$
et $(\ux\otimes\ux')$ la famille $(x_i \otimes x'_j)$ de $nn'$
\elts de~$\gB$. Montrer l'\egt

\snic {
\Disc_{\gB\sur\gk}(\ux\otimes\ux') =
\Disc_{\gA\sur\gk}(\ux)^{n'} \, \Disc_{\gA'\sur\gk}(\ux')^n.
}

}
\end {exercise}
%--- end -exercise-----------------------------------------

%--- Exercise{CyclicNormalBasis}-------------
\begin{exercise}\label{CyclicNormalBasis}
 {(Base normale d'une extension cyclique)}
\\
{\rm  
Soit~$\gL$ un \cdiz, $\sigma \in \Aut(\gL)$ d'ordre $n$
et~$\gK = \gL^\sigma$ le corps des invariants sous
$\sigma$. Montrer qu'il existe $x \in \gL$
tel que $\big(x, \sigma(x), \cdots, \sigma^{n-1}(x)\big)$
soit une~$\gK$-base de~$\gL$; on parle alors
de \emph{base normale} de~$\gL\sur\gK$ (définie par $x$).
}

\end {exercise}
%--- end -exercise-----------------------------------------

%--- Exercise{exoHomographieOrdre3}-------------
\begin{exercise}\label{exoHomographieOrdre3} \label{NOTAAn} {(Homographie d'ordre $3$ et \eqn \uvle
de groupe de Galois~$\rA_3$)}  
{\rm On note
$\rA_n$ le sous-groupe des permutations paires de~$\Sn$.
 Soit~$\gL=\gk(t)$ où~$\gk$ est un \cdi et $t$ une \idtrz. 
\begin {enumerate}
\item
Vérifier que $A = \crmatrix {0 &-1\cr 1 & 1\cr}$ est d'ordre $3$ dans
$\PGL_2(\gk)$ et expliquer la provenance de cette matrice.
\end{enumerate}
On note $\sigma\in \Aut_\gk\big(\gk(t)\big)$ l'\auto d'ordre $3$ associé à
$A$ (voir le \pbz~\ref{exoLuroth1}, on~a~$\sigma(f)=f({-1\over t+1})$), et $G = \gen
{\sigma}$.
\begin{enumerate}\setcounter{enumi}{1}
\item Calculer $g=\Tr_G(t)$ et montrer que~$\gk(t)^G = \gk(g)$.

\item
Soit $a$ une \idtr sur~$\gk$ et $f_a(T) = T^3 - aT^2 - (a+3)T - 1
\in \gk(a)[T]$.  Montrer que $f_a$ est \irdz, de groupe de Galois $\rA_3$. 

\item
Montrer que le \pol $f_a(X)$ est un \gui {\pol générique
de groupe de Galois $\rA_3$} au sens suivant: si~$\gL/\gK$ est une 
extension galoisienne de groupe de Galois $\rA_3$ ($\gL$ un \cdiz), il existe un \elt primitif
de~$\gL/\gK$ dont le \polmin est $f_\alpha(X)$ pour une certaine
valeur de $\alpha \in \gK$.

\end {enumerate}

}
\end {exercise}
%--- end -exercise-----------------------------------------

%--- Exercise{exoGroupAlgebra}-------------
\begin{exercise}
\label{exoGroupAlgebra} {(Algèbre d'un groupe commutatif fini)}\\
{\rm  
Soit~$\gk$ un anneau commutatif, $G$ un groupe commutatif d'ordre $n$
et~$\gA = \gk[G]$ l'\emph{\alg du groupe $G$}, i.e.~$\gA$ admet $G$
comme~$\gk$-base et le produit dans~$\gA$ de deux \elts de $G$
est leur produit dans $G$.\footnote{La \dfn fonctionne aussi pour \emph{l' \alg $\gk[M]$ d'un \mo $M$}.}\index{algèbre!d'un monoïde} 
\begin{enumerate}
\item Déterminer $\Ann(\rJ\iAk)$, son image par $\mu\iAk$ et
la forme bi\lin tracique sur~$\gA$.
\item Montrer que \propeq
\begin {itemize}
\item
$n$ est \iv dans~$\gk$.

\item
$\gA$ est \stez.

\item
$\gA$ est \spbz.
\end {itemize}

\item Montrer que~$\gk[G]$ est une \alg de Frobenius.
\end{enumerate}

}
\end {exercise}
%--- end -exercise-----------------------------------------

%:--- Exercise{exo1Frobenius}-------------
\begin{exercise}\label{exo1Frobenius}
{(Une \alg monogène finie est une \alg de Frobenius)}\\
{\rm  Soit $f = X^n + a_{n-1} X^{n-1} + \cdots + a_0 \in \gk[X]$ et~$\gA = \aqo{\kX}{f} = \gk[x]$.  On considère la forme \lin
$\lambda : \gA \to \gk$ définie par $x^{n-1} \mapsto
1$ et $x^i \mapsto 0$ pour $i < n-1$.
   On va montrer que $\lambda$ est
dualisante et que $\Tr\iAk = f'(x)\centerdot\lambda$. 
\\
\`A cet effet, on adjoint une \idtr $Y$. Le \sys
$(1, x, \ldots,x^{n-1})$ est  une base de~$\gA[Y]\sur{\gk[Y]}$, 
on note $\wi\lambda : \gA[Y] \to \gk[Y]$ l'extension de $\lambda$
et l'on définit l'application 
$\gk[Y]$-\lin  $\varphi : \gA[Y] \to \gk[Y]$, par $\varphi(x^i)=Y^i$ pour $i \in\lrb{0.. n-1}$.
\begin {enumerate}
\item 
Montrer que:
$\qquad\forall g \in \gA[Y],\quad f(Y) \wi\lambda(g) = \varphi\big((Y - x)g\big)\qquad(*)$

\item 
On définit la base (triangulaire de Horner) $(b_0, \ldots, b_{n-1})$ de
$\gA\sur\gk$ par

\snic {\arraycolsep2pt
\begin{array}{ccc} 
b_0  &  = & x^{n-1} + a_{n-1} x^{n-2} + \cdots + a_2x + a_1,  \\[1mm] 
b_1   & =   & x^{n-2} + a_{n-1} x^{n-3} + \cdots + a_3x + a_2  \,,
 \end{array}
 }

%\sni
et ainsi de suite: $b_i = x^{n-i-1} + \cdots + a_{i+1}$ et $b_{n-1} = 1$.
On a:
$$
f'(Y) = \frac{f(Y) - f(x)}{Y - x}
= \frac{f(Y)}{Y - x} = b_{n-1}Y^{n-1} + \cdots + b_1Y + b_0.
$$
Montrer en appliquant l'\egt $(*)$
à $g_i = x^if'(Y)$, que $(b_0\centerdot\lambda, \ldots, b_{n-1}\centerdot\lambda)$
est la base duale de $(1, x, \ldots, x^{n-1})$.
Conclure.

\item 
Montrer que $\Tr\iAk = f'(x)\centerdot\lambda$.
\end {enumerate}

}

\end {exercise}
%--- end -exercise-----------------------------------------

%--- Exercise{exoFrobeniusAlgExemples}-------------
\begin{exercise}\label{exoFrobeniusAlgExemples}
 {(Algèbres de Frobenius: exemples et contre-exemples \elrsz)}\\
{\rm  
Dans tout l'exercice,~$\gk$ est un anneau commutatif.

\emph {1.}
Soit $f_1$, \ldots, $f_n \in \gk[T]$ des \polusz. Montrer
que la \klg quotient 
$\gk[X_1, \ldots, X_n]\sur{\gen {f_1(X_1), \ldots, f_n(X_n)}}$
est libre de rang fini, de Frobenius.

\emph {2.}
Soit~$\gA = \gk[X,Y]\sur{\gen{X,Y}^2} = \gk[x,y]$. Décrire
$\Asta$ comme \Amo \pfz; en déduire que~$\gA$ n'est pas
une \alg de Frobenius.

\emph {3.}
Question analogue à la précédente avec~$\gA = \gk[X,Y]\sur{\gen
{X,Y}^n}$ pour $n \ge 2$ et~$\gB = \gk[X,Y]\sur{\gen {X^2,XY^{n+1},Y^{n+2}}}$
\hbox{pour $n \ge 0$}.

}
\end {exercise}
%--- end -exercise-----------------------------------------

%--- Exercise{exoAlgMonogeneJAnnJ}-------------
\begin{exercise}\label{exoAlgMonogeneJAnnJ} {(L'\id $\rJ\iAk$ pour une \klg monogène~$\gA$)}\\
{\rm  
Soit~$\gA = \gk[x]$ une \klg monogène et $\env{\gk}{\gA}
= \gA\te_\gk\gA$ son \alg enveloppante. 
On pose $y = x\te 1$, $z = 1\te x$, de sorte que $\env\gk\gA =
\gk[y,z]$.
On sait que $\rJ\iAk=\gen{y-z}$.
On suppose $f(x) = 0$ pour un $f \in \gk[X]$ non \ncrt \mon et
l'on considère %$f^\Delta \in \gk[Y,Z]$ 
le \pol \smq 
%défini par
$f^\Delta(Y,Z) = \big(f(Y) - f(Z)\big) /(Y-Z)$. Il vérifie 
l'\egt $f^\Delta(X,X) =
f'(X)$.

\emph{1.} 
Soit $\delta = f^\Delta(y,z)$. Montrer que $\delta \in \Ann(\rJ\iAk)$
et que $\delta^2 = f'(y)\delta = f'(z)\delta$.

\emph{2.} 
On suppose que $1 \in \gen {f, f'}$.\\
\emph{2a.} 
 Montrer que~$\gA$
est \spbz: expliciter l'\idm de séparabilité. \\
\emph{2b.} Montrer que  $\rJ\iAk=\gen{f^\Delta(y,z)}$
et que $f^\Delta(y,z) = f'(y)\vep\iAk = f'(z)\vep\iAk$.

Remarque:~$\gA$ n'est pas \ncrt \stfez.

}

\end {exercise}
%--- end -exercise-----------------------------------------

%--- Exercise{exoJacInversibleAlgSeparable}-------------
\begin{exercise}\label{exoJacInversibleAlgSeparable}
 {(Intersection complète, jacobien, bezoutien et \sptz)}\\
{\rm  
Dans cet exercice, le nombre d'\idtrs est égal au nombre de \polsz.  
On définit le \emph{bezoutien} de $(f_1, \ldots, f_n) $ où les $f_i\in \gk[\uX] = 
\gk[\Xn]$ par:

\snic{
\beta_{\uY,\uZ}(\uf) = \det \BZ_{\uY,\uZ}(\uf)
}
de sorte que $\beta_{\uX,\uX}(\uf) = \J_\uX(\uf)$.
\index{bezoutien!déterminant --- d'un \sypz}

On désigne par $\gA = \gk[\xn]$ une \klg \tfz, $\env{\gk}{\gA} =
\gk[\uy,\uz]$ son \alg enveloppante; on suppose que $f_i(\ux) = 0$ pour tout
$i$.

\emph {1.}
Dans le cas où $\J_\ux(f_1, \ldots, f_n)\in\Ati$, fournir une preuve
directe du fait que $\gA$ est une \alg \spbz.

\emph {2.}
On définit dans $\env{\gk}{\gA}$:

\snic {
\vep = \J_\uy(\uf)^{-1} \beta_{\uy,\uz}(\uf) = \beta_{\uy,\uz}(\uf) \J_\uz(\uf)^{-1}
.}

\sni
Vérifier que $\beta_{\uy, \uz}(\uf)$ et $\vep$ sont des \gtrs de $\Ann(\rJ\iAk)$
et que $\vep$ est l'\idst de $\gA$.

\emph {3.}
Donner des exemples.
}

\end {exercise}
\begin{exercise}\label{exoSepAlgDedekindLemma} 
{(Séparation des morphismes sur une \alg\spbz)}\\
{\rm  
Soient~$\gk$ un anneau commutatif et~$\gA$,~$\gB$ deux \klgs avec~$\gA$ \spbz.
Pour une fonction quelconque $f : \gA \to \gB$, on définit
$\Ann_\gB(f) = \Ann_\gB \gen {f(\gA)}$.

\emph {1.}
Montrer qu'à tout morphisme $\varphi \in \Hom_\gk(\gA,\gB)$ est attaché un
couple de familles finies $(a_i)_{i \in I}$, $(b_i)_{i \in I}$, avec $a_i \in
\gA$, $b_i \in \gB$, vérifiant les \prts suivantes% de commutation
:
\begin{itemize}
\item $\sum_i b_i\varphi(a_i) = 1$
\item $\sum_i \varphi(a)b_i \te a_i = \sum_i b_i \te aa_i$ pour tout $a \in \gA$.
\end{itemize}

\emph {2.}
Si au morphisme $\varphi' \in \Hom_\gk(\gA,\gB)$ est attaché le couple
de familles $(a'_j)_j$, $(b'_j)_j$, montrer que 

\snic {
\sum_i b_i\varphi'(a_i) = \sum_j b'_j\varphi(a'_j)
,}

%\sni
et que ce dernier \eltz, noté $e$, est un \idm de~$\gB$ ayant 
la \prt suivante de \gui {séparation des morphismes}:

\snic {
\Ann_\gB(\varphi - \varphi') = \gen {e}_\gB, \qquad
\gen {\Im(\varphi - \varphi')}_\gB  = \gen {1-e}_\gB
.}

%\sni
\emph {3.}
Soient $\varphi_1$, \ldots, $\varphi_n \in \Hom_\gk(\gA,\gB)$ et,
pour $i$, $j\in\lrbn$, $e_{ij} =
e_{ji}$ l'\idm défini par $\Ann_\gB(\varphi_i - \varphi_j) = \gen
{e_{ij}}_\gB$; en particulier, $e_{ii} = 1$.  On dit qu'une matrice~$A \in
\MM_{n,m}(\gB)$ est une \emph {matrice d'\evn de Dedekind} 
pour les~$n$ morphismes $\varphi_1$, \ldots, $\varphi_n$ 
si chaque colonne de $A$ est de
la forme $\tra {[\, \varphi_1(a) \; \cdots \;\varphi_n(a)\,]}$ pour un
$a \in \gA$ (dépendant de la colonne).  Montrer l'existence d'une matrice
d'\evn de Dedekind dont l'image contient les vecteurs $\tra
{ [\,e_{1i} \; \cdots \; e_{ni}\,]}$.  En particulier, si
$\Ann_\gB(\varphi_i - \varphi_j) = 0$ pour $i \ne j$, une telle matrice
est surjective.

}
\end {exercise}
%--- end -exercise-----------------------------------------

%--- Exercise{exoArtinAnOtherProof}-------------
\begin{exercise}\label{exoArtinAnOtherProof}
{(Une autre \dem du point \emph{2} du \tho d'Artin)}\\
{\rm
Le contexte est celui du théorème~\ref{thA}: on suppose que $(\gk,\gA,G)$
est une \aG et l'on veut montrer l'existence de $a_1$, \ldots, $a_r$, $b_1$, \ldots,
$b_r \in \gA$ tels que pour tout $\sigma \in G$ on~ait:

\snic{\sum_{i=1}^r a_i \sigma(b_i)=
\formule{1\;\;\mathrm{si}\;\sigma=\Id\\
0\;\;\mathrm{sinon}.}
}

%\sni
Pour $\tau \in \G$, $\tau \ne \Id$, montrer qu'il existe $m_\tau$
et $x_{1,\tau}$, \ldots, $x_{m_\tau,\tau}$, $y_{1,\tau}$, \ldots, $y_{m_\tau,\tau}$ 
dans~$ \gA$ tels que:

\snic {
\sum_{j=1}^{m_\tau} x_{j,\tau} \tau(y_{j,\tau}) = 0, \qquad
\sum_{j=1}^{m_\tau} x_{j,\tau} y_{j,\tau} = 1
.}

Conclure.

}
\end {exercise}

%--- Exercise{exoGalExemples}-------------
\begin{exercise}\label{exoGalExemples} 
{(Algèbres galoisiennes: quelques exemples \elrsz)}\\
{\rm  
On note $(e_1, \ldots, e_n)$ la base canonique de~$\gk^n$.  
On fait agir $\rS_n$  sur~$\gk^n$ par permutation des \coosz: $\sigma(e_i) = e_{\sigma(i)}$ pour $\sigma \in \rS_n$.

\emph {1.}
Soit $G \subset \rS_n$ un sous-groupe transitif de cardinal $n$.
\begin {enumerate}
\item [\emph{a.}]
Montrer que $(\gk, \gk^n, G)$ est une \aGz.
\item [\emph{b.}]
Donner des exemples.
\end {enumerate}

\emph {2.}
Soit~$\gB = \gk(e_1 + e_2) \oplus \gk(e_3 + e_4) \subset \gk^4$ et $G = \gen
{(1,2,3,4)}$. Déterminer $\Stp_{\rS_4}(\gB)$ et $H = \Stp_G(\gB)$.  
Est-ce que l'on a~$\gB = (\gk^4)^H$?

\emph {3.}
Soit $(\gk, \gA, G)$ une \aGz; le groupe $G$ opère naturellement
sur~$\gA[X]$. 
\begin {enumerate}
\item [\emph{a.}]
Montrer que $(\gk[X], \gA[X], G)$ est une \aGz.
\item [\emph{b.}]
Soit~$\gB = X\gA[X] + \gk$ ($\gB$ est donc constitué des \pols
de~$\gA[X]$ dont le \coe constant est dans~$\gk$). Alors,~$\gB$ 
est une sous-\klg de~$\gA[X]$ qui n'est pas de la forme
$\gA[X]^H$ sauf dans un cas particulier.
\end {enumerate}

}
\end {exercise}
%--- end -exercise-----------------------------------------

%--- Exercise{lemPaquesFerrero}-------------
\begin{exercise}
\label{lemPaquesFerrero}
{\rm  Soient~$\gk\subseteq\gB\subseteq\gC$ avec~$\gB$ \ste sur~$\gk$ et~$\gC$ \stfe
sur~$\gk$. On suppose que  $\rg_\gk(\gB)=\rg_\gk(\gC)$
(i.e.  $\gC$ et~$\gB$ ont même \polmu sur~$\gk$). Alors~$\gB=\gC$.
 
}
\end{exercise}
%--- end -exercise-----------------------------------------

%--- Exercise{exoCorGalste}-------------
\begin{exercise}
\label{exoCorGalste}
{\rm  En vous basant sur l'exercice~\ref{lemPaquesFerrero}
démontrer la correspondance galoisienne (\thref{thCorGalGen})
entre les sous-groupes finis de $G$ et les sous-\klgs \stes de~$\gA$
lorsque $\gA$ est connexe. 
}
\end{exercise}
%--- end -exercise-----------------------------------------

%--- exercise{exoIdeauxGlobInvariants}-------------
\begin{exercise}\label{exoIdeauxGlobInvariants} 
{(Algèbres galoisiennes: idéaux globalement invariants)}\\
{\rm  
Soit $(\gA, \gB, G)$ une \aGz; on dit qu'un idéal $\fc$ 
de~$\gB$ est \emph {globalement invariant} si $\sigma(\fc) 
= \fc$ pour tout $\sigma \in G$.

\emph {1.}
Montrer que $\fc$ est engendré par des \elts invariants,
i.e. par des \elts de~$\gA$.

\emph {2.}
De manière plus précise, on considère les deux transformations
entre idéaux de~$\gA$ et idéaux de~$\gB$: $\fa \mapsto \fa\gB$ et
$\fc \mapsto \fc \cap \gA$.  Montrer qu'elles établissent une
correspondance bijective croissante entre \ids de~$\gA$ et
\ids de~$\gB$ globalement invariants.

}
\end {exercise}
%--- end -exercise-----------------------------------------

%: pbs
%%%%%%%%%%%%%%%%%%%%%%%%%%%%%%%%%%%%%%%%%%%%%%%%%%%%%%%%%%%%%%%%%%%%%%%%%%%

%--- problem{exoLuroth1}-------------
\begin{problem}
\label{exoLuroth1} \label{NOTAPGL} (\Tho de Lüroth)\\
{\rm  Soit~$\gL=\gk(t)$ où~$\gk$ est un \cdi et $t$ une \idtrz.
Si $g=u/v\in\gL$ est une fraction non constante %écrite sous forme 
\ird ($u,v\in\gk[t]$, étrangers),
on définit la \emph{hauteur} de $g$ (par rapport à $t$) par: $\hauteur_t(g)\eqdefi \max\big(\deg_t(u),\deg_t(v)\big)$.
\begin{enumerate}
\item \emph{(Partie directe du \tho de Lüroth)} On pose~$\gK=\gk(g)\subseteq\gL$.
Montrer que~$\gL/\gK$ est une extension \agq de degré $d=\hauteur(g)$.
Plus \prmtz, $t$ est \agq sur~$\gK$ et son \polmin est, à un facteur multiplicatif près dans~$\gK\eti$, égal à $u(T)-g v(T)$.
Ainsi, tout \coe non constant de $\Mip_{\gK,t}(T)$,
$a\in\gK=\gk(g)$ s'écrit  $a=\fraC{\alpha g+\beta}{\gamma g+\delta}$ avec $\alpha\delta-\beta\gamma\in\gk\eti$, et~$\gk(a)=\gk(g)$.  
\item Soit un \elt arbitraire $f\in\gL$.
Donner une formule explicite utilisant les résultants pour exprimer~$f$ comme~$\gK$-\coli de $(1,t,\ldots,t^{d-1})$.
\item Si $h$ est un autre \elt de~$\gL\setminus\gk$ 
montrer que

\snic{\hauteur\big(g(h)\big)=\hauteur(g)\hauteur(h).}
%
%\sni
Montrer que tout~$\gk$-\homo d'algèbre~$\gL\to\gL$ s'écrit $f\mapsto f(h)$ pour un $h\in\gL\setminus\gk$. En déduire une description précise de $\Aut_\gk(\gL)$ au moyen des fractions de hauteur $1$.

\item On note $\PGL_n(\gA)$ le groupe quotient $\GLn(\gA)/\Ati$ (où $\Ati$
est identifié au sous-groupe des homothéties \ivs via $a\mapsto a\In$).  
\`A une matrice

\snic{A = \cmatrix {a &b\cr c&d\cr} \in \GL_2(\gA),}

%\sni  
on associe le~$\gA$-\autoz\footnote{Pour un anneau arbitraire~$\gA$,
l'anneau~$\gA(t)$ est le \gui{localisé de Nagata} de~$\gA[t]$ obtenu en inversant les \pols primitifs.}

\snic{\varphi_A : \gA(t) \to \gA(t),\quad t \mapsto {at + b\over ct + d}. }

%\sni
On a  $\varphi_A \circ \varphi_B =
\varphi_{BA}$ et $\varphi_A = \Id \Leftrightarrow A =
\lambda\I_2 \;(\lambda \in \Ati)$. Ainsi $A \mapsto \varphi_A$ définit
un \homo injectif $\PGL_2(\gA)\eo \to \Aut_\gA\big(\gA(t)\big)$.
\\
 Montrer que 
dans le cas d'un \cdi on obtient un \isoz.
\item \emph{(Partie réciproque du \tho de Lüroth)} 
Soient $g_1$, \ldots, $g_r\in\gL\setminus \gk$. Montrer que~$\gk(g_1,\ldots,g_r)=\gk(g)$ pour un $g$ convenable. Il suffit de traiter le cas~$n=2$. On montre que
$\gL$ est \stf sur~$\gK_1=\gk(g_1,g_2)$, on doit alors avoir~$\gK_1=\gk(g)$
pour n'importe quel \coe non constant $g$ de~$\Mip_{\gK_1,t}(T)$.
\\
NB.  Puisque~$\gL$ est un~$\gk(g_1)$-\evc de dimension finie,
tout sous-corps de~$\gL$ contenant strictement~$\gk$  est, en \clamaz, \tfz,  
donc de la forme~$\gk(g)$.
Notre formulation de la partie réciproque du \tho de Lüroth donne la signification \cov de cette affirmation.%
\index{hauteur!d'une fraction rationnelle}%
\index{Luroth@Lüroth!\tho de ---} 

\end{enumerate}
}
\end{problem}
%--- end -problem-----------------------------------------

%--- problem{exoBuildingFrobAlgebra}-------------
\begin{problem}\label{exoBuildingFrobAlgebra}
{(Opérateurs différentiels et \algs de Frobenius)}\\
{\rm  
Dans les premières questions,~$\gk$ est un anneau commutatif. La
\emph{dérivée de Hasse} d'ordre $m$ d'un \pol de $\kX$ se
définit \fmt par $f^{[m]} = {1 \over m!}f^{(m)}$. De même, pour
$\alpha \in \NN^n$, on définit $\partial^{[\alpha]}$ sur $\kuX = \gk[X_1,
\ldots, X_n]$ par:
\index{Hasse!dérivée de ---}
\index{derivee@dérivée de Hasse}
$$
\partial^{[\alpha]} f = {1\over\alpha!} 
{\partial^\alpha f \over \partial X^\alpha}
\quad \hbox {avec} \quad \alpha! = \alpha_1! \ldots \alpha_n!,\quad
f \in \kuX 
.$$
On a alors $\partial^{[\alpha]} (fg) = \sum_{\beta+\gamma=\alpha}
\partial^{[\beta]}(f) \, \partial^{[\gamma]}(g)$. On note $\delta^{[\alpha]}
: \kuX \to \gk$ la forme \lin $f \mapsto \partial^{[\alpha]}(f)(0)$. Ainsi,
 $f = \sum_\alpha \delta^{[\alpha]}(f)X^\alpha$. On
en déduit, en notant $\alpha\le\beta$ pour $X^\alpha \divi X^\beta$:

\snic {
X^\alpha \centerdot \delta^{[\beta]} = \formule{
\delta^{[\beta-\alpha]} \hbox{ si } \alpha\le\beta
\\
0 \hbox{ sinon,}}
\qquad
\partial^{[\alpha]}(X^\beta) = \formule{
X^{\beta-\alpha} \hbox{ si } \alpha\le\beta
\\
0 \hbox{ sinon}.}
}

%\sni
Soit $g = \sum_\beta b_\beta X^\beta\in \kuX$.  On évalue en $(\uze)$ le \emph{\pol \dilz}~\hbox{$\sum_\beta b_\beta
\partial^{[\beta]}$}, on obtient une forme \lin $\delta_g : \kuX \to \gk$,
$\delta_g = \sum_\beta b_\beta\delta^{[\beta]}$, puis un \id $\fa_g$ de~$\kuX$:

\snic {
\fa_g = \sotq{f \in \kuX} {f\centerdot \delta_g = 0} \eqdf {\rm def}
\sotq{f \in \kuX} {\delta_g(fu) = 0\ \forall u \in \kuX}
.}

%\sni
On obtient ainsi une \klg de Frobenius $\kuX\sur{\fa_g}$
(avec $\delta_g$ dualisante).

\emph {1.}
Soit $f=\sum_\alpha a_\alpha X^\alpha$, $g=\sum_\beta b_\beta X^\beta$.  On
note $\partial_f : \kuX \to \kuX$ l'opérateur \dil associé à~$f$,
i.e. $\partial_f = \sum_\alpha a_\alpha \partial^{[\alpha]}$. Vérifier
la relation suivante entre l'opérateur $\partial_f$ et la forme \lin
$\delta_g$:

\snic {
\sum_\gamma (f \centerdot \delta_g)(X^\gamma) X^\gamma =
\partial_f(g) = \sum_{\alpha \le \beta} a_\alpha b_\beta X^{\beta-\alpha}
.}

%\sni
En déduire que $f\centerdot \delta_g = 0 \iff \partial_f(g) = 0$.

Il faut maintenant noter que la loi $f*g = \partial_f(g)$ munit le groupe
additif $\kuX$ d'une structure de $\kuX$-module (car en particulier
$\partial_{f_1f_2} = \partial_{f_1} \circ \partial_{f_2}$). 
Mais comme~$X^\alpha * X^\beta = X^{\beta-\alpha}$ ou $0$, certains auteurs utilisent~$X^{-\alpha}$ au lieu de $X^{\alpha}$: ils munissent~$\kuX$ d'une structure
de~$\gk[\uX^{-1}]$-module. D'autres auteurs permutent~$\uX$ et~$\uX^{-1}$:
ils  munissent~$\gk[\uX^{-1}]$ d'une structure de~$\kuX$-module de façon
à ce que l'\id $\fa_g$ (annulateur de $g \in \gk[\uX^{-1}]$) soit un \id
d'un anneau de \pols avec des \idtrs à exposants~$\ge 0$.  Dans ce dernier
formalisme, un \polz~$f$ avec des \idtrs à exposants~$\ge 0$ agit donc sur un
\pol $g$ ayant des \idtrs à exposants $\le 0$ pour fournir un \pol $f*g$
ayant des \idtrs à exposants $\le 0$ (en supprimant les \moms contenant
un exposant $> 0$).  Ainsi, si $g = X^{-2} + Y^{-2} + Z^{-2}$, l'\id $\fa_g$
de~$\gk[X,Y,Z]$, contient par exemple $XY$, $X^2 - Y^2$ et tout \pol\hmg de
degré~$\ge 3$.

\emph {2.}
Soit $d \ge 1$. \'{E}tudier le cas particulier de la somme de Newton $g = \sum_i
X_i^{-d}$, \cad  $\delta_g : f \mapsto \sum_i {1 \over d!} {\partial^d f \over
\partial X_i^d}(0)$, somme des composantes sur $X_1^d, \ldots, X_n^d$.

\smallskip  Dans la suite, on fixe $g = \sum_\beta b_\beta X^\beta$, ou
selon les go\^uts, $g = \sum_\beta b_\beta X^{-\beta}$.

\emph {3.}
Montrer que l'on a une inclusion $\fb \subseteq \fa_g$ pour un \id $\fb = \gen {X_1^{e_1}, \cdots, X_n^{e_n}}$ avec des entiers $e_i \ge
1$. En particulier, $\kuX\sur{\fb}$ est un \kmo libre 
de rang fini et~$\kuX\sur{\fa_g}$ est un \kmo \tfz.

\emph {4.}
Définir une application~$\gk$-\lin $\varphi : \kuX\sur\fb \to \kuX$
telle que $\Ker\varphi = \fa_g\sur\fb$. On peut donc calculer
$\fa_g$ si l'on sait résoudre les \syss \lins sur~$\gk$.

\emph {5.}
On suppose que~$\gk$ est un \cdi et donc~$\gA := \kuX\sur{\fa_g}$ est un
\kev de dimension finie.  Montrer que $(\gA, \delta_g)$ est une \klg de
Frobenius.

}

\end {problem}
%--- end -problem-----------------------------------------

%--- problem{exoTh90HilbertAdditif}-------------
\begin{problem}\label{exoTh90HilbertAdditif} {(Le théorème 90 d'Hilbert, version additive)}\ihi \\
{\rm
Soit $(\gk,\gA,G)$ une \aG où $G = \gen {\sigma}$ est cyclique
d'ordre $n$.

 \emph {1.}
En considérant un
\elt $z \in \gA$ de trace 1, on montrera que:

\snic {
\gA = \Im(\Id_\gA - \sigma) \oplus \gk z, \qquad
\Im(\Id_\gA - \sigma) = \Ker\Tr_G.
}

%\sni
En conséquence  $\Im(\Id_\gA -
\sigma)$ est un \kmo \stl de rang $n-1$.
On pourra utiliser la famille d'\endos $(c_i)_{i \in \lrb{0..n}}$:

\snic {
c_0 = 0,\  c_1(x)=x,\ c_2(x) = x+\sigma(x),\ \ldots,\
c_i(x) = \sum_{j=0}^{i-1} \sigma^j(x),\ \ldots
}

\sni \emph {2.}
Pour qu'un $x\in\gA$ soit de la forme $y - \sigma(y)$, il faut et il
suffit que $\Tr_G(x) = 0$.

 \emph {3.}
Plus \gnltz, soit $(c_\tau)_{\tau\in G}$ une famille dans~$\gA$. Montrer qu'il existe un \elt $y$ tel que $c_\tau = y-\tau(y)$ \ssi la
famille vérifie la condition de cocycle additif suivante, pour tous $\tau_1$, $\tau_2\in G$: $c_{\tau_1\tau_2} =
\tau_1(c_{\tau_2}) + c_{\tau_1}$.

 \emph{4.}
On suppose que $n$ est un nombre premier $p$ et que $p=0$ dans~$\gk$. Montrer
l'existence d'un $y \in \gA$ tel que $\sigma(y) = y+1$. \\
En déduire que $(1,
y, \dots, y^{p-1})$ est une~$\gk$-base de~$\gA$ et que le \polcar de $y$ est
de la forme $Y^p - Y - \lambda$ avec~$\lambda \in \gk$.  \\
On a donc~$\gA =
\gk[y] \simeq \aqo{\gk[Y]}{Y^p - Y - \lambda}$ (extension d'Artin-Schreier).

 \emph{5.}
Donner une réciproque au point précédent.
}
\end {problem}
%--- end -problem-----------------------------------------

%--- Problem{exoHeitmannGaloisExemple}-------------
\begin{problem}\label{exoHeitmannGaloisExemple} 
 {(Algèbres galoisiennes: étude d'un exemple)}
{\rm  
On considère un anneau~$\gB$ dans lequel $2$ est \ivz, avec  $x$, $y
\in \gB$ et $\sigma \in \Aut(\gB)$ d'ordre $2$ vérifiant $x^2 + y^2 = 1$,
$\sigma(x) = -x$ \hbox{et $\sigma(y) = -y$}.  On peut prendre comme exemple l'anneau
$\gB$ des fonctions continues sur le cercle unité $x^2 + y^2 = 1$ 
et pour~$\sigma$ l'involution $f \mapsto \{(x,y) \mapsto f(-x,-y)\}$.  On note~$\gA =
\gB^{\gen {\sigma}}$ (sous-anneau des \gui {fonctions paires}).

\emph {1.}
Montrer que $(\gA, \gB, \gen\sigma)$ est une \aGz.
\\
En conséquence,~$\gB$ est un \Amo\prc 2.

\emph {2.}
Soit $E = \gA x + \gA y$ (sous-module des \gui {fonctions
impaires}). \\
Vérifier que~$\gB = \gA \oplus E$ et que $E$ est un \Amo\prc 1.

\emph {3.}
On pose $x_1 = 1$, $x_2 = x$, $x_3 = y$ de sorte que $(x_1, x_2, x_3)$
est un \sgr du \Amoz~$\gB$. Expliciter $y_1$, $y_2$, $y_3 \in \gB$
comme dans le lemme~\ref{lemArtin}, i.e. $\big((x_i)_{i \in \lrb {1..3}},
(y_i)_{i \in \lrb {1..3}}\big)$ est un \stycz. \\
En déduire 
une \mprn $P \in \MM_3(\gA)$ de rang 2 avec~$\gB \simeq_\gA \Im P$.

\emph {4.}
Soit $R = \crmatrix {x &-y\cr y & x\cr} \in \SL_2(\gB)$. Montrer
que cette \gui {rotation} $R$ induit un \iso de \Amos entre $E^2$ et~$\gA^2$:

\snic {
\cmatrix {f\cr g\cr} \mapsto R \cmatrix {f\cr g\cr} =
\cmatrix {xf-yg\cr yf+xg\cr}
.}

%\sni
En conséquence (question suivante), $E\te_\gA E \simeq \gA$; 
vérifier que $f \te g \mapsto fg$ réalise un \iso de \Amos
de $E \te_\gA E$ sur~$\gA$.

\emph {5.} Pour un \Amo $M$, on note $M^{2\te}=M\te_\gA M$,  $M^{3\te}=M\te_\gA M\te_\gA M$, etc\dots\, Soit $E$ un \Amo ($\gA$ quelconque) vérifiant
$E^n \simeq \gA^n$ pour un certain $n \ge 1$. Montrer que $E$ est
un \Amo\prcz~1 et que $E^{n\te} \simeq \gA$.

\emph {6.}
Soit $\fa$ l'\id de~$\gA$ défini par $\fa = \gen {xy, x^2}$.  Vérifier que
$\fa^2 = x^2\gA$ (donc si $x$ est \ndzz, $\fa$ est un \id\iv de~$\gA$), que
$\fa\gB$ est principal et enfin, que $\fa$, vu comme sous-\Amo de~$\gB$, est
égal à $xE$.

\emph {7.}
Soit~$\gk$ un anneau non trivial avec $2\in\gk\eti$ et $\gB =
\aqo{\gk[X,Y]}{X^2 + Y^2 - 1}$. On écrit $\gB=\gk[x,y].$
\\
On peut appliquer ce qui
précède en prenant $\sigma$  défini par
$\sigma(x) = -x$, $\sigma(y) = -y$.  On suppose que $\alpha^2+\beta^2=0\Rightarrow \alpha=\beta=0$ dans~$\gk$ (par exemple si~$\gk$ est un \cdi et $-1$ n'est pas un carré
dans~$\gk$).
\begin {enumerate}
\item [\emph {a.}]
Montrer que $\Bti = \gk^{\times}$; illustrer l'importance de l'hypothèse
\gui{de réalité} faite sur~$\gk$.

\item [\emph {b.}]
Montrer que $\fa$ n'est pas principal et donc $E$ n'est pas un
\Amo libre. En déduire que~$\gB$ n'est pas un \Amo libre.
\end {enumerate}

\emph {8.}
Soit~$\gB$ l'anneau des fonctions continues (réelles) sur le cercle unité
$x^2 + y^2 = 1$ et $\sigma$ l'involution $f \mapsto \{(x,y) \mapsto f(-x,-y)\}$.
Montrer que $\fa$ n'est pas principal et que~$\gB$ n'est pas un \Amo libre.

}

\end {problem}
%--- end -problem-----------------------------------------

%%%%%%%%%%%%%%%%%%%%%%%%%%%%%%%%%%%%%%%%%%%%%%%%%%%%%%%%%%%%%%%%%%%%%%%%%%%
% fin des exos
%:   solutions d'exos
\sol

%%%%%%%%%%%%%%%%%%%%%%%%%%%%%%%%%%%%%%%%%%%%%%%%%%%%%%%%%%%%%%%%%%%%%%%%%%%
\exer{exothEtalePrimitif}
On a~$\gB=\Kxn$, avec $\dex{\gB:\gK}=m$. On va faire un calcul qui montre que 
la \Klgz~$\gB$ est
 monogène ou contient un \idm $e\neq 0,1$. 
Dans le deuxième cas,~$\gB\simeq \gB_1\times \gB_2$, avec $\dex{\gB_i:\gK}=m_i<m $,
 $m_1+m_2=m$, ce qui permet de conclure par \recu sur $m$.
\\
Si l'on est capable de traiter le cas $n=2$, on a gagné, car~$\gK[x_1,x_2]$
est étale sur~$\gK$, donc ou bien on remplace~$\gK[x_1,x_2]$ par~$\gK[y]$ pour un
 certain  $y$,
ou bien on trouve un \idm $e\neq 0,1$ dedans. La \dem du point \emph{1} du \tho 
montre qu'une \Klg étale~$\gK[x,z]$ est monogène si~$\gK$ 
contient une suite infinie d'\elts distincts.
Elle utilise un \pol $d(a,b)$ qui, évalué dans~$\gK$ doit donner un \elt \ivz.
Si  l'on n'a pas d'information sur l'existence d'une suite infinie d'\elts distincts
de~$\gK$, 
on énumère les entiers de~$\gK$
jusqu'à obtenir $\alpha,\beta$ dans~$\gK$ avec $d(\alpha,\beta)\in\gK\eti$, 
ou à conclure que la \cara est égale à un nombre premier $p$. 
 On énumère ensuite
les puissances des \coes de~$f$ et de $g$ (les \polmins de $x$ et $z$ sur~$\gK$) jusqu'à obtenir $\alpha,\beta$ dans~$\gK$ avec $d(\alpha,\beta)\in\gK\eti$, 
ou à conclure que le corps~$\gK_0$ engendré par les \coes de~$f$
et $g$ est un corps fini. Dans ce cas,~$\gK_0[x,z]$
est une~$\gK_0$-\alg finie réduite. 
C'est un anneau réduit fini, donc ou bien c'est un corps fini, 
de la forme~$\gK_0[\gamma]$,  et~$\gK[x,z]=\gK[\gamma]$,
ou bien il contient un \idm $e\neq 0,1$. 

\smallskip 
\rem \Llec pourra vérifier que la transformation de preuve que l'on a fait subir
au cas \gui{$\gB$ est un \cdiz} est exactement la mise en {\oe}uvre de la
machinerie \lgbe \elr des anneaux \zedrsz. 
En fait la même machinerie s'applique aussi pour le \cdiz~$\gK$ et fournit le
résultat suivant: une \alg \ste sur un anneau \zedrz~$\gK$ (\dfnz~\ref{defdualisante}) est un produit fini
de \Klgs \stesz.    
\eoe

%%%%%%%%%%%%%%%%%%%%%%%%%%%%%%%%%%%%%%%%%%%%%%%%%%%%%%%%%%%%%%%%%%%%%%%%%%%

\exer{exoKnEltPrimitif}
\emph{1.} On écrit $x=(\xn)=\sum_{i=1}^nx_ie_i$ et l'on identifie~$\gA$ à un sous-anneau de~$\gB$ par $1\mapsto(1,\ldots,1)$.
En écrivant $e_i \in \gA[x]$, on obtient que les \elts $x_i - x_j$ sont \ivs pour  $j
\ne i$. Réciproquement, si $x_i - x_j$ est \iv pour tout~$i \ne j$,
on a~$\gB = \gA[x]=\gA \oplus \gA x\oplus \cdots \oplus \gA x^{n-1} $ (interpolation de Lagrange, \deter de
Vandermonde).

 \emph{2.} Si et seulement si $\#\gA \ge n$.
%%%%%%%%%%%%%%%%%%%%%%%%%%%%%%%%%%%%%%%%%%%%%%%%%%%%%%%%%%%%%%%%%%%%%%%%%%%

%%%%%%%%%%%%%%%%%%%%%%%%%%%%%%%%%%%%%%%%%%%%%%%%%%%%%%%%%%%%%%%%%%%%%%%%%%%

%%%%%%%%%%%%%%%%%%%%%%%%%%%%%%%%%%%%%%%%%%%%%%%%%%%%%%%%%%%%%%%%%%%%%%%%%%%
\exer{exoIdmChangBase} ~$\!\!\!$
\emph{1} et \emph{2.} Si $(a_1,\ldots,a_\ell)$ est une base de~$\gB$ sur~$\gK$,  c'est aussi une base
de~$\gB(v)$ sur~$\gK(v)$.

\emph{3.}
Soit $b/p$ un \idm de~$\gB(v)$: on a $b^2=bp$. Si $p(0)=0$, alors $b(0)^2=0$, et
puisque~$\gB$ est réduite, $b(0)=0$. On peut alors diviser $b$ et $p$ par
$v$. Ainsi, on peut supposer que $p(0)\in\gK\eti$. En divisant $b$  et $p$ par
$p(0)$ on est ramené au cas où $p(0)=1$. On voit alors que $b(0)$ est
\idmz. On le note $b_0$ et l'on pose~$e_0=1-b_0$. \'Ecrivons $e_0b=vc$. On
multiplie l'\egt $b^2=bp$ par $e_0=e_0^2$ et l'on obtient $v^2c^2=vcp$. Donc
$vc(p-vc)=0$, et puisque le \pol $p-vc$ est de terme constant $1$, donc
régulier, cela donne $c=0$. Donc $b=b_0b$.  Raisonnons un moment modulo
$e_0$: on a $b_0\equiv1$ donc $b$ est primitif et l'\egt $b^2=bp$ se simplifie
en $b\equiv p \mod e_0$. Ceci donne l'\egt $b=b_0b=b_0p$ dans~$\gB(v)$ et donc
$b/p=b_0$.

%%%%%%%%%%%%%%%%%%%%%%%%%%%%%%%%%%%%%%%%%%%%%%%%%%%%%%%%%%%%%%%%%%%%%%%%%%%

\exer{exoBiquadratique} ~\\
\emph{1.} Classique: c'est $\ZZ[\sqrt d]$ si $d \equiv 2$ ou $3 \bmod 4$
et $\ZZ[{1 + \sqrt d \over 2}]$ si $d \equiv 1 \bmod 4$.

\emph{2.} On a~$\gA = \ZZ[\sqrt a]$.
On a $\beta^2 = {a+1 \over 2} + \sqrt a \in \gA$, donc $\beta$ est entier sur
$\gA$, puis sur~$\ZZ$. En fait $(\beta^2 - {a+1 \over 2})^2 = a$ et $\beta$ est
racine de $X^4 - (a+1) X^2 + ({a-1 \over 2})^2$.
On trouve ainsi $(\sigma\tau)(\beta) = \beta - \sqrt 2$.

\emph{3.}
On trouve $(\sigma\tau)(z) = r - (s+u)\sqrt 2 + u\beta$ puis $z +
(\sigma\tau)(z) = 2r + u\sqrt {2a}$. Ce dernier \elt de $\QQ(\sqrt {2a})$ est
entier sur $\ZZ$, donc dans $\ZZ[\sqrt {2a}]$ car $2a \equiv 2 \bmod 4$.
\\
D'où $u \in \ZZ$ (et $2r \in \ZZ$). On remplace $z$ par $z - u\beta$ qui est
entier sur $\ZZ$, \linebreak 
i.e.~$z = r + s\sqrt 2 + t \sqrt a$. On a $\sigma(z) = r-s\sqrt 2
+ t\sqrt a$, $\tau(z) = r+s\sqrt 2 - t\sqrt a$; en utilisant $z+\sigma(z)$ et
$z+\tau(z)$, on voit que $2r$, $2s$, $2t \in \ZZ$. Utilisons:

\snuc {
z\sigma(z) = x + 2rt\sqrt a,\ z\tau(z) = y + 2rs\sqrt 2,
\hbox { avec } x = r^2 - 2s^2 + at^2,\  y = r^2 + 2s^2 - at^2
.}

%\sni
On a donc $x$, $y \in \ZZ$ puis $x + y = 2r^2 \in \ZZ$, $x - y = 2at^2 - (2s)^2
\in \ZZ$, donc $2at^2 \in \ZZ$. De $2r$, $2r^2 \in \ZZ$, on déduit $r \in \ZZ$. De même, de $2t$, $2at^2 \in \ZZ$ (et du fait que $a$ est impair),
on tire $t \in \ZZ$. Et puis enfin $s \in \ZZ$. Ouf!

Gr\^ace à  la $\ZZ$-base de~$\gB$, on obtient $\Disc_{\gB/\ZZ}
= 2^8 a^2$.

\emph{4.}
On a:
%
%\snic {
$\gB = \ZZ \oplus \ZZ\sqrt a \oplus \ZZ\sqrt 2 \oplus \ZZ \beta =
\gA \oplus E$   avec  $E = \ZZ\sqrt 2 \oplus \ZZ \beta
.
$\\%}
%
%\sni
On a aussi $2E = \sqrt2\,\fa$ avec $\fa = 2\,\ZZ  \oplus \ZZ(\sqrt a -1) =
\gen {2, \sqrt a-1}_\gA$. Ceci prouve d'une part que $E$ est un~$\gA$-module,
et d'autre part qu'il est isomorphe à l'idéal~$\fa$ de~$\gA$. En
conséquence, $E$ est un \Amo\ptf de rang $1$. L'écriture~$\gB =\gA\oplus E$
certifie que~$\gB$ est un \Amo\ptfz, écrit comme somme directe d'un
$\gA$-module libre de rang $1$ et d'un \ptf de rang $1$. En général,
l'idéal $\fa$ n'est pas principal, donc $E$ n'est pas un \Amo libre.
Voici un petit échantillon de valeurs de $a \equiv 3 \bmod 4$;
on a souligné quand l'idéal $\fa$ est principal:

\snic {
-33,\
-29,\
-21,\
-17,\
-13,\
-5,\
\und {-1},\
\und {3}, \
\und {7},\
\und {11},\
15,\
\und {19},\
\und {23},\
\und {31},\
35
.}

%\sni
Dans le cas où $\fa$ n'est pas principal,~$\gB$ n'est pas un~$\gA$-module libre: sinon,
$E$ serait stablement libre de rang 1, donc libre
(voir l'exercice~\ref{exoStabLibRang1}).
Enfin, on a toujours $\fa^2 = 2\gA$ (voir la suite), ou encore
$\fa \simeq \fa^{-1} \simeq \fa\sta$.\\ En conséquence~$\gB \simeq_\gA
\gB\sta$. Justification de $\fa^2 = 2\gA$; toujours dans le même
contexte ($a \equiv 3 \bmod 4$ donc~$\gA = \ZZ[\sqrt {a}\,]$), on
a pour $m \in \ZZ$:

\snic {
\gen {m, 1+\sqrt a} \gen {m, 1-\sqrt a} = \pgcd(a-1,m)\,\gA
.}

%\sni
En effet, l'idéal de gauche est engendré par $\big(m^2, m(1 \pm \sqrt a), 1-a\big)$,
tous multiples du pgcd. Cet idéal contient $2m = m(1 + \sqrt a)
+ m(1 - \sqrt a)$, donc il contient l'\elt $\pgcd(m^2, 2m, 1-a) = \pgcd(m, 1-a)$, (l'\egt est
 due à $a \equiv 3 \bmod 4$). Pour $m = 2$, on a
$\gen {2, 1+\sqrt a} = \gen {2, 1-\sqrt a} = \fa$ et $\pgcd(a-1, 2) = 2$.

%%%%%%%%%%%%%%%%%%%%%%%%%%%%%%%%%%%%%%%%%%%%%%%%%%%%%%%%%%%%%%%%%%%%%%%%%%%

\exer{exoTensorielDiscriminant}
On voit~$\gB = \gA\otimes_\gk\gA'$ comme \Algz, extension des scalaires à
$\gA$ de la \klgz~$\gA'$; elle est libre de rang $n'$. On dispose donc d'une
tour d'\algs libres~$\gk \to \gA \to \gB$ et la formule de transitivité du
discriminant fournit:

\snic {
\Disc_{\gB\sur\gk}(\ux\otimes\ux') = \Disc_{\gA\sur\gk}(\ux)^{n'}
\cdot \rN_{\gA\sur\gk}\big(\Disc_{\gB\sur\gA}(1 \otimes \ux')\big)
.}

%\sni
Mais $\Disc_{\gB\sur\gA}(1 \otimes \ux') = \Disc_{\gA'\sur\gk}(\ux')$. Comme c'est
un \elt de~$\gk$, sa norme~$\rN_{\gA\sur\gk}$ vaut~$\Disc_{\gA'\sur\gk}(\ux')^n$.  En fin de compte on obtient l'\egt

\snic {
\Disc_{\gB\sur\gk}(\ux\otimes\ux') =
\Disc_{\gA\sur\gk}(\ux)^{n'} \, \Disc_{\gA'\sur\gk}(\ux')^n.
}

%%%%%%%%%%%%%%%%%%%%%%%%%%%%%%%%%%%%%%%%%%%%%%%%%%%%%%%%%%%%%%%%%%%%%%%%%%%
\exer{CyclicNormalBasis} 
On va utiliser le résultat classique d'\alg\lin suivant. Soit~$E$ un \Kev de dimension finie  et $u\in\End_\gK(E)$. Si $d$ le degré du \polmin de $u$,  il
existe $x \in E$ tel que les \elts  
$x$, $u(x)$, \ldots, $u^{d-1}(x)$
soient~$\gK$-\lint indépendants.\\
Ici $\dex{\gL : \gK} = n$, et $\Id_\gL$, $\sigma$, \dots, $
\sigma^{n-1}$, sont~$\gK$-\lint indépendants, donc le
\polmin de $\sigma$ est $X^n - 1$, de degré $n$. On applique le 
résultat ci-dessus.

%%%%%%%%%%%%%%%%%%%%%%%%%%%%%%%%%%%%%%%%%%%%%%%%%%%%%%%%%%%%%%%%%%%%%%%%%%%

\exer{exoHomographieOrdre3}
\emph{1.} $A$ est la matrice compagne du \pol $X^2 - X + 1 = \Phi_6(X)$, 
donc $A^3=-\I_2$  dans
$\GL_2(\gk)$ et $A^3 = 1$ dans $\PGL_2(\gk)$. 
%On peut également
%considérer $\ZZ[j]$ ($j$ racine cubique de l'unité); la matrice de la
%multiplication par $-j$ dans la $\ZZ$-base $(1, -j)$ de $\ZZ[j]$ est $A$.

\emph{2.}
On sait par le \tho d'Artin que~$\gk(t)/\gk(t)^G$ est une extension galoisienne de groupe de Galois
$\rA_3$.
Le calcul donne 

\snic{g = t+\sigma(t)+\sigma^2(t)= {t^3 - 3t - 1 \over t(t+1)}.}

%\sni
On a évidemment $g\in\gk(t)^G$ et $t^3 - gt^2 - (g+3)t - 1=0$. 
Donc, (partie directe du \tho de Lüroth, \pbz~\ref{exoLuroth1}) $\dex{\gk(t):\gk(g)}=3$,
et~$\gk(t)^G=\gk(g)$.

 \emph{3.}
Puisque~$\gk(a)\simeq\gk(g)$ et $f_g(t)=0$, 
l'extension~$\gk(a)\to\aqo{\gk[T]}{f_a}$ est une photocopie de 
l'extension~$\gk(g)\to\gk(t)$.

 \emph{4.}
Soit $\sigma$ un \gtr  de $\Aut(\gL/\gK)$.  Cette question revient à
dire que l'on  peut trouver un $t\in\gL\setminus \gK$ 
tel que $\sigma(t)={-1\over t+1}~(*)$.
Puisque $t$ doit être de norme $1$, on le cherche de la forme 
$t={\sigma(u)\over u}$.
Le calcul montre alors que $(*)$ est satisfaite à condition que $u\in\Ker(\Tr_G)$. Il reste à montrer qu'il existe un $u\in\Ker(\Tr_G)$ tel que ${\sigma(u)\over u}\notin\gK$.
Cela revient à dire que la restriction de $\sigma$ à $E =\Ker(\Tr_G)$
n'est pas une homothétie. Or $E\subseteq\gL$ est un sous-\Kev de dimension $2$,
stable par $\sigma$.
D'après
l'exercice~\ref{CyclicNormalBasis}, le \Kevz~$\gL$ admet un \gtr pour
l'\endo $\sigma$.
Cette \prt d'\alg \lin reste vraie pour tout sous-espace stable par $\sigma$.

%%%%%%%%%%%%%%%%%%%%%%%%%%%%%%%%%%%%%%%%%%%%%%%%%%%%%%%%%%%%%%%%%%%%%%%%%%%
\exer{exoGroupAlgebra} 
Les \elts $g \te h$ forment une~$\gk$-base de $\env{\gk}{\gA}$.  \\
Soit $z =
\sum_{g,h} a_{g,h} g\te h$ avec $a_{g,h} \in \gk$.  Alors, $z \in
\Ann(\rJ\iAk)$ \ssi on~a~$g' \cdot z = z \cdot g'$ pour tout $g' \in G$. On obtient
$a_{g,h} = a_{1,gh}$, donc $z$ est combinaison~$\gk$-\lin des $z_k
\eqdf {\rm def} \sum_{gh = k} g \otimes h$.  \\
Réciproquement, on voit que
$z_k \in \Ann(\rJ\iAk)$ et l'on a $z_k = k \cdot z_1 = z_1 \cdot k$. \\
Donc
$\Ann(\rJ\iAk)$ est le \kmo engendré par les $z_k$, et c'est le \Amo (ou
l'\id de $\env{\gk}{\gA}$) engendré par $z_1 = \sum_g g\te g^{-1}$.
\\
L'image par $\mu\iAk$ de $\Ann(\rJ\iAk)$ est l'\id $n\gA$.\\
Concernant la trace, on a $\Tr\iAk(g) = 0$ si $g \ne 1$% et $\Tr\iAk(1) = n$
.
Donc $\Tr\iAk\!\big(\sum_g a_g g\big) = na_1$. 
\\
Si $a = \sum_g a_g g$ et $b = \sum_g b_g g$, alors
$\Tr\iAk(ab) = n\sum_g a_g b_{g^{-1}}$.
\\
Les \eqvcs du point \emph{2} sont donc claires, et dans le cas où $n\in\gk\eti$, l'\idst est
$n^{-1}\sum_g g\te g^{-1}$.

\emph {3.}
Soit $\lambda : \gk[G] \to \gk$ la forme \lin \gui{\coo sur $1$}. %\\
% Alors, p
 Pour
$g$, $h \in G$, on~a~$\lambda(gh) = 0$ si $h \ne g^{-1}$, et~$1$ sinon. 
%Ceci prouve que
Donc, $\lambda$ est dualisante et %que 
$(g^{-1})_{g \in G}$ est la base duale de
$(g)_{g \in G}$ relativement à $\lambda$.  
On a $\Tr_{\gk[G]\sur\gk} =
n\cdot\lambda$.

%%%%%%%%%%%%%%%%%%%%%%%%%%%%%%%%%%%%%%%%%%%%%%%%%%%%%%%%%%%%%%%%%%%%%%%%%%%
\exer{exo1Frobenius}
\emph{1.} 
Il suffit de le faire pour $g \in \{1, x, \ldots, x^{n-1}\}$, qui est une base de
$\gA[Y]$ sur $\gk[Y]$.  
%\\
Le membre droit de $(*)$ avec $g=x^{i}$ est 

\snic{h_i=\varphi\big((Y-x)x^i\big) =
\varphi(Yx^i -x^{i+1}) = Y^{i+1} - \varphi(x^{i+1}).}

%\sni
 Si $i < n-1$, on a
$\varphi(x^{i+1}) = Y^{i+1}$, donc $h_i=0$. 
 Pour $i = n-1$, on a

\snic {
\varphi(x^n) = -\varphi(a_0 + a_1x + \cdots + a_{n-1}x^{n-1}) = 
-(a_0 + a_1Y + \cdots + a_{n-1}Y^{n-1}),
}

%\sni
et $h_n(Y)=f(Y)$, ce qui permet de conclure.

 \emph{2.}
Pour $i < n$, on a 

\snic{ 
\begin{array}{ccc} 
f(Y)\wi\lambda\big(x^i f'(Y)\big) =
\varphi\big((Y-x)x^i f'(Y)\big) =  \varphi\big(x^i f(Y)\big) = Y^if(Y), \hbox{ i.e.}  \\[1mm] 
\wi\lambda\big(x^i f'(Y)\big) = \sum_{j < n} \lambda(x^i b_j) Y^j = Y^i. 
\end{array}
}

%\sni
Donc $(b_j\centerdot\lambda)(x^i)=\lambda(x^i b_j) = \delta_{ij}$. Ainsi, $\lambda$ est
dualisante.

 \emph{3.}
Pour deux bases duales $(e_i)$, $(\alpha_i)$, on a
$\Tr\iAk = \sum e_i\centerdot\alpha_i$. Avec les deux bases duales $(1, x, \ldots, x^{n-1})$ et
$(b_0\centerdot\lambda, b_1\centerdot\lambda, \ldots, b_{n-1}\centerdot\lambda)$ on obtient:

\snic {
\Tr\iAk =  b_0\centerdot\lambda + x b_1\centerdot\lambda + \cdots + 
x^{n-1} b_{n-1}\centerdot \lambda = f'(x)\centerdot\lambda.
}

%%%%%%%%%%%%%%%%%%%%%%%%%%%%%%%%%%%%%%%%%%%%%%%%%%%%%%%%%%%%%%%%%%%%%%%%%%%

\exer{exoFrobeniusAlgExemples} 
\emph {1.}
La \klgz~$\gA := \gk[X_1, \ldots, X_n]\sur{\gen {f_1(X_1), \ldots, f_n(X_n)}}$
est le produit tensoriel des $\aqo{\gk[X_i]}{f_i(X_i)}$ qui sont des \algs de
Frobenius, donc~$\gA$ est une \alg de Frobenius.
Précision avec $d_i = \deg(f_i)$. La \klg $\gA$ est libre  de rang
$d_1 \cdots d_n$,  les monômes $x^{\alpha}=x_1^{\alpha_1} \cdots x_n^{\alpha_n}$ avec $\alpha_i < d_i$ pour tout $i$ forment une~$\gk$-base, et la forme
$\gA$-\lin  \gui{\coo sur $x_1^{d_1-1} \cdots x_n^{d_n-1}$} est
dualisante.

\emph {2.}
Soient $\delta_0$, $\delta_x$, $\delta_y$ les trois formes \lins sur~$\gk[X,Y]$
définies par

\snic {
\delta_0(f) = f(0), \quad \delta_x(f) = f'_X(0), \quad \delta_y(f)=f'_Y(0).
}

%\sni
Elles passent au quotient sur~$\gA$ et définissent une~$\gk$-base
de $\Asta$, base duale de la~$\gk$-base $(1, x, y)$ de~$\gA$. On~a:

\snic {
x \centerdot \delta_x = y \centerdot \delta_y = \delta_0, 
}

%\sni

et donc $\Asta = \gA\centerdot\delta_x + \gA\centerdot\delta_y$.  Montrons
que $G = \cmatrix {x\cr -y\cr}$ est une \mpn de $\Asta$ pour $(\delta_x, \delta_y)$% (comme
%\Amoz), $G : \gA \to \gA^2$, est 
. 
Il faut voir que pour $u% = u(x,y)
$, $v% = v(x,y)
$ dans~$\gA$ on a l'implication

\snic {
u \centerdot \delta_x + v \centerdot \delta_y = 0 \;\;\Longrightarrow\;\;
\cmatrix {u\cr v\cr} \in \gA \cmatrix {x\cr -y\cr}.
}

%\sni
En multipliant $u \centerdot \delta_x + v \centerdot \delta_y = 0$ par $x$, on
obtient $u \centerdot \delta_0 + (xv) \centerdot \delta_y = 0$; on évalue en
$1$ et l'on fait $x := 0$ pour obtenir $u(0,y) = 0$, i.e.  $u \in \gA x$. De
même, $v \in \gA y$. Si l'on écrit $u = xr$, $v = ys$, on obtient $r
\centerdot \delta_0 + s \centerdot \delta_0 = 0$, i.e. $r + s = 0$, ce que
l'on voulait.

L'\idd $\cD_1(G) = \gen {x,y}$ est non nul, de carré nul,
donc il ne peut pas être engendré par un \idmz. En conséquence, le \Amo
$\Asta$ n'est pas  \proz. A fortiori, il n'est pas libre.

%%%%%%%%%%%%%%%%%%%%%%%%%%%%%%%%%%%%%%%%%%%%%%%%%%%%%%%%%%%%%%%%%%%%%%%%%%%
\exer{exoAlgMonogeneJAnnJ} 
\emph {1.}
On a $(y-z)f^\Delta(y,z) = 0$ donc $\delta := f^\Delta(y,z) \in \Ann(\rJ)$.
On sait alors que pour $\alpha \in \env{\gk}{\gA}$, on a $\alpha\delta =
\mu\iAk(\alpha) \cdot \delta = \delta \cdot \mu\iAk(\alpha)$. On applique ce
résultat
à $\alpha = \delta$ en remarquant que $\mu\iAk(\delta) = f'(x)$.

 \emph {2.}
On écrit
$f(Y) - f(Z) = (Y-Z)f'(Y) - (Y-Z)^2 g(Y,Z)$, ce qui donne dans
l'\alg $\env{\gk}{\gA}$ l'\egt $(y-z)f'(y) = (y-z)^2g(y,z)$. 
On écrit l'\egt $1=uf+vf'$ dans~$\AX$.
Alors $f'(y)v(y)=1$, donc $y-z = (y-z)^2v(y)g(y,z)$. \\
Lorsque $a = a^2b$, l'\elt $ab$ est
\idm et~$\gen {a} = \gen {ab}$. Donc~$\rJ=\gen{e}$ avec l'\idm $e = (y-z)v(y)g(y,z)$.
\\
On a $f^\Delta(Y,Z) = f'(Y) - (Y-Z)g(Y,Z)$, donc

\snic {
f^\Delta(y,z) = f'(y) - (y-z)g(y,z) = f'(y)\big(1  - (y-z)v(y)g(y,z)\big) = 
f'(y) (1-e). 
}
%%%%%%%%%%%%%%%%%%%%%%%%%%%%%%%%%%%%%%%%%%%%%%%%%%%%%%%%%%%%%%%%%%%%%%%%%%%

\exer{exoJacInversibleAlgSeparable} 
\emph {1.}
On note $f'_{ij}=\partial f_i/\partial X_j$ et l'on
écrit dans~$\gk[\uY,\uZ]$:

\snic {
f_i(\uY) - f_i(\uZ) - \sum_j (Y_j - Z_j) f'_{ij}(\uY) =: - g_i(\uY,\uZ) \in \gen {Y_1-Z_1, \ldots,
Y_n-Z_n}^2.
}

\sni Dans $\env{\gk}{\gA}$, on obtient, en posant $A=\JJ_\uy(f_1, \ldots, f_n)$:

\snuc {
A \cmatrix {y_1-z_1\cr \vdots\cr y_n-z_n\cr} =
\cmatrix {g_1(\uy,\uz) \cr \vdots\cr g_n(\uy,\uz)\cr}
\quad \hbox {avec} \quad
g_i(\uy,\uz) \in \gen {y_1-z_1, \ldots, y_n-z_n}^2 = \rJ\iAk^2.
}

En inversant $A$, on obtient $y_i - z_i \in \rJ\iAk^2$, i.e.
$\rJ\iAk = \rJ\iAk^2$.

\emph {2.}
Comme $\mu\iAk\big(\beta_{\uy,\uz}(\uf)\big) = \J_\ux(\uf)$, on
a $\mu\iAk(\vep) = 1$. \\
Comme $\vep \in \Ann(\rJ\iAk)$, $\vep$
est le \gtr \idm de $\Ann(\rJ\iAk)$. \\
Enfin,~$\beta_{\uy, \uz}(\uf)$, qui
est associé à $\vep$, est aussi un \gtr de $\Ann(\rJ\iAk)$.

\emph {3.}
Soient $f_1$, \ldots, $f_n$  dans $\gk[\uX]$ et
$\delta = \J_{\uX}(\uf)$. Inversons $\delta$ 
avec une \idtrz~$T$. Alors, dans $\gk[\uX,T]$ on obtient
$$
\JJ_{\uX,T}(\uf,\delta T-1) =
\bordercmatrix [\lbrack\rbrack]{
             &\partial_{X_1}  &\cdots         &\partial_{X_n} &\partial_T\cr
\quad f_1    &                &               &       &0 \cr
\ \quad\vdots&                &\JJ_{\uX}(\uf) &       &\vdots\cr
\quad f_n    &                &               &       &0\cr
\delta T-1   &\star           &\cdots         &\star  &\delta\cr
}
$$ 
et donc $\J_{\uX,T}(\uf,\delta T-1) = \delta^2$.  Notons 

\snic{\gA =
\aqo{\gk[\uX]}{\uf}\;\hbox{  et  }\;\gB = \gA[\delta^{-1}] =
\aqo{\gk[\uX,T]}{\uf,1-\delta T}.}

Alors le jacobien du \sys $(\uf,\delta T-1)$
qui définit $\gB$ est \iv dans $\gB$ et donc $\gB$ est une \alg\spbz.

%%%%%%%%%%%%%%%%%%%%%%%%%%%%%%%%%%%%%%%%%%%%%%%%%%%%%%%%%%%%%%%%%%%%%%%%%%%
%\exer{exoBézoutienJacobienIdm} 
%\emph {1.}
%Il suffit de montrer que $(y_j - z_j)\beta(\uy, \uz) = 0$ pour tout $j$.
%Ceci découle aussitôt de l'\egt $(\star)$.
%
%\emph {2.}
%Comme $\mu\iAk\big(\beta(\uy,\uz)\big) = \J_\ux(f_1, \ldots, f_n)$, on
%a $\mu\iAk(\vep) = 1$. \\
%Comme $\vep \in \Ann(\rJ\iAk)$, $\vep$
%est le \gtr \idm de $\Ann(\rJ\iAk)$. \\
%Enfin,~$\beta(\uy, \uz)$, qui
%est associé à $\vep$, est aussi un \gtr de $\Ann(\rJ\iAk)$.

%%%%%%%%%%%%%%%%%%%%%%%%%%%%%%%%%%%%%%%%%%%%%%%%%%%%%%%%%%%%%%%%%%%%%%%%%%%
\exer{exoSepAlgDedekindLemma} 
La~$\gB$-\algz~$\gB\te_\gk\gA$ est \spbz. On a une transformation (\prt \uvle de l'\edsz)

\snic{\Hom_\gk(\gA,\gB) \to \Hom_\gB(\gB\te_\gk\gA,\gB)$, $\psi \mapsto
\ov\psi,}

%\sni
définie par $\ov\psi(b\otimes a) = b\psi(a)$.

\emph {1.}
On considère alors l'\idm $\vep_{\ov\varphi} \in
\gB\te_\gk\gA$ du lemme~\ref {lemIdmHomSpb}, et on l'écrit
sous la forme
$\vep_{\ov\varphi} = \sum_{i \in I} b_i\te a_i$.

\emph {2.}
Découle directement du lemme~\ref {lemIdmHomSpb}: l'\idm
$e$ n'est autre que $e_{\{\ov\varphi, \ov{\varphi'}\}}$.

\emph {3.}
Puisque la juxtaposition horizontale de matrices d'évaluation de Dedekind
est une matrice d'évaluation de Dedekind, il suffit de montrer qu'il en
existe une, disons $A_1$, dont l'image contient le vecteur 
$v := \tra{[ \,e_{11} \;\cdots \; e_{n1}\,]}$. \\
 Soit $\big((a_j)_{j \in \lrbm}, (b_j)_{j \in \lrbm}\big)$ le
couple attaché à $\varphi_1$.
On met en colonne $j$ de~$A_1$
le vecteur $\tra {[\,\varphi_1(a_j) \; \cdots \; \varphi_n(a_j)\,]}$. On
a alors $A_1 \tra {[ \,b_1 \; \cdots \; b_m\,]} = v$.

%%%%%%%%%%%%%%%%%%%%%%%%%%%%%%%%%%%%%%%%%%%%%%%%%%%%%%%%%%%%%%%%%%%%%%%%%%%

%%%%%%%%%%%%%%%%%%%%%%%%%%%%%%%%%%%%%%%%%%%%%%%%%%%%%%%%%%%%%%%%%%%%%%%%%%%
\exer{exoArtinAnOtherProof}
Par hypotèse, pour chaque $\tau\in G \setminus \{\Id\}$ il existe $n_\tau\in\NN$ et des \elts $x_{1,\tau}$, \ldots,
$x_{n_\tau,\tau}$, $y_{1,\tau}$, \ldots, $y_{n_\tau,\tau} \in \gA$ tels que 
$1 = \sum_{j=1}^{n_\tau} x_{j,\tau}\big(y_{j,\tau} - \tau(y_{j,\tau})\big)$.
On pose $s_\tau = \sum_{j=1}^{n_\tau} x_{j,\tau} \tau(y_{j,\tau})$ de sorte 
que~$\sum_{j=1}^{n_\tau} x_{j,\tau}y_{j,\tau} = 1+s_\tau$, puis on définit $x_{n_\tau+1,\tau} =
-s_\tau$ et $y_{n_\tau+1,\sigma} = 1$. Alors, en posant $m_\tau=n_\tau+1$:

\snic {
\sum_{j=1}^{m_\tau} x_{j,\tau} \tau(y_{j,\tau}) = s_\tau - s_\tau = 0, \qquad
\sum_{j=1}^{m_\tau} x_{j,\tau} y_{j,\tau} = 1+s_\tau - s_\tau = 1.
}

%\sni
On fixe un $\sigma\in G$ et on fait le produit, on obtient

\snic{
\prod_{\tau \in G \setminus \{\Id\}}
  \sum_{j=1}^{m_\tau} x_{j,\tau} \sigma(y_{j,\tau})=
\formule{1\;\;\mathrm{si}\;\sigma=\Id\\
0\;\;\mathrm{sinon}.}
}

%\sni
Le développement du produit fournit deux familles $(a_i)$ et $(b_i)$
indexées par le même ensemble (chaque $a_i$ est un produit de certains $x_{j,\tau}$ et
$b_i$ est le produit des~$y_{j,\tau}$  correspondants) vérifiant:

\snic{\sum_{i=1}^r a_i \sigma(b_i)=
\formule{1\;\;\mathrm{si}\;\sigma=\Id\\
0\;\;\mathrm{sinon}.}
}

%%%%%%%%%%%%%%%%%%%%%%%%%%%%%%%%%%%%%%%%%%%%%%%%%%%%%%%%%%%%%%%%%%%%%%%%%%%

\exer{exoGalExemples} 
\emph {1.}
Comme $G$ agit transitivement sur $\lrb {1..n}$, on a $(\gk^n)^G = \gk$.  De
plus, $G$ étant de cardinal $n$, une permutation $\sigma \in G\setminus\so{\Id}$ %autre que $\Id$ 
n'a aucun point fixe.  On en déduit que $\sum_{\sigma \in G}
e_i\sigma(e_i) = 0$ si $\sigma \in G \setminus \{\Id\}$, et $1$ sinon.\\
 En prenant
$x_i = y_i = e_i$, les conditions du lemme~\ref{lemArtin} sont satisfaites et
$(\gk, \gk^n, G)$ est une \aGz.

L'application $G \to \lrb{1..n}$, $\sigma \mapsto \sigma(1)$, est une
bijection. L'action de $G$ sur $\lrb {1..n}$ est \ncrt isomorphe à l'action
de $G$ sur lui-même par translations.  Si $n$ est fixé, on peut prendre 
pour $G$ le groupe engendré par un $n$-cycle.

\emph {2.}
On a $\Stp_{\rS_4}(\gB) = \gen {(1,2), (3,4)}$ et $H = \Stp_G(\gB) = \{\Id\}$;
donc $(\gk^4)^H = \gk^4$.

\emph {3.}
Le premier point est immédiat. Supposons~$\gB = \gA[X]^H$ et soit $a \in
\gA$.\\
Alors $aX \in \gB$, donc $aX$ est invariant par $H$, i.e. $a$ est
invariant par $H$. \\
Bilan~:~$\gA = \Ae H$ donc $H = \{\Id\}$ puis~$\gA[X] =
X\gA[X] + \gk$ i.e.~$\gA = \gk$ et $G = \{\Id\}$. Hormis ce cas très
particulier,~$\gB$ n'est pas de la forme~$\gA[X]^H$.

%%%%%%%%%%%%%%%%%%%%%%%%%%%%%%%%%%%%%%%%%%%%%%%%%%%%%%%%%%%%%%%%%%%%%%%%%%%

\exer{lemPaquesFerrero} 
On suppose \spdg que~$\gB$ et~$\gC$ sont libres de rang $n\in\NN$: il suffit en effet de vérifier la conclusion après
\lon en des \eco  et l'on dispose du \tho de structure locale des \mptfsz.
Si $n=0$  alors~$\gk$ est trivial,
on peut donc supposer $1\leq n$.
 On considère une base $\cC=(c_1,\ldots,c_n)$
de~$\gC$  et une base $\cB=(b_1,\ldots,b_n)$ de~$\gB$ (sur~$\gk$),
et l'on écrit la matrice
$B\in\Mn(\gk)$ de $\cB$ sur $\cC$. Le fait que les $b_i$ forment une base implique
que $B$ est injective, i.e. $\delta=\det B$ est \ndz
(\thref{prop inj surj det}). En outre,
$\delta\gC\subseteq\gB$.
\\
Comparons  $\Tr\iBk(x)$ et $\Tr\iCk(x)$
pour un $x\in\gB$. Considérons~$\gk'=\gk[1/\delta]\supseteq\gk$.
Les deux~$\gk'$-\algs obtenues par \edsz,~$\gB[1/\delta]$ et~$\gC[1/\delta]$, sont les mêmes, et la trace se comporte bien par \edsz,
donc $\Tr\iBk(x)$ et $\Tr\iCk(x)$ sont égales parce qu'elles sont égales dans~$\gk'$.
Mais alors

\snic{\disc \gB\sur\gk=\disc\iBk(\cB)=\disc\iCk(\cB)=
\delta^{2}\disc\iCk(c_1,\ldots,c_n).}

%\sni
Enfin puisque $\disc \gB\sur\gk$ est \ivz, $\delta$ \egmt et~$\gB=\gC$.

%%%%%%%%%%%%%%%%%%%%%%%%%%%%%%%%%%%%%%%%%%%%%%%%%%%%%%%%%%%%%%%%%%%%%%%%%%%

\exer{exoCorGalste}
Tout d'abord, notons que puisque~$\gk$ est connexe, tous les \mptfs sur~$\gk$
sont de rang constant. Rappelons aussi que la correspondance galoisienne est déjà établie lorsque~$\gk$ est un \cdiz.
\\
Nous devons montrer que si~$\gk\subseteq\gB\subseteq\gA$ avec~$\gB$ \stez, alors

\snic{ \gB=\gC\eqdefi\Fix\big(\Stp(\gB)\big).}

%\sni
D'après le lemme~\ref{lemPaquesFerrero}, il suffit de montrer que~$\gB$ et~$\gC$ ont même rang. En \clama on conclut en notant qu'après \eds à n'importe quel
corps,~$\gB$ et~$\gC$ ont même rang puisque la correspondance galoisienne est  établie pour les corps.
\\
Cet argument de \clama fournit par relecture dynamique une \prcoz.
Ceci est lié au \nst formel (\thref{thNSTsurZ}).
\perso{fin de solution à rédiger!!!!}

%%%%%%%%%%%%%%%%%%%%%%%%%%%%%%%%%%%%%%%%%%%%%%%%%%%%%%%%%%%%%%%%%%%%%%%%%%%

\exer{exoIdeauxGlobInvariants} 
Soient $(x_i)$, $(y_i)$ deux \syss d'\elts de~$\gB$ comme dans le lemme~\ref{lemArtin}.

\emph {1.}
On sait que pour $x \in \gB$, $x = \sum_i \Tr_G(xy_i) x_i$. Si
$x \in \fb$, alors $xy_i \in \fb$, et comme~$\fb$ est globalement
invariant, $\Tr_G(xy_i) \in \fb$. \\
Bilan~: $\fb$ est engendré
par les \elts invariants $\Tr_G(xy_i)$ pour $x \in \fb$.

\emph {2.}
Soit $\fa$ un idéal de~$\gA$; il est clair que $\fa\gB$ est globalement
invariant.\\
 Il faut voir que $\fa\gB \cap \gA = \fa$. Cela vient du fait que
$\gA$ est facteur direct dans~$\gB$ (comme~\Amoz).  En effet, soit~$\gB = \gA
\oplus E$, donc $\fa\gB = \fa \oplus \fa E$. Si $x \in \fa\gB \cap \gA$, on
écrit $x = y+z$ avec $y \in \fa$ et $z \in \fa E \subseteq E$; on a alors
$x$, $y \in \gA$, donc $z \in \gA$, et comme $z \in E$, $z = 0$. En
conséquence, $x = y \in \fa$.

Réciproquement, si $\fb \subseteq \gB$ est globalement invariant,
il faut voir que $(\fb\cap \gA)\gB = \fb$; mais c'est ce qui a été
montré dans la question précédente.

%%%%%%%%%%%%%%%%%%%%%%%%%%%%%%%%%%%%%%%%%%%%%%%%%%%%%%%%%%%%%%%%%%%%%%%%%%%
%: sols pb
%%%%%%%%%%%%%%%%%%%%%%%%%%%%%%%%%%%%%%%%%%%%%%%%%%%%%%%%%%%%%%%%%%%%%%%%%%%

\prob{exoBuildingFrobAlgebra} 
De manière \gnlez, la forme \lin $\delta_g$ passe au quotient modulo l'\id
$\fa_g$ qu'elle définit. De plus, si $\delta_g(\ov u\,\ov v) = 0$ sur~$\gA =
\kuX\sur{\fa_g}$ pour tout $\ov v$, alors $\delta_g(u v) = 0$ pour tout $v$,
donc $u \in \fa_g$, i.e. $\ov u = 0$. Donc $\Ann_\gA(\delta_g) = 0$.

Pour $i \in \lrbn$, on note $\delta_i^m = \delta_{X_i^m}$ (\coo sur
$X_i^m$). \\
En particulier, $\delta_i(f) = {\partial f \over \partial
X_i}(0)$. Et l'on définit $\delta_0 : \kuX \to \gk$ par $\delta_0(f) = f(0)$.

\emph {1.}
Calcul facile.

\emph {2.}
On vérifie que $f*g = 0$ \ssi $f_m*g = 0$ pour toute
composante \hmg $f_m$ de~$f$, autrement dit l'\id $\fa_g$ est
\hmg (c'est toujours le cas si~$g$ est \hmgz).\\
 Il est clair aussi que pour $i
\ne j$, $X_iX_j * g = 0$, et pour $|\alpha| > d$, $X^\alpha * g = 0$. 
\linebreak 
Si $f = \sum_i
a_iX_i^m + \cdots$ est \hmg de degré $m \le d$, on a $f * g = \sum_i
a_i X_i^{-(d-m)}$. 
\\
 Si $m < d$, on a donc $f*g = 0$ \ssi $a_i = 0,\;\Tt i$, 
 \cade \hbox{si $f \in \gen {X_iX_j, i \ne j}$}.
\\
 Si $m = d$, on a $f*g = 0$ \ssi
$\sum_i a_i = 0$, \cade \hbox{si $f \in \gen {X_iX_j, i \ne j} +
\gen {X_i^d - X_1^d, i \in \lrb {2..n}}$}, car $\sum_i a_iX_i^d = \sum_i a_i(X_i^d - X_1^d)$. 
\\
Bilan: on a obtenu un \sgr de $\fa_g$ constitué de $n(n-1) \over 2$ \pogs de degré $2$ et de $n-1$ \pogs de degré $d$:

\snic {
\fa_g = \gen {X_iX_j, i < j} + \gen {X_i^d - X_1^d, i \in \lrb{2..n}}
.}

%\sni
On pose~$\gA = \kuX\sur\fa_g = \gk[x_1, \ldots, x_n]$. Alors:

\snic {
1,\quad x_1,\ \ldots,\ x_n,\quad x_1^2,\ \ldots,\ x_n^2, 
\quad\ldots\quad
x_1^{d-1},\ \ldots,\ x_n^{d-1}, \quad x_1^d
}

%\sni
est une~$\gk$-base de~$\gA$ de cardinal $(d-1)n+2$.  La~$\gk$-base duale
de $\Asta$ est:

\snic {
\delta_0,\quad \delta_1,\ \ldots,\ \delta_n,\quad \delta_1^2,\ \ldots,\ \delta_n^2, 
\quad\ldots\quad
\delta_1^{d-1},\ \ldots,\ \delta_n^{d-1}, \quad \delta_g
}

%\sni
et l'on a:

\snic {
x_i^m \centerdot\delta_g = \delta_i^{d-m} \hbox { pour $m \in \lrb{1..d-1}$},
\quad 
x_i^d \centerdot\delta_g = \delta_0
.}

%\sni
Donc $\Asta = \gA\centerdot\delta_g$ et $\delta_g$ est dualisante.

\emph {3.}
Si l'on prend $e_i$ strictement plus grand que l'exposant de $X_i$ dans
l'ensemble des monômes de $g$, on a $X_i^{e_i} * g = 0$.

\emph {4.}
Soit $f \in \kuX$. On a vu que $f\centerdot \delta_g = 0$
\ssi $\partial_f(g) = 0$. L'application~$\gk$-\lin $\kuX \to \kuX$, $f
\mapsto \partial_f(g)$, passe au quotient modulo $\fb$ pour
définir une application~$\gk$-\lin $\varphi$.

\emph {5.}
L'application~$\gk$-\lin~$\gA \to \Asta$, $f \mapsto f\centerdot \delta_g$,
est injective et comme~$\gA$ et~$\Asta$ sont \hbox{des~$\gk$-\evcsz} de même
dimension finie, c'est un \isoz.
%%%%%%%%%%%%%%%%%%%%%%%%%%%%%%%%%%%%%%%%%%%%%%%%%%%%%%%%%%%%%%%%%%%%%%%%%%%

\prob{exoTh90HilbertAdditif}
\emph{1.}  On pose comme par magie $\theta(x) =
\sum_{i=0}^{n-1} \sigma^i(z) c_i(x)$ (merci Hilbert\ihiz). On va vérifier que:

\snic {
\sigma\big(\theta(x)\big) = \theta(x) + \Tr_G(x)z - x
\quad \hbox {ou encore} \quad
x = (\Id_\gA - \sigma)\big(\theta(x)\big) + \Tr_G(x)z.
}

%\sni
Donc,
$(\Id_\gA - \sigma) \circ \theta$ et $x \mapsto
\Tr_G(x)z$ sont deux \prrs \orts de somme~$1$, \hbox{d'où $\gA =
\Im(\Id_\gA - \sigma) \oplus \gk z$}. Pour la \vfnz, notons $c_i$ pour
$c_i(x)$ et~$y = \theta(x)$. On a $\sigma(c_i) = c_{i+1} - x$, $c_n =
\tr_G(x)$ et

\snic {
\begin {array} {rcl}
\sigma(y) &=& \sum_{i=0}^{n-1} (c_{i+1} - x) \sigma^{i+1}(z) =
\sum_{i=0}^{n-1} c_{i+1} \sigma^{i+1}(z) -
\sum_{i=0}^{n-1} x\sigma^{i+1}(z)   \\[1mm]
&=& (y + \Tr_G(x)z) - x\Tr_G(z) = y + \Tr_G(x)z - x .
\end {array}
}

%\sni
Puisque $\Tr_G(z)=1$, $z$ est une base de~$\gk z$
(si $az=0$, \hbox{alors $0=\Tr_G(az)=a$}), donc $\Im(\Id_\gA - \sigma)$ est
bien \stl de rang $n-1$.

 \emph{2.} Il est clair que $\Im(\Id_\gA - \sigma) \subseteq \Ker\Tr_G$. L'autre inclusion résulte du point précédent.

 \emph{3.}
\Llec fera les \vfns en posant $y =
\sum_\tau c_\tau \tau(z)$. Il y a un lien avec la question \emph{1}: pour $x$
fixé avec $\Tr_G(x) = 0$, la famille $\big(c_i(x)\big)$ est un
$1$-cocycle additif à condition d'identifier $\lrb{0..n-1}$ et $G$ via $i
\leftrightarrow \sigma^i$.

\emph {4.}
L'\elt $-1$ est de trace nulle, d'où l'existence de $y\in\gA$ tel que
$-1=y-\sigma(y)$. On a alors, pour tout $i \in \ZZ$, $\sigma^i(y) = y+i$%
%. Comme l'application $i \mapsto i\cdot 1_\gk$ induit une injection 
% $\FF_p \hookrightarrow \gk$, 
, et $\sigma^j(y) - \sigma^i(y) = j-i$ est \iv pour $i \not\equiv j
\bmod p$.
\\ 
Notons $y_i = \sigma^i(y)$, $(i \in \lrb{0..p-1})$. La matrice
de Vandermonde de $(y_0, y_1, \ldots, y_{p-1})$ est \iv et par suite
$(1, y, \ldots, y^{p-1})$ est \hbox{une~$\gk$-base} de~$\gA$. On note~$\lambda = y^p
- y$. Alors $\lambda\in\gk$ puisque:

\snic {
\sigma(\lambda) = \sigma(y)^p - \sigma(y) = (y+1)^p - (y+1)
= y^p  - y = \lambda.
}

%\sni
Le \polcar de $y$ est $(Y-y_0) (Y-y_1)\cdots (Y-y_{p-1})$ et
ce \pol est égal à $f(Y) = Y^p - Y - \lambda$ (car les
$y_i$ sont racines de~$f$ et $y_i - y_j$ est \iv pour
$i \ne j$).

\emph {5.}
Soit~$\gk$ un anneau avec $p=_\gk0$. Fixons $\lambda \in \gk$ et posons
$\gA = \aqo{\gk[Y]}{f} = \gk[y]$, où~$f(Y) = Y^p-Y-\lambda$. Alors, $y+1$ est
racine de~$f$, et l'on peur définir $\sigma \in \Aut(\gA/\gk)$ par
$\sigma(y) = y+1$. L'\elt $\sigma$ est d'ordre $p$ et \llec vérifiera que
$(\gk,\gA, \gen {\sigma})$ est une \aGz.

%%%%%%%%%%%%%%%%%%%%%%%%%%%%%%%%%%%%%%%%%%%%%%%%%%%%%%%%%%%%%%%%%%%%%%%%%%%

\prob{exoHeitmannGaloisExemple} 
\emph {1.}
Considérons l'\idz~$\gen {x-\sigma(x), y-\sigma(y)} \eqdefi \gen {2x,2y}$.
Puisque $2$ est \ivz, c'est l'\idz~$\gen {x,y}$, et il contient $1$ car $x^2 + y^2=1$.
Ainsi,~$\gen {\sigma}$ est séparant.

\emph {2.}
Pour tout $f\in \gB$, on a $f = (xf)x + (yf)y$. Si~$f$ est impaire i.e.  si
$\sigma(f) = -f$, on~a~$xf$,~$yf \in \gA$, donc $f \in \gA x + \gA y$ et $E =
\sotq {f\in\gB}{\sigma(f) = -f}$.  
L'\egtz~$\gB = \gA \oplus E$ découle de
l'\egt $f = \big(f+\sigma(f)\big)/2 + \big(f-\sigma(f)\big)/2$ pour $f \in \gB$.\\  
\emph{Autre \demz.} On
sait qu'il existe $b_0 \in \gB$ de trace $1$ et que le noyau de la forme \lin
$\gB \to \gA$ définie par $b \mapsto \Tr(b_0b)$ est un \sul de~$\gA$
dans~$\gB$. Ici on peut prendre $b_0 = 1/2$, on retrouve $E$ comme \sulz.

\emph {3.}
Il s'agit de trouver  $y_1$, $y_2$, $y_3 \in \gB$ tels que $\sum_{i=1}^3
x_i\tau(y_i)=1$ pour $\tau = \Id$, $0$ sinon.  On remarque que:
$$\preskip.4em \postskip.4em
1\cdot 1 + x\cdot x + y\cdot y = 2, \quad\hbox{et}\quad 
1\cdot \sigma(1) + x\cdot \sigma(x) + y\cdot \sigma(y) = 0,
$$
%\sni
d'où une solution en prenant $y_i = x_i/2$. En posant
\smashbot{$X = \Cmatrix{.3em} {1 & x & y\cr {1\over 2} &-{x\over 2} &-{y\over 2}\cr}$}, 
on~a~\hbox{$X \tra {X} = \I_2$} et $\tra{X} X = P$ avec $P = 
\Cmatrix{.4em} {1 & 0 &0\cr 0 &x^2 &xy\cr 0 & xy &y^2 \cr}$. La matrice
$P$ est un \prr de rang $2$ dont l'image est isomorphe au \Amoz~$\gB$.\\
Note: on en déduit que $E$ est isomorphe à l'image du projecteur
\smash{$\Cmatrix{.3em} {x^2 &xy\cr  xy &y^2}$} et 
que~$\gB\te_\gA E$ est isomorphe à~$\gB$ comme \Bmoz.

\emph {4.}
Facile

\emph {5.}
L'isomorphie $E^n \simeq \Ae n$ prouve que $E$ est un module \prcz~1.
En appliquant $\Al{n}$, on obtient $E^{n\te} \simeq \gA$.\\
NB: pour plus de détails voir la section~\ref{sec ptf loc lib},
la \dem de la proposition~\ref{prop puissance ext},
 l'\egt (\iref{eqVik}) \paref{eqVik} et  l'\egt (\iref{eqfactPicStab2}) \paref{eqfactPicStab2}.

\emph {6.}
L'\egt $1 = x^2 + y^2$ implique $\fa^2 = \gen {x^2 y^2, x^3y, x^4} = x^2 \gen {y^2, xy, x^2} = x^2\gA
.$
\\
Et $\fa\gB = xy\gB + x^2\gB = x(y\gB + x\gB) = x\gB$. Dans~$\gB$,
$\fa = x(y\gA + x\gA) = xE$. Donc si~$x$ est \ndzz, $\fa \simeq_\gA E$
via la multiplication par $x$.

\emph {7a.} On a~$\gk[x]\simeq\kX$ et~$\gA=\gk[x^2,xy,y^2]$.
On regarde~$\gB$ comme un~$\gk[x]$-module libre de rang $2$, de base $(1,y)$,
et l'on note $\rN : \gB \to \gk[x]$ la norme. Pour $a$, $b \in \gk[x]$ on obtient:
$$\preskip-.4em \postskip.4em 
\rN(a+by) = (a+by)(a-by) = a^2 + (x^2-1)b^2
. 
$$
Comme $\rN(x)=x^2$, $x$ est \ndz (lemme~\ref{lemIRAdu} point \emph{2}).
Par ailleurs, $a + by\in\Bti$  \ssi $a^2 + (x^2-1)b^2 \in \gk\eti$.  
Supposons $b$ de degré formel $m\geq 0$ et~$a$ de degré formel $n\geq 0$.
 Alors, $(x^2-1)b^2 = \beta^2 x^{2m+2} + \dots$ et $a^2 = \alpha^2 x^{2n} + \dots$ 
Puisque $a^2 +(x^2-1)b^2\in\gk\eti$, on obtient:
\begin{itemize}
%%
%\item  si $b=0$, $a^2\in\gk[x]\eti$ 
%
\item   si $n>m+1$, $\alpha^2=0$ donc $\alpha=0$ et $a$ 
peut être réécrit en degré formel $<n$,
\item si $n<m+1$, $\beta^2=0$ donc $\beta=0$ et  
\begin{itemize}
\item si $m=0$, $b=0$ et $a=\alpha\in\gk\eti$ ou,
\item  si $m>0$, $b$ peut être réécrit en degré formel $<m$,
\end{itemize}
\item si $n=m+1$ (ce qui implique $n>0$),  
$\alpha^2 + \beta^2= 0$ donc  $\alpha=\beta=0$ et $a$ 
peut être réécrit en degré formel $<n$.
\end{itemize}
  On conclut par \recu sur $m+n$ 
que si $a + by\in\Bti$, alors $b = 0$ et $a \in \gk\eti$.

 On notera que si $-1= i^2$ dans~$\gk$, alors $(x+iy)(x-iy) = 1$
et l'on obtient un \iv $x + iy$ qui n'est pas une constante.

\emph {7b.}
Montrons que $\fa$ n'est pas principal. Comme $\fa \simeq_\gA E$, il
s'en suivra que $E$ n'est pas un \Amo libre. Et $\gB$ n'est pas libre non plus,
car sinon $E$ serait \stl de rang $1$, donc libre.\\
Supposons $\fa = a\gA$
avec $a \in \gA$. En étendant à~$\gB$, on obtient $\fa\gB = a\gB$. Mais
on a vu que $\fa\gB = x\gB$, et $x$ étant \ndzz, $x = ua$ avec
$u \in \gB\eti = \gk\eti$. Ceci entraînerait que~$x \in \gA$,
ce qui n'est pas le cas car $\gk$ est non trivial.

\emph {8.}
On reprend la preuve de la question précédente pour montrer que
$\fa$ n'est pas principal, mais ici~$\gB\eti$ n'est plus constitué
uniquement des constantes, par exemple la fonction (continue) $(x,y) \mapsto
x^2 + 1$ est \ivz. \`A l'endroit où $x = ua$ \hbox{et $u \in \gB\eti$}, on
raisonne comme suit. Puisque $u$ est un \elt \iv de $\gB$, sa valeur absolue est minorée par un \elt $>0$, et $u$ est de signe strict constant. 
Comme $x$ est impaire et $a$ paire,  $a$ et $x$ sont identiquement nulles: contradiction.

%%%%%%%%%%%%%%%%%%%%%%%%%%%%%%%%%%%%%%%%%

% fin des solutions d'exos

%:  ---- Section*{references}-----------
\penalty-2500	
\Biblio

Une étude \cov des \algs associatives (non \ncrt commutatives) \stfes sur un \cdi se trouve dans \cite[Richman]{ri82}
et dans \cite[Chapitre~IX]{MRR}.

La proposition~\ref{propIdemMini} se trouve dans \cite{MRR}
qui introduit la terminologie de \emph{corps \splz ment factoriel}.
Voir aussi \cite[Richman]{ri81}.

Le lemme~\ref{lemSqfDec} de \fcn sans carrés sur un \cdi parfait admet une \gnn
subtile sous forme d'un \gui{\algo de \fcn \splz} sur un \cdi arbitraire:
voir \cite[th IV.6.3, p. 162]{MRR} et~\cite[Lecerf]{Lecerf-Factsep}.

Les notions  d'\aG et d'\alg \spb ont été introduites par Auslander \& Goldman dans \cite[1960]{AG}. L'essentiel de la théorie des \aGs se trouve dans
l'article \cite[1968]{CHR} de Chase, Harrison \& Rosenberg. 
Un livre qui expose cette théorie est \cite{DI}. 
Presque tous les arguments dans \cite{CHR} sont déjà de nature \elr et \covz.

Le résultat donné dans l'exercice~\ref{lemPaquesFerrero} est d\^u à Ferrero et Paques dans \cite{FP}.

Le problème  \ref {exoBuildingFrobAlgebra} s'inspire du chapitre 21
(Duality, Canonical Modules, and Gorenstein Rings) de \cite{Eis} et en 
particulier des exercices 21.6 et 21.7.

\newpage \thispagestyle{CMcadreseul}
\incrementeexosetprob

%:        %%%%%%%%%%%%%%%%%%%%%%%%%%%%%%%%%%%%
%:        %%%%%%%%%%%%%%%%%%%%%%%%%%%%%%%%%%%%
%---- Chapitre  {Théorie de Galois}------------
\chapter[La méthode dynamique]{La méthode dynamique 
%Corps de racines, théorie de Galois, \nstz
}
\label{ChapGalois}
%--------------------
\vspace{-1.2cm}
{\LARGE \bf
\hspace*{1cm} \nst\\[.3cm] 
\hspace*{1cm} Corps de racines\\[.3cm] 
\hspace*{1cm} Théorie de Galois}

\vspace{1.2cm}
\minitoc

\newcommand{\mpN}{mise en position de Noether }

%: Intro
\Intro
\pagestyle{CMExercicesheadings}

La première section de ce  chapitre donne des versions \covs \gnles du \nst 
pour un \syp sur un \cdi (on pourra comparer les \thos \rref{thNstNoe}, \rref{thNstClassCof1} et \rref {thNstClassCof2},  aux \thos \rref{thNstfaibleClass}
et \rref{thNstClass}).
Nous avons \egmt indiqué un \tho de mise en position de \iNoe simultanée
(\thref{thNoetSimult}).

Il s'agit là d'un exemple significatif d'une reformulation d'un résultat de \clama \emph{dans un cadre plus \gnlz}: les \clama admettent que tout corps possède une clôture \agqz. Cela leur permet de ne pas se poser le \pb de la signification exacte du \nst de Hilbert lorsque l'on n'a pas à sa disposition une telle clôture \agqz.\ihi
Mais la question se pose vraiment et nous apportons une réponse tout à fait raisonnable: la clôture \agq n'est pas vraiment \ncrz, plutôt que chercher les zéros d'un \syp dans une clôture \agqz, on peut les chercher dans des \algs finies sur le corps donné au départ.

\smallskip  
Nous nous attaquons ensuite à un autre \pbz: celui d'interpréter \cot le discours classique sur la clôture \agq d'un corps. 
Le \pb pourrait sembler être surtout celui de l'utilisation du lemme de Zorn \ncr à la construction
de l'objet global. En fait, un \pb plus délicat se pose bien avant, au moment de la construction du corps de racines d'un \pol individuel. 

Le \tho de \clama disant que
tout \pol \spl de $\KT$ possède un \cdr \stf
sur $\gK$ (auquel cas la théorie de Galois s'applique),
n'est valable d'un point de vue \cof que sous des
hypothèses concernant la possibilité de factoriser les \pols \spls
(cf.
\cite{MRR} et dans cet ouvrage le \thref{thResolUniv} d'une part et le corolaire \ref{propIdemMini} d'autre part).
Notre but ici est de donner une théorie de Galois \cov
pour un \pol \spl arbitraire en l'absence de telles
hypothèses.

La contrepartie est que l'on ne doit pas considérer le \cdr d'un
\pol comme un objet usuel \gui{statique}, mais comme un objet
\gui{dyna\-mi\-que}. Ce phénomène est inévitable, car il faut gérer
l'ambiguïté qui résulte de l'impossibilité de connaître le
groupe de Galois d'un \pol par une méthode infaillible.
Par ailleurs, le dépaysement produit par cette mise en perspective
dynamique n'est qu'un exemple de la méthode \gnle dite
d'\evn paresseuse: \emph{rien ne sert de trop se fatiguer pour
connaître toute la vérité quand une vérité partielle est
suffisante pour les enjeux du calcul en cours.}

\smallskip 
Dans la section \ref{subsecDyna}, nous donnons une approche heuristique
de la méthode dynamique, qui constitue une pierre angulaire des
nouvelles méthodes en \alg \covz.

La section \ref{secBoole} consacrée aux \agBs est une courte introduction aux \pbs qui vont devoir
être gérés dans le cadre d'une \adu sur un \cdi lorsqu'elle n'est
pas connexe.

La section \ref{secadu} continue la théorie de l'\adu
déjà commencée en section \ref{sec0adu}. Sans
supposer le \pol \spl l'\adu a de nombreuses \prts intéressantes
qui sont conservées quand on passe à un \gui{quotient de Galois}.
En faisant le résumé de ces \prts nous avons été amenés à introduire
la notion d'\emph{\apGz}.

La section \ref{subsecCDR} donne une approche \cov et
 dynamique du \cdr d'un \pol sur un \cdiz,
sans hypothèse de séparabilité pour le \polz.

La théorie  de Galois dynamique d'un  \pol \spl sur un \cdi est
développée dans la section~\ref{secThGB}.

\penalty-2500
\medskip
Le chapitre présent peut être lu
\imdt après les sections \ref{secGaloisElr} et~\ref{sec2GaloisElr} 
sans passer par les chapitres  \ref{chap mpf} et
\ref{chap ptf0}
si l'on se limite pour l'\adu au cas des \cdis 
(ce qui simplifierait d'ailleurs
certaines \demsz).
Il nous a paru cependant naturel de développer 
les questions relatives à l'\adu dans un cadre plus \gnlz, 
 ce qui nécessite la notion de \mrc sur un anneau commutatif arbitraire.

\pagebreak
\pagestyle{CMheadings}

%--- SECTION{Le \nst sans clôture \agqz   secNstSCA
\section{Le \nst sans clôture \agqz}
\label{secNstSCA}
%-----------------------------------------

Il nous a semblé logique, dans ce chapitre consacré à la question \gui{comment récupérer \cot les résultats de \clama qui se basent sur l'existence
d'une clôture \agqz, même lorsque celle-ci fait dé\-faut?},
de reprendre le \nst et la mise en position de \Noe (\thref{thNstfaibleClass}) dans ce nouveau cadre.

%:  subsec Le cas d'un corps de base infini
\subsec{Le cas d'un corps de base infini}

Nous affirmons que le \tho \ref{thNstfaibleClass} peut être recopié quasiment mot à mot, simplement en supprimant la référence à un \cac qui contienne $\gK$.

On ne voit plus \ncrt les zéros du \syp considéré 
dans des extensions finies du \cdi $\gK$,
mais on construit des \Klgs non nulles \stfes (i.e., qui sont des \Kevs de dimension finie) et qui 
rendent compte de ces zéros:  en \clama les zéros se trouvent dans les corps quotients de ces \Klgsz,
de tels corps quotients existent facilement en application du principe du tiers exclu puisqu'il suffit de considérer un \id strict qui soit de dimension maximale en tant que \Kevz. 

%\newpage
%:     Theorem{thNstfaibleClassSCA}
\begin{theorem}\label{thNstfaibleClassSCA}\emph{(\nst faible et mise en position de \iNoez, 2)}\\
Soit $\gK$  un \cdi infini et $(\lfs)$ un \syp dans l'\alg $\KuX=\KXn$ ($n\geq 1$).\\
Notons $\ff=\gen{\lfs}_\KuX$ \hbox{et $\gA=\KuX\sur\ff$}. \\
$\bullet$ \emph{(\nst faible)} 
\begin{itemize}
\item Ou bien $\gA=\so 0$, \cad $1\in \gen{\lfs}$.
\item Ou bien il existe un quotient non nul de $\gA$ qui est une \Klg \stfez. 
\end{itemize}
$\bullet$ \emph{(Position de \Noez)} Plus \prmtz, on a un entier
$r\in\lrb{-1..n}$ bien défini avec les \prts suivantes.
\begin{enumerate}

\item 
Ou bien $r=-1$ et $\gA=\so0$. 
Dans ce cas, le \sys  $(\lfs)$ n'admet de zéro dans aucune \Klg non triviale.

\item 
Ou bien $r=0$, et $\gA$ est une \Klg \stfe non nulle 
 (en particulier, l'\homo naturel $\gK\to\gA$ est injectif). 

\item 
 Ou bien $r\geq 1$, 
 et il existe un  \cdv $\gK$-\lin (les nouvelles variables sont notées $\Yn$)
 qui satisfait les \prts suivantes.
\begin{itemize}
\item [$\bullet$] On a $\ff\,\cap\,\gK[\Yr]=\so0$. Autrement dit, l'anneau  $\gK[\Yr]$
s'iden\-tifie à un sous-anneau du quotient $\KuX\sur\ff$.
\item [$\bullet$] Chaque $Y_j$ pour $j\in\lrb{r+1..n}$ est entier sur $\gK[\Yr]$ modulo $\ff$ et l'anneau $\gA$ est un $\gK[\Yr]$-\mpfz.
\item [$\bullet$] Il existe un entier $N$
 tel que pour chaque $(\alpha_1,\ldots,\alpha_r)\in\gK^r$, l'\alg quotient 
$\aqo\gA{Y_1-\alpha_1,\ldots,Y_r-\alpha_r}$ est un \Kev non nul
de dimension finie $\leq N$.
\item [$\bullet$] On a des \itfs $\ff_j\subseteq\gK[Y_1,\ldots,Y_j]$
$(j\in\lrb{r..n})$
avec les inclusions et \egts suivantes.
\[ 
\begin{array}{ll} 
 \gen{0}=\ff_r\subseteq\ff_{r+1} \subseteq \ldots\subseteq\ff_{n-1}\subseteq\ff_n=\ff   \\[1mm] 
\ff_j\subseteq\ff_\ell\,\cap\,\gK[Y_1,\ldots,Y_j] &(j<\ell,\ j,\ell \in\lrb{r..n})   \\[1mm] 
\rD(\ff_j) = \rD\left(\ff_\ell\,\cap\,\gK[Y_1,\ldots,Y_j]\right) &(j<\ell,\ j,\ell  \in\lrb{r..n})   
 \end{array}
\]

\end{itemize}
\end{enumerate}
\end{theorem}
%--------- fin theorem ---------------------------------------------- 
%
\begin{proof}
On raisonne essentiellement comme dans la \dem du \thref{thNstfaibleClass}.
Pour simplifier nous gardons les mêmes noms de variables à chaque étape de la construction.
On pose $\ff_n=\ff$.
\begin{itemize}
\item Ou bien $\ff=0$, et $r=n$ dans le point \emph{3.} 
\item Ou bien il y a un \pol non nul parmi les $f_i$, on fait un \cdv \lin qui le rend \mon en la  dernière variable,  et l'on calcule l'\id résultant $\fRes_{X_n}(\ff_n)=\ff_{n-1}\subseteq\gK[X_1,\ldots,X_{n-1}]\cap \ff_n$.
%:2018 ajout rappel d'où vient l'idéal résultant
(cf. \thref{thElimAff}).
Puisque $\ff_n\cap \gK[X_1,\ldots,X_{n-1}]$ et $\ff_{n-1}$ ont même nilradical, ils sont simultanément nuls.
\item  Si  $\ff_{n-1}=0$,  le point \emph{3}
ou \emph{2} est vérifié avec $r=n-1$.  
\item  Sinon, on itère le processus. 
\item  Si le processus s'arrête avec
$\ff_r=0$, $r\geq0$,   le point \emph{3} ou \emph{2}
est vérifié avec cette valeur de~$r$.  
\item   Sinon, $\ff_0=\gen{1}$ et le calcul a permis de construire $1$ comme \elt de~$\ff$.
\end{itemize}
Il nous reste à vérifier deux choses. 
\\
Tout d'abord, que $\gA$
est un $\gK[\Yr]$-\mpfz. Il est clair que c'est un module \tfz, le fait qu'il est 
\pf est donc donné par le \thref{propAlgFinPresfin}.  
\\
Ensuite, que lorsque l'on  spécialise les $Y_i$ ($i\in\lrbr$) en des $\alpha_i\in\gK$, \hbox{le \Kev} obtenu
est \pf (donc de dimension finie) et non nul. 
Le \thref{factSDIRKlg} sur les changements d'anneau de base nous donne le fait que,
après spécialisation, l'\alg reste un \mpfz, 
donc que le \Kev obtenu est bien de dimension finie.
Il faut voir qu'il est non nul. %\\ 
Or on constate que, en supposant les \cdvs déjà
faits au départ, tous les calculs faits dans $\gK[\Yn]$ se spécialisent,
\cad restent inchangés si l'on remplace les \idtrs $Y_1$, \dots, $Y_r$
par les scalaires $\alpha_1$, \dots, $\alpha_r$. Et la conclusion $\ff\;\cap\;\gK[\Yr]=\so0$ est remplacée par le même résultat spécialisé en les $\alpha_i$, 
\cad \prmt ce que nous voulons. 
\\
On peut obtenir la même conclusion sous la forme plus savante que voici. 
Cette spécialisation est un changement d'anneau de base $\gK[\Yr]\to\gK$.\perso{Si on ne le donnait que sous la forme savante, l'argument pour la spécialisation non nulle pourrait sembler un peu subtil.
Sous la forme concrète \gui{il est clair que tous les calculs se spécialisent},
il est en fait très simple. 
}  
%:2018  référence au tho  thElimAff  plutôt qu'au lemme LemElimAff
On applique le point 
\emph{4} du \tho d'\eli \agq \ref{thElimAff} avec
%\emph{1c} du lemme d'\eli \gnl \ref{LemElimAff} avec
$$\preskip.4em \postskip.4em 
\gk=\gK[\Yr], \;\gC=\gA\;\; \hbox{et}\; \;\gk'=\gK. 
$$
L'\id d'\eli
et l'\id résultant dans $\gk$ sont nuls, donc après \eds l'\id résultant reste nul dans~$\gK$. 
\\
Donc, la même chose vaut pour l'\id d'\eliz, et l'\homo naturel $\gK\to\aqo\gA{Y_1-\alpha_1,\ldots,Y_r-\alpha_r}$ est injectif.

 Expliquons pour terminer pourquoi l'entier $r$ est bien défini. Tout d'abord
le cas $r=-1$ est le seul cas où $\gA=\so0$, ensuite pour $r\geq 0$, un calcul montre que $r$ est le nombre maximum d'\elts \agqt indépendants sur~$\gK$ dans~$\gA$. 
\end{proof}

\rems
1)
On a utilisé des \ids résultants $\fRes(\fb)$ (\thref{thElimAff})
à la place des \ids $\fR(g_1,\ldots,g_s)$ 
(avec $g_1$ \mon et $\gen{g_1,\ldots,g_s}=\fb$),
introduits au lemme~\ref{lemElimPlusieurs}. Mais le lemme~\ref{lemElimParametre} 
montre que ces derniers feraient aussi bien l'affaire.\\
2) 
Pour n'importe quel \homo $\gK[\Yr]\to\gB$, lorsque $\gB$ est une~\Klg réduite, 
le dernier argument dans la \dem du \tho fonctionne, de sorte que l'on sait que 
$\gB\subseteq\gB\otimes_{\gK[\Yr]}\gA$.    \\
3)
Le dernier item du point \emph{3} rappelle le fonctionnement de la preuve par \recuz, laquelle construit les \itfs $\ff_j$ pour aboutir à la mise en position de \iNoez. Cela donne aussi une certaine description des \gui{zéros} du \syp
(plus délicate que dans le cas où l'on a un \cac $\gL$ qui contient $\gK$,
et où l'on décrit les zéros à \coos dans $\gL$, 
comme dans le \thref{thNstfaibleClass}).  
\eoe

\medskip Il nous reste à lever la restriction introduite par la
considération d'un \cdi $\gK$ infini.
Pour ceci nous avons besoin d'un lemme de \cdv un peu plus \gnlz, qui utilise
une astuce de Nagata.

%:  subsec Changements de variables
\subsec{Changements de variables}

%:     Definition{defiCDV}
\begin{definition}\label{defiCDV}
On appelle \ix{changement de variables} dans l'anneau de \pols   $\kuX=\kXn$  un
\auto $\theta$ de cette \klgz. Si les~$\theta(X_i)$  sont notés $Y_i$,
les $Y_i$  sont appelés \emph{les nouvelles variables}. Chaque~$Y_i$ est un \pol
en les $X_j$, et chaque $X_i$ est un \pol en~les~$Y_j$. 
\end{definition}

Le plus fréquemment utilisés sont
les \gui{changements de variables \linsz}, dans lesquels on inclut, malgré leur nom, les translations et toutes les transformations affines.

\medskip 
\comm Un \cdv non \linz, comme par exemple 
$$
{(X,Y)\mapsto(X+Y^2,Y),}
$$
ne respecte pas la \gmt au sens intuitif de la chose. Par exemple une droite est transformée en une parabole: la \gmt \agq du plan affine n'est pas une extension de la \gmt affine,
elle est directement en contradiction avec elle!
C'est seulement dans le cadre des espaces projectifs que l'on retrouve ses petits: les \autos du plan projectif, du point de vue de la \gmt \agqz,
sont \ncrt \linsz, et la notion de \gui{droite} reprend ses droits.
\eoe

%-% ENTRE NOUS
\entrenous{ Je ne pense pas que cela vaut le coup d'introduire la notation
et le commentaire qui suivent.

%:     Notation{notaEspAff}
%\begin{notation}\label{notaEspAff}
{\bf Notation.} {\rm  On note $\An(\gk)$ l'espace affine $\gk^n$ en \gaqz. 
  }
%\end{notation}
%--------- fin notation ---------------------------------------------- 

\rm
 En \gaq un \ixc{morphisme}{entre espaces affines} 
 $\psi:\Am(\gk)\to\An(\gk)$ est donné par
 n'importe quelle application  de $\gk^m$ vers $\gk^n$
 définie par $n$ \pols dans $\kXm$, ou ce qui revient
 au même, par un $\gk$-morphisme
 
\snic{\wi \psi:\gk[\Yn]\to\kXm.}

%\sni
 On dit que $\wi \psi$ est le \ix{co-morphisme} de $\psi$.
 \\
   Si $\rho:\gk\to\gA$ est une \klg un tel morphisme $\psi$ définit
 aussi un morphisme de $\Am(\gA)$ vers $\An(\gA)$: son co-morphisme est obtenu
 par \edsz.
\sf}
%-% Fin ENTRENOUS

%:  subsubsec pols pseudo unitaires
\bni {\bf\Pols pseudo unitaires}~
\rdb\label{polpseudunit}

\smallskip Soit $\gk$ un anneau connexe. Un \pol  dans  $\kT$ est dit
\ixc{pseudo unitaire}{polynome@\pol ---} (en la variable $T$) s'il s'écrit
$\som_{i=0}^pa_kT^{k}$ avec $a_p$ \ivz.

En \gnlz, sans supposer $\gk$ connexe, un \pol dans $\kT$ est dit
\emph{pseudo unitaire} (en la variable $T$)  s'il existe un \sfio $(e_0,\ldots
,e_r)$ tel que,
pour chaque $j$,
en passant à $\gk[1/e_j]=\gk_j$, le \pol s'écrit
$\som_{k=0}^{j}a_{k,j}T^{k}$ avec $a_{j,j}$ \iv dans $\gk_j$.

Un \pol  dans  $\kXn=\kuX$ est dit
\emph{pseudo unitaire en la variable} $X_n$ s'il est pseudo unitaire
comme \elt de $\gk[X_1,\ldots ,X_{n-1}][X_n]$.\label{polpseudounitaire}%
\index{polynome@\pol!pseudo unitaire}

 NB: voir aussi la notion de \pol \lot \mon dans l'exercice~\ref{exoPolLocUnitaire}.

 Rappelons qu'un \pol de $\kXn$ est dit \emph{primitif} lorsque ses \coes
engendrent l'\id $\gen{1}$.
Rappelons aussi que si $\gk$ est réduit, on a l'\egt $\kXn\eti=\gk\eti$ (lemme~\ref{lemGaussJoyal}).

%:    factZedPrimFin
\begin{fact}\label{factZedPrimFin}
Soient $\gK$ un anneau \zed réduit et $P\in \gK[T]$. \Propeq
\begin{enumerate}
  \item [--] Le \pol $P$ est \ndzz.
  \item [--] Le \pol $P$ est primitif.
  \item [--] Le \pol $P$  est pseudo unitaire.
  \item [--] L'\alg quotient $\aqo{\gK[T]}{P}$  est finie sur $\gK$.
\end{enumerate}
\end{fact}
\begin{proof}
Les \eqvcs sont claires dans le cas des \cdisz. Pour obtenir le résultat \gnl on peut appliquer la machinerie \lgbe \elr des anneaux \zeds réduits 
\paref{MethodeZedRed}.\imlgz
\end{proof}
%

%  subsubsec un lemme
\penalty-2500
\subsubsection*{Un lemme simple et efficace}

%:    Lemma{lemNoether}---------------
\begin{lemma}
\label{lemNoether} \emph{(Lemme de \cdvs à la Nagata)}\\
Soit $\gK$  un anneau \zed réduit et  $g\in\KuX=\KXn$ un \elt \ndzz.
\begin{enumerate}
\item Il existe
un \cdv tel que, en appelant $Y_1$, $\ldots $, $Y_{n}$ les nouvelles
variables, le \pol $g$ devient pseudo unitaire en $Y_n$.
En conséquence
la \Klg $\aqo{\KuX}{g}$ est finie sur $\gK[Y_1,\ldots ,Y_{n-1}]$.
\item Lorsque $\gK$ est un \cdi infini, on peut prendre un changement \lin de variables.
\item Le résultat s'applique aussi pour une famille finie de \pols   \ndzs de $\KuX$ (ils peuvent être rendus simultanément pseudo \mons par un même
\cdvz).
\end{enumerate}
\end{lemma}
%--- end-lemma-----------------------------------------
%-----------------begin proof------------------
\begin{proof}
Pour le cas d'un \cdi infini voir le lemme \ref{lemCDV}.
\\
 Dans le cas \gnl
on peut supposer que $\gK$ est un \cdi et l'on fait un \cdv \gui{à la Nagata}. Par exemple avec trois variables,  si le
\pol $g$ est de degré $< d$ en chacune des variables  $X$, $Y$,~$Z$, on
fait le \cdv $X\mapsto X,\,Y\mapsto Y+X^d,\,Z\mapsto
Z+X^{d^2}.$
Alors, vu comme \elt de $\gK[Y,Z][X]$, $g$ est devenu pseudo unitaire en~$X$.

 Le point \emph{3} est laissé \alecz.
\end{proof}
%-----------------end proof------------------

%:  subsec Le cas général
\subsec{Le cas \gnl}

En raisonnant comme pour le  \tho \ref{thNstfaibleClassSCA} et en utilisant les
\cdvs du lemme précédent on obtient la forme \gnle du \nst 
faible et de la mise en position de \Noe en \comaz.

%:     Theorem{thNstNoe}
\begin{theorem}\label{thNstNoe}\emph{(\nst faible et mise en position de \iNoez, 3)}\\
Avec les mêmes hypothèses que dans le \thref{thNstfaibleClassSCA} mais 
%:HHH 11 legere modif, pas importante
%sans supposer que le \cdi $\gK$ est infini,
en supposant seulement que le \cdi $\gK$ est non trivial, 
on a les mêmes conclusions, à ceci près que le \cdv n'est pas \ncrt \linz.
\end{theorem}
%--------- fin theorem ---------------------------------------------- 

%:     Definition{defiDimsyp}
\begin{definition}\label{defiDimsyp} 
On considère le cas $1\notin\gen{\lfs}$ du \tho précédent. 
\begin{enumerate}
\item On dit que le \cdv (qui éventuellement ne change rien du
tout) a mis l'\id $\ff$ en \emph{position de \Noez}.
\item L'entier $r$ qui intervient dans la mise en position de \Noe 
est appelé la \ixc{dimension}{d'une
\vrt affine} \emph{du \sypz}, ou de la \vrt définie par le \sypz, ou de l'\alg quotient $\gA$.
Par convention l'\alg nulle est dite de dimension~$-1$.%
\index{dimension!d'un \syp sur un \cdiz}\index{dimension!d'une \apf sur un \cdiz}
\end{enumerate}
\end{definition}
%: remarques
\rems 1) Il est clair d'après le \tho %de mise en position de \Noe 
que $r=0$ \ssi l'\alg quotient est finie non nulle, ce qui implique
(lemme \ref{lemZrZr1}) que c'est un anneau \zed non trivial.  
\\
%:HHH 11 reciproque legerment amelioree, sans importance reelle
Inversement, si $\gA$ est \zed et $\gK$ non trivial,  le lemme~\ref{lemZrZr2}
montre que l'anneau $\gK[\Yr]$ est \zedz, ce qui implique que $r\leq 0$
(\hbox{si $r>0$,} alors une \egt $Y_r^m\big(1+Y_rQ(\Yr)\big)=0$ implique que $\gK$ est trivial).
Il n'y a donc pas de conflit avec la notion d'anneau \zedz.
Notons cependant que  l'\alg nulle  est encore un anneau \zedz.
 
 2) Le lien avec la \ddk sera fait dans le \thref{thDKAG}.
 
 3) Une version \gui{non \noeez} du \tho précédent pour 
 un anneau \zedr $\gK$ est donnée en exercice~\ref{exothNst1-zed}.
\eoe

%:   Theorem{thNoetSimult}------
\begin{theorem}
\label{thNoetSimult} \emph{(Mise en position de \iNoe simultanée)}
\\
Soient $\ff_1$, \ldots, $\ff_k$ des \itfs de $\KuX=\KXn$. %Alors:
\begin{enumerate}\itemsep0pt
\item  Il existe des entiers $r_1, $\ldots, $r_k\in\lrb{-1..n}$ et  un \cdv tels que, en appelant $Y_1$, \ldots, $Y_{n}$ les nouvelles variables, on ait
pour chaque $j\in\lrbk$ la situation suivante.
\\
Si $r_j=-1$, alors $\ff_j=\gen{1}$, sinon
%-----------------begin item------------------
\begin{enumerate}
\item 
$\gK[Y_1,\ldots , Y_{r_j}]\,\cap\, \ff_j=\so{0}$,
\item 
pour $\ell>r_j$, $Y_\ell$ est entier modulo 
$\ff_j$ sur $\gK[Y_1,\ldots , Y_{r_j}]$.
\end{enumerate}
%-----------------end item------------------
Lorsque $\gK$ est infini, on peut prendre un changement \lin de
variables.

\item Si $\gen{1}\neq\rD(\ff_1)\supset \rD(\ff_2)\supset \cdots \supset \rD(\ff_k)$ avec les dimensions $r_j$ strictement croissantes, on peut intercaler des radicaux d'\itfs de sorte que la suite 
des dimensions obtenue soit $0$, $1$, \ldots, $n$.  
%
%\item 
%
\end{enumerate}
\end{theorem}
%--- end-theorem-----------------------------------------
NB: dans le point \emph{1}, on dit que le \cdv (qui éventuellement ne change rien du
tout) a mis simultanément les \ids $\ff_1,\ldots,\ff_k$ en position de \Noez.

\begin{proof}  \emph{1.}
La même \dem que pour le \tho précédent fonctionne compte tenu du fait qu'un
\cdv peut rendre simultanément \mons en la dernière variable un nombre fini de \pols non~nuls.
\\
\emph{2.} Posons $\gA_i=\gK[X_1,\ldots,X_i]$.
Supposons par exemple que $\ff_1$ soit de dimension $2$ et $\ff_2$ de dimension $5$.  Nous devons intercaler des \ids de dimensions~$3$ et $4$. 
Nous supposons \spdg que les $\ff_i$ sont en position de \Noe par rapport à $\Xn$.
\\
On a par hypothèse $\gA_2\,\cap\,\ff_1=0$, avec des \polus 

\snic{h_3\in\gA_2[X_3]\,\cap\,\ff_1$, $h_4\in\gA_2[X_4]\,\cap\,\ff_1$, \ldots, $h_n\in\gA_2[X_n]\,\cap\,\ff_1.}

%\sni
On a alors les inclusions suivantes,

\snac{\fh_1=\ff_2+\gen{h_5,h_4} \supseteq \fh_2=\ff_2+\gen{h_5}
\supseteq \ff_2 \et \rD(\ff_1)\supseteq \rD(\fh_1)\supseteq \rD(\fh_2)\supseteq \rD(\ff_2),}

%\sni
avec $\fh_1$ de dimension 3 et $\fh_2$ de dimension 4, tous deux en position de \Noe
par rapport à $(\Xn)$. 
\end{proof}
%

%-% ENTRE NOUS
\entrenous{Le \tho précédent devrait être exploité dans le chapitre sur la \ddk pour donner sous forme \cov la \prt caténaire des anneaux de \pols
sur les \cdis (donc aussi pour les anneaux de fonctions intègres en \gaqz.) 
}
%-% Fin ENTRENOUS

%:  subsec Le \nst proprement dit
\subsec{Le \nst proprement dit}
\vspace{3pt}
Dans les \thos \ref{thNstfaibleClassSCA} (\cdi infini) et \ref{thNstNoe} 
(\cdi arbitraire) le \nst est sous forme faible, \cad
qu'est démontrée l'\eqvc entre d'une part,
\\
$\bullet$ le \syp %considéré 
n'a de zéro dans aucune \Klg finie non nulle,\\
 et d'autre part,
\\
$\bullet$   l'\alg quotient correspondante est nulle.

Le \nst \gnl dit à quelle condition un \pol s'annule aux zéros d'un \sypz. Ici, puisque nous n'avons pas de \cac
à notre disposition, nous considérerons les zéros dans les \Klgs
 finies et nous obtenons deux \nsts selon que nous considérons seulement
 les \Klgs réduites ou pas.
 
 Ces deux \thos généralisent d'un point de vue \cof (avec des~\gui{ou bien}
 explicites) le \nst classique énoncé sous la forme du \thref{thNstClass}.
% Notons que le deuxième énoncé, avec multiplicités, est 

%:     Theorem{thNstClassCof1}
\begin{theorem}\label{thNstClassCof1} \emph{(\nst classique, version \cov \gnlez)}\\
Soit $\gK$  un \cdi et %des \pols 
$f_1$, \dots, $f_s$, $g$  dans $\KXn$. 
Considérons l'\alg quotient $\gA=\aqo\KuX\lfs$. 
\begin{enumerate}
\item Ou bien il existe un quotient non nul $\gB$ de $\gA$ qui est une \Klg finie réduite avec $g\in\gB\eti$ (a fortiori $g\neq0$ dans $\gB$).
\item Ou bien $g$ est nilpotent dans $\gA$ (autrement dit, il existe un entier $N$ tel que $g^N\in\gen{\lfs}_\KuX$).
\end{enumerate}
\end{theorem}
%--------- fin theorem ---------------------------------------------- 
%
\begin{proof}
On utilise l'astuce de Rabinovitch. On introduit une \idtr supplémentaire
$T$ et l'on remarque que $g$ est nilpotent dans~$\gA$ 
\ssi l'\alg quotient~$\gA'$ pour le \syp
$(\lfs,1-gT)$ est nulle. On termine avec le \nst faible: si $\gA'\neq 0$,
on trouve un quotient non nul $\gB'$ de $\gA'$ qui est un \Kev de dimension
finie. Comme $g$ est \iv dans $\gA'$, il l'est aussi dans $\gB'$ et dans $\gB=\gB'\red$, et comme $\gB\neq 0$, $g\neq 0$ dans $\gB$.
\end{proof}
%

%:     Theorem{thNstClassCof2}
\begin{theorem}\label{thNstClassCof2} \emph{(\nst avec multiplicités)}\\
Soit $\gK$  un \cdi et %des \pols 
$f_1$, \dots, $f_s$, $g$  dans $\KXn$. 
Considérons l'\alg quotient $\gA=\aqo\KuX\lfs$. 
\begin{enumerate}
\item Ou bien il existe un quotient $\gB$ de $\gA$ qui est un \Kev de dimension finie  avec $g\neq0$ dans $\gB$.
\item Ou bien $g=0$ dans $\gA$ (autrement dit, $g\in\gen{\lfs}_\KuX$).
\end{enumerate}
\end{theorem}
%--------- fin theorem ---------------------------------------------- 
%
\begin{Proof}{\Demo utilisant les \bdgsz. }
Si dans la mise en position de \Noe on a $r=0$, le résultat est clair.
Le point délicat est lorsque~$r\geq 1$.
On suppose l'\id en position de \Noez. On considère un ordre d'\eli
pour les variables $(\Yr)$, puis la forme normale de $g$ par rapport à la \bdg de $\ff$.  Pour que \gui{tout reste en l'état} après
une spécialisation $Y_i\mapsto \alpha_i=\ov{Y_i}$ 
dans un anneau quotient  $\gL$ de $\gK[\Yr]$,
il suffit que les \coes de tête dans la \bdg de $\ff$ et dans la forme normale de $g$ (ces \coes sont des \elts de $\gK[\Yr]$) se spécialisent en des \elts \ivs de $\gL$. Si l'on  dispose  de suffisamment d'\elts distincts dans $\gK$ 
pour trouver des $\alpha_i$ convenables dans $\gK$  on peut prendre $\gL=\gK$, sinon
on considère le produit $h$ de tous les \coes de tête considérés précédemment, et l'on remplace $\gK[\Yr]$ par un quotient~$\gL$ non nul, \stf sur $\gK$, dans lequel $h$
est inversible (ceci est possible par le \thref{thNstClassCof1}, appliqué 
pour $h$ sans aucune équation $f_i$). La~solution à notre \pb est alors donnée par l'\alg

\snic{\gB=\gL\otimes_{\gK[\Yr]}\gA,}

%\sni
qui est un quotient de $\gA$  \stf sur $\gK$.
\end{Proof}
%

%:  subsec Module des \rdls
\subsec{Module des syzygies}
Une autre conséquence importante du lemme de \cdvz~\ref{lemNoether}  est le \tho suivant.

%:   Theorem{thpolcohfd}------
\begin{theorem}
\label{thpolcohfd}
Soit $\gK$ un anneau \zedr discret.
\begin{enumerate}
\item Toute \Klg \pf est un \cori et fortement discret.
\item En conséquence tout \mpf sur une telle \alg est \coh et fortement discret.
%
%\item
%
\end{enumerate}
\end{theorem}
%--- end-theorem---------------------------
%-----------------begin proof------------------
\begin{proof}
Nous montrons le premier point
pour $\KXn$ dans le cas où $\gK$ est un \cdiz. 
Le cas des anneaux \zeds s'en déduit par la technique habituelle (machinerie \lgbe \elr \num2). Ensuite le point \emph{2} est une conséquence du \thref{propCoh2}.\imlg 
\\
Nous faisons une preuve par \recu sur $n$,
le cas $n=0$ étant clair. Nous supposons $n\geq 1$ et nous notons $\gB=\KXn$.
Nous devons montrer qu'un \itf arbitraire $\ff=\gen{f_1,\ldots ,f_s}$ est \pf et
détachable.\\
Si $\ff=0$ c'est clair, dans le cas contraire on peut supposer en appliquant le lemme
\ref{lemNoether} que $f_s$ est unitaire en $X_n$ de degré $d$.
Si $s=1$, l'annulateur de~$f_1$ est nul, et donc aussi le module des \syzys pour
$(f_1)$. Et l'\id $\ff$ est détachable grâce à la division euclidienne par
rapport à~$X_n$. \\
Si $s\geq 2$, notons $\gA=\gK[X_1,\ldots ,X_{n-1}]$. L'anneau $\gA$ est \coh
\fdi par \hdrz. Notons $R_i$ la \syzy qui correspond à l'\egt $f_if_s-f_sf_i=0$ 
($i\in\lrb{1..s-1}$).  
Modulo les \syzysz~$R_i$ on peut réécrire les $X_n^kf_i=g_{k,i}$,
pour~$k\in\lrb{0..d-1}$ et~$i\in\lrb{1..s-1}$ comme des vecteurs dans le \Amo libre
$L\subseteq\gB$ engendré par $(1,X_n, \ldots ,X_n^{d-1})$. Modulo les \syzys
$R_i$ toute \syzy pour $(f_1,\ldots ,f_s)$ à \coes dans~$\gB$ se réécrit
comme une \syzy pour
\perso{Il semble que le contenu du \tho \ref{thpolcohfd} soit partie intégrante
de la preuve originelle du \nst par Hilbert, à moins que ce soit la preuve du
\tho "de la base". \`A vérifier et à signaler.
En tout cas, notre preuve du \nst n'utilise pas ce \thoz.}

\snic{V=(g_{0,1},\ldots ,g_{d-1,1},\ldots ,g_{0,s-1},\ldots ,g_{d-1,s-1})\in
L^{d(s-1)}}

à \coes dans $\gA$. Comme $L$ est un \Amo libre,
il est \coh \fdiz. On a en particulier un nombre fini
de $\gA$-\syzys pour $V$ qui les engendrent toutes.
Appelons les $S_1$, \ldots, $S_\ell$.
Chaque $\gA$-\syzy $S_j$  pour~$V$ peut être lue comme une
$\gB$-\syzy $S'_j$ pour  $(f_1,\ldots ,f_{s})$. Finalement, les \syzys $R_i$
et $S'_j$ engendrent le \Bmo des \syzys pour  $(f_1,\ldots ,f_{s})$.\\
Concernant le \crc \fdiz, on raisonne de la même manière. Pour tester
si un \elt de $\gB$ est dans $\ff$ on commence par le diviser par~$f_s$ par
rapport à $X_n$. On obtient alors un vecteur dans le \Amo $L$ dont il faut
tester s'il appartient au sous-module engendré par les~$g_{i,j}$.
%-% perso
\perso{On voudrait avoir aussi, pour pas trop cher si possible, un résultat
trivial en \clamaz, et que l'on  peut attraper \cot avec les \bdgsz:
un \itf de $\KXn$ possède une base finie ou dénombrable sur $\gK$,
et il admet un \Kev supplémentaire qui possède une base finie ou
dénombrable sur $\gK$.
Le même style de \dem que pour le \tho précédent devrait pouvoir fonctionner.
}
%-% Fin perso
\end{proof}
%-----------------end proof------------------

\goodbreak
%--- SUBsection{subsecLGIdeps}--------
\section{La méthode dynamique}
\label{subsecDyna}
%-------------------------------------

\begin{flushright}
{\em Je ne crois pas aux miracles.
}\\
Un mathématicien \cofz.
\end{flushright}

En \clama les preuves d'existence sont rarement expli\-cites.
Deux obstacles essentiels apparaissent chaque fois que l'on essaie de rendre une
telle preuve explicite.

Le premier obstacle est l'application du principe du tiers exclu. Par exemple, si
vous considérez la preuve que tout \pol univarié sur un corps $\gK$ admet une
\dcn en facteurs premiers, vous avez une sorte d'\algo dont
l'ingrédient essentiel est: si $P$ est irréductible c'est bon, si $P$ se
décompose en un produit de deux facteurs de degré $\geq 1$, c'est bon aussi,
par \hdrz. Malheureusement la disjonction qui sert à faire fonctionner la preuve
\gui{$P$ est irréductible ou $P$ se décompose en un produit de deux facteurs de
degré $\geq 1$} n'est pas en \gnl explicite. Autrement dit, même si un
corps est défini de manière \covz, on ne peut être certain que cette
disjonction puisse être explicitée par un \algoz. Nous nous trouvons ici en
présence d'un cas typique où le principe du tiers exclu \gui{pose problème},
car l'existence d'un facteur irréductible ne peut pas faire l'objet d'un
\algo \gnlz.

Le deuxième obstacle est l'application du lemme de Zorn, qui permet de
\gnr au cas non dénombrable les raisonnements par \recu usuels dans
le cas dénombrable.

Par exemple dans le Modern Algebra de van der Waerden le second écueil est
évité en se limitant aux structures \agqs dénombrables.

\penalty-2500
\ss Nous avons cependant deux
faits d'expérience désormais bien établis:
%-------begin item---
\begin{itemize}
\item  Les résultats concrets \emph{universels} démontrés par les méthodes
abstraites douteuses ci-dessus n'ont jamais été contredits. On a même très
souvent réussi à en fournir des  \prcos incontestables.
Cela signifierait que même si les méthodes abstraites sont quelque part
fautives ou contradictoires, elles n'ont jusqu'à présent été utilisées
qu'avec suffisamment de discernement.
\item  Les résultats concrets existentiels démontrés par les méthodes
abstraites douteuses n'ont pas non plus été infirmés. Bien au contraire, ils
ont souvent été confirmés par des \algos démontrés \cotz%
\footnote{Sur ce deuxième point, notre affirmation est moins nette.
Si nous revenons à
l'exemple de la \dcn d'un \pol en facteurs premiers, il est impossible de
réaliser le résultat de manière \algq sur certains corps.}.
\end{itemize}
%-----------------end item------------------

\ss Face à cette situation un peu paradoxale: les méthodes abstraites sont a
priori douteuses, mais elles ne nous trompent pas fondamentalement quand elles
donnent un résultat de nature concrète, il y a deux réactions possibles.

Ou bien l'on croit que les méthodes abstraites sont fondamentalement justes parce
qu'elles reflètent une \gui{réalité}, une sorte d'\gui{univers cantorien
idéal} dans lequel se trouve la vraie sémantique des \mathsz.
C'est la position du
réalisme platonicien, défendue par exemple par G\"odel.

Ou bien l'on pense que les méthodes abstraites sont vraiment sujettes à caution.
 Mais alors, à moins de croire que les
\maths relèvent de la magie ou du miracle, il faut expliquer pourquoi les \clama
se trompent si peu. Si l'on ne croit ni à Cantor, ni aux miracles, on est  conduit
à penser que les preuves abstraites de résultats concrets contiennent \ncrt des
\gui{ingrédients cachés} suffisants pour construire les preuves concrètes
correspondantes.

Cette possibilité de certifier \cot des résultats concrets obtenus par des
méthodes douteuses, si l'on arrive à la réaliser de manière assez
systématique, est dans le droit fil du programme de Hilbert.\ihi

La méthode dynamique en \alg constructive est une méthode \gnle
de décryptage des preuves abstraites des \clama lorsqu'elles utilisent
des objets \gui{idéaux} dont l'existence repose sur des principes non
\cofsz: le tiers exclu et l'axiome du choix.
L'ambition de cette nouvelle méthode
est de \gui{donner une sémantique \cov pour les \clama usuellement
pratiquées}.

Nous remplaçons les objets abstraits des \clama par des spécifications
incomplètes mais concrètes
de ces objets. C'est la contrepartie \cov des objets abstraits.
Par exemple un \emph{\ipp fini} (notion qui sera introduite en section
\ref{subsecMoco}) est donné par un nombre fini d'\elts dans l'\id et
un nombre fini d'\elts dans son complémentaire.
Cela constitue une spécification
incomplète mais concrète d'un \idepz.

Plus \prmtz, la méthode dynamique vise à donner une interprétation
systématique de preuves classiques qui utilisent des objets abstraits
en les relisant comme des  \prcos
au sujet de contreparties \covs de ces objets abstraits.

Cela se situe dans le même esprit que certaines techniques développées en  calcul formel. Nous pensons ici à  l'\gui{\evn
paresseuse}, ou l'\gui{\evn dynamique},
\cad  l'\evn paresseuse gérée de
manière arborescente, comme dans le système D5 \cite{D5} qui réalise de
manière très innocente ce tour de force:
calculer de manière sûre dans la
clôture algébrique d'un corps arbitraire, alors même
que l'on sait que cet objet
(la clôture algébrique)
ne peut pas être construit en toute \gntz.

Dans le chapitre présent
une spécification incomplète du \cdr d'un \pol
\spl sur un corps $\gK$ sera donnée par une \Klg $\gA$ et un groupe
fini d'\autos $G$ de cette \algz.
Dans $\gA$ le \pol se décompose en facteurs \lins
de sorte qu'un \cdr est un quotient de $\gA$,
et $G$ est une approximation du groupe de Galois en un sens convenable
(en particulier, il contient une copie du groupe de Galois).
Nous expliquerons comment calculer avec une telle approximation
sans jamais se tromper: quand une bizarrerie se manifeste, on sait comment
faire pour améliorer l'approximation en cours et faire disparaître
la bizarrerie.

%:  subsection{Galois d'un pol en clama}
\subsect{Corps de racines et théorie de Galois  en \clamaz}{Théorie de Galois classique}\label{GalPolCla}

Dans ce paragraphe nous indiquons un exposé
possible du \cdr d'un \pol arbitraire et de la théorie de Galois d'un \pol \spl en \clamaz.
Ceci permet de comprendre les \gui{détours} que nous serons obligés de faire pour
avoir une théorie pleinement \covz.

\smallskip
Si $f$ est un \poluz, on travaille avec l'\adu de $f$, $\gA=\Adu_{\gK,f}$ dans laquelle
$f(T)=\prod_i(T-x_i)$, avec $\Sn$ comme groupe d'\autos 
(voir la section~\ref{sec0adu}).

Cette \alg étant un \Kev de dimension finie, tous les \ids sont eux-mêmes des
\Kevs de dimension finie et
l'on a le droit de considérer un
\id strict $\fm$ de dimension maximum comme \Kev (tout ceci en application du principe du tiers exclu).
Cet \id est automatiquement un \idemaz.
L'\alg quotient $\gL=\gA/\fm$ est alors un \cdr pour $f$.
Le groupe $G=\St(\fm)$ opère sur $\gL$ et le corps fixe de $G$, $\gL^G=\gK_1$,
possède les deux \prts suivantes:
\begin{itemize}
\item $\gL/\gK_1$ est une extension galoisienne avec $\Gal(\gL/\gK_1)\simeq G$.
\item  $\gK_1/\gK$ est une extension obtenue par adjonctions successives de racines $p$-ièmes,
où $p=\car(\gK)$.
\end{itemize}

\smallskip
En outre, si $\gL'$ est un autre \cdr pour $f$ avec $f=\prod_i(T-\xi_i)$ 
dans~$\gL'[T]$, on a un unique \homo de \Klgs $\varphi:\gA\to\gL'$
vérifiant les \egts
$\varphi(x_i)=\xi_i$ pour $i\in\lrbn$. On peut alors montrer que~$\Ker\varphi$,
qui est un \idema de $\gA$, est \ncrt conjugué de~$\fm$ sous l'action de $\Sn$.
Ainsi le \cdr est unique, à \iso  près (cet \iso est non unique si
$G\neq\so{\Id}$).

Enfin, lorsque $f$ est \splz, la situation est simplifiée parce que l'\adu est étale, et $\gK_1=\gK$.

\smallskip
La démarche précédente est possible d'un point de vue \cof si le corps $\gK$
est \splz ment factoriel et si le \pol $f$ est \splz,  car alors,
puisque  l'\adu $\gA$ est étale, elle se décompose en un produit fini de
corps étales sur $\gK$ (corolaire~\ref{propIdemMini}).

Mais lorsque le corps n'est pas \splz ment factoriel, on est face à un obstacle
a priori rédhibitoire, et l'on ne peut pas espérer obtenir de
 manière systématique et \algq un \cdr qui soit \stf sur~$\gK$.

Si la \cara est finie et si le \pol n'est pas \splz,
on a besoin de propriétés de factorisation
plus fortes pour construire un \cdr (la question est délicate,
et très bien exposée dans \cite{MRR}).

%:  subsection{Contourner l'obstacle}
\subsect{Contourner l'obstacle de façon paresseuse}{Contourner l'obstacle}\label{GalPolComa}

Ce qui est \gnlt proposé en calcul formel c'est, par exemple dans le cas d'un \pol \splz,
 à tout le moins d'éviter de calculer une résolvante \uvle $R$ (comme dans le \thref{thResolUniv}) dont le degré,~$n!$,
rend rapidement les calculs impraticables.

Ici, nous nous situons dans le cadre le plus \gnl possible,
et nous évitons tout recours à la factorisation des \pols qui
peut s'avérer impossible, ou qui, lorsqu'elle est possible, risque de co\^uter trop cher.

L'idée est d'utiliser l'\adu $\gA$,
ou bien  un quotient de Galois $\aqo{\gA}{1-e}$, avec un \idm  galoisien $e$
(voir \paref{defQuoDeGal}) comme
substitut pour $\gL$. Ce \gui{\cdr dynamique} peut être géré
sans trop de \pbs parce que chaque fois qu'il se passe quelque chose
d'étrange, qui manifeste que le substitut à $\gL$ n'est pas
entièrement satisfaisant, on est
capable de \gui{réparer immédiatement la bizarrerie} en calculant un \idm
galoisien qui raffine le précédent, et dans la nouvelle approximation
du \cdrz, la chose étrange a disparu.

Pour développer ce point de vue nous aurons besoin de mieux connaître l'\aduz, et la section \ref{sec1adu} est consacrée à cet objectif.

Par ailleurs, nous étudierons dans la section \ref{subsecCDR} une version
dynamique et \cov du \cdr d'un \pol non \ncrt \splz.

%%%%%%%%%%%%%%%%%%%%%%%%%%%%%%%%%%%%%%%%%%%%%%%%%%%%%%%%%%%%%%%%%%%%%%%%%%%
%%%%%%%%%%%%%%%      algebres de Boole \label{secBoole}   %%%%%%%%%%%%%%%%%
%%%%%%%%%%%%%%%%%%%%%%%%%%%%%%%%%%%%%%%%%%%%%%%%%%%%%%%%%%%%%%%%%%%%%%%%%%%
%-- SECTION{L'\adu pour un \pol unitaire}--secadu
\section{Introduction aux \agBs} \label{secBoole}

 Un \ix{treillis} est un ensemble $\gT$ muni d'une relation d'ordre $\leq$
pour laquelle il existe un \elt minimum, noté $0_\gT$,  un \elt maximum,
noté $1_\gT$, et toute paire d'\elts $(a,b)$ admet une borne supérieure,
notée  $a\vu b$, et une  borne inférieure, notée $a\vi b$.
Une application d'un treillis vers un autre est
appelé un \emph{\homo de  treillis} si elle respecte les lois
$\vu$ et $\vi$ ainsi que les constantes $0$ et $1$.
Le treillis est appelé un
\ixx{treillis}{distributif} lorsque chacune des deux lois $\vu$ et $\vi$ est
distributive par rapport à l'autre.

Nous ferons une étude succincte de la structure de \trdi et de structures
qui s'y rattachent au chapitre~\ref{chapTrdi}.

%:     propdef defiBoole
\begin{propdef}\label{defiBoole}
\emph{(Algèbres de Boole)}\index{algèbre!de Boole}\index{Boole!algèbre de ---}
\begin{enumerate}
\item Par \dfn un anneau $\gB$ est une \emph{\agBz} \ssi tout \elt
est idempotent. En conséquence $2=_\gB 0$ (car $2=_\gB 4$).
\item On peut définir sur $\gB$
une relation d'ordre $x\preceq y$ par: $x$ est multiple
de~$y$, \cad $\gen{x}\subseteq\gen{y}$.
Alors, deux \elts arbitraires admettent une borne inférieure,
leur ppcm $x\vi y=xy$, et une borne supérieure, leur pgcd $x\vu y=x+y+xy$.
On obtient ainsi un \trdi avec~$0$ pour \elt minimum et $1$ pour \elt
maximum.
\item Pour tout $x\in\gB$, l'\elt $x'=1+x$ est l'unique \elt qui vérifie
les \egts $x\vi x'=0$ et $x\vu x'=1$, 
on l'appelle \emph{le complément} de $x$.
\end{enumerate}
\index{complement@complément!(dans une \agBz)}
\end{propdef}

{\it Conflit de notation.} On se retrouve ici avec un conflit de notation. En effet, la \dve dans un anneau  conduit à une notion
de pgcd, qu'il est usuel de noter $a\vi b$, car il est pris pour une borne inférieure ($a$ divise $b$ étant compris comme \gui{$a$ plus petit que $b$} au sens de la \dvez), en conflit avec le pgcd
des \elts dans une \agBz, qui est une borne supérieure. Cela tient à ce que la relation
d'ordre a été renversée pour que les \elts $0$ et $1$ de l'\agB
soient bien le minimum et le maximum dans le treillis.
Ce conflit, inévitable, apparaîtra plus fort encore lorsque l'on considérera l'\agB 
des \idms d'un anneau $\gA$.
\eoe

\smallskip 
Bien que tous les \elts d'une \agB soient \idms nous garderons la terminologie
de \gui{\sfioz\footnote{Il serait plus naturel de dire: \sys fondamental d'\elts
\ortsz.}} pour une famille finie $(x_i)$ d'\elts deux à deux \orts (\cadz~$x_ix_j=0$
pour $i\neq j$) dont la somme fait $1$.
Cette convention est d'autant plus justifiée que nous nous préoccuperons surtout de l'\agB qui apparaît naturellement en \alg commutative: celle des \idms
d'un anneau $\gA$.

%%%%%%%%%%%%%%%%%%%%%%%%%%%%%%%%%%%%%%%%%%%%%%%%%%%%%%%%%%%%%%%%%%%%%%%%%%%
%:--- SUBsec{Algebres de Boole discretes}
\subsec{Algèbres de Boole discrètes}
\label{subsec1AGBDiscretes}

%:     Proposition{propBoolFini}
\begin{proposition}\label{propBoolFini} \emph{(Toute \agB discrète
se comporte dans les calculs comme l'\alg des parties détachables d'un ensemble
fini)}\\
Soit $(r_1,\ldots,r_m)$  une famille finie dans une \agB $\gB$.\\
Nous posons $s_i=1-r_i$ et,
pour une partie finie $I$ de $\{1,\ldots ,m\}$, nous notons
$r_I=\prod_{i\in I}r_i\prod_{j\notin I}s_j$. %Alors:
\begin{enumerate}
\item  Les $r_I$ forment un \sfio et ils engendrent la même \agB
que les $r_i$.
\item Supposons  que $\gB$ soit discrète. Alors, s'il y a exactement $N$
\eltsz~$r_I$ non nuls, la sous-\agB engendrée par les $r_i$
est isomorphe à l'\alg des parties finies d'un ensemble à $N$ \eltsz.
\end{enumerate}
\end{proposition}

Comme corolaire on obtient le fait suivant et le \tho de structure fondamental
qui le résume.
Rappelons que l'on note $\Pf(S)$ l'ensemble des parties finies d'un ensemble $S$.

Dans une \agB discrète un \elt $e$ est appelé un \ix{atome} s'il vérifie l'une des \prts \eqves suivantes.
\begin{itemize}
\item $e$ est minimal parmi les \elts non nuls.
\item  $e\neq 0$ et pour tout $f$, $f$ est \ort ou supérieur à $e$.
\item  $e\neq 0$ et pour tout $f$, $ef=$ $0$ ou $e$,
ou encore $ef=0$ ou $e(1-f)=0$.
\item $e\neq 0$ et
une \egt $e=e_1+e_2$ avec $e_1e_2=0$ implique $e_1=0$ ou $e_2=0$.
\end{itemize}
On dit aussi que $e$ est \emph{indécomposable}. Il est clair qu'un \auto
d'une \agB discrète conserve l'ensemble des atomes  et que pour deux atomes $e$ et $f$,
on a $e=f$ ou $ef = 0$.\index{indecomp@indécomposable!\elt --- dans une \agBz}

%:     Theorem{corpropBoolFini}
\begin{theorem}\label{corpropBoolFini}\label{factagb} \emph{(\Tho de structure)}
\begin{enumerate}
\item Toute \agB finie
est isomorphe à l'\alg des parties détachables d'un ensemble
fini. 

\item 
Plus \prmtz, pour une \agB $C$ \propeq
%-----------------begin enum------------------
\begin{enumerate}
\item $C$ est finie.
\item $C$ est discrète et \tfz.
\item L'ensemble $S$ des atomes est fini, et  $1_C$ est la somme de cet ensemble.
\end{enumerate}
%-----------------end enum------------------
Dans un tel cas $C$ est isomorphe à l'\agB $\Pf(S)$.
\end{enumerate}
\end{theorem}

%:--- SUBsec{Algèbre de Boole des \idms d'un anneau commutatif}
\subsect{Algèbre de Boole des \idms d'un anneau \\commutatif}{Algèbre de Boole des \idms d'un anneau}
\label{subsecAGBIDMSAC}
%:     fact{propB(A)}
\begin{fact}\label{propB(A)}\relax
Les \idms d'un anneau $\gA$ forment une \agBz, notée $\BB(\gA)$, avec les lois $\vi$, $\vu$, $\lnot$ et $\oplus$ données par

\snic{r\vi s=rs$, $\;r\vu s =r+s-rs\;$ , $\;\lnot\ r=1-r\;$ et  $\;r\oplus s=(r-s)^2.}

%\sni

%En outre
%on a $(r-s)^3 = r-s$ et donc  $\gen {r\oplus s}_\gA = \gen {r-s}_\gA$.\\
Si $\gA$  est  une \agBz, $\BB(\gA)=\gA$.
Si $\varphi :\gA\to\gB$ est un morphisme d'anneaux, sa restriction à
$\BB(\gA)$ donne un morphisme $\BB(\varphi ):\BB(\gA)\to\BB(\gB)$.
\end{fact}
%---- end fact ----------------------------------------------
%
\begin{proof}
Il suffit de montrer que si l'on munit l'ensemble $\BB(\gA)$ des lois $\oplus$
et $\times$ on obtient une \agB avec $0_\gA$ et $1_\gA$ comme \elts neutres. Les calculs sont laissés \alecz.
\end{proof}

\medskip Le \thrf{corpropBoolFini} a la conséquence immédiate suivante.

%     Fact{factZerloc2}
\begin{fact}\label{factZerloc2}
\Propeq
\begin{enumerate}
\item L'\agB des \idms $\BB(\gA)$ est finie.
\item L'anneau $\gA$ est un produit fini d'anneaux connexes non triviaux.
\end{enumerate}
\end{fact}
\begin{proof}
Il suffit de montrer que \emph{1} implique \emph{2.}
Si $e$ est un atome de $\BB(\gA)$, l'anneau~$\gA[1/e]$ est non trivial et connexe.
Si $\BB(\gA)$ est finie, l'ensemble fini~$A$ de ses
atomes forme un \sfio de $\gA$, et l'on a  un \iso canonique $\gA\to\prod_{e\in A}\gA[1/e]$.
\end{proof}

\rem Si $\BB(\gA)$ a un seul \eltz, $\gA$ est trivial et le produit
fini est un produit vide. Ceci s'applique aussi pour le corolaire
suivant.
\eoe

%     Corollary{factZerloc2}
\begin{corollary}\label{lemZerloc2}
\Propeq
%-----------------begin enum------------------
\begin{enumerate}
\item $\BB(\gA)$ est finie et  $\gA$ est \zedz.
\item   $\gA$ est un produit fini d'anneaux locaux \zeds non triviaux.
\end{enumerate}
%-----------------end enum------------------
\end{corollary}

%%%%%%%%%%%%%%%%%%%%%%%%%%%%%%%%%%%%%%%%%%%%%%%%%%%%%%%%%%%%%%%%%%%%%%%%%%%
%:--- SUBsec{Elements galoisiens dans une \agB}
\subsec{\'Eléments galoisiens dans une \agB}

%:     Definition{defIdmGalAgb}-----
\begin{definition}
\label{defIdmGalAgb}~
\index{G-alg@$G$-algèbre de Boole}
\index{Boole!G-alg@$G$-algèbre de ---}
\index{G-alg@$G$-algèbre de Boole!transitive}
\index{transitive!G-alg@$G$-algèbre de Boole}
\begin{enumerate}
\item Si $G$ est un groupe qui opère sur une \agB $C$, on dit que 
le couple $(C,G)$
est une $G$-\agBz.
\item Un \elt $e$ d'une $G$-\agB $C$ est dit \ixc{galoisien}{\elt --- dans une \agBz}  si son orbite sous $G$ est un \sfioz.
\item Une $G$-\agB est dite \emph{transitive} si $0$ et $1$
sont les seuls \elts fixés par $G$.
\end{enumerate}
\end{definition}
%--- end-definition------------------------------------

Nous étudions maintenant le cas où le groupe est fini et l'\alg discrète.

%--- Fact{factAduIdmA}------------
\begin{fact}
\label{factAduIdmA}
Soit  $G$ un groupe fini et $C$ une $G$-\agB transitive, discrète et non
triviale. Soit $e\neq 0$ dans  $C$, et $\so{e_1,\ldots ,e_k}$ l'orbite de $e$ sous~$G$.
\Propeq
%-----------------begin enum------------------
\begin{enumerate}
\item \label{E1} L'\elt $e$ est  galoisien.
\item \label{E2} Pour tout $i>1$, $e_1e_i=0$.
\item \label{E21} Pour tout $\sigma\in G$, $e\sigma(e)=e$  ou $0$.
\item \label{E22} Pour tous $i\neq j\in\so{1,\ldots ,k}$, $e_ie_j=0$.
\end{enumerate}
%-----------------end enum------------------
\end{fact}
%--- end-fact-----------------------------------------
%-----------------begin proof------------------
\begin{proof}
Le point \emph{\ref{E1}} implique clairement les autres. Les points \emph{\ref{E2}} et
\emph{\ref{E22}} sont facilement \eqvs et impliquent le point \emph{\ref{E21}}.
Le point \emph{\ref{E21}} signifie que pour tout $\sigma$, $\sigma(e)\geq e$ ou
$\sigma(e)e=0$. Si l'on a  $\sigma(e)\geq  e$ pour un certain $\sigma$,
alors on obtient 

\snic{e\leq \sigma(e)\leq \sigma^2(e)\leq \sigma^3(e)\leq \ldots\,, }

%\sni
ce
qui donne $e=\sigma(e)$  en considérant un $\ell$ tel que $\sigma^\ell=1_G$.
Donc, le point~\emph{\ref{E21}} implique le point~\emph{\ref{E2}}.
Enfin si le point~\emph{\ref{E22}} est vérifié, la somme de l'orbite est un
\eltz~$>0$
fixé par~$G$ donc égal à~$1$.
\end{proof}
%-----------------end proof------------------

%:     Lemma{lemIdmGalAgb}
\begin{lemma}\label{lemIdmGalAgb} \emph{(Rencontre de deux \elts galoisiens)}\\
Soit $G$ un groupe fini et une $G$-\agB $C$, discrète et non
triviale. \'Etant donnés deux \elts galoisiens  $e$, $f$ dans $(C,G)$, nous noterons

\snic{G.e=\so{e_1,\ldots,e_m}$, $E=\St_G(e)$, et
$F=\St_G(f).}
\begin{enumerate}
\item Il existe $\tau\in G$ tel que $f\tau(e)\neq0$.
\item
Si $e\leq f$, alors $E\subseteq F$ et $f=\sum_{i:e_i\leq
f}e_i=\sum_{\sigma\in F/E}\sigma(e)$.
\end{enumerate}
On suppose $C$ transitive et $ef\neq0$. On obtient les résultats suivants.
\begin{enumerate}\setcounter{enumi}{2}
\item
L'\elt $ef$ est galoisien, de stabilisateur $E\,\cap\, F$, et l'orbite $G.ef$ est
constituée des \elts non nuls de $(G.e)(G.f)$.  En particulier, $G.ef$
engendre la même sous-\agB de $C$ que $G.e\,\cup\, G.f$.
\item Si  $E\subseteq F$, alors $e\leq f$.
\end{enumerate}
\end{lemma}
\begin{proof}
\emph{1.} En effet, $f=\sum_ife_i$.

\emph{2.}
De manière générale, pour $x' = \sigma(x)$ où $x \ne 0$ vérifie $x
\le f$, montrons:
$$x' \le f \;\;\buildrel {[a]} \over \Longrightarrow\;\;
fx' \ne 0 \;\;\buildrel {[b]} \over \Longrightarrow\;\;
\sigma(f) = f \;\;\buildrel {[c]} \over \Longrightarrow\;\;
x' \le f.\leqno(\star)
$$
On obtient $[a]$ en multipliant $x' \le f$ par $x'$, $[b]$ en multipliant $x'
\le \sigma(f)$ (déduite de $x \le f$) par $f$ et en utilisant $f$
galoisien, et enfin $[c]$ en appliquant~$\sigma$ à $x \le f$. Les assertions
de $(\star)$ sont donc des \eqvcsz. 
On en déduit~$\St_G(x) \subseteq
\St_G(f)$. Si de plus, $1 = \sum_{x' \in G.x} x'$, alors:

\snic {
f = \sum_{x' \in G.x} fx' = \sum_{x'\in G.x | x' \le f} x' =
\sum_{\sigma \in F/\St_G(x)} \sigma(x).
}

%\sni
Ceci s'applique à $x = e$.

\emph{3.} Notons  $G.f=\so{f_1,\ldots,f_p}$.
Pour $\sigma\in G$ il existe $i$, $j$ tels que 
$$\preskip.4em \postskip.4em 
e\,f\,\sigma(ef) = e\,f\,e_i\,f_j, 
$$
 qui est égal à $ef$ si $\sigma\in E\cap F$ et à $0$ sinon. D'après le fait \ref{factAduIdmA}, $ef$ est donc un \elt galoisien de stabilisateur 
 $E \cap F$.  Supposons maintenant $e_if_j\neq0$. Alors, d'après le point~\emph{1},
il existe $\tau\in G$ tel que $\tau(ef) e_if_j\neq0$. Puisque~$e$ et~$f$
sont galoisiens, ceci implique $\tau(e)=e_i$ et $\tau(f)=f_j$, donc~$e_if_j
\in G.ef$.

\emph{4.} Résulte \imdt de \emph{3.}
\end{proof}

L'exemple paradigmatique d'application du prochain \tho est le suivant.
On considère  un anneau connexe non trivial $\gk$, $(\gk,\gC,G)$ est une \aG ou 
prégaloi\-sienne (cf. \dfn \ref{defADG1}) \hbox{et  $C=\BB(\gC)$}.

%:    Theorem{lemAduIdmA} struct 1
\begin{theorem}
\label{lemAduIdmA} \emph{(\Tho de structure galoisien, 1)} 
Soit  $G$ un groupe fini et  $C$ une $G$-\agB transitive, discrète et non
triviale.
%-----------------begin enum------------------
\begin{enumerate}
% 1
\item \label{lemAduIdmA-1}\relax \emph{(Structure des $G$-\agBs finies transitives)}
\\
L'\alg $C$ est finie \ssi il existe un atome $e$. Dans
ce cas la structure de $(C,G)$ est entièrement \caree par  $E=\St_G(e)$.
\\
Plus \prmtz, l'\idm $e$ est galoisien, $G.e$ est
l'ensemble des atomes, {$C\simeq\Pf(G.e)$,}
$G$ opère sur $G.e$ comme sur $G/E$, et sur
$C$  comme \hbox{sur $\Pf(G/E)$}. En particulier, $|C|=2^{\idg{G:E}}.$
\\
On dira que $e$ est un \emph{\gtr galoisien de~$C$}.
% 2
\item \label{lemAduIdmA-2}\relax Toute famille finie d'\elts de $C$
engendre une sous-$G$-\alg finie.
% 3
\item \label{lemAduIdmA-3}\relax L'\agB $C$ ne peut avoir plus que $2^{|G|}$
\eltsz.
% 4
\item \label{lemAduIdmA-4}\relax Soient $e$ et $f$ des \elts galoisiens, $E=\St_G(e)$
et $F=\St_G(f)$.
\begin{enumerate}
% a
\item Il existe $\sigma\in G$ tel que $f\sigma(e)\neq0$.
% b
\item Si $ef\neq0$, $ef$ est un \gtr galoisien de la sous-$G$-\agB de $G$ engendrée par $e$ et $f$, et  
 $\St_G(ef)=E\cap F$.
% c
\item Si $e\leq f$ (i.e. $fe=e$), alors $E\subseteq F$ et $f=\sum_{\sigma\in
F/E}\sigma(e)$.
% d
\end{enumerate}
% 5
\item \label{lemAduIdmA-5}\relax \emph{(\Carn des \elts galoisiens dans une sous $G$-\alg finie)}
\\
Soit $e$ un \elt galoisien  et $f$ une somme de
$r$ \elts de $G.e$, dont~$e$. Soit $E=\St_G(e)$  et $F=\St_G(f)$. Alors \propeq
%-----------------begin enum------------------
\begin{enumerate}
% a
\item \label{lemAduIdmA-5a}$f$ est galoisien.
% b
\item \label{lemAduIdmA-5b}$E\subseteq F$ et  $f=\sum_{\sigma\in
F/E}\sigma(e)$.
% c
\item \label{lemAduIdmA-5c}$\abs{F}=r\times \abs{E}$.
% d
\item \label{lemAduIdmA-5d}$\abs{F}\geq r\times\abs{E} .$
\end{enumerate}
%-----------------end enum------------------
%
\end{enumerate}
%-----------------end enum------------------
\end{theorem}
%--- end-theorem-----------------------------------------

%-----------------begin proof------------------
\begin{proof}
\emph{\ref{lemAduIdmA-1}.} Si $C$ est finie il existe un atome.
Si $e$ est un atome, pour tout $\sigma\in G$, on a $e\,\sigma(e)=0$ ou $e$, donc
$e$ est galoisien (fait \ref{factAduIdmA}). Le reste en découle en tenant
compte du \thref{factagb}.
%:HHH 11 au dessus la reference est maintenant theorem au lieu de fait

\emph{\ref{lemAduIdmA-2}.} On considère la sous-\agB $C'\subseteq
C$ engendrée par les orbites des \elts de la famille finie donnée.
$C'$ est \tf et discrète donc finie. En conséquence ses atomes
 forment un ensemble fini $S=\so{e_1,\ldots,e_k}$ et $C'$ est isomorphe
à l'\agB des parties finies de $S$: 

\snic{C'=\sotq{\sum_{i\in
F}e_i}{F\in\cP_k}.}

%\sni
Clairement, $G$ opère sur $C'$. Pour
$\sigma\in G$, $\sigma (e_1)$ est un atome, donc~$e_1$
est galoisien (fait~\ref{factAduIdmA}~\emph{\ref{E21}}) et $(e_1,\ldots,e_k)$ est son orbite.

\emph{\ref{lemAduIdmA-3}.} Résulte de \emph{\ref{lemAduIdmA-1}} et \emph{\ref{lemAduIdmA-2}}.

\emph{\ref{lemAduIdmA-4}.}  Déjà vu dans le lemme \ref{lemIdmGalAgb}.

\emph{\ref{lemAduIdmA-5}.} On écrit $\sigma_1=1_G$,
$G.e=\so{\sigma_1.e,\ldots,\sigma_k.e}$ avec $k=\idg{G:E}$, 
ainsi que~$f=\sigma_1.e+\cdots +\sigma_r.e$.
\\
\emph{\ref{lemAduIdmA-5a}} $\Rightarrow$ \emph{\ref{lemAduIdmA-5b}.} On applique le
point \emph{\ref{lemAduIdmA-4}.}
\\
\emph{\ref{lemAduIdmA-5b}} $\Rightarrow$ \emph{\ref{lemAduIdmA-5a}.} Pour $\tau\in F$,
 $\tau.f=f$. 
 \\
Pour
$\tau\notin F$, $F.e\cap (\tau F).e=\emptyset$, et donc~$f\,\tau (f)=0$.
\\
\emph{\ref{lemAduIdmA-5b}} $\Rightarrow$ \emph{\ref{lemAduIdmA-5c}.} On a
$F.e=\so{1_G.e,\sigma_2.e,\ldots,\sigma_r.e}$, et puisque $E$ est le
stabilisateur de $e$,  on obtient $\abs{F}=r\times \abs{E}$.
\\
\emph{\ref{lemAduIdmA-5d}} $\Rightarrow$ \emph{\ref{lemAduIdmA-5b}.}
On a $F=\sotq{\tau}{\tau\so{\sigma_1.e,\ldots,
\sigma_r.e}=\so{\sigma_1.e,\ldots,\sigma_r.e}}$. D'où l'inclusion $F.e\subseteq
\so{\sigma_1.e,\ldots,\sigma_r.e}$, et %l'on peut écrire
$F.e=\so{\sigma_1.e,\ldots,\sigma_s.e}$ avec $s\leq r\leq k$. Le
stabilisateur de $e$ pour l'action de $F$ sur $F.e$ est égal
à $E\cap F$.
Donc
$$\abs{F}=\abs{F.e}\,\abs{E\cap F}= s\,\abs{E\cap F} \leq r\,\abs{E\cap
F}\leq r\,\abs{E}.$$
Donc si $\abs{F}\geq r\,\abs{E}$, on a $\abs{F.e}=r$ et $\abs{E}=\abs{E\cap
F}$, \cad $E\subseteq F$ \linebreak 
et $F.e=\so{\sigma_1,\ldots ,\sigma_r}$.
\end{proof}
%-----------------end proof------------------
%\newpage

%-% PERSO
\perso{Le point \emph{\ref{lemAduIdmA-5}} comporte quelques subtilités et il serait
bon de mettre un commentaire adéquat qui explique pourquoi on se fatigue
à ce point. Notamment il faudrait renvoyer à un endroit  où on utilise ce résultat. L'\algo suivant
éclaire un peu la chose}
%-% Fin PERSO

%:     Algorithme{algidmgal}-------
\begin{algor}[Calcul d'un \elt galoisien et de son
stabilisateur.] \label{algidmgal}
\Entree $e$: \elt non nul d'une \agB $C$; $G$: groupe fini d'\autos de $C$;
$S_e=\St_G(e)$.
\\ \# \, \emph{On suppose que $0$ et $1$  sont les seuls points fixes pour l'action de $G$  sur $C$.}
\Sortie $e_1$: un \elt galoisien de $C$ tel que $G.e_1$ engendre la même \agB que $G.e$; $H$: le sous-groupe
stabilisateur de $e_1$.
\Varloc $h$: dans $C$; $\sigma$: dans $G$;
%\\ 
$L$: liste d'\elts de $\ov{G/S_e}$;\\
$E$: ensemble correspondant d'\elts de $G/S_e$;
\\ \# \, $G/S_e$  est l'ensemble des classes
à gauche modulo $S_e$. 
\\ \# \, $\ov{G/S_e}$ est un \sys de représentants des classes
à gauche modulo $S_e$
\Debut
\hsu $E \aff \emptyset$; $L \aff [\,]$; $e_1\aff e$; 
\hsu \pur{\sigma}{\ov{G/S_e}}
\hsd $h\aff e_1\sigma(e)$;
\hsd \sialors{h\neq 0} $e_1\aff h$; $L\aff L\bullet [\sigma]$;
$E \aff E \cup \so{\sigma S_e}$
\hsd \finsi;
\hsu \finpour;
\hsu $H\aff \St_G(E)$\quad  \# \, $H =\sotq{\alpha\in G} 
 {\forall\sigma\in L,    \alpha\sigma\in \bigcup_{\tau\in L}\tau S_e}$.
\Fin
\end{algor}
%--- fin Algorithme ------------------

Sous les hypothèses du \thref{lemAduIdmA} on peut calculer un \elt
galoisien~$e_1$ tel que $G.e_1$
 et $G.e$ engendrent  la même \agBz, au moyen de l'\algo
\vref{algidmgal}.
%:HHH 11 le commentaire ci-après est prématuré 
%On pourra penser au cas où  $C=\BB(\gC)$,  $(\gk,\gC,G)$ est 
%une \apG  et $\gk$ est connexe. 
%En outre, on peut \egmt calculer le stabilisateur de $e_1$ (la
%\emph{nouvelle approximation du groupe de Galois}) \gui{sans sortir du
%groupe}.

%-----------------begin proof------------------
\begin{Proof}{Correction de l'\algoz. }
Nous avons noté $G/S$ un \sys de représentants des classes à gauche
modulo $S$. \'Ecrivons $e_1=e \sigma_2(e)\cdots \sigma_r(e)$ où les~$\sigma_i$ sont tous les $\sigma$ qui ont passé avec succès le test
$h\neq 0$ dans l'\algo (et $\sigma_1=\Id$). Nous voulons montrer que $e_1$
est un atome de $C'$ (l'\agB engendrée par $G.e$),
ce qui revient à dire que pour tout
$\sigma\in G/S$ on a $e_1\sigma(e)=e_1$ ou $0$ (puisque $C'$ est engendrée
par les $\tau(e)$). Or $\sigma$ a été testé par l'\algoz, donc ou bien
$\sigma$ est l'un des $\sigma_i$, auquel cas  $e_1\sigma(e)=e_1$,
ou bien $g\sigma(e)=0$ pour un \idm $g$  qui divise $e_1$, et a fortiori
$e_1\sigma(e)=0$.\\
Montrons que le stabilisateur $H$ de $e_1$ vérifie bien la condition
requise. On a $e_1=\prod_{\tau\in L}\tau(e)$, et pour $\sigma \in G$
on a les \eqvcsz:

\snic{
\sigma\in \bigcup_{\tau\in L}\tau S \;\Longleftrightarrow\;
e_1\sigma(e)=e_1\;\Longleftrightarrow\; e_1\leq
\sigma(e),\qquad\mathrm{et}}

\snic{
\sigma\notin \bigcup_{\tau\in L}\tau S \;\Longleftrightarrow\; e_1\sigma(e)=0.}

%\sni
Pour $\alpha\in G$ on a
$\alpha(e_1)=\prod_{\tau\in L}\alpha\big(\tau(e)\big)$. C'est un \elt de l'orbite de
$e_1$, il est égal à $e_1$ \ssi $e_1\leq \alpha(e_1)$, \ssi
$e_1\leq \alpha\big(\sigma(e)\big)$ pour chaque $\sigma$ in $L$. Enfin, pour un
$\sigma$
arbitraire dans $G$, $e_1\leq \alpha\big(\sigma(e)\big)$  \ssi
$\alpha \sigma $ est dans $\bigcup_{\tau\in L}\tau S$.
\end{Proof}
%-----------------end proof------------------

On notera que l'\elt $e_1$ obtenu comme résultat du calcul dépend de
l'ordre dans lequel est énuméré l'ensemble fini $G/S$ et qu'il n'y a
pas d'ordre naturel (intrinsèque) sur cet ensemble.

\medskip \exl
On peut se demander s'il existe un rapport  entre le stabilisateur~$S$ de $e$ et le stabilisateur $H$ d'un \elt galoisien $e_1$ associé à
$e$.
Voici un exemple qui montre qu'il n'y a pas de rapport étroit, avec
$G=\rS_6$ opérant sur $\Adu_{\QQ,f}$ avec le \pol $f(T)=T^6-4T^3+7$.
On considère   l'\idm 
$e=1/6(x_5^3x_6^3 - 2x_5^3 - 2x_6^3 + 7)$
que l'on calcule  partir d'une factorisation du
\polmin de l'\elt $x_5+x_6$  (cf. proposition \ref{propRvRel1}). \\
On trouve
$\St(e)=S=\gen{(1432),(12),(56)}\simeq \rS_4\times \rS_2$ avec $|S|=48$, et

\snac{\St(e_1)=H=\gen{(24),(123456)}=(\gen{(13),(135)}\times
\gen{(24),(246)})\rtimes \gen{(14)(25)(36)}
}

%\sni
 avec $H\simeq (\rS_3\times \rS_3)\rtimes
\rS_2$, $|H|=72$, et $S\cap H=\gen{(24),(1234)(56)}$ diédral
d'ordre 8.
\perso{cela serait bien de donner un autre exemple aboutissant a un sous-groupe $H$ complètement différent (par exemple $H=S$).}\\
En bref, $H$ (ni même la classe de conjugaison de $H$  dans $G$)
ne peut  être calculé à partir de $S$ seulement. En effet, la liste $L$  de
classes à gauche sélectionnée par l'\algo
ne dépend pas seulement du sous-groupe $S$  de~$G$ mais aussi
de la façon dont $G$  opère sur $C$.
\eoe

%\newpage

%%%%%%%%%%%%%%%%%%%%%%%%%%%%%%%%%%%%%%%%%%%%%%%%%%%%%%%%%%%%%%%%%%%%%%%%%%%
%%%%%%%%%%%%%%%      L'\adu pour un \pol unitaire         %%%%%%%%%%%%%%%%%
%%%%%%%%%%%%%%%%%%%%%%%%%%%%%%%%%%%%%%%%%%%%%%%%%%%%%%%%%%%%%%%%%%%%%%%%%%%
%-- SECTION{L'\adu pour un \pol unitaire}--secadu
%:2012 j'ai supprimé la , avant (2) dans le titre
\section{L'\adu  (2)}
%{L'\adu \\ pour un \pol unitaire sur un anneau commutatif (2)} 
\label{sec1adu}
\label{secadu}

Voici un petit guide de lecture pour la fin de ce chapitre.

Dans la section \ref{secGaloisElr}, nous avons vu que si $\gk$ est un corps
discret infini, si $f$ est \spl et si l'on est capable de décomposer une
résolvante de Galois en produit de facteurs \irdsz,
alors l'\adu $\gA$ est isomorphe à $\gL^r$, avec $\gL$  un \cdr pour $f$
\linebreak 
 et
$r=\idg{\Sn:G}$, où $G$ est un sous-groupe de $\Sn$ qui s'identifie à
$\Gal(\gL/\gk)$. En outre, $\dex{\gL:\gk}=\abs G$.

Nous allons voir que cette situation idéale peut servir de
ligne directrice
pour une approche paresseuse de la construction d'un \cdrz.
Ce qui remplace la \fcn complète d'une résolvante de Galois,
c'est la découverte ou la construction d'un \idm galoisien.
Alors, on a une situation analogue à la situation
idéale décrite auparavant: $\gA\simeq\gB^r$, où
$\gB$ est un 
{quotient de Galois} de $\gA$, muni d'un groupe d'\autos
qui s'identifie à un sous-groupe $G$ de $\Sn$, avec $\dex{\gB:\gk}=\abs G$
et $r=\idg{\Sn:G}$.

\medskip

\Grandcadre{
Dans toute la section \ref{secadu},  $\gk$ est un anneau commutatif, \\
$f=T^n+\sum_{k=1}^n(-1)^{k} s_kT^{n-k} \in \kT$ est  \mon de
degré $n$,\\
et $\gA =\Adu_{\gk,f}$ est l'\adu
de $f$ sur $\gk$.
}

\smallskip 
Rappelons que l'\adu

\snic{\gA =\Adu_{\gk,f}=\aqo{\gk[\uX]}{S_1-s_1,\ldots,S_n-s_n}=\kuX/\cJ(f)}

%\sni
(où les $S_i$ sont les \pols \smqs \elrs en les $X_i$)
est l'\alg qui résout le
\pb \uvl lié à la \dcn du \pol $f$ en un produit de
facteurs $T-\xi_j$ (cf. fait \ref{factEvident}).
Le \kmo $\gA=\Adu_{\gk,f}$ est libre, et une base est formée par les \gui{\momsz}
$x_1^{d_1}\cdots x_{n-1}^{d_{n-1}}$ tels que pour $k\in\lrb{0..n-1}$, on ait
$d_k\leq n-k$ (voir fait \ref{factBase}). Nous noterons cette base~$\cB(f)$.
\label{notaBaseadu}

Par changement d'anneau de base, on obtient
le fait important suivant (à distinguer du fait \ref{factAduAdu}).
%--- Fact{factChangeBase}---------
\begin{fact}
\label{factChangeBase}\emph{(Changement d'anneau de base pour une \aduz)}
Soit $\rho:\gk\to\gk_1$ une \klgz. Notons $f^{\rho}$ l'image de $f$
dans $\gk_1[T]$. Alors, l'\alg  $\rho\ist(\Adu_{\gk,f})=\gk_1\otimes_\gk \Adu_{\gk,f}$, est
naturellement isomorphe à  $\Adu_{\gk_1,f^{\rho}}$.
\end{fact}
%--- end-fact-----------------------------------------

%%%%%%%%%%%%%%%%%%%%%%%%%%%%%%%%%%%%%%%%%%%%%%%%%%%%%%%%%%%%%%%%%%%%%%%%%%%
%:--- SUBsec{Quotients de Galois de l'adu}  subsecIdGal
\penalty-2500
\subsec{Quotients de Galois des \apGsz} \label{subsecIdGal}

\medskip Si $\gC$ est une \klgz, on note $\Aut_\gk(\gC)$ son groupe d'\autosz.

Nous donnons maintenant une \dfn qui permet d'insérer l'\adu
dans un cadre un peu plus \gnl et utile.

%:    Definition{defADG1}----- defi alg pregaloisienne
\begin{definition}
\label{defADG1} \emph{(Algèbres prégaloisiennes)}\\
Une \emph{\apGz} est donnée par un triplet  $(\gk,\gC,G)$ où
%-----------------begin enum------------------
\begin{enumerate}
\item $\gC$ est une \klg avec $\gk\subseteq\gC$, $\gk$ facteur direct
dans $\gC$,
\item $G$ est un groupe fini de $\gk$-\autos de $\gC$,
\item $\gC$ est un \kmo \pro de rang constant $|G|$,
\item pour tout $z\in\gC$, on a $\rC {\gC/\gk}(z)(T)=\rC {G}(z)(T)$.
\end{enumerate}
%-----------------end enum------------------
\index{algèbre!prégaloisienne}
\end{definition}
%--- end-definition------------------------------------
\rem Rappelons que d'après le lemme~\ref{lemIRAdu},
si $\gB$ est une \klg \stfe et fidèle, alors $\gk$ (identifiée à son image dans $\gB$) est facteur direct dans $\gB$. En conséquence
le point \emph{1} ci-dessus résulte du point
\emph{3.}  \eoe

\smallskip 
\exls  1) D'après ce que l'on sait déjà sur l'\adu (section \ref{sec0adu})
et d'après le lemme \ref{lemPolCarAdu},
pour tout \pol unitaire $f$, le triplet $(\gk,\Adu_{\gk,f},\Sn)$ est une \apGz.

2) Le \tho d'Artin \ref{thA} montre que toute \aG
est une \apGz.
\eoe

\medskip  \Llec se reportera \paref{defQuoDeGal} pour les \dfns d'\idm galoisien,
d'\id galoisien et de quotient de Galois.

%:    Theorem{thADG1Idm} struct 2
\begin{theorem}
\label{thADG1Idm} \emph{(\Tho de structure galoisien, 2)}\\
Considérons une \apG $(\gk,\gC,G)$. Soit $e$ un \idm galoisien de~$\gC$, 
et~$\so{e_1,\ldots ,e_m}$ son orbite sous $G$. Soit~$H$ le stabilisateur
\linebreak 
 de $e=e_1$ et $r=|H|$, de sorte que
$rm=|G|$. Posons
$\gC_i={\gC}[1/e_i]$ \hbox{pour $(i\in\lrbm)$}. Soit enfin
$\pi:\gC\to\gC_1$ la projection canonique.
\begin{enumerate}
\item \label{thADG1Idm1} $(\gk,\gC_1,H)$ est une \apG (autrement dit  un quotient de Galois d'une \apG est une \apGz).
\item \label{thADG1Idm2} Les $\gC_i$ sont des \klgs deux à deux
isomorphes, et $\gC\simeq \gC_1^m$.
\item \label{thADG1Idm3} L'\alg $\gC_1$ est un \kmo \prc $r=|H|$. La restriction de $\pi$ à $\gk$,
et même à $\gC^G$, est injective. Et $\gk$ (identifié à son image
dans $\gC_1$) est facteur direct dans~$\gC_1$.
% 4
\item \label{thADG1Idm4} Le groupe $H$ opère sur $\gC_1$ et $\gC_1^H$ est
canoniquement isomorphe à $\gC^G$:
plus \prmtz, $\gC_1^H=\pi(\gC^H)=\pi(\gC^G)$.
% 5
\item \label{thADG1Idm5} Pour tout $z\in\gC_1$,
$\rC {\gC_1/\gk}(z)(T)=\rC {H}(z)(T)$.
% 6
\item \label{thADG1Idm6} Soit $g_1$ un \idm galoisien de  $(\gk,\gC_1,H)$,
$K$ son stabilisateur dans~$H$, $g'\in e_1\gC$ tel que $\pi(g')=g_1$. Alors,
$g'$ est un \idm galoisien de~$(\gk,\gC,G)$, son stabilisateur est $K$, et
l'on a un \iso canonique~$\aqo{\gC_1}{1-g_1}\simeq\aqo{\gC}{1-g'}$.
% 7
\item \label{thADG1Idm7} Si $(\gk,\gC,G)$ est une \aGz, alors $(\gk,\gC_1,H)$
\egmtz.
\end{enumerate}
%\index{quotient de Galois!d'une \apGz}
%-----------------end enum------------------
\end{theorem}
%--- end-theorem-------------------

%-----------------begin proof------------------
\begin{proof} Le point \emph{\ref{thADG1Idm1}} est une synthèse
partielle des points \emph{\ref{thADG1Idm2}},
\emph{\ref{thADG1Idm3}}, \emph{\ref{thADG1Idm4}}, \emph{\ref{thADG1Idm5}.}

Le point \emph{\ref{thADG1Idm2}} est évident.
 La première affirmation du point
\emph{\ref{thADG1Idm3}} en est une consé\-quence \imdez.  Soit
$\tau_1=\Id,\tau_2, \ldots, \tau_m$ un \sys de représentants pour
$G/H$, avec $\tau_i(e_1)=e_i$. Montrons que
la restriction de $\pi$ à $\gC^G$ est injective: si $a\in\gC^G$ et
$e_1a=0$,
alors, en transformant par les $\tau_j$, tous les~$e_ja$ sont nuls, et donc
aussi leur somme, égale à $a$.
Enfin $\pi(\gk)$ est facteur direct dans $\gC_1$ par le lemme~\ref{lemIRAdu}.

\emph{\ref{thADG1Idm4}.} 
Montrons d'abord $\gC_1^H=\pi(\gC^H)$.  Soit $z\in\gC_1^H$ et $u\in\gC$ tel
que $\pi(u) = z$.  Puisque $z\in\gC_1^H$, pour $\sigma\in H$,
$\sigma(u)\equiv u$ mod $\gen{1-e_1}$, i.e. $e_1\sigma(u)=e_1 u$
ou encore, puisque $\sigma(e_1) = e_1$, $\sigma(e_1u)=e_1 u$.  En posant $y =
e_1u$, on voit que~$y$ est $H$-invariant et $\pi(y) = z$.
\\
Montrons maintenant que $z\in\pi(\gC^G)$. On pose 

\snic{v=\sum_i\,\tau_i(y)
=\sum_i\,\tau_i(e_1y) =\sum_i\,e_i\tau_i(y).}

%\sni
Comme $\pi(e_i)=\delta_{1i}$, on
a $\pi(v)=\pi(y)$. L'\elt $v$ ainsi construit est indépendant du \sys de
représentants pour $G/H$. En effet, si $(\tau'_i)$ est un autre
\sys de représentants, quitte à réordonner les indices, on peut supposer
que $\tau'_iH = \tau_iH$, et donc, $y$ étant $H$-invariant, $\tau'_i(y) =
\tau_i(y)$.\\
 Pour $\sigma\in G$,
les $(\sigma \circ \tau_i)$ forment un \sys de représentants pour~$G/H$,
donc~$\sigma(v) = v$: l'\elt $v$ est $G$-invariant.

\emph{\ref{thADG1Idm5}.} %5
On a une \dcn $\gC=\gC'_1\oplus\cdots \oplus \gC'_m$, où $\gC'_j = e_j\gC$
est un~\hbox{\kmo} \ptf de rang $r$ et la restriction $\pi : \gC'_1 \to \gC_1$ est un
\iso de \kmosz. Pour tout $y \in \gC$, on~a:

\snac {
\rC{\gC/\gk}(y)(T) = \prod\limits_{j=1}^m \rC{\gC'_j/\gk}(e_jy)(T)  \; \hbox{ et }\;
\rC{G}(y)(T) = \prod\limits_{j=1}^m \prod\limits_{\tau\in H} \big(T - (\tau_j\circ \tau)(y)\big).
}

%\sni
Soit $y$ l'unique \elt de $\gC'_1$ tel que $\pi(y) = z$. L'\egt de gauche
donne 

\snic{\rC{\gC/\gk}(y)(T) = T^{(m-1)r} \rC{\gC_1/\gk}(z)(T).}

%\sni
Ensuite,
appliquons $\pi$ à l'\egt de droite en notant que $(\tau_j\circ\tau)(y)
\in \gC'_j$ (utiliser $y = e_1y$ et appliquer $\tau_j\circ\tau$). On obtient
alors: 

\snic{\rC{G}(y)(T) = T^{(m-1)r} \rC{H}(z)(T).}

%\sni
D'où $\rC{\gC_1/\gk}(z)(T) = \rC{H}(z)(T)$.

\emph{\ref{thADG1Idm6}.} %6
En tenant compte du fait que la restriction de
$\pi$ à $e_1\gC$ est un \iso on a $g'^2=g'=g'e_1$.
De même pour $\sigma\in H$ on~a: $\sigma(g')=g'$ si~$\sigma\in K$, ou
$g'\sigma(g')=0$ si $\sigma\notin K$.
Enfin pour $\sigma\in G\setminus H$, $e_1\sigma(e_1)=0$, 
et donc~$g'\sigma(g')=0$.
Ceci montre que $g'$ est un \idm galoisien de $\gC$ avec pour stabilisateur~$K$.
L'\iso canonique est immédiat.

\emph{\ref{thADG1Idm7}.} 
D'après le point \emph{4}, $\gk$ est l'ensemble des
points fixes. Il reste à voir que $H$ est séparant.
Si $\sigma\in H=\St(e)$ est distinct de l'\idtz, on a des
$a_i$ et des $x_i\in\gC$ tels que $\sum_ia_i(\sigma(x_i)-x_i)=1$.
Cette \egt reste vraie si on localise en $e$.
\end{proof}
%-----------------end proof------------------

%%%%%%%%%%%%%%%%%%%%%%%%%%%%%%%%%%%%%%%%%%%%%%%%%%%%%%%%%%%%%%%%%%%%%%%%%%%
%:--- SUBsec{Quand l'\agB d'une \adu est discrète}
\subsect{Cas où l'\agB d'une \alg de\\ \dcn \uvle est discrète}{Cas où l'\agB de $\Adu_{\gk,f}$ est discrète}
\label{subsec2AGBDiscretes}

Il est souhaitable que l'on puisse tester l'\egt
de deux \idmsz\hbox{ $e_1$, $e_2$} dans l'\adu $\gA$, 
ce qui revient à savoir tester~$e=0$ pour un \idm arbitraire 
de $\gA$ (comme dans tout groupe additif).
Or $e\gA$ est un \kmo \ptf et $e=0$ \ssiz $\rR {e\gA}(X)=1$
(\thrf{th ptf sfio} point \emph{\iref{remRang3}}).
Comme le \polmu $\rR {e\gA}$ peut être
calculé explicitement, on peut tester l'\egt
de deux \idms dans $\gA$ \ssi on peut tester l'\egt
de deux \idms dans $\gk$. L'argument ci-dessus fonctionne dans un cadre un
peu plus \gnl et l'on obtient le résultat suivant.

%--- Fact{factagbdisc}-----------
\begin{fact}
\label{factagbdisc}
Si $\BB(\gk)$ est une \agB discrète, il en va de même pour~$\BB(\gA)$.
Plus généralement, si  $\gC$ est une \klg \stfez, et si $\BB(\gk)$ est discrète, 
alors  $\BB(\gC)$ est discrète.
\end{fact}
%--- end-fact-----------------------------------------

%--- Fact{factIdemStable}-------
\begin{fact}
\label{factIdemStable}
Si $(\gk,\gC,G)$ est une \apGz,  tout \idm $e$ de~$\gC$  fixé par $G$ est un
\elt de~$\gk$.
\end{fact}
%--- end-fact-----------------------------------------
%-----------------begin proof------------------
\begin{proof}
Le \polcar $\rC {G}(e)=(T-e)^{|G|}$ est dans $\kT$
donc son \coe constant, qui est égal à $\pm e$ est dans~$\gk$.
\end{proof}
%-----------------end proof------------------

%:     Fact{factconnexe}-----------
\begin{fact}
\label{factconnexe}
Soit $(\gk,\gC,G)$  une \apG avec $\gk$ connexe et non trivial, alors:
%-----------------begin enum------------------
\begin{enumerate}
\item \label{factconnexe.1}\relax $0$ et $1$ sont les seuls \idms de $\gC$ fixés
par $G$,
\item \label{factconnexe.2}$\BB(\gC)$ est discrète,
\item \label{factconnexe.3}tout atome de $\BB(\gC)$ est un \idm galoisien,
\item \label{factconnexe.4}deux atomes sont conjugués sous $G$,
\item \label{factconnexe.5}un \idm $e\neq0$ est galoisien \ssi son orbite sous $G$ est formée d'\elts 2 à 2 \ortsz,
\item \label{factconnexe.6}si $f$ est un \idm $\neq0$, l'\id $\gen{1-f}$
est galoisien \ssi son orbite sous $G$ est formée d'\ids 2 à 2 \comz.
\end{enumerate}
%-----------------end enum------------------
\end{fact}
%--- end-fact-----------------------------------------
%-----------------begin proof------------------
\begin{proof}
Les points \emph{\ref{factconnexe.1}} et \emph{\ref{factconnexe.2}}
résultent clairement des faits
\emph{\ref{factIdemStable}} et~\emph{\ref{factagbdisc}.}

\emph{\ref{factconnexe.3}.} Si $e$ est un atome, $\sigma(e)$ aussi,
donc $\sigma(e) = e$ ou $e\sigma(e) = 0$.
Ainsi deux \elts de l'orbite de $e$ sont orthogonaux, donc la somme de
l'orbite de $e$ est un \idm non nul fixé par $G$: il est égal à~$1$.

\emph{\ref{factconnexe.4}.} Si $e'$ est un autre atome, il est égal
la somme des $e_ie'$, où $e_i$ parcourt l'orbite de $e$. Et comme les
$e_i$ sont des atomes, chacun des $e_ie'$ est nul ou égal à $e_i$.

\emph{\ref{factconnexe.5}.} Voir le fait \ref{factAduIdmA}.

\emph{\ref{factconnexe.6}.} Découle de \emph{\ref{factconnexe.5}}
puisque $\gen{1-f,1-f'}=\gen{1-ff'}$ pour des \idmsz~$f$ et~$f'$.
\end{proof}
%-----------------end proof------------------

Le  \thref{factagb}
implique que l'\agB $\BB(\gC)$  est
finie \ssi   les \idms in\dcps forment un ensemble fini (ils sont \ncrt 2 à 2
\ortsz) et s'ils engendrent $\BB(\gC)$.

\medskip
\comm \label{ensborn}
Un ensemble $X$ est dit \emph{borné}
 si l'on connaît
un entier $k$ qui majore le nombre de ses \eltsz. \Cad plus \prmt
si pour toute famille finie $(b_i)_{i\in\lrb{0..k}}$ dans $X$,
on a $b_i=b_j$ pour deux indices distincts.
En \clama ceci implique que l'ensemble est fini, mais du point de vue \cof
bien des situations distinctes peuvent se présenter.%
\index{borné!ensemble ---}%
\index{ensemble!borné}
\\
Une situation fréquente est celle  d'une \agB $C$ bornée
et discrète pour laquelle on ne connaît pas d'atome de manière s\^ure.
Les \itfs de $C$, tous principaux, s'identifient aux \elts de $C$, donc
 $C$ s'identifie à son propre treillis de Zariski{\footnote{Pour un anneau commutatif $\gk$, $\Zar\gk$
est l'ensemble des radicaux d'\itfs de $\gk$
(section \ref{secZarAcom}). C'est un treillis distributif.
En \clamaz, $\Zar\gk$ s'identifie au treillis des ouverts
quasi-compacts de l'espace spectral $\Spec\gk$ (section~\ref{secEspSpectraux}).}}
$\Zar C$.
Par ailleurs, en \clama les atomes sont en bijection avec les \ids premiers
(tous maximaux) de $C$ via $e\mapsto \gen{1-e}$. Ainsi l'ensemble des atomes
de $C$ (supposé borné) s'identifie à $\Spec C$. On retrouve donc
dans ce cas particulier le fait \gnl suivant: le treillis de Zariski
est la version \covz, maniable et \gui{sans point} du spectre de Zariski,
espace topologique dont les points peuvent s'avérer inaccessibles d'un
point de vue \cofz. Mais cette situation, bien que familière, est peut
être plus troublante dans le cas d'un espace topologique discret et
borné. Il s'agit typiquement d'un espace compact dont on n'a pas une bonne
description via un sous-ensemble énumérable dense, donc qui n'entre pas
dans le cadre des espaces métriques compacts à la Bishop (cf.~\cite{B67,BB85}).
\eoe

\medskip
Voici un corolaire du \tho de structure galoisien \ref{lemAduIdmA} dans le contexte des \apGsz.

%--- Proposition{corlemAduIdmA}----------
\begin{proposition}
\label{corlemAduIdmA}
Soit  $(\gk,\gC,G)$  une \apG  avec $\gk$ connexe. Pour un \idm $h$ de $\gC$
\propeq
%-----------------begin enum------------------
\begin{enumerate}
\item $h$ est un \idm galoisien.
\item $\gC[{1/h}]$ est un \kmo \pro de rang égal à $\St_G(h)$.
\item $\gC[{1/h}]$ est un \kmo \pro de rang inférieur ou égal à
$\St_G(h)$.
\end{enumerate}
%-----------------end enum------------------
\end{proposition}
%--- end-corollary------------------------------------
%-----------------begin proof------------------
\begin{proof}
On utilise le \thrf{lemAduIdmA}. D'après le point \emph{\ref{lemAduIdmA-2}}
de ce \tho on peut supposer qu'il existe un \idm galoisien $e$ tel que $h$
soit égal à une somme $e_1+\cdots +e_r$ d'\elts de l'orbite $G.e$.  On a des \isos de
\kmos $e\gC\simeq \gC[{1/e}]$ et $\gC\simeq (e\gC)^{\abs{G.e}}$, donc
$e\gC$ est \pro de rang constant $\idg{G:G.e} = \abs{\St_G(e)}$. On en déduit
que  le \kmo 
$$
\gC[{1/h}]\simeq h\gC = e_1\gC\oplus\cdots\oplus e_r\gC\simeq (e\gC)^r
$$
est \pro de rang $r\times \abs{\St_G(e)}$.  On applique alors le point
\emph{\ref{lemAduIdmA-5}} du \thref{lemAduIdmA} avec
$f = h$.
\\
Donc, le point \emph{2} (resp.\ le point \emph{3}) ici signifie la même
chose que le point~\emph{5c} (resp.\ le point~\emph{5d}) dans le \thref{lemAduIdmA}.
\end{proof}
%-----------------end proof------------------

%:--- SUBsection{Discriminant}-----
\subsec{Discriminant%, \mdiz
}
\label{subsecDiscrim}

Rappelons que dans $\gA=\Adu_{\gk,f}$ on a
$\disc(f)= \prod\nolimits_{1\leq i<j\leq n}(x_i-x_j)^2$
et~$\Disc\iAk =\disc(f)^{n!/2}$.

Dans le \tho suivant, on parle du $\gA$-\mdi $\Om{\gk}{\gA}$ de la \klg $\gA$.
Il suffit en fait de savoir
que le \mdi d'une \apf est isomorphe au conoyau de la
transposée de la matrice jacobienne
du \syp qui définit l'\algz. 
Pour plus de précisions sur ce sujet voir les \thosz~\ref{thDerivUniv} et \ref{thDerivUnivPF}.

%:     Theorem{lemAdu1}
\begin{theorem}\label{lemAdu1}
 Soit $J$ le jacobien du système de $n$ équations à $n$
inconnues définissant l'\adu $\gA=\Adu_{\gk,f}$.
%-----------------begin enum------------------
\begin{enumerate}
\item
\begin{enumerate}
\item \label{i1lemAdu1}
On a $J=\prod_{1\leq i<j\leq n}(x_i-x_j)$ dans $\gA$.
\item \label{i2lemAdu1}
 On a $J^2=\disc(f)\in\gk$.
\end{enumerate}
\item En particulier, \propeq
%-----------------begin enum------------------
\begin{enumerate}\itemsep0pt
\item \label{i3alemAdu1}
  $\Disc\iAk $ est \iv (resp.\,\ndzz) dans  $\gk$.
\item  $\disc(f)$ est \iv (resp.\,\ndzz) dans  $\gk$.
\item  $J$ est \iv (resp. \ndzz) dans  $\gA$.
\item Les $x_i-x_j$ sont \ivs (resp.\,\ndzsz) dans $\gA$.
\item \label{i3elemAdu1}
 $x_1-x_2$ est \iv  (resp.\,\ndzz) dans $\gA$.
\item \label{i3flemAdu1}
 $\Om{\gk}{\gA}=0$ (resp.\,$\Om{\gk}{\gA}$ est un \Amo \gui{de torsion},
i.e. annulé par un \elt \ndzz).
\item \label{i4lemAdu1}
 $\Sn$ est un groupe séparant
pour $\gA$ (resp.\,pour  $\Adu_{\Frac(\gk),f}$).
\end{enumerate}
%-----------------end enum------------------

%
\item \label{i5lemAdu1} Les \eqvcs analogues sont valables pour tout quotient de Galois
de l'\aduz.

\end{enumerate}
%-----------------end enum------------------
 \end{theorem}
%--- end-theorem
%-----------------begin proof------------------
\begin{proof}
Le point \emph{1a.}
est facile par \recu sur $n$, avec le signe exact si l'on considére le
\sys qui nous a servi pour la \dfn de l'\aduz.
Voici par exemple le calcul pour $n=4$
\snucc{\arraycolsep2pt
%--------------------begin array---------------
\begin{array}{rclcl}
J&=   &
\dmatrix{1&1&1&1\cr
\som_{i\neq 1}x_i&\som_{i\neq 2}x_i & \som_{i\neq 3}x_i&
\som_{i\neq 4}x_i  \cr
\som_{i,j\neq 1}x_ix_j&\som_{i,j\neq 2}x_ix_j & \som_{i,j\neq 3}x_ix_j&
\som_{i,j\neq 4}x_ix_j  \cr
x_2x_3x_4&x_1x_3x_4 &  x_1x_2x_4 &x_1x_2x_3  }
\\[8mm]
%\end{array}
%\begin{array}{rclcl}
\phantom{J}&  = &
\dmatrix{1&0&0&0\\[.5mm]
{\som_{i\neq 1}x_i} & x_1-x_2 & x_1-x_3&x_1-x_4  \cr
{\som_{i,j\neq 1}x_ix_j} & (x_1-x_2)\dsp{\sum_{i\neq 1,2}x_i} &
(x_1-x_3)\dsp{\sum_{i\neq 1,3}x_i} & (x_1-x_4)\dsp{\sum_{i\neq 1,4}x_i}
\cr
x_2x_3x_4&(x_1-x_2) x_3x_4 &  (x_1-x_3)x_2x_4 &(x_1-x_4)x_2x_3  }
   \\[10mm]
%\end{array}
%\begin{array}{rclcl}
\phantom{J}& = &(x_1-x_2)(x_1-x_3)(x_1-x_4)
\dmatrix{
 1 & 1&1  \cr
x_3+x_4 &
x_2+x_4 & x_2+x_3  \cr
 x_3x_4 &  x_2x_4 &x_2x_3  }
\end{array}
%---------------------end array--------------
}%
%\sni
etc\ldots

 On en déduit le point \emph{1b}, puis l'\eqvc des points \emph{2a} à \emph{2e}.

\emph{2f.} Puisque $\Om{\gk}{\gA}$ est un \Amo isomorphe au
conoyau de la transposée de la matrice jacobienne,
on obtient que $\Ann(\Om{\gk}{\gA})$
et $J\gA$ ont même nilradical (lemme \ref{fact.idf.ann}). Enfin l'\elt $J$ est \ndz (resp. \ivz) \ssi
l'\id $\sqrt{J\gA}$ contient un \elt \ndz (resp. contient~$1$).

\emph{2g.} Supposons $f$ \spl (resp.\,\ndzz), si $\sigma\in \Sn$ est distinct de $\Id_\gA$, il y a un $i\in\lrbn$ tel que $x_{\sigma i}\neq x_i$.
Puisque $x_{\sigma i}- x_i$ est \iv (resp.\,\ndzz),~$\sigma$ est  séparant (resp.\,séparant
une fois que l'on  inverse le \discriz).
Pour la réciproque, considérons par exemple la transposition~$\sigma$
qui échange~$1$ et~$2$. On a clairement $\gen{g-\sigma(g)\vert g\in\gA}=\gen{x_1-x_2}$. Ceci permet de conclure.

\emph{3.} Clair puisque l'\adu est toujours isomorphe à une puissance de
n'importe lequel de ses
quotients de Galois.
\end{proof}
%-----------------end proof------------------

%%%%%%%%%%%%%%%%%%%%%%%%%%%%%%%%%%%%%%%%%%%%%%%%%%%%%%%%%%%%%%%%%%%%%%%%%%%
%:--- SUBsection{points fixes}-----
\subsec{Points fixes}
\label{subsecptsfix}

Nous notons $\di(f)=\prod_{i<j\in\lrbn}(x_i+x_j)\in\gk$.

Il est clair que $\di(f)$ est congru modulo 2 à $\prod_{i<j\in\lrbn}(x_i-x_j)$,
ce qui donne $\gen{2,\di(f)^2}=\gen{2,\disc(f)}$.

%:    Theorem{theoremAdu1}            points fixes
\begin{theorem} \label{theoremAdu1} \emph{(Algèbre de \dcn \uvle et points fixes)}\\
Posons $\fa\eqdefi \Ann_{\gk}(\gen{2,\di(f)})$.
Alors: 

\snic{\Fix(\Sn)\subseteq\gk+\fa\gA.}

%\sni
En particulier, si $\fa=0$ et a fortiori si
$\Ann_{\gk}(\gen{2,\disc(f)})=0$,  on obtient 

\snic{\Fix(\Sn)=\gk.}
\end{theorem}
%--- end-theorem-----------------------------------------
%-----------------begin proof------------------
\begin{proof}
Il suffit  de démontrer la
première affirmation. \\
Voyons le cas où  $n=2$ avec $f(T) = T^2 - s_1T + s_2$. 
\\
Un \elt $z=c+dx_1\in\gA$ (avec $c$, $d\in\gk$) est
invariant par $\mathrm{S}_2$  \ssi $d(x_1-x_2)=d(s_1-2x_1)=0$, ou encore si
$ds_1=2d=0$, or on~a~$\di(f)=s_1$.
\\
On procède ensuite par récurrence sur $n$.
On reprend pour les modules de Cauchy les notations de la section~\ref{sec0adu}.
Pour $n>2$ on considère l'anneau~$\gk_1=\gk[x_1]\simeq\aqo{\gk[X_1]}{f(X_1)}$ et
le \pol $g_2(T)=f_2(x_1,T) $ qui est dans $\gk_1[T]$. On identifie  $\Adu_{\gk_1,g_2}$
avec $\Adu_{\gk,f}$ (fait~\ref{factAduAdu}).
Pour passer de l'écriture d'un \elt $y\in\gA$ sur la base $\cB(g_2)$ ($\gA$
vu comme $\gk_1$-module) à son écriture sur la base $\cB(f)$ ($\gA$ vu
comme \kmoz), il suffit d'écrire chaque coordonnée, qui est un \elt de
$\gk_1$, sur la $\gk$-base $(1,x_1,\ldots,x_1^{n-1})$ de $\gk_1$.
%\\
Notons aussi que $\di(f)=(-1)^{n-1}g_2(-x_1)\di(g_2)$ par un calcul direct.
Donc, si nous posons $\fa_1=\Ann_{\gk_1}(\gen{2,\di(g_2)})$,
nous obtenons $\fa_1\gA\subseteq \fa\gA$ et~$\fa_1\subseteq \fa\gk_1$.
\\
Passons à la \recu proprement dite.
\\
Soit $y\in \gA$ un  point fixe de $\Sn$, et regardons le comme un \elt de l'\adu  $\Adu_{\gk_1,g_2}$.
Puisque $y$ est invariant par~$\mathrm{S}_{n-1}$, on a %par \hdr que 
$y\in\gk_1+\fa_1\gA$, et donc $y\equiv h(x_1) \mod \fa_1\gA$
pour un $h\in\kX$.
A fortiori $y\equiv h(x_1) \mod \fa\gA$.
Il reste à voir que $h(x_1)\in\gk+\fa\gA$.
Puisque~$y$ est invariant par~$\Sn$, on obtient en permutant~$x_1$ et~$x_2$ la congruence

\snic{\qquad\qquad\qquad h(x_1)\equiv y \equiv h(x_2)
\quad\mod \fa\gA\qquad\qquad\qquad(*)}

%\sni
\'Ecrivons $h=\sum_{i=0}^{n-1}+c_{i}X^{i}\in\kX$.
On note que $h(x_1)$ est une écriture réduite sur la base canonique
$\cB(f)$. Concernant $h(x_2)$, pour obtenir l'écriture réduite,
nous devons remplacer dans le terme $c_{n-1}x_2^{n-1}$, $x_2^{n-1}$ par son
écriture sur la base canonique, qui résulte de $f_2(x_1,x_2)=0$.
\\
Cette réécriture fait apparaître le terme $-c_{n-1}x_1^{n-2}x_2$,
et ceci implique d'après~$(*)$
que~$c_{n-1}\in \fa$.
Mais alors, $h(x_2)-c_{n-1}x_2^{n-1}$ et $h(x_1)-c_{n-1}x_1^{n-1}$
sont des écritures
réduites de deux \elts égaux modulo $\fa\gA$.
Donc, les $c_i$ pour $i\in\lrb{1..n-2}$ sont dans $\fa$,
et l'on a vu que $c_{n-1}\in\fa$.
\end{proof}
%-----------------end proof------------------

\rem
Dans le cas $n=2$, l'étude faite ci-dessus 
montre que dès que~$\fa\neq 0$, l'anneau~$\Fix(\mathrm{S}_2)= \gk\oplus
\fa\,x_1=\gk+\fa\gA$ contient strictement $\gk$.
\\
Un calcul dans le cas~$n=3$ donne la même réciproque: si $\fa\neq0$,
l'anneau~$\Fix(\mathrm{S}_3)$ contient strictement~$\gk$. On trouve en effet
un \elt

\snic{v=x_1^2 x_2+  s_1 x_1^2 + (s_1^2+s_2) x_1 + s_2 x_2\neq 0}

%\sni
(une de ses
coordonnées sur $\cB(f)$ est égale à $1$) tel que~$\Fix(\rS_3)=
\gk\oplus \fa\,v$.
Par contre pour~$n\geq 4$, la situation se complique.
\eoe

\medskip
On obtient comme corolaire le \tho suivant.

%:     Theorem{thAduAGB}
\begin{theorem}\label{thAduAGB}
Si $f$ est un \pol \spl de $\kT$, l'\adu $\Adu_{\gk,f}$,
ainsi que tout quotient de Galois, est une \aGz.
\end{theorem}

\begin{proof}
D'après le \tho de structure \ref{thADG1Idm} (point \emph{7}) il suffit de montrer que $\Adu_{\gk,f}$
est galoisienne. Or on vient de démontrer la condition sur les points fixes,
et la condition sur les \autos séparants a été donnée dans le
\thrf{lemAdu1}.
\end{proof}

D'après le \tho d'Artin \ref{thA}, et vu le \tho précédent,
nous savons que toute \adu pour un \pol \splz, ou tout quotient de Galois d'une telle \klgz, se diagonalise elle-même. 
Nous examinons cette question plus en détail dans le paragraphe qui suit.
Même pour ce qui concerne le résultat précis que nous venons de citer,
il est intéressant de voir fonctionner la chose de façon \gui{concrète}
pour une \aduz.

%:--- SUBsection{Separabilite}-----
\subsec{Séparabilité}
\label{subsecsep}

Lorsque le \pol $f\in \gk[T]$ est \splz,  son \adu $\gA = \Adu_{\gk,f} = \gk[\xn]$ est \stez, d'après le fait~\ref{factDiscriAdu}. Le
\tho suivant est alors un simple rappel du \thref{prop2EtaleReduit} concernant
les \ases dans le cadre présent.

%:    Theorem{theoremAdu2}------
\begin{theorem}
\label{theoremAdu2} On suppose $f$ \splz.
%-----------------begin enum------------------
\begin{enumerate}
\item
Le nilradical $\rD_\gA(0)$ est l'\id engendré par $\rD_\gk(0)$. En particulier, si
$\gk$ est réduite, $\gA$ est réduite.
\item 
Pour toute \alg réduite $\gk\vers{\rho}\gk'$,
l'\alg $\rho\ist(\gA)\simeq\Adu_{\gk',\rho(f)}$ est réduite.
\end{enumerate}
%-----------------end enum------------------
\end{theorem}
%--- end-theorem-----------------------------------------
\subsubsection*{Diagonalisation d'une \aduz}

%:    Theorem{theoremAdu3}- \din d'une adu
\begin{theorem}
\label{theoremAdu3} \emph{(Diagonalisation  d'une \aduz)} 
Soit $\varphi:\gk\to\gC$ une \alg dans laquelle \emph{$f$ se factorise
complètement}, \cad $\varphi(f)=\prod_{i=1}^n(T-u_i)$.
Supposons aussi que $f$ est \spl sur~$\gC$, i.e. que les $u_i-u_j$
sont \ivs  pour $i \ne j$.\\
Notons $\gC\te_\gk\gA\simeq\Adu_{\gC,\varphi(f)}$,
et, pour $\sigma\in\Sn$,  $\phi_\sigma: \gC\te_\gk\gA \to\gC$ l'unique
\homo de \Clgs qui envoie chaque $1_\gC\otimes x_i$ sur $u_{\sigma i}$.
\\
Soit $\Phi: \gC\te_\gk\gA\to\gC^{n!}$ le $\gC$-\homo défini par
$y\mapsto\big(\phi_\sigma(y)\big)_{\sigma\in\Sn}$.
%-----------------begin enum------------------
\begin{enumerate}
\item 
$\Phi$ est un \isoz: $\gC$ diagonalise $\gA.$
\item 
Plus \prmtz, dans $\gC\te_\gk\gA$, notons $x_i$ à la place de
$1_\gC\te x_i$, $u_i$ à la place de $u_i\te 1_{\gA}$ (conformément
à la structure  
de \Clg de  $\gC\te_\gk\gA$) et posons 
$g_\sigma=\prod_{j\neq \sigma i} (x_i-u_j)$.  
Alors, 

\snic{\phi_\sigma(g_\sigma)=\pm \varphi\big(\disc(f)\big)=\pm\disc\big(\varphi(f)\big),
}

%\sni

et $\phi_\sigma(g_\tau)=0$ pour $\tau\neq \sigma$, de sorte que
si l'on pose $e_\sigma=g_\sigma/\phi_\sigma(g_\sigma)$, les~$e_\sigma$
forment le \sfio correspondant à l'\isoz~$\Phi$.
\item En outre, $x_ie_\sigma=u_{\sigma i}e_\sigma$, de sorte que la
base $(e_\sigma)$ du \Cmo $\gC\te_\gk\gA$ est une base diagonale commune pour
les multiplications par les~$x_i$.
\end{enumerate}
%-----------------end enum------------------
En particulier, lorsque $f$ est \splz, l'\alg enveloppante

\snic{\env\gk\gA=\gA\otimes_\gk\gA\simeq\Adu_{\gA,f}}

%\sni
est isomorphe
canoniquement à $\Ae {n!}$: $\gA$ se diagonalise elle-même.
\end{theorem}
%--- end-theorem-----------------------------------------
%-----------------begin proof------------------
NB: on prendra garde cependant à noter
$\Adu_{\gA,f}=\gA[u_1,\ldots,u_n]$
puisque les $x_i$ sont déjà pris comme \elts de~$\gA$.
\begin{proof}
\emph{ 1.}
Les deux \algs sont en tant que \Cmos isomorphes à $\gC^{n!}$
et $\Phi$ est une \Cli dont il suffit de démontrer la surjectivité.
La surjectivité résulte par le \tho chinois de ce que les
$\Ker\phi_\sigma$ sont deux à deux \comz: $\Ker\phi_\sigma$ contient
$x_i-u_{\sigma i}$,  $\Ker\phi_\tau$ contient
$x_i-u_{\tau i}$, 
\linebreak 
donc  $\Ker\phi_\sigma+\Ker\phi_\tau$ contient les
$u_{\sigma i}-u_{\tau i}$, et il y  a au moins un indice $i$ pour lequel
$\sigma i\neq \tau i$, ce qui donne $u_{\sigma i}-u_{\tau i}$ \ivz.

\emph{2.} Le \sfio correspondant à
l'\iso $\Phi$ est l'unique solution du \sli
$\phi_\sigma(e_\tau)=\delta_{\sigma,\tau}$ (où $\delta$ est le symbole de
Kronecker). \\
Or les \egts $\phi_\sigma(g_\sigma)=\pm \varphi\big(\disc(f)\big)$ et
$\phi_\sigma(g_\tau)=0$ sont faciles. 

\emph{3.} Fixons $i$. 
L'\egt $x_ig_\sigma=u_{\sigma i}g_\sigma$ résulte de ce que
dans $g_\sigma$ il y a déjà le produit des $x_i-u_j$ pour $j\neq \sigma i$,
donc $(x_i-u_{\sigma i})g_\sigma$ est multiple de $\varphi(f)(x_i)$, qui est nul. 
\end{proof}
%-----------------end proof------------------

\rem 
En fait, de manière \gnlez, $\Phi$ est une \ali dont on peut calculer le
\deter par rapport aux bases naturelles:  le carré de ce \deter est une
puissance de $\varphi\big(\disc(f)\big)$ 
et l'on trouve ainsi que~$\Phi$ est un \iso \ssi
$\varphi\big(\disc(f)\big)$ est \iv dans~$\gC$. Pour ceci, et pour une \gui{réciproque complète}, voir l'exercice~\ref{exoIdentiteDiscriminantale}.
\eoe

\perso{il ne semble pas que l'on  arrive à faire découler le
\thrf{theoremAdu2} du \ref{theoremAdu3}.}

\medskip Le \tho précédent implique le résultat suivant: si $\gA$
est une \adu pour un \pol \spl toute \Alg diagonalise~$\gA$.
Nous en donnons maintenant une \gnn pour un quotient de Galois de $\gA$.

%:    Theorem{theoremAdu4} \din d'un quotient de Galois d'une adu
\begin{theorem}
\label{theoremAdu4} \emph{(Diagonalisation d'un quotient de Galois d'une
\aduz)}
Soient $e$ un \idm galoisien de $\gA$, 

\snic{\gB=\aqo{\gA}{1-e}=\gk[\yn]\;\hbox{ et }\;G=\St_\Sn(e),}

%\sni
(on a noté $y_i=\pi(x_i)$ la classe de $x_i$ dans $\gB$).
Soit $\phi:\gB\to\gC$ un \homo d'anneaux. On note $u_i=\phi(y_i)$.
On considère la \Clg  

%\snic{\gC\otimes_\gk\gB\simeq\aqo{\Adu_{\gC,f}}{1-\phi(e)}}
\snic{\phi\ist(\gB)\simeq\gC\otimes_\gk\gB\simeq\aqo{\Adu_{\gC,f}}{1-\phi(e)}}

%\sni
obtenue
à partir de la \klg $\gB$ par \edsz.
Pour~$\sigma\in G$ notons $\phi_\sigma:\gC\otimes_\gk\gB\to\gC$ l'unique
\homo de \Clgs qui envoie chaque $1_\gC\otimes y_i$ sur $u_{\sigma i}$.
Soit $\Phi:\gC\otimes_\gk\gB \to\gC^{|G|}$ l'\homo de~\Clgs défini par
$z\mapsto\big(\phi_\sigma(z)\big)_{\sigma\in G}$.
%-----------------begin enum------------------
\begin{enumerate}
\item Si $\phi\big(\disc(f)\big)\in\gC\eti$, $\Phi$ est un \isoz, donc $\gC$ diagonalise $\gB$.
\item En particulier, si $f$ est \splz, $\gB\otimes_\gk\gB$ est isomorphe canoniquement à~$\gB^{|G|}$, i.e. $\gB$ se diagonalise elle-même.
\end{enumerate}
%-----------------end enum------------------
\end{theorem}
%--- end-theorem-----------------------------------------
%-----------------begin proof------------------
\begin{proof}
Les deux $\gC$-\algs sont des \Cmos \pros de rang constant  $|G|$
et $\Phi$ est une \Cli dont il suffit de démontrer la surjectivité.
Dans $\gC\otimes_\gk\gB$ nous notons $y_i$ à la place de  $1_\gC\otimes y_i$ et  $u_i$
à la place de  $u_i\otimes 1_{\gB}$.
La surjectivité résulte par le \tho chinois de ce que les
$\Ker\phi_\sigma$ sont deux à deux \comz: $\Ker\phi_\sigma$ contient
$y_i-u_{\sigma i}$,  $\Ker\phi_\tau$ contient
$y_i-u_{\tau i}$, donc  $\Ker\phi_\sigma+\Ker\phi_\tau$ contient les
$u_{\sigma i}-u_{\tau i}$. Or il y  a au moins un indice $i$ pour lequel
$\sigma i\neq \tau i$ et $u_{\sigma i}-u_{\tau i}$ est \iv parce que
$\phi\big(\disc(f)\big)$ est le produit des $(u_j-u_k)^2$ pour $1\leq j<k\leq n$.
\end{proof}
%-----------------end proof------------------

%:--- SUBsection{Structure triang}--subsecidGTri---
\subsec{Structure triangulaire des idéaux galoisiens}
\label{subsecidGTri}

Nous démontrons dans ce paragraphe le \thrf{thidGTri} qui  implique que la structure de l'\id $\cJ(f)$, qui est une
structure \gui{triangulaire} (au sens de Lazard) lorsque l'on considère les
modules de Cauchy comme \gtrsz, reste une structure triangulaire pour tous les
\ids galoisiens de l'\adu dans le cas d'un \pol \spl sur un \cdiz.

%:     Lemma{lemBaseDisc}
\begin{lemma}\label{lemBaseDisc}
Soit $\gk'$ une \klg qui est un \mptf de rang constant $m$,
$x \in \gk'$ et $r(T) \in \gk[T]$ le \polcar de~$x$ sur~$\gk$. 
Si $\disc(r) \in \gk\eti$, alors $\gk' = \gk[x]$
et $(1, x, \ldots, x^{m-1})$ est \hbox{une $\gk$}-base de $\gk'$. 
\end{lemma}
%--------- fin lemma ---------------------------------------------- 

\begin{proof}
Le cas où $\gk'$ est libre de rang $m$ a été prouvé en
\ref {propdiscTra}. Dans le cas \gnlz, on considère un \sys d'\eco de $\gk$
tel que chaque \lon fasse de $\gk'$ un \kmo libre de rang $m$.
\end{proof}

%:    Theorem{thidGTri}------------
\begin{theorem}
\label{thidGTri}
Soit  $(\gk,\gC,G)$ une \aG avec:
%-----------------begin item------------------
\begin{itemize}
\item   $\gC=\kxn\simeq\kXn\sur\fa$,
\item  $G$ opère sur $\{\xn\}$ et
\item  les $x_i-x_j$ sont \ivs pour $i\neq j$.
\end{itemize}
%-----------------end item------------------
Un exemple typique de cette situation: $\gC$ est un quotient de Galois de l'\adu
d'un \pol \splz.
\\
On pose

\snic{
\begin{array}{l} 
G_i=\sotQ{\sigma\in G}{\sigma(x_k)=x_k,\, k\in \lrb{1..i}}
\hbox{ pour }i \in \lrb{0..n} \hbox{ (donc } G_0 = G \hbox{)}, \\[1mm]
 r_i(T)=\prod\nolimits_{\sigma\in G_{i-1}/G_i}\big(T-\sigma(x_i)\big)
\quad  \hbox{ pour }  i \in\lrbn,
\end{array}
}

%\sni
où $G_{i-1}/G_i$ désigne un \sys de représentants des classes à
gauche. On note $d_i=\idg{G_{i-1}:G_i}$.
\\
On a alors les résultats suivants.
%-----------------begin item------------------
\begin{enumerate}
\item  $\gk[x_1,\ldots ,x_i]=\Fix(G_i)$ et $G_i=\Stp(\gk[x_1,\ldots ,x_i]).$
\item  Le \pol $r_i(T)$ est \mon à \coes
dans $\gk[x_1, \ldots, x_{i-1}]$,  de degré $d_i$.
On note $R_i(X_1,\ldots ,X_i) \in \gk[X_1,\ldots ,X_i]$
un \polu en $X_i$ de degré $d_i$
tel que $R_i(x_1,\ldots ,x_{i-1},X_i)=r_i(X_i)$.
\item  %L'idéal
$\fa_i=\fa\cap\gk[X_1,\ldots ,X_i]$ est engendré
par $R_1(X_1)$, $\ldots$, $R_i(X_1,\ldots ,X_i)$.
\end{enumerate}
%-----------------end item------------------
En conséquence chaque \alg $\gk[x_1,\ldots ,x_i]$ est à la fois un
$\gk[x_1,\ldots ,x_{i-1}]$-module libre de rang $d_i$ et un \kmo
libre de rang $\idg{G:G_i}$, et chacun des \ids $\fa_i$ est un \id
triangulaire (au sens de Lazard) de $\gk[X_1,\ldots ,X_i]$.
\end{theorem}
%--- end-theorem-----------------------------------------
%-----------------begin proof------------------
\begin{proof}
Le groupe $G_1$ est un groupe séparant d'\autos de l'anneau~$\gC$.  On note
$\gk_1=\gC^{G_1}$.  On sait que $\gC$ est un
$\gk_1$-module \pro de rang constant $|G_1|$ et que
$\gk[x_1]\subseteq\gk_1$. En outre, $\gk_1$ est facteur direct dans~$\gC$, donc
est un \kmo \pro de rang constant $d_1=\deg_T(r_1).$
\\
L'idéal $\fa_1$ est formé par tous les $R\in\gk[X_1]$ tels que $R(x_1) =
0$.\\
 Donc, $R\big(\sigma(x_1)\big)=0$ pour tout $\sigma\in G/G_1$, autrement dit $R$ est multiple de
chaque~$T- \sigma(x_1)$. Et puisque les $x_i-x_j$ sont \ivsz, $R$ est multiple de~$r_1$.  Ainsi
$\fa_1=\gen{r_1(X_1)}$ et $\gk[x_1] \simeq \aqo{\gk[X_1]}{r_1(X_1)}$.
\\
La proposition \ref{propAdiagAH} donne l'\egt

\snic{\rC{\gk_1/\gk}(x_1)(T) = \prod_{\sigma \in G/G_1}\big(T - \sigma(x_1)\big)=r_1(T).}

%\sni
Ceci implique que le \polcar $\rC{\gk_1/\gk}(x_1)(T)$ est \splz, et le lemme
\ref{lemBaseDisc} dit que $(1,x_1,\dots,x_1^{d_1-1})$ est une base de $\gk_1$.
\\
Ainsi  $\gk[x_1]=\gk_1=\Fix(G_1)$ et $(\gk[x_1],\gC,G_1)$ est une
\aGz. 

 Alors, $\gC=\gk_1[x_2,\ldots ,x_n]$
avec $G_1$ qui opère sur $\{x_2,\ldots ,x_n\}$ et les $x_i-x_j$
\ivsz.
Tout le raisonnement précédent fonctionne à l'identique en
remplaçant $\gk$ par $\gk_1$, $G$ par $G_1$, $x_1$ par $x_2$
et $G_1$ par $G_2$. 
On termine donc par \recuz.
 \end{proof}
%-----------------end proof------------------

%%%%%%%%%%%%%%%%%%%%%%%%%%%%%%%%%%%%%%%%%%%%%%%%%%%%%%%%%%%%%%%%%%%%%%%%%%%
%%%%%%%%%%%%%%%                                         %%%%%%%%%%%%%%%%%%%
%%%%%%%%%%%%%%%           Corps des racines             %%%%%%%%%%%%%%%%%%%
%%%%%%%%%%%%%%%                                         %%%%%%%%%%%%%%%%%%%
%%%%%%%%%%%%%%%%%%%%%%%%%%%%%%%%%%%%%%%%%%%%%%%%%%%%%%%%%%%%%%%%%%%%%%%%%%%
%\newpage
%--- Section{Corps des racines}-----
\section{Corps de racines d'un \pol sur un corps discret}
\label{subsecCDR}

Nous donnons dans cette section une approche \cov et dynamique du \cdr d'un \polu sur un \cdiz, en l'absence d'\algo de
factorisation des \polsz.

\Grandcadre{
Dans la section \ref{subsecCDR},  $\gK$ est un \cdi non trivial,
\\  $f$ est un \pol
 unitaire de degré $n$ et\\
  $\gA=\Adu_{\gK,f}=\KXn\sur{\cJ(f)}=\Kxn$.
}

\smallskip Les quotients de l'\adu $\gA$ sont des \Klgs finies,
donc ce sont des anneaux \zedsz.

%:--- SUBsection{Quotients réduits}-
\subsect{Quotients de Galois \gui{réduits}
de  l'\alg  de \\ \dcn \uvlez}{Bons quotients de l'\aduz}

Nous avons mis des guillemets à \gui{réduits} parce que, a priori,
on ne parle pas d'un quotient de Galois \emph{qui est} réduit, mais d'un quotient de Galois \emph{que l'on} réduit (en tuant les nilpotents).

Vu le fait \ref{factDiscriAdu} si le \pol $f$ est \spl l'\adu est
étale, donc réduite, et tout \id engendré par un \idm est égal à son
nilradical (puisque l'anneau quotient est réduit).
On peut alors remplacer dans les énoncés qui suivent chaque idéal
$\DA(1-e)=\sqrt{\gen{1-e}}$
par l'\id $\gen{1-e}$.

\smallskip 
Dans le lemme qui suit on sait par hypothèse que $\gB$ est \stfe sur~$\gK$,
mais on ne connaît pas \ncrt une base de~$\gB\red$ comme~\Kevz.
Le but est alors de donner une description \gui{assez satisfaisante} de~$\gB\red$
comme quotient de l'\aduz.

%--- Lemma{lemAquQuoRed}----------
\begin{lemma}
\label{lemAquQuoRed}
Soit $\gB$ une \Klg \stfez. On suppose que $f$ se décompose totalement dans
$\gB\red$ et que $\gB\red$ est engendrée par les zéros correspondants de
$f$.  Alors, il existe un \idm $e$ de $\gA=\Adu_{\gK,f}$ tel que
$\gB\red\simeq\gA\sur{\DA(1-e)}$.
\end{lemma}
%--- end-lemma-----------------------------------------
%-----------------begin proof------------------
\begin{proof}
Soient $y_1$, \dots, $y_n\in\gB$ tels que $f(T)=\prod_i(T-\ov{y_i})$ dans $\gB\red$.  Il existe un unique \homoz~$\lambda :\KXn\to\gB$ qui envoie les~$X_i$ sur les~$y_i$. Notons $\fb$ l'\id (\tfz) de $\gB$ engendré par $\lambda\big(\cJ(f)\big)$. On
a alors~$\fb\subseteq \DB({0})$, et $\gB' :=\gB\sur{\fb}$ est une~\Klg
strictement finie vérifiant~$\gB\red\simeq \gB'\red$. On obtient
ainsi un diagramme 

\snic {
\xymatrix @R = 0.4cm {
\KuX\ar[d]_{\lambda}\ar@{->>}[r] &\gA\ar[d]^{\varphi}\ar@{->>}[dr]^{\psi} 
\\
\gB \ar[r]         &\gB'\ar[r]   &\;\gB'\red = \gB\red 
\\
}}

%\sni
dans lequel $\varphi$ est l'unique \homo qui envoie $x_i$ sur la classe
de $y_i$. Puisque $\gB'$ est \stfez, $\Ker\varphi$ est un \itf de $\gA$, et il
existe $d \ge 0$ tel que $(\Ker\varphi)^d = (\Ker\varphi)^{d+1}$ donc\perso{Il serait intéressant d'obtenir
une borne raisonnable pour l'exposant $d$. 
}
$(\Ker\varphi)^d$ est engendré par un \idm $1-e$. Et l'on peut conclure
car d'une part, $\psi$ est surjectif, et d'autre part,
$\Ker\psi = \DA(\Ker \varphi) = \DA\big((\Ker \varphi)^d\big) = \DA(1-e)$.
\end{proof}
%-----------------end proof------------------

\rem
Notez que $\psi$ est surjectif, mais a priori $\Ker\psi$
n'est pas un \itf de $\gA$. Symétriquement, a priori $\varphi$ n'est pas
surjectif, mais $\Ker\varphi$ est un \itf de $\gA$.
\eoe

\medskip 
En \clama un \cdr (\dfn \ref{defCorpsdesRacines}) 
pour un \pol unitaire $f$ sur un corps
discret $\gK$ est obtenu comme quotient de l'\adu $\gA$ par un
\idemaz. Un tel \id existe: prendre un \id strict qui soit un~\Kev de dimension maximale,
d'après le principe du tiers exclu. 

En \coma on obtient le \tho plus précis qui suit.

%--- Theorem{thUnici}------
\begin{theorem}
\label{thUnici} ~
%-----------------begin enum------------------
\begin{enumerate}
\item \label{i1thUnici}
\Propeq
%----- begin enum ----
\begin{enumerate}
\item \label{i1athUnici} Il existe dans $\gA=\Adu_{\gK,f}$  un \idm
in\dcpz~$e$.
\item \label{i1bthUnici} Il existe une extension $\gL$ de $\gK$ qui est un \cdr de
$f$ et qui s'écrit $\gB\red$ où $\gB$ est une \Klg \stfez.
\item \label{i1cthUnici} L'\agB $\BB(\gA)$ est finie.
\end{enumerate}
%-----------------end enum------------------
\item \label{i2thUnici} Dans ce cas tout \cdr de $f$ est isomorphe à
$\gA\sur{\DA(1-e)}$, et il est  discret.
\end{enumerate}
%-----------------end enum------------------
\end{theorem}
%--- end-theorem-------------------
%-----------------begin proof------------------
\begin{proof}
L'\eqvc de \emph{1a} et \emph{1c} vaut dans le cadre \gnl des \agBs
(\thrf{lemAduIdmA}). 
Il est clair que \emph{1a} implique \emph{1b} Inversement si l'on a un 
\cdrz~$\gL=\gB\sur{\DB(0)}$, où~$\gB$ est une~\Klg \stfez, le lemme \ref{lemAquQuoRed}
fournit un \idm $e$, et celui-ci est
in\dcp parce que $\gL$ est connexe.
\\
 Voyons le point \emph{\ref{i2thUnici}}. Soit $\gM$ un \cdr pour~$f$.
 \'Ecrivons

\snic{f(T)=\prod_{i=1}^n(T-\xi_i)$ dans $\gM.}

%\sni
Par la \prt \uvle de $\Adu_{\gK,f}$, il existe un unique \homo de~\Klgs $\varphi
:\gA\to\gM$ tel que $\varphi (x_i)=\xi_i$ pour $i\in\lrbn$.  Soit
$(e_\ell)_{\ell=1,\ldots ,k}$ l'orbite de $e$. C'est un \sfioz, donc $\big(\varphi(e_\ell)\big)_{\ell=1,\ldots ,k}$ \egmtz, et, puisque $\gM$ est un \cdiz, cela implique qu'il y a
un $j$ pour lequel $\varphi(e_\ell)=\delta_{j,\ell}$ (symbole de Kronecker).
\\
Alors, $\gen{1-e_j}\subseteq\Ker\varphi$, donc $\gM$ est un quotient de
$\gA\sur{\DA(1-e_j)}$, qui est un \cdiz.  Comme $\gM$ est non
trivial, cela implique $\gM\simeq \gA\sur{\DA(1-e_j)}$.  Enfin les
$\gA\sur{\DA(1-e_\ell)}$ sont deux à deux isomorphes.
\end{proof}
%------ end proof -----

\comm
Dans \cite{MRR}, il est montré que tout \cdi énumérable
possède une clôture algébrique. Cependant, un \cdr pour
$f$, qui existe donc, ne possède pas \ncrt une base finie comme \Kevz, au
sens des \comaz. Et l'on ne connaît pas de \tho d'unicité \cof pour un
tel \cdrz.
On peut décrire  comme suit une procédure analogue à celle de
\cite{MRR} pour obtenir un \cdr pour $f$.
Tout d'abord on construit une énumé\-ration~$(z_m)_{m\in\NN}$ de l'\aduz.
Ensuite on construit une suite d'\itfs $(\fa_m)$ de
$\gA$ en posant $\fa_0=0$, et  $\fa_{m+1}=\fa_m+\gen{z_m}$ si
$\fa_m+\gen{z_m}\neq \gen{1}$,  et $\fa_{m+1}=\fa_m$ sinon
(le test fonctionne car on peut calculer une base du \Kev $\fa_m+\gen{z_m}$). 
Alors,
l'\id $\bigcup_m\fa_m$ est un \idema de $\gA$, et le quotient est un \cdrz, qui est discret. Notre point de vue est légèrement
différent.
Nous ne partons pas a priori d'un corps énumérable, et même dans le
cas d'un corps énumérable, nous ne privilégions pas une
énumération par rapport à une autre. Nous nous contentons plutôt
de répondre aux questions concernant le \cdr au fur et à
mesure qu'elles se posent, comme on va le voir dans le \tho qui suit.
\eoe  

\medskip
Le \tho suivant 
explique comment contourner la difficulté
que pose la non existence du \cdr en \comaz. Le \cdr
 de $f$ est remplacé par une \gui{approximation}
donnée sous forme d'un quotient réduit $(\aqo{\gA}{1-e})\red$ de
l'\aduz, avec $e$ un \idm galoisien.
\\
On s'appuie sur le fait suivant qui est
déjà établi dans le cadre  \gnl
des anneaux \zeds  (lemme \ref{lemme:idempotentDimension0}).
Nous en rappelons une preuve directe.

\emph{Pour tout $y\in\gA=\Adu_{\gK,f}$, il existe un \idm $e_y\in\gK[y]\subseteq\gA$ tel que~$y$ est \iv
modulo $1-e_y$ et nilpotent modulo $e_y$.\label{factIdmAduPolmin}
}
\begin{proof}
Soit $P(T)$  le  \polmin de $y$. Il existe un \elt \ivz~$v$ de $\gK$
tel que $vP(T)=T^k\big(1-TR(T)\big)$ avec $k\geq0$. L'\idm $e_y$ est~$\big(yR(y)\big)^k$.
\end{proof}
%

%--- Theorem{thdivzeridm}-----
\begin{theorem}
\label{thdivzeridm} \emph{(Gestion dynamique d'un \cdrz)}\\
Soit $(z_i)_{i\in I}$ une famille finie d'\elts de $\Adu_{\gK,f}=\gA$. Il
existe un \idm galoisien $e$ de $\gA$ tel qu'en posant $\gB=\aqo{\gA}{1-e}$
chaque $\pi(z_i)$ est nul ou \iv dans l'\alg quotient $\gB\red$ (ici, $\pi:\gA\to\gB\red$ est la
projection canonique).
\end{theorem}
%--- end-theorem------------------------------------
%-----------------begin proof------------------
\begin{proof}
Pour chaque $i\in I$ il y a un \idm $g_i\in\gA$ tel que $z_i$ est \iv modulo
$1-g_i$ et nilpotent modulo $g_i$.  Appliqué à la famille des $g_i$ le
\thrf{lemAduIdmA} donne un \idm galoisien $e$, tel que pour chaque~$i$, $1-e$
divise $g_i$ ou $1-g_i$. Donc, dans l'\alg quotient $\gB=\aqo{\gA}{1-e}$
 chaque~$\pi(z_i)$ est nilpotent ou \ivz.
\end{proof}
%-----------------end proof------------------

\rems 1) 
\Llec peut s'inquiéter du fait que l'on ne dispose pas a priori d'un \sys
\gtr fini de l'\id $\DA(1-e)$. En conséquence l'\alg finie $\gB\red$ n'est
pas \ncrt un \Kev de dimension finie au sens
\cofz. En fait les nilpotents peuvent être gérés eux aussi de façon
dynamique. On a dans $\gB=\aqo{\gA}{1-e}$ un test de nilpotence et si un \elt
$x$ nilpotent est mis en évidence, on peut quotienter $\gB$ par l'\idz~$\fa$
engendré par l'orbite de $x$ sous l'action de $G=\St_\Sn(e)$.  Alors,
$\gB\sur{\fa}$ est de dimension finie et $G$ opère sur~$\gB\sur{\fa}$.

 2) Dans le \thrf{thdivzeridm}  on peut avoir intérêt à saturer la famille
$(z_i)_{i\in I}$ par l'action de $\Sn$ de façon à rendre manifestes
dans $\gB$ tous les \gui{cas de figure} possibles.
\eoe

%%%%%%%%%%%%%%%%%%%%%%%%%%%%%%%%%%%%%%%%%%%%%%%%%%%%%%%%%%%%%%%%%%%%%%%%%%%
%: --- subsection{Unicité du \cdr}-
\subsec{Unicité du \cdr}

Le \tho d'unicité du \cdr admet une  version \cov
\gui{opératoire} (qui fonctionne à tout coup, même si l'on ne dispose
pas d'un \idm in\dcp dans l'\aduz) sous la forme suivante.

%--- Theorem{propUnicite}------
\begin{theorem}
\label{propUnicite} \emph{(Unicité du \cdrz, version dynamique)}\\
Soient deux \Klgs \stfes $\gB_1$, $\gB_2$ non nulles pour lesquelles 
le \pol $f$ se
décompose en produit de facteurs \lins dans $(\gB_1)\red$ et
$(\gB_2)\red$. On suppose en outre que $(\gB_i)\red$ est
engendrée par les zéros correspondants de $f$. Alors, il existe un \idm
galoisien $e$ de $\gA$ tel que, avec l'\alg $\gB=\aqo{\gA}{1-e}$, on a deux
entiers $r_i$ tels que $(\gB_i)\red \simeq \gB\red^{r_i}$.
\end{theorem}
%--- end-theorem-----------------------------------------
%-----------------begin proof------------------
\begin{proof}
Le lemme \ref{lemAquQuoRed} donne des \idms $e_1$, $e_2 \in \gA$ tels que

\snic{(\gB_i)\red \simeq \gA\sur{\DA(1-e_i)}\quad (i=1,\,2)}

%\sni
Le \thref{lemAduIdmA} point
\emph{2}   donne un \idm galoi\-sien $e$ \hbox{et 
 $r_1$, $r_2\in\NN$}
tels que $\aqo\gA{1-e_i}\simeq \gB^{r_i}$.  Donc $(\gB_i)\red\simeq
\gB\red^{r_i}.$
\end{proof}
%-----------------end proof------------------

\entrenous{ On aimerait ajouter qqch du style suivant
%: entrenous --- SUBsubsection{Plus d'information sur le cas non \splz}-
\subsubsection*{Plus d'information sur le cas non \splz}

\hum{Il serait bon de faire subir
à un nilpotent arbitraire dans une \apG un traitement du même style que
celui qui a été fait aux \idms dans le \thref{lemAduIdmA}.
En attendant voici quelques petites choses avec l'\adu qui semblent
raisonnables.

Mais peut être le fin mot de l'affaire est plutôt du coté des \gui{\fcns \splsz} qui sont par exemple très bien traitées dans \cite{MRR} chapitre 6.}

\rm
Lorsque $f$ n'est pas \splz, $\disc\,f=0$ et $z=x_1-x_2$
est un \elt non \iv de~$\gA$.

On considère alors l'\idm $e_z\neq 1$ de~$\gA$.

Si $e_z=0$ tous les $x_i-x_j$ sont nilpotents. En passant au quotient
par l'\id qu'ils engendrent, on obtient une \Klg isomorphe à
$\gK[x_1]\simeq\aqo{\gK[X_1]}{f(X_1}$.

Alors, le \cdr de $f$
est isomorphe à $\gK[x_1]\red$. C'est une extension purement 
radicielle de $\gK$?????????

\smallskip Si $e_z\neq 0$ on considère un \idm galoisien
$e$ correspondant (\thrf{lemAduIdmA} et \algo \ref{algidmgal}).
Dans le quotient de Galois $\aqo{\gA}{1-e}$, certains des $\pi(x_i-x_j)$
sont nilpotents, d'autres ($\pi(z)$ par exemple) \ivsz.
On passe au quotient par l'\id engendré par les $\pi(x_i-x_j)$ nilpotents.
Soit $\gC$ la \Klg obtenue. On remplace  $G=\St(e)$ par le groupe de
permutation
des $x_i$ \gui{restants} et l'on obtient un groupe  $H$
de $\gK$-\autos de $\gC$. Soit l'\alg $\gL=\gC^H$.

Alors  $(\gL,\gC,H)$ est une \aGz, et
$\gL$ est une \gui{approximation} d'une extension purement radicielle de
$\gK$, plus \prmt $\gL\red$ est un corps, c'est une extension
purement radicielle de $\gK$?????
\\ \sf}
%: fin entrenous

%--- SECTION{Théorie de Galois}- secThGB
\section[Théorie de Galois d'un \pol \splz]{Théorie de Galois d'un \pol \spl sur un corps discret} %par en bas
\label{secThGB}

\medskip 
\Grandcadre{
Dans la section \ref{secThGB},  $\gK$ est un \cdi non trivial,
\\  $f$ est un \pol \spl unitaire de degré $n$ et $\gA=\Adu_{\gK,f}$.\\
Nous soulignons le fait que $f$ \emph{n'est pas supposé irréductible}.
}

\medskip 
Rappelons que pour un corps \splz ment factoriel, tout \pol \spl possède un corps
de racines (corolaire \ref{propIdemMini}), unique à \auto près (\thref{propUnicCDR}).
Nous sommes intéressés maintenant par le cas où le corps \emph{n'est pas}
\splz ment factoriel (ou même par le cas où la \fcn des \pols \spls est trop co\^uteuse).

Nous donnons ici comme promis la version \cov et dynamique de la théorie de Galois
d'un \pol \spl sur un \cdiz.

%: --- SUBsubsection{Existence et unicité}------------
\subsect{Existence et unicité du \cdrz, statique,\\ dynamique}{Existence et unicité du \cdrz}

Le fait \ref{factDiscriAdu} (ou le corolaire \ref{corcorlemEtaleEtage}) nous assure que $\gA$ est une \Klg étale.
Il en va de même
pour ses quotients de Galois.
Le  \thref{thUnici} se relit comme suit.

%:    THO{thUnici}%bis
\penalty-2500
\THO{thUnici}%bis
{\emph{(\Pol \splz: quand un \cdr existe et est une extension \stfez)}
\label{thExistence}
%-----------------begin enum------------------
\begin{enumerate}
\item \label{thExistence1}
\Propeq
%-----------------begin enum------------------
\begin{enumerate}
\item \label{thExistence1a} Il existe dans $\gA=\Adu_{\gK,f}$  un \idm
in\dcp $e$.
\item \label{thExistence1a'}\relax Il existe une extension \stfe  $\gL$ de $\gK$
qui est un \cdr de $f$.
\item \label{thExistence1c} L'\agB $\BB(\gA)$ est finie.
\end{enumerate}
%-----------------end enum------------------
\item \label{thExistence2}
Dans ce cas tout \cdr de $f$ est une extension galoisienne de~$\gK$, isomorphe 
à~$\gA[{1/e}]$.
\end{enumerate}
%-----------------end enum------------------

}

\medskip
Le point \emph{\ref{thExistence2}} résulte aussi du fait que si
un \cdr existe et est \stf sur $\gK$, deux \cdr sont isomorphes (\thrf{propUnicCDR}).

Le \tho d'unicité \ref{propUnicite} se relit comme suit.

%:    Tho{propUnicite}bis
\THO{propUnicite}%bis
{\emph{(Unicité du \cdr d'un \pol \splz, version dynamique)} \'Etant données deux \Klgs \stfes $\gB_1$ et $\gB_2$ non nulles dans
lesquelles $f$ se décompose en produit de facteurs \lins et qui sont
engendrées par les zéros correspondants de $f$, il existe un
quotient de Galois $\gB=\gA[1/e ]$ 
de l'\adu et deux entiers $r_i$ tels que
$\gB_1\simeq\gB^{r_1}$ et~$\gB_2\simeq\gB^{r_2}$.

}

%:--- SUBsec Theoreme de structure  ..
\subsect{Structure des quotients de Galois de l'\adu}{Quotients de Galois de l'\adu}

Pour la suite de la section \ref{secThGB} nous fixons les notations suivantes.
%--- Notations{notas1.3}---------
\begin{notations}
\label{notas1.3}\relax \emph{(Contexte d'un quotient de Galois)}\\
Soit $e$ un \idm galoisien de $\gA=\Adu_{\gK,f}$, $\fb=\gen{1-e}_\gA$.
On note

\snic{\gB=\gA\sur\fb=\gA[{1/e}] ,\;  \pi=\pi_{\gA,\fb}: \gA\to\gB,\;\hbox{ et }\;G=\St_\Sn(e).}

%\sni
On note $(e_1,\ldots ,e_m)$
l'orbite de $e$ sous $\Sn$. Chaque \Klg
$\gA[{1/e_i}]$ est isomorphe à~$\gB$. Le groupe $G$ opère sur~$\gB$.  
\end{notations}
%--- end-notation-----------------------------------------

Notons que pour $y\in\gB$, le \pol $\Mip_y(T)$ est \spl (parce que $\gB$ est
étale sur $\gK$). En outre $y$ est \iv \ssi $\Mip_y(0)\neq 0$.  Notons aussi
qu'un \itf de $\gB$ (différent de $\gen{1}$) est un \id galoisien \ssi son
orbite sous $G$ est formée d'\ids deux à deux
\com (tout \itf est engendré par un \idmz, et fait~\ref{factconnexe}).

Le \tho de structure \ref{thADG1Idm} se décline comme suit,
compte tenu des \thrfs{thdivzeridm}{thAduAGB}.

%--- Thm{thStruc3}--------
\begin{theorem}
\label{thStruc3} \emph{(\Tho de structure galoisien, 3)}
 Dans le contexte~\ref{notas1.3} on obtient les résultats suivants.
\begin{enumerate}
\item \label{thADG3Idm3} \label{thADG3Idm2}
$(\gK,\gB,G)$ est un quotient de Galois de
$(\gK,\gA,\Sn)$. \\
En particulier, $\gB$ est un \Kev de dimension
finie $\abs{G}$  et pour tout $y\in\gB$,
$\rC {\gB/\gK}(y)(T)=\rC G(y)(T)$. En outre, $\Fix(G)=\gK$.
\item \label{thADG3Idm1} 
On a un \iso de \Klgs $\gA\simeq \gB^m$.
\item \label{thADG3Idm4} Si $\gB$ est connexe, c'est un \cdr
pour $f$ et une extension galoisienne de $\gK$ avec $G$ comme 
groupe de Galois.
\item \label{thADG3Idm8} 
Soit $(y_i)$ une famille finie d'\elts de $\gB$.  Il existe un \idm galoisien
$e_\gB$ de $(\gK,\gB,G)$, tel que dans $\gB[{1/e_\gB}]$, chacun des $y_i$ 
est nul ou \ivz.
\item \label{thADG3Idm7} 
La restriction $\pi : e\gA \to \gB$ est un \iso $\gK$-\lin et établit une
correspondance biunivoque entre les \idms galoisiens
de $(\gK,\gA,\Sn)$ contenus dans $e\gA$ et ceux de $(\gK,\gB,G)$.
Les stabilisateurs et quotients résiduels sont préservés;
\cad que si $e_\gA \in e\gA$ et $e_\gB \in \gB$ sont deux \idms galoisiens
qui se correspondent, alors $\St_\Sn(e_\gA) = \St_G(e_\gB)$ et 
$\gA[{1/e_\gA}] \simeq \gB[{1/e_\gB}]$.
\end{enumerate}
%-----------------end enum------------------
\end{theorem}
%--- end-theorem-----------------------------------------
NB: dans la suite on donne les énoncés uniquement pour la situation
relative, la situation absolue est en effet le cas particulier où $e=1$.

%--- Lemma{Thidgal}-- --

\begin{lemma}
\label{Thidgal} \emph{(Résolvante et \polminz)} 
  Contexte
\ref{notas1.3},   $y\in\gB$.
%-----------------begin enum------------------
\begin{enumerate}
\item \label{Thidgal5a} $\Rv_{G,y}(T)$ est à \coes dans $\gK$.
\item \label{Thidgal5b} $\Mip_y$ divise $\Rv_{G,y}$ qui divise une puissance
de $\Mip_y$.
\item \label{Thidgal5c} $\rC {\gB/\gK}(y)(T)= \rC G(y)(T) =
\Rv_{G,y}(T)^{|\St_G(y)|}$.
\end{enumerate}
\end{lemma}
%--- end-lemma-----------------------------------------
%-----------------begin proof------------------
\begin{proof}
\emph{\ref{Thidgal5a}.} Conséquence du point
\emph{\ref{thADG3Idm3}} dans le \tho de structure. 
 
\emph{2.} On en déduit  que  $\Mip_{y}$
divise $\Rv_{G,y}$, car $\Rv_{G,y}(y)=0$. 
Et comme chaque~$y_i$ annule  $\Mip_y$,
le produit des $T-y_i$ 
%:HHH annule remplacé par divise
divise une puissance de~$\Mip_y$.
 
\emph{\ref{Thidgal5c}.} La deuxième \egt est évidente, et la
première est dans le point~\emph{\ref{thADG3Idm3}} du \tho de structure.
\end{proof}
%-----------------end proof------------------

%:--- subsection{Ou se passent les calculs}------------
\subsec{Où se passent les calculs}

Rappelons que $f$ est un \polu \spl de $\KT$ avec $\gK$ un \cdi non trivial.

\Grandcadre{
On note $\gZ_0$ le sous-anneau de $\gK$ engendré par les \coes de $f$ et  
\\par $1/\!\disc(f)$. On note
$\gZ$ la \cli de $\gZ_0$ dans $\gK$.
}

Nous mettons ici en évidence
 que \gui{tous les calculs se passent, et tous les résultats s'écrivent, dans l'anneau $\gZ$}, comme
cela résulte des \thosz~\ref{thIdmEtale} et \ref{thidGTri}\footnote{Il s'ensuit que si $\gK$ est un corps \gnl (voir section \ref{secAloc1}),
les questions de calculabilité se discutent en fait entièrement
dans $\Frac(\gZ)=\Frac(\gZ_0)\te_{\gZ_0}\gZ={(\gZ_0\sta)^{-1}}\gZ$.
Et~$\Frac(\gZ)$ est discret si $\gZ_0$ est lui-même un anneau discret.
Comme $\gZ_0$ est un anneau \tfz, il est certainement, en \clamaz, un anneau
effectif (on dit encore calculable) avec test d'\egt explicite,
au sens de la théorie de la récursivité
via les machines de Turing.
\\
Mais ce dernier résultat n'est pas une approche vraiment satisfaisante
de la réalité du calcul. Il s'apparente en effet aux résultats de \clama du style
\gui{tout nombre réel récursif admet un développement en fraction continue récursif}. Théorème manifestement faux d'un point de vue pratique, puisque pour le mettre en \oe uvre, il faut d'abord savoir si le nombre est rationnel
ou pas.}.
\\
Ces \thos nous donnent dans le cadre présent les
points~\emph{1},~\emph{2} et~\emph{4} du \tho qui suit.
Quant au point~\emph{3}, c'est une conséquence \imde du point \emph{2.}

%:     Theorem{thZsuffit}
\begin{theorem}\label{thZsuffit} \emph{(Le sous-anneau $\gZ$ de $\gK$
est bien suffisant)}
\begin{enumerate}
\item Soit $\gZ_1$ un anneau intermédiaire entre $\gZ$ et $\gK$ (par exemple $\gZ_1=\gZ$).\\ 
Alors, les \adus 

\snic{\Adu_{\gZ_0,f}\subseteq\Adu_{\gZ,f}\subseteq\Adu_{\gZ_1,f}\subseteq\Adu_{\gK,f}}

%\sni
sont des \aGs (relativement à leurs anneaux de base, et avec le groupe $\Sn$).
\item Tout \idm de $\gA=\Adu_{\gK,f}$ est dans $\Adu_{\gZ,f}$: ses  coordonnées  sur la base $\cB(f)$ sont dans $\gZ$.
\item  Les \emph{théories de Galois de $f$}
sur $\gZ$, sur $\gZ_1$, sur $\Frac(\gZ)$ et sur $\gK$ sont \emph{identiques}, au sens suivant.
\begin{enumerate}
\item Tout quotient de Galois de $\Adu_{\gZ_1,f}$ est obtenu par \eds à 
$\gZ_1$ à partir d'un quotient de Galois de $\Adu_{\gZ,f}$.
\item Tout quotient de Galois de $\Adu_{\Frac(\gZ),f}$ est obtenu par \eds à 
$\Frac(\gZ)$ à partir d'un quotient de Galois de $\Adu_{\gZ,f}$.
Cette \eds revient d'ailleurs à passer à l'anneau total des
fractions du quotient de Galois.
\item Tout quotient de Galois de $\Adu_{\gK,f}$ est obtenu par \eds à 
$\gK$ à partir d'un quotient de Galois de $\Adu_{\Frac(\gZ),f}$.
\end{enumerate} 
\item Soit $e$ un \idm galoisien de $\gA$ et $\gZ_1$ le sous-anneau de $\gZ$
engendré par $\gZ_0$ et les coordonnées de $e$ sur $\cB(f)$.
Alors, la structure triangulaire de l'\id de $\gZ_1[\Xn]$ engendré par $1-e$
et les modules de Cauchy est explicitée au moyen de \pols à \coes dans $\gZ_1$.\\
Pour \ceux qui connaissent les \bdgsz: la \bdg (pour un ordre monomial lexicographique) de l'\id qui définit l'approximation correspondante du \cdr de $f$ est formée de \polus à \coes dans $\gZ_1$.
\end{enumerate}

\end{theorem}

Notez les simplifications dans les cas particuliers suivants. Si $\gK=\QQ$ 
\linebreak 
et $f
\in \ZZ[T]$ \monz, alors $\gZ=\gZ_0=\ZZ[1/\!\disc(f)]$. De même, pour $q$
puissance d'un premier et $\gK$ le corps de fractions rationnelles $\gK =
\FF_q(u)$ on a, si~$f\in \FF_q[u][T]$ est \monz,
$\gZ=\gZ_0=\FF_q[u][1/\!\disc(f)]$.

\medskip
\rems 1) L'étude expérimentale suggère que non
seulement \gui{$\gZ$ est suf\-fisant},
mais qu'en fait tous les résultats  de calculs (\coes d'un \idm
sur $\cB(f)$, \bdg d'un \id galoisien) n'utilisent comme \denos que des \elts
dont le carré divise le \discri de $f$. \perso{Help! La théorie des nombres
ne pourrait elle pas venir à la rescousse?}

 2) \emph{Théorie de Galois absolue d'un \polz.}
\'Etant donné un \pol \spl $f\in\KT$, plutôt que de considérer $\gK$ et la
\cli $\gZ$ de~$\gZ_0$ dans $\gK$,
on peut considérer $\gK'=\Frac(\gZ_0)$
et la \cli $\gZ'$ de~$\gZ_0$ dans $\gK'$. \eoe

%%%%%%%%%%%%%%%%%%%%%%%%%%%%%%%%%%%%%%%%%%%%%%%%%%%%%%%%%%%%%%%%%%%%%%%%%%%
%: --- SUBSUBsec Methodes modulaires
\subsec{Changement d'anneau de base, méthode modulaire}

Puisque tout se passe dans $\gZ$, on peut regarder ce qui se passe
après une \eds $\varphi:\gZ\to \gk$ arbitraire.

Il se peut par exemple que $\gk$ soit un \cdi \gui{simple}
et que l'on sache calculer $\Gal_\gk\big(\varphi(f)\big)$,
\cad identifier un \idm indécomposable~$e'$ dans $\Adu_{\gk,\varphi(f)}$.
Ce groupe sera \ncrt (isomorphe à) un sous-groupe du
groupe de Galois inconnu $\Gal_\gZ(f)=\Gal_\gK(f)$.

Supposons que l'on ait calculé un quotient de Galois $\gB$ de $\Adu_{\gZ,f}$ avec un groupe $G\subseteq\Sn$. 

Si $e$ est l'\idm galoisien de $\Adu_{\gZ,f}$ correspondant à $\gB$, on peut se ramener au cas où 
$\varphi(e)$ est une somme de conjugués de $e'$ et où 

\snic{\Gal_\gk\big(\varphi(f)\big)= H=\St_G(e')\subseteq G}

%\sni
Comme ceci est vrai de tout quotient de Galois de $\Adu_{\gZ,f}$, on obtient
une double inclusion
\begin{equation}
\label{eqmetMod} H\subseteq \Gal_{\gZ,f}\subseteq G
\end{equation}
 à ceci près que le groupe $\Gal_{\gZ,f}$ est seulement défini
à conjugaison prés, et qu'il peut a priori demeurer à tout jamais inconnu.

Ce type de renseignements, \gui{le groupe de Galois de $f$ sur $\gK$, à conjugaison près, contient $H$} ne relève pas de la méthode dynamique que nous avons exposée, car celle-ci fait une démarche en sens contraire: donner
des informations du type \gui{le groupe de Galois de $f$ sur $\gK$, à conjugaison près, est contenu dans~$G$}.

Il est donc intéressant a priori d'utiliser en parallèle les deux méthodes,
dans l'espoir de déterminer complètement $\Gal_\gK(f)$.

Le fait d'avoir remplacé le corps $\gK$ par un sous-anneau est important de ce point de vue car on dispose de beaucoup plus de morphismes d'\eds de source
$\gZ$ que de source $\gK$.

\smallskip En particulier, il est souvent utile de travailler modulo $\fp$, un \idema de $\gZ$.
Cette méthode est alors appelée une \emph{méthode modulaire}.

Elle semble avoir été inventée par Dedekind pour la déter\-mi\-nation du groupe de Galois de $f$ sur $\QQ$ lorsque $f\in\ZZ[T]$.
Notons que dans ce cas un \idema de $\gZ=\gZ_0=\ZZ[1/\disc(f)]$ est 
donné par un nombre premier $p$ qui ne divise pas $\disc(f)$.

%-% ENTRE NOUS
\entrenous{Le mieux serait sans doute ici de rajouter un exemple
particulièrement pertinent.

Mais si l'on pouvait comprendre mieux comment tirer parti des informations
obtenues modulo un premier, cela mériterait une vraie section
}
%-% Fin ENTRENOUS

%%%%%%%%%%%%%%%%%%%%%%%%%%%%%%%%%%%%%%%%%%%%%%%%%%%%%%%%%%%%%%%%%%%%%%%%%%%
%: --- SUBsec Th Galois par en bas paresseuse
\subsec{Théorie de Galois paresseuse}
%-----------------------------------------

Le  \tho de structure \vref{thStruc3} et le lemme  \ref{Thidgal} 
(qui donne quelques précisions) sont les
 résultats \cofs théoriques qui permet\-tent une \evn paresseuse
du \cdr et du groupe de Galois d'un \pol \splz.

Notez bien que le terme \gui{paresseux} n'est absolument pas péjoratif. Il
indique simplement que l'on peut travailler en toute confiance dans le corps
des racines d'un \pol \spl sur $\gK$, même en l'absence de tout
\algo de factorisation des \pols sur $\gK$. En effet, toute anomalie avec
l'\alg $\gB$, l'approximation \gui{en cours} du \cdr de $f$, par
exemple la présence d'un diviseur de zéro non nul,  peut être
exploitée pour améliorer significativement notre connaissance du groupe
de Galois et du \cdrz. Un \idm  galoisien strictement multiple de l'\idm
$e$ \gui{en cours} peut en effet être calculé.
Dans la nouvelle \aGz, qui
est un quotient de la précédente, tous les calculs faits auparavant
restent valables, par passage au quotient. Par ailleurs, le nombre
d'améliorations significatives qui peuvent se produire ainsi n'excède
pas la longueur maximum d'une chaîne de sous-groupes de $\Sn$.

Nous développons donc une variante \gui{galoisienne} du système D5
 qui fut le premier système de calcul formel à calculer de
manière systématique et sans risque d'erreur dans la clôture \agq
d'un corps en l'absence d'\algo de factorisation des \pols (voir \cite[Duval\&al.]{D5}).

Ici, contrairement à ce qui se passe avec le système D5, l'aspect
dynamique des choses ne consiste pas à \gui{ouvrir des
branches de calcul séparées} chaque fois qu'une anomalie se présente,
mais à améliorer à chaque fois l'approximation du \cdr (et
de son groupe de Galois) que constitue le quotient de Galois en cours de
l'\aduz.

%: --- SUBSUBsec L'algo de base
\subsubsec{L'\algo de base}

On peut réécrire comme suit l'\algo \ref{algidmgal} dans la situation
présente, lorsque l'on dispose d'un \elt $y$ ni nul ni \iv dans le quotient
de Galois $\gB=\gA\sur\fb$.

%:     Algorithme{alg2idmgal}-------
\begin{algoR}[Calcul d'un \id galoisien et de son
stabilisateur]\label{alg2idmgal}
\Entree $\fb$: \id galoisien d'une \adu $\gA$ pour un \pol \splz, $\fb$
est donné par un \sgr fini; $y$: \elt \dvz dans %le quotient de Galois
$\gB=\gA\sur\fb$; $G=\St_\Sn(\fb)$;
$S_y=\St_G(y)$.
\Sortie $\fb'$: un \id galoisien de $\gB$ contenant $y$;
$H$:  $\St_G(\fb')$.
\Varloc $\fa$: \itf de $\gB$;
$\sigma$: \elt de $G$;
\\$L$: liste d'\elts de $\ov{G/S_y}$;
\hsu \# \, $\ov{G/S_y}$ est un \sys de représentants des classes
à gauche modulo $S_y$\\
$E$: ensemble correspondant d'\elts de $G/S_y$;
\hsu \# \, $G/S_y$  est l'ensemble des classes
à gauche modulo $S_y$. 
\Debut
\hsu $E \aff \emptyset$; $L \aff [\,]$;  $\fb' \aff \gen{y}$;
\hsu \pur{\sigma}{\ov{G/S_y}}
\hsd $\fa \aff \fb' + \gen {\sigma(y)}$;
\hsd \sialors{1\notin\fa} $\fb' \aff \fa$;  
$L \aff L\bullet [\sigma]$; $E \aff E \cup \so{\sigma S_y}$
\hsd \finsi;
\hsu \finpour;
\hsu $H\aff \St_G(E)$
\quad \# \, $H = \sotq{\alpha\in G} 
  {\forall\sigma\in L, \alpha\sigma\in \bigcup_{\tau\in L}\tau S_y}$.
\Fin
\end{algoR}
%--- fin Algorithme ------------------

L'\id $\fb$ est donné par un \sgr fini, et $G=\St_\Sn(\fb)$.  Soient~$e$
l'\idm de $\gB$ tel que $\gen{1-e}_\gB=\gen{y}_\gB$, et $e'$ un \idm galoisien
tel que $G.e$ et $G.e'$ engendrent la même \agBz.  
\\
On cherche à calculer
l'\id galoisien $\fc$ de $\gA$ qui donne le nouveau quotient de Galois
$\gA\sur\fc \simeq \gB\sur{\fb'}$, où $\fb'=\gen{1-e'}_\gB$, i.e. $\fc=\pi_{\gA,\fb}^{-1}(\fb')$.

Dans l'\algo \ref{algidmgal} on fait le produit de $e$ par un nombre maximum
de conjugués, en évitant d'obtenir un produit nul.

Ici, on ne calcule pas $e$, ni $\sigma(e)$, ni $e'$, car l'expérimentation
montre que souvent, le calcul de $e$ est très long (cet \idm occupe souvent
beaucoup d'espace mémoire, nettement plus que $e'$).  On raisonne alors
avec les \ids correspondants $\gen{1-e}=\gen{y}$ et
$\gen{1-\sigma(e)}=\gen{\sigma(y)}$.  Il s'ensuit que dans l'\algo le produit
des \idms est remplacé par la somme des \idsz.

Par ailleurs, comme on ne calcule pas $e$, on remplace $\St_G(e)$ par $\St_G(y)$,
qui est contenu dans $\St_G(e)$, en \gnl de façon stricte.  Néanmoins,
l'expérience montre que, bien que $\ov{G/S_y}$ soit plus grand, l'ensemble du
calcul est plus rapide.  Nous laissons le soin \alec de montrer que la
dernière affectation dans l'\algo fournit bien le groupe $\St_G(\fb')$
voulu, i.e. que le sous-groupe $H$ de~$G$, stabilisateur de $E$ dans $G/S_y$,
est bien égal à $\St_G(\fb')$.

%%%%%%%%%%%%%%%%%%%%%%%%%%%%%%%%%%%%%%%%%%%%%%%%%%%%%%%%%%%%%%%%%%%%%%%%%%%
%: --- SUBSUBsec Quand une resolvante rel
\subsubsec{Quand une résolvante relative se factorise}

Souvent une anomalie dans un quotient de Galois de l'\adu correspond à la
constatation qu'une résolvante relative se factorise. Nous traitons donc
ce cas en toute
\gntz, pour le ramener au cas de la présence d'un diviseur de
zéro non nul.

%:    Proposition{propRvRel1}----
\begin{proposition}
\label{propRvRel1} \emph{(Quand une résolvante relative se factorise)}\\
Dans le contexte \ref{notas1.3} soit $y\in\gB$ et $G.y=\so{y_1,\ldots,y_r}$.
%-----------------begin enum------------------
\begin{enumerate}
\item \label{RvRenum5}  Si $\Mip_y=R_1R_2$ avec $R_1$ et $R_2$ de degrés $\geq 1$,
$R_1(y)$ et $R_2(y)$ sont des \dvzs non nuls, et il existe un
\idm $e$ tel que $\gen{e}=\gen{R_1(y)}$ et
$\gen{1-e}=\gen{R_2(y)}$.
\item \label{RvRenum6}  Si    $\deg(\Mip_y)<\deg(\Rv_{G,y})$, alors l'un des
$y_1-y_i$ divise zéro (on peut donc construire un  \idm $\neq 0,1$ de
$\gB$).
\item \label{RvRenum7}  Si  $P$ est un diviseur strict de $\Rv_{G,y}$
dans $\gK[T]$, alors au moins un des deux cas suivants se présente:
%-----------------begin item------------------
\begin{itemize}
\item  $P(y)$ est un diviseur de zéro non nul, on est ramené au point {\ref{RvRenum5}.}
\item un \elt ${y_1-y_i}$ est un diviseur de zéro non nul, on est ramené au point {\ref{RvRenum6}.}
\end{itemize}
%-----------------end item------------------
\end{enumerate}
%-----------------end enum------------------
\end{proposition}
%--- end-proposition----------------------------------------
%-----------------begin proof------------------
\begin{proof}
\emph{\ref{RvRenum5}.}   Puisque $\Mip_y$ est \splz, $R_1$ et $R_2$ sont \comz.
Avec une relation de Bézout $U_1R_1+U_2R_2=1$, posons  $e=(U_1R_1)(y)$ et $e'=1-e$. On~a~$ee'=0$, donc~$e$
et~$e'$ sont des \idmsz. On a aussi \imdt 

\snic{eR_2(y)=e'R_1(y)=0$,
$eR_1(y)=R_1(y)$ et $e'R_2(y)=R_2(y).
}

%\sni
 Donc  $\gen{e}=\gen{R_1(y)}$ et
$\gen{1-e}=\gen{R_2(y)}$.

\emph{\ref{RvRenum6}.}  La preuve qui montre que sur un anneau intègre
un \pol unitaire de degré $d$ ne peut avoir plus que $d$ racines
distinctes se relit comme suit. 
\\
Sur un anneau arbitraire, si un \pol $P$
unitaire de degré $d$ admet des zéros $(a_1,\ldots ,a_d)$ avec tous les
$a_i-a_j$ \ndzs pour $i \ne j$, alors on~a~$P(T)=\prod(T-a_i)$. Donc si $P(t)=0$ et $t$
distinct des $a_i$, %alors 
deux au moins des $t-a_i$ sont \dvzs non nuls. 
On applique ceci au \polmin $\Mip_y$ qui a plus de zéros dans $\gB$
que son degré (ce sont les $y_i$). Ceci donne un $y_j-y_k$ \dvzz,
et par un $\sigma\in G$ on transforme  $y_j-y_k$ en un $y_1-y_i$.

\emph{\ref{RvRenum7}.}  Si $P$ est multiple de $\Mip_y$,
on est ramené au point
\emph{\ref{RvRenum6}}.\\
 Sinon,  $\pgcd(\Mip_y,P)=R_1$
est un diviseur strict de $\Mip_y$, et $R_1\neq1$
car on~a~$\pgcd\big((\Mip_y)^k,P\big)=P$ pour $k$ assez grand. Donc  $\Mip_y=R_1R_2$, avec~$\deg(R_1)$ et $\deg(R_2)\geq1$.
On est ramené au point~\emph{\ref{RvRenum5}}.
\end{proof}
%-----------------end proof------------------

On en déduit le corolaire suivant qui \gns
le point \emph{\ref{thADG3Idm8}} dans le \tho de structure \vref{thStruc3}.
%:     Theorem{corRvRel1}---------
\begin{theorem}
\label{corRvRel1} Dans le contexte \ref{notas1.3}
soit $(u_j)_{j\in J}$ une famille finie  dans $\gB$. Il existe un \id
galoisien $\fc$ de $\gB$ tel que, en notant $H=\St_G(\fc)$, $\gC=\gB\sur{\fc}$, et $\beta:\gB\to\gC$  la
\prn canonique, on~a:
%-----------------begin enum------------------
\begin{enumerate}
\item Chaque  $\beta(u_j)$ est nul ou \ivz.
\item Dans $\gC$, $\Mip_{\beta(u_j)}(T)=\Rv_{H,\beta(u_j)}(T)$.
\item Les $\Mip_{\beta(u_j)}$ sont deux à deux égaux ou étrangers.
\end{enumerate}
%-----------------end enum------------------
\end{theorem}
%--- end-theorem------------------------------------

\rem
Dans le \tho précédent  on a parfois intérêt à saturer la famille
$(u_j)_{j\in J}$ par l'action de $G$ (voire de $\Sn$ en remontant les $u_j$
dans $\gA$) de façon à rendre manifestes dans $\gC$ tous les \gui{cas
de figure} possibles.
\eoe

\medskip
\exl Nous reprenons l'exemple de la \paref{exemple1Galois}. On demande
à {\tt Magma} ce qu'il pense de l'\elt {\tt x5 + x6}. Trouvant que
la résolvante est de degré~$15$ (sans avoir besoin de la calculer)
alors que le \polmin est de degré~$13$,
il se met en peine de réduire la bizarrerie et obtient un quotient de Galois
de l'\adu de degré $48$ (le \cdr de degré $24$ n'est pas encore atteint)
avec le groupe correspondant. Le calcul est presque instantané.
Voici le résultat.

%:     verbatim
{\small\label{exemple2Galois}
\begin{verbatim}
y:=x5+x6;
MinimalPolynomial(y);
  T^13 - 13*T^12 + 87*T^11 - 385*T^10 +
     1245*T^9 - 3087*T^8 + 6017*T^7 - 9311*T^6 + 11342*T^5 - 10560*T^4 +
     7156*T^3 - 3284*T^2 + 1052*T - 260
//nouvelle algebre galoisienne, calculee a partir de deg(Min)<deg(Rv) :
Affine Algebra of rank 6 over Rational Field
Variables: x1, x2, x3, x4, x5, x6
Quotient relations:
  x1 + x2 - 1,
  x2^2 - x2 + x4^2 - x4 + x6^2 - x6 + 3,
  x3 + x4 - 1,
  x4^4 - 2*x4^3 + x4^2*x6^2 - x4^2*x6 + 4*x4^2 - x4*x6^2 + x4*x6 -
        3*x4 + x6^4 - 2*x6^3 + 4*x6^2 - 3*x6 - 1,
  x5 + x6 - 1,
  x6^6 - 3*x6^5 + 6*x6^4 - 7*x6^3 + 2*x6^2 + x6 - 1
Permutation group acting on a set of cardinality 6
Order = 48 = 2^4 * 3
    (1, 2)
    (3, 5)(4, 6)
    (1, 3, 5)(2, 4, 6)
\end{verbatim}
}

Nous renvoyons en exercices quelques cas particuliers de la situation examinée
dans la proposition \ref{propRvRel1}.
Chaque fois le but est d'obtenir des informations plus précises sur ce qui se passe lorsque l'on réduit la bizarrerie constatée.
Voir les exercices \ref{propGaloisJordan}, \ref{propRvRelDec} et
\ref{propRvRelDecMin}.

%%:     Proposition{propGaloisJordan}
%\begin{proposition}\label{propGaloisJordan}
%\end{proposition}
%%%%%%%%%%%%%%%%%%%%%%%%%%%%%%%%%%%%%%%%%%
%\begin{proof}
%\end{proof}
%
%%:     Proposition{propRvRelDec}
%\begin{proposition}\label{propRvRelDec}
%\end{proposition}
%%
%\begin{proof}
%\end{proof}
%%

%%:     Proposition{propRvRelDecMin}
%\begin{proposition}\label{propRvRelDecMin}
%\end{proposition}
%%
%\facile
%%%%%%%%%%%%%%%%%%%%%%%%%%%%%%%%%%%%%%%%%%%%%%%%%%%%%%%%%%%%%%%%%%%%%%%%%%%
%:   subsubsec Quand la structure triangulaire manque
\subsubsec{Quand la structure triangulaire manque}

Considérons des \elts $\alpha_1,\ldots ,\alpha_\ell$ de $\gB$ et les \Klgs emboîtées
  $$\gK\subseteq\gK_1=\gK[\alpha_1]\subseteq\gK_2= \gK[\alpha_1,\alpha_2] \subseteq\cdots
\subseteq\gK_\ell= \gK[\alpha_1,\ldots ,\alpha_\ell]\subseteq \gB.$$

Pour $i=2,\ldots ,\ell$ la structure de $\gK_i $
comme $\gK_{i-1}$-module peut être explicitée par différentes
techniques.
Si $\gB$ est un \cdr pour $f$, tous les $\gK_i$
sont des corps et chacun des modules est libre.

Si l'un de ces modules n'est pas libre, alors on peut construire un
\idmz~$\neq0,1$ dans $\gB$
en utilisant la même technique que pour la \dem du \tho de structure
\vref{th1Etale}, point \emph{2b}.

\smallskip
Il peut s'avérer efficace d'utiliser la technique des bases de Gr\"obner,
avec l'\id qui définit $\gB$ comme quotient de $\KXn$.
On introduit  des noms de variables $a_i$ pour les $\alpha_i$
et l'on choisit un ordre lexicographique
avec $a_1< \cdots <a_\ell<X_1< \cdots <X_n$.

Si $\gB$ est un corps la \bdg doit avoir
une structure triangulaire. \`A chaque $\alpha_i$
doit correspondre un et un seul \pol dans la \bdgz, $P_i(a_1,\ldots, a_i)$
unitaire en $a_i$.

Si cette structure triangulaire n'est pas respectée pour la variable
$a_i$,
nous sommes certains que $\gK_{i-1}$ n'est pas
un corps, et nous pouvons expliciter un diviseur de zéro dans
cette \Klgz.

En fait soit $P(a_1,\ldots ,a_i)$ un \pol qui apparaît dans
la \bdg et qui n'est pas unitaire en $a_i$. Son \coe dominant en tant que
\pol en $a_i$ est un \pol $Q(a_1,\ldots ,a_{i-1})$ qui
donne \ncrt un \elt diviseur de zéro  $Q(\alpha_1,\ldots ,\alpha_{i-1})$
dans l'\alg \zede
$\gK_{i-1}\simeq \gK[a_1,\ldots ,
a_{i-1}]\sur{\fa}$,
où $\fa$ est l'\id engendré par les premiers \polsz, en les variables
$a_1,\ldots ,a_{i-1}$, qui apparaissent dans la \bdgz.
Sinon, on pourrait multiplier $P$ par l'inverse de $Q$ modulo $\fa$,
et réduire le résultat modulo~$\fa$,
et l'on obtiendrait un \pol unitaire en~$a_i$ qui précède $P$
pour l'ordre lexicographique, et qui rendrait la présence de~$P$ inutile.

%-% ENTRE NOUS
\entrenous{%\medskip
\exl ?
\eoe
 
}
%-% Fin ENTRENOUS

%%%%%%%%%%%%%%%%%%%%%%%%%%%%%%%%%%%%%%%%%
%:section: Exercices
%: newpage
\newpage
\Exercices{

%--- Exercise{exoChapGalLecteur}-------------
\begin{exercise}
\label{exoChapGalLecteur}
{\rm  Il est recommandé de faire les \dems non données, esquissées,
laissées \alecz,
etc\ldots\, On pourra notamment traiter les cas suivants.
\begin{itemize}
\item \label{exopropBoolFini} Démontrer les propositions \ref{defiBoole},
\ref{propBoolFini} et le \thrf{corpropBoolFini}.
\item Expliquer le fait \ref{factChangeBase}.
%
%%
%\item

%%
%\item

\end{itemize}
}
\end{exercise}
%--- end -exercise-----------------------------------------

%--- Exercise{exothSteCdiClass}-------------
\begin{exercise}
\label{exothSteCdiClass} (Structure des \algs finies sur un corps, version classique, version \cov dynamique)
\\
{\rm \emph{1.} Montrer en \clama  le résultat suivant.
\\
\emph{Toute \alg finie  sur un corps est un produit
fini d'\algs locales finies. 
}
\\
\emph{2.} Expliquez pourquoi on ne peut espérer en
obtenir une \prcoz, même
en supposant que le corps est discret.
\\
\emph{3.} Proposer une version \cov du résultat précédent.
 
}
\end{exercise}
%--- end -exercise-----------------------------------------

%--- Exercise{exothNst1-zed}----------
%:2018 énoncé modifié
\begin{exercise} \label{exothNst1-zed}
{\rm  Montrer que la machinerie \lgbe \elr \num2 \paref{MethodeZedRed}
appliquée à la \dem du \thrf{thNstNoe} donne le résultat
suivant, \eqv au \tho \ref{thNstNoe} dans le cas d'un \cdi
non trivial.\imlg

%:    Theorem{thNst1-zed}------
\smallskip   \textbf{Théorème \ref{thNstNoe} bis}
\label{thNst1-zed}\relax (\nst faible et mise en position de \iNoez, cas des anneaux \zedrsz) 
{\it Soit $\gK$ un anneau \zedrz, $\ff=\gen{f_1,\dots,f_q}$ un \itf de $\gB=\KuX=\KXn$ et~$\gA=\KuX/\ff$ l'\alg quotient. Alors,
il existe un \sfio $(e_{-1},e_{0},\ldots,e_{n})$ de $\gK$ tel que, 
en notant 
$$
\gK_r=\gK[1/e_r],\;\gB_r=\gK_r\otimes_{\gK}\gB\simeq \gK_r[\uX] \hbox{ et }\gA_r=\gA[1/e_r]=\gK_r\otimes_{\gK}\gA\simeq \gB_r\sur{\ff\, \gB_r,}
$$ 
on ait les résultats suivants.
%-----------------begin item------------------
\begin{enumerate}
\item  $\gA_{-1}=0$, i.e. $1\in\ff\,\gB_{-1}$.
\item  $\gK_0\cap \ff\,\gB_{0}=0$ et $\gA_0$ est un $\gK_0$-module quasi libre fidèle. 
\item  Pour $r=1$, $\ldots$, $n$ on a un \cdv tel que, en appelant $Y_1,\ldots ,Y_{n}$ les nouvelles variables,
\begin{itemize}
\item [$\bullet$]  %$Y_1,\ldots , Y_{r}$ sont \agqt indépendantes sur $\gK_r$ et 
$\gK_r[Y_1,\ldots , Y_{r}]\cap \ff\,\gB_{r}=0$, autrement dit l'\alg  
$\gK_r[Y_1,\ldots , Y_{r}]$ peut être considérée comme une sous-$\gK_r$-\alg de $\gA_r$;
\item [$\bullet$] $\gA_r$ est un \mpf sur $\gK_r[Y_1,\ldots ,Y_{r}]$;
\item [$\bullet$] il existe un entier $N$
 tel que pour chaque $(\alpha_1,\ldots,\alpha_r)\in{\gK_r}^{\!\!r}$, 
 la $\gK_r$-\alg 
$\gD_r=\aqo{\gA_r}{Y_1-\alpha_1,\ldots,Y_r-\alpha_r}
$
 est un  $\gK_r$-module quasi libre fidèle engendré par au plus  $N$ \eltsz.
\end{itemize}
\end{enumerate}
%-----------------end item------------------
En particulier, la \Klg $\gA$ est un \mpf sur  la sous-\alg \gui{\polle} $\prod_{r=0}^n\gK_r[Y_1,\ldots ,Y_{r}]$.
\\
Pour $r\geq 1$, on dit que le \cdv (qui éventuellement ne change rien du
tout) a mis l'\id en \emph{position de Noether}.\\
Enfin, le \sfio qui intervient ici ne dépend pas du \cdv
qui met l'\id en position de Noether.
}
}
\end{exercise}
%--- end-exercise-----------------------------------------

%--- Exercise{exoMatMag3}-------------
\begin{exercise}\label{exoMatMag3}
{(Matrices magiques et \alg commutative)}
\\
{\rm  
On fournit dans cet exercice une application de l'\alg commutative
à un problème combinatoire; le \crc libre qui intervient
dans la mise en position de \Noe de la question \emph {2} est
un exemple de la \prt Cohen-Macaulay en terrain gradué.
Une \emph {matrice magique} de taille $n$ est une matrice de $\Mn(\NN)$
dont la somme de chaque ligne et de chaque colonne est la
même. L'ensemble de ces matrices magiques est un 
sous-\mo additif de $\Mn(\NN)$; on admettra ici que c'est le
\mo engendré par les $n\,!$ matrices de permutation.
On s'intéresse au dénombrement des matrices
magiques de taille $3$ de somme $d$ fixée. 
Voici les 6 matrices de permutation de $\MM_3(\NN)$:

\snic {
\begin {array} {rcl}
P_1 = \cmatrix {1&0&0\cr 0&1&0\cr 0&0&1},\
P_2 = \cmatrix {0&0&1\cr 1&0&0\cr 0&1&0},\
P_3 = \cmatrix {0&1&0\cr 0&0&1\cr 1&0&0}
\\
\noalign {\smallskip}
P_4 = \cmatrix {0&1&0\cr 1&0&0\cr 0&0&1},\
P_5 = \cmatrix {0&0&1\cr 0&1&0\cr 1&0&0},\
P_6 = \cmatrix {1&0&0\cr 0&0&1\cr 0&1&0}
\\
\end {array}
}

%\sni
Elles sont liées par la relation $P_1+P_2+P_3 = P_4+P_5+P_6$. 
Soit l'anneau de \pols à neuf \idtrs $\gk[(x_{ij})_{i,j \in \lrb{1..3}}]$
où $\gk$ est un anneau quelconque. 
On identifie une matrice $M = (m_{ij}) \in
\MM_3(\NN)$ au \mom $\prod_{i,j} x_{ij}^{m_{ij}}$, noté $\ux^M$;
par exemple~\smashtop{$\ux^{P_1} = x_{11}x_{22}x_{33}$}.

\emph {1.}
Soient $U_1$, \ldots, $U_6$ six \idtrs sur $\gk$ et $\varphi : \gk[\uU]
\twoheadrightarrow \gk[\ux^{P_1}, \ldots, \ux^{P_6}]$ défini par $U_i
\mapsto \ux^{P_i}$. On veut montrer que $\Ker\varphi$ est \hbox{l'\id
$\fa = \gen {U_1U_2U_3 - U_4U_5U_6}$}.
\vspace{-1pt}
\begin {enumerate}\itemsep0pt
\item [\emph {a.}]
Montrer, pour $a$, $b$, $c$, $d$, $e$, $f \in \NN$ et $m = \min(a,b,c)$, que:

\snic {
U_1^a U_2^b U_3^c U_4^d U_5^e U_6^f \equiv
U_1^{a-m} U_2^{b-m} U_3^{c-m} U_4^{d+m} U_5^{e+m} U_6^{f+m} \bmod \fa
}

\item [\emph {b.}]
On note $\fa^\bullet$ le sous \kmo de $\gk[\uU]$ ayant pour base les \moms
non divisibles par $U_1U_2U_3$. Montrer que $\gk[\uU] = \fa\oplus\fa^\bullet$
et que $\Ker\varphi = \fa$.
\item [\emph {c.}]
En déduire que le nombre $M_d$ de matrices magiques de taille 3 et de 
somme~$d$ est égal à
${d+5 \choose 5} - {d+2 \choose 5}$ en convenant de ${i \choose j} = 0$ pour
$i < j$.
\end {enumerate}

\emph {2.}
Soit $\gB = \gk[\uU]\sur\fa = \gk[\uu]$. 
\vspace{-1pt}
\begin {enumerate}\itemsep0pt
\item [\emph {a.}]
Définir une mise en position de \Noe $\gA = \gk[v_2,v_3,u_4,u_5,u_6]$ de
$\gB$ \hbox{où $v_2$, $v_3$} sont des formes \lins en $\uu$, de façon à ce que
$(1, u_1, u_1^2)$ soit une $\gA$-base \hbox{de $\gB = \gA \oplus \gA u_1
\oplus \gA u_1^2$}.

\item [\emph {b.}]
En déduire que le nombre $M_d$ est aussi égal à ${d+4 \choose 4} + {d+3
\choose 4} + {d+2 \choose 4}$ (formule de MacMahon, qui donne en passant, une
\idt entre \coes binomiaux).
\end {enumerate}

\emph {3.}
On suppose que $\gk$ est un \cdiz. On veut montrer que l'anneau $\gB$, vu
comme l'anneau $\gk[\ux^{P_1}, \ldots, \ux^{P_6}]$ est \icl (voir aussi
le \pb \ref {exoFullAffineMonoid}).  Soit $E \subset
\MM_3(\ZZ)$ le sous-\ZZmo des matrices magiques (\dfn analogue) et le
sous-anneau $\gB_{11} \subset \gk[x_{ij}^{\pm 1}, i,j \in \lrb{1..3}]$

\snic {
\gB_{11} = \gk[\ux^{P_1}, \ux^{P_6}]
[\ux^{\pm P_2}, \ux^{\pm P_3}, \ux^{\pm P_4}, \ux^{\pm P_5}]
}

%\sni
de sorte que $\gB_{11} \subset \gk[\,\ux^M \,\vert\, M \in E,\ m_{11} \ge 0\,]$.
\vspace{-1pt}
\begin {enumerate}\itemsep0pt
\item [\emph {a.}]
Vérifier que $\gB$ et $\gB_{11}$ ont même corps de fractions, qui
est le corps de fractions $\gk(E)$, corps de fractions rationnelles
sur $\gk$ à 5 \idtrsz.
\item [\emph {b.}]
Montrer que $\gB_{11}$ est \iclz.
\item [\emph {c.}]
Pour $i$, $j \in \lrb{1..3}$, définir un anneau $\gB_{ij}$ analogue
à $\gB_{11}$ et en déduire que $\gB$ est \iclz.
\end {enumerate}

}

\end {exercise}
%--- end -exercise-----------------------------------------

%--- Exercise{exoNonCyclique}----
\begin{exercise}
\label{exoNonCyclique}
{\rm  Donner une preuve directe (pas par
l'absurde) que si un \cdi possède deux \autos qui engendrent un
groupe fini non cyclique, le corps contient un $x\neq 0$ dont toutes les
puissances sont distinctes, \cad qui n'est pas une racine de l'unité.
}
\end{exercise}
%--- end-exercise-----------------------------------------

%-% ENTRE NOUS
\entrenous{ 

%--- Exercise{exoPohstZassenhaus}-------------
%\begin{exercise}
%\label{exoPohstZassenhaus}
{\rm  Faire un exo a partir de la preuve \elr dans \cite{PZ}
pour: si $\gk$ est réduit et si $\disc(f)$ est \ndz
alors l'\adu est réduite ?

}
%\end{exercise}
%--- end -exercise-----------------------------------------
}
%-% Fin ENTRENOUS

%--- Exercise{exoIdentiteDiscriminantale}-------------
\begin{exercise}\label{exoIdentiteDiscriminantale}
{(Une \idt \gui {discriminantale})}\\
{\rm  
Soit $n \ge 1$. On note $E$ l'ensemble des $\alpha \in \NN^n$
tels que $0 \le \alpha_i < i$ pour $i \in \lrbn$; c'est un
ensemble de cardinal $n!$ que l'on ordonne par la \gui{numération
factorielle}, i.e. $\alpha \preceq \beta$ si $\sum_i \alpha_i i! 
\le \sum_i \beta_i i!$. On ordonne le groupe \smq $\rS_n$
par l'ordre lexicographique, $\In$ étant la plus petite
permutation. On considère $n$ \idtrs sur $\ZZ$ et on
définit une matrice $M \in \MM_{n!}(\ZZ[\ux])$, indexée
par $\rS_n \times E$:

\snic {
M_{\sigma, \alpha} = \sigma(\ux^\alpha), \qquad 
\sigma \in \rS_n, \quad \alpha\in E
}

%\sni
Ainsi pour $n = 3$:

\snic {
M = \cmatrix {
1& x_2& x_3& x_2x_3&  x_3^2& x_2x_3^2\cr
1& x_3& x_2& x_2x_3&  x_2^2& x_2^2x_3\cr
1& x_1& x_3& x_1x_3&  x_3^2& x_1x_3^2\cr
1& x_3& x_1& x_1x_3&  x_1^2& x_1^2x_3\cr
1& x_1& x_2& x_1x_2&  x_2^2& x_1x_2^2\cr
1& x_2& x_1& x_1x_2&  x_1^2& x_1^2x_2\cr
}}

%\sni
\emph {1.}
Montrer que $\det(M) = \delta^{n!/2}$ avec $\delta = \prod_{i<j}(x_i-x_j)$.

\emph {2.}
On note $s_1$, \ldots, $s_n \in \ZZ[\ux]$ les $n$ fonctions \smqs\elrsz, $F(T)$
le \pol \uvl $F(T) = T^n - s_1T^{n-1} + \cdots + (-1)^n s_n$, et $U \in
\MM_{n!}(\ZZ[\ux])$ la matrice tracique, indexée par $E \times E$, de terme
$\Tr_{\rS_n}(\ux^{\alpha+\beta})$, $\alpha, \beta\in E$.
\\  
Soit $f \in \gk[T]$ un \polu de degré $n$, $\gA=\Adu_{\gk,f}$.\\
Retrouver l'\egt $\Disc_\gk \gA = \disc(f)^{n!/2}$ (fait \ref{factDiscriAdu}). 
Et réciproquement?

\emph {3.}
Revisiter le \tho \ref{theoremAdu3}.
}

\end {exercise}
%--- end -exercise-----------------------------------------

%--- Exercise{exoTnADU}-------------
\begin{exercise}\label{exoTnADU}
 {(L'\adu du \pol $f(T) = T^n$)}\\
{\rm  
Soient $f(T) = T^n$ et $\gA
=\Adu_{\gk,f} = \gk[x_1, \ldots, x_n]$.  Décrire la structure de~$\gA$.
}
\end {exercise}
%--- end -exercise-----------------------------------------

%--- Exercise{exoNilIndexInversiblePol}-------------
\begin{exercise}\label{exoNilIndexInversiblePol}
{(\Pols \ivs et indices de nilpotence)} \\
{\rm  
On propose ici une version quantitative du résultat figurant dans le point
\emph {4}  du lemme \ref{lemGaussJoyal}. Soient $\gk$ un anneau commutatif,
$f$, $g \in \kX$ vérifiant $fg = 1$ \hbox{et $f(0) = g(0) = 1$}. On
écrit $f = \sum_{i=0}^n a_i X^i$, $g = \sum_{j=0}^m b_j X^j$.  Montrer que 

\snic {
a_1^{\alpha_1} a_2^{\alpha_2} \cdots a_n^{\alpha_n} 
b_1^{\beta_1} b_2^{\beta_2} \cdots b_m^{\beta_m}  = 0
\quad \hbox {si}\quad \sum_i i\alpha_i + \sum_j j\beta_j > nm
}

%\sni
En particulier, pour $i \ge 1$, on a $a_i^{\lceil (nm+1)/i\rceil}
= 0$ et par suite $a_1^{nm + 1} = 0$.

}

\end {exercise}
%--- end -exercise-----------------------------------------

%--- Exercise{exoADUpthRoot}-------------
\begin{exercise}\label{exoADUpthRoot}
{(L'\adu du \pol $f(T) = T^p - a$ en \cara $p$)}\\
{\rm  
Soit $p$ un nombre premier, $\gk$ un anneau dans lequel $p\cdot 1_\gk=0$ et $a \in \gk$.  \\ 
On note $f(T) = T^p - a \in \gk[T]$, $\gA =\Adu_{\gk,f} = \gk[x_1, \ldots,
x_p]$, $\gk[\alpha] = \aqo{\gk[T]}{f}$, de sorte que $T-a=(T-\alpha)^p$.  Soit $\varphi : \gA \twoheadrightarrow
\gk[\alpha]$ le $\gk$-morphisme $x_i \mapsto \alpha$. 
Expliciter l'\id $\Ker\varphi$ et décrire la
structure de la \klg $\gA$.
\\
NB: si $\gk$  est un \cdi et $a$ n'est pas une puissance
$p$-ième dans $\gk$, d'après l'exercice \ref{exoPrimePowerRoot}, le \pol
$f(T)$ est \ird et $\gk[\alpha]$ est un
corps de \dcn de $f$ sur $\gk$.

}

\end {exercise}
%--- end -exercise-----------------------------------------

%--- Exercise{exoA5GaloisGroup}-------------
\begin{exercise}\label{exoA5GaloisGroup}
 {(Le trinôme $T^5+5bT\pm 4b$ où $b= 5a^2-1$,
de groupe de Galois $\rA_5$)}\\
{\rm  
On considère un trinôme $T^5 + bT + c$. On va particulariser
$b$, $c$ de façon à ce que son discriminant soit un carré
et obtenir un \pol \ird de groupe de Galois $\rA_5$ comme
illustration de la méthode modulaire. 
\\
On utilise
l'\egt $\disc_T(T^5 + bT + c) = 4^4b^5 + 5^5c^4$ (voir le \pbz~\ref{exoDiscriminantsUtiles}).

\emph {1.} Pour forcer le \discri à être un carré dans $\ZZ$,
expliquer pourquoi ce qui suit est raisonnable: $b \aff 5b$, $c \aff 4c$, puis  $f_a(T) = T^5 + 5(5a^2
-1)T \pm 4(5a^2 - 1)$. Le \discri est alors le carré $2^8 5^6 a^2(5a^2 - 1)^4$.

\emph {2.}
On prend $a=1$ et on obtient $f_1(T) = T^5 + 20T \pm 16$ dans $\ZZ[T]$.  En
examinant les \fcns de $f_1$ modulo $3$ et $7$, montrer que $f_1$ est
\ird de groupe de Galois $\rA_5$.  En déduire que pour $a \equiv 1 \bmod
21$, $f_a$ est \ird de groupe de Galois $\rA_5$. Montrer qu'il en est de même de
$f_a$ vu comme \pol à \coes dans le corps de fractions rationnelles
$\QQ(a)$.

}

\end {exercise}
%--- end -exercise-----------------------------------------

%--- Exercise{propGaloisJordan}-------------
\begin{exercise}
\label{propGaloisJordan}
{\rm
\emph{(Quand une résolvante admet un zéro dans le corps de base)}\\
Dans le contexte \ref{notas1.3}
soit $y\in\gB$,  $G.y=\so{\yr}$ et $g(T)=\Rv_{G,y}(T)$.
\begin{enumerate}
\item \label{i1propGaloisJordan}
On suppose que $a\in\gK$ est un zéro simple de $g$. 
\begin{enumerate}
\item \label{i1apropGaloisJordan} $\fc=\gen{y-a}_\gB$ est un \id galoisien de $(\gK,\gB,G)$%, c'est le plus petit \id galoisien contenant $y-a$
.
\item \label{i1bpropGaloisJordan}  Si $\beta:\gB\to{\gC=\gB/\fc}$ est
la \prn canonique, et si $H=\St_G(\fc)$ est la nouvelle approximation
du groupe de Galois, alors $\beta(y_1)=a$ et pour $j\neq 1$, $\Rv_{H,y_j}$
divise $g(T)/(T-a)$ (comme d'habitude on identifie~$\gK$ à un sous-corps de $\gB$ et $\beta(\gK)$ à un sous-corps de $\gC$).
\end{enumerate}
\item \label{i2propGaloisJordan}
On suppose que $a\in\gK$ est un zéro de $g$ avec la multiplicité $k$.
\begin{enumerate}
\item \label{i2apropGaloisJordan}
Il existe $j_2$, \ldots, $j_k\in\lrb{2..r}$ tels \hbox{que $\fc=\gen{y_1-a,y_{j_2}-a,\ldots,y_{j_k}-a}$} est un \id galoisien minimal parmi ceux qui contiennent $y-a$.\\
 On pose $j_1=1$. 
 Pour $j\neq j_1$, \ldots, $j_k$, $y_j-a$ est \iv modulo $\fc.$
\item \label{i2bpropGaloisJordan}
Soient $\beta:\gB\to{\gC=\gB/\fc}$ 
la \prn canonique, et  $H=\St_G(\fc)$. \\
Alors $\beta(y_{j_1})=\cdots=\beta(y_{j_k})=a$, et pour $j\neq j_1$, \ldots, $j_k$, la résolvante $\Rv_{H,y_j}$ divise $g(T)/(T-a)^k$.
\end{enumerate}
\item \label{i3propGaloisJordan}
On suppose que $\fc$ est un \id galoisien de $\gB$ et que $\St_G(y)$ contient
$\St_G(\fc)$, alors $g(T)$ admet un zéro dans $\gK$.
\end{enumerate}
\rem Le point \emph{1} %de l'exercice~\ref{propGaloisJordan}
justifie la \gui{méthode de Jordan} pour le calcul du groupe de Galois. Voir \paref{Jordan}.
\eoe

}
\end{exercise}
%--- end -exercise-----------------------------------------

%--- Exercise{propRvRelDec}-------------
\begin{exercise}
\label{propRvRelDec}
{(Quand on connaît la décomposition en facteurs premiers
 d'une résol\-vante \splz)}
{\rm Dans le contexte \ref{notas1.3} soit $y\in\gB$ et $G.y=\so{y_1,\ldots,y_r}$.
% et $H=\St(y)$.
\\
On suppose que  $\Rv_{G,y}=\Mip_y=R_1\cdots R_\ell$, avec les $R_i$ \irds
et $\ell>1$.
%On note $d_i=\deg(R_i)$.
Calculer  un  \idm galoisien $e$ de $\gB$, avec les \prts suivantes, dans lesquelles on note  $(\gK,\gC,H)$ le quotient de Galois correspondant et
 $\beta:\gB\to\gC$
la \prn canonique.
\begin{enumerate}
\item Pour chaque $i\in\lrbr$, le \pol $\Mip_{\beta(y_i)}$
est égal à l'un des $R_j$.
\item Le groupe $H$ opère sur  $\so{\beta(y_1),\ldots,\beta(y_r)}$. 
\item Les orbites
sont de longueurs $d_1=\deg(R_1),\ldots,d_\ell=\deg(R_\ell)$.
\item Cette situation se reproduit dans tout quotient de Galois de $(\gK,\gC,H)$.
\end{enumerate}
\rem  L'exercice~\ref{propRvRelDec}
est la base de la \gui{méthode de McKay-Soicher} pour le calcul du groupe de Galois. Voir \paref{soicher}.
\eoe

}
\end{exercise}
%--- end -exercise-----------------------------------------

%--- Exercise{propRvRelDecMin}-------------
\begin{exercise}
\label{propRvRelDecMin}
{\rm
\emph{(Quand un \polmin divise strictement une résol\-vante)}\\
Dans le contexte \ref{notas1.3} soit $y\in\gB$ et $G.y=\so{y_1,\ldots,y_r}$.\\
On suppose que  $g(T)=\Rv_{G,y}(T)\neq\Mip_y(T)$.
Soit $(\gK,\gC,H)$ un quotient de Galois  (avec la \prn canonique
$\beta:\gB\to\gC$) dans lequel chaque $\beta(y_i)$ admet un \polmin
égal à sa résolvante.\\
Montrer que pour les différents zéros $\beta(y_j)$ de $g_1(T)=\Mip_{\beta(y_1)}(T)$
dans $\gC$, les fibres $\beta^{-1}\big(\beta(y_j)\big)$
ont toutes le même nombre d'\eltsz, disons $n_1$. \\
En outre,  $g_1^{n_1}$
divise $g$ et $g/g_1^{n_1}$ est comaximal avec $g_1$.
}
\end{exercise}
%--- end -exercise-----------------------------------------

%%%%%%%%%%%%%%%%%%%%%%%%%%%%%%%%%%%%%%%%%%%%%%%%%%%%%%%%%%%%%%%%%%%%%%%%%%%
}% fin des exos
%:  solutions

\sol{

%%%%%%%%%%%%%%%%%%%%%%%%%%%%%%%%%%%%%%%%%%%%%%%%%%%%%%%%%%%%%%%%%%%%%%%%%%%

\exer{exothSteCdiClass} \emph{1.} 
Cela résulte du fait qu'un anneau \zed connexe est local et du fait
que, par le principe du tiers exclu, on connaît les \idms in\dcps de
l'\algz, lesquels forment un \sfioz.

 \emph{2.} Dans le cas d'une \alg $\aqo\KX f$ avec $f$ \splz, trouver les \idms revient à factoriser le \polz. Mais il n'existe pas d'\algo \gnl de \fcn d'un \pol \splz.

 \emph{3.} Une version 
\cov   consiste  à affirmer que,
pour ce qui concerne un calcul, on peut toujours \gui{faire comme si}
le résultat  (démontré au moyen du tiers exclu) était vrai.
Cette \emph{version dynamique} s'exprime comme suit. \\
\emph{Soit $\gK$ un anneau \zed (cas particulier: un \cdiz).\\
Soit $(x_i)_{i\in I}$ une famille finie d'\elts dans une  \Klg entière  $\gB$ (cas particulier:  une \Klg finie).
\\
Il existe 
un \sfio $(e_1,\ldots,e_n)$ tel que dans chaque
composante $\aqo{\gB}{1-e_j}$, chaque $x_i$ est nilpotent ou \ivz.} \\
On démontre ce résultat comme suit: le lemme \ref{lemZrZr1} 
nous dit que $\gB$ est \zedz; on conclut par le lemme de
scindage \zedz~\ref{thScindageZed}.

%%%%%%%%%%%%%%%%%%%%%%%%%%%%%%%%%%%%%%%%%%%%%%%%%%%%%%%%%%%%%%%%%%%%%%%%%%%
%:2018 ajout du corrigé
\exer{exothNst1-zed}
Lorsque $\gK$ est un \cdiz, l'\algo de \mpN est contrôlé en taille par une fonction $\varphi(n,q,d)$ où  $n$ est le nombre de variables, $q$ le nombre \pols \gtrs de $\ff$ et $d$ une majoration $d$ du degré total de ces \polsz. Lorsque l'on reproduit cet \algo pour un anneau \zedr en scindant l'anneau courant $\gL$ en deux chaque fois que l'on doit certifier \gui{$a=0$ ou $a$ \iv dans $\gL$} pour un \elt $a$ de $\gK$\footnote{Si $\gL$ est l'anneau courant et si $h$ est l'\idm tel que $\gen{h}=\gen{a}$ dans $\gL$, les deux branches sont $\gL[1/h]$ et $\gL[1/(1-h)]$.}, 
chaque branche de l'\algo est majorée en taille par $\varphi(n,q,d)$ et le nombre de branches est fini. On obtient en fin de compte un \sfio $(h_j)_{j\in J}$ de $\gK$ et pour chaque $\gK_j=\gK[1/h_j]$, en notant $\gA_j=\gA[1/h_j]$ l'\algo a produit un des résultats suivants: 
\begin{itemize}
\item  $\gA_j=0$  (on définit alors $r_j=-1$),
\item ou bien $\gA_j$ est un module libre de rang fini non nul sur $\gK_j$ (on définit alors $r_j=0$), 
\item ou bien on a  un \cdv tel que $\gA_j$ a une \mpN avec $r_j$
\elts \agqt indépendants ($r_j\in\lrbn$).
\end{itemize}
 On définit enfin  $e_r$ égal à la somme de tous les $h_j$ pour lesquels $r_j=r$. On obtient comme cela les résultats
énoncés dans le \thoz\footnote{Les \cdvs des $\gK_j[\uX]$ peuvent bien ne pas se ressembler pour différents $j$. Comme les $h_j$ forment un \sfioz, ces \cdvs mis tous ensemble fournissent un \cdv de $\KuX$.}.\\
Notons que certains $h_j$ ou certains $e_r$ peuvent être nuls, mais cela n'affecte pas les résultats tels qu'ils sont formulés, car l'anneau nul a toutes les \prts voulues: si l'on sait que $e_r=0$ pour un certain $r\in\lrbzn$ le point correspondant dans l'énoncé peut être supprimé car sans intérêt, mais il n'est pas faux pour autant. Naturellement, la situation est un peu plus confortable si $\gK$ est discret, car on peut alors supprimer toute branche donnant un anneau nul dès qu'elle apparait au moment d'un scindage.

%%%%%%%%%%%%%%%%%%%%%%%%%%%%%%%%%%%%%%%%%%%%%%%%%%%%%%%%%%%%%%%%%%%%%%%%%%%

\exer{exoMatMag3} 
On vérifie facilement que le $\ZZ$-module des relations \lins entre
les matrices $P_1$, \ldots, $P_6$ est engendré par $(1,1,1, -1,-1,-1)$.
On utilisera aussi que le nombre de \moms de degré $d$
en $n$ variables est ${d + n-1 \choose n-1}$.

\emph {1a.}
Soit $S(Y,Z) = \sum_{i+j=m-1} Y^i Z^j$, donc $Y^m - Z^m = (Y-Z)S(Y,Z)$.
Dans cette \egtz, on fait  $Y = U_1U_2U_3$, $Z = U_4U_5U_6$.
On obtient le résultat demandé en multipliant par $U_1^{a-m} U_2^{b-m} U_3^{c-m} U_4^d U_5^e U_6^f$.

\emph {1b.}
On a clairement $\fa\subseteq\Ker\varphi$. L'\egt $\gk[\uU] =
\fa + \fa^\bullet$ résulte du point \emph{1a}.  Il suffit donc de
voir que $\Ker\varphi \cap \fa^\bullet = \{0\}$, i.e. que la restriction de
$\varphi$ à $\fa^\bullet$ est injective.  Comme $\varphi$ transforme un
\mom en un \momz, il suffit de voir que si deux \moms $U_1^a \cdots
U_6^f$ et $U_1^{a'} \cdots U_6^{f'} \in \fa^\bullet$ ont même image par
$\varphi$, ils sont égaux. \\
On a $(a,b,c, \ldots, f) = (a',b',c',
\ldots, f') + k(1,1,1,-1,-1,-1)$ avec $k \in \ZZ$, et \hbox{comme $\min(a,b,c) =
\min(a',b',c') = 0$}, on a $k = 0$, ce qui donne l'\egt des deux \momsz.

\emph {1c.}
Le nombre $M_d$ cherché est la dimension sur $\gk$ de la composante
\hmg de degré $d$ de $\gk[\ux^{P_1}, \ldots, \ux^{P_6}]$ ou encore (via
$\varphi$) celle de $\fa^\bullet_d$. \\
Mais on a aussi $\gk[\uU] =
\fb\oplus\fa^\bullet$ où $\fb$ est l'\id (monomial) engendré par les
\moms divisibles par $U_1U_2U_3$ (en quelque sorte, $\fb$ est un \id
initial de $\fa$).  \\
On a donc $\gk[\uU]_d = \fb_d \oplus \fa^\bullet_d$ et

\snic {
\dim_\gk \gk[\uU]_d = {d+5 \choose 5}, \quad
\dim_\gk \fb_d = {d+5-3 \choose 5}, \quad
M_d = \dim_\gk \fa^\bullet_d = {d+5 \choose 5} - {d+2 \choose 5}
}

%\sni
\emph {2a.}
On définit $V_2, V_3$ par $U_2 = U_1 + V_2$, $U_3 = U_1 + V_3$. \\
Le \pol $U_1U_2U_3 - U_4U_5U_6$, vu dans
$\gk[U_1, V_2, V_3, U_4, U_5, U_6]$ devient unitaire en $U_1$ de degré 3.
On laisse au lecteur le soin de vérifier les autres détails.

\emph {2b.}
Le nombre cherché est aussi $M_d = \dim_\gk \gB_d$. Mais on a

\snic{\gB_d =
\gA_d \oplus \gA_{d-1} u_1 \oplus \gA_{d-2} u_1^2 \simeq 
\gA_d \oplus \gA_{d-1} \oplus \gA_{d-2}.}

%\sni
Il suffit d'utiliser
le fait que $\gA$ est un anneau de \pols sur $\gk$ à 5 \idtrsz.
\`A titre indicatif, pour $d = 0$, $1$, $2$, $3$, $4$, $5$, $M_d$
vaut $1$, $6$, $21$, $55$, $120$, $231$.

\emph {3a.}
Le $\ZZ$-module $E$ est libre de rang 5: 5 matrices quelconques parmi
$\{P_1, \ldots, P_6\}$ en constituent une $\ZZ$-base.

\emph {3b.}
Puisque $P_1 + P_2 + P_3 = P_4 + P_5 + P_6$, on a:

\snic {
\gB_{11} = 
\gk[\ux^{P_1}][\ux^{\pm P_2}, \ux^{\pm P_3}, \ux^{\pm P_4}, \ux^{\pm P_5}] =
\gk[\ux^{P_6}][\ux^{\pm P_2}, \ux^{\pm P_3}, \ux^{\pm P_4}, \ux^{\pm P_5}]
.}

%\sni
On voit alors que $\gB_{11}$ est un localisé de 
$\gk[\ux^{P_1}, \ux^{P_2}, \ux^{P_3}, \ux^{P_4}, \ux^{P_5}]$, qui
est un anneau de \pols à 5 \idtrs sur $\gk$, donc \iclz.

\emph {3c.}
On définit $\gB_{ij}$ de façon à ce qu'il soit contenu dans $\gk[\,\ux^M
\,\vert\, M \in E,\ m_{ij} \ge 0\,]$. Par exemple, pour $(i,j) = (3,1)$, les matrices
$P_k$ ayant un coefficient nul en position $(3,1)$ sont celles autres que
$P_3$, $P_5$, ce qui conduit à la \dfn de $\gB_{31}$:

\snic {
\gB_{31} = \gk[\ux^{P_3}, \ux^{P_5}]
[\ux^{\pm P_1}, \ux^{\pm P_2}, \ux^{\pm P_4}, \ux^{\pm P_6}].
}

%\sni
On a alors l'\egt $\gB = \bigcap_{i,j} \gB_{ij}$, et comme les $\gB_{ij}$
sont tous \icl de même corps de fractions $\Frac\gB$, l'anneau $\gB$
est \iclz.

%%%%%%%%%%%%%%%%%%%%%%%%%%%%%%%%%%%%%%%%%%%%%%%%%%%%%%%%%%%%%%%%%%%%%%%%%%%

\exer{exoIdentiteDiscriminantale} 
\emph {2.}
On écrit $U = \tra{M} M$ et on prend le \deterz.\\
Cela donne $\Disc_\gk \gA = \disc(f)^{n!/2}$ à partir de $\det(M) = \delta^{n!/2}$. 
Réciproquement,
puisqu'il s'agit d'\idas dans $\ZZ[\ux]$, l'\egt $(\det M)^2=(\delta^{n!/2})^2$
implique $\det M=\pm\delta^{n!/2}$. 

\emph {3.}
Dans le \tho \ref{theoremAdu3}, ne supposons pas $f$ \spl sur  $\gC$.
Par hypothèse, \hbox{on a  $\varphi(f)(T) =
\prod_{i=1}^n(T- u_i)$}.  Avec $\gA = \gk[\xn] =
\Adu_\gk(f)$, on a alors un morphisme de $\gC$-\algs $\Phi : \gC\te_\gk\!\gA \to
\gC^{n!}$ qui réalise $1 \te x_i \mapsto (u_{\sigma(i)})_{\sigma\in
\rS_n}$. \\
La $\gk$-base canonique $\cB(f)$ de $\gA$ est une $\gC$-base de
$\gC\te_\gk\gA$ et  la matrice de $\Phi$ pour cette base (au départ) et pour
la base canonique de $\gC^{n!}$ (à l'arrivée) est la matrice $M$ ci-dessus
où $x_i$ est remplacé par $u_i$.  On en déduit que  
$\Phi$ est un \iso \ssi $\varphi\big(\disc(f)\big)\in\gC\eti$, i.e. si
$f$ est \spl sur $\gC$.

 Finalement faisons seulement l'hypothèse qu'une \alg $\varphi:\gk\to \gC$  
diagonalise~$\gA$. Cela signifie que l'on donne $n!$
caractères $\Adu_{\gC,\varphi(f)}\to\gC$ qui mis ensemble
donnent un \iso de $\gC$-\algs de $\Adu_{\gC,\varphi(f)}$ sur $\gC^{n!}$.
\\
 Puisqu'il existe un caractère $\Adu_{\gC,\varphi(f)}\to\gC$, le \pol $\varphi(f)(T)$  se factorise complètement dans $\gC$.
\\
Enfin le \discri de la base canonique de $\Adu_{\gC,\varphi(f)}$ est 
$\varphi\big(\disc(f)\big)^{n!/2}$
et le \discri de la base canonique de $\gC^{n!}$ est $1$. 
Donc, $f$ est \spl sur $\gC$.

%%%%%%%%%%%%%%%%%%%%%%%%%%%%%%%%%%%%%%%%%%%%%%%%%%%%%%%%%%%%%%%%%%%%%%%%%%%

\exer{exoTnADU} 
On a $\gA = \aqo {\gk[\Xn]}{S_1, \ldots, S_n}$ où $S_1$, \ldots, $S_n$ sont
les $n$ fonctions \smqs \elrs de $(\Xn)$; l'\id $\gen {S_1, \ldots, S_n}$ étant
\hmgz, la \klg $\gA$ est graduée (par le degré). On note $\gA_d$ sa
composante \hmg de degré $d$ et $\fm = \gen {\xn}$; on a donc
$\gA = \gA_0 \oplus \gA_1 \oplus \gA_2 \oplus \dots$ avec $\gA_0 = \gk$
et:

\snic {
\fm^d = \gA_d \oplus \gA_{d+1} \oplus \dots, \quad
\fm^d = \gA_d \oplus \fm^{d+1}.
}

%\sni
Puisque $x_i^n = 0$, on a $\fm^{n(n-1)+1} = 0$, donc $\gA_d = 0$ pour $d \ge
n(n-1)+1$. On rappelle la base  $\cB(f)$ de $\gA$, formée des $x_1^{\alpha_1} \ldots x_n^{\alpha_n}$ avec $0 \le \alpha_i <
n-i$. Pour tout~$d$, la composante \hmg $\gA_d$ de degré $d$ est un \kmo
libre dont une base est l'ensemble des $x_1^{\alpha_1}
\ldots x_n^{\alpha_n}$ avec $0 \le \alpha_i < n-i$ et $|\alpha| = d$.  Le cardinal de cette base est le \coe de degré
$d$ dans le \pol $S(t) \in \ZZ[t]$:

\snic {
S(t) = 1 (1+t) (1 + t + t^2) \cdots (1 + t + \cdots + t^{n-1}) =
\prod_{i=1}^n \frac {t^i-1} {t-1}.
}

%\sni
En effet, un multi-indice $(\alpha_1, \ldots, \alpha_n)$ tel que $0 \le
\alpha_i < n-i$ et $|\alpha| = d$ s'obtient en choisissant un \mom
$t^{\alpha_n}$ du \pol $1 + t + \cdots + t^{n-1}$, un \mom $t^{\alpha_{n-1}}$
du \pol $1 + t + \cdots + t^{n-2}$ et ainsi de suite, le produit de ces \moms
étant~$t^d$. On obtient ainsi la série d'Hilbert-Poincaré $S_\gA(t)$ de
$\gA$:

\snic{
S_\gA(t) \eqdefi \som_{i=0}^{\infty}\dim_\gk \gA_d\ t^d 
\eqdf {\rm ici} \som_{0 \le \alpha_i < n-i} t^{|\alpha|} = S(t).
}

%\sni
Le \pol $S$ est un \polu de degré $e= 1 + \cdots + n-1= n(n-1)/2$.
On a $S(1) = n!$, conforme à $S(1) = \dim_\gk\gA$.

\Deuxcol{.8}{.15}
% premiere colonne
{Variante. On pose $\gB = \gk[S_1, \ldots, S_n] \subset \gC = \gk[\Xn]$. \\
Alors $\gC$
est un $\gB$-module libre de base les $\uX^\alpha = X_1^{\alpha_1}
\cdots X_n^{\alpha_n}$ avec $0 \le \alpha_i < n-i$. Cette base est au dessus
de la base $\cB(f)$ de~$\gA$ sur $\gk$ si l'on considère que l'on a un
diagramme commutatif où chaque flèche verticale est une réduction
modulo $\gen {S_1, \ldots, S_n}$.  
}
% deuxieme colonne
{~

$\xymatrix @R = 0.5cm @C = 0.8cm
{\gB\ar@{->>}[d]\ar@{->}[r] &\gC\ar@{->>}[d]\\
\gk\ar@{->}[r] &\gA \\}$ 
}

On écrit $\gC = \bigoplus_{\alpha} \gB
\uX^\alpha$, avec le décalage $S_{\gB\uX^\alpha}(t) =
t^{|\alpha|}S_\gB(t)$, on a l'\egt suivante entre les séries d'Hilbert-Poincaré:

\snic {
S_\gC = S_\gA\ S_\gB  \qquad \hbox {avec} \qquad
S_\gA = \som_{0 \le \alpha_i < n-i} t^{|\alpha|}.
}

%\sni
Or il est facile de voir que

\snic {
S_\gC(t) = \fraC {1} {(1-t)^n}, \;
S_\gB(t) = \prod_{d=1}^n \fraC {1} {1 - t^d}, 
\; \hbox { et donc } \,
S_\gA(t) = \frac {S_\gC} {S_\gB} = \prod_{d=1}^n \fraC {1-t^d} {1-t},
}

%\sni
ce qui donne de nouveau le résultat pour $S_\gA$.

Passons maintenant aux puissances de l'\id $\fm$.  \\
Soit $\varphi :
\gA\twoheadrightarrow \gk$ le caractère $x_i \mapsto 0$ de noyau $\fm = \gen
{\xn}$. \\
On a $\gA=\kxn=\gk\oplus\fm$, $\fm\subseteq\DA(0)\subseteq \Rad(\gA)$
et pour $z \in \gA$, 

\snic{z\in\Ati\iff\varphi(z)\in\gk\eti\iff z\in\gk\eti\oplus\fm.}

%\sni
On a $\DA(0)=\rD_\gk(0)\oplus\fm$, $\Rad(\gA)=\Rad(\gk)\oplus\fm$.\\
Puisque $\fm \subseteq \Rad(\gA)$ est \tfz, on a 

\snic{\fm^d = \fm^{d+1} \iff \fm^d
= 0}

%\sni
(lemme \ref{lemLocaliseFini}), ce qui
équivaut, puisque $\fm^d = \gA_d \oplus \fm^{d+1}$, à $\gA_d = 0$. On en
déduit que $\fm^{e+1} = 0$.

Une remarque: si $\gk$ est local, alors $\gA$ \egmtz, et
$\Rad\gA=\varphi^{-1}(\Rad\gk)$.

%%%%%%%%%%%%%%%%%%%%%%%%%%%%%%%%%%%%%%%%%%%%%%%%%%%%%%%%%%%%%%%%%%%%%%%%%%%

\exer{exoNilIndexInversiblePol} 
On considère l'anneau de \pols $\gC =\gk[\an,b_1,\dots,b_m]$, on pose $f(X) = 1 + \sum_{i=1}^n a_i X^i$, $g(X) = 1 + \sum_{j=1}^m b_j
X^j$, et $\fc=\rc_\gC(fg - 1)$. On affecte
à  $a_i$ le poids $i$ et à  $b_j$ le poids $j$. Le
\coe de degré $k$ de $fg-1$ est \hmg de degré $k$, donc l'\id
$\fc$ est \hmgz. \\
On note $\gC' = \gC\sur\fc$. Cette \klg $\gC'$ est graduée
via le poids ci-dessus et on doit montrer que $\gC'_d = 0$ pour $d >
nm$. Il est clair que $\gC'_d = 0$ pour $d$ assez grand. 
On va déterminer la série d'Hilbert-Poincaré $S_{\gC'}$ de $\gC'$ (qui 
est ici un \polz):

\snic {
S_{\gC'}(t) \eqdefi \sum_{d \ge 0} \dim_\gk \gC'_d\, t^d =
{ \prod_{d=1}^{n+m} (1-t^d) \over \prod_{i=1}^{n} (1-t^i)\prod_{j=1}^{m} (1-t^j)}\;.
}

%\sni
Pour démontrer cette \egtz, on réalise $\gC$ et $\gC'$ d'une autre manière. \\
On
considère $n+m$ \idtrs $(\Xn, \Ym)$, et on note $(\an)$ les 
fonctions \smqs \elrs de $(\Xn)$, et $(b_1,\dots,b_m)$ les fonctions \smqs \elrs de 
$(\Ym)$. Puisque

\snuc {
\prod_{i=1}^n (T+X_i) \prod_{j=1}^m (T+Y_j) =
(T^n + a_1 T^{n-1} + \cdots + a_n) (T^m + b_1 T^{m-1} + \cdots + b_m), 
}

%\sni
on voit, en posant $a_0=b_0 = 1$, que $\sum_{i+j=d} a_ib_j$ est la
$d$-ième fonction \smq \elr de $(\Xn, \Ym)$. Comme $(\an, b_1,\dots,b_m)$
sont \agqt indépendants sur $\gk$, on peut considérer que
$\gC$ est la sous-\alg graduée suivante:

\snic {
\gC = \gk[\an, b_1,\dots,b_m]  \subset \gD = \gk[\Xn, \Ym],
}

%\sni
et que l'\id $\fc$ de $\gC$ est engendré par les $n+m$ sommes
$\sum_{i+j=d} a_ib_j$, qui sont les  fonctions \smqs \elrs de $(\Xn, \Ym)$.
\\
L'\alg $\gD$ est libre sur $\gC$ de rang $n!m!$,
comme pour une double \aduz.  Plus \prmtz, voici des bases. \\
Les $\uX^\alpha =
X_1^{\alpha_1} \cdots X_n^{\alpha_n}$ pour $0 \le \alpha_i < n-i$ forment une
base de $\gk[\uX]$ sur $\gk[\ua]$. \\
Les $\uY^\beta = Y_1^{\beta_1}
\cdots Y_m^{\beta_m}$ avec $0 \le \beta_j < m-j$ forment une base de
$\gk[\uY]$ sur $\gk[\ub]$. 
\\
Ainsi, les $\uX^\alpha \uY^\beta$
forment une base de $\gD = \gk[\uX,\uY]$ sur $\gC = \gk[\ua,\ub]$.\\
Enfin, par l'\eds $\gC\to\gC' = \gC\sur\fc$, les $\ux^\alpha \uy^\beta$ forment une base \hbox{de $\gD' = 
\gD\sur{\fc\gD} = \gk[\ux,\uy]$ sur $\gC'$}.

\Deuxcol{.77}{.18}
% premiere colonne
{
On dispose d'un diagramme commutatif où chaque flèche verticale est une réduction
modulo $\fc$.  Il s'agit de déterminer la série d'Hilbert-Poincaré
$S_{\gC'}$ de $\gC'$ sachant que l'on connait celles \hbox{de $\gD'$, $\gC$} et~$\gD$
(car $\gC$ et $\gD$ sont des anneaux de \polsz, \hbox{et $\gD'$} est l'\adu de $T^{n+m}$
sur~$\gk$).
}
% deuxieme colonne
{\vspace{-3mm}
$$\xymatrix @R = 0.5cm @C = 0.8cm
{\gC\ar@{->>}[d]\ar@{->}[r] &\gD\ar@{->>}[d]\\
\gC'\ar@{->}[r] &\gD' \\}$$ 
}

On peut conduire les calculs de la manière simple 
suivante. \\
On écrit $\gD = \bigoplus_{\alpha, \beta} \gC\, \uX^\alpha\uY^\beta$,
donc 

\snic {  
S_{\gD}(t) = F(t) S_{\gC}(t)  \quad \hbox {avec} \quad
F(t) = \sum_{\alpha, \beta} t^{|\alpha| + |\beta|} =
\sum_{\alpha} t^{|\alpha|} \sum_{\beta} t^{|\beta|},
}

%\sni
et l'on a également $S_{\gD'}(t) = F(t) S_{\gC'}(t)$. On a vu
dans l'exercice \ref{exoTnADU} que:

\snic {
F(t) = \prod_{i=1}^n \fraC {1-t^i} {1-t}\;\prod_{j=1}^m \fraC {1-t^j} {1-t},
\qquad
S_{\gD'}(t) = \prod_{d=1}^{n+m} \fraC {1-t^d} {1-t}.
}

%\sni
On note alors $S_d(t) = (1-t^d)/(1-t)$. C'est un \pol de degré $d-1$ 
\hbox{et $S_d(1) = d$}. On a donc obtenu
$$
S_{\gC'}(t) = {S_1 S_2 \cdots S_{n+m} \over 
S_1 S_2 \cdots S_{n}\,  S_1 S_2 \cdots S_{m}}, 
$$
avec

\snic {
\deg S_{\gC'} = {(n+m-1)(n+m) - (n-1)n - m(m-1) \over 2} = nm.
}

%\sni
Ainsi, comme souhaité, $\gC'_k = 0$ pour $k > nm$. \\
\`A noter que
$\dim_\gk \gC' = S_{\gC'}(1) = {n+m \choose n}$. 

%%%%%%%%%%%%%%%%%%%%%%%%%%%%%%%%%%%%%%%%%%%%%%%%%%%%%%%%%%%%%%%%%%%%%%%%%%%
\exer{exoADUpthRoot} 
Pour chaque $i\in\lrbp$ la restriction $\varphi : \gk[x_i] \to\gk[\alpha]$ est un \isoz.  
Considérons l'\id 

\snic{\fm = \gen {x_i - x_j, i,j \in\lrbp}= \gen {x_1 - x_i, i \in\lrb{2..p}}.}

%\sni
Alors $\gA=\gk[x_1]\oplus \fm$, d'où $\fm=\Ker\varphi$.
\\
En fait on peut voir $\gA$ comme l'\adu $\Adu_{\gk[x_1],g}$
pour le \pol $g(T)=f(T)/(T-x_1)=(T-x_1)^{p-1}$ sur l'anneau $\gk[x_1]$
 ce qui
nous ramène à l'exercice \ref{exoTnADU}.
En particulier: 

\snic{\fm^{1+(p-1)(p-2)/2}=0$, 
$\DA(0)=\rD_{\gk[x_1]}(0)\oplus \fm$ et $\Rad(\gA)=\Rad(\gk[x_1])\oplus \fm.}

%%%%%%%%%%%%%%%%%%%%%%%%%%%%%%%%%%%%%%%%%%%%%%%%%%%%%%%%%%%%%%%%%%%%%%%%%%%

\exer{exoA5GaloisGroup} 
\emph {1.}
L'opération $b \aff 5b$, $c \aff 4c$ a pour but de remplacer $4^4b^5 +
5^5c^4$ par $4^4 5^5 (b^5 + c^4)$; en imposant $c = \pm b$, on obtient $4^4
5^5 b^4 (b + 1)$ qu'il est facile de rendre carré en imposant $5(b + 1) =
a^2$.  Pour éviter le dénominateur $5$, on impose plutôt $5(b + 1) =
(5a)^2$, i.e. $b = 5a^2-1$.

\emph {2.}
Pour $a \in \QQ\sta$, le \pol $f_a(T) \in \QQ[T]$ est \splz.  Modulo les
petits nombres premiers on trouve les \dcns de $f_1(T) = T^5 + 20 T + 16\vep$,
\hbox{avec $\vep \in \{\pm 1\}$}, en facteurs \irds suivantes:
\[ 
\begin{array}{rccl} 
\mod &  2 &:& T^5 \\[1mm] 
\mod &  3 &:& f_1(T) \\[1mm] 
\mod &  5 &:& (T + \vep)^5\\[1mm] 
\mod &  7 &:& (T+2\vep)(T+3\vep)(T^3 + 2\vep T^2 + 5T + 5\vep)  
 \end{array}
\] 
Le résultat modulo $3$ prouve que $f_1(T)$ est \ird sur $\ZZ$.  Son groupe
de Galois~$G$ est un sous-groupe transitif de $\rA_5$ qui contient un
$3$-cycle (vu la réduction modulo $7$). Ceci implique $G = \rA_5$. En effet, un
sous-groupe transitif de $\rS_5$ contenant un $3$-cycle est
égal à $\rS_5$ ou $\rA_5$. En ce qui 
concerne $\QQ(a)$ comme corps de base, le \pol $f_a(T)$
est \ird dans $\QQ[a][T]$ puisque sa réduction modulo
$a = 1$ l'est dans $\QQ[T]$. Donc il est \ird dans $\QQ(a)[T]$.
En utilisant le fait que son discriminant est un carré et la réduction
modulo $a=1$, on obtient que son groupe de Galois est $\rA_5$.
\\
\Llec pourra se poser la question suivante: est-ce que pour tout $a \in \ZZ\setminus \{0\}$, le \pol $f_a(T)$
est \ird de groupe de Galois $\rA_5$?

 \emph{Expérimentation posssible}.\\
Voici la répartition des types de permutation
des sous-groupes transitifs de $\rS_5$% 
%où l'on rappelle que
%le type d'une permutation est la suite des longueurs de ses orbites
%rangée par ordre décroissant
. 
\\
Pour les 7 types qui apparaissent
dans $\rS_5$, on utilise les notations suivantes:

\snuc {
t_1 = (1^5),\, t_2 =(2,1^3),\, t_{22} = (2^2,1),\, t_3 = (3,1^2),\, 
t_{3,2} = (3,2),\, t_4 = (4,1),\ t_5 = (5).
}

%\sni
Ainsi $t_{22}$ est le type des double-transpositions, $t_3$ celui des $3$-cycles,
etc{\dots} La table annoncée:
\vskip2pt
\snic {
\begin {array} {cccccc}
\hline\\[-2.5mm]
G &\rC{5} &\ASL_1(\FF_5) &\AGL_1(\FF_5) &\rA_5 &\rS_5
   \\[1mm]\hline\\[-2.5mm]
\#G &5 &10 &20  &60 &120
   \\[1mm]\hline\\[-2.5mm]
&t_1^1\ t_5^4 &t_1^1\ t_{22}^5\ t_5^4 &t_1^1\ t_{22}^5\ t_4^{10}\ t_5^4 &
    t_1^1\ t_{22}^{15}\ t_3^{20}\ t_5^{24} &
    t_1^1\ t_2^{10}\ t_{22}^{15}\ t_3^{20}\ t_{32}^{20}\ t_4^{30}\ t_5^{24}
\\[1mm]\hline
\end {array}
}
\vskip2pt
%\sni

Par exemple sur la dernière ligne, sous $\rA_5$, $t_1^1\ t_{22}^{15}\ t_3^{20}\ t_5^{24}$ signifie que $\rA_5$ contient l'identité, $15$ double-transpositions, $20$
$3$-cycles et $24$ $5$-cycles ($1 + 15 + 20 + 24 = 60$).
\\
\Llec pourra tester expérimentalement le \tho de densité de Cebotarev
à l'aide d'un \sys de calcul formel. Il faut examiner 
la \fcn de $f_1(T)$ modulo \gui{beaucoup} de premiers $p$ et
comparer la répartition obtenue des types de \fcn à celle des
types de permutation de $\rA_5$.  
\\
L'auteur de
l'exercice a considéré les $120$
% --- deux fois le cardinal de $\rA_5$ ---
premiers nombres premiers --- autres que $2$ et $5$ qui divisent $\disc(f_1)$ ---
et son logiciel a trouvé la répartition suivante:

\snic {
t_{22}^{33}\ t_3^{38}\ t_5^{49}
}

%\sni

Ceci signifie que l'on a trouvé $33$ fois une \fcn de type $t_{22}$ (2
facteurs \irds de degré $2$, 1 facteur \ird de degré $1$), $38$ fois une
\fcn de type $t_{3}$ et $49$ fois une \fcn de type $t_{5}$ (pas de
factorisation de type~$t_1$). Répartition
à comparer avec celle de $\rA_5$. Quant au type $t_1$,
le plus petit premier $p$ pour lequel $f_1(T) \bmod p$ est totalement
décomposé est $p = 887$. Enfin, en traitant $1200$ premiers au lieu
de $120$, on trouve la répartition:

\snic {
t_1^{16}\ t_{22}^{304}\ t_3^{428}\ t_5^{452}
}

%%%%%%%%%%%%%%%%%%%%%%%%%%%%%%%%%%%%%%%%%%%%%%%%%%%%%%%%%%%%%%%%%%%%%%%%%%%
\exer{propGaloisJordan}
\emph{1a.} Il faut montrer que $\gen{y_1-a}+\gen{y_j-a}=\gen{1}$ pour $j\in \lrb{2..r}$. Par exemple dans le quotient $\aqo{\gB}{y_1-a,y_2-a}$ le \pol $g(T)=\prod(T-y_j)$ a deux facteurs égaux à $T-a$ ce qui implique $g'(a)=0$. Comme
$g'(a)$ est \iv par hypothèse (ce qui reste vrai dans un quotient), on a
bien $0=1$ dans le quotient.

 \emph{1b.} On voit facilement que $H=\St(y_1)$.
Donc $H$ opère sur $\so{\beta(y_2),\ldots,\beta(y_r)}$. \\
Or $g(T)/(T-y_1)=\prod_{j=2}^r(T-y_j)$ dans $\gB$, donc $g(T)/(T-a)=\prod_{j=2}^r\big(T-\beta(y_j)\big)$
dans $\gC$.\perso{ce n'est pas clair si l'orbite de $\beta(y_2)$
sous $H$ contient tous les $\beta(y_j)$ avec $j\neq1$}

 \emph{2a.} Il est clair que $y_1-a$ est un \dvz dans
$\gB$. Un \id galoisien minimal $\fc$ contenant $\gen{y_1-a}$
est obtenu en ajoutant le plus possible de conjugués de $\gen{y_1-a}$ sous la condition de ne pas atteindre l'\id $\gen{1}$.
L'\id $\fc$ est donc de la forme $\gen{y_j-a\,|\,j\in J}$
pour une partie $J$ de $\lrbr$. Il reste à voir que les $j$ tels que
$y_j-a\in\fc$ sont au nombre de $k$.
Or pour tout indice $j$, l'\elt $y_j-a$
est nul ou \iv modulo $\fc$. Puisque $g(T)=\prod_j\big(T-\beta(y_j)\big)$,
et puisque $a$ est un zéro de multiplicité $k$ de $g$, le nombre de
$j$ tels que $\beta(y_j)=a$ est égal à $k$ \linebreak 
(écrire $g(a)=g'(a)=\cdots=g^{(k-1)}(a)=0$ et $g^{(k)}(a)$ \ivz).

 \emph{2b.} On raisonne comme pour
\emph{1b.}

 \emph{3.} Le quotient de Galois
$\gC=\gB/\fc$ est obtenu avec son groupe $H=\St_G(\fc)$.
Par hypothèse $\ov {y_1}\in\Fix(H)$ donc
$\ov{y_1}\in\gK$. Notons $a$ l'\elt de $\gK$ en question.
Dans~$\gC$ on a $g(T)=\prod_j(T-\ov{y_j})$, donc ${g(a)}=0$.
Enfin $\gK $ s'identifie à son image dans~$\gC$.

\medskip
\exl Voici un exemple avec $\deg f=6$. On demande à {\tt Magma} de calculer le \polmin de $y=x_4+x_5x_6$, puis de le factoriser. Si $g$ est le premier facteur,
$z=g(y)$ est un \dvzz. On lance l'\algo \ref{alg2idmgal} avec~$z$.
On obtient donc les nouvelles approximations
du \cdr et du groupe de Galois en traitant la bizarrerie \gui{$z$
est \dvzz}, mais on peut constater a posteriori que $z$ est de multiplicité $6$ dans sa résolvante et que $\gen{z}$ est galoisien.
}

%\vspace{-.2cm}
%:     verbatim
%{\small
\label{exemple3Galois}
{\footnotesize
\begin{verbatim}
f:= T^6 - 3*T^5 + 4*T^4 - 2*T^3 + T^2 - T + 1;
y:=x4+x5*x6; pm:=MinimalPolynomial(y);
  T^60 - 46*T^59 + 1035*T^58 - 15178*T^57 + 163080*T^56 + ... + 264613
Factorization(pm);
    <T^6 - 4*T^5 + 8*T^4 - 6*T^3 + T + 1, 1>,
    ...
z:=Evaluate(T^6 - 4*T^5 + 8*T^4 - 6*T^3 + T + 1,y);
20*x4^3*x5^3*x6^3 - 15*x4^3*x5^3*x6^2 - 15*x4^3*x5^2*x6^3 +
  11*x4^3*x5^2*x6^2 + 2*x4^3*x5^2*x6 + 2*x4^3*x5*x6^2 + x4^3*x5*x6 - ...
// z divise 0, on calcule le nouveau quotient de Galois
Affine Algebra of rank 6 over Rational Field
Variables: x1, x2, x3, x4, x5, x6
Quotient relations:
  x1 + x2 + x3 - x6^5 + 2*x6^4 - x6^3 - x6^2 - 1,
  x2^2 + x2*x3 - x2*x6^5 + 2*x2*x6^4 - x2*x6^3 - x2*x6^2 - x2 + x3^2 -
      x3*x6^5 + 2*x3*x6^4 - x3*x6^3 - x3*x6^2 - x3 + x6^5 - 2*x6^4 +
      x6^3 + x6^2,
  x3^3 - x3^2*x6^5 + 2*x3^2*x6^4 - x3^2*x6^3 - x3^2*x6^2 - x3^2 +
      x3*x6^5 -  2*x3*x6^4 + x3*x6^3 + x3*x6^2 - x6^5 + 2*x6^4 - x6^3 -
      x6^2 + 1,
  x4 + x5 + x6^5 - 2*x6^4 + x6^3 + x6^2 + x6 - 2,
  x5^2 + x5*x6^5 - 2*x5*x6^4 + x5*x6^3 + x5*x6^2 + x5*x6 - 2*x5 -
      x6^4 + 2*x6^3 - x6^2 - x6,
  x6^6 - 3*x6^5 + 4*x6^4 - 2*x6^3 + x6^2 - x6 + 1
Permutation group G2 acting on a set of cardinality 6
Order = 72 = 2^3 * 3^2
    (1, 4)(2, 5)(3, 6)
    (1, 2)
    (2, 3)
Degree(MinimalPolynomial(z)); 55
#Orbite(z,G); 60
\end{verbatim}
}
{\small

%%%%%%%%%%%%%%%%%%%%%%%%%%%%%%%%%%%%%%%%%%%%%%%%%%%%%%%%%%%%%%%%%%%%%%%%%%%
\exer{propRvRelDec}
On note que les $y_i-y_j$ pour $i\neq j$ sont \ivsz, et que ceci reste vrai
dans tout quotient de Galois.

%%%%%%%%%%%%%%%%%%%%%%%%%%%%%%%%%%%%%%%%%%%%%%%%%%%%%%%%%%%%%%%%%%%%%%%%%%%

}% fin des solutions d'exos

%:   ---- Section*{references}-----------
\Biblio

Le \thref{thpolcohfd} sur le fait qu'un anneau de \pols sur un anneau \zedr
est \coh \fdi admet une extension remarquable pour les \adps \cohs \fdisz: voir \cite{Yengui} et \cite{DVY2015}. 

Les versions que nous avons données du \nst \gui{sans  côture \agqz} se trouvent 
sous une forme voisine
dans \cite[VIII.2.4,VIII.3.3]{MRR}.

La difficulté intrinsèque du \pb de l'\iso de deux clôtures \agqs d'un corps
est illustrée dans \cite[Sander, Theorem~26]{sand}, qui montre que, en présence du tiers exclu mais en l'absence d'axiome du choix dépendant, il est impossible de démontrer dans $\ZF$ que deux  clôtures \agqs de $\QQ$ sont isomorphes.

Le traitement de la théorie de Galois basé sur les quotients de Galois de
l'\adu remonte au moins à Jules Drach \cite[1898]{Drach} et à Ernest Vessiot \cite[1904]{Vessiot}. Voici un extrait de l'introduction de ce dernier article, qui parle dans le langage de l'époque des quotients de Galois de l'\aduz:

\gui{\'Etant donnée une \eqn \agqz, que l'on considère comme remplacée par le \sys $(S)$ des relations entre les racines $\xn$ et les \coesz,
on étudie d'abord le \pb fondamental suivant: \emph{Quel parti peut-on tirer de la connaissance de certaines relations $(A)$ entre $\xn$, en n'employant que des opérations rationnelles?} Nous montrons que l'on peut déduire du \sys $(S,A)$ un \sys analogue, dont le \sys $(S,A)$ admet toutes les solutions, et qui est, comme nous le disons, \emph{automorphe}: ce qui veut dire que ses diverses solutions se déduisent de l'une quelconque d'entre elles par les substitutions d'un groupe $G$, qui est dit \emph{le groupe associé au \sysz}, ou simplement le \emph{groupe du \sysz}. On remarquera que $S$ est déjà un \sys automorphe, ayant le groupe \gnl pour groupe associé. Dès lors, si l'on se place du point de vue de Galois,
 \dots\ 
on voit que l'on peut se limiter à ne considérer que des \syss $(S,A)$ rationnels et automorphes.}

\medskip 
L'\adu est traitée de manière assez détaillée dans le chapitre 2
du livre   \cite[1989]{PZ}.

Parmi les bons exposés modernes qui exposent toute la théorie classique de Galois, on peut citer \cite{Tignol} et~\cite{CoxGal}.

La \gui{théorie de Galois dynamique} exposée en détail dans ce chapitre
est présentée dans \cite[Díaz-Toca]{DiJSC} et \cite[Díaz-Toca\&al.]{DiLQ,DiL10}. 

Concernant le \thrf{theoremAdu1} sur les points fixes de $\Sn$ dans l'\aduz, le cas \gui{$f$ \splz}  fait partie du folklore. On le trouve avec
une preuve voisine de celle donnée ici
dans la thèse de Lionel Ducos \cite{Du}. Nous en avons donné une autre preuve dans le \thref{thResolUniv} pour le cas des \cdisz.
Le raffinement que nous donnons se trouve dans \cite{DiLQ}, il est inspiré de \cite{PZ} (voir le \tho 2.18 page~46, le corolaire 3.6 page 49
 et la remarque qui le suit, page~50).

Le \thrf{thidGTri}, publié dans \cite{DiLQ} sous une hypothèse restrictive,  \gns
un résultat donné séparément dans le cas de l'\adu sur un corps par
L. Ducos \cite{Duc} et par P. Aubry et  A. Valibouze~\cite{AV}. Notre méthode
de preuve se rapproche plus de celle de L. Ducos, mais elle est différente
car le cadre est plus \gnlz: nous avons à la base un
anneau commutatif arbitraire.

Une  version voisine du \thref{theoremAdu3} se trouve dans
\cite[lemme~II.4.1]{Du}.

\smallskip \rdb\label{Jordan}\label{soicher}
Concernant les méthodes explicites de calcul de groupes de Galois sur $\QQ$
récemment développées en calcul formel on pourra consulter 
\cite[Geissler\&Klüners]{GeKl}.

La méthode modulaire, popularisée par van der Waerden, 
est due à Dedekind
(lettre  adressée à  Frobenius le 18 juin 1882, voir \cite[Brandl]{Brandl}).

Les méthodes de Stauduhar \cite{Sta} et Soicher-McKay \cite{SoM} sont basées sur  des calculs de résolvantes et sur la connaissance des sous-groupes transitifs des groupes $\Sn$. Ceux-ci ont été tabulés jusqu'à $n=31$ \cite[Hulpke]{Hulpke}. 
Dans la plupart des \algos existants
le calcul détermine le groupe de Galois d'un \pol \irdz, sans calculer le corps des racines.\\
Voir cependant \cite[Klüners\&Malle]{KlMa} et \cite[Valibouze\&al.]{AV,ORV,Val}.

Citons par ailleurs le résultat remarquable \cite[Landau\&Miller]{LaMi} 
de calculabilité en temps \poll
concernant la résolubilité par radicaux.

\smallskip
Alan Steel \cite{steel,steel2} s'est inspiré de D5 pour implémenter
une très performante clôture algébrique \gui{dynamique} de $\QQ$ en {\tt Magma}. L'efficacité tient à ce qu'il n'utilise pas d'\algo de factorisation des \pols de $\ZZ[X]$, ni de représentation des extensions finies au moyen d'\elts primitifs.
Il utilise néanmoins des \algos de factorisation modulo $p$ pour contrôler le processus. Le processus est dynamique dans la mesure où la clôture construite progressivement dépend des questions de l'utilisateur.
L'auteur ne donne cependant pas (et il ne pourrait pas le faire dans le cadre qu'il se fixe) une implémentation
du corps de racines d'un \pol (disons \spl pour simplifier) sur un corps
\gui{général}. 

\smallskip 
Pour le \sys de calcul formel {\tt Magma}, voir \cite[Bosma\&al.]{Bosma,Cannon}.

\newpage \thispagestyle{CMcadreseul}
\incrementeexosetprob

%:        %%%%%%%%%%%%%%%%%%%%%%%%%%%%%%%%%%%%
%:        %%%%%%%%%%%%%%%%%%%%%%%%%%%%%%%%%%%%
%---- Chapitre  {Modules plats}------------     
\chapter{Modules plats} 
\label{chap mod plats}\relax
%--------------------
\minitoc
\newpage	
\Intro
%-----------------------------------------

\begin{flushright}
{\em Chers \eltsz, \\
si vous n'êtes pas libres, \\
ce n'est pas de ma faute.}\\
Un module plat. \\
\end{flushright}

La platitude est une notion fondamentale de l'\alg commutative,
introduite par Serre dans~\cite{Serregaga}.

Dans ce chapitre nous introduisons les notions de module plat, d'\alg plate et d'\alg \fpte et démontrons quelques unes des \prts essentielles de ces objets.

Un anneau intègre dont les \itfs sont plats est appelé un \ddpz.
C'est une autre notion fondamentale de l'\alg commutative. 
Elle sera seulement introduite ici et sera développée dans 
le chapitre~\ref{ChapAdpc}.

\goodbreak
%--- Sec {Premieres proprietes}--   
\section{Premières propriétés}\label{secPlatDebut} 

\vspace{5pt}
\subsec{Définition, et \prts de base}
Nous donnons une \dfn de nature \elr et développerons 
plus loin le rapport avec l'exactitude du foncteur $M\otimes \bullet$. 
%--- Definition{def.plat}-------
\begin{definition} 
\label{def.plat}\relax On considère un \Amo $M$.
\begin{enumerate}
\item Une \emph{\syzy dans $M$} est donnée par $L\in \Ae {1\times n}$ et $X\in M^{n\times 1}$ qui vérifient $LX=0$.
\item On dit que \emph{la \syzy $LX=0$ s'explique dans~$M$} si l'on trouve un vecteur $Y\in M^{m\times 1}$ et une matrice $G\in \Ae {n\times m}$ qui vérifient:
%---  equation eqdef.plat -------
\begin{equation}\label{eqdef.plat}\relax
LG=0\quad {\rm et } \quad GY = X\,.
\end{equation}
%---------------------end equation--------------
%
\item Le \Amo $M$ est appelé un {\em  \mplz}\index{module!plat}\index{plat!module} si  toute \syzy  
dans $M$ s'explique dans $M$.
(En langage intuitif: s'il y a une \syzy entre \elts de $M$ ce n'est pas la faute 
au module.) 
\end{enumerate}
\end{definition}
%--- end-definition---------------

\rems 1) Dans les points \emph{1} et \emph{2} le symbole $0$ est précisé  implicitement par le contexte. En \emph{1} il s'agit de $0_M$, tandis qu'en \emph{2}
%:2018    $0_{\Ae {m\times 1}}$  remplacé
il s'agit de $0_{\Ae {1\times m}}$. 

2) Dans le point \emph{2}, lorsque l'on dit que la \syzyz~\hbox{$LX=0$} s'explique dans~$M$, on signifie
que l'explication \gui{ne touche pas à $L$}.
Par contre,  les \egts données par l'\eqn matricielle $LG=0$ ont lieu dans $\gA$ et non dans $M$.
\eoe

\medskip 
\exls 1) Si $M$ est libre fini\footnote{Ou plus \gnlt si $M$ est librement engendré par un ensemble discret, \cad $M\simeq \gA^{(I)}=\bigoplus_{i\in I}\gA$ avec $I$ discret. Pour une autre \gnnz, 
voir l'exercice~\ref{propfreeplat}.}, il est plat: si $LX=0$,
on écrit~$X=GY$ avec un vecteur
colonne $Y$ qui forme une base, et $LX=0$ implique~$LG=0$.

 2)
Si $M=\bigcup_{i\in I}M_i$ avec $\forall i,j\in I,\;\exists k\in I,\;
M_k\supseteq M_i\cup M_j$ (on dit alors que~$M$ est \emph{réunion filtrante} des $M_i$),
et si chaque $M_i$ est plat, alors $M$ est~plat.
\index{reunion@réunion!filtrante}
\index{filtrante!réunion ---}

 3)
Soit $a$ un \elt \ndz dans  $\gA$, $M$ un \Amo et $u\in M$ tels que~$au=0$.
Si cette \syzy s'explique dans $M$, on écrit $u=\sum_ia_iu_i$ ($a_i\in\gA$,
$u_i\in M$) avec les $aa_i=0$, donc $u=0$. 
Ainsi dans un module plat, tout \elt annulé par un \elt \ndz est nul.  

 4) (Suite) %\rdb \label{NOTAReg}
Le \emph{sous-module de torsion} d'un module $M$ est le module 

\snic{N=\sotq{x\in M}{\exists a\in \Reg(\gA),\;ax=0},} 

%\sni
où $\Reg(\gA)$ désigne le filtre des \elts
\ndzs de $\gA$. Ce sous-module de torsion est le noyau du morphisme
d'\eds à $\Frac\gA$ pour le module $M$.
Le sous-module de torsion d'un module  plat est réduit à $0$.
\\
\emph{Lorsque l'anneau $\gA$ est intègre}, on dit qu'un module est \emph{sans torsion} si son sous-module de torsion
est réduit à $0$. Sur un anneau de Bézout intègre, ou plus \gnlt sur un \ddpz, un module est 
plat \ssi il est sans torsion (exercice \ref{exoPlatsLecteur} et  \thref{thPruf} point \emph{2b}). 
\\
Nous donnons plus loin une \gnn de la notion de module sans torsion pour
un anneau commutatif arbitraire (\dfnz~\ref{def.locsdz}).% 
\index{torsion!sous-module de ---}%
\index{torsion!module sans ---}%
\index{module!sans torsion}

 5) Nous verrons (proposition \ref{prop.itfplat}) 
qu'un \itf plat $\fa$ est \lopz, ce qui implique
 $\Al 2 \fa=0$ (\thrf{propmlm}). Ainsi, lorsque~$\gA$ est un anneau intègre
 non trivial et $\gB=\gA[x,y]$, l'\id $\fa=\gen{x,y}$ est un
 exemple de \Bmo sans torsion, 
 mais pas plat (puisque $\Al 2_\gB \fa=\gA$ 
 d'après l'exemple \paref{belexemple}).
 En fait, la relation \smashbot{$[\,y\;-x\,]\cmatrix{x\cr y}=0$} 
 ne s'explique pas dans $\fa$, mais dans~$\gB$.   
\eoe

%:    Proposition{propPlat1}  plusieurs relations 

\smallskip 
La proposition qui suit dit que l'\gui{explication} qui est donnée pour la \syzyz~$LX=0$ dans la \dfn d'un \mpl 
s'étend à un nombre fini de \syzysz.
%: propPlat1
\begin{proposition} 
\label{propPlat1} 
Soit $M$  un  \Amo plat. 
On considère une famille de $k$ \syzysz, écrites sous la forme  
$LX=0$,  où~$L\in \Ae {k\times n}$ et~$X\in M^{n\times 1}$. Alors, on peut trouver un entier~$m$, un vecteur $Y\in M^{m\times 1}$ et une matrice~$G$ dans~$ \Ae {n\times m}$ satisfaisant les \egts
%-----------------begin $$----------------
$$\preskip.3em \postskip.0em
 GY = X\quad {\rm et } \quad LG=0.
$$
%-----------------end $$------------------
\end{proposition}
%--- end-proposition-------------------
%-----------------begin proof------------------
\begin{proof}
Notons $L_1$, \ldots, $L_k$ les lignes de $L$. 
La \syzyz~$L_1X=0$  est expliquée 
par deux matrices $G_1$ et $Y_1$ et par deux \egtsz~$X = G_1Y_1$ et~$L_1G_1=0$. 
La \syzy $L_2X=0$ se réécrit $L_2G_1Y_1=0$  \cad $L'_2Y_1=0$. 
Cette \syzy s'explique 
sous la forme $Y_1 = G_2Y_2$ et $L'_2G_2=0$.\\ 
Donc $X = G_1Y_1 = G_1G_2Y_2$. Avec $L_1G_1G_2=0$  et $L_2G_1G_2=L_2'G_2=0$.
Le vecteur colonne $Y_2$ et la matrice $H_2=G_1G_2$ expliquent donc les deux
\syzys  $L_1X=0$  et  $L_2X=0$. \\
Il ne reste qu'à itérer le processus.
\end{proof}
%-----------------end proof------------------

%:--- Theorem{propPlat1}-----------
Une reformulation de la proposition~\ref{propPlat1} dans le langage 
catégorique 
est le \tho suivant. La \dem  est un exercice de traduction laissé \alecz.

%: thPlat1
\begin{theorem} 
\label{thPlat1} \emph{(\Carn des modules plats, 1)}
\perso{des dessins?}\\
Pour un \Amo $M$ \propeq
\begin{enumerate}
\item Le module $M$ est plat.
\item Toute \ali  d'un \mpf $P$ vers $M$ se 
factorise par un module libre de rang fini.
\end{enumerate}
\end{theorem}
%--- end-theorem-----------------------------------------
%
%

%:    Theorem{cor pf plat ptf}---  
\begin{theorem} 
\label{cor pf plat ptf}\relax 
Un \Amo 
$M$ est \pf et plat \ssi il est \ptfz.

%-----------------end item------------------
\end{theorem}
\begin{proof}
La condition est \ncr d'après la remarque qui suit.
Elle est suffisante, car l'\idt de $M$ se factorise par un \Amo
libre $L$ de rang fini. Alors, la composée $L\to M \to L$ est une
\prn d'image isomorphe à~$M$.
\end{proof}
%
%--- Lemma{lemPS}--------------

Il est \imd que le \Amo $M\oplus N$ est plat \ssi les modules~$M$ et~$N$ sont plats.

La proposition qui suit donne un peu mieux (voir aussi le \thref{propPlatQuotientdePlat} et l'exercice~\ref{propfreeplat}).

%:     Proposition{propSuExPlat}

\begin{proposition}\label{propSuExPlat}
Soit $N\subseteq M$ deux \Amosz. Si $N$ et $M/N$
sont plats, alors $M$ est plat.  
\end{proposition}
\begin{proof}
On écrit $\ov x$ pour l'objet $x$ (défini sur $M$) vu modulo $N$.
Considérons une \syzy $LX=0$ dans $M$. Puisque $M/N$ est plat, on obtient
$G$ sur $\gA$ et $Y$ sur $M$ tels que $LG=0$  et $G\ov Y=\ov X$.
Considérons le vecteur $X'=X-GY$ sur $N$. On a $LX'=0$, et puisque 
$N$ est plat, on obtient $H$ sur $\gA$ et $Z$ sur~$N$ tels que $LH=0$
et $HZ=X-GY$.\\
 Ainsi la matrice $\blocs{.6}{.5}{.8}{0}{$G$}{$H$}{}{}$
et le vecteur  $\blocs{.4}{0}{.6}{.5}{$Y$}{}{$Z$}{}$ 
expliquent la relation~$LX=0$.
\end{proof}
%

%--- Fact{fact.plat} -----------
\begin{fact} 
\label{fact.plat}\relax Soit $S$ un \mo de l'anneau $\gA$.
\begin{enumerate}
\item Le localisé $\gA_S$ est plat comme \Amoz. 
\item Si $M$ est un $\gA$-\mplz, alors $M_S$ est plat comme 
\Amo et comme $\gA_S$-module.
\end{enumerate}

\end{fact}
%--- end-fact--------------------
%
\begin{proof} Il suffit de montrer le point \emph{2.}
Si l'on a une \syzy $LX=0$ dans le \Amo $M_S$, on écrit $X=X'/s$ et l'on a une 
\syzy $uLX'=0$ dans $M$ (avec $u,s\in S$). On trouve donc $Y'$ sur $M$ et $G$
sur $\gA$ tels que
$GY'=X'$ dans $M$ et  $uLG=0$ dans $\gA$. Ceci implique pour $Y=Y'/(su)$ l'\egt
$uGY=X$ dans $M_S$, de sorte que $uG$ et $Y$ expliquent la relation $LX=0$ dans $M_S$. \Demo analogue si l'on part d'une \syzy dans $M_S$ vu comme $\gA_S$-module. 
\end{proof}
%

%:  --- Prc loc glob conc {plcc.plat} - 
\subsec{Principe local-global}

La platitude est une notion locale au sens suivant.

\begin{plcc} 
\label{plcc.plat}  
{\em  (Pour les modules plats)}\\
Soient $S_1$, $\ldots$, $S_r$  des \moco d'un anneau $\gA$, et soit $M$ \hbox{un \Amoz}.
%
%\vspace{-3pt}
\begin{enumerate}
\item Une \syzy $LX=0$ dans~$M$ s'explique dans~$M$ \ssi elle s'explique
dans chacun des $M_{S_i}$. 
\item Le module $M$ est plat sur $\gA$ \ssi chacun des $M_{S_i}$ est plat
sur~$\gA_{S_i}$.
\end{enumerate}
\end{plcc}
%--- end-plcc-----------------------------------------
%-----------------begin proof------------------
\begin{proof} Il suffit de démontrer le premier point.
Le \gui{seulement si} est donné par le fait \ref{fact.plat}. Voyons l'autre 
implication. Soit $LX=0$ une \syzy entre \elts de~$M$ (où~$L\in \Ae {1\times n}$ 
et $X\in M^{n\times 1}$). 
On cherche un entier 
$m\in\NN,$ un vecteur $Y\in M^{m\times 1}$ et une matrice~$G\in \Ae {n\times m}$ qui 
vérifient l'\eqrf{eqdef.plat}.
On a  une solution~$(m_i,Y_i,G_i)$  pour~$(\ref{eqdef.plat})$ dans chaque 
localisé $\gA_{S_i}$. 
\\
On peut
 écrire~$Y_i=Z_i/s_i$,~$G_i=H_i/s_i$ avec~$Z_i\in M^{m_i\times 1}$,
 $H_i\in \Ae {n\times m_i}$ et des~$s_i\in S_i$ convenables. 
On a alors~$u_iH_iZ_i=v_iX$ dans $M$ et~$u_iLH_i=0$ dans $\gA$ pour 
certains~$u_i$ et~$v_i\in S_i$. 
On écrit~$\sum_{i=1}^{r} b_i v_i =1$ dans~$\gA$. 
On prend pour~$G$  la matrice obtenue en juxtaposant 
en ligne les matrices~$b_i u_iH_i$, 
et pour $Y$ le vecteur obtenu en superposant en colonne les vecteurs~$Z_i$.
On obtient $ GY = \sum_{i=1}^rb_i v_iX=X$ dans $M$,  et $LG=0$  dans~$\gA$. 
\end{proof}
%-----------------end proof------------------

%:  --- Prc loc glob abs{plca.plat}---
Le principe correspondant en \clama est le suivant.

\begin{plca} 
\label{plca.plat}  
{\em  (Pour les modules plats)}
%\vspace{-3pt}
\begin{enumerate}

\item Une \syzy $LX=0$ dans $M$ s'explique dans~$M$ \ssi elle s'explique
dans $M_{\fm}$ pour tout \idemaz~$\fm$. 
\item Un \Amo   $M$  
 est plat \ssi pour tout \idemaz~$\fm$, le module $M_{\fm}$ est %un module 
 plat sur $\gA_{\fm}$.
\end{enumerate}
\end{plca}
%--- end-plca-----------------------------------------
%-----------------begin proof------------------
\begin{proof}
Il suffit de montrer le premier point. Or le fait qu'une \syzy $LX=0$ puisse être expliquée est  une \prt \carf (\dfn \ref{defiPropCarFini}). 
 On applique donc le fait~\ref{factPropCarFin} qui nous permet de passer du \plgc au \plga correspondant.
\end{proof}
%-----------------end proof------------------

%:  --- Autres \carns de la platitude -
%\penalty-5000 
\subsec{Autres \carns de la platitude}

Nous allons maintenant considérer des 
\emph{\syzys sur $M$ à \coes dans un autre module $N$} et
nous allons montrer que lorsque $M$ est plat,
toute \syzy à \coes dans n'importe quel module $N$ s'explique dans~$M$.

%:     Definition{defiRDLGen}
\begin{definition}\label{defiRDLGen}
Soient $M$ et $N$ deux \Amosz.  \\
Pour $L=[\,a_1\,\cdots\,a_n\,]\in N^{1\times n}$ et $X=\tra[\,x_1\,\cdots\,x_n\,]\in M^{n\times 1}$,
on note 
$$\preskip.4em \postskip.4em 
L\odot X\;\eqdefi\;\som_{i=1}^na_i\te x_i\;\in N\te M. 
$$
\begin{enumerate}
\item Si $L\odot X=0$ on dit que l'on a une \syzy entre les $x_i$ à \coes 
dans~$N$.
\item On dit que \emph{la \syzy $L\odot X=0$ s'explique dans~$M$} si l'on a $Y\in M^{m\times 1}$ et une matrice $G\in \Ae {n\times m}$ qui vérifient:
%---  equation eqdef2.plat -------
\begin{equation}\label{eqdef2.plat}\relax
LG =_{N^ {1\times m}}0 \quad {\rm et } \quad X=_{M^{n\times 1}}GY\,.
\end{equation}
%---------------------end equation--------------
%
\end{enumerate}
\end{definition}

\rem 1) Lorsque l'on dit que la \syzyz~\hbox{$LX=0$} s'explique dans $M$, on signifie
que l'explication \gui{ne touche pas à $L$}.

2) On notera que de manière \gnle l'\egt $L\odot G Y=LG\odot Y$ est assurée pour toute matrice $G$ à \coes dans $\gA$ parce que $a\otimes \alpha y=a\alpha\otimes y$ lorsque~$a\in N$, $y\in M$ et~$\alpha\in\gA$.
\eoe

%:     Proposition{propPlatRdl}
\begin{proposition}\label{propPlatRdl}
Soient $M$ et $N$ deux \Amosz. 
\\
Si $M$ est un \Amo plat
toute \syzy  à \coes dans $N$
 s'explique dans~$M$.  
\end{proposition}
\begin{proof}
On suppose donnée une \syzy $L\odot X=0$ avec~$L=[\,a_1\,\cdots\,a_n\,]\in N^{1\times n}$ et $X=\tra[\,x_1\,\cdots\,x_n\,]\in M^{n\times 1}$.

\emph{Cas où $N$ est libre de rang fini.} La proposition \ref{propPlat1} donne le résultat.

\emph{Cas où $N$ est \pfz.} \\
On écrit $N=P/R=\gA^k\sur{(\gA c_1+\cdots\gA c_r)}$.
Les $a_i$ sont donnés par des~$b_i$ de~$P$. La relation $L\odot X=0$
signifie que  $\sum_i b_i\otimes x_i\in R\otimes M\subseteq P\otimes M$, i.e. 
que l'on a une \egt
$$\preskip.4em \postskip.4em \ndsp
\sum_i b_i\otimes x_i+\sum_\ell c_\ell\otimes z_{\ell}=0 
$$ 
dans $P\otimes M$. On constate alors que lorsque l'on explique dans $M$ cette \syzy
(portant sur les $x_i$ et les $z_\ell$) à \coes dans le module libre $P$,  
on explique par la même occasion la \syzy  $L\odot X=0$
à \coes dans~$N$.

\emph{Cas d'un \Amo $N$ arbitraire.} 
\\
Une relation $L\odot X=\sum_i a_i\otimes x_i=0$ provient d'un calcul fini, dans lequel n'interviennent qu'un nombre fini d'\elts de $N$ et de relations entre ces \eltsz.
Il existe donc un \mpf $N'$, une \ali $\varphi: N'\to N$ et des $b_i\in N'$ tels que d'une part $\varphi(b_i)=a_i$ ($i\in\lrbn$), et d'autre part $\sum_i b_i\otimes x_i=0$ dans $N'\otimes M$.
On constate alors que lorsque l'on explique dans $M$ cette \syzy
%(portant sur les $x_i$) 
à \coes dans $N'$ (lequel est un \mpfz),  
on explique par la même occasion la \syzy  $L\odot X=0$
à \coes dans~$N$. 
\end{proof}
% 

%:     theorem{thplatTens}-------
\begin{theorem} \emph{(\Carn des modules plats, 2)}
\label{thplatTens} \\
Pour un \Amo $M$  \propeq
%-----------------begin enum------------------
\begin{enumerate}
\item \label{i1thplatTens} 
Le module $M$ est plat.
\item \label{i2thplatTens} 
  Pour tout \Amo $N$, toute \syzy entre \elts de $M$ à \coes dans
$N$ s'explique dans~$M$.
\item \label{i3thplatTens} 
  Pour tout \itf $\fb$ de $\gA$ l'application canonique 
$\fb\otimes_\gA M\alb\rightarrow M$ est injective (ceci établit donc un \iso de
$\fb\otimes_\gA M $ sur~$\fb M$).
\item \label{i4thplatTens} 
  Pour tous \Amos  $N\subseteq N'$, l'\ali canonique 
$$\preskip.4em \postskip.2em 
N\otimes_\gA M\rightarrow N'\otimes_\gA M 
$$  
est injective.
\item \label{i5thplatTens} 
  Le foncteur $\bullet\otimes M$ préserve les suites exactes.
\end{enumerate}
%-----------------end enum------------------ 
\end{theorem}
%--- end-theorem-----------------------------------------
%-----------------begin proof------------------
\begin{proof}
L'implication \emph{\ref{i5thplatTens}}  $\Rightarrow$ \emph{\ref{i3thplatTens}} est triviale. 

\emph{\ref{i4thplatTens}}  $\Rightarrow$ \emph{\ref{i5thplatTens}.}
Les suites exactes courtes sont préservées par le foncteur~$\bullet\otimes M$. 
Or toute suite exacte se décompose en suites exactes courtes (voir \paref{sexaseco}).

 \emph{\ref{i1thplatTens}} $\Leftrightarrow$ \emph{\ref{i3thplatTens}.} D'après le lemme du tenseur nul \ref{lem-tenul}.

 \emph{\ref{i1thplatTens}} $\Rightarrow$ \emph{\ref{i2thplatTens}.} 
C'est la proposition~\ref{propPlatRdl}.

 \emph{\ref{i2thplatTens}} $\Leftrightarrow$ \emph{\ref{i4thplatTens}.}
D'après le lemme du tenseur nul \ref{lem-tenul}.
\end{proof}
%-----------------end proof------------------

Le \tho précédent admet quelques corolaires importants. 
%:     Corollary{corPlatTens}
\begin{corollary}\label{corPlatTens} \emph{(Produit tensoriel)}\\
 Le produit tensoriel de deux \mpls est un \mplz.
\end{corollary}
\begin{proof}
Utiliser le point \emph{\ref{i4thplatTens}} du
\tho \ref{thplatTens}.
\end{proof}
%
%:     Corollary{cor2PlatTens}
\begin{corollary}\label{cor2PlatTens} \emph{(Autres constructions de base)}\\
 Les puissances tensorielles, extérieures et \smqs d'un \mpl sont des \mplsz.
\end{corollary}
\facile
%

%:     Corollary{cor3PlatTens}
\begin{corollary}\label{cor3PlatTens} \emph{(Intersection)}\\
Soient $N_1,\ldots, N_r$ des sous-modules d'un module $N$ et soit $M$ un \mplz.
Puisque~$M$ est plat, pour tout sous-module $N'$ de $N$, identifions $N'\otimes M$
avec son image dans $N\otimes M$.
Alors on a l'\egt 
$$\preskip.4em \postskip.4em\ndsp 
\left(\bigcap_{i=1}^r N_i\right) \otimes M= \bigcap_{i=1}^r (N_i\otimes M). 
$$ 
\end{corollary}
\begin{proof}
La suite exacte 
$%\preskip.0em \postskip.2em\ndsp 
0\to \bigcap_{i=1}^r N_i\to N\to \bigoplus_{i=1}^r(N\sur {N_i}) 
$
est préservée par le produit tensoriel avec $M$
et le module $(N\sur {N_i})\te M$ s'identifie à $(N\te M)\sur{(N_i\te M)}$.
\end{proof}
%

%:     Corollary{corPlatEds}
\begin{corollary}\label{corPlatEds} \emph{(\Edsz)}
Soit $\rho:\gA\to\gB$ une \algz. Si~$M$ est un \Amo plat, alors
$\rho\ist (M)$ est un \Bmo plat. 
\end{corollary}
\begin{proof} 
On note que 
pour un \Bmoz~$N$, on a 
$$\preskip.4em \postskip.4em 
N\te_\gB\rho\ist (M)\simeq N\te_\gB\gB\te_\gA M\simeq N\te_\gA M. 
$$
On applique alors  le point \emph{4} du \tho \ref{thplatTens}.
Notez que le dernier produit tensoriel est muni d'une structure de \Bmo via~$N$.
\end{proof}

\medskip 
\rem Nous venons d'utiliser sans le dire une forme \gnee d'associativité du produit tensoriel dont nous laissons la \dem \alecz.
C'est la suivante.\\
Tout d'abord on dit qu'un groupe abélien $P$ est un $(\gA,\gB)$-bimodule
s'il est muni de deux lois externes qui en font respectivement un \Amo et \hbox{un
\Bmoz}, et si ces deux structures sont compatibles au sens suivant: pour tous $a\in\gA$, $b\in\gB$ et $x\in P$, on a $a(bx)=b(ax)$.\\
Dans un tel cas, si $M$ est un \Bmoz, alors le produit tensoriel $M\te_\gB P$
peut \hbox{lui même} être muni d'une structure de $(\gA,\gB)$-bimodule
en posant, pour $a\in\gA$, $a(x\te y)=_{M\te_\gB P}x\te ay$.\\
De même, lorsque  $N$ est un \Amoz,  le produit tensoriel $P\te_\gA N$
peut \hbox{lui même} être muni d'une structure de $(\gA,\gB)$-bimodule
en posant, \hbox{pour $b\in\gB$}, $b(y\te z)=_{P\te_\gA N}by\te z$.

\emph{Lemme d'associativité.}
Sous ces hypothèses, il existe une unique application  $(\gA,\gB)$-\lin $\varphi:(M\te_\gB P)\te_\gA N\to M\te_\gB (P\te_\gA N)$ qui
vérifie 
$$\preskip.4em \postskip.4em 
\varphi\big((x\te y)\te z\big)=x\te (y\te z) 
$$
pour tous
$x\in M$, $y\in P$, $z\in N$. Enfin
$\varphi$ est un \isoz.
\eoe

%:subsec{Quotient plat d'un module plat}
\subsec{Quotients plats}

%:     theorem{propPlatQuotientdePlat}
\begin{theorem}\label{propPlatQuotientdePlat} \emph{(Quotients plats)}\\
Soit $M$ un \Amoz, $K$ un sous-module et $N=M/K$, avec la suite exacte
$$\preskip.2em \postskip.4em
0 \to K\vers\imath M\vers \pi N \to 0 
$$ 
\begin{enumerate}
\item Si $N$ est plat, pour tout module $P$, la suite
$$\preskip.2em \postskip.2em
0 \to K\te P\vers{\imath_P} M\te P\vers {\pi_P} N\te P \to 0 
$$ 
est exacte  ($\imath_P=\imath\te\Id_P$, $\pi_P=\pi\te\Id_P$).
\item Si $N$ et $M$ sont plats, $K$ est plat.
\item Si $N$ et $K$ sont plats, $M$ est plat.
\item Si $M$ est plat, \propeq
\begin{enumerate}
\item $N$ est plat.
\item Pour tout \itf $\fa$, on a $\fa M \cap K = \fa K$.
\item Tout \itf $\fa$ donne une suite exacte
$$\preskip.2em \postskip.2em
0 \to K/\fa K\vers{\imath_\fa} M/\fa M\vers {\pi_\fa} N/\fa N\to 0. 
$$ 
\end{enumerate}
\end{enumerate}
\end{theorem}
%--------- fin theorem ----------------------------------------------
%
NB: le point \emph{3} a déjà fait l'objet de  la proposition \ref{propSuExPlat}, nous en donnons ici une autre \demz, laissant \alec
le soin de les comparer.
\begin{proof} 
\emph{1.} \emph{Cas où $P$ est \tfz}. On écrit $P$ comme quotient d'un module libre fini $Q$ avec une suite exacte courte
$$0 \to R\vers a Q\vers p P \to 0. 
$$ 
On considère alors le diagramme commutatif suivant dans lequel 
toutes les suites horizontales et verticales sont exactes parce que
$N$ et $Q$ sont plats

\smallskip 
\centerline{\small
\xymatrix@C=3em@R=1.6em {
&   &     & 0\ar[d]
\\
 &K\te R\ar[d]_{a_K} \ar[r]^{\imath_R} 
&M\te R\ar[d]_{a_M}\ar[r]^{\pi_R} &
N\te R\ar[d]_{a_N}\ar[r]&0
\\
0\ar[r]&K\te Q\ar[d]_{p_K} \ar[r]^{\imath_Q} 
&M\te Q\ar[d]_{p_M}\ar[r]^{\pi_Q} &
N\te Q\ar[r]&0
\\
 &K\te P\ar[d] \ar[r]^{\imath_P} &M\te P\ar[d]  & & 
\\
&0  &  0  &  
\\
}}

On doit montrer que $\imath_P$
est injective. Il s'agit d'un cas particulier du lemme du serpent,
on peut le démontrer par une \gui{chasse dans le diagramme}. \\
On suppose $\imath_P(x)=0$. On écrit $x=p_K(y)$ et $v=\imath_Q(y)$.
On a $p_M(v)=0$, donc on écrit $v=a_M(z)$. 
\\
Comme $\pi_Q(v)=0$, on a 
$a_N(\pi_R(z))=0$, donc $\pi_R(z)=0$. \\
Donc on écrit $z=\imath_R(u)$ et
l'on a 

\snic{\imath_Q(a_K(u))=a_M(\imath_R(u))=a_M(z)=v=\imath_Q(y),}

\snii
et comme
$\imath_Q$ est injective, $y=a_K(u)$, d'où $x=p_K(y)=p_K(a_K(u))=0$.

\smallskip 
\emph{Cas \gnlz.} Une possibilité est d'écrire $P$ comme quotient d'un module plat $Q$ (voir à ce sujet l'exercice \ref{propfreeplat}) auquel cas la \dem précédente est inchangée. On peut aussi se passer de cette construction un peu lourde comme suit. Montrons que $\imath_P$
est injective. Soit $x=\sum_ix_i\te y_i\in K\te P$ tel que 
$\imath_P(x)=_{M\te P}0$, i.e. $\sum_ix_i\te y_i=_{M\te P}0$.
\\
D'après la \dfn du produit tensoriel, il existe un sous-\mtf $P_1\subseteq P$ tel que l'on a aussi $\sum_ix_i\te y_i=_{M\te P_1}0$. D'après le cas déjà examiné, on a 
$\sum_ix_i\te y_i=_{K\te P_1}0$, et ceci implique $\sum_ix_i\te y_i=_{K\te P}0$.

%\pagebreak	 
\emph{2} et \emph{3.} Soit $\fa$ un \itf arbitraire. Puisque $N$ est plat, on a d'après le point \emph{1} 
un diagramme commutatif avec suites exactes

\smallskip 
\centerline{\small
\xymatrix@C=3em@R=1.6em {
  &  & &0\ar[d] & 
\\
0\ar[r]&\fa\te K\ar[d]^{\varphi_K} \ar[r]^{\imath_\fa} &\fa\te M\ar[d]^{\varphi_M}\ar[r]^{\pi_\fa} &\fa\te N\ar[d]^{\varphi_N}\ar[r] & 0
\\
0\ar[r] &K\ar[r]^{\imath} &M\ar[r]^\pi &N\ar[r] & 0
\\
}.}

\smallskip Si $M$ est plat, $\varphi_M$ est injective, donc aussi $\varphi_M\circ \imath_\fa$, puis $\varphi_K$. On conclut par le point \emph{3} du \thref{thplatTens} que $K$ est plat.

Si $K$ est plat, $\varphi_K$ est injective et une petite chasse dans le diagramme montre que $\varphi_M$ est injective. Soit $x\in\fa\te M$ avec $\varphi_M(x)=0$. \\
Comme $\varphi_N(\pi_\fa(x))=0$, on a $\pi_\fa(x)=0$
et l'on peut écrire $x=\imath_\fa(y)$. 
\\
Alors $\imath(\varphi_K(y))=\varphi_M(x)=0$, donc $y=0$, donc $x=0$.

\emph{4a} $\Rightarrow$ \emph{4b.} Puisque $M$ et $N$ sont plats, $K$ l'est \egmt et la ligne du haut du diagramme précédent donne la suite exacte
$$\preskip.2em \postskip.4em
0 \to \fa K
\vvers{{\imath\,\frt{\fa\! K}}} 
\fa M
\vvers {{\pi\,\frt{\fa\! M}}} \fa N\to 0.\eqno(+) 
$$
Or le noyau de  $\pi\frt{\fa\! M}$ est par \dfn $\fa M\cap K$.

\emph{4b} $\Leftrightarrow$ \emph{4c.} La suite 
$$\preskip.2em \postskip.4em
0 \to K/\fa K\vers{\imath_\fa} M/\fa M\vers {\pi_\fa} N/\fa N\to 0 
$$ 
est obtenue à partir de la suite exacte $0\to K \to M \to N$
par \eds à $\gA/\fa$. Dire qu'elle est exacte revient à dire que $\imath_\fa$ est injective. Or un \elt
$\ov x\in K/\fa K$ est envoyé sur $0$ \ssi \hbox{on a $x\in \fa M\cap K$}.

\emph{4b} $\Rightarrow$ \emph{4a.}
Puisque $\fa K=\fa M\cap K$ la suite $(+)$ est exacte. On considère
le diagramme commutatif suivant avec suites exactes, pour lequel il nous faut montrer que $\varphi_N$ est injective.

\smallskip 
\centerline{\small
\xymatrix@C=3em@R=1.6em {
  &  & 0\ar[d]& & 
\\
&\fa\te K\ar[d]^{\varphi_K} \ar[r]^{\imath_\fa} &\fa\te M\ar[d]^{\varphi_M}\ar[r]^{\pi_\fa} &\fa\te N\ar[d]^{\varphi_N}\ar[r] & 0
\\
0\ar[r] &\fa K\ar[d]\ar[r]^{\imath\,\frt{\fa\! K}} &\fa M\ar[d]\ar[r]^{\pi\,\frt{\fa\! M}} &\fa N\ar[r]\ar[d] & 0
\\
& 0& 0&0
\\
}}
 
\smallskip 
Ceci résulte d'une petite chasse dans le diagramme. 
\\
Si $\varphi_N(x)=0$, on écrit $x=\pi_\fa(y)$. Comme  $\pi\,\frt{\fa\! M}(\varphi_M(y))=0$, on \hbox{a  $z\in\fa K$} tel que $\varphi_M(y)=\imath\,\frt{\fa\! K}(z)$, on écrit
$z=\varphi_K(u)$, avec $\varphi_M(\imath_\fa(u))=\varphi_M(y)$.
\\ Et comme $\varphi_M$ est injective, $y=\imath_\fa(u)$ et $x=\pi_\fa(y)=0$.
\end{proof}
%

%:     Corollary{corpropPlatQuotientdePlat}
\begin{corollary}\label{corpropPlatQuotientdePlat} \emph{(Une \alg plate)}
Soit $f\in\AuX = \gA[\Xn]$ et~$\Aux=\aqo\AuX{f}$. Alors, le \Amo $\Aux$
est plat \ssiz $\rc(f)^2 = \rc(f)$, \cad \ssi l'idéal $\rc(f)$ est engendré
par un \idmz.
\end{corollary}
%--------- fin corollary ---------------------------------------------- 
%
\begin{proof}
%D'après la proposition \ref{propPlatQuotientdePlat} l
Le \Amo
 $\Aux$ est plat \ssi pour tout \itf $\fa$ de $\gA$ on a
 {
\fbox{$\gen {f} \cap \fa[\uX] = f\fa[\uX]$}.
}

Si $\Aux$ est plat, on obtient en particulier pour $\fa = \rc(f)$, 
\hbox{que $\rc(f)^2 = \rc(f)$} 
(car~$f \in \gen {f} \cap \fa[\uX]$).

Réciproquement, supposons $\rc(f)^2 = \rc(f)$ et montrons que
$\Aux$ est plat. L'\idm $e$ tel que $\gen{e}=\gen{\rc(f)}$
scinde l'anneau en deux composantes. Dans la première on a $f=0$, et le résultat est clair. Dans la seconde, $f$ est primitif.  
On suppose maintenant $f$ primitif. 
\\
D'après le lemme de
Dedekind-Mertens\footnote{En fait, il s'agit d'une variante, avec essentiellement la même \demz, que nous laissons \alecz.}, pour tout \Amoz~$M$ l'\Ali
$M[\uX] \vvers {\times f} M[\uX]$ est injective. Appliqué à
$M = \gA/\fa$, cela donne l'encadré. En effet, écrivons $M[\uX] = \AuX/\fa[\uX]$
et 
supposons 
\linebreak 
que $g \in \gen {f} \cap \fa[\uX]$. Alors $g = fh$ pour un $h\in \AuX$, et $\overline h$ est dans
le noyau de $\AuX/\fa[\uX] \vvers {\times f} \AuX/\fa[\uX]$, donc
$\overline h = 0$, i.e. $h \in \fa[\uX]$, et $g \in f\fa[\uX]$. 
\end{proof}
%

%
%%%%%%%%%%%%%%%%%%%%%%%%%%%%%%%%%%%%%%%%%%%%%%%%%%%%%%%%%%%%%%%%%%%
%--- Sec {Modules plats tf}--    
\section{Modules plats \tfz}

%:     Lemma{lem.plat-tf}---------  
Dans le cas des  \mtfsz, la platitude est une 
\prt de nature plus \elrz.

\begin{lemma} 
\label{lem.plat-tf}\relax
On considère un \Amo $M$ \tfz, et  $X\in M^{n\times 1}$ un vecteur colonne
dont les \coos $x_i$ engendrent 
$M.$ 
   Le module $M$ est plat \ssi 
pour toute \syzy  $LX=0$  (où  $L\in \Ae {1\times n}$), on peut trouver 
 deux matrices $G$, $H\in \Mn(\gA)$ qui satisfont les \egts
$$\preskip.4em \postskip.4em 
H+G = \I_n,\;\; LG=0\;\; {\rm et } \;\; HX = 0.
$$
En particulier, un module monogène $M=\gA y$ est plat \ssi 

\snic{\forall a\in\gA, (\,ay=0\;\Longrightarrow\;\exists s\in\gA,\;\;as=0\et sy=y\,)\,.
} 
\end{lemma}
%--- end-lemma--------------------
\rem La symétrie entre $L$ et $X$ dans l'énoncé n'est qu'apparente:
le module $M$ est engendré par les \coos de $X$, alors que l'anneau $\gA$ n'est pas engendré (comme sous-module) par les \coos de $L$.  
\eoe
%-----------------begin proof------------------
\begin{proof}
On ramène une \syzy arbitraire $L'X'=0$ à une \syzy $LX=0$   en exprimant 
$X'$ en fonction de $X$. A priori on devrait écrire $X$ sous forme $G_1Y$ avec 
$LG_1=0$. 
\\
Comme $Y=G_2X$, on prend $G=G_1G_2$ et $H=\I_n-G$.
\end{proof}
%-----------------end proof------------------

\rem Pour les modules monogènes, en posant $t=1-s$, on obtient des conditions
sur $t$ plutôt que sur $s$:
$$\preskip.3em \postskip.4em
a=at \qquad   \hbox{et} \qquad  ty=0,
$$
ce qui implique que l'annulateur $\fa$ de $y$ vérifie $\fa^2=\fa$.
En fait, d'après le \thref{propPlatQuotientdePlat},
$\gA\sur\fa$ est plat sur $\gA$ \ssi
pour tout \itf $\fb$ on a l'\egt $\fa\cap\fb=\fa\fb$.
\eoe

%:   Proposition{propPlat2}-----    

\ms Nous pouvons donner une \gnn  du  lemme \ref{lem.plat-tf} exactement dans le 
style de la proposition \ref{propPlat1}.
\begin{proposition} 
\label{propPlat2} 
Soit $M$ un \Amo plat \tfz, et  $X\in M^{n\times 1}$ un vecteur colonne 
qui engendre~$M.$  
Soit une famille de $k$ \syzys écrites sous la forme  $LX=0$  où~$L\in 
\Ae {k\times n}$ et 
$X\in M^{n\times 1}$. Alors, on peut trouver une matrice~$G\in \Mn(\gA)$ qui 
vérifie les \egts

\snic{ LG=0\;\; \hbox{ et } \;\; GX = X.
}

\end{proposition}
%--- end-proposition-------------------
%-----------------begin proof------------------
\begin{proof}
Identique à la preuve de la proposition \ref{propPlat1}.
\end{proof}
%-----------------end proof------------------

%\medskip 
Un substitut \cof pour la \prt selon laquelle
tout \evc sur un corps admet une base (vraie seulement
en \clamaz) est le fait que tout \evc sur un corps discret est plat.
\Prmt on a le résultat suivant.

%: --- propEvcPlat
\begin{theorem}\label{propEvcPlat} \Propeq
\begin{enumerate}
  \item \label{i1propEvcPlat}  Tout \Amo $\aqo\gA a$ est plat.
  \item \label{i2propEvcPlat}  Tout \Amo est plat.
  \item \label{i3propEvcPlat}  L'anneau $\gA$ est \zed réduit.
\end{enumerate}
\end{theorem}
\begin{proof}
\emph{\ref{i1propEvcPlat}} $\Rightarrow$ \emph{\ref{i3propEvcPlat}}. 
Si $\aqo{\gA}{a}$ est plat, alors $\gen{a} = \gen{a}^2$ et
si c'est vrai pour tout $a$, c'est que $\gA$ est \zed réduit.
\\
\emph{\ref{i3propEvcPlat}} $\Rightarrow$ \emph{\ref{i2propEvcPlat}}. Traitons d'abord le cas d'un corps discret. 
\\
On considère une \syzy $LX=a_1x_1+\cdots+a_nx_n=0$
pour des \elts $x_1$, \ldots, $x_n$ d'un \Amo $M$. Si tous les~$a_i$ sont nuls
la relation est expliquée avec $Y=X$ et~$G=\I_n$: $LG=0$ et~$GY=X$.
Si un des~$a_i$ est \ivz, par exemple $a_1$, posons $b_j=-a_1^{-1}a_j$
pour $j\neq1$. On a $x_1=b_2x_2+\cdots+b_nx_n$ et $a_1b_j+a_j=0$ pour $j>1$.
La \syzy est expliquée par $Y=\tra [\, x_2\; \cdots\; x_n\, ] $ et par la matrice $G$ suivante,
car $LG=0$ et~$GY=X$:

\snic{G=\cmatrix{
b_2&b_3&\cdots&b_n\cr
1& 0&\cdots&0\cr
0&\ddots& &\vdots\cr
\vdots& &\ddots&0\cr
0&\cdots &0&1
}.
}

%\sni

Pour un anneau \zed réduit, on applique la machinerie
\lgbe \elr \num2 qui nous ramène au cas d'un corps discret.\imlgz
\end{proof}
%%%%%%%%%%%%%%%%%%%%%%%%%%%%%%%%%%%%%%%%%
NB: ceci justifie la terminologie \gui{absolument plat} pour \zed réduit.

%:     Lemma{lem.plat2-tf}---------  

\begin{lemma} 
\label{lem.plat2-tf}\relax
Même contexte que dans le lemme \ref{lem.plat-tf}.
 Si $\gA$ est un \alo et  $M$ est plat, on obtient sous 
l'hypothèse $LX=0$ l'alternative suivante:
le vecteur $L$ est nul, ou l'un des $x_i$ dépend \lint des autres (il peut donc être supprimé 
dans la liste des \gtrs de~$M$).
\end{lemma}
%--- end-lemma--------------------
%
\begin{proof}
C'est un \gui{truc du \deterz}. On note que
$\det(G)=\det(\I_n-H)$  s'écrit $1+\sum_{i,j}b_{i,j} h_{i,j}$. Donc $\det(G)$ ou 
l'un des $h_{i,j}$ est \ivz. Dans le premier cas $L=0$; dans le deuxième cas,  puisque $HX=0$,
l'un des vecteurs $x_i$ s'exprime en fonction des autres.
\end{proof}

La même \dem dans le cas d'un anneau arbitraire donne le résultat
suivant.

%:     Lemma{lem.plat3-tf}---------  

\begin{lemma} 
\label{lem.plat3-tf}\relax
Même contexte que dans le lemme \ref{lem.plat-tf}.
\\
 Si   $M$ est plat et $LX=0$, il existe des \eco $s_1$, \dots, $s_\ell$ tels que sur chacun des anneaux $\gA [1/s_j]$ on a
 $L=0$, ou l'un des $x_i$ est une \coli des autres.
\end{lemma}
%--- end-lemma--------------------

En \clamaz, le  lemme \ref{lem.plat2-tf} implique le fait suivant. 
%:     Fact{factplatlocalclass}
\begin{factc}\label{factplatlocalclass}
Un \mpl \tf sur un \alo est libre et une base peut être extraite de n'importe 
quel \sgrz. 
\end{factc}

Et à partir du lemme \ref{lem.plat3-tf}, on obtient ce qui suit.

%:     Fact{factplattfclass}
\begin{factc}\label{factplattfclass}
Un \mpl \tf sur un anneau intègre est \ptfz. 
\end{factc}

Voici une version \cov du fait$\etl$~\ref{factplatlocalclass}.

%:     Proposi{prop.mtf loc plat1}---
\begin{proposition} 
\label{prop.mtf loc plat1}\relax
Soit  $\gA$ un \alo et $M$ un \Amo \tf plat engendré par $(\xn)$. 
Supposons que  $M$ soit \fdi ou que l'existence de \syzys non triviales soit 
explicite dans $M$.
Alors, $M$ est librement engendré par une suite finie $(x_{i_1},\ldots,x_{i_k})$ 
(avec $k\geq 0$).
\end{proposition}
%--- end-proposition----------------------------------------
%-----------------begin proof------------------
\begin{proof}
Supposons d'abord que $M$ soit \fdiz, on peut alors trouver une suite finie 
d'entiers  
$1\leq i_1<\cdots<i_{k}\leq n$ (où $k\geq 0$) telle qu'aucun des~$x_{i_\ell}$ ne 
soit une \coli des autres, et~$(x_{i_1},\ldots,x_{i_k})$ engendre~$M$. 
Pour simplifier les notations, on suppose donc désormais que $k=n$, i.e., aucun 
des~$x_i$ n'est \coli des autres. Le lemme~\ref{lem.plat2-tf}  nous dit alors 
que toute \syzy entre les~$x_i$ est triviale.
\\
Supposons maintenant que l'existence de \syzys non triviales soit explicite dans 
$M$, \cad que pour toute famille d'\elts de $M$, on sache dire s'il y a une \syzy 
non triviale entre ces \elts et en fournir une le cas échéant. Alors, en 
utilisant le lemme \ref{lem.plat2-tf} on peut supprimer un à un les \elts 
superflus dans la famille $(x_i)$ sans changer le module $M$, jusqu'à ce qu'il 
ne reste qu'une sous-famille sans \syzy non triviale (un cas limite est fourni par 
la partie vide lorsque le module est nul).
\end{proof}
%-----------------end proof------------------
\comm Notez que la preuve utilise l'hypothèse \gui{$M$ est \fdiz}, ou \gui{l'exis\-tence 
de \syzys non triviales est explicite dans $M$}
uniquement avec des familles extraites du système \gtr $(x_i)$.
Par ailleurs, chacune de ces hypothèses est trivialement vraie en \clamaz. 
\eoe

\medskip 
Voici maintenant une version \cov du fait$\etl$~\ref{factplattfclass}.

%:     Proposi{prop.mtf integre plat}---
\begin{proposition} 
\label{prop.mtf integre plat}\relax
Soit  $\gA$ un anneau intègre et $M$ un \Amo \tf plat engendré par $(\xn)$. 
Supposons que  pour toute partie finie~$J$ de $\lrbn$ l'existence de \syzys non triviales entre $(x_j)_{j\in J}$ soit 
explicite dans $M$ (autrement dit, en passant au corps de fractions on obtient un \evc de dimension finie).
Alors, $M$ est \ptfz.
\end{proposition}
%--- end-proposition----------------------------------------
%-----------------begin proof------------------
\begin{proof}\perso{On écrit la \dem dans le cas où $\gA$
est non trivial. Dans le cas \gnlz, les \gui{ou bien} doivent être
replacés par de simples \gui{ou}, et \gui{{$\neq 0$}} doit être
remplacé par \gui{\ndzz}. }
On suppose \spdg $\gA$  non trivial, en 
utilisant le lemme~\ref{lem.plat3-tf} on obtient l'alternative suivante. 
Ou bien $(\xn)$ est une base, ou bien après \lon en des \eco le module est engendré 
\linebreak 
par $n-1$
des $x_j$. On conclut par \recu sur $n$: en effet, les \syzys après \lon en $s$
avec $s\neq 0$ sont les mêmes que celles sur $\gA$. 
\\
Notez que pour $n=1$,  ou bien $(x_1)$ est une base, ou bien $x_1=0$.
\end{proof}
%-----------------end proof------------------

%%%%%%%%%%%%%%%%%%%%%%%%%%%%%%%%%%%%%%%%%%%%%%%%%%%%%%%%%%%%%%%%%%%
%--- Sec {Idx principaux plats      
\section{Idéaux principaux plats}

Un anneau  $\gA$ est dit {\em  \sdzz}  
si l'on a:\index{sans diviseur de zéro!anneau ---}
\index{anneau!sans diviseur de zéro}
%-----begin equation- eqSDZ -- 
\begin{equation}\label{eqSDZ}
\forall a,b\in \gA\quad \big(ab=0 \; \Rightarrow \; (a=0\;{\rm ou}\;b=0)\big)
\end{equation}
%---------------------end equation--------------

Un anneau intègre (en particulier un corps discret) est \sdzz.
Un anneau discret \sdz est intègre.
Un anneau non trivial est intègre \ssi il est discret et \sdzz.

%--- Lemma{lemIdps}-------------  
\begin{lemma} \emph{(Quand un \id principal est plat)}
\label{lemIdps} ~
%-----------------begin item------------------
\begin{enumerate}
\item   Un \idpz, ou plus \gnlt un \Amo monogène
 $\gA  a$, est un \mpl \ssi 

\snic{\forall x\in\gA\quad \big(xa=0\;\Rightarrow\; \exists 
z\in\gA\;(za=0\;\;\hbox{et} 
\;\; xz=x)\big).
}

\item  Si $\gA$  est local, un \Amo  $\gA a$ %$\gen{a}$ 
est plat \ssi

\snic{\forall x\in\gA\quad \big(xa=0\;\Rightarrow\;(x=0 \;\;\hbox{ou}\;\; 
a=0)\big).
}

\item   Soit $\gA$  un \aloz, si $\gA$ est discret, ou si l'on  a un test pour 
répondre à la question \gui{$x$ est-il \ndzz?}, 
alors, un idéal $\gen{a}$ est plat \ssi $a$ est nul ou \ndzz.
%:2012 remplacé  x par  a  ligne au-dessus
\item   Pour un anneau local, $\gA$ \propeq
%-----------------begin item------------------
\begin{enumerate}
\item    Tout \idp est plat.
\item   L'anneau est \sdzz.
\end{enumerate}
%-----------------end item------------------
\end{enumerate}
%-----------------end item------------------
\end{lemma}
%--- end-lemma--------------------
%-----------------begin proof------------------
\begin{proof}
Le lemme \ref{lem.plat-tf} donne le point \emph{1}. Le calcul pour le point~\emph{2} en résulte,
car~$z$ ou $1-z$ est \ivz. La suite est claire. 
\end{proof}
%-----------------end proof------------------  

On a de même les \eqvcs suivantes.

%:     Lemma{corlemIdps}-----------  
\begin{lemma} 
\label{corlemIdps} 
Pour un anneau $\gA$, \propeq
%-----------------begin enum------------------
\begin{enumerate}
\item Tout \idp de $\gA$ est plat. 
\item Si $xy=0$,  on a
       $\Ann\,x+\Ann\,y=\gA.$
\item Si $xy=0$, il existe des \moco $S_{i}$ tels que dans chacun des 
localisés $\gA_{S_{i}}$, $x$ ou $y$ devient nul. 
\item Si $xy=0$, il existe $z\in\gA$ avec $zy=0$ et $xz=x$. 
\item Pour tous $x,\,y \in \gA, \quad \Ann\,xy =  \Ann\,x+\Ann\,y $.
\end{enumerate}
%-----------------end enum------------------
\end{lemma}
%--- end-lemma------------------------------------

%:     Definition{def.locsdz} et modules sans torsion 

La \prt pour un anneau d'être \sdz se comporte mal par 
\rcm et celle pour un module d'être plat se comporte bien par \lon et 
\rcmz. Cela justifie la \dfn suivante.

\begin{definition} 
\label{def.locsdz}~
\begin{enumerate}
\item  Un anneau $\gA$ est dit {\em \lsdzz} 
lorsqu'il vérifie les \prts \eqves 
du lemme~\ref{corlemIdps}.
\item Un \Amo $M$ est dit \emph{sans torsion} lorsque tous ses
sous-modules monogènes sont plats (voir le lemme \ref{lemIdps}). 
%:HHH reference plutot au lemme \ref{lemIdps}
%
\end{enumerate}
\index{localement!anneau --- sans diviseur de zéro}%
\index{anneau!localement sans diviseur de zéro}%
\index{torsion!module sans ---}%
\index{sans torsion!module ---}%
\index{module!sans torsion}%
\end{definition}
%--- end-definition------------------------------------

\rems
\\
1) Le sous-module de torsion d'un module sans torsion est réduit à $0$.
Notre \dfn est donc un peu plus contraignante que celle, plus usuelle, 
qui dit qu'un module est
sans torsion lorsque son module de torsion est réduit à $0$.
On notera que les deux \dfns coïncident 
lorsque l'anneau $\gA$ est \qiz.
 
2) Tout sous-module d'un module sans torsion est sans torsion, ce qui n'est
pas le cas en \gnl lorsque l'on remplace \gui{sans torsion} par \gui{plat}.
 
3) Un anneau \lsdz est réduit. 
 
4) Dans la littérature de langue anglaise, 
on trouve parfois l'appellation \gui{pf-ring} (principal ideals are flat) pour un  anneau \lsdzz.\index{pf-ring} 
 
5) Un anneau local est \lsdz \ssi il est \sdzz. 
 
6) Le corps des réels  {\em n'est pas} \sdz (\emph{ni} \lsdzz): c'est un anneau local pour lequel on ne sait pas réaliser 
explicitement l'implication (\ref{eqSDZ}) \paref{eqSDZ}.

7) En \clama un anneau est \lsdz \ssi il devient intègre après \lon en tout \idep (exercice \ref{exoClamlsdz}). 
\eoe

%:      Lemma{lem.platsdz} --------  
\begin{lemma} 
\label{lem.platsdz}\relax Soit $\gA$ un anneau \lsdz  
et $M$ un module plat sur $\gA$.
%-----------------begin enum------------------
\begin{enumerate}
\item  Le module $M$ est sans torsion.
\item  L'annulateur $(0:y)$ de n'importe quel  $y\in M$ est \idmz.
\end{enumerate}
%-----------------end enum------------------
\end{lemma}
%--- end-lemma--------------------
%-----------------begin proof------------------
\begin{proof}
\emph{1.} Supposons $ay=0$, $a\in\gA$, $y\in M$.
Puisque $M$ est plat on a des \elts $x_i$ de $M$, des \elts $b_i$ de $\gA$, et une \egt $y=\sum_{i=1}^n b_i x_i$
dans~$M$, avec $ab_i=0$ ($i\in\lrbn$) dans $\gA$.  
\\
Pour chaque $i$, puisque $ab_i=0$, il existe $c_i$ tel que $ac_i=a$ et $c_ib_i=0$.
 On pose  $c=c_1\cdots c_{n}$. Alors, $a=ca$ et $cy=0$.
 
\emph{2.} En effet, lorsque $ay=0$, alors $a=ca$ avec $c\in(0:y)$.
\end{proof}
%-----------------end proof------------------

%:     Fact{factLsdzCo}----------------

En utilisant le point \emph{2} du lemme \ref{lem.platsdz} et le 
fait qu'un \itf \idm est engendré par un \idm (lemme
\ref{lem.ide.idem}) on obtient le résultat qui suit.

\begin{fact}\label{factLsdzCo}
Soit  $\gA$ un anneau dans lequel l'annulateur de tout \elt est \tfz.
%-----------------begin enum------------------
\begin{enumerate}
\item $\gA$ est \lsdz \ssi il est \qiz.
\item $\gA$ est \sdz \ssi il est intègre.
\end{enumerate}
%-----------------end enum------------------
En particulier, un \cori \lsdz est \qiz.
\end{fact}
%--- end-fact--------------------

On notera que le point \emph{2} est évident en \clamaz, où 
l'hypothèse \gui{l'annulateur de tout 
\elt est \tfz} est superflue.

%%%%%%%%%%%%%%%%%%%%%%%%%%%%%%%%%%%%%%%%%%%%%%%%%%%%%%%%%%%%%%%%%%%
%--  Sec Idéaux plats \tf ---    
\section{Idéaux plats \tfz} 
\label{secIplatTf}
%-----------------------------------------

 On étudie maintenant la platitude pour les \itfsz.
En \clamaz, la proposition suivante est un corolaire \imd de la proposition 
\ref{prop.mtf loc plat1}. En \comaz, il est \ncr de fournir une nouvelle preuve, 
qui donne des informations \algqs de nature différente de celles 
données dans la preuve  de la proposition \ref{prop.mtf loc plat1}. En effet, on 
ne fait plus les mêmes hypothèses concernant le \crc discret des 
choses.
%:     Proposi{prop.itf plat local}  
\begin{proposition} 
\label{prop.itf plat local}{\em (Idéaux \tf plats sur un \aloz)}\\
Soit $\gA$ un \aloz, $x_1$, \ldots, $x_n\in\gA$ et $\fa=\gen{x_1,\ldots,x_n}$.
%-----------------begin item------------------
\begin{enumerate}
\item  Si $\fa$ est principal, il est 
engendré par l'un des $x_j$. \emph{(Bézout toujours trivial sur un anneau local).}
\item  Si $\fa$ est plat, il est 
principal, engendré par l'un des~$x_j$.
\item  Supposons que $\gA$ soit discret, ou que l'on  ait un test pour 
répondre à la question \gui{$x$ est-il \ndzz?}. 
Alors, un \itf est plat \ssi il est libre de rang $0$ ou~$1$.
\end{enumerate}
%-----------------end item------------------
\end{proposition}
%--- end-proposition---------------------------
%-----------------begin proof------------------
\begin{proof}
\emph{1.}  On a $\fa=\gen{x_1,\ldots,x_n}=\gen{z}$, $z=a_1x_1+\cdots+a_nx_n$, 
$zb_j=x_j$,
\linebreak 
donc $z(1-\sum_ja_jb_j)=0$. 
Si $1-\sum_ja_jb_j$ est \ivz, $\fa=0=\gen{x_1}$.
Si~$a_jb_j$ est \iv $\fa=\gen{x_j}$.

\smallskip 
\emph{2.}  On considère la \syzy $x_2x_1+(-x_1)x_2=0$. Soit 
\smash{$G=\cmatrix{ a_1 & \ldots & a_n \cr b_1 & \ldots & b_n}$} une matrice telle que  
$G\cmatrix{ x_1 \cr\vdots \cr x_n}=\cmatrix{ x_1 \cr  x_2}$  
et  
$[\,x_2\;-x_1\,]\,G= [\, 0\,\,\, 0\,]$.  \\
Si  $a_1$ est \ivz, l'\egt  $a_1x_2=b_1x_1$ 
montre que 
$\fa=\gen{x_1,x_3,\ldots,x_n }$.
\\ 
\hbox{Si $1-a_1$} est \ivz, l'\egt $a_1x_1+\cdots+a_nx_n=x_1$ montre que\linebreak 
l'on~a~ 
$\fa=\gen{x_2,x_3,\ldots,x_n }$. \\
On termine par \recu sur~$n$.
\\
\emph{3.}   Résulte de \emph{2} et du  lemme \ref{lemIdps}, point~\emph{3}.
\end{proof}
%-----------------end proof------------------

Rappelons qu'un \itf $\fa$ d'un anneau $\gA$ est dit {\em \lopz} s'il 
existe des \moco $S_1$, \ldots, $S_n$ de $\gA$ tels que chaque $\fa_{S_j}$ est 
principal dans $\gA_{S_j}$.
La proposition qui suit montre que \emph{tout \itf plat est \lopz}.
Sa \dem est directement issue de celle donnée dans le cas local.

%:   -- Proposition{prop.itfplat}--
\begin{proposition} 
\label{prop.itfplat}\relax 
{\em (Idéaux de type fini plats sur un anneau quelconque)}\\
Tout \itf plat est \lopz. Plus \prmtz,
\linebreak 
 si
 $\fa=\gen{x_1,\ldots,x_n}\subseteq\gA$, \propeq
%-----------------begin item------------------
\begin{enumerate}
\item  L'idéal  $\fa$ est un \mplz.
\item  Après \lon en des \moco convenables, l'idéal~$\fa$  est plat et principal.
\item  Après \lon en des \eco convenables, l'idéal  $\fa$  est 
plat et principal,
engendré par l'un des~$x_i$.
\end{enumerate}
%-----------------end item------------------
\end{proposition}
%--- end-proposition----------------------------------------
%-----------------begin proof------------------
\begin{proof}
On a évidemment \emph{3} $\Rightarrow$ \emph{2}. On a \emph{2} $\Rightarrow$ \emph{1}
par le \plgref{plcc.plat}. Pour montrer \emph{1} $\Rightarrow$ \emph{3} on reprend la 
preuve du point \emph{2} de la proposition~\ref{prop.itf plat local}.
On considère la \syzy $x_2x_1+(-x_1)x_2=0$. Soit 
\smashbot{$G=\cmatrix{ a_1 & \ldots & a_n \cr b_1 & \ldots & b_n}$} une matrice telle que 
%:2018  \cmatrix{ x_1 \cr  x_2}
$G\cmatrix{ x_1 \cr\vdots \cr x_n}=\cmatrix{ x_1 \cr  x_2}$  
%\linebreak 
et  
$[\,x_2\;-x_1\,]\,G= [\, 0\,\,\, 0\,]$.  Avec le localisé $\gA[1/a_1]$  l'\egt $a_1x_2=b_1x_1$ 
montre que $\fa=_{\gA[1/a_1]}\gen{x_1,x_3,\ldots,x_n }$. Avec le localisé 
$\gA[1/(1-a_1)]$  l'\egt $a_1x_1+\cdots+a_nx_n=x_1$ montre que 
$\fa =_{\gA[1/(1-a_1)]} \gen{x_2,x_3,\ldots,x_n }$. 
On termine par \recu sur~$n$.
\end{proof}
%-----------------end proof------------------

%:   subsec{Anars}
\subsec{Anneaux \aris et \adpsz}

La \dfn suivante des \adpsz, basée sur la platitude, 
est due à Hermida et S\'anchez-Giralda
\cite{HS}.

%-- Definition{defArit}---------- 
\begin{definition} 
\label{defArit} \emph{(Anneaux \arisz)} %\\
Un anneau $\gA$ est dit \ixd{arithmétique}{anneau} si tout \itf est
\lopz.
%-----------------end enum------------------
\end{definition}
%--- end-definition------------------------------------

%:      propdef{prop.itfplat 2}-----     
\begin{propdef} 
\label{prop.itfplat 2}\emph{(Anneaux de Prüfer)}\\
\Propeq
%-----------------begin item------------------
\begin{enumerate}
\item [1a.] Tout \itf de $\gA$ est plat. % anciennement a1
\item [1b.] Tout \id de $\gA$ est plat. % anciennement a2
\item [1c.] Pour tous \itfs $\fa$ et $\fb$ de $\gA$, l'\ali
canoni\-\hbox{que $\fa\te\fb\to\fa\fb$} est un \isoz. 
\item [2a.] L'anneau $\gA$ est \lsdz et  \ariz. % anciennement b1
\item [2b.] L'anneau $\gA$ est réduit et  \ariz. % anciennement b2
\end{enumerate}
Un anneau  vérifiant ces \prts est appelé \ixx{anneau}{de Prüfer}.
\index{Prufer@Prüfer!anneau de ---}
%-----------------end item------------------
\end{propdef}
%--- end-corollary--------------------------

%-----------------begin proof------------------
\begin{proof}
L'\eqvc entre  \emph{1a} et \emph{1c} est donnée par le \thref{thplatTens} (point \emph{\ref{i3thplatTens}}). L'\eqvc de \emph{1a} et \emph{1b} est \imdez. On sait déjà que  
\emph{1a} $\Rightarrow$ \emph{2a}, et l'implication \emph{2a} $\Rightarrow$ 
\emph{2b}  est claire. 

\emph{2b} $\Rightarrow$ \emph{2a.}
Soient $x$, $y$ tels que $xy=0$. Il existe $s$, $t$ avec $s+t = 1$,
$sx \in \gen {y}$ et~$ty \in \gen{x}$. Donc $sx^2 = 0$ et $ty^2 = 0$
puis ($\gA$ réduit) $sx = ty = 0$.

\emph{2a} $\Rightarrow$ \emph{1a}. 
Après des \lons convenables, l'\id devient principal,
et donc plat, puisque l'anneau est \lsdzz. On termine par le \plgref{plcc.plat} 
pour les modules plats. 
\end{proof}
%-----------------end proof------------------

%:   subsec{Principe local-global}
\subsec{Principe local-global}

Différentes notions introduites précédemment sont locales au 
sens du \plgc suivant. Les preuves sont basées sur le \plg de base et laissées \alecz.\iplg

%:  --- Principe local-global concret{plcc.arith}---
\begin{plcc} 
\label{plcc.arith}\relax \emph{(Anneaux \arisz)}\\
Soient $S_1$, $\ldots$, $S_n$  des \moco d'un anneau $\gA$ et $\fa$ un \id de~$\gA$.   
On a les \eqvcs suivantes.
%-----------------begin enum------------------
\begin{enumerate}
\item  L'\id $\fa$ est \lop \ssi chacun des $\fa_{S_i}$ est \lopz.
\item  L'anneau $\gA$ est \lsdz \ssi chacun des $\gA_{S_i}$ est \lsdzz.
\item  L'anneau $\gA$ est \ari \ssi chacun des $\gA_{S_i}$ est \ariz.
\item  L'anneau $\gA$ est un \adp \ssi chacun des~$\gA_{S_i}$ est un \adpz.
\end{enumerate}
%-----------------end enum------------------
\end{plcc}
%--- end-plcc-----------------------------------------

%:   subsec{Machinerie \lgbe  des \anars}
\subsec{Machinerie \lgbe  }

Un ensemble ordonné $(E,\leq)$ est dit \ixc{totalement ordonné}{ensemble ---}
si pour tous $x$, $y$ on a $x\leq y$ ou $y\leq x$. 
A priori on ne le suppose pas discret et l'on n'a donc pas de test pour l'in\egt stricte.

 Pour les \alosz,  la proposition~\ref{prop.itf plat local} donne le résultat suivant.

%--- Lemme{lemme:anarloc}------------  
\begin{lemma}\label{lemme:anarloc}\relax \emph{(Anneaux \aris locaux)}
%-----------------begin enum------------------
\begin{enumerate}
\item  Un anneau $\gA$ est local et \ari \ssi pour 
\linebreak tous~$a$,~$b\in\gA$,
on~a:~$a \in b\gA$ ou~$b \in a\gA$. De manière \eqvez, tout
\itf est principal et l'ensemble des \itfs
est totalement ordonné pour l'inclusion.
\item Soit $\gA$ un \anar local. Pour deux \ids arbitraires~$\fa$ et~$\fb$,
si~$\fa$ n'est pas contenu dans $\fb$, alors 
   $\fb$ est contenu dans $\fa$.
Donc en \clamaz, \gui{l'ensemble} de tous les \ids est totalement ordonné
pour l'inclusion.
\end{enumerate}
%-----------------end enum------------------
\end{lemma}
%--- end-Lemma-----------------------------------------

Ainsi, les \alos \aris sont la même chose que les anneaux de Bézout locaux.
Ils ont déjà été étudiés dans la section \ref{secBézout} \paref{secpfval}.

La facilité à démontrer des \prts
pour les \anars tient en grande partie à la machinerie \lgbe suivante.

%:   Machinerie loc glob  des an ariths
\mni{\bf Machinerie \lgbe des \anars}
\label{MetgenAnar}\imla
\\
\emph{Lorsque l'on doit prouver une \prt concernant un \anar
et qu'une famille finie d'\elts $(a_i)$ de l'anneau intervient dans le calcul,
on commence par démontrer le résultat dans le cas local.
On peut donc supposer que les \ids $\gen{a_i}$ 
sont totalement ordonnés par inclusion. 
Dans ce cas la preuve est en \gnl très simple.
Par ailleurs, puisque l'anneau est \ariz, on sait que l'on peut se ramener
à la situation précédente après \lon en un nombre fini
d'\ecoz. On peut donc conclure si la \prt à démontrer
obéit à un \plgcz.}

%\newpage 
\medskip Voici une application de cette machinerie.
%:2015
%:     Proposition{propIddsAnar}
\begin{proposition}\label{propIddsAnar} \emph{(Idéaux déterminantiels  sur un \anarz)}
\\
Soit $\gA$ un \anar \cohz, $M$ une matrice $\in\Ae{n\times m}$. 
\\ Notons $\fd_k=\cD_{\gA,k}(M)$ les \idds de $M$\footnote{On peut se limiter à $k\in\lrb{0..p+1}$ avec $p=\inf(m,n)$.}. 
Il existe des \itfs $\fa_1,\ldots,\fa_p$ vérifiant
$$\preskip.4em %\postskip.4em 
\fd_1=\fa_1  ,\; \fd_2=\fd_1\fa_1\fa_2,\; 
\fd_3=\fd_2\fa_1\fa_2\fa_3,\; \ldots
 $$
\end{proposition}
%--------- fin proposition -------------------------------------

\vspace{.02em}
\begin{proof} Soit $\fb_k=(\fd_k:\fd_{k-1})$ pour tout $k$, puis $\fc_k=\fb_1\cap\cdots\cap\fb_k$ et $\fa_k=(\fc_k:\fc_{k-1})$ pour $k\geq 1$. Ce sont tous les \itfsz. On a $\fb_1=\fd_1$, puis $\fb_k\fd_{k-1}=\fd_k$ et $\fa_k\fc_{k-1}=\fc_k$ pour~$k\geq 1$ parce que l'anneau est \ari \cohz. La suite des \ids $(\fc_k)_{k\geq 1}$ est décroissante par \dfnz.
\\
La proposition résulte de l'\egt \fbox{$\fc_k\fd_{k-1}=\fd_k$}
(claire pour $k=  1$).
\\ 
Si $\gA$ est un \anar local la matrice admet une forme réduite de Smith
(proposition \ref{propPfVal}). Notons $p=\inf(m,n)$.
\\
L'\algo qui produit la forme réduite de Smith dans le cas local et la
machinerie \lgbe des \anars précédente nous fournissent un
\sys d'\eco $(s_1,\ldots,s_r)$ tel que, sur chaque anneau $\gA[1/s_i]$, la matrice $M$
admet une forme réduite de Smith avec la sous-matrice diagonale $\Diag(c_1,c_2,\dots,c_p)$ et $c_1\mid c_2\mid \dots\mid c_{p}$. En outre, pour $k\geq 1$, $\fd_k=\gen{c_1\cdots c_k}$.
\\
Il suffit de prouver l'\egt encadrée après \lon en ces \ecoz.
Le fait que $\fb_k\fd_{k-1}=\fd_k$ implique $\fc_k\fd_{k-1}\subseteq \fd_k$.
\\
Il faut montrer l'inclusion réciproque. Plus \prmtz, montrons pour tout $k\geq 1$ que $c_k\in\fc_k$, ce qui implique $\fc_k\fd_{k-1}\supseteq \fd_k$.
On a  $c_k\fd_{k-1}=\fd_k$, donc
$c_k\in\fb_k$. Par ailleurs $c_k$ est multiple des $c_i\in\fb_i$ pour  $i\leq {k-1}$, donc $c_k\in\fb_1\cap\cdots\cap\fb_{k-1}$.
On a donc bien $c_k\in\fc_k$.
\end{proof}
\rem Si $\gA$ est un \ddp (un \anar intègre), on peut voir que la suite des \ids $\fb_k$ est décroissante jusqu'à ce que~$\fd_r$ s'annule. Si l'on prend pour $\gA$ un produit de \ddpsz, on peut constater qu'il est en \gnl faux que la suite des \ids $\fb_k$ soit décroissante.
\eoe

Nous reviendrons plus longuement sur les \anars et 
les \adps dans le chapitre~\ref{ChapAdpc}.

%%%%%%%%%%%%%%%%%%%%%%%%%%%%%%%%%%%%%%%%%%%%%%%%%%%%%%%%%%%%%%%%%%%
%%%%%%%%%%%%%%%%%%%%%%%%%%%%%%%%%%%%%%%%%%%%%%%%%%%%%%%%%%%%%%%%%%%
%%%%%%%%%%%%%%%%            Algèbres plates          %%%%%%%%%%%%
%%%%%%%%%%%%%%%%%%%%%%%%%%%%%%%%%%%%%%%%%%%%%%%%%%%%%%%%%%%%%%%%%%%
%%%%%%%%%%%%%%%%%%%%%%%%%%%%%%%%%%%%%%%%%%%%%%%%%%%%%%%%%%%%%%%%%%%
\section{Algèbres plates}
\label{secAlplates}

En langage intuitif, une \Alg $\gB$ est plate  lorsque les
\slis sur $\gA$ sans second membre n'ont \gui{pas plus} de solutions dans 
$\gB$ que dans~$\gA$, et elle est fidèlement plate si cette affirmation
est vraie \egmt des \slis avec second membre. Précisément on adopte les \dfns suivantes. 

%:     Definition{defiAlgPlate}
\begin{definition}\label{defiAlgPlate} Soit $\rho:\gA\to\gB$ une \Algz.
\begin{enumerate}
\item $\gB$ est dite \emph{plate (sur $\gA$)} 
lorsque toute \rde $\gB$-\lin entre 
\elts de $\gA$ est une combinaison $\gB$-\lin de  
\rdes $\gA$-\lins entre ces mêmes \eltsz.
Autrement dit,
pour toute forme \lin $\psi:\gA^n\to\gA$, 
on réclame que 
$$\preskip.4em \postskip.4em 
\Ker \rho\ist(\psi) =\gen{\rho (\Ker \psi )}_\gB. 
$$
On dira aussi que \emph{l'\homo d'anneaux $\rho$ est plat}.%
\index{plate!algèbre ---}\index{algèbre!plate}%
\index{plat!homomorphisme d'anneaux ---}
\item Une \Alg \emph{plate} $\gB$ est dite \emph{\fptez} si 
pour toute forme \lin $\psi:\gA^n\to\gA$ et tout $a\in\gA$, 
lorsque l'équation $\psi(X)=a$ 
admet une solution dans 
$\gB$ (i.e. $\exists X\in\gB^n,\,\big(\rho\ist(\psi)\big)(X)=\rho(a)$), 
alors elle admet une solution dans $\gA$. 
\\
On dira aussi que \emph{l'\homo d'anneaux $\rho$ est \fptz}.
\end{enumerate}
\index{fidèlement plate!algèbre ---}
\index{fidèlement plat!homomorphisme d'anneaux ---}
\index{algèbre!fidèlement plate}
\end{definition}

Pour une \Alg \fptez, en considérant le cas où $n=1$ 
\linebreak 
et $\psi=0$, 
on voit que $\rho(a)=0$ implique $a=0$. Ainsi, $\rho$ est un \homo injectif.
On dit donc que $\gB$ est une \emph{extension \fptez} de $\gA$.
On peut alors identifier $\gA$ à un sous-anneau de $\gB$ et la condition
sur l'équation \lin avec second membre est plus simple à formuler:
c'est exactement la même équation que l'on cherche à résoudre dans 
$\gA$ ou~$\gB$.

%:     Fact{factplatplat}
\begin{fact}\label{factplatplat} ~ \\
Une \Alg $\gB$ est plate \ssi $\gB$ est un \Amo plat. 
\end{fact}
\begin{proof}
Exercice de traduction laissé \alecz.
\end{proof}

{\bf Exemples fondamentaux. }
En voici dans le lemme suivant.

%:     Lemma{lemExAlgPlFiPl}
\begin{lemma}\label{lemExAlgPlFiPl}~
\begin{enumerate}
\item Un morphisme de \lon $\gA\to S^{-1}\gA$ donne une \Alg
plate.
\item Si $S_1$, \ldots, $S_n$ sont des \moco de $\gA$ et si $\gB=\prod_i\gA_{S_i}$
l'\homo \gui{diagonal} canonique $\rho:\gA\to\gB$ donne une \alg
\fptez.
\item Si $\gk$  est \zedrz, toute \klg $\gL$ est plate.
%Si en outre $\gK$ est facteur direct dans $\gL$, $\gL$ est une extension
%\fptez.
%
\end{enumerate}
\end{lemma}
\begin{proof}
\emph{1.} Voir le fait~\ref{fact.transporteur} ou les faits  \ref{factplatplat}
et~\ref{fact.plat}.
\\
\emph{2.} Cela résulte du \plg de base (on pourrait même dire que \und{c'est}
le \plg de base).\iplg
\\
\emph{3.} Résulte de  \ref{factplatplat}
et de ce que tout \Kmo est plat (\thref{propEvcPlat}). 
\end{proof}

\rems Concernant le point \emph{3} du lemme précédent.
 
1) 
 Il semble difficile
 de remplacer dans l'hypothèse $\gk$ par un corps (de Heyting) que l'on ne 
 suppose pas \zedz.
 
2)
Voir le \thref{thSurZedFidPlat} pour la question \fptez. 
\eoe

%:     proposition{lemAlgPlate}

\medskip 
Dans la proposition qui suit, analogue des propositions
\ref{propCoh1} (pour les \corisz) et \ref{propPlat1} (pour les \mplsz), on passe d'une \eqn à un \sys d'\eqnsz.
Pour alléger le texte,
 on fait comme si l'on avait une inclusion $\gA\subseteq\gB$
(même si $\gB$ est seulement supposée plate), autrement
dit on ne précise pas que quand on passe dans $\gB$, 
tout doit être transformé au moyen de l'\homo $\rho:\gA\to\gB$.

\begin{proposition}\label{lemAlgPlate} 
Soit $M\in\Ae{n\times m}$,  $C\in\Ae{n\times 1}$ et $\gB$ une \Alg plate. 
\begin{enumerate}
\item Toute solution
dans $\gB$ du \sli \hmg $MX=0$ est combinaison $\gB$-\lin 
de solutions dans~$\gA$.
\item Si en outre $\gB$ est \fptez, et si le \sys $MX=C$ admet une solution
dans $\gB$, il admet une solution dans~$\gA$.
\end{enumerate}
\end{proposition}
\begin{proof} Les \dfns des \Algs plates et \fptes concernent les \slis avec une seule équation.
Pour résoudre un \sli \gnl on applique la technique usuelle: on commence par résoudre la première équation,
puis on porte la solution \gnle de la première équation dans la seconde,
et ainsi de suite.
\end{proof}
%

%:     Proposition{propIdsAlgPlate}
\begin{proposition}\label{propIdsAlgPlate}~ \\
Soit $\gA\vers{\rho}\gB$ une \Alg plate et $\fa,\,\fb$ deux \ids de~$\gA$.
\begin{enumerate}
\item L'\Bli naturelle $\rho\ist(\fa)\to \rho(\fa)\gB$ est un \isoz.
\end{enumerate}
Dans la suite on identifie $\rho\ist(\fc)$ avec  l'\id $\rho(\fc)\gB$
pour tout \id $\fc$ de $\gA$.
\begin{enumerate} \setcounter{enumi}{1}
\item On a $\rho\ist(\fa\cap \fb)=\rho\ist(\fa)\cap\rho\ist(\fb)$. 
\item Si en outre $\fa$ est \tfz, on a 
$\rho\ist(\fb: \fa)=\big(\rho\ist(\fb):\rho\ist(\fa)\big)$. 
\end{enumerate}
\end{proposition}
\begin{proof}
Les deux premiers points résultent des faits analogues concernant
les \mpls (\thref{thplatTens} point \emph{\ref{i4thplatTens}} et corolaire~\ref{cor3PlatTens}).
\\
\emph{3.} Si $\fa=\gen{\an}$, alors $\fb:\fa=\bigcap_i(\fb:a_i)$, 
donc vu le point~\emph{2} on est ramené au cas d'un \idp $\gen{a}$. 
On considère alors la suite exacte  

\snic{0\to\,\fb: a\,\llra\gA\vvers a\gA\sur\fb,}

%\sni
on fait le produit tensoriel par $\gB$, et l'on obtient la suite exacte
(utiliser la platitude et le fait~\ref{factSuitExTens})

\snic{0\to\rho\ist(\fb: a)\llra\gB\vvvers {\rho(a)}\gB\sur{\rho\ist(\fb)},}

%\sni
qui donne le résultat voulu. 
\end{proof}
%

%:     theorem{thExtPlat}-------   
\begin{theorem} 
\label{thExtPlat} 
Soit $\rho:\gA\rightarrow \gB$  une \algz.  
\Propeq
%-----------------begin enum------------------
\begin{enumerate}
\item \label{i1thExtPlat} $\gB$ est une \Alg plate.
\item \label{i2thExtPlat} $\gB$ est un \Amo  plat.
\item \label{i3thExtPlat} Pour tout \Amo plat $M$, le \Amo  $\rho\ist (M)$  est plat.
\item \label{i4thExtPlat} Pour tout \itf $\fa$ de $\gA$, l'\Ali canonique

\snic{\gB\te_\gA\fa\simeq\rho\ist(\fa)\to\fa\gB}  

%\sni
est un \isoz.
\item \label{i5thExtPlat} Pour tous \Amos $N\subseteq M$ l'\Bli $\rho\ist(N)\to\rho\ist(M)$  est injective.
\item  \label{i6thExtPlat} Pour toute \Ali $\psi:M\to P$,  
l'\Bli naturelle $\rho\ist\big(\Ker(\psi)\big)\lora\Ker\big(\rho\ist (\psi)\big)$
est un \isoz.
\item \label{i7thExtPlat} Pour toute suite exacte de \Amos $M\vers{f} N\vers{g} P$ 
la suite  

\snic{\rho\ist(M)\vvvers{\rho\ist(f)} \rho\ist(N)\vvvers{\rho\ist (g)} \rho\ist (P)}
 
%\sni
est une suite exacte de \Bmosz.
\end{enumerate}
\end{theorem}
%--- end-theorem-------------------
Le point \emph{\ref{i5thExtPlat}} permet d'identifier  $\rho\ist(P)$ à un 
sous-\Bmo de $\rho\ist(Q)$ chaque fois que l'on a deux \Amos $P\subseteq Q$
et que $\gB$ est plate sur~$\gA$.
\begin{proof} 
\Llec vérifiera que les \eqvcs sont claires d'après ce que l'on sait déjà
(fait \ref{factplatplat},  \thref{thplatTens},
corolaire~\ref{corPlatTens}). On notera que la proposition~\ref{lemAlgPlate} donne le point \emph{\ref{i6thExtPlat}} dans le cas d'une \ali entre modules libres de rang fini.
\end{proof}

%:     Proposition{propPlateHom}
%\medskip 
La proposition suivante \gns les propositions \ref{fact.hom egaux}
et \ref{fact.homom loc pf}.
 
\begin{proposition}\label{propPlateHom}
Soient $\rho:\gA\to\gB$ une \Alg plate et $M$, $N$ 
\linebreak 
des \Amosz.
Si $M$ est \tf (resp. \pfz), l'\Bli naturelle 

\snic{\rho\ist\big(\Lin_\gA(M,N)\big)\to \Lin_\gB\big(\rho\ist(M),\rho\ist(N)\big) 
}

%\sni
est injective (resp. est un \isoz). 
\end{proposition}
%------------------------
%
\begin{proof}
On considère une suite exacte 

\snic{\qquad K\lora \gA^k\lora M\to 0,\qquad(*) 
}

%\sni
correspondant
au fait que $M$ est \tf (si $M$ est \pf le module $K$ est lui aussi libre de rang fini).  \\
Notons
$M_1=\rho\ist(M)$, $N_1=\rho\ist(N)$ et $K_1=\rho\ist(K)$.
Nous avons tout d'abord la suite
exacte

\snic{\qquad K_1\lora \gB^k\lora M_1\to 0.\qquad(**) 
}

%\sni
Ensuite on obtient les suites exactes ci-dessous. La première vient de 
$(*)$, la dernière vient de $(**)$ et la seconde résulte de la première
par \eds puisque $\gB$ est plate sur~$\gA$. 

\[ \arraycolsep2pt
\begin{array}{cccccccccc} 
0  &\to& \Lin_\gA(M,N) &\to& \Lin_\gA(\gA^k,N)\simeq N^k &\to&\Lin_\gA(K,N)
\\[3mm] 
0  &\to&\rho\ist\big(\Lin_\gA(M,N)\big)&\to&\rho\ist(\Lin_\gA(\gA^k,N)\simeq N_1^k   &\to&   \rho\ist\big(\Lin_\gA(K,N)\big)\\[1mm] 
 &&\dar{} &&\dar{} &&\dar{}  \\[1mm] 
0  &\to& \Lin_\gB(M_1,N_1)&\to&\Lin_\gB(\gB^k,N_1)\simeq N_1^k&\to&\Lin_\gB(K_1,N_1) 
 \end{array}
\]
 En outre, on a des \Blis \gui{verticales} naturelles de la deuxiè\-me
vers la troisième suite exacte, et les diagrammes commutent.
La deuxième flèche verticale est un \iso (l'identité de $N_1^k$ après
les identifications canoniques).
Ceci implique que la première flèche verticale (l'\Bli qui nous intéresse)
est injective.
\\
Si $M$ est \pf et si $K\simeq \gA^{\ell}$, les deux 
\Bmos à droite sont isomorphes
à $N_1^{\ell}$ et la flèche verticale
correspondante est un \isoz. Ceci implique que la première 
flèche verticale est un \isoz.
\end{proof}

Rétrospectivement la \dem donnée pour la proposition 
 \ref{fact.homom loc pf} semble bien compliquée.
La nouvelle \dem donnée ici dans un cadre plus \gnl 
est conceptuellement plus simple.

%-% ENTRE NOUS
\entrenous{1. Se pose sérieusement la question de jeter à la poubelle,
ou renvoyer en exo, la \dem \gui{compliquée} de la proposition \ref{fact.homom loc pf}.

2. Dans la proposition précédente il y a peut être des choses à dire en plus lorsque $\gB$ est \fptez. Si oui le signaler dès maintenant.

3. Peut être aussi des choses en rapport avec la \cohc (rappel: prop. \ref{propAliCoh})? 

4. Eisenbud signale son tho 2.13 et ses exos 3.3, 4.11, 4.12, 4.13 et 19.4 comme applications surprenantes de  la proposition \ref{propPlateHom}
}
%-% Fin ENTRENOUS

%%%%%%%%%%%%%%%%%%%%%%%%%%%%%%%%%%%%%%%%%%%%%%%%%%%%%%%%%%%%%%%%%%%
\section{Algèbres fidèlement plates}
\label{secAlfidplates}

On a déjà dit que si $\gA\vers\rho\gB$ est une \alg \fptez, $\rho$
est injectif. 
\perso{autres \prts analogues? Radical de Jacobson?}
Il est clair \egmt que $\rho$  \emph{réfléchit les
unités}, \cad que
\index{reflechit les@réfléchit les unités!\homo qui ---}
$$\preskip.4em \postskip.4em 
\rho(a)\in\gB\eti\;\Longrightarrow\;a\in\Ati. 
$$
 Voici maintenant quelques \prts \carasz. Dans la suite on retiendra
l'\eqvc des points \emph{1}, \emph{2a},  \emph{3a} et \emph{4.} 

%:     theorem{thExtFidPlat}-------   
\begin{theorem} 
\label{thExtFidPlat} \emph{(Caractérisation des \algs \fptesz)}\\
Soit $\rho:\gA\rightarrow \gB$ une \alg plate.  
\Propeq
%-----------------begin enum------------------
\begin{enumerate}
\item [1.] %\label{i1thExtFidPlat} 
L'\alg $\gB$ est  \fptez.
\item [2a.] %\label{i2athExtFidPlat} 
L'\homo $\rho$ est injectif, et en identifiant $\gA$
à un sous-anneau de $\gB$, pour tout \itf $\fa$ de $\gA$ on a

\snic{ \fa\gB\,\cap\,\gA\,=\,\fa.}

\item [2b.] %\label{i2bthExtFidPlat}
Même chose avec un \id arbitraire de $\gA$.
\item [3a.] %\label{i3athExtFidPlat} 
Pour tout \itf $\fa$ de $\gA$ on a l'implication

\snic{ 1_\gB\in \rho\ist(\fa)\;\Longrightarrow\;1_\gA\in\fa.} 
\item [3b.]  %\label{i3bthExtFidPlat} 
Pour tout \itf $\fa$ de $\gA$, si  
$\rho\ist(\gA\sur\fa)=0$, alors $\gA\sur\fa=0$.
\item [3c.] %\label{i3cthExtFidPlat}  
Pour tous \Amos $N\subseteq M$, si $\rho\ist(N)=
\rho\ist(M)$, alors~$N=M$.
\item [3d.] %\label{i3dthExtFidPlat} 
Pour tout \Amo $M$, si  $\rho\ist(M)=0$, alors~$M=0$.
\item [3e.] %\label{i3ethExtFidPlat} 
Pour tout \Amo $M$ l'\Ali naturelle $M\to\rho\ist(M)$ est injective.
\item [4.] %\label{i4thExtFidPlat}
L'\eds de $\gA$ à $\gB$ \emph{réfléchit les suites exactes}.
\\
Autrement dit,  une suite arbitraire de \Amos 
$$\preskip.4em \postskip.4em 
N\vers{f} M\vers{g} P
$$ 
est exacte si la suite de \Bmos  
$$\preskip.4em \postskip.4em 
\rho\ist (N)\vvvers{\rho\ist (f)} \rho\ist (M)\vvvers{\rho\ist (g)}
 \rho\ist (P) 
$$
est exacte.
\end{enumerate}
\end{theorem}
%--- end-theorem-------------------
%
\begin{proof}
Le point \emph{1} implique que $\rho$ est injectif.
Une fois ceci acquis,  \emph{2a} est une simple reformulation de \emph{1},
et  \emph{2a} est facilement \eqv à~\emph{2b.}

\emph{3a} $\Rightarrow$ \emph{1}.
On commence par remarquer que l'implication est encore valable si l'on remplace l'\itf $\fa$ par un \id arbitraire $\fc$. En effet, \hbox{si $1\in\rho\ist(\fc)$}
on aura \egmt $1\in\rho\ist(\fc')$ pour un \itf $\fc'$ contenu dans $\fc$.
\\
Soit maintenant $\fa=\gen{a_1,\ldots,a_n}$ et $c\in\gA$. 
L'équation $\sum_ia_ix_i=c$
admet une solution \ssi $c\in\fa$, i.e. $1\in(\fa:c)_\gA$. 
Puisque $\gB$ est plate, on a
$\big(\rho\ist(\fa):\rho(c)\big)_\gB=\rho\ist(\fa:c)$ (proposition \ref{propIdsAlgPlate}).
Si $\sum_i\rho(a_i)y_i=\rho(c)$ admet une solution dans $\gB$,
alors $1\in\big(\rho\ist(\fa):\rho(c)\big)_\gB$, donc l'hypothèse \emph{3a}
implique que $1\in(\fa:c)$, i.e. que 
$\sum_ia_ix_i=c$
admet une solution dans~$\gA$.

 Les implications  \emph{3e} $\Rightarrow$ \emph{3d} $\Rightarrow$ \emph{3b} sont triviales.

\emph{3d} $\Rightarrow$ \emph{3c}. On considère le module $M/N$.
Le module $\rho\ist(N)$ 
s'identifie à un sous-module de $\rho\ist(M)$ et $\rho\ist(M\sur N)$ 
s'identifie à $\rho\ist(M)\sur{\rho\ist(N)}$. On conclut.

\emph{3c} $\Rightarrow$ \emph{3d}. On prend $N=0$.

\emph{3a} $\Leftrightarrow$ \emph{3b}. Mêmes raisonnements. 

\emph{1} $\Rightarrow$ \emph{3e}.
 On identifie $\gA$ à un sous-anneau de $\gB$. \\
 Soit $x\in M$ tel que $1\te x=0$ dans $\rho\ist(M)$. 
 Puisque $\gB$ est un \Amo plat, cette \syzy s'explique
dans le \Amo $\gB$: il existe $u_1$, \dots, $u_n\in\gB$ et $a_1$, \dots, $a_n\in\gA$ tels que 
$\sum_ia_iu_i=1$ et $a_ix=0$ pour $i\in\lrbn$. 
L'équation en les $y_i$, $\sum_ia_iy_i=1$, admet une 
solution dans~$\gB$, donc elle en admet une dans $\gA$. D'où~$x=0$.

 \emph{4} $\Rightarrow$ \emph{3d.}
On fait $N=P=0$ dans la suite $N\to M\to P$.
Elle est exacte après \eds à $\gB$, donc elle est exacte.

\emph{1} $\Rightarrow$ \emph{4.}
On suppose que la suite de \Bmos est exacte. On doit montrer
que la suite de \Amos est exacte. Tout d'abord $g\circ f=0$, car l'\Bli
$P\to\rho\ist(P)$ est injective, et les diagrammes commutent. Ensuite puisque $\gB$ est plate,
on peut identifier $\rho\ist(\Ker g)$ avec $\Ker\rho\ist(g)$
et~$\rho\ist(\Im f)$ avec $\Im\rho\ist(f)$. On est ramené au point~\emph{3c.}
\end{proof}
Vu le \thref{propEvcPlat}, on obtient comme conséquence de la \carnz~\emph{2a} le \tho suivant.

%:     Theorem{thSurZedFidPlat}
\begin{theorem}\label{thSurZedFidPlat}
Toute extension d'un \cdi ou d'un anneau \zedr est \fptez. 
\end{theorem}
%--------- fin theorem ---------------------------------------------- 
%
\begin{proof}
Par hypothèse, on a $\gk\subseteq\gA$ avec $\gk$ \zedrz. On sait que l'extension est plate
par le  \thref{propEvcPlat}. On doit montrer que si $\fa$ est un \itf
de $\gk$, alors $\fa\gA\cap\gk=\fa$. Or $\fa=\gen{e}$ pour un \idm $e$;
l'appartenance d'un \elt $x$ à un \id $\gen {e}$ ($e$ \idmz) étant \caree 
par l'\egtz~\hbox{$x = xe$}, elle est indépendante de l'anneau. Dit autrement,
pour  $e$ \idm d'un anneau $\gB\subseteq \gB'$, on a toujours $e\gB' \cap \gB = e\gB$.
\end{proof}

Comme cas particulier de la \carn \emph{3a} on obtient le corolaire suivant. 

%:     Corollary{corpropFidPlatLoc}
\begin{corollary}\label{corpropFidPlatLoc}
Soit $\rho:\gA\to\gB$ un \homo plat entre \alosz. Il est \fpt \ssi il \emph{réfléchit les unités},
i.e. \hbox{si $\rho^{-1}(\gB\eti)=\Ati$}. 
\end{corollary}

Un \homo entre anneaux locaux 
qui réfléchit les unités 
est  appelé un \emph{\homo local}.%
\index{local!homomorphisme ---}%
\index{homomorphisme!local}

Les \dems des deux faits qui suivent résultent de considérations
simples sur la préservation et sur la \gui{réflection} des suites exactes.
Les détails sont laissés \alecz.

%:     Fact{factAlgPlate}
\begin{fact}\label{factAlgPlate} \emph{(Transitivité)}
Soit $\gB$ une \Alg et $\gC$ une \Blgz.
\begin{enumerate}
\item Si $\gB$ est plate sur $\gA$ et $\gC$ plate sur $\gB$, alors
$\gC$ est plate sur~$\gA$.
\item Si $\gB$ est \fpte sur $\gA$ et $\gC$ \fpte sur $\gB$, alors
$\gC$ est \fpte sur~$\gA$.
\item Si $\gC$ est \fpte sur $\gB$ et plate sur $\gA$, alors
$\gB$ est plate sur~$\gA$. 
\item Si $\gC$ est \fpte sur $\gB$ et  sur $\gA$, alors
$\gB$ est \fpte sur~$\gA$. 
\end{enumerate}
\end{fact}

%:     Fact{factAlgPlate2}
\begin{fact}\label{factAlgPlate2} \emph{(Changement d'anneau de base)}\\
Soit $\gB$  et $\gC$ deux \Algsz, et $\gD=\gB\otimes_\gA\gC$.
\begin{enumerate}
\item Si $\gC$ est plate sur $\gA$, $\gD$ est plate sur~$\gB$.
\item Si $\gC$ est \fpte sur $\gA$, $\gD$ est \fpte sur~$\gB$.
\end{enumerate}
\end{fact} 
%

%:     PrincipeLocGlob{plcc}
\begin{plcc}\label{plcc.algfptes}\relax \emph{(Localisation en bas, algèbres plates)}\\
Soient $\rho:\gA\to\gB$ une \alg et $S_1$, $\ldots$, $S_r$  des \moco de~$\gA$.\perso{Localisation en haut. Voir le \plgref{plcc2.apf}.
Comme c'est assez compliqué il semble qu'il vaut mieux attendre la chapitre 15}
\begin{enumerate}
\item L'\alg $\gB$ est plate sur $\gA$ \ssi pour chaque~$i$, 
$\gB_{S_i}$ est plate sur~$\gA_{S_i}$.
\item L'\alg $\gB$ est \fpte sur $\gA$ \ssi pour chaque~$i$, 
l'\alg $\gB_{S_i}$ est \fpte sur~$\gA_{S_i}$. 
\end{enumerate}
\end{plcc}
\begin{proof}
On introduit l'\Alg \fpte $\gC=\prod_i\gA_{S_i}$ qui donne par \eds la
\Blg \fpte $\gD=\prod_i\gB_{S_i}$. Il reste à appliquer
les faits \ref{factAlgPlate} et~\ref{factAlgPlate2}.
\end{proof}
%

%:     theorem{propFidPlatTf}
Le \tho suivant 
 \gns les \plgcs qui affir\-ment le \crc local (au sens \cofz)
  de certaines \prts de finitude pour les modules.

\begin{theorem}\label{propFidPlatTf}
Soit $\gA\vers \rho\gB$ une \Alg \fptez.
\\
Soit $M$  un~\Amo \hbox{et $M_1=\rho\ist(M)\simeq\gB\otimes_\gA M$}.
\begin{enumerate}
\item Le \Amo $M$ est plat \ssi le \Bmo $M_1$ est plat.
\item Le \Amo $M$ est \tf \ssi le \Bmo  $M_1$ est \tfz.
\item Si le \Bmo $M_1$ est \cohz, le \Amo $M$ est \cohz.
\item Le \Amo $M$ est \pf \ssi le \Bmoz~\hbox{$M_1$} est \pfz.
\item Le \Amo $M$ est \ptf \ssi le \Bmoz~\hbox{$M_1$} est \ptfz.
\item Si le \Bmo $M_1$ est \noez, le \Amo $M$ est \noez.
\end{enumerate}
\end{theorem}
\begin{proof} Dans les points \emph{1}, \emph{2}, \emph{4}, \emph{5}, on sait déjà que n'importe quelle \eds préserve la \prt concernée.  
 Il nous reste donc à prouver les
réciproques.

\emph{1.} On considère une suite exacte $N\vers f Q \vers g P$ de \Amosz.
On veut montrer qu'elle est exacte après tensorisation par $M$.
On sait qu'elle est exacte après tensorisation par $\gB\te M$. Or $\gB\te \bullet$
réfléchit les suites exactes.

\emph{2.} On considère des \elts $y_i\in\rho\ist(M)$ ($i\in\lrbn$) qui engendrent
ce module. Ces \elts sont fabriqués comme combinaisons $\gB$-\lins
d'une famille finie d'\elts  $1\te x_j$  ($x_j\in M$, $j\in\lrbm$).
Cela implique que l'\Ali $\varphi:\Ae m\to M$ qui envoie la base canonique sur~$(x_j)_{j\in\lrbm}$
est surjective après tensorisation par $\gB$. Or $\gB$ est \fptez,
donc $\varphi$ est surjective.

\emph{3.} Soit $N=\gA x_1+\cdots+\gA x_n$ un sous-\mtf de $M$.
On considère  l'\Ali surjective $\Ae n\to N$ correspondante, on note $K$ son noyau.
La suite exacte $0\to K\to \Ae n\to N\to 0$ donne par \eds une suite exacte
(ceci parce que $\gB$ est plate).
Puisque $\rho\ist(M)$ est \cohz, $\rho\ist(K)$ est \tfz. 
Il reste à appliquer le point~\emph{2.}

\emph{4.} Même raisonnement qu'au point~\emph{3.}

\emph{5.} Un module est \ptf \ssi il est plat et \pfz.

\emph{6.} 
On considère une suite croissante $(N_k)_{k\in\NN}$ de sous-\mtfs de~$M$ et l'on étend les scalaires à $\gB$. Deux termes consécutifs $\rho\ist(N_\ell)$
et~$\rho\ist(N_{\ell+1})$ sont égaux. Puisque $\gB$ est \fptez, on a aussi
les \egts $N_\ell=N_{\ell+1}$.
\end{proof}
%

%:     theorem{propFidPlatPrAlg}
Le \tho suivant  \gns les \plgcs qui affirment le \crc local 
(au sens \cofz)  de certaines \prts de finitude pour les \algsz.

\begin{theorem}\label{propFidPlatPrAlg}~
\\
\Deuxcol{.63}{.3}{
Soit  une \Alg \fpte $\rho:\gA\to\gB$ et une \Alg $\varphi:\gA\to\gC$. 
\\
On note $\gD=\rho\ist(\gC)$
la \Blg \fpte obtenue par extension des scalaires.}
{\vspace{-5mm}
\hspace{5mm}\xymatrix @R=18pt @C=12pt
{
\gA  \ar@{->}[d]_\rho\ar@{->}[rr]^{\varphi} 
   &\ar@{~>}[d]^{\rho\ist}
   &\gC\ar@{->}[d]%^{\varphi\ist(\rho)} 
   \\                      
\gB  \ar@{->}[rr]_{\rho\ist(\varphi)} && \gD  \\
}
}

\vspace{-1mm}
Pour que $\gC$ possède une des \prts ci-dessous en tant qu'\Alg il faut et suffit que $\gD$ possède la même \prt comme \Blgz:
\begin{itemize}
\item finie (comme module),
\item \pf comme module, 
\item \stfez, 
\item plate, 
\item \fptez, 
\item \stez,
\item \spbz,
%
%\item \fntz, 
%
\item \tf (comme \algz), 
\item \pf (comme \algz). 
%
%\item  
%
%\item  
\end{itemize}
\end{theorem} 
\begin{proof} 
Les trois premières \prts sont des \prts de modules et relèvent donc du \thrf{propFidPlatTf}.

\emph{Algèbres plates, \fptesz}.
On applique les faits \ref{factAlgPlate} et~\ref{factAlgPlate2}.

\emph{Algèbres \stesz}. On a déjà l'\eqvc pour le \crc
\stfz. Si $\gB$ est libre sur $\gA$ on utilise le fait que le \discri  
se comporte bien par \edsz, et l'on conclut en utilisant le fait qu'une 
extension \fpte réfléchit les unités.
\\
Dans le cas \gnl on se ramène au cas libre par \lon en des \ecoz, 
ou bien l'on invoque le \thref{thAlgStfSpbSte}: 
une \alg \stfe est \spb \ssi elle est \stez.

\emph{Algèbres \spbsz}. On regarde le diagramme commutatif 
dans le fait \ref{factSpbEds} (attention, les noms changent). 
La flèche verticale de droite est obtenue
par \eds \fpte à partir de celle de gauche. Elles sont donc simultanément surjectives.  

\emph{Algèbres \tfz}. Le fait d'être \tf ou \pf est préservé par n'importe quelle \edsz. Voyons la réciproque. 
\\
On identifie $\gA$ à un sous-anneau de $\gB$ et
$\gC$ à un sous-anneau de $\gD$.\\
Posons $\gA_1=\varphi(\gA)$ et $\gB_1=\rho\ist(\varphi)(\gB)$.
Puisque $\gD=\gB\otimes _\gA\!\gC$ est \tf sur $\gB$, et puisque tout \elt de $\gD$
s'écrit comme combinaison $\gB$-\lin d'\elts de $\gC$, on peut écrire
$\gD=\gB_1[\xm]$ avec des $x_i\in\gC\subseteq\gD$. Ceci donne une suite exacte
$$\preskip.4em \postskip.4em 
\gB[\Xm]\vvvvers{
\rho\ist(\varphi) ,\, X_i\mapsto x_i}
\gD\lora 0. 
$$
 On va montrer que $\gC=\gA_1[\xm]$. En effet, la suite exacte ci-dessus est obtenue par \eds \fpte à partir de la suite
$$\preskip.4em \postskip.4em 
\gA[\Xm]\vvvvers{
\varphi ,\, X_i\mapsto x_i}
\gC\lora 0. 
$$ 

\emph{Algèbres \pfz}. 
\\
Commençons par une remarque \gnle \elr mais utile
sur les \algs quotients $\kuX\sur\fa$. On peut  voir $\kuX$
comme le \kmo libre ayant pour base la famille des \moms $(X^{\ual})_{\ual\in{\NN^m}}$. Si $f\in\fa$, on obtient alors l'\egt
$$\preskip.4em \postskip.4em 
f\cdot\kuX=\som_\ual (X^\ual f)\cdot\gk. 
$$
  Donc l'\id $\fa$ est le sous-\kmo de $\kuX$ engendré par tous les  
$X^\ual f$, où~$\ual$ parcourt $\NN^m$ 
et  $f$ parcourt un \sgr de~$\fa$.  
\\
Reprenons alors la \dem en continuant avec les mêmes notations
que dans le point précédent.\\ 
Supposons que $\gD=\gB_1[\xm]\simeq\aqo{\gB[\uX]}{\lfs}$.
Dans la suite on regarde une \eqn $f_j=0$ comme une \syzy entre les \moms présents dans $f_j$.\\
Puisque le \Bmo $\gD$ est obtenu par \eds plate à partir du 
\Amo $\gC$, la \rde $\gB$-\linz~$f_j$  est 
une combinaison $\gB$-\lin
de \rdes $\gA$-\lins $f_{j,k}$ (entre les mêmes \moms vus dans~$\gC$).
\\
Chaque \egt $f_{j,k}(\ux)=0$  peut aussi être lue comme une \rde $\gA$-\agq (un relateur) entre les $x_i\in\gC$.
\\
Considérons alors le sous-\Amo de $\gA[\uX]$ engendré par tous les $X^\ual f_{j,k}$.
Par \eds de $\gA$ à $\gB$ la suite de \Amos
$$\preskip.4em \postskip.4em 
0\to \som_{j,k,\ual}(X^\ual  f_{j,k})\cdot\gA\to\AuX\to\gC\to 0\eqno(*) 
$$
donne la suite exacte de \Bmos
$$\preskip.4em \postskip.4em 
0\to \som_{j,k,\ual}(X^\ual  f_{j,k})\cdot\gB\to\BuX\to\gD\to 0. 
$$
En effet, $\sum_{j,k,\ual}(X^\ual  f_{j,k})\cdot\gB=\sum_{j,k}f_{j,k}\cdot \BuX=\sum_{j}f_{j}\cdot \BuX=\fa$. Donc,
puisque l'extension est \fptez,  la suite $(*)$ est elle-même
exacte. Enfin, puisque $\sum_{j,k,\ual}(X^\ual  f_{j,k})\cdot\gA=\sum_{j,k}  f_{j,k}\cdot\AuX$,
$\gC$ est une \Alg \pfz.
\end{proof}
%

%%%%%%%%%%%%%%%%%%%%%%%%%%%%%%%%%%%%%%%%%
%:section: Exercices   
\Exercices{ 

%--- Exercise{exoPlatsLecteur}-------------
\begin{exercise}
\label{exoPlatsLecteur}
{\rm  Il est recommandé de faire les \dems non données, esquissées,
laissées \alecz,
etc\ldots
\, On pourra notamment traiter les cas suivants.
\begin{itemize}
\item \label{exooPlatsLecteur0} 
Sur un anneau de Bézout intègre, un module est 
plat \ssi il est sans torsion.
\item \label{exooothPlat1} 
Montrer le \thrf{thPlat1}.

\item \label{exocorlemIdps} 
Montrer le lemme \ref{corlemIdps}.
\item \label{exothExtPlat} 
Montrer le fait  \ref{factplatplat} et le \thref{thExtPlat}.
\item Montrer les faits \ref{factAlgPlate} et~\ref{factAlgPlate2}.
%
%\item \label{exoZerRed} 

\end{itemize}
}
\end{exercise}
%--- end -exercise-----------------------------------------

%--- Exercise{exothPlat1}-------------
\begin{exercise}
\label{exothPlat1} 
{\rm Soit  $\pi:N\to M$ une \ali surjective.
\\
 1. \emph{(Cas particulier du \thref{propPlatQuotientdePlat})}
 Si $M$ est plat,  pour tout  \mpf $P$, l'\ali naturelle
$$\preskip.3em \postskip-.2em 
\Lin_\gA(P,\pi):\Lin_\gA(P,N)\to \Lin_\gA(P,M) 
$$
est surjective.
\\
2. Supposons que  $N=\Ae{(I)}$,  module 
libre sur un ensemble discret $I$. Si la \prt précédente est
vérifiée, $M$ est plat. 

\comm En \comaz, un module arbitraire $M$ n'est pas 
\ncrt un quotient d'un module  $N=\Ae{(I)}$ comme ci-dessus, 
mais cela est vrai dans le cas où $M$ est  discret, en prenant $I=M$.
Si l'on n'exige pas $I$ discret, on peut voir l'exercice~\ref{propfreeplat}.
\eoe
}
\end{exercise}
%--- end -exercise-----------------------------------------

%--- Exercise{exoplatFitt}------- 
\begin{exercise} 
\label{exoplatFitt} 
{\rm  Soit $M$ un \Amo \tfz.  Si $M$ est plat  ses \idfs sont \idmsz. 
\perso{Montrer que la réciproque n'est pas vraie: à rajouter si l'on en est certain}

}
\end{exercise}
%--- end-exercise-----------------------------------------

%--- Exercise{exoClamlsdz}-------
\begin{exercise} 
\label{exoClamlsdz} 
{\rm  Montrer en \clama qu'un anneau est \lsdz \ssi il devient intègre après 
\lon en n'importe quel \idepz.
} 
\end{exercise}
%--- end-exercise-----------------------------------------

%-% ENTRE NOUS
\entrenous{ Un exo du type suivant serait bienvenu:
%--- Exercise{exoPlatSanTor}-------------
%\begin{exercise}
%\label{exoPlatSanTor}
{  Un exemple de module plat avec un sous-module monogène non plat
(i.e., plat n'implique pas toujours sans torsion, 
comparer avec le lemme~\ref{lem.platsdz}) 
}
%\end{exercise}
%--- end -exercise-----------------------------------------
}
%-% Fin ENTRENOUS

%--- Exercise{exoClamAri}-------
\begin{exercise} 
\label{exoClamAri} 
{\rm  Montrer en \clama qu'un anneau est \ari \ssi il devient un anneau de Bézout après \lon en n'importe quel \idepz.
} 
\end{exercise}
%--- end-exercise-----------------------------------------

%--- Exercise{exoIlops0}-------------
\begin{exercise}
\label{exoIlops0}
{\rm  L'image d'un \id \lop  par un \homo d'anneaux est un \id \lopz.
Le résultat analogue pour les \ids \ivs n'est pas toujours vrai. 
 }
\end{exercise}
%--- end -exercise-----------------------------------------

%--- Exercise{exoIlops2}--------- 
\begin{exercise} 
\label{exoIlops2} 
{\rm 
Si $\fa = \gen{x_1,\dots, x_k}$ est \lopz, alors
   $\fa^n = \gen{x_1^n,\dots, x_k^n}$.
   \\
Calculer une \mlp pour $(x_1^n,\dots, x_k^n)$
à partir d'une \mlp pour $(x_1,\dots, x_k)$. 
\\
Expliciter
l'appartenance de $x_1^{n_1}\cdots x_k^{n_k}\in\gen{x_1^n,\dots, x_k^n}$ lorsque $n=n_1+\cdots+n_k$.
} 
\end{exercise}
%--- end-exercise-----------------------------------------

%--- Exercise{exoMLP}-------------
\begin{exercise}
\label{exoMLP}
{\rm \'Etant donnés $n$ \elts dans un \anar donner un \algo
qui construit une \mlp pour ces \elts à partir de \mlps
pour uniquement des couples d'\eltsz. 
 
}
\end{exercise}
%--- end -exercise-----------------------------------------

%--- Exercise{exoIlops}-------------
\begin{exercise}
\label{exoIlops}
{\rm  On considère deux \itfs $\fa$ et $\fb$ d'un anneau $\gA$,
engendrés respectivement par $m$ et $n$ \eltsz.
Soient $f$, $g\in\AX$ de degrés $m-1$ et $n-1$ avec $\rc(f)=\fa$ et $\rc(g)=\fb$.
 
\emph{1.}
Montrer que si $\fa$ est \lopz, on a $\fa\fb=\rc(fg)$ de sorte que $\fa\fb$ est engendré par
$n+m-1$ \elts (localiser et utiliser le corolaire~\ref{corlemdArtin}~\emph{4}).
 
\emph{2.}
Montrer que si $\fa$ et $\fb$ sont \lopsz, $\fa\fb$ est \lopz. Expliquez comment construire
une \mlp pour les \coes de $fg$ à partir de deux \mlpsz, respectivement pour les
\gtrs de $\fa$ et de ceux de $\fb$.
}
\end{exercise}
%--- end -exercise-----------------------------------------

%--- Exercise{exogcdlcm}--------- 
\begin{exercise} 
\label{exogcdlcm} 
{\rm On s'intéresse à l'\egt éventuelle 
%---  equation eqexogcdlcm --------
\begin{equation}\label{eqexogcdlcm}\preskip.3em \postskip.2em
\fa\,\fb =(\fa\cap\fb)(\fa+\fb)
\end{equation}
%---------------------end equation--------------
pour deux \itfs $\fa$ et $\fb$ d'un anneau $\gA$.
  
\emph{1.}  Montrer que l'\egt est vérifiée si $\fa+\fb$ est \lopz.
Si en outre $\fa$ et $\fb$ sont \lopsz, alors $\fa\cap \fb$ est \lopz.
 
\emph{2.}  Supposons $\gA$ intègre. Montrer que si l'\egt est vérifiée lorsque
$\fa$ et $\fb$ sont des 
\idps alors l'anneau est \ariz. 
 
\emph{3.}  Montrer que \propeq
%-----------------begin item------------------
\begin{itemize}
\item $\gA$ est un \adpz.
\item $\gA$ est \lsdz et l'\eqrf{eqexogcdlcm} est vérifiée pour les \idpsz.
\item $\gA$ est \lsdz et l'\eqrf{eqexogcdlcm} est vérifiée pour les \itfsz. 
\end{itemize}
%-----------------end enum------------------
} 
\end{exercise}
%--- end-exercise-----------------------------------------

%--- Exercise{exoSECSci}-------------
\begin{exercise}
\label{exoSECSci} {\rm  (Voir aussi l'exercice \ref{exoPgcdPpcm}) 
 Soient $\fa$, $\fb$, $\fc$ des \itfsz. 
 
\emph{1.} Si $\fa+\fb$ est \lopz, alors $(\fa:\fb)+(\fb:\fa)=\gen{1}$.
 
\emph{2.} Si $(\fa:\fb)+(\fb:\fa)=\gen{1}$, alors
\hsu 
\emph{a.}~
$(\fa+\fb):\fc = (\fa:\fc)+(\fb:\fc)$;
\hsu 
\emph{b.}~ $\,\fc :(\fa~\cap~ \fb)=(\fc:\fa)+(\fc:\fb)$;
\hsu 
\emph{c.}~ $(\fa+\fb)(\fa~\cap~ \fb) = \fa~\fb$;
\hsu
\emph{d.}~ $\fc\,(\fa~\cap~ \fb)= \fc~ \fa~ \cap~ \fc~ \fb$;
\hsu 
\emph{e.}~ $\,\fc+(\fa~\cap~ \fb)=(\fc+\fa)\cap (\fc+\fb)$;
\hsu 
\emph{f.}~ $\,\fc~\cap\,(\fa+ \fb)=(\fc~\cap~ \fa)+(\fc~\cap~ \fb)$;
\hsu 
\emph{g.}~ La suite exacte ci-après {\mathrigid 2mu (où $\delta(x)=(x, -x)$ et $\sigma(y, z)=y+z$)} est scindée:
$$\preskip.2em \postskip-.2em 
0 \longrightarrow \fa~\cap~ \fb \vers{\delta} \fa \times \fb 
\vers{\sigma} \fa+\fb
\longrightarrow 0. 
$$
}
\end{exercise}
%--- end -exercise-----------------------------------------

%--- Exercise{exoGaussien}-------------
\begin{exercise}
\label{exoGaussien} (Anneaux gaussiens)
{\rm Un anneau $\gA$ est dit \emph{gaussien} lorsque pour tous \pols $f$, $g\in\AX$,
on a l'\egt $\rc(fg)=\rc(f)\rc(g)$.
 
\emph{1.} Tout \anar est gaussien (voir l'exercice \ref{exoIlops}).

\emph{2.} Un anneau intègre gaussien est un \adpz. 

\emph{3.} Un anneau réduit gaussien est un \adpz.  Un anneau \qi gaussien est un \adpc  (voir le \thrf{th.adpcoh}).
}
\end{exercise}
%--- end -exercise-----------------------------------------

%--- Exercise{exononanar}--------- 
\begin{exercise} 
\label{exononanar} (Un anneau utile pour les contre-exemples)\\
{\rm Soit $\gK$ un corps discret non trivial et $V$ un \Kev de dimension 2.
On considère la \Klg $\gA=\gK\oplus V$  définie par
 $x$, $y\in V \Rightarrow xy=0$. Montrer que tout \elt de $\gA$ est \iv ou nilpotent
 (i.e. $\gA$ est local \zedz), 
et que l'anneau est \coh mais pas \ariz. Cependant tout \itf qui contient un \elt 
\ndz est égal à $\gen{1}$, a fortiori il est \ivz. 
} 
\end{exercise}
%--- end-exercise-----------------------------------------

%:--- Exercise{exoVascon}-------------
\begin{exercise}
\label{exoVascon}
{\rm  Soit $\gA$ un \alo \coh %discret et 
\dcdz. \\
On note $\fm=\Rad\gA$
et on suppose que $\fm$ est plat sur $\gA$.

\emph{1.} Montrer que $\gA$ est intègre.

\emph{2.} Montrer que $\gA$ est un \advz.

NB: on ne suppose pas  $\gA$  non trivial. 
}
\end{exercise}
%--- end -exercise-----------------------------------------

%--- Exercise{exoQuotientPlat}-------------
\begin{exercise}\label{exoQuotientPlat}
{(Quotient plat d'un module plat: une preuve directe)} \\
{\rm
Fournir une preuve directe de l'implication suivante de le
\thref {propPlatQuotientdePlat}: soit~$M$ un \Amo plat et $K$ un 
sous-module de $M$ vérifiant $\fa M \cap K = \fa K$ pour
tout \id $\fa$ de type fini; alors $M/K$ est plat.
}

\end {exercise}
%--- end -exercise-----------------------------------------

%--- problem{propfreeplat}-------------
\begin{exercise}
\label{propfreeplat} (Tout module est quotient d'un module plat)\\
{\rm  Cet exercice commence par un long texte introductif.
On précise dans la \dfn qui suit la construction d'une somme directe de modules et celle d'un module libre dans le cas d'un ensemble d'indices non \ncrt discret. Ceci nous permet de montrer que tout module est quotient d'un module  plat (en fait un module libre, pas \ncrt \pro d'un point de vue \cofz!).

%:     definition{deffree}
{\bf Définition.} \label{deffree}
Soit $I$ un ensemble arbitraire et $(M_i)_{i\in I}$ une famille \hbox{de \Amosz}. On suppose que si $i=_Ij$, alors $M_i$ et $M_j$ sont le même ensemble\footnote{Pour la notion
générale de famille d'ensembles indexée par un ensemble arbitraire, \llec peut consulter \cite[page 18]{MRR}; la construction de la somme directe d'une famille arbitraire de \Amos se trouve pages 54 et 55.}.
On définit la \emph{somme directe $\bigoplus_{i\in I}M_i$} comme un ensemble 
quotient de l'ensemble des
sommes formelles finies $\oplus_{k\in \lrbn}x_{i_k}$, où $i_k\in I$ \hbox{et $x_{i_k}\in M_{i_k}$} pour chaque~\hbox{$k\in\lrbn$}: plus \prmtz, une telle somme formelle est définie comme étant  une famille $(k,i_k, x_{k})_{k\in \lrbn}$, où $x_k\in M_{i_k}$ pour chaque $k$.
\\
La relation d'\eqvc qui définit l'\egt sur $\bigoplus_{i\in I}M_i$ est la relation d'\eqvc engendrée par les \gui{\egtsz} suivantes:
\begin{itemize}
\item associativité et commutativité des sommes formelles: on peut réordonner la famille comme l'on veut;
\item si $i_k=_Ii_\ell$ alors les deux termes $(k,i_k,x_{k})$  et $(\ell,i_\ell,x_{\ell})$ peuvent être remplacés par le seul terme $(k,x_k+x_\ell)$
en \gui{contractant la liste}; on réécrira ceci abusivement comme suit: si $i=_Ij$ alors $(i,x_i)\oplus(j,x_j)=(i,x_i+x_j)$;
\item tout terme $(k,0_{M_{i_k}})$ peut être supprimé. 
\end{itemize}
L'addition sur $\bigoplus_{i\in I}M_i$ est définie par concaténation,
et la loi externe est définie par $a\cdot \oplus_{k\in \lrbn}x_{i_k}=\oplus_{k\in \lrbn}ax_{i_k}$.\\
Enfin, le \emph{\Amo librement engendré par $I$} est le module $\bigoplus_{i\in I}\gA$
et il est \hbox{noté $\Ae{(I)}$}.

 La somme directe résout le \pb \uvl correspondant, ce que l'on peut schématiser par le dessin suivant  pour une famille $(\varphi_i)_{i\in I}$
d'\alis des $M_i$ dans un module arbitraire $N$.

\smallskip \centerline{\small
\xymatrix @C=1.5cm @R=.6cm %%@R=1.3cm @C=1.8cm
          {
                &&& M_i 
                \\
N\ar@{<-}[urrr]^{\varphi_i} \ar@{<-}[drrr]_{\varphi_j}\ar@{<-}[ddrrr]_{\varphi_k}
\ar@{<--}[rr]^(.6){\varphi!} && M\ar@{<-}[ur]_{\jmath_i}\ar@{<-}[dr]^{\jmath_j} \ar@{<-}[ddr]_(.5){\jmath_k}&& M=\bigoplus_{i\in I}M_i \\
  &&& M_j 
          \\
&&& M_k 
    \\
}}
 
\smallskip 
Le \Amo   $\Ae {(I)}$  résout le \pb \uvl correspondant, ce que l'on peut schématiser par le dessin suivant pour une famille $x=(x_i)_{i\in I}$
dans un module arbitraire $N$.

\vspace{-.5em}
\smallskip \centerline{\small
\xymatrix @C=1.5cm @R=.6cm %%@R=1.3cm @C=1.8cm
          {
                &&& I \ar[dl]^{\jmath}
                &x(i)=x_i\\
N\ar@{<-}[urrr]^{x} 
\ar@{<--}[rr]^(.6){\psi!} &&\Ae{(I)}   &&\\
}}

\smallskip Notons qu'en conséquence, \emph{si $(x_{i})_{i\in I}$ est un \sgr arbitraire du module $N$, ce dernier est isomorphe à un quotient de $\Ae{(I)}$}.

Soit $I$ un ensemble arbitraire et $(M_i)_{i\in I}$ une famille de \Amosz.
Montrer que le module $\bigoplus_{i\in I}M_i$ est plat \ssi chacun des modules $M_i$
est plat.
En particulier le module libre $\Ae{(I)}$ est plat.
}
\end{exercise}
%--- end -problem-----------------------------------------

}% fin des exos
%:  solutions
\sol

\exer{exoPlatsLecteur}
\emph{Sur un anneau de Bézout intègre $\gZ$, un module $M$ est 
plat \ssi il est sans torsion}.

On sait que la condition est \ncrz. Voyons qu'elle est suffisante.\\
On considère une \syzy $LX$ dans $M$ avec $L=\dex{a_1\;\cdots\;a_n}$
et $X=\tra{\dex{x_1\;\cdots\;x_n}}$.
Si les $a_i$ sont tous nuls, on a $L\,\In=0$ et $\In X= X$, ce qui explique $LX=0$ dans $M$.\\
Sinon,  on écrit
$\sum_i a_i u_i=g$ et $gb_i=a_i$, où $g$ est le pgcd des $a_i$. 
\\
On a $g(\sum_ib_ix_i)=0$, et puisque $M$ est sans torsion $\sum_ib_ix_i=0$.
\\
La matrice $C=\big((u_ib_j)_{i,j\in\lrbn}\big)=UB$ avec $B=\fraC1 g L$, est une \mlp \hbox{pour $(\an)$}. 
On pose $G=\In-C$, on a $CX=0$ \hbox{et $LC=L$}, \hbox{donc $LG=0$}
\hbox{et $GX=X$}, ce qui explique $LX=0$ dans $M$.

\exer{exothPlat1}

\emph{1.} Soit $\mu:P\to M$ une \aliz. On sait (\thref{thPlat1}) que $\mu$ se factorise via un module $L$ libre de type fini: $\mu=\lambda\circ \psi$.

\snic {
\xymatrix {
P\ar[d]_{?} \ar[r]^{\psi}\ar[dr]^{\mu} &L\ar[d]^{\lambda} &
\\
N \ar[r]^{\pi} &M\ar[r] &0 
\\
}}

Comme $L$ est libre, on peut écrire $\lambda=\pi\circ \nu$ pour une \ali
$\nu:L\to N$. Et donc $\mu=\pi\circ \varphi$ pour $\varphi= \nu\circ \psi$.

\emph{2.} Si la \prt est satisfaite avec $N=\gA^{(I)}$, où $I$ est un ensemble discret, on considère une \ali arbitraire $\mu:P\to M$ avec $P$ \pfz. On écrit $\mu=\pi\circ \varphi$ pour une \ali $\varphi:P\to N$. Il existe alors une partie finie $I_0$
de $I$ telle que pour chaque \gtr $g_j$ de $P$, $\varphi(g_j)$ ait ses \coos nulles en dehors de $I_0$. Ceci montre que l'on peut factoriser $\mu$ via le module libre de rang fini $\gA^{(I_0)}$. Donc par le \thref{thPlat1}, $M$ est plat.

%%%%%%%%%%%%%%%%%%%%%%%%%%%%%%%%%%%%%%%%%%%%%%%%%%%%%%%%%%%%%%%%%%%

\exer{exoplatFitt}\\
On considère un \mtf $M$ avec un \sgr $(\xn)$. 
On pose $X=\tra{\lst{x_{1}\;\cdots\; x_{n}}}$. 
Pour $k\in\lrbzn$ et $k+r=n$,
un \gtr typique de~$\cF_k(M)$  s'écrit $\delta=\det(L)$ où
$L\in\MM_r(\gA)$ et $LY=0$, pour un vecteur colonne extrait de $X$: $Y=\tra{\lst{x_{i_1}\;\cdots\;x_{i_r}}}$.
\\
On doit montrer que $\delta\in\cF_k(M)^{2}$. En fait on va montrer que $\delta\in\delta\,\cF_k(M)$.
\\
On suppose \spdg 
que $(i_1,\dots,i_r)=(1,\dots,r)$. On applique la proposition \ref{propPlat2}.
On a donc une matrice $H\in\MM_{r,n}$ avec $HX=Y$ et $LH=0$.\\
Soit \halfsmashbot{$H'=\I_{r,r,n}=\blocs{.9}{.7}{.9}{0}{$\I_{r}$}{$0$}{}{ }$}, et $K=H'-H$.
On a 

\snic{KX=Y-Y=0\;\hbox{  et  }\;LK=LH'=\blocs{.9}{.7}{.9}{0}{$L$}{$0$}{}{ }\,.}

\snii
Soit $K'$ la matrice formée par les $r$ premières colonnes de $K$. Alors $L=LK'$ \hbox{et $\det(L)=\det(L)\det(K')$}. Et puisque $KX=0$, on a $\det(K')\in\cF_k(M).$  

%%%%%%%%%%%%%%%%%%%%%%%%%%%%%%%%%%%%%%%%%%%%%%%%%%%%%%%%%%%%%%%%%%%

\exer{exoClamlsdz}
Supposons l'anneau $\gA$ \lsdzz. 
\\
Soit $\fp$ un \idep et $xy=0$ dans $\gA_\fp$.
Il existe $u\notin\fp$ tel que $uxy=0$ dans $\gA$. Soient $s$ et $t\in\gA$ tels que $s+t=1$, $sux=0$ et $ty=0$ dans $\gA$. Les \elts $s$ et~$t$ ne peuvent être tous deux dans $\fp$ (sinon $1\in\fp$). Si $s\notin\fp$, alors
puisque $sux=0$, on obtient $x=_{\gA_\fp}0$. Si $t\notin\fp$, alors
puisque $ty=0$, on obtient $y=_{\gA_\fp}0$. \\
Ainsi $\gA_\fp$ est un anneau intègre.
\\
Supposons maintenant que tout localisé $\gA_\fp$ en tout \idema $\fp$
soit intègre et supposons que $xy=_\gA0$. Pour un \idema $\fp$ arbitraire on a $x=_{\gA_\fp}0$ \hbox{ou $y=_{\gA_\fp}0$}. 
Dans le premier cas soit $s_\fp\notin\fp$ tel que
$s_\fp x=_\gA 0$. Sinon soit  $t_\fp\notin\fp$ tel que
$t_\fp y=_\gA 0$. la famille des $s_\fp$ ou $t_\fp$ engendre l'idéal 
$\gen{1}$ (car sinon tous les $s_\fp$ ou $t_\fp$ seraient dans un \idemaz).
\\
Il y a donc des $s_i$ en nombre fini vérifiant $s_ix=0$  (dans $\gA$) et des
$t_j$ en nombre fini vérifiant $t_jy=0$, avec une \eqn $\sum_ic_is_i+\sum_jd_jt_j=1$. 
\\
On prend $s=\sum_ic_is_i$, $t=\sum_jd_jt_j$ et l'on obtient
%les \egts 
$sx=ty=0$ et $s+t=1$.

%%%%%%%%%%%%%%%%%%%%%%%%%%%%%%%%%%%%%%%%%%%%%%%%%%%%%%%%%%%%%%%%%%%

\exer{exoClamAri}
On commence par un rappel: d'après le point \emph{3} du \thref{propmlm},
un \id $\gen{a,b}$ d'un anneau $\gA$ est \lop \ssi on peut trouver
$s$, $t$, $u$, $v\in\gA$ tels que $s+t=1$, $sa=ub$ et $tb=va$.
\\
Supposons l'anneau $\gA$ \ariz. 
\\
Soit $\fp$ un \idepz. Pour $a$, $b\in\gA_\fp$ on veut montrer que $a$ divise $b$ ou $b$ divise~$a$ (voir le lemme \ref{lemBezloc}). On peut \spdg prendre $a$ et~$b$ dans $\gA$. Soit alors $s$, $t$, $u$, $v$ comme ci-dessus.
 Les \elts $s$ et~$t$ ne peuvent être tous deux dans $\fp$ (sinon $1\in\fp$). Si $s\notin\fp$, alors
 $a=_{\gA_\fp} s^{-1}ub$ donc $b$ divise $a$ dans~$\gA_\fp$. Si $t\notin\fp$, alors $a$ divise $b$ dans~$\gA_\fp$. 
\\
Supposons maintenant que tout localisé $\gA_\fp$ en tout \idema $\fp$
soit un anneau de Bézout local et soient $a,\,b\in\gA$. 
\\
Pour un \idema $\fp$ arbitraire, on a: $b$ divise $a$ ou $a$ divise $b$ dans~$\gA_\fp$. 
Dans le premier cas soit $s_\fp\notin\fp$ et $u_\fp\in\gA$ tels que
$s_\fp a=_\gA u_\fp b$. Sinon soit  $t_\fp\notin\fp$ et $v_\fp$ tels que
$t_\fp b=_\gA v_\fp a$. la famille des $s_\fp$ ou $t_\fp$ engendre l'idéal 
$\gen{1}$ (car sinon tous les~$s_\fp$ ou $t_\fp$ seraient dans un \idemaz).
\\
Il y a donc des $s_i$, $u_i$ en nombre fini vérifiant $s_ia=u_ib$  (dans $\gA$) et des
$t_j$, $v_j$ en nombre fini vérifiant $t_jb=v_ja$, avec une \eqn $\sum_ic_is_i+\sum_jd_jt_j=1$. 
\\
On prend $s=\sum_ic_is_i$, $u=\sum_ic_iu_i$, $t=\sum_jd_jt_j$, $v=\sum_jd_jv_j$ et l'on obtient
les \egts 
$s+t=1$, $sa=ub$ et $tb=va$.
\\
Pour un \id avec un nombre fini de \gtrsz, on peut faire un raisonnement analogue, ou se reporter au résultat de l'exercice~\ref{exoMLP}.

%%%%%%%%%%%%%%%%%%%%%%%%%%%%%%%%%%%%%%%%%%%%%%%%%%%%%%%%%%%%%%%%%%%

\exer{exoIlops0} \\
L'image de l'\idp $\gen{60}$ de $\ZZ$ par l'\homo 
$\ZZ\to\ZZ/27\ZZ$ est l'\idz~$\gen{3}$ qui ne contient pas d'\elt \ndzz,
et qui n'est donc pas inversible. En fait, comme $\ZZ/27\ZZ$-module, l'\id
$\gen{3}$ n'est même pas \pro (son annulateur~$\gen{9}$ n'est pas \idmz).\\
Lorsque $\rho:\gA\to\gB$ est une \alg plate, l'image d'un \id $\fa\subseteq \gA$ est isomorphe à $\rho\ist(\fa)\simeq\gB\otimes_\gA\fa$. Donc si $\fa$ est \ivz, comme il est \pro de rang $1$, son image est aussi un module \pro de rang $1$.

%%%%%%%%%%%%%%%%%%%%%%%%%%%%%%%%%%%%%%%%%%%%%%%%%%%%%%%%%%%%%%%%%%%

\exer{exoIlops2} \\
On note d'abord qu'un produit d'\ids \lops est toujours \lot principal,
car après des \lons \come convenables, chaque \id devient principal, et  leur produit \egmtz.
\\
On se contente ensuite du cas $\fa=\gen{a,b}$ et de l'exemple
$\gen{a^{4},b^{4}}$. Il sera clair que la technique de calcul se \gns facilement.\\
On part avec $sa=ub$, $tb=va$ et $s+t=1$. Donc  $s^{4}a^{4}=u^{4}b^{4}$ et $t^{4}b^{4}=v^{4}a^{4}$. Puisque $\gen{s^{4},t^{4}}=\gen{1}$ (ce qui s'obtient en écrivant $1=(s+t)^{7}$), on obtient bien
que l'\id $\gen{a^{4},b^{4}}$ est \lopz. 
\\
Montrons ensuite par exemple que  $a^{2}b^{2}\in\gen{a^{4},b^{4}}$.\\
On écrit $s^{2}a^{2}=u^{2}b^{2}$ et $t^{2}b^{2}=v^{2}a^{2}$. Donc $s^{2}a^{2}b^{2}=u^{2}b^{4}$ et $t^{2}a^{2}b^{2}=v^{2}a^{4}$.\\
 Enfin $1=(s+t)^{3}=s^{2}(s+3t)+t^{2}(t+3s)$. Donc finalement

\snic{a^{2}b^{2}=(t+3s)v^{2}a^{4}+(s+3t)u^{2}b^{4}.}

%%%%%%%%%%%%%%%%%%%%%%%%%%%%%%%%%%%%%%%%%%%%%%%%%%%%%%%%%%%%%%%%%%%

\exer{exoMLP} On n'a pas besoin de supposer que l'anneau est \ariz.\\
On va montrer que si dans un anneau $\gA$ chaque couple $(a_i,a_j)$ admet une \mlpz,
il en va de même pour le $n$-uplet $(\an)$.
\\
Notons que ceci est à rapprocher de la \dem par Dedekind du \thref{th1IdZalpha}, dans laquelle il n'est question que d'\ids \ivsz, car sur un anneau intègre les \ids \ivs sont exactement les \ids \lops non nuls.
\\
Notons aussi que le résultat est a priori clair: par \lons \come successives, on obtiendra un \idp au bout de chaque branche d'un arbre de calcul a priori très grand. Ceci montrera que l'\id $\gen{\an}$ est 
toujours engendré par l'un des $a_i$ après des \lons en des \ecoz.
Mais ce que nous avons en vue ici, c'est plutôt un calcul pratique de la \mlpz.
\\
On procède par \recu sur $n$. 

Montrons l'étape de \recu pour le passage de $n=3$ à $n+1=4$.\label{corexoPruf3}\\
On considère $a_1$, $a_2$, $a_3$, $a_4\in\gZ$. \\
Par \hdr on a une matrice
\halfsmashtop{$C=\cmatrix{x_1&x_2&x_3\cr y_1&y_2&y_3\cr z_1&z_2&z_3\cr }$} qui convient
pour $(a_1,a_2,a_3)$, et des matrices $\cmatrix{c_{11}&c_{14}\cr d_{11}&d_{14}}$, $\cmatrix{c_{22}&c_{24}\cr d_{22}&d_{24}}$, $\cmatrix{c_{33}&c_{34}\cr d_{33}&d_{34}}$ 
qui conviennent respectivement pour  $(a_1,a_4)$,  $(a_2,a_4)$ et  
$(a_3,a_4)$.
Alors on va vérifier que la \und{trans}p\und{osée} de la matrice suivante convient
pour $(a_1,a_2,a_3,a_4)$:
$$
\cmatrix{c_{11}x_1&c_{22}y_1&c_{33}z_1&d_{11}x_1+d_{22}y_1+d_{33}z_1 
     \cr c_{11}x_2&c_{22}y_2&c_{33}z_2&d_{11}x_2+d_{22}y_2+d_{33}z_2  
     \cr c_{11}x_3&c_{22}y_3&c_{33}z_3& d_{11}x_3+d_{22}y_3+d_{33}z_3 
     \cr c_{14}x_1
        &c_{24}y_2 
        & c_{34}z_3
        & d_{14}x_1+d_{24}y_2+d_{34}z_3}
$$
Tout d'abord, on doit vérifier que la trace de la matrice est égale à $1$, i.e.
 
\snic{t=c_{11}x_1+c_{22}y_2+c_{33}z_3+d_{14}x_1+d_{24}y_2+d_{34}z_3=1,}

\snii or $c_{11}+d_{14}=1=c_{22}+d_{24}=c_{33}+d_{34}$ donc $t=x_1+y_2+z_3=1$.
\\
Et l'on doit vérifier que chacune des lignes de la 
matrice transposée est propotionnelle à $\lst{a_1\;a_2\;a_3\;a_4}$.
Deux cas se présentent. Tout d'abord, l'une des trois premières,
par exemple 
la ligne $\lst{c_{11}x_1\;c_{11}x_2\;c_{11}x_3\;c_{14}x_1}$. Il faut vérifier les deux types d'\egts suivantes

\snic{a_1c_{11}x_2=a_2c_{11}x_1, \quad \hbox{et}\quad a_1c_{14}x_1=a_4c_{11}x_1.}

Pour la première \egt on utilise $a_2x_1=a_1x_2$ et pour la seconde
$a_1c_{14}=a_4c_{11}$.\\
Enfin on doit vérifier que $[\,a_1\;a_2\;a_3\;a_4\,]$ est 
proportionnelle à la transposée de

\snac{\cmatrix{d_{11}x_1+d_{22}y_1+d_{33}z_1\\
d_{11}x_2+d_{22}y_2+d_{33}z_2\\
d_{11}x_3+d_{22}y_3+d_{33}z_3\\
d_{14}x_1+d_{24}y_2+d_{34}z_3}.}

\snii Ceci résulte d'une part de la proportionnalité de $[\,a_1\;a_2\;a_3\,]$
à chacune des lignes~\hbox{$[\,x_i\;y_i\;z_i\,]$}, et d'autre part de
la  proportionnalité des lignes $[\,a_i\;a_4\,]$ aux lignes $[\,d_{i1}\;d_{i4}\,]$.

Notons pour terminer que le passage de $n-1$ à $n$ (pour n'importe quel $n>2$) est tout à fait analogue.
%%%%%%%%%%%%%%%%%%%%%%%%%%%%%%%%%%%%%%%%%%%%%%%%%%%%%%%%%%%%%%%%%%%

\exer{exoIlops} On écrit $\fa=\gen{\am}$, $\fb=\gen{\bn}$.
On peut supposer \hbox{que  $f=\sum_{k=1}^{m}a_kX^{k-1}$} et $g=\sum_{h=1}^{n}b_hX^{h-1}$.

\emph{1.} Soit $F$ une \mlp pour $(\am)$. Si $\rc(f)=\fa$, on a des \eco $s_i$ (la diagonale de $F$) et des \pols $f_i\in\AX$ (donnés par les lignes de $F$) qui satisfont les \egts $s_if=a_if_i$ dans $\AX$. En outre,  le \coe \hbox{de $X^{i-1}$} dans $f_i$ est égal à $s_i$, donc $\rc(f_i)\supseteq \gen{s_i}.$
\\ 
En posant $\gA_i=\gA\big[\fraC1{s_i}\big]$, on a $\rc(f_i)=_{\gA_i}\! \gen{1}$ et les \egts

\snic{s_i\rc(fg)=\rc(a_if_ig)=a_i\rc(f_ig)=_{\gA_i}a_i\rc(g)=_{\gA_i}\rc(a_if_i)\rc(g)=s_i\rc(f)\rc(g) }

\snii (la troisième \egt vient du  corolaire~\ref{corlemdArtin}~\emph{4}
car $\rc(f_i)=_{\gA_i}\! \gen{1}$).\\
D'où l'\egt $\rc(fg)=\rc(f)\rc(g)=\fa\fb$ car elle est vraie dans chaque $\gA_i$.
 
\emph{2.} Si $g$ est aussi \lop on obtient de la même manière
$t_jb=b_jg_j$ \hbox{dans $\AX$}, avec $\rc(g_j)\supseteq \gen{t_j}$ et les $t_j$ \com dans $\gA$.
On a donc 

\snic{s_it_j\rc(fg)=_{\gA_{ij}}a_ib_j\rc(f_ig_j)=_{\gA_{ij}}\gen{a_ib_j}.}

Ceci nous dit que l'\id $\rc(fg)=\fa\,\fb$ devient principal après $mn$ \lons \comez. Comme cet \id admet $m+n-1$ \gtrs (les \coes de~$fg$)
il y a une \mlp pour ces \gtrsz.\\ 
Pour la calculer,
on peut utiliser la \dem de l'implication \emph{1} $\Rightarrow$  \emph{3} dans le \thref{propmlm}. 
Cette \dem est assez simple, de même que le calcul qu'elle sous-tend. Mais si l'on examine en détail ce qui va se passer, on s'aperçoit que dans la \dem ci-dessus on a utilisé le lemme de Gauss-Joyal:   sur l'anneau $\gA_{ij}$, on a $1\in\rc(f_i)\rc(g_j)$ car $1\in\rc(f_i)$ et $1\in\rc(g_j)$. Ce lemme admet plusieurs \dems \elrs (voir \ref{lemGaussJoyal} et \ref{corlemdArtin}), mais aucune ne donne une formule simple qui permette de fournir la \coli des \coes de $fg$ égale à $1$, à partir des deux \colis des \coes de~$f$ et de ceux de $g$.
\\
Merci \alec qui nous indiquera un calcul direct court,
par exemple dans le cas où l'anneau est intègre à \dve explicite\footnote{Notons que 
dans le cas d'un anneau intègre à \dve explicite, une \mlp est connue
à partir de ses seuls \elts diagonaux, ce qui peut simplifier les calculs.}. 

%%%%%%%%%%%%%%%%%%%%%%%%%%%%%%%%%%%%%%%%%%%%%%%%%%%%%%%%%%%%%%%%%%%
\exer{exogcdlcm} On écrit $\fa=\gen{\an}$, $\fb=\gen{\bbm}$.
\\
On utilisera le résultat de l'exercice \ref{exoMLP} qui montre que si
tout \id à deux \gtrs est \lopz, alors tout \itf est \lopz. 

\emph{1.} Dans l'exercice \ref{exoPgcdPpcm} point \emph{4} on a montré que  $1 \in (\fa : \fb) + (\fb :\fa)$, $\fa \cap \fb$ est \tf et  $\fa\fb = (\fa
\cap \fb)(\fa + \fb)$.\\
Si $\fa+\fb$ est \lopz, il y a un \sys d'\eco tel qu'en inversant l'un quelconque d'entre eux, l'\id est engendré par un $a_k$ ou un $b_\ell$. Mais si  $\fa+\fb=\gen{a_k}\subseteq \fa$, on a $\fb\subseteq \fa$, donc $\fa\cap\fb=\fb$, \lop par hypothèse.
Ainsi $\fa\cap\fb$ est \lop car il est \lop après \lon en des \ecoz.

\emph{2.} Si l'anneau est intègre et si  $(\fa+\fb)(\fa\cap\fb)=\fa\fb$
pour $\fa=\gen{a}$ et $\fb=\gen{b}$ (\hbox{où $a,b\neq 0$}), on obtient que $\gen{a,b}(\fa\cap\fb)=\gen{ab}$, donc $\gen{a,b}$ est \iv (et aussi~\hbox{$\gen{a}\cap\gen{b}$}
par la même occasion). Lorsque c'est vérifié pour tous $a,b\neq 0$, l'anneau est \ariz. 

\emph{3.} La seule implication délicate consiste à montrer que
si $\gA$ est \lsdz et si $(\fa+\fb)(\fa\cap\fb)=\fa\fb$ lorsque
$\fa=\gen{a}$ et $\fb=\gen{b}$ alors l'anneau est \ariz, autrement dit
tout \id $\gen{a,b}$ est \lopz. 

Si $\gen{a,b}(\fa\cap\fb)=\gen{ab}$, on écrit $ab=au+bv$ avec $u$ et $v\in\fa\cap\fb$:

\snic{u=ax=by,\;v=az=bt,\;\;\hbox{d'où} \;\;au+bv=ab(y+z)=ab.}

Puisque l'anneau est \lsdzz, de l'\egtz~\hbox{$ab(y+z-1)=0$}, on déduit trois \lons \come dans lesquelles on obtient respectivement~\hbox{$a=0$}, $b=0$ et $1=y+z$. Dans les deux premiers cas $\gen{a,b}$
est principal. Dans le dernier cas $\gen{a,b}$ est \lop (localiser en $y$ ou en~$z$).

%%%%%%%%%%%%%%%%%%%%%%%%%%%%%%%%%%%%%%%%%%%%%%%%%%%%%%%%%%%%%%%%%%%

\exer{exoSECSci} On écrit $\fa=\gen{\an}$, $\fb=\gen{\bbm}$.

\emph{1.} Démontré dans le  point \emph{4} de l'exercice \ref{exoPgcdPpcm}.
%Si $\fa+\fb$ est \lopz, alors $(\fa:\fb)+(\fb:\fa)=\gen{1}$.
%On doit résoudre le \sli 
%
%\snic{s+t=1,\quad  sa_i=\sum_jx_{ij}b_j \,(i\in\lrbn), \quad tb_j=\sum_iy_{ji}a_i
%\,(j\in\lrbm),}
%
%\snii
%avec les inconnues $s$, $t$, $x_{ij}$, $y_{ji}$.
%On considère des $u_i$ et $v_j$ \com tels que 
%
%\snic{u_i(\fa+\fb)=\gen{u_ia_i}$
%et $v_j(\fa+\fb)=\gen{v_jb_j}.}
%
%\snii
%Si on inverse $u_i$ ou $v_j$, une solution
%du \sys est évidente, avec $t=1$, $s=0$ dans le premier cas,  $s=1$, $t=0$
%dans le second cas.
%On conclut avec le \plg de base.

\emph{2.} On suppose maintenant $(\fa:\fb)+(\fb:\fa)=\gen{1}$, i.e.,
on a $s$, $t\in\gA$ avec 

\snic{s+t=1$, $s\fa\subseteq \fb$, $t\fb\subseteq \fa.}

%ce qui revient à donner une solution du \sli décrit au point \emph{1.}

\emph{2a.}~
$(\fa+\fb):\fc = (\fa:\fc)+(\fb:\fc)$. Dans cette \egt comme dans les suivantes
(jusqu'à~\emph{2f}), une inclusion n'est pas évidente (ici c'est $\subseteq $). Prouver l'inclusion non évidente revient à résoudre un \sli (ici, étant donné un $x$ tel \hbox{que $x\fc\subseteq \fa+\fb$}, on cherche $y$ et $z$ tels que $x=y+z$, $y\fc\subseteq \fa$ et $z\fc\subseteq \fb$). \\
On peut donc utiliser le \plg de base avec les \ecoz~$s$ et~$t$.\\
Lorsqu'on inverse $s$, on obtient $\fa\subseteq \fb$, et lorsqu'on inverse $t$, on obtient $\fb\subseteq \fa$. Dans les deux cas l'inclusion souhaitée devient triviale.

Pour mémoire: 
\emph{2b.}\, $\,\fc :(\fa~\cap~ \fb)=(\fc:\fa)+(\fc:\fb)$.
\,
\emph{2c.}\, $(\fa+\fb)(\fa~\cap~ \fb) = \fa~\fb$.
\\
\emph{2d.}\, $\fc\,(\fa~\cap~ \fb)= \fc~ \fa~ \cap~ \fc~ \fb$.
\,
\emph{2e.}\, $\,\fc+(\fa~\cap~ \fb)=(\fc+\fa)\cap (\fc+\fb)$.
\\
\emph{2f.}\, $\,\fc~\cap\,(\fa+ \fb)=(\fc~\cap~ \fa)+(\fc~\cap~ \fb)$.

\emph{2g.} La suite exacte courte ci-après (où $\delta(x)=(x, -x)$ et $\sigma(y, z)=y+z$) est scindée:
$$\preskip-.4em \postskip.4em 
0 \longrightarrow \fa~\cap~ \fb \vers{\delta} \fa \times \fb 
\vers{\sigma} \fa+\fb
\longrightarrow 0. 
$$
On veut définir $\tau:\fa+\fb\to\fa\times \fb$ telle que $\sigma\circ \tau=\Id_{\fa+\fb}$.\\
Si $\fa\subseteq \fb$, on peut prendre $\tau(b)=(0,b)$ pour tout $b\in\fb=\fa+\fb$.
Si $\fb\subseteq \fa$, on peut prendre $\tau(a)=(a,0)$ pour tout $a\in\fa=\fa+\fb$. \\
Dans le premier cas cela implique $s\tau(a_i)=\big(0,\sum_jx_{ij}b_j\big)$ et $s\tau(b_j)=(0,sb_j)$.\\
Dans le second cas cela implique $t\tau(b_j)=\big(\sum_iy_{ji}a_i,0\big)$ et $t\tau(a_i)=(ta_i,0)$.\\
On essaie donc de définir $\tau$ par la formule suivante qui coïncide
avec les deux précédentes dans les deux cas particuliers.

\snic{\tau(a_i)=\big(ta_i,\sum_jx_{ij}b_j\big),\;\;\tau(b_j)=\big(\sum_iy_{ji}a_i,sb_j\big) .}

Pour que cette tentative réussisse, il faut et suffit que lorsque
$\sum_i\alpha_ia_i=\sum_{j}\beta_jb_j$, on ait l'\egt
$$
\preskip-.2em \postskip.4em \ndsp
\sum_i\alpha_i\big(ta_i,\sum_jx_{ij}b_j\big)=\sum_{j}\beta_j\big(\sum_iy_{ji}a_i,sb_j\big) . 
$$
Ceci résulte pour la première \coo du calcul suivant
(même chose pour la deuxième \cooz).

\snic{\sum_i\alpha_ita_i=t\sum_i\alpha_ia_i=t\sum_{j}\beta_jb_j=\sum_{j}\beta_jtb_j=\sum_{j}\beta_j\sum_iy_{ji}a_i.
}

Enfin l'\egt $\sigma\circ \tau=\Id_{\fa+\fb}$ est satisfaite parce qu'elle
l'est en restriction à $\fa$ et $\fb$ (calcul \imdz).

%%%%%%%%%%%%%%%%%%%%%%%%%%%%%%%%%%%%%%%%%%%%%%%%%%%%%%%%%%%%%%%%%%%

\exer{exoGaussien}  
 \emph{1.} Démontré dans l'exercice \ref{exoIlops}.

\emph{2.} Soient $a$, $b$, $c$, $d\in\gA$. On pose $\fa=\gen{a,b}$\\
On considère $f=aX+b$ et $g=aX-b$, on obtient $\gen{a,b}^{2}=\geN{a^{2},b^{2}}$, i.e., 
$$
\preskip.2em \postskip.4em 
ab=ua^{2}+vb^{2}. 
$$
En considérant $f=cX+d$ et $g=dX+c$, on obtient $\gen{c,d}^{2}=\geN{c^{2}+d^{2},cd}$. Autrement dit $c^{2}$ et $d^{2}\in\geN{c^{2}+d^{2},cd}$.  
\\
On pose $\fb=\gen{ua,vb}$. On a $ab\in\fa\fb$. Il suffit de montrer que $\fa^{2}\fb^{2}=\geN{a^{2}b^{2}}$ car cela implique $\fa$ \iv (on traite le cas $a$, $b\in\Atl$).
Or on a

\snic{a^{2}b^{2}\in\fa^{2}\fb^{2}=\geN{a^{2},b^{2}}\geN{u^{2}a^{2},v^{2}b^{2}}.}

Il nous faut donc montrer que $u^{2}a^{4}$ et $v^{2}b^{4}\in\geN{a^{2}b^{2}}.$
Posons $u_1=ua^{2}$ et  $v_1=vb^{2}$. On a $u_1+v_1=ab$ et $u_1v_1\in\geN{a^{2}b^{2}}$. Donc aussi $u_1^{2}+v_1^{2}\in\geN{a^{2}b^{2}}$.\\
Et puisque $u_1^{2}\in\geN{u_1^{2}+v_1^{2},u_1v_1}$, on obtient bien
$u_1^{2}\in\geN{a^{2}b^{2}}$ (même chose \hbox{pour $v_1^{2}$}).

\emph{3.} Les \egts du point \emph{2} sont toutes satisfaites.
\\
Montrons d'abord que l'anneau est \lsdzz. 
\\
On suppose $cd=0$, puisque
$c^{2}\in\geN{c^{2}+d^{2},cd}$, on a $c^{2}=x(c^{2}+d^{2})$, i.e.

\snic{xd^{2}=(1-x)c^{2}.}

\snii
On en déduit que $xd^{3}=0$, et comme $\gA$ est réduit, $xd=0$. De même $(1-x)c=0$.

Voyons maintenant que l'anneau est \ariz. On part de $a$, $b$ arbitraires
et on veut montrer que $\gen{a,b}$ est \lopz.
D'après le point \emph{2} on a un \id $\fc$ tel que $\gen{a,b}\fc=\geN{a^{2}b^{2}}$. On a donc $x$ et $y$ avec

\snic{\gen{a,b}\gen{x,y}= \geN{a^{2}b^{2}} \;\hbox{ et } ax+by=a^{2}b^{2}.}

On écrit $ax=a^{2}b^{2}v$ et $by=a^{2}b^{2}u$.
De l'\egt $a(ab^{2}v-x)=0$, on déduit deux \lons \comez,
dans la première $a=0$, dans la seconde $x=ab^{2}v$.
On suppose donc \spdg que $x=ab^{2}v$ et, symétriquement $y=ba^{2}u$.
Ceci donne 
$$
\preskip-.4em \postskip.4em 
\gen{a,b}\gen{x,y}=ab\gen{a,b}\gen{au,bv}=\geN{a^{2}b^{2}}. 
$$
On peut encore supposer \spdg que $\gen{a,b}\gen{au,bv}=\gen{ab}$.
\\
On a aussi $ax+by=a^{2}b^{2}(u+v)$
et puisque $ax+by=a^{2}b^{2}$, on suppose \spdg que $u+v=1$. 
\\
Puisque $a^{2}u=abu'$, on suppose \spdg que $au=bu'$. 
\\
Symétriquement
$bv=av'$, et puisque $u+v=1$, $\gen{a,b}$ est \lopz.

%%%%%%%%%%%%%%%%%%%%%%%%%%%%%%%%%%%%%%%%%%%%%%%%%%%%%%%%%%%%%%%%%%%
\exer{exoVascon}
\emph{1.} Soit $a\in\gA$ et $a_1$, \dots,  $a_n\in \gA$ qui engendrent $\fa=\Ann(a)$.
 Si l'un des $a_i$ est dans $\Ati$, on obtient $a=0$ et $\fa=\gen{1}$.
Il reste à traiter le cas où tous les~$a_i$ sont dans $\fm$.
 Soit $b$ l'un des $a_i$. Puisque $\fm$ est plat et $b\in\fm$, la relation $ab=0$ nous donne des \eltsz~\hbox{$c_1$, \dots, $c_m\in\fa$} et $b_1$, \dots, $b_m\in\fm$ \hbox{avec
  $b=\sum_{i\in\lrbm}c_ib_i$}. Donc $b\in\fa\fm$, ce qui donne $b=\sum_{i\in\lrbn}a_iz_i$
pour des $z_i\in\fm$.
D'où une \egt matricielle
$$
[\,a_1\;\cdots\;a_n\,]=M\,[\,a_1\;\cdots\;a_n\,]\quad \hbox{ avec } M\in\Mn(\fm).
$$
Ainsi $[\,a_1\;\cdots\;a_n\,] (\In-M)=[\,0\;\cdots\;0\,]$ avec $\In-M$ \ivz, donc $\fa=0$.

\emph{2.}  On considère $a$, $b\in\gA$, on doit démontrer que l'un divise l'autre. 
Si l'un des deux est \ivz, l'affaire est entendue.
Il reste à examiner le cas où $a$ et $b\in\fm$. 
\\ 
On considère
une matrice 
$$
P=\cmatrix{a_1&\cdots&a_n\cr b_1&\cdots&b_n}
$$ 
dont les colonnes 
engendrent le module $K$ noyau de $(x,y)\mapsto bx-ay$. En particulier on a $a_ib=b_ia$ pour 
chaque $i$. Si l'un des $a_i$ ou $b_i$ est \ivz, l'affaire est entendue.
Il reste à examiner le cas où les $a_i$ et $b_i$ sont dans $\fm$.
\\
 Soit $(c,d)$ l'un des $(a_i,b_i)$. Puisque $\fm$ est plat et $a$, $b\in\fm$, la relation $cb-da=0$ donne 
$$\preskip-.4em \postskip.0em
\cmatrix{c\cr d}=\cmatrix{c_1&\cdots&c_m\cr d_1&\cdots&d_m} \cmatrix{y_1\cr\vdots\cr y_m} \quad \hbox {avec les }y_i\in\fm\; \hbox{ et les }\cmatrix{c_j\cr d_j}\in K.
$$
En exprimant les \smashbot{$\cmatrix{c_j\cr d_j}$} comme \colis des colonnes de $P$
on obtient
$$\preskip-.2em \postskip.3em
\cmatrix{c\cr d}=P \cmatrix{z_1\cr\vdots\cr z_n} \quad \hbox {avec les }z_i\in\fm.
$$
D'où ensuite $P=PN$ avec une matrice $N\in\Mn(\fm)$, puis $P=0$. Ceci implique
que $(a,b)=(0,0)$, et $a$ divise $b$ (en fait, dans ce cas, $\gA$ est trivial).
\perso{on peut se demander si l'hypothèse  \gui{\dcdz} est vraiment \ncrz.}
  
%%%%%%%%%%%%%%%%%%%%%%%%%%%%%%%%%%%%%%%%%%%%%%%%%%%%%%%%%%%%%%%%%%%
\exer{exoQuotientPlat}
Soient $a_i \in \gA$, $x_i \in M$ vérifiant $\sum_{i=1}^n a_i x_i \equiv
0 \bmod K$, relation que l'on doit expliquer.  On pose $\fa = \gen {\ua}$ de
sorte que $\fa K = \sum_i a_iK$; on peut écrire, puisque $\sum_i
a_ix_i \in \fa M \cap K = \fa K$, une \egt $\sum_i a_ix_i =
\sum a_iy_i$ avec \hbox{les $y_i \in K$}. On a donc, avec $z_i = x_i-y_i$,
la relation $\sum_i a_iz_i = 0$ dans $M$. Puisque $M$ est plat, cette
relation produit un certain nombre de vecteurs de $M$,
disons 3 pour simplifier, notés $u$, $v$, $w$ et 3 suites
de scalaires $\underline\alpha =(\alpha_1, \ldots,\alpha_n)$,
$\underline\beta =(\beta_1, \ldots,\beta_n)$ et
$\underline\gamma =(\gamma_1, \ldots,\gamma_n)$, le
tout vérifiant:

\snic {
(z_1, \ldots, z_n) = (\alpha_1, \ldots,\alpha_n)\,u +
 (\beta_1, \ldots,\beta_n)\,v + (\gamma_1, \ldots,\gamma_n)\,w
}

\snii
et $\scp {\ua}{\underline\alpha} = \scp {\ua}{\underline\beta} = 
\scp {\ua}{\underline\gamma} = 0$.

\snii
Puisque $z_i \equiv x_i \bmod K$, on obtient notre explication
convoitée dans $M/K$:

\snic {
(x_1, \ldots, x_n) \equiv (\alpha_1, \ldots,\alpha_n)\,u +
 (\beta_1, \ldots,\beta_n)\,v + (\gamma_1, \ldots,\gamma_n)\,w
\;\mod K.
}

%%%%%%%%%%%%%%%%%%%%%%%%%%%%%%%%%%%%%%%%%%%%%%%%%%%%%%%%%%%%%%%%%%%

\exer{propfreeplat} \emph{(Tout module est quotient d'un module plat)}\\
On suppose \spdg que $I$ est finiment énuméré. Autrement dit $I=\so{i_1,\dots,i_n}$.
On note $M=\bigoplus_{i\in I}M_i$.

\emph{Supposons d'abord que les modules $M_i$ sont plats}, et considérons une \syzy
dans $M$

\snic{0=\sum_{\ell\in\lrbm}a_\ell x_\ell=\sum_{\ell\in\lrbm}a_\ell(\oplus_{k\in\lrbn} x_{k,\ell})=\oplus_{k\in\lrbn}y_k}

avec $y_k=\sum_{\ell\in\lrbm}a_\ell x_{k,\ell}$ et les $x_{k,\ell}\in M_{i_k}$.\\
Par \dfn de l'\egt dans $M$, puisque $\oplus_{k\in\lrbn}y_k=0$, on est dans
(au moins) l'un des deux cas suivants:
\begin{itemize}
\item tous les $y_k$ sont nuls,
\item deux indices sont égaux dans $I$: $i_k=_Ii_h$ pour $h$ et $k$ distincts \hbox{dans $\lrbn$}. 
\end{itemize} 
Le premier cas se traite comme celui d'une somme directe sur un  $I$ fini.
Le deuxième cas se ramène au premier par \recu sur $n$.

\emph{Supposons maintenant que $M$ est plat} et considérons par exemple 
l'indice~\hbox{$i_1\in I$} et une \syzy $\sum_{\ell\in\lrbm}a_\ell x_\ell=0$ dans $M_{i_1}$.
On explique cette \syzy dans $M$ en écrivant 

\snic{x_\ell=_M\sum_{j\in \lrbp} g_{\ell,j}z_{j}$ avec $\sum_{\ell\in\lrbm} a_\ell g_{\ell,j}=_\gA0$ pour chaque $j.}

On réécrit $z_{j}=\oplus_{k\in\lrbn} y_{j,k}$ avec $y_{j,k}\in M_{i_k}$, ce qui donne

\snic{x_\ell=_M\oplus_{k\in\lrbp}\sum_{j\in \lrbp}g_{\ell,j}y_{j,k}=_M\oplus_{k\in\lrbn}y_{\ell,k}.}

Par \dfn de l'\egt dans $M$,  on est dans
(au moins) l'un des deux cas suivants:
\begin{itemize}
\item pour chaque $\ell$, on a  $x_\ell=\sum_{j\in \lrbn} g_{\ell,j}y_{\ell,1}$ dans $M_{i_1}$,
\item on a  dans $I$: $i_1=_Ii_h$ pour un $h\neq 1$ \hbox{dans $\lrbn$}. 
\end{itemize} 
Dans le premier cas on a dans $M_{i_1}$ les \egts qui nous conviennent.
\\
Le deuxième cas se ramène au premier par \recu sur $n$. 
%

%%%%%%%%%%%%%%%%%%%%%%%%%%%%%%%%%%%%%%%%%%%%%%%%%%%%%%%%%%%%%%%%%%%

%%%%%%%%%%%%%%%%%%%%%%%%%%%%%%%%%%%%%%%%%%%%%%%%%%%%%%%%%%%%%%%%%%%
%}% fin des solutions d'exos 
%:   ---- Section*{references}-----------    
\Biblio

Les \adps intègres ont été introduits par H. Prüfer en 1932 dans \cite{Prufer}.
Leur place centrale en théorie multiplicative des idéaux est mise en valeur
dans le livre de référence sur le sujet \cite{Gil}. 
Bien qu'ils aient été introduits de manière très concrète comme les anneaux intègres dans lequel tout \itf non nul est inversible, 
cette \dfn a fait souvent place dans 
la littérature moderne à la suivante, purement abstraite, qui ne fonctionne qu'en présence de principes non \cofs (tiers exclu et axiome du choix): 
la \lon en n'importe quel \idep donne un \advz. 

Les \anars sont introduits par L. Fuchs en 1949 dans~\cite{Fuchs}. 
 
Dans le cas d'un anneau non intègre, la \dfn 
que nous avons adoptée pour les \adps est due à Hermida et S\'anchez-Giralda
\cite{HS}. 
C'est celle qui nous a paru la plus naturelle, vue l'importance centrale 
du concept de platitude en \alg commutative. Un autre nom pour ces 
anneaux, dans la littérature est \emph{anneau de dimension globale faible 
inférieure ou égale à un}, ce qui est plutôt inélégant. 
Par ailleurs, on trouve souvent dans la littérature un anneau de Prüfer 
défini comme un anneau
dans lequel tout \id contenant un \elt \ndz est inversible. 
Ce sont donc presque des \anarsz, mais le comportement des \ids ne contenant
pas d'\elt \ndz semble tout à fait aléatoire (cf. exercice~\ref{exononanar}).

Un exposé assez complet sur les \anars et les \adps écrit dans le style
des \coma se trouve dans \cite[Ducos\&al.]{dlqs} et~\cite[Lombardi]{lom99}.

Un survey très complet sur les variations de la notion d'\adp intègre quand on
supprime l'hypothèse d'intégrité est donné dans~\cite[Bazzoni\&Glaz]{BG}, y compris les anneaux gaussiens (exercice~\ref{exoGaussien}). 

\newpage \thispagestyle{CMcadreseul}
\incrementeexosetprob

%:        %%%%%%%%%%%%%%%%%%%%%%%%%%%%%%%%%%%%
%:        %%%%%%%%%%%%%%%%%%%%%%%%%%%%%%%%%%%%

%---- Chapitre  {Anneaux locaux}------------
\chapter{Anneaux locaux, ou presque}
\label{chap Anneaux locaux}\relax
\minitoc
%\subsection*{Introduction}
%\addcontentsline{toc}{section}{Introduction}
%-----------------------------------------

%%%%%%%%%%%%%%%%%%%%%%%%%%%%%%%%%%%%%%%%%%%%%%%%%%%%%%%%%%%%%%%%%%%%%%%% 
%%%%%%%%%%%%%%%%%%%%%%%%%%%%%%%%%%%%%%%%%%%%%%%%%%%%%%%%%%%%%%%%%%%%%%%% 
%%%%%%%%%%%%%%%                                            %%%%%%%%%%%%% 
%%%%%%%%%%%%%%%   section Quelques \dfns constructives     %%%%%%%%%%%%% 
%%%%%%%%%%%%%%%                                            %%%%%%%%%%%%% 
%%%%%%%%%%%%%%%%%%%%%%%%%%%%%%%%%%%%%%%%%%%%%%%%%%%%%%%%%%%%%%%%%%%%%%%%
%%%%%%%%%%%%%%%%%%%%%%%%%%%%%%%%%%%%%%%%%%%%%%%%%%%%%%%%%%%%%%%%%%%%%%%%

%--- Sec{Quelques dfns constr} --- secAloc1
\section{Quelques \dfns constructives}
\label{secAloc1}
%-------------------------------

En \clama un anneau local est souvent défini com\-me un
anneau possédant un seul \idemaz.
Autrement dit les \elts non inversibles forment un idéal.
Cette deuxième \dfn a l'avantage d'être plus simple (pas de quantification
sur l'ensemble des idéaux).
Elle se prête cependant relativement mal à un
traitement algorithmique à cause de la négation contenue dans \gui{\elts non
inversibles}. C'est la raison pour laquelle on adopte en \coma la \dfn
donnée \paref{eqAloc}: si la somme de deux \elts est \ivz, l'un des deux est \ivz.

Nous nous trouvons maintenant
 dans l'obligation d'infliger \alec classique quelques
\dfns peu usuelles pour \luiz, dans la même lignée que la \dfn d'\aloz.
Qu'\il se rassure, sur d'autres planètes, dans d'autres systèmes solaires,
c'est sans doute la situation symétrique qui se produit.
Les \maths y ont toujours été \covs et l'on vient à peine d'y découvrir
l'intérêt du point de vue cantorien abstrait. \stMF
Un\e auteur\e dans le nouveau style est en train d'expliquer que pour \lui il est
beaucoup plus simple de voir un anneau local comme un anneau possédant un seul
\idemaz. \Llec fera-t-\il l'effort de \la suivre?

%:  subsec{Radical de Jacobson, anneaux locaux, corps}
\subsec{Radical de Jacobson, anneaux locaux, corps}
%-----------------------------------------

Rappelons que pour un anneau $\gA$ nous notons $\Ati$ le groupe multiplicatif
des \elts inversibles, encore appelé groupe des unités.

Un \elt $x$ d'un anneau $\gA$ est dit \ix{noninversible} (en un seul mot) s'il
vérifie{\footnote{Nous utilisons ici
une version légèrement affaiblie de la négation. 
Pour une \prtz~$\sfP$ portant sur des \elts de l'anneau $\gA$
ou d'un \Amo $M$, 
nous considérons la \prt $\sfP':=(\sfP\Rightarrow 1=_\gA0)$.
C'est la négation de $\sfP$ lorsque l'anneau n'est pas trivial.
Il arrive pourtant souvent qu'un anneau construit dans une preuve puisse
être trivial sans qu'on le sache.
Pour faire un traitement entièrement
\cof de la preuve classique usuelle dans une telle situation (la preuve classique
exclut le cas de l'anneau trivial par un argument ad hoc) notre version affaiblie  de la négation s'avère alors en \gnl utile.
Un \cdi (en un seul mot), ne vérifie pas \ncrt l'axiome
des ensembles discrets $\forall x,y\;\big(x=y\;\hbox{ou}\;\lnot(x=y)\big)$
mais il vérifie sa version affaiblie:

\snic{\forall\, x,\,y,\;\big(x=y\;\hbox{ ou }\;(x=y)'\big),
}

\vspace{2pt}
\noindent  
puisque si $0$ est inversible, alors $1=0$.
\label{footnoteNegation}}} l'implication suivante

\snic{x\in\Ati\; \Rightarrow \; 1=_\gA0.
}

\smallskip
Dans l'anneau trivial  l'\elt $0$ est à la fois \iv et non\ivz.

Pour un anneau commutatif arbitraire, l'ensemble des \elts  $a$ de $\gA$  qui
vérifient
\rdb
%---  equation eqDefRadJac --------
\begin{equation}\label{eqDefRadJac}
\forall x\in \gA~~~ 1+ax\in\Ati
\end{equation}
%---------------------end equation--------------
est appelé le \ixx{radical}{de Jacobson} de $\gA$. Il sera noté $\Rad(\gA)$.
C'est un \id parce que si $a$, $b\in\Rad\gA$, on peut écrire, pour $x\in\gA$:
%:2012 remplace x, y, a  par  a, b, x pour correspondre avec eqDefRadJac
\index{radical de Jacobson!d'un anneau}

\snic{1+(a+b)x=(1+ax)(1+(1+ax)^{-1}bx),}

%\sni
qui est produit de deux \elts \ivsz.

Dans un \alo le radical de Jacobson
est égal à l'ensemble des \elts
noninversibles (\llec est invité\e
à en écrire la \prcoz). \index{Jacobson!radical de ---}
En \clama le radical de Jacobson est \care comme suit.

%--- Lemma{lemcRadJ}-----------------
\begin{lemmac}
\label{lemcRadJ}
Le radical de Jacobson est égal à l'intersection des \idemasz. 
\end{lemmac}
%--- end-lemma-----------------------------------------
%-----------------begin proof------------------
\begin{proof}
Si $a\in\Rad\gA$ et $a\notin\fm$ avec $\fm$ un \idemaz, on a $1\in\gen{a}+\fm$
ce qui signifie que pour un $x$, $1+xa\in\fm$, donc $1\in\fm$: contradiction.\\
Si $a\notin\Rad\gA$, il existe $x$ tel que $1+xa$ est non inversible. Donc il
existe un idéal strict contenant $1+xa$. Par le lemme de Zorn il existe un
\idemaz~$\fm$ contenant $1+xa$, et $a$ ne peut pas être dans $\fm$ car sinon
on aurait $1=(1+xa)-xa\in\fm$.\\
\Llec pourra remarquer que la preuve dit en fait ceci: un
\elt $x$ est dans l'intersection des \idemas \ssi est vérifiée l'implication:
$\gen{x,y}=\gen{1}
\Rightarrow  \gen{y}=\gen{1}$.
\end{proof}
%-----------------end proof------------------

\rem Nous avons raisonné avec un anneau non trivial. Si l'anneau est trivial l'intersection de l'ensemble (vide)
des \idemas est bien égale à~$\gen{0}$.
\eoe

\medskip \rdb\label{corpsdeHeyting}
Un \ixx{corps}{de Heyting}, ou simplement  un \ix{corps}, est par \dfn un \alo
dans lequel tout \elt noninversible est nul, autrement dit un \alo dont le radical
de Jacobson est réduit à~$0$.\index{Heyting!corps de ---}

En particulier, un \cdiz,
donc aussi l'anneau trivial,  est un corps.
Les nombres réels forment un corps qui {\em n'est pas} un corps
discret{\footnote{Nous utilisons la négation en
italique pour indiquer que l'affirmation correspondante, ici ce serait  \gui{$\RR$ est un \cdiz},
n'est pas prouvable en \comaz.}}. Même remarque pour le corps $\QQ_p$ des
nombres $p$-adiques ou celui des
séries formelles de Laurent $\gk(\!(T)\!)$ lorsque $\gk$ est un \cdiz.

\Llec vérifiera qu'un corps est un \cdi \ssi il est \zedz.

\smallskip
Le quotient d'un \alo par son radical de Jacobson est un corps, appelé
{\em  corps résiduel de l'\aloz}.
\index{corps!résiduel d'un anneau local}

%:     Lemma{lemZeDRaD}
\begin{lemma}\label{lemZeDRaD}
Si $\gA$ est \zedz, $\Rad\gA=\DA(0)$.
\end{lemma}
\begin{proof}
L'inclusion $\Rad\gA\supseteq\DA(0)$ est toujours vraie. Si maintenant
$\gA$ est \zed et $x\in\Rad\gA$, puisque l'on a une \egt $x^\ell(1-ax)=0$, il est clair que $x^{\ell}=0$.
\end{proof}
%

%:     Lemma{lemRadAX}
\begin{lemma}\label{lemRadAX}
Pour tout $\gA$, $\Rad(\AX)=\DA(0)[X]$.
\end{lemma}
\begin{proof}
Si $f\in\Rad(\AX)$, alors $1+Xf(X)\in\AX\eti$. On conclut avec le lemme~\ref{lemGaussJoyal}~\emph{\iref{i4lemPrimitf}.} 
\end{proof}
%

%:     Fact{fact1Rad}
\begin{fact}\label{fact1Rad} Soit $\gA$ un anneau et $\fa$ un \id
contenu dans $\Rad\gA$.
\begin{enumerate}
\item \label{i1fact1Rad}
$\Rad\gA=\pi_{\gA,\fa}^{-1}(\Rad\big(\gA\sur{\fa})\big)\supseteq\DA(\fa)$.
\item \label{i2fact1Rad}
$\gA$ est local \ssi $\gA\sur{\fa}$ est local.
\item \label{i3fact1Rad}
$\gA$ est local  et $\fa=\Rad\gA$ \ssi $\gA\sur{\fa}$ est un corps.
\end{enumerate}
%\end{enumerate}
\end{fact}

Le fait qui suit décrit une construction qui force un \mo à s'inverser
et un idéal à se radicaliser (pour plus de détails voir le paragraphe \gui{Dualité dans les anneaux commutatifs} \paref{secIDEFIL} et suivantes, et la section~\ref{subsecMoco}).

%:     Fact{fact2Rad}
\begin{fact}\label{fact2Rad}
Soient $U$ un \mo et $\fa$ un \id de $\gA$.
Nous considérons le \mo $S=U+\fa$. Notons $\gB=S^{-1}\gA$ et $\fb=\fa\gB$.
\begin{enumerate}
\item L'\id $\fb$ est contenu dans $\Rad\gB$.
\item L'anneau $\gB\sur{\fb}$ est isomorphe à $\gA_U\sur{\fa\gA_U}$.
\end{enumerate}
\end{fact}

 Par \dfn un  \ixx{anneau}{local \dcdz}
est  un \alo dont le corps
résiduel est un \cdiz.
Un tel anneau $\gA$ peut être \care par l'axiome suivant
%---  equation eqDefAlrd --------
\begin{equation}\label{eqDefAlrd}
\forall x\in \gA \qquad x\in \Ati  \;\; {\rm  ou} \;\;
    1+x\gA\,\subseteq\,  \Ati
\end{equation}
%---------------------end equation--------------
(\llec est invité\e à en écrire la \prcoz).

Par exemple l'anneau des entiers $p$-adiques,
 quoique {\em non} discret, est \dcdz.

 On obtient  un \alo \emph{non}
\dcd en prenant $\gK[u]_{1+\gen{u}}$, où $\gK$ est un corps \emph{non} discret
(par exemple le corps des séries formelles~$\gk(\!(t)\!)$,
où $\gk$ est un \cdiz).

\medskip \comm
 La différence un peu subtile qui sépare les \alos des
\alos \dcds se retrouve, en permutant l'addition et la multiplication,
dans la différence qui sépare les anneaux \sdz des anneaux intègres.
\\
En \clama un anneau \sdz est intègre; les deux notions n'ont cependant pas le
même contenu \algqz, c'est la raison pour laquelle on les distingue en
\comaz. \eoe

%--- Definition{defiresidzed}-----------
\begin{definition}\label{residzed}
%\label{defipseudolocal}~
\index{anneau!résiduellement \zedz}
\index{residuellement zer@résiduellement \zedz!anneau ---}
 Un anneau $\gA$ est dit  \emph{\plcz}
lorsque $\gA\sur{\Rad\gA}$ est \zedz. Même chose pour \emph{\rdt connexe}.
%-----------------end enum-----------
\end{definition}
%--- end-definition-------------------

%\smallskip
Puisqu'un corps est \zed \ssi c'est  un \cdiz, un \alo
est \dcd \ssi il est \plcz.

\medskip
\comm En \clama un anneau $\gA$ est dit semi-local si $\gA\sur{\Rad\gA}$
 est isomorphe à
un produit fini de \cdisz. Ceci implique que c'est un anneau
\plcz.
En fait l'hypothèse de finitude présente dans la notion d'anneau
semi-local est rarement décisive. La plupart des \thos de la littérature concernant
les anneaux semi-locaux s'appliquent aux anneaux \plcsz,
voire aux  \algbs (section \ref{secAlocglob}).
Sur une \dfn possible d'anneau semi-local en \coma
voir les exercices \ref{exo1semilocal} et~\ref{exo2semilocal}.
\eoe

%:  subsec{Idéaux premiers, maximaux}
\subsec{Idéaux premiers, maximaux}
%-----------------------------------------
En \comaz, un idéal d'un anneau $\gA$ est appelé un \emph{idéal maximal}
lorsque
l'anneau quotient est un corps{\footnote{Nous avons
jusqu'à maintenant utilisé
la notion d'\idema uniquement dans le cadre des preuves en \clamaz.
Une \dfn \cov s'imposait à un moment ou un autre.
En fait nous n'utiliserons que rarement
cette notion en \comaz. En règle \gnle elle est avantageusement
remplacée par la considération du radical de Jacobson, par exemple dans le cas  des \alosz.}}. Un \id  est appelé un \emph{\idepz}
lorsque l'anneau quotient est \sdzz.%
\index{ideal@idéal!maximal}%
\index{ideal@idéal!premier}%
\index{premier!idéal --- d'un anneau commutatif}%
\index{maximal!idéal ---}

 Ces \dfns coïncident avec les
\dfns usuelles si l'on se situe en \clamaz,  à ceci
près que nous tolérons l'anneau trivial comme corps
et donc l'\id $\gen{1}$ comme  \idema et comme \idepz.

Dans un anneau non trivial, un \id
 est strict,
maximal et détachable \ssi l'anneau quotient est un \cdi
non trivial, il est strict,
premier et détachable \ssi l'anneau quotient est un anneau intègre
non trivial.

\medskip \comm \label{CommIdeps}\stMF
Ce n'est pas sans une certaine appréhension que nous décrétons
l'\id $\gen{1}$ à la fois premier et maximal.
Cela nous obligera à dire \gui{\idep strict} ou \gui{\idema strict}
pour parler des \ideps et \idemas \gui{usuels}.
Fort heureusement ce sera très rare.
\perso{ce fut une décision difficile à prendre, fin 2006,
elle semble  convenablement motivée}

Nous pensons en fait qu'il y a eu \emph{une erreur de casting au départ}.
Imposer à un corps ou à un anneau intègre d'être non trivial, chose qui semblait éminemment raisonnable a priori, a conduit inconsciemment les
mathéma\-ticie\nz s à transformer de nombreux raisonnements
\cofs en raisonnements par l'absurde. Pour démontrer qu'un \id
construit au cours d'un calcul est égal à $\gen{1}$, on a pris l'habitude
de faire le raisonnement suivant: si ce n'était pas le cas, il serait contenu dans un \idema et le quotient serait un corps, dans lequel on aboutit
à la contradiction $0=1$.
%Si l'on ne tient pas compte du fait que l'utilisation
%(problématique à cause du lemme de Zorn) de l'\idemaz,
Ce raisonnement s'avère être un raisonnement
par l'absurde pour l'unique raison que l'on a commis l'erreur de casting:
 on a interdit à l'anneau trivial d'être un corps.
 Sans cette interdiction, on présenterait le raisonnement comme un raisonnement
 direct sous la forme suivante: montrons que tout \idema de l'anneau quotient contient $1$.
Nous reviendrons
sur ce point dans la section~\ref{subsecLGIdeMax}.

Par ailleurs, comme nous utiliserons les \ideps et les \idemas essentiellement
à titre heuristique, notre transgression de l'interdit usuel
concernant l'anneau trivial n'aura pratiquement aucune conséquence
pour la lecture.
En outre, \llec pourra remarquer que cette convention inhabituelle
n'oblige pas à modifier la plupart des résultats établis spécifiquement
en \clamaz, comme  le \plgaz\etoz~\ref{plca.basic},
 le fait\eto \ref{factPropCarFin} ou le lemme\eto \ref{lemcRadJ}:
il suffit par exemple\footnote{Le fait\eto \ref{factMoco} pourrait \egmt
être traité selon le même schéma, en supprimant d'ailleurs la
restriction au cas non trivial.} pour la \lon en un \idep $\fp$ de la
définir comme la \lon en le filtre 

\snic{S\eqdefi \sotq{x\in\gA}{x\in\fp \Rightarrow 1\in\fp}.}

%\sni
\looseness-1
Sur le fond nous pensons que les \maths sont
plus pures et plus élégantes lorsque l'on évite d'utiliser la négation
(cela interdit radicalement les raisonnements par l'absurde
par exemple). C'est pour cette raison que l'on ne trouvera dans cet
ouvrage aucune \dfn qui utilise la négation\footnote{Si c'était le cas,
ce serait dans un cadre où la négation équivaut à une
affirmation positive, parce que la \prt considérée est décidable.}.
\eoe

%:  subsec{Radical de Jacobson et unités}
%\subsec{Radical de Jacobson et unités dans une extension entière}
\subsect{Radical de Jacobson et unités dans une extension\\ entière}%
{Radical de Jacobson et unités dans une extension entière}%-----------------------------------------

%:     Theorem{thJacUnitEntieres}
\begin{theorem}\label{thJacUnitEntieres} Soit $\gk\subseteq\gA$ avec $\gA$ entier sur $\gk$.
\begin{enumerate}
\item Si $y\in\Ati$, alors $y^{-1}\in\gk[y]$.
\item $\gk\eti=\gk\cap\Ati$.
\item $\Rad\gk=\gk\cap\Rad\gA$ et l'\homo $\gA\to\gA\sur{\Rad(\gk)\gA}$ 
réfléchit les unités\footnote{Rappelons que l'on dit qu'un \homo $\rho:\gA\to\gB$
réfléchit les unités lorsque~$\rho^{-1}(\gB\eti)=\Ati$.}.
%un \elt $y\in\gA$ est \iv \ssi $\ov y$ est \iv dans $\gA\sur{\Rad(\gk)\gA}$.
\end{enumerate}
\end{theorem}
\begin{proof}
\emph{1.} Soit $y,z\in\gA$ tels que $yz=1$. On a une \rdi pour $z$: $z^n=a_{n-1}z^{n-1}+\cdots+a_0$ ($a_i\in\gk$).
En multipliant par $y^n$ on obtient $1=yQ(y)$ donc $z=Q(y)\in\gk[y]$.

\emph{2.} En particulier, si $y\in\gk$ est \iv dans $\gA$, son inverse $z$ est dans~$\gk$.

\emph{3.} Soit $x\in \gk\cap\Rad\gA$, pour tout $y\in\gk$, $1+xy$ est inversible dans
$\gA$ donc aussi dans~$\gk$. Ceci donne l'inclusion $\Rad\gk\supseteq\gk\cap\Rad\gA$.\\
Soit $x\in\Rad\gk$ et $b\in\gA$. Nous voulons montrer que $y=-1+xb$
est \ivz. On écrit une \rdi pour $b$:

\snic{b^n+a_{n-1}b^{n-1}+\cdots+a_0=0,}

%\sni
on multiplie par $x^n$ et l'on remplace $bx$ par $1+y$. Il vient un \pol
en $y$ à \coes dans $\gk$:
$y^n+\cdots+(1+a_{n-1}x+\cdots+a_0x^n)=0$. Donc,~$yR(y)=1+xS(x)$ est \iv
dans $\gk$, et $y$ est \iv dans $\gA$.\\
Soit maintenant $y\in\gA$ qui est \iv modulo $\Rad(\gk)\gA$. A fortiori
il est \iv modulo $\Rad\gA$, donc il est \ivz.
\end{proof}

%:     Theorem{thJacplc}
\begin{theorem}\label{thJacplc}
Soit $\gk\subseteq\gA$ avec $\gA$ entier sur $\gk$.
\vspace{-3pt}
\begin{enumerate}
\item $\gA$ est \zed \ssi $\gk$ est \zedz.
\item $\gA$ est \plc \ssi $\gk$ est \plcz. Dans ce cas
 $\Rad\gA=\DA\big(\Rad(\gk)\gA\big)$.
\item Si $\gA$ est local, $\gk$ \egmtz.
\end{enumerate}
\perso{point 2: Le lying over nous dit aussi que
$\Rad\gk=\gk\cap\Rad(\gk)\gA$. Aurait-on  $\Rad\gA=\sqrt{\Rad(\gk)\gA}$
dès que $\gA$ est entière sur $\gk$?}

\end{theorem}
\begin{proof} \emph{1.} Déjà connu (lemmes \ref{lemZrZr1} et \ref{lemZrZr2}).

\emph{2.} Par passage au quotient, le morphisme entier $\gk\to\gA$ donne un morphisme entier $\gk\sur{\Rad\gk}  \to \gA\sur{\Rad\gA}$,  qui est injectif
parce que $\Rad\gk=\gk\cap\Rad\gA$ (\tho \ref{thJacUnitEntieres}). Donc, les deux anneaux sont simultanément \zedsz.
Dans ce cas, notons $\fa = \Rad(\gk)\,\gA \subseteq \Rad\gA$. On a un
morphisme~entier 

\snic{\gk\sur{\Rad\gk}  \to \gA\sur{\fa},}

%\sni
donc $\gA\sur{\fa}$ est \zedz, de sorte que son radical de Jacobson
est égal a son radical nilpotent (lemme~\ref {lemZeDRaD}), i.e.
$\Rad(\gA)\sur\fa = \DA(\fa)\sur\fa$, et donc~$\Rad\gA = \DA(\fa).$

 \emph{3.} Résulte du \tho \ref{thJacUnitEntieres}, point \emph{2.}
\end{proof}
%

%%%%%%%%%%%%%%%%%%%%%%%%%%%%%%%%%%%%%%%%%%%%%%%%%%%%%%%%%%%%%%%%%%%%%%%%%
%%%%%%%%%%%%%%%%%%%%%%%%%%%%%%%%%%%%%%%%%%%%%%%%%%%%%%%%%%%%%%%%%%%%%%%%%
%%%%%%%%%%%%%%%%%%%%%%%%%%%%%%%%%%%%%%%%%%%%%%%%%%%%%%%%%%%%%%%%%%%%%%%%%
%%%%%%%%%%%%%%%%%%%%%%%%%%%%%%%%%%%%%%%%%%%%%%%%%%%%%%%%%%%%%%%%%%%%%%%%%
%%%%%%%%%%%%%%%%%%%%%%%%%%%%%%%%%%%%%%%%%%%%%%%%%%%%%%%%%%%%%%%%%%%%%%%%%
%--- Sec{Quatre lemmes importants}---secAloc2
\newpage	
\section{Quatre lemmes importants}
\label{secAloc2}
%-------------------------------
%:  --- Lemme de Nakayama --
Tout d'abord nous donnons quelques variantes du \gui{truc du \deterz} souvent
appelé \gui{lemme de Nakayama}.
Dans ce lemme la chose importante à souligner est que le module $M$ est \tfz.
\index{Lemme de Nakayama}
\index{truc du \deterz}

\CMnewtheorem{lemNak}{Lemme de Nakayama}{\itshape}

%:     Lemma{lemNaka}
\begin{lemNak}\label{lemNaka}  
\emph{(Le truc du \deterz)}%
\index{Nakayama!Lemme de ---} \\
Soient $M$ un \Amo \tf et $\fa$ un idéal de $\gA$.
%-----------------begin enum------------------
\begin{enumerate}
\item  Si $\fa\,M=M$, il existe $x\in\fa$ tel que $(1-x)\,M=0$.
\item  Si en outre $\fa\subseteq\Rad(\gA)$, alors $M=0$.
\item  Si $N\subseteq M$,  $\fa\,M+N=M$ et $\fa\subseteq\Rad(\gA)$, alors $M=N$.
\item  Si $\fa\subseteq\Rad(\gA)$ et  $X\subseteq M$ engendre $M/\fa M$ comme
$\gA/\fa$-module, alors~$X$ engendre $M$ comme \Amoz.
\end{enumerate}
%-----------------end enum------------------ 
\end{lemNak}
%--------- fin lemma ---------------------------------------------- 

%-----------------begin proof------------------
\begin{proof}
Nous montrons le point \emph{1} et laissons les autres en exercice, comme
consé\-quences
faciles.
Soit $V\in M^{n\times 1}$ un vecteur colonne formé avec des \gtrs de
$M$.
L'hypothèse signifie qu'il existe une matrice
$G\in\Mn(\fa)$ vérifiant  $GV=V$.
Donc $(\I_n-G)V=0$, et en prémultipliant par la matrice cotransposée de
$\I_n-G$, on obtient $\det(\I_n-G)V=0$. Or $\det(\I_n-G)=1-x$ avec $x\in\fa$.
\end{proof}
%-----------------end proof------------------

Les \mptfs sont localement libres au sens (faible) suivant: ils deviennent libres
lorsque l'on  localise en un \idepz. Prouver ceci revient à montrer le {\em  lemme
de la liberté locale} (ci-après) qui affirme qu'un \mptf sur un anneau local
est libre.

%%%%%%%%%%%%%%%%%%%%%%%%%%%%%%%%%%%%%%%%%
\CMnewtheorem{lemlilo}{Lemme de la liberté locale}{\itshape}
%:     Lemma{lelilo}
\begin{lemlilo}\label{lelilo}\index{Lemme de la liberté locale}%
Soit $\gA$ un anneau local. Tout \mptf sur $\gA$ est libre de rang
fini.
De manière \eqvez: toute matrice $F\in\GAn(\gA)$  est semblable à une \mprn
standard
%--------------------begin{array---------------
$$\I_{r,n}=\cmatrix{
\I_r&   &0_{r,n-r}  \cr
0_{n-r,r}&   &0_{n-r}
}.$$
%---------------------end{array-------------- 
\end{lemlilo}
%--------- fin lemma ---------------------------------------------- 

%\mni {\bf  Lemme de la liberté locale.}   \label{lelilo}
%\index{Lemme de la liberté locale}
%{\em
%Soit $\gA$ un anneau local. Tout \mptf sur $\gA$ est libre de rang
%fini.
%De manière \eqvez: toute matrice $F\in\GAn(\gA)$  est 
%semblable à une \mprn
%standard}
%%--------------------begin{array-------------
%$$\I_{r,n}=\cmatrix{
%\I_r&   &0_{r,n-r}  \cr
%0_{n-r,r}&   &0_{n-r}
%}.$$
%%---------------------end{array--------------
%%----fin du Lemme de la liberté locale ---

\medskip \rem
La formulation matricielle implique évidemment la première formulation, plus
abstraite.
Inversement si $M\oplus N=\Ae n$, dire que $M$ et $N$ sont libres (de rangs $r$ et
$n-r$) revient à dire qu'il y a une base de $\Ae n$
dont les $r$ premiers \elts forment une base de $M$ et les $n-r$ derniers une base
de $N$, en conséquence la \prn sur $M$ parallèlement à $N$ s'exprime sur
cette base par la matrice~$\I_{r,n}$.
\eoe

%----begin{proof------------------
\begin{Proof}{Première \demz, (preuve classique usuelle). }
Nous notons $x\mapsto \ov{x}$ le passage au corps résiduel. Si $M\subseteq
\Ae n$ est l'image d'une matrice de \prnz~$F$ et si $\gk$ est le corps résiduel
on considère une base de  $\gk^n$ qui
commence par des colonnes de $\ov{F}$ ($\Im \ov{F}$  est un sous-\evc de
dimension $r$) et se termine par des colonnes de $\I_n-\ov{F}$ ($\Im  (\I_n-
\ov{F})=\Ker \ov{F})$. En considérant les colonnes correspondantes de  $\Im
{F}$ et  $\Im (\I_n- {F})=\Ker F$ on obtient un relèvement de la base
résiduelle en  $n$ vecteurs dont le \deter est \rdt \ivz, donc \ivz.
Ces vecteurs forment une base de~$\Ae n$ et sur cette base
il est clair que la \prn admet pour matrice~$\I_{r,n}$.
\\
Notez que dans cette preuve on extrait un \sys libre maximal
 parmi les colonnes d'une matrice
à \coes dans un corps. Cela se fait usuellement
par la méthode du pivot de  Gauss.
Cela réclame donc que le corps résiduel soit discret.
\end{Proof}
%----end{proof------------------
%----begin{proof------------------
\begin{Proof}{Deuxième \demz, (preuve par Azumaya). }
 Contrairement à la précédente cette preuve ne suppose  pas que l'\alo
soit \dcdz.
Nous allons diago\-naliser la matrice $F$.
La preuve fonctionne avec un anneau local non \ncrt commutatif.\\
Appelons $f_1$ le vecteur colonne $F_{1..n,1}$ de la matrice $F$,
$(e_1,\ldots ,e_n)$ la base canonique de~$\Ae n$
et $\varphi$ l'\ali représentée par $F$.\\
-- Premier cas, $f_{1,1}$ est inversible. Alors,
$(f_1,e_2,\ldots ,e_n)$ est une base de $\Ae n$. Sur cette
base, l'\ali
$\varphi$ admet une matrice
%-----------------begin $$------------------
$$
G=\cmatrix{
    1   &      L      \cr
    0_{n-1,1}&   F_1      }.
$$
%-----------------end $$------------------
En écrivant $G^2=G$, on obtient $F_1^2=F_1$ et $LF_1=0$.
On définit alors la matrice $P=\cmatrix{
    1   &      L      \cr
    0_{n-1,1}&   \I_{n-1}      }$ et l'on obtient les \egtsz:
%-----------------begin $$------------------
$$
\begin{array}{rcl}
PGP^{-1}& =  &
\cmatrix{
    1   &      L      \cr
    0_{n-1,1}&   \I_{n-1}      }
\cmatrix{
    1   &      L      \cr
    0_{n-1,1}&   F_1      }
\cmatrix{
    1   &      -L      \cr
    0_{n-1,1}&   \I_{n-1}      }   \\[4mm]
& =  &\cmatrix{
    1   &      0_{1,n-1}      \cr
    0_{n-1,1}&   F_1      }.
\end{array}
$$
%-----------------end $$------------------
-- Deuxième cas, $1-f_{1,1}$ est inversible.
On applique le calcul précédent à la matrice
$\In-F$, qui
est donc semblable à une matrice
%-----------------begin $$----------------
$$\preskip.4em \postskip.4em
A=\cmatrix{
    1   &      0_{1,n-1}      \cr
    0_{n-1,1}&   F_1      },
$$
%-----------------end $$------------------
avec $F_1^2=F_1$, ce qui signifie que $F$ est semblable  à une
matrice
%-----------------begin $$----------------
$$\preskip.4em \postskip.0em
\In-A=\cmatrix{
    0   &      0_{1,n-1}      \cr
    0_{n-1,1}&   H_1      },
$$
%-----------------end $$------------------
avec $H_1^2=H_1$.\\
On termine la preuve par \recu sur $n$.
\end{Proof}
%----end{proof------------------

%--- Comment{lelilo}----------
\comm \label{comment lelilo}\relax
Du point de vue classique, tous les ensembles sont discrets, et
l'hypothèse correspondante est superflue dans la première preuve.
La deuxième preuve doit être considérée comme supérieure
à la première car son contenu \algq est plus \uvl que
celui de la première (qui ne peut être rendue complètement explicite que lorsque  l'\alo est
\dcdz). \eoe
%---end-comment--------------------

%:  --- Lemme de l'application \lnl

\medskip \rdb
Le lemme suivant peut être considéré comme une variante
du lemme de la liberté locale.

\CMnewtheorem{lemnllo}{Lemme de l'application \lnlz}{\itshape}
%:     Lemma{lelnllo}
\begin{lemnllo}\label{lelnllo}%
\index{Lemme de l'application localement@Lemme de l'application \lnlz}%
Soit $\gA$ un anneau local et~$\psi $ une \ali entre \Amos libres de rang fini.
\Propeq 
%-----------------begin enum------------------
\begin{enumerate}
\item L'\ali $\psi$ est \nlz.
\item L'\ali $\psi$ est \lnlz.
\item L'\ali $\psi$ a un rang fini $k$.
\end{enumerate}
%-----------------end enum------------------ 
\end{lemnllo}
%--------- fin lemma ---------------------------------------------- 

%\mni {\bf  Lemme de l'application \lnlz.}   \label{lelnllo}
%\index{Lemme de l'application localement@Lemme de l'application \lnlz}
%{\em
%Soit $\gA$ un anneau local et~$\psi $ une \ali entre \Amos libres 
%de rang fini.
%\Propeq 
%%-----------------begin enum------------------
%\begin{enumerate}
%\item $\psi$ est \nlz.
%\item $\psi$ est \lnlz.
%\item $\psi$ a un rang fini $k$.
%\end{enumerate}
%%-----------------end enum------------------
%}
%%----fin du Lemme de l'application \lnl---

%-----------------begin proof------------------
\begin{proof}
\emph{2} $\Rightarrow$ \emph{3.} L'\egt
 $\psi\,\varphi\,\psi=\psi$ implique que les \idds de $\psi$ sont \idmsz.
D'après le lemme \ref{lem.ide.idem} ces \ids sont engendrés par des \idmsz.
Comme un \idm d'un \alo est \ncrt égal à $0$ ou $1$, que $\cD_0(\psi
)=\gen{1}$ et $\cD_{r}(\psi )=\gen{0}$ pour $r$ assez grand, il existe un entier $k\geq 0$ tel que
$\cD_k(\psi )=\gen{1}$ et $\cD_{k+1}(\psi )=\gen{0}$.

\emph{3} $\Rightarrow$ \emph{1.}  Par hypothèse
 $\cD_k(\psi )=\gen{1}$, donc les mineurs d'ordre $k$ sont
\com et comme l'anneau est local un des mineurs d'ordre $k$ est \ivz.
Comme $\cD_{k+1}(\psi )=\gen{0}$, le
résultat est alors une conséquence du lemme de la liberté
\ref{lem pf libre}.
\end{proof}
%-----------------end proof------------------

Notez que la terminologie d'application \lnl est en partie justifiée par le
lemme précédent. Notez aussi que le \thref{theoremIFD} peut être
considéré comme plus \gnl que le lemme précédent.

%:  --- Lemme du nb de gtrs local   {lemnbgtrlo}
\CMnewtheorem{lemnbgl}{Lemme du nombre de \gtrs local}{\itshape}
%:     Lemma{lem}
\begin{lemnbgl} 
\label{lemnbgtrlo}\index{Lemme du nombre de \gtrs local}~\\
Soit  $M$ un \Amo \tfz.
%-----------------begin item------------------
\begin{enumerate}
\item Supposons $\gA$  local. 
\begin{enumerate}
\item [a.\phantom{*}] Le module $M$  est engendré par $k$
\elts \ssi
son \idf $\cF_k(M)$ est égal à $\gA$.
\item [b.\phantom{*}]  
Si en outre $\gA$ est \dcd et $M$ \pfz, le module admet une matrice de \pn dont tous les \coes sont
dans l'idéal maximal $\Rad\gA$.
\end{enumerate}
\item En \gnlz, pour $k\in\NN$, \propeq
%-----------------begin item------------------
\begin{enumerate}
\item [a.\phantom{*}] $\cF_k(M)$ est égal à $\gA$.
\item [b.\phantom{*}] Il existe des \eco $s_j$ tels que après extension des
scalaires
à chacun des $\gA[1/s_j]$, $M$ est engendré par $k$ \eltsz.
\item [c.\phantom{*}] Il existe des \moco $S_j$ tels que
chacun des  $M_{S_j}$ est engendré par $k$ \eltsz.
\item [d*.] Après \lon en n'importe quel \idepz,  $M$ est engendré par
$k$ \eltsz.
\item [e*.] Après \lon en n'importe quel \idemaz,  $M$ est engendré
par $k$ \eltsz.
\end{enumerate}
\end{enumerate} 
\end{lemnbgl}
%--------- fin lemma ---------------------------------------------- 
%-----------------begin proof------------------
\begin{proof} Il suffit de prouver les \eqvcs  pour un \mpf
 en raison du fait~\ref{facttfpf}.
\\
 Supposons $M$ engendré par $q$ \elts et notons $k'=q-k$.
\\
\emph{1.}
La condition est toujours \ncrz, même si l'anneau n'est pas local.
Soit une matrice de \pn $A\in\gA^{q\times m}$ pour $M$.
Si l'anneau est local et si  $\cF_k(M)=\gA$, puisque les mineurs d'ordre $k'$
sont \comz, l'un d'entre eux est \ivz. 
Par le lemme du mineur \ivz~\ref{lem.min.inv}, la matrice $A$ est \eqve à une matrice
%--------------------begin pmatrix---------------

\snic{
\cmatrix{
   \I_{k'}   &0_{k',m-k'}      \cr
    0_{k,k'}&     A_1},}

%\sni
et donc, la matrice $A_1\in\Ae {k\times (m-k')}$ est aussi une
 \mpn de~$M$.
Enfin, si l'anneau est \dcdz, on peut diminuer le nombre de \gtrs jusqu'à ce que la
matrice de \pn correspondante ait tous ses \coes dans le radical.

 \emph{2.}
\emph{a} $\Rightarrow$ \emph{b}.
La même preuve montre que l'on peut prendre pour $s_j$ les mineurs d'ordre $k'$
de~$A$.
\\
\emph{b} $\Rightarrow$ \emph{c}.
Immédiat.
\emph{c} $\Rightarrow$ \emph{a}.
 Dire que  $\cF_k(M)=\gA$ revient à résoudre le \sli
$\sum_\ell x_\ell s_\ell =1$, où les inconnues sont les $x_\ell $ et où les
$s_\ell $ sont les mineurs d'ordre $k'$ de la matrice  $A$. On peut donc
appliquer le principe local-global de base.
\\
\emph{a} $\Rightarrow$ \emph{d}.
Résulte du point \emph{1}.
\emph{d} $\Rightarrow$ \emph{e}.
Trivial.
\\
\emph{e} $\Rightarrow$ \emph{a}.
Ceci ne peut être prouvé qu'en \clama (d'où l'étoile
que nous avons mise à \emph{d} et \emph{e}).
On raisonne par l'absurde en prouvant la
contraposée. Si $\cF_k(M) \neq \gA$ soit
$\fp$ un \idema strict
contenant~$\cF_k(M)$. Après \lon en $\fp,$ on obtient
 $\cF_k(M_\fp)\subseteq \fp\gA_\fp \neq \gA_\fp$, et donc~$M_\fp$ n'est pas
engendré par $k$ \eltsz.
\end{proof}
%-----------------end proof------------------

%--- Comment----------
\comm  Ce lemme donne la \emph{vraie signification} de l'\egtz~$\cF_k(M)=\gA$:
on peut dire que $\cF_k(M)$ \gui{mesure} la possibilité pour le module
d'être localement engendré par $k$ \eltsz. D'où la \dfn qui suit.
\\
Voir aussi les exercices \ref{exoAutresIdF}, \ref{exoNbgenloc}
et~\ref{exoVariationLocGenerated}.
\eoe

%----Definition{deflocgenk}----------
\begin{definition}
\label{deflocgenk}
Un \mtf sera dit \ixd{localement engendré par~$k$ \eltsz}{module}
lorsqu'il vérifie les \prts \eqves du point \emph{2} du lemme du nombre de \gtrs local.\index{localement!module --- engendré par $k$ \eltsz}
\end{definition}
%--- end-definition------------------------------------

%%%%%%%%%%%%%%%%%%%%%%%%%%%%%%%%%%%%%%%%%%%%%%%%%%%%%%%%%%%%%%%%%%%%%%%%
%%%%%%%%%%%%%%%%%%%%%%%%%%%%%%%%%%%%%%%%%%%%%%%%%%%%%%%%%%%%%%%%%%%%%%%%%
%%%%%%%%%%%%%%%%%%%%%%%%%%%%%%%%%%%%%%%%%%%%%%%%%%%%%%%%%%%%%%%%%%%%%%%%%
%%%%%%%%%%%%%%%%%%%%%%%%%%%%%%%%%%%%%%%%%%%%%%%%%%%%%%%%%%%%%%%%%%%%%%%%%
%%%%%%%%%%%%%%%%%%%%%%%%%%%%%%%%%%%%%%%%%%%%%%%%%%%%%%%%%%%%%%%%%%%%%%%%%
%%%%%%%%%%%%%%%%%%%%%%%%%%%%%%%%%%%%%%%%%%%%%%%%%%%%%%%%%%%%%%%%%%%%%%%%%

\penalty-2500
\section{Localisation en \texorpdfstring{$1+\fa$}{1 + a}}
\label{secLoc1+fa}\relax

\medskip 
\Grandcadre{Soient $\fa$ un \id de $\gA$, $S := 1 + \fa$,  $\jmath:
\gA \to \gB:=\gA_{1+\fa}$ \\ l'\homo canonique, et
$\fb := \jmath(\fa)\gB$.}

\medskip 
On note que $\fb$  s'identifie à $S^{-1}\fa$
(fait \ref{fact.sexloc}) et  que $1+\fb\subseteq\gB\eti$ 
(fait \ref{fact2Rad}).
%:     Lemma{lemLoc1+a}
\begin{lemma}\label{lemLoc1+a}\relax
 {(Quotient de puissances de $\fa$ dans le localisé 
$\gA_{1+\fa}$)}
\\
Sous les hypothèses précédentes on a les résultats suivants. 
\begin{enumerate}
\item $\Ker\jmath\subseteq\fa$, $\gB = \jmath(\gA) + \fb$ et l'\homo canonique $\gA\sur\fa\to\gB\sur\fb$
est un \isoz. 
\item La \lon en $1+\fa$ est la même que la \lon en  $1+\fa^n$ ($n\geq1$),
donc $\Ker\jmath\subseteq\fa^n$, $\gB = \jmath(\gA) + \fb^n$ et $\gA\sur{\fa^n}\simeq\gB\sur{\fb^n}$.
\item  Pour
tous $p$, $q \in \NN$, $\jmath$ induit un \iso
$\fa^p\sur{\fa^{p+q}} \simarrow \fb^p\sur{\fb^{p+q}}$
\hbox{de \Amosz}.
\end{enumerate}
\end{lemma}
%--------- fin lemma ---------------------------------------------- 
%
\begin{proof}
\emph{1.} L'inclusion $\Ker\jmath\subseteq\fa$ est \imdez.
\\
Le fait que l'\homo $\gA\sur\fa\to\gB\sur\fb$ est un \iso tient à ce que l'on résout deux \pbs \uvls \eqvsz: dans le premier on doit annuler les \elts de 
$\fa$, dans le second, on doit en plus inverser les \elts de $1+\fa$, mais inverser~$1$ ne co\^ute rien. Enfin la surjectivité de ce morphisme signifie exactement que
$\gB = \jmath(\gA) + \fb$.

\emph{2.} Les \mos  $1+\fa$ et $1+\fa^n$ sont \eqvs  car $1-a$ divise $1-a^n$.

\emph{3.} Notons que $\fb^q=S^{-1}\fa^q=\fa^q\gB$.
En multipliant  $\gB = \jmath(\gA) + \fb^q$ par~$\fa^p$, on obtient 
$\fb^p = \jmath(\fa^p) + \fb^{p+q}$.
Donc, l'application  $\jmath$ induit une surjection \hbox{de \Amosz}
$\fa^p \twoheadrightarrow \fb^p\sur{\fb^{p+q}}$. 
Il reste à voir que son noyau est~$\fa^{p+q}$.
\hbox{Si $x \in \fa^p$} vérifie $\jmath(x) \in \fb^{p+q}$, 
il existe~$s \in 1+\fa$ tel que $sx \in \fa^{p+q}$,
et comme~$s$ est \iv modulo $\fa$, il l'est modulo $\fa^{p+q}$,
et donc $x \in \fa^{p+q}$.
\end{proof}

\vspace{2pt}
%:     Lemma{lemLocaliseFini}
\CMnewtheorem{lemlofi}{Lemme du localisé fini}{\itshape}
\begin{lemlofi}\label{lemLocaliseFini}%
\index{Lemme du localisé fini} 
Si $\fa$ est un \itf et $n\in\NN\etl$, on a les \eqvcs

\snic{\fb^n = \fb^{n+1}
\iff \fb^n = 0 \iff \fa^n = \fa^{n+1}.}

%\sni
Dans ce cas, 
\begin{enumerate}
\item on a $\;\fa^n = \Ker \jmath=\gen{1-e}$ avec $e$ \idmz, de sorte que

\snic{\gB=\gA_{1 + \fa} = \gA[1/e]=\aqo{\gA}{1-e},}

\item
si en outre $\gA$ est une \klgz,
alors $\gA\sur\fa$ est finie sur $\gk$ \ssi $\gB$ est finie sur~$\gk$.
\end{enumerate} 
\end{lemlofi}
%--------- fin lemme localise fini 

\vspace{2pt}

\begin{proof}
%\emph {1.}  
Si $\fb^n = \fb^{n+1}$, alors $\fb^n$ est \idm \tfz, donc $\fb^n=\gen{\vep}$ 
avec~$\vep$ \idmz.  
Mais puisque $\varepsilon \in \fb$, l'\idm $1 - \varepsilon$
est \ivz, donc égal \` a $1$, i.e. $\varepsilon = 0$, donc $\fb^n=0$.
 La troisième \eqvc provient de l'\iso $\fb^n\sur{\fb^{n+1}}\simeq\fa^n\sur{\fa^{n+1}}$
(lemme~\ref{lemLoc1+a}).

 \emph {1.} 
Puisque $\fa^n$ est \idm \tfz, $\fa^n=\gen{1-e}$ avec $e$ \idmz.  Le reste
découle ensuite du fait~\ref{fact.loc.idm}.

\emph {2.} 
Si $\gB$ est un $\gk$-\mtfz, il en est de même de $\gA\sur\fa \simeq
\gB\sur\fb$. Réciproquement, supposons 
que $\gA\sur\fa$ soit un $\gk$-\mtf et considérons la filtration de
$\gB$ par les puissances de $\fb$:

\snic {
0 = \fb^n \subseteq \fb^{n-1} \subseteq\cdots\subseteq \fb^2 
\subseteq \fb \subseteq \gB.
}

%\sni
Alors, chaque quotient $\fb^i\sur{\fb^{i+1}}$ est un $\gB\sur\fb$-\mtfz,
ou encore un $\gA\sur\fa$-\mtfz, et par suite un $\gk$-\mtfz. On en déduit
que $\gB$ est un $\gk$-\mtfz.
\end {proof}
%

%:     Lemma{lemLocalisezeddim}
\CMnewtheorem{lemlozed}{Lemme du localisé \zedz}{\itshape}
\begin{lemlozed} \index{Lemme du localisé zéro-dimensionnel} 
\label{lemLocalisezeddim}~\\
Soit $\fa$ un \itf de $\gA$ tel que le localisé $\gB=\gA_{1+\fa}$ soit \zedz.  Alors, il existe un entier $n$ et un \idm $e$ tels que 
$$\preskip.2em \postskip.3em 
\fa^n =\gen{1-e}\quad\hbox{  et  }\quad\gA_{1 + \fa} = \gA\big[\frac 1{e}\big]=\aqo{\gA}{1-e}. 
$$
Si, en outre, $\gA$ est une \klg \tf avec $\gk$ \zed (par exemple un \cdiz), alors $\gB$ est finie
sur $\gk$. 
\end{lemlozed}
%--------- fin lemme localise zerodim --------------------------------------- 
%
\begin{proof} On applique le lemme du localisé fini: puisque $\gB$ est \zed et $\fb$ \tfz, il existe un entier $n$ tel que
$\fb^n = \fb^{n+1}$.
\\
On termine avec le \nst faible \ref{thNst0} {car $\gB=\aqo\gA{1-e}$} est  \hbox{une \klgz} \tfz.   
\end {proof}

\rem
Soit $\fa$ un \itf d'un anneau $\gA$ tel que le localisé~$\gA_{1+\fa}$ soit
\zedz.  %%Nous avons vu (lemme du localisé \zed \paref{lemLocalisezeddim}),
L'application naturelle $\gA \to \gA_{1 + \fa}$ est donc surjective de noyau
$\bigcap_{k \ge 0} \fa^k=\fa^m$ avec $m$ tel que $\fa^m = \fa^{m+1}$. En outre,
$\fa^m$ est engendré par un \idm $1-e$ et
$\gA_{1+\fa}=\gA[1/e]$. On a alors:

\snic {
\bigcap_{k \ge 0} \fa^k = \big(0 : (0 : \fa^\infty)\big).
}

%\sni
Cette remarque peut-être utile pour le calcul. Supposons que $\gA =
\kuX\sur\ff$ où $\kuX = \gk[\Xn]$ est un anneau de \pols à $n$ \idtrs sur
un \cdi $\gk$ et $\ff = \gen {f_1, \ldots, f_s}$ un \itfz.  Soit $\fa$ un \itf
de $\kuX$ et $\ov\fa$ son image dans $\gA$. Alors, si $\gA_{1 + \ov\fa}$
est \zedz, la composée $\kuX \to \gA_{1 + \ov\fa}$ est surjective et son
noyau s'exprime de deux manières:

\snic {
\bigcap_{k \ge 0} (\ff + \fa^k) = \big(\ff : (\ff : \fa^\infty)\big).
}

%\sni
La formule de droite peut se révéler plus efficace en
calculant $(\ff : \fa^\infty)$ de la façon suivante:

\snic {\qquad\qquad\qquad\qquad\;
(\ff : \fa^\infty) = \bigcap_{j=1}^r (\ff : g_j^\infty)
\hbox { pour } \fa = \gen {g_1, \ldots, g_r}.
$\eoe$ }

\comm
En \clama un \idep $\fa$ de l'anneau~$\gA$ est dit \emph{isolé} s'il est à la fois minimal
et maximal dans l'ensemble des \ideps de $\gA$. 
Autrement dit s'il ne se compare à aucun autre \idep pour la relation d'inclusion. Dire que $\fa$ est maximal revient
à dire que~$\gA/\fa$ est \zedz. Dire que $\fa$ est minimal revient à dire
que~$\gA_S$ est \zedz, où $S=\gA\setminus \fa$. Mais si $\fa$ est supposé maximal, le \moz~$S$ est le saturé de $1+\fa$.\\
Inversement supposons que $\gB=\gA_{1+\fa}$ soit \zedz. Alors $\gA/\fa$ est 
\egmt \zed puisque  $\gA/\fa\simeq\gB/\fa\gB$.
Ainsi, lorsque~$\fa$ est en outre \tfz, on se retrouve avec un cas particulier
du lemme du localisé \zed \ref{lemLocalisezeddim}.  
Il est à noter que les \ideps isolés  dans la littérature interviennent en \gnl dans le contexte
d'anneaux \noes et que donc en \clama il sont automatiquement \tfz.
\eoe

%%%%%%%%%%%%%%%%%%%%%%%%%%%%%%%%%%%%%%%%%%%%%%%%%%%%%%%%%%%%%%%%%%%%%%%%%
%%%%%%%%%%%%%%%%%%%%%%%%%%%%%%%%%%%%%%%%%%%%%%%%%%%%%%%%%%%%%%%%%%%%%%%%%
%%%%%%%%%%%%%%%%%%%%%%%%%%%%%%%%%%%%%%%%%%%%%%%%%%%%%%%%%%%%%%%%%%%%%%%%%
%%%%%%%%%%%%%%%%%%%%%%%%%%%%%%%%%%%%%%%%%%%%%%%%%%%%%%%%%%%%%%%%%%%%%%%%%
%%%%%%%%%%%%%%%%%%%%%%%%%%%%%%%%%%%%%%%%%%%%%%%%%%%%%%%%%%%%%%%%%%%%%%%%%

%--- Section{secExloc} ------------------------
\penalty-2500
\section{Exemples d'anneaux locaux en géométrie algébrique}
%-----------------------------------------
\label{secExlocGeoAlg}

On se propose ici d'étudier
dans quelques cas \gui{l'\alg locale en un zéro d'un \sypz}.
Nous fixons le contexte suivant pour toute la section~\ref{secExlocGeoAlg}.

\Grandcadre{$\gk$ est un anneau, $\uf = f_1,\ldots,f_s\in\kuX= \kXn$,\\[1mm]
$\gA=\aqo{\kuX}{\uf}=\kxn$,  \\[1mm]
$(\uxi)=(\xin)\in\gk^n$  est un zéro du \sysz,\\[1mm]
$\fm_\uxi= \gen{x_1-\xi_1,\ldots,x_n-\xi_n}_\gA$ est
l'\id du point $\uxi$,\\[1mm]
$J(\uX)=\JJ_\uX(\uf)$ est la matrice jacobienne du \sysz.} 

Rappelons  que   $\gA=\gk\oplus\fm_\uxi$ (proposition \ref{prdfCaracAlg}).
Plus \prmtz, on a
avec l'\evn en $\uxi$ une suite exacte scindée
de \kmos
$$\preskip.0em \postskip.4em 
0\to\fm_\uxi\to\gA\vvvvvers{\ov g\,\mapsto\, g(\uxi)}\gk\to0, 
$$
et deux \homos de \klgs 
$$\preskip.3em \postskip.3em 
\gk\to\gA\to\gk 
$$
qui se composent en donnant $\Id_\gk$.

Rappelons aussi (\thref{thidptva}) que $\fm_\uxi$ est un \Amo \pf
(la \mpn est donnée explicitement).

%:--- Subsec  AlgLocZer
\subsec{Algèbre locale en un zéro}
\label{subsecAlgLocZer}
%-----------------------------------------

Dans la suite nous parlons du point $\uxi$ de $\gk^{n}$, mais il serait plus correct de dire \gui{le point $(\uxi)$}.
Dans la \dfn suivante la terminologie \emph{\alg locale en $\uxi$}
ne doit pas  prêter à confusion:
on ne prétend pas qu'il s'agisse d'un anneau local; on
mime simplement la construction de l'\alg locale donnée
dans le cas où~$\gk$ est un corps.

%:     definition{notaAlgLocEnxi}
\begin{definition}\label{notaAlgLocEnxi}
\emph{(Algèbre locale en un zéro d'un \sypz)}
\\
L'anneau $\gA_{1+\fm_\uxi}$  est appelé \emph{l'\alg locale en $\uxi$
du \syp  $(\uf)$}. On utilise aussi la notation abrégée $\gA_\uxi$
à la place de $\gA_{1+\fm_\uxi}$.% 
\index{algebre locale@\alg locale!en un zéro d'un \sypz}
\end{definition}

Nous notons  $\xi:\gA\to\gk$  le \crc d'\evn en $\uxi$. Il se factorise par le localisé \hbox{en $1+\fm_\uxi$} 
et l'on obtient un \crc
$\gA_{\uxi}\to\gk$. \\
On a donc $\gA_{\uxi}=\gk\oplus\fm_\uxi\gA_{\uxi}$
et des \isos canoniques
$$\preskip.2em \postskip.4em 
\gA_{\uxi}\big/{(\fm_\uxi\gA_{\uxi})}\simeq \gA\big/{\fm_\uxi} \simeq \gk. 
$$
%:     Fact{factAlgLocEnxi}
\begin{fact}\label{factAlgLocEnxi}
 \emph{(Si $\gk$ est un \cdiz, l'\alg $\gA_{\uxi}$ est un \aloz)}
\begin{enumerate}
\item Soit $\gk$  un \alo et soit $\fp=\Rad\gk$. On pose 
 $\fM=\fp\gA+\fm_\uxi$ et $\gC=\gA_{1+\fM}$.
Alors, $\gC$ est un \alo avec $\Rad(\gC)=\fM \gC$
\hbox{et $\gC\sur{\Rad\gC}\simeq \gk\sur{\fp}$}.
 
%\penalty-2500
\item Si $\gk$ est un \cdiz, on a les résultats suivants.
\begin{enumerate}
\item L'anneau $\gA_{\uxi}$ est un \alo avec $\Rad\gA_{\uxi}=\fm_\uxi\gA_{\uxi}$  et son corps résiduel est (canoniquement isomorphe à)~$\gk$.
\item Les anneaux $\gA$ et $\gA_{\uxi}$ sont  \noes \cohsz, et $\gA$  est \fdiz.
\item $\bigcap_{r\in\NN}\bigl(\fm_\uxi\gA_{\uxi}\bigr)^r=0.$
\end{enumerate}
\end{enumerate}
\end{fact}
\begin{proof}
\emph{1.} On a $\gC\sur{\fM\gC} \simeq \gA\sur{\fm_\xi} = \gk/\fp$ d'après le
point \emph {2} du fait \ref{fact2Rad}. On termine en utilisant le
point \emph {3} du fait \ref{fact1Rad}.
 
\emph{2a.} Résulte de \emph{1.}
 
\emph{2b.} L'anneau $\gA$ est \coh \fdi d'après le \thref{thpolcohfd}.
On en déduit que  $\gA_{\uxi}$ est \cohz.
\\
Pour la \noet on renvoie à \cite[VIII.1.5]{MRR}.
 
\emph{2c.} Vu les points \emph{2a} et \emph{2b}, il s'agit d'un cas particulier du \tho d'intersection de Krull (\cite[VIII.2.8]{MRR}).
\end{proof}
%

%-% ENTRE NOUS
\entrenous{

1) On se dit que $\gA_{\uxi}$ doit être \fdiz. 
En fait il suffit prouver qu'il est discret, mais cela a l'air un peu subtil.

Certainement les gens des bases standard savent faire cela. 

Ahh, l'\alg locale ...!

2)
On se demande bien si \emph{2c)} ne serait pas vrai avec pour $\gk$
un anneau commutatif arbitraire. Sinon je serai curieux d'avoir un contre-exemple. 
}
%-% Fin ENTRENOUS

%%%%%%%%%%%%%
\subsubsection*{Espace tangent  en un zéro}%
\index{tangent!espace ---}\index{espace tangent}

%\smallskip 
Dans la suite nous notons $\partial_jf$ pour $\Dpp{f}{X_j}$. Ainsi la matrice jacobienne du \sysz, que nous avons notée $J=J(\uX)$, se visualise comme suit:
$$
\bordercmatrix [\lbrack\rbrack]{
    & X_1                     & X_2                     &\cdots  & X_n \cr
f_1 & \partial_1 {f_1} &\partial_2 {f_1}  &\cdots  &\partial_n {f_1} \cr
f_2 & \partial_1 {f_2} &\partial_2 {f_2}  &\cdots  &\partial_n {f_2} \cr
\;\vdots & \vdots                  &                         &        & \vdots              \cr
f_i & \vdots                  &                         &        & \vdots              \cr
\;\vdots & \vdots                  &                         &        & \vdots              \cr
f_s & \partial_1 {f_s} &\partial_2 {f_s}  &\cdots  &\partial_n {f_s} \cr
}~=~J.
$$
La congruence ci-dessous est \imdez,
pour $f \in \gk[\uX]$:
%  equation label {eqModM2}
\begin{equation}\label {eqModM2}
f(\uX) \equiv f(\uxi) +
\som_{j=1}^n (X_j - \xi_j)\, \partial_j f(\uxi)
\;\mod \gen {X_1-\xi_1, \ldots, X_n-\xi_n}^2
\end{equation}
%  end-equation
En spécialisant $\uX$ en $\ux$ on obtient dans $\gA$ la congruence fondamentale:
%  equation label {eq2ModM2}
\begin{equation}\label {eq2ModM2}
f(\ux) \equiv f(\uxi)+ 
\som_{j=1}^n  \,(x_j - \xi_j) \,\partial_j f(\uxi)
 \,\;\mod {\fm_\uxi}^2
\end{equation}
%  end-equation
Nous laissons \alec le soin de vérifier que le noyau de
$J(\uxi)$  ne dépend que
de l'\id $\gen{\lfs}$ et du point $\uxi$. C'est un sous-\kmo de $\gk^n$ qui peut être appelé \emph{l'espace tangent en $\uxi$
au schéma affine sur $\gk$ défini par $\gA$}. Nous le noterons
$\rT_\uxi(\gA\sur\gk)$ ou $\rT_\uxi$.

Cette terminologie est raisonnable en \gaq (i.e., lorsque $\gk$ est un \cdiz), 
au moins dans le cas où $\gA$
est intègre: on a une \vrt définie
comme intersection d'hypersurfaces $f_i=0$, et l'espace tangent en $\uxi$ à la \vrt
est l'intersection des espaces tangents aux hypersurfaces qui la définissent. 

Dans cette même situation (\cdi à la base), le zéro $\uxi$ du \syp 
est appelé un \emph{point lisse} ou \emph{régulier}, ou \emph{non singulier} (du schéma affine ou encore de la \vrt correspondante)
lorsque la dimension de l'espace tangent en $\uxi$ est égale à la dimension\footnote{Si $\gA$ est intègre, cette dimension
ne dépend pas de $\uxi$ et peut être définie 
via une mise en position de Noether. Dans le cas \gnlz, 
il faut considérer la \ddk de l'anneau~$\gA_{\uxi}$.}
de la \vrt au point $\uxi$.
Un point qui n'est pas régulier est appelé \emph{singulier}.

\smallskip Nous donnons maintenant une interprétation plus abstraite de l'espace tangent, en termes d'espace de \dvnsz. Ceci fonctionne avec à la base un anneau commutatif $\gk$ arbitraire.

Pour une \klg $\gB$ et un \crc $\xi:\gB\to\gk$ on définit une \emph{$\gk$-\dvn au point $\xi$ de $\gB$} comme une forme $\gk$-\lin $d:\gB\to\gk$ qui vérifie la règle de Leibniz, i.e. en notant $f(\xi)$ pour $\xi(f)$:%
\index{derivation@dérivation!en un point (un caractère) d'une \algz}
$$
d(fg)=f(\xi)d(g)+g(\xi)d(f).
$$
Ceci implique en particulier $d(1)=0$ (écrire $1=1\times 1$), et donc $d(\alpha)=0$ pour $\alpha\in\gk$.
Nous noterons $\DBxk$ le \kmo des $\gk$-\dvns de~$\gB$ au point $\xi$.
\\
Cette notation est légèrement abusive. 
En fait si l'on note $\gk'$ l'anneau $\gk$ muni de la structure de \Bmo donnée par $\xi$, la notation de la \dfnz~\ref{defiDeriv} serait $\Der \gk\gB{\gk'}$, d'ailleurs muni de sa structure \hbox{de \Bmoz}.

\smallskip  
Nous allons voir que l'espace tangent en $\uxi$ à $\gA$ et le \kmo des $\gk$-\dvns de~$\gA$ en~$\xi$  sont naturellement isomorphes. 

%:     Proposition{propTangent}
\begin{proposition}\label{propTangent}
\emph{ ($\rT_\uxi(\gA\sur\gk)$, $\DAxk$,
%\dvns au point $\uxi$, 
et $(\fm_\uxi/{\fm_\uxi}^2)\sta$)}\\
On note $\fm=\fm_{\uxi}$ et l'on rappelle la notation $\rT_\uxi(\gA\sur\gk)=\Ker J(\uxi)$.
\begin{enumerate}
\item Pour $u=(\un)\in\gk^n$, notons $D_u : \gk[\uX]\to\gk$  la forme $\gk$-\lin  
définie par
$$\preskip-.20em \postskip.4em\ndsp 
D_u(f) = \sum_{j=1}^n {\partial_j f}(\uxi)\; u_j. 
$$
%\sni
C'est une \dvn au point $\uxi$, on a $u_j = D_u(X_j) = D_u(X_j - \xi_j)$,
\linebreak 
et l'application 
$$\preskip-.20em \postskip.4em
u \mapsto D_u, \; \gk^n\to \DkXxk
$$
%\sni
 est un \iso $\gk$-\linz.

\item Si $u \in \Ker J(\uxi) \subseteq\gk^n$, alors $D_u$ passe au quotient modulo
$\gen {f_1, \ldots, f_s}$ et fournit une $\gk$-\dvn au point $\uxi$, $\Delta_u :
\gA\to\gk$. \\
On a
$u_j = \Delta_u(x_j) = \Delta_u(x_j - \xi_j)$, et l'application 
$$\preskip.4em \postskip.4em
u \mapsto \Delta_u, \; \Ker J(\uxi)\to \DAxk
$$
 est
un \iso $\gk$-\linz.

\item En outre,  $\Delta_u(\fm^2) = 0$ et l'on obtient, par restriction à $\fm$ et passage au quotient modulo $\fm^2$, 
une forme $\gk$-\lin $\delta_u:\fm\sur{\fm^2} \to \gk$. 
On construit ainsi une \kli $u \mapsto \delta_u$
de $\Ker J(\uxi)$ dans $(\fm\sur{\fm^2})\sta$.

\item Réciproquement, à $\delta \in (\fm\sur{\fm^2})\sta$, on
associe $u \in \gk^n$ défini par 

\snic {u_j = \delta\big((x_j-\xi_j) \bmod\fm^2\big).}

%\sni
Alors, $u$ appartient à $\Ker J(\uxi)$.

\item Les deux applications définies en 3. et 4., 

\snic{\Ker J(\uxi) \to (\fm\sur{\fm^2})\sta\quad$ et
$\quad(\fm\sur{\fm^2})\sta \to \Ker J(\uxi),}

%\sni
sont des \isos $\gk$-\lins
réciproques l'un de l'autre.

\end{enumerate}
\end {proposition}

\begin {proof}
\emph {1.}
Simple \vfn laissée \alecz.

\emph {2.}
Pour n'importe quel $u\in \gk^n$, on vérifie facilement
que l'ensemble

\snic {
\sotq {f \in \gk[\uX]}{D_u(f) = 0 \hbox { et } f(\uxi) = 0}
}

%\sni
est un \id de $\gk[\uX]$. Si $u\in\Ker J(\uxi)$, on a $D_u(f_i) = 0$
par définition \hbox{(et $f_i(\uxi) = 0$)}; on en déduit
que $D_u$ est nulle sur $\gen {f_1, \ldots, f_s}$.

\emph {3.}
Pour voir que $\Delta_u(\fm^2) = 0$, on utilise $\Delta_u(fg) =
f(\uxi)\Delta_u(g) + g(\uxi)\Delta_u(f)$ et $f(\uxi) = g(\uxi) =
0$ pour $f$, $g \in \fm$.

\emph {4.}
La congruence (\ref{eq2ModM2}) pour $f = f_i$ est $\sum_{j=1}^n (x_j - \xi_j)
{\partial_j f_i}(\uxi) \in \fm^2$. Ceci donne en appliquant
$\delta$, l'\egt $\sum_{j=1}^n u_j {\partial_j f_i}(\uxi) =
0$, i.e. $u \in \Ker J(\uxi)$.

\emph {5.}
Soit $\delta \in (\fm\sur{\fm^2})\sta$ et $u \in \Ker J(\uxi)$ l'\elt
correspondant; il faut montrer que $\delta_u = \delta$, ce qui revient à
vérifier, pour $f \in \fm$:

\snic {
\delta(f \bmod\fm^2) = \sum_{j=1}^n {\partial_j f}(\uxi)
\delta\big((x_j-\xi_j)\bmod\fm^2\big)
,}

%\sni
mais ceci découle de $(\ref{eq2ModM2})$.

Réciproquement, soit $u \in \Ker J(\uxi)$ et $v \in \Ker J(\uxi)$ l'\elt
correspondant à $\delta_u$; il faut voir que $v = u$; cela revient
à vérifier $\delta_u\big((x_j - \xi_j)\bmod\fm^2\big) = u_j$,
\egt qui a déjà été constatée.
\end {proof}

\rem Notons que la \dfn que nous avons donnée de l'espace tangent $\rT_\uxi(\gA\sur\gk)$,  naturelle et intuitive, fait voir celui-ci
comme un sous-module de $\gk^n$, où $n$ est le nombre de \gtrs de la \klg \pf $\gA$. Il faut donc lui préférer la \dfn plus 
abstraite~$\Der\gk\gA\xi$, ou~$\fm_\uxi/{\fm_\uxi}^2$,
qui est plus intrinsèque, puisqu'elle ne dépend que de la \klgz~$\gA$
et du \crc $\xi:\gA\to\gk$, sans tenir compte de la \pn choisie pour~$\gA$ (en fait seule la structure de la localisée~$\gA_\uxi$ intervient).  
\eoe

%%%%%%%%%%%%%
\subsubsection*{Espace  cotangent en un zéro}
%\addcontentsline{toc}{subsection}{}

 De manière \gnlez,  on a aussi la notion duale d'\emph{espace cotangent en $\uxi$}.
Nous le définirons ici %provisoirement 
comme le conoyau de
la transposée $\tra{\,J(\uxi)}$. 
En fait, il s'agit d'un \kmo
qui est intrinsèquement attaché à l'\alg $\gA$ et au \crcz~$\xi$,
car il peut aussi être défini de manière formelle comme \gui{l'espace des \diles au point $\uxi$}. Nous ne développerons pas ce point ici. 

Le \tho fondamental qui suit implique que l'espace tangent
est canoniquement isomorphe au dual de l'espace cotangent
(le fait \ref{factDualReflexif}~\emph{2} appliqué à $\tra J$
donne $(\Coker\tra J)\sta\simeq \Ker J$ 
puisque $\tra{(\tra J)}=J$). 
Par contre, lorsque l'on travaille avec un anneau arbitraire $\gk$, l'espace cotangent n'est pas \ncrt isomorphe au dual de l'espace tangent.

Lorsqu'un \Bmo $M$ admet une \mpn $W$ sur un \sgr $(\yn)$, si $\fb$ est un \id de
$\gB$, par le changement d'anneau de base  $\pi_{\gB,\fb}:\gB\to\gB\sur\fb$, on obtient le~$\gB\sur\fb$-module $M\sur{\fb M}$ avec la \mpn $W\mod\fb$
sur le \sgr $(\ov{y_1},\ldots,\ov{y_n})$. 
\\
Avec  le \Amo $M=\fm_\uxi$ et l'\id $\fb=\fm_\uxi$, on obtient  pour \mpn du \kmo $\fm_\uxi/{\fm_\uxi}^2$ sur %le \sgr 
$(\ov{x_1-\xi_1},\ldots,\ov{x_n-\xi_n})$, la matrice~$\ov W=W\mod\fm_\uxi$, avec la matrice $W$ donnée dans le \thref{thidptva}. Celle-ci, à des colonnes nulles près, est la matrice $\tra{J(\uxi)}$. Le \tho qui suit
dit la même chose sous une forme précise.

%:     Theorem{thCotangent}
\begin{theorem}\label{thCotangent}
 \emph{(Espace cotangent en $\uxi$ et $\fm_\uxi \big/ {{\fm_\uxi}^2}$)}
Soit  $(e_i)_{i\in\lrbn}$ la base canonique
de $\gk^n$. Considérons l'\kli
$$\preskip.2em \postskip.2em 
\varphi : \gk^n \twoheadrightarrow {\fm_\uxi}/{{\fm_\uxi}^2},\quad
e_j \mapsto (x_j - \xi_j) \bmod {\fm_\uxi}^2. 
$$
Alors, $\varphi$ induit un \iso de \kmos $\;\Coker\tra{J(\uxi)} \simarrow 
{\fm_\uxi}/{{\fm_\uxi}^2}$.
\\
Ainsi, on a un \iso canonique $\Coker\tra{J(\uxi)} \simarrow 
{\fm_\uxi}\gA_\uxi/{({{\fm_\uxi}\gA_\uxi})^2}$. 
\end{theorem}
%--------- fin theorem ---------------------------------------------- 

\begin{proof} On suppose \spdg que $\uxi=\uze$ et
on utilise les notations du \tho \ref{thidptva}. La \mpn de $\fm_\uze$ pour le
\sgr $(\xn)$ est la matrice $W = [\,R_\ux\,|\,U\,]$ avec
$U(\uze)=\tra{J(\uze)}$. Comme la matrice $R_\ux\mod\fm_\uze$ est nulle, on obtient le résultat annoncé.
\\
 La dernière assertion est donnée par le lemme \ref{lemLoc1+a}~\emph{3}. 
\end{proof}

%:     Definition{defiEspCot}
\begin{definition}\label{defiEspCot}
On définit \emph{l'espace cotangent
en $\uxi$} comme étant le \kmo $\fm_\uxi\gA_\uxi/{({\fm_\uxi\gA_\uxi})^2}$, 
pour lequel seule intervient la structure de l'\alg locale en $\uxi$.
 \end{definition}
%--------- fin definition ---------------------------------------------- 

 Dans la suite de la section \ref{secExlocGeoAlg},
 nous étudions quelques exemples d'\algs locales en des zéros de \sypsz,
sans supposer que l'on a \ncrt à la base un \cdiz: $\gk$ est seulement un anneau commutatif.
Nous cherchons ici seulement à illustrer la situation \gmq en nous libérant si cela se peut de l'hypothèse \gui{\cdiz}, mais sans viser à donner le cadre le plus \gnl possible.

%:--- Subsec  subsecZedLocPtIsoleSimple
\subsec{Anneau local en un point isolé}
\label{subsecZedLocPtIsole}
%-----------------------------------------

L'idée qui guide ce paragraphe provient de la \gaq où l'anneau local
en $\uxi$ est \zed \ssi le point $\uxi$ est un zéro isolé, et où le zéro
isolé  est simple \ssi l'espace tangent est réduit à $0$.

%:HHH susbsub supprimé
%%--- SUBsubsection*{tho general}--------
%\subsubsection*{Un \tho général}
%%-----------------------------------------

%:HHH modifications assez importantes à partir de maintenant dans 
% tout le paragraphe \gui{Anneau local en un point isolé}

%:     Theorem{thJZS}--------------
\begin{theorem}
\label{thJZS} \emph{(Un zéro isolé simple)}\\
Dans le contexte décrit au début de la section \ref{secExlocGeoAlg}, \propeq
\begin{enumerate}
\item Le morphisme naturel $\,\gk\to \gA_{\uxi}\,$  est un \iso (autrement dit, l'\id $\fm_\uxi$ est nul dans $\gA_{\uxi}$). En bref, on écrit $\gk=\gA_{\uxi}$.
\item La matrice $\tra J(\uxi)$  est surjective, i.e. $1\in\cD_n\big(J(\uxi)\big)$.
\item L'espace cotangent en $\uxi$ est nul, i.e. $\fm_\uxi={\fm_\uxi}^2$.
\item L'\id $\fm_\uxi$ est engendré par un \idm  $1-e$  de  $\gA$.
Dans ce cas  les morphismes  naturels  $\gk\to\gA[1/e]\to\gA_{\uxi}$ sont des \isosz.
\item Il existe $g\in\gA$ tel que $g(\uxi)=1$ et $\gA[1/g]=\gk$.

\end{enumerate}

%\sni
Si en outre $\gk$ est un \cdi (ou un anneau \zedrz), on a aussi l'\eqvc avec la \prt suivante.
\begin{enumerate}\setcounter{enumi}{5}
\item L'espace tangent $\rT_\uxi$ est nul.
%
%\item 
%
\end{enumerate}
 
\end{theorem}
%--- end-theorem----------------------------------------
Voici comment on peut décrire la situation précédente en langage plus intuitif: l'\alg locale en $\uxi$ est une \gui{composante connexe
de $\gA$} (i.e., la \lon en $\uxi$ est la même que la \lon
en un \idmz~$e$) \gui{réduite à un point simple} (i.e., cette \klg%,
  %qui correspond aux zéros du \syp obtenu en rajoutant l'équation $e=1$,
 est isomorphe à~$\gk$).
 En termes de \vgqz, le point \emph{5} signifie qu'il y a un ouvert de Zariski contenant le point $\uxi$ dans lequel la \vrt est réduite à ce point.
  
%-----------------begin proof------------------
\begin{proof}
\emph{1} $\Leftrightarrow$ \emph{3.}  Par le lemme du localisé fini~\ref{lemLocaliseFini} avec $n=1$.
 
\emph{2} $\Leftrightarrow$ \emph{3.}  Par le \thref{thCotangent}.
 
\emph{3} $\Leftrightarrow$ \emph{4.} Par le lemme de l'\itf \idmz~\ref{lem.ide.idem}. \\
On obtient alors les \isos voulus par le fait~\ref{fact.loc.idm}, et donc   le point~\emph{5} avec $g=e$.
 
\emph{5} $\Rightarrow$ \emph{1.} 
L'\egt $g(\uxi)=1$
signifie \hbox{que $g\in1+\fm_\uxi$}. Ainsi l'anneau $\gA_\uxi$ est un localisé de $\gA[1/g]=\gk$, et il est égal à $\gk$ puisque $\gA_{\uxi}=\gk\oplus\fm_\uxi\gA_{\uxi}$.

\emph{3} $\Leftrightarrow$ \emph{6.} (Cas d'un \cdiz.) Puisque l'espace tangent est le dual du cotangent, \emph{3} implique toujours \emph{6}
Sur un \cdi une matrice est surjective \ssi sa transposée est injective,
ceci donne l'\eqvc de \emph{3} et \emph{6} (en considérant la matrice ${J(\uxi)}$). 
\end{proof}
%-----------------end proof------------------
\rem La différence entre le cas $s$ (nombre d'\eqnsz) $=n$ (nombre d'\idtrsz) et le cas $s>n$  n'est guère visible
dans le \tho précédent, mais elle est importante: si l'on perturbe
un \sys avec $s=n$ et si le corps de base est \acz, un zéro simple continue d'exister, 
légèrement perturbé. 
Dans le cas $s>n$ une perturbation fait en \gnl
disparaître le zéro. Mais ceci est une autre histoire, car il faut définir en \alg la notion de perturbation.
\eoe

Voici pour le cas d'un \cdi un résultat 
dans le même style que le \thref{thJZS}, mais plus \gnl et plus 
précis. Cela peut être vu \egmt comme une version locale
du \tho de Stickelberger (\thosz~\ref{thSPolZed} et \ref{thStickelberger}).
On notera cependant que, contrairement à ce qui se passe pour
le \tho de Stickelberger, la \dem du \thref{thJZScdi} ne fait pas intervenir
le \nst ou la mise en position de \Noez. Cependant, un \cdv à la Nagata
intervient dans l'appel au \thref{thNst0} pour l'implication \emph{7} $\Rightarrow$ \emph{8.}

%:     Theorem \label{thJZScdi}
\begin{theorem}\label{thJZScdi} \emph{(Zéro isolé)}
On suppose que $\gk$ est un \cdiz.  \Propeq
\begin{enumerate}
\item  L'\alg \smash{$\gA_\uxi$} est finie sur $\gk$.
\item  L'\alg  \smash{$\gA_\uxi$} est entière sur $\gk$.
\item  L'\alg  \smash{$\gA_\uxi$} est \zedez.
\item  L'\id \smash{$\fm_\uxi$ est nilpotent dans  $\gA_\uxi$}.
\item  Il existe \smash{$ r\in\NN$ tel que $\,{\fm_\uxi}^r={\fm_\uxi}^{r+1}$}
\item  Il existe $ r\in\NN$ tel que l'\id ${\fm_\uxi}^r$ est engendré par un \idm $1-e$,
le morphisme   $\gA \to \gA_\uxi$ est surjectif, \hbox{et
${\aqo\gA {1-e}\simeq\gA_\uxi\simeq\gA[1/e]}$}.
\item  Il existe $g\in \gA$  tel que \smash{$g(\uxi)=1$ et $\gA[1/g]=\gA_\uxi$}.
\item  Il existe $g\in \gA$  tel que $g(\uxi)=1$ et $\gA[1/g]$ est local \zedz.
\item  Il existe $h\in \gA$  tel que $h(\uxi)=1$ et $\gA[1/h]$ est finie sur $\gk$.
\end{enumerate}
Dans ce cas,  $\gA_\uxi$ est \stfe sur $\gk$, $(\gA_\uxi)\red=\gk$,
et si $m=[\gA_{\uxi}:\gk]$,  pour tout  $\ell \in\gA_\uxi$, on a $\rC{\gA_\uxi/\gk}(\ell )(T)=\big(T-\ell (\uxi)\big)^m$.
\end{theorem}
%--------- fin theorem ---------------------------------------------- 
%
\begin{proof}
Le lemme du localisé fini~\ref{lemLocaliseFini}, appliqué avec $\fa=\fm_\uxi$,
montre que \emph{4} équivaut à \emph{5} et implique \emph{1.}
 
 \emph{3} $\Rightarrow$ \emph{4.} Par le lemme du localisé \zed \ref{lemLocalisezeddim}.
 
On a \emph{1} $\Rightarrow$ \emph{2}, et puisque $\gk$ est un \cdiz, \emph{2} $\Rightarrow$ \emph{3}.
 
Ainsi les points \emph{1} à \emph{5} sont \eqvsz.
 
Le point \emph{5} implique que $\fm_\uxi^{r}$ est \idmz. Donc \emph{5} $\Rightarrow$ \emph{6} par le lemme de l'\itf \idm \ref{lem.ide.idem} et le fait~\ref{fact.loc.idm}. 

On note que $e\in1+\fm_\uxi^{r}\subseteq 1+\fm_\uxi$,
donc $e(\uxi)=1$. Donc \emph{6} implique~\emph{7} \hbox{avec $g=e$}.
 
\emph{7} $\Rightarrow$ \emph{8.} L'\alg $\gA[1/g]=\gA_\uxi$ est locale et \tfz, on conclut par le \thref{thNst0}.
 
\emph{8} $\Rightarrow$ \emph{9.} Prendre $h=g$.
 
\emph{9} $\Rightarrow$ \emph{1.} 
Parce que $\gA_\uxi$ est un localisé de $\gA[1/h]$.

Dans ce cas  $\gA_\uxi$ est \stfe sur $\gk$ car c'est une \alg finie et \pf
(\thref{propAlgFinPresfin}).
\\
Enfin l'\egt  
$\rC{\gA_\uxi/\gk}(\ell )(T)=\big(T-\ell (\uxi)\big)^m$
vient de ce que  $\ell -\ell (\uxi)$ est dans~$\fm$, donc est nilpotent dans~$\gA_\uxi$,
donc admet~$T^m$ comme \polcarz.
\end{proof}
\hum{
Insérer un contre exemple pour justifier l'hypothèse que $\gk$ est un corps dans le \thref{thJZScdi}?
}

%On pose la \dfn suivante, sans hypothèse sur $\gk$.
%:     Definition{defiZerIsole}
\begin{definition}\label{defiZerIsole} \emph{(Zéro isolé d'un \syp sur un anneau)}
\begin{enumerate}
\item Le zéro $\uxi$ du \sys  est un \emph{zéro isolé simple} 
(ou \emph{zéro simple}) si $\gA_{\uxi}=\gk$.%
\index{zero iso@zéro isolé simple!d'un \sypz}%
\index{simple!zero iso@zéro isolé ---}
\item Le zéro $\uxi$
du \sys  est un \emph{zéro isolé} si  $\gA_{\uxi}$ est finie sur $\gk$.
\index{zero iso@zéro isolé!d'un \sypz}
\item Si en outre $\gk$ est un \cdiz, la dimension de  $\gA_{\uxi}$ 
comme \kev est appelée la \emph{multiplicité} du zéro isolé~$\uxi$.%
\index{multiplicite@multiplicité!d'un zéro isolé (cas des corps)}
\end{enumerate}
\end{definition}
%--------- fin definition ----------------------------------------------
\rem
Le point \emph{1} est une abréviation par laquelle on entend
 \prmt que les \homos canoniques
$\gk\to\gA_{\uxi}\to\gk$ sont des \isosz.\\
Dans le point \emph{3} on voit que sur un \cdiz, un zéro isolé est  
simple 
\ssi il 
est de multiplicité~$1$. 
\eoe
%

%:  le tho de StickelBerger a ete mis en chap 4

%\newpage

%:--- Subsec{subsecLocPtLisse1} -------------  
\subsec{Anneau local en un point non singulier d'une courbe 
\lot intersection complète} 
\label{subsecLocPtLisse1}
%-----------------------------------------

On considère toujours le  contexte défini au début de la section~\ref{secExlocGeoAlg}, et l'on suppose $s=n-1$.
Autrement dit on a maintenant 

\smallskip \centerline{\fbox{un \sys de $n-1$ \eqns \polles à $n$ inconnues}}

\vspace{.1em} 
et l'on s'attend à ce que la \vrt correspondante soit \gui{une courbe}.

Nous allons voir que si le zéro $\uxi$ de la courbe est non singulier au sens intuitif que
l'espace cotangent au point $\uxi$ est un \kmo \pro de rang $1$, alors la situation \gui{locale}
est conforme à ce à quoi on s'attend, \cad ce à quoi nous ont habitué les points non singuliers des courbes en \gmt \dilez.

%:  thm: \label{thPointLisseCourbeIC}
\begin {theorem} {\emph{(L'\id d'un point non singulier d'une courbe 
\lot intersection complète)}}
       \label{thPointLisseCourbeIC}
Lorsque $s=n-1$ \propeq
\begin {enumerate}

\item
Le point $\uxi$ est non singulier au sens que $J(\uxi)$ est une matrice de rang $n-1$
sur~$\gk$. 
\item
L'espace cotangent en $\uxi$, \smash{$\fm_\uxi/{\fm_\uxi}^2$}, est un \kmo \pro de rang~$1$. 
\item
L'\id $\fm_\uxi$ est un \Amo \pro de rang $1$. 
\item
L'\id \smash{$\fm_\uxi\gA_\uxi$} est un $\gA_\uxi$-module \pro de rang $1$. 
\item
L'\id \smash{$\fm_\uxi\gA_\uxi$} est un $\gA_\uxi$-module libre de rang~$1$. 
\item
L'espace cotangent en $\uxi$, \smash{$\fm_\uxi/{\fm_\uxi}^2$}, est un \kmo libre de rang~$1$.%
\end{enumerate}
\end {theorem}
%--------- fin theorem ---------------------------------------------- 
%
\begin{proof} On rappelle que pour un anneau $\gB$, un \Bmo $M$ et un \idz~$\fb$ de~$\gB$ on obtient par \eds $\gB\sur\fb\otimes_\gB M\simeq  M\sur{\fb M}$.
En particulier, si $\fc$ est un \id de $\gB$ on obtient 
$(\gB\sur\fb)\otimes_\gB \fc\simeq  \fc\sur{\fb \fc}$.
\\ 
Mais l'\Bli surjective
naturelle $\fb \,\otimes\, \fc\to\fb \fc$ n'est pas toujours un \iso (c'est le cas si l'un des deux \ids est plat).

\emph{1} $\Leftrightarrow$ \emph{2.} 
En effet, $\tra J(\uxi)$ est une \mpn de l'espace cotangent.

 \emph{3} $\Rightarrow$ \emph{4.}  
En effet,  le $\gA_\uxi$-module $\fm_\uxi\gA_\uxi$ est obtenu à partir du
\Amo $\fm_\uxi$ par \eds de~$\gA$ \smash{à~$\gA_\uxi$}.

 \emph{4} $\Rightarrow$ \emph{2} et  \emph{5} $\Rightarrow$ \emph{6.}
En effet, le \kmo $\fm_\uxi/{\fm_\uxi}^2\simeq \fm_\uxi\gA_\uxi/({\fm_\uxi\gA_\uxi})^2$ est obtenu à partir du
$\gA_\uxi$-module $\fm_\uxi\gA_\uxi$ par \eds de~$\gA_\uxi$ à~$\gk\simeq \gA_\uxi/{\fm_\uxi\gA_\uxi}$ (voir le rappel du début).

\emph{2} $\Leftrightarrow$ \emph{3.} Cela résulte de la considération
de la \mpn de~$\fm_\uxi$ comme \Amo donnée au \thref{thidptva}
et du lemme~\ref{lemDnRz}. 
\\
Pour simplifier l'exposé nous traitons le cas $n=4$ 
avec $\uxi=\uze$. %(on peut se ramener à ce cas par translation).
\\
On a quatre variables $X_i$ et trois \pols 
$$\preskip.2em \postskip.3em 
\begin{array}{rcl} 
 f_1(\uX) & =  &  X_1a_1(\uX)+X_2a_2(\uX)+X_3a_3(\uX)+X_4a_4(\uX), \\[1mm] 
 f_2(\uX) & =  & X_1b_1(\uX)+X_2b_2(\uX)+X_3b_3(\uX)+X_4b_4(\uX) , \\[1mm] 
 f_3(\uX) & =  & X_1c_1(\uX)+X_2c_2(\uX)+X_3c_3(\uX)+X_4c_4(\uX). 
\end{array}
$$
Une \mpn de~${\fm_\uze}$ sur  $(x_1,x_2,x_3,x_4)$  est 

\snic{W(\ux)=
\cmatrix{
x_2&x_3&0& x_4&0&0&a_1(\ux)&b_1(\ux)&c_1(\ux)
\cr
-x_1&0&x_3& 0&x_4&0&a_2(\ux)&b_2(\ux)&c_2(\ux)
\cr
0&-x_1&-x_2&0&0&x_4&a_3(\ux)&b_3(\ux)&c_3(\ux)
\cr
0&0&0&-x_1&-x_2&x_3&a_4(\ux)&b_4(\ux)&c_4(\ux)
},
}

%\sni
ou encore $W(\ux)=[\,R_\ux\mid U(\ux)\,]$ avec

\snic{ U(\ux)=
\cmatrix{
a_1(\ux)&b_1(\ux)&c_1(\ux)
\cr
a_2(\ux)&b_2(\ux)&c_2(\ux)
\cr
a_3(\ux)&b_3(\ux)&c_3(\ux)
\cr
a_4(\ux)&b_4(\ux)&c_4(\ux)
}\;
$ et $\;\tra J(\uze)=U(\uze).
}

%\sni
On veut montrer que %les matrices 
$W(\ux)$
(\mpn du \Amoz~${\fm_\uze}$) \hbox{et $W(\uze)$} (\mpn du \kmo ${\fm_\uze}/{\fm_\uze}^2$)
sont simultanément de rang $n-1=3$.\\
Reportons nous au lemme \ref{lemDnRz}.
%\\
 Le point \emph{3} donne l'\egt $\cD_4\big(W(\ux)\big)=0$
(car~$\cD_4\big(U(\ux)\big)=0$). Et puisque 
%
%\snic{
$U(\uze)=U(\ux)\mod {\fm_\uze}$,
%}
%
%\sni
le point \emph{2} donne l'\eqvc

\snic{1\in\cD_{\gA,3}\big(W(\ux)\big) \iff 1\in \cD_{\gk,3}\big(U(\uze)\big)\iff 1\in\cD_{\gk,3}\big(W(\uze)\big).}

%\sni
\emph{1} $\Rightarrow$ \emph{5.}
On reprend les notations précédentes avec $n=4$ et $\uxi=\uze$. Puisque la matrice
 $\tra J(\uze)=U(\uze)$ est de rang $n-1$, il existe $\lambda_1$, \ldots, $\lambda_4\in\gk$ tels que 

\snic{\det \big(V(\uze)\big)= 1 $,\quad où \quad $V(\ux)=\cmatrix{
a_1(\ux)&b_1(\ux)&c_1(\ux)&\lambda_1
\cr
a_2(\ux)&b_2(\ux)&c_2(\ux)&\lambda_2
\cr
a_3(\ux)&b_3(\ux)&c_3(\ux)&\lambda_3
\cr
a_4(\ux)&b_4(\ux)&c_4(\ux)&\lambda_4
}\;.
}

%\sni
On en déduit que $\det \big(V(\ux)\big)\in 1+\fm_\uxi$, et donc $V(\ux)\in\GL_4(\gA_\uxi)$. Or 

\snic{[\,x_1\;x_2\;x_3\;x_4\,]\,V=[\,0\;\;0\;\;0\;\;y\,]\;$ avec $\;y=\sum_i\lambda_ix_i.}

%\sni
Ceci montre que $\gen{x_1,x_2,x_3,x_4}=\gen{y}$ dans $\gA_\uxi$.
Enfin $y$ est \ndz puisque le module~$\fm_\uxi$ est de rang $1$.
\end{proof}

%
%-% ENTRE NOUS
\entrenous{Dans le \tho précédent, on est vraiment surpris que cela remonte jusqu'à l'anneau global et pas seulement à l'anneau local.
}
%-% Fin ENTRENOUS

On notera $M^{\otimes_\gB r}$ la puissance tensorielle $r$-ième du \Bmo $M$.

%:     Theorem{th2PtlisseCourbeIC}
\begin{theorem}\label{th2PtlisseCourbeIC}
On suppose satisfaites les \prts \eqves du \thref{thPointLisseCourbeIC}, 
on note $\Omega$ l'espace cotangent $\fm_\uxi/{\fm_\uxi}^2$ et l'on considère
un \elt $p$  de  $\fm_\uxi$
qui est une $\gk$-base de  $\Omega$.
\begin{enumerate}
\item Pour chaque $r>0$, l'\kli naturelle $\Omega^{\otimes_\gk r}\to{\fm_\uxi}^r/{\fm_\uxi}^{r+1}$ est un \isoz.
En d'autres termes, la \klg graduée 
$$
\preskip.4em \postskip.4em 
\gr_{\fm_\uxi}(\gA):=\gk\oplus \bigoplus\nolimits_{r\geq 1}{\fm_\uxi}^r/{\fm_\uxi}^{r+1} 
$$ 
associée au couple $(\gA,\fm_\uxi)$ est (naturellement) isomorphe à l'\alg \smqz~$\gS_\gk(\Omega)$ du \kmo $\Omega$, elle-même isomorphe à $\gk[X]$ parce que~$\Omega$ est libre de rang~$1$. 
\item Si $\gk$ est un \cdi non trivial, $\gA_\uxi$ est un \emph{\adv discrète} au sens suivant:   tout \elt non nul de \smash{$\gA_\uxi$} s'écrit
de manière unique sous forme $up^{\ell}$ avec $u\in\Ati\!$ et  $\ell\geq0$.%
\index{anneau!de valuation discrète}\index{valuation!anneau de --- discrète}
\index{valuation!discrète}%
\end{enumerate}
\end{theorem}
%--------- fin theorem ---------------------------------------------- 
%
\begin{proof} 
On note $\fm=\fm_\uxi$. On remarque aussi que pour un \kmo \pro
de rang $1$, l'\alg \smq est égale à l'\alg tensorielle.

\emph{1.} On a un \iso naturel $\fm^{\otimes_\gA r}\simarrow \fm^r$ parce que $\fm$ est plat. Par l'\eds $\gA\to\gA\sur\fm=\gk$, 
les \Amos $\fm$ et  $\fm^r$ donnent les \kmos $\fm/\fm^2$ et  $\fm^r/\fm\fm^r=\fm^r/\fm^{r+1}$. 
\\
Puisque l'\eds commute avec le produit tensoriel,
on en déduit que l'\homo naturel $\left(\fm/\fm^2\right)^{\otimes_\gk r}\to \fm^r/\fm^{r+1}$ est un \iso de \kmosz.\\
 Puisque le \kmo $\fm/\fm^2$ admet la $\gk$-base $p  \mod \fm^2$,
 le \kmo  $\fm^r/\fm^{r+1}$ admet la base $p^r\mod \fm^{r+1}$. D'où un \iso
 de \klgs

\snic{\kX\simarrow\bigoplus_{r\in\NN}{\fm_\uxi}^r/{\fm_\uxi}^{r+1} =\gS_\gk(\Omega),}

%\sni
donné par $X\mapsto p$. En pratique, vue la filtration

\snic{\fm^r\subset\cdots\subset\fm^2\subset\fm\subset\gA,}

%\sni
dont tous les quotients sont des \kmos libres de rang $1$, le quotient~$\gA\sur{\fm^r}$ admet pour $\gk$-base $(1,p\dots,p^{r-1})$, avec pour $\ell<r$ 
le sous-\kmoz~$\fm^{\ell}\sur{\fm^r}$ qui admet la base  $(p^{\ell},\dots,p^{r-1})$.

\emph{2.} D'après le fait \ref{factAlgLocEnxi}~\emph{2} nous obtenons le résultat gr\^ace au calcul suivant: si $x\in\gA_\uxi$ est non nul, il est non nul dans un $\gA_{\uxi}/{\fm^r}$. Vue la filtration précédente
il existe 
un $\ell$ minimum tel que $x\in \fm^{\ell}$. Si $x\equiv ap^{\ell}\mod\fm^{\ell+1}$ avec $a\in\gk\eti$, on écrit $x=p^\ell (a+vp)$ avec $v\in\gA$ et $u=a+vp$
est \iv dans~$\gA_\uxi$. 
\end{proof}
%
 
%-% ENTRE NOUS
\entrenous{les choses seraient plus simples si l'on savait que  $\gA_\uxi$ est discret,
car il serait intègre, et ce serait un \adv en notre sens, ce qui n'est pas prouvé
ci-dessus puisque lorsque $x$ et $y$ sont dans de trop grandes puissances de $\fm$, ou nuls, sans qu'on le sache, alors on ne sait pas décider que l'un divise l'autre.
}
%-% Fin ENTRENOUS

% --- Subsec*{subsecLocPtLisse2} -------------  
%\subsubsection*{Exemple: un point lisse d'une courbe dans l'espace} 
%\label{subsecLocPtLisse2}
%-----------------------------------------

% 
%\textbf
\EXL{La courbe monomiale $t \mapsto (x_1=t^4, x_2=t^5, x_3=t^6)$}
\label{exlcourbemonomiale}

Pour $n_1$, $n_2$, $n_3 \in \NN^*$  premiers dans leur
ensemble, on définit la courbe monomiale $(x_1 = t^{n_1}, x_2 = t^{n_2},
x_3 = t^{n_3})$, plongée dans l'espace
affine de dimension $3$.
\\
Par \dfnz, l'\id de cette courbe paramétrée est, pour un anneau~$\gk$, le noyau du morphisme $\gk[X_1, X_2, X_3] \to
\gk[T]$ défini par $X_i \mapsto T^{n_i}$.
\\
On peut montrer que cet \id  est toujours défini sur $\ZZ$ et
engendré par trois \gtrsz. Ici on a choisi (voir le commentaire à la fin) le cas
particulier où~$(n_1,n_2,n_3) = (4,5,6)$, cas pour lequel deux relateurs suffisent:

\snic{x_1^3 = x_3^2\quad$ et $\quad x_2^2 = x_1x_3.}

%\sni
(Laissé en exercice \alecz.) On note

\snic{\gA = \gk[x_1,x_2,x_3]=\aqo{\gk[X_1,X_2,X_3]}{X_1^3- X_3^2,X_2^2 - X_1X_3}
}

%\sni
l'anneau de la courbe. Pour $t_0 \in \gk$, on
considère le point 

\snic{(\uxi) = (\xi_1,
\xi_2, \xi_3) = (t_0^{4}, t_0^{5}, t_0^{6})  ,}

%\sni
avec son \id $\fm = \gen {x_1-\xi_1, x_2-\xi_2, x_3-\xi_3}_\gA$. 
La condition pour que le point $\uxi$ soit non singulier, au sens que la matrice jacobienne $J$ évaluée en $\uxi$ est de rang $2$, 
est donnée par $t_0 \in \gk^\times$, car $\cD_2(J)=\gen{4t_0^{11}, 5t_0^{12},6t_0^{13}}$.
On suppose désormais $t_0\in\gk\eti$. 
Une \mpn  de $\fm$  pour
le \sgr $(x_1-\xi_1, x_2-\xi_2, x_3-\xi_3)$ est donnée par:

\snic{W = \cmatrix {
x_2 - \xi_2  &x_3 - \xi_3  &0            &x_1^2 + \xi_1x_1 + \xi_1^2 & -x_3\cr
-x_1 + \xi_1 &      0      &x_3 - \xi_3  &0                &x_2 + \xi_2 \cr
0            &-x_1 + \xi_1 &-x_2 + \xi_2 &-x_3 - \xi_3     &  -\xi_1\cr
}.}

%\sni
Nous savons qu'elle est de rang $2$.
  On constate  que $W_2, W_3 \in \gen {W_1, W_5}$.
On obtient donc une nouvelle
\mpn plus simple $V$ avec les seules colonnes $W_1, W_4, W_5$.
On rappelle d'une part que pour $B \in \Ae {n \times m}$, on a $(\Ae n/\Im
B)\sta \simeq \Ker \tra B$ (fait \ref{factDualReflexif}); d'autre part (exercice \ref{exoMatriceCorangUn}), que
pour une matrice $A \in \Mn(\gA)$ de rang $n-1$, on a $\Ker A = \Im
\wi{A}$ facteur direct dans $\Ae n$. En appliquant ceci à $B = V$
et $A = \tra \,V$, on obtient:

\snic{\fm\sta \simeq (\Ae 3/\Im V)\sta \simeq \Ker\tra V = \Im\tra {\,\wi V}
$ avec $\Im\tra {\,\wi V}$ facteur direct dans $\Ae 3.}

%\sni
On réalise ainsi explicitement le \Amo  $\fm\sta$, de rang constant $1$,  comme facteur direct dans $\Ae 3$. 
\eoe

%\begin {comment}
\medskip
\comm
De manière \gnle un sous-\mo  $M$ de $(\NN,+,0)$ a un complément $G$ fini \ssi il est engendré par une liste d'entiers premiers entre eux dans leur ensemble
(par exemple avec la courbe monomiale ci-dessus on définit $ M= n_1\NN + n_2\NN + n_3\NN$ engendré par $\so{n_1,n_2,n_3}$).
On dit que les entiers de $G$, sont les \emph{trous}  du \mo $M$. \\
Leur nombre $g := \#G$ est appelé le \emph{genre} de $M$.
\\
On a toujours $[\,2g, \infty\,[
\;\subseteq M$.  Les \mos $M$ pour lesquels $2g-1 \in G$
 sont dits \emph{\smqsz}. Cette terminologie rend compte du fait que, dans ce cas,
l'intervalle $\lrb{0 .. 2g-1}$ contient autant de trous que de non-trous, et qu'ils sont échangés par la symétrie $x \mapsto (2g-1) - x$.
\\
Par exemple, pour $a$, $b$ premiers entre eux, le \mo
${a\NN + b\NN}$
est \smq de genre $g = {(a-1)(b-1) \over 2}$.  On sait \carar de
manière combinatoire les \mos $n_1\NN + n_2\NN + n_3\NN$ qui sont
\smqsz. On démontre que c'est le cas \ssi l'\id de
la courbe  $(x_1=t^{n_1}, x_2=t^{n_2}, x_3=t^{n_3})$ est
engendré par 2 \eltsz. Par exemple $4\NN + 5\NN + 6\NN$ est \smqz, de genre
$4$, et ses trous sont $\so{1, 2, 3, 7}$.
\eoe

%-% ENTRE NOUS
\entrenous{

sections à rajouter éventuellement, si l'on a au moins des exemples pertinents

\smallskip sous-section  Anneau local en d'autres points d'une courbe plane ?

sous-section  Anneau local en un point non singulier d'une surface ?

\smallskip section Exemples en théorie des nombres ?
}
%-% Fin ENTRENOUS

%%%%%%%%%%%%%%%%%%%%%%%%%%%%%%%%%%%%%%%%%%%%%%%%%%%%%%%%%%%%%%%%%%%%%%%%%
%%%%%%%%%%%%%%%%%%%%%%%%%%%%%%%%%%%%%%%%%%%%%%%%%%%%%%%%%%%%%%%%%%%%%%%%%
%%%%%%%%%%%%%%%%%%%%%%%%%%%%%%%%%%%%%%%%%%%%%%%%%%%%%%%%%%%%%%%%%%%%%%%%%
%%%%%%%%%%%%%%%%%%%%%%%%%%%%%%%%%%%%%%%%%%%%%%%%%%%%%%%%%%%%%%%%%%%%%%%%%
%%%%%%%%%%%%%%%%%%%%%%%%%%%%%%%%%%%%%%%%%%%%%%%%%%%%%%%%%%%%%%%%%%%%%%%%%
\section{Anneaux \dcpsz}
\label{secRelIdm}

Les anneaux qui sont isomorphes à des produits finis d'\alos jouent
un rôle important dans la théorie classique des \alos henséliens
(par exemple dans \cite{Ray} ou \cite{Laf}).
De tels anneaux sont appelés des \emph{anneaux décomposés}
et un \alo est dit hensélien (en \clamaz) si toute extension finie est
un anneau décomposé.

Nous donnons dans cette section un début de l'approche \cov
pour la notion d'anneau décomposé. En fait, comme nous voulons éviter les
\pbs de \fcnz, nous allons introduire la notion, \cot plus pertinente,
d'anneau \dcpz.

Tout commence avec cette remarque simple mais importante: dans un
anneau commutatif les \idms sont toujours \gui{isolés}.

%:     Lemma{lemIdmIsoles}
\begin{lemma}\label{lemIdmIsoles}
Dans un anneau commutatif $\gA$
deux \idms égaux modulo $\Rad\gA$
sont égaux.
\end{lemma}
\begin{proof}
On  montre que l'\homo $\BB(\gA)\to\BB(\gA\sur{\Rad\gA}\!)$ est injectif: si un \idm $e$ est dans $\Rad\gA$, $1-e$ est \idm et inversible,
donc égal à $1$.
\end{proof}

\rem Ceci n'est plus du tout vrai en non commutatif: les \idms
d'un anneau de matrices carrées $\Mn(\gA)$ sont les \mprnsz; sur un corps
on obtient, par exemple en fixant le rang à $1$, une \vrt
connexe de dimension $>0$ sans aucun point isolé (si $n\geq2$).
\eoe

%: --- Subsubsec*{\'Eléments décomposables}
\penalty-2500
\subsec{\'Eléments décomposables}
%-------------------

%:     Definition{defiRelIdm}
\begin{definition}\label{defiRelIdm} Soit $\gA$ un anneau et $a\in\gA$.
 L'\elt $a$ est dit \emph{\dcpz\footnote{Il faut faire attention que cette terminologie entre en conflit
 avec la notion d'\idm in\dcp dans la mesure où tout \idm est un \elt \dcp de l'anneau.}} s'il existe
un \idm $e$ tel que:

\snic{\formule{
a \mod \gen{1-e} %\,\in \,(\aqo{\gA}{1-e})\eti
\,\mathrm{est\;\iv\;dans} \,\aqo{\gA}{1-e}
\et
\\[.5mm]
 a \mod \gen{e}\,\in \,\Rad(\aqo{\gA}{e}).}
 }
\end{definition}
\index{decomposable@décomposable!element@élément --- dans un anneau}

Rappelons en soulignant les analogies qu'un \elt $a$ possède un
quasi inverse \ssi il existe un \idm $e$ tel que:

\snic{\formule{
a \mod \gen{1-e} \,\mathrm{est\;\iv\;dans} \,\aqo{\gA}{1-e} \et
\\[.5mm]
 a \mod \gen{e} = 0\;\mathrm{dans} \;\aqo{\gA}{e},}
 }

%\sni
et qu'un \elt $a$ a pour annulateur un \idm  \ssi il existe un \idm $e$ tel que:

\snic{\formule{
a \mod \gen{1-e} \,
\mathrm{est\;r\E egulier\;dans} \,\aqo{\gA}{1-e} \et
\\[.5mm]
 a \mod \gen{e} = 0\;\mathrm{dans} \;\aqo{\gA}{e}.}
 }

%:     Proposition{prop1DecEltAnneau}
\begin{proposition}\label{prop1DecEltAnneau}
 Un \elt $a$ de $\gA$ est \dcp \ssi il existe $b$ tel que
\begin{enumerate}
\item  $b(1-ab)=0$,
\item  $a(1-ab)\in\Rad\gA$.
\end{enumerate}
En outre, l'\elt $b$ vérifiant ces conditions est unique,
et $ab=e$ est l'unique \idm de $\gA$ vérifiant  $\gen{a}=\gen{e}\mod \Rad\gA$.
\end{proposition}
%%%%%%%%%%%%%%%%%%%%%%%%%%%%%%%%%%%%%%%%%
%
\begin{proof} Supposons $a$ \dcpz.
Alors, dans le produit $\gA = \gA_1 \times \gA_2$, avec~$\gA_1 = \aqo\gA{1-e}$
et~$\gA_2 = \aqo\gA{e}$, on a $e = (1,0)$, $a = (a_1, a_2)$, avec~$a_1\in\gA_1\eti$ et $a_2 \in \Rad(\gA_2)$.
On pose $b = (a_1^{-1}, 0)$, et l'on
a bien

\snic{b(1-ab)=(b,0)-(b,0)(1,0)=0_\gA\;$ et $\;a(1-ab)=(0,a_2)\in\Rad\gA.}

%\sni En outre si $b=(b_1,b_2)$ vérifie $b(1-ab)=0$ et $a(1-ab)\in\Rad\gA$, on obtient 
%$(b_1,b_2)=(b_1^2a_1,b_2^2a_2)$ donc \dots 

%\sni
Supposons qu'un \elt $b$ vérifie

\snic{\formule{
 b(1-ab)=0 \et
\\[.5mm]
 a(1-ab)\in\Rad\gA.}
 }

%\sni
Alors, l'\elt $ab=e$ est un \idm et $a$ est \iv modulo $1-e$.
Par ailleurs, modulo $e$ on a $a=a(1-e)$ qui est dans $\Rad\gA$,
donc $a\mod e$ est dans  $\Rad(\aqo{\gA}{e})$.

 Voyons l'unicité. Si $b(1-ab)=0$ et $a(1-ab)\in\Rad\gA$, alors $e=ab$
est un \idm  tel que $\gen{a}=\gen{e}\mod \Rad\gA$. Cela le caractérise comme
\idm de  $\gA\sur{\Rad\gA}$, donc comme \idm de $\gA$.
Les \egts $be=b$  et $ba=e$ impliquent que $\big(b+(1-e)\big)\big(ae+(1-e)\big)=1$.
L'\elt $b+(1-e)$ est donc  déterminé de manière unique:
c'est l'inverse de $ae+(1-e)$.
Par suite, l'\elt $b$ est lui-même déter\-miné de manière unique.
  \end{proof}
%

%:     Definition{defiAdcp}
\begin{definition}\label{defiAdcp}
 On dit que l'anneau $\gA$ est \emph{\dcpz} si tout \elt est
\dcpz.%
\index{decomposable@décomposable!anneau ---}%
\index{anneau!decomposable@décomposable}
\end{definition}

%:     Fact{factDCP}
\begin{fact}\label{factDCP}~
\begin{enumerate}
\item \label{i1factDCP} Un produit d'anneaux est \dcp \ssi chacun des facteurs est \dcpz.
\item Un anneau \zed est \dcpz. Un \alo \dcd est \dcpz. Un anneau \dcp connexe est local \dcdz.
\item \label{i4factDCP} La structure d'anneau \dcp est purement équationnelle
(elle peut être définie au moyen de lois de compositions soumises
à des axiomes universels).
\end{enumerate}
\end{fact}

\begin{proof}
\emph{\ref{i4factDCP}}. On ajoute aux lois des anneaux commutatifs deux lois

\snic{a\mapsto b \;\hbox{ et  }\; (a,x)\mapsto y,}

%\sni
avec les axiomes $b=b^2a$ et  $\big(1+x(a^2b-a)\big)y=1$.
D'où $a^2b-a\in\Rad\gA$.

\emph{\ref{i1factDCP}.}  Résulte du point \emph{\ref{i4factDCP}}.
\end{proof}

\rem Si l'on note $b=a\esh$, alors $(a\esh)\esh=b\esh=a^2b$ et
 $\big((a\esh)\esh\big)\esh=a\esh$. En outre, $(a\esh)\esh$ et $a\esh$
 sont quasi inverses l'un de l'autre.
 \eoe

%%%%%%%%%%%%%%%%%%%%%%%%%%%%%%%%%%%%%%%%%

%: --- Subsubsec*{Relèvement des \idmsz}
\subsec{Relèvement des \idmsz}
%-------------------

%:     Definition{defi2RelIdm}
\begin{definition}\label{defi2RelIdm} Soit $\gA$ un anneau.
\begin{enumerate}
\item On dit que \emph{l'anneau $\gA$ relève les \idmsz} si l'\homo naturel 

\snic{\BB(\gA)\to\BB(\gA\sur{\Rad\gA}\!)}

%\sni
est bijectif, autrement dit
si tout \idm
du quotient $\gA\sur{\Rad\gA}$ se relève en un \idm de $\gA$.
\item On dit que l'anneau $\gA$ est \emph{décomposé} s'il est
\dcp et si $\BB(\gA)$ est bornée.
\end{enumerate}
\index{anneau!decompose@décomposé}
\index{decompose@décomposé!anneau ---}
\index{anneau!qui relève les \idmsz}
\end{definition}
%

%:     Proposition{propDecEltAnneau}
\begin{proposition}\label{propDecEltAnneau}~
\Propeq
\begin{enumerate}
\item  $\gA$ est \plc et relève les \idmsz.
\item $\gA$ est \dcpz.
\end{enumerate}
\end{proposition}

\begin{proof}
\emph{1} $\Rightarrow$ \emph{2.}
Puisque $\gA\sur{\Rad\gA}$ est \zed réduit, il existe un
\idm $e$ de  $\gA\sur{\Rad\gA}$ tel que $\gen{a}=\gen{e}\mod \Rad\gA$.
Cet \idm se relève en un \idm de $\gA$, que nous continuons
d'appeler $e$.\\
 L'\elt $a+(1-e)$ est \iv dans  $\gA\sur{\Rad\gA}$,
donc dans $\gA$. Donc,~$a$ est \iv dans $\aqo{\gA}{1-e}$.
Enfin, puisque  $\gen{a}=\gen{e}\mod \Rad\gA$, on obtient $a\in\Rad(\aqo{\gA}{e}).$

 \emph{2} $\Rightarrow$ \emph{1.} 
Notons $\pi:\gA\to\gA\sur{\Rad\gA}$ la \prn canonique.
Tout \elt $a$ de~$\gA$
vérifie $\gen{\pi(a)}=\gen{\pi(e)}$ pour un \idm $e$ de $\gA$. 
Le quotient est donc \zedz. 
Montrons que $\gA$ relève les \idmsz. \\
Si~$\pi(a)$ est  \idm
et si $e$ est l'\idm tel que $\gen{\pi(a)}=\gen{\pi(e)}$, 
alors~$\pi(a)=\pi(e)$.
\end{proof}

\comm Il est maintenant facile de voir qu'en \clama un anneau est décomposé
\ssi il est isomorphe à un produit fini d'\alosz. \eoe

%%%%%%%%%%%%%%%%%%%%%%%%%%%%%%%%%%%%%%%%%%%%%%%%%%%%%%%%%%%%%%%%%%%%%%%%%
%%%%%%%%%%%%%%%%%%%%%%%%%%%%%%%%%%%%%%%%%%%%%%%%%%%%%%%%%%%%%%%%%%%%%%%%%
%%%%%%%%%%%%%%%%%%%%%%%%%%%%%%%%%%%%%%%%%%%%%%%%%%%%%%%%%%%%%%%%%%%%%%%%%
%%%%%%%%%%%%%%%%%%%%%%%%%%%%%%%%%%%%%%%%%%%%%%%%%%%%%%%%%%%%%%%%%%%%%%%%%
%%%%%%%%%%%%%%%%%%%%%%%%%%%%%%%%%%%%%%%%%%%%%%%%%%%%%%%%%%%%%%%%%%%%%%%%%

%--- Sec{Anneau \lgb} --- secAlocglob
\section{Anneau \lgb}
\label{secAlocglob}
%-------------------------------

Nous introduisons dans cette section une notion qui \gns à la fois
celle d'\alo et celle d'anneau \zedz.
Ceci éclaire un certain nombre de faits communs à ces deux classes d'anneaux,
comme par exemple celui que les \mptfs sont quasi libres.

%: --- Subsubsec*{Définition et \plgc}
\subsec{Définitions et \plgc}
%-------------------

%:     Definition{defalgb}-----------
\begin{definition}
\label{defalgb}\label{defipseudolocal}~
\index{anneau!local-global}
\index{polynome@\pol!primitif par valeurs}
\index{primitif par valeurs!\pol ---}
%-----------------begin enum------------------
\begin{enumerate}
\item
On dit qu'un \pol $f\in\AXn$ \emph{représente (dans $\gA$)
l'\elt $a\in\gA$} s'il existe
$\ux\in\Ae n$ tel que $f(\ux)=a$.
\item
On dit qu'un \pol $f\in\AXn$ est \emph{primitif par valeurs}
si les valeurs de $f$ engendrent l'idéal
$\gen{1}$ (les variables étant évaluées dans~$\gA$).
\item Un anneau $\gA$ est dit \ixc{local-global}{anneau ---}
 si tout \pol \ppv  représente un \ivz.
\end{enumerate}
%-----------------end enum-----------
\end{definition}
%--- end-definition-------------------

\rem Tout \pol primitif par valeurs est primitif, donc si un anneau
possède la \prt que tout \pol primitif représente un \ivz, c'est un anneau \lgbz. Ceci correspond à une \dfn dans la littérature 
(anneau fortement U-\irdz) qui a précédé celle
d'\algbz. 
\eoe

%--- Fact{factlgb1}-----------------
\begin{fact}
\label{factlgb1} ~
%-----------------begin enum------------------
\begin{enumerate}
\item \label{i1factlgb1} Un anneau $\gA$ est \lgb  \ssi $\gA\sur{\Rad(\gA)}$ est \lgbz.
\item Un produit fini d'anneaux est \lgb \ssi chacun des anneaux est \lgbz.
\item Un \alo est \lgbz.
\item \label{i4factlgb1} Un anneau \plc est \lgbz.
\item \label{i5factlgb1} Un quotient d'un anneau \lgb
(resp. \plcz) est \lgb (resp. \plcz).
\item \label{i6factlgb1}
Soit $\gA$ un anneau réunion filtrante croissante de sous-anneaux $\gA_i$,
i.e. pour tous $i, j$, il existe $k$ tel que $\gA_i \cup \gA_j \subseteq
\gA_k$. Alors, si chaque $\gA_i$ est \lgbz, il en est de même de $\gA$.
\end{enumerate}
%-----------------end enum------------------
\end{fact}
%--- end-fact-----------------------------------------
%-----------------begin proof------------------
\begin{proof}
Nous laissons les trois premiers points en exercice.

\emph{\ref{i4factlgb1}.} Vu le point \emph{\ref{i1factlgb1}}, il suffit de traiter    le cas d'un anneau \zedrz. Ce cas se ramène 
   au cas (évident) d'un \cdi par la machinerie \lgbe \elrz.

\emph{\ref{i5factlgb1}.} Voyons le cas \lgb (l'autre cas est évident).
Soit $\gA$ un anneau \lgbz, $\fa$ un \id et
$f\in\gA[\uX]$ un \pol \ppv dans $\gA\sur\fa$.
Il y a donc des valeurs $p_1$, \ldots, $p_m$ de $f$  et un $a\in\fa$ tels
que $\gen{p_1,\ldots ,p_m,a}=\gen{1}$. Le \pol $g(\uX,T)=Tf(\uX)+(1-T)a$
est donc \ppvz. Puisque $\gA$ est \lgbz,
il y a une valeur $tf(\ux)+(1-t)a$  de $g$ qui est inversible. La valeur
$f(\ux)$ est donc inversible modulo $\fa$.

\emph{\ref{i6factlgb1}.}
Soit $P\in\AXn$ primitif par valeurs: $1 = uP(\ux) + vP(\uy) +\dots$. \\ 
En considérant $u$, $\ux$,
$v$, $\uy$, \dots\, et les \coes de $P$, on voit qu'il y a un sous-anneau
$\gA_i$ tel que $P \in \gA_i[\uX]$ et tel que $P$ soit primitif par valeurs
sur~$\gA_i$.  Ainsi, $P$ représente un \iv sur $\gA_i$, a fortiori sur $\gA$.
\end{proof}
%-----------------end proof------------------

Pour un \pol les \prts de représenter un inversible ou d'être \ppv sont \carfz, comme indiqué dans le lemme qui suit.

%:    Lemma{lemAlgb0}
\begin{lemma}\label{lemAlgb0}
Soit $S$ un \mo de $\gA$ et un \pol $f\in\gA[X_1, \ldots, X_m]$.
\begin{enumerate}
\item  Le \pol $f$ représente un \iv dans $\gA_S$ \ssi il existe $s\in S$
tel que $f$ représente un \iv dans~$\gA_s$.
\item  Le \pol $f$ est \ppv dans $\gA_S$ \ssi il existe $s\in S$
tel que $f$ est \ppv dans~$\gA_s$.
\end{enumerate}
\end{lemma}
\begin{proof}
Nous montrons seulement le point \emph{1.} Soit $F(\uX,T) \in \gA[\uX,T]$
l'homogénéisé de $f(\uX)$ en degré assez grand.
L'hypothèse équivaut à
l'existence de $\ux \in \Ae m$ et $t$, $u \in S$ tels que $F(\ux,t)$ divise $u$ dans
$\gA$. En posant $s = tu$,  les \elts $t$ et $F(\ux,t)$ sont \ivs dans
$\gA_s$ donc $f$ représente un \iv dans $\gA_s$.
\end{proof}
%

%:     Lemma{lemAlgb1}
\begin{lemma}\label{lemAlgb1}
Soit $s \in \gA$ et $\fb$ un idéal de $\gA$ avec $1 \in \gen {s} + \fb$.  
\begin{enumerate}
\item Supposons que $f$ représente un \iv dans $\gA_s$. Il existe $\uz \in \Ae m$
tel que $1 \in \gen {f(\uz)} + \fb$.
\item Si $f$ est \ppv dans $\gA_s$ il existe un nombre fini d'\eltsz~$\uz_j$, ($j\in\lrbk$), dans~$\Ae m$
tels que  $1 \in \gen {f(\uz_j)\mid j\in\lrbk} + \fb$.
\end{enumerate}
\end{lemma}
\begin{proof}
\emph{1.} Soit $F(\uX,T) \in \gA[\uX,T]$ l'homogénéisé de~$f(\uX)$ en degré $d$
assez grand. L'hypothèse est que $F(\ux, t)$ divise $u$ dans~$\gA$ pour un $\ux \in
\Ae m$ et $t, u \in s^\NN$. Il existe $a$ tel que $ta \equiv 1 \bmod\fb$
donc:
$$
\preskip.4em \postskip.4em 
a^d F(\ux, t) = F(a\ux, at) \equiv F(a\ux, 1) = f(a\ux) \bmod \fb, 
$$
d'où  $a^du \in \gen {f(\uz)} + \fb$ avec $\uz = a\ux$.
Mais $1 \in \gen {a^du} + \fb$ donc $1 \in \gen {f(\uz)} + \fb$.

 On peut présenter le même argument \gui{sans calcul} comme suit.
\\
On a
$\gA_s\sur{(\fb\gA_s)} \simeq (\gA\sur\fb)_s$. Puisque $1 \in \gen {s} +
\fb$, $s$ est \iv dans $\gA\sur\fb$, et donc $\gA_s\sur{(\fb\gA_s)} \simeq
\gA\sur\fb$. Puisque $f$ représente un \iv dans $\gA_s$, a fortiori
il représente un \iv dans $\gA_s\sur{(\fb\gA_s)} \simeq \gA\sur\fb$, i.e. $f$
représente un \iv modulo $\fb$. 

\emph{2.} Calculs similaires.
\end{proof}

Nous allons utiliser dans la suite un \plgc un peu subtil que
nous énonçons sous forme d'un lemme. Voir aussi l'exercice \ref{exopropAlgb1}.

%:     Lemma{lempropAlgb1}
\begin{lemma}\label{lempropAlgb1}
Soient $S_1$, $\dots$, $S_n$ des \moco de $\gA$ et un \pol $f\in\gA[X_1, \ldots, X_m]$. \Propeq
\begin{enumerate}
\item  Le \pol $f$ est \ppvz.
\item  Dans chacun des anneaux $\gA_{S_i}$,  le \pol $f$ est \ppvz.
\item [3.]$\!\!\!\!$\eto Pour tout \idema $\fm$  de $\gA$, $f$ représente
un \iv dans~$\gA\sur\fm$.
\end{enumerate}
En particulier,   si $f$ représente un \iv dans chaque localisé $\gA_{S_i}$,  $f$ est primitif par valeurs.
\end{lemma}

\begin {proof}
Les implications \emph{1} $\Rightarrow$ \emph{2}  $\Rightarrow$ \emph{3\eto}$\!$
sont \imdesz. L'implication  \emph{3\eto}  $\!\!\Rightarrow$ \emph{1} est facile
en \clamaz.

  Voici maintenant une \dem directe et \cov de
\emph{2} $\Rightarrow$ \emph{1.} Il s'agit du décryptage de la
\dem classique de \emph{3.\eto}  $\!\!\Rightarrow$ \emph{1}, 
en utilisant la méthode qui sera
expliquée dans la section~\ref{subsecLGIdeMax}.
Pour simplifier les notations mais \spdgz, nous allons montrer le cas particulier
où $f$ représente un \iv dans chaque localisé $\gA_{S_i}$.
\\
On dispose donc d'\eco $(s_1, \ldots, s_n)$ tels que dans chaque
localisé~$\gA_{s_i}$, le \polz~$f$ représente un \iv (lemme \ref{lemAlgb0}).
\\
En appliquant le lemme \ref{lemAlgb1} on obtient successivement, pour $k = 0, \ldots, n$, \[1 \in \gen {f(\und
{z_1}), \ldots, f(\und {z_k}),\ s_{k+1}, \ldots, s_n}.\]  Au bout de $n$
étapes: $1 \in \gen {f(\und {z_1}), \ldots, f(\und {z_n})}$.
\end {proof}

%:     Proposition{propAlgb1}
\begin{proposition}\label{propAlgb1}
 \Propeq
\begin{enumerate}
\item L'anneau $\gA$ est \lgbz.
\item Pour   tout \pol $f\in\AXn$, s'il existe un \sys d'\eco $(s_1,\ldots,s_k)$
tel que $f$ représente un \iv dans chaque~$\gA_{s_i}$, alors $f$ représente un \ivz.
\item Pour  tout \pol $f\in\AXn$, s'il existe des \moco $S_i$ tel que
 $f$  est \ppv dans chaque $\gA_{S_i}$, alors~$f$ représente un \ivz.
\end{enumerate}
\end{proposition}

\begin{proof}
Vus les lemmes \ref{lemAlgb0} et \ref{lempropAlgb1}, 
il suffit de montrer que si $f$ est \ppv il existe des
\eco tels que $f$ représente un \iv dans chaque localisé.
On écrit la chose en une variable pour simplifier les notations.
On obtient $x_1$, \dots, $x_r \in\gA$ tels que $1\in\gen{f(x_1), \ldots, f(x_r)}$.  Soit~$s_i =
f(x_i)$: le \pol $f$ représente un \iv dans $\gA_{s_i}$.
\end{proof}

\rdb
D'après le lemme de Gauss-Joyal (\ref{lemGaussJoyal}) les \pols primitifs
forment un filtre $U\subseteq\AX$.
On appelle \ix{localisé de Nagata} l'anneau
 $\gA(X)=U^{-1}\gA[X]$.

%:     Fact{factLocNagata}
\begin{fact}\label{factLocNagata}
On utilise la notation ci-dessus.
\begin{enumerate}
  \item $\gA(X)$ est fidèlement plat sur $\gA$.
  \item $\gA(X)$ est un  \algbz.
\end{enumerate}
\end{fact}
\begin{proof} 
\emph{1.} Il est clair que $\gA(X)$ est plat sur $\gA$ (\lon de $\AX$). On utilise alors la \carn \emph{3a}
dans le \thref{thExtFidPlat}. \hbox{Soit $\fa=\gen{\an}$} un \itf de~$\gA$ tel \hbox{que
$1\in\fa\gA(X)$}. Nous devons montrer \hbox{que $1\in\fa$}. L'hypothèse donne~\hbox{$f_1$, \ldots, $f_n\in\AX$} tels que le  \pol $f=\sum_ia_if_i$ est  primitif, i.e., $1\in\rc_\gA(f)$. Or l'\id $\rc_\gA(f)$ est contenu dans $\fa$.

\emph{2.} Nous procédons en trois étapes.\\
a) Montrons d'abord que tout \pol primitif  $P(T)\in\gB[T]$
où $\gB:=\gA(X)$
représente un \elt \ivz.
En effet, soit $P(T) = \sum_i  Q_iT^i$ un tel \polz,
on peut supposer \spdg que les $Q_i$ sont dans $\gA[X]$. On a des \pols $B_i$
tels que $\sum_i  B_i(X)Q_i(X)$ est primitif. A fortiori les \coes des $Q_j$ sont \comz. Alors, pour
$k > \sup_i\big(\deg_X(Q_i)\big)$, puisque $P(X^k)$ a pour \coes tous les \coes des
$Q_j$ (astuce de \KRAz), c'est un
\pol primitif de $\gA[X]$, \cad un \elt inversible de $\gB$.\\
b) Montrons la même \prt pour un nombre quelconque de variables.
On considère un \pol primitif $Q(\Ym)\in\BuY$. Par l'astuce de Kronecker
en posant $Y_j=T^{n^j}$ avec $n$ assez grand, on obtient
un \pol $P(T)$ dont les \coes sont ceux de $Q$,
ce qui nous ramène au cas précédent.\\
c) Enfin considérons un \pol $Q$ à $m$ variables sur $\gB$
primitif par valeurs. Alors,  $Q$ est primitif et
on peut appliquer le point~b).
\end{proof}

%: --- Subsubsec*{Propriétés \lgbes remarquables}
\subsec{Propriétés \lgbes remarquables}
%-------------------

%:    plcc{thlgb1}----------------
\begin{plcc}
\label{thlgb1} Soient $S_1$, $\ldots $, $S_r$ des \moco d'un \algb $\gA$.
%-----------------begin enum------------------
\begin{enumerate}
\item Si deux matrices de $\Ae {m\times n}$ sont \eqves sur chacun des
$\gA_{S_i}$, alors elles sont \eqvesz.
\item Si deux matrices de $\Mn(\gA)$ sont semblables sur chacun des
$\gA_{S_i}$, alors elles sont semblables.
\end{enumerate}
%-----------------end enum------------------
\end{plcc}
%--- end-plcc-----------------------------------------
%-----------------begin proof------------------
\begin{proof}
\emph{1.} Soient $F$ et  $G$  les matrices, par hypothèse il existe un
\sys d'\eco  $(s_1, \ldots , s_r)$  et des matrices
$U_1$, \ldots, $U_r$, $V_1$, \ldots, $V_r $ telles que pour chaque~$i$ on a
     $U_i F = G V_i$ et
     $\det(U_i) \det(V_i) = s_i$.
%:HHH un peu plus propre, il y avait confusion entre x et x: maintenant \ux et \ual 
Introduisons des \idtrs  $(x_1, \ldots , x_r)=(\ux)$,
et considérons les matrices   

\snic{U=U(\ux) = x_1\,U_1+\cdots +x_r\,U_r$ et $V=V(\ux) = x_1\,V_1+\cdots
+x_r\,V_r.
}

%\sni

    On a  $U F = G V,$
et  $\det(U) \det(V)$ est un \pol en les $x_i$ qui vérifie les
hypothèses de la \dfn \ref{defalgb}:  il suffit d'évaluer  $(x_1, \ldots , x_r)$
successivement en  $(1,0,\ldots ,0)$,
\ldots, $(0,\ldots ,0,1)$.
Donc il existe un $\ual\in\Ae r$ tel que l'\eltz~$\det\big(U(\ual)\big) \det\big(V(\ual)\big)$ est \ivz.

\emph{2.} La même \demz, avec $U_i=V_i$ et $U=V$, fonctionne.
\end{proof}
%-----------------end proof------------------

On a le corolaire suivant.

%:    plcc{thlgb2}-----
\begin{plcc}
\label{thlgb2} Soient $S_1$, $\ldots$, $S_r$ des \moco d'un \algb $\gA$.
%-----------------begin enum------------------
\begin{enumerate}
\item Si deux \mpfs sont isomorphes sur chacun des~$\gA_{S_i}$, alors ils sont isomorphes.
\item Tout \mptf est quasi libre.
\end{enumerate}
%-----------------end enum------------------
\end{plcc}
%--- end-plcc-------------------------
%-----------------begin proof------------------
\begin{proof}
\emph{1.} On considère des \mpns et l'on caractérise le
fait que les modules sont isomorphes par l'\eqvc de matrices associées
(lemme~\ref{lem pres equiv}).
On applique alors le point~\emph{1} du \plgrf{thlgb1}.

\emph{2.} On applique le point \emph{1}
On considère un module quasi libre qui a les mêmes \idfsz,
on sait que les deux modules deviennent libres
après \lon en des \eco (et le rang est chaque fois le même
parce qu'ils ont les mêmes \idfsz).
\end{proof}
%-----------------end proof------------------

Signalons aussi les principes suivants.

%:    plcc{thlgb3}----------
\begin{plcc}
\label{thlgb3} Soit un \algb $\gA$.
%-----------------begin enum------------------
\begin{enumerate}
\item Soient $S_1$, $\ldots$, $S_r$ des \mocoz, $M$ un \mpf et $N$
un \mtfz. Si $N$ est un quotient de $M$ sur chacun des~$\gA_{S_i}$, alors $N$ est un quotient de $M$.
\item Un module localement engendré par $m$ \elts est  engendré par $m$
\eltsz.
\end{enumerate}
%-----------------end enum------------------
\end{plcc}
%--- end-plcc-----------------------------------------

%
\begin{proof}
Il suffit de montrer le point \emph{1} car un module est engendré par $m$
\elts \ssi c'est un quotient d'un module libre de rang $m$.
\\
Nous continuerons la \dem après  les deux  lemmes suivants.
\end{proof}

%:     Lemma{lem1thlgb3}
\begin{lemma}\label{lem1thlgb3}
Soit $M$ un \Amo \pfz, $N$ un \Amo \tfz, $S$ un \mo de $\gA$ et $\varphi:M_S\to N_S$ une \Ali surjective. 
\begin{enumerate}
\item Il existe $s\in S$ et  $\psi\in\Lin_\gA(M,N)$ tels que $s\varphi=_{\gA_S}\psi_S$.
\item Il existe $v\in S$ tel que $vN\subseteq\psi(M)$.
\item Il existe une matrice $Q$ de \syzys satisfaites par les \gtrs de~$N$ telle que, en considérant le module $N'$ admettant~$Q$ 
pour \mpnz, l'application~$\psi$  se décompose comme suit

\snic{M\vers{\theta}N'\vers{\pi}N,}

%\sni
($\pi$ est la \prn canonique), avec en outre $vN'\subseteq\theta(M)$
(a fortiori~$\theta_S$ est surjective). 
\end{enumerate}
\end{lemma}
%--------- fin lemma ---------------------------------------------- 
%
\begin{proof}
Le point  \emph{1} est une reformulation de la proposition \ref{fact.homom loc pf} (qui affirme un tout petit peu plus, dans un cas plus \gnlz). Le point \emph{2} en découle facilement.
\\
\emph{3.} On a $N=\gA y_1+\cdots+\gA
y_n$, et $M=\gA x_1+\cdots+\gA x_m$, avec une \mpn $P$.  
\\
Pour que la \fcn par $\theta$ existe, il suffit que parmi les colonnes de la matrice $Q$ on trouve
les \syzys \gui{images des colonnes de $P$ par $\psi$} (ce sont des \syzys
entre les $y_k$ une fois que l'on a exprimé les $\psi(x_j)$ en fonction des $y_k$). \\
Pour que  $vN'\subseteq\theta(M)$, il suffit que  parmi les colonnes de la matrice $Q$ on trouve les \syzys exprimant que les $vy_k$ sont dans $\gA\psi(x_1)+\cdots+\gA\psi(x_m)$ (une fois que l'on a exprimé les $\psi(x_j)$ en fonction des $y_k$).
\end{proof}
%

%:     Lemma{lempfthlgb3}
\begin{lemma}\label{lempfthlgb3}
Le \plgc \ref{thlgb3} est correct si $N$ est lui-même un \mpfz. 
\end{lemma}
%--------- fin lemma ---------------------------------------------- 
%
\begin{proof} L'hypothèse donne une \ali surjective $\varphi_i:M_{S_i}\to N_{S_i}$. Par les points \emph{1} et \emph{2} du lemme \ref{lem1thlgb3} on a $s_i,v_i\in S_i$ et une \ali $\psi_i:M\to N$ tels que $s_i\varphi_i=(\psi_i)_{S_i}$ et $v_iN\subseteq\psi_i(M)$.
Chaque \ali $\psi_i$ est représentée
 par deux matrices $K_i$ et $G_i$ qui font
commuter les diagrammes convenables (voir la section \ref{secCatMpf}).

\smallskip 
\centerline{\xymatrix {
\gA^{p} \ar[r]^{P} \ar[d]_{K_i} & \Ae m \ar[d]^{G_{i}}
\ar@{->>}[r]^{\pi_M}
                   & M \ar[d]^{\psi_i} \\
\gA^{q} \ar[r]_{Q}                    & \Ae n \ar@{->>}[r]_{\pi_{N}}
                   & N \\}
}

On considère $r$ inconnues  $a_i$  dans $\gA$
et l'application  $\psi = \sum a_i \psi_i$
correspondant aux matrices  $K=\sum a_i  K_i$ et   $G=\sum a_i  G_i$.
$$\preskip.2em \postskip.2em 
{\xymatrix {
\gA^{p} \ar[r]^{P} \ar[d]_{K} & \Ae m \ar[d]^{G}
\ar@{->>}[r]^{\pi_M}
                   & M \ar[d]^{\psi} \\
\gA^{q} \ar[r]_{Q}                    & \Ae n \ar@{->>}[r]_{\pi_{N}}
                   & N \\}
}
$$

Le fait que  $\psi$  soit surjective signifie que la matrice
\smashtop{$H=\blocs{.7}{.8}{.6}{0}{$G$}{$Q$}{}{}$} est surjective, \cad que
$\cD_n(H)=\gen{1}$.  
On introduit donc les \idtrs $c_\ell$ pour fabriquer une
\coli arbitraire des mineurs maximaux $\delta_\ell$ de la matrice $H$. 
Cette \coli $\sum_\ell  c_\ell\delta_\ell$ est un \pol en les $a_i$ et
$c_\ell$. Par hypothèse, ce \pol représente $1$ sur chacun des
$\gA\big[\frac 1{s_iv_i}\big]$, donc, puisque l'anneau
est \lgbz, il représente un  \iv (proposition~\ref{propAlgb1}).
\end{proof}
%

%%%%%%%%%%%%%%%%%%%%%%%%%%%%%%%%%%%%%%%%%%%%%%%%%%%%%%%%%%%%%%%%%%%%%%%%
%-----------------begin Proof------------------
%\begin{Proof}{Fin de la \dem du \plgc \ref{thlgb3}.}
\emph{Fin de la \dem du \plgc \ref{thlgb3}.}
\\ 
 On a $N=\gA y_1+\cdots+\gA y_n$, et $M=\gA x_1+\cdots+\gA x_m$, avec une \mpn $P$.  
Pour chaque $i\in\lrbr$ on applique le lemme \ref{lem1thlgb3}
avec le \mo $S_i$ et l'\ali surjective $\varphi_i:M_{S_i}\to N_{S_i}$ donnée 
dans l'hypothèse. On 
obtient une \ali $\psi_i:M\to N$, une matrice~$Q_i$ de \syzys satisfaites par les $y_k$,
une \aliz~$\theta_i:M\to N'_i$ (où~$N'_i$ est le \mpf correspondant à~$Q_i$),
des \elts $s_i$, $v_i\in S_i$ avec $s_i \varphi_i=(\psi_i)_{S_i}$, enfin~$\psi_i$  se factorise via~$\theta_i:M\to N'_i$ avec~$v_iN'_i\subseteq\theta_i(M)$.  
\\
On considère alors le module $N'$ \pf correspondant à la matrice de relations $Q$ obtenue en juxtaposant les matrices~$Q_i$, de sorte que~$N'$ est un quotient
de chaque~$N'_i$.
\\
Comme $N$ est un quotient de $N'$, on a ramené le \pb au cas où $N$
est lui-même \pfz, cas qui a été traité dans le lemme~\ref{lempfthlgb3}.
\eop
%\end{Proof}
%-----------------end Proof------------------

%:HHH 11 petit changement dans les noms section sous section

%: --- Subsec{Systèmes congruentiels}
\subsec{Systèmes congruentiels}
\label{subsecSysCong}
%-----------------------------------------

Une  \prt de stabilité importante des \algbs
est la stabilité par extension entière.
%C'est sans doute cette \prt qui légitime vraiment le fait de
%prendre le temps d'étudier en détail ces anneaux.

%:     Theorem{thLgbExtEnt}
\begin{theorem}\label{thLgbExtEnt}
Soit $\gA\subseteq\gB$ avec $\gB$ entier sur $\gA$.
Si $\gA$ est \lgbz, alors  $\gB$ \egmtz. \perso{\ssi?}
\end{theorem}

La \dem est renvoyée \paref{propSysCong}, après un détour par les anneaux congruentiels.

%:     Definition{defiSysCong}
\begin{definition}\label{defiSysCong}
Une partie $C$ d'un anneau $\gA$ est appelée un
\emph{système congru\-entiel}
si elle vérifie la \prt suivante: si $s_1+s_2=1$
dans $\gA$ et si~$c_1$, $c_2\in C$, alors il existe $c\in C$ tel que
$c\equiv c_1\mod s_1$ et $c\equiv c_2\mod s_2$.
\index{systeme cong@\sys congruentiel}
\index{congruentiel!système ---}
\end{definition}

\rems 1) Il revient au même de dire: si $\fa_1$ et $\fa_2$
sont deux \ids \com de $\gA$ et si $c_1$, $c_2\in C$, alors il existe $c\in C$ tel que
$c\equiv c_1\mod \fa_1$ et $c\equiv c_2\mod \fa_2$.\\
2) L'\elt $c'=c_2s_1+c_1s_2$ est le  candidat naturel  pour $c\in\gA$ vérifiant 
les congruences $c\equiv c_1\mod s_1$ et $c\equiv c_2\mod s_2$. 
On doit donc avoir un \eltz~$c$ de~$C$ tel que $c\equiv c'\mod s_1s_2$.
\eoe

\medskip \rdb 
\exl \label{defiSuslinSet}
Soit $(\ub) = (\bn)$ une suite dans un anneau $\gB$. L'\ix{ensemble de
Suslin de $(\bn)$} est la partie suivante de $\gB$:\index{Suslin!ensemble de --- d'une suite finie}

\snic {
\Suslin(\ub) = \sotq{u_1b_1 + \cdots + u_nb_n} {
(\un)\hbox{ est $\EE_n(\gB)$-complétable}},
} 

%\sni
($(\un)$ est la première ligne d'une matrice de $\EE_n(\gB)$).

Si l'un des $u_i$ est \ivz, alors $u_1b_1 + u_2 b_2
+ \cdots + u_n b_n \in \Suslin(\ub)$ et l'on a donc $\{\bn\} \subseteq
\Suslin(\bn) \subseteq \gen {\bn}$.

Montrons que l'ensemble $\Suslin(\ub)$ est toujours  congruentiel. 
\\
En effet, pour  $E$, $F \in \EE_n(\gB)$ et deux \elts comaximaux $s$, $t$
de $\gB$, il existe $G \in \EE_n(\gB)$ vérifiant $G \equiv E \bmod s$ et
$G \equiv F \bmod t$.

 Soient $f$, $g_1$, \ldots, $g_n  \in \AX$ avec $f$ \monz, et $\gB=\aqo\AX f$. 
Alors l'ensemble de Suslin de $(\ov {g_1},\dots,\ov {g_n})$ joue un
rôle important dans l'étude des  \vmds \polls (cf. lemme \ref
{lemSuslin1}).
\eoe

%:     Fact{fact1SysCong}
\begin{fact}\label{fact1SysCong}
Pour tout \pol $P\in\AXn$ l'ensemble $V_P$ des valeurs de $P$
est un \sys congruentiel  ($V_P=\sotq{P(\ux)}{\ux\in\Ae n}$).
\end{fact}
\begin{proof}
Soient $s$, $t$ deux \eco et $\ux$, $\uy$  dans
$\Ae n$. Un chinois nous donne un $\uz\in\Ae n$ tel que $\uz\equiv \ux\mod s$ et  $\uz\equiv \uy\mod t$. Alors, on~a~$P(\uz)\equiv P(\ux)\mod s$ et  $P(\uz)\equiv P(\uy)\mod t$.
\end{proof}
%

%:     Fact{fact2SysCong}
\begin{fact}\label{fact2SysCong}
Soit $C$  un \sys congruentiel. Si $\fa_1$, \ldots, $\fa_\ell$
sont des \ids deux à deux \com et si $c_1$, \ldots, $c_\ell\in C$,
alors il existe $c\in C$ tel que~$c\equiv c_j\mod \fa_j$ pour $j\in\lrb{1..\ell}$.
\end{fact}
\begin{proof}
Il s'agit de la preuve usuelle du \tho chinois, adaptée à la situation présente.
On raisonne par \recu sur $\ell\geq2$. L'initialisation est par \dfnz.
Si $\ell>2$ on considère les \ids deux à deux \com $\fa_1,\ldots,\fa_{\ell-2}$
et $\fa_{\ell-1}\fa_{\ell}$.
Soit $e\in C$ tel que $e\equiv c_{\ell-1}\mod\fa_{\ell-1}$
et  $e\equiv c_\ell\mod\fa_\ell$. \\
Par \hdrz, on trouve $c$ dans $C$ tel que
$c\equiv c_k\mod\fa_k$ pour~$k\in\lrb{1..\ell-2}$ et~$c\equiv e\mod\fa_{\ell-1}\fa_{\ell}$. A fortiori, $c\equiv c_{\ell-1}\mod\fa_{\ell-1}$
et~$c\equiv c_\ell\mod\fa_\ell$.
\end{proof}
%

%:     Fact{fact3SysCong}
\begin{fact}\label{fact3SysCong}
Soient $C$  un \sys congruentiel, $w_1$, \ldots, $w_n$ des \elts de~$C$
et~$(e_1,\ldots,e_n)$  un \sfioz.
Alors, l'\elt $w=e_1w_1+\cdots+e_nw_n$ est dans~$C$.
\end{fact}
\begin{proof}
On a $w\equiv w_i \mod 1-e_i$, et les $\gen{1-e_i}$
sont 2 à 2 \comz, mais $w$ est l'unique \elt vérifiant ces  
 congruences 
puisque $\bigcap_i\gen{1-e_i}=\gen{0}$. Il reste à appliquer le fait précédent.
\end{proof}
%
%:     Definition{defiAnneauCongruentiel}
\begin{definition}\label{defiAnneauCongruentiel}
Un anneau $\gA$ est dit \ixc{congruentiel}{anneau ---} si tout
\sys congruentiel  qui engendre l'\id $\gen{1}$ contient un \elt\ivz.\index{anneau!congruentiel}
%-% PERSO
\perso{Dans la \dfn \ref{defiAnneauCongruentiel}, la quantification ``pour tout \sys congruentiel''
n'est pas légitime, il faudrait se restreindre à un ``ensemble'' de \syss congruentiels
}
%-% Fin PERSO
\end{definition}
%:     Lemma{lemAnneauCongruentiel}
\begin{lemma}\label{lemAnneauCongruentiel}~
\begin{enumerate}
\item \label{i3lemAnneauCongruentiel} Soit $\fa\subseteq\Rad\gA$. Alors, l'anneau $\gA$ est congruentiel \ssi l'anneau $\gA\sur{\fa}$ est congruentiel.
\item \label{i1lemAnneauCongruentiel} Tout anneau \plc est congruentiel.
\item  \label{i2lemAnneauCongruentiel} Tout anneau congruentiel est \lgbz.
\end{enumerate}
 \end{lemma}
\begin{proof}
\emph{\ref{i3lemAnneauCongruentiel}.}
On utilise le fait que des \elts sont \com (resp. \ivsz)
dans $\gA$ \ssi ils sont \com (resp. \ivsz) dans~$\gA\sur{\fa}$.
 
\emph{\ref{i1lemAnneauCongruentiel}.}
Supposons $\gA$ \plcz. Il suffit de montrer que~$\gA\sur{\Rad\gA}$ est congruentiel.
 Soit $W$ un \sys congruentiel de $\gA\sur{\Rad\gA}$ tel que~$\gen{W}=\gen{1}$.
Soient $w_1$, \ldots, $w_n\in W$ avec $\gen{w_1,\ldots,w_n}=\gen{1}$.
Il existe  un \sfio  $(e_1,\ldots,e_n)$ tel que l'on ait~$\gen{e_1w_1+\cdots+e_nw_n}=\gen{1}$
(lemme \ref{lemZerRed}  point~\emph{5}).
On conclut avec le fait~\ref{fact3SysCong} que~$W$ contient
l'\elt\iv $e_1w_1+\cdots+e_nw_n$.
 
\emph{\ref{i2lemAnneauCongruentiel}.}
Supposons $\gA$ congruentiel et soit $P$ un \pol primitif par valeurs. Puisque les valeurs de $P$ forment un \sys congruentiel,
une valeur de $P$ est \ivz.
\end{proof}
%

%% --- Subsubsec*{Retour à l'extension entière d'un \alo \dcd}
%\subsubsection*{Retour à l'extension entière d'un \alo \dcd}
%%-----------------------------------------

%:HHH 11 suite de petit changement dans les noms section/ sous section
%: --- Subsec{Stabilité par extension entière}
\subsec{Stabilité par extension entière}
%-------------------

%:HHH modifications ci après
Comme corolaire \imd du lemme \ref{lemAnneauCongruentiel} on a le résultat qui suit.
%:     Proposition{CorSysCong}\label{corpropSysCong}
\begin{corollary}\label{CorSysCong}
Soit $\gB$ une \alg \stfe sur un \cdiz~$\gA$
et $W$ un \sys congruentiel dans $\gB$ tel que~$\gen{W}=\gen{1}_\gB$.
\\
Alors,  l'ensemble $\rN\iBA(W)$ contient un \elt \ivz.
\end{corollary}
\begin{proof}
On sait que $\gB$ est \zedz, donc il est congruentiel 
(lemme~\ref{lemAnneauCongruentiel}).
Puisque $W$ est congruentiel et engendre l'\id $\gen{1}$, il contient un \elt \ivz. Enfin, la norme d'un \elt \iv est \ivz.
\end{proof}
%

%:HHH modifié
%:     proposition{propSysCong}
\begin{proposition}\label{propSysCong}
Soit $\gB$ une \alg \stfe sur un anneau~$\gA$
et $W$ un \sys congruentiel dans $\gB$. Si~$1\in\gen{W}$,
alors,   \hbox{$1\in\gen{\rN\iBA(W)}$}.
\end{proposition}
\begin{proof} \emph{1.} Un \sys congruentiel reste congruentiel par passage à un anneau quotient. 
Si on lit la conclusion du corolaire \ref{CorSysCong} 
sous la forme~\hbox{$1\in\gen{\rN\iBA(W)}$} (plus faible), on constate
qu'il est sous une forme adéquate pour subir
la machinerie  \cov à \idemas qui sera expliquée \paref{MethodeIdemax}
dans la section~\ref{subsecLGIdeMax}, et qui sert à démontrer qu'un \id contient $1$. On obtient donc
le résultat souhaité.\imlma
\end{proof}

\rems \\
1) En \clama on dirait ceci: si 
 $1\notin \gen{\rN\iBA(W)}_\gA$, cet \id serait contenu dans un \idema $\fm$  de $\gA$.
Mais le corolaire~\ref{CorSysCong}, appliqué avec le \cdi $\gA/\fm$
et l'\alg \stfe $\gB/\fm\gB$, montre que c'est impossible. \\
La machinerie
\cov  à \idemas a justement pour but de décrypter ce type de preuve abstraite et de  la transformer en un
\algo qui construit $1$ comme \elt de~$\gen{\rN\iBA(W)}_\gA$ à partir des hypothèses.\imlma

2) Comme exemple, si $(\ub) = (b_1, \ldots, b_q)$  est un \sys d'\eco
dans $\gB$, on a $1 \in \gen{\rN_{\gB/\gA}(w) \mid w \in \Suslin(\ub)}_\gA$,
puisque l'ensemble $\Suslin(\ub)$ est congruentiel. \\
Mais on se gardera de croire  \hbox{que $1 \in \gen {\rN_{\gB/\gA}(b_1), \ldots,
\rN_{\gB/\gA}(b_q)}_\gA$}. 
\\
Une instance célèbre de cette \prt   est
un résultat d\^u à Suslin portant sur les vecteurs \pollsz, donné dans le lemme
\ref {lemSuslin1}.  Dans ce lemme, $\gB$ est de la forme $\aqo{\gA[X]}{v}$
avec $v \in \gA[X]$ \poluz. Un décryptage complet sera fourni dans la
\dem du lemme en question.
\eoe

\begin{Proof}{\Demo du \thrf{thLgbExtEnt}. }
Traitons en  premier le cas où~$\gB$ est  libre de rang fini,
disons $\ell$, sur $\gA$.  
Soit $P\in\BXn$ un \pol primitif par valeurs.  Nous voulons  un
$\ub\in\gB^n$ avec $P(\ub)$ \ivz.  Nous considérons le \sys congruentiel $W$ des valeurs de
$P$. Par hypothèse \hbox{on a $1\in\gen{W}$}. La proposition~\ref{propSysCong} nous
dit alors que~$\gen{\rN\iBA(W)}_\gA=\gen{1}_\gA$.\\
Mais $\rN\iBA\big(P(b_1,\ldots,b_n)\big)$
est un \pol à $n\ell$ variables dans $\gA$ si l'on exprime chaque $b_i\in \gB$
sur une $\gA$-base de $\gB$, et $\gA$ est \lgbz, donc il existe
$\ub\in\gB^n$ tel que $\rN\iBA\big(P(\ub)\big)$ est \ivz, et cela implique
que $P(\ub)$ est \ivz.\\
Dans le cas \gnl où $\gB$ est seulement supposé entier sur $\gA$,
considérons dans $\gB$ les sous-\Algs $\gB_i$ \tfz; $\gB$ en est la
réunion filtrante croissante. Puisque $\gB$ est entier sur $\gA$, $\gB_i$
 \egmtz, donc est un quotient d'une \Alg qui est un \Amo libre de rang
fini.  D'après le premier cas, et en vertu du point
\emph{\ref {i5factlgb1}} du fait \ref{factlgb1}, chaque $\gB_i$ est \lgbz.
Enfin d'après le dernier point du fait \ref{factlgb1},
$\gB$ est \lgbz.
\end{Proof}
%

%%%%%%%%%%%%%%%%%%%%%%%%%%%%%%%%%%%%%%%%%%%%%%%%%%%%%%%%%%%%%%%%%%%%%%%%%
%%%%%%%%%%%%%%%%%%%%%%%%%%%%%%%%%%%%%%%%%%%%%%%%%%%%%%%%%%%%%%%%%%%%%%%%%
%%%%%%%%%%%%%%%%%%%%%%%%%%%%%%%%%%%%%%%%%%%%%%%%%%%%%%%%%%%%%%%%%%%%%%%%%
%%%%%%%%%%%%%%%%%%%%%%%%%%%%%%%%%%%%%%%%%%%%%%%%%%%%%%%%%%%%%%%%%%%%%%%%%
%%%%%%%%%%%%%%%%%%%%%%%%%%%%%%%%%%%%%%%%%%%%%%%%%%%%%%%%%%%%%%%%%%%%%%%%%

%:section: Exercices
\Exercices

%--- Exercise{exoNilRad}------
\begin{exercise}
\label{exoNilRad}
{\rm  Démontrer en \clama que le nilradical d'un anneau est égal à
l'intersection de ses \idepsz.
}
\end{exercise}
%--- end-exercise-----------------------------------------

%--- Exercise{exoRadJacSat}-------------
\begin{exercise}
\label{exoRadJacSat}
{\rm Si $\fa$ est un \id de $\gA$ on note $\JA(\fa)$ son
\emph{radical de Jacobson}, \cad l'image réciproque
de $\Rad(\gA\sur\fa)$ par la projection canonique $\gA\to\gA\sur\fa$.
 Soit~$\fa$ un \id de $\gA$. Montrer que $\JA(\fa)$ est le plus grand \id $\fb$
 tel que le \moz~$1+\fb$ soit contenu dans le saturé de $1+\fa$.
} 
\end{exercise}
%--- end-exercise-----------------------------------------

%--- Exercise{exoRadJacLoc}------
\begin{exercise}
\label{exoRadJacLoc}
{\rm  Démontrer en \coma que le radical de Jacobson d'un \alo 
coïncide avec l'ensemble des \elts noninversibles. 
Et que c'est l'unique \id $\fa$ vérifiant:
\begin{itemize}
\item $\fa$ est maximal
\item $1\in\fa$ implique $1=0$.
\end{itemize}
}
\end{exercise}
%--- end-exercise-----------------------------------------

%--- Exercise{exoRadNonCom}-------------
\begin{exercise}
\label{exoRadNonCom}
{\rm Soit $\gA$ un anneau non commutatif, $a,b\in\gA$.
 \emph{1.} Si $a$ admet un inverse à gauche $c$, alors $c$ est inverse à droite de $a$ \ssi $c$ est l'unique inverse à gauche.

 \emph{2.} Si $1-ab$ admet un inverse à gauche $u$, alors $1-ba$ admet aussi un inverse à gauche $v$. Idée: si $ab$ et $ba$ sont \gui{petits}, $u$ doit être égal à $1+ab+abab+\dots $,
et~$v$ égal à $1+ba+baba+\cdots =1+b(1+ab+abab+\cdots)a$. 

 \emph{3.}  Si pour tout $x$, $1-xa$ est \iv à gauche, alors pour tout $x$, $1-xa$ est \iv à droite. 

 \emph{4.}  \Propeq
\begin{itemize}
\item Pour tout $x$, $1-xa$ est \iv à gauche.
\item Pour tout $x$, $1-xa$ est \iv à droite. 
\item Pour tout $x$, $1-xa$ est \ivz.  
\item Pour tout $x$, $1-ax$ est \iv à gauche.
\item Pour tout $x$, $1-ax$ est \iv à droite. 
\item Pour tout $x$, $1-ax$ est \ivz.  
\item Pour tous $x,y$, $1-xay$ est \ivz.  
\end{itemize}
Les \elts $a$ qui vérifient ces \prts forment un \id bilatère
appelé radical de Jacobson de $\gA$.

}
\end{exercise}
%--- end -exercise-----------------------------------------

%--- Exercise{exoLocalFreenessLemma}-------------
\begin{exercise}\label{exoLocalFreenessLemma}
 {(Un lemme de liberté)}
{\rm
Soit $(\gA, \fm)$ un anneau local intègre de corps résiduel $\gk$, de
corps de fractions $\gK$. Soit $E$ un \Amo de type fini; on suppose
que le \kev $E/\fm E = \gk \otimes_\gA E$ et le \Kev
$\gK \otimes_\gA E$ ont même dimension $n$. Montrer que $E$ est un
$\gA$-module libre de rang $n$. \\
Mieux: si $(x_1, \ldots, x_n) \in E^n$ est une
base résiduelle, c'est une $\gA$-base de $E$.
}
\end {exercise}
%--- end -exercise-----------------------------------------

%--- Exercise{exoNakayamaAppli}------
\begin{exercise}\label{exoNakayamaAppli} {(Une application du lemme de Nakayama)}
\\
{\rm
Soit  $E$ un \Amo \pf et $a \in \Rad(\gA)$ un \elt $E$-régulier.
On suppose que le $\gA/a\gA$-module $E/aE$ est libre de rang $n$.
Montrer que $E$ est libre de rang $n$. Plus \prmtz, soient
$e_1$, \ldots, $e_n \in E$, si $(\overline {e_1}, \ldots, \overline {e_n})$
est une~$\gA/a\gA$-base de $E/aE$, alors $(e_1, \ldots, e_n)$
est une $\gA$-base de~$E$.
}
\end {exercise}
%--- end-exercise-----------------------------------------

%--- Exercise{exoItfmonogene}------
\begin{exercise}
\label{exoItfmonogene}
{\rm  Soit $\gA$  un \aloz. Si $\gen{b}=\gen{a}$,  il existe un \elt
inversible~$u$ tel que $ua=b$. Si $\fa=\gen{\xn}=\gen{a}$,  il
existe un indice $i$ tel que $\fa=\gen{x_i}$.
}
\end{exercise}
%--- end-exercise-----------------------------------------

%--- Exercise{exoZeroSimpleIntersComp}-------------
\begin{exercise}
\label{exoZeroSimpleIntersComp}
{\rm   \Demo directe détaillée du \thref{thJZS} lorsque $n=s$.  
}
\end{exercise}
%--- end -exercise-----------------------------------------

%--- Exercise{exopropZerdimLib}------
\begin{exercise}
\label{exopropZerdimLib}
{\rm  On reprend certains points du \thref{propZerdimLib}, en supposant maintenant
que l'anneau $\gA$ est  \plcz. \Llec est invité\e à fournir des \dems
indépendantes des résultats obtenus pour les \algbsz.

 \emph{1.} Tout \Amo \ptf
est quasi libre.

 \emph{2.} Toute matrice  $G\in\gA^{q\times m}$ de rang $\geq k$ est \eqve
à une matrice

\snic{
\cmatrix{
    \I_{k}   &0_{k,m-k}      \cr
    0_{q-k,k}&  G_1      },}

%\sni
avec $\cD_r(G_1)=\cD_{k+r}(G)$ pour tout $r\geq 0$.
Les matrices sont  \elrt \eqves si $k<\sup(q,m)$.

 \emph{3.} Tout \mpf \lot engendré par $k$ \elts est engendré par $k$ \eltsz.
}
\end{exercise}
%--- end-exercise-----------------------------------------

%%--- Exercise{exoMnMn+1}-------------
%\begin{exercise}
%\label{exoMnMn+1}\relax
%{(Le \Amo $\fp^p/\fp^{p+q}$ dans le cas d'un \idema $\fp$)}\\
%{\rm Soit $\fp$ un \idema strict détachable d'un anneau $\gA$.
%\perso{Exo \ref{exoMnMn+1}. Je n'ai pas bien réfléchi, mais il me semble que l'on a besoin de
%l'hypothèse $\fp$  détachable. \`A vérifier.
%Il semble que cet exercice mériterait d'être développé.
%Que se passe-t-il si au lieu de supposer $\fp$ maximal on suppose $\gA\sur{\fp}$
%\zedz?}
%Le complé\-mentaire de $\fp$ dans $\gA$ est donc un filtre premier $S$,
%qui est aussi le saturé du \mo $1+\fp$, et l'on note
%$\gA_{\fp}$ le localisé $S^{-1}\gA$. On note $\fm$ l'idéal de $\gA_{\fp}$
%engendré par $\fp$ (ou plus exactement par l'image de $\fp$ dans $\gA_{\fp}$).
%Soient $p,\,q$ des entiers $\geq 0$. Montrer que l'\homo canonique de $\gA$ dans
%$\gA_{\fp}$ envoie $\fp^p$ dans $\fm^p$ et que l'on en déduit un \iso de \Amos de
% $\fp^p/\fp^{p+q}$ sur $\fm^p/\fm^{p+q}$.
%}
%\end{exercise}
%%--- end-exercise-----------------------------------------

%--- Exercise{exoSLnEn}------------
\begin{exercise}\label{exoSLnEn} 
(Si $\gA$ est local, $\SLn(\gA)=\En(\gA)$)
\\
{\rm  Soit $\gA$  un \aloz. Montrer que toute matrice $B\in\SLn(\gA)$  est produit de matrices \elrs (autrement dit,~$B$ est \elrt
\eqve à la matrice~$\In$). On pourra s'inspirer de la preuve du lemme de
la liberté locale. Voir aussi l'exercice~\ref{exoLgb3}.
}
\end{exercise}
%--- end-exercise-----------------------------------------

%--- Exercise{exoNbgenloc}-----------
\begin{exercise}
\label{exoNbgenloc}
{\rm  \emph{1.}
Démontrer qu'un \Amo \tf $M$ est localement engendré par $k$ \elts
(\dfn \ref{deflocgenk})
\ssi $\Vi_\Ae {k+1}M=0$. On pourra s'inspirer du cas $k=1$  traité dans
le \thref{propmlm}.
 
\emph{2.} En déduire que l'annulateur $\Ann\big(\Vi_\Ae {k+1}M\big)$ et l'\idfz~$\cF_k(M)$ ont même radical.
}
\end{exercise}
%--- end-exercise-----------------------------------------

%--- Exercise{exoVariationLocGenerated}-------------
\begin{exercise}\label{exoVariationLocGenerated}
{(Variation sur le thème localement engendré)}\\
{\rm
Soit $M$ un \Amo \tfz, avec deux \sgrs $(\xn)$ \hbox{et $(\yr)$} avec $r \le n$.  On veut
expliciter une famille $(s_I)$ de ${n \choose r}$ \ecoz, indexée par les $I
\in \cP_{r,n}$, telle que $s_IM \subseteq \gen {(x_i)_{i\in I}}$. Notez que 
sur chaque localisé $\gA[s_I^{-1}]$, le module $M$ est engendré par les $(x_i)_{i\in I}$.

\emph {1.}
Soient $A$ et $B\in \Mn(\gA)$. 
\vspace{-3pt}
\begin {enumerate}\itemsep=0pt
\item [a.]
Expliciter l'appartenance:

\snic {
\det(A+B) \in \cD_{n-r}(B) + \cD_{r+1}(A).
}
\item [b.]
En déduire que $1 \in \cD_{n-r}(\In-A) + \cD_{r+1}(A)$. 
\item [c.]
En particulier, si $\rg(A) \le r$, alors $\rg(\In-A) \ge n-r$.
\item [d.]
Soient $a_1$, \dots,  $a_n \in \gA$, $\pi_I = \prod_I a_i$, $\pi'_J = \prod_J (1-a_j)$.\\
Montrer que les $(\pi_I)_{\#I=r+1}$ et $(\pi'_J)_{\#J=n-r}$ 
forment un \sys de ${n+1 \choose r+1}$ \elts \comz.
\end {enumerate}

\emph {2.}
Prouver le résultat annoncé en début d'exercice en
explicitant la famille $(s_I)$.

\emph {3.}
Soit $E$ un \Amo \tf localement engendré par~$r$ \eltsz.  Pour n'importe
quel \sgr $(\xn)$, il existe des \ecoz~$t_j$ tels que chacun des
localisés $E_{t_j}$ est engendré par $r$ \elts parmi les~$x_i$.

\emph {4.}
Soit $E = \gen {\xn}$ un \Amo \tf et $A \in \Mn(\gA)$ vérifiant
$\ux\,A = \ux$ avec $\rg(A) \le r$. Montrer que $E$
est localement engendré par~$r$ \eltsz.
\'{E}tudier une réciproque.
}

\end {exercise}
%--- end -exercise-----------------------------------------

%--- Exercise{exoMorAdcp}-------------
\begin{exercise}
\label{exoMorAdcp}
{\rm Si $\gA$ et $\gB$ sont deux anneaux \dcps on dit qu'un \homo
d'anneaux $\varphi:\gA\to\gB$ est un \emph{morphisme d'anneaux \dcpsz} 
si, pour tous $a$, $b\in\gA$ vérifiant  $b(1-ab)=0$ et $a(1-ab)\in\Rad\gA$,
on a dans~$\gB$, \hbox{avec $a'=\varphi(a)$} et $b'=\varphi(b)$,   $b'(1-a'b')=0$ et $a'(1-a'b')\in\Rad\gB$ (cf. proposition~\ref{prop1DecEltAnneau}).
\begin{itemize}
\item [\emph{1.}] Montrer que $\varphi$ est un morphisme  d'anneaux \dcps \ssiz
$\varphi(\Rad\gA)\subseteq\Rad\gB$.
\item [\emph{2.}]  \'Etudier les morphismes injectifs et surjectifs
 des anneaux \dcpsz. En d'autres termes, préciser les notions
 de sous-anneau \dcp (en un seul mot) et d'anneau \dcp quotient.
\end{itemize}
\index{morphisme!d'anneaux dec@d'anneaux \dcpsz} 
}
\end{exercise}
%--- end -exercise-----------------------------------------

%--- Exercise{exoMachDec}-------------
\begin{exercise}
\label{exoMachDec} (Machinerie locale-globale \elr des anneaux \dcpsz)
\\
{\rm
Le fait de pouvoir scinder systématiquement en deux composantes un anneau \dcp
conduit à la méthode \gnle suivante.
\\
 {\it La plupart des \algos qui fonctionnent avec les \alos \dcds peuvent être modifiés de manière à fonctionner avec les anneaux
\dcpsz, en scindant l'anneau
en deux composantes chaque fois que l'\algo écrit pour les \alos \dcds
utilise le test
\gui{cet \elt est-il \iv ou dans le radical?}.
Dans la première composante l'\elt en question
est \ivz, dans la seconde il est dans le radical.}
\\
 En fait on a rarement l'occasion d'utiliser cette machinerie \elrz,
la principale raison étant qu'une machinerie \lgbe plus \gnle
(mais moins \elrz) s'applique avec un anneau arbitraire, comme cela sera
expliqué
dans la section \ref{secMachLoGlo}.\imlb

}
\end{exercise}
%--- end -exercise-----------------------------------------

%:HHH exo modifié

%--- Exercise{exopropAlgb1}-------------
\begin{exercise}
 \label{exopropAlgb1} (\Pol représentant \lot un \ivz, lemme \ref{lempropAlgb1})\\
{\rm Le point \emph{3} de cet exercice donne une version renforcée du lemme \ref{lempropAlgb1}. L'approche utilisée ici est due à Lionel Ducos.
\\
Soient $\gA$ un anneau, $d \in \NN$ et $e = d(d+1)/2$.

\emph{1.}
Ici,  $s$ est une \idtr sur $\ZZ$.  Construire
$d+1$ \pols $a_i(s) \in \ZZ[s]$ \hbox{pour $i \in \lrb{0..d}$}, vérifiant pour tout
$P \in \AuX = \gA[\Xn]$ de degré $\le d$:
$$
s^e P(s^{-1}\uX) = a_0(s)P(s^0\uX) + a_1(s)P(s^1\uX) + \cdots + a_d(s)P(s^d\uX). 
\leqno(\star_d)
$$

\emph{2.}
Pour $s \in \gA$, $\ux \in \gA^n$ et $P\in\AuX$ de degré total $\leq d$, montrer que:

\snic{
s^eP(\ux/s) \in\gen{P(\ux),P(s\ux), \ldots,P(s^{d}\ux)} \subseteq \gA.
}

\emph{3.} 
Soit $S$ un \mo et $P\in\AuX$. On suppose que~$P$ représente un \iv dans 
$\gA_S$. Montrer que $S$ rencontre l'\id engendré
par les valeurs de $P$.

} 
\end {exercise}
%--- end-exercise-----------------------------------------

%--- Exercise{exoLgb2}-------------
\begin{exercise}
 \label{exoLgb2} 
 {\rm 
%:HHH  reference a un autre exo
 (Voir aussi l'exercice \ref{exoUAtoUB})
Soit $\gA$ un \algb et $M$ \hbox{un \Amoz}.
 
 \emph{1.} Pour tout \id $\fa$,  l'\homo canonique
$\Ati\to(\gA\sur{\fa})\eti$ est surjectif.
 
 \emph{2.} Si $x$, $y\in M$ et $\gA x=\gA y$, il existe un
\iv $u$ tel  
que~$x=uy.$

} \end{exercise}
%--- end-exercise-----------------------------------------

%--- Exercise{exoLgb3}-------------
\begin{exercise}\label{exoLgb3} 
(Si $\gA$ est \lgbz, $\SLn(\gA)=\En(\gA)$)\\
{\rm
Soit $\gA$ un \algbz, et $(a_1,\ldots,a_n)$ un \vmd  ($n\geq2$).
 
 \emph{1.} Montrer qu'il existe $x_2$, \ldots, $x_n$ tels que $a_1+\sum_{i\geq2}x_ia_i\in\Ati$.
 
 \emph{2.} En déduire (pour $n\geq2$) que tout \vmd se transforme en le vecteur $(1,0,\ldots,0)$ par manipulations \elrsz.
 
 \emph{3.} En déduire que $\SLn\gA=\En\gA$.
} \end{exercise}
%--- end-exercise-----------------------------------------

%--- Exercise{exo1semilocal}-------------
\begin{exercise}
\label{exo1semilocal} (Anneaux \slgbsz, 1)%
\index{anneau!semi@\slgbz}%
\index{semi@\slgbz!anneau ---}%
\index{anneau!semi@\smlz}%
\index{semi@\smlz!anneau ---}
\\
{\rm\emph{1.}
Pour un anneau $\gB$,  \propeq
\begin {itemize}
\item [~\emph{a.}]
Si $(x_1, \ldots, x_k)$  est \umdz, il existe un \sys d'\idms  \orts
$(e_1, \ldots, e_k)$  tel que $e_1 x_1+\cdots + e_k x_k$ soit \ivz.

\item [~\emph{b.}]
Sous la même hypothèse, il existe un scindage $\gB \simeq \gB_1 \times
\cdots \times \gB_k$ tel que la composante de $x_i$ dans $\gB_i$ soit
\iv pour $i \in \lrbk$.

\item [~\emph{c.}]
Même chose que dans \emph{a}, mais avec $k = 2$.

\item [~\emph{d.}]
Pour tout $x\in\gB$, il existe un \idm $e\in\gB$ tel que $x+e$ soit \ivz.
\end {itemize}
Notez qu'au point \emph{a}, $(e_1,\ldots,e_k)$ est un \sfio puisque
$1\in\gen{e_1,\ldots,e_k}$.\\
Les anneaux vérifiant ces \prts \eqves ont été appelé  \gui{clean rings} dans \cite[Nicholson]{Nic}. Nous les appellerons les anneaux \emph{cleans}.%
\index{anneau!clean}\index{clean!anneau}

\emph{2.}
Les anneaux cleans sont  stables par quotient et par produit
fini. Tout anneau local est clean.

\emph{3.}
Si $\gB\red$ est clean, il en va de même pour~$\gB$. En déduire
qu'un anneau \zed est clean.

\emph{4.}
Si $\gB\red$ est clean,
 $\gB$ relève les \idms de $\gB/\Rad\gB$.

On dit qu'un anneau $\gA$ est \emph{\slgbz} si l'anneau $\gB=\gA\sur{\Rad\gA}$
est clean.  On dit qu'il est
\emph{\smlz} s'il est \slgb et si $\BB(\gA\sur{\Rad\gA}\!)$ est une \agB
bornée.
}
\end{exercise}
%--- end -exercise-----------------------------------------

%--- Exercise{exo2semilocal}-------------
\begin{exercise}
\label{exo2semilocal} ~(Anneaux \slgbsz, 2) 
{\rm Suite de l'exercice \ref{exo1semilocal}.

 \emph{1.}  
Un anneau local est \smlz. 
 
 \emph{2.}  
Un anneau \slgb et résiduellement connexe est local.
 
 \emph{3.} 
Un anneau \plc est \slgbz.
 
 \emph{4.}  
Un anneau \slgb est \lgbz.
 
 \emph{5.}  
Les anneaux \slgbs sont stables par quotient et par produit fini.
 
 \emph{6.} En \clamaz, un anneau est \sml \ssi il possède un
nombre fini d'\idemasz.  
}
\end{exercise}
%--- end -exercise-----------------------------------------

%%:     Fact{propRelIdm}
%\begin{fact}\label{propRelIdm} ~
%Un anneau $\gA$ \slgb relève les \idmsz.
%\end{fact}

%%
%\begin{proof}
%Soit $a\in\gA$ un \elt \idm dans  $\gA\sur{\Rad\gA}$  et $b=1-a$.
%Puisque $\gen{a,b}=1$ il existe deux \idms \orts $e$ et $f$ dans
%$\gA$ tels que $ae+bf$ est inversible. Puisque $\gen{e,f}=1$ on a $f=1-e$.
%Dans le quotient on a donc 4 \idms $a,b,e,f$ avec $a+b=1$, $ab=0$, $e+f=1$ et $ef=0$. Ainsi $ae$, $bf$, $af$, $be$
%forment un \sfioz. Comme $ae+bf$ est un \idm \iv cela donne $ae+bf=1$,
%d'où $af=be=0$.
%Finalement dans le quotient $a=e$ et $b=f$.
%\end{proof}
%%

%--- Exercise{exoNagatalocal}-------------
\begin{exercise}
\label{exoNagatalocal} 
(\Prts du localisé de Nagata)
{\rm Voir aussi l'exercice \ref{exoPrufNagata}.\\ 
Soit $\gA$ un anneau et
$U\subseteq\AX$ le \mo des \pols primitifs.\\
Notons $\gB=U^{-1}\AX=\gA(X)$ le localisé de Nagata de $\AX$.
 
 \emph{0.} %\label{i0exoNagatalocal}
Donner une \dem directe du fait que $\gB$ est \fpt sur $\gA$.
 
 \emph{1.} %\label{i1exoNagatalocal}
$\gA\cap\gB\eti=\Ati$.
 
 \emph{2.} %\label{i2exoNagatalocal}
$\Rad\gA= \gA\cap\Rad\gB$ et $\Rad\gB=U^{-1}(\Rad\gA)[X]$.
 
 \emph{3.} %\label{i3exoNagatalocal}
$ \gB\sur{\Rad\gB}\simeq  (\gA\sur{\Rad\gA}\!)(X)$.
 
 \emph{4.} %\label{i4exoNagatalocal}
Si $\gA$ est local (resp. local \dcdz), alors $\gB$ est local (resp. local \dcdz).
 
 \emph{5.} %\label{i5exoNagatalocal}
Si $\gA$ est un corps (resp. un \cdiz), alors $\gB$ est un corps (resp. un \cdiz).

} \end{exercise}
%--- end-exercise-----------------------------------------

%--- Exercise{exo2Nagata}-------------
\begin{exercise}
 \label{exo2Nagata} (Localisé de Nagata en plusieurs \idtrsz)\\
 {\rm  Soit $U$ l'ensemble des \pols primitifs de $\gA[X,Y]$.
 
 \emph{1.} Montrer que $U$ est un filtre. 
\\  
On note $\gA(X,Y)=U^{-1}\gA[X,Y]$, on l'appelle le \emph{localisé de Nagata de $\gA[X,Y]$}.
 
 \emph{2.} Montrer que l'application canonique $\gA[X,Y]\to\gA(X,Y)$ est injective
 et que l'on~a un \iso naturel $\gA(X,Y)\simarrow \gA(X) (Y)$.
 
 \emph{3.} Généraliser les résultats
 de l'exercice \ref{exoNagatalocal}.
 } \end{exercise}
%--- end-exercise-----------------------------------------

%%%%%%%%%%%%%%%%%%%%%%%%%%%%%%%%%%%%%%%%%%%%%%%%%%%%%%%%%%%%%%%%%%%%%%%%%%%

%--- exercise{exoAlgMon}-------------
\begin{exercise}\label{exoAlgMon}
{(Algèbre d'un \mo et \ids binomiaux)}\index{algèbre!d'un monoïde}\\
{\rm Soit $(\Gamma,\cdot,1_{\Gamma})$ un \mo commutatif noté multiplicativement, et $\gk$ un anneau commutatif. 
\\
\hbox{L'\emph{\alg de $(\Gamma,\cdot,1_{\Gamma})$ sur $\gk$}}, notée $\gk[(\Gamma,\cdot,1_\Gamma)]$ ou plus simplement $\kGa$, est
formée à partir du \kmo libre sur $\Gamma$ (si $\Gamma$ n'est pas supposé discret, voir l'exercice~\ref{propfreeplat}). Si $\gk$ est non trivial, on identifie tout \elt $\gamma$ de $\Gamma$ à son image dans le module libre. En cas de doute concernant $\gk$, on devrait noter $1_{\gk}\gamma$ au lieu de $\gamma$ cet \elt de~$\kGa$.
\\
La loi produit $\times$ de $\kGa$ est obtenue en posant $\gamma\cdot\gamma'=\gamma\times \gamma'$ et en prolongeant par $\gk$-bilinéarité.   
Notons que $1_{\gA}1_\Gamma=1_{\gk[\Gamma]}$. En pratique, on identifie $\gk$
à un sous-anneau de $\kGa$, et l'on identifie les trois $1$ ci-devant.
\begin{enumerate}
\item Vérifier que la \klg $\kGa$, considérée avec l'application  

\snic{\iota_{\gk,\Gamma}:\Gamma\to\kGa, \, \gamma\mapsto 1_{\gk}\gamma,}

\snii
donne la solution
du \pb \uvl résumé dans le dessin ci-dessous 

\vspace{-1.5em}
\PNV{\Gamma}{\iota_{\gk,\Gamma}}{\psi}{\kGa}{\theta}{\gL}{\mos commutatifs}{morphismes de \mos}{\klgsz}

\vspace{-1.8em}
On dit pour résumer que \emph{$\kGa$ est la \klg librement engendrée par
le \mo multiplicatif $\Gamma$}.
\end{enumerate}

Lorsque la loi de $\Gamma$ est notée additivement, on note $X^{\gamma}$
l'\elt de $\kGa$ image de $\gamma\in\Gamma$ de sorte que l'on a maintenant l'écriture naturelle $X^{\gamma_1}X^{\gamma_2}=X^{\gamma_1+\gamma_2}$.
\\
Par exemple, lorsque $\Gamma=\NN^{r}$ est le \mo additif librement engendré par
un ensemble à $r$ \eltsz, on peut voir les \elts de $\NN^{r}$ comme des multiexposants et $\kGa=\gk[(\NN^{r},+,0)]\simeq \kXr$. Ici  $X^{(m_1,\dots,m_r)}=X_1^{m_1}\cdots X_r^{m_r}$.
\\
Lorsque $\Gamma=(\ZZ^{r},+,0)$, on peut encore voir les \elts de $\ZZ^{r}$ comme des multiexposants et $\gk[\ZZ^{r}]\simeq \gk[\Xr,\fraC1{X_1},\dots,\fraC1{X_r}]$, l'anneau des \pols de~Laurent.
 
Supposons maintenant que $(\Gamma,\cdot,1)$ soit un \mo donné par \gtrs et relations.
Notons~$G$ l'ensemble des \gtrsz.  
\\
Les relations sont de la forme $\prod_{i\in I}  g_i^{k_i}=\prod_{j\in J}  h_j^{\ell_j}$ pour des familles finies 

\snic{(g_i)_{i\in I}$ et $(h_j)_{j\in J}$ dans $G$,
et $(k_i)_{i\in I}$ et $(\ell_j)_{j\in J}$ dans $\NN.}

\snii
Une telle relation peut être codée par le couple $\big((k_i,g_i)_{i\in I},(\ell_j,h_j)_{j\in J}\big)$.
\\
Si l'on espère contrôler les choses, il vaut mieux que $G$ et l'ensemble des relations soient énumérables et discrets. Du point de vue du calcul 
la place centrale est occupée par les présentations finies.
\\
{\bf Notation.} Pour visualiser une \pn finie, par exemple avec $G=\so{x,y,z}$ et des relations $xy^{2}=yz^{3}$, $xyz=y^{4}$ on écrira en notation multiplicative

\smallskip \centerline{\qquad \phantom{$(*)$}\fbox{$\Gamma=_\mathrm{MC}\scp{x,y,z}{xy^{2}=yz^{3},xyz=y^{4}}$},\qquad $(*)$}

 et  en notation additive

\smallskip \centerline{\fbox{$\Gamma=_\mathrm{MC}\scp{x,y,z}{x+2y=y+3z,x+y+z=4y}$}.}

L'indice $\mathrm{MC}$ est mis pour \gui{\mo commutatif}.  

\begin{enumerate} \setcounter{enumi}{1}
\item Montrer que $\kGa\simeq\gk[(g)_{g\in G}]/\fa$, où $\fa$ est l'\id engendré par les différences de \moms $\prod_{i\in I}{g_i}^{k_i}-\prod_{j\in J}
{h_j}^{\ell_j}$
(pour les relations  $\prod_{i\in I}  g_i^{k_i}=\prod_{j\in J}  h_j^{\ell_j}$ données dans la \pn de $\Gamma$). Un tel \id est appelé un \emph{\id binomial}.
Avec l'exemple $(*)$ ci-dessus, on peut donc écrire

\smallskip \centerline{\qquad \phantom{$(*)$}\fbox{$\kGa=_{\gk\mathrm{-\algsz}}\scp{x,y,z}{xy^{2}=yz^{3},xyz=y^{4}}$}.\qquad $(**)$}

\smallskip 
Autrement dit,  $\Gamma=_\mathrm{MC}\scp{truc}{muche}$ implique $\kGa=_{\gk\mathrm{-\algsz}}\scp{truc}{muche}$.
\end{enumerate}

}

\end {exercise}
%--- end -exercise-----------------------------------------

%:sinotenglish
\sinotenglish{

%--- exercise{exoMonfinMon}-------------
\begin{exercise}\label{exoMonfinMon}
{(Monoïde fini monogène)}\\ 
{\rm 
\emph{1.} On suppose que $(\Gamma,\cdot,1)=_\mathrm{MC}\scp{x}{x^{5}=x^{7}}$. Combien $\Gamma$ contient-il d'\eltsz? Montrer d'abord, que $e=x^{6}$ est \idmz, puis que $\kGa\simeq \gk[\epsilon]\times \gk[\alpha]$ où $\epsilon$ et~$\alpha$ sont soumis aux contraintes respectives $\epsilon^{5}=0$ et $\alpha^{2}=1$.

\emph{2.} Généraliser ces résultats en remplaçant $5$ et $2$
par deux entiers  $m$ et $r$ strictement positifs arbitraires.
En particulier $\Gamma$ contient exactement deux \idmsz, $1$ et~$e$, ce dernier égal à $x^{Nr}$
pour $N$ assez grand.

\emph{3.} Tout \mo fini engendré par un seul \elt est un groupe cyclique ou un \mo du type précédent.
}

\end {exercise}
%--- end -exercise-----------------------------------------

%--- exercise{exoMonIdm}-------------
\begin{exercise}\label{exoMonIdm}
{(Idempotent dans un \mo commutatif)} \\
{\rm Soit $e$ un idempotent dans un \mo commutatif $(\Gamma,\cdot,1_{\Gamma})$ et $\gk$ un anneau. 

\emph{1.} On note $\Gamma_{e=1}$, ou si le contexte s'y prête $\Gamma_{1}$, le \mo $\Gamma/(e=1)$. 
\\ On note $e\Gamma=\sotq{ex}{x\in\Gamma}$. La partie $e\Gamma$ est stable pour la multiplication par n'importe quel \elt de $\Gamma$. On dit que $e\Gamma$ est \emph{l'\idp de $\Gamma$ engendré par~$e$}.
\begin{itemize}
\item [\emph{a.}] Montrer que l'application composée $e\Gamma\to\Gamma\to\Gamma_1$ est une bijection et un morphisme pour la multiplication.
On peut identifier les deux \mos $e\Gamma$ et~$\Gamma_1$ au moyen de cet \isoz.
\item [\emph{b.}]  Montrer que le morphisme de \klg $\kGa\to \gk[\Gamma_1]$ est un morphisme de passage au quotient par l'\id $\gen{e-1}$. 
Ceci permet d'identifier $\aqo\kGa{e-1}$ \hbox{et $\gk[\Gamma_1]$}.
\item [\emph{c.}] Montrer l'\eqvc $x\in e\Gamma\iff ex=x$. Ainsi lorsque $\Gamma$ est discret, $e\Gamma$ est une partie détachable de $\Gamma$. 
%%
%\item [\emph{c.}]  En pratique  on identifie $\gk[\Gamma_1]$ avec le \kmo libre sur l'ensemble $e\Gamma$. Donner la table de multiplication pour cette base.
\end{itemize}

\smallskip \rem On peut identifier $\Gamma_{1}$ et l'\id $e\Gamma$ de $\Gamma$. Cependant $e\Gamma$ n'est pas un sous-\mo de $\Gamma$ car l'\elt neutre est $e$ et non pas $1_{\Gamma}$. Nous sommes ici dans la même situation que pour un \idm dans un anneau commutatif.
\eoe

\emph{2.} On note $\Gamma'_e$ le \mo $\Gamma/(e=ex,\,\forall x \in \Gamma)$.  
\\
On dit que $\Gamma'_e$ est \emph{le \mo quotient obtenu en forçant $e$
à être absorbant.}
\begin{itemize}
\item [\emph{a.}] Montrer que l'\egt de $\Gamma'_e$ est obtenue
comme la plus fine relation d'\eqvcz~$\sim$ sur $\Gamma$ qui vérifie
$ex\sim ey$ pour tous $x$, $y\in\Gamma$.   
\item [\emph{b.}]  Montrer que le morphisme de \klgs obtenu en composant

\snic{\kGa\to \gk[\Gamma'_e]\to\aqo{\gk[\Gamma'_e]}{e}}

est un morphisme de passage au quotient par l'\id $\gen{e}$. 
\\
Autrement dit on peut identifier $\aqo\kGa{e}$ et~$\aqo{\gk[\Gamma'_e]}{e}$. 
\item [\emph{c.}] Montrer que  $\aqo{\gk[\Gamma'_e]}{e-1}\simeq \gk$. Ainsi $\gk[\Gamma'_e]\simeq \gk \times \aqo{\gk[\Gamma'_e]}{e}\simeq \gk \times \aqo{\gk[\Gamma]}{e}$.  
\item [\emph{d.}]  On suppose maintenant que  $e\Gamma$ est une partie détachable de $\Gamma$. 
\begin{itemize}
\item Montrer qu'une $\gk$-base de la \klg $\aqo{\gk[\Gamma'_e]}{e}$ est l'ensemble $\Gamma\setminus e\Gamma$.
\item  Donner la table de multiplication pour cette base.
\item  On note $\Gamma_{0}={(\Gamma\setminus e\Gamma)}\cup\so e$ (réunion disjointe). Donner sur $\Gamma_0$
une structure de \mo pour laquelle on a  un \iso naturel
$\Gamma'_e\simarrow \Gamma_0$.
\end{itemize}
\end{itemize}

\smallskip 
\rems
1) Les \mos \pf sont discrets. Mais ce n'est pas un résultat facile. 
\\
2) Lorsque $e\Gamma$ n'est pas détachable, l'ensemble $\Gamma\setminus e\Gamma$ n'est pas toujours explicitement une base de $\aqo\kGa e$ comme \kmoz.
\\
Par exemple si $\Gamma=_{\mathrm{MC}}\scp x {x^{2}=x^{3}, \sfP \Rightarrow x=x^{2}}$, où $\sfP$ est une conjecture non résolue, considérons $e=x^{2}$, de sorte que
$\aqo\kGa e$, en tant que \kmoz, hésite entre $\gk$ et $\gk^{2}$.
%égal à $\gk\times  x\gk$ sans qu'on sache si le deuxième facteur est nul. 
\\
Si l'\kli  naturelle  $\iota:\gk^{\Gamma\setminus e\Gamma}\to \aqo\kGa e$ était un \isoz, il faudrait donner
 $\iota^{-1}(\ov x)$, qui est égal à $x$ si $\sfP$ est fausse et à $0$ si $\sfP$ est vraie. 
\\
3) En notation additive, un \idm $e$ est un \elt vérfiant l'\eqn $2e=e$.
Dans ce cas, il est légitime de noter $\Gamma_{e=0}$ ce qui était noté
 $\Gamma_{e=1}$ en notation multiplicative. Alors qu'en notation multiplicative, on pense un \elt absorbant comme le $0$ d'un anneau commutatif, en notation additive il faut penser un \elt absorbant comme étant égal à $+\infty$.
\eoe
}

\end {exercise}
%--- end -exercise-----------------------------------------

%--- exercise{exoMonExempleabcd}-------------
\begin{exercise}\label{exoMonExempleabcd}
{(Exemple d'étude d'un \mo commutatif fini)}\\ 
{\rm Soient $a$, $b$, $c$, $d \in\NN$ vérifiant $ad - bc \ne 0$.
%Soit $\gk$ un anneau commutatif arbitraire.
\\
On considère le \mo commutatif  $(\Gamma,\cdot,1)$ engendré
par deux \elts $x$, $y$ soumis aux relations $x^a = y^b, \; x^c = y^d$. Avec la notation donnée dans l'exercice~\ref{exoAlgMon}, on a donc $\Gamma=_\mathrm{MC}\scp{x,y}{x^a = y^b,x^c = y^d}$.
Vue la symétrie de la situation on peut supposer $ad>bc$. 
On pose $\Delta=ad-bc$ et $m=bc$.

\emph{1.} Montrer que $x^{m}=x^{ad}$ et $y^{m}=y^{ad}$. En déduire que $\Gamma$ est borné. \\
Si $b$ ou $c=0$, $\Gamma$ est un groupe fini. Sinon calculer un \idm $e\neq 1$. 
%\\
%On a donc $\kGa\simeq \aqo\kGa{e-1}\times \aqo\kGa{e}$ avec 
%$\aqo\kGa{e-1}=\gk[\Gamma_{e=1}]$. 

\emph{2.} Notons $\Gamma_1$ le \mo $\Gamma_{e=1}$ obtenu en ajoutant la relation 
$e=1$. Montrer que~$\Gamma_1$ est un groupe fini d'ordre $\Delta$. Il est cyclique \ssi $\gcd(a,b,c,d)=1$. Sinon, en notant $g$ ce pgcd, et $\delta=\Delta/g$, on obtient $\Gamma_1\simeq(\ZZ/g\ZZ \times \ZZ/ \delta\ZZ,+,0)$,
\cade $\Gamma_1=_\mathrm{MC}\scp{u,v}{u^{g}=1,v^{\delta}=1}$. 

\emph{3.} Quelle est la structure du sous-\mo de $\Gamma$ engendré par $x$?

\emph{4.} On va montrer que $\#\Gamma=ad$. Pour ceci on met en place un processus de réécriture des \moms $x^{k}y^{\ell}$. 
\\
Une relation d'ordre total $\preceq$ sur $\NN^{r}$ est appelée un \ix{ordre monomial} si elle raffine l'ordre partiel naturel (i.e. $m\preceq m+p$ pour
$m,p\in\NN^{r}$) et si elle est compatible avec l'addition:  $m\preceq n$
implique $m+p\preceq n+p$ pour  $m$, $n$, $p\in\NN^{r}$.
\\
On définit sur $\NN^2$ un ordre monomial $\preceq$ satisfaisant 
$(0,b)\prec(a,0)$ et $(c,0)\prec(0,d)$ (où si l'on préfère 
$ y^{b}\prec x^{a}$ et $x^{c}\prec y^{d}$) comme suit. 
\\
On pose  $\ell_1(u,v)=du+cv$, $\ell_2(u,v)=bu+av$ et

\snic{\alpha \preceq \alpha' \iff \big( {\ell_1}({\alpha}) <  {\ell_1}({\alpha'}) 
\hbox{ ou }\big( {\ell_1}({\alpha}) =  {\ell_1}({\alpha'})\hbox{ et } {\ell_2}({\alpha}) \le  {\ell_2}({\alpha'})\big)\big).}

Le procédé de réécriture consiste à réécrire $x^{k}y^{\ell}$
sous forme $x^{k-a}y^{\ell+d}$ si $k\geq a$, et sinon sous forme $x^{k+c}y^{\ell-d}$ si $\ell\geq d$. On remplace de cette manière un \mom \gui{formel} par un autre qui lui est égal dans $\Gamma$. Le processus s'arrête lorsque $k<a$ et $\ell<d$. Il y a donc $ad$ \moms possibles à la fin du processus.\\
Ici on demande de démontrer  trois choses.
\begin{itemize}
\item  [\emph{a.}] On a bien défini un ordre monomial pour lequel  $ y^{b}\prec x^{a}$ et $x^{c}\prec y^{d}$.
\item [\emph{b.}] Quel que soit le \mom de départ, le processus de réécriture termine  après un nombre fini d'étapes.
\item [\emph{c.}] Les  \moms $x^{k}y^{\ell}$ pour $k\in\lrb{0..a-1}$ et $\ell\in\lrb{0..d-1}$ sont des \elts deux à deux distincts de $\Gamma$.
\end{itemize}
Ce qui montre que $\#\Gamma=ad$.

\smallskip 
\emph{5.} En déduire que $\#\Gamma'_e=m+1$ ($\Gamma'_e$ est obtenu en forçant $e$ à devenir absorbant, voir l'exercice \ref{exoMonIdm}). 
%Ainsi $\aqo\kGa e$ est une
%\klg libre de base $\Gamma_0\setminus\so e$. Si $\gk$ est un \cdiz, c'est un \alo \zedz. 

\emph{6.} Montrer que les seuls \idms de $\Gamma$ sont $1$ et $e$
(distincts si $m>0$).
}
\end {exercise}
%--- end -exercise-----------------------------------------

%--- exercise{exoMonExemple}-------------
\begin{exercise}\label{exoSyst_abcd}
{(\'Etude d'un \sys \zedz, intersection de deux courbes planes affines)}\\
{\rm On étudie dans cet exercice l'intersection de deux courbes monomiales
sur un anneau de base arbitraire $\gk$. Ces deux courbes ont un \gui{cusp} à l'origine et, en langage savant, l'intersection est prise \gui{au sens des schémas}.\\
On va décrire la composante à l'origine de cette intersection et son
\copz. Les deux \klgs correspondantes sont des \kmos libres
de rang fini\footnote{Chacune des deux \algs $\gA_0$ et $\gA_1$ est donc \gui{\zede sur $\gk$} au sens de la \dfnz~\ref{defiKdimMor}.}. Lorsque $\gk$ est un \cdiz, ce sont deux anneaux \zedsz, le premier est local et le second est \smlz.
\\
Soient $a$, $b$, $c$, $d \in\NN\etl$ vérifiant $ad - bc \ne 0$. On veut montrer
que le \sys

\snic {
x^a = y^b, \quad x^c = y^d,
}

\snii
est \zed et calculer sa multiplicité en l'origine.  \\
Soient $f$, $g \in \gk[X,Y]$
définis par $f = X^a - Y^b$, $g = X^c - Y^d$ et:

\snic {
\gA = \gk[x,y] = \aqo {\gk[X,Y]} {f,g}, \qquad
\gA_0 = \gk[x,y]_{1 + \gen {x,y}}
}

\snii
De manière précise, on montrera ici que $\gA$ est libre sur
$\gk$ de rang $\max(ad,bc)$ et l'on donnera une base monomiale explicite.
En ce qui concerne la composante en l'origine, \hbox{la \klg
$\gA_0$} est libre de rang $m=\min(ad,bc)$ avec aussi une base
monomiale explicite. Rappelons que ce rang $m$ est appelé la multiplicité
du zéro~$(0,0)$ lorsque 
$\gk$ est un \cdiz.
\\
En fait vue la symétrie de la situation on supposera \spdg que~\fbox{$ad>bc$}. On posera $m=bc$ et $\Delta=ad-bc$.

En fait la solution de l'exercice présent est presque entièrement donnée 
par l'étude du \mo $\Gamma$ de l'exercice \ref{exoMonExempleabcd} et par les considérations
développées dans les exercices \ref{exoAlgMon}, \ref{exoMonfinMon} et \ref{exoMonIdm}.

\emph {1.}
Donner une preuve directe simple que $\gA$ est un \kmo \tfz. Donner les zéros de $\gA$
dans $\gk$ lorsque $\gk$ est un \cdi contenant les racines $\Delta$-èmes de l'unité.

\emph {2.}
D'après l'exercice \ref{exoAlgMon},   l'\alg $\gA$ est isomorphe à $\kGa$,
où $\Gamma$ est le \mo étudié dans l'exercice \ref{exoMonExempleabcd}: c'est donc un \kmo libre de base $\Gamma$ de rang~$ad=\#\Gamma$. 
\\
Plus \prmt $\kGa$
est un \kmo libre qui admet comme base les \momsz~\hbox{$x^{k}y^{\ell}$}
pour pour $k\in\lrb{0..a-1}$ et $\ell\in\lrb{0..d-1}$.\\
Soit $e$  l'\idm de $\Gamma$ dans l'exercice \ref{exoMonExempleabcd}, on a $\gA\simeq \aqo\gA e \times \aqo \gA{1-e}$. 
\\
D'après les exercices \ref{exoMonIdm}
et \ref{exoMonExempleabcd}, $\aqo \gA{1-e}\simeq\gk[\Gamma_{e=1}]$ où
$\Gamma_{e=1}\simeq e\Gamma$ est un groupe d'ordre $\Delta$.
Et $\aqo \gA{e}$ admet pour base les $m$ \moms de $\Gamma\setminus e \Gamma$.
En outre $x$ et $y$ sont nilpotents dans $\aqo \gA{e}$.

\emph {3.}  
La composante à l'origine est par \dfn $\gA_0=\gA_{1+\gen{x,y}}$.\\
Montrer qu'il s'agit de l'\alg $\gA[\fraC 1 {1-e}]\simeq \aqo\gA e$. 
D'après l'exercice \ref{exoMonIdm}, c'est donc un \kmo libre de base $\Gamma\setminus e\Gamma$, et son  rang est $m$ d'après l'exercice \ref{exoMonExempleabcd}.

\emph {4.}  Donner une base de $\gA_1=\gA[\fraC 1 e]\simeq \aqo\gA{1-e}$ sur $\gk$ et calculer le discriminant~$\Disc_{\gA_1/\gk}$.

\emph {5.}  On note $F$, $G \in \gk[X,Y,Z]$ les \pols homogénéisés de $f$ et $g$. 
\'Etudier l'intersection des deux courbes 
$\{F=0\}$ et $\{G=0\}$ dans $\PP^2(\gk)$ (ici $\gk$ est un \cdi contenant les racines $\Delta$-èmes de l'unité).

}
\end {exercise}
%--- end -exercise-----------------------------------------
}
%: fin sinotenglish

% fin des exos
%%%%%%%%%%%%%%%%%%%%%%%%%%%%%%%%%%%%%%%%%%%%%%%%%%%%%%%%%%%%%%%%%%%%%%%%%

%:  solutions d'exos
\penalty-2500
\sol 
%%%%%%%%%%%%%%%%%%%%%%%%%%%%%%%%%%%%%%%%%%%%%%%%%%%%%%%%%%%%%%%%%%%%%%%%

\exer{exoRadNonCom}
\emph{1.} Si $c$ est inverse à droite et à gauche c'est l'unique inverse à gauche car $c'a=1$ implique $c'=c'ac=c$. \\
Inversement, puisque $ca=1$, on a $(c+1-ac)a=ca+a-aca=1$. 
Donc $c+1-ac$ est un inverse à gauche, et s'il y a unicité, $1-ac=0$.

 \emph{2.} On vérifie que $v=1+bua$ convient.

\emph{3.} Si $u(1-xa)=1$, alors $u=1+uxa$, donc il est \iv à gauche. Ainsi $u$ est \iv à droite et à gauche, et $1-xa$ \egmtz.

%%%%%%%%%%%%%%%%%%%%%%%%%%%%%%%%%%%%%%%%%%%%%%%%%%%%%%%%%%%%%%%%%%%%%%%%

\exer{exoLocalFreenessLemma}
Soient $x_1$, \ldots, $x_n\in E$ tels que $(\ov {x_1}, \ldots, \ov {x_n})$ soit une
$\gk$-base de~$E/\fm E$. D'après Nakayama, les $x_i$ engendrent $E$. Soit $u
: \Ae n \twoheadrightarrow E$ la surjection %réalisant 
$e_i \mapsto x_i$. En
étendant les scalaires à $\gK$, on obtient une surjection $U : \gK^n
\twoheadrightarrow \gK\otimes_\gA E$ entre deux \evcs de même dimension $n$,
donc un \isoz. 
\\
\Deuxcol{.7}{.15}{ Puisque $\Ae n\hookrightarrow \gK^n$, on en déduit que
$u$ est injective.
En effet, si $y \in \Ae n$ vérifie $u(y) = 0$, alors
$1 \otimes u(y) = U(y) = 0$ dans $\gK\otimes_\gA E$, donc $y = 0$, cf.\,le
diagramme ci-contre:}
{$
\xymatrix @R = 7pt {
\Ae n \ar@{(->}[d] \ar@{->>}[r]_u & E\ar[d] \\
\gK^n \ar@{->>}[r]^{U~~}              & \gK\otimes_\gA E \\
}
$}

%\vspace{-2mm}
Bilan: $u$ est un \iso et $(x_1, \ldots, x_n)$ est une $\gA$-base de $E$.

%%%%%%%%%%%%%%%%%%%%%%%%%%%%%%%%%%%%%%%%%%%%%%%%%%%%%%%%%%%%%%%%%%%%%%%%
\exer{exoNakayamaAppli}
D'après  Nakayama, $(e_1,
\ldots, e_n)$ engendre le $\gA$-module $E$.\\
Soit $L = \Ae n$ et $\varphi : L \twoheadrightarrow E$, l'\ali (surjective) qui transforme la base canonique de $L$
en $(e_1, \ldots, e_n)$. Par hypothèse, $\overline\varphi
: L/aL \to E/aE$ est un \isoz.
Montrons que $\Ker\varphi = a\Ker\varphi$.
Soit $x \in L$ avec $\varphi(x) = 0$; \hbox{on a
$\overline\varphi(\overline x) = 0$}, donc $\overline x = 0$, i.e.
$x \in aL$, disons $x = ay$ avec $y \in L$. Mais
$0 = \varphi(x) = a\varphi(y)$ et $a$ étant $E$-régulier,
$\varphi(y) = 0$. On a bien $\Ker\varphi \subseteq a\Ker\varphi$.
Puisque $E$ est \pfz, $\Ker\varphi$ est \tfz, et l'on peut appliquer Nakayama à l'\egt $\Ker\varphi
= a\Ker\varphi$. On obtient $\Ker\varphi = 0$: $\varphi$
est un \isoz.

%%%%%%%%%%%%%%%%%%%%%%%%%%%%%%%%%%%%%%%%%%%%%%%%%%%%%%%%%%%%%%%%%%%%%%%%
\exer{exoVariationLocGenerated} 
\emph {1a, 1b, 1c.}
L'idée est de développer $\det(A+B)$ comme fonction multi\lin
des colonnes de $A+B$. Le résultat est une somme de $2^n$ \deters de matrices
obtenues en mélangeant des colonnes $A_j$, $B_k$ de $A$ et $B$.
On écrit:

\snic {
\det(A_1+B_1, \dots, A_n+B_n) = \sum_{2^n} \det(C_1, \dots, C_n) 
\quad \hbox {avec $C_j = A_j$ ou $B_j$.}
}

Pour  $J\in\cP_n$, notons $\Delta^{\rm col}_J$ le \deter lorsque $C_j =
B_j$ \hbox{pour $j \in J$} \hbox{et $C_j = A_j$} sinon.  Avec cette notation, on a donc:

\snic {
\det(A+B) = \sum_J \Delta^{\rm col}_J.
}

%\sni
Si $\#J \ge n-r$, alors $\Delta^{\rm col}_J \in \cD_{n-r}(B)$; sinon
$\#\ov {J} \ge r+1$ et donc $\Delta^{\rm col}_J \in \cD_{r+1}(A)$.

\emph {1d.}
Considérer $A = \Diag(\an)$.

\emph {2.}
On écrit $\ux=\uy\,U$ avec $U \in \gA^{r\times n}$, $\uy=\ux\,V$ avec $V \in
\gA^{n\times r}$. \\
On pose $A = VU$, $B = \In-A$. On a \framebox
[1.1\width][c]{$\ux\,B = 0$} et $\rg(B) \ge n-r$ puisque $\rg(A) \le r$.
L'\egt encadrée montre, pour $I \in \cP_{r,n}$ et $\nu$ mineur
 de $B$ sur les lignes de $\ov I$, l'inclusion 
$\nu M \subseteq \gen{(x_i)_{i\in I}}$.  Et c'est gagné car $1 \in \cD_{n-r}(B)$.

\Prmtz, notons $\Delta^{\rm row}_J$ le \deter  de la matrice \gui {mixte} dont
les {\em lignes} d'indice $i \in J$ sont les {\em lignes} correspondantes de
$B$ et les {\em lignes} d'indice $i \in \ov J$ sont celles de $A$.  Pour $J
\supseteq \ov I$, $\Delta^{\rm row}_J$ est une \coli de mineurs 
 de $B$ sur les lignes de $\ov I$. On pose donc

\snic {
s_I = \sum_{J \mid J \supseteq \ov I} \Delta^{\rm row}_J.
}

Alors d'une part, $s_IM \subseteq \gen {(x_i)_{i\in I}}$, et d'autre part, puisque
$\rg(B) \ge n-r$:

\snic {
1 = \sum_{I \in \cP_{r,n}} s_I.
}

\emph {3.}
Clair en utilisant le lemme des \lons successives (fait \ref{factLocCas}).

\emph {4.}
Si une matrice $A \in \Mn(\gA)$ existe comme indiquée, la preuve
du point \emph {2} s'applique \hbox{avec $B=\In-A$}.

La réciproque est problématique car la contrainte $\rg(A) \le r$ n'est pas
\lin en les \coes de $A$.  On y arrive pourtant pour $r=1$ par un autre moyen,
 (voir le \tho \ref{propmlm}).

%%%%%%%%%%%%%%%%%%%%%%%%%%%%%%%%%%%%%%%%%%%%%%%%%%%%%%%%%%%%%%%%%%%%%%%%

\exer{exoMorAdcp}
\emph{1.}
La condition $b'(1-a'b')=0$ s'obtient par 
 $\varphi\big(b(1-ab)\big)=0$. \\
Supposons que $\varphi$ est un morphisme d'anneaux \dcps et montrons \hbox{que $\varphi(\Rad\gA)\subseteq\Rad\gB$}: soit $a \in \Rad\gA$, alors $b = 0$ (par unicité de $b$), \hbox{donc $b' = 0$} et  $a'=a'(1-a'b')\in\Rad\gB$.\\
Réciproquement, supposons $\varphi(\Rad\gA)\subseteq\Rad\gB$. Si
$a$, $b\in\gA$ vérifient $b(1-ab)=0$ et $a(1-ab)\in\Rad\gA$, alors 
$\varphi\big(a(1-ab)\big)=a'(1-a'b')\in\Rad\gB$.

%%%%%%%%%%%%%%%%%%%%%%%%%%%%%%%%%%%%%%%%%%%%%%%%%%%%%%%%%%%%%%%%%%%%%%%%
\exer{exopropAlgb1}
\emph{1.} 
Il faut et il suffit que les $a_i$ vérifient l'\egt $(\star_d)$ pour 
les \moms de degré total $\leq d$. Soit $M=M(\uX)=\uX^{\alpha}$
un tel \mom  avec $\abs\alpha=j\leq d$. 
Puisque $M(s^{r}\uX)=s^{rj}M$, on veut réaliser

\snic {
s^e s^{-j}\,M = a_0(s) \,M + a_1(s) s^{j}\,\uX + \cdots + a_d(s)s^{dj}\,M, 
}

\cad après simplification par $M$ et multiplication par $s^j$:

\snic {
s^e = a_0(s)s^j + a_1(s)s^{2j} + \cdots + a_d(s)s^{(d+1)j} =
\sum_{i=0}^{d} a_i(s)(s^j)^{i+1}.
}

%\sni
Introduisons  le \pol $F(T) \in \ZZ[s][T]$
défini par $F(T) = T\,\sum_{i=0}^{d} a_i(s)T^{i}$. \\
Alors $\deg_T F \le d+1$ et
$F$ réalise l'interpolation $F(0) = 0$ et $F(s^j) = s^e$ pour
les $j \in
\lrb{1..d}$. Or un \pol $F \in \ZZ[s][T]$ qui satisfait cette interpolation
est le suivant:
$$
F(T) = s^e - (s^0-T)(s^1-T)(s^2-T) \cdots (s^d-T). \leqno (\#_d)
$$
Marche arrière toutes. Considérons le \pol défini par
l'\egt $(\#_d)$.
Il est de degré $d+1$ en $T$, nul en $T=0$, donc il s'écrit

\snic {
F(T) = T\, \sum_{i=0}^{d} a_i(s)T^{i}, \;\;\;\; \hbox {avec}\;\;
a_0(s),\, \ldots,\, a_d(s) \in \ZZ[s].
}

%\sni
Ces \pols $a_i(s)$ ont la \prt voulue.

%\sni
\emph{2.}
L'appartenance demandée se déduit de l'\egt $(\star_d)$
en évaluant $\uX$ en $\ux$.

\emph{3.} Supposons que $P$ soit de degré total $\leq d$.
Le fait que $P(\ux/s) \in (\gA_S)\eti$ signifie, dans~$\gA$, que $y = s^e
P(\ux/s)$ divise un \elt $t$ de $S$.  D'après le point \emph{2}, $y$
est dans l'\id engendré par les valeurs de~$P$; il en est de même de $t$.

%%%%%%%%%%%%%%%%%%%%%%%%%%%%%%%%%%%%%%%%%%%%%%%%%%%%%%%%%%%%%%%%%%%%%%%%
\exer{exoLgb2}{
\emph{1.} Soit $b\in\gA$ \iv modulo $\fa$. Il existe $a\in\fa$ tel
que $1\in\gen{b,a}$. Le \pol $aT+b$ prend les valeurs \come $a$, $a+b$,
donc il représente un \iv $b'=at+b$. Alors, $b'\equiv b \mod \fa$
avec $b'$ \ivz.\\
\emph{2.} On écrit $x=ay$, $y=bx$, donc $(1-ab)x=0$. \\
Puisque $b$
est \iv modulo $1-ab$, il existe $u\in\Ati$  tel que
$u\equiv b \mod {1-ab}$. Alors $ux=bx=y$.
}
%%%%%%%%%%%%%%%%%%%%%%%%%%%%%%%%%%%%%%%%%%%%%%%%%%%%%%%%%%%%%%%%%%%%%%%%

\exer{exo1semilocal}~\\
Pour deux \idms \orts $e$, $e'$, on a $\gen {ex,e'x'} =
\gen {ex + e'x'}$.  
\\
Donc pour  $(e_1, \ldots, e_k)$, on a $\gen
{e_1x_1 + \cdots + e_kx_k} = \gen {e_1x_1,\dots,e_kx_k}$.
\\
Ainsi, $e_1x_1 + \cdots + e_kx_k$ est \iv
\ssi $e_1x_1$, \ldots, $e_kx_k$ sont \comz. En conséquence, dans le
contexte  \emph{1.a}, soient $y_i \in \gen {x_i}$ avec $(y_1,\ldots, y_k)$ \com (a fortiori $(x_1, \ldots, x_k)$ \comz); si des
\idms $(e_1, \ldots, e_k)$ fonctionnent pour $(y_1, \ldots, y_k)$, ils
fonctionnent aussi pour $(x_1, \ldots, x_k)$.  Quitte à remplacer $x_i$ par
$u_ix_i$. On pourra donc supposer $\sum x_i = 1$.
\\
Pour deux \idms $e$, $e'$, on a $e \perp e'$ \ssi $1-e$, $1-e'$ sont \comz.

 \emph{1.}
\emph {c} $\Rightarrow$ \emph {d.} En prenant $x_1 = x$, $x_2 = 1+x$, $e = e_2
= 1-e_1$, on a $e_1x_1 + e_2x_2 = x+e$.\\
\emph {d}
$\Rightarrow$ \emph {c.} On peut supposer $1 = -x_1 + x_2$; on pose $x = x_1$.
\\
Alors, $e + x = (1-e)x + e(1+x) = (1-e)x_1 + ex_2$.
\\
\emph{a} $\Leftrightarrow$ \emph {b.} Obtenu facilement en posant $\gB_i = \aqo{\gB}{1-e_i}$.
\\
\emph {c} (ou \emph {d}) $\Rightarrow$ \emph {a.} Par \recu sur~$k$. 
On peut supposer $1 = \sum_i x_i$: il existe un
\idm $e_1$ tel que $e_1x_1 + (1-e_1)(1-x_1)$ soit \ivz.  On a $1 \in
\gen {x_2, \ldots, x_k}$ dans le quotient $\aqo{\gB}{e_1}$ qui possède aussi
la \prt \emph {d}; donc, par \recuz, il existe $(e_2, \ldots, e_k)
$ dans $\gB$ formant un \sfio dans le quotient $\aqo{\gB} {e_1}$ avec $e_2x_2 +
\cdots + e_kx_k$ \iv dans $\aqo{\gB}{e_1}$. Alors, $(e_1, (1-e_1)e_2,
\ldots, (1-e_1)e_k)$ est un \sfio de $\gB$ et $e_1x_1 + (1-e_1)e_2x_2 + \cdots
+ (1-e_1)e_kx_k$ est \iv dans~$\gB$.

 \emph{2.}
Pas de difficulté.

 \emph{3.}
Soit $x \in \gB$; il existe un \idm $e \in \gB\red$, tel que $e + \ov x$ soit
\iv dans~$\gB\red$. On relève $e$ en un \idm $e' \in \gB$. Alors, $e'
+ x$ relève $e + \ov x$ donc est \ivz.  Soit $\gB$ un anneau \zedz;
quitte à remplacer $\gB$ par~$\gB\red$, on peut supposer $\gB$ réduit; si
$x \in \gB$, il existe un \idm $e$ tel que $\gen {x} = \gen {1-e}$; alors $e +
x$ est \ivz.

 \emph{4.}
Soit $a\in\gB$ un \elt \idm dans $\gB\sur{\Rad\gB}$ et $b=1-a$.  \\
Puisque
$\gen{a,b}=1$, il existe deux \idms \orts $e$ et $f$ dans $\gB$ tels \hbox{que
$ae+bf$} est \ivz. Puisque $\gen{e,f}=1$, on a $f=1-e$.  Maintenant, on
raisonne dans le quotient.  Le \sys  $(ae, bf, af, be)$
est un \sfioz. Comme
$ae+bf$ est  \ivz, on a $ae+bf=1$, \hbox{d'où $af=be=0$}.  Finalement
(dans le quotient) $a=e$ et $b=f$.
%%%%%%%%%%%%%%%%%%%%%%%%%%%%%%%%%%%%%%%%%%%%%%%%%%%%%%%%%%%%%%%%%%%%%%%%

\exer{exo2semilocal}
Un anneau $\gA$ est local \ssi $\gA\sur{\Rad\gA}$ l'est;
un anneau $\gA$ est \slgb \ssi $\gA\sur{\Rad\gA}$ l'est.

\emph {1.}
Un anneau local vérifie le point \emph {1d} de
l'exercice précédent avec $e=0$ ou $e=1$. 

\emph {2.}
$\gA\sur{\Rad\gA}$ est connexe, semi-local donc local
(utiliser le point \emph {1d} de l'exercice
précédent sachant que  $e=0$
ou $1$); donc $\gA$ est local.

\emph {4.}
Se démontre pour l'anneau résiduel et résulte de la constatation
suivante.\\
Si $f$ est un \pol en $n$ \idtrs et $(e_1,\ldots, e_k)$ un
\sfioz, alors pour $(\xk)$ dans $\gA^n$, puisque l'\homo d'\evn commute aux
produits directs, on a l'\egt\\
\centerline{$f(e_1 x_1+\cdots +e_k
x_k) = e_1 f(x_1)+\cdots +e_k f(x_k)$.}
 
\emph {6.}
Un anneau $\gA$ possède un nombre fini d'\idemas \ssi c'est le cas de
$\gA\sur{\Rad\gA}$.  En \clamaz, $\gA\sur{\Rad\gA}$ est un produit fini de
corps.

%%%%%%%%%%%%%%%%%%%%%%%%%%%%%%%%%%%%%%%%%%%%%%%%%%%%%%%%%%%%%%%%%%%%%%%%
\exer{exoNagatalocal} ~\\
Dans la suite $f=\sum_ib_iX^i\in\AX$ et
$g=\sum_ic_iX^i\in U$, avec $1=\sum_ic_iu_i$.

\sni
\emph{0}. Soit $\uT$ un jeu d'\idtrs sur $\gA$ et $\gA(\uT)$ le localisé de
Nagata.  On sait \hbox{que $\gA(\uT)$} est plat sur $\gA$ et l'on montre que tout \sli sur $\gA$ qui admet une
solution sur $\gA(\uT)$ admet une solution sur $\gA$.
\\
Soit donc le \sli  $A x = b$ avec $A\in\Ae{n\times m}$ et $b \in \gA^n$. Supposons l'existence d'une
solution sur $\gA(\uT)$; elle s'écrit $P/D$ avec
$P \in \gA[\uT]^m$ \hbox{et $D \in \gA[\uT]$} un \pol primitif. On a donc
$A\,P = D\,b$ sur $\gA[\uT]$.
\\
\'Ecrivons $P = \sum_\alpha x_\alpha \uT^\alpha$ avec $x_\alpha
\in \gA^m$ et $D = \sum_\alpha a_\alpha \uT^\alpha$ où les
$a_\alpha \in \gA$ sont \comz. L'\egt $A\,P = D\,b$ donne
$A\, x_\alpha = a_\alpha b$ pour chaque $\alpha$.
\\
Si $\sum u_\alpha a_\alpha = 1$,
le vecteur $x = \sum_\alpha u_\alpha x_\alpha$ est
solution du \sys $Ax = b$.

\sni\emph{1}. Soit $a\in\gA$ \iv dans $\gB$. Il existe $f,\,g$
tels que $af=g$, donc $a$ et $f$ sont primitifs: $a\in\Ati$.

\sni\emph{2}. Montrons $\Rad\gA\subseteq\Rad\gB$. Soit $a\in\Rad\gA$, on veut montrer
que $1+a (f/g)$ est \iv dans $\gB$,
\cad que $g+af\in U$. On souhaite $1\in\gen{(c_i+ab_i)_i}$;
or cet \id contient $\sum_iu_i(c_i+ab_i)=1+az\in \Ati.$ \\
On sait donc que $\Rad\gB\supseteq U^{-1}(\Rad\gA)[X]$.
Soit $h=\sum_{i=0}^na_iX^{i}$. Montrons \hbox{que $h\in\Rad\gB$} implique
$a_n\in\Rad\gA$. 
\\
On en déduira par \recu que $h\in(\Rad\gA)[X]$.
\\
On considère $a\in\gA$, on prend $f=a$ et $g=X^n-a(h-a_nX^n)$.
Il est clair \hbox{que $g\in U$}, donc $g+fh=(1+aa_n)X^n$ doit être \iv dans $\gB$,
i.e., $1+aa_n$ doit être dans~$\Ati$.
%%%%%%%%%%%%%%%%%%%%%%%%%%%%%%%%%%%%%%%%%%%%%%%%%%%%%%%%%%%%%%%%%%%%%%%%

%%%%%%%%%%%%%%%%%%%%%%%%%%%%%%%%%%%%%%%%%%%%%%%%%%%%%%%%%%%%%%%%%%%%%%%%%%%

%%%%%%%%%%%%%%%%%%%%%%%%%%%%%%%%%%%%%%%%%%%%%%%%%%%%%%%%%%%%%%%%%%%%%%%%%%%
%: correction exoAlgMon

\exer{exoAlgMon}
\emph{(Algèbre d'un \mo et \ids binomiaux)}\\
\emph{1.} Tout d'abord on vérifie que $\kGa$ est bien une \klg et $\iota_{\gk,\Gamma}$ un morphisme de \mosz. Ensuite, si $\alpha:\Gamma\to\gA$
est un morphisme de \mosz, il y a a priori une unique manière de le prolonger en un morphisme $\wi\alpha$ de \klgs de $\kGa$ vers~$\gA$: il faut poser $\wi\alpha\big(\sum_{\gamma\in I}a_{\gamma}\gamma\big)=\sum_{\gamma\in I}a_{\gamma}\alpha(\gamma)$ (ici, $I$ est une partie finiment énumérée de~$\Gamma$). 
\\
On vérifie alors que $\wi\alpha$ est bien un morphisme de \klgsz. \Llec est invité\e à vérifier tous les détails lorsque $\Gamma$ n'est pas supposé discret, en se basant sur l'exercice~\ref{propfreeplat}.

\emph{2.} Il s'agit d'un résultat \gnl d'\alg \uvlez, car on est ici dans le cadre de structures \agqs purement équationnelles. Pour obtenir une \klg au moyen de \gtrs et relations données par des \egts de \momsz, on peut d'abord construire le \mo défini de la même manière, puis l'\alg librement engendrée par ce \moz.
\\
Si l'on ne veut pas invoquer un résultat aussi \gnlz, on peut simplement constater que les procédures de calculs dans $\kGa$ avec 
$\Gamma=_\mathrm{MC}\scp{truc}{muche}$ sont identiques à celles dans $\gA=_{\gk\mathrm{-\algsz}}\scp{truc}{muche}$.
 
%%%%%%%%%%%%%%%%%%%%%%%%%%%%%%%%%%%%%%%%%%%%%%%%%%%%%%%%%%%%%%%%%%%%%%%%%%%

%:sinotenglish
\sinotenglish{

%: correction exoMonfinMon

\exer{exoMonfinMon}
\emph{(Monoïde fini monogène)}\\
\emph{1.} On a $x^{5}=x^{7}=x^{5+2}$, donc $x^{5}=x^{5+2n}$ et  $x^{6}=x^{6+2n}$
 pour tout $n\in\NN$.
\\
On vérifie que l'ensemble $\lrb{0..6}$ muni de la loi d'addition correspondante est bien un \moz. Donc $\Gamma=\so{1,x,\dots,x^{6}}$ admet exactement $7$ \elts distincts.  
\\
En outre $x^{6}=x^{12}$ donc $e=x^{6}$ est un \idmz.
\\
La \klg $\kGa$ est isomorphe à $\aqo\kGa e\times \aqo\kGa {1-e}$.
\\
On a $\aqo\kGa e\simeq\aqo\kX{X^{6},X^{7}-X^{5}}=\aqo\kX{X^{6},X^{5}}=\aqo\kX{X^{5}}$.
\\
Et $\aqo\kGa {1-e}\simeq\aqo\kX{X^{6}-1,X^{7}-X^{5}}=\aqo\kX{X^{2}-1}$, car $\alpha=\ov X$ est \iv (ou si l'on préfère $\pgcd({X^{6}-1,X^{7}-X^{5})}={X^{2}-1}$). 

\emph{2.} Dans le cas \gnlz, on recopie la \dem du cas particulier. 
\\
La suite $x^{n}$ est périodique de période $r$ à partir du rang $m$:
$x^{m+k}=x^{m+k+nr}$
pour~$k\in\lrb{0..r-1}$ et pour tout~$n\in\NN$.
\\
On vérifie que l'ensemble $\lrb{0..m+r-1}$ muni de la loi d'addition correspondante  

\snic{(a,b)\mt\formule{a+b&\hbox{ si  }a+b<m,\\[.3em]
m+(a+b-m \mod r)&\hbox{ sinon,}}}

est bien un \moz. Donc $\Gamma=\so{1,x,\dots,x^{m+r-1}}$ admet exactement $m+r$ \elts distincts.
\\
Pour l'\idmz, on considère l'unique \elt de la forme $\ell r$ dans $\lrb{m..m+r-1}$. On a alors $x^{\ell r}=x^{(\ell+1) r}$ et donc $x^{\ell r}=x^{2\ell r}$. Soit $e$ cet \idmz. On pourra se rappeler  que $e=x^{Nr}$ pour tout $N\geq \ell$.
\\  
Dans $\aqo\kGa {1-e}$, $x$ est \iv donc $x^{r}=1$. Et ceci implique bien $e=1$. Ceci montre l'\iso $\aqo\kGa {1-e}\simeq \aqo\kX{X^{r}-1}$. 
\\
Si l'on préfère:  $\gen{X^{\ell r}-1,X^{m}(X^{r}-1)}=\gen{X^{r}-1}$.\\  
Enfin $\aqo\kGa {e}=\aqo\kX{X^{\ell r},X^{m}(X^{r}-1)}=\aqo\kX{X^{m}}$ car $m\leq \ell r$.
\\
En bref $\kGa\simeq \gk[\epsilon]\times \gk[\alpha]$ avec $\epsilon^{m}=0$ et $\alpha^{r}=1$.

\rem La seule chose un peu étonnante est que le résultat final, qui semble évident une fois constatée la périodicité de la suite $x^{n}$ à partir du rang $m$, demande un petit effort pour être démontré. \eoe

\emph{3.} Puisque l'ensemble $\sotq{x^{n}}{n\in\NN}$ est fini, la suite $x^{n}$ est périodique  à partir d'un certain rang, qui est le premier $m$ tel que $x^{m}$ soit égal à l'un de ses successeurs stricts. Si $m=0$ on obtient le groupe cyclique d'ordre $r$.
Dans le cas contraire, on retrouve la description du \mo donné 
au point~\emph{2.}

%%%%%%%%%%%%%%%%%%%%%%%%%%%%%%%%%%%%%%%%%%%%%%%%%%%%%%%%%%%%%%%%%%%%%%%%%%%
%: correction exoMonIdm

\exer{exoMonIdm}
\emph{(Idempotent dans un \mo commutatif)} \\
\emph{1a.} On note d'abord que  $e\Gamma$ est un \mo avec $e$ pour \elt neutre. L'application naturelle $\psi:e\Gamma\to \Gamma_1$, $ex\mt \ov {ex}=\ov x$  est la composée de deux morphismes pour la multiplication, donc l'est aussi. C'est un morphisme de \mos car $\ov e=\ov 1$. 
Par ailleurs  l'application
$\imath:\Gamma\to e\Gamma$, $x \mt ex$ est un morphisme de \mosz.
On a donc un unique morphisme $\theta: \Gamma_1\to e\Gamma$ 
vérifiant $\theta (\ov x)=ex$ pour tout $x$. 
Il reste à vérifier que $\psi$ et $\theta$ sont des bijections réciproques.
\\
Les deux \egts $\psi(\theta (\ov x))=\ov x$ et $\theta (\psi(ex))=ex$ résultent des \dfnsz.

\rem On aurait pu aussi montrer directement que $\imath$ résout le \pb \uvl
du passage au quotient par la relation $e=1$ dans le \mo $\Gamma$.
\eoe

\emph{1b.} Comme dans le point \emph{2} de l'exercice \ref{exoAlgMon}, 
ll s'agit ici d'un résultat \gnl d''\alg \uvlez. Dans les deux cas 
on considère la \klg obtenue librement à partir du \mo $\Gamma$ tout en imposant la relation $e=1$. On peut effectuer ces deux constructions 
dans l'ordre que l'on veut. L'important ici est que la relation~\hbox{$e=1$} puisse s'écrire dans le langage des \mosz.
\\
Si l'on ne veut pas invoquer un résultat aussi \gnlz, on peut simplement constater que les procédures de calculs dans $\kGa $ modulo $e=1$   sont identiques à celles dans~$\gk[\Gamma_{e=1}]$.

\emph{1c.} Si $x\in e\Gamma$, on écrit $x=ez$, donc $ex=e^{2}z=ez=x$.
Inversement %l'\egt 
$ex=x$ implique $x\in e\Gamma$.

\smallskip 
\emph{2a.} Il suffit de voir que la plus fine relation d'\eqvcz~$\sim$ sur $\Gamma$ qui vérifie~\hbox{$ex\sim ey$} pour tous $x$, $y\in\Gamma$ est compatible avec la multiplication.
Or  la relation $a\sim b$ est obtenue en liant  $a$ à $b$ par des suites de couples  $(a,x)$, \dots, $(u,v)$, \dots, $(z,b)$ où chacun des couples $(u,v)$ est certifié
par une \egt $u=_\Gamma v$ ou par deux \egtsz~$u=_\Gamma eu'$ \hbox{et $v=_\Gamma ev'$}. Par \recu sur la longueur de la suite on en déduit que   $a\sim b$ \ssi $a=_\Gamma b$ ou $\exists a',b'\,(a=_\Gamma a'e,\,b=_\Gamma b'e)$. Dans les deux cas on voit que pour un $x$
arbitraire, on a $ax\sim bx$. Notez que la \dem n'utilise pas que~$e\Gamma$ soit  détachable ou $\Gamma$ discret. 

\emph{2b.} Notons $\psi$ le morphisme naturel considéré $\kGa\to \aqo{\gk[\Gamma'_e]}{e}$. Il est composé de deux morphismes surjectifs et il annule $e$. Pour voir que c'est un morphisme de passage au quotient par $e$ il faut maintenant vérifier que $\Ker\psi\subseteq \gen{e}$.
\\
Soit $x=\sum_{\gamma\in I} u_{\gamma}\gamma$ (avec $I$ une partie finiment énumérée de $\Gamma$) un \elt de~$\kGa$ vérifiant $\psi(x)=0$.
Notons $\ov\gamma$ l'image de $\gamma$ dans $\Gamma'_e$. 
\\
On sait que
$\sum_{\gamma\in I} u_{\gamma}\ov\gamma\in \ov e\, \gk[\Gamma'_e]$, i.e.

\snic{\sum_{\gamma\in I} u_{\gamma}\ov\gamma=\ov e \sum_{\eta\in J} v_{\eta}\ov\eta= \sum_{\eta\in J} v_{\eta}\ov {e\eta}.}

La conclusion est claire si $\Gamma$ est discret  (donc $e\Gamma$ détachable), ce qui est souvent le cas dans les applications.
Donnons néanmoins la \dem dans le cas \gnlz, un peu plus sophistiquée.
L'\egt $\sum_{\gamma\in I} u_{\gamma}\ov\gamma=\sum_{\eta\in J} v_{\eta}\ov {e\eta}$ est calculée à partir de l'\egt de $\Gamma'_e$. Dans chacune des deux sommes on a le droit de regrouper des termes selon l'\egt de $\Gamma'_e$,
de supprimer un terme affecté d'un \coe nul
et l'on doit aboutir à deux sommes formellement identiques.
Or on a vu que 

\snic{\ov a=\ov b \hbox{ dans } \Gamma'_e \iff \big(a=b \hbox{ ou }\exists a',b'\,(a=a'e,\,b=b'e)\big)\hbox{ dans } \Gamma.}

Donc, pour produire l'\egt $\sum_{\gamma\in I} u_{\gamma}\ov\gamma=\sum_{\eta\in J} v_{\eta}\ov {e\eta}$ dans $\gk[\Gamma'_e]$, les $u_{\gamma}\gamma$  
ou bien doivent être égaux à des termes $u_{\gamma} \gamma 'e$, ou bien
doivent se réduire à $0$ dans~$\kGa$ après des regroupements utilisant des \egts dans~$\Gamma$.
Cela montre bien que~$x\in e\Gamma$.

\emph{2c.}    Calcul \imdz.

\emph{2d.}    Calculs \imdsz. Pour calculer un produit $\gamma_1\cdot\gamma_2$ 
dans $\gk[\Gamma\setminus e\Gamma]$ on teste \hbox{si $\gamma_1\gamma_2\in e\Gamma$} dans $\Gamma$.
Si la réponse est positive  $\gamma_1\cdot\gamma_2=0$, sinon $\gamma_1\cdot\gamma_2=\gamma_1\gamma_2$. 

%%%%%%%%%%%%%%%%%%%%%%%%%%%%%%%%%%%%%%%%%%%%%%%%%%%%%%%%%%%%%%%%%%%%%%%%%%%
%: correction exoMonExemple

\exer{exoMonExempleabcd}
\emph{(Exemple d'étude d'un \mo fini)}

\emph{1.} L'\egt $x^a = y^b$ élevée à la puissance $d$ 
donne $x^{ad} = y^{bd}= x^{bc}$, i.e. $x^{m}=x^{m+\Delta}$.
Même chose pour $y$.
Tout \elt de $\Gamma$ s'écrit donc $x^{k}y^{\ell}$ avec $k,\ell\in\lrb{0..\Delta-1}$. Ainsi~$\Gamma$ est borné. 
Si $m=0$, $x$ et $y$ sont \ivsz, donc $\Gamma$ est un groupe.
\\
Si $m>0$, l'exercice \ref{exoMonfinMon}
nous dit que $x^{\ell \Delta}$ est \idm si $\ell \Delta\geq m$.
On a le même résultat avec $y$, et comme $x^a = y^b$,
on a $x^{\ell \Delta}=x^{a\ell\Delta}=y^{b\ell\Delta}=y^{\ell \Delta}$.  On obtient donc deux fois le même \idmz, $e=x^{\ell \Delta}$.
\\
Pour voir que $e\neq 1$ si $m>0$ on peut passer au quotient en imposant $x=y$ \hbox{et $x^{2}=x$}. On vérifie que le \mo quotient est $\scp{x}{x^{2}=x}$ dans lequel $e\neq 1$.

\emph{2.} En forçant $e=1$ on force l'\ivt de $x$ et $y$ donc de tous
les \elts de~$\Gamma_1$. Ainsi $\Gamma_1=_\mathrm{GAb}\scp{x,y}{x^a = y^b,x^c = y^d}$, où l'indice GAb est mis pour \gui{groupe abélien}. En notation additive on considère donc le conoyau dans $\ZZ^{2}$ de la matrice $M=\Cmatrix{.4em}{a&c\cr b&d}$. Ainsi $\Gamma_1$,
que l'on peut identifier à~$e\Gamma$, est un groupe abélien fini tel que décrit dans l'énoncé. 

\emph{3.}  On sait que le sous-\mo de $\Gamma$
engendré par $x$  est décrit par
une seule \eqn $x^{m_x}=x^{m_x+r_x}$. L'entier $r_x$ divise $\Delta$ et peut être calculé
en examinant la réduction de Smith pour la matrice $M$.
\\
En fait, des \mlrs de colonnes sur $M$ nous donnent

\snic{M'=\Cmatrix{.3em}{*&\fraC \Delta h\cr h& 0}$ avec $h=\pgcd(b,d).}

\snii
Ainsi l'ordre  $r_x$ de $\ov x= ex$ dans $\Gamma_1\simeq e\Gamma$ est égal
à  $ \Delta/{\pgcd(b,d)}$.
\\
De même l'ordre $r_y$  de $\ov y= ey$ dans $\Gamma_1$ est égal
à  $ \Delta/{\pgcd(a,c)}$.
\\ Si l'un des deux pgcd est égal à $1$, on obtient que $\Gamma_1$
est engendré par $\ov x$ ou $\ov y$.
\\
Quant à l'entier $m_x$, il est $\leq m$ et peut être déterminé en calculant les formes normales (avec le calcul expliqué dans le 
point~\emph{4}) de
$x^{k}$ et $x^{k+r_x}$ pour les valeurs successives de $k$.
Il semble difficile de donner une réponse absolument \gnle
à la question posée.

\smallskip 
\rem Définissons les deux intervalles $I_0$, $J_0$:

\snic {
I_0=\lrb{0..c-1},\quad J_0=\lrb{0..b-1} .
}

\snii
L'auteur de l'exercice a constaté sur tous les exemples qu'il a traités que
$\cB\setminus e\cB$ est égal à $\sotQ{x^iy^j}{ i \in I_0,\  j \in J_0 }$ ce qui revient à dire que l'on a une réunion disjointe
$\cB=e\Gamma\cup \sotQ{x^iy^j}{ i \in I_0,\  j \in J_0 }$, ou encore

\snic {
\cB = \sotQ{ ex^iy^j}{ i\in\lrb{0..r_x-1},\ j\in\lrb{0.. r_y-1} } \ \cup\
\sotQ{x^iy^j}{ i \in I_0,\  j \in J_0 }
}

(le nombre  $\Delta=\#(e\Gamma)$ est en \gnl strictement plus petit que $r_xr_y$).\\ 
Il serait intéressant de
prouver cette \prt du \mo $\Gamma$.
\eoe

\smallskip 
\emph{4a.} Vérification laissée \alecz.

\emph{4b.} Pour chaque réécriture, le \mom réécrit est strictement
plus petit pour l'ordre monomial $\preceq$. Or vue la \dfn de $\preceq$ via un produit lexicographique,
une suite strictement décroissante pour $\preceq$ s'arrête en un temps fini.

\emph{4c.} Le \sys de réécriture mis au point comporte seulement deux règles et il est \emph{\lot confluent}
au sens suivant: si un \mom $x^{k}y^{\ell}$ peut être réécrit selon chacune des deux règles, sous formes  $x^{k'}y^{\ell'}$ et $x^{k''}y^{\ell''}\!$, alors les deux \moms réécrits peuvent être eux-mêmes réécrits (éventuellement en plusieurs étapes)
en deux \moms identiques. Ici, c'est particulièrement simple
car lorsque les deux règles peuvent s'appliquer séparément, elles
le peuvent successivement et elles commutent.\\
Ceci dit, quand deux \moms sont-ils égaux dans $\Gamma$? 
\\
Notons $\sim$
la plus fine relation d'\eqvc dans $\NN^{2}$ qui soit compatible avec l'addition et qui vérifie $(a,0)\sim(0,b)$ et $(c,0)\sim(0,d)$.
On a alors $(k_{1},\ell_{1})\sim(k_{2},\ell_{2})$ \ssi  $x^{k_1}y^{\ell_1}=x^{k_2}y^{\ell_2}$ dans $\Gamma$. On voit donc que $\sim$ est la clôture \smq et transitive de la relation \gui{$(k,\ell)$ se réécrit en $(k',\ell')$
selon l'une des deux règles de réécriture}.
\\
C'est un résultat bien connu que pour un \sys de réécriture confluent,
 deux \elts sont \eqvs selon la relation $\sim$ associée \ssi ils peuvent
se réécrire (éventuellement en plusieurs étapes) en deux termes égaux.
Voir par exemple la proposition II.3 page 52 dans {\sc Lalement R.} \emph{Logique Réduction Résolution.} Masson 1990.
\\ 
Cela implique donc ici que deux \moms finaux (i.e. qui ne peuvent pas être réécrits) formellement distincts sont deux \elts distincts  dans $\Gamma$. Ainsi on peut qualifier les formes finales  de formes normales.
\\
Et l'on a bien une bijection $(k,\ell)\mt x^{k}y^{\ell}$ de $\lrb{0..a-1}\times \lrb{0..d-1}$ sur $\Gamma$.   

\emph{5.} De manière \gnlez, lorsqu'un \mo commutatif $\Gamma$ est fini,
l'exercice \ref{exoMonIdm} montre que pour un \idm $e$, on a
l'\egt  $\#(\Gamma'_{e})+\#(\Gamma_{e=1})=1+\#(\Gamma)$, car on a des bijections
$\Gamma_{e=1}\longleftrightarrow e\Gamma$ et $\Gamma'_{e}\longleftrightarrow (\Gamma\setminus e\Gamma)\cup \so e$.

\emph{6.} Soit $f=x^{k}y^{h}$ (avec $k\in\lrb{0,a-1}$ et  $h\in\lrb{0,d-1}$)
un \idm de $\Gamma$. \hbox{Si $f\neq 1$}, on a $k$ ou $h$ $>0$. 
Alors pour $N$ assez grand $f=f^{N\Delta}=x^{kN\Delta}y^{hN\Delta}=e$
d'après le point \emph{2} de l'exercice \ref{exoMonfinMon}. 
%%%%%%%%%%%%%%%%%%%%%%%%%%%%%%%%%%%%%%%%%%%%%%%%%%%%%%%%%%%%%%%%%%%%%%%%%%%
%: correction exoSyst_abcd

\exer{exoSyst_abcd} \emph{(\'Etude d'un \sys \zedz, intersection de deux courbes planes affines)}

\emph {1.}
L'\egt $x^a = y^b$ élevée à la puissance $d$ 
donne $x^{ad} = y^{bd}= x^{bc}$, donc~$x$ est entier sur $\gk$ (car $ad > bc$). Même chose pour $y$. Notons que $x^{ad}-x^{bc}=x^{m}(x^{\Delta}-1)$.
\\
Lorsque $\gk$ est un \cdi contenant le groupe $U_{\Delta}$ des racines $\Delta$-èmes de l'unité, les zéros de $\gA$ dans $\gk$ sont d'une part~$(0,0)$, avec une multiplicité qui reste à déterminer, et d'autre part des couples~$(\xi,\zeta)$ avec $\xi$ et $\zeta$ dans $U_\Delta$ reliés par les relations
$\xi^{a}=\zeta^{b}$ et $\xi^{c}=\zeta^{d}$. Ce sont des zéros isolés. Ils sont simples si la \cara de $\gk$ ne divise pas $\Delta$.  On établira la réciproque plus loin.
Ce sont là \gui{tous les zéros de~$\gA$} au sens où il n'y en a pas d'autres
dans une extension quelconque de $\gk$. 

\emph{2.} Tout est écrit dans l'énoncé. 

\emph{3.} Si $e=x^{N}=y^{N}$ dans $\Gamma$, alors  $\gen{x,y}^{2N}=\gen{e}$ dans $\kGa$. 
Localiser en $1+\gen{x,y}$ est la même chose que localiser en $1+\gen{x,y}^{2N}$  
(lemme~\ref{lemLoc1+a}).
Enfin inverser les \elts de $1+\gen{e}$ avec $e$ \idm revient à inverser $1-e$
(fait \ref{fact.loc.idm}).

\emph{4.} D'après l'exercice \ref{exoMonIdm}, $\gA_1 \simeq \gk[\Gamma_1]$. D'après l'exercice \ref{exoMonExempleabcd},   $\Gamma_1=e\Gamma$ est un groupe abélien fini  d'ordre~$\Delta$. On peut voir $\gk[\Gamma_1]$
comme le \kmo de base $e\Gamma$. 
En s'inspirant de
l'exercice \ref{exoGroupAlgebra}, il est alors facile de déterminer la matrice
tracique~$T$ de $\gA_1$ dans la base $e\Gamma=\Gamma_1$: 

\smallskip 
\centerline{on a $\Tr(g_1g_2) = \# (\Gamma_1) = \Delta$ si
$g_1g_2 = 1_G$, et $\Tr(g_1g_2) = 0$ sinon,} 

\smallskip 
donc $T = \Delta P_\sigma$, où $P_\sigma$ est la
matrice de la permutation~$\sigma:g\mt g^{-1}$ de~$\Gamma_1$. Et le \discri de $\gA_1$ dans la base $\Gamma_1$ \hbox{est $\det(T)
= \vep(\sigma)\Delta^\Delta$}.
\\
En particulier lorsque $\Delta\in\gk\eti$, l'\alg $\gA_1$ est \ste sur
$\gk$. C'est le cas lorsque $\gk$ est un \cdi dans lequel $\Delta.1_\gk\neq 0$.

\emph{5.} On cherche les zéros projectifs éventuels du système sur la droite $Z=0$. On a déjà déterminé  les zéros dans le plan affine $\gk^{2}$ et la somme de leurs multiplicités et égale à $ad$. 
On note ensuite que puisque $ad>bc$, on a $a>  b$ ou $d> c$.

Si   $a>  b$ et $d> c$, on a $F=X^{a}-Y^{b}Z^{a-b}$ et $G=X^{c}Z^{d-c}-Y^{d}$. En faisant $Z=0$, on trouve le \sys $X^{a}=Y^{d}=0$, donc aucun zéro projectif. Même chose si l'une des deux in\egts est large.
La somme des multiplicités des zéros dans~$\PP^{2}(\gk)$ est donc $ad=\deg(f)\deg(g)$.

Si $a>b$ et $d<c$, on a $F=X^{a}-Y^{b}Z^{a-b}$ et $G=X^{c}-Y^{d}Z^{c-d}$.
Ceci donne pour $Z=0$, $X^{a}=X^{c}=0$ et correspond au zéro projectif $(0:1:0)$. En déshomogénéisant avec $Y=1$ 
on trouve le \sys $x^{a}=z^{a-b}$, $x^{c}=z^{c-d}$, qui admet $(0,0)$
comme zéro de multiplicité $\inf(a(c-d),(a-b)c)=ac-ad$ \hbox{(car $ad>bc$)}.
La somme des multiplicités des zéros dans $\PP^{2}(\gk)$ est donc égale à~$ad+(ac-ad)=ac=\deg(f)\deg(g)$.
 
Si $a<b$ et $d>c$, on trouve de même le zéro projectif $(1:0:0)$ avec la multiplicité~\hbox{$bd-ad$} et la somme des multiplicités des zéros dans $\PP^{2}(\gk)$
est encore égale à $\deg(f)\deg(g)$.

}
%: fin sinotenglish

%%%%%%%%%%%%%%%%%%%%%%%%%%%%%%%%%%%%%%%%%%%%%%%%%%%%%%%%%%%%%%%%%%%%%%%%

% fin des solutions d'exos

%%%%%%%%%%%%%%%%%%%%%%%%%%%%%%%%%%%%%%%%%%%%%%%%%%%%%%%%%%%%%%%%%%%%%%%%%
%%%%%%%%%%%%%%%%%%%%%%%%%%%%%%%%%%%%%%%%%%%%%%%%%%%%%%%%%%%%%%%%%%%%%%%%%
%%%%%%%%%%%%%%%%%%%%%%%%%%%%%%%%%%%%%%%%%%%%%%%%%%%%%%%%%%%%%%%%%%%%%%%%%
%%%%%%%%%%%%%%%%%%%%%%%%%%%%%%%%%%%%%%%%%%%%%%%%%%%%%%%%%%%%%%%%%%%%%%%%%
%%%%%%%%%%%%%%%%%%%%%%%%%%%%%%%%%%%%%%%%%%%%%%%%%%%%%%%%%%%%%%%%%%%%%%%%%
%:   ---- Section*{references}-----------
\Biblio

\Llec trouvera sans doute un peu arbitraire
notre volonté de donner à l'anneau trivial toutes
les \prts imaginables, notamment à travers notre utilisation
d'une version affaiblie de la négation
(cf. la note \ref{footnoteNegation} \paref{footnoteNegation}).
Nous espérons \la* convaincre de l'utilité pratique
d'une telle convention par les exemples.
Sur le bon usage de l'anneau trivial, voir~\cite[Richman]{Ri88}.

La \gui{preuve par Azumaya} du lemme de la liberté locale \ref{lelilo}
 est extraite  de la preuve du \tho d'Azumaya III.6.2 dans \cite{MRR},
pour le cas qui nous occupe ici.
Autrement dit, nous avons donné le contenu
\gui{matriciel} de la preuve du lemme de la liberté locale dans \cite{MRR}.

Les courbes monomiales (exemple \paref{exlcourbemonomiale}) sont traitées  dans \cite{Kun}, chapitre V,
exemple~3.13.f.

Les anneaux décomposés jouent
un rôle important dans la théorie classique des \alos henséliens
par exemple dans les ouvrages \cite{Ray} ou \cite{Laf}.

Un \algb est parfois appelé \gui{ring with many units} dans la
littérature de langue anglaise. 
Les \algbs ont notamment été étudiés dans \cite[Estes \& Guralnick]{EG}. 
D'autres \gui{anneaux avec beaucoup d'unités} sont apparus bien avant, sous la terminologie \gui{unit-irreducible rings} (voir par exemple \cite{vdK}). Ce sont les anneaux $\gA$ pour lesquels est vérifiée la \prt suivante: si deux
\pols de $\AX$ représentent un \ivz, alors leur produit représente aussi un \ivz. Ont \egmt été introduits les  anneaux \gui{primitifs} ou \gui{strongly U-irreducible} qui sont les anneaux pour lesquels est vérifiée la \prt suivante: tout \pol primitif représente un \ivz. Ce sont des \algbs particuliers. Dans la \dem du fait \ref{factLocNagata} on a montré qu'un localisé de Nagata est toujours \gui{primitif}.%
%:HHH index
\index{primitif!anneau ---}\index{anneau!primitif}

Concernant le localisé de Nagata $\gA(X)$, vu le fait \ref{factLocNagata}
et les bonnes \prts des \algbsz,
il n'est pas étonnant que cet anneau  joue un rôle crucial  pour 
la solution uniforme des \slis avec paramètres sur un \cdi et 
plus \gnlt pour les calculs uniformes \gui{en temps raisonnable} sur des anneaux
commutatifs arbitraires
(voir~\cite[Díaz-Toca\&al.]{DiGL,DiGLQ}).

\newpage \thispagestyle{CMcadreseul}
\incrementeexosetprob

%:        %%%%%%%%%%%%%%%%%%%%%%%%%%%%%%%%%%%%
%:        %%%%%%%%%%%%%%%%%%%%%%%%%%%%%%%%%%%%
%---- Chapitre  {Modules ptfs}---chap ptf1----
\chapter{Modules \ptfsz, 2}
\label{chap ptf1}\label{ChapThStrBa}
%--------------------
%\minitoc
\minitoc

\subsection*{Introduction}
\addcontentsline{toc}{section}{Introduction}
%-----------------------------------------

Nous poursuivons ici l'étude des \mptfs commencée dans  le chapitre \ref{chap ptf0}.

Dans la section
\ref{sec ptf loc lib} nous reprenons la question de la
\carn  des \mptfs comme modules \lot libres, \cad du \tho de structure locale.

La section \ref{subsecCalRang} est consacrée à l'anneau des rangs sur~$\gA$.
Dans la théorie usuelle en \clama le rang d'un \mptf est défini
comme une fonction localement constante sur le spectre de Zariski.
Nous donnons ici une théorie \elr du rang qui ne fait pas appel aux \idepsz.

Dans la section \ref{secAppliLocPtf}
nous donnons quelques applications simples du \tho
de structure locale.

La section \ref{secGrassman} est une introduction aux grassmanniennes.

Dans la section \ref{subsecClassifMptfs} nous introduisons le \pb \gnl de
la classification complète des \mptfs sur un anneau $\gA$ fixé.
Cette classification est un \pb fondamental et difficile, qui n'admet pas de
solution \algq \gnlez.

La section \ref{secAppliIdenti}
présente un exemple non trivial pour
lequel cette classification peut être obtenue.

%  Section {sec ptf loc lib}
\section{Les \mptfs sont localement libres}
\label{sec ptf loc lib}\relax
%------------------
Nous reprenons la théorie des \mptfs
après la section~\ref{subsec det ptf}. Nous demandons cependant
\alec d'oublier ce que lui a appris la section \ref{secMPTFlocLib}:
la \carn par les \idfsz, le \tho de structure locale \vref{prop Fitt ptf 2} et
les considérations sur le rang liées aux \idfs
ainsi que le \thrf{corth.ptf.sfio} dont la \dem dépend du \tho
de structure locale.

En fait, tous les résultats des sections \ref{subsec det ptf}
et \ref{sec ptf loc lib} pourraient être obtenus par des arguments de
\lon en des \eco puisque nous avons déjà obtenu le \tho de structure
locale des \mptfs (\thrfs{theoremIFD}{prop Fitt ptf 2}) par des méthodes d'\alg
extérieure.

Nous pensons néanmoins que le point de vue \gui{plus global}
développé dans ce chapitre est intéressant en soi, et, d'une certaine
manière, plus simple, comme le met en évidence la \dem
\elr du \tho matriciel
\ref{th matproj} qui résume (et précise) tous les \thos de structure
antérieurs. Là aussi, l'\alg extérieure
est un outil indispensable, mais il semble mieux utilisé, de façon
moins envahissante.

%:  Subsec{Complements puissances exterieures
\subsect{Compléments sur les
puissances extérieures d'un \mptfz}{Compléments sur les puissances extérieures}
\label{subsec comp ext}\relax
%------------------

Le lemme suivant est \imdz.
%--- Lemma{lem ptf RM2}----------
\begin{lemma}
\label{lem ptf RM2}\relax
Soit $P$ un \Amo libre de rang $h$ et $\varphi\in\End(P)$ un \endo \digz, avec une matrice semblable à $\Diag(\lambda_1, \ldots, \lambda_h)$, alors pour
le \polfon de $\varphi$ on obtient:

\snic{\rF{\varphi}(X)\eqdefi\det(\Id_{P[X]}+X\varphi)=(1+\lambda_1X)\,\cdots\,(1+\lambda_hX).}
\end{lemma}
%--- end-lemma--------

Nous établissons maintenant le résultat crucial.
%:   Propos{prop puissance ext}-
\begin{proposition}
\label{prop puissance ext}\relax \emph{(Puissances extérieures)}\\
Soit $P$  un \mptfz.
%-----------------begin enum------------------
%
\begin{enumerate}
\item La puissance extérieure
$k$-ième de $P$, notée $\Al{k}P$, est aussi un \mptfz.
Si  $P=\Im(F)$ pour  $F\in\GA(\gA)$, le module $\Al{k}P$ est (isomorphe à) l'image de la
\mprn $\Al kF$.
\item Si $\varphi$ est un \endo de $P$, le \polfon $\rF{\Al{k}\!\varphi}(X)$
ne dépend  que de $k$ et du \pol $\rF{\varphi}(X)$.
En particulier, le \polmu de $\Al{k} P$ ne
dépend que de  $k$ et du \polmu de~$P$.
\item
\begin{enumerate}
\item Si $P$ est de rang constant $h< k$, le module $\Al{k}P$ est nul.
\item Si $P$ est  de rang constant $h\ge k$, le module $\Al{k}P$ est de rang
constant~$h \choose k$.
\item Dans ce cas, si $\varphi$ est un \endo dont le \polfon s'écrit
$\rF{\varphi}=(1+\lambda_1X)\,\cdots\,(1+\lambda_hX)$, on~a
$$\rF{\Al{k}\!\varphi}(X)=
\prod_{1\leq i_1<\cdots<i_k\leq h}(1+\lambda_{i_1}\cdots\lambda_{i_k}X) .$$
\end{enumerate}
\item Si une \mprn  $F$ a pour image un \mrc $k$, alors $\cD_{k+1}(F)=0$.
\end{enumerate}
%-----------------end enum------------------
\end{proposition}
%--- end-proposition----------------
%---------begin proof----------
\begin{proof}
\emph{1.} Soient $M$ et $N$ deux \Amos et considérons les premières puissances
extérieures de leur somme directe $M\oplus N$. En examinant le \pb universel
que résout la puissance extérieure $k$-ième d'un module, nous obtenons les
\isos canoniques
%------------begin array-------
$$\arraycolsep2pt\begin{array}{rcl}
\Al{2}(M\oplus N) &\simeq& \Al{2}M \oplus(M\otimes N) \oplus \Al{2}N \\[1mm]
\Al{3}(M\oplus N)&\simeq&\Al{3}M \oplus \big(\,\Al{2}M\otimes N\big)  \oplus
\big(M\otimes \Al{2}N\big)  \oplus \Al{3}N, \end{array}$$
%-------------end array--------------
et plus \gnlt
%------begin equation--eqeqVik-----------
\begin{equation}\label{eqVik}
\Al{m}(M\oplus N)\,\,\simeq\quad
\bigoplus_{k=0}^m \Big(\big(\,\Al{k}M\big)\otimes \big(\,\Al{m-k}N \big)\Big)
\end{equation}
%---------------------end equation--------------
(avec $\Al0M=\gA$ et $\Al1M=M$). En particulier, si $P\oplus Q\simeq \Ae m $,
$\Al kP$ est en facteur direct dans $\Al k\Ae m \simeq\gA^{m \choose k}$.
On voit aussi que si  $P=\Im(F)$ pour une \mprn $F$, $\Al kP$
est (isomorphe à)
l'image de la \mprn $\Al kF$, car cette matrice représente l'\idt sur
$\Al kP$ et $0$ sur tous les autres facteurs de la somme directe.

 \emph{2.} On peut supposer $P = \Im(F)$, où $F \in \GA_n(\gA)$, et $n \ge k$.
\\
On a donc $P \oplus Q = \gA^n$ avec $Q = \Ker(F)$.
L'\endo $\varphi$
se prolonge en un \endo $\varphi_1:\Ae n \rightarrow \Ae n $, nul sur~$Q$,  de matrice~$H$
avec~$FHF=H$, et l'on a $\rF{\varphi}(X)=\rF{\varphi_1}(X)=\det(\I_n+XH)$.
Alors, on voit que~$\Al k\!\varphi_1$ est un
prolongement de  $\Al k\!\varphi$, nul sur
les termes distincts de $\Al kP$ dans la somme directe explicitée dans
la \dem du point \emph{1.} La
matrice de~$\Al k\!\varphi_1$ n'est autre que~$\Al kH$.
On veut donc montrer que $\det\big(\I_{n \choose k}+X\;\Al kH\big)$
ne dépend
que de~$k$ et de~$\det(\I_n+XH)$. Nous sommes donc ramenés au cas d'un
module libre. Et ce cas a été traité dans la
proposition~\ref{propPolCarPuissExt}.
\label{Preuvepoint2prop puissance ext}\relax

 \emph{3.}
Ce point résulte du précédent, puisque le \gui{cas \pro de rang
$k$} peut se déduire du cas \gui{libre de rang $k$}.
Notez que les points \emph{a} et \emph{b} disent  tous deux que lorsque $P$
est de rang constant $h$,  $\Al kP$ est de rang
constant~$h \choose k$ (qui est égal à $0$ si $h<k$).
On ne les a séparés
que dans le but de donner le résultat sous forme plus visible.

 \emph{4.}
Cela équivaut au fait que $\Al{k+1}P$ est nul, ce qui est le point~\emph{3a.}
\end{proof}
%---------end proof----------
%--- rem{rem puissance ext}-

\rems  (Conséquences de la proposition \ref{prop puissance ext}.)
\label{rem puissance ext}\relax

 1) %La convention naturelle est ${n \choose m}=0$ si $m>n$. Le point 3 a)
%est alors un cas particulier du point 3 b).
Posons $\rR{P}(X)=r_0+r_1X+\cdots+r_nX^n$. Chaque $r_hP$
est un \mrc $h$ sur $A[1/r_h]$, ce qui donne pour $k>0$  comme conséquence du point \emph{3}:

\snic{\rR{\Al k(r_hP)}(X)
=X^{h \choose k} \, \hbox{ sur }\,\gA[1/r_h].}

%\sni
En écrivant $P=\bigoplus_hr_hP$ et $\gA=\prod_h \gA[1/r_h]$ on obtient
\[\arraycolsep2pt
\begin{array}{rcl}
 \rR{\Al k\!P}(X) & =  &  r_0+\cdots+r_{k-1}+ r_kX+\cdots+
r_{k+j}X^{k+j \choose k} + \cdots+
r_{n}X^{n \choose k} \\
  & =  & \som_{h=0}^n r_{h}X^{h\choose k} .
\end{array}
\]
On a aussi par convention
$\Al0P=\gA$ et donc aussi $\rR{\Al0\!P}(X)=X$ (pour que la formule précédente s'applique il faut convenir que ${n \choose 0}=1$ pour tout~$n\geq0$).

 2) Si l'on note $\Vi P$ l'\alg extérieure de $P$, \llec montrera par un  calcul analogue que

\snic{\rR{\Vi \!P}(X)= r_0 X +r_{1}X^2 +\cdots+ r_k X^{2^k}  + \cdots+
r_{n}X^{2^n}
.}

%\sni
3) On peut calculer  $\rF{\Al k(\varphi)}$ à partir de  $\rF{\varphi}$ de
la manière suivante. \\
Puisque $\rF{\varphi}(0)=1$ et $\deg(\rF{\varphi})\leq n$,
si $\psi$ est l'\endo de $\Ae n $ ayant pour matrice la matrice compagne $C$  de
$X^n\rF{\varphi}(-1/X)$, on obtient $\rF{\varphi }=\rF{\psi}$. Donc
$$ \rF{\Al k\!\varphi}= \rF{\Al k\!\psi} =
\det \big(\I_{n \choose k}+X\,\Al kC\big)
\eqno\eoq  $$

\medskip Des remarques précédentes on déduit la proposition suivante.
%--- Propos{prop rang ext carac}
\begin{proposition}
\label{prop rang ext carac}\relax
Soit $P$ un \mptfz, et $k\leq h$ deux entiers $>0$.
\Propeq
\vspace{-2pt}
%---------begin item----------
\begin{enumerate}
\item  Le module $P$ est de rang constant $h$.
\item  Le module  $\Vi P$ est de rang constant $2^h$.
\item  Le module  $\Al k P$ est de rang constant $h \choose k$.
\end{enumerate}
%---------end item----------
Avec $h=0$, les \prts 1. et 2. sont \eqvesz.
\end{proposition}
%--- end-proposition----------------

%:  Subsection{subsec rang constant}----
\penalty-2500
\subsec{Cas des modules de rang constant}
\label{subsec rang constant}\relax
%-------------------
%--- Th{th rg const loc libre}--
\begin{theorem}
\label{th rg const loc libre}\relax
Soit $P$ un \Amrc $h$ avec~$n$ \gtrsz,
(isomorphe à l') image d'un \prr $F\in \GAn(\gA)$. Alors
les~$n\choose h$ mineurs principaux $(s_i)$ d'ordre $h$ de $F$ vérifient:
%---------begin item----------
\begin{itemize}
\item [--] $\sum_i s_i = 1$,  et
\item [--] chaque $\gA_{s_i}$-module  $P_{s_{i}}$ est libre de rang $h$, la matrice $F$
vue comme matrice à \coes dans $\gA_{s_i}$ est semblable à la \mprn standard $\I_{h,n}$.
\end{itemize}
%---------end item----------
\end{theorem}
%--- end-theorem-----------------
%---------begin proof----------
\begin{proof}
La somme des mineurs
principaux $s_i$ d'ordre $h$ de $F$ est égale à $1$ puisque~\hbox{$\det(\I_n+XF)=(1+X)^h$}.\\
Par ailleurs, puisque tout mineur d'ordre $h+1$ est nul (proposition \ref{prop puissance ext}), on peut appliquer le lemme de la liberté \ref{lem pf libre}
à chaque localisé~$P_{s_{i}}$, lequel est isomorphe à l'image de la
matrice $F$ vue comme matrice à \coes dans~$P_{s_{i}}$ (d'après la
proposition~\ref{propPtfExt}).
\end{proof}
%---------end proof----------
%--- rem{rem rg const loc libre}

\rem
\label{rem rg const loc libre}\relax
Dans le \tho précédent, il se peut que $s_i$ soit nilpotent pour certaines
valeurs de $i$, donc que $\gA_{s_i}$ soit trivial. Le fait de ne pas exclure ces
\lons nulles est inévitable lorsque l'on ne dispose pas d'un test pour
savoir si un \elt de $\gA$ est ou n'est pas nilpotent. Ceci justifie la convention
naturelle donnée dans la remarque \paref{conven rgc}.
\eoe

%:  Subsection{subsec cas general}------
\subsec{Cas général}
\label{subsec cas general}\relax
%------------------
%--- Theorem{th ptf loc libre}--
\begin{theorem}
\label{th ptf loc libre}\relax
Soit $P$ un \Amo \ptf avec $n$ \gtrsz. Alors pour chaque \idm $\ide_h(P)$ il
existe $n\choose h$ \elts $(s_{h,i})$ de $\gA$ avec les \prts suivantes:
%---------begin item----------
\begin{enumerate}
\item [--] $\sum_i s_{h,i} = \ide_h(P),$
\item [--] chaque $\gA_{s_{h,i}}$-module  $P_{s_{h,i}}$ est libre de rang $h$.
\end{enumerate}
%---------end item----------
En particulier, pour tout \mptf à $n$ \gtrsz, il existe $2^n$ \eco
$v_\ell$ tels que chaque $P_{v_\ell}$ soit libre.
\end{theorem}
%--- end-theorem-----------------
%---------begin proof----------
\begin{proof}
On localise d'abord en inversant $\ide_h(P)$ pour se ramener au \tho \ref{th rg const loc
libre}. On localise ensuite un peu plus conformément à ce dernier \thoz.
Le fait \ref{factLocCas} concernant les \lons successives s'applique.
\end{proof}
%---------end proof----------

Le \tho suivant résume les \thos \vref{th rg const loc libre} et 
\ref{th ptf loc libre}, et la réciproque donnée par le \plgrf{plcc.cor.pf.ptf}.

%--- Theorem{thptfloli}----
\begin{theorem}
\label{thptfloli}
Un \Amo  $P$ est \ptf \ssi il existe des \eco $s_1$, \ldots, $s_\ell$ tels que chaque
$P_{s_i}$ est libre sur $\gA_{s_i}$.
Il est \pro de rang $k$ \ssi il existe des \eco $s_1$, \ldots, $s_\ell$ tels que
chaque $P_{s_i}$ est libre de rang $k$ sur $\gA_{s_i}$.
\end{theorem}
%--- end-theorem-----------------------------------------

Une forme pratique du \thoz~\ref{th ptf loc libre} est sa forme matricielle. 
%:HHH supprime 
%En
%réalité, les \thos \vref{th decomp ptf}, \vref{th rg const loc libre} et
%\ref{th ptf loc libre} n'ont pas de signification calculatoire supérieure à
%celle du \tho matriciel suivant, qui est même légèrement plus précis.

%:     Theorem{th matproj}--------
\begin{theorem}
\emph{(Forme matricielle explicite des \thos \ref{th.ptf.loc}
et~\vref{th.ptf.idpt})}~\label{th matproj}\label{propPTFDec}\\
 Soit  $\gA$ un anneau, $F\in \Mn(\gA)$ avec $F^2=F$ et $P$ le \mptf
image de $F$ dans $\Ae n $.  On définit les \elts $r_h$ de $\gA$ pour
$h\in\lrb{0..n}$ par les \egtsz:

\snic{\rR{P}(1+X):=\det(\In+XF),\quad \rR{P}(X)=:r_0+r_1X+\cdots+r_nX^n.
}

%\sni
On a les résultats suivants.
%---------begin item----------
\begin{enumerate}
\item  La famille $(r_h)_{h=0,\ldots,n}$ est un \sfio de $\gA$.
\item  Pour $h\in\lrb{0..n-1}$ et $u$  mineur d'ordre $h+1$ de $F$, on
a $ r_hu=0 $.
\item  Si les  $t_{h,i}$ sont les mineurs principaux d'ordre $h$  de  $F$,
on obtient en posant~$s_{h,i}=r_ht_{h,i}$:
%---------begin item----------
\begin{enumerate}\itemsep0pt
\item [--] la somme (pour $h$ fixé) des $s_{h,i}$ est égale à $r_h$,
\item [--] chaque $\gA_{s_{h,i}}$-module  $P_{s_{h,i}}$ est libre de rang $h,$
\item [--] la matrice $F$  est semblable sur $\gA_{s_{h,i}}$ à la matrice
$\I_{h,n},$
\item [--] les $s_{h,i}$ sont \comz, \prmt $\som_{h,i}s_{h,i}=1$.
\end{enumerate}
%---------end item----------
\end{enumerate}
%---------end item----------
\end{theorem}
%--- end-theorem-----------------

%--- rem{rem th matproj}-----

\rem
\label{rem th matproj}\relax
Le \tho \ref{th matproj} résume les \thos \ref{th decomp ptf}, \ref{th rg const
loc libre} et  \ref{th ptf loc libre} qui l'ont précédé. 
Il est même légèrement plus précis. Il n'est donc pas
inintéressant d'en donner une preuve purement matricielle qui concentre toutes
les preuves précédentes, d'autant plus qu'elle est particulièrement \elrz.

%---------begin proof----------
\begin{Proof}{Preuve matricielle du \tho matriciel. }
\emph{1.}
Cela résulte de $\rR{P}(1)=1$ (évi\-dent) et de
$\rR{P}(XY)=\rR{P}(X)\rR{P}(Y)$ qui se voit comme suit:
%--------------------begin array---------------
$$\preskip3pt\postskip3pt
\arraycolsep2pt\begin{array}{rclc}
\rR{P}(1+X)\rR{P}(1+Y)&  = &\det(\In+XF)\det(\In+YF)   & =     \\[1mm]
\det\big((\In+XF)(\In+YF)\big)&  = & \det(\In+(X+Y)F+XYF^2)  & =     \\[1mm]
\det(\In+(X+Y+XY)F)&  = & \rR{P}\big((1+X)(1+Y)\big).  &
\end{array}$$
%---------------------end array--------------

\noi \emph{2.} La matrice $r_hF$ a pour \polfon $\det(\In+r_hXF)$. Dans l'anneau
$\gA_{r_h}$, on a $1=r_h$ et

\snic{\det(\In+r_hXF)=\det(\In+XF)=\rR{P}(1+X)=(1+X)^h.}

%\sni
%\looseness-1
En se plaçant dans l'anneau $\gA_{r_h}$ on est donc ramené
à démontrer le point~\emph{2} pour le cas 
où $r_h=1$ et $\det(\In+XF)=(1+X)^h$,
ce que nous supposons  désormais.
Nous devons montrer que les mineurs d'ordre $h+1$ de $F$  sont tous nuls.
Les mineurs d'ordre~$h+1$ sont les \coes de la matrice~$\Al{h+1}F=G$.
Puisque~$F^2=F$, on a aussi~$G^2=G$.
Par ailleurs, pour n'importe quelle matrice carrée~$H$,
le \polcar de~$\Al{k}H$
ne dépend que de~$k$ et du \polcar de~$H$ (proposition
\ref{propPolCarPuissExt}).
Appliquant ceci pour calculer le \polcar de~$G$, nous pouvons remplacer~$F$ par la
matrice $\I_{h,n}$ qui a même \polcar que $F$. Comme
la matrice $\Al{h+1}\I_{h,n}$ est
nulle, son \polcar est $X^{h+1 \choose n}$,
donc, par Cayley-Hamilton, la matrice $G$ est nilpotente, et comme elle est
\idmez, elle est nulle.

 \emph{3.} Résulte de \emph{1}, \emph{2} et du lemme de la liberté \ref{lem pf libre}.
\end{Proof}
%---------end proof----------

%%%%%%%%%%%%%%%%%%%%%%%%%%%%%%%%%%%%%%%%%%%%%%%%%%%%%%%%%%%%%%%%%%%%%%%%%%%
%:  subsec{secDecEqdim} ---------

\subsec{Modules de rang constant: quelques précisions}
\label{secDecEqdim}
%--------------------------------

%--- Propos{prop rgconstant local}
Les deux résultats suivants sont désormais faciles et nous laissons leur
preuve en exercice.

\begin{proposition}
\label{prop rgconstant local 2}\label{prop rgconstant local}\relax
\emph{(Modules \prcsz)}\\
Pour un \Amo $P$  \propeq
%-----------------begin enum------------------
\begin{enumerate}
\item $P$ est \prcz~$h$.
\item \label{i2prop rgconstant local}\relax
Il existe des \eco $s_i$ de $\gA$ tels que chaque $P_{s_i}$ est libre de
rang $h$ sur~$\gA_{s_i}$.
\item $P$ est \ptf et pour tout \elt $s$ de $\gA$, si $P_s$ est libre sur $\gA_s$,
il est de rang~$h$.
\item $P$ est \pfz, $\cF_h(P)=\gen{1}$ et  $\cF_{h-1}(P)=0.$
\item $P$ est isomorphe à l'image d'une \mprn de rang~$h$.
\end{enumerate}
%-----------------end enum------------------
\end{proposition}
%--- end-proposition----------------
En outre, si  $P$ engendré par $n$ \elts le nombre des \eco \ncrs dans
le point \emph{\ref{i2prop rgconstant local}.} est
majoré par  ${n\choose h}.$

%--- Propos{prop sfio unic}
\begin{proposition}
\label{prop sfio unic}{\em  (Localisés de rang constant et unicité du \sfioz)}\\
Soit $P$ un \Amo \ptfz. Posons $r_h=\ide_h(P)$. Soit $s$ un \elt de $\gA$.
%-----------------begin enum------------------
\begin{enumerate}
\item Le localisé $P_s$ est \pro de rang $h$ \ssiz $r_h/1=1$ dans~$\gA_s$, \cad si $r_hs^m=s^m$ dans $\gA$ pour un certain
exposant~$m$.
\item Si $s$ est un \idmz, cela signifie que $r_h$ divise $s$, ou encore que $1-r_h$ et~$s$ sont deux \idms \ortsz.
\item Enfin, si $ (s_0,\ldots,s_n) $ est un \sfio tel que chaque $P_{s_h}$ soit de rang $h$ sur $\gA_{s_h}$, alors $ r_h=s_h$
pour chaque $h\in\lrbzn$.
\end{enumerate}
%-----------------end enum------------------
\end{proposition}

Dans la proposition qui suit nous faisons  le lien entre notre \dfn et la \dfn usuelle (en \clamaz)
d'un module \pro de rang $k$. La preuve de cette \eqvc n'est
cependant pas \cov (et ne peut pas l'être).

%--- Prop{prop ptfrangconstant}----
\begin{proposition}\label{prop ptfrangconstant}\relax
Soit $k$ un entier naturel, $P$ un \mptf sur un anneau $\gA$ non trivial et
$\fa$ un \id contenu dans $\Rad\gA$.  Alors \propeq
%---------begin item----------
\begin{enumerate}
\item [1.\phantom{$^*$}] $P$ est de rang $k$, i.e., $\rR{P}(X)=X^k$
\item [2.$^*$] Pour tout \idema $\fm$ de $\gA$, l'espace vectoriel obtenu à
partir de $P$ en étendant les scalaires au corps résiduel $\gA/\fm$ est
de dimension~$k$.
\item [3.\phantom{$^*$}]  $\rR{P}(X)\equiv X^k$ modulo $\fa[X]$.
\end{enumerate}
%---------end item----------
\end{proposition}
%--------------
%----begin{proof----------
\begin{proof}
D'un point de vue classique, l'implication  \emph{2} $\Rightarrow$ \emph{3} est \imdez; il suffit de se rappeler que l'intersection des \idemas est
le radical de Jacobson de $\gA$.
Notez que d'un point de vue \cofz, la condition \emph{2} est a priori trop faible,
par manque d'\idemasz.
\\
Par ailleurs, \emph{1} implique trivialement  \emph{2} et \emph{3}. 
\\
Réciproquement, si
$\rR{P}(X)=X^k$ modulo $\fa[X]$, puisque les \idms sont toujours isolés
(lemme~\ref{lemIdmIsoles}), l'\egt a lieu dans $\AX$.
\end{proof}
%----end{proof----------

%:      Th {propRgConstant2}----
\begin{theorem}
\label{propRgConstant2} \emph{(Modules de rang constant $k$ comme sous-modules
de~$\Ae k$)}\\
Supposons que sur $\Frac \gA$ tout \mrc $k$ soit libre.
Alors tout $\gA$-\mrc $k$ est isomorphe à un sous-module
de $\Ae k$.
\end{theorem}
%--- end-theorem----------------------------------------
%-----------------begin proof------------------
\begin{proof}
D'après le lemme \dlg \ref{propIsoIm}, on peut supposer que le module est
image d'un \prr $F\in\GAn(\gA)$ de rang $k$ et qu'il existe une matrice~$P$
dans $\GLn(\Frac \gA)$ telle que $PFP^{-1}=\I_{k,n}$. On \hbox{a $P=Q/a$} 
\hbox{avec $Q\in\Mn(\gA)$} et~$a\in\Reg\gA$; ainsi $\det Q=a^n\det P$ est aussi \ndz dans~$\gA$.  
On définit une matrice~$Q_1$ par:
$$
Q\cdot F\,=\,\I_{k,n}\cdot Q\,=
\blocs{.8}{.6}{.8}{.6}{$\I_k$}{$0$}{$0$}{$0$}
\cdot  Q\,=\,
\blocs{1.4}{0}{.8}{.6}{$Q_1$}{}{$0$}{}\;.
$$
Or l'image de $Q\cdot F$ est isomorphe à l'image de $F$ parce que $Q$ est
injective, et l'image de $\I_{k,n}\cdot Q$ est clairement isomorphe à
l'image de~$Q_{1}= Q_{1..k, 1..n}$.
\end{proof}
%-----------------end proof------------------

\rem 1) Le \tho précédent
s'applique aux anneaux \qis et plus \gnlt à tout anneau $\gA$
%:HHH ci dessous Frac A plutot que Frac Ared
tel que $\Frac\gA$ soit \zedz,
ou même simplement \lgbz.
C'est par exemple le cas des anneaux réduits 
\noes \cohs  \fdis (voir le \pbz~\ref{exoAnneauNoetherienReduit}).  
On démontre en \clama que pour tout anneau \noez~$\gA$, $\Frac\gA$ est \rdt \zedz,
donc on peut lui appliquer le \thoz. On ne connaît pas d'analogue \cof
de ce \thoz. 
\\
2) Pour plus de précisions concernant le cas $k=1$ voir le
\thref{propRgConstant3}.

%:  Subsection{cas generique}----------
\subsec{Cas générique} \label{subsec cas generique}\relax
%-------------------

Qu'est-ce que nous appelons le cas \gnqz, concernant un module \pro
à $n$ \gtrsz?
Nous considérons l'anneau
$$\Gn =\ZZ[(f_{i,j})_{i,j\in\lrbn}]/\cGn ,\label{NOTAGN}$$
où les $f_{i,j}$ sont des \idtrsz,  $F$ est la matrice $(f_{i,j})_{i,j\in\lrbn}$ et $\cGn$
est l'\id défini par les $n^2$ \syzys obtenues en écrivant $F^2=F$.
\`A \coes dans cet anneau $\Gn$, nous avons la matrice $F$ dont l'image dans
$\Gn^{n}$ est ce qui mérite d'être appelé
{\em  le module \pro générique à $n$ \gtrsz}.

Reprenons les notations du \tho \ref{th matproj} dans ce cas particulier. 
\\
Dire que
$r_hr_k=0$ dans $\Gn$ (pour $0\leq h\not= k\leq n$) signifie que l'on a une appartenance 
$$
r_h({F})r_k({F})
\in \cGn\qquad (*)
$$
 dans
l'anneau $\ZZ[(f_{i,j})_{i,j\in\lrbn}]$.
Cela implique une \ida qui permet d'exprimer cette
appartenance. Cette \ida est naturellement valable dans tous
les anneaux commutatifs.
Il est donc clair que si l'appartenance $(*)$ est vérifiée dans le cas
\gnqz, elle implique  $r_hr_k=0$ pour n'importe quelle \mprn sur un anneau commutatif arbitraire.

La même chose vaut pour les \egts $r_hu=0$ lorsque $u$ est un mineur d'ordre
$h+1$.

En résumé: si le \tho \ref{th matproj} est vérifié dans le cas
\gnqz, il est vérifié dans tous les cas.
Comme souvent, nous constatons donc que des \thos importants d'\alg
commutative ne font rien d'autre qu'affirmer l'existence de certains types
particuliers d'\idasz.

%   section{L'anneau des  rangs}-
%:HHH ajout de [L'anneau des rangs généralisés $\HO(\gA)$
\section[L'anneau des rangs généralisés $\HO(\gA)$]
{\texorpdfstring{Le semi-anneau $\HOp(\gA)$, et l'anneau des rangs généralisés $\HO(\gA)$}
 {H0+(A),  et l'anneau des rangs généralisés H0(A)}}
\label{subsecCalRang}
%-----------------

Pour un module libre, en passant du rang $k$ au \polmu $X^k$, on passe de la
notation additive à la notation multiplicative.
Pour un \mptf \gnlz, on peut considérer en sens inverse
un \gui{rang \gnez} du module, qui est  le loga\-rithme (purement formel) en base
$X$ de son \polmuz. Bien qu'il s'agisse là d'un simple jeu de notation, il
s'avère que les calculs avec les rangs en sont facilités.
Expliquons comment cela fonctionne.

%  subection{Le semi-anneau des rangs}
\subsection*{Le semi-anneau des rangs}
\label{SemiAnneauRangs}

Nous dirons qu'un \pol $R(X)=r_0+r_1X+\cdots +r_nX^n$ est \ixc{multiplicatif}{polynome@\pol ---} lorsque $R(1)=1$ et $R(XY)=R(X)R(Y)$.
Il revient au même de dire que les $r_i$ forment un \sfioz, ou encore que $R(X)$ est le \polmu d'un \mptfz.

%:        Notation{notaHO+}
\begin{notation}
\label{notaHO+}\relax
{\rm   On note $\HOp (\gA)$ l'ensemble des classes d'\iso des modules
%:HHH petit ajout
%quasi libres sur $\gA$, et $[P]$ la classe d'un tel module $P$
quasi libres sur $\gA$, et $[P]_{\HOp (\gA)}$ (ou $[P]_\gA$, ou même
$[P]$) la classe d'un tel module $P$
% fin ajout
dans  $\HOp (\gA)$. L'ensemble $\HOp (\gA)$ est muni d'une structure
de \emph{semi-anneau\index{semi-anneau}\footnote{Ceci signifie que la structure
est donnée par une addition, commutative et associative, une multiplication,
commutative, associative et distributive par rapport à l'addition,
avec un neutre $0$ pour l'addition et un neutre $1$ pour la multiplication.
Par exemple $\NN$ est un semi-anneau.}
pour les lois héritées de $\oplus$ et $\otimes$: $[P\oplus Q]=[P]+[Q]$ \hbox{et
$[P\otimes Q]=[P]\cdot[Q]$}. Pour un \idm $e$ on notera aussi $[e]$
à la place de~$[e\gA]$, lorsque le contexte est clair.
L'\elt neutre pour la multiplication est~$[1]$.}
}
\end{notation}
%--- end-notation-----------------------------------------

Tout module quasi libre $P$ est isomorphe à un unique module\footnote{
On a aussi (exercice \ref{exoSfio}) $P\simeq e_1\gA \oplus  
\cdots \oplus e_n\gA$ avec $e_k=\sum_{j=k}^nr_j$, en outre $e_k$ divise~$e_{k+1}$ pour $1 \le k < n$.}

\snic{(r_1\gA)\oplus (r_2\gA)^2 \oplus \cdots \oplus
(r_n\gA)^n,}

%\sni
où les $r_i$ sont des \idms \ortsz, puisqu'alors $\ide_{i}(P)=r_i$.
On a donc $[P]=\sum_{k=1}^n k\,[r_k]$
et son \polmu est  

\snic{\rR{P}(X)=r_0+r_1X+\cdots +r_nX^n}

%\sni
avec $ r_0 =  1 - (r_1 +\cdots + r_n)$.

Mais alors que $\rR{P\oplus Q}=\rR{P}\rR{Q}$, on a $[P\oplus Q]=[P]+[Q]$: ceci assure le passage de la notation multiplicative à la notation additive.
Ainsi le \gui{logarithme en base $X$} du \pol multiplicatif $r_0+r_1X+\cdots +r_nX^n$
est défini comme l'\elt $\sum_{k=1}^n k\,[r_k]$ de $\HOp (\gA)$.

%:     Definition{defiRang}
\begin{definition}\label{defiRang}
Si $M$ est un \Amo \ptf on appelle \emph{rang (\gnez) de $M$} et l'on note
$\rg_\gA(M) $ ou $\rg(M)$ l'unique
\elt de $\HOp(\gA)$ qui a le même \polmu que 
lui.\index{rang!(généralisé) d'un \mptfz}
\end{definition}
\label{NOTAHO+}\relax

Ainsi si $\rR{M}(X)=r_0+r_1X+\cdots +r_nX^n$, alors $\rg(M)=\sum_{k=1}^n k\,[r_k]$.

Le module nul est \care par $\rg(M)=0$ (\thrf{th ptf sfio}).

Si $\gA$ est non trivial, alors $[1]\neq[0]$ et $\NN$ s'identifie
au sous-semi-anneau de~$\HOp(\gA)$ engendré par $[1]$ au moyen de l'injection $n\mapsto n\,[1]$.
La \dfn ci-dessus n'entre donc pas en conflit avec la notion de rang pour les \mrcsz,
définie auparavant.

Notez aussi que lorsque $\gA$ est trivial on a $\HOp (\gA)={0}$: ceci
est bien conforme à la convention selon laquelle le module nul sur l'anneau trivial
a pour rang n'importe quel entier, puisque dans $\HOp (\gA)$, $k=0$, ou
si l'on préfère, les deux \polmus $1$ et $X^k$ sont égaux sur l'anneau trivial.

\medskip
\rem Une règle de calcul pratique portant sur les rangs est
 la suivante:

\snic{[r]+[r']=[r+r'] \;\;\;\mathrm{si}\;\;\; rr'=0,}

%\sni
\cad plus \gnlt
\begin{equation}\label{eqsomHO}
[r]+[r']=[r\vu r']+ [r\vi r']=[r\oplus r']+2\,[r\vi r']
\end{equation}
où les lois
$\vu,$ $\vi$ et $\oplus$ sont celles de l'\agB des \idms de l'anneau:
$r\oplus r'=r+r'-2rr',$ $r\vu r'=r+r'-rr'$ et $r\vi r'=rr'$. Notez que
les deux \idms $r\oplus r'$ et $r\vi r'$ sont  \ortsz, de somme~$r\vu r'$, et que la signification
de l'\egt (\ref{eqsomHO}) est donnée  par les \isos suivants
$$
{r\gA\oplus r'\gA \;\simeq\; (r\vu r')\gA\oplus (r\vi r')\gA \;\simeq\;
(r\oplus r')\gA\oplus \big((r\vi r')\gA\big)^2.} \eqno \eoq
$$

%  subection{Notation exponentielle}
\subsection*{Notation exponentielle}
\label{ExpoRangs}

Remarquons que $a^n$ est le résultat de l'évaluation
du \pol multiplicatif~$X^n$
au point $a$: $X^n(a)=a^{\log_X(X^n)}$.

Ainsi, pour un \pol multiplicatif $R(X)=\sum_{k=0}^ne_kX^k$,
dont le logarithme en base $X$ est l'\elt $r=\sum_{k=0}^nk\,[e_k]$, on
adopte les notations légitimes suivantes: 
\begin{itemize}
\item $a^r=\sum_{k=0}^ne_ka^k=R(a)$,
\item et pour un \Amo $M$,  $M^r=\bigoplus_{k=0}^ne_kM^k$. 
\end{itemize}
Ce n'est pas une fantaisie: on a bien 
\begin{itemize}
\item $a^{r+r'}=a^r a^{r'}$, $a^{r r'}=(a^r)^{r'}$,
\item  $M^{r+r'}\simeq M^r \times M^{r'}$ et $M^{r r'}\simeq M^r \otimes M^{r'}\simeq (M^r)^{r'}$,
\end{itemize}
pour  $r,r'$ arbitraires dans $\HOp\gA$.

%  subection{Symetrisation}
\subsection*{Symétrisation}
\label{SymetrisationRangs}\rdb

Le \mo additif  $\HOp(\gA)$ est régulier: au choix d'après
le lemme de McCoy (corolaire \ref{corlemdArtin}), ou l'un
des deux \thos d'unicité \ref{prop unicyc} et \rref{prop quot non iso},
ou enfin le point \emph{\iref{th ptf sfio item reg}.} du \thrf{th ptf sfio}.

Le semi-anneau $\HOp (\gA)$  peut donc être considéré comme
un sous-semi-anneau de l'anneau obtenu en le symétrisant.
Cet anneau s'appelle l'\ixx{anneau des rangs}{(\gnesz) de \mptfsz} sur $\gA$,
et l'on le note  $\HO(\gA)$.
%:        Notation{notaHO}
%\begin{notation}
\label{notaHO}
%{\rm   On note $\HO(\gA)$ l'anneau des rangs (\gnesz) des \mptfs sur $\gA$.
%En abrégé on parlera de l'\emph{anneau des rangs sur~$\gA$}.
%}
%\end{notation}
%--- end-notation-----------------------------------------

Tout \elt de  $\HO(\gA)$ s'écrit sous forme $\som_{k\in J} k \,[r_k]$
où les $r_k$ sont des \idms deux à deux orthogonaux et $J$
est une partie finie de $\ZZ\setminus\so{0}$.
\\
L'écriture est unique au sens suivant:
si  $\som_{k\in J} k \,[r_k]=\som_{k\in J'}
k \,[r'_k],$ \hbox{alors  $r_k=r'_k$} si $k\in J\cap J',$ et les autres sont nuls.

%:HHH supprime, on parle de HO un peu partout
%Nous reparlerons de $\HO$ \paref{secGKO}.

%  subection{MultiplicationRangs}
\subsection*{Multiplication des rangs}
\label{MultiplicationRangs}

On a défini sur $\HOp\gA$ une multiplication, comme la loi héritée du
produit tensoriel. Ceci implique que pour deux \idms $e$ et $e'$
\hbox{on a $[e]\cdot [e']=[ee']$}. Les autres calculs de produits s'en déduisent par distributivité. D'où le fait suivant.

%:     Fact{fact0ProdRangs}
\begin{fact}\label{fact0ProdRangs}
L'\elt $1$ est le seul \iv de $\HOp(\gA)$. 
\end{fact}
\begin{proof}
Si $r=\sum_k k[r_k]$ et $s=\sum_k k[s_k]$, alors $rs=\sum_kk(\sum_{i,j,ij=k}[r_is_j])$. Par unicité de l'écriture, si $rs=1=1[1]$, alors
$r_1s_1=1$ donc $r_1=s_1=1$.
\end{proof}

On peut se demander quelle est la loi correspondante sur
les \pols multiplicatifs: \llec se convaincra qu'il s'agit de la loi

\snic{\big(R(X),R'(X)\big)\mapsto R\big(R'(X)\big)=R'\big(R(X)\big).}

%\sni
On a aussi le fait suivant qui découle de la proposition \ref{prop Produit tensoriel} à venir.

%:     Fact{factProdRangs}
\begin{fact}\label{factProdRangs}
Si $P$ et $Q$ sont deux \mptfsz, alors $P\otimes Q$ est un \mptf et
$\rg(P\otimes Q)=\rg( P)\cdot\rg(Q)$.
\end{fact}

%   subection{Relation d'ordre}
\subsection*{Relation d'ordre sur les rangs}
\label{Relation-dordreRangs}\relax

La relation d'ordre naturelle associée à la structure de \mo de $\HOp\gA $
est décrite dans la proposition suivante.

%:     PropDefi{def rang inferieur}--
\begin{propdef}
\label{def rang inferieur}~
\begin{enumerate}
\item Pour $s$, $t\in\HO\gA$ on définit $s\leq t$ par: $\exists r\in\HOp\gA $, $s+r= t$.
\item Cette relation fait de $\HO\gA$ un \emph{anneau 
ordonné}\index{anneau!ordonné}\footnote{Cela signifie que $\geq$ est une relation d'ordre partiel \emph{compatible}
avec les lois $+$ et $\times$: 
\begin{itemize}
\item $1\geq 0$,
\item $x\geq 0$ et $y\geq 0$
impliquent $x+y$ et $xy\geq 0$,
\item $x\geq y \Leftrightarrow x-y\geq 0$.
\end{itemize}
}, et $\HOp\gA$ est la partie positive de $\HO\gA$.
\item \label{i3def rang inferieur}\relax
Soient $P$  et $Q$ des \mptfsz, \propeq
\begin{enumerate}
\item  $\rg(P) \leq \rg(Q)$.
\item $\rR{P}$ divise $\rR{Q}$ dans $\AX $.
\item $\rR{P}$ divise $\rR{Q}$ dans $\BB(\gA)[X]$.
\item Pour tout $s\in \gA$, si $P_s$ et $Q_s$ sont libres, alors le rang de $P_s$
est inférieur ou égal à celui de $Q_s$.
\item Pour tous  $k>i, \; \; \ide_k(P)\cdot \ide_i(Q)=0.$
\item Pour tout $ k,\; \;  \ide_k(P)\cdot\sum_{i\geq k}\ide_i(Q)=\ide_k(P).$
\end{enumerate}
\end{enumerate}
\end{propdef}
%--- end-Propdefi -----------

\exl
Supposons que
$P\oplus R=Q$ et que l'on connaisse les rangs de
$P$ et $Q$, on demande de calculer le rang de $R$.\\
On a $\rg{P}=\sum_{i=0}^n i \,[r_i]$ et $\rg{Q}=\sum_{j=0}^m j \,[s_j]$.
On écrit
\[\arraycolsep2pt
\begin{array}{ccccccc}
 \rg{P} & = \big(\sum_{i=0}^n i \,[r_i]\big)\big(\sum_{j=0}^m \,[s_j]\big) & = &
 \sum_{i,j} i \,[r_is_j] & \leq\\[1mm]
 \rg{Q} & = \big(\sum_{j=0}^m j \,[s_j]\big)\big(\sum_{i=0}^n \,[r_i]\big) & = &
 \sum_{i, j} j \,[r_is_j]. \\
\end{array}
\]

\penalty-2500
Les $r_is_j$ forment un \sfio et l'on obtient par soustraction, sans avoir
à réfléchir, les \egts 

\snac{\quad \rg(Q)-\rg(P)=\rg(R)=\som_{i\leq j}
(j-i) \,[r_is_j] =\som_{k=0}^m k \big(\som_{j-i=k}\,[r_is_j]\big). \hfill\eoq}

\medskip
Dans la suite, nous laissons définitivement tomber le mot \gui{\gnez}
lorsque nous parlons de rang d'un \mptfz.

\medskip
\rem En \clamaz,  $\HO(\gA)$ est souvent défini
comme l'anneau des fonctions localement constantes
(i.e., continues) de $\SpecA$ vers~$\ZZ$. Un commentaire plus
détaillé sur ce sujet se trouve \paref{comHOclassique}.\eoe

%  subection{Autres utilisations du rang}
\goodbreak
\subsection*{Autres utilisations du rang}
\label{AutresRangs}

% :     notation
\begin{notations}\label{notaRangs} ~
\begin{enumerate}
\item Si $\varphi\in\Lin_{\gA}(P, Q)$ avec $P$, $Q$ \ptfsz, et si $\Im\varphi$ est facteur direct dans $Q$ on notera $\rg\varphi$ pour $\rg(\Im\varphi)$.
\item Si $p(X)$ est un \pol pseudo unitaire de $\AX$, on peut définir
son degré $\deg p$ comme un \elt de $\HOp (\gA)$.
\end{enumerate}
\end{notations}

Pour le point \emph{1.} on a
$\Ker \varphi$ qui est facteur direct dans $P$ et l'on obtient les
\gnns d'\egts bien connues dans le cas
des \evcs sur un corps discret: 

\snic{\rg(\Ker \varphi)+\rg \varphi= \rg P\;\hbox{ et} \;\rg(\Ker \varphi)+\rg Q=\rg(\Coker \varphi)+\rg P\,.}

%:HHH petit ajout
%\sni
En outre, en cas de modules libres, et pour un rang $r\in\NN$, on retrouve bien la notion de rang d'une matrice définie en~\ref{defRangk}. 

 Concernant le point \emph{2.}, notons que pour deux \pols pseudo unitaires~$p$ et~$q$
on a l'\egt $\deg pq=\deg p+\deg q$.
\\
Cette notion de degré s'étend de manière naturelle aux \pols \lot \mons
définis dans l'exercice~\ref{exoPolLocUnitaire}.

%  Sec{Qqs applis}---secAppliLocPtf--
\section{Quelques applications du \tho de structure locale}
\label{secAppliLocPtf}
%-----------------------------------------

Dans cette section nous envisageons des résultats concernant les \mptfs
et des \alis entre ceux-ci.

Vu le \tho de structure locale pour les \mptfsz, et puisque le \deter 
et les \pols corrélatifs se comportent bien par
changement d'anneau de base (fait \ref{fact.det loc}), on a de manière presque automatique tous les résultats souhaités au moyen de la \dem donnée dans l'encadré suivant. 

\Grandcadre{\DebP\,\, Dans le cas de modules libres, le résultat est facile à
établir. \quad \eop}

Nous ne le mentionnerons pas toujours dans cette section.

NB:  si dans l'hypothèse figure une \ali \lnl entre deux modules
différents, on est ramené par le \tho de structure locale
au cas d'une \ali \nlz.

\smallskip La \dem fonctionne chaque fois que le résultat à
établir est vrai \ssi il est vrai après \lon en des \ecoz.

\medskip
\rem Si l'on doit démontrer un résultat qui, dans le cas
de modules libres,
se résume à des \idas on peut en outre, supposer que les
endomorphismes sont \digsz.
%(comme nous l'avons fait implicitement \paref{Preuvepoint2prop puissance ext}
%dans la preuve du point 2
%de la proposition \ref{prop puissance ext} 
%lorsque nous avons invoqué la proposition~\ref{propPolCarPuissExt}).
L'argument est ici différent du
\tho de structure locale. C'est que pour vérifier une \ida
il suffit de le faire sur un ouvert de Zariski de l'espace des
paramètres, et une matrice \gnq est \dig d'après la proposition~\ref{propMatGenDiag}.

%:  Subsec{Trace}---------------
\subsect{Trace d'un \endo et nouvelle écriture du \polfon}{\Pol fondamental}
\label{subsec trace}\relax
%------------------
Rappelons que si $M$ et $N$ sont deux \Amosz, on note $\theta_{M,N}$ l'\ali  canonique

\snic{\theta_{M,N}:M\sta\otimes_\gA N\to\Lin_\gA(M,N), \; (\alpha \otimes y)\mapsto
(x\mapsto \alpha(x)y).}

%\sni
Rappelons aussi les résultats suivants (fait \ref{factDualPTF} et
proposition~\ref{propAliPtfs}).

\mni{\it Soit $P$ un \mptfz. \label{lem iso cano}\relax
%-----------------begin enum------------------
\begin{enumerate}
\item  $\theta_{P,N}$ est un \iso  de $P\sta\otimes_\gA N$ dans $\Lin_\gA(P,N)$.
\item  $\theta_{N,P}$ est un \iso  de $N\sta\otimes_\gA P$ dans $\Lin_\gA(N,P)$.
\item  L'\homo canonique $P\rightarrow P^{\star\star}$  est un \isoz.
\item  L'\homo canonique

\snic{\varphi \mapsto \tra{\varphi}~;~\Lin_\gA(N,P)\,\to\, \Lin_\gA(P\sta,N\sta),}

%\sni
est un \isoz.
\end{enumerate}
%-----------------end enum------------------
}

\medskip Si $P$ est un \mptfz, rappelons que la \emph{trace} de l'\endo
$\varphi$ de $P$ (notée $\Tr_P(\varphi)$) est le \coe en $X$ du \polfon
$\rF{\varphi}(X)$.
Elle peut être \egmt définie à partir de l'\ali  naturelle

\snic{\tr_P:P\sta\otimes_\gA P\rightarrow \gA~:~\alpha \otimes y\mapsto \alpha(y)
,}

%\sni
et de l'\iso  canonique $\theta_{P}:P\sta\te_\gA P\to\End(P)$, comme suit:
%% PERSO
\perso{Conflit de notation lorsque l'on  a affaire à une \Alg \stfe?
Peut-être faut-il le signaler? En cas de besoin $\Tr_P$ devrait être
noté $\Tr_{P,\gA}$, cela ferait
une claire distinction entre $\Tr\iBA$ et $\Tr_{\gB,\gA}$.
En fait $\Tr\iBA(a)=\Tr_{\gB,\gA}(\mu_{\gB,a})$.}
%% Fin PERSO
$$\Tr_P=\tr_P \circ\, {\theta_{P}}^{-1}.
$$
(\Llec pourra constater que les deux \dfns coïncident dans le cas d'un module libre, ou se reporter au fait \ref{factMatriceEndo}.)

Lorsque $P$ et $Q$ sont \ptfsz, la trace permet aussi de définir une dualité
canonique entre $\Lin_\gA(P,Q)$ et $\Lin_\gA(Q,P)$ au moyen de l'application
bi\lin $(\varphi,\psi)\mapsto \Tr(\varphi \circ\psi)=\Tr(\psi\circ\varphi)$. Cette
dualité peut aussi
être définie par l'\iso canonique $(P\sta\otimes_\gA Q)\sta\simeq P\otimes_\gA
Q\sta$.

%--- Propos{trace polfon}--
\begin{proposition}
\label{prop trace polfon}\relax
Soit $\varphi$  un \endo d'un \mptf $P$ à $n$ \gtrsz.
Les \coes du \polfon de $\varphi$ sont donnés par
%---------begin $$--------
$$
\rF{\varphi}(X) = 1+\som_{h\in\lrbn} \Tr\big(\,\Al h\varphi\,\big) \,X^h\,.
$$
%---------end $$----------
\end{proposition}
%--- end-proposition----------------

%--- Propos{prop trace surjective}
\begin{proposition}
\label{prop trace surjective}\relax
Si $P$ est un  \mptf fidèle, alors l'\Ali trace $\Tr_P:\End(P)\to \gA$
est surjective.
\end{proposition}
%--- end-proposition----------------
%NB: Rappelons que $P$ est fidèle \ssi $\rg P \geq 1$.

%\newpage
%:  Subsec{Produit tensoriel}-------
\subsec{Produit tensoriel}
\label{subsec Produit tensoriel}\relax
%------------------
%--- Proposi{prop Produit tensoriel}
\begin{proposition}
\label{prop Produit tensoriel}\relax
On considère deux $\gA$-\mptfs $P$ et $Q$, $\varphi$ et $\psi$ des
endomorphismes de $P$ et $Q$.  Le module $P\otimes_\gA Q$  est un \mptfz.
%-----------------begin enum------------------
\begin{enumerate}
\item On a l'\egt 

\snic{\det ({\varphi \otimes
\psi})=(\det {\varphi})^{\rg Q}\,(\det {\psi})^{\rg P}
\eqdefi \rR Q(\det \varphi) \rR P(\det \psi).}

\item Le \polfon $\rF{\varphi\otimes \psi}(X)$ de  $\varphi\otimes_\gA \psi$ ne
dépend que de  $\rg({P})$, de~$\rg({Q})$, de $\rF{\varphi},$ et de $\rF{\psi}$.
\item Si
$\rF{\varphi}=(1+\lambda_1X)\,\cdots\,(1+\lambda_mX)$, et
$\rF{\psi}=(1+\mu_1X)\,\cdots\,(1+\mu_nX)$,
on a l'\egt
$\rF{\varphi\otimes \psi}(X) =
\prod_{i,j}(1+\lambda_{i}\mu_{j}X)$.

\item
En particulier,  $\rg(P\otimes Q)=\rg(P)\rg(Q).$
\end{enumerate}
%-----------------end enum------------------
\end{proposition}
%--- end-proposition----------------

Remarquez que la dernière \egt se
réécrit

\snic{\ide_h(P\otimes Q)= \sum_{ jk=h}\,\ide_j(P)\ide_{k}(Q).
}

\smallskip
Notez aussi que la proposition précédente pourrait être démontrée
\gui{direc\-te\-ment} sans passer par le \tho de structure locale, avec une \dem
calquée sur celle qui a été faite pour les puissances extérieures
(proposition~\ref{prop puissance ext}).

%%%%%%%%%%%%%%%%%%%%%%%%%%%%%%%%%%%%%%%%%%%%%%%%%%%%%%%%%%%%%%%%%%%%%%%%%%%
%:  Subsec{Rangs et \alisz}--
\subsec{Rangs et \alisz}
\label{subsecSURJISO}
%------------------

%--- Propos{epi rang constant}
\begin{proposition}
\label{prop epi rang constant}\label{prop mono rang}\relax
Soit $\varphi :P\rightarrow Q$ une \ali entre \mptfsz.
%-----------------begin enum------------------
\begin{enumerate}
\item Si $\varphi$ est surjective, alors $P\simeq \Ker\varphi\oplus Q$. 
Si en outre $\rg(P)= \rg(Q)$, alors~$\varphi$ est un
\isoz.
\item Si $\varphi$ est injective, alors $\rg(P)\leq \rg(Q)$.
\end{enumerate}
%-----------------end enum------------------
\end{proposition}
%--- end-proposition----------

\begin{proof}
Dans le point \emph{2}, il suffit de prouver l'in\egt
après \lon en un \elt $s$ qui rend les
deux modules libres. Comme la \lon préserve l'injectivité,
on peut  conclure d'après le cas des modules libres
(voir le corolaire~\ref{corprop inj surj det} et la remarque qui suit).
\end{proof}

%--- Corollary{corEpiRC}--------
\begin{corollary}
\label{corEpiRC}
Soit $P_1 \subseteq P_2 \subseteq P$ avec $P_1$ facteur
    direct dans $P$. 
    \\
Alors $P_1$ est facteur direct dans $P_2$.
\\
    En conséquence, si les modules sont  \ptfsz, on
    a l'\eqvc
%    Si $P_1\oplus Q_1 = P =P_2\oplus Q_2$ avec $P_1\subseteq P_2$ on~a:

\snic{
 \rg(P_1)=\rg(P_2)\;%\Longleftrightarrow \;\rg(Q_1)=\rg(Q_2)\;
 \Longleftrightarrow
\;P_1=P_2.}

%\sni
 Si en outre  $P_1\oplus Q_1  =P_2\oplus Q_2=P$,
on a les \eqvcs

\snic{
 \rg(P_1)=\rg(P_2)\;\Longleftrightarrow \;\rg(Q_1)=\rg(Q_2)\;
 \Longleftrightarrow
\;P_1=P_2.}
\end{corollary}
%--- end-corollary------------------------------------

%:  Subsec{Formules de transitivité}
\subsec{Formules de transitivité}
\label{subsecTransPtf}

%:--- Notation{notaTraceDetCarAlg}-------
\begin{notation}
\label{notaTraceDetCarAlg}
{\rm
Soit une \Alg $\gB$, \stfe sur $\gA$. Alors on note
$[\gB:\gA]=\rg_\gA(\gB)$.
}
\end{notation}
%--- end-notation-----------------------------------------

Rappelons que d'après le fait \ref{factTransptf},
si  $\gB$ est \stfe sur $\gA$, et si~$P$ est un \Bmo \ptfz,
alors~$P$ est aussi un  \Amo \ptfz.

%:HHH ajouts et modif
Lorsque l'on prend pour $P$ un module quasi libre sur $\gB$, en considérant
son rang sur
$\gA$ cela définit
un \homo 
du groupe additif $\HO\gB$ vers
le groupe additif $\HO\gA$.
Cet \homo est appelé \emph{\homo de restriction} et il est noté~$\rRs{\gB/\!\gA}$.\index{restriction!homomorphisme de ---}
On obtient ainsi un foncteur contravariant d'une sous-catégorie des anneaux commutatifs vers celle des groupes abéliens.
Il s'agit de la catégorie dont les morphismes sont les $\rho:\gA\to\gB$
qui font de~$\gB$ une \asf sur $\gA$.

Par ailleurs, $\HO$ définit un foncteur covariant de la catégorie des anneaux commutatifs vers celle des semi-anneaux, 
puisque par \edsz, un module quasi-libre donne un module quasi-libre.

Comme $\HO(\gC)$ est complètement \care par $\BB(\gC)$ (pour une formulation catégorique, voir l'exercice~\ref{exoAnneaudesrangs}), les
points \emph{1} et \emph{2} du fait suivant décrivent complètement les deux foncteurs dont nous venons de parler. 

%:     Fact{factHORestriction}
\begin{fact}\label{factHORestriction}
Soit $\rho:\gA\to\gB$ une \algz. 
\begin{enumerate}
\item Pour $e\in\BB(\gA)$, on a $\HO(\rho)([e]_\gA)=\big[\rho(e)\big]_\gB$ dans $\HO\gB$. \\
En particulier $\HO(\rho)$ est injectif (resp., surjectif, bijectif) \ssi
la restriction de $\rho$ à $\BB(\gA)$ et  $\BB(\gB)$ est injective (resp., surjective, bijective). 
\end{enumerate}
On suppose maintenant que $\gB$ est \stfe sur $\gA$.
\begin{enumerate}\setcounter{enumi}{1}
\item Pour $e\in\BB(\gB)$, 
$\Rs\iBA([e]_\gB)=\rg_\gA(e\gB)$.
Et $\Rs\iBA(1)=[\gB:\gA]$.
\item Si un \Bmo $P$ est à la fois quasi libre sur $\gA$ et $\gB$, on obtient simplement
%
%\snic{
$\rRs{\gB/\!\gA}([P]_{\gB})=[P]_{\gA}$.
%}
%\item 
%
\end{enumerate} 
\end{fact}
%--------- fin fact ---------------------------------------------- 
\rem Si $\gA$ est connexe et contient $\ZZ$, on peut faire semblant de considérer $\HO(\gA)\simeq\ZZ$ comme un sous-anneau de $\gA.$ Dans le point \emph{2} ci-dessus on voit alors que
{$\Rs\iBA([e]_\gB)=\big[\Tr\iBA(e)\big]_\gA$}
(il suffit de considérer le cas où $e\gB$ est libre en facteur direct d'un libre dans~$\gB$). 
\eoe

Le lemme suivant généralise le \thref{Th.transitivity} 
(qui s'occupait du cas libre).

%:     Lemma{lem1TransPtf}
\begin{lemma}\label{lem1TransPtf} \emph{(Formules de transitivité pour la trace et le \deterz)} \\
Soit $\gB$  une \Alg \stfe
et $P$ un \Bmo \ptfz. Soit $u_\gB:P\to P$ une \Bliz, que nous notons
$u_\gA$ lorsque nous la regardons comme une \Aliz.
Alors on a les \egts fondamentales:

\snic{\det_\gA(u_\gA)=\rN\iBA\big(\det_\gB(u_\gB)\big)\;\hbox{ et }\;\Tr(u_\gA)=\Tr\iBA\big(\Tr(u_\gB)\big).}
\end{lemma}
\begin{proof}
Quitte à localiser en des \eco de $\gA$ nous pouvons supposer
que $\gB$ est un \Amo libre, de rang $k$.
Nous écrivons

\snic{P\oplus N=L\simeq\gB^n\simeq\Ae{ nk},}

\sni(le dernier \iso est un \iso de \Amosz).
%Ceci montre que $P$ est un \Amo \ptfz.
Nous considérons $v=u\oplus\Id_N\in\End_\gB(L)$.
Alors, par \dfn du \deterz, on obtient les \egts $\det_\gB(u_\gB)=\det_\gB(v_\gB)$
et $\det_\gA(u_\gA)=\det_\gA(v_\gA)$. On peut donc appliquer
la formule de transitivité du \thrf{Th.transitivity}.\\
Le raisonnement pour la trace est similaire.
\end{proof}

%:     corollary{lem2TransPtf}
\begin{corollary}\label{lem2TransPtf} Soit $\gA\vers{\rho}\gB$  une \alg \stfez,
 $P$ un \Bmo \ptf et $u_\gB\in\End_\gB(P)$.
\begin{enumerate}
\item $\rC{u_\gA}(X) = \rN_{\gB [X]/\!{\AX }} \big(
\rC{u_\gB}(X)\big)$.
\item  $\rF{u_\gA}(X) = \rN_{\gB [X]/\!{\AX }} \big(
\rF{u_\gB}(X)\big)$.
\item  En particulier, les \polmus de $P$ sur $\gA$ et $\gB$ sont reliés par

\snic{\rR{P_\gA}(X)= \rN_{\gB[X]/\!{\AX }}\big(\rR{P_\gB}(X)\big) }
\item  L'\homo de restriction vérifie

\snic{\rRs{\gB/\!\gA}\big(\rg_\gB(P)\big)=\rg_\gA(P).}
\item  Si $P$ est un \Amo \ptfz, alors

\snic{\rg_\gB\big(\rho\ist(P)\big)=\HO(\rho)\big(\rg_\gA(P)\big),\et
%}
%
%\snic{
\rg_\gA\big(\rho\ist(P)\big)=[\gB:\gA]\,\rg_\gA(P).}

\end{enumerate}
\end{corollary}
\begin{proof}
Les points \emph{1}, \emph{2}, \emph{3} résultent du lemme précédent.
\\
\emph{4.}
Le point \emph{3} nous dit que le \polmu de $P$ sur $\gA$
ne dépend que du \polmu de $P$ sur $\gB$. On peut donc supposer
$P$ quasi libre sur $\gB$ et on applique la \dfn de l'\homo $\Rs\iBA$.
%En bref le point 4 est la réécriture du point 3 en langage additif.
\\
Le point \emph{5} est laissé \alecz.
\end{proof}

Un autre corolaire est donné par le \tho suivant.

%:     theorem{corthTransPtf}
\begin{theorem}\label{corthTransPtf}
Soit   $\gB$ une \Alg \stfe  et  $\gC$ une \Blg \stfez.
Alors  $\gC$ est une \Alg \stfe et

\snic{
[\gC:\gA]=\rRs{\gB/\!\gA}([\gC:\gB]).}

%:HHH precision ci dessous
%\sni
En particulier, si $\HO(\gA)$ s'identifie à un sous-anneau de 
$\HO(\gB)$, et si le rang de $\gC$ sur $\gB$
est un \elt de $\HO(\gA)$, on a

\snic{
[\gC:\gA]= [\gB:\gA]\, [\gC:\gB].}
\end{theorem}

%:  Subsection{pros de rang 1}------
\subsec{Modules \pros de rang 1}
\label{subsec ptf rang 1}\relax
%------------------

%--- Fact{factRANG1}-------
\begin{fact}
\label{factRANG1}
Une matrice $F\in\Mn(\gA)$ est \idme et de rang $1$ \ssi $\Tr(F)=1$ et
$\Al2F=0$.
\end{fact}
%--- end-fact-----------------------------------------
\facile

%--- Proposition{propRANG1}------
\begin{proposition}
\label{propRANG1}
Soit $P$ un $\gA$-\mrcz~$1$.
%-----------------begin enum------------------
\begin{enumerate}
\item Les \homos canoniques

\snic{\gA\rightarrow \End(P)$, $\,a\mapsto \mu_{P,a}\,$  et
$\,\End(P)\rightarrow \gA$, $\varphi \mapsto \Tr(\varphi)}

%\sni
sont deux \isos
réci\-pro\-ques.
\item Pour tout $\varphi \in \End(P)$, on a  $\det(\varphi)=\Tr(\varphi)$.
\item L'\homo canonique
$P\sta\otimes_\gA P\rightarrow \gA$ est un \isoz.

\end{enumerate}
%-----------------end enum------------------
\end{proposition}
%--- end-proposition----------------------------------------

%--- Proposition{th ptrg1}---------
\begin{proposition}
\label{th ptrg1}\relax
Soient $M$ et $N$ deux \Amosz. \\
Si $N\te_\gA M$ est isomorphe à $\gA$, alors $M$
est un module \pro de rang $1$ et  $N$  est isomorphe à~$M\sta$.
\end{proposition}
%--- end-proposition-----------------
%-----------------begin proof------------------
\begin{proof}
Notons $\varphi$ un \iso de $N\otimes_\gA M$  sur $\gA$. Soit $u=\sum_{i=1}^n c_i\otimes a_i$ l'\elt de $N\otimes M$ tel que $\varphi(u)=1$.
On a deux \isos de $N\otimes M\otimes M$ vers $M$, construits à partir de
$\varphi$.

\snic{c \otimes a  \otimes b\mapsto \varphi(c\otimes a)\, b \quad  
\hbox{et} \quad
c \otimes a\otimes b\mapsto \varphi(c\otimes b)\, a .}

%\sni
Ceci donne un \iso  $\sigma :M\rightarrow M$ vérifiant

\snic{\sigma\big(\varphi(c\otimes a)\, b\big)=\varphi(c\otimes b)\,a \;\hbox{ pour tous }c \in N,\; a,b\in M, \hbox{ d'où}}

\snic{\sigma(x)=
\sigma\big(\sum_i \varphi(c_i\te a_i) x\big) =
\sum_i \varphi(c_i\otimes x) a_i,\hbox{ et }
x=\sum_i \varphi(c_i\te x) \sigma^{-1}(a_i).}

%\sni
Ceci montre que $M$ est \ptfz, avec le \syc $\big((\un),(\aln)\big)$, où $u_i=\sigma^{-1}(a_i)$ et
$\alpha_i(x)=\varphi(c_i\te x)$.
\\
De même, $N$ est \ptfz. 
Mais $1=\rg({N\otimes M})=\rg({N})\rg({M})$,
 %les modules  
donc $M$ et $N$  sont  de rang~1 (fait \ref{fact0ProdRangs}).
%\\
Enfin $N\otimes M\sta\te M\simeq N\simeq M\sta$.
\end{proof}
%-----------------end proof------------------

%%%%%%%%%%%%%
\penalty-2500
\section{Grassmanniennes}
\label{secGrassman}

%:--- Subsection {L'anneau générique Gn}
\subsec{Les anneaux génériques
\texorpdfstring{$\Gn$ et $\Gnk$}{Gn et Gn,k}}
\label{subsubsec AGBR}\relax
%------------------

Nous avons défini l'anneau $\Gn=\Gn(\ZZ) =\ZZ[(f_{ij})_{i,j\in\lrbn}]/\cGn$
\paref{subsec cas generique}.

En fait la construction est fonctorielle et l'on peut définir $\Gn(\gA)$
pour tout anneau commutatif $\gA$:  $\Gn(\gA)=\gA[(f_{ij})_{i,j\in\lrbn}]/\cGn\simeq\gA\otimes_\ZZ\Gn.$
Notons~$r_k=\ide_k(\Im\,F)$ où $F$ est la matrice $(f_{i,j})$ dans  $\Gn(\gA)$.

Si nous imposons en outre que le rang soit égal à
$k$, nous introduisons \hbox{l'idéal $\cGnk=\cGn+\gen{1-r_k}$} et nous obtenons l'anneau

\snic{\Gnk=\ZZ[F]\sur{\cGnk}\simeq \Gn[1/r_k]\simeq\aqo{\Gn}{1-r_k}.}

%\sni
Nous avons aussi la version relativisée à $\gA$:
$$
 \Gnk(\gA)=\gA[F]/\cGnk\simeq \Gn(\gA)[1/r_k] \simeq
\gA\otimes_\ZZ\Gnk.
$$

L'anneau $\Gn(\gA)$ est isomorphe au produit  des $\Gnk(\gA).$

\Grandcadre{Dans le paragraphe présent consacré à $\Gnk$ on pose $h=n-k$.}

Si $\gK$ est un corps, l'anneau $\Gnk(\gK)$ peut être considéré comme
l'anneau des \coos de la \vrt affine $\GAnk(\gK)$ dont les
points sont les paires~$(E_1,E_2)$ de sous-espaces de $\gK^n$ vérifiant les \egts
$\dim(E_1)=k$ et~$\gK^n=E_1\oplus E_2$.

En \gmt \agqz, il y a quelques arguments massue pour  affirmer que
l'anneau  $\Gnk(\gK)$ a toutes les bonnes \prts que l'on puisse imaginer,
ceci en relation avec le fait que la \vrt  $\GAnk(\gK)$ est un espace
homogène pour une action du groupe linéaire.

Nous allons retrouver ces résultats \gui{à la main} et en nous affranchissant
de l'hypothèse \gui{$\gK$ est un corps}.

En utilisant les \lons convenables en les mineurs principaux
d'ordre~$k$ de
la matrice $F=(f_{ij})$ (la somme de ces mineurs est égale à 1 
dans~$\Gnk(\gA)$), nous allons établir quelques \prts essentielles du
foncteur~$\Gnk$.

%:   Theorem{thBnk}--------
\begin{theorem} \emph{(Le foncteur $\Gnk$)}
\label{thBnk}
%---------begin enum----------
\begin{enumerate}
\item  \label{thBnk1} Il existe 
des \eco $\mu_i$ de l'anneau $\Gnk(\gA)$ tels que chaque localisé
$ \Gnk(\gA)[1/\mu_i]$  est isomorphe à l'anneau  

\snic{\gA[(X_j)_{j\in\lrb{1..2hk}}][1/\delta]}

%\sni
pour un certain $\delta$ qui vérifie $\delta(\uze)=1$. 
\item \label{thBnka} L'\homo naturel $\gA\rightarrow \Gnk(\gA)$ est injectif.
\item  \label{thBnkf} Si $\varphi \,:\gA\rightarrow \gB$ est un \homoz, le noyau
de $\Gnk(\varphi)$ est engendré par $\Ker\,\varphi$. En particulier, si
$\varphi$ est injectif,  $\Gnk(\varphi)$ l'est \egmtz.
\end{enumerate}
%---------end enum----------
\end{theorem}
%--- end-theorem-----------------

%:     Corollary{corthBnk}
\begin{corollary}\label{corthBnk}
Soit $\gK$ un corps discret et $\gA$ un anneau.
%---------begin enum----------
\begin{enumerate}
\item  \label{corthBnkc} L'anneau  $\Gnk(\gK)$ est intègre, \iclz,
\cohz, \noe ré\-gulier, de \ddk $2kh$.
\item  \label{thBnkd} Si $\gA$ est un anneau
intègre (resp. réduit,  \qiz,  \lsdzz,  normal, 
 \coh \noez,  \coh \noe régulier) il en va de même pour $\Gnk(\gA)$.
\item  \label{corthBnke} La \ddk de  $\Gnk(\gA)$ est égale à celle de
$\gA[X_1,\ldots ,X_{2hk}]$.
\item  \label{corthBnkb} L'anneau $\Gnk=\Gnk(\ZZ)$ est intègre, \iclz, \coh \noez, régulier, de dimension (de Krull) $2k h +1$.
\end{enumerate}
%---------end enum----------
\end{corollary}

\rdb
\comm Nous avons utilisé dans le corolaire la notion d'\anor et celle
 de \ddk
que nous n'avons pas encore définies (voir sections \ref{subsecEntiers} et \ref{secDefConsDimKrull}).
Enfin un \cori est dit \emph{régulier} lorsque tout \mpf
admet une résolution projective finie (pour cette dernière notion voir le \pbz~\ref{exoFossumKumarNori}).%
\index{regulier@régulier!anneau cohérent ---}%
\index{anneau!coherent@cohérent régulier}
%:2018 ajout de anneau cohérent regulier dans l'index 
\eoe

%---------begin proof----------
\begin{Proof}{Esquisse de la \demz. }
Si l'on rend un  mineur principal
d'ordre $k$ de~$F$ \ivz, alors l'anneau $\Gnk(\gA)$ devient isomorphe à un
localisé d'un anneau de \pols sur $\gA$, donc hérite de toutes les
\prts agréables de~$\gA$.
Pour l'intégrité il y a une subtilité en plus, car ce n'est pas une
\prt locale.
\end{Proof}
%---------end proof----------

Nous développons maintenant l'esquisse ci-dessus.
Pour le cas d'un \cdi nous commençons par le résultat suivant.
%:   Lemma{lemPSES1}--------
\begin{lemma}
\label{lemPSES1}
Soit $\gK$ un \cdi et $(E_1,E_2)$ une paire de sous-espaces
\supls de dimensions $k$ et $h$ dans $\gK^n$. \\
On suppose que la matrice
$\bloc{ \I_k}{L}{C}{\I_h}$ a ses $k$ premières colonnes qui engendrent $E_1$
et ses $h$ dernières colonnes qui engendrent $E_2$.
%Alors:
%---------begin enum----------
\begin{enumerate}
\item Les matrices $L$ et $C$ sont entièrement déterminées par la paire
$(E_1,E_2)$.
\item La matrice $\I_k-LC$  est \iv (on note $V$ son inverse).
\item La \mprn sur $E_1$ \paralm à $E_2$ est égale
à 

\snic{F=\bloc{ V}{-V\,L}{C\,V}{-C\,V\,L}.}
\end{enumerate}
%---------end enum----------
\end{lemma}
%--- end-lemma-----------------
%---------begin proof----------
\begin{proof}
L'unicité est claire. Soit
$F=\bloc{ V}{L'}{C'}{W}$ la matrice de la projection considérée.
Elle est \caree par l'\egt

\snic{F\,\bloc{ \I_k}{L}{C}{\I_h}=\bloc{ \I_k}{0}{C}{0}\,,}

%\sni
\cad encore 

\snic{V+L'C=\I_k$, $\,VL+L'=0$, $\,C'+WC=C\,$  et  $\,C'L+W=0,}

%\sni
ce qui équivaut à
$$
%\snic{
L'=-VL, \;W=-C'L, \;C'(\I_k-LC)=C\,\hbox{ et }\,V\,(\I_k-LC)=\I_k,%}
$$
%\sni
ou encore:  $(\I_k-LC)^{-1}=V$ , $C'=CV$,  $L'=-VL$, et~$W=-CVL$.
\end{proof}
%---------end proof----------

Ceci se \gns au cas d'une \mprn de rang $k$ sur un
anneau commutatif arbitraire de la manière suivante, qui est une variante
commune au lemme de la liberté et au lemme de la liberté locale.

%:   Lem de la liberte le 2eme {leli2}
\CMnewtheorem{lelib2}{Deuxième lemme de la liberté}{\itshape}
\begin{lelib2}\label{leli2} ~\\
Soit $F$ un \prr dans $\GAn(\gA)$; on rappelle que $k + h = n$.
\begin{enumerate}
\item Si  $\rg(F)\leq k$ et  si un
mineur principal d'ordre $k$ est \ivz, alors la matrice
 $F$ est semblable à une
\mprn standard $\I_{k,n}.$
\item Plus \prmtz, supposons que \smashbot{$F=\bloc{V}{L'}{C'}{W}$} avec
$V \in \GL_k(\gA)$.
Posons
$$\preskip.0em \postskip.4em 
B=\bloc{V}{-L'}{C'}{\I_{h}-W}. 
$$
Alors, en posant
$L=V^{-1}L'$ et 
 $C=-C'V^{-1}$, la matrice $B$ est \ivz, d'inverse
$\bloc{\I_k}{L}{C}{\I_{h}}$. En outre, on a les \egts

\snic{\begin{array}{c}
B^{-1}\,F\,B=\I_{k,n},\, W=C'V^{-1}L',\,  V=(\I_k-LC)^{-1},   \\[1mm]
\det V=\det(\I_h-W )\, \hbox{ et } \,
\I_{h}-W=(\I_{h}-CL)^{-1}.
\end{array}
}
\item Réciproquement, si  $L\in \Ae{ k{\times}h}$,
  $C\in \Ae{h{\times}k}$ et si $\I_k-LC$  est \iv
d'inverse $V$, alors
la matrice 

\snic{F=\bloc{V}{V\,L}{-C\,V}{-C\,V\,L}}

%\sni
est une \prn de rang $k$: c'est la
\prn sur le sous-module libre~$E_1$ engendré par les $k$ premières colonnes de
$\bloc{\I_k}{L}{C}{\I_h}$, \paralm au sous-module libre $E_2$ engendré
par les $h$ dernières colonnes de cette matrice.
\end{enumerate}

\end{lelib2}
%--------- fin lelib2 ---------------------------------------------- 

%----begin{proof----------
\begin{proof} Voir l'exercice \ref{exoleli2} et sa solution.
\perso{L'\egt $\det\,V=\det\,W_1$ est démontrée dans l'exo sous
l'hypothèse que $\det\,V$ est \ivz,
 mais il est probable qu'elle est vraie
pour toute \mprn de rang $k$.
}
\end{proof}
%----end{proof----------

Ce que l'on a gagné par rapport au premier lemme de la liberté
\ref{lem pf libre}, c'est que~$F$ est semblable à $\I_{k,n}$ au lieu d'être simplement \eqve
(cependant, voir l'exercice \ref{exo2.4.1}). Et
surtout, les précisions obtenues ici nous seront utiles.
%-% PERSO
\perso{Est-ce que c'était donc tant la peine de se fatiguer? Est-ce que l'exercice \ref{exoleli2}
doit être complété?
}
%-% Fin PERSO

Le lemme précédent se reformule de la manière suivante, plus abstraite,
mais essentiellement \eqve (quoique moins précise).
%:   Lemma{lemthBnk}---------
\begin{lemma}
\label{lemthBnk} \emph{(L'anneau $\Gnk(\gA)$ est presque un anneau de \polsz)}
\\
On considère la matrice \gnq $F=(f_{ij})_{i,j\in\lrbn}$ dans l'anneau $\Gnk(\gA)$. Soit $\mu=\det\big((f_{ij})_{i,j\in\lrbk}\big)$ son mineur principal dominant d'ordre $k$.
\\
Soit par ailleurs $\gA[L,C]$ l'anneau des \pols en $2kh$ indéterminées,
vues comme des \coes de deux matrices $L$ et $C$ de types respectifs
$k{\times}h$ et~$h{\times}k$. Enfin notons $\delta=\det(\I_k-LC)\in\gA[L,C]$.
\\
Alors les anneaux localisés $\Gnk(\gA)[1/\mu]$ et
$\gA[L,C][1/\delta]$ sont naturellement isomorphes.
\end{lemma}
%--- end-lemma-----------------
%---------begin proof----------
\begin{proof} Notons \smashtop{$F= \bloc{V}{L'}{C'}{W}$} avec $V\in \MM_k(\gA)$. \\
Lorsque l'on inverse  $\mu=\det(V)$, on obtient
$V\in\GL_k(\gA[1/\mu])$.  On applique le  point \emph{2} du lemme  \ref{leli2}. 
Avec les matrices   $L=V^{-1}L'$  et  $C=-C'V^{-1}$ on obtient $\delta=\det(\I_k-LC)\in\Ati$. Ceci définit un \homo d'\algs de $\gA[L,C][1/\delta]$ vers
$\Gnk(\gA)[1/\mu]$.\\
Dans l'autre sens: à $L$ et $C$ avec $\delta$ \iv
on fait correspondre la  
matrice  $F=\bloc{V}{V\,L}{-C\,V}{-C\,V\,L}$
(avec $V=(\I_k-LC)^{-1}$). \\
L'\homo correspondant va de
 $\Gnk(\gA)[1/\mu]$ vers $\gA[L,C][1/\delta]$.\\
En composant ces morphismes on trouve l'\idt dans les deux cas.
\end{proof}
%---------end proof----------

%---------begin proof----------
\begin{Proof}{Démonstration du \thrf{thBnk}. }
\\
\emph{\ref{thBnk1}.} Ce point se déduit du lemme pré\-cédent
puisque la somme des mineurs principaux d'ordre $k$ de $F$ est égale à
$1$ dans $\Gnk(\gA)$.

\emph{\ref{thBnka}.} Considérons le $\gA$-\homo
$\psi\,:\gA[(f_{i,j})]\rightarrow \gA$ de spécialisation en $\I_{k,n}$ défini par $\psi(f_{i,j})=1$ si $ i=j\in\lrbk$ et $=0$ sinon. \\
Il est clair que $\psi\big(\cGnk(\gA)\big)=0$. Ceci prouve que
$\gA\cap\cGnk(\gA)=0$ car si $a$ est dans cette intersection, $a=\psi(a)=0$.

\emph{\ref{thBnkf}.} Le noyau de
$\varphi_{L,C}\,:\gA[L,C]\rightarrow\gB[L,C]$ (l'extension naturelle de $\varphi$)
est engendré par le noyau de $\varphi$. La \prt reste vraie après
\lonz. Puis elle reste vraie en recollant des \lons en des \mocoz.
Donc dans notre cas on recolle en disant que
 $\Ker\,\Gnk(\varphi)$ est engendré par~$\Ker\,\varphi$.
\end{Proof}
%---------end proof----------
\begin{Proof}{\Demo du corolaire \ref{corthBnk}. }
\\
\emph{\ref{thBnkd}.} Mise à part la question de l'intégrité cela résulte du point \emph{\ref{thBnk1}} du \thref{thBnk},
car toutes les notions considérées sont stables par $\gA\leadsto \gA[X]$
et relèvent du \plg de base. Pour l'intégrité, cela se déduit du résultat dans le cas d'un \cdiz: si $\gA$ est intègre et $S=\Reg(\gA)$,
alors $\gK=\Frac\gA=\gA_S$ est un \cdi et le point \emph{\ref{thBnkf}} du \tho
\ref{thBnk} permet de conclure.\iplg

\emph{\ref{corthBnke}.} Vu le
\plgc pour la \ddk (voir \paref{thDdkLoc}),
il nous suffit de montrer que $\gA[L,C]$ et $\gA[L,C][1/\delta]$
ont la même dimension, ce qui résulte du lemme \ref{lemLocMemeKdim} ci-après.

\emph{\ref{corthBnkc}.}  Compte tenu des points \emph{\ref{thBnkd}} 
et \emph{\ref{corthBnke}} il reste à montrer que $\Gnk(\gK)$ est intègre. Pour
cela on se rappelle que $\SLn(\gK)$ opère transitivement sur $\GAnk(\gK)$,
ce qui signifie
que toute \mprn de rang $k$ et d'ordre $n$
peut s'écrire sous la forme $S\cdot\I_{k,n}\cdot S^{-1}$ avec $S\in\SLn(\gK)$.
Introduisons l'anneau des \coos de la \vrt  $\SLn(\gK)\subseteq\Mn(\gK)$:
$$\preskip.4em \postskip.0em
\Sln(\gK)=\aqo{\gK[(s_{i,j})_{i,j\in \lrbn}]}{1-\det\,S}.
$$
\`A l'application surjective
$$\preskip.2em \postskip.2em 
\theta_\gK:\SLn(\gK)\rightarrow\GAnk(\gK)\,:\,
S\mapsto S\cdot\I_{k,n}\cdot S^{-1}, 
$$
correspond le $\gK$-\homo

\snic{\wi{\theta}_\gK:\Gnk(\gK)\to\Sln(\gK),}

%\sni
qui envoie chaque  $f_{i,j}$ sur le \coe
$i,j$ de la matrice  $S\cdot\I_{k,n}\cdot S^{-1}$.
\\
Il est bien connu que $\Sln(\gK)$ est intègre, et il suffit donc de montrer que~$\widetilde{\theta}_\gK$ est injectif.
Comme $\theta_\gL$ est surjectif pour toute extension finie $\gL$ de $\gK$, tout
\elt de $\Ker\,\widetilde{\theta}_\gK$
est nilpotent (par le \nstz\footnote{Nous faisons ici une \dem \cov en supposant que
$\gK$ est contenu dans un \cdacz. On pourrait l'adapter au cas \gnlz.}).
Or $\Gnk(\gK)$ est réduit, donc $\widetilde{\theta}_\gK$ est injectif.

\emph{\ref{corthBnkb}.} Résulte des autres points (pour la \ddkz, il faut aussi le \thref{corthValDim}).
\end{Proof}

%:   lemLocMemeKdim
\begin{lemma}\label{lemLocMemeKdim}
Avec les notations précédentes
l'anneau
$\gA[L,C][1/\delta]$ est une extension entière monogène d'un anneau de \pols sur $\gA$ à $2kh$
\idtrsz.
En conséquence
$
\Kdim \gA[L,C][1/\delta] = \Kdim\gA[X_1,\ldots,X_{2kh}].
$
\end {lemma}
%%%%%%%%%%%%%%%%%%%%%%%%%%%%%%%%%%%%%%%%%
\begin {proof} On pose $L = (l_{ij})_{i
\in\lrbk,j \in\lrbh}$, $C = (c_{ij})_{ i\in\lrbh,j \in\lrbk}$.
Le \pol $\delta$ est de degré $2m$ avec $m =\min(h,k)$ et contient le \mom

\snic{(-1)^m l_{11} \ldots l_{mm} c_{11} \ldots c_{mm}.}

%\sni
Le localisé $\gA' =
\gA[L,C][1/\delta]$ peut être réalisé en adjoignant une \idtrz~$t$: $\gA'
= \aqo{\gA[L, C, t]}{t\delta - 1}$. On peut mettre le \pol $g=t\delta - 1$ en position de \iNoez.
En effet, avec le \cdv
$$
l'_{ii} = l_{ii} + t, \; c'_{ii} = c_{ii} + t, \  i \in\lrbm, \quad
l'_{ij} = l_{ij},\; c'_{ij} = c_{ij} \hbox { si } i\neq j,
$$
le \pol $g$ devient, au signe près, \mon en $t$. Donc $\gA'$ est une extension entière monogène de  $\gA[L',C']$.
On conclut avec le \thref{cor2thKdimMor}.
\end {proof}
%%%%%%%%%%%%%%%%%%%%%%%%%%%%%%%%%%%%%%%%%

Nous allons étudier maintenant les espaces tangents aux grassmanniennes.
Nous avons besoin pour ceci de définir le concept lui-même.

Nous commençons donc par une introduction heuristique à des notions
catégoriques et fonctorielles abstraites. \Llec non famili\er avec le
langage des catégories doit survoler cette introduction, dans laquelle nous ne donnons pratiquement pas de \demsz, et 
simplement essayer de se convaincre à partir des exemples donnés 
que la notion d'espace tangent
à un foncteur  en un point est somme toute assez raisonnable,
ce qui lui permettra de voir ensuite la belle application de ce concept 
aux grassmanniennes.

%\penalty-4500
\subsec{Schémas affines, espaces tangents}

\subsubsec{\nst et équivalence de deux catégories}
\label{subsecNstMorphismes}

Soit $(\uf)=(\lfs)$  un \syp dans $\kXn=\kuX$, et notons $\gA=\kxn=\kux$
l'\alg quotient correspondante.

 Nous avons vu \paref{ZerosCrucial}
 l'identification cruciale 
%:2012 les \mathrm{-alg} en indice sont inutiles
\Grandcadre{$\Hom_{\gk%\mathrm{-alg}
} (\gA,\gk)=\cZ(\uf,\gk)\subseteq \gk^n$} 
 entre les zéros sur $\gk$  du \syp $(\uf)$ et les \crcs  de l'\algz~$\gA$.
Si $\gk$ est réduit, on a évidemment $\Hom_{\gk%\mathrm{-alg}
} (\gA,\gk)=\Hom_{\gk%\mathrm{-alg}
} (\Ared,\gk)$.%
 \index{variété algébrique!sur un \cacz}
 \Grandcadre{Supposons maintenant que $\gk$ soit un \cac discret.}
 Un tel ensemble de zéros $\cZ(\uf,\gk)\subseteq \gk^n$ est alors appelé une \emph{\vgq  sur~$\gk$}. 
 
 Soient $\gA$ et
 $\gB$ deux \klgs quotients correspondant à deux \syps $(\uf)$ et $(\ug)$ 
 dans $\kuX=\kXn$.
 Le \nst (corolaire~\ref{corthNstClass}) 
 nous dit que les deux \algs réduites $\gA\red$ et $\gB\red$ sont égales \ssi elles ont la même 
 \vrt de zéros dans $\gk^n$:
\Grandcadre{$\cZ(\uf,\gk)=\cZ(\ug,\gk)\iff \rD_{\kuX}(\uf)=\rD_{\kuX}(\ug)\iff \gA\red=\gB\red$}
Cette constatation est la première étape dans la mise au point de l'\eqvc entre la catégorie des \klgs \rpfs d'une part, et celle des \vgqs sur $\gk$ d'autre part.
 
Pour que l'\eqvc soit complète, nous devons traiter aussi les morphismes.
Nous faisons pour cela une étude préliminaire concernant
l'\algz~$\Ared$.

Nous remarquons que tout \elt $p$ de $\kuX$ définit une fonction \polle  $\gk^n\to\gk,\;\uxi\mapsto p(\uxi)$, et qu'un \elt de $\Ared$ définit (par restriction) une fonction $\cZ(\uf,\gk)\to\gk$: en effet, si $p\equiv q\mod\rD_{\kuX}(\uf)$, une puissance \hbox{de $p-q$} est dans l'\id $\gen{\uf}$, donc les restrictions des fonctions \polles $p$ et $q$ à~$\cZ(\uf,\gk)$ sont égales.
Mais dans le cas où $\gk$ est un \cacz, nous avons la réciproque: si 
les restrictions de $p$ et $q$ à $\cZ(\uf,\gk)$ sont égales,
 $p-q$ s'annule sur $\cZ(\uf,\gk)$, et par le \nstz, une puissance \hbox{de $p-q$} est dans \hbox{l'\id $\gen{\uf}$}. 

\rdb
Ainsi, $\Ared$ peut être interprétée comme une \alg de fonctions sur 
la \vgq qu'elle définit, à savoir $A=\cZ(\uf,\gk)=\Hom_{\gk%\mathrm{-alg}
} (\gA,\gk)$. La structure de \klg de $\Ared$ est bien celle de cette \alg de fonctions.
Ces fonctions $\cZ(\uf,\gk)\to \gk$ sont appelées les \emph{\frgsz}. 
\index{fonction régulière}\index{reguliere@régulière!fonction ---}%

De la même manière, si $\gA=\kxn$ et $\gC=\gk[\ym]$ sont les \algs quotients correspondant à deux \syps 
$$ 
%\snic{
(\uf) \hbox{  dans } \kXn \hbox{ et } (\uh) \hbox{ dans } \kYm,
$$
\rdb%\sni
si $A=\cZ(\uf,\gk)\subseteq \gk^n$ et $C=\cZ(\uh,\gk)\subseteq \gk^m$ sont les \vgqs correspondantes, on définit une \emph{\aregz} de~$A$ vers~$C$ comme la restriction à $A$ et $C$ d'une application \pollez~\hbox{$\varphi:\gk^n\to\gk^m$} qui envoie~$A$ dans~$C$.%
\index{application régulière}\index{reguliere@régulière!application ---}
\\
Les \aregs sont, par \dfnz, \emph{les morphismes de $A$
vers~$C$ dans la catégorie des \vrts \agqs sur $\gk$}.
On notera  $\Mor_\gk(A,C)$ l'ensemble de ces morphismes.
 
L'application  $\varphi$ ci-dessus est donnée par un \sys $(F_1,\dots,F_m)$ dans~$\kuX$, ou encore, par l'\homo $F:\kuY\to\kuX,\;Y_j\mapsto F_j$.\\
Notons $\varphi_1:A\to C$ la restriction de $\varphi$;  si $\gamma:C\to \gk$ est une \frgz, alors la composée $\gamma\circ \varphi_1:A\to\gk$
est une \frgz, et l'application $\psi_1:\gamma\mapsto \gamma\circ \varphi_1$
peut être vue comme une application de~$\gC\red$ vers~$\Ared$. En fait, cette application n'est autre que l'\homo qui provient de $F$ par passage aux quotients. 

Dans l'autre sens, on peut voir que tout \homo $\psi_1:\gC\red\to\Ared$
provient d'un \homo $\psi:\gC\to\gA$, et que $\psi$ définit une \areg $\varphi:A\to C$, parfois appelée le \emph{co-morphisme} de $\psi$.
Cela se passe de la manière suivante: via les identifications  
$A=\Hom_{\gk%\mathrm{-alg}
} (\gA,\gk)$ \hbox{et $C=\Hom_{\gk%\mathrm{-alg}
} (\gC,\gk)$},
on a simplement l'\egt $\varphi(\uxi)=\uxi\circ \psi$ (ce qui fait de~$\varphi$
la \gui{transposée} de $\psi$).\index{co-morphisme}

Finalement,  $\Mor_\gk(A,C)$, s'identifie
 naturellement  
 à  $\Hom_{\gk%\mathrm{-alg}
 }(\gC\red,\gA\red)$, 
 identification que nous écrivons sous la forme d'une \egtz:
\Grandcadre{$\Mor_\gk(A,C)=\Hom_{\gk%\mathrm{-alg}
}(\gC\red,\gA\red)$.}
Notons cependant que le sens des flèches est renversé.

Considérons comme cas particulier le cas où $A$ est
la \vgq réduite à un point, associée à l'\algz~$\gk$, correspondant au \syp vide sur \emph{l'\alg \polle sans variable $\gk$}. Si l'on préfère, on peut voir ici $\gk$ comme le quotient $\aqo\kX X$, correspondant au point~$\so 0$, sous-\vrt de la \vgq $V=\gk$
associée à l'\alg $\kX$.
\\
Dans ces conditions, l'\egt encadrée ci-dessus
admet comme cas particulier~\hbox{${C=\Mor_\gk(\so{0},C)=
\Hom_{\gk%\mathrm{-alg}
}(\gC\red,\gk).}$}
La boucle est bouclée!

Le bilan de cette étude est le suivant: on peut entièrement réduire 
la considération des \vgqs  sur un \cac à l'étude des 
\klgs \rpfsz: il s'agit d'une interprétation en termes finis
(\syps finis sur~$\gk$ pour les objets aussi bien que pour les morphismes) d'objets a priori un peu plus mystérieux, et certainement plus infinis. En~termes catégoriques: on peut avantageusement remplacer la catégorie des \vrts \agqs sur $\gk$ par la catégorie
opposée à celle des \klgs \rpfsz. Il y a une \eqvc  naturelle entre ces deux catégories.

%: subsubsec{Schémas affines}  
\subsubsec{Schémas affines}
\label{subsecSchAff}
 
Maintenant on fait un grand saut dans l'abstraction. 
On admet tout d'abord que les \vrts 
peuvent avoir des multiplicités. Par exemple l'intersection d'un cercle et d'une droite doit être un point double, et non pas seulement un point, lorsque la droite est tangente au cercle. En conséquence, il est parfois néfaste de se limiter aux \klgs réduites.

On admet aussi qu'à la base on n'a pas \ncrt un \cac mais un anneau commutatif arbitraire. Auquel cas les points de la \vrt sur $\gk$ ne sauraient en \gnl caractériser ce que l'on a envie de considérer comme une \vrt \agq abstraite
définie sur $\gk$ (en autorisant les multiplicités).
Par exemple le cercle abstrait est certainement représenté par la \ZZlg

\snic{\ZZ[x,y]=\aqo{\ZZ[X,Y]}{X^2+Y^2-1},}

%\sni
mais ce ne sont pas ses points sur $\ZZ$
qui vont nous donner beaucoup d'information. Bien au contraire, ce sont ses points sur toutes les \ZZlgsz, \cad sur tous les anneaux commutatifs, qui nous importent.
De même un \emph{cercle double} abstrait  est certainement représenté par la \ZZlg

\snic{\ZZ[x',y']=\aqo{\ZZ[X,Y]}{(X^2+Y^2-1)^2},}

%\sni
mais on ne saurait distinguer un cercle simple d'un cercle double si l'on ne considère que les points sur les anneaux réduits (les anneaux sans multiplicité).

Nous voici donc en état de définir la catégorie des \emph{schémas affines sur l'anneau commutatif $\gk$}. Cela pourrait être simplement la catégorie 
opposée à la catégorie des \klgsz: celle dont les objets sont les \klgs et dont les flèches sont les \homos de \klgsz.

Mais il est une description \eqve nettement plus parlante (et élégante?):
\emph{un schéma affine sur l'anneau commutatif $\gk$ est connu 
lorsque l'on connaît ses zéros sur toutes les \klgsz}. 
Autrement dit, la \klgz~$\gA$ définit un schéma
affine qui n'est rien d'autre que le foncteur $\Hom_{\gk%\mathrm{-alg}
}(\gA,\bullet)$ de la catégorie des \klgs vers la catégorie des ensembles.

Et un \homo de \klgs $\gB\to\gA$ définit une transformation naturelle du foncteur
$\Hom_{\gk%\mathrm{-alg}
}(\gA,\bullet)$ vers le foncteur  $\Hom_{\gk%\mathrm{-alg}
}(\gB,\bullet)$: les transformations naturelles des foncteurs sont \gui{dans le bon sens},
\cad des zéros de $\gA$ vers les zéros de $\gB$.

Si l'on ne veut pas partir trop haut dans l'abstraction, on peut se limiter
aux \klgs \pfz, ce qui est bien assez pour faire de la très belle \gmt \agq
abstraite (i.e., non limitée à la \gmt \agq sur des \cdisz).

%: subsubsec{Espace tangent en un point à un foncteur}  
\subsubsec{Espace tangent en un point à un foncteur}
\label{subsubsecTanFonct}
%:HHH index
\index{tangent!espace ---}\index{espace tangent}

Rappelons tout d'abord la notion d'espace tangent à un \syp en un zéro du \sys
introduite à  la section \ref{secExlocGeoAlg}.

Prenons l'exemple de la sphère comme schéma affine défini sur $\QQ$.
Ce schéma est associé à la \QQlg $\gA=\QQ[x,y,z]=\aqo{\QQ[X,Y,Z]}{X^2+Y^2+Z^2-1}$.
Si $\uxi=(\alpha,\beta,\gamma)\in\QQ^3$ est un zéro de $\gA$ sur $\QQ$, \cad un point rationnel de la sphère, nous lui avons associé
\begin{itemize}
\item l'\id $\fm_\uxi=\gen{x-\alpha,y-\beta,z-\gamma}_\gA$,
\item  l'\alg locale $\gA_\uxi=\gA_{1+\fm_\uxi}$, et 
\item l'espace tangent $\rT_\uxi(\gA/\QQ)
%:HHH ajout \simeq\Der\QQ\gA\uxi
\simeq\Der\QQ\gA\uxi$, 
\end{itemize}
lequel est un \Qev canoniquement isomorphe à
$(\fm_\uxi/{{\fm_\uxi}^2})\sta$ ou encore à $(\fm_\uxi\gA_\uxi/{{\fm_\uxi\gA_\uxi}^2})\sta$.

De manière plus intuitive mais \eqve
(proposition \ref{propTangent}), un vecteur tangent à la sphère en $\uxi$ est simplement donné par un $(u,v,w)\in\QQ^3$ qui vérifie $u\alpha+v\beta+w\gamma=0$, \cad en posant $f=X^2+Y^2+Z^2-1$, 

\snic{u\Dpp f X(\uxi)+v\Dpp f Y(\uxi)+w\Dpp f Z(\uxi)=0.}

%\sni
Voici maintenant une nouvelle manière de voir cet espace tangent, 
que nous exprimons en termes du 
schéma affine correspondant, \cad du foncteur
$\Hom_\QQ(\gA,\bullet)=\cZ(f,\bullet)$. Nous devons pour ceci introduire de manière formelle un infinitésimal que nous notons $\vep$, \cad considérer la \QQlg
$\QQ[\vep]=\aqo{\QQ[T]}{T^2}$ ($\vep$ est la classe de $T$ modulo $T^2$).

Le point $\uxi$ est vu comme un \crc de
$\gA$, i.e. comme l'\elt $\wi\uxi:g\mapsto g(\uxi)$ \linebreak 
de $\Hom_\QQ(\gA,\QQ)$. 
Nous nous demandons alors quels sont les \elts $\lambda$ de l'ensemble
 $\Hom_\QQ(\gA,\QQ[\vep])$
qui \gui{relèvent $\wi\uxi$}, au sens que lorsque l'on compose avec l'évaluation de $\vep$ en $0$, de $\QQ[\vep]$ vers $\QQ$, on retombe sur $\wi\uxi$.

\centerline{\xymatrix {
             &\QQ[\vep]\ar[d]^{\vep := 0}  \\
\gA\ar@{-->}[ur]\ar[r]^{\wi\uxi} &\QQ                   \\
}} 

Un tel \elt est a priori donné par un zéro de $f$
sur $\QQ[\vep]$ qui redonne~$\uxi$ lorsque l'on  évalue $\vep$ en $0$, 
\cad un triplet %de la forme
$(\alpha+a\vep,\beta+b\vep,\gamma+c\vep)$,\linebreak 
 avec $f(\alpha+a\vep,\beta+b\vep,\gamma+c\vep)=0$ dans $\QQ[\vep]$. Mais ceci signifie exactement que $(a,b,c)$
est un vecteur tangent à la sphère en $\uxi$. 

Il ne s'agit quant au fond que de la banale constatation selon laquelle
\gui{la \dile est la partie \lin de l'accroissement de la fonction}:

\snic{f(\uxi+\vep V)=f(\uxi)+\vep \,\rd f(\uxi)(V)  \mod \vep^2.}

%\sni
Ce zéro $\uxi+\vep(a,b,c)$ de $\gA$ dans $\gk[\vep]$ définit un \homo
% l'\homo correspondant 
$\gA\to\gk[\vep]$ via 
$x\mapsto \alpha+a\vep$, $y\mapsto \beta+b\vep$, $z\mapsto \gamma+c\vep$.
\\
Cet \homo envoie   $g$ sur
$g(\uxi)+a\Dpp g X(\uxi)+b\Dpp g Y(\uxi)+c\Dpp g Z(\uxi)$, puisque

\snic{g\big(\uxi+\vep (a,b,c)\big)=g(\uxi)+\vep \,\rd g(\uxi)(a,b,c)  \mod \vep^2.}

\smallskip 
\Llec pourra vérifier que ce petit calcul que nous venons de faire sur un
petit exemple fonctionne pour n'importe quel zéro de n'importe quel \syp
basé sur n'importe quel anneau commutatif.

Il faut cependant au moins rajouter comment on peut interpréter en termes
du foncteur $\Hom_\gk(\gA,\bullet)$ la structure de \kmo sur l'espace tangent en un zéro d'un \syp sur un anneau $\gk$. 

Ici aussi, contentons-nous de notre petit exemple. 

Dans la catégorie des \QQlgsz, le produit fibré de la \gui{flèche de restriction} 

\snic{\QQ[\vep]\to\QQ$, $\vep\mapsto0}

%\sni
avec elle-même est l'\alg

\snic{\QQ[\vep]\times_\QQ\QQ[\vep]\simeq \QQ[\vep_1,\vep_2] 
\quad \hbox{avec }  \vep_1^2=\vep_1\vep_2=
\vep_2^2=0,}

%\sni
munie des deux \homos 
\gui{de \prnz} 

%:HHH un peu plus lisible comme cela
\snic{\QQ[\vep_1,\vep_2]\vers{\pi_1}\QQ[\vep],\;\vep_1\mapsto \vep,\; \vep_2\mapsto 0\quad$ et $\quad\QQ[\vep_1,\vep_2]\vers{\pi_2}\QQ[\vep],\;\vep_2\mapsto \vep,\; \vep_1\mapsto 0,~~~~~~}

%\sni
et de la flèche \gui{de restriction}

\snic{\QQ[\vep_1,\vep_2]\to\QQ,\quad\vep_1\mapsto 0,\;\vep_2\mapsto 0.}

%\sni
Il y a en outre un \homo naturel \gui{d'addition}

\snic{\QQ[\vep_1,\vep_2]\to\QQ[\vep],\quad\vep_1\mapsto\vep,\; \vep_2\mapsto\vep,}

%\sni
qui commute avec les restrictions.

Lorsque l'on donne deux zéros
$\uxi+\vep V_1$ et $\uxi+\vep V_2$ de $\gA$ dans $\QQ[\vep]$, vue la \prt \cara 
du produit fibré dans la catégorie des \QQlgsz, les deux \homos correspondants $\gA\to\QQ[\vep]$ se factorisent de manière unique pour donner
un \homo de  $\gA$ vers $\QQ[\vep_1,\vep_2]$ \gui{produit fibré des deux},
qui correspond au zéro $\uxi+\vep_1 V_1+\vep_2 V_2$ de $\gA$
dans $\QQ[\vep_1,\vep_2]$. 

Enfin en composant cet \homo {produit fibré} avec l'\homo d'addition $\QQ[\vep_1,\vep_2]\to\QQ[\vep]$, on obtient
l'\homo correspondant au zéro $\uxi+\vep (V_1+V_2)$. La boucle est donc bouclée, l'addition des vecteurs tangents a été décrite en termes purement catégoriques.

\entrenous{un joli dessin semblerait utile pour faire comprendre ce qui se passe ici}

Résumons nous. Dans le cas du foncteur qui est un schéma affine défini par un \syp sur 
un anneau $\gk$ avec son \alg quotient $\gA$, il y a une identification canonique entre $\rT_\uxi(\gA\sur\gk)$
et l'ensemble des points de~$\gA$ sur $\gk[\vep]$ qui relèvent $\uxi$, lorsque
l'on identifie $\uxi$ et $\uxi+\vep V$ aux \elts correspondants de $\Hom_\gk(\gA,\gk)$ et $\Hom_\gk(\gA,\gk[\vep])$.
En outre, la structure de~\kmo dans la deuxième interprétation est donné
par l'\gui{addition} fournie par l'\homo 

\snic{\gk[\vep_1,\vep_2]\simeq\gk[\vep]\times_\gk \gk[\vep]\to\gk[\vep],\quad
%\vep\mapsto \vep_1+\vep_2
\vep_1\mapsto\vep,\;\vep_2\mapsto\vep,}

%\sni
($\vep^2=\vep_1^2=\vep_2^2=\vep_1\vep_2
=0$). 

Notons que la \gui{loi externe}, multiplication par le scalaire $a$, provient, elle, de l'\homo 

\snic{\gk[\vep]\vers{\lambda_a}\gk[\vep],\quad b+\vep c\mapsto b+\vep ac.}

%\sni
Le mécanisme formel d'addition ainsi décrit pourra fonctionner avec n'importe quel autre foncteur qui voudra bien, lui aussi, transformer les produits fibrés (dans la catégorie des \QQlgsz) en produits fibrés 
(dans la catégorie des ensembles).

Ainsi la notion d'espace tangent en un point à un foncteur\footnote{Foncteur de la catégorie des \klgs vers la catégorie des ensembles.}
se \gns aux autres schémas sur un anneau $\gk$,
car ce sont 
\gui{de bons foncteurs}. I.e. les schémas de Grothendieck
(que nous ne définirons pas ici)
sont de bons foncteurs. Et les foncteurs 
grassmanniennes (qui ont déjà été définis) sont de tels schémas.

\subsec{Espaces tangents aux grassmanniennes}

%: subsubsec{Projecteurs et rangs}  %%%%%%%%%%%%
\subsubsec{Projecteurs et rangs}
\label{subsecPrRg}

Deux faits faciles avant d'entrer dans le vif du sujet.
On considère un module~$E$.
Deux \prrs $\pi_1$, $\pi_2 : E \to E$  sont dits
\ixc{orthogonaux}{projecteurs ---} s'ils vérifient $\pi_1 \circ
\pi_2 = \pi_2 \circ \pi_1 = 0$.

% :     Fact{lemSomProjOrt}
\begin{fact}\label{lemSomProjOrt}
Si $\pi_1$, $\pi_2 : E \to E$ sont des
\prrs \orts d'images $E_{1}$ \linebreak 
et $E_{2}$, alors
 $\pi_1 + \pi_2$ est un \prr et son
image est $E_1 \oplus E_2$. En conséquence, lorsque
$E$ est un \mptfz,  on obtient

\snic{\rg({\pi_{1}+\pi_{2}}) = \rg({E_{1}\oplus E_{2}}) = \rg {E_1} + \rg {E_2} .}
\end{fact}

% :     Fact
\begin {fact} \label{lemProjProj}
Soient $\pi_1$, $\pi_2 \in \End_\gA(E)$ deux \prrs d'images $E_1$
et $E_{2}$.
Alors l'\Ali

\snic{
\Phi : \End(E)\to\End(E),\; \varphi \mapsto \pi_2\circ\varphi\circ\pi_1,}

%\sni
est un \prr dont l'image est isomorphe à
$\Lin_{\gA}(E_1, E_2)$. En conséquence, lorsque
$E$ est un \mptfz,  on obtient l'\egt

\snic{\rg\Phi  = \rg E_1 \cdot \rg E_2.}
\end {fact}

%: subsubsec{Grassmannienne affine}  %%%%%%%%%%%%
\subsubsec{Grassmannienne affine}
\label{subsecGrassAff}

Ce paragraphe est consacré à la détermination
de l'espace tangent en un point au foncteur $\gA \mapsto \GAn(\gA)$.
Rappelons que l'acronyme $\GA$ est mis pour ``Grassmannienne Affine''.
L'interprétation \gmq d'un point $P$ de  $\GAn(\gA)$ est donnée par le
couple ordonné $(E,F)=(\Im P,\Ker P)$ de sous-modules en somme directe dans
$\Ae n $.

Plus \gnltz, si $\gk $ est un anneau donné en référence (en
géométrie usuelle ce serait un corps discret)
et si $M$ est un \kmo \ptf fixé,
on peut considérer la catégorie des \klgs
et le \linebreak 
foncteur $\gA \mapsto \GA_M(\gA)$, où
$\GA_M(\gA)$ désigne l'ensemble des couples ordonnés~$(E,F)$
de sous-modules en somme directe dans le module étendu~$\gA\otimes _{\gk}M$, que nous noterons~$M_{\gA}$.
Un tel couple peut être représenté par la projection~$\pi:M_{\gA}\to M_{\gA}$ sur $E$ \paralm à $F$.
La grassmannienne affine~$\GA_M(\gA)$ peut donc être vue comme
le sous-ensemble des \elts \idms dans~$\End_{\gA}(M_{\gA})$.
C'est ce point de vue que nous adoptons dans la suite.

Pour étudier l'espace tangent on doit considérer
l'\Alg $\gA[\varepsilon]$
où $\varepsilon$ est
l'\elt générique de carré nul.
Nous donnons tout d'abord l'énoncé pour la grassmannienne
usuelle $\GAn(\gA)$.

%%%%%%%%%%%%%%%%%%%%%%%%%%%%%%%%%%%%%%%%%
%:     theorem{prop1TanGrassmann}
\begin{theorem}\label{prop1TanGrassmann} \emph{(Espace tangent à une grassmanienne affine)}\\
Soit $P\in \GAn(\gA)$  un \prr d'image $E$ et de noyau~$F$.
Pour  $ H \in \MM_n(\gA)$
on a l'\eqvc suivante.

\snic{P + \varepsilon H\in \GAn(\gA[\varepsilon]) \quad  \Longleftrightarrow
\quad  H = HP + PH .}
%

%\sni
Associons à $P$
l'\Ali $\wh P :
\MM_n(\gA) \to \MM_n(\gA)$ définie par

\snic{\widehat P(G) = P\,G\,(\In-P) + (\In-P)\,G\,P.}

%\sni
On a les résultats suivants.
\begin{enumerate}
\item [--] Les \Alis 

\snic{\pi_1 : G\mapsto P\,G\,(\In-P)\;$ et
$\;\pi_2 : G\mapsto (\In-P)\,G\,P}

%\sni
sont des \prrs \ortsz.
En particulier,  $\widehat P$ est un \prrz.
\item [--] Pour $H \in \MM_n(\gA)$, on a $H
= PH + HP$ \ssi $H\in\Im\wh P$.
\item [--] Le module $\Im\wh{P}$ est canoniquement  isomorphe
à $\Lin_{\gA}(E,F)\oplus\alb \Lin_{\gA}(F,E)$.
En particulier, $\rg(\Im\wh P)= 2 \rg E  \cdot \rg  F $.
\end{enumerate}
En résumé, l'espace tangent en le $\gA$-point $P$ au foncteur $\GAn$
est canoniquement isomorphe au \mptf  $\Im\wh{P}$ (via $H \mapsto P+\varepsilon H$), lui-même canoniquement  isomorphe
à $\Lin_{\gA}(E,F)\oplus\alb \Lin_{\gA}(F,E)$.
\end{theorem}
%%%%%%%%%%%%%%%%%%%%%%%%%%%%%%%%%%%%%%%%%
\begin {proof}
Le premier point est immédiat. Notons $V_{P}$ le sous-module
des matrices~$H$ qui vérifient $H=HP+PH$. Ce module est
canoniquement isomorphe à l'espace tangent que nous cherchons.
Un calcul simple montre que~$\pi_1  $ et~$\pi_2  $ sont des \prrs
\ortsz. Donc $\widehat P$ est un \prrz.
L'\egt suivante
est claire:  $P\widehat P(G) + \widehat P(G)P = \widehat P(G)$.
Donc $\Im \wh P \subseteq V_{P}$.
Par ailleurs, si~$H = PH+HP$, on a $PHP = 0$, donc $\wh P(H) = PH + HP=H$.
Ainsi~$V_{P} \subseteq \Im \wh P$.
En bref $V_{P}=\Im \wh P=\Im\pi_{1} \oplus \Im\pi_{2}$: on conclut en
appliquant le fait~\ref{lemProjProj}.
\end{proof}
%%%%%%%%%%%%%%%%%%%%%%%%%%%%%%%%%%%%%%%%%

Nous donnons maintenant l'énoncé \gnl (la
preuve est identique).

%%%%%%%%%%%%%%%%%%%%%%%%%%%%%%%%%%%%%%%%%
%:     Proposition{prop2TanGrassmann}
\begin{proposition}\label{prop2TanGrassmann}
Soit $\pi\in \GA_M(\gA)$  un \prr d'image $E$ et de noyau~$F$.
Pour  $ \eta \in \End_{\gA}(M_{\gA})$
on a l'\eqvc

\snic{\pi + \varepsilon \eta\in \GA_M(\gA[\varepsilon]) \quad  \Longleftrightarrow
\quad  \eta = \pi\eta + \eta\pi .}
%

%\sni
On  associe à $\pi$
l'\Ali $\wh \pi :
\End(M_{\gA}) \to \End(M_{\gA})$ définie par
$
\widehat \pi(\psi) = \pi\,\psi\,(\I-\pi) + (\I-\pi)\,\psi\,\pi
$. Alors
\begin{enumerate}
\item [--] Les \alis $\pi_1 : \psi\mapsto \pi\, \psi\,(\I-\pi )$ et
$\pi_2 : \psi\mapsto (\I-\pi )\,\psi\, \pi $ sont des \prrs \ortsz.
En particulier,  $\widehat \pi $ est un \prrz.
\item [--] Une \Ali $\eta \in \End(M_{\gA})$ vérifie $\eta
= \pi\eta + \eta\pi$ \ssi $\eta\in\Im\wh \pi$.
\item [--] Le module $\Im\wh{\pi}$ est canoniquement  isomorphe
à $\Lin_{\gA}(E,F)\oplus\alb \Lin_{\gA}(F,E)$.
En particulier, $\rg(\Im\wh \pi )= 2 \rg E  \cdot \rg  F $.
\end{enumerate}
En résumé l'espace tangent en le $\gA$-point $\pi$ au foncteur $\GA_M$
est canoniquement isomorphe au \mptf  $\Im\wh{\pi}$
(via $\eta \isosim \pi+\varepsilon \eta$), lui même canoniquement  isomorphe
à $\Lin_{\gA}(E,F)\oplus\alb \Lin_{\gA}(F,E)$.
\end{proposition}
%%%%%%%%%%%%%%%%%%%%%%%%%%%%%%%%%%%%%%%%%

%\newpage
%:  subsubsec{Grassmannienne projective}  %%%%%%%%%%%%
\subsubsec{Grassmannienne projective}
\label{subsecGrassProj}

Ce paragraphe est consacré à la détermination
de l'espace tangent en un point au foncteur $\gA \mapsto \GGn(\gA)$,
où $\GGn(\gA)$ désigne l'ensemble des sous-modules en facteur direct
dans $\Ae{n}$.
%Naturellement l'acronyme $\GG$ est mis pour ``Grassmannienne''.

%%%%%%%%%%%%%%%%%%%%%%%%%%%%%%%%%%%%%%%%%
%:     Fact{fact1GrassProj}
\begin{fact}\label{fact1GrassProj} \emph{(L'espace des
\prrs qui ont la même image qu'un \prr fixé)}\\
Soit $P\in \GGn(\gA)$ un \prr d'image $E$. Notons $\Pi_E$
l'ensemble des \prrs qui ont $E$ pour image, et $V=\Ae n $.
Alors $\Pi_E$
est un sous-espace affine de~$\Mn(\gA)$, ayant pour
 \gui{direction} le \Amo  \ptfz~$\Lin_\gA(V/E,E)$ (naturellement
identifié à un sous-\Amo de $\Mn(\gA)$).
On  précise ce résultat de la manière suivante.
\begin{enumerate}
\item \label{i1fact1GrassProj}
Soit $Q \in \GGn(\gA)$ un autre \prrz. \\
Alors $Q\in \Pi_E$ \ssi
 $PQ = Q$ et $QP = P$. \\
Dans ce cas, la différence $N = Q-P$
vérifie les \egts $PN = N$ et~$NP = 0$, et donc $N^2 = 0$.
\item \label{i2fact1GrassProj}
 Réciproquement, si $N \in
\Mn(\gA)$ vérifie $PN = N$ et $NP = 0$ (auquel cas $N^2 = 0$),
alors $Q := P+N$ est dans $\Pi_E$.
\item \label{i5fact1GrassProj}
En résumé, l'ensemble $\Pi_E$ s'identifie au \Amo
$\Lin_\gA(V/E, E)$ via l'application affine

\snic{\Lin_\gA(V/E, E)\to \Mn(\gA),\;\varphi\mapsto P+j\circ\varphi\circ \pi,}

%\sni
où $j:E\to V$ est l'injection canonique et $\pi:V\to V/E$
la \prn canonique.
\end{enumerate}
%
%%%%%%%%%%%%%
Informations supplémentaires.
\begin{enumerate}\setcounter{enumi}3
\item \label{i3fact1GrassProj}
Si $Q \in \Pi_E$, $P$ et $Q$ sont
conjugués dans $\Mn(\gA)$.
Plus \prmtz, en posant $N = Q-P$, on a $(\In + N)(\In - N) = \In$
et
$(\In - N) P (\In + N) = Q$.
\item \label{i4fact1GrassProj}
Si $Q \in \Pi_E$, alors
pour tout $t \in \gA$, on~a: $tP + (1-t)Q\in \Pi_E$.
\end{enumerate}

\end{fact}
%%%%%%%%%%%%%%%%%%%%%%%%%%%%%%%%%%%%%%%%%
\begin{proof}
\emph{\ref{i1fact1GrassProj}.}
$N^2 = 0$ comme on le voit en multipliant $PN = N$ par $N$ à gauche. 

\emph{\ref{i5fact1GrassProj}.}
Les conditions
$PN = N$ et $NP = 0$ sur la matrice $N$
équivalent aux inclusions  $\,\Im N\subseteq E=\Im P$
et $E\subseteq \Ker N$. 
\\
Les matrices $N$
de ce type forment un \Amo $\wi E$ qui s'identifie au  
\linebreak
module $\Lin_\gA(\Ker P, \Im P)$
 \gui{par restriction du domaine et de l'image}.
\\
 De manière plus intrinsèque, ce module $\wi E$  s'identifie aussi
 à $\Lin_\gA(V/E, E)$ via l'\ali
 $\Lin_\gA(V/E, E) \to \Mn(\gA),\;\varphi\mapsto j\circ\varphi\circ \pi$,
qui est injective et admet $\wi E$ pour image.

\emph{\ref{i3fact1GrassProj}.} $(\In - N)\, P \,(\In + N) = P\,(\In + N) = P + PN = P + N = Q$.
\end{proof}
%

%:     Fact{fact2GrassProj}
\begin{fact}\label{fact2GrassProj}
Soit $E\in\GGn(\gA)$ et $E'\in\GGn(\gA[\varepsilon])$ qui donne $E$
par la spécialisation $\varepsilon\mapsto 0$ (autrement dit $E'$
est un point de l'espace tangent en $E$ au foncteur $\GGn$). Alors
$E'$ est isomorphe au module obtenu à partir de $E$ par \edsz:
$E'\simeq \gA[\varepsilon]\otimes _\gA E$.
\end{fact}
\begin{proof}
D'après le \thref{propComparRedRed},
un \mptf $M$ sur un anneau $\gB$ est caractérisé, à \iso près,
par sa \gui{réduction} $M\red$ (i.e., le module obtenu par \eds à $\gB\red$).
Or $E'$ et $\gA[\varepsilon]\otimes _\gA E$
ont même réduction $E\red$ à $(\gA[\varepsilon])\red\simeq\Ared$.
\perso{demander une preuve directe en exercice?}
\end{proof}
%
%%%%%%%%%%%%%%%%%%%%%%%%%%%%%%%%%%%%%%%%%
%:     theorem{prop3TanGrassmann}
\begin{theorem}\label{prop3TanGrassmann} \emph{(Espace tangent à une grassmanienne projective)}\\
Soit $E\in \GGn(\gA)$ un sous-\Amo en   facteur direct dans $\Ae{n}=V$.
%On considère le \Amo $\Lin_\gA(E,V/E)$.\\
Alors l'espace tangent en le $\gA$-point $E$ au foncteur $\GGn$
est canoniquement isomorphe à $\Lin_\gA(E,V/E)$.
Plus \prmtz, si $\varphi\in\Lin_\gA(E,V/E)$ et si l'on note
$$
E_\varphi=\sotq{x+\varepsilon h}{x\in E,\,h\in V,\,h\equiv\varphi(x)\mod E}
,$$
alors $\varphi\mapsto E_\varphi$ est une bijection du module $\Lin_\gA(E,V/E)$
sur l'ensemble des matrices $E'\in \GGn(\gA[\varepsilon])$ qui donnent $E$
lorsque l'on spécialise $\varepsilon$ en $0$.
\end{theorem}
%%%%%%%%%%% \end{theorem} %%%%%%%%%%%%%%%
%
%-% perso
\perso{Moralement le \tho devrait résulter du fait \ref{fact1GrassProj}
sans nouveau
calcul. Surtout si l'on donne une version plus explicite du
fait \ref{fact2GrassProj}. Un mystère pour moi: pourquoi $\Lin_\gA(V/E,E)$
la première fois et $\Lin_\gA(E,V/E)$ la deuxième? Henri
}
%-% Fin perso
%%%%%%%%%%%%  \begin {proof}  %%%%%%%%%%%%
\begin {proof}
Soient $E \in \GG_n(\gA)$ et $\varphi \in \Lin_\gA(E, V/E)$.
\\
\emph{Montrons d'abord que
$E_\varphi$ est dans $\GG_n(\gA[\varepsilon])$
et au dessus de~$E$.}
Fixons une matrice $P \in \GA_{n}(\gA)$ vérifiant $E = \Im P$. On a donc $V = E \oplus\Ker P$ et un isomorphisme $V/E \simeq\Ker P \subseteq V$.
On\vadjust{\break} peut donc relever l'\aliz~$\varphi$  en une matrice $H \in \MM_{n}(\gA) = \End(V)$
\vadjust{\nobreak}conformément au diagramme
$$\preskip-0.5ex\xymatrix @C=28pt @R = 16pt {
&V \ar@{->}[r]^{ H} \ar@{->>}[d]    & V  \\
&E \ar@{->}[r]^{ \varphi}          & V/E \ar@{(->}[u] &
}$$
%\sni
La matrice $H$ vérifie $PH = 0$ et $H(\I_{n}-P) = 0$, i.e. $HP = H$.
\\
%Pour prouver $E_\varphi \in\GG_n(\gA[\varepsilon])$ et 
%$E_\varphi$ au dessus de~$E$, i
Il suffit de montrer que $P + \varepsilon H$ est un \prr
d'image $E_\varphi$. \\
Pour l'inclusion $\Im(P + \varepsilon H) \subseteq
E_\varphi$, soit $(P + \varepsilon H)(y + \varepsilon z)$ avec $y, z \in V$:

\snic {
(P + \varepsilon H)(y + \varepsilon z) = Py + \varepsilon (Hy + Pz) =
Py + \varepsilon (HPy + Pz) = x + \varepsilon h
,}

%\sni
avec $x = Py \in E$, $h = Hx + Pz$. Puisque $x \in E$, on a $\varphi(x) = Hx$,
et \linebreak 
donc $h \equiv \varphi(x) \bmod E$.  Pour l'inclusion réciproque, soit $x +
\varepsilon h \in E_\varphi$ et montrons que $(P + \varepsilon H)(x +
\varepsilon h) = x + \varepsilon h$:

\snic {
(P + \varepsilon H)(x + \varepsilon h) = Px + \varepsilon(Hx + Ph)
.}

%\sni
Comme $x \in E$, on a $Px = x$. Il faut voir que $Hx + Ph = h$,  mais $h$ est
de la forme $h = Hx + y$ avec $y \in E$, donc $Ph = 0 + Py = y$ et l'on a bien
l'\egt $h = Hx + Ph$.
\\
Enfin, il est clair que $P + \varepsilon H$ est un \prrz:

\snic {
(P+\varepsilon H) (P + \varepsilon H) = P^2 + \varepsilon(HP + PH) =
P +\varepsilon H.
}

%%%%%%%%%%%%%%%%%%%%%%%%%%%%%%%%%%%%%%%%%%%%%%%%%%%%%%%%

%\sni
\emph{Montrons la surjectivité de $\varphi \mapsto E_\varphi$}. Soit
$E' \subseteq \gA[\varepsilon]^{n}$, facteur direct, au dessus de~$E$. Alors
$E'$ est l'image d'un \prr $P + \varepsilon H$ et l'on~a:

\snic {
(P + \varepsilon H)(P + \varepsilon H) = P^2 + \varepsilon (HP + PH)
\quad  \hbox {donc} \quad P^2 = P, \quad HP + PH = H
,}

%\sni
ce qui donne $PHP = 0$ (multiplier $HP + PH = H$ par $P$ à droite, par
exemple). On voit donc que~$P$ est un \prr d'image~$E$ 
(car $E'$, \linebreak 
pour $\varepsilon := 0$, c'est~$E$). On remplace $H$ par $K = HP$, qui
vérifie:

\snic {
KP = (HP)P = K, \qquad PK = P(HP) = 0.
}

%\sni
Ceci ne change pas l'image de $P + \varepsilon H$, i.e. $\Im(P + \varepsilon H) =
\Im(P + \varepsilon K)$. Pour le voir, il suffit de (et il faut) montrer que:

\snic {
(P + \varepsilon H)(P + \varepsilon K) = P + \varepsilon K, \qquad
(P + \varepsilon K)(P + \varepsilon H) = P + \varepsilon H
.}

%\sni
\`A gauche, on obtient $P + \varepsilon(HP + PK) = P + \varepsilon(HP + 0) = P +
\varepsilon K$; à droite, $P + \varepsilon(KP + PH) = P +
\varepsilon(K + PH) = P + \varepsilon(HP + PH) = P + \varepsilon H$.
\\
La matrice $K$ vérifie $KP = K$, $PK = 0$, et représente une
\ali $\varphi : E \to \Ae {n}/E$ avec $E' = {\rm
Im}(P + \varepsilon K) = E_\varphi$.

%%%%%%%%%%%%%%%%%%%%%%%%%%%%%%%%%%%%%%%%%%%%%%

%\sni
\emph{Prouvons l'injectivité de $\varphi \mapsto E_\varphi$}. Supposons donc
$E_\varphi = E_{\varphi'}$. On
fixe un \prr $P \in \GG_n(\gA)$ d'image~$E$ et l'on code $\varphi$
par~$H$, $\varphi'$ par~$H'$ avec:

\snic {
HP = H, \quad PH = 0, \qquad H'P = H', \quad PH' = 0
.}

%\sni
Comme $P + \varepsilon H$ et $P + \varepsilon H'$ ont même image,
on a les \egts

\snic {
(P + \varepsilon H)(P + \varepsilon H') = P + \varepsilon H' \quad\hbox{et}\quad 
(P + \varepsilon H')(P + \varepsilon H) = P + \varepsilon H
.}

%\sni
L'\egt de droite donne $H = H'$, donc $\varphi = \varphi'$.
\end {proof}

%%%%%%%%%%%%%%%%%%%%%%%%%%%%%%%%%%%%%%%%%%%%%%%%%%%%%%%%%%%%%%

\rem
La \prn $\GA_{n}\to\GGn$ %$\matrix {\GA_{n}\cr \downarrow\cr \GGn\cr}$ 
%est celle qui
associe à~$P$ son image $E = \Im P$.  Voici comment s'organisent les
espaces tangents et la \prn (avec $F = \Ker P$):

%\vspace{-1mm}

\snac{\xymatrix @C = 0.4cm
{
\rT_P(\GA_{n},\gA)\ar[d] \ar@{-}[r]^(0.39){_\sim}
&
\Lin_\gA(E,F)\oplus\Lin_\gA(F,E)\ar@{-}[r]^(0.42){_\sim} &
\{H \in \MM_{n}(\gA) \mid H = HP + PH \}
     \ar[d]^{ H \mapsto K = HP}
\\
\rT_E(\GGn, \gA) \ar@{-}[r]^(0.45){_\sim} &
\Lin_\gA(E, \Ae {n}\!/E) \ar@{-}[r]^(0.33){_\sim} &
\{K \in \MM_{n}(\gA) \mid KP = K,\, PK = 0\}
\\
}}
\vspace{-.2cm}\hfill\eoe
%$$\quad
%\xymatrix @C = 4cm {
%\rT_P(\GA_{n}, \gA) \simeq \{H \in \MM_{n}(\gA) \mid H = HP + PH \}
%   \ar[d]^{\textstyle H \mapsto K = HP} \\
%\rT_E(\GGn, \gA) \simeq \Lin_\gA(E, \Ae {n}/E) \simeq
%\{K \in \MM_{n}(\gA) \mid KP = K,\ PK = 0\}\\
%} \quad\matrix{ \cr  \cr \cr  \cr \cr  \cr \cr \eoq}
%$$

%   Section{subsecClassifMptfs}----
\section[Groupes de Grothendieck et de Picard]{Classification des \mptfsz,
groupes de Grothendieck et de Picard}
\label{subsecClassifMptfs}
%-----------------------------------------

Nous attaquons ici le \pb \gnl de la classification complète des \mptfs
sur un anneau $\gA$ fixé.

Cette classification est un \pb fondamental mais difficile, qui n'admet pas de
solution \algq \gnlez.

Nous commençons par poser quelques jalons pour le cas où tous les \mrcs sont
libres.

Nous donnons dans les sous-sections suivantes une toute petite introduction à
des outils
classiques qui permettent d'appréhender le \pb \gnlz.

%:--- Subsection{Quand les mprcs sont libres}
\subsec{Quand les \mrcs sont libres}\label{subsecMrcLibre}
%-----------------------------------------

Commençons par une remarque \elrz.
%--- Fact{factRgcstLib}-------------
\begin{fact}
\label{factRgcstLib}
Un \Amo \pro de rang $k$ est libre \ssi il est engendré par $k$ \eltsz.
\end{fact}
%--- end-fact-----------------------------------------
%-----------------begin proof------------------
\begin{proof}
La condition est clairement \ncrz. Supposons maintenant le
module engendré par $k$ \eltsz.
Le module est donc image d'une
matrice de \prn $F\in\Mk(\gA)$. Par hypothèse $\det(\I_k+XF)=(1+X)^k$. En
particulier, $\det F=1$, donc $F$ est \ivz, et puisque $F^2=F$, cela donne
$F=\I_k$.
\end{proof}
%-----------------end proof------------------

Voici une autre remarque facile.
%--- Fact{factRgcstLib2}-------------
\begin{fact}
\label{factRgcstLib2}
Tout \Amrc est libre \ssi tout \Amo \pro est quasi libre.
\end{fact}
%--- end-fact-----------------------------------------
%-----------------begin proof------------------
\begin{proof}
La condition est clairement suffisante. Si tout \Amrc est libre
et si $P$ est \pro soit $(r_0,\ldots ,r_n)$ le \sfio correspondant.
Alors $P_k=r_kP\oplus (1-r_k)\Ae k$ est un \Amo \pro de rang $k$ donc libre.
Soit une base $B_k$, la \gui{composante} $r_kB_k$ est dans $r_kP$, et $r_kP\simeq
\left(r_k\gA\right)^k$. Puisque $P$ est la somme directe des $r_kP$, il est bien
quasi libre.
\end{proof}
%-----------------end proof------------------

%:     Proposition{thZerdimLib}
\begin{proposition}\label{thZerdimLib}
Tout \mrc sur un \algb est libre. 
\end{proposition}
%--------- fin proposition ---------------------------------------------- 

%-----------------begin proof------------------
\begin{proof}
Déjà vu dans le \thrf{thlgb2}.
\end{proof}
%-----------------end proof------------------

%--- Theorem{thBézoutLib}----------
\begin{theorem}
\label{thBézoutLib}
Tout \mptf
sur un anneau de Bézout intègre est libre.
Tout \mptf de rang constant sur un anneau de Bézout \qi est libre.
\end{theorem}
%--- end-theorem-----------------------------------------
%-----------------begin proof------------------
\begin{proof}
Voyons le cas intègre. Une matrice de \pn du module peut être ramenée à
la forme  $\cmatrix{T&0\cr0&0}$ où $T$ est triangulaire
avec des \elts \ndzs sur la diagonale (voir l'exercice \ref{exoBézoutstrict}).
Comme les \idds de cette matrice sont \idms le \deter $\delta$ de $T$ est
un \elt \ndz qui vérifie $\delta\gA=\delta^2\gA$. Ainsi $\delta$ est \iv et la matrice de \pn est \eqve à
$\cmatrix{\I_k&0\cr0&0}$.\\
Pour le cas \qi on applique la machinerie locale-globale \elr expliquée
\paref{MethodeQI}.\imlgz
\end{proof}
%-----------------end proof------------------

Signalons un autre cas important:  $\gA=\BXn$ où $\gB$ est un
anneau de Bézout intègre.
Ceci est une extension remarquable du \tho de Quillen-Suslin, due à 
Bass (pour $n=1$), puis Lequain et Simis \cite{LS}.
Le \tho sera démontré dans la section~\ref{sec.Etendus.Valuation}.
\perso{Annoncer ici tout ce qui sera démontré (constructivement) dans le livre?}

%:--- Subsection{GKO et autres}---------
\subsec{\texorpdfstring{$\GKO(\gA)$, $\KO(\gA)$, $\KTO(\gA)$,  et
$\Pic(\gA)$} {GK0(A), K0(A) et Pic(A)}}
\label{secGKO}
%--------------------------------------

On note $\GKO \gA$ l'ensemble des classes d'\iso de \mptfs sur $\gA$.
C'est un semi-anneau pour les lois héritées de $\oplus$ et $\otimes$.
Le  \textsf{G}  de $\GKO$ est en hommage à Grothendieck. \label{NOTAGKO}

Tout \elt de $\GKO \gA$ peut être représenté par une matrice \idme à
\coes dans~$\gA$. Tout \homo d'anneaux $\varphi:\gA\to\gB$ induit
un \homo $\GKO \varphi : \GKO \gA\to \GKO \gB$.
Ceci fait de $\GKO$ un foncteur covariant de la catégorie des
anneaux commutatifs vers la catégorie des semi-anneaux.
On a $\GKO (\gA_{1}\times\gA_{2})  \simeq  \GKO \gA_{1} \times \GKO \gA_{2}$.
Le passage d'un module \pro à son dual définit un \auto
involutif de  $\GKO \gA$.

Si $P$ est un \Amo \ptf on peut noter $[P]_{\GKO \gA}$ l'\elt
de $\GKO \gA$ qu'il définit.

Le sous-semi-anneau de  $\GKO \gA$ engendré par $1$ (la classe du \mptf $\gA$)
est isomorphe à $\NN$, sauf dans le cas où $\gA$ est l'anneau trivial.
Comme sous-semi-anneau de  $\GKO \gA$ on a aussi celui  engendré par les classes d'\iso des
modules $r\gA$ où $r\in\BB(\gA)$,  isomorphe à $\HOp (\gA)$.
On obtient facilement l'\iso $\HOp (\gA)\simeq\GKO\big(\BB(\gA)\big)$. Par ailleurs, le rang
définit un \homo surjectif de semi-anneaux  $\GKO \gA\to \HOp (\gA)$,
et les deux \homos $\HOp (\gA)\to\GKO \gA\to \HOp (\gA)$
se composent selon l'identité.

Le \ix{groupe de Picard} $\Pic \gA$ est le sous-ensemble de $\GKO \gA$ formé par les classes d'\iso des \mrcs 1.
D'après les propositions \ref{propRANG1} et \ref{th ptrg1}
il s'agit du groupe des \elts \ivs
du semi-anneau  $\GKO \gA$
(l'\gui{inverse} de $P$ est le dual
de $P$). \label{NOTAPic}

%:HHH rajout rdb
\rdb
Le \mo additif (commutatif) de $\GKO \gA$ n'est pas toujours régulier.
Pour obtenir un groupe, on symétrise le \mo additif $\GKO\gA$ et l'on obtient
le \ix{groupe de Grothendieck} que l'on note $\KO \gA$.\label{NOTAK0}

%:HHH petit rajout: , ou $[P]_{\gA}$, ou même $[P]$
La classe du \mptf $P$ dans $\KO \gA$ se note $[P]_{\KO (\gA)}$, ou~$[P]_{\gA}$, ou même $[P]$ si le contexte le permet. 
Tout \elt de $\KO \gA$ s'écrit sous forme $[P]-[Q]$.
Plus \prmtz, il peut se représenter sous les deux formes
\begin{itemize}
\item [$\bullet$] [projectif] - [libre] d'une part,
\item [$\bullet$] [libre] - [projectif] d'autre part.
\end{itemize}
En effet:

\snic{[P] - [Q] = [P \oplus P'] - [Q \oplus P']
          = [P \oplus Q'] - [Q \oplus Q'],}

%\sni
avec au choix $P \oplus P'$ ou
$Q \oplus Q'$ libre.

 Le produit défini dans $\GKO \gA$ donne par passage au quotient un
produit dans $\KO \gA$, qui a donc une structure d'anneau
commutatif\footnote{Lorsque
l'anneau $\gA$ n'est pas commutatif, il n'y a plus de structure
multiplicative sur $\GKO \gA$. Cela explique que la terminologie
usuelle soit celle de groupe de Grothendieck
et non d'anneau de Grothendieck.}.

Les classes de deux \mptfs $P$ et $P'$ sont égales dans~$\KO \gA$ \ssi
il existe un entier $k$ tel que $P\oplus\Ae k\simeq P'\oplus\Ae k$.
On dit dans ce cas que $P$ et $P'$ sont \ixc{stablement isomorphes}{modules ---}.

Deux modules quasi libres stablement isomorphes sont isomorphes, de sorte
que $\HO \gA $ %(le symétrisé de $\HOp\gA$) 
s'identifie à un sous-anneau de $\KO\gA$.
%:HHH rajout pas de conflit
Et lorsque $P$ est quasi libre, il n'y a pas conflit entre les deux notations $[P]_\gA$ (ci-dessus et \paref{notaHO+}).

\rdb\label{NOTAKTO}
Deux \mptfs stablement isomorphes $P$ et $P'$ ont même rang
puisque $\rg(P\oplus\Ae k)=k+\rg(P)$. 
%:HHH inutile et $\HO\gA$ est un groupe additif.
En conséquence, le rang (\gnez) des \mptfs définit un \homo
surjectif d'anneaux
$\rg_\gA:\KO \gA\to\HO \gA$.
On note $\KTO \gA$ son noyau. Les deux \homosz~$\HO \gA\to\KO \gA\to\HO \gA$
se composent selon l'identité, autrement dit l'application $\rg_\gA$
est un \crc de la $\HO(\gA)$-\alg $\KO\gA$ et l'on peut écrire

\snic{\KO (\gA)=\HO (\gA)\oplus \KTO (\gA).}

\smallskip
La structure de  l'\id $\KTO \gA$ de $\KO \gA$
concentre une bonne part du mystère des classes d'\iso stable des \mptfsz,
puisque~$\HO \gA$ ne présente aucun mystère (il est complètement
décrypté par~$\BB(\gA)$). 
%:HHH simplifié
%En fait, comme déjà remarqué pour $\KO\gA$ tout \elt de $\KTO \gA$
%peut s'écrire sous forme $[P]-[\gA^p]=[P]-p$ 
% et l'on obtient le résultat suivant (voir le \pb \ref{exoLambdaGammaK0}).
Dans ce cadre le résultat suivant peut être utile (cf. \pbz~\ref{exoLambdaGammaK0}).

%:     Proposition{propKTO}
\begin{proposition}\label{propKTO}
L'\id  $\KTO \gA$ est le nilradical de $\KO \gA$.
\end{proposition}
%--------- fin proposition ---------------------------------------------- 
%

Notons enfin que si $\rho:\gA\to\gB$ est un \homo d'anneaux, on
obtient des \homos corrélatifs

\snic{
\KO \rho:\KO \gA\to\KO \gB, \quad \KTO \rho:\KTO \gA\to\KTO \gB\quad \hbox{et}
\quad \HO \rho:\HO \gA\to\HO \gB.}

%\sni
Et  $\KO$, $\KTO$ et $\HO$ sont des
foncteurs.

%:--- subsection{subsecPicGp}---------
\subsec{Le groupe de Picard}
\label{subsecPicGp}

Le groupe de Picard n'est pas affecté par le passage aux classes d'\iso stable, en raison du fait suivant.

%--- Fact{factPicStab}---------------
\begin{fact}
\label{factPicStab}
Deux \mrcs $1$  stablement isomorphes sont isomorphes.
En particulier, un module stablement libre de rang $1$ est libre.
Plus \prmtz, pour un module $P$ \prc $1$  on a
\begin{equation}
\label{eqfactPicStab0}
P\simeq\Al{k+1}(P\oplus\Ae k) .
\end{equation}
En particulier, $\Pic\gA$ s'identifie à un sous-groupe de $(\KO\gA)\eti$.
\end{fact}
%--- end-fact-----------------------------------------
%-----------------begin proof------------------
\begin{proof}
Voyons l'\isoz: cela résulte des \isos \gnls  donnés dans la preuve
de la proposition \ref{prop puissance ext} (\eqrf{eqVik}).
Pour des \Amos arbitraires
$P$, $Q$, $R$, \ldots, la considération de la \prt universelle
qui définit les puissances extérieures conduit à:
%--------------------begin array---------------
$$
\preskip.3em \postskip.3em
\arraycolsep2pt\begin{array}{rcl}
\Al2(P\oplus Q)& \simeq  & \Al2P\oplus(P\te Q)\oplus \Al2Q,
\\[1.4mm] \mathrigid 1mu
\Al3(P\oplus Q\oplus R)& \simeq  & \Al3P\oplus\Al3Q\oplus\Al3R\oplus\big(\Al2P\te Q\big)
\oplus \cdots \oplus (P\te Q\te R)
 ,
\end{array}
$$
%---------------------end array--------------
avec la formule \gnle suivante en convenant de $\Al{0}(P_i) = \gA$
%------begin equation--eqfactPicStab-----------
\begin{equation}\label{eqfactPicStab}
\Al k\big(\bigoplus\nolimits_{i=1}^mP_i\big)\simeq
\bigoplus_{\som_{i=1}^{m} k_i= k}
\Big(
\big(\, \Al{k_1}P_1 \big)\te\cdots \te
\big(\, \Al{k_m}P_m \big)
\Big).
\end{equation}
%---------------------end equation--------------
En particulier, si $P_1$, \ldots, $P_r$ sont des \mrcsz~1 on obtient
%------begin equation--eqfactPicStab2-----------
\begin{equation}
\label{eqfactPicStab2}
\Al r (P_1\oplus\cdots\oplus P_r) \simeq P_1\otimes \cdots\otimes  P_r.
\end{equation}
%%%%%%%%%%%%%%%%%%%%%%%%%%%%%%%%%%%%%%%%%
Il reste à appliquer ceci avec la somme directe
 $P\oplus\Ae k=P\oplus\gA\oplus \cdots \oplus\gA$.
 L'\iso de l'\eqrf{eqfactPicStab0} est alors obtenu avec l'\Ali
 $P\to \Al{k+1}(P\oplus\Ae k)$,
 $x\mapsto x\vi 1_1 \vi 1_2 \vi \cdots \vi 1_k$,
 où l'indice représente la position dans la somme directe
 $\gA\oplus \cdots \oplus\gA$.\\
 La dernière affirmation est alors claire puisque l'on 
 vient de montrer que l'application $\GKO\gA\to\KO\gA$, restreinte à
 $\Pic\gA$, est injective.
\end{proof}
%-----------------end proof------------------

\rem \Llec pourra comparer le résultat précédent et sa \dem
avec l'exercice~\ref{exoStabLibRang1}.
\eoe

\smallskip 
On en déduit le \tho suivant.
%--- Theorem{thPicKTO}--------------
\begin{theorem}
\label{thPicKTO} \emph{($\Pic \gA$ et $\KTO \gA$)} %\\
Supposons que tout $\gA$-\mrc $k+1$ ($k\geq1$) soit isomorphe à un module $\Ae k\oplus Q$. Alors l'application de $(\Pic \gA,\times )$ dans $(\KTO \gA,+)$
définie par
$$[P]_{\Pic \gA}\mapsto [P]_{\KO\! \gA} - 1_{\KO\! \gA} $$
est un \iso de groupes.
En outre, $\GKO\gA=\KO\gA$ et sa structure est entièrement
connue à partir de celle de $\Pic\gA$.\perso{Il semble qu'en \gnl
l'application ne soit pas un \homo de groupes.}
\end{theorem}
%--- end-theorem-----------------------------------------
\begin{proof}
L'application est injective par le fait \ref{factPicStab},
et surjective par hypothèse. C'est un \homo de groupe
parce que $\gA\oplus(P\otimes Q) \simeq P\oplus Q$,
\egmt en vertu du fait \ref{factPicStab},
puisque 

\smallskip 
\centerline{$\Al2\big(\gA\oplus(P\otimes Q)\big)\simeq P\otimes Q
\simeq \Al2(P\oplus Q).$}
\end{proof}

Notez que la loi de $\Pic \gA$ est héritée du produit tensoriel
tandis que celle de~$\KTO \gA$ est héritée de la somme directe.
Nous verrons au chapitre \ref{chapKrulldim} que l'hypothèse du \tho
est vérifiée pour les
anneaux \ddkz~$\leq1$.

\medskip \rdb
\comm \label{comHOclassique}
On a vu dans la section \ref{subsecCalRang} comment la structure
de $\HO(\gA)$ découle directement de celle de l'\agB $\BB(\gA)$.
\\
Du point de vue des \clama l'\agB $\BB(\gA)$ est l'\alg des
ensembles ouverts et fermés dans $\Spec\gA$ (l'ensemble des \ideps de $\gA$
muni d'une topologie convenable, cf. chapitre \ref{chapKrulldim}).
Un \elt de  $\BB(\gA)$ peut
donc être vu comme la fonction \cara d'un ouvert-fermé de $\Spec\gA$.
Alors la manière dont on construit $\HO(\gA)$ à partir de $\BB(\gA)$
montre que  $\HO(\gA)$ peut être vu comme
 l'anneau des fonctions à valeurs entières, \colis
  entières des \elts dans $\BB(\gA)$.
Il s'ensuit que $\HO(\gA)$ s'identifie à l'\alg des
fonctions localement constantes, à valeurs entières, sur $\Spec\gA$.
Toujours du point de vue des \clama le rang (\gnez) d'un \Amo \ptf $P$ peut être
vu comme la fonction (à valeurs dans $\NN$) définie sur $\Spec\gA$,
de la manière suivante:
à un \idep $\fp$ on associe le rang du module libre $P_\fp$ (sur un \alo
tous les \mptfs sont libres). Et l'anneau $\HO(\gA)$ est bien obtenu simplement en
symétrisant le semi-anneau $\HOp (\gA)$ des rangs de \Amos \ptfsz.
\eoe

%%%%%%%%%%%%%%%%%%%%%%%%%%%%%%%%%%%%%%%%%%%%%%%%%%%%%%%%%%%%%%%%%%%%%%%%%%%
\subsubs{Groupe de Picard et groupe des classes d'idéaux}
Considérons le \mo multiplicatif
des \emph{\ifrs \tfz} de l'anneau $\gA$, formé par les
sous-\Amos \tf de l'anneau total des fractions $\Frac\gA$.
Nous noterons ce \mo $\Ifr\gA$. \label{NOTAIfr}

Plus \gnlt un \emph{\ifrz} de $\gA$ est un sous-\Amo $\fb$ de $\Frac\gA$
tel qu'il existe  $b$ \ndz dans $\gA$ vérifiant $b\,\fb\subseteq\gA$.

En bref on peut voir $\Ifr\gA$ comme le \mo obtenu à partir de celui
des \itfs de $\gA$ en forçant l'inversibilité des \idps
engendrés par des \elts \ndzsz.%
\index{fractionnaire!idéal ---}%
\index{ideal@idéal!fractionnaire}

Un \id $\fa\in\Ifr\gA$ est
parfois dit \emph{entier} s'il est contenu dans $\gA$,
auquel cas c'est un \itf de $\gA$
au sens usuel.

Un \id $\fa$ arbitraire de $\gA$  est \iv comme \id de $\gA$ (au sens de la \dfnz~\ref{defiiv}) \ssi c'est un \elt \iv dans le \mo $\Ifr\gA$.
Inversement tout \id de  $\Ifr\gA$  \iv  dans ce \mo s'écrit $\fa/b$, où
$b\in\gA$ est \ndz et $\fa$ est un \id \iv de $\gA$.
Les  \elts \ivs de $\Ifr\gA$ forment un groupe, le \emph{groupe des \ifrs
\ivs de $\gA$}, que nous noterons $\Gfr\gA$.%
\index{groupe des idéaux fractionnaires inversibles}

En tant que \Amo un \ifr \iv est \prc 1.\label{pageclassgroup}
Deux \ids \ivs sont isomorphes en tant que \Amos
s'ils sont égaux modulo le sous-groupe des \idps \ivs (i.e. engendrés
par un \elt \ndz de $\Frac\gA$). On note $\Cl \gA $ le groupe quotient, que l'on appelle \emph{groupe des classes d'\ids inversibles}, ou parfois simplement le \emph{groupe des classes} de $\gA$, et l'on
obtient une application naturelle bien définie $\Cl \gA\to\Pic\gA$.%
\index{groupe des classes!d'un anneau $\gA$}%
\index{groupe des classes (d'idéaux inversible)}%
\index{classe!d'ideaux inversibles@d'\ids\ivsz}%
\index{classe!d'ideaux@d'\idsz}%
\label{NOTAGfr}

Par ailleurs, considérons un \id $\fa$ entier et \ivz.
Puisque $\fa$ est plat, l'application naturelle $\fa\te_\gA\fb\to\fa\fb$
est un \isoz, ceci pour n'importe quel \idz~$\fb$
(\thref{thplatTens}). Ainsi, l'application~$\Cl \gA\to\Pic\gA$
est un \homo de groupes, et c'est clairement un \homo injectif, donc~$\Cl\gA$ s'identifie à un sous-groupe de $\Pic\gA$.

Ces deux groupes sont souvent identiques comme le montre le \tho
suivant, qui résulte des considérations précédentes
et du \thref{propRgConstant2}.

%:      Th {propRgConstant3}----
\begin{theorem}
\label{propRgConstant3} \emph{(Modules de rang constant $1$ comme \ids
de $\gA$)}\\
Supposons que sur $\Frac \gA$ tout module \pro de rang $1$ soit libre.
\begin{enumerate}
\item Tout \Amo \pro de rang $1$
est isomorphe à un \id \iv de~$\gA$.
\item Tout \id \pro de rang $1$ est \ivz.
\item Le groupe des classes d'\ids \ivs est naturellement isomorphe au groupe de
Picard.
\end{enumerate}
 \end{theorem}
%--- end-theorem----------------------------------------
%
\begin{proof}
Le \tho \ref{propRgConstant2} montre que tout module \pro de rang $1$ est isomorphe à un \id  $\fa$. Il reste donc à voir qu'un tel \id est \ivz.
Puisqu'il est \lop il suffit de montrer qu'il contient un \elt \ndzz.
Pour cela on considère un  \id entier $\fb$ isomorphe à l'inverse de $\fa$
dans $\Pic\gA$. Le produit de ces deux \ids est isomorphe à leur produit
tensoriel (parce que $\fa$ est plat) donc c'est un module libre,
donc c'est un \idp engendré par un \elt \ndzz.
\end{proof}

NB: concernant la comparaison de $\Pic\gA$ et $\Cl\gA$ on trouvera un résultat plus \gnl en exercice \ref{exoPicAPicFracA}.

%%%%%%%%%%%%%%%%%%%%%%%%%%%%%%%%%%%%%%%%%%%%%%%%%%%%%%%%%%%%%%%%%%%%%%%%%%%
%--- subsection{subsecComparRed}---------
\subsect{Les semi-anneaux \texorpdfstring{$\GKO(\gA)$, $\GKO(\Ared)$
et $\GKO(\gA\sur{\Rad \gA})$} {GK0(A), GK0(Ared) et GK0\big(A/Rad(A)\big)}}{Semi-anneaux $\GKO(\gA)$, $\GKO(\Ared)$ et $\GKO(\gA\sur{\Rad \gA})$}
\label{subsecComparRed}
%-----------------------------------------------------

Dans ce paragraphe nous utilisons $\Rad \gA$, le radical de Jacobson de $\gA$,
qui est défini \paref{eqDefRadJac}. Nous comparons les \mptfs définis sur
$\gA,$  ceux définis sur $\gA'=\gA\sur{\Rad \gA }$
et  ceux définis sur $\Ared$.

L'\eds de $\gA$ à $\gB$ transforme un
\mptf défini sur $\gA$ en un \mptf sur $\gB$. Du point de vue
matrices de \prnz, cela correspond à considérer la  matrice
transformée par l'\homo $\gA\to\gB$.

%--- Proposition{propComparRedJac}---
\begin{proposition}
\label{propComparRedJac}
L'\homo naturel $\GKO(\gA)\to\GKO(\gA\sur{\Rad \gA})$ est injectif,
ce qui signifie que si deux \hbox{\mptfs $E$, $F$}  sur~$\gA$
sont isomorphes sur $\gA'=\gA\sur{\Rad \gA}$,
ils le sont \egmt sur $\gA$. De manière plus précise, si
deux matrices \idmes $P$, $Q$ de même format sont
conjuguées sur $\gA'$, elles le
sont \egmt sur $\gA$, via un \iso qui relève l'\iso
de conjugaison résiduel.
\end{proposition}
%--- end-proposition----------------------------------------
%-----------------begin proof------------------
\begin{proof}
On note $\ov{x}$ l'objet $x$  vu modulo $\Rad \gA$. Soit $C\in\Mn(\gA)$
une matrice telle que $\ov{C}$ conjugue $\ov{P}$ à $\ov{Q}$. Puisque $\det C$
est \iv modulo $\Rad \gA$, $\det C$ est \iv dans $\gA$ et
$C\in\GL_n(\gA)$. On a donc $\ov{Q}=\ov{C\,P\,C^{-1}}$. Quitte à remplacer $P$
par $C\,P\,C^{-1}$ on peut supposer $\ov{Q}=\ov{P}$ et $\ov C=\In$.
Dans ce cas on cherche une matrice inversible $A$ telle que $\ov A=\In$ et $APA^{-1}=Q$.
\\
On remarque que $QP$ code une \Ali
de $\Im P$ vers $\Im Q$ qui donne \rdt l'\idtz.
De même
$(\In-Q)(\In-P)$ code une \Ali de  $\Ker P$ vers $\Ker Q$ qui donne
\rdt l'\idtz. En s'inspirant du lemme
\dlg \ref{propIsoIm}
ceci nous conduit à la matrice $A=QP+(\In-Q)(\In-P)$ qui réalise $AP=QP=QA$  et $\ov{A}=\In$, donc $A$ est \iv et $APA^{-1}=Q$.\\
Pour deux \mptfs \rdt isomorphes~$E$ et~$F$ on utilise le
lemme  \dlg
qui permet de réaliser $\ov{E}$ et~$\ov{F}$ comme images de matrices \idmes de
même format conjuguées.
\end{proof}
%-----------------end proof------------------

Pour ce qui concerne la réduction modulo les nilpotents, on obtient en
plus la possibilité de relever tout \mptf en raison du corolaire
\ref{corIdmNewton}. D'où le \tho qui suit.

%:    theorem{propComparRedRed}---
\begin{theorem}
\label{propComparRedRed}
L'\homo naturel $\GKO(\gA)\to\GKO(\Ared)$ est un \isoz. De manière plus
précise, on a les résultats suivants.
%-----------------begin item------------------
\begin{enumerate}
\item 
\begin{enumerate}
\item Toute matrice \idme sur $\Ared$ se relève en
une matrice \idme sur~$\gA$.

\item   Tout \mptf sur $\Ared$ provient d'un \mptf sur~$\gA.$

\end{enumerate}
\item

\begin{enumerate}
\item Si deux matrices \idmes  de même format sont conjuguées sur
$\Ared$, elles le sont \egmt sur $\gA$, via un \iso qui relève l'\iso
de conjugaison résiduel.

\item   Deux \mptfs sur $\gA$ isomorphes sur $\Ared$ le sont \egmt sur $\gA.$

\end{enumerate}
\end{enumerate}
%-----------------end item------------------
\end{theorem}
%--- end-proposition----------------------------------------

%%%%%%%%%%%%%%%%%%%%%%%%%%%%%%%%%%%%%%%%%
%--- Subsection{Carre Milnor}---------
\subsec{Le carré de Milnor}
Un carré commutatif (dans une catégorie arbitraire) du style suivant
$$
\xymatrix @C=1.2cm{
A\,\ar[d]^{i_1}\ar[r]^{i_2}   & \,A_2\ar[d]^{j_2}   \\
A_1\,\ar[r]    ^{j_1}    & \,A'  \\
}
$$
est appelé un \ix{carré cartésien} s'il définit $(A,i_1,i_2)$ comme la
limite projective de $(A_1,j_1,A')$, $(A_2,j_2,A')$.
Dans une catégorie équationnelle
on peut prendre \[A=\sotq{(x_1,x_2)\in A_1\times A_2}{j_1(x_1)=j_2(x_2)}.\]

\Llec vérifiera par exemple qu'étant donné $\gA\subseteq \gB$ et $\ff$ un \id
de $\gA$ qui est aussi un \id de $\gB$ (autrement dit, $\ff$ est contenu dans le conducteur de
$\gA$ dans $\gB$), on a un carré cartésien d'anneaux commutatifs,
défini ci-dessous:
$$
\xymatrix @C=1.2cm{
\gA\,\ar@{->>}[d]\ar[r]   & \,\gB\ar@{->>}[d]   \\
\gA\sur{\ff}\,\ar[r]    & \gB\sur{\ff} \\
}
 $$

Soit  un \homo  $\rho:\gA\to\gB$, $M$ un \Amo et $N$ un \Bmoz. Rappelons
qu'une \Ali $\alpha
:M\to N$ est  un \emph{morphisme d'\edsz}
(cf.\ \dfn \ref{defAliAliExtScal}) \ssi l'\Bli
naturelle  $\rho\ist(M)\to N$ est un \isoz.

%:2012   Dans toute ce paragraphe !
Dans tout ce paragraphe nous considérons dans la catégorie des anneaux
commutatifs le \gui{carré de Milnor} ci-dessous
à gauche, noté $\cA$, dans lequel~$j_2$ est surjective,
$$
%\xymatrix @C=1.2cm{
%\gA\,\ar[d]^{i_1}\ar[r]^{i_2}   & \,\gA_2\ar@{->>}[d]^{j_2}   \\
%\gA_1\,\ar[r]    ^{j_1}    & \,\gA'  \\
%}
\xymatrix @C=1.2cm{
\ar@{}[dr]|*+<10pt>[o][F]{\cA}
\gA\,\ar[d]^{i_1}\ar[r]^{i_2}   & \,\gA_2\ar@{->>}[d]^{j_2}   \\
\gA_1\,\ar[r]    ^{j_1}    & \,\gA'  \\
}
\quad
\xymatrix @C=1.2cm{
M\,\ar[d]^{\psi_1}\ar[r]^{\psi_2}   & \,M_2\ar@{->>}[d]^{\varphi_2}   \\
M_1\,\ar[r]    ^{\varphi_1}    & \,M'  \\
}
\quad
\xymatrix @C=1.2cm{
E\,\ar[d]\ar[r]   & \,E_2\ar@{->>}[d]^{{j_2}\ist}   \\
E_1\,\ar[r]^{h\circ {j_1}\ist}    & \,E'  \\
}
 $$
\'Etant donnés un \Amo $M$, un $\gA_1$-module $M_1$, un  $\gA_2$-module $M_2$, \hbox{un $\gA'$-module $M'$}  et un carré cartésien de \Amos comme ci-dessus
au centre, ce dernier est dit \emph{adapté à $\cA$}, si les $\psi_i$ et
$\varphi_i$
sont des morphismes d'\edsz.

\'Etant donnés un $\gA_1$-module $E_1$, un  $\gA_2$-module $E_2$, et
un \iso \hbox{d'$\gA'$-modules}

\snic{h:{j_1}\ist(E_1)\to {j_2}\ist(E_2)=E',}

%\sni
nous notons
$M(E_1,h,E_2)=E$ (ci-dessus à droite) le \Amo limite projective de

\snic{\big(E_1,h\circ {j_1}\ist,{j_2}\ist(E_2)\big),\big(E_2,{j_2}\ist,{j_2}\ist(E_2)\big)}

\smallskip
Notons qu'a priori le carré cartésien obtenu n'est pas \ncrt adapté
à $\cA$.

%--- Theorem{thCarreMil1}----------
\begin{theorem}
\label{thCarreMil1} \emph{(\Tho de Milnor)}
%-----------------begin item------------------
\begin{enumerate}
\item On suppose que $E_1$ et $E_2$ sont \ptfsz, alors:
%-----------------begin enum------------------
\begin{enumerate}
\item $E$ est \ptfz,
\item le carré cartésien est adapté à $\cA$: les \homos naturels
${j_k}\ist(E)\to E_k$ ($k=1,2)$ sont des \isosz.
\end{enumerate}
%-----------------end enum------------------
\item Tout \mptf sur $\gA$ est obtenu (à \iso près) par ce procédé.
\end{enumerate}
%-----------------end item------------------
\end{theorem}
Nous aurons besoin du lemme suivant.
%--- Lemma{lemCarreMil1}---------
\begin{lemma}
\label{lemCarreMil1}
Soient $A\in\Ae{m\times n}$, $A_k= i_k(A)$ ($k = 1, 2$),
$A'=j_1(A_1)=j_2(A_2)$, $K=\Ker A\subseteq\Ae n $,  $K_i=\Ker A_i$ ($i=1,2$),
 $K'=\Ker A'$. Alors $K$ est la limite projective (comme \Amoz)
de $K_1\to K'$ et $K_2\to K'$.
\end{lemma}
%--- end-lemma-----------------------------------------
%-----------------begin proof------------------
\begin{proof}
Soient $x\in\Ae n $, $x_1 = {j_1}\ist(x) \in \gA_1^n$, 
$x_2 = {j_2}\ist(x) \in \gA_2^n$. Puisque $x\in K$ \ssi $x_i\in K_i$ pour $i=1,2$,
$K$ est bien la limite projective convoitée.
\end{proof}
%-----------------end proof------------------
\Llec remarquera que le lemme ne s'appliquerait
pas en \gnl pour les sous-modules images des matrices.
%-----------------begin proof------------------
\begin{Proof}{\Demo du \thref{thCarreMil1}. }
\emph{2.} Si $V\oplus W=\Ae n $, soit $P$ la matrice de
projection sur $V$ \paralm à $W$. Si $V_1$, $V_2$, $V'$ sont les modules
obtenus par \eds à $\gA_1$, $\gA_2$ et $\gA'$, ils
s'identifient aux
noyaux des matrices $P_1=i_1(\In-P)$,  $P_2=i_2(\In-P)$,
$P'=j_2(\In-P_2)=j_1(\In-P_1)$, et le lemme s'applique: $V$ est la limite
projective de $V_1\to V'$ et  $V_2\to V'$. L'\iso $h$ est alors $\Id_{V'}$.
Ce \gui{miracle} se produit gr\^ace à l'identification de ${j_i}\ist(V_i)$ et
$\Ker P_i$.

%\sni
\emph{1a.} Soit $P_i\in\MM_{n_i}(\gA_i)$ un \prr
d'image isomorphe à $E_i$ ($i=1,2$). On dispose d'un \iso d'$\gA'$-modules %,
%que l'on note encore $h$ 
 de $\Im\big(j_1(P_1)\big)\in\MM_{n_1}(\gA')$ 
 \hbox{vers $\Im\big(j_2(P_2)\big)\in\MM_{n_2}(\gA')$}. 
 Notons $n=n_1+n_2$. D'après le lemme \dlgz~\ref{propIsoIm} il existe une matrice $C\in{\En(\gA')}$
 réalisant la conjugaison

\snic{
\Diag(j_1(P_1),0_{n_2}) = C\,\,\Diag\big(0_{n_1},j_2(P_2)\big)\,\,C^{-1}.
}

%\sni
Puisque $j_2$ est surjective (ah ah !), $C$  se relève 
en une matrice $C_2\in\En(\gA_2)$. Posons

\snic{Q_1=\Diag(P_1,0_{n_2}),$ $\;Q_2=C_2\,\Diag(0_{n_1},P_2)\,C_2^{-1},}

%\sni
de sorte
que $j_1(Q_1)=j_2(Q_2)$ (pas mal, n'est-ce pas). Il existe alors une unique
matrice $Q\in\Mn(\gA)$ telle que $i_1(Q)=Q_1$ et $i_2(Q)=Q_2$. 
L'unicité de~$Q$  assure $Q^2=Q$, et le lemme précédent s'applique pour montrer que~$\Im Q$ est isomorphe à~$E$ (chapeau, M. Milnor!).

%\sni
\emph{1b.} Résulte du fait que $Q_k=i_k(Q)$ et $\Im Q_k\simeq \Im P_k\simeq E_k$
pour $k=1,2$.
\end{Proof}
%-----------------end proof------------------

Le fait suivant est purement catégorique et abandonné à la bonne volonté
\dlecz.
%--- Fact{lemCarMil2}-----------
\begin{fact}
\label{factCarMil2}
\'Etant donné deux carrés cartésiens adaptés à $\cA$
comme ci-dessous, il revient au même de se donner une \aliz~$\theta$ de~$E$ vers~$F$
ou de se donner trois \alis (pour les anneaux correspondants) $\theta_1:E_1\to
F_1$,
 $\theta_2:E_2\to F_2$ et  $\theta':E'\to F'$ qui rendent les carrés adéquats
commutatifs.
$$
\xymatrix@R=0.5cm@C=1.5cm{
E\ar[dd]\ar@{-->}[drr]_>>>>>>>>>>>>>>\theta\ar[r]   &
E_1\ar[dd]\ar[drr]^{\theta_1}             \\
            &            &  F\ar[r]\ar[dd]    & F_1\ar[dd] \\
E_2\ar[drr]_{\theta_2}\ar[r]        &  E'\ar[drr]^<<<<<<<<<<{\theta'}        &
\\
            &            &  F_2\ar[r]         & F'           \\
}
$$
\end{fact}
%--- end-fact-----------------------------------------

%--- Corollary{corthCarreMil1}------
\begin{corollary}
\label{corthCarreMil1} ~\\
On considère deux %\Amos \ptfs 
modules $E=M(E_1,h,E_2)$ et $F=M(F_1,k,F_2)$ comme dans le
\tho \ref{thCarreMil1}.
Tout \homo $\psi$ de $E$ dans $F$ est obtenu à l'aide de deux \homos de
$\gA_i$-modules $\psi_i:E_i\to F_i$ compatibles avec $h$ et $k$ en ce sens que
le diagramme ci-dessous est commutatif. L'\homo $\psi$ est un \iso \ssi
$\psi _1$ et $\psi _2$ sont des \isosz.
$$
\xymatrix @C=1.2cm{
{j_1}\ist(E_1)\,\ar[d]^{h}\ar[r]^{{j_1}\ist(\psi_1)}   & \, {j_1}\ist(F_1)
\ar[d]^{k}   \\
{j_2}\ist(E_2)\,\ar[r]^{{j_2}\ist(\psi_2)}        & \,{j_2}\ist(F_2)  \\
}
$$
\end{corollary}
%--- end-corollary------------------------------------

%  Section{secAppliIdenti}---------
\section[Identification de points
sur la droite affine]{Un exemple non trivial: identification de points
sur la droite affine}
\label{secAppliIdenti}

\vspace{3pt}

%:  subsec{Préliminaires}-----
\subsec{Préliminaires}

On considère un anneau commutatif $\gk$, la droite affine sur $\gk$
\emph{corres\-pond à} la \klg $\gk[t]=\gB$.
\'Etant donné $s$ points $\alpha_1$, \ldots, $\alpha_s$ de $\gk$ et des ordres de
multiplicité $e_1$, \ldots, $e_s\geq 1$,
on définit formellement une \klg $\gA$ qui représente le résultat
de l'identification de ces points avec les multiplicités données.
$$
\gA=\sotQ{ f\in\gB}{ f(\alpha_1)=\cdots =f(\alpha_s),\,\,%\&_{1\leq \ell<e_i 
f^{[\ell]}(\alpha_i)=0,\,\ell\in\lrb{1,e_i}, \,i\in\lrbs   }
$$
Dans cette \dfn $f^{[\ell]}$ représente la \emph{dérivée de Hasse} du \pol
$f(t)$
\cad $f^{[\ell]}=f^{(\ell)}\sur{\ell!}$ (formellement, car la \cara peut être finie). Les dérivées de Hasse permettent d'écrire une formule de Taylor pour
n'importe quel anneau $\gk$.
\index{derivee@dérivée de Hasse}

On pose $e=\sum_i e_i$, $x_0=\prod_i(t-\alpha_i)^{e_i}$ et $x_{\ell}=t^\ell\,x_0$
pour $\ell\in\lrb{0..e-1}$. On suppose $e>1$ sans quoi $\gA=\gB$.
Il est clair que les $x_\ell$ sont dans $\gA$.

On suppose aussi que les $\alpha_i-\alpha_j$ sont \ivs pour $i\neq j$.
On a alors par le \tho chinois un \homo surjectif

\snic{\varphi:\gB\to\prod_{i}\big(\vphantom{2^{2^2}}\aqo{\gk[t]}{(t-
\alpha_i)^{e_i}}\big)}

%\sni
dont le noyau est le produit des \idps  $(t-
\alpha_i)^{e_i}\gB$, \cad l'\id $x_0\gB$.

%--- Lemma{lemIdentiMil}-------------
\begin{lemma}
\label{lemIdentiMil}~
%-----------------begin enum------------------
\begin{enumerate}
\item $\gA$ est une \klg \tfz, plus \prmtz:
$\gA=\gk[x_0,\ldots ,x_{e-1}].$
\item $\gB=\gA\oplus \bigoplus_{1\leq \ell<e}\gk \,t^\ell$ en tant que \kmoz.
\item Le conducteur de $\gA$ dans $\gB$, $\ff=(\gA:\gB)$ est donné par
$$
\ff=\gen{x_0}_\gB=\gen{x_0,\ldots ,x_{e-1}}_\gA .
$$
\end{enumerate}
%-----------------end enum------------------
\end{lemma}
%--- end-lemma-----------------------------------------
%-----------------begin proof------------------
\begin{proof}
Soit $f\in\gB$, on l'écrit \gui{en base $x_0$},
 $f=r_0+r_1 x_0 + r_2 x_0^2+\cdots $ \linebreak
avec $\deg r_i<\deg x_0 =e$.
Pour $i\geq 1$, en écrivant $r_ix_0^i=(r_i\,x_0)x_0^{i-1}$ on voit que
$r_ix_0^i\in\gk[x_0,\ldots,x_{e-1}]$.
Ceci prouve que
$$\ndsp\gB=\gk[x_0,\ldots,x_{e-1}]+ \left(\bigoplus_{1\leq \ell<e}\gk
\,t^\ell\right).
$$
Soit $f\in\gA$ que l'on écrit $g+h$ dans la \dcn
précédente. On a donc~$h$ dans $\gA$, et si $\beta $ est la valeur commune des
$h(\alpha_i)$, on obtient 
\linebreak 
l'\egt  $\varphi (h-\beta)=0$. Donc $h-\beta \in x_0\gB$, et
puisque $h\in \bigoplus_{1\leq \ell<e}\gk \,t^\ell$
(le \kmo des \pols de degré $<e$ et sans terme constant), on \linebreak
obtient
$h-\beta=0$ puis $h=\beta=0$, donc $f\in\gk[x_0,\ldots ,x_{e-1}]$.
\\
En conclusion $\gA=\gk[x_0,\ldots ,x_{e-1}]$, les points \emph{1} et \emph{2} sont
démontrés.
\\
En  multipliant l'\egt du point \emph{2} par $x_0$ on obtient

\snic{x_0\gB=x_0\,\gA\;\oplus\;\bigoplus_{\ell\in\lrb{1..e-1}}x_\ell\,\gk,}

%\sni
puis l'\egt $x_0\,\gB=\gen{x_0,\ldots ,x_{e-1}}_\gA$,
ce qui implique $x_0\,\gB\subseteq\ff$. 
Soit 
\linebreak
enfin $f\in\ff$,
et donc $f\in\gA$, et $f=\lambda +g$ avec $\lambda \in\gk$ et $g\in
\gen{x_0,\ldots ,x_{e-1}}_\gA$. On en déduit que $\lambda \in\ff$, ce qui
implique $\lambda =0$: en effet,  $\lambda t\in\gA$, si $\beta$ est la valeur
commune des $\lambda \alpha_i$, on a $\varphi(\lambda t-\beta )=0$,
donc $\lambda t-\beta\in x_0\,\gB$, et puisque~$x_0$ est un \polu de degré
$\geq 2$, $\lambda =0$.
\end{proof}
%-----------------end proof------------------

%  SUBSUBsection{Un carré de Milnor}-----
\subsubsection*{Un carré de Milnor}

Dans la situation décrite dans le paragraphe précédent
on a le carré de Milnor suivant:
$$
\xymatrix @C=1.2cm{
\gA\,\ar@{->>}[d]\ar[r]   & \,\gB\ar@{->>}[d]^\varphi&{\kern-125pt}=\gk[t]   \\
\gk=\gA\sur{\ff}\,\,\ar@{>->}[r]^\Delta    & \gB\sur{\ff}&
{\kern-40pt}
\simeq \prod_i\left(\aqo{\gk[t]}{(t-\alpha_i)^{e_i}}\right)
\\
}
 $$

Dans la suite nous nous intéressons aux \Amrcsz~$r$ obtenus en recollant
le \Bmo  $\gB^r$ et le \kmo  $\gk^r$ à l'aide d'un~$(\gB\sur{\ff})$-\iso

\snic{h:\Delta \ist(\gk^r)\to\varphi \ist(\gB^r),}

%\sni
comme décrit avant le \thrf{thCarreMil1}.
\\
Nous avons noté $M(\gk^r,h ,\gB^r)$  un tel \Amoz.

En fait, $\Delta \ist(\gk^r)$ et $\varphi \ist(\gB^r)$ s'identifient tous les deux
à $(\gB\sur{\ff})^r$, et l'\iso $h$ s'identifie à un \elt de

\snic{\GL_r(\gB\sur{\ff})\simeq \prod_{i=1}^s\GL_r(\aqo{\gk[t]}{(t-
\alpha_i)^{e_i}}).}

%\sni
Nous utiliserons ces identifications dans la suite sans plus les
mentionner, et,
pour des raisons de commodité, nous coderons $h^{-1}$ (et non pas $h$) par les
$s$ matrices $H_i$ correspondantes (avec $H_i\in \GL_r(\aqo{\gk[t]}{(t-
\alpha_i)^{e_i}})$).   Et le module  $M(\gk^r,h ,\gB^r)$  sera noté
$M(H_1,\ldots ,H_s)$.

Dans le cas où les \mrcs sur $\gk$ et $\gB=\gk[t]$ sont toujours libres, le \tho de
Milnor  affirme que l'on obtient ainsi (à \iso près) tous les \mrcs $r$
sur~$\gA$.

Dans le paragraphe qui suit nous donnons une description
com\-plète de la catégorie des \mrcs sur $\gA$ obtenus par de tels recollements,
dans un cas particulier. Celui où toutes les multiplicités
sont égales à $1$.

%: subsec{sans multiplicites}---
\subsec{Identification de points sans multiplicités}

On applique maintenant les conventions précédentes en supposant que
les multiplicités $e_i$ sont toutes égales à~1.
%--- Theorem{thIdenPtsSimples}----
\begin{theorem}
\label{thIdenPtsSimples}
Avec les conventions précédentes.
%-----------------begin enum------------------
\begin{enumerate}
\item Le module $M(H_1,\ldots ,H_s)$, (avec $H_i\in \GL_r(\aqo{\gk[t]}{t-
\alpha_i})\simeq\alb \GL_r(\gk)$) s'identifie au sous-\Amo  $M'(H_1,\ldots ,H_s)$
de $\gB^r$ constitué des \eltsz~$f$ de~$\gB^r$ tels que

\snic{\forall 1\leq i< j\leq s,\;\;H_i\cdot f(\alpha_i)=H_j\cdot f(\alpha_j).}

%\sni
En particulier,   $M'(H_1,\ldots
,H_s)=M'(HH_1,\ldots ,HH_s)$ si $H\in\GL_r(\gk)$.
\item Soient, pour $i\in\lrbs$,  
\[ 
\begin{array}{ccc} 
G_i\in\alb
\GL_{r_1}(\aqo{\gk[t]}{t-\alpha_i})\simeq \GL_{r_1}(\gk)&
\hbox{et}
\\[1mm] 
H_i\in\alb \GL_{r_2}(\aqo{\gk[t]}{t-\alpha_i})\simeq
\GL_{r_2}(\gk).    
\end{array}
\]
Une \Ali $\phi$  de  $M(G_1,\ldots ,G_s)$
%(avec $G_i\in\alb
%\GL_{r_1}(\aqo{\gk[t]}{t-\alpha_i})\simeq \GL_{r_1}(\gk)$)
vers
 $M(H_1,\ldots ,H_s)$
% (avec $H_i\in\alb \GL_{r_2}(\aqo{\gk[t]}{t-\alpha_i})\simeq
%\GL_{r_2}(\gk)$)
peut être codée par une matrice $\Phi\in \gB^{r_2\times r_1}$
vérifiant, pour  $1\leq i< j\leq s$,
%------begin equation--eqthIdenPtsSimples-----------
\begin{equation}\label{eqthIdenPtsSimples}
H_i\cdot \Phi(\alpha_i)\cdot G_i^{-1} = H_j\cdot \Phi(\alpha_j)\cdot G_j^{-1}.
\end{equation}
%---------------------end equation--------------
 Une telle matrice envoie  $M'(G_1,\ldots ,G_s)$ dans $M'(H_1,\ldots ,H_s)$.
L'\Ali $\phi$ est un \iso \ssi $r_1=r_2$ et les $\Phi(\alpha_i)$ sont \ivsz.
\end{enumerate}
%-----------------end enum------------------
\end{theorem}
%--- end-theorem-----------------------------------------
%-----------------begin proof------------------
\begin{proof}
Le premier point n'a pas d'incidence sur les
résultats qui suivent, et il est laissé \alecz. Le deuxième point est un
conséquence immédiate du lemme \ref{lemCarreMil1}
et du corolaire~\ref{corthCarreMil1}.
\end{proof}
%-----------------end proof------------------

Dans le \tho qui suit on suppose que:
%-----------------begin item------------------
\begin{itemize}
\item  $\gk$ est réduit,
%\item  les \mrcs sur $\gk$ sont tous libres,
\item  les \mrcs sur $\gk[t]$ sont tous libres,
\item  les matrices carrées de déterminant $1$ sont produits de matrices
\elrsz, i.e. $\SL_n(\gk)=\En(\gk)$ pour tout $n$.
\end{itemize}
%-----------------end item------------------

Par exemple $\gk$ peut être un \cdiz, un anneau \zed réduit ou un anneau
euclidien intègre.
Notez aussi que si les \mrcs sur $\gk[t]$ sont libres, 
c'est a fortiori vrai des \mrcs sur $\gk$.

%--- Theorem{thIdenPtsSimples2}----
\begin{theorem}
\label{thIdenPtsSimples2}
Pour $a\in \gk$ notons $J_{r,a}\eqdefi\Diag(1,\alb\ldots ,1,a)
\in\MM_r(\gk)$. Sous les hypothèses  précédentes on obtient la classification complète
des  \mrcs sur l'anneau $\gA$ (on utilise les notations et conventions précédentes).
\begin{enumerate}
\item Les modules
de rang constant  $M(H_1,\ldots ,H_s)$ et  $M(G_1,\ldots ,G_s)$ sont isomorphes
\ssi $\det(H_j^{-1}\cdot H_1)=\det(G_j^{-1}\cdot G_1)$ pour tout~$j$.
\item Tout \Amrc $r$ est isomorphe à un unique module
\[M_r(a_2,\ldots ,a_s)\eqdefi M(\I_r,\alb J_{r,a_2},\ldots ,J_{r,a_s}),\] où les $a_i$ sont dans $\gk^\times$. En outre:
$$
M_r(a_2,\ldots ,a_s)   \simeq   \Ae{r-1}\oplus  M_1(a_2,\ldots ,a_s).
$$
\item Enfin la structure de $\GKO\gA$ est précisée par
%--------------------begin array---------------
$$\arraycolsep2pt\begin{array}{rcl}
%M_r(a_2,\ldots ,a_s)&  \simeq & \Ae{r-1}\oplus  M_1(a_2,\ldots ,a_s)    \\[1mm]
M_1(a_2,\ldots ,a_s)\otimes M_{1}(b_2,\ldots ,b_s)&  \simeq &  M_1(a_2b_2,\ldots
,a_sb_s)    \\[1mm]
M_1(a_2,\ldots ,a_s)\oplus M_{1}(b_2,\ldots ,b_s)&  \simeq &  \gA\oplus
M_1(a_2b_2,\ldots ,a_sb_s )
\end{array}$$
%---------------------end array--------------
En particulier, $\Pic(\gA)\simeq (\gk\eti)^{s-1}$.
\end{enumerate}

\end{theorem}
%--- end-theorem-----------------------------------------
%-----------------begin proof------------------
\begin{proof}
%Voyons le premier point. On doit examiner dans quelles conditions les modules
%de rang constant  $M(H_1,\ldots ,H_s)$ et  $M(G_1,\ldots ,G_s)$ sont isomorphes.
\emph{1.} En cas d'\iso toutes les matrices dans les équations (\ref{eqthIdenPtsSimples})
sont \ivsz, et il revient au même de demander

\snic{  H_j^{-1}\cdot H_1\cdot \Phi(\alpha_1) \cdot G_1^{-1}\cdot G_j =
\Phi(\alpha_j) }

%\sni
pour $j\in\lrb{2..s}$. Puisque $\Phi=\Phi(t)$ est \ivz, son \deter
est un \elt \iv de $\gk[t],$ donc de $\gk,$ et tous les
$\det \Phi(\alpha_i)$
sont égaux à~$\det \Phi$.
En conséquence les deux modules ne peuvent être isomorphes que si

\snic{  \det(H_j^{-1}\cdot H_1)=\det(G_j^{-1}\cdot G_1)}

%\sni
pour tout $j$
(ceci prouve en particulier l'unicité de la suite $a_2,\ldots ,a_s$ lorsque
$M_r(a_2,\ldots ,a_s)$ est isomorphe à un \mrc donné).
Inversement si cette condition est satisfaite, on peut trouver une matrice
\elr $\Phi$ qui réalise les conditions ci-dessus. Il suffit en effet d'avoir

\snic{
\Phi(\alpha_1)=\I_r$ et $\Phi(\alpha_j)= H_j^{-1}\cdot H_1 \cdot G_1^{-1}\cdot
G_j,}

%\sni
 ce que l'on obtient en appliquant le lemme qui suit.
\\
La fin est laissée \alecz. Rappel%
%\thrf{thPicKTO}
:  si $Q=P_1\oplus P_2\simeq \gA\oplus P$ (les $P_i$ \prcs 1), on a
$P\simeq \Vi_\Ae 2Q\simeq P_1\te_\gA P_2$.
\end{proof}
%-----------------end proof------------------
%--- Lemma{lemthIdenPtsSimples2}-----
\begin{lemma}
\label{lemthIdenPtsSimples2}
Soient $\alpha_1$, \ldots, $\alpha _s$ dans un anneau commutatif $\gk$ avec les \hbox{différences $\alpha_i-\alpha_j$} \ivs pour $i\neq j$. 
\\
\'Etant données
 $A_1$, \ldots, $A_s\in \EE_r(\gk)$, il existe une
\hbox{matrice $A\in\EE_r(\gk[t])$} telle \hbox{que $A(\alpha_i)=A_i$} pour chaque~$i$.
\end{lemma}
%--- end-lemma-----------------------------------------
%-----------------begin proof------------------
\begin{proof}
Si  une matrice $A\in\EE_r(\gk[t])$  s'évalue en  $s$
matrices $A_1$, \ldots, $A_s$, et une matrice $B\in\EE_r(\gk[t])$  s'évalue en
$s$ matrices $B_1$, \ldots, $B_s$, alors  $AB$ s'évalue \hbox{en  
$A_1B_1$,} \ldots, $A_sB_s$. En conséquence, il suffit de prouver le lemme lorsque
les~$A_i$ sont toutes
égales à $\I_r$ sauf une qui est une matrice \elrz.
Dans ce cas on peut faire une interpolation à la Lagrange puisque les
\hbox{\elts $\alpha_i-\alpha_j$} sont \ivsz.
\end{proof}
%-----------------end proof------------------

%-% ENTRE NOUS
\entrenous{Il faudra prévoir un jour de traiter
avec autant de détails \gui{Un point avec multiplicité}
ainsi que \gui{La droite projective privée d'un nombre fini de points}
}
%-% Fin ENTRENOUS

%: junk
\junk{

%--- Subsec{Un point avec multiplicite} --
\subsec{Un point avec multiplicité}

\hum{à écrire}

%  Section{secAppliDrProj}---------
\section{La droite projective privée d'un nombre fini de points}
\label{secAppliDrProj}
Cette section est basée sur l'article de Zachs \cite{} qui s'occupe de \ldots
et sur le livre de Lorenzini \cite{}.

}
%: fin junk

%:section: Exercices
\penalty-2500
\Exercices

%--- Exercise{exo3Lecteur}-------------
\begin{exercise}
\label{exoChapPtf2Lecteur}
{\rm  Il est recommandé de faire les \dems non données, esquissées,
laissées \alecz,
etc\ldots
\  On pourra notamment traiter les cas suivants.
\begin{itemize}
\item \label{exopropdiverschap6}
Démontrer les propositions \ref{prop rgconstant local}
et~\ref{prop sfio unic}.
\item \label{exodef rang inferieur}\relax
Vérifier les \eqvcs dans la
proposition \ref{def rang inferieur}~\emph{\ref{i3def rang inferieur}}.

\item Vérifier le corolaire \ref{lem2TransPtf}.
\item \label{exolemSomProjOrt} Vérifier les faits \ref{lemSomProjOrt}
et \ref{lemProjProj}
%
%\item \label{}

%
\end{itemize}
}
\end{exercise}
%--- end -exercise-----------------------------------------

%--- Exercise{exoleli2}-------------
\begin{exercise}\label{exoleli2}
{\rm  
  Vérifier les calculs dans le deuxième lemme de la liberté 
  locale~\ref{leli2}.  
}
\end{exercise}
%--- end -exercise-----------------------------------------

%--- Exercise{exo7.1}-----------
\begin{exercise}
\label{exo7.1} (Formule magique pour diagonaliser une \mprnz) \relax \\
{\rm Soit un entier $n$ fixé.
Si $\alpha\in\cP_{n}$ (ensemble des parties finies de
$\lrbn$), on considère le \prr canonique
obtenu à partir de $\In$ en annulant les \elts diagonaux dont
l'indice n'appartient pas à $\alpha$. Nous le notons $\I_{\alpha}$.
Soit $F\in\GAn(\gA)$ un \prrz, on va expliciter une
famille $(F_{\alpha})$ indexée par $\cP_{n}$ de matrices  vérifiant
les \gui{conjugaisons}
$$
FF_{\alpha}=F_{\alpha}\I_{\alpha}\eqno(\dag)
$$
ainsi que l'\ida
$$
\som_{\alpha}\det F_{\alpha}=1\eqno(\ddag)
$$
Ce résultat fournit une nouvelle
méthode uniforme pour expliciter la liberté locale  d'un \mptfz:
on prend les \lons en les \eco $\det(F_\alpha)$, puisque sur l'anneau
$\gA[1/\det(F_\alpha)]$ on a $F_{\alpha}^{-1}FF_{\alpha}=\I_{\alpha}$.\\
Nous allons voir que ceci est réalisé par la
famille définie comme suit:
$$
F_{\alpha}=F\,\I_{\alpha}+(\In-F)(\In-\I_{\alpha}).
$$
 Par exemple si $\alpha=\lrbk$,
on a les \dcns par blocs suivantes
$$
\I_{\alpha}=\I_{k,n}=\bloc{\I_{k}}{0}{0}{0}, \quad F=\bloc{F_{1}}{F_{2}}{F_{3}}{F_{4}}, \quad F_{\alpha}=\bloc{F_{1}}{-F_{2}}{F_{3}}{\I_{n-k}-F_{4}}.
$$
\begin{enumerate}
\item Montrer $(\ddag)$. Indication: pour deux matrices carrées $A$ et $B$
d'ordre $n$, on développe le \deter $\det(A+B)$ comme fonction multi\lin
 des colonnes de $A+B$. On obtient une somme de $2^n$ \deters de matrices
obtenues en mélangeant des colonnes de $A$ et des colonnes de $B$. On applique
cette remarque avec $A=F$ et $B=\In-F$.
\item Si $f$ et $e$ sont deux \idms dans un anneau non \ncrt commutatif,
et si l'on note $f*e=fe+(1-f)(1-e)$, montrer que $f(f*e)=fe=(f*e)e$.
Avec $f=F$ et $e=\I_{\alpha}$, on obtient $f*e=F_\alpha$
ce qui donne l'\egt $(\dag)$ ci-dessus.
\item Nous étudions maintenant quelques \egts qui font intervenir $\det F_\alpha$.
On note $\beta$ le \cop de $\alpha$
\begin{itemize}
  \item Montrer que $(1-2f)(1-e-f)=(1-e-f)(1-2e)=f*e$
  \item Montrer que $(1-2e)^2=(1-2f)^2=1$.
  \item Avec $f=F$ et $e=\I_{\alpha}$, on obtient
  $(\det F_\alpha)^2=\big(\det(\I_{\beta}-F)\big)^2$.
  \item Vérifier que $(1-e)f(1-e)+e(1-f)e=(e-f)^2$.
  \item Si l'on note $\mu_\alpha$ le mineur principal extrait de $F$ sur les
  indices appartenant à $\alpha$, et $\mu'_\beta$ 
  le mineur principal extrait de $\I-F$ sur les
  indices appartenant à $\beta$
  montrer que $(\det F_\alpha)^2=\mu_\alpha\mu'_\beta$.\\
  Indication:
  pour l'exemple ci-dessus avec  $f=F$ et $e=\I_{\beta}$ l'\egt
  dans le point précédent donne
$$
\bloc{F_{1}}{0}{0}{\I_{n-k}-F_{4}}= (\I_{\beta}-F)^2
$$
\end{itemize}
\end{enumerate}
}
\end{exercise}
%--- end-exercise-----------------------------------------
NB. Cette méthode uniforme de \din des \mprns donne un raccourci
pour le lemme de la liberté locale et pour le \tho de structure
qui affirme qu'un \mptf est localement libre au sens fort.
Nous avons pris la peine de démontrer ce \tho de structure deux fois.
Une fois par les \idfs dans le chapitre \ref{chap ptf0},
une autre fois de façon plus structurelle,
dans le chapitre présent. Nous espérons que
\llec ne nous en voudra pas de lui avoir fait subir des preuves nettement
moins \elrs dans le cours que celle de l'exercice \ref{exo7.1}.
C'est que les formules magiques
sont certainement une bonne chose, mais qu'elles cachent parfois la
signification profonde de preuves plus élaborées.

%--- Exercise{exo2Diag}-------------
\begin{exercise}\label{exo2Diag} (Généralisation de l'exercice précédent à la \din de matrices annulant un \pol \spl scindé)
\\
{\rm  Soient $a$, $b$, $c\in\gA$ tels que $(a-b)(a-c)(b-c)\in\Ati$, i.e., le \pol 

\snic{f(T)=(T-a)(T-b)(T-c)}

%\sni
est \splz, et $A\in\Mn(\gA)$ une matrice telle que $f(A)=0$. On considère les \pols de Lagrange $f_a(T)=\fraC{(T-b)(T-c)}{(a-b)(a-c)}$, \ldots\, qui vérifient
$f_a+f_b+f_c=1$. On pose $A_a=f_a(A), \,A_b=f_b(A), \,A_c=f_c(A)$.
\begin{enumerate}\itemsep0pt
\item Montrer que $AA_a=aA_a$, i.e.,
tout vecteur colonne $C$ de $A_a$ vérifie $AC=aC$.
\item En déduire que si une matrice $P$ a pour vecteurs colonnes des vecteurs colonnes de $A_a$ ou $A_b$ ou $A_c$, alors $AP=PD$, où $D$ est une matrice diagonale avec
pour \elts diagonaux, $a$, $b$ ou $c$.
\item En écrivant $1=\det(\In)=\det(A_a+A_b+A_c)$ et en utilisant la multilinéarité du \deter comme fonction des vecteurs colonnes, montrer qu'il existe~$3^n$ matrices
$P_i$ qui vérifient:
\begin{itemize}
\item $\sum_i\det(P_i)=1$
\item  Dans $\gA[1/\det(P_i)]$, la matrice $A$ est semblable à une matrice diagonale
 avec pour \elts diagonaux, $a$, $b$ ou $c$.
\end{itemize}
\item Si le \polcar de $A$ est égal à $(T-a)^m(T-b)^p(T-c)^q$,
montrer que de nombreuses matrices $P_i$ sont nulles et que la somme $\sum_i\det(P_i)=1$
peut être restreinte à une famille de matrices indexée par un ensemble fini
ayant $\fraC{(m+p+q)!}{m!p!q!}$ \eltsz.
\end{enumerate}
}
\end{exercise}
%--- end -exercise-----------------------------------------

%--- Exercise{exoJacobienneP2=P}-------------
\begin {exercise} \label {exoJacobienneP2=P}
       (Jacobienne du système $P^2 - P = 0$)\\
{\rm
Soit l'application $\Mn(\gA) \to \Mn(\gA)$ définie par $P
\mapsto P^2 - P$. Sa \dile en un point
$P\in\GAn(\gA)$ est

\snic{\varphi_P : \Mn(\gA) \to \Mn(\gA),\,
  H \mapsto HP + PH - H.}

%\sni
Si l'on identifie  $\Mn(\gA)$ et
$\Ae {n^2}$, $\varphi_P $ est donné par la matrice jacobienne au
point~$P$ des $n^2$
\eqns $P ^2 - P  = 0$.\\
En considérant

\snic{\gA = \Gn(\ZZ) = \aqo{\ZZ[(X_{ij})_{i,j \in\lrbn}]} {P ^2 - P }$  avec $P  = (X_{ij}),}

%\sni
d'après le  \thref{prop1TanGrassmann}, l'espace
tangent du schéma affine $\GAn$ au point~$P$ s'identifie canoniquement
à

\snic{\Ker \varphi_P = \sotq{H \in \Mn(\gA)}{HP + PH = H} = \Im \pi_P,}

%\sni
 où
$\pi_P \in \GAn(\gA)$ est le \prr défini par

\snic{\pi_P(H) = PH(\In-P) + (\In-P)HP= PH+HP-2PHP.}

%\sni
Ceci amène à étudier les relations entre $\varphi_P$ et
$\pi_P$.
Illustrer ce qui est affirmé concernant la matrice jacobienne et
l'identification de $\Mn(\gA)$ et $\Ae {n^2}$ pour $n = 2$.
En \gnl montrer les \egts

\snic{\varphi_P \circ \pi_P= \pi_P \circ \varphi_P=0$,
 $\;(\varphi_P)^2= \In - \pi_P$,  $\;(\varphi_P)^3 = \varphi_P,}

\snic{\Ker \varphi_P = \Ker (\varphi_P)^2 = \Im \pi_P\;$
 et
$\;\Im \varphi_P = \Im (\varphi_P)^2 = \Ker \pi_P.}
}
\end {exercise}
%--- end-exercise---------------------

%--- Exercise{exopropfideles}--------
\begin{exercise}\label{exopropfideles}
{\rm Démontrer la \carn locale suivante des \mptfs fidèles.
Pour un \Amo $P$, \propeq
%---------begin item----------
\begin{enumerate}
\item [$(a)$] $P$ est \ptf et fidèle.
\item [$(b)$] Il existe des \eco $s_i$ de $\gA$ tels que  chaque $P_{s_i}$ est
libre de rang $h\geq 1$ sur~$\gA_{s_i}=\gA[1/s_i]$.
\item [$(c)$] $P$ est \ptf et pour tout \elt $s$ de $\gA$, si $P_s$ est libre sur l'anneau
$\gA_s$, il est de rang~$h\geq 1$.
\end{enumerate}
%---------end item----------
}
\end{exercise}
%--- end-exercise-----------------------------------------

%--- Exercise{exoRangphi}-------------
\begin{exercise}
\label{exoRangphi}
{\rm Soit $\varphi:P\to Q$ une \Ali entre \mptfs et $r\in\HOp\gA$.
Exprimer $\rg(P)\leq r$ et  $\rg(P)\geq r$ en termes des \idds d'une \mprn
ayant pour image $P$.
 }
\end{exercise}
%--- end -exercise-----------------------------------------

%--- Exercise{exoP1FracRat}-------------
\begin{exercise}\label{exoP1FracRat}
{(Droite projective et fractions rationnelles)}\\
{\rm  
\emph {1.}
Soit $\gk$ un anneau, $P$, $Q \in \gk[u,v]$ deux \pogs de degrés
$p$ et~$q$. On définit:

\snic {
g(t) = P(t,1), \quad \wi g(t) = P(1,t), \quad
h(t) = Q(t,1), \quad \wi h(t) = Q(1,t)
}

\begin {enumerate}
\item [\emph {a.}]
Montrer que $\Res(g,p,h,q) = (-1)^{pq}\Res(\wi g,p, \wi h,q)$,
 que l'on notera $\Res(P,Q)$.
\item [\emph {b.}]
Montrer l'inclusion:

\snic {
\Res(P,Q) \gen {u,v}^{p+q-1} \subseteq \gen {P,Q}
}
\end {enumerate}

\emph {2.}
On rappelle que $\GA_{2,1}(\gk)$ est 
la partie de $\GA_{2}(\gk)$ formée par les
\prrs de rang~1; on a une projection $F \mapsto \Im F$ de $\GA_{2,1}(\gk)$ sur
$\PP^1(\gk)$.
\\
  Lorsque $\gk$ est un \cdi et $f \in
\gk(t)$ une fraction rationnelle,  on associe à~$f$ le \gui{morphisme}, noté
encore $f$, $\PP^1(\gk) \vers{f} \PP^1(\gk)$, qui réalise $t \mapsto f(t)$
(pour l'inclusion usuelle $\gk\subseteq\PP^1(\gk)$).\\
Comment \gnr à un anneau $\gk$ quelconque?
\\
Expliquer comment on peut relever ce morphisme $f$ en une application \pollez,
schématisée ci-dessous en pointillés:

\snic {
\xymatrix {
\GA_{2,1}(\gk)\ar[d]\ar[dr]\ar@{-->}[r]  &\GA_{2,1}(\gk)\ar[d]\\
\PP^1(\gk)\ar[r]^{f}                     &\PP^1(\gk) \\
}}

%\sni
\emph {3.}
Traiter les exemples $f(t) = t^2$, $f(t) = t^d$ et  $f(t) = (t^2 + 1)/t^2$.
Comment se relève une homographie $f(t) = {at + b \over ct +d}$ \,($ad-bc\in\gk\eti$)?
}

\end {exercise}
%--- end -exercise-----------------------------------------

%--- Exercise{exoConiqueFondamentale}-------------
\begin{exercise}\label{exoConiqueFondamentale}
{(La conique fondamentale ou le plongement de Veronese $\PP^1 \to \PP^2$)}\\
{\rm  
%On rappelle que $\GA_{n+1,1}(\gk)$ est la partie de $\GA_{n+1}(\gk)$ 
%formée par les \prrs de rang~1; on a une projection $F \mapsto
%\Im F$ de $\GA_{n+1,1}(\gk)$ sur $\PP^n(\gk)$.
%\\
Lorsque $\gk$ est un \cdiz, le plongement de Veronese $\PP^1(\gk) \to \PP^2(\gk)$
est défini par:

\snic {
(u : v) \mapsto (X=u^2 : Y=uv : Z=v^2)
.}

%\sni
Son image est la \gui {conique fondamentale} de $\PP^2$ d'\eqn 

\snic{\left| \matrix
{X & Y\cr Y &Z\cr}\right| = XZ - Y^2 = 0.}

%\sni
De manière analogue à
l'exercice \ref {exoP1FracRat} (voir aussi le \pb \ref{exoVeroneseMorphism}), montrer que l'on peut relever le morphisme de
Veronese en une application \pollez, schématisée ci-dessous en pointillés:

\snic {
\xymatrix {
\GA_{2,1}(\gk)\ar[d]_{F \mapsto \Im F}\ar[dr]\ar@{-->}[r]  &\GA_{3,1}(\gk)\ar[d]^{F \mapsto \Im F}\\
\PP^1(\gk)\ar[r]^{\rm Veronese}             &\PP^2(\gk) \\
}}

Votre relèvement obtenu doit s'appliquer à un anneau $\gk$ quelconque.
}

\end {exercise}
%--- end -exercise-----------------------------------------

%--- Exercise{exoProjecteurCorangUn}-------------
\begin{exercise}
\label{exoProjecteurCorangUn} (Matrices de \prn de corang $1$) \, {\rm Soit $n\geq2$.
\begin{enumerate}
\item  [\emph{1.}] Soit $P \in \GA_{n,n-1}(\gA)$. Montrer que
 $P + \wi{P} = \In$.

\item [\emph{2.}] Si $P \in \GAn(\gA)$ vérifie $P +
\wi {P} = \In$, alors $P$ est de rang  $n-1$.

\item [\emph{3.}] Si $P \in \Mn(\gA)$ vérifie $\det(P) = 0$ et $P + \wi {P} = \In$, alors $P\in \GA_{n,n-1}(\gA)$.

\end{enumerate}

}
\end{exercise}
%--- end -exercise-----------------------------------------

%--- Exercise{exoMatriceCorangUn}-------------
\begin {exercise}\label {exoMatriceCorangUn}
{\rm
Dans cet exercice, $A \in \Mn(\gA)$ est une matrice de {\it corang 1}, i.e.
de rang $n-1$. En utilisant l'exercice \ref {exoProjecteurCorangUn},
montrer les points suivants.
\begin {enumerate}

\item [\emph{1.}]
$\Im A = \Ker \wi {A}$ (module projectif de rang $n-1$).

\item [\emph{2.}]
$\Im \wi {A} = \Ker A$ (module projectif de rang $1$).

\item [\emph{3.}]
$\Im \tra A = \Ker \tra {\wi {A}}$ (module projectif de rang $n-1$).

\item  [\emph{4.}]
$\Im \tra {\wi {A}} = \Ker \tra A$ (module projectif de rang $1$).

\item  [\emph{5.}]
Les modules projectifs de rang 1, $\Ae n  / \Im A$ et $\Ae n  / \Im
\tra A$, sont duaux l'un de l'autre. En résumé, à partir
d'une matrice $A$ de corang $1$, on construit deux
modules projectifs de rang $1$ duaux l'un de l'autre:

\snic{\begin {array}{c}
\Ae n /\Im A = \Ae n /\Ker \wi {A} \simeq \Im \wi {A} = \Ker A,
\\[1mm]
\Ae n /\Im \tra {A} = \Ae n  / \Ker \tra{\wi{A}} \simeq
\Im \tra{\wi{A}} = \Ker \tra{A}.
\end {array}
}
\end {enumerate}
}
\end {exercise}
%--- end -exercise-----------------------------------------

%--- Exercise{exoIntersectionSchemasAffines}-------------
\begin{exercise}\label{exoIntersectionSchemasAffines}
{(Intersection de deux schémas affines sur $\gk$)}\\
{\rm  
Cet exercice se situe dans le cadre informel des schémas affines
sur un anneau $\gk$ \gui{définis} \paref{subsecSchAff}.
Soient $\gA=\kxn$, $\gB=\gk[\yn]$ deux \klgs quotients correspondant
à deux \syps $(\uf)$, $(\ug)$ dans $\kXn$. 
\\ Notons $A$,
$B$ les schémas affines correspondants.
Le schéma intersection $A\cap B$ est défini comme
associé à la \klg $\aqo\kuX{\uf,\ug} \simeq \gA\otimes_{\kuX}\gB$
(notez que le produit tensoriel est pris sur $\kuX$).

\Deuxcol{.55}{.4}
{\small 
\gui{Justifier} cette définition en vous appuyant sur le dessin ci-contre.
\\
Dans un repère \gui{euclidien},
le dessin comprend l'ellipse $\big(\fraC x a\big)^2 + y^2 = 1$, i.e.  $f(x,y) = 0$ avec $f
= x^2 + a^2y^2 - a^2$, et le cercle $g(x,y) = 0$ avec $g = (x-c)^2 + y^2 -
(c-a)^2$.  }
{~

\centerline{\includegraphics[width=4cm]{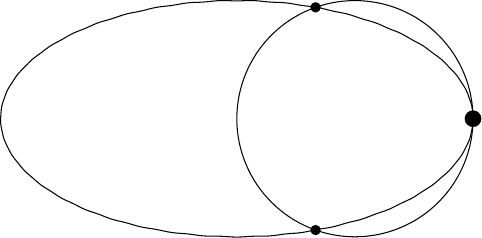}}
}
}
\end {exercise}
%--- end -exercise-----------------------------------------

%--- Exercise{exoPolPseudoUnitaire}-------------
\begin{exercise}\label{exoPolPseudoUnitaire}
{(\Pols pseudo \monsz)}\\
{\rm  
On rappelle
qu'un \pol $p(t) = \sum_{k\ge 0} a_kT^k\in \kT$ est dit pseudo \mon s'il
existe un \sfio $(e_0,\ldots ,e_r)$ tel que sur chaque~$\gk[1/e_j]$,~$p$
est un \pol de degré $j$ avec son \coe de degré $j$ \iv
 (voir \paref{polpseudunit}). Un tel \pol est
primitif et cette notion est stable par produit et morphisme.

\emph {1.}
%:HHH geN ci dessous
Vérifier que $a_k = 0$ pour $k > r$ et que $\geN {(1 - \som_{j>k} e_j)a_k}  =
\gen {e_k}$ pour $k \in \lrb{0..r}$.  En particulier, $\gen {a_r} = \gen
{e_r}$ et les $e_k$ sont uniques ou plutôt le \pol $\sum_k e_k X^k$ est
unique (on peut ajouter des \idms nuls).

\emph {2.}
Soit $P = \aqo{\gA[T]}{p}$. Montrer que $P$ est un $\gA$-\mptf dont
le \pol rang est $\rR{P}(X) = \sum_{k=0}^r e_kX^k$; on a \egmt
$\deg p = \sum_{k=1}^r k[e_k]$ (cf. le point \emph {2}
de \ref {notaRangs}).
Dans le même ordre d'idées, voir l'exercice~\ref {exoPolLocUnitaire}.

}

\end {exercise}
%--- end -exercise-----------------------------------------

%--- Exercise{exoPolLocUnitaire}-------------
\begin{exercise}\label{exoPolLocUnitaire}
 {(\Pols \lot \monsz)}\\
{\rm
\emph {1.}
Soit $\fa \subseteq \gA[T]$ un \id tel que $\gA[t] = \gA[T]\sur\fa$
soit un \Amo libre de rang~$n$. Soit $f \in \gA[T]$  le \polcar
de $t$ dans $\gA[t]$. Montrer que $\fa = \gen {f}$. En particulier,
$1, t, \ldots, t^{n-1}$
est une $\gA$-base de $\gA[t]$.

\emph {2.}
Résultat analogue en remplaçant l'hypothèse \gui{$\gA[T]\sur\fa$ est un \Amo libre
de rang~$n$} par  \gui{$\gA[T]\sur\fa$ est un \mrc $n$}.

Un \pol $f \in \gA[T]$ de degré $\leq r$ est dit \emph{\lot\monz} s'il existe un \sfio
$(e_0,\ldots,e_r)$ tel que $f$ soit unitaire de degré $d$
dans $\gA[1/e_d][T]$ pour chaque $d\in\lrb{0..r}$.%
\index{polynome@\pol!localement unitaire}%
\index{localement!polynome@\pol --- unitaire}
Ainsi,  pour chaque $d\in\lrb{0..r}$, le \pol   $f_d:= e_df$
est \mon de degré~$d$ modulo $\gen{1-e_d}$.
Il est clair que cette \dfn ne dépend pas du degré
formel $r$ choisi pour $f$, et que sur un anneau connexe, un \pol \lot \mon est \monz.

\emph {3.}
Caractériser un \pol \lot\mon à l'aide de ses \coesz.

\emph {4.}
Le \polcar d'un \endo d'un \mptf $M$ est  \lot\mon et le \sfio
correspondant est donné par les $\ide_i(M)$.

\emph {5.}
Soient $S_1$, \ldots, $S_m$ des \moco de $\gA$. Montrer que si $f$ est \lot\mon
(par exemple \monz) sur chaque $S_i^{-1}\gA$, il l'est sur $\gA$.

\emph {6.}
Si $f \in \gA[T]$ est \lot\monz, montrer que l'anneau $\gA[t] = \aqo{\gA[T]}{f}$
est un \Amo quasi libre et que $f$ est le \polcar de $t$.

%\sni
\emph {7.}
Réciproquement, si $\fa \subseteq \gA[T]$ est un \id tel que
l'anneau $\gA[t] =
\gA[T]\sur\fa$ soit un \Amo \ptfz, alors $\fa = \gen {f}$.
En particulier, si une \Alg monogène est un \Amo \ptfz, c'est
un \Amo quasi libre.

%\sni
\emph {8.}
Pour $g\in\AT$  \propeq
\begin{itemize}
\item $g$ peut s'écrire $uf$ avec $u\in\Ati$ et $f$ \lot \monz.
\item $g$ est pseudo \monz.
\item $\aqo{\AT}{g}$ est un \Amo \ptfz.
\end{itemize}

%\sni
\emph {9}\eto$\!\!$.
Démontrer en \clama qu'un \pol est \lot \mon \ssi il devient \mon
après \lon en n'importe quel \idepz.
}
\end {exercise}
%--- end -exercise-----------------------------------------

%--- Exercise{exoOneRankPtfIdeal}-------------
\begin{exercise}\label{exoOneRankPtfIdeal}
{(Modules \ivs et \mrcs 1)}\\
{\rm  
On propose une petite variation autour du \tho \ref{propRgConstant3}.
\\
\emph {1.}
Soient deux anneaux commutatifs $\gA \subseteq \gB$. 
Les sous-\Amos de $\gB$ forment un \mo multiplicatif, d'\elt neutre $\gA$.
Montrer qu'un sous-\Amo $M$ de~$\gB$   \emph{inversible} dans ce \mo  est \tf et que pour tout sous-\Amoz~$M'$ de~$\gB$ l'\homo canonique $M\otimes_\gA M' \to M.M'$ est un \isoz.  En
conséquence, les sous-\Amos de $\gB$ inversibles sont des \Amrcs 1.
\\
\emph {2.}
Soit $S\subseteq\Reg(\gA)$ un \mo  et $\fa$ un \id \lopz.  On suppose que $S^{-1}\fa$  est un \id
\iv de $S^{-1}\gA$; montrer que $\fa$ est un \id \iv de $\gA$. C'est le cas
par exemple si $S^{-1}\fa$ est un $S^{-1}\gA$-module libre.

}

\end {exercise}
%--- end -exercise-----------------------------------------

%--- Exercise{exoPicAPicFracA}-------------
\begin{exercise}\label{exoPicAPicFracA}
{(La suite exacte avec $\Pic\gA$ et $\Pic\gK$, où $\gK = \Frac\gA$)}
\\
{\rm  
Soit $\gA$ un anneau et $\gK = \Frac\gA$. Définir des morphismes 
naturels de groupes:

\snic {
1 \rightarrow \Ati \to \gK\eti \to \Gfr(\gA) \to \Pic\gA \to \Pic\gK
,}

%\sni
et montrer que la suite obtenue est exacte. En conséquence, on a une suite
exacte 

\snic{1 \to \Cl(\gA) \to \Pic\gA \to \Pic\gK.}

%\sni
Si $\Pic\gK$ est trivial, on
obtient un \iso $\Cl(\gA) \simeq \Pic\gA$, et l'on retrouve ainsi le \thref{propRgConstant3}.

}

\end {exercise}
%--- end -exercise-----------------------------------------

%--- Exercise{exoAnneaudesrangs}-------------
\begin{exercise}
 \label{exoAnneaudesrangs}
 {\rm
Montrer que $\HO\gA$ est l'anneau \gui{engendré par}
$\BB(\gA)$, l'\agB des \idms de $\gA$, au sens des foncteurs adjoints.\\
Plus \prmtz, si $B$ est une \agBz, l'anneau $\wi B$ librement engendré par $B$
est donné avec un \homo d'\agBs $\eta_B:B\to\BB(\wi B)$ tel que pour tout anneau $\gC$
l'application décrite ci-après soit une bijection:

\Deuxcol{.6}{.38}
{
\[ 
\begin{array}{ccc} 
 \Hom_{\mathrm{Anneaux}}(\wi B,\gC)\lora \Hom_{\mathrm{Alg.\, de\, Boole}}\big(B, \BB(\gC)\big)  \\[1mm] 
 \varphi\longmapsto \BB(\varphi)\circ \eta_B     
 \end{array}
\]
}
{\xymatrix @R = 0.2cm {
\wi {B}\ar[dd]_\varphi  &B\ar[dd]%^{\BB(\varphi)\circ \eta_B}
  \ar[dr]^{\eta_B} \\
\ar@{~>}[r]
&  & \BB(\wi B)\ar[dl]^{\BB(\varphi)}  \\
\gC  & \BB(\gC) \\
}
}

%\sni
Montrer alors que $\wi{\BB(\gA)}\simeq\HO\gA$.
}
\end{exercise}
%--- end-exercise-----------------------------------------

%--- Exercise{exoHOAclama}-------------
\begin{exercise}
 \label{exoHOAclama}
 {\rm
Démontrer en \clama que  $\HO(\gA)$ est canoniquement
isomorphe à l'anneau des fonctions localement constantes
(i.e., continues) de $\SpecA$ vers $\ZZ$.
}
\end{exercise}
%--- end-exercise-----------------------------------------

%--- Exercise{exoFonctDet1}----------
\begin{exercise}\label{exoFonctDet1}  (Le déterminant comme foncteur)\\
{\rm
On a défini le \deter d'un \endo d'un \mptfz. Nous allons voir que plus
\gnlt on peut définir le \deter comme un foncteur de la catégorie
des modules \pros vers celle des modules \pros de rang 1.
Sans doute, la \dfn la plus simple du \deter d'un \mptf est la suivante:

\sni{\it Définition:
%---------begin item----------
\vspace{-3pt}
\begin{enumerate}\itemsep0pt
\item [$(a)$] Soit $M$  un $\gA$-\mptf engendré par $n$ \eltsz.
\\ On note  
$r_h=\ide_h(M)$ ($h\in\lrbzn$) et $M\ep{h}=r_hM$. On définit $\det(M)$ par

\snic{ \det(M) := r_0\gA\oplus M\ep{1}\oplus \Al2M\ep{2}\oplus \cdots\oplus
 \Al nM\ep{n}
.}

%\sni
Nous utiliserons
aussi la notation suggestive $\det(M) =\Al{\rg(M)}M$ en utilisant le rang
$\rg(M)=\sum_{k=1}^{n}k\,[\ide_k(M)]\in\HO\gA$.
\item [$(b)$] Si $\varphi ~:M\rightarrow N$ est un \homo de
$\gA$-\mptfsz, avec $s_h=\ide_h(N)$, on définit $\det(\varphi)$ comme un \homo
de $\det(M)$ dans $\det(N)$ envoyant  $\Al hM\ep{h}$ dans
$\Al hN\ep{h}$ par $x\mapsto s_h (\Al h \!\varphi)(x)$.
\end{enumerate}
%---------end item----------
}

%\sni
On notera que lorsque $x\in\Al hM\ep{h}$ on a $x=r_hx$. 
%---------begin item----------
\vspace{-3pt}
\begin{enumerate}\itemsep0pt
\item  Le module $\det(M)$ est un \mrc 1, et l'on a les \egts $r_h\det(M)=\det(M)_{r_h}=\Al hM\ep{h}$.
Plus \gnltz, pour tout \idm $e$,
on a $e\det(M)=\det(M_e)$.
\item  La \dfn précédente fournit un foncteur qui commute avec la
\lon et transforme les sommes directes en produits tensoriels.
%:HHH une phrase en plus
En déduire que le foncteur $\det$ induit un morphisme surjectif de $(\KO\gA,+)$ sur $\Pic\gA$.
\item  Un \homo entre \mptfs est un \iso \ssi son \deter est un \isoz.
\item  Pour un \endo d'un \mptfz, la nouvelle \dfn du \deter coïncide
avec l'ancienne si l'on identifie $\End(L)$ avec~$\gA$ lorsque $L$ est un \mrcz~1.
\end{enumerate}
%---------end item----------
}
\end{exercise}
%--- end-exercise-----------------------------------------

%--- Exercise{exoFonctDet3}----------
\begin{exercise}
\label{exoFonctDet3}
{\rm
\`A \iso près, le foncteur \deter est le seul foncteur de
la catégorie des \Amos \ptfs dans elle-même qui possède
les \prts suivantes:
\begin{itemize}
\item il  transforme toute flèche $\varphi:\gA\rightarrow \gA$ en elle-même,
\item  il  transforme les
sommes directes en produits tensoriels,
\item  il commute à l'\eds pour tout changement de base \smash{$\gA\vers\alpha\gB$}.
\end{itemize}
\perso{Il semblerait que l'exo \ref{exoFonctDet3}
ait quelques rapports
avec le \deter d'un complexe et les preuves télégraphiques des
\prts de ce \deter données dans le \emph{Algèbre locale,
multiplicités} de J.-P. Serre.
}
}
\end{exercise}
%--- end-exercise-----------------------------------------

%--- Exercise{exoIDDPTF1}----------
\begin{exercise}
\label{exoIDDPTF1} (Idéaux déterminantiels d'une \ali entre \mptfsz)\, 
{\rm
Soit $\varphi : M\rightarrow N$ un \homo entre \mptfsz.
\'Ecrivons $M\oplus M'\simeq \Ae m $, $N\oplus N'\simeq \Ae n $,  et
\hbox{prolongeons $\varphi$ en}  

\snic{\psi:M\oplus M'\to N\oplus N'\;$ 
avec $\;\psi(x+x')=\varphi(x)\;  ( x\in M , \, x'\in M' ).}

%\sni
Alors, pour
chaque entier $h$, l'\idd $\cD_h(\psi)$ ne 
dépend que de~$h$ et de $\varphi$. On le
note $\cD_h(\varphi)$ et on l'appelle \emph{l'\idd d'ordre $h$ de~$\varphi$}.
}%
\index{ideaux determ@idéaux déterminantiels!d'une \ali (\mptfsz)}%
\index{ideal@idéal!déterminantiel}
\end{exercise}
%--- end-exercise-----------------------------------------

%:     Notation{notaRgfi}
\begin{notation}\label{notaRgfi}
 {\rm  Soit $r=\sum_{k=1}^{n}k\,[r_k]\in \HOp (\gA)$.
En application de l'exercice précé\-dent,
on appelle \emph{\idd de type $r$ pour  $\varphi$} et l'on note
$\cD_r(\varphi)$ l'idéal

\snic{r_0\gA+r_1\cD_1(\varphi)+\cdots + r_n\cD_n(\varphi).
}

%\sni
Les notations $\rg(\varphi)\geq k$ et  $\rg(\varphi)\leq k$
pour les \alis entre modules libres de rang fini se généralisent
comme suit aux \alis entre \mptfsz: on
note $\rg(\varphi)\geq r$  si $\cD_r(\varphi)=\gen{1}$,
$\rg(\varphi)\leq r$  si~$\cD_{1+r}(\varphi)=\gen{0}$, et
$\rg(\varphi)= r$ si $\rg(\varphi)\leq r$ et $\rg(\varphi)\geq r$.
\\
NB: voir l'exercice \ref{exoLocSimpPtf}.%
\index{rang!d'une application linéaire}

}
\end{notation}

%--- Exercise{exoIDDPTF2}----------
\begin{exercise}
\label{exoIDDPTF2} 
{\rm  (Suite de l'exercice \ref{exoIDDPTF1}) \, Soit $r\in\NN\etl$.
%-----------------begin enum------------------
\begin{itemize}
\item  [\emph{1.}]  Si $  M\vvers{\varphi}  N\vvers{\varphi'} L$ sont des \alis entre \mptfsz, on~a:
$\cD_r(\varphi'\varphi)\subseteq \cD_r(\varphi')\cD_r(\varphi).$
\item [\emph{2.}]  Si $S$ est un \mo de $\gA$, alors $\big(\cD_r(\varphi)\big)_S=\cD_{r}(\varphi_S)$.
\item [\emph{3.}]  Pour tout $s\in \gA$ tel que $M_s$ et $N_s$  sont libres, on a
$\big(\cD_r(\varphi)\big)_s=\cD_{r}(\varphi_s)$. 
%où $k_s=\rg(P_s)\in\NN$.
En outre, cette
\prt caractérise l'\id $\cD_r(\varphi)$.

\end{itemize}
%-----------------end enum------------------
Soit $r=\sum_{k\in\lrbn} k [r_k]\in\HOp\gA$. 
\begin{itemize}\setcounter{enumi}{3}
\item [\emph{4.}] Reprendre les points
précédents de l'exercice dans ce nouveau cadre.

\end{itemize}

}
\end{exercise}
%--- end-exercise-----------------------------------------

\entrenous{pour les deux exos précédents, on peut se demander
ce qui se \gns aux \mpfsz}

%--- Exercise{exoLocSimpPtf}-------------
\begin{exercise}
 \label{exoLocSimpPtf}
 {\rm  (Avec les notations \ref{notaRgfi}) \,
 Soit $\varphi:M\to N$ une \ali entre \Amos \ptfsz.
\Propeq
\\
 \emph{1.} $\varphi$ est \lnlz.
\\
 \emph{2.} $\varphi$ a un rang bien défini
dans  $\HOp(\gA)$.
\\
 \emph{3.} Après \lon en des \eco les modules sont libres
et l'\ali est \nlz.
 } \end{exercise}
%--- end-exercise-----------------------------------------

%--- Exercise{exoCoMatCoRang1}-------------
\begin{exercise}
\label{exoCoMatCoRang1}
{\rm
Soit $A \in \Ae {n \times m}$; si $A$ est de rang $m-1$, on
va expliciter un \sys fini de \gtrs du sous-module $\Ker A
\subseteq \Ae n$ sans utiliser ni test d'\egt ni test d'appartenance.
En fait, sous la seule hypothèse (plus faible) $n \ge m-1$,
on définit de manière
uniforme une matrice $A' \in \Ae {m \times N}$ avec $N = {n \choose m-1}$ qui
est \gui {une sorte de comatrice de $\gA$}. Cette matrice satisfait
à $\Im A' \subseteq \Ker A$ dès que~$A$ est de rang $\leq m-1$, avec
\egt lorsque $A$ est de rang $m-1$.

%\sni
On peut définir $A' \in \Ae {m \times N}$ via l'\alg
extérieure: on voit $A$ comme une \ali $u : \Ae m \to \Ae n$
et l'on considère 

\snic{u' = \Al {m-1}(\tra u) :
\Al {m-1}(\Ae n) \to \Al {m-1}(\Ae m).}

\snii
Dans les bases canoniques,
$\Al {m-1}(\Ae n) = \Ae N$ et $\Al {m-1}(\Ae m) = \Ae m$, donc $u'$ est
représenté par une matrice $A' \in \Ae {m \times N}$.  Pour expliciter
cette matrice $A'$, on ordonne l'ensemble des $N = {n \choose m-1}$ parties
$I$ de $\lrb{1..n}$ de cardinal $m-1$ de façon à ce que leurs
\cops soient rangés de manière croissante pour l'ordre
lexicographique; les colonnes de $A'$ sont indexées par cet ensemble de
parties, de la manière suivante:

\snic {
a'_{j,I} = (-1)^{k_I+j} \det(A_{I, \{1..m\} \setminus \{j\}}), \qquad
\hbox {$k_I$ étant le numéro de $I$.}
}

%\sni
Par exemple, si $m = 2$, alors $N = n$, et $A'=\cmatrix {
a_{n,2}  & -a_{n-1,2} & \cdots & \pm a_{1,2} \cr
-a_{n,1} & a_{n-1,1} & \cdots & \mp  a_{1,1}
}.$

\emph{1.}
Pour $i \in \lrb{1..n}$, on a $(AA')_{i,I} = (-1)^{k_I+1} \det(A_{\{i\} \cup
I, \{1..m\}})$.  En particulier, si~$\cD_m(A) = 0$, alors $AA' = 0$.

\emph{2.}
Si $n = m$, alors $A' = \wi{A}$ (la co-matrice de $A$).

\emph{3.}
Si $A$ est de rang $m-1$, alors $\Im A' = \Ker A$; en particulier, $A'$ est de
rang $1$.

\emph{4.}
Tout module \stl de rang 1 est libre. On pourra comparer
avec le fait \ref{factPicStab} et avec l'exercice \ref{exoStabLibRang1}.

\emph{5.} Si $B$ est une matrice vérifiant $ABA=A$, alors $P=BA$ est une
\mprn vérifiant $\Im(\In-P)=\Ker P=\Ker A$. Ceci fournit une autre manière
de répondre à la question: donner un \sys fini de \gtrs de $\Ker A$.
Comparer cette autre solution à celle de l'exercice présent.
Pour le calcul de la matrice $P$, on pourra utiliser la méthode
expliquée dans la section \ref{secCramer} (\thref{propIGCram}). 
Une autre méthode, nettement plus économique, se trouve dans \cite[Díaz-Toca\&al.]{DiGLQ} (basé sur \cite[Mulmuley]{Mul}).

}

\end {exercise}
%--- end -exercise-----------------------------------------

%--- Exercise{exoVarProj}-------------
\begin{exercise}
\label{exoVarProj} (\Pols\hmgs et $\Pn(\gk)$)
\\ 
{\rm Soit un \syp \hmg $(\lfs)=(\uf)$ dans $\gk[\Xzn]%=\kuX
$. On
cherche à définir les zéros de $(\uf)$ dans $\Pn(\gk)$. 
Soit $P$ un point de $\Pn(\gk)$, i.e. 
\linebreak
un \kmo \pro de rang $1$ en facteur direct dans
$\gk^{n+1}$. Montrer que si un \sgr de $P$ annule $(\uf)$, alors tout \elt de $P$
annule~$(\uf)$.
 
}
\end{exercise}
%--- end -exercise-----------------------------------------

%--- Exercise{exoGLnTangent}-------------
\begin{exercise}\label{exoGLnTangent}
 {(Espace tangent à  $\GLn$)}\\
{\rm Déterminer l'espace tangent en un point au foncteur $\gk \mapsto \GLn(\gk)$.
}

\end {exercise}
%--- end -exercise-----------------------------------------

%--- Exercise{exoSLnTangent}-------------
\begin{exercise}\label{exoSLnTangent} {(Espace tangent à  $\SLn$)}
\\
{\rm  
Déterminer l'espace tangent en un point au foncteur $\gk \mapsto \SLn(\gk)$.
}

\end {exercise}
%--- end -exercise-----------------------------------------

%--- Exercise{exoTangentJ0ConeNilpotent}-------------
\begin {exercise} \label {exoTangentJ0ConeNilpotent}
       (Espace tangent en $J_0$ au cône nilpotent)  \,  {\rm Soit $\gk$ un anneau.\\
On note $(e_{ij})_{i,j\in\lrbn}$ la base canonique de $\Mn(\gk)$ et $J_0
\in \MMn(\gk)$ la matrice de Jordan standard. Par exemple, pour $n = 3$, 
$J_0 =\cmatrix {0 &1 &0\cr 0 & 0& 1\cr 0& 0& 0\cr}$.

\emph{1.}
On définit $\varphi : \Mn(\gk) \to \Mn(\gk)$ par
$\varphi(H) = \sum_{i+j = n-1} J_0^i H J_0^j$.\\ 
Déterminer $\Im \varphi$.

\emph{2.}
Déterminer un \supl de $\Im \varphi$ dans $\Mn(\gk)$,
puis $\psi : \Mn(\gk) \to \Mn(\gk)$ vérifiant $\varphi \circ
\psi \circ \varphi = \varphi$. 
Montrer que  $\Ker\varphi$ est libre de rang $n^2-n$ et donner
une base de ce module. 

\emph{3.}
On considère le foncteur $\gk \mapsto \{N \in \Mn(\gk) \ | \
N^n = 0\}$. Déterminer l'espace tangent en   $J_0$ à ce foncteur.

}
\end {exercise}
%--- end -exercise-----------------------------------------

%--- Exercise{exoprop1TanGrassmann}----------
\begin{exercise}\label{exoprop1TanGrassmann} (Complément pour le \thref{prop1TanGrassmann})
{\rm On note $\gA[\vep]=\AT/\geN{T^2}$. \\
Soient $P$, $H \in \MM_n(\gA)$.
Montrer que la matrice $P + \varepsilon H$ est \idme \ssi

\snic{P^2 = P  \  $ et $ \  H = HP + PH.}

%\sni
Généraliser à un anneau non commutatif arbitraire
avec un \idm $\varepsilon$ dans le centre de l'anneau.

\comm L'exemple de l'anneau $\Mn(\gA)$ montre que dans le cas non commutatif la situation
pour les \idms est assez différente de celle dans le cas commutatif où $\BB(\gA)=\BB(\Ared)$
(corolaire~\ref{corIdmNewton})
et où les \idms sont \gui{isolés} (lemme \ref{lemIdmIsoles}).}\eoe
\end{exercise}
%--- end-exercise-----------------------------------------

%: sinotenglish
\sinotenglish{

%:--- Exercise{exoMonomialSyzygies}-------------
\begin{exercise}\label{exoMonomialSyzygies} {(Syzygies entre \momsz)}
\\
{\rm  
Soit $\kuX = \gk[\Xn]$ un anneau de \pols (où $\gk$ est un anneau quelconque)
et $s$ \moms $m_1$, \ldots, $m_s$ de $\kuX$. On note $m_i \vi m_j$ le pgcd de
$(m_i, m_j)$ et 
%
%\snic {
%\displaystyle 
${m_{ij} = {m_j \over m_i \vi m_j}}$
 {de sorte que} 
\framebox [1.1\width][c]{$m_{ij}\cdot m_i = m_{ji}\cdot m_j$}.
%}

\snii
On note $(\vep_1, \ldots, \vep_s)$ la base canonique de $\kuX^s$. 
\\
L'encadré
fournit des \syzys $r_{ij}=m_{ij}\vep_i - m_{ji}\vep_j$ pour $(m_1, \ldots, m_s)$.
Montrer que ces \syzys engendrent le module des \syzys pour
$(m_1, \ldots, m_s)$, i.e. le noyau de l'application $\kuX$-\lin $\vep_i \mapsto m_i$
de $\kuX^s$ vers $\kuX$.
}
\end{exercise}

}
%: fin sinotenglish

%%%%%%%%%%%%%%%%%%%%%%%%%%%%%%%%%%%%%%%%%%%%%%%%%%%%%%%%%%%%%%%%%%%%%%%%%%%
%: problemes

%--- Problem{exoAnneauCercle1}-------------
\begin{problem}\label{exoAnneauCercle1} {(L'anneau du cercle)}
\\
{\rm
Soit $\gk$ un \cdi de \cara $\ne 2$, $f(X,Y) = X^2 + Y^2 - 1 \in
\gk[X,Y]$. C'est un \pol \ird et lisse, i.e. $1 \in \gen {f, \Dpp{f}{X}, \Dpp{ f}{ Y}}$ (explicitement, on~a~$-2 = 2f -
X\Dpp{f}{X} - Y\Dpp{f}{Y}$). \\ 
Il est
donc licite de penser que l'anneau $\gA = \gk[X,Y]/\gen {f} = \gk[x,y]$ est un \adp
intègre. Ceci sera prouvé dans le \pbz~\ref{exoArithInvariantRing} (point~\emph{4}). 
\\ 
On note $\gK$ son corps
des fractions et l'on pose $t = {y \over x-1} \in \gK$.
\begin {enumerate}\itemsep0pt
\item
Montrer que $\gK = \gk(t)$; justifier \gmqt comment trouver $t$
(paramétrage d'une conique ayant un point $\gk$-rationnel)
et expliciter $x$, $y$ en fonction de~$t$.
\item
Soit $u = (1+t^2)^{-1}$, $v = tu$.  Vérifier que la \cli de
$\gk[u]$ dans $\gK = \gk(t)$~est

\snic {
\gk[x,y] = \gk[u, v] =
\sotq {h(t)/(1 + t^2)^s} {h \in \gk[t],\  \deg(h) \le 2s}
.}

%\sni
En particulier, $\gA = \gk[x,y]$ est \iclz.  Expliquer en quoi le $\gk$-cercle
$x^2 + y^2 = 1$ est la droite projective $\PP^1(\gk)$ privée
\gui {du $\gk$-point} $(x,y) = (1, \pm i)$.

\item
Si $-1$ est un carré dans $\gk$, montrer que $\gk[x,y]$
est un localisé $\gk[w,w^{-1}]$ (pour un $w$ à expliciter)
d'un anneau de \pols sur $\gk$, donc un anneau de Bézout.

\item
Soit $P_0 = (x_0,y_0)$ un $\gk$-point du cercle $x^2 + y^2 = 1$ et $\gen
{x-x_0, y-y_0} \subseteq \gA$ son \idemaz.  Vérifier que $\gen {x-x_0,
y-y_0}^2$ est un \idp de \gtr $xx_0 + yy_0 - 1$.
Interprétation \gmq de $xx_0 + yy_0 - 1$?

\item
%\sni
\mbox{\parbox[t]{.6\linewidth}{Ici $(x_0, y_0) = (1,0)$. Décrire les calculs permettant d'expliciter la
matrice (de projection)\hsd
 $~P = {1\over 2} \cmatrix {1-x & -y\cr -y & 1+x}$ 
\\
comme \mlp pour le couple $(x-1, y)$. La suite exacte:

\snic {
\Ae 2 \vvvers{\I_2-P} {\Ae 2}\vvvvers{(x-1, y)} \gen {x-1, y}\to 0
}

permet de réaliser l'idéal (\ivz) du point~$(1,0)$ comme l'image
du \prr $P$ de rang~$1$.
\\
Commenter le dessin ci-contre qui en est la
contrepartie \gmq (fibré vectoriel en droites sur le cercle).}%
\hspace*{.05\linewidth}%
\parbox[t]{.35\linewidth}{\begin{center}
\includegraphics*[width=3.5cm]{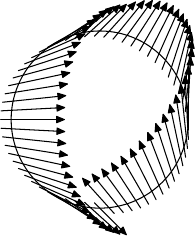} $\qquad\qquad$
\end{center}}}
\item
On suppose que $-1$ n'est pas un carré dans $\gk$
et l'on voit $\gk[x,y]$ comme \linebreak 
une $\gk[x]$-\alg
libre de rang $2$, de base $(1,y)$.
Expliciter la norme et vérifier, pour
$z = a(x) + b(x)y \ne 0$, l'\egt

\snic {
\deg \rN_{\gk[x,y]/\gk[x]} (z) = 2\max(\deg a, 1 + \deg b).
}

 En particulier, $\deg \rN_{\gk[x,y]/\gk[x]} (z)$ est pair.  En déduire
le groupe %des \ivs 
$\gk[x,y]^{\!\times}$, le fait que $y$ et $1 \pm
x$ sont \irds dans $\gk[x,y]$, et que l'idéal~$\gen {x-1, y}$ du point~$(1,0)$ n'est pas principal (i.e. le fibré en droites ci-dessus n'est pas
trivial).
\end {enumerate}
}
\end {problem}
%--- end -problem-----------------------------------------

%--- problem{exoLambdaGammaK0}-------------
\begin{problem}\label{exoLambdaGammaK0}
{(Les opérations $\lambda_t$ et $\gamma_t$ sur $\KO(\gA)$)} \\
{\rm  
Si $P$ est un \mptf sur $\gA$, on note, pour $n\in\NN$, $\lambda^n (P)$
ou~$\lambda^n ([P])$ la
classe de~$\Al{n}P$ dans~$\KO(\gA)$ et l'on a l'\egt fondamentale 

\snic{\qquad\lambda^n(P\oplus Q) =
\sum_{p+q = n}\lambda^p(P)\lambda^q(Q).\qquad(*)}

%\sni
 On définit \egmt
le \pol $\lambda_t(P) \in \KO(A)[t]$ par 
$\lambda_t(P) = \sum_{n\ge 0} \lambda^n(P) t^n$. C'est un \pol de terme
constant $1$ que l'on regarde dans l'anneau des séries formelles
$\KO(A)[[t]]$. Alors

\snic{\lambda_t(P) \,\in\, 1 + t\KO(\gA)[[t]]\;\subseteq\;(\KO(\gA)[[t]])\eti.}

%\sni
D'après $(*)$
on a $\lambda_t(P \oplus Q) = \lambda_t(P)\lambda_t(Q)$, ce qui permet
de prolonger $\lambda_t$ en un morphisme $(\KO(\gA), +) \to
(1 + t\KO(\gA)[[t]], \times)$. Ainsi si $P$, $Q$ sont deux \mptfsz,
pour $x = [P] - [Q]$, on a par \dfn:
$$
\lambda_t(x) = {\lambda_t(P) \over  \lambda_t(Q)} =
{1 + \lambda^1(P)t + \lambda^2(P)t^2 + \cdots \over
1 + \lambda^1(Q)t + \lambda^2(Q)t^2 + \cdots}
$$
série que l'on notera $\sum_{n\ge 0} \lambda^n(x)t^n$, avec $\lambda^0(x) =
1$, $\lambda^1(x) = x$. 

%\sni
Grothendieck a \egmt défini sur $\KO(\gA)$ une autre opération $\gamma_t$
par l'\egt

\snic{\gamma_t(x) = \lambda_{t/(1-t)}(x),}

%\sni
pour $x\in\KO(\gA)$. Ceci est licite
car le sous-groupe multiplicatif $1 + t\KO(\gA)[[t]]$ est stable par la
substitution $t \leftarrow t/(1-t)$. 
Cette substitution $t \leftarrow t/(1-t)$ laisse 
invariant le terme en $t$, on note 

\snic{\gamma_t(x)\,=\,1+tx+t^2\big(x+\lambda^2(x)\big)+\cdots\,=\, \sum_{n\ge 0}
\gamma^n(x)t^n.}

%\sni
\emph {1.} Donner $\lambda_t(p)$ et $\gamma_t(p)$ pour $p\in\NN\etl$.
Soit  
$x \in\KTO\gA$. 
Montrer que $\gamma_t(x)$ est un \pol  $t$. 
En utilisant $\gamma_{t}(-x)$, en déduire que $x$ est nilpotent.

\emph {2.}
Montrer que $\KTO(\gA)$ est le nilradical de l'anneau $\KO(\gA)$.

 On a $\rg\big(\lambda^n(x)\big)=\lambda^n(\rg x)$ et l'on dispose ainsi d'une 
série formelle $\rgl_t(x)$
à \coes dans $\HO\gA$ définie par 

\snic{\rgl_t(x) = \lambda_t(\rg x)=\sum_{n\ge 0}
\rg\big(\lambda^n(x)\big)t^n.}

%\sni
Si $x\in\HO\gA$ cela donne simplement $\rgl_t(x)=\lambda_t(x)$.
 
\emph {3.}  Si $x =
[P]$, on rappelle que $(1+t)^{\rg x} = \rR{P}(1+t) =
\rF{\Id_P}(t)\in\BB(\gA)[t]$.
Montrer que, lorsque l'on identifie $\BB(\gA)$ à $\BB(\HO\gA)$
en posant $e=[e\gA]=[e]$ pour~$e\in\BB(\gA)$, on obtient $\rgl_t(x) = 1+t\rg x=(1+t)^{\rg x}$ si $0\leq \rg x \leq 1$.\\
Montrer ensuite que $\rgl_t(x) = (1+t)^{\rg x}$ pour tout $x \in \KO(\gA)$.

\emph {4.}
On définit $\rgg_t(x) = \gamma_t(\rg x) =
\sum_{n\ge 0} \rg\big(\gamma^n(x)\big)t^n$.
\\
Montrer que $\rgg_t(x) = (1-t)^{-\rg x}$ pour tout $x \in \KO(\gA)$, ou encore
pour $x = [P]$ que $\rgg_t(x) = \rR{P}\big(1/(1-t)\big) = \rR{P}(1-t)^{-1}$. De plus,
si $0\leq\rg x \le 1$, on obtient l'\egt $\rgg_t(x) = 1+xt/(1-t) = 1 + xt + xt^2 + \dots$

\emph {5.}
Pour tout $x$ de $\KO\gA$,  ${\gamma_t(x)}{(1-t)^{\rg(x)}}$ est un \polz.

\emph {6.}
Montrer les formules de réciprocité entre $\lambda^n$ et
$\gamma^n$ pour $n \ge 1$:

\snic {
\gamma^n(x) = \sum_{p=0}^{n-1} {n-1 \choose p} \lambda^{p+1}(x),
\qquad
\lambda^n(x) = \sum_{q=0}^{n-1} {n-1 \choose q} (-1)^{n-1-q}\gamma^{q+1}(x).
}

}

\end {problem}
%--- end -problem-----------------------------------------

%--- Problem{exoApplicationProjectiveNoether}-------------
\begin{problem}\label{exoApplicationProjectiveNoether}
 {(L'application projective de \Noe et les 
\mrcs1 facteurs directs dans $\gk^2$)}\\
{\rm  
On fixe un anneau $\gk$, deux \idtrs $X$, $Y$ sur $\gk$ et un entier $n \ge1$. 
Etant donnés $x = (\xn)$, $y = (\yn)$ deux $n$-suites d'\elts de
$\gk$, on leur associe une $(n+1)$-suite $z = z(x,y) = (z_0, \ldots, z_n)$ de
la manière suivante:

\snic {
\prod_{i=1}^n (x_i X + y_i Y) = 
z_0X^n + z_1X^{n-1}Y + \cdots + z_{n-1}XY^{n-1} + z_nY^n
.}

%\sni
Ainsi, on a $z_0 = x_1\cdots x_n$, $z_n = y_1\cdots y_n$, 
et par exemple, pour $n = 3$:

\snic {
z_1 = x_1x_2y_3 + x_1x_3y_2 + x_2x_3y_1,\quad
z_2 = x_1y_2y_3 + x_2y_1y_3 + x_3y_1y_2
.}

%\sni
Pour $d \in \lrb{0..n}$, on vérifie facilement que $z_d(y,x) =
z_{n-d}(x,y)$ et que l'on a l'expression formelle suivante à l'aide des
fonctions \smqs\elrs de $n$ \idtrs $(S_0 = 1, S_1, \ldots, S_n)$:

\snic {
z_d = x_1\cdots x_n S_d(y_1/x_1, \ldots, y_n/x_n)
.}

%\sni
En particulier, $z_d$ est \hmg en $x$ de degré $n-d$, et
\hmg en $y$ de degré~$d$. On peut donner une \dfn directe
de $z_d$ de la manière suivante:

\snic {
z_d = \sum_{\#I = n-d} \prod_{i \in I} x_i \prod_{j \in \lrbn\setminus I} y_j 
.}

%\sni
Si $\gk$ est un \cdiz, on a une application $\psi : (\PP^1)^n =
\PP^1 \times \cdots \times \PP^1 \to \PP^n$, dite de \Noez, définie par:
$$
\psi : \big((x_1 : y_1), \dots, (x_n : y_n)\big) \mapsto (z_0 : \cdots : z_n)
\leqno (\star)
$$
On fait agir le groupe \smq $\rS_n$ sur le produit $(\PP^1)^n$ par permutation
des \coosz; alors l'application $(\star)$ ci-dessus, qui est $\rS_n$-invariante,
intervient en \gmt\agq pour mettre en isomorphie $(\PP^1)^n /\rS_n$ et
$\PP^n$. \perso{voir Shafarevich vol 2 exo 2 page 232}

\emph {1.}
Montrer que pour des points $P_1$, \ldots, $P_n$, $Q_1$, \ldots, $Q_n$ de
$\PP^1$,
on a 

\snic{\psi(P_1, \ldots, P_n) = \psi(Q_1, \ldots, Q_n) \iff 
(Q_1, \ldots, Q_n)$ est une permutation \hspace*{2.7cm} de $(P_1, \ldots, P_n).}

%\sni
On veut maintenant, $\gk$ étant un anneau quelconque, formuler
l'application $(\star)$ en termes de \kmrcs1.  
\\
De manière
précise, on pose $L = \gk X\oplus \gk Y \simeq \gk^2$, et l'on note 

\snic{S_n(L)
= \gk X^n \oplus \gk X^{n-1}Y\oplus \cdots \oplus \gk XY^{n-1} \oplus \gk Y^n
\simeq \gk^{n+1}}

%\sni
la composante \hmg de degré $n$ de $\gk[X,Y]$.  Si $P_1$,
\ldots, $P_n \subset L$ sont $n$ sous-\kmrcs1 facteurs directs, on veut
leur associer, de manière fonctorielle, un sous-\kmo $P = \psi(P_1, \ldots,
P_n)$ de $S_n(L)$, \prc1 et facteur direct.  Bien s\^ur, on doit avoir

\snic{\psi(P_1, \ldots, P_n) = \psi(P_{\sigma(1)}, \ldots, P_{\sigma(n)})}

%\sni
pour
toute permutation $\sigma \in \rS_n$.  De plus, si chaque $P_i$ est libre de
base $x_iX + y_iY$, alors $P$ doit être libre de base $\sum_{i=0}^n z_i
X^{n-i}Y^i$, de façon à retrouver $(\star)$.

\emph {2.}
Montrer que si chaque $(x_i, y_i)$ est \umdz, il en est de même
de $(z_0, \ldots, z_n)$.

\emph {3.}
Définir $\psi(P_1, \ldots, P_n) \subset S_n(L)$ à l'aide du
module 
$P_1 \te_\gk \cdots \te_\gk P_n$ et de l'application 
$\gk$-\lin $\pi : L^{n\te} \twoheadrightarrow S_n(L)$:

\snic {
\pi : \bigotimes_{i=1}^n (x_i X + y_iY) \longmapsto
\prod_{i=1}^n (x_i X + y_iY)
.}

%\sni
\emph {4.}
Soient $\gk[\uZ] = \gk[Z_0, \ldots, Z_n]$, $\gk[\uX, \uY] =
\gk[X_1,Y_1, \ldots, X_n,Y_n]$. \\
Que dire du $\gk$-morphisme $\varphi :
\gk[\uZ] \to \gk[\uX, \uY]$ défini par

\snic {
Z_d \longmapsto z_d = \sum_{\#I = n-d} \prod_{i \in I} X_i \prod_{j \in \lrbn\setminus I} Y_j 
~~~?}

%\sni
Note: $\varphi$ est le co-morphisme de $\psi$.

}

\end {problem}
%--- end -problem-----------------------------------------

%--- Problem{exoTh90HilbertMultiplicatif}-------------

\begin{problem}\label{exoTh90HilbertMultiplicatif} 
 {(Le \tho 90 multiplicatif d'Hilbert)}\\
{\rm
Soit $G$ un groupe fini agissant sur un anneau commutatif $\gB$; un
\emph{$1$-cocycle} de $G$ sur $\Bti$ est une famille $(c_\sigma)_{\sigma \in
G}$ telle que $c_{\sigma\tau} = c_\sigma \sigma(c_\tau)$; en conséquence,
$c_\Id = 1$. Pour tout \elt $b \in \Bti$, $(\sigma(b)b^{-1})_{\sigma \in G}$
est un $1$-cocycle appelé \emph{$1$-cobord}. \\
On note $\zcoho$ l'ensemble des
1-cocycles de $G$ sur $\Bti$; c'est un sous-groupe du groupe (commutatif) de
toutes les applications de $G$ dans $\Bti$ muni du produit à l'arrivée.
L'application $\Bti \to \zcoho$, $b \mapsto (\sigma(b)b^{-1})_{\sigma \in G}$,
est un morphisme; on note $\bcoho$ son image et l'on définit \emph{le premier
groupe de cohomologie de~$G$ sur $\Bti$} :

\snic {
\hcoho = \zcoho\sur\bcoho 
.}

%\sni
Enfin, on définit l'anneau (en \gnl non commutatif) $\tgaBG$
 comme étant \hbox{le \Bmoz} de base $G$,
muni du produit $(b\sigma) \cdot (b'\sigma') = b\sigma(b') \sigma\sigma'$.
Alors $\gB$ devient une $\tgaBG$-\alg via $(\sum_{\sigma} b_\sigma \sigma) \cdot
b = \sum_{\sigma} b_\sigma \sigma(b)$.
\\
On appelle $\tgaBG$  \emph{l'\alg tordue du groupe $G$} (twisted
group algebra of $G$).

%\sni
Soit $(\gA,\gB,G)$ une \aGz. Le but du \pb est d'associer
à tout 1-cocycle $c = (c_\sigma)_{\sigma \in G}$ un \Amrc $1$
noté $\gB_c^G$ et de montrer que $c \mapsto \gB_c^G$ définit un morphisme 
injectif de $\hcoho$ dans $\Pic(\gA)$. En particulier, si $\Pic(\gA)$
est trivial, alors tout $1$-cocycle de $G$ sur $\Bti$ est un cobord.

%\sni
\emph {1.}
Montrer que $\tgaBG \to \End_\gA(\gB)$, $\sigma \mapsto \sigma$ 
est un \iso de \Algsz.

%\sni
\emph {2.}
Soit $c \in \zcoho$. On définit $\theta_c : \tgaBG \to \tgaBG$
par $\theta_c(b\sigma) = bc_\sigma \sigma$. 
\begin {itemize}
\item [\emph {a.}]
Vérifier $\theta_c \circ \theta_{d} = \theta_{cd}$; en déduire que
$\theta_c$ est un $\gA$-\auto de $\tgaBG$.
\item [\emph {b.}]
Montrer que si $c \in \bcoho$, alors $\theta_c$ est un \auto intérieur.
\end {itemize}

\emph {3.}
Soit $c \in \zcoho$. On considère l'action de $\tgaBG$ sur $\gB$
\gui {tordue} par $\theta_c$, i.e. $z \cdot b = \theta_c(z)\, b$;
on note $\gB_c$ ce $\tgaBG$-module, $\gB_c^G$ l'ensemble des 
\elts de~$\gB$ invariants par $G$ (pour cette action tordue par $\theta_c$), et

\snic{
\pi_c = \sum_{\sigma \in G} c_\sigma\, \sigma \in \End_\gA(\gB)
.}

%\sni
Vérifier que $\gB_c^G$ est un sous-\Amo de $\gB$.  Montrer que $\pi_c$ est
une surjection de $\gB$ sur $\gB_c^G$ en en explicitant une section; en
déduire que $\gB_c^G$ est facteur direct dans $\gB$ (en tant que \Amoz).

\emph {4.}
On va montrer que pour tout $c \in \zcoho$, $\gB_c^G$ est un
\Amrcz~1.
\begin {enumerate}
\item [\emph {a.}]
Vérifier que $\gB_c^G \gB_{d}^G \simeq \gB_{cd}^G$ et $\gB_c^G \otimes_\gA
\gB_{d}^G \simeq \gB_{cd}^G$.
\item [\emph {b.}]
Montrer que si $c \in \bcoho$, alors $\gB_c^G \simeq \gA$.
Conclure.
\item [\emph {c.}]
Montrer que $c \mapsto \gB_c^G$ induit un morphisme 
injectif de $\hcoho$ dans~$\Pic(\gA)$. 
\end {enumerate}

\emph {5.}
Dans le cas où $\gA$ est un anneau \zed (par exemple un corps discret),
montrer que tout 1-cocycle $(c_\sigma)_{\sigma \in G}$ est le cobord d'un $b \in
\Bti$.

\emph {6.}
On suppose que $G$ est cyclique d'ordre $n$, $G = \gen {\sigma}$, 
et que $\Pic(\gA) = 0$. \\
Soit $x \in B$; montrer que $\rN\iBA(x)
= 1$ \ssi il existe $b \in \Bti$ tel \hbox{que $x = \sigma(b)/b$}.

}
\end {problem}
%--- end -problem-----------------------------------------

%--- problem{exoSegreMorphism}-------------
\begin{problem}\label{exoSegreMorphism} {(Le morphisme de Segre dans un cas particulier)}
\\
{\rm
Soit $\gA[\uX,\uY] = \gA[X_1, \ldots, X_n, Y_1, \ldots, Y_n]$. On considère l'idéal $\fa = \gen {X_iY_j - X_jY_i}$, i.e. l'idéal
$\cD_2(A)$, où $A$ est la matrice générique $\cmatrix {X_1 & X_2 &
\cdots & X_n\cr Y_1 & Y_2 & \cdots & Y_n\cr}$. On veut  montrer
que  $\fa$ est le noyau du morphisme:
$$
\varphi : \gA[\uX, \uY] \to \gA[T,U,\uZ] = \gA[T,U, Z_1, \ldots, Z_n],
\quad X_i\to TZ_i, \;\; Y_i\to UZ_i,
$$
où $T$, $U$, $Z_1$, \ldots, $Z_n$ sont des nouvelles \idtrsz.  Convenons
de dire qu'un \hbox{\mom $m \in \gA[\uX,\uY]$} est normalisé si $m$ est égal à $X_{i_1}\cdots X_{i_r} Y_{j_1} \cdots Y_{j_s}$ \linebreak 
{avec $1 \le i_1 \le \cdots \le i_r \le
j_1 \le \cdots \le j_s \le n$} (les indices de $\uX$ sont plus petits que
ceux de $\uY$). On note $\fa_{\rm nor}$ le sous-\Amo de
$\gA[\uX, \uY]$ engendré par les \moms normalisés.
\begin {enumerate}
\item
Si $m$, $m'$ sont  normalisés, montrer que
$\varphi(m) = \varphi(m') \Rightarrow m = m'$.
En déduire que $\Ker\varphi \cap \fa_{\rm nor} = \{0\}$.

\item
Montrer que l'on a une somme directe de \Amosz:
$\gA[\uX,\uY] = \fa \oplus \fa_{\rm nor}$

\item
En déduire que $\fa = \Ker\varphi$.
En particulier, si $\gA$ est réduit (resp. \sdzz), alors
$\fa$ est radical (resp. premier).

\end {enumerate}
}
\end {problem}
%--- end -problem-----------------------------------------

%%%%%%%%%%%%%%%%%%%%%%%%%%%%%%%%%%%%%%%%%%%%%%%%%%%%%%%%%%%%%%%%%%%
\comm
Le morphisme $\varphi$ induit, par co-morphisme, un morphisme entre espaces
affines
$$
\psi : \AA^2(\gA) \times \An(\gA) \to \MM_{2,n}(\gA) \simeq \AA^{2n}(\gA),\; \big((t,u),z\big)
\mapsto \cmatrix {tz_1 & \cdots & tz_n\cr uz_1 & \cdots & uz_n\cr}.
$$
Si $\gA$ est un corps, l'image de $\psi$ est
le lieu des zéros $\cZ(\fa)$, et   $\psi$
induit au niveau des espaces projectifs, une inclusion
$\PP^1(\gA) \times \PP^{n-1}(\gA) \to \PP^{2n-1}(\gA)$
(appelée \gui{plongement}).
\\
De manière plus \gnlez, en changeant totalement les notations,
avec des \idtrs $X_1$, \ldots, $X_n$, $Y_1$, \ldots, $Y_m$,
$Z_{ij}$, $ i \in\lrbn$, $j \in\lrbm$, considérons le
morphisme $\varphi : \gA[\uZ] \to \gA[\uX, \uY]$, $Z_{ij} \to X_iY_j$.
On montre \hbox{que $\Ker\varphi = \cD_2(A)$} où~\hbox{$A \in \MM_{n,m}(\gA[\uZ])$}
est la matrice générique.
Le morphisme $\varphi$ induit, par co-morphisme, un morphisme entre espaces affines 

\snic{\psi : \AA^n(\gA) \times \AA^m(\gA) \to \MM_{n,m}(\gA) \simeq \AA^{nm}(\gA), \;\big((x_i)_i,(y_j)_j\big)
\mapsto (x_iy_j)_{ij},}

%\sni
dont l'image est le lieu des zéros
$\cZ\big(\cD_2(A)\big)$. Si $\gA$ est un corps discret,  $\psi$ induit une injection 
$\PP^{n-1}(\gA) \times \PP^{m-1}(\gA) \to \PP^{nm-1}(\gA)$: c'est le
plongement de Segre. Cela permet de réaliser $\PP^{n-1} \times \PP^{m-1}$ \emph{comme une sous-\vgq projective de} $\PP^{nm-1}$  (en un sens précis que nous ne détaillons pas ici).
\\
Si $\gA$ est quelconque, soient $E \in \PP^{n-1}(\gA)$,
$F \in \PP^{m-1}(\gA)$; $E$ est donc un sous-module facteur direct dans $\Ae n$,
de rang~1; idem pour $F$. Alors $E\otimes_\gA F$ s'identifie canoniquement
à un sous-module de $\Ae n \otimes_\gA \Ae m \simeq \Ae {nm}$, facteur
direct, de rang~1. En posant $\psi(E, F) = E \otimes_\gA F$, on obtient ainsi
une application de~\hbox{$\PP^{n-1}(\gA) \times \PP^{m-1}(\gA)$} vers~\hbox{$\PP^{nm-1}(\gA)$}
qui \gui{prolonge} l'application précédemment définie:
si $x \in \Ae n$, $y \in \Ae m$ sont unimodulaires, il en est de même
de~\hbox{$x \otimes y \in \Ae n \otimes_\gA \Ae m$}, et en posant $E = \gA x$,
$F = \gA y$, on a $E \otimes_\gA F = \gA (x \otimes y)$.
\perso{Tout ceci s'éclairera avec \gui{le schéma $\PP^n$}, est-ce trop tôt pour y faire allusion?}

%--- problem{exoVeroneseMorphism}-------------
\begin{problem}\label{exoVeroneseMorphism} {(Le morphisme de Veronese dans un cas particulier)}
\\
{\rm
Soient $d \ge 1$, $\AuX = \gA[X_0, \ldots, X_d]$ et $\fa = \gen {X_iX_j -
X_kX_\ell, i+j = k+\ell}$. On va montrer que
l'idéal $\fa$ est le noyau du morphisme:

\snic {
\varphi : \AuX \to \gA[U,V], \qquad \varphi(X_i) = U^{d-i} V^i
.}

%\sni
où $U$, $V$ sont deux nouvelles \idtrsz. On définit
un autre idéal~$\fb$

\snic {
\fb = \gen {X_iX_j - X_{i-1}X_{j+1}, 1 \le i \le j \le d-1}
}

\begin {enumerate}
\item
Montrer que:

\snic {
\Ker\varphi \cap (\gA[X_0,X_d] + \gA[X_0,X_d]X_1 + \cdots +
\gA[X_0,X_d]X_{d-1}) = \{0\}
}

\item
Montrer que l'on a une somme directe de \Amosz:

\snic {
\AuX = \fb \oplus \gA[X_0,X_d] \oplus \gA[X_0,X_d]X_1
\oplus \cdots \oplus \gA[X_0,X_d]X_{d-1}
}

\item
En déduire que $\fa = \fb = \Ker\varphi$.
En particulier, si $\gA$ est réduit (resp. \sdzz), alors
$\fa$ est radical (resp. premier).

\end {enumerate}
}
\end {problem}
%--- end -problem-----------------------------------------

\comm
Plus \gnltz, soient $N = {n+d \choose d} = {n+d \choose n}$ et $n +
1 + N$ \idtrs $U_0$, \ldots, $U_n$, $(X_\alpha)_\alpha$,  où les indices $\alpha \in \NN^{n+1}$ sont tels
que $|\alpha| = d$. On dispose d'un morphisme $\varphi : \gA[\uX]
\to \gA[\uU]$, $X_\alpha \mapsto \uU^\alpha$ (le cas particulier
étudié ici est $n = 1 \mapsto N = d+1$); son noyau est
l'idéal

\snic {
\fa = \gen {X_\alpha X_\beta - X_{\alpha'} X_{\beta'},
\alpha + \beta = \alpha' + \beta'}
.}

%\sni

Par co-morphisme, $\varphi$ induit un morphisme entre espaces affines: 

\snic{\psi :
\AA^{n+1}(\gA) \to \AA^N(\gA), \;u = (u_0, \ldots, u_n) \mapsto
(u^\alpha)_{|\alpha| = d}.}

%\sni
Si $\gA$ est un corps discret, l'image de $\psi$ est le
lieu des zéros $\cZ(\fa)$ et l'on peut montrer que $\psi$ induit une injection
$\PP^{n}(\gA) \to \PP^{N-1}(\gA)$: c'est le plongement de Veronese de degré~$d$. 
\\
De manière encore plus \gnlez, soit $E$ un sous-module
facteur direct dans $\Ae {n+1}$, de rang~1.,
La composante \hmg de degré $d$ de l'\alg \smqz~$\gS_\gA(E)$, que l'on note $\gS_\gA(E)_d$, s'identifie à un sous-module
de $\gS_\gA(\Ae {n+1})_d \simeq \gA[U_0, \ldots, U_n]_{d}$ (composante \hmg de degré $d$), facteur
direct et de rang $1$. \\
Si l'on pose $\psi(E) = \gS_\gA(E)_d$, on \gui{prolonge}
ainsi l'application $\psi$ définie précédemment.
\eoe
%--- end - commentaire du probleme -----------------------------------------

%--- problem{exoVeroneseMatrix}-------------
\begin{problem}\label{exoVeroneseMatrix}
 {(Matrices de Veronese)}\\
{\rm
%:HHH   endo -> ali  
Soient deux anneaux de \pols $\kuX = \kXn$ et $\kuY = \kYm$.  \`A toute matrice $A
\in \gk^{m\times n}$, qui représente une \ali $\gk^n \to \gk^m$, on peut
associer (attention au renversement), un $\gk$-morphisme $\varphi_A : \kuY \to
\kuX$ construit de   la manière suivante: soient $X'_1$, \ldots, $X'_m$
les $m$ formes \lins de~$\kuX$ définies comme suit.

\snic {
\hbox{Si}\quad \cmatrix {X'_1\cr \vdots \cr X'_m} = A \cmatrix {X_1\cr \vdots \cr X_n},
\quad \hbox {alors }\quad
\varphi_A : f(\Ym) \mapsto f(X'_1,\ldots,X'_m).
} 

%\sni
Il est clair que $\varphi_A$ induit une \kli $A_d : \kuY_d \to \kuX_d$ entre
les composantes \hmgs de degré $d \ge 0$, et que la restriction \hbox{$A_1 :
\kuY_1 \to \kuX_1$} a pour matrice dans les bases $(\Ym)$ et $(\Xn)$, la \emph
{transposée} \hbox{de $A$}.  \hbox{Le \kmoz} $\kuX_d$ est libre de rang $n' = {n-1+d
\choose d} $; il possède une base naturelle, celle des
\moms de degré $d$, que l'on peut choisir d'ordonner par l'ordre
lexicographique,  avec $X_1 >
\cdots > X_n$.  
Idem pour $\kuY_d$ avec sa base de $m' = {m-1+d \choose m-1}$
\momsz.  On note $V_d(A) \in \gk^{m' \times n'}$ la \emph {transposée} de
la matrice de l'\endo $A_d$ dans ces bases (de sorte \hbox{que $V_1(A) = A$}) et l'on
dit que $V_d(A)$ est l'extension de Veronese de $A$ en degré $d$. 
\\
 Par
exemple, soit $n = 2$, $d = 2$, donc $n' = 3$; si $A = \cmatrix {a &b\cr c &
d\cr}$, on obtient la matrice $V_2(A) \in \MM_3(\gk)$ de la manière suivante

\snic {
\cmatrix {x'\cr y'\cr} = A \cmatrix {x\cr y\cr} =
\cmatrix {ax + by\cr cx + dy\cr}, \;
\cmatrix {x'^2\cr x'y'\cr y'^2\cr} = 
\cmatrix {a^2 & 2ab & b^2\cr ac & ad+bc & bd\cr c^2 & 2cd & d^2\cr}
\cmatrix {x^2\cr xy\cr y^2\cr} 
.}

%\sni
\emph {1.} 
Si $A, B$ sont deux matrices pour lesquelles le produit $AB$ a un sens,
vérifier les \egts $\varphi_{AB} = \varphi_B \circ \varphi_A$ et  $V_d(AB) =
V_d(A)V_d(B)$ pour tout $d \ge 0$.  
Vérifier \egmt que $V_d(\tra {A}) = \tra{V_d(A)}$.

\emph {2.}
Si $E$ est un \kmoz, le \emph{transformé $d$-Veronese de $E$} est le \kmo
$\gS_\gk(E)_d$, composante \hmg de degré $d$ de l'algèbre symétrique
$\gS_\gk(E)$.\\
Si $E$ est facteur  direct dans $\gk^{n}$, alors
$\gS_\gk(E)_d$ s'identifie à un sous-module de $\gS_\gk(\gk^n)_d \simeq
\gk[X_1, \ldots, X_n]_d$, facteur direct (voir aussi le \pbz~\ref
{exoVeroneseMorphism}).  Montrer que l'image par $V_d$ d'un \prr est un \prr et
que l'on a un diagramme commutatif:

\snic {
\xymatrix @C = 1.5cm{
\GA_n(\gk)\ar[d]_{\Im}\ar[r]^{V_d}         &\GA_{n'}(\gk)\ar[d]^{\Im}\\
\GG_n(\gk)\ar[r]^{d\rm -Veronese} &\GG_{n'}(\gk) \\
}
\qquad \hbox {avec} \quad n' = {n-1+d\choose d} = {n-1+d \choose n-1}
}

\emph {3.}
Montrer que si $A$ est un \prr de rang $1$, il en est de même
de $V_d(A)$.  Et plus \gnltz, si $A$ est un \prr de rang $r$, alors
$V_d(A)$ est un \prr de rang~${d+1-r \choose r-1}$.

}

\end {problem}
%--- end -problem-----------------------------------------

%--- Probleme{exoFossumKumarNori}-------------
\begin{problem}\label{exoFossumKumarNori}
{(Quelques exemples de résolutions projectives finies)}\\
{\rm  
\'Etant donnés $2n+1$ \elts $z$, $x_1$, \dots, $x_n$,  $y_1$, \dots, $y_n$,  d'un anneau $\gA$, on définit une
suite de matrices $F_k \in \MM_{2^k}(\gA)$, pour $k \in \lrb{0..n}$, de la
manière suivante:

\snic {
F_0 = \cmatrix {z}, \qquad F_k = 
\cmatrix {F_{k-1} & x_k\I_{2^{k-1}}\cr y_k\I_{2^{k-1}} &\I_{2^{k-1}} - F_{k-1}\cr}
.}

%\sni
Ainsi avec $\ov z = 1-z$:

\snic {
F_1 = \cmatrix {z & x_1\cr y_1 & \ov z}, \qquad
F_2 = \cmatrix {
z & x_1 & x_2 & 0 \cr
y_1 & \ov z &0  &x_2 \cr 
y_2 &0 &\ov z & -x_1 \cr 
0 & y_2 & -y_1 & z \cr}
.}

%\sni
\emph {1.}
Vérifier que $F_k^2 - F_k$ est la matrice scalaire de terme $z(z-1) +
\sum_{i=1}^k x_iy_i$. Montrer \egmt que $\tra{F_n}$ est semblable
à $\I_{2^n} - F_n$ pour $n \ge 1$.
En conséquence, si $z(z-1) + \sum_{i=1}^n x_iy_i = 0$,
alors $F_n$ est un \prr  de rang $2^{n-1}$.

\emph {2.}
On définit trois suites de matrices 

\snic{U_k, V_k \in \MM_{2^{k-1}}(\gA) \,
(k \in \lrb{1..n}),\quad G_k \in \MM_{2^k}(\gA)  \,(k \in \lrb{0..n}),}

%\sni
 de la
manière suivante: $U_1 = \cmatrix {x_1}$, $V_1 = \cmatrix {y_1}$, $G_0 =
\cmatrix {z}$ et:

\snic {
U_k = \cmatrix {U_{k-1} & x_k\I \cr y_k\I & -V_{k-1}}, \;
V_k = \cmatrix {V_{k-1} & x_k\I \cr y_k\I & -U_{k-1}}, \;
G_k = \cmatrix {z\I & U_k \cr V_k & \ov z\I}
.}

%\sni
Ainsi:

\snic {
U_2 = \cmatrix {x_1 & x_2\cr y_2 & -y_1}, \;
V_2 = \cmatrix {y_1 & x_2\cr y_2 & -x_1}, \;
G_2 = \cmatrix {
z & 0 & x_1 & x_2 \cr
0 & z & y_2  &-y_1 \cr 
y_1 &x_2 &\ov z & 0 \cr 
y_2 & -x_1 & 0 & \ov z \cr}
.}

%\sni
\begin{itemize}
\item [\emph {a.}]
Vérifier que $G_n$ et $F_n$ sont conjuguées par une matrice de permutation.

\item [\emph {b.}]
Vérifier que $U_kV_k$ est le scalaire $\sum_{i=1}^k x_iy_i$ et que $U_kV_k =
V_kU_k$.

\item [\emph {c.}]
Pour $n \ge 1$, si $z(z-1) + \sum_{i=1}^n x_iy_i = 0$, montrer
que $G_n$ (donc $F_n)$ est un \prr de rang $2^{n-1}$.
\end{itemize}

\emph {3.}
Soit $M$ un \Amoz. Une \emph {résolution projective finie} de $M$ est une
suite exacte de \mptfs $0 \rightarrow P_n \to \cdots \to P_1 \to P_0
\twoheadrightarrow M \to 0$; on dit que $n$ est \emph{la longueur
de la résolution}. Dans ce cas,  $M$ est
\pfz. 
\begin {itemize}
\item [\emph {a.}]
On considère deux résolutions projectives finies de $M$ que l'on
peut supposer de même longueur:

\snic {
\arraycolsep2pt\begin{array}{ccccccccccccccc}
0 &\rightarrow& P_n& \to & P_{n-1}& \to & \cdots & \to & P_1 & \to & P_0
& \to & M& \rightarrow& 0
\\
0 &\rightarrow& P'_n& \to & P'_{n-1}& \to & \cdots & \to & P'_1 & \to & P'_0
& \to & M&  \rightarrow& 0
\end {array}
.}

%\sni
En utilisant l'exercice \ref{exoSchanuelVariation}, montrer que l'on
a dans $\KO(\gA)$ l'\egt suivante:

\snic {
(\star) \qquad\qquad
\sum_{i=0}^n (-1)^i [P_i] = \sum_{i=0}^n (-1)^i [P'_i] 
.}

%\sni
\textit{Note.} L'exercice \ref{exoSchanuelVariation} fournit un résultat bien plus précis. \eoe

 \textit{Définition et notation.} Pour un module $M$ qui admet une résolution projective finie
on note $[M] \in \KO(\gA)$ la valeur commune de
$(\star)$ (même si $M$ n'est pas \ptfz). On définit alors 
\emph{le rang de $M$} comme celui de $[M]$
et l'on a $\rg M = \sum_{i=0}^n (-1)^i \rg P_i\in\HO\gA$.\index{rang!d'un module qui admet une résolution projective finie}

\item [\emph {b.}]
Soit $M$ un \Amo admettant une résolution projective finie; on
suppose \hbox{que $aM = 0$} avec $a \in \Reg(\gA)$. Montrer que
$\rg(M) = 0$ i.e. \hbox{que $[M] \in \KTO(\gA)$}.
\end {itemize}

%\sni
Si $\gk$ est un anneau quelconque, on définit l'anneau 

\snic{\gB_n =
\gk[z, \ux, \uy] = \aqo{\gk[Z, \Xn, \Yn]}{Z(Z-1) + \sum_{i=1}^n X_iY_i}}

%\sni
Ainsi $\gB_0 \simeq \gk\times\gk$.  On note $\fb_n$ l'\id
$\gen {z, \xn}$.

\emph {4.}
Montrer que les localisés $\gB_n[1/z]$ et $\gB_n[1/(1-z)]$ sont des localisés \elrs
(i.e., obtenus en inversant un seul \eltz) d'un anneau de \pols sur $\gk$ à $2n$ \idtrsz.
Montrer que $\aqo {\gB_n}{x_n} \simeq \gB_{n-1}[y_n]\simeq\gB_{n-1}[Y]$.

\emph {5.}
Pour $n=1$, définir une résolution projective du $\gB_1$-module
$\gB_1\sur{\fb_1}$ de longueur~$2$ et vérifier que
$[\gB_1\sur{\fb_1}] \in \KTO(\gB_1)$.

\emph {6.}
Pour $n=2$, définir une résolution projective du $\gB_2$-module
$\gB_2\sur{\fb_2}$ de longueur~$3$:

\snic {
0 \rightarrow \Im F_2 \to \gB_2^4 \to \gB_2^3 \vvvers{[z,x_1,x_2]} 
\gB_2 \twoheadrightarrow \gB_2\sur{\fb_2} \to 0
,}

%\sni
et vérifier que $[\gB_2\sur{\fb_2}] \in \KTO(\gB_2)$.

\emph {7.}
Expliciter une permutation $\sigma \in \rS_{2^n}$ telle que les $n+1$ premiers
\coes de la première ligne de la matrice $F'_n = P_\sigma F_n P_\sigma^{-1}$
soient $z$, $x_1$, \dots, $x_n$ ($P_\sigma$ est la matrice de la permutation $\sigma$).

\emph {8.}
Pour $n=3$, définir une résolution projective du $\gB_3$-module
$\gB_3\sur{\fb_3}$ de longueur~$4$:

\snic {
0 \rightarrow \Im (\I_8 -F'_3) \to \gB_3^8 \to \gB_3^7 \to 
\gB_3^4 \vvvvers{[z,x_1,x_2,x_3]}  \gB_3 \twoheadrightarrow \gB_3\sur{\fb_3} \to 0
,}

%\sni
et vérifier que $[\gB_3\sur{\fb_3}] \in \KTO(\gB_3)$.

\emph {9.}
Et en général?

}

\end {problem}
%--- end -problem-----------------------------------------

%:sinotenglish
\sinotenglish{

%:--- Probleme{exoPolynomialSyzygies}-------------
\begin{problem} \label{exoPolynomialSyzygies} {(Quand les \moms dominants sont 
premiers entre eux)}
\\
{\rm  
Dans ce \pbz, $\gk$ est anneau quelconque et $\kuX = \gk[\Xn]$ un anneau
de \pols avec un ordre monomial $\preceq$ fixé. \\ 
Si $m=\uX^{\alpha}$ est un \mom et~$f$
un \polz, la notation $f \preceq m$ signifie que $f$ est une
combinaison $\gk$-\lin de \moms $\preceq m$; convention analogue pour $f \prec m$.
Enfin $f \in \kuX$ est dit \emph{$\preceq$-unitaire} si  $f = m + r$
où $m$ est un \mom et $r \prec m$.

\snii
Soient  $f_1$, \ldots, $f_s \in \kuX$. On suppose que chaque $f_i$ est
$\preceq$-unitaire, de \mom dominant $m_i$ et enfin que $\pgcd(m_i,m_j) = 1$
pour $i \ne j$.  \\
On note $S$ le sous-$\gk$-module de $\kuX$ admettant base les \moms
$m$ tels que $m \notin \gen {m_1, \ldots, m_s}$:

\snic {
\displaystyle {S = \bigoplus_{m \notin \gen {m_1, \ldots, m_s}} \gk\,m}
}

\snii
Un des objectifs du \pb est de montrer que 
$$
\kuX = \gen {f_1, \ldots, f_s} \oplus S
\leqno (\star)
$$
La \dem ci-dessous fournit les moyens d'écrire tout \pol $f$ dans la
décomposition $\fa \oplus S$.

\snii
Pour ceux qui connaissent les \bdgsz: si $\gk$ est un corps, la
décomposition~$(\star)$ ci-dessus est classique; on dit parfois que les
\moms de $S$ (les \gui {\moms sous l'escalier}) sont standard: tout \pol $f$
est donc \eqv modulo $\fa$ à un et un seul \pol $\gk$-combinaison de \moms
standard, combinaison qualifiée de forme normale de $f$ modulo $\fa$
(relativement à l'ordre monomial~$\preceq$). Dans ce contexte, l'\id momonial
$\gen {m_1, \ldots, m_s}$ est l'\id initial de $\fa$ relativement à l'ordre
monomial~$\preceq$ (\id engendré par les \moms dominants des \pols
de $\fa$, on dit aussi souvent \gui{idéal de tête}).

\snii\emph {1.}
Pour un \mom $m$, on introduit deux $\gk$-modules $\fa_{\prec m}$ et
$\fa_{\preceq m}$ qui sont des sous-$\gk$-modules de l'\id $\fa = \gen {f_1,
\ldots, f_s}$,

\snic {
\fa_{\prec m} = \bigl\{ \sum_i g_i \,|\, g_i \in \gen{f_i}
\hbox { et } g_i \prec m \bigr\}
\ \subset \
\fa_{\preceq m} = \bigl\{ \sum_i g_i \,|\, 
g_i \in \gen{f_i} \hbox { et } g_i \preceq m \bigr\}
}

\snii
Il est clair que $\fa$ est la réunion des $\fa_{\prec m}$ a fortiori
des $\fa_{\preceq m}$.

\begin {itemize}
\item [\emph {a.}]
Il existe des $r_i \prec m_i$ tels que  $m_jf_i
- m_if_j = r_jf_i - r_if_j$.

\item [\emph {b.}]
Soient $m'_i, m'_j$ deux \moms tels que $m'_im_i = m'_jm_j$.\\
En notant $m$
le \mom $m'_im_i = m'_jm_j$, montrer que $m'_if_i - m'_jf_j \in \fa_{\prec m}$.

\item [\emph {c.}]
Plus \gnltz, pour une partie $I \subseteq \{1,\ldots,s\}$,
supposons disposer de \moms $(m'_i)_{i \in I}$ et d'un \mom $m$
tels que $m'_im_i = m$. Alors si
$(a_i)_{i \in I}$ est une famille finie d'\elts de $\gk$ de somme nulle,
on a $\sum_{i \in I} a_im'_if_i \in \fa_{\prec m}$.
\\
On pourra utiliser une \gui {transformation d'Abel}

\snic {
a_1 b_1 + \cdots + a_kb_k = \sum_{j=1}^{k-1} s_j(b_j-b_{j+1}),}

 avec  $s_j = \sum_{i=1}^j a_i$
 dès que  $\sum a_i = 0$.

\item [\emph {d.}]
Soit $g \in \gen {f_i}$ et un \mom $m$ tel que $g \preceq m$.  Selon que $m_i$
divise ou ne divise pas $m$, montrer que, ou bien $g \prec m$ ou bien $g -
am'f_i \in \fa_{\prec m}$ pour un $a \in \gk$ et un \mom $m'$ tel que $m'm_i =
m$.

\end {itemize}

%%%%%%%%%%%%%%%%%%%%%%%%%
\snii\emph {2.}
Soit $m$ un \momz.
Montrer que si $f \in \fa_{\preceq m}$ et $f \prec m$, alors
$f \in \fa_{\prec m}$. En déduire
que
$S \cap \fa_{\preceq m} \subset \fa_{\prec m}$. 
 Puis  $S \cap \gen {f_1, \ldots, f_s} = 0$.

%%%%%%%%%%%%%%%%%%%%%%%%%
\snii\emph {3.}
Montrer que $\kuX = \gen {f_1, \ldots, f_s} + S$.  Plus
\prmtz: si $f \in \kuX$ \hbox{et $f \preceq m$} pour un \mom $m$, alors $f \in
\fa_{\preceq m} + S$.  \\
Conclusion: $\kuX = \gen {f_1, \ldots, f_s} \oplus S$.  
\\
En outre $\gen {m_1, \ldots, m_s}$ est \emph{l'\id initial} de $\gen {f_1, \ldots, f_s}$
au sens suivant: si \hbox{pour $a \in \gk$} et $m$ \mom on a $am + h \in \gen {f_1,
\ldots, f_s}$ avec $h \prec m$, alors $am $ appartient à $\gen {m_1, \ldots, m_s}$.

%%%%%%%%%%%%%%%%%%%%%%%%%
\snii\emph {4.}
Montrer que le module des \syzys de $(f_1, \ldots, f_s)$ est engendré
par \hbox{les $f_i\vep_j - f_j\vep_i$} \hbox{ où $(\vep_1, \ldots, \vep_s)$} est la base
canonique de $\kuX^s$.

%%%%%%%%%%%%%%%%%%%%%%%%%
\snii\emph {5.}
Montrer que la suite de \moms $(m_1, \ldots, m_s)$ est \ndze et qu'il en
est de même de $(f_1, \ldots, f_s)$.

%%%%%%%%%%%%%%%%%%%%%%%%%
\snii\emph {6.}
On suppose $m_i \ne 1$ pour tout $i$.  Soit $M$ l'ensemble des \moms
appartenant à~$S$.  Montrer que $\kuX$ est un $\gk[f_1, \ldots, f_s]$-module
libre de base les $m \in M$ et que~\hbox{$(f_1, \ldots, f_s)$} sont $\gk$-\agqt
indépendants.

%%%%%%%%%%%%%%%%%%%%%%%%%
\snii\emph {7.}
Soient, dans $\gk[X,Y]$, les deux \pols $f_1 = X^2$ et $f_2 = XY + aY^2$
\hbox{avec $a \in \gk$}, tous les deux \hmgs de degré $2$ et unitaires
pour l'ordre lexi\-co\-gra\-phique~\hbox{$Y \prec X$}. 

\begin {itemize}
\item [\emph {a.}]
On a $a^2 Y^3 \in \gen {f_1,f_2}$, et  $Y^3 \in \gen {f_1,f_2}\iff a\in\gk\eti$.
\item [\emph {b.}]
On suppose $a$ \ndzz. Alors l'\id $\gen {f_1,f_2}$ est facteur direct
dans $\gk[X,Y]$ \ssi $a$ est \ivz.
\end {itemize}
}
\end{problem}

}
%: fin sinotenglish

%%%%%%%%%%%%%%%%%%%%%%%%%%%%%%%%%%%%%%%%%%%%%%%%%%%%%%%%%%%%%%%%%%%%%%%%%%%
% fin des exos

%:   solutions d'exos
\penalty-2500
\sol

%%%%%%%%%%%%%%%%%%%%%%%%%%%%%%%%%%%%%%%%%%%%%%%%%%%%%%%%%%%%%%%%%%%%%%%%%%%

%%%%%%%%%%%%%%%%%%%%%%%%%%%%%%%%%%%%%%%%%%%%%%%%%%%%%%%%%%%%%%%%%%%%%%%%%%%
\exer{exoleli2}{ 
On reprend à peu près la deuxième preuve du lemme de la liberté locale.
Notons $\varphi$ l'\ali qui a pour matrice $F$.
Appelons $f_j$ la  colonne $j$ de la matrice $F$, et
$(e_1,\ldots ,e_n)$ la base canonique de $\Ae n$.
Par hypothèse,
$(f_1,\ldots ,f_k,e_{k+1},\ldots ,e_n)$ est une base de $\Ae n$.
La matrice de passage correspondante est
$B_1=\bloc{V}{0}{C'}{\I_{h}}.$
Puisque $\varphi(f_i)=\varphi\big(\varphi(e_i)\big)=\varphi(e_i)=f_i$, par rapport à
cette  base,
$\varphi$ a une matrice du type $\bloc{\I_k}{X}{0}{Y}$.
Le calcul donne:
%---------begin $$----------
$$
B_1^{-1}=\bloc{V^{-1}}{0}{C}{\I_{h}},\qquad
G=B_1^{-1}\,F\,B_1=\bloc{\I_k}{L} {0}{-C'V^{-1}L'+W},
$$
%---------end $$----------
où $L=V^{-1}L'$, et $C=-C'V^{-1}$.\\
Puisque $\cD_{k+1}(G)=0$, on a
$G=\bloc{\I_k}{L} {0}{0}$, donc $W=C'V^{-1}L'$. 
\\
On
pose $B_2=\bloc{\I_k}{-L}{0}{\I_{h}}$, on a
$B_2^{-1}=\bloc{\I_k}{L}{0}{\I_{h}}$, puis
$B_2^{-1}\,G\,B_2=\I_{k,n}$.\\
Finalement
on obtient
$B^{-1}\,F\,B=\I_{k,n}$
avec
$$
B=B_1\,B_2=\bloc{V}{0}{C'}{\I_{h}}\cdot\bloc{\I_k}{-L}{0}{\I_{h}}=
\bloc{V}{-L'}{C'}{\I_h-W}
$$
et
$$
B^{-1}=B_2^{-1}\,B_1^{-1}=\bloc{\I_k}{L}{0}{\I_{h}}\cdot\bloc{V^{-
1}}{0}{C}{\I_{h}}=
\bloc{V^{-1}+LC}{L}{C}{\I_{h}}.
$$
L'\egt $F^2=F$ donne en particulier $V=V^2+L'C'$.\\
Donc
$\I_k=V\,(\I_k+L'C'V^{-1})=V\,(\I_k-LC)$, et finalement
$V^{-1}=\I_k-LC$. Donc comme annoncé $B^{-1}=\bloc{\I_k}{L}{C}{\I_h}$.\\
Avant de démontrer l'affirmation concernant $\I_h-W$ voyons la
réciproque. \\
La double \egt
$$
\bloc{\I_k}{L}{C}{\I_h}
= \bloc{\I_k-LC}{L}{0}{\I_h} \,\bloc{\I_k}{0}{C}{\I_h}=
 \bloc{\I_k}{L}{0}{\I_h}\,\bloc{\I_k}{0}{C}{\I_h-CL}
$$
montre que $\I_k-LC$ est \iv \ssi $\I_h-CL$ est \iv \ssi
 $\bloc{\I_k}{L}{C}{\I_h}$ est
\ivz. Cela donne aussi 
$$\det\bloc{\I_k}{L}{C}{\I_h}=\det(\I_k-LC)=\det(\I_h-
CL)\,.
$$
Le calcul donne alors
$$
{\bloc{\I_k}{L}{C}{\I_h}}^{-1}=\bloc{V}{-VL}{-CV}{\I_h+CVL},
$$
d'où
$$
{\bloc{\I_k}{L}{C}{\I_h}}^{-1} \cdot \bloc{\I_k}{0}{0}{0} \cdot
\bloc{\I_k}{L}{C}{\I_h}
=  \bloc{V}{VL}{-CV}{-CVL}\,,
$$
ce qui établit la réciproque.\\
Enfin l'\egt $B^{-1}\,F\,B=\I_{k,n}$ implique
$B^{-1}\,(\I_n-F)\,B=\I_n-\I_{k,n}$, ce qui donne
$$
\bloc{\I_k-V}{-L'}{-C'}{\I_h-W}
={\bloc{\I_k}{L}{C}{\I_h}}^{-1} \cdot \bloc{0}{0}{0}{\I_h}\cdot
\bloc{\I_k}{L}{C}{\I_h}
$$
et l'on se retrouve dans la situation symétrique,
donc $(\I_h-W)^{-1}=\I_h-CL$ \hbox{et $\det\,V=\det(\I_h-W)$}.
}

%%%%%%%%%%%%%%%%%%%%%%%%%%%%%%%%%%%%%%%%%%%%%%%%%%%%%%%%%%%%%%%%%%%%%%%%%%%

\exer{exoJacobienneP2=P}{
Noter $g$ (resp. $d$) la multiplication à gauche
(resp. à droite) par $P$. On a alors $g^2=g,\,d^2=d,\,gd=dg,\,\varphi=g+d-1$
et $\pi=g+d-2gd$.
}

%%%%%%%%%%%%%%%%%%%%%%%%%%%%%%%%%%%%%%%%%%%%%%%%%%%%%%%%%%%%%%%%%%%%%%%%%%%

\exer{exoP1FracRat} 
\emph {1a.}
La \gui{matrice de Sylvester \hmgz} $S$ est définie comme celle de l'\ali  $(A,B)\mapsto PA+QB$ sur les \hbox{bases $(u^{q-1},\ldots,v^{q-1})$} pour $A$ (\pog de degré $q-1$), 
$(u^{p-1},\ldots,v^{p-1})$ pour $B$ (\pog de degré $p-1$) et $(u^{p+q-1},\ldots,v^{p+q-1})$ pour $PA+QB$ (\pog de degré $p+q-1$).\\
En faisant $v=1$, on voit que $\tra S=\Syl(g,p,h,q)$, d'où $\det(S)=\Res(g,p,h,q)$. \\
En faisant $u=1$,
on voit que $\tra S$ est presque la matrice $\Syl(\wi g,p,\wi h,q)$: il faut renverser l'ordre des lignes, l'ordre des $q$ premières colonnes et l'ordre des $p$ dernières.
D'où le résultat annoncé car $(-1)^{\lfloor q/2\rfloor + \lfloor p/2\rfloor + \lfloor (p+q)/2\rfloor}
= (-1)^{pq}$.

\emph {1b.}
L'\egt $S \wi S=\Res(P,Q)\,\I_{p+q}$ signifie que,
si $k+\ell=p+q-1$,
 $u^kv^{\ell}\Res(P,Q)$ est une \coli des vecteurs colonnes de la matrice $S$.
Cela donne donc exactement l'inclusion requise, qui n'est en fin de compte que la version homogène de l'inclusion habituelle.

\emph {2.}
On écrit $f$ sous forme \ird $f = a/b$ avec $a$, $b \in \gk[t]$,
et l'on homogénéise~$a$ \hbox{et $ b$} en degré $d$ (maximum des degrés de $a$ et $b$)
pour obtenir deux \pogs $A$, $B \in \gk[u,v]$ de degré $d$. \\
Si $\gk$
est un anneau quelconque, on demande  que $\Res(A,B)$ soit \ivz.
Cela est \ncr pour que la fraction reste bien définie après toute \edsz.
Voyons alors que
le morphisme $f$ est d'abord défini au niveau des vecteurs \umdsz:

\snic {
(\xi : \zeta) \mapsto \big(A(\xi,\zeta) : B(\xi,\zeta)\big)
.}

%\sni
Ceci a bien un sens car si $1 \in \gen {\xi , \zeta}$, alors $1 \in \gen {A(\xi,\zeta),B(\xi,\zeta)}$ d'après le point~\emph{1b.}

Pour remonter au niveau $\GA_{2,1}(\gk)$, on prend deux nouvelles
\idtrsz~\hbox{$x$, $y$} en pensant à la matrice $\cmatrix {xu & yu\cr xv & yv\cr}$.
Comme $\gen {u,v}^{2d-1} \subseteq \gen {A,B}$, on peut écrire
$$(xu + yv)^{2d-1} = E(x,y,u,v)A(u,v) + F(x,y,u,v)B(u,v)
$$ 
avec
$E$ et $F$ \hmgs en $(x,y,u,v)$. \\
En fait, $E$ et $F$ sont bi-\hmgs en $\big((x,y), (u,v)\big)$,
de degré $2d-1$ en $(x,y)$, de degré $d-1$ en $(u,v)$. 
Comme $EA$ est bi-\hmgz, de bi-degré $(2d-1, 2d-1)$, il existe\footnote{Voir la justification ci-dessous.}
un \pog $\alpha'$ en 4 variables,
$\alpha' = \alpha'(\alpha, \beta, \gamma, \delta)$, 
tel que:

\snic {
EA = \alpha'(xu, yu, xv, yv), \quad \deg(\alpha')=2d-1
.}

%\sni
Même chose avec $FA$, $EB$, $FB$ pour produire $\beta'$, $\gamma'$, $\delta'$. On
considère alors les matrices:

\snic {
\xymatrix {\cmatrix {xu & yu\cr xv & yv\cr} \ar@{~>} [r] &
\cmatrix {\alpha & \beta\cr \gamma & \delta\cr}}, 
\qquad
\xymatrix {\cmatrix {EA & FA\cr EB & FB\cr} \ar@{~>} [r] &
\cmatrix {\alpha' & \beta'\cr \gamma' & \delta'\cr}
}.}

%\sni
Le relèvement cherché est alors
$\cmatrix {\alpha & \beta\cr \gamma & \delta\cr} \mapsto
\cmatrix {\alpha' & \beta'\cr \gamma' & \delta'\cr}$.
\\
Note: $\alpha'$, $\beta'$, $\gamma'$, $\delta'$ sont des \pogs
en $(\alpha, \beta, \gamma, \delta)$, de degré $2d-1$,
tels que:

\snic {
\left| \matrix {\alpha & \beta\cr \gamma & \delta\cr}\right|
\hbox { divise }
\left| \matrix {\alpha' & \beta'\cr \gamma' & \delta'\cr}\right|,
\qquad
\alpha + \delta - 1 \hbox { divise } \alpha' + \delta' - 1 
.}

%\sni
Justification de l'existence de $\alpha'$. \\
Cela repose sur le fait simple
suivant: $u^iv^jx^ky^\ell$ est un monôme en $(xu, yu,
xv, yv)$ \ssi $i+j = k+\ell$; en effet, si cette \egt est vérifiée,
il y a une matrice $\cmatrix {m & n\cr r & s\cr} \in \MM_2(\NN)$ telle
%:2012 lignes
que les sommes de lignes soient $(i,j)$ et les sommes de colonnes soient 
$(k,l)$.
Un schéma pour aider la lecture:

%:2012 \bordercmatrix[\lbrack\rbrack]
\snic {
\bordercmatrix[\lbrack\rbrack] {
   & k  & \ell \cr
i  & m  & n \cr
j  & r  & s \cr}
\qquad
\cmatrix {xu & yu\cr xv & yv\cr}
,}

%\sni
et alors

\snic {
u^iv^jx^ky^\ell = u^{m+n}v^{r+s}x^{m+r}y^{n+s}
= (xu)^m (yu)^n (xv)^r (yv)^s.
}

%\sni
On en déduit qu'un \pol bi-\hmg en $\big((x,y), (u,v)\big)$, de
bi-degré $(d,d)$, est l'\evn en $(xu, yu, xv, yv)$ d'un \pog de
degré $d$.

\emph {3.}
Pour $f(t) = t^2$, on obtient le relèvement:

\snic {
\cmatrix {\alpha & \beta\cr \gamma & \delta\cr} \mapsto
\cmatrix {\alpha^2(\alpha + 3\delta) & \beta^2(3\alpha + \delta)\cr 
\gamma^2(\alpha + 3\delta) & \delta^2(3\alpha + \delta)\cr}
.}

%\sni

Plus \gnltz,  on développe $(\alpha + \delta)^{2d-1}$
sous la forme $\alpha^d S_d(\alpha,\delta) + \delta^d S_d(\delta,\alpha)$,
et l'on obtient le relèvement:

\snic {
\cmatrix {\alpha & \beta\cr \gamma & \delta\cr} \mapsto
\cmatrix {\alpha^d S_d(\alpha,\delta) & \beta^d S_d(\delta, \alpha)\cr 
\gamma^d S_d(\alpha,\delta) & \delta^d S_d(\delta, \alpha)\cr}
.}

%\sni
Si $H = \cmatrix {a & b\cr c & d\cr}$, on obtient le relèvement 
suivant de $f(t) = {at + b \over ct +d}$:
$$\cmatrix {\alpha & \beta\cr \gamma & \delta\cr} \mapsto
H \cmatrix {\alpha & \beta\cr \gamma & \delta\cr} H^{-1}.$$

%%%%%%%%%%%%%%%%%%%%%%%%%%%%%%%%%%%%%%%%%%%%%%%%%%%%%%%%%%%%%%%%%%%%%%%%%%%

\exer{exoConiqueFondamentale} 
On procède comme dans l'exercice \ref {exoP1FracRat}, mais c'est plus simple
car puisque $\gen {u,v}^2 = \gen {u^2, uv, v^2}$, l'application $(u : v) \mapsto (u^2 : uv : v^2)$ est bien définie au niveau des
\vmdsz.%, comme restriction de: $\gk^2 \to \gk^3$, $\cmatrix {u\cr
%v\cr} \mapsto \cmatrix {u^2\cr uv\cr v^2\cr}$.

On introduit $(x, y)$ en pensant à la matrice $\cmatrix {\alpha & \beta\cr
\gamma & \delta\cr} \leftrightarrow \cmatrix {xu & yu\cr xv & yv\cr}$.
On développe $(xu + yv)^2 = x^2 u^2 + 2xy uv + y^2 v^2$, somme de trois termes
qui vont être les  termes diagonaux d'une matrice de $\GA_{3,1}(\gk)$, puis
on complète de façon à ce que chaque colonne soit le multiple ad-hoc
du vecteur $\tra[\,u^2\;uv\;v^2\,]$. Ce qui donne:

\snic {
\cmatrix {
x^2 u^2 & 2xy u^2 & y^2 u^2 \cr
x^2 uv  & 2xy uv  & y^2 uv \cr
x^2 v^2 & 2xy v^2 & y^2 v^2 \cr
}
\qquad
F = \cmatrix {
\alpha^2     & 2\alpha\beta  & \beta^2 \cr
\alpha\gamma & 2\alpha\delta  & \beta\delta \cr
\gamma^2     & 2\gamma\delta & \delta^2 \cr
}.}

%\sni
Le relèvement $\GA_{2,1}(\gk) \to \GA_{3,1}(\gk)$ est 
$\cmatrix {\alpha & \beta\cr \gamma & \delta\cr} \mapsto F$.
\\
On a bien s\^ur $\Tr(F) = (\alpha + \delta)^2 = 1$,
$\cD_2(F) \subseteq \gen {\alpha\delta - \beta\gamma} = 0$,
et $F$ est un \prr de rang~1.

%%%%%%%%%%%%%%%%%%%%%%%%%%%%%%%%%%%%%%%%%%%%%%%%%%%%%%%%%%%%%%%%%%%%%%%%%%%
\exer{exoProjecteurCorangUn}{
\emph{1.}
On fournit deux solutions pour cette question. La première
consiste à utiliser  l'expression de
l'adjointe en fonction de la matrice de départ; la seconde preuve utilise
la \lonz.

Pour $A \in \Mn(\gA)$ on a l'expression classique de
 $\wi{A}$ comme \pol en~$A$:

\snic{\wi A = (-1)^{n-1}Q(A) \;$  avec $\;XQ(X)=\rC{A}(X)-\rC{A}(0).}

%\sni
Appliquons ceci à un \prr $P$ de rang $n-1$.
Il vient 
$$\rC{P}(X) =  (X-1)^{n-1}X, \,Q(X)=(X-1)^{n-1}\hbox{ et }
 (P - \In)^{n-1}  = (-1)^{n-1} \wi P.$$
Puisque $(\In-P)^{n-1}=\In-P$, on obtient
$ P + \wi {P} = \In.$

Voici la preuve par \lonz.
D'après le \tho de structure locale des \mptfs (\thrf{prop Fitt ptf 2}
ou \thrf{th ptf loc libre}), il existe des
\lons \come telles que sur chaque localisé, $P$ est semblable à $\I_{r,n} $,
où l'entier $r$ dépend a priori de la \lonz. Ici, puisque $P$ est de
rang  $n-1$, on a $r=n-1$ ou $1=0$. Donc $P + \wi {P} = \In$ sur chaque localisé.
Et l'\egt est aussi vraie globalement d'après le \plg de base.\iplg

\emph{2.}
Voyons la \dem par \lons \comez.  Sur le localisé $\gA_s$, le \prr $P$ est
semblable à  $Q_s=\I_{r,n}$, où $r$ dépend de $s$.
On a $Q_s + \wi {Q_s} = \In$. 
\\
Si $r < n-1$, alors
$Q_s+\wi {Q_s}= \I_{r,n}$.  Si $r = n$, alors $Q_s+\wi {Q_s}  = 2\,\In$. 
\\
Bilan: si $r\neq n-1$, alors $1=0$ et le rang est aussi égal à $n-1$.
En conséquence sur tous les localisés $\gA_s$, le \prr $P$ est
de \hbox{rang $n-1$}, et donc globalement aussi.

\emph{3.}
Il suffit de multiplier $P + \wi {P} = \In$ par $P$ pour obtenir $P^2 = P$.
}

%%%%%%%%%%%%%%%%%%%%%%%%%%%%%%%%%%%%%%%%%%%%%%%%%%%%%%%%%%%%%%%%%%%%%%%%%%%
\exer {exoMatriceCorangUn}{
Il existe $B \in \Mn(\gA)$ telle que $ABA = A$, de sorte que $AB$
est un \prr de même image que~$A$, donc de rang $n-1$, et $BA$ un
\prr de même noyau que $A$, donc \egmt de rang $n-1$. On
définit $P$ et $Q \in \Mn(\gA)$ par
$AB = \In - P$, et $BA = \In - Q.$
\\
Ainsi $P$, $Q\in\GA_{1,n}(\gA)$, avec
$A = (\In-P)A = A(\In-Q)$.

\emph{1.}
On a $\det A = 0$, i.e. $\wi {A} A = A \wi {A} = 0$, donc $\Im A
\subseteq \Ker \wi {A}$.  \\
Ensuite
$\wi{AB} = \wi{\In - P}
= P$ (car $P\in\GA_{1,n}(\gA)$).
Et l'\egt  $\wi {B} \wi {A} = P$ prouve que:

\snic{\Ker \wi {A} \subseteq \Ker P = \Im(\In - P) = \Im A.}

%\sni
Conclusion: $\Ker \wi {A}  = \Im A = \Im(\In - P)$.

\emph{2.}
En raisonnant comme on point \emph{1}, on obtient $\Im \wi {A} \subseteq \Ker A = \Ker (BA)   = \Im Q$, puis $\wi {A} \wi {B} = \wi {BA}   = \wi {\In - Q}
= Q$, et $\Ker A = \Im \wi {A} = \Im Q$.

\emph{3.}
On applique le point \emph{1} à $\tra A$, donc $\Im \tra A = \Ker \tra
{\wi{A}}$.  Ensuite, on explicite le \prr ``de gauche'' (de rang~$1$)
associé à $\tra A$. On a:

\snic{\tra A \tra B \tra A = \tra A,
\quad \hbox {que l'on écrit} \quad
\tra (BA) \tra A = \tra A
\quad \hbox {avec} \quad
\tra (BA) = \In - \tra Q .}

%\sni
Ce \prr de gauche est donc $\tra Q$.

\emph{4.}
De même, le point \emph{2} donne $\Im \tra {\wi {A}} = \Ker \tra A$.  On
explicite le ``\prr de droite'' (de rang~1) associé à $\tra A$,
on obtient $\tra P$, d'où le résultat annoncé.

\emph{5.}
Enfin:

\snic{
\Ae n  / \Im A = \Ae n  / \Im(\In - P) \simeq \Im P, \qquad
\Ker \tra A = \Im \tra P,}

%\sni
donc les deux modules (\pros de rang $1$) sont bien duaux l'un de
l'autre.
\\
Remarque : on peut \egmt utiliser

\snic{\Ae n  / \Im \tra A = \Ae n  / \Im (\In - \tra Q) \simeq
\Im \tra Q, \qquad
\Ker A = \Im Q
,}

%\sni
pour voir que les deux modules (\pros de rang $1$) $\Ae n  / \Im \tra A$
et $\Ker A$ sont bien duaux l'un de l'autre.
}

%%%%%%%%%%%%%%%%%%%%%%%%%%%%%%%%%%%%%%%%%%%%%%%%%%%%%%%%%%%%%%%%%%%%%%%%%%%

\exer{exoIntersectionSchemasAffines}
Tout d'abord on remarque que les flèches surjectives $\kuX\vers{\pi_1}\gA$ et  $\kuX\vers{\pi_2}\gB$ dans la catégorie des \klgs \pf sont vues, du point de vue des schémas, comme
des \gui{inclusions} $A\vers{\iota_1}\gk^n$ et $B\vers{\iota_2}\gk^n$, \hbox{où $\gk^n$} est interprété comme le schéma affine correspondant à $\kuX$. La \dfn de l'intersection par produit tensoriel est donc en fait une \dfn comme somme amalgamée des deux flèches $\pi_1$ et $\pi_2$ dans la catégorie \hbox{des \klgsz} \pfz, 
ou comme produit fibré des deux flèches $\iota_1$ et $\iota_2$ dans la catégorie des
schémas affines sur $\gk$.    

Le centre de l'ellipse, le centre du cercle et point d'intersection double ont pour \coos
respectives $(0,0)$, $(c,0)$ et $(a,0)$.
Le calcul
des autres points d'intersection donne
$x = a(2ac + 1 - a^2)/(a^2 - 1)$ et $y^2 = 4ac(a^2 - ac - 1) / (a^2 - 1)^2$.

Du point de vue des \algs quotients on obtient

\snic {
\gA = \aqo{\gk[X,Y]}{f},\quad \gB = \aqo{\gk[X,Y]}{g}, \quad
\gC = \aqo{\gk[X,Y]}{f,g}.
}

%\sni
Ce qui donne  les morphismes

\snic {
\xymatrix @R=10pt @C=25pt{
              & \gK=\gk[X,Y]\ar[dl]_(.6){\pi_1} \ar[dd]\ar[dr]^(.6){\pi_2} \\
\gA\ar[dr]    &              &  \gB\ar[dl] \\
              & \gC = \gA\otimes_\gK \gB \\
}
\quad
\xymatrix @R=10pt @C=8pt{
                   & \gk^2   \\
\{f = 0\}\ar[ur]^{\iota_1}  &              & \{g = 0\}\ar[ul]_{\iota_2} \\
            &\{f = 0\} \cap \{ g = 0\}\ar[ul]\ar[uu]\ar[ur] \\
}
}

%\sni
Si $\gk$ est un \cdi   et si 
$4ac(a^2 - ac - 1) (a^2 - 1)\in \gk\eti$,
les \klgs $\gA$ et~$\gB$ sont intègres, mais pas $\gC$: on a un \iso

\snic{\gC\simarrow\gk[\zeta]\times \gk[\vep],\;\hbox{ où }\,\vep^2=0\;
\hbox{ et }\,\zeta^2=4ac(a^2 - ac - 1) /(a^2 - 1)^2.}

%\sni
L'\alg $\gC$ est un \kev de dimension 4,
correspondant au schéma affine formé par deux points de multiplicité 1 
(définis sur $\gk$ ou sur une extension quadratique de
$\gk$) et un point de multiplicité 2 (défini sur $\gk$).
%%%%%%%%%%%%%%%%%%%%%%%%%%%%%%%%%%%%%%%%%%%%%%%%%%%%%%%%%%%%%%%%%%%%%%%%%%%

\exer{exoPolPseudoUnitaire} 
On rappelle que pour un \idm $e$, on a $\gen {a,e} = \gen {(1-e)a + e}$; si
$e'$ est un autre \idm \ort à $e$, on a $\gen {\ov a} = \gen {\ov {e'}}$ dans
$\aqo{\gA}{e}$ \ssi $\gen {(1-e)a} = \gen {e'}$ dans $\gA$.

\emph {1.}
Pour $k > r$, on a $a_k = 0$ dans chaque composante, donc dans $\gA$.  L'\elt
$a_r$ est nul dans $\aqo{\gA}{e_r}$, \iv dans $\aqo{\gA}{1-e_r}$ donc $\gen
{a_r} = \gen {e_r}$.
\\
 De même dans $\aqo{\gA}{e_r}$, on a $\gen {\ov
{a_{r-1}}} = \gen {\ov {e_{r-1}}}$ donc $\gen {(1-e_r)a_{r-1}} = \gen
{e_{r-1}}$. Et ainsi de suite.

\emph {2.}
Localiser en chacun des $e_i$.

%%%%%%%%%%%%%%%%%%%%%%%%%%%%%%%%%%%%%%%%%%%%%%%%%%%%%%%%%%%%%%%%%%%%%%%%%%%
\exer{exoPolLocUnitaire}
\emph {1.}
Comme $f(t) = 0$, on a $f \in \fa$, d'où une \Ali
surjective $\aqo{\gA[T]}{f} \twoheadrightarrow \gA[T]\sur\fa$ entre deux
\Amos libres de même rang $n$: c'est un \iso
(proposition \ref{propDimMod1}), donc $\fa = \gen {f}$.

\emph {2.}
Le \polcar $f$ de $t$ est \mon de degré $n$ car
$\gA[t]$ est de rang constant $n$.
Comme $f(t) = 0$, on a $f \in \fa$, d'où une \Ali surjec-\linebreak 
tive
$\aqo{\gA[T]}{f} \twoheadrightarrow \gA[T]\sur\fa$, d'un
\Amo libre de rang $n$ sur un $\gA$-\mrc $n$;
c'est donc un \iso (proposition \ref{prop epi rang constant}),
\hbox{donc $\fa = \gen {f}$}.

\emph {3.}
Soit $f = \sum_{i=0}^r a_i T^i=\sum_{i=0}^r f_r$ un \pol \lot\mon
de degré formel $r$, avec le \sfio
$(e_0,\ldots,e_r)$, \hbox{et $fe_d=f_d$} \mon de degré $d$ modulo $\gen{1-e_d}$ pour chaque $d\in\lrb{0..r}$.
\\
Alors $a_r=e_r$ est \idmz. Ensuite $f-f_r=(1-e_r)f$
est \lot \mon de degré formel
$r-1$ et l'on peut terminer par \recu descendante sur $r$ pour calculer les $e_d$
à partir de $f$.
Si l'anneau est discret on obtient un test pour décider si un \pol donné est
\lot \monz: chacun des $e_d$ calculés successivement doit être \idm et la somme
des $e_d$ doit être égale à~$1$.

%%%%%%%%%%%%%%%%%%%%%%%%%%%%%%%%%%%%%%%%%%%%%%%%%%%%%%%%%%%%%%%%%%%%%%%%%%%

\exer{exoOneRankPtfIdeal} 
\emph {1.} Il existe un sous-\Amo $N$ de $\gB$ tel
que $M.N=\gA$.\\
On a $(\xn)$ dans $M$ et $(\yn)$ dans $N$ tels que $1=\sum_ix_iy_i$ et
$x_iy_j\in\gA$. On vérifie  que $M = \sum_i \gA x_i$ et $N =
\sum_i \gA y_i$.  Soit $\sum_kz_k\otimes z'_k$ dans $M\otimes_\gA
M'$. On a, en remarquant que $y_iz_k\in N.M = \gA$

\snic {\arraycolsep2pt
\begin{array}{rclcl} 
 \sum_kz_k\otimes z'_k &  = & \sum_{k,i}x_iy_iz_k\otimes z'_k  & =  & \sum_{k,i} x_i\left(y_iz_k\right)\otimes z'_k    \\[1mm] 
  & =  &  \sum_{k,i} x_i\otimes \left (y_iz_k\right)z'_k  &
=  & \sum_ix_i\otimes \left(y_i \sum_k  z_k z'_k\right),    
 \end{array}
}

%\sni
donc la surjection canonique $M\otimes_\gA M'\to M.M'$ est injective.

\emph {2.}
Il faut montrer que $\fa$ contient un \elt \ndz 
(lemme \ref{lemIdproj}~\emph{\iref{i6lemIdproj}}), ce qui est \imdz.

%%%%%%%%%%%%%%%%%%%%%%%%%%%%%%%%%%%%%%%%%%%%%%%%%%%%%%%%%%%%%%%%%%%%%%%%%%%

\exer{exoPicAPicFracA} 
Définir la suite va de soi; ainsi, l'application $\gK\eti \to \Gfr(\gA)$ est
celle qui à $x \in \gK\eti$ associe l'\id fractionnaire principal $\gA
x$. Pas de problème non plus pour vérifier que le composé de deux
morphismes consécutifs est trivial. 
\\
\emph{Exactitude en $\gK\eti$:} si $x\in
\gK\eti$ est tel que $\gA x = \gA$, alors $x\in\Ati$.  
\\
\emph{Exactitude en
$\Gfr(\gA)$:} si $\fa \in \Gfr(\gA)$ est libre, cela signifie qu'il est
principal i.e. de la forme $\gA x$ avec $x\in\gK\eti$.

Seule l'exactitude en $\Pic\gA$ est plus délicate.  De manière \gnlez, si
$P$ est \hbox{un \Amoz} \ptfz, alors l'application canonique $P \to \gK\te_\gA P$
est injective car $P$ est contenu dans un \Amo libre. Soit donc $P$ un
$\gA$-\mrc 1 tel que $\gK\te_\gA P \simeq \gK$. Alors $P$ s'injecte dans $\gK$
puis dans $\gA$ (multiplier par un dénominateur), i.e. $P$ est isomorphe à
un \id entier~$\fa$ de $\gA$. De même, le dual $P\sta$ est isomorphe à un
\id entier $\fb$ de $\gA$. \\
Et l'on a $\gA \simeq P\te_\gA P\sta \simeq \fa
\te_\gA \fb \simeq \fa\fb$, donc $\fa\fb$ est engendré par un \elt \ndz $x
\in \gA$. On a $x \in \fa$ donc $\fa$ est un \id \iv: on a trouvé un \id \iv
$\fa$ de $\gA$ tel que $\fa \simeq P$.

%%%%%%%%%%%%%%%%%%%%%%%%%%%%%%%%%%%%%%%%%%%%%%%%%%%%%%%%%%%%%%%%%%%%%%%%%%%

\exer{exoCoMatCoRang1}
\emph{1} et \emph{2.}
Immédiat.

\emph{3.}
On considère la courte suite $\Ae N \vers {A'}
\Ae m \vers {A}  \Ae n$; elle est exacte \lotz, donc
globalement.

\emph{4.}
Tout module \stl de rang 1 peut être donné sous la forme $\Ker A$
où $A
\in \Ae {n \times (n+1)}$ est une matrice surjective $\Ae {n+1} \vvers {A}  \Ae n$. 
Puisque $1 \in \cD_n(A)$, on applique la question \emph {3} avec $m = n+1$. On obtient
$A' \in \Ae {(n+1) \times 1}$ de rang $1$ avec $\Im A' = \Ker A$; donc la
colonne $A'$ est une base de $\Ker A$.

%%%%%%%%%%%%%%%%%%%%%%%%%%%%%%%%%%%%%%%%%%%%%%%%%%%%%%%%%%%%%%%%%%%%%%%%%%%

\exer{exoVarProj}
Soit $f\in\gk[\Xzn]$ un \pog de degré $m$ et (pour simplifier) 
$P=\gen{\ua,\ub,\uc}\subseteq\gk^{n+1}$
un  facteur direct de rang $1$. 
On suppose que $f(\ua)=f(\ub)=f(\uc)=0$
et l'on veut montrer que $f(\ux)=0$ si $\ux=\alpha\ua+\beta\ub+\gamma\uc$. 
La matrice
de $(\ua, \ub, \uc)$ est de rang $1$, donc les  $a_i$, $b_j$, $c_k$ sont \comz. Il suffit donc de prouver l'\egt après \lon en une de ces \coosz.
Par exemple sur $\gk[1/a_0]$ on a $\ux=(\alpha+\frac{b_0}{a_0}\beta+\frac{c_0}{a_0}\gamma)\ua=\lambda\ua$,
et \hbox{donc $f(\ux)=\lambda^mf(\ua)=0$}.

%%%%%%%%%%%%%%%%%%%%%%%%%%%%%%%%%%%%%%%%%%%%%%%%%%%%%%%%%%%%%%%%%%%%%%%%%%%

\exer{exoGLnTangent}
On considère la \klg $\gk[\vep]=\aqo{\kT}{T^2}$.\\ 
Soit  $A\in\GL_n(\gk)$ et $H \in \Mn(\gk)$. On a 
$A + \vep H=A(\In+\vep A^{-1}H)$. Et $\In+\vep M$ est inversible, d'inverse
$\In-\vep M$, pour tout $M\in\Mn(\gk)$. Donc 
 $A + \vep H \in \GLn(\gk)$ pour n'importe quel $H$.
Ainsi, l'espace tangent $\rT_A(\GL_n)$ est isomorphe à $\Mn(\gk)$.
\\
NB:  
$(A + \vep H)^{-1}=A^{-1}-\vep A^{-1}H A^{-1}$.

%%%%%%%%%%%%%%%%%%%%%%%%%%%%%%%%%%%%%%%%%%%%%%%%%%%%%%%%%%%%%%%%%%%%%%%%%%%

\exer{exoSLnTangent}
 On utilise la \klg $\gk[\vep]$ de l'exercice \ref{exoGLnTangent}.  
 Pour $A$, $H \in \Mn(\gk)$, on a $\det(A + \vep H) =
\det(A) + \vep \Tr(\wi{A}H)$. On en déduit

\snic {
\det(A + \vep H) = 1  \iff  (\det(A) = 1  \hbox { et } \Tr(\wi{A}H) = 0) 
.}

%\sni
On a donc, pour $A \in \SLn(\gk)$, 
$\rT_A(\SLn) = \sotQ {H \in \Mn(\gk)} {\Tr(\wi{A}H) = 0}$.\\
Montrons que $\rT_A(\SLn)$ est un $\gk$-module libre de rang
$n^2 -1$. \\
En effet, l'\auto $\gk$-\lin $H \mapsto AH$  de
$\Mn(\gk)$ transforme $\In$ en~$A$ et applique bijectivement
$\rT_{\In}(\SLn)$ sur $\rT_A(\SLn)$, comme on peut le vérifier en écrivant
$\Tr(H) = \Tr(\wi{A}\, AH)$.  Enfin $\rT_{\In}(\SLn)$ est le sous-$\gk$-module de
$\Mn(\gk)$ constitué des matrices de trace nulle (qui est bien libre de rang $n^{2}-1$).
\\
NB:  $H
\mapsto HA$ était aussi possible, car $\Tr(AH\, \wi{A}) = \Tr(\wi{A}\, AH) = \Tr(H)$.

%%%%%%%%%%%%%%%%%%%%%%%%%%%%%%%%%%%%%%%%%%%%%%%%%%%%%%%%%%%%%%%%%%%%%%%%%%%

\exer{exoTangentJ0ConeNilpotent}
\emph{1.}
On voit facilement que $\varphi(H)J_0 = 
J_0\varphi(H)$. Si $\gk$ était un corps, on
pourrait en déduire  que $\varphi(H)$
est un \pol en $J_0$. Le calcul direct
donne

\snic{\varphi(e_{ij}) = \cases {
0 & si $i < j$ \cr
J_0^{n-1-(i-j)} & sinon .}}

%\sni
En particulier, $\,
\varphi(e_{i1}) = J_0^{n-i}$.
On a donc $\Im \varphi = \bigoplus_{k=0}^{n-1} \gk J_0^{k}$.

%
%\sni
\emph{2.}
Pour $k \in \lrb{0.. n-1}$, la matrice $J_0^k$ a ses \coes nuls,
sauf ceux qui sont sur la $k$-ième sur-diagonale, tous égaux à 1.
On peut donc prendre comme \supl de $\Im\varphi$ le
sous-module engendré par les $e_{ij}$, avec $j < n$ (on omet
donc \hbox{les $e_{in}$} qui correspondent à la dernière position
des sur-diagonales des $J_0^k$). On  définit alors $\psi$
par

\snic{\psi(e_{ij}) = \cases {
0 & si $j < n$ \cr
e_{i1} & si $j = n$ \cr}\quad \hbox {ou encore} \quad
\psi(H) = H \tra {J_0^{n-1}}.}

%\sni
On vérifie facilement que $\psi(J_0^{n-i}) = e_{i1}$ pour $ i \in \lrbn$, puis
$(\varphi \circ \psi)(A) = A$ \hbox{si $A \in \Im\varphi$}, et enfin $\varphi \circ
\psi \circ \varphi = \varphi$. Par miracle, on a aussi
$\psi \circ \varphi \circ \psi = \psi$. \\
On a $e_{ij} - e_{i'j'} \in \Ker\varphi$ dès que $i'-j' = i-j$ ($i' \ge j'$,
$i \ge j$) et l'on obtient une base de $\Ker\varphi$ en considérant les
$n(n-1) \over 2$ matrices $e_{ij}$ avec $i < j$ et les $n(n-1) \over 2$
matrices $e_{i1} - e_{i+r,1+r}$, $r\in \lrb{1.. n-i}$, $i \in \lrb{1.. n-1}$.

%
%\sni
\emph{3.}
On utilise la \klg $\gk[\varepsilon]\simeq\aqo\kT{T^2}$.  
Pour $A, H \in \Mn(\gk)$, on a 

\snic{(A + \varepsilon H)^n = A^n + \vep\som_{i+j=n-1} A^i H A^j.}

%\sni
Pour $A = J_0$, on trouve que l'espace tangent \gui {au cône 
nilpotent} est $\Ker \varphi$ qui est un module libre de rang
$n^2 - n$ (c'est la dimension du cône nilpotent).

%: sol exer   exoMonomialSyzygies
\exer{exoMonomialSyzygies} \emph{(Syzygies entre \momsz)}\\
Soit~$r \in \kuX^s$; dire que~$r \in \sum_{i,j} \kuX (m_{ij}\vep_i - m_{ji}\vep_j)$
c'est dire qu'il existe des~$r_{ij} \in \kuX$ vérifiant~$r_{ii} = 0$,
$r_{ij} + r_{ji} = 0$ et

\snic {
r = \sum_{i,j} r_{ij}m_{ij}\vep_i = \sum_i \left( \sum_j r_{ij}m_{ij}\right) \vep_i.
}

\snii
Soit une relation~$\sum_i u_im_i = 0$ avec~$u_i \in \kuX$. En considérant la
composante sur un \mom quelconque fixé~$m$, on obtient un terme~$a_im'_i$ de
$u_i$ avec~\hbox{$a_i \in \gk$,~$m'_i$} \mom tel que~$m'_im_i = m$ et~$\sum_i a_im'_i
m_i = 0$, i.e.~$\sum_i a_i = 0$. Il est entendu \hbox{que~$a_i = 0$} si~$m_i \nedivi
m$ et l'on peut se limiter dans la suite aux~$m_i$ tels que~\hbox{$m_i \divi m$}. Il
suffit donc de montrer que~$\sum_i a_im'_i \vep_i$ est dans le module engendré
par les~\hbox{$m_{ij}\vep_i - m_{ji}\vep_j$}.

\snii
Puisque~$m'_im_i = m = m'_jm_j$, le \mom $m$ est divisible par~$m_i$, par
$m_j$ donc par leur ppcm~$m_i \vee m_j$; par conséquent, on peut écrire:

\snic {
m = q_{ij} (m_i \vee m_j),  \quad\hbox {et l'on a donc}\quad q_{ij} = q_{ji}.
}

\snii
Il vient:

\snic {
m'_im_i = q_{ij} (m_i \vee m_j), \quad\hbox {donc}\quad
m'_i =  {q_{ij} {m_i \vee m_j \over m_i} = 
q_{ij} {m_j \over m_i \vi m_j} = q_{ij} m_{ij}}.
}

\snii
D'autre part, puisque~$\sum_i a_i = 0$, il existe une matrice
anti\smqz~\hbox{$(a_{ij}) \in \MM_s(\gk)$} telle que 
$a_i = \sum_j a_{ij}$ (la somme sur la ligne~$i$ vaut~$a_i$).
Par exemple, pour~$s = 4$:

\snic {
\cmatrix {
0 & -a_2 & -a_3 &-a_4 \cr
a_2 & 0 & 0 & 0 \cr
a_3 & 0 & 0 & 0 \cr
a_4 & 0 & 0 & 0 \cr} \cmatrix {1\cr 1\cr 1\cr 1\cr} =
\cmatrix {a_1\cr a_2\cr a_3\cr a_4\cr}
}

\snii
On écrit alors (comme par magie en ayant gratté un peu sur son brouillon):

\snic {
a_im'_i = \sum_j a_{ij} m'_i = \sum_j a_{ij} q_{ij} m_{ij}. 
}

\snii Et il ne reste plus qu'à poser~$r_{ij} = a_{ij}q_{ij}$:
on a bien~$r_{ii} = 0$,~$r_{ij} + r_{ji} = 0$ \hbox{et
$\sum_i a_i m'_i \vep_i = \sum_{i,j} r_{ij} m_{ij} \vep_i$}.

%: sols pbs
%%%%%%%%%%%%%%%%%%%%%%%%%%%%%%%%%%%%%%%%%%%%%%%%%%%%%%%%%%%%%%%%%%%%%%%%%%%

\prob{exoAnneauCercle1} \emph{(L'anneau du cercle)}\\
\emph{1.}
De manière na\"ive: soit $f = f(x,y) \in \gk[x,y]$ une conique,
i.e. un \pol de degré 2, et $(x_0, y_0)$ un $\gk$-point de $\so{ f(x,y) = 0
}$.

\addhabille1\Habillage{\includegraphics{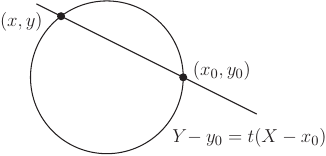}}{0}{0pt}
L'astuce classique de paramétrage consiste à définir $t$ par $y-y_0
$ $=$ $t(x-x_0)$ et, dans l'\eqn
$$f(x,y) = f\big(x, t_0 + t(x-x_0)\big) = 0,$$
à
chercher $x$ en fonction de $t$. Cette \eqn 
admet $x = x_0$ comme solution, d'où l'autre
solution sous forme rationnelle.
\endHabillage

Algébriquement parlant, on supppose $f$ \irdz, on pose $\gk[x,y] =
\gk[X,Y]/\gen {f}$ et l'on obtient $\gk(x,y) = \gk(t)$ avec $t = (y-y_0) /
(x-x_0)$.
Ici, \llec calculera les expressions de $x$, $y$ en fonction de $t$: $x =
{t^2 - 1 \over t^2 + 1}$, $y = {-2t \over t^2 + 1}$.
\\
\Gmqtz, les \elts de $\gk[x,y]$ sont exactement les fractions rationnelles
définies partout sur la droite projective $\PP^1(\gk)$ (paramétrée par
$t$) sauf peut-être au \gui {point} $t = \pm i$.

\emph{2.} %%%%%%%%%%%%%%%%%%%%%%%%%%%%%%%%%%%%%%%%%%%%%%%%%%%%%%%%%%%%%%%%%%
On a $x = 1-2u$, $y = -2v$, donc $\gk[x,y] = \gk[u,v]$. L'\egt $\gk[x,y]=\gk[u,v]$ n'est pas difficile et est laissée \alecz. Ce qui est plus
difficile, c'est de montrer \hbox{que $\gk[u,v]$} est la \cli de
$\gk[u]$ dans $\gk(t)$. On renvoie à l'exercice~\ref{exoAnneauOuvertP1}.
\\
\Gmqtz, les pôles de $x$ et $y$ sont $t = \pm i$, ce qui
confirme que $x$, $y$ sont entiers sur $\gk[(1+t^2)^{-1}] = \gk[u]$.
Algébriquement, on a $x = 1 - u$, $y^2 = -1-x^2 \in \gk[u]$, et
$x$, $y$ sont bien entiers sur $\gk[u]$.

\emph{3.}  %%%%%%%%%%%%%%%%%%%%%%%%%%%%%%%%%%%%%%%%%%%%%%%%%%%%%%%%%%%%%%%%%%
Si $i^2 = -1$, on a $(x+iy)(x-iy) = 1$.\\
En posant $w = x+iy$,
on a $\gk[x,y] = \gk[w, w^{-1}]$.

\emph{4.} %%%%%%%%%%%%%%%%%%%%%%%%%%%%%%%%%%%%%%%%%%%%%%%%%%%%%%%%%%%%%%%%%%
On applique la méthode standard en un point non singulier d'une courbe plane.
On écrit 
$$f(X,Y) - f(x_0,y_0) = (X-x_0)u(X,Y) + (Y-y_0)v(X,Y)
$$ 
avec
ici $u = X+x_0$, $v = Y+y_0$; la matrice $A = \cmatrix {y-y_0 & x+x_0\cr
x_0-x & y+y_0\cr}$ est donc une \mpn de $(x-x_0, y-y_0)$
avec $1 \in \cD_1(A)$. Explicitons l'appartenance $1 \in \cD_1(A)$:

\snic {
(-y_0)(y-y_0) + x_0(x+x_0) + x_0(x_0-x) + y_0(y+y_0) = 2
.}

%\sni
Ceci conduit à la matrice $B = {1\over 2} \crmatrix {-y_0 & x_0\cr
x_0 & y_0}$; celle-ci vérifie $ABA = A$ et la matrice $P$ cherchée
est $P = \I_2 - AB = \wi {AB}$:
{\footnotesize 
$$AB = {1\over 2} \cmatrix {y-y_0 & x+x_0\cr x_0-x & y+y_0\cr}
\crmatrix {-y_0 & x_0\cr x_0 & y_0} = {1\over 2}
\cmatrix {x_0x - y_0y + 1 & y_0x + x_0y \cr y_0x + x_0y &
          -x_0x + y_0y + 1}.$$}

\vspace{-.5em}
D'où l'expression \gnle de $P$:
$
P = {1 \over 2} \cmatrix {-x_0x + y_0y + 1 & -(y_0x + x_0y) \cr
-(y_0x + x_0y) & x_0x - y_0y + 1}
$,  pour   $x_0=1,\,y_0=0\,:
{1 \over 2} \cmatrix {1-x & -y\cr -y & 1+x}$.
Ainsi, $P$ est un \prr de rang 1, \mpn de $(x-x_0,
y-y_0)$.  Comme $P$ est \smqz, l'\egtz~(\iref{eqpmlm}) de la
proposition \ref{pmlm} a comme conséquence que $(x-x_0)^2 + (y-y_0)^2$ est
un \gtr de $\gen {x-x_0, y-y_0}$ avec $(x-x_0)^2 + (y-y_0)^2 =
-2(x_0x + y_0y - 1)$.

\Gmqtz, $xx_0 + yy_0 - 1 = 0$ est  la tangente au
cercle $x^2 + y^2 = 1$ au \hbox{point $P_0 = (x_0, y_0)$}.
Pour ceux qui connaissent les diviseurs: le diviseur des
zéros-pôles de cette tangente est le diviseur principal
$2P_0 - 2P_{t = \pm i}$, ce qui correspond au fait que
le carré de l'idéal $\gen {x-x_0, y-y_0}$ est principal.

%\sni
Variante I : on traite directement le cas du point $(x, y)=(1,0)$ (voir
question suivante) puis on utilise le fait que le cercle est un groupe pour passer du point $(1,0)$ à un point
quelconque $P_0 = (x_0, y_0)$.  Ainsi, on dispose de l'\auto \gui
{rotation}

\snic {
\cmatrix {x\cr y\cr} \mapsto
\crmatrix {x_0 & -y_0 \cr y_0 & x_0\cr} \cmatrix {x\cr y\cr}
\; \hbox {qui réalise} \;
\crmatrix {x_0 & -y_0 \cr y_0 & x_0\cr} \cmatrix {1\cr 0\cr} =
\cmatrix {x_0\cr y_0\cr}
.}

%\sni
On considère son inverse $R$,

\snic {
R = \crmatrix {x_0 & y_0 \cr -y_0 & x_0\cr}, \;
  R \cmatrix {x\cr y\cr} =\cmatrix {x'\cr y'\cr}, \;
 R \cmatrix {x_0\cr y_0\cr}= \cmatrix {1\cr 0\cr}
,}

%\sni
de sorte que :

\snic {
R \cmatrix {x-x_0\cr y-y_0\cr}=\cmatrix {x' - 1\cr y'\cr},
\quad \hbox{d'où} \quad
\gen {x'-1, y'} = \gen {x-x_0, y-y_0}
.}

%\sni
Comme $\gen {x'-1, y'}^2 = \gen {x'-1}$, on obtient
$\gen {x-x_0, y-y_0}^2 = \gen {x_0x + y_0y - 1}$.

%\sni
Variante II:
on fournit une autre justification de l'inversibilité de $\gen
{x-x_0, y-y_0}$ qui n'utilise pas directement le fait que le cercle est
lisse. On considère $\gk[x,y]$ comme une extension de degré 2 de
$\gk[x]$, en utilisant $(1, y)$ comme base.  On dispose d'un
$\gk[x]$-\auto $\sigma$ qui transforme $y$ en $-y$. \\
On considère la norme $\rN$ de
$\gk[x,y]$ sur $\gk[x]$. Pour $z = a(x) + b(x)y$, on a:

\snic {
\rN(z) = z\sigma(z) = (a + by)(a - by) = a^2 - (1-x^2)b^2 =
a^2 + (x^2-1)b^2
.}

%\sni
L'idée pour inverser $\gen {x-x_0, y-y_0}$ est de le multiplier par son
$\gk[x]$-conjugué. Montrons l'égalité suivante, certificat
de l'inversibilité de l'idéal $\gen {x-x_0, y-y_0}$:

\snic {
\gen { x-x_0, y-y_0}\ \gen {x-x_0, y+y_0} = \gen {x-x_0}
.}

%\sni
En effet, les générateurs du produit de gauche sont~:

\snic {
(x-x_0)^2, \quad (x-x_0)(y+y_0), \quad
(x-x_0)(y-y_0), \quad y^2 - y_0^2 = x_0^2 - x^2
.}

%\sni
D'où $\gen {x-x_0, y-y_0}\ \gen {x-x_0, y+y_0} =
(x-x_0) \gen { g_1, g_2, g_3, g_4}$ avec

\snic {
g_1 = x-x_0, \quad g_2 = y+y_0, \quad g_3 = y-y_0, \quad g_4 = x+x_0
.}

%\sni
Mais $\langle g_1, g_2, g_3, g_4\rangle$ contient ${g_4 - g_1\over 2} = x_0$
et ${g_2 - g_3\over 2} = y_0$ donc il contient $1 = x_0^2 + y_0^2$.

\emph{5.} %%%%%%%%%%%%%%%%%%%%%%%%%%%%%%%%%%%%%%%%%%%%%%%%%%%%%%%%%%
Par la brute force, en utilisant à droite que $1 \in \gen {x-1, x+1}$:

\snuc {
\gen {x-1, y}\,\gen {x-1, y} = \gen {(x-1)^2, (x-1)y, y^2} =
(x-1)\,\gen {x-1, y, -(x+1)} = \gen {x-1}.
}

%\sni
On divise cette \egt par $x-1$: $\gen {x-1,
y}\,\gen {1, {y\over x-1}}=\gen {1}$ et l'on pose~:

\snic {
x_1 = x-1, \quad x_2 = y, \; y_1 = 1, \; y_2 = {y \over x-1},
\; \hbox { de sorte que} \; x_1y_1 + x_2y_2 = -2,
}

%\sni
ce qui conduit à la matrice de projection $P$ de rang $1$:

\snic {
P = \frac{-1}2 \cmatrix {y_1\cr y_2\cr} [x_1, x_2] =
\frac{-1}2 \cmatrix {
x_1 y_1  &  x_2 y_1\cr
x_1 y_2  &  x_2 y_2\cr
} =
\frac{1}2 \cmatrix {
1 -x   &  -y\cr
-y  & 1+x\cr
}}

\emph{6.} %%%%%%%%%%%%%%%%%%%%%%%%%%%%%%%%%%%%%%%%%%%%%%%%%%%%%%%%%%%%%%%%%%%%%%
Notons $\rN = \rN_{\gk[x,y]/\gk}$.  Pour $a, b \in \gk[x]$, $\rN(a + by)
= a^2 + (x^2-1)b^2$.  L'\egt à prouver sur les
degrés est évidente si $a$ ou $b$ est nul. Sinon, on écrit, avec $n =
\deg a$ \hbox{et $m = 1 + \deg b$}, $a^2 =
\alpha^2 x^{2n} + \dots$, $(x^2 - 1)b^2 = \beta^2 x^{2m} + \dots$ ($\alpha,
\beta \in \gk\sta$). Le cas \hbox{où $2n \ne 2m$} est facile. Si $2n = 2m$, alors
$\alpha^2 + \beta^2 \ne 0$ (car $-1$ non carré dans $\gk$), et donc
le \pol $a^2 + (x^2 - 1)b^2$ est de degré $2n = 2m$.
\\
Si $a + by$ est \iv dans $\gA$, $\rN(a+by) \in \gk[x]^{\!\times} =
\gk\sta$; d'où $b = 0$ puis $a$ constant. Bilan:
$\gk[x,y]^{\!\times} = \gk\sta$. Ceci est spécifique au fait que $-1$
n'est pas carré dans $\gk$ car si $i^2 = -1$, l'\egt $(x + iy)(x-iy)
= 1$ montre l'existence d'\ivs autres que les constantes.
\\
Montrons que $y$ est \irdz. \\
Si $y = zz'$, alors $\rN(y) = \rN(z)\rN(z')$, i.e
$x^2 - 1 = (x-1)(x+1) = \rN(z)\rN(z')$.  Mais dans $\gk[x]$, $x \pm 1$ ne sont
pas associés à une norme (une norme non nulle est de degré pair). Donc
$\rN(z)$ ou $\rN(z')$ est une constante, i.e. $z$ ou $z'$ est \ivz.
De la même façon, $1 \pm x$ sont \irdsz.
\\
On va utiliser l'\egt

\snic {
y^2 = (1-x)(1+x), \; \hbox { analogue à } \,
2\cdot 3 = (1 + \sqrt {-5}) (1 - \sqrt {-5}) \hbox { dans }
\ZZ[\sqrt {-5}],
}

%\sni
pour voir que $\gen {x-1, y}$ n'est pas un \idp: une
\egt $\gen {x-1, y} = \gen {d}$ entra\^\i nerait $d \divi x-1$, $d \divi
y$, i.e. $d$ \ivz, i.e. $1 \in \gen {x-1, y}$, ce qui n'est pas.

%%%%%%%%%%%%%%%%%%%%%%%%%%%%%%%%%%%%%%%%%%%%%%%%%%%%%%%%%%%%%%%%%%%%%%%%%%%

%%%%%%%%%%%%%%%%%%%%%%%%%%%%%%%%%%%%%%%%%%%%%%%%%%%%%%%%%%%%%%%%%%%%%%%%%%%

\prob{exoLambdaGammaK0} \emph{(Les opérations $\lambda_t$ et $\gamma_t$ sur $\KO(\gA)$)}\\
\emph {1.} On a $\lambda_t(\gA)=\lambda_t(1)=1+t$ et $\gamma_t(1)=1/(1-t)$. \\
Donc $\lambda_t(p)=(1+t)^p$ et $\gamma_t(p)=1/(1-t)^p$ pour $p\in\NN\etl$.
\\
On écrit $x$ sous la forme $[P]-[\gA^p]=P-p$ 
pour un certain $p\in\NN\etl$, \hbox{avec $P$}  de rang constant $p$.
D'après la \dfn $\gamma_t([P]) = \sum_{n=0}^p \lambda^n(P) t^n/(1-t)^n$,
on a

\snic {
\gamma_t(x) =\frac{\gamma_t([P])}{\gamma_t(p)} =\sum_{n=0}^p \lambda^n(P) t^n(1-t)^{p-n}.
}

%\sni
Ainsi $\gamma_t(x)$ est un \pol de degré $\le p$ en $t$. 
\\
 Note:
$\gamma^p(x) = \sum_{n=0}^p \lambda^n(P) (-1)^{p-n}
= (-1)^p \sum_{n=0}^p \lambda^n(P) (-1)^{n}= (-1)^p \lambda_{-1}(P)$.
\\
On a $\gamma_{t}(x) \gamma_{t}(-x)  = 1$ et comme ce sont des \pols
de $\KO(\gA)[t]$, leurs \coes de degré $> 0$ sont nilpotents
(lemme \ref{lemGaussJoyal} et exercice \ref{exoNilIndexInversiblePol}).
En particulier l'\elt $x$, qui est le \coe de degré 1 de 
$\gamma_{t}(x)$, est nilpotent.

\emph {2.}
Soit $x \in \KO(\gA)$ nilpotent, alors $\rg x$ est un \elt nilpotent
de $\HO(\gA)$. Mais ce dernier anneau est réduit (en fait, \qiz);
donc $\rg x = 0$.

\emph {3.} Supposons $\rg x=[e]$ pour un \idm $e$.\\
On a $\Al n (e\gA)=0$ pour $n\geq2$, donc $\lambda_t([e])=1+[e]t$.
Par \dfn de %la notation exponentielle 
$a^r$ \hbox{pour $a\in\gB$} et $r\in\HO\gB$, on obtient
$(1+t)^{[e]}=(1-e)+e(1+t)=1+et$. \\
Par calcul direct on obtient aussi $\rR {e\gA}(t)=(1-e)+te$.\\
Enfin, on a  par convention
$\BB(\gA)\subseteq\HO\gA$
avec l'identification $e=[e]$.  
\\
On obtient ensuite
l'\egt \gnle pour $x = [P]$ en utilisant le
\sfio formé par les \coes de %du \pol rang 
$\rR{P}$ et en notant que
les deux membres sont des morphismes de $\KO(\gA)$ vers $1 +
t\KO(\gA)[[t]]$. \\
Notons aussi que $\lambda_t(p)=(1+t)^p$ pour $p\in\NN\etl$ est l'\egt voulue lorsque $\rg x\in\NN\etl$.

\emph {4.}
S'obtient à partir du point \emph{1} en remplaçant $t$
par $t/(1-t)$.

\emph {5.}
Un $x \in \KO(\gA)$ s'écrit $y+r$ avec $r=\rg x \in\HO\gA$ et $y\in\KTO\gA$.\\
Alors $\gamma_{t}(x)=\gamma_{t}(y)(1-t)^{-r}$.

\emph {6.}
On rappelle  les deux 
formules suivantes pour $d \ge 1$:

\snic {
{1 \over (1-t)^d} = \sum_{k\ge 0} {k+d-1 \choose d-1} t^k,
\qquad
(1-t)^{-d} = \sum_{k\ge 0} {-d \choose k} (-t)^k.
}

%\sni
Elles sont reliées par l'\egt 

\snic {
{k+d-1 \choose d-1}  = {k+d-1 \choose k}  =
{-d \choose k} (-1)^k.
}

%\sni
Par \dfnz,

\snic {
\gamma_t(x) = 1 + \sum_{d \ge 1} {\lambda^d(x) t^d \over (1-t)^d} =
1 + \sum_{d \ge 1, k \ge 0} \lambda^d(x) t^d 
{k+d-1 \choose d-1} t^k.
}

%\sni
Pour $n \ge 1$, le \coe $\gamma^n(x)$ de $t^n$ est:

\snic {
\sum_{k+d=n} \lambda^d(x) {k+d-1 \choose d-1} 
\quad \hbox {i.e. avec $p = d-1$} \quad
\sum_{p = 0}^{n-1} \lambda^{p+1}(x) {n-1 \choose p}. 
}

%\sni
L'autre \egt s'en déduit via l'\eqvc $\gamma_t =
\lambda_{t/(1-t)} \iff \lambda_t = \gamma_{t/(1+t)}$.

%%%%%%%%%%%%%%%%%%%%%%%%%%%%%%%%%%%%%%%%%%%%%%%%%%%%%%%%%%%%%%%%%%%%%%%%%%%

\prob{exoApplicationProjectiveNoether} \emph{(L'application projective de \Noe et les 
\mrcs1 facteurs directs dans $\gk^2$)}\\
\emph{1.}
Unicité de la factorisation à l'ordre près des facteurs et aux
inversibles près.

\emph{2.}
Le produit de \pols primitifs est un \pol primitif, cf. le lemme
\ref{lemGaussJoyal} (Gauss-Joyal du pauvre). On a le résultat plus
précis qui consiste en l'inclusion d'\idsz:

\snic {
\gen {x_1, y_1} \cdots \gen {x_n, y_n} \subseteq 
\rD_\gk(\gen {z_0, \ldots, z_n}).
}

%\sni

On peut le déduire du fait suivant: si $f$, $g$ sont deux \pols à une
\idtrz, le produit d'un \coe de $f$ et d'un \coe de $g$ est entier sur l'\id
engendré par les \coes du produit $fg$ (voir le lemme~\ref{lemthKroicl}), 
et en
particulier il est dans le radical de cet \idz.  
\\
On peut \egmt utiliser
l'approche qui suit: pour $I \subseteq \lrbn$, notons $I'$ son
complémentaire, $x_I = \prod_{i
\in I} x_i$, $y_I = \prod_{i \in I} y_i$. Pour $d = \#I$ et  $N = {n \choose d}$, on va montrer une \egtz:
$$
\prod\nolimits_{\#I = d}(T - x_Iy_{I'}) = T^N +
\som_{j=1}^{N} a_j T^{N-j}, \quad a_j \in \gen {z_0, \ldots, z_n}.
\leqno (\star')
$$
En faisant $T=x_Iy_{I'}$, on aura $(x_Iy_{I'})^N \in \gen {z_0, \ldots, z_n}$,
montrant ainsi l'inclusion d'\ids annoncée. Pour prouver $(\star')$, on
examine d'abord le cas où tous les $y_i$ sont égaux à $1$. On
écrit, en notant $S_1(x), \ldots, S_n(x)$ les fonctions
\smqs\elrs de $(\xn)$:

\snic {
\prod_{\#I = d}(T - x_I) = T^N + \sum_{j=1}^{N} b_j T^{N-j}, \qquad
b_j = f_j\big(S_1(x), \ldots, S_n(x)\big).
}

%\sni
Un examen attentif montre que $f_j$ est un \pol de degré $\le j$ en
$(S_1, \ldots, S_n)$.  
Remplaçons dans cette dernière \egt $x_i$ par
$x_i/y_i$ et multiplions par $(y_1\cdots y_n)^N$; on obtient, avec $U =
y_1\cdots y_nT$ et $s_i = S_i(x_1/y_1, \ldots, x_n/y_n)$:

\snic {
\prod_{\#I = d}(U - x_Iy_{I'}) = U^N + 
\sum_{j=1}^{N} (y_1\cdots y_n)^j f_j(s_1, \ldots, s_n) U^{N-j}.
}

%\sni

Soit $s_1^{\alpha_1} \cdots s_n^{\alpha_n}$ un monôme de $f_j(s_1, \ldots,
s_n)$; puisque $\sum_i \alpha_i \le \deg f_j \le j$, on obtient, en se souvenant
que $z_n = y_1\cdots y_n$, une \egtz:

\snic {
z_n^j s_1^{\alpha_1} \cdots s_n^{\alpha_n} =
z_n^{\alpha_0} (z_ns_1)^{\alpha_1} \cdots (z_ns_n)^{\alpha_n} =
z_n^{\alpha_0} z_{n-1}^{\alpha_1} \cdots z_0^{\alpha_n}
\hbox { avec } \alpha_0 = j - \sum_i \alpha_i.
}

%\sni
Puisque $j \ge 1$, l'un des exposants $\alpha_i$ ci-dessus n'est pas nul et
l'on a bien l'appartenance à $\gen {z_0, \ldots, z_n}$,
puis l'\egt $(\star')$.

\emph{3.}
Posons $E = P_1 \te_\gk \cdots \te_\gk P_n \subset L^{n\te}$; c'est un 
\mrcz~1. Montrons que la restriction de $\pi$ à $E$ est injective et que
$\pi(E)$ est facteur direct dans $S_n(L)$. Ceci prouvera bien que $\pi(E)$ est
un $\gk$-point de $\PP^n$.  On se ramène à l'aide d'un nombre fini de
\lons \come au cas où chaque $P_i$ est libre de base $x_iX +
y_iY$. Alors chaque $(x_i, y_i)$ est \umd et $\sum_{i=0}^n z_i
X^{n-i}Y^i$ est une base \umd de $\pi(E)$. Ceci prouve d'une part que
$\pi\frt E$ est injective (puisqu'elle transforme une base de $E$ en un
vecteur \umd de $S_n(L)$) et que $\pi(E)$ est facteur direct dans
$S_n(L)$.

\emph{4.}
Il semble que $\varphi$ soit injectif, i.e. que $(z_0, \ldots, z_n)$ sont
\agqt indépendants sur $\gk$. L'image par $\varphi$ est le sous-anneau
gradué $\gA = \gk[z_0, \ldots, z_n] \subset \gk[\uX,\uY]$ (la composante
\hmg d'un \elt de $\gA$ est dans $\gA$); si $f \in \gA$ est \hmg de degré
$m$, on a $m \equiv 0 \bmod n$, et pour $t_1, \ldots, t_n$ quelconques:

\snic {
f(t_1X_1, t_1Y_1, \ldots, t_nX_n, t_nY_n) = (t_1\ldots t_n)^{m/n}
f(X_1,Y_1, \ldots, X_n,Y_n)
.}

%\sni
Enfin, $\gA$ est invariant sous l'action du groupe symétrique $\rS_n$
qui agit sur $\gk[\uX, \uY]$ par

\snic {
\sigma \cdot f(X_1,Y_1, \ldots, X_n,Y_n) =
f(X_{\sigma(1)},Y_{\sigma(1)}, \ldots, X_{\sigma(n)},Y_{\sigma(n)})
.}

%\sni
Ces deux dernières propriétés caractérisent probablement $\gA$.
%%%%%%%%%%%%%%%%%%%%%%%%%%%%%%%%%%%%%%%%%%%%%%%%%%%%%%%%%%%%%%%%%%%%%%%%%%%

%%%%%%%%%%%%%%%%%%%%%%%%%%%%%%%%%%%%%%%%%%%%%%%%%%%%%%%%%%%%%%%%%%%%%%%%%%%

\prob{exoTh90HilbertMultiplicatif} \emph{(Le \tho 90 multiplicatif d'Hilbert)}\\
On fixe une fois pour toutes un \elt $b_0 \in \gB$ de trace 1.

\emph {1} et \emph {2.}
Pas de difficulté. Le fait que $\theta_c$ soit multiplicatif traduit
exactement le fait que $c$ est un 1-cocycle.

\emph {3.}
L'action de $G$ sur $\gB$ tordue par le 1-cocycle $c$ est $\sigma\cdot_c b =
c_\sigma\sigma(b)$; le fait que cela soit une action est exactement la condition
de $1$-cocyclicité de $c$. En effet:

\snic {
\tau\cdot_c(\sigma\cdot_c b) = \tau\cdot_c c_\sigma \sigma(b) = 
c_\tau\tau\big(c_\sigma \sigma(b)\big) = 
c_\tau\tau(c_\sigma)\, (\tau\sigma)(b) = 
c_{\tau\sigma}\, (\tau\sigma)(b) = (\tau\sigma) \cdot_c b
.}

%\sni
On remarquera que $\pi_c = \sum_\sigma c_\sigma\, \sigma$ est une sorte
de $G$-trace relativement à l'action de $G$ tordue par $c$.

On a donc $\gB_c^G = \sotq {b \in \gB} {c_\sigma \sigma(b) = b}$.  En
utilisant le fait que $c$ est un 1-cocycle, on trouve que $\tau \circ \pi_c =
c_\tau^{-1}\, \pi_c$; on en déduit que $c_\tau \tau(z) = z$ pour tout $z \in
\Im\pi_c$, i.e. $\Im\pi_c \subseteq \gB_c^G$. On définit $s : \gB_c^G
\to \gB$ par $s(b) = bb_0$. Alors $\pi_c \circ s = \Id_{\gB_c^G}$; en
effet, pour $b \in \gB_c^G$:

\snic {
\pi_c(b_0b) = \sum_\sigma c_\sigma\sigma(bb_0) =
\sum_\sigma c_\sigma\sigma(b)\sigma(b_0) = \sum_\sigma b\sigma(b_0) 
= b\Tr_{\gB\sur\gA}(b_0) = b
.}

%\sni
De l'\egt $\pi_c \circ s = \Id_{\gB_c^G}$, on déduit que $\pi_c$
est une surjection de $\gB$ sur $\gB_c^G$, que~$s$ est injectif
et que $\gB = s(\gB_c^G) \oplus \Ker\pi_c \simeq \gB_c^G \oplus \Ker\pi_c$.
En particulier, $\gB_c^G$ est \hbox{un $\gA$-\mptfz}.

Remarque. Voyons $s : b \mapsto b_0b$ dans
$\End_\gA(\gB)$, alors $(\pi_c \circ s)\big(\pi_c(z)\big) = \pi_c(z)$
pour tout $z \in \gB$, i.e. $\pi_c \circ s \circ \pi_c = \pi_c$.
En conséquence $\pi'_c \eqdf {\rm def} \pi_c \circ s =
\sum_\sigma c_\sigma \sigma(b_0\bullet)$ est un \prr;
on pourrait certainement calculer sa trace et trouver 1, ce qui
prouverait que $\pi'_c$ un \prr de rang $1$.

\emph {4.}
Soient $c$, $d$ deux 1-cocycles, $x \in \gB_c^G$, $y \in \gB_d^G$, donc
$c_\sigma\sigma(x) = x$, $d_\sigma\sigma(y) = y$; on vérifie facilement que
$xy \in \gB_{cd}^G$. \\
D'où une \Ali $\gB_{c}^G \otimes_\gA \gB_{d}^G \to
\gB_{cd}^G$, $x\otimes y \mapsto xy$, notée~$\mu_{c,d}$.
\\
Notons $(x_i)$, $(y_i)$ deux \syss d'\elts de $\gB$ comme dans le lemme
\ref{lemArtin} et posons $\varepsilon = \sum_i x_i \otimes y_i = 
\sum_i y_i \otimes x_i$ (\idm de séparabilité). On rappelle que
$\varepsilon \in \Ann(\rJ)$, ce qui se traduit par

\snic {
\forall\ b \in \gB\quad
\sum_i bx_i \otimes y_i = \sum x_i \otimes by_i  \quad\hbox{dans}\quad
\env\gA\gB \eqdf {\rm def} \gB\otimes_\gA \gB
.}

%\sni
On a aussi, pour $b$, $b' \in \gB$

\snic {
\Tr_{\gB\sur\gA}(bb') =\sum_i \Tr_{\gB\sur\gA}(bx_i)\Tr_{\gB\sur\gA}(b'y_i)
.}

%\sni
On va montrer que $z \mapsto (\pi_c \otimes \pi_d)(b_0z\varepsilon)$,
$\gB_{cd}^G \mapsto \gB_{c}^G \otimes_\gA \gB_{d}^G$ et $\mu_{c,d}$ sont
réciproques l'une de l'autre. Dans un sens:

\snic {
(\pi_c \otimes \pi_d)(b_0z\varepsilon) = \sum_i a_i \otimes b_i, 
\quad \hbox {avec} \quad
a_i = \sum_\sigma c_\sigma\sigma(b_0zx_i), \quad
b_i = \sum_\tau c_\tau\tau(y_i)
,}

%\sni
et l'on a

\snic {
\sum_i a_ib_i = \sum_{\sigma,\tau} \sigma(b_0z) c_\sigma d_\tau
\sum_i \sigma(x_i) \tau(y_i)
,}

%\sni
et comme la somme interne (sur $i$) vaut $1$ ou $0$, il reste,
pour $z \in \gB_{cd}^G$: 

\snic {
\sum\limits_i a_ib_i = \sum\limits_{\sigma} \sigma(b_0z) c_\sigma d_\sigma =
\sum\limits_{\sigma} \sigma(b_0) \sigma(z) (cd)_\sigma =
\sum\limits_{\sigma} \sigma(b_0) z = z \Tr_{\gB\sur\gA}(b_0) = z
.}

%\sni
Dans l'autre sens, soient $x \in \gB_c^G$ et $y \in \gB_d^G$. Alors,
puisque $\varepsilon \in \Ann(\rJ)$, on peut écrire:

\snic {
(\pi_c \otimes \pi_d)(b_0 xy\varepsilon) = \sum_i a_i \otimes b_i, 
\; \hbox {avec} \;
a_i = \sum_\sigma c_\sigma\sigma(b_0xx_i), \;
b_i = \sum_\tau d_\tau\tau(yy_i)
.}

%\sni
En utilisant 
$$c_\sigma \sigma(b_0xx_i) = c_\sigma \sigma(x) \sigma(b_0x_i)
= x\sigma(b_0x_i) \hbox{ et } d_\tau \tau(yy_i) = d_\tau
\tau(y)\tau(y_i) = y\tau(y_i),
$$
il vient
\[ 
\begin{array}{ccc} 
  \sum_i a_i \otimes b_i = 
\sum_i x\Tr_{\gB\sur\gA}(b_0x_i) \otimes y\Tr_{\gB\sur\gA}(y_i)  =   \\[1mm] 
(x \otimes y) \cdot
\big(\sum_i \Tr_{\gB\sur\gA}(b_0x_i)\Tr_{\gB\sur\gA}(y_i) \otimes 1\big) =
(x \otimes y) \cdot \big(\Tr_{\gB\sur\gA}(b_0) \otimes 1\big) = x\otimes y
. 
\end{array}
\]
Le point \emph {a} est prouvé. \\
Pour le point \emph {b}, soit un
1-cocycle, cobord de $b_1 \in \Bti$, $c_\sigma = \sigma(b_1)b_1^{-1}$. \\
Alors $b
\in \gB_c^G$ \ssi pour tout $\sigma$, $c_\sigma\sigma(b) = b$,
i.e. $\sigma(b_1b) = b_1b$ \cad $b_1b \in \gA$; donc $\gB_c^G = b_1^{-1}\gA$.
On en déduit que $\gB_c^G \otimes \gB_{c^{-1}}^G \simeq \gA$, 
\hbox{donc $\gB_c^G$}
est un \Amrc 1. \\
De plus $c \mapsto \gB_c^G$ induit un morphisme $\zcoho \to
\Pic(\gA)$.

Il reste à montrer que si $\gB_c^G$ est libre, i.e. $\gB_c^G = \gA b_1$
avec $b_1 \in \gB$ et $\Ann_\gA(b_1) = 0$, alors $c$ est un cobord.
Mais $\gB_{c^{-1}}^G$, étant l'inverse de $\gB_c^G$ est aussi
libre, $\gB_{c^{-1}}^G = \gA b_2$, et 
$\gB_c^G \gB_{c^{-1}}^G = \gB_1^G = \gA$.
On a donc
$\gA b_1b_2 = \gA$, puis $b_1$, $b_2$ sont \ivs dans~$\gB$ (et
$\gA b_2 = \gA b_1^{-1}$). Alors $c_\sigma^{-1} \sigma(b_2) =
b_2$, i.e. $c$ est le cobord de $b_2$.

\emph {5.}
Puisque $\gA$ est un anneau \zedz, $\Pic(\gA) = 0$ donc $\hcoho = 0$.

\emph {6.}
On pose $c_\tau = x\sigma(x) \cdots\sigma^{i-1}(x)$ avec $i \in \lrb{1..n}$ et
$\tau = \sigma^i$.\\
Ainsi, $c_\Id = \rN_{\gB\sur\gA}(x) = 1$, $c_\sigma = x$,
$c_{\sigma^2} = x\sigma(x)$. \\
C'est un 1-cocycle: $c_\sigma
\sigma(c_{\sigma^i}) = c_{\sigma^{i+1}}$, i.e. $c_\sigma\sigma(c_\tau) =
c_{\sigma\tau}$, puis $c_{\sigma^j}\sigma^j(c_\tau) = c_{\sigma^j \tau}$.

%%%%%%%%%%%%%%%%%%%%%%%%%%%%%%%%%%%%%%%%%%%%%%%%%%%%%%%%%%%%%%%%%%%%%%%%%%%
\prob{exoSegreMorphism} \emph{(Le morphisme de Segre dans un cas particulier)}\\
Il est clair que $\fa \subseteq \Ker\varphi$.

\emph{1.}
Soient $m = X_{i_1} \cdots X_{i_r} Y_{j_1} \cdots Y_{j_s}$, $m' = X_{i'_1}
\cdots X_{i'_{r'}} Y_{j'_1} \cdots Y_{j'_{s'}}$ avec
$$
1 \le i_1 \le \cdots \le i_r \le j_1 \le \cdots \le j_s \le n, \;
1 \le i'_1 \le \cdots \le i'_{r'} \le j'_1 \le \cdots \le j'_{s'} \le n.
$$
L'\egt $\varphi(m) = \varphi(m')$ fournit

\snic{
T^r U^s Z_{i_1} \ldots Z_{i_r} Z_{j_1} \ldots Z_{j_s} =
T^{r'} U^{s'} Z_{{i'}_{\!1}} \ldots Z_{{i'}_{\!r'}} Z_{j'_1} \ldots Z_{{j'}_{\!s'}}
.}

%\sni
Donc $r = r'$, $s = s'$ puis $i_k = i'_k$ et $j_\ell = j'_\ell$.
En définitive $m = m'$.
\\
Soient $s = \sum_\alpha a_\alpha m_\alpha$ une combinaison $\gA$-\lin
de \moms normalisés telle que $\varphi(s) = 0$. Comme les
\moms $\varphi(m_\alpha)$ sont deux à deux distincts, on
a $a_\alpha = 0$, i.e.~$s = 0$.

\emph{2.}
Puisque $X_iY_j \equiv X_jY_i \mod {\fa}$, on voit que tout
\mom est équivalent modulo~$\fa$ à un \mom
normalisé. Il vient donc $\gA[\uX,\uY] = \fa + \fa_{\rm nor}$.
Comme $\fa \subseteq \Ker\varphi$, la somme est directe d'après
la question précédente.

\emph{3.}
Soit $h \in \Ker\varphi$ que l'on décompose en $h = f + g$ avec $f \in \fa$,
$g \in \fa_{\rm nor}$.\\
Puisque $\fa \subseteq \Ker\varphi$, on a $g \in
\Ker\varphi$, donc $g = 0$.  Conclusion: $h = f \in \fa$, ce qui prouve
$\Ker\varphi \subseteq \fa$, puis $\Ker\varphi =\fa$.
%%%%%%%%%%%%%%%%%%%%%%%%%%%%%%%%%%%%%%%%%%%%%%%%%%%%%%%%%%%%%%%%%%%%%%%%%%%

%%%%%%%%%%%%%%%%%%%%%%%%%%%%%%%%%%%%%%%%%%%%%%%%%%%%%%%%%%%%%%%%%%%%%%%%%%%
\prob{exoVeroneseMorphism} \emph{(Le morphisme de Veronese dans un cas particulier)}\\
Il est clair que $\fb \subseteq \fa \subseteq \Ker\varphi$.

\emph{1.}
Soit $f$ dans l'intersection; $f$ s'écrit $f = f_0 + \sum_{i=1}^{d-1}
f_i X_i$ avec $f_i \in \gA[X_0,X_d]$; on écrit que $\varphi(f) = 0$:

\snic {
f_0(U^d,V^d) + f_1(U^d,V^d)U^{d-1}V + \cdots + f_{d-1}(U^d,V^d)UV^{d-1} = 0.
}

%\sni
Ceci est de la forme, dans $\gA[U][V]$, $h_0(V^d) + h_1(V^d)V + \cdots +
h_{d-1}(V^d)V^{d-1} = 0$; en examinant dans cette \egt les exposants de $V$
modulo $d$, on obtient $h_0 = h_1 = \cdots = h_{d-1} = 0$. Bilan: $f_i = 0$
puis $f = 0$.

\emph{2.}
On travaille modulo $\fb$ en posant 

\snic{\Aux = \AuX/\fb$, $\gB = \gA[x_0,x_d]
+ \gA[x_0,x_d]x_1 + \cdots + \gA[x_0,x_d]x_{d-1} \subseteq \Aux.}

%\sni
On va montrer
que $\gB$ est une sous-\Algz; comme elle contient les $x_i$, \hbox{c'est $\Aux$} tout entier. Il suffit de prouver que $x_ix_j \in \gB$
 pour $i \le j \in\lrb{1.. d-1}$, car les autres produits
sont dans $\gB$ par \dfn de $\gB$. On utilise les relations
$x_ix_j  = x_{i-1}x_{j+1}$ pour $i \le j \in\lrb{1.. d-1}$. On a $x_0x_k
\in \gB$ pour tout $k$; on en déduit $x_1x_j \in \gB$ pour tout
$ j \in\lrb{1.. d-1}$ et c'est encore vrai pour $j = d$ et $0$ par \dfn
de $\gB$. On en déduit ensuite $x_2x_j \in \gB$ pour $j \in\lrb{2.. d-1}$, et
ainsi de suite.
\\
L'\egt obtenue $\gB = \Aux$ s'écrit 

\snic{\AuX = \fb + \big(\gA[X_0,X_d]
\oplus \gA[X_0,X_d]X_1 \oplus \cdots \oplus \gA[X_0,X_d]X_{d-1}\big),}

%\sni
et le $+$
représente une somme directe d'après le point \emph{1} (puisque $\fb
\subseteq \Ker\varphi$).

\emph{3.}
Soit $h \in \Ker\varphi$ que l'on décompose en $h = f + g$ comme ci-dessus. \\
Puisque
$f\in\fb \subseteq \Ker\varphi$, on a $g \in \Ker\varphi$, donc $g = 0$.
Conclusion: $h = f \in \fb$, ce qui prouve $\Ker\varphi \subseteq \fb$, puis
$\Ker\varphi = \fb = \fa$.

%%%%%%%%%%%%%%%%%%%%%%%%%%%%%%%%%%%%%%%%%%%%%%%%%%%%%%%%%%%%%%%%%%%%%%%%%%%

\prob{exoVeroneseMatrix} \emph{(Matrices de Veronese)}\\
\emph {2.} 
Il est clair que $V_d(P)$ est un \prr si $P$ en est un, et le diagramme est
commutatif pour des raisons fonctorielles.  \\
On peut apporter la précision
suivante: si $P$, $Q \in \Mn(\gk)$ sont deux \prrs tels que $\Im P \subseteq \Im
Q$, alors $\Im V_d(P) \subseteq \Im V_d(Q)$. En effet, \hbox{on a $\Im P \subseteq
\Im Q$} \linebreak 
\ssi $QP = P$, et l'on en déduit que $V_d(Q)V_d(P) = V_d(P)$, \linebreak
i.e.  $\Im V_d(P) \subseteq \Im V_d(Q)$.

\emph {3.} 
Il suffit de le faire localement, i.e. de calculer $V_d(A)$ lorsque $A$ est un
\prr standard $\I_{r,n}$. Si $A=\Diag(a_1, \ldots, a_n)$, alors $V_d(A)$ est diagonale, de diagonale les $n'$
monômes $a^{\alpha} $ avec $\abs{\alpha}= d$. En particulier, pour $A = \I_{r,n}$, on voit que $V_d(A)$ est une projection standard, de rang 
le nombre de $\alpha$
tels que $\alpha_1 + \cdots + \alpha_r = d$, \cad ${d+1-r \choose
r-1}$.  Et $V_d(\I_{1,n}) = \I_{1,n'}$.
%%%%%%%%%%%%%%%%%%%%%%%%%%%%%%%%%%%%%%%%%%%%%%%%%%%%%%%%%%%%%%%%%%%%%%%%%%%

%:sinotenglish
\sinotenglish{

\prob{exoFossumKumarNori} \emph{(Quelques exemples de résolutions projectives finies)}\\
\emph {1.}
Le calcul de $F_k^2 - F_k$ se fait par \recu et ne pose pas de \pbz.
Pour la conjugaison ($n \ge 1$), on utilise

\snic {
\cmatrix {0& -\I\cr \I & 0} \cmatrix {A & B\cr C& D} \cmatrix {0& \I\cr -\I& 0} =
\cmatrix {D & -C\cr -B & A} 
.}

%\sni
Pour $\cmatrix {A & B\cr C& D} = F_n$, cela fournit une conjugaison
entre $F_n$ et $\I_{2^n} - \tra{F_n}$.\\
Lorsque $z(z-1) + \sum_{i=1}^n x_iy_i = 0$, les \prrs
 $F_n$ et $\I_{2^n} - {F_n}$ ont pour image des \mptfs $P$ et $Q$ avec $P\oplus Q\simeq \Ae{2^n}$ et $P\simeq Q\sta$. \\
 Donc $2\rg(P)=2^{n}$, et comme $mx=0\Rightarrow x=0$ pour $m\in\NN\etl$ et $x\in\HO \gA$, on obtient $\rg(P)=2^{n-1}$.

\emph {2.}
Le calcul de $U_kV_k$ et $V_kU_k$ se fait par \recuz.  Le fait que $F_n$ et
$G_n$ sont conjuguées par une matrice de permutation est laissé à la
sagacité \dlecz. \\
Par exemple, $G_2 = P_\tau F_2 P_\tau^{-1}$ pour $\tau
= (2,4,3) =(3,4)(2,3)$, et $G_3 = P_\tau F_3 P_\tau^{-1}$ \hbox{pour $\tau = (2, 4, 7, 5)(3,
6)=(3,6)(2,4)(4,7)(5,7)$}.  \\
En ce qui concerne le rang constant $2^{n-1}$ on peut invoquer le point \emph{1}, ou faire un calcul direct après \lon en $z$ et en $\ov z = 1-z$. 
%Examinons $G_n = \cmatrix {z\I
%& U\cr V & \ov z\I}$ avec $U = U_n$, $V = V_n$; sur le localisé $\gA_z$,
%$G_n$ possède (dans le coin nord-ouest) le mineur $z\I_{2^{n-1}}$ qui est
%\iv et on peut appliquer le lemme du mineur inversible (par exemple l'exercice
%\ref {exoTroisiemeLemmeLiberte}).  Ceci conduit à introduire sur $\gA$

%\snic {
%L = \cmatrix {\I & U\cr 0 &\I} \cmatrix {z\I & 0\cr -V & z\I}, \qquad
%R = \cmatrix {z\I & 0\cr V &z\I} \cmatrix {\I & -U\cr 0 & \I}
%}

%\sni On a alors les \egts (pour vérifier celle de droite, il faut
%utiliser $UV = VU = z\ov z \I$):

%\snic {
%LR = z^2\I, \qquad L G_n R = z^2 \cmatrix {\I & 0\cr 0 & 0}
%}

%\sni
%Ceci montre que sur le localisé $\gA_z$, la matrice $G_n$ est conjuguée
%à un \prr standard de rang $2^{n-1}$. Idem sur le localisé en $\ov z = 1-z$.

\emph {3a.}
Utilisation directe de l'exercice référencé.

\emph {3b.}
Soit $S$ le \mo $a^\NN$. On peut localiser une résolution projective finie
de $M$ sur $\gA$ pour en obtenir une sur $S^{-1}\gA$: 

\snic{0 \rightarrow S^{-1}P_n
\to \cdots \to S^{-1}P_1 \to S^{-1}P_0 \twoheadrightarrow S^{-1}M \to
0.}

%\sni
Comme $aM = 0$, on a $S^{-1} M = 0$, donc $\sum_{i=0}^n (-1)^i \rg
(S^{-1}P_i) = 0$. Mais le morphisme naturel $\HO(\gA) \to
\HO(S^{-1}\gA)$ est injectif. Donc $\sum_{i=0}^n (-1)^i \rg P_i = 0$.

\emph {4.}
Le localisé $(\gB_n)_z$ contient les  $y'_i= y_i/z$, et
puisque $z(1-z) = \sum_i x_iy_i$, \hbox{on a $1-z = \sum_i x_iy'_i$}.
%:HHH c'est $1 - \sum_i x_iy'_i\in(\gB_n)_z\eti$ et non pas  1 - \sum_i x'_iy'_i
Donc $z \in \gk[\xn, y'_1, \ldots, y'_n]$ et  $1 - \sum_i x_iy'_i\in(\gB_n)_z\eti$. On vérifie alors que 

\snic {
(\gB_n)_z =  \gk[\xn, y'_1, \ldots, y'_n]_{s}
\quad \hbox {avec} \quad s = 1 - \sum x_iy'_i
.}

%\sni
De même, $(\gB_n)_{1-z} = \gk[x'_1, \ldots, x'_n, \yn]_{1 - \sum x'_iy_i}$
avec $x'_i = x_i/(1-z)$.

\emph {5.}
Pour $n \ge 1$, tout \elt $a \in \{z, \xn\}$ est \ndz et $a(\gB_n/\fb_n)
= 0$.  Comme $F_1=\cmatrix {z & x_1\cr y_1 & \ov z}$ est un \prrz, on a $[z,x_1]
F_1 = [z,x_1]$. \Llec vérifiera que
$\Ker\,[z,x_1] = \Ker F_1 = \Im (\I_2 - F_1)$; d'où la suite exacte:

\snic {
0 \rightarrow \Im(\I_2-F_1) \to \gB_1^2
\vvers{[z,x_1]} \gB_1 \twoheadrightarrow \gB_1\sur{\fb_1} \to 0
.}

%\sni
On a bien $\rg (\gB_1\sur{\fb_1}) = 1-2+1 = 0$.

\emph {6.}
Soit $A$ la matrice constituée des 3 premières lignes de $\I_4 - F_2$:

\snic {
A = \cmatrix {1-z & -x_1 & -x_2 & 0\cr
-y_1 & z & 0 & -x_2\cr -y_2 & 0 & z & x_1\cr}
.}

%\sni
Il est clair que $A F_2 = 0$ et $[z,x_1,x_2]A = 0$. \Llec vérifiera que
la suite ci-dessous est exacte:

\snic {
0 \rightarrow \Im F_2 \to \gB_2^4 \vvers {A} \gB_2^3 \vvvvvers{[z,x_1,x_2]} 
\gB_2 \twoheadrightarrow \gB_2\sur{\fb_2} \to 0
.}

%\sni
On a bien $\rg (\gB_2\sur{\fb_2}) = 1-3+4-2 = 0$.

\emph {7.}
Immédiat vu la \dfn de $F_n$.

\emph {8.}
On considère la moitié haute de la matrice $\I_8 - F'_3$ et
on supprime sa dernière colonne (nulle) pour obtenir
une matrice $A$ de format $4 \times 7$. Soit $B$ la matrice
de format $7 \times 8$ obtenue en supprimant la dernière
ligne de $F'_3$.  Alors \llec courag\eUx vérifiera
l'exactitude de:

\snic {
0 \rightarrow \Im (\I_8-F'_3) \to \gB_3^8 \vvers {B} \gB_3^7 \vvers {A}
\gB_3^4 \vvvvvers{[z,x_1,x_2,x_3]}  \gB_3 \twoheadrightarrow \gB_3\sur{\fb_3} \to 0
.}

%\sni
On a $\rg(\gB_3\sur{\fb_3}) = 1-4+7-8+4 = 0$.

\emph {9.}
Il y a une suite exacte (on a posé $\gB = \gB_n$, $\fb = \fb_n$):

\snic {
L_{n+1} \vvers {A_{n+1}} L_n \vvers {A_n} L_{n-1} \vvers {A_{n-1}} \cdots
\lora L_2 \vvers {A_2} L_1 \vvers{A_1} L_0=\gB 
\twoheadrightarrow \gB\sur{\fb}
.}

%\sni
où $L_r$ est un module libre de rang $\sum_{i \in I_r} {n+1 \choose i}$ avec
$I_r = \sotq {i \in \lrb{0..r}}{i \equiv r \bmod 2}$.  En particulier, $L_1 =
\gB^{n+1}$ et $L_n = L_{n+1} = \gB^{2^n}$.  \\
Quant aux matrices $A_r$, on a
$A_1 = [z,\xn]$, et la matrice $A_r$ est extraite de~$F_n$ si $r$ est impair,
et extraite de $\I - F_n$ sinon. On a $A_{n+1} =
F_n$ pour $n$ pair, et $A_{n+1} = \I - F_n$ pour $n$ impair.
\\
En notant $P_{n+1} = \Im A_{n+1}$, le \Bmo $\gB\sur{\fb}$ admet une résolution projective
de longueur $n+1$ de type suivant:

\snic {
0 \rightarrow P_{n+1} \to L_n = \gB^{2^n} \vers {A_n} L_{n-1} \vvers {A_{n-1}} \cdots
\to L_2 \vers {A_2} L_1 \vers{A_1} L_0=\gB 
\twoheadrightarrow \gB\sur{\fb}
.}

%\sni
($P_{n+1}$ de rang constant $2^{n-1}$).\\
L'expression explicite du rang de $L_i$ confirme que $[\gB\sur{\fb}] \in
\KTO(\gB)$.  \\
On a $\rg L_{n-1} + \rg L_{0} = \rg L_{n-2} + \rg L_1 =
\cdots = 2^n$ (en particulier, si $n = 2m+1$, alors $\rg L_m =
2^{n-1}$).

 Note:
Si $\gk$ est un \cdiz, on peut montrer que $\KTO(\gB_n) \simeq \ZZ$
avec comme \gtr $[\gB_n\sur{\fb_n}]$.  On en déduit que l'\id $\KTO(\gB_n)$ est
de carré nul: de manière \gnlez, soit un anneau $\gA$ vérifiant
$\KTO(\gA)=\ZZ x \simeq \ZZ$, alors $x^2 = mx$ avec $m \in \ZZ$, donc
$x^{k+1} = m^k x$ pour $k \ge 1$, comme $x$ est nilpotent (voir le \pbz~\ref{exoLambdaGammaK0}), il y a un
$k\ge 1$ tel que $m^k x = 0$, donc $m^k = 0$, puis $m=0$ et $x^2 = 0$.

\prob{exoPolynomialSyzygies} \emph{(Quand les \moms dominants sont 
premiers entre eux)}
\\
\emph {1a.}
Soient $f_i = m_i - r_i$ avec $r_i \prec m_i$; alors

\snic{m_jf_i - m_if_j = (f_j+r_j)f_i - (f_i+r_i)f_j = r_jf_i - r_if_j.}

\snii\emph {1b.}
Par hypothèse $m'_im_i = m'_jm_j$; comme $\pgcd(m_i, m_j) = 1$, on a $m_i
\mid m'_j$.  On dispose donc d'un \mom $q$ défini par $q = m'_j/m_i =
m'_i/m_j$.  On écrit alors

\snic {
m'_if_i - m'_jf_j = q(m_jf_i - m_if_j) = q(r_jf_i - r_if_j) = qr_jf_i - qr_if_j,
}

\snii
et l'on a $qr_jf_i \prec qm_jm_i = m'_im_i = m$. De même $qr_if_j \prec m$.
\\
Bilan: $m'_if_i - m'_jf_j \in \fa_{\prec m}$.

\snii\emph {1c.}
On peut supposer $I = \{1,\ldots,k\}$; on écrit la transformation
d'Abel:

\snic {
\sum_{i \in I} a_im'_if_i = \sum_{j=1}^{k-1} 
s_j(m'_jf'_j - m'_{j+1}f'_{j+1}).
}

\snii
D'après la question précédente $m'_jf'_j - m'_{j+1}f'_{j+1}$ appartient
à $\fa_{\prec m}$, et il en est de même de leur somme $\sum_{i \in I}
a_im'_if_i$.

\snii\emph {1d.}
On écrit $g = qf_i \preceq m$; $q$ est donc une combinaison $\gk$-\lin de \moms
$m'$ tels que $m'm_i \preceq m$. Si $m_i \nedivi m$, on a $m'm_i \prec m$ donc
$g \prec m$. Si $m_i \divi m$, on écrit $m = m'm_i$; si $a \in \gk$ est le
\coe de $m'$ dans $q$, on a $q - am' \prec m'$.  
\\
Donc $g = qf_i = am'f_i +
(q - am')f_i$ avec $(q - am')f_i \prec m'm_i = m$.

%%%%%%%%%%
\sni\emph {2.}
Soit $f \in \fa_{\preceq m}$. On écrit $f = \sum_i g_i$ avec
$g_i \in \gen {f_i}$ et $g_i \preceq m$. On coupe cette somme
en deux:

\snic {
f = \sum_{j \notin I} g_j + \sum_{i \in I} g_i 
\quad \hbox {avec} \quad
I = \sotq {i \in \{1,\ldots,s\}} {m_i \hbox { divise } m}
.}

\snii
Si $m_j \nedivi m$, on a $g_j \prec m$ donc $\sum_{j \notin I} g_j \in
\fa_{\prec m}$. Pour les $i$ tels que $m_i \divi m$, on écrit $m = m'_im_i$ et
comme dans la question \emph {1d}, il y a un $a_i \in \gk$ tel \hbox{que $g_i -a_im'_if_i \in \fa_ {\prec m}$}.  Bilan:

\snic {
f = \hbox {un \elt de $\fa_{\prec m}$} + \sum_{i \in I} a_im'_i f_i. 
}

\snii
On utilise maintenant le fait que $f \prec m$. Ceci implique 
$\sum_{i \in I} a_i = 0$. D'après la question \emph {1c} on
a~$\sum_{i \in I} a_im'_i f_i \in \fa_{\prec m}$ donc 
$f \in \fa_{\prec m}$.

\snii
Soit $f \in S \cap \fa_{\preceq m}$, donc $f = \sum_i g_i$
avec~$g_i \in \gen {f_i}$ et~$g_i \preceq m$, a fortiori~$f \preceq m$. 
Si~$m_i \nedivi m$ pour chaque~$i$, alors~$g_i \prec m$ donc
$f \in \fa_{\prec m}$. Si~$m_i \mid m$ pour un indice~$i$ i.e.
si~$m \in \gen {m_1, \ldots, m_s}$, alors la composante de~$f$
sur~$m$ est nulle car~$f \in S$. Comme~$f \preceq m$, on a
$f \prec m$; d'après le début de cette question, on en
déduit que~$f \in \fa_{\prec m}$.

\snii
On a $S \cap \fa_{\preceq m} \subset \fa_{\prec m}$; et pour tout~$f \in
\fa_{\prec m}$, il existe un \mom $m' \prec m$ tel \hbox{que $f \in \fa_{\preceq
m'}$}.  On en déduit au bout d'un nombre fini d'étapes que~$f = 0$, \hbox{i.e.  $S \cap \fa_{\preceq m} = \{0\}$}. Comme~$\fa$ est la réunion des~$\fa_{\preceq
m}$, on a bien~$S \cap \fa = \{0\}$.

%%%%%%
\snii\emph {3.}
Il s'agit d'une technique de division d'un \pol par le \sys $(f_1, \ldots, f_s)$
au sens des \bdgsz. On utilise le fait que si~$f \prec m$ pour un \mom $m$,
alors~$f \preceq m'$ pour un \mom~$m' \prec m$.  Supposons le résultat à
montrer vrai pour les \pols $\prec m$ et montrons le pour~$f\preceq m$. On
écrit~$f = am + g$ avec~$a \in \gk$ et~$g \prec m$.  Si~$m \in \gen {m_1,
\ldots, m_s}$, alors~$m = m'm_i$ pour un~$i$ et on remplace~$f$ par~$f -
am'f_i \prec m$.  Si~$m \notin \gen {m_1, \ldots, m_s}$, alors~$m \in S$ a
fortiori~$am \in S$, et on remplace~$f$ par~$f - am \prec m$.

%%%%%%
\snii\emph {4.}
Soit $E$ le sous-module de $\kuX^s$ engendré par les~$f_i\vep_j -
f_j\vep_i$. Soit~$m$ un \momz, des \pols $u_i$ tels que
$\sum_i u_i f_i = 0$ et~$u_i f_i \preceq m$. Il suffit de montrer
qu'il existe des \pols $v_i$ tels que

\snic {
\sum_i u_i\vep_i = \hbox {un \elt de $E$} + \sum_i v_i\vep_i
\hbox { avec } v_if_i \prec m.
}

\snii
Soit $I = \sotq {i \in \{1,\ldots,s\}} {m_i \hbox { divise } m}$.
Pour~$i \notin I$, on a~$u_if_i \prec m$. Pour~$i \in I$, on
définit le \mom $m'_i = m/m_i$; si~$a_i$ est la composante de~$u_i$
sur~$m'_i$, on a~$(u_i - a_im'_i) \prec m'_i$. Alors
$\sum_{i\in I} a_im_im'_i = 0$ i.e.~$\sum_{i\in I} a_i = 0$.
En suivant la solution de l'exercice \ref {exoMonomialSyzygies},
on pose, pour~$i,j \in I$

\snic {
q_{ij} = m'_i/m_j = m'_j/m_i = m/(m_im_j).
}

\snii
Il existe une matrice antisymétrique $(a_{ij})_{I \times I}\in\gk^{I \times I}$  telle 
que $a_i = \sum_{j \in J} a_{ij}$ de sorte que

\snic {
\sum_{i\in I} a_im'_i\vep_i = \sum_{i \in I} 
\bigl(\sum_{j \in I} a_{ij} q_{ij} \und {m_j}\bigr) \vep_i
\qquad \hbox {(ici $m_{ij} = m_j/(m_i \vi m_j) = m_j$).}
}

\snii 
Pour $i \in I$, on définit $w_i = \sum_{j \in I} a_{ij} q_{ij} \und {f_j}$
de sorte que $\sum_{i\in I} w_i\vep_i \in E$ 
\linebreak
et $w_if_i \preceq q_{ij}m_jm_i = m$.
Enfin on pose

\snic {
v_i = \formule {u_i-w_i& \hbox {si $i \in I$}, \\[.2em] u_i& \hbox {sinon}.}
}

\snii On a $v_if_i \prec m$ (car pour $i \in I$, la composante de $v_if_i$ sur
$m$ est $a_i - \sum_{j\in J} a_{ij} = 0$).  Et c'est gagné puisque

\snic {
\sum_i u_i\vep_i = \sum_{i \in I} w_i\vep_i + \sum_i v_i\vep_i =
\hbox {un \elt de $E$} + \sum_i v_i\vep_i. 
}

%%%%%%
\sni
\emph {5.}
Le fait que la suite de \moms $(m_1, \ldots, m_s)$ soit \ndze et est laissé \alecz. Montrons le pour $(f_1, \ldots, f_s)$. Il
suffit de voir que 

\snic{uf_s \in \gen {f_1, \ldots, f_{s-1}}
\;\Rightarrow\; u \in \gen {f_1, \ldots, f_{s-1}}.}

Encore une fois on raisonne
par \recu sur l'ordre monomial en écrivant $u = am + v$ avec $v \prec m$. En
multipliant par $f_s$, on a donc un $w \prec mm_s$ tel que $amm_s + w \in \gen
{f_1, \ldots, f_{s-1}}$et donc $amm_s \in \gen {m_1, \ldots,
m_{s-1}}$. La suite $(m_1, \ldots, m_s)$ étant \ndzez, on a $am \in \gen {m_1,
\ldots, m_{s-1}}$, disons $am = qm_i$ avec $i < s$.  
\\
Alors $(u - qf_i) f_s
\in \gen {f_1, \ldots, f_{s-1}}$ avec $u - qf_i \prec m$. Donc, par \recuz,
$u - qf_i \in \gen {f_1, \ldots, f_{s-1}}$, d'où l'on tire $u \in \gen
{f_1, \ldots, f_{s-1}}$.

%%%%%%
\snii\emph {6.}
Soit $\gA$ le $\gk[f_1, \ldots, f_s]$-module engendré par les \moms
de~$M$. Il faut d'abord montrer que $\kuX = \gA$. On procède par \recu sur
l'ordre monomial en supposant que tout \pol $g \in \kuX$ tel que $g \prec m$
est dans $\gA$ et en montrant que c'est vrai pour un \pol $f \in \kuX$ tel que
$f \preceq m$. On écrit $f = \sum_i u_if_i + r$ avec $u_if_i \preceq m$ et $r
\in S \subset \gA$. On a $u_i \prec m$ (car $m_i \ne 1$) donc $u_i \in \gA$
par \recuz. Bilan: $f \in \gA$.

\snii
Pour $e = (e_1, \ldots, e_s) \in \NN^s$ notons $f^e = f_1^{e_1}\cdots
f_s^{e_s}$. On va montrer simultanément que les $m \in M$ forment une base
de $\kuX$ sur $\gk[f_1,\ldots, f_s]$ et que les $f_i$ sont $\gk$-\agqt
indépendants en prouvant, pour une famille $(a_{e,m})_{e\in\NN^s, m\in M}$
d'\elts de $\gk$ à support fini, l'implication:
$$
\sum_{e,m} a_{e,m} f^e\, m = \sum_m \bigl( \sum_e a_{e,m} f^e\bigr) m = 
\sum_e \bigl( \sum_m a_{e,m} m\bigr) f^e = 0 
\ \Rightarrow\  a_{e,m} = 0
\leqno (\star)
$$
La clef réside dans les différentes façons de regrouper les sommes;
ceci est lié au fait d'écrire que les $m \in M$ constituent un \sgr de
$\kuX$ sur $\gk[f_1, \ldots, f_s]$:

\snic {
\begin {array} {rcl} 
\kuX &=& \sum_{m \in M} \gk[f_1, \ldots, f_s]m = \sum_{m\in M} 
\bigl(\sum_{e\in \NN^s} \gk f^e\bigr)m 
\\[.4em]
&=&
\sum_{e\in \NN^s} f^e \sum_{m \in M} \gk m =
\sum_{e\in \NN^s} f^e S\,.
\end {array}
}

\sni Comme la suite $(f_1, \ldots, f_s)$ est \ndze et 
que $\gen {f_1, \ldots, f_s} \cap S = 0$, on se convainc que la dernière
somme est directe. Par exemple, avec $2$ \pols $f_1, f_2$ et des 
$s_{ij} \in S$, on veut:

\snic {
\big(s_{00}\cdot 1 + s_{10}\cdot f_1 + s_{01}\cdot f_2 + 
s_{20}\cdot f_1^2 + s_{11}\cdot f_1f_2 + s_{02}\cdot f_2^2 = 0
\big)\ \Rightarrow\ s_{ij} = 0.
}

\snii On a $s_{00} \in S \cap \gen {f_1,f_2}$ donc $s_{00} = 0$. Ensuite,
on raisonne modulo $f_1$, et en utilisant le fait que $f_2$ est
\ndz modulo $f_1$, il vient $s_{01} + s_{02} f_2 \equiv 0 \bmod f_1$, donc
$s_{01} = 0$ \hbox{puis $s_{02} = 0$}. On peut alors simplifier par $f_1$ et ainsi de
suite.

\snii
Bilan: on a bien justifié $(\star)$, ce qui prouve le résultat escompté.

%%%%%%
\sni\emph {7.}
On a $a^2Y^3 = Yf_1 + (aY-X)f_2$. Si $Y^3 \in \fa := \gen {f_1,f_2}$, comme
les \pols sont \hmgsz, on a $Y^3 = (\alpha X + \beta Y) f_1 + (\gamma X +
\delta Y) f_2$ avec $\alpha, \beta, \gamma, \delta \in \gk$.  L'examen de la
composante sur $Y^3$ donne $1 = \delta a$.

\snii 
Supposons $a$ \ndzz. Si $\fa$ est facteur direct, $\gk[X,Y] = \fa \oplus S$,
alors $Y^3 = f + r$ \hbox{avec $f \in \fa$} et $r \in S$. Comme $a^2Y^3 \in \fa$, on
a $a^2r = 0$ puis $r = 0$; donc $Y^3 \in \fa$ et $a$ \ivz. 
\\
Réciproquement,
si $a$ est \ivz, on considère $f_1$, $a^{-1} f_2 = Y^2 + a^{-1} XY$ et
l'ordre lexi\-co\-gra\-phique avec $X \prec Y$: les \moms dominants sont $X^2$,
$Y^2$, premiers entre eux et on peut appliquer l'étude précédente.

}
%: fin sinotenglish

%%%%%%%%%%%%%%%%%%%%%%%%%%%%%%%%%%%%%%%%%%%%%%%%%%%%%%%%%%%%%%%%%%%%%%%%%%%
% fin des solutions d'exos

%:   ---- Section*{references}-----------
\Biblio

Le \thrf{th rg const loc libre} précise  le \tho 2 dans
\cite{Bou} chap.~II~{\S}5. \perso{développer ce genre de références, préciser aussi le titre du livre concerné}

La section \ref{secAppliIdenti} est basée sur les articles 
\cite[Chervov\&Talalaev]{CT03,CT03bis} qui s'occupent de
\gui{\syss de Hitchin} sur les courbes singulières.

Le \pb \ref{exoLambdaGammaK0} est
inspiré d'un article non publié de R.G. Swan: \emph{On a theorem of
Mohan, Kumar and Nori}.

Le \pb \ref{exoTh90HilbertMultiplicatif} provient d'un exercice
du chapitre 4 de \cite{Jensen}.

Dans le \pb \ref{exoFossumKumarNori}, la matrice $F_k$ intervient dans l'article: \emph{Vector bundles over
Spheres are Algebraic}, {\sc R. Fossum}, Inventiones Math. {\bf 8}, 222--225 (1969).
L'anneau $\gB_n$ est un classique en K-théorie algébrique.

\newpage \thispagestyle{CMcadreseul}
\incrementeexosetprob

%:        %%%%%%%%%%%%%%%%%%%%%%%%%%%%%%%%%%%%
%:        %%%%%%%%%%%%%%%%%%%%%%%%%%%%%%%%%%%%

%---- Chapitre  {Treillis distributifs}------------
\chapter[Treillis distributifs, groupes réticulés]{Treillis distributifs Groupes réticulés}
\label{chapTrdi}
\minitoc

%: Intro
\Intro
\pagestyle{CMExercicesheadings}

Ce chapitre commence par une section introductive qui fixe le cadre
\agq formel des \trdis et des \agBsz.

\smallskip
Les \trdis sont importants en \alg commutative pour plusieurs raisons.

\smallskip
D'une part la théorie de la \dve a comme \gui{modèle idéal} la théorie\linebreak
 de la \dve des entiers naturels. Si l'on prend
pour relation \linebreak
d'ordre $a \preceq b$, la relation \gui{$a$ est multiple de $b$},
on obtient que $\NN$ est un \trdi avec: $0$ \elt minimum, $1$ \elt maximum,
le \linebreak
sup $a\vu b$ égal au pgcd et le inf $a\vi b$ égal au ppcm. Quelques belles \prts de la \dve dans $\NN$ s'expriment en termes modernes en disant que l'anneau
$\ZZ$ est un anneau de Bézout (voir les sections \ref{secApTDN} et 
\ref{secBézout}).
Les nombres idéaux en théorie des nombres ont été créés
par Kummer pour combler l'écart entre la théorie de
la \dve dans les anneaux de nombres et celle dans $\NN$.
Les anneaux de nombres ne sont pas en \gnl des anneaux de Bézout,
mais leurs \itfsz\footnote{Ce qui pour Kummer était \gui{le pgcd idéal de plusieurs nombres} a été remplacé en langage moderne par l'\itf
correspondant. Ce \emph{coup de force} d\^u à Dedekind a été
une des premières intrusions de l'infini \gui{actuel} en \mathsz.}
forment un \trdiz, et leurs \itfs non nuls forment la partie positive d'un \grl (voir la section
\ref{secGpReticules}) ce qui rétablit
le bon ordon\-nan\-cement des choses.
Les anneaux dont les \itfs forment un \trdi sont appelés des \anars
(traités ailleurs dans les sections \ref{secIplatTf} et \ref{secAnars}).
Leurs \ids \ivs forment aussi la partie positive d'un \grlz.
La théorie des anneaux à pgcd (section \ref{secAnnPgcd})
trouve aussi son cadre naturel
dans le contexte des \grlsz.

\smallskip
D'autre part les \trdis interviennent comme la contrepartie \cov
des espaces spectraux divers et variés qui se sont imposés comme
des outils puissants de l'\alg abstraite.
La discussion sur ce sujet est particulièrement éclairante quand on
considère le treillis de Zariski d'un anneau commutatif,
relativement peu connu,
 qui sert de contrepartie \cov au très célèbre spectre de Zariski.
Espace spectral que l'on pourrait croire indispensable à la théorie de la \ddk
 et à celle des schémas de Grothendieck.
Une étude systématique du treillis de Zariski sera donnée au chapitre \ref{chapKrulldim} concernant la dimension de Krull, avec une introduction
heuristique dans la section~\ref{secEspSpectraux}.
Dans la section~\ref{secZarAcom} nous mettons en place le treillis de Zariski d'un anneau commutatif $\gA$ essentiellement
en rapport avec la construction de la clôture \zede réduite $\Abul$
(\paref{secClotureZEDR}) de l'anneau. Cette construction
peut être vue comme une construction \paral à celle de l'\agB
engendrée par un \trdi (voir le \thref{thZedGenEtBoolGen}).
L'objet global $\Abul$ ainsi construit contient
essentiellement la même information que le produit des anneaux~$\Frac(\gA\sur\fp)$ pour tous les \ideps $\fp$ de $\gA$.
Ceci alors même que dans la situation
\gnle on n'a pas accès \cot aux \ideps d'un anneau
de manière individuelle.

\smallskip
Une autre raison de s'intéresser aux \trdis est la logique \cov
(ou intuitionniste) dans laquelle l'ensemble des valeurs de vérité de la logique classique, à savoir $\so{\Vrai,\Faux}$,
qui est une \agB à deux \eltsz, est remplacé par un \trdi
plus mystérieux\footnote{En fait les valeurs de vérité des \coma ne
forment pas un ensemble à proprement parler, mais une classe.
Néanmoins les connecteurs logiques \cofs agissent
sur ces valeurs de vérité avec les mêmes propriétés
\agqs que le $\vi$ le $\vu$ et le $\to$ des \agHsz. Voir la discussion
\paref{P(X)}.
\label{NoteValVer}}.
La logique \cov sera abordée dans l'annexe (voir \paref{chapPOM})
particulièrement dans les sections \iref{secAffirmerProuver} et \iref{secBHK}.
Dans la section~\ref{secEntRelAgH} du chapitre présent
nous mettons en place les outils qui servent
de cadre à une étude \agq formelle de cette logique:
les \entrels et les \agHsz.
Il est remarquable que Heyting ait défini ces \algs dans la
première tentative
de décrire la logique intuitionniste de façon formelle, et
qu'il n'y ait pas eu une virgule à rajouter depuis.
Par ailleurs, \entrels et \agHs ont leur utilité propre dans l'étude \gnle des
\trdisz. Par exemple il est parfois important de pouvoir dire que le treillis
de Zariski d'un anneau est une \agHz.

\newpage	

%--- Sec{Treillis distributifs et \agBsz}
\section{Treillis distributifs et \agBsz}
\label{secTrdis}\label{secAGB}
\pagestyle{CMheadings}

%-----------------------------------------

Dans un ensemble ordonné $X$ on note, pour un $a\in X$:
%---  equation eqda --------
\begin{equation}\label{eqda}
\dar a=\sotq{x\in X}{x\leq a},\quad \uar a=\sotq{x\in X}{x\geq a}.
\end{equation}
%---------------------end equation--------------

On appelle \emph{chaîne croissante} une liste finie $(a_0,\ldots,a_n)$
 d'\elts de $X$ rangés en ordre croissant. Le nombre $n$ est appelé la \emph{longueur} de la chaîne. Par convention la liste vide est une chaîne de longueur $-1$.
\index{chaine@chaîne!dans un ensemble ordonné}
\index{longueur d'une chaîne!dans un ensemble ordonné}

%--- Definition{deftrdi}---------
\begin{definition}
\label{deftrdi}~
\begin{enumerate}
\item Un \ix{treillis} est un ensemble $\gT$ muni d'une relation d'ordre $\leq$
pour laquelle toute famille finie admet une borne supérieure et une borne
inférieure. On note $0_\gT$ le minimum de $\gT$ (la borne supérieure de la famille
vide) et  $1_\gT$ le maximum de $\gT$. On note $a\vu b$ la borne supérieure de
$(a,b)$ et~$a\vi b$ sa borne inférieure.
  \item Une application d'un treillis vers un autre est
  appelé un \emph{\homo de  treillis} si elle respecte les lois~$\vu$ 
  et~$\vi$ ainsi que les constantes $0$ et $1$.
  \item Le treillis est appelé un
\ixx{treillis}{distributif} lorsque chacune des deux lois~$\vu$ et~$\vi$ est
distributive par rapport à l'autre.
\end{enumerate}
\end{definition}
%--- end-definition------------------------------------

Les axiomes des treillis peuvent être formulés avec des \egts universelles
concernant uniquement les deux lois $\vi$ et $\vu$ et les deux constantes $0_\gT$
et $1_\gT$. La relation d'ordre est alors définie par:
$a\leq_\gT b\equidef a\vi b=a$. Voici ces axiomes.
\[\arraycolsep2pt
\begin{array}{rclcrcl}
a \vu a   & = & a &\qquad\qquad& a \vi a   & = & a \\
a \vu b   & = &b\vu a &\quad\quad& a \vi b   & = &b\vi a \\
(a \vu b)\vu c   & =&a \vu(b\vu c)  &&  (a \vi b)\vi c   & =&a \vi(b\vi c) \\
(a \vu b)\vi a   & = &a && (a \vi b)\vu a   & = &a \\
a\vu 0_{\gT}   & =  &a & &a\vi 1_{\gT}   & =  &a
\end{array}
\]

On obtient ainsi une théorie purement
équationnelle, avec toutes les facilités afférentes.
Par exemple on peut définir un treillis par \gtrs et relations.
Même chose pour les \trdisz.

\rdb
Dans un treillis, une distributivité implique l'autre.
Supposons par exemple que $a\vi(b\vu c)=(a\vi b)\vu(a\vi c)$, pour tous
$a,b,c$. Alors l'autre \dit résulte du calcul suivant:
\label{DistriTrdi}
\[
\begin{array}{c}
 (a\vu b)\vi(a\vu c)= \big((a\vu b)\vi a\big)\vu \big((a\vu b)\vi c\big)=a \vu \big((a\vu b)\vi c\big)=
 \\[1mm]
 a \vu\big(( a\vi c)\vu ( b\vi c) \big)= \big(a \vu( a\vi c)\big)\vu ( b\vi c)=a\vu  ( b\vi c).
\end{array}
\]

Dans un treillis discret on a un test pour $a\leq b$, puisque cette relation
équivaut à $a\vi b=a$.

Les sous-groupes d'un groupe (ou les \ids d'un anneau commutatif) forment un
treillis pour l'inclusion, mais ce n'est pas en \gnl un \trdiz.

Un ensemble totalement ordonné\footnote{Rappelons que c'est un ensemble $E$ muni d'une relation d'ordre $\leq$ pour laquelle on a, pour tous $x$ \hbox{et $y\in E$},
$x\leq y$ ou $y\leq x$. Ceci n'implique pas que l'\egt soit décidable.} est un \trdi s'il possède un \elt maximum
et un \elt minimum.
On note ${\bf n}$ l'ensemble totalement ordonné
à $n$ \eltsz.
Une application entre deux treillis totalement ordonnés $\gT$ et $\gS$
est un \homo \ssi elle est croissante et~$0_\gT$  et $1_\gT$  ont pour images
$0_\gS$ et $1_\gS$.

Si $\gT$  et $\gT'$ sont deux \trdisz, l'ensemble $\Hom(\gT,\gT')$ des \homos
 de  $\gT$ vers $\gT'$ est muni
d'une structure d'ordre naturelle donnée par
%-----------------begin $$----------------
$$ \varphi \leq \psi \equidef  \forall x\in \gT\;
\;
\varphi(x) \leq \psi(x)
$$
%-----------------end $$------------------

Un produit cartésien de \trdis est un \trdi
(pour les lois $\vi$ et $\vu$ produits, ce qui donne
la relation d'ordre partiel produit).

\rdb\label{NOTAtrdiopp}
Pour tout \trdi $\gT$, si l'on remplace la relation d'ordre $x\leq_\gT y$ par la
relation symétrique  $y\leq_\gT x$  on obtient le \ixx{treillis}{opposé}
$\gT\eci$ avec échange de $\vi$ et $\vu$ (on dit parfois \emph{treillis dual}).

Si $A\in\Pfe(\gT)$  avec un \trdi $\gT$ on notera
%-----------------begin $$----------------
$$ \Vu A:=\Vu_{x\in A}x\qquad {\rm et}\qquad \Vi A:=\Vi_{x\in A}x.
$$
%-----------------end $$------------------

%: \subsec{Treillis quotients, \idsz, filtres}%%%%%%%%%%%%
\subsec{Treillis quotients, \idsz, filtres}

\index{treillis distributif!quotient}
Si une structure \agq est définie par des lois de composition
de diffé\-rentes arités et des axiomes qui sont des \egts universelles
(comme les groupes, les anneaux, les \trdisz), une structure quotient
est obtenue lorsque l'on a une relation d'équivalence et que les lois
de composition \gui{passent au quotient}.
Si l'on regarde la structure comme définie par \gtrs et relations
(ce qui est toujours possible), on obtient une structure quotient en rajoutant des relations.

Un \emph{treillis quotient $\gT'$ d'un treillis $\gT$} peut \egmt
être donné par une relation
binaire
$\preceq$ sur $\gT$ vérifiant les \prts suivantes:
%---  equation eqPreceq --------
\begin{equation}\label{eqPreceq}
\left.
%--------------------begin array---------------
\begin{array}{rcl}
a\leq b&  \Longrightarrow  & a\preceq b   \\
a\preceq b,\,b\preceq c&  \Longrightarrow  & a\preceq c   \\
a\preceq b,\,a\preceq c&  \Longrightarrow  & a\preceq b\vi c   \\
b\preceq a,\,c\preceq a&  \Longrightarrow  & b\vu c\preceq a
\end{array}
%---------------------end array--------------
\right\}
\end{equation}
%---------------------end equation--------------
La relation $\preceq$ induit alors une structure de treillis
sur l'ensemble quotient~$\gT'$ obtenu avec la nouvelle \egtz\footnote{Le fait
de ne pas changer d'objets quand on passe au quotient, mais de changer seulement la relation d'\egt sur les objets
est plus simple, mais aussi plus conforme à la tradition (Gauss)
et à l'implémentation sur machine. Sans doute le succès populaire
 des classes d'équivalence comme objets de l'ensemble quotient
est d\^u en bonne partie à l'heureux hasard qui fait que
dans le cas d'un groupe quotient $G/H$, en notation additive par exemple,
on a $(x+H)+(y+H)=(x+y)+H$ avec trois significations différentes
du symbole $+$. Les choses se passent pourtant moins bien pour les anneaux quotients, et par exemple $(3+7\ZZ)(2+7\ZZ)$ n'est pas égal à
$6+7\ZZ$, mais seulement contenu dedans.}
$$
(a,b\in\gT)\quad:\quad\quad a=_{\gT'}b \equidef (a\preceq b \;\mathrm{ et }\;
b\preceq a)
$$
Naturellement si $\gT$ est distributif, il en va de même pour $\gT'$.

Si $\varphi :\gT\rightarrow \gT'$ est un \homo de \trdisz,
$\varphi^{-1}(0)$
est appelé un \emph{\id de $\gT$}. Un \id $\fb $ de $\gT$ est une
partie de
$\gT$ soumise
aux  contraintes suivantes:
%---  equation eqIdeal --------
\begin{equation}\label{eqIdeal}
\left.
\begin{array}{rcl}
   & &  0 \in \fb    \\
x,y\in \fb & \Longrightarrow   &  x\vu y \in \fb    \\
x\in \fb ,\; z\in \gT& \Longrightarrow   &  x\vi z \in \fb    \\
\end{array}
\right\}
\end{equation}
%---------------------end equation--------------
(la dernière se réécrit $(x\in \fb ,\;y\leq x)\Rightarrow y\in
\fb $).
Un \emph{\idpz} est un \id engendré par un seul \elt $a$,
il est égal à $\dar a$.%
\index{ideal@idéal!d'un \trdiz}%
\index{ideal@idéal!principal (d'un \trdiz)}\label{NOTAdara}%
\index{principal!ideal@idéal --- d'un \trdiz}

L'\id $\dar a$, muni des lois $\vi$ et $\vu$ de $\gT$ est un \trdi
dans lequel l'\elt maximum est $a$.  L'injection canonique $\dar
a\rightarrow \gT$ n'est pas un morphisme de \trdis parce que
l'image de $a$ n'est pas égale à $1_\gT$.  Par contre
l'application surjective $\gT\rightarrow \dar a,\;x\mapsto x\vi a$
est un morphisme surjectif, qui définit donc $\dar a$ comme une
structure quotient.

\rdb \label{NOTAuara}
La notion opposée à celle d'\id est la notion de {\em
filtre}.  Le filtre principal engendré par $a$ est égal à $\uar a$.%
\index{filtre!d'un treillis distributif}%
\index{filtre!principal d'un treillis distributif}%
\index{principal!filtre --- d'un treillis distributif}

L'\emph{\id engendré} par une partie $J$ de $\gT$ est égal à
$$
\cI_\gT(J)=\big\{x\in\gT \;\big \vert \; \Ex J_0\in \Pfe(J),\,x\leq \Vu J_0\big\}.
$$
En conséquence \emph{tout \itf est principal}.
\label{NOTAidtrdi}

Si $A$ et $B$ sont deux parties de $\gT$ on note
%---  equation eqvuvi --------
\begin{equation}\label{eqvuvi}
A\vu B=\big\{ a\vu b \mid a\in A,\,b\in B\,\big\}  \; \hbox{ et } \; A\vi
B=\big\{ a\vi b \mid a\in A,\,b\in B\,\big\}.
\end{equation}
%---------------------end equation--------------
Alors l'\id engendré par deux \ids $\fa$ et $\fb$ est égal à
%---  equation eqSupId --------
\begin{equation}\label{eqSupId}
\cI_\gT(\fa\cup \fb) = \fa\vu\fb.
\end{equation}
%---------------------end equation--------------
L'ensemble des \ids de $\gT$ forme lui-même un \trdiz\footnote{En fait
il faut introduire une restriction pour obtenir vraiment un ensemble, de façon à ce que l'on ait un procédé bien défini de construction des \ids concernés.
Par exemple on peut considérer l'ensemble des \ids obtenus à partir des \idps par certaines opérations prédéfinies, comme les réunions et intersections dénombrables.
C'est le même \pb que celui indiqué dans la note~\vref{NoteValVer}.} pour
l'inclusion, avec pour borne inférieure de $\fa$ et $\fb$ l'\id 
%---  equation eqInfId --------
\begin{equation}\label{eqInfId}
\fa\cap \fb=\fa\vi\fb.
\end{equation}
%---------------------end equation--------------

Ainsi les opérations $\vu$ et $\vi$ définies en (\ref{eqvuvi})
correspondent au sup et au inf dans le treillis des \idsz.

 \rdb\label{NOTAfitrdi}
On notera
$\cF_\gT(S)=\sotq{\,x\in\gT}{\Ex S_0\in \Pfe(S),\,x\geq \Vi S_0}$
le filtre de $\gT$ engendré par le sous-ensemble $S$.
\\
Quand on considère le treillis des filtres, il faut
faire attention à ce que produit le renversement de la relation
d'ordre: $\ff\cap\ffg=\ff\vu\ffg$ est le inf des filtres~$\ff$ et $\ffg$,
tandis que leur sup est égal à $\cF_\gT(\ff\cup \ffg)=\ff\vi
\ffg$.

\rdb
 Le \emph{treillis quotient de $\gT$ par l'\id $\fa$}, noté
$\gT/(\fa=0)$ est défini comme le \trdi engendré par les \elts de
$\gT$ avec pour relations, les relations vraies dans $\gT$ d'une part,
et les relations $x=0$ pour les $x\in \fa$ d'autre part.  Il peut
aussi être défini par la relation de préordre suivante
%-----------------begin $$----------------
$$ %a\preceq b\quad\Longleftrightarrow\quad
\label{trquoideal}
a\leq_{\gT/(\fa=0)}b\;\;
\equidef \;\;
\exists x\in \fa \;\;a\leq  x\vu b.
$$
%-----------------end $$------------------
Ceci donne
%-----------------begin $$----------------
$$ a\equiv b\;\;\mod\;(\fa=0)\quad \Longleftrightarrow  \quad \exists
x\in \fa
\;\;\;a\vu x=b\vu x.
$$
%-----------------end $$------------------
En particulier, l'\homo de passage au quotient
$$\varphi:\gT\to\gT'=\gT/(\fa=0)$$
vérifie $\varphi^{-1}(0_{\gT'})=\fa$. Dans le cas du quotient par un \id principal $\dar a$ on obtient
$\gT/(\dar a=0)\simeq\uar a$ avec le morphisme $y\mapsto y\vu a$ de $\gT$
vers $\uar a$.

%:--- Proposition{propIdealFiltre}----
\begin{proposition}
\label{propIdealFiltre} Soit $\gT$ un \trdi et
$(J,U)$ un couple de parties de $\gT$.
On considère le quotient $\gT'$ de $\gT$ défini par les
relations $x=0$ pour les $x\in J$, et $y=1$ pour les $y\in U$. Alors
  l'in\egt $a\leq_{\gT'}b$ est satisfaite \ssi
il existe  $J_0\in\Pfe( J)$  et  $U_0\in\Pfe( U)$  tels que:
%---  equation eqpropIdealFiltre --------
\begin{equation}\label{eqpropIdealFiltre}
a \vi \Vi U_0 \; \leq_\gT\; b \vu \Vu J_0
\end{equation}
%---------------------end equation--------------
Nous noterons $\gT/(J=0,U=1)$ ce treillis quotient $\gT'.$
\end{proposition}
%--- end-proposition----------------------------------------

On voit sur l'exemple des ensembles totalement ordonnés
qu'une structure quotient d'un \trdi n'est pas en \gnl
\caree par les classes d'équivalence de $0$ et $1$.

\rdb
Soient $\fa$ un \id et $\ff$ un filtre de $\gT$.
On dit que $\fa$ est \emph{$\ff$-saturé} si l'on a
$$
(g\in \ff,\; x\vi g \in \fa) \Longrightarrow x\in \fa,
$$
on dit que $\ff$ est \emph{$\fa$-saturé} si l'on a
$$
(a\in \fa,\; x\vu a \in \ff) \Longrightarrow x\in \ff.
$$
Si $\fa$ est $\ff$-saturé et $\ff$ est $\fa$-saturé on dit que $(\fa,\ff)$
est un \emph{\paz} dans~$\gT$.%
\index{couple sature@\paz}%
\index{saturé!couple ---}%
\index{saturé!idéal $\ff$- ---}%
\index{saturé!filtre $\fa$- ---}

%:     Fact{factIdFiAsTrdi}
\begin{fact}\label{factIdFiAsTrdi}
Soit $\varphi:\gT \rightarrow \gT_1$ un \homo de \trdisz.
\linebreak 
L'idéal $\fa=\varphi^{-1}(0)$ et  le filtre
$\ff=\varphi^{-1}(1)$ forment un \paz.
Réciproquement, si $(\fa,\ff)$ est un \pa de $\gT$,
%et si $\gT'$ est le treillis quotient obtenu en forçant $\fa=0$ et $\ff=1$,
l'\homo de passage au quotient $\pi:\gT \rightarrow \gT/(\fa=0,\ff=1)$ vérifie $\pi^{-1}(0)=\fa$
et $\pi^{-1}(1)=\ff$.
\end{fact}

Lorsque $(\fa,\ff)$  est un \paz, on a
les \eqvcs

\snic{
1\in \fa\; \;\Longleftrightarrow\;\;  0\in \ff
\;\; \Longleftrightarrow\;\;  (\fa,\ff)=(\gT ,\gT )
}

%%%%%%%%%%%%%%%%%%%%%%%%%%%%%%%%%%%%%%%%%%%%%%%%%%%%%%%%%%%%%%%%%%%%
%:subsec{Les \agBsz}
\subsec{Les \agBsz}

Dans un \trdi un \elt $x'$ est appelé un \emph{complément}, ou un
\emph{complémentaire}, de $x$ si l'on
a $x\vi x'=0$ et $x\vu x'=1$. S'il existe le complément de $x$
est unique. Il est alors souvent noté  $\lnot x$.
\index{complement@complément!(dans un \trdiz)}

Rappelons que par \dfn un anneau $\gB$ est une \agB \ssi tout \elt
est idempotent. On définit alors une relation d'ordre $x\preceq y$ par: $x$ est multiple
de $y$, \cad $\gen{x}\subseteq\gen{y}$.\\
On obtient ainsi un \trdi
dans lequel tout \elt $x$ admet pour complément $x'=1+x$
(cf.~proposition~\ref{defiBoole}).

On a la réciproque suivante.

%:     propBooleTrdi
\begin{proposition}\label{defiBooleTrdi}
\emph{(Algèbres de Boole)}
\begin{enumerate}
\item Sur un \trdi dans lequel tout \elt $x$ admet
un complément, noté $\lnot x$, on peut définir une structure d'\agB
en posant
$$\preskip.0em \postskip.4em 
 xy=x\vi y\;\hbox{  et  }\;x\oplus y=(x\vi \lnot y)\vu (y\vi \lnot x).
$$
On retrouve  $x\vu y=x\oplus y\oplus xy$ et $\lnot x =1\oplus x$.
\item Tout \homo de \trdis entre deux \agBs
est un \homo d'\agBsz, et il respecte le passage au complémentaire.
\end{enumerate}
\end{proposition}

%:-% entrenous
\entrenous{
La proposition \ref{propBoolFini} affirme que tout \agB discrète se comporte dans les calculs comme l'ensemble ds parties finies d'un ensemble fini.

De même le principe de recouvrement par quotients 2.10 pour les \grls
peut être paraphrasé comme suit: dans les calculs, un \grl se
comporte toujours comme un produit fini de groupes totalement ordonnés.

De la même manière tout \trdi se comporte dans les calculs comme un
produit fini d'ensembles totalement ordonnés.
Mais nous n'avons pas donné l'énoncé analogue de la proposition \ref{propBoolFini} et du principe 2.10.

En \clamaz, ces résultats sont énoncés comme des \gui{\thos de représentation}, que l'on obtient à partir des énoncés \cofs au moyen de zornettes adéquates. 
}
%-% Fin entrenous

%\newpage
%: \subsec{Algèbre de Boole engendrée par un \trdiz}
\subsec{Algèbre de Boole engendrée par un \trdiz}

%:2015 Ce début de paragraphe, avant le lemme 1.5 est précisé

Commençons par quelques remarques sur les \elts qui possèdent
un complément
dans un \trdiz.
Si $a$ admet un complément $a'$, puisque $b=(b\vi a)\vu(b\vi a')$ pour tout $b\in\gT$,
l'\homo canonique
$$
\gT\to \gT\sur{(a=1)}\times \gT\sur{(a'=1)}
$$
est injectif. En outre, ce morphisme est surjectif parce que, pour $x,y\in \gT$, en posant $z=(x\vi a)\vu(y\vi a')$,
on a $z\vi a =x\vi a$, i.e. $z\equiv x \mod (a=1)$, et, de même, 
$z\equiv y \mod (a'=1)$.  Réciproquement, on a le résultat suivant
qui montre la similitude entre un \idm dans un anneau commutatif et un
\elt possédant un complément dans un \trdi
(voyez le fait~\ref{lemCompAnnComm}).

%:     Lemma{lemCompTrdi}
\begin{lemma}\label{lemCompTrdi} 
Pour tout \iso $\lambda:\gT\to\gT_1\times \gT_2$, il existe un (unique)
\elt $a\in\gT$ tel que:
\begin{enumerate}
\item $a$ possède un complément $\lnot a$,
\item l'\homo composé $\gT\to\gT_1$
identifie $\gT_1$ avec $\gT\sur{(a=0)}$ ainsi qu'avec $\gT\sur{(\lnot a=1)}$,
\item l'\homo composé $\gT\to\gT_2$
identifie $\gT_2$ avec $\gT\sur{(a=1)}$ ainsi qu'avec $\gT\sur{(\lnot a=0)}$.
\end{enumerate} 
\end{lemma}

\begin{proof}
L'\elt $a$ est donné par $\lambda(a)=(0_{\gT_1},1_{\gT_2})$.
\end{proof}

Lorsque deux \elts $a$ et $b$ possèdent des compléments
$\lnot a$ et $\lnot b$, les \ix{lois de Morgan} sont vérifiées:
\begin{equation}
\label{eqMorgan}
\lnot(a\vi b)=\lnot a \vu \lnot b \;\hbox{  et }\; \lnot(a\vu b)=\lnot a \vi \lnot b.
\end{equation}

\medskip
Par \dfnz, l'\emph{\agB librement engendrée par le \trdi $\gT$}
est donnée par un couple $(\Bo(\gT),\lambda)$, où $\Bo(\gT)$ est une \agBz, et où $\lambda:\gT\to\Bo(\gT)$ est un \homo de \trdis satisfaisant la \prt universelle suivante.
\\
\emph{Tout \homo $\psi$ de \trdis de $\gT$ vers une \agB $\gB$
se factorise de manière unique sous la forme $\varphi\circ \lambda$.}
\PNV{\gT}{\lambda}{\psi}{\Bo(\gT)}{\varphi}{\gB}{\trdisz}{}{\agBsz}

\vspace{-1.1em}
Puisque nous sommes dans le cadre de structures \agqs
purement équationnelles, cette \agB peut être construite à partir 
%\linebreak
de $\gT$ en rajoutant par force une loi unaire
$a\mapsto \lnot a$ et en imposant les 
axio-\linebreak
mes
$a\vi \lnot a=0$, $a\vu \lnot a=1$.

Autrement dit encore $\Bo(\gT)$ peut être définie comme une \agB
obtenue par \gtrs et relations.
Les \gtrs sont les \elts de~$\gT$ et les relations sont
celles qui sont vraies dans $\gT$: de la forme $a\vi b=c$ ou~$a\vu b= d$, sans oublier $0_{\Bo(\gT)}=0_\gT$ et $1_{\Bo(\gT)}=1_\gT$.

Cette description est cependant peu précise et nous allons
construire $\Bo(\gT)$ à la vitesse de la tortue pour y voir plus clair.

%:     Lemma{lemRajouCompl}
\begin{lemma}\label{lemRajouCompl}
Soit $\gT$ un \trdi et $a\in\gT$. Considérons le \trdi
$$\gT[a\bul]\eqdefi\gT\sur{(a=0)}\times \gT\sur{(a=1)}$$
et $\lambda_a:\gT\to\gT[a\bul]$
l'\homo canonique. %Alors
\begin{enumerate}
\item  L'\homo $\lambda_a$ est injectif et  $\lambda_a(a)=(0,1)$
admet $(1,0)$ pour com\-plément.
\item
Pour  tout \homo $\psi
:\gT\to\gT'$ tel que $\psi(a)$ admette un complément,
il existe un unique \homo $\varphi:\gT[a\bul]\to \gT'$
tel que~$\varphi\circ \lambda_a=\psi$.
\Pnv{\gT}{\lambda_a}{\psi}{\gT[a\bul]}{\varphi}{\gT'}{ }{}{$\psi(a)$ admet un complément}
\end{enumerate}
\end{lemma}

\begin{proof}
Le lemme \ref{lemCompTrdi} donne
$\gT'\simeq \gT'\sur{(\psi(a)=0)}\times \gT'\sur{(\psi(a)=1)}$, d'où
l'\homo $\varphi$ et l'unicité. L'injectivité de $\lambda_a$ ne saute pas aux yeux mais c'est un grand classique:
si $x\vi a=y\vi a$ et $x\vu a=y\vu a$, alors

\snic{x=(x\vu a)\vi x=(y\vu a)\vi x=(y\vi x)\vu(a\vi x).}

%\sni
De manière symétrique
$y=(y\vi x)\vu(a\vi y)$, donc $x=y$ puisque $a\vi x=a\vi y$.
\end{proof}

%:     Corollary{corlemRajouCompl}
\begin{corollary}\label{corlemRajouCompl}
Soient $a_1$, \ldots, $a_n\in\gT$.
\begin{enumerate}
\item Le treillis $\gT[a_1\bul][a_2\bul]\cdots[a_n\bul]$
est indépendant, à \iso unique près, de l'ordre des $a_i$.
Il sera noté $\gT[a_1\bul,a_2\bul,\ldots,a_n\bul]$.
\item Une description possible est la suivante:
$$
\gT[a_1\bul,a_2\bul,\ldots,a_n\bul]
\;\simeq\;
\prod\nolimits_{I\in\cP_n}{\gT\sur{\bigl( (a_i=0)_{i\in I}, (a_j=1)_{j\in \lrbn\setminus I}\bigr)}}.
$$
\item L'\homo naturel $\gT\to\gT[a_1\bul,a_2\bul,\ldots,a_n\bul]$ est injectif.
Il factorise de manière unique tout \homo $\psi$ de $\gT$ vers un \trdi
$\gT'$ tel que les $\psi(a_i)$ admettent un complément.
\end{enumerate}
\end{corollary}

%:     Theorem{thBoolGen}
\begin{theorem}\label{thBoolGen}
\emph{(Algèbre de Boole librement engendrée par un \trdiz)}
Pour tout \trdi $\gT$ il existe une \agBz, notée $\Bo(\gT)$, avec un
\homo $\lambda:\gT\to\Bo(\gT)$, qui
factorise de manière unique tout \homo $\psi :\gT\to\gB$ vers une \agBz.
Ce couple $(\Bo(\gT),\lambda)$ est unique à \iso unique près.
On a de plus les \prts suivantes.
\begin{enumerate}
\item [--] L'\homo  $\lambda$ est injectif.
\item [--] On a $\Bo(\gT)=\gT[(a\bul)_{a\in\gT}]$
\end{enumerate}
\end{theorem}
%
%:2015 La preuve est beaucoup plus détaillée
\begin{proof} On prend pour $\Bo(\gT)$
 la limite inductive (filtrante) des $\gT[a_1\bul,a_2\bul,\ldots,a_n\bul]$. L'injectivité de $\lambda$ et l'unicité de la \fcn résultent du lemme \ref{lemRajouCompl}. Il reste à voir que cette limite inductive est bien une \agBz. Tout d'abord
 c'est un \trdi comme limite inductive de \trdisz. Ensuite, chaque $a_i$ admet un complément dans $\gT[a_1\bul,a_2\bul,\ldots,a_n\bul]$, lequel donne bien un complément de $\lambda(a_i)$ dans la limite inductive, parce que tout morphisme de \trdis conserve les compléments qui existent.  
Enfin, les lois de Morgan, comme énoncées précédemment assurent par \recu  que tout \elt de $\Bo(\gT)$ admet un complément.
\end{proof}
%

%:2015 exemples plus nombreux et plus précis
\exls \label{exemplesthBoolGen}
1) Supposons que $\gT$ soit un treillis de parties détachables
d'un ensemble $E$, au sens que si $A$ et $B$ sont des \elts de $\gT$,
alors $A \cup B$ \hbox{et $A\cap B$} \egmtz, avec en outre $\emptyset$ et $E$ \elts de
$\gT$. Alors $\Bo(\gT)$ s'identifie à l'ensemble des combinaisons
booléennes finies d'\elts de $\gT$ (ce sont toutes des parties détachables de $E$), qui est une \agB de parties de $E$.
Voir à ce sujet l'exercice \ref{exoagGgen1} et son corrigé. 

2) Soit $\gT$ un ensemble totalement ordonné discret admettant un minimum~$0_\gT$ et un maximum $1_\gT$. Alors $\gT$ est un \trdi et il est isomorphe au treillis des parties de $\gT\setminus\so{1_\gT}$ qui sont de la forme $I_a=\sotq{x\in \gT}{x<a}$. L'\iso est donné par $a\mt I_a$. Par suite $\Bo(\gT)$ est formé par l'ensemble des parties qui sont réunions finies d'intervalles semi-ouverts~$[\,a_i,b_i\,[\,$.
L'écriture est unique si l'on réclame $a_1<b_1<a_2<\dots<b_n$ (la partie vide correspond à une réunion vide, et la partie pleine est $[\,0_\gT,1_\gT\,[\,$).\\
Un exemple simple avec $\gT$ infini est  le suivant. On considère $\gT=\NN\cup\so{+\infty}$, alors $\Bo(\gT)$ s'identifie à l'ensemble des parties finies ou cofinies de $\NN$ (données en tant que telles, bien évidemment).
\\
Notons que lorsque $\gT$ est un ensemble totalement ordonné \emph{non} discret, la description de $\Bo(\gT)$ est nettement plus délicate.
\eoe

\medskip
\comm En \clamaz, tout \trdi est isomorphe à un sous-treillis 
du treillis des
parties d'un ensemble. Cela fournit en suivant l'exemple~1) ci-dessus une \gui{construction} alternative
de l'\agB $\Bo(\gT)$.
\eoe

%%%%%%%%%%%%%%%%%%%%%%%%%%%%%%%%%%%%%%%%%%%%%%%%%%%%%%%%%%%%%%%%%
%--- Sec{Groupes réticulés}
\section{Groupes réticulés}
\label{secGpReticules}
%-----------------------------------------

%:  subsec{Premier pas}
\subsec{Premier pas}

Dans cet ouvrage nous nous limitons, pour les groupes ordonnés,
au cas des groupes commutatifs.

%:     Definition{defiGpReticule}
\begin{definition}\label{defiGpReticule}
On appelle \emph{groupe ordonné} un groupe abélien $G$ muni d'une relation
d'ordre partiel \emph{compatible} avec la loi de groupe, i.e., en notation additive,
$$\forall a,x,y\in G\qquad x\leq y\; \Longrightarrow\; a+x\leq a+y.
$$
Un groupe ordonné est dit \emph{réticulé} lorsque deux \elts arbitraires admettent une borne inférieure, que l'on notera $x\vi y$.
Si \ncrz, on précise la structure en écrivant $(G,0,+,-,\vi)$.
Un \emph{morphisme de \grlsz} est un \homo de groupe qui respecte la loi $\vi$.%
\index{groupe ordonné}\index{groupe réticulé}
\end{definition}

Un groupe abélien muni d'un ordre total compatible
(on dit un \emph{groupe totalement ordonné}) est un \grlz.
Les morphismes de groupes totalement ordonnés sont alors les \homos de groupes croissants.\index{totalement ordonné!groupe ---}

Un \emph{sous-\grlz} d'un \grl $G$ est par \dfn un sous-groupe stable pour
la loi de treillis $\vi$. Il ne suffit pas pour cela que la relation d'ordre induite sur le sous-groupe en fasse un treillis.
\index{sous-groupe réticulé}

Une idée directrice dans la théorie des \grls est qu'\emph{un \grl se comporte dans les calculs comme un produit de groupes totalement
ordonnés}. Ceci se traduira de manière \cov par le \prf \vref{prcfgrl}.

\medskip
\exls
1) (Attention, notation multiplicative!) L'ensemble $\QQ^{>0}$ des rationnels
strictement positifs est un \grl avec pour partie positive le \mo $(\NN^{>0},1,\times )$.
L'exemple de cette structure multiplicative est paradigmatique.
On a un \iso de \grls $\QQ^{>0}\simeq \ZZ^{(P)}$,\label{NotaZZP} où $P$ est l'ensemble des nombres premiers, $\ZZ^{(P)}=\bigoplus_{p\in P}\ZZ$ et l'ordre
est induit par l'ordre produit. Ceci n'est qu'une autre manière d'exprimer le \tho fondamental de l'arithmétique \gui{tout entier naturel s'écrit de manière unique comme produit de puissances de nombres premiers}.
C'est en voulant à tout prix faire ressembler la multiplication pour les
entiers des corps de nombres à la multiplication dans  $\NN^{>0}$
que les mathématiciens ont été amenés à inventer les \gui{nombres pgcd idéaux}.

\rdb
2)
Si $(G_i)_{i\in I}$ est une famille de \grls 
%:2015
indexée par un ensemble discret $I$%
, on définit
la \emph{\sdoz} de la famille, notée $\boxplus_{i\in I}G_i$%
\label{NotaSDirOr}, qui est un \grl avec comme groupe sous-jacent le groupe $\bigoplus_{i\in I}G_i$, la loi $\vi$ étant définie \coo par \cooz. 
%:2015 Si 
Lorsque $I$ est fini, par exemple
$I=\lrb{1..3}$, on notera $G_1\boxplus G_2\boxplus G_3 $.%
\index{somme directe o@\sdoz}
\\
Par exemple $\ZZ^{(P)}=\boxplus_{p\in P}\ZZ$. 
\\
On définit \egmt le produit $\prod_{i\in I}G_i$ de manière usuelle,
et c'est le produit dans la catégorie des \grlsz.
Si $I$ est un ensemble fini, les \grls  $\boxplus_{i\in I}G_i$
et $\prod_{i\in I}G_i$ sont naturellement isomorphes.

3) 
Si $(G_i)_{i\in I}$ est une famille de groupes totalement ordonnés discrets avec pour $I$ un ensemble totalement ordonné discret 
on définit la {somme lexicographique} de cette famille, c'est le groupe 
totalement ordonné discret~$G$ dont le groupe sous-jacent est  $\bigoplus_{i\in I}G_i$ et la relation d'ordre est l'ordre lexicographique: $(x_i)_{i\in I}< (y_i)_{i\in I}$ \ssi $x_{i_0}< y_{i_0}$ pour le plus petit indice $i_0$ tel que $x_{i_0}\neq y_{i_0}.$
\eoe

\medskip
Dans un \grl les translations sont des \autos de la structure d'ordre,
d'où la règle de \dit
\begin{equation}
\label{eq1grl}
{x+(a\vi b)\,=\,(x+a)\vi(x+b).}
\end{equation}

On voit aussi que la bijection $x\mapsto-x$ renverse l'ordre,
et donc que deux \elts arbitraires $x$, $y$ admettent aussi une borne
supérieure

\snic{x\vu y=-\big((-x)\vi(-y)\big),}

%\sni
avec $x+y-(x\vu y)=(x+y)+\big((-x)\vi(-y)\big)=(x+y-x)\vi (x+y-y)$, donc
\begin{eqnarray}
\label{eq2grl}
x+y&=&(x\vi y)+(x\vu y),\\
x+(a\vu b)&=&(x+a)\vu(x+b).\label{eq3grl}
\end{eqnarray}

Il manque cependant un \elt minimum et un \elt maximum pour obtenir
un treillis.

%:  subsec{Identités remarquables dans les \grlsz}
\subsec{Identités remarquables dans les \grlsz}

\Grandcadre{Dans ce paragraphe $G$ est un \grl \\
et $G^+$ le sous-\mo de $G$ formé des \elts positifs ou nuls.}

On note

\snic{x^+=x\vu 0, \quad x^-=(-x)\vu 0\; \et \;\abs{x}=x\vu(-x).}

%\sni
 On les appelle
respectivement la \ix{partie positive}, la \ix{partie négative} et la \ix{valeur absolue} de $x$.

%:     Theorem{th0GpRtcl}
\begin{theorem}\label{th0GpRtcl} \emph{(Distributivité dans les
\grlsz)}\\
Dans un \grl les lois $\vi$ et $\vu$ sont distributives
l'une par rapport à l'autre.
\end{theorem}
\begin{proof}
Il suffit de montrer $x\vu(y_1\vi y_2)= (x\vu y_1)\vi (x\vu y_2)$. En
translatant \hbox{par $-x$}, on se ramène à $x = 0$, i.e. à $(y_1
\vi y_2)^+ = y_1^+ \vi y_2^+$.  
\\
L'in\egt $(y_1 \vi y_2)^+ \le y_1^+ \vi y_2^+$ est \imdez.  
\\
Posons $y=y_1\vi y_2$.
L'\elt $y_i + y^+ - y$ est  $\ge
y_i$ et $\geq 0$, donc $\geq y_i^+$.
\\
D'où
$y_i^+ + y \le y_i + y^+$. Puis
$(y_1^+ + y) \vi (y_2^+ + y) \le (y_1 + y^+) \vi (y_2 + y^+)$, i.e.~$(y_1^+ \vi
y_2^+) + y \le (y_1 \vi y_2) + y^+$, i.e.~$y_1^+ \vi y_2^+ \le y^+$.
\end{proof}

Deux \elts $x$, $y$  sont dits \emph{disjoints} ou \emph{\ortsz}
si  $\abs{x}\vi\abs{y}=0$.\index{orthogonaux!éléments --- dans un \grlz}

%:     Lemma{lemx+x-}
\begin{lemma}\label{lemx+x-}\relax Soient $x$, $y\in G$.
\vspace{-1mm}
\begin{eqnarray}
\label{i1lemx+x-}\relax
&x=x^+ - x^-,\quad x^+ \perp x^-,\quad \abs{x}=x^++x^-=x^+\vu x^-\in G^+& \\
\label{i2lemx+x-}\relax
&x\leq y \; \iff\; x^+\leq y^+ \et y^-\leq x^-,\quad x=0\;\iff\;\abs{x}=0&%
\end{eqnarray}
\end{lemma}
\begin{proof} ({\ref{i1lemx+x-}}).
Tout d'abord $x^+-x=(x\vu 0)-x=(x-x)\vu\big(0+ (-x)\big)=x^-$. \\
Toujours par \dit on obtient

\snac{x^+ + x^-= (x\vu 0)+((-x)\vu 0)= (x-x)\vu (x+ 0)\vu\big(0+(-x)\big)\vu(0+ 0)=x^+\vu x^- .}

%\sni
 Enfin puisque $x^+ + x^-=(x^+\vu x^-)+(x^+\vi x^-)$, cela donne $x^+\vi x^-=0$.

 ({\ref{i2lemx+x-}}).
Laissé \alecz.
\end{proof}
%

%:     Lemma{lemGlem}
\begin{lemma}\label{lemGlem} \emph{(Lemme de Gauss)}
Soient $x$, $y$, $z\in G^+.$
\vspace{-1mm}
\begin{eqnarray}
\label{i1lemGlem}
(x\perp y \et x\leq y+z) &\;\Longrightarrow\;& x\leq z \\
\label{i2lemGlem}
x\perp y &\Longrightarrow& x\vi (y+z) = x\vi z\\
\label{i3lemGlem}
(x\perp y \et x\perp z) &\Longrightarrow& x\perp (y+z) \\
\label{i4lemGlem}
(x\perp y \et x\leq z \et y\leq z) &\Longrightarrow& x+y\leq z
\end{eqnarray}
\end{lemma}
\begin{proof}
{(\ref{i1lemGlem}).} On a $x\leq z+x$ parce que $z\geq0$
et $x\leq z+y$ par hypothèse, 
\linebreak
donc $x\leq (z+x)\vi(z+y)=z+(x\vi y)=z$.

{(\ref{i2lemGlem}).}
Soit $x' = x\vi (y+z)$. Il suffit de voir que $x' \le x\vi z$.
On a $x' \geq 0$, $x' \le x$ donc $x' \perp y$. On peut appliquer
le point précédent à l'in\egt $x' \le y+z$: elle fournit
$x' \le z$, ce que l'on voulait.

{(\ref{i3lemGlem}).}
Conséquence directe du point précédent.

{(\ref{i4lemGlem}).}
 Car $x+y=x\vu y$ et $x\vu y\leq z$.
\end{proof}
%

%:     Corollary{corlemGlem}
\begin{corollary}\label{corlemGlem} Soient $x$, $y$, $z\in G,\, n\in\NN\etl$.
\vspace{-1mm}
\begin{eqnarray}
\label{i1corlemGlem}
&(x=y - z\,\hbox{ et }\,y\geq 0\,\hbox{ et }\,z\geq 0 \,\hbox{ et }\, y\perp z) \,\Longleftrightarrow\, (y=x^+ \,\hbox{ et }\, z=x^-) &\;\;\;\;\;\;\;\\
\label{i11corlemGlem}
&(x\geq 0,\;y\geq 0,\;x\perp y)\;\Longrightarrow\;x\perp ny & \\
\label{i12corlemGlem}
& (nx)^+=nx^+,\; (nx)^-=nx^- ,\; \abs{nx}=n\abs{x}& \\
\label{i2corlemGlem}
&\;nx=0\;\Longrightarrow\;x=0 & \\
\label{i3corlemGlem}
&n(x\vi y)=nx\vi ny,\; \; \;n(x\vu y)=nx\vu ny &
\end{eqnarray}
\end{corollary}
\begin{proof}
{(\ref{i1corlemGlem}).} 
Il reste à montrer le sens $\Longrightarrow$. On a $x^++z=x^-+y$.
En appliquant le lemme de Gauss, on obtient $y\leq x^+$ (parce que $y\perp z$)
et $x^+\leq y$ (parce que $x^+\perp x^-$).

{(\ref{i11corlemGlem}).} 
 Résulte de {(\ref{i2corlemGlem}).} 

{(\ref{i12corlemGlem}).} 
 D'après  (\ref{i1corlemGlem}) et (\ref{i11corlemGlem})
puisque $nx=nx^+-nx^-$ et $nx^+\perp nx^-$.

{(\ref{i2corlemGlem}).} 
 D'après  (\ref{i12corlemGlem}) puisque l'implication est vraie pour $x\geq0$.

{(\ref{i3corlemGlem}).} 
 Les \elts $b=x\vu y$, $a=x\vi y$, $x_1=x-a$ et $y_1=y-a$
sont \cares par les relations suivantes:

\snic{x_1\geq 0,\;y_1\geq 0,\;x=x_1+a,\;y=y_1+a,\;x_1\perp y_1,\;a+b=x+y.}

On multiplie tout par $n$.
\end{proof}

\subsec{Congruences simultanées, principe de recouvrement par quotients}

%:     Definition{defiCongru}
\begin{definition}\label{defiCongru}
Si $a\in G$, on définit la \ixc{congruence modulo $a$}{dans un \grlz} comme suit

\snic{x\equiv y\; \mod\; a \equidef \exists n\in\NN^*,\;\abs{x-y}\leq n\,\abs{a}.}

On note $\cC(a)$ l'ensemble des $x$ congrus à $0$ modulo $a$.
\end{definition}

%:     Fact{factCongru}
\begin{fact}\label{factCongru}
L'ensemble $\cC(a)$ est un sous-\grl de $G$ et les lois du treillis passent au quotient dans $G/\cC(a)$. 
\\
Ainsi, l'application canonique $\pi_a:G\to G/\cC(a)$ est un morphisme de \grlsz,
et tout morphisme de \grls $G\to G'$ qui annule~$a$ se factorise par~$\pi_a$.
\end{fact}

La signification de la congruence $x\equiv 0 \mod a$ est donc que tout morphisme
de \grls $G\vers{\varphi} G'$ qui annule $a$ annule~$x$(\footnote{D'ailleurs, par calcul direct, si
 $\varphi(a)=0$, alors $\varphi(\abs{a})=\abs{\varphi(a)}=0$, et $\abs{\varphi(x)}=\varphi(\abs{x})\leq\varphi(n\abs{a})=n \varphi(\abs{a})=0$, donc~$\varphi(x)=0$.}).

%:     Lemma{lemChinoisGRL}
\medskip Le lemme suivant a la saveur d'un \tho des restes chinois
\ari (voir \thrf{thAnar} point \emph{5}) pour les \grlsz, mais seulement la saveur. Il est nettement plus simple.

\begin{lemma}\label{lemChinoisGRL} \emph{(Lemme des congruences simultanées)}\\
Soient $(\xn)$ dans $G^+$ et $(\an)$ dans  $G$.
\begin{enumerate}
\item Si sont satisfaites les in\egts
$ \abs{a_i - a_j} \leq x_i+x_j,\;i,j\in\lrbn,
$
 il existe un $a\in G$ tel que $ \abs{a - a_i} \leq x_i,\;i\in\lrbn
$. Par ailleurs:
\begin{itemize}
\item Si les $a_i$ sont dans $G^+$ on a une solution $a$ dans $G^+$.
\item Si $\Vi_ix_i=0$, la solution $a$ est unique.
\end{itemize}
\item De même, si 
$a_i\equiv a_j \mod x_i+x_j$ pour $i$, $j\in\lrbn$,
il \hbox{existe  $a\in G$} tel que $a \equiv a_i \mod x_i,\;i\in\lrbn$. Par ailleurs:
\begin{itemize}
\item Si les $a_i$ sont dans $G^+$ on a une solution $a$ dans $G^+$.
\item Si $\Vi_ix_i=0$, la solution $a$ est unique.
\end{itemize}
\end{enumerate}
\end{lemma}
\begin{proof} Il suffit de montrer le point \emph{1.} Voyons d'abord l'unicité.
Si l'on a deux solutions $a$ et $a'$ on aura $\abs{a-a'}\leq2x_i$ pour chaque $i$,
donc $\abs{a-a'}\leq2\Vi_ix_i$.\\
Passons à l'existence. Traitons le cas où les $a_i$ sont dans $G^+$.
Il s'agit en fait de montrer que les hypothèses impliquent l'in\egt
$\Vu_i(a_i-x_i)^+\leq \Vi_i(a_i+x_i)$.
Il suffit pour cela de vérifier que pour chaque $i$, $j$, on a $(a_i-x_i)\vu 0 \leq a_j+x_j$. Or $0\leq a_j+x_j$, et $a_i-x_i \leq a_j+x_j$ par hypothèse.
\end{proof}
%

%:     Lemma{lemRecfermebasique}
\begin{lemma}\label{lemRecfermebasique}
\'Etant donnés une famille finie $(a_j)_{j\in J}$ dans un \grl $G$
et une partie finie $P$ de $J\times J$, il existe une famille finie 
$(x_i)_{i\in I}$ dans~$G$ telle que
\begin{enumerate}
\item $\Vi_{i\in I}x_i=0$.
\item Modulo chacun des $x_i$, pour chaque $(j,k)\in P$, on a $a_j\leq a_k$ ou~$a_k\leq a_j$.
\end{enumerate}
\end{lemma}
%--------- fin lemma ---------------------------------------------- 
%
\begin{proof}
Posons $y_{j,k}=a_j-(a_j\vi a_k)$ et $z_{j,k}=a_k-(a_j\vi a_k)$.
On a $y_{j,k}\vi z_{j,k}=0$. Modulo $y_{j,k}$, on a $a_j=a_j\vi a_k$,
\cad $a_j\leq a_k$, et modulo $z_{j,k}$, on~a~$a_k\leq a_j$. 
\\
En développant par \dit la somme $0=\sum_{(j,k)\in P}(y_{j,k}\vi z_{j,k})$ on obtient un  $\,\Vi_{i\in I}x_i$, où chaque $x_i$ est une somme
$\sum_{j,k}t_{j,k}$, avec pour~$t_{j,k}$ l'un des deux \elts $y_{j,k}$ ou
$z_{j,k}$. Modulo un tel $x_i$ chacun des $t_{j,k}$ est nul (car ils sont $\geq 0$ et leur somme est nulle). On est donc bien dans la situation annoncée.
\end{proof}
%

%:    Principe de recouvrement par quotients
 Le principe ci-après est une sorte d'analogue, pour les \grlsz,
du \plg de base pour les anneaux commutatifs.\iplg

En fait il s'agit d'un simple cas particulier du point \emph{2} du lemme
\ref{lemChinoisGRL} lorsque les $a_i$ sont tous nuls: on applique l'unicité.  

\begin{prvq}\label{prcfgrl} \emph{(Pour les \grlsz)}
Soient $a$, $b\in G$, $x_1$, \ldots, $x_n\in G^+$
avec $\Vi_ix_i=0$. Alors  $a\equiv b \mod x_i$ pour chaque $i$
\ssi $a=b$.
\\
En conséquence, vu le lemme \ref{lemRecfermebasique},
pour démontrer une \egt $a=b$ on peut toujours supposer que les \elts
(en nombre fini) qui se présentent dans un calcul pour une \dem 
de l'\egt sont comparables,
si on en a besoin pour faire la \demz.
Ce principe s'applique aussi bien pour des in\egts que pour des \egts
puisque $a\leq b$ équivaut à $a\vi b=a$.
\end{prvq}

\rem En termes un peu plus abstraits, on aurait pu dire que le morphisme
canonique de \grls $G\to\prod_i G\sur{\cC(x_i)} $ est injectif.
Et le commentaire qui conclut le principe de recouvrement par quotients
peut être paraphrasé comme suit: dans les calculs, un \grl se
comporte toujours comme un produit de groupes totalement ordonnés.
\eoe

\medskip
Dans le \tho de Riesz qui suit on notera que les \gui{il existe} sont des
abréviations pour des formules explicites qui résultent de la \demz.
Ainsi ce \tho peut être vu comme une famille d'\idas dans $G$, sous
certaines conditions de signes (qui sont dans l'hypothèse).
Il est aussi possible de voir ce \tho comme une famille d'\idas \gui{pures} dans $G^+$, \cad
sans aucune condition de signe.
Dans ce cas il faut voir $G^+$ comme une structure \agq pour laquelle on rajoute l'opération
$x \dotdiv y\eqdefi x-(x\vi y)$ (bien définie sur $G^+$
malgré le fait qu'elle fasse appel à l'opération $-$ de $G$).

%:     Theorem{th1GpRtcl}
\begin{theorem}\label{th1GpRtcl} \emph{(\Tho de Riesz)}\\
Soient $G$ un \grl et $u$, $x_1$, \dots, $x_n$, $y_1$, \dots, $y_m$ dans $G^+$.
\begin{enumerate}
\item \label{i1th1GpRtcl}
Si $u\leq \sum_jy_j$, il existe $u_1$, \ldots, $u_m\in G^+$ tels que $u_j\leq y_j$
pour $j\in\lrbm$ et $u=\sum_ju_j$.
\item \label{i2th1GpRtcl}
Si $\sum_ix_i=\sum_jy_j$, il existe
 $(z_{i,j})_{i\in\lrbn,j\in \lrbm} $ dans $G^+$ tels que pour tous $i$, $j$ on ait: $ \sum_{k=1}^m z_{i,k}=x_i$ et $\sum_{\ell=1}^n z_{\ell,j}=y_j$.
\end{enumerate}
\end{theorem}
\begin{Proof}{\Demo \gui{directe}, mais astucieuse. }\\
\emph{\ref{i1th1GpRtcl}.} Il suffit de le prouver pour $m=2$ (\recu facile sur $m$). Si $u\leq y_1+y_2$, \linebreak 
on pose %$y=y_1+y_2$,
$u_1=u\vi y_1$ et $u_2=u-u_1$. Il faut prouver
$0\leq u_2\leq y_2$. \linebreak 
Or $u_2=%u-u_1=
u-(u\vi y_1)=u+\big((-u)\vu (-y_1)\big)=(u-u)\vu (u-y_1)\leq y_2$.

 \emph{\ref{i2th1GpRtcl}.}
Pour $n=1$ ou $m=1$ il n'y a rien à faire.
Pour $n=2$, c'est donné par le point \emph{\ref{i1th1GpRtcl}.}
Supposons donc $n\geq3$. Posons $z_{1,1}=x_1\vi y_1$, $x'_1=x_1-z_{1,1}$
et~$y'_1=y_1-z_{1,1}$. On a $x'_1+x_2+\cdots+x_n=y'_1+y_2+\cdots+y_m $.
\\
Puisque $x'_1\vi y'_1=0$, le lemme de Gauss donne $x'_1\leq y_2+\cdots+y_m$.
\\
Par le point \emph{\ref{i1th1GpRtcl}} on peut écrire $x'_1=z_{1,2}+\cdots+z_{1,m}$
avec les $z_{1,j}\leq y_j$, i.e.~$y_j=z_{1,j}+y'_j$ et $y'_j\in G^+$.
Donc $x_2+\cdots+x_n=y'_1+y'_2+\cdots+y'_m$. 
\\
Ceci nous permet donc une \recu sur~$n$.

\emph{\Demo par le principe de recouvrement par quotients.}\\
Il suffit de prouver le point \emph{2.}
En application du principe \ref{prcfgrl}, on peut supposer le groupe totalement ordonné.
Supposons par exemple $x_1\leq y_1$. On pose $z_{1,1}=x_1$, $z_{1,k}=0$ pour $k\geq2$.
On remplace $y_1$ par $y_1-x_1=y'_1$. On est ramené à résoudre le \pb
pour $x_2$, \ldots, $x_n$ et~$y'_1$, $y_2$, \ldots, $y_m$. De proche en proche, on fait ainsi diminuer $n+m$ jusqu'à ce que  $n=1$ \hbox{ou $m=1$}, auquel cas tout est clair.
\end{Proof}
%

%:     Fact{factGpRtcl}
\begin{fact}\label{factGpRtcl} \emph{(Autres identités dans les \grlsz)}\\
Soient $x$, $y$, $x'$, $y'$, $z$, $t\in G$, $n\in\NN$, $x_1$, \dots, $x_n\in G$.

\vspace{-.10em} 
\begin{enumerate}%\itemsep=1pt
\item \label{i1factGpRtcl} $x+y =\abs{x-y} +2(x\vi y)$
\item $(x\vi y)^+=x^+\vi y^+$, $(x\vi y)^-=x^-\vu y^-$, \\ $(x\vu y)^+=x^+\vu y^+$, $(x\vu y)^-=x^-\vi y^-$.
\item $2(x \vi y)^+ \leq (x+y)^+ \leq x^++y^+$.
\item $\abs{x+y} \leq \abs{x}+\abs{y}\;:\;$
$\abs{x}+\abs{y}=\abs{x+y}+2(x^+\vi y^-) +2( x^-\vi y^+)$.
\item $\abs{x-y} \leq \abs{x}+\abs{y}\;:\;$
$\abs{x}+\abs{y}=\abs{x-y}+2(x^+\vi y^+) +2( x^-\vi y^-)$.
\item $\abs{x+y}\vu\abs{x-y}=\abs{x}+\abs{y}$.
\item $\abs{x+y}\vi\abs{x-y}=\abS{\abs{x}-\abs{y}}$.

\item $\abs{x-y}=(x\vu y)-(x\vi y)$.
\item $\abs{(x\vu z)-(y\vu z)}+\abs{(x\vi z)-(y\vi z)}= \abs{x-y}.$
\item $\abs{x^+ - y^+} + \abs{x^- - y^-} = \abs{x-y}$.
\item \label{i11factGpRtcl} $x\leq z \;\Longrightarrow\; (x\vi y)\vu z= x\vi (y\vu z)$.
\item $x+y=z+t \;\Longrightarrow\; x+y=(x\vu z)+(y\vi t)$.
\item \label{i13factGpRtcl} $n\, x\geq \Vi_{k=1}^n (k y+(n-k)x)  \,\Longrightarrow\,x\geq y$.
\item $\Vu_{i=1}^nx_i =  \sum_{k=1}^n(-1)^{k-1}
        \big(\sum_{I\in \cP_{k,n}}\Vi_{i\in I}x_i\big)$.
\item  $x\perp y\,\Longleftrightarrow\, \abs{x+y}=\abs{x-y}
\,\Longleftrightarrow\, \abs{x+y}=\abs{x}\vu \abs{y}$.
\item  $x\perp y\,\Longrightarrow\, \abs{x+y}=\abs{x}
+\abs{y}=\abs{x}\vu \abs{y}$.
\item \label{i15bisfactGpRtcl} $(x\perp y,\,x'\perp y,\,x\perp y',\,x'\perp y',\,x+y=x'+y') \,\Longrightarrow\, (x=x', \, y=y')$.
\item  \label{i16bisfactGpRtcl} On définit $\Tri(\ux)=[\Tri_1(\ux), \Tri_2(\ux),\ldots,\Tri_n(\ux)]$, où 

\snic{\Tri_k(\xn) =\Vi_{I\in \cP_{k,n}}\left(\Vu_{i\in I}x_i\right) \quad(k\in\lrbn).}

%\sni
On a les résultats suivants.
\begin{enumerate}
\item $\Tri_k(\xn) =\Vu_{J\in \cP_{n-k+1,n}}\big(\Vi_{j\in J}x_j\big)$, $(k\in\lrbn)$.
\item $\Tri_1(\ux)\leq \Tri_2(\ux)\leq \cdots \leq \Tri_n(\ux).$
\item Si les $x_i$ sont deux à deux comparables, la liste $\Tr(\ux)$
est la liste des~$x_i$ rangée en ordre croissant (il n'est pas \ncr que
le groupe soit discret).
\end{enumerate}
\end{enumerate}
Supposons $u$, $v$, $w\in G^+$.
\vspace{-.25em} 
\begin{enumerate}\setcounter{enumi}{18}\itemsep=1pt
\item \label{i17factGpRtcl} $u\perp v\,\Longleftrightarrow\, u+v=\abs{u-v}$.
\item $(u+v)\vi w \leq (u\vi w)+(v\vi w)$.
\item $(x+y)\vu w \leq (x\vu w)+(y\vu w)$.
\item $v\perp w  \,\Longrightarrow\,(u+v)\vi w = u\vi w$.
\item $u\perp v  \,\Longrightarrow\,(u+v)\vi w = (u\vi w)+(v\vi w)$.
\end{enumerate}
\end{fact}
\begin{proof}
Tout ceci est à peu près \imd dans un groupe totalement ordonné,
en raisonnant cas par cas. On conclut avec le principe \ref{prcfgrl}.
\end{proof}
%

%-% ENTRE NOUS
\entrenous{Très probable: les points  \ref{i11factGpRtcl} et \ref{i16bisfactGpRtcl}
sont vrais dans tout \trdiz. Ils pourraient se démontrer au moyen
d'un principe similaire à celui utilisé ici.
}
%-% Fin ENTRENOUS

\rems
\\
1) Une implication comme par exemple

\snic{( u\vi v=0,\ u\geq0,\ v\geq0)\,\Longrightarrow\, u+v=\abs{u-v}}

%\sni
(voir le point \emph{\ref{i17factGpRtcl}})
peut être vue comme le résultat d'une \idt qui exprime,
pour un certain entier $n$, que
$n\abs{u+v-\abs{u-v}}$ est égal à une expression
qui combine $u^-$, $v^-$ et $\abs{u\vi v}$ au moyen des lois $\vu$, $\vi$
et $+$.
En fait,  l'\egt donnée au point \emph{\ref{i1factGpRtcl}} règle directement la question sans  hypothèse de signe sur $u$ et $v$:
$\abs{u+v-\abs{u-v}}=2\abs{u\vi v}$.

2) Il y a une différence importante entre les \idas usuelles, qui sont
en dernière analyse
des \egts entre \pols dans un anneau commutatif libre sur des \idtrsz, $\ZZ[\Xn]$, et les
\idas dans les \grlsz. Ces dernières sont certes des \egts entre expressions que l'on
peut écrire dans un \grl librement engendré par un nombre fini d'\idtrsz,
mais la structure d'un tel \grl libre est nettement
plus difficile à décrypter que celle d'un anneau de \polsz, dans lequel les objets ont une écriture normalisée. La comparaison de deux expressions dans
$\ZZ[\Xn]$ est \gui{facile} dans la mesure où on ramène chacune d'elle
à la forme normale. La t\^ache est beaucoup plus difficile dans les \grls libres,
pour lesquels il n'y a pas de forme normale unique (on peut ramener toute
expression à un sup de inf de \colis des \idtrsz, mais il n'y a pas unicité).
\eoe

%:  subsec{Décomposition partielle, \dcn complète}
\subsec{Décomposition partielle, \dcncz}

%:     Definition{defiDecPar}
\begin{definition}\label{defiDecPar}
Soit $(a_i)_{i\in I}$ une famille finie d'\elts $>0$ dans un \grl discret~$G$.
\begin{enumerate}
\item On dit que cette famille admet une \emph{\dcnpz}
si l'on peut trouver une famille finie $(p_j)_{j\in J}$ d'\elts $>0$ deux à deux \orts telle que chaque $a_i$ s'écrive $\sum_{j\in J}r_{i,j}p_j$
avec les $r_{i,j}\in\NN$.
La famille $(p_j)_{j\in J}$ est alors appelée une \emph{\bdpz} pour la famille~$(a_i)_{i\in I}$.
\item Une telle \dcnp est appelée une \emph{\dcncz} si les
$p_j$ sont \emph{\irdsz} (un \elt $q>0$ est dit \ird si une \egt $q=c+d$ dans $G^+$ implique $c=0$
ou~$d=0$).
\item Un \grl  est dit \emph{à \dcnpz}
 s'il est discret et si toute famille finie d'\elts $>0$ admet une \dcnpz.
\item Un \grl  est dit  \emph{à \dcncz} s'il est discret et si tout \elt $>0$ admet une \dcncz.
\item  Un \grl est dit \emph{à \dcnbz} lorsque pour \hbox{tout $x\geq 0$} il existe un entier $n$ tel que, lorsque $x=\sum_{j=1}^ny_j$ avec les $y_j\geq 0$, au moins l'un des $y_j$ est nul.
\item  Un \grl  est  dit \emph{\noez} si toute suite décroissante d'\elts $\geq 0$ admet deux termes consécutifs égaux.
\end{enumerate}%
\index{decomposition@\dcnz!partielle}%
\index{decomposition@\dcnz!complète}%
\index{decomposition@\dcnz!bornée}%
\index{irreductible@irréductible!\elt --- dans un \grlz}%
\index{noetherien@\noez!groupe réticulé ---}%
\index{groupe réticulé!a decomposition partielle@à décomposition partielle}%
\index{groupe réticulé!a decomposition bornee@à décomposition bornée}%
\index{groupe réticulé!a decomposition complete@à décomposition complète}%
\end{definition}

%\medskip 
\exls ~\\
Une famille vide, ou une famille d'\elts tous nuls, 
 admet la famille vide comme \bdpz.
 \\
Le \grl $\ZZ^{(\NN)}$  est à \dcncz.
\\
 Les
 \grls $\QQ^n$ ($n\geq 1$)  sont à \dcnp mais pas complète.
 \\ 
 Le \grl $\QQ[\sqrt 2]$
 n'est pas à \dcnp (considérer la famille finie $(1,\sqrt2)$). \\
 Le produit lexicographique $\ZZ\times \ZZ$
 n'est pas à \dcnpz. 
 \\
 Plus \gnlt un groupe totalement ordonné à \dcn
 partielle est isomorphe à un sous-groupe de $\QQ$.
\eoe

\medskip
Il est clair qu'un \grl à \dcnc est à \dcnb
et qu'un \grl à \dcnb est 
\noez.

 Dans un \grl à \dcnpz, deux \dcns partielles pour deux
familles finies de $G^+$ admettent un raffinement commun pour la réunion des deux familles: ici on entend qu'une \bdp $(q_1,\ldots,q_s)$ en raffine une autre
si elle est une \bdp pour cette autre.

%:     Proposition{propRaffinementCommun}
\begin{proposition}\label{propRaffinementCommun}
Dans un \grlz, si un \elt $>0$ admet une \dcncz, elle est unique à l'ordre près des facteurs.
\end{proposition}
\begin{proof}
Il suffit de montrer que si un \elt $q$ \ird est majoré par une
somme $\sum_ip_i$ d'\elts \irds il est égal à l'un d'eux.
\\
Or on a alors $q=q\vi \sum_ip_i$, et puisque $q\vi p_j=0$
ou  $p_j$, on peut conclure avec le lemme de Gauss (\egrf{i2lemGlem}).
\\  
Notez que l'on n'a pas besoin de supposer le groupe discret.
\end{proof}
%

%:     proposition{factgrldcntot}
\begin{proposition}\label{factgrldcntot}
Soit $G$ un \grl à \dcncz.
\begin{enumerate}
\item Les \elts \irds de $G^+$ forment une partie détachable $P$, et~$G$ est isomorphe 
à la somme directe \orte $\ZZ^{(P)}$. 
\item Le groupe $G$ est à \dcnb (et a fortiori \noez). 
\end{enumerate}
\end{proposition}

\begin{proof} \emph{1.}
Le test d'irréductibilité est donné par la \dcnc de l'\elt à tester.  
L'isomorphisme s'obtient à partir de l'unicité de la \dcnc (à l'ordre des facteurs près).
\\
\emph{2.} Soit $x\in G^+$. \'Ecrivons $x=\sum_{j\in J}n_jp_j$ avec les $p_j$ \irds et $n_j\in\NN$, et posons $n=\sum_jn_j$.
Alors si $x=\sum_{k=1}^{n+1}x_k$ avec des $x_k\geq 0$ on a \ncrt
l'un des $x_k$ nul (considérer la \dcn de chaque $x_k$ en somme d'\irdsz).
\end{proof}

En \clamaz, un \grl discret \noe
est à \dcncz. Ce résultat
 ne peut pas être obtenu \cotz. Néanmoins on obtient une \dcnpz.

%:     Theorem{th2GpRtcl}
\begin{theorem}\label{th2GpRtcl}
\emph{(Décomposition partielle sous condition \noeez)}
Un \grl $G$ discret et \noe  est à \dcnpz.
\end{theorem}

Pour la \demz, nous utiliserons le lemme suivant.

%--- lemme {decomp2}-----------------
\begin{lemma}\label{decomp2} (sous les hypothèses du \tho \ref{th2GpRtcl})\\
Pour  $a \in G^+$ et
$p_1$, \ldots, $p_m>0$  deux à deux \ortsz, on
peut trouver
des \elts deux à deux \orts $a_0$, $a_1$, \ldots, $a_m$ dans $ G^+ $ 
satisfaisant les \prts suivantes.
%-----------------begin item------------------
\begin{itemize}
\item [{1.}] $ a = \sum_{i=0}^{m}a_{i}.$
\item [{2.}] Pour tout $i \in \lrbm $, il
existe un entier $n_i \geq 0 $ tel que $ a_{i} \leq n_{i}p_i.$
\item [{3.}] Pour tout $ i \in \lrbm$,  $ a_0 \vi  p_i
=0.$
\end{itemize}
%-----------------end item------------------
\end{lemma}
\begin{proof}
Pour chaque $i$, on considère la suite croissante  $(a\vi np_i)_{n\in\NN}$
majorée par~$a$.
Il existe
$n_i$ tel que  $ a\vi n_ip_i =a\vi (n_{i}+1)p_i.$
On prend alors $a_i= a\vi n_ip_i$. Si~$a=a_i+b_i$, on a $b_i\vi p_i=0$
car $a_i\leq a_i+(b_i\vi p_i)\leq a\vi (n_{i}+1)p_i=a_i$.
Les~$a_i$ sont $\leq a$, deux à deux \orts et $\geq0$ donc $a\geq\Vu_ia_i=\sum_ia_i$.
Ainsi, on écrit dans $G^+$ $a = a_1 + \cdots + a_n + a_0$,  avec $a_{i} \leq n_{i}p_i$
pour $i \in \lrbm$.
%Enfin on a pour tout $i$,  $a_0 \vi p_i\leq b_i\vi p_i = 0$.
Enfin, on a $b_i = a_0 + \sum_{j \ne i} a_j$, donc $a_0 \le b_i$, puis
$a_0 \vi p_i \le b_i \vi p_i = 0$. Comme $a_i \le n_ip_i$, on a a fortiori
$a_0 \vi a_i = 0$.
\end{proof}
\begin{Proof}{\Demo du \tho \ref{th2GpRtcl}. }\\
Par \recu sur le nombre $m$ d'\elts de la famille. \\
$\bullet $ Supposons $m = 2 ,$ considérons les \elts $x_{1}$, $x_{2}$.
Pour les besoins de la notation, appelons les $a$, $b$.
Posons $L_1 = \left [a,b \right]$, $m_1=1$,
  $E_{1,a}= [1,0]$, et~$E_{1,b}= [0,1]$.
L'\algo procède par étapes, au début de l'étape $k$
on a un entier naturel $m_k$ et trois listes d'égale longueur: $L_k$,
une liste d'\elts $>0$ de $G$,  $E_{k,a}$ et  $E_{k,b}$ deux
listes d'entiers naturels.
A la fin de l'étape l'entier $m_k$ et les trois listes
sont remplacés par un nouvel entier et de nouvelles listes,
qui servent pour l'étape suivante (à moins que l'algorithme termine).
L'idée \gnle est la suivante: si $x$, $y$ sont deux termes
consécutifs de $L_k$ non \ortsz, on remplace dans $L_k$ le
\linebreak
segment  $(x,y)$ par le segment  $(x-(x \vi y), x \vi y, y-(x
\vi y))$ (en omettant le premier et/ou le dernier terme s'il est nul).
Nous noterons cette
procédure comme suit:

\snic{R:(x,y)\mapsto\;$
le nouveau segment (de longueur 1, 2 ou $3).}

%\sni
Notez que $x+y>\big(x-(x \vi y)\big)+ x \vi y+ \big(y-(x \vi y)\big)$.
\\
Nous devons définir un invariant de boucle. Précisément les
conditions vérifiées par l'entier $m_k$ et les trois listes  sont les
suivantes:
%-----------------begin item------------------
\begin{itemize}
\item $a$ est égal à la \coli des \elts de  $L_k$
affectés des \coes de  $E_{k,a}$,
\item $b$ est égal à la \coli des \elts de  $L_k$
affectés des \coes de  $E_{k,b}$,
\item si $L_k=[x_{k,1},\ldots ,x_{k,r_k}]$ les
\elts $x_{k,j}$ et $x_{k,\ell}$  sont \orts dès que
%-----------------begin item------------------
\begin{itemize}
\item  $j<m_k$ et $\ell\neq j$  ou
\item  $j\geq m_k$ et $\ell\geq j+2$
\end{itemize}
%-----------------end item------------------
\end{itemize}
%-----------------end item------------------
En bref, les $x_{k,j}$ sont deux à deux \ortsz, sauf peut-être certaines 
\linebreak
paires
$(x_{k,j},x_{k,j+1})$ avec $j\geq m_k$.
Ces conditions constituent \emph{l'invariant de boucle}.
Il est clair qu'elles sont (trivialement) vérifiées au départ. 
\\
L'\algo
termine à l'étape $k$ si les \elts de $L_k$
sont deux à deux \ortsz.
En outre, si l'\algo ne termine pas à l'étape $k$, on a
l'in\egt $\sum_{x\in L_k}x \,>\, \sum_{z\in L_{k+1}}z$, donc la condition de
chaîne décroissante assure la terminaison de l'\algoz.\\
Il nous reste à expliquer le déroulement d'une étape et à vérifier
l'invariant de boucle.
Pour ne pas manipuler trop d'indices, nous faisons un léger abus de notation et nous
écrivons
$L_k = \left [p_1,\ldots, p_n \right]$, %$m_k=i$,
$E_{k,a} =\left [\alpha_1,\ldots, \alpha_n \right]$ 
et~$E_{k,b} = \left [\beta_1,\ldots, \beta_n \right]$.\\
Le segment  $(x,y)$  de $L_k$ qui est traité par la procédure
$R(x,y)$ est le suivant:
%si $p_i\vi p_{i+1}\neq 0$ on prend
%$(x,y)=(p_i,p_{i+1})$ sinon
on considère
le plus petit indice $j$ (\ncrt $\geq m_k$) tel
que $p_j\vi p_{j+1}\neq 0$ et l'on prend
  $(x,y)=(p_{j},p_{j+1})$. Si un tel indice n'existe pas, les
\elts de~$L_k$ sont deux à deux \orts et l'\algo est terminé.
Dans le cas contraire on applique la procédure  $R(x,y)$
et l'on met à jour l'entier (on peut prendre $m_{k+1}=j$) et les trois
listes. 
\\
Par exemple en posant $q_j=p_j\vi p_{j+1}$,
$p'_j=p_j-q_j$ et $p'_{j+1}=p_{j+1}-q_j$, si~$p'_j\neq 0\neq
p'_{j+1}$, on aura:
%--------------------begin array---------------
$$\arraycolsep2pt
\begin{array}{rcl}
L_{k+1}& =  &  \left [p_1,\ldots,p_{j-1},p'_j,q_j,p'_{j+1},p_{j+2},\ldots ,
p_n \right]   \\
E_{k+1,a}& =  &  \left [\alpha_1,\ldots, \alpha_{j-1},\alpha_j,\alpha _j+\alpha
_{j+1},\alpha _{j+1},\alpha_{j+2},\ldots \alpha_n \right]   \\
E_{k+1,b}& =  &  \left [\beta_1,\ldots, \beta_{j-1},\beta_j,\beta _j+
\beta _{j+1},\beta _{j+1},\beta_{j+2},\ldots \beta_n \right]
\end{array}
$$
%---------------------end array--------------
On vérifie sans peine dans chacun des 4 cas possibles que l'invariant de
boucle est conservé.

$\bullet $  Si $m > 2$,  par \hdrz, on a pour 
$ (x_1, \ldots, x_{m-1})$  une \bdp  $(p_1,\ldots ,p_n) $.
En appliquant le lemme~\ref{decomp2} à $x_m$ et
$(p_1, \ldots ,p_n) $ on écrit $x_m = \sum_{i=0}^{n}a_i$.\\
Le cas de deux \elts nous donne pour chaque
$(a_i , p_i), \; i \in \lrbn$, une \bdp $S_i$. 
Finalement une \bdp pour $(x_1, \ldots ,x_m)$
est la concaténation des
$ S_i $ et de $ a_0 $.
\end{Proof}

\rem Il est facile de se convaincre que la \bdp calculée par l'\algo
est minimale: toute autre \bdp pour $(x_1, \ldots ,x_m)$ serait
obtenue en décomposant certains des \elts de la base précédente.

%%%%%%%%%%%%%%%%%%%%%%%%%%%%%%%%%%%%%%%%%%%%%%%%%%%%%%%%%%%%%%%%%
%  section{Anneaux à pgcd}
\section{Monoïdes à pgcd, anneaux à pgcd}
\label{secAnnPgcd}

Soit $G$ un \grlz.
Puisque $a\leq b$ \ssi $b\in a + G^+$, la structure d'ordre
est \caree par la donnée du sous-\mo $G^+$.
Vue l'\egt $x=x^+ - x^-$ le groupe  $G$  peut  être
obtenu par symétrisation du \mo $G^+$,
et il revient au même de parler de \grl ou de \mo vérifiant
certaines propriétés particulières (voir \thrf{lem1MonGcd}).

On aurait donc eu de bonnes raisons de commencer par la théorie
des objets du type \gui{partie positive d'un \grlz} plutôt que par
celle des \grlsz.
En effet, dans un \grl la relation d'ordre doit être donnée d'emblée
dans la structure, alors que dans sa partie positive, seule la loi du
\mo intervient, exactement comme dans la théorie multiplicative des
entiers naturels strictement positifs.

C'est donc uniquement des raisons de confort
dans les \dems qui nous ont fait choisir
de commencer par les \grlsz.

%--- subsec{Partie positive d'un \grlz}
\subsec{Partie positive d'un \grlz}

%:     Theorem{lem1MonGcd}
\begin{theorem}\label{lem1MonGcd}
Pour qu'un \mo commutatif $(M,0,+)$
soit la partie positive d'un \grlz,
il faut et suffit que les conditions \ref{i1lem1MonGcd},  \ref{i2lem1MonGcd}
et  \ref{i3lem1MonGcd} ci-dessous soient satisfaites. En outre, on peut
remplacer la condition  \ref{i3lem1MonGcd} par la condition  \ref{i4lem1MonGcd}.
\begin{enumerate}
\item \label{i1lem1MonGcd}
Le \mo est \emph{régulier}, i.e., $x+y=x+z\Rightarrow y=z$.
\item \label{i2lem1MonGcd}
La relation de préordre $x\in y+M$ est une relation d'ordre, 
autrement dit, on a: 
$x+y=0\Rightarrow x=y=0$.\\
On la note $y\leq_M x$, ou si le contexte est clair $y\leq x$.
\item \label{i3lem1MonGcd}
Deux \elts arbitraires admettent une borne supérieure,
i.e.,

\snic{\forall a,b\;\exists c\;\;\uar c=(\uar a)\cap (\uar b).}
\item \label{i4lem1MonGcd}
Deux \elts arbitraires admettent une borne inférieure,
i.e.,

\snic{\forall a,b\;\exists c\;\;\dar c= (\dar a)\cap (\dar b)}
\end{enumerate}
\index{regulier@régulier!monoïde ---}
\end{theorem}
\begin{proof}
A priori la condition \emph{\ref{i3lem1MonGcd}} pour un couple $(a,b)$
particulier est plus forte que
la condition \emph{\ref{i4lem1MonGcd}} pour la raison suivante: si $a$, $b\in M$, l'ensemble des \elts de $M$ inférieurs à $a$ et $b$
est contenu dans $X=\dar(a+b)$.
Sur cet ensemble $X$, l'application $x\mapsto a+b-x$ est une bijection qui renverse
l'ordre et échange donc  sup et inf quand ils existent.
Par contre dans l'autre sens, le inf
dans~$X$ (qui est le inf absolu) peut a priori être transformé
seulement en un sup pour la relation d'ordre restreinte
au sous-ensemble $X$, qui peut ne pas être une borne supérieure
dans $M$ tout entier.\\
Néanmoins, quand la condition \emph{\ref{i4lem1MonGcd}} est vérifiée pour tous $a$, $b\in M$, elle implique la condition~\emph{\ref{i3lem1MonGcd}}.
En effet, montrons que $m = a+b - (a \vi b)$ est le sup de $(a, b)$
dans $M$ en considérant un
$x \in M$ tel que $x \ge a$ et $x \ge b$. Nous voulons montrer que $x \ge m$,
i.e. en posant $y = x \vi m$, que $y \ge m$.
%En effet posons $m=a+b-a\vi b$ et considérons un $x\in M$ tel que
%$x\geq a$ et $x\geq b$. Considérons $y=x\vi m$. Nous voulons montrer $m\leq x$,
%\cad $m\leq y$.
Or $y$ est un majorant de $a$ et $b$, et $y\in X$.
Puisque $m$ est le sup de $a$ et $b$ dans $X$, on a bien $m\leq y$.\\
Le reste de la \dem est laissé \alecz.
\end{proof}

Le \tho précédent conduit à la notion de \mo à pgcd.
Comme cette notion est surtout utilisée
pour le \mo multiplicatif des \elts \ndzs d'un anneau commutatif,
nous passons en notation multiplicative, et nous acceptons que la relation
de divisibilité définie par le \mo ne soit qu'une relation de préordre,
de façon à tenir compte du groupe des unités.

%%%%%%%%%%%%%%%%%%%%%%%%%%%%%%%%%%%%%%%%%%%%%%%%%%%%%%%%%%%%%%%%%
%:--- subsec Monoïdes à pgcd
\subsec{Monoïdes à pgcd}

En notation multiplicative, un \mo commutatif $M$  est régulier lorsque, pour tous $a$, $x$, $y\in M$, l'\egt $ax=ay$ implique $x=y$.

%:     Definition{defiMoGcd}
\begin{definition}\label{defiMoGcd}
On considère un \mo commutatif, noté multiplicativement, $(M,1,\cdot)$.
On dit que \emph{$a$ divise $b$}
lorsque $b\in a\cdot M$, on dit aussi que~$b$ \emph{est multiple de $a$}, et l'on écrit $a\divi b$. Le \mo $M$ est appelé un \emph{\mo à pgcd} lorsque les deux \prts suivantes
sont vérifiées:%
\index{pgcd!mono@\mo à ---}\index{monoide@\moz!a pg@à pgcd}
\begin{enumerate}
\item $M$ est régulier.
\item Deux \elts arbitraires admettent un pgcd, i.e.,

\snic{\forall a,b,\;\exists g,
\;\forall x,
\quad (x \divi a \et x \divi b) \iff x\divi g.}
\end{enumerate}
\end{definition}

On note $U$ le groupe des \elts \ivs (c'est un sous-\moz), encore appelé \emph{groupe des unités}. Deux \elts $a$ et $b$ de $M$ sont dits
\ixc{associés}{elements@\elts --- dans un \moz} s'il existe un
\elt inversible $u$ tel que $ua=b$. Il s'agit d'une relation d'\eqvc (on dit \gui{la relation d'\emph{association}}) et la structure de \mo passe au quotient. On note $M/U$ le \mo quotient. C'est encore un \mo régulier,
et la relation de \dvez, qui était une relation de préordre sur $M$ devient une relation d'ordre sur $M/U$. \index{associ@association}

D'après le  \thrf{lem1MonGcd}, on obtient le résultat suivant.

%:     Theorem{thmoGCD}
\begin{theorem}\label{thmoGCD}
Avec les notations qui précèdent, un \mo commutatif
\ndz $M$ est un \mo à pgcd \ssi $M/U$
est la partie positive d'un \grlz.
\end{theorem}

En notation multiplicative, les \dcnsz, partielles ou complètes, s'appellent des
\emph{\fcnsz}. On parle alors de \emph{\bdfz} au lieu de \bdpz.

De même on utilise la terminologie suivante:
un \mo $M$ \emph{vérifie la condition de chaîne des diviseurs}
si $M/U$ est \noez, \cad si dans toute suite d'\elts $(a_n)_{n\in\NN}$ de $M$ telle
que $a_{k+1}$ divise $a_k$ pour tout $k$, il y a deux termes consécutifs
associés.%
\index{condition de chaîne des diviseurs}

Un \mo est dit \emph{à \fabz} si
$M/U$ est à \dcnbz, \cad si
pour chaque $a$ dans $M$ il existe un entier $n$ tel que
pour toute \fcn
$a=a_1\cdots a_n$ de $a$ dans $M$, l'un des $a_i$ est une unité.
Il est clair qu'un tel \mo vérifie la condition de chaîne des diviseurs.
\index{factorisation!bornée}

%%%%%%%%%%%%%%%%%%%%%%%%%%%%%%%%%%%%%%%%%%%%%%%%%%%%%%%%%%%%%%%%%
%:--- subsection{Anneaux a pgcd}
\subsec{Anneaux à pgcd}
\label{subsecAnnPgcd}

On appelle \emph{anneau à pgcd} un anneau commutatif pour lequel le
\mo multiplicatif des \elts \ndzs est un \mo à pgcd.
On définit de la même manière un anneau  \emph{à \fcn bornée}
ou \emph{qui vérifie la condition de chaîne des diviseurs}.%
\index{anneau!a pgcd@à pgcd}%
\index{factorisation!bornée}\index{anneau!a pgcd a factb@à pgcd à factorisation bornée}%

Un anneau intègre à pgcd pour lequel $\Reg(\gA)/\Ati$ est à \fap
s'appelle un \emph{anneau à pgcd à \fapz}. Rappelons qu'en particulier,
le \grl correspondant doit être discret, ce qui signifie ici que
$\Ati$ doit être une partie détachable de $\Reg(\gA)$.%
\index{factorisation!partielle}\index{anneau!a pgcd a factp@à pgcd à factorisation partielle}

Un anneau intègre à pgcd pour lequel $\Reg(\gA)/\Ati$ est à \fcn complète
s'appelle un \emph{anneau factoriel}. Dans ce cas on parle plutôt de \emph{\facz}.%
\index{factorisation!complète}%
\index{factorisation!totale}%
\index{factoriel}%
\index{anneau!factoriel}

\smallskip
Outre les résultats \gnls sur les \mos à pgcd (qui sont la traduction
en langage multiplicatif des résultats correspondants dans les \grlsz),
on établit quelques faits spécifiques aux anneaux à pgcd, car l'addition
intervient dans les énoncés. Ils pourraient être étendus
aux anneaux \qis sans difficulté.

%:     Fact{factBezGCD}
\begin{fact}\label{factBezGCD}~
\vspace{-.25em} 
\begin{enumerate}\itemsep=1pt%\partopsep=0pt
\item Un anneau intègre à pgcd dont le groupe des unités est détachable
et qui vérifie la condition de chaîne des diviseurs est à \fap
(\thrf{th2GpRtcl}).
\item Un anneau de Bézout est un anneau à pgcd.
\item Un anneau principal est un anneau à pgcd intègre
qui vérifie la condition de chaîne des diviseurs. 
Si le groupe des unités est détachable,
l'anneau est à \fapz.
\item Si $\gK$ est un \cdi non trivial,
$\KX$ est un anneau de Bézout, à \fabz, et le groupe des unités est détachable. En particulier, l'anneau $\KX$ est à \fapz.
\item Les anneaux $\ZZ$, $\ZZ[X]$ et $\QQ[X]$ sont  factoriels (proposition~\ref{propZXfactor}).
%
%\item
%
\end{enumerate}
\end{fact}
\facile

%:     Theorem{thAgcdNormal}
\begin{theorem}\label{thAgcdNormal}
Tout anneau à pgcd intègre est intégralement clos.
\end{theorem}
\begin{proof}
La preuve du lemme \ref{lemZintClos} peut être reprise mot à mot.
\end{proof}

Nous laissons \alec la \dem des  faits suivants (pour \ref{factAXiclgcd}
il faut utiliser le \tho de \KROz).
%:     Fact{factLocaliseGCD}
\begin{fact}\label{factLocaliseGCD}
Soit $\gA$  un anneau intègre à pgcd et  $S$ 
un \moz. Alors~$\gA_S$ est un anneau intègre à pgcd, et pour $a$, $b\in\gA$ un pgcd dans~$\gA$ est un pgcd dans~$\gA_S$.
\end{fact}

Nous dirons
qu'un sous-\mo $V$ d'un \mo $S$ est \emph{saturé} (dans $S$) si
$xy\in V$ et  $x$, $y\in S$ impliquent $x\in V$.
Dans la littérature, on trouve aussi: \emph{$V$ est factoriellement
clos dans $S$}.
Un \mo $V$ d'un anneau commutatif $\gA$ est donc saturé \ssi
il est saturé dans le \mo multiplicatif~$\gA$.%
\index{monoïde!saturé dans un autre}%
\index{saturé!sous-monoïde ---}%
\index{factoriellement clos!sous-monoïde --- }

%:     Fact{factSouspgcdsat}
\begin{fact}\label{factSouspgcdsat}
Un sous-\mo saturé $V$ d'un \mo à pgcd (resp.\,à \fcn bornée) $S$
est un \mo à pgcd (resp.\,à \fabz) avec les mêmes pgcd et ppcm
que dans $S$.
\end{fact}

%:     Fact{factAXiclgcd}
\begin{fact}\label{factAXiclgcd}~
\\
 Soit $\gA$ un anneau \icl  non trivial et $\gK$ son corps de fractions.
\\
Le \mo multiplicatif des \polus de $\AuX=\AXn$  s'identifie naturellement à un
sous-\mo saturé de $\KuX\etl\!\sur\gK\!\eti$.
\\
En particulier, le \mo multiplicatif des \polus de $\AuX$ 
est un \mo à pgcd à \fabz.
\end{fact}

%-% ENTRE NOUS
\entrenous{Sans doute le même résultat pour le \mo multiplicatif des \pols 
primitifs. En exo? 
}
%-% Fin ENTRENOUS

%%%%%%%%%%%%%%%%%%%%%%%%%%%%%%%%%%%%%%%%%%%%%%%%%%%%%%%%%%%%%%%%%
%: \subsubsec{Anneaux à pgcd de dimension $\leq1$}
\subsubsec{Anneaux à pgcd de dimension $\leq1$}

%:     Definition{defiDim1Integre}
\begin{definition}\label{defiDim1Integre}
Un anneau \qi $\gA$ est dit \emph{de dimension $\leq1$} si pour tout
\elt $a$ \ndz le quotient $\aqo{\gA}{a}$ est \zedz.
\index{dimension (de Krull) $\leq1$!anneau \qi de ---}
\end{definition}

\rem Sous l'hypothèse que $a$ est \ndzz, nous obtenons donc que pour tout
$b$, il existe $x$, $y\in\gA$ et $n\in\NN$ tels que

\snic{\;\;\qquad\qquad b^n(1+bx)+ay=0.\qquad\qquad\;\;(*)}

%\sni
Si nous ne faisons plus d'hypothèse sur $a$, nous pouvons considérer
l'\idm $e$ qui engendre $\Ann(a)$, et nous avons alors une \egt du type $(*)$,
mais en remplaçant $a$ par $a+e$, qui est \ndzz. Cette \egt donne,
après une multiplication par $a$ qui fait disparaître $e$:

\snic{\qquad\qquad a(b^n(1+bx)+ay)=0\qquad\qquad(+).}

%\sni
Nous obtenons ainsi une \egt conforme à celle donnée dans le chapitre~\ref{chapKrulldim} où apparaît une \dfn \cov de
la phrase \gui{$\gA$ est un anneau de \ddk $\leq 1$},
pour un anneau $\gA$
arbitraire (voir le point~\emph{\iref{i3corKrull}} de la proposition~\ref{corKrull}).
\eoe

%:     Lemma{lemDim1-1}
\begin{lemma}\label{lemDim1-1}\relax
\emph{(Une \fcn en dimension $1$)}
\begin{enumerate}
\item Soit dans un anneau $\gA$ deux \ids $\fa$, $\fb$ avec $\gA\sur\fa$ \zed et
$\fb$ \tfz.
 Alors on peut écrire:  

\snic{\fa=\fa_1\fa_2\;$ avec $\;\fa_1+\fb=\gen{1}$ et $\,\fb^n\subseteq \fa_2}

%\sni
pour un entier $n$ convenable. Cette
écriture est unique et l'on a 

\snic{\fa_1+\fa_2=\gen{1},$
   $\;\fa_2=\fa+\fb^n=\fa+\fb^{m}$ pour tout $m\geq n .}  

\item Le résultat s'applique si $\gA$ est \qi de dimension $\leq1$, $\fa$ est \ivz, et $\fb$ \tfz.
Dans ce cas $\fa_1$ et $\fa_2$ sont \ivsz.
En particulier, $\fa+\fb^n$ est \iv pour $n$ assez grand.
\end{enumerate}

\end{lemma}
%----------  end lemma -----------------------------
\begin{proof}
 Il suffit de prouver le point \emph{1.}
\\ 
\emph{Existence et unicité de la \fcnz}.
On considère un  triplet $(\fa_1,\fa_2,n)$ susceptible de vérifier les hypothèses. Puisque $\fa_1$ et $\fa_2$ doivent contenir $\fa$, on peut raisonner modulo $\fa$, et donc supposer $\gA$ \zed avec l'\egt $\fa_1\fa_2=\gen{0}$.
\\
 Notons que $\fa_1+\fb=\gen{1}$ implique $\fa_1+\fb^{\ell}=\gen{1}$ pour tout exposant $\ell\geq1$. En particulier, $\gA=\fa_1\oplus\fa_2=\fa_1\oplus \fb^m$ pour tout $m\geq n$. Ceci force, avec $e$ \idmz,
$\fa_1=\gen{1-e}$ et $\fa_2=\fb^m=\gen{e}$ pour $m$ tel que $\fb^m=\fb^{m+1}$
(voir le lemme \ref{lemfacile} et le point \emph{\iref{LID003}} du lemme \ref{lemme:idempotentDimension0}).
\end{proof}

\rem Le point \emph{2} est valable sans supposer $\gA$ \qiz. Cela deviendra
clair après la \dfn \cov \gnle de la dimension de Krull, puisque pour tout
\elt \ndz $a$, si $\gA$ est de dimension $\leq 1$, l'anneau $\aqo\gA a$ est \zedz. 
\eoe

%:     proposition{lemGCDLop}
\begin{proposition}\label{lemGCDLop}
Soit $\gA$ un anneau intègre à pgcd; alors
tout \id \lop est principal.
\end{proposition}
\begin{proof}
Soit $\fa = \gen {\an}$ \lop et $d = \pgcd(\an)$.  Montrons que $\fa = \gen
{d}$.  Il existe un \sys d'\eco  $(s_1, \ldots, s_n)$  avec $\gen {\an} = \gen
{a_i}$ dans $\gA_{s_i}$. Il suffit de voir que $\gen {\an} = \gen {d}$ dans chaque $\gA_{s_i}$ car cette \egtz,  vraie localement,  le sera globalement.
Or~$\gA_{s_i}$ reste un
anneau à pgcd, et les pgcds ne chan\-gent pas. Donc, dans~$\gA_{s_i}$, 
on obtient $\gen {\an} = \gen {a_i} = \gen {\pgcd(\an)} = \gen {d}$.
\end{proof}
%

%:     theorem{propGCDDim1}
\begin{theorem}\label{propGCDDim1}
Un anneau intègre à pgcd de dimension $\leq 1$ est un anneau de Bézout.
\end{theorem}
\begin{proof}
Puisque $\gen{a,b} = g \gen{a_1,b_1}$ avec $\pgcd(a_1,b_1) = 1$, il suffit de
montrer que $\pgcd(a,b) = 1$ implique $\gen{a,b} = \gen{1}$.  Or $\pgcd(a,b) =
1$ implique $\pgcd(a,b^n) = 1$ pour tout $n \ge 0$.  Enfin d'après le point
\emph {2} du lemme~\ref{lemDim1-1}, pour $n$ assez grand, $\gen {a, b^n}$ est
\iv donc \lopz,  et on conclut avec la proposition~\ref{lemGCDLop}.
\end{proof}
%

%%%%%%%%%%%%%%%%%%%%%%%%%%%%%%%%%%%%%%%%%%%%%%%%%%%%%%%%%%%%%%%%%
%: \subsubsec{Pgcd dans un anneau de \polsz}
\subsubsec{Pgcd dans un anneau de \polsz}

Si $\gA$ est un anneau à pgcd intègre et $f\in\AX$ on note $\G_X(f)$ ou $\G(f)$ un pgcd
des \coes de $f$ (il est défini à une unité près multiplicativement)
et on l'appelle le \ix{G-contenu} de $f$.
Un \pol dont le G-contenu est égal à $1$ est dit \ix{G-primitif}.

%:     Lemma{lemGcont}
\begin{lemma}\label{lemGcont}
Soit $\gA$ un anneau à pgcd intègre de corps de fractions
$\gK$ et~$f $ un \elt non nul de $\KX$.
\vspace{-.25em} 
\begin{itemize}\itemsep=0pt%\partopsep=0pt
\item On peut écrire $f=af_1$ avec $a\in\gK$
et $f_1$  G-primitif dans $\AX$.
\item Cette écriture est unique au sens suivant:
pour une autre écriture du même type $f=a'f_1'$,
il existe $u\in\Ati$ tel que $a'=ua$ et $f_1=uf_1'$.
\item  $f\in\AX$ \ssi $a\in\gA$, dans ce cas $a=\G(f)$.
\end{itemize}
\end{lemma}
\facile

%:   Proposition{propLG} ---------------
\begin{proposition}
\label{propLG} \emph{(Lemme de Gauss, un autre)} 
Soit $\gA$ un anneau à pgcd intègre et $f$, $g\in\AX$. Alors
$\G(fg)=\G(f)\G(g)$. En particulier, le produit de deux \pols G-primitifs
est un \polz~\hbox{G-primitif}.
\end{proposition}
%--- end-proposition----------------------------------------
%
\begin{proof} Notons $f_i$ et $g_j$ les \coes de $f$ et $g$.
Il est clair que $\G(f)\G(g)$ divise $\G(fg)$.
Par \dit le pgcd des $f_ig_j$ est
égal à $\G(f)\G(g)$, or  la proposition \ref{propArm} implique que
$\G(fg)$  divise les $f_ig_j$ donc leur pgcd.
\end{proof}
%

%:     Corollary{corpropLG}
\begin{corollary}\label{corpropLG}
Soit $\gA$ un anneau à pgcd intègre de corps de fractions
$\gK$ et $f$, $g\in\AX$. Alors $f$ divise $g$
dans $\AX$ \ssi $f$ divise~$g$ dans $\KX$ et $\G(f)$ divise $\G(g)$
dans $\gA$.
\end{corollary}
\begin{proof}
Le \gui{seulement si} résulte du lemme de Gauss.
Pour le \gui{si} nous pouvons supposer que $f$ est G-primitif.
Si $g=hf$ dans $\KX$, nous pouvons écrire $h=ah_1$ où $h_1\in\AX$
est G-primitif et $a\in\gK$. Par le lemme de Gauss, on a $fh_1$ G-primitif.
En appliquant le lemme \ref{lemGcont} à l'\egt $g = a(h_1f)$,
on obtient $a \in \gA$, puis $h\in\AX$.
\end{proof}
%

%:     Fact{factUnitesAX}
%\begin{fact}\label{factUnitesAX}

Rappelons que si $\gA$ est un anneau réduit, $\AX\eti=\Ati$
(lemme \ref{lemGaussJoyal} \emph{\iref{i4lemPrimitf}}).
En particulier, si $\gA$ est intègre non trivial
et si le groupe des unités de
$\gA$ est détachable,  il en va de même pour~$\AX$.

%:     Theorem{thAXgcd}
\begin{theorem}\label{thAXgcd}
Soit $\gA$ un anneau à pgcd intègre, de corps de fractions~$\gK$.
\begin{enumerate}
\item $\AXn$ est un anneau à pgcd intègre.
\item Si $\gA$ est à \fapz, il en va de même pour~$\AX$.
\item Si $\gA$ vérifie la condition de chaîne des diviseurs,
il en va de même pour~$\AX$.
\item Si $\gA$ est à \fabz, il en va de même pour~$\AX$.
\item Si $\AX$ est factoriel, il en va de même pour~$\AXn$ \emph{(Kronecker)}.
\end{enumerate}
\end{theorem}
\begin{proof}
\emph{1.} Il suffit de traiter le cas $n=1$. 
Soient $f$, $g\in\AX$. \\
\'Ecrivons $f=af_1$, $g=bg_1$, avec $f_1$ et $g_1$
G-primitifs. Soit $c=\pgcd_\gA(a,b)$ \hbox{et $h=\pgcd_\KX(f_1,g_1)$}. Nous pouvons
supposer \spdg \hbox{que $h$} est dans $\AX$ et qu'il est G-primitif. Alors,
en utilisant le corolaire~\ref{corpropLG}, on vérifie  que $ch$ est un pgcd de $f$ et $g$ dans $\AX.$
 
Les points \emph{2}, \emph{3} et \emph{4} sont laissés \alecz.
 
\emph{5.} Il suffit de traiter le cas $n=2$
et de savoir détecter si un \pol admet un facteur strict. 
On utilise l'astuce de \KRAz.
Pour tester le \polz~$f(X,Y)\in\gA[X,Y]$, supposé de degré $< d$ en $X$, 
on considère le \polz~$g(X)=f(X,X^d)$. Une \dcnc de $g(X)$ permet de savoir s'il existe un facteur strict de $g$ de la forme $h(X,X^d)$ (en regardant tous les facteurs stricts de $g$, à association près), ce qui correspond à un facteur strict de $f$. Pour quelques précisions voir l'exercice \ref{exoKroneckerTrick}.
\end{proof}
%

%:     Corollary{corthAXgcd}
\begin{corollary}\label{corthAXgcd}
Si $\gK$ est un \cdi non trivial, $\KXn$ est un anneau intègre à pgcd,
 à \fab et à \fapz. Le groupe des unités est $\gK\eti$. Enfin $\KXn$ ($n\geq 2$) est
 factoriel \ssi $\KX$ est factoriel. 

\end{corollary}

%%%%%%%%%%%%%%%%%%%%%%%%%%%%%%%%%%%%%%%%%%%%%%%%%%%%%%%%%%%%%%%%%
\section{Treillis de Zariski d'un anneau commutatif}
\label{secZarAcom}

\vspace{4pt}
%--- subsec{Généralités}
\subsec{Généralités}

Nous rappelons la notation $\DA(\fa)$ avec quelques précisions.
%: Notation{notaZA}--------------
\begin{notation}
\label{notaZA}\label{defZar}
{\rm  Si $\fa$ est un \id de $\gA$ on note $\DA(\fa)=\sqrt{\fa}$ le nilradical de~$\fa$. Si $\fa=\gen{x_1,\ldots ,x_n}$ on note $\DA(x_1,\ldots ,x_n)$ pour
$\DA(\fa)$. On \hbox{note $\ZarA$} l'ensemble des  $\DA(x_1,\ldots ,x_n)$ (pour $n\in\NN$
et $x_1$, \ldots, $x_n\in\gA$).
}
\end{notation}
%--- end-notation-----------------------------------------

On a donc $x\in\DA(x_1,\ldots ,x_n)$ \ssi une puissance de $x$ appartient à $\gen{x_1,\ldots ,x_n}$.

L'ensemble $\ZarA$ est ordonné par la relation d'inclusion.

%:    Fact{factZar}-------
\begin{fact}
\label{factZar}
$\ZarA$ est un \trdi avec
%--------------------begin array---------------
$$\begin{array}{rcl}
\DA(0)=0_{\ZarA}, && \DA(\fa_1)\vu\DA(\fa_2)=\DA(\fa_1+\fa_2), \\
\DA(1)=1_{\ZarA},  &\quad   &
\DA(\fa_1)\vi\DA(\fa_2)=\DA(\fa_1\,\fa_2).
\end{array}$$
%---------------------end array--------------
On l'appelle le \emph{treillis de Zariski de l'anneau $\gA$}.%
\index{treillis!de Zariski}
\end{fact}
%--- end-fact-----------------------------------------

En \clama $\DA(x_1,\ldots ,x_n)$ peut être vu comme un \oqc de $\SpecA$: l'ensemble des \ideps $\fp$ de $\gA$ tels que l'un au moins des $x_i$ n'appartienne pas à $\fp$. Et $\ZarA$ s'identifie au treillis des \oqcs de $\SpecA$.
Pour plus de détails sur ce sujet voir la section~\ref{secEspSpectraux}.

%:     Fact{fact1Zar}
\begin{fact}\label{fact1Zar}~
\begin{enumerate}
\item Pour tout morphisme $\varphi:\gA\to\gB$, on a un morphisme naturel
 $\Zar\varphi$ \hbox{de $\ZarA$} \hbox{vers $\Zar\gB$}, et l'on obtient ainsi un foncteur de la catégorie
des anneaux commutatifs vers celle des \trdisz.
\item Pour tout anneau $\gA$
l'\homo naturel $\ZarA\to\Zar\Ared$ est un \isoz,
de sorte que l'on peut identifier les deux treillis.
\item L'\homo naturel $\Zar(\gA_1\times \gA_2)\to\Zar\gA_1\times \Zar\gA_2$ est un \isoz.
\item Pour une \agB $\gB$, l'application $x\mapsto \DB(x)$
est un \iso de $\gB$ sur $\Zar\,\gB$.
\end{enumerate}
\end{fact}

%:     Fact{factZarABol}
\begin{fact}\label{factZarABol}
\Propeq
\begin{enumerate}
\item $\ZarA$ est une \agBz.
\item $\gA$ est \zedz.
\end{enumerate}
\end{fact}

\begin{proof} Rappelons qu'un \trdi \gui{est} une \agB \ssi tout \elt admet un complément
(proposition~\ref{defiBooleTrdi}).
\\
Supposons \emph{2.} Alors pour tout \itf $\fa$, il existe un \idm $e$
et un entier $n$ tels que $\fa^n=\gen{e}$. Donc $\DA(\fa)=\DA(e)$.
Par ailleurs, il est clair que $\DA(e)$ et $\DA(1-e)$ sont
complé\-mentaires dans $\ZarA$.
\\
Supposons \emph{1.} Soit $x\in\gA$ et $\fa$ un \itf de $\gA$ tel que $\DA(\fa)$
soit le complément de $\DA(x)$ dans $\ZarA$. Alors il existe $b\in\gA$
et $a\in\fa$ tels \hbox{que $bx+a=1$}. Comme $xa=x(1-bx)$ est nilpotent on obtient
une \hbox{\egt $x^n(1+cx)=0$}.
\end{proof}
%%%%%%%%%%%%%%%%%%%%%%%%%%%%%%%%%%%%%%%%%%%%%%%%%%%%%%%%%%%%%%%%%

%:     Fact{fact2Zar}
\begin{fact}\label{fact2Zar} \label{corthfactorisation2}
Soient $a\in\gA$ et $\fa\in\ZarA$.
\begin{enumerate}
\item L'\homo $\Zar\pi:\ZarA\to\Zar(\aqo{\gA}{a})$,
où $\pi:\gA\to \aqo{\gA}{a}$ est la projection canonique,
est surjectif,
et il permet d'identifier $\Zar(\aqo{\gA}{a})$ au treillis
quotient $\Zar(\gA)\sur{(\DA(a)=0)}$.
%
%\item
Plus \gnltz,
 $\Zar(\gA\sur{\fa})$ s'identifie à $\Zar(\gA)\sur{(\fa=0)}$.
\item  L'\homo $\Zar j:\ZarA\to\Zar(\gA[1/a])$,
où $j:\gA\to \gA[1/a]$ est l'\homo canonique,
est surjectif
et il permet d'identifier~$\Zar(\gA[1/a])$
au treillis quotient $\Zar(\gA)\sur{(\DA(a)=1)}$.
\item  Pour un \id $\fc$   et  un \mo $S$ de  $\gA$ on a un \iso naturel

\snic{\Zar(\gA_S\sur{\fc\gA_S})\;\simeq\; \Zar(\gA)\sur{(\fb=0,\ff=1)},}

%\sni
où $\fb$ est l'\id de $\ZarA$ engendré par les $\DA(c)$ pour $c\in\fc$,
et $\ff$ est le filtre de $\ZarA$ engendré par les $\DA(s)$ pour $s\in S$.
\end{enumerate}
\end{fact}

%Pour le dernier point on utilise le \tho de \fcn donné dans le fait~\ref{fact2Zar}.

%:--- Subsec{Dualité dans les anneaux commutatifs}--  secIDEFIL
\subsec{Dualité dans les anneaux commutatifs}\label{secIDEFIL} \perso{Je ne sais vraiment
pas où doit se trouver cette section dans le livre.}

%:   subsubsec{Annuler et inverser simultanément}
\subsubsec{Annuler et inverser simultanément}
\label{secAnnEtInv}

%Idéaux et filtres d'un anneau commutatif sont en rapport étroit
%avec les \ids et filtres de son treillis de Zariski.
Dans les \trdis on échange les rôles de $\vi$ et $\vu$ en passant au treillis
opposé, \cad en renversant la relation d'ordre.

Dans les anneaux commutatifs, une dualité féconde existe aussi entre l'addition et la multiplication,
plus mystérieuse, lorsque l'on  essaie d'échanger leurs rôles.

%Dans cette section les démonstrations sont laissées en exercice.

Rappelons qu'un \mo saturé est appelé un \emph{filtre}.
La notion de filtre est une notion duale de celle d'idéal, tout aussi
importante.

Les idéaux sont les images réciproques de $0$ par les \homosz, ils servent à passer au quotient,
\cad à annuler des \elts par force. Les filtres sont les
images réciproques du groupe des unités par les \homosz,
ils servent à localiser, \cad à rendre des \elts inversibles par force.

\'Etant donnés un idéal $\fa$ et un \mo $S$ de l'anneau $\gA$ on peut vouloir annuler les \elts de $\fa$ et inverser les \elts de $S$.
La solution de ce \pb est donnée par la considération de l'anneau
suivant.

%:     Definota{defiASa}
\begin{definota}\label{defiASa}
On note  (par abus) $\gA_{S}\sur{\fa}$ ou $S^{-1}\gA\sur{\fa}$ l'anneau
dont les \elts
sont donnés par les couples $(a,s)\in\gA\times S$,
 avec l'\egt $(a,s)=(a',s')$ dans
$\gA_{S}\sur{\fa}$ \ssi
il existe $s''\in S$  tels que $s''(as'-a's)\in \fa$
(on notera $a/s$ pour le couple $(a,s)$).
\end{definota}

Le fait que  $\gA_{S}\sur{\fa}$ ainsi défini répond au \pb posé
signifie que le \tho de \fcn suivant est vrai
(voir les faits analogues \ref{factUnivQuot} et~\ref{factUnivLoc}).

%:--- Fact{thfactorisation2}}-----
\begin{fact}
\label{thfactorisation2} \emph{(\Tho de factorisation)}\\
Avec les notations ci-dessus, soit
$\psi:\gA\to\gB$ un \homoz. Alors $\psi$ se factorise par $\gA_S\sur{\fa}$
\ssi $\psi(\fa)\subseteq \so{0}$ et $\psi(S)\subseteq \gB^\times$. Dans ce
cas, la \fcn est unique.
\end{fact}
%--- end-theorem-----------------------------------------

\pun{\gA}{\lambda}{\psi}{\gA_S\sur{\fa}}{\theta}{\gB}{$\psi(\fa)\subseteq \so{0}$ et $\psi(S)\subseteq \gB^\times$}

Naturellement on peut aussi résoudre le \pb en annulant d'abord $\fa$ puis en inversant (l'image de) $S$, ou bien en inversant d'abord $S$ puis
en annulant (l'image de) $\fa$. On obtient ainsi des \isos canoniques
$$
\gA_{S}\sur{\fa} \;\simeq \;\big(\pi_{\gA,\fa}(S)\big)^{-1}(\gA\sur{\fa})
\;\simeq\; (\gA_{S})\sur {\left(j_{\gA,S}(\fa)\gA_{S}\right)}.
$$

%:   subsec{Des \dfns duales} %%%%%%%%%%%%
\subsubsec{Des \dfns duales}

La dualité  entre idéaux et filtres est une forme de dualité entre l'addition et la multiplication.

Ceci se voit bien sur les axiomes respectifs qui servent à définir les
idéaux (resp. \idepsz) et les filtres (resp. filtres premiers):
%--------------------begin array---------------
$$\arraycolsep3pt\begin{array}{rclcrcl}
     \mathrm{Id\acute eal}   &\fa   &   & \qquad\qquad   & \mathrm{Filtre}
& \ff & \\[1mm]
       &\vdash   &  0\in\fa & \quad   &   & \vdash & 1\in\ff\\
   x\in\fa, \, y\in\fa     &\vdash   & x+y\in\fa  & \quad
  &x\in\ff, \, y\in\ff   & \vdash &  xy\in\ff\\
   x\in\fa      &\vdash   & xy\in\fa  & \quad
& xy\in\ff  & \vdash & x\in\ff\\[1mm]
     \mathrm{premier}   &  &   & \quad   & \mathrm{premier}  &   & \\   
xy\in\fa &\vdash   & x\in\fa \,\hbox{ ou }\, y\in\fa  & \quad   &   
x+y\in\ff&\vdash   & x\in\ff \,\hbox{ ou }\, y\in\ff
\end{array}$$
%---------------------end array--------------

Notez que selon la \dfn ci-dessus,  $\gA$ est à la fois un \idep et un
filtre premier de $\gA$.
Cette convention peut paraître étrange, mais il s'avère que c'est celle
qui est la plus pratique: un idéal est premier \ssi l'anneau quotient est \sdzz,
un filtre est premier \ssi le localisé est un \aloz.
Pour ce qui concerne les \ids nous nous sommes déjà expliqués
à ce sujet dans le commentaire \paref{CommIdeps}.
\index{filtre!premier}
\index{premier!filtre ---}

Nous adopterons la  \dfn suivante pour un \emph{filtre maximal}: le localisé est un \alo \zed (lorsque l'anneau est réduit: un \cdiz). 
En particulier, tout filtre maximal est premier.
Nous ferons usage de cette \dfn essentiellement à titre heuristique.
\index{filtre!maximal}\label{labfima}
\index{maximal!filtre ---}

Supposons maintenant l'anneau $\gA$ non trivial.
Alors un \id  strict détachable (resp. un filtre  strict détachable) est premier  \ssi son
\cop est un filtre (resp. un idéal).
Nous retrouvons dans ce cas le terrain familier en \clamaz.

De manière \gnle en \clama le \cop d'un \idep strict est un filtre premier strict
et vice versa, donc le \cop \idema strict est un filtre premier minimal,
et le \cop d'un \fima strict est un \idemiz. Les filtres premiers paraissent
donc plus ou moins inutiles et ont tendance à disparaître de la scène.

%\smallskip
%Rappelons que nous avons discuté nos ennuis terminologiques
%avec les \ideps et les \idemas dans le commentaire \paref{CommIdeps}.

%:   subsec{Couples saturés} %%%%%%%%%%%%
\subsubsec{Couples saturés}

Une bonne façon de comprendre la dualité est de traiter
simultanément idéaux et filtres.
Pour ceci nous introduisons
la notion de \emph{\paz}, analogue à celle que nous avons donnée
pour les \trdisz.

%--- Definition{defpa}-----------
\begin{definition}
\label{defpa}
Soient $\fa$  un idéal et $\ff$  un filtre.
On dit que $\fa$ est \emph{$\ff$-saturé} si l'on~a:
$$
(as\in\fa,\,s\in\ff)\Longrightarrow a\in\fa,
$$
on dit que $\ff$ est \emph{$\fa$-saturé} si l'on~a:
$$
(a+s\in\ff,\,a\in\fa)\Longrightarrow s\in\ff,
$$
si $\fa$ est $\ff$-saturé et $\ff$ est $\fa$-saturé on dit que $(\fa,\ff)$
est un \emph{\paz} dans~$\gA$.%
\index{saturé!couple ---}%
\index{saturé!idéal $\ff$- ---}%
\index{saturé!filtre $\fa$- ---}%
\end{definition}
%--- end-definition------------------------------------

Récapitulons les axiomes pour les \pas (notez que la dernière condition se
réécrit $\fa+\ff=\ff$).
%--------------------begin array---------------
$$\arraycolsep3pt\begin{array}{rclcrcl}
       &\vdash   &  0\in\fa & \quad   &   & \vdash & 1\in\ff\\
   x\in\fa, \, y\in\fa     &\vdash   & x+y\in\fa  &
  &x\in\ff, \, y\in\ff   & \vdash &  xy\in\ff\\
   x\in\fa      &\vdash   & xy\in\fa  &
& xy\in\ff  & \vdash & x\in\ff\\
   xy\in\fa,\,y\in\ff      &\vdash   & x\in\fa  &
& x+y\in\ff,\,y\in\fa  & \vdash & x\in\ff\\
\end{array}$$
%---------------------end array--------------

%--- Fact{factPaire}----------
\begin{fact}
\label{factPaire}~
\begin{enumerate}

\item Pour tout \homo $\varphi:\gA\to\gB$, le couple
$\big(\Ker\varphi,\varphi^{-1}(\gB^{\times})\big)$ est un \paz.

\item Réciproquement si $(\fa,\ff)$ est un \pa et si $\psi:\gA\to\gA_{\ff}\sur{\fa}=\gC$ désigne l'\homo canonique, on a $\Ker\psi=\fa$ et
$\psi^{-1}(\gC^{\times})=\ff$.

\item Soit $\varphi:\gA\to\gC$ un \homo et $(\fb,\ffg)$ un \pa de $\gC$, alors
$\big(\varphi^{-1}(\fb),\varphi^{-1}(\ffg)\big)$ est un \pa de $\gA$.

\end{enumerate}
\end{fact}
%--- end-fact-----------------------------------------

%--- Fact{actPaire2}-----------------
\goodbreak
\begin{fact}
\label{factPaire2}
Soit $(\fa,\ff)$  un \paz.
%Supposons $\fa\cap\ff=\emptyset$.
%-----------------begin enum------------------
\begin{enumerate}
\item  $\gA_{\ff}\sur{\fa}$ est local \ssi $\ff$ est un filtre premier
(\cad \ssi $\gA_{\ff}$ est local).
\item  $\gA_{\ff}\sur{\fa}$ est \sdz \ssi $\fa$ est un \id premier
(\cad \ssi $\gA\sur{\fa}$ est \sdzz).
\end{enumerate}
%-----------------end enum------------------
\end{fact}
%--- end-fact-----------------------------------------

%--- Definition{defRaffine}--------
\begin{definition}
\label{defRaffine}
Si $(\fa,\ff)$ et $(\fb,\ffg)$ sont deux \pas de $\gA$ on dit que
\emph{$(\fb,\ffg)$ raffine  $(\fa,\ff)$}
 et l'on écrit
$(\fa,\ff)\leq (\fb,\ffg)$ lorsque $\fa\subseteq\fb$ et $\ff\subseteq\ffg$.
\end{definition}\index{raffine}
%--- end-definition------------------------------------

Le lemme suivant décrit le \pa \gui{engendré} (au sens de la relation de
raffinement) par un couple de parties de $\gA$.
En fait il suffit de traiter le cas d'un couple formé par un \id et un \moz.

%--- Lemma{lempaireengendree}------
\begin{lemma}
\label{lempaireengendree}
Soit un \id $\fa$ et un \mo $\ff$ de $\gA$.
%-----------------begin enum------------------
\begin{enumerate}

\item Le \pa $(\fb ,\ffg )$ engendré par $(\fa ,\ff )$ est obtenu comme suit:
\[\fb =\sotq{x\in\gA}{\exists s\in \ff,\, xs\in \fa },
\hbox{ et }\,
\ffg =\sotq{y\in\gA}{\exists u\in\gA,\, uy\in \fa +\ff}.
\]

\item Si $\ff \subseteq\Ati$, alors
$\fb=\fa$ et $\ffg $ est le filtre obtenu en saturant le
\moz~\hbox{$1+\fa $}. Dans ce cas, $\gA_{\ffg }\sur{\fa }=\gA\sur{\fa }$.

\item Si $\fa =0$, alors
$\fb =\sotq{x\in\gA}{\exists s\in \ff,\,xs=0}=\sum_{s\in\ff}(0:s)$, et  $\ffg$ est le saturé de $\ff$.
Dans ce cas, $\gA_{\ffg }\sur{\fb }=\gA_{\ff }$.
Si en outre $\ff=s^{\NN} $,
$\fb =(0:s^\infty)$.
\end{enumerate}
%-----------------end enum------------------
\end{lemma}
%--- end-lemma-----------------------------------------

\rdb \label{NOTASatu}

Un cas important est celui du filtre obtenu par saturation d'un \mo $S$.
Nous introduisons la notation $\sat{S}$, ou, si \ncrz, $\satu S \gA$
pour ce filtre.

% :   subsec{Idéal et filtre incompatibles}
\subsubsection*{Idéal et filtre incompatibles}

Pour un \pa $(\fa,\ff)$ on a les équivalences suivantes.
%------begin equation--eqIncompatibles-----------
\begin{equation}\label{eqIncompatibles}
\fa=\gA\;\Longleftrightarrow\; 1\in\fa\;\Longleftrightarrow\;
0\in\ff\;\Longleftrightarrow\; \ff=\gA\;\Longleftrightarrow\;
\gA_\ff\sur{\fa}=\so{0}.
\end{equation}
%---------------------end equation--------------

Un idéal $\fa$ et un filtre $\ff$ sont dit  \emph{incompatibles} lorsqu'ils
engendrent la paire~\hbox{$(\gA,\gA)$}, \cad lorsque $0\in\fa+\ff$.
%Un \pa est incompatible \ssi il est égal à $(\gA,\gA)$.%

Un idéal $\fa$ et un filtre $\ff$ sont dit  \emph{compatibles} si
l'on \hbox{a $(0\in\fa+\ff\Rightarrow1=0)$}.
Si l'anneau est non trivial cela signifie aussi
$\fa\cap\ff= \emptyset$.
Dans ce cas on peut à la fois annuler les \elts de $\fa$ et rendre \ivs
les \elts de $\ff$ sans que l'anneau ne soit réduit à $0$.%
\index{compatibles!idéal et filtre ---}%
\index{compatible!\pa ---}%
\index{incompatible!\pa ---}%
\index{incompatibles!idéal et filtre ---}

%:     Fact{factPremComp}
\begin{fact}\label{factPremComp}
Soit $\fa$ un \id et $\ff$ un filtre compatibles.
\\
Si $\fa$ est premier, il est $\ff$-saturé, si $\ff$ est premier,
il est $\fa$-saturé.
\end{fact}

% :   subsec{Le treillis des \pas}
\subsubsection*{Le treillis des \pas}

%--- Fact{corlempaireengendree}-
\begin{fact}
\label{corlempaireengendree}
Les \pas de $\gA$  ont une structure de treillis
pour la relation de raffinement, avec:
%-----------------begin item------------------
\begin{enumerate}
\item [--] L'\elt minimum est $(\so{0},\Ati)$ et l'\elt maximum  $(\gA,\gA)$.
\item [--] $(\fa,\ff)\vu(\fb,\ffg)$ est le \pa engendré par $(\fa+\fb,\ff\,\ffg)$.
\item [--] $(\fa,\ff)\vi(\fb,\ffg)=(\fa\cap\fb,\ff\cap\ffg)$.
\end{enumerate}
%-----------------end item------------------
\end{fact}
%--- end-fact------------------------------------

%:   subsec{Idéaux et filtres dans un quotient localisé}
%%%%%%%%%%%%
\subsubsec{Idéaux et filtres dans un quotient localisé}

%--- Fact{propPaires}------
\begin{fact}
\label{propPaires} 
 Soit  $(\fa,\ff)$ un \pa de $\gA$ et $\pi:\gA\to\gB=\gA_\ff\sur{\fa}$
l'application canonique.
Alors:
%-----------------begin enum------------------
\begin{enumerate}
\item
L'application $(\fb,\ffg)\mapsto \big(\pi^{-1}(\fb),\pi^{-1}(\ffg)\big)$
est une bijection croissante (pour les relations de raffinement) entre d'une part,
les \pas de $\gB$, et d'autre part, les
\pas de $\gA$ qui raffinent $(\fa,\ff)$.
\item Si $(\fb,\ffg)$ est un \pa de $\gB$ l'application canonique
$$\gA_{\pi^{-1}(\ffg)}\sur{\pi^{-1}(\fb)} \;\longrightarrow\; \gB_{\ffg}\sur{\fb}$$
est un \isoz.
\item Dans cette bijection
%-----------------begin item------------------
\begin{enumerate}
\item [--] l'\id $\fb$ est premier \ssi $\pi^{-1}(\fb)$
est premier,
\item [--] tout \id premier de $\gA$ compatible avec $\ff$ et
contenant $\fa$ est obtenu,
\item [--] le filtre $\ffg$ est premier \ssi $\pi^{-1}(\ffg)$
est premier,
\item [--] tout filtre premier de $\gA$ compatible avec $\fa$
et contenant $\ff$ est obtenu.
\end{enumerate}
%-----------------end item------------------
%
\end{enumerate}
%-----------------end enum------------------
\end{fact}
%--- end-fact----------------------------------------

On en déduit la comparaison suivante qui est instructive sur la dualité entre
\ids et filtres.

%\newpage

\DeuxCol%debut DeuxCol
% debut premiere colonne
{
%--- Fact{factQuoIDFI}----------
\begin{fact}
\label{factQuoIDFI}
Soit $\fa$ un \id strict de $\gA$ et $\pi:\gA\to%\gB=
\gA/\fa$ l'\homo correspondant.
%-----------------begin enum------------------
\begin{enumerate}
\item L'application $\fb\mapsto\pi^{-1}(\fb)$ est une bijection
croissante  entre \ids de $\gA/\fa$ et \ids de $\gA$ contenant $\fa$.
Dans cette bijection les \ideps correspondent aux \idepsz.
\item L'application $\ffg\mapsto\pi^{-1}(\ffg)$ est une bijection
croissante
entre filtres de $\gA/\fa$ et filtres $\fa$-saturés  de $\gA$.
\item  Dans cette bijection les filtres premiers stricts de  $\gA/\fa$
correspondent exactement aux filtres premiers de $\gA$
com\-pa\-tibles avec~$\fa$.
\end{enumerate}
%-----------------end enum------------------
\end{fact}
%--- end-fact-----------------------------------------

}%fin premiere colonne debut deuxieme colonne
{
%--- Fact{factQuoFIID}-----
\begin{fact}
\label{factQuoFIID}
Soit $\ff$ un filtre strict de $\gA$ et $\pi:\gA\to%\gB=
\gA_\ff$ l'\homo correspondant.
%-----------------begin enum------------------
\begin{enumerate}
\item L'application $\ffg\mapsto\pi^{-1}(\ffg)$ est une bijection
croissante entre filtres de $\gA_\ff$ et filtres de $\gA$ contenant $\ff$.
Dans cette bijection les filtres premiers correspondent aux filtres premiers.
\item L'application $\fb\mapsto\pi^{-1}(\fb)$ est une bijection
croissante
entre \ids de $\gA_\ff$ et \ids $\ff$-saturés  de $\gA$.
\item  Dans cette bijection les \ideps stricts de  $\gA_\ff$
correspondent exactement aux \ideps de $\gA$
com\-pa\-tibles avec~$\ff$.
\end{enumerate}
%-----------------end enum------------------
\end{fact}
%--- end-fact-----------------------------------------

}%fin DeuxCol

%:   subsec{Principes de recouvrement fermé} %%%%%%%%%%%%
\subsec{Principes de recouvrement fermé}

La dualité entre \ids et filtres suggère qu'un principe dual
du \plg doit pouvoir fonctionner en \alg commutative.

Tout d'abord notons que les \ids de $\ZarA$ correspondent
bijectivement aux \ids
radicaux (i.e, égaux à leur nilradical) de
$\gA$ via: 

\snic{\fa$ (\id de $\ZarA$) $\mapsto$ $\sotq{x\in\gA}{\DA(x)\in\fa }.}

%\sni
En outre, les \ideps correspondent aux \idepsz.

Pour les filtres, cela ne se passe pas aussi parfaitement, mais pour un filtre
$\ff$ de $\gA$, l'ensemble $\sotq{\DA(x)}{x\in\ff}$ engendre un filtre de
$\ZarA$, et ceci donne une application injective qui est  bijective pour les filtres premiers.

Revenons au \plg et regardons ce que cela signifie  
dans le treillis
$\ZarA$. Lorsque l'on a des \moco $S_1$, \ldots, $S_n$ de $\gA$, cela correspond
à des filtres $\ff_i $ de $\ZarA$
(chacun engendré par les $\DA(s)$ pour les $s\in S_i$)
qui sont \gui{\comz}  en ce sens que $\bigcap_i\ff_i=\so{1_\ZarA}$.
Dans ce cas les \homos naturels

\snic{\gA\to\prod_i\gA_{S_i}\quad \hbox{ et }
\quad\ZarA\to\prod_i\ZarA\sur{(\ff_i=1)}}

%\sni
sont injectifs.

Par dualité, on dira qu'un \sys d'\ids $(\fa_1,\ldots,\fa_n)$ constitue
un \ix{recouvrement fermé} de $\gA$ lorsque $\bigcap_i\DA(\fa_i)=\so{0_\ZarA}$,
\cade  lorsque $\prod_i\fa_i\subseteq\DA(0)$.
Dans ce cas les \homos naturels

\snic{\gA/\DA(0)\to\prod_i\gA\sur{\DA(\fa_i)}\quad\hbox{ et }
\quad \ZarA\to\prod_i\ZarA\sur{(\DA(\fa_i)=0)}}

%\sni
 sont injectifs.

Nous dirons qu'une \prt $\sfP$ (concernant des objets reliés à
un anneau~$\gA$) vérifie le \gui{principe de recouvrement
fermé} lorsque:\\
\emph{chaque fois que des \ids $\fa_i$ forment un recouvrement fermé de $\gA$,
la \prt
$\sfP$ est vraie pour $\gA$ \ssi elle est vraie
après passage au quotient par chacun des $\fa_i$.
}

Par exemple on obtient facilement (voir aussi le lemme~\ref{lemNilpotProd}).

%%%%%%%%%%%%%%%%%%%%%%%%%%%%%%%%%%%%%%%%%
%:     prcf1
\begin{prcf}\label{prcf1} \emph{(\'Eléments nilpotents, \comz)}
On considère un recouvrement fermé $(\fa_1,\ldots,\fa_r)$ de l'anneau
$\gA$. Soient $x_1$, \dots, $x_n\in\gA$, $\fb$, $\fc$
deux \ids et $S$ un \moz.
\begin{enumerate}
\item Le \mo $S$ contient $0$  \ssi il contient $0$  modulo chaque~$\fa_i$.
\item On a $\fb\subseteq\sqrt{\fc}$   \ssi  $\fb\subseteq\sqrt{\fc}$
modulo chaque~$\fa_i$.
\item Les \elts $x_1$, \dots, $x_n$ sont \com \ssi ils sont \com modulo chaque~$\fa_i$.
\end{enumerate}
\end{prcf}
\begin{proof}
Il suffit de montrer le point \emph{2.} On suppose que $\DA(\fb)\leq \DA(\fc)\vu \DA(\fa_i)$,
donc $\DA(\fb)\leq \Vi_i \left(\DA(\fc)\vu \DA(\fa_i)\right) = \DA(\fc)\vu \left(\Vi_i\DA(\fa_i)\right)=\DA(\fc)$.
\end{proof}
\rem Par contre, il n'y a pas de \prf pour les solutions de \slisz.
Considérons en effet $u,\,v\in\gA$ tels \hbox{que $uv=0$}.
Le \sli (avec $x$ pour inconnue)  

\snic{ux=u,$  $vx=-v,}

%\sni
admet une solution modulo $u$ (à savoir $x=-1$) et une solution
modulo $v$ (à savoir $x=1$). Mais dans le cas de l'anneau
$\gA=\ZZ[u,v]=\aqo{\ZZ[U,V]}{UV}$ le \sli n'a pas de solution
dans $\gA$.
\eoe

%%%%%%%%%%%%%%%%%%%%%%%%%%%%%%%%%%%%%%%%%
%:     prcf2
\begin{prcf}\label{prcf2} \emph{(Modules \tfz)}\\ 
On considère un recouvrement fermé $(\fa_1,\ldots,\fa_r)$ de l'anneau
$\gA$.
On suppose \hbox{que $\prod_i\fa_i=0$} (c'est le cas si $\gA$ est
réduit).
Un \Amo $M$ est \tf  \ssi il est \tf
modulo chaque~$\fa_i$.
\end{prcf}
\begin{proof}
On suppose \spdg $r=2$. Soient $g_1$, \ldots, $g_k$ des \gtrs modulo $\fa_1$,
et  $g_{k+1}$, \ldots, $g_\ell$ des  \gtrs modulo $\fa_2$. Soit $x\in M$.
\linebreak 
On écrit
$x=\sum_{i=1}^k\alpha_ig_i+\sum_{j=1}^p\beta_jx_j$ avec $\alpha_i\in\gA$,
$\beta_j\in\fa_1$, $x_j\in M$.
\\ 
Chaque $x_j$ s'écrit comme une \coli de  $g_{k+1}$, \ldots, $g_\ell$
modulo $\fa_2$.  \hbox{Puisque $\fa_1\fa_2=0$}, on obtient $x$ comme \coli \hbox{de
$g_1$, \ldots, $g_\ell$}.
\end{proof}
%

%%%%%%%%%%%%%%%%%%%%%%%%%%%%%%%%%%%%%%%%%
%:     prcf3
\begin{prcf}\label{prcf3} \emph{(Modules \ptfsz)}
On considère un recouvrement fermé $(\fa_1,\ldots,\fa_r)$ de l'anneau
$\gA$, une matrice  $F\in\Ae{m\times n}$, $\fa$
un \itf
et $M$ un \mpfz.
\begin{enumerate}
\item La matrice $F$ est de rang $\geq k$  \ssi elle  est de rang $\geq k$
modulo chaque~$\fa_i$.
\end{enumerate}
Supposons $\bigcap_i\fa_i=0$ (c'est le cas si $\gA$ est
réduit). Alors:
\begin{enumerate} \setcounter{enumi}{1}
\item La matrice $F$ est de rang $\leq k$  \ssi elle  est de rang $\leq k$
modulo chaque~$\fa_i$.
\item L'\itf $\fa$ est engendré par un \idm  \ssi il  est engendré par un \idm
modulo chaque~$\fa_i$.
\item La matrice $F$ est \lnl  \ssi elle  est \lnl
modulo chaque~$\fa_i$.
\item Le module $M$ est \ptf  \ssi il  est \ptf
 modulo chaque~$\fa_i$.
\end{enumerate}
\end{prcf}
\begin{proof}
Le point \emph{1} résulte du \prf \ref{prcf1} en considérant
l'\idd d'ordre $k$. Le point \emph{2} vient de ce que si un \idd
est nul modulo chaque $\fa_i$, il est nul modulo leur intersection.
Le point~\emph{5} est une reformulation du point~\emph{4} 
qui est une conséquence du point~\emph{3.}\\
Montrons le point \emph{3.} On suppose \spdg $r=2$.
On utilise le lemme de l'\id engendré par un \idm (lemme~\ref{lem2ide.idem}).
On a

\snic{\fa+(0:\fa)_{\gA\sur{\fa_i}}=\gA\sur{\fa_i}\;\; (i=1,2). }

%\sni
Cela signifie $\fa+\fa_i+(\fa_i:\fa)=\gA$, et puisque $\fa_i\subseteq(\fa_i:\fa)$,
on a $1\in\fa+(\fa_i:\fa)$. En faisant le produit cela donne $1\in\fa+(\fa_1:\fa)(\fa_2:\fa)$ et puisque

\snic{(\fa_1:\fa)(\fa_2:\fa)\subseteq  (\fa_1:\fa)\cap(\fa_2:\fa) =((\fa_1\cap\fa_2):\fa)=(0:\fa), }

%\sni
 on obtient $1\in\fa+(0:\fa)$.
\end{proof}
%

%%%%%%%%%%%%%%%%%%%%%%%%%%%%%%%%%%%%%%%%%%%%%%%%%%%%%%%%%%%%%%%%%
%: \subsec{Anneau \zedr engendré par un anneau commutatif} %%%%%%%%%%%%
\subsec{Clôture \zede réduite d'un anneau commutatif}
\label{secClotureZEDR}

Commençons par des résultats concernant un sous-anneau $\gA$ d'un
anneau \zedrz. \Llec peut se reporter 
à l'étude des anneaux \zedrs   \paref{subsecAzedred}
et revisiter les \egts (\iref{eqQuasiInv}) pour la \carn d'un quasi inverse.

Si dans un anneau un \elt $c$
admet un quasi inverse, nous  notons \hbox{celui-ci $c\bul$},
et nous notons $e_c=c\bul c$ l'\idm associé à~$c$
qui satisfait les \egts $\Ann(c)=\Ann(e_c)= \gen {1-e_c}$.

%:     Lemma{lem1SousZedRed}
\begin{lemma}\label{lem1SousZedRed}
\emph{(Anneau engendré par un quasi inverse)}
\begin{enumerate}
\item Soit $a\in\gA\subseteq\gB$. On suppose que $\gA$ et $\gB$
sont réduits et que $a$ admet un quasi inverse dans $\gB$.
Alors

\snic{\gB\supseteq\gA[a\bul]\simeq\aqo{\gA[a\bul]}{1-e_a}\times \aqo{\gA[a\bul]}{e_a}=
\gA_1\times \gA_2.}

%\sni
En outre:
\begin{enumerate}
\item On a un \homo naturel bien défini $\gA[1/a]\to\gA_1$,
et c'est un \isoz.  En particulier, l'\homo naturel~\hbox{$\gA\to\gA_1$}
a pour noyau $\Ann_\gA(a)$.
\item
L'\homo naturel $\gA\to\gA_2$ est surjectif,
son  noyau est l'intersection \hbox{$\fa=\gA\cap e_a\gA[a\bul]$}
et  vérifie la double inclusion

\snic{\Ann_\gA\big(\Ann_\gA(a)\big)\supseteq \fa \supseteq \DA(a)\qquad(*).}
\end{enumerate}
En bref $\gA[a\bul]\simeq\gA[1/a]\times \gA\sur{\fa}$.
\item Inversement pour tout \id $\fa$ de $\gA$ vérifiant $(*)$, l'\elt $(1/a,0)$ est un quasi inverse
de (l'image de) $a$ dans l'anneau $\gC=\gA[1/a]\times \gA\sur{\fa}$ et l'\homo canonique
de $\gA$ dans $\gC$ est injectif.
\end{enumerate}
\end{lemma}
\begin{proof}
L'\iso $\gA[a\bul]\simeq \gA_1\times \gA_2$ signifie seulement que $e_a$ est un \idm dans $\gA[a\bul]$. Nous notons
%$j:\gA\to\gA[a\bul]$ et
$\pi_i:\gA[a\bul]\to\gA_i$ les \homos canoniques.

\emph{1b.} Soit $\mu$ l'\homo composé $\gA\lora \gA[a\bul] \lora \gA_2$.
Dans $\gA_2$, \linebreak 
on a \hbox{$a\bul=e_aa\bul=0$}, donc
$\gA_2=\gA\sur{\left(\gA\cap {e_a \gA[a\bul]}\right)}$. Ainsi~\hbox{$\fa=\gA\cap e_a\gA[a\bul]$}.
Dans~$\gA[a\bul]$, on a $a=e_aa$, donc $\mu(a)=\pi_2(a)=\pi_2(e_aa)=0$,  
et $a\in\fa$.
\\
Comme $\gB$ est réduit, les trois anneaux $\gA[a\bul]$, $\gA_1$ et $\gA_2$ le sont aussi.\linebreak 
Donc $\gen{a}\subseteq \fa$ implique $\DA(a)\subseteq \fa$.\\
Enfin,
$\fa\,\Ann_\gA(a) \subseteq \gen{e_a}\Ann_\gA(a) = 0$, donc
     $\fa\subseteq\Ann_\gA\big(\Ann_\gA(a)\big)$.

%Soit
%maintenant $x\in \fa$ et $y\in\Ann_\gA(a)$, alors $x=ze_a$ dans $\gA[a\bul]$ et $ay=0$ dans $\gA$. Donc $xy=za\bul ay=0$ dans $\gA[a\bul]\supseteq\gA$.
%Ainsi $x\in\Ann_\gA(y)$. En conclusion $\fa\subseteq\Ann_\gA(\Ann_\gA(a))$.

\emph{1a. }
Puisque $aa\bul=_{\gA_1}1$, on a un unique \homo $\lambda:\gA[1/a]\to\gA_1$
obtenu à partir de l'\homo composé $\gA\to \gA[a\bul] \to \gA_1$, et
$\lambda$ est clairement surjectif. Considérons un \elt $x/a^n$ de $\Ker\lambda$. Alors $\lambda(ax)=0$, donc $\pi_1(ax)=0$. Comme on a aussi $\pi_2(ax)=0$, on en déduit $ax=0$, \linebreak 
donc $x=_{\gA[1/a]}0$. Ainsi $\lambda$ est injectif.

\emph{2.} L'image de $a$ dans $\gC$ est $(a/1,0)$, donc $(1/a,0)$
est bien son quasi inverse. Soit maintenant $x\in\gA$ dont l'image dans $\gC$
est $0$.
D'une part $x=_{\gA[1/a]}0$, donc $ax=_\gA0$. D'autre part  $x\,\Ann_\gA(a)=0$ donc $x^2=_\gA0$, et $x=_\gA0$.
\end{proof}

\comm On voit
que la notation $\gA[a\bul]$ présente a priori une possible
ambiguïté, au moins lorsque $\DA(a)\neq\Ann_\gA\big(\Ann_\gA(a)\big)$. \eoe

\vspace{-.5mm}

%:     Lemma{lem2SousZedRed}
\begin{lemma}\label{lem2SousZedRed}
Si $\gA\subseteq\gC$ avec $\gC$ \zedrz, le plus petit sous-anneau \zed de $\gC$ contenant
$\gA$ est égal à $\gA[(a\bul)_{a\in\gA}]$. Plus \gnlt si $\gA\subseteq\gB$
avec $\gB$ réduit, et si tout \elt de $\gA$ admet un quasi inverse dans $\gB$,
alors le sous-anneau $\gA[(a\bul)_{a\in\gA}]$ de $\gB$ est \zedz.
En outre, tout \elt de $\gA[(a\bul)_{a\in\gA}]$ s'écrit sous forme 
$$\preskip.4em \postskip-.8em \ndsp 
\som_j{  a_jb_j\bul e_j}, \; \hbox{ avec}
$$
\begin{itemize}
\item les $e_j$ sont des \idms de $\gA[(a\bul)_{a\in\gA}]$
deux à deux \ortsz,
\item  $a_j,b_j\in\gA$ et $b_j{b_j\bul} e_j=e_j$ pour tout $j$, 
\end{itemize}
%\vspace{-1mm}
de sorte que $\big(\som_j{  a_jb_j\bul e_j}\big)\bul=\som_j{  a_j\bul b_j e_j}$.
\end{lemma}
NB. On prendra garde cependant que l'on n'a pas toujours $a_j{a_j\bul} e_j=e_j$. Il faut donc a priori remplacer $e_j$ par $e'_j=a_j{a_j\bul} e_j$ pour obtenir une écriture du même type que la précédente. On pourra aussi noter que tout \idm de $\gA[(a\bul)_{a\in\gA}]$ s'écrit sous forme $e_c\prod_i(1-e_{d_i})$ pour un $c$ et des $d_i\in\gA$.
\begin{proof}
Parmi les \elts de $\gB$, ceux qui s'écrivent comme somme de produits~$ab\bul$ avec $a$, $b\in\gA$ forment clairement un sous-anneau de $\gB$,
qui est donc égal à $\gA[(a\bul)_{a\in\gA}]$.
Par ailleurs, $ab\bul=ab\bul e_b$.
En considérant l'\agB engendrée par les $e_b$ qui interviennent dans une somme finie du type précédent, on en déduit que tout \elt de $\gA[(a\bul)_{a\in\gA}]$
peut s'écrire sous  la forme

\snic{\som_j{ \, \left(\som_i a_{i,j}b_{i,j}\bul\right)\,e_j},\;\;\hbox{avec}}

%\sni
\begin{itemize}
\item les $e_j$ sont des \idms dans $\gA[(a\bul)_{a\in\gA}]$
deux à deux orthogonaux,
\item  $a_{i,j},b_{i,j}\in\gA$, et $b_{i,j}b_{i,j}\bul e_j=e_j$, pour tous $i,j$.
\end{itemize}
On note que  $b_{i,j}\bul$ est l'inverse de  $b_{i,j}$
dans $\gA[(a\bul)_{a\in\gA}][1/e_j]$, et
on peut faire le calcul comme pour une somme de fractions
ordinaires $\sum_i a_{i,j}/b_{i,j}$.
Par exemple prenons pour simplifier un
terme avec une somme de 3 \elts 

\snic{  (a_1b_1\bul+a_2b_2\bul+a_3b_3\bul) e.}

%\sni
Puisque $b_2b_2\bul e=b_3b_3\bul e=e$, on a $a_1b_1\bul e=a_1b_2b_3(b_1b_2b_3)\bul e$,
et

\snic{(a_1b_1\bul+a_2b_2\bul+a_3b_3\bul) e = (a_1b_2b_3+a_2b_1b_3+a_3b_1b_2)(b_1b_2b_3)\bul e = d c\bul e,}

%\sni
qui admet pour quasi inverse $c d\bul e$.
\end{proof}

Rappelons que  $\gB\red$ désigne le quotient d'un anneau $\gB$ par son nilradical.

Dans le lemme qui suit on regarde ce qui se passe lorsque l'on rajoute de force
un quasi inverse à un \elt d'un anneau. C'est une opération voisine de la localisation, lorsque l'on rajoute de force un inverse d'un \eltz,
mais un peu plus délicate.

%: --- Lemma{lemZedClotBasicStep}-------
\begin{lemma}
\label{lemZedClotBasicStep} Soit $\gA$ un anneau et $a\in\gA$.
%-----------------begin enum------------------
\begin{enumerate}
\item Considérons l'anneau $\aqo{\gA[T]}{aT^2-T,a^2T-a}=\gA[a\efl  ]$
et notons $\lambda_a:\gA\to\gA[a\efl  ]$ l'\homo canonique ($a\efl$
désigne l'image de $T$). Alors
pour  tout \homo $\psi
:\gA\to\gB$ tel que $\psi(a)$ admette un quasi inverse
il existe un unique \homo $\varphi:\gA[a\efl  ]\to \gB$
tel \hbox{que $\varphi\circ \lambda_a=\psi$}.

\vspace{-1em}
{\Pnv{\gA}{\lambda_a}{\psi}{\gA[a\efl  ]}{\varphi}{\gB}{ }{}{$\psi(a)$ admet un quasi inverse}} 

\vspace{-1em}
\item  En outre, $aa\efl  $
est un \idm et $\gA[a\efl  ]\simeq \gA[1/a] \times \aqo{\gA}{a}$.
\item  Si $\gB$ est réduit on a une factorisation
unique via $(\gA[a\efl  ])\red$.
\end{enumerate}
Dans la suite on note $\gA[a\bul]$ pour l'anneau $(\gA[a\efl  ])\red$.
\label{NOTAAbul}
\begin{enumerate}\setcounter{enumi}{3}
\item
On a $\gA[a\bul] \simeq \Ared[1/a]\times \gA\sur{\DA(a)}$.
Si $\gA$ est réduit l'\homo canonique $\gA\to\gA[a\bul]$
est injectif.
\item \label{ite5lemZedClotBasicStep}
$\Zar(\gA[a\bul]) =\Zar(\gA[a\efl])$ s'identifie à
$(\ZarA)[\DA(a)\bul].$
\end{enumerate}
%-----------------end enum------------------
\end{lemma}
%--- end-lemma-----------------------------------------

\begin{proof}
Laissée \alecz. Le dernier point résulte du lemme \ref{lemRajouCompl}
et du fait~\ref{fact2Zar}.
\end{proof}

%:     Corollary{corlemZedClotBasicStep}
\begin{corollary}\label{corlemZedClotBasicStep}
Soient $a_1$, \ldots, $a_n\in\gA$.
\begin{enumerate}
\item L'anneau $\gA[a_1\bul][a_2\bul]\cdots[a_n\bul]$
est indépendant, à \iso unique près, de l'ordre des $a_i$.
Il sera noté $\gA[a_1\bul,a_2\bul,\ldots,a_n\bul]$.
\item Une description possible est la suivante:
$$\preskip.3em \postskip.3em
\gA[a_1\bul,a_2\bul,\ldots,a_n\bul]
\simeq
\big(\gA[T_1,T_2,\ldots,T_n]\sur{\fa}\big)\red
$$
avec $\fa=\gen{(a_iT_i^2-T_i)_{i=1}^n,(T_ia_i^2-a_i)_{i=1}^n}$.
\item Une autre description possible:
$$\preskip.3em \postskip.1em
\gA[a_1\bul,a_2\bul,\ldots,a_n\bul]
\simeq
\prod\nolimits_{I\in\cP_n}{\left(\aqo{\gA}{(a_i)_{i\in I}}\right)\red[1/\alpha_I]}
$$
avec $\alpha_I=\prod_{j\in \lrbn\setminus I}a_j$.
\end{enumerate}
\end{corollary}

%:     Theorem{thZedGen}
\begin{theorem}\label{thZedGen}
\emph{(Clôture \zede réduite d'un anneau
commu\-tatif)}
Pour tout anneau $\gA$ il existe un  anneau \zedr $\Abul$ avec un
\homo $\lambda:\gA\to\Abul$, qui
factorise de manière unique tout \homo $\psi :\gA\to\gB$ vers un anneau \zedrz.

\vspace{-1.2em}
\PNV{\gA}{\lambda}{\psi}{\Abul}{\varphi}{\gB}{anneaux commutatifs}{}{anneaux \zeds réduits}

\vspace{-1.5em}
Ce couple $(\Abul,\lambda)$ est unique à \iso unique près.
En outre:\index{cloture@clôture!zero@\zede réduite}
\begin{enumerate}
\item [--] L'\homo naturel $\Ared\to \Abul$ est injectif.
\item [--] On a $\Abul=\gA\red[(a\bul)_{a\in\Ared}]$
\end{enumerate}

\end{theorem}
\begin{proof}
Ceci est un corolaire des lemmes précédents.
On peut supposer $\gA$ réduit.
Le résultat d'unicité (corolaire \ref{corlemZedClotBasicStep}) permet de
 construire une limite inductive (qui mime une réunion filtrante)
basée sur les extensions du type~\hbox{$\gA[a_1\bul,a_2\bul,\ldots,a_n\bul]$},
et l'on conclut avec le lemme~\ref{lem2SousZedRed}.
\end{proof}

\comms 1) A priori,
puisque l'on  a affaire à des structures purement équa\-tionnelles,
 la clôture \zede réduite \uvle d'un anneau existe
et l'on pourrait la construire comme suit:
on rajoute \fmt l'opération unaire $a\mapsto a\bul$ et l'on force
$a\bul$ à être un quasi inverse de $a$. Notre preuve a permis de donner
en plus une description précise simplifiée de l'objet construit et de montrer l'injectivité dans le cas réduit.

2) En \clamaz, la clôture \zede réduite $\Abul$ d'un anneau $\gA$ peut être obtenue comme suit.
Tout d'abord on considère le produit $\gB=\prod_\fp \Frac(\gA\sur{\fp})$, où $\fp$
parcourt tous les \ideps de $\gA$. Comme $\gB$ est un produit de corps, il est \zed réduit.
Ensuite on considère le plus petit sous-anneau
\zed de $\gB$ contenant l'image de $\gA$ dans $\gB$ par l'\homo diagonal naturel.
\\
On comprend alors l'importance de la construction à la main que nous avons faite de $\Abul$.
Elle nous permet d'avoir accès de manière explicite
à quelque chose qui ressemble à \gui{l'ensemble de tous les} \ideps
de~$\gA$ (ceux des \clamaz) sans avoir besoin d'en construire un seul individuellement.
Le pari est que les raisonnements des \clama qui manipulent des \ideps
arbitraires non précisés de l'anneau $\gA$ (des objets en \gnl
inaccessibles) peuvent être relus comme
des arguments au sujet de l'anneau $\Abul$: un objet sans aucun mystère!
\eoe

\medskip\rdb
\exls \label{exlCLZED}~
\\
1) Voici une description de la clôture \zede réduite de $\ZZ$.\\
Tout d'abord, pour $n\in\NN\etl$ l'anneau $\ZZ[n\bul]$ est isomorphe à $ \ZZ[1/n]\times \prod_{p|n}\FF_p$,
où $p$ est mis pour \gui{$p$ nombre premier},
et $\FF_p=\ZZ\sur{p\ZZ}$. Ensuite, $\ZZ\bul$ est la
limite inductive (que l'on peut voir comme une réunion
croissante) des~$\ZZ[(n!)\bul]$.

2) Voici une description de  la clôture \zede réduite de
$\ZZ[X]$.\\
 Tout d'abord,  si $Q$ un \polu sans facteur carré,
 et si $n\in\NN\etl$ est multiple de $\disc(Q)$,
l'anneau $\ZZ[X][n\bul,Q\bul]$ est isomorphe à
$$
\ZZ[X,1/n,1/Q]\times \prod_{p\divi n}\FF_p[X,1/Q]
\times \prod_{P\divi Q}\aqo{\ZZ[X,1/n]}{P} \times
\prod_{p\divi n,R\divi Q}\!\!\!\aqo{\FF_p[X]}{R}~~
$$
avec $p$ mis pour \gui{$p$  premier}, $P$ mis pour \gui{$P$ irréductible
dans $\ZZ[X]$}, et~$R\divi Q$ mis pour \gui{$R$ \ird
dans $\FF_p[X]$ divise $Q$ dans $\FF_p[X]$}. 
\\
Ensuite, on passe à la limite inductive  des  %anneaux 
$\ZZ[X][u_n\bul,Q_n\bul]$ (ici, c'est une réunion croissante), où $Q_n$ est la partie sans carré du produit des
$n$ premiers \elts dans une énumération des \polus  sans facteur carré
de $\ZZ[X]$, et où~\hbox{$u_n=n!\disc(Q_n)$}. 
\\
Notez  que l'on obtient ainsi un anneau
par lequel se factorisent tous les \homos naturels $\ZZ[X]\to\Frac(\ZZ[X]\sur\fp)$
pour tous les \ideps $\fp$ de $\ZZ[X]$: un tel  $\Frac(\ZZ[X]\sur\fp)$
est en effet ou bien $\QQ(X)$, ou bien un $\aqo{\QQ[X]}{P}$, ou bien
un $\FF_p(X)$, ou bien un $\aqo{\FF_p[X]}{R}$.

3) L'anneau (\cot bien défini) $\RR\bul$ est certainement un des
objets les plus intrigants qui soient pour l'investigation du
monde \gui{sans tiers exclu} que constituent les \comaz.
Naturellement, en \clamaz, $\RR$ est \zed et $\RR\bul=\RR$. 
\eoe

%:     Theorem{thZedGenEtBoolGen}
\begin{theorem}\label{thZedGenEtBoolGen}
Pour tout anneau $\gA$ on a des \isos naturels 

\snic{\Bo(\ZarA)\simeq \BB(\Abul)\simeq\Zar(\Abul).}

\sni
\end{theorem}
\begin{proof}
Cela résulte du dernier point du lemme \ref{lemZedClotBasicStep}, et du fait que les deux constructions peuvent être vues
comme des limites inductives de \gui{constructions à un étage}
$\gE\leadsto\gE[a\bul]$ ($\gE$ un anneau ou un \trdiz).
\end{proof}

Notez que si l'on  adoptait la notation $\gT\bul$ pour $\Bo(\gT)$
on aurait la jolie formule $(\ZarA)\bul\simeq\Zar(\Abul)$.

%:     Proposition{propClZdrLoc}
\begin{proposition}\label{propClZdrLoc}
Soient $\gA$ un anneau, $\fa$ un \id et $S$ un \moz.  
\begin{itemize}
\item Les deux anneaux $(\gA\sur\fa)\bul$
et $\Abul\sur{\rD(\fa\Abul)}$ sont canoniquement isomorphes.
\item  Les deux anneaux $(\gA_S)\bul$
et $(\Abul)_S$ sont canoniquement isomorphes.
\end{itemize}
 
\end{proposition}
%--------- fin proposition ---------------------------------------------- 
%
\begin{proof} 
Notons que $(\Abul)_S$ est \zedr comme \lon d'un anneau \zedrz.  Et de même,
$\Abul\sur{\rD(\fa\Abul)}$ est \zedrz. \'Ecrivons la \dem pour les \lonsz.
Considérons les \homos naturels 

%%\snic{\gA\vers{\jmath}\gA_S\vers{\lambda}(\gA_S)\bul\eqdefi\gB\quad$ et 
%%$\quad\gA\vers{\theta}\Abul \vers{\imath}(\Abul )_S\eqdefi\gC.}

\snic{
\gA \to \gA_S \to (\gA_S)\bul  \quad\hbox {et}\quad
\gA  \to \Abul \to (\Abul)_S
.}

%\sni

L'\homo $\gA \to \Abul$ \gui {se prolonge} de manière unique en un
\homo $\gA_S \to (\Abul)_S$, et d'après la \prt \uvle de la clôture
\zede réduite, fournit un morphisme unique $(\gA_S)\bul \to (\Abul)_S$ qui
rend le diagramme ad-hoc commutatif. De même, le morphisme $\gA \to \gA_S$
donne naissance à un unique morphisme $\Abul \to (\gA_S)\bul$ qui se
prolonge en un morphisme $(\Abul)_S\to (\gA_S)\bul$.  En composant ces deux
morphismes, par unicité, on obtient deux fois l'identité.
\end{proof}
%

%%%%%%%%%%%%%%%%%%%%%%%%%%%%%%%%%%%%%%%%%%%%%%%%%%%%%%%%%%%%%%%%%
%  \subsec{Treillis distributifs et \entrelsz}%%%%%%%%%%%%
\section[Relations implicatives]{Treillis distributifs, \entrels et \agHsz}
\label{secEntRelAgH}
\vspace{4pt}
\subsec{Un nouveau regard sur les \trdisz}

Une règle particulièrement importante
pour les \trdisz, dite \emph{coupure}, est la
suivante
%-----------------begin $$----------------
\begin{equation}\label{coupure1}
 \bigl(\,x\vi a\; \leq\;  b\,\bigr)\quad\&\quad  \bigl(\,a\; \leq\; x\vu  b\,\bigr)
\quad \Longrightarrow \quad a \leq\;  b.
\end{equation}
%-----------------end $$------------------
Pour la démontrer on écrit $ x\vi a\vi b=x\vi a$  et
$a= a\vi(x\vu b)$ donc
%-----------------begin $$----------------
$$ a=(a\vi x)\vu(a\vi b)=(a\vi x\vi b)\vu(a\vi b)=a\vi b
$$
%-----------------end $$------------------

%--- Notation{notaVupVda}--------------
\begin{notation}
\label{notaVupVda}
{\rm
Pour un \trdi $\gT$ on  note $A \vda B$ ou $A \vdash_\gT B$ la relation définie comme suit  sur l'ensemble $\Pfe(\gT)$:
%-----------------begin $$----------------

\snic{A \vda B \; \; \equidef\; \; \Vi A\;\leq \;
\Vu B.}
}
\end{notation}
%--- end-notation-----------------------------------------
Notez que la relation  $A \vdash B$ est bien définie sur  $\Pfe(\gT)$ parce que les lois $\vi$  et $\vu$  sont associatives,
commutatives et idempotentes.
Notez que
$\; \emptyset  \vda \{x\}$ implique  $x=1 $ et
que $\; \{y\} \vda \emptyset$ implique $y=0$.
Cette relation vérifie les axiomes suivants, dans lesquels on
écrit $x$ pour $\{x\}$ et $A, B$  pour $ A\cup B$.
%--------------------begin array---------------
$$\arraycolsep3pt\begin{array}{rcrclll}
&    & a  &\vda& a    &\; &(R)     \\[1mm]
A \vda B &   \; \Longrightarrow \;  & A,A' &\vda& B,B'   &\; &(M)     \\[1mm]
(A,x \vda B)\;\;
% H % \land
\&
\;\;(A \vda B,x)  &  \; \Longrightarrow \; & A &\vda& B &\;
&(T).
\end{array}$$
%---------------------end array--------------
\rdb
On dit que la relation est \emph{réflexive}, \label{remotr} \emph{monotone} et
\emph{transitive}.
La troisième règle (transitivité) peut être vue comme une
\gnn de la règle (\ref{coupure1}) et s'appelle \egmt la
règle de \emph{coupure}.
\index{coupure}

Signalons aussi les deux règles suivantes dites \emph{de \ditz}:
%--------------------begin array---------------
$$\arraycolsep3pt\begin{array}{rcl}
(A,\;x \vda B)\;\& \;(A,\;y \vda B)  &  \;  \Longleftrightarrow  \; &
A,\;x\vu y \vda B  \\[1mm]
(A\vda B,\;x )\;\&\;(A \vda B,\;y)  &  \; \Longleftrightarrow \; &
A\vda B,\;x\vi y
\end{array}$$
%---------------------end array--------------

Une manière intéressante d'aborder la question des \trdis définis
par \gtrs et relations est de considérer la relation
$A \vda B$ définie sur l'ensemble $\Pfe(\gT)$ des parties
finiment énumérées d'un \trdiz~$\gT$.
En effet, si $S\subseteq \gT$ engendre $\gT$ comme treillis, alors la
connaissance de la relation
$\vda$ sur $\Pfe(S)$ suffit à caractériser sans ambigüité
le treillis $\gT$, car toute formule sur $S$ peut être réécrite, au
choix, en \gui{forme normale conjonctive} (inf de sups dans~$S$) ou \gui{normale
disjonctive} (sup de infs dans~$S$). Donc si l'on veut comparer deux
\elts du treillis engendré par~$S$ on écrit le premier en forme
normale disjonctive, le second en forme normale conjonctive, et l'on
remarque~que
%-----------------begin $$----------------
$$ \Vu_{i\in I}\big(\Vi A_i \big)\; \leq \; \Vi_{j\in J}\big(\Vu B_j
\big)
\quad \Longleftrightarrow\quad  \Tt i\in I,\ \Tt j \in J,\;\;
A_i \vda  B_j
$$
%-----------------end $$------------------

%--- Definition{defEntrel}-------------
\begin{definition}
\label{defEntrel}
Pour un ensemble $S$ arbitraire, une relation sur  $\Pfe(S)$  qui est
réflexive, monotone et transitive est
appelée une {\em \entrelz} (en anglais, {\em entailment relation}).
\end{definition}
%--- end-definition------------------------------------

Le théorème suivant est fondamental. Il dit que les
trois propriétés des \entrels sont exactement ce qu'il faut pour que
l'interprétation en forme de \trdi soit adéquate.

%:    Theorem{thEntRel1}----------------
\begin{theorem}
\label{thEntRel1} {\rm  (\Tho fondamental des \entrelsz)} \\
Soit un ensemble $S$  avec une \entrel
$\vdash_S$ sur $\Pfe(S)$. On considère le \trdi $\gT$ défini par
\gtrs et relations comme suit: les \gtrs sont les
\elts de $S$ et les relations sont les
%-----------------begin $$----------------
$$ A\; \vdash_\gT \;  B
$$
%-----------------end $$------------------
chaque fois que $A\; \vdash_S \; B$.  Alors, pour tous $A$,
 $B$ dans $\Pfe(S)$,  on a
%-----------------begin $$----------------
$$  A\; \vdash_\gT \;  B
\; \Longrightarrow \; A\; \vdash_S \;  B.
$$
%-----------------end $$------------------
\end{theorem}
%--- end-theorem-----------------------------------------

\smallskip 
\rem La relation $x\vdash_S y$ est a priori une relation de préordre, et non une relation d'ordre, sur~$S$. Notons $\ov x$ l'\elt $x$ vu dans l'ensemble ordonné~$\ov{S}$ associé à ce préordre, et pour une partie~$A$ de~$S$ notons $\ov A=\sotq{\ov x}{x\in A}$.
Dans l'énoncé du \tho on considère un \trdi $\gT$ qui donne sur $S$ la même \entrel que $\vdash_S$.
En toute rigueur, on aurait dû noter $ \ov A\; \vdash_\gT \; \ov B$
plutôt que $A\; \vdash_\gT \;  B$ pour tenir compte du fait que l'\egt dans $\gT$ est plus grossière que dans $S$. En particulier c'est $\ov{S}$, et non pas $S$, qui s'identifie à une partie de $\gT$.
\eoe

%-----------------begin proof------------------
\begin{proof}
On donne une description explicite du \trdi $\gT$. Les
\elts de $\gT$ sont représentés par ceux de
$\Pfe\big(\Pfe(S)\big)$, i.e. des $X$ de la forme:

\snic{X=\{A_1,\dots,A_n\}}

%:2015
(intuitivement $X$ représente $\Vi_{i\in\lrbn}\Vu {A_i}$).
On définit alors de manière inductive la relation
 $A\preceq Y$ pour $A\in \Pfe(S)$ et $Y\in \Pfe\big(\Pfe(S)\big)$
 comme suit:
\begin{itemize}
\item Si $B\in Y$ et $B\subseteq A$ alors $A\preceq Y$.
\item Si l'on a $A\vdash_S y_1,\dots,y_m$ et $A,y_j\preceq Y$ pour
$j=1$, \ldots, $m$ alors $A\preceq Y$.
\end{itemize}
%(intuitivement, $\Vi A\leq \Vu_{C\in Y} \left(\Vi C\right) $).\\
On montre facilement que si $A\preceq Y$ et $A\subseteq A'$ alors on
a $A'\preceq Y.$ On en déduit que $A\preceq Z$ si $A\preceq Y$ et $B\preceq Z$
pour tout $B\in Y$. On peut alors définir $X\leq Y$ par
\gui{$A\preceq Y$ pour tout $A\in X$}. On vérifie enfin que $\gT$ est 
un treillis distributif{\footnote{Plus \prmtz, comme $\leq$ est seulement un  préordre, on prend pour $\gT$  le quotient de
$\Pfe\big(\Pfe(S)\big)$ par la relation d'\eqvcz: $X\leq Y$ et $Y\leq X$.}}
pour les opérations ($0$-aires et binaires)
%:2015 permuter 0 et 1
\begin{equation}\label{eqentrel}
\left.\arraycolsep2pt
\begin{array}{rclcrcl}
1  & =  &  \emptyset  &\qquad\quad   &  0  &  = & \so{\emptyset}    \\[1mm]
X\vi Y  &  = & X\cup Y && X\vu Y  &  = & \sotq{ A \cup B} {A\in X,~B\in Y}
\end{array}
\right|
\end{equation}
Pour ceci on montre que si $C\preceq X$ et $C\preceq Y$, alors on a
$C\preceq X\vi Y$ par induction sur les preuves de $C\preceq X$ et
$C\preceq Y$.
\\
On remarque que si $A\vdash_S y_1,\dots,y_m$ et $A,y_j\vdash_S B$
pour tout $j$, alors on obtient $A\vdash_S B$ 
en utilisant $m$ fois la règle de coupure. 
Il en résulte que si l'on a $A\vdash_\gT B$,
c'est à dire $A\preceq \{\{b\}~|~b\in B\}$, alors on a $A\vdash_S B$.
\end{proof}
%-----------------end proof------------------

\rem Si on considère une partie $S$ d'un \trdi $T$, si $S$ engendre $T$
comme \trdi et si $S$ ne contient ni $0_T$ ni $1_T$ alors le treillis construit au \thref{thEntRel1} est canoniquement isomorphe à $T$.
Mais si $S$ contient $0_T$ le treillis construit ajoute à $T$ pour la disjonction vide (correspondant à $X=\so\emptyset$) un \elt minimum strictement plus petit que~$0_T$.
La remarque symétrique vaut dans le cas où $S$ contient $1_T$. \eoe

%:     Corollary{corthEntRel1}
\begin{corollary}\label{corthEntRel1} \emph{(Treillis distributifs \pfz)}
\begin{enumerate}
\item Un \trdi librement engendré par un ensemble $\,E\,$ fini est fini.
\item Un \trdi \pf est fini.
\end{enumerate}
\end{corollary}
\begin{proof}  %\cap au lieu de \cup
\emph{1.} On considère la relation implicative minimale sur $E$. 
%:2015
Elle est définie par  
$$\preskip-.2em \postskip.2em
(A,B\in \Pfe(E))\quad A\vdash_E B \;\equidef  \;\exists x\in A\cap B.\;\;\;\;\;\phantom{.}
$$
On considère alors le \trdi correspondant à cette \entrel via le
\thrf{thEntRel1}.
Il est isomorphe à un sous-ensemble %de $G(\Pfe(E))$
de $\Pfe\big(\Pfe(E)\big)$,
celui qui est représenté par les listes
$(A_1,\ldots,A_k)$ dans $\Pfe(E)$ telles que deux $A_i$ d'indices distincts sont incomparables pour l'inclusion. Les lois sont obtenues à partir 
de~(\ref{eqentrel}), en simplifiant les listes obtenues lorsqu'elles ne satisfont
pas le critère d'incomparabilité.

\emph{2.} Si l'on impose un nombre fini de relations entre les \elts de $E$, on doit passer à un treillis quotient du \trdi libre sur $E$. La relation d'\eqvc
engendrée par ces relations et compatible avec les lois de treillis
est décidable parce que la structure est définie en utilisant
seulement un nombre fini d'axiomes.
\end{proof}\rdb

\rems \label{remagBlibre}
1) Une autre preuve du point \emph{1} pourrait être la suivante.
L'\agB librement engendrée par le \trdi $\gT$ librement engendré par $E$
est l'\agB $\gB$ librement engendrée par $E$. Cette dernière peut être facilement décrite par les \elts de $\Pfe\big(\Pfe(E)\big)$, sans~aucun passage au quotient:
la partie $\{A_1,\dots,A_n\}$ représente intuiti\-vement
$\Vu_{i\in\lrbn}\left(\Vi {A_i}\vi \Vi {A'_i}\right)$, en désignant par $A'_i$ la partie de $E$ formée par les~$\lnot x$
pour les $x\notin A_i$. Donc $\gB$ possède $2^{2^{\#E}}$ \eltsz. Enfin
on a vu que $\gT$ s'identifie à un sous-\trdi de $\gB$ (\thrf{thBoolGen}).\\
2) La preuve donnée du point \emph{2} utilise un argument tout à fait \gnlz. Dans le cas des \trdis on peut plus \prmt se reporter à la description des quotients donnée \paref{eqPreceq}.
\eoe

%%%%%%%%%%%%%%%%%%%%%%%%%%%%%%%%%%%%%%%%%%%%%%%%%%%%%%%%%%%%%%%%%
%: \subsec{Dualité}
\subsect{Dualité entre \trdis finis et ensembles\\
ordonnés finis}{Dualité}\label{SpecTrdiFi}

Si $\gT$ est un \trdi on note
$\SpecT\eqdefi\Hom(\gT,\Deux)$. C'est un ensemble ordonné appelé
\emph{spectre (de Zariski) de $\gT$}.
Un \elt $\alpha$ de $\SpecT$ est \care par
son noyau. En \clama un tel noyau est appelé un \idepz.
Du point de vue \cof il doit être détachable.
Nous sommes intéressés ici par le cas où $\gT$ est fini, ce qui implique
que~$\SpecT$ est \egmt fini (au sens \cofz).%
\index{spectre!d'un \trdiz}%
\index{ideal@idéal!premier}%
\index{premier!idéal --- d'un treillis distributif}

Si $\varphi:\gT\to\gT'$ est un \homo de \trdis  et si $\alpha\in\Spec\gT'$,
alors  $\alpha\circ \varphi\in\SpecT$.
Ceci définit une application croissante de~$\Spec\gT'$ vers~$\SpecT$, 
notée $\Spec\alpha$, dite \gui{duale} de~$\varphi$.

Inversement soit $E$ un ensemble ordonné fini. On note $E\sta$ l'ensemble des
\emph{sections initiales} de $E$, i.e., l'ensemble des parties finies de $E$
stables pour l'opération $x\mapsto \dar x$. Cet ensemble, ordonné par la relation $\supseteq$, est un \trdi fini,  un sous-treillis du treillis $\Pf(E)\eci$ (le treillis opposé à $\Pf(E)$).

%:     Fact{factTrdiFiniDual}
\begin{fact}\label{factTrdiFiniDual}
Le nombre d'\elts d'un ensemble ordonné fini $E$ est égal à la longueur maximum d'une chaîne
strictement croissante d'\elts de~$E\sta$.
\end{fact}
\begin{proof} Il est clair qu'une chaîne strictement monotone d'\elts
de $E\sta$ (donc de parties finies de $E$) ne peut avoir plus que $1+\# E$ \eltsz. Sa \gui{longueur} est donc
$\leq \#E$. Concernant l'in\egt opposée,  
on la vérifie pour $E=\emptyset$ (ou pour un singleton), puis on fait une \recu sur $\#E$, en regardant un ensemble ordonné à $n$ \elts ($n\geq1$) comme un ensemble ordonné \hbox{à $n-1$}
\elts que l'on étend en rajoutant un \elt maximal.
\end{proof}

Si $\psi:E\to E_1$ est une application croissante entre ensembles ordonnés finis, alors pour tout $X\in E_1\sta$, $\psi^{-1}(X)$ est un \elt de
$E\sta$. Ceci définit un \homo $E_1\sta\to E\sta$ noté
$\psi\sta$, dit \gui{dual} de $\psi$.

%:     Theorem{thDualiteFinie}
\begin{theorem}\label{thDualiteFinie}\emph{(Dualité entre ensembles ordonnés finis et \trdis finis)}
\begin{enumerate}
\item Pour tout ensemble ordonné fini $E$ définissons $\nu_E:E\to\Spec(E\sta)$
par

\snic{\nu_E(x)(S)=0 \;\;\mathsf{si}\;\; x\in S,\,1\;\;\mathsf{sinon}.}

%\sni
Alors, $\nu_E$ est un \iso d'ensembles ordonnés. En outre,
pour toute application croissante  $\psi:E\to E_1$,
on a
$\nu_{E_1}\circ \psi=\Spec(\psi\sta)\circ \nu_E.$
\item Pour tout \trdi fini $\gT$ définissons $\iota_\gT:\gT\to(\SpecT)\sta$ par

\snic{\iota_\gT(x)=\sotq{\alpha \in \SpecT}{\alpha(x)=0}.}

%\sni
Alors, $\iota_\gT$ est un \iso de \trdisz. En outre,
pour tout morphisme $\varphi:\gT\to \gT'$,
on a
${\iota_{\gT'}\circ \varphi=(\Spec\varphi)\sta\circ \iota_\gT.}$

%
%\item En outre \ldots.
%
\end{enumerate}
\end{theorem}
\begin{proof}
Voir l'exercice \ref{exoTreillisDistributifFini}.
\end{proof}

En d'autres termes, les catégories des \trdis finis et
des ensembles ordonnés finis
sont anti\eqvesz. L'anti\eqvc est donnée par les foncteurs
contravariants $\Spec\bullet$
et $\bullet\sta$, et par les transformations naturelles~$\nu$ et~$\iota$
définies ci-dessus.

La \gnn de cette anti\eqvc de catégories pour le cas des \trdis
non \ncrt finis sera abordée brièvement \paref{SpecTrdi}.

%La \dem est laissée \alecz.

%%%%%%%%%%%%%%%%%%%%%%%%%%%%%%%%%%%%%%%%%%%%%%%%%%%%%%%%%%%%%%%%%
%:--- subsec{Algèbres de Heyting}
\subsec{Algèbres de Heyting}
\label{secagH}
%-----------------------------------------

%%: \subsec{Algèbres de Heyting}%%%%%%%%%%%%
%\subsec{Algèbres de Heyting}

Un \trdi $\gT$ est appelé un \ix{treillis implicatif}
ou une {\em \agHz}  lorsqu'il existe
une opération binaire  $\im$ vérifiant pour tous $a,\,b,\,c$:
\index{algèbre!de Heyting}\index{Heyting!algèbre de ---}
%---  equation eqAgHey --------
\begin{equation}\label{eqAgHey}
a\vi b \leq c \;\;\Longleftrightarrow \;\; a \leq  (b\im c)\,
\end{equation}
%---------------------end equation--------------
Ceci signifie que pour tous $b$, $c\in\gT$, l'\emph{\id transporteur}
\index{ideal@idéal!transporteur}\label{NOTATransp2}
\index{transporteur!d'un idéal dans un autre}
$$
(c:b)_\gT \eqdefi \sotq{x\in\gT}{x\vi b\leq c}
$$
est principal,
son \gtr
étant noté $b\im c$.
Donc si elle existe, l'opération~$\im$ est déterminée de
manière unique par la structure du treillis.
On définit alors la loi unaire  $\neg x = x\im 0$.
La structure d'\agH peut être définie comme purement
équationnelle en
donnant de bons axiomes, décrits dans le fait suivant.

%:     Fact{fact1AGH}
\begin{fact}\label{fact1AGH}
Un treillis $\gT$ (non supposé
distributif) muni d'une loi $\im$ est une \agH \ssi les axiomes
suivants sont vérifiés:

\vspace{-3mm}
%--------------------begin array---------------
$$\arraycolsep2pt\begin{array}{rcl}
a\im a&=   &1    \\
a\vi(a\im b)&=   &a\vi b    \\
b\vi(a\im b)&=   & b   \\
a\im(b\vi c)&=   &(a\im b)\vi(a\im c)
\end{array}$$
%---------------------end array--------------
\end{fact}

Notons aussi les faits importants suivants.

%:     Fact{fact2AGH}
\begin{fact}\label{fact2AGH} Dans une \agH on~a:
$$a\leq b \;\iff\;  a\im b =1
$$
\vspace{-4mm}
%--------------------begin array---------------
$$\arraycolsep2pt\begin{array}{rclcrcl}
a\im (b\im c) &=& (a\vi b) \im c, &&
         a\im b &\leq& \neg b\im \neg a,
\\
(a\vu b)\im c &=& (a\im c)\vi(b\im c), &&
 a&\leq    & \neg\neg a,
         \\
\neg\neg\neg a&=    &\neg a, &&
a\im b &\leq& (b\im c) \im (a\im c),~~
             \\
\neg(a\vu b)&=   & \neg a\vi \neg b,&&
\neg a\vu b&\leq    & a\im b.
\end{array}$$
%---------------------end array--------------
\end{fact}

  Tout \trdi fini est une \agHz, car tout \itf est principal.
 Un cas particulier important d'\agH est une \agBz.

 Un \emph{\homo d'\agHsz} est un \homo de \trdis $\varphi :\gT\to\gT'$
 tel que $\varphi(a\im b)=\varphi(a)\im\varphi(b)$ pour
tous $a$, $b\in\gT$.

Le fait suivant est immédiat.

%--- Fact{factQuoAgH}----------------
\begin{fact}
\label{factQuoAgH}
Soit $\varphi:\gT\to\gT_1$ un \homo de \trdisz, avec~$\gT$ et
$\gT_1$ des \agHsz. Notons $a\preceq b$ pour
$\varphi(a)\leq_{\gT_1}\varphi(b)$. Alors
$\varphi$ est un \homo d'\agHs \ssi on a pour tous $a$, $a'$, $b$, $b'\in\gT$:

\snic{
a\preceq a'\Longrightarrow (a'\im b)\preceq(a\im b) , \quad
\hbox{et}\quad
b\preceq b'\Longrightarrow (a\im b)\preceq(a\im b').
}
\end{fact}
%--- end-fact-----------------------------------------

%--- Fact{factQuoAgH2}---------------
\begin{fact}
\label{factQuoAgH2}
Si $\gT$ est une \agH tout quotient $\gT/(y=0)$ (\cad tout quotient
par un \id
principal) est aussi une \agHz.
\end{fact}
%--- end-fact-----------------------------------------
%-----------------begin proof------------------
\begin{proof}
Soit $\pi:\gT\to\gT'=\gT/(y=0)$ la projection canonique. On a
\[ \arraycolsep4pt
\begin{array}{ccccccccccccc} 
\pi(x)\vi\pi(a)\,\leq_{\gT'}\, \pi(b)  &  \Longleftrightarrow 
&  \pi(x \vi a)\,\leq_{\gT'}\, \pi(b)  &  \Longleftrightarrow  \\[1mm]
  x\vi a \,\leq\, b\vu
y  & \Longleftrightarrow   &  
x\,\leq\, a\im(b\vu y).  
 \end{array}
\]
 Or $y\,\leq\, b\vu y\,\leq\,a\im(b\vu y)$, 
 donc 
 
 \snic{\pi(x)\vi\pi(a)\,\leq_{\gT'}\,\pi(b)\;\Longleftrightarrow\;
x\,\leq\, \big(a\im(b\vu y)\big)\vu y,}
 
 %\sni
\cad $\pi(x)\leq_{\gT'}\pi\big(a\im(b\vu y)\big)$, ce
qui montre que $\pi\big(a\im(b\vu y)\big)$ vaut pour $\pi(a)\im\pi(b)$ dans
$\gT'$.
\end{proof}
%-----------------end proof------------------

%-% PERSO
\perso{

1.  Cependant il ne semble pas que $\pi$ soit en \gnl un \homo
d'\agHsz.

2.  Il serait bon d'avoir un exemple d'un treillis quotient d'une \agH
qui ne serait pas une \agHz.

3.  De manière \gnlez, il serait bon d'avoir des exemples
variés d'\agH non \noees à notre disposition.

}
%-% Fin PERSO

\rems
1) La notion d'\agH est reminiscente de la notion d'\cori en \alg commutative.  En effet, un \cori peut
être \care comme suit: l'intersection de deux \itfs est un
\itf et le transporteur d'un \itf dans un \itf est un \itfz.  Si l'on
\gui{relit} ceci pour un \trdi en se rappelant que tout \itf est
principal on obtient une \agHz. 

 2) Tout \trdi $\gT$ engendre une \agH de façon
naturelle.  Autrement dit on peut rajouter formellement un \gtr pour
tout \id $(b:c)$.  Mais si l'on part d'un \trdi qui se trouve être une
\agHz, l'\agH qu'il engendre est strictement plus grande.  Prenons par
exemple le treillis $\Trois$ qui est le \trdi libre à un \gtrz.
L'\agH qu'il engendre est donc l'\agH libre à un \gtrz.  Or
celle-ci est infinie (cf.  \cite{Johnstone}).  A~contrario le treillis
booléen engendré par $\gT$ (cf.\  \thrf{thBoolGen}) reste égal à~$\gT$
lorsque celui-ci est booléen. \eoe

%:2018 ajout d'une section
% --- sec{Constructions de \trdisz}
\section{Constructions de \trdisz}
\label{secconstrdi}
%-----------------------------------------

La description d'un \trdi via la \entrel qu'il induit sur un \sgr simplifie de nombreuses constructions usuelles en algèbre.

Étant données deux \entrels $\vdi1$ et $\vdi2$ sur un ensemble $S$, nous dirons que $\vdi1$ est \emph{plus forte que} $\vdi2$ si $A\vdi1 B$ implique $A\vdi2 B$. Nous dirons aussi que $\vdi2$ est \emph{plus lâche que} $\vdi1$. 

%: --- subsec Treillis distributifs quotients
\subsec{Treillis distributifs quotients}
Nous commençons par la construction de quotients d'un \trdiz.
Le contexte général est le suivant:
$\gT$ est un \trdiz, $S$ un \sgr de $\gT$, $a,b,x,y\in S$ et $A,B,X,Y\in\Pfe(S)$. On note $\vda$ la \entrel de~$\gT$ (restreinte à $\Pfe(S)$). 
On considère un treillis quotient $\gT'$. On note $\,\vdash'\,$ la \entrel de $\gT'$ restreinte à $\Pfe(S)$.
Notons que comme cas particulier important on a  $S=\gT$.

Le fait que $\gT'$ est un quotient de $\gT$ sous certaines relations imposées, et le fait que la structure de $\gT'$ est entièrement déterminée par la restriction \hbox{de $\vdi{\gT'}$} à $S$ implique que la \entrel  $\,\vdash'\,$ est la plus forte des \entrels sur $S$ plus lâches que $\vda$ et pour lesquelles  les conditions imposées sont réalisées.

Nous démontrons directement dans chaque cas que la \entrel construite résout bien le problème universel associé à la structure quotient.

%f
%:     lemma{lemquotrdi0}
\begin{lemma} \label{lemquotrdi0} 
On considère le treillis quotient $\gT'$ obtenu en forçant les relations $\Vi A=0$ et $\Vu B=1$.\\ 
Alors $X\,\vdash'\,Y$ \ssi $X,B\vda Y,A$.
\end{lemma}
%----------- fin lemma ----------------------------------------- 
%
\begin{proof}
On vérifie facilement que la relation $X\,\vdash'\, Y$ définie par $X,B\vda Y,A$ est une \entrel plus lâche que $\vda$ qui satisfait  $A\,\vdash'\, 0$ et $1\,\vdash'\,B$.
Il suffit donc de vérifier que $X\,\vdash'\, Y$ implique $X\vdi1Y$ si $\vdi1$
est une \entrel plus lâche que~$\vda$ qui satisfait $A\vdi1 0$ et $1\vdi1B$.
Or si $X\,\vdash'\, Y$ alors $X,B\vdi1Y,A$ et puisque $A\vdi1 0$ et $1\vdi1B$, cela donne $X\vdi1Y$.
\end{proof}

Le lemme précédent peut être considéré comme une variante de la 
proposition~\ref{propIdealFiltre}.
%c
%:     lemma{lemquotrdi1}
\begin{lemma} \label{lemquotrdi1}
On considère le treillis quotient $\gT'$ obtenu en forçant la relation $a\leq b$.  
Alors $X\,\vdash'\,Y$ \ssi $X\vda Y,a$ et   $X,b \vda Y$.
\end{lemma}
%--------- fin lemma ------------------------------- 
%
\begin{proof}
On vérifie facilement que la relation $X\,\vdash'\, Y$ définie par $X\vda Y,a$ 
\hbox{et $X,b \vda Y$} est une \entrel plus lâche que $\vda$ et qu'elle satisfait  $a\,\vdash'\, b$.
Il suffit donc de vérifier que $X\,\vdash'\, Y$ implique $X\vdi1Y$ si $\vdi1$
est une \entrel plus lâche que~$\vda$ qui satisfait $a\vdi1 b$.
Or si $X\,\vdash'\, Y$ \hbox{alors $X,b\vdi1Y$} \hbox{et $X\vdi1Y,a$}, et puisque $a\vdi1b$, cela donne après deux coupures $X\vdi1Y$.
\end{proof}
%

%:     proposition{lemquotrdi2}
\begin{proposition} \label{lemquotrdi2}
On considère le treillis  $\gT'$ quotient de $\gT$ obtenu en forçant la relation $A\,\vdash'\, B$. 
 \Propeq
\begin{enumerate}
\item $X\,\vdash'\,Y$
\item Pour chaque $a\in A$ et chaque $b\in B$, on a: $X,b\vda Y$ et $X\vda Y,a$.
\end{enumerate}
\end{proposition}
%--------- fin lemma ------------------------------- 
\begin{proof}
On vérifie facilement que la relation $X\,\vdash'\, Y$ définie par le point \emph{2} est une \entrel plus lâche que $\vda$ et qu'elle satisfait  $A\,\vdash'\, B$.
Il suffit donc de vérifier que $X\,\vdash'\, Y$ implique $X\vdi1Y$ si $\vdi1$
est une \entrel plus lâche que~$\vda$ qui satisfait $A\vdi1 B$.
Or si $X\,\vdash'\, Y$ alors $X,b\vdi1Y$ et $X\vdi1Y,a$ pour chaque $a\in A$ et $b\in B$, et puisque $A\vdi1B$, cela donne après de nombreuses coupures $X\vdi1Y$.
\end{proof}

%: --- subsec construction alternative  de l'\agB engendrée par un \trdi
\subsec{\AgB engendrée par un \trdi}

Une construction alternative  de l'\agB engendrée par un \trdi (obtenue dans le \thref{thBoolGen}) est donnée dans la proposition suivante.

%--- Proposition{propTrBoo}------------
\begin{proposition}
\label{propTrBoo}
Soit $\gT$ un \trdi avec un \sgrz~$S$. Le treillis booléen
librement engendré par $\gT$ peut être décrit comme le \trdi engendré par
l'ensemble $S_1=S\cup\overline{S}$ (où $\overline{S}$ est une copie de~$S$ disjointe de $S$) muni de la \entrel $\,\vdash'\,$ définie comme suit:
si $A,B,A',B'\in\Pfe(S)$  on a
%-----------------begin $$----------------
$$
A,\overline{B} \,\vdash'\, A',\overline{B'}\;\equidef\; A,B'\vda
A',B\;
{\rm dans} \;  \gT
$$
%-----------------end $$------------------
Si on note $\gT\bul$ ce treillis booléen, l'ensemble $\gT_1=\gT\cup\overline{\gT}$ s'injecte naturellement dans $\gT\bul$. En particulier le morphisme naturel de $\gT$ dans $\gT\bul$ est injectif, il identifie $\gT$ à son image dans~$\gT\bul$.
\end{proposition}
%--- end-proposition----------------------------------------
%-----------------begin proof------------------
\begin{proof}
Si l'on note $\ov x$ le complément de $x$ dans un \trdi booléen, il est
clair que le treillis booléen $\gT\bul$ engendré par $\gT$, qui existe en vertu de considérations générales liées à l'\alg universelle, est lui-même engendré par $S\cup\ov S$. La \entrel  sur $S\cup\ov S$  qui définit $\gT\bul$ est donc la plus forte des \entrels plus lâches que $\vda$ (restreinte à $S$) qui font du \trdi correspondant une \agB pour laquelle $x\in S$ et $\ov x$ sont des \elts complémentaires. \\
On vérifie facilement que la relation $A,\overline{B} \,\vdash'\, A',\overline{B'}$ définie dans l'énoncé est une \entrel plus lâche que $\vda$ (restreinte à $S$) et qu'elle fait du \trdi correspondant une \agBz:
le complément \hbox{de $x\in S$} est donné \hbox{par $\ov x\in \ov S$} car $x,\ov x \,\vdash'\, 0$ et $1\,\vdash'\, x,\ov x$ (le lemme \ref{lemagbgen} nous dit alors que le \trdi correspondant est une \agBz).
Il suffit donc de vérifier \hbox{que $A,\overline{B} \,\vdash'\, A',\overline{B'}$} implique $A,\overline{B} \vdi1 A',\overline{B'}$ si $\vdi1$
est une \entrel plus lâche que~$\vda$ qui fait du \trdi correspondant une \agB  pour laquelle $x\in S$ et $\ov x$ sont des \elts complémentaires.
Or si $A,\overline{B} \,\vdash'\, A',\overline{B'}$ alors \hbox{on a $A,B'\vdi1A',B$} et puisqu'on a affaire à une \agBz, comme $x\leq y$ équivaut à $x\vi\ov y=0$, cela équivaut à $A,B',\ov{A'},\ov B\vdi10$  \hbox{et à $A,\overline{B}\vdi1A',\overline{B'}$}.
\end{proof}
%-----------------end proof------------------

%l
%:     Lemma{lemagbgen}
\begin{lemma} \label{lemagbgen}
Soit $\gB$ un  \trdi avec un \sgr $G$.
Si tout \elt de $G$ possède un complément dans $G$, $\gB$ est une \agBz. 
\end{lemma}
%----------- fin lemma ----------------------------------- 
%
\begin{proof}
Un \elt arbitraire de $\gB$ s'écrit $\Vu_j a_j=\Vu_j \Vi A_j$ pour \hbox{des $A_j\in\Pfe(G)$}. Montrons d'abord que $a_j= \Vi A_j%=\Vi_{b\in A_j} b
$ admet $\Vu_{b\in A_j} \ov{b}$ pour complément. \\
Posons $B=\Vi A_j$. On a par distributivité 
$$\Vi B \vi \Vu_{b\in B} \ov b = \big(\Vi_{c\in B} c\big) \vi \Vu_{b\in B} \ov b=
\Vu_{b\in B} \big(\ov b\vi \Vi_{c\in B} c\big)=0$$
et
$$\Vi B \vu \big(\Vu_{b\in B} \ov b\big) = \Vi_{c\in B} c \vu \big(\Vu_{b\in B} \ov b\big)=
\Vi_{c\in B} \big(c\vu \Vu_{b\in B} \ov b\big)=1.$$
On montre ensuite de manière symétrique que $\Vu_j a_j$ admet $\Vi_j \ov{a_j}$
pour complément.
\end{proof}
%
%: --- subsec Somme directe de deux \trdis
\subsec{Somme directe de deux \trdis}

%--- Proposition{propSumtrdi}------------
\begin{proposition}
\label{propSumtrdi}
Soient $\gT_1$ et $\gT_2$ deux \trdis avec des \sgrsz~$S_1$ et $S_2$. Le treillis somme directe de $\gT_1$ et $\gT_2$ peut être décrit comme le \trdi engendré par
l'ensemble $S=S_1\cup S_2$ (où~$S_2$ est supposé disjoint de $S_1$) muni de la \entrel $\vdi0$ définie comme suit:
si $A_1,B_1\in\Pfe(S_1)$ et $A_2,B_2\in\Pfe(S_2)$  on a
%-----------------begin $$----------------
$$
A_1,A_2 \vdi0 B_1,B_2\;\equidef\; A_1\vdi1
B_1\;\;\mathrm{ou}\;\;
A_2\vdi2
B_2
$$
%-----------------end $$------------------
Si on note $\gT_0$ ce treillis, et si $\gT_1$ et $\gT_2$ sont non triviaux, les morphismes naturels de $\gT_1$ et $\gT_2$ dans $\gT$  sont injectifs.
\end{proposition}
%--- end-proposition----------------------------------------
%
\begin{proof}
Il est
clair que le treillis  $\gT_0$ somme directe  $\gT_1$ et $\gT_2$, qui existe en vertu de considérations générales liées à l'\alg universelle, est lui-même engendré par $S_1\cup S_2$. La \entrel $\vdi0$ sur $S_1\cup S_2$  qui définit $\gT_0$ est donc la plus forte des \entrels plus lâches que $\vdi1$ (restreinte à $S_1$) et~$\vdi2$ (restreinte à $S_2$).
On constate facilement que la relation $\vdi0$ donnée dans est une \entrel sur $S$ plus lâche que $\vdi1$  et $\vdi2$. Il suffit donc de démontrer que la relation $A_1,A_2 \vdi0 B_1,B_2$ implique
$A_1,A_2 \vdi3 B_1,B_2$ si~$\vdi3$ est une \entrel sur $S$ plus lâche que $\vdi1$  et $\vdi2$. Ceci résulte de façon immédiate de la monotonie des \entrelsz.
\end{proof}
%

%r
%:     Remark{rempropSumtrdi}
\begin{remark} \label{rempropSumtrdi} 
{\rm   Lorsque les treillis sont non triviaux, si l'on identifie $\gT_1$ et~$\gT_2$ à des sous-treillis de $\gT_0$
leur intersection est réduite à $\so{0,1}$. Si l'un des treillis est trivial (réduit à $0$),  il tue $\gT_0$ ainsi que l'image de $\gT_2$ dans~$\gT_0$.
Cette situation est analogue à celle de la catégorie des \Algsz. 
La somme directe de $\gB_1$ et $\gB_2$, représentée par $\gB_1\otimes_\gA\gB_2 $ est nulle si $\gB_1$ ou $\gB_2$ est nulle%, comme le montrent les \egts $0\otimes 0=1\otimes 0=0\otimes 1$
. 
}\end{remark}
%----------- fin remark ---------------------------------- 

%%%%%%%%%%%%%%%%%%%%%%%%%%%%%%%%%%%%%%%%%%%%%%%%%%%%%%%%%%%%%%%%%
%:section: Exercices
\Exercices

%--- Exercise{exo3Lecteur}-------------
\begin{exercise}
\label{exoTrdiLecteur}
{\rm  Il est recommandé de faire les \dems non données, esquissées,
laissées \alecz,
etc\ldots
\, On pourra notamment traiter les cas suivants. \perso{à compléter}
\begin{itemize}
\item \label{exotreillisquotient}
Montrer que les relations \pref{eqPreceq} \paref{eqPreceq} sont exactement ce qu'il faut
pour définir un treillis quotient.
\item \label{exopropIdealFiltre} Démontrer la proposition \ref{propIdealFiltre}.
\item \label{exocorlemRajouCompl}
Démontrer le corolaire \ref{corlemRajouCompl}.
\item
Démontrer les faits \ref{factBezGCD}, \ref{factLocaliseGCD},
\ref{factSouspgcdsat} et \ref{factAXiclgcd}.
\item Démontrer le fait \ref{fact1Zar} et  tous les faits numérotés entre \ref{fact2Zar} et \ref{factQuoFIID} (pour le fait~\ref{factZar} voir l'exercice
\ref{exoTreillisZariski}).
\item \label{exoexlCLZED} Démontrer ce qui est affirmé dans les exemples \paref{exlCLZED}.
\item \label{exofactAGH} Démontrer les faits \ref{fact1AGH}
et \ref{fact2AGH}.
\end{itemize}
}
\end{exercise}
%--- end -exercise-----------------------------------------

%--- Exercise{exoRegleCoupure}-------------
\begin{exercise}
\label{exoRegleCoupure}
{\rm
Soit $\gT$ un \trdi et $x\in\gT$. On a vu (lemme \ref{lemRajouCompl}) que

\snic{\lambda_x:\gT\to\gT[x\bul] \eqdefi\gT\sur{(x=0)}\times \gT\sur{(x=1)}}

%\sni
est injectif, ce qui signifie: si $y\vi x=z\vi x$ et $y\vu x=z\vu x$, 
alors $y =z$.
\\  
Montrer que l'on peut en déduire la règle de coupure (\ref{coupure1}).
}
\end {exercise}
%--- end -exercise-----------------------------------------

%--- Exercise{exo0GpRtcl}-------------
\begin{exercise}
\label{exo0GpRtcl}
{\rm  Soit $\gA$ un anneau intègre et $p$, $a$, $b \in \Reg(\gA)$, avec $p$ irréductible. On suppose que
$p \divi ab$, mais $p \nedivi a$, $p \nedivi b$. Montrer que
$(pa, ab)$ n'a pas de pgcd.
Montrer que dans $\ZZ[X^2,X^3]$ les \elts $X^2$ et $X^3$ admettent un pgcd,
mais pas de ppcm, et que les \elts $X^5$ et $X^6$ n'ont pas de pgcd.
}
\end{exercise}
%--- end -exercise-----------------------------------------

%-% ENTRE NOUS
\entrenous{ Un exo à écrire?

%--- Exercise{exo1GpRtcl}-------------
%\begin{exercise}
\label{exo1GpRtcl} \textbf{\dfn équationnelle des \grlsz}
{  Décrire divers \syss d'axiomes que doit vérifier
$\vi$ pour définir
une structure de \grl sur un groupe $(G,0,+,-)$
}
%\end{exercise}
%--- end -exercise----------------------------------------- 
}
%-% Fin ENTRENOUS

%--- Exercise{exo2GpRtcl}-------------
\begin{exercise}
\label{exo2GpRtcl} (Autre \dfn des \grlsz)\\
{\rm   Montrer que les axiomes que doit vérifier une partie $G^+$ d'un groupe
$(G,0,+,-)$ pour définir un ordre compatible réticulé sont:
\begin{itemize}
\item $G=G^+-G^+$,
\item $G^+\cap -G^+=\so{0}$,
\item $G^+ + G^+\subseteq G^+$,
\item $\forall a,b\;\exists c,\;\;c+G^+=(a+G^+)\cap (b+G^+)$.
\end{itemize}
}
\end{exercise}
%--- end -exercise-----------------------------------------

%--- Exercise{exoLemmeGaussAX}-------------
\begin{exercise}
\label{exoLemmeGaussAX}
{(Une autre preuve du lemme de Gauss)}\\
{\rm
Dans le contexte de la proposition \ref{propLG}, montrer que
$\G(fg)=\G(f)\G(g)$ à l'aide d'une \dem basée sur le lemme de \DKM
\ref{lemdArtin}.
}
\end{exercise}
%--- end -exercise-----------------------------------------

%--- Exercise{exoKroneckerTrick}-------------
\begin{exercise}\label{exoKroneckerTrick}
{(L'astuce de \KRAz)}
{\rm  
Soit $d$ un entier fixé $\ge 2$.
\\
\emph {1.}
Soit $\AuX_{<d} \subset \AuX = \gA[\Xn]$ le sous-\Amo constitué des
\pols $P$ tels que $\deg_{X_i} P < d$ pour tout $i \in \lrbn$, 
et $\gA[T]_{<d^n}\subset \gA[T]$ celui formé par les \pols $f \in \gA[T]$ de degré $< d^n$.
\\
Montrer que $\varphi : P(\Xn) \mapsto P(T, T^d, \ldots, T^{d^{n-1}})$ induit
un \iso de \Amos entre les \Amos $\AuX_{<d}$ et $\gA[T]_{<d^n}$.

\emph {2.}
On suppose $\AX$ factoriel. Soit $P \in \AuX_{<d}$ et $f = \varphi(P)
\in \AT_{d^n}$. Montrer que toute factorisation de $P$ dans
$\AuX$ peut être retrouvée par une procédure finie à partir de celles de $\varphi(P)$  dans $\AT$.
}
\end {exercise}
%--- end -exercise-----------------------------------------

%--- Exercise{exoTreillisZariski}-------------
\begin{exercise}
\label{exoTreillisZariski}
{\rm
Vérifier le fait \ref{factZar}, i.e.  $\ZarA$ est un \trdiz.
Montrer que ce \trdi peut être défini par \gtrs et relations comme suit. Les \gtrs sont  les symboles
$\rD(a)$, $a \in \gA$, avec le \sys de  relations:

\snic {
\rD(0) = 0,\quad  \rD(1) = 1,\quad
\rD(a + b) \le \rD(a) \vu \rD(b),  \quad \rD(ab) = \rD(a) \vi \rD(b).
}
}
\end {exercise}
%--- end -exercise-----------------------------------------

%--- Exercise{exoClosedCover1}-----------
\begin {exercise}\label{exoClosedCover1}
{\rm
Le contexte est celui du principe de recouvrement fermé~\ref{prcf2}. On considère un recouvrement fermé
de l'anneau $\gA$ par des \ids \hbox{$\fa_1$, \ldots, $\fa_r$}.  On ne suppose pas que
$\prod_i\fa_i=0$, mais on suppose que chaque $\fa_i$ est \tfz.  Montrer qu'un \Amo $M$  est \tf \ssi il est \tf modulo chaque~$\fa_i$.
}
\end {exercise}
%--- end-exercise-----------------------------------------

%--- Exercise{exoAbul}-------------
\begin{exercise}\label{exoAbul}
 (L'anneau $\Abul$) {\rm  On se situe en \clamaz.
\\
  Soit $\gA$ un anneau et $\varphi:\gA\to\Abul$ l'\homo naturel.

\emph{1.} Montrer  que l'application
$\Spec\varphi:\Spec\Abul\to\Spec\gA$ est une bijection et que pour
$\fq\in\Spec\Abul$, l'\homo naturel
$\Frac(\gA\sur{\varphi^{-1}(\fq)})\to\Abul\sur\fq$ est un \isoz.

\emph{2.} L'anneau 
$\Abul$ s'identifie au sous-anneau \zedr de 

\snic{\wi\gA\eqdefi\prod_{\fp\in\SpecA}\Frac(\gA\sur\fp)}

%\sni
engendré par (l'image de) $\gA$.
}

\end {exercise}

%--- Exercise{exoMinA}-------------
\begin{exercise}\label{exoMinA} (Idéaux premiers minimaux)
\\
{\rm  On se situe en \clamaz. Un \idep est dit minimal s'il est minimal parmi les \idepsz. On note $\Min\gA$ le sous-espace de $\SpecA$ formé par les \idemisz. 
Rappelons que l'on a défini un \emph{filtre maximal} comme un filtre dont le localisé est un \alo \zedz. Dans le point \emph{1} de cet exercice on fait le raccord avec la \dfn plus usuelle.

\emph{1.}
Montrer qu'un
filtre strict $\ff$ est maximal parmi les filtres stricts  \ssi pour tout $x\notin \ff$ il existe $a\in\ff$ tel que $ax$ est nilpotent.
Une autre \carn possible est que l'anneau localisé $\ff^{-1}\gA$ est local,  \zed et non trivial.
 En particulier, tout filtre strict maximal parmi les filtres stricts est premier. 

NB:  reformulation de la première \prt \cara pour l'\idep \copz:  un \idep $\fp$ est minimal \ssi
pour tout~\hbox{$x\in \fp$}, il existe $a\notin\fp$ tel que $ax$ est nilpotent.

\emph{2.} La notion duale du radical de Jacobson est le filtre intersection des filtres maximaux (\cad le \cop de la réunion des \ideps minimaux). Il peut être \care de la manière suivante en \clama
(comparez avec le lemme \ref{lemcRadJ} et sa preuve): c'est l'ensemble
des~\hbox{$a\in\gA$}  \gui{nilréguliers} au sens suivant:
%--------------------begin equation---------------
\begin{equation}\label{Eqnilreg}
\forall y \in\gA\quad\;ay  \;\mathrm{nilpotent}\;\Rightarrow \; y
\;\mathrm{nilpotent}.
\end{equation}
%---------------------end equation--------------
En particulier, dans un anneau réduit, c'est l'ensemble des \elts \ndzsz.
 
}
\end{exercise}
%--- end -exercise-----------------------------------------

%--- Exercise{exoFreeBooleAlgebra}-------------
\begin{exercise}\label{exoFreeBooleAlgebra} 
{(Algèbre de Boole librement engendrée par un ensemble fini)}
\\
{\rm
Soit $E = \{x_1, \ldots, x_n\}$ un ensemble fini.

\emph{1.}
Montrer que l'\agB $\gB$ librement engendrée par $E$ s'identifie à l'\alg
 \[\FF_2[X_1, \ldots, X_n]\sur\fa = \FF_2[x_1, \ldots, x_n]\] avec $\fa = \gen {(X_i^2 - X_i)_{i=1}^n}$.

\emph {2.}
Définir deux $\FF_2$-bases \gui {naturelles} de $\gB$, indexées par
$\Pf(E)$, l'une étant monomiale et l'autre un \sfioz.
Exprimer l'une en fonction de l'autre.
}
\end {exercise}
%--- end -exercise-----------------------------------------

%--- Exercise{exotrdifree}-------------
\begin{exercise}
\label{exotrdifree}
{\rm  Donner une description précise des \trdis librement engendrés par
des ensembles à $0$, $1$, $2$ et $3$ \eltsz. En particulier, préciser le
nombre de leurs \eltsz.
}
\end{exercise}
%--- end -exercise-----------------------------------------

%--- Exercise{exoTreillisDistributifFini}-------------
\begin{exercise}
\label{exoTreillisDistributifFini}
{\rm
On démontre le \thref {thDualiteFinie}.

\emph{1.} On utilise (comme dans le cours)
la structure d'ordre $\supseteq$ sur $E\sta$ (ensemble
des sections initiales de l'ensemble ordonné fini $E$).\\
Si $S_1, S_2 \in E\sta$, que valent $S_1 \vi S_2$, $S_1 \vu S_2$?

\emph{2.} Quelle est la structure d'ordre sur l'ensemble des
\ideps de $\gT$ correspondant à l'ordre qui a été défini
pour $\SpecT$?

\emph{3.}
Démontrer le point \emph{1} du \thoz.  \\
On commencera par vérifier que
pour $S \in E\sta$, $S$ engendre un \idep \ssi $S$ est de la forme $\dar x$
avec $x \in E$; puis que $\Ker \nu_E(x) = \cI_{E\sta}(\dar x)$.

\emph{4.}  Comment construire $E\sta$ à partir de $E$?
Traiter l'exemple suivant:

\snic {
E = \vcenter {\xymatrix @R = 5pt @C = 5pt{
   & &&d \\
c  &     & b\ar@{->}[ur] \\
   & a\ar@{->}[lu]\ar@{->}[ru] \\
}}
\qquad
a < b<d,\  a < c.
}

\'Etudier le cas $E$ totalement ordonné, et le cas $E$
ordonné par la relation d'égalité.

\emph{5.} Démontrer le point \emph{2} du \thoz.

\emph{6.}
Mêmes questions en considérant l'ordre opposé sur $E\sta$
et en adaptant l'ordre sur $\Spec(E\sta)$.
}
\end {exercise}
%--- end -exercise-----------------------------------------

%--- Exercise{exoIVpgcd}-------------
\begin{exercise}
\label{exoIVpgcd} {\rm
Soient $a$, $b$ non nuls dans un anneau intègre.
On suppose que l'\idz~\hbox{$\gen{a,b}$} est \iv et que $a$ et $b$ admettent un ppcm $m$% (autrement dit $\gen{a}\cap\gen{b}=\gen{m}$)
.\\
Montrer que $\gen{a,b}$ est un \idpz.}
\end{exercise}
%--- end -exercise-----------------------------------------

%--- Exercise{exoFactFini}-------------
\begin{exercise}
\label{exoFactFini} {(Un anneau factoriel avec seulement un nombre fini d'\elts \irdsz)}\\ {\rm
Montrer qu'un anneau factoriel avec seulement un nombre fini d'\elts \irds
est un anneau principal. }
\end{exercise}
%--- end -exercise-----------------------------------------

%--- Exercise{exoFactFini}-------------
\begin{exercise}
\label{exoPrincipalIntersecSoucorps} {(Une intersection intéressante)}\\
 {\rm Soit $\gk$ un corps. On considère l'intersection 
$$\gA = \gk(x,y)[z] \cap \gk(z,x+yz).$$
Ce sont deux  sous-anneaux de $\gk(x,y,z)$.
Le premier est principal, le second est un corps.
Montrer que $\gA = \gk[z,x+yz]$, isomorphe à $\gk[z,u]$.
Ainsi l'intersection n'est pas un anneau principal, ni même
un anneau de Bézout.}
\end{exercise}

%--- end -exercise-----------------------------------------

%:2015 nouvel exo
%--- Exercise{exoagGgen1}-------------
\begin{exercise}
\label{exoagGgen1} {(Algèbre de Boole engendrée par un treillis de parties détachables)}\\
 {\rm Démontrer ce qui est affirmé dans l'exemple 1) \paref{exemplesthBoolGen} après le \thref{thBoolGen}.}
\end{exercise}
%--- end -exercise-----------------------------------------

%:2018  nouveaux exos
 
%--- Exercise{exoQuotRIMP}-------------
\begin{exercise}\label{exoQuotRIMP}
 {(Un quotient de relations implicatives)}\\
 {\rm
 Soit $(S,\vdash)$ un ensemble avec une \entrel et $\am,\bn$ des \elts de $ S$.
 On considère la \entrel quotient $(S,\vdash')$ obtenue à partir de $\vda$ en forçant la  relation $\am \vdash' \bn$. Montrer que les relations $\cq \vda d$ pour
$\cq,d\in S$ sont inchangées pour $\vdash'$ \ssi est vérifiée l'implication suivante:
$$ (*)\hspace{2em}
\left. 
\begin{array}{c} 
\cq,\am,b_1\vda d\\[.3em]
\vdots\\[.3em]
\cq,\am,b_n\vda d 
\end{array}
\right \}
\;\Longrightarrow\;\cq,\am\vda d
$$  } 
 
\end{exercise}
%--- end -exercise-----------------------------------------

%--- problem{exoGpRtclQuotient}-------------
\begin{problem}
\label{exoGpRtclQuotient} (Groupes réticulés quotients, sous-groupes solides)\\
{\rm
Dans un ensemble ordonné $E$, si $a\leq b$, on appelle \emph{segment d'extrémités $a$ et $b$}
la partie $\sotq{x\in E}{a\leq x\leq b}$. On le note  $[a,b]_E$ ou $[a,b]$.
 Une partie $F$ de $E$ est dite
\ixc{convexe}{partie --- d'un ensemble ordonné} lorsqu'est satisfaite l'implication $a$, $b\in F\Rightarrow[a,b]\subseteq F$.
\\
Un sous-groupe $H$ d'un \grl est dit \ixc{solide}{sous-groupe ---} si c'est un sous-\grl convexe.\index{sous-groupe solide!d'un \grlz}  
On va voir que  cette notion est l'analogue pour les \grls de celle d'\id pour les anneaux. 

 \emph{1.}    Un sous-groupe $H$ d'un groupe ordonné $G$
 est convexe \ssi la relation d'ordre sur $G$ passe au quotient dans $G/H$, i.e. \prmt $G/H$ est muni d'une structure de groupe ordonné pour laquelle $(G/H)^{+}=G^{+}+H$. On dit aussi \emph{sous-groupe isolé} pour \gui{sous-groupe convexe d'un groupe ordonné}.\index{sous-groupe isolé!d'un groupe ordonné}\index{isole@isolé!sous-groupe ---}

 \emph{2.}  
Le noyau $H$ d'un morphisme de \grls $G\to G'$ est un sous-groupe
solide de $G$.

 \emph{3.} 
 Réciproquement, si $H$ est sous-groupe solide d'un \grl $G$, la loi $\vi$ passe au quotient, elle définit une structure de \grl sur $G/H$, et la
surjection canonique de~$G$ sur~$G/H$ est un morphisme de \grls qui factorise
tout morphisme de source $G$ qui s'annule sur $H$.

 \emph{4.}  
On a défini en \ref{defiCongru} le sous-\grl $\cC(x)$.\\  
Montrer que $\cC(x)\cap\cC(y)=\cC(\abs{x}\vi\abs{y})$, et que le sous-groupe
solide engendré par  $x_1$, \dots, $x_n\in G$ est égal à
$\cC(\abs{x_1}+\cdots+\abs{x_n})$.  En particulier, l'ensemble des sous-groupes
solides \emph{principaux}, i.e. de la forme $\cC(a)$, est \gui{presque} un
\trdi (il manque en \gnl un \elt maximum).

}
\end{problem}
%--- end -problem-----------------------------------------

%--- problem{exoSgpPolaire}-------------
\begin{problem}
\label{exoSgpPolaire} (Sous-groupes polaires, facteurs directs \ortsz)\\
{\rm
 \emph{1.} 
 Si $A$ est une partie quelconque d'un \grl $G$ on note

\snic{A\epr:=\sotq{x\in G\,}{\,\forall a\in A,\;\abs{x}\perp \abs{a}}.}

%\sni
Montrer que $A\epr$ est toujours un sous-groupe solide.
\\
Montrer que l'on  a comme d'habitude dans ce genre de situation:

\snic{ A\subseteq(A\epr)\epr, \; (A\cup B)\epr=A\epr\cap B\epr, \; A\subseteq B\Rightarrow B\epr\subseteq
A\epr\;\hbox{ et }\;{{{A\epr}\epr}\epr}=A\epr.}

%\sni
\emph{2.} 
Un sous-groupe solide $H$ d'un \grl est appelé un \ix{sous-groupe polaire}
lorsque ${{H\epr}\epr}=H$.
On dit encore \emph{une polaire} au lieu de \gui{un sous-groupe polaire}.
\\
   Un sous-groupe $H$ est dit \emph{facteur direct \ortz}
   lorsque  $G=H\oplus H\epr$ (somme directe de sous-groupes 
   dans un groupe abélien), auquel cas $G$ est naturellement 
   isomorphe à $H\boxplus H\epr$. 
   On dit aussi que $G$ est la \ix{somme directe \orte interne} de  $H$ 
   et $H\epr$ et l'on note (par abus) $G=H  \boxplus H\epr$.
\\
Montrer qu'un facteur direct \ort est toujours un sous-groupe polaire.
\\
Montrer que si $G=H\boxplus K$ (avec $H$ et $K$ identifiés à des sous-groupes de $G$) et si~$L$ est un sous-groupe solide, alors
$L=(L\cap H)\boxplus(L\cap K)$.

 \emph{3.} 
   De façon \gnlez, on dit que $G$ est la 
   \ixx{somme directe \orte interne} {d'une famille 
   de sous-\grlsz}  $(H_i)_{i\in I}$, indexée par un ensemble discret $I$, 
   lorsque l'on a $G=\sum_{i\in I}H_i$ et que les $H_i$ sont deux à deux \ortsz. 
   Dans ce cas, chaque~$H_i$ est un sous-groupe polaire de $G$ 
   et l'on a un \iso naturel de \grls $\boxplus_{i\in I}H_i\simeq G$.
   On écrit (par abus)  $G=\boxplus_{i\in I}H_i$.
\\
   On suppose qu'un \grl est somme directe \orte d'une famille 
   de sous-groupes polaires $(H_i)_{i\in I}$, ainsi que d'une autre 
   famille $(K_j)_{j\in J}$. 
   Montrer que ces deux \dcns admettent un raffinement commun.
\\
En déduire que si les composantes d'une \dcn en somme directe \orte
sont des sous-groupes non triviaux
\emph{in\dcpsz}, i.e.,
qui n'admettent pas de facteur direct \ort strict,
alors la \dcn est unique, à bijection
près de l'ensemble des indices.
}
\end{problem}
%--- end -problem-----------------------------------------

%--- problem{exoAutourGaussJoyal}-------------
\begin{problem}\label{exoAutourGaussJoyal} 
 {(Autour de Gauss-Joyal)\iJG}\\
{\rm
Soit $u : \gA \to \gT$ ($\gA$ est un anneau commutatif, $\gT$
un \trdiz) vérifiant:

\snic {
u(ab) = u(a) \vi u(b),\quad u(1) = 1_\gT, \quad
u(0) = 0_\gT, \quad u(a+b) \le u(a) \vu u(b).
}

%\sni
Pour $f = \sum_i a_i X^i \in \gA[X]$, on pose

\centerline {$u(f) = u\big(\rc(f)\big) \eqdefi \Vu_{i} u(a_i).$}

\emph{1.} 
 Montrer que \gui{c'est bien défini}, i.e.,
que $u(f)$ ne dépend que de $\rc(f)$.

 On veut prouver \gui
{de manière directe} (en particulier, sans utiliser le lemme~\ref
{lemGaussJoyal}), la version suivante du lemme de Gauss-Joyal:

\centerline{{\sf LGJ:} $\quad u(fg) = u(f)
\vi u(g)$.}

\emph{2.} 
Vérifier que  si $g = \sum b_jX^j\in \gA[X]$ le résultat équivaut à $u(a_ib_j) \le u(fg)$.

\emph{3.} 
Que dit {\sf LGJ} si $\gT = \so{\Vrai,\Faux}$
et $u(a) = (a \ne 0)$?

\emph{4.} 
En s'inspirant de la preuve classique du résultat de la question
précédente, démontrer~{\sf LGJ}.

\emph{5.} 
Que dit {\sf LGJ} si $\gT = \Zar\gA$ et $u(a) = \rD_\gA(a)$?

}
\end {problem}
%--- end -problem-----------------------------------------

%:--- problem{exoQiClot}----------
\begin{problem}
\label{exoQiClot} (Clôture \qi d'un anneau commutatif)\index{cloture@clôture!quasi intègre}
 \\
{\rm En vous inspirant de la clôture \zede réduite, donner une construction
de la clôture \qi $\Aqi $  d'un anneau commutatif arbitraire $\gA$.
\\
Il faut résoudre le \pb universel suivant: 
\Pnv{\gA}{\lambda}{\psi}{\Aqi }{\varphi}{\gB}{~}{\homos d'anneaux}{morphismes d'anneaux \qis}

\vspace{-3mm}où \emph{les morphismes d'anneaux \qis sont les \homos d'anneaux qui respectent la loi $a\mapsto e_a$} ($e_a$ est l'\idm vérifiant $\gen{1-e_a}=\Ann(a)$). Dans la suite on parlera de \emph{morphisme \qiz}.
\\
Une  clôture \qi d'un anneau $\gA$ existe \gui{a priori}, du simple fait que la théorie des anneaux \qis est purement équationnelle. En effet, pour n'importe quel \sys de \gtrs et de relations (une relation est une \egt entre deux termes construits à partir des \gtrsz, de $0$ et de $1$, en utilisant les lois $+,-,\times,  a\mapsto e_a$), il existe un anneau \qi \gui{le plus \gnl possible} correspondant à cette \pnz:
on prend sur l'ensemble des termes la plus petite relation d'\eqvc qui respecte les axiomes et qui mette dans la même classe d'\eqvc deux termes liés par une relation donnée au départ.
Dans ces conditions l'anneau $\Aqi$ est simplement l'anneau \qi engendré par les \elts de $\gA$ avec pour relations toutes les \egts $a+b=c$, $a\times b=d$, $a=-a'$ vraies dans $\gA$. 
\\
Mais on veut une description précise, comme pour la clôture \zede réduite.
\\
On pourra alors démontrer les résultats suivants.
\index{morphisme!d'anneaux qua@d'anneaux \qisz} 

\emph{1.} \emph{(Morphismes \qisz)}
\begin{itemize}
\item [\emph{a.}] Un morphisme $\varphi: \gA\to\gB$ est \qi \ssi il transforme tout \elt \ndz en
un \elt \ndzz.  
Dans ce cas, il se prolonge de manière unique en un morphisme
$\Frac(\varphi) : \Frac(\gA) \to \Frac(\gB)$.

\item [\emph{b.}] Un morphisme \qi  est injectif \ssi sa restriction à $\BB(\gA)$ est injective.
\item [\emph{c.}] Il existe des \homos injectifs entre anneaux \qis  qui ne sont pas \qisz.
\item [\emph{d.}] Tout \homo entre anneaux  \zedrs est  \qiz.
\item [\emph{e.}]  Si $\gA$ est \qiz, $\BB (\Frac\gA)$ s'identifie à $\BB(\gA)$
et l'injection $\gA\to\Frac(\gA)$ est un morphisme \qiz.
\end{itemize}

\emph{2.} 
 On a des \homos  naturels d'anneaux  $\Ared\to\Aqi\to\Frac(\Aqi)\to\Abul$.
 \\
Il sont tous injectifs et l'\homo naturel $\Frac(\Aqi)\to\Abul$ est un \isoz.

\emph{3.} Si $\gA\subseteq\gC$ avec $\gC$ \qiz, le plus petit sous-anneau \qi de $\gC$ contenant
$\gA$ est égal à $\gA[(e_a)_{a\in\gA}]$.

\emph{4.} 
 Si l'on identifie $\Ared$ à son image dans $\Abul$, 
 on peut identifier $\Aqi$ au sous-anneau de $\Abul$ engendré par $\Ared$
 et par les \idms $e_x$ pour $x\in\Ared$. 
 
\Grandcadre{Dans la suite on suppose \spdg que $\gA$ est réduit.}

\emph{5.} On se reporte au corolaire \ref{corlemZedClotBasicStep} pour la description des étages finis de la construction de $\Abul$. 
Vu le point \emph{4}, ceci nous donne une description des étages finis d'une construction possible de $\Aqi$.
\\
Pour $a_1$, \dots, $a_n\in \gA$, on a une
injection $\gA\to \gA[a_1\bul, \cdots, a_n\bul]=\gC$.\\ 
On note $e_i$ l'\idm $a_ia_i\bul$, $\gB=\gA[e_1, \ldots, e_n]\subseteq \gC$, et $e_I = \prod_{i\in I}(1-e_i) \prod_{j\notin
I}e_j$ pour $I \in \cP_n$.
Montrer les résultats suivants.
\begin{itemize}
\item [\emph {a.}] La famille
$(e_I)_{I\in \cP_n}$ est un \sfio de~$\gB$ et $\gen {1-e_I}_\gB = \gen
{(e_i)_{i\in I},\ (1-e_j)_{j\notin I}}_\gB$.

\item [\emph {b.}]
$\Ann_\gB(a_i) = \gen {1-e_i}_\gB$.

\item [\emph {c.}]
$\gA \cap \gen {e_i, \in I}_\gB = \DA(a_i, i \in I)$.

\item [\emph {d.}]
En notant $\fa'_I = (\DA(a_i, i \in I) : \prod_{j\notin I} a_j)$, on a $\gA
\cap \gen {1-e_I}_\gB = \fa'_I$ et un \iso $\gB \simeq \prod_{I\in\cP_n}
\gA\sur{\fa'_I}$.

\item [\emph {e.}]
L'anneau $\gC$ est un localisé \elr de l'anneau $\gB$:
$\gC=\gB_s$ avec $s \in \gB$ \ndzz.

\end{itemize}

En particulier, soit $a\in \gA$ et $\gA[e_a] \subseteq \gA[a\bul]$
avec $e_a = aa\bul$. \\
Alors, $\Ann_{\gA[e_a]}(a) =
\gen {1-e_a}$, $\gA[e_a] \simeq \gA\sur{\Ann_\gA(a)} \times
\gA\sur{\DA(a)}$, avec $e_a \leftrightarrow (1,0)$. Et $\gA[a\bul]$ est le
localisé $\gA[e_a]_s$ avec $s = 1-e_a + a$ \ndzz.
\\
Dans la suite on note $\gA_{\bra{\an}}$ pour $\gA[a_1a_1\bul,\ldots,a_na_n\bul]$

\emph{6.}
Soient $\varphi : \gA \to \gD$ un morphisme avec $\gD$ réduit, $a \in \gA$
et $b = \varphi(a)$. On suppose que $\Ann_\gD(b) = \gen {1-e_b}_\gD$ avec
$e_b$ \idmz. Montrer que l'on peut prolonger $\varphi$ en un morphisme de
$\gA_{\bra a}  \to \gD$ réalisant $e_a \mapsto e_b$.\\
Cependant, en général, pour $a$, $b \in \gA$, les anneaux $\gA_{\bra{a,b}}$
et $(\gA_{\bra a})_{\bra b}$ ne sont pas isomorphes.

\emph{7.} 
 Donner une description précise de $\ZZ_\mathrm{qi}$.\\ 
 Expliquer
pourquoi l'\homo $\ZZ_\mathrm{qi}\to(\ZZ_\mathrm{qi})_\mathrm{qi}$
n'est pas un \isoz.

\emph{8.} 
 (En \clamaz) Si $\gA$ est \qiz, et $\imath:\gA\to\Frac\gA$ est l'injection canonique, alors $\Spec \imath$ établit une bijection entre $\Spec (\Frac\gA)$
et~$\Min \gA$.

\emph{9.} 
 (En \clamaz) Pour tout anneau $\gA$, il y a une bijection naturelle entre
$\Min(\Aqi )$ et $\Spec \gA$.

\smallskip 
\comm
Bien que $\ZZ$ soit \qiz, $\ZZ_\mathrm{qi}$ n'est pas isomorphe à
$\ZZ$. Ceci se comprend en remarquant que la \prn naturelle $\ZZ\to\ZZ\sur{15\ZZ}$
n'est pas un morphisme  \qiz. 
Cette situation est différente de celle de la clôture \zede réduite: cela 
tient à ce que le quasi inverse $b$ d'un \eltz~$a$, quand il existe, est unique et simplement
défini par deux \eqnsz~\hbox{$ab^2=b$} et~\hbox{$a^2b=a$}, ce qui implique que tout \homo
d'anneaux respecte les quasi inverses.
\eoe

}
\end{problem}
%:--- end-problem-----------------------------------------

% fin des exos

%%%%%%%%%%%%%%%%%%%%%%%%%%%%%%%%%%%%%%%%%%%%%%%%%%%%%%%%%%%%%%%%%
%%%%%%%%%%%%%%%%%%%%%%%%%%%%%%%%%%%%%%%%%%%%%%%%%%%%%%%%%%%%%%%%%
%%%%%%%%%%%%%%%%%%%%%%%%%%%%%%%%%%%%%%%%%%%%%%%%%%%%%%%%%%%%%%%%%
%%%%%%%%%%%%%%%%%%%%%%%%%%%%%%%%%%%%%%%%%%%%%%%%%%%%%%%%%%%%%%%%%
%: solutions d'exos
\sol

%%%%%%%%%%%%%%%%%%%%%%%%%%%%%%%%%%%%%%%%%%%%%%%%%%%%%%%%%%%%%%%%%
\exer{exoRegleCoupure}
En effet, $(a \vi b) \vi x = a \vi x$ car $x \vi a \le b$.\\
Et $(a \vi b) \vu x = (a \vu x) \vi (b \vu x) = a \vu x$
car $a \vu x \le  b \vu x$ (puisque $a \le x \vu b$).\\
Donc $a \vi b = a$, i.e. $a \le b$.

%%%%%%%%%%%%%%%%%%%%%%%%%%%%%%%%%%%%%%%%%%%%%%%%%%%%%%%%%%%%%%%%%@
\exer{exo0GpRtcl} On note $a \sim b$ pour: $a$ et $b$ sont associés.
Montrons la forme suivante (qui est d'ailleurs plus forte si la \dve dans $\gA$ n'est pas explicite): si~$p$ irréductible, $p \divi ab$ 
et $(pa, ab)$ a un pgcd $d$, alors $p \divi a$ ou $p \divi b$.  
On a $p \divi pa$ et~$p \divi ab$, donc
$p \divi d$. Et aussi $a \divi pa$, $a \divi ab$, donc $a \divi d$. 
Soit $a' = d/a \in \gA$.  Comme~$d \divi pa$, on a $a' \divi p$. 
Mais $p$ étant irréductible, on a soit $a' \sim 1$, soit $a' \sim p$. \\
Dans le premier cas,~\hbox{$d \sim a$}, et comme $p \divi d$,  
on a $p \divi a$. 
Dans le second cas, on~a~\hbox{$d \sim ap$}, donc~\hbox{$ap \divi ab$}, i.e $p \divi b$.
\\
Dans $\ZZ[X^2, X^3]$, $X^2$ est irréductible, $X^2 \divi X^3 \cdot X^3$
mais $X^2 \nedivi X^3$, donc~\hbox{$X^2 \cdot X^3$} et $X^3 \cdot X^3$ n'ont pas de
pgcd. A fortiori ils n'ont pas de ppcm.
\\
Enfin: le pgcd de $X^2$ et $X^3$ dans $\ZZ[X^2,X^3]$ est $1$,
s'ils avaient un ppcm ce serait
donc $X^5$, mais $X^5$ ne divise pas $X^6$.

%%%%%%%%%%%%%%%%%%%%%%%%%%%%%%%%%%%%%%%%%%%%%%%%%%%%%%%%%%%%%%%%%@

\exer{exoLemmeGaussAX}
Notons $\G(\fa)$ le pgcd des \gtrs d'un \itf $\fa$. On constate facilement que
c'est bien défini. Ensuite, la \ditz~\hbox{$a(b\vi c)=ab\vi ac$} se généralise
sous la forme~\hbox{$\G(\fa)\G(\fb)=\G(\fa\fb)$} pour deux \itfsz~$\fa$ 
et~$\fb$. Enfin, pour deux \polsz~\hbox{$f$, $g\in\AX$}, \DKM dit que

\snic{\rc(f)^{p+1}\rc(g)=\rc(f)^{p}\rc(fg)$ pour $p\geq\deg g.}

%\sni
Comme~\hbox{$\G(f)=\G\big(\rc(f)\big)$} on obtient $\G(f)^{p+1}\G(g)=\G(f)^{p}\G(fg)$.
Et puisque ce sont des
\elts de l'anneau, on peut simplifier pour obtenir $\G(fg)=\G(f)\G(g)$.

%%%%%%%%%%%%%%%%%%%%%%%%%%%%%%%%%%%%%%%%%%%%%%%%%%%%%%%%%%%%%%%%%@

\exer{exoKroneckerTrick} 
\emph {1.}
Soit $\uX^\alpha = X_1^{\alpha_1} \cdots X_n^{\alpha_n}\in\AuX_{<d}$, alors:

\snic { 
\varphi(\uX^\alpha) = T^a  \quad \hbox {avec} \quad
a = \alpha_1 + \alpha_2 d + \cdots + \alpha_n d^{n-1}.
}

%\sni
On voit ainsi que $a < d^n$. La numération en base $d$ prouve
que $\varphi$ transforme \hbox{la $\gA$-base} de $\AuX_{<d}$ constituée des
$\uX^\alpha$ avec $\alpha_i < d$ en la $\gA$-base $(1, T, \ldots, T^{d^n-1})$
de $\gA[T]_{<d^n}$.

\emph {2.} Rappelons que $\AX\eti=\Ati=\AuX\eti$. Ici on suppose $\AX$ factoriel.
\\
Si $P = QR\in \AuX_{<d}$ alors $Q $ et $R \in \AuX_{<d}$, et
$\varphi(Q)$ et $\varphi(R) \in \gA[T]_{<d^n}$. \\
Comme
$\varphi(QR) = \varphi(Q)\varphi(R)$, et que $f=\varphi(P)$ n'a qu'un nombre fini de facteurs (dans $\AX\etl/\Ati$), il suffit de tester pour chaque facteur  $g(T)$ de $f(T)$ 
si $\varphi^{-1}(g)$ est un facteur de $P$. Ceci est possible car $\gA$ est supposé \dveez.

%%%%%%%%%%%%%%%%%%%%%%%%%%%%%%%%%%%%%%%%%%%%%%%%%%%%%%%%%%%%%%%%%

\exer{exoClosedCover1}{
On se ramène à $r = 2$. L'hypothèse $M$ \tf modulo $\fa_i$,
fournit un sous-module $M_i$ de $M$ de type fini tel que
$M = M_i + \fa_i M$. En reportant la valeur de $M$ dans le membre
droit, on obtient
$$
M = M_i + \fa_i M_i + \fa_i^2 M = M_i + \fa_i^2 M.
$$
En itérant, on obtient pour  $k \ge 1$, $M = M_i + \fa_i^k M$.
En reportant $M = M_2 + \fa_2^k M$ dans $M = M_1 + \fa_1^k M$, on
obtient $M = M_1 + M_2 + (\fa_1\fa_2)^k M$. Mais $\fa_1$, $\fa_2$
sont \tf et $\fa_1\fa_2 \subseteq \DA(0)$, donc il existe un
$k$ tel que $(\fa_1\fa_2)^k = \{0\}$, et par suite~\hbox{$M = M_1+M_2$}
est \tfz.
}

%%%%%%%%%%%%%%%%%%%%%%%%%%%%%%%%%%%%%%%%%%%%%%%%%%%%%%%%%%%%%%%%%

\exer{exoAbul} 
 On peut supposer $\gA$ réduit, sous-anneau de $\Abul$. 
\\ \emph{1.}
Soit $\fp$ un \idep de $\gA$; le morphisme canonique $\gA\to \gK=\Frac(\gA\sur\fp)$ se factorise à travers $\Abul$:

\centerline{\xymatrix @R=10pt @C=25pt{
\gA\ar[d]\ar[dr] \\
\Abul\ar@{-->}[r]_(.3){\pi_\fp}  & \Frac(\gA\sur\fp) \\
}}

Le morphisme $\pi_\fp$ est surjectif car pour $a \in \gA\setminus\fp$, on a
 $1/a = \pi_b(a\bul)$ dans $\gK$. Son noyau $\fq = \Ker\pi_p$
 est un \idema de $\Abul$; on a alors $\gA\sur\fp \subseteq
\gK \simeq \Abul\sur\fq$, donc la flèche naturelle $\gA\sur\fp
\to \Abul\sur\fq$ étant injective, $\fp = \fq\cap\gA$.
On dispose ainsi de deux transformations 

\snic{\Spec\Abul \to \Spec\gA$, $\fq
\mapsto \fq\cap\gA$, et $\Spec\gA \to \Spec\Abul$, $\fp \mapsto \Ker\pi_p,}

%\sni
qui sont réciproques l'une de l'autre.
En effet, si  $\fq\in\Spec\Abul$ et 
$\fp=\fq\cap\gA$, alors~\hbox{$\gK = \Abul\sur\fq$} (car $a\bul = 1/a$
pour $a \in \gA\setminus\fp$) donc $\Ker\pi_\fp = \fq$.
\\
\emph{2.} 
D'après le point \emph{1} l'\homo $\Abul \to \wi\gA$ qui factorise l'\homo naturel $\gA\to  \wi\gA$ est injectif, car son noyau est l'intersection de tous les \ideps de $\Abul$. On identifie $\gA\subseteq\Abul\subseteq\wi\gA$. Le lemme  \ref{lem2SousZedRed}
décrit le plus petit sous-anneau \zedr de $\wi\gA$ contenant $\gA$. On voit qu'il s'agit bien de $\Abul$ (d'après la construction de $\Abul$).
\\
\emph{Autre \demz, laissée \alecz.} On note $\gA_1$ le plus petit sous-anneau \zedr de $\wi\gA$ contenant $\gA$. On démontre alors que cet objet satisfait la \prt \uvle voulue.

%%%%%%%%%%%%%%%%%%%%%%%%%%%%%%%%%%%%%%%%%%%%%%%%%%%%%%%%%%%%%%%%%

\exer{exoMinA} \emph{1.} La première \carn des filtres stricts maximaux parmi les filtres stricts est \imdez: elle revient à dire que toute tentative de faire grandir~$\ff$ en lui rajoutant un \elt $x$ extérieur échoue, parce que le filtre engendré par~$\ff$ et $x$ contient $0$.
\\
Montrons ensuite qu'un filtre strict maximal parmi les filtres stricts est premier. Soient $x$, $y\in\gA$ avec $x+y\in\ff$. On veut montrer que $x\in\ff$ ou $y\in\ff$. Si~\hbox{$x\notin\ff$}, il existe $a\in\ff$ et $n\in\NN$ tel que $a^nx=0$,
donc $a^n(x+y)=a^ny\in\ff$ donc $y\in\ff$. 
\\
Montrons maintenant que le localisé est \zedz, \cad (puisque l'anneau est local) que tout \elt
non inversible est nilpotent. Un \elt non inversible dans le localisé est
un multiple d'un $x/1$ avec $x\notin \ff$. Il suffit de voir que~$x/1$ est nilpotent dans $\ff^{-1}\gA$, or il existe $a\in\ff$ tel que $ax$ est nilpotent
dans~$\gA$, et $a$ est \iv dans le localisé.
\\
Montrons pour terminer que si  $\ff^{-1}\gA$ est local \zed
et non trivial, alors $\ff$ est strict, maximal parmi les filtres stricts.
En effet, un $x\notin\ff$ n'est pas inversible, donc est nilpotent dans le localisé, ce qui signifie qu'il existe $a\in\ff$ tel que $ax$ est nilpotent dans $\gA$.

\emph{2.} 
Notons $S$ la partie définie par l'\eqvrf{Eqnilreg}.
Si $a\in S$ et $a\notin\ff$ avec $\ff$ un filtre maximal, on a $0\in a^\NN\ff$ ce
qui signifie que pour un $x\in\ff$ et $n\in\NN$, $xa^n=0$, donc, puisque $a\in S$,
$x$ est nilpotent: contradiction.\\
Si $a\notin S$, il existe $x$ non nilpotent tel que $xa$ est nilpotent. Donc il
existe un filtre strict contenant $x$. Par le lemme de Zorn il existe un filtre
maximal $\ff$ contenant $x$, et $a$ ne peut pas être dans $\ff$ car sinon $xa$
et donc $0$ serait dans $\ff$.

%%%%%%%%%%%%%%%%%%%%%%%%%%%%%%%%%%%%%%%%%%%%%%%%%%%%%%%%%%%%%%%%%
\exer{exoFreeBooleAlgebra} \emph {1.}
Résulte clairement de la définition d'une \agB comme anneau où tous les
\elts sont \idmsz, à condition de vérifier que l'objet construit est bien
une \agBz, ce qui n'offre pas de difficulté.  On notera que $\gB$ est
isomorphe~à

\snic {
\FF_2[x_1] \otimes_{\FF_2} \cdots  \otimes_{\FF_2} \FF_2[x_n],
}

%\sni
qui est la somme directe de $n$ \agBs librement engendrées par un seul
\gtr dans la catégorie des \agBsz.
En effet, la somme directe de deux \agBs $\gB$, $\gB'$ est
l'\agB $\gB\otimes_{\FF_2}\gB'$.

\emph {2.}
La $\FF_2$-base monomiale de $\gB$ est $(m_I)$ avec $m_I = \prod_{i \in I}
x_i$.  Elle est de cardinal~$2^n$, donc $\gB$ est de cardinal~$2^{2^n}$.  On
définit $e_I$ par $e_I = m_I \prod_{j \notin I} (1 + x_j)$; on vérifie
facilement que $(e_I)$ est un \sfioz, que~\hbox{$m_I e_J = e_J$} si~\hbox{$I \subseteq J$},
et $0$ sinon. \\
On a la même expression $e_I = \sum_{J \,|\, J
\supseteq I} m_J$ et $m_I = \sum_{J \,|\, J \supseteq I} e_J$ (ce qui
confirme que $(e_I)$ est une $\FF_2$-base de $\gB$).
\\
Par rapport à la description donnée dans le cours,
$x_1^{\varepsilon_1} \cdots x_n^{\varepsilon_n}$ correspond à
l'\elt suivant de $\Pf\big(\Pf(E)\big)$: $\so{\sotq{x_i}{\varepsilon_i = 1}}$.

%%%%%%%%%%%%%%%%%%%%%%%%%%%%%%%%%%%%%%%%%%%%%%%%%%%%%%%%%%%%%%%%%
%:2015
\exer{exotrdifree} Le \trdi librement engendré par $\emptyset$  est le treillis $\Deux$.\\
Le \trdi librement engendré par $\so{a}$ est $\so{0,a,1}$.\\
Le \trdi librement engendré par un ensemble $\so{a,b}$ ($a\neq b$) est formé par: $0,\,a\vi b,\, a,\, b, a\vu b,\, 1$.\\
Le \trdi librement engendré par un ensemble $\so{a,b,c}$ d'exactement trois \elts est formé par: 
\[\preskip.2em \postskip.4em 
\begin{array}{ccc} 
 0,\, 1,\, a,\,b,\,c,\,a\vu b,\,a\vu c,\,b\vu c,\,a\vu b\vu c,\,a\vi b,\,a\vi c,\,b\vi c,\,a\vi b\vi c,   \\[.3em] 
a\vi(b\vu c),\,b\vi(a\vu c),\,c\vi(a\vu b),\,(a\vu b)\vi(a\vu c),\,(a\vu b)\vi(b\vu c),\,\\[.3em] 
(a\vu c)\vi(b\vu c),\,(a\vu b)\vi(a\vu c)\vi(b\vu c). \end{array}
\]

%%%%%%%%%%%%%%%%%%%%%%%%%%%%%%%%%%%%%%%%%%%%%%%%%%%%%%%%%%%%%%%%%
\exer{exoTreillisDistributifFini}
\emph{1.} Par \dfn d'une section initiale l'intersection et la réunion de deux sections initiales en est une autre.   \\
Donc dans $E\sta$: $S_1 \vi S_2
= S_1 \cup S_2$, $S_1 \vu S_2 = S_1 \cap S_2$, $\emptyset =
1_{E\sta}$ et $E = 0_{E\sta}$.
\\
\emph{2.} Il revient au même de se donner $\alpha\in\SpecT$ ou l'\idep $\Ker\alpha$.
Ceci conduit à ordonner l'ensemble des idéaux premiers
de $\gT$ par la relation
$\supseteq$. En effet, si $\alpha, \beta : \gT \to \{0,1\}$ sont deux
morphismes, on a l'\eqvc

\snic {\alpha \le \beta \iff \Ker\alpha \supseteq
\Ker\beta.}

\emph{3.} On a:

\snuc {
\Ker\nu_E(x) = \sotq {S \in E\sta} {x \in S} =
\sotq {S \in E\sta} {\dar x \subseteq S} =
\sotq {S \in E\sta} {S \le \dar x}
,}

%\sni
i.e. $\Ker\nu_E(x) = \cI_{E\sta}(\dar x)$. On a bien les \eqvcsz:

\snuc {
x \le y \iff \dar x \subseteq \dar y \iff \dar y \le \dar x
\iff \cI(\dar y) \subseteq \cI(\dar x) \iff
\cI(\dar x) \le \cI(\dar y)
}

\mni
Par ailleurs, dans $E\sta$: $S_1 \vi S_2 \le \dar x \Rightarrow (S_1 \le \dar
x) \hbox { ou } (S_2 \le \dar x)$ (car la première in\egt signifie $\dar x
\subseteq S_1 \cup S_2$, i.e. $x \in S_1 \cup S_2$). Et comme
$\dar x \ne 1_{E^*} = \emptyset$,  $\dar x$ engendre un \idepz.  
Réciproquement, soit $\fp$ un \idep de $E\sta$.
\'Etant fini, il est principal: $\fp = \cI_{E\sta}(S)$ avec $S \ne 1_{E\sta}$,
i.e. $S$ non vide. Il faut montrer que $S$ est de la forme $\dar x$.
Si $S =\{x_1, \ldots, x_n\}$, on a $S = (\dar x_1)\cup\cdots\cup(\dar x_n)$, i.e.~\hbox{$(\dar x_1) \vi \cdots \vi (\dar x_n) = S$}.
Comme $S$ engendre un \idepz, il existe un~$i$ tel que $\dar x_i \le S$, i.e. $S \subseteq \dar x_i$,
puis $S = \dar x_i$.

 \emph{4.}
On détermine $E\sta$ en remarquant que toute section initiale est une
réunion de parties $\dar x$.  Le dessin du treillis $E\sta$ est le suivant:

\snic {
\xymatrix @R = 5pt @C = 4pt{
       & \emptyset\\
       & \dar a = \{a\}\ar@{-}[u] \\
\dar c =\{a,c\}\ar@{-}[ur] &&\dar b =\{a,b\}\ar@{-}[ul] \\
       &\{a,b,c\}\ar@{-}[ul]\ar@{-}[ur] &   & \dar d = \{a,b,d\}\ar@{-}[ul] \\
       &                              &\{a,b,c,d\}\ar@{-}[ul]\ar@{-}[ur] \\
}
}

%\sni
Si $E$ est totalement ordonné, alors $E\sta = \{\, \dar x \;\vert\; x \in E\,\}
\cup \{\emptyset\}$ est aussi totalement ordonné et $\#E\sta = 1 + \#E$.
Si $\gT$ est un ensemble fini totalement ordonné, alors~\hbox{$\Spec\gT =
\big\{\cI_{\gT}(a) \;\vert\;  a \in \gT \setminus \{1_\gT\}\big\}$}, et
$\#\Spec\gT = \#\gT - 1$. Si $E$ est ordonné par la relation d'égalité,
$E\sta = \cP(E)$ ordonné par $\supseteq$. Quant à $\Spec(E\sta)$, c'est
l'ensemble $\cI_{\cP(E)}(\{x\})$ avec $x \in E$ (qui est bien isomorphe
à $E$).

 \emph{5.}
\Llec vérifiera qu'en posant, pour
$a \in \gT$, $\wh {a} = \sotq {\fp \in \Spec\gT}{a \in \fp}$,
on obtient une section initiale, que toute section initiale
de $\Spec\gT$ est de cette forme, et enfin que $a \le b \iff
\wh {a} \le \wh {b}$.

\emph{6.}
On prend maintenant comme structure d'ordre $\subseteq$ sur $E\sta$ et sur
$\Spec \gT$.  \\
Alors $S_1 \vi S_2 =
S_1 \cap S_2$, $S_1 \vu S_2 = S_1 \cup S_2$, $\emptyset = 0_{E\sta}$, $E =
1_{E\sta}$.  Pour~\hbox{$x \in E$}, on pose~\hbox{$\wi x = E \setminus \uar x = \sotq {y \in E} {y \not\geq x}$}: cet \elt de $E\sta$ vérifie, pour~\hbox{$S \in E\sta$}, l'\eqvc $x \notin S \iff S \subseteq \wi x$.  
On a $\wi x \ne 1_{E\sta} = E$,
et $\wi x$ engendre un idéal premier du treillis $E\sta$: $S_1
\vi S_2 \le \wi x \Rightarrow S_1 \le \wi x \hbox { ou } S_2 \le \wi x$ 
(en effet, l'hypothèse est~\hbox{$\uar x \subseteq (E \setminus S_1) \cup (E \setminus S_2)$}, 
donc par exemple $x \notin S_1$, i.e.~\hbox{$S_1 \subseteq \wi {x}$)}.  On a
l'\eqvcz: $x \le y \iff \wi x \subseteq \wi y$.  On démontre que tout
idéal premier de $E\sta$ est de la forme $\wi x$, donc l'ensemble ordonné
$E$ est isomorphe, via~\hbox{$x \mapsto\cI_{E\sta} (\wi x)$},
à l'ensemble des idéaux premiers de $E\sta$, ordonné par inclusion.

\snic{
\xymatrix @R = 5pt @C = 4pt{
       & \{a,b,c,d\} \\
\wi c =\{a,b,d\}\ar@{-}[ur] &&\wi d =\{a,b,c\}\ar@{-}[ul] \\
       &\{a,b\}\ar@{-}[ul]\ar@{-}[ur] &   & \wi b = \{a,c\}\ar@{-}[ul] \\
       &                              &\{a\}\ar@{-}[ul]\ar@{-}[ur] \\
       &                              &\wi a = \emptyset \ar@{-}[u] \\
}
}

%%%%%%%%%%%%%%%%%%%%%%%%%%%%%%%%%%%%%%%%%%%%%%%%%%%%%%%%%%%%%%%%%

\exer{exoIVpgcd}
Puisque  $ \gen{a,b} $ est inversible  on a
     $s,$ $t$, $u$, $v$  avec  $sa = ub$,  $tb = va$ et $s+t = 1$.
\\
Puisque    $m$ est  ppcm de  $a$ et $b$ on peut écrire

\snic{m= ab' = ba'  \;\hbox{ et }\;   ab/m = g = b/b' = a/a'.}            
           
Ainsi $sa = mx = ab'x$ et $tb = m=ba'y$,  qui donnent $ s=b'x$ et $t=a'y$.
\\
Donc $b'x+a'y = 1$,
  $bx + ay = gb'x + ga'y = g$
et par suite $ \gen{a,b}=\gen{g} $.             

%%%%%%%%%%%%%%%%%%%%%%%%%%%%%%%%%%%%%%%%%%%%%%%%%%%%%%%%%%%%%%%%%

\exer{exoFactFini} \emph{(Un anneau factoriel avec seulement un nombre fini d'\elts \irdsz)} 
On note $(p_i)_{i\in I}$ la famille finie des \elts \irds distincts (à association près).
\\
On doit montrer que $\gA$ est un anneau de Bézout. Il suffit pour cela de montrer que si $a$ et $b\in\Atl$ ont pour pgcd $1$, alors $\gen{a,b}=\gen{1}=\gA$.
On écrit $$\ndsp a=\prod_{i\in A}p_i^{\alpha_i}, \;b=\prod_{j\in B}p_i^{\beta_j},
\hbox{avec les}  \;\alpha_i \,\hbox{et} \,\beta_j>0 \,\hbox{ et }\,A\cap B =\emptyset.$$ 
Soit $C=I\setminus(A\cup B)$ et $c=\prod_{k\in C}p_k$. On montre que $a+bc\in\Ati$. 
\\
En effet,
pour $i\in A$, $p_i$ divise $a$, donc il ne peut pas diviser~\hbox{$a+bc$}, sinon il diviserait~\hbox{$bc=(a+bc)-a$}. 
De même, pour $j\in B\cup C$, $p_j$ ne peut pas diviser~\hbox{$a+bc$}, sinon il diviserait $a=(a+bc)-bc$.
Ainsi $a+bc$ n'est divisible par aucun \elt \irdz.

%%%%%%%%%%%%%%%%%%%%%%%%%%%%%%%%%%%%%%%%%%%%%%%%%%%%%%%%%%%%%%%%%

\exer{exoPrincipalIntersecSoucorps} \emph{(Une intersection intéressante)}\\
On considère l'\homo d'\evn $$\varphi:\gk[z,u]\to \gk[z,x+yz] , \;z\mapsto z, \,u\mapsto x+yz.$$ 
Il est surjectif par construction. Il est injectif parce que, pour $f=f(z,u)$, en évaluant $\varphi(f)$  dans $\gk[x,y,z]$ on obtient
$\varphi(f)(x,0,z)=f(z,x)$. C'est donc bien un \isoz.
\\
Dans la suite on peut donc poser~\hbox{$u=x+yz$}, avec $\gk[z,x+yz]=\gk[z,u]$ où $z$ et $u$
jouent le rôle d'\idtrs distinctes. 
\\
Par ailleurs on remarque que $\gk[z,u][y]=\gk[x,y,z]$. Comme $\gk[z,u]$ est un anneau à pgcd, ceci implique que
deux \elts de $\gk[z,u]$ sont de pgcd  $1$ dans $\gk[z,u]$ \ssi ils sont de pgcd $1$ dans  $\gk[x,y,z]$.

Soit maintenant un \elt arbitraire $h\in\gA$ que l'on écrit sous forme d'une fraction \ird $f(z,u)/g(z,u)$ dans $\gk(z,u)$, et sous forme d'une fraction~\hbox{$a/b$} ($a\in\gk[x,y,z]$, $b\in\gk[x,y]$) en tant qu'\elt de $\gk(x,y)[z]$. Cette dernière fraction
peut elle-même être écrite sous forme \irdz, \cad que le pgcd de~$a$ et
$b$ dans $\gk[x,y,z]$ est égal à $1$. Par unicité de l'écriture d'une fraction sous forme réduite, on a donc une constante $\gamma\in\gk\etl$ telle
que $f(z,u)=\gamma a(x,y,z)$ et~\hbox{$g(z,u)=\gamma b(x,y)$}. 
\\
Il nous reste à montrer que le dénominateur $g(z,x+yz)$ est une constante. En faisant
$z=0$ dans l'\egt $g(z,x+yz)=\gamma b(x,y)$ on obtient $$g(0,x)=\gamma  b(x,y)=c(x).$$
Enfin, en faisant $(z,y)=(1,-x)$ dans l'\egt $g(z,x+yz)=c(x)$, on 
obtient~\hbox{$c(x)=g(1,0)$}.

%%%%%%%%%%%%%%%%%%%%%%%%%%%%%%%%%%%%%%%%%%%%%%%%%%%%%%%%%%%%%%%%%

%%%%%%%%%%%%%%%%%%%%%%%%%%%%%%%%%%%%%%%%%%%%%%%%%%%%%%%%%%%%%%%%%

\exer{exoagGgen1} \emph{(Algèbre de Boole engendrée par un treillis de parties détachables)}
On considère un sous-treillis $\gT$ de l'\agB  des parties détachables d'un ensemble $E$, \alg que nous notons $B_E$. 
Il est clair que les combinaisons booléennes d'\elts de $\gT$ sont des \elts de $B_E$ et qu'elles forment une sous-\agB $\gB$ de $B_E$. Cela implique que l'on a un morphisme d'\agB $\alpha:\Bo(\gT)\to \gB$ qui factorise le morphisme d'inclusion $\gT\to\gB$ (morphisme de \trdisz).\\
Vue la construction de $\gB$, le morphisme $\alpha$ est surjectif et il reste à démontrer qu'il est injectif. Pour cela, 
on considère des \elts $A_1,\dots,A_n,C_1,\dots,C_m$ de $\gT$ et deux combinaisons booléennes formelles $A$ et $C$, respectivement des~$A_i$ et des~$C_j$. On évalue~$A$ et $C$ dans $\gB$ et dans $\Bo(\gT)$. On suppose \hbox{que $\varphi(A)=\varphi(C)$} dans $\gB$ et l'on doit montrer que $\psi(A)=\psi(C)$ dans $\Bo(\gT)$. 
%\\
En fait, en utilisant la différence symétrique, on doit montrer que si $\varphi(A\oplus C)=\emptyset$, alors $\psi(A\oplus C)=0_{\Bo(\gT)}$.
\\
On note que toute combinaison booléenne formelle  se réécrit (en suivant les axiomes des \agBsz) en une forme normale disjonctive. Une telle forme normale (pour $A\oplus C$) est évaluée nulle dans une \agBz~$B$ \ssi chacun des disjoncts est évalué nul. Considérons par exemple un disjonct 
$$
A_1\vi \ov{A_2}\vi C_1\vi C_2\vi \ov{C_3}\vi \ov{C_4}=(A_1\vi C_1\vi C_2)\vi\ov{A_2\vu C_3\vu C_4}.
$$ 
Il s'évalue en $0$ via une \evn $\beta$ dans $B$ \ssi 
$$
\beta(A_1)\vi_B \beta(C_1)\vi_B \beta(C_2)\leq_B \beta(A_2)\vu_B \beta(C_3)\vu_B \beta(C_4).
$$
Or si cela est satisfait pour l'\evn $\varphi$ dans $\gB$, par \dfn de $\gT$, cela est satisfait dans $\gT$, et donc aussi dans $\Bo(\gT)$.
\\
Ceci termine la \demz.

%%%%%%%%%%%%%%%%%%%%%%%%%%%%%%%%%%%%%%%%%%%%%%%%%%%%%%%%%%%%%%%%%%%%
\exer{exoQuotRIMP} \emph{(Un quotient de relations implicatives)} \\ Définissons   $A=\so{a_1,\dots,a_m}$,  $X=\so{c_1,\dots,c_q,a_1,\dots,a_m}$, $Y=\so d$ et $B=\so{b_1,\dots,b_n}$ (parties de~$S$). 
La \entrel $\vdash'$ de l'énoncé s'obtient à partir de $\vda$ en forçant $A\vdash'B$. La proposition \ref{lemquotrdi2} nous dit que $X\vdash'd$
\ssi pour chaque $a\in A$ et chaque $b\in B$, on a: $X,b\vda d$ et $X\vda d,a$. Ici cela veut dire que pour chaque chaque $b\in B$, on a: $X,b\vda d$, ce qui est l'hypothèse dans
l'implication $(*)$ de l'énoncé.
Cela donne donc l'équivalence souhaitée.

%%%%%%%%%%%%%%%%%%%%%%%%%%%%%%%%%%%%%%%%%%%%%%%%%%%%%%%%%%%%%%%%%%%%

%%%%%%%%%%%%%%%%%%%%%%%%%%%%%%%%%%%%%%%%%%%%%%%%%%%%%%%%%%%%%%%%%
\prob{exoAutourGaussJoyal}
Le premier point est laissé \alecz. Notons $fg = \sum_k c_k X^k$.

 \emph{2.}
On a facilement $u(fg) \le u(f) \vi u(g)$. \\
En effet, $c_k = \sum_{i+j = k}
a_ib_j$, donc $u(c_k) \le \Vu_{i+j = k} u(a_ib_j) \le \Vu_{i} u(a_i)
 = u(f)$ (on~a utilisé $u(ab) \le u(a)$).
\\
Si l'on dispose de Gauss-Joyal, alors $u(a_ib_j) \le u(a_i)
\vi u(b_j) \le u(f) \vi u(g) = u(fg)$.  Réciproquement, si l'on sait
prouver $u(a_ib_j) \le u(fg)$ pour tous $i,j$, alors

\snic {
\Vu_{i,j} u(a_ib_j) \le u(fg),
\; \hbox {i.e. par \ditz} \;
\bigr(\Vu_{i} u(a_i)\bigl) \vi \bigr(\Vu_{j} u(b_j)\bigr)
\le u(fg),
}

%\sni
i.e. $u(f) \vi u(g) \le u(fg)$.

 \emph{3.}
Si $\gA$ est intègre, il en est de même de $\gA[X]$.

 \emph{4.}
Montrons par \recu descendante sur $i_0+j_0$ que $u(a_{i_0}b_{j_0}) \le u(fg)$.
C'est vrai si $i_0$ ou $j_0$ est grand car alors $a_{i_0}b_{j_0} = 0$.  
On écrit la \dfn du produit de deux \polsz:

\snic {
a_{i_0}b_{j_0} = c_{i_0+j_0}
- \sum\limits_{i+j = i_0+j_0 \atop i > i_0} a_ib_j
- \sum\limits_{i+j = i_0+j_0 \atop j > j_0} a_ib_j.
}

%\sni
On applique $u$ en utilisant d'une part $u(\alpha + \beta + \cdots) \le
u(\alpha) \vee u(\beta) \vee \dots$ et d'autre part $u(\alpha\beta) \le
u(\alpha)$ pour obtenir

\snic {(\star)\;:\;
u(a_{i_0}b_{j_0}) \le u(c_{i_0+j_0}) \vee
\Vu_{i > i_0} u(a_i) \vee \Vu_{j > j_0} u(b_j).}

%\sni
On dispose ainsi d'une in\egt $x \le y$ que l'on écrit $x
\le x \wedge y$. Autrement dit, dans $(\star)$, on réinjecte
$u(a_{i_0}b_{j_0})$ dans le membre droit, ce qui donne, en utilisant la
\ditz:

\snic {
u(a_{i_0}b_{j_0}) \le u(c_{i_0+j_0}) \vee
\Vu_{i > i_0} \big(u(a_i) \wedge u(a_{i_0}b_{j_0})\big) \vee
\Vu_{j > j_0} \big(u(b_j) \wedge u(a_{i_0}b_{j_0})\big).
}

%\sni
En utilisant $u(a_i) \wedge u(a_{i_0}b_{j_0}) \le u(a_i) \wedge u(b_{j_0})$ et
$u(b_j) \wedge u(a_{i_0}b_{j_0}) \le u(b_j) \wedge u(a_{i_0})$, et
(par définition) $u(c_{i_0+j_0}) \le u(fg)$, on majore
$u(a_{i_0}b_{j_0})$ par:

\snic {
u(fg) \vee \Vu_{i > i_0} u(a_{i}b_{j_0}) \vee
\Vu_{j > j_0} u(a_{i_0}b_{j}).
}

%\sni
Par \recu sur $i_0, j_0$, $u(a_{i}b_{j_0}) \le u(fg)$,
$u(a_{i_0}b_{j}) \le u(fg)$.\\ D'où $u(a_{i_0}b_{j_0}) \le u(fg)$.

 \emph{5.}
Dans ce cas: $a_ib_j \in \rD_\gA(c_k, k = 0, \ldots)$, ce qui est le lemme de Gauss-Joyal usuel.

%%%%%%%%%%%%%%%%%%%%%%%%%%%%%%%%%%%%%%%%%%%%%%%%%%%%%%%%%%%%%%%%%

\prob{exoQiClot} \emph{(Clôture \qi d'un anneau commutatif)}

Remarque préalable: si dans un anneau $\gA$, $\Ann(a) = \gen {e'_a}$ avec~$e'_a$
\idmz, alors $e'_a$ est l'unique $e'$ tel que:

\snic {
e' a = 0,\quad e'+a\hbox { est \ndzz} \quad\hbox{et}\quad  e' \hbox { est \idmz.}
}

En effet, $e' = e'e'_a$ (car $e'a=0$) et $(e'+a)e' = (e'+a)e'_a$ ($=e'$) d'où $e'=e'_a$.

\emph{1.} Soient  $\gA$, $\gB$  \qis et un morphisme \qi $\varphi:\gA\to\gB$.
\\
\emph{1a.} 
Si $a\in\gA$ est \ndzz, $e_a=1$ donc $e_{\varphi(a)}=1$ donc $\varphi(a)$ est
\ndzz. Inversement, soit $\psi:\gA\to\gB$ un \homo d'anneaux qui transforme
tout \elt \ndz en un \elt \ndzz.
Soit $a\in\gA$, $b=\psi(a)$ et $f=\psi(1-e_a)$.  \\
Alors $fb = \psi\big((1-e_a)a\big) = 0$,
$f+b = \psi(1-e_a + a)$ est \ndz et $f^2 = f$, et \hbox{donc $f = 1-e_b$}.
\\
\emph{1b.} Supposons $\varphi(x)=0$, alors $e_{\varphi(x)}=0$, i.e. $\varphi(e_x)=0$. Donc si $\varphi|_{\BB(\gA)}$ est injectif, on obtient $e_x=0$, i.e. $x=0$.
\\
\emph{1c.} On considère l'unique \homo  $\rho:\ZZ\to\prod_{n>0}\aqo\ZZ{2^n}$. Alors $\rho$ est injectif mais $\rho(2)$ n'est pas \ndzz.
\\
\emph{1d.} L'\homo conserve les quasi inverses, donc aussi les \idms associés
car $e_a=aa\bul$ si $a\bul$ est le quasi inverse de $a$.  
\\
\emph{1e.} 
Résulte \imdt du fait \ref{factQoQiZed}.

 \emph{2.} Puisque $\Aqi$ est réduit, il y a un unique \homo d'anneau $\Ared\to\Aqi$ qui factorise les deux \homos canoniques $\gA\to\Ared$ et $\gA\to\Aqi$.
Puisque $\Abul$  est \qiz, il y a un unique morphisme \qi $\Aqi\to\Abul$
 qui factorise les deux \homos canoniques $\gA\to\Aqi$ et $\gA\to\Abul$.
 Puisque le  morphisme $\Aqi\to\Abul$ transforme un \elt \ndz en un \elt \ndzz,
 et qu'un \elt \ndz dans un anneau \zed (réduit ou pas) est \ivz, il existe un unique \homo
 $\Frac(\Aqi)\to\Abul$ qui factorise les deux \homos canoniques $\Aqi\to\Frac(\Aqi)$ et $\Aqi\to\Abul$.
\\
 De la même manière, pour tout \homo $\gA\to\gB$ avec $\gB$ \zedrz,
 on obtient d'abord un unique morphisme \hbox{\qi $\Aqi\to\gB$} (qui factorise ce qu'il faut), puis un
 unique morphisme $\Frac(\Aqi)\to\gB$ qui factorise les deux \homos $\gA\to\Frac(\Aqi)$ et $\gA\to\gB$. \\
 Autrement dit, puisque $\Frac(\Aqi)$ est \zedrz, il résout le \pb \uvl de la clôture \zede réduite pour $\gA$. En conséquence l'\homo $\Frac(\Aqi)\to\Abul$
 que l'on a construit est un \isoz.  

 \emph{3.} Ce point est recopié  du lemme
\ref{lem2SousZedRed} qui concerne les anneaux \zeds réduits: \llec pourra aussi à très peu près recopier la \demz.

 \emph{4.} On note tout d'abord que l'\homo naturel $\Ared\to\Aqi$ est injectif parce que l'\homo $\Ared\to\Abul$ est injectif et qu'il y a \fcnz.
On peut donc identifier $\Ared$ à un sous-anneau de $\Aqi$, 
qui s'identifie lui-même à un sous-anneau de $\Frac(\Aqi)$
 que l'on identifie à $\Abul$.
Dans ce cadre $\Aqi$ contient \ncrt $\Ared$ ainsi que les \elts
$e_x=xx\bul$ pour les $x\in\Ared$ puisque le morphisme $\Aqi\to\Abul$ 
est \qi et injectif. 
\\ Notons $\gB$ le sous anneau de $\Abul$
engendré par $\Ared$ et les \idms $(e_x)_{x\in\Ared}$.
Il reste à voir que l'inclusion $\gB\subseteq\Aqi$ est une \egtz.
\\
Il est clair que $\Frac(\gB) = \Frac(\Aqi)$. D'autre part, comme $\gB$ est
\qiz, l'injection $\Ared \to \gB$ fournit un (unique) morphisme \qi $\varphi : \Aqi \to\gB$
tel que $\varphi(a) = a$ pour tout $a \in \Ared$.
Puisque le morphisme est \qiz, on en déduit que
que $\varphi(e_a) = e_a$ pour tout $a \in \Ared$, puis que $\varphi(b)=b$ pour tout~\hbox{$b\in\gB$}. Soit $x \in \Aqi$; on veut
montrer que $x \in \gB$; comme $x \in \Frac(\gB)$, il existe $b \in \gB$ \ndz
tel que $bx \in \gB$ donc $\varphi(bx) = bx$ \cad $b\varphi(x) = bx$; comme
$b$ est \ndz dans $\gB$, il l'est dans $\Frac(\gB)$, a fortiori dans $\Aqi$,
donc $x = \varphi(x) \in \gB$.

\emph{5a} et \emph{5b.} Facile.

\emph{5c.}
Puisque $a_j = a_je_j$, on a, pour $j \in I$, $a_j \in \gen{e_i, i \in I}_\gB
= \gen {e}_\gB$ 
avec $e$ l'\idm $1 - \prod_{i \in I}(1-e_i)$.  Mais dans un
anneau réduit, tout \idm engendre un \id radical:

\snic {
b^m \in \gen{e} \Rightarrow b^m(1-e) = 0 \Rightarrow b(1-e) = 0 
\Rightarrow b = be \in \gen {e}.
}

%\sni
Donc $\DA(a_i, i \in I) \subseteq \gen{e_i, i \in I}_\gB$.  \\
Montrons
maintenant que $\gA \cap \gen {e_i, \in I}_\gC \subseteq \DA(a_i, i
\in I)$. Soit $x \in \gA \cap \gen {e_i, \in I}_\gC$; en revenant à la
\dfn initiale de $\gC$, on a $x \in \gen{a_iT_i, i \in I}_{\gA[\uT]} + \fc$.
Travaillons sur l'anneau réduit $\gA' = \gA\sur{\DA(a_i, i \in I)}$; on a
alors 

\snic{\ov x \in \rD_{\gA'[\uT]}(a_kT_k^2 - T_k, a_k^2T_k - a_k, k \in
\lrbn).}

%\sni
Puisque $\gA' \to \gA'[\ov a_1\bul, \ldots, \ov a_n\bul]$ est
injectif, on a $\ov x = 0$ i.e. $x \in \DA(a_i, i \in I)$.

\emph{5d.} 
Notons $\pi$ le produit $\prod_{j\notin I} a_j$. Soit $x \in \gA\cap
\gen{1-e_I}_\gB$; puisque $\pi (1-e_j) = 0$ pour $j \notin I$, on a $\pi x \in
\gen{e_i, i \in I}_\gB$, donc, d'après \emph{5c)}, $\pi x \in
\DA(a_i, i \in I)$, \hbox{i.e.  $x \in \fa'_I = (\DA(a_i, i \in I) : \pi)$}.

Réciproquement, soit $x \in \fa'_I$; on écrit $x = \pi' x + (1-\pi')x$
avec $\pi' = \prod_{j \notin I} e_j$. \linebreak 
 On a $1-\pi' \in \gen{1-e_j, j \notin
I}$.  Quant à $\pi' x$, on remarque que dans $\gC$, $\gen{e_j}_\gC = \gen
{a_j}_\gC$, donc $\pi'x \in \gen{\pi x}_\gC \subseteq \rD_\gC(a_i, i \in I)
\subseteq \gen {e_i, i \in I}_\gC$. \\
Bilan: $x \in \gen {(e_i)_{i\in I},
(1-e_j)_{j \notin I}}_\gC = \gen {1-e_I}_\gC$. \\
 Mais $\gA \cap \gen{1-e_I}_\gC = \gA \cap \gen {1-e_I}_\gB$, donc $x \in \gen {1-e_I}_\gB$.

Enfin, $\gB$ est isomorphe au produit des $\gB\sur{\gen{1-e_I}_\gB}$
et $\gB\sur{\gen{1-e_I}_\gB} \simeq \gA\sur{\fa'_I}$.

\emph{5e.} 
Prendre $s = \sum_I e_I \prod_{j\notin I}a_j = \sum_I \prod_{i\in I}
(1-e_i) \prod_{j\notin I}a_j$: $s$ est l'unique \elt de $\gB$
qui vaut $\prod_{j\notin I}a_j$ sur la composante $e_I = 1$.

\emph{6.} Dans l'\iso $\gA[e_a] \simeq \gA\sur{\Ann_\gA(a)} \times
\gA\sur{\DA(a)}$, on a $e_a = (1,0)$ et donc $(\ov x, \ov y) = xe_a +
y(1-e_a)$.  On considère alors l'application 

\snic{\gA\times\gA \to \gD$, $(x,y)
\mapsto \varphi(x)e_b + \varphi(y) (1-e_b).}

%\sni
C'est un morphisme d'anneaux et
puisque $\gD$ est réduit, elle passe au quotient modulo
$\Ann_\gA(a) \times \DA(a)$.
\\ 
Comparons maintenant $\gA_{\bra{a,b}}$
et $(\gA_{\bra a})_{\bra b}$.
On trouve:

\snic {\arraycolsep2pt
\begin {array} {rcl}
\gA_{\bra{a,b}} &\simeq& 
\gA\sur{(0 : ab)} \times \gA\sur{(\rD(b) : a)} \times
\gA\sur{(\rD(a) : b)} \times \gA\sur{\rD(a,b)},
\\[1mm]
(\gA_{\bra a})_{\bra b} &\simeq&
\gA\sur{(0 : ab)} \times  \gA\sur{\rD\big((0 : a) + \gen {b}\big)} \times
\gA\sur{(\rD(a) : b)} \times \gA\sur{\rD(a,b)}.
\end {array}
}
Enfin, on note que  $\rD\big((0 : a) + \gen {b}\big)$ est contenu dans $(\rD(b) : a)$ mais
que l'inclusion peut être stricte. Prenons par exemple $\gA=\ZZ$,
$a = 2p$, $b = 2q$ où~$p$ et~$q$ sont deux premiers impairs disctints.
On utilise $(x : y) = x/\pgcd(x,y)$ pour  $x$, $y \in \ZZ$.\\
Alors $
\ZZ_{\bra{a,b}} \simeq \ZZ \times \ZZ\sur{q\ZZ} \times
\ZZ\sur{p\ZZ} \times \ZZ\sur{2\ZZ}$,
mais
$
(\ZZ_{\bra a})_{\bra b} \simeq \ZZ \times\ZZ\sur{2q\ZZ} \times
\ZZ\sur{p\ZZ} \times \ZZ\sur{2\ZZ}$. Dans le premier anneau, 
$\Ann(a)$ est engendré par  $(0,0,1,1)$. Dans le
second (le premier anneau en est un quotient), $\Ann(a)$
engendré par l'\idm $(0,q,1,1)$. 

\emph{8.} 
On rappelle (exercice \ref{exoMinA}) qu'un \idep $\fp$ d'un anneau $\gA$ est minimal \ssi pour tout
$x\in\fp$, il existe $s \in \gA\setminus\fp$ tel que $sx^n = 0$ pour un
certain $n$ (si $\gA$ est réduit, on peut prendre $n = 1$).
\\
D'abord, un \idemi de $\gA$ reste un \idep strict dans~$\Frac(\gA)$ (ce fait n'utilise pas $\gA$ \qiz), i.e.  $\fp \cap \Reg(\gA) =
\emptyset$: si $x \in \fp$, il existe  $s \notin
\fp$ et~\hbox{$n \in \NN$} tels que $sx^n = 0$, ce qui prouve que $x \notin \Reg(\gA)$.
\\
Réciproquement, pour $\fq$ \idep de $\Frac(\gA)$, prouvons que $\fp =
\fq\cap\gA$ est un \idemi de $\gA$. Soit $x \in \fp$; alors $x + 1-e_x$ est
\ndz dans $\gA$, donc \iv dans $\Frac(\gA)$, donc $1-e_x \notin \fp$. Alors $x(1-e_x)=xe_x(1-e_x)=0$: on a trouvé $s=1-e_x
\notin \fp$ tel que $sx = 0$.

\emph{9.}
D'après l'exercice \ref {exoAbul}, l'injection $\gA\to\Abul$ induit une
bijection $\Spec\Abul\to\Spec\gA$; mais $\Abul = \Frac(\Aqi)$ et $\Aqi$ est
\qiz.\\ 
Donc, d'après le point \emph{8} appliqué à $\Aqi$,
$\Spec\Abul$ s'identifie à $\Min(\Aqi)$, d'où la bijection naturelle entre
$\Spec\gA$ et $\Min(\Aqi)$.

%%%%%%%%%%%%%%%%%%%%%%%%%%%%%%%%%%%%%%%%%%%%%%%%%%%%%%%%%%%%%%%%%

% fin des solutions d'exos

%:   ---- Section*{references}-----------
\Biblio

Des livres de référence pour l'étude des treillis sont \cite{Birkhoff},
\cite{Grae} et~\cite{Johnstone}. Dans \cite{Johnstone} l'accent est
mis essentiellement sur les \trdisz, qui sont les objets qui nous intéressent
au premier chef. Ce livre présente la théorie des locales. La notion de \ix{locale} est
une \gnn de celle d'espace topologique. La structure
d'une locale est donnée par le \trdi de ses ouverts,
mais les ouverts ne sont plus \ncrt des ensembles de points.
C'est la raison pour laquelle une locale est parfois appelée un
\emph{espace topologique sans points} \cite[Johnstone]{Joh}.
L'auteur essaie en \gnl de donner des \prcos et  signale explicitement les 
\thos dont la preuve utilise l'axiome du choix.

En algèbre abstraite, les espaces spectraux sont omniprésents, au premier rang desquels on compte le spectre de Zariski d'un anneau commutatif.
Du point de vue
\cof ce sont des locales bien particulières qui \gui{manquent de points}.
Nous présenterons rapidement cette notion dans la section \iref{secEspSpectraux}
du chapitre \ref{chapKrulldim} consacré à la dimension de Krull.

Une \dem élégante du \thref{thAXgcd} (si $\gA$ un anneau à pgcd intègre
il en va de même pour $\AX$) se trouve dans \cite[th.\,IV.4.7]{MRR}.

L'origine des \entrels se trouve dans le calcul des séquents de Gentzen, qui
mit le premier l'accent sur la coupure (la règle $(T)$).
Le lien avec les \trdis a été mis en valeur dans \cite[Coquand\&al.]{cc,cp}.
Le \tho fondamental des relations implicatives \rref{thEntRel1} se trouve 
dans~\cite{cc}.

On trouve la terminologie de \emph{treillis implicatif} dans \cite{Curry},
celle d'\emph{\agHz} dans \cite{Johnstone}.

Un ouvrage de base pour théorie des \grls
et anneaux réticulés (non \ncrt commutatifs) est  \cite{BKW}.
Nous avons dit\footnote{Dans l'introduction et dans la remarque qui suit le principe de recouvrement fermé~\ref{prcfgrl}.} qu'une idée directrice dans la théorie des \grls est qu'un \grl  se comporte dans les calculs comme un produit de groupes totalement
ordonnés. Cela se traduit en \clama par le \tho de représentation
qui affirme que tout \grl (abélien) est isomorphe à un sous-\grl d'un produit
de groupes totalement ordonnés (\tho 4.1.8 dans le livre cité).

Les \grls qui sont des \Qevs constituent en quelque sorte la version
purement \agq de la théorie des espaces de Riesz. Tout bon livre sur les espaces de Riesz commence par développer les propriétés purement \agqs
de ces espaces, qui sont décalquées (avec des preuves très voisines, sinon
identiques) de la théorie des \grls (abéliens). Voir par exemple \cite{Zaanen}.

\perso{Concernant la structure des \grls libres ou \pf on peut consulter \ldots.
M. Anderson and T. Feil (1988) Lattice-Ordered Groups, D. Reidel, Dordrecht. ???

Et aussi algebres de Heyting libres???}

Dans les exercices de Bourbaki (Algèbre commutative, Diviseurs) un anneau
de Bézout intègre est appelé 
\ixx{anneau}{bezoutien}\index{bezoutien!anneau ---},
un anneau à pgcd intègre est appelé un 
\ixx{anneau}{pseudo-bezoutien}\index{pseudo-bezoutien!anneau ---}, et
un \ddp est appelé un \ixx{anneau}{pruferien}\index{pruferien!anneau ---}.
 
\Llec peut comparer le lemme \ref{lemquotrdi1} et sa démonstration avec le résultat
analogue donné dans le \tho 141 de \cite{Grae}.

Les constructions données dans les propositions \ref{propTrBoo} et \ref{propSumtrdi} se trouvent dans \cite{cc}.

\newpage \thispagestyle{CMcadreseul}
\incrementeexosetprob

%:        %%%%%%%%%%%%%%%%%%%%%%%%%%%%%%%%%%%%
%:        %%%%%%%%%%%%%%%%%%%%%%%%%%%%%%%%%%%%
%---- Chapitre  {Anneaux de Dedekind}------------   
\chapter{Anneaux de Prüfer et de Dedekind} 
\label{ChapAdpc}
%--------------------
\vskip-1em

\minitoc

\subsection*{Introduction} 
\addcontentsline{toc}{section}{Introduction}
%-----------------------------------------

Les \dfns usuelles d'\adk se prêtent
mal à un traitement \algqz.

Premièrement, la notion de \noet est délicate
(du point de vue \algqz). Deuxièmement les questions de \fcn
réclament en \gnl des hypothèses extrêmement fortes.
Par exemple, même si $\gK$  est un \cdi tout à fait explicite, il
n'y a pas de méthode \gnle
(valable sur tous les \cdisz) pour factoriser les \pols de~$\KX$.

Ainsi, un aspect essentiel de la théorie des \adksz, à
savoir que la clôture intégrale d'un \adk dans une
extension finie de son corps de fractions reste un \adkz,
ne fonctionne plus en toute \gnt (d'un point de vue
\algqz) si l'on exige la \fcn complète des \ids (voir par
exemple le traitement de cette question dans~\cite{MRR}).

Par ailleurs, même si une \fac est en
théorie
faisable (dans les anneaux d'entiers des corps de nombres par
exemple), on se heurte très rapidement à des \pbs d'une
complexité rédhibitoire comme celui de factoriser le \discri
(tâche en pratique impossible si celui-ci a plusieurs centaines de
chiffres).
Aussi Lenstra et Buchmann, \cite{LB}, ont-ils proposé de travailler
dans les anneaux d'entiers sans disposer d'une $\ZZ$-base.
Un fait \algq important est qu'il est toujours facile d'obtenir une
\emph{\fapz} pour une famille d'entiers naturels,
\cad une \dcn de chacun de ces entiers en produits de
facteurs pris dans une famille d'entiers premiers entre eux deux à deux 
(voir \cite[Bernstein]{Ber1}, et~\hbox{\cite[Bernstein]{Ber2}}
pour une mise en {\oe}uvre avec
les \ids de corps de nombres, voir aussi le \pbz~\rref{exoPlgb2}).

Un but de ce chapitre est de montrer la validité \gnle d'un
tel point de vue et de proposer des outils  dans ce cadre.

%-% ENTRE NOUS
\entrenous{   
On aimerait bien mettre:

Nous montrons en particulier que la possibilité de calculer des
\faps pour les familles finies d'\itfs se conserve par extensions entières normales.

Mais il nous manque ce \tho pour le moment}
%-% Fin ENTRENOUS

\smallskip Un rôle crucial et simplificateur dans la théorie
est joué par les
\anars (con\-for\-mé\-ment à une intuition de
Gian~Carlo~Rota~\cite{Rota}),
qui sont les anneaux dans lesquels le treillis des \ids est
distributif, et par les \emph{\mlpsz}, qui sont les matrices qui
explicitent la machinerie calculatoire des \itfs \lopsz,
de manière essentiellement \eqve à ce que Dedekind \cite{Ddk2} estimait
être une \prt fondamentale des anneaux d'entiers dans les corps de nombres
(voir \cite[Avigad]{Avi} et le point \emph{3$\,'.$} de notre proposition~\ref{propItfLocprinc}).

\smallskip La volonté de repousser le plus tard possible la mise en
place des hypothèses \noees nous a \egmt incité à 
développer un traitement \cof de nombreux points importants
de la théorie  dans un cadre plus simple et moins rigide que celui des
\adksz. C'est le contexte des anneaux qui ont les deux
\prts suivantes:
%-----------------begin item------------------
\begin{itemize}
\item  les \itfs sont \pros (ceci caractérise ce que nous
appelons un \emph{\adpcz}),
\item  la \ddk est inférieure ou égale 
à~$1$.
\end{itemize}
%-----------------end item------------------
Comme \llec pourra le constater, les preuves ne sont pas rendues 
plus compliquées,
bien au contraire, par cet affaiblissement des hypothèses.

\smallskip De même, nous avons été amenés à étudier les
\adpcs \gui{à \fapz} (dans le cas local, ce  sont les \advs dont le
groupe de valuation est isomorphe à un sous-groupe de~$\QQ$). Nous
pensons que ces anneaux constituent le cadre de travail naturel 
suggéré par Buchman et Lenstra~\cite{LB}.

\smallskip Enfin pour ce qui concerne les \adksz, nous
nous sommes libérés de l'hypothèse usuelle d'intégrité
(car elle se conserve difficilement d'un point de vue \algq par
extension \agqz) et nous avons laissé au second plan la \fac des
\itfs (pour la même raison) au profit du seul caractère
\noez. La \noet implique la \fap des familles d'\itfsz,
qui elle-même implique la dimension $\leq 1$ sous forme 
\covz.

%%%%%%%%%%%%%%%%%%%%%%%%%%%%%%%%%%%%%%%%%%%%%%%%%%%%%%%%%%%%%%%%%%%%%%%%%%%
%--- Sec {Anneaux \aris}---------   
\section{Anneaux \arisz} 
\label{secAnars}
%-----------------------------------------

Rappelons qu'un \anar est un anneau dont les \itfs sont \lops (voir la section \ref{secIplatTf}).
Nous commençons par quelques précisions concernant les \ids \lops dans
un anneau arbitraire.

%:   subsec{Matrice de \lon principale}
\subsec{Idéaux \lopsz, matrice de \lon principale}

Nous reprenons le \thref{propmlm} (énoncé pour les
modules \lmosz) dans le cadre des \ids \lopsz.

%:     Proposition{propItfLocprinc}--    
\begin{proposition} 
\label{propItfLocprinc} {\em (Idéaux  \lopsz)}\\
Soient $x_1$, \ldots, $x_n\in\gA$. 
\Propeq 
%-----------------begin item------------------
\begin{enumerate}
\item  L'idéal $\fa=\gen{\xn}$ est \lopz.
\item  Il existe $n$  \eco $s_i$ de $\gA$
tels que pour chaque $i$, après \lon en $s_i$,   $\fa$ devient
principal, engendré par $x_i$. 
\item  Il existe
une \mlp pour $(\xn)$, \cad une matrice $A = (a_{ij})\in\Mn(\gA)$
qui vérifie:
%---- equation {eqmlp} ----
\begin{equation}\label{eqmlp}\preskip.3em \postskip.3em
\left\{\arraycolsep2pt
 \begin{array}{rcl}
   \;\;\sum a_{ii}&=&1\\[1mm]
   \;\;a_{\ell j}x_{i}& =& a_{\ell i}x_{j} \qquad \forall i,j,\ell \in
 \lrbn
 \end{array}
\right.
\end{equation}
%---------------------end equation--------------
 NB. Le deuxième point  se lit comme suit:
 pour chaque ligne~$\ell$, les mineurs d'ordre~2 de 
la matrice $\cmatrix{
a_{\ell 1}&\cdots &a_{\ell n}\cr x_1&\cdots &x_n}$ sont nuls.

\item  $\Vi_\Ae 2(\fa)=0$.
\item  $\cF_1(\fa)=\gen{1}$.
\end{enumerate}
Dans le cas où   l'un des $x_k$ est \ndz l'existence de la matrice  $A$
dans le point 3 a la même signification que le point suivant.
\begin{enumerate}
\item [3$\,'$.] 
Il existe $\gamma_1$, \ldots, $\gamma_n$ dans $\Frac\gA$ tels que 
$\sum_i\gamma_ix_i=1$ et chacun des $\gamma_ix_j$ est dans $\gA$
 \emph{(formulation de Dedekind)}.
\end{enumerate}
%-----------------end item------------------
\end{proposition}
%--- end-proposition-------------------------
%
\begin{proof}
La seule chose nouvelle est la formulation \emph{3$\,'$.}
Si par exemple $x_1\in\Reg(\gA)$ et si l'on dispose de $A$, on pose
$\gamma_i=a_{i1}/x_1$. 
Réciproquement, si l'on dispose des $\gamma_i$, on pose $a_{ij}=\gamma_ix_j$.
\end{proof}

\medskip 
La proposition suivante reprend et précise la proposition~\ref{pmlm}.
Les résultats pourraient être obtenus de manière plus directe,
en utilisant la formulation de Dedekind,
lorsque l'un des $x_k$ est \ndzz.

%-- Proposition pilps1 ------------   
\begin{proposition}\label{pilps1}
Soit $\fa = \gen{\xn} $ un \id \lop de~$\gA$
et $A= (a_{ij})$ une \mlp pour $(\xn)$.
Nous avons les résultats suivants.
%--------begin item --------------------------
\begin{enumerate}
\item %1
$[\,x_1\; \cdots\;x_n\,]\; A = 
[\,x_1\; \cdots\;x_n\,]$.
\item %2
   Chaque $x_i$ annule $\cD_2(A)$ et $A^2 - A$.
\item %3
  Soit $\gA_i = \gA[1/a_{ii}]$, on a $\fa =_{\gA_i} \gen{x_i}$.
\item %4
 $\gen{\xn} \gen{a_{1j},\dots,a_{nj}} =
\gen{x_j}$.
\item %5
 Plus généralement, si $a = \sum\alpha_{i}x_i$ et
   $\tra{[\,y_1\; \cdots\;y_n\,]} = A~ \tra{[\,\alpha_1\; \cdots\;\alpha_n\,]}$, alors
$$\preskip-.4em \postskip.4em 
\displaystyle
\gen{\xn}\gen{\yn} = \gen{a}
. 
$$
 En outre, si $\Ann(\fa) = 0$, la matrice
$\tra{\!A}$ est une \mlp pour $(\yn)$.

\item %6
 En particulier, si $ \sum\alpha_{i}x_i = 0$ et
   $\tra{[\,y_1\; \cdots\;y_n\,]} = A~ \tra{[\,\alpha_1\; \cdots\;\alpha_n\,]}$, alors

\snic{\displaystyle
\gen{\xn}\gen{\yn} = 0
}
\item %7
 On considère la forme \lin
   ${\ux}: (\alpha_i)\mapsto \sum_i\alpha_ix_i$  associée à $(\xn)$,
   on note $\fN=\Ann\, \gen{\xn}$ et 
   $\fN^{(n)}$ le produit cartésien

\snic{\sotq{(\nu_1,\dots,\nu_n)} {\nu_i\in \fN,\,i\in\lrbn} \subseteq \Ae n.}

%\sni
Alors $\Ker{\ux}=\Im(\I_n-A)+\fN^{(n)}$.
\item %8
 Pour  $i\in\lrb{1..n-1}$ l'intersection
   $\gen{x_1,\dots,x_i}\cap\gen{x_{i+1},\dots,x_n}$ est
l'\itf
   engendré par les $n$ \coes du vecteur ligne 
   
\snic{
[\,x_1\;\cdots\;x_i\;0\;\cdots\; 0\,](\I_n-A) =
-[\,0\;\cdots\;0\;x_{i+1}\;\cdots\; x_n\,](\I_n-A).
}
\end{enumerate}
\end{proposition}

\begin{proof}
Le point \emph{3} est clair, les points \emph{4} et  \emph{6} sont des cas
particuliers de la première partie du point \emph{5}.\\
Les points \emph{1}, \emph{2} et  la première partie du point \emph{5} 
ont été montrés
pour les \mlmosz.

 \emph{5.} Il reste à montrer que, lorsque $\Ann(\fa)=0$, $\tra{A}$ est une \mlp pour~$(\yn)$. En effet, d'une part $\Tr(\tra{A})=1$, et
d'autre part, puisque $\fa\cD_2(A) = 0$, on a $\cD_2(A) = 0$,
ou encore $A_i \wedge A_j = 0$, $A_i$ étant 
la colonne $i$ de $A$. Comme le vecteur $y := \tra{[\,y_1\; \cdots\;y_n\,]}$ est dans~$\Im A$, on a aussi~$y \wedge A_j = 0$, ce qui traduit que~$\tra{A}$ est une \mlp pour~$(\yn)$.

\emph{7.} L'inclusion $\Ker\ux\subseteq\Im(\I_n-A)+\fN^{(n)}$
résulte du point \emph{6} et l'inclusion réciproque du~\emph{1.}

\emph{8.} Résulte de \emph{7} en remarquant que se donner un \elt $a$ de
\linebreak
{l'\id   $\fb=\gen{x_1,\dots,x_i}\cap\gen{x_{i+1},\dots,x_n}$} revient
à se donner
un \elt

\snic{(\alpha_1,\dots,\alpha_n)\in\Ker\ux~:~a=\alpha_1x_1+\cdots +\alpha_ix_i=-\alpha_{i+1}x_{i+1}-\cdots
-\alpha_nx_n.}

%\sni
Ainsi, $\fb$ est engendré par les \coes de  $[\,x_1\;\cdots\;x_i\;0\;\cdots\; 0\,] (\I_n-A)$.
\end{proof}

%:    Corollary{corpilps1}
\begin{corollary}\label{corpilps1}
Soit $\fa = \gen{\xn} $ un \itf de $\gA$.  
\begin{enumerate}
\item Si $\fa$ est \lopz, pour tout \itf $\fc$ contenu dans~$\fa$, il existe 
un \itf $\fb$ tel que $\fa\fb=\fc$.
\item Inversement, si $n=2$ et s'il existe un \id $\fb$ tel que 
$\gen{x_1}=\fa\fb$, alors $\fa$ est \lopz.
\item L'\id $\fa$  est un \mrc $1$ \ssi il est \lop et fidèle.
Dans ce cas, si $A$ est une \mlp pour $(\xn)$, c'est une \mprn 
de rang $1$ et $\fa\simeq \Im A$.
\item L'\id $\fa$ est \iv \ssi il est \lop et contient un \elt \ndzz.
\end{enumerate}
\end{corollary}
\begin{proof}
\emph{1}, \emph{3},  \emph{4.} Voir le lemme \ref{lemIdproj}, qui donne
des résultats un peu plus \gnlsz.
Ces points résultent aussi de la proposition précédente, points \emph{5} et \emph{7.} 
\\
\emph{2.} Dans $\fb$ on doit avoir $u_1$ et $u_2$ tels que 
d'une part $u_1x_1+u_2x_2=x_1$, \linebreak
donc $(1-u_1)x_1=u_2x_2$, et d'autre part $u_1x_2 \in \gen{x_1}$. 
 Lorsque l'on inverse l'\eltz~$u_1$,~$x_1$ engendre $\fa$, et lorsque l'on inverse
$1-u_1$, c'est $x_2$ qui engendre $\fa$.
\end{proof}
%

%:   subsec{Premières proprietes}
\subsec{Premières \prts}

Rappelons qu'un anneau est \coh \ssi d'une part l'intersection de deux \itfs
est un \itfz, et d'autre part l'annulateur de tout \elt est \tf
(\thrf{propCoh4}).  Par conséquent, en utilisant le point \emph{8} de la
proposition \ref{pilps1}, nous obtenons:

%--- Fact{factAnarCoh}------------
\begin{fact} 
\label{factAnarCoh} 
Dans un \anar l'intersection de deux \itfs est un \itfz. Un \anar est \coh \ssi 
l'annulateur de tout \elt est \tfz.
\end{fact}
%--- end-fact-----------------------------------------

Tout quotient et tout localisé d'un \anar est un \anarz.

Dans un anneau \fdiz, la relation de \dve est explicite.
On a la réciproque (remarquable) pour les \anarsz.

%--- proposition aritfdi  -----------
\begin{proposition}\label{aritfdi}
Un \anar est \fdi \ssi la relation de \dve est  explicite.
De manière plus précise, dans un anneau quelconque, si un \id $\gen {b_1,
\ldots, b_n}$ est \lop et si $A = (a_{ij})$ est une \mlp pour $(b_1,\dots,b_n)$,
on a l'\eqvc

\snic {
c \in \gen {b_1, \dots, b_n}  \iff  a_{jj}c \in \gen {b_j}
\hbox { pour tout } j.
}

%\sni
En particulier, on a $1\in\gen{b_1,\dots,b_n}$ \ssi pour tout
$j$, $b_j$ divise $a_{jj}$.
\end{proposition}

\begin{proof}
Si $a_{jj}c = u_jb_j$, alors $c = \sum_j u_jb_j$. Réciproquement, si
$c \in \gen {b_1, \dots, b_n}$, alors pour chaque $j$:

\snic {
a_{jj}c \in \gen {a_{1j},\dots,a_{nj}} \gen {b_1, \dots, b_n} = \gen {b_j}.
}
\end{proof}

Dans le \tho qui suit nous donnons quelques \carns possibles des \anarsz.
La \carn la plus simple des \anars est sans doute celle donnée
dans le point \emph{1b}. Puisqu'un \id $\gen{x,y}$ est \lop \ssi
il y a une \mlp pour $(x,y)$, la  condition \emph{1b}
signifie:
\Grandcadre{$\forall x,y\in\gA\;\;\exists u,a,b\in\gA,\quad\quad ux=ay,\;(1-u)y=bx$,}
ce qui est aussi exactement ce que dit le point \emph{2c}.

%:      Theorem{thAnar}-------------   
\begin{theorem} 
\label{thAnar} \emph{(\Carns des \anarsz)}\\
Pour un anneau $\gA$ \propeq
\begin{enumerate}
\item [1a.] $\gA$  est \ari (tout \itf est \lopz).
\item [1b.] Tout idéal $\fa=\gen{x_1,x_2}$ est \lopz.
\item [2a.] Pour tous \itfs $\fb\subseteq \fa$, il existe un \itf $\fc$ tel 
que $\fa\fc=\fb.$
\item [2b.] Pour tout idéal $\fa=\gen{x_1,x_2}$, il existe un \itf $\fc$ 
tel 
\linebreak
que~$\fa\fc=\gen{x_1}$.
\item [2c.] $\forall x_1,x_2\in \gA$ le \sli $BX=C$ suivant admet une 
solution:
%------- eqSLI -----
\begin{equation}\label{eqSLI}
[ \,B\mid C \,] =\cmatrix{ ~x_1 & x_2  &  0    &\vert&x_1 \cr
                                ~x_2 & 0    &  x_1  &\vert&~0~~}   
\end{equation}
%---------------------end equation--------------
\item [2d.] $\forall x_1,x_2\in \gA$ il existe $u\in \gA$ tel que
%-----------------begin $$----------------
$$ \gen{x_1}\; \cap \; \gen{x_2}\; =\; \gen{(1-u)x_1,ux_2}.
$$
%-----------------end $$------------------
%
\item [3.] Pour tous \itfs $\fa$ et $\fb$, la suite exacte courte 
ci-après est scindée:
%-----------------begin $$----------------
$$ 0\longrightarrow \gA/(\fa\cap \fb) \vers{\delta}\gA/\fa\times 
\gA/\fb\vers{\sigma}
\gA/(\fa+\fb)\longrightarrow 0 
$$
%-----------------end $$------------------
où $\delta : \ov x_{\fa\cap\fb} \mapsto (\ov x_\fa, \ov x_\fb)$ et 
$\sigma : (\ov y_\fa, \ov z_\fb) \mapsto \ov {(y-z)}_{\fa+\fb}.$
\item [4.] Pour tous \itfs $\fa$ et $\fb$, $(\fa:\fb)+(\fb:\fa)=\gen{1}$. 
\item [5.] {\em (\Tho des restes chinois, forme \ariz)}\\ 
Si $(\fb_k)_{k=1,\ldots,n}$ est une famille finie d'idéaux de $\gA$ 
et  $(x_k)_{k=1,\ldots,n}$ est une famille d'\elts de $\gA$
vérifiant $x_k\equiv x_\ell\; \mod\,\fb_k+\fb_\ell$ pour \linebreak
tous $k,\,\ell,$ alors 
il existe un $x\in\gA$ tel que  $x\equiv x_k\; \mod\,\fb_k$ 
pour tout~$k.$   
\item [6.] Le treillis des \ids de $\gA$ est un \trdiz. 
\end{enumerate}
\end{theorem}
%--- end-theorem-----------------------------------------
%-----------------begin proof------------------
\begin{proof}
\emph{1b} $\Rightarrow$ \emph{1a}.
Si l'on a un \itf avec $n$ \gtrsz, des \lons successives (chaque
fois en des \ecoz) le rendent principal.

Considérons le point \emph{2a} Soit $\fa=\gen{x_1,\ldots ,x_n}$  et  
$\fb=\gen{y_1,\ldots ,y_m}$.
Si $\fc$ existe, pour chaque $j=1,\ldots,m$ il existe des \elts $a_{i,j}\in \fc$ 
tels que 

\snic{\som_i a_{i,j} x_i =y_j.}

Par ailleurs, pour chaque $i,i',j$ on doit avoir  $a_{i,j}x_{i'}\in \fb$, ce qui 
s'exprime par l'existence d'\elts $b_{i,i',j,j'}\in \gA$ vérifiant 

\snic{\som_{j'}b_{i,i',j,j'}y_{j'}=a_{i,j}x_{i'}.
}

Réciproquement, si l'on peut trouver des $a_{i,j}$ et $b_{i,i',j,j'}\in \gA$ 
vérifiant les \eqns \lins 
ci-dessus (dans lesquelles les $x_i$ et $y_j$ sont des \coesz), alors l'idéal 
$\fc$ engendré par les $a_{i,j}$ vérifie bien $\fa\fc=\fb$.
Ainsi, trouver $\fc$ revient à résoudre un \sliz.

Il s'ensuit que pour prouver \emph{1a} $\Rightarrow$ \emph{2a}
on peut utiliser des 
\lons convenables:  les deux \ids $\fa$  et $\fb$ deviennent principaux, 
l'un étant inclus dans l'autre, auquel cas $\fc$ est évident.

 On vérifie facilement que les \prts 
\emph{1b}, \emph{2b}, \emph{2c} et  \emph{2d} sont \eqves (en tenant compte de la remarque précédente pour \emph{1b}).\\
Pour montrer que  \emph{1a} implique \emph{3}, \emph{4}, \emph{5} et  \emph{6}, on note que 
chacune des \prts considérées peut s'interpréter comme l'existence 
d'une solution d'un certain \sliz, et que cette solution est évidente lorsque 
les \ids qui interviennent sont principaux et totalement ordonnés pour 
l'inclusion.
  
Il reste à voir les réciproques.
 
\emph{3} $\Rightarrow$ \emph{2c}  et \emph{4} $\Rightarrow$ \emph{2c}. On considère dans \emph{3} ou
\emph{4} le cas où $\fa=\gen{x_1}$ \linebreak
et $\fb=\gen{x_2}$.
 
\emph{5}  $\Rightarrow$ \emph{1b}. 
Soient $a,b\in \gA$. Posons 

\snic{c=a+b, \; 
\fb_1=\gen{a},\;\fb_2=\gen{b}, \; \fb_3=\gen{c}, \; x_1=c, \;x_2=a\hbox{ et }x_3=b.}

 On a $\fb_1+\fb_2=\fb_1+\fb_3=\fb_3+\fb_2=\gen{a,b}$.\\
Les congruences 
$x_i\equiv x_k\; \mod\,\fb_i+\fb_k$ sont vérifiées, donc
il existe $u$, $v$, $w$ dans~$\gA$ tels que 

\snic{c+ua = a+vb = b+wc,}

 d'où

\snic{wb=(1+u-w)a,\; (1-w)a=(1+w-v)b.}

 Donc l'idéal $\gen{a,b}$ est \lopz.
 
\emph{6}  $\Rightarrow$ \emph{1b}. 
Prenons la \prt 
de \dit $\fa+(\fb\cap \fc)=(\fa+\fb)\cap (\fa+\fc)$, avec 
 $\fa=\gen{x}$, $\fb=\gen{y}$ et $\fc=\gen{x+y}$. On a donc $y\in 
\gen{x}+(\gen{y}\cap\gen{x+y})$, \cad qu'il existe $a$, $b$, $c$ tels que $y=ax+by$, 
$by=c(x+y)$. \\
D'où $cx=(b-c)y$ et $(1-c)y=(a+c)x$. Ainsi, $\gen{x,y}$ est \lopz.
\end{proof}
%-----------------end proof------------------

L'\iso $\gA/\fa\oplus \gA/\fb\simeq \gA/(\fa+\fb)\oplus \gA/(\fa\cap\fb)$
qui résulte du point~\emph{3} du \tho précédent admet la \gnn suivante.
%:     Corollary{corthAnar}
\begin{corollary}\label{corthAnar}
Soit $(\fa_i)_{i\in\lrbn}$  une famille d'\itfs d'un \anar $\gA$.
Posons
\[ 
\begin{array}{ccc} 
  \fb_1=\sum_{k=1}^n\fa_k,\; \fb_2=\sum_{1\leq j<k\leq n}(\fa_j\cap \fa_k),\;\dots  \\[2mm] 
\fb_r=\sum_{1\leq j_1<\cdots<j_r\leq n}(\fa_{j_1}\cap \cdots\cap \fa_{j_r}),\;\dots,\;
  \fb_n=\bigcap_{k=1}^n\fa_k.  
  \end{array}
\] 
Alors on a $\fb_n\subseteq\cdots\subseteq \fb_1$ avec un \iso
%-----------------begin $$----------------
$$\bigoplus\nolimits_{k=1}^n\gA/\fa_k \;\simeq \;\bigoplus\nolimits_{k=1}^n\gA/\fb_k.
$$
%-----------------end $$------------------
\end{corollary}
%--------- fin corollary ---------------------------------------------- 

En rapprochant ce résultat du \thref{prop unicyc} on obtient une classification
complète des \Amos de ce type.
On peut aussi comparer avec le fait~\ref{factGpRtcl}~\emph{\iref{i16bisfactGpRtcl}.}

%:     Corollary{cor2thAnar}
\begin{corollary}\label{cor2thAnar}
Soit $\gB$ une \Alg \fptez. Si $\gB$ est un \anar (resp. un \adpz, un \adpcz), alors $\gA$ \egmtz. 
\end{corollary}
%--------- fin corollary ---------------------------------------------- 
%
\begin{proof} Puisque $\gA\subseteq\gB$, si $\gB$ est réduit, $\gA$ \egmtz.
Le \tho \ref{propFidPlatTf}~\emph{3} implique que si $\gB$ est \cohz, $\gA$ \egmtz. Il reste à montrer le résultat pour \gui{\anarz}. \\
On considère  $x$, $y\in\gA$. On doit montrer qu'il existe $u$, $a$, $b\in\gA$ tels \linebreak
que $ux=ay$ et $(1-u)y=bx$.
Il s'agit en fait d'un \sli  à \coes dans $\gA$, avec les inconnues  $(u,a,b)$. 
Or ce \sys admet une solution dans
$\gB$ et~$\gB$ est \fpte sur $\gA$, donc il admet une solution dans $\gA$. 
\end{proof}
%

%:   subsec{Structure multiplicative des itfs}
\subsec{Structure multiplicative des \itfsz}

Rappelons que nous notons $\Ifr\gA$ le \mo multiplicatif des \ifrs \tf d'un anneau
arbitraire $\gA$ (voir \paref{NOTAIfr}). 

A priori une inclusion $\fa\subseteq\fb$ dans $\Ifr\gA$
n'implique pas l'existence d'un \ifr $\fc\in \Ifr\gA$ tel que $\fb\fc=\fa$. Mais ceci
est vérifié dans le cas des \anarsz.

\rdb
Pour $\fa$ et $\fb$ dans $\Ifr\gA$, on note \fbox{$\fa\div\fb=\sotq{x\in\Frac\gA}{x\fb\subseteq\fa}$}.\label{NOTAfadivfb} 

%:     Lemma{lemIfrCoh}
\begin{lemma}\label{lemIfrCoh} 
Soit $\gA$ un \emph{\coriz}. 
\begin{enumerate}
\item $\Ifr\gA$ est un treillis pour la relation d'inclusion, le sup est donné par la somme et le inf par l'intersection. 
\item  $\Ifr\gA$ est un \trdi \ssi l'anneau est \ariz.
\item  Au sujet des \elts \ivs de $\Ifr\gA$.
\begin{enumerate}
\item Si $\fa\,\fa'=\gA$ dans $\Ifr\gA$,  on a $\fa'\fc=\fc\div \fa$ et $\fa(\fc\div \fa)=\fc$ pour \hbox{tout $\fc\in\Ifr\gA$}. En particulier $\gA\div\fa$ est l'inverse de $\fa$.
\item Un \ifr $\fraC \fa a$ (où $\fa$ est un \itf de $\gA$) est \iv  dans  $\Ifr\gA$ \ssi  $\fa$ est un \id \ivz.
\item Si $\fa(\gA\div \fa)=\gA$, $\fa$ est \iv dans $\Ifr\gA$.
\end{enumerate}
\end{enumerate}
Soient  $\fa$, $\fb\in\Ifr\gA$ avec $b\in\fb\cap\Reg\gA$. 
On suppose que $\gA$ est \icl dans $\Frac\gA$.
\begin{enumerate}\setcounter{enumi}{3}
\item On a $\fa\div\fb\in\Ifr\gA$.  
\item Si en outre $\fa\subseteq\fb\subseteq\gA$, alors on a
$\fa\div\fb=\fa:\fb$. 
\end{enumerate}
\end{lemma}

\begin{proof}
Tout \elt de $\Ifr\gA$ s'écrit sous la forme $\fraC \fa a$ pour un \itfz~$\fa$ de~$\gA$ et un $a\in\Reg\gA$. En outre $\fraC \fa a\, \fraC \fb b =\fraC{\fa\,\fb}{ab}$. Enfin l'\elt neutre du \mo est $\gA=\gen{1}$.
Ceci montre les points \emph{1}, \emph{2} et \emph{3b.}

\emph{3a.} On a $\fa\fa'\fc=\fc $
donc $\fa'\fc\subseteq \fc\div\fa$ et $\fc=\fa\fa'\fc\subseteq\fa(\fc\div\fa)=\fc$. 
\\
Si $x\in\fc\div\fa$, i.e. $x\fa\subseteq \fc$, alors $x\gA=x\fa\fa'\subseteq \fa'\fc$, donc $x\in\fa'\fc$.

\emph{3c.} 
Avec $\fa=\gen{a_1,\dots,a_k}\subseteq \gA$, supposons que $\fa(\gA\div \fa)=\gA$. \\
Il existe $x_1$, \dots, $x_k\in (\gA\div \fa)$ tels que
$\sum_ix_ia_i=1$ et $x_ia_j\in\fa$ pour \hbox{tous $i,j$}. 
On peut écrire les
$x_i$ sous la forme $\fraC{b_i}c$ avec un même dénominateur~$c$.
On obtient $\sum_ia_ib_i=c$ et $a_ib_j\in\gen{c}$ pour tous $i,j$.
\\
Ainsi en posant $\fb=\gen{b_1,\dots,b_k}$ on obtient $\fa \,\fb=\gen{c}$.

\emph{5.} L'inclusion $\fa:\fb\subseteq \fa\div \fb$ est \imdez. 
Réciproquement, si un $x\in\gK$ vérifie~\hbox{$x\fb\subseteq \fa$}, nous devons montrer que $x\in\gA$.\\
Comme $\gA$ est \icl dans $\Frac\gA$, on applique le
point \emph{3} du fait~\ref{factEntiersAnn}, \hbox{avec $M=\fb$} \hbox{et $\gB=\Frac\gA$}, \hbox{car $x\fb\subseteq \fa\subseteq \fb$}.

\emph{4.} Résulte du point \emph{5} car
on se ramène au cas traité dans le point~\emph{5}, et dans un \coriz, le transporteur $\fa:\fb$ est \tf si $\fa$ et $\fb$  le sont.
\end{proof}

Le \tho suivant dit que la structure multiplicative
du \mo des \ids \ivs d'un \anar a toutes les \prts souhaitables. 
%:HHH morceau mis dans la preuve rajoutée:
%Cela résulte \imdt
%du corolaire \ref{corpilps1}, du \tho \ref{thAnar} et du lemme~\ref{lem1MonGcd}.

%:HHH petit ajout
Rappelons  que d'après le lemme \ref{lemIdproj}, un \itf est  \prc 1 \ssi il est \lop et fidèle.  

%:     Theorem{thiivanar}
\begin{theorem}\label{thiivanar}
Dans un \anar les \itfs fidèles 
forment un \mo multiplicatif qui est la partie positive d'un \grlz.
Les lois de treillis sont $\fa\vi\fb=\fa+\fb$ et $\fa\vu\fb=\fa\cap\fb$. 
\\
Les \ids \ivs (\cad les \itfs qui contiennent un \elt \ndzz) forment
la partie positive  d'un sous-\grl du précédent.  
\end{theorem}
\begin{proof}
Cela résulte 
du corolaire \ref{corpilps1}, du \tho \ref{thAnar} et du \thref{lem1MonGcd}.
\end{proof}

En fait les deux groupes sont confondus dès que $\gA$ est \qiz,
ou plus \gnlt lorsque les \mrcs $1$ sur $\Frac\gA$
sont libres (\thref{propRgConstant3}, point \emph{2}).

%%%%%%%%%%%%%%%%%%%%%%%%%%%%%%%%%%%%%%%%%%%%%%%%%%%%%%%%%%%%%%%%%%%%%%%%%%%
%--- Section{Elements entiers et localisation}     
\section{\'Eléments entiers et localisation} 
\label{subsecEntiers}
%-----------------------------------------

La \dfn suivante généralise la \dfn \ref{defEntierAnn0} dans deux directions.

%:    Definition{defPropACO}-------  
\begin{definition} 
\label{defPropACO} Soit $\varphi:\gA\to\gC$ un \homo entre anneaux 
commutatifs
et $\fa$ un \id de $\gA$. %Notons $\fa'=\varphi(\fa)\gC$ l'extension de $\fa$.
%-----------------begin item------------------
\begin{enumerate}
\item  Un \elt $x\in\gC$  est dit \ixc{entier}{element@élément --- 
sur un idéal} 
sur $\fa$ s'il existe un entier $k\geq 1$  tel que 
$$\preskip.2em \postskip.1em 
\;\;\;x^k=\varphi(a_1)x^{k-
1}+\varphi(a_2) x^{k-2}+\cdots +\varphi(a_k)\eqno(*)
$$
avec  les~$a_h\in\fa^h$.\\
Dans le cas où $\gC=\gA$, cela équivaut à $\big(\fa+\gen{x}\!\big)^k=\fa\big(\fa+\gen{x}\!\big)^{k-1}$.
\\
On dit aussi que l'\egt $(*)$  est \emph{une \rdiz}
de $x$ sur $\fa$.% 
\index{relation de dépendance!intégrale}

\item  Un idéal $\fa$ de $\gA$ est dit \ix{intégralement 
clos} dans $\gC$ si tout \elt de~$\gC$ entier sur $\fa$  est 
dans~$\varphi(\fa)$.%
\index{ideal@idéal!intégralement clos}
\item  L'anneau $\gA$ est dit \ixc{normal}{anneau ---} si tout idéal 
principal de $\gA$ est \icl dans $\gA$.%
\index{anneau!normal}
\end{enumerate}
%-----------------end item------------------
\end{definition}
%--- end-definition-------------------------

Dans tous les cas, un anneau normal est \icl dans son anneau total
de fractions. On a la réciproque partielle suivante.

%:     fact{lemNormalIcl}
\begin{fact}\label{lemNormalIcl}
Un anneau \qi est normal \ssi il est \icl dans son anneau total de fractions.
\end{fact}
%%%%%%%%%%%%%%%%%%%%%%%%%%%%%%%%%%%%%%%%%
%
\facile

Il est clair que tout anneau normal est réduit (car un nilpotent est 
entier sur $\gen{0}$). On a même un peu mieux.

%:     Lemma{lemiclplat}------------- 
\begin{lemma} 
\label{lemiclplat} 
Tout anneau normal est \lsdzz.
Plus \prmtz, on a pour tout anneau  $\gA$  les implications  
{1} $\Rightarrow$ {2} $\Rightarrow$ {3.}
%-----------------begin item------------------
\begin{enumerate}
\item  Tout \idp est \icl (i.e. $\gA$ est 
normal).
\item Pour tous  $x$, $y\in \gA$, si $x^2\in \gen{xy}$, alors $x\in \gen{y}.$
\item  Tout \idp est plat (i.e. $\gA$ est \lsdzz).
\end{enumerate}
%-----------------end item------------------
\end{lemma}
%--- end-lemma-----------------------------------------
%-----------------begin proof------------------
\begin{proof}
Notons que l'idéal $0$ est \icl \ssi l'anneau est réduit. 
On a évidemment \emph{1} $\Rightarrow$ \emph{2}, et \emph{2} 
implique que l'anneau est réduit.
Supposons  \emph{2} et  
soient  $x$, $y\in \gA$ tels que $xy=0$. On a $x^2 =x(x+y)$ 
\linebreak
donc  
$x\in \gen{x+y}$, e.g., $x=a(x+y)$. 
Alors $(1-a)x=ay$, 
$ay^2=(1-a)xy=0$.\\
Et puisque l'anneau est réduit, $ay=0$,
puis $(1-a)x=0$. 
\end{proof}
%-----------------end proof------------------

%--- Fact{factEntiers}------------- 
\begin{fact}\label{factEntiers}
Soient $x$ un \elt et $\fa$ un \id de $ \gA$.
Pour les \prts qui suivent on a {2} $\Rightarrow$ {1},
et {1} $\Rightarrow$ {2} si $\fa$ est fidèle et \tfz. 
\begin{enumerate}
\item
L'\elt $x$ est entier sur l'idéal~$\fa$. 
\item 
Il existe un \Amo fidèle $M$ \tf tel que $xM \subseteq \fa M$. 
\end {enumerate}
\end{fact}
%%%%%%%%%%%%%%%%%%%%%%%%%%%%%%%%
\begin{proof} (Comparer à la \dem du fait \ref{factEntiersAnn}.)\\
\emph{2} $\Rightarrow$ \emph{1}. 
On considère une
matrice $A$ à \coes dans $\fa$ qui repré\-sente~$\mu_{M,x}$ (la multiplication par $x$ dans $M$) sur 
un \sgr fini de~$M$. Si~$f$ est le \polcar de $A$, on a par le \tho de 
Cayley-Hamilton  $0=f(\mu_{M,x})=\mu_{M,f(x)}$, 
et puisque le module est fidèle, 
$f(x)=0$.\\ 
\emph{1} $\Rightarrow$ \emph{2}. Si l'on a une \rdi de degré $k$
de $x$ sur $\fa$ on prend $M = (\fa+\gen{x})^{k-1}$.
\end{proof}

%%%%%%%%%%%%%%%%%%%%%%%%%%%%%%%%
\rdb

Soit un idéal $\fa$ de $\gA$ et une \idtr $t$, la sous-\alg
$\gA[\fa t]$ de $\gA[t]$, \cad \prmt
$$
\gA[\fa t] = \gA \oplus \fa t \oplus \fa^2 t^2 \oplus \cdots
$$
est appelée l'\emph{\alg de Rees
de l'idéal $\fa$}.%
\index{algebre@\alg de Rees!de l'idéal $\fa$}% 
\label{NOTARees}

La \dem des deux faits suivants est laissée \alecz.
%--- Fact{fact2Entiers}------------- 
\begin{fact}\label{fact2Entiers}
Soit $\fa$ un \id de $ \gA$. 
\begin{enumerate}
\item Pour $x \in \gA$, \propeq
\begin{enumerate}
\item
L'\elt $x $ est entier sur l'\id $\fa$ de $\gA$. 
\item 
Le \pol $xt$ est entier sur la sous-\alg $\gA[\fa t]$ de $\gA[t]$. 
\end {enumerate}
\item De manière plus précise:
\begin{enumerate}
\item  Si $\ov \fa$ est l'ensemble des
\elts de $\gA$ entiers sur  $\fa$, alors la \cli de $\gA[\fa t]$ dans $\gA[t]$ est le sous-anneau
$\gA[\ov\fa t]$.%
\index{cloture@clôture intégrale!de l'idéal $\fa$ dans $\gA$}
\item En particulier,  $\overline \fa$ est un idéal de $\gA$, appelé
la \emph{\cli de l'idéal $\fa$ dans $\gA$}.
On le note $\Icl_\gA(\fa)$ ou $\Icl(\fa)$.%
\end{enumerate}
\end{enumerate}
\end{fact}
%

%:     Fact{fact3Entiers}
\begin{fact}\label{fact3Entiers}
Soient $\fa$ et $\fb$ deux \ids de $\gA$.
\begin{enumerate}
\item $\Icl\big(\Icl(\fa)\big)=\Icl(\fa)$.
\item $\fa\Icl(\fb)\subseteq\Icl(\fa)\Icl(\fb)\subseteq\Icl(\fa\fb).$
\end{enumerate}
\end{fact}

Nous revisitons maintenant deux résultats importants déjà établis.

Le point \emph{2c} du \tho de \KRO \ref{thKro} donne \prmt
le résultat suivant.

%:2012  le lemme ci dessous est un peu mieux dans la nouvelle formulation

%%:     Lemma{lemthKroicl}
%\begin{lemma}\label{lemthKroicl}\emph{(\Tho de \KROz, reformulé)}
%Si l'on a dans $\AT$
%
%\snic{f=\som_{i=0}^nf_iT^{n-i},\;\;  g=\som_{j=0}^m g_j T^{m-j} \;\; \mathit{et} \;\; h=fg=\som_{r=0}^{m+n} h_rT^{m+n-r},
%}
%
%%\sni
%alors, chaque $f_ig_j\in\Icl\big(\rc_\gk(h)\big)$,
%où $\gk$ est le sous-anneau engendré par les~$h_r$. 
%\end{lemma}

%:     Lemma{lemthKroicl}
\begin{lemma}\label{lemthKroicl}\emph{(\Tho de \KROz, reformulé)}\\
Supposons que l'on a dans $\AT$ une \egt

\snic{f=\som_{i=0}^nf_iT^{n-i},\;\;  g=\som_{j=0}^m g_j T^{m-j} \;\; \mathit{et} \;\; h=fg=\som_{r=0}^{m+n} h_rT^{m+n-r}.
}

Soit $\gk$ le sous-anneau de $\gA$ engendré par les~$f_ig_j$.
Alors, chaque $f_ig_j$ est entier sur l'\id $\rc_{\gk}(h)$ de $\gk$.
\end{lemma}

Notez que le point \emph{2c} du \tho de Kronecker \ref{thKro} nous dit \prmt
ceci:
\emph{il existe un \pog $R_{i,j}\in \ZZ[Y,H_0,\ldots,H_p]$
 (toutes les variables ont le même poids $1$), unitaire en $Y$, tel que}
 
 \snic{R_{i,j}(f_ig_j,h_0,\dots,h_p)=0.}

\smallskip 
Voici une nouvelle version du lemme \gui{lying over} (lemme~\ref{lemLingOver}).%
\index{lying over}

%:     Lemma{lemLingOver2}-------------- 
\begin{lemma} \emph{(Lying over, forme plus précise)}\label{lemLingOver2}\\
 Soit  $\gA\subseteq\gB$ avec $\gB$ entier sur $\gA$
et $\fa$ un idéal de $\gA$, alors $\fa\gB\cap\gA\subseteq\DA(\fa)$.
Plus \prmtz, tout \elt de $\fa\gB$ est entier sur $\fa$.
\end{lemma}
%--- end-lemma-----------------------------------------
%-----------------begin proof------------------
\begin{proof}
On reprend textuellement la preuve du lemme \ref{lemLingOver}.  Si
$x\in\fa\gB$, on~a~$x=\sum a_ib_i$, avec $a_i\in \fa,\;b_i\in \gB$.  Les $b_i$
engendrent une sous-\Algz~$\gB'$ qui est finie.  
Soit $G$ un \sys \gtr fini (avec $\ell$ \eltsz) du~\Amo $\gB'$.  
Soit $B_i\in\MM_\ell(\gA)$ une matrice
qui exprime la multiplication par~$b_i$ sur $G$.  
La multiplication par~$x$
est exprimée par la matrice $\sum a_iB_i$, qui est à \coes dans $\fa$.  Le
\polcar de cette matrice, qui annule~$x$ (parce que $\gB'$ est un \Amo
fidèle), a donc son \coe de degré $\ell-d$ dans~$\fa^d$.

On pourrait \egmt appliquer le fait \ref{factEntiers} en prenant
$M = \gB'$. En effet,  comme $x \in \fa\gB'$, on a bien $x\gB' \subseteq \fa\gB'$
et donc $x$ est entier sur $\fa$.
\end{proof}
%-----------------end proof------------------

Nous examinons maintenant les rapports entre \prts de type \gui{entier sur} 
et \lonsz.

%--- Fact{fact.loc.normal}--------- 
\begin{fact} 
\label{fact.loc.normal}\relax 
Soit  
$\fa$ un \id de $\gA$,  
$S$ un \mo de $\gA$ et $x\in\gA$.% et $s\in S$.
%-----------------begin item------------------
\begin{enumerate}
\item   L'\elt $x/1\in\gA_S$ est entier sur $\fa_S$ \ssi il existe 
$u\in S$ tel que
$xu$ est entier sur $\fa$ dans $\gA$.

\item  Si $\gA$ est normal, alors $\gA_S$ \egmtz.

\end{enumerate}
Soit $\gB\supseteq\gA$ une \alg \fptez.
\begin{enumerate}\setcounter{enumi}{2}
\item  Si $\gA'$ est la \cli de $\gA$ dans $\gB$,
alors $\gA'_{S}$ est la \cli de $\gA_S$ 
dans~$\gB_{S}$.

\item  Si $\gB$ est normal, alors $\gA$ est normal.
\end{enumerate}
%-----------------end item------------------
\end{fact}
%--- end-fact-----------------------------------------
%: preuve en exo
\penalty-2500
\begin{proof} On montre seulement le point \emph{1.}
Dans la \dem on confond un \elt de $\gA$ et son image dans $\gA_S$
pour alléger les notations.
Si une  \egt 
$x^k=a_1x^{k-1}+ a_2 x^{k-2}+ \cdots +a_k$ est réalisée dans 
 $\gA_S$ avec chaque $a_j\in(\fa\gA_S)^j$,
on obtient \gui{en réduisant toutes les fractions au même \deno
et en chassant le \denoz} 
une \egt 

\snic{sx^k=b_1x^{k-1}+ b_2 x^{k-2}+ \cdots +b_k}

%\sni
dans $\gA_S$ avec $s\in S$ et chaque $b_j\in\fa^j$. Ceci signifie une \egt dans $\gA$
après multiplication par un autre \eltz~$s'$ de~$S$. On peut
aussi bien multiplier par $s'^ks^{k-1}$ et l'on obtient avec~$u=ss'$ une \egt  

\snic{(xu)^k=c_1(xu)^{k-1}+ c_2 (xu)^{k-2}+ \cdots +c_k}

%\sni
dans~$\gA$ avec chaque $c_j\in\fa^j$.
\end{proof}

Le fait qu'un anneau est normal  est une notion locale, au sens 
suivant.

%:    Principe local global concret{plcc.normal}---
\begin{plcc} 
\label{plcc.normal}\relax \emph{(Anneaux normaux)}\\
Soient $S_1$, $\ldots$, $S_n$  des \moco d'un anneau $\gA$, $x\in \gA$ et $\fa$ un \id de $\gA$.   
%-----------------begin enum------------------
\begin{enumerate}
\item  L'\elt $x$ est entier sur $\fa$  \ssi il est entier sur chacun des~$\fa_{S_i}$.
\item  L'\id $\fa$ est \icl dans $\gA$ \ssi chacun des~$\fa_{S_i}$ est \icl
dans~$\gA_{S_i}$.
\item  L'anneau $\gA$ est normal \ssi chacun des $\gA_{S_i}$ est normal.
\end{enumerate}
%-----------------end enum------------------
\end{plcc}
%--- end-plcc-----------------------------------------
%-----------------begin proof------------------
%: preuve en exo
\begin{proof}
Il suffit de montrer le point \emph{1}, passage du local au global.
\perso{on devrait le laisser \alec?}
On obtient en appliquant le fait \ref{fact.loc.normal} pour chaque $i\in\lrbn$
un $s_i\in S_i$ tel que $s_ix$ est entier sur l'\id $\fa$ dans $\gA$.
On peut supposer que toutes les \rdis ont le même degré $k$.
\'Ecrivons ces \rdis
$$
%s_i^kx^k=\som_{i=1}^nc_{i,h}s_i^{n-h}x^{n-h}, \quad(c_{i,h}\in\fa^h), \qquad i\in\lrbn.
(s_ix)^k \in \som_{h=1}^k \fa^h (s_ix)^{k-h},
\qquad i\in\lrbn.
$$  
Une \coli de ces relations basée sur une \egt $\som_{i=1}^nb_is_i^k=1$
nous donne une \rdi de $x$ sur $\fa$ dans~$\gA$.
\end{proof}
%-----------------end proof------------------

Notons que puisque la \prt dans le point \emph{1} est \carfz, 
le lemme~\ref{factPropCarFin} nous dit que
le \plgc précédent est \eqv au \plga correspondant (dans lequel
intervient la \lon en n'importe quel \idema de $\gA$).

%%%%%%%%%%%%%%%%%%%%%%%%%%%%%%%%%%%%%%%%%%%%%%%%%%%%%%%%%%%%%%%%%%%%%%%%%%%
%%%%%%%%%%%%%%%%%%%%%%%%%%%%%%%%%%%%%%%%%%%%%%%%%%%%%%%%%%%%%%%%%%%%%%%%%%%
%--- Section{Anneaux de Prufer}  secAdP
\penalty-2500   
\section{Anneaux de Prüfer} 
\label{secAdP}
%-------------------------

Rappelons qu'un anneau est de Prüfer lorsque ses idéaux sont plats, ou
s'il est \ari et réduit, ou encore s'il est \ari et \lsdz   
(proposition~\ref{prop.itfplat 2}).

%:     Definition{defiAdv}
\begin{propdef}\label{defiAdv}
On appelle \ixx{anneau}{de valuation} un anneau $\gA$ vérifiant l'une des \prts \eqves suivantes.%
\index{valuation!anneau de ---}
\begin{enumerate}
\item $\gA$ est un anneau de Bézout local réduit.
\item $\gA$ est un  \adp local.
\item $\gA$ est réduit et vérifie: pour tous $a$, $b\in\gA,$ $a\divi b$ ou $b \divi a$.
\end{enumerate}
Si $\gK=\Frac\gA$, le groupe quotient $\gK\eti/\Ati$ est muni de la relation
d'ordre total $\ov x \divi \ov y$ définie par $\exists a\in\Reg(\gA)$, $y=ax$. Ce groupe totalement ordonné est appelé le \ixx{groupe}{de valuation} de $\gA$.%
\index{valuation!groupe de ---}    
\end{propdef}

En outre, $\gA$ est alors \sdzz.

%:2014  rajout de l'exemple ci-dessous
\medskip
\exl \label{exempleadvgroupe}
Soit $\gk$ un \cdi non trivial et $(\Gamma,\cdot,1_{\Gamma})$ un groupe totalement ordonné discret noté multiplicativement.
On fabrique une \klg qui est un \ddv avec $\Gamma$ pour groupe de valuation comme suit. On considère tout d'abord la \klg $\gA=\gk[\Gamma^{+}]$ décrite dans l'exercice~\ref{exoAlgMon}.
\\
Pour un \elt $a=\sum_{i}a_i\gamma_i$ de $\Atl$ on définit $v(a)$ comme le plus petit $\gamma_i$ qui intervient dans l'écriture de $a$ (on a pris les $\gamma_i$ deux à deux distincts et les $a_i\neq 0$). On vérife alors que $v(ab)=v(a)v(b)$, ce qui implique que $\gA$ est intègre. On pose aussi $v(0)=+\infty$.
Enfin notre \adv est le sous-anneau $\gV=\sotq{\fraC a b}{a\in\gA, b\in\Atl, v(a)\geq v(b)}$ de $\Frac\gA$.
\eoe

\medskip 
Nous donnons maintenant quelques autres \prts \caras des \adpsz, qui s'ajoutent
à celles que l'on peut obtenir à partir du \tho \ref{thAnar} pour les \anarsz.

%:    Theorem{thPruf}------------- 
\begin{theorem} 
\label{thPruf} \emph{(\Carns des \adpsz)}\\
Pour un anneau $\gA$ \propeq
%-----------------begin item------------------
\begin{enumerate}
%:2012 1a) remplace par 1a. etc
\item [1a.] $\gA$ est un \anar \lsdz (i.e., un \adpz).
\item [1b.] $\gA$ est  \lsdz et pour tous $x$, $y$ il existe $n\in\NN\etl$
et un \id $\fb$ tels que $\gen{x,y}\fb=\gen{x^n}$.
\item [2a.] Tout sous-module d'un \Amo plat est plat.
\item [2b.] $\gA$ est  \lsdz et tout module sans torsion est plat.
\item [3a.] Un \sli $BX=C$ arbitraire, 
       dès que les \idds de $[\,B \mid C\,]$ sont égaux à 
       ceux de $B$, admet une solution. 
\item [3b.] Même chose en se limitant à $B\in \Ae {2\times 3}$ et 
$C\in \Ae {2\times 1}$.
\item [4a.] Tout \id est \iclz.
\item [4b.] Tout \itf est \iclz. 
\item [4c.] Tout \id $\gen{x,y}$ est \iclz.
\item [4d.] $\gA$ est normal et pour tous $ x,y\in\gA$, on a $xy\in\gen{x^2,y^2}$.
\item [5a.] Si $\fa$, $\fa'$ et $\fc$ sont des \itfsz, on a l'implication:

\snic{ \fa+\fa'\subseteq\fc$, 
 $\;\fa\fc\subseteq\fa'\fc$ $\;\Longrightarrow\;$  $\fa\subseteq\fa'.}
\item [5b.] Si $\fa$, $\fa'$ et $\fc$ sont des \itfsz, on a l'implication:

\snic{ \Ann(\fa+\fa')\supseteq\Ann(\fc)$, 
  $\;\fa\fc\subseteq\fa'\fc$ $\;\Longrightarrow\;$  $\fa\subseteq\fa'.}

\end{enumerate}
\end{theorem}
%--- end-theorem-----------------------------------------
%-----------------begin proof------------------
\begin{proof} On s'occupe d'abord des \eqvcs entre \emph{1}, \emph{2}
et \emph{3}.
\\
Les implications \emph{1a} $\Rightarrow$ \emph{1b},  \emph{2a} $\Rightarrow$ \emph{1a} et  
\emph{3a} $\Rightarrow$ \emph{3b} sont évidentes.

\emph{1b} $\Rightarrow$ \emph{1a.} Résulte du lemme \ref{lemleszlop} 
ci-après.

\emph{3b} $\Rightarrow$ \emph{1a.} 
L'anneau est \ari parce que le \sli (\ref{eqSLI}) dans le 
\thrf{thAnar} admet une solution. 
En outre, l'anneau est réduit: si~$a^2=0$, le 
\sli $\so{\,ax=0,\,0x=a\,}$ admet une solution
car cela correspond~à 
$$  
   B=\cmatrix{a&0&0\cr0&0&0},\,\,C = \cmatrix{0\cr a}\,\,\,\,
\mathrm{avec}\,\,\cD_2([\,B\,\vert\, C\,])=\cD_2(B)=0\,\,!
$$

\emph{1a} $\Rightarrow$ \emph{3b}. Supposons tout d'abord l'anneau local.
Donc l'anneau est \sdz et tout \itf est principal. Alors, on peut conclure par le 
lemme \ref{lem.solsli} ci-après. Dans le cas \gnlz, la preuve du lemme 
peut être reproduite après des \lons en des \moco convenables, et 
comme il s'agit de résoudre un \sli le \plg de base s'applique.\iplg

\emph{2b} $\Rightarrow$ \emph{2a.} Un module plat est sans torsion
(lemme \ref{lem.platsdz}). 
Tout sous-module d'un module sans torsion est sans torsion, donc plat.

\emph{1a} $\Rightarrow$ \emph{2b.} 
Soit $M$ un module 
sans torsion sur un \adpz.
Nous voulons montrer qu'il est plat.
Supposons tout d'abord l'anneau local. \\
Soit $LX=0$ une \rdl avec 
$L=[\,a_1\;\cdots\; a_m\,]$ dans $\gA$ et $X\in M^{m\times 1}$. \Spdgz, on suppose \hbox{que $a_i=b_ia_1$} \hbox{pour $i>1$}. 
La \rdl se réécrit $a_1y=0$ \linebreak 
avec $y=x_1+b_2x_2+\cdots+b_mx_m$. 
Le sous-module monogène $\gA y$ est plat et l'anneau est local donc $a_1=0$ ou 
$y=0$. Dans le premier cas $L=0$. Dans le deuxième cas  $X=HX$ et $LH=0$ avec la 
matrice triangulaire $H$ suivante~:
$$ H=
\cmatrix{ 
   0 &   -b_2   &   -b_3   &  \ldots   &   -b_m   \cr 
   0 &    1     &    0     &  \ldots   &   0   \cr 
\vdots&    \ddots      &  \ddots    &    \ddots  &   \vdots   \cr 
\vdots &       &  \ddots    & \ddots    &   0   \cr 
   0 & \ldots &  \ldots    &  0   &    1  \cr 
}.
$$
Dans le cas d'un \adp arbitraire, on reprend le raisonnement précédent en 
utilisant les \lons (en des \ecoz) qui ren\-dent l'idéal  
$\gen{a_1,\ldots,a_m }$ engendré par l'un des $a_i$.

  On passe maintenant aux \eqvcs entre \emph{1}, \emph{4} et \emph{5}.
\\
Les implications \emph{4a} $\Leftrightarrow$ \emph{4b} $\Rightarrow$ \emph{4c} $\Rightarrow$ \emph{4d} et \emph{5b} $\Rightarrow$ \emph{5a} sont \imdesz.

 \emph{4d} $\Rightarrow$ \emph{1a.} 
L'anneau $\gA$ est \lsdz (lemme \ref{lemiclplat}). 
Il suffit donc de montrer que tout \id 
$\fa=\gen{x,y}$ est \lopz.
On a $xy=ax^2+by^2$, et $z=ax$ vérifie $z^2=zy-aby^2$. 
Donc,
puisque l'anneau est normal, $ax=a'y$ pour
un certain $a'$. 
De même, $by=b'x$ pour un certain~$b'$.
%Donc $\fa\fb=\gen{xy}$ où $\fb=\gen{xa,yb}$. 
Donc, $xy(1-a'-b')=0$.
Les  \elts $1-a'-b'$, $a'$ et $b'$ sont \comz. Lorsque l'on inverse $1-a'-b'$, on obtient 
$xy=0$, et  après deux nouvelles
\lonsz, $x=0$ ou $y=0$, donc $\fa$ est principal. Lorsque l'on inverse~$a'$, on obtient 
$\fa=\gen{x}$ car $a'y=ax$. Même chose lorsque l'on inverse $b'$.
  
 \emph{1a} $\Rightarrow$ \emph{4b}.   
Soit $x\in \gA$ entier sur un \itf $\fa$. On a pour un
certain~\hbox{$n\in\NN,$} $\fa (\fa +\gen{x})^n=(\fa +\gen{x})^{n+1}$. Puisque l'anneau est
\ariz, on a
un  \itf $\fb$ tel que  $(\fa +\gen{x})\fb=\gen{x}$. Donc en multipliant \hbox{par $\fb^n$}
on obtient
$x^n\fa =x^n(\fa +\gen{x})$ ce qui signifie qu'il existe un $y\in \fa$ tel que
$x^{n+1}=x^ny$ \cadz~\hbox{$x^{n}(y-x)=0$}. Puisque l'anneau est \lsdzz, cela
implique qu'après des \lons \come on a $x=0$ ou $y-x=0$, et dans chaque cas $x\in \fa$.

\emph{5a} $\Rightarrow$ \emph{4b.} 
Soit $x\in \gA$ entier sur un \itf $\fa$. On a pour un
certain~$n\in\NN,$ $\fa (\fa +\gen{x})^n=(\fa +\gen{x})^{n+1}$.
On applique plusieurs fois la \prt de simplification avec l'\id $\fc=\fa+\gen{x}$ et l'on obtient en fin de parcours~\hbox{$\fa+\gen{x}\subseteq\fa$}.

\emph{4b} $\Rightarrow$ \emph{5b}.
Soient $\fc ,\,\fa ,\,\fa' $ trois \itfs vérifiant l'hypothèse dans~\emph{5b}.
Soit~$x$ un \elt de $\fa$ et $X$ un vecteur colonne formé par un
\sgr de $\fc$. Puisque $x\fc \subseteq \fa'\fc $, il existe une
matrice $G \in \Mn(\fa')$ telle \hbox{que $xX=GX$}, i.e.
$(x\In - G) X=0$. Si $P$ est le \polcar de $G$, on a d'une part
$P(x)X=0$, et d'autre part $P(x) \in x^n + \fa'$.
\\ 
Donc $P(x)\in\Ann(\fc)\subseteq \Ann(\fa +\fa')$ et
$P(x)\in \fa +\fa'$. D'où $P(x)^2=0$, \hbox{puis
$P(x)=0$}. Ceci est une \rdi de $x$ sur $\fa'$. Donc $x\in \fa'$.
\end{proof}
%-----------------end proof------------------

%:     Lemma{lemleszlop}
\begin{lemma}\label{lemleszlop}
Dans un anneau \lsdzz, si \linebreak 
{l'on a $\gen{x,y}\fb=\gen{x^n}$} avec $n\geq1$, 
alors $\gen{x,y}$ est \lopz.
\end{lemma}
\begin{proof}
Il suffit de résoudre ce \pb après des 
\lons \comez. La \crc \lsdz de l'anneau va servir à fabriquer ces
\lonsz.

On a une \egt $\gen{u,v}\gen{x,y}=\gen{x^n}$ avec $x^n=ux+vy$, 
$ux=u_1x^n$, $vx=ax^n$ et $uy=bx^n$.  Il vient:

\snic {
(u_1y - bx)x^n  = 0, \quad 
(u_1x + ay - x)x^n = (ux + vy - x^n)x = 0.
}

On a donc des \lons \come dans lesquelles $x=0$
et le résultat est clair. Dans la dernière,
$u_1y = bx$ et $u_1x+ay = x$ i.e. $(1-u_1)x=ay$.
Ainsi, $\gen{x,y}$ est \lopz.
\end{proof}
%
%:      lemma} \label{lem.solsli}\relax
\begin{lemma}\label{lem.solsli}\relax
Soient $\gA$ un anneau arbitraire, $B\in \Ae {m\times n}$ et $C \in \Ae
{m\times 1}$. Le \sli $BX=C$ admet une solution dans $\Ae
{n\times 1}$ lorsque les conditions suivantes sont réalisées pour tout
$k\in\lrb{1..\inf(m,n)}$:
%-----------------begin item------------------
\begin{enumerate}
\item  
L'\idd $\cD_k(B)$ est de la forme $\delta_k\gA$, où $\delta_k$ est un mineur
d'ordre $k$.
\item  $\delta_k$ vérifie la condition~: 
$ \forall y\in \gA\; \; (y\delta_k=0\; \Rightarrow \; \big(\delta_k=0\; \lor\; y=0)\big)$.
\item  $\cD_k([\,B \mid C\,]) = \cD_k(B)$.
\end{enumerate}
%-----------------end item------------------
\end{lemma}
%--- end-lemma-----------------------------------------
%-----------------begin proof------------------
\begin{proof}
On commence avec $k=\inf(m,n)$. On écrit l'\idt à la Cramer 

\snic{\delta_k\times C=\delta_k\,\times $ (une \coli  des 
 colonnes  de  $B),}

%\sni

qui résulte de la nullité des \idds d'indice $k+1$ et du fait que 
$\cD_k([\,B \mid C\,])$ est engendré par $\delta_k$. Vu \emph{2}, on est 
dans l'un des deux cas suivants:\\ 
-- on peut simplifier  en divisant tout par $\delta_k$, donc $C \in \Im B$.
\\
-- $\delta_k=0$: ou bien $k=1$ auquel cas $C \in \Im B$ (car $B = C = 0$),
ou bien~$k \ge 2$, et l'on peut faire une \recu en remplaçant
$k$ par $k-1$.
\end{proof}
%-----------------end proof------------------

%%%%%%%%%%%%%%%%%%%%%%%%%%%%%%%%%%%%%%%%%%%%%%%%%%%%%%%%%%%%%%%%%%%%%%%%%%%
%:   subsec{Extensions  d'\adpsz}%%%
\subsec{Extensions  d'\adpsz}

Le fait qu'un anneau normal est \lsdz signifie que \lot il se comporte comme un anneau \sdzz. En fait la machinerie de \lons \come à l'{\oe}uvre
dans la \dfn d'un anneau \lsdz permet souvent de se ramener au cas intègre,
comme on a déjà pu le voir dans la \dem du lemme~\ref{lemleszlop}.

On a le \tho important suivant, qui est une \gnn du résultat analogue
obtenu en théorie des nombres (\thrf{th1IdZalpha}).

%       EXTENSION ENTIERE NORMALE D'UN PRUFER
%:     Theorem{thExtEntPruf}
\begin{theorem}\label{thExtEntPruf} \emph{(Extension entière normale
d'un \adpz)}.\\
Soient $\gA\subseteq\gB$ avec $\gB$ normal entier sur $\gA$ et $\gA$ de Prüfer.
Alors $\gB$ est un \adpz. 
\end{theorem}
\begin{proof}
On va montrer que tout \id  $\gen{\alpha,\beta}$   
est \lopz. 

 Voyons d'abord le cas 
d'un \id $\gen{a,\beta}$ avec $(a,\beta)\in\gA\times\gB$.
On peut alors 
reprendre presque mot à mot la
\dem \gui{à la Dedekind\footnote{Cela fonctionnerait aussi avec la \dem à la Kronecker.}} du \thoz~\ref{th1IdZalpha}.
\\
Soit $f \in \gA[X]$ \mon s'annulant en $\beta$. On écrit
$f(X)=(X-\beta)h(X)$
où~$h \in \gB[X]$. On a donc $f(a X)=(a X-\beta )h(a X)$, que l'on écrit 
$f_1=g_1h_1$.  
Soit $\fc=\rc_\gA(f_1)$,  $\fb=\rc_\gB(h_1)$ et
$\fa=\rc_\gB(g_1)=\gen{a,\beta}$.\\
Si $\deg(f)=n$, on a $a^n\in\fc$. 
Soit $\fc'$ un \itf de $\gA$ avec $\fc\fc'=a^n\gA$.
 \\
En utilisant le \tho de \KRO (reformulé dans le lemme \ref{lemthKroicl}),
on obtient 
$\fc\gB\subseteq\fa \fb \subseteq\Icl_\gB(\fc)$ et donc

\snic{a^n\gB=(\fc\gB) (\fc'\gB)\subseteq \fa \fb (\fc'\gB)\subseteq\Icl_\gB(\fc)(\fc'\gB)\subseteq\Icl_\gB(\fc\fc')
=\Icl_\gB(a^n)=a^n\gB.
}

%\sni
Donc $\fa \fb(\fc'\gB)=a^n\gB$ et $\fa$ est \lop d'après le lemme \ref{lemleszlop}.

Passons au cas \gnlz, avec $\alpha,\beta\in\gB$. 
Si $\gB$ est intègre, on peut supposer que~$\alpha\neq0$ et l'on trouve 
$\gamma\neq0$ dans $\gB$ tel que $\alpha\gamma=a\in\gA$, ce qui nous ramène au \pb déjà traité. 
\perso{se pose la question de dépenser de l'espace pour la version \gnlez,
renvoyer en exo la fin?}

Il reste à voir le cas, plus délicat, où l'on ne suppose pas $\gB$
intègre.  En fait on applique avec persévérance la recette des \lons
\come fournies par le \crc \lsdz de l'anneau, et cela marche.  On
écrit $p(\alpha)=0$ avec $p$ \mon dans $\AX$. On fait une \recu sur
$m=\deg(p)$. On a déjà traité le cas $m=0$. Passons de $m$ à~$m+1$.
On écrit $p(X)=Xq(X)+a$ avec $q$ \monz, $\deg(q)=m$ et l'on \linebreak 
pose $\wi\alpha =
q(\alpha)$. Puisque $\alpha\wi\alpha = -a \in \gA$, on sait trouver
$u,v\in\gB$ pour lesquels on~a~$\gen{u,v}\gen{\alpha\wi\alpha,\beta\wi\alpha} =
\gen{\alpha\wi\alpha}$, avec

\snic {
\alpha\wi\alpha = u\alpha\wi\alpha + v\beta\wi\alpha, \quad
u\beta\wi\alpha = u_1\alpha\wi\alpha, \quad
v\beta\wi\alpha = v_1\alpha\wi\alpha
.}

%\sni
Si l'on pouvait simplifier par $\wi\alpha$, on aurait
$\gen {u,v}\,\gen{\alpha, \beta} = \gen {\alpha}$. Les trois \egts
ci-dessus s'écrivent:

\snic {
(\alpha - u\alpha - v\beta)\,\wi\alpha = 0, \quad
(u\beta - u_1\alpha)\,\wi\alpha = 0 \quad
(v\beta - v_1\alpha)\,\wi\alpha = 0
.}

%\sni
\`A partir de ces \egtsz, on va trouver des \lons \come (comme dans la preuve
du lemme \ref{lemleszlop}).  Dans les unes $\wi\alpha = 0$, i.e. $q(\alpha)=0$
et l'on applique l'\hdrz.  Dans la dernière on a
$\gen{u,v}\gen{\alpha,\beta}=\gen{\alpha}$, ce qui montre que
$\gen{\alpha,\beta}$ est \lopz.
\end{proof}

\rem Cette \demz, comme d'ailleurs celle du lemme \ref{lemleszlop},
 est plus redoutable qu'il n'y paraît. Elle arrive à traiter
d'une seule manière le cas où $\alpha=0$, le cas où $\alpha$ est \ndzz, et \gui{tous les cas intermédiaires.}
\eoe

\medskip 
On a aussi le résultat facile suivant.

%:     Theorem{thSurAdp}
\begin{theorem}\label{thSurAdp}
Soit $\gA\subseteq\gB\subseteq\Frac\gA$. 
\begin{enumerate}
\item Si $\gA$ est \lsdzz, il en va de même pour $\gB$.
\item Si $\gA$ est \ariz, il en va de même pour $\gB$.
\item Si $\gA$ est de Prüfer, il en va de même pour $\gB$.
\end{enumerate}
\end{theorem}
\begin{proof}
Le point \emph{1} est laissé \alecz. 
\\
\emph{2.} 
Soient $x$, $y \in \gB$. Il existe $d \in \Reg(\gA)$  tel que  $x_1 = dx$, et $y_1 = dy$ sont dans $\gA$. Alors $d (x,y) = (x_1,y_1)$, et une \mlp dans $\gA$ pour $(x_1,y_1)$ est aussi une \mlp
pour~$(x, y)$.
%Montrons que si $a,b\in\gA$ avec $b$ \ndzz, alors l'anneau $\gA[z]$, avec $z=a/b$, est \ariz.
%Considérons un \id $\gen{p(z),q(z)}$, si $m=\deg(p)\geq n=\deg(q)$, on
%a $b^{m}\gen{p(z),q(z)}=\gen{p_1(a,b),q_1(a,b)}$, où $p_1=b^mp$ et $q_1=b^mq$
%sont \hmgs de degré $m$ à \coes dans $\gA$. Alors une \mlp dans $\gA$
%pour $p_1(a,b),q_1(a,b)$ est aussi une \mlp pour $p(z),q(z)$.
\end{proof}

Les deux \thos précédents sont reliés à deux résultats
classiques dans le cadre \noe  (cf.~\cite[page~17]{FJ}):

\mni
{\bf \Tho de Krull-Akizuki.} 
\emph{Si $\gA$ est un \adk et $\gL$ une extension
finie du corps de fractions de $\gA$, alors la \cli de~$\gA$
dans $\gL$ est un \adkz.}

\mni\rdb
{\bf \Tho de Grell-Noether.}\label{thGrellNoether} 
\emph{Si $\gA$ est un \adkz, alors tout anneau
compris entre $\gA$ et son corps de fractions est de Dedekind.}

\medskip 
Vue la \carn des \adks (en \clamaz) comme \adps \noes intègres, on voit que nous
avons établi les versions non-\noees et non-intègres de ces deux \thosz.

Nous démontrerons plus loin que dans les circonstances analogues, la \ddk de 
$\gB$ est toujours inférieure ou égale à celle de $\gA$, ce qui cette
fois-ci est lié à la \carn des \adks comme anneaux \icl
de dimension $\leq 1$ et \noesz.   
%%%%%%%%%%%%%%%%%%%%%%%%%%%%%%%%%%%%%%%%%%%%%%%%%%%%%%%%%%%%%%%%%%%%%%%%%%%
%%%%%%%%%%%%%%%%%%%%%%%%%%%%%%%%%%%%%%%%%%%%%%%%%%%%%%%%%%%%%%%%%%%%%%%%%%%
%%%%%%%%%%%%%%%%%%%%%%%%%%%%%%%%%%%%%%%%%%%%%%%%%%%%%%%%%%%%%%%%%%%%%%%%%%%
%%%%%%%%%%%%%%%%%%%%%%%%%%%%%%%%%%%%%%%%%%%%%%%%%%%%%%%%%%%%%%%%%%%%%%%%%%%

%%%%%%%%%%%%%%%%%%%%%%%%%%%%%%%%%%%%%%%%%%%%%%%%%%%%%%%%%%%%%%%%%%%%%%%%%%%
%--- Section{Anneaux de Prufer coh}  secAdPcoh   
\section{Anneaux de Prüfer \cohsz} 
\label{secAdPcoh}
%-------------------------

\vspace{3pt}
%%%%%%%%%%%%%%%%%%%%%%%%%%%%%%%%%%%%%%%%%%%%%%%%%%%%%%%%%%%%%%%%%%%%%%%%%%%
%: subsec{Premières \prts}%%%
\subsec{Premières \prtsz}

Rappelons que sur un anneau \qi un \itf est fidèle \ssi il contient un \elt \ndz
(voir le corolaire \ref{corlemQI}). En fait tout \itf contient un \elt qui
a le même annulateur que lui.
En particulier,  sur un anneau \qi un \itf \pro est \iv \ssi il est fidèle.

\smallskip 
Après avoir fourni des \carns des \adps (voir la proposition et
\dfnz~\ref{prop.itfplat 2} et le \thref{thPruf}), en voici pour les \adpcs;
\llec en trouvera d'autres dans l'exercice~\ref{exocaracPruCoh}.

%:     Th{th.adpcoh} ----------------
\begin{theorem}
\label{th.adpcoh}\relax \emph{(\Carns des \adpcsz)}\\
Pour un anneau $\gA$, \propeq
%-----------------begin item------------------
\begin{enumerate}
\item \label{i1th.adpcoh}\relax $\gA$ est un \adpcz.
\item  $\gA$ est un \anar \qiz.
\item  Tout \itf  est \proz.
\item  Tout idéal à deux \gtrs est \proz.
\item  $\gA$ est \qi et tout \id 
$\gen{a,b}$ avec $a\in\Reg\gA$ est \ivz.
\item \label{i6th.adpcoh}\relax  $\gA$ est \qi et tout \itf fidèle est un \mrcz~1.
\end{enumerate}
%-----------------end item------------------
\end{theorem}
\begin{proof} \emph {1} $\Leftrightarrow$ \emph {2.}  Utiliser le fait \ref{factLsdzCo}.
\\
\emph{3} $\Rightarrow$ \emph{4.} Trivial.
\\
\emph{4} $\Rightarrow$ \emph{2.}
Le \thref{thAnar} donne l'implication pour le caractère \lop des \idsz.
Par ailleurs un anneau est \qi \ssi les \idps sont projectifs.
\\
Les implications \emph {1} $\Rightarrow$ \emph {3},  \emph {5},  \emph {6} tiennent à la \carn des \ids \pros comme \ids \lops
dont l'annulateur est un \idm et celle des \ids \ivs comme \ids \lops contenant un \elt
\ndz (lemme~\ref{lemIdproj}, points \emph{2} et \emph{6}.).
\\
Pour les réciproques, on se rappelle qu'un \idp est \pro \ssi son annulateur est 
engendré par un \idm (lemme~\ref{lemIdpPtf}), 
et l'on peut voir le corrigé de l'exercice \ref{exocaracPruCoh}.
On peut aussi examiner ces réciproques 
dans le cas intègre, où elles sont claires,
et utiliser la machinerie \lgbe des anneaux \qisz.\imlgz
\end{proof}

Dans le cas local on obtient le résultat suivant (trivial en \clamaz, mais significatif d'un point de vue \cofz).

%:     Fact{factValCoh}
\begin{fact}\label{factValCoh}
Un \adv est \coh \ssi il est intègre. 
\end{fact}
\begin{proof}
Un \adp est \coh \ssi il est \qiz. Un \alo est connexe. Un anneau connexe
est intègre \ssi il est \qiz. 
\end{proof}

Dans ce cas $\gK=\Frac\gA$ est un \cdi et pour tout $x\in\gK\eti$, $x$ ou~$1/x$ est dans $\gA$. De manière \gnlez, on appelle \emph{\adv d'un \cdi $\gK$}
un sous-anneau vérifiant la \prt précédente. Et il est clair que c'est un \adv intègre.\index{valuation!anneau de --- d'un corps discret}\index{anneau!de valuation}

%:     factAdpc  
Les \prts de stabilité suivantes sont faciles.

\begin{fact}\label{factAdpc}  ~
\begin{enumerate}
\item %1
Un anneau \zed réduit est un \adpcz.

\item %2
Un localisé, un quotient réduit d'un \adpc par un \itf
est un \adpcz. 

\item %3
Un anneau est de Prüfer et \coh \ssi il a la même
\prt après \lon en des \mocoz.

\end{enumerate}
\end{fact}

%:     factAdpcdisc 
Un simple rappel ci-après: le point \emph{1} est valable pour les anneaux \qis et le point \emph{2} pour les \anarsz. 

\begin{fact}\label{factAdpcdisc}  Soit $\gA$ un \adpcz.
\begin{enumerate}
\item %1
$\gA$ est discret \ssi $\BB(\gA)$ est discret.
\item %2
$\gA$ est \fdi \ssi il est  \dveez. 
\end{enumerate}
\end{fact}

%%%%%%%%%%%%%%%%%%%%%%%%%%%%%%%%%%%%%%%%%%%%%%%%%%%%%%%%%%%%%%%%%%%%%%%%%%%
%: subsec{Noyau, image et conoyau d'une matrice}%%%
\subsec{Noyau, image et conoyau d'une matrice}

%:     theorem {ThImMat}
\begin{theorem}
\label{ThImMat} Soit $\gA $ un \adpcz. 
%-----------------begin enum------------------
\begin{enumerate}
\item %1
L'image d'une matrice $F\in\gA^{n\times m}$ est isomorphe à une somme directe de $n$ \itfsz.
\item %2
Tout sous-\mtf d'un \mptf est un \mptfz.
\item %3
Le noyau d'une \ali entre \mptfs est facteur direct (donc \ptfz).
\item %4
Tout \mpf est somme directe de son sous-module de
torsion (qui est \pfz) et d'un sous-\mptfz.
\item %5
Tout module \pro de rang $k\geq0$ est  isomorphe à une somme directe de
$k$ \ids \ivsz.
\item %6
Tout module \pro de rang  $\leq k$ est 
isomorphe à une somme directe de $k$ \itfsz.
\end{enumerate}
%-----------------end enum------------------
\end{theorem}
%--- end-theorem-----------------------------------------
NB: on ne demande pas que $\gA$ soit discret.
%-----------------begin proof------------------
\begin{proof}
On considère une \ali arbitraire
$\varphi:\Ae m\rightarrow \Ae n$.

 \emph{1.} 
 On traite le cas du module
 $M=\Im\,\varphi \subseteq\Ae n$. Soit
$\pi_n:\Ae n\rightarrow \gA$ la dernière
forme \cooz. L'\id $\pi_n(M)=\fa_n$ est \tf donc \proz,
et l'\ali surjective induite $\pi'_n:M\rightarrow \fa_n$ 
est scindée, et
$$
M\simeq \Ker \pi'_n\oplus \Im \pi'_n=(M\cap \Ae {n-1})\oplus \fa_n.
$$
On termine la preuve par \recu sur $n$:
$M\cap \Ae {n-1}$ est \tf puisque isomorphe à un quotient de $M$.
On obtient donc $M\simeq \fa_1\oplus \cdots \oplus \fa_n$.

 \emph{2.} Résulte \imdt de  \emph{1.} 

 \emph{3.}  Cela résulte de ce que l'image de l'\ali
est un \mptfz.

 \emph{4.} On traite le cas du module $N=\Coker \varphi$.

\emph{Voyons d'abord le cas où $\gA$ est local}, i.e., est un \adv
intègre. \\
La matrice de $\varphi$ se met en forme de Smith (proposition \ref{propPfVal}). Puisque l'anneau est intègre, $N$ est somme directe d'un module libre (correspondant aux \elts diagonaux nuls dans la réduite de Smith)
et d'un sous-module de torsion, lui-même somme directe de sous-modules 
$\aqo{\gA}{d_i}$ correspondant aux \elts diagonaux \ndzsz.

\emph{Voyons ensuite le cas où $\gA$ est intègre}.\\
Au moyen d'un nombre fini de \lons en des \ecoz, disons~$s_1$, $\ldots$, $s_r$, on  se ramène
à la situation du cas local (réduction de Smith de la matrice).
Puisque $\Ann_\gA(s_i)= \gen{0}$ ou $\gen{1}$, et puisque les \lons en $0$ sont inutiles, on peut supposer que les $s_i $ sont dans~$\Reg(\gA)$.\\ 
Notons $T$ le sous-module de torsion de $N$ et regardons ce qui se passe après
\lon en $S_i=s_i^{\NN}$. On constate facilement que le sous-module de torsion de $N_{S_i}$ est égal à $T_{S_i}$. 
Ainsi, $T$ est \pf parce qu'il est \pf après \lon en les $S_i$.
Il est facteur direct dans $N$ parce que $T_{S_i}$ est facteur direct
dans $N_{S_i}$ pour chaque $i$: l'injection canonique $T\to N$ admet un
inverse à gauche d'après le \plgrf{plcc.scinde}. 
Enfin, le module $N/T$, qui est \ptf après \lon en les $S_i$, est bien \ptfz.
\\  
Nous obtenons donc ce que nous souhaitions, avec un petit plus: le module~$T$
devient, après \lon en chacun des \elts $s_j$ d'un \sys comaximal $(s_1,\ldots, s_r)$, 
une somme directe de modules de torsion monogènes, i.e.
isomorphes à $\aqo{\gA[1/s_j]}{u_{k,j}}$, avec $u_{k,j}\in\Reg(\gA)$.

\emph{Voyons enfin le cas \gnlz, où $\gA$ est \qiz}.
\\
En partant de la \dem du cas intègre, la machinerie \lgbe \elr des anneaux \qis produit un \sfio $(e_1,\ldots,e_r)$ 
tel que le résultat soit acquis dans chacune des composantes
$e_iN$ (vue comme $\gA[1/e_i]$-module). Et cela donne \imdt le résultat 
global.\imlg

 \emph{5.} 
Dans le cas où $\gA$ est intègre, cela résulte du point~\emph{1}
puisque chaque \id dans la \dcn en somme directe est de rang $0$ ou $1$.
\\
On peut déduire le cas \gnl par la machinerie \lgbe \elrz.
Voici une autre \demz
\footnote{Plus savante ou moins savante, c'est difficile à dire. Cela dépend des go\^uts.}, indépendante de la \dem du point~\emph{1.}
Si $M$ est de rang constant $k \geq 1$, alors son dual $M\sta$
l'est \egmtz, leurs annulateurs sont nuls,
et il existe $\mu \in M\sta$ tel que $\Ann(\mu) = \gen{0}$
(voir le lemme \ref{lemQI}). Alors $\mu(M)$ est un \id \iv de $\gA $
car son annulateur est \egmt nul. De plus,
$M \simeq \Ker \mu \oplus \Im \mu$, ce qui prouve que $\Ker\mu$
est
\ptf de rang constant $k-1$. On termine par \recuz.

 \emph{6.} On considère $M$ comme somme directe de ses composantes de
rang constant, et l'on applique le point \emph{4} à chacune d'elles.
\end{proof}

%%%%%%%%%%%%%%%%%%%%%%%%%%%%%%%%%%%%%%%%%%%%%%%%%%%%%%%%%%%%%%%%%%%%%%%%%%%
%:   subsec{Extensions d'\adpcsz}
\subsec{Extensions d'\adpcsz}%
\index{algebrique@\agq!element prim@élément primitivement --- sur un anneau}\index{primitivement algébrique}
\label{subsecExtAdpC}

Un \elt $x$ d'une \Alg $\gB$ est dit \emph{primitivement \agq sur $\gA$} s'il annule un \pol
primitif de $\AX$. Après changement d'anneau de base, un \elt primitivement \agq reste primitivement
\agqz. La \prt pour un \elt d'être primitivement \agq est locale au sens suivant.

%:    Principe local global concret{plcc.agq}---
\begin{plcc} 
\label{plcc.agq}\relax 
\emph{(\'Eléments primitivement \agqsz)} 
Soient $S_1$, $\ldots$, $S_n$  des \moco d'un anneau $\gA$, $\gB$ une \Alg
et $x\in \gB$.   
Alors $x$ est primitivement  \agq sur $\gA$  \ssi 
il est primitivement \agq sur chacun des $\gA_{S_i}$.
\end{plcc}
%--- end-plcc-----------------------------------------
%
\begin{proof}
Il faut montrer que la condition est suffisante. On a des \eco $s_1$, $\ldots$, $s_n$
($s_i\in S_i$) et des \pols $f_i\in\AX$ tels que $s_i\in\rc(f_i)$ et~$f_i(x)=0$.
Si $d_i\geq\deg_X(f_i)+1$, on considère le \pol 

\snic{f=f_1+X^{d_1}f_2+X^{d_1+d_2}f_3+\cdots.}

%\sni
On a alors $f(x)=0$ et $\rc(f)=\sum_{i=1}^n\rc(f_i)=\gen{1}$.
\end{proof}
%

%:     Lemma{lemEmmanuel}
\begin{lemma}\label{lemEmmanuel} \emph{(Les entiers d'Emmanuel)} 
Soient $\gB$ un anneau et $\gA$ un sous-anneau, soient $\gA'$ la \cli de $\gA$
dans $\gB$ et $s$ un \elt de $\gB$ qui annule un \pol $f(X)=\sum_{k=0}^na_kX^k\in\AX$.\\
On note $g(X)=\sum_{k=1}^nb_kX^{k-1}$ le \pol $f(X)/(X-s)$. 
\begin{enumerate}
\item Les \elts $b_i$ et $b_is$ sont dans $\gA'$.
\item Dans $\gA'$ on obtient:
\vspace{-1mm}

\snic{
\gen{a_0,\ldots,a_n}=\rc(f)\subseteq\rc(g)+\rc(sg)=\gen{b_1,\ldots,b_n,b_1s,\ldots,b_ns}.}
\item Dans $\gA'[s]$ les deux \ids sont égaux. 
\end{enumerate}
\end{lemma}
\begin{proof}
Puisque $f(X)=(X-s)g(X)$, le \tho de \KRO nous dit que les~$b_i$ et $b_is$
sont entiers sur $\gA$. On a 

\snic{b_n=a_n$, $b_{n-1}=b_ns+a_{n-1}$, \ldots,
$b_{1}=b_2s+a_{1}$,
$0=b_1s+a_{0}.}

%\sni
Donc chaque $a_i\in\rc(g)+\rc(sg)$. Et, dans $\gA'[s]$, de proche en proche, on obtient
$b_n\in\rc(f)$, $b_{n-1}\in\rc(f)$, \ldots, $b_{1}\in\rc(f)$.
\end{proof}
%

%:     Theorem{th.2adpcoh}
\begin{theorem}\label{th.2adpcoh}\label{i7th.adpcoh}\relax 
\emph{(Une autre \carn des \adpcsz, voir aussi les exercices \ref{exocaracPruferC}
et \ref{exocaracPruCoh})}\\
Un anneau $\gA$ est un \adpc \ssi il est \qiz,
\icl dans $\Frac\gA$, et si tout \elt de $\Frac\gA$ est primitivement \agq sur~$\gA$. 
\end{theorem}
\begin{proof} 
Supposons que $\gA$ est un \adpcz. Il nous reste à montrer que tout \elt
de $\Frac\gA$ est primitivement \agq sur $\gA$.
Soit $x=a/b\in\Frac\gA$. 
Il y a une matrice $\cmatrix {s & u\cr v & t\cr} \in \MM_2(\gA)$, 
de \lon principale pour $(b,a)$, i.e. $s + t = 1$,
$sa = ub$ et $va = tb$.
\\
Ce qui donne
$sx-u=0$ et  $t=vx$. Ainsi, $x$ annule le \pol
pri-\linebreak
mitif~$-u+sX+X^2(t-vX)$, ou si l'on préfère $t-(u+v)X+sX^2$.

Voyons la réciproque.
Il suffit de faire la preuve dans le cas intègre.
On doit montrer que tout \id $\gen{a,b}$ est \lopz.
On suppose \spdg $a$, $b\in\Reg(\gA)$. L'\elt $s=a/b$ annule un \pol 
primitif $f(X)$. Puisque $\rc(f)=\gen{1}$ dans $\gA$, d'après le lemme
\ref{lemEmmanuel} (points~\emph{1} et~\emph{2}), on a des \elts $b_1$, $\ldots$, $b_n$, $b_1s$, $\ldots$, $b_ns$ \com dans
$\gA$. \\
On a alors $s\in\gA[1/b_i]$ et $1/s\in\gA[1/(b_is)]$: dans chacune des
\lons \comez, $a$ divise $b$ ou $b$ divise~$a.$
\end{proof}

Le \tho qui suit contient une nouvelle \dem de la stabilité des \adps intègres
par extension
entière et intégralement close (voir le \thref{thExtEntPruf}).
Elle semble d'une facilité déconcertante par rapport à celle donnée 
sans l'hypothèse de \cohcz.

%:     Theorem{factAdpIntExt}
\begin{theorem}\label{factAdpIntExt}
Si $\gB$ est un anneau \qi normal, et une extension entière
d'un \adpc $\gA$, alors $\gB$ est un \adpcz.
\end{theorem}
%%%%%%%%%%%%%%%%%%%%%%%%%%%%%%%%%%%%%%%%%
%
\begin{proof} Voyons d'abord le cas où \emph{$\gB$ est intègre et non trivial.}
Soit $s\in\Frac\gB$. Il suffit de montrer que $s$ est primitivement \agq sur $\gB$.
On a un \pol non nul $f(X)\in\AX$ tel que $f(s)=0$.
\\
\emph{Cas où $\gA$ est un anneau de Bézout.} On divise $f$ par $\rc(f)$
et l'on obtient un \pol primitif qui annule $s$.
\\
\emph{Cas d'un \adpz.} Après \lon en des \ecoz, l'\id $\rc(f)$ est engendré 
par un des \coes de $f$, le premier cas s'applique.

 Dans le cas \gnlz,  la machinerie
\lgbe \elr des anneaux \qis  nous ramène au cas intègre.\imlgz 
\end{proof}

Voici maintenant l'analogue de   la proposition \ref{propAECDN},
qui décrivait l'anneau d'entiers d'un corps de nombres.  Dans le cas où $\gA$ est un anneau de Bézout intègre,
on aurait pu reprendre presque mot pour mot les mêmes \demsz.
Notez aussi que le \thref{Thextent} étudie une situation du même style avec une hypothèse un peu plus faible.
Voir aussi le point \emph{1} du \pbz~\ref{exoLemmeFourchette}.

%:     Theorem{thAESTE}
\begin{theorem}\label{thAESTE}%\label{thPfpEE}
 \emph{(Anneau d'entiers dans une extension \agqz)}\\
Soit $\gA$ un \adpcz, 
$\gK=\Frac(\gA)$, $\gL\supseteq\gK$ une \Klg entière réduite et $\gB$ la \cli de $\gA$ dans $\gL$.
\begin{enumerate}
\item \label{i1thAESTE} $\Frac\gB=\gL=(\Reg\gA)^{-1}\gB$ et  $\gB$ est un \adpcz.
\item  \label{i2thAESTE} Si $\gL$ est \stfe sur $\gK$ et si  $\gA$ est \fdiz, $\gB$ est \fdiz.
\end{enumerate}
Si en outre $\gL$ est étale sur $\gK$, on obtient:
\begin{enumerate} \setcounter{enumi}{2}
\item  \label{i3thAESTE} Si  $\gA$ est \noez, il en va de même pour $\gB$.
\item  \label{i4thAESTE} Si  $\gA$ est un \adk (\dfn \ref{defDDK}), $\gB$ \egmtz.
\item  \label{i5thAESTE} Si $\gL=\Kx=\aqo\KX f$ avec $f\in\AX$ \mon et 
$\disc_X(f) \in\alb\Reg\gA$, %\linebreak 
alors $\,\fraC1\Delta\gA[x]\subseteq \gB\subseteq\gA[x]$ ($\Delta=\disc_X(f)$). \\
En particulier $\gA[x][\fraC1\Delta]=\gB[\fraC1\Delta].$
\item  \label{i6thAESTE} Si en outre $\disc_X(f)\in\Ati$, on a $\gB=\gA[x]$
 \ste sur $\gA$.
\end{enumerate}
 
\end{theorem}
%--------- fin theorem ---------------------------------------------- 
%
\begin{proof}
\emph{\ref{i1thAESTE}.} Conséquence directe du fait \ref{factReduitEntierQi} et du \thrf{factAdpIntExt}.

 \emph{\ref{i2thAESTE}.}
Puisque $\gB$ est un \adpz, il suffit de savoir tester la \dve dans $\gB$, \cad
l'appartenance d'un \elt de $\gL$ à $\gB$.  Soit~\hbox{$y \in \gL$} et~$Q
\in \KY$ son \polmin (\monz) sur $\gK$. Alors $y$ est entier sur~$\gA$
\ssi $Q \in \AY$: dans le sens non \imdz, \linebreak 
soit $P \in \AY$ \mon tel que
$P(y) = 0$, alors $Q$ divise $P$ dans $\KY$ et le lemme \ref{lem0IntClos}
implique que $Q\in\AY$. 
\\
Note: on aurait aussi bien pu utiliser le \polcarz, 
mais la \dem qui utilise le \polmin fonctionne dans un cadre plus \gnl
(il suffit que $\gL$ soit \agq sur $\gK$ et que l'on sache calculer les \polminsz).

 \emph{\ref{i5thAESTE}.} 
Dans le cas où $\gA$  est un anneau de Bézout intègre et $\gL$ un corps, on applique le \thref{Thextent}. 
Le résultat dans le cas \gnl est ensuite obtenu à partir de cette \dem en utilisant les machineries \lgbes des anneaux \qis et des \anarsz.\imlgz\imla
%\\
%Soit $z={h(x) \over \delta}$ un \elt de $\gB$, écrit avec $\delta\in\Reg\gA$, $\gen{\delta}+\rc(h)=\gen{1}$ \linebreak 
%et~$\deg_X(g)<n=\deg_X(f)$.
%On va montrer,  résultat plus précis, que~$\delta^2$ divise~$\disc_X(f)$.
%\\
%Le module $M=\gA[x]+z\gA\subseteq \gB$ est un sous-\Amo \tf de~$\frac 1 \delta\gA[x]$,
%et~$\frac 1 \delta\gA[x]$ est libre de rang $n$. Donc $M$ est un \Amo \ptfz, et puisque 
%$\gA$ est de Bézout intègre,
% $M$ est libre de rang $n$ sur~$\gA$.
% \\
%  Si $d$ est le \deter d'une matrice
%qui exprime la base naturelle $\cB_0$ de~$\gA[x]$ sur une base  $\cB_1$ de $M$,
%on obtient: 
%
%\snic{d^2\disc_{\gK[x]/\gK}(\cB_1)=\disc_{\gK[x]/\gK}(\cB_0)=\disc_X(f).}
%
%%\sni
%(Propositions \ref{defiDiscTra} et \ref{propdiscTra}.) Enfin  $\gen{d}=\gen{\delta}$ d'après le lem\-me~\ref{lemSousLibre}.   

 \emph{\ref{i3thAESTE}.} On fait la \dem sous les hypothèses du point \emph{5}. Ce n'est pas restrictif car d'après le \thref{thEtalePrimitif}, $\gL$ est un produit \hbox{de \Klgsz} étales monogènes.
% : il suffit donc de montrer le résultat lorsque~$\gL$ est  du type $\gK[y]=\aqo\KY F$, où $F$ est un \polu \spl de $\KY$, mais on a aussi
%pour un $a\in\Reg\gA$, $x=ay\in\gB$   
%et $\gL=\gK[x]=\aqo\KX f$, où~$f\in\AX$ est  \mon   avec $\disc_X(f)\in\Reg\gA$.
\\
Soit $\fb_1\subseteq \fb_2 \subseteq\cdots\subseteq \fb_n\subseteq \dots$
une suite d'\itfs de $\gB$ que l'on écrit $\fb_n=\gen{G_{n}}_\gB$ avec $G_n
\subseteq G_{n+1}$; on définit

\snic {
L_n = \disc_X(f) \cdot \left(\sum_{g \in G_n} \gA g\right) \subseteq \gA[x]
.}

%\sni
Alors $L_1\subseteq L_2\subseteq \cdots\subseteq L_n\subseteq \dots$ est une
suite de sous-\Amos \tf de $\gA[x]$. 
Or $\gA[x]$ est un \Amo libre de rang fini (égal à $\deg(f)$), donc \noez. 
On termine en notant que si $L_m=L_{m+1}$, alors $\fb_m=\fb_{m+1}$.

\emph{\ref{i4thAESTE}.} Résulte de \emph{\ref{i2thAESTE}} et \emph{\ref{i3thAESTE}.}

\emph{\ref{i6thAESTE}.}  Il est clair que $\gB=\gA[x]$.
\end{proof}
%

%-% ENTRE NOUS
\entrenous{Il semble que pour les points \emph{\ref{i3thAESTE}} et \emph{\ref{i4thAESTE}} l'hypothèse \ste pourrait être affaiblie en \stfez, mais

1) ce n'est pas certain, mais il y a quelque chose de ce style obtenu gr\^ace à la théorie des diviseurs

2) je ne sais pas faire sans la théorie des diviseurs 
}
%-% Fin ENTRENOUS

\rem
 Le \tho précédent s'applique dans deux cas importants dans l'histoire de l'\alg commutative. 
 \\
 Le premier cas est celui des anneaux d'entiers de corps de nombres, \linebreak 
avec $\gA=\ZZ$ et $\gB$ l'anneau d'entiers d'un corps de nombres (cas déjà examiné en section~\ref{secApTDN}).
 \\
  Le deuxième cas est celui des courbes \agqsz. On considère un \cdi
$\gk$, l'anneau principal $\gA=\gk[x]$ et un \pol $f(x,Y)\in\gk[x,Y]$ \mon en $Y$, \irdz, avec $\disc_Y(f)\neq0$. On note $\gK=\gk(x)$. 
\\
L'anneau $\gA[y]=\gk[x,y]=\aqo{\gk[x,Y]}f$ est intègre. 
La courbe plane $\cC$ d'équation $f(x,Y)=0$ peut avoir des points singuliers, auquel cas $\gA[y]$ n'est pas \ariz.
Mais la \cli $\gB$ de $\gA$ dans $\gK[y]=\aqo{\gK[Y]}f$ est bien
un \ddp (\thref{thcohdim1}), en fait un \dDkz. Le corps $\gK[y]$ est appelé le corps de fonctions de $\cC$.  L'anneau $\gB$ correspond à une courbe (qui n'est plus \ncrt plane) sans point singulier,
avec le même corps de fonctions que $\cC$.   
\eoe

%%%%%%%%%%%%%%%%%%%%%%%%%%%%%%%%%%%%%%%%%%%%%%%%%%%%%%%%%%%%%%%%%%%%%%%%%%%
%%%%%%%%%%%%%%%%%%%%%%%%%%%%%%%%%%%%%%%%%%%%%%%%%%%%%%%%%%%%%%%%%%%%%%%%%%%
%%%%%%%%%%%%%%%%%%%%%%%%%%%%%%%%%%%%%%%%%%%%%%%%%%%%%%%%%%%%%%%%%%%%%%%%%%%
%%%%%%%%%%%%%%%%%%%%%%%%%%%%%%%%%%%%%%%%%%%%%%%%%%%%%%%%%%%%%%%%%%%%%%%%%%%
%%%%%%%%%%%%%%%%%%%%%%%%%%%%%%%%%%%%%%%%%%%%%%%%%%%%%%%%%%%%%%%%%%%%%%%%%%%
%--- Section{Anneaux de Prufer coh}  secAdPcohDim1   
\section{Anneaux \qis de dimension \texorpdfstring{$\leq1$}{<=1}}
\label{subsecQiDim1}
%-------------------------

La plupart des \thos \gui{classiques} concernant les \dDks sont déjà
valables pour les \adpcs \ddi1, voire pour les \anarsz. 
Nous en démontrons un certain nombre dans cette section et la suivante.

Dans cette section les résultats concernent les anneaux \qis \ddi1.

Le \tho suivant  est un cas particulier du \gui{stable range} de Bass
dont nous donnerons des versions \gnles (\thosz~\ref{Bass0} et~\ref{Bass}).

%:     Th{thK1-SLE}----------
\begin{theorem}
\label{thK1-SLE}\relax
Soit $n\geq 3$ et $\tra[\,x_1 \;\cdots\;x_n \,]$ un \vmd sur un
anneau \qi $\gA$ \ddi1. Ce vecteur est la première colonne d'une matrice
de
   $\En(\gA)$. En particulier, $\SLn(\gA)$ est engendré par $\En(\gA)$
   et $\SL_2(\gA)$ pour $n\geq 3$.
Et pour $n\geq 2$ tout \vmd est la première colonne
d'une matrice de~$\SLn(\gA)$.
\end{theorem}
%--- end-theorem----------------------------
%-----------------begin proof---------------
\begin{proof}
L'annulateur de $\gen{\xn}$ est nul, donc
on peut par \mlrs transformer le vecteur 
$v= {\tra[\,x_1 \;\cdots\;x_n \,]}$ en un \vmd $ {\tra[\,y_1\;x_2 \;\cdots\;x_n \,]}$, avec~$y_1\in\Reg(\gA)$ (cf. lemme~\ref{lemQI}).
\\
Considérons 
l'anneau~$\gB=\aqo{\gA}{y_1}$. Cet
anneau est \zed et le vecteur $v$ 
devient égal à
${\tra[\,0\;x_2 \;\cdots\;x_n \,]}$ toujours \umdz.
\\
 Puisque $n\geq 3$,
on peut transformer dans $\gB$ par \mlrs ${\tra[\,x_2 \;\cdots\;x_n \,]}$ en ${\tra[\,1\;0 \;\cdots\;0 \,]}$ (exercise \ref{exoSLnEn}).
Regardons dans $\gA$ ce que l'on obtient alors:
${\tra[\,y_1\;1+ay_1 \;z_3\;\cdots\;z_n \,]}$,
d'où ensuite, toujours par \mlrs ${\tra[\,y_1\;1 \;z_3\;\cdots\;z_n \,]}$, puis ${\tra[\,1\;0 \;\cdots\;0 \,]}$.
\end{proof}
%-----------------end proof-----------------

%:     Theorem{th1-5}
Le \tho suivant \gns le résultat analogue déjà obtenu en
théorie des nombres (corolaire \ref{corpropZerdimLib}). Le  point \emph{1} concerne les \ids \ivsz.  Le  point \emph{2} s'applique à tous les \itfs d'un \adpc de dimension $\leq1$. Une \gnn est proposée dans le \thref{th1.5}.

\begin{theorem}\label{th1-5}\relax \emph{(\Tho un et demi)}%
%:HHH index
\index{un et demi!\Tho ---}
\\
Soit $\gA$ un anneau \qi de dimension $\leq 1$ et $\fa$ un \id \lop (donc \ptfz).
\begin{enumerate}
\item 
Si  $a\in\fa\cap\Reg(\gA)$, il existe $b\in\fa$ tel que $\fa=\gen{a^n,b}$ pour tout $n\geq1$. 
\item 
Il existe $a\in\fa$ tel que $\Ann(a)=\Ann(\fa)$. Pour un tel $a$ il existe $b\in\fa$ tel que $\fa=\gen{a^n,b}$ pour tout $n\geq1$. 
\end{enumerate}
\end{theorem}
\begin{proof}
La \dem du point \emph{1} est identique à celle du corolaire \ref{corpropZerdimLib}
qui donnait le résultat en théorie des nombres.\\
\emph{2.} Tout \itf $\fa$ contient un \elt $a$ tel que $\Ann(a)=\Ann(\fa)$
(corolaire \ref{corlemQI}). On passe au quotient $\aqo{\gA}{1-e}$ où $e$ est l'\idm tel que $\Ann(a)=\Ann(e)$ et l'on applique le point \emph{1.} 
\end{proof}
%

%:     proposition}\label{avant.dekinbe}\relax
\begin{proposition} \label{avant.dekinbe}\relax
Soit $\gA$ un anneau  \qi de dimension $\leq1$,
dont le radical de Jacobson 
contient un \elt \ndzz, 
et $\fa$ un \id \ivz. Alors $\fa$ est principal.
\end{proposition}

\begin{proof}
Soient $y \in \Rad(\gA)$
et $x \in \fa$ tous deux \ndzsz.
Alors $\fa\cap \Rad(\gA)$ contient $a=xy$ qui est \ndzz.
Par le \tho un et demi, il existe $z\in\fa$ tel que $\fa=\gen{a^2,z}$.
Donc $a=ua^2+vz$ ce qui donne $a(1-ua)=vz$ et puisque $a\in\Rad(\gA)$, $a\in\gen{z}$ donc $\fa=\gen{z}$.
\end{proof}

%--- Lemma{lemEvitCo}----------------
Nous revisitons maintenant le résultat classique suivant,
dans lequel nous allons nous débarrasser de l'hypothèse
\noeez:
{\em si $\gA$ est un anneau \noe intègre \ddi1 et $\fa$, $\fb$
   deux \ids avec $\fa$ \iv et $\fb\neq 0$, alors il existe $u\in\Frac(\gA)$
tel que  $u\,\fa\subseteq\gA$ \linebreak 
et $u\fa+\fb=\gen{1}$.}

%--- Lemma
\begin{lemma}
\label{lemEvitCo}
Soit $\gA$ un anneau \qi (par exemple un \adpcz)  \ddi1. Soit 
$\fa$ un \id \iv
de $\gA$ et $\fb$ un \id contenant un \elt \ndzz. 
Alors il existe un \elt $u$ \iv dans $\Frac(\gA)$
   tel que $u\fa\subseteq\gA$ et $u\fa+\fb=\gen{1}$.
\end{lemma}
%--- end-lemma-----------------------------------------

%-----------------begin proof------------------
\begin{proof}
Nous faisons la \dem dans le cas intègre, en laissant le soin \alec d'appliquer ensuite la machinerie \lgbe \elr des anneaux \qisz.
Pour lui faciliter la t\^ache, nous ne supposons pas $\gA$ non trivial et nous mettons \gui{\ndzz} lorsque dans le cas non trivial nous aurions mis 
\gui{non nul}.\imlg\\
On cherche $a$ et $b$ \ndzs tels que $\fraC{b}{a}\,\fa\subseteq\gA$,
\cad encore $b\,\fa\subseteq a\gA$, et~$\gA=\fraC{b}{a}\,\fa+\fb$. 
Si $c$ est un \elt \ndz de $\fb$, comme 
la condition devrait être aussi réalisée 
lorsque $\fb$ est l'idéal $c\gA$,
on doit trouver  $a$ et $b$ \ndzs tels que $b\,\fa\subseteq a\gA$ et
$\gA=\fraC{b}{a}\,\fa+c\gA$. 
Si l'on s'arrange pour que~$a\in\fa$, 
on aura~$b\in \fraC{b}{a}\fa$, 
%on a $b\in a\,\fa$ 
et il  suffit donc de réaliser les conditions $b\,\fa\subseteq a\gA$
\\ 
et~$\gA=\gen{b,c}$. C'est ce que nous allons faire. \\
Soit $c \in \fa\cap\fb$ un \elt \ndz (par exemple le produit
de deux \elts \ndzsz, l'un dans $\fa$ et l'autre dans $\fb$).
D'après le \tho un et demi, il existe un~$a\in\fa$ tel que 
$\fa=\gen{a,c^2}=\gen{a,c}$. 
Si $a=0$, l'\id $\fa=\gen{c}$ est \idm donc égal à $\gen{1}$,
et ce n'était donc pas la peine de se fatiguer\footnote{Notons cependant que nous ne sommes pas censés savoir d'avance si un \id
\iv de $\gA$ contient $1$, nous ne nous sommes donc pas fatigués complètement pour rien, le calcul nous a permis de savoir que $1\in\fa.$}: on pouvait choisir $b=a=1$. 
\\
On suppose donc $a$ \ndzz. 
Puisque $c\in\fa$, on a une \egt $c=\alpha a+\beta c^2$, ce qui donne
$c(1-\beta c)=\alpha a$. Posons $b=1-\beta c$ de sorte que 
$\gA=\gen{b,c}$. On obtient $b\,\fa=b\gen{a,c}=\gen{ba,bc}
=a\gen{b,\alpha}\subseteq a\gA$.
Si  $b$ est \ndzz, on a donc gagné, et si $b=0$,  alors $1\in\gen{c}$
et ce n'était
pas la peine de se fatiguer. 
\end{proof}
%-----------------end proof------------------

%%%%%%%%%%%%%%%%%%%%%%%%%%%%%%%%%%%%%%%%%
\begin {proposition}
\label{prop-a/ab} \relax
Soit $\fa$ un idéal \iv d'un anneau intègre $\gA$
de dimension~$\le 1$. Pour tout \id non nul $\fb$
de $\gA$, on a un \iso \hbox{de \Amosz} 
$\fa/\fa\fb \simeq \gA/\fb$.
\end {proposition}
%%%%%%%%%%%%%%%%%%%%%%%%%%%%%%%%%%%%%%%%%
\begin {proof}
D'après le lemme~\ref{lemEvitCo}, il existe un \id entier~$\fa'$ dans la
classe\footnote{Voir \paref{pageclassgroup}.} de~$\gA\div\fa$ tel \hbox{que $\fa' + \fb = \gA$}; \hbox{on a $\fa\fa' = x\gA$} \hbox{avec
$x \in \Reg\gA$}.  La
multiplication par~$x$, $\mu_x : \gA \to \gA$, induit un \iso 

\snic{\gA/\fb \simarrow x\gA/x\fb = \fa'\fa/\fa'\fa\fb.}

 Considérons maintenant
l'application canonique 

\snic{f : \fa'\fa \to \fa/\fa\fb}

\snii
qui à $y \in \fa'\fa\subseteq \fa$ associe la classe de $y$ modulo $\fa\fb$.
Montrons que $f$ est surjective:
en \hbox{effet, $\fa' + \fb = \gA \Rightarrow \fa'\fa + \fa\fb = \fa$},
 donc tout
\elt de $\fa$ est congru à un \elt de $\fa'\fa$
modulo $\fa\fb$.  Examinons enfin  $\Ker f=\fa'\fa \cap \fa\fb$. Puisque
$\fa$ est \ivz, $\fa'\fa \cap \fa\fb = \fa(\fa' \cap
\fb)$, et enfin $\fa' + \fb = \gA$ entraîne \hbox{que $\fa' \cap \fb = \fa'\fb$},
donc $\Ker f=\fa'\fa\fb$.  On a ainsi des \isos de
\Amos
$$\preskip-.4em \postskip.0em 
\gA/\fb \simeq x\gA/x\fb = \fa'\fa/\fa'\fa\fb \simeq \fa/\fa\fb, 
$$
d'où le résultat.
\end {proof}

%-% ENTRE NOUS
\entrenous{J'ai un \pb avec le résultat de la proposition \ref{prop-a/ab}.
Puisque $\fa$ est projectif il est plat, $\fa\fb\simeq\fa\otimes \fb$, et la suite exacte $0\to\fb\to\gA\to\gA\sur\fb\to0$ donne que $\fa\otimes \gA\sur\fb$ est isomorphe à $\fa\sur{\fa\fb}$, cela impliquerait donc (avec la proposition \ref{prop-a/ab}) que $\fa\otimes \gA\sur\fb\simeq \gA\sur\fb$ pour tout \id $\fb$ non nul de $\gA$. 

Peut-être en fait ce n'est pas si étonnant vu que $\gA\sur\fb$ est \zedz.
En effet, sur un anneau \zed un \id \pro de rang $1$ est toujours égal à $\gen{1}$!
Mézalor cela donne  une \dem alternative plus simple de la
proposition \ref{prop-a/ab}. 

Est-ce que la \dem du lemme \ref{lemEvitCo} ne pourrait pas elle-même être simplifiée???

Par ailleurs, le lemme \ref{lemEvitCo} labélisé {\tt lemEvitCo} a sans doute à voir avec le \thref{propEvitementConducteur} labélisé {\tt propEvitementConducteur}.

Il faudrait un commentaire sur ce sujet, et, qui sait?, une \gnn du 
\tho \ref{propEvitementConducteur}.
}
%-% Fin ENTRENOUS

%%%%%%%%%%%%%%%%%%%%%%%%%%%%%%%%%%%%%%%%%%%%%%%%%%%%%
\begin {corollary} \label{corprop-a/ab}\relax
Soient $\gA$ un anneau intègre avec $\Kdim\gA \le 1$, $\fa$ un \id
\iv et $\fb$ un \id non nul. On a alors une
suite exacte de \Amosz:
$$\preskip.2em \postskip.4em
0 \to \gA/\fb \to \gA/\fa\fb \to \gA/\fa \to 0
$$
\end {corollary}
%%%%%%%%%%%%%%%%%%%%%%%%%%%%%%%%%%%%%%%%%

%:     Lemma{lemRadJDIM1}
\begin{lemma} \label{lemRadJDIM1}
\emph{(Radical de Jacobson d'un anneau de dimension~$\leq 1$)}
\\ 
Soit $\gA$ un anneau intègre  de dimension $\leq 1$.  
\begin{enumerate}
\item Pour tout $a$ non nul dans $\gA$,  $\Rad(\gA)\subseteq \sqrt[\gA]{a\gA}$.
\item Pour $\fb$ \tf contenant $\Rad(\gA)$,  $\Rad(\gA)=\fb \big(\Rad(\gA):\fb\big)$.
\item Si $\Rad(\gA)$ est un \id \ivz, $\gA$ est un  domaine de Bézout.
\end{enumerate} 
\end{lemma}
%--------- fin lemma ---------------------------------------------- 
%
\begin{proof} On note~\hbox{$\fa=\Rad(\gA)$}.\\
\emph{1.}
Soit $x\in\fa$,  $\aqo\gA a$ est \zedz, 
donc il existe~\hbox{$y$, $z\in\gA$} \hbox{et $m\in\NN$} tels que $x^{m}(1+xz)=ay$. 
Comme  $x\in\Rad(\gA)$, on a $1+xz\in\Ati$, \hbox{donc $x^{m}\in a\gA$} \hbox{et
$x\in \sqrt[\gA]{a\gA}$}.

\emph{2.} Si $\fa=0$ c'est clair, sinon l'anneau $\gA/\fa$ est \zedrz, donc l'\itfz~$\fb$ est égal à un
\id $\gen{e}$ modulo $\fa$, avec $e$ \idm modulo $\fa$.
 Donc $\fb=\fb+\fa=\fa+\gen{e}$, puis $(\fa:\fb)=\fa+\gen{1-e}$, et enfin

\snic{\fb(\fa:\fb)=(\fa+\gen{e})(\fa+\gen{1-e})=\fa.}

\emph{3.} Soit $\fc_1$ un \itf non nul arbitraire.
On définit $\fb_1=\fc_1+\fa$ \hbox{et $\fc_2=(\fc_1:\fb_1)$}.  D'après le point \emph{2}, puisque $\fa$ est \ivz, $\fb_1$ \egmtz. 
Si $\fb_1\fb'=\gen{b}$ ($b$ \ndzz), tous les \elts de $\fc_1\fb'$ sont divisibles par $b$, on considère alors $\fd=\fraC 1 b \fc_1\fb'$, donc~\hbox{$\fd\fb_1=\fc_1$}
et~$\fd$ est \tfz. On a clairement $\fd\subseteq \fc_2$. Réciproquement
si $x\fb_1\subseteq \fc_1$ alors $bx=x\fb_1\fb'\subseteq b\fd$, \hbox{donc $x\in\fd$}. En bref $\fc_2=\fd$ et l'on a établi l'\egt $\fb_1\fc_2=\fc_1$, avec $\fc_2$ \tfz.
En itérant le processus on obtient une suite croissante d'\itfsz~$(\fc_k)_{k\in\NN}$ avec $\fc_{k+1}=(\fc_k:\fb_k)$ et $\fb_k=\fc_k+\fa$.
\\
En fait $\fc_2=\big(\fc_1:(\fc_1+\fa)\big)=(\fc_1:\fa)$, puis $\fc_3=(\fc_2:\fa)=(\fc_1:\fa^{2})$ et plus \gnltz~\hbox{$\fc_{k+1}=(\fc_1:\fa^{k})$}. \\
Soit $a\neq 0$ dans $\fc_1$. Par le point \emph{1}, $\fa\subseteq \sqrt{a\gA}$. Or $\fa$ est \tfz, donc l'inclusion~\hbox{$\fa\subseteq \sqrt{a\gA}$}  implique que pour un certain $k$, $\fa^{k}\subseteq a\gA\subseteq \fc_1$, \hbox{donc
$\fc_{k+1}=\gen{1}$}.
\\
Lorsque $\fc_{k+1}=\gen{1}$, on a $\fc_1=\prod_{i=1}^{k}\fb_i$, qui est \iv comme produit d'\ids \ivsz.
\\
On a montré que tout \itf non nul est \ivz, donc l'anneau est un \ddpz, et d'après la proposition~\ref{avant.dekinbe} c'est un anneau de Bézout.
\end{proof}

%--- Section{Anneaux de Prufer coh}  secAdPcohDim1   
\section[Anneaux de Prüfer \cohs de dimension \texorpdfstring{$\leq1$}{<=1}]{Anneaux de Prüfer \cohs de ~~~~~~~~~~~~~~~~dimension $\leq1$}
\label{secAdPcohDim1}
%-------------------------
%\label{subsecAdpCD1}

%\vspace{3pt}
\subsec{Quand un \adp est un anneau de Bézout}
%:     theorem{dekinbe}---------
Nous généralisons maintenant
un résultat classique souvent formulé ainsi\footnote{Voir la \dfn
\cov d'un anneau de Dedekind, page \pageref{defDDK}.}:
{\em un anneau de Dedekind
intègre ayant un nombre fini d'\idemas est un anneau principal.}

%     theorem \label{dekinbe} 
\begin{theorem}\label{dekinbe}
Soit $\gA$ un \adpc   \ddi1 et
dont le radical de Jacobson 
contient un \elt \ndzz. Alors $\gA$ est un anneau de Bézout.
\end{theorem}

\begin{proof}
Soit $\fb$ un \itfz. Il existe $b\in\fb$ tel que~\hbox{$\Ann\,\fb=\Ann\, b=\gen{e}$} avec $e$ \idmz. Alors $\fa=\fb+\gen{e}$ contient l'\elt \ndz $b+e$: il est \iv et $\fb=(1-e)\fa$.
Il suffit de montrer que  $\fa$
est principal. Or cela résulte de la proposition~\ref{avant.dekinbe}.
\end{proof}

Le \tho précédent et le suivant sont à comparer avec le 
\thref{propGCDDim1} qui affirme qu'un anneau intègre à pgcd de dimension
$\leq1$ est un anneau de Bézout.

%%%%%%%%%%%%%%%%%%%%%%%%%%%%%%%%%%%%%%%%%%%%%%%%%%%%%%%%%%%%%%%%%%%%%%%%%%%
\subsec{Une \carn importante}

%%%%%%%%%%%%%%%%%%       thcohdim1       %%%%%%%%%%%%%%%%%%%%%%%%%
 
Le résultat donné dans le \tho \ref{thcohdim1} ci-après est important: 
les trois machineries calculatoires de la normalité, de la \cohc
et de la dimension~1 se combinent pour fournir la 
machinerie de la \lon principale des \itfsz.

%:     Theorem{thcohdim1}
\begin{theorem}\label{thcohdim1}
Un anneau $\gA$, normal, \cohz, de dimension $\leq1$ est un \adpz. 
\end{theorem}
\begin{proof}
Commençons par remarquer que $(\gA\div\fa\fb) = (\gA\div\fa)\div\fb$.
\\
Puisque $\gA$ est \qiz, il suffit de traiter le cas intègre et de terminer
avec la machinerie \lgbe des anneaux \qisz.
Nous supposons donc que $\gA$ est intègre et nous montrons que tout
\itf $\fa$ contenant un \elt \ndz est \ivz.\imlg\\
Considérons $(\gA\div \fa) \in \Ifr\gA$
et $\fb=\fa (\gA\div \fa)$, qui est un \itf (entier) de $\gA$;
nous voulons montrer que $\fb=\gA$.
Montrons d'abord \linebreak 
que  $\gA\div \fb=\gA$.
Soit $y \in \gA\div \fb$, d'où $y (\gA\div \fa) \subseteq (\gA\div \fa)$.
Puisque $\gA\div \fa$ est un module fidèle (il contient 1) et \tfz, $y$ est
entier sur $\gA$ (cf. fait \ref{factEntiersAnn}) donc  $y\in\gA$ car $\gA$ est
normal.
\\
Par  \recuz, en utilisant $\gA\div \fb^{k+1} = (\gA\div \fb)\div
\fb^k$, on obtient $\gA\div \fb^k = \gA$ pour tout $k \geq 1$.

Fixons $x \in \fb$ un \elt\ndzz.  Par le lemme \ref{lemDim1-1}, il existe un
 $k \in\NN\sta$ tel que $\fb' := \gen{x}+\fb^k$ est \ivz.  En conséquence
$\fb' (\gA\div\fb') = \gA$. Enfin, comme $\fb^k\subseteq\fb'\subseteq\fb$, on a
$\gA\div \fb'=\gA$, d'où $\fb' =\gA$ puis $\fb=\gA$.
\end{proof}
%

%:2014  rajout de l'exemple
\medskip
\exl Outre l'exemple des \advs donné \paref{exempleadvgroupe}, 
qui peuvent avoir une \ddk arbitraire, il y a d'autres
exemples naturels de \ddps qui ne sont pas de dimension~$\leq 1$.
\\
On appelle \emph{anneau des \pols à valeurs entières} le sous-anneau de
$\QQ[X]$ formé par les \pols $f(X)$ tels que $f(x)\in\ZZ$
pour tout $x\in\ZZ$. On montre facilement que  c'est un \ZZmo libre admettant
pour base les \pols combinatoires $ {x} \choose{n}$ pour $n\in\NN$. 
L'\id engendré par les \pols $ {x} \choose{n}$ pour $n\geq 1$
n'est pas \tfz. On montre qu'un \pol à valeurs entières
peut être évalué en un entier $p$-adique arbitraire, ce qui fournit
un ensemble non dénombrable d'\idepsz. Cet anneau est un \ddp de dimension 2, mais la \dem de ce résultat n'est pas simple, surtout si l'on demande qu'elle soit \covz.
Voir à ce sujet  \cite[Ducos]{DucPVE} et \cite[Lombardi]{LomPVE}.
\eoe

%%%%%%%%%%%%%%%%%%%%%%%%%%%%%%%%%%%%%%%%%%%%%%%%%%%%%%%%%%%%%%%%%%%%%%%%%%%
\subsec{Structure des \mpfsz}

%:      Theorem{thPTFDed}----------
\begin{theorem}
\label{thPTFDed}
Soit $\gA$ un \adpc \ddi1. Tout module \pro $M$
de rang constant \hbox{$k \geq 1$} est isomorphe à 
$\Ae {k-1} \oplus \fa$, où $\fa$
est un \id \ivz. 
En particulier, il est engendré par $k+1$ \eltsz.
Enfin puisque $\fa\simeq \Al k M$,
la classe d'\iso de $M$ comme \Amo détermine celle de $\fa$.
\end{theorem}
%--- end-theorem-----------------------------------------
\begin{proof}
D'après le \thrf{ThImMat}, $M$ est une somme directe
de $k$ \ids \ivsz. Il suffit donc de traiter le cas $M\simeq \fa\oplus \fb$, avec
des \ids \ivs
   $\fa$ et~$\fb$.
Par le lemme \ref{lemEvitCo}, on peut trouver un \id $\fa_1$ tel que
   $\fa_1\simeq \fa$ (en tant que~$\gA$-modules) et $\fa_1+\fb=\gen{1}$
(comme \idsz). On a alors la suite exacte courte

\snic{
\displaystyle
\gen{0} \longrightarrow \fa_1 \fb = \fa_1\cap \fb\vers{\delta} \fa_1\oplus \fb
\vers{\sigma}
\fa_1+\fb=\gA\longrightarrow \gen{0},}

%\sni
où $\delta(x)=( x, -x)$ et $\sigma(x, y)=x+y$.
Enfin, puisque cette suite est scindée, 
on obtient
$M\simeq \fa\oplus \fb \simeq \fa_1\oplus \fb \simeq \gA\oplus (\fa_1 \cap \fb)
=\gA\oplus (\fa_1 \fb) $.
\end{proof}
%-----------------end proof------------------

Une conséquence \imde est le \tho de structure suivant.
%--- corolaire{corthPTFDed}-----------------
\begin{corollary}\label{corthPTFDed}
Soit $\gA$ un \adpc de dimension $\leq 1$. Tout  \mptf 
 est isomorphe à une somme directe

\snic{r_{1}\gA \oplus r_{2}\Ae 2  \oplus \cdots   \oplus r_{n}\Ae n\oplus \fa,}

%\sni
où 
les $r_{i}$ sont des \idms \orts (certains peuvent être nuls)
et $\fa$ est un \itfz.
\end{corollary}

%:     Proposition{propAriCohZed}
\begin{proposition}\label{propAriCohZed}
Soit $\gA$ un \anar  \zedz. Toute matrice admet une forme réduite de Smith.
En conséquence tout \Amo \pf    est isomorphe à une somme directe de
modules monogènes~$\aqo\gA{a_k}$. 
\end{proposition}
%--------- fin proposition ----------------------------------------------

%
\begin{proof} 
Si $\gA$ est local, c'est un anneau de Bézout local et la matrice admet une forme réduite de Smith (proposition \ref{propPfVal}), ce qui donne le résultat. 
En suivant la \dem du cas local, et en appliquant la machinerie \lgbe  des \anars \paref{MetgenAnar}, on produit une famille d'\eco $(s_1,\ldots,s_r)$
tels que le résultat est assuré sur chaque anneau $\gA[1/s_i]$.
\\
Puisque $\gA$ est \zedz, tout filtre principal est engendré par un \idm
(lemme \ref{lemme:idempotentDimension0}~\emph{\iref{LID002}}). On considère les \idms
$e_i$ correspondants aux $s_i$, puis un \sfio $(r_j)$ 
tel que chaque $e_i$ soit une somme  de certains $r_j$.
\\
L'anneau s'écrit comme un produit fini $\prod\gA_j$ avec la réduction de Smith sur chaque  $\gA_j$.
Le résultat est donc assuré. 
\end{proof}
\rems
Pour l'unicité de la \dcnz, voir le  \thref{prop unicyc}. Par ailleurs, la \dem montre que la réduction peut se faire avec des produits de matrices \elrsz.
Enfin une \gnn est proposée en exercice \ref{exoAnarlgb}.
\eoe

%:     Corollary{corpropAriCohZed}
\begin{corollary}\label{corpropAriCohZed}
Soit $\gA$ un \anar \ddi $1$. Tout \Amo \pf de torsion  est 
isomorphe à une somme
directe de  modules monogènes~$\aqo\gA{b,a_k}$ avec $b\in\Reg(\gA)$.
\end{corollary}
%--------- fin corollary ---------------------------------------------- 
%
\begin{proof}
Le module est annulé par un \elt $b\in\Reg(\gA)$. On le regarde
comme un $\aqo\gA b$-module et l'on applique la proposition~\ref{propAriCohZed}.
\end{proof}

On peut maintenant synthétiser les \thrfs{thPTFDed}{ThImMat}, et le corolaire  
\ref{corpropAriCohZed} comme suit.
Nous laissons le soin \alec de donner l'énoncé pour le cas \qi
(i.e. pour les \adpcsz).

%:     Theorem{thMpfPruCohDim}
\begin{theorem}\label{thMpfPruCohDim}\emph{(\Tho des facteurs invariants)}\\
Sur un \ddp   $\gA$ de dimension $\leq1$,  
tout \mpf est somme directe
\begin{itemize}
\item d'un \Amo \ptfz, nul ou de la forme $\Ae r\oplus\fa$ ($r\geq0$, $\fa$
un \id \ivz),
\item   et de son sous-module de torsion,
qui est isomorphe à une somme directe de modules monogènes~$\aqo\gA{b,a_k}$
avec $b\in\Reg(\gA)$.
\end{itemize}
 En outre: 
 \begin{itemize}
\item l'\id $\fa$ est uniquement déterminé par
le module, 
\item  on peut supposer que les \ids $\gen{b,a_k}$
sont totalement ordonnés pour la relation d'inclusion, 
et la \dcn du sous-module de torsion est alors unique
au sens précisé dans le \thref{prop unicyc}.
\end{itemize}
\end{theorem}
%--------- fin theorem ---------------------------------------------- 

%:HHH petit rajout
\rems
1) En particulier, le \tho de structure pour les \mpfs sur un anneau principal
(proposition \ref{propPfPID}) est valable pour tout anneau de Bézout intègre de dimension~$\leq1$.
\\
2) Pour un module de torsion $M$, les \ids $\gen{b,a_k}$
du \tho précédent sont les \ix{facteurs invariants} de $M$, conformément à la \dfn donnée au \thref{prop unicyc}. 
\eoe

%%%%%%%%%%%%%%%%%%%%%%%%%%%%%%%%%%%%%%%%%%%%%%%%%%%%%%%%%%%%%%%%%%%%%%%%%%%
\subsubsection*{Réduction de matrices}

%--- theoreme{mat33}-----------------
Le \tho suivant donne une forme réduite
pour une matrice colonne, à la Bézout. Il serait intéressant de le \gnr
à une matrice quelconque.

%%%%%%%%%%%%%%%%%%%%%%%%%%%%%%%%%%%%%%%%%
\begin{theorem}\label{mat33}
Soit $\gA$ un \adpc de dimension $\leq 1$
et $x_1$, \dots, $x_n\in \gA$.  
Il existe une matrice  $M\in\GL_n(\gA)$ telle que 

\snic{M\,\tra[\,x_1 \;\cdots\; x_n\,]=\tra[\,y_1 \;y_2 \;0 \;\cdots\; 0\,].}
\end{theorem}
%%%%%%%%%%%%%%%%%%%%%%%%%%%%%%%%%%%%%%%%%
\begin{proof}
Il suffit de traiter le cas où $n=3$.
\\
Si $e$ est un \idmz, alors $\GL_n(\gA) \simeq \GL_n(\gA_e) \times 
\GL_n(\gA_{1-e})$:
quitte à localiser en inversant l'\idm annulateur de $\gen{x_1,x_2,x_3}$ et 
son \copz, on peut donc supposer que
$\Ann(\gen{x_1,x_2,x_3})=\gen{0}$.
\\
Soit $A$ une \mlp pour $(x_1,x_2,x_3)$.\\
Le module $K=\Im(\I_3-A)$ est le noyau de la forme \lin associée
au vecteur ligne $X=[\,x_1 \;x_2\; x_3\,]$ et c'est un module \pro de rang $2$
en facteur direct dans $\Ae 3$. Le
\thrf{thPTFDed} nous dit que $K$ contient un sous-module libre de
rang $1$ en facteur direct dans $\Ae 3$, \cad un module~$\gA v$ où $v$ est un \vmd de  $\Ae 3$.
Par le \thrf{thK1-SLE}, ce vecteur est la dernière colonne
d'une matrice \ivz~$U$; le dernier \coe de~$XU$ est nul
et la matrice $M = \tra {U}$ convient.
\end{proof}

%-% ENTRE NOUS
\entrenous{Il nous manque toujours  une réduction en une sorte de forme se Smith
pour les matrices avec plusieurs lignes. 
Le \thref{thMpfPruCohDim} fournit presque le résultat, à ceci près que 
pour une \mpn il y a un peu plus de transformations \elrs autorisées.
Le \tho \ref{thMpfPruCohDim} ne nous donne pas d'argument décisif pour un anneau de Bézout  intègre \ddi1.

Il semblerait utile de comprendre en détail le papier \cite[Brewer\&Klinger]{BrKl}
dans lequel il semble que: 
\begin{itemize}
\item \emph{a.} le \thref{thMpfPruCohDim} est démontré dans une variante légèrement plus forte
\item \emph{b.}  un certain lien avec un \tho purement matriciel est établi
\item \emph{c.} la condition \gui{\ddi $1$} est affaiblie en: pour tout \elt \ndzz,
l'anneau quotient est \lgb (c'est ce que l'on trouve dans la littérature sous
l'appellation  des anneaux
\gui{presque \lgbsz}).
\end{itemize}

Noter aussi que \cite[Couchot]{Couc1} étend les résultats de \cite{BrKl}
des intègres aux \qisz, ce qui n'est pas vraiment un scoop, au moins pour nous.
L'avantage est que son papier est peut-être un peu mieux structuré
et que ses \dfns sont plus faciles à lire pour nous.

En fait il semble que l'exercice \ref{exoAnarlgb} donne exactement
ce que font ces gens là concernant \emph{c}. Il resterait donc surtout à comprendre le \emph{b.}

}
%-% Fin ENTRENOUS

%%%%%%%%%%%%%%%%%%%%%%%%%%%%%%%%%%%%%%%%%%%%%%%%%%%%%%%%%%%%%%%%%%%%%%%%%%%
%%%%%%%%%%%%%%%%%%%%%%%%%%%%%%%%%%%%%%%%%%%%%%%%%%%%%%%%%%%%%%%%%%%%%%%%%%%
%%%%%%%%%%%%%%%%%%%%%%%%%%%%%%%%%%%%%%%%%%%%%%%%%%%%%%%%%%%%%%%%%%%%%%%%%%%
%%%%%%%%%%%%%%%%%%%%%%%%%%%%%%%%%%%%%%%%%%%%%%%%%%%%%%%%%%%%%%%%%%%%%%%%%%%
%%%%%%%%%%%%%%%%%%%%%%%%%%%%%%%%%%%%%%%%%%%%%%%%%%%%%%%%%%%%%%%%%%%%%%%%%%%

%%%%%%%%%%%%%%%%%%%%%%%%%%%%%%%%%%%%%%%%%%%%%%%%%%%%%%%%%%%%%%%%%%%%%%%%%%%
%%%%%%%%%%%%%%%%%%%%%%%%%%%%%%%%%%%%%%%%%%%%%%%%%%%%%%%%%%%%%%%%%%%%%%%%%%%
%%%%%%%%%%%%%%%%%%%%%%%%%%%%%%%%%%%%%%%%%%%%%%%%%%%%%%%%%%%%%%%%%%%%%%%%%%%
%%%%%%%%%%%%%%%%%%%%%%%%%%%%%%%%%%%%%%%%%%%%%%%%%%%%%%%%%%%%%%%%%%%%%%%%%%%
%%%%%%%%%%%%%%%%%%%%%%%%%%%%%%%%%%%%%%%%%%%%%%%%%%%%%%%%%%%%%%%%%%%%%%%%%%%
%%%%%%%%%%%%%%%%%%%%%%%%%%%%%%%%%%%%%%%%%%%%%%%%%%%%%%%%%%%%%%%%%%%%%%%%%%%

%--- Section{Factorisation d'\itfs} secAdPfactpar    
\section{Factorisation d'\itfs} 
\label{secAdPfactpar}
%-------------------------

\vspace{3pt}

%:     SUBsec{subsecfactogene}-----------
\subsec{Factorisations \gnles} 
\label{subsecfactogene}
%-----------------------------------------

Dans un \anar \gnl il semble que l'on n'a pas de résultats 
de \fcn qui aillent au delà
de ce qui découle du fait que les \ids \ivs (\cad les \itfs contenant
un \elt \ndzz) forment un \mo à pgcd, et plus \prmt la partie positive d'un \grlz.

Par exemple le \tho de Riesz se relit comme suit.

%:     Theorem{thRieszAnar}
\begin{theorem}\label{thRieszAnar} \emph{(\Tho de Riesz pour les \anarsz)}\\
Soit $\gA$ un \anarz, $(\fa_i)_{i\in\lrbn}$ et $(\fb_j)_{j\in\lrbm}$ des \ids
\ivs tels que $\prod_{i=1}^n\fa_i=\prod_{j=1}^m\fb_j$.\\ Alors il existe des
\ids \ivs $(\fc_{i,j})_{i\in\lrbn,j\in\lrbm}$ tels que l'on ait pour tout
$i$ et tout $j$: 

\snic{\fa_i=\prod_{j=1}^m\fc_{i,j}\,$ et $\,\fb_j=\prod_{i=1}^n\fc_{i,j}.}  
\end{theorem}

%:     SUBsec{subsecfactodimun}-----------
\subsec{Factorisations en dimension 1} 
\label{subsecfactodimun}
%-----------------------------------------

%--- Théorème{thFactor}----------
\begin{theorem}
\label{thFactor}
Dans un \adpc \ddi1, on con\-si\-dère deux \itfs $\fa$ et $\fb$
avec $\fa$ \ivz. Alors on peut écrire:  

\snic{\fa=\fa_1\fa_2\;$ avec $\;\fa_1+\fb=\gen{1}$ et $\,\fb^n\subseteq \fa_2,}

%\sni
pour un entier $n$ convenable. Cette
écriture est unique et l'on a 

\snic{\fa_1+\fa_2=\gen{1},$
   $\;\fa_2=\fa+\fb^n=\fa+\fb^{n+1}.}
\end{theorem}
%--- end-theorem-----------------------------------------

%-----------------begin proof------------------
\begin{proof}
Ceci est un cas particulier du lemme \ref{lemDim1-1}.
\end{proof}
%-----------------end proof------------------

\rem On n'a pas  besoin de supposer les \ids détachables.

%--- Théorème{thFactor2}---------
\begin{theorem}
\label{thFactor2}
On considère dans un \adpc de dimension $\leq 1$ des \itfs $\,\fp_1$, $\dots$, $\fp_n$ deux à deux \comz, et un \id \iv $\fa$.\\
On peut écrire $\fa=\fa_0\cdot \fa_1\cdots\, \fa_n$ avec les \itfs 
$\,\fa_0$, $\dots$, $\fa_n$ deux à deux \com
et, pour $j \geq 1$, $\fp_j^{m_j}\subseteq \fa_j$ avec $m_j$ entier convenable.
\\
Cette écriture est unique et l'on a $\,\fa_j=\fa+\fp_j^{m_j}=\fa+\fp_j^{1+m_j}$.
\end{theorem}
%--- end-theorem-----------------------------------------
\begin{proof} Par \recu en utilisant
le \thrf{thFactor} avec $\fb \in \{\fp_1,\dots,\fp_n\}$.
\end{proof}

%:     SUBsec{subsecDedfactpar}-----------
\subsec{Anneaux de Prüfer à \fapz} 
\label{subsecDedfactpar}
%-----------------------------------------
Reprenons la \dfn des \dcnps  (donnée pour les \grlsz)
dans le cadre du \mo des \ids \ivs d'un \adpc $\gA$  (ce \mo est la partie positive du 
\grl formé par les \elts \ivs de $\Ifr(\gA)$).

%%%%%%%%%%%%%%%%%%%%%%%%%%%%%%%%%%%%%%
%--- Defi factorisation partielle --
\begin{definition}
Soit $F=(\fa_1,\dots,\fa_n) $ une famille finie
d'\ids \ivs dans un anneau $\gA$. On dit que $F$ admet une
\emph{\fapz} s'il existe une famille
   $P=(\fp_1,\dots,\fp_k)$ d'\ids \ivs deux à deux \com telle que tout
\id $\fa_j$  peut s'écrire sous la forme:
   $\fa_j = \fp_1^{m_{1j}}\cdots \fp_k^{m_{kj}} $ (certains des $m_{ij}$
peuvent être nuls).
%:HHH ci dessous  bdf  plutot que  bdp
On dit alors que $P$ est une \emph{\bdfz} pour la
famille $F$.%
\index{factorisation partielle!base de ---}
\end{definition}
%%%%%%%%%%%%%%%%%%%%%%%%%%%%%%%%%%%%%%

Pour que le \mo $\Ifr(\gA)$ soit discret il faut supposer que $\gA$ est \fdiz.
Ceci conduit à la \dfn suivante.

%%%%%%%%%%%%%%%%%%%%%%%%%%%%%%%%%%%%%%
%--- Defi Prufer a fact partielle defPrufFP
\begin{definition}
\label{defPrufFP}
Un anneau est appelé \emph{\adp à \fapz} si c'est un \adpc
\fdiz\footnote{D'après la proposition~\ref{aritfdi} un \anar est
\fdi \ssi la relation de \dve est  explicite.}
et si toute famille finie d'\ids \ivs admet une \fapz.%
\index{anneau!de Prüfer à \fapz}%
%:HHH ajout d'un index \index{factorisation!partielle}
\index{factorisation!partielle}\index{Prufer@Prüfer!anneau de --- à \fapz}
\end{definition}
%%%%%%%%%%%%%%%%%%%%%%

%--- Lemme{lemK1}--------------------
\begin{lemma}
\label{lemK1}
Un \adp à \fap est \ddi1.
\end{lemma}
%--- end-lemme-----------------------------------------

%-----------------begin proof------------------
\begin{proof}
On considère un \elt  $y$  \ndzz. On veut montrer que $\aqo{\gA}{y}$
est \zedz. \\
Pour cela on prend un $x$ \ndz et l'on veut trouver $a\in\gA$ et $n\in\NN$
tels que $x^n(1-ax)\equiv 0 \mod y$. 
%\\
%Nous remplaçons chaque \idp par un \gtr pour alléger les notations.
La \fap de   $(x,y)$ nous donne 

\snic{\gen{x} = \fp_{1}^{\alpha_1} \cdots \fp_{i}^{\alpha_i}\fq_{1}^{\beta_1} \cdots
\fq_{j}^{\beta_j} = \fa\fb, \hbox{  et }
\gen{y} = \fp_{1}^{\gamma_1} \cdots \fp_{i}^{\gamma_i} \fh_{1}^{\delta_1}
\cdots \fh_{k}^{\delta_k} = \fc\fd}

%\sni
avec tous les exposants non nuls.
Il existe $n\geq0$ tel que $\fa^n$ soit un multiple de $\fc$ ce qui donne
$\gen{x^n} = \fc\ffg$.
Comme $\gen{x} + \fd = 1$, il existe un~$a \in \gA$ tel que~$1 - ax \in \fd$. On a
donc $\gen{y}=\fc\fd\supseteq  \fc\ffg\fd = \gen{x^n}\fd \supseteq \gen{x^n(1 - ax)}$,
\cadz~$x^n(1-ax)\equiv 0 \mod y$.
\end{proof}
%-----------------end proof------------------

%:     SUBsec{subsecDedfactot}-----------
\subsec{Anneaux de Dedekind} 
\label{subsecDedfactot}
%-----------------------------------------

%:     definition{defDDK}
\begin{definition}\label{defDDK}
On appelle \emph{\adkz} un \adp \coh \fdi et \noez. 
Un \emph{\dDkz} est un \adk intègre (ou encore connexe).%
\index{anneau!de Dedekind}%  
\index{Dedekind!anneau de ---}  
\end{definition}

%:     theorem{prop2DDK}
\begin{theorem}\label{prop2DDK}
Un \adk est un \adp à \fapz, donc \ddi1. 
\end{theorem}
\begin{proof}
Le \thrf{th2GpRtcl}
donne le résultat de \fap dans le cadre des \trdis et l'on termine avec le
lemme \ref{lemK1}.
\end{proof}
%

%:     theorem{propDDK}
\begin{theorem}\label{propDDK} \emph{(\Carns des \adksz)}\\
Pour un anneau $\gA$ \propeq

\begin{itemize}
\item [1.] $\gA$ est un \adkz.
\item [2.] $\gA$ est \qiz, \ariz, \dvee et \noez.
\item [3.] $\gA$ est \qiz, normal, \ddi1,  \dveez, \coh et \noez.
\end{itemize}
\end{theorem}
\begin{proof}
Puisque $\gA$ est un anneau de Prüfer \coh \ssi il est \ari et \qiz, et puisque 
un anneau \ari est \fdi \ssi il est  \dveez, les points \emph{1}
et \emph{2} sont \eqvsz.
L'implication \emph{1} $\Rightarrow$  \emph{3} résulte du \thrf{prop2DDK},
et le \thrf{thcohdim1} donne la réciproque (il faut simplement ajouter \fdi et \noe dans l'hypothèse et la conclusion).
\end{proof}
% 

%:--- Defi factorisation totale --
\begin{definition}
Soit $\fa $ un \id d'un anneau $\gA$. On dit que $\fa$ admet une
\emph{\facz} s'il s'écrit $\fa = \fp_1^{m_{1}}\cdots \fp_k^{m_{k}} $ ($m_{i}>0$, $k>0$) avec les \ids $\fp_i$  maximaux stricts détachables (autrement dit, chaque 
\linebreak 
anneau $\gA\sur{\fp_i}$ est un \cdi non trivial).%
\index{factorisation totale! d'un \id dans un anneau}%
\index{factorisation!totale}%
\end{definition}

%:     Thdef{thdefDDKTOT}
\begin{thdef}\label{thdefDDKTOT}
Pour un anneau $\gA$ \qi \fdi non trivial, \propeq
\begin{enumerate}
\item Tout \idp $\gen{a}\neq\gen{1}$ avec $a\in\Reg\gA$ admet une \facz.
\item L'anneau $\gA$
est un \adkz, et tout \id \iv $\neq\gen{1}$ admet une \facz.
\end{enumerate}
Un tel anneau est appelé un \emph{\adk à \facz}.% 
\index{anneau!de Dedekind à \fcn totale}%  
\index{Dedekind! anneau de --- à \fcn totale}%  
\index{factorisation totale!anneau de Dedekind à ---}%
\index{factorisation!totale}%
\end{thdef}
%--------- fin theorem ---------------------------------------------- 
%
\begin{proof} Il faut montrer que \emph{1} implique \emph{2.}
On traite le cas intègre (le cas \qi s'en déduit facilement). 
\\
On se reporte à l'exercice \ref{exoDecompIdeal} et à sa correction. On voit que tout \itf
contenant un \elt \ndz est \ivz, et qu'il admet une \facz. Le \thref{th.adpcoh} nous dit alors que $\gA$
est un \adpcz. Il reste à voir qu'il est \noez. 
On considère un \itf et sa \fac  $\fa = \fp_1^{m_{1}}\cdots \fp_k^{m_{k}} $. Tout \itf  $\fb\supseteq\fa$ s'écrit $\fp_1^{n_{1}}\cdots \fp_k^{n_{k}}$
avec les $n_i\in\lrb{0..m_i}$. 
Toute suite croissante d'\itfs démarrant avec $\fa$
admet donc deux termes consécutifs égaux. 
\end{proof}

\rem L'exercice \ref{exoDecompIdeal} n'utilise aucun attirail théorique 
compliqué. 
Aussi il est possible d'exposer la théorie des \adks en commençant par
le \tho précédent, qui conduit rapidement aux résultats essentiels. Le principal inconvénient de cette approche est qu'elle est basée sur une \prt de \fac qui n'est pas \gnlt satisfaite du point de vue \cofz, même par les anneaux principaux, et 
qui ne s'étend pas en \gnl aux extensions entières.  
\eoe

\medskip 
Rappelons que nous avons déjà établi le \thref{thAESTE} 
concernant les extensions finies d'\adksz. 

Nous pouvons rajouter la précision suivante.

%:h2013 énoncé et preuve modifiés
%:     theorem{lemthAESTE}
\begin{theorem} \label{lemthAESTE} \emph{(Un calcul de \cliz)}\\
Soit $\gA$ un \adkz, 
$\gK=\Frac(\gA)$, $\gL\supseteq\gK$ une \Klg étale et $\gB$ la \cli de $\gA$ dans $\gL$.\\
Supposons que $\gL=\aqo\KX f$ avec $f\in\AX$ \mon et 
$\disc_X(f) \in\Reg\gA$ (ce qui n'est pas vraiment restrictif). Si $\gen{\disc_X(f)}$ admet une \facz, et si pour chaque \idema $\fm$ de cette \fcnz, le corps résiduel $\gA/\fm$ est parfait, alors $\gB$ est un \Amo \ptfz.
\end{theorem}
%--------- fin theorem ---------------------------------------------- 
%
\begin{proof} 
Comme $\gA$ est \qiz, il suffit de traiter le cas où $\gA$
est intègre (machinerie \lgbe \elr des anneaux \qisz), donc~$\gK$ est un \cdiz.
L'hypothèse $\gL=\aqo\KX f$ avec $f\in\AX$ \mon et 
$\disc_X(f) \in\Reg\gA$ n'est pas vraiment restrictive
car d'après le \thref{thEtalePrimitif}, $\gL$ est un produit de \Klgs étales monogènes.
On peut même supposer que $\gL$ est un corps étale sur $\gK$  (machinerie \lgbe \elr des anneaux \zedrsz).\\
On pose $\Delta=\disc_X(f)$.
D'après le point \emph{\ref{i5thAESTE}} du \thref{thAESTE} on a les inclusions

\snic{\gA[x]\subseteq\gB\subseteq\fraC 1 {\Delta}\,\gA[x].}

Ainsi $\gB$ est un sous-module du \Amo \tf $\fraC 1 {\Delta}\gA[x]$. D'après le \thref{ThImMat}, si $\gB$ est \tfz, il est \ptfz.   
\\
On a  $\gA[x,\fraC 1 {\Delta}]=\gB[\fraC 1 {\Delta}]$, donc $\gB$
est \tf après \lon en $\Delta^{\NN}$. Il reste à montrer que
$\gB$ est \tf après \lon en $S=1+\Delta\gA$. L'anneau $\gA_S$ est un anneau de Bézout (\thref{dekinbe}). Si $\fp_1$, \dots, $\fp_r$ sont les \idemas qui interviennent
dans la \fac de $\Delta$, les \mos $1+\fp_i$ sont \com dans 
 $\gA_S$, et il suffit de montrer que $\gB$ est \tf après \lon en chacun des
 $1+\fp_i$.
On est ainsi ramené au cas traité dans le lemme \ref{lemlemthAESTE}
qui suit. \end{proof}

Notez qu'un \dDk local $\gV$ est aussi bien un \ddv \noe \fdiz, ou encore un anneau principal local avec~$\gV\eti$ détachable. Dans le lemme suivant, on demande en outre que $\Rad\gV$ soit principal (ce qui est  automatique en \clamaz). Dans ce cas on dira que $\gV$ est un \emph{\adv discrète}, selon la terminologie classique, et un \gtr de $\Rad(\gV)$ est appelé \emph{une uniformisante}.%
\index{anneau!de valuation discrète}\index{valuation!anneau de --- discrète}
\index{valuation!discrète}\index{uniformisante} 

%:     Lemma{lemlemthAESTE}
\begin{lemma} \label{lemlemthAESTE}
Soit $\gV$  un \advd  avec $\Rad\gV=p\gV$. On suppose le corps résiduel $\gk=\aqo\gV p$  parfait. Soit  $f\in\VX$ un \polu \ird  \hbox{(donc 
$\Delta=\disc_X(f) \in\Reg\gV$)}.\\
On note $\gK=\Frac(\gV)$, $\gL=\Kx=\aqo\KX f$, et~$\gW$ la \cli de~$\gV$
dans~$\gL$. Alors $\gW$ est \tf sur $\gV$.
\end{lemma}
%--------- fin lemma ---------------------------------------------- 
%
\begin{proof} Puisque $\gk$ est parfait, d'après le lemme \ref{lemSqfDec}, pour tout \polu $f_i$ de $\VX$ on sait calculer la \gui{partie sans carré} de $\ov {f_i}$
($f_i$ vu modulo $p$), i.e. un \pol $\ov {g_i}$ \spl dans~$\gk[X]$ qui divise $\ov {f_i}$, et dont une puissance est multiple de $\ov {f_i}$.
\\
La stratégie est de rajouter des \elts $x_i\in\gW$ à $\gV[x]$ jusqu'au
moment où l'on obtient un anneau $\gW'$ dont le radical est un \id \ivz.
Lorsque ceci est réalisé, nous savons d'après le lemme \ref{lemRadJDIM1}
que $\gW'$ est un \ddpz, donc qu'il est \iclz, donc égal à $\gW$.\\
Pour \gui{construire} $\gW'$ (de type fini sur $\gV$) on va utiliser dans une \recu le fait suivant, initialisé avec $\gW_1=\gV[x]$ ($x_1=x$, $r_1=1$). 

{\it Fait. Soit $\gW_k=\gV[x_1,\dots,x_{r_k}]\subseteq \gW$, alors 

\snic{\Rad(\gW_k)=\geN{p,g_1(x_1),\dots,g_k(x_{r_k})},
}

où $\ov{g_i}$ est la partie sans carré de $\ov{f_i}$, $f_i$ \polmin sur $\gK$ de l'entier~$x_i$.}
\begin{Proof}{}
Le \thref{thJacplc} nous dit que $\Rad(\gW_k)=\rD_{\gW_k}(p\gW_k)$.
Cet \id est l'image réciproque de $\rD_{\gW_k/p\gW_k}(0)$
et l'on a 
$\gW_k/p\gW_k=\gk[\ov{x_1},\dots,\ov{x_{r_k}}]$. Comme les $g_i(x_i)$ sont nilpotents modulo $p$ par construction, ils sont dans le nilradical $\rD_{\gW_k}(p\gW_k)$. Il nous suffit maintenant de vérifier
que la \klg 

\snic{\Aqo{\gk[\ov{x_1},\dots,\ov{x_{r_k}}]}{\ov {g_1}(\ov{x_1}), \dots, 
\ov {g_{r_k}}(\ov{x_{r_k}})}}

est réduite. En fait $\gW_k$ est un sous-$\gV$-\mtf de $\fraC1\Delta\gV[x]$,
donc est libre fini sur $\gV$. En conséquence $\gW_k/p\gW_k$ est strictement finie sur $\gk$, et elle est étale  
parce qu'elle est engendrée par des \elts qui annulent des \pols \spls sur $\gk$ (\thref{corlemEtaleEtage}). 
\end{Proof}

Ceci étant vu, puisque $\gW$ est un \ddpz, nous savons inverser l'\itf
$\Rad(\gW_k)$ dans $\gW$. \\
Cela signifie calculer des \elts
$x_{r_{k}+1}$, \dots, $x_{r_{k+1}}$ de $\gW$ et un \itf $\ffg_k$ dans le nouvel anneau $\gW_{k+1}$ tels que l'\id produit $\ffg_k\Rad(\gW_k)$ soit principal (non nul).\\
Il se peut cependant que les \gtrs de $\Rad(\gW_k)$ n'engendrent pas
l'\id $\Rad(\gW_{k+1})$ de $\gW_{k+1}$, ce qui oblige à itérer le processus.\\
La suite croissante des $\gW_k$ est une suite croissante de
$\gV$-\mtfs contenus dans $\fraC1\Delta \gV[x]$, donc elle admet deux termes consécutifs égaux. Dans ce cas on a atteint le but prescrit.
\hfill\eop
\end{proof}

%:    Principe local global concret{plcc.ddk}---
\begin{plcc} 
\label{plcc.ddk}\relax \emph{(Anneaux de Dedekind)}\\
Soient $s_1$, $\ldots$, $s_n$  des \eco d'un anneau $\gA$.   
Alors
\begin{enumerate}
\item L'anneau $\gA$ est  \coh \noe \fdi   \ssi chacun des~$\gA_{s_i}$ est
 \coh \noe \fdiz.
\item L'anneau $\gA$ est un \adk  \ssi chacun des~$\gA_{s_i}$ est un \adkz.
\end{enumerate}
 
\end{plcc}
%--- end-plcc-----------------------------------------
%
\begin{proof} On sait déjà que le \plgc fonctionne pour les \adps  et pour les \coris avec des \mocoz. Il en va de même pour les anneaux ou modules \noes
(une \dem est donnée avec le \plgref{plcc.ptf}).
\\
Il reste  à examiner la \prt \gui{\fdiz} dans le cas d'\ecoz. 
Soit $\fa$ un \itf et~$x\in\gA$.
Il est clair que si l'on a un test pour~$x\in\fa\gA_{s_i}$ 
pour chacun des~$s_i$, cela fournit un test pour  $x\in\fa\gA$. La difficulté est dans l'autre sens: si $\gA$ est \fdi et si $s\in\gA$, alors $\gA[1/s]$ est \fdiz. Ce n'est pas vrai en \gnlz, mais c'est vrai 
pour les \coris \noesz. En effet, l'appartenance $x\in\fa\gA[1/s]$
équivaut à $x\in(\fa:s^{\infty})_\gA$. Or l'\id  $(\fa:s^{\infty})_\gA$
est la réunion de la suite croissante des \itfs $(\fa:s^{n})_\gA$, et dès
que $(\fa:s^{n})_\gA=(\fa:s^{n+1})_\gA$, la suite devient constante.
\end{proof}
%

% section 8
\section{Anneau intègre versus anneau \sdzz}
\label{secAnor}

%\section{Anneau intègre versus anneau \sdzz}
%\label{secAnor}

%: subsec{Motivation}
\subsec{Motivation}
La principale motivation de cette section est de fournir  une \prco du \tho suivant.

%:     Theorem{thAnormalAXnormal}
\begin{theorem} \label{thAnormalAXnormal}
Si $\gA$ est un \anor il en va de même pour $\AX$. 
\end{theorem}
%--------- fin theorem ---------------------------------------------- 

On a vu en comparant les \dems du \thref{thExtEntPruf} et du \thref{factAdpIntExt} que le cas \gui{normal intègre}
(ou normal \qiz)
se traite plus facilement que le cas \gui{normal tout court}. Pourtant un anneau normal est \lsdz et dans le cas local, qui sert de référence pour les preuves classiques, {après \lon en un \idep arbitraire}, 
la différence entre \gui{intègre} et \gui{\sdzz} n'est sensible qu'en \comaz.

Ceci pose un intéressant \pb de décryptage des \dems classiques qui utilisent un test d'\egt dans un contexte où un tel test fait défaut.

Nous donnons dans cette section des exemples significatifs qui montrent que la stratégie implicitement utilisée pour démontrer le \thref{thExtEntPruf} peut être  souvent appliquée avec succès, sinon toujours. 

La stratégie est la suivante. Tout d'abord on vérifie que la \dem
classique du \gui{cas intègre} (voire du cas local intègre) est suffisamment claire pour être rendue \covz. 

Ensuite on reprend la \dem précédente en la modifiant \gui{un petit peu} de façon à ce qu'elle puisse s'appliquer au cas \gui{\sdzz}.
Enfin, on traite le cas \lsdz par la technique usuelle qui consiste à ouvrir deux branches de \lons \come pour le calcul chaque fois que \ncrz, et à \gui{recoller} les résultats.

Ainsi on peut voir cette section comme  quelques exemples illustrant une \emph{machinerie \lgbe \elr des anneaux \lsdzz}.  

%: subsec{Un premier exemple}
\subsec{Un premier exemple}

On commence par un énoncé très simple, qui permet de voir les choses fonctionner.

%:     Lemma{lemSDZAAX}
\begin{lemma} \label{lemSDZAAX}
Si l'anneau $\gA$ est \lsdz il en va de même pour l'anneau $\AX$. 
\end{lemma}
%--------- fin lemma ---------------------------------------------- 
%
\begin{proof} On considère $f$ et $g$ deux \pols tels \hbox{que $fg=0$}.
\\
\Prmtz,  $f=\sum_{k=0}^na_kX^{k}$ et $g=\sum_{j=0}^mb_jX^{j}$.
\\
On va trouver $u$ et $v\in\AX$, \hbox{avec $u+v=1$}, $uf=0$ 
et $vg=0$. En fait, on va trouver $u$ et $v\in\gA$. 

On rappelle d'abord la \dem donnée dans le cas où l'anneau est non trivial et intègre.
\\ Si $f=0$, c'est bon.
Si un \coe de $f$ est non nul, le degré de~$f$ est un entier~\hbox{$d\geq 0$} bien défini.  On démontre par \recu descendante
sur $j$ que tous les~$b_j$ sont nuls.

Voici la preuve dans le cas \und{\sdzz} qui en résulte. On doit montrer que 
$f=0$ ou $g=0$ ($u=1$ ou $v=1$). 
\\
On  fait une \dem par \recu descendante sur $m+n$.\\
Si $n=m=0$, le résultat est clair (on est dans $\gA$).\\ 
Supposons $n+m\geq 1$.
Comme $a_nb_m=0$, on a  $a_n=0$ \hbox{ou $b_m=0$}.  On applique l'\hdrz.

Voici enfin la preuve dans le cas \und{\lsdzz} qui résulte
de la précédente (laquelle fonctionnait notamment dans le cas local).
\\
On fait une   \dem par \recu descendante sur $m+n$ (avec $m$ et $n\geq 0$). L'initialisation pour $n=m=0$ est \imdez. Voyons l'étape de \recuz.
\\
On a $a_nb_{m}=0$, donc il existe $s$, $t\in\gA$
tels que 

\snic{s+t=1 ,\;\;   sa_n=0\; \hbox{ et } \;tb_{m}=0.}

Si $n>0$, l'\hdr s'applique aux \pols $sf$ et $g$ et fournit~\hbox{$u_1,v_1\in\AX$} tels que $u_1+v_1=1$, $u_1sf=0$, $v_1g=0$. 
\\
Si $m>0$, l'\hdr s'applique aux \pols $f$ et $tg$ et fournit~\hbox{$u_2,v_2\in\AX$} tels que $u_2+v_2=1$, $u_2f=0$, $v_2tg=0$.
\\
 Ainsi si $n$ et $m$ sont $>0$, on pose $u=su_1+tu_2$ et $v=sv_1+tv_2$, et l'on obtient $u+v=1$, $uf=0$ et $vg=0$.
\\
 Enfin, si par exemple $n>0$ et $m=0$, on reprend le même calcul \hbox{avec $u_2=0$} \hbox{et $v_2=1$}: en effet $u_2f=0$ et $v_2tg=tg=tb_0=0$.
\end{proof}
%

%%%%%%%%%%%%%%%%%%%%%%%%%%%%%%%%%%%%%%%%%%%%%%%%%%%%%%%%%%%%%%%%%%%%%%%%%%%
%: subsec{Une version \gnee du lemme \ref{lemZintClos}}
\subsec{Une version \gnee du lemme \ref{lemZintClos}}

On cherche maintenant une \gnn du lemme \ref{lemZintClos} dont nous rappelons l'énoncé. 

\medskip {\bf Lemme \ref{lemZintClos}.} {\it Lorsque $\gK$ est un \cdiz, l'anneau $\KX$ 
est \iclz.
}

\medskip La \gnn  recherchée consiste à remplacer~$\gK$ par un anneau réduit~$\gA$.
On doit alors remplacer dans la conclusion \gui{anneau \iclz} par \gui{anneau normal}.
On a donc en vue une \prt du style: \gui{Si $a\in\AX$ est entier \hbox{sur $b\AX$}, \hbox{alors $a\in b\AX$}}. 

Nous introduisons pour ceci la notion de \emph{pseudo-reste pour la division
en puissances croissantes}.

Considérons $a$ et $b$ de degrés formels  $\leq m$, autrement dit 

\snic{a=\sum_{k=0}^ma_kX^{k}$ et $b=\sum_{j=0}^mb_jX^{j}.}

Nous ne pouvons pas diviser $a$ par $b$ dans $\AX$ en puissances croissantes, mais c'est possible lorsque $b_0=1$. Le reste à l'ordre $m$ est alors
égal au \deter d'une matrice \gui{de type Sylvester}.
 
Par exemple l'\egt de la division par puissances croissantes
à l'ordre $3$ est 

\snic{a(X)=b(X)q_3(X)+X^{3}s_3(X), \quad \deg(q_3)\leq 2.}

Et le reste  $X^{3}s_3(X)=r_3(X)$  est égal au \deter suivant
(dans lequel on a mis $\cdot$ à la place de $0$ pour plus de lisibilité)  
 $$
X^{3}s_3(X)=\dmatrix{ 
                          1 & b_{1}& b_{2}& b(X)\cr
                          \cdot & 1& b_{1}& Xb(X)\cr
                          \cdot & \cdot &1 & X^{2}b(X)\cr
                          a_0 &a_1 &a_2 & a(X)}
=
\dmatrix{ 
                          1 & b_{1}& b_{2}& b(X)\cr
                          \cdot & 1& b_{1}& Xb(X)\cr
                          \cdot & \cdot &1 & X^{2}b(X)\cr
                          \cdot &\cdot &\cdot & r_3(X)}.
$$
La deuxième matrice est en effet obtenue par \mlrs de lignes à partir de la première.

Lorsque $b_0$ n'est plus supposé égal à $1$, on définit le \emph{pseudo-reste à l'ordre $3$ de $a$ par $b$} comme le \deter suivant:
$$
\mathrm{Prsc}_X(a,b,3)=\dmatrix{ 
                          b_{0} & b_{1}& b_{2}& b(X)\cr
                          \cdot & b_{0}& b_{1}& Xb(X)\cr
                          \cdot & \cdot &b_{0} & X^{2}b(X)\cr
                          a_{0} &a_{1} &a_{2} & a(X)}.
$$
Lorsque $b_0$ est \ivz, la division est encore possible et l'on obtient, au lieu \hbox{de $\mathrm{Prsc}(a,b,3)=r_3$}, l'\egt $\mathrm{Prsc}(a,b,3)=b_0^{3}r_3$. Ainsi, lorsque l'anneau est intègre et $b_0\neq 0$, le \pol $\mathrm{Prsc}(a,b,3)$  est proportionnel au reste de la division par puissances croissantes à l'ordre $3$ opérée dans le corps de fractions.

Lorsque $a_0=b_0 c$ dans $\gA$, on a la version \gui{améliorée} suivante
$$
R_3(X)=\dmatrix{ 
                          1 & b_{1}& b_{2}& b(X)\cr
                          \cdot & b_{0}& b_{1}& Xb(X)\cr
                          \cdot & \cdot &b_{0} & X^{2}b(X)\cr
                          c &a_{1} &a_{2} & a(X)}.
$$
de sorte que $\mathrm{Prsc}(a,b,3)=b_0R_3$ et $R_3=b_0^{2}r_3$.
\\
En développant le \deter $R_3$ selon la dernière colonne, on 
obtient une \egt dans $\AX$

\snic{R_3(X)=b_0^{2} a(X)- Q_3(X) b(X),\quad \deg(Q_3)\leq 2.}

Notons enfin dans le cas des \cdisz, que si  $a$ et $b$
sont de degrés~\hbox{$\leq 2$} et \hbox{si $b_0\neq 0$}, alors $a$ est multiple de $b$ \ssi
le reste à l'ordre $3$ est nul, \ssi $R_3=0$.

On peut \gnr ces \dfns pour un reste à un ordre $m\geq 1$ arbitraire,
sur un anneau $\gA$ arbitraire, avec

\snic{a=\sum_{k=0}^ma_kX^{k}$ et $b=\sum_{j=0}^mb_jX^{j},}

et en hypotèse une \egt  $a_0=b_0 c$ dans $\gA$. On pose alors
$$
R_m(X)=\dmatrix{ 
               1 & b_{1}& b_{2}& \cdots&b_{m-1} & b(X)\cr
               0 & b_{0}& \ddots& &  & Xb(X)\cr
               \vdots & \ddots &\ddots  &\ddots&&\vdots \cr
               \vdots &  & \ddots &\ddots&&\vdots \\[.2em]
               0 & \cdots & \cdots &0&b_0&X^{m-1}b(X) \\[.3em]
               c &a_{1} &\cdots &\cdots&a_{m-1}& a(X)}.
$$
En développant ce \deter selon la dernière colonne, on obtient une \egt dans $\AX$ 

\smallskip \centerline{\fbox{$b_0^{m-1}a(X)=Q_m(X)b(X)+R_m(X),\quad \deg(Q_m)<m, \,R_m\in X^{m}\AX$}.}

\smallskip 
Le {lemme \ref{lemZintClos}} peut alors se relire comme suit.\\
 \emph{Soit $\gK$ un \cdiz, $c\in\gK$, $a=\sum_{k=0}^ma_kX^{k}$ et $b=\sum_{j=0}^mb_jX^{j}\in\KX$  avec $a(X)$ entier sur l'\id ${b(X)}\KX$ et $a_0=b_0c$. \hbox{Alors  $b_0 R_{m+1}=0.$}}

Notez que le résultat est bien correct pour toute valeur de $b_0$, \cad y compris lorsque $b_0=0$.

%:     Lemma{lem2ZintClos}
\begin{lemma} \label{lem2ZintClos} \emph{(Version du lemme \ref{lemZintClos} pour un anneau réduit)}\\
Soit $\gA$ un anneau réduit, $c\in\gA$, 
$a=\sum_{k=0}^ma_kX^{k}$ et $b=\sum_{j=0}^mb_jX^{j}\in\AX$  avec $a(X)$ entier sur l'\id ${b(X)}\AX$ et $a_0=b_0c$. \hbox{Alors  $b_0 R_{m+1}=0.$} 
\end{lemma}
%--------- fin lemma ---------------------------------------------- 
%
\begin{proof}
Cela résulte du lemme \ref{lemZintClos}, relu comme ci-dessus, par le \nst formel. En effet l'hypothèse selon laquelle $a$ et est entier sur l'\id $\gen{b}$
dans $\AX$ signifie une famille d'\egts \polles portant sur les \coes de $a$, ceux de $b$ et ceux d'autres \pols $e_j(X)$ donnés implicitement dans l'hypothèse
puisqu'on doit avoir
une \egt dans $\AX$

\snic{a^{p}=e_1ba^{p-1}+\cdots+ e_{p-1}b^{p-1}a+e_pb^{p}.}

Quant à la conclusion du lemme, elle signifie que certains \pols
en $c$ et en les \coes de $a$ et $b$ sont nuls. 
On est donc bien dans le cadre du \nst formel \ref{thNSTsurZ}.
\end{proof}
%

%:\subsec{\Dem du \thref{thAnormalAXnormal}
\subsec{\Demo du \thref{thAnormalAXnormal}}

Rappelons d'abord la \dem \cov donnée dans le cas d'un anneau \iclz,
\cad normal et intègre.

\medskip  {\bf \Tho \ref{thIntClosStab}.} {\it Si $\gA$ est normal intègre, il en est 
de même pour~$\AX$.}

\begin{proof}
Posons $\gK=\Frac\gA$. Si  un \elt $f$ de $\gK(X)$ est entier sur $\AX$, il est entier sur $\KX$,
donc dans $\KX$ car $\KX$ est intégralement clos (lemme \ref{lemZintClos}).
On conclut avec le lemme~\ref{lemPolEnt}: tous les \coes du \pol $f$ sont entiers sur $\gA$,
donc dans $\gA$.
\end{proof}

Voici ensuite  la \dem \gui{naturelle} du \thref{thAnormalAXnormal} en \clamaz.
\begin{Proof}{\Demo du \thref{thAnormalAXnormal} en \clamaz. } \\
Par \dfn un anneau est normal \ssi
il est \icl après \lon en un \idep arbitraire.
Soit $\fP$ un \idep de~$\AX$ \hbox{et $\fp=\gA\cap \fP$} sa trace sur $\gA$.
L'anneau $\AX_\fP$ est un localisé de~$\gA_\fp[X]$,  il suffit
donc de montrer que $\gA_\fp[X]$ est \iclz. Or $\gA_\fp$ est \iclz.
On termine avec le \thref{thIntClosStab}. 
\end{Proof}

Pour une \dem \covz, nous appliquerons la stratégie indiquée auparavant et démontrerons  d'abord la version \und{\sdzz} (notez que si l'on était en \clama cette deuxième version serait identique à version III-8.12!).

Nous commençons par un cas particulier, qui est une légère \gnn
du lemme \ref{lemPolEnt}.

%:     Lemma{lemAnorAXnor1}
\begin{lemma} \label{lemAnorAXnor1}
Si $\gA$ est normal \sdzz, et si $q\in\AX$ est entier sur $\beta\gA$ ($\beta\in\gA$)
alors $q\in {\beta}\AX$. 
\end{lemma}
%--------- fin lemma ---------------------------------------------- 
%
\begin{proof} Conséquence \imde de l'exercice \ref{exoEntierSurIX}.
\\
Voici une \dem alternative. Dans les deux cas on utilise le \tho de Kronecker.
\\
On considère l'anneau
 $\gB=\gA[\fraC 1 {\beta}]$ et  l'image $\gC$ de $\gA$ dans $\gB$. 
 \\
On a $\gC\simeq\gA/(0:\beta)$.
L'\eltz~$\fraC {q(X)} \beta$ de~$\BX$ est entier sur~$\gC[X]$,
donc chaque \coe  $\fraC {q_k}\beta$ de $\fraC {q(X)} \beta$ est entier sur~$\gC$ (lemme \ref{lemPolEnt}). 
\\
En multipliant la \rdi par la puissance convenable de
$\beta$, on obtient dans $\gC$ une
\egt $U_k({q_k},\beta)=0$ où $U_k(Y,B)$ est un \pog unitaire en $Y$.
On peut lire les \coes de $U_k$ dans~$\gA$, et l'\egt dans $\gC$
 donne une \egt $\beta U_k({q_k},\beta)=0$  dans~$\gA$. 
 \\
 Puisque $\gA$ est \sdz
on a

\snic{\beta=0\; \hbox{  ou }\; U_k({q_k},\beta)=0.}

\snii
Si $\beta=0$, $q$ est nilpotent donc nul, et $q\in\gen{\beta}$.\\
Dans le deuxième cas, puisque $\gA$ est normal et $q_k$ entier sur $\gen{\beta}$,
on \hbox{a $q_k=\beta v_k$} pour un certain $v_k$. En bref
on a un \pol $v$ tel que $q=\beta v$.
\end{proof}

Montrons maintenant la version \und{\sdzz} du \thref{thIntClosStab}. 
%:     Lemma{lemAnorAXnor2}
\begin{lemma} \label{lemAnorAXnor2}
Si $\gA$ est normal \sdzz, il en est 
de même pour $\AX$.
\end{lemma}
%--------- fin lemma ---------------------------------------------- 

%
\begin{proof} On a déjà montré (lemme \ref{lemSDZAAX}) que $\AX$ est \sdzz.
\\ 
 On considère $a$, $b\in\AX$ de degrés formels $\leq  m$, et l'on suppose que~$a$ est entier sur l'\idz~\hbox{${b}\AX$} (cela correspond dans la \dem du \thref{thIntClosStab} à la fraction~\hbox{$f= \fraC a b$} entière sur $\AX$). On doit montrer \hbox{que $a\in {b}\AX$}.
 On fait une \recu sur $m$. L'initialisation avec $m=0$ est claire.
 Supposons $m\geq 1$.
 \\
Notons $b_0=b(0)$ et $a_0=a(0)$. Il est clair que $a_0$ est entier sur $\gen{b_0}$ dans $\gA$. 
Puisque $\gA$ est normal, on a $a_0=b_0c$ pour un $c\in\gA$. 
D'après le lemme \ref{lem2ZintClos} on a $b_0R_{m+1}=0$. Donc $b_0=0$ ou $R_{m+1}=0$.
 
\hspace*{0em}--- Si $b_0=0$, on a $a_0=0$, $b=XB$ et $a=XA$ avec $A,B\in\AX$.
\\
Puisque $a$ est entier sur $b$, $A$ est entier sur $B$: on a une \rdi $a^{\ell}=\sum_{j=1}^{\ell} c_j b^{j} a^{\ell-j}$, et en divisant par $X^{\ell}$, on obtient l'\egt $A^{\ell}=\sum_{j=1}^{\ell} c_j B^{j} A^{\ell-j}$.
On peut alors appliquer l'\hdr avec $A$ et $B$.
 
\hspace*{0em}--- Si $R_{m+1}=0$, on a une \egt \fbox{$b_0^{m} a=b {q}$}
  dans l'anneau  $\gA[X]$ avec $q$ de degré formel $m$.
  \\
On pose $\beta=b_0^{m}$, ce qui donne \fbox{$\beta a(X)=b(X)q(X)$}, et l'on considère une \rdi de $a$ sur $\gen{b}$ dans $\AX$

\snic{a^{\ell}-\som_{j=1}^{\ell} c_j b^{j} a^{\ell-j}=0 \quad (c_j\in\AX).}

\snii
On multiplie par $\beta^{\ell}$ et l'on remplace $\beta a$ par $bq$, on obtient

\snic{
 b^{\ell}\big(q^{\ell}-\som_{j=1}^{\ell} c_j   \beta^{j} q^{\ell-j}\big)=0.}

Puisque $\AX$ est \sdzz, on a 

\snic{b=0\; \hbox{ ou }\;
 q^{\ell}-\som_{j=1}^{\ell} c_j   \beta^{j} q^{\ell-j}=0.}

\snii
Dans le premier cas $a$ est nilpotent donc nul, et $a\in\gen{b}$. 
\\
Dans le second cas, 
$q$ est entier sur $\gen{\beta}$. Le lemme \ref{lemAnorAXnor1} nous donne un $v$ tel que $q=\beta v$.
D'où finalement une \egt
$\beta(a-bv)=0$, ce qui permet de conclure car  $a=bv$ ou $\beta=0$
(cas déjà traité: $b_0=0$). 
\end{proof}

\begin{Proof}{\Demo \cov du \thref{thAnormalAXnormal}. }
On prend $a(X)$ et $b(X)$ de degrés formels $\leq m$.
Puisque $\gA$ est normal, il est \lsdzz.
On remplace chaque disjonction qui apparaît dans le raisonnement 
précédent lorsqu'un produit est nul
par deux \lons \come de l'\gui{l'anneau en cours} (au départ il s'agit de $\gA$) dans chacune desquelles la \dem peut se poursuivre. 
\\
\`A la fin, on a \hbox{des $u_i$} \com dans $\gA$; et dans chaque $\gA[\fraC 1{u_i}][X]$ on a une \egt   $a(X)=b(X)w_i(X)$ avec $w_i$ de degré formel~$m$. Ainsi le \sli
en les \coes de $w(X)$ ($\deg(w)\leq m$) qui signifie $a(X)=b(X)w(X)$ admet une solution locale. On conclut par le \plg de base.
 \end{Proof}

Notez que la \prco construit un \pol $w(X)$ vérifiant l'\egt $a(X)=w(X)b(X)$
sans faire appel à aucune hypothèse du type \gui{l'anneau $\AX$ possède un test de \dvez}.
Cette \dem ne se contente donc pas de dire
abstraitement \gui{cette chose est vraie} 
(à savoir $b$ divise $a$ dans $\AX$).
La \dem classique, elle, certifie la vérité de la conclusion
uniquement en un sens affaibli, puisque la construction du \pol $w$ n'y est
indiquée en aucune manière.

Nous devons noter que, contrairement à de très nombreuses \dems dans cet ouvrage, qui sont directement simples et élégantes sous forme \algqz, la \prco demande ici un effort \sul non négligeable par rapport à la \dem classique. En particulier nous avons d\^u faire appel au  \nst formel \ref{thNSTsurZ}.
% 

%-% ENTRE NOUS
\entrenous{Ce qui nous manque peut-être encore dans ce chapitre.

1) parler de la \fcn d'un \idep (\tfz) dans une extension finie normale
d'un \adpc

2) une belle section sur les courbes en affine
}
%-% Fin ENTRENOUS

%%%%%%%%%%%%%%%%%%%%%%%%%%%%%%%%%%%%%%%%%%%%%%%%%%%%%%%%%%%%%%%%%%%%%%%%%%%
%:section: Exercices
%\newpage	   
\Exercices{

%--- Exercise{exoDetTrick1}--------- 
\begin {exercise} \label {exoDetTrick1}
       (Encore un \gui {determinant trick})\\
{\rm 
Soit $E$ un $\gA$-module fidèle engendré par $n$ \elts et
$\fa \subseteq \fb$ deux idéaux de $\gA$ vérifiant $\fa E = \fb E$.
Montrer que $\fa\fb^{n-1} = \fb^n$.
}
\end {exercise}
%--- end -exercise-----------------------------------------

%--- Exercise{exoMatLoc22}-------------
\begin{exercise}\label{exoMatLoc22}
{(Matrices de \lon principale $2 \times 2$)}   
{\rm Soient $x, y \in\gA$.

\emph {1.} 
Montrer que l'\id $\gen {x,y}$ est \lop \ssi il existe une matrice~$B \in \MM_2(\gA)$ de
trace $1$ vérifiant $[\,x\; y\,]B = 0$; dans ce cas, $A = \wi B$ est une
\mlp pour $(x, y)$.

\emph {2.}
Soit $z \in \gA$; on suppose qu'il existe un \id $\fb$ tel que $\gen
{x,y}\fb = \gen {z}$.  Montrer qu'il existe $B \in \MM_2(\gA)$ tel que $z[x
\; y\,]B = 0$ et $z\big(1 - \Tr(B)\big) = 0$.

\emph {3.}
Déduire des questions précédentes une autre preuve 
du lemme \ref{lemleszlop}.

}

\end {exercise}
%--- end -exercise-----------------------------------------

%--- Exercise{exoPrufNagata}-------------
\begin{exercise}
\label{exoPrufNagata}\label{exoWhenNagataRingIsArith}
($\gA$ \ari $\Leftrightarrow$  $\gA(X)$ Bézout)  {\rm Voir aussi l'exercice~\ref{exoBézoutKdim1TransfertArX}.
\\
 Soit $\gA$ un anneau  et $\gA(X)$ le localisé de Nagata. 

 \emph{1.}
Montrer que pour $a,b\in\gA$, $a\mid b$
dans $\gA$ \ssi $a \mid b$ dans $\gA(X)$.

 \emph{2.} Si $\gA$ est un \anar et $f\in\AX$,
on a  dans $\gA(X)$

\snic{\gen{f}=\rc_\gA(f)\gA(X).}

%\sni
Montrer aussi que $\gA(X)$ est un anneau de Bézout.

 \emph {3.}
Soient $x$, $y \in \gA$. Montrer que si $\gen{x,y}$ est \lop
dans  $\gA(X)$, il est \lop dans $\gA$ (utiliser l'exercice \ref{exoMatLoc22}).
En particulier, si $\gA(X)$ est \ariz, il en est de même de $\gA$.
A fortiori, si $\ArX$ est \ari (notation \ref{notaA<X>}), il en est de même de $\gA$.

\emph {4.}
Conclure. 

 NB: sur l'anneau $\gA(X)$ voir le fait \ref{factLocNagata} 
et l'exercice~\ref{exoNagatalocal}. 
}
\end{exercise}
%--- end -exercise-----------------------------------------

%--- Exercise{exoanars}--------- 
\begin{exercise} 
\label{exoanars} (Quelques autres \prts \caras des \anarsz)
{\rm  
Pour un anneau $\gA$, \propeq
%-----------------begin item------------------
\begin{itemize}
\item [\phantom{.2}$(1)$] $\gA$ est un \anarz. 
\item [$(2.1)$] Pour tous \ids  $\fa$, $\fb$ et $\fc$ on a 
$\fa\cap(\fb+\fc)=(\fa\cap \fb)+(\fa\cap \fc).$ 
\item [$(2.2)$] Même chose en se limitant aux \idpsz. 
\item [$(2.3)$] Même chose  en se limitant au cas $\fb=\gen{x}$, $\fc=\gen{y}$ 
et $\fa=\gen{x+y}.$
\item [$(3.1)$] Pour tous \ids  $\fa$, $\fb$ et $\fc$ on a 
$\fa+(\fb\cap \fc)=(\fa+\fb)\cap (\fa+\fc).$
\item [$(3.2)$] Même chose en se limitant aux \idpsz. 
\item [$(3.3)$] Même chose  en se limitant au cas $\fa=\gen{x}$, $\fb=\gen{y}$ 
et $\fc=\gen{x+y}.$  
\item [$(4.1)$] Pour tous \itfs   $\fa$, $\fb$ et $\fc$ on a 
$(\fb+\fc):\fa=(\fb:\fa)+(\fc:\fa).$
\item [$(4.2)$] Même chose avec  $\fb$ et $\fc$  \idps et $\fa=\fb+\fc$.
\item [$(5.1)$] Pour tout idéal  $\fa$ et tous \itfs   $\fb$ et $\fc$ on a l'\egt
\\
\centerline{$\fa:(\fb\cap \fc)=(\fa:\fb)+(\fa:\fc).$}
\item [$(5.2)$] Même chose avec  $\fb$ et $\fc$  \idps et $\fa=\fb\cap \fc$. 
\end{itemize}
%-----------------end item------------------
} 
Indication: pour démontrer que les conditions sont \ncrs
on utilise la méthode \gnle expliquée \paref{MetgenAnar}.
\end{exercise}
%--- end-exercise-----------------------------------------

%--- Exercise{exoNormalClass}--------
\begin{exercise} 
\label{exoNormalClass} 
{\rm Démontrer en \clama qu'un anneau est normal \ssi il devient normal 
lorsque l'on  localise en un \idep arbitraire (rappelons que dans le cas intègre, 
normal signifie \icl dans son corps de fractions).
}
\end{exercise}
%--- end-exercise-----------------------------------------

%--- Exercise{FermetureAlgZariski}-------------
\begin{exercise}\label{FermetureAlgZariski}
 {(Fermeture \agqz: un \tho dû à Zariski)}\\
{\rm  
Soient $\gK \subseteq \gL$ deux \cdisz, $\gK'$ la fermeture \agq
de $\gK$ dans $\gL$. Alors la fermeture \agq de
$\gK(\Xn)$ dans $\gL(\Xn)$ est $\gK'(\Xn)$; résultat analogue
si l'on remplace fermeture \agq par fermeture \agq\spbz.
}

\end {exercise}
%--- end -exercise-----------------------------------------

%--- Exercise{exoNotAbsIntClos}-------------
\begin{exercise}\label{exoNotAbsIntClos}
{(Un défaut d'intégralité par extension des scalaires)}\\
{\rm  
Soit $\gk$ un \cdi de \cara $p \ge 3$, $a \in \gk$ et $f = Y^2 - f(X)
\in \gk[X,Y]$ avec $f(X) = X^p - a$.

\emph {1.}
Montrer que $Y^2 - f(X)$ est \ix{absolument \irdz}, \cad que pour tout surcorps $\gk'
\supseteq \gk$, le \pol $Y^2 - f(X)$ est \ird dans $\gk'[X,Y]$.

 On note $\gk[x,y] = \aqo {\gk[X,Y]}{Y^2 - f(X)}$ et $\gk(x,y) =
\Frac(\gk[x,y])$.

\emph {2.}
Montrer que $\gk$ est \agqt fermé dans $\gk(x,y)$ et
que pour toute extension \agq $\gk'$ de $\gk$, on a
$\gk' \te_\gk \gk(x,y) = \gk'(x,y)$.

\emph {3.}
On suppose que $a \notin \gk^p$. Montrer que $\gk[x,y]$
est \icl et que~$\gk(x,y)$ n'est pas un corps
de fractions rationnelles à une \idtr sur $\gk$.

\emph {4.}
On suppose $a \in \gk^p$ (par exemple $a = 0$). Montrer que $\gk[x,y]$ n'est
pas \icl et expliciter $t \in \gk(x,y)$ tel que $\gk(x,y) = \gk(t)$.

}

\end {exercise}
%--- end -exercise-----------------------------------------

%\newpage
%--- Exercise{exoAnneauOuvertP1}-------------
\begin{exercise}
\label{exoAnneauOuvertP1} (L'anneau
des fonctions sur la droite projective privée d'un nombre fini
de points)\\
{\rm 
On utilise dans cet exercice de manière informelle les notions de schéma affine et de droite projective qui ont déjà été discutées
dans les sections \ref{sec1Apf} et~\ref{secGrassman} (voir notamment les pages \pageref{subsecNstMorphismes} à \pageref{subsubsecTanFonct}).\\
Si $\gk$ est un \cdiz, la \klg  {des fonctions \polles
définies sur la droite affine $\AA^1(\gk)$} est $\gk[t]$. 
Si l'on pense à  $\AA^1(\gk) \cup \{\infty\}  = \PP^1(\gk) $, 
les \elts de~$\gk[t]$ sont alors les fractions rationnelles sur $\PP^1(\gk)$ qui sont \gui{définies partout, sauf peut-être en $\infty$}.
\\
Soient $t_1$, \ldots, $t_r$, des points de cette droite affine (on peut avoir $r=0$). 
\\
On munit $\AA^1(\gk) \setminus \{t_1,
\ldots, t_r\}$ (droite affine privée de $r$ points) d'une structure de \vrt
affine en forçant l'inversibilité du produit des $t-t_i$, i.e. en définissant

\snic {
\gB = \gk\big[t, (t-t_1)^{-1}, \ldots, (t-t_r)^{-1}\big] 
 \simeq 
\aqo{\gk[t,x]} {F(t,x)},
}

%\sni
avec $F(t,x) = (t-t_1) \cdots (t-t_r)\cdot x - 1$. Cette
\klg $\gB$  apparaît alors comme l'\alg des fractions rationnelles sur
$\PP^1(\gk)$ définies
partout sauf aux points~$\infty$ et~$t_i$.  C'est un \acl 
et même un anneau de Bézout (en effet, c'est un localisé de $\gk[t]$).
\\
De manière analogue, pour $n$ points  $t_1$, \ldots, $t_n$ de la droite affine
(avec $n \ge 1$ cette fois), on peut considérer la \klg 

\snic {
\gA = \gk\big[(t-t_1)^{-1}, \ldots, (t-t_n)^{-1}\big]\subseteq \gk(t).
}

%\sni
Cet anneau $\gA$ est un localisé de $\gk[(t-t_1)^{-1}]$ (qui est isomorphe à $\kX$) puisqu'en
\linebreak 
posant $v = (t-t_1)^{-1}$, on a $t-t_i = \big((t_1-t_i)v + 1\big)/v$.
Ainsi, 
$$\preskip.8pt
{\gA=\gk\big[v,\big((t_1-t_2)v+1\big)^{-1},\dots,\big((t_1-t_n)v+1\big)^{-1}\big]\subseteq \gk(v)=\gk(t)}.
\postskip2.5pt
$$
La \klg $\gA$ est donc un \acl (et même un anneau de Bézout). En notant $p(t) = (t-t_1) \cdots (t-t_n)$, on a aussi facilement l'\egt

\snic {
\gA = \gk[1/p, t/p, \ldots, t^{n-1}/p].
}

%\sni
La \klg $\gA$, constituée des fractions rationnelles $u/p^s$
avec $\deg(u) \le ns$, apparaît comme celle des fractions rationnelles
définies partout sur $\PP_1(\gk)$ (y~compris au point $t = \infty$) sauf
éventuellement aux points $t_i$.  En bref,
on peut convenir que~$\gA$ est la \klg des \gui{fonctions} définies
sur la droite projective privée des points $t_1$, \dots, $t_n$.

On étudie dans cet exercice un cas plus \gnl où $p$ est un
\polu de degré $n \ge 1$.

Soit $\gk$ un \cdiz, $p(t) = t^n + a_{n-1} t^{n-1} + \cdots + a_1 t + a_0\in\gk[t]$ $(n\geq1)$,
où~$t$ est une \idtrz. On pose $x_i = \fraC{t^i}p$.
 \\
Montrer que la \cli de $\gk[x_0]$ dans $\gk(t)$ est la \klg

\snic{\gA = \gk[x_0, \ldots, x_{n-1}] = 
\sotq {\fraC u {p^s}} {s\in\NN,\,u \in \gk[t],\  \deg(u) \le ns}.
}

%\sni
En outre, $\Frac(\gA) = \gk(t)$. 
}
\end{exercise}
%--- end -exercise-----------------------------------------

%--- Exercise{exoPresentationAlgOuvertP1}-------------
\begin{exercise}\label{exoPresentationAlgOuvertP1}
 {(Une présentation de l'\alg des fonctions
sur la droite projective privée d'un nombre fini de points)}\\
{\rm  
Le contexte est celui de l'exercice \ref {exoAnneauOuvertP1}, mais cette fois-ci $\gk$
est un anneau quelconque. On note $p = a_nt^n + \cdots + a_1t + a_0 \in
\gk[t]$ un \polu ($a_n = 1$) et

\snic {
\gA = \gk\big[\fraC1 p, \fraC t p, \ldots, \fraC{t^{n-1}}p\big].
}

%\sni

On pose $x_i = \fraC{t^i}p$ pour $i \in \lrb{0..n-1}$. On peut alors écrire
$\gA = \gk[\uX]/\fa$ \hbox{où $(\uX) = (X_0, \ldots, X_{n-1})$} et $\fa$ est l'\id
des relateurs pour $(x_0, \ldots, x_{n-1})$.  Il va se révéler commode de
définir $x_n$ par $x_n = \fraC{t^n}p$; on a donc $x_j = x_0t^j$ et
$$\preskip.5pt\ndsp
\sum_{i=0}^n a_i x_i = 1  \quad \hbox {ou encore} \quad
x_n = 1 - \sum_{i=0}^{n-1} a_i x_i.
\postskip2.5pt
$$
L'\egt de droite prouve que $x_n \in \gA$.

\emph {1.}
Vérifier que la famille $R$ suivante est constituée de relateurs entre
les $x_j$.%, $j \in \lrb{0..n}$:

\snic {
R \ : \quad
x_ix_j = x_kx_\ell  \qquad \hbox {pour $i+j = k+\ell$}, \qquad
0 \le i, j, k, \ell \le n.
}

%\sni
On définit la sous-famille $R_{\rm min}$, constituée de ${n(n-1) \over 2}$
relateurs.

\snic {
R_{\rm min}\ : \quad x_i x_j = x_{i-1}x_{j+1}, \qquad 1 \le i \le j \le n-1.
}

%\sni
\emph {2.}
Montrer que la famille $R_{\rm min}$ (donc $R$ aussi) engendre l'\id des
relateurs entre les $x_i$ pour $i \in \lrb{0..n-1}$. 
En d'autres terms, si nous notons $\varphi : \gk[\uX] \to \gk[t,1/p]$ le morphisme défini par
$X_i \mapsto x_i$ pour $i \in \lrb {0..n-1}$, cela signifie que 

\snic {
\Ker\varphi=\geN{X_iX_j - X_{i-1}X_{j+1}, \, 1 \le i \le j \le n-1} \;\; (\hbox{où }  X_n := 1 -
\sum_{i=0}^{n-1} a_iX_i)
.}

%\sni
On pourra faire intervenir le \kmo
$\gk[X_0] \oplus \gk[X_0]X_1 \oplus \cdots \oplus \gk[X_0]X_{n-1}$.

}

\end {exercise}
%--- end -exercise-----------------------------------------

%--- Exercise{exoEntiersEmmanuel}-------------
\begin{exercise}\label{exoEntiersEmmanuel} {(Les entiers d'Emmanuel)}
\\
{\rm
Donner une preuve directe du point \emph{1} du lemme  \ref{lemEmmanuel} sans 
utiliser le \tho de \KROz.}
\end {exercise}
%--- end -exercise-----------------------------------------

%--- Exercise{exoKroneckerTheorem}-------------
\begin{exercise}\label{exoKroneckerTheorem} 
{(Une autre \dem du \tho de \KROz)}\\
{\rm 
On considère les \pols
\[ 
\begin{array}{ccc} 
  f(T) = a_0T^n + \cdots + a_n, \quad g(T) = b_0T^m + \cdots + b_m \;\;\hbox{et}  \\[1mm] 
h(T)
= f(T)g(T) = c_0T^{n+m} + \cdots + c_{n+m}.
 \end{array}
\] Le
\tho de Kronecker \ref{thKro} affirme que chaque produit $a_ib_j$ est
entier sur l'anneau $\gA = \ZZ[c_0, \ldots, c_{n+m}]$.  
\\
Il suffit de traiter le cas où 
les $a_i$ et $b_j$ sont des \idtrsz. Alors dans un anneau contenant
$\ZZ[a_0,\dots,a_n,b_0,\dots,b_m]$ on~a:

\snic {
f(T) = a_0(T-x_1) \cdots (T-x_n), \qquad
g(T) = b_0(T-y_1) \cdots (T-y_m)
.}

%\sni
\emph {1.}
En utilisant les entiers d'Emmanuel (lemme \ref{lemEmmanuel},
avec la \dem donnée dans l'exercice \ref{exoEntiersEmmanuel},
indépendante du \tho de Kronecker), montrer que
pour tous~$I \subseteq \lrb{1..n}$, $J \subseteq \lrb{1..m}$, le produit
$a_0b_0 \prod_{i\in I}x_i \prod_{j\in J}y_j$ est entier sur $\gA$. 

\emph {2.}
Conclure.

}
\end {exercise}
%--- end -exercise-----------------------------------------

%%--- Exercise{exofactAdpIntExt}-------------
%\begin{exercise}
%\label{exofactAdpIntExt}
%{\rm  Montrer la variante suivante du lemme \ref{factAdpIntExt}: 
%Si $\gB$ est un anneau normal, et une extension entière
%d'un \adpc $\gA$, alors $\gB$ est un \adpcz.
% 
%}
%\end{exercise}
%%--- end -exercise-----------------------------------------

%--- Exercise{exoAnneauCoinceBézout}-------------
\begin{exercise}\label{exoAnneauCoinceBézout} 
{(Anneau intermédiaire $\gA \subseteq \gB \subseteq\Frac(\gA)$, cas Bézout)}
 \\
{\rm  
Soit $\gA$ un anneau de Bézout intègre de corps de fractions $\gK$ et $\gB$
un anneau intermédiaire $\gA \subseteq \gB \subseteq \gK$.  Montrer que
$\gB$ est un localisé de $\gA$ (donc un anneau de Bézout).
}
\end {exercise}
%--- end -exercise-----------------------------------------

%--- Exercise{exoGrellNoether}-------------
\begin{exercise}\label{exoGrellNoether}
{(Anneau intermédiaire, cas Prüfer)}\\
{\rm  
Dans cet exercice on généralise le résultat de l'exercice \ref{exoAnneauCoinceBézout}
au cas où $\gA$ est un \ddp et l'on
précise le \thref{thSurAdp}. Il s'agit donc d'une variation autour du \tho de 
Grell-\Noe (\paref{thGrellNoether}).

\emph {1.}
Soit $x \in \gK=\Frac\gA$. 
\begin{itemize}
\item [\emph {a.}]
Montrer qu'il existe $s \in \Reg(\gA)$  tel que $sx \in \gA$ et $1-s \in \gA
x$.

\item [\emph {b.}]
Soit $t \in \gA$ tel que $tx = 1-s$. Pour tout anneau  $\gA'$
intermédiaire entre~$\gA$ et~$\gK$, montrer que $\gA'[x] = \gA'_s \cap \gA'_t$.
En particulier, $\gA[x] = \gA_s \cap \gA_t$. En conséquence, $\gA[x]$ est
\iclz, et c'est un \adpz.
\end{itemize}

\emph {2.}
Montrer que toute sous-\Alg $\gB$ de $\gK$ de type fini est intersection d'un
nombre fini de localisés  de $\gA$ de la forme $\gA_s$
avec $s \in \gA$. En conséquence, $\gB$ est \iclz, et c'est un \adpz.

\emph {3.}
En déduire que tout anneau intermédiaire entre $\gA$ et $\gK$ est
de Prüfer.

\emph {4.}
Donner un exemple d'anneau \icl $\gA$,  avec
un anneau  $\gB$ intermédiaire  entre~$\gA$ et~$\Frac(\gA)$ qui n'est pas  \icl (en
particulier, $\gB$ n'est pas un localisé de $\gA$).

}

\end {exercise}
%--- end -exercise-----------------------------------------

%--- Exercise{exoPrimitivementAlg}-------------
\begin{exercise}\label{exoPrimitivementAlg}
 {(\^Etre primitivement \agqz)}\\
{\rm  
Soient $\gA = \aqo{\ZZ[A,B,U,V]}{AU+BV-1} = \ZZ[a,b,u,v]$ et
$\gB = \gA[1/b]$. \\
On pose $x = a/b$.
Montrer que $x$ est primitivement \agq sur $\gA$,
mais que~$y=2x$ ne l'est pas.

}

\end {exercise}
%--- end -exercise-----------------------------------------

%--- Exercise{exocaracPruferC}-------------
\begin{exercise}\label{exocaracPruferC}  {(\Carns des \adpcsz, 1)}
\\
{\rm  
Soit $\gA$  \qi et $\gK=\Frac\gA$.
\Propeq

%-----------------begin enum------------------
\vspace{-1pt}
\begin{itemize}
\item [\emph {1.}]
$\gA$ est de Prüfer.

\item [\emph {2.}]
$\gA$ est normal et $x \in \gA[x^2]$ pour tout $x \in \gK$.

\item [\emph {3.}]
Tout anneau  $\gA[y]$ où $y\in\gK$ est normal.

\item [\emph {4.}]
Tout anneau intermédiaire entre $\gA$ et $\gK$ est normal.

\item [\emph {5.}]
$\gA$ est normal et $x \in \gA+ x^2 \gA$ pour tout $x \in \gK$.
\end{itemize}

}
\end {exercise}
%--- end -exercise-----------------------------------------

%--- Exercise{exocaracPruCoh}-------------
\begin{exercise}
\label{exocaracPruCoh} (\Carns des \adpcsz, 2)\\
{\rm   Pour un anneau $\gA$ \qiz,
\propeq
%-----------------begin item------------------
\vspace{-1pt}
\begin{itemize}
\item [\emph {1.}] $\gA$ est un \adpz.
\item [\emph {2.}]  Tout \itf contenant un \elt \ndz est \ivz.
\item [\emph {3.}]  Tout \id $\fa=\gen{x_1,x_2}$ avec $x_1$, $x_2\in\Reg(\gA)$ est
\ivz.
\item [\emph {4.}]  Pour tous $a$, $b\in\gA,$ on a:
$\gen{a,b}^2=\gen{a^2,b^2}=\gen{a^2+b^2,ab}.$
\item [\emph {5.}]  Pour tous $f$, $g\in\gA[X]$, on a: $\rc(f)\rc(g)=\rc(fg)$.
\end{itemize}
 
}
\end{exercise}
%--- end -exercise-----------------------------------------

%--- Exercise{exoAnarlgb}-------------
\begin{exercise}\label{exoAnarlgb} (Une \gnn de la proposition \ref{propAriCohZed})\\
{\rm  Soit $\gA$ un \adpc  \lgb (par exemple \rdt \zedz). 
\vspace{-1pt}
\begin{itemize}
\item [\emph {1.}]
 Toute
matrice est \eqve à une matrice en forme de Smith (i.e. $\gA$ est un anneau de Smith).
\item  [\emph {2.}] Tout \Amo \pf est caractérisé, à \iso près, par ses \idfsz.  En fait il est isomorphe à une somme directe de
modules monogènes~$\gA\sur{\fa_k}$ avec
des \idps $\fa_1\subseteq\cdots\subseteq\fa_n$ ($n\geq0$).
\\ 
NB: on peut naturellement en déduire une \gnn 
analogue du corolaire~\ref{corpropAriCohZed}.
\end{itemize} 
}
\end{exercise}
%--- end -exercise-----------------------------------------

%--- Exercise{exoReductionIdeal}-------------
\begin{exercise}\label {exoReductionIdeal} 
  {(Idéal réduction d'un autre idéal)}
{\rm  
\vspace{-1pt}
\begin{itemize}
\item [\emph {1.}] Soient $E$ un $\gA$-module engendré par $n$ \eltsz, $b \in \gA$
et $\fa$ un idéal tels que~$bE \subseteq \fa E$.
Montrer qu'il existe $d= b^n + a_{1} b^{n-1} + \cdots + a_{n-1}b + a_n$, 
avec les~$a_i \in \fa^i$, qui annule $E$.
\end{itemize}
On dit qu'un idéal $\fa$ est une \emph{réduction} d'un idéal $\fb$
si $\fa \subseteq \fb$ et si $\fb^{r+1} = \fa\fb^r$ pour un certain
exposant $r$ (c'est alors vrai pour tous les exposants plus grands).
\begin {itemize} 
\item [\emph {2.}]
Soient $f$, $g \in \gA[\uX]$. Vérifier que $\rc(fg)$ est une
réduction de $\rc(f)\rc(g)$.
\item [\emph {3.}]
Dans $\gA[X,Y]$, montrer que $\fa = \gen {X^2, Y^2}$ est
une réduction de $\fb = \gen {X,Y}^2$.\\
Montrer 
que $\fa_1 = \gen {X^7, Y^7}$
et $\fa_2 = \gen {X^7, X^6Y + Y^7}$ sont des réductions
de \linebreak 
l'\idz~$\fb' = \gen {X^7, X^6Y, X^2Y^5, Y^7}$. Donner les
plus petits exposants possibles.
\item [\emph {4.}]
Soit $\fa \subseteq \fb$ deux idéaux avec $\fb$ de type fini.
Montrer que $\fa$ est une réduction de~$\fb$ \ssi $\Icl(\fa) = \Icl(\fb)$.
\end {itemize}
}
\end {exercise}
%--- end -exercise-----------------------------------------

%:--- Exercise{exolemNormalIcl}-------------
\begin{exercise}  \label{exolemNormalIcl} {(Anneau normal \qiz)}\\
{\rm  
Voici une légère \gnn du fait \ref{lemNormalIcl}.
D'après le \pb \ref{exoAnneauNoetherienReduit} l'hypothèse est satisfaite pour les anneaux \noes
\cohs réduits \fdis (en \clama ce sont les anneaux \noes réduits).

On considère un anneau $\gA$ réduit. On suppose que son anneau total des fractions est \zedz. 

\emph{1.} Si $\gA$ est normal, il est \qiz.

\emph{2.} L'anneau $\gA$ est normal \ssi il est \icl dans $\Frac\gA$.

}

\end{exercise}

%:--- Exercise{exoEntierSurIX}-------------
\begin{exercise}  \label{exoEntierSurIX} {(\Pol entier sur $\fa[X]$)}

{\rm  
Soient $\gA \subseteq\gB$ deux anneaux, $\fa$ un \id de $\gA$ et $\fa[X]$
l'\id de $\gA[X]$ constitué des \pols à \coes dans $\fa$. Pour
$F \in \gB[X]$, montrer que $F$ est entier \hbox{sur $\fa[X]$} \ssi chaque \coe de
$F$ est entier sur $\fa$.

}

\end{exercise}

%:--- Exercise{exosdirindec}-------------
\begin{exercise}  \label{exosdirindec} {(Modules indécomposables)}

{\rm  
On dit qu'un module $M$ est \emph{indécomposable} si les seuls sous-modules facteurs directs de~$M$ sont $0$ et $M$.
Le but de l'exercice est de démontrer que sur un \dDk à \facz, tout \mpf est somme directe d'un nombre fini de modules indécomposables, cette \dcn étant unique à l'ordre des termes près lorsque le module est de torsion.

\emph{1.} Soit $\gA$ un anneau et $\fa$ un \idz. Si le \Amo $M=\gA/\fa$ est somme directe de deux sous-modules $N$ et $P$ on a $N=\fb/\fa$, $P=\fc/\fa$
avec~$\fb\supseteq \fa$ et~$\fc\supseteq \fa$  \comz. \Prmtz, $\fb=\gen{b}+\fa$, $\fc=\gen{c}+\fa$, où $b$ et $c$ sont des \idms \cops modulo~$\fa$. 

\emph{2.} Soit $\gZ$ un \dDkz. 
 
\emph{2a.} Montrer qu'un \mrc 1 est indécomposable.
 
\emph{2b.} Montrer qu'un module cyclique  $\gZ/\fa$ avec $\fa$ \tfz, $\neq \gen{0},\,\gen{1}$ est indécomposable \ssi $\fa=\fp^{m}$ pour un \idema $\fp$ et un~$m\geq 1$.
 
\emph{2c.} En déduire que si $\gZ$ est à \facz, tout \mpf est somme directe d'un nombre fini de modules indécomposables.

\emph{3.} Lorsque le module est de torsion, montrer l'unicité de la  \dcn 
en un sens à préciser.

}

\end{exercise}

%%%%%%%%%%%%%%%%%%%%%%%%%%%%%%%%%%%%%%%%%%%%%%%%%%%%%%%%%%%%%%%%%%%%%%%%%%%
%:  pbs
%%%%%%%%%%%%%%%%%%%%%%%%%%%%%%%%%%%%%%%%%%%%%%%%%%%%%%%%%%%%%%%%%%%%%%%%%%%
%--- problem{exoArithInvariantRing}-------------
\begin{problem}\label{exoArithInvariantRing} \\
{(Sous-anneau d'invariants par un groupe fini et
\crc \ariz)}\\
{\rm  
 NB: voir aussi le \pb \ref{exoGaloisNormIdeal}.\\
\emph {1.}
Si $\gA$ est un \anorz, tout \id\lop %$\fa$ de $\gA$ 
est \iclz. En conséquence, si $f$, $g \in \AX$ avec $\rc(fg)$
\lopz, alors $\rc(f)\rc(g) = \rc(fg)$.

\emph {2.}
On suppose $\gA$ normal et $\gB\supseteq\gA$ entière sur $\gA$.
Si $\fa$ est un \itf de~$\gA$ \lopz, alors $\fa\gB\cap\gA = \fa$.

\emph {3.}
Soit $(\gB, \gA, G)$ où $G \subseteq \Aut(\gB)$ est un groupe fini
et $\gA =\Fix_\gB(G)= \gB^G$. Si $\fb$ est un \id de $\gB$, on note $\rN'_G(\fb) =
\prod_{\sigma \in G} \sigma(\fb)$ (c'est un idéal de $\gB$)\linebreak 
 et $\rN_G(\fb) =
\gA \cap \rN'_G(\fb)$. 

 On suppose $\gB$ normal et $\gA$ de Prüfer (donc normal). 
\begin {itemize}
\item [{\it a.}]
Pour $b \in \gB$, vérifier que $\rN'_G(b\gB)=\rN_G(b)\gB$  
et
$\rN_G(b\gB)=\rN_G(b)\gA$.

\item [{\it b.}]
Si $\fb$ est un \itf de $\gB$, montrer que $\rN_G(\fb)$ est un \itf de
$\gA$ et que
$\rN'_G(\fb) = \rN_G(\fb)\gB$. 
On pourra écrire $\fb = \gen {b_1, \ldots, b_n}$, introduire $n$
\idtrs $\uX = (X_1, \ldots, X_n)$ et le \pol normique $h(\uX)$:

\snic {
h(\uX) = \prod_{\sigma \in G} h_\sigma(\uX) 
\quad \hbox {avec} \quad
h_\sigma(\ux) = \sigma(b_1) X_1 + \cdots +  \sigma(b_n) X_n.
}

\item [{\it c.}]
Pour $\fb_1$, $\fb_2$ \itfs de $\gB$, on obtient $\rN_G(\fb_1\fb_2) = \rN_G(\fb_1) \rN_G(\fb_2)$.

\item [{\it d.}]
Un \itf $\fb$ de $\gB$ est \iv \ssi $\rN_G(\fb)$  est \iv dans $\gA$.
\end {itemize}

Note: on sait que $\gB$ est un \adp (\tho \ref{thExtEntPruf});
dans le cas où $\gB$ est intègre, la question \emph {3d}
en fournit une nouvelle preuve.

\emph {4.}
Soit $\gk$ un \cdi avec $2\in\gk\eti$, $f(X) \in \gk[X]$ un
\polu\spbz. Le \pol $Y^2 - f(X)
\in \gk[X,Y]$ est absolument \ird (voir l'exercice \ref{exoNotAbsIntClos}); 
on pose $\gk[x,y] =\aqo {\gk[X,Y]}{Y^2 - f(X)}$. 
Montrer que $\gk[x,y]$ est un \adpz.

}

\end {problem}
%--- end -problem-----------------------------------------

%--- problem{exoFullAffineMonoid}-------------
\begin{problem}\label{exoFullAffineMonoid}
 {(Sous-\mos pleins de $\NN^n$)}
\\
{\rm  
Soit $M\subseteq\NN^n$ un sous-\moz; pour  un anneau $\gk$, on
note $\gk[M]$ la \klg  du \mo $M$. C'est la sous-\klg de
$\gk[\NN^n] \simeq \gk[\ux] = \gk[\xn]$ engendrée par les monômes $\ux^m =
x_1^{m_1} \cdots x_n^{m_n}$ pour $m \in M$.
On dit que $M$ est un \emph {sous-\mo plein} de $\NN^n$ si pour
$m \in M$, $m' \in \NN^n$, on a $m+m' \in M \Rightarrow m' \in M$.

\emph {1.}
Le sous-groupe de $\ZZ^n$ engendré par $M$ est égal à $M-M$, 
et si $M$ est plein, alors $M = (M-M)\cap\NN^n$. Réciproquement,
si $L \subseteq \ZZ^n$ est un sous-groupe, alors le \mo $M = L\cap\NN^n$
est un sous-\mo plein de $\NN^n$.

\emph {2.}
Soit $M \subseteq \NN^n$ un sous-\mo plein et $\gk$ un \cdiz.
\begin {enumerate}
\item [\emph {a)}]
Soit $\gA=\gk[M] \subseteq \gB=\kux$. Montrer que si $a \in
\gA\setminus\{0\}$,  $b \in \gB$, et~$ab \in\gA $,  alors~$b \in\gA$.

\item [\emph {b)}]
Soient $\gA\subseteq\gB$ deux anneaux intègres vérifiant: si $a \in
\gA\setminus\{0\}$,  $b \in \gB$, et~$ab \in\gA $,  alors~$b \in\gA$.  

\begin {enumerate}
\item [\emph {i.}]
Montrer que $\gA = \gB\cap\Frac(\gA)$; en déduire que si $\gB$ est \iclz, il
en est de même de $\gA$.
\item [\emph {ii.}]
En particulier, si $M \subseteq \NN^n$ un sous-\mo plein, alors
$\gk[M]$ est \icl pour tout \cdi $\gk$.
\item [\emph {iii.}]
Plus \gnltz, si $\gB\subseteq\gC$ est intégralement
fermé dans $\gC$, alors $\gA$ est intégralement fermé
dans $\gC\cap\Frac(\gA)$.
\end {enumerate}
\end {enumerate}

\emph {3.}
Soit $M \subseteq \NN^n$ le sous-\mo des matrices magiques (voir l'exercice \ref
{exoMatMag3}); alors $\gk[M]$ est \icl pour tout \cdi $\gk$.

}

\end {problem}
%--- end -problem-----------------------------------------

%--- problem{exoBaseNormaleAlInfini}-------------
\begin{problem}\label{exoBaseNormaleAlInfini} 
 {(Base normale à l'infini)}\\
{\rm 
Un \adv intègre $\gB$ de corps de fractions $\gK$ est un \advd si $\gK\eti\!/{\gB\eti}\simeq \ZZ$ (\iso de groupes ordonnés).
On rappelle qu'une uniformisante est un \elt $b\in\gB$ tel que $v(b)=1$, \hbox{où
$v:\gK\eti\to\ZZ$} est l'application définie via l'\iso précédent.
Cette application $v$ s'appelle aussi une \emph{valuation du corps $\gK$}.
Tout \elt $z$ de $\gK\eti$ s'écrit alors $ub^{v(z)}$ pour \hbox{un $u\in\gB\eti$}.%
\index{valuation!d'un corps discret}

Soient $\gk$ un \cdiz, $t$ une \idtr sur $\gk$, $\gA =
\gk[t]$, $\gA_\infty = \gk[t^{-1}]_{\langle t^{-1} \rangle}$, et~$\gK = \Frac(\gA) =
\gk(t) = \Frac(\gA_\infty) = \gk(t^{-1})$.  Si $L$ est un \Kev de dimension
finie, on étudie dans ce \pb l'intersection d'un $\gA$-réseau de~$L$
et d'un $\gA_\infty$-réseau de~$L$ (cf. les \dfns question \emph {2}),
intersection qui est toujours un \kev de dimension finie.  Cette étude est,
dans la théorie des corps de fonctions \agqsz,
à la base de la détermination des espaces de Riemann-Roch, quand toutefois
certaines clôtures intégrales sont connues par des bases; comme
sous-produit, on obtient la détermination de la fermeture \agq de $\gk$ dans
une extension finie de $\gk(t)$.
\\
L'anneau $\gA_\infty$ est un \adv discrète; on note $v : \gK \to \ZZ \cup
\{\infty\}$ la valuation correspondante, définie par $v = -\deg_t$, et l'on
fixe $\pi = t^{-1}$ comme uniformisante.  Si $x = \tra {[x_1, \ldots, x_n]}$,
on pose $v(x) = \min_i v(x_i)$.  Ceci permet de définir une
\emph{réduction modulaire}

\snic{\gK^n
\setminus \{0\} \to \gk^n \setminus \{0\}$, $\quad x \mapsto \xi = \ov {x},}

%\sni
avec
$\xi_i = (x_i/\pi^{v(x)}) \mod \pi \in \gk$.  
\\
De manière \gnlez,
si $\gV$ est un \adv d'un \cdi $\gK$, de corps résiduel $\gk$,
on a une \emph{réduction}:

\snic {
\PP^{m}(\gK) \to \PP^{m}(\gk),  \qquad
(x_0 : \ldots : x_m) \mapsto (\xi_0 : \ldots : \xi_m)
\quad \hbox {avec $\xi_i = \ov {x_i/x_{i_0}}$}
,}

%\sni
où $x_{i_0} \,\vert\, x_i$ pour
tout $i$; l'\elt $(\xi_0: \ldots: \xi_n)\in\PP^{m}(\gk)$ est bien défini: il correspond à un \vmd de $\gV^{m+1}$. 
En bref on a un \gui{\isoz} $\PP^{m}(\gV)\simeq\PP^{m}(\gK)$ et une {réduction} $\PP^{m}(\gV)\to\PP^{m}(\gk)$.
\\
Ici le choix de l'uniformisante $\pi =
t^{-1}$ donne une \dfn directe, sans  passer au projectif, 
de la réduction  $\gK^n \setminus\{0\} \to \gk^n \setminus \{0\}$.

On dira qu'une matrice $A \in \GL_n(\gK)$ de colonnes $(A_1, \ldots, A_n)$
est \emph{$\gA_\infty$-réduite} si la matrice $\ov A \in \MM_n(\gk)$  est dans $\GL_n(\gk)$.

\emph {1.}
Soit $A \in \GL_n(\gK)$ de colonnes $A_1, \ldots, A_n$. Montrer
que $\sum_{j=1}^n v(A_j) \le v(\det A)$.

\emph {2.}
Soit $A \in \GL_n(\gK)$;  calculer $Q \in
\GL_n(\gA)$ telle que $AQ$ soit $\gA_\infty$-réduite.  Ou encore: soit $E
\subset \gK^n$ un $\gA$-réseau, i.e. un $\gA$-module libre de rang $n$;
alors $E$ admet une $\gA$-base $\gA_\infty$-réduite (une base $(A_1,
\ldots, A_n)$ telle que $(\ov {A_1}, \ldots, \ov {A_n})$ soit une $\gk$-base de
$\gk^n$). On pourra commencer par l'exemple $A = \cmatrix {\pi^2 &\pi\cr 1 &
1\cr}$.

\emph {3.}
Pour $P \in \GL_n(\gA_\infty)$, montrer les points suivants.
\begin {itemize}
\item [\emph {a.}]
$P$ est une $v$-isométrie, i.e. $v(Px) = v(x)$ pour tout $x \in \gK^n$.

\item [\emph {b.}]
Pour tout $x \in \gK^n \setminus \{0\}$, $\ov {Px} = \ov{P}\,\ov{x}$.

\item [\emph {c.}]
Si la matrice $A \in \GL_n(\gK)$ est $\gA_\infty$-réduite, il en est de
même de $PA$.
\end {itemize}

\emph {4.}
Soit $A \in \GL_n(\gK)$ triangulaire. Que signifie \gui{$A$ est
$\gA_\infty$-réduite}?

\emph {5.}
Soit $A \in \GL_n(\gK)$. Montrer qu'il existe $Q \in \GL_n(\gA)$,
$P \in \GL_n(\gA_\infty)$ et des entiers $d_i \in \ZZ$ tels que
$PAQ = \Diag(t^{d_1}, \ldots, t^{d_n})$; de plus, si l'on range
les $d_i$ par ordre croissant, ils sont uniques.

\emph {6.}
Soit $L$ un $\gK$-\evc de dimension $n$, $E \subset L$ un $\gA$-réseau,  et $E'\subset L$ un $\gA_\infty$-réseau. 
\begin {itemize}
\item [\emph {a.}]
Montrer qu'il existe une $\gA$-base $(e_1,\ldots, e_n)$ de $E$, une $\gA_\infty$-base $(e'_1, \ldots, e'_n)$ de $E'$ et des
entiers $d_1$, \ldots, $d_n \in \ZZ$ vérifiant $e'_i = t^{d_i} e_i$ pour $i\in \lrb{1..n}$.  De plus, les $d_i$ rangés par ordre croissant ne dépendent que de $(E, E')$.
\item [\emph {b.}]
En déduire que $E \cap E'$ est un \kev de dimension finie.\\
De manière
précise:

\snic {
E \cap E' = \bigoplus_{d_i \ge 0} \bigoplus_{j=0}^{d_i} \gk t^j e_i
,}

%\sni
et en particulier:

\snic {
\dim_\gk(E\cap E') = \sum_{d_i \ge 0} (1+d_i) = \sum_{d_i \ge -1} (1+d_i) 
}
\end {itemize}

\emph {7.}
On suppose que $\gL$ est une $\gK$-extension finie de degré $n$. On
définit des clôtures intégrales dans $\gL$: $\gB$ celle de $\gA$,
$\gB_\infty$ celle de $\gA_\infty$ et $\gk'$ celle de $\gk$.  On dit qu'une
base $(\ue) = (e_1, \ldots, e_n)$ de $\gB$ sur $\gA$ est \emph {normale à
l'infini} s'il existe $r_1$, \ldots, $r_n \in \gK^*$ tels que $(r_1e_1, \ldots,r_ne_n)$ est une $\gA_\infty$-base de $\gB_\infty$. Montrer que les \elts de la base $(\ue)$ \gui {entiers à l'infini}, 
i.e. qui appartiennent à $\gB_\infty$,
forment une $\gk$-base de l'extension~$\gk'$.

\emph {8.}
Soit $\gk=\QQ$, $\gL = \aqo{\gk[X,Y]}{X^2 + Y^2} = \gk[x,y]$, $\gA =
\gk[x]$. \\
Montrer que $(y+1, y/x)$ est une $\gA$-base de $\gB$ mais qu'elle
n'est pas normale à l'infini.  Expliciter une base normale à l'infini de
$\gB\sur\gA$.

}
\end {problem}
%--- end -problem-----------------------------------------

%--- Problem{exoHyperEllipticFunctionRing}-------------
\begin{problem}\label{exoHyperEllipticFunctionRing}
{(Anneau des fonctions d'une courbe hyper-elliptique
affine ayant un seul point à l'infini)}\\
{\rm  
On utilisera ici une notion de \emph{norme d'un \idz} dans le contexte suivant: $\gB$
étant une \Alg libre de rang fini $n$ et $\fb$ un \itf de $\gB$,
la norme de $\fb$ est l'\id 

\snic{\rN_{\gB\sur\gA}(\fb)=\rN(\fb) \eqdefi \cF_{\gA,0}(\gB\sur\fb)\subseteq\gA.}

%\sni
 Il est clair que pour $b \in \gB$,
$\rN(b\gB) = \rN_{\gB\sur\gA}(b)\gA$, que $\rN(\fa\gB) = \fa^n$ pour $\fa$
un \itf de $\gA$ et que $\fb_1 \subseteq \fb_2 \Rightarrow \rN(\fb_1)
\subseteq \rN(\fb_2)$.%
\index{norme!d'un idéal}

Soit $\gk$ un corps de \cara $\ne 2$, $f = f(X)\in \gk[X]$ un \polu \spl de
degré impair $2g+1$. Le \pol $Y^2 - f(X) \in \gk[X,Y]$ est absolument \ird;
on pose $\gB = \aqo {\gk[X,Y]}{Y^2 - f(X)} = \gk[x,y]$ et $\gA = \gk[x] \simeq
\gk[X]$. %On identifiera un \itf  de $\gA$ et son unique \gtruz.
L'anneau $\gB$ est intègre, c'est un \Amo libre de base $(1,y)$. 
Pour  $z=a + by$ avec $a$, $b \in \gA$, on note $\ov z=a-yb$, et $\rN =
\rN_{\gB\sur\gA}$: $\rN(z) = z\ov z = a^2 - fb^2$. 
\\
Le but du \pb est de paramétrer les \itfs non nuls de $\gB$, 
de montrer que $\gB$ est un \adp
et d'étudier le groupe $\Cl(\gB)$ des classes d'\ids \ivs de $\gB$.
\\
Si $\fb$ est un \itf de $\gB$, son \emph{contenu} est l'\idf $\cF_{\gA,1}(\gB\sur\fb)$. 
\\
\`A deux \elts $u$, $v \in \gA$ vérifiant $v^2 \equiv f \bmod u$, on associe
le sous-\Amo \hbox{de $\gB$}: $\fb_{u,v} = \gA u + \gA (y-v)$. On a
$u \neq  0$ parce que $f$ est \splz. On fera
parfois intervenir le \pol $w \in \gA$ tel que $v^2 - uw = f$ et l'on notera
$\fb_{u,v,w}$ au lieu de $\fb_{u,v}$ (même si $w$ est complètement
déterminé par $u,v$). 

\emph {1.}
Montrer que $\fb_{u,v}$ est un \id de $\gB$ et que $\fb_{u,v} = \gA u \oplus \gA
(y-v)$.  Réciproquement,
pour $u, v \in \gA$, si $\gA u + \gA (y-v)$ est un \id de $\gB$, alors $v^2
\equiv f \bmod u$.

\emph {2.}
Montrer que $\gA \to \gB\sur{\fb_{u,v}}$ induit un \iso $\gA\sur{\gA u}
\simeq \gB\sur{\fb_{u,v}}$; en conséquence, $\Ann_\gA(\gB\sur{\fb_{u,v}})=\gA u$.
%le \Amo $\gB\sur{\fb_{u,v}}$ est d'annulateur $u\gA$. 
En déduire \gui {l'unicité de $u$}:

\snic {u_1,\,u_2
\hbox { \mons et } \fb_{u_1,v_1} = \fb_{u_2,v_2} \;\; \Longrightarrow  \;\;u_1 = u_2.
}

%\sni

Vérifier \egmt que $\rN(\fb_{u,v}) = u\gA$ et que $v$ est unique modulo $u$:

\snic {
\fb_{u,v_1} = \fb_{u,v_2}  \iff  v_1 \equiv v_2 \mod u
.}

%\sni
\emph {3.}
Montrer que

\snic {
\fb_{u,v,w}\fb_{w,v,u} = \gen {y-v}_\gB, \qquad
\fb_{u,v}\fb_{u,-v} = \gen {u}_\gB.
}

%\sni

En conséquence, l'\id $\fb_{u,v}$ est \ivz.  \\
De plus, pour $u = u_1u_2$
vérifiant $v^2 \equiv f \bmod u$, on a $\fb_{u,v} = \fb_{u_1,v}\,\fb_{u_2,v}$.

\emph {4.}
Soit $\fb$ un \itf non nul de $\gB$. 
\begin {itemize}
\item [\emph {a.}]
Montrer qu'il existe deux \polus uniques $d$, $u \in \gA$ et $v \in \gA$\linebreak 
avec $v^2 \equiv f \bmod u$, tels que $\fb = d\,\fb_{u,v}$.  En conséquence,
$\fb$ est un \id\iv (donc $\gB$ est un \adpz).  De plus, $v$ est unique modulo
$u$, donc unique si l'on impose $\deg v < \deg u$.

\item [\emph {b.}]
En déduire que $\fb\ov\fb = \rN(\fb)\gB$ puis que la norme est
multiplicative sur les \idsz.

\item [\emph {c.}]
Montrer que $\gB\sur\fb$ est un \kev de dimension finie.\\
Montrer que $\dim_\gk
(\gB\sur\fb) = \dim_\gk (\gA\sur\fa)$ avec $\fa = \rN(\fb)$. Cet entier sera \linebreak 
noté $\deg(\fb)$. Vérifier que  $\deg(\fb_{u,v}) = \deg u$, que $\deg
(\fb) = \deg \rN(\fb)$, et enfin que $\deg$ est additif, i.e. que
$\deg(\fb_1\fb_2) = \deg(\fb_1) + \deg(\fb_2)$.
\end {itemize}

Soient $u$, $v \in \gA$ avec $v^2 \equiv f \bmod u$. On dit que \emph{le couple $(u,v)$ est
réduit} si~$u$ est \mon et \framebox [1.1\width]{$\deg v < \deg u \le g$}.
Par abus de langage, on dit aussi que $\fb_{u,v}$ est réduit.
Par exemple, si $(x_0, y_0)$ est un point de la courbe hyperelliptique
$y^2 = f(x)$, son \id $\gen {x-x_0, y-y_0}$ est un \id réduit
(prendre $u(x) = x-x_0$, $v = y_0$).

\emph {5.}
Montrer que tout \itf non nul de $\gB$ est associé à un idéal réduit
de $\gB$ (deux \ids $\fa$ et $\fa'$ sont dits \emph{associés} s'il existe deux \elts
\ndzs $a$ et $a'$ tels que $a\fa'=a'\fa$, on note alors $\fa\sim\fa'$).

\emph {6.}
Dans cette question, pour un \itf non nul $\fb$ de $\gB$, on désigne 
par~$\rN(\fb)$ le \polu \gtr de l'\id $\rN_{\gB/\gA}(\fb)$.  Soit
$\fb_{u,v}$ un \id réduit.
\begin {itemize}
\item [\emph {a.}]
Soit $z \in \fb_{u,v} \setminus \{0\}$ de sorte que $u = \rN(\fb_{u,v}) \divi
\rN(z)$, i.e. $\rN(z)/\rN(\fb_{u,v})$ est un \polz.  Montrer que
$$
\deg \left(\rN(z)/\rN(\fb_{u,v})\right) \ge \deg u,
$$
avec l'\egt \ssi $z \in \gk\eti u$. 

\item [\emph {b.}]
Soit $\fb'$ un \itf  de $\gB$ vérifiant $\fb' \sim
\fb_{u,v}$. Montrer que $\deg(\fb') \ge \deg(\fb_{u,v})$ avec l'\egt \ssi
$\fb' = \fb_{u,v}$. En résumé, dans une classe d'idéaux inversibles de
$\gB$, il y a donc un et un seul idéal de degré minimum: c'est l'unique
idéal réduit de la classe.
\end {itemize}
\emph {7a.}
Montrer que la courbe affine $y^2 = f(x)$ est lisse; de manière précise,
en posant $F(X,Y) = Y^2 - f(X) \in \gk[X,Y]$, montrer que $1 \in \gen {F,
F'_X, F'_Y}$; cela utilise uniquement le fait que $f$ est \spb et que la \cara
de $\gk$ est distincte de $2$, pas le fait que $f$ est de degré impair.

Si $\gk$ est \acz, on obtient ainsi une correspondance biunivoque entre les
points $p_0 = (x_0,y_0)$ de la courbe affine $y^2 = f(x)$ et les \advs
(discrète) $\gW$ de $\gk(x,y)$ contenant $\gB = \gk[x,y]$: à $p_0$, on
associe son anneau local $\gW$ et dans l'autre sens, à $\gW$ on associe le point
$p_0 = (x_0,y_0)$ tel \linebreak 
que $\gen {x-x_0, y-y_0}_\gB = \gB \cap \fm(\gW)$.

\emph {b.}
On étudie maintenant \gui{les points à l'infini de la courbe projective
lissifiée}, à l'infini relativement au modèle $y^2 = f(x)$. De manière
\agqz, il s'agit des \advs pour $\gk(x,y)$ ne contenant pas $\gB$ (mais
contenant $\gk$ bien entendu).  Soit l'\adv discrète $\gA_\infty =
\gk[x^{-1}]_{\gen {x^{-1}}}$. Montrer qu'il existe un et un seul anneau
$\gB_\infty$, $\gA_\infty \subseteq \gB_\infty \subseteq \Frac(\gB) =
\gk(x,y)$, ayant $\Frac(\gB)$ pour corps de fractions. Montrer que
$\gB_\infty$ est un \adv discrète, que $\gB_\infty\sur{\fm(\gB_\infty)}
\simeq \gA_\infty\sur{\fm(\gA_\infty)} \simeq \gk$ et que c'est le seul point
à l'infini.
}
\end {problem}
%--- end -problem-----------------------------------------

%--- problem{exoTrifolium}-------------
\begin{problem}\label{exoTrifolium}
{(Trifolium: \cli et paramétrisation)}
\\
{\rm  Soit $\gk$ un \cdi et 

\snic{F(X,Y) = (X^2+Y^2)^2 + \alpha X^2Y + \beta Y^3,}

%\sni
avec $\alpha \ne \beta$ dans $\gk$.

\Deuxcol{.50}{.45}
{
On étudie la courbe $F(x,y) = 0$, ses points singuliers,
son corps de fonctions 

\snic{\gL = \gk(x,y)}

%\sni

(on va montrer que $F$ est \irdz), son
anneau de fonctions $\gk[x,y]$, la \cli $\gB$ de $\gk[x,y]$
dans $\gL$~\ldots\ {etc}\ \dots \\
On note que $F(-X,Y) = F(X,Y)$ et donc
l'involution $(x,y) \mapsto (-x,y)$ laisse invariante la courbe $F(x,y) = 0$.
\\
Ci-contre, un exemple d'une telle courbe.}
{~\vspace{-5mm}

\includegraphics*[width=5cm]{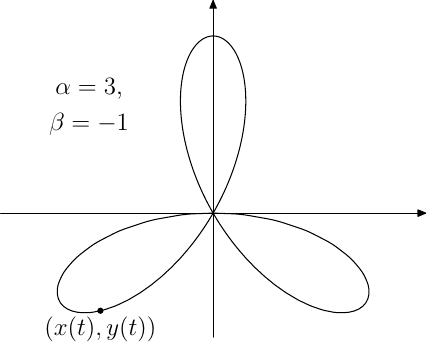}}

\vspace{2mm}
\emph {1.}
Montrer que $F$ est un \pol absolument \irdz. Plus \gnltz: soit $\gk$ un anneau intègre, $\gk[\uT]$ un anneau de \pols
à plusieurs \idtrs et $F \in \gk[\uT]$, $F = F_N + F_{N+1}$ avec 
$F_N$, $F_{N+1}$ \hmgs non nuls, de degrés respectifs $N$, $N+1$. Alors, dans toute
\fcn $F = GH$, l'un des deux \pols $G$ ou $H$ est \hmg; enfin, si
$\gk$ est un corps, alors $F$ est \ird \ssi $F_N$, $F_{N+1}$ sont premiers entre
eux.

\emph {2.}
Déterminer les points singuliers de la courbe $F = 0$.

On note $\gL = \gk(x,y)$ et $\gB$ la \cli de
$\gk[x,y]$ dans $\gL$.

\emph {3.}
Soit $t = y/x$ de sorte que $\gL = \gk(x,t)$.  
\begin {itemize}
\item [\emph {a.}]
Déterminer une \eqn primitive \agq de $t$ sur $\gk[x]$. \\
On note $G(X,T) =
a_4 T^4 + \cdots + a_1 T + a_0 \in \gk[X][T]$, avec $a_i=a_i(X) \in \gk[X]$,
un tel \pol primitif, vérifiant donc $G(x,t) = 0$. Vérifier que $(x,
t)=(0,0)$ est un point non singulier de la courbe $G = 0$.

\item [\emph {b.}]
Déterminer les entiers d'Emmanuel $b_4$, \ldots, $b_1$
associés à $(G, t)$ avec $\gA = \gk[x]$ comme anneau de base
(lemme \ref{lemEmmanuel}).
En déduire une \mlp pour $(x,y)$ et préciser l'\id $\fq$ de $\gB$
que $\gen{x}_\gB = \fq \gen {x,y}_\gB$.
\end {itemize}

\emph {4.}
Montrer que $\gL = \gk(t)$ et exprimer $x$, $y$ comme \elts de $\gk(t)$.  

\emph {5.}
Détermination de la \cli $\gB$ de $\gk[x,y]$ dans $\gL$.
\begin {itemize} 
\item [\emph {a.}]
Montrer que $\gB = \gk[g_0, g_1]$ avec $g_0 = 1/(1+t^2)$ et $g_1 = tg_0$.
Exprimer $x$, $y$ dans $\gk[g_0, g_1]$. Quelle est \gui {l'\eqnz} 
liant $g_0$  et $g_1$?
\item [\emph {b.}]
Montrer que $(1, y, b_3t, b_2t)$ est une $\gA$-base de $\gB$.
\item [\emph {c.}]
Vérifier que $\dim_\gk  \gB\sur{\gen{x,y}_\gB} = 3$.
\end {itemize}

\emph {6.}
On note $\gV$ l'\advz\footnote{Un sous-anneau $\gV$ d'un \cdi $\gL$ est appelé un \emph{\adv de $\gL$} si pour tout $x\in\gL\eti$ on a $x\in\gV$ ou $x^{-1}\in\gV$.} de $\gL$ défini par le point non singulier $(0, 0)$ de la courbe $
G = 0$. C'est le seul \adv de $\gL$ contenant $\gk$ et tel \hbox{que $x$, $t
\in \Rad\gV$} (et donc aussi $y \in \Rad\gV$). \\
On considère l'\idep $\fp_1 =
(\Rad\gV) \cap \gB$. Montrer que:

\snic {
\fp_1 = \gen {x, y, b_4t, b_3t, b_2t, b_1t} = \gen {g_0-1, g_1} 
\quad \hbox {et} \quad \gB\sur{\fp_1} = \gk
,}

%\sni
et vérifier que $\fp_1^2 = \gen {g_0-1, g_1^2}$.

\emph {7.}
Déterminer la \fcn dans $\gB$ de l'\id $\gen {x,y}_\gB$ en produit d'
\idepsz.  La réponse n'est pas uniforme en $(\alpha, \beta)$, contrairement
à la détermination de la \cli $\gB$ de $\gA$.

\emph {8.}
Reprendre les questions en supposant seulement que $\gk$ est
un anneau \icl et que $\beta-\alpha\in\gk\eti$.

}

\end {problem}
%--- end -problem-----------------------------------------

}% fin des exos
%:  solutions
\sol{

%%%%%%%%%%%%%%%%%%%%%%%%%%%%%%%%%%%%%%%%%%%%%%%%%%%%%%%%%%%%%%%%%%%%%%%%%%%
\exer{exoDetTrick1}{Il faut montrer l'inclusion $\fb^n\subseteq\fa\fb^{n-1}$.
Soient $(x_1, \ldots, x_n)$ un \sgr  de $E$, $X=\tra{ [\,{x_1\;\cdots\; x_n}\,]}$,  $b_1$, \ldots, $b_n \in \fb$ et  $B=\Diag(b_1, \ldots, b_n)$.
Puisque $b_ix_i \in \fa E$ ($i\in\lrbn$), il existe  $A \in \Mn(\fa)$ telle que
$B\,X = A \,X$. \\
Soit $C = B - A$, il vient $C \,X = 0$, 
 et puisque  $E$ est fidèle, $\det C  = 0$.
On développe ce \deterz, il vient $b_1 \cdots b_n + a = 0$ 
avec $a \in \fa\fb^{n-1}$ (car $\fa\subseteq\fb$).
}

%%%%%%%%%%%%%%%%%%%%%%%%%%%%%%%%%%%%%%%%%%%%%%%%%%%%%%%%%%%%%%%%%%%%%%%%%%%
\exer{exoMatLoc22}
\emph {1.}
Immédiat, car si $B = \cmatrix {b_{11}&b_{12}\cr b_{21}&b_{22}\cr}$,
alors $\wi B = \crmatrix {b_{22}&-b_{12}\cr -b_{21}&b_{11}\cr}$ et
$[x \; y\,]B = [x' \; y'\,]$ avec:

\snic {
x' = -\left|\matrix {-b_{21} & b_{11}\cr x & y\cr}\right|, \quad
y' = \left|\matrix {b_{22} & -b_{12}\cr x & y\cr}\right|
.}

%\sni
\emph {2.}
On a $u$, $v \in \fb$ avec $z = ux + vy$ et $ux$, $uy$, $vx$, $vy$
multiples de $z$, ce que l'on écrit
$\crmatrix {y\cr -x\cr} [v \; -u\,] = zB$.
%%\cmatrix {vy & -uy\cr -vx & ux\cr} = zB
%%\quad \hbox {avec} \quad B = \cmatrix {a & b\cr c & d\cr} 
Comme $[x \; y\,] \crmatrix {y\cr -x\cr} = 0$, on a
$[x \; y\,]zB = 0$; de \hbox{plus $\Tr(zB) = yv + xu = z$}.

\emph {3.}
Dans le lemme en question, $z = x^n$ et l'anneau est \lsdzz.  Les \egts
$x^n[x\; y\,]B = 0$ et $x^n\big(1 - \Tr(B)\big) = 0$ fournissent deux \lons\come de
$\gA$: une dans laquelle $x^n = 0$, auquel cas $x = 0$ car l'anneau $\gA$ et
son localisé sont réduits, et l'autre dans laquelle $[x\; y\,]B = 0$ et
$\Tr(B) = 1$. Dans chacune d'entre elles, $\gen {x,y}$ est \lop donc il l'est
dans $\gA$.
%%%%%%%%%%%%%%%%%%%%%%%%%%%%%%%%%%%%%%%%%%%%%%%%%%%%%%%%%%%%%%%%%%%%%%%%%%%

\exer{exoPrufNagata}  \emph{1.} En effet, $\gA(X)$ est \fpt sur $\gA$.

 \emph{2.}
Soit $f=\sum_{k=0}^n a_kX^k\in\AX$%, $\fa=\rc_\gA(f)$ et $\fA=\fa\gA(X)$.
%On peut aussi noter ces deux \ids $\gen{a_0,\ldots,a_n}$
. 
Pour chaque $k$, on a dans  $\gA$ une \egt 

\snic{\gen{a_0,\ldots,a_n}\gen{b_{0,k},\ldots,b_{n,k}}=\gen{a_k}
$  avec $ a_0b_{0,k}+\cdots +a_nb_{n,k}=a_k.}

%\sni
Considérons alors le \pol $g_k=\sum_{j=0}^n b_{j,k}X^{n-j}$. 
Tous les \coes de $fg_k$ sont dans $\gen{a_k}$. On peut donc écrire
$fg_k=a_kh_k$ avec le \coe de degré $k$ dans~$h_k$ égal à $1$.
Ceci implique que dans $\gA(X)$, $a_k\in\gen{f}$. Or on a $f\in\gen{a_0,\ldots,a_n}$ dans $\AX$. Ainsi, dans  $\gA(X)$, $\gen{f}=\gen{a_0,\ldots,a_n}$.
\\
On en déduit que  $\gA(X)$ est un anneau de Bézout, car pour
$f_0$, \ldots, $f_m\in\AX$ de degrés $<d$, une conséquence du résultat précédent est que dans $\gA(X)$:

\snic{\gen{f_0,\ldots,f_m}=\gen{f_0+X^df_1+\cdots+X^{dm}f_m}.}

%\sni
\emph{3.}
D'après l'exercice \ref {exoMatLoc22}, $(x,y)$ admet une \mlp sur~$\gB$ \ssi
il existe $B \in \MM_2(\gB)$ de trace $1$ vérifiant $[\,x\; y\,]B = [\,0\;
0\,]$%: dans ce cas, $\wi B$ est une \mlp pour $(x, y)$
.
\\
Soit donc $B \in \MM_2\big(\gA(X)\big)$ vérifiant $[\,x\; y\,] B = [\,0\; 0\,]$ et
$\Tr(B) = 1$.  \\
En multipliant les \coes de $B$ par un \deno commun,
on obtient des \elts $p$, $q$, $r$, $s$ de $\AX$ tels que $[\,x\; y\,] \cmatrix {p &q\cr r
&s\cr} = [\,0\; 0\,]$ et $p + s$ primitif.  
On a donc (avec $p=\sum_k p_kX^{k}$, \dots):  $[\,x\; y\,] \cmatrix {p_i &q_i\cr r_i &s_i\cr} = [\,0\; 0\,]$.
Comme $p+s$ est
primitif, on a des $u_i \in \gA$ tels que $\sum u_i(p_i+s_i) = 1$. Soit $B' = \sum_i u_i\cmatrix
{p_i&q_i\cr r_i &s_i}\in\MM_2(\gA)$: on obtient $[\,x\; y\,] B' = [\,0\;
0\,]$ avec $\Tr(B')=1$.

\emph{4.} $\ArX$ \ari $\Rightarrow$ $\gA$ \ari et 

\snic{ \gA $ \ari $ \,\iff\,\gA(X) $ \ari  $\,\iff\,\gA(X)$ Bézout$.}

%\sni
La dernière \eqvc résulte aussi du \plgref{thlgb3}.
En outre, le \mo de la \dve dans $\gA(X)$, i.e. $\gA(X)/\gA(X)\eti$, est isomorphe
au \mo des \itfs de $\gA$.

%%%%%%%%%%%%%%%%%%%%%%%%%%%%%%%%%%%%%%%%%%%%%%%%%%%%%%%%%%%%%%%%%%%%%%%%%%%
\exer{FermetureAlgZariski} 
On montre seulement le premier point. Il est clair que
$\gK'(\uX)$ est \agq sur $\gK(\uX)$. \\
Réciproquement, soit
$z \in \gL(\uX)$ \agq sur $\gK(\uX)$, il existe $a \in \gK[\uX]$
non nul tel que $az$ soit entier sur $\gK[\uX]$, a fortiori
sur $\gL[\uX]$. Comme $\gL[\uX]$ est un anneau à pgcd, 
on a $az \in \gL[\uX]$. Par ailleurs, on sait que la
\cli de $\gK[\uX]$ dans~$\gL[\uX]$ est $\gK'[\uX]$
(lemme \ref{lemPolEnt});
donc $az \in \gK'[\uX]$ puis $z = (az)/a \in \gK'(\uX)$.

%%%%%%%%%%%%%%%%%%%%%%%%%%%%%%%%%%%%%%%%%%%%%%%%%%%%%%%%%%%%%%%%%%%%%%%%%%%

\exer{exoNotAbsIntClos} 
\emph {1.}
Immédiat.

Dans la suite on va utiliser le fait que $(1,y)$ est une $\gk[x]$-base de $\gk[x,y]$;
c'est aussi une $\gk(x)$-base de $\gk(x,y)$ et l'extension 
$\gk(x,y)\sur{\gk(x)}$ est galoisienne de groupe $\gen {\sigma}$
où $\sigma : \gk(x,y) \to \gk(x,y)$ est le $\gk(x)$-\auto
involutif qui réalise $y \mapsto -y$.

\emph {2.}
Soit $z = u(x) + yv(x) \in\gk(x,y)$ \agq sur $\gk$. \\
Alors $z+\sigma(z) = 2u$
et $z\sigma(z) = u^2 - fv^2$ sont \agqs sur $\gk$ et dans $\gk(x)$ donc
dans $\gk$. D'où $u \in \gk$, $v = 0$ et $z = u \in \gk$.

\emph {3.}
Comme $a \notin \gk^p$, on voit facilement que $f(X)$ est \ird dans $\gk[X]$.
Montrons que $\gk[x,y]$ est la \cli de $\gk[x]$ dans~$\gk(x,y)$. \\
Soit $z = u(x) + yv(x) \in\gk(x,y)$ entier sur
$\gk[x]$. \\
Alors $z+\sigma(z) = 2u$ et $z\sigma(z) = u^2 - fv^2$ sont dans
$\gk(x)$ et entiers sur $\gk[x]$, donc dans $\gk[x]$. 
Ainsi $u$ et $fv^2 \in
\gk[x]$. En utilisant le fait que $f$ est \irdz, on voit que $v \in
\gk[x]$. Bilan: $z \in \gk[x,y]$.

\emph {4.}
Soit $\alpha=a^{1/p}\in \gk$, d'où $f(X) = (X-\alpha)^p$.
On pose $t = y/(x-\alpha)^{p-1\over 2}$. \\
Alors $t^2 = x-\alpha$,
donc $x \in \gk[t]$. Et $y = t(x-\alpha)^{p-1\over 2} = t^p \in \gk[t]$.\\
Donc $\gk [x,y]  \subseteq  \gk [t] $ et $\gk(x,y) = \gk(t)$. On voit que $t$ est entier sur $\gk[x]$,
mais \linebreak 
que $t\notin\gk[x,y] = \gk[x] \oplus \gk[x]y$. La \cli de $\gk[x]$ (ou celle de $\gk[x,y]$) dans~$\gk(x,y)$ est $\gk[t]$
(qui contient bien $x$ et $y$).

%%%%%%%%%%%%%%%%%%%%%%%%%%%%%%%%%%%%%%%%%%%%%%%%%%%%%%%%%%%%%%%%%%%%%%%%%%%

\exer{exoAnneauOuvertP1} Rappelons que $x_0=\fraC1 p$. L'\egt

\snic{\gk[x_0, \ldots, x_{n-1}] = \sotq {u/p^s} {u \in \gk[t],\  \deg(u) \le ns},}

%\sni
est facile en remarquant que $t^n x_0 \in \gk[x_0, \ldots, x_{n-1}]$ puisque

\snic{{t^n \over p} = 1 + {t^n - p \over p} \;\in \;\sum_{i=0}^{n-1}\gk\,{ {t^i}  \over p}.}

%\sni
\'Ecrivons que $t$ est \agq sur $\gk(x_0)$ comme racine en $T$ du \pol

\snic{p(T)x_0 - 1=x_0T^n + x_0a_{n-1}T^{n-1} + \cdots+x_0 a_1T + (x_0a_0 - 1).}

%\sni

Les \gui{entiers d'Emmanuel} (cf. lemme \ref{lemEmmanuel} ou exercice \ref{exoEntiersEmmanuel}) sont 
\[ 
\begin{array}{rcl} 
  x_0t, &x_0t^2 + x_0a_{n-1}t,   & x_0t^3 + x_0a_{n-1}t^2 + x_0a_{n-2}t,  \\[1mm] 
&\ldots,&x_0t^{n-1} + \cdots + x_0a_2t . 
\end{array}
\]
Ainsi, $t^k x_0$ est entier sur $\gk[x_0]$ pour $k \in \lrb
{0..n-1}$ et $\gk[x_0, \ldots, x_{n-1}] \subseteq \gA$.

Reste à voir que $\gA \subseteq \gk[x_0, \ldots, x_{n-1}]$.  On
utilise l'inclusion

\snic{\gk[x_0]\subseteq\gV_\infty := \gk[1/t]_{1+\gen{1/t}}.}

%\sni
Ce dernier anneau est constitué des fractions rationnelles de degré $\le 0$,
i.e. définies en $t=\infty$. Il est isomorphe à $\gk[y]_{1+\gen{y}}$ 
donc il est \iclz, et $\gA\subseteq\gV_\infty$. 
L'anneau $\gV_\infty$ est appelé \gui{l'anneau local du point
$t = \infty$}.
\\
Soit $z \in \gk(t)$ une fraction
rationnelle entière sur $\gk[x_0]$. En multipliant une \rdi de
$z$ sur $\gk[x_0]$ par $p^N$ avec $N$ assez grand, on obtient

\snic {
p^N z^m + b_{m-1} z^{m-1} + \cdots + b_1 + b_0 = 0, \qquad b_i \in \gk[t]
.}

%\sni
Ceci entraîne que $p^N z$ est entier sur $\gk[t]$ donc appartient
à $\gk[t]$ ($\gk[t]$ est \iclz). Par ailleurs, $z \in \gV_\infty$, i.e., $\deg
z \le 0$.  En définitive, $z$ est une fraction rationnelle de degré $\le
0$ dont le \deno divise une puissance de $p$, \linebreak 
donc $z \in \gk[x_0, \ldots,
x_{n-1}]$.

Enfin, on a $x_1 = tx_0$ donc $t = x_1/x_0 \in \Frac(\gA)$ puis
$\gk(t) = \Frac(\gA)$.

%%%%%%%%%%%%%%%%%%%%%%%%%%%%%%%%%%%%%%%%%%%%%%%%%%%%%%%%%%%%%%%%%%%%%%%%%%%

\exer{exoPresentationAlgOuvertP1} 
\emph {1.}
Immédiat.

\emph {2.}
On note $\fb$ l'\id engendré par les
$X_iX_j - X_{i-1}X_{j+1}$ pour $1 \le i \le j \le n-1$ et $E$
le \kmo:

\snic {
E = \gk[X_0] \oplus \gk[X_0]X_1 \oplus \cdots \oplus \gk[X_0]X_{n-1}
.}

%\sni
On va prouver que $E \cap \Ker\varphi = 0$ et que
$\gk[\uX] = E + \fb$. Comme $\fb \subseteq \Ker\varphi$, on
obtiendra $\gk[\uX] = E \oplus \fb$. Soit $y \in \Ker\varphi$
que l'on écrit $y = y_1 + y_2$ avec $y_1 \in E$ \linebreak 
et $y_2 \in \fb$.
En appliquant $\varphi$, on obtient $\varphi(y_1) = 0$, donc $y_1 = 0$,
puis $y = y_2 \in \fb$. On a obtenu $\Ker\varphi \subseteq \fb$, d'où
l'\egt $\Ker\varphi = \fb$.

$\bullet$ \emph{Justification de $E \cap \Ker\varphi = 0$}. Soit $f \in E$

\snic {
f = f_0(X_0) + f_1(X_0)X_1 + \cdots + f_{n-1}(X_0)X_{n-1}
.}

%\sni
On écrit que $\varphi(f) = 0$:

\snic {
\varphi(f) = f_0(1/p) + f_1(1/p)t/p + \cdots + f_{n-1}(1/p)t^{n-1}/p = 0
.}

%\sni
En multipliant chaque $f_i(1/p)$ par $p^N$,  $N$ assez grand, on
obtient $g_i(p)\in \gk[p]$:

\snic {
pg_0(p) + g_1(p)t + \cdots + g_{n-1}(p)t^{n-1} = 0
.}

%\sni

Mais $(1, t, \ldots, t^{n-1})$ est une base de $\gk[t]$ sur
$\gk[p]$, donc les $ g_{k} = 0$, puis
$f = 0$.

$\bullet$
\emph{Justification de $\gk[\uX] = E + \fb$}. \\
En notant $\gk[x_0, \ldots, x_{n-1}] =
\gk[\uX]/\fb$ et $E' = \gk[x_0] + \gk[x_0]x_1 + \cdots + \gk[x_0]x_{n-1}$,
cela revient à montrer que $\gk[\ux] = E'$. Or $E'$ contient $x_n := 1 - \sum_{i=0}^{n-1}
a_ix_i$. Il suffit donc de prouver
que $E'$ est un sous-anneau, ou encore \hbox{que $x_ix_j \in E'$} \hbox{pour
$i$, $j \in \lrb {0..n-1}$}.  Par \dfn, il contient $x_0^2$, $x_0x_1$, \dots, $x_0x_{n-1}$
donc aussi~$x_0x_n$. Mais $x_1x_j = x_0x_{j+1}$ pour $j \in \lrb {1..n-1}$,
donc $E'$ contient ces~$x_1x_j$ donc aussi~$x_1x_n$. Et en 
utilisant $x_2x_j = x_1x_{j+1}$ \hbox{pour $j \in \lrb {2..n-1}$}, on
voit que~$E'$ contient tous les~$x_2x_j$. Et ainsi de suite.

\rem
L'auteur de l'exercice a opéré ainsi pour un \cdi $\gk$: il a utilisé
une \idtr supplémentaire $X_n$ et a choisi sur $\gk[X_0, X_1, \ldots, X_n]$
l'ordre monomial gradué lexicographique inversé en ordonnant les \idtrs de
la manière suivante: $X_0 < X_1 < \cdots < X_n$.  On constate alors que
l'\id initial de l'\id $\gen {R_{\rm min}} + \gen {1 - \sum_{i=0}^n a_iX^i}$
est l'\id monomial engendré par les monômes:

\snic {
(\star) \qquad\qquad
X_n  \hbox { et } X_iX_j  \hbox { pour }  1 \le i \le j \le n-1
.}

%\sni
Le \kev engendré par les monômes non divisibles par un monôme de~$(\star)$
est le \kev $E = \gk[X_0] \oplus \gk[X_0]X_1 \oplus \cdots \oplus \gk[X_0]X_{n-1}$.
C'est celui qui apparaît dans le corrigé ci-dessus (dans lequel $\gk$
est un anneau quelconque, pas \ncrt un \cdiz).
\eoe

%%%%%%%%%%%%%%%%%%%%%%%%%%%%%%%%%%%%%%%%%%%%%%%%%%%%%%%%%%%%%%%%%%%%%%%%%%%

\exer{exoEntiersEmmanuel}
En multipliant l'\eqn initiale
par $a_n^{n-1}$, on obtient $a_n s$ entier sur~$\gA$. \'Ecrivons ensuite
l'\eqn initiale de la manière suivante:

\snic {
(a_n s + a_{n-1}) s^{n-1} + a_{n-2} s^{n-2} + \cdots + a_1 s + a_0 = 0,
\, \hbox { avec } b = b_{n-1}= a_n s + a_{n-1}
,}

%\sni
et considérons l'anneau $\gA[b]$. Ainsi,~$s$ annule un \pol de~$\gA[b][X]$ dont le \coe  dominant est~$b$; d'après ce qui précède, $bs$ est
entier sur~$\gA[b]$.  Mais~$b$ est entier sur~$\gA$ donc $bs = a_n s^2 +
a_{n-1} s$ est entier sur~$\gA$.  
\\
L'étape suivante consiste à écrire
l'\eqn initiale sous la forme:

\snic {
cs^{n-2} + a_{n-3} s^{n-3} + \cdots + a_1 s + a_0 = 0, \; \hbox { avec } c=  b_{n-2}= a_n s^2 + a_{n-1} s + a_{n-2}.
}

%%%%%%%%%%%%%%%%%%%%%%%%%%%%%%%%%%%%%%%%%%%%%%%%%%%%%%%%%%%%%%%%%%%%%%%%%%%

\exer{exoKroneckerTheorem}
\emph {1.}
On écrit $\lrb{1..n} \setminus I = \{i_1, i_2, \ldots\}$.  En utilisant le lemme
\ref{lemEmmanuel}, on voit que les \coes de $h_1(T) = h(T)/(T -
x_{i_1})$ sont entiers sur $\gA$, que ceux de $h_2(T) = h_1(T)/(T- x_{i_2})$
sont entiers sur $\gA[{\rm coeffs.\ de\ } h_1]$, donc entiers sur~$\gA$ et
ainsi de suite. Donc en posant $q(T) = \prod_{i'\notin I}(T - x_{i'})
\prod_{j'\notin J}(T - y_{j'})$, les \coes du \pol $h(T)/q(T)$ sont entiers
sur $\gA$. Le \coe constant de ce dernier \pol est 
$\pm a_0b_0 \prod_{i\in I}x_i \prod_{j\in J}y_j$.

\emph {2.}
Fonctions \smqs \elrsz: on a $a_i = \pm a_0 S_i(\ux)$, $b_j = \pm b_0 S_j(\uy)$,
donc~$a_ib_j$ est entier sur $\gA$.

%%%%%%%%%%%%%%%%%%%%%%%%%%%%%%%%%%%%%%%%%%%%%%%%%%%%%%%%%%%%%%%%%%%%%%%%%%%
\exer{exoAnneauCoinceBézout}
Soit $S \subseteq \gA\setminus \{0\}$ l'ensemble des \denos $b$ des
\elts de~$\gB$ écrits sous la forme $a/b$ avec $a, b \in \gA$, $b\ne
0$ et $1 \in \gen {a,b}$. C'est clairement un \moz. Pour montrer que $\gB =
\gA_S$, il suffit de vérifier que $S^{-1} \subseteq \gB$.\\
Soit $a/b \in \gB$
écrit de manière \irdz; il existe $u, v \in \gA$ tels que $1 = ua
+ vb$ si bien que $1/b = u(a/b) + v \in \gA\gB + \gA \subseteq \gB$.

%%%%%%%%%%%%%%%%%%%%%%%%%%%%%%%%%%%%%%%%%%%%%%%%%%%%%%%%%%%%%%%%%%%%%%%%%%%
\exer{exoGrellNoether} 
On veut montrer que tout anneau intermédiaire entre $\gA$ et $\gK$ est
de Prüfer. Tout \elt de $\gK$ est primitivement \agq sur n'importe quel anneau intermédiaire entre $\gA$ et $\gK$. Il reste à montrer que l'anneau intermédiaire est 
\icl pour pouvoir appliquer le \thref{th.2adpcoh}.

\emph {1.}
Si $x = a/b$, avec $a$, $b \in \gA$, il y a une matrice $\cmatrix {s & c\cr t & 1-s\cr}
\in \MM_2(\gA)$, de \lon principale pour $(b,a)$, i.e. %$s + t = 1$,
%:2018 une égalité en trop ! 
$sa = cb$ et $ta = (1-s)b$.  \\
Donc $x = c/s = (1-s)/t$ et $x \in \gA'_s \cap \gA'_t$. Réciproquement, si $x' \in \gA'_s \cap \gA'_t$,
il \hbox{y a $a'$, $b' \in \gA'$} \hbox{et $n$, $m \in \NN$} tels que $x' = a'/s^n =
b'/t^m$. Donc, pour $u$, $v \in \gA$, puisque $1/t = x/(1-s)$ on a:
$$
x' = \frac {a'}{s^n} = \frac{b' x^m}{(1-s)^m} =
\frac {ua' + vb' x^m}{us^n + v(1-s)^m}\;.
$$
Il suffit de prendre $us^n + v(1-s)^m = 1$ pour constater que $x' \in \gA'[x]$.

\emph {2.}
Soit $\gB \subseteq \gK$ une \Alg engendrée par $n$ \elts $(n\geq1)$. \\
On écrit $\gB = \gA'[x]$, où $\gA'$ est une
\Alg engendrée par $n-1$ \eltsz. D'après le point \emph{1},
il existe $s$, $t \in \gA$ tels que $\gA'[x] = \gA'_s \cap \gA'_t$.\\
Par \recuz, il existe $u_1$, \ldots, $u_k \in \gA$ tels que
$\gA' = \gA_{u_1} \cap \cdots \cap \gA_{u_k}$. \\
Alors,
$\gA'_s = \gA_{su_1} \cap \cdots \cap \gA_{su_k}$ et 
$\gA'_t = \gA_{tu_1} \cap \cdots \cap \gA_{tu_k}$, donc

\snic {
\gB = \gA_{su_1} \cap \cdots \cap \gA_{su_k} \cap
\gA_{tu_1} \cap \cdots \cap \gA_{tu_k}.
}

%\sni
\emph {3.}
Soit $\gB$ un anneau intermédiaire et $x \in \gK$ entier sur $\gB$.
Alors $x$ est entier sur une sous-\Alg de type fini, donc appartient
à cette sous-\Alg de type fini, donc à $\gB$, i.e. $\gB$ est \iclz.

\emph {4.}
Soient $x,y$ deux \idtrs sur un \cdi $\gk$ et $\gA = \gk[x,y]$.\\
Posons $\gB =\gk[x,y,(x^2+y^2)/xy]$. Alors $\gA$ est \icl mais pas $\gB$: en effet,
$x/y$ et $y/x$ sont entiers sur $\gB$ (leur somme et leur produit 
appartient à~$\gB$) mais $x/y$ et $y/x \notin\gB$ comme on le vérifie
facilement à l'aide d'un argument d'homogénéité.

%%%%%%%%%%%%%%%%%%%%%%%%%%%%%%%%%%%%%%%%%%%%%%%%%%%%%%%%%%%%%%%%%%%%%%%%%%%

\exer{exoPrimitivementAlg} 
On a $bx - a = 0$ avec $1 = ua + vb$.
\Llec vérifiera que \hbox{si $f(Y) \in \gA[Y]$} satisfait
à $f(y) = 0$, alors $f$ est multiple, dans $\gA[Y]$, 
de $bY - 2a$. Donc $\rc(f) \subseteq \gen {2a, b}$
et comme $1 \notin \gen {2a, b}$, $y$ n'est pas
primitivement \agq sur $\gA$.

%%%%%%%%%%%%%%%%%%%%%%%%%%%%%%%%%%%%%%%%%%%%%%%%%%%%%%%%%%%%%%%%%%%%%%%%%%%

\exer{exocaracPruferC} %On traite le cas intègre.
Les implications \emph {4} $\Rightarrow $ \emph {3} $\Rightarrow $ \emph {2}
et \emph {5} $\Rightarrow $ \emph {2} sont triviales. Le \thref
{thSurAdp} donne \emph {1} $\Rightarrow $ \emph {4} 
et le \tho \ref{thPruf}~\emph{4d}
\paref{thPruf} donne \emph {1}~$\Rightarrow$~\emph{5.}

  \emph
{2} $\Rightarrow $ \emph {1.} % $x \in \gA[x^2]$ entra\^\i ne que 
$x$ est  primitivement \agq sur $\gA$, on applique le \thref{th.2adpcoh}.

%%%%%%%%%%%%%%%%%%%%%%%%%%%%%%%%%%%%%%%%%%%%%%%%%%%%%%%%%%%%%%%%%%%%%%%%%%%
\exer{exocaracPruCoh}
On sait déjà que \emph{1} $\Rightarrow$ \emph{2} $\Rightarrow$ \emph{3}
et
\emph{1} $\Rightarrow$ \emph{5.} 

Montrons que \emph{3} implique que l'anneau est \ariz.
Considérons un \id à deux \gtrs arbitraire
$\fa=\gen{y_1,y_2}$ et soit $r_i$ l'annulateur \idm de $y_i$.
Considérons les \idms \ortsz:  $e=(1-r_1)(1-r_2),$
 $f=r_1(1-r_2)$, \hbox{et $g=r_2$}. On a $e+f+g=1$.
Si l'on inverse $f$ ou $g$, un des $y_i$ est nul et \hbox{l'\id $\fa$}
devient principal.
Pour voir ce qui se passe si l'on inverse $e$,
considérons les \elts \ndzs $x_1=(1-e)+ey_1$ \hbox{et $x_2=(1-e)+ey_2$}. L'\id $\fb=\gen{x_1,x_2}$ est \iv dans $\gA$.
Soient alors $u$, $v$, $w$ tels \hbox{que $ux_1=vx_2$} \hbox{et
$(1-u)x_2=wx_1$}. On multiplie par~$e$ et l'on obtient
 $uey_1=vey_2$ \hbox{et $(1-u)ey_2=wey_1$}, ce qui implique que
l'\id $\fa\gA_e=\gen{ey_1,ey_2}\gA_e$ est \lopz.

\emph{5} $\Rightarrow$ \emph{4.} Considérer d'abord $f=aX+b$, $g=aX-b$, puis
$f=aX+b$, $g=bX+a$.

\emph{4} $\Rightarrow$ \emph{3.} Soit $\fa=\gen{a,b}$, avec $a$ et $b$ \ndzsz.  Soient
$\alpha$, $\beta$ tels que $ab=\alpha a^2+\beta b^2,$ et soit $\fb=\gen{\alpha
a,\beta b}$.  On a $ab \in \fa\fb$, donc

\snic{
a^2b^2 \in\fa^2\fb^2= \gen{a^2,b^2}\gen{\alpha^2a^2,\beta^2b^2}.
}

%\sni
Montrons l'\egt $\gen{a^2b^2} = \fa^2\fb^2$, ce qui impliquera
$\fa$ \ivz.  \\
En posant $u=\alpha
a^2,$ $v=\beta b^2$, il suffit de montrer \hbox{que $u^2=\alpha^2a^4$} \hbox{et
$v^2=\beta^2b^4$} sont dans $\gen{a^2b^2}$.  Par \dfnz, $u+v = ab \in
\fa\fb$ et $uv \in \gen {a^2b^2}$.
\\
Donc $u^2 + v^2 = (u+v)^2 - 2uv \in \gen
{a^2b^2}$. Comme $u^2$, $v^2 \in\gen{u^2+v^2,uv}$, on a bien $u^2$, $v^2 \in \gen
{a^2b^2}$.

%%%%%%%%%%%%%%%%%%%%%%%%%%%%%%%%%%%%%%%%%%%%%%%%%%%%%%%%%%%%%%%%%%%%%%%%%%%

\exer{exoAnarlgb} On fait la \dem dans le cas intègre. Le cas \qi
s'en déduit par application de la machinerie \lgbe \elr usuelle.\imlg
\\
\emph{1.} Soit $M\in \Ae{n\times m}$, $p=\inf(m,n)$. 
La proposition \ref{propIddsAnar} nous donne des \ids \lops $\fa_i$
tels que

\snic{
\cD_{\gA,1}(M)=\fa_1  ,\; \cD_{\gA,2}(M)=\fa_1^2\fa_2,\; \cD_{\gA,3}(M)=\fa_1^3\fa_2^2\fa_3,\; \cD_{\gA,4}(M)=\fa_1^4\fa_2^3\fa_3^2\fa_4,\;\ldots }

%\sni
Puisque l'anneau est \lgbz, les \ids \lops $\fa_j$ sont principaux
(\plgref{thlgb3}). 
\\
Posons $\fa_j=\gen{a_j}$ et considérons la matrice
$M'\in \Ae{n\times m}$ en forme de Smith, dont les \elts diagonaux sont
$a_1$, $a_1a_2$, \ldots, $a_1a_2\cdots a_p$.
\\
 Comme dans la \dem de la proposition \ref{propIddsAnar}
l'\algo qui produit la forme réduite de Smith dans le cas local et la machinerie \lgbe  des \anars nous fournissent un
\sys comaximal $(s_1,\ldots,s_r)$ tel que, sur chaque $\gA[1/s_i]$, la matrice $M$
admet une forme réduite de Smith.
En comparant les \idds on voit que cette forme réduite peut toujours être prise
égale à $M'$ (ici intervient le fait que sur un anneau intègre,
deux \gtrs d'un \idp sont toujours associés).
\\
Ainsi, $M$ et $M'$ sont \eqves sur chaque $\gA[1/s_i]$.
On conclut par le  \plgref{thlgb1} qu'elles sont \eqvesz.

\emph{2.} Conséquence \imde du \emph{1.}

%-% ENTRE NOUS
\entrenous{Cela semble probable que cela marche pour tout \anar \lgbz. 
}
%-% Fin ENTRENOUS

%%%%%%%%%%%%%%%%%%%%%%%%%%%%%%%%%%%%%%%%%%%%%%%%%%%%%%%%%%%%%%%%%%%%%%%%%%%
\exer{exoReductionIdeal}
\emph{1.} On écrit $E = \gA x_1 + \cdots + \gA x_n$ donc $\fa E = \fa x_1 + \cdots
+ \fa x_n$. En utilisant $bE \subseteq \fa E$, on obtient une matrice $A \in
\Mn(\fa)$ telle que 

\snic{b \tra { [\,x_1\; \cdots\; x_n\,] } = A\tra {[\,x_1\; \cdots\; x_n\,]}.}

\snii
Il suffit alors de poser $d = \det(b\In - A)$.

 \emph{2.}
Si $\deg(g) \le m$, on sait que $\rc(f)^{m+1}\rc(g) = \rc(f)^m \rc(fg)$
(lemme~\ref{lemdArtin}).
En multipliant par $\rc(g)^m$, on obtient $\big(\rc(f)\rc(g)\big)^{m+1} =
\rc(fg)\big(\rc(f)\rc(g)\big)^m$.

 \emph{3.}
On a $\fb^2 = \fa\fb$, $\fb'^5 = \fa_1\fb'^4$ et $\fb'^4 = \fa_2\fb'^3$.

 \emph{4.}
Supposons $\fb^{r+1} = \fa\fb^r$. On applique la première question
avec $E = \fb^r$ et $b \in \fb$. On obtient $d = b^n + a_{1} b^{n-1}
+ \cdots + a_{n-1} b + a_n \in \Ann(\fb^r)$ avec $a_i \in \fa^i$.
\\
Comme $d \in \fb$ et $d \in \Ann(\fb^r)$, on a $d^{r+1} = 0$, qui
est une \rdi de $b$ sur $\fa$.\\
Pour la réciproque, soit $\fb$ entier sur $\fa$. Pour $b \in \fb$, en
écrivant une \rdi de $b$ sur $\fa$, on
obtient $n$ tel que $b^{n+1} \in \fa \fb^n$. Or si l'on a deux \ids $\fb_1,
\fb_2 \subseteq \fb$ avec $\fb_i^{n_i+1} \subseteq \fa\fb^{n_i}$, on
a $(\fb_1 + \fb_2)^{n_1+n_2+1} \subseteq \fa\fb^{n_1 + n_2}$.
%En effet, tout d'abord, $\fb_i^{m+1} \subseteq \fa\fb^{m}$ pour
%$m \ge n_i$; ensuite  $(\fb_1 + \fb_2)^{n_1+n_2+1} = \sum_{j+k = n_1+n_2+1}
%\fb_1^j\fb_2^k$; on a par exemple $j \ge n_1+1$ donc
%$\fb_1^j\fb_2^k \subseteq \fa\fb^{j-1}\fb_2^k \subseteq \fa\fb^{n_1+n_2}$.
\\
En utilisant un \sys fini \gtr de $\fb$, on obtient
un exposant $r$ avec l'inclusion~\hbox{$\fb^{r+1} \subseteq \fa\fb^r$}%
% donc $\fb^{r+1} = \fa\fb^r$
.

%%%%%%%%%%%%%%%%%%%%%%%%%%%%%%%%%%%%%%%%%%%%%%%%%%%%%%%%%%%%%%%%%%%%%%%%%%%

\exer{exolemNormalIcl}
On pose $\gK=\Frac\gA$.

\emph{1.} Soit $a\in\gA$ et $e_a$ l'\idm de $\gK$ tel que
$\Ann_\gK(a)=\Ann_\gK(e_a)$. L'\eltz~$e_a$ est entier sur~$\gA$, donc $e_a\in\gA$.
Et $\Ann_\gA(b)=\Ann_\gK(b)\cap\gA$ pour tout $b\in\gA$.

\emph{2. Implication directe.}  Le calcul est \imdz.

\emph{2. Implication réciproque.} Soit $a$ entier sur l'\idp $\gen{b}$ dans $\gA$.
\'Ecrivons la \rdi de $a$ sur $\gen{b}$.
$$
 a^n=b(u_{n-1}a^{n-1}+u_{n-2}ba^{n-2}+\cdots+u_0b^{n-1}).\eqno(*)
$$ 
On a $(1-e_b)a^n=0$, donc puisque $\gA$ est réduit $(1-e_b)a=0$.
On introduit l'\elt \ndz $b_1=b+(1-e_b)$. Alors l'\elt $c=a/b_1\in\gK$ est entier sur~$\gA$. En effet, l'\egt $(*)$ reste vraie en remplaçant $b$ par~$b_1$ 
et les~$u_i$ \hbox{par $e_bu_i$}, car sur la composante $e_b=1$ on obtient $(*)$ et sur la composante~\hbox{$e_b=0$} on obtient $0=0$.
\\
Donc $c$ est dans $\gA$, et $a=e_ba=e_bb_1c=bc$. 

%%%%%%%%%%%%%%%%%%%%%%%%%%%%%%%%%%%%%%%%%%%%%%%%%%%%%%%%%%%%%%%%%%%%%%%%%%%
\exer{exoEntierSurIX}
Soit $T$ une nouvelle \idtr sur $\gB$. Pour $b \in\gB$, on va utiliser
le résultat (similaire au fait \ref{fact2Entiers}): $b$ est entier sur l'\id $\fa$ \ssi $bT$ est
entier sur le sous-anneau
%
%\snic {
$\gA[\fa T] \eqdefi \gA\oplus \fa T \oplus \fa^2 T^2 \oplus \dots$
%} 
de $\gB[T]$.

\snii
Voyons le cas difficile. Soit $F\in\gB[X]$ entier sur
$\fa[X]$, on doit montrer que chaque \coe de $F$ est entier
sur $\fa$. On écrit une \rdiz:

\snic {
F^n + G_1F^{n-1} + \cdots + G_{n-1}F + G_n = 0, \quad
G_k=G_k(X) \in  (\fa[X])^k = \fa^k[X].
}

\snii
On a donc une \egt dans $\gB[X][T]$ avec des $Q_i\in\gB[X]$

\snic {
T^n + G_1T^{n-1} + \cdots + G_{n-1}T + G_n = 
(T - F) (T^{n-1} + Q_1T^{n-2} + \cdots + Q_{n-1}).
}

\snii
On remplace $T$ par $1/(TX)$ et on multiplie par $(TX)^n = TX \times
(TX)^{n-1}$. Ce qui donne

\snic {
1 + XTG_1 + \cdots + X^nT^nG_n = (1 - XTF) 
(1 + XTQ_1 + \cdots + X^{n-1}T^{n-1}Q_{n-1}).
}

\snii
On regarde maintenant cette \egt dans $\gB[T][X]$. \\
Si $b$ est un \coe de $F$,  $bT$
est un \coe en $X$ de $1 - XTF$ et $1$ est un \coe en $X$ de $1 + XTQ_1 +
\cdots + X^{n-1}T^{n-1}Q_{n-1}$. D'après le \tho de Kronecker, le produit
$bT = bT \times 1$ est entier sur l'anneau engendré par les \coes (en $X$)
du \pol $1 + XTG_1 + \cdots + X^nT^nG_n$. Mais le \coe en $X^k$ de ce
dernier \pol est dans $\gA[\fa T] = \gA \oplus \fa T \oplus \fa^2
T^2 \oplus \dots$ et donc $bT$ est entier sur $\gA[\fa T]$ et par suite $b$
est entier sur $\fa$.

%%%%%%%%%%%%%%%%%%%%%%%%%%%%%%%%%%%%%%%%%%%%%%%%%%%%%%%%%%%%%%%%%%%%%%%%%%%

\exer{exosdirindec} \emph{(Modules indécomposables)}\\
\emph{1.} Tout se passe modulo $\fa$. On considère donc l'anneau quotient
$\gB=\gA/\fa$. Alors le résultat est évident (lemme \ref{lemfacile}).

\emph{2a.} Si $M=N\oplus P$, $N$ et $P$ sont \prcs et la somme des rangs vaut $1$, donc l'un des deux est nul.

\emph{2b.} On se réfère au point \emph{1}. Si le module est décomposable,
on a $\fa\subseteq \fb$ et $\fc$ avec~$\fb$ et~$\fc$ \com \tfz. Ces \ids sont donc obtenus à partir de la \fac de $\fa$ comme deux produits partiels de cette \fcnz.
\\
Ainsi, on ne peut pas avoir $\fb$ et $\fc$ \com 
si la \fac de $\fa$  fait intervenir un seul \idemaz. 
\\
Dans le cas contraire la \fac de $\fa$ fournit deux \ids \comz~$\fb$ et $\fc$ tels \hbox{que
$\fb\fc=\fa$}. Donc $\fb+\fc=\gZ$ et $\fb\cap\fc=\fa$ ce qui donne $\gZ/\fa=\fb/\fa\oplus \fc/\fa$.
\\
En fait, si $\fa=\prod_{i=1}^{k}\fq_i=\prod_{i=1}^{k}\fp_i^{m_i}$ est la \fac de $\fa$, on obtient pr \recu sur $k$ que
  $\gZ/\fa=\bigoplus_{i=1}^{k}\fq_i/\fa$.

\emph{2c.} Résulte des considérations précédentes et du \tho
de structure des \mpfs sur un \dDkz.

\emph{3.} L'unicité peut s'énoncer comme suit: si $M$ s'écrit de deux manières comme somme de modules indécomposables, il y a un \auto de $M$
qui envoie les modules de la première \dcn sur ceux de la seconde. 
\\
Si un \mpf et de torsion  $M$ est décomposé en somme directe de modules indécomposables, chaque terme de la somme est lui-même \pf et de torsion. Et donc de la forme $\gZ/\fp^{m}$ d'après le point \emph{1.} 
\\
Par le \tho des restes chinois on se ramène au cas où un seul \idema intervient dans la somme directe, et l'unicité résulte alors du \thref{prop unicyc}.
\\
Notons aussi que dans le cas d'un anneau principal à \facz,
l'unicité est valable  pour la \dcn de tout \mpfz. 

%%%%%%%%%%%%%%%%%%%%%%%%%%%%%%%%%%%%%%%%%%%%%%%%%%%%%%%%%%%%%%%%%%%%%%%%%%%
%:  sols pbs
%%%%%%%%%%%%%%%%%%%%%%%%%%%%%%%%%%%%%%%%%%%%%%%%%%%%%%%%%%%%%%%%%%%%%%%%%%%

%%%%%%%%%%%%%%%%%%%%%%%%%%%%%%%%%%%%%%%%%%%%%%%%%%%%%%%%%%%%%%%%%%%%%%%%%%%

\prob{exoArithInvariantRing}  
Ci-après  le mot \gui{\lotz} signifie \gui{après \lon en des \ecoz}.

\emph {1.}
L'\id $\fa$ est \lopz, donc puisque $\gA$ est normal, \lot \iclz, donc
il est \icl (\plgref{plcc.normal}). On termine avec le lemme~\ref{lemthKroicl} (variante du \tho de \KROz).

\emph {2.}
Si $x \in \fa\gB \cap \gA$, alors $x$ est entier sur l'\id $\fa$ 
(lying over, lemme \ref{lemLingOver2}) donc
dans~$\fa$ d'après la question précédente.

\emph {3a.}
Si $a = \rN_G(b)$, on a $\rN_G(b\gB) = a\gB \cap \gA = a\gA$.

\emph {3b} et \emph{3c.}
L'\itf $\fa = \rc_\gA(h)$ est \lopz, donc $\rc_\gB(h) = \fa\gB$ est un
\id\lop de $\gB$. D'après la première question, on a~:

\snic {
\prod_\sigma \rc_\gB(h_\sigma) = \rc_\gB(h),
\quad \hbox {i.e.} \quad \rN'_G(\fb) = \fa\gB.
}

%\sni
Et d'après la question \emph {2}, $\fa =  \rN_G(\fb)$.
Ensuite on note que
$$\preskip.3em \postskip.3em 
\rN_G(\fb_1\fb_2)\gB = \rN'_G(\fb_1\fb_2) = \rN'_G(\fb_1)\rN'_G(\fb_2)
= \rN_G(\fb_1)\rN_G(\fb_2)\gB, 
$$
d'où le résultat en prenant l'intersection avec $\gA$.

\emph {3d.} Cela résulte du point 2. et des deux faits suivants.
\\
$\bullet$ Si $b\in\gB$ est \ndz alors $a = \rN_G(b)\in\gA$ est \ndz dans $\gA$:
en effet, c'est un produit d'\elts \ndzs dans $\gB$ donc il est \ndz dans $\gB$.
\\
$\bullet$  Si $a \in \gA$ est \ndz
dans $\gA$ alors il est \ndz dans $\gB$. Soit en effet $x \in \gB$ 
tel que $ax = 0$, on veut montrer que $x=0$. On considère le \pol

\snic{\rC{G}(x)(T)  =    \prod_{\sigma\in G}\big(T-\sigma(x)\big).}

%\sni
Comme $a\sigma(x)=0$ pour chaque $\sigma$, les \coes de $\rC{G}(x)(T)$
sont annulés par $a$ donc nuls, à l'exception du \coe dominant. 
Ainsi $x^{\abs{G}}=0$, or $\gB$ est normal donc réduit.

\emph {4.}
Soit $\gk(x,y) = \Frac \gk[x,y]$.  On va utiliser le fait que $(1,y)$ est une
$\gk[x]$-base de $\gk[x,y]$; c'est aussi une $\gk(x)$-base de $\gk(x,y)$ et
l'extension $\gk(x,y)\sur{\gk(x)}$ est galoisienne de groupe $\gen {\sigma}$
où $\sigma : \gk(x,y) \to \gk(x,y)$ est le $\gk(x)$-\auto involutif qui
réalise $y \mapsto -y$.
Montrons que $\gk[x,y]$ est la \cli de $\gk[x]$ dans
$\gk(x,y)$. Soit $z = u(x) + yv(x) \in\gk(x,y)$ entier sur
$\gk[x]$. Alors $z+\sigma(z) = 2u$ et $z\sigma(z) = u^2 - fv^2$ sont dans
$\gk(x)$ et entiers sur $\gk[x]$ donc dans $\gk[x]$. On a \hbox{donc $fv^2 \in
\gk[x]$}. En utilisant le fait que $f$ est \spbz, on voit que $v \in
\gk[x]$. Bilan: $z \in \gk[x,y]$. Donc $\gk[x,y]$ est \iclz.
On applique ce qui précède avec $\gA = \gk[x]$, $\gB = \gk[x,y]$, 
$G = \gen {\sigma}$.

%%%%%%%%%%%%%%%%%%%%%%%%%%%%%%%%%%%%%%%%%%%%%%%%%%%%%%%%%%%%%%%%%%%%%%%%%%%
\prob{exoFullAffineMonoid} 
\emph {2a)}
Posons $a = \sum_\alpha a_\alpha \ux^\alpha$, $b = \sum_\beta b_\beta
\ux^\beta$. 
\\
On doit montrer que $\beta \in M$ pour chaque $\beta$ tel que
$b_\beta \ne 0$. On peut supposer~$b$ non nul. Soient $a_\alpha\ux^\alpha$ le
monôme dominant de $a$ pour l'ordre lexicographique \hbox{et $b_\beta\ux^\beta$}
celui de $b$. Le monôme dominant de $ab$ est $a_\alpha
b_\beta\ux^{\alpha+\beta}$, donc $\alpha+\beta \in M$.
\\
Comme $\alpha\in M$ et
que $M$ est plein, on a $\beta\in M$.  On recommence ensuite en remplaçant
$b$ par $b' = b-b_\beta\ux^\beta$ qui satisfait $ab' \in \gk[\ux]$.
On obtient $b' \in \gk[\ux]$ et finalement $b \in \gk[\ux]$.

%%%%%%%%%%%%%%%%%%%%%%%%%%%%%%%%%%%%%%%%%%%%%%%%%%%%%%%%%%%%%%%%%%%%%%%%%%%

\prob{exoBaseNormaleAlInfini}
\emph {1.}
Si $A = (a_{ij})$, alors $\det A = \sum_{\sigma \in \rS_n} a_{\sigma(1)1} \cdots
a_{\sigma(n)n}$ et

\snic {
v(a_{\sigma(1)1} \cdots a_{\sigma(n)n}) \ge v(A_1) + \cdots + v(A_n)
.}

%\sni
On en déduit que $v(\det A) \ge v(A_1) + \cdots + v(A_n)$.

\emph {2.}
Pour la matrice donnée en exemple: on a $\det(A) = \pi^2 - \pi \ne 0$.
\\
Mais
$\ov A = \cmatrix {0 & 0\cr 1 & 1\cr}$ n'est pas inversible.  En réalisant
$A_1 \leftarrow A_1 - A_2$, on obtient l'\egt $A' = \cmatrix {\pi^2 - \pi &\pi \cr 0
& 1\cr}$ et cette fois $\ov {A'} = \cmatrix {1 & 0\cr 0 & 1\cr}$ est
inversible.

Voici la méthode générale: si $\det \ov A \ne 0$, $A$ est
$\gA_\infty$-réduite et c'est terminé. Sinon, il y a des $\lambda_1,
\ldots, \lambda_n \in \gk$, non tous nuls, tels que $\lambda_1 \ov {A_1} +
\cdots + \lambda_n \ov {A_n} = 0$. On considère une colonne $A_j$
avec $\lambda_j \ne 0$ et $v(A_j)$ minimum; pour simplifier, on
peut supposer que c'est $A_1$ et que $\lambda_1 = 1$ (quitte à 
diviser la relation par $\lambda_1$); on réalise alors l'\op
\elr:

\snic {
A_1 \leftarrow A'_1 = A_1 + \sum_{j=2}^n \lambda_j \pi^{v(A_1)-v(A_j)} A_j
.}

%\sni
Dans cette somme, en ne faisant intervenir que
les $A_j$ pour lesquels $\lambda_j \ne 0$, chaque exposant de $\pi$ est $\ge
0$. Il s'agit donc d'une \op $\gk[\pi^{-1}]$-\elr sur les colonnes, i.e.
$\gk[t]$-\elr, et  l'on ne change pas le $\gk[t]$-module engendré par les
colonnes. Par ailleurs, $v(A'_1) > v(A_1)$; en effet, (en se souvenant que
$\lambda_1 = 1$):

\snic {
A'_1 / \pi^{v(A_1)} = s \eqdf {\rm def} 
\som_{\lambda_j \ne 0} \lambda_j A_j / \pi^{v(A_j)}
,}

%\sni
et $v(s) > 0$ puisque par hypothèse 
$\sum_{\lambda_j \ne 0} \lambda_j \ov {A_j}= 0$. 
On itère ce processus qui finit par s'arrêter car à chaque
étape, la somme $\sum_j v(A_j)$ croit strictement tout en
étant bornée par $v(\det A)$, invariant par les opérations ci-dessus.

\emph {3.}
Soit $y = Px$, i.e. $y_i = \sum_j p_{ij} x_j$; on a $v(p_{ij}) \ge 0$,
$v(x_j) \ge v(x)$ donc $v(y_i) \ge v(x)$ puis $v(y) \ge v(x)$.
Par symétrie, $v(y) = v(x)$. Le reste ne pose pas plus de
difficultés.

\emph {4.}
$A$ est $\gA_\infty$-réduite \ssi tout coefficient diagonal (nécessairement
non nul) divise (au sens $\gA_\infty$) tous les coefficients de sa colonne.

\emph {5.}
Quitte à remplacer $A$ par $AQ$ avec $Q \in \GL_n(\gA)$ convenable, on peut
supposer que $A$ est $\gA_\infty$-réduite. On va réaliser des \ops $A
\leftarrow PA$ avec $P \in \GL_n(\gA_\infty)$ (i.e. considérer le
$\gA_\infty$-réseau engendré par les lignes de $A$), ce qui ne modifie pas
le \crc $\gA_\infty$-réduit de $A$.  Il existe $P \in
\GL_n(\gA_\infty)$ telle que $PA$ soit triangulaire supérieure et l'on
remplace $A$ par $PA$.  Soient $L_1, \ldots, L_n$ les lignes de $A$; on
réalise alors l'\op $\gA_\infty$-\elr

\snic {
L_1 \leftarrow L_1 - {a_{12} \over a_{22}} L_2
\qquad  \hbox {rappel : }  a_{22} \divi_{\gA_\infty} a_{12}
,}

%\sni
ce qui amène un $0$ en position $a_{12}$ (et la nouvelle matrice est
toujours triangulaire \hbox{et $\gA_\infty$-réduite)}. On continue pour annuler tous
les coefficients de la première ligne (sauf $a_{11}$); on peut ensuite
passer à la deuxième ligne et ainsi de suite de façon à obtenir une
matrice diagonale (en utilisant constamment le fait que dans une matrice
triangulaire $\gA_\infty$-réduite, chaque coefficient diagonal
$\gA_\infty$-divise tous les coefficients de sa colonne).  Comme $\gA_\infty$
est un \adv discrète, on peut faire en sorte que la matrice diagonale finale
obtenue soit $\Diag(\pi^{d_1}, \ldots,
\pi^{d_n})$ avec $d_i \in \ZZ$.

\emph {6a.}
Soit $\und\varepsilon$ une $\gA$-base de $E$, $\und{\varepsilon'}$ une
$\gA_\infty$-base de $E'$  et $A = \Mat_{\und\varepsilon,
\und{\varepsilon'}}(\Id_L)$.  Il existe alors $P \in \GL_n(\gA_\infty)$ et $Q
\in \GL_n(\gA)$ telles que $PAQ = \Diag(t^{-d_1}, \ldots, t^{-d_n})$.  Soient
$\ue$ et $\und{e'}$ définies par $\Mat_{\ue,\und\varepsilon}(\Id_L) = Q$,
$\Mat_{\und{\varepsilon'},\und{e'}}(\Id_L) = P$. 
\\
Alors $\ue$ est une
$\gA$-base de $E$, $\und{e'}$ une $\gA_\infty$-base de $E'$ et $e_i = t^{-d_i}
e'_i$.

\emph {6b.}
Puisque $t^je_i = t^{j-d_i} e'_i$, il est clair 
que $t^j e_i \in E \cap E'$ pour $0 \le j \le d_i$. 
Réci\-proquement, soit $y \in E \cap E'$ que l'on 
écrit 

\snic {
y = \sum_i a_ie_i = \sum_i a'_i t^{d_i} e_i, \quad\hbox{avec }
a_i \in \gA \hbox{ et } a'_i \in \gA_\infty
,}

%\sni
et donc $a_i = a'_i t^{d_i}$. \\
Si $d_i < 0$, on obtient $a_i = a'_i = 0$,
et si $a_i \ne 0$,  $0 \le \deg a_i \le d_i$. D'où la
$\gk$-base annoncée.

\emph {7.}
Tout d'abord $\gk' = \gB \cap \gB_\infty$, donc $\gB$ et $\gB_\infty$ sont des
$\gk'$-\evcsz. Montrons que chaque $r_i \in \gA_\infty$ et que de plus, si
$e_i \notin \gB_\infty$, alors  $v(r_i) > 0$, \hbox{i.e. $\deg(r_i) < 0$}. 
Si $e_i \in
\gB_\infty$, on a $e_i \in \gB \cap\gB_\infty = \gk'$, donc
aussi $e_i^{-1} \in \gk'$; par suite $r_i = e_i^{-1} (r_ie_i)
\in \gB_\infty$ donc $r_i \in \gB_\infty \cap \gK = \gA_\infty$.
\\
Si $e_i \notin \gB_\infty$, on écrit $e_i = r_i^{-1}(r_ie_i)$,
\egt qui prouve que $r_i^{-1} \notin \gA_\infty$ (ne pas oublier
que $r_ie_i \in \gB_\infty$) donc $v(r_i^{-1}) < 0$, i.e.
$v(r_i) > 0$.

Soit maintenant $c \in \gk'$ que l'on écrit dans la
$\gA$-base $(e_i)$ et la $\gA_\infty$-base $(r_ie_i)$

\snic {
c = \sum_i a_ie_i = \sum_i a'_i r_ie_i, \quad a_i \in \gA, \quad
a'_i \in \gA_\infty, \quad a_i = a'_i r_i
.}

%\sni
Pour les $i$ tels que $e_i \in \gk'$, comme $r_i \in \gA_\infty$, on a $a_i =
a'_ir_i \in \gA \cap \gA_\infty = \gk$. Reste à voir que pour $e_i \notin
\gk'$, $a_i = 0$; l'\egt $a_i = a'_ir_i$ et le fait que $a_i \in \gA$,
$a'_i \in \gA_\infty$ et $\deg(r_i) < 0$ entraînent alors $a_i = a'_i = 0$.
Bilan: les $e_i$ qui sont dans $\gk'$ forment une $\gk$-base de $\gk'$.

\emph {8.}
En posant $i = y/x$, on a $i^2 = -1$ et

\snic {
\cmatrix {1 & x\cr 0 & 1\cr} \cmatrix {1\cr i} = \cmatrix {y+1\cr i}
.}

%\sni
La matrice de gauche est de déterminant $1$, donc $(1,i)$ et
$(y+1,i)$ sont deux base du même $\gA$-module. Mais
$y+1$ n'est pas entier sur $\gA_\infty$ (car $x$ est entier
\hbox{sur $\gk[y] = \gk[y+1]$} et n'est pas entier sur $\gA_\infty$).
La base $(1, i)$ est normale a l'infini mais pas la base $(y+1,i)$.

%%%%%%%%%%%%%%%%%%%%%%%%%%%%%%%%%%%%%%%%%%%%%%%%%%%%%%%%%%%%%%%%%%%%%%%%%%%
\prob{exoHyperEllipticFunctionRing} 
\emph {1.}
Soit $z = y-v$, $(1,z)$ est une $\gA$-base de $\gB$ et $\gA u
\cap \gA z = \{0\}$. Pour montrer que $\fb_{u,v} = \gA u \oplus \gA z$ est un
\idz, il suffit de voir que $z^{2}\in\fb_{u,v}$. \\
Or $y^2 = (z+v)^2 = z^2 +
2vz + v^2$, i.e.  $z^2 + 2vz + uw = 0$.

\emph {2.}
Comme $(1,z)$ est une $\gA$-base de $\gB$ et $(u, z)$ une
$\gA$-base de $\fb_{u,v}$, on obtient l'\egt $\gA \cap \fb_{u,v} = u\gA$. D'autre part,
tout \elt de $\gB$ est congru modulo $z$ à un \elt de $\gA$, donc $\gA \to
\gB\sur{\fb_{u,v}}$ est surjective de noyau $u\gA$. \\
 La matrice $M$ de $(u,
y-v)$ sur $(1, y)$ est $M = \cmatrix {u & -v\cr 0 & 1}$ avec $\det(M) = u$,
ce qui donne $\rN(\fb_{u,v}) = u\gA$. On voit \egmt que le contenu de $\fb_{u,v}$ est
$1$.  Les autres points sont faciles.

\emph {3.}
On a $\fb_{u,v,w} = \gA u \oplus \gA z$, $\fb_{w,v,u} = \gA w
\oplus \gA z$. Le produit de ces deux \ids est engendré (en tant
qu'\id ou \Amoz) par les 4 \elts $uw$, $uz$, $wz$, $z^2$, tous multiples de $z$ (car
$z^2 + 2vz + uw = 0$). Il suffit donc de voir que

\snic {
z \in \gen {uw, uz, wz, z^2}_\gB = \gen {uw, uz, wz, 2vz}_\gB =
\gen {uw, uz, wz, vz}_\gB
.}

%\sni
Or $v^2 - uw = f$ est séparable, donc  $1 \in \gen {u,w,v}_\gA$, et $z \in \gen {uz,wz,vz}_\gB$.

Quant à $\fb_{u,-v}$ c'est $\gA u \oplus \gA \ov z$ avec $z\ov z = uw$ et
$z + \ov z = -2v$.  Le produit $\pi$ des deux \ids $\fb_{u,v}$ et
$\fb_{u,-v}$ est égal à $\gen{u^2, u\ov z, uz, z\ov z}$, avec $z\ov z = uw$, donc $\pi\subseteq \gen{u}$. Enfin $-2uv = uz + u\ov z \in \pi$ et donc
$\pi\supseteq \gen{uv, u^2, uw}=u\gen {v, u, w}=\gen{u}$.

Enfin, avec $u = u_1u_2$, on a $\fb_{u_1,v}\fb_{u_2,v} = \gA u + \gA u_1z +
\gA u_2z + \gA z^2$ clairement inclus dans $\gA u + \gA z = \fb_{u,v}$. Comme
$z^2 + 2vz + uw = 0$ on obtient

\snic {\mathrigid1mu
\gA u + \gA u_1z + \gA u_2z + \gA z^2 =
\gA u + \gA u_1z + \gA u_2z + \gA vz = 
\gA u + (\gA u_1 + \gA u_2 + \gA v)z
.}

%\sni
D'où $\fb_{u_1,v}\fb_{u_2,v} = \fb_{u,v}$; en effet, 
$1\in\gen{u_1,u_2, v}_\gA$ car $v^2 - u_1u_2w = f$ est
séparable, donc $\gen{u_1,u_2, v}_\gA z = \gA z$.

\emph {4a.}
Soit $\fb$ un \itf non nul de $\gB$. Comme \Amo libre de rang~$2$, il
admet  une $\gA$-base $(e_1, e_2)$ et nous notons $M = \cmatrix {a & b\cr 0 &d\cr}$ la matrice de $(e_1, e_2)$ sur $(1,y)$.
On écrit que $\fb$ est un \idz, i.e. $y\fb \subseteq \fb$:
l'appartenance $ye_1 \in \gA e_1 \oplus \gA e_2$ donne $a$ multiple de $d$
et l'appartenance $ye_2 \in \gA e_1 \oplus \gA e_2$ donne $b$ multiple de $d$.
En définitive, $M$ est de la forme $M = d \cmatrix {u & -v\cr 0 &1\cr}$
et l'on obtient $\fb = d\fb_{u,v}$. On voit que $\gen {d}_\gA$ est le
contenu de $\fb$, et $d$ est unique si l'on impose $d$ unitaire.

%%Soit $\fb$ un \itf non nul de $\gB$; c'est un \Amo libre de rang $2$; si l'on
%%note $M \in \MM_2(\gA)$ la matrice d'une $\gA$-base quelconque de $\fb$ dans
%%une $\gA$-base quelconque de $\gB$, alors l'\idd $\cD_1(M)$ est le $1$-Fitting
%%$\cF_1(\gB\sur\fb)$ du \Amo $\gB\sur\fb$; il est engendré par le pgcd
%%(unitaire) des 4 coefficients de $M$.  Par exemple, $\cF_1(\gB\sur{\fb_{u,v}})
%%= \gen {1}_\gA$

\emph {4b.}
On a vu que $\fb = d\fb_{u,v}$ donc $\ov\fb = d\fb_{u,-v}$ puis $\fb\ov\fb =
d^2u\gB$. \\
Mais on aussi $\rN(\fb) = d^2u\gA$ car $d\cmatrix {u & -v\cr 0 &
1\cr}$ est la matrice d'une $\gA$-base de~$\fb$ sur une $\gA$-base de $\gB$.
On en déduit que $\fb\ov\fb = \rN(\fb)\gA$. Alors, pour deux \ids non nuls
$\fb_1, \fb_2$ de $\gB$:

\snic {
\rN(\fb_1\fb_2)\gB = \fb_1\fb_2\ov {\fb_1\fb_2} =
\fb_1\ov{\fb_1} \fb_2\ov{\fb_2} = \rN(\fb_1)\rN(\fb_2)\gB
,}

%\sni
d'où $\rN(\fb_1\fb_2) = \rN(\fb_1)\rN(\fb_2)$ puisque les trois \ids
sont des \idps de $\gA$.

\emph {4c.}
Tout d'abord, si $\fb$ est un \itf non nul de $\gB$, il contient
un \elt \ndz $b$ et $a = \rN(b) = b\wi b$ est un \elt
\ndz de $\fb$ contenu dans $\gA$. On a alors une
surjection $\gB\sur{a\gB} \twoheadrightarrow \gB\sur{\fb}$
et comme $\gB\sur{a\gB}$ est un $\gk$-\evc de dimension finie, il
en est de même de $\gB\sur{\fb}$.

Si $d \in \gA \setminus \{0\}$, on a une suite exacte:

\snic {
0 \to \gB\sur{\fb'} \simeq d\gB\sur{d\fb'} \to \gB\sur{d\fb'} \to
\gB\sur{d\gB} \to 0
.}

%\sni
On en déduit $\deg(d\fb') = \deg(\fb') + \deg(d\gB) = \deg(\fb') + \deg(d^2)$.
En particulier, pour $\fb' = \fb_{u,v}$ et $\fb = d\fb_{u,v}$, on obtient 

\snic {
\deg(\fb) = \deg(u) + \deg(d^2) = \deg \rN(\fb)
.}

%\sni
Ceci montre que $\deg$ est additif.

\emph {5.}
On fournit d'abord un algorithme de réduction de $(u,v)$ vérifiant $v^2
\equiv f \bmod u$. Quitte à remplacer $v$ par $v \bmod u$, on peut supposer
$\deg v < \deg u$. Si $\deg u \le g$, alors, en rendant $u$
\monz, $(u,v)$ est réduit. Sinon, avec $v^2 - uw =
f$ montrons que $\deg w < \deg u$; cela permettra de considérer $\wi u := w$,
$\wi v := (-v) \bmod \wi u$, ayant la \prt $\fb_{u,v} \sim \fb_{\tilde u, \tilde v}$
et d'itérer le processus $(u,v) \leftarrow (\wi u, \wi v)$ jusqu'à
l'obtention de l'inégalité $\deg u \le g$.  Pour montrer $\deg u > g
\Rightarrow \deg w < \deg u$, on considère les deux cas suivants; ou bien
$\deg(uw) > 2g+1 = \deg f$, auquel cas l'\egt $f + uw = v^2$ fournit $\deg(uw)
= 2\deg v < 2\deg u$ donc $\deg w < \deg u$; ou bien $\deg(uw) \le 2g+1$,
auquel cas $\deg w \le 2g+1 - \deg u < 2g+1 - g$ donc $\deg w \le g < \deg u$.

Tout \id $\fb_{u,v}$ est donc associé à un \id réduit et comme tout \itf
non nul $\fb$ de $\gB$ est associé à un \id $\fb_{u,v}$, $\fb$ est donc
associé à un \id réduit.

\emph {6a.}
Soit $w$ vérifiant $v^2 - uw = f = y^2$; comme $(u,v)$ est réduit, on
a:

\snic {
\deg v < \deg u \le g < g+1 \le \deg w
\quad \hbox {et} \quad
\deg u + \deg w = 2g + 1
.}

%\sni
Posons $y' = y-v$ et $z = au + by'$ avec $a, b \in \gA$. 
\\
On a $y' + \ov {y'} = -2v$, $y \ov {y'} = -(y^2 - v^2) =
uw$, donc:

\snic {
\rN(z)= z\ov z = a^2u^2 + aub(y' + \ov {y'}) +
b^2 y'\ov {y'} = u(a^2u - 2vab + b^2w)
,}

%\sni
d'où 
$\rN(z)/\rN(\fb_{u,v}) = \rN(z)/u = a^2u - 2vab + b^2w$,
\pol dont il s'agit de minorer le degré. Notons le cas particulier $b = 0$
(donc $a \ne 0$) auquel cas $\rN(z)/u = a^2u$, de degré $2\deg a + \deg u
\ge \deg u$. On voit ici que l'\egt $\deg (\rN(z)/u) = \deg u$ est atteinte
\ssi $\deg a = 0$, i.e. \ssi $z \in \gk\eti u$. 
\\
Il y a aussi le cas
particulier $a = 0$ (donc $b \ne 0$) auquel cas $\rN(z)/u = b^2 w$,
qui est de degré
$2\deg b + \deg w > \deg u$.

Il reste donc à montrer que pour $a \ne 0$, $b \ne 0$, on a $\deg(\rN(z)/u)
> \deg u$. On introduit  $\alpha = \deg a \ge 0$, $\beta = \deg b \ge 0$
et:

\snic {\mathrigid1mu
d_1 = \deg(a^2u) = 2\alpha + \deg u, \
d_2 = \deg(vab) = \alpha + \beta + \deg v, \
d_3 = \deg(b^2w) = 2\beta + \deg w
.}

%\sni

On a $d_1 + d_3 \equiv \deg u + \deg w = 2g+1 \bmod 2$ donc $d_1 \ne d_3$. Aussi, $\alpha \ge \beta \Rightarrow d_1 > d_2$ et enfin  $\beta \ge
\alpha \Rightarrow d_3 > \max(d_1,d_2)$.  
\\
Si $d_3 > \max(d_1,d_2)$, alors
$\deg (\rN(z)/u) = d_3 \ge \deg w > \deg u$. Si $d_3 \le \max(d_1,d_2)$, alors
$\alpha > \beta$, donc $d_1 > d_2$, puis $d_1 > \max(d_2,d_3)$. 
\\
On a donc $\deg
(\rN(z)/u) = d_1 = 2\alpha + \deg u \ge 2 + \deg u > \deg u$.

\emph {6b.}
On a $\fb' = d\fb_{u_1,v_1}$ et $\deg(\fb') = 2\deg(d) + \deg(\fb_{u_1,v_1})$.
On peut donc supposer \hbox{que $d = 1$}.  On a $c$, $c_1 \in \gB \setminus \{0\}$
avec $c \fb_{u,v} = c_1\fb_{u_1,v_1}$, que nous notons $\fb$. On a $\rN(\fb) = u\rN(c) = u_1\rN(c_1)$.  Le degré minimum des
$\rN(z)/\rN(\fb)$ pour $z \in
\fb \setminus \{0\}$ est $\deg u$ et il est atteint uniquement pour $z \in
\gk\eti cu$.  \\
Pour $z = c_1u_1 \in \fb$, on a $\rN(z) = u_1^2 \rN(c_1)$
donc $\rN(z)/\rN(\fb) = %\big(u_1^2 \rN(c_1)\big) \big/ \big(u_1\rN(c_1)\big) 
\fraC{u_1^2 \rN(c_1)}{u_1\rN(c_1)}= u_1$.
On a donc $\deg u_1 \ge \deg u$, i.e. $\deg(\fb_{u_1,v_1}) \ge \deg(\fb_{u,v})$.
L'\egt n'est possible que \hbox{si $c_1u_1 \in \gk\eti cu$}. Dans ce cas,
$u\fb_{u_1,v_1} = u_1\fb_{u,v}$. Puisque le contenu de 
$u\fb_{u_1,v_1}$ est~$u$, et que celui de $u_1\fb_{u,v}$ est $u_1$,
l'\egt précédente entraîne $u = u_1$ puis
$v = v_1$.

\emph {7a.}
On a $F'_X(X,Y) = -f'(X)$, $F'_Y(X,Y) = 2Y$. \\
Comme 
$\car(\gk)\neq 2$, on obtient $f(X)\in\gen {F, F'_X, F'_Y}$, puis  $1\in\gen {F, F'_X, F'_Y}$.

\emph {7b.}
On réalise le \cdv $\gx = 1/x$ dans 

\snic {y^2 = f(x) = x^{2g+1} + a_{2g} x^{2g}
+ \cdots + a_1 x + a_0,}

et l'on multiplie par $\gx^{2g+2}$ pour obtenir:

\snic {
\gy^2 = \gx + a_{2g}\gx^2 + \cdots + a_0\gx^{2g+2} = \gx \big(1 + \gx h(\gx)\big)
\quad \hbox {avec} \quad \gy = y\gx^{g+1}
.}

%\sni

Bilan: le \cdv $\gx = 1/x$, $\gy = y/x^{g+1}$ donne $\gk(x) =
\gk(\gx)$ \hbox{et $\gk(x,y) = \gk(\gx, \gy)$}. Et $\gy$ est entier sur $\gk[\gx]$, a
fortiori sur $\gA_\infty$. 
\\
On pose $\gB_\infty =
\gk[\gx,\gy]_{\gen {\gx,\gy}}$; dans ce localisé, on a $\gen {\gx, \gy} =
\gen {\gy}$ puisque $\gx = {\gy^2 \over 1 + \gx h(\gx)}$.  Conclusion:
$\gB_\infty$ est un \adv discrète d'uniformisante $\gy$.

Enfin, soit $\gW$ un \adv pour $\gk(x,y)$ contenant $\gk$.  
\\
Si $x \in \gW$, alors
$\gk[x] \subset \gW$. Alors  $y$, entier sur $\gk[x]$, est dans $ \gW$, donc $\gB  \subset \gW$.\\
Si $x \notin \gW$,
on a $x^{-1} \in \fm(\gW)$, donc $\gA_\infty = \gk[x^{-1}]_{\gen {x^{-1}}}
\subset \gW$, et $\gW = \gB_\infty$.

%%%%%%%%%%%%%%%%%%%%%%%%%%%%%%%%%%%%%%%%%%%%%%%%%%%%%%%%%%%%%%%%%%%%%%%%%%%

\prob{exoTrifolium} 
On note $\vep$ l'\iv défini par \framebox[1.1\width][c]{$\vep = \beta-\alpha$}.

\emph {1.}
On décompose $G$ et $H$ en composantes \hmgs $G_i, H_j$:

\snic {
G = G_a + \cdots + G_b, \  a \le b, \quad
H = H_c + \cdots + H_d, \  c \le d
.}

%\sni
La composante \hmg basse de $GH$, de degré $a+c$, est $G_aH_c$ tandis
que la composante \hmg haute de $GH$, de degré $b+d$, est $G_bH_d$.
On en déduit que $a+c=N$, $b+d=N+1$; on ne peut pas avoir à la
fois $a < b$ et $c < d$ (car on aurait alors $a+c+2 \le b+d$, i.e.
$N+2 \le N+1$). Si $a=b$, alors $G$ est \hmgz, si $c=d$ c'est $H$.
Supposons $F_N$, $F_{N+1}$ premiers entre eux et soit une \fcn
$F = GH$; par exemple, $G$ est \hmg de degré $g$; on en
déduit que $H = H_{N-g} + H_{N+1-g}$ et que $F_N = GH_{N-g}$,
$F_{N+1} = GH_{N+1-g}$: $G$ est un facteur commun à $F_N, F_{N+1}$,
donc $G$ est \ivz. La réciproque est facile.

Les \pols $(X^2 + Y^2)^2$ et $\alpha X^2Y + \beta Y^3 = Y(\alpha X^2 + \beta
Y^2)$ sont premiers entre eux \ssi les \pols $X^2 + Y^2$ et $\alpha X^2 +
\beta Y^2$ le sont i.e. \ssi $\alpha \ne \beta$.

\emph {2.}
\Llec vérifiera que $(0,0)$ est le seul point singulier; on a
le résultat plus précis:

\snic {
\vep^2 X^5, \vep^2 Y^5 \in \gen {F, F'_X, F'_Y}
}

\emph {3.}
On pose $Y = TX$ dans $F(X,Y)$ et l'on obtient $F(X, TX) =
X^3 G(X,T)$ avec

\snic {
G(X, T) = XT^4 + \beta T^3 + 2X T^2 + \alpha T + X
%, \qquad  G(0,T) = T(\beta T^2 + \alpha)
.}

%\sni
Le \pol $G$ est primitif (en $T$) et $(x=0,t=0)$ est un point simple de
la courbe $G = 0$. Avec $a_4 = x$, $a_3 = \beta$, $a_2 = 2x$, $a_1 = \alpha$,
$a_0 = x$, on considère les entiers d'Emmanuel:
$$
\preskip.4em \postskip.4em 
b_4 = a_4,\quad  b_3 = a_3 + tb_4,\quad b_2 = a_2 + tb_3,\quad
b_1 = a_1 + tb_2 
. 
$$
Ainsi, $b_4 = x$, $b_3 = \beta+y$ et $b_2 = 2x + (\beta+y)y/x$.

Il est clair que $a_4$, $a_3$, \dots, $a_0 \in \sum_i \gA b_i + \sum_i \gA tb_i$.
Comme $a_3-a_1 = \vep$ est \ivz,  il y a des $u_i, v_i \in \gA$ tels que
$1 = \sum_i u_i b_i + \sum_i v_i tb_i$. On écrit formellement
(sans se soucier de la nullité d'un $b_i$):
$$\preskip.0em\postskip.2em \ndsp
t = {b_1t \over b_1} = \cdots = {b_4t \over b_4} =
{\sum_i v_ib_it \over \sum_i v_i b_i} = 
{\sum_i u_ib_it \over \sum_i u_i b_i}.
$$
Ainsi, $t= y/x = a/b = c/d$ avec $a$, $b$, $c$, $d \in \gB$ et $a+d = 1$. 
\\
Les
\egts $by = ax$, $dy = cx$,  $a+d = 1$ sont celles convoitées. 
\\
On obtient ainsi
$\fq\gen{x,y}_\gB = \gen {x}_\gB$ avec $\fq =\gen {d,b}_\gB$. Ici en posant:

\snic {
a = b_2t - b_4t, \quad  b = b_2 - b_4,\quad
c = b_3t - b_1t, \quad  d = b_3 - b_1 
,}

%\sni

on a $\vep = a+d$. En posant $g_0 = 1/(1+t^2)$, $g_1 = tg_0$, on
trouve $b = \vep g_1$, $d = \vep g_0$, donc $\fq = \gen
{g_0,g_1}_\gB$.  On va montrer (question \emph {5}) que $\gB = \gk[g_0,g_1]$,
donc $\gB\sur\fq = \gk$.

\emph {4.}
Une idée \gmq conduit à l'\egt $\gk(t) =
\gk(x,y)$. C'est  le paramétrage du trifolium. Le \pol définissant la courbe est de degré $4$
et l'origine est un point singulier de multiplicité $3$. Donc une droite
rationnelle passant par l'origine recoupe la courbe en un point rationnel.
Algébriquement, cela correspond au fait que le \pol $G(X,T)$ est de degré
$1$ en $X$:

\snic {
G(T,X) = (T^4 + 2T^2 + 1)X + \beta T^3 + \alpha T =
(T^2 + 1)^2 X + T(\beta T^2 + \alpha)
,}

%\sni
d'où:
$$
x = -{t(\beta t^2 + \alpha) \over (t^2 + 1)^2}, \quad
y = tx = -{t^2(\beta t^2 + \alpha) \over (t^2 + 1)^2}.
$$
En $t = 0$, on a $(x,y) = (0,0)$. Quelles sont les autres valeurs du
paramètre $t$ pour
lesquelles $\big(x(t), y(t)\big) = (0,0)$? \\
Il faut d'abord trouver les zéros de
$x(t)$, fraction rationnelle de hauteur~4. Il y a la valeur $t = \infty$, pour
laquelle $y(t) = -\beta$.  \\
Si $\alpha = 0$, on a seulement deux zéros de
$x$: $t=0$ (de multiplicité 3) et $t = \infty$ (de multiplicité 1).  
\\
Si
$\beta = 0$, on a seulement deux zéros de $x$: $t=0$ (de multiplicité 1)
et $t = \infty$ (de multiplicité 3).  
\\
Si $\beta \ne 0$, on a deux autres
zéros de $x$ (éventuellement confondus): $t = \pm \sqrt {-\alpha/\beta}$.
On peut rendre cela plus uniforme en faisant intervenir le caractère
quadratique de $-\alpha\beta$, voir la question \emph {7.}

 Remarque: dans tous les cas, en $t = \infty$, on a $(x,y) = (0,-\beta)$.

\emph {5.}
On sait d'après l'exercice \ref {exoAnneauOuvertP1} que $\gk[g_0,g_1]$ est
un \aclz, \cli de $\gk[g_0]$ dans $\gk(t)$. Pour obtenir un
$\gk$-relateur entre $g_0$ et $g_1$, on reporte $t = g_1/g_0$ dans
l'expression $g_0 = 1/(1+t^2)$, ce qui donne $g_0^2 - g_0 + g_1^2 = 0$ et
confirme que $g_1$ est entier sur $\gk[g_0]$.  En $t=0$, on a
$(g_0,g_1) = (1,0)$; ce point est un point non singulier de la courbe
$g_0^2 - g_0 + g_1^2 = 0$. En fait la conique  $C(g_0, g_1) = g_0^2 -
g_0 + g_1^2$ est lisse sur tout anneau puisque

\snic {
1 = -4C + (2g_0-1)\Dpp{C}{g_0} + 2g_1 \Dpp{C}{g_1}.
}

%\sni

Il en est de même de la conique homogénéisée notée encore $C$, $C =
g_0^2 - g_0g_2 + g_1^2$, qui vérifie $\gen {g_0,g_1,g_2}^2 \subseteq \gen {C,
\Dpp{C}{g_0}, \Dpp{C}{g_1}, \Dpp{C}{g_2}}$:
$$
%\snic {
g_0 = -\Dpp{C}{g_2}, \quad g_1^2 = C + (g_0-g_2)\Dpp{C}{g_2}, \quad
g_2 = -\Dpp{C}{g_0} - 2\Dpp{C}{g_2}.
%}
$$
%\sni
On dispose de $\PP^1 \to \PP^2$ défini par $(u : v) \mapsto
(g_0 : g_1 : g_2) = (u^2 : uv : u^2 + v^2)$ dont l'image est la conique
\hmg $C = 0$; à quelque chose près, il s'agit du plongement
de Veronese $\PP^1 \to \PP^2$ de degré $2$.

Par ailleurs, la \dcn en \elts simples fournit les
expressions suivantes de $x$, $y$, $b_3t$, $b_2t$ dans $\gk[g_0,g_1]$:

\snic {
\begin {array} {c}
x = \vep g_0g_1 - \beta g_1, \quad
y = \vep g_0^2 + (2\beta-\alpha)g_0 - \beta =  (g_0-1)(\beta - \vep g_0)
\\[2pt]
b_2t = 2y + (\beta+y)t^2  = -\alpha + \beta g_0 - \vep g_0^2, \quad
b_3t = (2\beta-\alpha)g_1 - \vep g_0g_1
.
\end {array}
}

%\sni
On y voit que $g_0$ est entier sur $\gk[y]$, donc entier sur $\gk[x]$;
comme $g_1$ est entier sur~$\gk[g_0]$, il l'est sur $\gk[x]$. On
vient d'obtenir l'\egt $\gB = \gk[g_0,g_1]$.

Considérons d'abord $\gk[y] \subset \gk[g_0] \subset \gk[g_0,g_1]$; il est
clair que $(1,g_0)$ est une base de $\gk[g_0]$ sur $\gk[y]$ et $(1,g_1)$ est
une base de $\gk[g_0,g_1]$ sur $\gk[g_0]$, donc $(1, g_0, g_1, g_0g_1)$ est
une base de $\gk[g_0,g_1]$ sur $\gk[y]$ (mais pas sur $\gA = \gk[x]$).  
\\
Montrons que $(1, y, b_3 t, b_2t)$ est une $\gA$-base, soit le \Amo $E $ engendré. En utilisant $y-b_2t = \vep(g_0-1)$ et
$x+b_3t = \vep g_1$, on voit que $g_0$, $g_1 \in E$. Enfin, $E$ contient $x +
\beta g_1 = \vep g_0g_1$, donc $g_0g_1 \in E$ et $E = \gk[g_0,g_1] = \gB$.

Un \id\iv $\fb$ de $\gB$ contient un \elt\ndz donc $\gB\sur\fb$ est un \kev de
dimension finie, ce qui permet de définir $\deg \fb$ par $\deg \fb =
\dim_\gk \gB\sur\fb$; on a alors (voir la proposition \ref {prop-a/ab} et son
corolaire \ref{corprop-a/ab}) $\deg(\fb\fb') = \deg(\fb) + \deg(\fb')$.  On
en déduit que $\deg \gen {x,y}_\gB = 4-1 = 3$.

\emph {6.}
On a $\fp_1 = \gen {g_0-1, g_1}$, donc pour montrer l'\egt $\fp_1^2 = \gen
{g_0-1, g_1^2}$, il suffit de voir que $g_0-1 \in \gen {(g_0-1)^2,
g_1^2}$. Cela résulte de l'\egt $1-g_0 = (1-g_0)^2 + g_1^2$ qui découle de
$g_0^2 - g_0 + g_1^2 = 0$.

\emph {7.}
On pose $X = UY$ dans $F(X,Y)$ et on obtient $F(UY,Y) = U^3H(U,Y)$
avec 

\snic {
H(U,Y) = YU^4 + (2Y + \alpha) U^2 + Y + \beta, \qquad
H(U,0) = \alpha U^2 + \beta
.}

%\sni

Ce \pol $H = a'_4U^4 + a'_2U^2 + a'_0$ est primitif en $U$ (on a $a'_2 = 2a'_4
+ \alpha$ \hbox{et $a'_0 = a'_4 + \beta$} donc $\epsilon = a'_0 - a'_2 +
a'_4$). Il vérifie $H(u,y) = 0$ avec $u = x/y$; l'entier~$b'_3$ d'Emmanuel
associé est $x$, et l'on a donc $b'_3u \in \gB$ avec:

\snic {
b'_3u = x^2/y = \vep g_0 - \beta - y
.}

%\sni
En $t$ racine de $\beta t^2 + \alpha = 0$, on a $g_0 = \beta/\vep$ et $g_1^2 =
-\alpha\beta / \vep^2$, ce qui rend naturel l'introduction de l'\id $\fa =
\gen {\vep g_0 - \beta, \vep^2 g_1^2 + \alpha\beta}$.  On vérifie l'\egt:

\snic {
\gen {y, x^2/y}_\gB = 
\gen {\vep g_0 - \beta,  \vep^2 g_1^2 + \alpha\beta}
.}

%\sni
On a alors $\gen {x,y}_\gB = \fp_1\fa$ et $\deg \fa = 2$.  Si $-\alpha\beta$
n'est pas un carré, alors $\fa$ est premier. Sinon, on a $\fa = \fp_2\fp_3$
avec $\fp_2, \fp_3$ s'exprimant avec les deux racines carrées de
$-\alpha\beta$. On a $\fp_2 = \fp_3$ \ssi les deux racines carrées sont
confondues; ceci arrive quand $\alpha\beta = 0$ par exemple ou en \cara $2$.
Enfin, pour $\alpha = 0$, on a $\fp_1 = \fp_2 = \fp_3$.

%%%%%%%%%%%%%%%%%%%%%%%%%%%%%%%%%%%%%%%%%%%%%%%%%%%%%%%%%%%%%%%%%%%%%%%%%%%

}% fin des solutions d'exos 

%:   ---- Section*{references}-----------
\penalty-2500    
\Biblio

\vspace{2pt}

Concernant la genèse de la théorie des \ids de corps de nombres développée 
par Dedekind, on peut lire les articles de H. Edwards \cite{Edw81} et de J. Avigad~\cite{Avi}.

Les \adps intègres ont été introduits par H. Prüfer en 1932 dans \cite{Prufer}.
Leur place centrale en théorie multiplicative des idéaux est mise en valeur
dans le livre de référence sur le sujet \cite{Gil}. 
Voir aussi les commentaires bibliographiques en fin du chapitre~\ref{chap mod plats}.

Dans la littérature classique un \adpc est souvent appelé un 
\ixx{anneau}{semihéréditaire} (selon le point \emph{3} dans le \thref{th.adpcoh}), ce qui n'est pas très joli. 
Ces anneaux sont signalés comme importants dans \cite{CE}.
La preuve du point \emph{1} du \thref{ThImMat} y est donnée, \covz, 
dans le chapitre 1, proposition~6.1.%
\index{semihéréditaire!anneau ---}%
\index{heredit@héréditaire!anneau ---}

Un \ixx{anneau}{héréditaire} est un anneau dans lequel 
 tout \id est \proz. Cette notion est mal définie en \coma à cause 
 de la quantification non légitime \gui{tout idéal}.
Un exemple d'un tel anneau non \noe est le sous anneau d'un produit
 dénombrable de corps $\FF_2$, formé par les suites qui sont ou bien presque partout nulles, ou bien presque partout égales à $1$.  
Le cas le plus intéressant est celui des \adpcs \noesz, que l'on décrit en \clama comme les anneaux dans lesquels tout \id est \ptfz.
Notre \dfn pour un \adk (libéré de la contrainte d'intégrité)
correspond exactement (en \clamaz) à la notion d'anneau héréditaire
\noez.

Des exposés assez complets sur les \anars et les \adps écrits dans le style des \coma se trouvent dans les articles \cite[Ducos\&al.]{dlqs} et~\cite[Lombardi]{lom99}.

Les \gui{entiers d'Emmanuel} du lemme \ref{lemEmmanuel}
sont très présents dans le mémoire de thèse 
d'Emmanuel Hallouin \cite{TheseHallouin}. 

Le \thref{th.2adpcoh} est de Gilmer  et Hoffmann \cite{GHo}.
Le \thref{dekinbe} pour le cas d'un \adp intègre est donné 
par Heitman et Levy dans~\cite{HeLe}.
Le \thref{thcohdim1} a été démontré en \clama par
Quentel dans~\cite{Que}. La \dem \cov est due à I. Yengui.

Le \thref{thPTFDed} est classique (\tho de Steinitz)
pour les \adksz.
Il a été \gne pour les \ddps possédant la 
\prt un et demi dans  \cite[Kaplansky]{Kap52} et
\index{un et demi!\Tho ---}
\cite[Heitmann\&Levy]{HeLe}. L'inspection détaillée de notre \dem montrerait 
d'ailleurs que 
l'hypothèse \gui{anneau \ddi1} pourrait être affaiblie en \gui{anneau possédant la 
\prt un et demi}.

On trouve le \thref{thMpfPruCohDim} 
(voir aussi l'exercice \ref{exoAnarlgb}) dans   \cite[Brewer\&Klinger]{BrKl}
pour le cas intègre. Il a été \gne au cas \qi dans \cite[Couchot]{Couc1}.

Le lemme \ref{lemRadJDIM1} et le \thref{lemthAESTE} sont dus à Claire Tête et Lionel Ducos. 

Le \pb \ref{exoBaseNormaleAlInfini} est basé sur l'article~\cite[Hess]{Hess02}.

Une démonstration alternative du \thref{thAnormalAXnormal} se trouve dans l'article \cite[Coquand \& Lombardi, 2016]{cl2016}.

\newpage \thispagestyle{CMcadreseul}
\incrementeexosetprob

%:        %%%%%%%%%%%%%%%%%%%%%%%%%%%%%%%%%%%%
%:        %%%%%%%%%%%%%%%%%%%%%%%%%%%%%%%%%%%%

%---- Chapitre  {Dimension de Krull}------------
\chapter{Dimension de Krull}
\label{chapKrulldim}
%--------------------

\vspace{-1em}
\minitoc

\subsection*{Introduction}
\addcontentsline{toc}{section}{Introduction}
%-----------------------------------------

Dans ce chapitre on introduit la \ddk dans sa version \cov \elr et on la compare à la notion classique correspondante.

On établit ensuite les premières \prts de cette dimension. La facilité avec laquelle on obtient la \ddk d'un anneau de \pols sur un \cdi montre que la version \cov de la \ddk peut être vue comme une simplification 
conceptuelle de la version classique usuelle. 

Nous appliquons ensuite le même type d'idées pour définir la \ddk
d'un \trdiz, celle d'un morphisme d'anneaux commutatifs, puis la dimension valuative
des anneaux commutatifs.

Nous établissons quelques \thos de base importants concernant ces notions.

Nous terminons en indiquant les versions \covs des notions classiques usuelles de
 Lying over, Going up, Going down et Incomparabilité, avec quelques applications.

%%%%%%%%%%%%%
\section{Espaces spectraux}
\label{secEspSpectraux}

Dans cette section, nous décrivons l'approche de la \ddk
en \clamaz.

Pour nous il s'agit avant tout d'une heuristique.
C'est la raison pour laquelle nous ne donnons aucune \demz.
Cela n'aura aucune incidence dans la suite de l'ouvrage.
En effet, l'aspect \cof des espaces spectraux est entièrement concentré dans
les \trdis obtenus par dualité. En particulier,
l'aspect \cof de la \ddk est entièrement concentré dans
la \ddk des  \trdis et elle peut être définie de manière
complètement indépendante des espaces spectraux.

Néanmoins l'heuristique donnée par les espaces spectraux est essentielle
à la compréhension du petit miracle qui va advenir avec
l'introduction des notions \covs  duales. Petit miracle dont on ne prendra pleinement conscience que dans les chapitres suivants, quand on verra tant de beaux \thos
abstraits se transformer en \algosz.

%: --- Subsec{Treillis de Zariski}-- subsecZarZar
\subsec{Treillis et spectre de Zariski}
\label{subsecZarZar}

Rappelons que l'on note $\DA(\fa)$ le nilradical de
l'\id $\fa$ dans l'anneau $\gA$ et que le treillis de Zariski $\ZarA$ est
l'ensemble des  $\DA(x_1,\ldots ,x_n)$ (pour~$n\in\NN$
et $x_1$, \ldots, $x_n\in\gA$).
On a donc $x\in\DA(x_1,\ldots ,x_n)$ \ssi une puissance de $x$ appartient à $\gen{x_1,\ldots ,x_n}$.
L'ensemble $\ZarA$, ordonné par la relation d'inclusion,
est un \trdi avec 

\snic{\DA(\fa_1)\vu\DA(\fa_2)=\DA(\fa_1+\fa_2)\;
\hbox{ et } \;\DA(\fa_1)\vi\DA(\fa_2)=\DA(\fa_1\,\fa_2).}

%: Definition{nota Spec(A)} ---------
\begin{definition}
\label{nota Spec(A)}\relax
%{\rm
On appelle  \ix{spectre de Zariski} de l'anneau $\gA$  et l'on note~$\SpecA$
l'ensemble des  idéaux premiers stricts
de~$\gA$. On le  munit de la topologie
possé\-dant pour base d'ouverts les
$\fD_\gA(a)=\sotq{\fp\in\SpecA}{a\notin\fp}$.\\
 On note $\fD_\gA(\xn)$
pour $\fD_\gA(x_1)\cup\cdots \cup\fD_\gA(x_n)$.

\noi Pour $\fp\in\SpecA$ et $S=\gA\setminus\fp$ on note $\Ap$ pour $\gA_S$
(l'ambigüité entre les deux notations contradictoires $\Ap$ et $\gA_S$ est
levée en pratique par le contexte).
%}
\end{definition}
%---end-definition-----------------

En \clamaz,  on obtient alors le résultat suivant.

%:     Theorem{thZarZar}
\begin{theoremc}\label{thZarZar} ~
\begin{enumerate}
\item Les \oqcs de
$\SpecA$ sont  les ouverts
$\fD_\gA(x_1,\ldots ,x_n)$.
\item L'application $\DA(x_1,\ldots ,x_n)\mapsto
\fD_\gA(x_1,\ldots ,x_n)$ est bien définie, c'est un \iso de \trdisz.
\end{enumerate}
\end{theoremc}

%: --- Sec{Spectre d'un \trdiz}-- subsecSpecT
\subsec{Spectre d'un \trdiz}
\label{subsecSpecT}

Le spectre de Zariski est l'exemple paradigmatique d'un \ix{espace
spectral}.
Les espaces spectraux ont été introduits par Stone \cite{Sto} en
1937. \index{spectral!espace ---}

Ils peuvent être \cares comme les espaces topologiques vérifiant les
\prts suivantes:
\begin{itemize}
\item l'espace est quasi-compact,
\item tout ouvert  est réunion d'\oqcsz,
\item l'intersection de deux \oqcs
est un \oqcz,
\item pour deux points distincts, il y a un ouvert  contenant l'un mais pas
l'autre,
\item tout fermé irréductible est l'adhérence d'un point.
\end{itemize}
Les \oqcs  forment alors un \trdiz, le sup et le inf étant la réunion et l'intersection. Une application continue entre espaces spectraux est dite
\ixc{spectrale}{application ---} si l'image réciproque de tout \oqc est un \oqcz.
Le résultat fondamental de Stone peut être énoncé comme suit.

\emph{En \clama la catégorie des  espaces spectraux et
applications spectrales  est anti\eqve à  la catégorie des \trdisz. }

\smallskip  Voici comment cela fonctionne.

Tout d'abord si $\gT$ est un \trdiz,
un  \emph{\idepz} est un \id $\fp$ qui vérifie
%--------------------begin array---------------
$$\begin{array}{llllll}
x\vi y \in \fp \;\Rightarrow\; (x\in\fp\;\mathrm{ou}\; y\,\in\fp) ,&   &
1_\gT\notin\fp%\;\Rightarrow\; %1_\gT=0_\gT ; \Faux.
\end{array}$$
%---------------------end array--------------

\rdb \label{SpecTrdi}
Le \emph{spectre} de $\gT$, noté $\SpecT$ est
alors défini comme l'espace dont les points sont les \ideps de
$\gT$ et dont une base d'ouverts est donnée par les
parties $\fD_{\gT}(a):=\sotq{\fp\in\SpecT}{a\notin\fp}$ pour $a\in\gT$.%
\index{spectre!d'un \trdiz}
\\
Si $\varphi:\gT\to\gV$ est un morphisme de \trdisz,
on définit l'application 

\snic{\Spec\varphi:\Spec\gV\to\SpecT,\quad \fp\mapsto\varphi^{-1}(\fp).}

%\sni
C'est une application spectrale et tout
ceci définit $\Spec$ comme foncteur contravariant.
\\
On montre  que les $\fD_{\gT}(a)$ sont tous les \oqcs de $\SpecT$.
En fait le \thoz\eto \ref{thZarZar} s'applique à tout \trdi $\gT$:
{\it \begin{enumerate}
\item Les \oqcs de
$\SpecT$ sont exactement les
$\fD_\gT(u)$.
\item L'application $u\mapsto \fD_\gT(u)$ est bien définie et c'est un \iso de \trdisz.
\end{enumerate}
}\label{NOTASpecT}
Dans l'autre sens, si $X$ est un espace spectral on note
$\OQC(X)$ le \trdi formé par ses \oqcsz.
Si $\xi:X\to Y$ est une application spectrale,
l'application 

\snic{\OQC(\xi):\OQC(Y)\to\OQC(X),\quad U\mapsto\xi^{-1}(U)}

%\sni
est un \homo de \trdisz.
Ceci définit $\OQC$ comme foncteur contravariant.

L'anti\eqvc de catégories qui était annoncée est
définie par les foncteurs~$\Spec$ et $\OQC$.
Elle \gns l'anti\eqvc donnée dans le cas fini au \thrf{thDualiteFinie}.

Notez que l'espace spectral vide correspond au treillis $\Un$, et
qu'un espace spectral réduit à un point correspond au treillis $\Deux$.

%: --- Subsec{Sous espace spectraux}
\subsec{Sous-espaces spectraux}
\label{subsecSesSpec}

Par \dfnz,
un sous-ensemble $Y$ d'un espace spectral $X$ est un \ix{sous-espace spectral}
si la topologie induite fait de $Y$ un espace spectral et si l'injection canonique
$Y\to X$ est spectrale.
\\
Cette notion est en fait exactement la notion duale de la notion de \trdi quotient.
Autrement dit une application spectrale $\alpha:Y\to X$ identifie $Y$
à un sous-espace spectral de $X$ \ssi l'\homo de \trdis
$\OQC(\alpha)$ identifie $\OQC(Y)$
à un \trdi quotient de~$\OQC(X)$.

Les sous-espaces fermés de $X$ sont spectraux et
correspondent aux quotients
par les \idsz. Plus \prmt un \id $\fa$ de $\OQC(X)=\gT$
définit le fermé~$\fV_\gT(\fa)=\sotq{\fp \in X}{\fa \subseteq \fp}$,
(à condition d'identifier les points de~$X$ avec les \ideps de $\OQC(X)$)
et l'on a alors un \iso canonique
$$\OQC(\fV_\gT(\fa)\big)\simeq\OQC(X)\sur{(\fa=0)}.$$
Les fermés irréductibles
correspondent aux \ideps de $\OQC(X)$.

Enfin les \oqcs correspondent aux quotients par des filtres principaux:
$$\OQC(\fD_\gT(u)\big)\simeq\OQC(X)\sur{(\uar u=1)}.$$

%: --- Subsec{Une approche heuristique pour la dimension de Krull}
\subsec{Une approche heuristique pour la \ddkz}
\label{subsecHKdim}

Notons par ailleurs que le spectre de Zariski d'un anneau commutatif s'identifie
de façon naturelle avec le spectre de son treillis de Zariski.

En \clamaz, la notion de \ddk peut être définie, pour un espace spectral
arbitraire~$X$, comme la longueur maximale
des chaînes strictement croissantes de fermés irréductibles.

Une manière intuitive d'appréhender cette notion
de dimension est la suivante. La dimension peut être \caree par
\recu en disant que d'une part, la dimension $-1$ correspond à l'espace vide, et
d'autre part, pour $k\geq 0$, un espace $X$ est de dimension $\leq k$ \ssi pour
tout \oqc $Y$, le bord   de $Y$ dans $X$ est de dimension $\leq k-1$
 (ce bord est fermé donc c'est un sous-espace spectral de $X$).

Voyons par exemple, pour un anneau commutatif $\gA$,
 comment on peut définir le bord de l'ouvert
$\fD_\gA(a)$ dans $\Spec\gA$. Le bord est l'intersection de l'adhérence de
$\fD_\gA(a)$ et du fermé \cop de $\fD_\gA(a)$,
que nous notons $\fV_\gA(a)$.
L'adhérence de  $\fD(a)$ c'est l'intersection de tous
les $\fV(x)$ qui contiennent  $\fD(a)$,
\cad tels que  $\fD(x)\cap\fD(a)=\emptyset$.
\\
Comme  $\fD(x)\cap\fD(a)=\fD(xa)$, et
comme on a $\fD(y)=\emptyset$ \ssi $y$ est nilpotent, on obtient une approche
heuristique de l'\id \gui{bord de Krull
de~$a$}, qui est l'idéal engendré par $a$ d'une part
(ce qui correspond à $\fV(a)$),
et par tous les $x$ tels que $xa$ est nilpotent d'autre part (ce qui
correspond à l'adhérence de $\fD(a)$).

%--- Sec{Definition \cov}--secDefConsDimKrull
\section[Une \dfn \covz]{Définition \cov et premières
consé\-quences}
\label{secDefConsDimKrull}

En \clamaz, la \ddk d'un anneau commutatif est définie
comme le maximum (éventuellement infini) des longueurs des chaînes
strictement croissantes d'\ideps stricts (attention, une chaîne 
$\fp_0\subsetneq \cdots\subsetneq \fp_\ell$ est dite de longueur $\ell$).
Puisque le \cop d'un \idep est un filtre premier, la \ddk est aussi le
maximum  des longueurs des chaînes strictement croissantes de filtres
premiers.

Comme cette \dfn est impossible à manipuler d'un point de vue \algqz, on la
remplace en \coma par une \dfn \eqve
(en \clamaz) mais de nature plus \elrz.

La quantification sur l'ensemble des \ideps de
l'anneau est alors remplacée par une quantification
sur les \elts de l'anneau et les entiers naturels.
Depuis cette découverte (de manière surprenante elle est très récente) 
les \thos qui font intervenir la \ddk ont pu rentrer à part entière dans 
le domaine des \coma et du Calcul Formel.

%:  Definition{defZar2} bords de Krull
\begin{definition}
\label{defZar2} Soient  $\gA$  un anneau commutatif, $x\in\gA$ et $\fa$ un \itfz.
%-----------------begin item------------------
\begin{enumerate}
\item [$(1)$] Le \ix{bord supérieur de Krull} de $\fa$
dans $\gA$ est l'anneau quotient
%---  equation eqBKAC --------
\begin{equation}\label{eqBKAC}
\gA_\rK^{\fa}:=\gA/\JK_\gA(\fa)  \quad \hbox{où} \quad
 \JK_\gA(\fa):=\fa+(\sqrt{0}:\fa).
\end{equation}
%---------------------end equation--------------
On note  $\JK_\gA(x)$  pour $\JK_\gA(x\gA)$ et $\gA_\rK^{x}$
pour $\gA_\rK^{x\gA}$. Cet anneau est appelé le \emph{bord
supérieur de $x$ dans $\gA$}. \\
On dira  que $\JK_\gA(\fa)$ est
\emph{l'\id bord de Krull de $\fa$ dans $\gA$.}%
\index{ideal@idéal!bord de Krull}
\item [$(2)$] Le \ix{bord inférieur de Krull} de $x$ dans $\gA$ est l'anneau
localisé
%---  equation eqBKAS --------
\begin{equation}\label{eqBKAS}
\gA^\rK_{x}:=\SK_\gA(x)^{-1}\!
\gA \quad\hbox{où}\quad \SK_\gA(x)=x^\NN(1+x\gA).
\end{equation}
%---------------------end equation--------------
On dira que  $\SK_\gA(x)$ est le
\emph{monoïde bord de Krull de $x$} dans~$\gA$.%
%:HHH index un peu modifié
\index{monoide@monoïde!bord de Krull}%
\index{bord de Krull!idéal ---}\index{bord de Krull!monoide@\mo ---}
\end{enumerate}
\end{definition}
%--- end-definition-------------------------

Rappelons qu'en \clama la \ddk
d'un anneau est $-1$ \ssi l'anneau n'admet
pas d'\idepz, ce qui signifie qu'il est trivial.

Le \tho suivant donne alors en \clama une \carn inductive \elr
de la \ddk d'un anneau commutatif.

%:  Theorem{thDKA} dimKrull
\begin{theoremc}
\label{thDKA} Pour un anneau commutatif  $\gA$
et un entier $k\geq 0$ \propeq
%-----------------begin enum------------------
\begin{enumerate}
\item  La \ddk de $\gA$ est $\leq k$.
\item  Pour tout $x\in \gA$ la \ddk de $\gA_\rK^{x}$ est
$\leq k-1$.
\item  Pour tout $x\in \gA$ la \ddk de $\gA^\rK_{x}$ est
$\leq k-1$.
\end{enumerate}
%-----------------end enum------------------
\end{theoremc}
%--- end-theorem-----------------------------------------
NB. Ceci est un \tho de \clama qui ne peut pas admettre de \prcoz. \eoe

Dans
la \dem qui suit, tous les \ids et filtres premiers ou maximaux sont
pris au sens usuel en \clamaz: ils sont stricts.
%-----------------begin proof------------------
\begin{proof}
Montrons d'abord l'\eqvc des points \emph{1} et \emph{3.}
Rappelons que les \ideps de $S^{-1}\gA$ sont de la forme $S^{-1}\fp$ où $\fp$
est un \idep
de $\gA$ qui ne coupe pas $S$ (fait \ref{factQuoFIID}).
L'\eqvc résulte alors clairement des deux affirmations suivantes.\\
(a) Soit $x\in\gA$, si $\fm$ est un idéal maximal de $\gA$ il coupe toujours
$\SK_\gA(x)$. En effet, si $x\in\fm$ c'est clair et sinon, $x$ est \iv
modulo $\fm$ ce qui signifie que $1+x\gA$ coupe $\fm$.\\
(b) 
Soient $\fa$ un \idz, $\fp$ un \idep avec $\fp\subset\fa$ et $x\in\fa\setminus\fp$;
si $\fp\cap \SK_\gA(x)$ est non vide, alors $1 \in \fa$.
En effet, soit $x^n(1+xy)\in\fp$; 
puisque $x\notin\fp$, on~a~$1+xy\in\fp\subset\fa$, ce qui donne
avec $x \in \fa$, $1\in\fa$.\\
Ainsi, si $\fp_0\subsetneq \cdots \subsetneq \fp_\ell$ est une chaîne avec
$\fp_\ell$ maximal, elle est raccourcie d'au moins son dernier terme quand on
localise en $\SK_\gA(x)$, et elle n'est raccourcie que de son dernier terme si
$x\in\fp_\ell\setminus\fp_{\ell-1}$.

L'\eqvc des points \emph{1} et \emph{2} se démontre de manière duale,
en remplaçant les \ideps par les filtres premiers.
Soit $\pi:\gA\to\gA/\fa$ la \prn canonique. On remarque que les filtres
premiers de $\gA/\fa$ sont exactement les $\pi(S)$, où $S$ est un filtre premier
de $\gA$ qui ne coupe pas $\fa$  (fait \ref{factQuoIDFI}).
Il suffit alors de démontrer les deux affirmations duales de (a) et (b) qui sont
les suivantes.\\
(a') Soit $x\in\gA$, si $S$ est un filtre maximal de $\gA$ il coupe toujours
$\JK_\gA(x)$. En effet, si $x\in S$ c'est clair et sinon,
puisque $S$ est maximal, $Sx^\NN$ contient~$0$, ce qui signifie qu'il
y a un entier $n$ et un \elt $s$ de $S$ tels que $sx^n=0$.
Alors $(sx)^n=0$ et $s\in (\sqrt{0}:x)\subseteq \JK_\gA(x)$.\\
(b') 
Soient $S'$ un filtre premier contenu dans un filtre $S$ et $x\in S\setminus
S'$. Si~$S'\cap \JK_\gA(x)$ est non vide, alors $S = \gA$. En effet, soit
$ax+b\in S'$ avec~$(bx)^n=0$. Alors, puisque $x\notin S'$, on a $ax\notin S'$
et, vu que $S'$ est premier,~$b\in S'\subseteq S$. Et comme $x\in S$, on obtient
$(bx)^n=0\in S$.
\end{proof}
%-----------------end proof------------------

En \coma on remplace la \dfn usuellement donnée en \clama par la \dfn plus \elr suivante.
%:  Definition{defDiKrull}-----
\begin{definition}
\label{defDiKrull}
La \ixc{dimension de Krull}{d'un anneau commutatif}
(notée $\Kdim$) \label{NOTAKdim}
d'un anneau commutatif $\gA$ est définie
par \recu comme suit:
%-----------------begin enum------------------
\begin{enumerate}
\item  $\Kdim\gA=-1$ \ssi $\gA$ est trivial.
\item  Pour $k\geq 0$,  $\Kdim\gA\leq k$ signifie: $\forall  x\in
\gA,\,\Kdim(\gA^\rK_{x})\leq k-1$.
\end{enumerate}
%-----------------end enum------------------
\end{definition}
%--- end-definition-------------------

Naturellement $\gA$ sera dit de dimension infinie \ssi pour tout entier $k\geq 0$
on a l'implication  $\Kdim\gA\leq k\Rightarrow 1=_\gA0$.

%:--- Lemma{lemZEdDef}---------
\medskip Le lemme suivant résulte immédiatement des \dfnsz.

\begin{lemma}
\label{lemZEdDef}
Un anneau est \zed \ssi il est \ddi$0$.
\end{lemma}
%--- end-lemma-----------------------------------------

Notez que la terminologie \gui{anneau \zedz} constitue donc
un léger abus de langage car affirmer
que la dimension est
inférieure ou égale à~$0$ laisse ouverte la
possibilité de dimension égale à
$-1$, ce qui signifie que l'anneau est trivial.

\medskip \exls \label{exlKdim} ~
\\
1) Si $x$ est nilpotent ou \iv dans $\gA$, l'\id et le \mo bords de
$x$ dans $\gA$ sont tous deux égaux à $\gA$. Les deux anneaux bords
sont triviaux.

 2)
Pour  $x\neq 0$, $1$, $-1$ dans $\ZZ$, les anneaux bords $\ZZ_\rK^x=\ZZ/x\ZZ$ et
$\ZZ_x^\rK=\QQ$  sont \zedsz.
On retrouve donc que   $\Kdim\ZZ\leq1$.

 3)
Soit $\gK$ un corps contenu dans un  \cac discret  $\gL$. Soient
$\fa$ un \itf de $\gK[X_1,\ldots ,X_n]$ et $\gA=\gK[X_1,\ldots ,X_n]/\fa$.
Soient $V$ la \vrt affine correspondant à $\fa$ dans  $\gL^n$, 
$W$ la sous-\vrt de~$V$ définie par $f$. Alors le \gui{bord de~$W$ dans~$V$}, 
défini comme l'intersection de~$W$ avec la clôture de Zariski de son 
\cop dans~$V$,
est la \vrt affine correspondant à l'anneau~$\gA_\rK^f$.
De manière abrégée:

\snic {
\hbox {bord}_V\, \cZ(f) = \cZ_V(\hbox {bord de $f$}).
}

\sni 4) Soit $\gA$ intègre et $k\geq 0$: $\Kdim\gA\leq k$ équivaut à
$\Kdim(\gA\sur{a\gA})\leq k-1$ pour tout $a$ \ndz (utiliser les \ids bords de Krull).

 5) Soit $\gA$ local \dcd et $k\geq 0$: $\Kdim\gA\leq k$ équivaut à
$\Kdim \gA[1/a]\leq k-1$ pour tout $a\in\Rad\gA$  
(utiliser les \mos bords de Krull).
\eoe

\medskip
\comms
1) L'avantage de la \dfn \cov de la \ddk par rapport à la \dfn usuelle est qu'elle est plus simple (pas de quantification sur l'ensemble des
\idepsz) et plus \gnle (pas besoin de supposer l'axiome du choix).
Cependant nous avons seulement défini  la phrase
\gui{$\gA$ est de \ddk $\leq k$}.

2)
Situons nous en \clamaz. La \ddk de $\gA$ peut être définie comme
un \elt de $\so{-1}\cup\NN\cup\so{+\infty}$ en posant

\snic{\Kdim\gA=\inf\sotq{k\in\ZZ,\,k\geq -1}{\Kdim\gA\leq k},}

%\sni
(avec $\inf\,\emptyset_\ZZ=+\infty$). Cette \dfn
basée sur la \dfn \cov \ref{defDiKrull}
est \eqve  à la \dfn donnée usuellement via les chaînes d'\ideps
(cf. le \thoc\ref{thDKA}).

3)
Du point de vue \cof la méthode précédente ne définit pas la \ddk de $\gA$ comme un \elt de
$\so{-1}\cup\NN\cup\so{+\infty}$.
En fait il s'avère que le
concept en question n'est en \gnl pas \ncr
(mais \llec doit nous
croire sur parole).
\perso{remarque  qui semble inévitable pour
répondre à l'angoisse bien naturelle de la lectrice}
\\
 Le point de vue le plus proche des \clama serait de
regarder $\Kdim\gA$ comme une partie de $\NN\cup\so{-1}$,
définie par

\snic{\sotq{k\in\ZZ,\,k\geq -1}{\Kdim\gA\leq k}.}

%\sni
On raisonne alors avec des parties finales (éventuellement vides) de
$\NN\cup\so{-1}$, la relation d'ordre est donnée par l'inclusion renversée, la borne
supérieure par l'intersection et la borne inférieure par la réunion.
\\
Cette approche trouve sa limite avec le \gui{contre-exemple}
du corps des nombres réels (voir le commentaire \paref{remDKRR}).
\eoe

\medskip
On utilise en \coma les \emph{notations} suivantes, pour se
rapprocher du langage classique:

%--- Notation{notaKdiminf}-------
\begin{notation}
\label{notaKdiminf}
{\rm  Soient $\gA$, $\gB$, $(\gA_i)_{i\in I}$,
$(\gB_j)_{j\in J}$ des anneaux commutatifs (avec $I$, $J$ finis).
%-----------------begin item------------------
\begin{enumerate}
\item [--]  $\Kdim\gB\leq \Kdim\gA$ signifie: $\forall \ell\geq -1\;
(\Kdim\gA\leq \ell\;\Rightarrow \Kdim\gB\leq \ell)$.
\item [--]  $\Kdim\gB= \Kdim\gA$ signifie:  $\Kdim\gB\leq  \Kdim\gA$ et
$\Kdim\gB\geq  \Kdim\gA$.
\item  [--] $\sup_{j\in J}\Kdim\gB_j\leq  \sup_{i\in I}\Kdim\gA_i$ signifie:
$$\preskip.1em \postskip.3em
\forall\ell\geq -1\quad  \big(\,\&_{i\in I} \, \Kdim\gA_i\leq
\ell \;\Rightarrow\;\&_{j\in J}\, \Kdim\gB_j\leq
\ell  \,\big).
$$
\item  [--] $\sup_{j\in J}\Kdim\gB_j=  \sup_{i\in I}\Kdim\gA_i$ signifie:
$$\preskip.1em \postskip.0em
\forall\ell\geq -1\quad \big(\, \&_{i\in I} \,\Kdim\gA_i\leq
\ell\;\Leftrightarrow\;\&_{j\in J}\,\Kdim\gB_j\leq
\ell  \,\big).
$$
\end{enumerate}
%-----------------end item------------------
}
\end{notation}
%--- end-notation-----------------------------------------

%%%%%%%%%%%%%%%%%%%%%%%%%%%%%%%%%%%%%%%%%%%%%%%%%%%%%%%%%%%%%%%%%%%%%%%%%%
%: subsec{Bords itérés, suites singulières}
\subsec{Bords itérés, \susisz, suites \copsz}

La \dfn \ref{defDiKrull} peut être réécrite en terme d'\idasz.
Pour cela, nous introduisons la notion   de \emph{\susiz}.

%:     definition{notaBordsIteres}
\begin{definition}\label{notaBordsIteres}
  Pour une suite  $(\ux)=(\xzk)$ dans $\gA$ on définit les \emph{bords de Krull itérés} de la manière suivante.

 1.  Une version \gui{itérée} du \mo
$\SK_\gA(x)$: l'ensemble
\begin{equation}\label{eqMonBordKrullItere}
\SK_\gA(\xzk):=
x_0^\NN(x_1^\NN\cdots(x_k^\NN (1+x_k\gA) +\cdots)+x_1\gA) + x_0\gA)
\end{equation}
est un \moz.
Pour une suite vide, on définit $\SK_\gA()=\so{1}$.%
\index{monoide@monoïde!bord de Krull itéré}%
\index{bord de Krull!monoide@monoïde --- itéré}%

 2. On définit deux variantes pour l'\id bord de Krull itéré.
\\
 --- 2a) L'\id $\JK_\gA(\xzk)=\JK_\gA(\ux)$ est défini comme suit:
\begin{equation}\label{eq1IdBordKrullItere}
\JK_\gA()=\so{0}, \;
\JK_\gA(x_0, \ldots, x_k) =
\big(\DA\big(\JK_\gA(x_0, \ldots, x_{k-1})\big) : x_k\big) + \gA x_k.%
\index{ideal@idéal!bord de Krull itéré}%
\index{bord de Krull!ideal@idéal --- itéré}%
\end{equation}
 --- 2b)  L'\id $\IK_\gA(\xzk)=\IK_\gA(\ux)$ est défini comme suit:
\begin{equation}\label{eqIdBordKrullItere}
\IK_\gA(\ux):=\sotq{y\in\gA}{\,0\in
x_0^\NN\big(\cdots\big(x_k^\NN (y+x_k\gA) +\cdots\big) + x_0\gA\big)}
\end{equation}
Pour une suite vide, on définit $\IK_\gA()=\so{0}$.

\end{definition}

On montrera (lemme \ref{fact2BordKrullItere}) que les deux \ids \gui{bord itéré} définis ci-dessus ont même nilradical.

%:  propdef{defSeqSing}-----
\begin{propdef}
\label{defSeqSing}~
\begin{enumerate}
\item Une suite $(x_0,\ldots ,x_k)$ dans  $\gA$ est dite
\ixd{singulière}{suite}
si $0\in\SK_\gA(\xzk)$, \emph{i.e.} si $1\in\IK_\gA(\xzk)$,
\cade s'il existe $a_0$, \ldots, $a_k\in \gA$ et
$m_0$, \ldots, $m_k\in \NN$ tels que
%-----equation--eqsing-------------
\begin{equation}\label{eqsing}
x_0^{m_0}(x_1^{m_1}(\cdots(x_k^{m_k} (1+a_k x_k) +
\cdots)+a_1x_1) + a_0x_0) =  0
\end{equation}
%---------------------end equation--------------
%
\item La \prt \gui{la suite $(x_0,\ldots ,x_k)$ est singulière dans  $\gA$}
est une \prt \carfz.
\item (Principe \lgb pour les suites singulières) Soient $S_1$, $\ldots$, $S_n$ des \moco de $\gA$, la suite $(x_0,\ldots ,x_k)$ est singulière dans~$\gA$ \ssi elle est singulière dns chacun des $\gA_{S_i}$
\end{enumerate}
\end{propdef}
%--- end-propdef------------------------------------
%
\begin{proof}
\emph{2}. Soit $S$ un \mo de $\gA$. Par la technique usuelle de calcul dans un anneau $\gA_S$, dire que la suite $(x_0,\ldots ,x_k)$ est singulière dans  $\gA_S$ revient à dire qu'il existe $s\in S$, $a_0$, \ldots, $a_k\in \gA$ et
$m_0$, \ldots, $m_k\in \NN$ tels que 
$$x_0^{m_0}(x_1^{m_1}(\cdots(x_k^{m_k} (s+a_k x_k) +
\cdots)+a_1x_1) + a_0x_0) =  0,$$
ce qui implique que, dans $\gA[1/s]$
$$\ndsp
x_0^{m_0}(x_1^{m_1}(\cdots(x_k^{m_k} (1+\fraC{a_k}s x_k) +
\cdots)+\fraC{a_1}s x_1) + \fraC{a_0}sx_0) =  0.$$

\emph{3}. Une fois les exposants $m_i$ fixés, on peut regarder l'\egt \pref{eqsing} comme une \eqn linéaire en les  $a_j$.
On note par ailleurs que si l'\eqn est satisfaite pour un \sys d'exposants, elle est aussi satisfaite pour un \sys d'exposants supérieurs.
Donc en prenant un \sys d'exposants qui majore chacun de ceux obtenus
séparément pour chaque  $\gA_{S_i}$, on obtient une unique équation \lin en les $a_j$ qui a une solution dans chaque  $\gA_{S_i}$.
On peut donc appliquer le \plg de base  \ref{plcc.basic}.\iplg 
\end{proof}

\rem En \clama on déduit des points \emph{2} et \emph{3} qu'une suite est singulière \ssi elle est singulière après \lon en tout \idemaz.
\eoe

%On peut aussi prendre tous les exposants $m_i$ égaux à un même $m$.

%:--- proposition{corKrull}------
\begin{proposition}
\label{corKrull}
Pour un anneau commutatif $\gA$  et un entier $k\geq 0$, \propeq
%-----------------begin enum------------------
\begin{enumerate}
\item La \ddk de $\gA$ est $\leq k.$
\item Pour tout $x\in \gA$ la \ddk de $\gA_\rK^{x}$ est
$\leq k-1$.
\item \label{i3corKrull} Toute suite $(\xzk)$ dans $\gA$  est \singz.
\item \label{i4corKrull} Pour tous $x_0$, \ldots, $x_k \in \gA$ il existe
$b_0$, \ldots, $b_k\in \gA$ tels que
%---  equation eqCG --------
\begin{equation}\label{eqCG}
\left.\arraycolsep2pt
\begin{array}{rcl}
\DA(b_0x_0)& =  &\DA(0)    \\
\DA(b_1x_1)& \leq  & \DA(b_0,x_0)  \\
\vdots~~~~& \vdots  &~~~~  \vdots \\
\DA(b_k x_k )& \leq  & \DA(b_{k -1},x_{k -1})  \\
\DA(1)& =  &  \DA(b_k,x_k )
\end{array}
\right\}
\end{equation}
%---------------------end equation--------------
%
\item  \label{i5corKrull} Pour tous $x_0$, \ldots, $x_k \in \gA$, en posant
$\pi_i = \prod_{j < i} x_j$ pour $i \in \lrb {0..k+1}$
(donc $\pi_0=1$), il existe $n\in\NN$ tel que

\snic{\pi_{k+1}^n \in \gen {
\pi_{k}^n x_k^{n+1},\
\pi_{k-1}^n x_{k-1}^{n+1},\ \ldots,\
\pi_1^n x_1^{n+1},\
\pi_0^n x_0^{n+1}}.\qquad\qquad\quad}%
\end{enumerate}
%-----------------end enum------------------
\end{proposition}
%--- end-proposition-----------------
Par exemple pour $k=2$ le point \emph{\ref{i4corKrull}} correspond au dessin suivant dans $\ZarA$.
$$
\SCO{\DA(x_0)}{\DA(x_1)}{\DA(x_2)}{\DA(b_0)}{\DA(b_1)}{\DA(b_2)}
$$

%-----------------begin proof------------------
\begin{proof}
Les \eqvcs pour la dimension $0$ sont \imdes par application des
\dfnsz.
\\
\emph{1} $\iff$ \emph{3.} Supposons l'\eqvc 
 établie pour la dimension $\leq k$ et pour tout anneau commutatif. On voit
alors que
$S^{-1}\gA$ est de dimension  $\leq k$ \ssi l'on  a:
\\
\emph{pour tous
$x_0$, \ldots, $x_k\in \gA$ il existe $a_0$, \ldots, $a_k\in \gA$, $s\in S$
et $m_0$, \ldots, $m_k\in \NN$ tels que}
%----- equation--eqsing3----
\begin{equation}\label{eqsing3}
x_0^{m_0}(x_1^{m_1}\cdots(x_k^{m_k} (s+a_k x_k) +
\cdots + a_1x_1) + a_0x_0) =  0.
\end{equation}
%---------------------end equation--------------
Notez que par rapport à l'\eqrf{eqsing}, un $s\in S$ a remplacé le $1$ au
centre de l'expression du premier membre.
\\
Il reste donc à remplacer $s$ par un \elt arbitraire de
$\SK_\gA(x_{k+1})$, \cad un \elt de la
forme $x_{k+1}^{m_{k+1}} (1+a_{k+1} x_{k+1})$.
\\
L'\eqvc entre \emph{2} et \emph{3} se prouve de manière analogue.
\\
\emph{3} $\Rightarrow$ \emph{4.} On prend
$b_k=1+a_k x_k$, puis
$b_{\ell -1} = x_\ell^{m_\ell} b_\ell+ a_{\ell -1}x_{\ell -1}$, successivement pour $\ell=k,$
\ldots, $1$.
\\
\emph{4} $\Rightarrow$ \emph{2.} Preuve par \recuz.
L'implication pour la dimension $\leq 0$ est claire. Supposons la chose établie pour
la dimension $< k $. Supposons la \prtz~\emph{4} et montrons que pour tout
$x_0$
la dimension de $\gB= \gA\ul{x_0}$  \hbox{est $<k$}.
\\
Par \hdrz, il suffit de trouver, pour tous  $x_1$, \ldots, $x_k$
des \hbox{\elts $b_1$, \ldots, $b_k$}  tels que
$$\left.\arraycolsep2pt
\begin{array}{rcl}
\rD_\gB(b_1x_1)& =  & \rD_\gB(0)  \\
\vdots~~~~& \vdots  &~~~~  \vdots \\
\rD_\gB(b_k x_k )& \leq  & \rD_\gB(b_{k -1},x_{k -1})  \\
\rD_\gB(1)& =  &  \rD_\gB(b_k,x_k ).
\end{array}
\right\}
$$
Or, par hypothèse, on a des \elts $b_0$, \ldots, $b_k$  tels que
$$
\left.\arraycolsep2pt
\begin{array}{rcl}
\DA(b_0x_0)& =  &\DA(0)    \\
\DA(b_1x_1)& \leq  & \DA(b_0,x_0)  \\
\vdots~~~~& \vdots  &~~~~  \vdots \\
\DA(b_k x_k )& \leq  & \DA(b_{k -1},x_{k -1})  \\
\DA(1)& =  &  \DA(b_k,x_k ).
\end{array}
\right\}
$$
et les in\egts avec $\DA$ impliquent les mêmes avec $\rD_\gB$. La
deuxième in\egt signifie que $(b_1x_1)^m \in \gen {b_0, x_0}$ (pour un
certain $m$); la première nous dit que $b_0x_0$ est nilpotent donc $\gen
{b_0,x_0} \subseteq \JK_\gA(x_0)$. Bilan: $b_1x_1$ est nilpotent
dans $\gB$.
\\
On pourrait aussi démontrer \emph{4} $\Rightarrow$ \emph{3}  par un calcul direct un peu fastidieux.
\\
\emph{3} $\Leftrightarrow$ \emph{5}. Dans la \dfn d'une \susiz,
on peut remplacer tous les exposants $m_i$ par leur maximum $n$. Une fois
ceci acquis,  le point \emph{5}  est une simple reformulation du point \emph{3.}
\end{proof}
%-----------------end proof------------------

Nous aurions donc pu donner une \dfn par \recu de la \ddk basée sur les bords
supérieurs $\gA_\rK^{x}$ plutôt que sur les bords inférieurs~$\gA^\rK_{x}$:
nous venons d'obtenir une \prco directe (sans utiliser le \thoc\vref{thDKA}) de
l'\eqvc entre les deux \dfns par \recu possibles.

\medskip \rem
Le système d'in\egts (\ref{eqCG}) dans le point \emph{\ref{i4corKrull}} de
la proposition~\ref{corKrull} établit une relation intéressante et
\smq entre les deux suites~$(b_0,\ldots ,b_k)$ et $(x_0,\ldots ,x_k)$. \\
Lorsque $k=0$,
cela signifie $\DA(b_0)\vi\DA(x_0)=0$ et $\DA(b_0)\vu\DA(x_0)=1$,
\cad
que les deux \elts  $\DA(b_0)$ et $\DA(x_0)$
sont complé\-ments l'un de l'autre
dans le treillis $\ZarA$.
Dans $\Spec\gA$ cela signifie que les ouverts de base
correspondants sont \copsz. \\
Nous introduisons donc
la terminologie suivante: lorsque les suites~$(b_0,\ldots ,b_k)$ 
et~$(x_0,\ldots ,x_k)$ vérifient les
in\egts (\ref{eqCG}) nous dirons que ce sont deux \emph{suites \copsz}.%
\index{suites complémentaires!dans un anneau commutatif}%
\index{complémentaires!suites --- (anneaux commutatifs)}%
\eoe

%%%%%%%%%%%%%%%%%%%%%%%%%%%%%%%%%%%%%%%%%%%%%%%%%%%%%%%%%%%%%%%%%%%%
%:    fact0BordKrullItere
\begin {fact}\label{fact0BordKrullItere}
Soient $(\ux) = (\xn)$, $(\uy) = (\ym)$ deux suites d'\elts de $\gA$,
$\gA \to \gA'$ un morphisme et $(\uxp)$  l'image de $(\ux)$ dans $\gA'$.
\begin {enumerate}
\item
On a les \eqvcsz:
$$\preskip.4em \postskip.4em
\exists z \in\IK_\gA(\ux) \cap \SK_{\gA}(\uy)   \iff
1 \in \IK_\gA(\ux, \uy) \iff 0 \in \SK_\gA(\ux,\uy).
$$

\item
Si $\gA \to \gA'$ est surjectif, l'image de $\SK_\gA(\ux)$ est
$\SK_{\gA'}(\uxp)$.

\item
Si $\gA' = S^{-1}\gA$, avec $S$  un \mo de $\gA$, alors $S^{-1}
\IK_{\gA}(\ux) = \IK_{\gA'}(\uxp)$.

\end {enumerate}
\end {fact}

%:    fact1BordKrullItere
\begin {fact}\label{fact1BordKrullItere}
Soient $\fa$ un \id de $\gA$, $Z \subseteq \gA$ une partie quelconque et
$x \in \gA$. 

\snic{x^\NN(Z + \gA x) \hbox { rencontre } \fa
\;  \iff \;
Z  \hbox { rencontre } (\fa : x^\infty) + \gA x.}
\end {fact}

%%%%%%%%%%%%%%%%%%%%%%%%%%%%%%%%%%%%%%%%%%%%%%%%%%%%%%%%%%%%%%%%%%%%

%:    lem1BordKrullItere
\begin{lemma}\label{lem1BordKrullItere} \emph{(Idéaux bords de Krull à la Richman)}\\
Pour une suite $(\ux) = (\xn)$ d'\elts de $\gA$, l'\id bord itéré
$\IK_\gA(\ux)$ peut être défini de manière récursive comme suit:
$$
\IK_\gA()=\so{0}, \qquad
\IK_\gA(x_1, \ldots, x_n) =
\big(\IK_\gA(x_1, \ldots, x_{n-1}) : x_n^\infty\big) + \gA x_n.
$$
Par exemple:
$$
\IK_\gA(x_1) = (0 : x_1^\infty) + \gA x_1, \;
\IK_\gA(x_1, x_2) =
\big(\big((0 : x_1^\infty) + \gA x_1\big) : x_2^\infty\big) + \gA x_2.
$$
\end {lemma}

\begin {proof}
On définit provisoirement
$$
N() = \{0\}, \quad
N(x_1, \ldots, x_n) = \big(N(x_1, \ldots, x_{n-1}) : x_n^\infty\big) + \gA x_n
$$
Prenons $n = 3$ pour fixer les idées. Alors, pour $y \in \gA$, on a
les \eqvcsz:
$$
\begin {array} {ll}
0 \in x_1^\NN\big(x_2^\NN\big(x_3^\NN(y + \gA x_3) + \gA x_2\big) + \gA x_1\big)
&\quad  \Longleftrightarrow \\[1mm]
x_2^\NN\big(x_3^\NN(y + \gA x_3) + \gA x_2\big) \hbox { rencontre } N(x_1)
&\quad  \Longleftrightarrow \\[1mm]
x_3^\NN(y + \gA x_3) \hbox { rencontre }
\big(N(x_1) : x_2^\infty\big) + \gA x_2 \eqdefi  N(x_1, x_2)
&\quad  \Longleftrightarrow \\[1mm]
y \in \big(N(x_1,x_2) : x_3^\infty\big) + \gA x_3 \eqdefi
  N(x_1, x_2, x_3),
\end {array}
$$
ce qui prouve que $\IK_\gA(x_1, x_2, x_3) = N(x_1, x_2, x_3)$.
\end {proof}

%%%%%%%%%%%%%%%%%%%%%%%%%%%%%%%%%%%%%%%%%%%%%%%%%%%%%%%%%%%%%%%%%%%%
%:    lem2BordKrullItere
\begin {lemma}\label{lem2BordKrullItere} \emph{(Enchaînement de bords itérés, coté idéal)}\\
Soient $(\ux) = (\xn)$ et $(\uy) = (\ym)$ deux suites d'\elts de~$\gA$. 
On pose $\gA' = \gA/\IK_\gA(\ux)$ et l'on note
$(\uyp) = (\ypm)$ l'image de $(\uy)$ dans~$\gA'$.
\begin {enumerate}
\item
Le noyau du morphisme canonique (surjectif) $\gA \to \gA'\big/\IK_{\gA'}(\uyp)$
est l'\id $\IK_\gA(\ux, \uy)$.
\item
On définit  $\gA_0 =
\gA$ et $\gA_i = \gA_{i-1}/\IK_{\gA_{i-1}}(x_i)$ pour $ i\in \lrbn$. Alors le
noyau du morphisme canonique (surjectif) $\gA \to \gA_n$ est l'\id
$\IK_\gA(\ux)$.
\end {enumerate}
\end {lemma}

\begin {proof}
Il suffit de montrer le premier point pour $n = 1$. Notons $x = x_1$.
\\
Soit $z \in \gA$  et $z'$ son image dans $\gA' = \gA/\IK_{\gA}(x)$.
On a les \eqvcsz:

\snic{\begin {array} {ll}
z=0 \quad\mathrm{dans}\quad \gA'/\IK_{\gA'}(\uyp)
&\quad  \Longleftrightarrow \\[1mm]
0 \in {y'_1}^\NN\big(\cdots\big({y'_m}^\NN (z' + y'_m\gA') +\cdots\big) + y'_1\gA'\big)
&\quad  \Longleftrightarrow \\[1mm]
{y_1}^\NN\big(\cdots\big({y_m}^\NN (z + y_m\gA) +\cdots\big) + y_1\gA\big)
\hbox { rencontre } \IK_\gA(x)
&\quad  \Longleftrightarrow \\[1mm]
0 \in x^\NN({y_1}^\NN\big(\cdots\big({y_m}^\NN (z + y_m\gA) +\cdots\big) + y_1\gA) + x\gA\big)
&\quad  \Longleftrightarrow \\[1mm]
z \in \IK_\gA(x, \uy).
\end {array}
}
\vspace{-.5em}
\end {proof}

%:    fact2BordKrullItere
\begin{fact}\label{fact2BordKrullItere}
Pour toute suite $(\ux)$ d'\elts de $\gA$, les idéaux $\IK_\gA(\ux)$ 
et~$\JK_\gA(\ux)$ ont même nilradical.
\end {fact}

\begin{proof}
Pour tout \id $\fa$ et tout $x \in \gA$, on vérifie facilement que la
racine de l'\id $(\fa : x^\infty)$ est $(\DA(\fa) : x)$.  En utilisant
$\DA(\fb + \fc) = \DA(\DA(\fb) + \fc)$, on en déduit que les idéaux $(\fa
: x^\infty) + \gA x$ et $(\DA(\fa) : x) + \gA x$ ont même racine.
Le résultat annoncé s'en déduit par \recu sur la longueur
de la suite $(\ux)$ en utilisant la \dfn récursive des deux
idéaux bord itérés.
\end{proof}

%%%%%%%%%%%%%%%%%%%%%%%%%%%%%%%%%%%%%%%%%%%%%%%%%%%%%%%%%%%%%%%%%%%%
%:    lem3BordKrullItere
\begin {lemma}\label{lem3BordKrullItere}
%:HHH j'ai introduit V=\SK_{\gA'}(x') pour alleger
Soient $S$ un \mo de $\gA$, $\gA' = S^{-1} \gA$, $x \in \gA$, $x'$
son image dans $\gA'$ et $V=\SK_{\gA'}(x')$.  Alors le morphisme canonique $\gA \to
V^{-1}\gA'$ est un morphisme de \lonz\footnote{Voir éventuellement la \dfn \ref{defHomloc}.} et induit un
\iso \hbox{de $T^{-1}\gA$} \hbox{sur $ V^{-1}\gA'$}, où $T$ est le \mo
$x^\NN(S + \gA x)$.
\end {lemma}

\begin {proof}
L'image dans $V^{-1}\gA'$ de l'\elt $s + ax \in S + \gA x$
est \iv puisque l'on peut écrire $s+ ax = s(1 + ax/s)$ (avec quelques abus de notations). D'où un morphisme (canonique) $\varphi :
T^{-1}\gA \to V^{-1}\gA'$. \\ 
Par ailleurs, puisque $S \subseteq T$,
on a un morphisme $\gA' \to T^{-1}\gA$. L'image par ce morphisme de $1 + xa/s
\in 1 + x\gA'$ est \iv car $1 + xa/s = (s + xa)/s$, d'où un
morphisme $\varphi' : V^{-1}\gA' \to T^{-1}\gA$.\\
On vérifie sans peine que $\varphi$ et $\varphi'$ sont inverses l'un de l'autre.
\end {proof}

%%%%%%%%%%%%%%%%%%%%%%%%%%%%%%%%%%%%%%%%%%%%%%%%%%%%%%%%%%%%%%%%%%%%
%:    cor1BordKrullItere
%:HHH rajout d'un titre
\begin{corollary}\label{cor1BordKrullItere} \emph{(Enchaînement de bords itérés, coté \moz)}\\
Soient $(\ux) = (\xn)$ et $(\uy) = (\ym)$ %deux suites 
dans
$\gA$, 
$\gA' = \SK_\gA(\uy)^{-1}\gA$, \hbox{et 
$(\uxp) = (\xpn)$} l'image de $(\ux)$ dans~$\gA'$.
\\
Alors, le morphisme  
$\gA \to \SK_{\gA'}(\uxp)^{-1}\gA'$ 
donne par \fcn  un \iso
$\SK_{\gA}(\ux,\uy)^{-1}\gA \simarrow \SK_{\gA'}(\uxp)^{-1}\gA'$.
\end{corollary}

%%%%%%%%%%%%%%%%%%%%%%%%%%%%%%%%%%%%%%%%%%%%%%%%%%%%%%%%%%%%%%%%%%%%%%%%%%
%: subsec{Suites regulières}
\subsec{Une suite régulière \gui{n'est pas}\ singulière}

Le point \emph{4} de la proposition qui suit implique qu'une suite régulière qui n'engendre pas
l'\id $\gen{1}$ est non singulière, ce qui explique le titre du paragraphe.

Un avantage des bords de Krull itérés à la Richman est que sur un anneau \noe \coh ce sont des \itfs (voir le lemme \ref{lem1BordKrullItere}). Un autre avantage est donné par le point \emph{1} dans la proposition qui suit.

%:     Proposition{lemRegsing}--------
\begin{proposition}
\label{lemRegsing} \emph{(Suites régulières et \ddkz)}\\
Soit $(\xn)$ une suite régulière dans $\gA$ et $(\yr)$
une autre suite. 
\begin{enumerate}
\item \label{i1lemRegsing} On a $\IK_\gA(\xn)=\gen{\xn}$.
\item \label{i2lemRegsing} La suite $(\xn,\yr)$ est singulière dans $\gA$ \ssi la suite
$(\yr)$  est singulière dans $\aqo\gA\xn$.
\item \label{i3lemRegsing}  L'implication suivante est satisfaite pour tout $k\geq -1$:

\snic{\Kdim \gA \leq n+k\;\Longrightarrow\;\Kdim \aqo\gA\xn \leq k.}

%\sni
Et si $1\notin\gen{\xn}$, on a ${n+\Kdim \aqo\gA\xn \leq \Kdim\gA}$. 
\item \label{i4lemRegsing} 
Si la suite $(\xn)$ est \egmt singulière,
on a $1\in\gen{\xn}$. \\
En conséquence si $\Kdim\gA\leq n-1$ toute suite régulière de longueur~$n$
engendre l'\id $\gen{1}$.
\end{enumerate}
\end{proposition}
%--- end-lemma-----------------------------------------
%
\begin{proof}
\emph{1.} Calcul \imd en tenant compte de la \dfn récursive donnée dans le lemme~\ref{lem2BordKrullItere} (point~\emph{2}).\\
\emph{2.} On applique le point \emph{1} du lemme~\ref{lem2BordKrullItere}.\\
\emph{3.} Résulte du point \emph{2.}\\
\emph{4.} Cas particulier du point \emph{2}, avec la suite $(\yr)$ vide.
\end{proof}
%

%%%%%%%%%%%%%%%%%%%%%%%%%%%%%%%%%%%%%%%%%%%%%%%%%%%%%%%%%%%%%%%%%%%%%%%%%%
%: subsec{Minorer la \ddk}
\subsec{Minorer la \ddk}

Il peut être confortable, voire parfois utile,
de définir la phrase \gui{$\gA$ est de \ddk $\geq k$}.

Tout d'abord $\Kdim\gA\geq 0$
doit signifier $1\neq0$. Et une possibilité, \linebreak 
pour $k\geq1$, sera de demander: \gui{il existe une suite
$(x_1,\ldots ,x_k)$ qui n'est pas singulière}.
Notez que du point de vue \cof cette affirmation est plus forte
 que la négation de
\gui{toute suite  $(x_1,\ldots ,x_k)$ est singulière}.

Un anneau a alors une \ddk bien
définie s'il existe un entier~$k$ tel que l'anneau soit à la fois de \ddk
$\geq k$ et de \ddk $\leq k$.

La chose ennuyeuse est le \crc trop compliqué de l'assertion

\centerline{\gui{la suite $(x_1,\ldots ,x_k)$ n'est pas singulière}.}

Il semble de toute manière ici impossible  d'éviter l'usage de la négation, car on ne voit pas comment on pourrait définir $\Kdim\gA\geq 0$
autrement que par la négation $1\neq0$. Même dans le cas où $\gA$ est un anneau \fdiz, la
phrase \gui{il existe une suite
$(\ux)=(x_1,\ldots ,x_k)$ telle que $1\notin\IK_\gA(\ux)$}
reste problématique d'un point de vue \cof car on doit tester tous les \elts de l'\id $\IK_\gA(\ux)$, qui n'est pas a priori \tfz. Dans le cas où l'anneau est en plus \noez, le lemme \ref{lem1BordKrullItere} montre que l'\id est \tf et l'on peut tester si la suite $(\ux)$ est singulière.

La \dfn de $\Kdim\gA\leq k$ qui correspond à une assertion de type $\forall\exists$ est plus sympathique que celle de  $\Kdim\gA\geq k$, qui correspond, lorsque l'anneau est discret,
à une assertion de type $\exists\forall$.
Cependant même la \dfn de $\Kdim\gA\leq k$ ne peut \gnlt pas être certifiée par un simple
calcul: il faut une preuve.

Notons que pour l'anneau $\RR$, si l'on  utilise la négation forte
(de \crc positif), pour laquelle $x\neq0$ signifie \gui{$x$ est \ivz},
pour définir la phrase $\Kdim\RR\geq k$,
alors il est absurde que $\Kdim\RR\geq1$. Mais on ne peut prouver \cot
que $\Kdim\RR\leq0$ (commentaire \paref{remDKRR}).

%%%%%%%%%%%%%%%%%%%%%%%%%%%%%%%%%%%%%%%%%%%%%%%%%%%%%%%%%%%%%%%%%%%%%%%%%%%
%%%%%%%%%%%%%%%%%%%%%%%%%%%%%%%%%%%%%%%%%%%%%%%%%%%%%%%%%%%%%%%%%%%%%%%%%%%
%%%%%%%%%%%%%%%%%%%%%%%%%%%%%%%%%%%%%%%%%%%%%%%%%%%%%%%%%%%%%%%%%%%%%%%%%%%
%%%%%%%%%%%%%%%%%%%%%%%%%%%%%%%%%%%%%%%%%%%%%%%%%%%%%%%%%%%%%%%%%%%%%%%%%%%
%--- Sec{Proprietes elem}---------
\section{Quelques \prts \elrs de la \ddkz}
\label{secKrullElem}
%-----------------------------------------

Les faits énoncés dans la proposition suivante sont faciles
(notez que l'on utilise les notations \ref{notaKdiminf}).

%: --- Proposition{propDdk0}-------
\begin{proposition}
\label{propDdk0} Soit $\gA$ un anneau, $\fa$ un \idz, $S$ un \moz.
\begin{enumerate}
\item Une \susi reste singulière
dans  $\gA\sur{\fa}$ et $\gA_S$.
\item $\Kdim \gA\sur{\fa} \leq \Kdim \gA$, $\Kdim \gA_S \leq \Kdim \gA$.
\item $\Kdim (\gA\times \gB) = \sup( \Kdim \gA, \Kdim\gB)$.
\item $\Kdim\gA=\Kdim\Ared$.
\item Si $a$ est \ndz dans $\Ared$ (a fortiori s'il est \ndz dans $\gA$), \linebreak 
alors $\Kdim \gA\sur{a\gA}\leq \sup(\Kdim\gA,0)-1$.
\item Si $a\in\Rad\gA$, alors $\Kdim \gA[1/a] \leq \sup(\Kdim\gA,0) -1$.
\end{enumerate}
 \end{proposition}
%--- end-proposition----------------------------------------

%:HHH dÂplacÂ, un peu modifiÂ et rajout d'un commentaire gÂomÂtrique
\exl On 
donne un anneau $\gB$ pour lequel $\Frac(\gB)$ est de \ddk $n>0$,
mais $\gB\red=\Frac(\gB\red)$ est \zedz.
\\
Considérons $\gB = \gA\sur{x\fm}$, où $\gA$ est local \dcdz,
$\fm=\Rad\gA$ et~$x \in \fm$.  L'anneau $\gB$ est local, $\Rad\gB=\fm' = \fm/x\fm$ et $\gB/\fm'=\gA/\fm$.\\ 
Si $\overline x = 0$, alors $x \in x\fm$,
i.e. $x(1-m) = 0$ avec $m \in \fm$, ce qui implique $x = 0$.
\\
Pour $y\in \fm$ on a $\ov y\,\ov x=0$. 
Donc si $\ov y \in\Reg \gB\cap\fm'$, on obtient $x=0$. 
\\
Or on a $\ov y\in\fm'$ ou $\ov y\in\gB\eti$, 
donc si $x\neq 0$ et~$\ov y \in\Reg \gB$,  on obtient $\ov y\in\gB\eti$. Autrement dit, si $x\neq 0$, $\gB=\Frac(\gB)$.
\\
Prenons  $\gA = \gk[\xzn]_{\gen {\xzn}}$ où $\gk$ est un \cdi non trivial, 
\hbox{et $x = x_0$}.
On a alors $\aqo\gA{x_0} \simeq \gk[x_1, \ldots, x_n]_{\gen {x_1, \ldots, x_n}}$ \hbox{et $\Kdim\aqo\gA{x_0} = n$}.
Comme $\ov {x_0}^{2}=0$ dans $\gB$, on a $\gB\red \simeq \aqo\gA{x_0}$ et donc $\Kdim\gB=n$.\\
Enfin $\Frac(\gB\red)=\gk(\xn)$ est un \cdiz, \zedz.
\\
Géométriquement: on a considéré l'anneau d'une variété \gui{avec multiplicités} consistant en un point immergé dans un hyperplan de dimension $n$, et on a localisé en ce point immergé. 
\\
NB. En \clamaz, si $\gC$ est \noe et réduit, $\Frac(\gC)$ est un produit
fini de corps, donc \zedz. Pour une version \cov on peut se reporter
au \pb \ref{exoAnneauNoetherienReduit} et à \cite[Coquand\&al.]{cls}.
\eoe

%: --- Plcc{thDdkLoc}-----------
\begin{plcc}\label{thDdkLoc}
\emph{(Pour la \ddkz)}\\
Soient $S_1$, $\ldots$, $S_n$ des \moco d'un anneau $\gA$ et $k\in\NN$.
\begin{enumerate}
\item Une suite est singulière dans $\gA$ \ssi elle est singulière dans
chacun des $\gA_{S_i}$.
\item L'anneau $\gA$ est \ddi$k$ \ssi les $\gA_{S_i}$ sont \ddi$k$.
\end{enumerate}
\end{plcc}
%--- end-plcc-----------------------------------------
On aurait pu écrire: $\Kdim\gA=\sup_{i}\Kdim\gA_{S_i}=\Kdim\prod_{i}\gA_{S_i}$.
%-----------------begin proof------------------
\begin{proof}
Le point \emph{2} résulte du point \emph{1}, qui a été démontré en \ref{defSeqSing}.
%-% PERSO
\perso{La preuve télégraphique suivante semble fautive.
Une suite $\xzk$ est singulière \ssi
le \mo $\SK(\xzk)$
 contient $0$. On conclut donc par le
 \plgrf{plcc.basic.monoides} qui concerne les \mosz.
}
%-% Fin PERSO
\end{proof}
%-----------------end proof------------------

\rem Comme la \prt pour une suite d'être singulière
est \carfz, on déduit en \clama du résultat précédent le \plga correspondant: \emph{un anneau est \ddi$k$ \ssi il est \ddi$k$ après \lon en tout \idemaz}. 
\eoe

De même on note que le  point \emph{1} dans le \plgc précédent s'applique en fait
toujours avec une famille d'\ecoz, ce qui correspond à un recouvrement
fini du spectre de Zariski par des ouverts de base.

Dans le cas d'un recouvrement fini par des fermés, le résultat
tient encore.

%:     prcf{thDdkRecFer}
\begin{prcf}\label{thDdkRecFer} \emph{(Dimension de Krull)}\\
Soit $\gA$ un anneau, $k$ un entier $\geq 0$, et  $\fa_1$, \ldots, $\fa_r$ des \ids de $\gA$.\\
On suppose tout d'abord que les $\fa_i$ forment un recouvrement fermé de 
$\gA$.
\begin{enumerate}
\item Une suite $(\xzk)$ est singulière dans $\gA$ \ssi elle est
singulière dans chacun des $\gA\sur{\fa_i}$.
\item L'anneau $\gA$ est \ddi$k$ \ssi chacun des $\gA\sur{\fa_i}$ est \ddi$k$.
\end{enumerate}
Plus \gnltz, sans hypothèse sur les $\fa_i$ on a
\begin{enumerate}\setcounter{enumi}{2}
\item L'anneau $\gA\big/\!\bigcap_i\fa_i$ est \ddi$k$ \ssi chacun des $\gA\sur{\fa_i}$ est \ddi$k$.\\
Ceci peut s'écrire sous forme abrégée 

\snic{\Kdim\gA\big/\!\prod_i\fa_i=\Kdim\gA\big/\!\bigcap_i\fa_i=\sup_{i}\Kdim\gA\sur{\fa_i}=\Kdim\prod_{i}\gA\sur{\fa_i}.}
\end{enumerate}

\end{prcf}
\begin{proof}
Il suffit de montrer le point \emph{1.} La suite $(\xzk)$ est singulière \ssi
le \mo $\SK(\xzk)$
 contient $0$. \\
 En outre, $\SK_{\gA\sur{\fa_i}\!}(\xzk)$
 n'est autre que $\SK_\gA(\xzk)$ vu modulo $\fa_i$. On conclut par le \prf \vref{prcf1}.
\end{proof}

%: --- Theorem{th1.5}-------------
\goodbreak
\begin{theorem} \emph{(\Tho un et demi)}%
%:HHH index
\index{un et demi!\Tho ---}
\label{th1.5}\relax
\begin{enumerate}
\item
\begin{enumerate}
\item Si $\gA$ est \zedz, ou plus \gnlt si $\gA$ est \lgbz, tout
module \lmo est monogène.
\item  Si $\gA$ est \zedz, tout \id \ptf est engendré par un \idmz.
\end{enumerate}
\item Soit $\gA$  \ddi$k$, soit $(x_1,\ldots ,x_k)$ une suite
régulière
et $\fb$ un \id \lop contenant $\fa=\gen{x_1,\ldots ,x_k}$ alors il existe
$y\in\fb$ tel que

\snic{\fb=\gen{y,x_1,\ldots ,x_k}=\gen{y}+\,\fb\fa=\gen{y}+\,\fa^m}

%\sni
pour n'importe quel exposant $m\geq 1$.
\item Soit $\gA$ tel que  $\gA\sur{\Rad\gA}$ est \ddi$k$, \linebreak 
soit $(x_1,\ldots
,x_k)$ une suite régulière dans  $\gA\sur{\Rad\gA}$ et  $\fb$ un \id \ptf de
$\gA$ contenant $\fa=\gen{x_1,\ldots ,x_k}$ alors il existe $y\in\fb$ tel que

\snic{\fb=\gen{y,x_1,\ldots ,x_k}=\gen{y}+\,\fb\fa=\gen{y}+\,\fa^m}

%\sni
pour n'importe quel exposant $m\geq 1$.
\end{enumerate}
\end{theorem}
%--- end-theorem-----------------------------------------
%-----------------begin proof------------------
\begin{proof}
Le point \emph{1a}  est un rappel (voir le point \emph{\iref{ite4propZerdimLib}}
du \thrf{propZerdimLib}
pour les anneaux \zeds  et le point \emph{2} du \thrf{thlgb3}
pour les \algbsz).
 
 \emph{1b.} Rappelons que dans un anneau arbitraire un \id \ptfz~$\fa$ a pour annulateur un \idmz~$h$. Dans $\aqo\gA h$, $\fa$ est fidèle, donc $\fa^k$ aussi,
pour tout $k\geq1$. Dans $\gA[1/h]$, $\fa=0$. Donc $\Ann(\fa^k)=\Ann(\fa)=\gen{h}$
\linebreak 
pour $k\geq1$.
\\
Dans le cas \zedz, puisqu'un \id \ptf  est \lopz, il est principal d'après le point \emph{1a}, notons le $\gen{x}$. On
sait que pour $k$ assez grand, $\gen{x}^k=\gen{e}$ avec $e$ \idmz.
D'après la remarque préliminaire $\Ann (x)=\Ann (e)= \gen{1-e} $.
Dans $\aqo\gA {1-e}$, $x$ est \ivz, donc $\gen{x}=\gen{1}$;  
dans $\aqo\gA {e}$, $x$ est nul; ainsi dans $\gA$, $\gen{x}=\gen{e}$.
 
\emph{3.} Résulte de \emph{2} par le lemme de Nakayama.
 
\emph{2.}
L'\id $\fb$ vu comme \Amoz, après \eds à $\gA\sur{\fa}$,
devient le module $\fb\sur{\fb\fa}$ et il reste \lmoz. Puisque l'anneau
quotient $\gA\sur{\fa}$
est \zedz,
 le point \emph{1a} nous dit que  $\fb\sur{\fb\fa}$ est engendré par un \elt $y$. Cela
signifie $\fb=\gen{y}+\,\fb\fa$ et les autres \egts suivent
\imdtz.
\end{proof}
%-----------------end proof------------------

\rem
Dans le cas de la dimension 1
et d'un \id \ivz, le point \emph{2} du \tho précédent est souvent
appelé
\gui{\tho un et demi}. Voir le corolaire \ref{corpropZerdimLib}
et le \thref{th1-5}.
\eoe

%%%%%%%%%%%%%%%%%%%%%%%%%%%%%%%%%%%%%%%%%%%%%%%%%%%%%%%%%%%%%%%%%%%%%%%%%%%
\section{Extensions entières}\label{secDDKExtEn}
% Rappelons %(cf. section \ref{subsecEntiers})
%que si  $\gA\subseteq\gB$, un \elt $x$ de $\gB$ est dit entier sur $\gA$ s'il annule
%un \pol unitaire à \coes dans $\gA$. Il revient au même de dire que $\gA[x]$
%est un \Amo \tfz. Si tout \elt de $\gB$ est entier sur $\gA$ on dit que $\gB$ est
%\emph{entier sur $\gA$} ou encore que $\gB$ est une \ix{extension entière} de
%$\gA$. La proposition suivante se formule en \clama en disant que la dimension de
%$\gA$ est inférieure ou égale à celle de $\gB$. Nous donnons une formulation
%\cov \eqvez.

%:     Proposition{propDKEXENT}---
\begin{proposition}
\label{propDKEXENT}
Soit des anneaux $\gA\subseteq\gB$ avec $\gB$ entier sur $\gA$.
Toute suite finie d'\elts de $\gA$ qui est singulière dans $\gB$ est 
singulière dans $\gA$.
En particulier, $\Kdim\gA\leq \Kdim\gB$. 
\end{proposition}
%--- end-proposition----------------------------------------
NB: l'in\egt opposée est prouvée un peu
plus loin (\thref{cor2thKdimMor}).
%-----------------begin proof------------------
\begin{proof} 
Supposons par exemple que la suite  $(x,y)\in\gA$ soit singulière dans $\gB$, i.e. 

\snic{\exists a,b\in\gB, \;\exists m,\ell\in\NN,\;  x^\ell\big(y^m(1+ay)+bx\big)=0.}

%\sni
On veut réaliser le même type d'\egtz,  avec
des \elts $a'$, $b'$ de $\gA$ au lieu des \elts $a$, $b$ dans $\gB$. L'idée intuitive
est de transformer l'\egt précédente par l'opération \gui{norme}.
Considérons des \polusz~$f$,~$g\in\AT$ qui annulent $a$ et $b$.
Soit $\gB_1=\aqo{\gA[T,T']}{f(T),g(T')}$.
Notons $\alpha$ et $\beta$ les classes de $T$ et $T'$ dans $\gB_1.$
Le sous-anneau
$\gA[a,b]$ de~$\gB$ est un quotient de $\gB_1=\gA[\alpha,\beta]$, via un 
$\gA$-\homo qui envoie~$\alpha$ et~$\beta$ sur~$a$ et~$b$.
En outre, $\gB_1$ est un module libre de rang fini sur $\gA$
ce qui permet de définir la norme et l'\elt cotransposé d'un 
\elt de~$\gB_1[X,Y]$ arbitraire.
Soit alors

\snic{U  (\alpha,\beta,X,Y)=X^\ell\big(Y^m(1+\alpha Y)+\beta X\big)\;$ et $\;V(X,Y)=\rN_{\gB_1[X,Y]/\gA[X,Y]}(U).}

%\sni
D'après le lemme \ref{lemNormeSuiteSing}, $V(X,Y)$ est de la forme

\snic{X^{p}\big(Y^q\big(1+A(Y) Y\big)+B(X,Y) X\big),}

%\sni
 avec $A\in\gA[Y]$, $B\in\gA[X,Y]$.
Soit par ailleurs $W(\alpha,\beta,X,Y)\in\gB_1[X,Y]$ l'\elt
cotransposé de $U(\alpha,\beta,X,Y)$. 
En spécialisant $X,Y,\alpha,\beta$ en $x,y$ dans $\gA$ et $a,b$
dans $\gB$, on obtient une \egt dans $\gB$

\snic{V(x,y)=x^{p}\big(y^q\big(1+A(y) y\big)+B(x,y) x\big)=U(a,b,x,y)W(a,b,x,y),}

%\sni
ce qui termine la \dem puisque $V(x,y)=0$ est une \egt dans~$\gA$:
notez que l'on a $U(a,b,x,y)=0$ dans $\gB$ mais que~$U(\alpha,\beta,x,y)$
n'est peut-être pas nul dans $\gB_1$.
\end{proof}

%:     Lemma{lemNormeSuiteSing}
\begin{lemma}\label{lemNormeSuiteSing}
Soit $\gC$ une \Alg  libre de rang fini sur $\gA$, $(\czn)$ dans $\gC$ et
$(\Xzn)=(\uX)$ une liste d'\idtrsz. Posons
$$
U(\uX) =X_0^{k_0}\big(X_1^{k_1}\big(\cdots(X_n^{k_n} (1+c_nX_n) +\cdots)+c_1X_1\big) + c_0X_0\big)\in\gC[\uX].
$$
Alors $V(\uX)\eqdefi\rN_{\gC[\uX]/\gA[\uX]}(U(\uX)\big)$ est de la forme
$$
V(\uX)=
X_0^{\ell_0}\big(X_1^{\ell_1}\big(\cdots(X_n^{\ell_n}
(1+a_nX_n) +\cdots)+a_1X_1\big) + a_0X_0\big)\in\gA[\uX],
$$
avec $a_n\in\gA[X_n]$,  $a_{n-1}\in\gA[X_n,X_{n-1}]$, \ldots,
$a_0\in\gA[\uX]$.
\end{lemma}
\begin{proof}
Tout d'abord la norme $\rN(1+c_nX_n)$ est un \pol  $h(X_n)\in\gA[X_n]$ qui vérifie
$h(0)=1$, donc qui s'écrit sous la forme $1+a_n(X_n)X_n$. Ensuite
on utilise la multiplicativité de la norme, et une \evn en $X_{n-1}=0$
pour montrer que $\rN\big(X_n^{k_n}(1+c_nX_n)+c_{n-1}X_{n-1}\big)$ est de la forme

\snic{X_n^{\ell_n}\big(1+a_n(X_n)X_n\big)+a_{n-1}(X_n,X_{n-1})X_{n-1}.}

%\sni
Et ainsi de suite. \Llec sceptique ou \llec
pointil\leux peut faire une preuve par \recu en bonne et due forme.
\end{proof}

%  --- Sec{anneaux geometriques}----
\section{Dimension des anneaux géométriques}
\label{secDimGeom}
%-----------------------------------------

%: subsec Anneaux de \pols sur un \cdiz
\subsec{Anneaux de \pols sur un \cdiz}

Un premier résultat important dans la théorie de
la \ddk est la dimension des anneaux de \pols
sur un \cdiz.
%: --- th{propDimKXY}---
\begin{theorem}
\label{propDimKXY}
Si $\gK$ est un \cdi non trivial, la \ddk de
l'anneau de \pols
$\gK[X_1,\ldots,X_\ell]$ est égale à $\ell$.
\end{theorem}
%--- end-proposition-----------------

Nous établissons d'abord le résultat suivant qui nécessite une \dfn précise.
Des \elts $x_1$, \ldots, $x_\ell$ d'une \Klg avec $\gK$ \zed sont dits 
\emph{\agqt dépendants sur $\gK$} s'ils annulent un \pol primitif\footnote{Ceci \gns pour $\ell$ \elts la notion d'\elt primitivement \agq
 introduite \paref{subsecExtAdpC}. Pour un anneau $\gK$ arbitraire, il serait plus raisonnable d'utiliser une terminologie plus contraignante du style \gui{\rde primitivement \agqz}. Il est clair aussi que le \plgref{plcc.agq} se généralise dans le cas de~$\ell$ \eltsz.}
 $f\in\gK[X_1,\ldots, X_\ell]$.

%--- Proposition{propAdepSing}--
\begin{proposition}
\label{propAdepSing}
Soit $\gK$ un \cdiz, ou plus \gnlt un anneau \zedz,
$\gA$ une \Klgz, et $x_1$, \ldots, $x_\ell
\in \gA$ \agqt dépendants sur $\gK$.
Alors la suite $(x_1,\dots,x_\ell)$ est singulière.
\end{proposition}
%--- end-proposition----------------------------------------
%-----------------begin proof------------------
\begin{proof}
On traite le cas d'un \cdiz. Le cas \gnl s'en déduit par application
de la machinerie \lgbe \elr \num2.\imlg
\\
Soit $Q(x_1,\ldots,x_\ell)= 0$ une \rde \agq sur
$\gK$.
Mettons un ordre lexicographique sur les \moms non nuls $\alpha_{p_1,\ldots
,p_\ell}\allowbreak x_1^{p_1}x_2^{p_2}\cdots x_\ell^{p_\ell}$ de $Q$, en accord
avec les  \gui{mots} $p_1\,p_2\,\ldots \,p_\ell$.
Nous pouvons supposer le \coe du plus petit \mom non nul égal à   $1$ (ici
on utilise l'hypothèse que le corps est discret, car on suppose que l'on peut
déterminer pour chaque $\alpha_{p_1,\ldots ,p_\ell}$ s'il est nul ou
\ivz).
Soit $x_1^{m_1}x_2^{m_2}\cdots x_\ell^{m_\ell}$ ce \momz.
En suivant l'ordre lexicographique, nous voyons que nous pouvons écrire
$Q$ sous la forme
%--------------------begin array---------------
$$\arraycolsep2pt\begin{array}{rcl}
Q& =   &  x_1^{m_1}\cdots x_\ell^{m_\ell}+
x_1^{m_1}\cdots x_\ell^{1+m_\ell}R_\ell  +x_1^{m_1}\cdots x_{\ell-1}^{1+m_{\ell-
1}}R_{\ell-1} \\[1mm]
&   & +\cdots+ x_1^{m_1}x_2^{1+m_2}R_2+ x_1^{1+m_1}R_1
\end{array}$$
%---------------------end array--------------
où $R_j\in \gK[x_k\,;\,k\geq j]$. Alors $Q=0$ est l'\egt voulue.
\end{proof}
%-----------------end proof------------------

\begin{Proof}{Preuve du \tho \ref{propDimKXY}. }
Nous notons d'abord que la suite $(X_1,\ldots,X_\ell)$ est régulière, ce qui
montre que la \ddk de
$\gK[X_1,\ldots,X_\ell]$ est  $\geq \ell$. On peut voir aussi directement qu'elle
est non singulière: dans  l'\egref{eqsing} avec $x_i=X_i$ le membre de gauche est
non nul (considérer le \coe de
$X_1^{m_1}X_2^{m_2}\cdots \allowbreak X_\ell^{m_\ell}$).

 Pour prouver que la dimension de
$\gK[X_1,\ldots,X_\ell]$  est $\leq \ell$,
il suffit, vu la proposition \ref{propAdepSing}, de montrer que $\ell+1$
\elts de $\gK[X_1,\ldots,X_\ell]$ sont toujours \agqt dépendants sur
$\gK$. Voici une preuve \elr de ce résultat classique. \label{prvalgdep}
Soient $y_1,\ldots ,y_{\ell+1}$ ces \eltsz,
et $d$ une borne sur leurs degrés. Pour un entier $k\geq 0$ considérons la
liste ${L}_k$ de tous \hbox{les
 $y_1^{\delta_{1}}\cdots
y_{\ell+1}^{\delta_{\ell+1}}$}   tels que $\sum_{i=1}^{\ell+1}
\delta_i\leq k$. Le nombre d'\elts de la liste  ${L}_k$  est égal à
${k +\ell+1}\choose{k }$: ceci est un \pol de degré $\ell+1$ en $k$. Les \elts
de  ${L}_k$ vivent dans l'\evc $E_{\ell,kd}$ des \elts de $\gK[X_1,\ldots,X_\ell]$
de degré $\leq k \,d$, qui a pour dimension
${d\,k +\ell}\choose{d\,k }$: ceci est un \pol de degré~$\ell$ en~$k$. Ainsi pour
$k$ assez grand,
le cardinal de  ${L}_k$ est plus grand que la dimension de l'\evc  $E_{\ell,kd}$
où se trouvent les \elts
de $L_k$, donc il y a une \rdl entre les \elts de  ${L}_k$. Ceci fournit une \rde \agq entre les~$y_i$.
\end{Proof}
%-----------------end proof------------------

%--- Commentaire
\comm   \label{remDKRR}
La preuve de la proposition \ref{propAdepSing} ne peut pas fournir  \cot le
même résultat pour le corps des réels $\RR$ (qui \emph{n'est pas}
discret). En fait il est impossible de réaliser pour $\RR$
le test de
zéro-dimensionnalité: 

\snic{\forall x\in \RR \;\;\exists a\in \RR\;\;\exists n\in \NN,\,
x^n\, (1-ax) = 0.}

%\sni
Cela signifierait en effet que pour tout nombre réel $x$, on sache
trouver un réel $a$ tel que $x(1-ax)=0$.
Si  l'on a trouvé un tel $a$, on obtient:
%-----------------begin item------------------
\begin{enumerate}
\item [--]  si $ax$ est \iv alors $x$ est \ivz,
\item  [--] si  $1-ax$ est \iv alors $x=0$.
\end{enumerate}
%-----------------end item------------------
Or l'alternative  \gui{$ax$ ou $1-ax$  \ivz} est explicite sur $\RR$.
Ainsi fournir le test de zéro-dimensionnalité revient à fournir le  test
pour \gui{$x$ est nul ou \ivz~?}.
Mais ceci n'est pas possible du point de vue \cofz.\\
Par ailleurs, on peut montrer qu'il est impossible
d'avoir une suite non singulière
de longueur $1$, si l'on  prend $y\neq 0$ au sens fort de
\gui{$y$ est \ivz} (dans la \dfn de \gui{non singulière}).
En effet, si l'on  a un $x$ tel que

\snic{\forall a\in \RR \;\;\forall n\in \NN,\;\; x^n\,(1-ax)\;
\mathrm{est}\;\mathrm{\iv},
}

%\sni
cela donne une contradiction: avec $a=0$ on obtient $x$  \ivz,
donc il existe un  $b$ tel que $1-bx=0$. \eoe

%: subsec{Un corolaire intéressant}
\subsec{Un corolaire intéressant}

%:     Lemma{lemahbonvraiment}
\begin{lemma}\label{lemahbonvraiment}
Un anneau engendré par $k$ \elts est de \ddk finie.
\end{lemma}
\begin{proof}
Puisque la dimension ne peut que diminuer par passage à un quotient,
il suffit de montrer que $\ZZXk$ est de \ddkz~$\leq2k+1$
(en fait cet anneau est de \ddk $k+1$ d'après le \thrf{corthValDim}). \\
Soit $(h_1,\ldots,h_{2k+2})$ une suite de $2k+2$ \elts dans
$\ZZXk=\ZZuX$. Nous devons montrer qu'elle est singulière.\\
La suite  $(h_1,\ldots,h_{k+1})$ est singulière dans
$\QQXk=\QQuX$. Cela signifie que l'\id bord itéré
$\IK_\QQuX(h_1,\ldots,h_{k+1})$ contient $1$. 
\\
En {chassant les \denosz}
on obtient que  $\IK_\ZZuX(h_1,\ldots,h_{k+1})$ contient un
entier~$d>0$.
Donc l'anneau $\gB=\ZZuX\big/{\IK_\ZZuX(h_1,\ldots,h_{k+1})}$ est un quotient de
l'anneau $\gC=(\aqo{\ZZ}{d})[\uX]$. Comme $\aqo{\ZZ}{d}$ est \zedz,
la \linebreak 
suite $(h_{k+2},\ldots,h_{2k+2})$ est  singulière dans $\gC$ (proposition \ref{propAdepSing}), autrement dit l'\id $\IK_\gC(h_{k+2},\ldots,h_{2k+2})$
contient $1$. A fortiori $\IK_\gB(h_{k+2},\ldots,h_{2k+2})$
contient $1$. Finalement l'anneau

\snic{\ZZuX\sur{\IK_\ZZuX(h_1,\ldots,h_{2k+2})}= \gB\sur{\IK_\gB(h_{k+2},\ldots,h_{2k+2})}}

%\sni
est trivial.
\end{proof}
%

%: subsec{Anneaux géométriques}
\subsec{Anneaux géométriques}\rdb

Nous appelons \emph{anneau \gmqz} un anneau $\gA$ qui est
une \Klg \pf avec $\gK$  un \cdi non trivial.\index{anneau!geome@géométrique}

Le \thref{thNstNoe} de mise en position
de \iNoe affirme qu'un tel anneau quotient est
 une extension entière finie
d'un anneau $\gB=\KYr$ contenu dans $\gA$ (ici, $Y_1$, \dots, $Y_r$ sont des
\elts de $\gA$ \agqt indépendants sur $\gK$).

%--- Theorem{thDKAG}------------
\begin{theorem}
\label{thDKAG}
Sous les hypothèses précédentes, la \ddk de l'anneau $\gA$ est égale à $r$.
\end{theorem}
%--- end-theorem-----------------------------------------
%-----------------begin proof------------------
\begin{proof}
La \ddk est $\leq r$ par application de la proposition~\ref{propAdepSing}. On peut
d'ailleurs donner une preuve du fait que $r+1$ \elts de $\gA$ sont
\agqt dépendants sur $\gK$ dans le même style que celle donnée
\paref{prvalgdep} pour une \alg de \polsz.\\
Enfin la \ddk  
est $\geq r$ d'après la proposition~\ref{propDKEXENT}.\\
Notons que le \thref{cor2thKdimMor} nous donne une autre \demz, via l'\egt $\Kdim\gA=\Kdim\gB$.
\end{proof}
%-----------------end proof------------------

%%%%%%%%%%%%%%%%%%%%%%%%%%%%%%%%%%%%%%%%%%%%%%%%%%%%%%%%%%%%%%%%%%%%%%%%%%%
%%%%%%%%%%%%%%%%%%%%                                 %%%%%%%%%%%%%%%%%%%%%%
%%%%%%%%%%%%%%%%%%%%    Dimension des treillis dist  %%%%%%%%%%%%%%%%%%%%%%
%%%%%%%%%%%%%%%%%%%%                                 %%%%%%%%%%%%%%%%%%%%%%
%%%%%%%%%%%%%%%%%%%%%%%%%%%%%%%%%%%%%%%%%%%%%%%%%%%%%%%%%%%%%%%%%%%%%%%%%%%
\section{Dimension de Krull des \trdisz}
\label{secDDKTRDIS}

Comme nous l'avons déjà signalé, la \ddk d'un anneau
commutatif $\gA$ n'est autre que la \ddk de
l'espace spectral $\SpecA$, du moins
en \clamaz.

En \coma on introduit la \ddk d'un \trdi $\gT$
de façon à ce qu'elle soit égale,
en \clamaz, à la  \ddk de son spectre $\SpecT$.
La \dem de cette \egt est à très peu près identique
à celle que nous avons donnée pour les anneaux commutatifs.
Nous ne la répéterons pas, puisque de toute manière, nous utiliserons
toujours la \ddk d'un \trdi via la \dfn \cov qui suit.

%:     Definition{defiDDKTRDI}
\begin{definition}\label{defiDDKTRDI}~
\index{dimension de Krull!d'un treillis distributif}%
\index{suites complémentaires!dans un treillis distributif}
\begin{enumerate}
\item Deux suites $(\xzn)$ et $(\bzn)$ dans
un \trdi $\gT$ sont dites \ixc{complémentaires}{suites --- (treillis distributifs)} si
%---  equation eqC2G --------
\begin{equation}\label{eqC2G}
\left.\arraycolsep3pt
\begin{array}{rcl}
 b_0\vi x_0& =  & 0    \\
 b_1\vi x_1& \leq  &  b_0\vu x_0  \\
\vdots~~~~& \vdots  &~~~~  \vdots \\
 b_n\vi  x_n & \leq  &  b_{n -1}\vu x_{n -1}  \\
 1& =  &   b_n\vu x_n
\end{array}
\right\}
\end{equation}
%---------------------end equation--------------
Une suite qui possède une suite \cop sera dite \ixd{singulière}{suite}.
\item Pour $n\geq0$ on dira que le \trdi $\gT$ est \emph{de \ddk inférieure ou égale à $n$} si toute suite $(\xzn)$ dans $\gT$ est singulière.
Par ailleurs,  on dira que le \trdi $\gT$ est de \ddk $-1$
s'il est trivial, \cad si $1_\gT=0_\gT$.
\end{enumerate}
\end{definition}
\hmodeHabillage{\hbox{\xyrowsp=3pt
$\SCO{x_0}{x_1}{x_2}{b_0}{b_1}{b_2}$}}{0}{-15pt}Par exemple pour $k=2$ le point {1} correspond au dessin suivant dans $\gT$.
%$$
%\SCO{x_0}{x_1}{x_2}{b_0}{b_1}{b_2}
%$$
\\On notera $\Kdim\gT\leq n$ lorsque la \ddk est $\leq n$.
\label{NOTAKdimTrdi}\\
Il est évident qu'un treillis a la même \ddk que le treillis opposé.
On voit aussi tout de suite qu'un treillis est \zed \ssi c'est une \agBz.
\'Egale\-ment: un ensemble totalement ordonné de $n$ \elts  a pour
\ddkz~$n-2$.
\endHabillage
%:     Fact{factDDKTRDI}
\begin{fact}\label{factDDKTRDI}
Soit $S$ une partie de $\gT$ qui engendre $\gT$ en tant que \trdiz. Alors
$\gT$ est de \ddk $\leq n$ \ssi toute suite $(\xzn)$ dans $S$
admet une suite \cop dans $\gT$.
\end{fact}
\begin{proof}
Illustrons les calculs sur un exemple suffisamment \gnl dans le cas~$n=4$. 
\\
On vérifie que si $(x_0,x_1,x_2,x_3,x_4)$ admet $(a_0,a_1,a_2,a_3,a_4)$
pour suite \copz, et si  $(x_0,x_1,y_2,x_3,x_4)$ admet $(b_0,b_1,b_2,b_3,b_4)$
pour suite \copz, alors  la suite $(x_0,x_1,x_2\vu y_2,x_3,x_4)$ admet la suite \copz~$(a_0\vu b_0,a_1\vu b_1,a_2\vi b_2,a_3\vi b_3,a_4\vi b_4)$. 
\\
Et dualement,
la suite $(x_0,x_1,x_2\vi y_2,x_3,x_4)$ admet la suite \copz~$(a_0\vu b_0,a_1\vu b_1,a_2\vu b_2,a_3\vi b_3,a_4\vi b_4)$.
\\ 
Le même calcul fonctionnerait pour un $x_i$ arbitraire (au lieu de $x_2$ ci-dessus) dans une suite finie arbitraire. Ainsi si chaque suite $(z_0,\ldots,z_n)$ dans $S$ admet une suite
\cop dans $\gT$, il en sera de même pour toute suite de $n+1$ termes dans le treillis engendré
par $S$.
\end{proof}
%

%:     fact{corfactDDKTRDI}
\begin{fact}\label{corfactDDKTRDI}
Un anneau commutatif a même \ddk que son treillis de Zariski.
\end{fact}
\begin{proof}
La \demz, basée sur le fait \ref{factDDKTRDI}, est laissée \alecz. 
Une autre \dem sera donnée plus loin
sous la forme du lemme~\ref{lemsutfidele}.
\end{proof}

On peut aussi accéder à la \ddk via les \ids bords de Krull comme pour les anneaux commutatifs.

%:     Definition{defiBordKrullTrdi}
\goodbreak
\begin{definition}\label{defiBordKrullTrdi}~
\begin{enumerate}
\item Le treillis $\gT\ul x=\gT/(\JK_\gT(x)=0)$, où
%---  equation eqbordsupTrdi --------
\begin{equation}\label{eqbordsupTrdi}
\JK_\gT(x)\,=\,\dar x \,\vu\, (0:x)_\gT
\end{equation}
%---------------------end equation--------------
est appelé \emph{le bord supérieur (de 
Krull) de $x$ dans $\gT$}. On dit aussi que l'\id $\JK_\gT(x)$ est \emph{l'\id bord de 
Krull de $x$ dans $\gT$.} \index{ideal@idéal!bord de Krull (\trdiz)}
\index{bord supérieur de Krull!(\trdiz)}
\item Plus \gnltz, pour une suite  $(\ux)$ dans $\gT$, l'\id bord de Krull itéré $\JK_\gT(\ux)$ est défini par \recu comme suit: $\JK_\gT()=\so{0},$ et
\begin{equation}\label{eq1IdBordKrullItereTrdi}
\JK_\gT(x_0, \ldots, x_k) =
\big(\JK_\gT(x_0, \ldots, x_{k-1}) : x_k\big)_\gT \vu \dar x_k.%
\index{ideal@idéal!bord de Krull itéré (\trdiz)}
\end{equation}
\end{enumerate}
\end{definition}
%--------- fin definition ---------------------------------------------- 

%:     Fact{factDdkTrdiBord}
\begin{fact}\label{factDdkTrdiBord}
Soient   $n\in\NN$  et $\gT$ un  \trdi. 
\begin{enumerate}
\item Une suite $(\xzn)$ dans $\gT$ est singulière \ssi l'\id bord itéré $\JK_\gT(\xzn)$ contient $1$.
\item On a $\Kdim\gT\leq n$ \ssi pour tout $x$, 
$\Kdim \gT\ul x\leq n-1$. %\index{dimension de Krull!d'un \trdiz}
\end{enumerate}
\end{fact}
%--------- fin fact ---------------------------------------------- 

%:     fact{propDDKagH}
\begin{fact}\label{propDDKagH}
Dans une \agHz,
 tout \id bord de Krull itéré est principal: $\JK_\gT(x)=\dar{\big(x\vu\lnot x\big)}$
et plus \gnltz,
%---- equation {} ----
\begin{equation}
   \!\!\!  \JK_\gT(\xzn) =\dar{\big(x_{n}\vu \big(x_{n}\im( \cdots (x_1 \vu (x_1 \im (x_0\vu \neg x_0)))\cdots)\big)\big)}
\end{equation}
 \end{fact}
%
%\begin{proof}
%\end{proof}
%

%:     Lemma{lemZarAHeyt}
\begin{lemma}\label{lemZarAHeyt}
Soient $\fa$, $\fb$ deux \itfs d'un anneau $\gA$. \\
Dans le treillis $\ZarA$, l'\elt
$\DA(\fa)\im\DA(\fb)$ existe \ssi l'\idz~$(\fb:\fa^{\infty})$  a même radical qu'un \itfz.
\end{lemma}
%--------- fin lemma ---------------------------------------------- 
%
\begin{proof}
Dans un \trdiz,
l'\elt $u\im v$  existe si l'\id $(v:u)$ est principal (son \gtr est alors noté $u\im v$). 
Or pour un \id $\fa$ \tfz,  $\big(\DA(\fb):\DA(\fa)\big)=\DA(\fb:\fa^{\infty})$. 
%:HHH raccourci
D'où le résultat annoncé.
%Donc $(\DA(\fb):\DA(\fa)\big)$ est principal dans $\ZarA$
%\ssi $(\fb:\fa^{\infty})$ a même radical qu'un \itfz.
\end{proof}
%
%:     Lemma{lem2ZarAHeyt}
\begin{lemma}\label{lem2ZarAHeyt}
Supposons que $\ZarA$ soit une \agHz. \\
Pour $(\xzn)$ dans $\gA$, on a l'\egt 

\snic{\DA\big(\JK_\gA(\xzn)\big)=\JK_\ZarA\big(\DA(x_0),\ldots,\DA(x_n)\big) .}
\end{lemma}
%--------- fin lemma ----------------------------------------------
\facile 

%:     Proposition{propNoetAgH}
\begin{proposition}\label{propNoetAgH}
Soit $\gA$ un \cori \noez.
\begin{enumerate}
\item Si $\fa$ et $\fb$ sont deux \itfsz,  l'\id $(\fb:\fa^{\infty})$ est \tfz.
\item $\ZarA$ est une \agHz, avec $\DA(\fa)\im\DA(\fb)=\DA(\fb:\fa^{\infty})$.  
\item Les \ids bords de Krull itérés
définis \paref{eq1IdBordKrullItere} ont même radical que des \itfsz.
\item Si en outre $\gA$ est  \fdiz, $\ZarA$ est discret
et l'on dispose d'un test pour décider si une suite dans
$\gA$ admet une suite \copz.
\end{enumerate}
\end{proposition}

\begin{proof} 
\emph{1.} Soient  $\fa$, $\fb\in \Zar\gA$.  
Notons $\fJ_k=(\fb:\fa^k)$.
Puisque $\gA$ est \cohz, chaque \id $\fJ_k$ est \tfz. Puisque  $\gA$ est \noez, la suite admet deux termes consécutifs égaux, 
par exemple d'indices $p$ et $p+1$, à partir desquels
il est clair qu'elle devient stationnaire.
On a alors $\fJ_p=(\fb:\fa^{\infty})$.
\\
\emph{2.} Conséquence de \emph{1} vu le lemme \ref{lemZarAHeyt}.
\\
\emph{3.} Résulte par \recu de \emph{2} vu le fait \ref{propDDKagH} et le lemme \ref{lem2ZarAHeyt}.
\\
\emph{4.}
Résulte de \emph{2}, du fait \ref{propDDKagH}  et du lemme \ref{lem2ZarAHeyt}.
\end{proof}

%%%%%%%%%%%%%%%%%%%%%%%%%%%%%%%%%%%%%%%%%%%%%%%%%%%%%%%%%%%%%%%%%%%%%%%%%%%
%%%%%%%%%%%%%%%%%%%%                                 %%%%%%%%%%%%%%%%%%%%%%
%%%%%%%%%%%%%%%%%%%%    Dimension des morphismes     %%%%%%%%%%%%%%%%%%%%%%
%%%%%%%%%%%%%%%%%%%%                                 %%%%%%%%%%%%%%%%%%%%%%
%%%%%%%%%%%%%%%%%%%%%%%%%%%%%%%%%%%%%%%%%%%%%%%%%%%%%%%%%%%%%%%%%%%%%%%%%%%
\section{Dimension des morphismes}
\label{secKdimMor}

%:  Définition
\subsec{Définition et premières \prtsz}

%:     Definition{defiKdimMor}
\begin{definition}\label{defiKdimMor}
Si $\rho:\gA\to\gB$ est un \homo d'anneaux, \emph{la \ddk du morphisme} $\rho$
est par \dfn la \ddk de l'anneau $\Abul\otimes_\gA\!\gB$
obtenu par le changement d'anneau de base qui remplace $\gA$ par sa clôture
\zede réduite $\Abul$ (définie dans le paragraphe \paref{secClotureZEDR}).
\end{definition}

\exls ~\\
1) Si $\gk$ est \zedz,  on a vu que $\Kdim\kXn\leq n$. On en déduit
que la \ddk du morphisme $\gA\to\AXn$ est~$\leq n$, avec \egt si $\gA$
est non trivial.
\\
2) Si $\gB$ est une \Alg entière, après
\eds l'\alg est entière sur $\Abul$, donc \zedez. Ainsi, le morphisme
$\gA\to\gB$ est \zedz.
\eoe

%:     Lemma{lem1KdimMor}
\begin{lemma}\label{lem1KdimMor}
Soit $\gB$  et $\gC$ deux \Algsz. 
Alors par \eds on obtient $\Kdim(\gC\to\gC\otimes_\gA\!\gB)\leq \Kdim(\gA\to\gB)$
dans les cas suivants.
\begin{enumerate}
\item $\gC$ est un quotient de $\gA$, ou un localisé de $\gA$, ou le quotient d'un localisé de $\gA$.
\item  $\gC$ est un produit fini d'anneaux du type précédent.
\item  $\gC$ est une limite inductive filtrante d'anneaux du type précédent.  
\end{enumerate}

\end{lemma}
\begin{proof} On utilise la constatation $\gC\bul\otimes_\gC(\gC\otimes_\gA\gB)\simeq\gC\bul\otimes_\gA\gB
\simeq \gC\bul\otimes_{\Abul}(\Abul\otimes_\gA\gB)$.
On vérifie ensuite que le foncteur $\gB\mapsto\gB\bul$
transforme un quotient en un quotient, un localisé en un localisé
(proposition \ref{propClZdrLoc}),
un produit fini en un produit fini, et une limite inductive filtrante 
en une limite inductive filtrante. Par ailleurs, l'\eds commute aussi avec toutes ces constructions. Enfin la dimension de Krull ne peut que diminuer par ces constructions.
\end{proof}

\rem{Il n'est pas vrai que $\gC\otimes_{\gA}\gB$
soit \zed dès que les trois anneaux le sont.
Par exemple on peut prendre $\gA$ un \cdi et~$\gB=\gC=\gA(X)$.
Alors $\Kdim(\gC\otimes_\gA\gB)=1$ (voir l'exercice \ref{exoKdimSomTr}).
 Il s'ensuit que l'\edsz, même dans le cas d'une extension \fptez, 
peut augmenter strictement la dimension de Krull des morphismes.
A contrario on a le \plgc suivant. 
\eoe}

%:     PrincipeLocGlob{plcc.KdimMor}
\begin{plcc}\label{plcc.KdimMor}\relax
Soient $S_1$, $\ldots$, $S_n$  des \moco d'un anneau $\gA$, 
un entier $k\geq-1$ et $\gB$ une \Algz.
La \ddk du morphisme $\gA\to\gB$ est $\leq k$ \ssi 
la \ddk de chacun des morphismes
$\gA_{S_i}\to\gB_{S_i}$ est $\leq k$.  
\end{plcc}
%--------- fin plcc ---------------------------------------------- 
%
\begin{proof}
Comme $(\gA_{S_i})\bul \simeq (\gA\bul)_{S_i}$
(proposition \ref{propClZdrLoc}), on obtient 

\snic{(\gA_{S_i})\bul\otimes _{\gA_{S_i}}\gB_{S_i}
\simeq(\gA\bul\otimes_{\gA}\gB)_{S_i},}

et l'on est ramené au \plgref{thDdkLoc}.
\end{proof}

\medskip 
Le but de cette section est de montrer,
pour un morphisme $\rho:\gA\to\gB$, l'in\egt
$$\preskip.0em \postskip.4em
\fbox{$1+\Kdim\gB\,\leq\, (1+\Kdim\gA)(1+\Kdim\rho)$}.
$$

Notons que pour $\Kdim\gA\leq0$ on a trivialement $\Kdim\gB=\Kdim\rho$.
 Nous traitons ensuite un cas simple mais non trivial pour y voir clair.
Le cas vraiment simple serait celui où $\gA$ est intègre et $\Kdim\gA\leq1$.
Comme la \dem est inchangée, nous supposerons seulement $\gA$ \qiz, ce qui
facilitera la suite.

%:     Proposition{prop1KdimMor}
\begin{proposition}\label{prop1KdimMor} 
Soit  $\rho:\gA\to\gB$ un morphisme, avec $\gA$ \qiz.\\
 Si $\Kdim \rho \leq n$ et $\Kdim\gA\leq 1$, alors $\Kdim\gB\leq 2n+1$.
\end{proposition}
\begin{proof}
Soit $\uh=(h_0,\ldots,h_{2n+1})$ une suite de $2n+2$ \elts dans
$\gB$. On doit montrer qu'elle est singulière.\\
Par hypothèse l'anneau $\Abul\otimes_\gA\gB$ est de \ddk $\leq n$.\\
 Soit $\gK=\Frac\gA$ l'anneau total des fractions, il est \zed réduit et engendré par
 $\gA$ comme anneau \zed réduit, donc c'est un quotient de $\Abul$.
On en conclut que  la suite  $(h_0,\ldots,h_{n})$ est singulière dans~$\wi\gB=\gK\otimes_\gA\!\gB$. \\
 Cela signifie que l'\id bord itéré
$\IK_{\wi\gB}(h_0,\ldots,h_{n})$ contient $1$, et {en chassant les \denosz}
que  $\IK_{\gB}(h_0,\ldots,h_{n})$ contient un $a\in\Reg(\gA)$.\\
Donc  $\gB_0=\gB\sur{\IK_\gB(h_0,\ldots,h_{n})}$
est un quotient de $\gB\sur{a\gB}=\gA\sur{a\gA}\otimes_\gA\!\gB$.
Puisque $a$ est \ndz et $\Kdim\gA\leq1$, le quotient $\gA\sur{a\gA}$
est \zedz, donc $(\gA\sur{a\gA})\red$ est un quotient de $\Abul$
et l'anneau $(\gB_0)\red$
est un quotient de  $\Abul\otimes_\gA\!\gB$.
On en déduit que la suite $(h_{n+1},\ldots,h_{2n+1})$ est singulière 
dans  $(\gB_0)\red$,
donc aussi dans $\gB_0$.
\\
Donc l'anneau $\gB\sur{\IK_{\gB}(\uh)}=
\gB_0\sur{
\IK_{\gB_0}(h_{n+1},\ldots,h_{2n+1})}$ est  trivial.
\end{proof}

Pour passer du cas \qi au cas \gnl on a envie de dire que tout anneau
réduit peut se comporter dans les calculs comme un anneau intègre à condition de remplacer $\gA$ par
$$
\gA\sur{\Ann_\gA (a) }\times \gA\sur{\Ann_\gA(\Ann_\gA( a) )}
$$
lorsqu'un \algo demande de savoir si l'annulateur de $a$ est égal à $0$
ou $1$.
La chose importante dans cette construction est que le principe de recouvrement fermé pour
les \susis s'applique puisque le produit des deux \ids
$\Ann_\gA (a)$
et $\Ann_\gA(\Ann_\gA (a) )$ est nul.

Ce type de  \dem sera sans doute plus facile à saisir quand
on sera familiarisé avec la machinerie \lgbe de base expliquée
\paref{MethodeIdeps}.\imlb
Ici nous ne procédons pas par \lons \come successives mais par
\gui{recouvrements fermés} successifs.

\smallskip
En fait nous n'allons pas introduire d'arbre de calcul dynamique en tant que tel,
nous construirons plutôt un objet universel
qui en tient lieu. Cet objet universel est une \gui{approximation finitaire \covz}
du produit de tous les quotients de $\gA$ par ses \idemisz, un objet
des \clama un peu trop idéal pour pouvoir être considéré \cotz,
du moins sous la forme que l'on vient de définir: en fait, si~$\gB$ est ce produit
et si $\gA_1$ est l'image naturelle de $\gA$ dans $\gB$, alors l'anneau 
\uvl que nous construisons devrait être égal à la clôture \qi de~$\gA_1$
dans $\gB$, du moins en \clamaz.

%:  Cloture \qi minimale d'un anneau reduit
\subsec{Clôture \qi minimale d'un anneau réduit}
Dans la suite  nous notons $a\epr$ l'\id annulateur de l'\elt $a$ lorsque le contexte
 est clair (ici le contexte est simplement l'anneau dans lequel on doit
 considérer~$a$). Nous utiliserons aussi la notation $\fa\epr$ pour l'annulateur d'un \idz~$\fa$.

Les faits énoncés ci-après sont \imdsz.
\begin{eqnarray}
\fa & \subseteq & (\fa\epr)\epr \label{ann1}\\
\fa\subseteq\fb & \Longrightarrow & \fb\epr \subseteq \fa\epr\label{ann2}\\
\fa\epr & = & \big((\fa\epr)\epr\big)\epr\label{ann3}\\
(\fa+\fb)\epr & = & \fa\epr \cap \fb\epr\label{ann4}\\
\fa\epr\subseteq\fb\epr & \Longleftrightarrow &(\fa+\fb)\epr=\fa\epr
\label{ann5}\\
\fa\epr\subseteq\fb\epr & \Longleftrightarrow & (\fb\epr)\epr \subseteq (\fa\epr)\epr\label{ann6}\\
(\fa\epr:\fb)&=& (\fa\fb)\epr \label{ann7}
\\
(\gA\sur{\fa\epr})\big/{{~\ov\fb~}\epr}& = &  \gA\sur{(\fa\fb)\epr} \label{ann8}
\end{eqnarray}

\rems 1) Un \id $\fa$ est un annulateur (d'un autre \idz) \ssi $\fa=(\fa\epr)\epr$.
\\
2) L'inclusion $\fa\epr+\fb\epr \subseteq (\fa\cap\fb)\epr$
peut être stricte, même si $\fa=\fa_1\epr$ et~$\fb=\fb_1\epr$.
Prenons par exemple $\gA=\ZZ[x,y]=\aqo{\ZZ[X,Y]}{XY}$, $\fa_1=\gen{x}$
et $\fb_1=\gen{y}$. Alors, $\fa=\fa_1\epr=\gen{y}$,  $\fb=\fb_1\epr=\gen{x}$,
$\fa\epr+\fb\epr=\gen{x,y}$, et~$(\fa\cap\fb)\epr=\gen{0}\epr=\gen{1}$.\eoe

\medskip
Si nous supposons $\gA$ réduit, nous avons en plus les résultats suivants.
\begin{eqnarray}
\sqrt{\fa\epr}\quad=\quad\fa\epr & = & (\sqrt{\fa})\epr\quad=\quad(\fa^2)\epr \label{ann9}\\
 (\fa\fb)\epr & = &  (\fa\cap \fb)\epr \label{ann10}
\\
\fa\epr\subseteq\fb\epr & \Longleftrightarrow & (\fa\fb)\epr = \fb\epr
\label{ann11}
\end{eqnarray}

%:     Lemma{lem20MorRc}
\begin{lemma}\label{lem20MorRc}
Soit   $\gA$ un  anneau réduit
et $a\in\gA$.
On définit
$$\gA_{\so{a}}\eqdefi\gA\sur{a\epr}\times \gA\sur{({a\epr})\epr}$$
et l'on note $\psi_a:\gA\to\gA_{\so{a}}$  l'\homo canonique. 
\begin{enumerate}
\item $\psi_a(a)\epr$ est engendré par l'\idm $(\ov 0,\wi 1)$,
donc $\psi_a(a)\epr=(\ov 1,\wi 0)\epr$.
\item $\psi_a$ est injectif (on peut identifier $\gA$ à un sous-anneau de $\gA_{\so{a}}$).
\item Soit $\fb$  un \id  dans $\gA_{\so{a}}$, alors l'\id $\psi_a^{-1}(\fb\epr)=\fb\epr\cap\gA$ est un \id annulateur dans $\gA$.
\item L'anneau $\gA_{\so{a}}$ est réduit.
\end{enumerate}
\end{lemma}
\begin{proof} \emph{1.}  On a $\psi_a(a)=(\ov a,\wi 0)$, où $\ov x$ est la classe modulo
$a\epr$ et $\wi x$ est la classe modulo
$(a\epr)\epr$. Si $c=(\ov y, \wi z)$, l'\egt $\psi_a(a)c=0$ signifie $\ov{ya}=0$,
i.e. $ya^2=0$, ou encore $ya=0$, \cad $\ov y=\ov 0$.
\\
\emph{2.} Si $xa=0$ et $xy=0$ pour tout $y\in a\epr$ alors $x^2=0$ donc $x=0$.
\\
\emph{3.} Notons $\psi_1:\gA \to\gA\sur{a\epr}$ et
$\psi_2:\gA \to\gA\sur{(a\epr)\epr}$ les deux \prnsz.
On~a~$\fb=\fb_1\times \fb_2$. Si $x\in\gA$ on a

\snic{\psi_a(x)\in\fb\epr \;\Longleftrightarrow\;\psi_1(x)\fb_1=0 \et \psi_2(x)\fb_2=0,}

%\sni
i.e.  $x\in\psi_1^{-1}(\fb_1\epr)
\cap \psi_2^{-1}(\fb_2\epr)$. L'\egt (\ref{ann8})
nous dit que chaque $\psi_i^{-1}(\fb_i\epr)$ est un \id annulateur.
On conclut avec l'\egt~(\ref{ann4}).
\\
\emph{4.} Dans un anneau réduit, tout \id annulateur $\fb\epr$
est radical: en effet, si~$x^2\fb=0$, alors $x\fb=0$. 
\end{proof}
%

%:     Lemma{lem3MorRc}
\begin{lemma}\label{lem3MorRc}
Soit   $\gA$ réduit
et $a,b\in\gA$. Alors avec les notations du lemme \ref{lem20MorRc}
 les deux anneaux $(\gA_{\so{a}})_{\so{b}}$ et $(\gA_{\so{b}})_{\so{a}}$ sont canoniquement isomorphes.
\end{lemma}
\begin{proof}
L'anneau $(\gA_{\so{a}})_{\so{b}}$ peut être décrit de manière \smq comme suit:

\snic{\gA_{\so{a,b}}=\gA\sur{(ab)\epr}\times \gA\sur{(ab\epr)\epr}\times \gA\sur{(a\epr b)\epr}\times \gA\sur{(a\epr b\epr)\epr},}

%\sni
et si $\psi:\gA\to\gA_{\so{a,b}}$ est l'\homo canonique,
il est clair que l'on~a~$\psi(a)\epr=(1,1,0,0)\epr$ et $\psi(b)\epr=(1,0,1,0)\epr$.
\end{proof}

\rem Le cas où $ab=0$ est typique: quand on le rencontre, on voudrait
bien scinder l'anneau en composantes où les choses sont \gui{claires}.
La construction précédente donne alors les trois composantes

\snic{\gA\sur{(ab\epr)\epr}, \; \gA\sur{(a\epr b)\epr}\, \hbox{ et } \,\gA\sur{(a\epr b\epr)\epr}.}

%\sni
Dans la première, $a$ est \ndz et $b=0$, dans la seconde
$b$ est \ndz et $a=0$, et dans la troisième $a=b=0$.
\eoe

\medskip
Le lemme suivant qui concerne les anneaux \qis est recopié du lemme~\ref{lem2SousZedRed} qui concernait les anneaux \zeds réduits (\llec pourra aussi à très peu près recopier la \demz).

%:     Lemma{lem2qi}
\begin{lemma}\label{lem2qi}
Si $\gA\subseteq\gC$ avec $\gC$ \qiz, le plus petit sous-anneau \qi de $\gC$ contenant
$\gA$ est égal à $\gA[(e_a)_{a\in\gA}]$,
où $e_a$ est l'\idm de $\gC$ tel que $\Ann_\gC(a)=\gen{1-e_a}_\gC$. Plus \gnlt si $\gA\subseteq\gB$
avec $\gB$ réduit, et si tout \elt $a$ de $\gA$ admet un annulateur dans $\gB$
engendré par un \idm $1-e_a$,
alors le sous-anneau $\gA[(e_a)_{a\in\gA}]$ de $\gB$ est \qiz.
\end{lemma}

%:     Theorem{thAmin}
\begin{thdef}\label{thAmin} \emph{(Clôture \qi minimale)}
\\
Soit $\gA$ un anneau réduit.
On peut définir un anneau $\Amin$
 comme limite inductive filtrante 
 en itérant la construction de base qui consiste à
 remplacer~$\gE$ (l'anneau \gui{en cours}, qui contient $\gA$) par

\snic{\gE_{\so{a}}\eqdefi\gE\sur{a\epr}\times \gE\sur{({a\epr})\epr}=\gE\sur{\Ann_\gE (a) }\times \gE\sur{\Ann_\gE(\Ann_\gE (a) )},}

%\sni
 lorsque $a$ parcourt $\gA$.
\begin{enumerate}
\item  Cet anneau $\Amin$ est \qiz, contient $\gA$ et est entier sur $\gA$.
\item Pour tout $x\in\Amin$,
$%\Ann_{\Amin}(x)
x\epr\cap\gA$ est un \id annulateur dans $\gA$.
\end{enumerate}
 Cet anneau $\Amin$ est appelé la \emph{clôture \qi minimale de $\gA$}.
 \\
Lorsque $\gA$ n'est pas \ncrt réduit,
  on prendra $\gA\qim\eqdefi (\Ared)\qim$.
\end{thdef}
\index{cloture@clôture!quasi intègre minimale}
\begin{proof} \emph{1.} D'après le lemme \ref{lem2qi}, il suffit de rajouter un \idm $e_a$ pour chaque $a\in\gA$
pour obtenir un anneau \qiz. La limite inductive est bien définie gr\^ace à la relation de commutation donnée par le lemme \ref{lem3MorRc}.
\\
Pour le point \emph{2} on note que $x$ est obtenu à un étage fini de la construction, et que $x\epr\cap\gA$ ne change plus à partir du moment où $x$
est atteint parce que les \homos successifs sont des injections. On peut donc faire appel au point \emph{3} du lemme~\ref{lem20MorRc}.
\end{proof}
\rem
On peut se demander si $\Amin$ ne pourrait pas être \care par une \prt \uvle liée au point \emph{2.}
\eoe

\medskip 
En \gui{itérant} la description de  $(\gA_{\so{a}})_{\so{b}}$
donnée dans
la \dem du lemme \ref{lem3MorRc} on obtient 
la description suivante de chaque anneau
obtenu à un étage fini de la construction de $\Amin$
(voir l'exercice \ref{exoAminEtagesFinis}).

%:     Lemma{lem4MorRc}
\begin{lemma}\label{lem4MorRc}
Soit $\gA$ un anneau réduit et $(\ua) = (\an)$ une suite de~$n$ \elts de
$\gA$.  Pour $I\in\cP_n$, on note $\fa_I$ l'\id

\snic {
\fa_I = \big(\prod_{i\in I} \gen{a_i}\epr \prod_{j\notin I} a_j\big)\epr
= \big(\gen{a_i, i \in I}\epr \prod_{j\notin I} a_j\big)\epr
.}

%\sni
Alors $\Amin$ contient l'anneau suivant, produit de $2^n$
anneaux quotients de~$\gA$ (certains éventuellement nuls):

\snic {
\gA_{\so\ua} = \prod_{I\in\cP_n} \gA\sur{\fa_I}.
} 
\end{lemma}
%--------- fin lemma ---------------------------------------------- 

%:     Fact{fact2Amin}
\begin{fact}\label{fact2Amin}~
\begin{enumerate}
\item Soit $\gA$ un anneau \qiz. 
\begin{enumerate}
\item $\Amin=\gA$.
\item $\AX$ est \qiz, et $\BB(\gA)=\BB(\AX)$.
%
%\item Donc $\Amin[X]=\AX=\AX\qim$, même chose $\AXn$.
%
\end{enumerate}

\item Pour tout anneau $\gA$ on a un \iso canonique

\snic{\Amin[\Xn]\simeq(\AXn)\qim.}
\end{enumerate}
\end{fact}
\begin{proof} \emph{1a.} Résulte de la construction de $\Amin.$
\\
\emph{1b.} Le résultat est clair pour les anneaux intègres.
On peut appliquer la machinerie \lgbe \elr \paref{MethodeQI}.
On pourrait aussi utiliser le lemme de McCoy,
 corolaire \ref{corlemdArtin}~\emph{2.}\imlg
\\
\emph{2.} On suppose \spdg l'anneau $\gA$ réduit.
Il suffit aussi de traiter le cas d'une variable.
Vu le lemme \ref{lem3MorRc} on peut \gui{commencer} la construction
de  $\AX\qim$ avec les constructions $\gE\leadsto \gE_{\so{a}}$ pour des $a\in\gA$.
Mais si $\gE=\gB[X]$ et $a\in\gA\subseteq\gB$ alors $\gE_{\so{a}}=\gB_{\so{a}}[X]$.
Ainsi $\Amin[X]$ peut être vu comme une première étape de la construction
de $\AX\qim$. Mais puisque d'après le point \emph{1} $\Amin[X]$ est \qi
et que pour un anneau \qi $\gC$ on a $\gC=\gC\qim$, la construction de $\AX\qim$
est terminée avec $\Amin[X]$.
\end{proof}
%

%\medskip
\comm
En pratique, pour utiliser l'anneau $\Amin$, on n'a besoin que des étages finis de la construction.
On peut noter cependant que même un seul étage de la construction est un peu mystérieux,
dans la mesure où les \ids $a\epr$ et $(a\epr)\epr$ sont difficiles
à maîtriser. C'est seulement dans le cas des \coris que l'on sait les décrire par des \sgrs finis.
En fait si l'anneau est \noez, la construction doit s'arrêter en un nombre fini d'étapes (au moins du point de vue des \clamaz), et elle remplace l'anneau par le produit de ses quotients par les \idemisz. Nous sommes ici dans une situation où la construction de $\Amin$ répondant aux standards des \coma semble plus compliquée que le résultat
en \clama (au moins si l'anneau est \noez). Néanmoins, puisque nous n'avons pas besoin de connaître les \idemisz, notre méthode est plus \gnle (elle ne nécessite pas le principe du tiers exclu). En outre, sa complication est surtout apparente. Quand nous utilisons $\gA\sur{a\epr}$ par exemple, nous faisons en fait des calculs dans $\gA$  en forçant $a$ à être \ndzz, \cad en annulant par force tout $x$ qui se présente et qui annule~$a$.
Quand nous utilisons $\gA\sur{(a\epr)\epr}$, la chose est moins facile, car a priori, nous avons besoin d'une preuve (et non du simple résultat d'un calcul)
pour certifier qu'un \elt $x$ est dans $(a\epr)\epr$.\\
C'est un fait que l'utilisation des \idemis dans un raisonnement de \clama peut en \gnl être rendue inoffensive (\cad \covz) par l'utilisation de $\Amin$ (ou d'un autre anneau universel du même type\footnote{$\Amin$ correspond à l'utilisation de tous les quotients par les \idemisz, $\Frac(\Amin)$ correspond à l'utilisation de tous les corps de fractions de ces quotients.}),
même si l'on  ne dispose pas d'autre moyen  pour \gui{décrire un \id
$(a\epr)\epr$} que celui d'appliquer la \dfnz.
\eoe

%:  Application
\subsec{Application}

%:     Corollary{corthAmin}
\begin{corollary}\label{corthAmin}
Soit  $\rho:\gA\to\gB$ un morphisme de \ddk finie. On \gui{étend les scalaires}
de $\gA$ à $\Aqim$: on obtient
$\gB'=\Aqim\otimes_\gA\!\gB$  et l'on note $\rho':\Aqim\to\gB'$ le morphisme naturel.
\\
Alors $\Kdim\Aqim=\Kdim\gA$, $\Kdim\gB'=\Kdim\gB$ et $\Kdim\rho'\leq\Kdim\rho$.
\end{corollary}
\begin{proof}
Les deux premiers points résultent du fait que dans la construction de l'anneau $\Aqim$,
%:2012 le  .  remplacé par  ,  ci-dessus
à chaque étape \elr
$$
\gE \quad \rightsquigarrow
\quad\gE\sur{\Ann_\gE (a) }\times \gE\sur{\Ann_\gE(\Ann_\gE a )},
$$
le produit des deux \ids est nul, ce qui se retrouve après tensorisation
par~$\gB$. Donc, le principe de recouvrement fermé pour la \ddkz~\vref{thDdkRecFer}
s'applique. 
Enfin l'in\egt $\Kdim\rho'\leq\Kdim\rho$ résulte du lemme~\ref{lem1KdimMor}.
\end{proof}

\rem De manière \gnle l'anneau $\Frac\Aqim$ semble un meilleur concept
que~$\Abul$ pour remplacer le corps de fractions dans le cas d'un anneau réduit non intègre.
Dans le cas où $\gA$ est \qiz, on a en effet
$\Aqim=\gA$, donc~$\Frac\Aqim=\Frac\gA$, tandis que $\Abul$ est en \gnl
nettement plus encombrant (comme le montre l'exemple $\gA=\ZZ$).
\eoe

%:     corollary {prop2KdimMor}
\begin{corollary}\label{prop2KdimMor} 
Soit  $\rho:\gA\to\gB$ un morphisme. \\
Si $\Kdim \rho \leq n$ et~$\Kdim\gA\leq 1$, alors
$\Kdim\gB\leq 2n+1$.
\end{corollary}
\begin{proof}
Cela résulte clairement de la proposition \ref{prop1KdimMor}
et du corolaire \ref{corthAmin}.
\end{proof}
%
%%%%%%%%%%%%%%%%%%%%%%%%%%%%%%%%%%%%%%%%%
\entrenous{\'Etrange que seule la dimension de $(\Aqim)\bul\otimes_{\gA}\gB$
soit en fin de compte pertinente ici? Non, ce que dit cette preuve \gui{raffinée}
du corolaire \ref{prop2KdimMor}, c'est que l'on  pourrait considérer d'une part la longueur maximale
des chaînes de premiers dans $\gB$ au dessus d'un premier minimal de $\gA$, d'autre part
 la longueur maximale
des chaînes de premiers dans $\gB$ au dessus d'un premier maximal de $\gA$
et faire la somme des deux.}
%%%%%%%%%%%%%%%%%%%%%%%%%%%%%%%%%%%%%%%%%

%:     Theorem{thKdimMor}
\begin{theorem}\label{thKdimMor}
Soit  $\rho:\gA\to\gB$ un morphisme.\\ 
Si $\Kdim \rho \leq n$ et $\Kdim\gA\leq m$, alors
$\Kdim\gB\leq mn+m+n$.
 \end{theorem}
%------------------------
%
\begin{proof}
On fait une \dem par \recu sur $m$. Le cas $m=0$ est trivial.
La preuve donnée pour $m=1$ dans le cas où $\gA$ est \qi
(proposition \ref{prop1KdimMor}), qui s'appuyait sur la dimension $0$
pour prouver le résultat en dimension $m=1$,
 s'adapte sans \pb pour passer de la dimension~$m$ à la dimension $m+1$.
 Nous recopions la \dem dans le cas où~$\gA$ est \qiz.\\
  Pour passer au cas d'un anneau arbitraire
 on utilise le corolaire~\ref{corthAmin}.\\
On suppose donc $\gA$ \qi et on considère  une suite $(\uh)=(h_0,\ldots,h_{p})$   dans
$\gB$ avec  $p=(m+1)(n+1)-1$. On doit montrer qu'elle est singulière.\\
Par hypothèse l'anneau $\Abul\otimes_\gA\!\gB$ est de \ddk $\leq n$.
L'anneau total des fractions $\gK=\Frac\gA$  est \zedrz, et il est engendré par~$\gA$ comme anneau \zedrz, donc c'est un quotient de~$\Abul$.
On en conclut que  la suite  $(h_0,\ldots,h_{n})$ est singulière dans
l'anneau~$\wi\gB=\gK\otimes_\gA\!\gB$. \\
 Cela signifie que l'\id bord itéré
$\IK_{\wi\gB}(h_0,\ldots,h_{n})$ contient $1$, et {en chassant les \denosz}
que  $\IK_\gB(h_0,\ldots,h_{n})$ contient un $a\in\Reg(\gA)$.
Donc l'anneau $\gB_0=\gB\sur{\IK_\gB(h_0,\ldots,h_{n})}$
est un quotient de~$\gB\sur{a\gB}=\gA\sur{a\gA}\otimes_\gA\!\gB$.
Puisque $a$ est \ndz et $\Kdim\gA\leq m$, le quotient $\gA\sur{a\gA}$
est de \ddk $\leq m-1$. L'\homo naturel
$\gA\sur{a\gA}\to\gB\sur{a\gB}$
reste de \ddk $\leq n$ (lemme \ref{lem1KdimMor}). Donc, par \hdrz, la suite $(h_{n+1},\ldots,h_{p})$ est singulière dans $\gB\sur{a\gB}$.
Donc la suite~$(h_{n+1},\ldots,h_{p})$ est singulière dans $\gB_0$.
\\
En conclusion, l'anneau $\gB\sur{\IK_{\gB}(\uh)}=
\gB_0\sur{\IK_{\gB_0}(h_{n+1},\ldots,h_{p})}$ est  trivial.
 \end{proof}

%:     Corollary{corthKdimMor}
\begin{corollary}\label{corthKdimMor}
Supposons  $\Kdim\gA\leq m$.  Alors

\snic{\Kdim\AXn\leq mn+m+n.}
\end{corollary}
\begin{proof}
On sait que si $\gK$ est \zed réduit, $\Kdim\KXn\leq n$.
Ainsi $\Kdim(\gA\to\AXn)\eqdefi\Kdim\Abul[\Xn]\leq n$. On applique le \thoz~\ref{thKdimMor}.
\end{proof}
%
%:HHH rajout
On dispose également d'une minoration de
$\Kdim \AXn$.
%:     Lemma{lemKdimAxn}
\begin{lemma}\label{lemKdimAxn}
Pour tout anneau $\gA$ non trivial et tout $n>0$ on a 

\snic{n+\Kdim\gA\leq \Kdim\AXn.}

%\sni
Plus \prmtz, l'implication suivante est satisfaite pour tout $k\geq -1$
et pour tout anneau:

\snic{\Kdim\AXn\leq n+k\;\Longrightarrow\;\Kdim\gA\leq k}
\end{lemma}
%--------- fin lemma ---------------------------------------------- 
%
\begin{proof} Conséquence \imde de la proposition \ref{lemRegsing}.
\end{proof}
%

%:     theorem  Corollary{cor2thKdimMor}
\begin{theorem}\label{cor2thKdimMor}
On considère une \alg $\rho:\gA\to\gB$.
\begin{enumerate}
\item Supposons que  $\gB$ est engendrée par des \elts primitivement  \agqs sur $\gA$,
alors $\Kdim \rho\leq0$ et donc $\Kdim\gB\leq\Kdim\gA$.
\item Si $\rho$ est injectif et $\gB$  entier sur $\gA$, alors $\Kdim\gB=\Kdim\gA$.
\end{enumerate}
 \end{theorem}

\begin{proof}
\emph{1.} Vu le fait \ref{factZedPrimFin}, l'anneau $\Abul\otimes_\gA\gB$ est \zedz,
autrement dit  $\Kdim\rho\leq0$. On conclut par le \thoz~\ref{thKdimMor}.
\\
\emph{2.} D'après le point \emph{1} et la proposition \ref{propDKEXENT}.
\end{proof}
Pour une \dem  plus directe de l'in\egt $\Kdim\gB\leq\Kdim\gA$, voir l'exercice~\ref{exoInclusionBordLionel}.

%%%%%%%%%%%%%%%%%%%%%%%%%%%%%%%%%%%%%%%%%%%%%%%%%%%%%%%%%%%%%%%%%%%%%%%%%%%
%%%%%%%%%%%%%%%%%%%%                                 %%%%%%%%%%%%%%%%%%%%%%
%%%%%%%%%%%%%%%%%%%%    Dimension valuative          %%%%%%%%%%%%%%%%%%%%%%
%%%%%%%%%%%%%%%%%%%%                                 %%%%%%%%%%%%%%%%%%%%%%
%%%%%%%%%%%%%%%%%%%%%%%%%%%%%%%%%%%%%%%%%%%%%%%%%%%%%%%%%%%%%%%%%%%%%%%%%%%

\penalty-2500
\section{Dimension valuative}
\label{secValdim}

\vspace{3pt}
%:  Dimension des \advsz
\subsec{Dimension des \advsz}

Rappelons qu'un \adv est un anneau réduit dans lequel on a, pour tous $a,\,b$:
$a$ divise $b$ ou $b$ divise $a$. Autrement dit c'est un \alo de Bézout et réduit.
Un \adv est un anneau normal, local et \sdzz.
Il est intègre \ssi il est \cohz.

%:HHH le morceau ci apres a ete mis dans le chapitre 8
%Un ensemble ordonné $(E,\leq)$ est dit \ixc{totalement ordonné}{ensemble ---}
%si pour tous $x$, $y$ on a $x\leq y$ ou $y\leq x$. A priori on ne le suppose pas discret et l'on n'a donc pas de test pour l'in\egt stricte.

Il est clair que le treillis de Zariski d'un \adv est un 
ensemble totalement ordonné.
%:     Fact{fact1SeqSingTD}
\begin{fact}\label{fact1SeqSingTD}
Dans un \trdi si une sous-suite
de $(\ux)=(\xn)$ est singulière,
la suite $(\ux)$ est singulière.
\end{fact}
\begin{proof}
On considère une \susi $(\yr)$, avec une suite \cop $(b_1,\ldots,b_r)$.
Rajoutons un terme $z$ à $(\yr)$.  Pour en obtenir une suite \copz,
on procède comme suit.
Si $z$ est mis à la fin, on rajoute
$1$ à la fin de $(b_1,\ldots,b_r)$. Si $z$ est mis au début, on rajoute $0$
au début de~$(b_1,\ldots,b_r)$. Si $z$ est intercalé entre $y_i$ et $y_{i+1}$
on intercale~$b_i$ entre~$b_i$ et~$b_{i+1}$.
\end{proof}
%
%:     Fact{fact2SeqSingTD}
\begin{fact}\label{fact2SeqSingTD}
Dans  un \trdiz,  si  $(\ux)=(\xn)$ et si l'on \hbox{a $x_1=0$}, ou~$x_n=1$, ou~$x_{i+1}\leq x_i$ pour un $i\in\lrb{1..n-1}$,
alors la suite $(\ux)$ est singulière.
\end{fact}
\begin{proof}
On applique le fait précédent en notant que $(0)$ et $(1)$ sont deux suites \cops et que la suite  $(x_i,x_{i+1})$ avec $x_{i+1}\leq x_i$ admet $(0,1)$
%:2012 rajout () autour de 0  1  et  0,1 ci-dessus
pour suite \copz.
\end{proof}

Un rappel: la signification \cov de la phrase
\gui{le nombre d'\elts de $E$ est borné par $k$} (ce que l'on note $\#E\leq k$)
 est que pour toute liste finie de $k+1$ \elts dans $E$,
il y en a deux égaux.

%:     Lemma{lemSeqSingTD}
\begin{lemma}\label{lemSeqSingTD}
Pour une suite  \emph{croissante} $(\ua)=(a_1,\ldots, a_{n})$ dans un treillis
totalement ordonné 
 \propeq
\begin{enumerate}
\item La suite est singulière.
\item $a_1=0$, ou $a_{n}=1$, ou il existe $i\in\lrb{1..n-1}$
tel que $a_{i}= a_{i+1}$.
\item Le nombre d'\elts dans $(0,\an,1)$ est borné par $n+1$
\end{enumerate}
\end{lemma}
\begin{proof}
\emph{1} $\Rightarrow$ \emph{2.}
Faisons le calcul pour le cas $n=3$ en laissant la \recu \alec sceptique.
On considère une suite \cop $(b_1,b_2, b_3)$.
On a
\[\arraycolsep3pt
\begin{array}{ccc}
1  & \leq  &  a_3\vu b_3 \\
a_3\vi b_3  & \leq  &  a_2\vu b_2 \\
a_2\vi b_2  & \leq  &  a_1\vu b_1 \\
a_1\vi b_1  & \leq & 0
 \end{array}
\]
Ainsi, $a_1=0$ ou $b_1=0$. \\
Si $b_1=0$, alors $a_1\vu b_1=a_1\geq a_2\vi b_2$.
Donc $a_2\leq a_1$ ou $b_2\leq a_1$. Dans le premier cas, $a_1=a_2$.
Dans le deuxième cas, $b_2\leq a_1\leq a_2$ donc $a_2\vu b_2=a_2$.
Ceci implique $a_3\leq a_2$
ou $b_3\leq a_2$. Dans le premier cas,~$a_2=a_3$.
Dans le deuxième cas, $b_3\leq a_2\leq a_3$, donc $a_3\vu b_3=a_3=1$.

 \emph{2} $\Rightarrow$ \emph{1.} D'après le fait \ref{fact2SeqSingTD}.

 \emph{3} $\Rightarrow$ \emph{2.} Si l'on  a deux \elts égaux dans une suite
croissante, alors il y a aussi deux \elts consécutifs égaux.
\end{proof}

\smallskip
Le \tho suivant donne une interprétation \cov précise et \elr de la \ddk
d'un ensemble totalement ordonné. Il résulte directement du fait \ref{fact2SeqSingTD} et du lemme \ref{lemSeqSingTD}.

%:     theorem{thKdimTDTO}
\begin{theorem}\label{thKdimTDTO}
Pour un \trdi totalement ordonné $\gT$, \propeq
\begin{enumerate}
\item $\gT$ est \ddi$n$.
\item Le nombre d'\elts de $\gT$ est borné par $n+2$ ($\#\gT\leq n+2$).
\item Pour toute suite croissante
$(\xzn)$ dans $\gT$, on a $x_0=0$, ou $x_n=1$, ou $x_{i+1}= x_i$ pour un $i\in\lrb{0,n-1}$.
\end{enumerate}
\end{theorem}

Notez que le \tho précédent s'applique au treillis de Zariski
d'un anneau de valuation.
Nous présentons maintenant deux faits fort simples et utiles concernant les \advsz.

%:     Fact{fact1ValRing}
\begin{fact}\label{fact1ValRing}
Soient dans un \adv des \elts $u_1$, \ldots, $u_m$ avec $\sum_iu_i=0$ (et~$m\geq2$).
Alors il existe $j\neq k$ et un \elt \iv $v$ tels que $\gen{u_1,\ldots,u_m}=\gen{u_j}=\gen{u_k}$ et $vu_j=u_k$.
\end{fact}
\begin{proof}
Tout d'abord il existe $j$ tel que $\gen{u_1,\ldots,u_m}=\gen{u_j}$. Soit
alors pour chaque $k$ un \elt $v_k$ tel que $u_k=v_ku_j$, avec
$v_j=1$. \\
On obtient l'\egt $u_j(1+\sum_{k\neq j}v_k)=0$.
Donc $u_j=0$ ou $1+\sum_{k\neq j}v_k=0$.
Si $u_j=0$,  on peut prendre tous les $v_k$ égaux à $1$.\\
Si $1+\sum_{k\neq j}v_k=0$, l'un des $v_k$ est \iv puisque $\gV$ est local.
\end{proof}
%

%:     Fact{fact2ValRing}
\begin{fact}\label{fact2ValRing}
Soit $\gV$ un \adv et une suite $(\an)$ dans~$\gV\etl$.
Pour des exposants $p_i$ tous $>0$,  posons $a=\prod_{i=1}^n a_i^{p_i}$.
Alors il existe un~$j\in\lrbn$ tel que $\rD_\gV(a)=\rD_\gV(a_j)$.
\end{fact}
\begin{proof}
%:HHH oups il fallait permuter i et j dans  $a_j$ divise $a_i$
Considérons un $j$ tel que $a_i$ divise $a_j$ pour tous les $i\in\lrbn$.
Alors $a_j$ divise $a$ qui divise $a_j^p$, où $p=\sum_{i=1}^n p_i$.
\end{proof}
Nous aurons besoin du lemme combinatoire suivant.

%:     Lemma{lemBornes}
\begin{lemma}\label{lemBornes}
Soient deux ensembles $E\subseteq F$. On suppose que pour
toute suite $(x_0,\ldots,x_m)$ dans $F$, une des deux alternatives suivantes a lieu:
\begin{itemize}
\item il existe $i<j\in\lrb{0..m}$ tels que $x_i=x_j$,
\item il existe $i\in\lrb{0..m}$ tel que $x_i\in E$.
\end{itemize}
Alors  $\#E\leq\ell$ implique $\#F\leq\ell+m$.
\end{lemma}
\begin{proof}
On considère une suite $(y_0,\ldots,y_{\ell+m})$ dans $F$. On doit montrer qu'il y a deux termes égaux. On considère les $m+1$ premiers termes. Ou bien il y en a deux égaux, et l'affaire est entendue, ou bien l'un des termes est dans $E$.
Dans ce cas, on supprime ce terme qui est dans $E$ de la suite 
$(y_0,\ldots,y_{\ell+m})$ 
et l'on considère les $m+1$ premiers termes de cette nouvelle suite. 
Ou bien il y en a deux égaux, et l'affaire est entendue, ou bien l'un des termes est dans $E$ \ldots\, Dans le cas pire, on poursuit le processus
jusqu'au bout et l'on obtient à la fin
$\ell+1$ termes dans $E$ et deux d'entre eux sont égaux.
\end{proof}
%

%:     theorem{th1Valdim}
\begin{theorem}\label{th1Valdim}
Soient $\gV$ un \adv intègre, $\gK$ son corps de fractions, $\gL\supseteq \gK$ un \cdi de degré de transcendance $\leq m$ sur~$\gK$,
et~$\gW\supseteq \gV$ un \adv de $\gL$. Alors $\Kdim\gW\leq \Kdim\gV+m$.
\end{theorem}
\begin{proof} 
On doit montrer que si $\,\Kdim\gV\leq n$, alors $\,\Kdim\gW\leq n+m$.
\\
Puisqu'il s'agit d'\advsz, on doit simplement montrer que

\snic{\#\Zar\gV\leq n+2\;$ implique $\;\#\Zar\gW\leq n+m+2.}

%\sni
(Voir le \thrf{thKdimTDTO}.)
Il suffit donc de montrer que les hypothèses du lemme \ref{lemBornes} sont satisfaites
pour les entiers $\ell=n+2$ et $m$, et pour les ensembles~$E=\Zar\gV$ et $F=\Zar\gW$.
\\
Soit $\gV'=\gW\cap\gK$. Puisque $\gV'$ est un localisé de $\gV$, on
a $\Kdim\gV'\leq\Kdim\gV$.  On se ramène ainsi au cas où $\gV=\gW\cap\gK$,
ce qui implique $\Zar\gV\subseteq\Zar\gW$.
\\
Soient maintenant $x_0$, \ldots, $x_{m}\in\Reg \gW$, noté $\gW\etl$. \\
 Considérons une \rde \agq
sur $\gK$ entre $(x_0,\ldots,x_{m})$. On peut supposer que les \coes du \pol
$P\in\gK[X_0,\ldots,X_m]$ qui donne cette \rde \agq sont dans
$\gV\cap\gK\eti=\gV\etl$. En notant, pour $p \in \NN^{m+1}$, $x^p =
x_0^{p_0}\cdots x_{m}^{p_m}$,  le fait~\ref{fact1ValRing} nous donne~$p$ et~$q$ distincts dans $\NN^{m+1}$ tels que~$ax^{p}$ et~$bx^{q}$ sont associés dans $\gW$, 
avec~$a$,~$b\in\gV\etl$. 
En simplifiant par~$x^{p \wedge q}$, on peut supposer $p \wedge q = 0$. 
Puisque~$a$ divise~$b$ ou~$b$ divise~$a$, on peut supposer~$b = 1$.  On
a donc $ax^p$ associé à $x^q$ dans~$\gW$.  Si $q = 0$, alors chaque $x_j$
qui figure dans $x^p$ (il y en a au moins un) est \iv dans $\gW$, i.e.
$\rD_\gW(x_j) = \rD_\gW(1)$.  Sinon, le fait \ref{fact2ValRing} appliqué à~$x^q$ nous donne  un~$x_j$ présent dans~$x^q$ tel que~$\rD_\gW(x^q) =
\rD_\gW(x_j)$; appliqué à $ax^p$, il nous fournit,
$\rD_\gW(ax^p) = \rD_\gW(a)$ ou $\rD_\gW(x_k)$ avec~$x_k$ présent dans~$x^p$;
on a donc~$\rD_\gW(x_j) = \rD_\gW(a)$, ou bien
$\rD_\gW(x_j) = \rD_\gW(x_k)$. La \dem est complète.
 \end{proof}
%

%%%%%%%%%%%%%%%%%%%%%%%%%%%%%%%%%%%%%%%%%%%%%%%%%%%%%%%%%%%%%%%%%%%%%%%%%%%
%: subsec{Dimension valuative d'un anneau commutatif}
\subsec{Dimension valuative d'un anneau commutatif}

%:     Definition{defiValdim}
\begin{definition}\label{defiValdim}~
\begin{enumerate}
\item Si $\gA$ est un anneau \qiz, la \ix{dimension valuative}
est définie comme suit. Soit $d\in\NN\cup\so{-1}$ et $\gK=\Frac\gA$, on dit que la dimension
valuative de $\gA$ est inférieure ou égale à $d$ et l'on écrit
$\Vdim\gA\leq d$ si pour toute suite finie $(\ux)$ dans $\gK$ on a
$\Kdim \Aux \leq d$.

\item Dans le cas \gnl on définit \gui{$\Vdim\gA\leq d$} par \gui{$\Vdim\Amin\leq d$}.
\end{enumerate}
\end{definition}
On a \imdtz:
\begin{itemize}
  \item $\Kdim\gA\leq\Vdim\gA$,
  \item $\Vdim\gA=-1$ \ssi $\gA$ est trivial,
  \item $\Vdim\gA\leq0$ \ssi $\Kdim\gA\leq0$,
  \item si $\gA$ est \qi alors 
\begin{itemize}
\item $\Kdim\gA=\Vdim\gA$ \ssi $\Kdim \gB  \leq \Kdim\gA$
pour tout anneau $\gB$ intermédiaire entre $\gA$ et~$\Frac\gA$,
\item si $\gB$ est intermédiaire entre $\gA$ et~$\Frac\gA$, on a $\Vdim\gB\leq\Vdim\gA$. 
\end{itemize}
\end{itemize}

\smallskip Le fait suivant résulte directement de la construction de $\Amin$.
%:     Fact{fact1Amin}
\begin{fact}\label{fact1Amin}
Si $\gA$ est un \anarz, alors $\Amin$ \egmtz.
\end{fact}

%:     Lemma{lem1Valdim}
\begin{lemma}\label{lem1Valdim}
Si $\gA$ est un \anarz, on a $\Kdim\gA = \Vdim\gA$.
\end{lemma}
\begin{proof}
Puisque $\Kdim\gA=\Kdim\gA\qim$, et puisque $\gA\qim$ est un \anar si $\gA$ en est un, 
il suffit de traiter le cas où $\gA$ est \qiz.
On applique alors le \thref{th.2adpcoh} qui dit que tout \elt de
$\Frac\gA$ est primitivement \agq sur $\gA$, et le \thref{cor2thKdimMor}
qui dit que dans un tel cas $\Kdim \gB  \leq \Kdim\gA$
pour tout anneau $\gB$ intermédiaire entre $\gA$ et~$\Frac\gA$.
 \end{proof}

\rem Voici une fin de preuve (cas où $\gA$ est un
\anar \qiz) moins savante.
On suppose d'abord que $\gA$ est local,
i.e. que c'est un \adv intègre.
Pour tout $x=a/b\in\Frac\gA$, on a l'alternative: $b$ divise $a$,
auquel cas $x\in\gA$, ou~$a$ divise $b$, i.e., $ac=b$ auquel cas $c$ est \ndz et $x=1/c$ de sorte que~$\gA[x]$ est un \adv localisé de $\gA$, donc
$\Kdim\gA[x]\leq\Kdim\gA$. On termine par \recu sur le nombre
d'\elts de $\Frac\gA$ que l'on rajoute à~$\gA$. 
Enfin, dans le cas \gnlz, on reprend la \dem précédente. On remplace l'alternative
\gui{$b$ divise $a$ ou $a$ divise $b$} par la création de deux \lons \come
de $\gA$. Dans la première $b$ divise $a$, dans la seconde $a$ divise $b$.
\eoe

%:     Lemma{lemVdimKdim}
\begin{lemma}\label{lemVdimKdim}
Soit $\gA$ un anneau intègre, $n\geq 1$ et $k\geq -1$.\\
Si $\Kdim\AXn\leq n+k$,
alors pour tous $x_1$, \dots, $x_n$ dans $\Frac\gA$, on a $\Kdim\Axn\leq k$.
\end{lemma}
%--------- fin lemma ---------------------------------------------- 
%
\begin{proof} On introduit les anneaux intermédiaires

\snic{\gB_0=\AXn, \;\gB_1=\gA[x_1,X_2,\dots,X_n],\;\dots,
\;\gB_n=\Axn.}

%\sni
%On note $d_i$ (un majorant de) la \ddk de $\gB_i$.
Pour $i\in\lrbn$, on note $\varphi_i$ l'\homo d'\evn $\gB_{i-1}\to\gB_i$
défini par $X_i\mapsto x_i$. Si $x_i=a_i/b_i$, le noyau $\Ker\varphi_i$
contient $f_i=b_iX_i-a_i$.\\ 
Soit $i\in\lrb{0..n-1}$. 
Puisque $b_{i+1}\in\Reg\gA[(x_j)_{1\leq j\leq i}]$, on a $f_{i+1}\in\Reg\gB_{i}$ (lemme de McCoy, corolaire \ref{corlemdArtin}).
Donc, d'après le point~\emph{5} de la proposition~\ref{propDdk0}, on a 
$\Kdim\aqo{\gB_{i}}{f_{i+1}}\leq \Kdim\gB_i-1$. 
%\\
Enfin, puisque $\gB_{i+1}$ est un quotient
de~$\aqo{\gB_{i}}{f_{i+1}}$, on obtient  $\Kdim\gB_{i+1} \le \Kdim\gB_i-1$.
\end{proof}

Dans la proposition suivante, comme nous le verrons un peu plus loin, les trois \prts sont en fait \eqves (\thref{thValDim} point~\emph{2}).

%:     Proposition{propVdimKdim}
\begin{proposition}\label{propVdimKdim}
Soit $\gA$ un anneau intègre et $n\geq 1$, on a pour les points suivants les implications  {1} $\Rightarrow$  {2} $\Rightarrow$  {3.}  
\begin{enumerate}
\item On a $\Kdim\AXn\leq 2n$.  
\item Pour tous $x_1$, \dots, $x_n$ dans $\Frac\gA$, on a $\Kdim\Axn\leq n$.
\item On a $\Vdim\gA\leq n$.
\end{enumerate}
 
\end{proposition}
%--------- fin proposition ---------------------------------------------- 
%
\begin{proof}
\emph{1}$ \Rightarrow$ \emph{2.} Cas particulier du lemme 
\ref{lemVdimKdim}.

\emph{2} $\Rightarrow$ \emph{3.}
On considère une suite arbitraire $(\yr)$ dans $\Frac\gA$, puis une
suite $(\xzn)$ arbitraire dans $\gB=\Ayr$.   On doit démontrer que la suite $(\xzn)$ est singulière
dans $\gB$. Il suffit de montrer qu'elle est singulière dans 
$\gC=\gA[\xzn]$, ou encore, que la suite $(\xn)$ est singulière dans 
$\gC/\IK_\gC(x_0)$. \\
On écrit $x_0=a_0/b_0$ avec $b_0\in\Reg\gA$.
Si $a_0=0$, c'est terminé. \\
Si $a_0$ est \ndzz, alors
$\IK_\gC(x_0)= x_0 \gC\supseteq  a_0\gC$. Donc $\gC/\IK_\gC(x_0)$ est un quotient  de $\aqo\gC{a_0}$ qui est égal à $\aqo\Axn{a_0}$, lequel est \ddk $\leq n-1$. Ainsi $\gC/\IK_\gC(x_0)$ est de \ddk $\leq n-1$, et la suite $(\xn)$ est singulière dans $\gC/\IK_\gC(x_0)$.
\end{proof}
%

%%%%%%%%%%%%%%%%%%%%%%%%%%%%%%%%%%%%%%%%%%%%%%%%%%%%%%%%%%%%%%%%%%%%%%%%%%%
%: subsec{Dimension valuative d'un anneau de \pols}
\subsec{Dimension valuative d'un anneau de \pols}

%:HHH petites modifs jusqu'a En \clama on a le résultat suivant: 

Le but de ce paragraphe est de démontrer l'\egtz: 

\snic{\fbox{$\Vdim\AXn = n + \Vdim\gA $}, }

%\sni
pour tout $n\geq 1$. On en déduit la même
\egt pour les dimensions de Krull dans le cas d'un \anarz.

 \underline{Par \dfnz}, cette \egt de dimensions signifie  l'\eqvc suivante\begin{equation}
\label{eqVdimAXn}
\fbox{$\forall k\geq -1,\;\;\Vdim\gA\leq k\iff \Vdim\AXn \leq  n +   k$}.
\end{equation}

Ainsi la première \egt encadrée ne colle pas vraiment pour l'anneau trivial
(il faudrait dire que la dimension de l'anneau trivial est $-\infty$ plutôt \hbox{que $-1$}).

\emph{Remarque préliminaire.} 
\'Etant donné que $\Vdim\gA=\Vdim\Amin$ par \dfnz, et que
$\Amin[\Xn]\simeq(\AXn)\qim$ (fait \ref{fact2Amin}), 
il suffit de traiter le cas où $\gA$
est \qiz, et par la machinerie \lgbe \elr des anneaux \qisz, il suffit de traiter le cas intègre.
Dans la suite du paragraphe, nous utiliserons donc parfois la phrase
salvatrice \gui{on peut \spdg supposer l'anneau intègre}, ou parfois,
si nous voulons donner une explication sur le fonctionnement de la machinerie \lgbe \elrz,  \gui{on peut \spdg supposer l'anneau \qiz}. 
\eoe

%:     Fact{fact}
\begin{fact}\label{factVdimAXnfacile}
Dans \pref{eqVdimAXn}, l'implication réciproque (de droite à gauche) est correcte. 
\end{fact}
%--------- fin fact ---------------------------------------------- 
%
\begin{proof}
On suppose \spdg $\gA$ intègre. On note $[\uX]=[\Xn]$. 
On suppose $\Vdim\AuX\leq n+k$. Soit $\gB=\gA[\yr]$,
avec $y_i\in\Frac\gA$ pour $i\in\lrbn$. 
 On veut démontrer que $\Kdim \gB \leq k$.
\\
Or $\BuX=\AuX[\yr]$ avec les $y_i$ dans $\Frac(\AuX)$.\\
Donc $\Kdim\BuX\leq n+k$, et par le lemme \ref{lemKdimAxn}, $\Kdim\gB\leq k$.
\end{proof}

On étudie maintenant l'implication directe difficile dans  (\ref{eqVdimAXn}).
En \clama on a le résultat suivant: \\
$(*)$ \emph{la dimension valuative d'un anneau intègre $\gA$ est aussi le
maximum des dimensions des anneaux de valuation contenant $\gA$ et contenus
dans son corps de fractions.}

Cette affirmation $(*)$ n'est plus vraie en \gnl d'un point de vue \cof (par manque d'\advsz), mais elle est une conséquence directe (en \clamaz) du corolaire
\ref{prop4ValDim}, qui est donc une version \cov
de $(*)$.

%%:     Proposition{prop2ValDim}
%\begin{proposition}\label{prop2ValDim}
%Soient $\gA\subseteq\gB$, $\fa$ un \id de $\gA$ et $\alpha,\beta\in\gB$
%tels que $\alpha\beta=1$. Supposons que $1\in\fa[\alpha]$ et $1\in\fa[\beta]$,
%alors $1\in \fa$.
%\end{proposition}
%%
%\begin{proof}
%On écrit les hypothèses sous forme $P(\alpha)=0$ et $Q(\beta)=0$. Si $Q$
%est de degré formel $m$ on considère $R(T)=T^mQ(1/T)$.
%En écrivant $\alpha^mQ(\beta)=0$ on obtient $R(\alpha)=0$.
%On vérifie que: le \coe constant de $P$ et le \coe dominant de
%$R$ sont dans $1+\fa$, tous les autres \coes étant dans $\fa$.
%Alors le résultant $\Res(R,m,P,n)$  de $P$ et $R$ est nul (car
%la matrice de Sylvester correspondante annule le vecteur colonne
%$\tra{[{\,\,\ldots\,\alpha^2\,\alpha\,1}}]$)
%et il est de la forme $1+a$ avec $a\in\fa$.
%\end{proof}

%:     lem2Valdim
\begin {lemma}\label{lem2ValDim}
Soient $x_0$, $x_1$, \ldots, $x_\ell$, $u$, $v$, $\alpha$ des indéterminées sur un
anneau $\gA$, $P_0(\alpha)$, \ldots, $P_\ell(\alpha) \in \gA[\alpha]$ et
$Q_0(\alpha^{-1})$, \ldots, $Q_\ell(\alpha^{-1}) \in \gA[\alpha^{-1}]$.  Pour des
$m_i$, $n_i \in \NN$, on définit $P = P(\alpha)$ et $Q = Q(\alpha^{-1})$ comme suit:

\snac{\arraycolsep2pt \begin {array} {rcl}
P &=&
x_0^{m_0}(x_1^{m_1}(\cdots (x_\ell^{m_\ell}
  (u + P_\ell(\alpha)x_\ell) + \cdots) + P_1(\alpha)x_1) + P_0(\alpha)x_0),
\\[1mm]
Q &=&
x_0^{n_0}(x_1^{n_1}(\cdots (x_\ell^{n_\ell}
  (v + Q_\ell(\alpha^{-1})x_\ell) + \cdots) + Q_1(\alpha^{-1})x_1) + Q_0(\alpha^{-1})x_0).
\end {array}
}

%\sni
Si $P$ est de degré formel $p$ (en $\alpha$), $Q$ de degré formel $q$ (en
$\alpha^{-1}$), on considère le résultant:
$$\preskip-.4em \postskip.2em
R = \Res_{\alpha}(\alpha^q Q,q, P,p )  \in \gA[x_0, \ldots, x_\ell, u, v].
$$
Alors, en posant $r_i = qm_i + pn_i$ et $w = u^q v^p$, $R$ est de la forme
$$\preskip.2em \postskip.4em \mathrigid 2mu 
R = x_0^{r_0}(x_1^{r_1}(\cdots (x_\ell^{r_\ell} (w + a_\ell x_\ell) + \cdots) + a_1x_1) + a_0x_0)
\quad \hbox {avec $a_i \in \gA[\ux,u,v]$}. 
$$
\end {lemma}
\begin {proof}
On note $\Res_{\alpha,q,p}(U,V)$
à la place de $\Res_{\alpha}(U,q, V,p )$,
 on suppose $n = 1$ et l'on note $x = x_0$, $y =
x_1$, de sorte que $P = x^{m_0} S$, $\alpha^q Q = x^{n_0} T$, avec:

\snic{S = y^{m_1} (u + P_1(\alpha) y) + P_0(\alpha)x, \quad
T = y^{n_1} (v\alpha^q + T_1(\alpha) y) + T_0(\alpha)x.}

%\sni
On obtient $R = x^{r_0} \Res_{\alpha,q,p}(T,S)$,  $r_0 = qm_0 +
pn_0$. En faisant $x := 0$ on~a:

\snic{\arraycolsep2pt
\begin {array} {rl}
\Res_{\alpha,q,p}(T,S)_{x := 0} &= \Res_{\alpha,q,p}(T_{x:=0} ,S_{x := 0})
\\
&=
\Res_{\alpha,q,p}(y^{n_1} (v\alpha^q + T_1(\alpha)y), y^{m_1} (u + P_1(\alpha)y)\big)
\\
&=
y^{r_1} \Res_{\alpha,q,p}(v\alpha^q + T_1(\alpha)y, u + P_1(\alpha) y),
\\
\end {array}
}

%\sni
avec $r_1 = qm_1 + pn_1$. En faisant $y := 0$ on~a:

\snic{
\Res_{\alpha,q,p}(v\alpha^q + T_1(\alpha)y, u + P_1(\alpha) y)_{y:=0} =
\Res_{\alpha,q,p}(v\alpha^q, u)  = u^q v^p,}

%\sni
ce qui donne le résultat annoncé.
\end {proof}

%:     proposition{prop3ValDim}
\begin{proposition}\label{prop3ValDim}
Soient $\gA\subseteq\gB$, $(\ux)=(x_0,\ldots,x_n)$ une suite dans $\gA$ et~$\alpha_0$,~$\beta_0$ dans $\gB$ tels que $\alpha_0\beta_0=1$. Supposons que la suite soit singulière dans~$\gA[\alpha_0]$ et~$\gA[\beta_0]$, alors elle est singulière dans~$\gA$.
\end{proposition}
\begin{proof}
On applique le lemme précédent en spécialisant $u$ et $v$ en $1$. Les \polsz~$P(\alpha)$ et $\alpha^q Q(\alpha^{-1})$ ayant la racine commune $\alpha_0$ dans $\gB$, leur résultant est nul
 (lemme~\ref{lem0Resultant}).
\end{proof}
%
%:     corollary{prop4ValDim}
\begin{corollary}\label{prop4ValDim}
Soient $a$ et $b$ des \elts \ndzs d'un anneau \qi $\gA$.
Alors $\Vdim\gA=\sup\big(\Vdim\gA[\fraC a b],\Vdim\gA[\fraC b a]\big)$.
\end{corollary}
\begin{proof}
Les in\egts $\Vdim\gA[\fraC a b]\leq \Vdim\gA$ et $\Vdim\gA[\fraC b a]\leq \Vdim\gA$ résultent de la \dfn de la dimension valuative.
\\
Supposons enfin que $\Vdim\gA[\fraC a b]\leq n$ et $\Vdim\gA[\fraC b a]\leq n$ pour un $n\in\NN$. Soit $(\xzn)$ une suite dans $\gA$. Elle est singulière 
dans $\Vdim\gA[\fraC a b]$ et~$\Vdim\gA[\fraC b a]$, donc elle est singulière dans $\gA$ par la proposition \ref{prop3ValDim}.
\end{proof}
% 

%:     Proposition{th0ValDim}
\begin{proposition}\label{th0ValDim}
Pour tout anneau $\gA$  et tout $n\geq 1$,
on a
$$\preskip.2em \postskip.4em
\Vdim\gA[\Xn] \leq n + \Vdim\gA.$$
\end{proposition}
\begin{proof}
On doit montrer que si $\Vdim\gA\leq k$ alors $\Vdim\gA[\Xm] \leq k+m$.
D'après le fait \ref{fact2Amin}, il suffit de traiter le cas % $m=1$ et 
$\gA$ \qiz. 
\\
On suppose d'abord $\gA$ intègre. On reprend la preuve du
 \thref{th1Valdim} et l'on utilise la méthode dynamique. 
 Chaque fois que l'on a une disjonction du type \gui{$a$ divise $b$ ou 
 $b$ divise $a$} on introduit les anneaux $\gC[ \fraC a b]$ et 
 $\gC[\fraC b a]$, où~$\gC$ est l'anneau \gui{en cours}.
\`A chaque feuille de l'arbre ainsi construit on a un anneau
$\gA[u_1,\ldots,u_\ell]\subseteq\Frac\gA$ dans lequel la suite considérée est singulière. On conclut par la proposition \ref{prop3ValDim} que la suite
est singulière dans l'anneau~$\gA$. 
 \\
Dans le cas où $\gA$ est \qi on peut faire appel à la machinerie \lgbe \elr
des anneaux \qisz. On peut aussi raisonner plus directement: la donnée de~$a$ et~$b$
produit la décomposition de \gui{l'anneau en cours} $\gC$ en un produit de 4 composantes. Dans trois d'entre elles, $a$ ou $b$ est nul et tout est facile. Dans la quatrième, $a$ et $b$ sont \ndzs et l'on est ramené au cas intègre.\imlgz
\end{proof}

Comme corolaires on obtient les  \thos suivants.

%:     Theorem{thValDim}
\begin{theorem}\label{thValDim} 
Pour un anneau $\gA$, on a les \eqvcs suivantes.
\begin{enumerate}
\item Si $n\geq 1$ et $k\geq -1$, alors
$$\preskip.2em \postskip.4em
\Vdim\gA\leq k\iff\Vdim\AXn\leq n+k. 
$$
Autrement dit, $\Vdim\AXn=n+\Vdim\gA$.
\item Si $n\geq 0$, alors
$$\preskip.2em \postskip.4em
\Vdim\gA\leq n\iff \Kdim \AXn \leq 2n.
$$
Dans le cas où $\gA$ est \qiz, c'est aussi \eqv à:\\
pour tous $x_1$, \dots, $x_n$ dans $\Frac\gA$, on a $\Kdim\Axn\leq n$.
\end{enumerate}
\end{theorem}
\begin{proof}
\emph{1.} Démontré dans le fait \ref{factVdimAXnfacile} et la proposition \ref{th0ValDim}.
\\
\emph{2.} Le cas $n=0$ a déjà été noté. Voyons le cas $n\geq 1$.
L'implication directe résulte du point \emph{1}  parce que $\Kdim\AXn\leq\Vdim\AXn$. L'implication
réciproque est donnée (dans le cas intègre, mais ce n'est pas restrictif) dans la proposition \ref{propVdimKdim}. 
\end{proof}
%

%:     Theorem{corthValDim}
\begin{theorem}\label{corthValDim}~
\begin{enumerate}
\item Si $\gA$ est un \anar de \ddk finie on a
\[
\Vdim \gA[\Xn]= \Kdim \gA[\Xn]\leq n+\Kdim \gA.
\] 
avec \egt si $\gA$ est non trivial.
\item $\Vdim \ZZ[\Xn]= \Kdim \ZZ[\Xn]=1+n$.
\item Tout anneau engendré par $n$ \elts est de
dimension valuative (donc de \ddkz) $\leq1+n$.
\item \label{i5corthValDim} Soit $\gA$ un anneau \qi engendré par $n$ \elts et
$\gB$ un anneau intermédiaire entre $\gA$ et $\Frac\gA$. Alors $\Vdim\gB\leq1+n$.
\end{enumerate}
\end{theorem}
\begin{proof}
Le point \emph{1} résulte du \tho plus \gnl \ref{cor0thValDim}
et le point  \emph{2} est un cas particulier.

\emph{3.} L'anneau $\gA$ est un quotient de $\ZZXn$, donc 
$\gA[Y_1,\dots,Y_{n+1}]$ est un quotient de $\ZZXn[Y_1,\dots,Y_{n+1}]$
qui est de \ddkz~$2n+2$ par le point \emph{2.} Donc $\Vdim\gA\leq n+1$
par le point \emph{2} du \thref{thValDim}.

\emph{4.} Conséquence du point \emph{3} puisque $\Vdim\gA\leq n+1$.
\end{proof}
%

%:     Theorem{cor0thValDim}
\begin{theorem} \label{cor0thValDim} 
 Pour un anneau $\gA$ de \ddk $\leq n$ ($n\geq 1$) \propeq
\begin{enumerate}
\item $\Vdim\gA=\Kdim\gA$.
\item Pour tout $k\geq 1$, $\Kdim(\AXk)\leq k+\Kdim\gA$%
%, avec \egt si $\gA$ est non trivial
.
\item $\Kdim(\AXn)\leq n+\Kdim\gA$%
%, avec \egt si $\gA$ est non trivial
%
\end{enumerate}
Si $\gA$ est non trivial, on peut remplacer $\leq $ par $=$ dans les points 2. et 3.\\
Si $\Vdim\gA=\Kdim\gA$, pour tout $k\geq 1$, on a l'\egt

\snic{\Kdim(\AXk)=\Vdim(\AXk).}
\end{theorem}
\begin{proof}
Notons que l'on ne suppose pas connue de manière exacte la \ddk de $\gA$.
\\
\emph{1} $\Rightarrow$ \emph{2.} On fixe un $k\geq 1$ et l'on doit montrer que pour tout $m\geq -1$, on a~$\Kdim\gA\leq m\Rightarrow\Kdim(\AXk)\leq m+k$.
\\
On a $\Vdim(\gA)\leq m$, donc $\Vdim(\AXk)\leq m+k$ d'après \ref{th0ValDim},
\hbox{donc $\Kdim(\AXk)\leq m+k$} car on a toujours $\Kdim\gB\leq \Vdim\gB$.

\emph{2} $\Rightarrow$ \emph{3.} C'est le cas particulier où $k=n$.
 
\emph{3} $\Rightarrow$ \emph{1.} On suppose $\Kdim\gA\leq m$  
et on doit montrer $\Vdim\gA\leq m$. \Spdg $0\leq m\leq n$. Si $m=n$ on conclut par le point \emph{2} du \thref{thValDim}.
Si $n= m+r$, on a $\Kdim(\AXn)\leq n+m$ par hypothèse. Comme $(X_{m+1},\dots,X_n)$ est singulière de longueur $r$,
 le point \emph{3} de la proposition \ref{lemRegsing} nous donne $\Kdim(\AXm)\leq n+m-r= 2m$ et on conclut par le point \emph{2} du \thref{thValDim}.

Le dernier point   est laissé \alecz.
\end{proof}
%

%%%%%%%%%%%%%%%%%%%%%%%%%%%%%%%%%%%%%%%%%%%%%%%%%%%%%%%%%%%%%%%%%%%%%%%%%%%
\section{Lying over, Going up et Going down}
\label{secGoingUp}

Nous sommes intéressés dans cette section pour comprendre en termes \cofs
certaines propriétés des anneaux commutatifs
et de leurs morphismes qui sont introduites en \clama  via les
notions  de spectre de Zariski ou de morphisme spectral entre spectres de Zariski (donné par un \homo d'anneaux).

En faisant fonctionner nos \dfns \covs
nous espérons obtenir des versions \covs de nombreux \thos  de \clamaz,
réellement utilisables en pratique.
En fait, c'est ce qui se passera de manière systématique dans cette section et dans les
chapitres suivants.

%-% ENTRE NOUS
\entrenous{La section manque d'un ou deux exemples d'applications.
}
%-% Fin ENTRENOUS

%:H2018  subsection* remplacés par subsec
%:subsec{Relèvement des \ideps (lying over)} 
\subsec{Relèvement des \ideps (lying over)}

En \clama on dit qu'un \homo $\alpha:\gT\to\gV$ de \trdis \gui{possède la \prt de relèvement des \idepsz} lorsque l'\homo dual $\Spec\alpha:\Spec\gV\to\Spec\gT$ est surjectif, autrement dit lorsque tout \idep de $\Spec\gT$
est l'image réciproque d'un \idep de $\Spec\gV$. Pour abréger on dit
aussi que le morphisme est \gui{lying over}.
Nous allons donner une \dfn \cot pertinente sans utiliser l'\homo dual. Pour l'\eqvc en \clama avec la \dfn via les spectres,
voir l'exercice~\ref{exoLYOV}.

%:H2018 reecriture du point 1
%:     Definition{defiLYO}
\begin{definition}\label{defiLYO}~
\begin{enumerate}
\item Un \homo $\alpha:\gT\to\gV$ de \trdis est dit \emph{lying over} lorsqu'il est injectif. Il revient au même de dire que $\alpha$ réfléchit les inégalités: pour $a$, $b\in\gT$,
$\alpha(a)\leq\alpha(b)$ implique $a\leq b$.
\item Un \homo $\varphi:\gA\to\gB$ d'anneaux
commutatifs est dit \emph{lying over} lorsque
l'\homo  $\Zar\varphi:\Zar\gA\to\Zar\gB$ est injectif.
\end{enumerate}
\end{definition}
\index{lying over!morphisme ---}

\rdb\label{remLY}
\rem On a aussi les formulations \eqves suivantes pour les morphismes lying over.
\begin{itemize}
\item Pour les \trdisz:
\begin{itemize}
\item Pour tout $b\in\gT$, $\alpha^{-1}(\dar \alpha(b)\big)=\dar b$.
\item Pour tout \id $\fa$ de $\gT$,  
$\alpha^{-1}\big(\cI_{\gV}(\alpha(\fa))\big)=\fa$.
\end{itemize}
\item Pour les anneaux commutatifs:
\begin{itemize}
\item Pour tout $x\in\gA$ et tout \itf $\fa$ de $\gA$ on a
l'implication

\centerline{$\varphi(x)\in\varphi(\fa)\gB\
\Longrightarrow\ x\in\DA(\fa).$}
\item  Pour tout \itf $\fa$ de $\gA$ on a
$\varphi^{-1}(\gen {\varphi(\fa)}) \subseteq \DA({\fa})$.
\item   Pour tout \id $\fa$ de $\gA$ on a
$\varphi^{-1}\big(\DB(\gen{\varphi(\fa)})\big)=\DA({\fa})$. \eoe
\end{itemize}
\end{itemize}

%:     Fact{factLOInt}
\begin{fact}\label{factLOInt}
Soit $\gB\supseteq\gA$ une extension d'anneaux. Si $\gB$ est entière ou \fpte  sur $\gA$,
le morphisme d'inclusion
$\gA\to\gB$ est lying over.
\end{fact}
%%%%%%%%%%%%%
%
\begin{proof}
Le premier cas   est une simple reformulation du lemme
\ref{lemLingOver} (lying over). Dans le deuxième cas, 
pour tout \itf $\fa$ de $\gA$, on~a~$\fa\gB\cap\gA=\fa$.
\end{proof}
%

%-% ENTRE NOUS
\rem {Il n'est pas vrai que $\varphi:\gA\to\gB$ soit lying over si
$\gB$ est entière sur~$\gA$ mais avec $\varphi$ non injectif. En effet,
prenons le cas où~$\gB$ est un quotient de $\gA$,
donc entière sur $\gA$. Si~$\fp$ est un \idep
de~$\gA$ strictement contenu dans $\Ker\varphi$ il n'est pas l'image réciproque
d'un \idep de $\gB$. Ou bien en termes constructifs, si deux \itfs de $\gA$
sont contenus dans~$\Ker\varphi$ il n'y a aucune raison pour que leurs
radicaux soient égaux, sauf si $\Ker\varphi\subseteq\DA(0)$.
%
%Il serait intéressant d'avoir une \carn exacte en termes ``simples''
%des morphismes lying over pour les anneaux commutatifs.
%
%\`A nilradical près, ce sont des inclusions, mais des inclusions
%d'un type particulier.
}
%-% Fin ENTRENOUS

%:subsec{Montée (going up)} %%%%%%%%%%%%
\subsec{Montée (going up)}
%:H2018 paragraphe allongé avant la définition
En \clama on dit qu'un \homo $\alpha:\gT\to\gV$ de \trdis \emph{possède la \prt de montée pour les chaînes d'\idepsz}, ou plus simplement qu'il est \emph{going up} lorsque  la \prt suivante est satisfaite.

\emph{Si $\fq\in\Spec\gV$ et $\alpha^{-1}(\fq)=\fp$, toute chaîne $\fp=\fp_1\subseteq\cdots\subseteq\fp_n$ d'\ideps de $\gT$ 
est l'image réciproque d'une  chaîne $\fq=\fq_1\subseteq\cdots\subseteq\fq_n$ d'\ideps de $\gV$.
}

\smallskip Naturellement on peut se limiter au cas $n=2$.
On voit alors que la \dfn peut se relire comme suit.

\emph{Si $\fq\in\Spec\gV$ et $(\Spec\alpha)(\fq)=\fp$, et si on note $$\alpha':(\gV/(\fq=0)\to \gT/(\fp=0)$$ le morphisme obtenu par factorisation,  
alors le morphisme $$\Spec\alpha':\Spec(\gV/(\fq=0))\to\Spec(\gT/(\fp=0))$$ est surjectif.} 

Ceci nous ramène au cas d'un morphisme lying over.

Les mêmes \dfns sont utilisées pour les morphismes d'anneaux commutatifs.

Voici des \dfns \covs en termes de \trdis et d'anneaux commutatifs.
Pour l'équivalence avec la \dfn classique en \clama voir l'exercice \ref{exodefiGoingup2}.
%:     Definition{defiGoingup}
\begin{definition}\label{defiGoingup}~
\begin{enumerate}
\item Un \homo $\alpha:\gT\to\gV$ de \trdis est dit \emph{going up}
lorsque   pour
tous $a,c\in\gT$ et $y\in\gV$ on a

\snic{\alpha(a)\leq\alpha(c)\vu y \quad\Longrightarrow\quad\exists x\in\gT\; (a\leq c \vu x\et \alpha(x)\leq y).}

\vspace{1mm}

\item Un \homo $\varphi:\gA\to\gB$ d'anneaux
commutatifs est dit \emph{going up} lorsque
l'\homo  $\Zar\varphi:\Zar\gA\to\Zar\gB$ est going~up.
\end{enumerate}
\end{definition}
\index{going up!morphisme ---}

\rems
1)
Pour le point \emph{1}, si $\fa=\alpha^{-1}(0_\gV)$ et $\gT_1=\gT\sur{(\fa=0)}$,
alors $\alpha$ est going up \ssi  $\alpha_1:\gT_1\to\gV$ est going~up.
\\
Pour le point \emph{2}, si $\gT=\ZarA$, alors $\gT_1\simeq\Zar(\varphi(\gA)\big)$. On en déduit en posant $\gA_1=\varphi(\gA)$,
que $\varphi$ est  going up \ssi  $\varphi_1:\gA_1\to\gB$ est going~up.

 2)
Pour les \trdisz, si $\alpha^{-1}(0)=0$ et si $\alpha$ est going up, alors il est lying over. Pour les anneaux commutatifs, si $\Ker\varphi\subseteq\DA(0)$
et si $\varphi$ est going up,  alors il est lying over.
\eoe

%:     Proposition{prop1Gup}
\begin{proposition}\label{prop1Gup}
Si $\gB$ est  une \Alg entière, le morphisme
$\gA\to\gB$ est going up.
\end{proposition}
\begin{proof}
D'après la remarque précédente on peut supposer $\gA\subseteq\gB$.
On sait alors que l'\homo est lying over, i.e. que
$\ZarA\to\Zar\gB$ est injectif, donc on identifie
$\ZarA$ à un sous-treillis de $\Zar\gB$.
On doit montrer qu'étant donnés $a_1$, \dots, $a_n$, $c_1$, \dots, $c_q$ dans $\gA$ et $y_1$, \dots, $y_p$ dans $\gB$
vérifiant

\snic{\DB(\ua)\leq\DB(\uc)\vu\DB(\uy),}

%\sni
on peut trouver une suite $(\ux) $ dans $\gA$ telle que

\snic{\DA(\ua)\leq\DA(\uc)\vu\DA(\ux)\;$ et  $\;
\DB(\ux)\leq \DB(\uy).}

%\sni
Soit $\fb=\DB(\uy)$,  $\fa=\fb\cap\gA$, $\gB_1=\gB\sur{\fb}$
et $\gA_1=\gA\sur{\fa}$. On
considère l'extension entière $\gB_1\supseteq\gA_1$.
L'hypothèse est maintenant
que ${\rD_{\gB_1\!}(\ua)\leq\rD_{\gB_1\!}(\uc).}$
\\
Par le lying over on sait que cela implique que
${\rD_{\gA_1\!}(\ua)\leq\rD_{\gA_1\!}(\uc).}$
Cela signifie que pour chaque $i\in\lrbn$ on a un $x_i\in\fa$
tel que $\DA(a_i)\leq\DA(\uc)\vu\DA(x_i)$. On a donc réalisé
le but recherché avec $(\ux)=(\xn).$
\end{proof}
%

%:subsec{Descente (going down)} %%%%%%%%%%%%
\subsec{Descente (going down)}
%:H2018 paragraphe allongé avant la définition
En \clama on dit qu'un \homo $\alpha:\gT\to\gV$ de \trdis \gui{possède la \prt de descente pour les chaînes d'\idepsz} 
ou plus simplement qu'il est \gui{going down} lorsque 
le morphisme opposé $\alpha\eci:\gT\eci\to\gV\eci$ est going up.
Autrement dit lorsque la \prt suivante est satisfaite.

\emph{Si $\fq\in\Spec\gV$ et $\alpha^{-1}(\fq)=\fp$, toute chaîne $\fp_1\subseteq\cdots\subseteq\fp_n=\fp$ d'\ideps de $\Spec\gT$
est l'image réciproque d'une  chaîne $\fq_1\subseteq\cdots\subseteq\fq_n=\fq$ d'\ideps de $\Spec\gV$.
}

On pourrait dire tout aussi bien ceci, en se rappelant qu'un filtre premier de $\gT$ est un \idep de $\gT\eci$.

\emph{Si $\ff\in\Spec\gV\eci$ et $\alpha^{-1}(\ff)=\ffg$, toute chaîne $\ff=\ff_1\subseteq\cdots\subseteq\ff_n$ de filtres premiers de $\gT$
est l'image réciproque d'une  chaîne $\ffg=\ffg_1\subseteq\cdots\subseteq\ffg_n$ de filtres premiers de $\gV$.
}

\smallskip Naturellement on peut se limiter au cas $n=2$. Et la \dfn \cov que nous proposons est la notion opposée à celle du going up: elle est obtenue en renversant les relations d'ordre. 

%:     Definition{defiGoingdown}
\begin{definition}\label{defiGoingdown}~
\begin{enumerate}
\item Un \homo $\alpha:\gT\to\gV$ de \trdis est dit \emph{going down} lorsque le même \homo pour les treillis
opposés $\gT\eci$ et $\gV\eci$ est going up.
Autrement dit pour
tous $a,c\in\gT$ et $y\in\gV$ on a

\snic{\alpha(a)\geq\alpha(c)\vi y \quad\Longrightarrow\quad\exists x\in\gT\; (a\geq c \vi x \et \alpha(x)\geq y).}

\vspace{1mm}
\item Un \homo $\varphi:\gA\to\gB$ d'anneaux
commutatifs est dit \emph{going down} lorsque
l'\homo  $\Zar\varphi:\Zar\gA\to\Zar\gB$ est going down.
\end{enumerate}
\end{definition}
\index{going down!morphisme ---}

\rems 1) La \dfn \emph{1} revient à dire que l'image par $\alpha$ de l'\id transporteur
 $(a:c)_\gT$
engendre l'\id $\big(\alpha(a):\alpha(c)\big)_{\gV}$. Si donc les \trdis sont des \agHs
cela signifie que l'\homo de treillis est aussi un \homo d'\agHsz.

On a aussi les mêmes remarques que pour le going up. 

2) Si $\ff=\alpha^{-1}(1_\gV)$ et
$\gT_2=\gT\sur{(\ff=1)}$, alors
$\alpha$ est going down \ssi  $\alpha_2:\gT_2\to\gV$ est going down.
\\
Ceci donne pour les anneaux commutatifs: si $S=\varphi^{-1}(\gB\eti)$ et $\gA_2=\gA_S$, alors $\varphi$ est  going down \ssi  $\varphi_2:\gA_2\to\gB$ est going down.

3)
Pour les \trdisz, si $\alpha^{-1}(1)=1$ et $\alpha$ est going down, alors il est lying over. Pour les anneaux commutatifs, si $\varphi^{-1}(\gB\eti)\subseteq\Ati$
et $\varphi$ est going down,  alors il est lying over.
\eoe

%:subsec{Conséquences pour la dimension de Krull}
\subsec{Conséquences pour la dimension de Krull}

%:     Proposition{propGupLY}
\begin{theorem}\label{propGupLY}
Si un \homo $\alpha:X\to Y$ (de \trdis ou d'anneaux commutatifs) est lying over et going up,
ou s'il est  lying over et going down, on a $\Kdim X\leq\Kdim Y$.
\end{theorem}
%%%%%%%%%%%%%%%%%%%%%%%%%%%%%%%%%%%%%%%%%
\rem C'est par exemple le cas lorsque
l'anneau $\gB$ est une extension entière de $\gA$. On retrouve ainsi la
proposition~\ref{propDKEXENT}. Pour les extensions plates, voir la proposition~\ref{prop1Gdown}. 
\eoe
%%%%%%%%%%%%%%%%%%%%%%%%%%%%%%%%%%%%%%%%%
\begin{proof}
Il suffit de traiter le cas  going up avec les treillis.\\
On suppose $\Kdim Y\leq n$
et
l'on considère une suite $(\azn)$ dans $X$. \linebreak 
On a dans $Y$ une suite $(y_0,\ldots,y_n)$  \cop
de $\alpha(\ua)$:

\snic{\alpha(a_0)\vi y_ 0\leq0,\;\ldots,\;\alpha(a_{n})\vi y_n\leq\alpha(a_{n-1})\vu y_{n-1}, \; 1\leq\alpha(a_n)\vu y_n.}

%\sni
On va construire une suite  $(\xzn)$  \cop
de $(\ua)$ dans $X$. \`A l'étage $n$, par going up, il existe $x_ n\in X$ tel que

\snic{1\leq a_n\vu x_n$ et $\alpha(x_n)\leq y_n.}

%\sni
 Ceci donne à l'étage $n-1$ l'in\egtz:
$\alpha(a_{n}\vi x_n)\leq \alpha(a_{n-1})\vu y_{n-1}.$
\\
  Par going up
il existe $x_ {n-1}\in X$ tel que

\snic{a_{n}\vi x_n\leq a_{n-1}\vu x_{n-1}\hbox{  et  }\alpha(x_{n-1})\leq y_{n-1}.}

%\sni
  On continue de la même manière jusqu'à l'étage $0$, où cette fois-ci il faut utiliser le lying over.
\end{proof}
%

%:     Lemma{lem1Gdown}
\begin{lemma}\label{lem1Gdown}
Pour qu'un \homo d'anneaux $\varphi:\gA\to\gB$ soit going down il faut et suffit que pour tous $c$, $a_1$, \ldots, $a_q\in\gA$ et $y\in\gB$ tels que~$\varphi(c)y\in\DB(\varphi(\ua)\big)$,
il existe des \elts $x_1$, \dots, $x_m\in\gA$ tels que:

\snic{\DA(c)\vi\DA(\ux)\leq\DA(\ua)\;\et\;\DB(y)\leq\DB(\varphi(\ux)\big).}
\end{lemma}
\begin{proof}
Nous avons remplacé dans la \dfn un \elt arbitraire $\DA(\uc)$ de~$\ZarA$ et un
\elt arbitraire $\DB(\uy)$ de $\Zar\gB$ par des \gtrsz~$\DA(c)$ et~$\DB(y)$.
Comme les \gtrs  $\DA(c)$ (resp. $\DB(y)$) engendrent $\ZarA$ (resp. $\Zar\gB$)
par sups finis, les règles de \dit impliquent que la restriction à ces \gtrs
est suffisante (calculs laissés \alecz).
\end{proof}
%

%:     Proposition{prop1Gdown}
\begin{proposition}\label{prop1Gdown}
Un \homo $\varphi:\gA\to\gB$ d'anneaux
commutatifs est going down dans les deux cas suivants.
\begin{enumerate}
\item $\gB$ est une \Alg plate.
\item $\gB\supseteq \gA$ est intègre et entier sur $\gA$, et $\gA$ est \iclz.
\end{enumerate}
\end{proposition}
\begin{proof} On se place dans les hypothèses du lemme \ref{lem1Gdown}, avec une \egt dans~$\gB$:
$$\preskip.4em \postskip.4em 
\varphi(c)^\ell y^\ell +\som_{i=1}^qb_i\varphi(a_i)=0\eqno(*) 
$$

\emph{1.}
On regarde $(*)$ comme une \rde $\gB$-\lin entre les \elts $c^\ell$, $a_1$, \ldots, $a_q$.
On exprime qu'elle est une combinaison $\gB$-\lin de \rdes $\gA$-\linsz.
% entre $a^\ell,\cq$.
\\
 Ces relations s'écrivent $x_jc^\ell+\sum_{i=1}^q u_{j,i}a_i=0$
pour $j\in\lrbm$, avec les~$x_j$ et les~$u_{j,i}$ dans $\gA$.
D'où $\DA(cx_j)\leq\DA(\ua)$, et 
$\DA(c)\vi\DA(\ux)\leq\DA(\ua)$. 
Enfin,~$y^\ell$
est une combinaison $\gB$-\lin des $\varphi(x_j)$, d'où $\DB(y)\leq\DB(\varphi(\ux)\big)$.

 \emph{2.}
D'après $(*)$, $(cy)^\ell \in \gen{\ua} \,\gB$. Par le lying over
\ref{lemLingOver2}, $(cy)^\ell$, et a fortiori~$cy$, est entier sur
$\gen{\ua}_\gA$. On écrit une \rdi pour~$cy$ sur l'\id $\gen{\ua}_\gA$ sous
la forme $f(cy)=0$ avec

\snic{f(X)=X^k+\sum_{j=1}^k\mu_jX^{k-j}
\qquad \hbox {où }\mu_j\in\gen{\ua}_{\!\gA}^{\,j}.}

%\sni
Par ailleurs, $y$ annule un \polu $g(X)\in\AX$. Considérons
dans~$(\Frac\gA)[X]$ le pgcd \mon $h(X)=X^{m}+x_1X^{m-1}+\cdots+x_m$ des deux
\pols $f(cX)$ et $g(X)$. Puisque $\gA$ est \iclz, le \tho de \KRO dit que
$x_j\in\gA$, et l'\egt $h(y)=0$ donne $y\in\DB(\ux)$. \\
Il reste à voir que
$cx_j\in\DA(\ua)$ pour $j\in\lrbm$. En remplaçant formellement~$X$ par~$Y/c$, 
on obtient que le
\pol 

\snic{h_c(Y)=Y^{m}+cx_1Y^{m-1}+\cdots+c^mx_m}

%\sni
divise~$f(Y)$ dans
$(\Frac\gA)[Y]$. Le \tho de \KRO (sous la forme du lemme~\ref{lemthKroicl}) nous dit que $cx_j\in\DA(\mu_1,\ldots,\mu_k)$.\\
Enfin, comme~$\DA(\mu_1,\ldots,\mu_k)\leq\DA(\ua)$, on a bien $cx_j\in\DA(\ua)$.
\end{proof}
%

%:subsec{Incomparabilité} 
\subsec{Incomparabilité}
En \clama on dit qu'un \homo $\alpha:\gT\to\gT'$ de \trdis \gui{possède la \prt d'incomparabilité} lorsque les fibres de l'\homo dual
$\Spec\alpha:\Spec\gT'\to\Spec\gT$ sont constituées d'\elts
deux à deux incomparables. %\\
Autrement dit, pour $\fq_1$
et $\fq_2$ dans~$\Spec\gT'$, si $\alpha^{-1}(\fq_1)=\alpha^{-1}(\fq_2)$ et $\fq_1\subseteq\fq_2$, alors $\fq_1=\fq_2$.

La \dfn \cov correspondante est que le morphisme $\gT\to\gT'$ est \zedz.

Nous avons déjà donné la \dfn de la dimension d'un morphisme
dans le cas des anneaux commutatifs. Une \dfn analogue peut être fournie pour les \trdisz, mais nous n'en aurons pas l'usage.

 La principale conséquence de la situation d'incomparabilité
pour un \homo $\varphi:\gA\to\gB$ est le
fait que $\Kdim\gB\leq\Kdim\gA$. Ceci est un cas particulier du \thref{thKdimMor}
avec l'important \thref{cor2thKdimMor}.

%%%%%%%%%%%%%%%%%%%%%%%%%%%%%%%%%%%%%%%%%%%%%%%%%%%%%%%%%%%%%%%%%%%%%%%%%%%
%%%%%%%%%%%%%%%%%%%%%%%                              %%%%%%%%%%%%%%%%%%%%%%
%%%%%%%%%%%%%%%%%%%%%%%               EXOS           %%%%%%%%%%%%%%%%%%%%%%
%%%%%%%%%%%%%%%%%%%%%%%                              %%%%%%%%%%%%%%%%%%%%%%
%%%%%%%%%%%%%%%%%%%%%%%%%%%%%%%%%%%%%%%%%%%%%%%%%%%%%%%%%%%%%%%%%%%%%%%%%%%
%: section Exos
\Exercices{

%--- Exercise{exoKdimLecteur}-------------
\begin{exercise}
\label{exoKdimLecteur}
{\rm  Il est recommandé de faire les \dems non données, esquissées,
laissées \alecz,
etc\ldots
\, On pourra notamment traiter les cas suivants.
\begin{itemize}
\item \label{exoDdk1} Démontrer la proposition \ref{propDdk0}.
\item \label{exobord}
Démontrer ce qui
est affirmé dans les exemples \paref{exlKdim}.
\item \label{exocorfactDDKTRDI} Démontrer le fait \ref{corfactDDKTRDI}.
\item \label{exopropDDKagH}
Démontrer les faits \ref{factDdkTrdiBord} et \ref{propDDKagH}.
\item \label{exopropNoetAgH} 
Vérifier les détails dans la \dem de la proposition
\ref{propNoetAgH}.
\item \label{exolem2qi} 
Démontrer le lemme
\ref{lem2qi} en vous inspirant de la preuve du lemme~\ref{lem2SousZedRed}.

%
%\item \label{exoprop4ValDim} Démontrer le corolaire \ref{prop4ValDim}.
%
%\item \label{exocorthValDim} Vérifier les affirmations
%du \thrf{corthValDim} comme corolaires du \tho \ref{thValDim}.
%
\item \label{exolem1Gdown} Vérifier les détails dans la \dem 
du lemme \ref{lem1Gdown}.
\end{itemize}
}
\end{exercise}
%--- end -exercise-----------------------------------------

%--- Exercise{exoIdFilPrem}-------------
\begin{exercise}
 \label{exoIdFilPrem}
 {\rm  Si $\ff$ est un filtre de l'anneau $\gA$, définissons son \emph{complément}~$\wi{\ff}$
 comme étant $\sotq{x\in\gA}{x\in\ff\Rightarrow0\in\ff}$.
 En particulier, on a toujours $0\in\wi\ff$, même si $0\in\ff$.
 De même, si $\fa$ est un \id de l'anneau $\gA$, définissons son \emph{complément} $\bar{\fa}$
 comme étant $\sotq{x\in\gA}{x\in\fa\Rightarrow1\in\fa}$. Montrer que si $\ff$ est filtre premier son complément est un \idz. Si en outre $\ff$ est détachable,
 alors $\fa$ est un \idep détachable. Montrer aussi les affirmations duales.

} \end{exercise}
%--- end-exercise-----------------------------------------

%--- Exercise{exoDimKX}--------------
\begin{exercise}
\label{exoDimKX}
{\rm  \emph{1.} Si la suite $(\Xn)$ est singulière dans l'anneau $\AXn$,
alors $\gA$ est
trivial.
\\
\emph{2.} Soit $k\in\NN$. Démontrer que si $\AX$ est un anneau  \ddi$k$ alors
$\gA$ est \ddi$k-1$.
Retrouver ainsi le point~\emph{1}.
}
\end{exercise}
%--- end-exercise-----------------------------------------

%--- Exercise{exoDimKX1n}------------
\begin{exercise}
\label{exoDimKX1n}
{\rm  Démontrer que si $\gK$ est un anneau de \ddk exactement égale à $0$ alors,  $\gK[X_1,\ldots ,X_n]$ est de \ddk exactement égale à $n$.
}
\end{exercise}
%--- end-exercise-----------------------------------------

%--- Exercise{exoPartitionUnite}-------------
\begin{exercise}\label{exoPartitionUnite}
 {(Partition de l'unité associée à un recouvrement
ouvert du spectre)}\\
{\rm  
Soit $\gA$ un anneau et $(U_i)_i$ un recouvrement ouvert de $\Spec(\gA)$. 
Montrer en \clama 
qu'il existe une famille $(f_i)_i$ d'\elts de $\gA$ avec $f_i = 0$ sauf
pour un nombre fini d'indices $i$ et

\snic {
(\star) \qquad\qquad
\DA(f_i) \subseteq U_i, \qquad   \sum_i f_i = 1
.}

%\sni
Remarque: ainsi, on remplace tout recouvrement ouvert de
$\Spec(\gA)$ par un \sys fini d'\elts de $\gA$ qui \gui {recouvrent} $\gA$
(puisque leur somme vaut $1$), sans \gui {perdre d'information} puisque
$(\star)$ confirme de nouveau que $(U_i)_i$ est un recouvrement.

}

\end {exercise}
%--- end -exercise-----------------------------------------

%--- Exercise{exoKdimGeom}-------------
\begin{exercise}
\label{exoKdimGeom}
{\rm  Pour une \apf $\gA$ sur un \cdi non trivial
appelons \gui{dimension de \Noe de $\gA$} le nombre de variables
\agqt indépendantes après une mise un position de \iNoez.
\\
\emph{1.} Soit $f\in\gA\supseteq\KYr=\KuY$ ($\gA$ entière sur $\KuY$). 
\\
\emph{1a.} Montrer que l'\id bord de $f$ contient un $g\in\KuY\setminus\so{0}$. 
\\
\emph{1b.} 
En déduire que l'anneau bord de Krull $\gA\sur{\JK_\gA(f)}$ est un quotient
d'une \apf dont la dimension de \Noe  
est $\leq r-1$.
\\
\emph{2.}  
En déduire une \dem directe de l'\egt des dimensions de Krull
et de \Noe des \apfs sur un \cdi non trivial.
}
\end{exercise}
%--- end -exercise-----------------------------------------

%--- Exercise{exolemLocMemeKdim}-------------
\begin{exercise}
\label{exolemLocMemeKdim}
{\rm  \emph{1.} Soient $\gK$  un \cdi non trivial, $\KuX=\KXn$ et~$f\in \KuX\setminus\so{0}$,
alors $\Kdim \KuX[1/f] =n$.
\\
\emph{2.} 
Plus \gnltz, donner une condition suffisante sur le \pol $\delta\in\AuX$ 
pour que l'on ait $\Kdim(\AuX[1/\delta])=\Kdim\AuX$  
(voir la \dem du lemme~\ref{lemLocMemeKdim}).
}
\end{exercise}
%--- end -exercise-----------------------------------------

%--- Exercise{exoPruferKdim}-------------
\begin{exercise}\label{exoPruferKdim} {(\Carn des \adps intègres de dimension $\leq1$)}
\\
{\rm
Soit $\gA$ un anneau \iclz.
\\
\emph{1.} Montrer que si $\Kdim\AX\leq2$, alors $\gA$ est un \adpz, en montrant 
que tout \elt de~$\Frac\gA$ est primitivement  \agq sur $\gA$.
\\
\emph{2.} Montrer que $\gA$ est un \adp \ddi1 \ssi $\Kdim\AX\leq2$.
\\
\emph{3.} Peut-on \gnr à un anneau normal?
}
\end {exercise}
%--- end -exercise-----------------------------------------

%--- Exercise{exoMultiplicativiteIdeauxBords}-------------
\begin{exercise}\label{exoMultiplicativiteIdeauxBords}
{(Une \prt de multiplicativité des \ids bord)}\\
{\rm  
\emph{1.}
Pour $a, b \in \gA$ et deux suites $(\ux)$, $(\uy)$ d'\elts de $\gA$,
montrer que

\snic {
\IK_\gA(\ux, a, \uy)\,\IK_\gA(\ux, b, \uy) \subseteq \IK_\gA(\ux, ab, \uy)
.}

\emph{2.}
En déduire que
$\IK_\gA(a_1b_1, \ldots, a_nb_n)$ 
 contient le produit  $\prod_{\uc} \IK_\gA(\uc)$,
 dans lequel la suite $(\uc) = (c_1, \ldots, c_n)$ parcourt l'ensemble des $2^n$ suites
telles que $c_i = a_i$ \linebreak 
ou $c_i = b_i$ pour chaque $i$.
}

\end {exercise}
%--- end -exercise-----------------------------------------

%--- Exercise{exoInclusionBordLionel}-------------
\begin{exercise}\label{exoInclusionBordLionel}
{(Idéaux bord et relations \agqsz)}
\\
{\rm  
\emph {1.}
On considère l'ordre lexicographique sur $\NN^n$. Soit $\alpha = (\alpha_1,
\ldots, \alpha_n) \in\NN^n$. Vérifier, pour 
$\beta > \alpha$, que

\snic {
\uX^\beta \in \gen {X_1^{1+\alpha_1},\ X_1^{\alpha_1} X_2^{1+\alpha_2},\
X_1^{\alpha_1} X_2^{\alpha_2} X_3^{1+\alpha_3},\ \cdots,\ 
X_1^{\alpha_1} X_2^{\alpha_2} \cdots X_{n-1}^{\alpha_{n-1}}X_n^{1+\alpha_n}}
.}

\emph {2.}
Soit $\gA$ un anneau réduit, $(\ux) = (\xn)$ une suite dans $\gA$ et $P = \sum_\beta a_\beta \uX^\beta$ dans~$\gA[\uX]$, qui annule $\ux$.
\begin {itemize}
\item [\emph {a.}]
Montrer, pour $\alpha\in\NN^n$, que
$a_\alpha \prod_{\beta < \alpha} \Ann(a_\beta) \subseteq \IK(\ux)$.

\item [\emph {b.}]
En déduire:

\snic {
\prod_\beta \IK(a_\beta) \subseteq \IK(\ux) + 
\prod_\beta \Ann(a_\beta)}
\end {itemize}

%\sni
\emph {3.}
Soient une \alg $\gA \to \gB$ avec $\gB$ réduit et $x \in \gB$ primitivement
\agq sur $\gA$: $\sum_{i=0}^d a_i x^i = 0$ avec $a_i \in \gA$ et $1 \in \gen
{a_i, i \in \lrb{0..d}}$. Déduire de la question précédente que
$\IK_\gB(x)$ contient l'image de $\prod_{i=0}^d \IK_\gA(a_i)$.

 \emph {4.}
En déduire une nouvelle preuve du \thref{cor2thKdimMor}:
si tout \elt de $\gB$ est primitivement \agq sur $\gA$, alors
$\Kdim\gB\leq\Kdim\gA$.

}

\end {exercise}
%--- end -exercise-----------------------------------------

%--- Exercise{exoExtEntiereIdealBord}-------------
\begin{exercise}\label{exoExtEntiereIdealBord}
 {(Extension entière de l'\id bord $\IK$)}\\
{\rm  
Soit $\gA \subseteq \gB$ une extension entière d'anneaux.

\emph {1.}
Si $\fa$ est un \id de $\gA$, $\fb$ un \id de $\gB$, montrer que

\snic {
\gA \cap (\fb + \fa\gB) \subseteq \rD_\gA(\fa + \gA\cap\fb)
.}

\emph {2.}
En déduire, pour $a_0$, \ldots, $a_d \in \gA$:

\snic {
\gA\cap\IK_\gB(a_0, \ldots, a_d) \subseteq 
\rD_\gA\big(\IK_\gA(a_0, \ldots, a_d)\big)
.}

\emph {3.}
Donner une nouvelle \dem du fait que $\Kdim\gA \le \Kdim\gB$, cf.
la proposition \ref{propDKEXENT} et le \thref{propGupLY}. Comparer à l'exercice~\ref{exoExtEntiereMonoideBord}.
}

\end {exercise}
%--- end -exercise-----------------------------------------

%--- Exercise{exoExtEntiereMonoideBord}-------------
\begin{exercise}\label{exoExtEntiereMonoideBord}
 {(Extension entière du \mo bord $\SK$)}\\
{\rm  
Soit $\gA \subseteq \gB$ une extension entière d'anneaux.

\emph {1.}
Soient $\fa$ un \id de $\gA$ et $S \subseteq \gA$ un \moz.
Montrer que

\snic {
S + \fa\gB \subseteq \satu {(S + \fa)}{\gB}
.}

\emph {2.}
En déduire, pour $a_0$, \ldots, $a_d \in \gA$:

\snic {
\SK_\gB(a_0, \ldots, a_d) \subseteq 
\satu{\SK_\gA(a_0, \ldots, a_d)}{\gB}
.}

\emph {3.}
Donner une nouvelle \dem du fait que $\Kdim\gA \le \Kdim\gB$.
}

\end {exercise}
%--- end -exercise-----------------------------------------

%--- Exercise{exoKdimSomTr}-------------
\begin{exercise}
\label{exoKdimSomTr}
{\rm  Soit $\gK$ un \cdi non trivial. On note $(\uX)$ pour $(\Xn)$  et~$(\uY)$ pour $(\Ym)$.
On note $\gA=\gK(\uX)\otimes _\gK\gK(\uY)$.
On se propose de déterminer la \ddk de $\gA$.
\begin{itemize}
\item [\emph{1.}] $\gA$ est la \lon de $\gK[\uX,\uY]$ en 
${S=(\KuX)\etl(\KuY)\etl}$.
Il est aussi une \lon de $\gK(\uX)[\uY]$ et de $\gK(\uY)[\uX]$. 
En conséquence
$\Kdim\gA\leq\inf(m,n)$.
\item [\emph{2.}] Supposons $n\leq m$. Montrer que la suite $(X_1-Y_1,\ldots,X_n-Y_n)$
est une suite régulière dans $\gA$.
\end{itemize}
 Conclure que $\Kdim\gA = \inf(n,m)$.
 
}
\end{exercise}
%--- end -exercise-----------------------------------------

%%%%%%%%%%%%%%%%%%%%%%%%%%%%%%%%%%%%%%%%%%%%%%%%%%%%%%%%%%%%%%%%%%%%%%%%%%%

%--- Exercise{exoDualiteBords}-------------
\begin{exercise}\label{exoDualiteBords}
{(Idéaux premiers, bords et dualité)}\\
{\rm  
Soit $\fp_0\subsetneq\fp_1\subsetneq \cdots \subsetneq\fp_{d-1}
\subsetneq\fp_d\subsetneq\gA$ une chaîne d'\ideps détachables avec 
$x_1\in\fp_1\!\setminus\fp_0$, $x_2\in\fp_2\!\setminus\fp_1$, \dots,
$x_{d}\in\fp_{d}\!\setminus\fp_{d-1}$, selon le schéma suivant:

\snic {
\def \foo {\ar@{.}[ld]|\notin\ar@{.}[rd]|\in}
\def \subnot {\subsetneq}
\xymatrix @R=15pt @C=8pt {
     &x_1\foo    &      &x_2\foo &&&
         &x_{d-1}\foo    &      &x_d\foo   \\
0\in\fp_0&\subnot &\fp_1 & \subnot & \fp_2 
&\subnot\cdots\subnot &\fp_{d-2}&\subnot& \fp_{d-1} &\subnot& \fp_d&\!\!\!\!\!\!\not\ni1 }}

\emph {1.}
Montrer que $\IK(x_1, \ldots, x_i) \subseteq \fp_i$ pour $i\in\lrb{0..d}$.  Donc $\IK(x_1, \ldots, x_d) \subseteq \fp_d$. De plus, si $x_{d+1}
\notin \fp_d$, alors $\IK(x_1, \ldots, x_d, x_{d+1}) \subseteq
\fp_d + \gA x_{d+1}$. En conséquence,  \linebreak 
si $x_{d+1}\notin \fp_d$ et 
$1 \in \IK(x_1, \ldots, x_d, x_{d+1}) $, alors $1 \in \fp_d + \gA x_{d+1}$.

\emph {2.}
On considère les filtres premiers \cops $\ff_i = \gA\setminus\fp_i$
pour $i\in\lrb{0..d}$. On a le schéma dual du  précédent:

\snic {
\def \foo {\ar@{.}[ld]|\notin\ar@{.}[rd]|\in}
\def \subnot {\subsetneq}
\xymatrix @R=15pt @C=8pt {
     &x_d\foo    &      &x_{d-1}\foo &&&
         &x_2\foo    &      &x_1\foo   \\
1\in\ff_d&\subnot &\ff_{d-1} & \subnot & \ff_{d-2} 
&\subnot\cdots\subnot &\ff_2&\subnot& \ff_1 &\subnot& \ff_0&\!\!\!\!\!\!\not\ni0 \\
}}

%\sni
Montrer que $\SK(x_{i+1}, \ldots, x_d) \subseteq \ff_i$ pour $i \in \lrb
{0..d}$. Donc $\SK(x_{1}, \ldots, x_d) \subseteq \ff_0$.  De plus,
si $x_{0} \notin \ff_0$, i.e. si $x_0 \in \fp_0$, alors $\SK(x_0, x_1, \ldots,
x_d) \subseteq x_0^\NN \ff_0$. En conséquence,  si $x_0\notin \ff_0$ et $0 \in \SK(x_0,x_1,\ldots, x_d)$, alors $ 0 \in x_0^\NN \ff_0$.

 NB:  $\fp_d + \gA x_{d+1}$ est l'\id engendré par
$\fp_d$ et $x_{d+1}$, dualement $x_0^\NN\ff_0$ est le \mo engendré par $\ff_0$ et $x_0$.

}

\end {exercise}
%--- end -exercise-----------------------------------------

%--- Exercise{exoEliminationEtBord}-------------
\begin{exercise}\label{exoEliminationEtBord}
{(\'Elimination et \ids bord dans les anneaux de \polsz)}\\
{\rm  
 Voici une \dem  détaillée  de l'in\egt 
$\Kdim\gA[T] \le 1 + 2\Kdim\gA$ (section~\ref {secKdimMor}), 
avec quelques précisions. 
\Spdg  $\gA$ est supposé réduit.

\emph {1.}
Soit $f \in \gA[T]$ un \pol tel que l'annulateur de chaque
\coe soit engendré par un \idmz.  Pour $g \in \gA[T]$, définir  $R 
\in \gA[X,Y]$ tel   que $\Ann(R) = 0$ \linebreak
et $R(f,g) = 0$: noter que le \pol
$\Res_T(f(T)-X, Y-g(T)\big)$ résout la question lorsque $f$ est \mon de degré
$\geq1$ (pourquoi?), et utiliser le lemme~\ref{lemQI}.

\emph {2.}
En utilisant l'exercice \ref {exoInclusionBordLionel}, montrer que si $R = \sum_{i,j} r_{ij}X^i Y^j$, on~a:

\snic {
\prod_{i,j} \IK_{\AT}(r_{ij}) \subseteq \IK_{\AT}(f,g).
}

\emph {3.}
En utilisant un anneau de type $\gA_{\so\ua}$ (lemme
\ref {lem4MorRc} et exercice \ref{exoAminEtagesFinis}),
retrouver l'in\egt $\Kdim\gA[T] \le 1 + 2\Kdim\gA$.

\emph {4.}
Montrer le résultat plus précis suivant: pour un anneau réduit
$\gA$ \hbox{et  $f$, $g \in \gA[T]$}, 
 l'\id  $\rD_\AT\left(\IK_{\gA[T]}(f,g)\right)$ contient un produit fini
d'\ids bord $\IK_\gA(a)$, $a \in \gA$.

\emph {5.}
Plus \gnltz: si $\gA[\uT] = \gA[T_1, \ldots, T_r]$ et $f_0$,
\ldots, $f_r \in \gA[\uT]$, alors la racine de l'\id bord $\IK_{\gA[\uT]}(f_0,
\ldots, f_r)$ contient un produit fini d'\ids bord $\IK_\gA(a_i)$, avec les~$a_i \in
\gA$. On en déduit de nouveau que $1 + \Kdim\gA[\uT] \le (1+r)(1 +
\Kdim\gA)$.

}

\end {exercise}
%--- end -exercise-----------------------------------------

%--- Exercise{exoIdealBordPolynomes}-------------
\begin{exercise}\label{exoIdealBordPolynomes}
{(Idéaux bord de \polsz)}
{\rm  
Suite de l'exercice \ref{exoEliminationEtBord}.
 \\
\emph {1.}
Soient $x$, $y \in \gB$ et $(z_j)$ une famille finie
dans $\gB$ vérifiant $\prod_j \IK(z_j) \subseteq \IK(x,y)$.
Montrer que pour $(\bn)$ dans $\gB$,
$\prod_j \IK(z_j,\bn) \subseteq \IK(x,y,\bn)$.

\emph {2.}
Soit $T$ une \idtr sur un anneau $\gA$.
\begin {itemize}
\item [\emph {a.}]
Pour $(\an)$ dans $\gA$, vérifier que $\IK_\gA(\an)\gA[T] = \IK_{\gA[T]}(\an)$. 

\item [\emph {b.}]
Montrer que l'\id bord de $2d$ \pols de $\gA[T]$
contient, à radical près, un produit d'\ids bord de $d$ \elts de $\gA$.
\\
En conséquence $\Kdim\gA < d \Rightarrow \Kdim\gA[T] < 2d$: ceci est une autre forme de
l'in\egt $\Kdim\gA[T] \le 1 + 2\Kdim\gA$.
\end {itemize}

\emph {3.}
Comment peut-on généraliser le premier point? le second?

}
\end {exercise}
%--- end -exercise----------

%--- Exercise{exoKdimEspanol}-------------
\begin{exercise}
\label{exoKdimEspanol} (Une autre \dfn de la \ddk des \trdisz, cf. \cite[Espa\~nol]{Espa08})
{\rm  Dans un ensemble ordonné, une suite $(\xzn)$ est appelée
une \emph{chaîne de longueur $n$} si l'on a $x_0\leq x_1\leq\cdots\leq x_{n}$. 
Dans un \trdi deux chaînes   $(\xzn)$
et  $(b_0,\dots,b_n)$ sont dites \emph{liées}, s'il existe une chaîne  
$(c_1,\dots,c_n)$ avec 
%---  equation eqEspanol --------
\begin{equation}\label{eqEspanol}
\left.\arraycolsep2pt
\begin{array}{rcccl}
 x_0\vi b_0& =  & 0    \\
 x_1\vi b_1& = & c_1 &=  &  x_0\vu b_0  \\
\vdots~~~~& \vdots &\vdots& \vdots & ~~~~  \vdots \\
 x_n\vi  b_n & = & c_n & =& x_{n -1}\vu b_{n -1}  \\
&& 1& =  &   x_n\vu b_n
\end{array}
\right\}
\end{equation}
%---------------------end equation--------------
On pourra comparer avec la \dfn \ref{defiDDKTRDI} pour les suites \copsz.
Noter aussi que si des suites $(\xzn)$, $(b_0,\dots,b_n)$ et $(c_1,\dots,c_n)$ satisfont les \eqns
(\ref{eqEspanol}), alors ce sont des chaînes.

\emph{1.} Si dans un \trdi on a $x \le y$ et 
$x\vu  a\geq  y \vi b$, alors on peut expliciter~$a'$ et $b'$ tels que

\snic {
x\vi a' =  x\vi a, \qquad  y \vu b' = y\vu b, \qquad
x\vu  a' =  y \vi b'
.}

%\sni
Donc à partir d'une configuration de gauche (en supposant toujours $x \le y$), 
on peut construire une configuration de droite:
$$
\left\{\,\arraycolsep2pt
\begin{array}{rcl}
x\vi a   & =  &  p  \\
x\vu  a & \geq  & y \vi b  \\
q & =  &   y\vu b
\end{array}
\right.
\qquad\qquad
\left\{\,\arraycolsep2pt
\begin{array}{rcl}
x\vi a'   & =  &  p  \\
 x\vu  a' & =  & y \vi b'  \\
 q& =  &   y\vu b'
\end{array}
\right.
$$

\emph{2.} Dans un \trdi une chaîne $(\xzn)$ possède une suite
\cop \ssi il existe une chaîne qui lui est liée.

\emph{3.} Pour un \trdi $\gT$ \propeq
\begin{enumerate}
\item [\emph{a.}] $\gT$ est de \ddk $\leq n$.
\item [\emph{b.}] Toute chaîne de longueur $n$ admet une suite \copz.  
\item [\emph{c.}] Toute chaîne de longueur $n$ admet une chaîne qui lui est liée.
\end{enumerate}
}
\end{exercise}
%--- end -exercise-----------------------------------------

%--- Exercise{exoAminEtagesFinis}-------------
\begin{exercise}\label{exoAminEtagesFinis}
 {(Quelques précisions sur les étages finis de $\Amin$)}
\\
{\rm  Soit $\gA$ un anneau réduit. Pour 
$\fa$, $\fb$ \ids de $\gA$ on note $\fa\diamond\fb = (\fa\epr\fb)\epr =
(\fa^{\perp\perp} : \fb)$.

\emph {1.}
Vérifier que $\gA\sur{\fa\diamond\fb}$ est un anneau réduit dans lequel
$\fa$ est nul et $\fb$ fidèle.

\emph{2.}
Vérifier que $(\gA\sur{\fa_1\diamond\fb_1})\sur{(\ov{\fa_2}\diamond\ov{\fb_2})}
\simeq \gA\sur{\fa_3\diamond\fb_3}$ avec $\fa_3 = \fa_1+\fa_2$,
$\fb_3 = \fb_1\fb_2$.

\emph {3.}
Soit $(\ua) = (a_1, \ldots, a_n)$ dans $\gA$.
Dans le lemme \ref{lem4MorRc} on a défini (pour 
$I\in\cP_n$):

\snic {
\fa_I = \gen {a_i, i \in I} \diamond \prod_{j\notin I}a_j
\qquad
\gA_{\so\ua} = \prod_{I\in\cP_n} \gA\sur{\fa_I}.
}

%\sni
Ainsi, modulo $\fa_I$, $a_i$ est nul pour $i \in I$ et \ndz pour $i \notin I$.
On notera~$\vep_i$  l'\idm de $\gA_{\so\ua}$ dont la \coo dans
$\gA\sur{\fa_I}$ est $1$ si $i \in I$ et $0$ si $i \notin I$.
\begin {itemize}
\item [\emph {a.}]
Vérifier que l'intersection (et a fortiori le produit) des \ids $\fa_I$ est
nulle; en conséquence, le morphisme $\gA \to \gA_{\so\ua}$ est 
injectif et $\Kdim\gA = \Kdim\gA_{\so\ua}$.

\item [\emph {b.}]
Vérifier que $\Ann_{\gA_{\so\ua}}(a_i) = \gen{\vep_i}_{\gA_{\so\ua}}$.
\end {itemize}

}
\end {exercise}
%--- end -exercise-----------------------------------------

%--- Exercise{exoAmin}-------------
\begin{exercise}
\label{exoAmin} (Quelques précisions sur $\Amin$)
\index{regulier@régulier!morphisme --- (d'anneaux)}
\index{morphisme!regulier@régulier (d'anneaux)}
\\
{\rm Voir le \pbz~\ref{exoQiClot} pour ce qui concerne $\Aqi$.
\\
  Un \homo d'anneaux $\gA\to\gB$ est dit \emph{\regz}
lorsque l'image de tout \elt \ndz est un \elt \ndzz.
\\
\fbox{Soit $\gA$ un anneau réduit.}

\vspace{-.5em}
\begin{enumerate}\itemsep0pt
\item Soit $\theta:\gA\to\gB$ un \homo \reg et $a\in\gA$. Si  $a\epr$ est engendré par un \idm $e$, alors $\theta(a)\epr$ est engendré par l'\idm $\theta(e)$. 
\\
En particulier, comme déjà noté dans le \pb \ref{exoQiClot}, un \homo entre anneaux \qis est \qi \ssi il est \regz.
\item L'\homo naturel  $\Aqi \to \Amin$  est \reg et surjectif.
\item Pour $a\in\gA$, l'\homo naturel $\psi_a:\gA\to\gA_{\so{a}}$ est \regz. 
\item L'\homo naturel $\psi:\gA\to\Amin$ est \reg et  
l'\homo naturel $\ZZ\to\ZZ_\mathrm{qi}$ n'est pas \regz.
\end{enumerate}
}
\end{exercise}
%--- end -exercise-----------------------------------------

%--- Exercise{exolemVdimKdim}-------------
% juste apres exoAmin
\begin{exercise}
\label{exolemVdimKdim}
{\rm  Expliciter la \dem du lemme \ref{lemVdimKdim} en termes 
de \susisz. 
}
\end{exercise}
%--- end -exercise-----------------------------------------

%--- Exercise{exothValDim}-------------
\begin{exercise}
\label{exothValDim} (Une \gnn du \thref{thValDim})\\
{\rm  
Pour $\gA\subseteq\gB$ et $\ell\in \NN$, si pour toute suite $(\ux)=(x_0,\ldots,x_\ell)$
dans $\gB$, on a un \pol primitif de $\gA[\uX]$ qui annule 
$(\ux)$, alors  $\Vdim\gB\leq\ell + \Vdim\gA$. 
 
}
\end{exercise}
%--- end -exercise-----------------------------------------

%--- Exercise{exoLyingOverClassique}-------------
\begin{exercise}
\label{exoLyingOverClassique} 
{(Morphisme lying over)}
\\
{\rm Démontrer ce qui est affirmé  dans la remarque qui suit la
\dfn du lying over \paref{remLY}.
}
\end{exercise}
%--- end -exercise-----------------------------------------

%--- Exercise{exoLYOV}-------------
\begin{exercise}\label{exoLYOV} {(Morphisme lying over, 2)}
\\
{\rm  Dans la catégorie des ensembles ordonnés finis, il est clair qu'un morphisme
est surjectif \ssi c'est un épimorphisme.
Cela correspond donc pour les \trdis duaux à un monomorphisme,
ce qui signifie ici  un homomorphisme injectif, \cad lying over.
\\
Donner une preuve en \clama de l'\eqvcz, pour un \homo $\alpha:\gT\to\gT'$
de \trdisz, entre: $\alpha$ est  lying over  d'une part,
et $\Spec \alpha:\Spec \gT'\to\Spec\gT$ est surjectif, d'autre part. 
\\
Idée. Utiliser le \emph{lemme de Krull}, qui se démontre par une zornette: \emph{si dans un \trdi on a un \id $\fa$
et un filtre $\ff$ qui ne se coupent pas, il existe un \idep contenant $\fa$
dont le \cop est un filtre contenant $\ff$.} Voir \egmt la remarque qui suit
le lemme \ref{lemLingOver}.  
}
\end{exercise}
%--- end -exercise-----------------------------------------

%--- Exercise{exodefiGoingup}-------------
\begin{exercise}
\label{exodefiGoingup} {(Morphisme going up, going down)}\\
{\rm  Démontrer ce qui est affirmé dans la remarque qui suit la
\dfn du going up \paref{defiGoingup} (utiliser la description
du treillis quotient $\gT\sur{(\fa=0)}$ donnée \paref{trquoideal}).
Même chose avec le going down.
 }
\end{exercise}
%--- end -exercise-----------------------------------------

%--- Exercise{exodefiGoingup2}-------------
\begin{exercise}
\label{exodefiGoingup2} {(Morphisme going up, going down, 2)}\\
{\rm  Démontrer en \clama que la \dfn \ref{defiGoingup} que nous avons donnée pour un morphisme
going up $\alpha:\gT\to\gV$ équivaut à la définition classique expliquée
juste avant. Plus précisément, pour un morphisme $\alpha:\gT\to \gV$ de \trdis \propeq
\begin{enumerate}
\item Pour tout \idep $\fq$ de $\gV$, en notant $\fp=\alpha^{-1}(\fq)$ le morphisme 
$$\alpha':\gT/(\fp=0)\to \gV/(\fq=0)$$ est injectif (i.e., lying over).
\item Pour tout idéal $\fa$ de $\gV$, en notant $\fb:=\alpha^{-1}(\fa)$, le morphisme $$\alpha':\gT/(\fb=0)\to \gV/(\fa=0)$$ est injectif. 
\item Pour tout $i\in \gV$, avec $\fb=\alpha^{-1}(\dar i)$ le morphisme $$\alpha':\gT/(\fb=0)\to \gV/(i=0)$$ est injectif. 
\item Pour
tous $a,b\in \gT$ et $i\in \gV$ on a
$$
\alpha(a)\vdi \gV \alpha(b),\, i \quad\Longrightarrow\quad\exists j\in \gT \;\;\; a \vdi \gT b ,\, j \;\;\; \hbox{et}\;\;\; \alpha(j)\leq_\gV i.
$$
\end{enumerate}
En fait on a \cot \emph{2} $\Leftrightarrow$ \emph{3} $\Leftrightarrow$ \emph{4} $\Rightarrow$ \emph{1} et seule l'implication \emph{1} $\Rightarrow$ \emph{4}
utilise le tiers exclu et le lemme de Zorn.  
 }
\end{exercise}
%--- end -exercise-----------------------------------------

%%%%%%%%%%%%%%%%%%%%%%%%%%%%%%%%%%%%%%%%%%%%%%%%%%%%%%%%%%%%%%%%%%%%%%%%%%%
%:  pbs

%%%%%%%%%%%%%%%%%%%%%%%%%%%%%%%%%%%%%%%%%%%%%%%%%%%%%%%%%%%%%%%%%%%%%%%%%%%

%--- problem{exoAnneauNoetherienReduit}-------------
\begin{problem}
\label{exoAnneauNoetherienReduit}
{(Annulateur d'un \id dans un anneau \noe réduit)}\\
{\rm
On considère un anneau
réduit  $\gA$ tel que toute suite croissante d'\ids de la
forme $\DA(x)$
possède deux termes consécutifs égaux.

\emph{1.}
Soit $\fa$ un idéal  de $\gA$ tel que l'on sache tester pour $y\in\gA$ si $\Ann(y)\fa=0$
(et en cas de réponse négative fournir le certificat correspondant).
\\
\emph{1a.} Si un $x \in \fa$ vérifie $\Ann(x)\fa \ne 0$,
 déterminer un $x' \in \fa$ tel que $\DA(x) \subsetneq \DA(x')$.
\\
\emph{1b.}
En déduire l'existence d'un $x \in \fa$ tel que $\Ann(x) = \Ann(\fa)$.

\emph{2.}
On suppose de plus
que tout \elt\ndz de $\gA$ est \ivz, et que pour tous $y$, $z$
on sait tester  si $\Ann(y)\Ann(z)=0$.
Montrer que $\Kdim\gA \le 0$.

\emph{3.}
Soit $\gB$ un anneau \noe \coh \fdiz. Montrer que 
 $\Frac(\gB\red)$ est un anneau \zedz.
\\
NB. En \clama $\gB$ admet un nombre fini d'\idemis 
$\fp_1$, \dots, $\fp_k$
et $\Frac(\gB\red)$ est isomorphe au produit fini de corps correspondant: $\Frac(\gA/\fp_1)\times \cdots\times \Frac(\gA/\fp_k)$. En \gnl cependant, on n'a pas accès aux
$\fp_i$ de façon \algqz.
}
\end {problem}
%--- end-problem-----------------------------------------

%--- problem{exoContractedInclusion}-------------
\begin{problem}
\label{exoContractedInclusion} (Lying over, going up, going down, exemples)
\\
{\rm
\emph {1.}
Soit une inclusion d'anneaux $\gA\subseteq\gB$ telle que, en tant que \Amoz,
$\gA$ soit facteur direct dans $\gB$. Montrer que $\fa\gB \cap \gA = \fa$
pour tout idéal $\fa$ de~$\gA$.  En particulier, $\gA \hookrightarrow
\gB$ est lying over.

\emph {2.}
Soit $G$ un groupe fini agissant sur un anneau $\gB$ avec
$\abs G 1_\gB$ inversible dans~$\gB$. On note $\gA = \gB^G$ le
sous-anneau des points fixes et l'on définit l'\emph{opérateur
de Reynolds} $R_G : \gB \to \gA$:
\index{operateur@opérateur de Reynolds}

\snic {
R_G(b) = {1 \over \abs G} \sum_{g \in G} g(b)
.}

%\sni
Vérifier que $R_G$ est un $\gA$-\prr d'image $\gA$; en
particulier, $\gA$ est facteur direct (comme \Amoz) dans $\gB$.

\emph {3.}
Soit $\gA \hookrightarrow\gB$ avec $\gA$ facteur direct (comme \Amoz)
dans $\gB$.  Fournir une preuve directe de $\Kdim\gA \le \Kdim\gB$.

\emph {4.}
Soient $\gk$ un \cdi non trivial et $\gA = \gk[XZ,YZ] \subset\gB = \gk[X,Y,Z]$.
Alors $\gA$ est facteur direct dans $\gB$, donc $\gA \hookrightarrow
\gB$ est lying over. Mais $\gA \hookrightarrow \gB$ n'est ni going up ni
going~down.

}
\end {problem}
%--- end -problem-----------------------------------------

%--- problem{exoChainesPotPrem}-------------
\begin{problem}
\label{exoChainesPotPrem} (Chaînes potentielles d'\idepsz)\index{chaine po@chaîne potentielle!d'ideaux pr@d'\idepsz}
\\
{\rm  Sur un anneau $\gA$ on appelle \emph{chaîne potentielle d'\idepsz},
ou encore \emph{chaîne potentielle}
une liste $[(I_0,U_0),\ldots,(I_n,U_n)]$, où les $I_j$ et $U_j$
sont des parties de $\gA$ (i.e., chaque $(I_j,U_j)$ est un \idep potentiel de $\gA$). Une chaîne potentielle est dite \emph{finie} si les $I_j$
et $U_j$ sont des parties finiment énumérées.
\\
Une chaîne potentielle est dite \emph{complète} si les conditions suivantes sont satisfaites:
\begin{itemize}
\item les $I_j$ sont des \ids et les $U_j$ sont des \mosz,
\item $I_0\subseteq I_1\subseteq \cdots\subseteq I_n$ et  $U_0\supseteq U_1\supseteq \cdots\supseteq U_n$,
\item $I_j+U_j=U_j$ pour chaque $j$.
\end{itemize}
On dit que la chaîne potentielle $[(I_0,U_0),\ldots,(I_n,U_n)]$
\emph{raffine} la chaîne $[(J_0,V_0),\alb\ldots,\alb(J_n,V_n)]$ si l'on  a les inclusions
$J_k\subseteq I_k$ et $V_k\subseteq U_k$ pour chaque $k$.

 \emph{1.} Toute chaîne potentielle engendre une chaîne potentielle 
complète (au sens de la relation de raffinement). Plus \prmtz, à partir de $[(I_0,U_0),\ldots,(I_n,U_n)]$, on construit successivement
\begin{itemize}
\item $\fa_j=\gen{I_j}$, $\fb_j=\sum_{i\leq j}\fa_i$ ($j\in\lrb{0..n}$),
\item $\ff_n=\cM(U_n)+\fb_n$(\footnote{Rappelons que $\cM(A)$ est le \mo engendré par la partie $A$.}), $\ff_{n-1}=\cM(U_{n-1}\cup \ff_n)+\fb_{n-1}$, \ldots,
$\ff_{0}=\cM(U_{0}\cup \ff_1)+\fb_{0}$.
\end{itemize}
Et l'on considère $[(\fb_0,\ff_0),\ldots,(\fb_n,\ff_n)]$.  

 \emph{2.} On dit qu'une chaîne potentielle $\cC$ \emph{collapse} si
dans  la chaîne
complète qu'elle engendre $[(\fb_0,\ff_0),\ldots]$ on a $0\in\ff_0$.
Montrer qu'une suite $(x_1,\ldots,x_n)$ est \sing \ssi la chaîne
potentielle $[(0,x_1),(x_1,x_2),\ldots,(x_{n-1},x_n),(x_n,1)]$ collapse.

 \emph{3.} En \clamaz, une chaîne potentielle 
$\cC$ de $\gA$ collapse
\ssi il est impossible de trouver  des \ideps   
$\fp_0\subseteq\fp_1\subseteq\cdots\subseteq\fp_n$, tels que 
%, en notant $\fh_j=\gA\setminus \fp_j$, la chaîne $[(\fp_0,\fh_0),\alb\ldots,\alb(\fp_n,\fh_n)]$ 
la chaîne $[(\fp_0,\gA\setminus \fp_0),\ldots,\alb(\fp_n,\gA\setminus \fp_n)]$ raffine
la chaîne $\cC$.

 \emph{4.} \'Etant donnée une chaîne potentielle $\cC=[(I_0,U_0),\ldots,(I_n,U_n)]$, on la \emph{sature} en  rajoutant dans $I_k$ (resp. dans $U_k$) tout
$x\in\gA$ qui, rajouté à $U_k$ (resp. à $I_k$) conduirait à un collapsus.
Ainsi une chaîne potentielle collapse \ssi sa saturée est $[(\gA,\gA),\ldots,(\gA,\gA)]$.
\\
Montrer que l'on obtient ainsi une  chaîne potentielle $[(J_0,V_0),\ldots,(J_n,V_n)]$ qui raffine la chaîne complète engendrée par~$\cC$.
\\
Montrer en \clama que $J_k$ est l'intersection des \ideps qui s'insèrent en position $k$ dans une chaîne d'\ideps qui raffine $\cC$
(comme dans la question précédente). Démontrer aussi l'affirmation duale pour $V_k$.   
 
}
\end{problem}
%--- end -problem-----------------------------------------

}% fin des exos

%%%%%%%%%%%%%%%%%%%%%%%%%%%%%%%%%%%%%%%%%%%%%%%%%%%%%%%%%%%%%%%%%%%%%%%%%%%
%:  solutions
\sol{
%%%%%%%%%%%%%%%%%%%%%%%%%%%%%%%%%%%%%%%%%%%%%%%%%%%%%%%%%%%%%%%%%%%%%%%%%%%

%%%%%%%%%%%%%%%%%%%%%%%%%%%%%%%%%%%%%%%%%%%%%%%%%%%%%%%%%%%%%%%%%%%%%%%%%%%
\exer{exoDimKX}
\emph{2.} On considère
une suite de longueur $k$ dans $\gA$, on
lui rajoute $X$  au début, et elle devient \sing
dans $\gA[X]$. On se
débarrasse ensuite de $X$ dans l'\egrf{eqsing} \paref{eqsing} correspondante.
\\
NB: on peut aussi invoquer le point \emph{3} de la proposition \ref{lemRegsing}.

%%%%%%%%%%%%%%%%%%%%%%%%%%%%%%%%%%%%%%%%%%%%%%%%%%%%%%%%%%%%%%%%%%%%%%%%%%%
\exer{exoDimKX1n}
On peut supposer $\gK$  réduit 
($\gK\red[X_1,\ldots ,X_n]=\KXn\red$ a même dimension que $\KXn$).
Deux possibilités s'offrent alors. La première est de
réécrire la preuve donnée dans le cas d'un
\cdi en utilisant l'exercice~\ref{exoZDpiv} et le \plgrf{thDdkLoc}.
La deuxième est d'appliquer la machinerie \lgbe \elr \num2.\imlg

%%%%%%%%%%%%%%%%%%%%%%%%%%%%%%%%%%%%%%%%%%%%%%%%%%%%%%%%%%%%%%%%%%%%%%%%%%%

\exer{exoPartitionUnite} 
On écrit chaque $U_i$ sous la forme
$U_i = \bigcup_{j \in J_i} \DA(g_{ij})$. Dire que les~$\DA(g_{ij})$  recouvrent $\Spec(\gA)$, 
signifie que
$1 \in \gen {\,\DA(g_{ij})\mid j \in J_i,\ i \in I}$,
d'où une \egt
$1 = \sum_{j,i} u_{ji} g_{ij}$,
les $u_{ji}$ étant nuls sauf un nombre fini d'entre eux (\hbox{i.e., $i\in I_0$}, $j\in J_i$,
$I_0$ et les $J_i$ finis).
\\ 
On pose $f_i = \sum_{j\in J_i} u_{ji}g_{ij}$. On obtient $\DA(f_i) \subseteq U_i$  car pour $\fp \in \DA(f_i)$, on \hbox{a $f_i \notin\fp$}, donc un indice $j$ tel 
que $g_{ij} \notin \fp$, i.e. $\fp \in \DA(g_{ij})
\subseteq U_i$. Et  $\sum_{i\in I_0} f_i = 1$.
%%%%%%%%%%%%%%%%%%%%%%%%%%%%%%%%%%%%%%%%%%%%%%%%%%%%%%%%%%%%%%%%%%%%%%%%%%%

\exer{exoKdimGeom} ~\\
\emph{1a.} On écrit une \rdi de $f$ sur $\KYr$ 

\snic{f^n+a_{n-1}f^{n-1}+\cdots+a_kf^k=0,}

%\sni
avec $n\geq1$, les $a_i\in\KYr$ et $a_k\neq0$. 
L'\egt $(a_k+bf)f^{k}=0$ montre que $a_k+bf\in(\DA(0):f)$
(même si $k=0$). Donc $a_k\in \JK_\gA(f)$.

%%%%%%%%%%%%%%%%%%%%%%%%%%%%%%%%%%%%%%%%%%%%%%%%%%%%%%%%%%%%%%%%%%%%%%%%%%%
\exer{exolemLocMemeKdim} 
\emph{1.} On écrit $\KuX[1/f]=\aqo{\gK[\uX,T]}{1-fT}$.\\
Alors, une mise en position de \iNoe du \pol non constant $1-fT$ nous ramène à une extension entière de $\gK[\Yn]$.
\\
\emph{2.} On écrit $\AuX[1/\delta]=\aqo{\gA[\uX,T]}{1-\delta T}$. On cherche à appliquer le  \thref{cor2thKdimMor} sur les extensions entières.
On veut d'une part que $\delta $ soit \ndzz, pour que l'\homo  $\AuX\to\AuX[1/\delta ]$ soit injectif, et d'autre part que l'on puisse faire une mise en position de \iNoe du
\pol $1-\delta T$, pour que~$\AuX[1/\delta ]$ soit entier sur un anneau
$\gA[\Yn]$. 
\\
La première condition signifie que l'\id $\rc(\delta )$ est fidèle
(McCoy, corolaire \ref{corlemdArtin}). 
\\
La deuxième condition est réalisée si l'on se trouve dans la même situation que pour le lemme~\ref{lemLocMemeKdim}:  
\begin{itemize}
\item $\delta$ est de degré formel $d$,
\item  un des \moms de degré $d$, portant sur un sous-ensemble de 
variables $(X_i)_{i\in I}$, a pour \coe
un \elt de $\Ati$,
\item  et c'est le seul \mom de degré $d$
en les variables $(X_i)_{i\in I}$ présent dans $\delta$.
\end{itemize}   
En effet, le \cdv \gui{$X'_i = X_i + T$ si $i \in I$, $X'_i = X_i$ sinon},
rend alors le \pol $1 - \delta T$ \mon en $T$ (à un \iv près).
Notons que dans ce cas le \pol $\delta$ est primitif et la première condition est \egmt réalisée.

%%%%%%%%%%%%%%%%%%%%%%%%%%%%%%%%%%%%%%%%%%%%%%%%%%%%%%%%%%%%%%%%%%%%%%%%%%%

%%%%%%%%%%%%%%%%%%%%%%%%%%%%%%%%%%%%%%%%%%%%%%%%%%%%%%%%%%%%%%%%%%%%%%%%%%%
\exer{exoPruferKdim} \emph {1.} On considère $s=a/b\in\Frac\gA$ avec $b$ \ndzz.\\
La suite $(bX-a,b,X)$ est \sing dans $\AX$. Cela donne une \egt 
dans~$\AX$ du type suivant

\snic{
(bX-a)^{k_1}\big(b^{k_2}\big(X^{k_3}\big(1+Xp_3(X)\big)+bp_2(X)\big)+(bX-a)p_1(X)\big)=0.}

%\sni
Puisque $\AX$ est intègre, on peut supprimer le facteur $(bX-a)^{k_1}$,
à la suite de quoi on spécialise $X$ en $s$. Il vient

\snic{
b^{k_2}\big(s^{k_3}\big(1+sp_3(s)\big)+bp_2(s)\big)=0,
}

%\sni
et puisque $b$ est \ndzz:

\snic{
s^{k_3}\big(1+sp_3(s)\big)+bp_2(s)=0.}

%\sni
Ainsi $s$ annule %les deux \pols 
$g(X)=X^{k_3}\big(1+Xp_3(X)\big)+bp_2(X)$
\hbox{et  $f(X)=bX-a$}.\\
Enfin, puisque le \coe de $X^{k_3}$ dans $g$ est de la forme $1+bc$, on obtient \hbox{que $1\in\rc(f)+\rc(g)=\rc(f+X^2g)$}.

\emph{2.} Résulte de \emph{1} et des résultats \gnls sur la dimension de
$\AX$, pour un anneau arbitraire et pour un \adpz.

\emph{3.} Il semble bien que oui.

%%%%%%%%%%%%%%%%%%%%%%%%%%%%%%%%%%%%%%%%%%%%%%%%%%%%%%%%%%%%%%%%%%%%%%%%%%%

\exer{exoMultiplicativiteIdeauxBords} 
\emph {1.}
Il suffit de montrer, pour deux \ids $\fa$, $\fb$ et deux \eltsz~\hbox{$u$, $v \in \gA$},
que:
$$\preskip.0em \postskip.4em
\begin{array} {c}
\big((\fa : u^\infty) + \gA u\big)\,\big((\fb : u^\infty) + \gA u\big) \subseteq
(\fa\fb : u^\infty) + \gA u \et
\\[1mm]
\big((\fa : u^\infty) + \gA u\big) \, \big((\fa : v^\infty) + \gA v\big) \subseteq
(\fa : (uv)^\infty) + \gA uv.
\end{array}
$$
La première inclusion découle de $(\fa : u^\infty)\,(\fb : u^\infty)
\subseteq (\fa\fb : u^\infty)$ et la seconde 
\hbox{de $(\fa : u^\infty) + (\fa :
v^\infty) \subseteq (\fa : (uv)^\infty)$}.

%%%%%%%%%%%%%%%%%%%%%%%%%%%%%%%%%%%%%%%%%%%%%%%%%%%%%%%%%%%%%%%%%%%%%%%%%%%

\exer{exoInclusionBordLionel} 
\emph {1.}
Comme $\beta > \alpha$, $\uX^\beta$ est multiple de l'un des
\moms suivants:

\snic {
X_1^{\alpha_1} X_2^{\alpha_2} \cdots X_{n-1}^{\alpha_{n-1}}X_n^{1+\alpha_n},
\;
X_1^{\alpha_1} X_2^{\alpha_2} \cdots X_{n-1}^{1+\alpha_{n-1}},
\; \dots,\; 
X_1^{\alpha_1}  X_2^{1+\alpha_2},\;
X_1^{1+\alpha_1}.
}

%\sni
\emph {2a.}
Soit $y\in \prod_{\beta < \alpha} \Ann(a_\beta)$; en posant $Q(\uX) =
yP(\uX)$, on a

\snic{Q(\uX) = ya_\alpha\uX^\alpha + \sum_{\beta>\alpha} ya_\beta
\uX^\beta\quad$ et $\quad Q(\ux) = 0.}

%\sni
Pour montrer que $ya_\alpha\in\IK(\ux)$, on peut
donc supposer que l'on a $y=1$ et~\hbox{$P(\uX) = a_\alpha\uX^\alpha +
\sum_{\beta>\alpha} a_\beta\uX^\beta$}.
En utilisant l'\egt $P(\ux) = 0$ et la première question,
on obtient $a_\alpha \in \IK(\ux)$.

\emph {2b.}
D'abord, puisque $\gA$ est réduit, on a $\IK(a) = \Ann(a) + \gA a,\;\forall a \in \gA$. Ensuite, on utilise la remarque suivante: soit $\fc$ un \id
et $2m$ \ids $\fa_1$, $\fb_1$, \ldots, $\fa_m$, $\fb_m$ tels que
$\fa_1\cdots\fa_{k-1}\fb_k\subseteq \fc$ pour tout $k \in \lrbm$. Alors
on obtient l'inclusion:

\snic {
(\fa_1+\fb_1) \cdots (\fa_m+\fb_m) \subseteq \fc + \fa_1\cdots\fa_m.
}

%\sni
En effet, par \recu sur $m$, si
$
(\fa_1+\fb_1) \cdots (\fa_{m-1}+\fb_{m-1}) \subseteq \fc + \fa_1\cdots\fa_{m-1}$,  on déduit:

\snic {
(\fa_1+\fb_1) \cdots (\fa_m+\fb_m) \subseteq \fc + 
\fa_1\cdots\fa_{m-1}\fa_m + \fa_1\cdots\fa_{m-1}\fb_m
\subseteq \fc  + 
\fa_1\cdots\fa_{m-1}\fa_m+\fc  
}

%\sni
d'où l'inclusion
annoncée.
\\
Appliquons cela à $\fc = \IK(\ux)$ et aux \ids $\fa_\beta = \Ann(a_\beta)$,
$\fb_\beta = \gA a_\beta$.  \\
Comme $\Ann(a_\beta) + \gA a_\beta =
\IK(a_\beta)$, on obtient l'inclusion désirée.

\emph {3.}
Application directe avec $n = 1$.

\emph {4.}
On peut supposer $\gA$, $\gB$ réduits quitte à remplacer $\gA \to \gB$
par $\gA\red \to \gB\red$ (tout $z\in\gB\red$ reste primitivement \agqz). On peut aussi
supposer $\gA \subseteq \gB$ quitte à remplacer~$\gA$ par son
image dans $\gB$. \\
Montrons alors l'implication $\Kdim\gA \le m \Rightarrow \Kdim\gB \le m$
par \recu sur $m$. Il suffit de montrer pour $x\in\gB$ arbitraire  que
$\Kdim(\gB\sur{\IK_\gB(x)}) \le m-1$; mais~$\IK_\gB(x)$ contient un
\id $\fa$ de $\gA$ produit fini d'\ids bord $\IK_\gA(a)$, $a\in \gA$.
\\
On a donc une \alg $\gA\sur\fa \to \gB\sur{\IK_\gB(x)}$   à laquelle on peut appliquer l'hypothèse
de \recu puisque $\Kdim \gA\sur\fa \le m-1$.

%%%%%%%%%%%%%%%%%%%%%%%%%%%%%%%%%%%%%%%%%%%%%%%%%%%%%%%%%%%%%%%%%%%%%%%%%%%
\exer{exoExtEntiereIdealBord} 
\emph {1.}
On utilise l'extension entière $\ov\gA = \gA\sur{\gA\cap\fb}
\hookrightarrow \ov\gB = \gB\sur{\fb}$. 
\\
Soit $a \in \gA \cap (\fb + \fa\gB)$;
le lying over (\ref{lemLingOver}) avec $\ov\gA \subseteq \ov\gB$,
donne $\ov a^n \in \ov\fa$, \hbox{i.e. $a^n \in \fa + \fb$} et comme
$a \in \gA$, $a^n \in \fa + \gA\cap\fb$.

\emph {2.}
Par \recu sur $d$. On pose

\snuc {
\fa = \IK_\gA(a_0, \ldots, a_{d-1}), \
\fa' = \IK_\gA(a_0, \ldots, a_d), \
\fb = \IK_\gB(a_0, \ldots, a_{d-1}), \
\fb' = \IK_\gB(a_0, \ldots, a_d)
.}

%\sni
On a donc par \dfn $\fa' = (\fa : a_d^\infty)_\gA + \gA a_d$ et
$\fb' = (\fb : a_d^\infty)_\gB + \gB a_d$. On veut montrer que
$\gA\cap\fb' \subseteq \rD_\gA(\fa')$. Le point \emph {1},
donne $\gA\cap\fb' \subseteq \rD_\gA(\fc)$ \hbox{avec $\fc = 
\gA a_d + \gA \cap (\fb : a_d^\infty)_\gB =
\gA a_d + (\gA\cap\fb : a_d^\infty)_\gA$}. 
\\
Par \recuz, $\gA\cap\fb
\subseteq \rD_\gA(\fa)$, donc
$$\preskip.0em \postskip.4em
\fc \subseteq \gA a_d + (\rD_\gA(\fa) : a_d^\infty)_\gA 
\subseteq  \rD_\gA(\gA a_d + (\fa : a_d^\infty)_\gA)
\eqdf {\rm def} \rD_\gA(\fa')
,
$$
d'où $\gA\cap\fb' \subseteq \rD_\gA(\fa')$. 
%%%%%%%%%%%%%%%%%%%%%%%%%%%%%%%%%%%%%%%%%%%%%%%%%%%%%%%%%%%%%%%%%%%%%%%%%%%

\exer{exoExtEntiereMonoideBord} 
\emph {1.}
Soit $t \in S + \fa\gB$; i.e. 
$t+s \in \fa\gB$ avec $s \in S$. 
\\
Alors $t+s$ est entier sur $\fa$, donc zéro d'un \polu 

\snic{P(X) \in X^n + \fa X^{n-1} + \cdots +
\fa^{n-1} X + \fa^n .}

On écrit que
$P(T+s) = TQ(T) + P(s)$ et l'on remarque que $P(s) \in s^n + \fa$. \\
Ainsi, $tQ(t) \in S + \fa$. 

\emph {2.}
Par \recu sur $d$.  Soit $V = \SK_\gB(a_0, \ldots, a_d) = a_0^\NN(
\SK_\gB(a_1, \ldots, a_d) + a_0\gB)$; la \recu fournit $\SK_\gB(a_1, \ldots,
a_d) \subseteq \satu{\SK_\gA(a_1, \ldots, a_d)}{\gB}$ donc

\snic {
V \subseteq 
a_0^\NN( \satu{\SK_\gA(a_1, \ldots, a_d)}{\gB} + a_0\gB) 
\subseteq
a_0^\NN \satu{(\SK_\gA(a_1, \ldots, a_d) + a_0\gB)}{\gB}
.}

%\sni
La première question fournit:

\snic {
V \subseteq a_0^\NN \satu{(\SK_\gA(a_1, \ldots, a_d) + a_0\gA)}{\gB}
\subseteq \satu{(a_0^\NN(\SK_\gA(a_1, \ldots, a_d) + a_0\gA)\big)}{\gB}
,}

%\sni
\cad $V \subseteq \satu{\SK_\gA(a_0, \ldots, a_d)}{\gB}$.

%%%%%%%%%%%%%%%%%%%%%%%%%%%%%%%%%%%%%%%%%%%%%%%%%%%%%%%%%%%%%%%%%%%%%%%%%%%

\exer{exoKdimSomTr} 
\emph{2.} L'anneau quotient $\aqo\gA{X_1-Y_1}$ peut être vu comme
la \lon de~$\gK[\Xn,Y_2,\ldots,Y_m]$ en 

\snic{S_1=(\KXn)\etl(\gK[X_1,Y_2,\ldots,Y_m])\etl.}

%\sni
Il est donc intègre. On décrit de la même façon les quotients successifs.

%%%%%%%%%%%%%%%%%%%%%%%%%%%%%%%%%%%%%%%%%%%%%%%%%%%%%%%%%%%%%%%%%%%%%%%%%%%

\exer{exoDualiteBords} 
\emph {1.}
Soit $\fa_i := \IK(x_1, \ldots, x_i)$, avec 
$\fa_{i+1} = (\fa_i : x_{i+1}^\infty) + \gA x_{i+1}$. Par \recuz,
 $\fa_i \subseteq \fp_i$:  $x_{i+1}\notin\fp_i$ donne $(\fp_i : x_{i+1}^\infty)
\subseteq \fp_i$, puis $\fa_{i+1} \subseteq \fp_i + \gA x_{i+1}$, donc $
\fa_{i+1}\subseteq \fp_{i+1}$. Le reste ne pose pas de problème.

\emph {2.}
En posant $S_i = \SK(x_{i+1}, \ldots, x_d)$, on a $S_d = 1$ et $S_{i-1} =
x_i^\NN(S_i + \gA x_i)$. De proche en proche on montre $S_i \subseteq \ff_i$,
en utilisant $x_i \in \fp_i$ et $x_i \in f_{i-1}$:

\snic {
S_{i-1} = x_i^\NN(S_i + \gA x_i) \subseteq x_i^\NN(\ff_i + \fp_i) = x_i^\NN \ff_i 
\subseteq x_i^\NN\ff_{i-1} \subseteq \ff_{i-1}
.}

%\sni
Le reste ne pose pas de problème.

%%%%%%%%%%%%%%%%%%%%%%%%%%%%%%%%%%%%%%%%%%%%%%%%%%%%%%%%%%%%%%%%%%%%%%%%%%%

\exer{exoEliminationEtBord}  
Si $f$ est \mon de degré $n \ge 1$, le \pol $R(X,Y)$ de l'énoncé est $Y$-\mon de degré $n$,
donc $\Ann(R) = 0$. 
\\
Et $R(f,g) = 0$ car
$R \in \gen {f(T)-X, Y-g(T)}_{\gA[T,X,Y]}$.

\emph {1.} 
Soit $f = \sum_{k=0}^n a_kT^k$.  D'après le lemme \ref{lemQI} il existe un
\sfio $(t_n, t_{n-1}, \ldots, t_0, t_{-1})$ tel que:
\\
-- dans la composante $t_k = 1$
pour $k \in\lrbzn$, on a $a_i = 0$ pour $i > k$ et $a_k$ \ndz; 
\\
-- dans la
composante $t_{-1} = 1$, on a $f = 0$, i.e. $t_{-1}f = 0$ et même $\Ann(f) =
\gen {t_{-1}}$.  
\\
Soit $m$ le degré formel de $g$.  Pour $1 \le k \le n$, on
pose:

\snic {
R_k(X,Y) = t_k\Res_T(t_kf(T) - X, k, Y-g(T), m).
}

%\sni
On définit $R_0(X,Y) = t_0(t_0f(T)-X)$ et $R_{-1}(X,Y) = t_{-1}X$.  Pour $k
\in \lrb{-1..n}$, on a $\Ann(R_k) = \gen {1-t_k}$ et $R_k(f,g) = 0$. Ainsi en
posant $R = \sum_{k=-1}^n R_k(X,Y)$, on a $\Ann(R) = 0$ et $R(f,g) = 0$.  

\emph {2.}
Application directe de l'exercice référencé.

\emph {3.}
Par \recu sur la dimension de Krull de $\gA$.  On peut remplacer $\gA$ par un anneau $\gA' :=
\gA_{\so\ua}$ de façon à ce que l'annulateur (dans $\gA'$) de chaque
\coe de $f$ soit engendré par un \idm (rappelons que $\Kdim\gA = \Kdim\gA'$). \\
Alors, si $\fa=\prod_{i,j}\IK_\gA(r_{ij})$, 
l'anneau $\gA[T]\sur{\IK_\AT(f,g)}$ est un quotient de
$(\gA\sur\fa)[T]$. Comme $\Kdim(\gA\sur\fa) < \Kdim\gA$, on obtient par \hdrz:

\snic {
\Kdim (\gA\sur\fa)[T] \le 1 + 2\Kdim(\gA\sur\fa) \le
1 + 2(\Kdim\gA - 1)$,\ \  puis $ \
}

\snic {
\Kdim \gA[T] \le 2 + \Kdim \gA[T]\big/{\IK_\AT(f,g)} \le 
2 + 1 + 2(\Kdim\gA - 1) = 1 +  2\Kdim\gA. 
}

%\sni
\emph {4.}
On conserve les notations des questions précédentes.  Chaque $\IK_{\gA'}(r_{ij})$
contient un produit fini d'\ids bord de $\gA$ (exercice \ref
{exoInclusionBordLionel}) donc le produit des $\IK_{\gA'}(r_{ij})$ contient un
\id $\fa$ de $\gA$, produit fini d'\ids bord de $\gA$. 
\\
Ainsi, $\fa \subset
\gA[T] \cap \IK_{\gA'[T]}(f,g) \subseteq \rD_{\gA[T]}(\IK_{\gA[T]}(f,g)\big)$
(exercice~\ref {exoExtEntiereIdealBord}).
%%%%%%%%%%%%%%%%%%%%%%%%%%%%%%%%%%%%%%%%%%%%%%%%%%%%%%%%%%%%%%%%%%%%%%%%%%%

\exer{exoIdealBordPolynomes}
\emph {1.}
Par \recu sur $n$, le cas $n=0$ étant l'hypothèse.  Ajoutons un \elt $b$
à $\bn$ et posons $\fb'_j = \IK(z_j, \bn, b)$.
\\
Par \dfn $\fb'_j = \gB b +
(\fb_j : b^\infty)$ avec $\fb_j = \IK(z_j, \bn)$; le produit des $\fb_j$ est
contenu dans $\IK(x,y,\bn)$ (par \recuz). En utilisant des inclusions du type
$(\fb : b^\infty)(\fb' : b^\infty) \subseteq (\fb\fb' : b^\infty)$, on
obtient:

\snic {\arraycolsep2pt
\begin {array}{rcl}
\prod_j \fb'_j &\subseteq& \gB b + \prod_j (\fb_j : b^\infty)
\subseteq \gB b + \left( \prod_j \fb_j : b^\infty\right)
\\
&\subseteq& \gB b + (\IK(x,y,\bn) : b^\infty) = \IK(x,y,\bn,b)
\\
\end{array}
.}

%\sni
\emph {2a.} Résulte du fait que pour deux \ids $\fa, \fb$ de $\gA$,
on a $(\fa : \fb)_\gA\gA[T] = (\fa : \fb)_{\gA[T]}$. 

\emph {2b.}
Pour deux \ids $\fa, \fb$, notons $\fa \Subset \fb$ pour  $\fa \subseteq\rD(\fb)$.  On raisonne par \recu sur $d$,
le cas $d = 1$ figurant dans l'exercice
\ref{exoEliminationEtBord}. 
\\
Considérons $2(d+1)$
\pols $p$, $q$,  $g_1$, \ldots, $g_{2d} \in \gA[T]$. Il existe des $a_j \in \gA$
tels que $\prod_j \IK_\gA(a_j) \Subset \IK_{\gA[T]}(p,q)$ (le cas $d = 1$).
D'après la première question:

\snic {
\prod_j \IK_{\gA[T]}(a_j, g_1, \ldots, g_{2d}) \Subset 
\IK_{\gA[T]}(p,q, g_1, \ldots, g_{2d})
.}

%\sni
Il suffit donc de montrer, pour $a \in \gA$, qu'un \id bord $\IK_{\gA[T]}(a,
g_1, \ldots, g_{2d})$ contient, à radical près, un produit d'\ids bord de
$d+1$ \elts de $\gA$.  
On pose $\ov\gA= \gA\sur{\IK(a)}$ et $\varphi : \gA[T] \to \ov\gA[T] \simeq
\gA[T]\sur{(\IK_{\gA[T]}(a)\big)}$ l'\homo de passage au quotient. Par \recuz, l'\id bord
$\IK_{\ov\gA[T]}(\ov {g_1}, \ldots, \ov{g_{2d}})$ contient, à radical près,
un produit $\prod_j \fa_j$ où chaque $\fa_j$ est un \id 
bord de $d$ \elts de~$\ov\gA$. En prenant l'image réciproque par $\varphi$, on obtient

\snic {
\prod_i \varphi^{-1}(\fa_i) \subseteq
\varphi^{-1}\big(\prod_i \fa_i\big) \Subset
\varphi^{-1}\big(\IK_{\ov\gA[T]}(\ov {g_1}, \ldots, \ov{g_{2d}})\big)
.}

%\sni
En utilisant le lemme \ref{lem2BordKrullItere}, on a d'une part

\snic {
\varphi^{-1}\big(\IK_{\ov\gA[T]}(\ov {g_1}, \ldots, \ov{g_{2d}})\big) =
\IK_{\gA[T]}(a, g_1, \ldots, g_{2d})
,}

%\sni
et d'autre part $\varphi^{-1}(\fa_i)$ est un \id bord de $d+1$ \elts de
$\gA$ (le premier \elt étant $a$).  Ceci montre que $\IK_{\gA[T]}(a,
g_1, \ldots, g_{2d})$ contient à radical près, un produit d'\ids bord de
$d+1$ \elts de $\gA$.

\emph {3.}
Si $\gA[\uT] = \gA[T_1, \ldots, T_r]$, l'\id bord de $(r+1)d$ \pols de
$\gA[\uT]$ contient, à radical près, un produit d'\ids bord de $d$ \elts
de $\gA$.  
\\
En conséquence, $\Kdim\gA < d \Longrightarrow \Kdim\gA[\uT] < (r+1)d$,
i.e. 

\snic{\Kdim\gA[\uT]+1\, \le \,(r+1)(\Kdim\gA +1) .}

%%%%%%%%%%%%%%%%%%%%%%%%%%%%%%%%%%%%%%%%%%%%%%%%%%%%%%%%%%%%%%%%%%%%%%%%%%%
\exer{exoKdimEspanol}
\emph{1.} 
On prend $a'=y\vi a$ et $b'=x\vu b\vu a'$.  Alors $x \vi a' = x \vi y \vi a =
x \vi a$ (car $x \le y$).  Puis $y \vu b' = y \vu x \vu b \vu a' = (x \vu a')
\vu (y \vu b) = y \vu b$ (la dernière \egt utilise 
$x \vu a' \le y$ qui découle de $x \le y$ et $a' \le y$, a fortiori
$x \vu a' \le y \vu b$).
\\
Reste à voir que $y \vi b' = x \vu a'$; on a l'\idt pour tous $y, b, z$, $y
\vi (b \vu z) = y \vi z'$ avec $z' = (y \vi b) \vu z$ que l'on utilise avec $z
= x \vu a'$.  Mais on a $y \vi b \le x \vu a'$ car l'hypothèse est $y \vi b
\le x \vu a$, donc $y \vi b \le (x \vu a) \vi y = (x \vi y) \vu (y \vi a)
\le x \vu a'$. Donc $z' = x \vu a'$ et $y \vi b' = y \vi (x \vu a')$.
Enfin, $y \vi (x \vu a')=x \vu a'$ car $x \vu a' \le y$
(en utilisant $x \le y$ et $a' \le y$).

\emph{2.} D'après \emph{1} par \recu sur $n$.

\emph{3.} Le point \emph{3a} implique le point \emph{3c} d'après le point \emph{2}. Le point \emph{3c} implique le point~\emph{3b} parce qu'une chaîne liée est un cas particulier de suite \copz. 
Pour voir que \emph{3b} implique  \emph{3a}, \hbox{soit $y_0,\dots,y_{n}$} une suite arbitraire. \\
On définit alors
$x_0=y_0$, $x_i=y_i\vu x_{i-1}$  ($i\in\lrbn$).
Soit $(a_0,\dots,a_{n})$ une suite complémentaire \hbox{de $(x_0,\dots,x_{n})$}.\\
On définit $b_0=a_0$  et $b_{i} = a_{i}\vu  x_{i-1}$  \hbox{pour $i\in\lrbn$}. On a alors $x_{i}\vu  a_i=y_i\vu b_i$ pour $i\in\lrb{0..n}$. Donc $0 = x_0\vi a_0=y_0\vi b_0$ et $1=x_{n}\vu a_{n}=y_n\vu b_n$. Voyons maintenant les inégalités intermédiaires. Pour $i\in\lrb{1..n  }$ 
on a $x_{i}\vi a_{i}\leq  x_{i-1}\vu a_{i-1}$, et donc
\[ 
y_{i} \vi a_{i} \leq  x_{i} \vi a_{i}  \leq  x_{i-1}\vu a_{i-1}= y_{i-1}\vu b_{i-1},
\]
d'où
\[ y_{i}\vi b_{i}=  y_{i}\vi (a_{i}\vu  x_{i-1})=
(y_{i}\vi a_{i})\vu (y_{i}\vi x_{i-1})\leq (y_{i}\vi a_{i})\vu x_{i-1}.
\]
Comme les deux derniers termes après $\leq $ sont majorés par $x_{i-1}\vu a_{i-1}= y_{i-1}\vu b_{i-1}$, on obtient bien l'inégalité $y_{i}\vi b_{i}\leq y_{i-1}\vu b_{i-1}$.

%%%%%%%%%%%%%%%%%%%%%%%%%%%%%%%%%%%%%%%%%%%%%%%%%%%%%%%%%%%%%%%%%%%%%%%%%%%

\exer{exoAminEtagesFinis} 
Tout d'abord, pour tout \id $\fc$, l'anneau $\gA\sur{\fc\epr}$ est  réduit.  \\
Montrons que $(\fa_1\epr\fa_2\epr)\epr
= (\fa_1+\fa_2)^{\perp\perp}$: l'\egt $\fa_1\epr\cap\fa_2\epr =
(\fa_1+\fa_2)\epr$ implique que les  \ids $\fa_1\epr\cap \fa_2\epr$, $\fa_1\epr\fa_2\epr$ et
$(\fa_1+\fa_2)\epr$ ont même racine donc même annulateur. 
\\
On en déduit
que

\snic {
(\fa_1\epr\fa_2\epr\fb)\epr = (\fa_1+\fa_2) \diamond \fb.
}

%\sni
En effet:

\snic {
(\fa_1\epr\fa_2\epr\fb)\epr = \big((\fa_1\epr\fa_2\epr)\epr : \fb\big) =
\big((\fa_1+\fa_2)^{\perp\perp} : \fb\big) = (\fa_1+\fa_2) \diamond \fb.
}

%\sni
\emph {1.}
Comme $\fa\epr\fb \subseteq \fa\epr$, on a $\fa\diamond\fb \supseteq
\fa^{\perp\perp} \supseteq \fa$. Soit $x\in \gA$ tel que dans le quotient 
on ait $\ov
x\,\ov\fb = 0$, i.e.  $x\fb \subseteq \fa\diamond\fb$, i.e.  $x\fb\fa\epr\fb =
0$. On a donc $x\fb\fa\epr = 0$, i.e.  $x\in \fa\diamond\fb$, i.e. $\ov x =
0$.

\emph {2.}
On a:

\snic {
(\gA\sur{\fa_1\diamond\fb_1})\sur{(\ov{\fa_2}\diamond\ov{\fb_2})} \simeq
\gA\sur{(\fa_1\epr\fa_2\epr\fb_1\fb_2)\epr} =
\gA\sur{\big((\fa_1+\fa_2)\diamond(\fb_1\fb_2)\big)}. 
}

%%%%%%%%%%%%%%%%%%%%%%%%%%%%%%%%%%%%%%%%%%%%%%%%%%%%%%%%%%%%%%%%%%%%%%%%%%%

\exer{exoAmin}
\emph{1.} Soit $y\in\gB$ et supposons $y\theta(a)=0$. Posons $e'=\theta(e)$.
On doit montrer
que $y=ye'$. Puisque $e+a$ est \ndzz, $e'+\theta(a)$ est \ndzz.
\hbox{Or $y(e'+\theta(a)\big)=ye'=ye'(e'+\theta(a)\big)$} parce que $e'$ est \idmz.

\emph{2.} L'\homo  $\Aqi \to \Amin$ provient de la \prt \uvle de $\Aqi$.
Il est surjectif parce que  $\Amin=\gA[(e_x)_{x\in\gA}]$ et que le morphisme
$\Aqi \to \Amin$ est \qiz.

\emph{3.} Soit $x$ \ndz dans $\gA$ et 
$u=(\ov y, \wi z)\in\gA_{\so{a}}=\gA\sur{a\epr}\times \gA\sur{({a\epr})\epr}$,
avec $ux=0$. On doit montrer que $u=0$, i.e. $\ov y=\ov 0$
et $\wi z=\wi 0$. 
\\
On a $xy\in a\epr$, \hbox{i.e. $xay=0$}, donc
$ay=0$, puis $\ov y=\ov 0$. 
\\
Pour voir que $\wi z=\wi0$ on considère un \eltz~$t$ 
arbitraire de $a\epr$ et l'on doit montrer que $zt=0$.
Or $\wi{xz}=\wi0$, donc $xzt=0$, puis $zt=0$.

\emph{4.} Si $a\in\gA$ est \ndzz, il reste \ndz aux étages finis de 
la construction de~$\Amin$  d'après le point \emph{3} %(voir  lemme~\ref{lem2qi}) 
et cela suffit pour qu'il soit \ndz dans $\Amin$.
Si l'\homo naturel $\ZZ\to\ZZ_\mathrm{qi}$ était \reg tous les \homos
de $\ZZ$ vers des anneaux \qis seraient \regs vue la \prt \uvle de $\ZZ_\mathrm{qi}$. Or la surjection $\ZZ\to\aqo\ZZ n$ n'est pas un  \homo \reg pour $n\geq2$.
Notez que l'argument s'applique à tout anneau $\gA$ pour lequel il existe un \elt
\ndz $x$ tel que $\aqo\gA x$ est \qi et non trivial.

%%%%%%%%%%%%%%%%%%%%%%%%%%%%%%%%%%%%%%%%%%%%%%%%%%%%%%%%%%%%%%%%%%%%%%%%%%%
\exer{exolemVdimKdim}
\'Ecrivons le calcul pour $n=k=2$.
\\
Soient $x_1={a_1\over b_1}$, $x_2={a_2\over b_2}\in\Frac\gA$ et $s=\big(P(x_1,x_2),(Q(x_1,x_2),(R(x_1,x_2)\big)$ une suite dans $\gA[x_1,x_2]$, avec $P$, $Q$, $R\in\gA[X_1,X_2]$.
On doit montrer que la suite~$s$ est singulière. 
On note $\gA_1=\gA[x_1]$. On sait que la suite

\snic{(b_1X_1-a_1,b_2X_2-a_2,P,Q,R)=(f_1,f_2,P,Q,R)}

%\sni
est \sing dans $\gA[X_1,X_2]$, ce qui donne une \egt

\snic{f_1^{m}(f_2^{m}(P^{m}(Q^{m}(R^{m}(1+AR)+BQ)+CP)+Df_2)+Ef_1)=0}

%\sni 
dans $\gA[X_1,X_2]$. Puisque $b_1\in\Reg\gA$, on a $f_1\in\Reg\gA[X_1,X_2]$ (lemme de McCoy, corolaire \ref{corlemdArtin}). On  simplifie donc l'\egt par $f_1^{m}$, puis on l'évalue dans $\gA_1[X_2]$ par le morphisme $X_1\mapsto x_1$. On obtient l'\egt suivante dans  $\gA_1[X_2]$:

\snic{f_2^{m}(p^{m}(q^{m}(r^{m}(1+ar)+bq)+cp)+df_2)=0,}

%\sni 
avec $p=P(x_1,X_2)$,  $q=Q(x_1,X_2)$, \dots, $d=D(x_1,X_2)$.\\
Puisque $b_2\in\Reg\gA_1$, on a $f_2\in\Reg\gA_1[X_2]$. On peut donc simplifier l'\egt par~$f_2^{m}$, puis l'évaluer dans $\gA[x_1,x_2]$ par le morphisme $X_2\mapsto x_2$.
On obtient une \egt qui dit que la suite $s$ est singulière.

%%%%%%%%%%%%%%%%%%%%%%%%%%%%%%%%%%%%%%%%%%%%%%%%%%%%%%%%%%%%%%%%%%%%%%%%%%%
\exer{exoLyingOverClassique} \emph{(Morphisme lying over)}
Notons \emph{a}, \emph{b} et \emph{c} les trois \prts pour les anneaux commutatifs. L'\eqvc de \emph{a} et \emph{b} est facile.
L'implication \emph{a} $\Rightarrow$
\emph{c} a été donnée en remarque après le lemme lying over \ref{lemLingOver}.

 \emph {c} $\Rightarrow$ \emph {a.}
En \clama $\DA({\fa})$ est l'intersection des \ideps
qui contiennent $\fa$. On veut donc montrer que  pour tout idéal premier $\fp$
tel \hbox{que $\fa\subseteq\fp$}, on a $\varphi^{-1}(\gen {\varphi(\fa)})
\subseteq\fp$. Soit $\fq$ un idéal premier de $\gB$ au dessus de $\fp$, \hbox{i.e.
$\varphi^{-1}(\fq)=\fp$}. Alors, $\gen {\varphi(\fa)} \subseteq
\gen {\varphi(\fp)} \subseteq \fq$, d'où
$\varphi^{-1}(\gen {\varphi(\fa)}) \subseteq\fp$.

%%%%%%%%%%%%%%%%%%%%%%%%%%%%%%%%%%%%%%%%%%%%%%%%%%%%%%%%%%%%%%%%%%%%

\exer{exoLYOV}   \emph{(Morphisme lying over, 2)}\\
Démontrons tout d'abord (en \clamaz) le lemme de Krull pour les \trdisz.
Étant donnés un deux parties $I$ et $F$ de $\gT$, on dit que le couple $(I,F)$ \emph{s'effondre} si $1=0$ dans le treillis quotient $\gT/(I=0,F=1)$. Comme conséquence, il ne peut y avoir de morphisme $\varphi:\gT\to\Deux$ envoyant $I$ sur $0_\gT$ et~$F$ sur $1_\gT$. Le lemme de Krull est la réciproque de la remarque précédente: si le couple ne s'effondre pas, il y a un morphisme $\varphi:\gT\to\Deux$ envoyant $I$ sur $0_\gT$ et $F$ sur $1_\gT$. 
\\
La proposition \ref{propIdealFiltre} nous dit qu'un couple $(J,G)$ s'effondre
\ssi il existe une partie finie $J_0$ de $J$ et une partie finie $G_0$ de $G$ 
telles \hbox{que $\Vi G_0\leq \Vu J_0$}.
Soit $x\in \gT$. On note que si les couples $(J\cup \so x,G)$ et $(J,G\cup \so x)$ s'effondrent alors il existe une partie finie $J_1$ de $J$ et une partie finie $G_1$ de $G$ telles que $$\Vi G_1\vi x\leq \Vu J_1\hbox{ et }\Vi G_1 \leq \Vu J_1\vu x.$$ En notant $g=\Vi G_1$ et $j=\Vu J_1$, on obtient $g\vi x\leq  j$ et
$g\leq  j\vu x$. Donc par coupure (voir \pref{coupure1} \paref{coupure1}), on a $g\leq j$, \hbox{et $(J,G)$} s'effondre également.\\  
On dit qu'un couple $(I_2,F_2)$ raffine un couple $(I_1,F_1)$ si $I_1\subseteq I_2$ et $F_1\subseteq F_2$.
On considère un couple 
$(J,G)$ qui est 
maximal, pour la relation de raffinement, parmi les couples qui raffinent 
$(I,F)$ et qui ne s'effondrent pas (le lemme de Zorn s'applique car l'effondrement ne concerne que des parties finies). Tout d'abord, il est clair que~$J$ est un idéal et $G$ 
un filtre car on ne change pas l'effondrement en remplaçant $J$ par l'idéal qu'il engendre et $G$ par le filtre qu'il engendre.
Si on n'avait pas $J \cup G=\gT$, on aurait un $x\in A\setminus (J\cup G)$ et d'après 
la remarque initiale  l'un des deux
couples $(J\cup \so{x},G)$ et $(J,G\cup\so x)$ doit ne pas s'effondrer. Et ceci contredit la maximalité.

Passons à la démonstration de l'équivalence convoitée.
On a les équivalences:
%-----------------begin $$----------------
$$\quad  a\not= b\quad \Longleftrightarrow\quad  a\vi b\not= a\vu b\quad
\Longleftrightarrow\quad
a\vu b\not\leq a\vi b
.$$
%-----------------end $$------------------
Supposons que $\Spec(\varphi)$ est surjectif. Si $a\not= b$ dans $\gT$, soit
$a'=\varphi(a)$, $b'=\varphi(b)$ et
  $\psi\in\Spec(\gT)$ tel que $\psi(a\vu b)=1$ et
 $\psi(a\vi b)=0$. Puisque  $\Spec(\varphi)$ est surjectif il existe
$\psi'\in\Spec(\gT')$ tel \hbox{que $\psi=\psi'\varphi$} donc  $\psi'(a'\vu b')=1$ 
et
 $\psi'(a'\vi b')=0$, donc $a'\vu b'\not\leq a'\vi b'$ et $ a'\not= b'$.\\
Supposons que  $\varphi$ est injectif. Identifions $\gT$ à un sous-\trdi 
de $\gT'$. \\
Si $\psi\in\Spec(\gT)$, \hbox{soient $I=\psi^{-1}(0)$} et $F=\psi^{-1}(1)$.
Alors $(I,F)$ ne peut pas s'effondrer dans~$\gT'$ car cela le ferait s'effondrer
dans $\gT$. Donc il existe $\psi'\in\Spec(\gT')$ tel que $\psi'(I)= 0$ 
\hbox{et $\psi'(F)= 1$}, ce qui signifie~\hbox{$\psi=\psi'\varphi$}.

%\exer{exodefiGoingup} \emph{(Morphisme going up, going down)}\\

\exer{exodefiGoingup2} \emph{(Morphisme going up, going down, 2)}\\
 Notons que dans le point \emph{3}, $j\in \alpha^{-1}(\dar i)$ signifie  $\alpha(j)\leq_\gV i$.
Le point \emph{4} est donc simplement la traduction du point \emph{3} en tenant compte de la description de l'inégalité~\hbox{$a\leq b$} dans un quotient d'un \trdi $S$ par un idéal prinicpal: $a\vdi{S/(j=0)}b$ équivaut à $a\vdi{S}b,\,j$ (point \emph{3} de la proposition~\ref{propIdealFiltre}).
\\
Le point \emph{3} est un cas particulier du point \emph{2}. Pour l'implication réciproque si \hbox{on a $\alpha(a)\vdi{\gV/(\fa=0)}\alpha(b)$}, il y a un $i\in \fa$ tel que
$\alpha(a)\vdi{\gV/(i=0)}\alpha(b)$. 
\\
Le point \emph{1} est un cas particulier du point \emph{2}.
Il reste à prouver \emph{1} $\Rightarrow$ \emph{4}, ce qui ne peut se faire qu'avec des arguments non \cofsz. 

\smallskip \noindent 
Les \elts $a,b,i$ satisfaisant $\alpha(a)\vdi \gV \alpha(b),\,i$ sont fixés et on cherche un $j$ convenable. Par hypothèse (point \emph{1}) pour tout \idep $\fq$ de $\gV$ contenant $i$, on a un $j_\fq\in \gT$ vérifiant: $a\vdi \gT b,\,j_\fq$ et $\alpha(j_\fq)\in \fq$. Soit $\ff$ le filtre de $\gV$ engendré par les~$\alpha(j_\fq)$. Si $\ff$ contient $i$, on obtient un nombre fini de $j_{\fq_k}$ dont la borne inférieure  $j$ satisfait~\hbox{$\alpha(j)\leq i$} et $a\vdi \gT b,\,j$ (par distributivité): on a gagné.
Si ce n'est pas le cas, une zornette nous donne un filtre $\ffg$ maximal parmi les filtres contenant $\ff$ et ne contenant pas $i$, et ce filtre est premier par construction (comme dans la démonstration du lemme de Krull donnée à l'exercice \ref{exoLYOV}). D'où une contradiction en considérant l'\idep $\gV\setminus \ffg$
qui contient~$i$ mais aucun des $\alpha(j_\fq)$. 

%%%%%%%%%%%%%%%%%%%%%%%%%%%%%%%%%%%%%%%%%%%%%%%%%%%%%%%%%%%%%%%%%%%%%%%%%%%
%:  sols pbs

%%%%%%%%%%%%%%%%%%%%%%%%%%%%%%%%%%%%%%%%%%%%%%%%%%%%%%%%%%%%%%%%%%%%%%%%%%%

\prob{exoAnneauNoetherienReduit}
\emph {1a.}
On a un $a \in \Ann(x)\fa$ non nul, en particulier $ax = 0$. 
\\
Montrons que $a \notin
\rD(x)$: si $a^n \in \gen {x}$, alors $a^{n+1} \in \gen {ax} = 0$, et donc $a
= 0$. 
\\
Donc $\rD(x) \subsetneq \rD(x,a) = \rD(ax, a+x) = \rD(a+x)$: on
prend $x' = a+x$ (qui est bien dans $\fa$).

\emph{1b.}
On pose $x_0 = 0$. Si $\Ann(x_0)\fa = 0$, \cad $\fa=0$, alors $\Ann(x_0) \subseteq \Ann(\fa)$,
\hbox{donc $\Ann(x_0) = \Ann(\fa)$}. Dans ce cas on pose $x_i=x_0$ pour tout $i\geq0$. Sinon, il y a un $x_1 \in \fa$ avec $\rD(x_0)
\subsetneq \rD(x_1)$.  Si $\Ann(x_1)\fa = 0$, alors $\Ann(x_1) \subseteq
\Ann(\fa)$, \hbox{donc $\Ann(x_1) = \Ann(\fa)$}. Dans ce cas on pose $x_i=x_1$ pour tout $i\geq1$. Sinon, il y a un $x_2 \in \fa$ avec
$\rD(x_1) \subsetneq \rD(x_2)$ \ldots\,\ldots
\\
On construit de cette manière une suite infinie
croissante d'\ids $\rD(x_i)$, qui est stationnaire
dès que deux termes consécutifs sont égaux,
auquel cas le \pb initial est résolu\footnote{La preuve de terminaison de l'\algo sous hypothèse \noee \cov
qui vient d'être donnée est un peu déroutante. Spontanément
on aurait préféré dire:
il faut bien que l'\algo termine un jour car sinon, on aurait une suite infinie
strictement croissante. L'ennuyeux dans ce dernier argument est qu'il est un argument par l'absurde.
Ici on a utilisé l'hypothèse \noee sous forme \cov
et cela nous a fourni le moyen de savoir a priori quand l'\algo terminera.
Ce point délicat renvoie à la discussion sur le principe de Markov (annexe
\paref{principeMarkov}).}.

\emph{2.}
Soit $y\in\gA$.
D'après l'hypothèse, on applique le point \emph{1} avec l'\id $\fa=\Ann(y)$
et l'on sait déterminer un
$x \in \Ann(y)$ tel que $\Ann(x) = \Ann(\Ann(y)\big)$, i.e. $xy =
0$ et $\Ann(x)\Ann(y) = 0$.
On a alors $(\Ann(y) \cap \Ann(x)\big)^2 \subseteq
\Ann(x)\Ann(y) = 0$, \hbox{donc $\Ann(x) \cap \Ann(y) = 0$} (l'anneau est
réduit). Montrons que $x+y$ est \ndz; supposons $z(x+y) =
0$. En multipliant par $y$, $zy^2 = 0$, donc $zy = 0$, puis $zx = 0$, donc $z \in
\Ann(x) \cap \Ann(y) = 0$. En conséquence, $x + y$ est \iv et cet \elt est
dans l'\id bord de $y$ puisque $x\in\Ann(y)$.

\emph{3.}
Pour tout anneau $\gC$, tout \elt \ndz de $\Frac(\gC)$ est \ivz.
Nous pouvons appliquer le résultat du point \emph{2} à l'anneau $\gC=\Frac(\gB\red)$.
\\
En effet, la première hypothèse à vérifier est que toute suite croissante d'\ids de la forme $\rD_\gC(x_n/y_n)$ ($x_n\in\gB\red$, $y_n\in\Reg(\gB\red)$) admet deux termes consécutifs égaux.
Or, dans $\gC$ on a l'\egt $\rD_\gC(x_n/y_n)=\rD_\gC(x_n)$, et l'on conclut par le fait que dans $\gB$, la suite croissante $\gen{\xzn}_\gB$ admet deux termes consécutifs égaux.\\
La deuxième hypothèse est que l'on  sache tester, pour $\fraC x u,\fraC y v\in\gC$, 

\snic{\Ann\big(\fraC x u\big)\Ann\big(\fraC y v\big)=0\;?}

%\sni
ce qui est la même chose que 
$\Ann(x)\Ann(y)=0$ dans $\gB\red$. 
Or,  dans $\Zar\gB$, on a l'\egt $\Ann_{\gB\red}(x)=\rD_\gB(x)\im \rD_\gB(0)$,
 et l'on sait que
$\Zar\gB$ est une \agH  discrète (proposition~\ref{propNoetAgH}). 

%%%%%%%%%%%%%%%%%%%%%%%%%%%%%%%%%%%%%%%%%%%%%%%%%%%%%%%%%%%%%%%%%%%%%%%%%%%

%%%%%%%%%%%%%%%%%%%%%%%%%%%%%%%%%%%%%%%%%%%%%%%%%%%%%%%%%%%%%%%%%%%%%%%%%%%

\prob{exoContractedInclusion}
\emph{1.} Soit $\pi : \gB \to \gA$ un $\gA$-\prr d'image $\gA$. \\
Soit
$a \in \fa\gB \cap\gA$, $a = \sum_i a_ib_i$ avec $a_i \in \fa$,
$b_i \in \gB$; donc $a = \pi(a) = \sum_i a_i\pi(b_i) \in \fa$.

 \emph{2.}
Il est clair que $R_G$ est $\gA$-\lin et que $R_G(a) = a$ pour
tout $a \in \gA$. Le reste en découle.

 \emph{3.}
Supposons $\Kdim\gB \le d$ et montrons $\Kdim\gA \le d$. 
\\
Soient $a_0$, \ldots, $a_d \in \gA$; comme $\Kdim\le d$, il existe $n \ge 0$ tels que:
$$\preskip.2em \postskip.2em
(a_0 \ldots a_d)^n \in \gen {c_{d}, c_{d-1}, \cdots, c_0}_\gB
\quad \hbox {avec} \quad
c_i = (a_0 \ldots a_{i-1})^n a_{i}^{n+1}.
$$
Mais $\gen {c_{d}, c_{d-1}, \cdots, c_0}_\gB \cap \gA = \gen {c_{d}, c_{d-1},
\cdots, c_0}_\gA$. Donc $\Kdim\gA \le d$.

 \emph{4.} (\Demo en \clamaz)
\\
On gradue $\gB$ par $\deg X=\deg Y=1$ et $\deg Z = -1$. Alors $\gA$ est la
composante \hmg de degré $0$, donc est facteur direct dans $\gB$.
\\
Soit $\fq' =  \gen{Z}_\gB$  (c'est un \idepz) et $\fp' := \gA \cap \fq' = \gen
{XZ,YZ}$.  \\
On pose $\fp =  \gen{XZ}_\gA$; c'est un \idep avec $\fp \subset
\fp'$ mais il n'existe pas d'\idep $\fq$ de $\gB$ contenu dans
$\fq'$ et au dessus de $\fp$. Ainsi $\gA\subseteq\gB$ n'est pas going~down.
\\
Soit $\fq = \gen {X,Y^2Z-1}_\gB$ (c'est un \idepz) et $\fp := \gA \cap
\fq = \gen {XY}$.  
\\
On pose $\fp' = \gen {XZ, YZ}_\gA$; c'est un \idep avec $\fp \subset \fp'$ mais il n'existe pas d'\idep $\fq'$
de $\gB$ contenant $\fq$ et au dessus de $\fp'$ (un \idep au dessus
de $\fp'$ doit contenir  $Z$, ou $X$ et $Y$). Ainsi $\gA\subseteq\gB$
n'est pas going~up.

%%%%%%%%%%%%%%%%%%%%%%%%%%%%%%%%%%%%%%%%%%%%%%%%%%%%%%%%%%%%%%%%%%%%%%%%%%%
}% fin des solutions d'exos

%:   ---- Section*{references}-----------
\Biblio

Un très bon exposé de la \ddk du point de vue des \clama
se trouve dans \cite{Eis}.

Les espaces spectraux ont été introduits par Stone \cite{Sto} en 1937.
La théorie des espaces spectraux est au c{\oe}ur du livre  \cite{Johnstone}.

 Un \tho important de Hochster \cite{Hoc} affirme que tout espace spectral
est homéomorphe au spectre d'un anneau commutatif.
Une version sans point du \tho de Hochster est: tout \trdi est isomorphe au
treillis de Zariski d'un anneau commutatif (pour une preuve non \cov  
voir~\cite[Banaschewski]{Ba}). Le point délicat est de savoir construire un anneau dont le treillis de Zariski est un ensemble ordonné fini donné.

La \dfn \cov de la \ddk des \trdis et anneaux commutatifs
remonte aux travaux d'André Joyal \cite{Joyal71,Joyal} et de Luis Espa\~nol~\cite{espThesis,esp,esp83,esp86,Es,Espa08}.
L'idée de Joyal était de construire pour chaque entier $\ell\geq1$,
 à partir du \trdi
$\gT$, un \trdi $\gT_\ell$,
qui vérifie une \prt universelle adéquate de façon
que, en \clamaz, les \ideps de $\gT_\ell$ s'identifient aux chaînes
$\fp_0\subseteq\cdots\subseteq\fp_\ell$ d'\ideps de~$\gT$
(les inclusions n'étant pas \ncrt strictes).
Ceci permet ensuite de définir $\Kdim\gT\leq\ell$ au moyen d'une
\prt reliant $\gT_{\ell+1}$ et $\gT_\ell$.
Enfin la \ddk d'un anneau commutatif peut être définie comme
celle de son treillis de Zariski. Pour plus de détails à
ce sujet, voir les \hbox{articles \cite[Coquand\&al.]{cl,cl1,clq2}}.

Le \thrf{thDKA} qui donne  une \carn inductive \elr
de la \ddk d'un anneau commutatif se trouve dans \cite[Coquand\&al.]{clr03}.
La \carn en terme d'\idas donnée dans la proposition \ref{corKrull}
se trouve dans \cite[Coquand\&Lombardi]{cl} et~\cite[Lombardi]{lom}.
\\
Bien que la \carn en terme de suites \cops soit déjà présente
pour les \trdis dans \cite{cl}, elle apparaît pour
les anneaux commutatifs seulement dans \cite[Coquand\&al.]{clq2}.

Des précisions supplémentaires sur le traitement de la \ddk dans les extensions entières se trouvent dans \cite[Coquand\&al.]{CDLQ06}

La dimension valuative des anneaux intègres a été introduite par Jaffard \cite{Jaffard2} (voir aussi \cite[Chap. 5, \S 30]{Gil}) et généralisée aux anneaux commutatifs par Cahen \cite{Cah90}.  Un traitement \cof très élégant est donné dans le cas intègre par
T. Coquand dans~\cite{coqval}. 

Le résultat de l'exercice \ref{exoInclusionBordLionel}
est d\^u à Lionel Ducos. Le \pb \ref{exoAnneauNoetherienReduit} est en rapport direct avec l'article \cite[Coquand\&al.]{cls}.
 Le \pb \ref{exoChainesPotPrem} est tiré des articles \cite[Brenner]{Brenner}, \cite[Coquand\&Lombardi]{cl1} et \cite{lom}.
 Une variante pour les \trdis se trouve dans~\cite{cl}. 

\newpage \thispagestyle{CMcadreseul}
\incrementeexosetprob

%:        %%%%%%%%%%%%%%%%%%%%%%%%%%%%%%%%%%%%
%:        %%%%%%%%%%%%%%%%%%%%%%%%%%%%%%%%%%%%
%---- Chapitre  Nombre de \gtrs d'un module
\chapter{Nombre de \gtrs d'un module}
\label{chapNbGtrs}
%--------------------

\vskip-1em

\minitoc

\subsection*{Introduction}
\addcontentsline{toc}{section}{Introduction}
%-----------------------------------------

Dans ce chapitre  on établit la version \elrz,  non
\noee et \cov de \gui{grands} \thos d'\alg
commutative.

Ces \thosz, dus dans leur version originale
à \KRNz, Bass, Serre, Forster et Swan, concernent
le nombre de \gtrs radicaux d'un \itfz,
le nombre de \gtrs d'un module, la possibilité de produire
un sous-module libre en facteur direct dans un module,
et la possibilité de simplifier des \isosz,
dans le style suivant: si $M\oplus N\simeq M'\oplus N$, alors $M\simeq M'$.

Un progrès décisif a été accompli par Heitmann \cite{Hei84}
qui a montré comment se débarrasser des hypothèses \noeesz.

Un autre progrès décisif a été accompli par T. Coquand
qui a montré dans plusieurs articles comment obtenir tous les résultats
classiques (parfois sous une forme plus forte) au moyen de 
 \dems à la fois \covs et \elrsz.

Les preuves données ici sont essentiellement celles de \cite[Coquand]{Coq3,coq07}
et de \cite[Coquand\&al.]{clq,clq2}.

%--- Sec{Kronecker}  secKroBass
\section[Le \tho de Kronecker et le stable range de Bass]{Le \tho de Kronecker et le stable range de Bass (versions non \noees de Heitmann)}
\label{secKroBass}
%-------------------------------------------

%:--- SUBsection{subsecKro}--- Kronecker ----
\subsec{Le \tho de Kronecker}
\label{subsecKro}
%-----------------------------------------

Le \tho de \KRNz\footnote{Il s'agit d'un autre \tho de Kronecker que celui
donné au chapitre \ref{chapGenerique}.} est usuellement énoncé
sous la forme suivante (\cite{Kronecker}): une
variété \agq dans $\CC^n$ peut toujours être définie par $n+1$
équations.

Pour Kronecker il s'agissait plutôt de remplacer un \sys d'équations
arbitraires dans $\QQXn$ par un \sys \gui{\eqvz} ayant au plus~$n+1$
équations. L'\eqvc de deux \syss vue par Kronecker
se traduit dans le langage actuel par
le fait que les deux \ids ont même nilradical.
C'est en utilisant le \nst que l'on
 obtient la formulation \gui{variétés \agqsz} ci-dessus.

Dans la version démontrée dans cette section, on donne la
formulation à la Kronecker en remplaçant l'anneau $\QQXn$
par un anneau de \ddk $\leq n$ arbitraire.

Le lemme suivant, bien que terriblement anodin, est une clef essentielle.
%: --- Lemma{gcd}----------
\begin{lemma}
\label{gcd}
Pour $u,v\in\gA$ on a
%--------------------begin array---------------
$$\begin{array}{c}
\DA(u,v) =   \DA(u+v,uv)   =   \DA(u+v)\vu\DA(uv)
\end{array}.$$
%---------------------end array--------------
En particulier, si $uv\in\DA(0)$, %\cad $\DA(uv)=0_{\ZarA}$ on obtient
alors $\DA(u,v)=\DA(u+v)$.
\end{lemma}
%--- end-lemma-----------------------------------------
%-----------------begin proof------------------
\begin{proof}
On a évidemment $\gen{u+v,uv}\subseteq\gen{u,v}$,
donc $\DA(u+v,uv) \subseteq\DA(u,v)$. Par ailleurs,
$u^2=(u+v)u-uv\in\gen{u+v,uv}$, donc $u\in\DA(u+v,uv)$.
\end{proof}
%-----------------end proof------------------

Rappelons que deux suites qui vérifient les inégalités
(\iref{eqCG})
dans la proposition~\ref{corKrull} sont dites \copsz.

%--- Lemma{lemKroH}----------
\begin{lemma}
\label{lemKroH}
Soit $\ell\geq 1$.
Si $(b_1,\ldots ,b_\ell)$ et  $(x_1\ldots ,x_\ell)$ sont deux suites
\cops dans $\gA$
alors pour tout $a\in\gA$ on~a:
$$\DA(a,b_1,\dots,b_\ell) = \DA(b_1+ax_1,\dots,b_\ell+ax_\ell),$$
\cad encore: $a\in \DA(b_1+ax_1,\dots,b_\ell+ax_\ell)$.
\end{lemma}
%--- end-lemma---------------------------------
%-----------------begin proof------------------
\begin{proof}
On a les in\egts
$$\arraycolsep2pt
\begin{array}{rcl}
\DA(b_1x_1)& =  &\DA(0)    \\
\DA(b_2x_2)& \leq  & \DA(b_1,x_1)  \\
\vdots~~~~& \vdots  &~~~~  \vdots \\
\DA(b_\ell x_\ell )& \leq  & \DA(b_{\ell -1},x_{\ell -1})  \\
\DA(1)& =  &  \DA(b_\ell,x_\ell ).
\end{array}
$$
On en déduit celles-ci
$$\arraycolsep2pt
\begin{array}{rcl}
\DA(ax_1b_1)& =  &\DA(0)    \\
\DA(ax_2b_2)& \leq  & \DA(ax_1,b_1)  \\
\vdots~~~~& \vdots  &~~~~  \vdots \\
\DA(ax_\ell b_\ell )& \leq  & \DA(ax_{\ell -1},b_{\ell -1})  \\
\DA(a)& \leq   &  \DA(ax_\ell,b_\ell ).
\end{array}
$$
On a donc d'après le lemme \ref{gcd}
$$\arraycolsep2pt
\begin{array}{rcl}
\DA(a)& \leq   & \DA(ax_\ell+b_\ell)\vu \DA(ax_\ell b_\ell)\\
\DA(ax_\ell b_\ell)& \leq &\DA(ax_{\ell-1}+b_{\ell-1})\vu \DA(ax_{\ell-1} b_{\ell-
1}) \\
\vdots~~~~& \vdots  &~~~~  \vdots  \\
\DA(ax_3b_3)& \leq  & \DA(ax_2+b_2)\vu \DA(ax_2 b_2)\\
\DA(ax_2b_2)& \leq  & \DA(ax_1+b_1)\vu \DA(ax_1 b_1) = \DA(ax_1+b_1).
\end{array}
$$
Donc finalement
$$\arraycolsep2pt
\begin{array}{rcl}
\DA(a)&\leq &\DA(ax_1+b_1)\vu \DA(ax_2+b_2)\vu \cdots  \vu \DA(ax_\ell+b_\ell)\\
&= &\DA(ax_1+b_1,ax_2+b_2,\ldots ,ax_\ell+b_\ell).
\end{array}
$$
\end{proof}
%-----------------end proof------------------

%: --- Theorem{thKroH}----------
\begin{theorem}
\label{thKroH} \emph{(\Tho de \KRNz-Heitmann,
avec la dimension de Krull, non \noez)}
\begin{enumerate}
\item Soit $n\geq 0$.
Si $\Kdim \gA <n$ et $b_1$, \dots, $b_n\in\gA$,  il existe
$x_1$, \dots, $x_n$ tels que pour tout $a\in\gA$,
$\DA(a,b_1,\dots,b_n) = \DA(b_1+ax_1,\dots,b_n+ax_n)$.
\item En conséquence, dans un anneau de dimension de Krull $\leq n$, tout
\itf a même nilradical qu'un \id engendré par au plus $n+1$ \eltsz.
\end{enumerate}
\end{theorem}
%--- end-theorem------------------------------

%-----------------begin proof-----------------
\begin{proof}
\emph{1.} Clair d'après le lemme
\ref{lemKroH} et la proposition \ref{corKrull}
(si $n=0$, l'anneau est trivial et $\DA(a)=\DA(\emptyset)$).\\
\emph{2.} Découle de \emph{1} car il suffit
d'itérer le processus. En fait, si $\Kdim \gA\leq n$ 
et~$\fa=\DA(b_1,\ldots,b_{n+r})$  ($r\geq 2$), on obtient en fin de compte

\snic{\fa=\DA(b_1+c_1,\ldots,b_{n+1}+c_{n+1})}

%\sni
avec les
$c_i\in\gen{b_{n+2},\ldots ,b_{n+r} }$.
\end{proof}
%-----------------end proof-------------------

%:--- SUBsec{subsecBass}-- Bass stable range
\subsec{Le \tho \texorpdfstring{\gui{stable range}}{"stable range"} \,de Bass, 1}
\label{subsecBass}
%-----------------------------------------

%: --- Theorem{Bass0}-------------
\begin{theorem}
\label{Bass0} \emph{(\Tho de Bass, avec la dimension de Krull,
sans \noetz)} 
Soit $n\geq 0$. Si $\Kdim \gA <n$, alors $\Bdim\gA<n$. \\
En abrégé: $\Bdim\gA\leq \Kdim \gA$. En particulier, si $\Kdim\gA<n$, tout \Amo stablement libre de rang $\geq n$ est libre (voir le \thref{corBass2}).
\end{theorem}
%--- end-theorem-----------------------------------------

%-----------------begin proof------------------
\begin{proof} Rappelons que  $\Bdim\gA<n$ signifie que pour 
tous $b_1$, \dots, $b_n\in\gA$,
il  existe des~$x_i$ tels que l'implication suivante soit satisfaite:

\snic{\forall
a\in\gA\quad (
1 \in\gen{a,b_1,\dots,b_n } \,\Rightarrow\, 1
\in\gen{b_1+ax_1,\dots,b_n+ax_n}).}

%\sni 
Cela résulte directement du premier point dans le \thoz~\ref{thKroH}.
\end{proof}
%-----------------end proof------------------

\vspace{.7em}

%%: --- Theorem{Bass0}-------------
%\begin{theorem}
%\label{Bass0} \emph{(\Tho de Bass, avec la dimension de Krull,
%sans n{\oe}thé\-rianité)} 
%Soit $n\geq 0$. Si $\Kdim \gA <n$, pour tous $b_1$, \dots, $b_n\in\gA$,
%il  existe $x_1$, \dots, $x_n$ tels l'implication suivante soit satisfaite:
%
%\snic{\forall
%a\in\gA\quad (
%1 \in\gen{a,b_1,\dots,b_n } \,\Rightarrow\, 1
%\in\gen{b_1+ax_1,\dots,b_n+ax_n}).}
%\end{theorem}
%%--- end-theorem-----------------------------------------
%
%%-----------------begin proof------------------
%\begin{proof}
%Cela résulte directement du premier point dans le \thoz~\ref{thKroH}.
%\end{proof}
%%-----------------end proof------------------

%:--- SUBsec{subsecKroloc}
\subsec{Le \tho de Kronecker local}
\label{Kroloc}
%-----------------------------------------

%:     propdef
\begin{propdef}\label{propdefdisjointes}
Dans un anneau on considère deux suites
$(a_0, \dots, a_d)$ et $(x_0, \dots, x_d)$ telles que
$$\arraycolsep2pt\left\{\begin{array}{rcl}
a_0 x_0 & \in & \rD(0) \\
a_1 x_1 & \in & \rD(a_0,x_0) \\
a_2 x_2 & \in & \rD(a_1,x_1) \\
a_3 x_3 & \in & \rD(a_2,x_2) \\
& \vdots & \\
a_d x_d & \in & \rD(a_{d-1},x_{d-1}) 
\end{array}\right.
$$
On dira que ces deux suites sont \ixc{disjointes}{suites ---}.
Alors pour $0 \leq i < d$, on a
\index{suites disjointes}

\snic{
\rD(a_0, \dots, a_i, x_0, \dots, x_i, a_{i+1} x_{i+1})
= \rD(a_0 + x_0, \dots, a_i + x_i).
}
\end{propdef}
%%%%%%%%%%%%%%%%%%%%%%%%%%%%%%%%%%%%%%%%%
\begin{proof}
Une inclusion est évidente. Pour démontrer l'inclusion réciproque, on
utilise les \egts $\rD(a_i,x_i) = \rD(a_i x_i,a_i+x_i)$.

Il vient alors successivement
{%\small
\vspace{1mm}\mathrigid2mu
$$ \hspace{0mm}
{\arraycolsep1pt %\parsep1pt
\begin{array}{rcccccl}
a_0 x_0 & \in & \rD(0) & = & \rD( ) & =  & \rD()\\
& & & & & & ~\VRT{~\supseteq~~} \\
a_0, x_0, a_1 x_1 & \in & \rD(a_0,x_0) & = & \rD(a_0 x_0, a_0 + x_0)
         & =  & \rD(a_0 + x_0)\\
& & & & & & ~~~~~~\VRT{~\supseteq~~} \\
a_1, x_1, a_2 x_2 & \in & \rD(a_1,x_1) & = & \rD(a_1 x_1, a_1+x_1)
        & \subseteq & \rD(a_0 + x_0, a_1 + x_1)\\
& & & & & & ~~~~~~~~~~~\VRT{~\supseteq~~} \\
a_2, x_2, a_3 x_3 & \in & \rD(a_2,x_2) & = & \rD(a_2 x_2, a_2+x_2)
        & \subseteq & \rD(a_0+x_0, a_1 + x_1, a_2 + x_2)\\
\vdots~~~~~~& \vdots & \vdots& \vdots &\vdots & \vdots & ~~~~~~~~~~~\vdots \\
a_i, x_i, a_{i+1} x_{i+1} & \in
& \rD(a_i,x_i) & = & \rD(a_i x_i, a_i+x_i)
       & \subseteq & \rD(a_0 + x_0, \dots, a_i + x_i).
\end{array}
}
$$
}
\end{proof}
%%%%%%%%%%%%%%%%%%%%%%%%%%%%%%%%%%%%%%%%%

Notons que les suites $(a_0, \dots, a_d)$ et $(x_0, \dots, x_d)$
sont \cops \ssi elles sont disjointes et $1\in\gen{a_d,x_d}$.

%:     theorem{thKroLoc}
\begin{theorem}\label{thKroLoc}
Soit $\gA$ un \alo \dcd \ddi$d$,  de radical de Jacobson $\fm$.
On suppose que~$\fm$ est \emph{radicalement \tfz},
i.e., qu'il existe $z_1$, \dots,$ z_n\in\gA$ tels \linebreak 
que $\fm=\DA(z_1, \dots, z_n)$.
Alors $\fm$ est radicalement engendré par $d$~\eltsz.%
\index{radicalement \tfz!ideal@\id ---}\index{ideal@idéal!radicalement \tfz}
\end{theorem}
%%%%%%%%%%%%%%%%%%%%%%%%%%%%%%%%%%%%%%%%%
\begin{proof}
Puisque $\Kdim \gA\leq d$ et $\fm$ est radicalement \tfz,  
le \tho de \KRN \ref{thKroH} nous dit que~$\fm= \rD(x_0,\dots,x_d)$.
En outre, il existe une \linebreak 
suite $(\ua) = (a_0,\dots,a_d)$ \cop
de $(\ux) = (x_0,\dots,x_d)$.
En particulier (suites disjointes), pour tout $i \leq d$, on a

\snic{
\rD(a_0, \dots, a_{i-1}, x_0, \dots, x_{i-1}, a_i x_i)
= \rD(a_0 + x_0, \dots, a_{i-1} + x_{i-1})
,}

%\sni
mais aussi (suites \copsz)
$
1 \in \gen{a_d, x_d}
$.
Ceci montre que $a_d$ est \iv puisque
$x_d \in \fm$. Soit $i$ le plus petit indice tel que $a_i$ soit \iv
(ici on utilise l'hypothèse que $\fm$ est détachable).\\
Il vient alors $a_0$, \dots, $a_{i-1} \in \fm$, puis

\snic{
\rD(x_0, \dots, x_{i-1}, x_i)
\subseteq \rD(a_0 + x_0, \dots, a_{i-1} + x_{i-1})
\subseteq \fm,}

%\sni
et enfin

\snic{
\begin{array}{c}
 \fm = \rD(x_0, \dots, x_{i-1}, x_i, x_{i+1},\dots, x_d)
\subseteq \hspace*{5cm}    \\[1mm]
\hspace*{4cm}%\qquad\qquad\qquad
\rD(\underbrace{a_0 + x_0, \dots, a_{i-1} + x_{i-1}, x_{i+1},
\dots, x_d}_{d \hbox{\scriptsize ~\eltsz}})
\subseteq \fm.
\end{array}
}
\end{proof}
%%%%%%%%%%%%   end{proof}  %%%%%%%%%%%%%%

\rem Pour une \gnn voir les exercices \ref{exoKroLocvar} et \ref{exoKroLocvarbis}.
\eoe

%--- SUBsec{subsecBassHeitmann}-- Bass stable range
\section{Dimension de Heitmann et \tho de~Bass}
\label{subsecDimHeit}
%-----------------------------------------

Nous allons introduire une nouvelle dimension, que nous appellerons la  dimension de Heitmann d'un anneau commutatif. Sa \dfn sera calquée sur la \dfn inductive de la dimension de~Krull, et nous la noterons $\Hdim$.
Auparavant, nous introduisons la dimension $\Jdim$ définie par Heitmann. 

%:  --- definotation{notaJA}--------------
\begin{definota}
\label{notaJA}\label{defHeit}~
\begin{enumerate}

\item [--] Si $\fa$ est un \id de $\gA$ on note $\JA(\fa)$ son
\ixc{radical de Jacobson}{d'un idéal}, \cad l'image réciproque
de $\Rad(\gA\sur\fa)$ par la projection canonique $\gA\to\gA\sur\fa$.%
\index{radical!de Jacobson}\index{Jacobson!radical de ---}

\item [--] Si
$\fa=\gen{\xn}$ on notera $\JA(\xn)$ pour $\JA(\fa)$. En
particulier,~$\JA(0)=\Rad \gA.$

\item [--] On note $\HeA$ l'ensemble des \ids  $\JA(\xn)$.
On l'appelle le \ixx{treillis}{de Heitmann} de l'anneau~$\gA$.

\item [--] On définit $\Jdim\gA$ comme égale à $\Kdim(\Heit\gA)$.
\end{enumerate}
\end{definota}
%--- end-notation-----------------------------------------

On a donc   $x\in\JA(\fa)$ \ssi pour tout $y\in\gA$, $1+xy$ est \iv
modulo~$\fa$. Autrement dit encore

\snic{x\in\JA(\fa)\iff 1+x\gA\subseteq \sat{(1+\fa)},}

%\sni
et $\JA(\fa)$ est le plus grand \id $\fb$ tel que  $1+\fb\subseteq \sat{(1+\fa)}$. \\
On a donc $\sat{\big(1+\JA(\fa)\big)}=\sat{(1+\fa)}$ et $\JA\big(\JA(\fa)\big)=\JA(\fa)$. \\
En particulier $\JA\big(\JA(0)\big)=\JA(0)$ et l'anneau $\gA/\Rad\gA$
a son radical de Jacobson réduit à $0$.

%:  --- Lemma{gcd2}----------
\begin{lemma}
\label{gcd2} ~
%-----------------begin enum------------------
\begin{enumerate}
\item Pour un \id  arbitraire $\fa$ on a  $\JA(\fa)=\JA\big(\DA(\fa)\big)=\JA\big(\JA(\fa)\big)$.\\
En conséquence,  $\HeA$ est un \trdi quotient de $\ZarA$.
\item Pour $u$, $v\in\gA$ on a

\snic{\JA(u,v)\; =  \; \JA(u+v,uv)  \;=\;  \JA(u+v)\vu\JA(uv).}

%\sni
En particulier, si $uv\in\JA(0)$,
alors $\JA(u,v)=\JA(u+v)$.
\end{enumerate}
%-----------------end enum------------------
\end{lemma}
%--- end-lemma-----------------------------------------
%-----------------begin proof------------------
\begin{proof}
On a $\fa\subseteq\DA(\fa)\subseteq\JA(\fa)$, donc $\JA(\fa)=\JA\big(\DA(\fa)\big)=\JA\big(\JA(\fa)\big)$. 
\\
L'\egt
$\DA(u,v)= \DA(u+v,uv)$ implique donc  $\JA(u,v) =\JA(u+v,uv)$.
\end{proof}
%-----------------end proof------------------

\rdb
\comm \label{NOTAJdim}
La $\Jdim$ introduite par Heitmann dans \cite{Hei84} correspond à l'espace
spectral $\Jspec \gA$ suivant: c'est le plus petit
sous-espace spectral de~$\Spec \gA$
contenant l'ensemble $\Max \gA$ des \idemas de $\gA$.
Cet espace peut être décrit comme l'adhérence de $\Max \gA$ dans
$\Spec \gA$
pour la topologie constructible, topologie ayant pour \sys \gtr d'ouverts les $\fD_\gA(a)$ et leurs \cops $\fV_\gA(a)$.
Heitmann remarque que la dimension utilisée dans le \tho de Swan ou dans le
\SSOz, à savoir la dimension de $\Max \gA$,
ne fonctionne bien que dans le cas où cet espace est \noez. En outre,
 dans ce cas, la dimension de~$\Max \gA$  est celle d'un espace spectral, l'espace $\jspec \gA$ formé par
les \ideps qui sont intersections d'\idemasz. Par contre, dans le cas \gnlz,
le sous-espace $\jspec \gA$ de $\Spec\gA$ n'est plus spectral, et donc, selon une remarque qu'il qualifie de philosophique,
 $\jspec \gA$ doit être remplacé par  l'espace
spectral qui s'offre naturellement comme solution de rechange, 
à savoir~$\Jspec \gA$. En fait, 
$\Jspec \gA$ s'identifie au spectre du \trdiz~$\HeA$
(voir \cite[\Tho 2.11]{clq2}).
Et les \oqcs de~$\Jspec \gA$ forment un treillis quotient de~$\ZarA$,
canoniquement isomorphe à~$\HeA$.
En \comaz, on définit donc $\Jdim\gA$ comme égale \hbox{à $\Kdim(\Heit\gA)$}. 
\eoe

\medskip La \dfn de la dimension de Heitmann qui est
donnée ci-après est assez naturelle, dans la mesure où elle mime 
la \dfn \cov de la dimension de Krull en remplaçant $\DA$ par $\JA$. 

%--- Definition{defHei2}--------
\begin{definition}
\label{defHei2} Soit  $\gA$  un anneau commutatif, $x\in\gA$ et $\fa$ un idéal \linebreak 
\tfz.
Le \emph{bord de Heitmann de $\fa$ dans $\gA$} est l'anneau quo-\linebreak 
tient
$\gA_\rH^\fa:=\gA\sur{\JH_\gA(\fa)}$ avec
         $$\JH_\gA(\fa):=\fa+(\JA(0):\fa).$$
Cet \id est appelé le \emph{l'\id bord de Heitmann de $\fa$ dans $\gA$}.
\\
On notera aussi  
$\JH_\gA(x):=\JH_\gA(x\gA)$ et $\gA_\rH^x:=\gA\sur{\JH_\gA(x)}$.%
\index{bord de Heitmann!anneau quotient, idéal}\index{ideal@idéal!bord de
Heitmann}
%-----------------end item------------------
\end{definition}
%--- end-definition-------------------------

%--- Definition{defDHA}--------------
\begin{definition}
\label{defDHA} La \ix{dimension de Heitmann} de $\gA$ est définie par \recu
comme suit:
%-----------------begin item------------------
\begin{enumerate}
\item [--] $\Hdim \gA=-1$ \ssi $1_\gA=0_\gA$.
\item [--] Soit $\ell\geq 0$,  on a l'\eqvcz:

\snic{\Hdim \gA\leq \ell$  $\;\Longleftrightarrow\;$ pour tout $x\in\gA$,
$\Hdim(\gA_\rH^x)\leq \ell-1.}

\end{enumerate}
\end{definition}
%--- end-definition-----------------------------------------

Cette dimension  est inférieure ou égale à la $\Jdim$ définie par Heitmann
dans~\cite{Hei84}, \cad la \ddk du \trdi $\Heit(\gA)$.

%--- Fact{factKdimHdim}--------------
\begin{fact}
\label{factKdimHdim} ~
%-----------------begin enum------------------
\begin{enumerate}
\item La dimension de Heitmann ne peut que diminuer par passage  à
un anneau quotient.
\item La dimension de Heitmann est toujours inférieure ou égale
à la dimension de Krull.
\item Plus précisément $\Hdim \gA\leq \Kdim\big(\gA/\JA(0)\big)\leq \Kdim\gA.$
\item Enfin $\Hdim \gA\leq0$ \ssi $\Kdim\big(\gA/\JA(0)\big)\leq 0$
(i.e.,  $\gA$ est \plcz).
\end{enumerate}
%-----------------end enum------------------
\end{fact}
%--- end-fact-----------------------------------------
%-----------------begin proof------------------
\begin{proof}
\emph{1.} Par \recu sur $\Hdim\gA$(\footnote{En fait par \recu sur $n$ si $\Hdim\gA\leq n$.}) en remarquant que pour tout $x\in\gA$, l'anneau
$(\gA/\fa)_\rH^x$ est un quotient de~$\gA_\rH^x$.
\\
\emph{2.} Par \recu sur $\Kdim\gA$ (en utilisant  \emph{1}) en remarquant que
$\gA_\rH^x$ est un quotient de~$\gA_\rK^x$.
\\ 
\emph{3} et \emph{4.} Soit $\gB=\gA/\JA(0)$. Alors $\rJ_\gB(0)=\gen{0}$, et l'on a
 $\gA_\rH^x\simeq \gB_\rH^{\ov x}=\gB_\rK^{\ov x}$ pour tout~$x\in\gA$.
\end{proof}
%-----------------end proof------------------

%:--- SUBsec{subsecBass}-- Bass stable range
\subsec{Le \tho \texorpdfstring{\gui{stable range}}{"stable range"} \,de Bass, 2}

%:HHH Enonce legerement modifie
%:  --- Theorem{Bass}-------------
\begin{theorem}
\label{Bass} \emph{(\Tho de Bass, avec la dimension de Heitmann, sans
\noetz)} 
Soit $n\geq 0$. Si $\Hdim \gA <n$, alors $\Bdim\gA<n$. 
\\
En abrégé: $\Bdim\gA\leq \Hdim \gA$. En particulier si $\Hdim\gA<n$, tout \Amo stablement libre de rang $\geq n$ est libre.
\end{theorem}
%--- end-theorem-----------------------------------------

%%:  --- Theorem{Bass}-------------
%\begin{theorem}
%\label{Bass} \emph{(\Tho de Bass, avec la dimension de Heitmann, sans
%\noetz)}\\
%Soit $n\geq 0$. Si $\Hdim \gA <n$ et $1 \in\gen{a,b_1,\dots,b_n } $ alors il
%existe $x_1$, \dots, $x_n$ tels que $1  \in\gen{b_1+ax_1,\dots,b_n+ax_n}$.
%\end{theorem}
%%--- end-theorem-----------------------------------------

%-----------------begin proof------------------
\begin{proof}
La même \dem donnerait le \thref{Bass0}
en remplaçant le bord  de Heitmann par le bord de Krull.
%:HHH rajout rappel sur la Bdim
Rappelons que  $\Bdim\gA<n$ signifie que pour 
tous $b_1$, \dots, $b_n\in\gA$,
il  existe des~$x_i$ tels l'implication suivante soit satisfaite:

\snic{\forall
a\in\gA\quad (
1 \in\gen{a,b_1,\dots,b_n } \,\Rightarrow\, 1
\in\gen{b_1+ax_1,\dots,b_n+ax_n}).}
 
%\sni
On rappelle que  $1\in\gen{L}$ équivaut  
à $1\in\JA(L)$ pour toute liste $L$.
On raisonne par \recu sur $n$.\\
 Lorsque $n=0$ l'anneau est trivial et $\JA(1)=\JA(\emptyset)$.\\
Supposons  $n\geq 1$. Posons $\fj=\JH_\gA(b_n)$.
L'\hdr nous donne
$x_1$, \ldots, $x_{n-1}\in\gA$ tels que
\begin{equation}
\label{eqBass1}
1\in \gen{b_1+x_1a,\ldots , b_{n-1}+ x_{n-1} a} \quad \mathrm{dans} \quad
\gA/\fj .
\end{equation}
Notons  $L=(b_1+x_1a,\ldots ,b_{n-1}+ x_{n-1} a)$. Un \elt
arbitraire de $\fj$ s'écrit sous la forme $b_{n}y+x$ avec $xb_n\in\JA(0)$. L'appartenance (\ref{eqBass1}) signifie donc
qu'il existe un~$x_n$ tel que
$$
x_nb_n\in\JA(0)\quad  \mathrm{et}\quad
1\in \gen{L,b_n,x_n}.
$$
A fortiori
\begin{equation}
\label{eqBass3}
 1\in \JA(L,b_n,x_n)=\JA(L,b_n)\vu\JA(x_n).
\end{equation}
Notons que par hypothèse $1\in\gen{a,b_1,\dots,b_n}=\gen{L,b_n,a}$.
Donc
\begin{equation}
\label{eqBass4}
1\in\JA(L,b_n,a)=\JA(L,b_n)\vu\JA(a).
\end{equation}
Comme $\JA(x_{n}a)=\JA(a)\vi\JA(x_{n})$,  (\ref{eqBass3}) et (\ref{eqBass4}) donnent par distributivité
$$
1\in \JA(L,b_n)\vu\JA(x_{n}a)=\JA(L,b_n,x_na)
.$$
Puisque $b_n\,x_n\,a\in\JA(0)$,  le lemme \ref{gcd2} donne
$\JA(b_n,x_na)=\JA(b_n+x_na)$, et donc 
$$
1\in\JA(L,b_n+x_na)=\JA(L,b_n,x_na),
$$
ce qui était le but recherché.
\end{proof}
%-----------------end proof------------------

%:HHH partie sur les stablement libres supprimée

%%%%%%%%%%%%%%%%%%%%%%%%%%%%%%%%%%%%%%%%%%%%%%%%%%%%%%%%%%%%%%%%%%%%%%%%%%%
%%%%%%%%%%%%%%%%%%%       Variante Kronecker          %%%%%%%%%%%%%%%%%%%%%
%%%%%%%%%%%%%%%%%%%%%%%%%%%%%%%%%%%%%%%%%%%%%%%%%%%%%%%%%%%%%%%%%%%%%%%%%%%

%:--- SUBsec{Variante Kronecker
\subsec{Variante \texorpdfstring{\gui{améliorée}}{"améliorée"} \,du \tho de Kronecker}

%--- Lemma{HenriLemma}----------
%\goodbreak
\begin{lemma}\label{HenriLemma} ~\\
Soient $a$, $c_1$, \ldots, $c_m$, $u$, $v$, $w \in \gA$ et notons $Z=(c_1,\ldots,c_m)$.
%-----------------begin enum------------------
\begin{enumerate}
\item
 $a \in \DA(Z)  \iff  1 \in \gen {Z}_{\gA[a^{-1}]}$.
\item 
$\big(w\in \Rad(\gA[a^{-1}])$ et $a\in \DA(Z,w)\big)$
$\Longrightarrow$ $a\in \DA(Z)$.
\item 
$\big(uv\in \Rad(\gA[a^{-1}])$ et $a\in \DA(Z,u,v)\big)$ $\Longrightarrow$ $a\in \DA(Z,u+v).$ 
\end{enumerate}
%-----------------end enum------------------
\end{lemma}
%--- end-lemma-----------------------------------------

%-----------------begin proof------------------
\begin{proof}
\emph{1.} Immédiat.

\emph{2.} Supposons $a\in \DA(Z,w)$ et travaillons dans l'anneau
$\gA[a^{-1}]$. \\
On a $1 \in \gen {Z}_{\gA[a^{-1}]} +
\gen{w}_{\gA[a^{-1}]}$, et comme $w$ est dans 
$\Rad(\gA[a^{-1}])$, cela implique que $1 \in \gen {Z}_{\gA[a^{-1}]}$, i.e.
$a\in \DA(Z)$.

\emph{3.}
Résulte du point \emph{2} car  $\DA(Z,u,v)=\DA(Z,u+v,uv)$ (lemme \ref{gcd}).
\end{proof}
%-----------------end proof------------------

\rem On peut se demander si l'idéal $\Rad\gA[a^{-1}]$ est le meilleur
possible. La réponse est oui. L'implication du point \emph {2} est
vérifiée (pour tout $Z$) en remplaçant $\Rad\gA[a^{-1}]$
par $\fJ$ \ssi $\fJ\subseteq \Rad\gA[a^{-1}]$. \eoe

%--- Lemma{thCor2.2Heit}-----------
\begin{lemma}
\label{thCor2.2Heit}\relax ~\\
Supposons que $\Hdim\gA [1/a]< n$, $L\in \Ae n$ et $\DA(b)\leq\DA(a)\leq \DA(b, L)$.
Alors il existe~$X\in \Ae n$ tel que $\DA(L+bX)=\DA(b, L)$, ce qui équivaut
à~$b \in\DA(L+bX)$, ou encore à $\,a \in\DA(L+bX)$.
En outre, nous pouvons prendre $X=aY$ avec~$Y\in \Ae n$.
\end{lemma}
%--- end-lemma-----------------------------------------
%-----------------begin proof------------------
\begin{proof}
\emph{Remarque préliminaire.}
Si $\DA(L+bX)=\DA(b, L)$, on a 
$a\in\DA(L+bX)$ car $a\in \DA(b,L)$. 
Réciproquement, si $a\in\DA(L+bX)$, on a $b \in
\DA(L+bX)$ (puisque $b\in\DA(a)$),
donc $\DA(L+bX)=\DA(b, L)$.

On raisonne par \recu sur $n$. Le cas $n=0$ est trivial.  
\\
On pose $L=(b_1,\ldots,b_n)$ et on commence par
chercher $X\in \Ae n$.  \\
Soient $\fj=\JH_{\gA[1/a]}(b_n)$ et
$\gA'=\gA/(\fj\cap\gA)$, où $\fj\,\cap\,\gA$ est mis pour
\gui{l'image réciproque de $\fj$ dans $\gA$}.
On a une identification $\gA[1/a]/\fj=\gA'[1/a]$.
\\
Comme $\Hdim\gA'[1/a]< n-1$, on peut appliquer l'hypothèse de \recu à
$\gA'$ et $(a,b,b_1,\ldots ,b_{n-1})$, en remarquant que $b_n = 0$ dans $\gA'$. 
On obtient alors $x_1$, \dots, $x_{n-1}\in \gA$ tels que,  en notant~$Z=(b_1+bx_1,\ldots, b_{n-1}+ bx_{n-1})$,  on ait $\rD(Z)=\rD(b, b_1,\dots,b_{n-1})$ dans $\gA'$. D'après la
remarque préliminaire, cette dernière \egt équivaut à
$a \in \rD_{\gA'}(Z)$, ce qui, d'après le lemme \ref{HenriLemma}~\emph{1}, signifie $1 \in \gen{Z}$ dans $\gA'[1/a]$,
i.e. $1 \in \gen{Z} + \fj$ dans $\gA[1/a]$. Par
\dfn du bord de Heitmann, cela veut dire qu'il existe $x_n$, que l'on peut
choisir dans $\gA$, tel que $x_nb_n\in\Rad\gA[1/a]$ et $1 \in \gen{Z,
b_n,x_n}_{\gA[1/a]}$. \\
 On a donc $a\in\DA(Z,b_n,x_n)$. Mais on a aussi $a\in\DA(Z,b_n,b)$, puisque

\snic {
\gen {Z, b_n, b} = \gen {b_1,, \ldots, b_{n-1}, b_n, b} 
\eqdefi {\gen{L, b}},
}

%\sni
et que $a \in \DA(L,b)$ par hypothèse.
Bilan: $a \in\DA(Z,b_n,x_n)$, $a \in\DA(Z,b_n,b)$ donc $a \in\DA(Z,b_n,bx_n)$.
L'application du lemme~\ref{HenriLemma}~\emph{3}  avec
$u = b_n$, $v = bx_n$ fournit $a \in\DA(Z,b_n+bx_n)$, i.e.  $a \in \DA(L +
bX)$ où $X = (x_1, \ldots, x_n)$.
\\
Enfin, si $b^p\in\gen{a}_\gA$, nous pouvons appliquer le résultat avec $b^{p+1}$ à la place de $b$
puisque $\DA(b)=\DA(b^{p+1})$. Alors $L+b^{p+1}X$ se réécrit $L+baY$.
\end{proof}
%-----------------end proof------------------

Pour $a\in\gA$, on a toujours  $\Hdim\gA[1/a]\leq 
\Kdim\gA[1/a]\leq
\Kdim\gA$. Par conséquent le \tho suivant améliore le \tho de 
\KRNz.

%--- Theorem{KroH2}-------------
\begin{theorem}
\label{KroH2} \emph{(\Tho de \KRNz,  dimension de Heitmann)}
\begin{enumerate}
\item Soit $n\geq 0$. Si  $a$, $b_1$, $\dots$, $b_n\in\gA$ et $\Hdim\gA[a^{-1}] <n$,
alors il existe~$x_1$,~$\dots$,~$x_n\in \gA$ tels que

\snic{\DA(a,b_1,\dots,b_n) =\DA(b_1+ax_1,\dots,b_n+ax_n).}
\item En conséquence, si $a_1$, \ldots, $a_r$, $b_1$, \dots, $b_n\in\gA$ et
$\Hdim\gA[1/a_i] <n$ pour $i\in\lrbr$, alors  il existe $y_1$, \dots, $y_n\in
\gen{a_1,\ldots,a_r}$ tels que

\snic{\DA(a_1,\ldots ,a_r,b_1,\dots,b_n) =\DA(b_1+y_1,\dots,b_n+y_n).}
\end{enumerate}
\end{theorem}
%--- end-theorem-----------------------------------------
%-----------------begin proof------------------
\begin{proof}
\emph{1.} Conséquence directe du lemme \ref{thCor2.2Heit} en faisant $a=b$. 

 \emph{2.} Se déduit de \emph{1} par \recu sur $r$: 

\snac{
\begin{array}{cclcl} 
  \fa&=&\DA(a_1,\ldots,a_r,b_1,\dots,b_n)& =  & \DA(a_1,\ldots ,a_{r-1},b_1,\dots,b_n) \vu
\DA(a_r)  \\[1mm] 
  & =  &  \fb\vu\DA(a_r), \quad \hbox{avec}&& \\[1mm] 
\fb  & =  & \DA(b_1+z_1,\dots,b_n+z_n) 
\end{array}
}

%\sni
où $z_1$, \dots, $z_n\in \gen{a_1,\ldots ,a_{r-1}}$, donc
$\fa=\DA(a_r,b_1+z_1,\dots,b_n+z_n)$, et l'on applique une nouvelle fois le
résultat.
\end{proof}
%-----------------end proof------------------

%--- Section{Splitting off et Forster-Swan}--------------
\section[Splitting-off et Forster-Swan]{Le Splitting-off de Serre, \\ le \tho de Forster-Swan, et \\
le \tho de simplification de Bass}
%-----------------------------------------
\label{secSOSFSa}

Dans cette section, nous expliquons quelles sont les \prts 
matricielles d'un anneau qui permettent de faire fonctionner
les \thos de Serre (Splitting-off) et de Forster-Swan 
(contrôle du nombre de \gtrs d'un \mtf en fonction du nombre
de \gtrs local).

Les sections suivantes consisteront à développer des résultats 
qui montrent que certains anneaux satisfont les \prts matricielles en question.
Les premiers anneaux qui sont apparus (gr\^ace à Serre et Forster)
étaient les anneaux \noes avec certaines \prts de dimension (la \ddk pour Forster et la dimension du spectre maximal pour Serre et Swan).
Plus tard Heitmann a montré comment se débarrasser de la \noet concernant la \ddkz, et a donné les idées directrices pour faire le même travail
avec la dimension du spectre maximal. 
En outre Bass a aussi introduit une \gnn dans laquelle il remplaçait
la \ddk par le maximum des dimensions de Krull pour les anneaux associés à une partition du spectre de Zariski en sous-ensembles constructibles.
Enfin, Coquand apporta une lumière \gui{définitive} sur ces questions
en généralisant les résultats et en les traitant de manière \cov gr\^ace à deux notions sous-jacentes aux preuves antérieures: la $n$-stabilité d'une part, la dimension de Heitmann d'autre part.
L'aspect purement matriciel des \pbs à résoudre a été clairement mis en évidence dans un article de synthèse d'Eisenbud-Evans. La section ici peut être considérée comme une approche non \noee et \cov de ces derniers travaux.

%:     Definition{defiSdimGdim}
\begin{definition}\label{defiSdimGdim}
Soit un anneau $\gA$ et un entier $n\geq 0$.
\begin{enumerate}
\item On écrit $\Sdim\gA< n$ si, pour toute matrice $F$ de rang $\geq n$,
il y a une \coli des colonnes qui est \umdz.\\
Autrement dit $1=\cD_n(F)\Rightarrow \exists X,\, 1=\cD_1(FX)$
\item On écrit $\Gdim\gA< n$
lorsque la \prt suivante est satisfaite. Pour toute matrice 
$F=[\,C_0\,|\,C_1\,|\,\dots\,|\,C_p\,]$  (les~$C_i$ 
sont les colonnes, et on note
$G=[\,C_1\,|\,\dots\,|\,C_p\,]$)
telle que $1=\cD_1(C_0)+\cD_n(G)$,
il y a une \coli des colonnes, avec le premier \coe égal à $1$, qui est \umdz. 
\end{enumerate}
\end{definition}
%--------- fin definition ---------------

Dans l'acronyme $\Sdim$, $\mathsf{S}$ fait allusion à \gui{splitting} ou à \gui{Serre} et est justifié par le \thref{thSerre}.
De même,  dans $\Gdim$, $\mathsf{G}$ fait allusion à \gui{générateurs} et est justifié par le \thref{thSwan}.

Les notations $\Sdim\gA<n$ et  $\Gdim\gA<n$ sont justifiées par les implications évidentes, pour tout $n\geq 0$, 

\snic{\Sdim\gA<n\Rightarrow\Sdim\gA<n+1$ et $\Gdim\gA<n\Rightarrow\Gdim\gA<n+1.}

Notez que $\cD_n(F) \subseteq \cD_1(C_0) + \cD_n(G)$, et par conséquent
 l'hypothèse pour~$F$ dans $\Sdim\gA< n$ implique 
l'hypothèse pour~$F$ dans $\Gdim\gA< n$.
Par ailleurs la conclusion dans $\Gdim\gA< n$ est plus forte.
Cela donne le point \emph{2} qui suit.

%:     Fact{factSDimGdimHdim}
\begin{fact}\label{factSDimGdimHdim}~
\begin{enumerate}
% 1
\item  $\Sdim\gA< 0\iff\Gdim\gA< 0\iff$ l'anneau $\gA$ est trivial.
% 2
\item Pour tout $n\geq 0$, on a $\;\Gdim\gA< n\,\Longrightarrow\,\Sdim\gA< n$.
\\
On note en abrégé $\Sdim\gA\leq \Gdim\gA$.
% 2
\item Si $\gB=\gA/\fa$, on a 
$\,\Sdim\gB\leq \Sdim\gA$ 
et $\,\Gdim\gB\leq \Gdim\gA$.
% 3
\item On a  $\Sdim\gA= \Sdim\gA/\!\Rad\gA$ et  $\Gdim\gA= \Gdim\gA/\!\Rad\gA$.
% 4
\item Si $\gA$ est $n$-stable (section \ref{secSUPPORTS}), alors $\Gdim\gA<n$ (\thref{matrix}).\\
 En abrégé, $\Gdim\gA\leq \Cdim\gA$. 
% 5
\item Si $\Hdim\gA<n$, alors $\Gdim\gA<n$ (\thref{MAINCOR}).\\
 En abrégé, $\Gdim\gA\leq \Hdim\gA$. 
\end{enumerate}
 
\end{fact}
%--------- fin fact ---------------------------------------------- 
%
\begin{proof}
Il faut seulement démontrer les points \emph{3} et \emph{4.}
Le point \emph{4} est clair parce qu'un \elt de $\gA$ est \iv dans $\gA$
\ssi il est \iv dans $\gA/(\Rad\gA)$.

\emph{3 pour $\Sdim$.} Soit $F\in\Ae{m\times r}$ avec $\cD_n(F)=1$ modulo $\fa$. Si $n>\inf(m,r)$ on obtient $1\in\fa$ et tout va bien.
Sinon, soit $a\in\fa$ tel que $1-a\in \cD_n(F)$. \\
On considère la matrice $H\in\Ae{(m+n)\times r}$ obtenue en superposant $F$ et la matrice
$a\In$ suivie de $r-n$ colonnes nulles. \\
On a $1-a^{n}\in \cD_n(F)$, donc
$1\in\cD_n(H)$. Une \coli des colonnes de $H$ est \umdz. La même \coli des colonnes de $F$ est \umd modulo $\fa$. 

\emph{3 pour $\Gdim$.} La même technique fonctionne, mais il suffit
ici de considérer la matrice $H\in\Ae{(m+1)\times r}$ obtenue
en insérant la ligne $[\,a\;0\;\cdots\;0\,]$ en dessous de $F$.
\end{proof}

La \dem du fait suivant aide à justifier la \dfn un peu étonnante 
choisie pour $\Gdim\gA<n$.
 
%:     Fact{factGdimBdim}
\begin{fact}\label{factGdimBdim}
Pour tout $n\geq 0$, on a $\Gdim\gA<n\Rightarrow\Bdim\gA
< n$. \\
En abrégé, $\Bdim\gA\leq \Gdim\gA$.
\end{fact}
%--------- fin fact ----------------------------------------------
%   
\begin{proof}
Par exemple avec $n=3$.
On considère $(a,b_1,b_2,b_3)$ avec $1=\gen{a,b_1,b_2,b_3}$.
On veut des $x_i$ tels que $1=\gen{b_1 + ax_1 ,  b_2 + ax_2,  b_3 + ax_3}$.
On considère la matrice
$F=\cmatrix{b_1 & a & 0 & 0\cr
b_2 & 0 & a & 0\cr
b_3&  0&  0 & a
}= [\,C_0\mid G\,]$ avec $G=a\I_3$.
On a  

\snic{1=\cD_1(C_0)+\cD_3(G)$,\quad  i.e. $1=\gen{b_1,b_2,b_3}+\gen{a^{3}},}

%\sni
car $1=\gen{b_1,b_2,b_3}+\gen{a}$.
En appliquant la \dfn de $\Gdim\gA<3$ à $F$, on obtient
un \vmd $\tra[\,b_1 + ax_1\; b_2 + ax_2\; b_3 + ax_3\,]$.
\end{proof}
% 

%:--- SUBsection{Serre}------------------
\subsec{Le \tho \SSOz}

La version suivante du \tho de Serre est relativement facile,
la partie délicate étant d'établir que $\Sdim  \gA<k$
pour un anneau $\gA$. Modulo les \thrfs{matrix}{MAINCOR} 
on obtient les vraies versions fortes du \thoz.
%:--- Theorem{thSerre}----------------
\begin{theorem}
\label{thSerre} \emph{(\Tho \SSOz, pour la $\Sdim$)}\\
Soit $k\geq 1$ et soit $M$  un \Amo  \pro de rang $\geq k$, ou plus \gnlt
 isomorphe à l'image d'une matrice de rang $\geq k$.\\
Supposons que $\Sdim  \gA< k$.
Alors $M\simeq N\oplus \gA$ pour un certain module~$N$ isomorphe à l'image d'une matrice de rang $\geq k-1$. 
\end{theorem}
%--- end-theorem-----------------------------------------

%-----------------begin proof------------------
\begin{proof}
Soit $F\in\Ae{n\times m}$ une matrice avec $\cD_{k}(F)=1$.
Par \dfnz, on a un vecteur
$u=\tra{[\,u_1\,\cdots \,u_n\,]}\in\Im F$ qui est \umd dans~$\Ae n$. Donc $\gA u$ est un sous-module libre
de rang $1$ en facteur direct dans~$\Ae n$, et a fortiori dans~$M$.
\Prmtz, si $P\in\GAn(\gA)$ est un \prr d'image~$\gA u$, on obtient $M=\gA u\oplus N$ avec 

\snic{N=\Ker(P)\cap M=(\In-P)(M)=\Im\big((\In-P)\,F\big).}

%\sni
Il nous reste à voir que  $(\In-P)\,F$ est de rang $\geq k-1$.
Quitte à localiser et à faire un changement de base, 
on peut supposer que $P$ est la projection standard $\I_{1,n}$.
Alors la matrice $G=(\In-P)\,F$ est la matrice $F$ dans laquelle on a remplacé sa première ligne par $0$, et il est clair que $\cD_k(F)\subseteq \cD_{k-1}(G)$.

\end{proof}
%-----------------end proof------------------

Ainsi, \prmtz, si $M$ est l'image de $F\in\Ae {n\times m}$ de rang $\geq k$,
on obtient une \dcn $M= N\oplus L$ où $L$ est libre de rang $1$ en facteur
direct dans~$\Ae n$ et $N$ isomorphe à l'image d'une matrice de rang $\geq
k-1$. Si maintenant~$F$ est de rang plus grand, on peut itérer le processus
et on a le corolaire suivant (avec la correspondance $h \leftrightarrow k-1$).
	
%:     Corollary{corthSerre}
\begin{corollary}\label{corthSerre}
Soit un anneau $\gA$ tel que $\Sdim \gA\le h$  et soit $M$ un \Amo 
isomorphe à l'image d'une matrice de rang $\geq h+s$.  Alors~$M$
contient en facteur direct un sous-module libre de rang~$s$.  \Prmtz, si~$M$
est l'image d'une matrice $F\in\Ae {n\times m}$ de rang $\geq h+s$, on a $M= N\oplus L$
où $L$ est libre de rang $s$ en facteur direct dans $\Ae n$, et $N$ est
l'image d'une matrice de \hbox{rang $\geq h$}.
\end{corollary}

%:HHH subsec Le \tho de Forster-Swan
\subsec{Le \tho de Forster-Swan} \label{subsecFoSw}

Rappelons qu'un \mtf $M$ est dit \lot engendré par $r$ \elts 
si $\cF_r(M)=1$.
Voir à ce sujet le lemme du nombre de \gtrs local \ref{lemnbgtrlo}.

Le \tho de Forster-Swan ci-dessous a d'abord été établi pour la \ddk ($\Kdim$
à la place de $\Gdim$). La version présentée ici est relativement facile, et la partie délicate est d'établir que $\Gdim  \gA\leq \Kdim  \gA$
pour tout anneau $\gA$. Modulo les \thrfs{matrix}{MAINCOR}
 on obtient les meilleures versions connues du \thoz, sous forme entièrement \covz. 

%:   --- Theorem{thSwan}---------------
\begin{theorem}
\label{thSwan} \emph{(\Tho de Forster-Swan pour la $\Gdim$)} Soit $k\geq 0$ \hbox{et $r\geq 1$}.
Si $\Gdim  \gA\leq   k$, ou même %seulement 
si $\Sdim  \gA\leq   k$
et $\Bdim  \gA\leq   k+r$, et
si un \Amo \tfz~$M$ est   \lot
engendré par $r$ \eltsz, alors il est engendré par $k+r$
\eltsz. 
Dans le premier cas, plus \prmtz, si $M=\gen{y_1 , \dots,  y_{k+r+s}}$,  on peut calculer 

\snic{z_1,\, \ldots, \,z_{k+r}  
\in \gen{y_{k+r+1},\ldots,y_{k+r+s}}}

 tels que
$M$  soit engendré par $(y_1+z_1, \ldots, y_{k+r}+z_{k+r})$.
%si $M$ est engendré   
%par $y_1$, \dots, $y_{k+r+s}$,  on peut calculer $z_1$, \ldots, $z_{k+r}$ 
%dans $\gen{y_{k+r+1},\ldots,y_{k+r+s}}$ tels que
%$M$  soit engendré par $y_1+z_1$, $\ldots$, $y_{k+r}+z_{k+r}$.
\end{theorem}
%--- end-theorem-------------------------------
%-----------------begin proof------------------
\begin{proof}
Puisque $M$ est \tf et $\cF_r(M)=1$, $M$ est le quotient d'un \mpf $M'$
vérifiant $\cF_r(M')=1$. On peut donc supposer que~$M$ est \pfz.
\\
On part d'un \sgr à plus que $k+r$ \elts et on va le remplacer par un \sgr de la forme annoncée avec un \elt de moins.
Soit donc $(y_0,y_1,\dots,y_p)$ un  \sgr de
$M$ \hbox{avec $p\geq k+r$}, et $F$ une \mpn de $M$ pour ce \sysz.
Alors par hypothèse $1 = \cF_{r}(M)= \cD_{p+1-r}(F)$, et puisque $p+1-r\geq k+1$
on a $1 = \cD_{k+1}(F)$.

\emph{Premier cas.}  
Notons $L_0$, \ldots, $L_p$ les lignes
de $F$. Nous appliquons la \dfn de $\Gdim\gA<k+1$ avec  
la matrice transposée de $F$ (qui est de \hbox{rang $\geq k+1$}).
Nous obtenons des  $t_i$ tels que la  
ligne $L_0+t_1L_1+\cdots +t_pL_p$ soit
\umdz.
Remplacer la ligne~$L_0$ par la ligne $L_0+t_1L_1+\cdots +t_pL_p$ revient à
remplacer le \sgr $(y_0,y_1,\dots,y_p)$  
par  

\snic{(y_0,y_1-t_1y_0,\dots,y_p-
t_py_0)=(y_0,y'_1,\dots,y'_p).}

%\sni
Puisque la nouvelle ligne $L_0$ est \umdz, une \coli
convenable des colonnes est de la forme $\tra{[\,1\;y_1\,\cdots \,y_p\,]}$. Cela signifie que l'on a~\hbox{$y_0+y_1y'_1+\cdots +y_py'_p=0$} dans $M$, 
et donc que $(y'_1,\dots,y'_p)$ engendre~$M$.

\emph{Deuxième  cas.}  Nous appliquons la \dfn de $\Sdim\gA<k+1$ avec  
la matrice  $F$. Nous obtenons une \coli de colonnes qui est \umdz, et nous rajoutons cette colonne en première position devant~$F$. 
Puis, en appliquant le fait \ref{corBass} avec $\Bdim\gA<k+r+1\leq p+1$, par \mlr de lignes, nous obtenons une nouvelle \mpn de~$M$ (pour un autre \sgrz) avec la  première colonne égale à $\tra[\,1\,0\,\cdots\,0\,]$. Cela signifie
que le premier \elt du nouveau \sgr est nul.
\end{proof}
%-----------------end proof------------------

Le \thrf{thSwan} est évidemment valable en remplaçant 
l'anneau $\gA$ par l'anneau $\gA/\Ann(M)$ ou  $\gA/\cF_0(M)$.
Nous proposons en \ref{thSwan2} un raffinement un peu plus subtil.

%:--- proposition{basic2}-------------
\begin{proposition}
\label{basic2} Notons $F=[\,C_0\,|\,C_1\,|\,\dots\,|\,C_p\,]\in\Ae{n\times (p+1)}$ (les~$C_i$ sont les colonnes) et
$G=[\,C_1\,|\,\dots\,|\,C_p\,]$, 
de sorte que~$F=[\,C_0\,\vert\,G\,]$.
\\
Si $1\in\cD_1(F)$ et si l'on a 
$\Gdim(\gA/\cD_{k+1}(F))<k$  pour $k\in\lrbq$,
alors il existe $t_1$, \dots, $t_p$ tels que le vecteur 
$C_0+t_1C_1+\cdots+t_pC_p$ est \umd modulo $\cD_{q+1}(F)$.
\end{proposition}
%--- end-theorem-----------------------------------------
%
\begin{proof}
On considère d'abord l'anneau
$\gA_2=\gA/\cD_2(F)$. Puisque $1=\cD_1(F)$ \hbox{et $\Gdim(\gA_2)<1$}, on obtient
\hbox{des $t_{1,i}$} et $C_{1,0}=C_0+t_{1,1}C_1+\cdots+t_{1,p}C_p$ 
tels que $1=\cD_1(C_{1,0})$ modulo~$\cD_2(F)$, \cad $1=\cD_1(C_{1,0}) +\cD_2(G)$. On change $F$ en $F_1$ en remplaçant $C_0$ par $C_{1,0}$
sans changer $G$. Notons que l'on a  $\cD_i(F_1)=\cD_i(F)$ pour tout $i$.
\\
On considère ensuite l'anneau  
$\gA_3=\gA/\cD_3(F_1)$ avec $\Gdim(\gA_3)<2$. 
\\
Puisque $1=\cD_1(C_{1,0}) +\cD_2(G)$, on obtient  
 $C_{2,0}=C_{1,0}+t_{2,1}C_1+\cdots+t_{2,p}C_p$ 
tel que $1=\cD_1(C_{2,0})$ modulo $\cD_3(F)$, \cad $1=\cD_1(C_{2,0}) +\cD_3(G)$. On change $F_1$ en $F_2$ en remplaçant $C_{1,0}$ par $C_{2,0}$
sans changer $G$. On a de nouveau $\cD_i(F_2)=\cD_i(F)$ pour tout $i$.
\\
On continue de la même manière jusqu'à 
obtenir un vecteur $C_{q,0}$ de la \hbox{forme $C_0+t_1C_1+\cdots+t_pC_p$} \umd modulo  $\cD_{q+1}(F)$.
\end{proof}
%

%:--- Theorem{thSwan2}---------------
\begin{theorem}
\label{thSwan2} \emph{(\Tho de Forster-Swan, plus \gnlz, pour la $\Gdim$)}\\
Soit  $M$  un module \tf sur $\gA$. Notons $\ff_k=\cF_k(M)$ ses \idfsz.
Supposons que $1\in\ff_m$ (i.e., $M$ est \lot engendré par $m$ \eltsz) et que pour  $k\in\lrb{0..m-1}$,  on a $\Gdim(\gA/\ff_k)<  m-k$. 
Alors $M$ est engendré par~$m$ \eltsz. 
%\\
Plus \prmtz, si $M=\gen {y_1,\dots,y_{m+s}}$,  on peut calculer des~$z_i $,  
\hbox{dans $\gen {y_{m+1},\dots,y_{m+s}}$} tels que
$M=\gen {y_1+z_1,\dots,y_{m}+z_{m}}$.
\end{theorem}
%--- end-theorem-----------------------------------------
\begin{proof}
Puisque $\ff_0$ annule $M$, on peut remplacer $\gA$ par $\gA/\ff_0$,
ou, ce qui revient au même, supposer que $\ff_0=\cF_0(M)=0$, ce que nous faisons désormais.\\
On considère un \sgr de $M$ avec plus que $m$ \elts et on va le remplacer par un \sgr de la forme annoncée avec un \elt de moins.
Soit donc $(y_0,y_1,\dots,y_p)$ un  \sgr de~$M$ \hbox{avec $p\geq m$}.

Lorsque le module est \pf on raisonne comme  pour le \thref{thSwan}.
\\
Soit $F$ une \mpn de $M$ pour le \sgr considéré.
On a $\ff_{k+1}=\cD_{p-k}(F)$, et en particulier $1\in\ff_p=\cD_1(F)$.
Les hypothèses de la proposition \ref{basic2} sont alors satisfaites avec $q=p$
pour la matrice transposée de $F$.
Si $L_0$, \ldots, $L_p$ sont les lignes
de $F$, nous obtenons des  $t_i$ avec  $L_0+t_1L_1+\cdots +t_pL_p$ 
\umd modulo $\cD_{p+1}(F)=\ff_0=0$. 
On conclut comme au \thref{thSwan}.

%\sni
Le raisonnement dans le cas où $M$ est seulement supposé \tf 
consiste à montrer que $M$ est le quotient d'un module \pf 
qui possède une \mpn  supportant avec succès la \dem de la proposition \ref{basic2}.
Notons $\uy=[\,y_0\;\cdots \;y_p\,]$.   
Toute \syzy entre les $y_i$ 
s'écrit   $\uy\,C=0$ pour un $C\in\gA^{p+1}$.
\\
L'\idf $\ff_{p+1-i}$ de~$M$ est l'idéal $\Delta_i$
somme des \idds $\cD_i(F)$, pour \hbox{les $F\in \gA^{(p+1)\times n}$} 
qui vérifient $\uy\,F=0$, \cad pour les matrices qui sont des \gui{matrices de relations
pour $(y_0, \ldots, y_p)$}.\\
D'après les hypothèses, on a $\Delta_1=1$ et  $\Gdim(\gA/\Delta_{k+1})<k$ \hbox{pour $k\in\lrbp$}.
\\
Le fait que $\Delta_1=1$ se constate sur une matrice de
\syzys $F_1$.  
\\
On considère   la
matrice $\tra F_1$  et l'anneau~$\gA_2=\gA/\Delta_2$.  
\hbox{Comme $\Gdim(\gA_2)<1$},  on obtient une \coli $C_{1,0}$ des colonnes de $\tra F_1$  \umd modulo~$\Delta_2$, 
\cade telle que  $1= \cD_1(C_{1,0})+ \Delta_2$.
\Prmtz, on obtient $C_{1,0}=\tra F_1 X_1$ avec $X_{1}=\tra[\,1\;x_{1,1}\;\cdots\;x_{1,p}\,]$.
\\
L'\egt $1= \cD_1(C_{1,0})+ \Delta_2$ fournit un \elt $a\in\Delta_2$ 
obtenu comme \coli d'un nombre fini de mineurs d'ordre 2 de matrices de \syzysz, et donc  $a\in\cD_2(F_2)$
pour une matrice de \syzys $F_2$. 
On considère alors la  matrice $F'_2=[\,F_1\,|\,F_2\,]$. 
Pour la transposée de $F'_2$ on obtient d'abord que la colonne 
$C_2=\tra {F'_2} \,X_1$ est \umdz. 
On remplace la première colonne de $\tra {F'_2}$
par $C_2$, ce qui donne une matrice $\tra {F''_2}$ qui convient pour les hypothèses de $\Gdim\gA_3<2$ (où $\gA_3=\gA/\Delta_3$), i.e. 
$1=\cD_1(C_{2})+\cD_2(F''_2)$. On obtient 
en définitive  une \coli $C_{2,0}$ des colonnes de $\tra {F'_2}$  \umd  modulo~$\Delta_3$, \cade telle que  $1= \cD_1(C_{2,0})+ \Delta_3$.
\Prmtz, $C_{2,0}=\tra {F'_2}\, X_2$ avec $X_{2}=\tra[\,1\;x_{2,1}\;\cdots\;x_{2,p}\,]$.
Et ainsi de suite.\\
On obtient au bout du compte une matrice de \syzys pour $\uy$, 
$$\preskip.4em \postskip.4em 
F=[\,F_1\,|\,\cdots\,|\,F_p\,] 
$$
et un vecteur $X_{p}=\tra[\,1\;x_{p,1}\;\cdots\;x_{p,p}\,]$
avec la \coli $\tra F\, X_p$ \umd (car \umd modulo $\Delta_{p+1}=\ff_0=0$).\\
On conclut comme au \thref{thSwan}.
\end{proof}

\comm Le \thrf{thSwan2} avec $\Hdim$ ou $\Kdim$ à la place de $\Gdim$ a pour conséquence facile en \clama des
énoncés beaucoup plus abstraits et qui ont l'air beaucoup
plus savants. Par exemple l'énoncé usuel du \tho de Forster-Swan\footnote{Corolaire 2.14 page 108 dans
\cite{Kun} ou \tho 5.8 page 36 dans \cite{Mat}.
En outre, les auteurs remplacent $\gA$ par $\gA\sur{\Ann(M)}$, ce qui ne mange pas
de pain.}
 (énoncé dans le cas où $\Max\gA$ est \noez)
utilise le maximum\footnote{Rappelons que $\jspec\gA$ désigne le sous-espace
de $\SpecA$ formé par les \ideps qui sont intersections d'\idemasz.}, 
pour $\fp\in \jspec\gA$ de $\mu_\fp(M)+\Kdim(\gA/\fp)$:  ici $\mu_\fp(M)$
est le nombre minimum de \gtrs de $M_\fp$.
Ce genre d'énoncé laisse croire
que les \ideps qui sont intersections
d'\idemas jouent un rôle essentiel dans le
\thoz. En réalité, ce n'est pas nécessaire de faire peur aux enfants
avec $\jspec\gA$.
Car ce \tho abstrait est exactement \eqv (dans le cas envisagé,
et en \clamaz) au \thref{thSwan2} pour la~$\Jdim$, qui dans le cas
envisagé est égale à la~$\Hdim$.
En outre, d'un strict point de vue pratique on voit mal comment avoir accès
au bien mystérieux maximum des  $\mu_\fp(M)+\Kdim(\gA/\fp)$.
Par contraste, les hypothèses du \thref{thSwan2} sont susceptibles
d'être certifiées par une \dem \covz, ce qui conduira dans ce cas à un
\algo permettant d'expliciter la conclusion.\eoe

%:HHH subsec{Le \tho de simplification de Bass
\subsec{Le \tho de simplification de Bass}
\label{secSFS}

%:     Definition{defiSimplifiMod}
\begin{definition}\label{defiSimplifiMod}
\'Etant donnés deux modules $M$ et $L$ on dira que $M$
\emph{est simplifiable pour $L$} si $M\oplus L\simeq N\oplus L$
implique   $M \simeq N$.%
\index{simplifiable!module ---}%
\index{module!simplifiable}
\end{definition}

%:     Lemma{lemSimplifiMod}
\begin{lemma}\label{lemSimplifiMod}
Soient $M$ et $L$ deux \Amosz.
Dans les énoncés suivants on a
  1 $\Leftrightarrow$ 2 et 3 $\Rightarrow$ 2.
\begin{enumerate}
\item $M$ est simplifiable pour $L$.
\item Pour toute décomposition
$M\oplus L=M'\oplus L'$ avec $L'\simeq L$, il existe un \auto $\sigma$
de $M\oplus L$ tel que $\sigma(L')=L$.
\item Pour toute décomposition
$M\oplus L=M'\oplus L'$ avec $L'\simeq L$, il existe un \auto $\theta$
de $M\oplus L$  tel que $\theta(L')\subseteq M$.
\end{enumerate}
\end{lemma}
\begin{proof} L'\eqvc de \emph{1} et \emph{2} est un jeu de photocopies.
\\
 \emph{1 $\Rightarrow$ 2.} On suppose $M\oplus L=M'\oplus L'$.
Puisque $L\simarrow L'$, on obtient 
un \iso $M\oplus L\simarrow M'\oplus L$,
donc $M\simarrow M'$. Et l'on obtient en faisant la somme un \iso 
$M\oplus L\simarrow M'\oplus L'$,
\cad un \auto de $M\oplus L$ qui envoie $L$ sur~$L'$.
\\
 \emph{2 $\Rightarrow$ 1.} On suppose $N\oplus L\simarrow M\oplus L$.
Cet \iso envoie $N$ sur~$M'$ et~$L$ sur~$L'$, de sorte que $M\oplus L=M'\oplus L'$.
Il y a donc un \autoz~$\sigma$ de $M\oplus L$ qui envoie $L$ sur $L'$, et disons
$M$ sur $M_1$. Alors:

\snic{N\simeq M'\simeq (M'\oplus L')/L'= (M\oplus L)/L' = (M_1\oplus L')/L'
\simeq M_1\simeq M.}

%\sni
 \emph{3 $\Rightarrow$ 2.} 
 Puisque $\theta(L')$ est facteur direct dans $M\oplus L$, il
est en facteur direct dans $M$, que l'on  écrit $M_1\oplus \theta(L')$.
Soit $\lambda$ l'\auto de $M\oplus L$ qui échange $L$ et $\theta(L')$
en fixant $M_1$.
Alors $\sigma=\lambda\circ \theta$ envoie $L'$ sur $L$.
\end{proof}

Rappelons qu'un \elt $x$ d'un module \emph{arbitraire} $M$ est dit \umd
lorsqu'il existe une forme linéaire $\lambda\in M\sta$ telle que
$\lambda(x)=1$. Il revient au même de dire que $\gA x$ est libre
(de base $x$) et en facteur direct dans $M$
(proposition~\ref{propSplittingOffAlgExt}).

%:--- Theorem{thBassCancel2}----------
\begin{theorem}
\label{thBassCancel2} \emph{(\Tho de simplification de Bass, pour la $\Gdim$)}\\
Soit $M$ un \Amo
\ptf de rang $\geq k$.
Si $\Gdim \gA< k$,
alors $M$ est simplifiable pour tout  \Amo \ptfz:
 si $Q$ est \ptf et 
$M\oplus Q\simeq N\oplus Q$, alors~$M\simeq N$.
\perso{on ne pourra sans doute pas affaiblir l'hypothèse
au cas où $M$ est isomorphe à l'image d'une matrice $F$ de rang $\geq k$}
\end{theorem}
%--- end-theorem-----------------------------------------

\begin{proof}
Supposons avoir montré que $M$ est simplifiable pour $\gA$. \\
Alors,
puisque $M\oplus\Ae\ell$ vérifie aussi l'hypothèse,
on montre par \recu sur $\ell$ que $M$ est simplifiable pour $\Ae{\ell+1}$.
Par suite  $M$ est simplifiable   pour tout facteur direct
dans $\Ae{\ell+1}$.\\
Enfin,  $M$ est simplifiable pour $\gA$ parce qu'il vérifie
le point \emph{3} du lemme~\ref{lemSimplifiMod} pour~$L=\gA$.
En effet, supposons que  $M=\Im F\subseteq\Ae n$,
où $F$ est une \mprn (de rang $\geq k$),
et soit
$L'$ facteur direct \hbox{dans $M\oplus \gA$}, isomorphe à $\gA$:
$L'=\gA(x,a)$ avec $(x,a)$ \umd dans  $M\oplus \gA$.
Puisque toute forme \lin sur $M$ se prolonge à $\gA^{n}$,
il  existe  une forme  $\nu\in (\Ae n)\sta$ telle que $1\in\gen{\nu(x),a}$.
D'après le lemme~\ref{forbass} ci-après, \hbox{avec $x=C_0$},
il existe \hbox{un $y\in M$} tel que $x'=x+a y$
est \umd dans $M$. Considérons une \hbox{forme $\mu\in M\sta$} telle que
$\mu(x')=1$. On définit alors un \auto $\theta$ \hbox{de $M\oplus \gA$} comme suit:
$$\theta=\cmatrix{1&0\cr -a\mu&1}\cmatrix{1&y\cr0&1}\;
i.e.,\; \cmatrix{m\cr b}\mapsto \cmatrix{m+by\cr \mu(x)b-a\mu(m)}.$$
Alors
$\theta(x,a)=(x',0)$, donc $\theta(L')\subseteq M$.
On conclut avec le lemme~\ref{lemSimplifiMod}.
\end{proof}

Dans le lemme suivant, qui termine la \dem
du \thrf{thBassCancel2},  nous reprenons les notations
de la proposition \ref{basic2}, la matrice $F=[\,C_0\,|\,C_1\,|\,\dots\,|\,C_p\,]$ étant celle du \tho précédent.
%:--- lemma {forbass}  %%%%%%%%%%%%%%%%%%%%%%%%%%%
\begin{lemma}\label{forbass}
Si  $\Gdim \gA <k$
et $\cD_k(F) = 1 =\DA(C_0 )\vu \DA(a)$, alors il
\hbox{existe
$t_1$, \dots, $t_p$} tels que 

\snic{1 = \DA(C_0 +at_1C_1+\cdots+at_pC_p).}
\end{lemma}

\begin{proof}
On considère la matrice
$[\,C_0\,|\,aC_1\,|\,\dots\,|\,aC_p\,] $,
obtenue en remplaçant $G$ par $aG$ dans $F$.
Comme $\DA(C_0)\vu \cD_k(G)=1=\DA(C_0)\vu\DA(a)$,
on a bien par \dit $\DA(C_0)\vu \cD_k(aG)=1$. 
On conclut avec  $\Gdim\gA<k$.
\end{proof}

%:
\subsec{Une \prt \cara simple pour $\Gdim\gA<n$}

Pour démontrer $\Gdim\gA<n$ pour un anneau $\gA$, il suffit de vérifier la conclusion (dans la \dfn de $\Gdim\gA<n$)
pour des matrices particulièrement simples.
Cela fait l'objet de la proposition qui suit.

%:     Proposition{propGdimGdim}
\begin{proposition}\label{propGdimGdim}
Pour un anneau $\gA$ on a $\,\Gdim\gA<n$ \ssi  
pour toute matrice  $V\in\MM_{n+1}(\gA)$ de la forme
$$\preskip.4em \postskip.4em
V=\cmatrix{
b&c_1&\cdots&\cdots&c_n\cr
b_1&a&0&\cdots&0\cr
\vdots&0&\ddots&\ddots&\vdots\cr
\vdots&\vdots&\ddots& \ddots&0\cr
b_n&0&\cdots&0&a\cr
}=[\,V_0\mid V_1 \mid \dots \mid V_n\,]
,$$
et pour tout $d \in \gA$ tel que $1 = \gen{b,a,d}$,
il existe  \hbox{des  $x_i\in \gA$} tels que
$$\preskip.3em \postskip.8em
1 = \cD_1(V_0+x_1V_1+\cdots+x_nV_n)+\gen{d}.
$$
\end{proposition}
%--------- fin proposition ---------------------------------------------- 
\rem Au lieu d'utiliser un \elt $d$ soumis à la contrainte $1=\gen{a,b,d}$,
on aurait pu utiliser un couple $(u,v)$ soumis à aucune contrainte et remplacer $d$ \hbox{par $1+au+bv$} dans la conclusion. Sous cette forme, il est particulièrement évident que si la condition ci-dessus est satisfaite pour l'anneau $\gA$, elle est satisfaite pour tout quotient de $\gA$. \eoe
\begin{proof}
Pour montrer que la condition est \ncrz, nous raisonnons avec l'anneau quotient 
$\gB=\aqo\gA d$ et nous considérons la matrice 

\snic{F=V=[\,V_0\mid V_1 \mid \dots \mid V_n\,].}

%\sni
Avec  les notations de la \dfn \ref{defiSdimGdim} on a $p=n$, $F=[\,C_0\mid G\,]$,  \hbox{et $C_i=V_i$} pour $i\in\lrbzn$. 
Puisque $1 = \gen{b,a,d}$ dans $\gA$,
on a 
$$\preskip.4em \postskip.4em 
1 = \gen{b,a^{n}}\subseteq \cD_1(C_0)+\cD_n(G)\hbox{  dans }\gB, 
$$ 
et l'hypothèse de la \dfn est satisfaite. Puisque $\Gdim\gB<n$, on obtient
des~$x_i$ dans $\gA$ tels que 
$$\preskip.4em \postskip.4em 
1=\cD_1(C_0+x_1C_1+\cdots+x_nC_n)\hbox { dans } \gB . 
$$
D'où la conclusion souhaitée dans $\gA$.

Pour démontrer la réciproque nous   procédons en deux étapes.
Rappelons tout d'abord que si la condition est vérifiée pour l'anneau $\gA$, elle est vérifiée pour tout quotient de $\gA$.
Nous allons en fait utiliser cette condition \hbox{avec $d=0$} (l'hypothèse 
sur $V$ devient alors du même type que celle qui sert à définir~\hbox{$\Gdim<n$}),
avec l'anneau $\gA$ et certains de ses quotients.

\emph{Première étape: le cas où la matrice $F$ possède $n+1$ colonnes, i.e. $p=n$.} 
%On note $m$ le nombre de lignes de $F$.
Avec $F\in\Ae{m\times (n+1)}$, on a par hypothèse une forme \lin
$\varphi_0:\Ae m \to\gA$ et une forme $n$-\lin alternée $\psi : (\Ae m)^n
\to\gA$ telles que 

\snic{1=\varphi_0(C_0)+\psi(C_1, \dots, C_n).}

%\sni
Pour $j\in\lrbn
$ notons $\varphi_j:\Ae m \to\gA$ la forme \lin 

\snic{X\mapsto \psi(C_1, \dots,
C_{j-1},X,C_j, \dots, C_n).}

%\sni
En notant $a = \psi(C_1, \ldots, C_n)$, on a alors:
\begin{itemize}
\item $\varphi_1(C_1)=\cdots=\varphi_n (C_n) = a$,
\item  $\varphi_i(C_j)=0$ si $1\leq i\neq j\leq n$
\end{itemize}
Nous considérons la matrice des $\varphi_i(C_j)$, nous obtenons:
$$
V = [\,V_0\mid \dots \mid V_n\,] :=
\Cmatrix{2pt}{
\varphi_0(C_0)&\varphi_0(C_1)&\cdots&\cdots&\varphi_0(C_n)\\[.3em]
\varphi_1(C_0)&a&\phantom{m}0\phantom{m}&\cdots&0\\[.3em]
\vdots&0&\ddots& \ddots&\vdots\\[.1em]
\vdots&\vdots&\ddots& \ddots&0\\[.3em]
\varphi_n(C_0)&0&\cdots&\phantom{m}0\phantom{m}&a
}
,$$
\cad $V = [\,\varphi(C_0)\mid \dots \mid \varphi(C_n)\,]$
en notant
$\varphi(Z)=\Cmatrix{2pt}{\varphi_0(Z)\cr\vdots\cr\varphi_n(Z)}$.\\
On peut appliquer l'hypothèse avec $d=0$. On trouve $x_1$, \dots, $x_n\in\gA$
tels  que le vecteur $V_0+x_1V_1+\cdots+x_nV_n$ est \umdz.
Ce vecteur est égal \hbox{à
$\varphi(C_0+x_1C_1+\cdots+x_nC_n)=\varphi(C)$}.
Puisque ce vecteur est \umd et que  $\varphi$ est
\linz, le vecteur $C$ est lui-même  \umdz.
%Notez que la conclusion ici est \gui{un peu meilleure} que celle qui sert 
%à définir $\Gdim\gA<n$ puisqu'on a obtenu des $x_i\in a \gA$. En fait,
% ceci n'est qu'un détail technique qui permet de traiter la deuxième 
% étape plus confortablement.

\emph{Deuxième étape: le cas général.}\\
Comme $1=\cD_1(C_0 )+\cD_n(G)$, on a une famille
$(\alpha_i)_{i\in \lrbq}$ de parties à~$n$
\elts de $\lrbp$ telle que $1 = \cD_1(C_0 )+ \som_i \cD_n(G_{\alpha_i})$,
où $G_{\alpha_i}$ est la matrice extraite de $G$  en considérant
uniquement les colonnes dont l'indice est dans~$\alpha_i$.
On note $C_{0,0}=C_0$ et $J_\ell = \som _{i>\ell} \cD_n(G_{\alpha_i})$. On applique alors le cas de la première étape successivement avec $\ell=1$, \ldots, $q$
pour obtenir 

\snic{1= \cD_1(C_{0,\ell})=\cD_1(C_{0,\ell-1}) + \cD_n(G_{\alpha_\ell})
\hbox{ dans }\gA/J_\ell}

%\sni
et donc $\cD_1(C_{0,q}) = 1$ dans~$\gA$.
\\
Notez que dans cette deuxième étape, nous utilisons le résultat de la première étape avec des anneaux quotients de $\gA$.
\end{proof}
%

%--- Section{Supports}--------------
\section{Supports et \texorpdfstring{$n$}{n}-stabilité}
%-----------------------------------------
\label{secSUPPORTS}

Dans la section \ref{secManipElemCol} nous établirons des \thos concernant
les \mlrs de matrices. Ils auront comme corolaires
de grands \thos dus à Serre, Forster, Bass et Swan.

Nous les donnerons en deux versions similaires mais néanmoins différentes.
Nous ne pensons pas qu'elles puissent être ramenées à une forme unique.

La première version est basée
sur la notion de $n$-stabilité%
%, notion dérivée de la \ddk
%que nous introduisons dans la section \ref{secSUPPORTS}
.
Cette version aboutit notamment à un résultat sophistiqué d\^u
à Bass dans lequel intervient une partition du spectre de Zariski en
un nombre fini de parties qui sont toutes de dimension petite
(plus petite que la \ddk de l'anneau). Ce résultat
sera utilisé dans le chapitre \ref{ChapMPEtendus}
pour démontrer le \thrf{thBassValu} de Bass
concernant les modules étendus.

La deuxième version utilise la dimension de Heitmann,
introduite dans la section \ref{subsecDimHeit},
inférieure ou égale à la \ddkz,
mais pour laquelle on ne connaît pas
d'analogue de la version sophistiquée de Bass.

La section \ref{secSUPPORTS} donne quelques préliminaires \ncrs pour
la première version basée sur la $n$-stabilité. 

%:  subsec  supports
\subsec{Supports, dimension, stabilité}

%:     Definition{defiSupportT}
\begin{definition}\label{defiSupportT}
Un  \emph{\sutz} sur un anneau $\gA$ dans un \trdiz~$\gT$ est une application
$D:\gA\to\gT$, qui vérifie les axiomes suivants:%
\index{support!sur un anneau commutatif}
\vspace{-2pt}
\begin{itemize}
\item [$\qquad\bullet$] $D(0_\gA)=0_\gT,\;\;\; D(1_\gA)=1_\gT$,
\item [$\qquad\bullet$]  $D(ab)=D(a)\vi D(b),$
\item [$\qquad\bullet$]  $D(a+b)\leq D(a)\vu D(b)$.
\end{itemize}
On notera $D(\xn)=D(x_1)\vu\cdots\vu D(x_n)$.
\end{definition}

Il est clair que $\DA:\gA\to\ZarA$ est un \sutz, appelé \emph{\sut de Zariski}.
Le lemme suivant montre que le \sut de Zariski est le \sut \gui{libre}.%
\index{support!de Zariski}

%:     Lemma{lemAnodin}
\begin{lemma}\label{lemAnodin}
Pour tout \sut $D$ on a:
\begin{enumerate}
\item \label{i1lemAnodin} $D(a^m)=D(a)$ pour $m\geq1$, $D(ax)\leq D(x),$ $D(a,b)=D(a+b,ab).$
\item \label{i2lemAnodin} $\gen{\xn}= \gen{\yr} $ implique $D(\xn)=D(\yr).$
\item \label{i3lemAnodin} $\DA(y)\leq\DA(\xn)$ implique $D(y)\leq D(\xn).$
\item \label{i4lemAnodin} Il existe un unique
\homo $\theta$ de \trdis qui fait commuter le diagramme suivant:
\Pnv{\gA}{\DA}{ D }{\ZarA}{\theta}{\gT}{ }{\suts}{\homos de \trdisz}
\end{enumerate}
\end{lemma}
\facile

Ainsi tout \sut $D:\gA\to\gT$  tel que $D(\gA)$ engendre $\gT$ en tant que
\trdi est obtenu en composant le \sut de Zariski avec un
passage au quotient $\ZarA\to\ZarA\sur{\sim}$ par une relation d'\eqvc
compatible avec la structure de treillis.

On notera $D(\fa) $  pour $D(\xn)$ si $\fa=\gen{\xn}$.
On dira qu'un vecteur $X\in\Ae n$ est \emph{$D$-\umdz} si $D(X)=1$.
\index{Dunimo@$D$-unimodulaire!vecteur}

%   subsec{Dimension d'un support
\subsubsection*{Dimension d'un \sutz, \tho de \KRNz}

%:     Definition{defi2SuitComp}
\begin{definition}\label{defi2SuitComp} 
\'Etant données deux suites $(\xzn)$ et~$(\bzn)$ dans~$\gA$ et un \sutz~$D$ sur  $\gA$, on dit que les deux suites
 sont  \emph{$D$-\copsz} si l'on a les in\egts suivantes:
%---  equation eqC1G --------
\begin{equation}\label{eqC1G}
\left.\arraycolsep2pt
\begin{array}{rcl}
 D(b_0x_0)& =  & D(0)    \\
 D(b_1x_1)& \leq  &  D(b_0,x_0)  \\
\vdots ~~~~& \vdots  &~~~~  \vdots \\
 D(b_n x_n )& \leq  &  D(b_{n -1},x_{n -1})  \\
 D(1)& =  &   D(b_n,x_n )
\end{array}
\right\}
\end{equation}
%---------------------end equation--------------
Le \sut $D$ est dit \emph{de \ddk $\leq n$} si toute suite $(\xzn)$ dans $\gA$
admet une suite $D$-\copz. On notera $\Kdim(D)\leq n$.%
\index{Dcomple@$D$-complémentaires!suites ---}%
\index{dimension de Krull!d'un support}%
\index{complémentaires!suites ---}%
\index{suites complémentaires!pour un support}
 \end{definition}

Par exemple pour $n=2$ les suites \cops correspondent au dessin suivant dans $\gT$.
$$
\SCO{D(x_0)}{D(x_1)}{D(x_2)}{D(b_0)}{D(b_1)}{D(b_2)}
$$

%\smallskip 
\rem Notons que $\Kdim\gA=\Kdim(\DA)$.
\eoe

\medskip
La preuve du lemme suivant peut être recopiée sur celle du lemme
\ref{lemKroH}
en remplaçant $\DA$ par $D$. Le \tho de Kronecker est ensuite une
conséquence directe.

%--- Lemma{lem2KroH}----------
\begin{lemma}
\label{lem2KroH}
Soit $\ell\geq 1$.
Si $(b_1,\ldots ,b_\ell)$ et  $(\xl)$ sont deux suites
$D$-\cops dans $\gA$,
alors pour tout $a\in\gA$ on~a:
$$ D(a,b_1,\dots,b_\ell) =  D(b_1+ax_1,\dots,b_\ell+ax_\ell),$$
\cad encore: $D(a)\leq  D(b_1+ax_1,\dots,b_\ell+ax_\ell)$.
\end{lemma}
%--- end-lemma-----------------------------------------

%: --- Theorem{th2KroH}----------
\begin{theorem}
\label{th2KroH} \emph{(\Tho de \KRNz, pour les \sutsz)}\\
Si $D$  est un \sut de \ddk $\leq n$, pour tout
\itfz~$\fa$
il existe un \id $\fb$ engendré par $n+1$ \elts tels que $D(\fa)=D(\fb)$.
\linebreak 
En fait, pour tous $b_1$, \ldots, $b_{n+r}$ ($r\geq 2$), il existe des  $c_j\in\gen{b_{n+2},\ldots ,b_{n+r} }$
tels que $D(b_1+c_1,\ldots,b_{n+1}+c_{n+1})=D(b_1,\ldots,b_{n+r})$.
\end{theorem}
%--- end-theorem-----------------------------------------

%   subsec{Supports fidèles
\subsubsection*{Supports fidèles}

Dans ce paragraphe on démontre en particulier que la \ddk d'un
anneau (que l'on sait déjà égale à la dimension de son support de Zariski) 
est égale à celle de son treillis de Zariski: on tient ici la promesse faite en \ref{corfactDDKTRDI}.

%:     Definition{defisutfidele}
\begin{definition}\label{defisutfidele}
Un \sut $D:\gA\to\gT$ est dit \emph{fidèle} si
$\gT$ est engendré par l'image de $D$ et si, pour tout $a\in\gA$
 et $L\in \Ae m$, l'in\egt $D(a)\leq D(L)$
implique l'existence d'un $b\in\gen{L}$ tel que
 $D(a)\leq  D(b)$.
\index{support!fidèle}
\index{fidèle!support ---}
\end{definition}

Par exemple le \sut de Zariski $\DA$
est toujours fidèle.

\medskip Soit $D:\gA\to\gT$ un support. Si l'image de $\gA$ engendre $\gT$,
puisqu'on a l'\egt $D(a_1)\vi\cdots \vi D(a_n)=D(a_1\cdots a_n)$, tout \elt
de $\gT$ peut s'écrire sous forme $D(L)$ pour une liste $L$ d'\elts de $\gA$.

%:     Lemma{lemsutfidele}
\begin{lemma}\label{lemsutfidele}
Si $D$ est fidèle et $\Kdim\gT<k$ alors $\Kdim(D)<k$.
En particulier la \ddk d'un anneau est égale à celle 
de son treillis de Zariski.
\end{lemma}
\begin{proof}
Soit la suite $(a_1,\ldots,a_k)$ dans $\gA$. Nous devons démontrer qu'elle
admet une suite $D$-\copz. \\
Puisque $\Kdim\gT<k$, la suite
$\big(D(a_1),\ldots,D(a_k)\big)$ possède une suite \cop $\big(D(L_1),\ldots,D(L_k)\big)$
dans $\gT$, avec pour $L_i$ des listes dans
$\gA$:
$$\arraycolsep2pt
\begin{array}{rcl}
 D(a_1)\vi D(L_1)& =  & D(0)    \\
 D(a_2)\vi D(L_2)& \leq  &  D(a_1,L_1)  \\
\vdots~~~~~~~~& \vdots  &~~~~  \vdots \\
 D(a_k)\vi D(L_k )& \leq  &  D(a_{k -1},L_{k -1})  \\
 D(1)& =  &   D(a_k,L_k).
\end{array}
$$
Puisque $D$ est fidèle, il existe $c_k$ dans $\gen{a_k,L_k}$
tel que $D(1)\leq D(c_k)$, ce qui donne $b_k\in\gen{L_k}$
tel que $D(1)\leq D(a_k,b_k)$.\\
Notons que l'on  a
$$
D(a_kb_k)=D(a_k)\vi D(b_k )\leq D(a_k)\vi D(L_k )\leq D(a_{k -1},L_{k -1}).
$$
Puisque $D$ est fidèle, on a $c_{k-1}\in\gen{a_{k-1},L_{k-1}}$
avec $D(a_kb_k)\leq D(c_{k-1})$, ce qui donne $b_{k-1}\in\gen{L_{k-1}}$
tel que $D(a_kb_k)\leq D(a_{k-1},b_{k-1})$.\\
Et ainsi de suite. Au bout du compte on a construit une suite
$(b_1,\ldots,b_k)$ qui est $D$-\cop de $(a_1,\ldots,a_k)$.
\end{proof}
%

%:HHH passe en subsec \Suts $n$-stables
%   subsec{\Suts $n$-stables
\subsubsection*{\Suts $n$-stables}

On abstrait maintenant la \prt décrite au lemme \ref{lem2KroH}
pour les suites \cops sous la forme suivante.

%:HHH deux defi fusionnées
%:     Definition{defisutnstab}
\begin{definition}\label{defisutnstab}~
\begin{enumerate}
\item Soit $n\geq 1$. Un \sut $D:\gA\to\gT$ est dit \emph{$n$-stable} lorsque, pour 
\linebreak 
tout $a\in\gA$
 et $L\in \Ae n$, il existe $X\in \Ae n$
tel que $D(L,a)=  D(L+aX)$, \cad $D(a)\leq  D(L+aX)$.%
\index{support!nstab@$n$-stable}%
\index{nstab@$n$-stable!support ---}
\item \label{defiAnneaunStable}
L'anneau $\gA$ est dit \emph{$n$-stable} si son \sut de Zariski
$\DA$ est $n$-stable.
%:HHH ajout
On notera $\Cdim\gA<n$ pour dire que $\gA$ est $n$-stable.%
\index{anneau!nstab@$n$-stable}%
\index{nstab@$n$-stable!anneau ---}
\item L'anneau $\gA$ est dit $0$-stable s'il est trivial.
\end{enumerate}

\end{definition}

Dans l'acronyme $\Cdim$, $\mathsf{C}$ fait allusion à \gui{Coquand}.

Naturellement si  $\Kdim(D)< n$ alors $D$ est $n$-stable.
En particulier, avec le support libre $\DA$, on obtient $\Cdim\gA\leq \Kdim\gA$.
Par ailleurs,  le \tho de Kronecker s'applique (presque par \dfnz)
à tout \sut $n$-stable.
%:HHH ajout

La notation $\Cdim\gA<n$ est justifiée par le fait que si $D$ est $n$-stable, il \hbox{est $(n+1)$-stable}. % pour $m>n$. 
Enfin, le point \emph{3} dans la \dfn a été donné pour plus de clarté,
mais il n'est pas vraiment \ncrz: en lisant le point \emph{1}
\hbox{pour $n=0$}, on obtient que pour tout $a\in\gA$, $D(a)\leq D(0)$.

%:HHH ajout
\medskip 
\exls ~\\ 
1) Un \advz, ou plus \gnlt un anneau $\gV$ qui vérifie \gui{$a\mid b$ ou $b\mid a$ pour tous $a$, $b$}, est $1$-stable, même en \ddk infinie.
Pour tout $(a, b)$ il suffit de trouver un $x$ tel que $\gen{a, b} = \gen{b + x a}$.
\hbox{Si $a = q b$}, on a $\gen{a, b} = \gen{b}$ et l'on prend $x=0$.
Si $b = q a$, on a $\gen{a, b} = \gen{a}$ et l'on prend $x=1-q$.

%\sni
2) Un  anneau de Bézout intègre est $2$-stable. Plus \gnltz, un anneau de Bézout strict   (cf. la section \ref{secBézout} \paref{secpfval} et
l'exercice~\ref{exoAnneauBézoutStrict})  est $2$-stable. 
%Autrement dit, pour $a$, $b_1$, $b_2 \in \gA$,
%il existe $x_1$, $x_2$ tels \hbox{que $a\in 
%\DA({b_1 + x_1a, b_2 + x_2a})$}.
Plus \prmtz, pour $a$, $b_1$, $b_2 \in \gA$,
il existe $x_1$, $x_2$ tels \hbox{que $a\in 
\gen {b_1 + x_1a, b_2 + x_2a}$}, i.e. 
$\gen {a, b_1, b_2} = \gen {b_1 + x_1a, b_2 + x_2a}$.
\\
En effet,  d'après la question \emph {1c} de
l'exercice, il existe $u_1$ et $u_2$ \com tels que $u_1{b_1} + u_2{b_2} = 0$.
On prend $x_1$, $x_2$ tels que $u_1x_1 + u_2x_2 = 1$ et l'on obtient l'\egt

$\qquad\qquad\quad 
a = u_1b_1 + a + u_2b_2 =  %u_1b_1 + u_1x_1a + u_2x_2a + u_2b_2 =  
u_1(b_1 + x_1a) + u_2(b_2 + x_2a)
$.
\eoe

%:     Fact{factStBdim}
\begin{fact}\label{factStBdim}
 On a toujours  $\,\Bdim\gA\leq \Cdim\gA$.
\end{fact}
%--------- fin fact ---------------------------------------------- 
%
\begin{proof}
Si $\gA$ est $n$-stable, alors $\Bdim \gA 
< n$: en effet, on applique la \dfn avec $(a,\an)$ dans $\gA$ vérifiant $1 \in \gen{a, \an}$.
\end{proof}
%

%:     Fait{lemsutnstab}
\begin{fact}\label{lemsutnstab}
Si $D$ est $n$-stable, pour tout $a\in\gA$
 et $L\in \Ae n$, il existe $X\in \Ae n$
tel que $D(L,a)=  D(L+a^2X)$, \cad $D(a)\leq  D(L+a^2X)$.
\end{fact}

En effet, $D(a)=D(a^2)$ et $D(L,a)=D(L,a^2)$.

%:HHH passe en subsec  Constructions et recollements de \sutsz
%:  subsec{Constructions de \sutsz
\subsec{Constructions et recollements de \sutsz}

%:     definition{lemsutHeit}
\begin{definition}\label{lemsutHeit} ~\\
L'application  $\JA:\gA\to\Heit\gA$ définit le
\emph{\sut de Heitmann}.
\index{support!de Heitmann}
\end{definition}

\rem
A priori $\Kdim\DA=\Kdim\gA\geq \Kdim\JA \geq \Jdim\gA$.
On manque d'exemples  qui montreraient que
les deux in\egts peuvent être strictes.
\eoe

%:     Lemma{lem2GaussJoyal}
\begin{lemma}\label{lem2GaussJoyal} \emph{(Variante du lemme de Gauss-Joyal
\ref{lemGaussJoyal})\iJG}\\
Si $D$ est un \sut sur $\gA$, on obtient un \sut $D[X]$ sur $\gA[X]$
en posant 
$$\preskip.2em \postskip.3em
D[X](f)=D\big(\rc(f)\big).$$
\end{lemma}
\begin{proof}
Le lemme \ref{lemGaussJoyal}  donne
$\DA\big(\rc(fg)\big)=\DA\big(\rc(f)\big)\vi \DA\big(\rc(g)\big)$.
\end{proof}
%
%:     Lemma{lemsutquo}
\begin{lemma}\label{lemsutquo} \emph{(\Sut et quotient)}
Soit $D:\gA\to\gT$  un \sut et $\fa$ un \itf de
$\gA$.
On obtient un \sut
$$\preskip.1em \postskip.4em
D/\fa:\gA\to\gT\sur{\fa}\eqdefi \gT\sur{(D(\fa)=0)},
$$
en composant
$D$ avec la projection $\Pi_{D(\fa)}:\gT\to\gT\sur{(D(\fa)=0)}$. 
\begin{enumerate}
\item $\DA/\fa$ est canoniquement isomorphe à $\rD_{\gA\sur{\fa}}\circ \Zar (\pi_\fa)$,
où $\pi_\fa$ est l'application canonique $\gA\to\gA\sur{\fa}$.
\item Si $D$ est fidèle, alors $D/\fa$ \egmtz.
%:HHH ajout
\item Si $D$ est $n$-stable, alors $D/\fa$ \egmtz.\\
En particulier $\Cdim\gA/\fa\leq \Cdim\gA$.
\end{enumerate}
\end{lemma}
\begin{proof}Rappelons que
$\Pi_{D(\fa)}(x)\leq\Pi_{D(\fa)}(y)\iff x\vu {D(\fa)}\leq y\vu {D(\fa)}$.
\\
\emph{1.} Résulte du fait \ref{fact2Zar}. \\
\emph{2.} Notons
$D'=D\sur{\fa}$. Soient $a\in\gA$ et $L$ un vecteur tels  
que~\hbox{$D'(a)\leq D'(L)$}.
On cherche un $b\in\gen{L}$  tel que~$D'(a)\leq D'(b)$. Par \dfn de~$D'$
on a~\hbox{$D(a)\leq D(L,\fa)$}, et puisque $D$ est fidèle, il existe $c\in\gen{L}+\fa$
tel que~\hbox{$D(a)\leq D(c)$}, ce qui donne un $b\in L$ tel que
$D(a)\leq D(b,\fa)$, autrement dit~$D'(a)\leq D'(b)$.\\
%:HHH ajout
\emph{3.} Soient $a \in \gA$, et $L \in \Ae n$. On cherche $X \in \Ae n$
tel que~\hbox{$D'(a) \leq  D'(L + aX)$}, i.e.
$D(a) \vu D(\fa)  \leq  D(L + aX) \vu D(\fa)$.
Or on a un $X$ qui convient pour $D$, i.e. $D(a) \leq  D(L + aX)$, 
donc il convient pour $D'$.
\end{proof}

De manière duale on a le lemme suivant.
%:     Lemma{lemsutloc}
\begin{lemma}\label{lemsutloc} \emph{(\Sut et \lonz)}
Soit $D:\gA\to\gT$  un \sut et un \elt $u$ de $\gA$.
On obtient un \sut

\snic{D[1/u]:\gA\to\gT[1/u]\eqdefi\gT\sur{(D(u)=1)},}

%\sni
en composant
$D$ avec $j_{D(u)}:\gT\to\gT\sur{(D(u)=1)}$.
\begin{enumerate}
\item $\DA[1/u]$ est canoniquement isomorphe à $\rD_{\gA[1/u]}\circ \Zar (\iota_u)$,
où $\iota_u$ est l'application canonique $\gA\to\gA[1/u]$.
\item Si $D$ est fidèle, alors $D[1/u]$ \egmtz.
%:HHH ajout
\item Si $D$ est $n$-stable, alors $D[1/u]$ \egmtz. \\
En particulier $\Cdim\gA[1/u]\leq \Cdim\gA$.
\end{enumerate}
\end{lemma}
\begin{proof}
Rappelons que $
j_{D(u)}(x)\leq j_{D(u)}(y)\iff x\vi {D(u)}\leq y\vi {D(u)}.$
\\
\emph{1.} Résulte du fait \ref{fact2Zar}. \\
\emph{2.} Notons
$D'=D[1/u]$. Soit $a\in\gA$ et $L$ un vecteur tels que $D'(a)\leq D'(L)$.
 Par \dfn de $D'$
on a $D(au)=D(a)\vi D(u)\leq D(L)$. Puisque $D$ est fidèle, il existe
$b\in\gen{L}$
tel que $D(au)\leq D(b)$, \cadz~$D'(a)\leq D'(b)$.\\
%:HHH ajout
\emph{3.} Comme pour le lemme \ref{lemsutquo} en remplaçant
$D/\fa$ et $\vu$ par $D[1/u]$  et~$\vi$.
\end{proof}
%

%:HHH deplace augmenté
%:     Lemma{lemPartitionSpec}
\begin{lemma}\label{lemPartitionSpec}\label{corlemPartitionSpec}~
\begin{enumerate}
\item Soit un \sut $D:\gA\to\gT$ et $b\in \gA$. 
\begin{enumerate}
\item $D/b$ et $D[1/b]$ sont $n$-stables \ssi $D$ est $n$-stable.
\item Si
$D$ est fidèle et si $\gT/b$ et $\gT[1/b]$ sont de \ddk $< n$,
alors $D$ est $n$-stable.
\end{enumerate}
\item Soit un anneau $\gA$ et $b\in\gA$. Alors $\aqo\gA b$ et $\gA[1/b]$ sont $n$-stables \ssi $\gA$ est $n$-stable. \\
De manière abrégée: $\Cdim\gA=\sup\big(\Cdim\aqo\gA b,\Cdim\gA[1/b]\big)$.
%
%\item 
%
\end{enumerate}

\end{lemma}
\begin{proof} Il suffit de montrer l'implication directe dans le point \emph{1a.}\\
Soient $a\in \gA$ et~$L\in\Ae n$. 
Puisque $D/b$ est $n$-stable, on a un $Y\in\Ae n$ tel que~\hbox{$D(a)\leq  D(L+aY)$}
dans ${\gT/\big(D(b)=0\big)}$, \cad dans~$\gT$:
$$
D(a)\leq D(b)\vu D(L+aY).\eqno(*)
$$
Ensuite on applique la $n$-stabilité de $D[1/b]$ avec $ab$ et $L+aY$
ce qui fournit un $Z\in\Ae n$ tel que $D(ab)\leq  D(L+aY+abZ)$
dans ${\gT/\big(D(b)=1\big)}$. 
\\
Dans~$\gT$, en posant $X = Y+bZ$, cela s'écrit:
$$D
(ab)\vi D(b) \leq D(L+aX), \quad \hbox {i.e.}\quad D(ab)\leq D(L+aX).
\eqno (\#)
$$
Mais on a $\gen{b,L+aX}=\gen{b,L+aY}$, donc $D(b,L+aX)= D(b,L+aY)$.
Les in\egts $(*)$ et $(\#)$ s'écrivent alors

\snic{D(a)\leq D(b)\vu D(L+aX) \et D(a)\vi D(b)\leq D(L+aX).}

%\sni
Ceci implique (par \gui{coupure}, cf. \paref{coupure1})
que $D(a)\leq  D(L+aX)$.
\end{proof}
%

%  subsec{Partitions constructibles du spectre de Zariski
\subsubsection*{Partitions constructibles du spectre de Zariski}

Une partie
\ix{constructible} de  $\SpecA$ est une combinaison
booléenne d'ouverts de base $\fD(a)$. En \clamaz, si l'on munit l'ensemble~$\SpecA$
de la \gui{topologie constructible} ayant pour base d'ouverts les parties constructibles, on obtient un espace spectral, le \ix{spectre constructible} \emph{de l'anneau~$\gA$}, que l'on peut identifier à $\Spec\Abul$. 

D'un point de vue \cofz, on a vu que l'on  peut remplacer $\SpecA$ (un objet un peu trop idéal) par le treillis~$\ZarA$ (un objet concret), isomorphe en \clama au treillis des \oqcs de $\Spec\gA$.
Lorsque l'on  passe de la topologie de Zariski à la topologie constructible en \clamaz, on 
passe de $\ZarA$ à
$\Bo(\ZarA)\simeq\Zar(\Abul)$ en \coma (pour ce dernier \isoz, voir le \thrf{thZedGenEtBoolGen}).

Hyman Bass s'est intéressé aux partitions constructibles du spectre de Zariski.
%:2015 partitions constructibles
Une étape \elr de la construction d'une telle partition consiste en le remplacement
d'un anneau $\gB$ par les deux anneaux $\aqo{\gB}{b}$ et~$\gB[1/b]$,
pour un \elt $b$ de $\gB$.
Une remarque importante qu'a faite Bass est que ces deux anneaux peuvent avoir
chacun une \ddk strictement plus petite que celle de~$\gB$,
alors que certaines \prts de l'anneau n'ont besoin, pour être vérifiées
dans $\gB$, que d'être vérifiées dans chacun de ses deux fils.
C'est le cas pour la $n$-stabilité du \sut libre. Telle est en tout cas l'analyse qu'a faite T.~Coquand de quelques 
pages de~\cite{BASS}.
\perso{citer les pages en question}

En \clamaz, de tout recouvrement de $\Spec\gA$
par des ouverts de la topologie constructible, on peut extraire un recouvrement fini, que l'on peut raffiner en une partition finie par des \oqcs (\cad des combinaisons booléennes finies d'ouverts de base~$\fD(a)$).
C'est beaucoup d'abstractions de haute volée, mais le résultat est 
extrêmement concret, et c'est ce résultat qui nous intéresse
en pratique. 

Nous définissons en \coma une \emph{partition constructible du spectre de Zariski} par sa version duale, qui est un \sfio dans l'\agBz~\hbox{$\Zar\Abul=\Bo(\ZarA)$}.\\
En pratique, un \elt de $\Zar\Abul$ est donné par une liste double
dans l'anneau $\gA$

\snic{(a_1,\dots,a_\ell;u_1,\dots,u_m)=(I;U)}

%\sni
qui définit l'\elt suivant de $\Zar\Abul$: 

\snic{\Vi_i{\lnot\DAbul(a_i)}\vi\Vi_j\DAbul(u_j)=\lnot\DAbul(a_1,\dots,a_\ell) \vi \DAbul(u), \hbox{ où }u=\prod_j u_j.}

%\sni
\`A cet \elt $(I;U)$, est associé l'anneau $(\aqo\gA I) [1/u]$(\footnote{En \clama $\Spec(\aqo\gA I) [1/u]=\bigcap_{a\in I}\fV(a)\cap\bigcap_{v\in U}\fD(v)$, où~$\fV(a)$ désigne le \cop de $\fD(a)$.}). Un \sfio de $\Bo(\ZarA)$ peut alors être obtenu
comme résultat d'une construction arborescente qui démarre avec
la liste double  $(0;1)$ et qui autorise le remplacement d'une liste $(I;U)$
par deux listes doubles $(I,a;U)$ et $(I;a,U)$ pour un $a\in\gA$.

%:HHH rappatriee plus tot     Definition{defiAnneaunStable}
%\begin{definition}\label{defiAnneaunStable}~
%\end{definition}

%:HHH ci apres legerement modifie, grÂce aux partitions plus générales

Le \tho crucial qui suit est un corolaire du point \emph{2} du lemme~\ref{corlemPartitionSpec}.
%:     Theorem{thPartitionSpec}
\begin{theorem}\label{thPartitionSpec}
On considère une partition constructible de~$\SpecA$,
décrite comme ci-dessus par une famille $(I_k;U_k)_{k\in\lrbm}$. 
On note $\fa_k$ l'\idz~$\gen{I_k}$ et $u_k$ le produit des \elts de $U_k$.
\begin{enumerate}
\item Si $D:\gA\to\gT$ est un \sutz, et si tous les $(D/\fa_k)[1/u_k]$
sont $n$-stables, alors $D$ est $n$-stable.
\item En particulier, si chaque anneau $\gA[1/{u_k}]/{\fa_k}$ est $n$-stable (par exemple si sa \ddk
 est $<n$), alors $\gA$ est $n$-stable.
\end{enumerate}
\end{theorem}
%
%\begin{proof}
%Résulte par \recu du corolaire~\ref{corlemPartitionSpec}.
%\end{proof}
%%

\rems ~

1) Le cas paradigmatique d'anneau $n$-stable
est  donné dans le \tho précédent
lorsque chaque anneau  $\gA[1/u_i]/{\fa_i}$ est de \ddkz~$<n$.

2) Toute partition constructible de $\SpecA$ peut être raffinée
en la partition décrite par  les $2^{n}$ couples \cops 
formés à partir d'une liste finie~\hbox{$(\an)$} dans $\gA$.

3) Des constructions arborescentes analogues apparaissent au chapitre \ref{chap gen loc}
dans le cadre du \plgc de base, mais ce sont d'autres  anneaux, des  localisés  notés $\gA_{\cS(I;U)}$, qui interviennent alors.
\eoe

%%%%%%%%%%%%%%%%%%%%%%%%%%%%%%%%%%%%%%%%%%%%%%%%%%%%%%%%%%%%%%%%%%%%%%%%%%%
%%%%%%%%%%%%%%%%%%%%%%%%%%%%%%%%%%%%%%%%%%%%%%%%%%%%%%%%%%%%%%%%%%%%%%%%%%%

%--- Section{Manipulations}--------------
\penalty-2500
\section{Manipulations \elrs de colonnes}
%-----------------------------------------
\label{secManipElemCol}

\vspace{3pt}
Dans cette section nous établissons des \thos analogues dans deux contextes différents.
Le premier utilise la stabilité d'un support, le deuxième utilise la
dimension de Heitmann.

\Llec peut visualiser l'essentiel des résultats du chapitre sur le dessin
suivant, en gardant en mémoire les \thos \ref{thSerre}, \ref{thSwan},  \ref{thSwan2} et~\ref{thBassCancel2}.

Une flèche telle que $\Sdim\lora\Gdim$ est mise pour $\Sdim\gA\leq \Gdim\gA$.

\newcommand\URRU[1]%
  {\ar@{}[urr]|{\rotatebox{20}{\hbox{$\lllra$}}}^{\rotatebox{20}{#1}}}

\newcommand\DRRD[1]%
  {\ar@{}[drr]|{\rotatebox{-25}{\hbox{$\lllra$}}}_{\rotatebox{-25}{#1}}}

\newcommand\DRRU[1]%
  {\ar@{}[drr]|{\rotatebox{-20}{\hbox{$\lllra$}}}^{\rotatebox{-20}{#1}}}

\newcommand\URRD[1]%
  {\ar@{}[urr]|{\rotatebox{15}{\hbox{$\lllra$}}}_{\rotatebox{15}{#1}}}

\vspace{5pt}
\Grandcadre{
\xymatrix @R=0.2cm @C=1cm{
\Sdim\DRRU{fait \ref{factSDimGdimHdim} \emph{2}}&&&&\Hdim\DRRU{fait \ref{factKdimHdim}}
\\
                                       &&\Gdim\URRU{th. \ref{MAINCOR}}
                                         \DRRD{th. \ref{matrix}}
                                       &&&&\Kdim
\\
\Bdim\URRD{fait  \ref{factGdimBdim}}   &&&&\Cdim\URRD{def. \ref{defisutnstab}}
\\
}
}

%\vspace{-2pt}

%: SUBSEC  Avec la stabilite d'un support
\subsec{Avec la stabilité d'un \sutz}
%-----------------------------------------

%%%%%%%%%%%%%%%%%%%%%%%%%%%%%%%%%%%%%%%%%%%%%%%%%%%%%%%%%%%%%%%%%%%%%%%%%%%

\Grandcadre{Dans ce paragraphe, $D:\gA\to\gT$ est un \sut fixé}

Nous fixons les notations suivantes, analogues à celles qui sont utilisées pour définir $\Gdim\gA<n$ dans la \dfnz~\ref{defiSdimGdim}.
%:--- Notation{nota0Matrix}-------------
\begin{notation}
\label{nota0Matrix}
{\rm Soit $F=[\,C_0\,|\,C_1\,|\,\dots\,|\,C_p\,]$  une matrice dans $\Ae{m\times (p+1)}$ (les~$C_i$ sont les colonnes) et
$G=[\,C_1\,|\,\dots\,|\,C_p\,]$, 
de sorte que~$F=[\,C_0\,\vert\,G\,]$.
}
\end{notation}
%--- end-notation-----------------------------------------

Remarquons que pour tout $n$ on a $\DA\big(C_0,\cD_n(F)\big)=\DA\big(C_0,\cD_n(G)\big)$, et a fortiori $D\big(C_0,\cD_n(F)\big)=D\big(C_0,\cD_n(G)\big)$.

%:--- Lemme{cor0main}-------------
\begin{lemma}
\label{cor0main}
On suppose que  $D$ est $n$-stable et on prend la notation \ref{nota0Matrix} avec   $m=p=n$. On note $\delta=\det(G)$.
Il existe $x_1$, \dots, $x_n$ tels que

\snic{D(C_0,\delta)\leq  D\big(C_0+\delta(x_1C_1+\cdots+x_nC_n)\big).}

\end{lemma}
%--- end-corollary------------------------------------

%-----------------begin proof------------------
\begin{proof}
Il suffit de réaliser $D(\delta)\leq  D\big(C_0+\delta(x_1C_1+\cdots+x_nC_n)\big)$, \cad 

\snic{D(\delta)\leq D(C_0+\delta G X)$ pour un $X\in\Ae n.}

%\sni
Soit $\wi{G}$ la matrice adjointe de $G$ et $L = \wi{G}C_0$.
Pour n'importe quel $X\in \Ae n$, on a $\widetilde {G}(C_0+\delta GX) = L+\delta^{2} X$, donc   $\DA(L+\delta^{2} X)\leq \DA(C_0+\delta GX)$, et a fortiori~\hbox{$D( L+\delta^{2} X) \leq  D(C_0+\delta GX)$}.
Puisque  $D$ est $n$-stable, d'après le fait~\ref{lemsutnstab},
on a un  $X\in \Ae n$ tel que
$D(\delta)\leq  D(L+\delta^2X)$.\\
Donc  $D(\delta) \leq   D(C_0+\delta GX)$, comme demandé.
\end{proof}
%-----------------end proof------------------

%:--- Theorem{matrix}------------
\begin{theorem}
\label{matrix}  
\emph{(\Tho de Coquand, 1: Forster-Swan et autres avec la $n$-stabilité)} 
On a  $\Gdim\gA\leq \Cdim\gA$.
En conséquence, les \thos  \SSOz, de Forster-Swan et de simplification de Bass (\ref{thSerre}, \ref{thSwan},  \ref{thSwan2},  
\ref{thBassCancel2}) s'appliquent avec la~$\Cdim$.
\end{theorem}
%--- end-theorem-----------------------------------------
%
\begin{proof} On suppose que $\Cdim\gA< n$ et on montre que $\Gdim\gA< n$.
On utilise la \carn de  $\Gdim\gA< n$ donnée dans la proposition \ref{propGdimGdim}.
Le lemme~\ref{cor0main}, avec le support  $D=\rD_{\gA/\! \gen{d}}$, nous dit que la \prt \eqve décrite en
\ref{propGdimGdim}  est satisfaite si  $\Cdim\aqo\gA d<n$.
On conclut en notant \hbox{que $\Cdim\aqo\gA d\leq \Cdim\gA$}. 
\end{proof}
%

%:--- Theorem{matrixC}------------
\begin{theorem}
\label{matrixC}  \emph{(\Tho de Coquand, 2: \mlrs de colonnes, support et $n$-stabilité)}
Avec les notations \ref{nota0Matrix}.\\
Soit $n\in\lrbp$. Si $D$ est $n$-stable
 il existe $t_1$, \dots, $t_p\in \cD_n(G)$ tels que

\snic{
 D\big(C_0,\cD_n(G)\big)\leq  D(C_0+t_1C_1+\cdots+t_pC_p).}
 
%%\sni
%En particulier, si $1=D\big(C_0,\cD_n(G)\big)$,  
%il existe $t_1$, \dots, $t_p\in \cD_n(G)$ tels que
%le vecteur $C_0+t_1C_1+\cdots+t_pC_p$
%est $D$-\umdz.
\end{theorem}
%--- end-theorem-----------------------------------------

\vspace{.1em} La \dem de ce \tho comme conséquence du lemme \ref{cor0main}
est analogue à la \dem de l'implication difficile dans la proposition~\ref{propGdimGdim}, dans un contexte légèrement différent. Le résultat est ici plus fort
car la proposition~\ref{propGdimGdim} ne s'intéresse qu'au cas particulier
donné dans le corolaire~\ref{t0basic}, avec en \hbox{outre $D=\DA$}.
%-----------------begin proof------------------
\begin{proof} 
On doit trouver $t_1$, \ldots, $t_p$ dans $\cD_n(G)$ tels que,
 pour tout mineur $\nu$ d'ordre $n$ de $G$, on ait $D(C_0,\nu)\leq
D(C_0+t_1C_1+\cdots+t_pC_p)$.  
\\ En
fait il suffit de savoir réaliser

\snic{D(C_0,\delta)\leq  D\big(C_0+\delta (x_1C_1+\cdots+x_pC_p)\big)}

%\sni
pour \emph{un} mineur $\delta$ d'ordre $n$ de $G$,
et comme on l'a déjà remarqué, il suffit pour cela que
$D(\delta)\leq  D\big(C_0+\delta (x_1C_1+\cdots+x_pC_p)\big)$.\\
En effet dans ce cas, nous
remplaçons $C_0$ par $C'_0=C_0+\delta (x_1C_1+\cdots+x_pC_p)$ dans~$F$ (sans
changer $G$), et nous pouvons passer à un autre mineur~$\delta'$ de~$G$ pour
lequel nous obtiendrons $x'_1$, \dots, $x'_p$ vérifiant

\snic{D(C_0,\delta,\delta')\leq D(C'_0,\delta')\leq
 D\big(C'_0+\delta'(x'_1C_1+\cdots+x'_pC_p)\big)=  D(C''_0)}

%\sni
avec $C''_0=C_0+t''_1C_1+\cdots+t''_pC_p$ et ainsi de suite.
\\
Pour réaliser l'in\egt

\snic{D(\delta)\leq  D\big(C_0+\delta (x_1C_1+\cdots+x_pC_p)\big)}

%\sni
  pour un
mineur $\delta$ d'ordre $n$ de $G$, on utilise le lemme~\ref{cor0main} avec la matrice extraite~$\Gamma$ correspondant au mineur~$\delta$,
et on se limite pour $C_0$ aux lignes de~$\Gamma$, 
ce qui nous donne un vecteur $\Gamma_0$.
On obtient \hbox{un $X\in\Ae n$} tel que

\snic{
D(\delta)\leq  D(\Gamma_0+\delta \Gamma X)\leq D(C_0+\delta G Z).
}

où $Z\in \gA^{p}$ est obtenu en complétant $X$ par des $0$.
\end{proof}
%-----------------end proof------------------

Toujours avec les notations \ref{nota0Matrix}, nous obtenons comme
corolaire le résultat suivant, qui signifie, lorsque $D=\DA$,
que $\Gdim\gA\leq \Cdim\gA$.
%:--- corollary{t0basic}-------------
\begin{corollary}
\label{t0basic}
Soit $n\in\lrbp$. Si $D$ est $n$-stable et si $1=D\big(C_0,\cD_n(G)\big)$, 
il existe $t_1$, \dots, $t_p$ tels que le vecteur $C_0+t_1C_1+\cdots+t_pC_p$
est $D$-\umdz. 
\end{corollary}
%--- end-theorem-----------------------------------------

%%%%%%%%%%%%%%%%%%%%%%%%%%%%%%%%%%%%%%%%%%%%%%%%%%%%%%%%%%%%%%%%%%%%%%%%%%%
%%%%%%%%%%%%%%%%%%%%%%%%%%%%%%%%%%%%%%%%%%%%%%%%%%%%%%%%%%%%%%%%%%%%%%%%%%%
%%%%%%%%%%%%%%%%%%%%%%%%%%%%%%%%%%%%%%%%%%%%%%%%%%%%%%%%%%%%%%%%%%%%%%%%%%%
%%%%%%%%%%%%%%%%%%%%%%%%%%%%%%%%%%%%%%%%%%%%%%%%%%%%%%%%%%%%%%%%%%%%%%%%%%%
%%%%%%%%%%%%%%%%%%%%%%%%%%%%%%%%%%%%%%%%%%%%%%%%%%%%%%%%%%%%%%%%%%%%%%%%%%%
%: SUBSEC  Avec la dimension de Heitmann
\subsec{Avec la dimension de Heitmann}
%-----------------------------------------
%:     Lemma{mainlemma}----------------
\begin{lemma}
\label{mainlemma}
On considère une matrice de la forme
$$
\cmatrix{
b_0&c_1&\cdots&\cdots&c_n\cr
b_1&a&0&\cdots&0\cr
\vdots&0&\ddots&\ddots&\vdots\cr
\vdots&\vdots&\ddots& \ddots&0\cr
b_n&0&\cdots&0&a\cr
}
,$$
dont nous notons les colonnes par $V_0$, $V_1$, \ldots, $V_n$.
\\
Si $\Hdim\gA < n$ et
$1 = \DA(b_0,a)$,
alors il existe $x_1$, \dots, $x_n\in a\gA$ tels que
$$\preskip.4em \postskip.6em
1 = \DA(V_0+x_1V_1+\cdots+x_nV_n).
$$
\end{lemma}
%--- end-lemma-----------------------------------------

\vspace{.1em}
%-----------------begin proof------------------
\begin{proof}
La preuve est par \recu sur $n$.  Pour $n=0$, c'est clair.  \\
Si $n>0$, soit $\fj=\IH_\gA(b_n)$.  On a
$b_n\in \fj$ et $\Hdim\gA/\fj< n-1$, donc par \hdrz, on peut
trouver $y_1$, \dots, $y_{n-1}\in \gA$ tels que
$$\preskip.4em \postskip.4em
1 = \rD(U_0+ay_1U_1+\cdots+ay_{n-1}U_{n-1})\quad \mathrm{dans}\;\gA/\fj,\eqno (\alpha)
$$
où $U_i$ désigne le vecteur $V_i$ privé de sa dernière \cooz.\\
Posons $U'_0=U_0+ay_1U_1+\cdots+ay_{n-1}U_{n-1}$,
on a
$\DA(U'_0,a)=\DA(U_0,a)$. %% et $\DA(X,a)=\DA(a,b_1,\ldots b_{n-1})$. 
L'\egt $(\alpha)$  signifie qu'il existe $y_n$ tel que $b_ny_n\in \JA(0)$ et
$$\preskip.4em \postskip.4em
1 = \DA(%X,
U'_0) \vu \DA(b_n,y_n).\eqno (\beta)
$$
Posons $V'_0 = V_0+ay_1V_1+\cdots+ay_{n-1}V_{n-1}+ay_nV_n$. Le lemme
est prouvé si $1 \in \DA(V'_0)$. 
Remarquons que $V'_0$ privé de sa dernière \coo est le
vecteur $U'_0 + a_ny_nU_n$ et que 
sa dernière \coo est $b_n + a^2 y_n$, d'où le jeu serré qui vient
avec $b_n$, $a$, $y_n$.
On a
$$\preskip.4em \postskip.4em
\DA(U_0'+ay_nU_n) \vu \DA(a) = \DA(U'_0,a) = \DA(U_0,a)
\supseteq\DA(b_0,a) = 1,\eqno (\gamma)
$$
et, d'après $(\beta)$,
$$\preskip.4em \postskip.4em
\DA(U_0'+ay_nU_n)\vu \DA(b_n,y_n) = \DA(U_0') \vu \DA(b_n,y_n) = 1.\eqno (\delta)
$$
Ensuite $(\gamma)$ et $(\delta)$ impliquent
$$\preskip.4em \postskip.4em 
\DA(U_0'+ay_nU_n) \vu \DA(b_n,a^2y_n) = 1 =
\JA(U_0'+ay_nU_n,b_n,a^2y_n),\eqno (\eta)
$$
et d'après le lemme \ref{gcd2}, puisque $b_na^2y_n\in \JA(0)$,
$$\preskip.4em \postskip.4em
1 = \JA(U_0'+ay_nU_n,b_n+a^2y_n),$$
\cad $1=\DA(V'_0)$.
\end{proof}
%-----------------end proof------------------

%%%%%%%%%%%%%%%%%%%%%%%%%%%%%%%%%%%%%%%%%%%%%%%%%%%%%%%%%%%%%%%%%%%%%%%%%%%
\incertain{
\rems \\
1) Lorsque l'on  remplace la dimension de Heitmann par celle de Krull,
le lemme précédent admet une version \gui{uniforme} dans laquelle
les $x_i$ ne dépendent pas de $a$ (comme dans le
lemme~\ref{lemKroH}).\\
2) Le lemme ci-dessus peut être vu comme une variante raffinée
du \thrf{Bass}.
\eoe
}
%%%%%%%%%%%%%%%%%%%%%%%%%%%%%%%%%%%%%%%%%%%%%%%%%%%%%%%%%%%%%%%%%%%%%%%%%%%

%:--- Theorem{MAINCOR}------------ \Tho de Coquand
\begin{theorem}
\label{MAINCOR} \emph{(\Tho de Coquand, 3: Forster-Swan et autres avec la
dimension de Heitmann)} 
On a  $\Gdim\gA\leq \Hdim\gA$. En conséquence,  les \thos  \SSOz, de Forster-Swan et de simplification de Bass s'appliquent avec la~$\Hdim$ (\thos \ref{thSerre}, \ref{thSwan},  \ref{thSwan2},  
\ref{thBassCancel2}).
\end{theorem}
%--- end-theorem-----------------------------------------
\begin{proof}
On utilise la \carn de  $\Gdim\gA< n$ donnée 
dans la proposition~\ref{propGdimGdim}.
Le lemme~\ref{mainlemma} nous dit que la \prt \eqve décrite en
\ref{propGdimGdim}  est satisfaite si  $\Hdim\aqo\gA d<n$.
Enfin, on note \hbox{que $\Hdim\aqo\gA d\leq \Hdim\gA$}.
\end{proof}

\REM{ finale}
 Tous les
\thos
d'\alg commutative que nous avons démontrés dans ce chapitre
se ramènent en fin
de compte à des \thos concernant les matrices et leurs
\mlrsz.
\eoe

%:section: Exercices
%\newpage	
\Exercices

%\let\S\relax
%--- Exercise{exo16Lecteur}-------------
\begin{exercise}
\label{exo16Lecteur}
{\rm   
Expliciter le calcul que donne la preuve du \thrf{thKroH} dans le \linebreak 
cas $n=1$.
}
\end{exercise}
%--- end -exercise-----------------------------------------

%--- Exercise{exoRegularSequence1}-------------
\begin{exercise}\label{exoRegularSequence1} {(Une \prt des suites \ndzesz)}
\\
 {\rm  
Soit $(\an)$ une suite \ndze de $\gA$  et $\fa = \gen {\an}$ ($n\geq1$).
\\
\emph {1.}
Montrer que $(\ov {a_1}, \ldots, \ov{a_n})$ est une $(\gA\sur\fa)$-base de
$\fa\sur{\fa^2}$.
\\
\emph {2.}
En déduire, lorsque $1\notin\fa,$ que $n$ est le nombre minimum de générateurs de l'\id $\fa$.
Par exemple, si $\gk$ est un anneau non trivial et $\gA = \gk[\Xm]$,
alors pour $n \le m$, le nombre minimum de générateurs de l'\id
$\gen {\Xn}$ est $n$.
}
\end {exercise}

%--- Exercise{exoNbGensIdeal}-------------
\begin{exercise}\label{exoNbGensIdealBis}
\label{exoNbGensIdeal} (Nombre de \gtrs  de $\fa/{\fa^2}$ et de
 $\fa$)
 \\ 
{\rm
Soit $\fa$ un \itf de $\gA$ avec 
$\fa/\fa^2 = \gen {\ov {a_1}, \cdots,\ov{a_n}}$. 
\\
 \emph {1.}
Montrer que $\fa$ est engendré par $n+1$ \eltsz.
\\
 \emph {2.}
Montrer que $\fa$ est \lot engendré par $n$ \elts au sens précis
suivant: il existe $s \in \gA$ tel que sur les deux localisés $\gA_s$ et
$\gA_{1-s}$, $\fa$ est engendré par $n$ \eltsz.
\\
 \emph {3.}
En déduire que si $\gA$ est \lgb (par exemple si $\gA$
est \rdt \zedz), alors $\fa$ est engendré par $n$ \eltsz.
}
\end {exercise}
%--- end -exercise----------

%--- Exercise{exoGensPolIdeal}-------------
\begin{exercise}
\label {exoGensPolIdeal}
{\rm
\emph{1.}
Soient $E$ un \Amo et $F$ un \Bmoz.
Si $E$ et $F$ sont engendrés par $m$ \eltsz, il
en est de même du $(\gA\times\gB)$-module
$E\times F$.
\\
 \emph{2.}
Soit $\fa \subseteq \gA[X]$ un idéal contenant un \pol $P = \prod_{i = 1}^s
(X - a_i)$ \splz.  On note $\ev_{a_i} :
\gA[X] \twoheadrightarrow \gA$ le morphisme d'évaluation qui spécialise $X$ en $a_i$. On suppose que
chaque $\fa_i := \ev_{a_i}(\fa)$ est engendré par $m$ \eltsz. Montrer que
$\fa$ est engendré par $m+1$ \eltsz.
\\
 \emph{3.}
Soit $\gK$ un corps discret et $V \subset \gK^n$ un ensemble fini. 
Montrer que l'\id

\snic{\fa(V) = \sotq{ f \in \KXn}{ \forall\ w \in V,\ f(w) = 0}}

%\sni
est engendré par $n$ \elts (notez que cette borne ne dépend pas de $\#V$
et que le résultat est clair pour $n=1$).

}
\end {exercise}
%--- end -exercise-----------------------------------------

%--- Exercise{exoCubiqueGaucheP3}-------------
\begin{exercise}\label{exoCubiqueGaucheP3}
 {(La cubique gauche de $\PP^3$, image de $\PP^1$
par le plongement de Veronese de degré $3$)} 
{\rm  
L'anneau de base $\gk$ est quelconque, sauf dans la première question où
c'est un \cdiz. On définit le morphisme de Veronese $\psi : \PP^1 \to \PP^3$
par

\snic {
\psi : (u : v) \mapsto (x_0 : x_1 : x_2 : x_3) \quad\hbox {avec}\quad 
x_0=u^3,\ x_1=u^2v,\  x_2=uv^2,\  x_3=v^3
.}

%\sni
\emph{1.}
Montrer que $\Im\psi = \cZ(\fa)$ où $\fa = \gen {D_1, D_2, D_3}=\cD_2(M)$ avec 

\snic{M=\cmatrix {X_0 & X_1 & X_2\cr X_1 & X_2 & X_3\cr},}

\snic{D_1 = X_1X_3-X_2^2,\; D_2 = -X_0X_3+X_1X_2\hbox{ et } 
D_3 = X_0X_2-X_1^2.}

%\sni
\emph{2.}
Montrer que $\fa$ est le noyau de $\varphi : \gk[X_0,X_1,X_2,X_3] \to
\gk[U,V]$, $X_i \mapsto U^{3-i}V^i$.  En particulier, si $\gk$ est intègre,
$\fa$ est premier et si $\gk$ est réduit, $\fa$ est radical.
On pourra montrer qu'en posant:

\snic {
\fa^\bullet = \gA \oplus \gA X_1 \oplus \gA X_2
\quad \hbox {avec} \quad \gA = \gk[X_0,X_3],
}

%\sni
on a
$\gk[X_0,X_1,X_2,X_3] = \fa + \fa^\bullet \; \hbox { et } \;
\ker\varphi \cap \fa^\bullet = 0$.

\emph{3.}
Montrer que $\fa$ ne peut pas être engendré par deux \gtrsz.

\emph{4.}
Expliciter un \pol \hmg $F_3$ de degré $3$ tel que 
$\DA  (\fa) = \DA  ({D_1,F_3})$. En particulier,
si $\gk$ est réduit, $\fa = \DA  ({D_1,F_3})$.
}

\end {exercise}
%--- end -exercise-----------------------------------------

%--- Exercise{exoKroLocvar}-------------
\begin{exercise}
\label{exo0KroLocvar}
{\rm
Montrer que si deux suites sont disjointes
(voir \paref{propdefdisjointes}) elles restent disjointes
lorsque l'on  multiplie une des suites par un \elt de l'anneau.
}
\end{exercise}
%--- end -exercise-----------------------------------------

\entrenous{ Exercice incertain.
%--- Exercise{exosutnstab}-------------
%\begin{exercise}
\label{exosutnstab} 
{%\rm
Soit $D:\gA\to\gT$ un \sut  $n$-stable.
\\
 1. Si $\varphi:\gT\to\gT'$
est un \homo de \trdisz, alors $\varphi\circ D$ est $n$-stable.
\\
 2. Si $\fb$ est un \id de $\gB$, si $\pi_\fb:\gB\to \gB\sur{\fb}\simeq\gA$ est l'\homo
canonique, alors $D\circ \pi_\fb$ est $n$-stable.
} 
%\end{exercise}
%--- end-exercise-----------------------------------------
} % fin entrenous

%--- Exercise{exoJPFurter}-------------
\begin{exercise}\label{exoJPFurter}
{(Transitivité de l'action de $\GL_2(\gk[x,y])$ sur
les \syss de deux \gtrs de $\gen {x,y}$)} 
{\rm  Le résultat de la question \emph{1} est d\^u à 
Jean-Philippe Furter, de l'Université de La Rochelle.
\\
Soient $\gk$ un anneau, $\gA = \gk[x,y]$ et
$p$, $q \in \gA$ vérifiant $\gen {p,q} = \gen {x,y}$.

\emph {1.}
Construire une matrice $A \in \GL_2(\gA)$ telle que $A \cmatrix {x\cr y}
= \cmatrix {p\cr q}$ et $\det(A) \in \gk^\times$.

\emph {2.}
On écrit $p = \alpha x + \beta y + \dots$, $q = \gamma x + \delta y + \dots$
avec $\alpha$, $\beta$, $\gamma$, $\delta \in \gk$. 
\begin {itemize}
\item [\emph {a.}]
Montrer que $\cmatrix {\alpha & \beta\cr \gamma &\delta} \in \GL_2(\gk)$.

\item [\emph {b.}]
Soit $G \subset \GL_2(\gA)$ l'intersection de $\SL_2(\gA)$
et du noyau de l'\homo \gui{réduction modulo $\gen {x,y}$} $\GL_2(\gA) \to \GL_2(\gk)$. Le sous-groupe
$G$ est  distingué dans $\GL_2(\gA)$.
Le sous-groupe $G\,\GL_2(\gk) = \GL_2(\gk)\,G$ de $\GL_2(\gA)$
opère transitivement sur les \syss de deux \gtrs de $\gen {x,y}$.
\end {itemize}

\emph {3.}
Soient $p = x + \sum_{i+j=2} p_{ij}x^iy^j$, $q = y + \sum_{i+j=2}
q_{ij}x^iy^j$. On a $\gen {x,y} = \gen {p,q}$
\ssi les \eqns suivantes sont satisfaites:

\snic { \arraycolsep2pt
\begin{array}{rclc} 
p_{20}p_{02} + p_{02}q_{11} + q_{02}^2  &  = &  
p_{20}p_{11} + p_{02}q_{20} + p_{11}q_{11} - p_{20}q_{02} + q_{11}q_{02} & = 
\\[1mm] 
p_{20}^2 + p_{11}q_{20} + q_{20}q_{02}  & =  &  0 
\end{array}
}

%\sni
\emph {4.}
Généraliser le résultat de la question précédente.

}

\end {exercise}
%--- end -exercise-----------------------------------------

%--- Exercise{exoSdimSmithRing}-------------
\begin{exercise}\label{exoSdimSmithRing}
{(Autour des anneaux de Smith et de la $\Sdim$)} \\
{\rm Pour les notions d'anneau de Bézout strict et de Smith, voir la
section \ref{secBézout} \paref{secpfval} et les
exercices~\ref{exoAnneauBézoutStrict} et~\ref{exoSmith}. L'exercice \ref{exoSmith} donne une solution directe du point \emph{5} ci-dessous.

\emph {1.}
Si $\gA$ est un anneau de Smith, on a $\Sdim\gA \le 0$.\\
En déduire $\Sdim\ZZ$, $\Bdim\ZZ$, $\Gdim\ZZ$ et $\Cdim\ZZ$.

Dans les questions \emph {2} et \emph {3}, l'anneau $\gA$ est quelconque.

\emph {2.}
Soit $A \in \MM_2(\gA)$ et $u \in \Ae2$ un vecteur \umdz.
Montrer que $u \in \Im A$ \ssi il existe $Q \in \GL_2(\gA)$
telle que $u$ soit la première colonne de $AQ$.

\emph {3.}
Soit $A \in \MM_2(\gA)$ de rang $\ge 1$. Alors $A$
est \eqve à une matrice diagonale \ssi 
$\Im A$ contient un vecteur \umdz.

\emph {4.}
Soit $\gA$ un anneau de Bézout strict. Montrer que 
$\Sdim\gA \le 0$ \ssiz$\gA$ est un anneau de
Smith.

\emph {5.}
En déduire qu'un anneau $\gA$ est de Smith \ssi
il est de Bézout strict et si pour $a$, $b$, $c$
\com la matrice $\cmatrix {a & b\cr 0 &c\cr}$
possède un vecteur \umd dans son image. Cette
dernière condition peut s'exprimer par 
la condition dite de Kaplansky:

\snic {
1 \in \gen {a,b,c}  \Rightarrow
\hbox { il existe $p,q$ tels que $1 \in \gen {pa, pb + qc}$.}
}

Remarque: on dispose de la \carn \elr: $\gA$ est de Bézout strict
\ssi pour $a, b\in \gA$, il existe $d$ et $(a',b')$ \com
avec $a=da'$ \hbox{et $b=db'$}. Si l'on ajoute la condition de
Kaplansky ci-dessus, on obtient une \carn \elr des
anneaux de Smith.

}

\end {exercise}
%--- end -exercise-----------------------------------------

%%%%%%%%%%%%%%%%%%%%%%%%%%%%%%%%%%%%%%%%%%%%%%%%%%%%%%%%%%%%%%%%%%%%%%%%%%%
% fin des exos

%:  solutions
\penalty-2500
\sol{

%%%%%%%%%%%%%%%%%%%%%%%%%%%%%%%%%%%%%%%%%%%%%%%%%%%%%%%%%%%%%%%%%%%%%%%%%%%
\exer{exo16Lecteur} %Point 1.
La preuve donnée dit ceci.
Puisque $\Kdim \gA\leq 0$, il existe un $x_1$ tel que $b_1x_1\in\DA(0)$ et
$1\in\DA(b_1,x_1)$. A fortiori  $b_1ax_1\in\DA(0)$ et $a\in\DA(b_1,ax_1)$. Le
lemme \ref{gcd} nous dit que $\DA(b_1,ax_1)=\DA(b_1+ax_1)$, donc
$a\in\DA(b_1+ax_1)$.

%%%%%%%%%%%%%%%%%%%%%%%%%%%%%%%%%%%%%%%%%%%%%%%%%%%%%%%%%%%%%%%%%%%%%%%%%%%

\exer{exoRegularSequence1}
\emph {1.}
Soient $b_1$, \dots, $b_n \in \gA$ tels que $\sum_i \ov{b_i} \ov{a_i} = 0$ dans
$\fa\sur{\fa^2}$. Autrement dit  $\sum_i b_ia_i = \sum_i c_ia_i$ avec $c_i \in
\fa$. D'après le lemme \ref{PetitLemmeAlterne}, il existe une matrice
alternée $M \in \Mn(\gA)$ telle que $[\,b_1-c_1 \; \cdots \;
b_n-c_n\,] =  [\,a_1 \; \cdots \; a_n\,]M$.\\
 D'où $b_i - c_i \in
\fa$, et donc $b_i \in \fa$. 

 \emph{Même chose, présentée de façon plus abstraite.} On
sait qu'une \mpn du \Amo $\fa$ pour le \sgr $(\an)$ est $R_{\ua}$.  Par le
changement d'anneau de base $\gA\to\gA\sur\fa$, cela donne une \mpn nulle
($R_\ua\mod\fa $)  du $\gA\sur\fa$-module $\fa/{\fa^2}$ pour
$(\ov{a_1},\ldots,\ov{a_n})$, ce qui signifie que ce \sys
est une base.

\emph {2.}
Si $(\yp)$ est un \sgr de l'\id $\fa$, $(\ov {y_1}, \ldots, \ov {y_p})$
est un \sgr du $(\gA\sur\fa)$-module $\fa/{\fa^2}$,  libre
de rang $n$. Donc, si $p<n $,  $\gA\sur\fa$ est trivial.

%%%%%%%%%%%%%%%%%%%%%%%%%%%%%%%%%%%%%%%%%%%%%%%%%%%%%%%%%%%%%%%%%%%%%%%%%%%
\exer{exoNbGensIdeal}
\emph{1.} En posant $\fb = \gen {a_1, \cdots, a_n}$, 
l'\egt $\fa/\fa^2 = \gen {\ov {a_1}, \cdots,\ov{a_n}}$ signifie que $\fa =
\fb + \fa^2$. On a alors $(\fa/\fb)^2 = (\fa^2 + \fb)/\fb = \fa/\fb$, et
l'\itf $\fa/\fb$ de~$\gA/\fb$ est \idmz, donc
engendré par un \idmz. Il existe donc un $e \in \fa$,
\idm modulo $\fb$, tel que $\fa = \fb + \gen {e}$: $\fa=\gen{a_1, \ldots, a_n, e}$.

\emph {2.}
Avec les mêmes notations on voit %alors 
que $(1-e)\fa \subseteq \fb + \gen {e^2 - e} \subseteq \fb$.\\
Donc dans $\gA_{1-e}$, $(\an)$ engendre $\fa$ tandis que dans $\gA_e$,
$1 \in \fa$.

\emph{Variante.}
On introduit $S = 1 + \fa$ et l'on travaille sur $\gA_S$:  $\fa\gA_S\subseteq\Rad(\gA_S)$ et donc, par Nakayama, un
\sgr de $\fa_S\sur{\fa_S^2}$ est aussi un \sgr de $\fa_S$. On a donc $\fa_S =
\fb_S$, d'où l'existence d'un $s \in S$ tel \hbox{que $s\fa \subseteq\fb$}. \hbox{Dans
$\gA_s$, $(\an)$} engendre $\fa$, tandis que dans $\gA_{1-s}$, $1 \in \fa$
($s\in1+\fa$, \hbox{donc $1-s \in \fa$}).

%%%%%%%%%%%%%%%%%%%%%%%%%%%%%%%%%%%%%%%%%%%%%%%%%%%%%%%%%%%%%%%%%%%%%%%%%%%
\exer{exoGensPolIdeal}
\emph{1.} 
\'Evident, et l'on en déduit  \emph{2} puisque
$\aqo{\fa}{P} \simeq \fa_1 \times \cdots \times \fa_s$.
On en déduit \emph{3} par \recu sur $n$. On peut remarquer que le \tho chinois utilisé dans le point \emph{2} est réalisé concrètement par l'interpolation à la Lagrange.
\\
 NB: voir aussi l'exercice \ref{exoGensIdealEnsFini}.

%%%%%%%%%%%%%%%%%%%%%%%%%%%%%%%%%%%%%%%%%%%%%%%%%%%%%%%%%%%%%%%%%%%%%%%%%%%

%%%%%%%%%%%%%%%%%%%%%%%%%%%%%%%%%%%%%%%%%%%%%%%%%%%%%%%%%%%%%%%%%%%%%%%%%%%

\exer{exoCubiqueGaucheP3} 
\emph {1.}
$\psi$ est \hmg de degré $3$. Soit $p = (x_0:x_1:x_2:x_3)$ \hbox{dans $
\cZ(\fa)$}. Si $x_0 \ne 0$, on est ramené à 
$x_0 = 1$, donc $(x_0, x_1, x_2, x_3) = (1, x_1, x_1^2, x_1^3) =
\psi(1 : x_1)$. Si $x_0 = 0$, alors
$x_1 = 0$, puis $x_2 = 0$, donc $p = \psi(0 : 1)$.

\emph {2.}
Soit $\kux = \kuX\sur\fa$ et $\ov\gA = \gk[x_0,x_3]$. Montrer l'\egt $\kuX = \fa
+ \fa^\bullet$ revient à montrer que $\kux = \ov\gA +
\ov\gA x_1 + \ov\gA x_2$.  On a les relations $x_1^3 = x_0^2x_3 \in
\ov\gA$, \hbox{et $x_2^3 = x_0x_3^2 \in \ov\gA$}, donc $\ov\gA[x_1,x_2]$ est le
$\ov\gA$-module engendré par les $x_1^ix_2^j$ pour \hbox{les $i$, $j \in \lrb{0..2}$}. Mais
on a aussi $x_1x_2 = x_0x_3$, $x_1^2 = x_0x_2$, $x_2^2 = x_1x_3$, ce qui
achève de montrer que $\ov\gA[x_1,x_2] = \ov\gA + \ov\gA x_1 + \ov\gA x_2$.

Soit $h = a + bX_1 + cX_2 \in \fa^\bullet$ vérifiant $\varphi(h) = 0$
($a, b, c \in \gA = \gk[X_0,X_3]$). On a donc
$$
a(U^3,V^3) + b(U^3,V^3)U^2V + c(U^3,V^3)UV^2 = 0
.$$
En posant $p(T) = a(U^3,T)$, $q(T) = b(U^3,T)U^2$, $r(T) = c(U^3,T)U$, on
obtient l'\egt $p(V^3) + q(V^3)V + r(V^3)V^2 = 0$, et un examen modulo $3$ des exposants
en~$V$ de $p$, $q$, $r$ fournit $p=q=r=0$. D'où $a=b=c=0$, i.e. $h = 0$.
Maintenant, \hbox{si $f \in \ker\varphi$}, en écrivant $f = g + h$ avec
$g \in \fa$, $h \in \fa^\bullet$, on obtient $h \in \ker\varphi \cap
\fa^\bullet = 0$, donc $f = g \in \fa$.

\emph {3.}
Posons $E = \fa/\!\gen{\uX}\fa$. C'est un $\aqo{\kuX}{\uX}$-module engendré
par les $d_i = \ov {D_i}$. Autrement dit  $E = \gk d_1 + \gk d_2 + \gk d_3$.  De plus,
$d_1$, $d_2$, $d_3$ sont $\gk$-\lint 
indépendants. En effet, si $ad_1 + bd_2 +cd_3 = 0$, alors $aD_1 + bD_2 + cD_3 \in \gen{\uX}\fa$,
qui pour des raisons d'homogénéité donne $aD_1 + bD_2 + cD_3 = 0$, puis
$a = b = c = 0$. Donc $E$ est libre de rang $3$ sur $\gk$.  Si $G$ est un \sgr
de $\fa$, \hbox{alors $\ov G$} est un \sgr du \kmo $E$, donc $\#\ov G \ge 3$, a
fortiori $\# G \ge 3$.

\emph {4.}
Posons
$F_3 = X_0D_2 + X_1D_3 = -X_0^2X_3 + 2X_0X_1X_2 - X_1^3  \in 
\gen {D_2, D_3}$.
 On a
$$\preskip.4em \postskip.4em
\begin {array} {c}
D_2^2 = -(X_3F_3 + X_1^2D_1) \in \gen {D_1, F_3}, \quad 
D_3^2 = -(X_1F_3 + X_0^2D_1) \in \gen {D_1, F_3}, \\[1mm]
D_2D_3 = X_0X_1D_1 + X_2F_3  \in \gen {D_1, F_3} \;\;
\hbox{puis}\\[1mm]
\gen {D_1, D_2, D_3}^2 \subseteq \gen {D_1, F_3} \subseteq
\gen {D_1, D_2, D_3}, \; \hbox {d'où} \;
\sqrt {\gen {D_1, D_2, D_3}} = \sqrt {\gen {D_1, F_3}}.
\end {array}
$$

%%%%%%%%%%%%%%%%%%%%%%%%%%%%%%%%%%%%%%%%%%%%%%%%%%%%%%%%%%%%%%%%%%%%%%%%%%%

\exer{exoJPFurter}
 \emph {1.}
Remarquons d'abord que pour $m_{ij} \in \gA=\gk[x,y]$, une \egt
$$
\cmatrix {m_{11} & m_{12}\cr m_{21} & m_{22}} \cmatrix {x\cr y} =
\cmatrix {0\cr 0}
$$ 
entraîne $m_{ij} \in \gen {x,y}$. Par ailleurs,
on va utiliser les \idts suivantes pour des matrices $2\times 2$:
$\det(A + B) = \det(A) + \det(B) + \Tr(\wi AB)$ et:

\snic {
\hbox {pour } H = \crmatrix {v\cr -u} [\,y\;-x\,], \quad
\Tr(\wi AH) = [\,u\;v\,]\,A\cmatrix {x\cr y} 
.}

%\sni
Par hypothèse, on a $A$, $B \in \MM_2(\gA)$ telles que 
$$A \cmatrix
{x\cr y} = \cmatrix {p\cr q}\hbox{  et  }B \cmatrix {p\cr q} = \cmatrix {x\cr y}$$
donc $(BA - \I_2) \cmatrix {x\cr y} = \cmatrix {0\cr 0}$.  Ainsi, modulo $\gen
{x,y} = \gen {p,q}$, on a $BA \equiv \I_2$. Donc $a = \det(A)(0,0) \in \gk\eti$
et l'on peut écrire, avec $u$, $v \in \gA$, $\det(A) = a + up + vq$.  On
pose $H = \crmatrix {v\cr -u} [\,y\;-x\,]$. On a $H\cmatrix {x\cr y} =
\cmatrix {0\cr 0}$, $\det(H) = 0$, et l'on corrige $A$ en $A' = A-H$.  Alors $A'
\cmatrix {x\cr y} = \cmatrix {p\cr q}$ et 
$$
\det(A') = \det(A) + \det(H) -
\Tr(\wi AH) = a + up + vq - [\,u\;v\,] \cmatrix {p\cr q} = a.
$$
\emph {2.}
On décompose $A$ en composantes \hmgs: $A = A_0 + A_1 + \dots$, et l'on
examine l'\egt $A\cmatrix {x\cr y} = \cmatrix {p\cr q}$. \\
L'examen de la composante
\hmg de degré $1$ donne $A_0 = \cmatrix {\alpha & \beta\cr \gamma &\delta}$,
et l'on sait que $\det(A) = \det(A_0) \in \gk^\times$.  On peut alors
écrire $A_0(A_0^{-1}A)\cmatrix {x\cr y} = \cmatrix {p\cr q}$ avec $A_0 \in
\GL_2(\gk)$ et $A_0^{-1}A \in G$.

\emph {3.}
On écrit $A\cmatrix {x\cr y} = \cmatrix {p\cr q}$ avec $A \in G$.
Pour des raisons de degré, on obtient une \egt $A = \I_2 + xB + yC$ avec $B$, $C \in \MM_2(\gk)$.
On a alors:
$$ 
\begin{array}{ccc} 
 \cmatrix {p\cr q} = A\cmatrix {x\cr y} = 
\cmatrix {x\cr y} + B\cmatrix {x^2\cr xy} + C\cmatrix {xy\cr y^2}
= \\[1.2em] 
 \cmatrix {x + b_{11}x^2 + (c_{11}  + b_{12})xy + c_{12}y^2 \\[.3em]
y + b_{21}x^2 + (c_{21} + b_{22})xy + c_{22}y^2} 
\end{array}
\eqno(\star)
$$
Par ailleurs, on remarque que le \coe de $\det(A)-1$ en $x^iy^j$ est un \pog de
degré $i+j$ en les \coes de $B$ et $C$:

\snic {
\det(A)-1 = \Tr(B)x + \Tr(C)y + \det(B)x^2 + \Tr(\wi BC)xy + \det(C)y^2
.}

%\sni

Si $\gk$ était un corps \acz, on pourrait tenir le discours suivant.
L'\egt $\det(A) = 1$ définit une sous-\vrt projective $V \subset
\PP^{8-1}$ ($2 \times 4$ \coes pour $(B, C)$); d'autre part $(\star)$ définit
un morphisme $V \to \PP^{6-1}$ ($6$ pour les \coes de $p-x$, $q-y$). L'image
de ce morphisme est l'ensemble $W$ défini par les \eqns de l'énoncé.

L'anneau $\gk$ étant quelconque, on examine attentivement les \eqns $(\star)$;
en utilisant $\Tr(B) = \Tr(C) = 0$, on peut exprimer $B$ et $C$
en fonction des \coes de $p$ et $q$

\snic {
B = \cmatrix {p_{20} & p_{11} + q_{02}\cr q_{20} & -p_{20}}, \qquad
C = \cmatrix {-q_{02} & p_{02}\cr p_{20}+q_{11} & q_{02}}
.}

%\sni
On construit ainsi une section $s : W \to G$ de l'application $(\star)$, et
en fait les trois \eqns de $W$ figurant dans l'énoncé sont, au signe près,
$\det(C)$, $\Tr(\wi BC)$ \hbox{et $\det(B)$}.

%%%%%%%%%%%%%%%%%%%%%%%%%%%%%%%%%%%%%%%%%%%%%%%%%%%%%%%%%%%%%%%%%%%%%%%%%%%
\exer{exoSdimSmithRing} 
\emph {1.}
Une matrice rectangulaire \gui {diagonale} de rang
$\ge 1$ possède dans son image un \vmd (ceci pour tout anneau).  Soit $A$ une matrice de rang $\geq 1$, si $\gA$ est de Smith, $A$ est \eqve à une matrice \gui
{diagonale}~$D$, donc $\Im D$ contient un vecteur \vmdz, et \egmt
$\Im A$.

On a donc $\Sdim\ZZ = 0$. Par ailleurs, $\Cdim\ZZ \le 1$
($\ZZ$ est $2$-stable car $\ZZ$ est un anneau de Bézout
intègre). Enfin, $\Bdim\ZZ > 0$ car $1 \in \gen {2,5}$
sans que l'on puisse trouver un $x \in \ZZ$ tel que
$1 \in \gen {2 + 5x}$. \\
Bilan: $\Bdim\ZZ = \Gdim\ZZ =
\Cdim\ZZ = 1$ mais $\Sdim\ZZ = 0$.

\emph {2.}
Si $u = Av$, alors $v$ est \umdz. \hbox{Donc $v = Q\cdot
e_1$} \hbox{avec $Q \in \SL_2(\gA)$} et $u$ est la première colonne de $AQ$. L'autre
sens est immédiat.

\emph {3.}
Supposons que $\Im A$ contienne un vecteur \umdz. D'après le point \emph{2}, on a~$A \sim B$ avec $B\cdot e_1$ \umdz. Donc l'espace des lignes de
$B$ contient un vecteur de la forme $[\,1\ *\,]$. 
Le point \emph{2} pour $\tra B$ nous donne:

\snic {
\tra {B} \sim \cmatrix {1 &*\cr * &*\cr} \sim \cmatrix {1 & 0\cr 0 &*\cr}, \hbox{ diagonale.}
}

Bilan: $A$ est \eqve à une matrice diagonale. L'autre sens est \imdz.

\emph {4.}
Soit $\gA$ de Bézout strict avec $\Sdim\gA \le 0$. On montre
que toute matrice triangulaire~$M \in \MM_2(\gA)$ est \eqve à
une matrice diagonale. 
\\
On peut écrire $M = dA$ avec
$A$ de rang $\ge 1$ (car $\gA$ est de Bézout strict).
\\
Puisque
$\Sdim\gA \le 0$, $\Im A$ contient un vecteur \umd donc
est \eqve à une matrice diagonale $D$. En définitive $M \sim dD$.

\emph {5.}
Facile maintenant.

%%%%%%%%%%%%%%%%%%%%%%%%%%%%%%%%%%%%%%%%%%%%%%%%%%%%%%%%%%%%%%%%%%%%%%%%%%%

}% fin des solutions d'exos

%:   ---- Section*{references}-----------
\Biblio

Si l'on s'en tient à l'aspect \cof des résultats, l'ensemble du chapitre est essentiellement d\^u à T. Coquand,
avec parfois l'aide des auteurs du livre que vous tenez entre les mains.
Il s'agit ici d'un succès remarquable de l'approche \cov de la théorie de 
la \ddkz. Sans cette approche, il était simplement impensable d'obtenir
sous forme \cov \gnle les \gui{grands} théorèmes classiques démontrés ici.
En outre, cette approche a guidé la mise au point d'une nouvelle dimension,
que nous appelons dimension de Heitmann, gr\^ace à laquelle ont pu être
encore un peu améliorés les remarquables résultats non \noes de Heitmann,
notamment la version \gnle non \noee du \SSO et du \tho de Forster-Swan.

\smallskip
Le \tho de \KRN est usuellement énoncé
sous la forme suivante: une
variété algébrique dans $\CC^n$ peut toujours être définie par $n+1$
\eqnsz.
\\
Il a été étendu au cas des anneaux \noes par van der Waerden \cite{vW}
sous la forme suivante: dans un anneau \noe de \ddkz~$n$,
tout \id  a même nilradical qu'un \id engendré par au plus $n+1$ \eltsz.
\\
La version de Kronecker a été améliorée par divers auteurs  dans
les articles \cite[Storch]{Stor} et \cite[Eisenbud\&Evans]{E2} qui ont montré que $n$ \eqns suffisent en \gnlz.
Une \prco de ce dernier \tho est dans \cite[Coquand\&al.]{cls}.
Par ailleurs, on ne sait toujours pas
si toute courbe dans l'espace complexe de dimension 3 est ou non intersection de
deux surfaces (voir \cite{Kun}, chapitre~5).

Le lemme \ref{thCor2.2Heit} est le corolaire 2.2 de Heitmann \cite{Hei84}, 
(chez nous, la  $\Hdim$ remplace la $\Jdim$) qui débouche sur le \thref{KroH2}
(amélioration par Heitmann du \tho de Kronecker).

Le \tho de Kronecker local \ref{thKroLoc} est d\^u à Lionel Ducos
\cite{Duc08}.

\smallskip
\emph{Note concernant le \gui{stable range}}. Le \thrf{Bass} est d\^u à Bass dans le cas \noe
(avec la dimension du spectre maximal, qui dans ce cas coïncide
avec la $\Jdim$ et la $\Hdim$) et à Heitmann dans le cas
non \noe avec la $\Jdim$. Le \thrf{Bass0} est une version non \noeez,
mais avec la \ddkz,
du \thoz~\ref{Bass}.

\smallskip
\emph{Note concernant la $\Jdim$.} Dans \cite{Hei84} Heitmann introduit la $\Jdim$ pour un anneau
non \ncrt \noe comme le bon substitut à la dimension du
spectre maximal $\Max\gA$. C'est la dimension de $\Jspec\gA$, le
 plus petit sous-espace spectral
de $\SpecA$ contenant $\Max\gA$.
Il établit le \tho \gui{stable range} de Bass pour cette dimension.
Par contre pour les \thos  de Serre et de Forster-Swan, il doit utiliser
une dimension ad hoc,  la borne supérieure des $\Jdim(\gA[1/x])$
pour $x\in\gA$. Comme cette dimension ad hoc est de toute manière
majorée par la \ddkz, il obtient en particulier
une version non \noee des
grands \thos cités avec la \ddkz.

\smallskip \emph{Note concernant les \thos de Serre et Forster-Swan.} Le \tho de Serre
est dans \cite[Serre]{Serre}.
Le \tho de Forster-Swan (version \noeez) est dans \cite[Forster]{Forster}
pour la \ddk et dans \cite[Swan]{Swan} pour la dimension du spectre maximal.
Des versions non-\noees pour la dimension de Krull sont dues à Heitmann
\cite{Hei76,Hei84}.
Enfin l'article d'Eisenbud-Evans \cite{Eisenbud} a beaucoup
aidé à clarifier les choses concernant les
\thos de Serre, Forster et Swan.
\\
Les  sections \ref{secSOSFSa} et \ref{secManipElemCol} (deuxième partie: dimension de Heitmann) s'inspirent des grandes lignes des articles~\cite[Eisenbud\&Evans]{Eisenbud} et \cite[Heitmann]{Hei84}. Elles donnent des versions \covs des \thos
de Serre (Splitting-off), Forster-Swan, et Bass (\tho de simplification).
Ceci améliore (même sans tenir compte de l'aspect \cof de la \demz) tous les \thos connus sur la question, en répondant positivement pour la dimension de Heitmann (et a fortiori pour la $\Jdim$),  à une question laissée ouverte par Heitmann.

\smallskip
\emph{Note concernant la $\Hdim$}. 
La dimension de Heitmann, notée $\Hdim$, a été introduite dans \cite{clq}
(voir aussi \cite{clq2}). C'est elle au fond qui fait fonctionner
les \dems dans l'article de Heitmann \cite{Hei84}.
Le fait qu'elle soit meilleure a priori
que la $\Jdim$ n'est pas le point essentiel.
C'est bien plutôt le fait que les \thos de Serre et de Forster-Swan
passent pour la $\Hdim$, et donc a fortiori pour la $\Jdim$, ce qui donne
la version non \noee complète de ces \thosz, laquelle avait été conjecturée
par Heitmann.
\\
Dans le cas d'un anneau \noez,
la $\Hdim$, la $\Jdim$ de Heitmann ainsi que la dimension du spectre maximal
 $\Max \gA$ qui
intervient dans les \thos de Serre et de Swan \cite{Swan} sont les mêmes
(cf. \cite{clq2,Hei84}).

\smallskip
\emph{Note concernant la $n$-stabilité}. 
La notion de \sut remonte à Joyal \cite{Joyal} et Espa\~ nol \cite{espThesis}, qui l'utilisent
pour donner une \carn constructive de la \ddk des anneaux commutatifs.
Elle est utilisée de manière systématique dans les articles
récents de T.~Coquand.
Dans la section \ref{secSUPPORTS} et
dans la première partie de la section \ref{secManipElemCol}
la notion de \sut $n$-stable est décisive.
Elle a été inventée par
 T.~Coquand \cite{coq07} pour mettre à jour le contenu \cof du discours
 de Bass sur les partitions finies  de $\SpecA$ dans \cite{BASS}.

\smallskip
La version du \tho de simplification de Bass pour la $\Hdim$ a d'abord
été démontrée par L. Ducos \cite{DQ}. La preuve que nous donnons
est plutôt basée sur \cite{clq2}.

\smallskip En ce qui concerne l'exercice \ref {exoJPFurter}, Murthy, dans
\cite {Murthy03}, a prouvé le résultat \gnl suivant. Soient $\gA =
\gk[\xm]$ un anneau de \pols ($\gk$ anneau commutatif) et $r \ge 1$ fixés.
Supposons, pour tout $n \in \lrbr$, que tout vecteur \umd de $\Ae n$ soit
complétable et considérons, \hbox{pour $n \le \inf(r,m)$}, l'ensemble des
\syss de $r$ \gtrs de l'\id $\gen {\xn}$ de $\gA$, comme par exemple
$(\xn, 0, \ldots, 0)$ où il y a $r-n$~zéros. Alors le \hbox{groupe $\GL_r(\gA)$}
opère transitivement sur cet ensemble (le résultat de Murthy est en fait
bien plus précis).

\newpage \thispagestyle{CMcadreseul}
\incrementeexosetprob

%:        %%%%%%%%%%%%%%%%%%%%%%%%%%%%%%%%%%%%
%:        %%%%%%%%%%%%%%%%%%%%%%%%%%%%%%%%%%%%
%---- Chapitre  {le principe local-global}-{ChapPlg}-

\chapter{Le \plgz}
\label{chapPlg} \label{chap gen loc}\relax
\minitoc

\subsection*{Introduction}
\addcontentsline{toc}{section}{Introduction}
%-----------------------------------------

Dans ce chapitre, nous discutons quelques méthodes importantes
directement reliées à ce qu'il est convenu d'appeler  le \plg
en \alg commutative.

Dans la section \ref{subsec loc glob conc}
nous le développons sous la forme de
\plgcsz. Il s'agit de dire que certains \prts
sont vraies globalement
dès qu'elles le sont localement.
Ici localement est pris au sens \cofz: après \lon en
un nombre fini de \mocoz.

Dans la section \ref{subsec loc glob abs}, nous établissons les \plgas correspondants, en
utilisant, de manière inévitable, des preuves non \covsz: ici localement est pris au sens
abstrait, \cad après \lon en n'importe quel \idepz.

Dans la section \ref{secColleCiseaux}, nous expliquons la construction d'objets
\gui{globaux} à partir d'objets de même nature définis uniquement
de manière locale.

\smallskip
Les sections \ref{secMachLoGlo}, \ref{subsecLGIdeMax} et \ref{subsecLGIdepMin} sont consacrées au \gui{décryptage dynamique et \cofz}
de méthodes utilisées en \alg abstraite.
Rappelons que nous avons présenté  dans la section \ref{subsecDyna}
la philosophie \gnle de cette méthode dynamique.

Dans la section \ref{secMachLoGlo}, nous discutons le décryptage \cof de
 méthodes abstraites qui rentrent dans un
cadre \gnl du type \gui{\plgz}. Nous donnons un énoncé \gnl
(mais inévitablement un peu informel) pour cela, et nous donnons des exemples
simples, qui pourraient être traités de manière plus directe.
Les exemples vraiment pertinents arriveront dans le chapitre \ref{ChapMPEtendus}.

Cette méthode dynamique est un outil fondamental de l'\alg \covz.
On aurait pu écrire l'ouvrage présent en commençant par
cette explication préliminaire et en utilisant de manière systématique
ce décryptage.
Nous avons préféré commencer par développer
tout ce qu'il était possible de faire de manière directe,
en établissant les
\plgcs qui permettent la plupart du temps d'éviter l'utilisation
du décryptage dynamique en tant que tel.
Bref, plutôt que mettre en avant la magie à l'{\oe}uvre dans l'\alg
classique nous avons préféré faire voir d'abord un autre type de
magie, à l'{\oe}uvre en \alg \covz, sous le slogan \gnlz:
\gui{pourquoi faire compliqué quand on peut faire simple?}.
%Confessons que
%nous pensons que la plupart des livres d'\alg classique adoptent
%le slogan directement contraire, sans l'afficher trop ouvertement.

%Dans la section \ref{chapApPlg}, nous donnons des exemples d'application
%des méthodes \covs \gnles exposées auparavant.

Dans la section \ref{subsecLGIdeMax}, nous analysons la méthode   d'\alg abstraite, qui consiste à \gui{aller voir ce qui se passe
lorsque l'on quotiente par un \idema arbitraire}.

Dans la section \ref{subsecLGIdepMin}, nous analysons la méthode
qui consiste à \gui{aller voir ce qui se passe
lorsque l'on localise en un \idemi arbitraire}.

Dans les sections \ref{secPlgcor} et \ref{secPlgprof2}, nous examinons
dans quelle mesure certains \plgs restent valides lorsque l'on remplace 
dans les énoncés les listes d'\eco 
par des listes de profondeur $\geq 1$ ou de profondeur $\geq 2$.

%--- SUBsection{subsecMoco}------------------
\section{Monoïdes comaximaux, recouvrements}
\label{subsecMoco}
%-----------------------------------------

Nous traiterons dans la section \ref{subsec loc glob conc} des versions concrètes de principes du type local-global.
Pour ces versions concrètes, la \lon peut être réclamée en un nombre fini
d'\eco (ou de \mocoz) de $\gA$:
\emph{si la \prt considérée est vraie après
\lon en un nombre fini d'\ecoz, alors elle est vraie.}

Nous  introduisons une \gnnz.
%--- Definition{def.moco} ------
\begin{definition}
\label{def.moco}\relax 
On dit que
{\em les \mos $S_1$, \ldots, $S_n$ de l'anneau $\gA$ recouvrent le \moz~$S$}
si $S$ est contenu dans le saturé de chaque $S_i$ et si un \id de $\gA$ qui
coupe
chacun des $S_i$  coupe toujours $S$, autrement dit si l'on~a:

\snic{\forall s_1\in S_1 \;\dots\; \forall s_n\in S_n \;\;
\exists$ $a_1$, \ldots, $a_n\in \gA\quad \sum_{i=1}^{n} a_i s_i \in S.
}
\end{definition}
%--- end-definition------------

Des \mos sont \com s'ils recouvrent le \moz~$\so{1}$.

En \clama (avec l'axiome de l'\idepz{\footnote{L'axiome de l'\idep affirme que tout idéal strict d'un anneau est contenu dans un \idepz. Il s'agit d'une version affaiblie de l'axiome du choix.
Dans la théorie
des ensembles classique $\ZF$,
l'axiome du choix équivaut à l'axiome de l'\idemaz,
qui affirme que tout idéal strict d'un anneau est contenu dans un
\idemaz. Il est un peu plus fort que l'axiome de l'\idepz.
Ce dernier équivaut au fait que toute
théorie formelle cohérente admet un modèle (c'est le \tho de compacité en
logique classique).
Dans la théorie des ensembles avec axiome du choix, l'axiome de l'\idep devient un \tho et s'appelle \gui{lemme de Krull}.% 
\index{Lemme de Krull}\index{axiome de l'\idepz}}})
on a la \carn donnée dans le lemme qui suit.
Pour un \mo $S$, nous notons $U_S$ la partie de $\Spec \gA$ définie~par

\snic{U_S=\sotq{\fp\in \Spec \gA}{\fp\cap S=\emptyset}.}

%\sni
Si $S$ est le \mo engendré par l'\elt $s$, nous notons $U_s$ pour $U_S$.
Du point de vue \cofz, $\Spec\gA$ est un espace topologique connu via ses
ouverts de base $U_s=\fD_\gA(s)$ mais dont les points sont souvent difficilement
accessibles.

Rappelons que l'on  note $\sat{S}$ le saturé du \mo $S$.

%: --- Lemmac{lemdefmoco}--------------
\goodbreak
\begin{lemmac}
\label{lemdefmoco}~
\begin{enumerate}
\item Pour tout \mo $S$ on a $\sat S = \bigcap_{\fp\in U_S}(\gA\setminus \fp)$.
En conséquence pour deux \mos $S$ et $T$, $\sat S\subseteq\sat{T}\Leftrightarrow U_{T}\subseteq U_S$.
\item $S_1$, \ldots, $S_n$ sont \com \ssi  $\Spec \gA=\bigcup_iU_{S_i}$.
\item $S_1$, \ldots, $S_n$  recouvrent le \mo $S$ \ssi $U_S=\bigcup_iU_{S_i}$.
\end{enumerate}
\end{lemmac}
%--- end-lemma-----------------------------------------
%-----------------begin proof------------------
\begin{proof} \emph{1.} Résulte du lemme de Krull:  si un \id $\fa$ ne coupe pas un \mo $S$, il existe un \idep $\fp$ tel que
$\fa\subseteq\fp$ et $\fp\cap S=\emptyset$.\\
\emph{2.} On peut supposer que $\gA$ n'est pas trivial.
Si les \mos sont \com et si $\fp$ est un \idep
n'appartenant à aucun des $U_{S_i}$, il y a dans chaque~$S_i$ un \elt $s_i$ de
$\fp$, donc par la \dfn des \mocoz,~$1\in\fp$, ce qui est absurde.
Réciproquement, supposons que l'on ait  $\Spec \gA=\bigcup_iU_{S_i}$, et soient $s_1\in
S_1$, \ldots, $s_n\in S_n$. Si $\gen{s_1,\ldots ,s_n}$ ne contient pas~$1$, il est
contenu dans un \idep $\fp$. Donc $\fp$ n'est dans aucun \hbox{des  $U_{S_i}$}, ce qui
est absurde.
\end{proof}
%-----------------end proof------------------

Le lemme ci-après est une variation sur le thème: \index{recouvrement}%
\emph{un recouvrement de recouvrements est un recouvrement}.
C'est aussi une \gnn du fait~\ref{factLocCas}. Les calculs
correspondants sont immédiats.
En \clama ce serait encore plus rapide via le
lemme\eto\ref{lemdefmoco}.
%: --- Lemma{lemAssoc}-------
%\goodbreak
\begin{lemma} \index{Lemme des \lons successives, 2}
\label{lemAssoc} \emph{(Lemme des \lons successives, 2)}
\begin{enumerate}
\item \emph{(Associativité)} Si les \mos $S_1$, \ldots, $S_n$ de l'anneau $\gA$
recouvrent le \mo $S$  et si chaque $S_\ell$ est recouvert par des \mos
$S_{\ell,1},\ldots ,S_{\ell,m_\ell}$, alors les $S_{\ell,j}$ recouvrent $S$.
\item \emph{(Transitivité)}
\begin{enumerate}
\item Soit  $S$ un \mo de  l'anneau $\gA$ et $S_1,\ldots
,S_n$ des \mos de  l'anneau $\gA_S$.
Pour $\ell\in\lrbn$ soit $V_\ell$ le \mo
de  $\gA$ formé par les numérateurs des \elts de  $S_\ell$. Alors
 les \mos $V_1,\ldots ,V_n$ recouvrent  $S$ \ssi les \mos  $S_1$, \ldots, $S_n$ sont
\comz.
\item Plus \gnlt soient
$S_0,\ldots ,S_n$ des \mos de l'anneau  $\gA_S$ et pour $\ell=0,\ldots,n$ soit
$V_\ell$ le \mo de  $\gA$ formé par les
numé\-ra\-teurs des \elts de  $S_\ell$.
Alors les \mos $V_1,\ldots ,V_n$ recouvrent  $V_0$ \ssi  $S_1$, \ldots, $S_n$
recouvrent $S_0$ dans $\gA_S$.
\end{enumerate}
\end{enumerate}
\end{lemma}
%--- end-lemma-----------------

%--- Definotation{nota mopf} ---
\begin{definota}
\label{nota mopf}\relax
Soient $U$ et $I$ des parties de l'anneau $\gA$. Nous notons $\cM(U)$ le \mo engendré par  $U$,
 et  $\,\cS(I,U)$ est le \moz:

\snic{\cS(I,U)= \gen {I}_\gA + \cM(U).}

%\sni
Le couple $\fq=(I,U)$ est encore appelé un \emph{\ippz}, et l'on note
(par abus) $\gA_{\fq}$ pour $\gA_{\cS(I,U)}$.
De la même manière on note:

\snic{
\cS(a_1,\ldots,a_k;u_1,\ldots,u_\ell) = \gen{a_1,\ldots ,a_k}_\gA + \cM(u_1,\ldots,u_\ell). 
}

%\sni
Nous disons qu'un tel \mo {\em admet une description finie}.
Le couple

\snic{(\so{a_1,\ldots,a_k},\so{u_1,\ldots,u_\ell})}

%\sni
est appelé
un \emph{\ipp fini}.\index{ideal@idéal!premier potentiel}%
\index{ideal@idéal!premier potentiel fini}
\end{definota}
%--- end-notation-----------------

Il est clair que pour $u=u_1\cdots u_\ell$, les \mos
$\cS(a_1,\alb\ldots,\alb a_k;u_1,\ldots,u_\ell)$ \hbox{et $\cS(a_1,\ldots,a_k;u)$} sont
\eqvsz, i.e. qu'ils ont même saturé.

\medskip \rem
L'\ipp  $\fq=(I,U)$ est fabriqué dans le but suivant: \emph{lorsqu'on localise
en $\cS(I,U)$, on obtient $U\subseteq \gA_\fq\eti$  et $I\subseteq \Rad(\gA_\fq)$.} \\
De même, pour tout \idep $\fp$ tel que $I\subseteq\fp$ et $U\subseteq \gA\setminus \fp$, on a $U\subseteq \gA_\fp\eti$  et $I\subseteq \Rad(\gA_\fp)$. Le couple $\fq=(I,U)$ représente donc une information
partielle sur un tel \idepz. Il peut être considéré comme une approximation de $\fp$. Ceci explique la terminologie d'\ipp
et la notation $\gA_{\fq}$.\\ 
On peut comparer les approximations de $\fp$
par des \ipps finis aux approximations d'un nombre réel par des intervalles rationnels.
\eoe

%: --- Lemma{lemRecouvre}---------
\begin{lemma} \index{Lemme des \lons successives, 3}
\label{lemRecouvre} \emph{(Lemme des \lons successives, 3)}\\
Soient $U$ et $I$ des parties de l'anneau $\gA$ et $a\in \gA$,
alors les \mos
$$\preskip.1em \postskip.4em
\cS(I;U,a)\eqdefi\cS(I,U\cup\so a)\;\hbox{  et  }\;\cS(I,a;U)\eqdefi\cS(I\cup\so a,U)
$$
recouvrent le \mo $\cS(I,U)$.\\
En particulier, les \mos 
$S=\cM(a)=\cS(0;a)$ et
$S'=\cS(a;1)= 1+a\gA$  sont  \comz.
\end{lemma}
%--- end-lemma-----------------
%---------begin proof----------
\begin{proof}
Soient $x\in \cS(I;U,a)$, $y\in \cS(I,a;U)$. 
Il faut voir que $\gen {x,y}$
rencon-\linebreak 
tre $\gen {I} + \cM(U)$, ou encore que $\gen {x,y} + \gen {I}$
rencontre $\cM(U)$.  \\
On a $k \ge 0$, $u, v \in \cM(U)$ et $z \in \gA$ tels que
$x \in ua^k + \gen {I}$ et $y \in v - az + \gen {I}.$
Modulo $\gen {x,y} + \gen
{I}$, $ua^k \equiv 0$, $v \equiv az$ donc $uv^k \equiv 0$, i.e. $uv^k \in \gen
{x,y} + \gen {I}$ avec $uv^k \in \cM(U)$.
%et l'on écrit $z^kx+u_1y_2y=u_1u_3+u_1j_3+j_1z^k=u_4+j_4$.
\end{proof}
%---------end proof----------

\comm
Le lemme précédent est fondamental. Il est la contrepartie \cov de la
constatation banale suivante en \clamaz:
après que l'on ait localisé en un \idep
tout \elt se retrouve être \iv ou bien dans le radical.
Quand on a affaire à ce genre d'argument dans une preuve classique, on peut la
plupart du temps l'interpréter \cot au moyen de ce lemme. Sa preuve est très
simple, à l'image de la banalité de la constatation faite dans la preuve
classique. Mais ici il y a un vrai calcul. On peut
d'ailleurs se demander si la preuve classique évite
ce calcul. Une analyse détaillée montre que non: il se trouve
dans la preuve du lemme\eto\ref{lemdefmoco}.
\eoe

\medskip Les exemples donnés dans le lemme suivant sont fréquents.
%:   --- lemma{exaMoco} ---------
\begin{lemma}
\label{exaMoco}
Soit $\gA$ un anneau,  $U$ et $I$ des parties de  $\gA$, et $S=\cS(I,U)$.
%---------begin item----------
\begin{enumerate}
\item  Si $s_1$, \ldots, $s_n\in \gA$  sont des \ecoz, les \mosz~$\cM(s_i)$
sont \comz.
 Plus \gnltz, si $s_1$, \ldots, $s_n\in \gA$ sont des \eco dans $\gA_S$, les \mos
$\cS(I;U,s_i)$ recouvrent le \mo $S$.

\item  Soient $s_1$, \ldots, $s_n\in \gA$. Les \mosz:
%
%\snic{
$$\begin{array}{l}
S_1=\cS(0;s_1), \,\,S_2=\cS(s_1;s_2),\,\,S_3=\cS(s_1,s_2;s_3),\, \ldots ,  \\[1mm]
S_n=\cS(s_1,\ldots,s_{n-1};s_n)\, \, \, \mathrm{et}\,  \,\,
S_{n+1}=\cS(s_1,\ldots,s_{n};1)
\end{array}
$$
%}
%
%\sni
sont \comz.\\
Plus \gnltz, les \mosz:
$$
%\snic{
\begin{array}{l}
V_1=\cS(I;U,s_1),\,\, V_2=\cS(I,s_1;U,s_2),\,\, V_3=\cS(I,s_1,s_2;U,s_3),\,\,
\ldots ,\\[1mm]
 V_n=\cS(I,s_1,\ldots,s_{n-1};U,s_n)\, \,\, \mathrm{et}\,  \, \,
V_{n+1}=\cS(I,s_1,\ldots,s_{n};U)
\end{array}
$$
%}
%
%\sni
recouvrent le \mo $S=\cS(I,U)$.

\item 
Si $S$, $S_1$, $\ldots$, $S_n\subseteq \gA$ sont des \moco et si  $a\in\gA$, alors les \mos
$\cS(I;U,a),$ $\cS(I,a;U)$, $S_1$, $\dots$, $S_n$ sont \comz.
\end{enumerate}
%---------end item----------
\end{lemma}
%--- end-example-----------------
%---------begin proof----------
\begin{proof}
Les points \emph{2} et \emph{3} résultent \imdt des
lemmes \ref{lemAssoc} et~\ref{lemRecouvre}. 
\\
\emph{1.} Le premier cas résulte du fait que pour $k_1$, \ldots, $k_n \ge 1$, 
on a,
pour $k$ assez grand, $\gen {s_1, \ldots, s_n}^k \subseteq \geN {s_1^{k_1},
\ldots, s_n^{k_n}}$; on peut prendre  $k = \sum_i
(k_i-1) + 1$.
\\
Pour le cas \gnlz, soient $t_1$, \ldots, $t_n$ avec $t_i \in \cS(I;U, s_i)$; on
veut montrer que $\gen {t_1, \ldots, t_n}$ rencontre $S = \cS(I,U)$. Par
\dfnz, il y a un $u_i \in \cM(U)$ et $k_i \ge 0$ tels que $t_i \in
u_is_i^{k_i} + \gen {I}$; en posant $u = u_1 \cdots u_n \in \cM(u)$,
on obtient $us_i^{k_i} \in \gen {t_i} + \gen {I} \subseteq \gen {t_1, \ldots, t_n} +
\gen {I}$. Donc pour $k$ assez grand:

\snic {
u\gen {s_1, \ldots, s_n}^k \subseteq u\geN {s_1^{k_1}, \ldots, s_n^{k_n}}
\subseteq \gen {t_1, \ldots, t_n} + \gen {I}.
}

%\sni
Mais comme $s_1$, $\ldots$, $s_n$ sont des \eco dans $\gA_S$, il y a un~\hbox{$s \in S$}
tel que $s \in \gen {s_1, \ldots, s_n}$; donc $us^k \in \gen {t_1, \ldots,
t_n} + \gen {I}$, \cade  $\gen {t_1, \ldots, t_n}$ rencontre $us^k + \gen
{I} \subseteq S$.
\end{proof}
%---------end proof----------

\penalty-2500
%--- Section{subsec loc glob conc}----
\section{Quelques \plgcsz}
\label{subsec loc glob conc}\relax
%-------------------
\vspace{4pt}
%:--- SUBsection{subsecLGSLI}-----------------
\subsec{Systèmes linéaires}
\label{subsecLGSLI}
%-----------------------------------------

Le \plgc suivant est une légère \gnn du \plgrf{plcc.basic}
(\plgc de base), qui
ne concernait que le point \emph{4} ci-dessous dans le cas de modules libres de rang fini.\iplg 
 En fait l'essentiel a déjà été donné dans le \plgrf{plcc.basic.modules} (\plgc pour les modules). Nous redonnons les \dems pour insister sur leur grande simplicité.

Soient $M_1$, $\dots$, $M_\ell$, $P$ des \Amosz. Nous disons qu'une applica-\linebreak 
tion
$\Phi:M_1\times \cdots \times M_\ell\to P$ est
\ixc{homogène}{application ---}
s'il existe des entiers $r_1$, $\ldots$, $r_\ell$ tels que l'on ait
identiquement
\hbox{$\Phi(a_1x_1,\ldots ,a_\ell x_\ell)=a_1^{r_1}\cdots a_\ell^{r_\ell}\Phi(x_1,\ldots
,x_\ell)$}. Dans un tel cas, l'application $\Phi$ \gui{passe aux \lonsz}:
elle
peut être étendue naturellement en une application 

\snic{\Phi_S:S^{-1}M_1\times
\cdots \times S^{-1}M_\ell\to S^{-1}P}

%\sni
pour n'importe quel \mo $S$.
Le prototype d'une application \hmg est une application donnée par des
\pols \hmgs en les
\coos lorsque les modules sont libres de rang fini.

%: --- Prc lgc {plcc.sli}
\begin{plcc}
%\label{plcc.ring}\relax
\label{plcc.sli}\relax
Soient $S_1$, $\dots$, $S_n$  des \moco de $\gA$, $M$, $N$, $P$ des \Amosz, $\varphi$, $\psi$
des \alis  de~$M$ dans~$N$, $\theta:N\to P$ une \aliz, et $x$, $y$ 
des \elts de~$N$. On note
$\gA_i$ pour $\gA_{S_i}$, $M_i$ pour $M_{S_i}$ etc.
Alors on a les \eqvcs suivantes.
%---------begin item----------
\begin{enumerate}
\item  Recollement concret des \egtsz:
$$\preskip.4em \postskip.4em 
x=y\quad   \hbox{ dans} \quad   N  \quad  \Longleftrightarrow
\quad
\forall i\in\lrbn\;\; x/1= y/1  \quad   \hbox{ dans} \quad   N_{i}. 
$$

\item   Recollement concret des \egts d'\alisz:
\[\preskip.4em \postskip.4em 
\begin{array}{ccc} 
 \varphi =\psi \quad   \hbox{ dans} \quad   \Lin_\gA(M,N)  \qquad
\Longleftrightarrow    \\[.3em] 
\quad
\forall i\in \lrbn\;\; \varphi/1= \psi/1  \quad   \hbox{ dans} \quad
\Lin_{\gA_{i}}(M_{i},N_{i}). \end{array}
\]

\item   Recollement concret des \elts \ndzsz:
\[\preskip.4em \postskip.4em 
\begin{array}{ccc} 
x \hbox{ est \ndz dans }  N \qquad \Longleftrightarrow   \\[.3em] 
\forall i\in \lrbn
\;\;x/1\; \hbox{ est \ndz dans } N_{i}.
\end{array}
\]

\item   Recollement concret des solutions de \slisz:
$$\preskip.4em \postskip.4em 
x\in\Im \varphi  \quad  \Longleftrightarrow \quad
\forall i\in \lrbn\;\; x/1\in\Im \varphi_{i} 
$$

%5
\item   Recollement concret des solutions de \slis sous conditions
\hmgsz. 
 Soit  $(\Phi_\ell)$  une famille finie d'applications \hmgs
$$\preskip.2em \postskip.2em \mathrigid 2mu
\Phi_\ell:\Lin_\gA(M,N)\times N\to Q_\ell,\hbox{ ou }\Phi_\ell:\Lin_\gA(M,N)\to
Q_\ell ,\hbox{
ou  }\Phi_\ell: N\to Q_\ell. 
$$
Alors:
\[\preskip-.2em \postskip.4em 
\begin{array}{ccc} 
\left(\big(\&_\ell \; \Phi_\ell(\varphi ,y)=0\big)\;\Rightarrow \;
 y\in\Im \varphi\right) \quad \quad \Longleftrightarrow     
                        \\[.3em] 
\forall i\in \lrbn\;
\left(\big(\&_\ell \; \Phi_\ell(\varphi ,y)
=_{Q_{\ell,i}}
0\big)\;\Rightarrow \;
 y/1\in\Im \varphi_{i} \right).
 \end{array}
\]
où l'on a noté $Q_{\ell,i}$ pour  $(Q_\ell)_{S_i}$.

\item   \label{plcc.sex}\relax Recollement concret des suites exactes.
La suite

\snic{
M\vvers{\varphi}N\vvers{\theta}P}

est exacte \ssi les suites
$$\preskip.2em \postskip.4em 
M_{i}\vvers{\varphi_{S_i}}N_{i}\vvers{\theta_{S_i}}P_{i} 
$$
sont exactes pour $i\in \lrbn$.

\item   Recollement concret de facteurs directs dans les \mpfsz.
Ici  $M$ est un sous-\mtf d'un \mpf $N$.
\[\preskip.2em \postskip.4em 
\begin{array}{ccc} 
M  \hbox{  est facteur direct dans  } N 
 \iff \\[.3em] 
\forall i\in \lrbn, \,  M_{i}\hbox{  est facteur direct dans }  N_{i}.
 \end{array}
\]
\end{enumerate}
\end{plcc}
%----begin{proof----------
\begin{proof}
Les conditions sont \ncrs en raison du fait \ref{fact.sexloc}.
Une vérification directe est d'ailleurs \imdez.
Prouvons que les conditions locales sont suf\-fisan\-tes.

\emph{1.} Supposons que $x/1=0$ dans chaque $N_{i}$. 
Pour des $s_i\in S_i$ convenables on~a donc~$s_ix=0$ dans $N$.
Comme $\sum_{i=1}^{n} a_i s_i =1$,  on obtient $x=0$ dans~$N$.

\emph{2.} Conséquence \imde de \emph{1.}

\emph{3.} Supposons que $x/1$ soit \ndz   dans chaque $N_{i}$. 
Soit $a\in \gA$ avec
$ax=0$ dans $\gA$, donc aussi $ax/1=0$ dans chaque~$N_{i}$. On a donc~$a/1=0$ dans chaque
$\gA_{i}$, donc aussi dans~$\gA$.

\emph{4.}
Supposons que l'équation $\varphi (z)=x$ admette une solution $z_i$  dans chaque~$M_{i}$. On peut écrire $z_i=y_i/s_i$ avec $y_i\in M$ et $s_i\in S_i$. On a donc~$u_i\varphi(y_i)=s_iu_ix$ dans $N$ avec $u_i\in S_i$. 
Comme $\sum_{i=1}^{n} a_i s_iu_i =1$, on pose~$z=\sum_{i=1}^{n} a_i u_iy_i$ et l'on obtient $\varphi(z)=x$ dans~$N$.

\emph{5.} C'est une simple variante de \emph{4},  l'homogénéité des $\Phi_\ell$
intervient pour que la \prt locale soit bien définie, et pour qu'elle
résulte de la \prt globale.

\emph{6.} C'est un cas particulier du point précédent.

\emph{7.} Soit $\rho:N\rightarrow N/M$ la \prn canonique. Le
module $N/M$ est \egmt un \mpfz.  Le module $M$ est facteur direct 
dans~$N$ \ssi $\rho$ est \iv à droite, on peut donc conclure par le
\plgrf{plcc.scinde}.
 \end{proof}
%----end{proof----------
\rdb

\rem  \label{remplgcsli}
On peut voir que le point \emph{5},
simple variante du point \emph{4}, implique tous
les autres comme cas particuliers. Par ailleurs, le point \emph{4} résulte du point~\emph{1} avec $y=0$ en considérant le module
$\big(N/\varphi(M)\big)_{S_i}\simeq N_{S_i}/\varphi_{S_i}(M_{S_i})$.
On aurait donc pu énoncer le point \emph{1}  comme seul principe de
base et en déduire les points \emph{2} à \emph{6} comme corolaires.
Enfin le point \emph{7} résulte aussi directement du point \emph{4}
(voir la \dem du \plgz~\ref{plcc.scinde})
\eoe

%:--- SUBsection{subsecLGFiMo}----------------
\subsec{Propriétés de finitude pour les modules}
\label{subsecLGFiMo}
%-----------------------------------------

Les \prts de finitude usuelle des modules ont un \crc local.
La plupart ont déjà été démontrées, nous récapitulons.
%:--- Prc lgc {plcc.ptf} ---
\begin{plcc}
\label{plcc.ptf}\relax
{\em  (Recollement concret de \prts de finitude pour les modules)}
Soient $S_1$, $\dots$, $S_n$  des \moco de $\gA$  et 
$M$ un \Amoz.   Alors on a les \eqvcs suivantes.
\begin{enumerate}\itemsep0pt\parsep0pt
\item  $M$ est  \tf \ssi chacun des $M_{S_i}$
est un $\gA_{S_i}$-\mtfz.
\item $M$ est \pf \ssi chacun des $M_{S_i}$ est un~$\gA_{S_i}$-\mpfz.
\item  $M$ est plat \ssi chacun des $M_{S_i}$ est un~$\gA_{S_i}$-\mplz.
\item \label{iptfplcc.ptf}\relax $M$ est \ptf \ssi chacun des $M_{S_i}$ 
est un~$\gA_{S_i}$-\mptfz.
% 5
\item $M$ est projectif de rang $k$ \ssi chacun des $M_{S_i}$ 
est un~$\gA_{S_i}$-module projectif de rang~$k$.
\item $M$ est cohérent \ssi chacun des $M_{S_i}$ est un
$\gA_{S_i}$-\comoz.
\item  $M$ est \noe \ssi chacun des $M_{S_i}$ est un
$\gA_{S_i}$-module \noez.
\end{enumerate}
\end{plcc}
%--- end-plcc-----------------
\begin{proof}
\emph{1.} Voir le \plgrf{plcc.tf}.

 \emph{2.} Voir le \plgrf{plcc.pf}.

 \emph{3.} Voir le \plgrf{plcc.plat}.

 \emph{4.} Voir le \plgrf{plcc.cor.pf.ptf}. On peut aussi utiliser le fait qu'un
\mpf est \pro \ssi il est plat
(et appliquer les points \emph{2} et \emph{3}).

 \emph{5.} Résulte du point \emph{4} et du fait que le \pol rang 
  peut-être calculé localement (il est égal à $X^k$ \ssi il est 
  égal à $X^k$ après \lon en des \mocoz).

 \emph{6.} Voir le \plgrf{plcc.coh}.

 \emph{7.} Nous faisons la preuve pour la \noet définie \cot à la
Richman-Seidenberg. Limitons nous au cas de deux \lons \come en
$S_1$ et $S_2$.
Considérons, une suite croissante
$(M_k)_{k\in\NN}$ de sous-\mtfs de $M$. 
 Elle admet une sous-suite infinie
$\big(M_{\sigma(k)}\big)_{k\in\NN}$, où~$\sigma(k)<\sigma(k+1)\, \forall \, k,$
avec:
$M_{\sigma(k)}=M_{\sigma(k)+1}$ après \lon en $S_1$ pour tout~$k$.
Considérons la suite infinie  $M_{\sigma(k)}$ vue dans $M_{S_2}$. Elle admet
deux termes consécutifs égaux  $M_{\sigma(k)}$ et~$M_{\sigma(k+1)}$.
Donc~$M_{\sigma(k)}$ et  $M_{\sigma(k)+1}$ sont égaux 
à la fois dans~$M_{S_1}$ et~$M_{S_2}$. Ils sont donc égaux dans $M$.
\end{proof}
%----end{proof----------

%:--- SUBsection{subsecLGA1}------------------
\subsec{Propriétés des anneaux commutatifs}
\label{subsecLGA1}
%-----------------------------------------

Nous rappelons quelques résultats déjà établis concernant le \crc
local de quelques \prts intéressantes pour les anneaux commutatifs, au
sens de la \lon en des \mocoz.
%--- Prc lgc {plcc.propaco} ---
\begin{plcc}
\label{plcc.propaco}\relax
{\em  (Recollement concret de \prts des anneaux commutatifs)} 
Soient $S_1$, $\dots$, $S_n$  des \moco et~$\fa$ un \id de~$\gA$.
Alors on a les \eqvcs suivantes.
%---------begin item----------
\begin{enumerate}\itemsep0pt\parsep0pt
\item  $\gA$ est  \coh \ssi chaque  $\gA_{S_i}$ est  \cohz.
\item  $\gA$ est  \lsdz \ssi chaque  $\gA_{S_i}$ est  \lsdzz.
\item  $\gA$ est  \qi \ssi chaque  $\gA_{S_i}$ est  \qiz.
\item  $\gA$ est réduit \ssi chaque  $\gA_{S_i}$ est réduit.
\item  L'\id $\fa$ est \lop \ssi chaque  $\fa_{S_i}$ est \lopz.
\item  $\gA$ est  \ari \ssi chaque  $\gA_{S_i}$ est  \ariz.
\item  $\gA$ est  de Prüfer \ssi chaque  $\gA_{S_i}$ est  de Prüfer.
\item  L'\id $\fa$ est \icl \ssi chaque  $\fa_{S_i}$ est \iclz.
\item  $\gA$ est normal \ssi chaque  $\gA_{S_i}$ est normal.
\item  $\gA$ est de \ddk $\leq k$ \ssi 
chaque  $\gA_{S_i}$ est de \ddkz~$\leq k$.
\item  $\gA$ est  \noe \ssi chaque  $\gA_{S_i}$ est  \noez.
\end{enumerate}
\end{plcc}
%--- end-plcc-----------------

Rappelons \egmt que, pour des \lons en des \ecoz, le \plgc s'applique aussi
pour les notions d'\adk et d'\cori \noe \fdi (\plgref{plcc.ddk}).

\penalty-2500
%:--- subsec{Principes local-global pour les algebres
\subsec{Principes \lgbs concrets pour les \algsz}
\vspace{3pt}
\subsubsection*{Localisation en bas}

%:      --- Principe locglob plcc.apf
\begin{plcc}\label{plcc.etale}\label{plcc.apf}\relax
Soient $S_1$, $\dots$, $S_n$ des \moco d'un anneau $\gk$ et $\gA$ une \klgz. Alors
\propeq
\begin{enumerate}
\item $\gA$ est \tf (resp. plate, \fptez, \pfz, finie, entière, \stfez, \spbz, \stez) sur~$\gk$.
\item Chacune des \algs $\gA_{S_i}$ est \tf (resp. plate, \fptez, \pfz, finie, entière, \stfez, \spbz, \stez) sur~$\gk_{S_i}$.
\end{enumerate}
De même si $\gA$ est \stfe et si $\lambda\in\Asta$, 
alors $\lambda$ est dualisante \ssi chacune des formes $\lambda_{S_i}$
est dualisante. 
\end{plcc}
\begin{proof} %On pose $\gA_i=\gA_{S_i}$ et $\gk_i=\gk_{S_i}$.
\emph{1} $\Leftrightarrow$ \emph{2.} 
On introduit la \klg \fpte $\prod_i\gk_{S_i}$. Il suffit alors d'appliquer le \thref{propFidPlatPrAlg}.

  La question de la forme dualisante (lorsque $\gA$ est \stfez) est une question
d'\iso de modules et relève des \plgcs pour les modules
(en tenant compte du fait~\ref{factSpbEds}). 
%
%\sni Enfin le cas des \ases résulte du cas des \asfs et de celui des
%formes dualisantes. 
\end{proof}

%-% ENTRE NOUS
\entrenous{Le cas des \algs de Frobenius semble mystérieux.

Pour les \aGs on a déjà énoncé le principe dans la section \ref{secAGTG},
mais cela serait peut être mieux de tout regrouper ici. Une idée de \dem pourrait en outre être donnée ici. 

Il semble aussi que le cas des \algs \fnts  est problématique,
par manque de finitude.
}
%-% Fin ENTRENOUS

\paragraph{Localisation en haut}~

Il y a aussi des \plgs qui correspondent à des \prts dites
\gui{locales dans~$\gA$}. Ici nous avons besoin de \lons en des \eco (les \moco ne suffisent pas).

%:      --- Principe locglob plcc2.apf
\begin{plcc}\label{plcc2.apf}~\\
Soit $\gA$ une \klg et $s_1$, \ldots, $s_m$ des \eco de $\gA$. Alors
\propeq
\begin{enumerate}
\item $\gA$ est \tf (resp. \pfz,  plate) sur~$\gk$.
\item Chacune des \algs $\gA_{s_i}$ est \tf (resp. \pfz,  plate) sur~$\gk$.
\end{enumerate}
\end{plcc}
\begin{proof}
Tout d'abord si $\gA=\kxn=\kXn\sur\fa$ et $s=S(\ux)$ (où~\hbox{$S\in\kuX$}), alors~$\gA_s=\gk[x_1,\ldots,x_n,t]$ avec $t=1/s$ dans $\gA_s$, ce qui donne aussi

\snic{\gA_s=\gk[X_1,\ldots,X_n,T]\sur{(\fa+\gen{TS(\uX)-1})}.}

%\sni
Ainsi la \prt d'être \tf
ou \pf est stable par \lon en un \elt (mais elle ne l'est pas pour une \lon en un \mo arbitraire).
\\
Concernant la platitude, comme $\gA_s$ est plate sur $\gA$, si $\gA$ est plate sur~$\gk$,~$\gA_s$ est plate sur $\gk$ (fait~\ref{factAlgPlate}).

Supposons maintenant que $\sum_is_iu_i=1$ dans $\gA$.
\\
 Voyons tout d'abord ce que l'on  obtient si chacune des \klgs $\gA_{s_i}$ est \tfz. On peut supposer que les \gtrs proviennent d'\elts de~$\gA$ 
 (en considérant la fraction correspondante de \denoz~$1$).
Faisons une seule liste $(\xn)$ avec tous ces \elts de $\gA$.
\Llec constatera alors par un petit calcul que $\gA$ est engendrée par

\snic{(\xn,s_1,\ldots,s_m,u_1,\ldots,u_m)=(y_1,\ldots,y_p),\hbox{  avec  }p=n+2m.}

%\sni
Voyons maintenant le cas où toutes les \algs $\gA_{s_i}$ sont \pfz.
On considère des \idtrs $Y_i$ correspondant à
la liste $(y_1,\ldots,y_p)$
définie ci-dessus. On écrit $s_i=S_i(\ux)$, $u_i=U_i(\ux)$ avec des \pols à \coes dans~$\gk$.
\\
Pour le \sgr commun $(\xn)$ que nous venons de considérer, et pour chaque
$i\in\lrbm$, nous avons un
\syp correspondant, disons $F_i$, dans  $\gk[\uX,Y_{n+i},T_i]$,
qui permet de définir l'\iso
$$\preskip.4em \postskip.4em 
\gk[\uX,Y_{n+i},T_i]\sur{\fa_i}\to \gA_{s_i}, 
$$
avec $\fa_i=\gen{F_i,Y_{n+i}-S_i(\uX),Y_{n+i}T_i-1}$.  Pour chaque
$f\in F_i$ il y a un exposant $k_f$ tel que $s_i^{k_f}f(\ux)=0$ dans $\gA$.
 On peut prendre tous les $k_f$ égaux, disons à~$k$.
\\
 On considère alors le
\syp suivant dans $\gk[Y_1,\ldots,Y_p]$,\linebreak  avec~$Y_j=X_j$ pour $j\in\lrbn$.
 On prend tout d'abord tous les $Y_{n+i}^{k}f(\uX)$ pour $f\in F_i$ et $i\in\lrbm$.\\
 Ensuite
 on écrit les relations $Y_{n+i}-S_i(\uX)$
 et $Y_{n+m+i}-U_i(\uX)$ pour les indices $i\in\lrbm$.
  Enfin, on prend la relation
 qui correspond à~$\sum_iu_is_i=1$, \cad $\sum_{i=1}^mY_{n+i}Y_{n+m+i}-1$.
\\
 \Llec fera le calcul pour se convaincre que l'on  a bien ainsi une description sans faille de la \klgz~$\gA$. Le contraire e\^ut été étonnant,
voire immoral, puisque l'on  a transcrit tout ce que l'on  pouvait savoir de la situation. L'important était que cela puisse s'exprimer par un \sys fini de relations sur un \sys fini d'\idtrsz.  En fait on a procédé exactement comme dans la \dem du
\plgrf{plcc.pf} pour les \mpfsz.

%:HHH
 Concernant la platitude, considérons $(\an)$
dans $\gk$ et $(\xn)$ dans $\gA$ tels que $\sum_ix_ia_i=0$. Nous voulons montrer que $(\xn)$ est combinaison $\gA$-\lin de \rdls  dans $\gk$. Nous savons que ceci est vrai après \lon en chacun des $s_k$. On a donc un exposant $N$ tel que
pour chaque $k$ on ait une \egt 
$$s_k^N(\xn)=\som_{j=1}^{p_j}b_{k,j}(x_{1,k,j},\dots,x_{n,k,j}),$$
 ($x_{i,k,j}\in\gk$, $b_{k,j}\in\gA$)
avec $\sum_ix_{i,k,j}a_i=0$.
On termine en prenant une combinaison $\gA$-\lin des $s_k^N$ égale à $1$.    
\end{proof}

%--- Section{subsec loc glob abs}-----
\section{Quelques \plgs abstraits}
\label{subsec loc glob abs}\relax
%-------------------

Un outil essentiel en \alg commutative classique est la \lon en
 un idéal premier. Cet outil est a priori difficile
à utiliser \cot parce que l'on ne sait pas fabriquer les \ideps qui interviennent dans les \dems classiques, et dont l'existence repose
sur l'axiome du choix. Cependant, on peut remarquer que ces \ideps sont
en \gnl utilisés à l'intérieur de \dems par l'absurde, et ceci
donne une explication du fait que le recours à ces objets 
\gui{idéaux} peut être contourné et même interprété \cot (voir section~\ref{secMachLoGlo}).

%:2015 rajout du soustitre
\subsect{Principe \lgb abstrait pour les \prts de \linebreak caractère fini}{Pour les \prts \carfz}
Le \plga en \alg commutative est un principe informel selon lequel certaines
\prts concernant les modules sur les anneaux commutatifs sont vraies si
et seulement si elles sont vraies après \lon en n'importe quel \idepz.

Nous rappelons maintenant quelques cas où le \plga s'applique en \clamaz, en
expliquant le lien avec les principes concrets correspondants.
%:  --- Pr lga{plcaring}

Une version abstraite du \plgc \ref{plcc.sli} est la suivante.
\begin{plca}
\label{plca.ring}\relax
Soient  $\varphi$, $\psi$ des \alis $M\to N$, $\theta$ une \ali $N\to P$, et $x$, $y$ des
\elts de~$N$.
Alors on a les \eqvcs suivantes.
%---------begin item----------
\begin{enumerate}
\item   Recollement abstrait des \egtsz:
$$\preskip.2em \postskip.4em 
x=y\quad   {\rm dans} \quad   N  \quad  \Longleftrightarrow
\quad
\forall \fp\in\Spec\gA\;\; x/1= y/1  \quad   {\rm dans} \quad   N_{\fp} 
$$

\item   Recollement abstrait des \egts d'\alisz:

\ss\centerline {$ \varphi =\psi \quad   {\rm dans} \quad   \Lin_\gA(M,N)  \qquad
\Longleftrightarrow $}

\ss\centerline {$\quad
\forall \fp\in\Spec\gA\; \; \varphi/1= \psi/1  \quad   {\rm dans} \quad
\Lin_{\Ap}(M_{\fp},N_{\fp})$.}

\item   Recollement abstrait des \elts \ndzsz:

\ss\centerline {$x$  est \ndz dans   $N \qquad \Longleftrightarrow$}

\ss\centerline {$\forall \fp\in\Spec\gA\;
\;\;x/1\;$   est \ndz dans   $N_{\fp}$.}

\item   Recollement abstrait des solutions de \slisz:
$$\preskip.4em \postskip.4em 
 x\in\Im \varphi  \quad  \Longleftrightarrow \quad
\forall \fp\in\Spec\gA\; \; x/1\in\Im \varphi_{\fp}. 
$$

\item   Recollement abstrait des solutions de \slis sous conditions
\hmgsz. Soit  $(\Phi_\ell)$  une famille finie d'applications \hmgs
$$\preskip.2em \postskip.2em \mathrigid 2mu
\Phi_\ell:\Lin_\gA(M,N)\times N\to Q_\ell,\hbox{ ou }\Phi_\ell:\Lin_\gA(M,N)\to
Q_\ell ,\hbox{
ou  }\Phi_\ell: N\to Q_\ell. 
$$
Alors:
\[\preskip-.4em \postskip.3em 
\begin{array}{ccc} 
\left(\big(\&_\ell \; \Phi_\ell(\varphi ,y)=0\big)\;\Rightarrow \;
 y\in\Im \varphi\right) \quad \quad \Longleftrightarrow \\[.3em] 
\forall \fp\in\Spec\gA\;
\left(\big(\&_\ell \; \Phi_\ell(\varphi ,y)
=_{Q_{\ell,\fp}}
0\big)\;\Rightarrow \;
 y/1\in\Im \varphi_{\fp} \right),
 \end{array}
\]
où l'on a noté $Q_{\ell,\fp}$ pour  $(Q_\ell)_{\fp}$.

\item   Recollement abstrait des suites exactes.
La suite
$$\preskip.2em \postskip.3em 
M\vers{\varphi}N\vers{\theta}P 
$$
est exacte \ssi la suite
$$\preskip.0em \postskip.3em 
M_{\fp}\vers{\varphi_{\fp}}N_{\fp}\vers{\theta_{\fp}}P_{\fp} 
$$
est exacte pour tout $ \fp\in\Spec\gA\; $.

\item  Recollement abstrait de facteurs directs dans les \mpfsz.
Ici  $M$ est un sous-\mtf d'un \mpf $N$.
\[\preskip-.4em \postskip.4em 
\begin{array}{ccc} 
M  \hbox{ est facteur direct dans   }  N  
 \iff  \\[.1em] 
\forall \fp\in\Spec\gA\;M_{\fp}   \hbox{ est facteur direct dans   }N_{\fp}
 \end{array}
\]
\end{enumerate}

\end{plca}
%--- end-plca-----------------
%-----------------begin proof------------------
\begin{Proof} {Démonstrations (non constructives). }
Les conditions sont \ncrs en raison du fait \ref{fact.sexloc}.
Une vérification directe est d'ailleurs \imdez.
Pour les récipro\-ques, nous supposons \spdg
que l'anneau $\gA$ est non trivial. Il suffit de traiter le point \emph{4}
(voir la remarque \paref{remplgcsli}).
En fait nous avons déjà établi le point \emph{6}, qui implique
le point \emph{4}, dans le \plga \rref{plca.basic.modules}, mais nous pensons
qu'il n'est pas inutile de redonner deux \dems classiques distinctes
(la seconde est celle donnée au chapitre \ref{chapSli})
et de comparer leur degré d'effectivité.

 {\em Première \demz}.
\\
Supposons  $x\notin\Im\varphi$, cela revient à dire que $x\neq 0$
dans $N/\varphi(M)$. Puisque pour un \idep $\fp$ on a
 $\big(N/\varphi(M)\big)_\fp\simeq N_\fp/\varphi_\fp(M_\fp)$, il suffit de prouver le
point \emph{1} avec $y=0$.  On raisonne par l'absurde en supposant~$x\neq 0$ dans~$N$. Autrement dit $\Ann_\gA(x)\neq
\gen{1}$, et il existe $\fp\in\Spec \gA$ qui contient~$\Ann_\gA(x)$.
Alors, %en notant~$x/1$ l'\elt correspondant de $N_\fp$, 
puisque~$\big(\Ann_\gA(x)\big)_\fp=\Ann_{\gA_\fp}(x/1)$, on obtient~$x\neq_{N_\fp}0$.

 {\em Deuxième \demz}.
\\
La \prt $x\in\Im\varphi$ est  \carfz, on peut donc appliquer le
fait\etoz~\ref{factPropCarFin} qui dit (en \clamaz) que pour une \prt 
\carfz, le \plgc (\lon en des \mocoz) est \eqv au \plga (\lon en tous les \idemasz).
\end{Proof}
%----end{proof----------

%--- Comment{comment plga1} ----
\comms
\label{comment plga1}~
\\
 1) Il ne semble pas que la deuxième preuve, de \crc trop \gnlz,
 puisse jamais être  rendue \covz. 
 La première preuve n'est pas non plus \gui{en \gnlz}
\covz, mais il existe des cas où elle l'est. Il suffit pour
cela que les conditions suivantes soient vérifiées, dans le cas du point
\emph{4.}
%-----------------begin item------------------
\begin{enumerate}
\item [--] Le module $N$ est \pf et le module $M$ \tfz.
\item [--] L'anneau $\gA$ est \coh \fdiz.
\item [--]
Pour tout \itf strict $\fa$ de $\gA$ on sait construire un \idep $\fp$
contenant~$\fa$.
\end{enumerate}
%-----------------end item------------------
Les deux dernières conditions sont vérifiées
lorsque $\gA$ est une \alg \pf sur $\ZZ$
ou sur un corps \gui{pleinement factoriel} (voir~\cite{MRR}).

 2)
Ceci nous permet, par exemple, de donner une autre \prco du
\tho matriciel~\ref{th matproj}.
Comme nous l'avons remarqué \paref{subsec cas generique},
il nous suffit de traiter le cas \gnq et de montrer certaines
\egtsz~\hbox{$r_ir_j=0$} et $r_hu=0$. Comme l'anneau $\Gn$ est une
\alg \pf sur $\ZZ$, nous pouvons montrer ces \egts en
appliquant le \rca des \egtsz. Nous sommes donc
ramenés au cas d'un localisé local de $\Gn$, et dans ce cas les
\egts sont vraies puisque le module est libre par application
du lemme de la liberté locale.

 3) En pratique, on peut comprendre le \plga \ref{plca.ring}
sous la forme intuitive suivante: pour démontrer un \tho
d'\alg commutative dont la signification est qu'un certain
\sli sur un anneau commutatif $\gA$ admet une solution, il suffit de traiter le
cas où l'anneau est local.
C'est un principe du même genre que le principe de Lefschetz: pour
démontrer un \tho d'\alg commutative dont la
signification est qu'une certaine \ida a lieu, il
suffit de traiter le cas où l'anneau est le corps des complexes (ou
n'importe quel sous-anneau qui nous arrange, d'ailleurs).
Cette remarque est développée dans la section~\ref{secMachLoGlo}.

 4)
Dans l'article \cite{Bass}, Hyman Bass fait le commentaire suivant concernant une
version \noee du \plga \ref{plca.ring}, point~\emph{7.}\\
{\em Aussi \elr que ce résultat puisse paraître, il ne semble
pas qu'aucune preuve puisse en être donnée sans utiliser, ou reconstruire pour
l'essentiel, le foncteur $\mathrm{Ext}^1.$}\\
Ce commentaire est étonnant, au vu du \crc tout à fait anodin de notre
preuve du principe concret correspondant, laquelle ne calcule rien qui ressemble
à un $\mathrm{Ext}^1.$
En fait, lorsque le but est de montrer le scindage d'une suite suite exacte courte,
il semble que l'efficace machinerie calculatoire des
$\mathrm{Ext}$ est souvent inutile, et qu'elle peut être court-circuitée par un
argument plus \elrz.

 5) Le \plga ci-dessus fonctionne aussi en utilisant uniquement la \lon en n'importe quel \idemaz, comme vu dans le \plga \rref{plca.basic.modules}. Mais ceci n'est pas vraiment utile car les \lons en les \idemas
sont les moins poussées (parmi les \lons en les \idepsz).
Il y a par contre des cas où le raisonnement classique utilise uniquement
des \lons en des \idemisz. Ce sont des \dems plus subtiles et plus difficiles à
décrypter \cotz. Nous en parlerons dans la section~\ref{subsecLGIdepMin}.

 6) 
Comme nous l'avons remarqué \paref{plcc.tf},
le \plga pour les \mtfs ne fonctionne pas.
%:2015 
Nous revenons sur cette question dans le paragraphe qui suit.
%Il en serait de même pour le \prcc des  \mpfs
%ou pour celui des  \comosz.
%Ceci dénote une certaine supériorité des \plgcs
%sur les \plgasz.
\eoe

%:2015 \Lon au voisinage de tout \idepz
\subsec{\Lon au voisinage de tout \idepz}
Le \plg informel en \clama dit que les bonnes \prts des anneaux ou des modules sont celles qui obéissent à la règle suivante:  

\smallskip 
$\bullet$ \emph{La \prt est satisfaite \ssi elle est satisfaite après \lon en n'importe quel \idepz}. 

\smallskip 
Il y a cependant des \prts qui méritent d'être qualifiées de bonnes, comme le fait pour un module d'être \tf ou \cohz, et qui n'obéissent pas à la règle ci-dessus, mais seulement à la règle
suivante:  

\smallskip 
$\bullet$  Forme variante d'un \plgaz. \emph{La \prt est satisfaite \ssi elle est satisfaite après \lon au voisinage de n'importe quel \idepz.} 

\smallskip 
Dans cette règle, \gui{après \lon  au voisinage de l'\idepz~$\fP$}%
%:2015 rajout index
\index{localisation!au voisinage d'un \idepz} 
signifie qu'il existe un $s\notin \fP$ tel que la \prt est satisfaite pour le changement d'anneau de base $\gA\to \gA[1/s]$. 

Dans la forme variante, la vérification de la \prt locale est plus délicate, car l'anneau $\gA[1/s]$ n'est pas local, mais l'implication du local au global est plus facile à établir, et plus souvent satisfaite.

En fait cette deuxième forme de principe local-global est \gui{la meilleure}.
On démontre facilement en \clama que cette règle est \eqve à la règle
la plus stricte que nous utilisons en \comaz, celle où nous faisons appel à des \eco plutôt qu'à des \mocoz. 

\smallskip 
$\bullet$  Forme \cov stricte d'un \plgz. \emph{La \prt est satisfaite \ssi elle est satisfaite après \lon en des \ecoz.} 

\smallskip 
Enfin, sans doute le plus important, c'est cette règle qui permet de déclarer légitime le passage  des \gui{bonnes} \prts des schémas affines aux schémas de Grothendieck.
Par exemple la \prt pour un module d'être \pf et \coh vérifie la deuxième forme du \plg et c'est ce qui légitime la \dfn des \gui{faisceaux de modules \cohsz\footnote{Les faisceaux \agqs \cohs de Serre, à l'origine de l'histoire, sont \lot \gui{\cohsz} au sens de Serre et Bourbaki, \cad \pf et \cohs dans la terminologie actuelle.}}.

%\penalty-2500
%--- Section{secColleCiseaux}--------------------
\section{Recollement concret d'objets}
\label{secColleCiseaux}
%-----------------------------------------

\vspace{4pt}
%:--- SUBsection{colle et ciseaux}----------
\subsec{La colle et les ciseaux} \label{subsecColleCiseaux}
%-----------------------------------------

Nous faisons ici une brève discussion concernant les méthodes de recollement en \gmt \dile et leurs traductions en \alg commutative.

\smallskip   
Tout d'abord nous examinons la possibilité de construire une \vrt lisse à partir de cartes locales, \cad par recollement d'ouverts  $U_i$ de $\RR^n$ au moyen de difféomorphismes (ou \isosz)
$\varphi_{ij}:U_{ij}\to U_{ji}$:   $U_{ij}$ est un ouvert de $U_i$ et $\varphi_{ji}=\varphi_{ij}^{-1}$.

{

\centerline{\includegraphics[width=11cm]{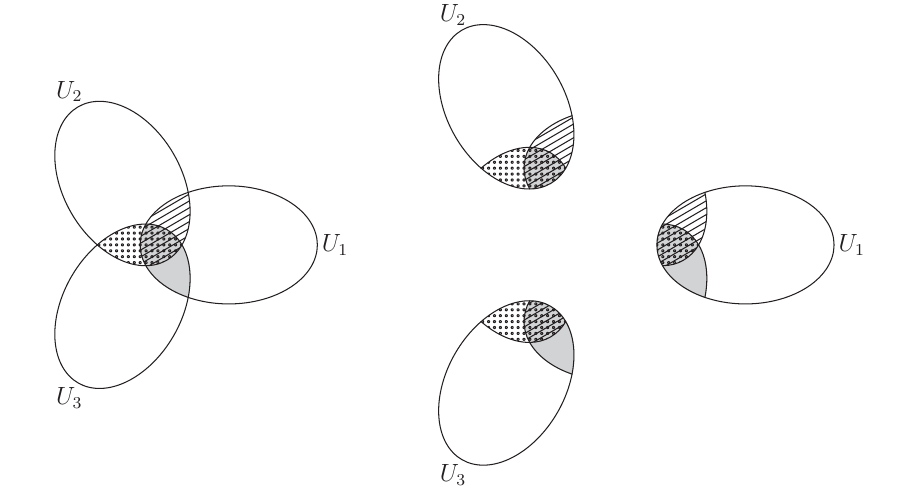}}
}

Nous allons considérer le cas simple où la \vrt est obtenue en recollant seulement un nombre fini d'ouverts de $\RR^n$.

Dans ce cas la condition à remplir est que les morphismes de \rcms doivent être \emph{compatibles entre eux trois par trois}. Cela signifie \prmt la chose suivante.
Pour chaque triplet d'indices distincts $(i,j,k)$ on considère l'ouvert
$U_{ijk}= U_{ij}\cap U_{ik}$ (avec donc $U_{ijk}=U_{ikj}$). La \cpbt signifie d'une part que, 
pour chaque $(i,j,k)$, la restriction $\Frt{\varphi_{ij}} {U_{ijk}}$ établit un \iso de $U_{ijk}$ sur $U_{jik}$, et d'autre part que  si l'on compose les \isos 

\snic{U_{ijk}\vvvvers{ \Frt{\varphi_{ij}} {U_{ijk}} }U_{jik}\quad$ et  $\quad
U_{jki}\vvvvers{ \Frt{\varphi_{jk}} {U_{jki}} }U_{kji} }

%\sni
on obtient l'\iso 
$U_{ijk}\vvvvers{ \Frt{\varphi_{ik}} {U_{ijk}} }U_{kij}$: 
$\Frt{\varphi_{ik}}{\bullet}=\Frt{\varphi_{jk}}{\bullet}\circ \Frt{\varphi_{ij}}{\bullet}\,$.

\smallskip 
Si l'on essaie de faire la même chose en \alg commutative, on va considérer des anneaux $\gA_i$ (correspondant aux anneaux $\Cin(U_i)$) et des \eltsz~\hbox{$f_{ij}\in\gA_i$}.
L'anneau $\Cin(U_{ij})$ correspondrait à $\gA_i[1/f_{ij}]$ et le morphisme de
recollement $\varphi_{ij}$ à un \iso $\omega_{ij}:\gA_i[1/f_{ij}]\to\gA_j[1/f_{ji}]$.
On devra \egmt formuler des conditions de compatibilité trois par trois.
On espère alors construire un anneau $\gA$ et des \elts $f_i\in\gA$, de telle sorte
que $\gA_i$ puisse s'identifier à $\gA[1/f_i]$, $f_{ij}$ à \gui{$f_j$ vu dans $\gA[1/f_i]$}, et~$\omega_{ij}$ à l'identité entre  $\gA[1/f_i][1/f_j]$
et  $\gA[1/f_j][1/f_i]$.

Malheureusement, cela ne fonctionne pas toujours bien. L'anneau $\gA$ censé recoller les $\gA_i$ n'existe pas toujours (cependant, s'il existe il est bien déterminé, à \iso unique près).

Le premier exemple de cet échec patent du recollement est l'espace projectif. L'espace projectif complexe $\Pn(\CC)$ est obtenu en recollant des cartes affines $\CC^n$, mais les anneaux de fonctions correspondants, isomorphes à~$\CC[\Xn]$ ne se recollent pas:
il n'y a pas de fonctions \polles définies sur $\Pn(\CC)$, à part les constantes.
Et en localisant l'anneau $\CC$ on ne risque pas d'obtenir l'anneau $\CC[\Xn]$.

Cela illustre le fait que la \gmt \agq est beaucoup plus rigide que la \gmt
$\Cin$. 

Ce phénomène désagréable est à l'origine de la création des schémas par Grothendieck, qui sont les objets abstraits obtenus \fmt en recollant des anneaux le long de morphismes de recollement lorsque les conditions de compatibilité trois par trois sont vérifiées, mais qu'aucun anneau ne veut bien réaliser le
recollement.

\smallskip  
Voyons maintenant la question du recollement des fibrés vectoriels définis
localement sur une \vrt lisse fixée $U$, recouverte par un nombre fini d'ouverts
$U_i$. On note $U_{ij}=U_i\cap U_j$. Le fibré $\pi:W\to U$ que l'on veut construire, dont toutes les fibres sont isomorphes à un \evcz~$F$ donné, est connu a priori seulement par ses restrictions~\hbox{$\pi_i:W_i\to U_i$}. 
Pour que l'on puisse recoller il faut donner des difféomorphismes de recollement~\hbox{$\psi_{ij}: W_{ij}\to W_{ji}$} où $W_{ij}=\pi_i^{-1}(U_{ij})$. 
Ces morphismes doivent tout d'abord respecter la structure d'\evc fibre par fibre. En
outre, là aussi, on a besoin de conditions de compatibilité trois par trois, analogues à celles que nous avons définies dans le premier cas.

Si maintenant on passe au cas analogue en \alg commutative, on doit partir d'un anneau
$\gA$ avec un \sys d'\eco $(f_1,\ldots,f_\ell)$. On note $\gA_i=\gA[1/f_i]$
et $\gA_{ij}=\gA[1/f_if_j]$.
Pour chaque indice $i$, on donne le \gui{module des sections du fibré $\pi_i:W_i\to U_i$}, \cad un \hbox{$\gA_i$-module $M_i$}. Les  $\psi_{ij}$ sont maintenant
représentés par des \isos de $\gA_{ij}$-modules 

\snic{\gA_{ij}\otimes_{\gA_i}M_i\vvers{\theta_{ij}} \gA_{ji}\otimes_{\gA_j}M_j
\simarrow M_{ij}=M_{ji}.}

%\sni
Nous allons voir dans les paragraphes qui suivent que cette fois-ci tout se passe bien:
si les conditions de compatibilité trois par trois sont satisfaites, 
on a bien un 
\Amo $M$ qui \gui{recolle} les $\gA_i$-modules $M_i$.

%:--- SUBsec{Un cas simple}-
\subsec{Un cas simple}\label{subsecLGFAC}

%:     Theorem{thRecolSousMod}
\begin{theorem}\label{thRecolSousMod}
Soit $\gA$ un anneau intègre de corps de fractions $\gK$, 
$N$ un \Amo sans torsion,  $S_1$, \ldots, $S_n$ des \moco de $\gA$ et
pour chaque $i\in\lrbn$ un sous-$\gA_{S_i}\!$-module $M_i$ de $S_i^{-1}N\subseteq \gK\te_\gA N$.  On suppose que pour chaque $i,j\in\lrbn$ on a $S_j^{-1}M_i=S_i^{-1}M_j$ (vus comme sous-\Amos de $\gK\te_\gA N$). On a les résultats suivants. 
\begin{enumerate}
\item Il existe un unique sous-\Amo $M$ de $N$ tel que l'on ait $S_i^{-1}M=M_i$
pour chaque $i\in\lrbn$.
\item Ce sous-module $M$ est égal à l'intersection des $M_i$.
\item Si les $M_i$ sont \tf (resp. \pfz, \cohsz, \ptfsz), il en va de même pour $M$.
\end{enumerate}
\end{theorem}

\begin{proof}
\emph{1} et \emph{2.} Soit $P=\bigcap_iM_i$. Tout d'abord $P\subseteq N$ parce qu'un
\elt de l'intersection s'écrit

\snic{\frac{x_1}{s_1}=\cdots=\frac{x_n}{s_n}=\frac{\som_ia_ix_i}{\som_ia_is_i}=\som_ia_ix_i\quad\mathrm{si}\;\;\som_ia_is_i=1\;\;\hbox{dans}\;\gA}

%\sni
(avec  $x_i\in N$, $s_i\in S_i$ pour $i\in\lrbn$).
\\Montrons que le module $P$ satisfait les conditions requises.
\\  
Tout d'abord $P \subseteq M_i$ donc $S_i^{-1} P \subseteq M_i$ pour chaque $i$. 
Inversement soit par exemple
$x_1\in M_1$, nous voulons voir qu'il est dans $S_1^{-1}P$.\\
Puisque $S_j^{-1}M_1=S_1^{-1}M_j$, il existe $u_{1,j}\in S_1$
tel que $u_{1,j}x_1\in M_j$. En posant $s_1=\prod_{j\neq1}u_{1,j}$,
on obtient bien $s_1x_1\in\bigcap_iM_i$.
\\
 Voyons maintenant l'unicité. \\
 Soit $Q$ un module satisfaisant les conditions requises. On a $Q \subseteq S_i^{-1}Q=M_i$ et ainsi
$Q\subseteq P$. Considérons alors la suite $Q \to P\to 0$.
Puisqu'elle est exacte après \lon en des \mocoz, elle est exacte
(\plgref{plcc.basic.modules}), i.e.,
l'\homo d'inclusion est surjectif: $Q=P$.

 Enfin le point \emph{3} résulte de \plgcs déjà établis. 
\end{proof}

Si l'on ne suppose pas l'anneau intègre et le module sans torsion, le \tho précédent est un peu plus délicat. Cela sera l'objet du 
\plgref{plcc.modules 2}. 

%:--- SUBsec{Recollement d'objets}-
\subsec{Recollement d'objets dans les modules}

%:--- Prc lgc {plcc.modules 1}

Soient $\gA$ un anneau commutatif et $(S_i)_{i\in\lrbn}$ des
\moco de~$\gA$.
Notons $ \gA_i:= \gA_{S_i}$ et  $ \gA_{ij}:= \gA_{S_iS_j}$ ($i\neq j$)
de sorte que $ \gA_{ij}= \gA_{ji}$. Notons~\hbox{$\alpha_i: \gA\to  \gA_i$} et
$\alpha_{ij}: \gA_i\to  \gA_{ij}$ les
\homos naturels. 

Dans la suite, des notations comme $(M_{ij})_{i<j\in\lrbn}$ et $(\varphi_{ij})_{i\neq j\in\lrbn})$ signifient que l'on a $M_{ij}=M_{ji}$ mais pas (a priori)
$\varphi_{ij}=\varphi_{ji}$.

\begin{plcc}
\label{plcc.modules 1}\relax 
{\em (Recollement concret d'\elts dans un module, et
d'\homos entre modules) }  
%-----------------begin enum------------------
\begin{enumerate}
\item \label{i1plcc.modules 1}\relax Soit un \elt $(x_i)_{i\in\lrbn}$ de  $\prod_{i\in\lrbn}  \gA_i$.
Pour qu'il existe un~\hbox{$x\in  \gA$} vérifiant $\alpha_i(x)=x_i$ dans chaque $ \gA_i$, il faut et
suffit que pour chaque~\hbox{$i<j$} on ait $\alpha_{ij}(x_i)=\alpha_{ji}(x_j)$ dans $ \gA_{ij}$. En outre, cet $x$
est alors déterminé de manière unique.
En d'autres termes l'anneau $ \gA$ (avec les \homos $\alpha_{i}$) est la limite projective du diagramme:

\snic{\big(( \gA_i)_{i\in\lrbn},( \gA_{ij})_{i<j\in\lrbn};(\alpha_{ij})_{i\neq j\in\lrbn}\big)}
$$\preskip-.4em \postskip.4em\hspace*{.1em}
\xymatrix @C=2em @R=1.5em
          {  &&& \gA_i \ar[r]^{\alpha_{ij}}\ar[ddr]^(.7){\alpha_{ik}}
                & \gA_{ij}\\
\gC\ar[urrr]^{\psi_i} \ar[drrr]_{\psi_j}\ar[ddrrr]_{\psi_k}
\ar@{-->}[rr]_(.6){\psi!} &&\gA \ar[ur]_(.4){\alpha_i}\ar[dr]^(.4){\alpha_j} \ar[ddr]_(.5){\alpha_k} &&\\
  &&& \gA_j \ar[uur]_(.7){\alpha_{ji}}\ar[dr]
          &\gA_{ik}\\
&&& \gA_k\ar[ur] \ar[r]_{\alpha_{kj}} 
   & \gA_{jk}\\
}
$$

\item \label{i2plcc.modules 1}\relax
Soit $M$ un \Amoz.
Notons $M_i:=M_{S_i}$ et  $M_{ij}:=M_{S_iS_j}$ ($i\neq j$)
de sorte que $M_{ij}=M_{ji}$. Notons $\varphi_i:M\to M_i$ et
$\varphi_{ij}:M_i\to M_{ij}$ les
\alis naturelles.
Alors le \Amo $M$ (avec les \alis $\varphi_i:M\to M_i$) est la limite projective du diagramme

\snic{\big((M_i)_{i\in\lrbn},(M_{ij})_{i<j\in\lrbn};(\varphi_{ij})_{i\neq j\in\lrbn}\big).}

%\sni Si les $M_i$ sont \tf (resp. \pfz, \cohsz, \ptfsz), 
%il en va de même pour $M$.
%
\item \label{i3plcc.modules 1}\relax
Soit  un autre module
$N$, posons $N_i:=N_{S_i}$, $N_{ij}:=N_{S_iS_j}$. 
 Soit  pour chaque $i\in\lrbn$ une application $\gA_i$-\lin
$\psi_i:M_i\rightarrow N_i$.
Pour qu'il existe une \Ali $\psi: M\rightarrow N$ vérifiant
$\psi_{S_i}=\psi_i$ pour chaque $i$, il faut et suffit que pour 
chaque $i<j$, les deux \alis $(S_j)^{-1}\psi_i$ et $(S_i)^{-1}\psi_j$  de $M_{ij}$ vers $N_{ij}$
soient égales. En outre,
l'\ali $\psi$ est alors déterminée de manière unique.
$$\preskip-.2em \postskip.4em\hspace*{3em}
\xymatrix @C=2em @R=1em
          {
    & M_i \ar[ddr]^(.3){\varphi_{ij}}\ar[rrr]^{\psi_{i}}&&&   
                  N_i \ar[ddr]^{\phi_{ij}}&\\
M\phantom{_{ij}} \ar[ur]^{\varphi_i}\ar[ddr]_{\varphi_j}\ar@{-->}[rrr]^(.7){\psi!}
&&& ~\phantom{_i} N\phantom{_j}\ar[ur]^{\phi_i}\ar[ddr]^(.3){\phi_j} \\
&& M_{ij}\ar@{-->}[rrr] &&&N_{ij}&\\
                 & M_j \ar[ur]^{\varphi_{ji}}\ar[rrr]^{\psi_{j}}&&&   
                  N_j \ar[ur]^{\phi_{ji}} 
}
$$
En d'autres termes le \Amo $\Lin_\gA(M,N)$ est la limite projective du diagramme formé par les $\Lin_{\gA_i}(M_i,N_i)$, les $\Lin_{\gA_{ij}}(M_{ij},N_{ij})$
et les \alis naturelles.

\end{enumerate}
%-----------------end enum------------------
\end{plcc}
%--- end-plcc-----------------

%---------begin proof----------
\begin{proof}
%:2012 suppression de la premiere phrase 
%Nous écrivons la \dem avec des \lons en des \ecoz, ce qui ne change rien.
%\\
 \emph{\ref{i1plcc.modules 1}.} Cas particulier de \emph{\ref{i2plcc.modules 1}.}
 
\emph{\ref{i2plcc.modules 1}.} Soit un \elt $(x_i)_{i\in\lrbn}$ de  $\prod_{i\in\lrbn}  M_i$. On doit montrer que
pour qu'il existe 
%:2012 $x\in M$ au lieu de $x\in \gA$
un $x\in M$ vérifiant $\varphi_i(x)=x_i$ dans chaque $ M_i$ il faut et
suffit que pour chaque $i<j$ on ait $\varphi_{ij}(x_i)=\varphi_{ji}(x_j)$ dans $ M_{ij}$. En outre, cet $x$
doit être unique.
\\
 La condition est clairement \ncrz. Voyons qu'elle est suffisante.
\\
Montrons l'existence de $x$. Il existe des $s_i\in S_i$ et des $y_i$ dans $M$ tels
que l'on ait $x_i=y_i/s_i$ dans chaque $M_i$.
Si $\gA$ est intègre, 
%:2012 rajout de la ligne en dessous
$M$ sans torsion et \hbox{les $s_i\neq 0$},
on a dans le module obtenu par \eds au corps
des fractions:
$$\preskip-.1em \postskip.4em\ndsp 
\;\;\;\frac{y_1}{s_1}=\frac{y_2}{s_2}=\cdots
=\frac{y_n}{s_n}=\frac{\som_ia_iy_i}{\som_ia_is_i}= \som_ia_iy_i=x\in M, 
$$
avec $\sum_ia_is_i=1$. Dans le cas \gnl on fait à peu près la même
chose.
Pour chaque couple $(i,j)$  avec $i\neq j$, le fait que $x_i/1=x_j/1$ dans $M_{ij}$
signifie que pour certains
$u_{ij}\in S_i$ et  $u_{ji}\in S_j$ on a
$s_j u_{ij} u_{ji} y_i = s_i u_{ij} u_{ji} y_j $.
Posons~\hbox{$u_i=\prod_{k\neq i}u_{ik}\in S_i$}. On a
$s_j u_{i} u_{j} y_i = s_i u_{i} u_{j} y_j $.
Soient $(a_i)$ des \elts de $\gA$ tels que $\sum_i a_i s_i u_i =1$.
Posons~\hbox{$x=\sum a_i  u_i y_i$}.
Nous devons montrer que $x/1=x_i$ dans $M_i$ pour chaque $i$. Par exemple pour
$i=1$, on écrit les \egts suivantes dans $M$:
%--------------------begin array---------------
$$\begin{array}{c}
s_1 u_1 x\,= \,s_1 u_1 \som_i a_i u_i y_i\, = \,
\som_i a_i s_1 u_1 u_i y_i\,=   \\[1mm]
\quad \som_i a_i s_i u_1 u_i y_1\,= \,
\big(\som_i a_i s_i u_i\big)u_1y_1 \,=\, u_1y_1.
\end{array}
$$
%---------------------end array--------------
Ainsi $s_1 u_1 x=u_1y_1$ dans $M$ et $x=y_1/s_1$ dans $M_{S_1}$.
\\ Enfin, l'unicité de $x$  résulte du \prcc des \egtsz.

\emph{\ref{i3plcc.modules 1}.} Les \alis composées $M\to M_i\to N_i$ sont compatibles avec les
\alis naturelles $N_i\to N_{ij}$. On conclut avec le fait que $N$ est la limite projective du diagramme des $N_i$ et $N_{ij}$ 
%:2012 point 2 et non pas 1
(point \emph{2}).
\end{proof}
%---------end proof----------

{\bf Un point délicat} (au sujet du point \emph{\ref{i3plcc.modules 1}}). Si $M$ est un \Amo \pf ou si $\gA$ est intègre et $M$ \tfz,
 les applications $\gA_i$-\lins naturelles $\Lin_\gA(M,N)_{s_i}\to \Lin_{\gA_i}(M_i,N_i)$ sont des \isos
 (voir les propositions \ref{fact.homom loc pf} et \ref{propPlateHom}).
 \\
 Dans le cas \gnlz, la notation $\psi_{s_i}$ est ambigu\"{e} car cela peut représenter, au choix, un \elt de $\Lin_{\gA_i}(M_i,N_i)$
ou un \elt de $\Lin_\gA(M,N)_{s_i}$. Et l'\ali naturelle $\Lin_\gA(M,N)_{s_i}\to \Lin_{\gA_i}(M_i,N_i)$
n'est a priori injective que si $M$ est \tfz. Cette ambiguïté peut être une source d'erreur. D'autant plus que $\Lin_\gA(M,N)$ apparaît alors comme limite projective de deux diagrammes essentiellement distincts. 
Celui basé sur les 
$\Lin_{\gA_i}(M_i,N_i)$ (le plus intéressant des deux) et celui basé sur les  $\Lin_\gA(M,N)_{s_i}$.
\eoe

{\bf Un exemple d'application de \rcm d'\eltsz.}  Vu que les \deters d'\endos de modules libres
se comportent bien par \lonz, vu le \tho qui affirme que les \mptfs sont
localement libres (au sens fort) et vu le \plgc précédent, on obtient la possibilité de \emph{définir} le \deter d'un \endo d'un \mptf
en utilisant seulement des \deters d'\endos entre
modules libres, après des \lons \come convenables.
Autrement dit, le fait suivant peut être établi indépendamment
de la théorie des \deters développée aux chapitres~\ref{chap ptf0} et~\ref{chap ptf1}.

%--- fact{prop carc loc det}--
%\begin{fact}
\sni
{\bf Fait.} \label{prop carc loc det} \relax
\emph{Pour un \endo $\varphi$ d'un \Amo \ptfz~$M$, il existe
un unique \elt $\det\varphi$ vérifiant
 la \prt suivante: si~$s\in \gA$ est tel que
le module $M_s$
soit libre, alors $(\det\varphi)_s=\det(\varphi_s)$ dans~$\gA_s$.}
\eoe
%\end{fact}
%--- end-fact----------------

%:--- SUBsec{Recollement de modules}-
\subsec{Recollement de modules}

Le principe de \rcm \ref{plcc.modules 2} qui suit précise sous quelles
conditions la limite projective d'un \sys de modules analogue rentre dans le cadre
indiqué dans le \plgrf{plcc.modules 1}.

%:--- Definition{defAliloc}-----------
\begin{definition}
\label{defAliloc}
Soient $S$ un \mo de $\gA$, $M$ un \Amo et $N$ un~\hbox{$\gA_S$-module}. Une \Ali $\alpha
:M\to N$ est appelée un \ix{morphisme de localisation en $S$} si
c'est un morphisme d'\eds de $\gA$ à $\gA_S$ pour $M$
(voir \paref{pageChgtBase}).%
\index{morphisme!de localisation en $S$ (modules)}%
%:2015 rajout index
\index{localisation!morphisme de --- (modules)}
\end{definition}
%--- end-definition------------------------------------

Autrement dit, si $\alpha:M\to N$ est un morphisme de \lon en $S$, 
et  
si~\hbox{$\beta_{M,S} :M\to M_S$} est l'\ali naturelle, l'unique 
application~\hbox{$\gA_S$-\lin} $\varphi
:M_S\to N$
telle que $\varphi \circ\beta_{M,S} =\alpha $ est un \isoz.
Un morphisme de \lon en $S$ peut être \care par les
conditions suivantes:
%-----------------begin item------------------
\begin{enumerate}
\item [--] $\forall x,x'\in M,\;\;\big(\alpha(x)=\alpha(x')\;\Longleftrightarrow\; \exists s\in
S,\;sx=sx'\big)$,
\item [--] $\forall y\in N,\;\exists x\in M,\;\exists s\in S,\;\;sy=\alpha(x)$.
\end{enumerate}
%-----------------end item------------------

%:--- Prc lgc {plcc.modules 2}
\begin{plcc}
\label{plcc.modules 2} {\em (Recollement concret de modules) }\\
 Soient $S_1$, $\dots$, $S_n$ des \moco de $\gA$. 
 \\
Notons $\gA_i=\gA_{S_i}$,
 $\gA_{ij}=\gA_{S_iS_j}$ et $\gA_{ijk}=\gA_{S_iS_jS_k}$.
On donne dans la catégorie des \Amos un diagramme
commutatif~$\fD$:

\snac{\big((M_i)_{i\in  I}),(M_{ij})_{i<j\in  I},(M_{ijk})_{i<j<k\in  I};(\varphi_{ij})_{i\neq j},(\varphi_{ijk})_{i< j,i\neq k,j\neq k}\big)}

\snii
comme dans la figure ci-après.

\smallskip {\small\hspace*{8em}{
$
\xymatrix @R=2em @C=7em{
 M _i\ar[d]_{\varphi _{ij}}\ar@/-0.75cm/[dr]^{\varphi _{ik}} &
     M _j\ar@/-1cm/[dl]^{\varphi _{ji}}\ar@/-1cm/[dr]_{\varphi _{jk}} &
        M _k\ar@/-0.75cm/[dl]_{\varphi _{ki}}\ar[d]^{\varphi _{kj}} &
\\
 M _{ij} \ar[rd]_{\varphi _{ijk}} & 
    M _{ik}  \ar[d]^{\varphi _{ikj}} & 
      M _{jk}  \ar[ld]^{\varphi _{jki}} 
\\
   &  M _{ijk} 
}
$
}}

\vspace{-.1em}
On fait les hypothèses suivantes. %
\begin{itemize}
\item Pour tous $i$, $j$, $k$ (avec $i<j<k$), $M_i$ est un $\gA_i$-module,  $M_{ij}$ est un~$\gA_{ij}$-module et $M_{ijk}$ est un~$\gA_{ijk}$-module.
Rappelons que selon nos conventions de notation on pose $M_{ji}=M_{ij}$, $M_{ijk}=M_{ikj}=\dots$

\item Pour $i\neq j$,  $\varphi_{ij}:M_i\to M_{ij}$ est un  \molo en $S_j$ (vu dans $\gA_i$).
\item Pour $i\neq k$, $j\neq k$ et $i<j$, $\varphi_{ijk}:M_{ij}\to M_{ijk}$ est un  \molo en $S_k$ (vu dans $\gA_{ij}$).
\end{itemize}
 Alors, en notant $\big(M,(\varphi_i)_{i\in\lrbn}\big)$  la limite projective du diagramme, chaque morphisme $\varphi_i:M\to M_i$ est
un  \molo en~ $S_i$.
En outre $\big(M,(\varphi_{i})_{i\in\lrbn}\big)$ est,
à \iso unique près, l'unique
\sys  qui rend le diagramme commutatif et qui fait de chaque $\varphi_i$
un  \molo en $S_i$.
\end{plcc}
%--- end-plcc-----------------

\begin{proof} 
Le premier point ne dépend pas du fait que les  $S_i$ sont \comz. En effet la construction d'une limite projective de \Amos pour un diagramme arbitraire est stable par \eds plate (parce qu'il s'agit du noyau
d'une \ali entre deux produits). 
\\
Or si l'on prend comme \eds le \moloz~\hbox{$\gA\to \gA_i$}, le diagramme se simplifie 
comme suit 

\centerline{%\small
$
%\hspace*{.7cm}
\xymatrix @R=.5cm @C=.8cm{
 M _i\ar[d]_{\varphi _{ij}}\ar[dr]^{\varphi _{ik}} &               
\\
 M _{ij} \ar[rd]_{\varphi _{ijk}} & 
    M _{ik}  \ar[d]^{\varphi _{ikj}}  
\\
   &  M _{ijk} 
}
$}

et il admet trivialement la limite projective $M_i$.

Pour démontrer l'unicité, nous raisonnons \spdg avec un \sys d'\eco $(s_1,\dots,s_n)$.
Soit $\big(N,(\psi_i)\big)$ un concurrent. \\
Puisque~$M$ est la limite projective du diagramme, il y a une unique  \Ali  $\lambda:N\to M$ telle que
$\psi_i=\varphi_i\circ \lambda$ pour tout $i$. En fait \hbox{on a $\lambda(v)=\big(\psi_1(v),\dots,\psi_n(v)\big)$}. Montrons que $\lambda$ est injective. Si $\lambda(v)=0$
tous les $\psi_i(v)$ sont nuls, et puisque $\psi_i$ est un  \molo en~$s_i$, il existe des exposants $m_i$ tels que $s_i^{m_i}v=0$. \\
Puisque les $s_i$
sont \comz, on a $v=0$. Comme $\lambda$ est injective on peut supposer $N\subseteq M$ et $\psi_i=\varphi_i\frt N$. 
Montrons que $N=M$. Soit $x\in M$. Comme~$\psi_i$ et~$\varphi_i$
sont deux  \molos en $s_i$, il y a un exposant 
$m_i$ tel que $xs_i^{m_i}\in N$. Puisque les $s_i$ sont \comz, $x\in N$.
\end{proof}

\rem Pour comprendre pourquoi la condition de comaximalité est vraiment \ncr pour l'unicité, examinons l'exemple \gui{trop simple} suivant.  
Avec l'anneau $\ZZ$, et l'unique \elt $s=2$, prenons pour $M$ \hbox{un $\ZZ[1/2]$-module} libre de base~$(a)$ (où~$a$ est un objet individuel arbitraire). \\
Pour y voir clair, nous notons $M'$ le \ZZmo $M$.\\
Considérons aussi le \ZZmo $N$ libre de base $(a)$. Considérons deux morphismes de \lon en $2^{\NN}$, $\varphi:M'\to M$ et $\psi:N\to M$. Ils  envoient tous deux $a$ sur $a$. Ainsi $M'$ et $N$ ne sont pas isomorphes comme
\ZZmos et l'unicité est en défaut.\\
Si l'on avait pris $s=1$ on pourrait définir deux morphismes de \lon 
en $1$ distincts, à savoir
$\phi_1:N\to N,\;a\mapsto a$, et $\phi_2:N\to N,\;a\mapsto -a$, et l'unicité serait assurée au sens demandé dans l'énoncé. \eoe

En pratique, on construit souvent un module en donnant des 
\hbox{$\gA_i$-modules~$M_i$} et en les recollant via leurs \lons 
$M_{ij}=M_i[1/s_j]$. Dans ce cas les modules $M_{ij}$ et $M_{ji}$  sont distincts,
et l'on doit donner pour chaque $(i,j)$ un \iso de $\gA_{ij}$-modules $\theta_{ij}:M_{ij}\to M_{ji}$.
Cela donne la variante suivante, dans laquelle on ne donne pas les modules $M_{ijk}$ en hypothèse,
mais où l'on indique les conditions de compatibilité que doivent satisfaire les $\theta_{ij}$. 
%pour permettre l'existence de modules~$M_{ijk}$ convenables. 

%: PLCC{plcc.modules 2}
\PLCC{plcc.modules 2}{ {\em (Recollement concret de modules) }
 Soient $S_1$, $\dots$, $S_n$ des \moco de $\gA$. 
 \\
Notons $\gA_i=\gA_{S_i}$,
 $\gA_{ij}=\gA_{S_iS_j}$ et $\gA_{ijk}=\gA_{S_iS_jS_k}$.
\\
On suppose donnés  des $\gA_i$-modules $M_i$ et l'on note

\snic{M_{j\ell}=M_j[1/s_\ell] \hbox{ et } M_{jk\ell}=M_j[1/s_ks_\ell] \hbox{ pour tous } j,k,\ell\hbox{ distincts }\in\lrbn,}

de sorte que $M_{jk\ell}=M_{j\ell k}$, 
avec les morphismes de \lon  

\snic{\varphi_{j\ell}: M_{j}\to M_{j\ell}\;$
et $\;\varphi_{j\ell k}: M_{j\ell}\to M_{j\ell k}.}

On suppose donnés aussi des morphismes de $\gA_{ij}$-modules 
$\theta_{ij}:M_{ij}\to M_{ji}$.
\\
On note $\theta_{ij}^{k}:M_{ijk}\to M_{jik}$ le morphisme
de $\gA_{ijk}$-modules obtenu par \lon en $s_k$ à partir de $\theta_{ij}$.
{\small
$$\preskip.2em \postskip-.4em\hspace*{-.2em}
\xymatrix @R=1em @C=.15em{
&M_i\ar[ddl]_(.4){\varphi_{ij}}\ar[ddrrr]^(.4){\varphi_{ik}} 
&&&&
M_j\ar[ddlll]_(.4){\varphi_{ji}}\ar[ddrrr]^(.4){\varphi_{jk}} 
&&&&
M_k\ar[ddlll]_(.4){\varphi_{ki}}\ar[ddr]^{\varphi_{kj}} 
\\
\\
M_{ij}\ar@<0.4ex>[rr]^{\theta_{ij}} \ar[dddrr]_{\varphi_{ijk}}    
&&
M_{ji}\ar@<0.4ex>[ll]^{\theta_{ji}} \ar[ddrrr]^(.7){\varphi_{jik}}
&\hspace{3.5em}&
M_{ik}\ar@<0.4ex>[rr]^{\theta_{ik}}  \ar[dddll]_(.6){\varphi_{ikj}}
&&
M_{ki}\ar@<0.4ex>[ll]^{\theta_{ki}} \ar[dddrr]^(.6){\varphi_{kij}}
&\hspace{3.5em}&
M_{jk}\ar@<0.4ex>[rr]^{\theta_{jk}} \ar[ddlll]_(.7){\varphi_{jki}}
&&
M_{kj}\ar@<0.4ex>[ll]^{\theta_{kj}} \ar[dddll]^{\varphi_{kji}}
\\
\\
&\hspace{2em}&     
&&&
M_{jik}\ar@{-->}[drrr]^{\theta_{jk}^{i}}  
&&&  
&\hspace{2em}
\\
&&M_{ijk}\ar@{-->}[urrr]^{\theta_{ij}^{k}}     
&&&  
&&&
M_{kij}\ar@{-->}[llllll]^{\theta_{ki}^{j}}  
&
}
$$
}

On suppose enfin que les relations de compatibilité suivantes sont satisfaites: 
\begin{itemize}
\item $\theta_{ji}\circ \theta_{ij}=\Id_{M_{ij}}$ pour $i\neq j\in\lrbn$,
\item pour $i$, $j$, $k$ distincts dans $\lrbn$, en composant circulairement 

\snic{M_{ijk} \vvvers{\theta_{ij}^{k}} M_{jik}=M_{jki}
 \vvvers{\theta_{jk}^{i}} M_{kji}=M_{kij} \vvvers{\theta_{ki}^{j}}M_{ikj}}

\smallskip 
on doit obtenir l'identité de ${M_{ijk}}$.  

%
%\item  
\end{itemize}
Alors, si $\big(M,(\varphi_i)_{i\in\lrbn}\big)$ est la limite projective du diagramme

\snic{\big((M_i)_{i\in \lrbn}),(M_{ij})_{i\neq j\in \lrbn};(\varphi_{ij})_{i\neq j},(\theta_{ij})_{i\neq j}\big),}
 
\smallskip  chaque morphisme $\varphi_i:M\to M_i$ est
un  \molo en $S_i$. \\
En outre, $\big(M,(\varphi_{i})_{i\in\lrbn}\big)$ est,
à \iso unique près, l'unique
\sys qui rend le diagramme commutatif et qui fait de chaque $\varphi_i$
un  \molo en $S_i$.

}
\begin{proof}
Notez que le diagramme ci-dessus est commutatif par construction,
sauf éventuellement le triangle du bas en trait-tirets, chaque fois qu'il est possible de joindre deux modules par deux chemins différents: par
exemple~\hbox{$\varphi_{ij}\circ \varphi_{ijk}= \varphi_{ik}\circ \varphi_{ikj}$}
et $\theta_{ij}^{k}\circ \varphi_{ijk}=\varphi_{jik}\circ \theta_{ij}$. \\  
Nous devons ici nous convaincre que les conditions de compatibilité indiquées sont exactement ce qui est \ncr et suffisant pour se ramener à la situation décrite au \plgref{plcc.modules 2}.
\\
Pour cela, lorsque $i<j<k$ nous conserverons seulement $M_{ij}$, $M_{ik}$, $M_{jk}$ \hbox{et $M_{ijk}=M_{ikj}$}. 
\\
Ceci nous force à remplacer
\[ 
\begin{array}{rclcrcl} 
 \varphi_{ji} & :  & M_j\to M_{ji} & \;\;\hbox{par}\;\;&\gamma_{ji}=\theta_{ji}\circ \varphi_{ji} &: & M_j\to M_{ij}, \\[1mm] 
 \varphi_{ki}  & :  & M_k\to M_{ki}  & \;\hbox{par}\;&\gamma_{ki}=\theta_{ki}\circ \varphi_{ki}  &: & M_k\to M_{ik},  \\[1mm] 
\varphi_{kj}   & :  & M_k\to M_{kj}  & \;\hbox{par}\;& \gamma_{kj}=\theta_{kj}\circ \varphi_{kj} &: &  M_k\to M_{jk}, \\[1mm] 
 \varphi_{jki}   & :  & M_{jk}\to M_{jik}  & \;\hbox{par}\;&\gamma_{jki}=\theta_{ji}^{k}\circ \varphi_{jki}  &: & M_{jk}\to M_{ijk}.   
 \end{array}
\] 
Jusqu'ici, tout se passe sans encombre (en relation avec les modules à deux et trois indices que l'on a décidé de conserver): les carrés  $(M_i,M_{ij},M_{ijk},M_{ik})$ et~$(M_j,M_{ij},M_{ijk},M_{jk})$ sont commutatifs et les flèches sont des \molosz.
\\
C'est seulement avec les deux \molos $M_k\to M_{ijk}$ que l'on va voir le \pbz.

{\small
$$\hspace*{0mm}
\xymatrix @R=1.2em @C=.15em{
&M_i\ar[ddl]_(.4){\varphi_{ij}}\ar[ddrrr]^(.3){\varphi_{ik}} 
&&&&
M_j\ar[ddlllll]_(.3){\gamma_{ji}}\ar[ddrrr]^(.3){\varphi_{jk}} 
&&&&
M_k\ar[ddlllll]_(.3){\gamma_{ki}}\ar[ddl]^{\gamma_{kj}} 
\\
\\
M_{ij} \ar[ddrrr]_{\varphi_{ijk}}    
&&
&\hspace{2.5em}&
M_{ik}  \ar[ddl]_{\varphi_{ikj}}
&&
&\hspace{2.5em}&
M_{jk} \ar[ddlllll]^{ \gamma_{jki}}
&&
\\
\\
&&&M_{ijk}     
&&&  
&&&
}
$$
}

Ces deux \molos sont maintenant imposés, à savoir, celui qui passe par $M_{ik}$, qui doit être 

\snic{\varphi_{ikj}\circ \gamma_{ki}=\varphi_{ikj}\circ\theta_{ki}\circ \varphi_{ki}=\theta_{ki}^{j}\circ \varphi_{kij}\circ \varphi_{ki},}

et celui qui passe par $M_{jk}$, qui doit être 

\snic{\gamma_{jki}\circ \gamma_{kj}=\theta_{ji}^{k}\circ \varphi_{jki}\circ\theta_{kj}\circ \varphi_{kj} 
=\theta_{ji}^{k}\circ\theta_{kj}^{i}\circ\varphi_{kji}\circ \varphi_{kj}.}

Comme $\varphi_{kij}\circ \varphi_{ki}=\varphi_{kji}\circ \varphi_{kj}$, la fusion
 a réussi si l'on a $\theta_{ki}^{j}=\theta_{ji}^{k}\circ\theta_{kj}^{i}$.
 \\
 En fait la condition est aussi \ncr parce que \gui{tout \molo est un épimorphisme}: si $\psi_1\circ \varphi=\psi_2\circ \varphi$ avec $\varphi$ un \moloz, 
 alors $\psi_1=\psi_2$. 
\end{proof}

\penalty-2500
%:--- subsec{subsecRecolAnn}--------------
\subsec{Recollement d'\homos entre anneaux}
\label{subsecRecolAnn}
%-----------------------------------------

%:--- Definition{defHomloc}-----------
\begin{definition}
\label{defHomloc}
Soit $S$ un \mo de $\gA$.
%-----------------begin enum------------------
Un morphisme $\alpha  :\gA\to \gB$ est appelé un \ix{morphisme de
localisation en $S$} si tout morphisme~\hbox{$\psi:\gA\to\gC$} tel que
$\psi(S)\subseteq \gC\eti$ se factorise de manière unique par $\alpha$.
\index{morphisme!de localisation en $S$ (anneaux)}%
%:2015 rajout index
\index{localisation!morphisme de --- (anneaux)}
%-----------------end enum------------------
\end{definition}
%--- end-definition------------------------------------

%--- Remarque --------
\rem %\label{factHomloc}
Si $\alpha :\gA\to\gB$ est un \homo de \lonz, et 
si~\hbox{$S=\alpha^{-1}(\gB^\times)$}, alors
$\gB$ est canoniquement isomorphe à $\gA_S$.
Par ailleurs, un morphisme de \lon peut aussi être \care comme suit:
%-----------------begin item------------------
\begin{enumerate}
\item [--] $\forall x,x'\in \gA\;(\alpha(x)=\alpha(x')\;\iff\; \exists s\in
S\;sx=sx')$
\item [--] $\forall y\in \gB,\;\exists x\in \gA,\;\exists s\in S\;\;sy=\alpha(x)$.
\eoe
\end{enumerate}
%-----------------end item---------------------

\medskip Dans la théorie des schémas développée par
Grothendieck, les morphismes de \lon $\gA\to\gA[1/s]$ jouent un rôle prépondérant.

Nous avons déjà discuté au début de cette section \ref{secColleCiseaux} l'impossibilité de recoller en \gnl des anneaux, avec l'exemple de $\Pn(\CC)$, ce qui conduit à la \dfn des schémas.

La possibilité de définir une catégorie des schémas comme
\gui{recollements d'anneaux} repose en dernière analyse sur le \prcc suivant
pour les \homos entre anneaux.
La \dem du principe est très simple. La chose importante est que le morphisme
est défini uniquement à travers des \lons et que les conditions de \cpbt sont elles-mêmes décrites via des \lons plus poussées. 

%:--- Prc lgc {plcc.RecolHomAnn}
\begin{plcc}
\label{plcc.RecolHomAnn}\relax  {\em (Recollement de morphismes d'anneaux) } 
Soient $\gA$ et $\gB$ deux anneaux, $s_1$, $\ldots$, $s_n$ des \eco de $\gA$ et
$t_1$, $\ldots$, $t_n$ des \eco de $\gB$. Posons
$$\gA_i=\gA[1/s_i],\;
\gA_{ij}=\gA[1/s_is_j],\;  \gB_i=\gB[1/t_i]\; \mathit{et} \;\gB_{ij}=\gB[1/t_it_j].$$
Pour chaque $i\in\lrbn$, soit $\varphi_i:\gA_{i}\to\gB_{i}$ 
un \homoz. Supposons 
que les conditions de \cpbt suivantes soient satisfaites: pour
$i\neq j$
les deux \homos  $\beta_{ij}\circ \varphi_i:\gA_{i}\to\gB_{ij}$   
et $\beta_{ji}\circ\varphi_j:\gA_{j}\to\gB_{ij}$ 
se factorisent via  $\gA_{ij}$ et donnent un même \homo $\varphi_{ij}:\gA_{ij}\to\gB_{ij}$ 
(voir le diagramme).
$$\preskip-.4em \postskip.4em\hspace*{3em}
\xymatrix @C=2em @R=1em 
          {
    & \gA_i \ar[ddr]^(.3){\alpha_{ij}}\ar[rrr]^{\varphi_{i}}&&&   
                  \gB_i \ar[ddr]^{\beta_{ij}}&\\
\gA\phantom{_{ij}} \ar@{..>}[ur]^{\alpha_i}\ar@{..>}[ddr]_{\alpha_j}\ar@{-->}[rrr]^(.7){\varphi!}
&&& ~\phantom{_i} \gB\phantom{_j}\ar@{..>}[ur]^{\beta_i}\ar@{..>}[ddr]^(.3){\beta_j} \\
&& \gA_{ij}\ar@{.>}[rrr]^(.3){\varphi_{ij}}
&&&\gB_{ij}&\\
                 & \gA_j \ar[ur]^{\alpha_{ji}}\ar[rrr]^{\varphi_{j}}&&&   
                  \gB_j \ar[ur]^{\beta_{ji}} 
}
$$
Alors il existe un unique \homo $\varphi:\gA\to\gB$ tel que pour chaque $i$,
on ait $\varphi_i\circ \alpha_i=\beta_i\circ \varphi$.
\end{plcc}
%--- end-plcc-----------------

%-----------------begin proof------------------
\begin{proof}
Les conditions de \cpbt sont clairement \ncrsz.
Montrons qu'elles sont suffisantes. 
D'après le \plgrf{plcc.modules 1},  $\gB$ est la limite projective du
diagramme des $\gB_i$, $\gB_{ij}$ et $\beta_{ij}$. Les  conditions de \cpbt 
impliquent que l'on a aussi les \egts 

\snic{\beta_{ij}\circ (\varphi_i\circ \alpha_i)=
\beta_{ji}\circ (\varphi_j\circ \alpha_j)}

%\sni
qui sont les conditions pour assurer l'existence et l'unicité de $\varphi$.
\end{proof}
%-----------------end proof------------------

%:-% ENTRE NOUS
\entrenous{ 
%--- Section{subsecFidPlat}----
%\subsec{
Sous-section incertaine qui pourrait être utile?

\mni {\large \bf Généralisations aux extensions fidèlement plates}
%}
%\label{subsec.gene.fidplat}\relax
\label{subsecFidPlat}
%-------------------
}
%-% Fin ENTRENOUS

%--- Section{secMachLoGlo}---- secMachLoGlo -
\section{La machinerie locale-globale \cov de base} % secMachLoGlo
\label{secMachLoGlo}\imlb
%-----------------------------------------

\begin{flushright}
{\em Localisez donc en n'importe quel \idepz.
}\\
Une mathématicienne classique
\end{flushright}

Rappelons que nous avons présenté  dans la section \ref{subsecDyna}
la philosophie \gnle de la méthode dynamique en \alg \covz.
\\
Nous indiquons maintenant comment de nombreuses preuves
utilisant le \plg en \alg abstraite peuvent être décryptées
en des \prcos conduisant aux mêmes résultats sous forme explicite.

Dans la section \ref{subsecLGIdeMax} nous nous occuperons du décryptage de
preuves abstraites qui utilisent les quotients par tous les \idemas
et dans la section \ref{subsecLGIdepMin} nous nous occuperons du décryptage de
preuves abstraites qui utilisent des \lons en tous les \idemisz.

%: subsec{Décryptage de \dems classiques qui utilisent la \lon en tout \idepz}
\subsec{Décryptage de \dems classiques qui utilisent la \lon en tout \idepz}
Un argument de \lon typique fonctionne comme suit en \clamaz. 
Lorsque l'anneau est local une certaine \prt $\sfP$ est vérifiée 
en vertu d'une \dem assez concrète. 
Lorsque l'anneau n'est pas local, la même \prt est encore vraie (d'un point de 
vue classique non \cofz) car il suffit de la vérifier \lotz.
Ceci en vertu d'un \plgaz.

Nous examinons avec un peu d'attention la première preuve. Nous voyons
alors apparaître certains calculs qui sont faisables en vertu du
principe suivant:
$$\preskip-.2em \postskip.4em 
\forall x\in\gA\;\;\; x\in\Ati\; \hbox{ou}\; x\in
\Rad(\gA), 
$$
principe qui est appliqué à des \elts $x$ provenant de la
preuve elle-même.
Autrement dit, la preuve classique donnée dans le cas local nous fournit une \prco sous l'hypothèse d'un \alo \dcdz.
Voici maintenant notre décryptage dynamique \cofz.
Dans le cas d'un anneau arbitraire, nous
répétons la même preuve, en remplaçant chaque disjonction  \gui{$x$
est \iv ou $x$ est dans le radical}, par la considération des deux
anneaux $\gA_{\cS(I,x;U)}$ et $\gA_{\cS(I;x,U)}$, où $\gA_{\cS(I,U)}$ est la \lon \gui{courante} de
l'anneau $\gA$ de départ, à l'endroit de la preuve où l'on se trouve.
Lorsque la preuve initiale est ainsi déployée, on a construit à la fin
un certain nombre, fini parce que la preuve est finie, de localisés~$\gA_{S_i}$, pour lesquels la \prt est vraie.
D'un point de vue \cofz, nous obtenons au moins le résultat
\gui{quasi global}, \cad après \lon en des \mocoz,
en vertu du lemme \ref{lemRecouvre}.
On fait alors appel à un \plgc pour conclure.

Notre décryptage de la preuve classique est rendu possible par le fait que la \prt $\sfP$ étudiée est \carfz:
elle est conservée par \lonz, et si elle est vraie après
\lon en un \moz~$S$, elle est \egmt vraie  après
\lon en un $s\in S$  (voir la section \ref{secPLGCBasic} pages \pageref{defiPropCarFini} et suivantes, ainsi que la section \ref{secSalutFini}).

Le décryptage complet contient donc deux ingrédients essentiels.
Le premier est le décryptage de la preuve donnée dans le cas local
qui permet d'obtenir un résultat quasi global (parce que la \prt est \carfz). Le deuxième
est la \prco du \plgc correspondant au \plga utilisé en \clamaz.
Dans tous les exemples que nous avons rencontré, cette \prco n'offre
aucune difficulté parce que la \dem que nous trouvons dans la littérature classique donne déjà l'argument concret,
au moins sous forme télégraphique
(sauf parfois dans Bourbaki, lorsqu'il réussit à dissimuler habilement
les arguments concrets).

La conclusion \gnle est que les \dems classiques \gui{par \plgaz} sont
déjà \covsz, si l'on veut bien se donner la peine de les lire en détail.
C'est une bonne nouvelle, outre le fait que cela confirme que les
\maths ne sont le lieu d'aucun miracle surnaturel.

\smallskip La méthode indiquée ci-dessus donne donc,
comme corolaire du lemme \ref{lemRecouvre}, le principe \gnl de
décryptage suivant, qui
\emph{permet d'obtenir automatiquement une version \cov globale (ou au moins quasi globale) d'un \tho à partir de sa version locale.}

%:   Machinerie locale-globale
\mni {\bf Machinerie locale-globale à \idepsz.}\label{MethodeIdeps}\\
{\it  Lorsque l'on relit une \prcoz, donnée pour le
cas d'un \alo \dcdz, avec un anneau $\gA$ arbitraire,
que l'on
considère au départ comme $\gA=\gA_{\cS(0;1)}$ et qu'à chaque
disjonction (pour un \elt $a$ qui se présente au cours du calcul
dans le cas local)

\sni \centerline{$ a\in\Ati\; \hbox{ou }\; a\in \Rad(\gA)$,}

\sni on remplace l'anneau \gui{en cours} $\gA_{\cS(I,U)}$ par les deux anneaux
 $\gA_{\cS(I;U,a)}$ et~$\gA_{\cS(I,a;U)}$ (dans chacun desquels le calcul
peut se poursuivre), on obtient à la fin de la relecture, une
famille finie d'anneaux
$\gA_{\cS(I_j,U_j)}$ avec les \mosz~\hbox{${\cS(I_j,U_j)}$} \com et $I_j,\;U_j$
finis. Dans chacun de ces anneaux, le calcul a été poursuivi avec succès
et a donné le résultat souhaité.
}

\medskip   On notera que si \gui{l'anneau en cours} est $\gA'=\gA_{\cS(I;U)}$
 et si la disjonction porte sur 
$$\preskip-.4em \postskip.4em 
b\in{\gA'}^{\times}\; \hbox{ ou }\; b\in \Rad(\gA'), 
$$ 
avec $b=a/(u+i)$, $a\in\gA$, $u\in\cM(U)$ et
$i\in \gen{I}_\gA$, alors il faut considérer  les localisés  $\gA_{\cS(I;U,a)}$
et  $\gA_{\cS(I,a;U)}$.

\smallskip Dans la suite nous parlerons de la  machinerie \lgbe à \ideps
comme de la \gui{machinerie \lgbe de base}.\imlb

%:   Subsection{Applications du plg}-{chapApPlg}-
\subsect{Exemples d'applications de la machinerie \lgbe de base}{Exemples d'applications}
\label{chapApPlg}\label{chapComplements}
%-------------------

Nous donnons dans ce paragraphe deux exemples significatifs.
\Llec pourra aussi consulter  l'exercice \ref{exoGCDplg}.
%:      subsubsec{Premier exemple}
\subsubsec{Premier exemple}
On veut démontrer le résultat suivant.

%:        Lemma{lemPrimUnitaire}
\begin{lemma}\label{lemPrimUnitaire}
Soit $f\in\gA[X]$ un \pol primitif et $r\in\gA$ un \elt \ndz avec $\Kdim\gA\leq1$.
Alors l'\id $\gen{f,r}$ contient un \poluz.
\end{lemma}
\begin{proof}
On commence par montrer le lemme dans le cas où $\gA$ est un \alo \dcdz.
On peut écrire $f=f_{1}+f_{2}$ avec $f_{1}\in(\Rad\gA)[X]$ et $f_{2}$ pseudo
\monz.
Par ailleurs, pour tout $e\in\Rad \gA$ on a une \egtz~\hbox{$r^m(e^m(1+ye)+zr)=0$},
donc $r$ divise $e^m$. Par suite $r$ divise une puissance de $f_{1}$, disons
d'exposant $N$.
On a ${f_{2}}^N=(f-f_{1})^N\in\gen{f,f_{1}^N}\subseteq\gen{f,r}$.
Alors,~$f_2^N$ fournit le \polu  cherché.

Pour un anneau arbitraire on reprend la preuve précédente de manière dynamique. Par exemple si $f=aX^2+bX+c$, on explicite le
raisonnement précédent sous la forme suivante.

Ou bien $a$ est \ivz, ou bien il est dans le radical. 
Si $a$ est \ivz, alors on prend $f_{2}=f,\,f_{1}=0$.

Sinon, ou  bien $b$ est \ivz, ou bien il est dans le radical.
Si $b$ est \ivz, alors on prend  $f_{2}=bX+c, \,f_{1}=aX^2$.

Sinon, ou bien $c$ est \ivz, ou bien il est dans le radical.
Si  $c$ est \ivz, alors on prend  $f_{2}=c, \,f_{1}=aX^2+bX$.

Sinon $\gen{1}=\gen{a,b,c}\in\Rad \gA$ donc l'anneau est trivial.

Voir plus loin le dessin de l'arbre des \lons successives. Les \moco se trouvent aux feuilles de l'arbre, le dernier contient $0$ et n'intervient pas dans le calcul.

%\sni
Terminons en indiquant comment on construit un \polu dans
l'\id $\gen{f,r}$ de ${\gA_{\cS(I,U)}[X]}$ à partir de deux \polus $g$ et $h$ dans les \ids $\gen{f,r}$
de ${\gA_{\cS(I,y;U)}[X]}$ et ${\gA_{\cS(I;y,U)}[X]}$.
On a d'une part 

\snic{sg=sX^m+g_{1}$ avec $\deg g_1 < m$, 
$s\in\cS(I,y;U)$ et $sg\in\gen{f,r}_{\AX},}

%\sni
et d'autre part 

\snic{th=tX^n+h_{1}$ avec $\deg h_1 < n$, 
$t\in\cS(I;y,U)$ et $th\in\gen{f,r}_{\AX}.}

\vspace{2pt}
\snac{\xymatrix {
        &\cS(0;1)\ar[ld]\ar[rd] \\
\dessus {\cS(0;a)} {a^\NN} &    &\dessus{\cS(a;1)} {1 + \langle a\rangle} \ar[ld]\ar[rd] \\
        &\dessus {\cS(a;b)} {b^\NN  + \langle a\rangle}&  &
                                    \dessus {\cS(a,b ; 1)} {1 + \langle a, b\rangle}
                                           \ar[ld]\ar[rd] \\
        &    &\dessus {\cS(a,b;c)} {c^\NN + \langle a,b\rangle}&&
                 \dessus {\cS(a,b,c;1)} {1 + \langle a,b,c\rangle} \\
}
}

%\sni
\medskip 
Les \pols $sX^ng$ et $tX^mh$ de degré formel $n+m$ ont pour \coes \fmt
dominants $s$ et $t$. En prenant $us+vt\in \cS(I,U)$, le travail est terminé
avec $usX^ng+vtX^mh$.
\end{proof}
%\eoe

%:      subsubsec  Deuxième exemple. Un résultat quasi global
\subsubsect{Deuxième exemple: un résultat quasi global
obtenu à partir d'une \dem donnée pour un \alo}{Deuxième exemple: un résultat quasi global}
\label{quasiglobaldynamique}

{\it Relecture dynamique du lemme de la liberté locale.}
La relecture dynamique de la \gui{preuve par Azumaya}
\paref{lelilo} du lemme de la liberté locale
donne une nouvelle preuve du \tho qui affirme que les \mptfs sont localement libres,
avec la formulation précise suivante.

\textit{Si $F\in\Mn(\gA)$ est une \mprnz, il existe $2^n$ \eco
$s_{i}$ tels que sur chaque $\gA_{s_{i}}$, la matrice est semblable
à une \mprn standard. Plus \prmtz, pour chaque $k=0,\ldots,n$
il y a $n \choose k$ \lons où la matrice est semblable
à $\I_{k,n}$.}

\rdb
Rappelons d'abord (voir le dessin un peu plus loin) comment se présente l'arbre du calcul pour un \alo
avec une matrice $F$ dans $\MM_3(\gA)$.

%\sni
Au point 1 le calcul démarre par le test \gui{$f_{11}$ ou $1-f_{11}$
est \ivz} (notez que la disjonction n'est en \gnl pas exclusive,
et le test doit seulement certifier que l'une des deux possibilités a bien lieu).
Si le test certifie que~$f_{11}$ est \ivz, on suit la branche de gauche,
on va en $2$ où l'on fait un changement de base qui permet de ramener la
matrice à la forme $\bloc{1}{0}{0}{G}$ avec $G\in\MM_2(\gA)$ et $G^2=G$.
Si le test certifie que $1-f_{11}$ est \ivz, on suit la branche de droite,
on va en $3$ où l'on fait un changement de base qui permet de ramener la
matrice à la forme $\bloc{0}{0}{0}{H}$ avec   
$H^2=H$.

\def \AAA {\ar@{-}[dll]\ar@{-}[drr]}
\def \BBB {\ar@{-}[dllll]\ar@{-}[drrrr]}
\snuc{\hspace*{-3mm}
\xymatrix  @C = 1pt @R = 1.2cm{
&&&&&&&&&&&&&&& 1\ar@{-}[dllllllll]\ar@{-}[drrrrrrrr]
\\
&&&&&&&\BBB2 &&&&&&&&&&&&&&&&\BBB3
\\
&&&\AAA4 &&&&&&&&\AAA5 &&&&&&&&\AAA6 &&&&&&&&\AAA7
\\
&8 &&&&9 &&&&10 &&&&11 &&&&12 &&&&13 &&&&14 &&&&15
} \label{arbrebinaire}
}

\medskip 
Si l'on est arrivé en $2$, on teste l'\elt $g$ en position $1,1$ dans $G$.
Selon le résultat, on se dirige en $4$ ou $5$ pour faire un changement
de base qui nous ramène à l'une des deux formes $\cmatrix{1&0&0\cr0&1&0\cr0&0&a}$  avec $a^2=a$, ou $\cmatrix{1&0&0\cr0&0&0\cr0&0&b}$ avec $b^2=b$.
Si l'on est arrivé en $3$, on teste l'\elt $h$ en position $1,1$ dans~$H$.
Selon le résultat, on se dirige en $6$ ou $7$ pour faire un changement
de base qui nous ramène à l'une des deux formes $\cmatrix{0&0&0\cr0&1&0\cr0&0&c}$
avec $c^2=c$, ou $\cmatrix{0&0&0\cr0&0&0\cr0&0&d}$ avec $d^2=d$.

Dans tous les cas, on termine avec un test d'inversibilité qui
certifie que l'\idm est égal à $1$ ou à $0$, ce qui donne l'une des 8
\mprns diagonales possible (avec uniquement des $0$ et $1$ sur la diagonale).

Si l'on relit ce calcul de manière dynamique avec un anneau arbitraire, on
obtient les \lons \come suivantes. \\
Au départ en 1, on a l'anneau~\hbox{$\gA_{1}=\gA$}. En $2$ et $3$ on a
les \lons \come $\gA_{2}=\gA_{1}[1/f_{11}]$
et $\gA_{3}=\gA_{1}[1/(1-f_{11})]$.
En $4$ et $5$ on a les \lons
\come de $\gA_{2}$ suivantes: $\gA_{4}=\gA_{2}[1/g]$
et $\gA_{5}=\gA_{2}[1/(1-g)]$. En $6$ et $7$ on a les \lons
\come de $\gA_{3}$ suivantes: $\gA_{6}=\gA_{3}[1/h]$
et $\gA_{7}=\gA_{3}[1/(1-h)]$. \\
On passe au dernier étage. En $8$ et $9$ on crée les
 \lons \come de $\gA_{4}$ suivantes: $\gA_{8}=\gA_{4}[1/a]$
 et $\gA_{9}=\gA_{4}[1/(1-a)]$.  En $10$ et $11$ on crée les
 \lons \come de $\gA_{5}$ suivantes: $\gA_{10}=\gA_{5}[1/b]$
  et~\hbox{$\gA_{11}=\gA_{5}[1/(1-b)]$} etc.

Au bout du compte, en considérant les \denos
 $d_{i}$ ($i=8,\ldots,15$)
des fractions créées dans les différentes branches (par exemple $d_{11}$
est le \deno dans $\gA$ de la fraction $1/f_{11}(1-g)(1-b)$,
où $g\in\gA_{2}$ et~\hbox{$b\in\gA_{5}$}), on obtient huit \elts comaximaux de $\gA$.
Et pour chacune des \lons on obtient la réduite diagonale
correspondante de la matrice de départ~$F$.

Autrement dit, la relecture dynamique de la preuve donnée dans le
cas d'un \alo crée huit \ecoz, là où la preuve classique abstraite
nous enjoindrait de localiser en tous les \idemasz, ce qui risquerait
de prendre un certain temps.

%--- section{subsecLGIdeMax}--------
\section{Quotienter par tous les \idemas}
\label{subsecLGIdeMax}
%-------------------------------------
\begin{flushright}
{\em Un anneau qui n'a pas d'\idemas est réduit à $0$.
}\\
Un mathématicien classique
\end{flushright}

On trouve dans la littérature un certain nombre de preuves dans lesquelles
l'auteur\e démontre un résultat en considérant \gui{le passage au quotient
par un \idema arbitraire}. L'analyse de ces preuves montre que le résultat
peut être compris comme le fait qu'un anneau obtenu à partir de
constructions plus ou moins compliquées est en fait réduit à $0$.
Par exemple, si l'on veut démontrer qu'un idéal $\fa$ de $\gA$ contient $1$,
on raisonne par l'absurde, on considère un \idema $\fm$ qui contiendrait $\fa$,
et l'on trouve une contradiction en faisant un calcul dans le corps
résiduel $\gA\sur{\fm}$.

Cela revient à appliquer le principe donné en exergue:
un anneau qui n'a pas d'\idemas est réduit à $0$.

Le fait de présenter le raisonnement comme une preuve par l'absurde
est le résultat d'une déformation professionnelle.
Car prouver qu'un anneau est réduit à $0$ est un fait de nature
concrète (on doit prouver que $1=0$ dans l'anneau considéré), et non pas une absurdité. Et le calcul fait dans le corps~$\gA\sur{\fm}$
ne conduit à une absurdité que parce que l'on  a décidé un jour que
dans un corps, il est interdit que $1=0$. Mais le calcul n'a rien à voir avec une telle interdiction.
Le calcul dans un corps utilise le fait que tout \elt
est nul ou \ivz, mais pas le fait que cette disjonction serait
exclusive.

En conséquence, la relecture dynamique de la preuve par l'absurde
en une \prco est possible selon la méthode suivante.
Suivons le calcul que l'on nous demande de faire comme si l'anneau $\gA\sur{\fa}$
était vraiment un corps. Chaque fois que le calcul exige de savoir si
un \elt $x_{i}$ est nul ou \iv modulo $\fa$, parions sur $x_{i}=0$
et rajoutons le à $\fa$. Au bout d'un certain temps, on trouve
que $1=0$ modulo l'idéal construit. Au lieu de perdre courage devant
une telle absurdité, voyons le bon coté des choses. Nous venons
par exemple de constater que $1\in\fa+\gen{x_{1},x_{2},x_{3}}$. Ceci est
un fait positif et non une absurdité. Nous venons en fait
de calculer un inverse~\hbox{$y_{3}$} de
$x_{3}$ dans $\gA$ modulo $\fa+\gen{x_{1},x_{2}}$.
Nous pouvons donc examiner le calcul que nous demande de faire
la preuve classique lorsque $x_{1},x_{2}\in\fm$ et $x_{3}$ est \iv
modulo $\fm$. \`A ceci près que nous n'avons pas besoin de~$\fm$
puisque nous venons d'établir que $x_{3}$ est \iv
modulo $\fa+\gen{x_{1},x_{2}}$.

Contrairement à la stratégie qui correspondait à la \lon
en n'importe quel \idepz, nous n'essayons
pas de déployer tout l'arbre du calcul
qui semble se présenter à nous.
Nous n'utilisons que des quotients,
et pour cela nous suivons systématiquement la branche \gui{être nul}
(modulo $\fm$)
plutôt que la branche \gui{être \ivz}. Ceci crée
des quotients successifs de plus en plus poussés.
Lorsqu'une soi-disant contradiction apparaît, \cad lorsqu'un
calcul a abouti à un certain résultat de nature positive,
nous revenons en arrière en profitant de l'information que nous venons
de récolter: un \elt a été certifié \iv dans le quotient
précédent.

Par exemple, avec un arbre déployé du type de
celui de la \paref{arbrebinaire}
et en prenant pour contexte \gnl l'anneau $\gA\sur{\fa}$,
si à chaque fois la branche de droite correspond à $x=0$ et celle de gauche
à $x$ \ivz,
il faut commencer par suivre le chemin $1\to3\to7\to15$ et considérer les
quotients successifs. En $15$ le calcul
nous a donné un résultat positif qui nous permet de remonter en $7$
pour suivre la branche $7\to14$. En~$14$ un résultat positif
nous permet de remonter au point $3$ (par le chemin $14\to 7\to 3$)
en sachant que l'\eltz~$a_3$ qui produit la disjonction en ce point est en fait \ivz.
Nous pouvons alors suivre le calcul proposé pour la branche $3\to6\to13$.
En~$13$ la preuve classique nous donne une soi-disant contradiction,
en fait un résultat positif dans le quotient considéré en $6$.

On aura suivi en fin de compte le chemin 

\snic{1\to 3\to 7\to 15\to 7\to 14 \to 3\to 6\to 13\to 6\to 12\to 
 1\to 2\to 5\to 11\to 5\to 10\to 2\to 4\to 9\to 4 \to 8\to 1.}

%\sni
On aura calculé uniquement dans des quotients
de $\gA\sur\fa$ et le résultat final est
que $1=0$ dans $\gA\sur\fa$,
\cad que $\fa=\gA$, qui était le but poursuivi.

Notez que lors du premier passage au point $7$, on travaille avec l'anneau
$\gA_{1,3,7}=\gA\sur{(\fa+\gen{a_1,a_3,a_7})}$. En arrivant en $15$, on apprend que
cet anneau est trivial donc que $a_7$ est \iv
dans $\gA_{1,3}=\gA\sur{(\fa+\gen{a_1,a_3})}$. En~$14$, on apprend que $\gA_{1,3}$ est trivial, i.e., que $a_3$ est \iv dans l'anneau \hbox{$\gA_{1}=\gA\sur{(\fa+\gen{a_1})}$}.
On part donc vers le point $6$ avec l'anneau
$\gA_{1}$ et un inverse de~$a_3$ en mains \ldots\, Ainsi lors des différents passages à un même point
nous ne travaillons pas avec le même anneau, car nous accumulons des informations
au fur et à mesure des calculs.

L'argument de passage au quotient
par tous les \idemas de $\gA\sur{\fa}$ (supposé par l'absurde
non réduit à $0$), qui semblait un peu magique, est ainsi
remplacé par
un calcul bien concret, donné en filigrane
par la preuve classique.
Résumons la discussion précédente.

%:   Machinerie locale-globale idemax
\goodbreak
\mni {\bf Machinerie locale-globale à \idemasz.}\label{MethodeIdemax}
%:HHH index 
\imlma\\
{\it Pour relire une preuve classique qui démontre par l'absurde
qu'un anneau $\gA$ est trivial en supposant le contraire, puis en considérant
un \idemaz~$\fm$ de cet anneau, en faisant un calcul dans le corps résiduel
et en trouvant la contradiction $1=0$, procéder comme suit. Premièrement
s'assurer que la preuve devient une \prco que $1=0$ sous l'hypothèse
supplémentaire que $\gA$ est un corps discret. Deuxièmement, supprimer
l'hypothèse supplémentaire et suivre pas à pas la preuve précédente
en privilégiant la branche $x=0$ chaque fois que la disjonction \gui{$x=0$
ou $x$ \ivz} est requise pour la suite du calcul.
Chaque fois que l'on prouve $1=0$ on a en fait montré que dans
l'anneau quotient précédemment construit, le dernier \elt
à avoir subi le test était  \ivz, ce qui permet
de remonter à ce point pour suivre la branche \gui{$x$ \ivz}
conformément à la preuve proposée pour le cas \iv
(qui est maintenant certifié).
Si la preuve considérée est suffisamment uniforme (l'expérience montre
que c'est toujours le cas), le calcul obtenu dans son ensemble est fini
et aboutit à la conclusion souhaitée.
}

\medskip \exl \\
Le lemme crucial suivant était le seul ingrédient vraiment non
\cof dans la solution par \Sus du \pb de Serre. Nous exposerons cette solution
\paref{subsecSuslin} et suivantes (voir notamment la \dem du \thrf{th4SusQS}).
Ici, nous donnons la \dem du lemme crucial par Suslin en \clamaz, puis
son décryptage \cofz.

%:     Lemma{lemSuslin1}
\begin{lemma}\label{lemSuslin1}
Soit $\gA$ un anneau, $n$ un entier $\geq2$ et $U=
\tra[\,v_1\;\cdots\;v_{n}\,] $ un vecteur \umd dans $\gA[X]^{n\times 1}$ avec
$v_{1}$ \monz. \\
Notons $V=\tra[\,v_2\;\cdots\;v_{n}\,]$.
Il existe des matrices $E_1$, \ldots, $E_\ell\in
\EE_{n-1}(\gA[X]) $, telles que, en notant $w_i$ la première \coo du vecteur
$E_i V$, l'\id $\fa$ ci-après contient 1:

\snic{\fa=\gen{\Res_X(v_1,w_1), \Res_X(v_1,w_2), \ldots,\Res_X(v_1,w_\ell) }_\gA.}
\end{lemma}
\begin{proof}
Si $n=2$, on a $u_1v_1+u_2v_2=1$ et puisque $v_1$ est \monz, $\Res(v_1,v_2)\in\Ati$:

\snac{\Res(v_1,v_2)\Res(v_1,u_2)=\Res(v_1,u_2v_2)=\Res(v_1,u_2v_2+u_1v_1)=\Res(v_1,1)=1.}

%\sni
 Si $n\geq3$, soit $d_1=\deg v_1$. On suppose \spdg que les $v_i$ sont des \pols formels de degrés $d_i<d_1$ ($i\geq2$). On a au départ des \pols $u_i$
 tels que $u_1v_1+\cdots+u_nv_n=1$.

\sni{\it \Demo classique de Suslin.} On montre que pour tout \idemaz~$\fm$,
on peut trouver une matrice $E_\fm\in \EE_{n-1}(\gA[X])$ telle que,
en notant~$w_\fm$ la première \coo
de $E_\fm V$ on ait $1\in\gen{\Res_X(v_1,w_\fm)}$ modulo~$\fm$.
Pour cela on se place sur le corps $\gk=\gA\sur{\fm}$.
En utilisant l'\algo d'Euclide, le pgcd $w_\fm$
des $v_i$ ($i\geq2$) est la première \coo d'un
vecteur obtenu par manipulations \elrs  sur $V$. On relève la matrice \elr
qui a été calculée dans $\EE_{n-1}(\gk[X])$ en une matrice
 $E_\fm\in\EE_{n-1}(\AX)$. Alors,
 puisque $v_1$ et $w_\fm$ sont premiers entre eux, le
 résultant $\Res_X(v_1,w_\fm)$ est non nul dans le corps $\gA/\fm$.

{\it \Demo \cov (par décryptage).}\\
Nous faisons une preuve par \recu sur le plus petit
des degrés formels~$d_i$, que nous notons $m$ (rappelons que $i\geq2$).
Supposons pour fixer les idées que ce soit $d_2$. \\
Initialisation: si $m=-1$,  $v_2=0$ et par une transformation \elr
on met $u_3v_3+\cdots+u_nv_n$ en position 2, ce qui nous ramène au cas $n=2$.
\\
 Récurrence: de $m-1$ à $m$.
Soit $a$ le \coe de $v_2$ de degré $m$ et
$\gB $  l'anneau  $\aqo{\gA}{a}$. Dans cet anneau l'\hdr est vérifiée.
Ainsi, on a des matrices  $E_1$, \ldots, $E_\ell\in \EE_{n-1}(\gB[X]) $, telles que, en notant $\wi{w_i}$ la première \coo
de $E_iV$, on a l'\egt

\snic{\gen{\Res_X(v_1,\wi{w_1}), \Res_X(v_1,\wi{w_2}), \ldots,\Res_X(v_1,\wi{w_\ell}) }_\gB = \gen{1}.}

%\sni
Ceci signifie, en relevant les matrices dans  $\EE_{n-1}(\AX)$
sans les changer de nom, et en notant $w_i$ la première \coo de
$E_iV$ que l'on~a:

\snic{\gen{a,\Res_X(v_1,w_1), \Res_X(v_1,w_2), \ldots,\Res_X(v_1,w_\ell)}_\gA=\gen{1}.}

%\sni
Considérons alors
$\fb=\gen{\Res_X(v_1,w_1), \Res_X(v_1,w_2), \ldots,\Res_X(v_1,w_\ell) }_\gA$,
\linebreak 
et \hbox{$\gC=\gA\sur{\fb}$}. Puisque $a$ est \iv dans $\gC$,
on peut par une manipulation \elr remplacer $v_3$ par un \pol
$v'_3=v_3-qv_2$ avec $\deg v_3'\leq m-1$.
On applique l'\hdr avec l'anneau $\gC$, on a
des matrices \elrs $E'_1$, \ldots, $E'_q\in \EE_{n-1}(\gC[X]) $ que l'on relève dans 
$\EE_{n-1}(\AX)$ sans les changer de noms. Si $w'_1$, \ldots, $w'_q$ sont les \pols correspondants (pour chaque $j$, $w'_j$ est la première \coo de $E'_jV$), 
on obtient

\snic{
1\in\gen{\Res_X(v_1,w_1), \ldots,\Res_X(v_1,w_\ell),\Res_X(v_1,w'_1), \ldots,\Res_X(v_1,w'_q)  }_\gA.
}
\end{proof}

\comm Voyons maintenant pourquoi cette élégante
 \dem est bien un décryptage de celle de Suslin selon la méthode indiquée
 auparavant. Posons $a_2=u_2v_2+\cdots+u_nv_n$.

 Lorsque l'on veut traiter sur un corps discret le vecteur $V$ par l'\algo d'Euclide,
 on doit faire des divisions. Une division dépend du degré du dividende
 (le \pol
 par lequel on divise). Dans le décryptage dynamique, on a donc des tests
 à faire sur les \coes du dividende pour déterminer son degré.
 Si l'on choisit de commencer par la division de $v_3$ par $v_2$,
 la méthode indiquée demande donc de considérer en premier
 le cas où $v_2$ est identiquement nul. Notez que cela correspond à l'initialisation de la \recuz.

 Soit $\fa_1=\gen{(v_{2,i})_{i\in\lrb{0..d_2}}}$ l'\id engendré par les
\coes de $v_2$.
\\
Si $v_2$ est identiquement nul, on a le résultant
 $r_1=\Res(v_1,a_2)=\Res(v_1,w_1)$ (\ivz) avec $w_1$
 qui est du type première cooordonnée  de
 $E_1V$ pour une matrice  $E_1\in\EE_{n-1}$
 explicite.
\\
Naturellement, ceci n'est vrai que modulo $\fa_1$,
ce qui donne $\fa_1+\gen{r_1}=\gen{1}$.
Soit $\fa_2=\gen{(v_{2,i})_{i\in\lrb{1..d_2}}}$. On a établi
que $\fa_2+\gen{r_1}+\gen{v_{2,0}}=\gen{1}$.

On raisonne maintenant modulo $\fb_2=\fa_2+\gen{r_1}$. Puisque $v_{2,0}$
est \iv et $v_2=v_{2,0}$, on peut réduire à $0$ le vecteur $v_3$
par manipulations \elrs puis mettre en position 3 un \elt égal à $a_2$
modulo $\fb_2$,
puis le ramener en position $2$.
On a donc une matrice   $E_2\in\EE_{n-1}$ avec~$w_2$ première \coo
de $E_2V$ et $\Res(v_1,w_2)=r_2$ \iv dans $\gA/\fb_2$, \cad
$\fa_2+\gen{r_1}+\gen{r_2}=\gen{1}$.
Soit $\fa_3=\gen{(v_{2,i})_{i\in\lrb{2..d_2}}}$. On vient d'établir
\linebreak 
que $\fa_3+\gen{r_1,r_2}+\gen{v_{2,1}}=\gen{1}$.

On raisonne maintenant modulo $\fb_3=\fa_3+\gen{r_1,r_2}$.
Puisque $v_{2,1}$ est \iv et $\fa_3=0$, on peut réduire à une constante
le vecteur $v_3$
par manipulations \elrs (correspondant à la division de $v_3$ par $v_2$),
puis l'emmener en position $2$.
Nous nous retrouvons dans la situation
précédemment étudiée (où $v_2$ était réduit à une constante).
On sait donc calculer deux nouvelles matrices \elrs $E_3$ et $E_4$ telles
que, en notant $w_3$ et $w_4$ leurs premières \coosz,
et $r_i=\Res(v_1,w_i)$, on obtient $\fa_3+\gen{r_1,r_2,r_3,r_4}=\gen{1}$.
\\
Soit $\fa_4=\gen{(v_{2,i})_{i\in\lrb{3..d_2}}}$. On a établi
que $\fa_4+\gen{r_1,r_2,r_3,r_4}+\gen{v_{2,2}}=\gen{1}$.

On raisonne maintenant modulo $\fb_4=\fa_4+\gen{r_1,r_2,r_3,r_4}$.
Puisque $v_{2,2}$ est \iv et $\fa_4=0$, on peut réduire au degré $1$
le vecteur $v_3$
par manipulations \elrs (correspondant à la division de $v_3$ par $v_2$),
puis l'emmener en position $2$.  Nous nous retrouvons dans la situation
précédemment étudiée (où $v_2$ était de degré $1$).
\,\ldots\,\ldots\, On obtient $\fa_4+\gen{r_1,r_2,\ldots,r_8}=\gen{1}$.
\\
Soit $\fa_5=\gen{(v_{2,i})_{i\in\lrb{4..d_2}}}$.
On a établi
que $\fa_5+\gen{r_1,r_2,\ldots,r_8}+\gen{v_{2,3}}=\gen{1}$.

Et ainsi de suite \, \ldots \,\ldots

L'important est que les inverses de \coes dominants de $v_2$ successifs
qui apparaissent dans l'\algo sont toujours calculés en tant qu'\elts de l'anneau et non pas par un procédé de \lonz.
\`A chaque fois ils ne sont \ivs que modulo un certain \id spécifié,
mais ce n'est pas grave, l'\id grandit en incorporant les résultants
autorisés mais diminue en expulsant les intrus que sont les \coes de $v_2$.
\eoe

%--- SUBsection{subsecLGIdepMin}--------
\penalty-2500
\section{Localiser en tous les \ideps minimaux}
\label{subsecLGIdepMin}
%-------------------------------------
\begin{flushright}
{\em Un anneau qui n'a pas d'\ideps minimaux est réduit à $0$.
}\\
Une mathématicienne classique
\end{flushright}

\Llec est maintenant mis\e à contribution pour se convaincre de
la justesse de la méthode suivante, en remplaçant
dans la section précédente  l'addition par la multiplication et le passage au quotient par la \lonz.

%:   Machinerie locale-globale idemin
\mni {\bf Machinerie locale-globale à \idemisz.}\label{MethodeIdemin}\\
{\it
Pour relire une preuve classique qui démontre par l'absurde
qu'un anneau~$\gA$ est trivial en supposant le contraire, puis en considérant
un \idemi de cet anneau, en faisant un calcul dans l'anneau localisé
(qui est local et \zedz, donc un corps dans le cas réduit)
et en trouvant la contradiction $1=0$, procéder comme suit.
\\
Premièrement
s'assurer que la preuve devient une \prco de l'\egt $1=0$ sous l'hypothèse
supplémentaire que $\gA$ est local et \zedz. Deuxièmement, supprimer
l'hypothèse supplémentaire et suivre pas à pas la preuve précédente
en privilégiant la branche \gui{$x$ \ivz} chaque fois que la disjonction \gui{$x$
nilpotent ou $x$ \ivz} est requise pour la suite du calcul.
Chaque fois que l'on prouve $1=0$ on a en fait montré que dans
l'anneau localisé précédemment construit, le dernier \elt
à avoir subi le test était nilpotent, ce qui permet
de remonter à ce point pour suivre la branche \gui{$x$ nilpotent}
conformément à la preuve proposée pour le cas nilpotent
(qui est maintenant certifié).
Si la preuve considérée est suffisamment uniforme (l'expérience montre
que c'est toujours le cas), le calcul obtenu dans son ensemble est fini
et aboutit à la conclusion souhaitée.

}

\medskip \exl
Un exemple assez spectaculaire est donné dans le chapitre suivant
avec le décryptage \cof d'une preuve abstraite du \tho de Traverso
concernant les anneaux seminormaux.

%%%%%%%%%%%%%%%%%%%%%%%%%%%%%%%%%%%%%%%%%%%%%%%%%%%%%%%%%%%%%%%%%%%%%%%%%%%
\penalty-2500
\section{Principes local-globals en profondeur $1$}
\label{secPlgcor}

%: subsec Un \plg relatif à la régularité
%\subsec{Un \plg relatif à la régularité}

Jusqu'à maintenant les différentes variantes du \plg étaient basées
sur les familles d'\ecoz, \cad les familles finies qui engendrent l'\id $\gen{1}$.
Une notion plus faible est suffisante pour les questions de 
régularité: il s'agit des familles finies qui engendrent un \id fidèle, ou plus \gnlt un \id \Ergz. On dit que ce sont des familles de \prof $\geq 1$.
Dans la section suivante, on examinera ce qu'on appelle les familles de \prof $\geq 2$.

%:     Definition{defiCoreg1}
\begin{definition} \label{defiCoreg1} 
\begin{enumerate}
\item Une famille finie $(\an)$ d'un anneau $\gA$ est appelé un \ix{système d'\elts coréguliers} si l'\id $\gen{\an}$ est fidèle\footnote{\`A ne pas confondre avec la notion de suite corégulière introduite par Bourbaki, comme notion duale de celle de \srgz.}.
\\
On dit aussi que \emph{l'\id $\fa$, ou la liste $(\an)$, est de profondeur $\geq 1$}, et l'on note
ceci sous la forme $\Gr_\gA(\an)\geq 1$.%
\index{coreguliers@\corz!elem@\elts ---}%
\index{profondeur!famille finie de --- $\geq 1$}
\item Soit $E$ un \Amoz.
\begin{itemize}
\item On dit qu'un \elt $a\in\gA$ est  \emph{\Ergz} (ou  \emph{\ndz pour~$E$}) si:  

\snic{\forall x\in E,\; 
(ax=0\;\Longrightarrow\; x=0).}

\item Une famille finie $(\an)$ est dire une fois \emph{\Ergez} si:

\snic{\forall x\in E,\; 
\big((a_1x=0,\,\dots,\,a_nx=0)\;\Longrightarrow\; x=0\big).}

On dit aussi que les $a_i$ sont \emph{\cor pour $E$}. 
\\
On note
ceci sous la forme $\Gr_\gA(\an,E)\geq 1$.
\item Un \itf  $\fa\subseteq \gA$ est dit \emph{\Ergz} si un (tout) \sgr de $\fa$
est une fois \Ergz.
On dit aussi que  \emph{la profondeur de~$E$
relativement à $\fa$ est supérieure ou égale à $1$}, et l'on note
ceci sous la forme $\Gr_\gA(\fa,E)\geq 1$.%
\index{E-reg@\Ergz!ideal@\id ---}\index{E-reg@\Ergz!elem@\elt ---}
\end{itemize}
\end{enumerate}
\end{definition}
%--------- fin definition ----------------------------------------------

Ainsi $\Gr_\gA(\ua)\geq 1$ signifie $\Gr_\gA(\ua,\gA)\geq 1$.
Dans la suite on donnera souvent uniquement l'énoncé
avec $\Gr_\gA(\ua,E)\geq 1$.
%:     factdefiCoreg1
\begin{fact} \label{factdefiCoreg1}~
\begin{itemize}
\item Le produit de deux \itfs \Ergs est \Ergz.
\item Si $\fa\subseteq \fa'$ avec $\fa$ \Ergz, alors $\fa'$ est \Ergz.
\end{itemize}
\end{fact}

%:     Lemma{lemCoreg1}
\begin{lemma} \label{lemCoreg1} \emph{(Astuce $(a,b,ab)$ pour la profondeur $1$)}\\
On suppose que les \ids $\gen{a,c_2, \dots, c_n}$ et $\gen{b,c_2, \dots, c_n}$ sont \Ergsz. Alors
l'\idz~$\gen{ab,c_2, \dots, c_n}$ est \Ergz. 
\end{lemma}
%--------- fin lemma ---------------------------------------------- 
%
\begin{proof}
Soit $x\in E$ tel que $abx= c_1x=\cdots=c_nx=0$.\\
Alors~$abx=c_1bx=\cdots=c_nbx=0$,
donc~$bx=0$, donc~$x=0$.
\end{proof}

On a le corolaire  \imd suivant\footnote{On aurait pu aussi remarquer que pour $q$ assez grand, l'\id $\gen{\an}^{q}$, qui \hbox{est \Ergz}, est contenu dans l'\id $\gen{a_1^p,\dots,a_n^p}$.}. 
%
%:     Lemma{lemCoreg2}
\begin{lemma} \label{lemCoreg2}
Soient $\gen{\an}$ un \id \Erg et des $p_i\in\NN$. \\
Alors l'\id $\geN{a_1^{p_1},\dots,a_n^{p_n}}$
est \Ergz. 
\end{lemma}
%--------- fin lemma ---------------------------------------------- 

On peut comparer le \plg suivant aux  points \emph{1} et \emph{3} du \plgref{plcc.sli}.

Notons que l'affirmation \gui{$\fb$ est \Ergz} est stable par \lon 
lorsque~$\fb$ est \tfz. Ceci donne l'implication dans le sens direct pour le point \emph{2\ref{i3plcc.regularite}} dans le \plg qui suit.

%: --- Prc lgc {plcc.regularite}
\begin{plcc}
%\label{plcc.reg1} 
\label{plcc.regularite}
\emph{(Localisations en \prof $\geq 1$)}\\
Soient  $b$,  $a_1$, \dots, $a_n\in\gA$,   et~$\fb$ un
\itfz. On note $\gA_i=\gA[1/a_i]$.
\begin{enumerate}
\item On suppose que les~$a_i$ sont \corz.  
\begin{enumerate}
\item  On a $x=0$ dans $\gA$   \ssi $x=0$ dans chaque $\gA_i$.
\item  L'\elt $b$ est \ndz \ssi il est  \ndz dans chaque $\gA_i$.
\item   L'\id $\fb$ est fidèle \ssi il est  fidèle dans chaque~$\gA_i$.
\end{enumerate}%

\item Soit $E$ un \Amoz, on note $E_i=E[1/a_i]$. \\
On suppose que l'\id $\gen{\an}$ est \Ergz.
\begin{enumerate}
\item \label{i1plcc.regularite} On a $x=0$ dans $E$   \ssi $x=0$ dans chaque $E_i$.
\item \label{i2plcc.regularite} L'\elt $b$ est \Erg \ssi il est  $E_i$-\ndz pour chaque~$i$.
\item  \label{i3plcc.regularite} L'\id $\fb$ est \Erg \ssi il est  $E_i$-\ndz pour chaque~$i$.
\end{enumerate}%
\end{enumerate}
\end{plcc}
\begin{proof} Il suffit de traiter le point \emph{2.}
\\
\emph{2\ref{i1plcc.regularite}.} Si $x=0$ dans  $E_i$ il y a un exposant $k_i$ tel que $a_i^{k_i}x=0$ dans~$E$. On conclut par le lemme \ref{lemCoreg2}  (avec le module~$\gA x$) que~$x=0$.

\noindent \emph{2\ref{i3plcc.regularite}.} Supposons que $\fb$ est $E_i$-\ndz pour chaque~$i$, et~$\fb\,x=0$.
Alors~$x=0$ dans chaque~$E_i$, donc~$x=0$ par le point \emph{2\ref{i1plcc.regularite}.}
\end{proof}

On utilisera souvent de manière implicite le lemme suivant, qui est une variante du lemme~\ref{factLocCas} énoncé pour les \syss d'\ecoz.

%-- Fact{factLocCasreg}----------------
\begin{fact} \emph{(Lemme des \lons \core successives)}
\label{factLocCasreg}\index{Lemme des \lons successives, profondeur 1}\\
Si $\Gr_\gA(s_1, \ldots, s_n,E)\geq 1$ et si pour chaque $i$,
on a des \elts
$s_{i,j}$, \hbox{$(j\in\lrb{1..k_i})$}
\cor pour $E\left[1/s_i\right]$,
alors les $s_is_{i,j}$ sont \cor pour~$E$.
\end{fact}
%--- end-fact-----------------------------------------
%
\begin{proof}
Soit $\fb$ l'\id engendré par les $s_is_{i,j}$. D'après le point \emph{2c}
du \plgref{plcc.regularite}, il suffit de démontrer qu'il
est \Erg après \lon en des \elts \cor pour $E$. Les $s_i$ conviennent. 
\end{proof}
%

%: subsec{\Tho de McCoy et variantes}
\subsec{Un \tho de McCoy}

Comme application du \plgref{plcc.regularite}, nous donnons une nouvelle \dem d'un \tho de McCoy (\ref{prop inj surj det} point~\emph{2}).

%:     Proposition{propInjIdd}
\CMnewtheorem{thoMccoy}{\Tho de McCoy}{\itshape}
\begin{thoMccoy} \label{propInjIdd} 
Une matrice $M\in\Ae{m\times n}$ %représente une \ali 
est injective \ssi
l'\idd $\cD_n(M)$ est fidèle.  
\end{thoMccoy}
%--------- fin proposition ---------------------------------------------- 
%
\begin{proof}
L'implication \gui{si} est simple. Montrons que si la matrice~$M$ est injective,
l'\idz~$\cD_n(M)$ est fidèle. On fait un \recu sur le nombre de colonnes.
Puisque~$M$ est injective, les \coes de la première colonne
(qui représente l'image du premier vecteur de base), engendrent un \id fidèle.
Par le \plgrf{plcc.regularite}, il suffit donc de démontrer que~$\cD_n(M)$ est fidèle sur l'anneau~$\gA_a=\gA[1/a]$, où~$a$ est un \coe  de la première colonne. 
\\
Sur cet anneau il est clair que la matrice $M$ est \eqve à une matrice de la forme 
$\blocs{.4}{.6}{.4}{.9}{$1$}{$0$}{$0$}{$N$}$. 
En outre~$N$ est injective
donc par \hdr l'\id $\cD_{n-1}(N)$ est fidèle sur~$\gA_a$. 
Enfin~$\cD_{\gA_a,n-1}(N)=\cD_{\gA_a,n}(M)$.    
\end{proof}

\rems ~\\ 
1) La \dem donne aussi que si~$m<n$ et~$M$ est injective, alors l'anneau
est trivial. En effet à chaque étape de \recuz, quand on remplace~$M$
par~$N$ la différence~$m-n$ reste constante. Donc si~$m<n$ on obtient \gui{à l'initialisation} une application injective de~$\Ae0$ dans~$\Ae{n-m}$
ce qui implique~$1=0$ dans~$\gA$. Ceci est conforme à l'énoncé \gnl
du~\thref{propInjIdd}, car pour $m<n$, $\cD_n(M)=0$, et si $0$
est un \elt \ndzz, l'anneau est trivial.
 
2) On trouve souvent dans la littérature le \tho de McCoy énoncée comme suit, sous forme contraposée (en apparence).
\\
 \emph{Si l'idéal n'est pas fidèle, l'application n'est pas injective}. 
 \\
 Ou encore de manière plus précise. 
\\
\emph{Si un \elt $x\in\gA$ non nul
annule $\cD_n(M)$, il existe un vecteur colonne non nul $C\in\Ae{m\times 1}$
tel que $MC=0$}.
\\
Malheureusement, cet énoncé ne peut être démontré qu'avec la logique classique, et l'existence du vecteur $C$ ne peut pas résulter d'un \algo \gnlz. Voici un contre-exemple, bien connu des numériciens. Si $M$ est une matrice à \coes réels avec $m<n$, on ne sait pas produire un vecteur non nul dans son noyau tant que l'on ne connaît pas le rang de la matrice. Par exemple pour $m=1$ et $n=2$, on donne deux réels $(a,b)$, et l'on cherche un couple $(c,d)\neq (0,0)$ tels que $ac+bd=0$. Si le couple~$(a,b)$ est a priori indiscernable du couple $(0,0)$, il est impossible de fournir un couple~$(c,d)$ convenable tant que l'on n'a pas élucidé si
$\abs a+\abs b$ est nul ou non.
\\
Des variantes \covs de la contraposée sont proposées dans les exercices \ref{exoMcCoyContr1} et \ref{exoMcCoyContr2}.
\eoe

%%%%%%%%%%%%%%%%%%%%%%%%%%%%%%%%%%%%%%%%%%%%%%%%%%%%%%%%%%%%%%%%%%%%%%%%%%%
\section{Principes local-globals en profondeur $2$}\label{secPlgprof2}

%:     Definition{defiProf2}
\begin{definition} \label{defiProf2} Soient $a_1$, \dots, $a_n\in\gA$ et $E$
un \Amoz.
\begin{itemize}
\item La liste  $(\ua)=(\an)$  est dite \emph{de \prof $\geq 2$} si elle est de \profz~\hbox{$\geq 1$} et si, pour toute liste $(\ux)=(\xn)$ dans $\gA$ proportionnelle\footnote{Rappelons que cela signifie que les \deters $\Dmatrix{.2em}{a_i&a_j\\ x_i&x_j}$ sont tous nuls.} à $(\ua)$, il existe $x\in\gA$ tel que $(\ux)=x (\ua)$. 
\\
On note ceci sous la forme $\Gr_\gA(\ua)\geq 2$ ou $\Gr(\ua)\geq 2$.
\item   La liste $(\ua)=(\an)$ est dite \emph{$2$ fois \Ergez} si
$\Gr_\gA(\ua,E)\geq 1$ et si, pour toute liste $(\ux)=(\xn)$ dans $E$ proportionnelle  à $(\ua)$
il existe $x\in E$ tel que $(\ux)=(\ua)x$. \\
On note ceci sous la forme $\Gr_\gA(\an,E)\geq 2$ ou $\Gr(\ua,E)\geq 2$. 
\\
On dit aussi\footnote{Eisenbud parle de la \prof de $\fa$ sur $E$, et Matsumura de la $\fa$-\prof de~$E$. La terminologie adoptée ici est celle de Bourbaki.} que \emph{la profondeur de~$E$
relativement à $(\an)$ est supérieure ou égale à $2$}.
\end{itemize}
\index{profondeur!famille finie de --- $\geq 2$}
\end{definition}
%--------- fin definition ----------------------------------------------

\rem La notation $\Gr(\fa,E)$ est prise dans \cite{Nor}. Dans ce merveilleux livre,
Northcott definit le \gui{true grade} à la Hochster comme le bon substitut
non \noe pour la profondeur.
\eoe

\medskip 
\exls 1) Dans un anneau intègre une liste $(a,b)$ avec $a,b\in\Reg(\gA)$ est de \prof $\geq 2$
\ssi $\gen{a}\cap\gen{b}=\gen{ab}$, i.e. $ab$ est le ppcm de $a$ et $b$
au sens de la \dvez.
\\
2) Dans un anneau intègre à pgcd une liste $(\an)$ est de \profz~$\geq 2$ \ssi $1$ est le pgcd de la liste.
\\
3) Si $n=1$ et la liste est réduite au seul terme $a$, $\Gr(a,E)\geq 2$ signifie que tout $y\in E$ s'écrit $y=ax$, i.e. $aE=E$. \\
En particulier
$\Gr_\gA(a)\geq 2$ signifie $a\in\Ati$.
\\
4) Toute liste d'\eco est de \prof $\geq 2$ (d'après le \plg de base). 
\eoe

\medskip Il est clair que $\Gr(\ua)\geq 2$ signifie $\Gr(\ua,\gA)\geq 2$.
Cela dispense dans la suite de dédoubler les énoncés: on les présente avec $\Gr(\ua,E)\geq 2$ pour un module~$E$ arbitraire chaque fois que c'est possible.

%:  propdef lem1prof2
\begin{propdef} \label{lem1prof2}~\\
Soient $(\ua)=(\an)$ et $(\ub)=(\br)$ dans $\gA$ et $E$
un \Amoz. Si $\Gr_\gA(\ua,E)\geq 2$ et $\gen{\ua}\subseteq \gen{\ub}$, alors
$\Gr_\gA(\ub,E)\geq 2$. \\
En conséquence, on dit qu'\emph{un \itf $\fa$ est $2$ fois \Ergz} si tout \sgr fini de $\fa$ est $2$ fois \Erg (il suffit de le vérifier pour un seul).
On note ceci sous la forme $\Gr_\gA(\fa,E)\geq 2$.
\end{propdef}
\begin{proof} Il suffit de montrer les deux faits suivants:
\begin{itemize}
\item $\Gr(\ua,E)\geq 2\Rightarrow\Gr(\ua,b,E)\geq 2$. 
\item Si $c\in\gen{\ua}$ et $\Gr(\ua,c,E)\geq 2$ alors \hbox{$\Gr(\ua,E)\geq 2$}. 
\end{itemize}
Cela montre en effet d'abord qu'on peut remplacer un \sgr d'un \itf par un autre  sans changer \gui{la \prof $\geq 2$} et ensuite que lorsqu'on
remplace $\fa$ par un \itf plus grand, la \prof $\geq 2$ se conserve.
\\
Voyons le premier point. On a une liste $(\xn,y)$ dans $E$ proportionnelle 
à $(\an,b)$. On trouve un $x$ (d'ailleurs unique) tel que $(\ux)=(\ua)x$.
On doit montrer que $bx=y$. Or $a_iy=bx_i$ et $bx_i=ba_ix$ pour $i\in\lrbn$. 
\\Donc $a_i(y-bx)=0$ et l'on conclut que $y=bx$ parce que $\Gr(\ua,E)\geq 1$.
\\
Le deuxième point est laissé \alecz.
\end{proof}
%

%%:     lemma  lem2prof2
%\begin{lemma} \label{lem2prof2}
%Soient $(\ua)=(\an)$ et $(\ub)=(\br)$ dans $\gA$ et $E$
%un \Amoz. On note $(\ua \star \ub)$ la famille finie des $a_ib_j$.
%\\
%Si $\Gr_\gA(\ua,E)\geq 2$ et $\Gr_\gA(\ub,E)\geq 2$ alors $\Gr_\gA(\ua\star\ub,E)\geq 2$.
%\\
%En termes d'\itfsz: \\
%--- si $\Gr_\gA(\fa,E)\geq 2$ et $\Gr_\gA(\fb,E)\geq 2$ alors $\Gr_\gA(\fa \fb,E)\geq 2$. 
%\end{lemma}
%%
%\begin{proof}
%On considère une famille $(x_{ij})$ dans $E$ proportionnelle à la famille $(a_ib_j)$.\\
%Pour chaque $i\in \lrbn$ et $j$, $k\in\lrbr$, on a $a_i \Dmatrix{.2em}{x_{ij}&x_{ik}\\ b_{j}&b_k}=0$. Puisque la liste~$(\ua)$ est une fois \Ergez, on a $\Dmatrix{.2em}{x_{ij}&x_{ik}\\ b_{j}&b_k}=0$ pour tous $j$, $k\in\lrbr$. \\
%On a donc un $y_i\in E$
%qui vérifie $x_{ij}=y_ib_j$ pour $j\in\lrbr$. 
%\\
%Ainsi pour $i,\ell \in\lrbn$ et $j\in\lrbr$ on a $\Dmatrix{.2em}{b_jy_{i}&b_jy_{\ell}\\ b_{j}a_i& b_ja_\ell}=b_j^{2}\Dmatrix{.2em}{y_{i}&y_{\ell}\\ a_i&a_\ell}=0$. Comme la famille des $b_j^{2}$ est une fois \Ergez, on obtient $\Dmatrix{.2em}{y_{i}&y_{\ell}\\ a_i&a_\ell}=0$. 
%On a donc un $y\in E$ pour lequel $y_i=a_iy$ pour tout $i\in\lrbn$.
%\\
%Finalement pour cet $y$, on a $x_{ij}=ya_{i}b_j$ pour tous $i,j$. 
%\end{proof}
%%

%:     Lemma{lemtrickprof2}
\begin{lemma} \label{lemtrickprof2} \emph{(Astuce $(a,b,ab)$ pour la profondeur $2$)}\\
On suppose que les listes $({\an,a})$ et $({\an,b})$ sont deux fois \Ergesz. Alors
la liste~$({\an,ab})$ est deux fois \Ergez.
\\
En conséquence, si  $\Gr_\gA(\an,E)\geq 2$, alors $\Gr_\gA(a_1^m,\dots,a_n^m,E)\geq 2$ pour $m>0$. 
\end{lemma}
%--------- fin lemma ----------------------------------------------
%
\begin{proof} On sait déjà que $({\an,ab})$ est une fois \Ergez.
\\
Soit $(\xn,y)$ une liste dans $E$ proportionnelle à $({\an,ab})$.
La suite $(x_1b,\dots,x_nb,y)$ est proportionnelle à $({\an,a})$.
Il existe donc un~$z\in E$ tel que 
$$
x_1b=a_1z,\,\dots,\,x_nb=a_nz,\,y=az
$$
Cela implique que la liste $(\xn,z)$ est proportionnelle à $({\an,b})$. Il existe donc un  $x\in E$ tel que 
$$
x_1=a_1x,\,\dots,\,x_n=a_nx,\,z=bx\;\hbox{ et a fortiori }y=abx
$$

\vspace{-2em}
\end{proof}
% 
%:     PrincipeLocGlob{plcc1Profondeur2}
\begin{plcc} \label{plcc0Profondeur2} \label{plcc1Profondeur2} \emph{(Pour la \dve et les anneaux \iclz, \lons en profondeur 2)}\\
On considère une famille $(\us)=(s_1,\dots,s_n)$ dans $\gA$  avec $\Gr_\gA(\us,E)\geq 2$. On note $\gA_{i}=\gA[\fraC1{s_i}]$ et $E_i=E[\fraC1{s_i}]$.
\begin{enumerate}
\item Soit $a\in\gA$ un \elt \Erg et $y\in E$. Alors $a$ \gui{divise} $y$ 
dans~$E$ \ssi  $a$ divise $y$ après \lon en chaque $s_i$.
\item Soit $(\bbm)$ dans $\gA$. Alors $\Gr_{\gA}(\bbm,E)\geq 2$ \ssi $\Gr_{\gA_i}(\bbm,E_i)\geq 2$ pour chaque $i$.
\item \label{I2plccAnnDivl} Supposons $\gA$ intègre et $\Gr_\gA(\us)\geq 2$, alors $\gA$ est \icl \ssi chaque anneau~$\gA_i$ est \iclz.
\end{enumerate}
\end{plcc}
\begin{proof} \emph{1.} Supposons que $a$ divise $y$ après \lon en $s_i$. On a $ax_i=u_iy$ dans~$E$ pour un $u_i=s_i^{n_i}$ et un $x_i\in E$. La liste des $u_i$ est $2$ fois \Erge (lemme \ref{lemtrickprof2}). 
\\
On a $au_jx_i=u_iu_jy=au_ix_j$ et comme $a$ est \Ergz, $u_jx_i=u_ix_j$. Donc on a un $x\in E$ tel que $x_i=u_ix$ pour chaque $i$.
Ceci donne $u_iax=u_iy$ et comme $\Gr(u_1,\dots,u_n,E)\geq 1$, on obtient $ax=y$.

\emph{2.} 
Considérons dans $\gA$ une suite $(\xm)$ proportionnelle à $(\bbm)$.
On cherche un $x\in E$ tel que $x_\ell=xc_\ell$ pour tout $\ell \in\lrbm$.
Dans chaque $E_i$ on trouve un  $y_i$ tel que $x_\ell=y_i c_\ell$ pour tout $\ell \in\lrbm$. 
Cela signifie qu'on a un $u_i\in s_i^{\NN}$ et un $z_i\in E$ tels que $u_ix_\ell=z_i c_\ell$ dans $E$ pour tout $\ell \in\lrbm$. Il nous suffit
de montrer qu'il existe un $z\in E$ tel que  $z_i=u_i z$ pour chaque $i$,
car alors $u_i(x_\ell-z c_\ell)=0$ pour chaque~$i$ (et les $u_i$ sont
\cor pour~$E$). Il suffit donc de montrer que les $z_i$ forment une famille proportionnelle aux~$u_i$, i.e. que $u_iz_j=u_jz_i$ pour tous $i,j\in\lrbn$.
Or on sait que les $c_\ell$ sont \cor pour $E$ (d'après le \plgref{plcc.regularite}). Donc il suffit de montrer que l'on a les \egts
$u_iz_jc_\ell=u_jz_ic_\ell$, or les deux membres sont égaux à $u_iu_jx_\ell$.

\emph{3.} Soient $x$ et $y$ dans $\gA$ avec $y$ entier sur l'\id $x\gA$.
Ceci reste vrai pour chaque localisé $\gA_i$, lequel est \iclz. Donc $x$ divise $y$ dans  chaque  $\gA_i$. Donc par le point \emph{1} avec $E=\gA$, $x$ divise $y$ dans $\gA$. 
\end{proof}

%
%-- Fact{lelosuccprof}----------------
\begin{fact} \emph{(Lemme des \lons successives, avec la profondeur 2)}
\label{lelosuccprof}\index{Lemme des \lons successives, profondeur 2}
%:2015 rajout d'un index
\\
Si $\Gr_\gA(s_1, \ldots,s_n,E)\geq 2$
 et si pour chaque $i$
on a une liste
$(s_{i,1},\ldots ,
s_{i,k_i})
$ dans $\gA$
qui est  $2$ fois  $E[1/s_i]$-\ndzez,
alors le \sys des $s_{i}s_{i,j}$ est  $2$ fois \Ergz.
\end{fact}
\begin{proof}
D'après le \plgrf{plcc1Profondeur2}, il suffit de vérifier que les~$s_is_{ij}$ sont~$2$ fois \Ergs
après \lon en des \elts qui forment une liste $2$ fois \Ergez. Cela fonctionne avec la liste des $s_i$. 
\end{proof}
%

%:     lemma  lem2prof2
\begin{lemma} \label{lem2prof2}
Soient $(\ua)=(\an)$ et $(\ub)=(\br)$ dans $\gA$ et $E$
un \Amoz. On note $(\ua \star \ub)$ la famille finie des $a_ib_j$.
\\
Si $\Gr_\gA(\ua,E)\geq 2$ et $\Gr_\gA(\ub,E)\geq 2$ alors $\Gr_\gA(\ua\star\ub,E)\geq 2$.
\\
En termes d'\itfsz: \\
--- si $\Gr_\gA(\fa,E)\geq 2$ et $\Gr_\gA(\fb,E)\geq 2$ alors $\Gr_\gA(\fa \fb,E)\geq 2$. 
\end{lemma}
\begin{proof}
D'après le \plgref{plcc0Profondeur2}, il suffit de montrer que la famille des~$a_ib_j$ est 2 fois \Erge après \lon en chacun des $a_i$.
Or, lorsqu'on localise en $a_1$ par exemple, la suite des $a_1b_j$
engendre le même \id que la suite des $b_j$, et cet \id est 2 fois \Ergz. 
\end{proof}

%:subsec{Recollements en profondeur 2}
\subsec{Recollements en profondeur 2}

La \dfn suivante permet de simplifier un peu la rédaction de certaines \demsz.
%:     Definition{defiProfMon2}
\begin{definition} \label{defiProfMon2}  
\emph{(Sustème de \mos $2$ fois~\hbox{\Ergz})}\\
 Un \sysz~$(S_1, \dots, S_n)=(\uS)$ de \mos de $\gA$ est dit  \emph{$2$ fois \Ergz} si pour tous~$s_1\in S_1,\,\dots,\,s_n\in S_n$, on a~$\Gr_\gA(s_1,\dots,s_n,E)\geq 2$.
\end{definition}
%--------- fin definition ------------------------------------------

Le cas le plus important est le \sys des \mos $(s_1^{\NN},\dots,s_n^{\NN})$ 
lorsque $\Gr_\gA(s_1,\dots,s_n,E)\geq 2$.

Nous reprenons maintenant le \plgref{plcc.modules 1} en remplaçant l'hypothèse selon laquelle les \mos sont \com par une hypothèse plus faible 
(\sys de \mos deux fois \ndzz).

Le contexte est le suivant.
On considère  $(\uS)=(S_i)_{i\in\lrbn}$ un \sys de \mosz.
\\
Nous notons $ \gA_i:= \gA_{S_i}$ et  $ \gA_{ij}:= \gA_{S_iS_j}$ ($i\neq j$)
de sorte que $ \gA_{ij}= \gA_{ji}$. 
\\Nous notons~\hbox{$\varphi_i: \gA\to  \gA_i$} et
$\varphi_{ij}: \gA_i\to  \gA_{ij}$ les
\homos naturels. 
\\
Dans la suite des notations comme $(E_{ij})_{i<j\in\lrbn}$ et $(\varphi_{ij})_{i\neq j\in\lrbn})$ signifient que l'on a $E_{ij}=E_{ji}$ mais pas (a priori)
$\varphi_{ij}=\varphi_{ji}$.

%: plcc.modules1bis 
\begin{plcc}
\label{plcc.modules1bis}  {\em (Recouvrement un module par des localisés
en \prof 2) }  
%-----------------begin enum------------------
On considère le contexte décrit ci-dessus.
\begin{enumerate}
\item \label{i1plcc.modules1bis} 
On suppose $(\uS)$  deux fois \ndzz. On considère un \elt $(x_i)_{i\in\lrbn}$ de  $\prod_{i\in\lrbn}  \gA_i$.
Pour qu'il existe un~\hbox{$x\in  \gA$} vérifiant $\varphi_i(x)=x_i$ dans chaque $ \gA_i$, il faut et
suffit que pour chaque~\hbox{$i<j$} on ait $\varphi_{ij}(x_i)=\varphi_{ji}(x_j)$ dans $ \gA_{ij}$. En outre, cet $x$
est alors déterminé de manière unique.
\\
En d'autres termes l'anneau $ \gA$ (avec les \homos $\varphi_{i}$) est la limite projective du diagramme:

\snic{\big(( \gA_i)_{i\in\lrbn},( \gA_{ij})_{i<j\in\lrbn};(\varphi_{ij})_{i\neq j\in\lrbn}\big)}

\item \label{i2plcc.modules1bis} 
Soit $E$ un \Amoz.  On suppose $(\uS)$ deux fois \Ergz.
\\Notons $E_i:=E_{S_i}$ et  $E_{ij}:=E_{S_iS_j}$ ($i\neq j$)
de sorte que $E_{ij}=E_{ji}$. 
\\Notons $\varphi_i:E\to E_i$ et
$\varphi_{ij}:E_i\to E_{ij}$ les
\alis naturelles.
Alors le couple $\big(E,(\varphi_{i})_{i\in\lrbn}\big)$ donne la limite projective du diagramme
suivant dans la catégorie des \Amosz:

\snic{\big((E_i)_{i\in\lrbn},(E_{ij})_{i<j\in\lrbn};(\varphi_{ij})_{i\neq j\in\lrbn}\big)}

\vspace{-.5em}
$$
\xymatrix @C=3.5em @R=1.5em %%@R=1.3cm @C=1.8cm
          {
                &&& E_i \ar[r]^{\varphi_{ij}}\ar[ddr]_(.3){\varphi_{ik}}
                & E_{ij}\\
F\ar[urrr]^{\psi_i} \ar[drrr]^{\psi_j}\ar[ddrrr]_{\psi_k}
\ar@{-->}[rr]^(.6){\psi!} &&E \ar[ur]_{\varphi_i}\ar[dr]^{\varphi_j} \ar[ddr]_(.5){\varphi_k} &&\\
  &&& E_j \ar[uur]_(.7){\varphi_{ji}}\ar[dr]
%^{\varphi_{}} 
          &E_{ik}\\
&&& E_k \ar[ur]\ar[r]_{\varphi_{kj}} 
   & E_{jk}\\
}
$$
\end{enumerate}
%-----------------end enum------------------
\end{plcc}
%--- end-plcc-----------------

%---------begin proof----------
\begin{proof}
 \emph{\ref{i1plcc.modules1bis}.} Cas particulier de \emph{\ref{i2plcc.modules1bis}.}

\emph{\ref{i2plcc.modules1bis}.} Soit un \elt $(x_i)_{i\in\lrbn}$ de  $\prod_{i\in\lrbn}  E_i$. On doit montrer l'\eqvc suivante: 
il existe un $x\in  E$ 
vérifiant $\varphi_i(x)=x_i$ dans chaque $ E_i$ \ssi pour chaque $i<j$ on a $\varphi_{ij}(x_i)=\varphi_{ji}(x_j)$ dans $E_{ij}$. En outre, cet $x$ doit être unique.
\\
 La condition est clairement \ncrz. Voyons qu'elle est suffisante.
\\
Montrons l'existence de $x$. Notons tout d'abord qu'il existe des $s_i\in S_i$ et des $y_i$ dans $E$ tels
que l'on ait $x_i=y_i/s_i$ dans chaque $E_i$.
\\
Si $\gA$ est intègre, 
$E$ sans torsion et \hbox{les $s_i\neq 0$},
on a dans l'\evc obtenu par \eds au corps
des fractions les \egtsz

\snic{\frac{y_1}{s_1}=\frac{y_2}{s_2}=\cdots
=\frac{y_n}{s_n},}

et vue l'hypothèse concernant les $s_i$ il existe un  $x\in E$
tel que $xs_i=y_i$ pour chaque $i$. 
\\
Dans le cas \gnl on fait à peu près la même
chose.
\\
Pour chaque couple $(i,j)$  avec $i\neq j$, le fait que $x_i/1=x_j/1$ dans $E_{ij}$
signifie que pour certains
$u_{ij}\in S_i$ et  $u_{ji}\in S_j$ on a
 $s_j u_{ij} u_{ji} y_i = s_i u_{ij} u_{ji} y_j $.
Pour chaque $i$, soit  $u_i\in S_i$ un multiple commun des $u_{ik}$ (pour $k\neq i$). 
\\
On a
alors $(s_j u_{j}) (u_{i}  y_i) = (s_i u_{i}) (u_{j} y_j) $.
Ainsi le vecteur des $u_{i}  y_i$ est proportionnel au vecteur des $s_i u_{i}$.
Puisque le \sys $(\uS)$ est deux fois \Ergz, il existe \hbox{un $x\in E$} tel que 
$u_{i}  y_i= s_i u_{i} x$ pour tout $i$, ce qui donne les \egtsz~\hbox{$\varphi_i(x)=\fraC{u_{i}  y_i}{s_i u_{i}}=\fraC{  y_i}{s_i} =x_i$}.
\\
Enfin cet $x$ est unique parce que les $S_i$ sont $E$-\corz.
\end{proof}
%---------end proof----------

Voici maintenant une variante du \plgref{plcc.modules 2}.
Cette variante apparaît cette fois-ci comme une réciproque du \plg précédent.

%:--- Prc lgc {plcc.modules2bis}
\begin{plcc}
\label{plcc.modules2bis}\relax {\em (Recollement concret de modules) }
\\
 Soit $(\uS)=(S_1, \dots, S_n)$ un \sys de \mos de $\gA$. \\
 On note $\gA_i=\gA_{S_i}$,
$\gA_{ij}=\gA_{S_iS_j}$ et $\gA_{ijk}=\gA_{S_iS_jS_k}$.
Supposons donné dans la catégorie des \Amos un diagramme commutatif

\snic{\big((E_i)_{i\in \lrbn}),(E_{ij})_{i<j\in \lrbn},(E_{ijk})_{i<j<k\in \lrbn};(\varphi_{ij})_{i\neq j},(\varphi_{ijk})_{i< j,i\neq k,j\neq k}\big)}

(comme dans la figure ci-après)
avec les \prts suivantes. 
\begin{itemize}
\item Pour tous $i$, $j$, $k$ (avec $i<j<k$), $E_i$ est un $\gA_i$-module,  $E_{ij}$ est un~$\gA_{ij}$-module et $E_{ijk}$ est un~$\gA_{ijk}$-module.
Rappelons que selon nos conventions de notation on pose $E_{ji}=E_{ij}$, $E_{ijk}=E_{ikj}=\dots$

\item Pour $i\neq j$,  $\varphi_{ij}:E_i\to E_{ij}$ est un  \molo en $S_j$ (vu dans $\gA_i$).
\item Pour $i\neq k$, $j\neq k$ et $i<j$, $\varphi_{ijk}:E_{ij}\to E_{ijk}$ est un  \molo en $S_k$ (vu dans $\gA_{ij}$).
\end{itemize}
Alors, si $\big(E,(\varphi_i)_{i\in\lrbn}\big)$ est la limite projective du diagramme, on a  les résultats suivants.
\begin{enumerate}
\item 
Chaque morphisme $\varphi_i:E\to E_i$ est
un  \molo en~$S_i$.  
\item 
Le \sys $(\uS)$ est deux fois \Ergz.  
\item  Le \sys $\big(E,(\varphi_{i})_{i\in\lrbn}\big)$ est,
à \iso unique près, l'unique
\sys  $\big(F,(\psi_{i})_{i\in\lrbn}\big)$ avec les $\psi_i\in\Lin_\gA(F,E_i)$
vérifiant les points suivants:
\begin{itemize}
\item le diagramme est commutatif,
\item  chaque $\psi_i$
un  \molo en $S_i$,
\item  le \sys $(\uS)$ est deux fois \Frgz.
\end{itemize}
\end{enumerate}

\smallskip {\small\hspace*{6em}{
$
\xymatrix @R=2em @C=7em{
 E_i\ar[d]_{\varphi _{ij}}\ar@/-0.75cm/[dr]^{\varphi _{ik}} &
     E_j\ar@/-1cm/[dl]^{\varphi _{ji}}\ar@/-1cm/[dr]_{\varphi _{jk}} &
        E_k\ar@/-0.75cm/[dl]_{\varphi _{ki}}\ar[d]^{\varphi _{kj}} &
\\
 E_{ij} \ar[rd]_{\varphi _{ijk}} & 
    E_{ik}  \ar[d]^{\varphi _{ikj}} & 
      E_{jk}  \ar[ld]^{\varphi _{jki}} 
\\
   &  E_{ijk} 
}
$
}}
\end{plcc}
%--- end-plcc-----------------
%
\begin{proof} \emph{1.} Cette \prt est valable sans aucune hypothèse sur le \sys de \mos considéré (voir la \dem du \plgref{plcc.modules 2}).

\emph{2.}  On considère des $s_i\in S_i$ et une suite $(\bmx_i)_{i\in\lrbn}$ dans $E$ proportionnelle à $(s_i)_{i\in\lrbn}$. Notons $\bmx_i=(x_{i1},\dots,x_{in})$. La proportionalité des deux suites signifie que $s_ix_{jk}=s_jx_{ik}$
dans $E_k$ pour tous $i$, $j$, $k$. On pose $\bmx=(\fraC{x_{ii}}{s_i})_{i\in\lrbn}$. On vérifie ensuite que $s_i\bmx=\bmx_{i}$: i.e., que
$s_i\fraC{x_{jj}}{s_j}=x_{ij}$ dans chaque $E_j$. En effet,
cela résulte de l'\egt de proportionalité  $s_ix_{jk}=s_jx_{ik}$ pour $k=j$.

\emph{3.} 
Puisque~$E$ est la limite projective du diagramme, il y a une unique \Ali  $\psi:F\to E$ telle que
$\psi_i=\varphi_i\circ \psi$ pour tout~$i$.
\\ 
En fait on a $\psi(y)=\big(\psi_1(y),\dots,\psi_n(y)\big)$.
\\ Montrons d'abord que $\psi$ est injective. \hbox{Si $\psi(y)=0$}
tous les $\psi_i(y)$ sont nuls, et puisque $\psi_i$ est un  \molo en~$S_i$, il existe des  $s_i\in S_i$ tels que $s_iy=0$. 
Puisque  $(\uS)$ est un \sys \Frgz, on a $y=0$. 
\\
Comme $\psi$ est injective on peut supposer $F\subseteq E$ et $\psi_i=\varphi_i\frt F$. \\
Dans ce cas montrer que $\psi$ est bijective revient à montrer que $F=E$. 
\\
Soit $\bmx\in E$. Comme~$\psi_i$ et~$\varphi_i$
sont deux  \molos en $S_i$, il y a des  
$u_i\in S_i$ tels que $u_i\bmx\in F$. Puisque  $(\uS)$ est deux fois \Frgz, 
et que la suite des $u_i\bmx$ est proportionnelle à la suite des
$u_i$, il existe un $y\in F$ tel \hbox{que $u_i\bmx=u_iy$} pour tout $i$, donc $y=\bmx\in F$.
\end{proof}
%

%%%%%%%%%%%%%%%%%%%%%%%%%%%%%%%%%%%%%%%%%%%%%%%%%%%%%%%%%%%%%%%%%%%%%%%%%%%
%:section: Exercices
%\pagebreak
%\newpage	
\Exercices

%--- Exercise{exoLocMon1}--------
\begin{exercise}
\label{exoLocMon1}
{\rm  
Soient $S_1$, \ldots, $S_n$, $S$ des \mos de $\gA$ tels que $S$ est contenu dans le saturé de chaque $S_i$. \Propeq

 \emph{1.} 
Les $S_i$ recouvrent $S$.

 \emph{2.} 
Les $S_i$ sont \com dans $\gA_S$.
}
\end{exercise}
%--- end-exercise------------------

%--- Exercise{exoLocMon2}--------
\begin{exercise}
\label{exoLocMon2}
{\rm Soit $I$ un \id et $U$ un \mo de $\gA$.
Posons $S = \cS(I , U)$. %Alors:

 \emph{1.} 
 Dans $\gA_S$, le monoïde $\cS(I; U,a)$ est
\eqv à $\cS(I ; a) = I + a^\NN$.

 \emph{2.} 
 Dans $\gA_S$, le monoïde $\cS(I,a ; U)$  est \eqv à $\cS(a ; 1) = 1 +\lra a$.
 }
\end{exercise}
%--- end-exercise-----------------------------------------

%--- Exercise{exoLocMon3}--------
\begin{exercise}
\label{exoLocMon3}
{\rm  Donner une \dem du lemme \ref{lemRecouvre}
basée sur les deux exercices précédents.}
\end{exercise}
%--- end-exercise------------------

%--- Exercise{exoComCom}-------------
\begin{exercise}
 \label{exoComCom}  {\rm Soit $\gA$ un anneau.
 
\emph{1.} Pour  $a_1$, \ldots, $a_n$ dans
$\gA$, si $a_1 \cdots a_n \in \Rad\gA$, les \mos $1 + \gen {a_i}$ sont \comz.

 \emph{2.}
 Si $\fa_1$, \ldots, $\fa_\ell$ sont des \ids de $\gA$,
 les \mosz~\hbox{$1+\fa_i$} recouvrent le \moz~\hbox{$1+\prod_i\fa_i$}.
} 
 \end{exercise}
%--- end-exercise-----------------------------------------

%--- Exercise{exoIdepDyna}-------------
\begin{exercise}
\label{exoIdepDyna} (Conformément à la \dfn des \idepsz, si un produit 
de facteurs est dans un \idep potentiel, 
on peut ouvrir des branches de calcul dans chacune 
desquelles un au moins des facteurs est dans
le nouvel \idep potentiel)

\noindent 
{\rm  On reprend les notations de la \dfn \ref{nota mopf}. 
On considère deux parties $I$ et~$U$ de~$\gA$ et le \mo 
correspondant $\cS(I,U)$. Soient $a_1$, \dots, $a_k\in\gA$ 
pour lesquels on a
$$\preskip-.4em \postskip.4em\ndsp 
\prod_{i=1}^ka_i\in \gen{I}_{\gA_{\cS(I,U)}}. 
$$
\emph{1.} Montrez que les \mos $\cS(I\cup\so{a_i},U)$ recouvrent le \mo $\cS(I,U)$.

\emph{2.}  Si  l'on a $a_i-a_{j}\in \cS(I,U)$,
alors  $a_{j}$  est \iv dans $\gA_{\cS(I\cup\so{a_i},U)}$.
 
\emph{3.}  Supposons que pour chaque $j\in\lrbk$,  on a
un \auto de $\gA$ qui fixe le \mo $\sat{\cS(I,U)}$ et qui envoie $a_1$ sur $a_j$, et que chacun des $\gA_{\cS(I\cup\so{a_i},U)}$
est trivial, alors
$\gA_{\cS(I,U)}$ est trivial. 
}
\end{exercise}
%--- end -exercise-----------------------------------------

%--- Exercise{exoSLequiv}-------------
\begin {exercise}\label{exoMonoidesComax1}
{\rm
Soit $S = (S_1, \ldots, S_n)$ une famille de \mos\ de $\gA$.

  \emph{1.}
On considère la famille $S'$ obtenue à partir de $S$ en répétant
chaque $S_i$ un certain nombre de fois (au moins une):
$$
S' = (S_1, S_1, \ldots, S_2, S_2, \ldots, S_n, S_n, \ldots)
$$
Montrer que $S$ est une famille de \mocoz\ de $\gA$ si et seulement si il en
est de même de $S'$.

  \emph{2.}
On considère une seconde famille $U = (U_1, \ldots, U_m)$ de \mos\ de
$\gA$. On suppose que pour chaque $i \in\lrb{1..n}$, il existe $j
\in\lrb{1..m}$ tel que $S_i \subseteq U_j$ et pour chaque $j \in \lrb{1..m}$ il
existe $i \in\lrb{1..n}$ tel que $U_j \supseteq S_i$. Montrer que si $U$ est une
famille de \mocoz\ de $\gA$, il en est de même de $S$.
}
\end {exercise}
%--- end -exercise-----------------------------------------

%--- Exercise{exoKroLocvar}-------------
\begin{exercise}
\label{exoKroLocvar} (Variation sur le \tho de \KRN local \paref{thKroLoc})\\
{\rm Pour résoudre l'exercice, on peut remarquer que le résultat souhaité
est un énoncé \gui{quasi global} que l'on peut obtenir par relecture
de la \dem du \tho de \KRN local.
\\
 Soient $x_0$, \dots, $x_d\in\gA$ et $\fa = \DA(x_0,\dots,x_d)$. Si $\Kdim \gA \leq d$ et
$\Kdim \gA/\fa \leq 0$,
il existe des \elts $s_0$, \dots, $s_d \in \gA$ et
des \ids $\fb_0$, \dots, $\fb_d$, chacun engendré par~$d$ \eltsz,
tels que\footnote{s.f.i.o.: \sfioz.}

\snic{
(s_0, \dots, s_d) \, \hbox{ est un s.f.i.o. de } \gA/\fa
\;\;\;  \hbox{et} \;\;\;
\forall i,~~s_i \fa \subseteq \sqrt{\fb_i} \subseteq \fa}

%\sni
(localement,
$\fa$ est maximal et radicalement engendré par $d$~\eltsz).
}
\end{exercise}
%--- end -exercise-----------------------------------------

%--- Exercise{exoKroLocvarbis}-------------
\begin{exercise}
\label{exoKroLocvarbis}
(Deuxième variation sur le \tho de \KRN local)\\
{\rm
Soit $\gA$ un anneau et
un \id $\fa$ \tfz.
Si $\Kdim \gA/\fa \leq 0$ et $\Kdim \gA_{1+\fa}\leq d$, \linebreak 
il existe des \elts $s_0,$ \dots, $s_d \in \gA$ et
des \ids $\fb_0$, \dots, $\fb_d \subseteq \fa$, chacun engendré par~$d$ \eltsz,
tels que

\snic{
(s_0, \dots, s_d)\, \hbox{ est un s.f.i.o. de } \gA/\fa
\;\;\;  \hbox{et} \;\;\;
\forall i,~~s_i \fa \subseteq \sqrt{\fb_i}. }
 }
\end{exercise}
%--- end -exercise-----------------------------------------

%--- Exercise{exoFonctDet2}----------
\begin{exercise}
\label{exoFonctDet2}
{\rm
Vu le \prcc des modules (\plgc \ref{plcc.modules 2}), et vu l'\iso canonique 

\snic{\big(\Lin_\gA(M,N)\big)_S\rightarrow \Lin_{\gA_S}(M_S,N_S)}

%\sni
dans le cas de \mpfs (proposition \ref{fact.homom loc pf}), on a des \carns
locales pour le \deter d'un \mptf et celui d'un \homo entre \mptfs
(cf. exercice \ref{exoFonctDet1}).

 \emph{1.} 
 Le module $\det(M)$ est \care à \iso unique près par
la \prt suivante: si $s\in \gA$ est tel que $M_s$ est libre, alors
$\det(M)_s\simeq\det(M_s)$, avec des \isos compatibles lorsque l'on fait une
\lon plus poussée{\footnote{Cela signife \prmtz: si $s''=ss'$,
alors l'\iso $(\det(M))_{s''}\simeq\det(M_{s''})$ est donné par la \lon
de l'\iso $(\det(M))_s\simeq\det(M_s)$.}}.

 \emph{2.} 
 Si $\varphi ~:M\rightarrow N$ est un \homo de
 $\gA$-\mptfsz, l'\homo $\det(\varphi)$  est \care par la \prt
suivante: si $s\in \gA$ est tel que $M_s$ et $N_s$ sont libres, alors
$\det(\varphi)_s=\det(\varphi_s)$ (modulo les \isos canoniques).
}
\end{exercise}
%--- end-exercise-----------------------------------------

%--- Exercise{exoRecolle2PTF}--------
\begin{exercise}
\label{exoRecolle2PTF} ~{\rm Soit $n\geq3$,   $s_1$, $s_2$ deux \eco de $\gA$.
On se propose de recoller concrètement deux \mptfs $P_1$ et $P_2$
définis respectivement sur $\gA_{s_1}$ et $\gA_{s_2}$ qui ont des extensions
à  $\gA_{s_1s_2}$
isomorphes. En utilisant le lemme \dlgz, on peut 
supposer qu'ils sont images de \mprns conjugées $F_1$ et $F_2$ sur  $\gA_{s_1s_2}$
au moyen d'un produit de matrices \elrsz.

 \emph{1.} Soit $E\in\En(\gA_{s_1s_2})$.
Montrer qu'il existe $E_1\in\En(\gA_{s_1})$ et $E_2\in\En(\gA_{s_2})$
tels que $E=E_1E_2$ sur $\gA_{s_1s_2}$.

 \emph{2.} Soient $F_1\in\Mn(\gA_{s_1})$ et $F_2\in\Mn(\gA_{s_2})$ deux
\mprns conjuguées sur $\gA_{s_1s_2}$ au moyen d'une matrice
$E\in\En(\gA_{s_1s_2})$. Que faire?
}
\end{exercise}
%--- end-exercise-----------------------------------------

%:--- Exercise{exoMcCoyContr1}-------------
\begin{exercise}
\label{exoMcCoyContr1} {(\Tho de McCoy contraposé, version pénible)}\\
{\rm 
Soient $\gA$ un anneau discret non trivial et
une matrice $M\in\Ae{m\times n}$.
\begin{enumerate}
\item Si $\cD_n(M)$ est fidèle, $M$ est injective.
\item Si l'on connaît un entier $k<n$ et un $x\in\gA$ non nul, tels
que 

\snic{x\cD_{k+1}(M)=0\hbox{  et  }\cD_k(M)\hbox{  est fidèle},}

\snii alors on peut construire 
un vecteur non nul dans le noyau de $M$.  
\end{enumerate}
 
}
\end{exercise}
%--- end -exercise-----------------------------------------

%:--- Exercise{exoMcCoyContr2}-------------
\begin{exercise}
\label{exoMcCoyContr2} (\Tho de McCoy contraposé, version digeste)\\
{\rm Soient $\gA$ un anneau \coh discret non trivial et
une matrice $M\in\Ae{m\times n}$.
\begin{enumerate}
\item Ou bien $\cD_n(M)$ est fidèle, et $M$ est injective.
\item Ou bien on peut construire  dans le noyau de $M$
un vecteur avec au moins une \coo dans $\Atl$.  
\end{enumerate} 
}
\end{exercise}
%--- end -exercise-----------------------------------------

%%:--- Exercise{exotrickprof2}-------------
%\begin{exercise}
%\label{exotrickprof2}
%{\rm On remarque que les \dfns de la profondeur $1$ et de la profondeur~$2$ sont données en des termes qui ne font pas intervenir
%la structure additive de l'anneau considéré, mais seulement sa structure multiplicative, \cad le \mo $(\gA,\smalltimes,1)$.  
%
%Comme l'énoncé du lemme \ref{lemtrickprof2} ne fait pas non plus usage de la structure additive, on peut espérer une \dem purement multiplicative de ce lemme. L'inspection de la \dem donnée dans le cours montre que ce n'est pas le cas. On propose donc \alec de trouver une \dem du lemme  \ref{lemtrickprof2}
%qui fonctionne pour n'importe quel \moz.
% 
%}
%\end{exercise}
%%--- end -exercise-----------------------------------------

%: notenglish

%:--- Exercise{exoDDMcCoy}-------------
\begin{exercise} 
\label{exoDDMcCoy} (Un autre résultat de McCoy: le lemme de McCoy)
{\rm Le résultat suivant est pour l'essentiel une reprise du 
corolaire \ref{corlemdArtin}.
\emph{Soient $E$ un \Amo et $f$  un \pol de $\AuY$. \Propeq
\begin{enumerate}
\item L'\elt $f$ est $E[\uY]$-\ndzz, i.e. $(0_{E[\uY]}:f)_{E[\uY]}=0_{E[\uY]}$.
\item L'\id $\rc(f)$ est \Ergz, i.e. $\big(0_{E}:\rc(f)\big)_{E}=0_{E}$.
\item Pour tout $x\in E$, $\,x\, f =0\;\Rightarrow\; x=0$, i.e. $(0_{E[\uY]}:f)_{E}=0_{E}$. 
\end{enumerate}
Un tel \pol est appelé un \emph{\pol de Kronecker attaché à~$\fa=\rc_\gA(f)$}.}%
\index{polynome de Kro@\pol de Kronecker!attaché à l'\id $\fa$}

Les points \emph{2} et \emph{3} ont la même signification, et
il est clair que le point \emph{3} est un cas particulier du point \emph{1.}
Ce qui est vraiment l'objet du lemme est l'implication~\hbox{\emph{3} $\Rightarrow$
\emph{1.}}  Il suffit de démontrer le lemme dans le cas des \pols en une seule variable
(le cas \gnl peut s'en déduire en utilisant l'astuce de \KRAz).

\emph{a.} Donner une \dem inspirée de celle du corolaire \ref{corlemdArtin}.

\emph{b.} Donner une \dem directe.
 
}
\end{exercise}
%--- end -exercise-----------------------------------------

%--- Exercise{exoGCDplg}-------------
\begin{exercise}\label{exoGCDplg}
{(Un exemple d'application de la machinerie \lgbe \cov de base)} 
{\rm Nous reprenons l'exercice \ref{exoPgcdNst} et nous proposons une solution
complète fondée sur le \plg de base.
Il s'agit ici de \gnr au cas d'un anneau commutatif arbitraire un résultat
utile en théorie des \cdisz: \emph{si l'on divise un \pol $f(x)$ par le pgcd de~$f$ et $f'$, on obtient un \pol \splz.} 
Pour un anneau arbitraire, on devra supposer que le pgcd de~$f$ et $f'$
existe en un sens fort. 

Les points \emph{1}  et \emph{2a} sont des rappels du points \emph{2} et \emph{3a} de l'exercice \ref{exoPgcdNst}.

\emph{1.} Soit $\gK$ un corps discret, 
$x$ une \idtrz, $f\in\Kx$ un \pol non nul de degré $n\geq 0$, \hbox{$h=\pgcd(f,f')$} et
${f_1}=f/h$. 
Alors $\Res_x({f_1},f_1')\in\gK\eti$, ou, ce qui revient au même,  $1\in\gen{f_1,f_1'}\subseteq \Kx$.
Si en outre $\deg(f)=n$ et  $n!\in\gK\eti$, alors $f$ divise~${f_1}^{n}$.

\emph{2.} Soit $\gk$ un anneau commutatif, 
$x$ une \idtrz, et $f\in\kx$ primitif de degré
formel $n\geq 0$. On suppose que \fbox{l'\id $\gen{f,f'}$ est engendré par un \pol $h$} (\ncrt primitif). %
\begin{enumerate}
\item [\emph{a.}] Montrer qu'il existe des \pols  $u$, $v$, $f_2$, ${f_1}\in\kx$,  satisfaisant
les \egts 

\snic{u{f_1}+vf_2=1 \;\;\hbox{ et }\;\;\cmatrix{u&v\cr-f_2&{f_1}}\cmatrix{f\cr f'}=\cmatrix{h\cr 0}.}

En outre, on a $f_1h=f$ et $f_2h=f'$, de sorte que $f_1$ est \ncrt primitif.

\item [\emph{b.}] On suppose que $\gk$ est un \alo \dcdz. Montrer que $1\in\gen{f_1,f_1'}\subseteq \kx$. Si en outre  $n!\in\gk\eti$,
montrer que $f$ divise $f_1^{n}$.
\item [\emph{c.}] Montrer les mêmes résultats qu'au point \emph{2b}, mais pour un anneau commutatif arbitraire. Utiliser  la machinerie \lgbe de base.
\end{enumerate}

\emph{3.} Question subsidiaire.
Donner une \dem directe du point \emph{2c}
qui n'utilise pas la machinerie \lgbe de base.

}

\end {exercise}
%--- end -exercise-----------------------------------------

%: fin notenglish

%: pb

%--- problem{exoChasserIdeauxPremiers1}-------------
\begin{problem}\label{exoChasserIdeauxPremiers1}
{(\'Eviter les \ids premiers)}\\
{\rm  
 Dans ce \pbz, on examine comment décrypter \cot une \dem classique
qui utilise comme outil de base \gui{aller voir ce qui se passe dans les corps $\Frac(\gA\sur\fp\!)$ pour tous les \ideps $\fp$ de $\gA$}.

  \emph{1.}
Soit $t$ une \idtrz, $a$, $b$, $c_1$, \ldots, $c_n \in \gA$ tels que
$(at+b, c_1, \ldots, c_n)$ soit un \vmd sur $\gA[t,t^{-1}]$.  On veut
montrer que $ab \in \DA(c_1, \ldots, c_n)$.  La \dem suivante, typique en \clamaz,
utilise le principe du tiers exclu et l'axiome du choix.  
Si $ab \notin \DA(c_1, \ldots, c_n)$,
il existe un \id premier $\fp$ avec $c_i \in \fp$ pour $i \in \lrbn$ et
$ab \notin \fp$. Sur le corps $\gK = \Frac(\gA\sur\fp\!)$, puisque $\ov a$ est
non nul, l'\eqn $at + b = 0$ a une unique solution $t = -\ov b/\ov a$, qui est
non nulle car $\ov b$ est non nul; on peut alors définir un morphisme
$\varphi : \gA[t,t^{-1}] \to \gK$ par $t \mapsto -\ov b/\ov a$; $\varphi$
transforme le \vmd $(at+b, c_1t, \ldots, c_nt)$ en le vecteur nul de
$\gK^{n+1}$. Absurde.
\\
 Qu'en pensez vous?

 \emph{2.} Si $\gB$ est un anneau réduit décrire les
unités de $\gB[t,1/t]$. 
\\
 On pourra montrer que si $p$, $q\in\gB[t]$ vérifient $pq=t^m$ (avec
$p=\sum_k p_kt^k$ et~\hbox{$q=\sum_k q_kt^k$}), alors
$1 \in \rc(p)$,  $1 \in \rc(q)$ et:

\snic {
p_kp_\ell=q_kq_\ell=0 \hbox { si } k\neq\ell,\quad
p_kq_\ell=0 \hbox { si } k+\ell\neq m, \quad
p_i = q_i=0 \hbox { si } i>m.
}

En conséquence, pour $k\in\lrb{0..m}$, $\gen {p_k}$ est engendré par un
\idm $e_k$. On dispose alors d'un \sfio $(e_0,\ldots,e_m)$ dans $\gB$ tel que
$\gen{e_k}=\gen{p_k}=\gen{q_{m-k}}$ pour $k\in\lrb{0..m}$ et:

\snic{e_kp=e_kp_kt^k, \; e_kq=e_kq_{m-k}t^{m-k}\;$ et $\;e_k=e_kp_kq_{m-k}.}

 Le résultat est clair lorsque l'anneau est intègre, donc le réflexe en \clama est d'utiliser des \idepsz.
Une solution possible pour décrypter \cot ce raisonnement est d'utiliser le \nst formel (\thref{thNSTsurZ}).

}
\end {problem}
%--- end -problem----------

% fin des exos
%:  solutions

\sol

%%%%%%%%%%%%%%%%%%%%%%%%%%%%%%%%%%%%%%%%%%%%%%%%%%%%%%%%%%%%%%%%%%%%%%%%%%%
\exer{exoLocMon1}
$\emph {2} \Rightarrow \emph {1.}$ Soient $s_1, \ldots, s_n$ avec
$s_i \in S_i$. On veut des $b_i \in \gA$ tels que~\hbox{$b_1s_1 + \cdots + b_ns_n
\in S$}. Le fait que les $S_i$ soient \com dans $\gA_S$ fournit un $s \in S$ et
des $a_i \in \gA$ tels que $a_1s_1 + \cdots + a_ns_n = s$ dans $\gA_S$ ; il
existe donc \linebreak 
un $t \in S$ tel que, dans $\gA$, $(ta_1)s_1 + \cdots + (ta_n)s_n
= ts \in S$.

%%%%%%%%%%%%%%%%%%%%%%%%%%%%%%%%%%%%%%%%%%%%%%%%%%%%%%%%%%%%%%%%%%%%%%%%%%%
\exer{exoLocMon2}~\\
\emph{1.} Un \elt $s$ de $\cS(I; U,a)$
est de la forme $x + ua^k$. 
On voit que $s$ divise dans~$\gA_S$
l'\elt $xu^{-1} + a^k \in \cS(I ;
a)$.

\emph{2.} Un \elt $s$ de $\cS(I,a ; U)$ est de la forme $x
+ ya + u$. On pose $x' =u^{-1}x$ \linebreak 
et $y' = u^{-1}$.
Alors $s$ divise dans $\gA_S$
l'\elt  $x' + y'a + 1$, qui  divise $1 + y''a$,  \linebreak 
où $y'' = (1 +
x')^{-1} y'$.

%%%%%%%%%%%%%%%%%%%%%%%%%%%%%%%%%%%%%%%%%%%%%%%%%%%%%%%%%%%%%%%%%%%%%%%%%%%
\exer{exoLocMon3}
Posons $S = \cS(I,U)$, $S_1 = \cS(I; U,a)$, $S_2 = \cS(I,a ; U)$.  
%Il est clair que $S \subseteq S_i$, $i = 1, 2$.  Il suffit donc de 
Il faut vérifier que $S_1$
et $S_2$ sont \com dans $\gA_S$.  Dans $\gA_S$, $S_1$ est \eqv à
$I + a^\NN$, et~$S_2$ est \eqv  à $1 + \langle a\rangle$.  Utilisons l'\idt suivante:

\snic{y^k(x + a^k) + \bigl(\som_{j < k} y^j a^j\bigr)(1 - ya) =
y^k(x + a^k) + 1 - y^ka^k = 1 + y^kx.}

%\sni
Appliquée a $x \in I$, elle prouve que $x + a^k$ et $1 - ya$
sont \com dans $\gA_S$ (puisque $1 + y^kx \in 1 + I$
et que $I$ est contenu dans le radical de $\gA_S$).

%%%%%%%%%%%%%%%%%%%%%%%%%%%%%%%%%%%%%%%%%%%%%%%%%%%%%%%%%%%%%%%%%%%%%%%%%%%

%%%%%%%%%%%%%%%%%%%        exoComCom       %%%%%%%%%%%%%%%%%%%%%%
\exer{exoComCom} \emph{1.} Pour $j\in\lrbn$  soit $b_j=1-a_jx_j$ dans le \moz~\hbox{$1+a_j\gA$}. Soit~\hbox{$a=\prod_ia_i$}.
On doit montrer que l'\id $\fm=\gen{b_1,\dots,b_n}$ contient $1$.
\\
Or, modulo  $\fm$  on a $1=a_jx_j$, donc $1=a\prod_ix_i=ax$. Ainsi $1-ax\in\fm$, \hbox{mais $1-ax\in\Ati$} car $a\in\Rad\gA$.

\emph{2.}
Il est clair que $S=1+\prod_i\fa_i\subseteq1+\fa_j=S_j $ pour chaque $j$. On doit 
donc vérifier (exercice \ref{exoLocMon1}) que les $1+\fa_j$ vus dans $\gA_S$
sont \comz. Or le produit $\prod_i\fa_i$, vu dans $\gA_S$, est dans $\Rad\gA_S$.
Donc il suffit d'appliquer le point \emph{1.}

%%%%%%%%%%%%%%%%%%%%%%%%%%%%%%%%%%%%%%%%%%%%%%%%%%%%%%%%%%%%%%%%%%%%%%%%%%%

%%%%%%%%%%%%%%%%%%%        exoIdepDyna       %%%%%%%%%%%%%%%%%%%%%%

\exer{exoIdepDyna} L'hypothèse signifie que l'on a un $u\in\cM(U)$ et un $j\in\gen{I}_\gA$ tels que 
%
%\snic{
$(u+j)\,\prod_{i=1}^ka_i\in\gen{I}_\gA$,
%}
%
 ou encore: $u\,\prod_{i=1}^ka_i\in\gen{I}_\gA$.

\sni
\emph{1.} 
\emph{Première solution, par calcul direct.}
\\On considère des $x_i\in S_i=\cS(I\cup\so{a_i},U)$ et on cherche une \coli qui se trouve dans  $S=\cS(I,U)$.
Pour chaque $i$ on écrit
$$\preskip.2em \postskip.2em 
x_i=u_i+j_i+a_iz_i\hbox{  avec  }u_i\in \cM(U),\;j_i\in\gen{I}_\gA\hbox{  et  }z_i\in\gA. 
$$
Dans le produit
$$\preskip-.2em \postskip.3em\ndsp 
u\,\prod_{i=1}^k(x_i-(u_i+j_i))=u\;\prod_{i=1}^ka_iz_i\in\gen{I}_\gA , 
$$
on réécrit le premier membre sous forme

\snic{\som_{i=1}^kc_ix_i\pm u\,\prod_{i=1}^k(u_i+j_i),}

%\sni 
et l'on obtient, en faisant passer $\pm u\,\prod_{i=1}^k(u_i+j_i)$
dans le second membre, l'appartenance souhaitée $\som_{i=1}^kc_ix_i\in \cS(I,U)$.

%\sni 
\emph{Deuxième solution, conceptuelle.}
\\ 
Il est clair que $S \subseteq S_i$. Il suffit donc (exercice \ref
{exoLocMon1}) de montrer que les  $S_i$ sont \com dans $\gA_S$.
Dans $\gA_S$, (exercice~\ref {exoLocMon2}) les \mos
$S_i$ et $1 + \gen {a_i}$ ont même saturé. Il suffit donc de voir
que les \mos $1 + \gen {a_i}$ sont comaximaux dans $\gA_S$.
Par ailleurs, on sait que, vu dans $\gA_S$, $I$ est contenu dans
$\Rad(\gA_S)$. On applique donc le point \emph{1} de l'exercice~\ref{exoComCom}.

%\sni 
\emph{2.} Clair puisque $a_{j}\in \cS(I,U)+\gen{a_i}\subseteq \cS(I\cup\so{a_i},U)$.

%\sni 
\emph{3.} Si l'un des $\gA_{S_i}$ est trivial, 
ils le sont tous car ils sont
deux à deux isomorphes. Puisque les $S_i$ recouvrent $S$, 
$\gA_S$ est lui même trivial.

%%%%%%%%%%%%%%%%%%%%%%%%%%%%%%%%%%%%%%%%%%%%%%%%%%%%%%%%%%%%%%%%%%%%%%%%%%%
\exer{exoMonoidesComax1}
 \emph{1.}
Il suffit de le montrer pour $S' = (S_1, S_1, S_2, \ldots, S_n)$.  \\
Supposons la famille $S$ comaximale. Soient $s'_1$, $s''_1 \in S_1$, et $s_i \in S_i$ pour $i \in
\lrb{2..n}$. Les \elts $s'_1s''_1$, $s_2$, \ldots, $s_n$ sont \comz, et
puisque $s'_1s''_1 \in \gen {s'_1, s''_1}$, il en est de même de $s'_1$, $
s''_1$, $s_2$, $\ldots$, $s_n$. Dans l'autre sens, supposons $S'$ comaximale et
soient $s_i \in S_i$ pour $i \in \lrb{1..n}$; alors $s_1$, $s_1$, $s_2$, $\ldots$, $s_n$ sont comaximaux donc il en est de même de $s_1$, $s_2$, $\ldots$, $s_n$.

  \emph{2.}
En répétant certains des $S_i$ et des $U_j$, on obtient deux familles
$S'$, $U'$ de \mos\ de $\gA$, indexées par le même intervalle $\lrb{1..p}$
et vérifiant $S'_k \subseteq U'_k$ \linebreak 
pour $k \in \lrb {1..p}$. Puisque $U$ est
comaximale, il en est de même de $U'$ donc de $S'$ puis de~$S$.

%%%%%%%%%%%%%%%%%%%%%%%%%%%%%%%%%%%%%%%%%%%%%%%%%%%%%%%%%%%%%%%%%%%%%%%%%%%

\exer{exoKroLocvar}
Comme $\Kdim \gA\leq d$,
il existe une suite $(\underline{a}) = (a_0,\dots,a_d)$ \cop
à $(\underline{x}) = (x_0,\dots,x_d)$.
Donc (suites disjointes), pour tout $i \leq d$, on~a

\snic{\rD(a_0, \dots, a_{i-1}, x_0, \dots, x_{i-1}, a_i x_i)
= \rD(a_0 + x_0, \dots, a_{i-1} + x_{i-1}).}

%\sni
Comme $\Kdim \gA/\fa \leq 0$ et $\fa=\DA(\fa)$, on a \egmt
$\gA = \gA a_i + (\fa : a_i)$ pour tout~$i$. On construit alors le triangle :

\snuc{\arraycolsep2pt\begin{array}{rl}
& \gA \\
  = & \gA a_0 + (\fa : a_0) \\
  = & \gA a_0 + (\fa : a_0)a_1 + (\fa : a_0)(\fa : a_1) \\
  & \qquad \vdots \\
  = & \gA a_0 + (\fa : a_0)a_1 + \cdots + (\fa : a_0)\cdots(\fa : a_{d-1}) a_d
 + (\fa : a_0)\cdots(\fa : a_d). 
\end{array}
}

%\sni
Maintenant, on écrit
$$
%\snic{
1 = b_0 a_0 + b_1 a_1 + \cdots + b_d a_d + t
%}
$$
%\sni
avec $b_i \in (\fa : a_0)\cdots(\fa : a_{i-1})$ et
$t \in (\fa : a_0)\cdots(\fa : a_d)$. Pour $i \leq d$, on a d'une part

\snic{
\begin{array}{c}
b_i \gen{x_0, \dots, x_{i-1}}
\subseteq \rD\big( b_i(a_0+x_0), \dots, b_i(a_{i-1}+x_{i-1}) \big)
  \hspace*{3cm}    \\[1mm] \hspace*{3cm}
\subseteq \rD(  b_ia_0,x_0, \dots, b_ia_{i-1},x_{i-1} ) \subseteq \rD( \fa ) =\fa,
\end{array}
}

%\sni
et d'autre part

\snic{
\begin{array}{c}
b_i a_i x_i \in b_i \rD(a_0+x_0,\dots, a_{i-1}+x_{i-1})
  \hspace*{5cm}    \\[1mm] \hspace*{4cm}
\subseteq \rD\big( b_i(a_0+x_0), \dots, b_i(a_{i-1}+x_{i-1}) \big)
\subseteq \fa.
\end{array}
}

%\sni
Posons maintenant $s_i = b_i a_i$. On arrive ainsi à

\snic{
s_i \gen{x_0, \dots, x_{i-1}, x_i}
\subseteq \rD\big(b_i(a_0+x_0),\dots, b_i(a_{i-1}+x_{i-1}) \big) \subseteq \fa
,}

%\sni
puis

\snuc{
\begin{array}{c}
s_i \gen{x_0, \dots, x_{i-1}, x_i, x_{i+1}, \dots, x_d}
\subseteq
  \hspace*{6cm}    \\[1mm] \hspace*{3cm}
\rD\big( \underbrace{b_i(a_0+x_0), \dots, b_i(a_{i-1}+x_{i-1}),
     s_i x_{i+1}, \dots, s_i x_d}_{\hbox{\scriptsize engendrent }\fb_i
\hbox{\scriptsize ~(def.) }} \big)
\subseteq \fa.
\end{array}
}

%\sni
Il existe donc un idéal $\fb_i$ engendré par $d$ \elts vérifiant

\snic{
s_i \fa \subseteq \rD(s_i \fa) \subseteq \rD( \fb_i ) \subseteq \fa.}

%\sni
On termine la \dem
en utilisant $1 \in \gen{a_d,x_d}$, il vient

\snic{
t \in t \gen{a_d, x_d} \subseteq \gen{ta_d, x_d} \subseteq \fa,
}

%\sni
si bien que la somme des $s_i$ vaut $1 \mod \fa$.
Par ailleurs, pour $i>j$,

\snic{
s_i s_j \in b_i a_j \gA \subseteq \fa,
}

%\sni
ce qui permet de conclure que $(s_0,\dots,s_d)$ est un \sfio de~$\gA/\fa$.

%%%%%%%%%%%%%%%%%%%%%%%%%%%%%%%%%%%%%%%%%%%%%%%%%%%%%%%%%%%%%%%%%%%%%%%%%%%
\exer{exoKroLocvarbis}
Pour commencer, le \tho de \KRN donne que $\sqrt{\fa}$ est
radicalement engendré par $d+1$ \eltsz.
Ensuite on applique le résultat de l'exercice \ref{exoKroLocvar} à $\sqrt{\fa}$ dans le localisé
$(1+\fa)^{-1} \gA$.
Cela fournit des~$s_i$ formant un \sfio modulo~$\sqrt{\fa}$ et
des~$\fb_i \subseteq \sqrt{\fa}$.
Quitte à prendre des multiples de puissances des~$s_i$, on peut imposer
que $(s_0, \dots, s_d)$ est un \sfio de~$\gA/\fa$.
Quitte à prendre des puissances des
\gtrs des~$\fb_i$, on peut imposer $\fb_i \subseteq \fa$.

%%%%%%%%%%%%%%%%%%%%%%%%%%%%%%%%%%%%%%%%%%%%%%%%%%%%%%%%%%%%%%%%%%%%%%%%%%%
\exer{exoRecolle2PTF} \emph{1.} \cite[page 208 proposition 1.14]{Lam06}.

\emph{2.} On a $F_1E=EF_2$, on écrit $E=E_1E_2$, donc  $F_1E_1E_2=E_1E_2F_2$
et 

\snic{\wi{E_1}F_1E_1=_{\gA_{s_1s_2}}E_2F_2\wi{E_2}.}

%\sni
La matrice $\wi{E_1}F_1E_1$ (resp. $E_2F_2\wi{E_2}$) est une \mprn  sur  $\gA_{s_1}$ (resp. sur~$\gA_{s_2}$) car $\wi{E_1}E_1=\In$ (resp. $\wi{E_2}E_2=\In$). Par le \plg de \rcm des \elts dans un module (ici $\Mn(\gA)$),
il existe une unique matrice~\hbox{$F\in\Mn(\gA)$} qui est égale à
$\wi{E_1}F_1E_1$ sur $\gA_{s_1}$ et à $E_2F_2\wi{E_2}$ sur $\gA_{s_2}$.
Pour vérifier que $F^2=F$, il suffit de le vérifier sur  $\gA_{s_1}$
et $\gA_{s_2}$. Soit $P=\Im F\subseteq\Ae n$.
\\
Par construction, pour $i=1,2$ on obtient

\snic{P_{s_i}= _{\gA_{s_i}^n}\Im F_{s_i}\simeq _{\gA_{s_i}^n}\Im F_i\simeq _{\gA_{s_i}^n}P_i.}

\exer{exoMcCoyContr1} \emph{(\Tho de McCoy contraposé, version pénible)}\\
\emph{1.} Déjà vu.\\
\emph{2.}  On a $x\neq 0$,  $\cD_k(M)$ est fidèle et l'anneau est discret, 
donc il existe un mineur~$\mu$ d'ordre $k$ de $M$ tel que $x\mu\neq 0$.
Supposons par exemple que $\mu$ soit le mineur nord-ouest et notons $C_1$, \dots, $C_{k+1}$ les premières colonnes de $M$, notons~$\mu_i$ ($i\in\lrbk$) les \deters convenablement signés des matrices extraites sur les lignes $\lrbk$
et les colonnes précédentes, sauf la colonne d'indice $i+1$. Alors
les formules de Cramer donnent l'\egt $\sum_{i=1}^{k}x \mu_iC_i+x\mu C_{k+1}=0$. Comme $x\mu\neq 0$, ceci donne un vecteur non nul dans le noyau de $M$.

\comm En \clamaz, si $\cD_n(M)$ n'est pas fidèle, comme $\cD_0(M)=\gen{1}$
est fidèle, il existe un $k<n$ tel que $\cD_k(M)$ est fidèle 
et~$\cD_{k+1}(M)$ n'est pas fidèle. Toujours en \clamaz, si $\cD_{k+1}(M)$ n'est pas fidèle, il existe un $x\neq 0$ tel que $x\cD_{k+1}(M)=0$.
Pour que ces choses deviennent explicites, il faut par exemple disposer d'un test pour la fidélité des \itfsz, en un sens assez fort. 
\eoe

%%%%%%%%%%%%%%%%%%%%%%%%%%%%%%%%%%%%%%%%%%%%%%%%%%%%%%%%%%%%%%%%%%%%%%%%%%%

\exer{exoMcCoyContr2}\emph{(\Tho de McCoy contraposé, version digeste)}\\
Puisque les \idds sont des \itfsz, et l'anneau est \cohz, leurs annulateurs
sont \egmt des \itfsz. Et l'on peut tester la nullité d'un \itf parce que l'anneau est discret. Les hypothèses de l'exercice \ref{exoMcCoyContr1} sont donc satisfaites.\\
NB: l'alternative \gui{\emph{1} ou \emph{2}} est exclusive car l'anneau est non nul, ceci justifie le \gui{ou bien, \dots, ou bien} de l'énoncé.

%%%%%%%%%%%%%%%%%%%%%%%%%%%%%%%%%%%%%%%%%%%%%%%%%%%%%%%%%%%%%%%%%%%%%%%%%%%

%%%%%%%%%%%%%%%%%%%%%%%%%%%%%%%%%%%%%%%%%%%%%%%%%%%%%%%%%%%%%%%%%%%%%%%%%%%

%:sinotenglish
\sinotenglish{

\exer{exoDDMcCoy}~
\emph{a.} 
Comme indiqué dans le corolaire \ref{corlemdArtin}, 
c'est une conséquence facile du lemme de \DKM \ref{lemdArtin}. \\
En fait nous utilisons la variante \gui{pour les modules} du lemme de Dedekind-Mertens: \emph{Pour $f\in\AT$ et $g\in E[T]$ avec $m\geq\deg g$ on a}

\snic{\rc(f)^{m+1}\rc(g)=\rc(f)^m\rc(fg).}

\snii Ici le contenu $\rc(g)$ du \pol $g\in E[T]$ est le sous-\Amo de $E$
engendré par les \coes de $g$.
La variante pour les modules est une conséquence du lemme lui-même.
En effet ce lemme peut être vu comme une famille 
d'\idas reliant les \coes de $f$
et de $g$. Ces \idas sont toutes \lins par rapport aux \coes de~$g$.
Or toute \ida \lin par rapport à certaines des \idtrs peut être évaluée en spécialisant les \idtrs dans un anneau arbitraire $\gA$, sauf 
certaines des \idtrs \gui{\linsz} qui peuvent être spécialisées en des \elts d'un \Amo arbitraire~$E$. 

\emph{b.} 
On a $f\in\AT$ et $g\in E[T]$. On suppose $fg=0$ et on veut montrer $g=0$.
On fait une \recu  sur le degré formel $m$ de $g$. Pour  
$m=0$ le résultat est clair. Passons de $m-1$ à $m$. Appelons~$f_i$ les \coes de $f$ et $g_j$ ceux de $g$.  
On va montrer que $f_ig$ est nul pour chaque~$i$. Cela implique alors 
que tous les $f_ig_j$ sont nuls, et puisque $\Ann_E(\rc(f))=0$ que tous les 
$g_j$ sont nuls.
Pour montrer que~$f_ig$ est nul on fait une \recu descendante sur~$i$, 
qui s'initialise avec~$i=n+1$ sans \pb ($n$ est le degré formel de~$f$). Supposons donc avoir montré 
que~$f_ig=0$ pour tous les~$i>i_0$ et montrons le pour~$i_0$.
On a 
  $$\ndsp\preskip.4em \postskip.4em 
      fg = \left(\som_{i\leq i_0}f_iX^i\right) \,g 
  $$
Le \coe de degré~$m+i_0$ de ce \pol est égal à~$g_mf_{i_0}$ et
donc~$g_mf_{i_0}=0$. Donc le \pol $\tilde{g}=f_{i_0}g$ est de degré 
$\leq m-1$ et, évidemment,~$f\tilde{g}=f_{i_0}fg=0$. On peut donc 
appliquer l'\hdr avec~$\tilde{g}$, on conclut~$f_{i_0}g=0$, et on a gagné.

%%%%%%%%%%%%%%%%%%%%%%%%%%%%%%%%%%%%%%%%%%%%%%%%%%%%%%%%%%%%%%%%%%%%%%%%%%%

\exer{exoGCDplg} \emph{1} et \emph{2a.}
Voir la solution des points \emph{2} et \emph{3a} de l'exercice \ref{exoPgcdNst}.

\emph{2b}. Dans le point \emph{2a}, le degré formel de $f_1$ est le même que celui de $f$. Sur un \cdi $\gK$, un \pol $f_1$ non nul de degré $\leq n$ vérifie $1\in\gen{f_1,f'_1}$ \ssiz~$1\in M$, où $M$ est le sous-\Kevz\footnote{En fait, si $1\in\gen{f_1,f'_1}$,
alors l'inclusion est une \egtz.} 
$$
\preskip.4em \postskip.0em 
M=\!\!\sum_{k\in\lrb{0..n-2}}\!\!\!X^{k}f_1\,\gK\;+\!\!\!\!\sum_{\ell\in\lrb{0..n-1}}\!\!\!\!X^{\ell}f'_1\,\gK\;\subseteq\; \sotq{g\in\Kx}{\deg(g)\leq 2n-1}\simeq\gK^{2n}. 
$$
Si maintenant $\gk$ est un \alo \dcdz, on peut considérer le
sous-\kmo  $M\subseteq \sotq{g\in\gk[x]}{\deg(g)\leq 2n-1}\simeq\gk^{2n}$ analogue,
ainsi que le sous-\kmo $N=\gk+M$. On montre que $N=M$, donc que $1\in\gen{f_1,f'_1}$, en appliquant le lemme de Nakayama. Le sous-module $M$ de $N$
est égal à $N$ parce qu'il lui est \rdt égal: \rdtz, on est ramené à la situation d'un \cdiz, car le \pol $f_1$ reste primitif, donc est \rdt non nul.

Si en outre  $n!\in\gk\eti$,
on veut montrer que $f$ divise $f_1^{n}$. On procède de la même manière. Dans le cas d'un \cdi $\gK$, avec $f$ et $f_1$ non nuls de degrés $\leq n$, $f$ divise $f_1^{n}$ \ssi $f_1^{n}\in P$, où
$$\preskip.4em \postskip.4em 
P=\sum_{k\in\lrb{0..n^{2}-n}}\!\! X^{k}f\,\gK\;\subseteq\;  \sotQ{g\in\Kx}{\deg(g)\leq n^{2}}.
$$
Pour le cas d'un \alo \dcdz, on termine de la même manière avec le lemme de Nakayama.

\emph{2c}. Puisqu'il s'agit de démontrer des \egts $E=F$, où $E\subseteq F$,
pour des \kmos $E$ et $F$, la machinerie \lgbe de base s'applique: il suffit de démontrer une telle \egt après \lon en des \mocoz.
Ces \moco sont fournis par la relecture de la \dem donnée dans le cas d'un \alo \dcdz.

\emph{3.} Merci \alec qui nous indiquera une solution plus \elr de l'exercice.

%%%%%%%%%%%%%%%%%%%%%%%%%%%%%%%%%%%%%%%%%%%%%%%%%%%%%%%%%%%%%%%%%%%%%%%%%%%
}
%: fin sinotenglish

%:   sol pb
%%%%%%%%%%%%%%%%%%%%%%%%%%%%%%%%%%%%%%%%%%%%%%%%%%%%%%%%%%%%%%%%%%%%%%%%%%%

%%%%%%%%%%%%%%%%%%%%%%%%%%%%%%%%%%%%%%%%%%%%%%%%%%%%%%%%%%%%%%%%%%%%%%%%%%%

\prob{exoChasserIdeauxPremiers1}
\emph{1.}
Il n'y a pas de miracle: un certificat pour $ab \in \DA(c_1,
\ldots, c_n)$ peut être obtenu à partir d'un certificat d'unimodularité
de $(at+b, c_1t, \ldots, c_nt)$ dans $\gA[t,t^{-1}]$.
\\ 
En remplaçant $\gA$ par $\gA_1=\gA\sur{\DA(c_1, \ldots, c_n)}$, on se ramène
à $\gA$ réduit \hbox{et $c_i = 0$}. L'hypothèse est alors $at+b$ \iv dans
$\gA[t,t^{-1}]$, et le résultat à montrer est~\hbox{$ab = 0$} (résultat \smq en
$a,b$ tout comme l'hypothèse). \\
On a $(at+b)g(t) = t^e$ pour un $g \in
\gA[t]$ et un $e \in \NN$, donc le \pol $at + b$ est primitif. Pour montrer $ab =
0$, il suffit de localiser en $a$ puis en $b$. Sur le localisé en $a$, on
prend $t = -b/a$ dans $(at+b)g(t) = t^e$, on obtient $(-b/a)^e = 0$.
Ainsi, on~a~$b = 0$, puis $ab = 0$. Par symé\-trie, on obtient dans $\gA_b$, $a =
0$ donc $ab=0$. \\
En fait, si $ua+vb=1$, l'anneau $\gA_1$ est cassé en deux par l'\idm
$ua$. Dans la première composante, $at+b=a$ avec $a$ \ivz, dans la seconde,
 $at+b=b$ avec $b$ \ivz.
 
 \emph{2.} En \clamaz: si l'on passe au quotient par un \idep le résultat est clair.
Par continuité, le spectre est partitionné en un nombre fini d'ouverts correspondant au \sfio convoité. 
\\
 Une \dem \cov est donnée dans \cite[Yengui]{Ye0}. \Llec pourra aussi
 s'inspirer de la \dem du point \emph{1}.
\\
 Une méthode que l'on peut utiliser de manière systématique  consiste à faire appel au \nst formel (\thref{thNSTsurZ}).
\\ 
Dans le cas présent, on note que le \pb revient à démontrer que les $p_kp_\ell$ sont nuls pour $k\neq\ell$
et que les $p_{m+r}$ sont nuls pour $r>0$. Une fois ceci constaté, puisque les
$p_k$ sont \comz, on obtient  un \sfio $(e_0,\ldots,e_m)$ tel que $e_kp=e_kp_kt^k$
pour tout $k$, ce qui permet de conclure.
\\ 
La philosophie est la suivante: si l'on prend tous les \coes du \pb comme des \idtrs
sur $\ZZ$,  l'hypothèse revient à passer au quotient par le radical $\fa$ d'un \itfz, qui représente les hypothèses. Le but est alors de démontrer que les conclusions sont \egmt dans~$\fa$. Pour cela il suffit de vérifier que c'est bien le cas lorsque l'on évalue le \pb dans  un corps fini arbitraire.
\\
 Ici les \idtrs sont $p_0$, \ldots, $p_n$, $q_0$, \ldots, $q_n$.\\
 Pour
$p = \sum_{k=0}^n p_kt^k$ et $q = \sum_{k=0}^n q_kt^k$;
on définit le \pol

\snic{\sum r_jt^j\eqdefi pq-t^m}

%\sni
(avec les  $r_j\in\ZZ[p_0,\ldots,p_n,q_0,\ldots,q_n]$)  et l'\id $\fa$ est $\rD(r_0,\ldots,r_{2n})$. 
On va montrer
que les  $p_kp_\ell$ et $q_kq_\ell\in\fa$ si $k\neq\ell$, 
que les $p_kq_\ell\in\fa$ pour $k+\ell\neq m$
et que les  $p_{m+r}$ et $q_{m+r}\in\fa$ pour $r>0$.
\\
 Or cela résulte directement du point  \emph{2} dans le \nst formel
(ou alors du point \emph{4} dans le corolaire~\ref{corthNSTsurZ}).
\\
 En termes \gmqsz: si $n\geq m$, 
la \vrt des zéros de $pq-t^m$ sur un corps $\gK$ est un espace formé de 
$m+1$ copies de $\gK\eti$ isolées les unes des autres; 
sur un anneau réduit la réponse est fondamentalement la même, 
mais les composantes isolées dans le cas des corps font ici 
apparaître un \sfioz.

 Ainsi, le \nst formel (\tho \ref{thNSTsurZ}) fournit une méthode \cov pour décrypter les \algos cachés dans certains raisonnements des \clamaz,  lorsque l'argument
consiste à aller voir ce qui se passe dans tous les $\Frac(\gA\sur\fp\!)$ pour tous les \ideps de $\gA$.

%%%%%%%%%%%%%%%%%%%%%%%%%%%%%%%%%%%%%%%%%%%%%%%%%%%%%%%%%%%%%%%%%%%%%%%%%%%
% fin des solutions d'exos

%:   ---- Section*{references}-----------
\Biblio

La méthode dynamique telle qu'elle est expliquée
dans la machinerie \lgbe à \ideps (section \ref{secMachLoGlo}) consiste pour l'essentiel à mettre à
plat les calculs qui sont impliqués par la \emph{méthode
de l'évaluation dynamique} donnée dans \cite[Lombardi]{Lom97},
héritière de la méthode dynamique 
mise en {\oe}uvre dans \cite[Coste\&al.]{clr}
pour des preuves du type \nstz, elle-même héritière de l'évaluation
dynamique à la D5 en calcul formel \cite[Duval\&al.]{D5}.\imlb
Par rapport à ce qui est proposé dans \cite{clr,Lom97}, 
la différence dans le chapitre présent
est surtout que nous avons évité  la référence à la logique formelle.

En \clamaz, on trouve le \prcc des \mptfs 
(point \emph{\ref{iptfplcc.ptf}} du \plgref{plcc.ptf}) par exemple dans \cite[proposition~2.3.5 et lemme~3.2.3]{Kni} (avec une \dem
presque entièrement \covz) et dans \cite[règle~1.14 du chapitre~IV]{Kun}.

Le traitement constructif du lemme de \Sus \ref{lemSuslin1} est d\^u à Ihsen Yengui \cite{Y1}, qui donne la clé de la machinerie \lgbe à \idemasz.

La machinerie \lgbe à \idemis est due à Thierry Coquand \cite[On seminormality]{coq}.

Les exercices \ref{exoKroLocvar} et \ref{exoKroLocvarbis} sont dus à
Lionel Ducos
\cite{Duc08}.\perso{penser à demander la référence à Lionel}

La méthode dynamique a été appliquée pour le calcul de \gui{bases de Gr\"obner dynamiques} par Yengui dans \cite{Y3}.

\newpage \thispagestyle{CMcadreseul}
\incrementeexosetprob

%:        %%%%%%%%%%%%%%%%%%%%%%%%%%%%%%%%%%%%
%:        %%%%%%%%%%%%%%%%%%%%%%%%%%%%%%%%%%%%
%---- Chapitre {Modules projectifs étendus}
\chapter[Modules projectifs étendus]{Modules projectifs étendus}
\label{ChapMPEtendus}
%------------------
\vskip-1em

\minitoc

\subsection*{Introduction}
\addcontentsline{toc}{section}{Introduction}
%-----------------------------------------
Dans ce chapitre nous établissons de manière \cov
quelques résultats importants
concernant les situations où les \mptfs sur un anneau de \pols sont étendus
depuis l'anneau de base.

Nous traitons notamment le \tho de Traverso-Swan (section \ref{sec.Traverso}), le recollement à la Vaserstein-Quillen (section \ref{subsecQPatch}), les \thos de Horrocks (section \ref{sec.Horrocks}), le \tho de Quillen-\Sus (section \ref{sec.QS}),
et dans la section \ref{sec.Etendus.Valuation}, le \tho de Bass \ref{thBass.Valuation} et le \tho
de Lequain-Simis \ref{thLSValu}. 

%%%%%%%%%%%%%%%%%%%%%%%%%%%%%%%%%%%%%%%%%%%%%%%%%%%%%%%%%%%%%%%%%%%
%%%%%%%%%%%%%%%%%%%%%%%%%%%%%%%%%%%%%%%%%%%%%%%%%%%%%%%%%%%%%%%%%%%
%%%%%%%%%%%%%%%%%%%%%       MODULES ETENDUS           %%%%%%%%%%%%%
%%%%%%%%%%%%%%%%%%%%%%%%%%%%%%%%%%%%%%%%%%%%%%%%%%%%%%%%%%%%%%%%%%%
%%%%%%%%%%%%%%%%%%%%%%%%%%%%%%%%%%%%%%%%%%%%%%%%%%%%%%%%%%%%%%%%%%%

%---- Section{Modules etendus}---
\section{Modules étendus} \label{sec etendus}

\'Etant donnée une \alg $\gA\vers{\rho}\gB$, l'\eds de
$\gA$ à $\gB$ transforme un
module $M$  sur $\gA$ en un module
$\rho\ist(M)\simeq \gB\otimes_{\gA}M$ sur $\gB$.
Rappelons qu'un \Bmo isomorphe à un tel module $\rho\ist(M)$ est dit
étendu depuis~$\gA$. On dit aussi qu'il provient du \Amo $M$ par extension
des scalaires.

Dans le cas d'un \mpfz, du point de vue
des \mpnsz, cela correspond à considérer la  matrice
transformée par l'\homo $\rho$.

Une condition \ncr pour qu'un \mpf soit étendu est que
ses idéaux de Fitting soient de la forme $\rho(\fa_{i}) \gB$
pour des \itfs $\fa_{i}$ de $\gA$.
Cette condition est réalisée pour les \mptfs \ssi les
\idms de $\gB$ sont tous images d'\idms de~$\gA$.

%:   subsec{Le problème de l'extension}
\subsec{Le problème de l'extension}

Pour les \mptfsz, le \pb suivant se pose naturellement au vu du morphisme  $\GKO \rho:\GKO \gA\to \GKO \gB$.

\smallskip  \PB{\num1} Tout \mptf sur $\gB$
provient-il d'un \mptf sur $\gA$?
Ou encore: $\GKO\rho$ est il surjectif?

\medskip Rappelons que $\GKO \Ared=\GKO \gA$ et  $\GKO \gB\red=\GKO \gB$,
de sorte que le \pb de l'extension des \mptfs peut être restreint au cas des anneaux réduits.
Par ailleurs, si $\HO\rho:\HO\gA\to\HO\gB$ n'est pas surjectif, la réponse
au \pb \num1 est négative \gui{pour une mauvaise raison}
et le \pb suivant est alors plus naturel.

\smallskip \PB{\num2}  Tout \mrc sur $\gB$
provient-il d'un \mptf sur $\gA$?

\medskip Pour les \mpfs la \gnn naturelle du \pb précédent est alors la suivante.

\smallskip \PB{\num3} Tout \mpf  sur $\gB$ dont les \idfs sont extensions d'\itfs de $\gA$
provient-il d'un \mpf sur $\gA$?

%%%%%%%%%%%%%%%%%%%%%%%%%%%%%%%%%%%%%%%%%%%%%%%%%%%%%%%%%%%%%%%%%%%
%%%%%%%%%%%%%%%%%%%%%%%%%%%%%%%%%%%%%%%%%%%%%%%%%%%%%%%%%%%%%%%%%%%
%%%%%%%%%%%%%%%%%%%%%%%%%%%%%%%%%%%%%%%%%%%%%%%%%%%%%%%%%%%%%%%%%%%
%%%%%%%%%%%%%%%%%%%%%%%%%%%%%%%%%%%%%%%%%%%%%%%%%%%%%%%%%%%%%%%%%%%
%:   subsec{Cas des anneaux de pols}
\subsec{Cas des anneaux de \polsz}

Soit $\gB=\AXr=\AuX$.
Si $(\ua)\in\Ae r$ nous noterons $\ev_\ua$ l'\homo d'évaluation
en $\ua$:

\snic{\ev_\ua\,:\,\gB\to\gA,\; p\mapsto p(\ua).}

%\sni
Les deux \homos $\gA\vers{j}\gB\vvers{\ev_\ua~}\gA$
 se composent selon l'identité.

La plupart de ce qui suit dans ce paragraphe pourrait être écrit dans le
cadre plus \gnl d'une \Alg $\gB$ qui possède un charactère (cf. proposition~\ref{prdfCaracAlg}).
Pour les anneaux de \pols
 on obtient les résultats suivants (avec une notation intuitive évidente pour $M(\uX)$).

%:       factEtPol
%:HHH rajout Avec $gB=\AuX$.
\begin{fact}\label{factEtPol} Avec $\gB=\AuX$.
\begin{enumerate}
  \item Un \Bmo  $M=M(\uX)$ est étendu \ssi il est
isomorphe à $M(\uze)$.
  \item En particulier, si $M$ est \pf avec une
\mpn   $G(\uX)\in\gB^{q\times m}$, vu le lemme \ref{lem pres equiv},
$M$ est
étendu depuis $\gA$ \ssi  les matrices
$H(\uX)$ et $H(\uze)$, où $H$ est dessinée ci-dessous, sont \eqves
sur l'anneau $\gB$
$$H(\uX)\;=\kern-7pt\raise5pt\hbox{$\begin{array}{c|p{35pt}|p{23pt}|p{23pt}|p{35pt}|}
\multicolumn{1}{c}{} & \multicolumn{1}{c}{$m$} & \multicolumn{1}{c}{$q$} &
\multicolumn{1}{c}{$q$} & \multicolumn{1}{c}{$m$} \\
\cline{2-5}
\vrule height20pt depth13pt width0pt \rlap{\hskip6cm$q$}\; & \hfil $G(\underline{X})$  \hfil &\hfil
0\hfil &\hfil $0$\hfil &\hfil $0$\hfil \\
\cline{2-5}
\vrule height20pt depth13pt width0pt \rlap{\hskip6cm$q$}\; & \hfil $0$  &\hfil
${\rm I}_{q}$\hfil & \hfil $0$\hfil &\hfil$0$\hfil \\
\cline{2-5}
\end{array}$}
$$
\vskip0pt
\end{enumerate}
\end{fact}

\rem D'après le lemme \ref{lem pres equiv} lorsque les matrices
$H(\uX)$ et $H(\uze)$ sont \eqvesz, elles sont \elrt \eqvesz.
\eoe

\medskip
Concernant les \mptfs on obtient
des \homos de semi-anneaux qui se composent selon l'identité:
$$\GKO\gA\vvvers{\GKO j}\GKO\AuX\vvvvers{\GKO\ev_\ua}\GKO\gA.$$
En conséquence $\GKO j$ est injectif, et la phrase \gui{tout \mptf sur $\AuX$ est étendu depuis $\gA$} signifie que $\GKO j$ est un \isoz, ce que l'on écrit sous forme abrégée \gui{$\GKO\gA=\GKO\AuX$}.

De même pour les anneaux de Grothendieck:

\snic{\KO\gA\vvers{\KO j}\KO\AuX\vvvvers{\KO\ev_\ua}\KO\gA,\quad \hbox{avec}
\quad \KO(\ev_\ua)\circ \KO( j)=\Id_{\KO\gA}.}

\medskip On a par ailleurs les résultats \elrs suivants,
dans lesquels chaque \egt a la signification qu'un morphisme naturel
est un \isoz.

%:       fact A->A[X]
\begin{fact}\label{fact A->A[X]} Avec $\gB=\AuX$.
\begin{enumerate}
\item \label{i1fact A->A[X]}
  $\rD_{\gB}(0)=\DA(0)\gB$ (un \pol est nilpotent \ssi tous ses \coes le sont). En particulier, $\gB\red=\Ared[\uX]$.
\item \label{i2fact A->A[X]}
\label{item3fact A->A[X]} Si $\gA$ est réduit, $\Bst=\Ati$.
  Plus \gnltz, $\Bst=\Ati+\DA(0)\gen{\uX}$.
\item \label{i3fact A->A[X]}
 $\BB(\gA)=\BB(\AuX)$ et $\HO \gA =\HO \AuX $.
\item \label{i4fact A->A[X]}  $\GKO \gA =\GKO \Ared $.
\item \label{i5fact A->A[X]}  $\GKO \gA =\GKO \AuX \iff \GKO \Ared =\GKO \Ared[\uX]$.%
\item \label{i6fact A->A[X]}   $\Pic \gA =\Pic \gB \iff \Pic \Ared =\Pic \Ared[\uX]$.
\end{enumerate}
\end{fact}
\begin{proof}
\emph{\ref{i1fact A->A[X]}} et \emph{\ref{i2fact A->A[X]}.} Voir le lemme \ref{lemGaussJoyal}.\\
\emph{\ref{i3fact A->A[X]}.} On doit montrer que tout \pol \idm
est constant. Cela se fait (en une variable) 
par \recu sur le degré formel du \polz.\\
\emph{\ref{i4fact A->A[X]}.}  C'est le \thrf{propComparRedRed}.\\
\emph{\ref{i5fact A->A[X]}}  et  \emph{\ref{i6fact A->A[X]}.} Résultent des points \emph{\ref{i1fact A->A[X]}} et \emph{\ref{i4fact A->A[X]}.}    
\end{proof}
%

%%%%%%%%%%%%%%%%%%%%%%%%%%%%%%%%%%%%%%%%%%%%%%%%%%%%%%%%%%%%%%%%%%%
%%%%%%%%%%%%%%%%%%%%%%%%%%%%%%%%%%%%%%%%%%%%%%%%%%%%%%%%%%%%%%%%%%%
%%%%%%%%%%%%%%%%%%%%%%                                 %%%%%%%%%%%%%%%%%%%%
%%%%%%%%%%%%%%%%%%%%%%           TRAVERSO              %%%%%%%%%%%%%%%%%%%%
%%%%%%%%%%%%%%%%%%%%%%                                 %%%%%%%%%%%%%%%%%%%%
%%%%%%%%%%%%%%%%%%%%%%%%%%%%%%%%%%%%%%%%%%%%%%%%%%%%%%%%%%%%%%%%%%%
%%%%%%%%%%%%%%%%%%%%%%%%%%%%%%%%%%%%%%%%%%%%%%%%%%%%%%%%%%%%%%%%%%%

%---- Section{Traverso}---------
\section[\Tho de Traverso-Swan]{\Tho de Traverso-Swan, anneaux semi\-normaux}
\label{sec.Traverso} 
%------------------

Cette section est consacrée à l'étude des anneaux $\gA$ pour lesquels
l'\homo naturel de $\Pic\gA$ vers $\Pic\AXr$ est un \iso (i.e., les
\mrcs $1$ sur $\AXr$ sont tous étendus depuis $\gA$).
La réponse est donnée par le \tho de Traverso-Swan-Coquand
(\cite{Tra,Swan80,coq}):

%:     THo{(Traverso-Swan-Coquand)}
\THo{(Traverso-Swan-Coquand)\\}
{
\Propeq
\begin{enumerate}
\item L'anneau $\Ared$ est seminormal (\dfn \ref{defiseminormal})
\item L'\homo naturel $\Pic\gA\vvers{\Pic j}\Pic\AX$ est un \isoz.
\item %Pour tout entier $n\geq1$,
$\Tt r\geq1$, l'\homo naturel $\Pic\gA\to\Pic\AXr$ est un \isoz.
\item %Pour un entier $n\geq1$,
$\Ex r\geq1$, l'\homo naturel $\Pic\gA\to\Pic\AXr$ est un \isoz.
\end{enumerate}
}

\medskip
On montrera  \emph{1  $\Rightarrow$ 3} et \emph{2  $\Rightarrow$ 1}.
Comme corolaire, $\gA$ est seminormal \ssi $\AX$ est seminormal.

%:   subsec{Préliminaires}
\subsec{Préliminaires}

Rappelons tout d'abord le résultat suivant
(voir la proposition \ref{propImProjLib}).

%--- Lemma{lempropImProjLib}---------
\begin{lemma}
\label{lempropImProjLib}
 Une matrice de \prn de rang $1$, $P$, a son image libre \ssi il
existe un vecteur colonne $C$ et un vecteur ligne $L$ tels que $LC=1$
et $CL=P$. En outre, $C$ et $L$  sont uniques, au produit par une
unité près, sous la seule condition que $CL=P$.
\end{lemma}
%--- end-lemma-----------------------------------------

%:HHH reecrit, et fact supprimé
%Par ailleurs, $\Ared[\uX]=(\AuX)\red$ (fait \ref{fact A->A[X]}) et
%le morphisme naturel de~$\Pic\gA$ vers $\Pic\Ared$ est un \iso
%(\thrf{propComparRedRed}).
%D'où le fait suivant.
%
%%--- fact{corpropComparRed}----
%\begin{fact}
%\label{corpropComparRed}
%Le morphisme naturel $\Pic\gA\to \Pic\AuX$ est un \iso \ssi le 
% morphisme naturel  $\Pic\Ared\to \Pic\Ared[\uX]$ est un \isoz.
%\end{fact}
%%--- end-fact------------------------------------

Par ailleurs rappelons que le morphisme naturel $\Pic\gA\to \Pic\AuX$ est un \iso \ssi le 
 morphisme naturel  $\Pic\Ared\to \Pic\Ared[\uX]$ est un \iso
  (fait \ref{fact A->A[X]}~\emph{\ref{i6fact A->A[X]}.}).

Les deux \homos de groupe

\snic{\Pic\gA\vvvers{\Pic j}\Pic\AuX\vvvers{\Pic\ev_\uze}\Pic\gA}

%\sni
se composent selon l'identité. Le premier est injectif, le
second surjectif.
Ce sont des \isos \ssi le premier est surjectif, \ssi
le second est injectif. \perso{la structure de groupe n'intervient
pas dans le raisonnement, cependant  l'injectivité du second
est plus facile à caractériser dans le cas des groupes.}
\\
Cette dernière \prt signifie: toute matrice carrée $P(\uX)$ \idme
de rang 1 sur $\AuX$ qui vérifie \gui{$\Im\big(P(\uze)\big)$ est libre},
vérifie elle-même  \gui{$\Im(P\big(\uX)\big)$ est libre}.
\\
En fait, si  $\Im\big(P(\uze)\big)$ est libre,  la matrice
%diagonale par blocs
$\Diag(P(\uze),0_1)$ est semblable à une matrice de \prn
standard $\I_{1,n}=\Diag(1,0_{n-1,n-1})$ (lemme \dlg \ref{propIsoIm}). D'où le lemme suivant.

%:--- Lemma{lemPicPic1}-----------
\begin{lemma}
\label{lemPicPic1} \Propeq
%-----------------begin enum------------------
\begin{enumerate}
\item L'\homo naturel $\,\Pic\gA\to\Pic\AuX$ est un \isoz,
\item Pour toute
matrice  $M(\uX)=(m_{i,j})
%(\uX))_{i,j\in 1,\ldots ,n}
\in\GAn(\AuX)$
telle que $M(\uze)=\I_{1,n}$,
il existe
$f_1$, \ldots, $f_n$, $g_1$, \ldots, $g_n\in\AuX$ tels que $m_{i,j}=f_ig_j$
pour \hbox{tous $i$, $j$}.
\end{enumerate}
%-----------------end enum------------------
\end{lemma}
%--- end-lemma-----------------------------------------

Notez que l'hypothèse $M(\uze)=\I_{1,n}$ implique $\rg(M)=1$ parce que
l'\homo $\HO(\AuX)\to\HO(\gA)$ est un \isoz.

%--- Convention{convPicPic}--------
\mni{\bf Convention.}
\label{convPicPic}
Nous abrégerons la phrase \gui{le morphisme naturel de~$\Pic\gA$ vers~$\Pic\AuX$ est un \isoz} en écrivant:
 \gui{$\Pic\gA=\Pic\AuX$}.
%--- end-convention------------------------------------

%:     Lemma{lemPicPic2}
\begin{lemma}\label{lemPicPic2}
Soient $\gA\subseteq\gB$ des anneaux réduits
et  $f_1$, \ldots, $f_n$,  $g_1$, \ldots, $g_n$
dans $\BuX$ qui vérifient les \prts suivantes:
$$\leqno{(*)} \quad\quad
\formule{f_1(\uze)=g_1(\uze)=1,  f_i(\uze)=g_i(\uze)=0 \;\;(i=2,\ldots,n),
\\[1mm]
m_{ij}\eqdefi f_ig_j \in\AuX  \;\;  (i,j=1,\ldots,n),
\\[1mm]
\sum_if_ig_i=1.}
$$
Sous ces hypothèses, la matrice $M:=(m_{ij})$ est une \mprn de rang $1$, $M(\uze)=\I_{1,n}$,
et \propeq
\begin{enumerate}
\item Le module  $\Im M$ est libre sur $\AuX$, i.e., étendu depuis $\gA$.
\item Les $f_i$ et les $g_i$ sont dans $\AuX$.
\item Le \pol $f_1$ est dans $\AuX$.
\end{enumerate}
\end{lemma}

\begingroup\def\mou{}
\begin{proof}
\emph{3 $\Rightarrow$ 2.} Les $g_j$ s'obtiennent à partir
de $f_1$ et des $m_{1j}$ en faisant des divisions par puissances croissantes,
car le \coe constant de $f_1$ est égal à $1$. De même, on obtient
ensuite les $f_i$ à partir de $g_1$ et des $m_{i1}$.
L'implication réciproque est triviale.\\
\emph{2 $\Leftrightarrow$ 1.} D'après le lemme \ref{lempropImProjLib}, le \pb est
de trouver des $f_i$ et $g_j$ convenables à partir de la matrice $(m_{ij})$.
Or ces $f_i$ et $g_j$ existent dans $\BuX$, et la condition $f_1(\uze)=1$ force
leur unicité parce que les anneaux sont réduits (donc les \ivs dans l'anneau des \pols sont des constantes).
\end{proof}

Les lemmes \ref{lemPicPic1} et \ref{lemPicPic2} impliquent le résultat suivant.

%:     Corollary{corlemPicPic2}
\begin{corollary}\label{corlemPicPic2}
Soit $\gA\subseteq\gB$ deux anneaux réduits
avec  $\Pic\gB=\Pic\BuX$.  \Propeq
\begin{enumerate}
\item $\Pic\gA=\Pic\AuX$.
\item Si des \pols $f_1$, \ldots, $f_n$,  $g_1$, \ldots, $g_n$
dans $\BuX$  vérifient les conditions~$(*)$ du lemme~\ref{lemPicPic2},
alors les $f_i$ et les $g_i$ sont dans $\AuX$.
\item Si des \pols $f_1$, \ldots, $f_n$,  $g_1$, \ldots, $g_n$
dans $\BuX$  vérifient les conditions~$(*)$, alors
$f_1\in\AuX$.
\end{enumerate}

\end{corollary}
%
% --- Subsubsection{Anneaux seminormaux}--
\subsec{Anneaux seminormaux}

Un anneau intègre $\gA$ est dit \emph{seminormal} si, chaque fois \hbox{que $b^2=c^3\neq 0$}, 
l'\elt $a=b/c$ de $\Frac(\gA)$ est en fait dans $\gA$. 
Dans ce cas, $a^3=b$ \hbox{et $a^2=c$}.

%     Definition{defiseminormal}
\begin{definition}\label{defiseminormal}
Un anneau quelconque $\gA$ est dit \ixc{seminormal}{anneau ---} si chaque fois
que $b^2=c^3$, il existe $a\in\gA$ tel que $a^3=b$ et $a^2=c$.%
\index{anneau!seminormal}
\end{definition}

%--- Fact{factRed1}--------------
\begin{fact}
\label{factRed1}~
1. Un anneau seminormal est réduit.\\
2. Dans un anneau réduit,  $x^2=y^2$ et  $x^3=y^3$ impliquent \hbox{$x=y$}.
\end{fact}
%--- end-fact-----------------------------------------
%
\begin{proof}
\emph{1.} Si $b^2=0$,  alors $b^2=0^3$,
d'où $a\in\gA$ avec $a^3=b$ et $a^2=0$, donc $b=0$.\\
\emph{2.} Dans tout anneau,  $(x-y)^3=4(x^3-y^3)+3(y^2-x^2)(x+y)$.
\end{proof}
En conséquence l'\elt $a$ dans la \dfnz~\ref{defiseminormal}
est toujours unique. En outre,
$\Ann(b)=\Ann(c)=\Ann(a)$.

%--- Fact{fact1seminormal}--------------
\begin{fact}
\label{fact1seminormal}
Tout anneau normal est seminormal.
\end{fact}
%--- end-fact-----------------------------------------
\begin{proof}
Un anneau est normal lorsque tout \idp est \iclz.
Un tel anneau est réduit et \lsdzz: si $uv=0$, il existe
$s$ tel que $su=(1-s)v=0$ (lemme~\ref{lemiclplat}). Soient $b$ et $c$ 
tels \hbox{que $b^3=c^2$}, alors~$c$ est entier sur l'idéal $\gen{b}$, d'où un $x$ tel que
$c=xb$,  
\hbox{d'où $b^3=c^2=x^2b^2$} \hbox{et $b^2(x^2-b)=0$}. Donc il existe $s$ tel que
$s(x^2-b)=0$ et $b^2(1-s)=0$. Ceci donne  $b(1-s)=0$, puis $(sx)^2=s^2b=sb=b$.
En posant $a=sx$, il vient $a^2=b$, $a^3=bsx=bx=c$.
\end{proof}
\endgroup

% --- Subsec{Schanuel}--
\subsubsection*{La condition est nécessaire: l'exemple de Schanuel}

%:     Lemma{lemSchaSeminor}
\begin{lemma}\label{lemSchaSeminor}
Si $\gA$ est réduit et $\Pic\gA=\Pic\AX$, alors
$\gA$ est seminormal.
\end{lemma}

\begin{proof}
Soient $b,c\in\gA$ avec $b^2=c^3$. Soit $\gB=\gA[a]=\gA+a\gA$
un anneau réduit contenant $\gA$, avec $a^3=b, \,a^2=c$.
\\
Considérons les \pols $f_i$ et $g_j$ ($i,j=1,2$) définis comme suit:

\snic{
f_1=1+aX,\; f_2=g_2=cX^2 \; \hbox{et} \;g_1=(1-aX)(1+cX^2).}

%\sni
On a $f_1g_1+f_2g_2=1$, $f_1(0)=g_1(0)=1$, $f_2(0)=g_2(0)=0$. Et
chaque produit $m_{ij}=f_ig_j$ est dans $\AX$.
On applique le  lemme \ref{lemPicPic2}: l'image de la matrice $(m_{ij})$
est libre \ssi $f_1\in\AX$, i.e.  $a\in\gA$.
\end{proof}

\noi NB. Pour $\gB$ on peut prendre $\left(\aqo{\gA[T]}{T^2-c, T^3-b}
\right){\!}\red$. Si un \elt $a$ convenable est déjà présent dans $\gA$, on obtient
par unicité $\gB=\gA$.

%: --- Subsection{Cas des anneaux intègres}--
\subsec{Cas des anneaux intègres}
Nous traitons d'abord les anneaux à pgcd,
puis  les anneaux normaux et enfin les anneaux seminormaux.
% --- Subsection{anneau a pgcd}--
\subsubsection*{Cas d'un anneau à pgcd}

Rappelons qu'un anneau (intègre) à pgcd est un anneau dans lequel
deux \elts arbitraires admettent un plus grand commun diviseur,
\cad une borne supérieure pour la relation de divisibilité.
Rappelons aussi que si~$\gA$ est un anneau à pgcd, il en va de
même pour l'anneau des \polsz~$\AuX$.

%:--- Lemma{lemPicGcd}--------------
\begin{lemma}
\label{lemPicGcd}
Si $\gA$ est un anneau intègre à pgcd, alors $\Pic\gA=\so{1}$.
\end{lemma}
%--- end-lemma-----------------------------------------

%-----------------begin proof------------------
\begin{proof}
On utilise la \carn donnée dans le lemme \ref{lempropImProjLib}.\\
Soit $P=(m_{ij})$ une matrice \idme de rang 1.
Puisque $\sum_i m_{ii}=1$, on peut supposer que $m_{1,1}$ est
régulier.
Soit $f$ le pgcd des \elts de la première ligne. On écrit $m_{1j}=fg_j$
avec le pgcd des $g_j$ égal à 1.  L'\egtz~\hbox{$m_{1,1}m_{ij}=m_{1j}m_{i1}$}  donne, en simplifiant par $f$, $g_1m_{ij}=m_{i1}g_j$.
Ainsi,~$g_1$ divise tous les $m_{i1}g_j$, donc aussi leur pgcd
$m_{i1}$.
On écrit $m_{i1}=g_1f_i$. Puisque~\hbox{$g_1f_1=m_{1,1}=fg_1$}, cela donne
$f_1=f$. Enfin,  $m_{1,1}m_{ij}=m_{1j}m_{i1}$ donne
l'\egtz~\hbox{$f_1g_1m_{ij}=f_1g_jg_1f_i$}, puis $m_{ij}=f_ig_j.$
\end{proof}
%-----------------end proof------------------
 On a alors le corolaire suivant.

%:     proposition{thcorlemPicGcd}
\begin{proposition}\label{thcorlemPicGcd}
Si $\gA$ est un corps discret ou un anneau \zed réduit,
alors $\Pic\gA=\Pic\AuX=\so{1}$.
\end{proposition}
\begin{proof}
Le lemme \ref{lemPicGcd} donne le résultat pour les corps discrets.
Il suffit ensuite d'appliquer la machinerie \lgbe \elr \num2~\paref{MethodeZedRed}.\imlgz
\end{proof}

% --- Subsec{anneau intègre normal}--
\penalty-2500
\subsubsection*{Cas d'un anneau intègre normal}

%:--- Lemma{lemIntegclos}-------------
\begin{lemma}
\label{lemIntegclos}
Si $\gA$ est intègre normal, alors
$\Pic\gA=\Pic\AuX$.
\end{lemma}
%--- end-lemma-----------------------------------------
%-----------------begin proof------------------
\begin{proof}
On utilise la \carn donnée au corolaire \ref{corlemPicPic2}~\emph{3},
avec ici $\gA\subseteq\gK$, le corps de fractions de $\gA$. Soient
$f_i$ et $g_j$, $(i,j\in\lrbn)$ les \pols convenables de $\KuX$.
Alors, puisque $f_1g_1=m_{1,1}\in\AuX$ et  $g_1(\uze)=1$, vu le \tho de
\KRO \ref{thKro},
les \coes de $f_1$ sont entiers sur l'anneau engendré
par les \coes de~$m_{1,1}$. Ainsi $f_1 \in \gA[X]$.
\end{proof}
%-----------------end proof------------------

\rem De même que pour la proposition \ref{thcorlemPicGcd}, on peut étendre
le résultat du lemme \ref{lemIntegclos} au cas d'un anneau réduit
$\gA$ intégralement fermé dans un anneau \zed réduit $\gK\supseteq\gA$.
\eoe

% --- Subsec{intègre seminormal}--
\subsubsection*{Cas d'un anneau intègre seminormal}

%:--- proposition{propIntSemin}-------
\begin{proposition}
\label{propIntSemin}
Si $\gA$ est intègre et seminormal, alors $\Pic\gA=\Pic\AuX$.
\end{proposition}
%--- end-proposition----------------------------------------
%-----------------begin proof------------------
\begin{Proof}{Début de la \demz. }
Comme dans la \dem du lemme~\ref{lemIntegclos},
on a au départ des \pols $f_1(\uX)$, \ldots, $f_n(\uX)$,
$g_1(\uX)$, \ldots, $g_n(\uX)$ dans~$\KuX$ 
qui vérifient les conditions $(*)$ du lemme \ref{lemPicPic2}.
 On appelle $\gB$ le sous-anneau de $\gK$ engendré
par $\gA$ et par les \coes des $f_i$ et des $g_j$, ou encore,
cela revient au même,  engendré
par $\gA$ et par les \coes de $f_1$. Alors,
vu le \tho de \KROz, $\gB$ est une extension finie de
$\gA$. Notre but est de montrer que
$\gA=\gB$.
On note  $\fa$ le conducteur de $\gA$ dans $\gB$, \cad
l'ensemble~\hbox{$\sotq{x\in\gB}{x\gB\subseteq\gA}$}. C'est à la fois un
\id de $\gA$ et $\gB$. Notre but est maintenant de montrer
$\fa=\gen{1}$, \cade que $\gC=\gA\sur{\fa}$ est trivial.
\end{Proof}

Nous commençons par deux lemmes.

%:--- Lemma{lemIntSemin1}-------------
\begin{lemma}
\label{lemIntSemin1}
Si $\gA\subseteq\gB$, $\gA$ seminormal et $\gB$ réduit, alors le
conducteur~$\fa$ de $\gA$ dans $\gB$ est un \id radical de $\gB$.
\end{lemma}
%--- end-lemma-----------------------------------------
%-----------------begin proof------------------
\begin{proof}
On doit montrer  que si $u\in\gB$ et $u^2\in\fa$, alors $u\in\fa$. Soit
donc $c\in\gB$, on doit montrer que  $uc\in\gA$. On sait que
$u^2c^2$ et
 $u^3c^3=u^2(uc^3)$ sont dans $\gA$ puisque $u^2\in\fa$.
Puisque $(u^3c^3)^2=(u^2c^2)^3$, on a un $a\in\gA$ tel que~\hbox{$a^2=(uc)^2$} et $a^3=(uc)^3$. Comme $\gB$ est réduit, on obtient
$a=uc$, et donc $uc\in\gA$.
\end{proof}
%-----------------end proof------------------

\rem La \emph{clôture seminormale} d'un anneau
$\gA$ dans un  anneau réduit~$\gB\supseteq \gA$ est obtenue en partant de $\gA$ et en rajoutant
les \elts $x$ de $\gB$  tels que $x^2$ et $x^3$ sont dans l'anneau
préalablement construit. Notez que par le fait \ref{factRed1}, $x$ est uniquement déterminé par la donnée de~$x^2$ et~$x^3$.
La preuve du lemme précédent peut alors être interprétée comme
une démonstration de la variante suivante.%   
\index{seminormale!clôture --- dans un suranneau réduit}
\eoe

%:--- Lemma{lemIntSemin1bis}-------------
\begin{lemma}
\label{lemIntSemin1bis}
Soient $\gA\subseteq\gB$ réduit, $\gA_{1}$ la
clôture  seminormale de  $\gA$ dans~$\gB$, et
$\fa$ le conducteur de $\gA_{1}$ dans $\gB$.
Alors, $\fa$ est un \id radical de~$\gB$.
\end{lemma}
%--- end-lemma-----------------------------------------

%:--- Lemma{lemIntSemin2}-------------
\begin{lemma}
\label{lemIntSemin2}
Soient $\gA\subseteq\gB$,  $\gB=\gA[\cq]$ réduit fini sur $\gA$ et
$\fa$ le conducteur de $\gA$ dans $\gB$. On suppose que $\fa$ est un
\id radical, alors il est
égal à $\sotq{x\in\gA}{xc_1,\ldots ,xc_q\in\gA}$.
\end{lemma}
%--- end-lemma-----------------------------------------
%-----------------begin proof------------------
\begin{proof}
En effet, si $xc_i\in\gA$, alors  $x^\ell c_i^\ell\in\gA$ pour tout
$\ell$, et donc pour un $N$ assez grand
$x^N y\in\gA$ pour tout $y\in\gB$, donc $x$ est dans le nilradical de
$\fa$ (si $d$ majore les degrés des équations de dépendance
intégrale des $c_i$ sur $\gA$, on pourra prendre $N=(d-1)q$).
\end{proof}

%: fin de la demo du \tho \ref{propIntSemin}
\begin{Proof}{Fin de la \dem de la proposition \ref{propIntSemin}. }\\
Nous la donnons
d'abord  en \clamaz.
Le raisonnement classique naturel procéderait par l'absurde:
 l'anneau $\gC$ est trivial parce que sinon, il posséderait un \idemi
 et la \lon en cet \idemi mènerait à une contradiction.\\
Pour éviter le \crc non \cof du raisonnement par l'absurde,
nous allons localiser en un \fimaz, en rappelant notre \dfn \gui{sans négation}
selon laquelle un filtre
est  maximal \ssi l'anneau localisé est un \alo \zedz.
Autrement dit nous tolérons pour les  \fimas d'un anneau, non seulement
les \cops des \idemis  mais aussi le filtre engendré par $0$
qui donne par \lon l'anneau  trivial.
En \clama un anneau  est alors trivial
\ssi  son seul \fima est l'anneau tout entier (autrement dit, le filtre
engendré par $0$).
\\ 
Insistons sur le fait que c'est seulement dans l'affirmation précédente
que se situe le \crc \gui{classique} du raisonnement.
Car la preuve de ce qui suit est parfaitement \covz: si $S$ est un \fima de $\gC$,
alors~\hbox{$0\in S$} (donc $S=\gC$).
\\ 
On considère l'inclusion $\gC=\gA\sur{\fa}\subseteq \gB\sur{\fa}=\gC'$.
Soit $S$ un \fima de~$\gC$, et $S_1$ le \fima correspondant de $\gA$
(l'image réciproque de $S$ par la \prn canonique).
Puisque $S$ est un \fimaz, et puisque~$\gC$ est réduit,
$S^{-1}\gC=\gL$ est un \alo \zed réduit, \cad un corps discret,
contenu dans l'anneau réduit  $S^{-1}\gC'=\gL'$.\\
Si  $x$ est un objet défini sur $\gB$, notons $\ov{x}$ ce qu'il
devient après le changement de base  $\gB\to\gL'$.
 Puisque $\gL$ est un \cdiz, $\gL[\uX]$
 est un anneau intègre à pgcd, et les $\ov{f_i}$ et $\ov{g_j}$
 sont dans $\gL[\uX]$.
 Cela signifie qu'il existe $s\in S_1$  tel que  $sf_1\in\AX$.
 D'après le lemme \ref{lemIntSemin2}, ceci
implique que $s\in\fa$. Ainsi $\ov s=0$  et $\ov s\in S$.
\end{Proof}
%-----------------end proof------------------

La démonstration donnée ci-dessus pour la proposition  \ref{propIntSemin}
est finalement assez simple.
Elle n'est cependant pas totalement
\cov et elle semble ne  traiter que le cas intègre.

\begin{Proof}{Démonstration \cov de la proposition \ref{propIntSemin}. }\\
Nous reprenons la \dem donnée en \clama en consi\-dérant que
le \fima $S$ de $\gC$ est un objet purement générique qui
nous guide dans la \prcoz.
\\ 
Imaginons que l'anneau $\gC$ soit un corps discret, \cad que l'on  ait déjà
fait la \lon en un \fimaz. 
\\ 
Alors, des \pols ${F_i}$ et ${G_j}$ de $\gC[\uX]$
vérifiant ${F_i}{G_j}=\ov{m_{ij}}$
et $F_1(\uze)=1$
sont calculés  à partir des $\ov{m_{ij}}$
selon un \algo que l'on déduit des \prcos données
auparavant pour le cas des corps discrets (lemme~\ref{lemPicGcd}).
L'unicité de la solution force alors l'\egt $F_1=\ov{f_1}$, ce qui montre que~$\ov{f_1}\in\gC[\uX]$, et donc que $\gC$ est trivial.
\\ 
Cet \algo utilise la disjonction \gui{$a$ est nul ou $a$ est
inversible},
pour les \elts $a\in\gC$ qui sont produits par l'\algo à partir des
\coes des \pols $\ov{m_{i,j}}$.  Comme $\gC$ est seulement un anneau réduit,
sans test d'\egt à $0$ ni test d'inversibilité, l'\algo
pour les corps discrets, si on l'exécute avec $\gC$, doit être remplacé
par un arbre dans lequel on ouvre deux branches chaque fois
qu'une question \gui{$a$ est-il nul ou inversible?} est posée par
l'\algoz.
\\ 
Nous voici en face d'un arbre, gigantesque, mais fini.
Disons que systémati\-quement on a mis la branche \gui{$a$ inversible}
à gauche, et la branche \gui{$a=0$} à droite.
Regardons ce qui se passe dans la branche d'extrême gauche. 
\\ 
On a inversé successivement $a_1$, \ldots , $a_p$
et l'on a obtenu un $s$ qui montre  que l'anneau
$\gC[1/(a_1\cdots a_p)]$ est trivial.
\\
 \emph{Conclusion: dans l'anneau $\gC,$ on a l'\egt $a_1\cdots
a_p=0$.}
\\
 Remontons d'un cran. 
\\ 
Dans l'anneau $\gC[1/(a_1\cdots a_{p-1})]$, nous savons que $a_p=0$.
\\ 
La branche de gauche n'aurait pas d\^u être ouverte.
Regardons le calcul dans la branche  $a_p=0$.
\\ 
Suivons à partir de là
la branche d'extrême gauche.
\\ 
On a inversé  $a_1$, \ldots , $a_{p-1}$, puis, disons $b_1,\ldots
,b_k$ (éventuellement, $k=0$).
Nous obtenons un $s$ qui montre  que  l'anneau
$\gC[1/(a_1\cdots a_{p-1} b_1\cdots  b_k)]$ est trivial.
\\
 \emph{Conclusion: dans l'anneau $\gC,$ on a l'\egt
$a_1\cdots a_{p-1} b_1\cdots  b_k=0$.}
\\ 
  Remontons d'un cran: nous savons que $b_k=0$ (ou, si $k=0$, $a_{p-1}=0$)
  dans l'anneau qui était là juste avant le dernier branchement:
à savoir l'anneau  $\gC[1/(a_1\cdots a_{p-1} b_1\cdots  b_{k-1})]$  (ou, si $k=0$, $\gC[1/(a_1\cdots a_{p-2})]$).
  La branche de
gauche n'aurait pas d\^u être ouverte.
Regardons le calcul dans la branche~\hbox{$b_k=0$} (ou, si $k=0$, la branche $a_{p-1}=0$)\,\ldots
\\
 \emph{Et ainsi de suite.} Quand on poursuit le processus
jusqu'au bout,
on se retrouve à la racine de l'arbre avec l'anneau
$\gC=\gC[1/1]$ qui est trivial.
\end{Proof}
%-----------------end proof------------------

En utilisant le lemme \ref{lemIntSemin1bis} à la place du lemme
\ref{lemIntSemin1} on obtiendra le résultat suivant, plus précis
que la proposition~\ref{propIntSemin}.

%:--- proposition{propIntSeminBis}-------
\begin{proposition}
\label{propIntSeminBis}
Si $\gA$ est un anneau intègre et $P$ un module
\pro de rang $1$ sur $\AuX$ tel que $P(\uze)$ est libre,
il existe $c_{1}$, \ldots, $c_{m}$ dans le corps de fractions de $\gA$ tels que:
\begin{enumerate}
\item $c_{i}^2$ et $c_{i}^3$ sont dans $\gA[(c_{j})_{j<i}]$ pour $i=1,\ldots,m$,
\item $P$ est libre sur $\gA[(c_{j})_{j\leq m}][X]$.
\end{enumerate}
\end{proposition}
%--- end-prop----------------------------------------
\rem En fait, seul intervient le corps de fractions du sous-anneau engendré
par les \coes présents dans une matrice de \prn 
dont l'image est isomorphe à $P$.
\eoe

%:----   subsec{Cas général} %%%%%%%%%%%%
\subsec{Cas général}
\label{subsecTravSwanGeneral}

%:--- Proposition{thseminormalCoq}----------
\begin{proposition}
\label{thseminormalCoq}\emph{(Coquand)}
Soit $\gA\subseteq\gK$ avec $\gK$ réduit.
\begin{enumerate}
\item
\'Etant donnés  $f$ et $g\in\KuX^{n}$ qui vérifient les conditions $(*)$ du lemme~\ref{lemPicPic2},
on peut construire  $c_{1}$, \ldots, $c_{m}$ dans  $\gK$ tels que:
\begin{enumerate}
\item [--]  $c_{i}^2$ et $c_{i}^3$ sont dans $\gA[(c_{j})_{j<i}]$ pour $i\in\lrbm$,
\item [--]  $f$ et $g$ ont leurs \coos dans $\gA[(c_{k})_{k\in\lrbm}][\uX]$
\end{enumerate}
\item
Si   $\Pic\gK=\Pic\KuX$ et si $P$ est un module
\pro de rang $1$ sur $\AuX$, il existe $c_{1}$, \ldots, $c_{m}$ dans  $\gK$ tels que:
\begin{enumerate}
\item [--]  $c_{i}^2$ et $c_{i}^3$ sont dans $\gA[(c_{j})_{j<i}]$ pour $i\in\lrbm$,
\item [--]  $P\simeq P(\uze)$  sur $\gA[(c_{k})_{k\in\lrbm}][\uX]$.
\end{enumerate}
\end{enumerate}
\end{proposition}
%-----------------begin proof------------------
\begin{proof}
La \dem de la proposition \ref{propIntSemin}, ou de sa variante plus précise~\ref{propIntSeminBis},
est en fait une \dem du point \emph{1} ci-dessus.
Le point \emph{2} s'en déduit facilement.
\end{proof}
%-----------------end proof-----------

%:     Theorem{thTSC}
\begin{theorem}\label{thTSC}\emph{(Traverso-Swan-Coquand)}\\
Si $\gA$ est un anneau seminormal, alors $\Pic\gA=\Pic\AuX$.
\end{theorem}
\begin{proof}
On le déduit de la proposition précédente en utilisant
le fait qu'il existe un suranneau $\gK$ de $\gA$ tel que $\Pic\gK=\Pic\KuX$.
En effet, tout anneau réduit est contenu dans un anneau \zed réduit
 (\thrf{thZedGen} ou~\ref{thAmin})
$\gK$, lequel vérifie $\Pic\gK=\Pic\KuX=\so1$ (proposition~\ref{thcorlemPicGcd}).
\end{proof}
%

%%%%%%%%%%%%%
\subsubsection*{Un calcul direct menant au résultat}

Comme souvent lorsque l'on  essaie d'implémenter sur machine
un \tho \cof qui a une preuve élégante, on est amené à trouver certains
raccourcis dans les calculs qui donnent en définitive une solution
plus simple. Mais cette solution cache en partie,
sinon le mécanisme profond de la preuve initiale, du moins
la démarche de la pensée qui a élaboré la preuve.
Voir par exemple comment l'exercice \ref{exo7.1} trivialise la preuve
du \tho de structure locale des \mptfsz.

C'est ce qui s'est produit avec la proposition \ref{thseminormalCoq}
qui a finalement été réalisée par un \algo de nature assez \elr 
dans \cite[Barhoumi\&Lombardi]{BL07}, basé sur la théorie de l'\id résultant 
(cf. section \ref{subsecIdealResultant}) et des modules
sous-résultants.

%%%%%%%%%%%%%%%%%%%%%%%%%%%%%%%%%%%%%%%%%%%%%%%%%%%%%%%%%%%%%%%%%%%%
%%%%%%%%%%%%%%%%%%%%%%%%%%%%%%%%%%%%%%%%%%%%%%%%%%%%%%%%%%%%%%%%%%%%
%%%%%%%%%%%%%%%%%%%%%%                                 %%%%%%%%%%%%%%%%%%%
%%%%%%%%%%%%%%%%%%%%%%       QUILLEN PATCHING          %%%%%%%%%%%%%%%%%%%
%%%%%%%%%%%%%%%%%%%%%%                                 %%%%%%%%%%%%%%%%%%%
%%%%%%%%%%%%%%%%%%%%%%%%%%%%%%%%%%%%%%%%%%%%%%%%%%%%%%%%%%%%%%%%%%%%
%%%%%%%%%%%%%%%%%%%%%%%%%%%%%%%%%%%%%%%%%%%%%%%%%%%%%%%%%%%%%%%%%%%%

%---- Section{Recollement à la Quillen-Vaserstein}---
\section{Recollement à la Quillen-Vaserstein}
\label{subsecQPatch}

Nous exposons dans cette section ce qu'en anglais on appelle  le Quillen patching.
C'est un résultat profond qui pourrait sembler a priori
un peu trop abstrait (utilisation abusive d'\idemasz) mais
qui s'avère plein de bon sens \cofz.

Les démonstrations que nous donnons sont (pour l'essentiel)
 recopiées de~\cite{Kun}.
Nous avons remplacé
la \lon en n'importe quel idéal maximal par
la \lon en des \mocoz.

%--- Lemma{lem0PrepVaser}-------------
\begin{lemma}
\label{lem0PrepVaser}
Soit $S$ un \mo de l'anneau $\gA$ et $P\in\AX$ un \pol  tel que
$P=_{\gA_S[X]}0$ et
$P(0)=0$. Alors, il existe $s\in S$  tel que
$P(sX)=0$.
\end{lemma}
%--- end-lemma-----------------------------------------
%-----------------begin proof------------------
\facile
%-----------------end proof------------------

Voici une légère variante.

%--- Fact{lem01PrepVaser}-------------
\begin{fact}
\label{lem01PrepVaser}
Soit $S$ un \mo de l'anneau $\gA$ et $P\in\gA_S[X]$ un \pol  tel que
  $P(0)=0$. Alors, il existe $s\in S$  et $Q\in\AX$ tels \hbox{que
$P(sX)=_{\gA_S[X]}Q$}.
\end{fact}
%--- end-fact-----------------------------------------
%
%\facile
%

%--- Lemma{lem1PrepVaser}-------------
\begin{lemma}
\label{lem1PrepVaser}
Soit $S$ un \mo de l'anneau $\gA$. On considère trois matrices à \coes
dans $\AX$, $A_1,\,A_2,\,A_3$ telles que le produit $A_1\,A_2$ a le
même format que $A_3$. Si $A_1\,A_2=_{\gA_S[X]}A_3$ et
$A_1(0)\,A_2(0)=A_3(0)$, il existe $s\in S$  tel que
$A_1(sX)\,A_2(sX)=A_3(sX)$.
\end{lemma}
%--- end-lemma-----------------------------------------
%-----------------begin proof------------------
\begin{proof}
Appliquer le lemme \ref{lem0PrepVaser} aux coefficients de la matrice  $A_1\,A_2-A_3$.
\end{proof}
%-----------------end proof------------------

%--- Lemma{lem2PrepVaser}-------------
\begin{lemma}
\label{lem2PrepVaser}
Soit $S$ un \mo de l'anneau $\gA$ et 
$C(X)\in\GL_p(\gA_S[X])$.  Il existe $s\in S$  et $U(X,Y) \in
\GL_p(\gA[X,Y])$ tels que $U(X,0)=\I_p$, et, sur l'anneau $\gA_S[X,Y]$,
$U(X,Y)=C(X+sY)C(X)^{-1}$.
\end{lemma}
%--- end-lemma-----------------------------------------
%-----------------begin proof------------------
\begin{proof}
Posons $E(X,Y)=C(X+Y)C(X)^{-1}$. Notons $F(X,Y)=E(X,Y)^{-1}$. On a
$E(X,0)=\I_p$, donc $E(X,Y)=\I_p+E_1(X)Y+\cdots+E_k(X)Y^k$. Pour un $s_1\in
S$, les ${s_1}^jE_j$ peuvent se réécrire \gui{sans dénominateur}.
On obtient ainsi une matrice $E'(X,Y)\in \MM_p(\gA[X,Y])$ telle que
$E'(X,0)=\I_p$ et, sur  $\gA_S[X,Y]$,
$E'(X,Y)=E(X,s_1Y)$. On procède de même avec $F$ (et l'on peut choisir un $s_1$ commun).
On a alors $E'(X,Y)F'(X,Y)=\I_p$ \linebreak 
dans   $\MM_p(\gA_S[X,Y])$
et $E'(X,0)F'(X,0)=\I_p$. \\
En appliquant le lemme \ref{lem1PrepVaser}
dans lequel on remplace $X$ par $Y$ et $\gA$ par~$\AX$, on obtient
un $s_2\in S$ tel  
que  $E'(X,s_2Y)F'(X,s_2Y)=\I_p$. \\
D'où le résultat
souhaité avec $U=E'(X,s_2Y)$ et $s=s_1s_2$.
\end{proof}
%-----------------end proof------------------

%--- Lemma{lem3PrepVaser}-------------
\begin{lemma}
\label{lem3PrepVaser}
Soit $S$ un \mo de  $\gA$ et $G\in \AX^{q\times m}$.
Si $G(X)$ et $G(0)$ sont \eqves sur $\gA_S[X]$, il existe
$s\in S$ tel que $G(X+sY)$ et $G(X)$ sont \eqves sur $\gA[X,Y]$.
\end{lemma}
%--- end-lemma-----------------------------------------
%-----------------begin proof------------------
\begin{proof}
\'{E}crivons $G=C\,G(0)\,D$ avec $C\in\GL_q(\gA_S[X])$ et
 $D\in\GL_m(\gA_S[X])$. On a donc
\[\preskip.4em \postskip.2em \arraycolsep2pt
\begin{array}{rcl} 
 G(X+Y) &  = &  C(X+Y)G(0)D(X+Y) \\[1mm] 
  &  = &   C(X+Y)C(X)^{-1}G(X)D(X)^{-1}D(X+Y).
 \end{array}
\]
En appliquant le lemme \ref{lem2PrepVaser}, on obtient $s_1\in S$,
$U(X,Y)\in\GL_q(\gA[X,Y])$ et $V(X,Y)\in\GL_m(\gA[X,Y])$,
tels que 
\[\preskip.4em \postskip.2em 
\begin{array}{rcll} 
U(X,0)=\I_q  &,   & V(X,0)=\I_m \;,  \\[1mm] 
  \hbox{et sur l'anneau } \gA_S[X,Y]&  :   \\[1mm] 
 U(X,Y)=C(X+s_1Y)C(X)^{-1} & \hbox{et}  &  V(X,Y)=D(X)^{-1}D(X+s_1Y).
 \end{array}
\] 
Donc
\[\preskip-.4em \postskip.4em 
\begin{array}{ccc} 
  G(X)    =    U(X,0)G(X)V(X,0),   \\[1mm] 
  \hbox{et sur l'anneau } \gA_S[X,Y]\quad :  \hspace{4cm}~ \\[1mm] 
  G(X+s_1Y)   =   U(X,Y)G(X)V(X,Y).  
 \end{array}
\]
 En appliquant le lemme \ref{lem1PrepVaser}
(comme dans le lemme \ref{lem2PrepVaser}), on obtient $s_2\in S$
tels que $G(X+s_1s_2Y)=U(X,s_2Y)G(X)V(X,s_2Y)$.\\
D'où le résultat avec $s=s_1s_2$.
\end{proof}
%-----------------end proof------------------

%:     plcc{thPatchV}-------------
\begin{plcc}
\label{thPatchV} \emph{(Recollement de Vaserstein)}\\
Soit $G$ une matrice  sur $\AX$ et $S_1$, $\ldots$, $S_n$ des \moco de~$\gA$.
\begin{enumerate}
\item Les matrices $G(X)$ et $G(0)$ sont \eqves sur $\AX$ \ssi  elles sont
\eqves sur $\gA_{S_i}[X]$ pour chaque $i$.
\item Même résultat pour \gui{l'\eqvc à gauche}: deux matrices
$M$ et $N$ de même format sur un anneau commutatif sont dites \emph{\eqves à gauche}
s'il existe une matrice carrée \iv $H$ telle que $H\,M=N$.
\end{enumerate}
\index{matrices!equiva@équivalentes à gauche}
\index{equivalentes@équivalentes à gauche!matrices}
\end{plcc}
%--- end-plcc-----------------------------------------
%-----------------begin proof------------------
\begin{proof}
\emph{1.} On vérifie %sans peine 
que l'ensemble des $s\in\gA$
tels que la matrice $G(X+sY)$
soit \eqve à $G(X)$ sur $\gA[X,Y]$ forme un
idéal de $\gA$. En appliquant le lemme \ref{lem3PrepVaser},
cet idéal contient un \elt $s_i$ dans $S_i$ pour chaque $i$, donc il
contient 1,
et  $G(X+Y)$ est \eqve à $G(X)$. Il reste à faire $X=0$.\\
\emph{2.} Dans toutes les \dems précédentes, on peut remplacer
l'\eqvc par l'\eqvc à gauche.
\end{proof}
%-----------------end proof------------------

%:     plcc{thPatchQ}--------
\begin{plcc}
\label{thPatchQ} \emph{(Recollement de Quillen)}\index{Quillen!recollement de ---}\\
Soit $M$ un \mpf sur $\AX$ et $S_1$, $\ldots$, $S_n$ des \moco de $\gA$.
Alors, $M$ est un module étendu depuis $\gA$ \ssi chaque  $M_{S_i}$
est étendu depuis~$\gA_{S_i}$.
\end{plcc}
%--- end-plcc-----------------------------------------
%-----------------begin proof------------------
\begin{proof} C'est un corolaire du théorème précédent
car l'\iso entre les modules $M(X)$ et $M(0)$ s'exprime
par l'\eqvc de deux matrices $H(X)$ et  $H(0)$
construites à partir d'une \mpn $G$ de $M$
(voir le fait~\ref{factEtPol}).
\end{proof}
%-----------------end proof------------------

\comm
La formulation originale de Quillen\index{Quillen!recollement de ---}, \eqve au
\plgrf{thPatchQ} en \clamaz,
est la suivante: \emph{si $M_\fm$ est étendu depuis $\gA_\fm$
après \lon en tout
\idema $\fm$, 
%:HHH un peu mieux ci dessous
% alors il est étendu.
alors $M$ est étendu depuis $\gA$.
}
\\
Pour faire façon d'une \dem classique basée sur le \rcm
de Quillen dans la formulation originale,
nous devrons faire appel à la machinerie \lgbe de base
expliquée dans la section~\ref{secMachLoGlo}.\imlbz
\eoe

%%%%%%%%%%%%%%%%%%%%%%%%%%%%%%%%%%%%%%%
%---- Subsection{secRoitman}---
%%%%%%%%%%%%%

\subsec{Un \tho de Roitman}
\label{secRoitman}
Ce paragraphe est consacré à la \dem du \tho suivant,
qui constitue une sorte de réciproque du \tho de \rcm de Quillen.

%:HHH j'ai numerote ce theoreme
%:     Theorem{thRoitman}
\begin{theorem}\label{thRoitman} \emph{(\Tho de Roitman)}\\
Soit $r$ un entier $\geq1$ et $\AuX=\AXr$. Si tout $\AuX$-module \ptf est étendu depuis $\gA$, alors toute \lon $\gA_S$ de $\gA$ vérifie la même \prtz. 
\end{theorem}
%--------- fin theorem ---------------------------------------------- 

\subsubsection*{En une variable}
%------------------
%:     Lemma{lemRoitman}
\begin{lemma}\label{lemRoitman}
Si tout $\AX$-module \ptf est étendu depuis~$\gA$, alors toute \lon $\gA_S$ de $\gA$ vérifie la même \prtz.
\end{lemma}
\begin{proof}
\emph{Cas spécial: 
$\gA_S$ est un \alo \dcdz.} 
\\
Notons  $\rho:\AX\to\gA_S[X]$ le morphisme naturel.
Soit $M \in\GAn(\gA_S[X])$.
Puisque $\gA_S$ est local,
$M(0)$ est semblable à un \prr standard $\I_{k,n}$.
Nous pouvons donc supposer \spdg que $M(0)=\I_{k,n}$,
i.e. que $M(X)=\I_{k,n}+M'(X)$ avec
$M'(X) \in\Mn(\gA_S[X])$ et $M'(0) = 0$. 
\\
Soit~$v$ le \gui{produit des
\denosz} dans les \coes des entrées de~$M'(X)$.
Puisque $M'(0) = 0$, on a une matrice $N' = N'(X)\in\Mn(\AX)$ telle  
que $M'(vX)=_{\gA_S[X]}N'(X)^{\rho}$ et $N'(0)=0$.
\\
Avec $N(X)=\I_{k,n}+N'(X)$  
on obtient~\hbox{$N(0)=\I_{k,n}$} et $M(vX)=N(X)^{\rho}$.
Puisque $M^2=_{\gA_S[X]}M$, on a \hbox{un $s\in S$} tel que
$s(N^2-N)=0$.
\\
Comme $(N^2-N)(0)=0$, on écrit $N^2-N=XQ(X)$.\\
Maintenant $sXQ(X)=0$ implique $sQ(X)=0$.
A fortiori $sQ(sX)=0$, \linebreak 
donc $N(sX)^2=N(sX).$
Mais les \mptfs sur  $\AX$ sont étendus depuis $\gA$, donc la \mprn $N(sX)$
a un noyau et une image isomorphes au noyau et à l'image de
$N(0)=\I_{k,n}$. Donc $N(sX)$ est semblable à $\I_{k,n}$: il existe
$G=G(X)\in\GL_n(\AX)$ telle que
$$
G^{-1}(X)N(sX)G(X)=\I_{k,n}.
$$
En posant $H(X)=G(X)^{\rho}\in\GL_n(\gA_S[X])$, on obtient sur $\gA_S[X]$
l'\egt  \[H^{-1}(X)M (svX)H(X)=\I_{k,n}\]
et par conséquent
$$
H^{-1}(X/sv)M(X)H(X/sv)=_{\gA_S[X]}\I_{k,n}
$$
avec $H(X/sv)\in\GL_n(\gA_S[X])$.

 \emph{Cas \gnlz.} Soit $P$ un $\gA_S[X]$-\mptf arbitraire.\\ 
Posons $\gB=\gA_S$.
Comme d'habitude, $P(0)$ note le \Bmo  obtenu par \eds via
le morphisme $\ev_0:\BX\to\gB$.
On applique la machinerie \lgbe de base \paref{MethodeIdeps}
à la \dem \cov que nous venons de donner dans le cas spécial.
On obtient des \moco $V_1$, \ldots, $V_m$ de $\gB$ avec
$P\simeq_{\gB_{V_i}}P(0)$. On conclut avec le \rcm de Quillen: 
$P\simeq_{\gB}P(0)$.\imlbz
\end{proof}

\medskip \REM{. Pour implémenter l'\algo correspondant à
cette \demz} En fait la seule \prt particulière que nous avons utilisée dans
la \dem du cas spécial, c'est que les $\gA_S$-\mptfs sont libres.
Donc la mise {\oe}uvre de la machinerie \lgbe de base est ici très \elrz.
Elle consiste à construire
des \lons \come pour lesquelles 
la matrice $M(0)$ devient semblable à une \mprn standard, 
et à faire
tourner l'\algo donné par la \dem du cas spécial dans chacune de ces \lonsz.
Naturellement, on termine avec
l'\algo correspondant à la \dem \cov du \rcm de Quillen.\imlbz
\eoe
\subsubsection*{En plusieurs variables}
%------------------
\begin{Proof}{\Demo du \tho de Roitman \ref{thRoitman}. }
On raisonne par \recu sur $r$. Le cas $r=1$ est déjà traité.
Passons de $r\geq1$ à $r+1$. 
On considère un \mo $S$ d'un anneau $\gA$. Notons $(\Xr)=(\uX)$.\\  
On a $\gA_S[\uX,Y]=\gA[\uX,Y]_S=(\gA[Y]_S)[\uX]$.
Soit $P$ un $\gA_S[\uX,Y]$-\mptfz. D'après l'\hdr appliquée avec l'anneau $\gA[Y]$, $P$ est étendu depuis $\gA[Y]_S=\gA_S[Y]$, i.e.,
$P(\uX,Y)$ est isomorphe à $P(\uze,Y)$ comme $\gA_S[\uX,Y]$-module. Et d'après
le cas $r=1$ appliqué avec l'anneau $\gA$, $P(\uze,Y)$ est étendu
depuis~$\gA_S$.
\end{Proof}
%
%:HHH changement la question ouverte est maintenant fermée
\subsubsection*{Une question longtemps ouverte résolue par la négative}
%------------------

Il s'agit de la question suivante.

\emph{Si tout $\AX$-\mptf est étendu depuis $\gA$, est-il toujours vrai que
pour n'importe quel $r$ tout
$\AXr$-\mptf est étendu depuis $\gA$?}

Une réponse négative est donnée dans \cite[Corti\~nas\&al., (2011)]{CHWW}. 

\penalty-2500
%%%%%%%%%%%%%%%%%%%%%%%%%%%%%%%%%%%%%%%%%%%%%%%%%%%%%%%%%%%%%%%%%%%
\subsubsection*{Un principe \lgb à la Roitman}

%%%%%%%%%%%%%%%%%%%%%%%%%%%%%%%%
%:principe loc glob {plgcetendus}
\begin{plcc} \label{plgcetendus}
Soient $n$ et $r>0$. On considère la \prt suivante pour un anneau $\gA$. 
{\sf P}$_{n,r}(\gA)\,:$ tout \mrc $ r$ sur $\AXn$ est étendu depuis $\gA$.\\
Soient $S_1$, \ldots, $S_k$ des \moco d'un anneau  $\gA$.
Alors $\gA$ satisfait la \prt {\sf P}$_{n,r}$ \ssi chacun des  $\gA_{S_i}$ la satisfait.\\
En particulier $\gA$ est seminormal \ssi chacun des  $\gA_{S_i}$ est seminormal.
\end{plcc}
%--------- fin lemma ---------------------------------------------- 
%
\begin{proof}
La condition est \ncr d'après le  \tho de Roitman \ref{thRoitman},  
dont la \dem reste valable si l'on se limite aux \mrcs  $r$.
\\
La condition est suffisante d'après le recollement de Quillen~\ref{thPatchQ}. 
\end{proof}

%%%%%%%%%%%%%%%%%%%%%%%%%%%%%%%%%%%%%%%%%%%%%%%%%%%%%%%%%%%%%%%%%%%%
%%%%%%%%%%%%%%%%%%%%%%%%%%%%%%%%%%%%%%%%%%%%%%%%%%%%%%%%%%%%%%%%%%%%
%%%%%%%%%%%%%%%%%%%%%%                                 %%%%%%%%%%%%%%%%%%%
%%%%%%%%%%%%%%%%%%%%%%           HORROCKS              %%%%%%%%%%%%%%%%%%%
%%%%%%%%%%%%%%%%%%%%%%                                 %%%%%%%%%%%%%%%%%%%
%%%%%%%%%%%%%%%%%%%%%%%%%%%%%%%%%%%%%%%%%%%%%%%%%%%%%%%%%%%%%%%%%%%%
%%%%%%%%%%%%%%%%%%%%%%%%%%%%%%%%%%%%%%%%%%%%%%%%%%%%%%%%%%%%%%%%%%%%

%---- Section{Horrocks}---------
\section{Le \tho de Horrocks}
\label{sec.Horrocks} 
%------------------

Le lemme suivant est un cas particulier de la proposition~\ref{propFiniBon1}~\emph{4}.
%     lemma  lemLocptfLib
\begin{lemma}\label{lemLocptfLib}
Soit $S$ un \mo de $\gA$ et $P$, $Q$ des \Amos \ptfs tels que
$P_{S}\simeq Q_S$. Alors, il existe $s\in S$ tel que
$P_{s}\simeq Q_s$.
\end{lemma}

%     notation  notaA<X>
\begin{notation}\label{notaA<X>} 
On note $\ArX$ l'anneau $S^{-1}\AX$, où $S$ est le \mo des \polus de $\AX$.\label{NOTAAlraX}
\end{notation}

%:     theorem Horrocks local  thHor0
\begin{theorem}\label{thHor0}\emph{(\Tho de Horrocks local)}\\
Soit $\gA$ un \alo \dcd et $P$ un \mptf sur $\AX$.
Si $P_S$ est libre sur $\ArX$, alors $P$ est libre sur $\AX$ (donc étendu
depuis~$\gA$).
\end{theorem}
%---end-theorem-------------------

 Nous reprenons la preuve de \cite[Nashier \& Nichols]{NaNi} qui est  presque \covz,
telle qu'exposée dans \cite{Lam06} ou \cite{IRa}.

Nous avons besoin de quelques résultats préliminaires.

%:     Lemma{lemLocLocCom}
\begin{lemma}\label{lemLocLocCom}
Soit $\gA$ un anneau, $\fm=\Rad\gA$  et $S\subseteq\AX$
le \mo des \polusz. Les \mos $S$ et $1+\fm[X]$ sont \comz.
\end{lemma}
\begin{proof}
Soit $f(X)\in S$ et $g(X)\in 1+\fm[X]$.
Le résultant $\Res_X(f,g)$ appartient à l'\id $\gen{f,g} $ de $\AX$.
Puisque $f$ est \monz,
le résultant subit
avec succès la spécialisation
$\gA\to\gA\sur\fm$. Donc
$\Res_X(f,g)\equiv\Res_X(f,1)=1  \mod \fm$.
\end{proof}
%
%:     Lemma{lem0thHor0}
\begin{lemma}\label{lem0thHor0}
Soit $\gA\subseteq\gB$,  $s\in\Reg(\gB)$, et $P$, $Q$
deux $\gB$-\mptfs avec $sQ\subseteq P\subseteq Q$. Si $\gB$ et
$\aqo{\gB}{s}$ sont des $\gA$-modules \pros (non \ncrt \tfz), alors
 il en va de même pour \hbox{le \Amoz~$Q/P$}.
\end{lemma}
\begin{proof}
Puisque $s$ est \ndz et que $Q$ et $P$ sont des sous-modules d'un module
libre, la multiplication
par $s$ (notée $\mu_s$) est injective dans $P$ et dans $Q$.
On a les suites exactes de \Amos suivantes.
\[\arraycolsep2pt
\begin{array}{ccccccccccccccccc}
0 &\rightarrow& Q& \vvers{\mu_s} &P & \lora & P/sQ& \rightarrow& 0   \\[1mm]
0 &\rightarrow& sQ/sP& \llra &  P/sP& \lora & P/sQ& \rightarrow& 0   
\end{array}
\]
Le \Amo $P$ est \pro par transitivité, le $\gB/s\gB$-module $P/sP$
est \proz, donc par transitivité   $P/sP$ est un \Amo \proz.
On peut alors appliquer le lemme de Schanuel (lemme \ref{lemScha}):
$(P/sP)\oplus Q\simeq (sQ/sP) \oplus P$ comme \Amosz. Puisque
$Q$ est un \Amo \proz, il en va de même pour $sQ/sP$. Mais
puisque $\mu_s$ est injective, $sQ/sP$ est isomorphe à $P/Q$.
\end{proof}
%

%:     Lemma{lem1thHor0}
\begin{lemma}\label{lem1thHor0}\emph{(Murthy \& Pedrini, \cite{MuPe})}
\\ 
Soient $\gA$ un anneau, $\gB=\AX$, $S$ le \mo des \polus de~$\AX$$,
P$, $Q$ deux \mptfs sur $\gB$, et
$f\in S$.
\begin{enumerate}
\item Si $fQ \subseteq P\subseteq Q$, alors $P$ et $Q$ sont stablement isomorphes.
\item Si $P_S\simeq Q_S$,  alors $P$ et $Q$ sont stablement isomorphes.
\item Si en plus $P$ et $Q$ sont de rang $1$,
alors $P\simeq Q$.
\end{enumerate}
\end{lemma}
\begin{proof}
\emph{1.} Puisque $f$ est \monz, l'\Alg $\aqo{\gB}{f}$ est un \Amo libre de rang $\deg f$. Le $\aqo{\gB}{f}$-module $Q/fQ$ est \ptf sur $\gA$.  D'après le lemme précédent, le \Amo $M=Q/P$ est 
 \proz. Et il est \tf sur $\aqo{\gB}{f}$,
donc sur $\gA$. Donc $M[X]$ est un \Bmo \ptfz.
Nous avons deux suites exactes ($\mu_X$ est la multiplication par $X$)
\[\arraycolsep2pt
\begin{array}{ccccccccccccccccc}
0 &\rightarrow& P& \llra &Q& \lora& M& \rightarrow& 0 ,  \\[1mm]
0 &\rightarrow& M[X]& \vvers{\mu_X} &M[X]& \lora & M& \rightarrow& 0.   
\end{array}
\]
D'après le lemme de Schanuel (lemme \ref{lemScha}) on a
$P\oplus M[X]\simeq Q\oplus M[X]$ comme \Bmosz. Puisque $M[X]$ est \ptf
sur $\gB$, $P$ et $Q$ sont stablement isomorphes.

\emph{2.} On sait que $P_f\simeq Q_f$ pour un $f\in S$.\\
Par hypothèse,  on a $F\in\GAn(\gB)$  et $G\in\GAn(\gB)$ 
avec $P\simeq \Im F$,
et~\hbox{$Q\simeq \Im G$}. 
Nous savons que $F'=\Diag(F,0_m)$ et $G'=\Diag(G,0_n)$
sont conjuguées sur $\gB_f$
(lemme \dlgz~\ref{propIsoIm}). Ceci signifie qu'il
existe une matrice $H\in\MM_{m+n}(\gB)$
telle que $HF'=G'H$ et $\det(H)=\delta$ divise une puissance
de $f$. On a alors $P_1=\Im(HF')\subseteq\Im G'$.
Puis, en postmultipliant par $\wi H$, $(HF')\wi H=\delta G'$,
ce qui implique $\delta \Im G'\subseteq \Im(HF')$.
Puisque $H$ est injective, on a $P_1\simeq P$, et par ailleurs
$\Im G'=Q_1\simeq Q$. On
peut conclure d'après le point \emph{1} puisque
 $\delta Q_1\subseteq P_1\subseteq Q_1$.

\emph{3.} Les modules
$P$ et $Q$ sont de rang $1$ et stablement isomorphes, donc
 isomorphes (fait \ref{factPicStab}).
\end{proof}

\begin{Proof}{\Demo du \thref{thHor0}. } ~\\
Notations:  $\fm=\Rad\gA$, $\gk=\gA/\fm$ (\cdiz),
$\gB=\AX $, $n=\rg(P)$  ($n\in\NN$ car $\gB$ est connexe), $U=1+\fm [X]$,
et $\ov{E}$ l'objet $E$ réduit modulo $\fm$.
%On considère un \pol \mon  $f$  dans $\gB$ tel que $P_f$ est libre.

\emph{1.} On montre par \recu sur $n$ que l'on  a
un \iso $P\simeq P_1\oplus\gB^{n-1}$. Pour $n=1$ c'est trivial.

%:     Lemma{lem2thHor0}
\textbf{Petit lemme} {(voir la \dem plus loin)}\\
\emph{Il
existe  $z$, $y_2$, \ldots, $y_n$,  $z_2$, \ldots, $z_n$ dans $P$
tels que $(z,y_2, \ldots,y_n)$ est une base de $P_S$ sur $\gB_S$
et $(\ov{z}, \ov{z_2}, \ldots, \ov{z_n})$ est une base de $\ov{P}$ sur $\ov{\gB}=\gk[X]$.}

Le $\gB_U$-module $P_U$ est libre de base
$(z,  z_2, \ldots, z_n)$: en effet, $\fm\subseteq\Rad\gB_U$,
et modulo $\fm$, $(\ov{z}, \ov{z_2}, \ldots, \ov{z_n})$
engendre  $\ov{P}=\ov{P_U}$, donc $(z,  z_2, \ldots, z_n)$
engendre~$P_U$ par le lemme de Nakayama. Enfin un \mptf de rang~$n$ 
engendré par $n$ \elts est libre.\\
On pose $P'=P/\gB z$. Les deux modules $P'_U$ et $P'_S$
sont libres. Les \mosz~$U$ et $S$ sont \com
(lemme \ref{lemLocLocCom}), donc $P'$ est \ptf sur~$\gB$,
d'où $P\simeq P'\oplus \gB z$.
Par \hdrz, $P'\simeq P_1\oplus \gB^{n-2}$, ce qui donne
$P\simeq P_1\oplus \gB^{n-1}$

\emph{2.} L'\iso $P\simeq P_1\oplus \gB^{n-1}$ avec $P_1$ de rang $1$
donne par \lon que $(P_1)_S$ est stablement libre.
On applique le point \emph{3} du lemme \ref{lem1thHor0}: on obtient
que $P_1$ est libre.
\end{Proof}

\begin{Proof}{\Demo du petit lemme. }\\
Soit, dans le module $P$,  une $\gB_S$-base $(\yn)$  de~$P_S$.
Il existe une base  $(\ov{z_1}, \ov{z_2}, \ldots, \ov{z_n})$
de $\ov P$, avec les $z_i$ dans~$P$ telle que $\ov{y_1}\in\kX\,\ov{z_2}$
(en divisant~$\ov{y_1}$ par le pgcd de ses \coesz, on obtient un \vmdz,
et sur un anneau de Bézout, tout \vmd est complétable).
On cherche $z$ sous la forme~\hbox{$z_1+X^ry_1$}.
Il est clair que, pour n'importe quel $r$,   $(\ov{z},
\ov{z_2}, \ldots, \ov{z_n})$ est une base de $\ov P$.
Puisque $(\yn)$  est une base de~$P_S$
sur~$\gB_S$, il existe~\hbox{$s\in S$} tel que $sz_1=\sum_{i=1}^nb_iy_i$, avec
les $b_i$ dans~$\gB$. \\
Alors, $sz=(b_1+sX^r)y_1+\sum_{i=2}^nb_iy_i$,
et pour $r$ assez grand,   $b_1+sX^r$ est un \poluz:
$(z, y_2, \ldots, y_n)$ est une base de $P_S$ sur $\gB_S$.
\end{Proof}
On donne maintenant la version globale.

%:     theorem Horrocks global   thHor
\begin{theorem}\label{thHor} \emph{(\Tho de Horrocks global)}\\
Soit $S$ le \mo des \polus de $\AX$ et
$P$ un \mptf sur $\AX$. Si $P_S$ est étendu depuis
$\gA$, alors $P$  est étendu depuis~$\gA$.
\end{theorem}
%---end-theorem-------------------
%-----------------begin proof------------------
\begin{proof}
Nous appliquons la machinerie \lgbe
de base \paref{MethodeIdeps} avec la \dem
\cov du \thrf{thHor0}.\imlb Nous obtenons une
famille finie de \moco de $\gA$, $(U_i)_{i\in J}$,   avec chaque localisé $P_{U_i}$
étendu depuis  $\gA_{U_i}$.
On conclut avec le \rcm de Quillen\index{Quillen!recollement de ---} (\plgc \ref{thPatchQ}).
\end{proof}
%-----------------end proof------------------

Cet important \tho peut être complété par le résultat subtil suivant,
qui ne semble pas pouvoir être étendu aux \mpfsz.

%:     Theorem{th2Bass}
\begin{theorem}\label{th2Bass} \emph{(Bass)}\\
Soient $P$ et $Q$ deux \Amos \ptfsz. S'ils sont isomorphes après extension
des scalaires à $\ArX$, ils sont isomorphes.
\end{theorem}
\begin{proof}
Nous raisonnons avec des \mprns et des similitudes entre ces matrices
qui correspondent à des \isos entre les modules images. Implicitement donc,
nous utilisons de manière systématique le lemme \dlg \ref{propIsoIm},
sans le mentionner.
\\ 
On démarre avec  $F$ et $G$ dans $\GAn(\gA)$, conjuguées sur
l'anneau $\ArX$. Les \mptfs sont  $P \simeq \Im F$ et $Q \simeq \Im G$. Nous avons donc
une matrice~\hbox{$H\in\Mn(\AX)$}, avec $\det(H)\in S$ (\mo des \polusz), et
$HF=GH$.\\
En posant $Y=1/X$, pour $N$ assez grand,
 la matrice $Y^N H=H'$ est dans~\hbox{$\Mn(\AY)$},
avec $\det(H')=Y^r\big(1+Yg(Y)\big)=Y^rh(Y)$ où $h(0)=1$, et
évidemment $H'F=GH'$. Autrement dit,~$F\sim G$~sur l'anneau $\gA[Y]_{Yh}$. \\
Les \elts $Y$ et $h$ sont \comz,
donc, par application du \tho de \rcm des modules
(\plgc \ref{plcc.modules 2}), il existe un $\AY$-module $M$ tel que
$M_Y$ est isomorphe %à $\AY_Y\te_\gA P$ 
à~\gui{$P$ étendu à $\AY_Y$}, 
et~$M_h$ est isomorphe à \gui{$Q$ étendu à $\AY_h$}. Et~$M$ est  \ptf
puisqu'il a deux \lons \come qui sont des \mptfsz.
Ceci nous fournit une \mprn $E$ à \coes dans $\AY$ telle
que $E\sim F$ sur $\AY_Y$ et $E\sim G$ sur $\AY_h$. Puisque $Y$ est un \poluz, le \tho de Horrocks nous dit que $\Im E$ provient par \eds d'un \Amo \ptf $M'$.
En conséquence, pour tous $a$, $b\in \gA$ les matrices \gui{évaluées}
$E(a)$ et $E(b)$ sont conjuguées sur $\gA$
(leurs images sont toutes deux isomorphes à $M'$).\\
Enfin $F\sim E(1)$ et $G\sim E(0)$ sur $\gA$, donc $F \sim G$ sur $\gA$.
\end{proof}

\rem Pour \la mathématicie\n qui désire implémenter l'\algo
sous-jacent à la \dem précédente, on suggérera d'utiliser
des \mpns (des \mptfs considérés) plutôt  que des \mprns
(dont les images sont isomorphes à ces modules). Cela évitera en particulier
d'avoir à utiliser de manière répétée une implémentation
du lemme \dlgz. 
\eoe

\medskip
Nous terminons cette section avec un corolaire du
lemme \ref{lem1thHor0}. Ce \tho est à comparer avec le \thref{th2QUILIND}.
\perso{Je n'ai pas de référence pour ce \thoz}

%:     theorem{thStabLibPol}
\begin{theorem}\label{thStabLibPol} \emph{(Induction de Quillen\index{Quillen!induction de ---} concrète, cas stablement libre)}\\
Soit $\cF$ une classe  d'anneaux qui satisfait les \prts suivantes.
\perso{donner un nom et rajouter une référence dans la biblio}
\hsu
\emph{1.} Si $\gA\in\cF$, alors $\gA\lra{X} \in\cF$.\hsu
\emph{2.} Si $\gA\in\cF$, tout
  $\gA$-module de rang constant est stablement libre.\\
Alors, pour  $\gA\in\cF$ et  $r\in\NN$, tout  $\AXr$-module de rang constant  est stablement libre.
\end{theorem}
\begin{proof}
On fait une preuve par \recu sur $r$, le cas $r=0$ est clair.\\
Nous passons de $r-1$ à $r$ ($r\geq 1$). Soit $\gA$ un anneau dans la classe $\cF$,
et $P$ un \mrc sur $\AXr$.\\
On note $\gB=\gA[(X_i)_{i<r}]$, $\gC=\gA[X_r]$, et $V$ est le \mo des \polus de $\gA[X_r]$.
Ainsi 
$\AXr\simeq \gB[X_r]\simeq \gC[(X_i)_{i<r}].$\\
%:HHH   V^{-1}\gC  au lieu de  V^{-1}\gL
L'anneau $\gA\lra{X_r}=V^{-1}\gC$  est dans
la classe $\cF$.\\
Le $\gA\lra{X_r}[(X_i)_{i<r}]$-module $P_V$, qui
est \prcz,
est \stl par \hdrz. \\
Si $S$ est le \mo des \polus de $\gB[X_r]$, on a $V\subseteq S$, et donc~$P_S$ est
\stl sur l'anneau~\hbox{$S^{-1}\gB[X_r]$}. Par le point \emph{2} du lemme \ref{lem1thHor0}, $P$ est \stlz.
\end{proof}

%:     Corollary{corthStabLibPol}
\begin{corollary}\label{corthStabLibPol}
Si $\gK$ est un \cdiz, tout \mptf sur $\KXr$ est stablement libre.
\end{corollary}
\begin{proof}
On applique le résultat précédent avec la classe $\cF$ des \cdisz:
 si~$\gK$ est un \cdiz, alors $\gK\lra{X}=\gK(X)$ est aussi un \cdiz.
\end{proof}
%

%%%%%%%%%%%%%%%%%%%%%%%%%%%%%%%%%%%%%%%%%%%%%%%%%%%%%%%%%%%%%%%%%%%
%%%%%%%%%%%%%%%%%%%%%%%%%%%%%%%%%%%%%%%%%%%%%%%%%%%%%%%%%%%%%%%%%%%
%%%%%%%%%%%%%%%%%%%%%%                                 %%%%%%%%%%%%%%%%%%%%
%%%%%%%%%%%%%%%%%%%%%%       QUILLEN SUSLIN            %%%%%%%%%%%%%%%%%%%%
%%%%%%%%%%%%%%%%%%%%%%                                 %%%%%%%%%%%%%%%%%%%%
%%%%%%%%%%%%%%%%%%%%%%%%%%%%%%%%%%%%%%%%%%%%%%%%%%%%%%%%%%%%%%%%%%%
%%%%%%%%%%%%%%%%%%%%%%%%%%%%%%%%%%%%%%%%%%%%%%%%%%%%%%%%%%%%%%%%%%%

%---- Section{Quillen-Suslin}---------
\section{Solution de la conjecture de Serre}
\label{sec.QS} 
%------------------

Dans cette section nous exposons plusieurs solutions \covs au \pb de Serre,
dans lequel  $\gK$ est un corps discret.

\Grandcadre{Les \mptfs sur $\KXr$ sont  libres}

%%%%%%%%%%%%%%%%%%%%%%                                 %%%%%%%%%%%%%%%%%%%%
%%%%%%%%%%%%%%%%%%%%%%       A LA QUILLEN              %%%%%%%%%%%%%%%%%%%%
%%%%%%%%%%%%%%%%%%%%%%                                 %%%%%%%%%%%%%%%%%%%%

%:   --- SUBsection{Quillen}-------------------------
\subsec{\`A la Quillen\index{Quillen}}
\label{subsecQuillen}
%-----------------------------------------

La solution par Quillen du \pb de Serre est basée sur le \tho
de Horrocks local et sur l'\ix{induction de Quillen} suivante (voir~\cite{Lam06}). 

%:     Induction de Quillen abstraite
\CMnewtheorem{IQa}{Induction de Quillen\index{Quillen!induction de --- abstraite} abstraite}{\itshape}
\begin{IQa}\label{InductionQabs}~\\
Soit $\cF$ une classe  d'anneaux qui satisfait les \prts suivantes.
%-----------------begin item------------------
\begin{description}
\item [{\rm (Q1)}]
Si $\gA\in\cF$, alors $\ArX \in\cF$.
\item [{\rm (Q2)}]
Si $\gA\in\cF$, alors $\gA_{{\fm}}\in\cF$ pour tout \idema ${\fm}$ de $\gA$.
\item [{\rm (Q3)}]
Si $\gA\in\cF$ est local, tout
  $\AX$-\mptf  est étendu depuis $\gA$ (i.e., libre).
\end{description}
%-----------------end item------------------
Alors, pour tout $\gA\in\cF$ et tout $r\geq1$, tout \mptf sur~$\AXr$ est étendu depuis $\gA$. 
\end{IQa}
%--------- fin lemma ---------------------------------------------- 

En fait, les \prts (Q1), (Q2) et (Q3) sont d'abord utilisées par Quillen pour obtenir le cas $r=1$, en utilisant le \tho de Horrocks local et le recollement de Quillen.
La partie \gui{preuve par \recuz} est basée sur le cas $r=1$, sur (Q1)
et sur le \tho de Horrocks (local ou global).

Dans la suite nous isolons cette preuve par \recuz,
que nous qualifions d'induction de Quillen \gui{concrète}.
Nous remplaçons (Q3) par une version plus forte (q3) qui est le cas $r=1$.

Dans un commentaire postérieur, nous expliquons comment nous pouvons, 
en fait, remplacer d'une certaine manière (q3) par (Q3) sans pour autant
perdre le \crc \cof de la \demz.  

%\penalty-2500
\subsubsection*{La preuve par \recu proprement dite}

%:     Theorem{thQUILIND}
\begin{theorem}\emph{(Induction de Quillen\index{Quillen!induction de --- concrète} concrète)} \label{thQUILIND}~\\
Soit $\cF$ une classe  d'anneaux qui satisfait les \prts suivantes.
\begin{description}
  \item[{\rm (q1)}] Si $\gA\in\cF$, alors $\ArX \in\cF$.
%\item[{\rm (q2)}] Si $\gA\in\cF$ et $S$ est un \mo de $\gA$ alors $\gA_S \in\cF$.
  \item[{\rm (q3)}] Si $\gA\in\cF$, tout $\AX$-\mptf  est étendu depuis~$\gA$.
\end{description}
Alors, pour tout $\gA\in\cF$ et tout $r\geq1$, tout \mptf sur~$\AXr$ est étendu depuis $\gA$.
\end{theorem}
%%%%%%%%%%%%%%%%%%%%%%%%%%%%%%%%%%%%%%%%%
\begin{proof}
Passons de $r\geq1$ à $r+1$. On considère un $\gA[\Xr,Y]$-\mptf $P=P(\Xr,Y)=P(\uX,Y)$. 
 On note
\begin{enumerate}
  \item [--] $P(\uX,0)$ le $\AuX$-module obtenu par l'\homo $Y\mapsto 0$,
  \item [--] $P(\uze,Y)$ le $\gA[Y]$-module obtenu par l'\homo
  $\uX\mapsto \uze$,
  \item [--] $P(\uze,0)$ le \Amo  obtenu par l'\homo $\uX,Y\mapsto \uze,0$.
\end{enumerate}
On doit montrer que $P(\uX,Y)\simeq P(\uze,0)$ sur $\gA[\uX,Y]$.
\\ 
On appelle $S$ le \mo des \polus de $\gA[Y]$, qui est contenu
dans le \mo $S'$ des \polus de $\AuX[Y]$.
On a alors:
\begin{enumerate}
  \item  $P(\uX,Y)\simeq P(\uze,Y)$ sur $\gA\lra{Y}[\uX]=\gA[\uX,Y]_{S}$ 
  par \hdr puisque  $\gA\lra{Y}\in \cF$,
  \item  a fortiori $P(\uX,Y)\simeq P(\uze,Y)$ sur
  $\gA[\uX]\lra{Y}=\gA[\uX,Y]_{S'}$,
  \item  $P(\uze,Y)\simeq P(\uze,0)$ sur $\gA[Y]$ par le cas $r=1$,
  \item  $P(\uze,0)\simeq P(\uX,0)$ sur $\AuX$ par \hdrz,
  \item  en combinant 2, 3 et 4, on a $P(\uX,Y)\simeq P(\uX,0)$ sur
   $\gA[\uX]\lra{Y}$,
  \item  donc, par le \tho de Horrocks global,  $P(\uX,Y)\simeq P(\uX,0)$ sur l'anneau
  $\gA[\uX,Y]$,
  \item  on combine ce dernier \iso avec l'\iso entre $P(\uX,0)$ et $P(\uze,0)$ sur
 l'anneau $\gA[\uX]$ obtenu par \hdrz.
\end{enumerate}
\vspace{-1em}\end{proof}

%:     Corollary{corthQUILIND}
\begin{corollary}\label{corthQUILIND} \emph{(\Tho de Quillen-Suslin, preuve
de Quillen\index{Quillen})}\\
Si $\gK$ est un corps discret (resp. un anneau \zedz),
tout \mptf sur $\KXr$ est libre (resp. quasi libre).
\end{corollary}

\begin{proof}
L'induction de Quillen concrète s'applique avec la classe $\cF$ des \cdisz:
on note que  
$\KX$ est un anneau de Bézout intègre, donc les \mptfs sur $\KX$ sont libres, et a fortiori étendus. 
On passe aux anneaux \zedrs par la machinerie \lgbe \elr \num2. 
Enfin, pour les anneaux \zedsz, on utilise l'\egt $\GKO(\gA)=\GKO(\Ared)$.\imlgz
\end{proof}
\rems
1) On rappelle qu'un anneau \zed est connexe \ssi il est local.
Si $\gK$ est un tel anneau, tout \mptf sur $\KXr$ est libre.
\\
 2) L'induction de Quillen concrète s'applique aux
anneaux de Bézout intègres de \ddk $\leq 1$ 
(voir l'exercice \ref{exoBézoutKdim1TransfertArX}) et plus \gnlt aux \adps de
dimension $\leq1$ (voir le \thref{thMaBrCo}). Ceci généralise le cas des \dDks obtenu par Quillen. Pour le cas des anneaux \noes réguliers de \ddk $\leq2$ (que nous ne traiterons pas dans cet ouvrage),
voir \cite{Lam06}. 
\eoe

\subsubsection*{(Q3) versus (q3)}

L'induction de Quillen abstraite (qui ne fournit pas de résultat sous forme
\covz) présente l'avantage d'utiliser une hypothèse (Q3) plus faible que
l'hypothèse (q3) utilisée dans l'induction concrète. 
Nous expliquons maintenant comment nous pouvons récupérer \cot la mise,
même pour l'hypothèse~(Q3).

\sni{\it Le cas libre.}

 Dans le cas où la classe $\cF$ est telle que les \mptfs sont libres,
 on remarque que l'hypothèse (q3) est en fait inutile. En effet,
soient~$P$ un $\AX$-\mptf et $S$ le \mo des \polus de $\AX$.
Alors, d'après (q1) le $\ArX$-module $P_S$ est libre, donc étendu
depuis~$\gA$. Mais alors, d'après le \tho de Horrocks global, le module $P$ est étendu
depuis~$\gA$. Autrement dit nous avons démontré la version suivante,
adaptée au cas libre, et particulièrement simple.

%:     Theorem{th2QUILIND}
\begin{theorem}\label{th2QUILIND} \emph{(Induction de Quillen\index{Quillen!induction de --- concrète, cas libre} concrète, cas libre)}\\
Soit $\cF$ une classe  d'anneaux qui satisfait les \prts suivantes.
\begin{description}
  \item[{\rm (q0)}] Si $\gA\in\cF$, tout $\gA$-\mptf  est libre.
  \item[{\rm (q1)}] Si $\gA\in\cF$, alors $\ArX \in\cF$.
%\item[{\rm (q2)}] Si $\gA\in\cF$ et $S$ est un \mo de $\gA$ alors $\gA_S \in\cF$.
\end{description}
Alors, pour tout $\gA\in\cF$ et tout $r\geq1$, tout \mptf 
sur~$\AXr$  est libre.
\end{theorem}
%%%%%%%%%%%%%%%%%%%%%%%%%%%%%%%%%%%%%%%%%
%-% ENTRE NOUS
\entrenous{1) il semble que cela fonctionne aussi avec quasi libre à la place de
libre, car on doit avoir $\BB(\ArX)=\BB(\gA)$. 

2) Cela serait intéressant de voir si l'on est capable de faire fonctionner
\cot cette induction de Quillen concrète pour la classe des corps non \ncrt discrets. Le (q0) est clair, le (q1) nettement moins. 
}
%-% Fin ENTRENOUS

\sni{\it Le cas \gnlz.}

On aura remarqué que la \prt  (Q2) n'intervient pas dans  l'induction de Quillen concrète: cette hypothèse est
rendue inutile par l'hypothèse (q3).

La \prt  (Q2) intervient cependant lorsque nous voulons remplacer (q3) par (Q3), qui est une hypothèse a priori plus faible que (q3).

Nous pensons que cet affaiblissement de l'hypothèse est toujours possible en pratique, sans pour autant perdre le \crc \cof du résultat. Mais comme ceci est basé sur la machinerie \lgbe de base 
(machinerie \lgbe à \idepsz), et comme cette dernière est une méthode de
\dem et non pas un \tho à proprement parler, nous n'avons pas pu formuler
notre induction concrète directement avec (Q3), car nous voulions un \tho
en bonne et due forme.

Venons en à l'explication du remplacement de l'hypothèse forte (q3) par l'hypothèse faible (Q3). 

Nous reprenons l'hypothèse (Q2) sous la forme plus \gnle suivante.

{\rm (q2)} Si $\gA\in\cF$ et $S$ est un \mo de $\gA$, alors $\gA_S \in\cF$.

 Nous supposons que (Q3) est satisfaite sous la forme suivante: 
sous l'hypothèse que $\gA$ est un \alrd dans la classe $\cF$ 
on a une \dem \cov du fait que tout
\mptf $P$ sur $\AX$ est étendu, ce qui se traduit par un \algo de calcul
(pour l'\iso entre $P$ et $P(0)$) basé sur les \prts de la classe $\cF$ 
et sur la disjonction
$$
a\in\Ati  \quad\hbox{ou}\quad a\in\Rad(\gA)
$$
pour les \elts $a$ qui se présentent au cours de l'\algoz.
Dans ces conditions la machinerie \lgbe de base s'applique.\imlb
En conséquence pour un \mptf $P$ sur $\AX$ pour un anneau $\gA\in\cF$ arbitraire, la \dem donnée dans le cas local \dcdz, suivie pas à pas, nous fournit des \moco $S_1$, \ldots, $S_\ell$ tels que pour chacun
d'entre eux, le module $P_{S_i}$ (sur $\gA_{S_i}[X]$) 
est étendu depuis~$\gA_{S_i}$. Notez que pour que cette méthode fonctionne, la classe d'anneaux considérée doit vérifier (q2), et que l'on peut se limiter aux \lons en des \mos
$\cS(\an;b)$. Il ne reste alors qu'à appliquer le recollement de Quillen\index{Quillen!recollement de ---}
(\plgc \ref{thPatchQ}) pour obtenir le résultat souhaité: 
le module $P$ est étendu depuis $\gA$.

%%%%%%%%%%%%%%%%%%%%%%                                 %%%%%%%%%%%%%%%%%%%
%%%%%%%%%%%%%%%%%%%%%%       A LA SUSLIN               %%%%%%%%%%%%%%%%%%%
%%%%%%%%%%%%%%%%%%%%%%                                 %%%%%%%%%%%%%%%%%%%

%:  --- SUBsection{Suslin}-------------------------
\subsec{\`A la Suslin, Vaserstein ou Rao}
\label{subsecSuslin}
%-----------------------------------------

La solution par \Sus de la conjecture de Serre consiste à montrer que
tout module stablement libre sur $\KXr$ est libre (Serre avait déjà
démontré que tout \mptf sur $\KXr$ est stablement libre),
autrement dit que le noyau de toute matrice surjective est libre,
ou encore que tout \vmd est la première colonne d'une matrice \iv
(cf. fait \ref{factStablib} et proposition~\ref{corpropStabliblib}).

 Si \rdb\label{NOTAfGg} $\cG$ est un sous-groupe de $\GL_n(\gA)$ et $A$, $B\in \Ae{n\times 1}$, nous
noterons~\hbox{$A\sims{\cG }B$} pour dire qu'il existe une
matrice $H\in\cG $ telle que $HA=B$.
Il est clair qu'il s'agit d'une relation d'\eqvcz.

Rappelons qu'un \vmd $f\in\Ae {n\times 1}$ est dit complétable s'il est le
premier vecteur colonne d'une matrice $G\in\GL_n(\gA)$. Cela revient à dire~que l'on a  
$$\preskip-.3em \postskip.4em 
 f\;\sims{\GLn(\gA)}\;\tra{[\,1\;0\;\cdots\;0\,]}.
$$
Le but dans ce paragraphe est donc d'obtenir une \dem \cov
du \tho suivant.

%:     Theorem{thSuslinQS}
\begin{theorem}\label{thSuslinQS} \emph{(\Susz)}\\
Tout \vmd $f$ à \coos
dans $\KXr=\KuX$ (où $\gK$ est un corps discret) 
est complétable.
%on a
%$$
%\cmatrix{f_1(\uX)\cr f_2(\uX)\cr\vdots\cr f_n(\uX)}\;\;\;\sims{\GL_n(\KuX)}
%\;\;\;\cmatrix{1\cr0\cr\vdots\cr0}.
%$$
\end{theorem}

Nous donnerons trois \dems distinctes, par ordre chronologique.

%%%%%%%%%%%%%%%%%%%%%%       A LA SUSLIN               %%%%%%%%%%%%%%%%%%%%

\subsubsection*{Première \dem}

Nous suivons ici de très près la \dem originale de Suslin.
Nous devons seulement nous débarrasser d'une utilisation non \cov
d'un \idema générique, et nous avons déjà fait ce travail
lorsque nous avons donné une \dem \cov du lemme de Suslin \ref{lemSuslin1}
au chapitre~\ref{chapPlg}.

%\newpage

%:     fact{factM2}
\begin{fact}\label{factM2}
Soient $M,N\in\MM_2(\gA)$. On a $\Tr(M)\,\I_2=M+\wi M$
et

\snic{\det(M+N)=\det (M) + \Tr(\wi{M}\,N) + \det( N).}
\end{fact}
\begin{proof}
Pour les matrices dans $\MM_2(\gA)$, l'application $M\mapsto\wi M$ est \linz, donc

%\vspace{-1mm}
\snac{\arraycolsep2pt
\begin{array}{rcl}
 \det (M+N)\, \I_2\;=\; (\wi M +  \wi N)(M+N)&=&
\wi M M + (\wi M N + \wi N M)+ \wi N N\\
  & = & \big(\det (M) + \Tr(\wi{M}\,N) + \det( N)\big)\,\I_2.
\end{array}}
\vspace{-1mm}
\end{proof}
%
%:     Lemma{lem02SusQS}
\begin{lemma}\label{lem02SusQS}
Soit  $B\in\MM_2(\gA),\,H=H(X)\in\MM_2(\AX)$, $\gB$
une \Alg et $x\in\gB$.
On pose    $C(X)=B+XH$.
On suppose $\det C=\det B=a$.
En notant $S=\I_2+ x\wi H(ax)\,B $, on a alors  $S\in\SL_2(\gA)$ et $S\wi B=\wi C(ax)$.
\end{lemma}
\begin{proof}
Le fait \ref{factM2} donne
$\det (C)=\det (B) + X\,(\Tr(\wi{H}\,B) + X \det H)$,
et donc
$$\preskip.3em \postskip.3em
E(X)=\Tr(\wi{H}\,B) + X \det H=0.
$$
Posons $H_1=H(ax)$ et $C_1=C(ax)$.\\
On a alors $S\wi B=\wi B+ x\wi{H_1} B\wi B=\wi B+ax\wi{H_1}=
\wi{C_1}$ et
$$\preskip.3em \postskip-1em 
\begin{array}{rcl}
\det(S)&=&1+x\Tr(\wi{H_1} B)+\det(x\wi{H_1}B)\\
&=&1+x\Tr(\wi{H_1} B)+x^2a\det(H_1) = 1+x E(ax)=1.
\end{array} 
$$
\end{proof}
%

%:     Lemma{lem2SusQS}
\begin{lemma}\label{lem2SusQS} \emph{(Lemme de \Susz)}\index{Lemme de Suslin}\\
Soient $u,\,v\in\AX$, $a\in\gA\,\cap\,\gen{u,v}$, $\gB$
une \Alg et $b,\,b'\in\gB$. \\
Si $
b\equiv b' \mod a\gB$,  alors $
\Cmatrix{2pt}{u(b)\cr v(b)}  \sims{\SL_2(\gB)}
\cmatrix{u(b')\cr v(b')}$.
\end{lemma}
\begin{proof}
Soient $p, q\in\AX$ tels que $up+vq=a$ et $x\in\gB$ tel que $b'=b+ax$.\\
Considérons la matrice $M=\Cmatrix{3pt}{p&q\cr-v&u}\in\MM_2(\AX)$.
On applique le lemme~\ref{lem02SusQS}~avec les  matrices $B=M(b)$ et $C(X)=M(b+X)$.\\
Notez que la première colonne de $\wi B$ est $\Cmatrix{2pt}{u(b)\cr v(b)}$
et que la première colonne de  $\wi C(ax)$ est $\Cmatrix{2pt}{u(b')\cr v(b')}$.
\end{proof}
%

%:     Lemma{lem3SusQS}
\begin{lemma}\label{lem3SusQS}
Soient $f\in \AX^{n\times 1}$, $\gB$ une \Alg   et $\cG$
un sous-groupe de~$\GLn(\gB)$, alors l'ensemble
$$
\fa=\sotQ{a\in\gA}{\Tt b,\,b'\in\gB, \big((b\equiv b' \mod a\gB) \;\Rightarrow\;f(b)\sims{\cG}f(b')\big)}
$$
est un \id de $\gA$.
\end{lemma}
\facile
%

%:     Theorem{th4SusQS}
\begin{theorem}\label{th4SusQS}
Soient $n\geq2$, $f$ un \vmd de $\AX^{n\times 1}$ avec~$f_1$ \monz,
 $\gB$ une \Algz, et $\cG\subseteq \GLn(\gB)$ le sous-groupe engendré
 par $\En(\gB)$ et $\,\SL_2(\gB)
 (\footnote{$\SL_2(\gB)$ est plongé dans $\GLn(\gB)$ par l'injection  $A\mapsto \Diag(A,\I_{n-2})$.})$.
 Alors, pour tous $b,\,b'\in\gB$, on a $f(b)\sims{\cG}f(b')$.
\end{theorem}
\begin{proof}
Il nous suffit de montrer que l'\id $\fa$ défini au lemme \ref{lem3SusQS}
contient $1$.
Pour une matrice \elr $E=E(X)\in\EE_{n-1}(\AX)$, nous considérons
le vecteur
$$
\cmatrix{g_2\cr\vdots\cr g_n}\;=\;E\,\cmatrix{f_2\cr\vdots\cr f_n}.
$$
Nous allons montrer que le résultant $a=\Res_X(f_1,g_2)$, qui est bien
défini puisque $f_1$ est \monz, est un \elt de $\fa$. Nous aurons donc terminé en invoquant le lemme de Suslin \ref{lemSuslin1}.\\
Montrons donc que $a\in\fa$. Nous utilisons juste le fait que $a\in\gen{f_1,g_2}\cap\gA$.
On prend $b,\,b'\in\gB$ avec $b\equiv b'\mod a\gB$.
On veut aboutir à $f(b)\sims{\cG}f(b')$.
Notons que pour $i\geq2$ on~a:
\begin{equation}\label{eqth4SusQS}
\begin{array}{ccc}
g_i(b')-g_i(b)\in \gen{b'-b}\subseteq \gen{a}  \subseteq \gen{f_1(b),g_2(b)},
\\[1mm]
\hbox{i.e., } ~~~~
g_i(b')\in g_i(b) + \gen{f_1(b),g_2(b)}.
\end{array}
\end{equation}
On a alors une suite d'\eqvcs
$$\hss \mathrigid1.5mu
\Cmatrix{.18em}{f_1(b)\cr f_2(b)\cr f_3(b)\cr\vdots\cr f_n(b)} \sims{E(b)}
\cmatrix{f_1(b)\cr g_2(b)\cr g_3(b)\cr\vdots\cr g_n(b)} \sims{\En(\gB)}
\cmatrix{f_1(b)\cr g_2(b)\cr g_3(b')\cr\vdots\cr g_n(b')} \sims{\SL_2(\gB)}
\cmatrix{f_1(b')\cr g_2(b')\cr g_3(b')\cr\vdots\cr g_n(b')} \sims{E(b')^{-1}}
\cmatrix{f_1(b')\cr f_2(b')\cr f_3(b')\cr\vdots\cr f_n(b')}
\hss.$$
La seconde est donnée par l'\eqrf{eqth4SusQS}, la troisième
par le lemme \ref{lem2SusQS} appliqué à $u=f_1$ et $v=g_2$.
\end{proof}
%

%:     Corollary{corth4SusQS}
\begin{corollary}\label{corth4SusQS}
Soient $n\geq2$, $f$ un \vmd de $\AX^{n\times 1}$ avec~$f_1$ \mon
et $\cG$ le sous-groupe de~$\GLn(\AX)$ engendré
 par $\En(\AX)$ et $\,\SL_2(\AX)$. Alors  $f\sims{\cG}f(0)$.
\end{corollary}
\begin{proof}
Dans le \thrf{th4SusQS}, on prend $\gB=\AX$, $b=X$ et~$b'=0$.
\end{proof}
%

%:     Corollary{corth5SusQS}
\begin{corollary}\label{corth5SusQS}
Soient $\gK$ un corps discret,  $n\geq2$, $f$ un \vmd de~$\KuX^{n\times 1}$,
où $\KuX=\KXr$, et~$\cG\subseteq \GLn(\KuX)$ le sous-groupe engendré
 par $\En(\KuX)$ et $\,\SL_2(\KuX)$. Alors  $f\sims{\cG} \tra{[\,1\;0\;\cdots\;0\,] }$.
\end{corollary}
\begin{proof}
Si $f_1=0$, on transforme facilement
par des \mlrs le vecteur $f$ en $\tra{[\,1\;0\;\cdots\;0\,]}$.
Sinon, un \cdv permet de transformer $f_1$ en un \pol pseudo \mon en $X_r$
(lemme~\ref{lemNoether}).
Nous pouvons donc supposer $f_1$ \mon en $X_r$,
nous appliquons le corolaire \ref{corth4SusQS} avec l'anneau
$\gA=\gK[X_1,\ldots,X_{r-1}]$,
et nous obtenons  $f\sims{\cG} f(X_1,\ldots,X_{r-1},0)$.
On conclut par \recu sur~$r$.
\end{proof}

On a bien obtenu le \thrf{thSuslinQS}, en fait avec une précision
intéressante sur le groupe~$\cG$.

%%%%%%%%%%%%%%%%%%%%%%       A LA VASERSTEIN          %%%%%%%%%%%%%%%%%%%%

\subsubsection*{Deuxième \dem}

Nous suivons maintenant de près une \dem de Vaserstein~\cite {Vaserstein2}
telle qu'elle est exposée dans~\cite{Lam06} mais en
utilisant des arguments constructifs.

De manière plus \gnle nous sommes intéressés par la possibilité de trouver dans la classe
d'\eqvc d'un vecteur défini sur $\AX$ un vecteur défini sur
$\gA$, en un sens convenable.

Nous utiliserons  le lemme suivant.
%--- Lemma{lemfij}----------------
\begin{lemma}
\label{lemfij}
Soit $\gA $ un anneau  et $f(X)=\tra{[\,f_1(X)\;\cdots\;f_n(X)\,]}$ un \vmd dans
$\AX^{n{\times}1}$, avec $f_1$ \mon de degré $\geq1$. \\
Alors, l'idéal $\fa=\rc(f_2)+\cdots +\rc(f_n)$  contient~$1$.
\end{lemma}
%--- end-lemma-----------------------------------------
%-----------------begin proof------------------
\begin{proof} 
On a: $1= u_1f_1$ dans $\gA\sur\fa$. Cette \egt dans l'anneau $(\gA\sur\fa)[X]$,
avec~$f_1$ \mon de degré $\geq1$ implique que $\gA\sur\fa$
est trivial (par \recu sur le degré formel de~$u_1$).
\end{proof}
%-----------------end proof------------------

\penalty-2500
%:     Theorem{th2HorrocksLocal}
\begin{theorem}\label{th2HorrocksLocal}\emph{(Petit \tho de Horrocks local)}\\
Soit un entier $n\ge 3$, $\gA $ un \alo \dcd et
 un \vmd dans $\AX^{n{\times}1}$: $f(X)=\tra{[\,f_1(X)\;\cdots\;f_n(X)\,]}$, avec
$f_1$ \monz. Alors
$$
f(X)=\cmatrix{ f_1   \cr \vdots\cr \vdots \cr f_n}
\sims{\En(\AX)}
\cmatrix{ 1 \cr 0\cr \vdots \cr 0}\sims{\En(\gA)}
\cmatrix{ f_1(0)   \cr \vdots\cr \vdots \cr f_n(0)}
.
$$
 \end{theorem}

%-----------------begin proof------------------
\begin{proof}
Soit  $d$ le degré de $f_1$. Par \mlrs de lignes,
on ramène les \pols $f_2$, \ldots, $f_n$ à être de degrés $< d$. Notons $f_{i,j}$
le \coe de~$X^j$ dans $f_i$. Le vecteur
$\tra{[\,f_1(X)\;\cdots\;f_n(X)\,]}$ reste \umdz. Si $d=0$, c'est terminé. Sinon
vu le lemme \ref{lemfij} et puisque l'anneau est local, l'un des  $f_{i,j}$
pour  $i\in\lrb{2..n}$ est une unité. Supposons par exemple que $f_{2,k}$
est \ivz. On va voir que l'on peut trouver deux \pols $v_1$ et
$v_2$ tels que le \pol $g_2=v_1f_1+v_2f_2$ soit \mon de degré $d-1$.
Si $k=d-1$, cela marche avec $v_1=0$ et $v_2$ constant.
Si $k<d-1$, considérons la disjonction suivante

\snic {   f_{2,d-1}\in\Ati\; \vee\; f_{2,d-1}\in
\Rad(\gA ).}

%\sni
Dans le premier cas, on est ramené à  $k=d-1$. Dans le deuxième
cas
le  \pol $q_2=Xf_2-f_{2,d-1}f_1$ est  de degré $\le d-1$ et vérifie:
$q_{2,k+1}$ est une unité. On a gagné un cran: il suffit 
d'itérer le processus.
\\ 
Nous avons donc maintenant $g_2=v_1f_1+v_2f_2$  de degré $d-1$ et
\monz. On peut donc diviser $f_3$ par $g_2$ et l'on obtient $g_3=f_3-g_2q$
de degré $<d-1$ ($q\in\gA $), donc le \pol

\snic{h_1=g_2+g_3=f_3+g_2(1-q)=f_3+(1-q)v_1f_1+(1-q)v_2f_2}

%\sni
 est \mon de degré $d-1$.
Ainsi, par une  \mlr de lignes on a pu remplacer
%:HHH plus beau avec \tra{[\,f_1\;f_2\;f_3\,]} et \tra{[\,f_1\;f_2\;h_1\,]
$\tra{[\,f_1\;f_2\;f_3\,]}$ par $\tra{[\,f_1\;f_2\;h_1\,]}$ avec $h_1$ \mon
de degré~$d-1$.
Nous pouvons donc par une suite de  \mlrs de lignes ramener 
$\tra{[\,f_1(X)\;\dots\;f_n(X)\,]}$, avec $ f_1 $ \mon
de degré $ d $,
 à

\snic{\tra{[\,h_1(X)\;\dots\;h_n(X)\,]} $ avec $ h_1 $ \mon
de degré $ d-1.}

%\sni
On obtient le résultat souhaité par \recu sur~$d$.
\end{proof}
%-----------------end proof------------------

 {\bf Terminologie.} Nous considérons un \sys de \pols formels
$(f_i)$ avec $\deg f_i=d_i$.
On appelle alors \gui{\id de tête du \sys $(f_i)$} l'idéal des
 \coes formellement dominants des~$f_i$.

\penalty-2500
%:     Theorem{th2HorrocksGlobal}
\begin{theorem}\label{th2HorrocksGlobal}
\emph{(Petit \tho de Horrocks global)}\\
Soit un entier $n\ge 2$, $\gA $ un anneau et
$f\in \AX^{n{\times}1}%=\tra{(f_1(X),\ldots,f_n(X))}
$ un \vmdz.
On suppose que l'\id de tête des $f_i$ contient $1$.
Alors
$$
f(X)=\cmatrix{ f_1   \cr \vdots \cr f_n}
\sims{\GLn(\AX)}
\cmatrix{ f_1(0) \cr \vdots \cr f_n(0)}=f(0).
$$
\end{theorem}
%
%NB: 
%Si $f_i$ est \mon on se ramène à ce cas
%par des \mlrs de lignes
%qui font descendre en dessous de $\deg f_i$ les degrés des autres
%\coos du vecteur.
\begin{proof}
Le cas $n=2$ est à part: si $u_1f_1+u_2f_2=1$, l'\egt 

\snic{\cmatrix{u_1&u_2\cr -f_2&f_1}
\cmatrix{f_1\cr f_2}=\cmatrix{1\cr 0}}

%\sni
 donne la matrice cherchée, dans $\SL_2(\AX)$.
\\ 
Pour $n\geq3$, nous appliquons la machinerie \lgbe
de base \paref{MethodeIdeps} avec la \dem
\cov du \thrf{th2HorrocksLocal}.\imlb Nous obtenons une
famille finie de \mocoz, $(S_i)_{i\in J}$  dans $\gA$, de
telle sorte que pour chaque $i$ on a
$
f(X)\sims{\En(\gA_{S_i}[X])} f(0).
$
On conclut avec le \rcm de Vaserstein
pour les \eqvcs de matrices à gauche (point~\emph{2} du \plgrf{thPatchV}).
\end{proof}

 {\bf Conclusion.}
On vient d'obtenir une variante (légèrement plus faible) du
corolaire \ref{corth4SusQS}. Et ceci donne la \dem du \tho
de Suslin~\ref{thSuslinQS}
de la même manière que dans la première solution.

\mni\comm Le petit \tho de Horrocks global peut aussi être obtenu
comme conséquence du \gui{grand} \tho de Horrocks global \vref{thHor}.
\\
On pose $P=\Ker\tra{f(X)}$. En localisant en $f_1$, $P$ devient libre.
\\
Le \tho de Horrocks global nous dit que $P$ est libre, ce qui signifie
que~\hbox{$f(X)\sim \tra{[\,1\;0\;\cdots\;0\,]}$} sur~$\GLn(\AX)$.
\eoe

%%%%%%%%%%%%%%%%%%%%%%       A LA RAO           %%%%%%%%%%%%%%%%%%%%

\subsubsection*{Troisième \dem}

Nous suivons maintenant de près une \dem de Rao.
Nous n'aurons cette fois-ci pas besoin de \recu sur le nombre de variables
pour aboutir au \tho de Suslin.

%:     Lemma{lemUMD}
\begin{lemma}\label{lemUMD}
On considère un vecteur
$x=(x_1,\ldots ,x_n)\in\Ae n$ et $s\in\gA$. Si $x$ est \umd sur
$\aqo{\gA}{s}\!$ et sur $\gA[1/s]$, il est \umdz.
 \end{lemma}
\begin{proof}
Posons $\fa=\gen{x_1,\ldots,x_n}$. On a $s^r\in\fa$ (pour un certain
$r$) et~\hbox{$1-as\in\fa$} (pour un certain $a$). On écrit $1=a^rs^r+(1-as)(1+as+\cdots) \in \fa$.
\end{proof}
%

%:     Lemma{lem1Rao}
\begin{lemma}\label{lem1Rao}
Soient un entier $n\ge 2$, $\gA $ un anneau, et
 un \vmd dans $\AX^{n{\times}1}$: $f=\tra{[\,f_1(X)\;\cdots\;f_n(X)\,]}$.
Pour chaque $f_i$ de degré formel $d_i$, nous notons $f_i\sta$ le \emph{\pol
formel réciproque} $X^{d_i}f_i(1/X)$. Nous notons $f\sta(X)=\tra{[\,f_1\sta(X),\;\cdots\;f_n\sta(X)\,]}$.\\
Si $f\sta(0)$ est \umdz, il en va de même pour~$f\sta$.
\end{lemma}
\begin{proof}
D'après le lemme \ref{lemUMD}, il suffit de vérifier que $f\sta(0)$
est \umd (c'est vrai par hypothèse) et que $f\sta$ est \umd sur $\gA[X,1/X]$,
ce qui vient de l'\egt $\sum_iu_i(1/X)X^{-d_i}f_i\sta=1$
 (où $\sum_iu_if_i=1$ dans~$\AX$).
\end{proof}
%

%:     Theorem{th1Rao}
\begin{theorem}\label{th1Rao}\emph{(\Tho de Rao, \cite{Rao85b})}\\
Soit un entier $n\ge 2$, $\gA $ un anneau, et
$f=\tra{[\,f_1(X)\;\cdots\;f_n(X)\,]}$ un \vmd dans $\AX^{n{\times}1}$,
avec $1$ dans l'\id de tête des $f_i$. Alors:
$$
f\sims{\GLn(\AX)} f(0)\sims{\GLn(\gA)}f\sta(0)\sims{\GLn(\AX)} f\sta.
$$
Si en outre l'un des $f_i$ est \monz, on a $f\sims{\GLn(\AX)}
\tra[\,1\;0\;\cdots\;0\,]$.
\end{theorem} 
\begin{proof}
On sait que $f\sim f(0)$ par le petit \tho de Horrocks global,
on en déduit $f\sim f(1)$. En outre,  $f\sta(0)$
est \umd donc $f\sta$ est \umd
(d'après le lemme \ref{lem1Rao}).
Par ailleurs, $1$ est dans l'\id de tête des $f_i\sta$,
ce qui permet d'appliquer à $f\sta$
le petit \tho de Horrocks global.\\
On conclut:   $f\sim f(0)\sim f(1)= f\sta(1)\sim f\sta$.
\end{proof}

\comm
 Le même résultat est valable en remplaçant $\GLn$
par $\En$, mais la \dem est nettement plus délicate (voir le
\thrf{th2Rao}).
\eoe

\mni {\bf Conclusion.} On obtient alors le \tho de Suslin \rref{thSuslinQS}
comme suit.
On prend pour $\gA $ l'anneau $\gK[X_1,\ldots,X_{r-1}]$
et l'on fait un \cdv qui rend l'un des \pols pseudo \monz.
\\ 
Ainsi, 
\begin{itemize}
\item  d'une part, la solution est beaucoup
plus \gui{efficace} que dans les deux pre\-mières \dems puisqu'il n'y a plus maintenant de \recu sur le nombre de variables,
\item  et d'autre part, le \tho est beaucoup plus \gnlz.
\end{itemize}

%\vspace{15mm} 
% ~

%%%%%%%%%%%%%%%%%%%%%%%%%%%%%%%%%%%%%%%%%%%%%%%%%%%%%%%%%%%%%%%%%%%%
%%%%%%%%%%%%%%%%%%%%%%%%%%%%%%%%%%%%%%%%%%%%%%%%%%%%%%%%%%%%%%%%%%%%
%%%%%%%%%%%%%%%%%%%%%%                               %%%%%%%%%%%%%%%%%%%%%
%%%%%%%%%%%%%%%%%%%%%%       Bass.Valuation          %%%%%%%%%%%%%%%%%%%%%
%%%%%%%%%%%%%%%%%%%%%%                               %%%%%%%%%%%%%%%%%%%%%
%%%%%%%%%%%%%%%%%%%%%%%%%%%%%%%%%%%%%%%%%%%%%%%%%%%%%%%%%%%%%%%%%%%%
%%%%%%%%%%%%%%%%%%%%%%%%%%%%%%%%%%%%%%%%%%%%%%%%%%%%%%%%%%%%%%%%%%%%
%%%%%%%%%%%%%%%%%%%%%%%%%%%%%%%%%%%%%%%%%%%%%%%%%%%%%%%%%%%%%%%%%%%%

%\newpage
%\penalty-10000
%---- Section{sec.Etendus.Valuation}--
\section[Modules \pros étendus depuis les \anarsz]
{Modules \pros étendus depuis les \advs ou \aris}  \label{sec.Etendus.Valuation} 
%------------------

Rappelons qu'un \adv est un anneau réduit dans lequel on a, pour tous $a,\,b$:
$a$ divise $b$ ou $b$ divise $a$. C'est un anneau normal, local et \sdzz.
\perso{Il est intègre \ssi il est \cohz, mais cela semble sans importance ici.
}

Nous commençons par un résultat utile concernant les \advs et la \ddk
(on peut aussi consulter l'exercice~\ref{exoPrufNagata}).

%:     Lemma{lemV-V(X)}
\begin{lemma}\label{lemV-V(X)} 
Si  $\gA$ est un \advz, alors $\gA(X)$ \egmtz.
Si~$\gA$ est un \adv de \ddk finie, alors $\gA(X)$ a même \ddkz.
\end{lemma}
\begin{proof}
Si $\gA$ est un \adv tout $f\in\AX$ s'écrit sous forme $f=ag$ avec $a\in\gA$ et $g\in\AX$ qui admet un \coe égal à 1.
En particulier,~$g$ est \iv dans $\gA(X)$.
Si  $F_1=a_1g_1/u_1$ et  $F_2=a_2g_2/u_2$ sont deux \elts
arbitraires de $\gA(X)$ (avec $a_i\in\gA$ et $g_i,u_i$ primitifs dans
$\AX$), alors $F_1$ divise $F_2$ dans $\gA(X)$
\ssi $a_1$ divise $a_2$ dans $\gA$. Donc \gui{la \dve est identique
dans $\gA$ et $\gA(X)$} et $\gA(X)$ est un \advz.
En outre, puisque les \itfs sont principaux, l'\homo canonique $\ZarA\to\Zar\gA(X)$ est un \iso de \trdis (NB: ce sont des ensembles totalement ordonnés),
ce qui implique que la \ddk est la même. 
%:HHH argument trop comliqué le lecteur est censé se rappeler 
%que la dimension d'un anneau est celle de son Zar
%(par exemple en application du lemme \ref{lemSeqSingTD}
%ou de la section~\ref{secDDKTRDIS}).
\end{proof}
%

%%%%%%%%%%%%%
%:--- subsec{En une variable}
\subsec{En une variable} \label{subsec.Bass.Valuation} 

Ce paragraphe est consacré pour l'essentiel à la \dem
\cov du \tho de Bass suivant.

%:HHH j'ai numerote ce theoreme

%:     Theorem{thBass.Valuation}
\begin{theorem}\label{thBass.Valuation} 
 Si $\gV$ est un \adv de
\ddk finie,  tout $\VX$-\mptf est libre.
\end{theorem}
%--------- fin theorem ---------------------------------------------- 

Nous démontrerons en fait des variantes un peu plus fortes:
on peut se débarrasser de l'hypothèse sur la \ddkz,
et l'on a une version avec des \anarsz.

Nous commençons par un exemple simple.
%    --- SUBsection{subsecZXlibre}------
\subsubsection*{Un exemple simple}
\label{subsecZXlibre}
%-----------------------------------------

%:     Proposition{propZXlibre}
\begin{proposition}\label{propZXlibre}
Tout \mptf sur $\ZZX$ est libre.
\end{proposition}
\begin{proof}
Soit $M$ un $\ZZX$-\mptfz.
Notons tout d'abord que si $M$ est de rang $1$, il est libre parce que
$\ZZX$ est un anneau à pgcd (lemme~\ref{lemPicGcd}). \\
Supposons maintenant que $M$ est de rang $r>1$.
Si nous étendons les scalaires à $\QQX$,
le module devient libre.
Il existe donc un entier $d>0$ tel que
$M$ devient libre sur $\ZZ[1/d][X]$. 
Si $d=1,$ il n'y a rien à faire.
Sinon, soient~\hbox{$p_1$, \ldots, $p_k$}  les facteurs
premiers de~$d$. \\
Les \mos $d^{\NN}$, $1+p_1\ZZ$, \ldots,  $1+p_k\ZZ$
sont \com (exemple fondamental \paref{explfonda}). Il nous suffit donc
de montrer que les modules $M_{1+p_i\ZZ}$ sont libres (donc étendus),
car alors le \tho de \rcm de Quillen\index{Quillen!recollement de ---}  %(\thrf{thPatchQ})
implique que $M$ est étendu depuis $\ZZ$, donc libre. \\
Notons $p$ l'un quelconque
des $p_i$.  Puisque $\ZZ_{1+p\ZZ}[X]$ est $2$-stable
(lemme ci-après), en application du splitting off de Serre
(\thrf{thSerre}), on obtient
$M_{1+p\ZZ}\simeq \ZZ_{1+p\ZZ}[X]^{r-1}\oplus N$, avec $N$ 
un $\ZZ_{1+p\ZZ}[X]$-\mrcz~$1$.
D'après la remarque initiale (qui s'applique en remplaçant~$\ZZ$ par $\ZZ_{1+p\ZZ}$), $N$ est libre, donc $M$ est libre.
\end{proof}
%
%:     Lemma{lemZ1+pZ}
\begin{lemma}\label{lemZ1+pZ} 
L'anneau $\ZZ_{1+p\ZZ}[X]$ est $2$-stable.
\end{lemma}
\begin{proof}
On considère la partition de $\Spec(\ZZ_{1+p\ZZ}[X])$ attachée à $\so{p}$:
plus \prmtz,
l'anneau $\ZZ_{1+p\ZZ}[X]$ est remplacé par les deux anneaux

\snic{\ZZ_{1+p\ZZ}[X][1/p]\simeq\QQX$ et
$\aqo{(\ZZ_{1+p\ZZ}[X])}{p}\simeq \FF_{p}[X],}

%\sni
qui sont de \ddk $1$.\\
Le \thrf{thPartitionSpec} nous dit alors que $\ZZ_{1+p\ZZ}[X]$  est $2$-stable.
\end{proof}

\rem
En fait le recours aux facteurs premiers de $d$, bien qu'intuitivement
naturel, introduit une complication inutile.
En effet, les \mosz~$d^{\NN}$ et $1+d\ZZ$ étant \comz, il suffit de
démontrer que $M_{1+d\ZZ}$ est libre. Comme $\ZZ_{1+d\ZZ}[X]$
est un anneau à pgcd, le raisonnement précédent s'applique
si l'on sait montrer que  $\ZZ_{1+d\ZZ}[X]$ est $2$-stable. Or
la preuve du lemme \ref{lemZ1+pZ} fonctionne en remplaçant $p$ par $d$,
 car $\ZZ_{1+d\ZZ}[X][1/d]\simeq\QQX$, et~$\aqo{\ZZ_{1+d\ZZ}[X]}{d}\simeq (\aqo{\ZZ}{d}\!)[X]$ qui sont de \ddk $1$
 ($\aqo{\ZZ}{d}$ est \zedz).
\eoe

%    --- SUBsection{subsec2ZXlibre}------
\subsubsection*{Un exemple plus élaboré}
\label{subsec2ZXlibre}
%-----------------------------------------

Au vu de la remarque précédente nous laissons \alec la \dem
de la \gnn qui suit.

%:     Proposition{propBezIntlibre}
\begin{proposition}\label{propBezIntXlibre}
Soient $\gA$ un anneau intègre de \ddkz~$\leq 1$, un \elt
  $d$ de $\Reg(\gA)$, et $M$ un
$\AX$-\mptfz.
\begin{enumerate}
\item $\gA_{1+d\gA}[1/d]=\Frac\gA$ est \zedz.
\item $\aqo{\gA_{1+d\gA}}{d}\simeq\aqo{\gA}{d}$ est \zedz.
\item $\gA_{1+d\gA}[X]$ est $2$-stable.
\item
\begin{enumerate}
\item  Si $\gA$ est un anneau de Bézout,  $M$ est libre.
\item  Si $\gA$ est seminormal,  $M$ est étendu depuis $\gA$.
\end{enumerate}
\end{enumerate}
\end{proposition}

%    --- SUBsection{subsec1VXlibre}------
\subsubsection*{Un exemple en \ddk finie $>0$}
\label{subsec1VXlibre}
%-----------------------------------------

Soit $\gV$ un \adv intègre avec des \elts $a_1$, \ldots, $a_k$. On suppose:
$$\preskip-.1em \postskip.4em 
\rD_\gV(a_1)< \rD_\gV(a_2)< \cdots < \rD_\gV(a_k). 
$$
La partition
en constructibles de $\Spec\gV$ associée à cette famille contient seulement $k+1$ \eltsz:
$$\preskip.4em \postskip.4em 
 \rD_\gV(a_1),\; \rD_\gV(a_2)\setminus \rD_\gV(a_1),\; \dots,\;
\rD_\gV(a_k)\setminus \rD_\gV(a_{k-1}),\;
\rD_\gV(1)\setminus \rD_\gV(a_{k}),
$$
qui correspondent aux anneaux

\snic{\gV[1/a_1]$, $(\aqo{\gV}{a_1}\!)[1/a_2]$, \ldots, $(\aqo{\gV}{a_{k-1}})[1/a_k]$
et $\aqo{\gV}{a_k}.}

%\sni
Supposons maintenant que ces anneaux sont tous \emph{\zedsz}. Alors, on a pareillement une partition en $k+1$ constructibles de $\Spec\VX$ et les
anneaux correspondants

\snic{\gV[1/a_1][X]$, $(\aqo{\gV}{a_1})[1/a_2][X]$, \ldots,
$(\aqo{\gV}{a_{k-1}})[1/a_k][X]$ et $(\aqo{\gV}{a_k})[X]}

%\sni
 sont tous de \ddkz~$\leq 1$.
Le \tho \ref{thPartitionSpec} nous dit alors que~$\gV[X]$ est $2$-stable.
Donc si $M$ est un $\VX$-\mrc $r$, par le splitting off de Serre (version $\Cdim$), on obtient 
$M\simeq \VX^{r-1}\oplus N$, avec~$N$ de rang constant $1$.

Si $\gV$ est en plus un anneau seminormal (resp. un anneau à pgcd),
alors  $N$
est étendu depuis $\gV$ (resp. alors $N$ est libre), donc $M$ est étendu depuis $\gV$
(resp.  $M$ est libre).

\smallskip Ainsi le résultat \gui{$\VX$ est $2$-stable} est vérifié  lorsque $\gV$ est un domaine de
valuation de \ddkz~$k$ pour lequel on a une connaissance suffisamment
précise du groupe de valuation: on connaît $a_1$, \dots, $a_k$ tels que
$\rD_\gV(0)<\rD_\gV(a_1)< \rD_\gV(a_2)< \cdots < \rD_\gV(a_k)< \rD_\gV(1)$.

En \clama (avec le principe du tiers exclu mais sans utiliser
les \ideps ni l'axiome du choix) on obtient donc déjà le
\tho de Bass souhaité pour les \advs de \ddk finie.

Cependant le résultat n'est pas de nature \algq
si l'on ne sait pas calculer des \elts $a_i$ convenables.

Cette difficulté va être contournée de manière dynamique.

%    --- SUBsection{subsec2VXlibre}------
\subsubsection*{\Demo \cov du \tho de Bass}
\label{subsec2VXlibre}
%-----------------------------------------

Le plus important est d'établir le \tho suivant.

%:     Th{thVX2stab}
\begin{theorem}\label{thVX2stab}%~\\
Si $\gV$ est un \advz, %l'anneau
$\VX$ est $2$-stable.
\end{theorem}

On commence par le lemme suivant (la \dem du \tho est reportée
\paref{DemothVX2stab}).

%:     Lemma{lemVX3stab}
\begin{lemma}\label{lemVX3stab}
Soient $\gV$  un \adv et $\gV'$ le sous-\adv de $\gV$
engendré par une famille finie d'\elts de $\gV$.
Alors, $\gV'[X]$ est $2$-stable.
\end{lemma}
\begin{proof}
Notons $\gV_1$ le sous-anneau de $\gV$ engendré par la famille finie.
Notons
$$
\gV'=\sotq{c/b}{c,b\in\gV_1,\;b\;\mathrm{r\E egulier}\;\mathrm{divise}\; c\; \mathrm{dans}\; \gV}\subseteq\Frac(\gV_1).
$$
On voit facilement que $\gV'$ est un \advz.
On sait que $\gV\!_1$ est de \ddk finie
(lemme~\ref{lemahbonvraiment}), disons $\Kdim(\gV_1)\leq m$. Montrons que l'on a aussi
$\Kdim(\gV')\leq m.$
Le \thref{thKdimTDTO} nous dit le contrat à remplir pour cela. On doit considérer une suite $(y_0,\dots,y_m)$ avec
$$\preskip.4em \postskip.4em 
\rD_{\gV'}(y_0)\leq \rD_{\gV'}(y_1)\leq \dots\leq  \rD_{\gV'}(y_{m}).
$$
On peut écrire $y_k=x_k/b$ pour un même dénominateur $b\in\Reg(\gV\!_1)$ et des  $x_k\in \gV\!_1$.
On a maintenant
$$\preskip.4em \postskip.4em 
\rD_{\gV'}(x_0)\leq \rD_{\gV'}(x_1)\leq \dots\leq  \rD_{\gV'}(x_{m}),
$$
Introduisons une suite $(a_0,\dots,a_m)$ \cop de celle des $x_i$ dans $\gV\!_1$.
A fortiori elle est \cop dans $\gV'$. Puisque $\Zar\gV'$ est totalement ordonné, le lemme \ref{lemSeqSingTD} nous dit que $\rD_{\gV'}(x_0)=\rD_{\gV'}(0)$, ou $\rD_{\gV'}(x_{m})=\rD_{\gV'}(1)$, ou $\rD_{\gV'}(x_i)=\rD_{\gV'}(x_{i+1})$ pour un $i\in\lrb{0..m-1}$. Et l'on en déduit la même chose 
pour les  $\rD_{\gV'}(y_i)$: le contrat est rempli.

Soient $\ell_1$, $\ell_2$ et $a$ dans $\gV'[X]$.
Nous posons $L=(\ell_1,\ell_2)$ et $Q=(q_1,q_2)$.
Nous cherchons $q_1$, $q_2\in\gV'[X]$
qui vérifient $\rD_{\gV'[X]}(a,L)=\rD_{\gV'[X]}(L+aQ)$.
\\
Si $\gV'$ était un \cdiz, on disposerait d'un
\algo pour calculer $Q$ à partir de $L$.
En exécutant cet \algoz, nous utiliserions  le test \gui{$y=0$ ou $y$ \ivz?}
pour des \elts $y\in\gV_1$ qui se présentent au cours du calcul
(en effet, dans le cas où $\gV'$ est un corps discret,
un $y/z$ dans $\gV'$ est nul si $y$ est nul, \iv si $y$
est \ivz, $z$ ayant été déjà certifié \ivz).\\
Nous pouvons transformer l'\algo de façon  dynamique en
remplaçant chaque test \gui{$y=0$ ou $y$ \ivz?}
par le scindage de \gui{l'anneau $\gA$ en cours},
 qui donne les  deux anneaux $\gA[1/y]$ et
${\gA}\sur{\DA(y)}$.\\
Au départ $\gA=\gV'$.
Comme dans $\gV'$ les \elts sont comparables pour la \dvez, 
tous les anneaux introduits peuvent être ramenés à la forme
standard ${\gV'}\sur{\rD_{\gV'}( y_i)}\![1/y_{i-1}]$ ($i\in\lrb{2..k}$)
pour une famille finie $(y_i)_{i\in\lrbk}$ de $\gV_1$,
avec $y_{i-1}$ divise $y_i$ dans $\gV'$ pour $i\in\lrb{2..k}$.
\\
Ici nous pourrions avoir l'impression d'avoir gagné dans la mesure
où nous pourrions dire:
nous appliquons maintenant le lemme \ref{lemPartitionSpec}.
\\
Mais en lisant la \dem
de ce lemme, nous voyons que lors d'un scindage  $\gB\mapsto(\gB[1/b],\aqo{\gB}{b}\!)$, d'abord les données $L$
et $a$ produisent un~$Q$ pour  $\aqo{\gB}{b}$, puis $L+aQ$ et $ab$ produisent
un $R$ pour $\gB[1/b]$, le résultat final étant que $Q+bR$ convient pour $L$ et $a$
dans $\gB$.\\
Ainsi la dynamique de notre \algo transformé doit être
mieux contrôlée\footnote{Sinon, le lemme pourrait en fait
être démontré sans aucune hypothèse sur $\gV$.}.
Ce qui nous sauve la mise, c'est que dans notre utilisation dynamique du lemme \ref{lemPartitionSpec}, les calculs qui démarrent avec $L$ et $a$ restent entièrement
dans $\gV'\subseteq\Frac(\gV_1)$. En conséquence, nous pouvons être certains de ne pas entrer dans une boucle infinie où le nombre d'anneaux ${\gV'}\sur{\rD_{\gV'}( y_i)}\![1/y_{i-1}]$
croîtrait indéfiniment, ce qui empêcherait la terminaison
de l'\algoz.
En effet, puis\-que $\gV'$ est de \ddk  $\leq m$, on dispose
d'une procédure qui, étant donnée une famille
finie $(y_i)$ comme ci-dessus,
permet de raccourcir la famille des $\rD_{\gV'}(y_i)$
à au plus $m$ \elts gr\^ace au lemme~\ref{lemSeqSingTD}.\\
Considérons en effet les $m+1$ premiers termes consécutifs
dans la suite des~$y_i$, nous savons  que l'une des trois situations suivantes se produit
 \begin{itemize}
\item $\rD_{\gV'}(y_{1})=0$, auquel cas l'anneau
${\gV'}\sur{\rD_{\gV'}( y_2)}\![1/y_{1}]$ est trivial et la liste est raccourcie en supprimant $y_1$,
\item $\rD_{\gV'}(y_{m+1})=1$, auquel cas l'anneau
${\gV'}\sur{\rD_{\gV'}( y_{m+1})}\![1/y_{m}]$ est trivial et la liste est raccourcie en supprimant $y_{m+1}$,
\item pour un $i\in\lrb{2,m+1}$, on a l'\egt $\rD_{\gV'}(y_{i-1})=\rD_{\gV'}(y_{i})$, auquel cas l'anneau
${\gV'}\sur{\rD_{\gV'}( y_{i})}\![1/y_{i-1}]$ est trivial et la liste est raccourcie en supprimant $y_{i}$.
\end{itemize}
\vspace{-1em}
\end{proof}

\rem  Ainsi, une fois $\gV_1$ fixé, l'anneau $\gV'$ se comporte,
pour ce qui concerne la $2$-stabilité de $\gV'[X]$ comme l'anneau
de \ddk  \gui{finie $>0$ mais entièrement contrôlée}
qui était donné dans
le paragraphe précédent: la suite des $y_i$, limitée
à $m$ termes, se comporte comme la suite des $a_i$
du paragraphe précédent,
à ceci près que les $y_i$ sont produits de façon dynamique
par l'exécution de l'\algo alors que les $a_i$ étaient donnés au départ.
\eoe

\begin{Proof}{\Demo du \tho \ref{thVX2stab}. } \label{DemothVX2stab}
Soient $\ell_1,\ell_2$ et $a$ dans $\VX$.
Nous cherchons $q_1,q_2\in\gV[X]$ vérifiant
$\rD_{\VX}(a,L)=\rD_{\VX}(L+aQ)$
(avec $L=(\ell_1,\ell_2)$ et $Q=(q_1,q_2)$).
On applique le lemme \ref{lemVX3stab} avec la famille finie constituée
par les \coes de $\ell_1$, $\ell_2$ et $a$. On trouve $q_1,q_2$ dans
$\gV'[X]\subseteq\VX$.
\end{Proof}
%

%:     Theorem{thBassValu}
\begin{theorem}\label{thBassValu}\emph{(Bass-Simis-Vasconcelos)}
Si $\gV$ est un \advz,  tout $\VX$-module \ptf est libre.
\end{theorem}
\begin{proof}
Soit  $M$ un \mptf sur $\VX$. Puisque $\VX$ est connexe, $M$ a un rang constant $r\in\NN$.
Puisque $\VX$ est $2$-stable, le splitting off de Serre nous donne que $M\simeq \VX^{r-1}\oplus N$, où $N$
est un $\VX$-\mrc $1$.
Il reste à montrer que $N\simeq\VX$.\\
Si $\gV$ est intègre nous terminons comme ceci:
puisque $\VX$ est un anneau à pgcd, $N\simeq \VX$.
En \gnl nous pouvons dire: $\gV$ est normal, donc tout \mrc 1 sur
$\VX$ est étendu depuis $\gV$.
Or $\gV$ est local,
en conclusion $N$ est libre sur $\VX$.
\end{proof}
%

%    --- SUBsection{subsec4VXlibre}------
\subsubsection*{Le cas des \anarsz}
\label{subsec4VXlibre}
%-----------------------------------------
%:     Theorem{thBassAri}
\begin{theorem}\label{thBassAri}\emph{(Bass-Simis-Vasconcelos)}
Si $\gA$ est un \anarz,  tout $\AX$-module \ptf est étendu depuis $\gA$.
\end{theorem}
\begin{proof}
Tout d'abord, puisque $\GKO(\gA)=\GKO(\Ared)$ et $\AX\red=\Ared[X]$, il
suffit de faire la \dem dans le cas réduit, \cad pour les \adpsz.
\\
On considère un $\AX$-\mptf $M$.\\
En \clama on appliquerait le \tho de \rca  de Quillen:
un \mptf sur $\AX$ est étendu parce qu'il est étendu si on localise
en un \idep arbitraire de $\gA$ (l'anneau devient un \advz).
\\
En \comaz, on relit la \prco donnée dans le cas local (pour le \thrf{thBassValu}) en appliquant la machinerie \lgbe
de base.\imlb
\\
Précisément, supposons que dans le cas local
(i.e., pour un \advz) on utilise
la disjonction \gui{$a$ divise $b$ ou $b$ divise $a$}. Puisque l'on  est
avec un \adpz, on connaît $u$, $v$, $s$, $t$ tels que \hbox{$s+t=1$, $sa=ub$} et~$tb=va$.
Si $\gB$ est l'anneau \gui{en cours}, on considère les deux \lons \come
$\gB[1/s]$ et $\gB[1/t]$. Dans la première, $a$ divise $b$, et dans la seconde,~$b$ divise $a$.
\\
En fin de compte on obtient une famille finie $(S_i)$ de \moco de $\gA$
telle qu'après \lon en l'un quelconque des $S_i$, le module~$M$ devient libre,
donc étendu. On conclut avec le \rcm de Quillen (\plgc \ref{thPatchQ}).
\end{proof}

\rems 1) On n'a pas eu besoin de supposer que l'\adv était \dcd pour faire
fonctionner la \dem \cov des \thos \ref{thVX2stab}, \ref{thBassValu}
et \ref{thBassAri}. Cela se traduit notamment par le fait que dans la
dernière \demz, les \moco sont basés sur la disjonction (dans un \aloz)
\gui{$s$ ou $1-s$ est \ivz} et sont directement donnés par des \ecoz.

 2)
Dans ce type de passage du local au global, pour être certain
que l'\algo termine, il faut s'assurer que la version donnée
dans le cas local est \gui{uniforme}, cela veut dire que son exécution
se fait en un nombre d'étapes qui est borné par une fonction des paramètres
discrets de l'entrée: la taille de la matrice et les degrés de ses
\coesz. C'est bien le cas ici, modulo la preuve du lemme
\ref{lemahbonvraiment}. Notons que le fait que l'\algo dans le cas local
n'utilise pas de test d'\egt à $0$ nous simplifie beaucoup la vie
pour apprécier la validité de sa mise en {\oe}uvre dynamique
dans le passage du local au global.
\eoe

%-% ENTRE NOUS
\entrenous{La remarque 2) ci-dessus produit un certain mal de tête.
On se demande si c'est vraiment possible de traiter localement
un \adp comme un \adv quand il est question de
borner a priori la \ddk à la simple vue des données du \pbz.

Dans le cas de l'\adv on se restreint à un sous-\adv 
engendré par un nombre fini d'\elts ce qui implique
que la \ddk est contrôlée une fois pour toutes.

Dans le cas d'un \anar les \lons \come se font au prix d'utiliser
des nouveaux \elts qui proviennent de  \mlps et donc le premier
réflexe serait plutôt de dire que l'on  ne contrôle plus la \ddkz.

En fait lorsque l'anneau est intègre, on voit assez bien ce qui va
se produire. On a le corps $\gL$ engendré par les \coes des \coes d'une
\mpn du module. Ce corps est de degré de transcendance fini, donc tous ses \advs
ont une \ddk bornée une fois pour toutes.

Certes, les \lons \come utilisent des \elts en dehors de $\gL$. Mais les calculs 
qui transforment une matrice $H(X)$ (comme dans le lemme \ref{factEtPol})
en la matrice $H(0)$ ont lieu à chaque fois dans $\gL$
(avec un sous-\adv générique de $\gL$). On 
peut donc bien suivre pas à pas la \dem du cas \gui{\adv de \ddk finie}.

C'est seulement au moment du recollement à la Vaserstein que l'on  utilisera
les \elts dans $\gA$ mais hors de $\gL$ qui ont servi à produire
les \lons \comez.

Notons que le lemme \ref{lemahbonvraiment} (\ddk finie
pour un anneau \tfz) peut être remplacé
par un résultat plus fort, le \thrf{corthValDim} par exemple
(\ddk finie pour un suranneau dans son corps de fraction d'un anneau
intègre \tfz),
si cela nous rend les choses plus faciles à comprendre.

Dans la \dem du lemme \ref{lemVX3stab} 
je pense que l'on  peut bien comprendre ce qui se passe 
et voir que le recours
au lemme \ref{lemahbonvraiment} est vraiment suffisant,
avec des arguments du même type que ceux ci-dessus.

Je signale que le même \pb se pose pour le passage du local au
global dans Lequain-Simis.
}
%-% Fin ENTRENOUS

%%%%%%%%%%%%%%%%%%%%%%%%%%%%%%%%%%%%%
%%%%%%%%%%%%%%%%%%%%%%%%%%%%%%%%%%%%%%%
%---- Section{Lequain-Simis}--sec.LS-
%%%%%%%%%%%%%
%:--- subsection{En plusieurs variables}
\subsec{En plusieurs variables}
\label{sec.LS} 
%------------------
Ce paragraphe est consacré à la \dem \cov du \tho de Lequain-Simis suivant.

%:     Theorem{%:     Theorem   {LSValu}
\THo{(Lequain-Simis)}
{Si $\gA$ est un \anarz, tout \mptf sur $\AXr$ est étendu depuis
$\gA$.
}

\medskip 
{\bf Une comparaison dynamique entre les anneaux
$\gA(X)$ et $\ArX$} \label{dyn}

\smallskip 
Dans le \tho suivant, nous démontrons que pour un anneau $\gA$
\ddi$d$, l'anneau $\ArX$ se comporte dynamiquement comme l'anneau $\gA(X)$
ou comme une \lon d'un anneau  $\gA_S[X]$
 pour un \mo $S$   de $\gA$ avec $\Kdim\gA_S\leq d-1$.

%: ---Thm{compa}  R(X) et R<X>
\begin{theorem} \emph{(Comparaison dynamique de $\gA(X)$ avec $\ArX$)}
\label{compa} \\
Soit  un anneau $\gA$,
  $f=\sum_{j=0}^ma_jX^j \in \AX$ un \pol primitif,
et,   \linebreak 
pour $j\in\lrbm$,  $S_j=\cS_\gA^\rK(a_j)=a_j^{\NN}(1+a_j\gA)$ (le \mo bord de Krull de~$a_j$ dans $\gA$).\\
Alors, les \mos
$f^{\NN}$, $S_1$, \ldots, $S_m$ sont \com dans $\ArX$.\\
En particulier,  si $\Kdim \gA$ et $d\geq 0$, chaque anneau $\ArX_{S_j}$ est une \lon d'un $\gA_{S_j}[X]$ avec
$\Kdim\gA_{S_j}\leq d-1$.
\end{theorem}

\begin{proof}
Pour $x_1$, \ldots, $x_m \in \gA$ et $n$, $d_1$, \ldots, $d_m \in \NN$, on doit
montrer que les \elts suivants de $\gA[X]$

\snic{
f^n,\; a_m^{d_m}(1 - a_mx_m),\; \dots, \;
a_1^{d_1}(1 - a_1x_1),}

%\sni
engendrent un \id de $\gA[X]$ qui contient un \pol \monz. On
raisonne par \recu sur $m$; c'est évident pour $m = 0$ car $a_m =
a_0$ est \ivz. \\
Pour $m \ge 1$ et $j\in \lrb{1..m-1}$, posons 

\snic{a = a_m, \;\;x = x_m, \;\;d = d_m\;\;\hbox{et}\;\;
a'_j = a_j^{d_j}(1 - a_jx_j).}

%\sni
Considérons le
quotient
$\gB=\aqo{\gA}{a^d(1 - ax)}$; il faut montrer que la famille

\snic{\cF = (f^n, \; a'_{m-1},\; \dots,\; a'_1 )
}

%\sni
engendre un \id de $\BX$ qui contient un \poluz.  
\\
Puisque
$a^d(1 - ax) = 0$, $e = a^d x^d$ est un \idm et $\gen{e} = \gen
{a^d}$. \\
Notons $\gB_e \simeq \aqo{\gB}{1-e}$ et $\gB_{1-e} \simeq \aqo{\gB}{e}$.
Il suffit de montrer que $\gen {\cF}_{\gB_e[X]}$ et $\gen {\cF}_{\gB_{1-e}[X]}$ contiennent un
\poluz.
\\
 Dans $\gB_e[X]$, c'est \imd car $a$ est \ivz.  Dans $\gB_{1-e}[X]$, on a $a^d = 0$.  \'Ecrivons $f = aX^m +
r$ avec $r = \sum_{j=0}^{m-1} a_j X^j$.  
Dans $\gB$, pour tout exposant~$\delta$, les \elts de  $(a^\delta, a_{m-1}, \ldots, a_1, a_0)$ sont
comaximaux. Pour $\delta = d$, on en déduit que dans $\gB_{1-e}[X]$, le
\pol $r$ est primitif. Puisque $r = f - aX^m$ et $a^d = 0$, on a $r^d
\in \gen {f}$ donc $r^{dn} \in \gen {f^n}$. \\
On applique l'\hdr au \pol
$r \in \gB_{1-e}[X]$ de degré (formel) $m-1$:
l'\id $\gen{r^{dn}, a'_{m-1}, \cdots, a'_1}$ de $\gB_{1-e}[X]$ contient un \poluz; il en est donc de même de l'\id
 $\gen{f^n, a'_{m-1}, \cdots, a'_1}$.
\end{proof}

\rem Le \tho précédent semble tombé du ciel comme par miracle.
En fait il est le résultat d'une histoire un peu compliquée.
Dans l'article~\cite{ELY07} était démontré le \tho suivant,
en commençant
par le cas spécial d'un \alo \dcdz,
puis en généralisant à un anneau arbitraire au moyen de
la machinerie \lgbe de base.\imlb
\eoe

%\penalty-2500 
{\bf Théorème.} \emph{Soit  un anneau $\gA$ tel que $\Kdim \gA\leq d\in\NN$.
Soit  $f \in \AX$ un \pol primitif.
Il existe des \moco $V_1$, \ldots, $V_\ell $ de $\ArX$ tels que pour chaque $i\in \lrbl$, ou bien  $f$ est \iv dans  $\ArX_{V_i}$,  ou
bien~$\ArX_{V_i}$ est une \lon d'un $\gA_{S_i}[X]$ avec
$\Kdim\gA_{S_i}<d$.} 

En explicitant l'\algo contenu dans la preuve de ce \thoz, on a obtenu
le \thoz~\ref{compa}.
\eoe

%%%%%%%%%%%%%%%%%%%%%%%%%%%%%%%%%%%%%%%%%
%:HHH rajout rdb et label
\rdb \label{MachDynAxAx}
\medskip 
{\bf Machinerie dynamique avec $\ArX$ et $\gA(X)$}

\smallskip 
Le \tho précédent permet de mettre en {\oe}uvre une machinerie
dynamique d'un nouveau type.

On suppose que l'on  a établi un \tho pour les \advs de \ddk $\leq n$.
On veut le même \tho pour les anneaux~$\ArX$ lorsque $\gA$ est un \adv
de \ddk $\leq n$.

On suppose aussi que la propriété à démontrer  est stable par
\lon et qu'elle relève d'un \plgcz.

On fait une preuve par \recu sur la \ddkz. Lorsque la \ddk est nulle,
$\gA$ est un corps discret et l'on a $\ArX=\gA(X)$, qui est aussi un \cdiz,
donc le \tho s'applique.

Voyons le passage de la dimension $k$ à la dimension $k+1$ ($k<n$).
Remarquons que $\gA(X)$ est un \adv de même \ddk que~$\gA$
(lemme \ref{lemV-V(X)}). Nous supposons $\Kdim\gA\leq k+1$.
Nous avons une \prco du \tho pour les \advs de \ddk $\leq n$,
en particulier elle fonctionne pour $\gA(X)$.
Nous essayons de faire fonctionner cette \dem (i.e., cet \algoz)
avec $\ArX$ au lieu de~$\gA(X)$.
Cette preuve utilise le fait que dans $\gA(X)$ les \pols primitifs de~$\AX$
sont inversibles. Chaque fois que la preuve initiale utilise l'inverse
d'un tel \pol $f$, nous faisons appel au \thrf{compa}, qui remplace l'anneau
\gui{en cours} par des \lons \comez. Dans la première \lon le \pol $f$
a été inversé, et la preuve peut se poursuivre comme si $\ArX$ était
$\gA(X)$. Dans chacune des autres \lons on a remplacé $\ArX$ par un
localisé d'un anneau $\gA_{S_i}[X]$ avec $\Kdim \gA_{S_i}\leq k$,
et, \emph{si nous avons de la chance}, l'\hdr permet de conclure.

Au bout du compte on a prouvé le \tho pour des \lons \come de $\ArX$.
Puisque la conclusion relève d'un \plgcz, on a démontré le \tho
pour $\ArX$.

%\penalty-2500
\subsubsection*{Application au \tho de Maroscia et Brewer{\&}Costa}

La machinerie dynamique expliquée au paragraphe précédent
s'applique pour le premier des  résultats suivants.
\begin{enumerate}
  \item [(i)] \emph{Si $\gA$ est un \adv avec $\Kdim\gA \leq 1$,
alors $\ArX$ est un \adp avec $\Kdim\ArX \leq 1$.} \\
En effet, il suffit de vérifier la conclusion \lot
(ici, après \lon de $\ArX$ en des \mocoz).
Or le \tho \ref{compa} nous permet de scinder l'anneau $\ArX$
en composantes qui se comportent (pour le calcul à faire),
soit comme $\gA(X)$, soit comme un localisé d'un~$\KX$ où
$\gK$ est \zed réduit. Dans les deux cas on obtient un \adp
de \ddk $\leq1$.
  \item [(ii)] \emph{Si $\gA$ est un \adp avec $\Kdim\gA \leq 1$,
  alors $\ArX$ \egmtz.}\\
En effet, il suffit de vérifier la conclusion \lot
(ici, après \lon de $\gA$ en des \mocoz). On applique la machinerie \lgbe 
des  \anars à la \dem du point~(i): l'anneau $\gA$ subit des \lons \comez,
dans chacune desquelles il se comporte comme un \advz.\imla
\end{enumerate}

\rdb
Comme conséquence on obtient une version particulière du
\tho de Lequain-Simis
en utilisant l'induction de Quillen\index{Quillen!induction de ---} concrète (\thrf{thQUILIND}).

%:     Theorem{thMaBrCo}
\begin{theorem}\label{thMaBrCo} \emph{(Maroscia, Brewer{\&}Costa)}\\
 Si $\gA$ est un \anar avec $\Kdim\gA\leq1$,
tout \mptf sur $\AXr$ est étendu depuis~$\gA$.
\end{theorem}
\begin{proof}
Puisque $\Ared[\uX]=\AuX\red$ et $\GKO(\gB)=\GKO(\gB\red)$, il suffit de traiter le cas réduit, i.e., le cas des \adpsz.\\
Vérifions que la classe des \adps de \ddkz~\hbox{$\leq 1$}
satisfait les hypothèses du \thrf{thQUILIND}.
La première condition est le point~(ii) ci-avant que nous venons
de démontrer. \\
La deuxième condition est que les \mptfs sur $\AX$
sont étendus depuis $\gA$. C'est le \tho de Bass-Simis-Vasconcelos.
\end{proof}

\rdb
\subsubsection*{L'induction de Lequain-Simis} \label{LSinduction}

Dans le but de \gnr le \tho de Quillen-Suslin aux domaines de Prüfer,
et constatant que cette classe d'anneaux n'est pas stable pour
le passage de $\gA$ à $\ArX$,  Lequain et
Simis \cite{LS} ont trouvé un moyen habile pour contourner
la difficulté en démontrant un nouveau \tho d'induction
\gui{à la Quillen}, convenablement modifié.

\CMnewtheorem{ILSa}{Induction de Lequain-Simis abstraite}{\itshape}
%:     Induction de Lequain-Simis abstraite
\begin{ILSa}\label{LSindabs} ~\\
Soit $\cF$ une classe  d'anneaux qui satisfait les \prts suivantes.
%-----------------begin item------------------
\begin{description}
\item [{\rm (LS1)}]
Si $\gA \in \cF$, tout \idep $\fp$ non maximal de $\gA$
a une hauteur finie\footnote{I.e., $\Kdim(\gA_\fp)<\infty$.}.
\item [{\rm (LS2)}]
Si $\gA \in \cF $, alors $ \AX_{{\fp}[X]}  \in \cF $ pour tout \idep
${\fp}$  de $\gA$.
\item [{\rm (LS3)}]
Si $\gA \in \cF $, alors $\gA_{\fp}\in\cF$ pour tout \idep  ${\fp}$ de $\gA$.
\item [{\rm (LS4)}]
Si $\gA \in \cF $ est local,  tout \mptf
sur $\AX$ est~libre.
\end{description}
%-----------------end item------------------
Alors, pour tout $\gA\in\cF$ et tout $r\geq1$, 
tout \mptf  sur~$\AXr$ est étendu depuis $\gA$. 
\end{ILSa}
%--------- fin fact ---------------------------------------------- 

Notez ici  que si $\gA$ est local avec $\Rad\gA={\fm}$,
alors $\gA(X)=\AX_{{\fm}[X]}$.

Nous proposons une \gui{variation \covz} sur le thème de
l'induction de Lequain-Simis. Il s'agit d'une
application importante de notre comparaison dynamique
entre $\gA(X)$ et $\ArX$. Cette induction \cov \gui{à la Lequain-Simis}
est due à I.~Yengui.

%: --- Constructive induction theorem \label{induYengui}
\begin{theorem}\emph{(Induction de Yengui)}\label{induYengui} 
\\ Soit $\cF$ une classe  d'anneaux
commutatifs de \ddk finie (non \ncrt bornée) qui satisfait les \prts suivantes.
\begin{description}
\item [{\rm (ls1)}] Si $\gA \in \cF$, alors $\gA(X) \in \cF$.
\item [{\rm (ls2)}] Si $\gA \in \cF$, alors $\gA_S \in \cF$ pour tout \mo $S$ de $\gA$.
\item [{\rm
(ls3)}] Si $\gA \in \cF$, alors tout $\AX$-\mptf est étendu depuis~$\gA$.
\end{description}
%-----------------end item------------------
Alors, pour tout $\gA\in\cF$ et tout $r\geq1$, 
tout \mptf  sur~$\AXr$ est étendu depuis $\gA$.
\end{theorem}
NB: (ls1) remplace (LS2), (ls2) remplace (LS3) et (ls3) remplace (LS4).
\begin{proof}
En raison du fait \ref{fact A->A[X]}~\emph{\ref{i5fact A->A[X]}},
nous nous limitons au cas des anneaux réduits.
Nous raisonnons par \recu double sur le nombre $r$ de
variables et sur la \ddk $d$ de~$\gA$.
\\
 L'initialisation pour $r=1$ ($d$ arbitraire)
est donnée par (ls3), et \hbox{pour $d=0$} (avec $r$ arbitraire)
c'est le \tho de Quillen-Suslin.
\\
 Nous supposons le résultat prouvé en $r$ variables
pour les anneaux dans~$\cF$.
Nous considérons le cas de $r+1$
variables et nous faisons une preuve par \recu sur
(un majorant $d$ de) la \ddk d'un anneau~\hbox{$\gA\in\cF$}.
Soit donc un anneau $\gA$ de \ddk $\leq d+1$.
 Soit $P$ un \mptf sur $\gA[\Xr,Y]=\gA[\uX,Y]$.
Soit~$G=G(\uX,Y)$ une \mpn de $P$ à \coes  
dans~\hbox{$\gA[\uX,Y]$}.
Soit $H(\uX,Y)$ la matrice construite à partir de $G$
comme dans le fait~\ref{factEtPol}.
\\
 En utilisant l'\hdr pour $r$ et (ls1), nous obtenons que les matrices
$H(\uX,Y)$ et $H(\uze,Y)$ sont \elrz ment \eqves sur~$\gA(Y)[\uX]$.
Cela signifie qu'il existe des matrices $Q_{1}$, $R_{1}$ sur
$\gA[\uX,Y]$ telles que
$$\preskip-.2em \postskip.4em
\begin{array}{c}
Q_{1}H(\uX,Y)  =  H(\uze,Y)R_{1} \label{eqeq}\\[1mm]
\hbox{avec}\quad
\det(Q_{1})   \hbox{ et }   \det(R_{1}) \;\hbox{ primitifs  dans } \gA[Y].
\end{array}
$$
Nous montrons maintenant que $H(\uX,Y)$ et $H(\uze,Y)$ sont \eqves \linebreak 
sur~$\ArY[\uX]$.
D'après le \rcm de Vaserstein il suffit de montrer qu'elles sont \eqves sur
$\ArY_{S_{i}}[\uX]$ pour des \mocoz~$S_{i}$ de $\ArY$.
\\
 Nous considérons le \pol primitif
$f=\det(Q_{1})\det(R_{1})\in\gA[Y]$, et nous appliquons le \tho \ref{compa}.
Si $f$ est de degré formel $m$, nous obtenons des \mos $(S_i)_{i\in \lrbm}$ de
$\gA$ tels que les \mos $V=f^{\NN}$ et $(S_i)_{i\in \lrbm}$ 
sont \com  dans $\ArY$. En outre, $\Kdim\gA_{S_i}\leq d$ pour $i\in \lrbm$.
\\
 Pour le localisé en $V$, $\det(Q_{1})$ et
 $\det(R_{1})$ sont \ivs dans $\ArY_V$. Ceci implique que
  $H(\uX,Y)$ et $H(\uze,Y)$ sont \eqves sur $\ArY_V[\uX]$.
\\
Pour un localisé en $S_i$ ($i\in \lrbm$), par \hdr sur $d$  
et en utilisant (ls2),  $H(\uX,Y)$ et $H(\uze,0)$ sont \eqves sur
$\gA_{S_i}[\uX,Y]$.
A~fortiori $H(\uX,Y)$ et $H(\uze,Y)$ sont \eqves sur $\gA_{S_i}[\uX,Y]$, donc aussi  
sur $\ArY_{S_i}[\uX]$, qui est une \lon de
$\gA_{S_i}[Y][\uX]=\gA_{S_i}[\uX,Y]$.

 Ainsi, nous avons rempli le contrat et nous obtenons des matrices \ivsz~$Q$ et $R$ sur $\ArY[\uX] \subseteq \AXY$  telles~que
%\perso{le petit raisonnement de la fin semble inutilement compliqué à première vue, mais je ne vois pas comment le simplifier}
$$
Q\,H(\uX,Y)=H(\uze,Y)\,R.
$$
  Par ailleurs, nous savons par  (ls3) que $H(\uze,0)$ et $H(\uze,Y)$
sont \eqves sur
$\gA[Y]\subseteq \AXY$, et, par \hdr sur $r$,
que $H(\uze,0)$ et~$H(\uX,0)$ sont \eqves sur $\AuX\subseteq \AXY$.
En conclusion $H(\uX,0)$ et~$H(\uX,Y)$
sont \eqves sur $\AXY$.
Donc par le \tho de Horrocks global, $P$ est étendu depuis $\AuX$.\\
 Enfin, par \hdr sur $r$,
$P(\uX,0)$ est étendu depuis~$\gA$.
\end{proof}

\rem Nous avons demandé dans (ls2) que la classe $\cF$ soit stable par
\lon pour n'importe quel \moz. En fait dans la \dem interviennent seulement
des \lons en des \mos bords de Krull, ou par inversion d'un unique \elt (ceci de manière itérée).
\eoe
%%%%%%%%%%%%%%%%%%%%%%%%%%%%%%%%%%%%%%%%%%%%%%%%%%%%%%%%%%%%%%%%%%%

\subsubsection*{Lequain-Simis en dimension finie}

%: --- corollary {corLSValu}  (Lequain-Simis
\begin{corollary}~\label{corLSValu}\\
Si $\gA$ est un \anar de \ddk finie,
tout \mptf sur $\AXr$ est étendu depuis $\gA$.
\end{corollary}
%--------------------------------------------------
\begin{proof}
On montre que la classe des \anars de dimension finie satisfait l'induction
de Lequain-Simis concrète. La condition (ls1) est donnée 
dans l'exercice \ref{exoPrufNagata}, (ls3) par le \tho de Bass-Simis-Vasconcelos, et~(ls2) est clair.
\end{proof}

%%%%%%%%%%%%%%%%%%%%%%%%%%%%%%%%%%%%%%%%%%%%%%%%%%%%%%%%%%%%%%%%%%%
\subsubsection*{Lequain-Simis local sans hypothèse de dimension}

%: --- corollary {corLSValu}  (Lequain-Simis
\begin{corollary}\label{corLSValu2} Si $\gV$ est un \advz,
tout \mptf sur $\VXr$ est étendu depuis $\gV$ (i.e., libre).
\end{corollary}
%--------------------------------------------------
%
\begin{proof}
 Soit  $M$ un \mptf sur $\VXr$.
Nous devons montrer que $M$ est libre.
Soit $F=(f_{ij})\in\GAq(\VXr)$ une matrice dont l'image
est isomorphe au module $M$. Soit
$\gV\!_1$ le sous-anneau de $\gV$
engendré par les \coes des \pols $f_{ij}$ et $\gV'$ le sous-\adv de $\gV$
engendré par $\gV\!_1$. Le point \emph{\iref{i5corthValDim}} du
\thref{corthValDim} nous dit que tout anneau compris entre $\gV\!_1$ et $\Frac\gV\!_1$,
en particulier $\gV'$, est de \ddk
finie. On applique le corolaire~\ref{corLSValu}.
\end{proof}

\subsubsection*{\Tho de Lequain-Simis \gnlz}

%:    th   thLSValu
\begin{theorem}\label{thLSValu} \emph{(Lequain-Simis)}
Si $\gA$ est un \anarz, tout \mptf sur $\AXr$ est étendu depuis
$\gA$.
\end{theorem}
\begin{proof}
Cela résulte du corolaire \ref{corLSValu} (le cas local) avec la même \dem que pour déduire le \thref{thBassAri} du \thoz~\ref{thBassValu}.
\end{proof}

\entrenous{ très incertain \ldots
%:   --- SUBsection{subsecBreCos}------
%\subsubsec{Amélioration de Brewer \& Costa}

{\bf Amélioration de Brewer \& Costa}
\label{subsecBreCos}
%-----------------------------------------

Il serait question de démontrer dans le cas intègre que si la clôture
intégrale de $\gA$ est un anneau seminormal, alors les \mptfs sur $\AXr$
sont tous étendus depuis $\gA$ \ssi $\gA$ est seminormal.

Dans le \tho de Brewer\&Costa proprement dit, $\gA$ doit en plus vérifier
une hypothèse mystérieuse difficile à décrypter.
}

%%%%%%%%%%%%%%%%%%%%%%%%%%%%%%%%%%%%%%%%%%%%%%%%%%%%%%%%%%%%%%%%%%%
%%%%%%%%%%%%%%%%%%%%%%%%%%%%%%%%%%%%%%%%%%%%%%%%%%%%%%%%%%%%%%%%%%%
%%%%%%%%%%%%%%%%%%%%%%                               %%%%%%%%%%%%%%%%%%%%
%%%%%%%%%%%%%%%%%%%%%%       Conjectures             %%%%%%%%%%%%%%%%%%%%
%%%%%%%%%%%%%%%%%%%%%%                               %%%%%%%%%%%%%%%%%%%%
%%%%%%%%%%%%%%%%%%%%%%%%%%%%%%%%%%%%%%%%%%%%%%%%%%%%%%%%%%%%%%%%%%%
%%%%%%%%%%%%%%%%%%%%%%%%%%%%%%%%%%%%%%%%%%%%%%%%%%%%%%%%%%%%%%%%%%%
%%%%%%%%%%%%%%%%%%%%%%%%%%%%%%%%%%%%%%%%%%%%%%%%%%%%%%%%%%%%%%%%%%%

\vspace{5pt}
%---- Section{sec.BQHConj}--
\section*{Conclusion: quelques conjectures}
\addcontentsline{toc}{section}{Conclusion: quelques conjectures}
  \label{sec.BQHConj} 
%------------------

La solution du \pb de Serre a naturellement conduit à poser quelques
conjectures sur de possibles \gnnsz.

Nous citerons les deux plus célèbres et renvoyons à \cite[chap.V,VIII] {Lam06} pour des informations détaillées sur le sujet.

\smallskip  La première, et la plus forte, est la \emph{conjecture des anneaux de Hermite},
qui peut être énoncée sous deux formes \eqvesz, une locale et une globale,
vu le principe de recollement de Quillen.
Rappelons qu'un anneau est appelé \gui{anneau de Hermite} lorsque les modules \tf \stls sont libres, ce qui revient à dire que les \vmds sont complétables. 

 {\bf (H)} Si $\gA$ est un anneau de Hermite, alors $\AX$ \egmtz.

 {\bf (H')} Si $\gA$ est un \alo \dcdz, alors $\AX$ est un anneau de Hermite.

Le \gui{stable-range} de Bass donne une première approche du \pb (voir
la proposition \ref{propStabliblib}, le corolaire \ref{corBass} et le \thref{corBass2}). Des cas particuliers sont traités par exemple dans \cite[Roitman]{Roi86} et \cite[Yengui]{Ye3,Ye4}, qui traite le cas $n=1$ de la conjecture
suivante: sur un anneau~$\gA$ de dimension de Krull $\leq 1$, les $\AXn$-modules stablement 
libres sont libres.

\smallskip  La deuxième est la \emph{conjecture de Bass-Quillen\index{Quillen}}. 

Un \cori est 
appelé un \emph{anneau régulier} si tout \mpf admet une résolution projective finie (pour la \dfn et un exemple de résolution projective finie, voir le \pb
\ref{exoFossumKumarNori}). Ici aussi il y a une version locale et une version globale \eqvesz.

 {\bf (BQ)} Si $\gA$ est un anneau \noe \coh régulier\footnote{Naturellement en \clama l'hypothèse \gui{\cohz} est superflue.}, alors les
\mptfs sur $\AXn$ sont étendus depuis $\gA$.

 {\bf (BQ')} Si $\gA$ est un \alo \dcd \noe \coh régulier\footnote{Naturellement en \clama les hypothèses \gui{\cohz} et \gui{\dcdz} sont superflues.}, alors les
\mptfs sur $\AXn$ sont libres.

En fait, puisque $\gA$ \noe régulier implique $\AX$ \noe régulier, 
il suffirait de démontrer le cas $n=1$. 
Des résultats partiels ont été obtenus.
Par exemple, la conjecture est démontrée en dimension de Krull $\leq2$, 
pour~$n$ arbitraire
(mais on ne dispose pas pour le moment de \dem \covz).
On peut a priori \egmt envisager une version non \noee pour les  \coris
réguliers de \ddk $\leq k$ fixé.

%%%%%%%%%%%%%%%%%%%%%%%%%%%%%%%%%%%%%%%%%%%%%%%%%%%%%%%%%%%%%%%%%%%
%%%%%%%%%%%%%%%%%%%%%%%%%%%%%%%%%%%%%%%%%%%%%%%%%%%%%%%%%%%%%%%%%%%
%%%%%%%%%%%%%%%%%%%%%%%%%%%%%%%%%%%%%%%%%%%%%%%%%%%%%%%%%%%%%%%%%%%
%%%%%%%%%%%%%%%%%%%%%%%%       EXOS              %%%%%%%%%%%%%%%%%%%%%%%%
%%%%%%%%%%%%%%%%%%%%%%%%%%%%%%%%%%%%%%%%%%%%%%%%%%%%%%%%%%%%%%%%%%%
%%%%%%%%%%%%%%%%%%%%%%%%%%%%%%%%%%%%%%%%%%%%%%%%%%%%%%%%%%%%%%%%%%%
%%%%%%%%%%%%%%%%%%%%%%%%%%%%%%%%%%%%%%%%%%%%%%%%%%%%%%%%%%%%%%%%%%%
%%%%%%%%%%%%%%%%%%%%%%%%%%%%%%%%%%%%%%%%%%%%%%%%%%%%%%%%%%%%%%%%%%%

%:section: Exercices
\Exercices{

%-% ENTRE NOUS
\entrenous{Les deux  exos qui suivent semblent en rapport avec Quillen-Suslin
mais le rapport mérite sans doute une explication 
}
%-% Fin ENTRENOUS

%--- Exercise{exolem1SusQS}-------------
\begin{exercise}
 \label{exolem1SusQS}
 {\rm
Soit $\fA$ un \id de $\AX$ contenant un \polu et $\fa$ un \id de $\gA$.
Alors $\gA\,\cap\, (\fA+\fa[X])$ est contenu dans $\DA\big((\gA\,\cap\,\fA)+\fa\big)$.
En particulier, si $1\in \fA+\fa[X]$, alors $1\in(\gA\,\cap\,\fA)+\fa$.
} \end{exercise}
%--- end-exercise-----------------------------------------

%%%%%%%%%%%%%%%%%%%%%%%%%%%%%%%%%%%%%%%%%

%--- Exercise{exolemLocLocCom}-------------
\begin{exercise}
 \label{exolemLocLocCom} (top-bottom lemma)
 {\rm Soit $\gA$ un anneau,  $\fm=\Rad\gA$.
 \begin{enumerate}\itemsep0pt
\item Soit  $S\subseteq\AX$
le \mo des \polusz. Les \mos $S$ et $1+\fm[X]$ sont \comz.
\item Soit $U\subseteq \AX$ le \mo
$\sotq{X^n+\sum_{k<n}a_kX^k}{n\in\NN,a_k\in\fm\; (k<n)}$.
Les \mos $U$ et $1+\fm+X\AX$ sont \comz.
\end{enumerate}

} \end{exercise}
%--- end-exercise-----------------------------------------

%--- Exercise{exoSLequiv}-------------
\begin{exercise}
\label{exoSLequiv}
{\rm  Le but de l'exercice est de montrer un résultat analogue
 au \rcm de Vaserstein (\plgrf{thPatchV}) dans lequel on remplace $\GLn$ par~$\SLn$.
 \begin{enumerate}\itemsep0pt
\item Soit un anneau $\gB$ et un \mo $S$ de $\gB$.
\begin{enumerate}\itemsep0pt
\item Soit $P\in\BY$ tel que $P(0)=0$ et $P=0$ dans $\gB_S[Y]$. Montrer
qu'il existe $s\in S$ tel que $P(sY)=0$.
\item Soit $H\in\Mn(\BY)$ telle que $H(0)\in\SLn(\gB)$
et $H\in\SLn(\gB_S[Y])$. Montrer qu'il existe $s\in S$ tel que
$H(sY)\in\SLn(\BY)$.
\end{enumerate}
\item Montrer le lemme \ref{lem2PrepVaser} en remplaçant
%dans la conclusion (mais pas dans l'hypothèse)
$\GL$ par $\SL$.
\item Montrer le \plgrf{thPatchV} en remplaçant
$\GL$ par $\SL$.
\end{enumerate}
} \end{exercise}
%--- end-exercise-----------------------------------------

%%%%%%%%%%%%%%%%%%%%%%%%%%%%%%%%%%%%%%%%%

%--- Exercise{exoCasParticHorrocks}-------------
\begin{exercise}
 \label{exoCasParticHorrocks}
{\rm   Soit $\gA$ un \alo \dcd et soit $\fb \subseteq \gA[X]$ un \id \iv
contenant un \poluz. On veut montrer que $\fb$ est un \idpz.}

Ceci constitue un cas particulier du théorème de Horrocks local
(\thrf{thHor0}): en effet,
d'une part $\fb$ est un $\gA[X]$-module projectif, et d'autre part, si $f \in\fb$ est un \poluz, alors en localisant en $f$, $\fb_f =
\gA[X]_f$, et donc, par le \tho de Horrocks local, $\,\fb$ est un $\gA[X]$-module libre. Cet exercice donne une \dem indépendante de celle du cours.
Dans le cas particulier étudié ici, on apporte
la précision que $\fb$ est engendré par un \poluz.

 {\rm  Soit $\gA$ un anneau, on note $\fm=\Rad\gA$ et $\gk = \gA/\fm$.
Soit $\fb \subseteq \AX$ un \id contenant un \poluz.
On note $\ov a$ la réduction de $a$ modulo $\fm$.
\begin{enumerate}\itemsep0pt
\item
%
% ,$\gA \to \gA/\fa$ et $\gA[X] \to (\gA/\fa)[X]$
Montrer que tout \polu de $\ov\fb\subseteq\kX$ peut être relevé en un \polu de~$\fb$.
\end{enumerate}
On suppose maintenant $\gA$ local \dcdz.
\begin{enumerate} \setcounter{enumi}{1}\itemsep0pt
\item Montrer l'existence d'un \polu $f \in \fb$ tel
que $\ov \fb = \langle \ov f\rangle$ dans $\gk[X]$
et donc $\fb = \langle f\rangle + \fb \cap \fm[X]$.
\end{enumerate}
On suppose maintenant que l'\id $\fb$ est \ivz.
\begin{enumerate} \setcounter{enumi}{2}\itemsep0pt
\item
%En utilisant le fait que $\fb$ est \ivz,
Montrer que $\fb \cap \fm[X] = \fb\fm[X]$%
% : on pourra considérer un \id $\fb'$ tel que $\fb\fb' = \fb
%\cap \fm[X]$
.

\item On considère l'anneau $\aqo{\AX}{f}$.
Montrer que $\fm (\aqo{\fb}{f}) = \aqo{\fb}{f}$.
En déduire que $\fb = \gen{f}$.
\end{enumerate}
On propose une \gnnz.
\begin{enumerate} \setcounter{enumi}{4}
\item La \dem fonctionne-t-elle avec un anneau $\gA$ \plcz?

\end{enumerate}

} \end{exercise}
%--- end-exercise-----------------------------------------

%%%%%%%%%%%%%%%%%%%%%%%%%%%%%%%%%%%%%%%%%

%--- Exercise{exoBézoutKdim1TransfertArX}-------------
\begin{exercise}\label{exoBézoutKdim1TransfertArX}
{(\Tho de Brewer{\&}Costa: cas des anneaux
de Bézout intègres de dimension $\le 1$)}
{\rm  Voir aussi l'exercice \ref{exoPrufNagata}
et le \thref{thMaBrCo}.
 
     Soit  $\cF$  la classe des anneaux intègres de Bézout de dimension $\le
1$, et $\gA\in\cF$. 
 
\emph {1.}
Montrer que $\Kdim\ArX \le 1$ (utiliser l'exercice~\ref{exoMultiplicativiteIdeauxBords}).

\emph {2.}
En déduire que $\ArX$ est un anneau de Bézout.
 
\emph {3.}
La classe $\cF$  satisfait les hypothèses du \thref {th2QUILIND} (induction de Quillen concrète, cas libre). Ainsi,  tout $\AXr$-\mptf est libre.

}

\end {exercise}
%--- end -exercise-----------------------------------------

%--- Exercise{exoseminorlgb}-------------
\begin{exercise}
\label{exoseminorlgb} (Principe \lgb pour les anneaux seminormaux)\\
{\rm  On donne une \dem directe du principe \ref{plgcetendus} dans le cas particulier des anneaux \lsdz seminormaux.
\emph{1.} Dans un anneau \lsdzz, si $xc=b$ et $b^2=c^3$, alors
il existe $z$ tel que $zc=b$ et $z^2=c$, donc $z^3=b$.

\emph{2.} Soient $S_1$, \ldots, $S_n$ des \moco d'un anneau  $\gA$.
On suppose que chacun des $\gA_{S_i}$ est \lsdz et seminormal. 
Montrer que~$\gA$ est \lsdz et seminormal.
 
}
\end{exercise}
%--- end -exercise-----------------------------------------

%--- Exercise{exoKdimBounded}-------------
\begin{exercise}\label{exoKdimBounded}
{(Anneaux vérifiant certaines des conditions de la section 
\gui{Un exemple en \ddk finie $>0$} \paref{subsec1VXlibre})}\\
{\rm
Soient $a_1$, \ldots, $a_k \in \gA$, $a_0 = 0$,
$a_{k+1} = 1$, et les anneaux $\gA_1$, \ldots, $\gA_{k+1}$ suivants

\snic{
\gA_i = \big(\aqo{\gA}{a_{i-1}}\!\big)[1/a_i] \quad \hbox {pour $i \in \lrb{1..k+1}$}
}

\snii
On va montrer que si chaque $\gA_i$ est \zedz, alors $\Kdim\gA \le k$.
Le même résultat vaut avec $\gA_i = \big(\gA/\rD_\gA(a_{i-1})\!\big)[1/a_i]$.

\snii
\emph {1.}
Soit $a \in \gA$. Si $\Kdim \gA[1/a] \le n$ et
$\Kdim \aqo{\gA}{a} \le m$, alors $\Kdim\gA \le n+m+1$.

\snii
\emph {2.}
En déduire le résultat annoncé.

}

\end {exercise}
%--- end -exercise-----------------------------------------
%%%%%%%%%%%%%%%%%%%%%%%%%%%%%%%%%%%%%%%%%

%%%%%%%%%%%%%%%%%%%%%%%%%%%%%%%%%%%%%%%%%

}% fin des exos
%:  solutions

\penalty-2500
\sol{

%%%%%%%%%%%%%%%%%%%%%%%%%%%%%%%%%%%%%%%%%

\exer{exolem1SusQS}
{Posons $\gB=\gA\sur{\gA\,\cap\,\fA}$, $\gB'=\AX\sur{\fA}$, $\fb=\ov{\fa}$,
$\fb'=\fb\,\gB'$. \\
L'anneau $\gB'$ est une extension entière de
$\gB$. On applique le lying over (\ref{lemLingOver}).\\
\emph{Une autre solution}.
Soit $f\in\fA$ \monz. Soit $a\in\gA\,\cap\, (\fA+\fa[X])$,
il existe $g\in\fA$ tel que $g\equiv a \mod \fa$. Alors $\Res(f,g)\equiv\Res(f,a)\mod\fa$. Mais $\Res(f,a)=a^{\deg f}$ et $\Res(f,g)\in\fA\,\cap\,\gA$.
}

%%%%%%%%%%%%%%%%%%%%%%%%%%%%%%%%%%%%%%%%%%%%%%%%%%%%%%%%%%%%%%%%%%%

\exer{exolemLocLocCom}Utiliser le résultant.

%%%%%%%%%%%%%%%%%%%%%%%%%%%%%%%%%%%%%%%%%%%%%%%%%%%%%%%%%%%%%%%%%%%

\exer{exoSLequiv}
\emph{1b.} On pose $P(Y)=1-\det\big(H(Y)\big)$ et l'on applique le point \emph{1a}.

\emph{2.} Le lemme  \ref{lem2PrepVaser} nous fournit une matrice $U(X,Y)\in
\GL_r(\gA[X,Y])$ telle que 

\snic{U(X,0)=\I_r$ et, sur $\gA_S[X,Y]$,
$U(X,Y)=C(X+sY)C(X)^{-1}.}

%\sni
D'après le point \emph{1}, il existe  $t\in S$
tel que $U(X,tY)\in \SL_r(\gA[X,Y])$.\\
On pose $V(X,Y)=U(X,tY)$
et l'on remplace $s$ par $st$. 

\emph{3.} Le lemme  \ref{lem3PrepVaser} subit avec succès le remplacement
de $\GL$ (implicite dans le mot \gui{\eqvez}) par $\SL$.
Même chose ensuite pour le \rcm de Vaserstein.

%%%%%%%%%%%%%%%%%%%%%%%%%%%%%%%%%%%%%%%%%%%%%%%%%%%%%%%%%%%%%%%%%%%

\exer{exoCasParticHorrocks}
\emph{1.}
Montrons d'abord le résultat suivant: si l'on a $g$, $f \in \fb$ avec
$\ov g$ \mon de degré $r$ et $f$ \mon de degré $r+1$, alors
$\ov g$ peut être relevé en un \polu de $\fb$
(de degré $r$). On écrit $g = aX^{r +
\delta} + \dots $, avec $\delta \in \NN$ et l'on montre par \recu sur
$\delta$ que $\ov g$ peut être relevé en un \polu de
$\fb$. Si $\delta = 0$, on a $a \equiv 1 \bmod \fm$ (car $\ov g$ est
\monz), donc $a$ est \iv et le \polu $a^{-1}g
\in \fb$  relève $\ov g$.  Si $\delta
\ge 1$, on a $a \in \fm$ (car $\ov g$ est \monz), et l'on considère $h
= g - aX^{\delta-1}f \in \fb$. Il est de la forme $bX^{r + \delta-1} + \dots$,
et il vérifie~\hbox{$\ov h = \ov g$}. On applique l'\hdrz.
\\ 
Il suffit ensuite de montrer que pour tout $g \in \fb$ tel que $\ov g$
est \mon de degré~$r$, l'\id $\fb$ contient un \polu de degré
$r+1$. Par hypothèse,~$\fb$ contient un \polu $f$. Si $
\deg(f) \le r+1$, alors le résultat est clair. Si~\hbox{$n =
\deg(f) > r+1$}, alors le \pol $X^{n-(r+1)}\ov g$ est
\mon de degré $n-1$, et d'après la première étape,
%(appliquée à $(n,n-1)$ au lieu de $(r+1,r)$),
$\fb$ contient un \polu de
degré $n-1$. On~conclut par \recu sur $n-r$.
%Une \recu \imde sur $n$ montre que $\fb$ contient
%un \polu de degré $r+1$.

 \emph{2.}
L'idéal $\ov \fb$ est un \itf de $\gk[X]$,
%qui est un anneau de Bézout intègre,
donc $\ov\fb$ est principal.  Comme $\fb$
contient un \polu
on peut prendre le \gtr $\ov h$ \mon et l'on le
relève en un \polu de $\fb$ d'après la question
précédente.

 \emph{3.} Soit $f$ \mon dans $\fb$, et $\fb_1$ l'\id qui vérifie $\fb\fb_1=\gen{f}$.\\
On considère  $\fb'=\fb_1(\fb\cap \fm[X])/f$
(c'est un \id de $\AX$).
Alors $\fb\fb' = \fb\cap \fm[X]$.
On a $f\fb'\subseteq\fm[X]$ et $f$ \mon donc $\ov{\fb'}=0$,
\cad $\fb' \subseteq \fm[X]$. En multi\-pliant par $\fb$, on obtient $\fb
\cap \fm[X] \subseteq \fb\fm[X]$, donc $\fb \cap \fm[X] = \fb\fm[X]$.

 \emph{4.} On a

\snic{
\fm (\aqo{\fb}{f}) = \aqo{\fc}{f}
\; \hbox {avec} \;
\fc = \fm\fb + \gen{f} = \fm[X]\fb +  \gen{f} =
\fm[X] \cap \fb + \gen{f} = \fb.}

%\sni
Le $\aqo{\AX}{f}$-module $\aqo{\fb}{f}$ est \tf et
comme $f$ est \monz, $\aqo{\AX}{f}$ est un \Amo \tfz.
On en déduit que $\aqo{\fb}{f}$ est un \Amo \tfz.
Par le lemme de Nakayama on obtient $\aqo{\fb}{f} =
0$, i.e. $\fb = \gen{f}$.

%%%%%%%%%%%%%%%%%%%%%%%%%%%%%%%%%%%%%%%%%%%%%%%%%%%%%%%%%%%%%%%%%%%

\exer{exoBézoutKdim1TransfertArX} 
\emph {1.}
Il faut montrer que pour $f$, $g \in \gA[X]$, on a $1 \in
\IK_{\ArX}(f, g)$.  Puisque $\gA$ est de Bézout intègre, tout
\pol de $\gA[X]$ est le produit d'un \elt de $\gA$ par un \pol
primitif. D'après l'exercice \ref {exoMultiplicativiteIdeauxBords}, il
suffit de montrer que $1 \in \IK_{\ArX}(f, g)$, soit lorsque $f$ ou $g$ est
primitif, soit lorsque $f$ et $g$ sont des constantes $a$, $b$. Dans ce dernier
cas, cela découle, puisque $\Kdim\gA \le 1$, 
de $1\in \IK_{\gA}(a,b) \subseteq \IK_{\ArX}(a,b)$.
\\ 
On suppose donc que $f$ ou $g$ est primitif, par exemple $f$. Il suffit de montrer \linebreak 
que $1 \in \IK_{\ArX}(f, g)$ après \lon en des
\mocoz. Or le \thref {compa} fournit des \mos bord $S_j$ dans $\gA$ tels
que $f^\NN$ et les $S_j$ sont \com dans $\ArX$. Pour la \lon en
$f^\NN$, il est clair que $1 \in \IK(f,g)$.  
\\
 Quant à $S_j^{-1}\ArX$, c'est
une \lon de $\gA_{S_j}[X]$ avec $\gA_{S_j}$ \zedz, ce qui donne 
 $\Kdim \gA_{S_j}[X] \le 1$. Donc $1 \in \IK(f,g)$ dans $\gA_{S_j}[X]$, 
et a fortiori dans le localisé $S_j^{-1}\ArX$.

\emph {2.}
L'anneau $\gA[X]$ est un anneau intègre à pgcd, donc il en
est de même de son localisé~$\ArX$. Comme $\Kdim \ArX \le 1$,
 le \thref{propGCDDim1} nous dit que $\ArX$
est un anneau de Bézout.

\emph {3.} On a démontré la \prt (q1) et l'on sait 
déjà que la \prt (q0) est satisfaite (\thref{thBézoutLib}).

%%%%%%%%%%%%%%%%%%%%%%%%%%%%%%%%%%%%%%%%%

\exer{exoseminorlgb}On peut commencer par donner une \dem directe que toute \lon d'un anneau seminormal est encore un anneau seminormal (implicitement supposé dans l'énoncé).
Soit  en effet $\gA$ seminormal et $S$ un \moz. Supposons  $(\fraC x s)^2=(\fraC y s)^3$ dans $\gA_S$. On a donc pour un $t\in S$, on a $tsx^2=ty^3$ dans $\gA$. D'où $t^6s^4x^2=t^6s^3y^3$. Or $\gA$ est seminormal, on a donc un $z\in\gA$
tel que $t^3s^2x=z^3$ et $t^2sy=z^2$, ce qui donne dans $\gA_S$ les \egts $\fraC x s=(\fraC z {st})^3$ et $\fraC y s=(\fraC z {st})^2$.

\emph{1.}  On a $x^2c^2=b^{2}=c^3$, donc $c^2(x^2-c)=0$,
donc $c(x^2-c)=0$. 
\\
Soient alors $s$, $t$ tels que $s+t=1$, $sc=0$ et $t(x^2-c)=0$. Posons
$z=tx$. \\
On a $tc=c$, $z^2=t^2c=c$ et $zc=xtc=xc=b$.

\emph{2.} 
On suppose que chacun des  $\gA_{S_i}$ est \lsdz et seminormal. Donc $\gA$ est \lsdzz. 
 \\
 Soient $b$, $c\in\gA$  avec $b^2=c^3$. Si les  
$\gA_{S_i}$ sont seminormaux, il existe $x_i\in\gA_{S_i}$ tels que $x_i^2=c$ et
$x_i^3=b$, et donc $x_ic=b$. Ceci implique qu'il existe $x\in\gA$ 
tel \hbox{que $xc=b$}.
On conclut par le  point \emph{1}.

NB: Il y a des anneaux seminormaux qui ne sont pas \lsdzz: par exemple $\gk[x,y]$ avec $xy=0$
où $\gk$ est un corps discret.

%%%%%%%%%%%%%%%%%%%%%%%%%%%%%%%%%%%%%%%%%%%%%%%%%%%%%%%%%%%%%%%%%%%

\exer{exoKdimBounded} 
\\
\emph {1.}
Soient $n+m+2$ \elts de $\gA$, $(\ux) = (x_0, \ldots, x_n)$, $(\uy) = (y_0, \ldots,y_m)$.  
En considérant le \mo bord itéré de $(\uy)$ dans $\aqo{\gA}{a}$, on obtient
que $\SK_\gA(\uy)$ contient un multiple de $a$, disons $ba$. En considérant
l'\id bord itéré de $(\ux)$ dans $\gA[1/a]$, on obtient que $\IK_\gA(\ux)$ contient
une puissance de $a$, disons $a^e$.
\\
 Alors $(ba)^e \in \IK_\gA(\ux) \cap
\SK_\gA(\uy)$, donc $1 \in \IK_\gA(\ux,\uy)$ d'après le
fait~\ref{fact0BordKrullItere}, point~\emph{1}.

\snii
\emph {2.}
En utilisant la question précédente, on montre par \recu sur $i \in
\lrb{0..k+1}$, que l'on a $\Kdim\gA[1/a_i] \le i-1$; pour $i = k+1$, on
obtient $\Kdim\gA \le k$.

%%%%%%%%%%%%%%%%%%%%%%%%%%%%%%%%%%%%%%%%%%%%%%%%%%%%%%%%%%%%%%%%%%%

%%%%%%%%%%%%%%%%%%%%%%%%%%%%%%%%%%%%%%%%%%%%%%%%%%%%%%%%%%%%%%%%%%%

}% fin des solutions d'exos

%:   ---- Section*{references}-----------
\Biblio

Carlo Traverso a démontré dans  \cite{Tra} le \tho qui porte son nom pour un anneau \noe réduit
$\gA$ (avec une restriction supplémentaire).
Pour le cas intègre sans hypothèse
\noee on peut
consulter \cite[Querré]{Querre}, \cite[Brewer\&Costa]{BC2} et \cite[Gilmer\&Heitmann]{GH}. 
%:H
Le cas le plus \gnl est donné par \cite[Swan]{Swan80}.

Le \tho de Traverso-Swan sur les anneaux seminormaux a été décrypté du
point de vue \cof par Coquand dans \cite{coq}.
Le décryptage a commencé par la \dem \elr de la proposition \ref{propIntSemin} telle qu'elle est donnée ici.
Cette \dem est une simplification (assez spectaculaire)
des preuves existantes dans la littérature. Il fallait ensuite
contourner l'argument de la considération d'un \idemi pour obtenir
une \prco complète du résultat. Il est remarquable que par la même
occasion, le cas d'un anneau non intègre ait pu être traité sans plus d'effort, contrairement à ce qui se passe avec la \dem de Swan dans~\cite{Swan80}. Pour une explication détaillée de \cite{coq} voir \cite[Lombardi\&Quitté]{LQ06}.
%:HHH 
Pour un \algo \gui{simple} qui réalise le \tho 
dans le cas univarié, voir~\cite[Barhoumi\&Lombardi]{BL07}. 
Une \dem directe, dans le même esprit, pour l'implication \gui{$\gA$ seminormal implique $\AX$ seminormal} se trouve
dans \cite[Barhoumi]{Barh09}.

Le \tho de Roitman \ref{thRoitman} se trouve dans \cite{Roi}.

En ce qui concerne l'historique de la résolution du problème
de Serre sur les anneaux de \polsz, le lecteur pourra consulter
le chapitre~III de l'ouvrage de Lam~\cite{Lam06} ainsi
que l'exposé de Ferrand à Bourbaki~\cite{Ferrand}.

Les \dems originales du \tho de Quillen-Suslin (solution du \pb de Serre)
se trouvent dans\index{Quillen}\index{Suslin} \cite[Quillen]{Qu} et \cite[Suslin]{Sus}. Les \thos de Horrocks ont leur source dans \cite[Horrocks]{Hor}.

Le \gui{Quillen patching} qui apparaît dans \cite{Qu} est parfois
appelé \plg de Quillen. Un survol remarquable des applications de ce principe
et de ses extensions se trouve dans \cite[Basu\&al.]{BaRaKha}.
\`A lire \egmtz~\cite[Rao\&Selby]{RaSe}.

L'anneau $\ArX$ a joué un grand rôle dans la solution du problème de Serre
par Quillen et dans ses \gnns successives (\thos de Maroscia et Brewer{\&}Costa, et
de Lequain\&Simis).
L'anneau $\gA(X)$ s'est avéré un outil efficace pour plusieurs
résultats d'\alg commutative. On pourra consulter l'article \cite[Glaz]{Glaz1}
pour une bibliographie assez complète concernant ces deux anneaux.

Le livre de Lam \cite{Lam06} (qui fait suite à \cite{Lam}) 
est une mine d'or concernant les modules \pros étendus. 
Il contient notamment plusieurs preuves des \thos de Horrocks (local et global), avec tous les détails et toutes les références \ncrsz,
au moins du point de vue des \clamaz.

Le \thrf{thHor} de Horrocks global a été démontré \cotz,  tout d'abord 
(pour  une variante un peu plus faible) dans l'article \cite[Lombardi\&Quitté]{LQ02}, puis dans \cite[Lombardi,Quitté\&Yengui]{LQY05}.
La version que nous donnons \paref{thHor} reprend ce dernier article en précisant
tous les détails. Elle s'appuie sur les livres de Kunz et Lam.

Le \thrf{thBassValu} de Bass-Simis-Vasconcelos (\cite{BASS,SV})
a été décrypté du
point de vue \cof par Coquand dans~\cite{coq07}.

Concernant le \tho de Maroscia et Brewer{\&}Costa (\thref{thMaBrCo}), voir
les articles originaux \cite{BC,Ma}. 
 On en trouve une \dem \cov dans \cite{LQY05}.
Ce \tho est légèrement antérieur au \tho de Lequain\&Simis.
Ce dernier  a été décrypté du
point de vue \cof essentiellement par I. Yengui~\cite{BLY06,ELY07}.

De nombreux \algos pour le \tho de Quillen-Suslin (cas des corps) ont été proposés en calcul formel, en \gnl basés sur la \dem de \Susz.
 
Le \tho de Quillen-Suslin a été étudié du point de vue de sa
complexité \algq dans \cite[Fitchas\&Galligo]{FiGa} et \cite[Caniglia\&al.]{CCDHKS} (pour des \algos efficaces, mais semble-t-il non encore implémentés).

Un nouvel \algoz, simple et efficace, pour le \tho de \Sus (compléter un \vmd
(contenant un \poluz)  de $\AX$) est donné dans \cite[Lombardi\&Yengui]{YL04} et amélioré dans \cite[Mnif\&Yengui]{MnY}.  

%-% ENTRE NOUS
\entrenous{Les références  
Vaserstein pour le
petit \tho de Horrocks local \rref{th2HorrocksLocal} ou le global \rref{th2HorrocksGlobal} sont-elles bien dans la biblio?
en tout cas il y a les références \cite{Vaserstein1} et 
\cite{Vaserstein3} ne sont citées nullepart
}
%-% Fin ENTRENOUS

\newpage \thispagestyle{CMcadreseul}
\incrementeexosetprob

%:        %%%%%%%%%%%%%%%%%%%%%%%%%%%%%%%%%%%%
%:        %%%%%%%%%%%%%%%%%%%%%%%%%%%%%%%%%%%%
%---- Chapitre  {Tho de stabilité de Suslin} 
\chapter[\Tho de stabilité de Suslin]{\Tho de stabilité de \Susz, le cas des corps} 
\label{ChapSuslinStab}
\minitoc

%: Intro ------------------
\Intro

%\incertain{

Dans ce chapitre, nous donnons un traitement entièrement \cof du \tho
de stabilité de Suslin pour le cas des \cdisz. 
%%%%%%%%%%%%%%%%%%%%%%%%%%%%%%%%%%%%%%%%%%%%%%%%%%%%%%%%%%%%%%%%%%%
%%%%%%%%%%%%%%%%%%%%%%%%%%%%%%%%%%%%%%%%%%%%%%%%%%%%%%%%%%%%%%%%%%%
%%%%%%%%%%%%%%%%%%%%                                     %%%%%%%%%%
%%%%%%%%%%%%%%%%%%%%        groupe elementaire           %%%%%%%%%%
%%%%%%%%%%%%%%%%%%%%                                     %%%%%%%%%%
%%%%%%%%%%%%%%%%%%%%%%%%%%%%%%%%%%%%%%%%%%%%%%%%%%%%%%%%%%%%%%%%%%%
%%%%%%%%%%%%%%%%%%%%%%%%%%%%%%%%%%%%%%%%%%%%%%%%%%%%%%%%%%%%%%%%%%%
%\newpage	
\section{Le groupe \elrz} \label{secGpEn}

\vspace{4pt}
%:   subsec{Transvections}
\subsec{Transvections}
\label{subsecTransvections} 

Concernant le groupe \elr $\En(\gA)$ rappelons
qu'il est engendré par les matrices \elrs $\rE^{(n)}_{i,j}(a)=\rE_{i,j}(a)$. 

Si l'on note $(e_{ij})_{1 \leq i
, j \leq n}$  la base canonique de $\Mn(\gA)$, 
on~a:
$$\preskip.4em \postskip.4em 
\rE_{i,j}(a) = \I_n + a e_{ij},
\quad 
\dsp\rE_{i,j}(a)\,e_k =
\cases {e_k               &si $k \ne j$    \cr  
        e_j + a e_i &si $k = j$\cr}
        \quad (i\neq j), 
$$
avec par exemple
$$\preskip.0em \postskip.4em
\rE_{2,3}(a) = \cmatrix {
1 & 0 & 0 & 0\cr
0 & 1 & a & 0\cr
0 & 0 & 1 & 0\cr
0 & 0 & 0 & 1\cr
}.
$$

Pour $i$ fixé (resp.\ pour $j$ fixé) les matrices $\rE_{i,j}(\bullet)$
commutent, et forment un sous-groupe de $\EE_n(\gA)$ isomorphe à
$(\Ae {n-1},+)$. Par exemple
$$\preskip.4em \postskip-.2em
\rE_{2,1}(a)\cdot \rE_{2,3}(b)\cdot \rE_{2,4}(c) = \cmatrix {
1 & 0 & 0 & 0\cr
a & 1 & b & c\cr
0 & 0 & 1 & 0\cr
0 & 0 & 0 & 1\cr
}
$$
et
$$\preskip.2em \postskip.4em 
\cmatrix {
1 & 0 & 0 & 0\cr
a & 1 & b & c\cr
0 & 0 & 1 & 0\cr
0 & 0 & 0 & 1\cr
} \cdot\cmatrix {
1 & 0 & 0 & 0\cr
a' & 1 & b' & c'\cr
0 & 0 & 1 & 0\cr
0 & 0 & 0 & 1\cr
}=
 \cmatrix {
1 & 0 & 0 & 0\cr
a+a' & 1 & b+b' & c+c'\cr
0 & 0 & 1 & 0\cr
0 & 0 & 0 & 1\cr
}.
$$

Plus \gnlt soit $P$ un \Amo \ptfz. On dira qu'un couple $(\lambda,w)\in P\sta\times P$
est \emph{\umdz} si $\lambda(w)=1$. 
Dans ce cas $w$ est un \elt \umd de $P$, 
$\lambda$  est un \elt \umd de $P\sta$ et l'\Ali $\theta_P(\lambda\otimes w):P\to P$
définie par $x\mapsto\lambda(x)w$ est la projection sur $L=\gA w$ \paralm
à $K=\Ker \lambda$, représentée sur $K\times L$ par la matrice%
\index{couple unimodulaire}%
\index{unimodulaire!couple ---}

\snic{\bloc{0_{K\to K}}{0_{L\to K}}{0_{K\to L}}{1_{L\to L}}=
\bloc{0_{K\to K}}{0}{0}{\Id_L}.
}

\smallskip  Si $u\in K$, l'\Ali $\tau_{\lambda,u}:=\Id_P+\theta_P(\lambda\otimes u)$,
$x\mapsto x+\lambda(x)u$ est appelée une \ix{transvection}, 
elle est représentée  sur $K\times L$ par la matrice

\snic{ \bloc{1_{K\to K}}{(\lambda\otimes u)\frt{L}}{0_{K\to L}}{1_{L\to L}}
= \bloc{\Id_{K}}{(\lambda\otimes u)\frt{L}}{0}{\Id_L}.
}

 Par exemple, si $P=\Ae n$, une
matrice \elr définit une transvection.

\rdb\smallskip\label{NOTAtransvec}  
On note $\GL(P)$ le groupe des \autos
\lins de $P$ et $\SL(P)$ le sous-groupe des \endos de \deter $1$.
Le sous-groupe de $\SL(P)$ engendré par les transvections
sera noté $\wi{\EE}(P)$.
L'application affine 

\snic{u\mapsto \tau_{\lambda,u},\;\Ker\lambda\to \End_\gA(P)}

fournit un \homo du groupe $(\Ker\lambda,+)$ dans le groupe $\wi{\EE}(P)$. 

Dans le cas où $P=\Ae {n}$, si  $\lambda$ est une forme coordonnée, on trouve que la matrice de
la transvection est un produit de matrices \elrsz. Par exemple
avec le vecteur $u=\lst{u_1\,u_2\,u_3\,0}$:
$$\preskip.4em \postskip.4em 
\bloc{\I_{3}}{u}{0}{1}\;=\;
\cmatrix{1&0&0&u_1\cr
0&1& 0&u_2\cr
0& 0 &1&u_{3}\cr
0&0 & 0 &1}\; =\;\dsp\prod_{i=1}^{3}\;\rE_{i,4}(u_i). 
$$
Notez cependant qu'a priori $\En(\gA)$ est seulement \emph{contenu
dans} $\wi{\EE}(\Ae n)$. Ceci montre que le groupe \elr est a priori
dépourvu de signification \gmq claire. Comme point crucial,
 $\En(\gA)$ n'est pas a priori stable par $\GLn(\gA)$-conjugaison.

%:   subsec{Matrices spéciales}
\subsec{Matrices spéciales}
\label{subsecMatSpec} 

%Nous parlons maintenant uniquement des groupes $\En(\gA)$.

Soient $u = \cmatrix {u_1\cr \vdots\cr u_n} \in \Ae {n \times 1}$ et $v =
\cmatrix {v_1 & \cdots & v_n} \in \Ae {1 \times n}$ auxquels on associe la
matrice $\I_n + uv \in \Mn(\gA)$. On va fournir des résultats
précisant l'appartenance de cette matrice au groupe \elr
$\En(\gA)$.  
\\
Puisque $\det(\I_n + uv) = 1 + \tr(uv) = 1 + vu$, il est
impératif de réclamer l'\egt $vu \eqdefi v_1u_1 + \cdots + v_nu_n
= 0$.  Dans ce cas, on a $(\I_n + uv)(\I_n - uv) = \I_n$.

Les  transvections admettent pour matrices les matrices de ce type,
avec $v$ \umdz.
En outre, l'ensemble de ces matrices $\I_n + uv$ (avec $vu = 0$) est
un ensemble stable par $\GLn(\gA)$-conjugaison. 
\\
Par exemple pour $A \in \GLn(\gA)$, on obtient 
$A\, \rE_{ij}(a)\, A^{-1} = \I_n +
a uv$, où $u$ est la colonne $i$ de $A$ et $v$ la ligne $j$ de~$A^{-1}$.

Prenons garde cependant 
que si l'on ne suppose pas $v$ \umd ces matrices ne représentent
en \gnl pas des transvections. Si ni $u$, ni $v$ n'est \umd
la matrice ne représente même pas a priori un \elt de~$\wi\EE(\Ae n)$.

%%%%%%%%%%%%%%%%%%%%%%%%%%%%%%%%%%%%%%%%%
%:     Lemma{lemE01Rao}
\begin{lemma}\label{lemE01Rao}
Supposons  $u \in \Ae {n \times 1}$, $v \in \Ae {1 \times n}$ et $vu
= 0$. \\
Alors 
$\cmatrix {\I_n + uv & 0\cr 0 & 1\cr} \in \EE_{n+1}(\gA)$.
\end{lemma}
%%%%%%%%%%%%%%%%%%%%%%%%%%%%%%%%%%%%%%%%%

%%%%%%%%%%%%%%%%%%%%%%%%%%%%%%%%%%%%%%%%%
\begin {proof}
On a une suite d'opérations \elrs à droite (la première utilise l'\egt $vu=0$):

\snic{\arraycolsep2pt\begin{array}{cccccccccccccc}
\cmatrix {\I_n + uv & 0\cr 0 & 1\cr}  & \vers{\alpha}  &   
\cmatrix {\I_n + uv & -u\cr 0 & 1\cr} & \vers{\beta}\\[3mm]
   \cmatrix {\I_n  & -u\cr v & 1\cr}  & \vers{\gamma} & 
\cmatrix {\I_n & 0\cr v & 1\cr}  & \vers{\delta}  &  \cmatrix {\I_n & 0\cr 0 & 1}. 
\end{array}
}

%\sni
%qui
%transforme la matrice $\cmatrix {\I_n + uv & 0\cr 0 & 1\cr}$
%en $\I_{n+1}$.
Ceci implique
 $\cmatrix {\I_n + uv & 0\cr 0 & 1\cr} = 
   \delta^{-1} \gamma^{-1} \beta^{-1} \alpha^{-1}$,
\cad
$$\preskip.4em \postskip.4em 
\cmatrix {\I_n + uv & 0\cr 0 & 1\cr} =\cmatrix {\I_n & 0\cr v & 1\cr} \cdot
\cmatrix {\I_n & -u\cr 0 & 1\cr} \cdot
\cmatrix {\I_n & 0 \cr -v & 1\cr} \cdot
\cmatrix {\I_n & u\cr 0 & 1\cr}. 
$$
\end {proof}
%%%%%%%%%%%%%%%%%%%%%%%%%%%%%%%%%%%%%%%%%

Un vecteur colonne $u$ sera dit \emph{spécial} si au moins une de ses 
\coos
est nulle.
\perso{Remords, dommage que cela soit faux: 
Si $vu=0$ et si la \coo \num$i$
de $u$ est nulle, on peut remplacer $v$ par $v'$  en remplaçant 
la \coo  \num$i$
de $v$ par $1$. Alors $v'$ est \umd et $\In+uv=\In+u'v$ est la matrice
d'une transvection, que nous appellerons \emph{transvection spéciale}.}
%:finperso
Si $vu=0$ et si  $u$ est spécial nous dirons que  $\In+uv$
est une \emph{matrice spéciale}.\index{matrice!spec@spéciale}
%:HHH index

%%%%%%%%%%%%%%%%%%%%%%%%%%%%%%%%%%%%%%%%%
%:     Corollary{corlemE01Rao}
\begin{corollary}\label{corlemE01Rao}
Soient $u \in \Ae {n \times 1}$ et $v \in \Ae {1 \times n}$ vérifiant $vu
= 0$. Si $u$ est spécial, alors $\I_n + uv \in \En(\gA)$.
Autrement dit toute matrice spéciale est dans~$\En(\gA)$.
\end{corollary}
%%%%%%%%%%%%%%%%%%%%%%%%%%%%%%%%%%%%%%%%%

%%%%%%%%%%%%%%%%%%%%%%%%%%%%%%%%%%%%%%%%%
\begin {proof}
On peut supposer que $n\geq2$ et $u_n = 0$. \\
\'Ecrivons $u = \cmatrix {\mathring
{u} \cr 0\cr}$, $v = \cmatrix {\mathring {v} & v_n}$, avec $\mathring {u} \in
\Ae {(n-1) \times 1}$ et $\mathring v \in \Ae {1 \times (n-1)}$. Alors
$$
\I_n + uv = 
\cmatrix {\I_{n-1} + \mathring{u}\,\mathring{v} & v_n\mathring{u}\cr 0 & 1\cr} =
\cmatrix {\I_{n-1} & v_n\mathring {u}\cr 0 & 1\cr} 
\cmatrix {\I_{n-1} + \mathring {u}\,\mathring{v} & 0\cr 0 & 1\cr}. 
$$
Puisque $\mathring {v}\,\mathring{u} = vu = 0$, le lemme \ref{lemE01Rao}
s'applique et
$\I_n + uv \in \En(\gA)$.
\end {proof}
%%%%%%%%%%%%%%%%%%%%%%%%%%%%%%%%%%%%%%%%%

Les matrices spéciales se \gui
{relèvent} facilement d'un localisé $\gA_S$ à $\gA$
lui-même. Précisément, on obtient ce qui suit.

%%%%%%%%%%%%%%%%%%%%%%%%%%%%%%%%%%%%%%%%%
%:     Fact{factE02Rao}
\begin{fact}\label{factE02Rao}
Soit $S \subseteq \gA$ un \moz, $u \in \gA_S^{n \times 1}$, $v \in
\gA_S^{1 \times n}$ avec $vu = 0$ et $u$ spécial.
Alors il existe $s \in S$, $\wi{u} \in \Ae {n \times 1}$, $\wi{v} \in \Ae {1
\times n}$ avec $\wi{v}\wi{u} = 0$, $\wi{u}$ spécial et~\hbox{$u = \wi{u}/s$}, $v = \wi{v}/s$ sur $\gA_S$.
\end{fact}
%%%%%%%%%%%%%%%%%%%%%%%%%%%%%%%%%%%%%%%%%

%%%%%%%%%%%%%%%%%%%%%%%%%%%%%%%%%%%%%%%%%
\begin {proof}
Par \dfnz, $u = u'/s_1$, $v = v'/s_1$ avec $s_1 \in S$, $u' \in
\Ae {n \times 1}$ et $v' \in \Ae {1 \times n}$. L'égalité $vu = 0$ fournit
un $s_2 \in S$ tel que $s_2v'u' = 0$, et $u_i = 0$ fournit un $s_3 \in S$ tel
que $s_3u'_i = 0$. Alors $s = s_1s_2s_3$, $\wi{u} = s_2s_3u'$ et $\wi{v} =
s_2s_3v'$ remplissent les conditions requises.
\end {proof}
%%%%%%%%%%%%%%%%%%%%%%%%%%%%%%%%%%%%%%%%%

%%%%%%%%%%%%%%%%%%%%%%%%%%%%%%%%%%%%%%%%%
%:     Theorem{lemE03Rao}  
\begin{theorem}\label{lemE03Rao}
Si $n \ge 3$, alors $\wi\EE(\Ae n)=\En(\gA)$.
En particulier, $\En(\gA)$ est stable par $\GLn(\gA)$-conjugaison.
\\
 {\emph{Précisions:} Soient $u \in \Ae {n \times 1}$, $v \in \Ae {1 \times n}$ avec $vu
= 0$ et $v$  \umdz. Alors, on peut écrire $u$ sous la forme $u = u'_1 + u'_2 + \cdots + u'_N, $  
avec $vu'_k = 0$ et chaque $u'_k$ a au plus deux composantes non nulles.
La matrice $\I_n + uv$ s'écrit alors comme un produit de matrices
spéciales:
$$\preskip.4em \postskip.4em
\I_n + uv =  (\I_n + u'_1v)\  (\I_n + u'_2v)\ \cdots\ (\I_n + u'_Nv) 
$$
et par conséquent, elle appartient à $\En(\gA)$.
}
\end {theorem}
%%%%%%%%%%%%%%%%%%%%%%%%%%%%%%%%%%%%%%%%%

%%%%%%%%%%%%%%%%%%%%%%%%%%%%%%%%%%%%%%%%%
\begin{proof} La base canonique de $\Ae n$ est notée $(e_1, \ldots, e_n)$.
On a $a_1$, \ldots, $a_n$ dans~$\gA$ tels que $a_1v_1 + \cdots + a_nv_n = 1$.
\\
Pour $i \leq j$, définissons $a_{ij} \in \gA$ par $a_{ij} = u_i a_j - u_j a_i$.
Alors:

\snic{u = \som_{i < j} a_{ij} (v_j e_i - v_i e_j) =
\som_{i \leq j} a_{ij} (v_j e_i - v_i e_j).}

%\sni
En effet, pour $k$ fixé, le coefficient de $e_k$ dans la somme de droite est
$$\preskip.4em \postskip.4em 
{\arraycolsep2pt\begin{array}{lllllll}
\dsp 
\sum_{j \ge k} a_{kj} v_j - \sum_{i < k} a_{ik} v_i  &  = &  
\dsp
\sum_{j \ge k}\, (u_ka_j - u_ja_k) v_j - \sum_{i < k}\, (u_ia_k - u_ka_i) v_i
 \\[5mm]
  & =  & 
  \dsp 
  u_k \som_{j = 1}^n a_jv_j - a_k \som_{j = 1}^n u_jv_j\;=\;
   u_k. 
  \end{array}} 
$$
Pour $i < j$, on définit alors $u'_{ij} \in \Ae {n \times 1}$ par $u'_{ij} =
a_{ij}(v_j e_i - v_i e_j)$. Il est clair que $u'_{ij}$ a au plus deux
composantes non nulles et que $vu'_{ij} = 0$. 
%Enfin, si $u', u'' \in \Ae {n
%\times 1}$ satisfont à $vu' = vu'' = 0$, on a $\I_n + (u'+u'')v = (\I_n +
%u'v) (\I_n + u''v)$ d'où la \dcn en produit de $\I_n + uv$.
\end {proof}
%%%%%%%%%%%%%%%%%%%%%%%%%%%%%%%%%%%%%%%%%

%%%%%%%%%%%%%%%%%%%%%%%%%%%%%%%%%%%%%%%%%%%%%%%%%%%%%%%%%%%%%%%%%%%
%%%%%%%%%%%%%%%%%%%%%%%%%%%%%%%%%%%%%%%%%%%%%%%%%%%%%%%%%%%%%%%%%%%
%%%%%%%%%%%%%%%%%%%%                                     %%%%%%%%%%%%%%%
%%%%%%%%%%%%%%%%%%%%        Symbole de Mennicke          %%%%%%%%%%%%%%%
%%%%%%%%%%%%%%%%%%%%                                     %%%%%%%%%%%%%%%
%%%%%%%%%%%%%%%%%%%%%%%%%%%%%%%%%%%%%%%%%%%%%%%%%%%%%%%%%%%%%%%%%%%
%%%%%%%%%%%%%%%%%%%%%%%%%%%%%%%%%%%%%%%%%%%%%%%%%%%%%%%%%%%%%%%%%%%

%---- Section{Symbole de Mennicke}--- 
\section{Le symbole de Mennicke} \label{secMennicke} 

%%%%%%%%%%%%%%%%%%%%%%%%%%%%%%%%%%%%%%%%%
\begin {lemma} \label {lemMennicke1}
Soient des \elts $a, b$ \com dans $\gA$. Alors la classe d'\eqvc dans
$\;\SL_3(\gA)/\EE_3(\gA)$ de la matrice $\cmatrix {a & b & 0\cr c & d & 0\cr 0 &
0 & 1\cr}$ ne dépend pas du choix de $c$ et $d$ vérifiant $1 = ad - bc$.\\
On notera $\meck {a}{b}$ l'élément de $\;\SL_3(\gA)/\EE_3(\gA)$ ainsi
obtenu. On l'appelle le \ix{symbole de Mennicke} de~$(a, b)$.
\end {lemma}

%%%%%%%%%%%%%%%%%%%%%%%%%%%%%%%%%%%%%%%%%
\begin {proof}
Soient $A = \cmatrix {a&b\cr c&d\cr}$, 
$A' = \cmatrix {a&b\cr c'&d'\cr}$ avec $ad - bc = ad' - bc' = 1$. Alors

\snic{
AA'^{-1} = \cmatrix {a&b\cr c&d\cr}\cmatrix {d'&-b\cr -c'&a\cr} =
\cmatrix {1 & 0\cr cd' - c'd& 1\cr} ,  
}

et $\bloc{A}{0_{2,1}}{0_{1,2}}{1}{\bloc{A'}{0_{2,1}}{0_{1,2}}{1}}^{-1}=\bloc{AA'^{-1}}{0_{2,1}}{0_{1,2}}{1}$ est dans $ \EE_3(\gA)$.
\end {proof}
%%%%%%%%%%%%%%%%%%%%%%%%%%%%%%%%%%%%%%%%%

%%%%%%%%%%%%%%%%%%%%%%%%%%%%%%%%%%%%%%%%%
\begin {proposition} \label {propMennicke1}
Le symbole de Mennicke vérifie les \prts suivantes.
\begin {enumerate}
\item
Si $a \in \Ati$, alors $\meck {a}{b} = 1$ pour tout~$b \in \gA$.
\item
Si $\gen{1} = \gen {a,b}= \gen {a',b}$ alors
$1 \in \gen {aa',b}$ et $\meck {aa'}{b} = \meck {a}{b}\meck{a'}{b}$.
\item
Si $1 \in \gen {a,b}$, alors $\meck {a}{b} = \meck {b}{a} = \meck{a+tb}{b}$
pour tout~$t \in \gA$.

\end {enumerate}
\end {proposition}
%%%%%%%%%%%%%%%%%%%%%%%%%%%%%%%%%%%%%%%%%

%%%%%%%%%%%%%%%%%%%%%%%%%%%%%%%%%%%%%%%%%
\begin {proof}
\emph{1.} 
La matrice $\cmatrix {a & b\cr 0 & a^{-1}\cr}$ appartient
à $\EE_2(\gA)$.

 \emph{2.} 
On a :
$$\preskip.0em \postskip.0em
\cmatrix {a & b & 0\cr c & d & 0\cr 0 &0 & 1\cr} \sims{\EE_3(\gA)} 
\cmatrix {a & 0 & b\cr 0 & 1 & 0\cr c &0 & d\cr},
$$
et
$$\preskip.0em \postskip.4em
\cmatrix {a' & b & 0\cr c' & d' & 0\cr 0 &0 & 1\cr} \sims{\EE_3(\gA)} 
\cmatrix {a' & 0 & -b\cr c' & 0 & -d'\cr 0 &1 & 0\cr} \sims{\EE_3(\gA)} 
\cmatrix {a' & 0 & -b\cr c' & 0 & -d'\cr 0 &1 & a\cr}.
$$
Le produit $\meck {a}{b}\meck{a'}{b}$ est représenté par le
produit des matrices de droite, i.e. par
$$
\cmatrix {aa' & b & 0\cr c' & 0 & -d'\cr ca' &d & 1\cr} \sims{\EE_3(\gA)} 
\cmatrix {aa' & b & 0\cr * & * & 0\cr ca' &d & 1\cr} \sims{\EE_3(\gA)} 
\cmatrix {aa' & b & 0\cr * & * & 0\cr 0 &0 & 1\cr},
$$
et donc $\meck {a}{b}\meck{a'}{b} = \meck {aa'}{b}$.

 \emph{3.} 
Si $ad - bc = 1$, alors $\cmatrix {a &b\cr c&d\cr} \sims{\EE_2(\gA)} 
\cmatrix {-b &a\cr -d&c\cr},$ et donc
$$
\meck {a}{b} = \meck{-b}{a} = \meck{-1}{a}\meck{b}{a} = \meck{b}{a}.
$$
Enfin, $\cmatrix {a &b\cr c&d\cr} \sims{\EE_2(\gA)} 
\cmatrix {a + tb & b\cr c+td & d\cr},$ donc $\meck {a}{b} = \meck{a+tb}{b}$.
\end {proof}

%--- Lemma{lemII3.6}------------
\begin{lemma}\label{lemII3.6}\relax \emph{(Version locale)}\\ Soit $\gA$ un \alo \dcd
et $f$, $g\in \AX $ \com avec $f$ \monz. Alors  on~a: 
$$\meck{f}{g}=\meck{f(0)}{g(0)}=1.$$
\end{lemma}
%--- end-lemma-----------------------------------------

\begin{proof}%(cf. \cite{GM})~\\
\'Ecrivons $af+bg=1$. Notons pour commencer que l'on  
peut diviser $b$ par~$f$ et que l'on  obtient alors une \egt
$a_1f+b_1g=1$ avec $\deg(b_1)<\deg(f)$, et donc, puisque $f$ est \monz,
$\deg(a_1)<\deg(g)$. 
Nous supposerons donc \spdg que $\deg(b)<\deg(f)$ et $\deg(a)<\deg(g)$.\\
Soit $r$ le reste de la division euclidienne de $g$ par $f$. 
Alors $\meck{f}{g}=\meck{f}{r}$. 
En particulier, si $\deg(f)=0$ on a terminé. Dans le cas contraire, on peut supposer  
$\deg(g)< \deg(f)$ et l'on raisonne par \recu sur $\deg(f)$. 
Puisque~$\gA$ 
est local \dcdz,  $g(0)\in\Ati$  ou $g(0)\in\fm=\Rad\gA$.\\
Supposons tout d'abord $g(0)$ inversible. Alors
%-----------------begin $$----------------
$$
 \meck{f}{g}=\meck{f-g(0)^{-1}f(0)g}{g},
$$
%-----------------end $$------------------
si bien que nous pouvons supposer $f(0)=0$ et $f=Xf_1$.
Alors 
%-----------------begin $$----------------
$$
 \meck{Xf_1}{g}=\meck{X}{g}\meck{f_1}{g}=\meck{X}{g(0)}\meck{f_1}{g}
=\meck{f_1}{g}
$$
%-----------------end $$------------------
et la preuve est terminée par \recu puisque $f_1$ est \monz.
\\
Supposons maintenant que $g(0)$ est dans  $\fm$.  Comme 
$a(0)f(0)+b(0)g(0)=1$, on a  $a(0)f(0)\in 1+\fm\subseteq \Ati$, et 
donc $a(0)\in\Ati$. Or
%-----------------begin $$----------------
$$ \left[\matrix{ 
  \phantom-f  &  g   &   0   \cr 
  -b  &  a   &   0   \cr
  \phantom-0  &  0   &   1
}\right]\equiv 
\left[\matrix{ 
  f-b  &  g+a   &   0   \cr 
  -b  &  a   &   0   \cr
  0  &  0   &   1
}\right]\;\;\mod\; \EE_3(\AX ),
$$
%-----------------end $$------------------
donc
%-----------------begin $$----------------
$$ 
\meck{f}{g}=\meck{f-b}{g+a},
$$
%-----------------end $$------------------
avec $f-b$ unitaire,  $\deg(f-b)=\deg (f)$, $\deg(g+a)<\deg(f)$ et 
$(g+a)(0)$ dans~\hbox{$\fm+\Ati=\Ati$}. On est donc ramené au cas 
précédent.
\end{proof}
%-----------------end proof------------------
Notre machinerie  \lgbe de base \vref{MethodeIdeps}, appliquée à la \dem locale précédente, donne le lemme 
quasi global suivant.\imlb  

%:--- Lemma{lemII3.6bis}----- quasi global ----
\begin{lemma} 
\label{lemII3.6bis}\relax \emph{(Version quasi globale)}\\
Soit $\gA$ un anneau  et $f$, $g\in \AX $ \com avec $f$ \monz. 
Alors, il existe dans $\gA$ un \sys d'\eco  $(s_i)$  tels que 
dans chaque localisé $\gA[1/s_i]$, on ait l'\egt des symboles de 
Mennicke suivante:
 $$
 \meck{f}{g}=\meck{f(0)}{g(0)}=1.
 $$
\end{lemma}
%--- end-lemma-----------------------------------------

%%%%%%%%%%%%%%%%%%%%%%%%%%%%%%%%%%%%%%%%%%%%%%%%%%%%%%%%%%%%%%%%%%%
%%%%%%%%%%%%%%%%%%%%%%%%%%%%%%%%%%%%%%%%%%%%%%%%%%%%%%%%%%%%%%%%%%%
%%%%%%%%%%%%%%%%%%%%                                     %%%%%%%%%%%%%%%
%%%%%%%%%%%%%%%%%%%%     Complements sur les \vmds       %%%%%%%%%%%%%%%
%%%%%%%%%%%%%%%%%%%%                                     %%%%%%%%%%%%%%%
%%%%%%%%%%%%%%%%%%%%%%%%%%%%%%%%%%%%%%%%%%%%%%%%%%%%%%%%%%%%%%%%%%%
%%%%%%%%%%%%%%%%%%%%%%%%%%%%%%%%%%%%%%%%%%%%%%%%%%%%%%%%%%%%%%%%%%%
  
%---- Section{Complements sur les \vmds \polsz}--- 
\section{Vecteurs \umds \pollsz} 
\label{secCompVmdsPols}

%:     Notation{notaGLIdeal}
\begin{notation}\label{notaGLIdeal} ~\\
{\rm Si $\fb$ est un \id de $\gB$, on note $\GLn(\gB,\fb)$ le sous-groupe de
 $\GLn(\gB)$ noyau du morphisme naturel $\GLn(\gB)\to\GLn(\gB\sur{\fb})$.
 Notation analogue pour $\SLn$. 
 Mais attention pour le groupe $\En$!
 On note~\hbox{$\En(\gB, \fb)$} le sous-groupe normal engendré par les $\rE_{ij}(b)$ avec
 $b \in \fb$. 
 }
\end{notation}

Le groupe~\hbox{$\En(\gB, \fb)$} est un sous-groupe du noyau de $\En(\gB)\to\En(\gB\sur{\fb})$, et en \gnlz, c'est un sous-groupe strict.  Cependant, dans le cas où $\gB =
 \gA[X]$ \hbox{et $\fb = \gen {X}$}, les deux groupes coïncident. Ce résultat 
est l'objet du lemme suivant.

%%%%%%%%%%%%%%%%%%%%%%%%%%%%%%%%%%%%%%%%%
%:     Lemma{lemE04Rao}
\begin{lemma}\label{lemE04Rao}
Le groupe $\En(\AX, \gen{X})$ est le noyau de l'\homo
canonique $\En(\gA[X])\to \En(\gA[X]\sur{\gen {X}}\!) = 
\En(\gA)$.
Il est engendré par les matrices du type $\gamma\, \rE_{ij}(Xg)\, \gamma^{-1}$
avec $\gamma \in \En(\gA)$ et $g \in \gA[X]$.
\end{lemma}
%%%%%%%%%%%%%%%%%%%%%%%%%%%%%%%%%%%%%%%%%

%%%%%%%%%%%%%%%%%%%%%%%%%%%%%%%%%%%%%%%%%
\begin {proof}
Soit $H$ ce noyau. On va utiliser la \dcn suivante, valide dans tout groupe, d'un
produit $\alpha_1\beta_1 \alpha_2\beta_2  \cdots
\alpha_m\beta_m$, par exemple avec $m=3$:
$$
\bigl( \alpha_1\beta_1\alpha_1^{-1} \bigr)\,
\bigl( (\alpha_1\alpha_2)\beta_2(\alpha_1\alpha_2)^{-1} \bigr)\,
%\bigl( (\alpha_1\alpha_2\alpha_3)\beta_3(\alpha_1\alpha_2\alpha_3)^{-1} \bigr)
\bigl( (\alpha_1\alpha_2\alpha_3)\beta_3(\alpha_1\alpha_2\alpha_3)^{-1} \bigr)
\, (\alpha_1\alpha_2\alpha_3).
$$
Soit donc $E = E(X) \in H$, $E = \prod_{i=1}^m
\rE_{i_k,j_k}(f_k)$ avec $f_k \in \AX$. \\
On écrit $f_k = c_k + Xg_k$
avec $c_k = f_k(0) \in \gA$ et
$$
\rE_{i_k,j_k}(f_k) = \alpha_k \beta_k, \quad \hbox {avec} \quad
\alpha_k = \rE_{i_k,j_k}(c_k), \quad
\beta_k = \rE_{i_k,j_k}(Xg_k). 
$$
On termine en appliquant  la \dcn donnée ci-dessus et en utilisant l'\egtz~\hbox{$\alpha_1\cdots\alpha_m = E(0) = \I_n$}.
\end {proof}

%%%%%%%%%%%%%%%%%%%%%%%%%%%%%%%%%%%%%%%%%
%:     Proposition{propE0Rao}
\begin{proposition}\label{propE0Rao}
Soient $n \ge 3$, $s \in \gA$ et $E = E(X) \in \En(\gA_s[X],
\gen{X})$. Il existe $k \in \NN$ et $E' = E'(X) \in \En(\AX,
\gen{X})$ vérifiant $E'(X) = E(s^k X)$ sur~$\gA_s[X]$.
\end {proposition}
%%%%%%%%%%%%%%%%%%%%%%%%%%%%%%%%%%%%%%%%%
\begin {proof}
On peut supposer $E = \gamma \rE_{ij}(Xg) \gamma^{-1}$
avec $\gamma \in \En(\gA_s)$ et $g \in \gA_s[X]$. En notant
$u \in \gA_s^{n \times 1}$ la colonne $i$ de $\gamma$ et
$v \in \gA_s^{1 \times n}$ la ligne $j$ de $\gamma^{-1}$, on~a:
$$
E(X) = \gamma \rE_{ij}(Xg) \gamma^{-1} = \I_n + (Xg) uv,
\quad vu = 0, \quad \hbox {$v$ unimodulaire}.
$$
Le \thrf{lemE03Rao} permet d'écrire:
$
u = u'_1 + u'_2 + \cdots + u'_N $ 
 {avec $vu'_k = 0$ et~\hbox{$u'_k \in \gA_s^{n \times 1}$} 
a au plus deux composantes non nulles}.
On a donc
$$
E(X) = (\I_n + (Xg) u'_1v)\  (\I_n + (Xg) u'_2v)\ \cdots\ (\I_n + (Xg) u'_Nv). 
$$
En utilisant une méthode analogue au fait \ref{factE02Rao}, 
on vérifie facilement qu'il
existe $k \in \NN$, $\wi {g} \in \AX$, $\wi u_k \in \Ae {n \times 1}$ et
$\wi {v} \in \Ae {1 \times n}$ tels que l'on ait sur $\gA_s$ les \egts
$
g = \wi{g}/s^k $, $   u'_k = \wi u_k /s^k$,  $ 
v = \wi{v}/s^k$,  $\wi{v} \wi u_k = 0$,  
et $\wi u_k$ a au plus deux composantes non nulles.
On pose alors :
$$
E'(X) = (\I_n + (X\wi{g}) \wi u_1\wi{v})\  
(\I_n + (X\wi{g}) \wi u_2 \wi{v})\ \cdots\ (\I_n + (X\wi{g}) \wi u_N \wi{v}) .
$$
D'après le corolaire \ref{corlemE01Rao}, 
chaque $\I_n + (X\wi{g}) \wi u_k \wi{v}$ appartient
à $\En(\AX)$.  On a donc $E'(X) \in \En(\AX)$, $E'(0) = \I_n$ et
$E'(s^{3k}X) = E(X)$ sur $\gA_s[X]$.
\end {proof}
%%%%%%%%%%%%%%%%%%%%%%%%%%%%%%%%%%%%%%%%%

%
%:     Lemma{lemE1Rao}
\begin{lemma}\label{lemE1Rao}
Soit un entier $n\geq3$, $s\in\gA$ et $E=E(X)\in\En(\gA_s[X])$. 
\\ 
Il existe un entier $k\geq0$ tel que pour tous
$a , \, b\in\gA$ congrus modulo $s^k$, la matrice $E^{-1}(aX)E(bX)$ est dans 
l'image de l'\homo naturel 

\snic{\En(\AX,\gen{X})\longrightarrow\En(\gA_s[X],\gen{X}).}
\end{lemma}
NB: en bref, mais de manière moins précise,
 si $a$ et $b$ sont suffisamment \gui{proches}, la matrice 
$E^{-1}(aX)\,E(bX)$ n'a plus de \denoz.
%---  end lemma ---------------------------------------------
\begin{proof}
On introduit deux nouvelles \idtrs $T$, $U$ et l'on pose 
$$
E'(X,T,U)=E^{-1}\big((T+U)X\big)\,E(TX).
$$ 
On a $E'(X,T,0)=\In$. 
On applique la proposition \ref{propE0Rao} avec $F=E'$ en prenant $\gA[X,T]$
au lieu de $\gA$ et $U$ au lieu de $X$: il existe une matrice~$G$ dans $\En(\gA[X,T,U],\gen{U})$ et un entier $k\geq0$ tels que 

\snic{E'(X,T,s^kU)=G(X,T,U) \;\hbox{ dans }\;  \En(\gA_s[X,T,U],\gen{U}).} 

%\sni
Donc $G(X,T,U)=E^{-1}\big((T+s^kU)X\big) \, E(TX)$ sur $\gA_s$, et
si $b=a+s^kc$:
$$\preskip.3em \postskip.4em 
E^{-1}(aX)\,E(bX)=G(X,a,c) \; \hbox{ sur } \; \gA_s . 
$$
On a $G(0,T,U)=\In$ sur $\gA_s$, mais pas \ncrt sur  $\gA$.
On pose: 

\snic{H(X,T,U)=G^{-1}(0,T,U)\,G(X,T,U).}

%\sni
On a alors $H(0,T,U)=\In$ sur $\gA$
et $H(X,T,U)=G(X,T,U)$ sur $\gA_s$.
On obtient donc

\snic{E^{-1}(aX)\,E(bX)=H(X,a,c) \;\;\; \mathrm{dans} \;\;
\; \En(\gA_s[X],\gen{X}),}

%\sni
avec $H(X,a,c)\in \En(\AX,\gen{X})$.
\end{proof}
%

%:     Lemma{lemE2Rao}
\begin{lemma}\label{lemE2Rao}
Soit un entier $n\geq3$, $s\in\gA$ et 
$$\preskip.4em \postskip.3em
E=E(X)\in\GLn(\gA[X])\cap \En(\gA_s[X]).
$$ 
Il existe un entier $k\geq0$ tel que pour tous
$a , \, b\in\gA$ congrus modulo $s^k$, la matrice $E^{-1}(aX)E(bX)$ 
est dans~$\En(\gA[X],\lra X)$. 
\end{lemma}
\facile
%

%:     Lemma{lemE3Rao}
\begin{lemma}\label{lemE3Rao}
Soit un entier $n\geq3$, $s$, $t$ \com dans $\gA$ et 
$$\preskip.4em \postskip.4em
E\in\GLn(\gA[X],\lra X)\cap \En(\gA_s[X]) \cap \En(\gA_t[X])
.$$
Alors~$E\in \En(\gA[X])$.
\end{lemma}
\begin{proof}
D'après le lemme \ref{lemE2Rao}, il existe un $k$
tel que   pour tous
$a , \, b\in\gA$ congrus modulo $s^k$, ou modulo $t^k$, 
la matrice $E^{-1}(aX)E(bX)$ est dans 
$\En(\gA[X],\lra X)$. Soit $c\in\gA$ tel que $c\equiv 0 \mod s^k$ 
et $c\equiv 1 \mod t^k$. \\
Alors on écrit  $E=E^{-1}(0\cdot X)\,E(c\cdot X)\,E^{-1}(c\cdot X)\,E(1\cdot X)$. 
\end{proof}
%

%

%%%%%%%%%%%%%%%%%%%%%%%%%%%%%%%%%%%%%%%%%%%%%%%%%%%%%%%%%%%%%%%%%%%
%%%%%%%%%%%%%%%%%%%%%%%%%%%%%%%%%%%%%%%%%%%%%%%%%%%%%%%%%%%%%%%%%%%
%%%%%%%%%%%%%%%%%%%%                                     %%%%%%%%%%%%%%%
%%%%%%%%%%%%%%%%%%%%     Principe \lgb de Rao            %%%%%%%%%%%%%%%
%%%%%%%%%%%%%%%%%%%%                                     %%%%%%%%%%%%%%%
%%%%%%%%%%%%%%%%%%%%%%%%%%%%%%%%%%%%%%%%%%%%%%%%%%%%%%%%%%%%%%%%%%%
%%%%%%%%%%%%%%%%%%%%%%%%%%%%%%%%%%%%%%%%%%%%%%%%%%%%%%%%%%%%%%%%%%%

%---- Section{Principes \lgb de Suslin et Rao}--- 
\section{Principes \lgbs  de Suslin et Rao} 
\label{secPlgbRao}

Maintenant nous démontrons \cite[lemme I 5.9 page 26]{GM}.
%:     Theorem{lem159}
\begin{theorem}\label{lem159}
Soit $n\geq 3$ et $A = A(X)\in \GL_n(\AX)$.
\begin{enumerate}
\item 
Si $A(0) = \In$, alors l'ensemble
$$\preskip.2em \postskip.0em
\fa=\sotq{s \in \gA}{  A\in \En(\gA_s[X]) }
$$ 
est un idéal de $\gA$. 

\item 
L'ensemble
$
\fa = \sotQ {s\in\gA}{A(X)\sims{\En(\gA_s[X])} A(0)}
$
est un idéal de $\gA$.
\end{enumerate}
\end{theorem}

\begin{proof}
Les deux formulations sont \eqvesz; on démontre la seconde
à partir de la première en considérant la matrice $A(X)A(0)^{-1}$.
\\ 
\emph{1.}  Il est clair que $s \in \fa
\Rightarrow as \in \fa$ pour tout $a \in \gA$.  Soient maintenant $s$, $t$
dans $\fa$. On doit montrer que $s+t\in\fa$, ou encore que $1\in\fa\gA_{s+t}$. En bref, on suppose que $s$ et $t =1-s$ sont dans
$\fa$, et l'on doit montrer que $1\in\fa$.
\\ 
Par \dfnz, on a $A \in\En(\gA_{s}[X])$ et $A \in\En(\gA_{t}[X])$; d'après le
lemme~\ref{lemE3Rao}, on a $A \in\En(\gA[X])$, i.e.  $1 \in \fa$.
\end {proof}

Ce lemme aurait pu être écrit sous la forme du \plgc suivant  (à très peu près \cite[lemme I 5.8]{GM}).

%--- Prc loc glob conc {plcc.159}     
\begin{plcc} 
\label{plcc.159}\relax Soient $n\geq 3,$  $S_1$, $\ldots$, $S_k$  des 
\moco de $\gA$ et $A\in \GL_n(\AX)$, avec $A(0)=\In$. Alors:  
$$\preskip.3em \postskip.4em
A\in \En(\AX )\;\;\;\;\Longleftrightarrow\;\;\;\; 
\hbox{ pour }\;i\in\lrbk,\quad A\in \En(\gA_{S_i}[X]).
$$
\end{plcc}
%---------------------------------

Le \tho suivant reprend \cite[corolaire II 3.8]{GM}.
%:--- theorem{II3.8}-----------  
\begin{theorem}\label{II3.8}\relax 
\emph{(Version globale du lemme \ref{lemII3.6})}
\\ 
Soient $n\geq 3,$    et $f$, $g\in \AX $ 
\comz, avec $f$ \monz. 
Alors, on a l'\egt des symboles de Mennicke suivante:
 $\meck{f}{g}=\meck{f(0)} {g(0)}.$ 
\end{theorem}
%--- end-theorem------------------------------------
%-----------------begin proof------------------
\begin{proof}
\'Ecrivons $af-bg=1$. Soit $B=\cmatrix {f & g & 0\cr b & a & 0\cr 0 &
0 & 1\cr}$. \\
L'\egt $\meck{f}{g}=\meck{f(0)}{g(0)}$ signifie:
$A=BB(0)^{-1}\in \EE_3(\AX )$. On a évidemment $A(0)=\I_3$. Le \plgc 
\ref{plcc.159}
nous dit qu'il suffit de vérifier l'assertion après \lon en des \ecoz~$(s_i)$. Et le lemme~\ref{lemII3.6bis}  
a construit une telle famille.
\end{proof}
%-----------------end proof------------------

%: corollary KuXMennickeTrivial
\begin{corollary} 
\label{KuXMennickeTrivial}
\emph{(Trivialité du symbole de Mennicke sur $\KuX$)}\\
Soit $\gK$ un corps discret, et  $f$, $g\in \KuX$ \comz.
Alors $\meck{f}{g} = 1$.
\end{corollary}

\begin{proof}
On raisonne par \recu sur le nombre  $r$ de variables dans $\uX$. \\
Le cas $r=0$, i.e.  $\KuX = \gK$ découle de
$\EE_3(\gK) = \SL_3(\gK)$ ($\gK$ est un corps discret).  Pour $r \ge 1$, on
suppose \spdg que $f$ est non nul. Un \cdv permet de transformer $f$ en un \pol
pseudo \mon en $X_r$ (lemme \ref{lemNoether}), disons $f = ah$ avec $a\in\gK^*$
et $h$ \mon en~$X_r$. Alors, en posant $h_0 = h(X_1, \ldots, X_{r-1}, 0)$ et
$g_0 = g(X_1, \ldots, X_{r-1}, 0)$, qui sont dans $\gK[X_1, \ldots, X_{r-1}]$,
on a $\meck{f}{g} = \meck{h}{g} = \meck{h_0}{g_0}$.
\end{proof}

Dans la fin de cette section, les résultats sont démontrés dans le
cas d'un anneau intègre. Ils sont en fait vrais pour un anneau arbitraire. 
Pour le cas \gnlz, il faut se reporter
à \cite{Rao85a,Rao85b,Rao85c}.

%:     Theorem{th3Rao}
\begin{theorem}\label{th3Rao}
Soit $n\geq3$,  $\gA $ un anneau intègre et 
$f(X)$ un \vmd dans $\AX^n$,
alors l'ensemble
$$\preskip.1em \postskip.2em
\fa=\sotQ{s\in\gA}{f(X)\sims{\En(\gA_s[X])} f(0)}
$$
est un \idz.
\end{theorem}

On exprime la même chose dans le \plgc suivant.

%:    Principe \lgb concret de Rao
\CMnewtheorem{plgcRao}{Principe \lgb concret de Rao}{\itshape}
%:     Lemma{lem}
\begin{plgcRao}\label{plgc-Rao}\relax
Soit $n\geq3$,  $\gA $ un anneau intègre,~$f(X)$ un \vmd dans 
$\AX^n$, et  $S_1$, $\ldots$, $S_k$ des \moco de $\gA$. 
\Propeq 
\begin{enumerate}
\item ${f(X)}\sims{\En(\AX)} f(0)$.
\item ${f(X)}\sims{\En(\gA_{S_i}[X])} f(0)$ pour chaque $i$.
\end{enumerate}
\end{plgcRao}
%--------- fin lemma ---------------------------------------------- 

%
\begin{Proof}{\Demo du \thrf{th3Rao}. }\\
On doit montrer que l'ensemble
$$\preskip.0em \postskip.1em 
 \fa=\sotQ{s\in\gA}{ f(X)\sims{\En(\gA_s[X])} f(0)}
$$
est un \idz. Puisque tous les calculs dans $\gA_s$ sont valables
dans $\gA_{sa}$, on~a:~$s\in\fa$ implique
$ as \in \fa$. Soit maintenant
$s_1$ et $s_2$ dans $\fa$. On doit montrer que $s_1+s_2\in\fa$, ou encore
que~$1\in\fa\gA_{s_1+s_2}$. En bref, on suppose
que $s_1$ et $s_2=1-s_1$ sont dans $\fa$, et l'on doit montrer que~$1\in\fa$.
\\ 
Par \dfnz, pour $i=1$, $2$, on a une matrice $E_i=E_i(X)\in\En(\gA_{s_i}[X])$
telle que $E_if(X)=f(0).$ Comme $E_i(0)f(0)=f(0)$, quitte à remplacer~$E_i$ par $E_i^{-1}(0)E_i$, on peut supposer que $E_i(0)=\In$.
\\ 
On introduit $E=E_1E_2^{-1}\in \En(\gA_{s_1s_2}[X],\gen{X})$, ce qui
donne un entier $k\geq0$ satisfaisant la conclusion du lemme \ref{lemE1Rao}
pour la matrice $E$ et pour les deux \lons $\gA_{s_1}\to \gA_{s_1s_2}$ et $\gA_{s_2}\to \gA_{s_1s_2}$.
\\
Soit $c\in \gA$ avec $c\equiv 1\mod s_1^k$ et $c\equiv 0\mod s_2^k$.  
On a donc deux matrices $E'_1\in \En(\gA_{s_1}[X],\lra X)$,
 $E'_2\in \En(\gA_{s_2}[X],\lra X)$ qui vérifient:
\begin{enumerate}
\item [--]  $E^{-1}(cX)E(X)=E'_2$ 
sur $\gA_{s_1s_2}$ (puisque $c\equiv 1 \mod s_1^k$),
\item [--]  $E(cX)=E(cX)E(0\cdot X)=E'_1$ sur $\gA_{s_1s_2}$ (puisque $c\equiv 0 \mod s_2^k$).
\end{enumerate}
On obtient  $E=E'_1E'_2=E_1E_2^{-1}$ sur  $\gA_{s_1s_2}$,
et  ${E'_1}^{-1}E_1=E'_2E_2$ sur  $\gA_{s_1s_2}$.
Puisque ${E'_1}^{-1}E_1=F_1$ est définie sur $\gA_{s_1}$, que $E'_2E_2=F_2$
est définie sur~$\gA_{s_2}$, qu'elles sont égales sur $\gA_{s_1s_2}$,
et que~$s_1$ et $s_2$ sont \comz, il existe une unique matrice $F\in\Mn(\AX)$
qui donne $F_1$ sur $\gA_{s_1}$ et $F_2$ sur~$\gA_{s_2}$.
Il faut encore vérifier que $F\in\En(\AX)$ et~$Ff=f(0)$.
\\ 
Le premier point résulte du lemme \ref{lemE3Rao}.
Pour vérifier $Ff=f(0)$, on va supposer $\gA$ intègre,
ce qui légitime les \egts suivantes sur $\gA$ 
\[ \preskip.3em \postskip-.3em
\arraycolsep2pt
\begin{array}{cccccccc} 
Ff  & =  & {E'_1}^{-1}E_1f &= & {E'_1}^{-1}f(0) & =&\\[.25em] 
E^{-1}(cX)f(0)  & =  & E_2(cX)E_1^{-1}(cX)f(0) &= &E_2(cX)f(cX) &=&f(0).  
 \end{array}
\]
\end{Proof}
NB: dans cette \dem
la dernière vérification est le seul endroit où nous avons besoin
de supposer l'anneau intègre.

%:     Theorem{th2Rao}
\begin{theorem}\label{th2Rao}
Soit $n\geq3$, $\gA $ un anneau et 
$f=\tra{\big(f_1(X),\ldots,f_n(X)\big)}$ un \vmd dans $\AX^n$,
avec $1$ dans l'\id de tête des $f_i$. Alors: 
$$
f\;\sims{\En(\AX)} \;f(0)\;\sims{\En(\gA)}\;f\sta(0)\;\sims{\En(\AX)}\; f\sta.
$$
Si l'un des $f_i$ est \monz, on a $f\sims{\En(\AX)}  
 \tra{[\,1\;0\;\cdots\;0\,]}$. 
\end{theorem}
\begin{proof}
Le petit \tho de Horrocks local (\thrf{th2HorrocksLocal})
et le \plg de Rao donnent la première \eqvcz. Ensuite 
on recopie la \dem du \tho de Rao (\thrf{th1Rao}) 
en remplaçant~$\GLn$ par $\En$.
\end{proof}
%

%:      cor TransitiviteEn
\begin{corollary} 
\label{TransitiviteEn}
\emph{(Transitivité de $\En$ pour $n \ge 3$)}
\\
Si $\gK$ est un \cdi et $\KuX = \gK[\Xr]$, 
alors $\En(\KuX)$ agit transitivement sur l'ensemble
des \vmds de $\KuX^n$\linebreak pour~$n\geq 3$.
\end{corollary}

\begin{proof}
On raisonne par \recu sur $r$. Le cas $r=0$ découle
du fait que $\gK$ est un corps discret. \\
Soit $r \ge 1$ et  $f =
\tra{[\,f_1(\uX)\;\cdots\;f_n(\uX)\,]}$ un \vmd de~$\KuX^n$.  
Notons $\gA = \gK[X_1, \ldots, X_{r-1}]$.
L'un des $f_i$
est non nul et un \cdv permet de le transformer  
en un \pol pseudo \mon en~$X_r$ (lemme~\ref{lemNoether}).  Avec $f_i$ \mon en~$X_r$, nous appliquons le \thrf {th2Rao} pour obtenir

\snic{f\;\sims{\En(\KuX)}\; f(X_1, \ldots, X_{r-1},0).}

%\sni
Ce dernier vecteur
est un \vmd de $\Ae n$. On applique  l'\hdrz. 
\end{proof}
%

%%%%%%%%%%%%%%%%%%%%%%%%%%%%%%%%%%%%%%%%%%%%%%%%%%%%%%%%%%%%%%%%%%%
%%%%%%%%%%%%%%%%%%%%%%%%%%%%%%%%%%%%%%%%%%%%%%%%%%%%%%%%%%%%%%%%%%%
%%%%%%%%%%%%%%%%%%                                       %%%%%%%%%%%%%%%
%%%%%%%%%%%%%%%%%%           Le cas des corps            %%%%%%%%%%%%%%%
%%%%%%%%%%%%%%%%%%                                       %%%%%%%%%%%%%%%
%%%%%%%%%%%%%%%%%%%%%%%%%%%%%%%%%%%%%%%%%%%%%%%%%%%%%%%%%%%%%%%%%%%
%%%%%%%%%%%%%%%%%%%%%%%%%%%%%%%%%%%%%%%%%%%%%%%%%%%%%%%%%%%%%%%%%%%

Enfin, la preuve que le \thrf{II3.8} implique le \tho de stabilité de 
Suslin est simple et \covz, comme  dans \cite{GM}.

%\penalty-2500
%:     Theorem{thStabSuslCorps}
\begin{theorem}\label{thStabSuslCorps}
 \emph{(Théorème de stabilité de \Susz, cas des corps discrets)}\\
Soit $\gK$ un corps discret.  Pour $n\geq 3$, on a 
$\SL_n(\KuX)= \En(\KuX)$.
\end{theorem}

\begin{proof}
Montrons le résultat préliminaire suivant.\\
\emph{Pour $A \in
\GL_n(\KuX)$, il existe $P$, $Q\in \En(\KuX)$ telles que} 

\snic{P\,A\,Q \in \GL_2(\KuX)
\subseteq  \GL_n(\KuX) (\footnote{L'inclusion $\GL_r \hookrightarrow \GL_n$
est définie comme d'habitude par $B \mapsto \Diag(B,\I_{n-r})$.}).}

%\sni
En effet,
considérons la dernière ligne de $A$. C'est un \vmdz, donc (corolaire~\ref{TransitiviteEn}),
il existe $Q_n \in \En(\KuX)$ telle que la dernière
ligne de~$A\,Q_n$ soit $[\,0\;\cdots\;0\;1\,]$.
%:HHH je change l'argument, ce n'est pas: De meme 
%De même, il existe
%$P_n \in \En(\KuX)$ telle que la première colonne de $P_n(A\,Q_n)$ soit
%$\tra{[\,0\;\cdots\;0\;1\,]}$, 
D'où ensuite $P_n \in \En(\KuX)$ telle que la dernière colonne de $P_n(A\,Q_n)$ soit
$\tra{[\,0\;\cdots\;0\;1\,]}$, i.e.  $P_n\,A\,Q_n \in \GL_{n-1}(\KuX)$. 
\\
 En itérant,
on  trouve des matrices $P$, $Q \in \En(\KuX)$ de la forme 

\snic{P = P_3\cdots P_n$, $\;Q =
Q_n\cdots Q_3,}

%\sni
telles que $P\,A\,Q \in \GL_2(\KuX)$.
\\ 
Si de plus $A \in \SLn(\KuX)$, on obtient $P\,A\,Q \in \SL_2(\KuX)
\hookrightarrow \SL_3(\KuX)$. 
\\
On peut~alors considérer son image dans
$\SL_3(\KuX)\sur{\EE_3(\KuX)}$.
\\
Comme le symbole de Mennicke correspondant
vaut $1$ (corolaire~\ref{KuXMennickeTrivial}), on obtient $P\,A\,Q \in \EE_3(\KuX)$, et
en fin de compte $A \in \EE_n(\KuX)$.
\end{proof}
%

%%%%%%%%%%%%%%%%%%%%%%%%%%%%%%%%%%%%%%%%%%%%%%%%%%%%%%%%%%%%%%%%%%%
%%%%%%%%%%%%%%%%%%%%%%%%%%%%%%%%%%%%%%%%%%%%%%%%%%%%%%%%%%%%%%%%%%%
%%%%%%%%%%%%%%%%%%%%                                     %%%%%%%%%%%%%%%%%%
%%%%%%%%%%%%%%%%%%%%     Miracles susliniens             %%%%%%%%%%%%%%%%%%
%%%%%%%%%%%%%%%%%%%%                                     %%%%%%%%%%%%%%%%%%
%%%%%%%%%%%%%%%%%%%%%%%%%%%%%%%%%%%%%%%%%%%%%%%%%%%%%%%%%%%%%%%%%%%
%%%%%%%%%%%%%%%%%%%%%%%%%%%%%%%%%%%%%%%%%%%%%%%%%%%%%%%%%%%%%%%%%%%

\incertain{  
%---- Section{Miracles susliniens}--- 
%\section{Miracles susliniens}
\mni 
{\bf \Large 17.5 ~~~ Miracles susliniens} 
\label{secMiracleSus}
\hum{très incertain}
}
%%%%%%%%%%%%%%%%%%%%%%%%%%%%%%%%%%%%%%%%%%%%%%%%%%%%%%%%%%%%%%%%%%%
%%%%%%%%%%%%%%%%%%%%%%%%%%%%%%%%%%%%%%%%%%%%%%%%%%%%%%%%%%%%%%%%%%%
%%%%%%%%%%%%%%%%%%%%%%%%%%%                        %%%%%%%%%%%%%%%%%%%%%%%%
%%%%%%%%%%%%%%%%%%%%%%%%%%%        Exercices       %%%%%%%%%%%%%%%%%%%%%%%%
%%%%%%%%%%%%%%%%%%%%%%%%%%%                        %%%%%%%%%%%%%%%%%%%%%%%%
%%%%%%%%%%%%%%%%%%%%%%%%%%%%%%%%%%%%%%%%%%%%%%%%%%%%%%%%%%%%%%%%%%%
%:   Exercices   
\Exercices{

%%%%%%%%%%%%%%%%%%%%%%%%%%%%%%%%%%%%%%%%%
%--- Exercise{exoChap18-1}-------------
\begin{exercise}
\label{exoChap18-1} 
{\rm  
Soient $U \in \Ae {n \times m}$ et $V \in \Ae {m \times n}$.
\begin{itemize}
\item [\emph{1.}]
Vérifier, pour $N \in \MM_n(\gA)$, que

\snic{(\I_m - VNU) (\I_m + VU) = \I_m + V\big( \I_n - N(\I_n + UV)\big)U
}

En déduire que si $\I_n + UV$ est inversible d'inverse $N$, alors
$\I_m + VU$ est inversible d'inverse $\I_m - VNU$.

\item [\emph{2.}]
En déduire que $\I_n + UV$ est inversible si et seulement si $\I_m + VU$
l'est et établir des formules symétriques pour leurs inverses. 

\item [\emph{3.}]
Montrer  que $\det(\I_n + VU) = \det(\I_m + UV)$ dans tous les cas.

\item [\emph{4.}]
On suppose que $\I_m + VU$ est inversible. Montrer l'appartenance suivante
due à Vaserstein.

\snic{\cmatrix {\I_n + UV & 0\cr 0 & (\I_m + VU)^{-1} \cr} \in \EE_{n+m}(\gA).
}

\smallskip  Que se passe-t-il lorsque $VU = 0$? 

\end{itemize}
}
\end{exercise}
%--- end -exercise-----------------------------------------

%%%%%%%%%%%%%%%%%%%%%%%%%%%%%%%%%%%%%%%%%
\begin {exercise}\label{exoMennicke}
{\rm 
Avec les notations du lemme \ref {lemMennicke1}, vérifier que
la matrice $A'^{-1} A$
est de la forme $\I_2 + uv$ avec $u, v \in \Ae {2 \times 1}$, $vu = 0$ et $v$
\umdz.
}
\end {exercise}

%%%%%%%%%%%%%%%%%%%%%%%%%%%%%%%%%%%%%%%%%
\begin {exercise}\label{exoMennicke2}
{\rm 
Soient $a$, $b$, $u$, $v \in \gA$ vérifiant $1 = au + bv$.  Montrer, en utilisant
uniquement les \prts du symbole de Mennicke figurant dans la
proposition~\ref {propMennicke1}, que $\meck {a}{b} = \meck {u}{v} = \meck
{a-v}{b+u}$.
}
\end {exercise}

%%%%%%%%%%%%%%%%%%%%%%%%%%%%%%%%%%%%%%%%%
\begin {exercise}\label{exoStabFreeCoRang1}
{\rm 
Un \Amo $E$ \stl de rang $r$ est dit \emph{de type $t$} \linebreak
si~$E \oplus \Ae t \simeq \Ae {r+t}$. On s'intéresse ici aux
relations entre d'une part les classes d'\iso des modules \stls de rang
$n-1$, de type $1$, et d'autre part le $\GLn(\gA)$-ensemble $\Um_n(\gA)$ constitué des
\vmds de $\Ae n$.
\begin {enumerate}\itemsep0pt

\item
Soit $x \in \Um_n(\gA)$. Vérifier que le module $\Ae n / \gA x$ est
\stl de rang $n-1$, de type $1$, et que pour $x' \in \Um_n(\gA)$,
on a $\Ae n / \gA x \simeq \Ae n / \gA x'$ \ssi $x \sims{\GL_n(\gA)}  x'$. 
Montrer que l'on obtient ainsi une (première)
correspondance bijective: $x \longleftrightarrow \Ae n / \gA x$

\vbox {\advance \vsize by -4cm
\Grandcadre{$
\displaystyle{\Um_n(\gA) \over \GL_n(\gA)} \buildrel {(1)} \over \simeq
{\hbox {modules \stls de rang $n-1$, de type $1$} 
\over \hbox {\iso}}
$}}

Quels sont les \vmds qui correspondent à un module
libre?

\item
Soit $x \in \Um_n(\gA)$. Montrer que $x^\perp \eqdefi  \Ker
\tra {x}$ est un module \stl de rang $n-1$, de type $1$, et que
pour $x' \in \Um_n(\gA)$, on a $x^\perp \simeq x'^\perp$ si et seulement si $x
\sims{\GL_n(\gA)} x'$.  Vérifier que l'on obtient ainsi une
(deuxième) correspondance bijective: $x \longleftrightarrow x^\perp$

\vbox {\advance \vsize by -4cm
\Grandcadre{$
\displaystyle{\Um_n(\gA) \over \GL_n(\gA)} \buildrel {(2)} \over \simeq
{\hbox {modules \stls de rang $n-1$, de type $1$} 
\over \hbox {\iso}}
$}}

\item
Si $E$ est \stl de rang $r$, de type $t$, il en est de même de
son dual $E\sta$. Pour $t = 1$, décrire l'involution de
$\Um_n(\gA)/\GLn(\gA)$ induite par l'involution $E \leftrightarrow E\sta$.

\item
Soient $x$, $x'$, $y \in \Ae n$ tels que $\tra {x}y = \tra {x'}y = 1$. Pourquoi
a-t-on $x \sims{\GL_n(\gA)}  x'$? Expliciter $g \in \GL_n(\gA)$
tel que $gx = x'$, $g$ de la forme $\I_n + uv$ avec 
$vu = 0$ et~$v$ \umdz. En déduire que pour
$n \ge 3$, $g \in \En(\gA)$, et donc $x \sims{\En(\gA)}  x'$.

\end {enumerate}

}
\end {exercise}

%%%%%%%%%%%%%%%%%%%%%%%%%%%%%%%%%%%%%%%%%
\begin {exercise}\label{exoStabFree2}
       (Modules \stls de type 1 autoduaux)
{\rm 
\begin {enumerate}\itemsep0pt
\item
Soient $a$, $b \in \gA$, $x = (x_1, \ldots, x_n) \in \Ae n$ avec $n \ge 3$ et
$ax_1 + bx_2$ inversible modulo $\gen {x_3, \ldots, x_n}$ (en
particulier, $x$ est \umdz). \\
On pose $x' = (-b,a, x_3, \ldots, x_n)$.
Expliciter  $z \in \Ae n$ tel que $\scp {x}{z} = \scp {x'}{z} = 1$.
En déduire, pour $G = \GLn(\gA)$ (ou mieux pour $G = \En(\gA)$) que

\snic{
x \,\sims G \, x'\, \sims G \, (a,b, x_3, \ldots, x_n).
}

\item
Soient $x, y \in \Ae 4$ tels que $\scp {x}{y} = 1$.  Montrer que $x \sims
{\EE_4(\gA)}  y$. En particulier, le module \stl
$x^\perp = \Ker \tra {x}$ est isomorphe à son dual.

\item
Question analogue à la précédente en remplaçant $4$ par
n'importe quel nombre pair $n \ge 4$. 

\end {enumerate}
}

\end {exercise}

%%%%%%%%%%%%%%%%%%%%%%%%%%%%%%%%%%%%%%

%%%%%%%%%%%%%%%%%%%%%%%%%%%%%%%%%%%%%%%%%

%%%%%%%%%%%%%%%%%%%%%%%%%%%%%%%%%%%%%%%%%%%%%%%%%%%%%%%%%%%%%%%%%%%
%%%%%%%%%%%%%%%%%%%%%%%%%%%%%%%%%%%%%%%%%%%%%%%%%%%%%%%%%%%%%%%%%%%
%%%%%%%%%%%%%%%%%%%%%%%%%%%%%%%%%%%%%%%%%%%%%%%%%%%%%%%%%%%%%%%%%%%
}% fin des exos
%:  solutions
%\penalty-2500
\sol

%%%%%%%%%%%%%%%%%%%%%%%%%%%%%%%%%%%%%%%%%
\exer{exoChap18-1}{
\emph{2}. On établit les formules 

\snic{N = (\I_n + UV)^{-1} = \I_n - UMV, $ $  M = (\I_m + VU)^{-1} = \I_m - VNU}

%\sni
\emph{4}.
On sait que $\I_n + UV$ est \ivz; on note $N = (\I_n + UV)^{-1}$,
$M  = (\I_m + VU)^{-1}$. On a donc $N + UVN = \I_n = N + NUV$ et
$M + VUM = \I_m = M + MVU$. \\
On réalise les opérations élémentaires suivantes:

\snic{\cmatrix {\I_n + UV & 0\cr 0 & M\cr} \cmatrix {\I_n  & -NU\cr 0 & \I_m\cr} =
\cmatrix {\I_n + UV & -U\cr 0 & M\cr},}

\snic{\cmatrix {\I_n + UV & -U\cr 0 & M\cr} \cmatrix {\I_n & 0\cr V & \I_m\cr} =
\cmatrix {\I_n & -U\cr MV & M\cr},
}

%\sni
puis

\snic{
\cmatrix {\I_n & -U\cr MV & M\cr} \cmatrix {\I_n & U\cr 0 & \I_m\cr} =
\cmatrix {\I_n & 0\cr MV & MVU + M\cr} = \cmatrix {\I_n & 0\cr MV & \I_m\cr},
}

%\sni
et enfin

\snic{
\cmatrix {\I_n & 0\cr MV & \I_m\cr} \cmatrix {\I_n & 0\cr -MV & \I_m\cr} =
\cmatrix {\I_n & 0\cr 0 & \I_m\cr}.
}

%\sni
On a donc explicité des matrices $\alpha$, $\beta$, $\gamma$, $\delta \in
\EE_{n+m}(\gA)$ telles que

\snic{\cmatrix {\I_n + UV & 0\cr 0 & (\I_m + VU)^{-1}\cr}
\,\alpha\,\beta\,\gamma\,\delta = \I_{n+m},}

%\sni
d'où

\snic{
\cmatrix {\I_n + UV & 0\cr 0 & (\I_m + VU)^{-1}\cr} = 
\delta^{-1}\, \gamma^{-1}\, \beta^{-1}\, \alpha^{-1} =}

\snic{
\cmatrix {\I_n & 0\cr MV & \I_m\cr} 
\cmatrix {\I_n & -U\cr 0 & \I_m\cr} 
\cmatrix {\I_n & 0 \cr -V & \I_m\cr} 
\cmatrix {\I_n & NU\cr 0 & \I_m\cr} .
}

%\sni
Dans le cas particulier où $VU = 0$, on a montré que 

\snic{
\cmatrix {\I_n + UV
& 0\cr 0 & \I_m \cr} \in \EE_{n+m}(\gA).}

}

%%%%%%%%%%%%%%%%%%%%%%%%%%%%%%%%%%%%%%%%%
\exer{exoMennicke} 
En utilisant $ad' = 1+bc'$, $ad = 1+bc$, on obtient pour $A'^{-1}A$

\snuc{
\cmatrix {d'&-b\cr -c'&a\cr} \cmatrix {a&b\cr c&d\cr} =
\cmatrix {ad'-bc & bd'-bd\cr ac-ac'&ad-bc'\cr} =
\cmatrix {1 + b(c'-c) & b(d'-d)\cr a(c-c') &1+b(c-c')\cr}}

%\sni
En remplaçant $b(c'-c)$ par $a(d'-d)$, on voit que 

\snuc{A'^{-1}A = \I_2 + uv$
avec $u=\cmatrix {d'-d\cr c-c'\cr}$, $v=\cmatrix {a & b}$,  $vu = 0$, 
et $v  \hbox{ \umdz}.
}

%%%%%%%%%%%%%%%%%%%%%%%%%%%%%%%%%%%%%%%%%
\exer{exoMennicke2}{
On a $\meck {au}{b} = \meck {a}{b}\meck {u}{b}$.\\
Mais $au = 1-bv$ donc $\meck
{au}{b} = \meck {1-bv}{b} = \meck {1}{b} = 1$. Bilan: $\meck {a}{b}\meck
{u}{b} = 1$. \\
De la même manière, $\meck {u}{b}\meck {u}{v} = 1$, donc
$\meck {a}{b} = \meck {u}{v}$. \\
Enfin, $(a-v)u +
(b+u)v = 1$, donc $\meck {a-v} {b+u} = \meck {u} {v}$.
}

%%%%%%%%%%%%%%%%%%%%%%%%%%%%%%%%%%%%%%%%%
\exer{exoStabFreeCoRang1}{
 \emph{1.} 
Soit $y \in \Ae n$ tel que $\tra {y} x = 1$.\\
On a $\Ae n = \gA x \oplus \Ker
\tra {y}$ et donc $\Ae n/\gA x \simeq \Ker \tra {y}$ est \stlz.  \\
{Si
$x \sims{\GL_n(\gA)} x'$}, il est clair que $\Ae n/\gA x \simeq
\Ae n/\gA x'$. \\
 Réciproquement soit $\varphi : M = \Ae n/\gA x \to M' =
\Ae n/\gA x'$ un \isoz.  \linebreak 
On a $\Ae n \simeq M \oplus \gA x 
\simeq M' \oplus \gA x'$. 
On définit $\psi : \gA x \to \gA x',\; 
ax \mapsto ax'$. \linebreak 
Alors $\varphi \oplus \psi$ vu dans
$\GLn(\gA)$ transforme $x$ en $x'$, donc $x \sims{\GL_n(\gA)} x'$.
\\
 Un \vmd $x \in \Ae n$ fournit un module libre $\Ae n / \gA x$
\ssi $x$ fait partie d'une base de $\Ae n$.

 \emph{2.}
Posons $M = x^\perp$, $M' = x'^\perp$ et supposons $M \simeq M'$.  En
désignant par $\mathring M \subseteq (\Ae n)\sta$ l'\ort de $M \subseteq
\Ae n$, on a $\mathring M = \gA\tra {x}$ et $\mathring M' = \gA\tra {x'}$. Si
$\scp {x}{y} = 1$, $\scp {x'}{y'} = 1$, on a $\Ae n = \gA y \oplus M =
\gA y' \oplus M'$, d'où un \auto de $\Ae n$
transformant $M$ en $M'$ (envoyer $y$ sur $y'$), puis par
dualité, un \auto $u$ de $(\Ae n)\sta \simeq \Ae n$ transformant
$\gA\tra{x}$ en $\gA\tra{x'}$. On en déduit $u(\tra{x}) = \varepsilon\tra{x'}$
avec $\varepsilon \in \Ati$. Alors, $\varepsilon^{-1}\tra{u}$
transforme $x$ en $x'$.

 \emph{3.}
Si $G = E \oplus F$, alors $G\sta \simeq E\sta \oplus F\sta$; avec $G = \Ae {r+t} \simeq G\sta$, $F = \Ae {r} \simeq F\sta$, on obtient
le résultat.
 L'involution induite sur $\Um_n(\gA)/\GLn(\gA)$ est la suivante: à la classe
modulo $\GLn(\gA)$ de $x \in \Um_n(\gA)$, on associe la classe modulo
$\GLn(\gA)$ d'un \elt $y \in \Um_n(\gA)$ qui satisfait $\scp {x}{y} = 1$. Naturellement,
il y a plusieurs~$y$ qui conviennent mais leur classe modulo $\GLn(\gA)$ est
bien définie.

 \emph{4.}
On a $\Ae n = \gA y \oplus x^\perp = \gA y \oplus x'^\perp$ d'où $x^\perp
\simeq x'^\perp \simeq \Ae n / \gA y$ donc $x \sims{\GLn(\gA)} 
x'$.  Pour déterminer $g \in \GLn(\gA)$ réalisant $gx = x'$, on utilise
$\Ae n = \gA x \oplus y^\perp = \gA x' \oplus y^\perp$. De manière
\gnlez, soit $G = E \oplus F = E' \oplus F$; pour expliciter un
\auto de $G$ qui envoie $E$ sur $E'$, on procède comme suit. Soit
$\pi$ la \prn sur $E$, $\pi'$ celle sur $E'$ et $p = \I_G - \pi$, $p' =
\I_G - \pi'$. \\
Les \prrs $p$ et $p'$ ont même image $F$. Notons $h
= p'-p = \pi-\pi'$. \\
On obtient $h^2 = 0$ et $(\I_G - h) p (\I_G + h) = p'$, ou encore
$(\I_G - h) \pi (\I_G + h) = \pi'$. \hbox{Donc $\I_G - h$} est un \auto de
$G$ transformant $\Im\pi = E$ en $\Im\pi' = E'$. Ici $E = \gA x$, $E' = \gA
x'$, $F = y^\perp$, donc
$$
\pi(z) = \scp {z}{y}x, \quad \pi'(z) = \scp {z}{y}x', \quad
h(z) = \scp {z}{y} (x - x').
$$
L'\auto cherché de $\Ae n$ qui transforme $x$ en $x'$ 
est donc
$$
\I_n - h : z \mapsto z + \scp {z}{y} (x' - x) 
\quad \hbox {i.e.} \quad
\I_n - h = \I_n + uv 
$$
avec $u = x'-x \in \Ae {n \times 1}$, $v = \tra {y} \in \Ae {1 \times n}$;
on a bien $vu = 0$ et $v$ \umdz.

}

%%%%%%%%%%%%%%%%%%%%%%%%%%%%%%%%%%%%%%%%%
\exer{exoStabFree2}{
\emph{1.}
La clef du \pb se trouve dans la double \egt suivante pour un $u$ dans~$\gA$:
$z_1 = u(a+x_2)$, $z_2 = u(b-x_1)$, ce qui implique 

\snic{
z_1x_1 + z_2x_2 = u (ax_1 + bx_2) = z_1 b + z_2 (-a).} 

%\sni
Soit $u$ tel que $u(ax_1 + bx_2) + z_3x_3 + \cdots + z_nx_n = 1$ et $z =
(z_1, z_2, z_3, \ldots, z_n)$. On a alors $\scp {z}{x} = \scp
{z}{x'} = 1$. D'après l'exercice \ref{exoStabFreeCoRang1}, $x \sims G  x'$.  
Comme $(b, -a) \sims{\EE_2(\gA)}  (a,b)$, on a $x
\sims G (a,b, x_3, \ldots, x_n)$.

 \emph{2.}
Comme $x_1y_1 + x_2y_2 + x_3y_3 + x_4y_4 = 1$, on a

\snic{(x_1,x_2,x_3,x_4)\, \sims G \,
(y_1,y_2,x_3,x_4)\, \sims G  \,
(y_1,y_2,y_3,y_4)} 

%\sni
Le reste de la question en découle aussitôt.

 \emph{3.}
Méthode analogue à la question précédente.
}

%:   ---- Section*{references}-----------    
\Biblio
La section \ref{secMennicke} et la \dem du \thref{thStabSuslCorps} 
suivent de très près l'exposé de \cite{GM}.
Pour l'essentiel nous avons seulement transformé quelques arguments \lgbs
abstraits en arguments concrets via une utilisation de la machinerie \lgbe à \ideps expliquée dans 
la section \ref{secMachLoGlo}.\imlb

La section \ref{secCompVmdsPols} est directement inspirée de 
\cite[chapitre VI, section 2]{Lam06}.

\newpage \thispagestyle{CMcadreseul}
\incrementeexosetprob

\appendix
\let\showchapter\oldshowchapter
\let\showsection\oldshowsection
%%%%%%%%%%%%%%%%%%%%%%%%%%%%%%%%%%%%%%%%%%%%%%%%%%%%%%%%%%%%%%%%%%%%%%%%%%
%%%%%%%%% A-annexe1PTF

%%%%%%%%%%%%%% Principes d'omniscience %%%%%%%%%%
\def\thechapter{}
%\def\chaptername{\appendixname}
%\def\partname{Annexes}

%\part*{Annexes}

\def\chaptername{Annexe}
\refstepcounter{chapter}

\chapter*{Annexe\\[.6em] Logique constructive}
\addcontentsline{toc}{chapter}{\hskip-0.8emAnnexe. Logique constructive}
\mtcaddchapter
\markboth{Annexe}{Logique constructive}
\label{chapPOM}

\minitoc

\newpage	
\Intro

Cette annexe est consacrée à l'exposition de quelques concepts de base
des \coma dans le style de Bishop, illustré par les trois ouvrages 
fondateurs \cite{B67,BB85,MRR}.

Par logique constructive, nous entendons la logique des \comaz.

\setcounter{section}{0}

%--- Section Objets de base Ensembles Fonctions
\section{Objets de base, Ensembles, Fonctions}
%------------------------------------------------------------------
Entiers naturels et constructions sont deux notions primitives.
Elles ne peuvent pas être définies.

D'autres notions primitives sont liées au langage usuel et difficiles à
situer précisément. Par exemple l'\egt du nombre $2$ en deux
occurrences distinctes.

La formalisation d'un morceau de \maths peut être utilisée pour mieux
comprendre ce que l'on est en train d'y faire. Mais pour parler à propos d'un
formalisme il faut comprendre beaucoup de choses qui sont du même genre de
complexité que les entiers naturels.
Ainsi, le formalisme est seulement un outil et il ne peut pas remplacer les
intuitions et les expériences de base (par exemple les entiers naturels,
les constructions): si puissant que soit un ordinateur, il ne comprendra
jamais \gui{ce qu'il fait}, ou encore, comme le disait René Thom, \gui{Tout
ce qui est rigoureux est insignifiant}.

%---- subsubsec Ensembles -----------------
\subsubsec{Ensembles}
%------------------------------------------------------------------

\noi Un \emph{ensemble} $(X,=_X,\neq_X)$ est défini en disant:

\noi --- comment on peut construire un \elt de l'ensemble (nous
disons que nous avons défini un  \emph{préensemble} $X$)
\index{preensemble@préensemble}

\noi --- quelle est la signification de l'\emph{\egtz} pour deux
\elts de l'ensemble (nous avons à montrer que c'est bien une
relation d'\eqvcz)

\noi --- quelle est la signification de la
\ix{distinction}{\footnote{Cette terminologie \emph{n'est pas} un
hommage à Pierre Bourdieu. Tous comptes faits, nous préférons
\emph{distinction} à \emph{non-\egtz}, qui présente
l'inconvénient d'une connotation négative, et à \emph{in\egtz}
qui est plutôt utilisé dans le cadre des relations d'ordre. Pour les
nombres réels par exemple, c'est l'\egt et non la distinction qui
est une assertion négative.}} pour deux \elts de l'ensemble (on
dit alors que les \elts sont \emph{discernables} ou
\emph{distincts}). Nous avons à montrer les \prts suivantes:

\hspace*{5mm}-- $ \; (x\neq_X y \; \land \; x=_Xx' \; \land \; y=_Xy')\; \Rightarrow \;
x' \neq_X y'$,

\hspace*{5mm}-- $ \; x\neq_X x$  est impossible,

\hspace*{5mm}-- $ \; x\neq_X y\; \Rightarrow \; y \neq_X x$.

\ss Ordinairement, on laisse tomber l'indice $X$ pour les symboles $=$ et
$\neq $. Si la distinction n'est pas précisée, elle est implicitement
définie comme signifiant l'absurdité de l'\egtz.

\ss Une relation de distinction est appelée une relation de \emph{séparation} si
elle vérifie la propriété de \ix{cotransitivité}  suivante (pour
trois \elts $x,y,z$ de $X$ arbitraires):\index{separation@séparation}

\hspace*{5mm}-- $\; \; x\neq_X y\; \Rightarrow \; (x \neq_X z\; \lor\;  y\neq_X z)$

\noi Une relation de séparation $\neq_X$ est dite \emph{étroite} si
$x=_X y$ équivaut à l'absurdité de  $x\neq_X y$. Dans un ensemble avec
une séparation étroite, la distinction est souvent plus importante que
l'\egtz. \index{separation@séparation!étroite}

\ss Un ensemble  $(X,=_X,\neq_X)$ est dit \emph{discret} \label{discret} si
l'on a  $$ \forall x,y\in X\; (x=_X y \lor x\neq _X y).$$  Dans ce cas la
distinction est une séparation étroite et elle équivaut à
l'absurdité de l'\egtz.
\index{ensemble!discret}

%---- subsubsec  -----------------
\subsubsec{Les entiers naturels}

\noi  L'ensemble $\NN=\so{0,1,2,\ldots}$ des entiers naturels
est considéré comme bien défini a priori. Notez
cependant que \cot il s'agit d'un \ix{infini potentiel} et pas d'un
\ix{infini actuel}. On entend par l'idée d'infini potentiel que
l'infinitude de $\NN$ est appréhendée comme une notion essentiellement
négative: on n'a jamais fini d'épuiser les entiers naturels. Au
contraire, la sémantique de $\NN$ en \clama est celle d'un infini achevé,
qui existe \gui{quelque part} au moins de manière purement idéale.

Un entier naturel peut être codé d'une manière usuelle. La comparaison
de deux entiers donnés sous forme codée peut être faite de manière
s\^ure. En bref, l'ensemble des entiers naturels est un ensemble discret et
la relation d'ordre est \emph{décidable}:
%-----------------begin $$----------------
$$\preskip.2em \postskip.1em
\forall n,m\in\NN\; \; \; (n<m\; \; \lor\; \;  n=m\; \; \lor\; \;  n>m)$$
%-----------------end $$------------------

%---- subsubsec  -----------------
\subsubsec{Ensembles de couples}

\noi Quand deux ensembles sont définis, leur \emph{produit cartésien} est
également défini, de manière naturelle: la fabrication des couples
d'objets est une construction \elrz.  L'\egt et la distinction
sur un produit cartésien sont définis de manière naturelle.

%---- subsubsec Fonctions -----------------
\subsubsec{Fonctions}
%------------------------------------------------------------------

\noi  L'ensemble $\NN^\NN$ des suites d'entiers naturels dépend de la notion
primitive de construction. Un \elt de $\NN^\NN$ est une construction
qui prend en entrée  un \elt de $\NN$ et donne en sortie un
\elt de $\NN$.  L'\egt de deux \elts dans $\NN^\NN$ est
l'\emph{\egt extensionnelle}:
$$\preskip.4em \postskip.4em
(u_n)=_{\NN^\NN}(v_n) \quad \mathrm{ signifie} \quad  \forall n\in\NN\;\;
u_n=v_n.
$$
Ainsi, l'\egt entre deux \elts de $\NN^\NN$ demande a priori une
infinité de \gui{calculs \elrsz}, en fait l'\egt réclame une
preuve.\index{extensionnelle!égalité ---}

\noi La distinction de deux \elts de $\NN^\NN$ est la relation de
\emph{distinction extensionnelle}:
$$\preskip-.2em \postskip.4em
(u_n)\neq_{\NN^\NN}(v_n) \quad {\eqdef} \quad  \exists n\in\NN\;\;  u_n\neq
v_n.
$$
Ainsi, la distinction de deux \elts de $\NN^\NN$ peut être
constatée par un simple calcul.\index{extensionnelle!distinction ---}

%:--- example{exo neq dans N^N}--------
\begin{example}
\label{exo neq dans N^N}\relax
\emph{
La distinction de $\NN^\NN$ est une relation de séparation
étroite. }
\end{example}
%--- end-example-----------------------------------------

 L'argument diagonal de Cantor est \cofz. Il montre  que  $\NN^\NN$ est
\emph{beaucoup plus compliqué} que $\NN$. D'un point de vue \cofz, $\NN$ et
$\NN^\NN$ sont seulement des  infinis potentiels: cela n'a pas de
signification de dire qu'un infini potentiel est \emph{plus grand} qu'un
autre.

\mni \emph{Digression.}
Quand vous dites \gui{Je vous donne une suite d'entiers naturels}, vous
devez prouver que la construction $\; n\mapsto u_n\; $ que vous proposez
fonctionne pour n'importe quelle entrée $n$. Par ailleurs, quand vous dites
\gui{Considérons une suite arbitraire de nombres naturels $(u_n)_{n\in\NN}$},
la seule chose
que vous savez avec certitude est que pour tout $n\in\NN$, vous avez
$u_n\in\NN$, et que cet $u_n$ est non ambigu: vous pouvez par exemple
concevoir la suite comme donnée par un oracle.
En fait, vous pourriez a priori demander, de
manière symétrique, ce qu'est exactement la construction $\; n\mapsto u_n$,
et une preuve que cette construction
fonctionne pour toute entrée $n$.

Mais, dans le constructivisme à la Bishop, on ne fait aucune hypothèse
précise concernant \gui{ce que sont les constructions légitimes de $\NN$
vers $\NN$}, ni non plus sur \gui{qu'est-ce \prmt qu'une preuve
qu'une construction marche?}. Ainsi nous sommes dans une situation
dissymétrique.

Cette dissymétrie a la conséquence suivante. Tout ce que vous prouvez a
un contenu  calculatoire. Mais tout ce que vous prouvez est également
valide  d'un point de vue classique. Les \clama
pourraient voir les \coma comme parlant seulement d'objets \cofsz. Et les
\coma  de Bishop sont certainement intéressées au premier chef par les
objets \cofs (cf. \cite{bi}). Mais en fait, les \coma à la Bishop font des
preuves constructives qui marchent pour n'importe quel type d'objets
\mathsz\footnote{\ldots~s'il existe des objets
\maths non \cofsz.}. Les \thos que l'on trouve dans \cite{BB85} et \cite{MRR}
sont valables en \clamaz, mais ils supportent aussi l'interprétation
constructive russe (dans laquelle tous les objets \maths sont des mots d'un
langage formel que l'on pourrait fixer une fois pour toutes) ou encore la
philosophie intuitionniste de Brouwer, qui a une composante nettement
idéaliste.    \eoe

\ss Après cette digression revenons à nos moutons: les fonctions. De
manière \gnlez,  une \emph{fonction} $f:X\rightarrow Y$ est une
construction qui prend en entrée un $x\in X$ et une preuve que $x\in X$, et
donne en sortie un $y\in Y$ et une preuve que $y\in Y$. En outre, cette
construction doit être \emph{extensionnelle}:
$$
x=_Xx'\Rightarrow f(x)=_Yf(x')\qquad  \mathrm{et} \qquad  f(x)\neq_Yf(x')
\Rightarrow x\neq_Xx'
$$
Quand $X$ et $Y$ sont des ensembles bien définis, on considère (dans les
\coma à la Bishop) que l'ensemble $\cF(X,Y)$ des fonctions~\hbox{$f:X\rightarrow Y$} est aussi bien défini. Pour l'\egt et la distinction on prend
les définitions extensionnelles usuelles.

\ss Une fonction $f:X\rightarrow Y$ est \emph{injective} si elle vérifie
$$f(x)=_Yf(x')\Rightarrow x=_Xx'\qquad \mathrm{et} \qquad x\neq_Xx'
\Rightarrow f(x)\neq_Yf(x') $$

%---- subsubsec Ensembles finis, enum, deno
\subsubsec{Ensembles finis, bornés, \enums et \denbs}
%------------------------------------------------------------------

\noi Nous donnons maintenant un certain nombre de définitions \cot
pertinentes en relation avec les concepts d'ensembles finis, infinis et
\denbs en \clamaz.\rdb
%-----------------begin item------------------
\begin{itemize}
\item  Un ensemble est \emph{fini} s'il y a une bijection entre cet ensemble
et l'ensemble des entiers $<n$ pour un certain entier $n$.%
\footnote{C'est la \dfn donnée \paref{Deux mots}, dans le paragraphe \gui{Deux mots sur les ensembles finis}.}
\item  Un ensemble $X$ \emph{finiment \enumz} si l'on a une application
surjective d'un segment $[0,n[\,$ de $\NN$ sur $X$ (c'est la \dfn donnée \paref{Deux mots}).%
\index{ensemble!finiment énumérable}\index{finiment énumérable!ensemble ---}%
\item  Un ensemble $X$ est dit \emph{\enumz} si l'on a donné un
moyen de l'énumé\-rer en lui laissant la possibilité d'être vide\footnote{Page \pageref{Deux mots}, on a donné la \dfn pour les ensembles non vides.}, ce qui se passe en pratique comme suit. 
On donne un $\alpha\in\so{0,1}^\NN$ et une opération $\varphi$ qui satisfont les deux assertions suivantes:
\begin{itemize}
\item si $\alpha(n)=1$ alors $\varphi$ construit à partir de l'entrée $n$
un \elt de $X$,
\item  tout \elt de $X$ est construit de cette façon.
\end{itemize}
\index{ensemble!enumerable@énumérable}\index{enumerable@énumérable!ensemble ---}%
\item  Un ensemble est dit \emph{\denbz} s'il est \enum  et discret.
\index{ensemble!dénombrable}\index{denombrable@dénombrable!ensemble ---}
\item  Si $n$ est un entier non nul, on dit qu'un ensemble \emph{possède
au plus $n$ \eltsz} si pour toute famille $(a_i)_{i=0,\ldots,n}$ dans
l'ensemble il existe des entiers $h$ et $k$ ($0\leq h<k\leq n$) tels que
$a_h=a_k$.
\item  Un ensemble $X$ est \emph{borné en nombre} (\emph{borné} tout
court s'il n'y a pas d'ambiguïté) s'il existe un entier $n$ non nul tel
que $X$ ait au plus $n$ \elts (\dfn donnée~\paref{ensborn}).
\index{borné!ensemble ---}
\index{ensemble!borné}
\item  Un ensemble $X$  est \emph{faiblement fini} si pour toute suite
$(u_n)_{n\in\NN}$ dans $X$ il existe $m$ et $p>m$ tels que $u_m=u_p$.
\index{ensemble!faiblement fini}\index{faiblement fini!ensemble ---}
\item  Un ensemble $X$ est \emph{infini} s'il existe une application
injective $\NN\rightarrow X$.%
\index{ensemble!infini}\index{infini!ensemble ---}
\end{itemize}
%-----------------end item------------------

%--- example{exo inf deno1}-----------
\begin{example}
\label{exo inf deno1}\relax
{\rm Un ensemble infini et \denb peut être mis en bijection
avec~$\NN$.
}
\end{example}
%--- end-example-----------------------------------------

\rdb
%:---- subsubsec Parties d'un ensemble -----
\subsubsec{Parties d'un ensemble}\label{P(X)}\relax
%------------------------------------------------------------------

\noi Une partie d'un ensemble $(X,=_X,\neq_X)$ est définie par une
propriété $P(x)$ \emph{portant sur les \elts de $X$}, c.-à-d.
vérifiant 
$$\forall x,y\in X\;\big( \; (\; x=y\; \land \; P(x)\;)\; \;
\Longrightarrow\; \;  P(y)\; \big).
$$
Un \elt de la partie $\sotq{x\in X}{P(x)}$ est donné par un
couple $(x,p)$ où~$x$ est un \elt de $X$ et $p$ est une preuve que
$P(x)$({\footnote{Par exemple, un nombre réel $\geq 0$ est \emph{un peu
plus qu}'un nombre réel.}}). Deux \prts concernant les
\elts de $X$ définissent la même partie lorsqu'elles sont
\eqvesz.

On peut aussi présenter les choses de la manière suivante, qui, bien que
revenant au même, fait un peu moins mal à la tête au nouveau venu.
Une partie de $X$ est donnée par un couple $(Y,\varphi)$  où  $Y$ est un
ensemble et $\varphi$ est une fonction injective de $Y$ dans
$X$({\footnote{Par exemple on peut définir les nombres réels $\geq 0$
comme ceux qui sont donnés par des suites de Cauchy de rationnels $\geq
0$.}}).
Deux couples $(Y,\varphi)$ et $(Y',\varphi')$ définissent la même partie
de $X$ si l'on a
$$
\forall y\in Y\; \exists y'\in Y'\; \varphi(y)=\varphi'(y')\;\;
\hbox{ et }\;\;
\forall y'\in Y'\; \exists y\in Y\; \varphi(y)=\varphi'(y').
$$

En \coma on considère que les
parties de $X$ ne forment pas un ensemble,
mais une \emph{classe}.\index{classe!(versus ensemble)} 
Cette classe n'a pas clairement le statut d'un ensemble (au sens
donné plus haut).
L'intuition est la suivante: les ensembles sont des classes suffisamment
bien définies pour que l'on puisse quantifier universellement ou existentiellement
sur leurs \eltsz. Pour cela, il faut que le procédé de construction
des \elts soit clair.

\ss Rappelons qu'une partie $Y$ de $X$ est dite \emph{détachable} lorsque l'on a un test
pour \gui{$x\in Y~?$} lorsque $x\in X$. Les parties détachables de $X$
forment un ensemble qui s'identifie à $\{0,1 \}^X$.

Constructivement, on ne connaît aucune partie détachable de $\RR$,
hormis $\emptyset$ et $\RR$: \emph{il n'y a pas de trou dans le continu sans
la logique du tiers exclu.}

\mni\rem  Une variante \cot intéressante  pour \gui{une partie $Y_1$ de $X$} est
obtenue en considérant un couple $(Y_1,Y_2)$ de parties de $X$ qui
vérifient les deux propriétés suivantes
$$\forall x_1\in Y_1\; \forall  x_2\in Y_2\;  \; x_1\neq_X x_2\qquad
\mathrm{et}\qquad
\forall x\in X\; \lnot( x\notin Y_1\; \land\;  x\notin Y_2)
$$
Le \emph{\copz} est alors donné par le couple  $(Y_2,Y_1)$, ce
qui rétablit une certaine symétrie.\eoe

%:---- La classe des parties d'un ensemble ----
\mni\emph{La classe des parties d'un ensemble}~
\rdb
\\
Notons  $\rP(X)$ la classe des parties de l'ensemble $X$.
Si l'on admettait  $\rP(\{0\})$ comme un ensemble, alors 
 $\rP(X)$ serait \egmt un ensemble et il y aurait
une bijection naturelle entre  $\rP(X)$  et
$\cF\big(X,\rP(\{0\})\big)=\rP(\{0\})^X$.

Ceci montre que toute la difficulté avec
l'ensemble des parties est concentrée sur la classe $\rP(\{0\})$, \cad
 la classe des \emph{valeurs de vérité}. En \clamaz, on admet que cette classe est un ensemble à deux \eltsz, c'est le
\emph{principe du tiers exclu} \TEMz:
$$
\rP(\{0\})=\{\{0\},\emptyset  \}
$$
(la classe des valeurs de vérité se réduit à l'ensemble
 $\{\Vrai,\Faux \}$) et l'on n'a évidemment plus aucun \pb avec
$\rP(X)$.
%En fait il s'agit plutôt de l'art de l'esquive: cachez ce
%\pb que je ne saurais voir.

%\ss Il ne semble pas que l'on connaisse un seul \tho \mathe pertinent de \clama
%dont l'étude du point de vue \cof nécessite le recours à  l'ensemble
%$\rP(\{0\})$.

%--- Section Affirmer signifie prouver ---------
\section{Affirmer signifie prouver} \label{secAffirmerProuver}
%------------------------------------------------------------------
En  \coma la vérité est aussi le résultat d'une construction. Si $P$
est une assertion \mathz, nous écrirons \gui{$\; \vda P\; $} pour
\gui{nous avons une preuve de $P$}.

\ss Les assertions \elrs peuvent être testées par des calculs simples.
Par exemple, la comparaison de deux entiers naturels.
Quand une assertion signifie une infinité d'assertions \elrs (e.g., la
conjecture de Goldbach\footnote{Tout nombre pair $\geq4$ est somme de deux nombres premiers.}), les \coma considèrent qu'elle n'est pas a priori
\gui{vraie ou fausse}. A fortiori, les assertions ayant une complexité
logique encore plus grande ne sont pas considérées (d'un point de vue
\cofz) comme ayant a priori la valeur de vérité $\Vrai$ ou $\Faux$.

Ceci ne doit pas être nécessairement considéré comme une position
philosophique concernant la vérité. Mais c'est s\^urement une position
\mathe concernant les assertions \mathsz. En fait, cette
position est nécessaire pour avoir une signification calculatoire pour
tous les \thos qui sont prouvés de manière \covz.

\mni \emph{Digression carrément philosophique.} Cette position est
également à distinguer de la position qui consiste à dire qu'il y a
certainement différents univers \maths possibles, par exemple
l'un dans lequel l'hypothèse du continu{\footnote{L'hypothèse du
continu est, dans la théorie des ensembles classiques, l'affirmation qu'il
n'y a pas de cardinal strictement compris entre celui de $\NN$ et celui de
$\RR$, autrement dit, que toute partie infinie de $\RR$ est équipotente à
$\NN$ ou à $\RR$.}}  est vraie, un autre dans lequel elle est fausse. Cette
position est naturellement parfaitement défendable (Cantor, et sans
doute G\"odel, l'auraient  refusée au nom d'un réalisme platonicien des
Idées), mais elle intéresse peu les \coma à la Bishop qui ont pour
objet d'étude une abstraction de l'univers concret des calculs finis, avec
l'idée que cette abstraction doit correspondre d'aussi près que possible
à la réalité qu'elle veut décrire.
Ainsi, l'hypothèse du continu est plutôt dans ce cadre considérée
comme vide de signification, car il est vain de vouloir comparer des
infinis potentiels selon leur taille. Si l'on désire les comparer selon leur
complexité, on
s'aperçoit bien vite qu'il n'y a aucun espoir de mettre une vraie
relation d'ordre total dans ce fouillis. En conséquence, l'hypothèse du
continu ne semble rien d'autre aujourd'hui qu'un jeu des spécialistes de
la théorie formelle $\ZF$.
Mais chacun et chacune est bien libre de croire Platon, ou même  Cantor,
ou Zermelo-Frankel, ou encore, pourquoi pas, de croire  en la multiplicité
des mondes. Personne ne pourra jamais lui prouver qu'il a tort.
Et rien ne dit par ailleurs que le jeu $\ZF$ ne s'avérera pas un jour
vraiment utile, par exemple pour comprendre certains points subtils des
\maths qui ont une signification concrète.  \eoe

%--- Section Connect et quantificateurs --------
\section{Connecteurs et quantificateurs} \label{secBHK}
%-----------------------------------------

Ici nous donnons l'explication \gui{Brouwer-Heyting-Kolmogorov}  pour la
signification constructive des symboles logiques usuels. Ce sont seulement
des explications informelles, pas des \dfnsz\footnote{Pour le point de vue de Kolmogorov plus \prmt sur \gui{la logique des \pbsz} on pourra
consulter  \cite[Kolmogorov]{Kol} et \cite[Coquand]{CoqK}.}.

Il s'agit d'explications \gui{détaillées}, 
pour ce qui concerne les connecteurs logiques et les \qtfsz,
concernant ce que l'on entend par le slogan \gui{affirmer
signifie prouver}. Quand on écrit $\vda P$ on sous-entend que l'on dispose d'une preuve \cov de $P$. Nous expliciterons ceci en donnant un nom, par exemple
$p$, à cet objet \mathe qu'est la preuve de $P$.
Les explications concernent alors ces objets particuliers $p$, mais tout ceci
reste informel.

\mni \textbf{Conjonction:} $\vda P \; \land\;  Q$  signifie: \gui{$\vda P$  et $\vda Q$}
(comme pour la logique classique). En d'autres termes: une preuve de $P \;
\land\;  Q$ est un couple $(p,q)$ où $p$ est une preuve de $P$ et $q$ une
preuve de~$Q$.

\mni \textbf{Disjonction:} $\vda P \; \lor\;  Q$  signifie: \gui{$\vda P$  ou $\vda Q$}
(ce qui ne marche pas avec la logique classique).
En d'autres termes: une preuve de
 $P \, \lor\,  Q$ est un couple $(n,r)$ avec $n\in\{0,1 \}$.
Si $n=0$, $r$ doit être une preuve de $P$, et si~$n=1$, $r$ doit être une
preuve de~$Q$.

\mni \textbf{Implication:} $\vda P \; \Rightarrow \;  Q$  a la signification
suivante: \\
une preuve de $\; P \; \Rightarrow \;  Q\; $ est une
cons\-truction $p\mapsto q$ qui transforme toute preuve $p$ de $P$ en une
preuve $q$ de~$Q$.

\mni \textbf{Négation:} $\lnot P $  est une abréviation de $P \; \Rightarrow \;
0=_\NN1$.

\mni \textbf{Quantificateur universel:} (similaire à l'implication). \emph{Une
\qtn est toujours une \qtn sur les objets d'un ensemble
défini au préalable.}
Soit $P(x)$ une propriété concernant les objets $x$ d'un ensemble $X$.
\\
Alors   $\vda \forall x\in X\; \; P(x)$ a la signification suivante:
nous avons une construction $(x,q)\mapsto p(x,q)$ qui prend en entrée 
n'importe quel couple $(x,q)$, où 
$x$ est un objet et  $q$ est  une preuve que $x\in X$, et donne en sortie  une preuve~\hbox{$p(x,q)$} de l'assertion~$P(x)$.

 Pour une \qtn sur $\NN$, on estime que la donnée d'un entier
$x$ (sous-entendu, sous forme standard) suffit à prouver que $x\in\NN$: la
partie $q$ dans le couple $(x,q)$ ci-dessus peut être omise.

%--- example{exo log0}----------------
\begin{example}\label{exo log0}\relax
\emph{Supposons que les \prts $P$ et $Q$ dépendent d'une variable $x\in \NN$.
Alors une preuve de
$\,  \forall x\in \NN \;\big(P(x) \, \lor\,  Q(x)\big)\, $ est une cons\-truction
$\NN\ni x\mapsto \big(n(x),r(x)\big)$, où $n(x)\in\{0,1 \}$:
 si $n(x)=0$, $r(x)$
est une preuve de $P(x)$, et si $n(x)=1$, $r(x)$ est une preuve
de~$Q(x)$.
}\eoe
\end{example}
%--- end-example--------------------------------

\mni \textbf{Quantificateur existentiel:} (similaire à la disjonction) \emph{Une
\qtn est toujours une \qtn sur les objets d'un ensemble
défini au préalable.} Soit $P(x)$ une \prt concernant les objets
$x$ d'un ensemble~$X$. Alors    $\vda \exists x\in X\;  P(x)$
 a la signification
suivante: une preuve de $\exists x\in X\;  P(x)$ est un triplet $(x,p,q)$ où $x$ est un objet, 
$p$ est une preuve de $x\in X$,  et  $q$ une  preuve  de~$P(x)$.

%--- example{exo log1}----------------
\begin{example}
\label{exo log1}\relax
\emph{Soit $P(x,y)$ une propriété concernant les entiers naturels~$x$  et $y$.
Alors l'affirmation
$$\vda \forall x\in\NN\; \; \exists y\in\NN\;\;  P(x,y)$$
signifie: voici un couple $(u,p)$ où $u$ est une construction $u:x\mapsto y=u(x)$ de $\NN$ vers $\NN$ et
$p$ est une preuve de $\vda \forall x\in\NN\; P\big(x,u(x)\big)$.
}\eoe
\end{example}
%--- end-example--------------------------------

%--- example{exo log2}----------------
\begin{example}
\label{exo log2}\relax (Logique des propositions)\\
{\rm La classe des valeurs de vérité en \coma est
une \agHz.\\
NB: $\rP(\so{0})$ étant une classe et non un ensemble
on entend simplement par là que les connecteurs $\vi$, $\vu$ et $\im$
et les constantes $\Vrai$ et $\Faux$ satisfont les axiomes des \agHsz.

En particulier, soient $A,\; B,\; C$ des propriétés \mathsz. On a les
\eqvcs suivantes.
%-----------------begin item------------------
\begin{itemize}
\item [$\vda$]  $\big((A\Rightarrow C)\; \land\; (B\Rightarrow C)\big)\; \;
\Longleftrightarrow \; \; \big((A\; \lor \; B)\; \Rightarrow \; C\big) $
\item [$\vda$]  $\big(A\Rightarrow (B\Rightarrow C)\big)\; \; \Longleftrightarrow \; \;
\big((A\; \land \; B)\; \Rightarrow \; C\big)$
\item [$\vda$]  $\lnot (A\; \lor B\; ) \; \; \Longleftrightarrow \; \;
        (\lnot A \land \lnot B)$
\item [$\vda$]  $(A\Rightarrow B)\; \; \Longrightarrow \; \;
(\lnot B\; \Rightarrow  \; \lnot A)$
\item [$\vda$]  $\lnot\lnot\lnot A\; \; \Longleftrightarrow \; \;\lnot A$
\end{itemize}
%-----------------end item------------------
Si en outre on a $ \vda A\; \lor \lnot A$ et $\vda B\; \lor \lnot B$, alors
on~a:
%-----------------begin item------------------
\begin{itemize}
\item [$\vda$] $\lnot \lnot A\; \; \Longleftrightarrow \; \; A$
\item [$\vda$] $\lnot (A\; \land B\; ) \; \; \Longleftrightarrow \; \;
        (\lnot A \lor \lnot B)$
\item [$\vda$] $(A\Rightarrow B) \; \; \Longleftrightarrow \; \;(\lnot A \lor B)$
\eoe\end{itemize}
%-----------------end item------------------
}
\end{example}
%--- end-example--------------------------------

%--- Remark-
\rem
\label{rem propneg}\relax
Puisque  $\lnot\lnot\lnot A\,\Leftrightarrow \,\lnot A$, une propriété
$C$ est équiva\-len\-te à une propriété $\lnot B$ (pour une certaine
propriété $B$ non encore précisée) \ssi
$\lnot\lnot C\,\Rightarrow \,C$. Ainsi, on peut définir en \coma le
concept de \emph{propriété négative}. En \clamaz, le concept n'a pas
d'intérêt puisque toute propriété est négative. En \comaz,
il faut prendre garde que $\Vrai$ est aussi une propriété négative:
puisque $\Faux \Rightarrow \Faux$, $\lnot \Faux$ est vrai.  \eoe

%--- Section para meca -------------------------
\section{Calculs mécaniques}\label{AnnexeCalculsMec}

Nous discutons ici un point qui est souvent mal apprécié par
les mathématiciens classiques. Une fonction de $\NN$ vers $\NN$ est donnée
par une construction. Les constructions usuelles correspondent à des
programmes algorithmiques  qui peuvent tourner sur un ordinateur
\gui{idéal}{\footnote{Un ordinateur disposant de tout l'espace et de tout
le temps nécessaire au calcul envisagé.}}. Ceci conduit à la notion de
\emph{calculs mécaniques}. Une fonction  $f\in\NN^\NN$ obtenue par un tel
calcul mécanique est appelée une \emph{fonction récursive}.
\\
Le sous-ensemble $\Rec\subset \NN^\NN$ formé par les fonctions récursives
peut alors être décrit de manière plus formelle comme nous allons l'expliquer
maintenant.

Rappelons qu'une \emph{fonction primitive récursive} est une fonction
$\NN^k\rightarrow \NN$ qui peut être définie par composition ou par
récurrence simple à partir de fonctions primitives récursives déjà
définies (nous commençons avec les fonctions constantes et l'addition
$+$).
Appelons $\Prim_2$ l'ensemble des fonctions primitives récursives
$\NN^2\rightarrow \NN$. On vérifie sans peine que $\Prim_2$ est un ensemble
énumérable.
\\
Une fonction $\beta\in\Prim_2$ peut être pensée comme simulant l'exécution
d'un programme de la manière suivante. Pour une entrée $n$ nous calculons~\hbox{$\beta(n,m)$} pour
$m=0$, $1$, $\ldots$ jusqu'à ce que   $\beta(n,m)\neq 0$ (intuitivement: jusqu'à ce que le programme arrive à l'instruction {\sf Stop}). Alors, la fonction $\alpha\in\Rec$ calculée par le
\gui{programme} $\beta\in\Prim_2$ est: $f:n\mapsto \beta(n,m_n)-1$ où
$m_n$ est la première valeur  de $m$ telle que $\beta(n,m)\neq 0$.

Ainsi, nous obtenons une application surjective d'un sous-ensemble $Rec$ de
$\Prim_2$ sur $\Rec$, et   $\Rec$ peut être identifié au préensemble
$Rec$ muni de l'\egt et de la distinction convenables. Cela signifie
que  $\Rec$ est défini comme un \gui{quotient}({\footnote{Puisque  $\Rec$
est l'image de  $Rec$ par une application surjective.}}) d'un sous-ensemble
d'un ensemble énumérable. Les \elts de la partie $Rec$ de
$\Prim_2$ sont définis par la condition suivante:
%-----------------begin $$----------------
$$\beta\in Rec \; \eqdef \;(*)\;:\;  \forall n\in \NN \;\; \exists m\in
\NN\;\;\; \beta(n,m)\neq 0 
$$
%-----------------end $$------------------
D'un point de vue classique, pour n'importe quel $\beta\in\Prim_2$,
l'assertion~$(*)$ ci-dessus est vraie ou fausse dans l'absolu, en
référence à la logique du tiers exclu (ou, si l'on préfère, à
l'infinité actuelle de $\NN$): la notion de calcul mécanique peut ainsi
être définie sans référence aucune à une notion primitive de
construction.

D'un point de vue \cof par contre, l'assertion $(*)$  doit être prouvée,
et une telle preuve est elle-même une construction. Ainsi \emph{la notion
de calcul mécanique dépend de la notion de construction, qui ne peut pas
être définie}.

Signalons pour terminer ce paragraphe que le constructivisme russe à la
Markov admet comme principe fondamental l'\egt  $\Rec=\NN^\NN$,
principe parfois appelé \textbf{Fausse Thèse de Church}\index{These@Thèse de Church!Fausse ---}. Voir
\cite{Be,BR} et \cite[Richman]{ri2}. La vraie  \textbf{Thèse de Church}\index{These@Thèse de Church} est qu'aucun
système de calcul  automatique ne pourra jamais calculer d'autres
fonctions que les fonctions récursives: on pourra améliorer les
performances des ordinateurs, mais aucun système de calcul automatique ne
pourra dépasser ce qu'ils savent calculer \gui{en principe} (\cad s'ils disposent du temps et de l'espace nécessaire).
La vraie  Thèse de Church est extrêmement vraisemblable, mais
elle n'est évidemment susceptible d'aucune preuve.

%--- Section Principes d'omniscience -----------
\penalty-2500
\section{Principes d'omniscience}
%--------------------------------------------------
On appelle \emph{principe d'omniscience} un principe qui, bien que vrai en
\clamaz, pose manifestement \pb en \comaz, car il suppose une
connaissance a priori de ce qui se passe avec un infini potentiel. Le mot
omniscience vaut donc ici pour \gui{prescience de l'infini potentiel}. Les
principes d'omniscience ont en \gnl des contre-exemples durs dans les
\coma russes. Ils ne peuvent cependant pas être démontrés faux dans
les \coma à la Bishop, car elles sont compatibles avec les \clamaz.
%---- subsubsec LPO -----------------------
\subsubsec{Le Petit Principe d'Omniscience}

\noi Soit $\alpha=(\alpha_n)\in\{0,1\}^\NN$ une {\it suite binaire}, i.e.,
une construction qui donne pour chaque entier naturel (en entrée) un
\elt de
$\{0,1\}$ (en sortie). Considérons les assertions suivantes:
$$\arraycolsep2pt
\begin{array}{rcl}
P(\alpha) & :~ & \alpha_n=1 \mathrm{\;  pour \; un\;  } n,\\[1mm]
\lnot P(\alpha) & : & \alpha_n=0 \mathrm{\; pour\;  tout\;  }n,\\[1mm]
P(\alpha)\vee \lnot P(\alpha) & : & P(\alpha) \mbox{ ou }\lnot P(\alpha),\\[1mm]
\forall \alpha \;\; \big(P(\alpha)\vee \lnot P(\alpha)\big) & : & \mathrm{pour\;
toute \;suite
\;binaire\; } \alpha, \; P(\alpha)\mathrm{\;ou\; }\lnot P(\alpha).
\end{array}
$$

Une preuve \cov de $P(\alpha)\vee \lnot P(\alpha)$ devrait fournir un
algorithme qui ou bien montre que $\alpha_n=0$ pour tout $n$, ou bien
calcule un entier naturel $n$ tel que $\alpha_n=1$.

Un tel algorithme est beaucoup trop performant, car il permettrait de
résoudre de manière automatique un grand nombre de conjectures
importantes. En fait nous savons que si un tel algorithme existe, il n'est
certainement pas \gui{mécaniquement calculable}: un programme qui tourne
sur machine ne peut s\^urement pas accomplir un tel travail même lorsque l'on
impose la limitation sur l'entrée $\alpha$ qu'elle soit une suite binaire
primitive récursive explicite. Cette impossibilité est un grand \tho
d'informatique théorique, souvent indiqué sous l'appellation \gui{\tho de
l'arrêt des programmes}.

%--- Thm{arrêt des programmes}-----------

\mni
{\bf Théorème de l'arrêt des programmes} (On ne peut pas tout savoir)\\
{\it Sous trois formes immédiatement \eqvesz:
%-----------------begin item------------------
\begin{itemize}
\item  On ne peut pas assurer automatiquement la terminaison des programmes:  
il n'existe pas de programme \,$T$\, qui puisse tester si un programme arbitraire \,$P$\,
finira par aboutir à l'instruction Stop.
\item  Il n'existe pas de programme qui puisse tester si une suite primitive
récur\-sive  arbitraire est identiquement nulle.
\item Il n'existe pas de programme \,$U$\, qui prenne en entrée deux
entiers, donne en sortie un booléen, et qui énumère toutes les suites
binaires programmables  (la suite \,$n\mapsto U(m,n)$\, est la
\,$m$-ième suite énumérée par \,$U$).
\end{itemize}
%-----------------end item------------------
 }
%--- end-Thm---------------------------

\medskip
Non seulement ce \thoz, sous sa dernière formulation, ressemble au \tho de
Cantor qui affirme que l'on ne peut pas énumérer l'ensemble des suites
binaires, mais la preuve, très simple, est essentiellement la même.

\smallskip
Bien que le \tho précédent n'interdise pas a priori l'existence d'une
procédure effective mais non mécanisable pour résoudre de manière
systématique ce type de \pbsz, il confirme l'idée intuitive selon
laquelle il faudra toujours faire preuve de nouvelle inventivité pour
progresser dans notre connaissance du monde \mathz.

\smallskip Ainsi, d'un point de vue \cofz, nous rejetons  le {\it Limited
Principle of Omniscience}.
%-----------------begin item------------------
\begin{description}
\item \LPOz:  Si $(\alpha_n)$ est une suite binaire, alors ou bien il existe
un $n$
tel que $\alpha_n=1$, ou bien $\alpha_n=0$ pour tout $n$.
\end{description}
%-----------------end item------------------
%
Le voici sous forme plus concentrée.
\begin{description}
\item \LPOz: \hfill $\forall \alpha \in\NN^\NN, \; \;
(\alpha\not= 0 \;\lor \;\alpha= 0)$\hfill~
\end{description}
%-----------------end item------------------

Nous appellerons \emph{\prt \elrz} une \prt \eqve  à

\snic{\exists n\; \alpha(n)\not=0}

%\sni
pour un certain $\alpha\in\NN^\NN$.

\medskip Le principe \LPO a de nombreuses formes équivalentes. En voici quelques unes.
%-----------------begin item------------------
\begin{enumerate}
\item  Si $A$ est une propriété \emph{\elrz}, on a
$ \; A \lor \lnot A $.
\item Toute suite dans $\NN$ est ou bien bornée, ou bien non bornée.
\item Toute suite décroissante dans $\NN$ est constante à partir d'un
certain rang.
\item D'une suite bornée dans $\NN$ on peut extraire une sous-suite infinie
constante.
\item Toute partie \enum de $\NN$ est détachable.
\item Toute partie \enum de $\NN$ est ou bien finie, ou bien infinie.
\item Pour toute suite double d'entiers $\beta :\NN^2\rightarrow \NN$ on~a:
$$\preskip.4em \postskip.0em \forall n\; \exists m\;\;  \beta(n,m)=0\quad \lor\quad
\exists n\; \forall m\; \; \beta(n,m)\neq 0$$
\item Tout sous-groupe détachable de $\ZZ$ est engendré par un seul
\eltz.
\item Tout sous-groupe de $\ZZ^p$ engendré par une suite infinie est \tfz.
\item $\forall x\in\RR$, $(\; x\not= 0 \; \lor\;  x= 0\; )$.
\item $\forall x\in\RR$, $(\; x> 0 \; \lor\;  x= 0 \; \lor\;  x<0\; )$.
\item Toute suite bornée monotone dans $\RR$  converge.
\item D'une suite bornée dans $\RR$ on peut extraire une sous-suite
convergente.
\item Tout nombre réel est ou bien rationnel ou bien irrationnel.
\item Tout sous-\evc \tf de \,$\RR^n$\, admet une base.
\item Tout espace de Hilbert séparable admet \\
-- ou bien une base
hilber\-tienne  finie \\
-- ou bien une base hilbertienne  \denbz.
\end{enumerate}
%-----------------end item------------------

%---- subsubsec LLPO ----------------------
\subsubsec{Le Mini Principe d'Omniscience}\rdb

\noi  Un autre principe d'omniscience, plus faible, \LLPO  (Lesser Limited
Principle of Omniscience) est le suivant.
%-----------------begin item------------------
\begin{description}\label{LLPO}
\item \LLPOz:  Si $A$ et $B$  sont deux propriétés \elrsz, on a

\snic{\lnot (A\; \land\;  B) \quad \Longrightarrow \quad(\lnot
A\;  \lor\;  \lnot B)}
\end{description}
%-----------------end item------------------
\rdb
Ce principe \LLPO  a de nombreuses formes équivalentes.
%-----------------begin item------------------
\begin{enumerate}
\item $\forall\, \alpha, \,\beta$ suites croissantes $\in\NN^\NN$, si $\, \forall
n\,
     \alpha(n)\beta(n)=0$, alors $\, \alpha=0$ ou~\hbox{$\beta=0$}.
\item \label{recinsep}$\forall \alpha,\beta\in\NN^\NN$, si $\; \forall
n,m\in\NN\;  \alpha(n)\neq
      \beta(m)\;  $ alors $ \; \exists \gamma\in\NN^\NN\;  $ tel que
$$\preskip.4em \postskip.0em 
\forall n,m\in\NN\quad
    \big(\gamma(\alpha(n))=0\; \land\;  \gamma(\beta(m))=1\big) 
$$
\item $\forall \alpha\in\NN^\NN$, $\exists k\in \{0,1\}$, ($\; \exists n\;
     \alpha(n)=0\; \Rightarrow \; \exists m\;
     \alpha(2m+k)=0)$.
\item  $ \forall x \in \RR \quad (\; x\leq 0\;  \lor\;  x\geq 0\; )$
 (ceci permet de faire de nombreuses preuves par dichotomie avec les nombres
réels.)
\item $\forall x,y \in \RR \quad (\; xy=0\; \Rightarrow\;  (\; x= 0 \; \lor\;
y=0\; )\; )$.
\item L'image d'un intervalle $[a,b]\subset \RR$ par une fonction réelle
uniformément continue est un interval\-le~$[c,d]$.
\item Une fonction réelle uniformément continue sur un espace
métrique compact atteint ses bornes.
\item \kl (une des versions du lemme de K\"onig) Tout arbre infini explicite
à embranchements finis possède une branche infinie. \label{K1LLPO}
\end{enumerate}
%-----------------end item------------------

Il est connu que si un algorithme existe pour le troisième item il ne peut
pas être  \gui{mécaniquement calculable} (i.e., récursif): on peut
construire $\alpha$ et $\beta$ mécaniquement calculables vérifiant
l'hypothèse, mais pour lesquels aucun $\gamma$ mécaniquement calculable
ne vérifie la conclusion.
De même, l'arbre singulier de Kleene est un arbre récursif  \denb infini
à embranchements finis qui ne possède aucune branche infinie
récursive. Ceci donne un \gui{contre-exemple récursif} pour \klz.

\smallskip Nous allons montrer maintenant\footnote{Comme pour toutes les \dems dans cette annexe, elle est informelle et l'on ne précise pas dans quel cadre formel elle pourrait être écrite. \Llec repérera dans cette \dem une utilisation d'une construction
par \recu qui relève en fait de l'axiome du choix dépendant, \gnlt considéré comme non problématique en \comaz.} l'\eqvc
\kl $\Leftrightarrow$ \LLPOz.

Un arbre infini explicite à embranchements finis peut être décrit par
un ensemble $A\subset \Lst(\NN)$ de listes d'entiers vérifiant les
propriétés suivantes (les quatre premières correspondant à la notion
d'arbre explicite à embranchements finis).
%-----------------begin item------------------
\begin{itemize}
\item  La liste vide $[\,]$ représente la racine de l'arbre, elle
appartient à $A$,
\item  un $a=[a_1,\ldots,a_n ]\in A$ représente à la fois un noeud de
l'arbre et le chemin qui mène de la racine jusqu'au noeud,
\item  si $[a_1,\ldots,a_n ]\in A$ et $n\geq 1$, alors $[a_1,\ldots,a_{n-
1}]\in A,$
\item  si $a=[a_1,\ldots,a_n ]\in A$ alors les  $x\in\NN$ tels que
$[a_1,\ldots,a_n,x ]\in A$ forment un segment $\sotq{x\in\NN}{x<\mu(a)}$
où $\mu(a)$ est donné explicitement en fonction de $a$: les branches
issues de $a$ sont numérotées $0,\ldots,\mu(a)-1$.
\item  Pour tout $n\in\NN$ il y a au moins un $[a_1,\ldots,a_n]\in A$
(l'arbre est explicitement infini).
\end{itemize}
%-----------------end item------------------
Ainsi la partie $A$ de
$\Lst(\NN)$ est   détachable (c'est en fin de compte ce que signifie ici le
mot \gui{explicite}).  Et $A$ est \denbz.

%-----------------begin proof------------------
%
\begin{Proof}{\Demo de \kl $\Leftrightarrow$ \LLPOz.}\\
Nous prenons pour \LLPO la variante donnée dans l'item 1.\\
Supposons \klz. Soit $\alpha,\beta\in\NN^\NN$ comme dans l'item 1.
Considérons l'arbre suivant. Après la racine on ouvre deux branches qui
se poursuivent indéfiniment sans jamais créer de nouveaux
embranchements, jusqu'à ce que~\hbox{$\alpha(n)\neq 0$} ou $\beta(n)\neq 0$ (si
jamais cela se produit). Si cela se produit avec~\hbox{$\alpha(n)\neq 0$}, on
arrête la branche de gauche et l'on continue celle de droite. Si c'est avec~$\beta(n)=0$, on fait le contraire. Donner explicitement une branche infinie
dans cet arbre revient à certifier d'avance que  $\alpha=0$ ou $\beta=0$.

 Inversement supposons  \LLPOz. Considérons un arbre infini  explicite
à embranchements finis. Supposons \spdg  que l'arbre est binaire: au
delà d'un noeud il y a au plus deux branches. Nous prouvons par
récurrence que nous pouvons sélectionner jusqu'à la profondeur $n$ un
chemin qui aboutit à un noeud $K_n$ en dessous duquel l'arbre est infini.
Ceci est vrai pour $n=0$ par hypothèse. Si cela est vrai pour $n$, il y a
au moins une branche en dessous du noeud $K_n$ sélectionné. S'il y en a
deux, considérons les suites  $\alpha_n$ et $\beta_n\in\NN^\NN$ définies
comme suit:\\
--- $\alpha_n(m)=0$ si il y a au moins une branche de longueur $m$ en
dessous de~$K_n$ partant sur la droite, sinon $\alpha_n(m)=1$ \\
--- $\beta_n(m)=0$ si il y a au moins une branche de longueur $m$ en dessous
de~$K_n$ partant sur la gauche, sinon $\beta_n(m)=1$.\\
Par \hdr les suites $(\alpha_n)_{n\in\NN}$ et $(\beta_n)_{n\in\NN}$ sont  croissantes et leur produit
est nul.  On applique l'item~1 de \LLPO: l'une des deux suites est nulle et cela nous
donne le moyen de sélectionner le chemin vers la droite ou
celui vers la gauche.
\end{Proof}
%-----------------end proof------------------

%---- subsubsec  PTE ----------------------
\subsubsec{Le Principe du Tiers Exclu}
%------------------------------------------------------------------

\noi  Le Principe du Tiers Exclu (\TEMz) affirme que $P\vee \lnot P$ est
vrai pour toute proposition $P$. Ce principe d'omniscience extrêmement
fort implique \LPOz. Il suppose de manière implicite que des ensembles
tels que $\NN$ ou $\NN^\NN$ ou même nettement plus compliqués, sont des
\emph{infinis actuels}. Il implique également que tout ensemble $X$ est
discret  si l'on définit
$x\neq_Xy$ comme signifiant~\hbox{$\lnot (x=_Xy)$}.

%--- Section Principes problematiques ----------
\section[Principes problématiques \ldots]{Principes problématiques en
\coma}
\label{para pripro}\relax
%------------------------------------------------------------------
Nous entendons par \emph{principe problématique} un principe qui,
quoi\-que vérifié en pratique si l'on fait des \coma dans le style de
Bishop, est indémontrable \cotz. En \clamaz, ces principes peuvent être connus
comme vrais ou connus comme faux. 
\\
Par exemple, en pratique, chaque fois
qu'un $\alpha \in\NN^\NN$ est bien défini \cotz, il peut être calculé
par un programme. \\
Autrement dit, en pratique, la \textbf{Fausse Thèse de Church},
que l'on peut écrire  sous la forme \fbox{$\Rec=\NN^\NN$}, est vérifiée en \comaz. Mais elle ne peut pas être démontrée dans le cadre minimaliste
des \coma à la Bishop, qui sont compatibles avec les \clamaz.
Car la Fausse Thèse de Church est un
principe faux en \clamaz, en vertu d'un argument de cardinalité. Par contre, les \coma russes le prennent comme un
axiome fondamental.

\ss Nous allons ici examiner (brièvement) seulement deux principes
pro\-blé\-ma\-tiques, tous deux vrais en \clamaz. \rdb

%:---- Le Principe de Markov -----
\subsubsec{Le Principe de Markov} \label{principeMarkov}

\noi  Le \emph{Principe de Markov}, \MPz, est le suivant:
$$ \forall x \in \RR \quad (\lnot x= 0 \Rightarrow  x\neq 0).
$$

Affirmer  \MP  revient à dire: pour toute suite binaire $\alpha$, s'il est
impossible que tous ses termes soient nuls, alors il doit y avoir un terme
non nul.

C.-à-d. encore: si $A$ est une propriété \elr alors $\lnot\lnot A
\Rightarrow A$.

L'école  \cov russe admet \MPz. En fait,
pour un $\alpha\in\NN^\NN$, il semble impossible de donner une preuve
constructive  de $\lnot (\alpha=0)$ sans trouver un $n$ tel que
$\alpha(n)\neq 0$. Ainsi \MP est valide d'un point de vue pratique dans le
constructivisme à la Bishop.  Remarquons  aussi que \LPO  implique
clairement~\MPz.

%:---- Principes de Continuité Uniforme --
\subsubsec{Principes de continuité uniforme}

\noi  Le principe de continuité uniforme  affirme que toute fonction
ponctuellement continue sur un espace métrique compact est uniformément
continue. Il est équivalent à la même affirmation dans un cas
particulier, qui est elle-même très proche de l'une des formes
classiques du lemme de K\"onig. Les principes problématiques suivants
semblent intéressants à étudier dans leurs relations mutuelles,
d'autant plus qu'ils apparaissent très souvent en analyse classique% (cf.
%la section \ref{ssec LPO KL})
.

%-----------------begin description--------
\begin{description}
\item [\UCp]  Toute fonction ponctuellement continue $f:X\rightarrow Y$,
avec
$X$ espace métrique compact et $Y$ espace métrique, est uniformément
continue.
\item [\UC]  Toute fonction ponctuellement continue $f:\{0,1\}^\NN\rightarrow
\NN$ est
uniformément continue.
\item [\Mini]  Toute fonction réelle $>0$  uniformément continue sur un
espace métrique compact est minorée par un réel  $>0$.
\item [\Minim]  Toute fonction réelle $>0$  uniformément continue sur un
intervalle compact $[a,b]$ est minorée par un réel  $>0$.
\item [\Minip]  Toute fonction réelle $>0$  ponctuellement continue sur un
espace métri\-que compact est minorée par un réel  $>0$.
\item  [\KL]  Un arbre binaire explicite $A$ qui ne possède pas de
branche infinie (i.e., $\forall \alpha\in \{0,1\}^\NN\; \exists m\in \NN\;
\alpha |^m\notin A$) est fini.
%\item  [\KLp] Un arbre binaire énumérable qui ne possède pas de
%branche
%infinie a une profondeur bornée{\footnote{Il revient au même de dire
%qu'il est borné en nombre.}}.
\end{description}
%-----------------end description------------

Dans la formulation \KLz, on voit que ce principe est apparemment voisin de
\LLPO (voir la dernière forme équivalente  \kl citée \paref{K1LLPO}). En
fait, on peut montrer qu'il est une conséquence de \LPOz%
% (voir la section \ref{ssec LPO KL})
. Mais il ne s'agit pas d'un principe d'omniscience.
D'ailleurs, il n'implique pas \LLPOz. En \comaz, \LLPO est manifestement faux en
pratique, tandis que \KL est vérifié en pratique. Car chaque fois
que l'on sait prouver \cot qu'un arbre à embranchements finis n'a pas de
branche infinie, on sait également prouver qu'il est fini.

%:section: Exercices
\penalty-2500
\Exercices{

\Oui{\setcounter{exercise}{0}}
%--- Exercise{exo}-------------
\begin{exercise}
\label{exoAnnexe1}
{\rm  Donnez des \dems pour les exemples
 \iref{exo neq dans N^N}, \iref{exo inf deno1}, \iref{exo log0}, \iref{exo log1}  
 et~\iref{exo log2}.
 }
\end{exercise}
%--- end -exercise-----------------------------------------

%--- exercise{exo inf deno2}-----------
\begin{exercise}
\label{exo inf deno2} 
{\rm
Expliquez pourquoi les  notions d'ensemble fini, finiment énu\-mérable,
borné en nombre, faiblement fini, énumérable et borné en nombre, ne peuvent pas être
identifiées en \comaz. Expliquez pourquoi ces notions coïncident si
l'on admet le principe du tiers exclu.
}
\end{exercise}
%--- end-exercise-----------------------------------------

%--- exercise{exoA22}-------------------
\begin{exercise}
\label{exoA22}
\emph{Démontrer quelques unes des équivalences signalées pour \LPOz.
}
\end{exercise}
%--- end-exercise-----------------------------------------

%--- exercise{exoA23}-------------------
\begin{exercise}
\label{exoA23}
\emph{Démontrer quelques unes des équivalences signalées pour \LLPOz.
}
\end{exercise}
%--- end-exercise-----------------------------------------

}% fin des exos
%:  solutions

%\sol{

%
%}% fin des solutions d'exos

%:  ---- Section*{references}-----------
\vspace{5pt}
\Biblio
\vspace{3pt}

La polémique sur la nature et l'usage de l'infini en \maths
a été très vive au début du 20\ieme\  siècle:
voir par exemple Hilbert \cite[1926]{Hi2}\ihiz, Poincaré \cite[1909]{Po}, H. Weyl \cite[1918]{W1}, 
 \cite[1951]{Bro} et \cite[1987]{INS}).
Le débat a semblé en un premier temps se terminer à l'avantage du point
de vue représenté par la logique classique.
En fait, depuis les années 60 et notamment la parution du livre de Bishop,
 les deux points de vue sont nettement moins opposés qu'il ne pouvait
 paraître. \\
 Quelques références intéressantes sur ce thème: \cite[1962]{Loren}, \cite[Fred Richman, 1990]{Ri90}, \cite[2007]{Dow2} et \cite[Per Martin-L\"of, 2008]{ML2008}.

La logique constructive s'appelle souvent \gui{logique intuitionniste}.
Elle a été mise au point en tant que \sys formel par A.\ Heyting.

On trouve des exposés agréables de tels \syss formels
(en confrontation avec les \syss correspondant à la logique classique avec Tiers Exclu) dans les livres \cite[1962]{Loren} et \cite[2001]{DNR}. 
Le livre  plus difficile \cite[2012]{SW} est une référence essentielle pour une étude plus avancée des rapports entre logique classique, logique intuitionniste et logique minimale. 
Le petit livre \cite[1995]{Dow} donne quant à lui une présentation informelle
intéressante.

Concernant la discussion sur les rapports entre effectivité et récursivité
voir \cite[Coquand]{CoqRec}, \cite[Heyting]{HeyRec} et \cite[Skolem]{SkoRec}.

Le livre \cite[1985]{Be} fait une étude systématique de nombreux principes
 problé\-ma\-tiques en \comaz.
 Pour l'arbre singulier de Kleene voir \cite[page 68]{Be} et \cite[1965]{KV}.

La mise au point et la comparaison de \syss formels pouvant servir de cadre aux \coma pratiquées dans \cite{B67} ou \cite{MRR} est depuis longtemps
un sujet très actif de recherche.
On notera l'influence prépondérante de la théorie \cov
des types \tsbf{CTT} de Per Martin-L\"of,  
%: 2015\cite[1984]{PML}, \cite{ML}, 
\cite{ML1973,ML} et \cite[1984]{PML}, 
et de la théorie \tsbf{CZF} de Peter Aczel et Michael Rathjen (\cite[Aczel]{Acz} et \cite{AcRa}). Voir aussi les développements récents dans \cite[2014]{HoTT} et la page web de Thierry Coquand: \url{http://www.cse.chalmers.se/~coquand/}.

Citons \egmt le beau livre \cite[1998]{fef} qui se situe dans la lignée des
propositions d'Hermann Weyl.

Pour une discussion sur le \gui{Fan Theorem}
voir~\cite[Coquand]{CoqFT}.

L'étude systématique de la comparaison  (en \comaz) de principes  d'omniscience (tels que \LPO ou \LLPOz), ainsi que celle de principes
 problématiques (tels que \MP ou \FANz)  a pris récemment
un essor important. On pourra consulter à ce sujet \cite[Berger\&al.]{Berg,BeBr,BeIs} et \cite[Ishihara]{Ishi,Ishi2,Ishi3}.

\newpage \thispagestyle{CMcadreseul}
\incrementeexosetprob

 %%%%%%%%%%%%%%%%%%%%%%%%%%%%%%

\pagestyle{CMExercicesheadings}
%!TEX root =  ACMC-A.tex

%%@@@@@@@@@------- table des theoremes ----------
\finincrementeexosetprob
\cleardoublepage \thispagestyle{CMcadreseul}
\refstepcounter{chapter}
\addcontentsline{toc}{chapterbis}{Tables des th\'eor\`emes}
\addtocontents{toc}{\vskip-0.8em}

\catcode`\@=11
\@makesindexhead{Tables des th\'eor\`emes}
\perso{compil\'e le \today}
\pagestyle{CMExercicesheadings}

\setlength\parskip{1pt}

\markboth{Tables des \thosz}{Tables des \thosz}

% Commandes
% \newcommand \TTD[2]{%
% \sni\mbox{\parbox[t]{.775\textwidth}{#2}%
% \hspace{.05\textwidth}% 
% \parbox[t]{.175\textwidth}{#1}}}

\newcommand \TTT[3]{%
%\sni
\noindent 
\mbox{%
\parbox[b]{.80\textwidth}{\leftskip10pt\parindent-\leftskip\strut#3\dotfill\par}%
\hspace{.01\textwidth}% 
\parbox[b]{.12\textwidth}{\hfill#1}%
\hspace{.01\textwidth}% 
\parbox[b]{.06\textwidth}{\hfill#2}%
}%
}

\newcommand \TTref[2]{\TTT{\ref{#1}}{\pageref{#1}}{#2}\par}

\newcommand \TTpa[3]{\TTT{\ref{#3}}{\pageref{#1}}{#2}\par}

\newcommand \chap[1]{

\goodbreak
\mni\textbf{#1}}

%%%%%%%%%%%%%%%%%%%%%%%%%%%%%%%%%%%
\small
%:  Méthodes dynamiques
\begin{center}
    {\large\bf Méthodes dynamiques}\label{machineries}
\end{center}

\vspace{2mm}
%\vskip2pt

\TTT{~}{page}{Nom} \par\vspace{1mm}

\TTpa{MethodeQI}{Machinerie \lgbe \elr des anneaux \qisz}{subsecAnneauxqi}
\TTpa{MethodeZedRed}{Machinerie \lgbe \elr des anneaux \zeds réduits}{secKrull0dim}
\TTpa{subsecDyna}{La méthode dynamique}{subsecDyna}
\TTpa{MetgenAnar}{Machinerie \lgbe  des \anars}{secIplatTf}
\TTpa{MethodeIdeps}{Machinerie \lgbe de base (\`a \idepsz)}{secMachLoGlo}
\TTpa{MethodeIdemax}{Machinerie dynamique \`a \idemasz}{subsecLGIdeMax}
\TTpa{MethodeIdemin}{Machinerie dynamique \`a \idemisz}{subsecLGIdepMin}
\TTpa{MachDynAxAx}{Machinerie dynamique avec $\ArX$ et $\gA(X)$}{sec.Etendus.Valuation}

\vspace{5mm}
%:  Principes local-globals concrets

\begin{center}
    {\large\bf Principes local-globals concrets}
\end{center}
\vspace{2mm}
\TTref{plcc.basic}{Principe local-global concret de base}
\TTref{PrTransfertBasic}{Principe de Transfert de base}
\TTref{plcc.coh}{Modules \cohsz}
\TTref{plcc.tf}{Modules \tfz}
\TTref{plccRangMat}{Rang d'une matrice}
\TTref{fact.lnl.loc}{Applications \lins \lnlsz} 
\TTref{plcc.basic.modules}{Suites exactes de modules} 
\TTref{plcc.basic.monoides}{Pour les \mosz} 
\TTref{plcc.entier}{\'Eléments entiers sur un anneau}
\TTref{plcc.scinde}{Propriétés d'\alis entre \mpfsz}
\TTref{plcc.pf}{Modules \pfz}
\TTref{plcc.aqi}{Anneaux \qisz}
\TTref{plcc.cor.pf.ptf}{Modules \ptfsz}
\TTref{plcc.aGs}{Alg\`ebres galoisiennes}
\TTref{plcc.plat}{Modules plats}
\TTref{plcc.arith}{Anneaux \lsdzz, \arisz, de Pr\"ufer, \ids \lopsz}
\TTref{plcc.algfptes}{Alg\`ebres plates ou \fptesz, \lon en bas}
\TTref{thlgb1}{Matrices semblables, ou équivalentes (anneau \lgbz)}
\TTref{thlgb2}{Modules \pf isomorphes (anneau \lgbz)}
\TTref{thlgb3}{Modules quotients (anneau \lgbz)}
\TTref{plcc.normal}{Anneaux normaux et \ids \iclz}
\TTref{plcc.agq}{\'Eléments primitivement \agqsz}
\TTref{plcc.ddk}{Anneaux de Dedekind}
\TTref{defSeqSing}{Suites singulières}
\TTref{thDdkLoc}{Dimension de Krull des anneaux}
\TTref{plcc.KdimMor}{Dimension de Krull des morphismes}\rdb
\TTref{plcc.sli}{Suites exactes et généralisations}
\TTref{plcc.ptf}{Propriétés de finitude pour les modules}
\TTref{plcc.propaco}{Propriétés des anneaux commutatifs}
\TTref{plcc.apf}{Propriétés de finitude pour les \algsz, \lon en bas}
\TTref{plcc2.apf}{Propriétés de finitude pour les \algsz, \lon en haut}
\TTref{plcc.modules 1}{Recollement concret d'\elts dans un module, 
 ou d'\homos entre modules}
\TTref{plcc.modules 2}{Recollement concret de modules}
\TTref{plcc.RecolHomAnn}{Recollement concret d'\homos d'anneaux}
\TTref{plcc.regularite}{Principe \lgb pour l'\egt en profondeur 1}
\TTref{plcc0Profondeur2}{Principes \lgbs en profondeur 2}
\TTref{plcc.modules1bis}{Recollement concret d'\elts dans un module
en \prof 2}
\TTref{plcc.modules2bis}{Recollement concret de modules en \prof 2}
\TTref{thPatchV}{Recollement de Vaserstein: matrices équivalentes sur $\gA[X]$}
\TTref{thPatchQ}{Recollement de Quillen: modules étendus (Quillen patching)}
\TTref{plgcetendus}{Principe local-global \`a la Roitman}
%:2017 rectification d'un anglicisme
\TTref{plcc.159}{Recollement concret dans le groupe \elrz}
\TTref{plgc-Rao}{Principe \lgb concret de Rao}
\label{plgvarios}

\vspace{5mm}
%

%:  Principes de recouvrement fermé

\begin{center}
    {\large\bf Principes de recouvrement fermé}
\end{center}

\vspace{2mm}

%\TTT{\num}{page}{Nom} \par\vspace{3mm}
\TTref{prcfgrl}{Pour les \grlsz}
\TTref{prcf1}{\'Eléments nilpotents, \comz}
\TTref{prcf2}{Modules \tfz}
\TTref{prcf3}{Rang d'une matrice,  \mptfsz}
\TTref{thDdkRecFer}{Dimension de Krull}
% 
% \TTref{prcf}{}

\vspace{5mm}
%

%:  Stabilité par \eds
%\pagebreak
\begin{center}
    {\large\bf Stabilité par \eds}
\end{center}

\vskip2mm

%\TTT{\num}{page}{Nom} \par\vspace{3mm}
\TTref{propPfExt}{Modules \tf et \pfz, produits tensoriels, puissances \smqs et extérieures, \alg extérieure}
\TTref{fact.idf.change}{Idéaux de Fitting}
\TTref{propPtfExt}{Modules \ptfsz}
\TTref{fact.det loc}{Déterminant, \polcarz, \polfonz, \polmuz,
\endo cotransposé} 
\TTref{factEdsAlg}{Alg\`ebres \tfz, \pfz, \stfesz}
\TTref{factEdsDualisante}{Formes dualisantes, alg\`ebres de Frobenius}
\TTref{factEdsEtale}{Alg\`ebres  \stesz}
\TTref{factSpbEds}{Alg\`ebres \spbsz}
\TTref{factEdsAG}{Automorphismes séparants}
\TTref{corAGextsca}{Alg\`ebres galoisiennes}
\TTref{factChangeBase}{Alg\`ebre de \dcn \uvlez}
\TTref{corPlatEds}{Modules plats}
\TTref{propFidPlatPrAlg}{Réciproques dans le cas des extensions \fptesz}
%:HHH supprime ci dessous
%\TTref{propComparRedRed}{$\GKO(\gA)\simeq\GKO(\Ared)$}
% \TTref{propOmExt}{Module des \dilesz}
% \TTref{}{}
% \TTref{}{}
% \TTref{}{}
% \TTref{}{}

\vspace{5mm}
%
%:  Théor\`emes
\begin{center}
\bigskip    {\large\bf Théor\`emes}
\end{center}

%\TTT{\num}{page}{Nom} %\par\vspace{3mm}
%: 2

\CHAP{\nameref{chapSli}}

%:HHH ajout
\TTref{plcc.basic}{Principe local-global concret de base}
\TTref{lemGaussJoyal}{Lemme de Gauss-Joyal}
\TTref{propCoh4}{Caractérisation des \comosz}
\TTref{fact.sfio}{Syst\`eme fondamental d'\idms \ortsz}
%:HHH l'ordre des deux qui suivent a été changé
\TTref{lem.ide.idem}{Lemme de l'\itf \idmz}
\TTref{restes chinois}{\Tho des restes chinois, forme générale (pour la forme \ari voir le \thref{thAnar})}
\TTref{lem.min.inv}{Lemme du mineur inversible}
\TTref{lem pf libre}{Lemme de la liberté}
\TTref{lemCram}{Formule de Cramer généralisée}
\TTref{propIGCram}{Formule magique \`a la Cramer} 
\TTref{propFactDirRangk}{Sous-\mtfs en facteur direct dans un module libre} 
\TTref{prop inj surj det}{Crit\`eres d'injectivité et de surjectivité} 
\TTref{theoremIFD}{Matrices \lnlsz}
%:HHH 2 rajouts
\TTref{Th.transitivity}{Formules de transitivit\'{e} pour la trace et le \deterz}
\TTref{thTransDisc}{Formule de transitivit\'{e} pour les \discrisz}

%: 3
\CHAP{\nameref{chapGenerique}}

\TTref{thSymEl}{Polyn\^omes \smqs \elrsz}
\TTref{lemdArtin}{Lemme de Dedekind-Mertens}
\TTref{thKro}{\Tho de \KRO (1)}
\TTref{propUnicCDR}{Unicité du \cdr  (cas \stfz)}
\TTref{thPrIs}{\Tho de prolongement des \isosz}
\TTref{thGaloiselr}{Correspondance galoisienne}
%\TTref{}{}
\TTref{thResolUniv}{Construction d'un corps de racines}
\TTref{LemElimAffBasic}{Lemme d'élimination de base}
\TTref{thIntClosStab}{Anneau de \pols intégralement clos}
\TTref{thClotAlgQ}{Corps de racines, \tho de l'\elt primitif}
\TTref{th1IdZalpha}{Tout \itf non nul d'un corps de nombres est \ivz}
\TTref{th2IdZalpha}{Structure multiplicative des \itfs d'un corps de nombres}
\TTref{propEvitementConducteur}{\Tho de Dedekind, \ids qui évitent le conducteur}
\TTref{thNstfaibleClass}{\nst faible et mise en position de \Noez, voir aussi
le \tho \ref{thNstNoe}}
\TTref{thNstClass}{\nst classique}
\TTref{thNSTsurZ}{\nst sur $\ZZ$, \nst formel}
\TTref{corthNSTsurZ}{\nst sur $\ZZ$, \nst formel, 2}
\TTref{thNewtonQuad}{Méthode de Newton}
\TTref{corIdmNewton}{Lemme des \idms résiduels}

%: 4
\CHAP{\nameref{chap mpf}}

\TTref{lem pres equiv}{Matrices qui présentent le m\^eme module}
\TTref{propsregpf}{Un \id engendré par une suite réguli\`ere est \pfz}
\TTref{thidptva}{L'\id d'un point est un \mpfz}
\TTref{propCoh2}{Cohérence et \pn finie (voir aussi la proposition \ref{propAliCoh})}
\TTref{prop unicyc}{Somme directe de modules cycliques (unicité)}
\TTref{prop quot non iso}{Un quotient isomorphe est un quotient par $0$}
\TTref{thScindageQi}{Lemme de scindage \qiz}
\TTref{thScindageZed}{Lemme de scindage \zedz}
\TTref{thZerDimRedLib}{Le paradis des anneaux \zedrsz}
\TTref{thSPolZed}  {Syst\`eme \poll \zed sur un \cdiz}
\TTref{thStickelberger}  {\Tho de Stickelberger (\sys \zedz)}
\TTref{fact.idf.ann}{Lemme du premier idéal de Fitting}
\TTref{LemElimAff}{Lemme d'élimination général}
\TTref{thElimAff} {\Tho d'\eli \agqz, idéal résultant}

%: 5
\CHAP{\nameref{chap ptf0}}

\TTref{propdef ptf}{Modules \ptfsz}
\TTref{prop pf ptf}{Matrice de \pn d'un \mptfz}
\TTref{corlemScha}{Lemme de Schanuel}
\TTref{propIsoIm}{Lemme \dlgz}
\TTref{propZerdimLib}{Lemme de la liberté \zedez: point \emph{\iref{ite2propZerdimLib}.} du \thoz}
\TTref{propZerdimLib}{\Tho de la base incompl\`ete: point \emph{\iref{ite5propZerdimLib}.} du \thoz}
\TTref{corBass2}{\Tho de Bass, modules stablement libres}
\TTref{prop Fitt ptf 2}{\Tho de structure locale des \mptfsz. 
Voir aussi les \thos \ref{theoremIFD},
\ref{corth.ptf.sfio}, \ref{th ptf loc libre} et~\ref{th matproj}}
%\TTref{prop Fitt ptf 2}{Caractérisation des \mptfs par 
%leurs \idfsz. Voir aussi le \thref{theoremIFD}}
\TTref{factLocCas}{Lemme des localisations successives, 1}
\TTref{propmlm}{Modules \tf \lmosz, voir aussi  \ref{pmlm}}
\TTref{propdef det ptf}{Déterminant d'un \endo d'un \mptfz}
\TTref{th ptf sfio} {Le \sfio associé \`a un \mptfz}
\TTref{lem calculs}{Calculs explicites: \deterz, \polcarz, \&ct}
\TTref{th decomp ptf}{Décomposition d'un \mptf en somme 
directe de modules de rang constant}

%: 6
\CHAP{\nameref{chap AlgStricFi}}

\TTref{th1Etale}{\Tho de structure des \Klgs étales, 1}
\TTref{cor1lemEtaleEtage}{\'Eléments \spls dans une \Klgz}
\TTref{corlemEtaleEtage}{Caractérisation des \Klgs étales}
\TTref{thEtalePrimitif}{\Tho de l'\elt primitif}
\TTref{th2Etale}{\Tho de structure des \Klgs étales, 2}
\TTref{thClsep}{Cl\^oture \spbz}
\TTref{thNormalEtaleGalois}{Caractérisation des extensions galoisiennes}
\TTref{thCGSynthese}{Correspondance galoisienne, synth\`ese}
\TTref{factSDIRKlg}{Somme directe dans la catégorie des \klgsz}
\TTref{lemLingOver}{Lying over: voir aussi le lemme \vref{lemLingOver2}}
\TTref{thNst0}{Un \nst faible}
\TTref{propAlgFinPresfin} {\klgs qui sont des \kmos \pfz}
\TTref{Thextent}{Extension enti\`ere et \iclz e d'un anneau \iclz}
\TTref{propTraptf}{Transitivité pour les \algs \stfesz}
\TTref{factCarDua}{\Carn des formes dualisantes dans le cas \stfz}
\TTref{factCarAste}{\Carn des \asesz}
\TTref{prop2EtaleReduit}{Une \klg \ste est réduite}
\TTref{thIdmEtale}{Idempotents et extension des scalaires dans les \ases}
\TTref{thAlgSteIdmSpb}{Idempotent de \spt d'une \asez}
\TTref{thSepIversen}{\Prts \caras des \klgs \spbsz}
\TTref{thAlgStfSpbSte}{Une \alg \spb \stfe est \stez}
\TTref{thSepProjFi}{Propriété de finitude des \algs \spbsz}
\TTref{corthSepProjFicdi}{Sur un \cdiz, une \alg \spb \pf est \stez.}
\TTref{thDA}{Lemme de Dedekind}
\TTref{thA}{\Tho d'Artin, version  \aGsz}
\TTref{corAGextsca}{Extension des scalaires pour les \aGsz}
\TTref{thAGACar}{\Carn des \aGsz}
\TTref{corAGAlibreCar}{\Carn des \aGs libres}
\TTref{corDAexplicite2}{Correspondance galoisienne, version  \aGsz}
\TTref{thCorGalGen}{Correspondance galoisienne, \aGs connexes}
\TTref{corAGQuo}{Quotient de Galois d'une \aGz}
\TTref{exoLuroth1}{\Tho de L\"uroth (exercice)}{}
%\TTref{}{}
%\TTref{}{}
%\TTref{}{}

%: 7
\CHAP{\nameref{ChapGalois}}

\TTref{thNstfaibleClassSCA}{\nst faible et mise en position de \Noez, 2}
\TTref{thNstNoe}{\nst faible et mise en position de \Noez, 3}
\TTref{thNoetSimult}{Mise en position de \Noe simultanée}
\TTref{thNstClassCof1}{\nst classique, version \cov \gnlez}
\TTref{thNstClassCof2}{\nst avec multiplicités}
\TTref{thpolcohfd}{Une \apf sur un \cdi est un \cori  \fdiz}
\TTref{corpropBoolFini}{\Tho de structure des \agBs finies}
\TTref{lemAduIdmA}{\Tho de structure galoisien (1), $G$-\algs de Boole}
\TTref{thADG1Idm}{\Tho de structure galoisien  (2), quotients de Galois d'une \apGz}
\TTref{lemAdu1}{Alg\`ebre de \dcn \uvle et séparabilité. Voir aussi le \thref{theoremAdu2}}
\TTref{theoremAdu1}{Alg\`ebre de \dcn \uvle et points fixes}
\TTref{thAduAGB}{L'\adu comme \aGz}
\TTref{theoremAdu3}{Diagonalisation d'une \aduz, voir aussi le \thref{theoremAdu4}}
\TTref{thidGTri}{Base triangulaire de l'\id définissant une \aGz}
\TTref{thUnici}{Unicité éventuelle du \cdr d'un \pol \splz}
\TTref{thdivzeridm}{Gestion dynamique d'un \cdrz, voir aussi le \thref{corRvRel1}}
\TTref{propUnicite}{Unicité du \cdrz, version dynamique}
\TTref{thStruc3}{\Tho de structure galoisien (3), quotients de Galois de l'\adu d'un \pol \spl sur un \cdiz}
\TTref{thZsuffit}{O\`u se passent les calculs: le sous-anneau $\gZ$ de $\gK$
est bien suffisant}
\TTref{exothNst1-zed}{\nst et mise en position de Noether, cas des anneaux \zeds réduits (exercice)}
%\TTref{}{}
%\TTref{}{}
%\TTref{}{}

%: 8
\CHAP{\nameref{chap mod plats}}

\TTref{thPlat1}{Caractérisation des modules plats, 1}
\TTref{cor pf plat ptf}{Caractérisation des \mptfs par la platitude}
\TTref{thplatTens}{Caractérisation des modules plats, 2}
\TTref{propPlatQuotientdePlat}{Quotients plats}
\TTref{thExtPlat}{Caractérisation des \algs plates}
\TTref{thExtFidPlat}{Caractérisation des \algs \fptesz}
\TTref{thSurZedFidPlat}{Toute extension d'un \cdi est \fptez}
\TTref{propFidPlatTf}{Extensions \fptes et \prts de finitude des modules}
\TTref{propFidPlatPrAlg}{Extensions \fptes et \prts de finitude des \algsz}

%: 9
\CHAP{\nameref{chap Anneaux locaux}}

\TTref{thJacUnitEntieres}{Radical de Jacobson et unités d'une extension enti\`ere}
\TTref{thJacplc}{Propriétés locales d'extensions enti\`eres}
\TTref{lemNaka}{Lemme de Nakayama (le truc du déterminant)}
\TTref{lelilo}{Lemme de la liberté locale}
\TTref{lelnllo}{Lemme de l'application \lnlz}  
\TTref{lemnbgtrlo}{Lemme du nombre de \gtrs local}
\TTref{lemLocaliseFini}{Lemme du localisé fini}
\TTref{lemLocalisezeddim}{Lemme du localisé \zedz}
\TTref{thCotangent}{Espace cotangent en $\uxi$ et $\fm_\uxi / {{\fm_\uxi}^2}$}
\TTref{thJZS}{Zéro simple}
\TTref{thJZScdi}{Zéro isolé simple}
\TTref{thPointLisseCourbeIC}{L'\id d'un point non singulier d'une courbe \lot intersection compl\`ete. Voir aussi le \thref{th2PtlisseCourbeIC}}
\TTref{thLgbExtEnt}{Extension enti\`ere d'un \algbz}
%\TTref{}{}

%: 10
\CHAP{\nameref{chap ptf1}}

\TTref{th rg const loc libre}{Les modules de rang constant sont \lot libres}
\TTref{th ptf loc libre}{Les modules \ptfs sont \lot libres;
voir aussi les \thos \ref{thptfloli} et \vref{th matproj}}
\TTref{propRgConstant2}{Modules de rang constant $k$ comme sous-modules
de $\Ae k$}
\TTref{corthTransPtf}{\Algs \stfesz: formule de transitivité pour les rangs}
\TTref{thBnk}{Le foncteur $\Gnk$}
\TTref{leli2}{Deuxi\`eme lemme de la liberté}  
\TTref{prop1TanGrassmann}{Espace tangent \`a une grassmannienne,
voir aussi le \thref{prop3TanGrassmann}}
\TTref{thBézoutLib}{Tout \mrc sur un anneau de Bézout \qi est libre}
\TTref{thPicKTO}{$\Pic \gA$ et $\KTO \gA$}
\TTref{propRgConstant3}{Groupe de Picard et groupe de classes}
\TTref{propComparRedRed}{$\GKO(\gA)\simeq\GKO(\Ared)$}
\TTref{thCarreMil1}{Carré de Milnor}
\TTref{thIdenPtsSimples2}{Une classification compl\`ete de $\GKO(\gA)$;
voir aussi le \thref{thIdenPtsSimples}}

%: 11
\CHAP{\nameref{chapTrdi}}
\nobreak

\TTref{thBoolGen}{Alg\`ebre de Boole librement engendrée par un \trdiz}
\TTref{th0GpRtcl}{Distributivité dans les \grlsz}
\TTref{th1GpRtcl}{\Tho de Riesz (\grlsz)}
\TTref{th2GpRtcl}{\Tho de \dcn partielle sous condition \noeez}
\TTref{thAgcdNormal}{Un anneau \`a pgcd int\`egre est \iclz.}
\TTref{propGCDDim1}{Un anneau int\`egre \`a pgcd de dimension $\leq 1$ est un anneau de Bézout}
\TTref{thAXgcd}{Anneaux \`a pgcd int\`egres: $\gA$ et $\AX$}
\TTref{thZedGen}{Cl\^oture \zede réduite d'un anneau commu\-tatif}
\TTref{thEntRel1}{\Tho fondamental des \entrelsz}
\TTref{thDualiteFinie}{\Tho de dualité entre \trdis finis et ensembles ordonnés finis}
%\TTref{}{}
%\TTref{}{}

%: 12
\CHAP{\nameref{ChapAdpc}}\label{tdtChapAdpc}

\TTref{thAnar}{Caractérisations des \anarsz, restes chinois}
\TTref{thiivanar}{Structure multiplicative des \ids \ivs dans un \anarz}
\TTref{thPruf}{Caractérisations des \adpsz}
\TTref{thExtEntPruf}{Extension enti\`ere normale d'un \adpz}
\TTref{thSurAdp}{Suranneau d'un \adpz}
\TTref{th.adpcoh}{Caractérisations des \adpcsz }
\TTref{ThImMat}{Modules \pf sur un \adpcz}
\TTref{th.2adpcoh}{Une autre \carn des \adpcsz}
\TTref{factAdpIntExt}{Extension finie d'un \adpc (avec le \tho \ref{thAESTE})}
\TTref{thK1-SLE}{$\SL_3=\EE_3$ pour un anneau \qi de dimension $\leq 1$}
\TTref{th1-5}{\Tho un et demi: anneaux \qis de dimension~$\leq1$}
\TTref{dekinbe}{Un \adpc de Bézout}
\TTref{thcohdim1}{Un anneau normal, \cohz, de dimension $\leq1$ est un \adpz}
\TTref{thPTFDed}{Modules \pro de rang $k$ sur un \ddp de dimension~$\leq1$}
\TTref{thMpfPruCohDim}{\Tho des facteurs invariants:  \mpfs sur un \ddp de dimen\-sion~$\leq1$}
\TTref{mat33}{Réduction d'une matrice ligne }
\TTref{thRieszAnar}{\Tho de Riesz pour les \anarsz}
\TTref{thFactor}{Factorisation d'\itfs sur un \adpc de dimension~$\leq1$
(voir aussi le \tho \ref{thFactor2})}
\TTref{prop2DDK}{Un \adk est \`a \fapz}
\TTref{propDDK}{Caractérisations des \adksz }
\TTref{thdefDDKTOT}{Anneaux de Dedekind \`a \fcn totale}
\TTref{lemthAESTE}{Un calcul de \cli (\adkz)}
%: sinotenglish
\sinotenglish{\TTref{thAnormalAXnormal}
{Si $\gA$ est un \anor il en va de m\^eme pour $\AX$}}
%

%: 13
\CHAP{\nameref{chapKrulldim}}

\TTref{thZarZar}{La dualité entre spectre et treillis de Zariski}
\TTref{thDKA}{Caractérisation \elr de la \ddkz}
\TTref{th1.5}{\Tho un et demi (un autre)}
\TTref{propDimKXY}{Dimension de Krull d'un anneau de \pols sur un corps}
\TTref{thDKAG}{Dimension de Krull et mise en position de \Noez}
\TTref{thAmin}{Cl\^oture \qi minimale d'un anneau}
\TTref{thKdimMor}{Dimension de Krull d'un morphisme}
\TTref{corthKdimMor}{Dimension de Krull d'un anneau de \polsz}
\TTref{cor2thKdimMor}{Dimension de Krull d'une extension enti\`ere}
\TTref{thKdimTDTO}{Dimension de Krull d'un ensemble totalement ordonné}
\TTref{th1Valdim}{Dimension de Krull d'une extension d'\advsz}
\TTref{thValDim}{Dimension valuative d'un anneau de \polsz}
\TTref{corthValDim}{Dimension de Krull et \anarsz}
\TTref{propGupLY}{Going up, Going down et \ddkz}
%\TTref{}{}

%: 14
\CHAP{\nameref{chapNbGtrs}}

\TTref{thKroH}{\Tho de \KRN (2), non \noez, pour la \ddkz}
\TTref{Bass0}{\Tho \gui{stable range} de Bass, non \noez}
\TTref{thKroLoc}{\Tho de Kronecker, version locale}
\TTref{Bass}{\Tho \gui{stable range} pour la dimension de Heitmann}
\TTref{KroH2}{\Tho de Kronecker, variante Heitmann}
\TTref{thSerre}{\Tho \SSO pour la $\Sdim$}
\TTref{thSwan}{\Tho de Forster-Swan pour la $\Gdim$}
\TTref{thSwan2}{\Tho de Forster-Swan \gnlz, pour la $\Gdim$}
\TTref{thBassCancel2}{\Tho de simplification de Bass, pour la $\Gdim$}
\TTref{th2KroH}{\Tho de Kronecker, pour les \sutsz}
\TTref{thPartitionSpec}{Partition constructible du spectre de Zariski 
et $k$-stabilité}
\TTref{matrix}{\Tho de Coquand, 1: Forster-Swan et autres
avec la $n$-stabilité}
\vspace{3pt}
\TTref{matrixC}{\Tho de Coquand, 2: manipulations \elrs de colonnes
avec la $n$-sta\-bi\-lité}
\vspace{2pt}
%\TTref{t0basic}{Vecteurs \umds et $n$-stabilité}
\TTref{MAINCOR}{\Tho de Coquand, 3: Forster-Swan et autres
avec la dimension de Heitmann}
%\TTref{}{}

%: 15
\CHAP{\nameref{chapPlg}}

Machineries dynamiques et  \plgs  variés sont indiqués pages \pageref{machineries}
et \pageref{plgvarios}.
%\TTref{lemAssoc}{Lemme des localisations successives, 2}
%\TTref{lemRecouvre}{Lemme des localisations successives, 3}
%\TTpa{MethodeIdeps}{Machinerie \lgbe de base (\`a \idepsz)}
%\TTpa{MethodeIdemax}{Machinerie dynamique \`a \idemasz}
%\TTpa{MethodeIdemin}{Machinerie dynamique \`a \idemisz}

%: 16
\CHAP{\nameref{ChapMPEtendus}}

\TTref{thTSC}{\Tho de Traverso-Swan-Coquand}
\TTref{thRoitman}{\Tho de Roitman}
\TTref{thHor0}{\Tho de Horrocks local}
\TTref{thHor}{\Tho de Horrocks global}
\TTref{th2Bass}{\Tho de Bass}
%:HHH modif
\TTref{thStabLibPol}{Induction de Quillen concr\`ete, cas  stablement libre}
\TTref{InductionQabs}{Induction de Quillen abstraite}
\TTref{thQUILIND}{Induction de Quillen concr\`ete}
\TTref{corthQUILIND}{\Tho de Quillen-Suslin, preuve
de Quillen}
%:HHH rajout
\TTref{th2QUILIND}{Induction de Quillen concr\`ete, cas libre}
\TTref{thSuslinQS}{\Tho de Suslin}
\TTref{th4SusQS}{\Tho de Suslin (un autre)}
\TTref{th2HorrocksLocal}{Petit \tho de Horrocks \`a la Vaserstein
 (et \tho \ref{th2HorrocksGlobal})}
\TTref{th1Rao}{\Tho de Rao}
\TTref{thBass.Valuation}{\Tho de Bass (un autre)}
\TTref{thVX2stab}{La $2$-stabilité de $\VX$, pour le \thref{thBass.Valuation}}
\TTref{thBassValu}{\Tho de Bass-Simis-Vasconcelos (et \tho \ref{thBassAri})}
\TTref{compa}{Comparaison dynamique de $\gA(X)$ avec $\ArX$}
%:HHH permuter les deux ci dessous
\TTref{thMaBrCo}{\Tho de Maroscia \& Brewer-Costa}
\TTref{LSindabs}{Induction de Lequain-Simis abstraite}
\TTref{induYengui}{Induction de Yengui}
\TTref{thLSValu}{\Tho de Lequain-Simis}

\newpage \thispagestyle{CMcadreseul}  

%\rdb
\refstepcounter{chapter}

%---     ---------------
\cleardoublepage \thispagestyle{CMcadreseul}

\catcode`\@=11
\@makesindexhead{Index des notations}
\refstepcounter{chapter}
\addtocontents{toc}{\vskip-0.8em}
\addcontentsline{toc}{chapterbis}{Index des notations}
\markboth{Index des notations}{Index des notations}
\perso{compil\'e le \today}
\pagestyle{CMExercicesheadings}

\newbox\toto
\newbox\tata
\newlength\largeurtoto
\newlength\largeurtata

\newcommand \ttt[3]{\setbox\tata=\hbox{#1}%
\ifdim\wd\tata<.15\textwidth\relax%
\setlength{\largeurtoto}{.78\textwidth}%
\setlength{\largeurtata}{.15\textwidth}%
\else%
\setlength{\largeurtoto}{.93\textwidth}%
\addtolength{\largeurtoto}{-\wd\tata}%
\setlength{\largeurtata}{\wd\tata}%
\fi%
\setbox\toto=\hbox{\parbox[b]{\largeurtoto}{\leftskip10pt\parindent-\leftskip\strut#2\dotfill\par}}%
\smallskip\noindent\mbox{%
\parbox[b][\ht\toto][t]{\largeurtata}{#1}%
\parbox[b]{\largeurtoto}{\leftskip10pt\parindent-\leftskip\strut#2\dotfill\par}%
\hspace{.02\textwidth}%
\parbox[b]{.05\textwidth}{\hfill#3}%
}%
\par}

\newcommand
\NOTA[3]{\ttt{$#1\phantom{_{S}^{1}}$}{#2}{\pageref{#3}}

\vspace{-3pt}}

\begin{flushright}{page}\end{flushright}%

%%%%%%%%%%%%%%%%%%%%%%%%%%%%%%%%%%%%%%%%%%%%%%%%%%%%%%%%%%%%%%%%%%%%%%%%%%%
%:   Exemples
%:HHH nouveaux
\CHAP{\nameref{chapMotivation}}
%%%%%%%%%%%%%%%%%%%%%%%%%%%%%%%%%%%%%%%%%%%%%%%%%%%%%%%%%%%%%%%%%%%%%%%%%%%
\NOTA{\Der{\RR}{\gB}{M}}{le \Bmo des \dvns de $\gB$ dans $M$}{Leibniz2}
\NOTA{\mathrm{Der}(\gB)}{le \Bmo des \dvns de $\gB$}{Leibniz2}
\NOTA{\Om{\RR}{\gB}}{le \Bmo des  \diles (de K\"ahler) de~$\gB$, voir aussi \paref{thDerivUniv}}{dilesK}
%\NOTA{ }{ }{NOTAgAfa}

%%%%%%%%%%%%%%%%%%%%%%%%%%%%%%%%%%%%%%%%%%%%%%%%%%%%%%%%%%%%%%%%%%%%%%%%%%%
%:   Syst\`emes d'équations linéaires
\CHAP{\nameref{chapSli}}
%%%%%%%%%%%%%%%%%%%%%%%%%%%%%%%%%%%%%%%%%%%%%%%%%%%%%%%%%%%%%%%%%%%%%%%%%%%
\NOTA{\pi_{\gA,\fa}}{l'\homo canonique $\gA\to\gA\sur{\fa}$}{NOTAgAfa}
\NOTA{\Ati}{le groupe multiplicatif des \elts inversibles de $\gA$}{NOTAAst}
\NOTA{\gA_{S}}{(ou encore $S^{-1}\gA$) le localisé de $\gA$ en $S$}{NOTAAst}
\NOTA{\sat{S}}{le saturé du \mo $S$}{NotaSatmon}
\NOTA{j_{\gA,S}}{l'\homo canonique $\gA\to\gA_{S}$}{NOTAAst}
\NOTA{\gA[1/s]}{(ou encore $\gA_{s}$) le localisé de $\gA$ en $s^\NN$}{NOTAA[1/s]}
\NOTA{(\fb:\fa)_\gA}{le transporteur de l'idéal $\fa$ dans l'idéal $\fb$}%
{NOTATransp}
\NOTA{(P:N)_\gA}{le transporteur du module $N$ dans le module $P$}%
{NOTATransp}
\NOTA{\Ann_\gA(x)}{l'annulateur de l'\elt $x$}{NOTAAnn}
\NOTA{\Ann_\gA(M)}{l'annulateur du module $M$}{NOTAAnn}
\NOTA{(N:\fa)_M}{$\sotq{x\in M}{\fa x\subseteq N}$}{NOTAAnn2}
\NOTA{{(N:{\fa^\infty})}_M}{$\sotq{x\in M}% (N:\fa^\infty)_M
{\exists n\,\, \fa^n x\subseteq N}$}{NOTAAnn2}
%:HHH Reg, introduit au debut supprime plus loin
\NOTA{\Reg \gA }{\mo des \elts \ndzs de $\gA$}{NOTATotFrac}
\NOTA{\Frac\gA}{anneau total des fractions de $\gA$}{NOTATotFrac}
\NOTA{\Ae{m\times p}}{(ou $\MM_{m,p}(\gA)$) matrices \`a $m$ lignes et
$p$ colonnes}{NOTAmatrices}
\NOTA{\Mn(\gA)}{$\MM_{n,n}(\gA)$}{NOTAmatrices}
\NOTA{\GLn(\gA)}{groupe des matrices \ivsz}{NOTAmatrices}
\NOTA{\SLn(\gA)}{groupe des matrices de \deter $1$}{NOTAmatrices}
\NOTA{\GAn(\gA)}{\mprnsz}{NOTAmatrices}
\NOTA{\DA(\fa)}{(ou encore $\sqrt{\fa}$) nilradical de l'\id $\fa$ de $\gA$}{NOTADA}
\NOTA{\gA\red}{$\gA\sur{\DA({0})}$: anneau réduit associé \`a $\gA$}{NOTADA}
\NOTA{\rc_{\gA,\uX}(f)}{(ou $\rc(f)$) idéal de $\gA$, contenu du \pol $f\in\AuX$}{NOTADA}
\NOTA{\rg_\gA(M)}{rang d'un module libre, voir aussi les \gnns
aux \mptfs pages \pageref{ModStabLibre}, \pageref{Nota2rang} et \pageref{defiRang}}{Nota1rang}
\NOTA{\Adj\,B\phantom{\wi{B}}}{(ou encore $\wi{B}$) matrice cotransposée de $B$}{NOTACotrans}
\NOTA{\cD_k(G)}{\idd d'ordre $k$ de la matrice $G$}{defIdDet}
\NOTA{\cD_k(\varphi)}{\idd d'ordre $k$ de l'\ali $\varphi$,
 voir aussi \paref{exoIDDPTF1}}{fact.idd prod}
\NOTA{{\rg(\varphi)\geq k}}{notation
qui se comprend avec la \dfn \ref{defRangk},
voir aussi la notation \ref{notaRgfi}}{defRangk}
\NOTA{{\rg(\varphi)\leq k}}{m\^eme chose}{defRangk}
\NOTA{\rE^{(n)}_{i,j}(\lambda)}{(ou $\rE_{i,j}(\lambda)$) matrice \elrz}{NOTAEn}
\NOTA{\En(\gA)}{groupe \elrz}{NOTAEn}
\NOTA{\I_{k}}{matrice identité d'ordre $k$}{lem.min.inv}
\NOTA{0_{k}}{matrice carrée d'ordre $k$}{lem.min.inv}
\NOTA{0_{k,\ell}}{matrice nulle de type $k\times \ell$}{lem.min.inv}
\NOTA{\I_{k,q,m}}{matrice simple standard}{NOTAIkqm}
\NOTA{\I_{k,n}}{\mprn standard}{NOTAIkqm}
\NOTA{A_{\alpha,\beta}}{matrice extraite}{NOTAextraite}
\NOTA{\Adj_{\alpha,\beta}(A)}{voir la notation \ref{notaAdjalbe}}{notaAdjalbe}
\NOTA{\cP_{\ell}}{ensemble des parties finies de
$\{1,\ldots ,\ell\}$}{notaAdjalbe}
\NOTA{\cP_{k,\ell}}{parties \`a $k$ \eltsz}{notaAdjalbe}
\NOTA{\GAnk(\gA)}{sous ensemble de $\GAn(\gA)$: \mprns de rang $k$}{defiGrassmanniennes}
\NOTA{\GGnk(\gA)}{grassmannienne projective sur $\gA$}{defiGrassmanniennes}
\NOTA{\GGn(\gA)}{grassmannienne projective sur $\gA$}{defiGrassmanniennes}
\NOTA{\Pn(\gA)}{espace projectif de dimension $n$ sur $\gA$}{defiGrassmanniennes}
\NOTA{\Diag(a_1,\ldots,a_n)}{%~~~~~~~~~
matrice carrée diagonale}{NOTADiag}
\NOTA{\Tr(\varphi)}{trace de $\varphi$ (\endo de $\Ae n$),
voir aussi \paref{NOTAPolcar}}{NOTA1Polcar}
\NOTA{\rC\varphi(X)}{\polcar de $\varphi$ (idem), voir aussi \paref{NOTAPolcar}}{NOTA1Polcar}
\NOTA{\dex{\gB:\gA}}{$\rg_\gA(\gB)$, voir aussi  \paref{subsecTransRang} et \ref{notaTraceDetCarAlg}}{notaCTrN}
\NOTA{\Tr_{\gB/\!\gA}(a)}{trace de (la multiplication par) $a$, voir aussi \vref{defi2STF}}{notaCTrN}
\NOTA{\rN_{\gB/\!\gA}(a)}{norme de $a$, voir aussi \vref{defi2STF}}{notaCTrN}
\NOTA{\rC{\gB/\!\gA}(a)}{\polcar de (la multiplication par) $a$, voir aussi \vref{defi2STF}}{notaCTrN}
\NOTA{\Gram_\gA(\varphi,\ux)}{%~~~~~
matrice de Gram de  $(\ux)$ pour $\varphi$}{defiGram}
\NOTA{\gram_\gA(\varphi,\ux)}{%~~~~~
\deter de Gram de  $(\ux)$ pour $\varphi$}{defiGram}
\NOTA{\disc_{\gB/\!\gA}(\ux)}{%~~~~~
\discri de la famille
$(\ux)$}{defiDiscTra}
\NOTA{\Disc_{\gB/\!\gA}}{discriminant d'une extension libre}{defiDiscTra}
\NOTA{\Lin_\gA(M,N)}{\Amo d'\alisz}{NOTAAlis}
\NOTA{\End_\gA(M)}{$\Lin_\gA(M,M)$}{NOTAAlis}
\NOTA{M\sta}{module dual de $M$}{NOTAAlis}
\NOTA{\AuX_d}{sous-\Amo de $\AuX$ des \pogs de degré $d$}{NOTAAXd}
%%%%%%%%%%%%%%%%%%%%%%%%%%%%%%%%%%%%%%%%%%%%%%%%%%%%%%%%%%%%%%%%%%%%%%%%%%%
%:   Methode coeffs indetermines
\CHAP{\nameref{chapGenerique}}
%%%%%%%%%%%%%%%%%%%%%%%%%%%%%%%%%%%%%%%%%%%%%%%%%%%%%%%%%%%%%%%%%%%%%%%%%%%
\NOTA{\Pf(E)}{ensemble des parties finies de $E$}{NOTAPfPfe}
\NOTA{\Pfe(E)}{ensemble des parties finiment énumérées de $E$}{NOTAPfPfe}
\NOTA{\Hom_\gA(\gB,\gB')}{ensemble des \homos d'\Algs de $\gB$ vers $\gB'$ }{def0Alg}
\NOTA{\mu_{M,b}}{(ou $\mu_b$)  $y\mapsto by$, $\in\End_\gB(M)$ ($b\in \gB$, $M$ un \Bmoz)}{NOTAmux}
\NOTA{\cJ(f)}{idéal des relateurs symétriques}{definotaAdu}
\NOTA{\Adu_{\gA,f}}{\adu de $f$ sur $\gA$}{definotaAdu}
\NOTA{\disc_{X}(f)}{discriminant du \polu  $f$ de $\AX$}{eqDiscri}
\NOTA{\Tsc_g(f)}{transformé de Tschirnhaus de $f$ par $g$}{defiTschir}
\NOTA{\Mip_{\gK,x}(T)}{ou $\Mip_{x}(T)$, \polmin \mon de $x$ (sur le corps $\gK$)}{defiSTF}
\NOTA{G.x}{orbite de $x$ sous $G$}{NOTAStStp}
\NOTA{{G.x=\so{x_1,\ldots,x_k}}}{%~~~~~~~~~~~~
orbite énumérée sans répétition avec $x_1=x$ }{NOTAStStp}
\NOTA{\St_G(x)}{(ou $\St(x)$) sous-groupe stabilisateur du point $x$}{NOTAStStp}
\NOTA{\Stp_G(F)}{(ou $\Stp(F)$) stabilisateur point par point de la partie $F$}{NOTAStStp}
\NOTA{\idg{G:H}}{indice du sous-groupe $H$ dans le groupe $G$: $\#(G/H)$}{NOTAStStp}
\NOTA{\Fix_E(H)}{(ou encore $E^H$) partie de $E$ formée des points fixes de $H$}{NOTAStStp}
\NOTA{\sigma\in G/H}{on prend un $\sigma$ dans chaque classe \`a gauche modulo $H$}{NOTAStStp}
\NOTA{\rC{G}(x)(T)}{$=\prod_{\sigma\in G}(T-\sigma(x))$}{NOTAStStp}
\NOTA{\rN_G(x)}{$=\prod_{\sigma\in G}\sigma(x)$}{NOTAStStp}
\NOTA{\Tr_G(x)}{$=\som_{\sigma\in G}\sigma(x)$}{NOTAStStp}
\NOTA{\Rv_{G,x}(T)}{résolvante de $x$ (relativement \`a $G$)}{NOTAStStp}
\NOTA{\Aut_\gA(\gB)}{groupe des $\gA$-\autos de $\gB$}{NOTAAutAB}
\NOTA{\Gal(\gL/\gK)}{idem, pour une extension galoisienne}{defiCorGal}
\NOTA{\cG_{\gL/\gK}}{sous-groupes finis de $\Aut_\gK(\gL)$}{defiCorGal}
\NOTA{\cK_{\gL/\gK}}{sous-$\gK$-extensions \stfes de $\gL$}{defiCorGal}
\NOTA{\Gal_\gK(f)}{groupe de Galois du \pol \spl $f$}{NOTAGalKf}
%\NOTA{(G:H)}{indice du groupe $H$ dans le groupe $G$: $\#(G/H)$}{thGaloiselr}
\NOTA{\Syl_X(f,p,g,q)}{%~~~~~~~
matrice de Sylvester de $f$ et $g$
en degrés $p$ et $q$}{secRes}
\NOTA{\Res_X(f,p,g,q)}{%~~~~~~~
résultant des \pols $f$ et $g$
en degrés $p$ et $q$}{secRes}
\NOTA{\car(\gK)}{\cara d'un corps}{NOTACarK}
\NOTA{\Adj\iBA (x)}{ou $\wi x$: \elt cotransposé, voir aussi \paref{eqeltcotransp}}{eqelt0cotransp}
\NOTA{(\gA:\gB)}{conducteur de $\gA$ dans $\gB$}{defiConducteur}
\NOTA{\fR_X(f,g_1,\ldots,g_r)}{}{lemElimPlusieurs}
\NOTA{\JJ_\uX (\uf)}{matrice jacobienne d'un \sypz}{secNewton}
\NOTA{\J_{\uX}(\uf)}{jacobien d'un \sypz}{secNewton}
\NOTA{\idg{L:E}_\gA}{indice d'un sous-module \tf dans un module libre}{exolemSousLibre}

%%%%%%%%%%%%%%%%%%%%%%%%%%%%%%%%%%%%%%%%%%%%%%%%%%%%%%%%%%%%%%%%%%%%%%%%%%%
%:   Modules \pf
\CHAP{\nameref{chap mpf}}
%%%%%%%%%%%%%%%%%%%%%%%%%%%%%%%%%%%%%%%%%%%%%%%%%%%%%%%%%%%%%%%%%%%%%%%%%%%
\NOTA{R_\ua}{matrice des relateurs triviaux}{secRelTrivSeqReg}
\NOTA{\scp{\ux}{\uz}}{$\sum_{i=1}^n x_iz_i$}{NOTAProdScal}
\NOTA{\fm_\uxi}{$\gen{x_1-\xi_1,\ldots,x_n-\xi_n}_\gA$: idéal du zéro $\uxi$}
{secExempleGeo}
\NOTA{M\otimes_\gA N}{produit tensoriel de deux \Amosz}{propPftens}
\NOTA{\Al k_\gA M}{puissance extérieure $k$-i\`eme de $M$}{propPfPex}
\NOTA{\gS_\Ae k M}{puissance symétrique $k$-i\`eme de $M$}{propPfPex}
\NOTA{\rho\ist(M)}{\Bmo obtenu \`a partir du  \Amo $M$  par l'\eds $\rho:\gA\rightarrow \gB$}{NOTArhosta}
\NOTA{\cF_n(M)}{ou $\cF_{\gA,n}(M)$:  $n$-i\`eme \idf du \Amo \tf $M$}{def ide fit}
\NOTA{\fRes_{X}(\ff)}{\id résultant de $\ff$
(avec un \polu dans $\ff$)}{thElimAff}
\NOTA{\cK_n(M)}{$n$-i\`eme \id de Kaplanski du \Amo  $M$}{exoAutresIdF}
%\NOTA{}{}{}

%%%%%%%%%%%%%%%%%%%%%%%%%%%%%%%%%%%%%%%%%%%%%%%%%%%%%%%%%%%%%%%%%%%%%%%%%%%
%:   Modules \ptfsz, 1
\CHAP{\nameref{chap ptf0}}
%%%%%%%%%%%%%%%%%%%%%%%%%%%%%%%%%%%%%%%%%%%%%%%%%%%%%%%%%%%%%%%%%%%%%%%%%%%
\NOTA{\theta_{M,N}}{\Ali naturelle  $M\sta\te_\gA N\to \Lin_\gA(M,N)$}{NOTAthetaMN}
\NOTA{\theta_{M}}{\Ali naturelle  $M\sta\te_\gA M\to \End_\gA(M)$}{NOTAthetaM}
\NOTA{\Diag(M_1,\ldots,M_n)}{%
matrice carrée diagonale par blocs}{Notadiagblocs}
\NOTA{\Bdim\gA<n}{stable range (de Bass) inférieur ou égal \`a $n$}{defiStableRange}
\NOTA{\det\varphi}{\deter de l'\endo $\varphi$ d'un \mptfz }{NOTAdet}
\NOTA{\rC\varphi(X)}{\polcar de $\varphi$ \ldots\, (idem) }{NOTAPolcar}
\NOTA{\wi{\varphi}}{\endo cotransposé de $\varphi$ \ldots\, (idem) }{NOTACotransp}
\NOTA{\rF{\varphi}(X)}{\polfon de  $\varphi$, i.e.,
$\det(\Id_{P}+X\varphi)$}{NOTAPolfon}
\NOTA{\Tr_P(\varphi)}{trace de l'\endo $\varphi$}{NOTAPolfon}
\NOTA{\rR{P}(X)}{\polmu du \mptf $P$}{NOTAPolmu}
\NOTA{\ide_{h}(P)}{l'\idm  associé \`a l'entier $h$ et au module
\pro $P$}{NOTAide}
\NOTA{P\ep{h}}{composant du module $P$ en rang~$h$}{NOTAep}
%\NOTA{}{}{}

%%%%%%%%%%%%%%%%%%%%%%%%%%%%%%%%%%%%%%%%%%%%%%%%%%%%%%%%%%%%%%%%%%%%%%%%%%%
%:   Alg\`ebres \STFES
\CHAP{\nameref{chap AlgStricFi}}
%%%%%%%%%%%%%%%%%%%%%%%%%%%%%%%%%%%%%%%%%%%%%%%%%%%%%%%%%%%%%%%%%%%%%%%%%%%
\NOTA{\rC{\gB/\!\gA}(x)(T)}{\polcar de (la multiplication par) $x$}{defi2STF}
\NOTA{\rF{\gB/\!\gA}(x)(T)}{\polfon de (la multiplication par) $x$}{defi2STF}
\NOTA{\rN\iBA(x)}{norme de $x$: \deter de la multiplication par $x$}{defi2STF}
\NOTA{\Tr\iBA(x)}{trace de (la multiplication par) $x$}{defi2STF}
\NOTA{a\centerdot\alpha}{$\alpha\circ \mu_a:x\mapsto \alpha(ax)$}{definotaAsta}
\NOTA{\Adj\iBA(x)}{ou $\wi x$: \elt cotransposé}{eqeltcotransp}
\NOTA{[\gB:\gA]}{$\rg_\gA(\gB)$, voir aussi pages \pageref{notaCTrN} et \pageref{notaTraceDetCarAlg}}{subsecTransRang}
\NOTA{\Phi_{\gA/\gk,\lambda}}{$\Phi_\lambda(x,y) = \lambda(xy)$}{defdualisante}
\NOTA{\phi\otimes\phi'}{produit tensoriel de formes bi\linsz}{NOTAptfb}
\NOTA{\env\gk\gA}{$\gA\otimes_\gk\gA$, \alg enveloppante de $\gA/\gk$}{NotaAGenv}
\NOTA{\rJ\iAk}{idéal de $\env\gk\gA$}{eqDiAk}
\NOTA{\Delta\iAk}{$\Delta(x)=x\te1-1\te x$}{eqDiAk}
\NOTA{\mu\iAk}{$\mu\iAk \left(\sum_ia_i\te b_i\right) =\sum_ia_ib_i$}{eqmuAk}
\NOTA{\Der{\gk}{\gA}{M}}{le \Amo des \dvns de $\gA$ dans $M$}{defiDeriv}
\NOTA{\mathrm{Der}(\gA)}{le \Amo des \dvns de $\gA$}{defiDeriv}
\NOTA{\Om{\gk}{\gA}}{le \Amo des  \diles (de K\"ahler) de~$\gA$}{thDerivUniv}
\NOTA{\vep\iAk}{\idm qui engendre $\Ann(\rJ\iAk)$, s'il existe}{thSepIversen}
\NOTA{\LIN_\gk\!(\gA,\gA)}{\Amo des \klis de $\gA$ dans $\gA$}{NotaAGAL}
\NOTA{\PGL_n(\gA)}{groupe quotient $\GLn(\gA)/\Ati$}{NOTAPGL}
\NOTA{\rA_n}{sous-groupe des permutations paires de~$\Sn$}{NOTAAn}
%%%%%%%%%%%%%%%%%%%%%%%%%%%%%%%%%%%%%%%%%%%%%%%%%%%%%%%%%%%%%%%%%%%%%%%%%%%
%:   Théorie de Galois
\CHAP{\nameref{ChapGalois}}
%%%%%%%%%%%%%%%%%%%%%%%%%%%%%%%%%%%%%%%%%%%%%%%%%%%%%%%%%%%%%%%%%%%%%%%%%%%
\NOTA{\BB(\gA)}{\agB des \idms de $\gA$}{propB(A)}
\NOTA{\cB(f)}{base \gui{canonique} de l'\aduz}{notaBaseadu}
%\NOTA{}{}{}
%%

%%%%%%%%%%%%%%%%%%%%%%%%%%%%%%%%%%%%%%%%%%%%%%%%%%%%%%%%%%%%%%%%%%%%%%%%%%%
%:   Modules plats
%\CHAP{\nameref{chap mod plats}}
%%%%%%%%%%%%%%%%%%%%%%%%%%%%%%%%%%%%%%%%%%%%%%%%%%%%%%%%%%%%%%%%%%%%%%%%%%%
%\NOTA{}{}{}
%\NOTA{}{}{}

%%%%%%%%%%%%%%%%%%%%%%%%%%%%%%%%%%%%%%%%%%%%%%%%%%%%%%%%%%%%%%%%%%%%%%%%%%%
%:   Anneaux locaux
\CHAP {\nameref{chap Anneaux locaux}}
\NOTA{\Rad(\gA)}{radical de Jacobson  de $\gA$}{eqDefRadJac}
%\NOTA{\J_\uX (\uf)}{jacobien d'un \syp}{subsecZedLocPtIsoleSimple}
\NOTA{\gA(X)}{localisé de Nagata de $\gA[X]$}{factLocNagata}
%:h2013
\NOTA{\Suslin(\bn)}{ensemble de Suslin de $(\bn)$}{defiSuslinSet}
\NOTA{\gk[G]}{\alg d'un groupe, ou d'un \moz}{exoAlgMon}

%\NOTA{}{}{}
%\NOTA{}{}{}

%%%%%%%%%%%%%%%%%%%%%%%%%%%%%%%%%%%%%%%%%%%%%%%%%%%%%%%%%%%%%%%%%%%%%%%%%%%
%:   Modules \ptfsz, 2
\CHAP{\nameref{chap ptf1}}
%%%%%%%%%%%%%%%%%%%%%%%%%%%%%%%%%%%%%%%%%%%%%%%%%%%%%%%%%%%%%%%%%%%%%%%%%%%
\NOTA{\Gn}{$\Gn =\ZZ[(f_{i,j})_{i,j\in \lrbn}]/\cGn$}{NOTAGN}
\NOTA{\cGn}{relations obtenues en écrivant $F^2=F$}{NOTAGN}
\NOTA{\HOp(\gA)}{semi-anneau des rangs des \Amos quasi libres}{notaHO+}
\NOTA{[P]_{\HOp (\gA)}}{ou $[P]_\gA$, ou
$[P]$: classe d'un \Amo quasi libre dans $\HOp(\gA)$}{notaHO+}
\NOTA{\rg_\gA(M)}{rang (\gnez) du \Amo \ptf $M$}{NOTAHO+}
\NOTA{\HO\gA}{anneau des rangs sur $\gA$}{notaHO}
\NOTA{\dex{\gB:\gA}}{$\rg_\gA(\gB)$, voir aussi pages \pageref{notaCTrN} et \pageref{subsecTransRang}}{notaTraceDetCarAlg}
\NOTA{\Gn(\gA)}{$\Gn\otimes_\ZZ\gA$}{subsubsec AGBR}
\NOTA{\cGnk}{$\cGn+\gen{1-r_k}$, avec  (dans $\Gn$) $r_k=\ide_k(\Im\,F)$}{subsubsec AGBR}
\NOTA{\Gnk}{$\Gnk =\ZZ[(f_{i,j})_{i,j\in \lrbn}]/\cGnk$
ou encore $\Gn[1/r_k]$}{subsubsec AGBR}
\NOTA{\GAnk(\gA)}{\gui{sous-\vrtz} de $\GAn(\gA)$: \prrs de rang $k$}{subsubsec AGBR}
\NOTA{\GKO\gA}{semi-anneau des classes d'\iso de \mptfs sur $\gA$}{NOTAGKO}
\NOTA{\Pic\gA}{groupe des classes d'\iso  des \mrcs 1 sur $\gA$}{NOTAPic}
\NOTA{\KO\gA}{anneau de Grothendieck de $\gA$}{NOTAK0}
\NOTA{[P]_{\KO (\gA)}}{ou $[P]_\gA$, ou
$[P]$: classe d'un \Amo \ptf dans $\KO (\gA)$}{NOTAK0}
\NOTA{\KTO\gA}{noyau de l'\homo rang $\rg:\KO \gA\to\HO \gA$}{NOTAKTO}
\NOTA{\Ifr\gA}{mono\"{\i}de des \ifrs \tf de l'anneau $\gA$}{NOTAIfr}
\NOTA{\Gfr\gA}{groupe des \elts \ivs de $\Ifr\gA$}{NOTAIfr}
\NOTA{\Cl \gA}{groupe des classes d'\ids \ivs (quotient de $\Gfr\gA$
par le sous-groupe des \idps \ivsz)}{NOTAIfr}

%%%%%%%%%%%%%%%%%%%%%%%%%%%%%%%%%%%%%%%%%%%%%%%%%%%%%%%%%%%%%%%%%%%%%%%%%%%
%:   Treillis distributifs
\CHAP{\nameref{chapTrdi}}
%%%%%%%%%%%%%%%%%%%%%%%%%%%%%%%%%%%%%%%%%%%%%%%%%%%%%%%%%%%%%%%%%%%%%%%%%%%
\NOTA{\dar a}{$\sotq{x\in X}{x\leq a}$, voir aussi \vpageref{NOTAdara}}{eqda}
\NOTA{\uar a}{$\sotq{x\in X}{x\geq a}$, voir aussi \vpageref{NOTAuara}}{eqda}
\NOTA{\gT\eci}{treillis opposé du treillis $\gT$}{NOTAtrdiopp}
\NOTA{\cI_\gT(J)}{idéal engendré par $J$ dans le \trdi $\gT$}{NOTAidtrdi}
\NOTA{\cF_\gT(S)}{filtre engendré par $S$ dans le \trdi $\gT$}{NOTAfitrdi}
\NOTA{{\gT/(J=0,U=1)}}{treillis quotient particulier}{propIdealFiltre}
\NOTA{\Bo(\gT)}{\agB engendrée par le \trdi $\gT$}{thBoolGen}
\NOTA{\ZZ^{(P)}}{somme directe \orte de copies de $\ZZ$, indexée par $P$}{NotaZZP}
\NOTA{\boxplus_{i\in I}G_i}{somme directe \orte de groupes ordonnés}{NotaSDirOr}
\NOTA{\cC(a)}{Sous-groupe solide engendré par $a$ (dans un \grlz)}{defiCongru}
\NOTA{\DA(x_1,\ldots ,x_n)}{$\DA(\gen{x_1,\ldots ,x_n})$: un \elt de
$\ZarA$}{notaZA}
\NOTA{\ZarA}{treillis de Zariski de $\gA$}{notaZA}
\NOTA{\gA_{S}\sur{\fa}}{(ou encore $S^{-1}\gA\sur{\fa}$) on inverse les \elts de $S$
et on annule les \elts de $\fa$}{defiASa}
\NOTA{\satu S \gA}{ou $\sat{S}$: le filtre obtenu en saturant le \mo $S$ dans $\gA$}{NOTASatu}
\NOTA{\Abul}{anneau \zed réduit engendré par $\gA$}{thZedGen}
\NOTA{A \vda B}{$\Vi A\;\leq \;\Vu B$: relation implicative}{notaVupVda}
\NOTA{\SpecT}{spectre du \trdi fini $\gT$, voir aussi \paref{SpecTrdi}}{SpecTrdiFi}
\NOTA{(b:a)_\gT}{le transporteur de $a$ dans $b$ (\trdisz)}%
{NOTATransp2}
\NOTA{\gA_\mathrm{qi}}{cl\^oture \qi de $\gA$}{exoQiClot}
\NOTA{\Min\gA}{sous-espace de $\SpecA$ formé par les \idemisz}{exoMinA}
%\NOTA{}{}{}

%%%%%%%%%%%%%%%%%%%%%%%%%%%%%%%%%%%%%%%%%%%%%%%%%%%%%%%%%%%%%%%%%%%%%%%%%%%
%:   Anneaux de Pr\"ufer et de Dedekind
\CHAP{\nameref{ChapAdpc}}
%%%%%%%%%%%%%%%%%%%%%%%%%%%%%%%%%%%%%%%%%%%%%%%%%%%%%%%%%%%%%%%%%%%%%%%%%%%

\NOTA{\fa\div\fb}{$\sotq{x\in\Frac\gA}{x\fb\subseteq\fa}$}{NOTAfadivfb}
\NOTA{\gA[\fa t]}{\alg de Rees de l'\id $\fa$ de $\gA$}{NOTARees}
\NOTA{\Icl_\gA(\fa)}{\cli de l'\id $\fa$ dans $\gA$}{fact2Entiers}

%%%%%%%%%%%%%%%%%%%%%%%%%%%%%%%%%%%%%%%%%%%%%%%%%%%%%%%%%%%%%%%%%%%%%%%%%%%
%:   Dimension de Krull
\CHAP{\nameref{chapKrulldim}}
%%%%%%%%%%%%%%%%%%%%%%%%%%%%%%%%%%%%%%%%%%%%%%%%%%%%%%%%%%%%%%%%%%%%%%%%%%%
\NOTA{\SpecA}{spectre de Zariski de l'anneau $\gA$}{nota Spec(A)}
\NOTA{\fD_\gA(\xn)}{%~~~~~~~~
ouvert quasi compact de $\SpecA$}{nota Spec(A)}
\NOTA{\SpecT}{spectre du \trdi $\gT$}{SpecTrdi}
\NOTA{\fD_\gT(u)}{ouvert quasi compact de $\SpecT$}{SpecTrdi}
\NOTA{\OQC(\gT)}{\trdi des ouverts quasi compacts de $\SpecT$}{SpecTrdi}
\NOTA{\JK_\gA(x)}{$\gen{x}+(\DA(0):x)$: idéal bord de Krull de $x$ dans $\gA$}{defZar2}
\NOTA{\JK_\gA(\fa)}{$\fa+(\DA(0):\fa)$: idéal bord de Krull de $\fa$ dans $\gA$}{defZar2}
\NOTA{\gA_\rK^x}{$\gA\sur{\JK_\gA(x)}$: (anneau) bord supérieur de $x$ dans $\gA$}{defZar2}
\NOTA{\SK_\gA(x)}{$x^\NN(1+x\gA)$: \mo bord de Krull de $x$ dans $\gA$}{defZar2}
\NOTA{\gA^\rK_{x}}{$(\SK_\gA(x))^{-1}\gA$: (anneau) bord inférieur de $x$ dans $\gA$}{defZar2}
%\NOTA{\Kdim}{dimension de Krull, voir aussi la \dfn \rref{defiDDKTRDI})}{NOTAKdim}
\NOTA{\Kdim\gA\leq r}{la dimension de Krull de l'anneau $\gA$ est  $\leq r$}{NOTAKdim}
\NOTA{\Kdim\gA\leq\Kdim\gB}{}{notaKdiminf}
\NOTA{\SK_\gA(\xzk)}{\mo bord de Krull itéré}{notaBordsIteres}
\NOTA{\JK_\gA(\xzk)}{\id bord de Krull itéré}{notaBordsIteres}
\NOTA{\IK_\gA(\xzk)}{\id bord de Krull itéré, variante}{notaBordsIteres}
%:HHH rajout
\NOTA{\Kdim\gT\leq r}{la \ddk du \trdi $\gT$ est $\leq r$}{defiDDKTRDI}
\NOTA{\JK_\gT(x)}{$\dar x \,\vu\, (0:x)_\gT$: idéal bord de Krull de $x$ dans 
le \trdiz~$\gT$}{defiBordKrullTrdi}
\NOTA{\gT_\rK^x}{$\gT\sur{\JK_\gT(x)}$: (treillis) bord supérieur de $x$}{defiBordKrullTrdi}
\NOTA{\JK_\gT(\xzk)}{\id bord de Krull itéré dans un \trdiz}{eq1IdBordKrullItereTrdi}
%:HHH rajout
\NOTA{\Kdim \rho}{\ddk du morphisme $\rho$}{defiKdimMor}
\NOTA{\gA_{\so{a}}}{$\gA\sur{a\epr}\times \gA\sur{({a\epr})\epr}$}{lem20MorRc}
\NOTA{\Amin}{cl\^oture \qi minimale de $\gA$}{thAmin}
\NOTA{\Vdim\gA}{dimension valuative}{defiValdim}
%\NOTA{}{}{}
%%%%%%%%%%%%%%%%%%%%%%%%%%%%%%%%%%%%%%%%%%%%%%%%%%%%%%%%%%%%%%%%%%%%%%%%%%%
%:   Nombre de \gtrs d'un module
\CHAP{\nameref{chapNbGtrs}}
%%%%%%%%%%%%%%%%%%%%%%%%%%%%%%%%%%%%%%%%%%%%%%%%%%%%%%%%%%%%%%%%%%%%%%%%%%%
\NOTA{\JA(\fa)}{radical de Jacobson de l'\id $\fa$  de $\gA$}{notaJA}
\NOTA{\JA(x_1,\ldots ,x_n)}{%~~~~~~~
$\JA(\gen{x_1,\ldots ,x_n})$: un \elt de
$\HeA$}{notaJA}
\NOTA{\HeA}{treillis de Heitmann de $\gA$}{defHeit}
\NOTA{\Jdim}{dimension du J-spectre de Heitmann}{NOTAJdim}
\NOTA{\Max\gA}{sous espace de $\SpecA$ formé par les \idemasz}{NOTAJdim}
\NOTA{\Jspec\gA}{$\Spec(\HeA)$: J-spectre de Heitmann}{NOTAJdim}
\NOTA{\JH_\gA(x)}{$\gen{x}+(\JA(0):x)$: idéal bord de Heitmann (de $x$ dans $\gA$)}{defHei2}
\NOTA{\gA_\rH^x}{$\gA\sur{\JH_\gA(x)}$: l'anneau bord de Heitmann de $x$}{defHei2}
\NOTA{\Hdim}{dimension de Heitmann}{defDHA}
\NOTA{\Sdim\gA<n}{}{defiSdimGdim}
\NOTA{\Gdim\gA<n}{}{defiSdimGdim}
\NOTA{\Cdim\gA<n}{l'anneau $\gA$ est $n$-stable}{defiAnneaunStable}

%%%%%%%%%%%%%%%%%%%%%%%%%%%%%%%%%%%%%%%%%%%%%%%%%%%%%%%%%%%%%%%%%%%%%%%%%%%
%:   Le principe local-global
\CHAP{\nameref{chapPlg}}
%%%%%%%%%%%%%%%%%%%%%%%%%%%%%%%%%%%%%%%%%%%%%%%%%%%%%%%%%%%%%%%%%%%%%%%%%%%
\NOTA{\cM(U)}{le \mo engendré par l'\elt ou la partie $U$ de $\gA$}{nota mopf}
\NOTA{\cS(I,U)}{$\sotq{v\in\gA}{\exists u\in\cM(U)\;\exists a\in \gen{I}_\gA,
\;  v=u+a }$}{nota mopf}
\NOTA{\cS(a_1,\ldots,a_k;u_1,\ldots,u_\ell)}
{$\cS(\so{a_1,\ldots,a_k},\so{u_1,\ldots,u_\ell})$}{nota mopf}
%\NOTA{}{}{}

%%%%%%%%%%%%%%%%%%%%%%%%%%%%%%%%%%%%%%%%%%%%%%%%%%%%%%%%%%%%%%%%%%%%%%%%%%%
%:   Modules étendus
\CHAP{\nameref{ChapMPEtendus}}
%%%%%%%%%%%%%%%%%%%%%%%%%%%%%%%%%%%%%%%%%%%%%%%%%%%%%%%%%%%%%%%%%%%%%%%%%%%
\NOTA{\gA\lra{X}}{localisé de $\gA[X]$ en les \pols unitaires}{NOTAAlraX}
\NOTA{A\sims{\cG }B}{il existe une
matrice $H\in\cG $ telle que $HA=B$}{NOTAfGg}

%%%%%%%%%%%%%%%%%%%%%%%%%%%%%%%%%%%%%%%%%%%%%%%%%%%%%%%%%%%%%%%%%%%%%%%%%%%
%:   \Tho de stabilité de Suslin
\CHAP{\nameref{ChapSuslinStab}}
\NOTA{\GL(P)}{groupe des \autos\lins de $P$}{NOTAtransvec}
\NOTA{\wi{\EE}(P)}{sous-groupe  de $\GL(P)$ engendré par les transvections}{NOTAtransvec}
\NOTA{\meck {a}{b}}{symbole de Mennicke}{lemMennicke1}
\NOTA{\GLn(\gB,\fb)}{noyau de $\GLn(\gB)\to\GLn(\gB\sur{\fb})$}{notaGLIdeal}
\NOTA{\En(\gB, \fb)}{sous-groupe normal de $\En(\gB)$ engendré par les $\rE_{ij}(b)$ avec
 $b \in \fb$.}{notaGLIdeal}
%\NOTA{}{}{}
%%%%%%%%%%%%%%%%%%%%%%%%%%%%%%%%%%%%%%%%%%%%%%%%%%%%%%%%%%%%%%%%%%%%%%%%%%%
%:   Annexe
\CHAP {Annexe: logique \covz}
\NOTA{\rP(X)}{la classe des parties de $X$}{P(X)}

\newpage
\thispagestyle{CMcadreseul}
\cleardoublepage 
\rdb

%: %@@@@@@@@@------ index --------------------
% \chapter*{Index des termes}

\addtocontents{toc}{\vskip-0.8em}
\addcontentsline{toc}{chapterbis}{Index des termes}
\markboth{Index des termes}{Index des termes}

\printindex


\begin{thebibliography}{155}



\bibitem[Abdeljaoued \& Lombardi]{Jou} {\sc Abdeljaoued A., Lombardi H.}
{\em Méthodes Matricielles.
Introduction \`a la Complexité Algébrique}.
Springer, (2003). \perso{compilé le \today}

\bibitem[Aczel \& Rathjen]{AcRa} {\sc Aczel P.,  Rathjen M.}
{\em Notes on Constructive Set Theory}.
\url{http://www1.maths.leeds.ac.uk/~rathjen/book.pdf}.

\bibitem[Adams \& Loustaunau]{Lou} {\sc Adams W., Loustaunau P.}
{\em An Introduction to Gr\"obner Bases}.
American Mathematical
Society, (1994).


\bibitem[Apéry  \& Jouanolou]{AJ} {\sc Apéry F., Jouanolou J.-P.}
 {\em \'Elimination. Le cas d'une variable.}
 {Hermann, (2006).}

\bibitem[Atiyah \& Macdonald]{Atiyah} {\sc Atiyah  M.F., Macdonald
I.G.}
 {\em Introduction to Commutative Algebra.}
 {Addison Wesley, (1969).}


\bibitem[Basu, Pollack \& Roy]{BPR} {\sc Basu S., Pollack R.,  Roy
M.-F.} {\em Algorithms in real algebraic Geometry.} Springer, (2006).


\bibitem[Bass]{BASS} {\sc Bass H.} {\em Algebraic $K$-theory.}
W. A. Benjamin, Inc., New York-Amsterdam, (1968).

\bibitem[Beeson]{Be} {\sc Beeson M.} {\em Foundations of Constructive
Mathematics.}
Springer-Verlag, (1985).

%\bibitem[Ben-Israel \& Greville]{BIG} {\sc Ben-Israel A., Greville T.}
%{\em Generalized Inverses: Theory and Applications}.
%New York: Wiley (1977).

\bibitem[Bhaskara Rao]{Bha} {\sc Bhaskara Rao K.}
{\em The Theory of Generalized Inverses over a Commutative Ring}.
Taylor \& Francis. Londres, (2002).

\bibitem[Bigard, Keimel \& Wolfenstein]{BKW} {\sc Bigard A., Keimel K., Wolfenstein S.} {\em Groupes et anneaux réticulés}. Springer LNM 608,
(1977).

%\bibitem[Bini \& Pan]{BP} {\sc Bini D., Pan V.} {\em Polynomial and matrix
%computations}. Progress in Theoretical Computer Science,
%Birkh\"auser (1994).


\bibitem[Birkhoff]{Birkhoff} {\sc Birkhoff G.}
 {\em Lattice theory.}
{Third edition.
American Mathematical Society Colloquium Publications, Vol. XXV
American Mathematical Society, Providence, R.I., (1967).}



\bibitem[Bishop]{B67} 
{\em Foundations of Constructive Analysis.}
McGraw Hill, (1967). Reprint avec une préface de Michael Beeson.
2012. IshiPress New York and Tokyo.

\bibitem[Bishop \& Bridges]{BB85} {\sc Bishop E., Bridges D.}
 {\em  Constructive Analysis}.
Springer-Verlag, (1985).

%\bibitem[Bochnak, Coste \& Roy]{BCR} {\sc Bochnak J., Coste M., Roy M.-F.}
% {\em Géométrie algébrique réelle}.
% Springer Verlag, (1987).

\bibitem[Bourbaki]{Bou} {\sc Bourbaki. }{\it   Alg\`ebre Commutative.}
Hermann, (1961-2002).

\bibitem[Bridges \& Richman]{BR} {\sc Bridges D., Richman F.}
{\em  Varieties of Constructive Mathematics.}
London Math. Soc. LNS 97. Cambridge University Press, (1987).

\bibitem[Brouwer]{Bro} {\sc Brouwer L.} {\em Brouwer's Cambridge Lectures on
Intuitionism, 1951} (Van Dalen ed.) Cambridge University Press, (1981).


\bibitem[Burris \& Sankappanavar]{BuSa}  {\sc Burris S., Sankappanavar H.}
{\em A  Course  in Universal  Algebra}.
Springer, (1981).

%: C
\bibitem[Cartan \& Eilenberg]{CE} {\sc  Cartan H., Eilenberg S.} 
{\it  Homological algebra}.
Princeton University Press, (1956).

\bibitem[COCOA]{KrRo} {\sc Kreuzer  M., Robbiano  L.}
 {\it  Computational commutative algebra.}
 Springer Verlag, Berlin. Vol. 1 (2000), Vol. 2 (2005)

\bibitem[Cohn]{Cohn} {\sc Cohn P.}
{\it Basic Algebra. Groups, rings and fields}. (2nd edition)
Springer Verlag, (2002).

\bibitem[Cox]{CoxGal} {\sc Cox D.}
{\it Galois theory}. Wiley-Interscience, (2004).


\bibitem[Cox, Little \& O'Shea]{CLS} {\sc Cox D., Little J, O'Shea D.}
{\it Ideals, Varieties, and Algorithms}. (2nd edition)
Springer Verlag UTM, (1998).

\bibitem[CPMPCS]{PS2016}  
{\em Concepts of proof in mathematics, philosophy, and computer science.  (Based on the Humboldt-Kolleg, Bern, Switzerland, September 9--13, 2013).}
{\sc Eds:   {Probst} D., {Schuster} P.} Berlin: De Gruyter, (2016).

\bibitem[CRA]{CRA} {\it Commutative ring theory and applications.} {\sc  Eds: Fontana M., Kabbaj S.-E., Wiegand S.}
Lecture notes in pure and applied mathematics vol 231. M. Dekker, (2002).

\bibitem[Curry]{Curry}
{\sc Curry H. B.} {\em Foundations of mathematical logic}
McGraw-Hill Book Co., Inc., New York-San Francisco,
Calif.-Toronto-London, (1963).

\bibitem[David, Nour \& Raffalli]{DNR} {\sc David R., Nour K., Raffalli C.}
{\it Introduction \`a la logique}. Dunod, (2001).




\bibitem[Demeyer \& Ingraham]{DI}
{\sc Demeyer F., Ingraham E.} {\it  Separable algebras over
commutative rings}. Springer Lecture Notes in Mathematics 181, (1971).

\bibitem[D\'{\i}az, Lombardi \& Quitté]{DiLQ14} {\sc D\'{\i}az-Toca G., Lombardi H., Quitté}
{\em Modules sur les anneaux commutatifs.}
Calvage\&Mounet, (2014).

\bibitem[Dowek1]{Dow} {\sc Dowek G.}
{\it La logique}. Flammarion. Collection Dominos, (1995).

\bibitem[Dowek2]{Dow2} {\sc Dowek G.}
{\it Les métamorphoses du calcul. Une étonnante histoire de mathématiques}. Le Pommier, (2007).


\bibitem[Edwards89]{Ed} {\sc Edwards H.} {\it Divisor Theory}. Boston, MA:
Birkh\"auser, (1989).


\bibitem[Edwards05]{Ed2} {\sc Edwards H.} {\it Essays in Constructive Mathematics}. Springer Verlag, (2005).


\bibitem[Eisenbud]{Eis} {\sc Eisenbud D.}
{\it Commutative Algebra with a view toward
Algebraic Geometry}.
Springer Verlag, (1995).

\bibitem[Ene \& Herzog]{EnHe} {\sc Ene, V., Herzog, J.}
{\it  Gr\"obner bases in commutative algebra}.
Graduate Studies in Mathematics \num130,
American Mathematical Society
(2012).

\bibitem[Elkadi \& Mourrain]{ElkMo} {\sc Elkadi M., Mourrain B.}
{\it  Introduction \`a la résolution des syst\`emes polynomiaux}.
Collection Mathématiques \& Applications, \num59, Springer Verlag, Berlin (2007).


\bibitem[Feferman]{fef} {\sc Feferman S.}
{\it  In the Light of Logic}.
Oxford University Press, (1998).


\bibitem[Frege-G\"odel]{vH} {\sc van Heijenoort J.} (ed.),
{\it  Frome Frege to G\"odel: a source book in mathematical logic}.
Harvard University Press, Cambridge, Massachussets
(1967). (troisi\`eme réimpression en 2002).

\bibitem[Freid \& Jarden]{FJ} {\sc Freid M. D., Jarden M.} {\em Field Arithmetic}.
Springer-Verlag. (1986).


%: G


\bibitem[von zur Gathen \alb\& Gerhard]{vzGaGe} {\sc von zur Gathen J. Gerhard J.}
{\it  Modern computer algebra}.
Cambridge University Press, Cambridge, (2003). 


\bibitem[Gilmer]{Gil} {\sc Gilmer R.}
{\it Multiplicative Ideal Theory.} Queens papers in pure and applied Math,
vol. 90, (1992).

\bibitem[Glaz]{Gla} {\sc Glaz S.} {\it Commutative Coherent Rings}.
Lecture Notes in Math., vol. 1371, Springer Verlag,
Berlin-Heidelberg-New York, second edition, (1990).

\bibitem[Gr\"atzer]{Grae} {\sc Gr\"atzer G.} {\it Lattice Theory: foundation.}
Birkh\"auser/Springer Basel AG, Basel, (2011).


\bibitem[Gupta \& Murthy]{GM} {\sc Gupta S., Murthy M.}
 {\it  Suslin's work on linear groups over
polynomial rings and Serre conjecture}. ISI Lecture Notes \num8.
The  Macmillan Company of India Limited, (1980).

\bibitem[HoTT]{HoTT} 
 {\it  Homotopy Type Theory:  Univalent Foundations of Mathematics.}  \url{http://homotopytypetheory.org/} (2014).

\bibitem[Infini]{INS} {\sc Toraldo di Francia G.} (ed.),
{\it  L'infinito nella scienza}.
Istituto della Enciclopedia Italiana, Rome,
(1987).

\bibitem[Ireland \& Rosen]{IR} {\sc Ireland K., Rosen M.}
{\it A classical introduction to modern number theory}.
Graduate Texts in Mathematics, vol. 84, Springer-Verlag,
Berlin-Heidelberg-New York, (1989).

\bibitem[Ischebeck \& Rao]{IRa} {\sc  Ischebeck F., Rao R.}
{\it Ideals and Reality. Projective modules and
number of generators of ideals}.
 Springer Monograph in Mathematics,
Berlin-Heidelberg-New York, (2005).


%\bibitem[Iversen]{Ive} {\sc Iversen B.}
%{\it  Generic local structure in commutative  algebras}.
%Springer Lecture Notes in Mathematics \num310 (1973).
\bibitem[Jaffard]{Jaffard2} 
{\sc Jaffard, P.}
{\em Théorie de la dimension dans les anneaux de polynômes}
Gauthier-Villars, Paris,
(1960).

\bibitem[Jensen, Ledet \& Yui]{Jensen} 
{\sc Jensen C., Ledet A., Yui N.}
{\em Generic Polynomials, Constructive
Aspects of the Inverse Galois Problem}
Cambrigde University Press, 
MSRI Publications 45, (2002).

\bibitem[Johnstone]{Johnstone} {\sc Johnstone P.}  {\it  Stone spaces}.
Cambridges studies in advanced mathematics \num3.
Cambridge University Press, (1982).

\bibitem[Kaplansky]{Kapl} {\sc Kaplansky I.}  {\it  Commutative rings}.
Boston, Allyn and Bacon, (1970).

\bibitem[Kleene \& Vesley]{KV}
{\sc Kleene  S.C., Vesley R.}  {\it  The Foundations of intuitionistic mathematics}. Amsterdam (North-Holland), (1965).

 
\bibitem[Knapp, 1]{Knap1} {\sc Knapp A.} {\it  Basic algebra}. Birkh\"auser, (2006).

\bibitem[Knapp, 2]{Knap2} {\sc Knapp A.} {\it  Advanced algebra}. Birkh\"auser, (2007).



\bibitem[Knight]{Kni} {\sc Knight J.}  {\it    Commutative Algebra}.
London Mathematical Society LNS \num 5. Cambridge University
Press, (1971).

		



\bibitem[Kunz]{Kun} {\sc Kunz E.} {\it  Introduction to Commutative
Algebra and Algebraic Geometry}. Birkh\"auser, (1991).

%\bibitem[Kunz-2]{Kun2} {\sc Kunz E.} {\it  K\"ahler Differentials}. 
%Advanced Lectures in Mathematics,
%Friedr. Vieweg \& Sohn, (1986).

%: L
\bibitem[Lafon \& Marot]{Laf} {\sc Lafon J.-P., Marot J.}
{\it Alg\`ebre locale}.
Hermann, Paris, (2002).

\bibitem[Lakatos]{La} {\sc Lakatos I.} {\em Preuves et réfutations}. 
Version fran\c{c}aise, Hermann (1984).


\bibitem[Lam]{Lam} {\sc  Lam T.Y.}
{\it Serre's conjecture}.
Lecture Notes in Mathematics, Vol. 635.
Springer Berlin Heidelberg New York, (1978).

\bibitem[Lam06]{Lam06} {\sc  Lam T.Y.}
{\it Serre's Problem on Projective Modules}.
Springer Berlin Heidelberg New York, (2006).

\bibitem[Lancaster \& Tismenetsky]{Lan} {\sc Lancaster P., Tismenetsky M.}
{\em The Theory of Matrices, 2/e} Academic Press, (1985).

\bibitem[Lawvere \& Rosebrugh]{LaRo} {\sc Lawvere W., Rosebrugh R.} 
{\em Sets for Mathematics} Cambridge University Press, (2003).

\bibitem[Lorenzen]{Loren} {\sc Lorenzen P.} 
{\em Métamathématique}. Traduit de l'allemand par J. B. Grize, Gauthier-Villars, Paris, Mouton, Paris La Haye (1967), édition originale 1962.

\bibitem[Mac Lane]{MACL} {\sc Mac Lane, S.} {\em Categories for the Working Mathematician}. Second edition, Springer, (1998).


\bibitem[Martin-L\"of]{PML} {\sc Martin-L\"of P.} {\it Intuitionistic type theory}. Notes by Giovanni Sambin. Studies in Proof Theory. Lecture Notes, 1. Bibliopolis, Naples, (1984).

\bibitem[Matsumura]{Mat} {\sc Matsumura H.} {\it Commutative ring theory}.
Cambridge studies in advanced mathematics \num8.
Cambridge University Press, (1989).

\bibitem[MITCA]{MITCA} {\it  Multiplicative Ideal Theory in Commutative Algebra: A tribute to the work of Robert Gilmer.} {\sc Eds: Brewer J., Glaz G., Heinzer W.,  Olberding B.} Springer, (2006)

\bibitem[MRR]{MRR} {\sc Mines R., Richman F., Ruitenburg W.} {\it  A Course in
Constructive Algebra.} Universitext. Springer-Verlag, (1988).

\bibitem[Mora]{Mora} {\sc Mora T.} {\it   Solving Polynomial Equation Systems I: The Kronecker-Duval Philosophy.}   Cambridge University Press, (2003)

\bibitem[Northcott]{Nor} {\sc Northcott D.} 
{\it Finite free resolutions.} Cambridge
tracts in mathematics No 71. Cambridge University Press, (1976).

\bibitem[PFCM]{PFCM}{\sc Crosilla L.,
Schuster P., eds.} {\it From Sets and Types to Analysis and Topology: Towards
Practicable Foundations for Constructive Mathematics}. 
Oxford University Press, (2005).

\bibitem[Pohst \& Zassenhaus]{PZ}  {\sc Pohst,  Zassenhaus}
{\em Algorithmic algebraic number theory  (Encyclopedia of Mathematics and its Applications).}
Cambridge University Press, (1989).




%: R

\bibitem[Rao \& Mitra]{RM}  {\sc Rao C.,  Mitra S.}
{\em Generalized Inverses of Matrices and its Applications.}
John Wiley \& Sons, (1971).

\bibitem[Raynaud]{Ray}  {\sc Raynaud M.}
{\em Anneaux locaux henséliens.}
Springer Lecture Notes in Mathematics \num169, (1970).

\bibitem[SINGULAR]{Singular}  {\sc Greuel G.-M., Pfister G.}
{\em A Singular Introduction to Commutative Algebra.}
Springer (2002). {\url{http://www.singular.uni-kl.de/}}

\bibitem[Schwichtenberg \& Wainer]{SW}  {\sc Schwichtenberg H., Wainer S.}
{\em Proofs and Computations.}
Perspectives in Logic. Association for Symbolic Logic and Cambridge University Press, (2012).


\bibitem[Stacks-Project]{Stacks}  {\sc Stacks-Project.} Ouvrage collectif.
\url{http://stacks.math.columbia.edu}

\bibitem[TAPAS]{CCS} {\sc Cohen A., Cuypers H., Sterk H. }(eds)
{\it Some Tapas of Computer Algebra}.  Springer Verlag, (1999).

\bibitem[Tignol]{Tignol} {\sc Tignol J.-P.} {\it Galois' theory of algebraic equations.} World Scientific Publishing Co., Inc., River Edge, NJ, (2001). 

%\bibitem[Tchebotarev]{Tchebotarev} {\sc Tchebotarev N.}. 1934
%{\it Osnovy teorii Galua}.
%1950 {\it Grundz\"uge des Galois'shen Theorie}. P. Noordhoff.
\bibitem[Yengui]{Yengui}{\sc Yengui I.} {\it Constructive commutative algebra. Projective modules over polynomial rings and dynamical Gr\"obner bases.} Springer LNM \num 2138 (2015).

\bibitem[Zaanen]{Zaanen} {\sc Zaanen A.} {\it Introduction to Operator
Theory in Riesz Spaces}.
Springer Verlag, (1997).
~\par\medskip
%:  articles

\goodbreak
\noi{\bf\Large Articles}
%: A

\bibitem{Acz} {\sc Aczel P.} {\it  Aspects of general topology in constructive set theory.}  Ann. Pure Appl. Logic, {\bf 137},  (2006),   3--29.

\bibitem{AV} {\sc Aubry P.,  Valibouze A.} {\it  Using Galois Ideals for
Computing Relative Resolvents.}
J. Symbolic Computation, {\bf 30}, (2000), 635--651.

\bibitem{AG} {\sc Auslander M., Goldman, O.} {\it  The Brauer group of a commutative ring.}
Trans. Amer. Math. Soc., {\bf 97}, (1960), 367--409.


\bibitem{Avi}  {\sc Avigad J.} \emph{Methodology and metaphysics in the
development of \hbox{Dedekind's} theory of ideals}. Dans: José Ferreir\'os and Jeremy Gray, editors, The Architecture of Modern Mathematics, Oxford University Press, (2006), 159--186.

%: B
\bibitem{Ba} {\sc Banaschewski B.}
{\em Radical ideals and coherent frames.}
Comment. Math. Univ. Carolin. {\bf 37}  2, (1996),  349--370.

\bibitem{Barh09} {\sc Barhoumi S.}
 {\it Seminormality and polynomial rings.}
Journal of Algebra {\bf 322} (2009), 1974--1978.

\bibitem{BL07} {\sc Barhoumi S., Lombardi H.}
 {\it An Algorithm for the Traverso-Swan theorem on seminormal rings.}
Journal of Algebra {\bf 320} (2008), 1531--1542.

\bibitem{BLY06} {\sc Barhoumi S., Lombardi H., Yengui I.}
{\it  Projective modules over polynomial rings:
a constructive approach.}  Math. Nachrichten  {\bf 282}  (2009), 792--799.


\bibitem{BaRaKha} {\sc Basu R., Rao R., Khanna R.} 
{\it On Quillen's Local Global Principle.}
Contemporary Mathematics, Commutative Algebra and Algebraic Geometry, Volume 390, (2005), 17--30.

\bibitem{Bass} {\sc Bass H.} 
{\it Torsion free and projective modules.}
Trans. Amer. Math. Soc. {\bf 102}, (1962), 319--327.

\bibitem{BG} {\sc Bazzoni S., Glaz S.} {\it Pr\"ufer rings}.
dans \cite{MITCA}, 55--72.


\bibitem{Berg} {\sc Berger J.}
{\it  Constructive Equivalents of the Uniform Continuity Theorem}.
Journal of Universal Computer Science  {\bf 11} \num12 (2005), 1878--1883.

\bibitem{BeBr} {\sc Berger J., Bridges D.} {\it
A fan-theoretic equivalent of the antithesis of Specker's theorem}.
Proc. Koninklijke Nederlandse Akad. Wetenschappen.
Indag. Math.   {\bf 18} \num2 (2007), 195-202.

\bibitem{BeIs} {\sc Berger J., Ishihara H.} {\it
Brouwer's fan theorem and unique existence in constructive analysis}.
Math. Logic Quarterly  {\bf 51} (2005), 360--364.

\bibitem{Ber1}  {\sc Bernstein D.} {\it  Factoring into coprimes in
essentially
linear time}.  Journal of Algorithms {\bf 54} (2005), 1--30.

\bibitem{Ber2} {\sc Bernstein D.}
{\it  Fast ideal arithmetic via lazy localization.}
Cohen, Henri (ed.), Algorithmic number theory. Second international
symposium,
ANTS-II, Talence, France, May 18-23, 1996. Proceedings. Berlin:
Springer.
Lect. Notes Comput. Sci. \num1122 (1996), 27--34.


\bibitem{bi} {\sc Bishop, E.} {\it  Mathematics as a numerical language}.
in  Intuitionism and Proof Theory. Eds. Myhill, Kino, and Vesley,
North-Holland, Amsterdam, (1970)

\bibitem{BoSc} {\sc Boniface J., Schappacher N.} {\it "Sur le concept de nombre en mathématique": cours inédit de Leopold Kronecker \`a Berlin (1891).}
Rev. Histoire Math.  {\bf 7} (2001), 206--275.

\bibitem{Bosma} {\sc Bosma W., Cannon J., Playoust C.} {\it The Magma algebra system. I. The user language.}
J.~Symbolic Comput.  {\bf 24} (1997), 235--265.

\bibitem{Brandl} {\sc Brandl R.} {\it Integer polynomials that are reducible modulo all primes.}
Amer. Math. Month. {\bf 93}, no. 4 (1986), 286--288.

\bibitem{Brenner} {\sc Brenner H.} {\it Lifting chains of prime ideals.}
J. Pure Appl. Algebra {\bf 179} (2003), 1--5.

\bibitem{BC} {\sc Brewer J., Costa D.}
{\it  Projective modules over some non-Noetherian polynomial
rings.} J. Pure Appl. Algebra {\bf 13} (1978), no. 2, 157--163.

\bibitem{BC2} {\sc  Brewer J., Costa  D.} {\it Seminormality and projective modules over polynomial rings.}
J. Algebra {\bf 58}, no. 1 (1979), 208--216.

\bibitem{BrKl} 
 {\sc  Brewer J., Klinger L.} {\it Pole assignability and the invariant factor theorem in Pr\"ufer domains 
and Dedeking domains}. J. Algebra {\bf 111} (1987), 536--545.

\bibitem{LB} {\sc Buchmann J.,  Lenstra H.}  {\it Approximating rings of
integers in number fields}.  J. Th\'{e}or. Nombres Bordeaux,  {\bf 6} (2)
(1994),   221--260.

%: C

\bibitem{Cah90}
{\sc Cahen, P.-J.},
{\it  Construction B, I, D et anneaux localement ou r\'esiduellement de Jaffard. (B, I, D construction and locally or residually Jaffard rings).},
Archiv der Mathematik,
{\bf 54}, (1990), 125--141.

\bibitem{CCDHKS} {\sc Caniglia L., Cortinas G., Dan{\'o}n S.,
Heintz J., Krick T., Solern{\'o} P.} 
{\it  Algorithmic Aspects of Suslin's Proof of Serre's Conjecture}. 
Computational Complexity {\bf 3} (1993), 31--55. 

\bibitem{Cannon} {\sc Cannon J.,  Bosma W.} {\it Handbook of Magma functions.}
Version 2.14, Oct. 2007, 4400 pages.

\bibitem{cc}  {\sc Cederquist J., Coquand T.} {\em Entailment relations and
Distributive Lattices} Logic Colloquium '98 (Prague), 127--139,
Lect. Notes Log., 13. Assoc. Symbol. Logic, Urbana, (2000).

\bibitem{CHR}  {\sc Chase S., Harrison D., Rosenberg A.}
{\em Galois theory and Galois cohomology of commutative rings.}
 Mem. Amer. Math. Soc. {\bf 52} (1965), 15--33.


\bibitem{CT03} {\sc Chervov A., Talalaev D.}
{\it Hitchin systems on singular curve I}. 
Theor. Math. Phys. {\bf 140} (2004), 1043--1072.

\bibitem{CT03bis} {\sc Chervov A., Talalaev D.}
{\it Hitchin systems on singular curve II. Glueing subschemes}.
Int. J. Geom. Meth. Mod. Phys {\bf 4} (2007), 751--787.

%\bibitem{CFGK} {\sc Chistov A., Fournier H., Gurvits L., Koiran P.}
%{\em  Vandermonde Matrices, NP-Completeness, and Transversal Subspaces}.
%preprint 2002.\perso{donner la référence exacte}

\bibitem{CoqK} {\sc Coquand T.} {\it La contribution de Kolmogorov
en logique intuitionniste}.
dans: L'héritage de Kolmogorov en mathématiques. Charpentier~E.,
Lesne~A., Nikolski~N. (eds). Belin, Paris (2004).

\bibitem{CoqFT} {\sc Coquand T.} {\it About Brouwer's fan theorem.}
Revue internationale de philosophie, {\bf 230} (2004),
 483--489.


\bibitem{Coq3} {\sc Coquand T.}
{\it Sur un théor\`eme de Kronecker concernant les variétés algébriques}.
C.~R. Acad. Sci. Paris, Ser. I {\bf 338} (2004),  291--294.


\bibitem{coq} {\sc Coquand T.}
{\it On seminormality}. Journal of Algebra, {\bf 305} (1),
(2006), 585--602.


\bibitem{coq07} {\sc Coquand T.}
{\it A refinement of Forster's theorem.}
Preprint (2007).

\bibitem{coqval} {\sc Coquand T.}
{\it Space of valuations}. Annals of Pure and Applied Logic,
{\bf 157}  (2009), 97--109.


%: paru
\bibitem{CoqRec} {\sc Coquand T.}
{\it Recursive functions and constructive mathematics.}
p. 159--167 dans:
Bourdeau M., Dubucs J. (Eds.),
Calculability and Constructivity. Historical and Philosophical Aspects. Logic, Epistemology and the Unity of Science, Vol. 34. Springer (2014).


\bibitem{CDLQ02} {\sc Coquand T., Ducos L., Lombardi H., Quitté C.}
{\it L'idéal des coefficients du produit de deux polyn\^omes}.
Revue des Mathématiques de l'Enseignement Supérieur, {\bf 113} (3), (2003),
25--39.

\bibitem{CDLQ06} {\sc Coquand T., Ducos L., Lombardi H., Quitté C.}
{\it  Constructive Krull Dimension. I: Integral Extensions.}
Journal of Algebra and Its Applications.
{\bf 8}  (2009), 129--138.

\bibitem{CL05} {\sc Coquand T., Lombardi H.} {\it  A logical approach to abstract algebra.}
(survey) Math. Struct. in Comput. Science {\bf 16} (2006), 885--900. 


\bibitem{cl} {\sc Coquand T., Lombardi H.}
{\it Hidden constructions in abstract algebra (3)
Krull dimension of distributive lattices and commutative rings}, dans \cite{CRA},
 477--499.

\bibitem{cl1} {\sc Coquand T., Lombardi H.}
{\it Constructions cachées en alg\`ebre abstraite (3) 
Dimension de Krull, Going Up, Going Down}. Rapport technique (2001)
\url{http://hlombardi.free.fr/publis/GoingUpDownFrench.pdf}
(version anglaise \url{http://hlombardi.free.fr/publis/GoingUpDown.pdf}).

\bibitem{cl2016} {\sc Coquand T., Lombardi H.}
{\it Some remarks on normal rings},   dans \cite{PS2016} 141--149.

\bibitem{clq} {\sc Coquand T., Lombardi H., Quitté C.}
{\it Generating non-Noetherian modules constructively}.
Manuscripta mathematica, {\bf 115}  (2004), 513--520.


\bibitem{clq2} {\sc Coquand T., Lombardi H., Quitté C.}
{\it Dimension de Heitmann des \trdis et des anneaux commutatifs}.
Publications Mathé\-matiques de Besan\c{c}on.
Théorie des nombres (2006). 51 pages.
Version corrigée:
\url{http://hlombardi.free.fr/publis/AHeitmann.html}

\bibitem{clr03} {\sc Coquand T., Lombardi H., Roy M.-F.} {\it
An elementary characterization of Krull dimension}.
dans \cite{PFCM}, 239--244.

\bibitem{cls} {\sc Coquand T., Lombardi H., Schuster P.}
{\it A nilregular element property}. Archiv der Mathematik,
{\bf 85} (2005), 49--54.

\bibitem{cp} {\sc Coquand T., Persson H.}
{\it Valuations and Dedeking Prague theorem}.
J. Pure Appl. Algebra, \textbf{155} (2001), 121--129.

\bibitem{CHWW} {\sc Corti\~nas G., Haesemayer C., Walker M.E. and Weibel C.}
{\it A negative answer to a question of Bass}.
Proc. AMS, \textbf{139} (2011), 1187--1200.


\bibitem{clr} {\sc Coste M., Lombardi H., Roy M.-F.} {\it
Dynamical method in algebra: Effective Nullstellens\"atze}.
Annals of Pure and Applied Logic,
{\bf 111}, (2001), 203--256.

\bibitem{Couc1} {\sc Couchot F.} {\it
Finitely presented modules over semihereditary rings}.
Communications in Algebra,  {\bf 35}  (9),  (2007) 2685--2692.

\bibitem{Ddk1} {\sc Dedekind R.} {\it
\"Uber einen arithmetischen Satz von Gauss}.
Mitt. dtsch. math. Ges. Prag.   (1892), 1--11.

\bibitem{Ddk2} {\sc Dedekind R.} {\it
\"Uber die Begr\"undung der IdealTheorie}.
Nachr. K. Ges. Wiss. G\"ottingen   (1894), 272--277.

\bibitem{D5} {\sc Della Dora J., Dicrescenzo C., Duval D.} {\it About a new
method for computing in algebraic number fields}.
In Caviness B.F. (Ed.)   EUROCAL '85. Lecture Notes in
Computer Science 204,  289--290. Springer (1985).

%: D

\bibitem{DiJSC} {\sc D\'{\i}az-Toca G.}
{\em  Galois theory, splitting fields and computer algebra}. 
J. Symbolic Computation {\bf 41} \num11, (2006),  1174--1186.

\bibitem{DiGL} {\sc D\'{\i}az-Toca G., Gonzalez-Vega L., Lombardi H.}
{\em  Generalizing Cramer's Rule: Solving uniformly linear systems
of equations}. SIAM Journal on Matrix Analysis and Applications
{\bf 27} \num3, (2005), 621--637.

\bibitem{DiGLQ} {\sc D\'{\i}az-Toca G., Gonzalez-Vega L., Lombardi H. \&
Quitté C.}
{\em Modules projectifs de type fini, \alis croisées et inverses généralisés}. Journal of Algebra
{\bf 303} \num2, (2006), 450--475.

\bibitem{DiL09} {\sc D\'{\i}az-Toca G., Lombardi H.} 
{\it  A polynomial bound on the number of comaximal localizations 
needed in order to make free a projective module}.
Linear Algebra and its Application. {\bf  435}, (2011), 354--360.

\bibitem{DiLQ} {\sc D\'{\i}az-Toca G.,  Lombardi H., Quitté C.}
{\it  L'alg\`ebre de décomposition universelle.}
Proceedings du colloque TC2006 (Granada), 169--184.

\bibitem{DiL10} {\sc D\'{\i}az-Toca G.,  Lombardi H.} 
{\it  Dynamic Galois Theory}. Journal of Symbolic Computation.
{\bf  45}, (2010), 1316--1329.

%:2015  intégration
\bibitem{Drach} {\sc Drach J.}
{\em Essai sur la théorie générale de l'intégration 
et sur la classification des transcendantes}. 
Ann. Sci. Ec. Norm. Sup {\bf 3} \num15, (1898), 243--384.

\bibitem{Du} {\sc Ducos L.} 
{\it  Effectivit\'{e} en th\'{e}orie de Galois. Sous-r\'{e}sultants}. 
Université de Poitiers, Th\`ese doctorale. Poitiers (1997).

\bibitem{Duc} {\sc Ducos L.}
{\it  Construction de corps de décomposition
gr\^ace aux facteurs de résolvantes.  (French) [Construction of
splitting fields in favour of resolvent factors]}. Communications
in Algebra {\bf 28} \num2, (2000), 903--924.

\bibitem{DQ} {\sc Ducos L.}
{\it Vecteurs unimodulaires et syst\`emes générateurs.}
Journal of Algebra {\bf 297}, (2006), 566--583.


\bibitem{Duc08} {\sc Ducos L.}
{\it Sur la dimension de Krull des anneaux \noesz.}
Journal of Algebra {\bf 322}, (2009), 1104--1128.

%\bibitem{DucPVE} {\sc Ducos L.}
%{\it Polyn\^omes \`a valeurs enti\`eres : un anneau de Pr\"fer de dimension 2.} (2011) To appear in Communications in Algebra.
%
% \bibitem{dlqs}
%{\sc Ducos L., Lombardi H., Quitt\'e C., Salou M.}
%{\it Th\'{e}orie algorithmique des anneaux arithm\'etiques,
%des anneaux de Pr\"ufer et des anneaux de Dedekind}.
%Journal of Algebra. {\bf 281}, (2004), 604--650.
\bibitem{DucPVE} {\sc Ducos L.}
{\it Polyn\^omes à valeurs enti\`eres : un anneau de Pr\"fer de dimension 2.} (2011) To appear in Communications in Algebra. 

 
\bibitem{dlqs}{\sc Ducos L., Lombardi H., Quitté C., Salou M.}
{\it   Th\'{e}orie algorithmique des anneaux arithmétiques,
des anneaux de Pr\"ufer et des anneaux de Dedekind}.
Journal of Algebra {\bf 281}, (2004), 604--650.
 
\bibitem{DVY2015} {\sc Ducos L., Valibouze A., Yengui I.}
{\it   Computing syzygies over $V [X_1, \ldots, X_k]$, $V$ a valuation domain}.
Journal of Algebra {\bf 425}, (2015), 133--145.


%: E

%\bibitem{La} {\sc Ekedahl E., Laskov D.} {\it  Splitting algebras, symmetric
%functions ans Galois Theory.} Journal of Algebra and its Applications, {\bf
%4} (1), (2005), 59--76.

\bibitem{Edw81}
{\sc Edwards H.}
{\it  The genesis of ideal theory.}
Arch. Hist. Exact Sci. {\bf 23} \num4 (1980/81), 321--378.

\bibitem{Eisenbud} {\sc Eisenbud D., Evans E., Jr}.
{\em Generating modules efficiently: theorems from algebraic $K$-theory}.
J. Algebra {\bf 27} (1973), 278--305.

\bibitem{E2} {\sc Eisenbud D., Evans E., Jr.}
{\em Every algebraic set in $n$-space is the intersection of $n$ hypersurfaces}.
Inventiones math. {\bf 19} (1973), 107--112.


\bibitem{ELY07} {\sc Ellouz A., Lombardi H., Yengui I.}
{\it  A constructive comparison of the rings
$\mathbf{R}(X)$ and $\mathbf{R} \langle X \rangle$
and application to the Lequain-Simis Induction 
Theorem.} 
Journal of Algebra. {\bf 320}  (2008), 521--533.


\bibitem{espThesis} {\sc Espa\~nol L.} {\em Dimensi\'on en \'algebra constructiva}. Th\`ese doctorale. Université de Zaragoza,
Zaragoza, (1978).

\bibitem{esp} {\sc Espa\~nol L.} {\em Constructive Krull dimension of
lattices.}
Rev.\ Acad.\ Cienc.\ Zaragoza (2) {\bf 37} (1982), 5--9.

\bibitem{esp83} {\sc Espa\~nol L.} {\em Le spectre d'un anneau dans l'alg\`ebre constructive et applications \`a la dimension.}
Cahiers de topologie et géométrie différentielle catégorique. {\bf 24} \num2 (1983), 133--144.




\bibitem{esp86} {\sc Espa\~nol L.} {\em Dimension of Boolean Valued Lattices and Rings.}
Journal of Pure and Applied Algebra {\bf 42} (1986), 223--236.

\bibitem{Es} {\sc Espa\~nol L.}
{\em The spectrum lattice of Baer rings and polynomials.}
Categorical algebra and its applications.
(Louvain-La-Neuve, 1987), 118--124,
Lecture Notes in Math., 1348,
Springer, Berlin-New York, (1988).

\bibitem{Espa08} {\sc Espa\~nol L.} \emph{Finite chain calculus in distributive lattices and elementary Krull dimension.} Contribuciones cientificas en honor de Mirian Andres Gomez.
Eds. L. Lamban, A. Romero y J. Rubio, 
Servicio de Publicaciones, Universidad de La Rioja, Logrono, Spain, (2010).

\bibitem{EG} {\sc Estes R., Guralnick R.}
{\em Module equivalences: local to global when primitive polynomials
represent units}.
J. of Algebra {\bf 77} (1982), 138--157.

\bibitem {Ferrand} {\sc Ferrand D.}
{\it Les modules projectifs de type fini sur un anneau
de polyn\^omes sur un corps sont libres}.
Sém. Bourbaki, exposé {\bf 484}, (1975-1976), 202--221.

\bibitem{FP} {\sc Ferrero M., Paques A.}
{\em Galois theory of commutative rings revisited}.
Contributions to Algebra and Geometry {\bf 38} (1997), 399--410.

\bibitem{FiGa} {\sc Fitchas N., Galligo A.} {\it  Nullstellensatz effectif et 
Conjecture de Serre (Théor\`eme de Quillen-Suslin) pour le Calcul
 Formel}. 
Math. Nachr. {\bf 149}  (1990), 231--253. 

\bibitem{FL} {\sc Fontana  M., Loper A.}
{\em An historical overview of Kronecker function rings, Nagata rings and related star and semistar operations}.  dans   \cite{MITCA},  169--187.

\bibitem{Forster} {\sc Forster O.}
{\em \"{U}ber die Anzahl der Erzeugenden eines Ideals in einem Noetherschen Ring.}
{Math. Z. {\bf 84} (1964), 80--87.}

\bibitem{Fuchs} {\sc Fuchs L.}
{\em \"{U}ber die  Ideale arithmetischer ringe.}
{Math. Helv. {\bf 23} (1949), 334--341.}



%: G


\bibitem{Gauss} {\sc Carl Friedrich Gauss} {\it Demonstratio nova altera theorematis omnem functionem algebraicam 
rationalem integram unius variabilis in factores reales primi vel secundi gradus resolvi posse}.
Comm. Recentiores (Gottingae), {\bf 3}  (1816), 107--142. 
Reproduit dans Werke III, 31--56.
Traduction anglaise: \url{http://www.monad.me.uk/misc/gauss-web.php}
sur la page web de Paul Taylor. \url{http://www.monad.me.uk/}

\bibitem{GeKl} {\sc Geissler K., Kl\"uners J.} {\it Galois Group Computation for Rational Polynomials}.
J. Symbolic Computation {\bf 30} (2000), 653--674.

\bibitem {GillmanHenriksen} {\sc Gillman L., Henriksen M.}
{\it Some remarks about elementary divisor rings}.
Trans. Amer. Soc. {\bf 82}, (1956) 362--365

\bibitem{GH} {\sc Gilmer  R., Heitmann  R.} {\it On Pic R[X] for R seminormal}.
 J. Pure Appl. Algebra {\bf 16}
(1980), 251--257.


\bibitem{GHo} {\sc Gilmer  R., Hoffmann, J.} {\it A characterization of Pr\"ufer
domains in terms of polynomials}.
Pacific J. Math. {\bf 60} (1)
(1975), 81--85.

%\bibitem{GMa} {\sc Glass A., Marra V.}
%{\it Embedding finitely generated abelian lattice-ordered groups:
%Higman's theorem and a realisation of $\pi$.}
%J. London Math. Soc. (2) {\bf 68} (2003), no. 3, 545--562.


\bibitem{Glaz1} {\sc Glaz S.} {\it Finite conductor properties
of $ \gR(X)$ and $\gR\lra X$}.  dans: Proceeding of conference in honor to J.
Huckaba's retirement, Missouri, Dec. 1999. Marcel Dekker Lecture
Notes.


\bibitem{Glaz2} {\sc Glaz, S., Vasconcelos W.} {\em  Gaussian polynomials.}
Marcel Dekker Lecture Notes 186 (1997), 325--337.


\bibitem{Gold}{\sc Goldman  O.}
{\em Determinants in projective modules}. Nagoya Math. J.
{\bf 18} (1961), 27--36.

%:HHH non utilis
%\bibitem{GV}{\sc Greenberg  B., Vasconcelos W.}
%{\em Coherence of polynomial rings}. Proc. Amer. Math. Soc.
%{\bf 54} (1976), 59--64.


%: H
\bibitem{ValuHallouin} {\sc Hallouin E}.
{\em Parcours initiatique \`a travers la théorie des valuations.} Rapport technique. Université de Poitiers, (1996). \url{http://www.picard.ups-tlse.fr/~hallouin/eh-valuation.ps}

\bibitem{TheseHallouin} {\sc Hallouin E}.
{\em Calcul de fermeture intégrale en dimension 1 et factorisation intégrale.}
Th\`ese. Université de Poitiers, (1998). \url{http://www.picard.ups-tlse.fr/~hallouin/eh-these.ps}

\bibitem{Hei76} {\sc Heitmann R}.
{\em Generating ideals in Pr\"ufer domains.}
Pacific J. Math. {\bf 62} (1976), 117--126.

\bibitem{Hei84} {\sc Heitmann R.}
{\em Generating non-Noetherian modules efficiently.}
Michigan Math. {\bf 31} 2 (1984), 167--180.


\bibitem{HeLe} {\sc Heitmann R., Levy L.}
{\it 1 1/2 and 2 generator ideals in Pr\"{u}fer domains.}
Rocky Mountain J. Math.  {\bf 5} 3  (1975), 361--673.


\bibitem{HS}  {\sc Hermida J., S\'anchez-Giralda T.}
{\it  Linear Equations over Commutative Rings and Determinantal Ideals.}
Journal of Algebra {\bf 99} (1986), 72--79.

\bibitem{Hess02}  {\sc Hess F.}
{\it  Computing Riemann-Roch space in algebraic function fields.}
Journal of Symbolic Computation {\bf 33} (2002), 425--445.

\bibitem{HeyRec} {\sc Heyting A.}
{\it  After thirty years.} 
In: 1962 Logic, Methodology and Philosophy of Science (Proc. 
1960 Internat. Congr.) pp. 194--197 Stanford Univ. Press, Stanford, Calif.

\bibitem{Hi2} {\sc Hilbert D.}
{\em \"Uber das Unendliche}.
Math. Annalen {\bf 95}  (1926), 161--190.
(Sur l'infini) traduction anglaise dans \cite{vH} 367--392.

\bibitem{Hoc} {\sc Hochster M.}
{\em Prime ideal structure in commutative rings.}
Trans. Amer. Math. Soc. {\bf 142} (1969), 43--60.

\bibitem{Hor} {\sc Horrocks G.}
{\it  Projective modules over an extension of a local ring}.
Proc. Lond. Math. Soc. {\bf 14} (1964), 714--718.

\bibitem{Hulpke} {\sc Hulpke A.}
{\it  Konstruktion transitiver Permutationsgruppen}.
Dissertation, Rheinisch-Westf\"alische 
Technische Hochschule, Aachen, Germany.
(1996).  

%: I

\bibitem{Ishi} {\sc Ishihara H.} {\it
Constructive reverse mathematics: compactness properties}. dans
\cite{PFCM}, 245--267.


\bibitem{Ishi2} {\sc  Ishihara H.} {\it
Weak K\"onig lemma implies Brouwer's fan theorem: a direct proof}.
Notre Dame J. Formal Logic  {\bf 47} (2006), 249--252.

\bibitem{Ishi3} {\sc  Ishihara H.} {\it
Reverse mathematics in Bishop's constructive mathematics}.
Philosophia Scientiae, Cahier Spécial  {\bf 6} (2006), 43--59.

\bibitem{JL} {\sc  Jacobsson C., L{\"o}fwall C.} {\it
Standard Bases for General Coefficient Rings and a New Constructive Proof 
of {H}ilbert's Basis Theorem}.
J. Symb. Comput.  {\bf 12} (1991),  {337--372}.


\bibitem{Joh} {\sc Johnstone, P.} {\it
The art of pointless thinking: a student's guide to the category of locales.} Category theory at work (Bremen, 1990), 85--107,
Res. Exp. Math., 18, Heldermann, Berlin, 1991.


\bibitem{Joyal71} {\sc Joyal A.} {\em Spectral spaces and distibutive lattices.}
Notices AMS {\bf 18} (1971), 393.

\bibitem{Joyal} {\sc Joyal A.} {\em Le théor\`eme de Chevalley-Tarski.}
Cahiers de topologie et géometrie différentielle catégorique, 1975.

%: K
\bibitem{vdK} {\sc van der Kallen W.}  {\it The K2 of rings with many units.}
 Ann. Sci. \'E.N.S.  4\ieme série, {\bf 10}, (1977), 
473--515.


\bibitem{Kap} {\sc Kaplansky I.}  {\it    Elementary divisors and
modules.}
Transactions of the AMS {\bf 66}, (1949), 464--491.

\bibitem{Kap52}{\sc  Kaplansky I.} {\it Modules over Dedekind Rings and
Valuation Rings.}
Trans. Amer. Math. Soc. {\bf 72},
(1952), 327--340.

\bibitem{KlMa}{\sc  Kl\"uners J., Malle G.} 
{\it Explicit Galois realization of transitive groups of degree up to 15.}
J. Symbolic Comput. {\bf 30} (6), (2000), 675--716.

\bibitem{Kol} {\sc Kolmogorov A.}  {\it  Zur Deutung der intuitionistischen Logik.}  Math. Zeitschr., {\bf 35} (1932) 58--65.

\bibitem{Kro1} {\sc Kronecker L.}
{\it  Zur Theorie der Formen h\"oherer Stufen} Ber. K. Akad. Wiss. Berlin
   (1883), 957--960. (Werke 2, 417--424).

\bibitem{Kronecker} {\sc Kronecker L.}
{\it Grundz\"uge einer arithmetischen Theorie der algebraischen 
Gr\"ossen.}
{J. reine angew. Math.} {\bf 92}, (1882) 1--123.
Réimprimé dans {\it Leopold Kronecker's Werke}, II, 237--387.

%: L

\bibitem{LaMi} {\sc Landau, S., Miller, G.}  {\it  Solvability by radicals is in polynomial time.} 
J. Comput. Syst. Sci.  {\bf 30} (1985), 179--208.


\bibitem{Lecerf-Factsep} {\sc Lecerf, G.} 
{\it Fast separable factorization and applications.}
 Applicable Algebra in Engineering, Communication and Computing,  {\bf 19} (2) (2008), 135--160.


\bibitem{LS} {\sc Lequain, Y., Simis, A.}
{\it Projective modules over $R[X_1,... ,X_n]$,
$R$ a Pr\"ufer domain.}
J. Pure Appl. Algebra {\bf 18} (2) (1980),  165--171.

%\bibitem{Lom} {\sc Lombardi H.} {\it  Relecture constructive de la
%théorie d'Artin-Schreier.}  Annals of Pure and Applied
%Logic  {\bf 91}, (1998), 59--92.

\bibitem{Lom97}  {\sc Lombardi H. } {\it  Le contenu constructif
d'un principe local-global avec une application \`a la
structure d'un module projectif de type fini.}
Publications Mathé\-matiques de Besan\c{c}on.
Théorie des nombres. Fascicule   (1997), 94--95 \& 95--96.

\bibitem{lom99} {\sc Lombardi H.} {\it Platitude, localisation et anneaux de
Pr\"ufer: une approche constructive}.
64 pages.
Publications Mathématiques de Besan\c{c}on.
Théorie des nombres. Années 1998-2001.


\bibitem{lom}  {\sc Lombardi H.} {\em Dimension de Krull,
Nullstellens\"atze  et \'Evaluation dynamique}.
Math. Zeitschrift,  {\bf 242}, (2002),  23--46.

\bibitem{LomPVE}  
{\it  Un anneau de Pr\"ufer.}   
  Third International Meeting on Integer-Valued Polynomials.
Actes des rencontres du CIRM,  {\bf 2}  
(2010). \url{http://acirm.cedram.org/cgi-bin/browse}


\bibitem{LQ02}  {\sc Lombardi H., Quitté C.}
{\it Constructions cachées en alg\`ebre abstraite (2)
Le principe local global}.
dans \cite{CRA} 461--476.


\bibitem{LQ06}  {\sc Lombardi H., Quitté C.}
{\it  Seminormal rings (following Thierry Coquand).} 
Theoretical Computer Science. {\bf 392}, (2008), 113--127.


\bibitem{LQY05} {\sc Lombardi H., Quitté C., Yengui I.}
{\it  Hidden constructions in abstract
algebra (6) The theorem of Maroscia, Brewer and Costa.}
Journal of Pure and Applied Algebra.
{\bf 212} 7 (2008),  1575--1582.



\bibitem{YL04} {\sc Lombardi H., Yengui I.}
{\it  Suslin's algorithms for reduction of
unimodular rows.}
Journal of Symbolic Computation {\bf 39} (2005), 707--717.

%: M N

\bibitem{Ma} {\sc Maroscia P.} {\it  Modules projectifs sur
certains anneaux de polyn\^omes.} C.R.A.S. Paris {\bf 285} série A
(1977), 183--185.

%:2015
\bibitem{ML1973} {\sc Per Martin-L\"of.} 
{\it  An intuitionistic theory of types: Predicative part.}
In H. E. Rose and J. C. Shepherdson, editors, Logic Colloquium Ô73, 
pages 73--118. North Holland, (1975).

\bibitem{ML} {\sc Martin-L\"of P.}
{\it  An intuitionistic theory of types}.   127--172,  in: Twenty-five years of constructive type theory (Venice, 1995),
Oxford Logic Guides, 36, Oxford Univ. Press, New York, 1998.

\bibitem{ML2008} {\sc Per Martin-L\"of}
{\it  The Hilbert-Brouwer controversy resolved?}  
dans: One hundred years of intuitionism (1907-2007), (Cerisy),  (Mark Van Atten \& al., editors) Publications des Archives Henri Poincaré, 
Birkh\"auser Basel, (2008), pp. 243--256.

\bibitem{Mer} {\sc Mertens F.} {\em
\"Uber einen algebraischen Satz}. Ber. K. Akad. Wiss. Wien
 (1892).\perso{pages \`a préciser, semble introuvable}

\bibitem{MnY}  {\sc Mnif A., Yengui I.}
{\it An algorithm for  unimodular completion over Noetherian rings.}
J. Algebra {\bf 316} (2007), 483--498. 


\bibitem{Mul} {\sc Mulmuley K.} {\em A fast parallel algorithm
to compute the rank of a matrix over an arbitrary field}.
Combinatorica,  {\bf 7}/1, (1987), 101--104.


\bibitem{Murthy03} {\sc Murthy M.} {\em Generators of a general ideal}.
in: A tribute to C. S. Seshadri (Chennai, 2002), Trends in Math.,
Birkh\"auser, Basel, (2003), 379--384.

\bibitem{MuPe} {\sc Murthy M., Pedrini C.} {\em $K_0$ and $K_1$ of polynomial rings}.
in Algebraic K-Theory II,
Lecture Notes in Math. 342,  (1973), 109--121.

\bibitem{NaNi} {\sc Nashier B., Nichols W.}
{\it  Ideals containing monics.}
 Proc. Amer. Math. Soc.
{\bf 99} (1987), 634--636.

\bibitem{Nic} {\sc Nicholson W.} {\it  Lifting idempotents and exchange rings.} Trans. Amer. Math. Soc. {\bf 229} (1977), 269--278.

\bibitem{Nor2} {\sc Northcott D.}
{\it  A generalization of a theorem on the content of polynomials.}
Proc. Cambridge Philos. Soc.
{\bf 55} (1959), 282--288.

%: O P

\bibitem{ORV} {\sc Orange  S., Renault G., Valibouze A.}
{\it Calcul efficace de corps de décomposition.}
Rapport technique LIP6 2003/005.
 

\bibitem{Per1} {\sc  Perdry H.}
{\it Strongly Noetherian rings and constructive ideal theory.}
J. Symb. Comput.
{\bf 37} (2004), 511--535.

\bibitem{Per2} {\sc  Perdry H.}
{\it   Lazy bases: a minimalist constructive theory of Noetherian rings.}
 Math. Log. Quart.
{\bf  54} (2008), 70--82.


\bibitem{Po} {\sc  Poincaré H.}
{\em La logique de l'infini}.
Revue de Métaphysique et de Morale {\bf 17}, 461--482, (1909)
réédité dans {\em Derni\`eres pensées}, Flammarion (1913).



%\bibitem{Pra}  {\sc Prasad K.}
%{\em Generalized Inverses of Matrices over Commutative Rings.}
%Linear Algebra Appl. {\bf 211} (1994), 35--52.

%\bibitem{PB}  {\sc Prasad K.,  Bapat R.}
%{\em The generalized Moore-Penrose inverse.}
%Linear Algebra Appl. {\bf 165} (1992), 59--69.

\bibitem{Prufer}  {\sc Pr\"ufer H.}
{\em Untersuchunger uber teilbarkeitseigenschaften in korpen.}
Angew. Mat. {\bf 168} (1932), 1--36.


%: Q

\bibitem{Que} {\sc Quentel Y.}
{\it  Sur une caractérisation des anneaux de valuation de hauteur~1.}
C.\,R.\,Acad.\,Sci., Paris, Ser. A {\bf 265} (1967), 659--661.

\bibitem{Querre} {\sc Querré J.}
{\em Sur le groupe de classes de diviseurs.}
C. R. Acad. Sci. Paris, {\bf 284} (1977), 397--399.


\bibitem{Qu} {\sc Quillen D.} {\it   Projective modules over
polynomial rings.}    Invent. Math. {\bf 36} (1976), 167--171.

%: R

%\bibitem{Rai}  {\sc Rais M.}
%Th\`ese (1970): {\it Distributions homog\`enes sur des espaces de matrices},
%Bulletin de la Société Mathématique de France, supplément au
%numéro de Juin 1972, mémoire \num30.


\bibitem{Rao85a}  {\sc Rao R.} {\it On projective $R_{f_1,\ldots,f_t}$-modules.} Amer. J. Math. {\bf 107} (1985), 387--406.

\bibitem{Rao85b}  {\sc Rao R.} {\it An elementary transformation of a special unimodular vector to its top coefficient vector.} 
Proc. Amer. Math. Soc. {\bf 93} (1985), 21--24.


\bibitem{Rao85c}  {\sc Rao R.} {\it A note on the Serre dimension of polynomial rings.} J. Pure Appl. Algebra
{\bf 38} (1985), 87--90.

\bibitem{RaSe}  {\sc Rao R., Selby J.} {\it Quillen-Suslin theory revisited.} J. Pure Appl. Algebra
{\bf 211} (2007),  541--546.


\bibitem{ric74} {\sc Richman F.} {\it  Constructive aspects of Noetherian
rings.} Proc. Amer. Mat. Soc. {\bf 44} (1974), 436--441.


\bibitem{ri81} {\sc Richman F.} {\it  Seidenberg's condition P.}
in Constructive Mathematics. Springer LNM 873 (1981), 1--11.

\bibitem{ri82} {\sc Richman F.} {\it  Finite dimensional algebras over discrete fields.}
L.\,E.\,J. Brouwer centenary symposium, Troelstra and van Dalen eds.,
North-Holland Pub. Co. (1982), 397--411.

\bibitem{ri2} {\sc Richman F.} {\it  Church Thesis without tears.}
Journal of Symbolic Logic, {\bf 48} (3) (1983), 797--803.

\bibitem{Ri88}  {\sc Richman F.} {\it Non trivial uses of trivial
rings.} Proc. Amer. Math. Soc., {\bf 103} (1988),
1012--1014.

\bibitem{Ri90}  {\sc Richman F.} {\it Intuitionism as generalization.}
Philosophia Mathematica, {\bf 5} (1990), 124--128.

%:h2013 rajout
\bibitem{Ri98} F. Richman.   \textit{The regular element property}. 
Proc. Amer. Math. Soc. {\bf 126} 7 (1998), 2123--2129.

\bibitem{Roi}  {\sc Roitman M.}
{\it  On projective modules over polynomial rings.}
Journal of Algebra {\bf 58} (1979),  51--63.

\bibitem{Roi86}  {\sc Roitman M.} {\it  On stably extended projective modules over polynomial rings} Proc. Amer. Math. Soc. {\bf 97} (1986), 
585--589.

\bibitem{Rota} {\sc Rota Gian Carlo}
{\it  The many lives of lattice theory.}  Notices Amer. Math. Soc.
{\bf 44} 11 (1997), 1440--1445.


%: S



\bibitem{sand} {\sc Sander T.} {\it  Existence and uniqueness of the real closure of an ordered field without Zorn's Lemma.} J. Pure and Applied Algebra
{\bf 73} (1991), 165--180.


%\bibitem{sei70} {\sc Seidenberg A.} 
% {\it  Construction of the integral closure of a finite integral domain.}
% Rend. Sem. Mat. e Fis.  Milano {\bf 40} (1970), 100--120.

\bibitem{sei74b} {\sc Seidenberg A.} {\it  What is Noetherian?}
 Rend. Sem. Mat. e Fis.  Milano {\bf 44} (1974), 55--61.

\bibitem{sei84} {\sc Seidenberg A.} {\it  On the Lasker-Noether decomposition theorem.}
Amer. J. Math  {\bf 106} (1984), 611--638.


\bibitem{Serregaga} {\sc Serre J.-P.}
{\it  Géométrie algébrique et géométrie analytique.}
Ann. Inst. Fourier Grenoble  {\bf 6} (1955-1956), 1--42.


\bibitem{Serre} {\sc Serre J.-P.}
{\it   Modules projectifs et espaces fibrés \`a fibre vectorielle}.
Séminaire P. Dubreil, Année 1957/1958.

\bibitem{SV} {\sc Simis A., Vasconcelos W.} {\it
Projective modules over $R[X]$, $R$ a valuation ring, are free.} Notices. Amer. Math.
Soc. {\bf 18} (5), (1971).

\bibitem{SkoRec} {\sc Skolem T.}
{\it   A critical remark on foundational research.} Norske Vid. Selsk. Forh., Trondheim 
{\bf 28} (1955), 100--105.

\bibitem{SoM} {\sc Soicher L., McKay J.} {\it Computing Galois groups over the rationals.} 
J. Number Theory, {\bf 20}, 
(1985), 273--281.



\bibitem{Sta} {\sc Stauduhar R.}
{\em The determination of Galois groups.}
Math. Comp. {\bf 27},  (1973), 981--996.

\bibitem{steel} {\sc Steel A.} {\it A New Scheme for Computing with Algebraically Closed
Fields}. Lecture Notes In Computer Science \textbf{2369}.
Proceedings of the 5th International Symposium on Algorithmic
Number Theory, (2002), 491--505.


\bibitem{steel2} {\sc Steel A.} {\it Computing with algebraically closed fields}.
Journal of Symbolic Computation. \textbf{45}, 342--372, (2010). 

\bibitem{Sto} {\sc Stone M. H.}
{\em Topological representations of distributive
lattices and Brouwerian logics.}
Cas. Mat. Fys. {\bf 67},  (1937), 1--25.


\bibitem{Stor} {\sc Storch U.}
{\em Bemerkung zu einem Satz von M. Kneser.}
Arch. Math. {\bf 23},  (1972), 403--404.

\bibitem{Sus} {\sc Suslin A.} {\em Projective modules over
polynomial rings are free. (Russian)}.   Dokl. Akad. Nauk SSSR
{\bf 229} no. 5 (1976), 1063--1066.

\bibitem{Sus77} {\sc Suslin A.} {\em On the structure of the special linear group over polynomial rings. (Russian)}.   Izv. Akad. Nauk. SSSR Ser. Mat.
{\bf 41}  (1977), 235--252.
English translation Math. USSR Izvestija, {\bf 11}, no. 2, 221--238.

\bibitem{Sus77b} {\sc Suslin A.} {\em Stably Free Modules. (Russian)}.   
Mat. Sb. (N.S.) {\bf 102} (1977), 537--550. English translation: 
Math. USSR Sb. {\bf 31}, 479--491.


\bibitem{Swan62} {\sc Swan R.}
{\em Factorization of Polynomials over Finite Fields.}
Pacific Journal of Mathematics {\bf 12}, no. 3, (1962), 1099--1106.

\bibitem{Swan} {\sc Swan R.}
{\em The Number of Generators of a Module.}
Math. Z. {\bf 102} (1967), 318--322.

\bibitem{Swan80} {\sc Swan R.}
{\em On Seminormality.}
Journal of Algebra, {\bf 67} (1980), 210--229.


%: T - Z

\bibitem{Ten} {\sc Tennenbaum J. B.}
{\em A constructive version of Hilbert's basis theorem.}
 Dissertation, University of California San Diego, 
(1973).

\bibitem{Tra} {\sc Traverso C.}
{\em Seminormality and the Picard group}. 
Ann. Scuola Norm. Sup. Pisa, {\bf 24} (1970), 585--595.


%\bibitem{Tu1} {\sc Turnbull H.}
%{\em On differentiating a matrix}, Proc. Edinb. Math. Soc. {\bf 2} vol. 1, part 2, (1928), 111--128.

%\bibitem{Tu2} {\sc Turnbull H.}
%{\em The theory of determinants, matrices and invariants}, Dover (1960).

\bibitem{Val} {\sc Valibouze A.}
\emph{Sur le corps des racines d'un polyn\^ome}.
Acta Arithmetica {\bf 131} (1), (2008), 1--27.

\bibitem {Vaserstein1}{\sc Vaserstein L.N. ({\rm with} A.A. Suslin)}
{\em Serre's problem on projective modules over 
polynomial rings and algebraic K-theory}.
Funk. An. {\bf 8} (1974), 65--66 = 
Funct. Anal. Appl. {\bf 8}, 148--150. 

\bibitem {Vaserstein2} {\sc Vaserstein L.N.}
{\em Serre's problem on projective modules over polynomial rings
after Suslin and Quillen}. 
(1976), Unpublished notes.

\bibitem {Vaserstein3}{\sc Vaserstein L.N. ({\rm with} A.A. Suslin)}
{\em Serre's problem on projective modules over 
polynomial rings and algebraic K-theory}.
Izv. Akad. Nauk SSSR Ser. Mat. {\bf 40} (1976), 993--1054 = 
Math. USSR Izv. {\bf 10}, 937--1001. 

\bibitem{Vessiot} {\sc Vessiot E.}
\emph{Sur la théorie de Galois et ses diverses généralisations}.
Ann. Sci. E. N. S. 3\`eme série {\bf 21}, (1904), 9--85.

\bibitem{vW} {\sc Van der Waerden.}
% {\em }
Review Zentralblatt f\"ur Math {\bf 24}, (1941), 276.

\bibitem{W1}  {\sc Weyl H.}
{\it  Das Kontinuum, Kritische Untersuchungen \"uber die Grundlagen der Analysis.}
Veit, Leipzig (1918). Traduction italienne
{\it  Il Continuo. Indagine critiche sui fondamenti dell' Analisi.} par A. B. Veit
Riccioli, Bibliopolis, Naples (1977).
Traduction anglaise  {\it  The Continuum. A critical examination of the foundations
of Ana\-lysis.}  par  S. Polard et T. Bole. Thomas Jefferson Press, University
Press of America (1987).
En fran{\c c}ais: {\it Le continu et autres écrits}. Traduits
et commentés par Jean Largeault. Librairie Vrin (1994).




%\bibitem{Y2} {\sc Yengui I.},
%\emph{}  

\bibitem{Ye0} {\sc Yengui I.}
\emph{An algorithm for the divisors of monic polynomials over a
commutative ring}.
Math. Nachr.  {\bf 260} (2003), 93--99.

%:HHH mise a jour Y3
\bibitem{Y3} {\sc Yengui I.}
\emph{Dynamical Gr\"obner bases}.
Journal of Algebra  {\bf 301} (2006), 447--458.
Corrigendum: \cite{Y3c}

%:HHH nouvelle ref Y3c
\bibitem{Y3c} {\sc Yengui I.}
 Corrigendum to  \emph{Dynamical Gr\"obner bases}  [J. Algebra 301 (2) (2006) 447--458] and to  \emph{Dynamical Gr\"obner bases over Dedekind rings} [J. Algebra 324 (1) (2010) 12--24].  
Journal of Algebra {\bf 339} (2011), 370--375.

\bibitem{Y1} {\sc Yengui I.}
\emph{Making the use of maximal ideals constructive}. 
Theoretical Computer Science, {\bf 392}, (2008) 174--178.


\bibitem{Ye3} {\sc Yengui I.} {\it The Hermite ring conjecture in dimension
one}. Journal of Algebra {\bf 320} (2008), 437--441. 


\bibitem{Ye4}  {\sc Yengui I.} {\it Stably free modules over $R[X]$ of rank $> \dim R$ are free}. Mathematics of Computation {\bf 80} (2011), 1093--1098.

\end{thebibliography}
\end{document}